\documentclass[11pt]{aart}
\usepackage{series-format}

\newif\ifwignerheadings

\usepackage{amsmath,amssymb,amsfonts,latexsym,amsthm,amsxtra,amscd}

\usepackage{mathrsfs,bm,mathtools}

\usepackage{tabularx,booktabs,array,enumitem,float,needspace}
\usepackage{tikz}
\usetikzlibrary{arrows.meta,calc,positioning,decorations.pathreplacing,shapes.geometric}

\tikzset{
 mg site/.style={circle,draw=black,fill=white,inner sep=0pt,minimum size=3.8pt},
 mg occurrence/.style={circle,draw=black,fill=black,inner sep=0pt,minimum size=3.4pt},
 mg block/.style={diamond,draw=black,fill=white,inner sep=0pt,minimum size=5.5pt},
 mg kernel/.style={draw=black,line width=.45pt},
 mg insertion/.style={draw=black,line width=1pt},
 mg incidence/.style={draw=black,densely dashed,line width=.45pt},
 mg arrow/.style={draw=black,line width=.45pt,-{Stealth[length=3.6pt,width=3pt]}},
 mg text/.style={font=\footnotesize,inner sep=1.5pt},
 mg label/.style={font=\footnotesize,fill=white,inner sep=1.5pt}
}

\usepackage{bbm,stmaryrd,longtable,multicol,textcomp,etoolbox}
\SetSymbolFont{stmry}{bold}{U}{stmry}{m}{n}
\usepackage[english]{babel}

\usepackage{hyperref}
\hypersetup{colorlinks,linkcolor={red!50!black},citecolor={blue!50!black},
 urlcolor={blue!80!black},pdftitle={Operator limits of random matrices: Edge universality and sharp approximations},
 pdfauthor={Hyung-Gil Kwak}}
\usepackage[capitalise,nameinlink]{cleveref}
\usepackage{bookmark}

\theoremstyle{plain}
\newtheorem{theorem}{Theorem}[subsection]
\newtheorem{proposition}[theorem]{Proposition}
\newtheorem{lemma}[theorem]{Lemma}
\newtheorem{corollary}[theorem]{Corollary}
\theoremstyle{definition}
\newtheorem{definition}[theorem]{Definition}
\newtheorem{assumption}[theorem]{Assumption}
\newtheorem{remark}[theorem]{Remark}
\AddToHook{env/lemma/begin}{\crefalias{theorem}{lemma}}
\AddToHook{env/proposition/begin}{\crefalias{theorem}{proposition}}
\AddToHook{env/assumption/begin}{\crefalias{theorem}{assumption}}
\AddToHook{env/remark/begin}{\crefalias{theorem}{remark}}
\AddToHook{env/definition/begin}{\crefalias{theorem}{definition}}
\AddToHook{env/corollary/begin}{\crefalias{theorem}{corollary}}
\numberwithin{equation}{subsection}
\numberwithin{figure}{section}

\newcommand{\localheading}[1]{\par\smallskip\Needspace{3\baselineskip}\noindent\textit{#1.}\enspace\ignorespaces}

\newcommand{\me}{\mathrm{e}}
\newcommand{\ii}{\mathrm{i}}
\newcommand{\dd}{\mathrm{d}}
\newcommand{\cum}{\operatorname{cum}}

\renewcommand{\le}{\leqslant}

\renewcommand{\ge}{\geqslant}
\renewcommand{\epsilon}{\varepsilon}

\newcommand{\s}{\mathsf}                        %

\newcommand{\q}[1]{[#1]}

\newcommand{\h}[1]{\{#1\}}

\newcommand{\norm}[1]{\lVert #1 \rVert}

\newcommand{\ang}[1]{\langle #1 \rangle}

\newcommand{\N}{\mathbb{N}}

\newcommand{\C}{\mathbb{C}}

\newcommand{\E}{\mathbb{E}}

\newcommand{\tr}{\operatorname{tr}}
\newcommand{\Tr}{\operatorname{Tr}}
\newcommand{\RNum}[1]{\uppercase\expandafter{\romannumeral #1\relax}}

\let\originalleft\left
\let\originalright\right
\renewcommand{\left}{\mathopen{}\mathclose\bgroup\originalleft}
\renewcommand{\right}{\aftergroup\egroup\originalright}

\renewcommand{\b}[1]{\boldsymbol{\mathrm{#1}}}
\renewcommand{\r}{\mathrm}

\newcommand{\PathCoeff}{\mathcal C_{\beta}}
\newcommand{\GeoCoeff}{C_{\beta}^{\mathrm e}}

\newtheorem*{theorem*}{Theorem}
\newtheorem*{lemma*}{Lemma}
\newtheorem*{proposition*}{Proposition}
\newtheorem*{corollary*}{Corollary}
\newtheorem*{definition*}{Definition}
\newtheorem*{assumption*}{Assumption}
\newtheorem*{remark*}{Remark}

\newtheorem*{conjecture*}{Conjecture}

\newtheorem*{example*}{Example}

\newtheorem*{question*}{Question}

\newcommand{\eqtag}[1]{\refstepcounter{equation}\tag{\theequation}\label{A:eq:#1}}

\newcommand{\HS}{\mathfrak S_2}
\newcommand{\HA}{\mathcal H^{\mathrm A}_\beta}
\newcommand{\HB}{\mathcal H^{\mathrm B}_{\beta,a}}
\newcommand{\KA}{\mathsf K^{\mathrm A}}
\newcommand{\KB}{\mathsf K^{\mathrm B}}
\newcommand{\KS}{\mathsf K^{\mathrm S}}
\newcommand{\RA}{\mathscr R^{\mathrm A}}
\newcommand{\RW}{\mathscr R^{\mathrm W}}
\newcommand{\Pp}{\mathbb P}

\definecolor{darkred}{rgb}{0.9,0,0.3}
\definecolor{darkblue}{rgb}{0,0.3,0.9}
\usetikzlibrary{decorations.markings,backgrounds}

\hypersetup{bookmarksnumbered=true,bookmarksdepth=1,pdfauthor={Hyung-Gil Kwak}}
\title{Operator limits of random matrices:\\Edge universality and sharp approximations}
\SeriesAuthor{Hyung-Gil Kwak}{hyunggilkwak581@gmail.com}

\hypersetup{bookmarksnumbered=true,bookmarksdepth=1}

\AddToHook{env/theorem/before}{\par\Needspace{5\baselineskip}}
\AddToHook{env/lemma/before}{\par\Needspace{5\baselineskip}}
\AddToHook{env/proposition/before}{\par\Needspace{5\baselineskip}}
\AddToHook{env/corollary/before}{\par\Needspace{5\baselineskip}}
\AddToHook{env/definition/before}{\par\Needspace{5\baselineskip}}
\AddToHook{env/assumption/before}{\par\Needspace{5\baselineskip}}
\AddToHook{env/remark/before}{\par\Needspace{5\baselineskip}}

\makeatletter
\newcommand{\manualstatementlabel}[2]{%
  \phantomsection
  \begingroup
  \def\@currentlabel{#2}%
  \label{#1}%
  \endgroup}
\makeatother

\begin{document}
\maketitle
\begin{abstract}
\fontsize{10pt}{11pt}\selectfont
We establish operator universality at the spectral edge of general Wigner matrices by showing that their suitably rescaled Lanczos tridiagonalizations, constructed from
deterministic initial vectors, converge to the stochastic Airy operator with Dirichlet boundary condition.

The result requires control of the entire tridiagonal structure on the scale relevant to the limiting differential operator. The Lanczos coefficients for non-Gaussian Wigner
matrices are strongly dependent. We analyze their fluctuations by combining multi-resolvent local laws with quantitative control of Krylov subspaces and a comparison of the
resulting Jacobi coefficients with those of the Gaussian ensembles. This approach identifies the limiting Brownian potential directly from the original matrix entries, without
invoking Dyson Brownian motion in the step of Green function comparison. The resolvent limit yields joint convergence of edge eigenvalues, embedded eigenfunctions, and spectral
measures.

For Gaussian, Laguerre and circular $\beta$-ensembles we construct sharp couplings to the stochastic Airy, Bessel and Sine operators, preserving the finite-matrix laws. For fixed
$\beta>0$, their squared Hilbert--Schmidt inverse errors have almost-sure orders $n^{-1/3}$, $n^{-1}$ and $n^{-1}$, respectively, with $a>0$ in the hard-edge upper bound. We also
identify the positive leading order coefficient for the circular ensemble, whose second derivative jumps at $\beta=2$.
\end{abstract}

\clearpage
\papercontents{11}{13.2}
\clearpage

\section{Introduction}\label{sec:introduction}

Universality in random matrix theory describes spectral statistics whose limits do not depend on the entry distributions. For Wigner matrices it is known for local eigenvalue
statistics and many eigenvector observables, with related results for covariance matrices and invariant ensembles; see \cite{erdos2017dynamical,benaych2018lectures} and the
references below. We study an operator version of this question: whether a rescaled matrix, expressed in suitable coordinates, converges to a universal random differential
operator.

 Wigner introduced random matrices as statistical models for complicated Hamiltonians. He also asked explicitly whether large matrix ensembles approximate operators on Hilbert
 space and how ensembles of such operators should be constructed \cite{wigner1955characteristic,wigner1967physics}. The connection between matrices and operators is also central to
 random band matrices, whose localization problem relates random matrix models to disordered Hamiltonians such as the Anderson model \cite{anderson1958absence,fyodorov1991scaling}.
 An operator limit specifies the limiting equation, its domain and boundary condition, and a topology for convergence on the comparison Hilbert space.

We compare the resolvents on a common Hilbert space and derive eigenvalue and eigenvector limits. The spectral scaling is essential: the global semicircle scale and the scale of
individual eigenvalues can lead to different operator limits. We treat the upper spectral edge of Wigner matrices, where the limiting operator is the stochastic Airy operator and
convergence is in the operator norm topology for the embedded resolvents.

Our main theorem derives this random differential operator from non-invariant Wigner ensembles, including Bernoulli matrices. Lanczos tridiagonalization from any deterministic unit
vector provides the comparison coordinates. After rescaling, the resolvents converge in distribution, locally uniformly in operator norm at nonreal parameters. The entry
assumptions are the Lee--Yin off-diagonal tail condition and uniform integrability of the standardized diagonal squares \cite{lee2014necessary}; an infinite fourth off-diagonal
moment is allowed. The main difficulty is that no explicit joint law is available for the induced recurrence coefficients. We obtain their limit by comparison through the
independent entries of the original matrix.

We also estimate the accuracy of approximation while preserving the finite dimensional matrix law. For the Gaussian, Laguerre and circular $\beta$-ensembles, we construct couplings
with the stochastic Airy, Bessel and Sine operators. Their squared Hilbert--Schmidt inverse errors have almost sure orders $n^{-1/3}$, $n^{-1}$ and $n^{-1}$, respectively, in the
parameter ranges stated below. The edge results extend to the specified polynomial potentials. We prove matching lower bounds under a rank constraint at the soft edge and under a
constraint on Brownian observations at the hard edge, with both approximation classes specified below. For the circular ensemble, we determine the leading order coefficient of the
error for the constructed coupling.

The Gaussian soft edge connects the two parts of the paper. The tridiagonal representation of Dumitriu and Edelman \cite{dumitriu2002matrix} has independent coefficients, whose
continuum limit is the stochastic Airy operator \cite{edelman2007random,ramirez2011beta}. At $\beta=1,2$ this is the distribution of the tridiagonal matrix obtained by applying the
Lanczos recurrence to the GOE and GUE. We use this operator limit as the reference in the comparison for Wigner matrices and first analyze the coupling error in the Gaussian case.
The qualitative and quantitative results share weighted resolvent estimates, while the proof for Wigner matrices also requires comparison of growing Gram matrices and localization
of the entire spectral subspace near the edge.

\subsection{Operator universality for Wigner matrices}
\label{sec:intro-wigner-universality}

An operator limit requires coordinates as well as a spectral scaling. For a finite self-adjoint matrix $A_n$ and an isometry $J_n$ into a fixed Hilbert space, the objects to
compare are
\[
 J_n(A_n-z)^{-1}J_n^*
 \quad\text{and}\quad
 (\mathcal H-z)^{-1}.
\]
The norm of their difference controls all unit vectors in this comparison space. Spectral projections on isolated intervals can then be recovered by contour integration. The choice
of $J_n$ matters: information about eigenvectors in these coordinates need not give their projections onto fixed vectors in the original matrix basis.

Let $H_N=(h_{ij})$ be real symmetric or complex Hermitian, with $\beta=1$ or $2$, respectively, and independent upper-triangular entries. The off-diagonal variables
$\xi_{ij}=\sqrt N\,h_{ij}$ have a common law, independent of $N$, satisfying
\[
 \E\xi_{12}=0,\qquad
 \E|\xi_{12}|^2=1,\qquad
 \E\xi_{12}^2=2/\beta-1,
 \qquad
 \lim_{s\to\infty}s^4\Pp(|\xi_{12}|>s)=0.
\]
The real diagonal variables $\sqrt N\,h_{ii}$ are centered, have uniformly bounded second moments, and have uniformly integrable squares. The tail condition is necessary and
sufficient for edge eigenvalue universality in the setting of Theorem~1.2 of Lee and Yin~\cite{lee2014necessary}. We use it as a sufficient hypothesis for the operator limit; we do
not prove its necessity in this topology. Assumption~\ref{W:ass:wigner-lee-yin} gives the complete entry conditions.

Fix a deterministic sequence of unit vectors $w_N$. Orthogonalizing $w_N,H_Nw_N,H_N^2w_N,\ldots$ gives the Lanczos basis up to the first invariant subspace. At a breakdown, the
construction continues on its orthogonal complement using the first nonzero projected coordinate vector. Repeating this rule gives a full tridiagonal representation $T_N(w_N)$,
including invariant subspaces invisible to the spectral measure at $w_N$. Set
\[
 A_N=N^{2/3}(2I-T_N(w_N)),\qquad h_N=N^{-1/3},
\]
and embed the $k$th coordinate as $\mathcal E_Ne_k=h_N^{-1/2}\mathbf1_{[(k-1)h_N,kh_N)}$ in $L^2(\mathbb R_+)$. The limiting operator is
\begin{equation}\label{eq:unified-airy}
 \HA=-\frac{\dd^2}{\dd x^2}+x+\frac2{\sqrt\beta}B'(x),
 \qquad x\ge0,
\end{equation}
where $B$ is standard Brownian motion. The white noise potential is interpreted through its quadratic form, with the Dirichlet condition at zero. Theorem~\ref{W:thm:main-operator}
gives
\begin{equation}\label{eq:unified-wigner}
 \RW_N(\cdot):=\mathcal E_N(A_N-\cdot)^{-1}\mathcal E_N^*
 \ \Longrightarrow\ (\HA-\cdot)^{-1}
\end{equation}
in $C_{\mathrm{loc}}(\mathbb C\setminus\mathbb R; \mathcal K(L^2(\mathbb R_+)))$, with uniform operator norm on each compact spectral set. The convergence holds for every
deterministic sequence of unit vectors; equivalently, the distance between the laws tends to zero uniformly over deterministic starting vectors. The vector is fixed before the
matrix is sampled. For general Wigner matrices the theorem gives convergence in distribution.

This representation also identifies the orthogonal polynomials involved in the limit. If $H_N\xi_j=\mu_j\xi_j$ with an orthonormal eigenbasis, the spectral measure at $w_N$ is
\[
 \mu_{H_N,w_N}=\sum_{j=1}^N|w_N^*\xi_j|^2\delta_{\mu_j}.
\]
Before breakdown, the Lanczos coefficients are the three term recurrence coefficients of the orthonormal polynomials for this measure. Consequently \eqref{eq:unified-wigner}
describes the scaled recurrence at degrees of order $N^{1/3}$ by a differential operator. The continued Lanczos matrix also contains the rest of the spectrum. This is the operator
interpretation discussed in~\cite{virag2014operator}, here applied to a random spectral measure arising from independent matrix entries.

Under their variance-profile assumptions and uniform subexponential decay, Erd\H{o}s, Yau and Yin proved edge eigenvalue universality when the first two entry moments match; see
Theorem~2.4 of~\cite{erdos2012rigidity}. Under analogous variance-profile and tail assumptions, Knowles and Yin proved joint edge eigenvalue and eigenvector universality in
Theorem~1.6 of~\cite{knowles2013eigenvectors}, assuming matching first two moments and level repulsion for one of the two ensembles. For the present operator limit, a bounded
interval contains order $N^{1/3}$ steps of orthogonalization. Entry differentiation therefore acts on the inverse of a Gram matrix whose dimension grows on this scale. The
derivative estimates must retain this dimension dependence.

The cumulative recurrence gives the required derivative estimate: differentiating its sum leaves a boundary square and differences of adjacent off-diagonal coefficients. This
cancellation would be lost in separate bounds on the individual coefficient derivatives. The total variation of those coefficients replaces a sum of order $N^{1/3}$ unrelated
terms. At derivative order $r$, the resulting factor is $\rho^{r-1}$ for a Gram matrix of dimension $\rho$; the calculation in Section~\ref{W:sec:proof-outline} explains why this
is sufficient for the remainders of the expansions at $\rho\asymp N^{1/3}$.

The third moments require a separate cancellation. The derivative coefficient and its resolvent factors are evaluated at the same matrix before summing over entries. The common
entry law then gives a deterministic matrix of trace norm $O(N\log(2N))$. Together with the entry normalization and the polynomial degree, this gives
$N^{-3/2}N^{1/3}N^{1+o(1)}=N^{-1/6+o(1)}$; see Section~\ref{W:sec:intro-third-moment-comparison}. This bounds the error for the expectations used in the comparison. The operator
theorem gives convergence without a rate.

Local coefficient convergence identifies the differential expression. To obtain compactness in operator norm, we also prove
\begin{equation}\label{eq:unified-capture}
 \lim_{L\to\infty}\limsup_{N\to\infty}\sup_{\|w\|=1}
 \mathbb P\!\left(\|Q_{N,L}\Pi_N(E)\|>\eta\right)=0,
 \qquad E\in\mathbb R,\quad\eta>0.
\end{equation}
Here $\Pi_N(E)$ is the spectral projection of the whole scaled Lanczos matrix onto $(-\infty,E]$, and $Q_{N,L}$ projects onto coordinates beyond the first
$\lfloor LN^{1/3}\rfloor$; see \eqref{W:eq:operator-coefficient-tail-notation}. The estimate must detect a normalized eigenfunction lying entirely in a later invariant block, even
when its spectral weight at $w_N$ is zero. We therefore compare an \emph{unnormalized} trace of the low energy tail in the tridiagonal coordinates. A separate cubic trace
calculation gives its comparison with the GOE or GUE. Theorem~\ref{W:thm:operator-probabilistic} combines this localization with local coefficient convergence. Its final
application is in Section~\ref{W:sec:main-assembly}.

The auxiliary estimates for bounded entries, Gaussian entries, and entries whose distributions are convex combinations of these laws use the isotropic local semicircle law in
Theorem~2.12 of~\cite{bloemendal2014isotropic} and the spectral deviation bounds of~\cite{erdosxu2023small}. The original truncated entries require separate conditional row
estimates. Products of resolvents arise from explicit entry differentiation; we prove their bounds in the comparison argument. Section~\ref{W:sec:intro-entry-comparison} specifies
the hypotheses and scales of the local laws used in these estimates.

\paragraph{Consequences in the original matrix basis}
The operator construction also has conclusions that can be stated without Lanczos coordinates. For each deterministic $a\in\mathbb R$,
Corollaries~\ref{W:cor:real-energy-inverse-convergence} and~\ref{W:cor:real-spectral-distance} imply
\[
 \begin{split}
 N^{-2/3}\bigl\|((2-aN^{-2/3})I-H_N)^{-1}\bigr\|
 &\Longrightarrow\| (\HA-a)^{-1}\|,\\
 N^{2/3}\operatorname{dist}(2-aN^{-2/3},\sigma(H_N))
 &\Longrightarrow\operatorname{dist}(a,\sigma(\HA)).
 \end{split}
\]
The inverse is assigned value zero on the event that the finite matrix is singular; the probability of this event tends to zero. The limiting spectral distance is almost surely
strictly positive. Thus these statements control the size of the inverse at a prescribed real energy, using the whole finite spectrum.

For $F\in C_c(\mathbb R)$, the convergence of the spectral measure in Corollary~\ref{W:cor:ev-starting-vectored-measure} yields
\[
 Nw_N^*F\bigl(N^{2/3}(2I-H_N)\bigr)w_N
 \Longrightarrow\sum_{j\ge1}F(\Lambda_j)|\varphi_j'(0)|^2.
\]
Here $\Lambda_j$ and $\varphi_j$ are the ordered eigenvalues and $L^2$-normalized eigenfunctions of $\HA$. Only finitely many terms on the right are nonzero. Together with the
joint eigenvalue and projection limits, this describes spectral observables at the prescribed vector in the original matrix space. The displayed quantities depend on $H_N$ and, in
the last formula, $w_N$; their definitions require no tridiagonal representation.

The quantitative results also have consequences for entire spectral sequences. Theorem~\ref{A:res:20.6} identifies the order $n^{-1/3}$ of the squared $\ell^2$ distance between the
complete reciprocal spectra in the coupling of the tridiagonal matrices to the limiting operator, with the stated matching convention. This distance is defined by the spectra
themselves. It controls the full reciprocal spectrum.

\subsection{Sharp approximation of the Airy, Bessel and Sine operators}
\label{sec:intro-sharp-approximations}

The classical $\beta$-ensembles admit one-dimensional representations with explicit coefficient laws. At the soft edge, Ram\'irez, Rider and Vir\'ag~\cite{ramirez2011beta} proved
the limit given by the stochastic Airy operator, and Krishnapur, Rider and Vir\'ag~\cite{krishnapur2016universality} extended it to Jacobi matrices with uniformly convex polynomial
potentials. At the hard edge, the inverse kernel construction of Ram\'irez and Rider~\cite{ramirez2009diffusion} gives the limit given by the stochastic Bessel operator. Rider and
Waters~\cite{rider2019universality} established operator convergence for their polynomial potential models with $\beta\ge1$. For the circular ensemble, Valk\'o and
Vir\'ag~\cite{valko2017sine,valko2020operator} constructed the limiting Sine$_\beta$ operator and a quantitative coupling of its inverse. We estimate the approximation errors in
these operator representations.

Fix the parameters of one model. We construct its finite coefficient arrays at all sizes, with their exact laws, on a probability space carrying the Brownian path of the limiting
operator. The arrays may depend on one another across sizes. We write $\|\cdot\|_{\HS}$ for the unnormalized Hilbert--Schmidt norm. The following three statements concern its
square.

At the soft edge, let $\KA_n$ be the inverse of the scaled Gaussian or general potential Jacobi matrix, expressed on the comparison space, and put $\KA=(\HA)^{-1}$. For every fixed
$\beta>0$, and for each polynomial potential with $V''\ge c_V>0$ and the stated regular soft edge,
\begin{equation}\label{eq:unified-airy-rate}
 \|\KA_n-\KA\|_{\HS}^2=\Theta_\omega(n^{-1/3}).
\end{equation}
Here $u_n=\Theta_\omega(v_n)$ means that, on an event of probability one, there are $0<c_\omega\le C_\omega<\infty$ and $n_0(\omega)<\infty$ such that
$c_\omega v_n\le u_n\le C_\omega v_n$ for every $n\ge n_0(\omega)$. On this event, zero belongs to the resolvent set of the limit and, after increasing $n_0(\omega)$ if necessary,
each scaled finite matrix is invertible. The Gaussian map is a full isometry. For the general potential the map $\mathcal J_n$ of \eqref{A:eq:1.1} satisfies
$\mathcal J_n^*\mathcal J_n=P_{I_n}^*P_{I_n}$, the projection onto the retained coordinates. The full matrix is inverted before this compression is applied. On the same event, the
two-sided estimates hold uniformly on each nonempty compact subset of $\mathbb C\setminus\mathbb R$, with constants $c_{\omega,K},C_{\omega,K}$ and a starting index $n_0(\omega,K)$
(Theorems~\ref{A:res:1.1}--\ref{A:res:1.2}).

At the hard edge, let $\KB=(\HB)^{-1}$ be the positive inverse of the stochastic Bessel operator on $L^2(0,1)$, and let $\KB_n$ be the isometrically embedded inverse of the scaled
beta-Laguerre or polynomial potential bidiagonal matrix. For fixed $\beta>0$ and $a>0$,
\begin{equation}\label{eq:unified-bessel-rate}
 0<\liminf_{n\to\infty}n\|\KB_n-\KB\|_{\HS}^2
 \le\sup_{n\ge1}n\|\KB_n-\KB\|_{\HS}^2<\infty
 \quad\text{almost surely}.
\end{equation}
In the polynomial case we require $x\mapsto V(x^2)$ to be uniformly convex, with the normalization and coordinate change of
Section~\ref{sec:hard-edge-models-brownian-information-and-bessel-limits}. The same two-sided order holds locally uniformly on compact subsets of the limiting resolvent set. At
size $n$, the coefficient array uses only that size's Brownian samples at deterministic logarithmic times and auxiliary randomness independent of the entire path. The limiting
operator exists for $a>-1$; the coupled upper estimate uses $a>0$. The range $0<\beta<1$ requires separate estimates for the singular initial coefficient densities.

For the circular ensemble, let $\KS_n$ and $\KS$ be the inverses of $\mathrm{Circ}_{\beta,n}$ and $\mathrm{Sine}_\beta$, conjugated to the same unweighted $L^2$ space. For each
fixed $\beta>0$,
\begin{equation}\label{eq:unified-sine-rate}
 n\|\KS_n-\KS\|_{\HS}^2
 \longrightarrow \mathcal C_\beta(Z)\in(0,\infty)
 \qquad\text{almost surely},
\end{equation}
where $\mathcal C_\beta(Z)$ is an explicit functional of the limiting hyperbolic Brownian path. The terminal boundary point is the same in both operators. The circular differential
operator represents the periodic extension of $n$ eigenangles and has an inverse of infinite rank. Theorem~1 of Valk\'o and Vir\'ag~\cite{valko2020operator} bounds the squared
Hilbert--Schmidt error by $C_\omega(\log n)^6/n$ for every $n\ge n_0(\omega)$, on an event of probability one. Our construction uses a different coupling of the same operator laws
and identifies the positive leading order coefficient in \eqref{eq:unified-sine-rate}. The coefficient is specific to this coupling; we do not prove a lower bound over all
couplings.

Thus the error norms are of order $n^{-1/6}$ at the soft edge and $n^{-1/2}$ at the hard edge and for the circular ensemble. The next two subsections explain the coefficient
estimates needed for the upper bounds and the distinct reasons for these error scales.

\subsection{Coefficient estimates and operator errors}
\label{sec:intro-coefficient-comparisons}

The inverse formula determines how accurately the coefficients must be coupled. We estimate this error separately from the discretization of the limiting path. An operator built
from samples of that path, with the same normalization and boundary data, provides the intermediate comparison. Denoting its inverse or resolvent by $\mathsf R_n^B$, and the finite
and limiting objects by $\mathsf R_n$ and $\mathsf R$, gives
\begin{equation}\label{eq:unified-decomposition}
 \mathsf R_n-\mathsf R
 =\bigl(\mathsf R_n-\mathsf R_n^{B}\bigr)
  +\bigl(\mathsf R_n^{B}-\mathsf R\bigr).
\end{equation}
The first difference measures the coefficient coupling error while preserving the finite ensemble law. The second measures the error from discretizing the Brownian path. Each
comparison includes its interval restrictions and boundary corrections.

For the tridiagonal matrices at the soft edge, the resolvent identity places a diagonal error between two resolvents. A coefficient difference that is a discrete derivative has a
better estimate after summation by parts. Section~\ref{A:sec:earlier-results} displays the two bounds: estimating the derivative coordinate by coordinate can give an error of order
$h^{-1}$, while estimating its primitive gives order $h$. This determines how the conditional Gaussian field is coupled to Brownian increments. Its partial sum variance fixes the
coefficient multiplying the Brownian increments, and a factorization of its covariance separates those increments from a discrete derivative. The remaining errors satisfy the
weighted square sum bound for diagonal perturbations. The coefficient estimates establish the quadratic form bounds needed for the resolvent comparison before that comparison is
applied.

The finite difference approximation driven by the Brownian increments must also be compared with the differential operator on the half-line. With lattice spacing $h\asymp n^{-1/3}$
and retained length $L\asymp h^{-2}$, the proof gives squared error bounded by $C_\omega(h+L^{-1/2})=O_\omega(h)$. The estimate separates the discretization error $h$ from the
truncation error $L^{-1/2}$. This is the content of Theorem~\ref{A:res:1.3}; Section~\ref{A:sec:outline-of-the-proof} explains how it is combined with the coefficient comparison
and the boundary estimates.

The finite matrix inverse at the hard edge contains products of bidiagonal coefficient ratios over intervals of order $n$. Their logarithms are partial sums. The weighted estimate
\eqref{B:eq:1-55} therefore requires bounds on the cumulative logarithmic errors and on the local coefficient errors, together with a maximum bound that controls the exponential
factors during interpolation. Bounds on individual ratios alone do not control their products over a growing interval. For dependent coefficients, the Gaussian approximation and
the diffusion comparison are estimated using weighted sums of squared errors in the partial sums. The remaining coordinates are then sampled from their joint conditional law. The
integral identity in Section~\ref{B:subsection:intro-green-kernel-comparison} verifies that this recovers the original coefficient distribution while retaining the Gaussian vector
already coupled to the Brownian samples.

For the circular ensemble, the Green kernel is bilinear in its two boundary columns. In the square of its first variation, integration in the unchanged variable supplies the weight
on the perturbed column. The estimate therefore remains useful when one boundary solution grows near the singular boundary. Applying it to products of the coupled transfer matrices
gives $\|\KS_n-\mathsf K_n^{\mathrm S,B}\|_{\HS}=o_\omega(n^{-1/2})$, where $\mathsf K_n^{\mathrm S,B}$ is formed from the Brownian samples. The product calculation estimates the
centered increment errors and their means separately, using the same terminal boundary line. This leaves the sampling error as the source of the leading order coefficient in
\eqref{eq:unified-sine-rate}.

\subsection{Lower bounds and the leading order coefficient}
\label{sec:intro-error-scales}

At the soft edge the reciprocal spectrum already forces the rate. The eigenvalues of the stochastic Airy operator satisfy $\Lambda_j\asymp_\omega j^{2/3}$; see
Lemma~\ref{A:res:weyl} and~\cite{virag2014operator} (Theorem~6.1). For $\RA(z)=(\HA-z)^{-1}$ and nonreal $z$, the singular value tail gives
\begin{equation}\label{eq:unified-rank}
 \inf_{\operatorname{rank}F\le n}
 \|F-\RA(z)\|_{\HS}^2
 =\sum_{j>n}s_j(\RA(z))^2
 \asymp_{\omega,z}n^{-1/3}.
\end{equation}
This minimum is attained by spectral truncation, without a restriction on the approximating law or range. The coupled matrices attain the same order while keeping their exact laws
and identifying maps. The coefficient of the unrestricted minimum is computed in Section~\ref{A:sec:rank-lower-bound-and-optimality}. A leading order coefficient for the matrix
coupling is not identified.

The inverse of the stochastic Bessel operator is trace class. Its rank-$n$ spectral truncations therefore have squared Hilbert--Schmidt error $o(n^{-1})$
(Proposition~\ref{B:prop:bessel-best-finite-rank-approximation}). For the coupling in \eqref{eq:unified-bessel-rate}, the restriction comes instead from the available Brownian
observations. Conditional on the values at consecutive observation times, the intervening bridges remain random. After restriction to an interior rectangle and localization of the
sampled data, their linear contribution has conditional covariance with trace at least $c/n$ and operator norm at most $C/n^2$. The covariance trace bounds the mean squared
Gaussian fluctuation from below. The smaller operator norm prevents a few directions from accounting for that fluctuation. The resulting small ball estimate is uniform over the
translate chosen from the observations. Its failure probability and that of the nonlinear remainder are summable, so the lower bound holds almost surely.
Theorem~\ref{B:thm:brownian-information-lower-bound} applies for $a>-1$ to each fixed sequence of Hilbert--Schmidt approximations using only the samples at its own size and
auxiliary randomness independent of the whole path. Its probability one event may depend on the approximation sequence.

For the circular model we compute the error for the specified sampling times and coupling. This calculation does not minimize over other approximations. The first variation of the
Green kernel has a conditional Gaussian term and a mean measurable with respect to the sampled values. Both the covariance trace and the squared mean enter its conditional second
moment. The position of the sampling time within each operator interval also contributes to the coefficient. The samples and terminal point are fixed before applying the
concentration inequality for the squared norm of the conditional Gaussian term. A covariance operator norm smaller than $n^{-1}$, together with a conditional second moment of order
$n^{-1}$, gives summable localized tails. The coefficient calculation also estimates the path dependent weights and cross terms. The proof therefore requires both concentration and
a separate identification of the conditional mean.

To transfer the leading coefficient to the finite ensemble, set $D_n=\KS_n-\mathsf K_n^{\mathrm S,B}$ and $E_n=\mathsf K_n^{\mathrm S,B}-\KS$. The exact identity is
\[
 n\|D_n+E_n\|_{\HS}^2
 =n\|E_n\|_{\HS}^2
  +2n\operatorname{Re}\langle D_n,E_n\rangle_{\HS}
  +n\|D_n\|_{\HS}^2.
\]
Since $\|D_n\|_{\HS}=o_\omega(n^{-1/2})$ and $\|E_n\|_{\HS}=O_\omega(n^{-1/2})$, Cauchy--Schwarz gives
\[
 \bigl|n\|D_n+E_n\|_{\HS}^2-n\|E_n\|_{\HS}^2\bigr|
 \le 2n\|D_n\|_{\HS}\|E_n\|_{\HS}+n\|D_n\|_{\HS}^2=o_\omega(1).
\]
An $O(n^{-1/2})$ bound on $D_n$ alone would leave both contributions at the scale of the leading order coefficient.

\subsection{Spectral consequences and organization}
\label{sec:intro-spectral-consequences}

The operator limits give convergence of isolated spectral projections in the comparison coordinates. Ordered edge eigenvalues also require the lower spectral bounds established in
the corresponding arguments. For the full matrix with a general potential under the soft edge scaling, projections are compared using the enlarged isometry $\widehat{\mathcal J}_n$
and the inclusion $f\mapsto(f,0)$ of Corollary~\ref{A:res:full-ordered}; compression by $\mathcal J_n$ alone need not produce an orthogonal projection.

For Wigner matrices, convergence of the boundary derivative identifies the spectral weights of finitely many edge eigenvectors at one deterministic direction. After independent
signs or phases are assigned, these give the Gaussian projection laws. Applying the limit for one deterministic direction to deterministic linear combinations gives the joint
projections of one fixed eigenvector onto finitely many deterministic directions with convergent Gram matrices (Corollaries~\ref{W:cor:intro-boundary-weights},
\ref{W:cor:ev-one-direction-gaussian} and~\ref{W:cor:ev-many-directions}). The quantitative inverse estimates also give global bounds for reciprocal eigenvalues and fixed interval
spectral estimates, with their stated gap conditions.

Section~\ref{sec:random-matrix-models-and-main-results} gives the complete models, normalizations and main statements. The \hyperlink{proof-strategy}{proof strategy} following that
section separates the argument for Wigner matrices from the quantitative coupling arguments. The common resolvent identity, the tridiagonal convergence criterion and the model
specific weighted estimates are in Section~\ref{sec:resolvent-estimates-and-operator-convergence}. The quantitative arguments continue through the Brownian approximations of
Section~\ref{sec:discrete-approximations-of-the-limiting-operators}, the coefficient couplings of Section~\ref{sec:coupling-matrix-coefficients-to-brownian-motion}, and the rate
and coefficient calculations of Section~\ref{sec:sharp-convergence-rates-and-the-limiting-error-coefficient}. Section~\ref{sec:the-circular-comparison-scheme} reduces the circular
estimate to the coefficient comparison and the error from sampling the Brownian path, and Section~\ref{S:sec:proof-of-main-theorems} combines them to prove the circular asymptotic.

The estimates for the GOE and GUE used in the entry comparison are collected in Section~\ref{sec:finite-gaussian-estimates-for-the-wigner-comparison}.
Section~\ref{sec:green-function-comparison-for-wigner-matrices} proves the local coefficient limit. Section~\ref{sec:spectral-localization-and-the-wigner-operator-limit} proves
localization and concludes with the operator theorem in Section~\ref{W:sec:main-assembly}. Section~\ref{sec:spectral-consequences-and-dependence-on-parameters} contains the
spectral consequences and parameter results.

\clearpage
\section{Random matrix models and main results}\label{sec:random-matrix-models-and-main-results}

The four models share an operator formulation, but their coordinates and normalizations are different. These choices determine both the limiting operator and the meaning of the
error, so each model is introduced together with its identifying map. The symbols $A_n$, $J_n$ and $\mathcal H$ are therefore local to the model under discussion; the superscripts
$\mathrm A$, $\mathrm B$ and $\mathrm S$ distinguish them in the introduction. Both soft edge representations have the Dirichlet boundary condition at zero.

The operator comparisons use $(A-z)^{-1}$. In expansions in the Wigner entries we instead write $R_H(z)=(zI-H)^{-1}$, with the sign convention specified in
Section~\ref{W:sec:notation-and-scales}. The comparisons at the soft and hard edges use scalar products linear in the first variable. The complex canonical systems use
$\langle f,g\rangle_R=\int f^*Rg$, linear in the second variable; their probabilistic estimates are on the real Hilbert space. The inverse kernel at the hard edge acts by
integration in its first variable, $(Tf)(t)=\int K_T(s,t)f(s)\,\dd s$; the canonical kernel acts by integration in its second variable, $(Kf)(s)=\int K(s,t)f(t)\,\dd t$.

All Hilbert--Schmidt norms in the approximation theorems are unnormalized. For entry comparison, $\|A\|_{2,N}=N^{-1/2}\|A\|_{\mathrm{HS}}$ and
$\langle A\rangle=N^{-1}\operatorname{Tr}A$ have the normalization stated in Section~\ref{W:sec:notation-and-scales}. A bare norm has the meaning fixed in the current argument; in
particular, the error estimates for the circular ensemble use the norm of $\mathscr K=\mathfrak S_2(\mathcal H)$. For the forms of the scalar differential and difference operators,
inner products are linear in the first argument. The treatment of the canonical systems uses $\langle f,g\rangle_R=\int f^*Rg$, linear in the second argument, as specified with its
canonical system; its Gaussian field spaces are real. In each part a quadratic form is the diagonal of the specified sesquilinear or real symmetric bilinear form, with its domain
stated where it is introduced. Subscripts on constants record allowed parameter dependence, and $\omega$ allows dependence on the sampled outcome.

The operator limit for Wigner matrices is stated in Theorem~\ref{W:thm:main-operator}. The soft edge inverse and resolvent rates are Theorems~\ref{A:res:1.1} and~\ref{A:res:1.2};
the hard edge results are Theorems~\ref{B:result:1-1} and~\ref{B:result:1-3}; the circular asymptotic and the corresponding approximation by a sampled path are
Theorems~\ref{S:thm:main} and~\ref{S:thm:reference-main}. The parameter ranges and probability quantifiers are specified in each statement. The circular theorem gives a pathwise
asymptotic for its coupling; the theorem for Wigner matrices gives convergence in distribution in the operator topology.

\subsection{Wigner matrices and their operator limit}
\label{sec:wigner-matrices-and-the-airy-limit-for-a-deterministic-starting-vector}
\label{W:sec:introduction}

We construct the comparison coordinates by Lanczos iteration from a deterministic unit vector. At a breakdown, we continue in the orthogonal complement of the invariant subspace,
obtaining a tridiagonal representation of the whole matrix. We state the entry assumptions, define the embedding and prove convergence at the soft edge.

\subsubsection{Wigner matrices and their Lanczos representation}
\label{W:sec:wigner-model-and-main-result}

The entry assumptions allow distributions for which fourth moments need not exist, while keeping the variance normalization that places the upper spectral edge at two. Throughout
the treatment of Wigner matrices, $\beta\in\{1,2\}$, $\mathbb F_1=\mathbb R$ and $\mathbb F_2=\mathbb C$. The matrix $H_N=(h_{ij})$ is real symmetric for $\beta=1$ and complex
Hermitian for $\beta=2$. Its standardized entries are written as $X_N=\sqrt N H_N=(x_{ij})$, with $\xi_{ij}=x_{ij}$ for $i<j$ and $d_{N,i}=x_{ii}$. The precise assumptions
distinguish the off-diagonal law from the diagonal laws.

\begin{assumption}
\label{W:ass:wigner-lee-yin}
The upper-triangular entries of $H_N$ are independent. The variables $\xi_{ij}$, $i<j$, have a common law independent of $N$, with
\begin{equation}
 \mathbb E\xi_{12}=0,\qquad \mathbb E|\xi_{12}|^2=1,\qquad
 \mathbb E\xi_{12}^2=2/\beta-1.
 \label{W:eq:off-diagonal-normalization}
\end{equation}
The real variables $d_{N,i}$ are centered and satisfy
\begin{equation}
 \sup_{N,i}\mathbb E d_{N,i}^2<\infty,
 \label{W:eq:diagonal-second-moment}
\end{equation}
\begin{equation}
 \lim_{s\to\infty}\sup_{N,i}
 \mathbb E[d_{N,i}^2\mathbf1_{\{|d_{N,i}|>s\}}]=0.
 \label{W:eq:diagonal-uniform-integrability}
\end{equation}
The off-diagonal tail satisfies
\begin{equation}
 \lim_{s\to\infty}s^4\mathbb P(|\xi_{12}|>s)=0.
 \label{W:eq:lee-yin-tail}
\end{equation}
\end{assumption}

The diagonal laws may depend on $N$ and $i$. The off-diagonal condition \eqref{W:eq:lee-yin-tail} is the criterion of Lee and Yin \cite{lee2014necessary}; it implies finiteness of
every moment of order less than four. An infinite fourth moment is permitted. The tail estimates needed for truncation are proved in Lemma~\ref{W:lem:offdiag-tail-moments}. We keep
the diagonal assumptions \eqref{W:eq:diagonal-second-moment}--\eqref{W:eq:diagonal-uniform-integrability} separate throughout the comparison.

Fix a deterministic starting sequence
\begin{equation}
 w_N\in\mathbb F_\beta^N,\qquad \|w_N\|=1.
 \label{W:eq:deterministic-starting-vector}
\end{equation}
Both coordinate vectors and vectors with entries of order $N^{-1/2}$ are included. The vector is fixed before $H_N$ is sampled. Uniform probability estimates below have the form
$\sup_{\|w\|=1}\mathbb P(\mathcal B_N(w))$, with $w$ deterministic inside each probability. The equivalence with estimates for every deterministic starting sequence is proved in
Section~\ref{W:sec:measurability-and-truncation}.

For a Hermitian $H$ and a unit vector $w$, the Krylov subspaces are
\begin{equation}
 \mathcal K_r(H,w)=\operatorname{span}\{w,Hw,\ldots,H^{r-1}w\},
 \qquad \mathcal K_0(H,w)=\{0\}.
 \label{W:eq:lanczos-krylov-subspace}
\end{equation}
The vector $w$ is cyclic when $\mathcal K_N(H,w)=\mathbb F_\beta^N$. The recurrence below is Lanczos tridiagonalization \cite{lanczos1950iteration} (\S VII, equations~(63)--(75)).
We specify its continuation when the Krylov subspace becomes invariant. Put $q_0=0$, $q_1=w$, and $b_0=0$. Having constructed $q_1,\ldots,q_k$, set
\begin{equation}
 a_k=q_k^*Hq_k,\qquad
 r_k=Hq_k-a_kq_k-b_{k-1}q_{k-1},\qquad
 \widehat P_k=\sum_{j=1}^kq_jq_j^*.
 \label{W:eq:lanczos-recurrence}
\end{equation}
If $k<N$ and $r_k\ne0$, put $b_k=\|r_k\|$ and $q_{k+1}=r_k/b_k$. If $k<N$ and $r_k=0$, put $b_k=0$ and
\begin{equation}
 j_k=\min\{j:\|(I-\widehat P_k)e_j\|>0\},\qquad
 q_{k+1}=\frac{(I-\widehat P_k)e_{j_k}}
                    {\|(I-\widehat P_k)e_{j_k}\|}.
 \label{W:eq:intro-lanczos-continuation}
\end{equation}
The equality $r_k=0$ is called a breakdown. The preceding span is then invariant, and \eqref{W:eq:intro-lanczos-continuation} starts a new block in its orthogonal complement. The
construction gives an orthonormal basis $U=(q_1,\ldots,q_N)$ and a real symmetric tridiagonal matrix with nonnegative off-diagonal entries; the proof is in
Proposition~\ref{W:prop:lanczos-construction}.

For $H=H_N$ write $U_N=U_N(w_N)$ and define the Lanczos tridiagonal matrix by
\begin{equation}
 T_N(w_N)=U_N(w_N)^*H_NU_N(w_N).
 \label{W:eq:lanczos-tridiagonal-matrix}
\end{equation}
Before the first breakdown, the first $r$ vectors span $\mathcal K_r(H_N,w_N)$. Afterwards the Krylov subspace stays fixed, whereas the continued Lanczos basis includes the
remaining invariant blocks. Lemma~\ref{W:lem:lanczos-first-block} establishes these properties of the Lanczos construction. The continuation is Borel by
Section~\ref{W:sec:measurability-and-truncation}, so $T_N$ is a random matrix for every entry law in the assumption.

Unitary conjugation gives
\begin{equation}
 \sigma(T_N)=\sigma(H_N)
 \label{W:eq:spectrum-preserved}
\end{equation}
with multiplicity. For the decreasing eigenvalues $\mu_1^{(N)}\ge\cdots\ge\mu_N^{(N)}$ of $H_N$, set
\begin{equation}
 \mathcal A_N=N^{2/3}(2I-T_N).
 \label{W:eq:scaled-lanczos-operator}
\end{equation}
Its increasing eigenvalues are
\begin{equation}
 \lambda_{j,N}=N^{2/3}(2-\mu_j^{(N)}).
 \label{W:eq:scaled-eigenvalue-identity}
\end{equation}

\subsubsection{The operator limit}
\label{W:sec:intro-limiting-operator}

The Lanczos representation supplies the coordinates; an embedding into a fixed Hilbert space is still needed to compare operators of different sizes. At the soft edge the spatial
scale is $h_N=m_N^{-1}$, where $m_N=N^{1/3}$. The corresponding embedding uses piecewise constant functions:
\begin{equation}
 (\mathcal E_Nv)(x)=h_N^{-1/2}v_k,
 \qquad (k-1)h_N\le x<kh_N,\quad 1\le k\le N,
 \label{W:eq:piecewise-constant-embedding}
\end{equation}
with value zero for $x\ge Nh_N$. For $v\in\mathbb F_\beta^N$ and $f\in L^2(\mathbb R_+)$,
\[
 \|\mathcal E_Nv\|_2^2=\sum_{k=1}^N|v_k|^2,\qquad
 (\mathcal E_N^*f)_k=h_N^{-1/2}
             \int_{(k-1)h_N}^{kh_N}f(x)\,\dd x.
\]
Thus $\mathcal E_N^*\mathcal E_N=I$, and the embedded resolvent is
\begin{equation}
 \mathscr R_N(z)=\mathcal E_N(\mathcal A_N-z)^{-1}\mathcal E_N^*,
 \qquad z\in\mathbb C\setminus\mathbb R.
 \label{W:eq:embedded-resolvent}
\end{equation}
We write $\mathscr R_N^{(w)}$ when the starting vector is displayed.

Let $B$ be standard Brownian motion on $[0,\infty)$. The stochastic Airy operator with Dirichlet boundary condition has the formal expression
\begin{equation}
 \mathcal H_\beta=-\frac{\dd^2}{\dd x^2}+x+\frac2{\sqrt\beta}B'(x).
 \label{W:eq:sao-formal}
\end{equation}
We use the inner product $\langle f,g\rangle=\int f\overline g$, linear in its first argument. A sesquilinear form is linear in its first argument and conjugate linear in its
second, and its quadratic form is the diagonal value $\mathfrak q[f]=\mathfrak q[f,f]$. The operator in \eqref{W:eq:sao-formal} is defined by closing the quadratic form
\begin{equation}
 \mathfrak q_\beta[f]
 =\int_0^\infty|f'|^2\,\dd x+\int_0^\infty x|f|^2\,\dd x
  -\frac2{\sqrt\beta}\int_0^\infty B(x)(|f(x)|^2)'\,\dd x,
 \qquad f\in C_c^\infty(0,\infty).
 \label{W:eq:sao-form}
\end{equation}
On this initial domain the associated Hermitian sesquilinear form is
\[
 \mathfrak q_\beta[f,g]
 =\int f'\overline{g'}+\int x f\overline g
      -\frac2{\sqrt\beta}\int B(f\overline g)'.
\]
Its closure has form domain $H_0^1(0,\infty)\cap L^2((0,\infty),x\,\dd x)$ almost surely. The representation theorem associates a unique self-adjoint operator with this closed
semibounded form; see Theorems~VI.2.1 and~VI.2.6 of \cite{kato1995perturbation}. Its spectrum is bounded below and its resolvent is compact \cite{ramirez2011beta}. These properties
are also proved from the pathwise form estimates in Lemmas~\ref{W:lem:operator-limit-operator} and~\ref{W:lem:operator-brownian-form}. Put
\begin{equation}
 \mathscr R_\beta(z)=(\mathcal H_\beta-z)^{-1}.
 \label{W:eq:sao-resolvent}
\end{equation}

Both resolvents take values in the compact operators on $\mathcal H=L^2(\mathbb R_+)$. We use
\begin{equation}
 \mathscr Y=C_{\rm loc}(\mathbb C\setminus\mathbb R;\mathcal K(\mathcal H))
 \label{W:eq:operator-path-space}
\end{equation}
with locally uniform operator norm convergence. For this topology use the compact exhaustion
\begin{equation}
 K_j=\{z:|z|\le j,\ |\operatorname{Im}z|\ge j^{-1}\},\qquad j\ge2,
 \label{W:eq:resolvent-compact-exhaustion}
\end{equation}
and the metric
\begin{equation}
 d_{\mathscr Y}(S,T)=\sum_{j=2}^\infty2^{-j}
       \left(1\wedge\sup_{z\in K_j}\|S(z)-T(z)\|\right).
 \label{W:eq:resolvent-path-metric}
\end{equation}
Lemma~\ref{W:lem:resolvent-path-space} proves completeness and separability. In particular, convergence in $\mathscr Y$ means uniform operator norm convergence on each compact
subset of $\mathbb C\setminus\mathbb R$.

\begin{theorem}
\label{W:thm:main-operator}
\label{W:thm:wigner-norm-resolvent}
Under Assumption~\ref{W:ass:wigner-lee-yin}, for every deterministic sequence of unit vectors $w_N$,
\begin{equation}
 \mathscr R_N\Longrightarrow\mathscr R_\beta\quad\text{in }\mathscr Y.
 \label{W:eq:main-resolvent-convergence}
\end{equation}
\end{theorem}

In terms of the operators, we write this conclusion as
\begin{equation}
 N^{2/3}(2I-T_N(w_N))\Longrightarrow\mathcal H_\beta
 \quad\text{in the norm-resolvent topology after embedding}.
 \label{W:eq:main-scaled-operator-limit}
\end{equation}
There are copies of the resolvent paths with the same individual laws for which, almost surely,
\begin{equation}
 \sup_{z\in K}\|\mathscr R_N(z)-\mathscr R_\beta(z)\|\longrightarrow0
 \quad\text{for every compact }K\Subset\mathbb C\setminus\mathbb R.
 \label{W:eq:skorokhod-uniform-resolvent}
\end{equation}
This follows from the Polish space formulation and is included in the proof in Section~\ref{W:sec:main-assembly}. The isometry from the original matrix coordinates is
$\mathcal E_NU_N(w_N)^*$, with $U_N$ fixed by the recursion above.

The conclusion is uniform over deterministic starting vectors in the following sense. Let $d_{\rm BL}$ be the bounded Lipschitz metric on probability measures on
$(\mathscr Y,d_{\mathscr Y})$, namely
\[
 d_{\rm BL}(\nu,\nu')=
 \sup_{\|F\|_\infty\le1,\,\operatorname{Lip}(F)\le1}
       \left|\int F\,\dd\nu-\int F\,\dd\nu'\right|.
\]
Here $\operatorname{Lip}(F)$ is computed with respect to $d_{\mathscr Y}$. For the laws of the embedded resolvents,
\begin{equation}
 \sup_{\|w\|=1}d_{\rm BL}\bigl(\mathcal L(\mathscr R_N^{(w)}),
                                      \mathcal L(\mathscr R_\beta)\bigr)
 \longrightarrow0.
 \label{W:eq:intro-uniform-distribution}
\end{equation}
The sequential argument establishing this uniform formulation is given in Section~\ref{W:sec:measurability-and-truncation}. In \eqref{W:eq:intro-uniform-distribution}, each law
corresponds to a vector $w$ fixed before the matrix is sampled.

\subsubsection{Local coefficients and compactness}
\label{W:sec:local-coefficients-and-full-spectral-localization}

In the quadratic form on a bounded interval, summation by parts pairs the cumulative potential with differences of the squared test function. We therefore prove local convergence
of the cumulative potential.

Write $(T_N)_{kk}=a_{k,N}$ and $(T_N)_{k,k+1}=b_{k,N}\ge0$. For the quadratic form we introduce the auxiliary coefficient at zero $\bar b_{0,N}=1$ and put $\bar b_{k,N}=b_{k,N}$
for $1\le k<N$. This boundary convention is separate from the initialization $b_0=0$ of the Lanczos recursion. Define
\begin{equation}
 \widetilde Y_{N,0}=0,\qquad
 \widetilde Y_{N,k}=m_N\sum_{r=1}^k
       (2-a_{r,N}-\bar b_{r-1,N}-b_{r,N}),\qquad 1\le k<N.
 \label{W:eq:operator-exact-integrated-potential}
\end{equation}
Let $\widetilde Y_N^{\rm lin}$ be its linear interpolation at $kh_N$, constant after $(N-1)h_N$. The step interpolation has value $\widetilde Y_{N,k}$ on $[kh_N,(k+1)h_N)$ and the
same constant extension. For $N=1$ both functions are zero. We write $\widetilde Y_N^+=\widetilde Y_N$ and define $\widetilde Y_N^-$ by replacing $a_{k,N}$ with $-a_{k,N}$.

For $L>0$ and $E\in\mathbb R$, set
\begin{equation}
 \begin{split}
 B_N(L)&=\max_{1\le k\le\min\{N-1,\lceil m_NL\rceil+1\}}
                                  |b_{k,N}-1|,\\
 P_{N,L}&=\sum_{k=1}^{\min\{N,\lfloor m_NL\rfloor\}}e_ke_k^*,
 \qquad Q_{N,L}=I-P_{N,L},\\
 \Pi_N(E)&=\mathbf1_{(-\infty,E]}(\mathcal A_N).
 \end{split}
 \label{W:eq:operator-coefficient-tail-notation}
\end{equation}
An empty maximum is zero. The spectral projection refers to the entire tridiagonal matrix. These cumulative sums are the integrated potentials used in \cite{ramirez2011beta} and
\cite{krishnapur2016universality}. Let
$Y_\beta(x)=x^2/2+2B(x)/\sqrt\beta$.

\begin{theorem}
\label{W:thm:integrated-potential-and-offdiagonal-coefficients}
Under Assumption~\ref{W:ass:wigner-lee-yin}, for every deterministic sequence of unit vectors $w_N$,
\begin{equation}
 \widetilde Y_N^{\rm lin}\Longrightarrow Y_\beta
                  \quad\text{in }C_{\rm loc}([0,\infty)).
 \label{W:eq:tightness-integrated-potential-limit}
\end{equation}
For every fixed $L,\varepsilon>0$,
\begin{equation}
 \sup_{\|w\|=1}\mathbb P(B_N(L)>\varepsilon)\longrightarrow0.
 \label{W:eq:tightness-offdiagonal-coefficients}
\end{equation}
Both reflected interpolations are tight in $C([0,T])$ for each finite $T$, uniformly over deterministic starting vectors.
\end{theorem}

The proof in Section~\ref{W:sec:proof-local-coefficients} compares the integrated potential with the corresponding GOE or GUE potential. Both reflected sums are needed: the exact
recursion in Theorem~\ref{W:thm:tightness-offdiagonal-coefficients} recovers the off-diagonal coefficients from them.

A low energy eigenfunction supported farther along the tridiagonal matrix can prevent operator norm convergence while leaving the coefficients in every fixed initial interval
unchanged. The following localization estimate controls such directions.

\begin{theorem}
\label{W:thm:localization-spectral-projection}
Under Assumption~\ref{W:ass:wigner-lee-yin}, for every $E\in\mathbb R$ and $\eta>0$,
\begin{equation}
 \lim_{L\to\infty}\limsup_{N\to\infty}\sup_{\|w\|=1}
       \mathbb P(\|Q_{N,L}\Pi_N(E)\|>\eta)=0.
 \label{W:eq:localization-spectral-projection}
\end{equation}
\end{theorem}
The projection $Q_{N,L}$ removes the first $\lfloor N^{1/3}L\rfloor$ Lanczos coordinates. Thus \eqref{W:eq:localization-spectral-projection} controls tails in the tridiagonal
coordinates. We use this bound to prove compactness. Separating local coefficient convergence from compactness follows the operator approach of \cite{ramirez2011beta} and
\cite{krishnapur2016universality}.

The estimate includes every invariant block after a breakdown. It is proved in Section~\ref{W:sec:proof-full-localization} by comparing the unnormalized trace described in
Section~\ref{W:sec:intro-localization-statistic}.

To see how the two theorems imply the operator limit, let $K\Subset\mathbb C\setminus\mathbb R$ be nonempty, set $d_K=\inf_{z\in K}|\operatorname{Im}z|$ and
$e_K=\sup_{z\in K}\operatorname{Re}z$, and take $E>e_K$. The spectral theorem gives
\begin{equation}
 \sup_{z\in K}\|Q_{N,L}(\mathcal A_N-z)^{-1}\|
 \le d_K^{-1}\|Q_{N,L}\Pi_N(E)\|+(E-e_K)^{-1}.
 \label{W:eq:intro-resolvent-localization}
\end{equation}
The first term is controlled by localization in the tridiagonal coordinates, and the second can be made small by increasing $E$. Local coefficient convergence identifies the
limiting quadratic form. The deterministic criterion in Section~\ref{W:sec:convergence-of-tridiagonal-operators} combines these facts; the complete proof of
Theorem~\ref{W:thm:main-operator} is given in Section~\ref{W:sec:main-assembly}.

\subsubsection{Spectral consequences}
\label{W:sec:resolvent-matrix-elements-and-spectral-consequences}

 Isolated spectral projections follow from resolvent convergence, whereas the ordering of the lowest eigenvalues and the weights at the starting vector also use lower spectral and
 boundary estimates. Throughout the corollaries below, Assumption~\ref{W:ass:wigner-lee-yin} and the deterministic unit vectors $w_N$ remain fixed. The proofs are given in
 Section~\ref{W:sec:comparison-of-the-original-ensembles-and-spectral-consequences}.

Write
\begin{equation}
 \Lambda_1<\Lambda_2<\cdots
 \label{W:eq:sao-eigenvalues}
\end{equation}
for the simple eigenvalues of $\mathcal H_\beta$ and choose real eigenfunctions with
\begin{equation}
 \mathcal H_\beta\varphi_j=\Lambda_j\varphi_j,\qquad
 \|\varphi_j\|_2=1.
 \label{W:eq:sao-eigenfunctions}
\end{equation}
Simplicity follows from the first-order differential system in Lemma~\ref{W:lem:continuous-integrated-potential-eigenfunctions}. Fix
$f_1,\ldots,f_p,g_1,\ldots,g_q\in L^2(\mathbb R_+)$, with $p,q$ independent of $N$, and use scalar products linear in the first argument. Write
\[
 \mathcal M_N(z)=(\langle\mathscr R_N(z)f_i,g_j\rangle)_{i,j},
 \qquad
 \mathcal M_\beta(z)=(\langle\mathscr R_\beta(z)f_i,g_j\rangle)_{i,j}.
\]

\begin{corollary}\label{W:cor:intro-resolvent-matrix-elements}
Under Assumption~\ref{W:ass:wigner-lee-yin},
\begin{equation}
 \mathcal M_N\Longrightarrow\mathcal M_\beta
 \quad\text{in }C_{\rm loc}(\mathbb C\setminus\mathbb R;
                                         \mathbb C^{p\times q}).
 \label{W:eq:intro-resolvent-matrix-elements}
\end{equation}
\end{corollary}

For functions supported in $[0,L]$, the corresponding vectors in the original matrix coordinates belong to the span of the first $\lceil m_NL\rceil$ Lanczos vectors; their exact
expression is \eqref{W:eq:intro-original-coordinate-resolvent}.

Define the embedded finite inverse to be zero when $\mathcal A_N$ is singular, as in Corollary~\ref{W:cor:wigner-edge-inverse}.

\begin{corollary}\label{W:cor:intro-real-inverse}
The singularity probability tends to zero, and $\mathcal E_N\mathcal A_N^{-1}\mathcal E_N^* \Longrightarrow\mathcal H_\beta^{-1}$ in operator norm, jointly with
\eqref{W:eq:main-resolvent-convergence}.
\end{corollary}

At a simple eigenvalue $\lambda_{j,N}$ of $\mathcal A_N$, choose a normalized eigenvector $v_{j,N}$. The embedded rank-one projection is
$\mathcal E_Nv_{j,N}v_{j,N}^*\mathcal E_N^*$. The limiting projection maps $f$ to $\langle f,\varphi_j\rangle\varphi_j$. The lower bound on the spectrum in
Proposition~\ref{W:prop:spectral-lower-bound-in-probability} is also needed to identify the ordered eigenvalues. It excludes eigenvalues escaping to $-\infty$, which nonreal
resolvent convergence alone would not exclude.

\begin{corollary}\label{W:cor:intro-ordered-edge-spectrum}
For each fixed $k$, the first $k$ eigenvalues of $\mathcal A_N$ and their embedded spectral projections converge jointly in distribution to those of $\mathcal H_\beta$, with
operator norm for the projections. In particular,
\[
 (N^{2/3}(\mu_j^{(N)}-2))_{j=1}^k
                    \Longrightarrow(-\Lambda_j)_{j=1}^k.
\]
\end{corollary}

The weights in the next statement depend on boundary derivatives. Their convergence is proved in Lemma~\ref{W:lem:ev-boundary-derivative}; $L^2$ convergence of the eigenfunctions
alone is insufficient.

Let $\xi_{j,N}$ be unit eigenvectors of $H_N$ with the decreasing ordering of its eigenvalues. Measurable orthonormal choices on a multiple eigenspace do not affect the following
limit; by the preceding corollary the first finitely many eigenvalues are simple with probability tending to one.

\begin{corollary}\label{W:cor:intro-boundary-weights}
For fixed $r,k$,
\begin{equation}
 \left((N^{2/3}(2-\mu_j^{(N)}))_{j=1}^r,
                         (N|w_N^*\xi_{j,N}|^2)_{j=1}^k\right)
 \Longrightarrow
 \left((\Lambda_j)_{j=1}^r,(|\varphi_j'(0)|^2)_{j=1}^k\right).
 \label{W:eq:intro-original-eigenvector-weights}
\end{equation}
The limiting weights are independent $\operatorname{Gamma}(\beta/2,2/\beta)$ variables, with shape and scale in that order, independent of the limiting eigenvalues.
\end{corollary}

Corollary~\ref{W:cor:ev-many-directions} gives projections of one fixed edge eigenvector onto finitely many deterministic directions, with an independent uniform sign or phase
assigned to the eigenvector.

\subsubsection{Outline of the comparison}
\label{W:sec:proof-outline}

The operator criterion requires the local coefficient limit and localization of the full low energy spectral subspace. We prove these first for the GOE when $\beta=1$ and the GUE
when $\beta=2$, then compare the relevant expectations in the independent matrix entries. Write $H_N^{\mathrm G}$ for the Gaussian matrix, with the normalization in
\eqref{W:eq:Gaussian-matrix-density}. Invariance makes the law of $T_N(w)$ independent of the deterministic starting vector. The coefficients are independent and satisfy
\[
 a_{k,N}\sim N(0,2/(\beta N)),\qquad
 b_{k,N}\sim(\beta N)^{-1/2}\chi_{\beta(N-k)}.
\]
Section~\ref{W:sec:gaussian-ensembles-and-airy-operator} proves this representation and the Gaussian process and localization bounds. We then compare these quantities by replacing
the independent matrix entries.

\paragraph{The integrated potential}
\label{W:sec:three-estimates-in-the-comparison}

For a Hermitian $H$ and a deterministic unit vector $w$, let $\mathcal L(f)=w^*f(H)w$. Use the fixed polynomials $p_j(x)=U_j(x/2)$ of
Section~\ref{W:sec:the-gram-representation-of-the-jacobi-coefficients}, and form $G_\rho=(\mathcal L(p_ip_j))_{0\le i,j<\rho}$. On $\kappa I\le G_\rho\le\kappa^{-1}I$,
orthogonalization gives smooth functions of the moments. At the required dimension $\rho\asymp m_N=N^{1/3}$, we must track the dependence of the derivative bounds on
$\rho$.

Let $q_j$ denote the orthonormal polynomials with positive leading coefficient. In this calculation $q_j$ is a polynomial; before breakdown, $q_j(H)w$ is the $(j+1)$st Lanczos
vector. For a fixed real moment direction $v$ satisfying $v(1)=0$, differentiated orthogonality gives
\[
 Dq_j[v]=-\tfrac12v(q_j^2)q_j-\sum_{i<j}v(q_iq_j)q_i.
\]
Indeed, differentiating $\mathcal L(q_iq_j)=0$ determines the coefficient of $q_i$ for $i<j$, since $Dq_i$ has degree at most $i$; differentiating $\mathcal L(q_j^2)=1$ gives the
diagonal coefficient. The three term recurrence then gives
\[
 Da_{j+1}[v]=b_{j+1}v(q_jq_{j+1})-b_jv(q_{j-1}q_j),\qquad
 Db_j[v]=\tfrac12b_jv(q_j^2-q_{j-1}^2).
\]
The term involving $q_{-1}$ is zero. In particular, summing the first identity leaves $b_Kv(q_{K-1}q_K)$. The interior mixed terms disappear before absolute values are taken.

The complete potential includes the off-diagonal contribution. For $1\le K\le\rho-2$ and $\sigma\in\{-1,1\}$, put
\[
 \mathcal T_K^\sigma
 =\sigma\sum_{j=1}^Ka_j+2\sum_{j<K}(b_j-1)+(b_K-1).
\]
Thus $\widetilde Y_{N,K}^\sigma=-m_N\mathcal T_K^\sigma$ at matrix moments. Summing the two variation identities with these weights gives
\begin{equation}
 D\mathcal T_K^\sigma[v]
 =v\!\left(\tfrac12b_K(q_K+\sigma q_{K-1})^2
 +\sum_{j=1}^{K-1}(b_j-b_{j+1})q_j^2-b_1q_0^2\right).
 \label{W:eq:intro-potential-variation}
\end{equation}
The interior coefficients $b_j-b_{j+1}$ vanish for the free Jacobi matrix. In general, their total variation controls the contribution of the interior squares.

The gain in this first variation also controls higher derivatives. Lemma~\ref{W:lem:tightness-directions} gives $\rho^{r-1}$ at order $r$ in the norm specified there. For
$\Theta=m_N\sum_l\lambda_l\mathcal T_{K_l}^{\sigma_l}$, a partition $\pi$ of $r$ derivative labels into $k$ nonempty sets contributes the dimension factor
\begin{equation}
 m_N^k\prod_{A\in\pi}\rho^{|A|-1}
 =m_N^k\rho^{r-k}\asymp m_N^r.
 \label{W:eq:overview-phase-partition}
\end{equation}
Thus all terms in the derivative of $e^{-\mathrm i\Theta}$ have the same power of $m_N$. The full product formula and its cutoff estimates are proved in
Proposition~\ref{W:prop:comparison-global}; the ensuing sum over entries is \eqref{W:eq:summed-fourth-remainder-terms}.

The bounds for the Gram matrix and the polynomial coefficients enter through smooth cutoffs. The comparison must also prove that these bounds hold under the distribution of the
Wigner matrix. We therefore compare the expectations of the cutoff functions under the original ensemble and under the GOE or GUE. Proposition~\ref{W:prop:comparison-global}
includes their derivatives as well as those of the observable. Applying the same calculation to the smooth functions of the modulus of continuity gives tightness; applying it to
the characteristic functions identifies the local limit. The size limit is taken with the cutoffs fixed before they are removed.

\paragraph{The third moment terms}
\label{W:sec:intro-third-moment-comparison}

A Taylor expansion in an off-diagonal entry contains a cubic term, since its third moments need not match those of the GOE or GUE. The first comparison uses bounded standardized
entries that match the first three off-diagonal moments of the normalized truncated law; their construction is in Lemma~\ref{W:lem:cubic-defect}. In the subsequent comparison with
the GOE and GUE, the cubic terms are summed with their signs. At first, each summand uses a different matrix with its selected entry removed. These summands cannot yet be combined
into one matrix product. Lemma~\ref{W:lem:cubic-entry-restoration} rewrites these terms with all resolvents evaluated at the same matrix, estimating the changes in both the
resolvents and the derivative coefficient. The following trace norm calculation is applied only after that replacement.

Write $\Delta\nu_N$ for the bounded off-diagonal law minus the off-diagonal law of the GOE or GUE, and set $\gamma_{3,N}=\int\zeta|\zeta|^2\,\Delta\nu_N(\dd\zeta)$. The common
off-diagonal distribution gives the deterministic Hermitian matrix $B_{3,N}$ with zero diagonal and upper-triangular entries $\gamma_{3,N}$. Lemma~\ref{W:lem:cumulant-trace-norm}
proves
\[
 \|B_{3,N}\|_1\le C|\gamma_{3,N}|N\log(2N).
\]
For real entries, $B_{3,N}=\gamma_{3,N}(\mathbf1\mathbf1^*-I)$ and the trace norm is $2(N-1)|\gamma_{3,N}|$. The trace norm estimate uses the constant triangular entries supplied
by the common entry law.

Let $B_{3,N}=\sum_l\lambda_lv_lv_l^*$ be a deterministic spectral decomposition and $R=(zI-H)^{-1}$. The spectral polar decomposition of $R$ and Cauchy--Schwarz give
\begin{equation}
 |w^*RB_{3,N}Rw|
 \le(w^*|R|w)\sum_l|\lambda_l|\,v_l^*|R|v_l.
 \label{W:eq:intro-signed-third-sum}
\end{equation}
The deterministic vector estimates for $|R|$ in Lemma~\ref{W:lem:cumulant-event} make the right side $N^{1+o(1)}$ on the event used in the comparison. The remaining coefficient and
contour factors contribute $m_NN^{o(1)}$. With the entry normalization, the resulting scale is
\[
 N^{-3/2}m_NN^{1+o(1)}=N^{-1/6+o(1)}.
\]
Theorem~\ref{W:thm:cumulant-complete} also includes the other cubic terms and the errors introduced by evaluating all resolvents at the same matrix.

The fourth-order and diagonal remainders are estimated conditionally on the other entries. Section~\ref{W:sec:truncation-and-resolvent-estimates} proves the required truncated
moment bounds under Assumption~\ref{W:ass:wigner-lee-yin}; the fourth moment enters through its value after truncation. The interpolation samples each independent entry from a
convex combination of the two distributions being compared.

\paragraph{Localization of the whole low energy subspace}
\label{W:sec:intro-localization-statistic}

Put $A=m_N^2(2I-H_N)$, $r=\min\{N,\lfloor m_NL\rfloor\}$, and $\widehat P_r=U_NP_{N,L}U_N^*$. For $\theta_E(x)=(1+e^{x-E})^{-2}$, define
$C_{r,E}=\operatorname{Tr}[(I-\widehat P_r)\theta_E(A)]$. Since $\theta_E\ge\tfrac14\mathbf1_{(-\infty,E]}$, functional calculus gives
\begin{equation}
 \|(I-\widehat P_r)\Pi_E\|^2
 \le4\operatorname{Tr}[(I-\widehat P_r)\theta_E(A)]
 =4\|(I-\widehat P_r)\theta_E(A)^{1/2}\|_{\rm HS}^2,
 \label{W:eq:intro-projection-trace-bound}
\end{equation}
where $\Pi_E=\mathbf1_{(-\infty,E]}(A)$. To obtain the inequality, multiply the functional calculus inequality on both sides by $I-\widehat P_r$ and use that the norm of a positive
matrix is at most its trace. In an orthonormal eigenbasis
$A\psi_j=\lambda_j\psi_j$,
\begin{equation}
 C_{r,E}=\sum_{j=1}^N\theta_E(\lambda_j)
                              \|(I-\widehat P_r)\psi_j\|^2.
 \label{W:eq:intro-spectral-mass}
\end{equation}
An eigenfunction with $\lambda_j\le E$ orthogonal to the first $r$ Lanczos vectors contributes at least $1/4$, whether or not its spectral weight at $w_N$ is zero. Dividing the
trace by $N$ would hide this single direction in a quantity tending to zero. This explains both the normalization and the need to use the whole spectral subspace.

The comparison uses a smooth version of this projection. Its agreement with the projection onto the initial vectors of the completed Lanczos basis is proved on the events
controlling the spectrum and the Gram matrix; the probabilities of their complements are estimated separately. The unnormalized trace has its own cubic term, treated in
Lemma~\ref{W:lem:localization-trace-cubic}. The Gaussian expectation is bounded by $C_{\beta,E}/L$. For fixed $L,E$ and cutoff parameter $\kappa$ for the eigenvalues of the Gram
matrix, the entry comparison errors tend to zero as $N\to\infty$. We then remove the cutoff in the eigenvalues of the Gram matrix by taking $\kappa$ to zero, and finally let
$L\to\infty$ to obtain localization of the full spectral projection. The comparison is made at each fixed $L$ before this last limit. The inequality with all exceptional
probabilities is proved in Section~\ref{W:sec:proof-full-localization}; the final application of the operator criterion is in Section~\ref{W:sec:main-assembly}.

\subsubsection{Background and relation to earlier work}
\label{W:sec:relation-to-earlier-work}

Wigner's moment calculation identified the global semicircle density for a random symmetric matrix ensemble \cite{wigner1955characteristic}. Dyson's symmetry classification
distinguishes real, complex, and quaternionic ensemble types \cite{dyson1962threefold}. Universality in random matrix theory predicts that microscopic spectral statistics are
determined by a small number of structural features of the ensemble. For Wigner matrices, the Wigner--Dyson--Mehta conjecture identifies the symmetry class as the parameter
governing local eigenvalue statistics after the appropriate rescaling \cite{mehta2004random}. The development of this conjecture has produced universality theorems for bulk and
edge eigenvalues, followed by joint eigenvalue and eigenvector laws and Gaussian limits for eigenvector projections
\cite{erdos2012bulk,erdos2012rigidity,bourgade2016fixed,knowles2013eigenvectors,bourgade2017momentflow}. These results suggest a common operator description of microscopic spectral
data: after choosing coordinates and an embedding, matrices in the same universality class should approach the same random operator. We establish this form of universality at the
upper spectral edge of real symmetric and complex Hermitian Wigner matrices.

The spectral theorem gives a useful formulation of the question. If a self-adjoint matrix $A_N$ has orthonormal eigenvectors $v_{j,N}$ and eigenvalues $\lambda_{j,N}$, and $J_N$ is
an isometry into a fixed Hilbert space, then
\[
 J_N(A_N-z)^{-1}J_N^*
 =\sum_{j=1}^N\frac{(J_Nv_{j,N})\otimes(J_Nv_{j,N})}{\lambda_{j,N}-z},
 \qquad z\notin\mathbb R,
\]
where $u\otimes u$ denotes the rank-one operator $f\mapsto\langle f,u\rangle u$. An operator limit controls the spectral locations and their projections in one topology. Operator
norm gives uniform control over the unit ball of the comparison space; convergence of the associated spectral subspaces then follows by contour integration. The extra analytic
requirements are tightness of these operator-valued objects and control of spectral tails and tails in the tridiagonal coordinates. This is the operator limit viewpoint developed
for the classical ensembles in \cite{edelman2007random,ramirez2011beta,virag2014operator}.

For Wigner matrices we choose the coordinates by Lanczos iteration from a deterministic unit vector $w_N$. Orthogonalization of $w_N,H_Nw_N,H_N^2w_N,\ldots$ yields the successive
coordinates and a tridiagonal representation $T_N(w_N)$. A fixed continuation at breakdown completes this to a representation of the whole matrix. Put $h_N=N^{-1/3}$ and let
$\mathcal E_Ne_k$ be the normalized indicator of $[(k-1)h_N,kh_N)$. Our main theorem states that the paths
\[
 z\longmapsto
 \mathcal E_N\bigl(N^{2/3}(2I-T_N(w_N))-z\bigr)^{-1}\mathcal E_N^*
\]
converge in distribution, locally uniformly in operator norm on $\mathbb C\setminus\mathbb R$, to the resolvent of
\[
 \mathcal H_\beta=-\frac{\mathrm d^2}{\mathrm dx^2}+x+
                   \frac{2}{\sqrt\beta}B',\qquad x>0,
\]
with Dirichlet boundary condition at zero. Here $B$ is standard Brownian motion and $\beta=1,2$ records the symmetry class. Theorem~\ref{W:thm:main-operator} gives the precise
assumptions and topology. Its hypotheses permit an infinite fourth off-diagonal moment and a nonidentically distributed diagonal array with uniformly integrable squares.

We derive a random differential operator from non-invariant Wigner matrices, including Bernoulli ensembles. The proof controls the induced recurrence through the independent
entries of the original matrix. Its limiting operator law is independent of the entry distribution and of the deterministic starting sequence. Operator limits for invariant
$\beta$-ensembles with general potentials \cite{krishnapur2016universality,rider2019universality} and the different global scaling considered in \cite{banerjee2025operator} are
discussed below.

The Gaussian tridiagonal model indicates both the coordinates and the scales. Near the first coordinate, its off-diagonal coefficients have a slowly varying mean close to one,
while the centered diagonal and off-diagonal coefficients have fluctuations of order $N^{-1/2}$ \cite{dumitriu2002matrix}. Under the scaling $k\sim N^{1/3}x$, the second difference
becomes a second derivative, the mean variation produces the linear potential $x$, and the cumulative fluctuations converge to Brownian motion. A bounded interval therefore
involves order $N^{1/3}$ steps of orthogonalization. The extension to Wigner matrices requires comparison estimates whose constants remain controlled at this growing degree.

We compare the cumulative potential by differentiating its sum. Cancellation between neighboring terms gives the degree dependence required when comparing the independent matrix
entries. A separate comparison of an unnormalized trace controls the mass of the entire low energy spectral subspace outside the first Lanczos coordinates. These two estimates
identify the limiting differential expression and give operator norm compactness. Boundary estimates then recover the spectral weights of the starting vector in the original matrix
coordinates.

\subparagraph{Eigenvalue and eigenvector universality} The Wigner--Dyson--Mehta conjecture compares local eigenvalue statistics with those of the Gaussian ensemble of the same
symmetry class. Their exact spectral formulas and the role of symmetry are developed in \cite{mehta2004random,forrester2010log}; general introductions to Wigner matrices and the
universality problem include \cite{anderson2010introduction,erdos2017dynamical,benaych2018lectures}. For bulk eigenvalues, Johansson \cite{johansson2001universality} proved
universality for a class of Gaussian divisible Hermitian Wigner matrices. The work of Erd\H{o}s, Ram\'irez, Schlein and Yau \cite{erdos2010sine} and the local relaxation flow of
Erd\H{o}s, Schlein and Yau \cite{erdos2011relaxation} developed the dynamical approach. Tao and Vu \cite{tao2011local} established the four moment comparison theorem. Bulk
universality for generalized Wigner matrices \cite{erdos2012bulk}, gap universality \cite{erdos2015gap}, and fixed energy universality \cite{bourgade2016fixed} successively refined
the comparison of local spectral statistics.

At the soft edge, Tracy and Widom \cite{tracy1994airy,tracy1996orthogonal} identified the limiting largest eigenvalue laws for the classical Gaussian ensembles. Soshnikov's moment
method \cite{soshnikov1999edge} established their universality for Wigner matrices under symmetry and moment assumptions. Further comparison results include Tao and Vu
\cite{tao2010edge}, the rigidity and edge comparison estimates of Erd\H{o}s, Yau and Yin \cite{erdos2012rigidity}, and the necessary and sufficient tail criterion of Lee and Yin
\cite{lee2014necessary}. The tail condition used here is that criterion. Our diagonal assumptions allow an array of laws with uniformly integrable squares; the required conditional
estimates are proved in Section~\ref{W:sec:truncation-and-resolvent-estimates}.

Eigenvector universality gives complementary information about the spectral projections. Knowles and Yin \cite{knowles2013eigenvectors} proved joint comparison of eigenvalues and
eigenvectors, including edge comparison under subexponential tail assumptions, two moment matching, and level repulsion. Tao and Vu \cite{tao2012eigenvectors} extended their
comparison theorem to eigenvector entries and resolvent statistics. The eigenvector moment flow of Bourgade and Yau \cite{bourgade2017momentflow} proves universality of eigenvector
projections and local quantum unique ergodicity. Subsequent results include joint Gaussian limits for bulk eigenvector projections of real symmetric generalized Wigner matrices
\cite{marcinek2022normality}, fluctuations of eigenvector mass \cite{benigni2022fluctuations}, and eigenstate thermalization estimates
\cite{cipolloni2021thermalization,cipolloni2024edge}. These theorems motivate an operator universality principle in a specified microscopic regime. Its formulation includes the
choice of coordinates, the limiting Hilbert space, the boundary condition, and the topology on resolvents. In the present theorem that topology is locally uniform operator norm
convergence in distribution. The ordering of the lowest eigenvalues also uses a lower spectral bound, and the boundary weight limits require convergence of boundary derivatives.

\subparagraph{From Gaussian matrices to a differential operator} Dumitriu and Edelman \cite{dumitriu2002matrix} constructed explicit tridiagonal and bidiagonal models for the
Hermite and Laguerre $\beta$-ensembles. For the GOE and GUE, orthogonal or unitary invariance identifies the tridiagonal model with Lanczos iteration from any fixed unit vector. In
our normalization its coefficients are independent and have laws
\[
 a_{k,N}^{G}=\frac{g_k}{\sqrt{\beta N}},\quad g_k\sim N(0,2),
 \qquad
 b_{k,N}^{G}=\frac{\chi_{\beta(N-k)}}{\sqrt{\beta N}}.
\]
Edelman and Sutton \cite{edelman2007random} interpreted the classical matrix models as difference approximations to stochastic differential operators. Ram\'irez, Rider and Vir\'ag
\cite{ramirez2011beta} constructed the stochastic Airy operator and proved the soft edge limit through convergence and compactness of the discrete quadratic forms. Operator
convergence and its spectral interpretation are surveyed in \cite{virag2014operator}.

The scaling can be read from the Gaussian coefficients. If $k\asymp m_N$ and $m_N=N^{1/3}$, their means satisfy
\[
 \mathbb E b_{k,N}^{G}
 =\sqrt{1-k/N}+O(N^{-1})
 =1-\frac{k}{2N}+O(k^2/N^2+N^{-1}).
\]
For $x=kh_N$ and $h_N=m_N^{-1}$, the deterministic second difference is
\[
 h_N^{-2}\bigl(2f(x)-f(x-h_N)-f(x+h_N)\bigr)\longrightarrow-f''(x).
\]
At the same time, $m_N^2(2-\mathbb E b_{k-1,N}^{G}-\mathbb E b_{k,N}^{G})\to x$ for $k/m_N\to x>0$. The centered coefficient errors accumulate on this scale. With the auxiliary
form coefficient $\bar b_{0,N}^{G}=1$, the polygonal interpolation of
\[
 \widetilde Y_{N,k}^{G}
 =m_N\sum_{j=1}^{k}
     (2-a_{j,N}^{G}-\bar b_{j-1,N}^{G}-b_{j,N}^{G})
\]
converges locally to $Y_\beta(x)=x^2/2+2B(x)/\sqrt\beta$. The diagonal contribution has limiting variance $2x/\beta$; the off-diagonal contribution has the same limiting variance,
since interior coefficients occur twice in the cumulative sum. Their independence gives the coefficient $2/\sqrt\beta$ of the Brownian term.
Section~\ref{W:sec:convergence-of-gaussian-coefficient-process} gives the complete process calculation, including its boundary terms.

Summation by parts identifies the limiting quadratic form. For compactly supported smooth $f$ with $f(0)=0$, it is
\[
 \mathfrak q_\beta[f]
 =\int_0^\infty |f'|^2\,\mathrm dx
   -\int_0^\infty Y_\beta(x)(|f(x)|^2)'\,\mathrm dx.
\]
Substituting $Y_\beta(x)=x^2/2+2B(x)/\sqrt\beta$ gives
\[
 \mathfrak q_\beta[f]
 =\int_0^\infty(|f'|^2+x|f|^2)\,\mathrm dx
   -\frac{2}{\sqrt\beta}\int_0^\infty B(x)(|f(x)|^2)'\,\mathrm dx.
\]
This construction gives a pathwise meaning to the white noise potential and retains the boundary condition selected by the first matrix coordinate. Coercivity and localization
control the passage from bounded intervals to the half-line. The resulting operator has discrete spectrum; the negatives of its eigenvalues have the upper soft edge point process
law. Its eigenfunctions also determine the orthogonal projections onto its eigenspaces \cite{ramirez2011beta,krishnapur2016universality}. For a general cyclic pair $(H,w)$, the
same Lanczos construction orthogonalizes polynomials in its spectral measure
\[
 \mu_{H,w}=\sum_j|w^*\xi_j|^2\delta_{\mu_j},
 \qquad H\xi_j=\mu_j\xi_j,
 \qquad
 w^*(H-z)^{-1}w=\int\frac{\mu_{H,w}(\mathrm d\lambda)}{\lambda-z}.
\]
This identity explains the connection with eigenvector weights. The continued construction specified in this paper additionally represents invariant subspaces with zero weight at
$w$.

\subparagraph{Other regimes and forms of operator convergence} The same program has produced several random differential operators. Ram\'irez and Rider
\cite{ramirez2009diffusion,ramirez2011erratum} obtained the stochastic Bessel operator at the hard edge of the Laguerre ensembles. Holcomb and Moreno Flores \cite{holcomb2012edge}
established soft- and hard edge limits for the $\beta$-Jacobi ensemble. The sparse circular models of Killip and Nenciu \cite{KillipNenciu2004} and the recurrence analysis of
Killip and Stoiciu \cite{KillipStoiciu2009} led to a corresponding bulk theory. Valk\'o and Vir\'ag developed the Brownian carousel \cite{valko2009continuum}, identified its point
process with the spectrum of the Sine$_\beta$ operator \cite{valko2017sine}, and constructed quantitative operator couplings for the circular $\beta$-ensemble
\cite{valko2020operator}.

Gaussian perturbations also have operator interpretations. Bloemendal and Vir\'ag \cite{bloemendal2013limits,bloemendal2016limits} encoded finite-rank spikes through boundary
conditions and matrix-valued stochastic Airy operators; hard edge spikes were treated in \cite{ramirez2017spiking}. Gorin and Shkolnikov \cite{gorin2018stochastic} obtained the
stochastic Airy semigroup from powers of tridiagonal matrices. Representations using reflected Brownian motion and convergence of semigroups generated by more general tridiagonal
matrices appear in \cite{gaudreaulamarre2019spiked,gaudreaulamarre2020convergence}. Recent solvable block tridiagonal models give further matrix-valued operator limits
\cite{rider2026block}. Across these constructions, an explicit finite-dimensional representation makes the coefficient fluctuations accessible and specifies the limiting
differential expression and boundary data.

\subparagraph{General potentials and non-invariant Wigner matrices} Operator universality already extends beyond quadratic invariant potentials. Krishnapur, Rider and Vir\'ag
\cite{krishnapur2016universality} proved norm resolvent convergence to the stochastic Airy operator for the Jacobi matrices of Dyson $\beta$-ensembles with uniformly convex
polynomial potentials. Their starting point is the explicit joint density of the recurrence coefficients, whose dependence is analyzed through the potential. Rider and Waters
\cite{rider2019universality} established convergence to the stochastic Bessel operator for a class of hard edge models with general potentials. General tridiagonal convergence
theorems \cite{gaudreaulamarre2020convergence} and scaling limits for discrete random Schr\"odinger operators \cite{kritchevski2012scaling} give further examples with non-Gaussian
coefficient laws. The common feature of these results is direct access to a one-dimensional matrix or recurrence model.

The independent entries of a non-invariant Wigner matrix supply a different starting point. Its Lanczos coefficients are dependent nonlinear functions of the entire matrix and of
the starting vector. Theorem~\ref{W:thm:main-operator} obtains their soft edge operator limit under Assumption~\ref{W:ass:wigner-lee-yin}. Banerjee~\cite{banerjee2025operator}
constructed an operator limit for Wigner matrices at the global spectral scale on a Hilbert-space extension of $L^2[0,1]$, under the bounds on the growth of moments in Theorem~5
there. The convergence is formulated through matrix powers and inner products, and Theorem~6 identifies semicircular spectral measures at the embedded deterministic unit functions.
The present scaling resolves the upper edge on the $N^{-2/3}$ scale and produces the stochastic Airy operator on $L^2(\mathbb R_+)$.

\subparagraph{Comparison and the present proof} Green function comparison at the edge was developed by Erd\H{o}s, Yau and Yin; see Section~6.1 and Theorem~6.3
of~\cite{erdos2012rigidity}. Their comparison permits nonzero third moments under the stated subexponential decay and second-moment hypotheses. Here the observables contain the
inverse of a growing Gram matrix and order $N^{1/3}$ recurrence coefficients. The derivative identity \eqref{W:eq:intro-potential-variation}, its higher order bounds, and the
separate trace comparison control this degree dependence. We make the probabilistic estimates in the independent Wigner entries and compare with the GOE and GUE through product
measures of entry laws.

We use the isotropic local semicircle law \cite{bloemendal2014isotropic} for the bounded moment interpolating ensembles, together with conditional row estimates for the truncated
arrays. Entry perturbations are supported on at most two coordinates, and the products are handled by finite-rank contractions. The precise hypotheses and spectral scales of the
published local laws used in these estimates are recorded in Section~\ref{W:sec:multi-resolvent-identities-and-entry-expansions}. The boundary recurrence identifies the original
coordinate weights once the joint coefficient and operator limit is available.

Sections~\ref{sec:soft-edge-ensembles-and-quantitative-airy-limits}, \ref{sec:hard-edge-models-brownian-information-and-bessel-limits},
and~\ref{sec:circular-ensembles-canonical-systems-and-the-sine-coefficient} study quantitative operator approximation in the soft edge, hard edge, and circular models. The present
theorem supplies universality of the soft edge operator limit for the Wigner matrix distributions in Assumption~\ref{W:ass:wigner-lee-yin}, with convergence in distribution in the
specified resolvent topology.

Section~\ref{W:sec:convergence-of-tridiagonal-operators} proves the operator convergence criterion. Its application to the two probabilistic theorems stated above is given in
Section~\ref{W:sec:main-assembly}, after the localization proof. Section~\ref{W:sec:orthogonal-polynomials-and-spectral-statistics} develops the moment formulas and derivative
estimates, and Section~\ref{W:sec:gaussian-ensembles-and-airy-operator} establishes the required estimates for the GOE and GUE.
Sections~\ref{W:sec:truncation-and-resolvent-estimates} and~\ref{W:sec:green-function-comparison} prove truncation and comparison for the local coefficients.
Section~\ref{W:sec:localization-of-the-low-energy-spectral-subspace} proves localization of the full low energy subspace.
Section~\ref{W:sec:comparison-of-the-original-ensembles-and-spectral-consequences} contains the spectral consequences and boundary estimates.

\subsection{Gaussian and polynomial ensembles at the soft edge}
\label{sec:soft-edge-ensembles-and-quantitative-airy-limits}
\label{A:sec:introduction}

For the theorem for Wigner matrices the objective is a common limiting operator. In the classical $\beta$-ensembles, the coefficient laws allow a further question: how accurately
can the finite matrix approximate one realization of that operator while retaining its prescribed distribution? At the soft edge, we answer this question for the Gaussian model and
the polynomial Jacobi models. The comparison includes eigenfunctions and resolvents as well as eigenvalues, in the operator formulation discussed in \cite{virag2014operator}.

We study this question at a regular soft edge. The limiting random differential operator is
\[
 \mathcal H_\beta=-\frac{\mathrm d^2}{\mathrm dx^2}+x+c_\beta B'(x),
 \qquad c_\beta=\frac2{\sqrt\beta},
\]
on $L^2(\mathbb R_+)$, with Dirichlet boundary condition at zero and $B$ a standard Brownian motion. Its form domain and operator domain are specified in
Section~\ref{A:sec:intro-limiting-operator}. Ram\'irez, Rider, and Vir\'ag \cite{ramirez2011beta} obtained the Gaussian $\beta$-ensemble limit from the tridiagonal model of
Dumitriu and Edelman \cite{dumitriu2002matrix}. Krishnapur, Rider, and Vir\'ag \cite{krishnapur2016universality} extended the operator limit to Jacobi matrices with uniformly
convex polynomial potentials. We ask how accurately these matrices can approximate $\mathcal H_\beta$ when their coefficient laws and their comparison maps are fixed before the
coupling is chosen.

We compare inverses and resolvents in Hilbert--Schmidt norm, $\|A\|_{\mathfrak S_2}^2=\operatorname{Tr}(A^*A)$. Since this norm includes the entire reciprocal spectrum, the high
eigenvalues of the limiting operator impose a lower bound on the approximation error. Write $R_\beta(z)=(\mathcal H_\beta-z)^{-1}$. The ordered eigenvalues of $\mathcal H_\beta$
satisfy $\Lambda_k\sim(3\pi k/2)^{2/3}$ almost surely (Lemma~\ref{A:res:weyl}; see also \cite{virag2014operator} (Theorem~6.1)). For a nonempty compact
$K\Subset\mathbb C\setminus\mathbb R$, let $s_k$ be the decreasing singular values of $R_\beta(z)$. Then
\[
 \inf_{\operatorname{rank}A\le n}
       \|A-R_\beta(z)\|_{\mathfrak S_2}^2
   =\sum_{k>n}s_k^2
   \asymp_{B,K}\sum_{k>n}k^{-4/3}
   \asymp n^{-1/3},\qquad z\in K.
\]
Subsection~\ref{A:sec:rank-lower-bound-and-optimality} proves the equality and the uniformity in $z$.

The spectral truncation attaining this minimum depends on the limiting eigenfunctions. We obtain the same order for the matrix laws and comparison maps specified below. We couple
all sizes to one Brownian motion while preserving each finite coefficient law. The result holds for every fixed $\beta>0$, for the Gaussian model and for the polynomial Jacobi laws
at a regular soft edge.

The comparison recovers convergence of each fixed group of eigenvalues and of its spectral projections, and also controls quantities involving the whole spectrum. For a fixed real
shift outside the limiting spectrum, Theorem~\ref{A:res:20.6} proves that the optimal squared $\ell^2$ matching distance between the full reciprocal spectra, with unmatched values
paired with zero, has order $n^{-1/3}$. This quantity depends only on the eigenvalues and is independent of the coordinates used to compare the operators. The full matrix results
in Corollary~\ref{A:res:full-ordered} also imply convergence of $\operatorname{Tr}f(\mathcal A_{n,V,\beta})$ for continuous compactly supported $f$. These conclusions concern every
matrix eigenvalue, including those associated with coordinates outside the range of the comparison map.

\subsubsection{The matrix models and their normalization}\label{A:sec:intro-models-scaling}
The Gaussian and polynomial models can be written with the same Jacobi coefficient density. We specify the deterministic normalization for each model. The parameter $\beta>0$ is
fixed throughout.

Let $V$ be a real polynomial of even degree with $V''\ge c_V>0$. The $\beta$-ensemble with potential $V$ has eigenvalue density
\[
 \frac1{Z_{n,V,\beta}}
 \exp\left\{-\frac{\beta n}{2}\sum_{k=1}^nV(\lambda_k)\right\}
 \prod_{j<k}|\lambda_j-\lambda_k|^\beta
 \prod_{k=1}^n\mathrm d\lambda_k.
\]
Let $T(a,b)$ be the real symmetric tridiagonal matrix with diagonal $a_1,\ldots,a_n$ and positive off-diagonal entries $b_1,\ldots,b_{n-1}$. Its coefficient law is
\[
 \mathrm d\mu_{n,V,\beta}(a,b)
 =\frac1{\mathcal Z_{n,V,\beta}}
   \exp\left\{-\frac{\beta n}{2}\operatorname{Tr}V(T(a,b))\right\}
   \prod_{k=1}^{n-1}b_k^{\beta(n-k)-1}
   \prod_{k=1}^n\mathrm da_k\prod_{k=1}^{n-1}\mathrm db_k.
 \eqtag{jacobi-law}
\]
Under \eqref{A:eq:jacobi-law}, the eigenvalues of $T(a,b)$ have the preceding density \cite{krishnapur2016universality}, with the potential in that reference replaced by $V/2$.
Write $T_{n,V,\beta}=T(a,b)$.

For $V(t)=t^2/2$,
\[
 \operatorname{Tr}V(T)=\frac12\sum_{k=1}^na_k^2
                              +\sum_{k=1}^{n-1}b_k^2.
\]
The density therefore factorizes, with independent coefficients
\[
 a_k\sim N\left(0,\frac2{\beta n}\right),\qquad
 b_k\sim\frac1{\sqrt{\beta n}}\chi_{\beta(n-k)}.
\]
These are the Dumitriu--Edelman matrices $T_{n,\beta}^{\mathrm{DE}}$ \cite{dumitriu2002matrix}. For higher degree potentials, the trace couples neighboring coefficients. The
Gaussian model admits a direct coupling of independent coefficients; for the general law we construct the coupling conditionally on exterior coefficients.

To locate and scale the edge for a general $V$, let $J(a,b)$ be the Jacobi operator on $\ell^2(\mathbb Z)$ with constant coefficients, $(J(a,b)u)_j=bu_{j-1}+au_j+bu_{j+1}$. Define
\[
 W_V(a,b)=[V(J(a,b))]_{00}
       =\frac1{2\pi}\int_0^{2\pi}V(a+2b\cos\theta)\,\mathrm d\theta.
 \eqtag{11.4}
\]
The Fourier symbol $a+2b\cos\theta$ gives both this identity and the upper edge $a+2b$ of the constant operator. At macroscopic index $s$, the leading local energy in
\eqref{A:eq:jacobi-law} also contains the logarithmic term from the off-diagonal density. For $0\le s<1$, let $(a^\dagger(s),b^\dagger(s))$ minimize
\[
 F_{V,s}(a,b)=W_V(a,b)-2(1-s)\log b,
 \qquad a\in\mathbb R,\quad b>0.
 \eqtag{local-energy}
\]
This minimizer exists, is unique, and is analytic in $s$; see Lemmas~\ref{A:res:11.1}--\ref{A:res:11.2} and Proposition~\ref{A:res:11.3}. Put
\[
 a_\star=a^\dagger(0),\qquad b_\star=b^\dagger(0),\qquad E(s)=a^\dagger(s)+2b^\dagger(s),\qquad
 E_\star=E(0),\qquad \tau_\star=-E'(0)>0.
 \eqtag{edge-data}
\]
The last inequality is our regular edge assumption. The spacing and the matrix normalization are chosen so that the second difference coefficient and the slope of the limiting
potential are both one. Explicitly,\label{A:par:edge-scaling}
\[
 h_{n,V}=\left(\frac{b_\star n}{\tau_\star}\right)^{-1/3},
 \qquad
 \mathcal A_{n,V,\beta}
       =\frac{h_{n,V}^{-2}}{b_\star}(E_\star I-T_{n,V,\beta}).
 \eqtag{edge-scaling}
\]
With $h=h_{n,V}$, the factor $b_\star^{-1}$ makes the constant coefficient difference term equal to $-\Delta_{h,D}$. For $s=x/(nh)$ in a fixed compact subinterval of $[0,1)$, the
expansion $E_\star-E(s)=\tau_\star s+O(s^2)$ gives
\[
 \frac{h^{-2}}{b_\star}(E_\star-E(s))
   =\frac{\tau_\star}{b_\star nh^3}x+O(h^2x^2)
   =x+O(h^2x^2).
\]
The finite-$n$ coefficient calculation, including the center and logarithmic shifts, is in Subsection~\ref{A:sec:the-local-minimizer-and-edge-scaling}. In the Gaussian case
$W_V(a,b)=a^2/2+b^2$, so $a^\dagger(s)=0$, $b^\dagger(s)=\sqrt{1-s}$, and $(E_\star,b_\star,\tau_\star)=(2,1,1)$. This gives the Gaussian normalization
$\mathcal A_{n,\beta}^{\mathrm{DE}} =n^{2/3}(2I-T_{n,\beta}^{\mathrm{DE}})$.

For the comparison on $L^2(\mathbb R_+)$, use the isometry
\[
 (\iota_h^{(m)}u)(x)
   =h^{-1/2}u_k\quad((k-1)h\le x<kh,\ 1\le k\le m),
 \qquad
 (\iota_h^{(m)}u)(x)=0\quad(x\ge mh).
 \eqtag{step-embedding}
\]
Write $x_k=(k-\tfrac12)h$ for the interval centers. The conditional coefficient analysis for general $V$ is performed away from the original boundaries. Accordingly, the comparison
map omits a logarithmic number of coordinates at the left and a fixed number at the right. Set $R_V=\deg V$ and
\[
 r_n=\lceil A_0\log n\rceil,\qquad s_\beta=\lceil\beta^{-1}\rceil-1+2R_V+10,\qquad
 I_n=\{r_n+1,\ldots,n-s_\beta\},\qquad M_n=|I_n|.
\]
Here $A_0$ is a sufficiently large deterministic constant, fixed in Subsection~\ref{A:sec:boundary-reduction}. Let $P_{I_n}$ restrict to these coordinates and let $U_{r_n}$ reindex
them from $1$. Define
\[
 \mathcal J_n=\iota_{h_{n,V}}^{(M_n)}U_{r_n}P_{I_n}.
 \eqtag{1.1}
\]
When $r_n+s_\beta\ge n$, take $I_n$ empty and $\mathcal J_n=0$. For the other sizes, $\mathcal J_n^*\mathcal J_n=P_{I_n}^*P_{I_n}$ and
$\operatorname{rank}\mathcal J_n=M_n=n-O(\log n)$. Thus the operator in the theorem is the compressed full resolvent $\mathcal J_n(\mathcal A_{n,V,\beta}-z)^{-1}\mathcal J_n^*$.

The left separation occupies $r_n$ sites, with physical length $r_nh\to0$. At the right, the removed interval contains the indices with $\beta(n-k)-1<0$ and all coefficients
interacting with them. The interaction range is at most $R_V$, as follows from the trace expansion in Subsection~\ref{A:sec:the-local-minimizer-and-edge-scaling}.
Subsection~\ref{A:sec:boundary-reduction} estimates both changes by Schur complements. These comparisons leave the law \eqref{A:eq:jacobi-law} of the full matrix unchanged.

\subsubsection{The limiting operator}\label{A:sec:intro-limiting-operator}
The white noise in the formal expression for the stochastic Airy operator is interpreted through a quadratic form, rather than as a pointwise multiplication operator. This also
fixes the domain and the boundary condition used in the comparison.

We work on complex Hilbert spaces, with inner products linear in the first argument. A sesquilinear form $\mathfrak q$ is linear in its first argument and conjugate linear in its
second. It is Hermitian when $\mathfrak q[f,g]=\overline{\mathfrak q[g,f]}$; its associated quadratic form is $\mathfrak q[f]=\mathfrak q[f,f]$. On real vectors, the restriction of
a Hermitian form with real coefficients is a symmetric bilinear form. Whenever $\mathcal E^*$ denotes the dual of an energy space below, it is the conjugate dual, with pairing
$\langle F,v\rangle$ conjugate linear in $v$.

Let $B$ be standard Brownian motion and $c_\beta=2/\sqrt\beta$. The common form domain is
\[
 \mathcal D_{\mathrm A}
 =\left\{f\in H_0^1(\mathbb R_+):
                 \int_0^\infty x|f(x)|^2\,\mathrm dx<\infty\right\}.
\]
The deterministic Hermitian sesquilinear form $\mathfrak a$ is defined on $\mathcal D_{\mathrm A}\times\mathcal D_{\mathrm A}$. Initially define the sesquilinear form
$\mathfrak w_B$ on compactly supported pairs in this domain:
\[
 \begin{aligned}
 \mathfrak a[f,g]&=\int_0^\infty f'\overline{g'}\,\mathrm dx
                    +\int_0^\infty xf\overline g\,\mathrm dx,\qquad
 \mathfrak w_B[f,g]=-\int_0^\infty B(f\overline g)'\,\mathrm dx,\\
 \mathfrak q_\beta[f,g]&=\mathfrak a[f,g]+c_\beta\mathfrak w_B[f,g].
 \end{aligned}
 \eqtag{airy-form}
\]
The perturbation $\mathfrak w_B$ has a continuous extension to $\mathcal D_{\mathrm A}$ and is infinitesimally bounded relative to $\mathfrak a$, almost surely. The estimates in
Section~\ref{A:sec:galerkin-approximation-and-the-resolvent-error} prove that $\mathfrak q_\beta$ is closed and bounded below on this domain. The associated operator is specified
by
\[
 \begin{split}
 \operatorname{dom}\mathcal H_\beta
 &=\{f\in\mathcal D_{\mathrm A}:\text{ there is }g\in L^2(\mathbb R_+)
       \text{ with }\mathfrak q_\beta[f,v]=\langle g,v\rangle
       \text{ for every }v\in\mathcal D_{\mathrm A}\},\\
 \mathcal H_\beta f&=g.
 \end{split}
\]
The trace $f(0)=0$ is imposed by $H_0^1(\mathbb R_+)$. The quadratic form interprets $B'$ by integration against the continuous path $B$. Write
$R_\beta(z)=(\mathcal H_\beta-z)^{-1}$ and let $\Lambda_1<\Lambda_2<\cdots$ be the eigenvalues of $\mathcal H_\beta$, in increasing order. Compactness and simplicity are proved in
Section~\ref{A:sec:atomlessness-real-resolvents-and-inverse-coupling} and the form estimates cited above.

\subsubsection{Main results}\label{A:sec:intro-main-results}
The following statements cover both normalizations. In the general model the comparison map compresses the full inverse to the retained coordinates.

The matrices, comparison maps and spacings in the two cases are\label{A:par:notation}
\[
 (\mathscr A_n,\mathscr I_n,h_n)
 =\begin{cases}
  (\mathcal A_{n,\beta}^{\mathrm{DE}},\iota_{n^{-1/3}}^{(n)},n^{-1/3}),
       &\text{Gaussian model},\\
  (\mathcal A_{n,V,\beta},\mathcal J_n,h_{n,V}),
       &\text{general Jacobi model}.
 \end{cases}
\]
In both cases $h_n\asymp n^{-1/3}$. The model, $\beta$, the potential where applicable, and the comparison maps are fixed throughout each coupling. The Hilbert--Schmidt norm is
unnormalized.

\begin{theorem}\label{A:res:1.1}
There is a coupling of one standard Brownian motion $B$ and all the matrices $T_n$, each with its original coefficient law, such that almost surely $\mathcal H_\beta$ is invertible
and, for all sufficiently large $n$, $\mathscr A_n$ is invertible and
\[
 c_\omega h_n
 \le \|\mathscr I_n\mathscr A_n^{-1}\mathscr I_n^*
                  -\mathcal H_\beta^{-1}\|_{\mathfrak S_2}^2
 \le C_\omega h_n,
 \qquad 0<c_\omega\le C_\omega<\infty.
\]
\end{theorem}

The matrices at different sizes need not be independent. The embedding for the Gaussian ensemble uses all coordinates; $\mathcal J_n$ uses $n-O(\log n)$ of them. In both cases the
inverse is taken before applying the comparison map. For the general model, $\mathcal J_n$ compresses the full inverse to the retained coordinates. The squared norm in
Theorem~\ref{A:res:1.1} has order $n^{-1/3}$; the norm itself has order $n^{-1/6}$.

The same coupling gives a locally uniform estimate for the nonreal resolvents.
\begin{theorem}\label{A:res:1.2}
Almost surely, for every nonempty compact $K\Subset\mathbb C\setminus\mathbb R$ and all sufficiently large $n$,
\[
 c_{\omega,K}h_n
 \le\|\mathscr I_n(\mathscr A_n-z)^{-1}\mathscr I_n^*
                       -R_\beta(z)\|_{\mathfrak S_2}^2
 \le C_{\omega,K}h_n
 \qquad(z\in K),
\]
where $0<c_{\omega,K}\le C_{\omega,K}<\infty$.
\end{theorem}
The event in this statement is common to all $K$; the constants and the threshold in $n$ may depend on $K$.

The lower bound holds for every approximation of rank at most $n$, regardless of its law. Subsection~\ref{A:sec:rank-lower-bound-and-optimality} computes the leading order
coefficient of this unrestricted minimum. The matrix coupling attains this optimal order. Its leading order coefficient remains undetermined; other comparison errors may also
contribute at order $h_n$.

The operator comparison also yields spectral information on the retained coordinates. Together with the common lower spectral bound and the estimates for the omitted coordinates,
it controls the ordered eigenvalues of the full finite matrix. Write $\lambda_{1,n}\le\cdots\le\lambda_{n,n}$ for those of $\mathscr A_n$. Theorem~\ref{A:res:21.1} and
Corollary~\ref{A:res:full-ordered} give, for each fixed $m$,
\[
 \max_{k\le m}|\lambda_{k,n}-\Lambda_k|
       \le C_{\omega,m}h_n^{1/2}
\]
almost surely for all sufficiently large $n$. The corresponding spectral projections converge, and
\[
 \sum_{j=1}^n\delta_{\lambda_{j,n}}
       \longrightarrow\sum_{k\ge1}\delta_{\Lambda_k}
\]
almost surely in the vague topology. These conclusions are proved in Subsection~\ref{A:sec:spectral-consequences}. The rate for a fixed eigenvalue is a consequence of the operator
bound; we do not prove that this eigenvalue rate is optimal.

If $\theta_{1,n}\ge\cdots\ge\theta_{n,n}$ are the eigenvalues of $T_{n,V,\beta}$, then
\[
 \lambda_{k,n}=b_\star^{-1/3}\tau_\star^{-2/3}n^{2/3}
                          (E_\star-\theta_{k,n}).
\]
Thus $b_\star^{-1/3}\tau_\star^{-2/3}n^{2/3} (\theta_{1,n}-E_\star)$ converges to $-\Lambda_1$, whose law is the Tracy--Widom $\beta$ distribution \cite{ramirez2011beta}.

\subsubsection{The comparison with the differential operator}\label{A:sec:intro-operator-embeddings}
The upper bound uses a difference operator whose potential is formed from the increments of the Brownian motion driving $\mathcal H_\beta$. This gives an intermediate operator for
both the coefficient coupling and the approximation of the differential operator. For $N\ge1$, impose $u_0=u_{N+1}=0$ and define
\[
 (\mathcal B_{N,h,\beta}^Bu)_k =h^{-2}(2u_k-u_{k-1}-u_{k+1})+(x_k+\omega_{k,h})u_k,\qquad
 \omega_{k,h} =c_\beta\frac{B(kh)-B((k-1)h)}h.
 \eqtag{brownian-difference}
\]
In the later form notation, $\mathcal B_{N,h,\beta}^B=-\Delta_{h,D}+X_h+M_{\omega_h}$. For fixed $h$, the noise coefficients are independent Gaussians with variance $4/(\beta h)$.

Let $h_n,N_n$ be deterministic sequences with $h_n\asymp n^{-1/3}$ and $N_nh_n\asymp h_n^{-2}$. For a continuous function $b$ on an interval $I$, write
$\operatorname{osc}_I b=\sup_{s,t\in I}|b(s)-b(t)|$. Let $\mathbb W$ be Wiener measure on $C_0(\mathbb R_+)=\{b\in C(\mathbb R_+):b(0)=0\}$. For $\alpha>0$, put
\[
 \mathcal O_\alpha(b)=\sup_{k\in\mathbb N_0}(1+k)^{-\alpha}
 \bigl(\operatorname{osc}_{[k,k+1]}b+|b(k+1)-b(k)|\bigr).
\]
For $r\in\mathbb N_0$, define
\[
\begin{aligned}
d_r& =2^{-r},\qquad T_r(b)=\sum_{k\ge0}(1+kd_r)^{-3/2} \operatorname{osc}_{[kd_r,(k+1)d_r]}(b)^2,\\
Q(b)& =\sup_{r\in\mathbb N_0}T_r(b),\qquad \Omega_B =\{b:Q(b)<\infty\}\cap \bigcap_{\alpha\in\mathbb Q\cap(0,1/2)} \{b:\mathcal O_\alpha(b)<\infty\}.
\end{aligned}
\]
Oscillations of continuous functions can be computed at rational points; hence these functions are Borel measurable and $\Omega_B$ is a Borel set.

Lemma~\ref{A:res:A.6} proves $\mathbb W(\Omega_B)=1$. On this set the form \eqref{A:eq:airy-form} defines $\mathcal H_\beta$. The set is independent of the deterministic spacings
and boundary points in the next theorem.
\begin{theorem}\label{A:res:1.3}
For every $B\in\Omega_B$, every compact $K\Subset\mathbb C\setminus\mathbb R$, and all sufficiently large $n$,
\[
 \sup_{z\in K}
 \left\|\iota_{h_n}^{(N_n)}
       (\mathcal B_{N_n,h_n,\beta}^B-z)^{-1}
       (\iota_{h_n}^{(N_n)})^*-R_\beta(z)\right\|_{\mathfrak S_2}^2
 \le C_{B,\beta,K}h_n.
 \eqtag{1.4}
\]
\end{theorem}
The constant may depend on the constants in the two assumptions on $h_n,N_n$. This comparison must hold on a growing interval: its length is $L=(N_n+1)h_n\asymp h_n^{-2}$.
Truncation gives squared error $C_{B,K}L^{-1/2}$, and approximation of the integrated potential followed by restriction of the form to piecewise linear functions gives
$C_{B,\beta,K}h_n$. Thus the growing interval puts both errors at the required scale; the proof is in
Subsection~\ref{A:sec:finite-difference-approximation-of-the-stochastic-airy-operator}.

For the Gaussian model, a rotation of independent normal variables followed by a chi quantile map couples the scalar coefficients to these increments; see
Subsection~\ref{A:sec:coupling-the-dumitriu-edelman-fluctuations-to-brownian-increments}. For the general Jacobi model, conditional Gaussian approximation on a regular interval
gives coordinate $L^p$ error $O(n^{-1})$, uniformly over the admissible exterior values (Proposition~\ref{A:res:A.14A}). The kernel in Lemma~\ref{A:res:A.11D} gives the exact
conditional coefficient law. Integrating against the original exterior marginal gives the full matrix law.

The covariance of the Gaussian approximation is the inverse of the Hessian at the conditional minimizer. The variance of long partial sums of its scalar coefficients determines the
normalization of the Brownian term. Lemma~\ref{A:res:15.1} computes this quantity and gives the normalization $4/\beta$. The factorization in Theorem~\ref{A:res:15.3} decomposes
the correlated Gaussian field into Brownian increments, a discrete derivative, and remainders. The weighted resolvent estimates treat these terms differently: summation by parts
transfers the derivative to the test functions, while the remaining diagonal error must satisfy a weighted square sum estimate. This distinction determines the accuracy required of
the coupling. Section~\ref{A:sec:outline-of-the-proof} gives the identities and the estimates at their respective scales.

To pass from nonreal resolvents to the inverse, we prove atomlessness at zero and then use a rational resolvent identity. In the general model this identity applies first to the
isometrically embedded principal block; comparisons using Schur complements then give the bound for the full inverse.
Subsection~\ref{A:sec:atomlessness-real-resolvents-and-inverse-coupling} gives both comparisons.

\subsubsection{Related results}\label{A:sec:intro-related-results}
The tridiagonal convergence criterion of Ram\'irez, Rider, and Vir\'ag \cite{ramirez2011beta} derives the low lying spectral limit from convergence of the integrated potential and
suitable form bounds. Krishnapur, Rider, and Vir\'ag \cite{krishnapur2016universality} prove a coupled operator limit for uniformly convex polynomial Jacobi laws, stated in norm
resolvent topology. These works establish qualitative operator limits on a common probability space.

For the inverse itself, Li, Valk\'o, and Xu \cite{li2026laguerre} (Theorem~2) prove Hilbert--Schmidt convergence in probability for the Gaussian tridiagonal model. Their lower edge
normalization is equivalent in law to the Gaussian upper edge normalization used here. Theorems~\ref{A:res:1.1} and~\ref{A:res:1.2} quantify this inverse comparison: they give an
almost sure squared error of order $n^{-1/3}$, locally uniformly for nonreal resolvents, and treat the polynomial Jacobi laws under the stated regular edge assumption. The lower
bound holds for every approximation of rank at most $n$.

An earlier quantitative operator coupling is due to Valk\'o and Vir\'ag \cite{valko2020operator}. For the circular $\beta$-ensemble and the Sine$_\beta$ operator, their conjugated
inverses satisfy
\[
 \|\mathtt r\mathrm{Circ}_{\beta,n}
              -\mathtt r\mathrm{Sine}_\beta\|_{\mathfrak S_2}^2
       \le\frac{(\log n)^6}{n}
\]
almost surely for all sufficiently large $n$. Their symbol $\mathtt r$ denotes the inverse after unitary identification with an operator on unweighted $L^2$. The circular operator
has the periodically extended finite spectrum, and its inverse has infinite rank. The lower bound in Subsection~\ref{A:sec:rank-lower-bound-and-optimality} applies to
approximations of finite rank.

The other parts of this work concern the hard edge, the circular bulk, and Wigner matrices. In Section~\ref{sec:hard-edge-models-brownian-information-and-bessel-limits}, the
squared inverse error has order $n^{-1}$ for beta-Laguerre and polynomial bidiagonal models in the stated parameter ranges; the matching lower bound requires access only to
Brownian observations at the model's deterministic sampling times and auxiliary randomness independent of the Brownian path. In
Section~\ref{sec:circular-ensembles-canonical-systems-and-the-sine-coefficient}, $n$ times the squared circular inverse error converges almost surely to an explicit positive
functional of the Brownian path. The result for Wigner matrices in Section~\ref{sec:wigner-matrices-and-the-airy-limit-for-a-deterministic-starting-vector} establishes an operator
limit for models whose tridiagonal coefficient laws have no explicit formula. The proof of the soft edge approximation rate is independent of these companion results. For the
broader operator approach to random matrices, see \cite{virag2014operator}.

Section~\ref{A:sec:outline-of-the-proof} identifies the analytic estimates and coefficient bounds needed for the upper rate.
Section~\ref{A:sec:inverse-coupling-and-spectral-consequences} proves the inverse and resolvent theorems. The spectral consequences are established in
Subsection~\ref{A:sec:spectral-consequences}, and Section~\ref{A:sec:parameter-dependence-and-further-directions} treats uniform moment bounds and converging parameter sequences.

\subsection{Hard edge models and coupling theorems for their inverses}
\label{sec:hard-edge-models-brownian-information-and-bessel-limits}
\label{B:sec:reorganized-1}

At the hard edge we use the inverse of the bidiagonal matrix. Its product formula gives an integral kernel, through which we compare the coefficients with those of the limiting
differential operator. The bidiagonal representation of Dumitriu and Edelman \cite{dumitriu2002matrix} and the stochastic operator description of Edelman and Sutton
\cite{edelman2007random} provide the starting point for this comparison.

For covariance models, the relevant limit is at the lower spectral edge. Ram\'irez and Rider \cite{ramirez2009diffusion,ramirez2011erratum} identified the stochastic Bessel
operator as the limit of beta-Laguerre matrices. Rider and Waters \cite{rider2019universality} extended this operator limit to polynomial potentials under a convexity assumption.
We ask how accurately a whole sequence of such matrices can be coupled to one realization of the limiting operator.

For every fixed $\beta>0$ and $a>0$, our coupling preserves the coefficient law at each size and attains squared Hilbert--Schmidt error of order $n^{-1}$ for the positive inverses,
almost surely. The estimate recovers eigenvalue and eigenspace convergence with a rate. It also controls quantities defined without reference to the chosen embedding: the complete
sequence of reciprocal eigenvalues in $\ell^2$, traces of inverse powers of order at least two, and regularized determinants together with their derivatives. These consequences
compare the finite matrix with the limiting operator beyond any fixed collection of eigenvalues; their precise statements appear in
Section~\ref{B:sec:spectral-consequences-and-further-questions}.

The lower bound reflects the information available to the coupling. At size $n$, the matrix depends on Brownian values at specified deterministic logarithmic times and on
additional randomness independent of the whole path. For each fixed sequence of approximations using those observations, the bridges between successive times give the matching
lower bound almost surely. The statement and proof appear in Theorem~\ref{B:thm:brownian-information-lower-bound}. To attain that order while preserving the matrix law, we control
partial sums of logarithmic coefficient ratios in the weighted norm dictated by the inverse kernel. Subsection~\ref{B:subsection:intro-green-kernel-comparison} explains this
reduction.

\subsubsection{The finite beta-Laguerre matrices}
\label{B:subsection:the-finite-matrix-models}
We use the following sign and indexing convention for the bidiagonal matrices.

Fix $\beta>0$ and $a>-1$. The finite operators act on $\mathbb C^n$, and their embedded inverses act on $L^2(0,1)$.

For positive coefficients $x_1,\ldots,x_n,y_1,\ldots,y_{n-1}$, let $B=B(x,y)$ be lower bidiagonal, with
\begin{equation}
 B_{kk}=x_k,\qquad B_{k+1,k}=-y_k;
 \qquad B_{ij}=0\quad\text{otherwise}.
 \label{B:eq:1-3}
\end{equation}
Since $\det B=\prod_kx_k>0$, the matrix $BB^{\mathsf T}$ is positive definite: $v^{\mathsf T}BB^{\mathsf T}v=|B^{\mathsf T}v|^2>0$ for real $v\ne0$. In the beta-Laguerre model the
coefficients at size $n$ are independent, with distributions
\begin{equation}
 \sqrt{n\beta}\,X_{n,k}^{\mathrm L}\sim\chi_{\beta(k+a)}
       \quad(1\le k\le n),\qquad
 \sqrt{n\beta}\,Y_{n,k}^{\mathrm L}\sim\chi_{\beta k}
       \quad(1\le k<n).
 \label{B:eq:intro-laguerre-coefficients}
\end{equation}
Here $\chi_\nu^2\sim\Gamma(\nu/2,2)$, with shape and scale in that order. Dependence between different sizes is specified by the coupling. Put
\begin{equation}
 B_{n,\beta,a}^{\mathrm L}=B(X_n^{\mathrm L},Y_n^{\mathrm L}),\qquad
 \mathcal A_{n,\beta,a}^{\mathrm L}
       =n^2B_{n,\beta,a}^{\mathrm L}(B_{n,\beta,a}^{\mathrm L})^{\mathsf T}.
 \label{B:eq:1-4}
\end{equation}
This is the beta-Laguerre representation of Dumitriu--Edelman \cite{dumitriu2002matrix}, with the order of the coefficients and scaling used in \cite{ramirez2009diffusion}.

The comparison space is fixed by the isometry
\begin{equation}
 \mathscr E_ne_k=\sqrt n\,\mathbf1_{[(k-1)/n,k/n)}
 \label{B:eq:flat-cell-embedding}
\end{equation}
The indicators have disjoint supports of length $1/n$, so $\mathscr E_n^*\mathscr E_n=I$. Define the embedded inverse and resolvent by
\begin{equation}
 \begin{aligned}
 \mathcal R_{n,\beta,a}^{\mathrm L}
   &=\mathscr E_n(\mathcal A_{n,\beta,a}^{\mathrm L})^{-1}\mathscr E_n^*,\\
 \mathcal G_{n,\beta,a}^{\mathrm L}(z)
   &=\mathscr E_n(\mathcal A_{n,\beta,a}^{\mathrm L}-z)^{-1}\mathscr E_n^*.
 \end{aligned}
 \label{B:eq:1-6}
\end{equation}
\subsubsection{The limiting inverse and boundary conditions}
\label{B:subsection:the-stochastic-bessel-inverse}

The limiting inverse is the adjoint product of a triangular integral operator whose coefficients are functions of a Brownian path. We construct this integral operator before
identifying its differential inverse. The construction uses standard Brownian motion $W$ on $[0,\infty)$ and initially allows $\beta>0$ and $a>-1$. In logarithmic time,
\begin{equation}
 \mathcal Y(s)=\beta^{-1/2}W(-\log s),\qquad 0<s\le1.
 \label{B:eq:1-29}
\end{equation}
All integral operators below act on the complex Hilbert space $L^2(0,1)$. We integrate the first kernel variable: $(Tf)(t)=\int_0^1K_T(s,t)f(s)\,\dd s$. Define
\begin{equation}
 \begin{gathered}
 p^B(s)=s^{a/2}e^{\mathcal Y(s)},\qquad
 q^B(s)=s^{-a-1}e^{-2\mathcal Y(s)},\\
 K_{\beta,a}(s,t)
 =p^B(s)\sqrt{q^B(t)}\,\mathbf1_{\{s<t\}}
 =\frac1{\sqrt t}\left(\frac st\right)^{a/2}
     e^{\mathcal Y(s)-\mathcal Y(t)}\mathbf1_{\{s<t\}}.
 \end{gathered}
 \label{B:eq:intro-brownian-coefficients}
\end{equation}
This is the inverse factor in \cite{rider2019universality}, in logarithmic time. To compare with the stochastic integral used in that paper, let $b$ be Brownian motion on $[0,1]$
and set $W(u)=\int_{e^{-u}}^1s^{-1/2}\,\dd b_s$. Its covariance is
\[
 \int_{e^{-\min\{u,v\}}}^1\frac{\dd s}{s}=\min\{u,v\}.
\]
The Gaussian fourth moment of an increment is $3|u-v|^2$. Thus the process has a continuous version which is Brownian motion, and
$\mathcal Y(s)-\mathcal Y(t)=\beta^{-1/2}\int_s^t r^{-1/2}\,\dd b_r$.

Let $\mathcal K_{\beta,a}$ be the integral operator with this kernel. The sublinear growth estimate for the Brownian path and its application in Proposition~\ref{B:result:2-5} give
$\mathcal K_{\beta,a}\in\mathfrak S_2$ almost surely for every fixed $a>-1$. Define
\begin{equation}
 \mathcal R_{\beta,a}=\mathcal K_{\beta,a}^*\mathcal K_{\beta,a},\qquad
 \mathcal H_{\beta,a}=\mathcal R_{\beta,a}^{-1},\qquad
 \operatorname{dom}\mathcal H_{\beta,a}
       =\operatorname{ran}\mathcal R_{\beta,a}.
 \label{B:eq:1-1}
\end{equation}
Proposition~\ref{B:result:2-5} proves that $\mathcal R_{\beta,a}$ is positive, injective and trace class, with dense range. Hence its inverse on the displayed domain is a positive
self-adjoint operator. The kernel used in the comparison is
\begin{equation}
 R_{\beta,a}(s,t)=p^B(s)p^B(t)\int_{s\vee t}^1q^B(r)\,\dd r.
 \label{B:eq:1-2}
\end{equation}
Its differential domain can be stated without differentiating $p^B$ or $q^B$. For $f\in L^2(0,1)$ put $y_f=f/p^B$ and $v_f=(q^B)^{-1}y_f'$, with the derivative defined when $y_f$
is locally absolutely continuous. Then
\[
 \begin{split}
 \operatorname{dom}\mathcal H_{\beta,a}
 =\{f\in L^2(0,1):{}&y_f,v_f\in\operatorname{AC}([\epsilon,1])
                       \quad(0<\epsilon<1),\\
                 &(p^B)^{-1}v_f'\in L^2(0,1),\quad
                   y_f(1)=0,\quad v_f(0+)=0\}.
 \end{split}
\]
On this domain $\mathcal H_{\beta,a}f=-(p^B)^{-1}v_f'$. The derivatives and the identity defining $v_f$ hold almost everywhere. The boundary value at zero is the trace of $v_f$. At
one, $p^B(1)=1$, so $y_f(1)=0$ is equivalent to $f(1)=0$. Theorem~\ref{B:result:2-5a} proves that this domain equals $\operatorname{ran}\mathcal R_{\beta,a}$; the generic
coefficient version is \eqref{B:eq:bessel-domain-conditions}. Its eigenvalues, denoted by $0<\Lambda_1<\Lambda_2<\cdots$, are simple and tend to infinity by
Corollary~\ref{B:cor:bessel-simple-spectrum}. We write $\mathcal G_{\beta,a}(z)=(\mathcal H_{\beta,a}-z)^{-1}$ for $z\in\rho(\mathcal H_{\beta,a})$.

For a deterministic partition $0=z_{n,0}<z_{n,1}<\cdots<z_{n,n}=1$, write
\begin{equation}
 \mathcal N_n=\sigma\{W(-\log z_{n,k}):1\le k\le n\}.
 \label{B:eq:intro-observation-field}
\end{equation}
The matrix of size $n$ will be measurable with respect to $\mathcal N_n\vee\mathcal U$, where $\mathcal U$ is independent of $\sigma\{W(u):u\ge0\}$. One choice of $W$ and
$\mathcal U$ serves all sizes; the partition depends on the model. The index starts at $k=1$, so every observed time is finite. The upper bound below assumes $a>0$; the operator
construction and the lower bound have the range $a>-1$.

\subsubsection{The beta-Laguerre coupling}
\label{B:subsection:the-beta-laguerre-coupling-theorem}

The lower bound depends on the information used to construct the approximation. We therefore specify the joint distribution of the matrix coefficients and the Brownian
observations, as well as the finite matrix distribution. The matrices and embedded inverses are those of Section~\ref{B:subsection:the-finite-matrix-models}, and the partition used
for observation in \eqref{B:eq:intro-observation-field} is $z_{n,k}=k/n$. A row denotes the whole coefficient vector $(X_n^{\mathrm L},Y_n^{\mathrm L})$; the auxiliary sigma-field
$\mathcal U$ is independent of the entire Brownian path $W$.

\begin{theorem}[Beta-Laguerre ensemble]\label{B:result:1-1}
Fix $\beta>0$ and $a>0$. There is a coupling of the matrices in \eqref{B:eq:intro-laguerre-coefficients} with $W$ and an independent auxiliary sigma-field $\mathcal U$, such that
each row is $\mathcal N_n\vee\mathcal U$-measurable and, almost surely,
\begin{equation}
 0<\liminf_{n\to\infty}
       n\|\mathcal R_{n,\beta,a}^{\mathrm L}-\mathcal H_{\beta,a}^{-1}\|_{\mathfrak S_2}^2
 \le\sup_{n\ge1}
       n\|\mathcal R_{n,\beta,a}^{\mathrm L}-\mathcal H_{\beta,a}^{-1}\|_{\mathfrak S_2}^2
 <\infty.
 \label{B:eq:1-33}
\end{equation}
\end{theorem}

Section~\ref{B:sec:the-beta-laguerre-coupling} constructs the coupling and proves the theorem. Quantiles couple each logarithmic chi ratio to the corresponding Brownian increment;
an independent gamma variable then recovers the pair of coefficients with the joint law in \eqref{B:eq:intro-laguerre-coefficients}. The common Brownian path and the auxiliary
variables in this construction determine the dependence between matrix sizes.

\subsubsection{Polynomial potentials}
\label{B:subsection:the-rider-waters-general-potential-model}
\label{B:subsection:profile-normalization-and-embeddings}
\label{B:subsection:the-rider-waters-embedding}
\label{B:subsection:the-general-potential-coupling-theorem}

A general potential changes both the coefficient distribution and its deterministic profile. We specify the change of variable and the matrix scale before stating the coupling
theorem. Let $V(\lambda)=\sum_{m=1}^dg_m\lambda^m$ be a real polynomial satisfying
\begin{equation}
 \frac{\dd^2}{\dd x^2}V(x^2)\ge c_*>0\quad(x\in\mathbb R),
 \qquad c_0=2V'(0)>0.
 \label{B:eq:1-8}
\end{equation}
The coefficient distribution is
\begin{equation}
 \begin{split}
 \dd\mathbb P_{n,V,\beta,a}(x,y)
 ={}&\frac1{Z_{n,V,\beta,a}}
       e^{-n\beta\operatorname{Tr}V(BB^{\mathsf T})}
       \prod_{k=1}^n x_k^{\beta(k+a)-1}
       \prod_{k=1}^{n-1}y_k^{\beta k-1}\,\dd x\,\dd y,
 \end{split}
 \label{B:eq:1-10}
\end{equation}
on $(0,\infty)^{2n-1}$, as in \cite{rider2019universality}. Its normalizing constant is finite and positive for $a>-1$. For $V(\lambda)=\lambda/2$ it is exactly the independent law
\eqref{B:eq:intro-laguerre-coefficients}. Both calculations are given with the coefficient density in Section~\ref{B:subsection:the-polynomial-hamiltonian}. Write $B_n=B(X_n,Y_n)$
when $(X_n,Y_n)$ has distribution $\mathbb P_{n,V,\beta,a}$.

The stationary equations at constant coefficients determine the profile $\phi:[0,1]\to[0,\infty)$ through
\begin{equation}
 t=\sum_{m=1}^d m\binom{2m}{m}g_m\phi(t)^{2m}.
 \label{B:eq:1-12}
\end{equation}
Lemma~\ref{B:lem:bessel-coordinate-regularity} proves uniqueness and shows that $\phi(t)/\sqrt t$ extends to a positive smooth function on $[0,1]$. In particular,
\begin{equation}
 \phi(t)=c_0^{-1/2}\sqrt t\,(1+O(t)),\qquad t\downarrow0.
 \label{B:eq:1-13}
\end{equation}
Set
\begin{equation}
 L_V(t)=\int_0^t\frac{\dd s}{\phi(s)},\qquad
 L_V=L_V(1),\qquad
 \theta_V(t)=\left(\frac{L_V(t)}{L_V}\right)^2.
 \label{B:eq:1-15}
\end{equation}
The coordinate $\theta_V$ is the one used in \cite{rider2019universality}. Differentiating its definition gives
\begin{equation}
 \theta_V'(t)=\frac{2L_V(t)}{L_V^2\phi(t)}
             =\frac{2\sqrt{\theta_V(t)}}{L_V\phi(t)},
 \qquad \theta_V'(t)^2\phi(t)^2=\frac4{L_V^2}\theta_V(t).
 \label{B:eq:1-16}
\end{equation}
The same lemma gives $\theta_V\in C^2([0,1])$ and $0<c\le\theta_V'\le C$, $|\theta_V''|\le C$. Integrating \eqref{B:eq:1-13} in \eqref{B:eq:1-15} yields
\begin{equation}
 L_V(t)=2\sqrt{c_0t}\,(1+O(t)),\qquad
 \theta_V(t)=d_Vt(1+O(t)),\qquad
 d_V=\theta_V'(0)=\frac{4c_0}{L_V^2}.
 \label{B:eq:1-18}
\end{equation}

To place the operator in this coordinate, define
\begin{equation}
 \begin{aligned}
 (U_{\theta_V}f)(t)&=\sqrt{\theta_V'(t)}f(\theta_V(t)).
 \end{aligned}
 \label{B:eq:1-22}
\end{equation}
The substitution $z=\theta_V(t)$ gives
\[
 \int_0^1|U_{\theta_V}f(t)|^2\,\dd t
       =\int_0^1|f(z)|^2\,\dd z,
 \qquad
 (U_{\theta_V}^{-1}g)(z)
 =\frac{g(\theta_V^{-1}(z))}
        {\sqrt{\theta_V'(\theta_V^{-1}(z))}}.
\]
Thus $U_{\theta_V}$ is unitary. Define the isometry
\begin{equation}
 \mathscr J_{n,V}=U_{\theta_V}^{-1}\mathscr E_n,
 \qquad
 (\mathscr J_{n,V}e_k)(z)
 =\frac{\sqrt n\,\mathbf1_{[\theta_V((k-1)/n),\theta_V(k/n))}(z)}
        {\sqrt{\theta_V'(\theta_V^{-1}(z))}}.
 \label{B:eq:1-24}
\end{equation}

The same coordinate change fixes the matrix scale. The factor kernel in \cite{rider2019universality} is
\[
 K_V(s,t)=\frac1{\sqrt{\phi(s)\phi(t)}}
    \left(\frac{\theta_V(s)}{\theta_V(t)}\right)^{a/2+1/4}
    e^{\mathcal Y(\theta_V(s))-\mathcal Y(\theta_V(t))}
    \mathbf1_{\{s<t\}}.
\]
The substitution in the preceding unitary formula and \eqref{B:eq:1-16} give, as an identity of integral kernels,
\begin{align*}
 K_{U_{\theta_V}\mathcal K_{\beta,a}U_{\theta_V}^{-1}}(s,t)
 &=\sqrt{\theta_V'(s)\theta_V'(t)}
                 K_{\beta,a}(\theta_V(s),\theta_V(t))\\
 &=\frac2{L_V\sqrt{\phi(s)\phi(t)}}
   \frac{(\theta_V(s)\theta_V(t))^{1/4}}{\theta_V(t)^{1/2}}
   \left(\frac{\theta_V(s)}{\theta_V(t)}\right)^{a/2}
   e^{\mathcal Y(\theta_V(s))-\mathcal Y(\theta_V(t))}
   \mathbf1_{\{s<t\}}\\
 &=\frac2{L_V}K_V(s,t).
\end{align*}
Consequently, if $\mathcal K_V$ denotes the operator with kernel $K_V$, then
\begin{equation}
 \mathcal H_{\beta,a}^{-1}
 =\frac4{L_V^2}U_{\theta_V}^{-1}
                     \mathcal K_V^*\mathcal K_VU_{\theta_V}.
 \label{B:eq:intro-factor-conjugation}
\end{equation}
The factor multiplier is $2/L_V$ and the inverse multiplier is $4/L_V^2$. Accordingly the finite matrix is scaled by
\begin{equation}
 c_V=\frac{L_V^2}{4},\qquad
 \mathcal A_{n,V,\beta,a}=c_Vn^2B_nB_n^{\mathsf T}.
 \label{B:eq:1-20}
\end{equation}
For $V(\lambda)=\lambda/2$, the formulas give
\begin{equation}
 \phi(t)=\sqrt t,\quad L_V(t)=2\sqrt t,\quad L_V=2,
 \quad\theta_V(t)=t,\quad c_V=1,\quad\mathscr J_{n,V}=\mathscr E_n.
 \label{B:eq:1-21}
\end{equation}

For this model set
\begin{equation}
 \begin{aligned}
 \mathcal R_{n,V,\beta,a}
 &=\mathscr J_{n,V}\mathcal A_{n,V,\beta,a}^{-1}\mathscr J_{n,V}^*,\\
 \mathcal G_{n,V,\beta,a}(z)
 &=\mathscr J_{n,V}(\mathcal A_{n,V,\beta,a}-z)^{-1}\mathscr J_{n,V}^*.
 \end{aligned}
 \label{B:eq:1-26}
\end{equation}
The two choices of Brownian observation times are thus
\begin{equation}
 u_{n,k}^{\mathrm L}=-\log(k/n),\qquad
 u_{n,k}^{V}=-\log\theta_V(k/n),\qquad 1\le k\le n.
 \label{B:eq:1-30}
\end{equation}
Use $z_{n,k}=\theta_V(k/n)$ in \eqref{B:eq:intro-observation-field}. The potential $V$, as well as $\beta$ and $a$, is fixed in the next theorem.

\begin{theorem}[Polynomial potentials]\label{B:result:1-3}
Fix $\beta>0$, $a>0$, and a polynomial $V$ satisfying \eqref{B:eq:1-8}. There is a coupling of matrices with row distributions $\mathbb P_{n,V,\beta,a}$ with $W$ and an independent
auxiliary sigma-field $\mathcal U$, such that each row is $\mathcal N_n\vee\mathcal U$-measurable and, almost surely,
\begin{equation}
 0<\liminf_{n\to\infty}
       n\|\mathcal R_{n,V,\beta,a}-\mathcal H_{\beta,a}^{-1}\|_{\mathfrak S_2}^2
 \le\sup_{n\ge1}
       n\|\mathcal R_{n,V,\beta,a}-\mathcal H_{\beta,a}^{-1}\|_{\mathfrak S_2}^2
 <\infty.
 \label{B:eq:1-38}
\end{equation}
\end{theorem}
The proof in Section~\ref{B:subsection:proof-general-potential-theorem} uses the conditional coefficient distributions constructed in
Sections~\ref{B:sec:general-potentials-and-convex-reference-measures}--\ref{B:sec:construction-of-the-coupling-for-all-matrix-sizes}. At $V(\lambda)=\lambda/2$ the marginal
coefficient laws and the comparison maps agree with those of Theorem~\ref{B:result:1-1}. The conditional construction specifies its own joint law with $W$.

\subsubsection{Resolvents and spectral consequences}
\label{B:subsection:intro-resolvents-and-tightness}
\label{B:subsection:the-inverse-formulation}
\label{B:subsection:spectral-consequences}

The resolvent is a rational function of the positive inverse. The finite-dimensional formula also includes the isometry into the common Hilbert space.

The notation $O_\omega(a_n)$ means a bound by a finite random constant times $a_n$, for all sufficiently large $n$. For either coupling abbreviate its embedded inverse by $R_n$,
put $R=\mathcal H_{\beta,a}^{-1}$, and write $G_n(z),G(z)$ for the corresponding embedded finite and limiting resolvents. Set $e_n=\|R_n-R\|_{\mathfrak S_2}$ and
$Z=\sup_{n\ge1}ne_n^2$. The preceding theorems give $Z<\infty$ and $\liminf_n\sqrt n\,e_n>0$ almost surely. Thus the norm error has order $n^{-1/2}$.
Section~\ref{B:subsection:inverse-identities-and-embedded-resolvents} proves
\begin{equation}
 \begin{gathered}
 G_n(z)=R_n(I-zR_n)^{-1},\qquad G(z)=R(I-zR)^{-1},\\
 G_n(z)-G(z)=(I-zR_n)^{-1}(R_n-R)(I-zR)^{-1}.
 \end{gathered}
 \label{B:eq:intro-resolvent-identity}
\end{equation}
On the orthogonal complement of the finite-dimensional comparison space, $R_n$ and $G_n(z)$ are zero and $I-zR_n$ is the identity. The next two corollaries are proved in
Section~\ref{B:subsection:resolvent-rates-and-tightness}.

\begin{corollary}\label{B:cor:intro-resolvent-rates}
On the event of either coupling theorem, every nonempty compact $K\subset\rho(\mathcal H_{\beta,a})$ lies in the finite resolvent set for all sufficiently large $n$, and
\begin{equation}
 \frac{c_{\omega,K}}n
 \le\inf_{z\in K}\|G_n(z)-G(z)\|_{\mathfrak S_2}^2
 \le\sup_{z\in K}\|G_n(z)-G(z)\|_{\mathfrak S_2}^2
 \le\frac{C_{\omega,K}}n
 \label{B:eq:intro-resolvent-rates}
\end{equation}
for positive finite constants and all sufficiently large $n$.
\end{corollary}

\begin{corollary}\label{B:cor:intro-tightness}
The sequence $(ne_n^2)$ is tight. For each deterministic nonempty compact $K\subset\mathbb C\setminus(0,\infty)$,
\begin{equation}
 \sup_{n\ge1}n\sup_{z\in K}\|G_n(z)-G(z)\|_{\mathfrak S_2}^2<\infty
 \quad\text{almost surely};
 \label{B:eq:intro-resolvent-supremum}
\end{equation}
in particular the normalized squared resolvent errors are tight.
\end{corollary}

Use $\mathfrak m\in\{\mathrm L,V\}$ to distinguish the two models. Let $\lambda_{j,n}^{\mathfrak m}$ be the increasing eigenvalues of its scaled matrix. The symbols
$\Pi_{j,n}^{\mathfrak m}$ and $\Pi_j$ denote the corresponding embedded finite and limiting rank-one spectral projections. Explicitly, for a real normalized eigenvector $v_{j,n}$
of $\mathcal A_n^{\mathfrak m}$ and eigenfunction $\varphi_j$ of $\mathcal H_{\beta,a}$, put $u_{j,n}=\mathscr J_n^{\mathfrak m}v_{j,n}$. Then
$\Pi_{j,n}^{\mathfrak m}f=\langle f,u_{j,n}\rangle u_{j,n}$ and $\Pi_jf=\langle f,\varphi_j\rangle\varphi_j$, where $\mathscr J_n^{\mathrm L}=\mathscr E_n$ and
$\mathscr J_n^V=\mathscr J_{n,V}$. Section~\ref{B:sec:spectral-consequences-and-further-questions} constructs measurable normalized choices. The three corollaries below are deduced
in Section~\ref{B:subsection:bessel-spectral-probability-rates}, after the estimates for reciprocal eigenvalues, spectral projections, and determinants.

\begin{corollary}\label{B:cor:intro-eigenvalue-rates}
For either coupling, almost surely,
\begin{equation}
 \sum_{j=1}^n\left|\frac1{\lambda_{j,n}^{\mathfrak m}}-\frac1{\Lambda_j}\right|^2
       +\sum_{j>n}\Lambda_j^{-2}\le e_n^2\le Z/n.
 \label{B:eq:1-60}
\end{equation}
For every fixed $m$, choose $n_0(\omega,m)\ge m$ such that $e_n\le(2\Lambda_m)^{-1}$ for $n\ge n_0(\omega,m)$. For these
$n$,
\begin{equation}
 \max_{j\le m}|\lambda_{j,n}^{\mathfrak m}-\Lambda_j|
       \le2\Lambda_m^2e_n=O_\omega(n^{-1/2}).
 \label{B:eq:1-61}
\end{equation}
The errors for fixed indices are also $O_{\mathbb P}(n^{-1/2})$.
\end{corollary}

\begin{corollary}\label{B:cor:intro-projection-rates}
For each fixed $j$,
\begin{equation}
 \|\Pi_{j,n}^{\mathfrak m}-\Pi_j\|_{\mathfrak S_2}
       =O_\omega(n^{-1/2}),\qquad
  \|\Pi_{j,n}^{\mathfrak m}-\Pi_j\|_{\mathfrak S_2}=O_{\mathbb P}(n^{-1/2}).
 \label{B:eq:1-62}
\end{equation}
The same rates hold for the corresponding real normalized eigenfunctions in $L^2$, after minimizing over their sign.
\end{corollary}

The regularized determinants in the next corollary are
\[
 \begin{aligned}
 D_n^{\mathfrak m}(z)&=\det{}_2(I-zR_n)
   =\prod_{j=1}^n(1-z/\lambda_{j,n}^{\mathfrak m})
                          \exp(z/\lambda_{j,n}^{\mathfrak m}),\\
 D(z)&=\det{}_2(I-zR)
   =\prod_{j=1}^\infty(1-z/\Lambda_j)\exp(z/\Lambda_j).
 \end{aligned}
\]
Section~\ref{B:subsection:regularized-determinants} proves local uniform convergence of the infinite product and describes its zeros.

\begin{corollary}\label{B:cor:intro-determinant-rates}
For every deterministic nonempty compact $K\subset\mathbb C$,
\begin{equation}
 \sup_{z\in K}|D_n^{\mathfrak m}(z)-D(z)|
       =O_\omega(n^{-1/2}).
 \label{B:eq:1-63}
\end{equation}
The normalized errors are tight. Both conclusions hold for every fixed derivative order.
\end{corollary}

\subsubsection{The lower bound}
\label{B:subsection:information-used-by-the-coupling}
\label{B:subsection:parameter-range}

The lower bound depends on which values of the Brownian path are observed at each size. For either partition in \eqref{B:eq:1-30}, let
\begin{equation}
 \mathcal F_n=\mathcal N_n\vee\mathcal U_n,\qquad
 \mathcal U_n\ \text{independent of }\sigma(W(u):u\ge0).
 \label{B:eq:1-32}
\end{equation}
The auxiliary fields may be dependent for different $n$. For every sequence of $\mathcal F_n$-measurable Hilbert--Schmidt operators $T_n$,
Theorem~\ref{B:thm:brownian-information-lower-bound} gives
\begin{equation}
 \liminf_{n\to\infty}
 n\|T_n-\mathcal H_{\beta,a}^{-1}\|_{\mathfrak S_2}^2>0
 \quad\text{almost surely}.
 \label{B:eq:intro-observation-lower-bound}
\end{equation}
The event may depend on the sequence. This result holds for $a>-1$; $T_n$ may have any distribution subject to the stated measurability. The observation field in this lower bound
is fixed separately at each size. On continuous path space, for the uniform partitions,
\[
 \bigvee_{m\ge1}\mathcal N_m=\sigma\{W(u):u\ge0\}.
\]
Indeed, for each $u\ge0$ there are integers $1\le k_j\le m_j$ with $k_j/m_j\to e^{-u}$. Then $-\log(k_j/m_j)\to u$ and $W(-\log(k_j/m_j))\to W(u)$. Thus using the observations from
all sizes gives access to the complete path. The lower bound uses $\mathcal N_n\vee\mathcal U_n$ at size $n$.

The scale follows from the conditional variance between consecutive observations. On a time interval of length $h$, the Brownian bridge $B_h$ has covariance
$\mathbb E B_h(s)B_h(t)=s\wedge t-st/h$. Consequently
\[
 \mathbb E\int_0^h B_h(t)^2\,\dd t
 =\int_0^h(t-t^2/h)\,\dd t=\frac{h^2}{6}.
\]
On a fixed bounded interval of logarithmic time there are order $n$ cells of length order $n^{-1}$. Their total mean squared bridge fluctuation is therefore of order $n^{-1}$. To
obtain \eqref{B:eq:intro-observation-lower-bound} for the Green kernel, Section~\ref{B:subsection:brownian-information-lower-bound} conditions its linearization on the Brownian
observations. The conditional covariance bound and the estimate for the nonlinear remainder are proved separately in Propositions~\ref{B:prop:conditional-covariance-bounds}
and~\ref{B:prop:bridge-remainder-bound}. These estimates give summable failure probabilities; Borel--Cantelli proves that the lower bound holds almost surely.

The unrestricted spectral truncation has smaller error. Let $R^{[n]}$ be the rank-$n$ spectral truncation of $R=\mathcal H_{\beta,a}^{-1}$, and put $\mu_j=\Lambda_j^{-1}$. Since
$R$ is trace class,
\[
 \begin{aligned}
 n\|R-R^{[n]}\|_{\mathfrak S_2}^2
 &=n\sum_{j>n}\mu_j^2
 \le n\mu_{n+1}\sum_{j>n}\mu_j\\
 &\le\operatorname{Tr}R\sum_{j>n}\mu_j\longrightarrow0.
 \end{aligned}
\]
Proposition~\ref{B:prop:bessel-best-finite-rank-approximation} identifies this spectral truncation as a best approximation among operators of rank at most $n$. Such an
approximation may use the whole Brownian path. Theorems~\ref{B:result:1-1} and~\ref{B:result:1-3} attain the lower bound within the class \eqref{B:eq:1-32}, while retaining the
original matrix distributions.

The upper bound uses $a>0$ in two places. The weighted kernel estimate requires $1<\lambda<\min\{a+1,2\}$. The rows near zero are bounded by $d^{1+\nu_a-\eta}$, where
$\nu_a=\min\{a,1\}$ and $0<\eta<\nu_a$. For $d_n\asymp\varepsilon_n\asymp(\log n)^A/n$, the three truncation terms in Theorem~\ref{B:thm:inverse-truncation-bound} obey
\[
 \begin{aligned}
 nd_n^{1+\nu_a-\eta}
 &=O\!\left(n^{-(\nu_a-\eta)}(\log n)^{A(1+\nu_a-\eta)}\right)\longrightarrow0,\\
 n\varepsilon_n^3&=O\!\left(n^{-2}(\log n)^{3A}\right)\longrightarrow0,\qquad
 n\varepsilon_n^2=O\!\left(n^{-1}(\log n)^{2A}\right)\longrightarrow0.
 \end{aligned}
\]
The estimates at the singular boundary therefore determine the range $a>0$ of the upper bound. The extension to $-1<a\le0$ is discussed in
Section~\ref{B:subsection:bessel-spectral-parameter-range}. All constants in the coupling theorems may depend on the fixed model parameters.

\subsection{Circular ensembles and the leading order coefficient}
\label{sec:circular-ensembles-canonical-systems-and-the-sine-coefficient}
\label{S:sec:introduction}
 For the circular $\beta$-ensemble, the Dirac operators of Valk\'o and Vir\'ag~\cite{valko2017sine,valko2020operator} give finite and limiting operators on a common Hilbert space.
 We determine the leading coefficient of their inverse approximation error.

For $\beta>0$, the circular $\beta$-ensemble is the distribution of $n$ points $\me^{\ii\theta_{j,n}}$ on the unit circle with density
\[
 \frac1{Z_{n,\beta}}\prod_{1\le j<k\le n}
       |\me^{\ii\theta_j}-\me^{\ii\theta_k}|^\beta
       \prod_{j=1}^n\frac{\dd\theta_j}{2\pi},
 \qquad 0\le\theta_j<2\pi.
\]
Its rescaled periodic extension
\begin{equation}\label{S:eq:intro-periodic-spectrum}
 \{n(\theta_{j,n}+2\pi m):1\le j\le n,\ m\in\mathbb Z\}
\end{equation}
converges locally in distribution to Sine$_\beta$, with intensity $1/(2\pi)$ \cite{KillipStoiciu2009,valko2009continuum,Nakano2014,valko2017sine}.

Valk\'o and Vir\'ag~\cite{valko2017sine,valko2020operator} associated random Dirac operators with these point processes. The finite operator is defined by a hyperbolic random walk
and the limiting operator by a hyperbolic Brownian path. After unitary conjugation to unweighted coordinates, their inverses act on the same Hilbert space. The distance between
these inverses controls reciprocal eigenvalues and isolated eigenspaces, as well as the convergence of the point processes.

We determine the leading squared Hilbert--Schmidt error under a coupling that preserves the exact circular operator law at every size. The coefficient is an explicit positive
functional of the limiting Brownian path. To find it, we first couple the circular inverse to the inverse obtained by sampling that path at the times defined below, with an error
smaller than $n^{-1/2}$. The leading error is therefore already present in sampling the Brownian path. We compute the contribution of the oscillations between successive sampling
times to the Green kernel and prove that the resulting squared norm is asymptotic to its conditional second moment.

\subsubsection{The finite and limiting operators}
\label{S:sec:intro-operator-representations}
Both the circular operator and its limit are defined by a path in the upper half-plane and a boundary direction. The matrix factor below identifies their inverses with operators on
a common unweighted space.

For a path $z(t)=x(t)+\ii y(t)$ in the upper half-plane, put
\[
 X(t)=\mathcal A(z(t)),\qquad
 \mathcal A(x+\ii y)=\frac1{\sqrt y}
       \begin{pmatrix}1&-x\\0&y\end{pmatrix},\qquad
 J=\begin{pmatrix}0&-1\\1&0\end{pmatrix}.
\]
We use the differential expression
\begin{equation}\label{S:eq:intro-canonical-expression}
 \mathcal D_Xf=2R^{-1}Jf',\qquad R=X^{\mathsf T}X,
\end{equation}
on $L^2_R=L^2([0,1),R(t)\,\dd t;\mathbb R^2)$. Given $q\in\mathbb R$, set $u_0=(-1,0)^{\mathsf T}$ and $u_1=(q,1)^{\mathsf T}$, so that $u_0^{\mathsf T}Ju_1=1$. The domain consists
of locally absolutely continuous functions $f$ such that $f,\mathcal D_Xf\in L^2_R$ and
\[
 f(0)^{\mathsf T}Ju_0=0,\qquad
 \lim_{t\uparrow1}f(t)^{\mathsf T}Ju_1=0.
\]
For the paths below, $\mathcal D_X$ with this domain is self-adjoint and has a Hilbert--Schmidt inverse. At a regular boundary the second condition is the usual boundary value. For
the Brownian path it is an additional condition when $\beta>2$ and is automatic on the maximal domain when $\beta\le2$. The operator domains and Green formulas are established in
Theorems~\ref{S:thm:3.2} and~\ref{S:thm:sine-selfadjoint-domain}.

Multiplication by $X$ is a unitary map $\mathcal U_X$ from $L^2_R$ to $\mathcal H=L^2([0,1);\mathbb R^2)$, since $|Xf|^2=f^{\mathsf T}Rf$. For the canonical systems we compare
$\mathcal U_X\mathcal D_X^{-1}\mathcal U_X^{-1}$ using the displayed factor $\mathcal A$. This choice of $\mathcal A$ fixes the unitary map for each canonical operator. In the
convention of~\cite{valko2017sine,valko2020operator}, the Hamiltonian is $R/2$, so the differential expression is the same as \eqref{S:eq:intro-canonical-expression}.

For the circular operator, start a hyperbolic walk at $b_0=\ii$. For $0\le j\le n-2$, take independent variables $\zeta_j\sim\operatorname{Beta}(1,\beta(n-j-1)/2)$ and move from
$b_j$ in an independent uniform direction through distance $2\operatorname{arctanh}\sqrt{\zeta_j}$. Conditional on the walk, choose $b_n$ from harmonic measure at $b_{n-1}$; at
$x+\ii y$ this measure has density $y/[\pi((q-x)^2+y^2)]$. The step path $z(t)=b_{\lfloor nt\rfloor}$ and boundary point $q=b_n$ define $\mathrm{Circ}_{\beta,n}$. For $n=1$ only
$b_1$ is sampled from harmonic measure at $b_0$. The circular operator representation of~\cite{valko2020operator}, with the normalization in
Proposition~\ref{S:prop:operator-normalization}, identifies its spectrum with \eqref{S:eq:intro-periodic-spectrum}. In particular, this differential operator has infinite periodic
spectrum, despite having only $n$ eigenvalues in one period.

Let $Z_\theta=x_\theta+\ii y_\theta$ be standard hyperbolic Brownian motion started at $\ii$, with generator $\tfrac12y^2(\partial_x^2+\partial_y^2)$. Define
\[
 \tau_\beta(t)=\frac4\beta\log\frac1{1-t},\qquad
 X_\beta(t)=\mathcal A(Z_{\tau_\beta(t)}),\qquad
 q_\infty=\lim_{\theta\to\infty}x_\theta.
\]
The terminal point is finite almost surely by Lemma~\ref{S:lem:5.0}. Taking $X=X_\beta$ and $q=q_\infty$ defines $\mathrm{Sine}_\beta$, whose spectrum has the Sine$_\beta$ law
\cite{valko2009continuum,valko2017sine}.

Write $\mathtt r\mathrm{Circ}_{\beta,n}$ and $\mathtt r\mathrm{Sine}_\beta$ for the conjugated inverses.
\subsubsection{Sampling the Brownian path}
\label{S:sec:intro-brownian-samples}
The leading coefficient depends on the sampling times. For large $k$, a circular step with radial parameter $\beta k/2$ is compared with a Brownian increment of duration
$4/(\beta k+1)$. This choice matches $\beta k/2$ to the parameter $2/h-1/2$ in the short time radial approximation of~\cite{valko2020operator}. Accordingly, put
\[
 \vartheta_{n,j}=\sum_{l=j}^{n-1}\frac4{\beta l+1},\qquad
 I_{n,k}=[1-(k+1)/n,1-k/n),\quad 0\le k<n.
\]
Here $1\le j\le n$, and an empty sum is zero. Define $\mathrm{Br}_{\beta,n}$ using the same boundary lines as $\mathrm{Sine}_\beta$ and the step coefficient
\[
 X_{\beta,n}^{\mathrm{br}}(t)
      =\mathcal A(Z_{\vartheta_{n,k+1}}),\qquad t\in I_{n,k}.
\]
This is a regular system for each $n$, with conjugated inverse $\mathtt r\mathrm{Br}_{\beta,n}$. Its first coefficient is the identity, since $\vartheta_{n,n}=0$. It depends only
on $Z$ and its terminal point $q_\infty$.

\subsubsection{The leading order coefficient}
\label{S:sec:intro-coefficient-and-coupling}
The leading order coefficient is determined by the Green kernel. Use the distance $r=1-t$ from the singular boundary and set $\theta_\beta(r)=(4/\beta)\log(1/r)$ and
$\mathsf X_\omega(r)=X_\beta(1-r,\omega)$. The boundary solutions at spectral parameter zero become the columns
\[
 a_\omega(r)=\mathsf X_\omega(r)u_0,\qquad
 c_\omega(r)=\mathsf X_\omega(r)u_1,
 \qquad u_1=(q_\infty,1)^{\mathsf T}.
\]
Write $p_\omega=|a_\omega|^2$ and $q_\omega=|c_\omega|^2$. In these coordinates the inverse kernel is
\[
 K_X(r,s)=\tfrac12a(r)c(s)^{\mathsf T}\quad(s<r),\qquad
 K_X(s,r)=K_X(r,s)^{\mathsf T}.
\]
We suppress $\omega$ in the next calculation. If $X$ is perturbed to $\me^{\epsilon U}X$, with $U$ bounded, traceless and supported away from zero, differentiation gives
\[
 (\mathscr D_XU)(r,s)
 =\tfrac12\bigl[(U(r)a(r))c(s)^{\mathsf T}
                     +a(r)(U(s)c(s))^{\mathsf T}\bigr],\qquad s<r.
\]
Consequently the two square terms in $\|\mathscr D_XU\|_{\mathfrak S_2}^2$ have sum
\begin{equation}\label{S:eq:intro-local-square-terms}
 \frac12\int_0^1\bigl[A(r)|U(r)a(r)|^2
                         +C(r)|U(r)c(r)|^2\bigr]\,\dd r,
\end{equation}
where
\[
 A_\omega(r)=\int_0^r q_\omega(s)\,\dd s,\qquad
 C_\omega(r)=\int_r^1p_\omega(s)\,\dd s.
\]
The factor $1/2$ in \eqref{S:eq:intro-local-square-terms} comes from squaring the kernel factor $1/2$ and integrating both triangles. The cross term involving different positions
is estimated after summation. The full identity for the Hilbert--Schmidt inner product and the nonlinear remainder are proved in
Section~\ref{S:sec:stability-and-fr-echet-expansion-of-the-inverse-kernel}.

The derivatives of $\mathcal A$ with respect to $x$ and $y$ give the matrices
\[
 \mathsf A_1=\begin{pmatrix}0&-1\\0&0\end{pmatrix},\qquad
 \mathsf A_2=\frac12\begin{pmatrix}-1&0\\0&1\end{pmatrix}.
\]
Indeed, $(y\partial_x\mathcal A)\mathcal A^{-1}=\mathsf A_1$ and $(y\partial_y\mathcal A)\mathcal A^{-1}=\mathsf A_2$. For a fixed vector $v$ and $g\sim N(0,hI_2)$,
\[
 \mathbb E\left|\sum_{\mu=1}^2g_\mu\mathsf A_\mu v\right|^2
 =h\mathfrak h_{\mathrm{aff}}(v),\qquad
 \mathfrak h_{\mathrm{aff}}(v)=\frac14v_1^2+\frac54v_2^2,
 \qquad v\in\mathbb R^2.
\]
This is a positive definite quadratic form on $\mathbb R^2$; its associated symmetric bilinear form is $(v_1w_1+5v_2w_2)/4$. Applying this covariance calculation to the two squares
in \eqref{S:eq:intro-local-square-terms} gives the weight
\begin{equation}\label{S:eq:intro-gamma}
 \Gamma_{\beta,\omega}(r)
 =\frac12\left[A_\omega(r)\mathfrak h_{\mathrm{aff}}(a_\omega(r))
             +C_\omega(r)\mathfrak h_{\mathrm{aff}}(c_\omega(r))\right].
\end{equation}
For the columns $a=(-y^{-1/2},0)^{\mathsf T}$ and $c=((q_\infty-x)y^{-1/2},y^{1/2})^{\mathsf T}$, the weight also has the form
\begin{equation}\label{S:eq:intro-gamma-lengths}
 \Gamma_{\beta,\omega}(r)
 =\frac{A_\omega(r)p_\omega(r)+C_\omega(r)q_\omega(r)}8
            +\frac{C_\omega(r)}{2p_\omega(r)}.
\end{equation}

The second factor depends on the position of the sample relative to its operator interval. Put $b_\beta=\tfrac12-\beta^{-1}$ and $r_k=k/n$. For $1\le k<n$, Lemma~\ref{S:lem:4.1}
gives
\[
 \frac{\vartheta_{n,k+1}-(4/\beta)\log(n/(k+1))}
      {(4/\beta)\log(1+1/k)}
 =b_\beta(1-r_k)+O_\beta(k^{-1}).
\]
The denominator is the Brownian time elapsed across that operator interval. Define
\begin{equation}\label{S:eq:intro-kappa}
 \kappa_\beta(r)=\int_0^1|s-b_\beta(1-r)|\,\dd s
 =\begin{cases}
   \frac12-b_\beta(1-r),&0<\beta\le2,\\[1mm]
   \frac12-b_\beta(1-r)+b_\beta^2(1-r)^2,&\beta>2.
  \end{cases}
\end{equation}
For $v\le0$, $\int_0^1|s-v|\,\dd s=1/2-v$; for $0\le v\le1$, the same integral is $v^2/2+(1-v)^2/2$. These two cases apply to $v=b_\beta(1-r)$ in \eqref{S:eq:intro-kappa}, and show
that $1/4\le\kappa_\beta\le C_\beta$. Proposition~\ref{S:prop:sampling-offset} gives the corresponding integrated variance:
\begin{equation}\label{S:eq:intro-integrated-time-variance}
 \int_0^1|\theta_\beta((k+s)/n)-\vartheta_{n,k+1}|\,\dd s
       =\frac4{\beta k}\kappa_\beta(r_k)+O_\beta(k^{-2}).
\end{equation}
The factor $1/n$ from the operator interval, together with $1/k=1/(nr_k)$ in this formula, explains the squared error scale $n^{-1}$ after summing over $k$. The leading order
coefficient is
\begin{equation}\label{S:eq:intro-coefficient}
 \PathCoeff(\omega)=\frac4\beta\int_0^1
          \kappa_\beta(r)\frac{\Gamma_{\beta,\omega}(r)}r\,\dd r.
\end{equation}
Identifying this integral requires controlling the dependence between the Brownian increments and the random weight, as well as the cross terms in the squared kernel. We describe
these points below.

To verify integrability in \eqref{S:eq:intro-coefficient}, set $a_\beta=\min\{1,2/\beta\}$. Lemma~\ref{S:lem:Gamma-convergence} applies the boundary estimates to give
$0\le\Gamma_{\beta,\omega}(r)\le C_{\beta,\omega,\zeta}r^\zeta$ for every fixed $0<\zeta<a_\beta$. Hence
\begin{equation}\label{S:eq:intro-coefficient-tail}
 0\le\int_0^\delta\kappa_\beta(r)
          \frac{\Gamma_{\beta,\omega}(r)}r\,\dd r
 \le\frac{C_{\beta,\omega,\zeta}}\zeta\delta^\zeta.
\end{equation}
The columns never vanish, and $A_\omega,C_\omega$ are positive on $(0,1)$. Thus $\Gamma_{\beta,\omega}>0$ there, proving positivity of the integral.
Proposition~\ref{S:prop:coefficient-positive-measurable} also records its measurability as a function of $Z$. At $\beta=2$, $\kappa_2=1/2$. Its change of formula comes from the
sampling position passing into the interval. This change depends on the chosen sampling times. The limit-point/limit-circle classification concerns the domain of the Sine operator.
Section~\ref{S:sec:the-critical-value-beta-2} discusses the distinction.

\subsubsection{The coupling and sampling theorems}
\label{S:sec:intro-main-results}
The following theorems identify the leading coefficient for approximation by Brownian samples and for a coupling with the circular $\beta$-ensemble. Throughout, $\beta>0$ is fixed,
and the finite walk, the limiting and sampled operators, and $\PathCoeff$ are those defined above.

Valk\'o and Vir\'ag~\cite{valko2020operator} constructed a coupling for which, almost surely for all sufficiently large
$n$,
\begin{equation}\label{S:eq:intro-previous-bound}
 \|\mathtt r\mathrm{Circ}_{\beta,n}
             -\mathtt r\mathrm{Sine}_\beta\|_{\mathfrak S_2}^2
       \le\frac{(\log n)^6}{n}.
\end{equation}
Our result identifies a leading order coefficient under a different coupling of the same operator laws.

\begin{theorem}\label{S:thm:main}
For each fixed $\beta>0$, there is a coupling of $Z$ and circular walks $b^{(n)}$ of the law above, with $b_n^{(n)}=q_\infty$ for every $n$, such that, almost surely,
\[
 n\|\mathtt r\mathrm{Circ}_{\beta,n}
             -\mathtt r\mathrm{Sine}_\beta\|_{\mathfrak S_2}^2
       \longrightarrow\PathCoeff(\omega)\in(0,\infty).
\]
\end{theorem}
The coefficient is given by \eqref{S:eq:intro-coefficient}. Thus
\[
 \|\mathtt r\mathrm{Circ}_{\beta,n}
             -\mathtt r\mathrm{Sine}_\beta\|_{\mathfrak S_2}
       =\frac{\sqrt{\PathCoeff(\omega)}+o_\omega(1)}{\sqrt n}.
\]
The construction preserves each finite walk law and places all sizes on the same probability space. The positive limit establishes the exact order for this coupling. Optimality
over all couplings would require a separate lower bound. Each circular inverse has infinite rank. Section~\ref{S:sec:proof-strategy} gives the reduction, and
Section~\ref{S:sec:proof-of-main-theorems} proves both approximation theorems.

The same functional is the limit for the sampled operator when only the Brownian path is used.

\begin{theorem}\label{S:thm:reference-main}
For each fixed $\beta>0$, almost surely under the law of $Z$,
\[
 n\|\mathtt r\mathrm{Br}_{\beta,n}
             -\mathtt r\mathrm{Sine}_\beta\|_{\mathfrak S_2}^2
       \longrightarrow\PathCoeff(\omega).
\]
\end{theorem}
In the coupling of Proposition~\ref{S:prop:7.2}, the remaining comparison satisfies
\begin{equation}\label{S:eq:intro-finite-comparison-rate}
 \sqrt n\,\|F_n^{\mathrm{full}}\|_{\mathfrak S_2}\longrightarrow0,
 \qquad F_n^{\mathrm{full}}
 =\mathtt r\mathrm{Circ}_{\beta,n}-\mathtt r\mathrm{Br}_{\beta,n},
\end{equation}
by Proposition~\ref{S:prop:circular-brownian-inverse-comparison}. Write $\varepsilon_{\beta,n}^{\mathrm{samp}} =\mathtt r\mathrm{Br}_{\beta,n}-\mathtt r\mathrm{Sine}_\beta$. Then
\begin{equation}\label{S:eq:intro-two-comparisons}
 \mathtt r\mathrm{Circ}_{\beta,n}-\mathtt r\mathrm{Sine}_\beta
       =F_n^{\mathrm{full}}+\varepsilon_{\beta,n}^{\mathrm{samp}}.
\end{equation}
After multiplication of the norm by $\sqrt n$, the first term tends to zero, while the squared norm of the second has a nonzero limit. An $O(n^{-1/2})$ bound for
$F_n^{\mathrm{full}}$ would leave both its square and its inner product with $\varepsilon_{\beta,n}^{\mathrm{samp}}$ in the leading order coefficient. We therefore use the
little-$o$ estimate \eqref{S:eq:intro-finite-comparison-rate}.

The finite operator has the prescribed circular $\beta$-ensemble distribution, and its terminal boundary line agrees with that of the limiting operator. Conditional on $q_\infty$,
we couple the radial and angular variables of the two increments by quantiles. The resulting relative matrix errors have smaller means than their individual fluctuations.
Estimating their products before taking absolute values preserves this gain. The estimates for the inverse kernel then turn it into \eqref{S:eq:intro-finite-comparison-rate};
separate bounds near the singular boundary permit removal of the cutoff. These are the tasks of
Sections~\ref{S:sec:coupling-the-finite-transfer-increments-to-brownian-increments}--\ref{S:sec:comparison-of-the-finite-and-sampled-brownian-operators}.

\subsubsection{Conditional second moments}
\label{S:sec:intro-local-variance}
The preceding covariance calculation freezes the boundary columns. In the actual problem those columns depend on the same Brownian path as the increments being squared, and the
terminal column also depends on $q_\infty$. We cannot treat these weights as deterministic or obtain the coefficient by averaging the bridge variance alone.

Both the conditional mean and conditional variance contribute, as the following scalar calculation shows. For Brownian motion $W$ on $[0,h]$, write $\Delta=W_h$ and
$b(s)=h^{-1/2}(W_{hs}-sW_h)$. Then $b$ is a standard Brownian bridge independent of $\Delta$, so
\begin{equation}\label{S:eq:intro-scalar-bridge-energy}
 \int_0^1\mathbb E[W_{hs}^2\mid\Delta]\,\dd s
       =\frac{\Delta^2}{3}+\frac h6.
\end{equation}
Indeed, the conditional second moment is $s^2\Delta^2+hs(1-s)$. The first term is the square of the conditional mean; the second is the conditional variance. Both are of order $h$
after averaging over $\Delta$, and their sum has mean $h/2$. Keeping only the bridge variance would therefore give the wrong coefficient even in this example.

For the operators, condition on one row of samples and the terminal point,
\[
 \mathscr S_n=\sigma(q_\infty,Z_{\vartheta_{n,j}}:1\le j\le n),
 \qquad \mathscr K=\mathfrak S_2(\mathcal H).
\]
With the kernel on the discarded boundary intervals included in the remainder, Sections~\ref{S:sec:the-error-from-sampling-the-brownian-path}
and~\ref{S:sec:the-conditional-gaussian-operator} construct
\[
 \varepsilon_{\beta,n}^{\mathrm{samp}}
       =G_n+B_n+\widehat R_n^{\mathrm{br}},\qquad
 G_n\mid\mathscr S_n\sim\mathcal N_{\mathscr K}(0,Q_n),\qquad
 \sqrt n\,\|\widehat R_n^{\mathrm{br}}\|_{\mathscr K}\longrightarrow0.
\]
Here $B_n$ and $Q_n$ are $\mathscr S_n$-measurable. It follows that
\[
 \mathbb E[\|G_n+B_n\|_{\mathscr K}^2\mid\mathscr S_n]
       =\operatorname{tr}Q_n+\|B_n\|_{\mathscr K}^2.
\]
Theorem~\ref{S:thm:trace-conditional-complete} and Proposition~\ref{S:prop:coefficient-riemann-limit} identify the limit of $n$ times this quantity as
\eqref{S:eq:intro-coefficient}. The calculation uses weighted quadratic variation for the increments between sampled positions. Products involving separated intervals require
additional estimates because the terminal boundary data depend on the future path.

We next estimate the fluctuation about the conditional second moment. For a real Gaussian vector in $\mathscr K$, the exact variance is
\[
 \operatorname{Var}(\|G_n+B_n\|_{\mathscr K}^2\mid\mathscr S_n)
       =2\operatorname{tr}(Q_n^2)+4\langle Q_nB_n,B_n\rangle.
\]
We bound the operator norm of $Q_n$ on a scale smaller than $n^{-1}$ and combine it with the estimate for the conditional second moment. Section~\ref{S:sec:covariance-estimates}
proves, for every $\epsilon>0$,
\[
 \operatorname{tr}Q_n+\|B_n\|_{\mathscr K}^2=O_\omega(n^{-1}),
 \qquad
 \|Q_n\|_{\mathrm{op}}
       =O_{\beta,\omega,\epsilon}(n^{-1-a_\beta+\epsilon}).
\]
Taking $0<\epsilon<a_\beta$, the conditional Gaussian exponential bound gives summable tail probabilities after localization and hence almost sure concentration along the full
sequence. This is Theorem~\ref{S:thm:covariance-concentration}. The small operator remainder and \eqref{S:eq:intro-finite-comparison-rate} then leave the same coefficient for the
sampled and circular inverses. The complete deduction, including the passage from the auxiliary bridge construction to the distribution of the Brownian path, is given in
Section~\ref{S:sec:proof-of-main-theorems}.

\subsubsection{Spectral consequences}
\label{S:sec:intro-spectral-consequences}
The inverse comparison controls the reciprocal eigenvalues globally. Passing back to the eigenvalues of the differential operators requires control away from zero in the reciprocal
spectrum; comparisons of individual eigenspaces also retain the corresponding spectral gaps.

Let $\mu_j^+(A)$ and $\mu_j^-(A)$ be the positive and negative eigenvalues of a compact self-adjoint operator, ordered respectively in decreasing and increasing order, with
multiplicity and with finite lists extended by zeros. For $A_n=\mathtt r\mathrm{Circ}_{\beta,n}$ and $A=\mathtt r\mathrm{Sine}_\beta$, the signed compact self-adjoint
Hoffman--Wielandt inequality \cite{bhatia1994hoffman} (equation~(27), p.~490) gives
\begin{equation}\label{S:eq:intro-reciprocal-spectral-bound}
 \sum_{j\ge1}\sum_\pm|\mu_j^\pm(A_n)-\mu_j^\pm(A)|^2
 \le\|A_n-A\|_{\mathfrak S_2}^2
 =\frac{\PathCoeff(\omega)+o_\omega(1)}n.
\end{equation}
The sum in \eqref{S:eq:intro-reciprocal-spectral-bound} depends only on the two spectra. It compares the full sequences of reciprocal eigenvalues, including those outside any fixed
spectral interval. Fixed direct eigenvalues and isolated spectral projections converge at rate $O_\omega(n^{-1/2})$, with constants depending on the spectral gap of the limiting
operator. Theorems~\ref{S:thm:spectral-isolated-projections} and~\ref{S:thm:spectral-fixed-window} give the precise statements, including pairing on bounded spectral intervals
whose boundary points avoid the limiting spectrum. Section~\ref{S:sec:spectral-consequences} proves expansions of resolvents and spectral projections in the common Hilbert space
specified above. Their norm estimates are unchanged by applying a common unitary transformation to both operators. The comparison between different operators uses the stated
identifications of their Hilbert spaces. The upper estimates alone do not identify a leading order coefficient for an individual eigenvalue error.

\subsubsection{Related work}
\label{S:sec:intro-related-work}
The exact circular models of Killip and Nenciu~\cite{KillipNenciu2004}, the circular bulk limit of Killip and Stoiciu~\cite{KillipStoiciu2009}, and the Brownian carousel of Valk\'o
and Vir\'ag~\cite{valko2009continuum} provide the finite and limiting descriptions underlying the operator representation. The coupling defined through stopping times
in~\cite{valko2020operator} controls the inverse through path comparison. We compare the exact finite law with Brownian values at deterministic times. The comparison error is
smaller than the leading sampling error, whose coefficient can then be computed from the first variation of the Green kernel.

Operator couplings have also been used to study other observables. Lambert~\cite{Lambert2021} used the coupling of~\cite{valko2020operator} to obtain a mesoscopic central limit
theorem for Sine$_\beta$. The Dirac representation gives the stochastic zeta function~\cite{valko2022zeta} and Palm measures~\cite{VVPalm2023}. Li and Valk\'o proved operator
limits and convergence of normalized characteristic polynomials for circular Jacobi and real orthogonal ensembles~\cite{li2022circularjacobi}, and spectral limits and convergence
of normalized characteristic polynomials for rank-one truncations of circular ensembles~\cite{li2025truncated}. Corrections to marginal statistics are a different question:
Forrester and Shen~\cite{ForresterShen2026} obtain expansions in powers of $n^{-2}$ for bulk correlations and spacing statistics at $\beta=1,2,4$.

For beta ensembles, the tridiagonal models of Dumitriu and Edelman~\cite{dumitriu2002matrix} led to the stochastic operator interpretation of Edelman and
Sutton~\cite{edelman2007random}. The limit given by the stochastic Airy operator was proved by Ram\'irez, Rider, and Vir\'ag~\cite{ramirez2011beta} and extended to uniformly convex
polynomial potentials by Krishnapur, Rider, and Vir\'ag~\cite{krishnapur2016universality}. Gorin and Shkolnikov~\cite{gorin2018stochastic} obtained the stochastic Airy semigroup;
related transfer equations appear in~\cite{kritchevski2012scaling}. Parameter limits connect the stochastic Bessel operator to the stochastic Airy operator~\cite{dumaz2021hardsoft}
and the canonical system associated with the stochastic Airy operator at high energy to the limiting canonical system in the bulk~\cite{painchaudpaquette2025softbulk}
(Theorems~1--2). After a change of variables, the coefficient matrices in the latter result converge vaguely in distribution; control of the boundary data also enters the
convergence of Weyl--Titchmarsh functions and spectral measures.

The result for the tridiagonal models in Section~\ref{sec:soft-edge-ensembles-and-quantitative-airy-limits} proves squared Hilbert--Schmidt order $n^{-1/3}$ for Gaussian beta
ensembles and Jacobi models with uniformly convex polynomial potential at a regular soft edge, with a lower bound for every approximation of rank at most $n$. The result for the
bidiagonal models in Section~\ref{sec:hard-edge-models-brownian-information-and-bessel-limits} proves order $n^{-1}$ for beta-Laguerre and polynomial bidiagonal models with
$x\mapsto V(x^2)$ uniformly convex, for fixed $\beta>0$ and $a>0$; its lower bound concerns Brownian observations at deterministic times and independent auxiliary randomness. Here
we identify a coefficient for a specified circular coupling. The result for Wigner matrices in
Section~\ref{sec:wigner-matrices-and-the-airy-limit-for-a-deterministic-starting-vector} proves an operator limit in distribution from a deterministic Lanczos starting vector under
its entry hypotheses. The circular argument is independent of these other results.

Section~\ref{S:sec:proof-strategy} explains the estimates required for the two asymptotics; their final proof is in Section~\ref{S:sec:proof-of-main-theorems}.
Sections~\ref{S:sec:green-kernels-for-canonical-systems} and \ref{S:sec:stability-and-fr-echet-expansion-of-the-inverse-kernel} establish the estimates for the Green kernel.
Section~\ref{S:sec:sampling-times-and-estimates-for-the-boundary-solutions} proves the boundary estimates and analyzes the sampling times.
Sections~\ref{S:sec:coupling-the-finite-transfer-increments-to-brownian-increments}--\ref{S:sec:comparison-of-the-finite-and-sampled-brownian-operators} construct and estimate the
circular coupling. Section~\ref{S:sec:the-error-from-sampling-the-brownian-path} constructs the Brownian bridge expansion.
Sections~\ref{S:sec:the-conditional-gaussian-operator}--\ref{S:sec:the-limiting-coefficient-of-the-error} treat the conditionally Gaussian random operator, its concentration, and
the leading order coefficient. Sections~\ref{S:sec:the-critical-value-beta-2} and \ref{S:sec:spectral-consequences} discuss $\beta=2$ and spectral consequences, respectively.

\clearpage
\section{Proof strategy}\label{sec:overall-proof-strategy}
\hypertarget{proof-strategy}{}

The proof for Wigner matrices starts from the independent entries of the original matrix. No explicit joint distribution of the recurrence coefficients is available, so we prove
their convergence by comparison. Two estimates are needed: convergence of the integrated potential on bounded intervals, and localization of the entire low energy spectral
subspace. Together they identify the limiting differential operator and give compactness in the operator norm topology. The coefficient limits, bounds on the Gram matrix, and
localization estimates for the GOE and GUE used in this comparison are proved in Section~\ref{sec:finite-gaussian-estimates-for-the-wigner-comparison}, independently of the
quantitative coupling theorems.

The coupling arguments use the specified joint distributions of the finite matrix coefficients. We construct these coefficients together with the Brownian path of the limiting
operator. An identity for the inverse or resolvent expresses the operator difference in terms of the coefficient differences. The kernel estimates determine the weighted norms used
for these differences. At the soft edge, summation by parts retains the discrete derivative in the coefficient error. At the hard edge, the inverse formula involves partial sums of
logarithmic ratios. For the circular ensemble, we apply the first derivative of the Green kernel to the linear interpolation and Gaussian bridge terms, and compute the resulting
conditional mean and covariance. We treat these three estimates after the Wigner argument.

\subsection*{Operator limit for Wigner matrices}
\hypertarget{proof-strategy-wigner}{}

Fix a deterministic unit vector $w_N$, and write $a_{k,N},b_{k,N}$ for the coefficients of its completed Lanczos matrix $T_N$. Put $m=N^{1/3}$, $h=m^{-1}$ and
$\mathcal A_N=m^2(2I-T_N)$, and let $\mathcal E_N$ be the normalized step embedding with cell length $h$. The quantity that identifies the limiting potential is
\[
 \widetilde Y_{N,k}
 =m\sum_{j=1}^k(2-a_{j,N}-\bar b_{j-1,N}-b_{j,N}),
 \qquad \bar b_{0,N}=1,\quad \bar b_{j,N}=b_{j,N}\ (j\ge1).
\]
To see why, take a compactly supported smooth function $f$ with $f(0)=0$ and set $f_j=f(jh)$. The potential part of the quadratic form at the vector $(h^{1/2}f_j)_j$ is
\[
 m\sum_j(2-a_{j,N}-\bar b_{j-1,N}-b_{j,N})|f_j|^2
 =-\sum_j\widetilde Y_{N,j}(|f_{j+1}|^2-|f_j|^2).
\]
Thus local convergence of the interpolated potential to $Y$, together with $b_{k,N}\to1$ on bounded intervals, gives the limiting form
\[
 \mathfrak q_Y[f]=\int|f'|^2-\int Y(|f|^2)'.
\]
For the GOE and GUE, the independent normal and chi coefficients give $Y(x)=x^2/2+2B(x)/\sqrt\beta$. Substitution yields the stochastic Airy operator with Dirichlet boundary
condition. We have to transfer this coefficient limit through the independent entries of $H_N$ and then control the operator outside bounded intervals.

\paragraph*{Cancellation in the integrated potential}
Suppress the index $N$. Write $\mathcal L_H(f)=w^*f(H)w$ and $p_j(x)=U_j(x/2)$. Orthogonalizing these polynomials expresses the recurrence coefficients through
\[
 G_\rho=(\mathcal L_H(p_ip_j))_{0\le i,j<\rho},\qquad \rho\asymp m.
\]
The growing degree is essential: a fixed limiting interval contains order $m$ coefficients. Smoothness of orthogonalization requires control of $G_\rho^{-1}$ as well as $G_\rho$.
For the GOE and GUE both norms are tight by Lemma~\ref{W:lem:gaussian-gram-conditioning}. We therefore differentiate with the Gram eigenvalues restricted to $[\kappa,\kappa^{-1}]$,
and include the smooth cutoff in the comparison itself.

Let $q_j$ be the orthonormal polynomials, with positive leading coefficient. For a variation of the moment functional in a direction $v$ with $v(1)=0$, differentiation of
orthogonality and of the recurrence gives
\[
 Da_{j+1}[v]=b_{j+1}v(q_jq_{j+1})-b_jv(q_{j-1}q_j),
 \qquad Db_j[v]=\tfrac12b_jv(q_j^2-q_{j-1}^2).
\]
Summing the first identity gives $D\sum_{j=1}^Ka_j[v]=b_Kv(q_{K-1}q_K)$. The interior terms cancel. For $\sigma=\pm1$ and $K\le\rho-2$, define
\[
 \mathcal T_K^\sigma=\sigma\sum_{j=1}^Ka_j
       +2\sum_{j<K}(b_j-1)+(b_K-1),
 \qquad \widetilde Y_{N,K}^\sigma=-m\mathcal T_K^\sigma.
\]
Combining the two variations yields
\[
 D\mathcal T_K^\sigma[v]
 =v\!\left[\tfrac12b_K(q_K+\sigma q_{K-1})^2
       +\sum_{j<K}(b_j-b_{j+1})q_j^2-b_1q_0^2\right].
\]
The interior squares are weighted by the variation of $b_j$. Bounding the original summands separately would lose this cancellation. At each fixed derivative order $r$, it saves
one factor of $\rho$: the degree factor is $\rho^{r-1}$. The other factors involve polynomial coefficient sums and coefficient variation; their cutoffs have thresholds $N^{o(1)}$.
For a fixed finite linear combination $\Theta=m\sum_l\lambda_l\mathcal T_{K_l}^{\sigma_l}$, a partition $\pi$ of $r$ derivative labels into $k$ sets consequently contributes the
factor
\[
 N^{o(1)}m^k\prod_{A\in\pi}\rho^{|A|-1}
 =N^{o(1)}m^k\rho^{r-k}=N^{r/3+o(1)}.
\]
Without this cancellation, the bound for each set of the partition would contain an additional factor $\rho$. Proposition~\ref{W:prop:comparison-global} gives the resulting bounds
in the specified moment norms, including derivatives of both cutoffs.

\paragraph*{Comparison through the matrix entries}
On the event where spectral regularization agrees with $H$, functional calculus turns a polynomial moment derivative into
\[
 D_H\mathcal L_H(f)[V]
 =\frac1{2\pi\mathrm i}\oint_{\mathcal C}
       f(z)w^*R(z)VR(z)w\,\mathrm dz,
 \qquad R(z)=(zI-H)^{-1},
\]
where $\mathcal C$ encloses the spectrum. A change in entry $(a,b)$ makes $V$ supported on the corresponding two coordinates. The product then contains the mixed resolvent entries
$w^*R(z)e_a$ and $e_b^*R(z)w$. Their dependence on $a,b$ is retained when the entry contributions are summed; an operator norm bound on every resolvent factor would discard this
information.

The moment assumptions permit truncation with probability $o(1)$ of changing $H_N$. An intermediate bounded entry law matches the first three off-diagonal moments of the normalized
truncation. In the comparison with the GOE or GUE, the cubic term remains because those moments need not vanish. Taylor expansion in one matrix entry produces products of
resolvents evaluated at matrices with different entries set to zero. Lemma~\ref{W:lem:cubic-entry-restoration} rewrites these products using resolvents evaluated at the same matrix
and also estimates the changes in the derivative coefficients.

The common off-diagonal law now allows the cubic terms to be summed before taking absolute values. The relevant deterministic Hermitian matrix $B_{3,N}$ has zero diagonal and
constant upper triangular entries equal to the difference of the third moments $\mathbb E\xi|\xi|^2$ of the two standardized laws. Its trace norm is $O(N\log N)$. Writing
$B_{3,N}=\sum_l\lambda_lv_lv_l^*$ and $R=(zI-H)^{-1}$ gives
\[
 |w^*RB_{3,N}Rw|
 \le(w^*|R|w)\sum_l|\lambda_l|\,v_l^*|R|v_l
 \le N^{1+o(1)}
\]
on the event used in the comparison. The entry normalization contributes $N^{-3/2}$ and the coefficient and contour estimates contribute $mN^{o(1)}$. The resulting scale is
\[
 N^{-3/2}mN^{1+o(1)}=N^{-1/6+o(1)}.
\]
This bounds the cubic contribution to the comparison after regularization. Theorem~\ref{W:thm:cumulant-complete} includes the other cubic terms. The fourth order and diagonal
remainders use conditional estimates under the stated moment assumptions.

Apply the comparison to characteristic functions and to smooth functions of the modulus of continuity. The first identifies finite dimensional distributions; the second gives
tightness. We also compare the expectations of the cutoff functions, so the cutoffs can be removed after $N\to\infty$; see Theorem~\ref{W:thm:tightness-main} and
Corollary~\ref{W:cor:tightness-Fourier}. The two choices of sign are both needed. With $e_k=b_k-1$, $e_0=0$, and
$\Delta\widetilde Y_{N,k}^\sigma=\widetilde Y_{N,k}^\sigma-\widetilde Y_{N,k-1}^\sigma$, their definitions give
\[
 e_k+e_{k-1}=-\tfrac h2
       (\Delta\widetilde Y_{N,k}^++\Delta\widetilde Y_{N,k}^-),
 \qquad
 \max_{k\le K}|e_k|
 \le\tfrac{Kh}{2}\sum_{\sigma=\pm1}
                     \max_{j\le K}|\Delta\widetilde Y_{N,j}^\sigma|.
\]
For $K=\lceil mL\rceil+1$, tightness of the continuous interpolations makes the maximum increments tend to zero. Since $Kh$ stays bounded, the off-diagonal coefficients converge
uniformly to one on every fixed interval.

\paragraph*{Control of the whole low energy subspace}
Local convergence of the coefficients does not imply localization of eigenvectors whose eigenvalues are bounded above. For example, setting $b_{N-1,N}=0$ and $a_{N,N}=2+m^{-2}$
leaves every fixed initial interval unchanged but gives $\mathcal A_Ne_N=-e_N$. The vector $\mathcal E_Ne_N$ is supported outside every fixed bounded interval for all sufficiently
large matrix sizes, while its resolvent image has norm $|-1-z|^{-1}$. We must therefore control the full spectral subspace, including any blocks created by a Lanczos breakdown.

Let $Q_{N,L}$ be the orthogonal projection that annihilates the first $\lfloor mL\rfloor$ coordinates and put $\Pi_N(E)=\mathbf1_{(-\infty,E]}(\mathcal A_N)$. For
$\theta_E(x)=(1+e^{x-E})^{-2}$, use
\[
 C_{N,L,E}=\operatorname{Tr}[Q_{N,L}\theta_E(\mathcal A_N)]
 =\sum_j\theta_E(\lambda_{j,N})\|Q_{N,L}v_{j,N}\|^2,
 \qquad \|Q_{N,L}\Pi_N(E)\|^2\le4C_{N,L,E}.
\]
A normalized eigenvector with eigenvalue at most $E$, supported outside the retained coordinates, contributes at least $1/4$, even if it has zero weight at $w_N$. This explains the
use of an unnormalized trace. We compare $C_{N,L,E}/(1+C_{N,L,E})$, multiplied by the cutoff on the eigenvalues of the Gram matrix. Its derivatives include a trace term without the
starting vector; its cubic contribution is estimated separately in Lemma~\ref{W:lem:localization-trace-cubic}. Coercivity and eigenvalue counting estimates for the GOE and GUE
bound its expectation by $C_{\beta,E}/L+o(1)$, uniformly in that cutoff, in Proposition~\ref{W:prop:gaussian-localization-expectation}. Taking $N\to\infty$, removing the cutoff,
and then letting $L\to\infty$ yields
\[
 \lim_{L\to\infty}\limsup_N
 \sup_{\|w\|=1}\mathbb P(\|Q_{N,L}\Pi_N(E)\|>\eta)=0.
\]
Each probability fixes $w$ before sampling the matrix.

Finally, for compact $K\Subset\mathbb C\setminus\mathbb R$, set $d_K=\inf_{z\in K}|\operatorname{Im}z|$ and $e_K=\sup_{z\in K}\operatorname{Re}z$. The spectral theorem gives, for
$E>e_K$,
\[
 \sup_{z\in K}\|Q_{N,L}(\mathcal A_N-z)^{-1}\|
 \le d_K^{-1}\|Q_{N,L}\Pi_N(E)\|+(E-e_K)^{-1}.
\]
Choose $E$ large and then $L$ large. This controls the tail, while local form estimates give compactness on the retained interval and identify the weak resolvent equation.
Theorem~\ref{W:thm:operator-probabilistic} combines these facts to give convergence of the embedded resolvents to $(\mathcal H_\beta-z)^{-1}$, locally uniformly in operator norm in
distribution.

\subsection*{Coupling rates}
\hypertarget{proof-strategy-coupling}{}

Fix the parameters of the ensemble. We construct all finite sizes and one limiting Brownian path on a common probability space, preserving the prescribed distribution at each size.
Different sizes may be dependent. We write $X_n=O_\omega(a_n)$ if, almost surely, $X_n\le C(\omega)a_n$ for all sufficiently large $n$, with a finite constant that may also depend
on the fixed parameters.

The comparison passes through an operator constructed from samples of the Brownian path. Writing $\mathsf R_n$, $\mathsf R_n^B$ and $\mathsf R$ for the finite, sampled and limiting
inverse or resolvent on the common Hilbert space, we use
\[
 \mathsf R_n-\mathsf R
   =(\mathsf R_n-\mathsf R_n^B)+(\mathsf R_n^B-\mathsf R).
\]
The inverse formula determines the norm in which the finite coefficients must be coupled to the samples. We then estimate the second difference from the Brownian path and control
the discarded boundary intervals. For the Airy and Bessel operators, these estimates give the upper bound, which is matched by a lower bound in the stated approximation class. For
the circular ensemble, the first difference is $o_\omega(n^{-1/2})$ in Hilbert--Schmidt norm, and the second determines the positive leading coefficient of the squared error. Each
construction also verifies the prescribed joint distribution of the coefficients at every matrix size.

\subsubsection*{Approximation of the stochastic Airy operator}

At the soft edge, put $h=h_n\asymp n^{-1/3}$. The coefficient comparison starts from the resolvent identity
\[
 R_1(z)-R_0(z)=-R_1(z)(H_1-H_0)R_0(z),
 \qquad R_a(z)=(H_a-z)^{-1}.
\]
The size of a coefficient difference alone does not determine its contribution to this product. Under the common form bounds, Propositions~\ref{A:res:3.16} and~\ref{A:res:3.17}
give the following estimates. For the diagonal perturbation $H_1-H_0=h^{-2}M_r$,
\[
 \|R_1-R_0\|_{\mathfrak S_2}^2
 \le C_{\omega,K}h^{-3}\sum_{j=1}^{N}(1+x_j)^{-7/2}|r_j|^2.
\]
For $H_1-H_0=M_{(e_j-e_{j-1})/h}$,
\[
 \|R_1-R_0\|_{\mathfrak S_2}^2
 \le C_{\omega,K}\sum_{j=0}^{N}h(1+x_j)^{-3/2}|e_j|^2.
\]
Here $M$ denotes diagonal multiplication, $x_j=(j-\tfrac12)h$ for $j\ge1$, $x_0=0$, and $z$ ranges over a fixed compact $K\Subset\mathbb C\setminus\mathbb R$. The second estimate
uses summation by parts, retaining the terms at both Dirichlet boundaries. For example, if $|e_j|^2\le Ch$, its right side is $O(h)$. Substituting $r_j=h(e_j-e_{j-1})$ into the
first estimate gives only $O(h^{-1})$. Thus a discrete derivative in the coefficient difference must be kept in that form. The same estimate applies to the first order forms
arising from the off-diagonal coefficients.

For independent Gaussian and chi coefficients, the coupling can be constructed directly: an orthogonal rotation combines each Brownian increment with an auxiliary normal variable,
and the chi quantile map gives the prescribed coefficient law. The resulting scalar discrepancy satisfies the first weighted estimate, while the off-diagonal fluctuations are
handled by the second. For polynomial potentials, the coefficients are dependent, and the principal difficulty is to obtain the same decomposition from their conditional Gaussian
approximation.

The covariance of this approximation is obtained from the inverse conditional Hessian. The Brownian coefficient is determined by the variance of long partial sums. At the edge that
variance per term is $4h^3/\beta$. To explain its use, fix the macroscopic parameter at its edge value and let $g(w)=\sum_{\ell\ge0}g_\ell w^\ell$ be the analytic factor of the
covariance symbol divided by $h^3$. Its coefficients decay exponentially and $g(1)=2/\sqrt\beta$. Put $\sigma=g(1)$, $p_\ell=-\sum_{r>\ell}g_r$, and, for independent standard
normal variables $(\xi_j)_{j\in\mathbb Z}$, put $e_j=h^{1/2}\sum_{\ell\ge0}p_\ell\xi_{j-\ell}$. Then
\[
 \begin{aligned}
 g(w)&=\sigma+(1-w)\sum_{\ell\ge0}p_\ell w^\ell,\\
 Z_j:=h^{3/2}\sum_{\ell\ge0}g_\ell\xi_{j-\ell}
      &=h^{3/2}\sigma\xi_j+h(e_j-e_{j-1}).
 \end{aligned}
\]
Taking $\xi_j=h^{-1/2}(B(jh)-B((j-1)h))$ for $j\ge1$ gives exactly the Brownian increment term required by the difference operator. The remaining term is a discrete derivative, to
which the second weighted estimate applies. Lemma~\ref{A:res:15.1} and Theorem~\ref{A:res:15.3} justify the variance identity and factorization. In the finite model the factors
vary with the index; that variation, the finite covariance comparison, and the conditional approximation produce further terms controlled by the first weighted bound. The common
form bounds are proved before applying these inequalities. Integrating the conditional construction against the joint distribution of the coefficients outside the conditioned block
preserves the prescribed joint coefficient distribution.

It remains to compare the Brownian difference operator with the differential operator. If the retained interval has length $L\asymp h^{-2}$, Theorem~\ref{A:res:1.3} gives squared
resolvent error
\[
 C_{\omega,K}\bigl(h+L^{-1/2}\bigr)=O_\omega(h).
\]
The two contributions are discretization of the integrated potential and restriction of the half-line. The discretization includes the mass matrix of the piecewise linear basis and
the change to the step embedding. For general potentials, the omitted coordinates are estimated by Schur complements before the comparison map is applied. The operator compared is
therefore the compressed full resolvent, with the full matrix inverted first (Theorem~\ref{A:res:17.6}). The coefficient and discretization estimates now give the upper bound. The
passage to inverses is treated in Section~\ref{A:sec:inverse-coupling-and-spectral-consequences}.

The matching lower bound requires no coupling construction. Since $\Lambda_j\sim(3\pi j/2)^{2/3}$,
\[
 \inf_{\operatorname{rank}F\le n}\|F-R_\beta(z)\|_{\mathfrak S_2}^2
 =\sum_{j>n}s_j(R_\beta(z))^2\asymp n^{-1/3}.
\]
Proposition~\ref{A:res:5.1} therefore gives sharpness under the rank constraint. The leading coefficient of the matrix coupling error is not identified.

\subsubsection*{Approximation of the stochastic Bessel operator}

At the hard edge the inverse has a product structure. On the cells of the normalized embedding, its kernel is
\[
 R_n^\Delta(i,j)=p_{n,i}p_{n,j}\sum_{\ell=i\vee j}^{n}q_{n,\ell}.
\]
The factors $p_{n,k}$ contain products of bidiagonal coefficient ratios over intervals of increasing length. Errors in individual ratios can accumulate in those products, so the
coupling must control their logarithmic partial sums. For the polynomial model the partition is $z_{n,k}=\theta_V(k/n)$; for beta-Laguerre it is $z_{n,k}=k/n$. The matrix
normalization and the coordinate change are those of the main theorems.

Retain cells in an interval $I$ separated from the boundaries, and set $p_{n,k}^B=p^B(z_{n,k})$, $Q_{n,k}^B=\int_{(z_{n,k-1},z_{n,k}]}q^B(s)\,\dd s$, and
$h_{n,k}=\log(z_{n,k}/z_{n,k-1})$. Define the logarithmic differences by
\[
 e^{-c_n}p_{n,k}=p_{n,k}^Be^{\delta_{n,k}},
 \qquad
 e^{2c_n}q_{n,k}=Q_{n,k}^Be^{-2\delta_{n,k}+\rho_{n,k}}.
\]
The scalar $c_n$ cancels from every product $p_{n,i}p_{n,j}q_{n,\ell}$. The quantity $\delta_{n,k}$ contains the cumulative ratio difference; $\rho_{n,k}$ contains the local
coefficient and cell integral differences. Interpolating these exponents in the inverse formula gives Theorem~\ref{B:thm:weighted-inverse-stability}. For $1<\lambda<\min\{a+1,2\}$,
it reduces the coefficient comparison to
\[
 \begin{gathered}
 \max_{k\in I}(|\delta_{n,k}|+|\rho_{n,k}|)\longrightarrow0,\\
 \sum_{k\in I}h_{n,k}z_{n,k}^{\lambda}
       (|\delta_{n,k}|^2+|\rho_{n,k}|^2)=O_\omega(n^{-1}).
 \end{gathered}
\]
The first condition bounds the exponential factors throughout the interpolation. Once it holds, the second bounds the squared Hilbert--Schmidt difference of the inverses obtained
from the restricted factors. These are separate requirements: small logarithmic errors alone do not give the rate, and a weighted average alone does not bound the exponentials.

For beta-Laguerre, quantile coupling gives centered logarithmic ratio errors of order $k^{-3/2}$ in finite moments and means of order $k^{-2}$. Their sums have the required bounds.
An independent Gamma variable reconstructs each pair of coefficients, so the ratio coupling preserves the joint distribution of all coefficients at each matrix size. For polynomial
potentials, a conditional Gaussian approximation is first coupled to the Brownian increments. Its spectral factorization produces a discrete difference $\Xi_{n,k}-\Xi_{n,k+1}$;
summing leaves a value at the upper index, which is included in $c_n$. The transport of the exact logarithmic law to this Gaussian law is estimated directly in a weighted partial
sum norm (Theorem~\ref{B:thm:logarithmic-ratio-coupling}). Conditional sampling of the remaining coordinates then recovers the original coefficient law on the same Brownian space.
Theorem~\ref{B:thm:simultaneous-coefficient-construction} performs this construction at every size, including the original singular exterior densities when $0<\beta<1$.

The estimates so far concern restricted triangular factors. Compressing their full positive products gives an additional positive operator of rank at most one, which must also be
estimated. With discarded lengths $d_n,\varepsilon_n=O((\log n)^A/n)$, the squared boundary contributions are bounded by
\[
 C_\omega\bigl(d_n^{1+\min\{a,1\}-\eta}
                 +\varepsilon_n^3+\varepsilon_n^2\bigr)
 =o(n^{-1}),\qquad 0<\eta<\min\{a,1\},
\]
as in Theorem~\ref{B:thm:inverse-truncation-bound}. This explains the condition $a>0$ in the upper bound. The Brownian sampling comparison has squared error $O_\omega(n^{-1})$ by
Theorem~\ref{B:thm:brownian-inverse-approximation}; changing the cell embedding to the embedding of the main theorem has norm $O_\omega(n^{-1})$. Combining these estimates proves
the upper bound for the full inverse.

For the lower bound, the construction uses only $\mathcal N_n=\sigma\{W(-\log z_{n,k}):1\le k\le n\}$ and auxiliary randomness independent of the entire path. The integrals
$Q_{n,k}^B$ above are used in the comparison, not as additional observations in the construction. Conditional on $\mathcal N_n$, Brownian bridges remain random between consecutive
sample times. After localization to an interior rectangle, the linearized inverse has conditional covariance with trace at least $c/n$ and operator norm at most $C/n^2$. A Gaussian
small ball estimate, uniform over the translate selected from the observations, and a separate bound on the nonlinear remainder give summable failure probabilities.
Theorem~\ref{B:thm:brownian-information-lower-bound} yields the almost sure lower bound for each fixed approximation sequence with this measurability. It holds for $a>-1$, and its
probability one event may depend on the sequence. The restriction to the observations at size $n$ is essential: unrestricted spectral truncations of the trace class inverse have
squared error $o(n^{-1})$. The coupling estimates determine the order $n^{-1}$; they do not identify a leading coefficient.

\subsubsection*{The circular ensemble}

Fix $\beta>0$. The circular result identifies the coefficient of the squared inverse error. We compare the circular operator with the operator obtained by sampling the limiting
Brownian path. Write $K_n$, $K_n^B$ and $K$ for the circular, sampled and limiting inverses in the common unweighted space. We prove
\[
 \sqrt n\,\|K_n-K_n^B\|_{\mathfrak S_2}\longrightarrow0,
 \qquad
 n\|K_n^B-K\|_{\mathfrak S_2}^2\longrightarrow\mathcal C_\beta(Z).
\]
The second comparison determines the coefficient. The first preserves it when passing to the circular operator, since Cauchy--Schwarz makes the corresponding mixed term tend to
zero after multiplication by $n$.

To see what must be estimated, use $r=1-t$ and let $X(r)$ be the matrix coefficient of the canonical system. Its two boundary columns are $a(r)=X(r)u_0$ and $c(r)=X(r)u_1$, where
$u_0=(-1,0)^{\mathsf T}$ and $u_1=(q_\infty,1)^{\mathsf T}$. For bounded traceless $U$ supported away from $r=0$, differentiate $X\mapsto e^{\epsilon U}X$ at $\epsilon=0$. On one
triangle the inverse kernel and its derivative are
\[
 \begin{aligned}
 K_X(r,s)&=\tfrac12a(r)c(s)^{\mathsf T},\\
 (DK_X[U])(r,s)
 &=\tfrac12\bigl[(U(r)a(r))c(s)^{\mathsf T}
                  +a(r)(U(s)c(s))^{\mathsf T}\bigr],\qquad s<r.
 \end{aligned}
\]
The other triangle is determined by symmetry. Squaring and integrating gives the two local terms
\[
 \frac12\int_0^1\bigl[A(r)|U(r)a(r)|^2
                         +C(r)|U(r)c(r)|^2\bigr]\,\dd r,
 \qquad
 A(r)=\int_0^r|c(s)|^2\,\dd s,
 \quad C(r)=\int_r^1|a(s)|^2\,\dd s,
\]
together with a mixed term involving two positions. These are the relevant weights: a boundary column may grow while its product with the integral of the other column remains
integrable. The same weights control the nonlinear part of the kernel expansion. Uniform bounds on both columns would discard this compensation near $r=0$.

The scale $n^{-1}$ follows from the two different time coordinates. The Brownian time is $\theta_\beta(r)=(4/\beta)\log(1/r)$ and the sample times are
$\vartheta_{n,j}=\sum_{l=j}^{n-1}4/(\beta l+1)$. An operator interval has length $1/n$, whereas its Brownian time has length approximately $4/(\beta k)$ when $r\approx k/n$. The
sample need not lie at the endpoint of that Brownian interval. Its limiting relative position is $b_\beta(1-r)$, where $b_\beta=1/2-1/\beta$. The integrated local variance is
consequently proportional to
\[
 \kappa_\beta(r)=\int_0^1|s-b_\beta(1-r)|\,\dd s.
\]
The affine coefficient $\mathcal A(x+\mathrm iy)$ has infinitesimal noise matrices sending $v$ to $(-v_2,0)^{\mathsf T}$ and $(-v_1,v_2)^{\mathsf T}/2$. Their squared lengths sum
to $\mathfrak h_{\mathrm{aff}}(v)=(v_1^2+5v_2^2)/4$. Applying this covariance to the two local terms gives
\[
 \Gamma(r)=\tfrac12\bigl[A(r)\mathfrak h_{\mathrm{aff}}(a(r))
                           +C(r)\mathfrak h_{\mathrm{aff}}(c(r))\bigr],
 \qquad
 \mathcal C_\beta(Z)=\frac4\beta\int_0^1
                   \kappa_\beta(r)\frac{\Gamma(r)}r\,\dd r.
\]
Indeed, multiplication of the interval length $1/n$ by its variance $4\kappa_\beta(r)/(\beta nr)$ produces the Riemann sum for this expression after scaling by $n$. Boundary
estimates give $\Gamma(r)\le C_{\beta,\omega,\zeta}r^\zeta$ for every fixed $0<\zeta<a_\beta:=\min\{1,2/\beta\}$. They both justify removing the boundary restriction and prove
finiteness of the integral. Its positivity follows from the nonvanishing columns.

Two points prevent the preceding covariance calculation from being a proof by itself. The weights depend on the same Brownian path as the increments, and the conditional mean
contributes at the same scale as the conditional variance. For example, if $W$ is Brownian motion and $\Delta=W_h$, then
\[
 \int_0^1\mathbb E[W_{hs}^2\mid\Delta]\,\dd s
       =\frac{\Delta^2}{3}+\frac h6.
\]
Keeping only the bridge variance $h/6$ would miss the squared conditional mean. In the operator calculation we condition on the sigma-field
$\mathscr S_n=\sigma(q_\infty,Z_{\vartheta_{n,j}}:1\le j\le n)$. The path expansion and the kernel derivative give, with the omitted boundary kernel included in $R_n$,
\[
 K_n^B-K=G_n+B_n+R_n,
 \quad G_n\mid\mathscr S_n\sim\mathcal N(0,Q_n),
 \quad \sqrt n\,\|R_n\|_{\mathfrak S_2}\longrightarrow0.
\]
Here $B_n$ is measurable with respect to the samples. Weighted quadratic variation, including estimates for products on distinct intervals, identifies the limit of
$n(\operatorname{tr}Q_n+\|B_n\|_{\mathfrak S_2}^2)$ as the integral above. This is the calculation in Theorem~\ref{S:thm:trace-conditional-complete}.

The covariance operator norm is of smaller order than the conditional second moment, which gives concentration of the squared norm. More precisely,
\[
 \operatorname{tr}Q_n+\|B_n\|_{\mathfrak S_2}^2=O_\omega(n^{-1}),
 \qquad
 \|Q_n\|_{\mathrm{op}}
       =O_{\beta,\omega,\epsilon}(n^{-1-a_\beta+\epsilon}).
\]
For $0<\epsilon<a_\beta$, the second bound is smaller than $n^{-1}$. Conditional Gaussian exponential bounds then give summable failure probabilities after localization, proving
almost sure convergence along all sizes. No independence between sizes is needed. Theorem~\ref{S:thm:covariance-concentration} carries out this step.

It remains to construct the circular ensemble increments from these samples with a smaller error. Conditional on $q_\infty$, radial and angular quantiles preserve the exact
increment distributions. In the relative transfer product, the centered increment errors have size $k^{-3/2}$ and their means have size $k^{-2}$. Centering before summation gives
$(\sum_{j\ge k}j^{-3})^{1/2}=O(k^{-1})$, the same order as the sum of the means. The preceding Brownian product is predictable in the chosen filtration. Applying the weighted
kernel bound and removing the boundary restriction yields $\sqrt n\|K_n-K_n^B\|_{\mathfrak S_2}\to0$. Proposition~\ref{S:prop:7.2} verifies the joint distribution of the increments
and terminal boundary point at every matrix size, and Proposition~\ref{S:prop:circular-brownian-inverse-comparison} gives this estimate.

Finally, the sample enters the Brownian interval when $\beta$ crosses $2$, changing $\kappa_\beta$. With one Brownian path fixed, the coefficient has the form
$F_\omega(\beta)+(\beta-2)_+^2H_\omega(\beta)$, where $F_\omega,H_\omega$ are analytic and $H_\omega>0$. Thus its second derivative has a positive jump. This conclusion concerns
the specified sampling convention and coupling; the inverse asymptotics hold for each fixed $\beta>0$.

\clearpage
\section{Resolvent estimates and operator convergence}\label{sec:resolvent-estimates-and-operator-convergence}

This section bounds operator errors in terms of coefficient errors. For Wigner matrices the issue is compactness and identification of the resolvent; for the quantitative couplings
it is the size of the error. A common inverse identity comes first, followed by the convergence criterion for tridiagonal operators and the estimates adapted to each model. These
retain additive potential errors, partial sums of logarithmic errors, or relative errors in transfer matrices, as appropriate. The boundary conditions are included in the
corresponding domain identifications.

\subsection{Common inverse and resolvent identities}\label{sec:common-inverse-and-resolvent-identities}
We express the resolvent difference in terms of the inverse difference on a common Hilbert space. On the orthogonal complement of the embedded finite space, the inverse is zero and
$I-zK$ is the identity.

Let $\mathscr H$ be a Hilbert space. For a bounded operator $K$ and a complex number $z$ with $I-zK$ invertible, set
\[
 \Phi_z(K)=K(I-zK)^{-1}.
\]
We first prove the algebraic identity for bounded operators. We then identify $\Phi_z(A^{-1})$ on the domain of an unbounded operator $A$, estimate the factors uniformly on compact
spectral sets, and account for the orthogonal complement of an embedded finite-dimensional space.

\begin{proposition}[Resolvent estimates from the inverse]\label{prop:shared-resolvent-transfer}
Let $K$ and $L$ be bounded operators and suppose that $I-zK$ and $I-zL$ are invertible. Then
\begin{equation}\label{eq:shared-resolvent-identity}
 \Phi_z(K)-\Phi_z(L)
 =(I-zK)^{-1}(K-L)(I-zL)^{-1}.
\end{equation}
For any unitarily invariant operator ideal norm $\|\cdot\|_{\mathfrak I}$ for which $K-L$ belongs to the ideal, one has
\begin{align}
 \|\Phi_z(K)-\Phi_z(L)\|_{\mathfrak I}
 &\le \|(I-zK)^{-1}\|\,\|(I-zL)^{-1}\|
          \|K-L\|_{\mathfrak I},\label{eq:shared-resolvent-upper}\\
 \|K-L\|_{\mathfrak I}
 &\le \|I-zK\|\,\|I-zL\|
          \|\Phi_z(K)-\Phi_z(L)\|_{\mathfrak I}.
          \label{eq:shared-resolvent-lower}
\end{align}
In particular, inverse estimates and their lower bounds transfer locally uniformly to resolvents wherever the displayed factors are uniformly bounded. The identity and estimates
also hold with the operator norm.
\end{proposition}
\begin{proof}
Set $A=I-zK$, $B=I-zL$ and $D=\Phi_z(K)-\Phi_z(L)$. Direct multiplication gives $AK=K-zK^2=KA$ and $BL=L-zL^2=LB$. Multiplication on both sides by the respective inverses yields
\[
 KA^{-1}=A^{-1}K,
 \qquad LB^{-1}=B^{-1}L.
\]
Consequently
\[
 ADB=AKA^{-1}B-ALB^{-1}B=KB-AL
       =K-zKL-L+zKL=K-L.
\]
Multiplication by $A^{-1}$ and $B^{-1}$ proves \eqref{eq:shared-resolvent-identity}. The ideal inequality $\|UXV\|_{\mathfrak I} \le\|U\|\|X\|_{\mathfrak I}\|V\|$ now gives
\[
 \|D\|_{\mathfrak I}\le\|A^{-1}\|\|K-L\|_{\mathfrak I}\|B^{-1}\|,
 \qquad
 \|K-L\|_{\mathfrak I}\le\|A\|\|D\|_{\mathfrak I}\|B\|.
\]
These are \eqref{eq:shared-resolvent-upper} and \eqref{eq:shared-resolvent-lower}. The first bound also proves that $D$ belongs to the ideal. Taking suprema of the displayed
factors over a compact spectral set gives the stated uniform bounds. Submultiplicativity proves the same bounds in operator norm.
\end{proof}

\localheading{Unbounded operators}\phantomsection\label{par:shared-operator-domain}
Let $A$ be self-adjoint on $\mathscr H$, suppose $0\in\rho(A)$, and put $K=A^{-1}$. Then $K\mathscr H=\operatorname{dom}A$. If $D_z=I-zK$ is invertible, $KD_z^{-1}$ takes values in
this domain and
\[
 (A-z)KD_z^{-1}f=f\quad(f\in\mathscr H),\qquad
 KD_z^{-1}(A-z)u=KD_z^{-1}D_zAu=u\quad(u\in\operatorname{dom}A).
\]
Thus $KD_z^{-1}=(A-z)^{-1}$. Conversely, if $z\in\rho(A)$, write $R=(A-z)^{-1}$. Applying $K$ to $(A-z)Rf=f$ and applying $R$ to $(A-z)Kf=f-zKf$ gives $R-K=zKR=zRK$. Hence
\[
 (I-zK)(I+zR)=I+z(R-K-zKR)=I,\qquad
 (I+zR)(I-zK)=I+z(R-K-zRK)=I.
\]
This proves
\[
 z\in\rho(A)\ \Longleftrightarrow\ I-zK\text{ is invertible},\qquad
 \Phi_z(K)=(A-z)^{-1},\qquad (I-zK)^{-1}=I+z(A-z)^{-1}.
\]
For nonreal $z$, self-adjointness therefore gives $\|(I-zK)^{-1}\|\le1+|z|/|\operatorname{Im}z|$.

\localheading{Compact spectral sets}\phantomsection\label{par:shared-compact-sets}
Suppose $K_n\to K$ in operator norm. Let $C\subset\mathbb C$ be nonempty and compact, with $I-zK$ invertible for every $z\in C$. For $z$ close to $z_0\in C$,
\[
 I-zK=(I-z_0K)\bigl[I-(z-z_0)(I-z_0K)^{-1}K\bigr].
\]
The second factor has a norm convergent Neumann series when $|z-z_0|\|(I-z_0K)^{-1}K\|<1$. This proves continuity of the inverse, and compactness gives finite constants
\[
 M_C=\sup_{z\in C}\|(I-zK)^{-1}\|,\qquad r_C=\sup_{z\in C}|z|.
\]
Choose $n_0$ so that $r_CM_C\|K_n-K\|\le1/2$ for $n\ge n_0$. With $F_{n,z}=z(I-zK)^{-1}(K_n-K)$, one has
\[
 I-zK_n=(I-zK)(I-F_{n,z}),\qquad
 (I-zK_n)^{-1}=\left(\sum_{j=0}^{\infty}F_{n,z}^j\right)(I-zK)^{-1}.
\]
The series converges uniformly for $z\in C$, since its tail after order $m$ has norm at most $2^{-m}$. Consequently $\sup_C\|(I-zK_n)^{-1}\|\le2M_C$. If $K_n-K\in\mathfrak I$, put
$M_* =\max\{\|K\|,\sup_{n\ge n_0}\|K_n\|\}<\infty$. Applying \eqref{eq:shared-resolvent-upper}--\eqref{eq:shared-resolvent-lower} at each $z\in C$ gives
\[
 \frac{\|K_n-K\|_{\mathfrak I}}{(1+r_CM_*)^2}
 \le\inf_{z\in C}\|\Phi_z(K_n)-\Phi_z(K)\|_{\mathfrak I}
 \le\sup_{z\in C}\|\Phi_z(K_n)-\Phi_z(K)\|_{\mathfrak I}
 \le2M_C^2\|K_n-K\|_{\mathfrak I}.
\]
The same bounds hold in operator norm. For $C=\{0\}$ the difference is exactly $K_n-K$, which proves the estimate also for this compact set. For a pathwise comparison, $n_0$ and
the constants may depend on the outcome and on $C$; this argument supplies no integrability bound on them.

\localheading{Isometric embeddings}\phantomsection\label{par:shared-embeddings}
If $J:\mathbb C^n\to\mathscr H$ is an isometry and $A$ is invertible, then $K=JA^{-1}J^*$ satisfies
\[
 \Phi_z(K)=J(A-z)^{-1}J^*.
\]
Indeed, put $P=JJ^*$ and $D=I_{\mathbb C^n}-zA^{-1}$. Since $J^*J=I_{\mathbb C^n}$ and $J^*(I-P)=0$,
\[
 I-zK=JDJ^*+(I-P),\qquad
 (I-zK)^{-1}=JD^{-1}J^*+(I-P).
\]
Multiplying the two right hand sides gives $P+(I-P)=I$ in both orders. These formulas also show that $I-zK$ is invertible exactly when $D$, equivalently $A-z$, is invertible. As
$D=(A-z)A^{-1}$, we obtain
\[
 K(I-zK)^{-1}=JA^{-1}D^{-1}J^*=J(A-z)^{-1}J^*.
\]
Thus the resolvent is zero on $(\operatorname{ran}J)^\perp$.

\localheading{Compressions}\phantomsection\label{par:shared-compression}
Let $P$ be an orthogonal projection and write a finite self-adjoint matrix on $\operatorname{ran}P\oplus\ker P$ as
\[
 A=\begin{pmatrix}A_{11}&A_{12}\\A_{21}&A_{22}\end{pmatrix}.
\]
Suppose that $A-z$ and $A_{22}-z$ are invertible. The second row of $(A-z)(u,v)^{\mathsf T}=(f,0)^{\mathsf T}$ gives $v=-(A_{22}-z)^{-1}A_{21}u$. Substituting this in the first row
yields
\[
 S_z u=f,\qquad S_z=A_{11}-z-A_{12}(A_{22}-z)^{-1}A_{21}.
\]
If $S_zu=0$, then $(u,-(A_{22}-z)^{-1}A_{21}u)$ is in $\ker(A-z)$, hence $u=0$. Thus $S_z$ is invertible, and
\[
 P(A-z)^{-1}P\big|_{\operatorname{ran}P}=S_z^{-1}.
\]
When $A_{11}-z$ is also invertible, the inverse difference identity gives
\[
 S_z^{-1}-(A_{11}-z)^{-1}
 =S_z^{-1}A_{12}(A_{22}-z)^{-1}A_{21}(A_{11}-z)^{-1}.
\]
For nonreal $z$, self-adjointness of $A$, $A_{11}$ and $A_{22}$ gives the required resolvents. The compressed full resolvent therefore contains the contribution of the eliminated
coordinates displayed above. The Schur complement estimates for the approximation of the limiting resolvent by tridiagonal matrix resolvents control that contribution before using
the inverse identity for the isometrically embedded principal block. The identification maps for the approximation of the limiting inverse by bidiagonal matrix inverses are
isometries on the whole finite space; their positive inverse and spectral gap bounds give the factors on the compact sets stated in the hard edge theorems.

\subsection{Convergence of tridiagonal operators}
\label{sec:a-norm-resolvent-criterion-for-jacobi-operators}
\label{W:sec:convergence-of-tridiagonal-operators}

We prove that local convergence of the coefficients and localization of the full low energy spectral subspace imply convergence of the resolvent in operator norm. The local
coefficients identify the weak equation, and localization gives compactness. The same assumptions yield a common lower spectral bound.

We prove the deterministic criterion before its probabilistic version, then verify the finite Lanczos representation and its measurability.
Section~\ref{W:sec:operator-limit-reduction} identifies the two estimates needed for Wigner matrices; the final deduction follows their proofs in Section~\ref{W:sec:main-assembly}.
Notation for entry comparison is collected in Section~\ref{W:sec:notation-and-scales}.

\subsubsection{The convergence criteria}
\label{W:sec:jacobi-scaling}
\label{W:sec:necessity-of-localization}

\label{W:sec:criteria-for-norm-resolvent-convergence}

The criterion below assumes local coefficient convergence and localization of the low energy spectral subspace. A common lower spectral bound follows from these assumptions.

Let $T_N$ be a real symmetric tridiagonal matrix with diagonal $a_{k,N}$ and nonnegative off-diagonal entries $b_{k,N}$. Set $m_N=N^{1/3}$, $h_N=m_N^{-1}$, and
$\mathcal A_N=m_N^2(2I-T_N)$. The hypotheses below concern the local coefficients and the spectral projections; the smallest eigenvalue of $\mathcal A_N$ is allowed to depend on
$N$ without a uniform lower bound. All scalar products are linear in their first argument, and all resolvents at nonreal parameters are taken on the complex Hilbert spaces.

Use the isometry $\mathcal E_N$ in \eqref{W:eq:piecewise-constant-embedding}. Integration over the disjoint cells gives $\|\mathcal E_Nv\|_2^2=\sum_{k=1}^N|v_k|^2$ and
$\mathcal E_N^*\mathcal E_N=I$. Its range projection is cell averaging:
\[
 (\mathcal E_N\mathcal E_N^*f)(x)
 =h_N^{-1}\int_{(k-1)h_N}^{kh_N}f(t)\,\dd t
 \quad((k-1)h_N\le x<kh_N),
\]
with value zero for $x\ge Nh_N$. For $f\in H^1(0,\infty)$, Cauchy--Schwarz on each cell gives
\[
 \|\mathcal E_N\mathcal E_N^*f-f\|_2^2
 \le h_N^2\|f'\|_2^2+\int_{Nh_N}^{\infty}|f|^2\,\dd x\longrightarrow0.
\]
The projections are contractions, and $H^1(0,\infty)$ is dense in $L^2(0,\infty)$, so the same convergence holds for every $f\in L^2$. For nonreal $z$, extending the finite
operator itself by zero gives
\[
 (\mathcal E_N\mathcal A_N\mathcal E_N^*-zI)^{-1}
 =\mathscr R_N(z)-z^{-1}(I-\mathcal E_N\mathcal E_N^*).
\]
The comparison in the theorem uses $\mathscr R_N$, with zero action on the complementary space.

The quadratic forms are expressed through linear interpolation at the nodes,
\begin{equation}
 \mathcal I_Nv(kh_N)=h_N^{-1/2}v_k,
 \qquad 0\le k\le N+1,\qquad v_0=v_{N+1}=0,
 \label{W:eq:piecewise-linear-embedding}
\end{equation}
with zero extension for $x\ge(N+1)h_N$.

The convergence assumptions concern the cumulative potential and the full spectral projections of $\mathcal A_N$. The auxiliary value at the left boundary is
\begin{equation}
 \bar b_{0,N}=1,\qquad \bar b_{k,N}=b_{k,N}\quad(1\le k<N).
 \label{W:eq:tightness-left-offdiagonal-coefficients}
\end{equation}
This auxiliary value is used in the quadratic form decomposition. The Lanczos recurrence retains its initialization $b_0=0$. Define
\begin{equation}
 \widetilde Y_{N,k}=m_N\sum_{r=1}^k
 (2-a_{r,N}-\bar b_{r-1,N}-b_{r,N}),\qquad 0\le k<N.
 \label{W:eq:tightness-exact-grid-potential}
\end{equation}
The empty sum is zero. Write $\widetilde Y_N^{\rm lin}$ for its linear interpolation and
\begin{equation}
 \widetilde Y_N^{\rm st}(x)
 =\widetilde Y_{N,\min\{N-1,\lfloor m_Nx\rfloor\}},\qquad x\ge0,
 \label{W:eq:tightness-exact-step}
\end{equation}
for its step interpolation. Both are constant after $(N-1)h_N$; for $N=1$ they are identically zero. The off-diagonal step function is
\begin{equation}
 \mathfrak b_N(x)=
 \begin{cases}
 \bar b_{k,N},&kh_N\le x<(k+1)h_N,\quad0\le k<N,\\
 1,&x\ge Nh_N.
 \end{cases}
 \label{W:eq:jacobi-step-coefficient}
\end{equation}
As in \eqref{W:eq:operator-coefficient-tail-notation}, let $B_N(L)$ be the maximum of $|b_{k,N}-1|$ over the off-diagonal coordinates up to $\min\{N-1,\lceil m_NL\rceil+1\}$. For
the coordinate and spectral projections use
\[
 P_{N,L}=\sum_{k=1}^{\min\{N,\lfloor m_NL\rfloor\}}e_ke_k^*,\qquad
 Q_{N,L}=I-P_{N,L},\qquad
 \Pi_N(E)=\mathbf1_{(-\infty,E]}(\mathcal A_N).
\]
The projection $\Pi_N(E)$ includes all invariant blocks of $T_N$, also when some $b_{k,N}$ vanish.

For a continuous real function $Y$ with $Y(0)=0$, define the Hermitian sesquilinear form
\begin{equation}
 \mathfrak q_Y[f,g]
 =\int_0^\infty f'\overline{g'}\,\dd x
       -\int_0^\infty Y(f\overline g)'\,\dd x,
 \qquad f,g\in C_c^\infty(0,\infty).
 \label{W:eq:operator-limit-form}
\end{equation}
The second integral is well-defined even when $Y$ is nowhere differentiable. Put
\[
 \mathcal V=H_0^1(0,\infty)\cap L^2((0,\infty),x\,\dd x),\qquad
 \|f\|_{\mathcal V}^2=\int_0^\infty(|f'|^2+(1+x)|f|^2)\,\dd x.
\]
We assume that there are $c_Y,D_Y>0$ and $C_Y\ge0$ such that
\begin{equation}
 c_Y\|f\|_{\mathcal V}^2
 \le\mathfrak q_Y[f]+C_Y\|f\|_2^2
 \le D_Y\|f\|_{\mathcal V}^2
 \quad(f\in C_c^\infty(0,\infty)).
 \label{W:eq:operator-limit-coercivity}
\end{equation}
Lemma~\ref{W:lem:operator-limit-operator} constructs the associated self-adjoint operator $\mathcal H_Y$ with Dirichlet boundary condition and compact resolvent.
Lemma~\ref{W:lem:operator-brownian-form} verifies \eqref{W:eq:operator-limit-coercivity} for $Y_\beta(x)=x^2/2+2B(x)/\sqrt\beta$ almost surely.

The deterministic coefficient assumptions are, for each finite $L$,
\begin{equation}
 \|\widetilde Y_N^{\rm lin}-Y\|_{C([0,L])}\longrightarrow0,
 \qquad B_N(L)\longrightarrow0.
 \label{W:eq:operator-deterministic-local}
\end{equation}
The spectral localization assumption is, for each fixed $E\in\mathbb R$,
\begin{equation}
 \lim_{L\to\infty}\limsup_{N\to\infty}
              \|Q_{N,L}\Pi_N(E)\|=0.
 \label{W:eq:operator-deterministic-localization}
\end{equation}

\begin{theorem}
\label{W:thm:operator-deterministic}
Under \eqref{W:eq:operator-limit-coercivity}, \eqref{W:eq:operator-deterministic-local}, and \eqref{W:eq:operator-deterministic-localization},
\begin{equation}
 \sup_{z\in K}\|\mathcal E_N(\mathcal A_N-z)^{-1}\mathcal E_N^*
                           -(\mathcal H_Y-z)^{-1}\|\longrightarrow0
 \label{W:eq:operator-deterministic-conclusion}
\end{equation}
for every compact $K\Subset\mathbb C\setminus\mathbb R$. There are constants $C<\infty$ and $N_0<\infty$, depending on the coefficient sequence, such that $\mathcal A_N\ge-CI$ for
every $N\ge N_0$.
\end{theorem}

For random $T_N$, the corresponding assumptions are
\begin{equation}
 \widetilde Y_N^{\rm lin}\Longrightarrow
 Y_\beta=\frac{x^2}{2}+\frac2{\sqrt\beta}B(x)
 \quad\text{in }C_{\rm loc}([0,\infty)),\qquad
 B_N(L)\xrightarrow{\mathbb P}0\quad(L<\infty),
 \label{W:eq:operator-probabilistic-local}
\end{equation}
and, for every $E\in\mathbb R$ and $\eta>0$,
\begin{equation}
 \lim_{L\to\infty}\limsup_{N\to\infty}
       \mathbb P(\|Q_{N,L}\Pi_N(E)\|>\eta)=0.
 \label{W:eq:operator-probabilistic-localization}
\end{equation}
Here $\mathscr R_N=\mathcal E_N(\mathcal A_N-\boldsymbol\cdot)^{-1}\mathcal E_N^*$ and $\mathscr R_\beta=(\mathcal H_\beta-\boldsymbol\cdot)^{-1}$.

\begin{theorem}
\label{W:thm:operator-probabilistic}
Under \eqref{W:eq:operator-probabilistic-local} and \eqref{W:eq:operator-probabilistic-localization},
\begin{equation}
 (\widetilde Y_N^{\rm lin},\mathscr R_N)
 \Longrightarrow(Y_\beta,\mathscr R_\beta)
 \quad\text{in }C_{\rm loc}([0,\infty))\times\mathscr Y,
 \label{W:eq:operator-probabilistic-conclusion}
\end{equation}
and
\begin{equation}
 \lim_{M\to\infty}\limsup_{N\to\infty}
       \mathbb P(\lambda_{\min}(\mathcal A_N)<-M)=0.
 \label{W:eq:operator-bottom-tightness}
\end{equation}
\end{theorem}

The random data in these assumptions are taken jointly, with their given dependence. Assumptions that hold uniformly over deterministic starting vectors yield uniform conclusions
in the sense of \eqref{W:eq:uniformity-in-the-starting-vector}; this is proved in Section~\ref{W:sec:measurability-and-truncation}.
Theorems~\ref{W:thm:integrated-potential-and-offdiagonal-coefficients} and \ref{W:thm:localization-spectral-projection} verify these assumptions for Wigner matrices. For
deterministic data, the local form identity gives uniform $H^1$ bounds on bounded intervals. Spectral localization bounds the remainder outside those intervals. The resulting
approximation of both $\mathscr R_N(z)$ and $\mathscr R_N(\bar z)$ by one finite-dimensional space proves compactness in operator norm. A limit satisfies
\eqref{W:eq:operator-limit-weak-equation}; the domain argument in Lemma~\ref{W:lem:weak-resolvent-equation} identifies it with the resolvent of $\mathcal H_Y$. The lower spectral
bound is deduced separately from the local coercivity estimate and localization at energy zero.

For random data we first prove tightness in the operator norm path space, then identify joint subsequential limits. The localization estimate
\eqref{W:eq:operator-probabilistic-localization} provides the tail control needed for tightness. The local weak equation then identifies each joint subsequential limit.

A single remote coordinate shows why \eqref{W:eq:operator-deterministic-localization} accompanies the local assumptions. Starting from any sequence of Jacobi matrices, set the last
off-diagonal entry to zero and set $a_{N,N}=2+m_N^{-2}$. The coefficients on a fixed interval $[0,L]$ are unchanged whenever $(N-1)h_N>L$. For the new scaled matrix,
$\mathcal A_Ne_N=-e_N$. The unit vector $g_N=\mathcal E_Ne_N$ is supported on $[(N-1)h_N,Nh_N)$, whose left end tends to infinity. For $f\in L^2$,
\[
 |\langle g_N,f\rangle|
 \le\|f\mathbf1_{[(N-1)h_N,\infty)}\|_2\longrightarrow0.
\]
Thus $g_N\rightharpoonup0$, while for fixed $L$ and large $N$,
\begin{equation}
 \|Q_{N,L}\Pi_N(0)\|=1,\qquad
 \|\mathscr R_N(z)g_N\|=|-1-z|^{-1}.
 \label{W:eq:calculation-escaping-state}
\end{equation}
For any compact operator $R$, weak convergence of $g_N$ implies $\|Rg_N\|\to0$: every norm convergent subsequence of $Rg_N$ has weak limit zero. Consequently
\[
 \liminf_N\|\mathscr R_N(z)-R\|\ge|-1-z|^{-1}>0.
\]
Replacing the eigenvalue $-1$ by $-M_N$ with $M_N\to\infty$ makes its contribution to a nonreal resolvent tend to zero, although it still belongs to $\Pi_N(0)$. Localization of the
full interval $(-\infty,E]$ will also provide the lower spectral bound in Lemma~\ref{W:lem:operator-bottom}.

\subsubsection{Discrete forms and interpolation}
\label{W:sec:discrete-quadratic-forms-on-finite-intervals}

The scaling and the choice of cumulative potential can be read from the quadratic form. Let $J_0$ have zero diagonal and unit off-diagonal entries. For
$f\in C_c^\infty([0,\infty))$ with $f(0)=0$, set $v_k=h_N^{1/2}f(kh_N)$ and use zero values beyond its support. Expansion of the quadratic form gives
\begin{equation}
 \langle m_N^2(2I-J_0)v,v\rangle
 =\sum_{k\ge0}h_N\left|
       \frac{f((k+1)h_N)-f(kh_N)}{h_N}\right|^2
 \longrightarrow\int_0^\infty|f'|^2\,\dd x.
 \label{W:eq:calculation-free-form}
\end{equation}
To verify the limit, on each cell the difference quotient equals $h_N^{-1}\int_{kh_N}^{(k+1)h_N}f'(t)\,\dd t$. Its difference from $f'(x)$ is at most $h_N\|f''\|_\infty$ on that
cell. On an interval of length $S$ containing the support this gives an $L^2$ error at most $h_N\sqrt S\|f''\|_\infty$; taking squared norms proves
\eqref{W:eq:calculation-free-form}.

The deterministic coefficients $a_k=0$ and $b_k=1-k/(2N)$ give, with the auxiliary value $\bar b_0=1$,
\begin{equation}
 m_N^2(2-\bar b_{k-1}-b_k)=(k-\tfrac12)h_N,\qquad
 m_N\sum_{j=1}^k(2-\bar b_{j-1}-b_j)=\frac{(kh_N)^2}{2}.
 \label{W:eq:calculation-linear-drift}
\end{equation}
Indeed, $2-\bar b_{j-1}-b_j=(2j-1)/(2N)$ and $\sum_{j=1}^k(2j-1)=k^2$. Thus the lattice scale $h_N=N^{-1/3}$ yields the potential $x$ and the integrated potential $x^2/2$.

The alternative cumulative sum
\begin{equation}
 Y_{N,k}=m_N\sum_{r=1}^k(2-a_{r,N}-2b_{r,N})
 \label{W:eq:tightness-convenient-grid}
\end{equation}
will also occur in the comparison of coefficients. Its relation to the exact potential is
\begin{equation}
 \widetilde Y_{N,k}-Y_{N,k}
 =m_N\sum_{r=1}^k(b_{r,N}-\bar b_{r-1,N})
 =m_N(\bar b_{k,N}-1),\qquad 0\le k<N.
 \label{W:eq:integrated-potential-difference}
\end{equation}
For any coefficient array $\sharp$ with $\bar b_{0,N}^{\sharp}=1$, the sum $\sum_{r=1}^k(b_{r,N}^{\sharp}-\bar b_{r-1,N}^{\sharp})$ telescopes to $\bar b_{k,N}^{\sharp}-1$. Hence
\begin{equation}
 \widetilde Y_{N,k}^{\sharp}-Y_{N,k}^{\sharp}
 =m_N(\bar b_{k,N}^{\sharp}-1).
 \label{W:eq:regularized-potential-difference}
\end{equation}
Linear interpolation is a convex combination of adjacent sampled values. Thus, for fixed $T$ and large $N$,
\begin{equation}
 \|\widetilde Y_N^{\rm lin}-Y_N^{\rm lin}\|_{C([0,T])}
 \le m_N\max_{0\le k\le\lceil m_NT\rceil}|\bar b_{k,N}-1|.
 \label{W:eq:potential-interpolation-difference}
\end{equation}
The step and $\sharp$ versions satisfy the same bound. Replacing $\widetilde Y_N$ by $Y_N$ requires the corresponding coefficient error multiplied by $m_N$ to tend to zero.

We derive the bilinear identity with only the test vector supported in an initial interval. This is the version needed for solutions of the full resolvent equation, which need not
vanish at the right end of that interval. Put $x_k=kh_N$ and let $K<N$. With $u\in\mathbb C^N$ and $v$ supported on $\{1,\ldots,K\}$, write $f_N=\mathcal I_Nu$ and
$g_N=\mathcal I_Nv$.

\Needspace{12\baselineskip}
\begin{lemma}
\label{W:lem:operator-mesh}
For these vectors,
\begin{equation}
 \langle\mathcal A_Nu,v\rangle
 =\int\mathfrak b_N f_N'\overline{g_N'}
       -\int\widetilde Y_N^{\rm st}(f_N\overline{g_N})'.
 \label{W:eq:operator-exact-bilinear}
\end{equation}
For every $v\in\mathbb C^N$,
\begin{equation}
 \|v\|_{\ell^2}^2
 =\|\mathcal I_Nv\|_2^2+\frac{h_N^2}{6}\|(\mathcal I_Nv)'\|_2^2,
 \qquad
 \|\mathcal E_Nv-\mathcal I_Nv\|_2^2
 =\frac{h_N^2}{3}\|(\mathcal I_Nv)'\|_2^2.
 \label{W:eq:operator-mass-interpolation}
\end{equation}
Consequently $\|\mathcal I_Nv\|_2^2\le\|v\|_{\ell^2}^2 \le3\|\mathcal I_Nv\|_2^2$. On a bounded interval,
\begin{equation}
 \|\mathcal E_Nv-\mathcal I_Nv\|_{L^2(0,L)}
 \le h_N\|(\mathcal I_Nv)'\|_{L^2(0,L+h_N)}.
 \label{W:eq:operator-local-interpolation-error}
\end{equation}
Under \eqref{W:eq:operator-deterministic-local}, $\widetilde Y_N^{\rm st}$, $\mathfrak b_N$, and $\mathfrak b_N^{-1}$ converge locally uniformly to $Y$, $1$, and $1$. The
reciprocal is defined locally for all sufficiently large $N$.
\end{lemma}
\begin{proof}
Set $u_0=v_0=0$ and $v_{K+1}=0$, but retain $u_{K+1}$. Expanding gives
\begin{equation}
 \langle\mathcal A_Nu,v\rangle ={}h_N^{-2}\sum_{k=0}^K\bar b_{k,N} (u_{k+1}-u_k)\overline{(v_{k+1}-v_k)} +h_N^{-1}\sum_{k=1}^K (\widetilde Y_{N,k}-\widetilde Y_{N,k-1})u_k\overline{v_k}.
 \label{W:eq:operator-discrete-bilinear}
\end{equation}
Indeed, the coefficient of $u_k\overline{v_k}$ is $h_N^{-2}(\bar b_{k-1,N}+b_{k,N})+ h_N^{-1}(\widetilde Y_{N,k}-\widetilde Y_{N,k-1}) =h_N^{-2}(2-a_{k,N})$; the adjacent
coefficients are $-h_N^{-2}b_{k,N}$. On the $k$th cell, $f_N'=h_N^{-3/2}(u_{k+1}-u_k)$, so the first sum is the derivative integral in \eqref{W:eq:operator-exact-bilinear}.

For an absolutely continuous $\varphi$ on $[0,x_{K+1}]$ vanishing at both boundary points, summation by parts gives
\begin{equation}
 \sum_{k=1}^K(\widetilde Y_{N,k}-\widetilde Y_{N,k-1})\varphi(x_k)
 =-\int_0^{x_{K+1}}\widetilde Y_N^{\rm st}(x)\varphi'(x)\,\dd x.
 \label{W:eq:tightness-stieltjes-identity}
\end{equation}
Apply this to $\varphi=f_N\overline{g_N}$. The product is quadratic on each cell, and the integral of its derivative is the difference of its boundary values. Since
$\varphi(x_k)=h_N^{-1}u_k\overline{v_k}$, this proves \eqref{W:eq:operator-exact-bilinear}, including the term involving $u_{K+1}$.

For the mass identity, write $f(x_k+t h_N)=(1-t)a+tb$ for $0\le t\le1$. Then
\[
 \int_{x_k}^{x_{k+1}}|f|^2
 =\frac{h_N}{3}(|a|^2+\operatorname{Re}(a\overline b)+|b|^2)
 =\frac{h_N}{2}(|a|^2+|b|^2)-\frac{h_N}{6}|b-a|^2.
\]
The trapezoidal terms sum to $\|v\|_{\ell^2}^2$. Since $\int_{x_k}^{x_{k+1}}|f'|^2=h_N^{-1}|b-a|^2$, their difference from the integral is $h_N^2\|f'\|_2^2/6$. On every cell,
including the terminal cell $[Nh_N,(N+1)h_N]$, $\mathcal E_Nv$ is the value of $\mathcal I_Nv$ at the right boundary. Hence
\[
 \mathcal E_Nv(x)-\mathcal I_Nv(x)
 =(x_{k+1}-x)(\mathcal I_Nv)'(x),\qquad x_k\le x<x_{k+1}.
\]
Integration gives the factor $1/3$ in \eqref{W:eq:operator-mass-interpolation}. Also $h_N^2\|(\mathcal I_Nv)'\|_2^2 =\sum_{k=0}^N|v_{k+1}-v_k|^2\le4\|v\|_{\ell^2}^2$, which proves
the mass inequalities. Summing only over cells meeting $[0,L]$ gives \eqref{W:eq:operator-local-interpolation-error}.

Finally, on any fixed $[0,L]$, write $\epsilon_N=\|\widetilde Y_N^{\rm lin}-Y\|_{C([0,L+1])}$. The adjacent successive increments are at most $2\epsilon_N+\omega_Y(h_N;L+1)$, where
$\omega_Y$ is the modulus of continuity of $Y$. It follows that
\begin{equation}
 \|\widetilde Y_N^{\rm lin}-\widetilde Y_N^{\rm st}\|_{L^\infty(0,L)}
 \le2\epsilon_N+\omega_Y(h_N;L+1)\longrightarrow0.
 \label{W:eq:interpolation-jump-bound}
\end{equation}
Convergence of the off-diagonal coefficients and their reciprocals follows from
\[
 \|\mathfrak b_N-1\|_{L^\infty(0,L)}\le B_N(L+1),\qquad
 \|\mathfrak b_N^{-1}-1\|_{L^\infty(0,L)}
 \le\frac{B_N(L+1)}{1-B_N(L+1)},
\]
whose denominator is positive for large $N$.
\end{proof}

When $u=v$ is supported through $K$, the same identity reads
\begin{equation}
 \langle\mathcal A_Nv,v\rangle
 =m_N^2\sum_{k=0}^K\bar b_{k,N}|v_{k+1}-v_k|^2
   +m_N\sum_{k=1}^K(\widetilde Y_{N,k}-\widetilde Y_{N,k-1})|v_k|^2.
 \label{W:eq:tightness-discrete-form-potential}
\end{equation}
The last kinetic term is $m_N^2b_{K,N}|v_K|^2$ even though $v_{K+1}=0$. It cancels the term $-m_N^2b_{K,N}|v_K|^2$ contained in the potential. Omitting this edge without changing
the potential would change the diagonal of the compressed matrix. The auxiliary coefficient at zero cancels in the same way.

For $L_{N,K}=(K+1)h_N$, let $\mathcal I_{N,K}v=\mathcal I_Nv|_{[0,L_{N,K}]}$ and $\mathfrak b_{N,K}=\mathfrak b_N|_{[0,L_{N,K})}$. These restrictions satisfy
\begin{equation}
 \langle\mathcal A_Nv,v\rangle
 =\int_0^{L_{N,K}}\mathfrak b_{N,K}|(\mathcal I_{N,K}v)'|^2
 -\int_0^{L_{N,K}}\widetilde Y_N^{\rm st}(|\mathcal I_{N,K}v|^2)'.
 \label{W:eq:tightness-finite-element-form}
\end{equation}
All these identities hold for any other coefficient array, with its own cumulative potential and the same auxiliary boundary convention.

\subsubsection{The limiting operator}
\label{W:sec:compactness-and-identification-of-the-limit}

The limiting potential need not have a derivative as a function. The form \eqref{W:eq:operator-limit-form} is therefore the definition of the operator. We first close this form and
verify its bounds for Brownian paths. We then show that an $L^2$ solution of the local weak equation belongs to the operator domain. This domain identification allows us to pass
from local weak limits to the full resolvent without imposing a boundary condition at a finite cutoff.

\begin{lemma}
\label{W:lem:operator-limit-operator}
Under \eqref{W:eq:operator-limit-coercivity}, the form $\mathfrak q_Y$ closes on $\mathcal V$ and has $C_c^\infty(0,\infty)$ as a form core. Its associated operator $\mathcal H_Y$
is self-adjoint, bounded below by $-C_Y$, and has compact resolvent.
\end{lemma}
\begin{proof}
Let $(f_j)$ be Cauchy for $\|f\|_{\mathcal V}^2=\|f'\|_2^2+\int_0^\infty(1+x)|f|^2\,\dd x$. Completeness of $H_0^1(0,\infty)$ gives a limit $f$ there, and completeness of
$L^2((1+x)\,\dd x)$ gives a limit $g$ in that space. Both convergences imply convergence in ordinary $L^2$, so $f=g$ almost everywhere. In particular $f\in\mathcal V$, and
\[
 \|f_j-f\|_{\mathcal V}^2
 =\|f_j'-f'\|_2^2+\|(1+x)^{1/2}(f_j-f)\|_2^2\longrightarrow0.
\]
Thus $\mathcal V$ is complete. Under the stated lower and upper bounds, the shifted form norm
\[
 \|f\|_{\mathfrak q,\kappa}^2=\mathfrak q_Y[f]+\kappa\|f\|_2^2
\]
is equivalent to this norm, for a fixed sufficiently large $\kappa$. A sequence Cauchy in the shifted form norm therefore converges in $\mathcal V$ and in that same form norm. This
proves closedness of the form extended from its core. For $f\in\mathcal V$, choose smooth cutoffs $\chi_R,\eta_\delta$ with $\chi_R=1$ on $[0,R]$, $\chi_R=0$ after $R+1$,
$|\chi_R'|\le2$, and $\eta_\delta=0$ on $[0,\delta]$, $\eta_\delta=1$ after $2\delta$, $|\eta_\delta'|\le2/\delta$. Then
\begin{equation}
 \|(1-\chi_R)f\|_{\mathcal V}^2
 \le C\int_R^\infty(|f'|^2+(1+x)|f|^2)\,\dd x\longrightarrow0.
 \label{W:eq:form-core-infinity}
\end{equation}
Since $f(0)=0$, $|f(x)|^2\le x\int_0^x|f'|^2$, and hence
\begin{equation}
 \|\eta_\delta'f\|_2^2
 \le\frac4{\delta^2}\int_0^{2\delta}x\int_0^x|f'(t)|^2\,\dd t\,\dd x
 \le8\int_0^{2\delta}|f'|^2\longrightarrow0.
 \label{W:eq:form-core-zero}
\end{equation}
Thus $\eta_\delta\chi_R f\to\chi_Rf$ in $\mathcal V$ as $\delta\downarrow0$. Mollifying its zero extension approximates it by $C_c^\infty(0,\infty)$ functions in $H^1$, and
therefore in $\mathcal V$, since the supports lie in one compact interval. This proves density of the proposed form core.

On that core put $a_Y[f,g]=\mathfrak q_Y[f,g]+C_Y\langle f,g\rangle$. Positivity and the Cauchy--Schwarz inequality for Hermitian forms give
$|a_Y[f,g]|\le D_Y\|f\|_{\mathcal V}\|g\|_{\mathcal V}$. It therefore extends uniquely to $\mathcal V\times\mathcal V$, with
\begin{equation}
 c_Y\|f\|_{\mathcal V}^2\le a_Y[f,f]\le D_Y\|f\|_{\mathcal V}^2.
 \label{W:eq:closed-form-equivalence}
\end{equation}
Equivalence of these norms proves that the extension is closed. The representation theorem for closed semibounded forms (Theorems~VI.2.1 and~VI.2.6 of \cite{kato1995perturbation})
associates a self-adjoint operator $\mathcal H_Y\ge-C_Y$ with the closed form $\mathfrak q_Y=a_Y-C_Y\langle\cdot,\cdot\rangle$. For $s>C_Y$ and $T=(\mathcal H_Y+s)^{-1}$,
\begin{equation}
 (\mathfrak q_Y+s\langle\cdot,\cdot\rangle)[TF,v]
 =\langle F,v\rangle\qquad(v\in\mathcal V),
 \label{W:eq:form-riesz-inverse}
\end{equation}
so testing with $v=TF$ yields
\begin{equation}
 \|TF\|_2\le(s-C_Y)^{-1}\|F\|_2,\qquad
 \|TF\|_{\mathcal V}^2\le[c_Y(s-C_Y)]^{-1}\|F\|_2^2.
 \label{W:eq:form-inverse-bound}
\end{equation}
The operator domain is
\begin{equation}
 \begin{split}
 \mathcal D(\mathcal H_Y)
 &=\{u\in\mathcal V:\exists F\in L^2,
    \ \mathfrak q_Y[u,v]=\langle F,v\rangle\ (v\in\mathcal V)\},\\
 \mathcal H_Yu&=F.
 \end{split}
 \label{W:eq:form-operator-domain}
\end{equation}
In particular every domain element has zero trace at the origin.

Let $P_{L,J}$ replace a function by its average on each of $J$ equal intervals in $[0,L]$ and by zero after $L$. The interval Poincar\'e inequality and the weighted tail estimate
give
\begin{equation}
 \|f-P_{L,J}f\|_2^2
 \le\bigl((L/J)^2+(1+L)^{-1}\bigr)\|f\|_{\mathcal V}^2.
 \label{W:eq:form-compact-approximation}
\end{equation}
Together with \eqref{W:eq:form-inverse-bound}, this shows that $T$ is a norm limit of the finite-rank operators $P_{L,J}T$, by taking first $L\to\infty$ and then $J/L\to\infty$.
Hence $T$ is compact, and the resolvent identity proves compactness of every resolvent.
\end{proof}

To control the Brownian term, average $B$ over intervals of length one. Both the error in this averaging and the derivative of the averaged path are bounded by local increments of
$B$.

\begin{lemma}
\label{W:lem:operator-brownian-form}
For each fixed $\beta>0$, almost every Brownian path gives $Y_\beta(x)=x^2/2+2B(x)/\sqrt\beta$ satisfying \eqref{W:eq:operator-limit-coercivity}. The associated operator is the
stochastic Airy operator with Dirichlet boundary condition.
\end{lemma}
\begin{proof}
Put $c_\beta=2/\sqrt\beta$ and $Y(x)=x^2/2+c_\beta B(x)$. We first construct an event independent of the test function. For each integer $k\ge0$, let
\[
 O_k=\sup_{s,t\in[k,k+2]}|B(t)-B(s)|.
\]
The inequality $O_k\le2\sup_{0\le t\le2}|B(k+t)-B(k)|$, the reflection principle, and a Gaussian tail bound imply
\[
 \mathbb P(O_k>v)\le4\exp(-v^2/16),\qquad v>0.
\]
Choose a deterministic $A$ with $A^2/16>1$. The series
\[
 \sum_{k\ge0}\mathbb P\{O_k>A\sqrt{\log(k+2)}\}
 \le4\sum_{k\ge0}(k+2)^{-A^2/16}
\]
converges. Borel--Cantelli and continuity on the finitely many remaining intervals give, on one event of probability one,
\[
 C_B:=\sup_{k\ge0}\frac{O_k}{\sqrt{\log(k+2)}}<\infty.
\]
Fix a path in this event. Set $\overline B(x)=\int_x^{x+1}B(s)\,\dd s$. This function is continuously differentiable, with $\overline B'(x)=B(x+1)-B(x)$. If $k\le x\le k+1$, all
values used in the average lie in $[k,k+2]$. Thus, after increasing a deterministic constant,
\[
 |B(x)-\overline B(x)|+|\overline B'(x)|
 \le C C_B\sqrt{\log(2+x)}.
\]
For $f\in C_c^\infty(0,\infty)$, integration by parts for the differentiable average gives
\[
 -\int_0^\infty B(|f|^2)'
 =-2\operatorname{Re}\int_0^\infty(B-\overline B)f'\overline f
       +\int_0^\infty\overline B'|f|^2.
\]
There is no boundary term because $f$ has compact support in $(0,\infty)$. For every $\varepsilon>0$, Young's inequality yields
\[
 \begin{split}
 \left|c_\beta\int B(|f|^2)'\right|
 &\le\frac\varepsilon2\|f'\|_2^2
  +\int\left(C\varepsilon^{-1}c_\beta^2C_B^2\log(2+x)
        +Cc_\beta C_B\sqrt{\log(2+x)}\right)|f(x)|^2\,\dd x.
 \end{split}
\]
Denote the nonnegative coefficient in parentheses by $v_\varepsilon(x)$. Since $v_\varepsilon(x)/x\to0$, there is a finite $R=R(\varepsilon,\beta,B)$ such that
$v_\varepsilon(x)\le\varepsilon x$ for $x\ge R$. On $[0,R]$ its maximum is finite. Consequently
\begin{equation}
 \left|c_\beta\int_0^\infty B(x)(|f(x)|^2)'\,\dd x\right|
 \le\varepsilon\left(\|f'\|_2^2+\int_0^\infty x|f(x)|^2\,\dd x\right)
       +C_{\varepsilon,\beta,B}\|f\|_2^2.
 \label{W:eq:operator-brownian-small-bound}
\end{equation}
The same event works for every $\varepsilon>0$: only the deterministic choice of $R$ and the resulting constant change. Polarization extends the noise form continuously from the
core to
\[
 \mathcal V=H_0^1(0,\infty)\cap L^2((0,\infty),x\,\dd x).
\]
The deterministic part of the form is $\|f'\|_2^2+\int x|f|^2$, because $-\int(x^2/2)(|f|^2)'=\int x|f|^2$ on the core. Taking $\varepsilon=1/2$ in the displayed estimate and
choosing $\kappa>C_{1/2,\beta,B}+1$ gives
\[
 \tfrac12\|f\|_{\mathcal V}^2
 \le\mathfrak q_Y[f]+\kappa\|f\|_2^2
 \le C_{\beta,B}\|f\|_{\mathcal V}^2.
\]
Here $\|f\|_{\mathcal V}^2=\|f'\|_2^2+\int(1+x)|f|^2$. These are the required lower and upper form bounds. The construction uses no differentiation of the Brownian path.
\end{proof}

We will identify limits through a local differential equation. For $u\in H^1_{\rm loc}([0,\infty))$ with $u(0)=0$, its weak form is
\begin{equation}
 \int u'\overline{\varphi'}-\int Y(u\overline\varphi)'
       -z\langle u,\varphi\rangle=\langle f,\varphi\rangle,
 \qquad\varphi\in C_c^\infty(0,\infty).
 \label{W:eq:operator-limit-weak-equation}
\end{equation}
The factor $\varphi$ gives every integral in this formula compact support.

\begin{lemma}
\label{W:lem:weak-resolvent-equation}
Assume \eqref{W:eq:operator-limit-coercivity}. If $f,u\in L^2$, $z\notin\mathbb R$, and $u$ satisfies \eqref{W:eq:operator-limit-weak-equation} with the local regularity and zero
trace specified above, then $u\in\mathcal D(\mathcal H_Y)$ and $u=(\mathcal H_Y-z)^{-1}f$.
\end{lemma}
\begin{proof}
On a bounded interval $I$,
\[
 \left|\int_I u'\overline{v'}-\int_I Y(u\overline v)'\right|
 \le(1+2\|Y\|_{L^\infty(I)})\|u\|_{H^1(I)}\|v\|_{H^1(I)}.
\]
Let $v$ be a compactly supported $H^1$ function with $v(0)=0$. For the cutoff $\eta_\delta$ of \eqref{W:eq:form-core-zero}, $|v(x)|^2\le x\int_0^x|v'|^2$ implies
\[
 \|v-\eta_\delta v\|_{H^1}^2
 \le 18\int_0^{2\delta}|v'|^2\,\dd x
       +\int_0^{2\delta}|v|^2\,\dd x\longrightarrow0.
\]
For fixed $\delta$, mollification of the zero extension approximates $\eta_\delta v$ in $H^1$ by smooth functions supported in a fixed compact subinterval of $(0,\infty)$. The
preceding form bound and $|z\langle u,v\rangle+\langle f,v\rangle| \le(|z|\|u\|_2+\|f\|_2)\|v\|_2$ permit both limits in the weak equation. It therefore holds for $v$. Take
$v=\chi_R^2u$, with $\chi_R$ as in \eqref{W:eq:form-core-infinity}. On its support $u\in H^1$ and $u(0)=0$, so this is an admissible test function. Expansion of the derivative
terms gives
\[
 |(\chi_Ru)'|^2
 =\chi_R^2|u'|^2+2\chi_R\chi_R'\operatorname{Re}(u'\overline u)
                         +|\chi_R'|^2|u|^2.
\]
The first two terms equal $\operatorname{Re}(u'\overline{(\chi_R^2u)'})$, and the two potential terms both contain $(\chi_R^2|u|^2)'$. Hence
\begin{equation}
 \mathfrak q_Y[\chi_Ru]
 =\operatorname{Re}\mathfrak q_Y[u,\chi_R^2u]
                              +\int|\chi_R'|^2|u|^2.
 \label{W:eq:continuous-cutoff-formula}
\end{equation}
The coercivity bound applies to $\chi_Ru$ by the core approximation on its support. The weak equation yields
\[
 c_Y\|\chi_Ru\|_{\mathcal V}^2
 \le\|f\|_2\|u\|_2+(|z|+C_Y+4)\|u\|_2^2.
\]
Since $\chi_R=1$ on $[0,R]$, the same constant bounds
\[
 c_Y\int_0^R(|u'|^2+(1+x)|u|^2)\,\dd x
 \le\|f\|_2\|u\|_2+(|z|+C_Y+4)\|u\|_2^2.
\]
Monotone convergence in the integral on the left gives $\int_0^\infty(|u'|^2+(1+x)|u|^2)<\infty$. The zero trace and the approximations
\eqref{W:eq:form-core-infinity}--\eqref{W:eq:form-core-zero} now give $u\in H_0^1(0,\infty)$, and hence $u\in\mathcal V$. For a compactly supported test function, the extended form
agrees with the local integrals by convergence in $H^1$ on its support. Moreover,
\[
 |\mathfrak q_Y[u,v]-z\langle u,v\rangle-\langle f,v\rangle|
 \le\bigl((D_Y+C_Y+|z|)\|u\|_{\mathcal V}+\|f\|_2\bigr)
                                      \|v\|_{\mathcal V}.
\]
Approximating any $v\in\mathcal V$ by the core extends the equation to all $v$. Formula \eqref{W:eq:form-operator-domain} gives $u\in\mathcal D(\mathcal H_Y)$ and
$(\mathcal H_Y-z)u=f$. If two solutions existed, their difference $w$ would satisfy $(\mathcal H_Y-z)w=0$; the imaginary part of its scalar product with $w$ gives
$(\operatorname{Im}z)\|w\|_2^2=0$. This proves uniqueness.
\end{proof}

\subsubsection{Local estimates for the full resolvent}
\label{W:sec:local-resolvent-estimates}

We localize a solution of the full resolvent equation by multiplying it by a cutoff. This keeps the boundary terms involving the solution at the right boundary of the chosen
interval. Only the self-adjoint resolvent bound is used here; the lower bound on the whole spectrum will be a consequence of this estimate and spectral localization.

Fix $K\Subset\mathbb C\setminus\mathbb R$, $L<\infty$, and suppose
\begin{equation}
 0<b_*\le\mathfrak b_N\le b^*,\qquad
 |\widetilde Y_N^{\rm st}|\le M\quad\text{on }[0,L+4].
 \label{W:eq:operator-local-bounds}
\end{equation}
Choose a real cutoff $\chi$ equal to one on $[0,L+1]$, zero after $L+2$, with $0\le\chi\le1$ and $\|\chi'\|_\infty\le2$. Multiplication of a vector by $\chi$ means multiplication
of its $k$th coordinate by $\chi(kh_N)$.

\begin{lemma}
\label{W:lem:operator-local-energy}
For sufficiently large $N$ and $u=(\mathcal A_N-z)^{-1}F$,
\begin{equation}
 \|\mathcal I_Nu\|_{H^1(0,L)}
       +\|\mathcal I_N(\chi u)\|_{H^1(0,\infty)}
 \le C_{K,L,b_*,b^*,M}\|F\|_{\ell^2},
 \qquad z\in K\cup\overline K.
 \label{W:eq:operator-local-H1}
\end{equation}
\end{lemma}

\begin{proof}
For real cutoff values $\chi_k=\chi(kh_N)$, expansion of one edge gives
\[
 \operatorname{Re}\{(a-b)\overline{(A^2a-B^2b)}\}
 =|Aa-Bb|^2-(A-B)^2\operatorname{Re}(a\overline b).
\]
The diagonal terms are identical on both sides. Therefore, for a cutoff supported before the terminal coordinate,
\begin{equation}
 \langle\mathcal A_N\chi u,\chi u\rangle
 =\operatorname{Re}\langle\mathcal A_Nu,\chi^2u\rangle
   +h_N^{-2}\sum_k\bar b_{k,N}(\chi_{k+1}-\chi_k)^2
                         \operatorname{Re}(u_{k+1}\overline{u_k}).
 \label{W:eq:discrete-cutoff-formula}
\end{equation}
Its last sum is bounded in absolute value by
\[
 b^*\|\chi'\|_\infty^2\sum_k|u_{k+1}u_k|
 \le\frac{b^*\|\chi'\|_\infty^2}{2}
                       \sum_k(|u_{k+1}|^2+|u_k|^2)
 \le4b^*\|u\|_{\ell^2}^2.
\]
Let $w=\mathcal I_N(\chi u)$. Its exact local form and Young's inequality give
\begin{equation}
 \langle\mathcal A_N\chi u,\chi u\rangle
 \ge \frac{b_*}{2}\|w'\|_2^2-
\frac{2M^2}{b_*}\|w\|_2^2.
 \label{W:eq:operator-local-coercive}
\end{equation}
The mass identity implies $\|w\|_2\le\|\chi u\|_{\ell^2}\le\|u\|_{\ell^2}$. Write $d_K=\inf_{z\in K}|\operatorname{Im}z|$. Self-adjointness gives
$\|u\|_{\ell^2}\le d_K^{-1}\|F\|_{\ell^2}$. Since $\mathcal A_Nu=F+zu$, \eqref{W:eq:discrete-cutoff-formula} gives the upper bound
\[
 \langle\mathcal A_N\chi u,\chi u\rangle
 \le\|F\|_{\ell^2}\|u\|_{\ell^2}
        +(\sup_{z\in K}|z|+4b^*)\|u\|_{\ell^2}^2.
\]
Combining the two bounds yields
\[
 \|w'\|_2^2 \le\frac2{b_*}\left[d_K^{-1} +\left(\sup_{z\in K}|z|+4b^*+\frac{2M^2}{b_*}\right)d_K^{-2}\right] \|F\|_{\ell^2}^2, \qquad \|w\|_2^2\le d_K^{-2}\|F\|_{\ell^2}^2.
\]
For large $N$, every cell meeting $[0,L]$ has both cutoff values at its ends equal to one, so $w=\mathcal I_Nu$ there. The constants are unchanged at $\bar z$, proving
\eqref{W:eq:operator-local-H1}.
\end{proof}

\Needspace{6\baselineskip}
\subsubsection{Spectral localization and lower bounds}
\label{W:sec:spectral-localization-compactness}

The local energy estimate controls oscillation on a bounded interval. To control the mass outside it, we split the spectrum at a fixed energy $E$. The proof commutes the spectral
projection $\Pi_N(E)$ through the resolvent and leaves the coordinate projection on the left.

\begin{lemma}
\label{W:lem:operator-resolvent-tail}
For $K\Subset\mathbb C\setminus\mathbb R$, set $d_K=\inf_{z\in K}|\operatorname{Im}z|$ and $e_K=\sup_{z\in K}\operatorname{Re}z$. If $E>e_K$, then
\begin{equation}
 \sup_{z\in K\cup\overline K}
 \|Q_{N,L}(\mathcal A_N-z)^{-1}\|
 \le d_K^{-1}\|Q_{N,L}\Pi_N(E)\|+(E-e_K)^{-1}.
 \label{W:eq:operator-spectral-tail-bound}
\end{equation}
\end{lemma}

Consequently \eqref{W:eq:operator-deterministic-localization} implies
\begin{equation}
 \lim_{L\to\infty}\limsup_N t_N(L;K)=0,\qquad
 t_N(L;K):=\sup_{z\in K\cup\overline K}
                \|Q_{N,L}(\mathcal A_N-z)^{-1}\|.
 \label{W:eq:operator-pathwise-resolvent-tail}
\end{equation}
The probabilistic localization assumption implies
\begin{equation}
 \lim_{L\to\infty}\limsup_N
       \mathbb P(t_N(L;K)>\eta)=0\qquad(\eta>0).
 \label{W:eq:operator-probability-resolvent-tail}
\end{equation}

\begin{proof}
Write $R_N(z)=(\mathcal A_N-z)^{-1}$ and $\Pi=\Pi_N(E)$. The spectral theorem gives $R_N(z)\Pi=\Pi R_N(z)$ and
\[
 \begin{aligned}
 \|Q_{N,L}R_N(z)\|
 &\le\|Q_{N,L}\Pi R_N(z)\|+
                 \|Q_{N,L}(I-\Pi)R_N(z)\|\\
 &\le |\operatorname{Im}z|^{-1}\|Q_{N,L}\Pi\|
                +\sup_{\lambda>E}|\lambda-z|^{-1}\\
 &\le d_K^{-1}\|Q_{N,L}\Pi\|+(E-e_K)^{-1}.
 \end{aligned}
\]
This bound is uniform on $K\cup\overline K$ and proves \eqref{W:eq:operator-spectral-tail-bound}. For $\eta>0$ choose $E>e_K+2/\eta$. Then
\[
 \{t_N(L;K)>\eta\}
 \subseteq\{\|Q_{N,L}\Pi_N(E)\|>d_K\eta/2\}.
\]
The double limit of its probability is zero by the localization hypothesis. For deterministic data the same inequality gives $\lim_{L\to\infty}\limsup_Nt_N(L;K)\le\eta/2$. Letting
$\eta\downarrow0$ proves the corresponding deterministic limit.
\end{proof}

\Needspace{6\baselineskip}

\begin{lemma}
\label{W:lem:operator-bottom}
Suppose \eqref{W:eq:operator-local-bounds} holds on a window containing $[0,L+4]$ and $\|Q_{N,L}\Pi_N(0)\|\le\theta<1$. Then, for large enough $N$ to fit that window,
\begin{equation}
 \lambda_{\min}(\mathcal A_N)
 \ge -\frac{2M^2}{b_*}-\frac{4b^*}{1-\theta^2}.
 \label{W:eq:operator-bottom-bound}
\end{equation}
\end{lemma}

Under the assumptions of Theorem~\ref{W:thm:operator-deterministic}, the resulting $C$ and $N_0$ satisfy $\mathcal A_N\ge-CI$ for all $N\ge N_0$. Under the assumptions of
Theorem~\ref{W:thm:operator-probabilistic} it gives \eqref{W:eq:operator-bottom-tightness}.

\begin{proof}
If $\lambda\le0$ is an eigenvalue, choose a normalized eigenvector $v$. It belongs to $\operatorname{Ran}\Pi_N(0)$, so $\|Q_{N,L}v\|\le\theta$. Use the cutoff from
Lemma~\ref{W:lem:operator-local-energy}. Since it equals one on the first $\lfloor m_NL\rfloor$ coordinates, $\|\chi v\|^2\ge1-\theta^2$. By \eqref{W:eq:discrete-cutoff-formula}
and \eqref{W:eq:operator-local-coercive}, after using $\|\mathcal I_N\chi v\|\le\|\chi v\|$ and discarding the positive derivative term,
\[
 -\frac{2M^2}{b_*}\|\chi v\|^2
 \le\langle\mathcal A_N\chi v,\chi v\rangle
 \le\lambda\|\chi v\|^2+4b^*.
\]
Division by $\|\chi v\|^2$ gives the bound for every $\lambda\le0$; positive eigenvalues already satisfy it. The projection contains the eigenvectors of all invariant blocks.

In the deterministic case, choose $L$ for which $\limsup_N\|Q_{N,L}\Pi_N(0)\|<1/2$. Local convergence on $[0,L+5]$ then supplies fixed $M$, $b_*=1/2$ and $b^*=2$ for all
sufficiently large $N$, proving a common lower bound.

For random data and $\varepsilon>0$, first choose $L$ so that $\limsup_N\mathbb P(\|Q_{N,L}\Pi_N(0)\|>1/2)<\varepsilon/3$. Tightness of $\widetilde Y_N^{\rm lin}$ on this fixed
larger interval provides $M_\varepsilon<\infty$ with
\[
 \limsup_N\mathbb P\bigl(\|\widetilde Y_N^{\rm lin}\|_{C([0,L+5])}
                           >M_\varepsilon\bigr)<\varepsilon/3.
\]
Define
\[
 \Omega_{N,\varepsilon}
 =\{\|Q_{N,L}\Pi_N(0)\|\le1/2\}
  \cap\{B_N(L+5)\le1/2\}
  \cap\{\|\widetilde Y_N^{\rm lin}\|_{C([0,L+5])}\le M_\varepsilon\}.
\]
The step potential is bounded by the supremum of the linear interpolation on the enlarged interval. On $\Omega_{N,\varepsilon}$, \eqref{W:eq:operator-local-bounds} thus holds with
$b_*=1/2$, $b^*=2$, and $M=M_\varepsilon$, for all sufficiently large matrix sizes, so that the enlarged interval lies within the finite partition. The union bound and
$B_N(L+5)\to0$ in probability give $\limsup_N\mathbb P(\Omega_{N,\varepsilon}^c)<\varepsilon$. On this event,
\[
 \lambda_{\min}(\mathcal A_N)\ge
              -4M_\varepsilon^2-\frac{32}{3}.
\]
Thus, for $R_\varepsilon=4M_\varepsilon^2+32/3+1$,
\[
 \limsup_N\mathbb P(\lambda_{\min}(\mathcal A_N)<-R_\varepsilon)
 \le\varepsilon.
\]
The probabilities are decreasing in the threshold $R$, so $\lim_{R\to\infty}\limsup_N\mathbb P(\lambda_{\min}<-R) \le\varepsilon$ for every $\varepsilon>0$.
\end{proof}

This conclusion is stronger than one can infer from nonreal resolvents alone. For example, on $\ell^2(\mathbb N)$ let $Ae_j=je_j$ and let $A^{(n)}$ agree with $A$ except that
$A^{(n)}e_n=-ne_n$. Then
\[
 \|(A^{(n)}-z)^{-1}-(A-z)^{-1}\|
 =\left|\frac1{-n-z}-\frac1{n-z}\right|\longrightarrow0
\]
locally uniformly for nonreal $z$, but $\lambda_{\min}(A^{(n)})=-n$. The localization assumption controls the whole lower spectral half-line and therefore excludes this escape. A
bound on a fixed bounded spectral interval alone does not provide that control.

\Needspace{6\baselineskip}
\subsubsection{Compactness and identification of the resolvent}
\label{W:sec:compactness-of-resolvent-paths}

We approximate the resolvents by operators on one fixed finite-dimensional space. The local $H^1$ bound controls high frequencies, and the preceding spectral estimate controls the
tail in the tridiagonal coordinates. Both the resolvent and its adjoint must be approximated:
\[
 \|R-FRF\|\le\|(I-F)R\|+\|(I-F)R^*\|
\]
for an orthogonal projection $F$. The two terms require finite-rank approximation of the images of the unit ball under both $R$ and $R^*$.

Fix $L\ge1$, put $S=L+3$, and let $F_{S,J}$ be the orthogonal projection onto the first $J$ functions $(2/S)^{1/2}\sin(\pi jx/S)$ on $(0,S)$, extended by zero after $S$. Assume
\eqref{W:eq:operator-local-bounds}.

\begin{lemma}
\label{W:lem:operator-finite-rank}
For all sufficiently large $N$,
\begin{align}
 &\sup_{z\in K}
 \bigl(\|(I-F_{S,J})\mathscr R_N(z)\|
          +\|(I-F_{S,J})\mathscr R_N(z)^*\|\bigr)
 \notag\\
 &\hspace{15mm}\le
 2t_N(L;K)+C_{K,L,b_*,b^*,M}
                  \left(h_N+\frac{S}{\pi(J+1)}\right).
 \label{W:eq:operator-finite-rank-bound}
\end{align}
\end{lemma}

Under \eqref{W:eq:operator-deterministic-local} and \eqref{W:eq:operator-pathwise-resolvent-tail}, the maps $\mathscr R_N|_K$ form a relatively compact subset of
$C(K;\mathcal K(L^2))$ in the uniform operator norm topology.

\begin{proof}
For $\|f\|_2\le1$ put $u=(\mathcal A_N-z)^{-1}\mathcal E_N^*f$ and $w=\mathcal I_N(\chi u)$, with the cutoff from Lemma~\ref{W:lem:operator-local-energy}. That lemma bounds
$\|w'\|_2$ uniformly for $\|f\|_2\le1$. Its value at zero is zero, and it vanishes after $L+2+h_N<S$ for $h_N<1$. Thus $w\in H_0^1(0,S)$. The triangle inequality and
\eqref{W:eq:operator-mass-interpolation} give
\[
 \|\mathcal E_Nu-w\|_2 \le\|\mathcal E_N((1-\chi)u)\|_2 +\|\mathcal E_N(\chi u)-\mathcal I_N(\chi u)\|_2 \le\|Q_{N,L}u\|_{\ell^2}+h_N\|w'\|_2.
\]
Indeed, $1-\chi$ vanishes on the coordinates of $P_{N,L}$, so $\|(1-\chi)u\|\le\|Q_{N,L}u\|$. The remaining difference is bounded by the interpolation estimate applied to the
supported vector $\chi u$. For $e_j(x)=(2/S)^{1/2}\sin(\pi jx/S)$ write $c_j=\langle w,e_j\rangle$. Parseval's identities for $H_0^1(0,S)$ give
\[
 \|(I-F_{S,J})w\|_2^2=\sum_{j>J}|c_j|^2
 \le\left(\frac{S}{\pi(J+1)}\right)^2
        \sum_{j>J}\left(\frac{\pi j}{S}\right)^2|c_j|^2
 \le\left(\frac{S}{\pi(J+1)}\right)^2\|w'\|_2^2.
\]
Taking the supremum over unit $f$ yields the first resolvent bound. The bound holds on $K\cup\overline K$ with the same constants. Substituting $\bar z$ and using
$\mathscr R_N(z)^*=\mathscr R_N(\bar z)$ gives the adjoint bound in \eqref{W:eq:operator-finite-rank-bound}.

For compactness, fix $\varepsilon>0$. Choose first $L$, then $J$, and finally $N_0$ so that the two errors in \eqref{W:eq:operator-finite-rank-bound} have sum at most $\varepsilon$
for $N\ge N_0$. Enlarge $F_{S,J}$ to a finite-rank projection $G$ containing $\operatorname{Ran}\mathcal E_N$ for $N<N_0$. Then, for every $N$,
\[
 \sup_{z\in K}\|\mathscr R_N(z)-G\mathscr R_N(z)G\|
 \le\sup_{z\in K}\bigl(\|(I-G)\mathscr R_N(z)\|
                     +\|(I-G)\mathscr R_N(z)^*\|\bigr)\le\varepsilon.
\]
The resolvent identity gives
\begin{equation}
 \|\mathscr R_N(z)\|\le d_K^{-1},\qquad
 \|\mathscr R_N(z)-\mathscr R_N(w)\|\le d_K^{-2}|z-w|.
 \label{W:eq:operator-resolvent-modulus}
\end{equation}
We prove compactness by constructing finite nets. Let $\mathcal F$ be the family of resolvent paths on the fixed compact set $K\subset\mathbb C\setminus\mathbb R$, and put
$d=\inf_{z\in K}|\operatorname{Im}z|>0$. The resolvent identity gives
\[
 \|R(z)-R(w)\|\le d^{-2}|z-w|,\qquad R\in\mathcal F,
\]
and the same bound holds after compression by an orthogonal projection. The estimates just proved, applied also at conjugate spectral parameters, give for every $\varepsilon>0$ a
finite-rank orthogonal projection $F$ such that
\[
 \sup_{R\in\mathcal F,\ z\in K}\|R(z)-FR(z)F\|<\varepsilon/4.
\]
Indeed, $R-FRF=(I-F)R+FR(I-F)$ and
$\|FR(I-F)\|\le\|(I-F)R^*\|$.
Choose a finite $\delta$-net $z_1,\ldots,z_l$ of $K$, where $2d^{-2}\delta<\varepsilon/4$. All tuples $(FR(z_i)F)_{i\le l}$ lie in a bounded subset of the finite-dimensional space
$\mathcal B(\operatorname{Ran}F)^l$. This subset has a finite $\varepsilon/8$-net for the maximum operator norm. For every net ball that meets the tuples, retain one path from the
family. If $R$ and its retained representative $T$ lie in the same ball, then
\[
 \max_{i\le l}\|F(R(z_i)-T(z_i))F\|<\varepsilon/4.
\]
For $z\in K$ choose $z_i$ within $\delta$. The two Lipschitz estimates imply
\[
 \|F(R(z)-T(z))F\|<\varepsilon/4+2d^{-2}\delta<\varepsilon/2.
\]
Adding the two compression errors gives $\sup_K\|R-T\|<\varepsilon$. Thus finitely many balls of every positive radius cover $\mathcal F$. The space $C(K;\mathcal K(L^2))$ is
complete in the supremum norm, so the closure of this totally bounded family is compact.

\end{proof}

\Needspace{6\baselineskip}

\begin{lemma}
\label{W:lem:operator-identification}
Suppose \eqref{W:eq:operator-deterministic-local} holds, and let $Y$ satisfy \eqref{W:eq:operator-limit-coercivity}. If, along a subsequence, $\mathscr R_N\to\mathscr R_*$ in
$C(K;\mathcal K(L^2))$, then $\mathscr R_*(z)=(\mathcal H_Y-z)^{-1}$ for every $z\in K$.
\end{lemma}

\begin{proof}
Fix $z\in K$ and $f\in L^2$. Put $u_N=(\mathcal A_N-z)^{-1}\mathcal E_N^*f$ and $u=\mathscr R_*(z)f$. Then $\mathcal E_Nu_N\to u$ in $L^2$. The interpolation and local energy
estimates imply, on every fixed interval $[0,L]$,
\[
 \mathcal I_Nu_N\rightharpoonup u\ \text{in }H^1(0,L),\qquad
 \mathcal I_Nu_N\to u\ \text{in }L^2(0,L).
\]
Continuity of the trace on $H^1(0,L)$ gives $u(0)=0$. For $\varphi\in C_c^\infty(0,\infty)$ take $v_{N,k}=h_N^{1/2}\varphi(kh_N)$. Its interpolants satisfy
$\mathcal I_Nv_N\to\varphi$ in $H^1$ and $\mathcal E_Nv_N\to\varphi$ in $L^2$, with supports in one fixed compact interval. The exact sesquilinear identity gives
\begin{equation}
 \begin{split}
 &\int\mathfrak b_N(\mathcal I_Nu_N)'
                         \overline{(\mathcal I_Nv_N)'}
 -\int\widetilde Y_N^{\rm st}
                         (\mathcal I_Nu_N\overline{\mathcal I_Nv_N})'\\
 &\hspace{1cm}-z\langle\mathcal E_Nu_N,\mathcal E_Nv_N\rangle
 =\langle f,\mathcal E_Nv_N\rangle.
 \end{split}
 \label{W:eq:operator-discrete-weak-equation}
\end{equation}
On the support of the test function, $\mathfrak b_N\to1$ and $\widetilde Y_N^{\rm st}\to Y$ uniformly. Write $f_N=\mathcal I_Nu_N$, $g_N=\mathcal I_Nv_N$, and
$Y_N=\widetilde Y_N^{\rm st}$ there. The potential term equals $-\langle f_N',Y_Ng_N\rangle-\langle Y_Nf_N,g_N'\rangle$. Here $f_N'\rightharpoonup u'$, while $Y_Ng_N\to Y\varphi$,
$Y_Nf_N\to Yu$, and $g_N'\to\varphi'$ strongly in $L^2$. Thus both potential terms converge. Weak convergence of one factor and strong convergence of the other also give
convergence of the kinetic term, and strong $L^2$ convergence treats the mass and forcing terms. The limiting identity is \eqref{W:eq:operator-limit-weak-equation}.
Lemma~\ref{W:lem:weak-resolvent-equation} therefore identifies $u=(\mathcal H_Y-z)^{-1}f$. Since $f,z$ were arbitrary, the limiting operator path is the resolvent of
$\mathcal H_Y$.
\end{proof}

\begin{proof}[Proof of Theorem~\ref{W:thm:operator-deterministic}]
Fix $K\Subset\mathbb C\setminus\mathbb R$. For each $L<\infty$, local coefficient convergence and Lemma~\ref{W:lem:operator-mesh} imply, for all sufficiently large $N$,
\[
 \tfrac12\le\mathfrak b_N\le2,\qquad
 \|\widetilde Y_N^{\rm st}\|_{L^\infty(0,L+4)}
 \le1+\|Y\|_{C([0,L+5])}.
\]
Lemma~\ref{W:lem:operator-local-energy} therefore applies with constants independent of $N$. The localization assumption and Lemma~\ref{W:lem:operator-resolvent-tail} give
\[
 \lim_{L\to\infty}\limsup_N t_N(L;K)=0.
\]
For the finite rank projection in Lemma~\ref{W:lem:operator-finite-rank}, its estimate is
\[
 \sup_{z\in K}\bigl(\|(I-F_{L+3,J})\mathscr R_N(z)\|
                 +\|(I-F_{L+3,J})\mathscr R_N(z)^*\|\bigr)
 \le2t_N(L;K)+C_{L,K}\left(h_N+\frac{L+3}{J+1}\right).
\]
Choose $L$, then $J$, and then $N$ to make the right side arbitrarily small. Together with \eqref{W:eq:operator-resolvent-modulus}, that lemma proves relative compactness in
$C(K;\mathcal K(L^2))$.

If this uniform convergence failed, some $\varepsilon>0$ and subsequence would satisfy
\[
 \sup_{z\in K}\|\mathscr R_{N_j}(z)-(\mathcal H_Y-z)^{-1}\|\ge\varepsilon.
\]
A further subsequence converges uniformly on $K$ to $\mathscr R_*$. Lemma~\ref{W:lem:operator-identification} applies to the unchanged local coefficient limits and gives
$\mathscr R_*(z)=(\mathcal H_Y-z)^{-1}$ for every $z\in K$. Hence
\[
 \varepsilon\le
 \sup_{z\in K}\|\mathscr R_{N_j}(z)-(\mathcal H_Y-z)^{-1}\|
 =\sup_{z\in K}\|\mathscr R_{N_j}(z)-\mathscr R_*(z)\|\longrightarrow0,
\]
a contradiction. Applying this conclusion to each set in \eqref{W:eq:resolvent-compact-exhaustion} proves \eqref{W:eq:operator-deterministic-conclusion}.
Lemma~\ref{W:lem:operator-bottom} supplies $C,N_0<\infty$ with $\mathcal A_N\ge-CI$ for $N\ge N_0$.
\end{proof}

\subsubsection{Dirichlet restrictions to bounded intervals}
\label{W:sec:forms-on-bounded-intervals}
\label{W:sec:interpolation-and-interval-convergence}

The preceding proof concerns the full resolvent. We also record the corresponding result for a finite interval with a Dirichlet condition at its right boundary. This requires a
separate argument: for a coordinate projection $P$, in general
\[
 P(\mathcal A_N-z)^{-1}P\big|_{\operatorname{Ran}P}
 \ne\bigl(P\mathcal A_NP\big|_{\operatorname{Ran}P}-z\bigr)^{-1}.
\]
The proof of Theorem~\ref{W:thm:operator-deterministic} applies the equation for the full resolvent before restricting its test functions to $\operatorname{Ran}P$.

Fix $L>0$. Set
\begin{equation}
 n_N(L):=\lfloor m_NL\rfloor,
 \qquad
 K_N(L):=n_N(L)-1,
 \qquad
 \ell_N:=n_N(L)h_N.
 \label{W:eq:quadratic-form-window-mesh}
\end{equation}
Then
\begin{equation}
 0\le L-\ell_N<h_N,
 \qquad
 \ell_N\longrightarrow L.
 \label{W:eq:interval-length-error}
\end{equation}
We suppress $L$ from $K_N(L)$ when no ambiguity is possible.

Let $\mathscr V_N^L\subset H_0^1(0,L)$ be the space of continuous functions which are linear on each interval $[x_k,x_{k+1}]$, $0\le k\le K_N$, vanish at $x_0=0$ and
$x_{K_N+1}=\ell_N$, and are extended by zero on $[\ell_N,L]$. For $v\in\mathbb C^{K_N}$, the map $\mathcal I_{N,L}$ is defined by $\mathcal I_{N,L}v(x_k)=h_N^{-1/2}v_k$ for
$1\le k\le K_N$, zero boundary values, and linear interpolation. Define the isometry
\begin{equation}
 \mathcal E_{N,L}v=h_N^{-1/2}\sum_{k=1}^{K_N}
               v_k\mathbf1_{[(k-1)h_N,kh_N)},\qquad
 \mathcal C_N=\mathcal E_{N,L}\mathcal I_{N,L}^{-1}.
 \label{W:eq:interval-embedding-maps}
\end{equation}
The first map takes values in $L^2(0,L)$ and satisfies $\|\mathcal E_{N,L}v\|_2=\|v\|_{\ell^2}$. For $f\in\mathscr V_N^L$, $\mathcal C_Nf$ equals the value of $f$ at the right end
of each full cell, including the last cell where that value is zero, and is zero on the remaining interval.

Extend the exact off-diagonal coefficient and potential data to the fixed interval $[0,L]$ by
\begin{align}
 \mathfrak b_{N,L}(x)
 &:={}
 \begin{cases}
  \bar b_{k,N},&x\in[x_k,x_{k+1}),\quad 0\le k\le K_N,\\
  1,&x\in[\ell_N,L],
 \end{cases}
 \label{W:eq:quadratic-form-extended-offdiagonal-coefficients}\\
 \mathcal Y_{N,L}(x)
 &:={}
 \begin{cases}
  \widetilde Y_{N,k},&x\in[x_k,x_{k+1}),\quad 0\le k\le K_N,\\
  \widetilde Y_{N,K_N},&x\in[\ell_N,L].
 \end{cases}
 \label{W:eq:quadratic-form-extended-potential}
\end{align}
Functions in $\mathscr V_N^L$ vanish on the last interval. The extension there allows us to state uniform convergence on the fixed domain $[0,L]$.

For this estimate, let $Y_N^{\rm lin}$ and $Y_N^{\rm st}$ denote arbitrary linear and step interpolations of the same real sampled values. Let $Y\in C([0,L])$. Suppose
\begin{equation}
 \eta_N:=\|Y_N^{\rm lin}-Y\|_{L^\infty(0,\ell_N)}
       +\|Y_N^{\rm lin}-Y_N^{\rm st}\|_{L^\infty(0,\ell_N)}\to0.
 \label{W:eq:interval-extension-hypotheses}
\end{equation}
Extend $Y_N^{\rm st}$ by its value on the last full cell to $[\ell_N,L]$. Let $b_N^{\rm ext}$ extend a real function $b_N$ by one on this interval, and suppose
\begin{equation}
 \delta_N=\|b_N-1\|_{L^\infty(0,\ell_N)}\to0.
 \label{W:eq:interval-coefficient-hypothesis}
\end{equation}
Write $\omega_Y(t)=\sup_{|x-y|\le t}|Y(x)-Y(y)|$ on $[0,L]$.

\begin{lemma}
\label{W:lem:extension-to-fixed-interval}
If $\ell_N\to L$, then
\begin{equation}
 \|Y_N^{\rm st}-Y\|_{L^\infty(0,L)}\le\eta_N+\omega_Y(L-\ell_N)\to0.
 \label{W:eq:interval-extension-convergence}
\end{equation}
For $\delta_N\le1/2$ one has $1/2\le b_N^{\rm ext}\le3/2$ and
\begin{equation}
 \|b_N^{\rm ext}-1\|_\infty=\delta_N,
 \qquad\|(b_N^{\rm ext})^{-1}-1\|_\infty
        \le\frac{\delta_N}{1-\delta_N}\le2\delta_N.
 \label{W:eq:interval-reciprocal-estimate}
\end{equation}
\end{lemma}

\begin{proof}
On $[0,\ell_N)$ the triangle inequality bounds the potential error by $\eta_N$. The step function is constant on its last cell. Taking a limit through that cell gives
$|y_N^{\rm term}-Y(\ell_N)|\le\eta_N$. Hence, for $x\in[\ell_N,L]$,
\[
 |y_N^{\rm term}-Y(x)|\le\eta_N+\omega_Y(L-\ell_N).
\]
This proves \eqref{W:eq:interval-extension-convergence}. The coefficient extension is one on the last interval. On the remaining interval, $|b_N^{\rm ext}-1|\le\delta_N$ and
$b_N^{\rm ext}\ge1-\delta_N$, so
\[
 |(b_N^{\rm ext})^{-1}-1|
 =\frac{|b_N^{\rm ext}-1|}{b_N^{\rm ext}}
 \le\frac{\delta_N}{1-\delta_N}.
\]
The two bounds in \eqref{W:eq:interval-reciprocal-estimate} follow.
\end{proof}

For the actual discrete potential, \eqref{W:eq:interpolation-jump-bound} verifies the step interpolation assumption from local uniform convergence of $\widetilde Y_N^{\rm lin}$.
The coefficient error is $\max_{1\le k\le K_N}|b_{k,N}-1|$, since $\bar b_{0,N}=1$. Thus the extensions used in the form converge uniformly to $Y$ and one.

For $f,g\in\mathscr V_N^L$, define
\begin{equation}
 \mathfrak q_{N,L}[f,g]=\int_0^L\mathfrak b_{N,L}f'\overline{g'}\,\dd x-\int_0^L\mathcal Y_{N,L}(f\overline g)'\,\dd x.
 \label{W:eq:quadratic-form-discrete-form}
\end{equation}
Its quadratic part is denoted by $\mathfrak q_{N,L}[f]$. The limiting form is
\begin{equation}
 \mathfrak q_{Y,L}[f,g]=\int_0^Lf'\overline{g'}\,\dd x-\int_0^LY(f\overline g)'\,\dd x,\qquad f,g\in H_0^1(0,L).
 \label{W:eq:quadratic-form-limit-form}
\end{equation}
For $Y=Y_\beta$, this is the form on a bounded interval of the stochastic Airy operator with Dirichlet boundary condition in \eqref{W:eq:sao-form}. We write
$\mathcal H_{Y,L}^{\mathrm D}$ for the associated self-adjoint operator and
\begin{equation}
 \mathcal H_{\beta,L}^{\mathrm D}
 :=\mathcal H_{Y_\beta,L}^{\mathrm D}.
 \label{W:eq:quadratic-form-window-sao}
\end{equation}

For the Dirichlet restriction argument below, assume
\begin{equation}
 \widetilde Y_N^{\rm lin}\longrightarrow Y\quad\text{in }C([0,L]),
 \qquad
 \epsilon_{b,N}=\max_{1\le k\le K_N}|b_{k,N}-1|\longrightarrow0,
 \qquad Y\in C([0,L];\mathbb R).
 \label{W:eq:quadratic-form-minimal-data}
\end{equation}
Only $N$ with $1\le K_N<N$ matter. Let $A_N^{\rm D}$ be the compression of $\mathcal A_N$ to these $K_N$ coordinates. The coefficient $b_{K_N,N}$ is still retained in the form and
canceled by its potential contribution, as explained after \eqref{W:eq:tightness-discrete-form-potential}. For $v,w\in\mathbb C^{K_N}$, put $f=\mathcal I_{N,L}v$ and
$g=\mathcal I_{N,L}w$. Restricting the bilinear identity and polarizing the mass identity of Lemma~\ref{W:lem:operator-mesh} gives
\begin{equation}
 \langle A_N^{\rm D}v,w\rangle=\mathfrak q_{N,L}[f,g],
 \label{W:eq:quadratic-form-exact-form}
\end{equation}
\begin{equation}
 \langle v,w\rangle=\langle f,g\rangle
 +\frac{h_N^2}{6}\langle f',g'\rangle=:\mathfrak m_N[f,g].
 \label{W:eq:quadratic-form-exact-mass}
\end{equation}
Likewise,
\begin{equation}
 \|\mathcal C_Nf-f\|_2^2=\frac{h_N^2}{3}\|f'\|_2^2,
 \qquad \|f\|_2^2\le\mathfrak m_N[f]\le3\|f\|_2^2.
 \label{W:eq:quadratic-form-exact-embedding}
\end{equation}
These formulas include the last cell, where the right boundary value is zero. The variational equation uses the discrete scalar product $\mathfrak m_N$. Its correction
$h_N^2\langle f',g'\rangle/6$ is estimated after obtaining the energy bound.

For $Y\in C([0,L])$, the limiting form has domain $H_0^1(0,L)$: with $M=\|Y\|_\infty$,
\[
 \left|\int_0^L Y(|f|^2)'\right|
 \le2M\|f\|_2\|f'\|_2
 \le\tfrac12\|f'\|_2^2+2M^2\|f\|_2^2.
\]
Its shifted form norm is equivalent to the $H^1$ norm. The finite-rank approximation in the proof of \eqref{W:eq:form-compact-approximation}, now on the fixed interval, shows that
$H_0^1(0,L)\hookrightarrow L^2(0,L)$ is compact. Thus $\mathcal H_{Y,L}^{\rm D}$ has compact resolvent.

\begin{proposition}
\label{W:prop:interval-resolvent-convergence}
Under \eqref{W:eq:quadratic-form-minimal-data}, for every $K\Subset\mathbb C\setminus\mathbb R$,
\begin{equation}
 \sup_{z\in K}\left\|
 \mathcal E_{N,L}(A_N^{\mathrm D}-z)^{-1}\mathcal E_{N,L}^*
              -(\mathcal H_{Y,L}^{\mathrm D}-z)^{-1}\right\|\to0.
 \label{W:eq:interval-resolvent-convergence}
\end{equation}
\end{proposition}
\begin{proof}
Let $b_h=\mathfrak b_{N,L}$ and $Y_h=\mathcal Y_{N,L}$. The extension estimate and the bound for adjacent increments give
\[
 \epsilon_N:=\|b_h-1\|_\infty+\|Y_h-Y\|_\infty\to0.
\]
Choose $M<\infty$ with $\|Y_h\|_\infty,\|Y\|_\infty\le M$ for all sufficiently large $N$. On $\mathscr V_N^L$,
\begin{align}
 |\mathfrak q_{N,L}[f,g]-\mathfrak q_{Y,L}[f,g]|
 &\le\|b_h-1\|_\infty\|f'\|_2\|g'\|_2
       +\|Y_h-Y\|_\infty
             (\|f'\|_2\|g\|_2+\|f\|_2\|g'\|_2)\notag\\
 &\le2\epsilon_N\|f\|_{H^1}\|g\|_{H^1},
 \label{W:eq:interval-form-difference}\\
 \mathfrak q_{N,L}[f]&\ge\tfrac12\|f'\|_2^2-2M\|f\|_2\|f'\|_2
       \ge\tfrac14\|f'\|_2^2-4M^2\|f\|_2^2.
 \label{W:eq:interval-uniform-coercivity}
\end{align}
Here $b_h\ge1/2$ for large $N$. Put $s=4M^2+1$ and write $\mathfrak m_N$ for the discrete scalar product in \eqref{W:eq:quadratic-form-exact-mass}. Since
$\mathfrak m_N[f]\ge\|f\|_2^2$,
\begin{equation}
 \mathfrak q_{N,L}[f]+s\mathfrak m_N[f]
 \ge\tfrac14\|f'\|_2^2+\mathfrak m_N[f]
 \ge\tfrac14\|f\|_{H^1}^2.
 \label{W:eq:interval-shifted-coercivity}
\end{equation}
The same lower bound holds for the limiting shifted form, with $\mathfrak m_N[f]$ replaced by $\|f\|_2^2$.

Suppose $F_N\rightharpoonup F$ in $L^2(0,L)$ and $\|F_N\|_2\le1$. Set
\[
 v_N=(A_N^{\mathrm D}+s)^{-1}\mathcal E_{N,L}^*F_N,
 \qquad f_N=\mathcal I_{N,L}v_N.
\]
The variational equation is
\begin{equation}
 \mathfrak q_{N,L}[f_N,g]+s\mathfrak m_N[f_N,g]
                   =\langle F_N,\mathcal C_Ng\rangle,
 \qquad g\in\mathscr V_N^L.
 \label{W:eq:interval-variational-equation}
\end{equation}
With $g=f_N$, its real part and $\|\mathcal C_Nf_N\|_2=\mathfrak m_N[f_N]^{1/2}$ give
\[
 \tfrac14\|f_N'\|_2^2+\mathfrak m_N[f_N]
 \le\mathfrak m_N[f_N]^{1/2}.
\]
Consequently $\mathfrak m_N[f_N]^{1/2}\le1$ and $\|f_N'\|_2\le2$. A subsequence therefore satisfies $f_N\rightharpoonup f$ in $H_0^1(0,L)$ and $f_N\to f$ in $L^2(0,L)$. The
interpolation identity also gives
\begin{equation}
 \|\mathcal E_{N,L}v_N-f_N\|_2
          =\frac{h_N}{\sqrt3}\|f_N'\|_2\to0.
 \label{W:eq:interval-linear-interpolation}
\end{equation}

Fix $\varphi\in C_c^\infty(0,L)$. For sufficiently large $N$ its support lies before $\ell_N$. Its piecewise linear interpolant $g_N\in\mathscr V_N^L$ satisfies
\[
 \|g_N'-\varphi'\|_2\le h_N\sqrt L\|\varphi''\|_\infty,
 \qquad
 \|g_N-\varphi\|_2\le h_N^2\sqrt L\|\varphi''\|_\infty.
\]
The first estimate follows by comparing the mean of $\varphi'$ on a cell with each value on that cell; integrating the derivative error from a sampling point proves the second. In
particular, $\mathcal C_Ng_N\to\varphi$ in $L^2$. By \eqref{W:eq:interval-form-difference}, $\mathfrak q_{N,L}[f_N,g_N]-\mathfrak q_{Y,L}[f_N,g_N]\to0$. To pass to the limit in the
latter expression, write its potential term as
\[
 -\langle f_N',Yg_N\rangle-\langle Yf_N,g_N'\rangle.
\]
Here $f_N'\rightharpoonup f'$, $Yg_N\to Y\varphi$, $Yf_N\to Yf$, and $g_N'\to\varphi'$ in the indicated weak or strong $L^2$ topologies. Also
\[
 |\mathfrak m_N[f_N,g_N]-\langle f_N,g_N\rangle|
 \le\frac{h_N^2}{6}\|f_N'\|_2\|g_N'\|_2\to0.
\]
The right side of \eqref{W:eq:interval-variational-equation} tends to $\langle F,\varphi\rangle$, since $F_N$ converges weakly and $\mathcal C_Ng_N$ converges strongly. Thus
\[
 \mathfrak q_{Y,L}[f,\varphi]+s\langle f,\varphi\rangle
                       =\langle F,\varphi\rangle.
\]
Given $g\in H_0^1(0,L)$, choose $\varphi_m\in C_c^\infty(0,L)$ with $\|\varphi_m-g\|_{H^1}\to0$. Since $Y$ is bounded on this interval, the difference of the two sides at $g$ and
$\varphi_m$ has absolute value at most
\[
 \bigl[(1+2\|Y\|_\infty+s)\|f\|_{H^1}+\|F\|_2\bigr]
       \|g-\varphi_m\|_{H^1}\longrightarrow0.
\]
The equality therefore holds for $g$. If $f,\widetilde f$ are two solutions, testing their difference by $f-\widetilde f$ gives $\|f-\widetilde f\|_{H^1}^2/4\le0$ by
\eqref{W:eq:interval-shifted-coercivity}. Hence $f=(\mathcal H_{Y,L}^{\mathrm D}+s)^{-1}F$. If the original sequence did not converge strongly to $f$, a subsequence would have
$L^2$ distance at least some $\epsilon>0$ from $f$. Its uniform $H^1$ bound supplies a further strongly convergent subsequence, and the preceding identification forces its limit to
be $f$, a contradiction.

Let
\[
 S_N=\mathcal E_{N,L}(A_N^{\mathrm D}+s)^{-1}\mathcal E_{N,L}^*,
 \qquad S=(\mathcal H_{Y,L}^{\mathrm D}+s)^{-1}.
\]
If $\|S_N-S\|\not\to0$, there are a subsequence, unit vectors $F_N$, and $\varepsilon>0$ such that $\|(S_N-S)F_N\|\ge\varepsilon$. Pass to a subsequence with
$F_N\rightharpoonup F$. The preceding argument gives $S_NF_N\to SF$. Compactness of $S$ gives $SF_N\to SF$, contradicting the inequality. Therefore $\|S_N-S\|\to0$.

For $z\notin\mathbb R$,
\[
 \mathcal E_{N,L}(A_N^{\mathrm D}-z)^{-1}\mathcal E_{N,L}^*
       =S_N[I-(z+s)S_N]^{-1},
 \qquad
 (\mathcal H_{Y,L}^{\mathrm D}-z)^{-1}=S[I-(z+s)S]^{-1}.
\]
On the complement of $\operatorname{Ran}\mathcal E_{N,L}$ the first inverse in brackets is the identity. On that range its norm is at most $1+|z+s|/|\operatorname{Im}z|$, by the
spectral theorem; the same bound holds for the limiting inverse. Direct multiplication gives
\[
 \begin{split}
 &S_N[I-(z+s)S_N]^{-1}-S[I-(z+s)S]^{-1}\\
 &\qquad=[I-(z+s)S_N]^{-1}(S_N-S)[I-(z+s)S]^{-1}.
 \end{split}
\]
In particular, with $M_{K,s}=1+\sup_{z\in K}|z+s|/\inf_{z\in K}|\operatorname{Im}z|$,
\[
 \sup_{z\in K}\left\|
 \mathcal E_{N,L}(A_N^{\rm D}-z)^{-1}\mathcal E_{N,L}^*
             -(\mathcal H_{Y,L}^{\rm D}-z)^{-1}\right\|
 \le M_{K,s}^2\|S_N-S\|\longrightarrow0.
\]
This proves \eqref{W:eq:interval-resolvent-convergence} without changing the discrete scalar product or the identifying map.
\end{proof}

\subsubsection{Resolvent paths and measurable dependence}
\label{W:sec:measurability-and-truncation}

The deterministic convergence result must now be formulated for random resolvents. We use the path space $\mathscr Y$ in \eqref{W:eq:operator-path-space} and prove measurable
dependence on the limiting potential. This is weaker than global norm resolvent continuity, which would require additional control at infinity.

For nonreal $z,z'$, the spectral theorem and the resolvent identity give
\[
 \|\mathscr R_N(z)\|\le |\operatorname{Im}z|^{-1},\qquad
 \|\mathscr R_N(z)-\mathscr R_N(z')\|
 \le \frac{|z-z'|}{|\operatorname{Im}z|\,|\operatorname{Im}z'|}.
\]
The same bounds hold for $\mathscr R_\beta$. The identity uses $\mathcal E_N^*\mathcal E_N=I$ when multiplying the embedded resolvents.

For the space $\mathscr Y$, use the exhaustion \eqref{W:eq:resolvent-compact-exhaustion} and metric \eqref{W:eq:resolvent-path-metric} defined in
Section~\ref{W:sec:intro-limiting-operator}. The next lemma proves completeness and separability and identifies the induced topology.

\begin{lemma}
\label{W:lem:resolvent-path-space}
The metric \eqref{W:eq:resolvent-path-metric} is complete and separable and induces locally uniform convergence in operator norm.
\end{lemma}
\begin{proof}
Let $P_n$ project onto the first $n$ elements of an orthonormal basis. For each compact $A$, strong convergence of $P_n$ is uniform on the relatively compact images of the unit
ball under $A$ and $A^*$. Consequently
\[
 \|A-P_nAP_n\|\le\|(I-P_n)A\|+\|(I-P_n)A^*\|\longrightarrow0.
\]
Finite matrices with rational real and imaginary entries therefore form a countable dense subset of $\mathcal K(\mathcal H)$. This space is complete because compact operators are
closed in operator norm.

For compact metric $K$, the space $C(K;\mathcal K(\mathcal H))$ is complete in the uniform norm. Fix a countable dense set of centers in $K$. For every finite cover by balls
$B(x_i,r_i)$ with such centers and positive rational radii, set
\[
 \phi_i(x)=\frac{(r_i-d(x,x_i))_+}
                  {\sum_j(r_j-d(x,x_j))_+}.
\]
The denominator is positive on $K$, and $\sum_i\phi_i=1$. The functions $\sum_i\phi_i A_i$, with $A_i$ rational finite matrices, form a countable family. Given
$F\in C(K;\mathcal K(\mathcal H))$ and $\varepsilon>0$, choose the radii small enough that $\|F(x)-F(x_i)\|<\varepsilon/2$ on each ball, and choose $\|A_i-F(x_i)\|<\varepsilon/2$.
Then
\[
 \sup_{x\in K}\left\|F(x)-\sum_i\phi_i(x)A_i\right\|
 \le\sup_{x\in K}\sum_i\phi_i(x)
          \bigl(\|F(x)-F(x_i)\|+\|F(x_i)-A_i\|\bigr)<\varepsilon.
\]
This proves separability. Restriction identifies $\mathscr Y$ with the closed subspace
\[
 \{(F_j)_{j\ge2}:F_{j+1}|_{K_j}=F_j\}
 \subset\prod_{j\ge2}C(K_j;\mathcal K(\mathcal H)).
\]
The metric \eqref{W:eq:resolvent-path-metric} is the product metric on this subspace. Each compact subset of $\mathbb C\setminus\mathbb R$ lies in some $K_j$, and each nonreal
point is interior to some $K_j$. The compatible functions therefore define exactly the continuous paths with the topology of local uniform convergence. A closed subspace of a
countable product of Polish spaces is Polish. Since the stated metric is the product metric on this subspace, this proves completeness, separability, and local uniform convergence
in the stated metric.
\end{proof}

Uniformity in a deterministic unit vector means
\begin{equation}
 \sup_{\|w\|=1}\mathbb P(F_N(w))\longrightarrow0
 \label{W:eq:uniformity-in-the-starting-vector}
\end{equation}
for each failure event under consideration. This is equivalent to the same conclusion for every deterministic sequence $w_N$. One implication follows by bounding
$\mathbb P(F_N(w_N))$ by the supremum. Conversely, if the supremum does not tend to zero, there are $\varepsilon>0$, indices $N_j\uparrow\infty$, and deterministic unit vectors
$w_{N_j}$ such that
\[
 \mathbb P(F_{N_j}(w_{N_j}))\ge\varepsilon/2\quad\text{for every }j.
\]
Completing these vectors to a deterministic sequence contradicts convergence for every such sequence. This quantifier is different from taking a union over all vectors: if $V_N$ is
uniform on $\{e_1,\ldots,e_N\}$ and $F_N(w)=\{|w^*V_N|^2>1/2\}$, then
\[
 \mathbb P(F_N(w))=\frac1N\#\{j:|w_j|^2>1/2\}\le\frac1N,
 \qquad
 \mathbb P\Bigl(\bigcup_{\|w\|=1}F_N(w)\Bigr)=1.
\]
The upper bound is attained at each coordinate vector.

Use the bounded Lipschitz metric on probability laws defined in Section~\ref{W:sec:intro-limiting-operator}. The same selection argument applies to laws: convergence to a fixed law
for every deterministic sequence $w_N$ is equivalent to convergence uniformly in $w$ in a bounded Lipschitz metric. Indeed, failure of uniform convergence supplies indices $N_j$
and deterministic choices $w_{N_j}$ for which the distance is bounded away from zero. Completing these choices to a sequence contradicts the assumed convergence. We will apply this
observation after the operator limit is proved.

Let $\mathcal D$ be the set of continuous real functions $Y$ with $Y(0)=0$ that satisfy \eqref{W:eq:operator-limit-coercivity} for some $c_Y,D_Y>0$ and $C_Y\ge0$.

\begin{lemma}
\label{W:lem:operator-measurability}
The set $\mathcal D$ is Borel in $C_{\rm loc}([0,\infty);\mathbb R)$, and $Y\mapsto(\mathcal H_Y-\boldsymbol\cdot)^{-1}$ is Borel from $\mathcal D$ into $\mathscr Y$.
\end{lemma}

\begin{proof}
A single countable family can be chosen before the potential is specified. For each integer $m\ge2$, the Hilbert space $H_0^1((1/m,m))$ is separable and $C_c^\infty((1/m,m))$ is
dense in it. Choose a countable dense subset of this smooth subspace and extend its elements by zero to the half-line. Take the union over $m$ and then all finite linear
combinations with rational real coefficients, or rational real and imaginary parts in the complex case. The resulting family is countable.

If $f\in C_c^\infty(0,\infty)$, choose $m$ with $\operatorname{supp}f\subset(1/m,m)$. The family from this particular interval approximates $f$ in $H^1$, with all approximants
supported in $[1/m,m]$. On this interval,
\[
 \|g\|_{H^1}^2\le\|g\|_{\mathcal V}^2
                 \le(1+m)\|g\|_{H^1}^2.
\]
The approximation is consequently also in $\mathcal V$. Since the compactly supported smooth functions are a form core, the full countable family has dense span in $\mathcal V$.

Enumerate it and perform Gram--Schmidt in $L^2$, omitting each zero remainder. Every retained vector is a finite linear combination of compactly supported smooth functions and
remains smooth with compact support. Every member of the original family is in the finite span of the retained vectors produced by its enumeration index. The resulting orthonormal
sequence therefore has the same dense span. The sequence is chosen independently of the potential.

Denote this orthonormal sequence by $(e_j)$.

Put
\[
 V_j=\operatorname{span}\{e_1,\ldots,e_j\},\quad
 A_j(Y)_{ab}=\mathfrak q_Y[e_b,e_a],\quad
 J_jc=\sum_{a=1}^jc_ae_a,
 \quad R_{Y,j}(z)=J_j(A_j(Y)-z)^{-1}J_j^*.
\]
If these basis functions are supported in $[0,L_j]$, then
\[
 |A_j(Y)_{ab}-A_j(\widetilde Y)_{ab}|
 \le\|Y-\widetilde Y\|_{C([0,L_j])}\|(e_b\overline{e_a})'\|_1.
\]
The resolvent identity consequently makes $Y\mapsto R_{Y,j}$ continuous into $\mathscr Y$.

For a fixed admissible potential choose $\kappa$ so that $\mathfrak a[f,g]=\mathfrak q_Y[f,g]+\kappa\langle f,g\rangle$ is an inner product equivalent to that of $\mathcal V$. Put
$S=(\mathcal H_Y+\kappa)^{-1}$ and $S_m=R_{Y,m}(-\kappa)$. Let $V_m$ be the span of the first $m$ vectors of the fixed form core, and let $u=Sf$ and $u_m=S_mf\in V_m$ solve
\[
 \mathfrak a[u,v]=\langle f,v\rangle\quad(v\in\mathcal V),\qquad
 \mathfrak a[u_m,v]=\langle f,v\rangle\quad(v\in V_m).
\]
Subtraction gives $\mathfrak a[u-u_m,v]=0$ on $V_m$. Hence for every $v\in V_m$,
\[
 \|u-u_m\|_{\mathfrak a}^2
 =\mathfrak a[u-u_m,u-v]
 \le\|u-u_m\|_{\mathfrak a}\|u-v\|_{\mathfrak a}.
\]
Density and the increasing trial spaces prove $S_mf\to Sf$ in $\mathcal V$ for each fixed $f$. Testing with $u_m$ also gives $\|S_mf\|_{\mathcal V}\le C_Y\|f\|_2$ uniformly in $m$.

To prove operator-norm convergence, take an arbitrary bounded sequence $f_m$ and pass to a subsequence with $f_m\rightharpoonup f$ in $L^2$. The preceding estimate and the compact
inclusion $\mathcal V\hookrightarrow L^2$ give a further subsequence such that $S_mf_m$ converges strongly in $L^2$ and weakly in $\mathcal V$ to some $u$. Testing its equation
against any fixed core vector, which belongs to $V_m$ for all sufficiently large $m$, shows $\mathfrak a[u,v]=\langle f,v\rangle$. Density yields $u=Sf$. Also $Sf_m\to Sf$ strongly
in $L^2$ because $S$ is compact. Thus $\|(S_m-S)f_m\|_2\to0$ along this further subsequence. Were $\|S_m-S\|$ not tending to zero, unit vectors witnessing a fixed positive lower
bound would contradict this conclusion. Therefore $S_m\to S$ in operator norm.

The finite inverse at a nonreal parameter is the bounded transform
\[
 S_m[I-(z+\kappa)S_m]^{-1},
\]
with the denominator equal to the identity on $V_m^\perp$. The corresponding transform of $S$ is $(\mathcal H_Y-z)^{-1}$. The inverse identity and a Neumann series transfer the
norm convergence to these resolvents. The shift $\kappa$ is used only to prove convergence for the fixed path: each unshifted finite matrix and its nonreal inverse were already
defined measurably without a path-dependent choice. Their pointwise operator-norm limits are therefore Borel. Convergence in the resolvent path space follows from the resolvent
identity and a countable compact exhaustion of the nonreal plane.

Let $\mathcal Q$ be the rational complex span of the chosen smooth family. The set of admissible potentials is
\[
 \mathcal D=\{Y:Y(0)=0\}\cap
 \bigcup_{\substack{c,D\in\mathbb Q_{>0}\\C\in\mathbb Q_{\ge0}}}
 \bigcap_{\varphi\in\mathcal Q}
 \{Y:c\|\varphi\|_{\mathcal V}^2
 \le\mathfrak q_Y[\varphi]+C\|\varphi\|_2^2
 \le D\|\varphi\|_{\mathcal V}^2\}.
\]
Each set inside the intersection is closed in the topology of local uniform convergence. Approximation with fixed compact support extends these inequalities to all smooth test
functions. Conversely, rational choices $0<c<c_Y$, $C\ge C_Y$, $D>D_Y+C-C_Y$ follow from the assumed bounds. Thus the displayed set is exactly the domain in the statement and is
Borel. On $\mathcal D$ the limiting resolvent path is a pointwise limit of the continuous maps $R_{Y,j}$ into the metric space $\mathscr Y$, and is therefore Borel. Its zero
extension on $\mathcal D^c$ is Borel as well.
\end{proof}

\Needspace{6\baselineskip}

\subsubsection{Convergence in distribution}
\label{W:sec:operator-convergence-in-distribution}

The local coefficient norms are bounded in probability. Together with the estimate for the tail in the tridiagonal coordinates, the finite-rank approximation proves tightness of
the full resolvent paths. We then apply the Skorokhod representation theorem to the joint law of the coefficients and resolvents and identify their almost sure limits by the
discrete weak equation.
\begin{lemma}
\label{W:lem:operator-probability-compactness}
Under \eqref{W:eq:operator-probabilistic-local} and \eqref{W:eq:operator-probabilistic-localization}, the laws of the embedded resolvent paths are tight in
$C_{\mathrm{loc}}(\mathbb C\setminus\mathbb R;\mathcal K(L^2))$.
\end{lemma}

\begin{proof}
Fix $K\Subset\mathbb C\setminus\mathbb R$ and $\varepsilon>0$. We construct a compact $\mathcal C_K\subset C(K;\mathcal K(L^2))$ with
$\inf_N\mathbb P(\mathscr R_N|_K\in\mathcal C_K)\ge1-\varepsilon$. For each $j\ge1$ choose small positive numbers $\delta_j$ and $p_j$, with $\delta_j=2^{-j-4}$ and
$p_j=\varepsilon2^{-j-3}$. By \eqref{W:eq:operator-probability-resolvent-tail}, choose $L_j\ge1$ and an index after which
\[
 \mathbb P(t_N(L_j;K)>\delta_j)\le p_j.
\]
The $L_j$ may be made increasing. On this fixed window, tightness of $\widetilde Y_N^{\mathrm{lin}}$ and convergence of $B_N(L_j+5)$ give a finite $M_j$ and a further index after
which
\[
 \mathbb P\left(\frac12\le\mathfrak b_N\le2,
         \ |\widetilde Y_N^{\mathrm{st}}|\le M_j
               \text{ on }[0,L_j+4]\right)\ge1-p_j.
\]
Let $C_j$ be the constant in \eqref{W:eq:operator-finite-rank-bound} for this window. Choose $J_j$ so that $C_j(L_j+3)/(\pi(J_j+1))\le\delta_j$, and then $N_j$ large enough for the
two probability bounds, the window inclusion, and $C_jh_N\le\delta_j$. Let $G_j$ be the projection onto a finite-dimensional space containing $\operatorname{Ran}F_{L_j+3,J_j}$ and
all $\operatorname{Ran}\mathcal E_N$ for $N<N_j$. Enlarging a projection only decreases the two approximation errors. For $N\ge N_j$, the preceding bounds and
\eqref{W:eq:operator-finite-rank-bound} show that
\[
 \sup_{z\in K}\bigl(\|(I-G_j)\mathscr R_N(z)\|
               +\|(I-G_j)\mathscr R_N(z)^*\|\bigr)\le2^{-j}
\]
except on an event of probability at most $2p_j$. For $N<N_j$ the error is identically zero, because both ranges lie in $\operatorname{Ran}G_j$. For every $N$, the probability that
one of these estimates fails is at most $\sum_{j\ge1}2p_j=\varepsilon/4$. Thus the same sequence of deterministic projections works with probability at least $1-\varepsilon$,
uniformly in $N$.

Let $\mathcal C_K$ consist of the continuous paths taking values in compact operators $F$ satisfying the norm and Lipschitz bounds \eqref{W:eq:operator-resolvent-modulus} and, for
every $j$,
\[
 \sup_{z\in K}\bigl(\|(I-G_j)F(z)\|+\|(I-G_j)F(z)^*\|\bigr)\le2^{-j}.
\]
These conditions are closed in $C(K;\mathcal K)$. Moreover,
\[
 \sup_{z\in K}\|F(z)-G_jF(z)G_j\|\le2^{-j}.
\]
Fix an approximation tolerance and choose $j$ so that $2^{-j}$ is smaller than this tolerance. The compressed paths take values in a bounded finite-dimensional matrix ball and have
the common Lipschitz bound \eqref{W:eq:operator-resolvent-modulus}. The compactness argument in Lemma~\ref{W:lem:operator-finite-rank} therefore gives a finite net for these
compressed paths, and the displayed compression error gives one for $\mathcal C_K$. The set is closed and the ambient space is complete, so it is compact. The preceding event
estimates imply
\[
 \inf_N\mathbb P(\mathscr R_N|_K\in\mathcal C_K)
 \ge1-\sum_{j\ge1}2p_j\ge1-\varepsilon.
\]
This proves tightness on $K$.

For the full space $\mathscr Y$, use the exhaustion \eqref{W:eq:resolvent-compact-exhaustion}. For each $K_j$ choose a compact set $\mathcal C_j$ with exceptional probability at
most $\varepsilon2^{-j-1}$, and put
\[
 \mathcal C=\{F\in\mathscr Y:F|_{K_j}\in\mathcal C_j\text{ for every }j\}.
\]
Every sequence in $\mathcal C$ has, by successive subsequence choices, a subsequence converging uniformly on each $K_j$. The limits agree on overlaps, define a continuous path on
$\mathbb C\setminus\mathbb R$, and belong to the closed sets $\mathcal C_j$. The subsequence thus converges in the metric of $\mathscr Y$ to an element of $\mathcal C$. This proves
compactness. Finally,
\[
 \sup_N\mathbb P(\mathscr R_N\notin\mathcal C)
 \le\sum_{j\ge2}\sup_N\mathbb P(\mathscr R_N|_{K_j}\notin\mathcal C_j)
 \le\varepsilon.
\]
This is tightness in $\mathscr Y$.
\end{proof}

\begin{proof}[Proof of Theorem~\ref{W:thm:operator-probabilistic}]
The resolvent paths are tight by Lemma~\ref{W:lem:operator-probability-compactness}. To retain all finite matrix identities when taking a joint limit, let $b_N^{\rm lin}$
interpolate $1,b_{1,N},\ldots,b_{N-1,N},1$ at $0,h_N,\ldots,Nh_N$, with value one afterwards, and regard all coefficients as the element
\[
 \mathbf t_N=(a_{1,N},b_{1,N},\ldots,b_{N-1,N},a_{N,N},0,\ldots)
 \in\mathbb R^{\mathbb N}.
\]
The local assumptions give $b_N^{\rm lin}\to1$ in probability in $C_{\rm loc}$. For a fixed index $k$, local convergence of the adjacent coefficients gives $b_{k,N}\to1$ in
probability. The definition of the integrated potential gives
\[
 a_{k,N}=2-\bar b_{k-1,N}-b_{k,N}
               -h_N(\widetilde Y_{N,k}-\widetilde Y_{N,k-1}),
 \qquad \bar b_{0,N}=1.
\]
Tightness of the polygonal potentials on any fixed interval containing the first $k$ cells implies tightness of both values in the last difference. Since $h_N\to0$, their
difference multiplied by $h_N$ tends to zero in probability. Thus $a_{k,N}\to0$ in probability. More explicitly, if $V_N$ denotes the supremum of the absolute value of the
polygonal potential on that fixed interval, then for $\varepsilon>0$ and $M>0$,
\[
 \mathbb P(2h_NV_N>\varepsilon)
 \le\mathbb P(V_N>M)+\mathbf1_{\{2h_NM>\varepsilon\}}.
\]
First take the upper limit in $N$ and then let $M\to\infty$.

Extend the finite arrays by any fixed deterministic convention and enumerate their coordinates as $\mathbf t_N=(a_{1,N},b_{1,N},a_{2,N},b_{2,N},\ldots)$. In the product metric
\[
 d(t,s)=\sum_{j\ge1}2^{-j}\min\{1,|t_j-s_j|\},
\]
the limiting array is $\mathbf t_*=(0,1,0,1,\ldots)$. For any $\varepsilon>0$, choose $J$ with $2^{-J}<\varepsilon/2$. Then
\[
 \mathbb P\{d(\mathbf t_N,\mathbf t_*)>\varepsilon\}
 \le\sum_{j=1}^J\mathbb P\{|(\mathbf t_N)_j-(\mathbf t_*)_j|>\varepsilon/2\}\longrightarrow0.
\]
The extension convention does not affect a fixed coordinate for all sufficiently large $N$.

Tightness of the marginals proves tightness of
\begin{equation}
 (\widetilde Y_N^{\rm lin},b_N^{\rm lin},\mathbf t_N,\mathscr R_N)
 \quad\text{in}\quad
 C_{\rm loc}(\mathbb R_+)^2\times\mathbb R^{\mathbb N}\times\mathscr Y.
 \label{W:eq:operator-joint-state}
\end{equation}
These are Polish spaces. Prokhorov's theorem and the Skorokhod representation theorem therefore give, from every subsequence, a further subsequence with almost sure joint
convergence. Its limiting potential $Y$ has the law of $Y_\beta$, its off-diagonal path is one, and its resolvent limit is a compact operator path $\mathscr R_*$. The form estimate
for the Brownian potential makes $Y$ admissible almost surely.

For each fixed $N$, the two interpolants are continuous functions of the finite row $\mathbf t_N$, and its resolvent map is continuous:
\[
 \sup_{z\in K}\|\mathcal E_N(m_N^2(2I-T)-z)^{-1}\mathcal E_N^*
 -\mathcal E_N(m_N^2(2I-T')-z)^{-1}\mathcal E_N^*\|
 \le m_N^2d_K^{-2}\|T-T'\|.
\]
Their graphs are Borel. The identities that define these functions, the zero coordinates after $2N-1$, and the nonnegative off-diagonal entries therefore hold for the copies with
probability one. Intersect over $N$. The copied resolvents are now the resolvents of the copied matrices, with exactly the same local coefficient identities.

The copied potential and off-diagonal paths converge locally uniformly. Lemma~\ref{W:lem:operator-mesh} gives the corresponding step function limits. Since the copied resolvents
already converge in norm, Lemma~\ref{W:lem:operator-identification}, applied on each $K_j$, identifies
\[
 \mathscr R_*(z)=(\mathcal H_Y-z)^{-1}\qquad(z\notin\mathbb R)
\]
almost surely. Localization was used to prove tightness; this identification uses only the local coefficient convergence and the already convergent resolvents. By
Lemma~\ref{W:lem:operator-measurability}, the right side is a Borel function of $Y$. Every subsequential joint law is therefore the law of $(Y_\beta,\mathscr R_\beta)$, proving
\eqref{W:eq:operator-probabilistic-conclusion}. Lemma~\ref{W:lem:operator-bottom} proves \eqref{W:eq:operator-bottom-tightness}.

If the assumptions hold uniformly over deterministic starting vectors, the argument applies to every deterministic sequence $w_N$. The selection argument in
Section~\ref{W:sec:measurability-and-truncation} then proves uniform convergence of the joint laws. The choices in the proof of Lemma~\ref{W:lem:operator-bottom} can also be made
uniformly: for each $\varepsilon>0$ they provide $R_\varepsilon<\infty$ with
\[
 \limsup_N\sup_{\|w\|=1}
 \mathbb P(\lambda_{\min}(\mathcal A_N(w))<-R_\varepsilon)
 \le\varepsilon.
\]
Letting $\varepsilon\downarrow0$ gives uniform lower spectral tightness.
\end{proof}

For Wigner matrices, Section~\ref{W:sec:green-function-comparison} will provide both local coefficient assumptions, and
Section~\ref{W:sec:localization-of-the-low-energy-spectral-subspace} will provide spectral localization. To apply the criterion, we first verify that the continued Lanczos
construction represents the full matrix and is measurable.

\subsubsection{Lanczos representation and truncation}
\label{W:sec:lanczos-tridiagonalization}

The first Krylov subspace may be a proper subspace of the matrix space. The deterministic continuation in \eqref{W:eq:intro-lanczos-continuation} restarts the recurrence whenever
$r_k=0$ and $k<N$, and includes all remaining invariant blocks. This gives the full Lanczos matrix needed for the operator criterion.

Fix a Hermitian $H$ and a unit vector $w$, and use the Krylov subspaces, recurrence, and continuation defined in
\eqref{W:eq:lanczos-krylov-subspace}--\eqref{W:eq:intro-lanczos-continuation}. We write $U=(q_1,\ldots,q_N)$ and $T=U^*HU$ when the matrix size is fixed. Scalar products are linear
in the first argument: $\langle u,v\rangle=v^*u$. The relation with orthogonal polynomials and the spectral measure at $w$ is described in \cite{krishnapur2016universality}.

The Krylov subspaces increase with $r$. Once $\mathcal K_{r+1}=\mathcal K_r$, multiplication by $H$ preserves $\mathcal K_r$, and every later Krylov subspace equals it.
Consequently there is a unique $1\le d\le N$ such that
\begin{equation}
 \dim\mathcal K_r=\min\{r,d\},\qquad r\ge0.
 \label{W:eq:lanczos-krylov-dimension}
\end{equation}

We use the convention that a Jacobi matrix has strictly positive off-diagonal entries, as in \cite{krishnapur2016universality}. The tridiagonal matrix $T$ above has nonnegative
off-diagonal entries. If $b_k=0$, its first $k$ coordinates and their orthogonal complement are invariant. Splitting at every such index writes $T$ as an orthogonal direct sum of
Jacobi matrices; a block of size one is permitted. The estimates below are stated for the full tridiagonal matrix.

\begin{proposition}
\label{W:prop:lanczos-construction}
The recurrence \eqref{W:eq:lanczos-recurrence} with \eqref{W:eq:intro-lanczos-continuation} produces an orthonormal basis. The matrix $T=U^*HU$ is real symmetric, tridiagonal, and
unitarily equivalent to $H$. Its zero off-diagonal entries separate invariant subspaces of $H$.
\end{proposition}
\begin{proof}
Suppose $q_1,\ldots,q_k$ are orthonormal and the recurrence has been verified through index $k-1$. For $i\le k-2$, self-adjointness gives
\[
 q_i^*Hq_k=(Hq_i)^*q_k
 =\bigl(\mathbf1_{\{i\ge2\}}b_{i-1}q_{i-1}
                     +a_iq_i+b_iq_{i+1}\bigr)^*q_k=0.
\]
For $k\ge2$ it also gives $q_{k-1}^*Hq_k=b_{k-1}$, whereas $q_k^*Hq_k=a_k$. Consequently
\begin{equation}
 q_i^*r_k=0\quad(1\le i\le k),\qquad
 r_k=(I-\widehat P_k)Hq_k.
 \label{W:eq:lanczos-residual-orthogonality}
\end{equation}
If $b_k>0$, the vector $r_k/b_k$ is therefore a unit vector orthogonal to all its predecessors. If $b_k=0$ and $k<N$, the normalized projected coordinate vector in
\eqref{W:eq:intro-lanczos-continuation} has these same properties. Such a vector exists, since
\[
 \sum_{j=1}^N\|(I-\widehat P_k)e_j\|^2
 =\operatorname{Tr}(I-\widehat P_k)=N-k>0.
\]
Its denominator is positive by the definition of $j_k$, and its numerator is orthogonal to $q_1,\ldots,q_k$. In either case,
\[
 Hq_k=\mathbf1_{\{k\ge2\}}b_{k-1}q_{k-1}+a_kq_k+b_kq_{k+1}.
\]
When $b_k=0$, the last term is zero regardless of the choice of the next vector. This proves the induction, beginning with $q_1=w$. At the last coordinate, the same scalar product
calculation gives
\[
 Hq_N-a_Nq_N-\mathbf1_{\{N\ge2\}}b_{N-1}q_{N-1}
                    \perp\operatorname{span}\{q_1,\ldots,q_N\}.
\]
The constructed vectors span the entire space, so the displayed vector is zero. The last recurrence therefore ends at $q_N$.

Taking scalar products of these identities with $q_i$ gives
\[
 q_i^*Hq_k=
 \begin{cases}
 a_k,&i=k,\\
 b_k,&i=k+1,\\
 b_{k-1},&i=k-1,\\
 0,&|i-k|>1.
 \end{cases}
\]
Thus $U^*U=UU^*=I$ and $U^*HU=T$ with the entries displayed above. If $b_k=0$, each $Hq_i$ for $i\le k$ belongs to $\operatorname{span}\{q_1,\ldots,q_k\}$. Its orthogonal
complement is also invariant, since for $v$ in that complement and $u$ in the span, $\langle Hv,u\rangle=\langle v,Hu\rangle=0$. Applying this observation at every zero $b_k$
writes $T$ as an orthogonal direct sum of scalar Jacobi matrices with positive interior off-diagonal entries. Finally, $\det(T-zI)=\det(U^*(H-zI)U)=\det(H-zI)$, so all eigenvalues
are preserved with multiplicity.
\end{proof}

Only the first invariant block is generated by $w$. In the next lemma, set $r_N=Hq_N-a_Nq_N-\mathbf1_{\{N\ge2\}}b_{N-1}q_{N-1}$, as in the final recurrence identity above.

\begin{lemma}
\label{W:lem:lanczos-first-block}
The smallest index $k$ for which $r_k=0$ is $d$ from \eqref{W:eq:lanczos-krylov-dimension}. For $1\le r\le d$,
\begin{equation}
 \operatorname{span}\{q_1,\ldots,q_r\}=\mathcal K_r(H,w),\qquad
 q_r=\pi_{r-1}(H)w,
 \label{W:eq:lanczos-polynomial-vectors}
\end{equation}
where $\pi_{r-1}$ has degree $r-1$ and positive leading coefficient $(b_1\cdots b_{r-1})^{-1}$. Moreover,
$d=\#\{\lambda\in\sigma(H):\mathbf1_{\{\lambda\}}(H)w\ne0\}$.
\end{lemma}
\begin{proof}
Set $\pi_{-1}=0$ and $\pi_0=1$. As long as $b_k>0$, define
\[
 \pi_k(x)=\frac{(x-a_k)\pi_{k-1}(x)
                    -\mathbf1_{\{k\ge2\}}b_{k-1}\pi_{k-2}(x)}{b_k}.
\]
Induction in \eqref{W:eq:lanczos-recurrence} yields $q_{k+1}=\pi_k(H)w$, and the leading coefficient is divided by $b_k$ at each step. The change of basis from $1,x,\ldots,x^{r-1}$
to $\pi_0,\ldots,\pi_{r-1}$ is triangular with nonzero diagonal. Hence its image under $p\mapsto p(H)w$ has the same span, proving \eqref{W:eq:lanczos-polynomial-vectors} before
breakdown. If breakdown occurred at $k<d$, invariance from Proposition~\ref{W:prop:lanczos-construction} would imply $\dim\mathcal K_{k+1}=k$, a contradiction. The vector $r_d$
belongs to $\mathcal K_{d+1}=\mathcal K_d$ and is orthogonal to $\mathcal K_d$, so it is zero.

For the dimension formula, write $v_\lambda=\mathbf1_{\{\lambda\}}(H)w$ and $S_w=\{\lambda:v_\lambda\ne0\}$. The vectors $(v_\lambda)_{\lambda\in S_w}$ are nonzero and mutually
orthogonal, and
\[
 H^jw=\sum_{\lambda\in S_w}\lambda^jv_\lambda.
\]
It follows that every Krylov subspace lies in their span. Conversely, for $\lambda\in S_w$ the polynomial
\[
 p_\lambda(x)=\prod_{\substack{\mu\in S_w\\\mu\ne\lambda}}
                        \frac{x-\mu}{\lambda-\mu}
\]
has degree $|S_w|-1$ and satisfies $p_\lambda(H)w=v_\lambda$. Thus $\mathcal K_{|S_w|}(H,w)=\operatorname{span}\{v_\lambda:\lambda\in S_w\}$, which has dimension $|S_w|$.
\end{proof}

If $d<r\le N$, the Krylov subspace $\mathcal K_r(H,w)$ still has dimension $d$, whereas $\widehat P_r=\sum_{j\le r}q_jq_j^*$ has rank $r$. It includes vectors from the subsequent
invariant blocks. In particular, the projection onto the first $r$ Lanczos coordinates agrees with the Krylov projection exactly when $r\le d$. For any monic polynomials
$p_0,\ldots,p_{r-1}$ of successive degrees, let $F_r=(p_0(H)w,\ldots,p_{r-1}(H)w)$. Then
\begin{equation}
 \operatorname{rank}(F_r^*F_r)=\min\{r,d\},\qquad
 \widehat P_r=F_r(F_r^*F_r)^{-1}F_r^*\quad(r\le d).
 \label{W:eq:lanczos-gram-projection}
\end{equation}
Indeed, the triangular change of polynomial basis identifies $\operatorname{Ran}F_r$ with $\mathcal K_r$. When $r\le d$, the matrix on the right is self-adjoint, its square equals
itself, and it is the identity on $\operatorname{Ran}F_r$; it is therefore the orthogonal projection onto that space. This is the identity used in
Section~\ref{W:sec:cholesky-representation-of-jacobi-coefficients}.

For example, if $H=\operatorname{diag}(\lambda_1,\ldots,\lambda_N)$ and $w=e_1$, then $r_1=0$, $d=1$, and the continuation produces $q_j=e_j$ for every $j$. Thus $T=H$, while the
spectral measure at $w$ is only $\delta_{\lambda_1}$. In general, with $\mu_H^w=\sum_{\lambda\in S_w}\|v_\lambda\|^2\delta_\lambda$,
\[
 \langle f(H)w,w\rangle
 =\langle f(T)e_1,e_1\rangle
 =\int f\,\dd\mu_H^w.
\]
The spectrum of $T$ contains the eigenvalues of all its invariant blocks. We therefore use the spectral projection of $T$ in the localization argument.

Using the same continuation preserves equality of the Lanczos matrices whenever the original matrices agree. The truncation error is therefore bounded by the probability that those
matrices differ.

For a matrix $H$ let $\mathscr R_N(H)$ be the embedded resolvent of its Lanczos tridiagonal matrix, with the fixed starting vector and continuation in
Section~\ref{W:sec:lanczos-tridiagonalization}.

\begin{proposition}
\label{W:prop:operator-equality-transfer}
The map $H\mapsto\mathscr R_N(H)$ is Borel measurable into $\mathscr Y$. For any coupled Hermitian matrices $H_N,\widehat H_N$ and bounded Borel function
$f:\mathscr Y\to\mathbb C$,
\begin{equation}
 \bigl|\mathbb Ef(\mathscr R_N(H_N))-
              \mathbb Ef(\mathscr R_N(\widehat H_N))\bigr|
 \le2\|f\|_\infty\mathbb P(H_N\ne\widehat H_N).
 \label{W:eq:operator-equality-transfer}
\end{equation}
\end{proposition}
\begin{proof}
We verify the measurability of the recurrence by induction. If the first $k$ vectors are Borel functions of $H$, so are $\widehat P_k$, $a_k$, and $r_k=(I-\widehat P_k)Hq_k$.
Dividing a nonzero vector by its norm is a continuous operation. At breakdown, with $v_{k,j}=(I-\widehat P_k)e_j$, the choice of $j$ is made on the disjoint Borel sets
\[
 B_{k,j}=\{r_k=0\}\cap\{\|v_{k,j}\|>0\}
                      \cap\bigcap_{i<j}\{\|v_{k,i}\|=0\}.
\]
They partition $\{r_k=0\}$ when $k<N$. Thus $q_{k+1}$ is Borel on each member of a finite Borel partition, proving the induction. The entries of $T=U^*HU$ are consequently Borel
functions of $H$. For Hermitian $A,B$ the resolvent identity gives
\[
 \sup_{z\in K}\|(A-z)^{-1}-(B-z)^{-1}\|
 \le d_K^{-2}\|A-B\|,\qquad
 d_K=\inf_{z\in K}|\operatorname{Im}z|>0.
\]
Hence the map from $T$ to its embedded resolvent path is continuous for fixed $N$. Composing proves the path valued measurability.

On $D_N=\{H_N=\widehat H_N\}$, induction in the Lanczos recurrence gives the same vectors and coefficients for both matrices. The deterministic continuation rule makes the same
choice at each breakdown. Hence $\mathscr R_N(H_N)=\mathscr R_N(\widehat H_N)$ on $D_N$. It follows that
\[
 \begin{split}
 &\bigl|\mathbb Ef(\mathscr R_N(H_N))-
                  \mathbb Ef(\mathscr R_N(\widehat H_N))\bigr|\\
 &\quad=\left|\mathbb E\left[
       \bigl(f(\mathscr R_N(H_N))-f(\mathscr R_N(\widehat H_N))\bigr)
                                             \mathbf1_{D_N^c}\right]\right|\\
 &\quad\le2\|f\|_\infty\mathbb P(D_N^c).
 \end{split}
\]
This proves \eqref{W:eq:operator-equality-transfer}.
\end{proof}

In particular, if $\mathbb P(H_N\ne H_N^{\rm tr})\to0$, convergence in distribution of $\mathscr R_N(H_N^{\rm tr})$ implies the same convergence for $\mathscr R_N(H_N)$: apply
\eqref{W:eq:operator-equality-transfer} to every bounded continuous function on the metric space $\mathscr Y$. Centering and variance normalization are treated by the entry
comparison estimates in the subsequent sections.

\subsubsection{Reduction to coefficient convergence and localization}
\label{W:sec:operator-limit-reduction}

For the Wigner matrices of Assumption~\ref{W:ass:wigner-lee-yin}, Theorem~\ref{W:thm:operator-probabilistic} reduces the operator limit to the local convergence in
\eqref{W:eq:operator-probabilistic-local} and the localization in \eqref{W:eq:operator-probabilistic-localization}.
Theorem~\ref{W:thm:integrated-potential-and-offdiagonal-coefficients} establishes the first condition, and Theorem~\ref{W:thm:localization-spectral-projection} establishes the
second. Their proofs are in Sections~\ref{W:sec:tightness-and-identification-of-coefficients} and~\ref{W:sec:proof-full-localization}, respectively. The final application is given
in Section~\ref{W:sec:main-assembly}, after both estimates.

\subsection{Reduction of the resolvent estimate for the tridiagonal models}
\label{sec:the-soft-edge-comparison-scheme}
\label{A:sec:outline-of-the-proof}

The upper bound has two contributions: the error in replacing the matrix coefficients by Brownian increments, and the error in replacing the resulting difference operator by
$\mathcal H_\beta$. We first explain the coefficient estimates required by the resolvent identity. The independent coefficients of the Gaussian $\beta$-ensemble give a direct
construction. For a polynomial potential, the covariance factorization must also retain a discrete derivative, whose primitive has a better bound than its coordinate differences.
At the end of this reduction, we combine these estimates with the discretization of the limiting operator and the lower bound of Proposition~\ref{A:res:5.1}, whose proof uses only
the limiting spectrum and the permitted rank.

Fix the model and $\beta>0$, and use the notation of Subsection~\ref{A:sec:intro-main-results}. Let $K\Subset\mathbb C\setminus\mathbb R$ be nonempty and compact; replacing it by
$K\cup\overline K$ allows us to use the same estimates for resolvents and their adjoints. Write $h=h_n$, take $N=n$ in the Gaussian case and $N=M_n$ in the general case, and set
$\iota_h=\iota_h^{(N)}$. For this resolvent comparison, abbreviate
\[
 S_n(z)=\mathscr I_n(\mathscr A_n-z)^{-1}\mathscr I_n^*,\qquad
 S_{N,h}^B(z)=\iota_h(\mathcal B_{N,h,\beta}^B-z)^{-1}\iota_h^*.
\]
The basic comparison is
\[
 \|S_n(z)-R_\beta(z)\|_{\mathfrak S_2}^2
 \le2\|S_n(z)-S_{N,h}^B(z)\|_{\mathfrak S_2}^2
   +2\|S_{N,h}^B(z)-R_\beta(z)\|_{\mathfrak S_2}^2.
 \eqtag{outline-comparison}
\]
Both terms will be bounded by $C_{\omega,K}h$, uniformly for $z\in K$ and all sufficiently large $n$. The same Brownian path is used in both comparisons and at every matrix size.

\subsubsection{Weighted perturbations}\label{A:sec:earlier-results}
\label{A:sec:rank-lower-bound-and-perturbation-estimates}
The resolvent identity distinguishes an ordinary diagonal error from an error that is a discrete derivative. Estimating both coordinate by coordinate would lose the gain from the
primitive of the latter.

Let $R_a(z)=(H_a-z)^{-1}$, $a=0,1$, where the finite operators satisfy \eqref{A:eq:3.8}, \eqref{A:eq:3.10}, and \eqref{A:eq:3.12} with common constants. The two estimates below
specify the accuracy required of the coefficient coupling. For a diagonal difference $H_1-H_0=h^{-2}M_r$, the resolvent identity and Proposition~\ref{A:res:3.16} give
\[
 \|R_1(z)-R_0(z)\|_{\mathfrak S_2}^2
 \le C_{\omega,K}h^{-3}\sum_{j=1}^N(1+x_j)^{-7/2}|r_j|^2.
\]
Since the prefactor is $h^{-3}$, a squared coefficient sum of order $O(h^4)$ suffices for squared resolvent error $O(h)$. The coefficients in this sum include the matrix scaling.

A discrete derivative can be estimated through its primitive. Use counting measure, as in Subsection~\ref{A:sec:weighted-resolvent-estimates-and-comparison-principles}, put
$u_0=u_{N+1}=0$, and write $(D_hu)_j=(u_{j+1}-u_j)/h$ and $(J_hu)_j=(u_{j+1}+u_j)/2$.

For a real edge sequence $e=(e_j)_{j=0}^N$, summation by parts gives
\[
 \sum_{j=1}^N\frac{e_j-e_{j-1}}h|u_j|^2
 =-\sum_{j=0}^Ne_j\frac{|u_{j+1}|^2-|u_j|^2}h
 =-2\operatorname{Re}\sum_{j=0}^Ne_j(D_hu)_j\overline{(J_hu)_j}.
\]
Consequently, if $r_j=h(e_j-e_{j-1})$, then $h^{-2}M_r=\Gamma_e$, with $\Gamma_e=-(D_h^*M_eJ_h+J_h^*M_eD_h)$. Proposition~\ref{A:res:3.17} therefore gives the estimate
\[
 \|R_1(z)-R_0(z)\|_{\mathfrak S_2}^2
 \le C_{\omega,K}\sum_{j=0}^Nh(1+x_j)^{-3/2}|e_j|^2.
\]
The same estimate holds for the first-order forms using $J_-$ or $J_+$. Unlike $\Gamma_e$, these need not be diagonal. Both Dirichlet edges are included in the estimates.

The difference between these estimates matters at the scales used below. Suppose, for example, that only $\mathbb E|e_j|^2\le Ch$ is available. Since
$\sum_j(1+x_j)^{-a}\le C_a h^{-1}$ for $a>1$, coordinatewise estimation of $r_j=h(e_j-e_{j-1})$ gives only
\[
 \mathbb E\!\left[h^{-3}\sum_{j=1}^N(1+x_j)^{-7/2}|r_j|^2\right]
 \le Ch^{-1},\qquad
 \mathbb E\!\left[\sum_{j=0}^Nh(1+x_j)^{-3/2}|e_j|^2\right]\le Ch.
\]
Here the first bound follows from $|e_j-e_{j-1}|^2\le2|e_j|^2+2|e_{j-1}|^2$ and a sum of order $h^{-1}$. Summation by parts instead retains $e_j$ with the factor $h$, giving the
second bound. This example explains the required decomposition; it is only a comparison of expectations. The coefficient constructions below supply the stronger tail estimates
needed for an almost sure rate.

The common form hypotheses follow from Proposition~\ref{A:res:8.2B} in the Gaussian case and Theorem~\ref{A:res:16.2} in the general case. The coordinate maxima and relative form
bounds make the resolvent estimates uniform over the intermediate operators. The weighted square sums then determine the rate.

\subsubsection{The Gaussian model}\label{A:sec:outline-gaussian-model}
For the Gaussian model, independence of the Jacobi coefficients makes the preceding decomposition explicit. 

In the Dumitriu--Edelman model, write $\delta_{j,n}=b_{j,n}-\mathbb E b_{j,n}$ for $1\le j<n$ and set $\delta_{0,n}=\delta_{n,n}=0$. Proposition~\ref{A:res:6.1} gives
\[
 \mathcal A_{n,\beta}^{\mathrm{DE}}
 =\overline H_{n,\beta}+h^{-2}M_{\zeta_n}+\Gamma_{\eta_n}^-,\qquad
 \zeta_{j,n}=-a_{j,n}-2\delta_{j,n},\qquad
 \eta_{j,n}=h^{-1}\delta_{j,n},
\]
where $\overline H_{n,\beta}=\mathbb E\mathcal A_{n,\beta}^{\mathrm{DE}}$ and $\Gamma_{\eta_n}^-$ uses $J_-$ in place of $J_h$. The variables $\zeta_{j,n}$ use disjoint coefficient
pairs as $j$ varies and are therefore independent. The two fields $\zeta_n$ and $\eta_n$ themselves share the variables $\delta_{j,n}$. Their operator errors are combined by a
triangle inequality, which requires no independence between these fields. We couple the first field to Brownian increments and apply the primitive estimate to the second.

Put
\[
 \zeta_{j,n}^B=\frac2{\sqrt\beta}h\bigl(B(jh)-B((j-1)h)\bigr).
\]
The coupling in Subsection~\ref{A:sec:coupling-the-dumitriu-edelman-fluctuations-to-brownian-increments} uses an orthogonal rotation of each normalized Brownian increment and an
independent normal variable, followed by the chi quantile map. The rotation preserves the joint standard Gaussian distribution. The quantile map then gives the exact chi laws, and
the diagonal Gaussian terms cancel in $\zeta_{j,n}-\zeta_{j,n}^B$. Equation~\eqref{A:eq:8.28} yields
\[
 \|\zeta_{j,n}-\zeta_{j,n}^B\|_{L^{2p}}
 \le C_{p,\beta}\frac{h^{3/2}}{\sqrt{1+n-j}},\qquad 1\le j\le n.
\]
Minkowski's inequality now gives
\[
 \left\|h^{-3}\sum_{j=1}^n(1+x_j)^{-7/2}
       |\zeta_{j,n}-\zeta_{j,n}^B|^2\right\|_{L^p}
 \le C_{p,\beta}\sum_{j=1}^n\frac{(1+x_j)^{-7/2}}{1+n-j}
 \le C_{p,\beta}h^2.
\]
For the last inequality, split at $\lfloor n/2\rfloor$. On the first part the denominator is of order $n$; on the second, $1+x_j\ge ch^{-2}$. Hence the two contributions are
bounded by $Cn^{-1}h^{-1}=Ch^2$ and $Ch^7(1+\log n)$, respectively. The latter bound includes $j=n$, where the centered chi variable in the coupling is zero.

Let $V_n\ge0$ be the scalar weighted error in this comparison, and retain its finite-size estimate $\|V_n\|_{L^p}\le C_p h_n^2$ for every fixed $p\ge1$, with deterministic
$h_n\asymp n^{-1/3}$. Fix $0<\varepsilon<1$ and choose $p$ with $p\varepsilon>3$. Markov's inequality gives
\[
 \mathbb P\{V_n>h_n^{2-\varepsilon}\}
 \le h_n^{-p(2-\varepsilon)}\mathbb E V_n^p
 \le C_p h_n^{p\varepsilon}
 \le C_{p,\varepsilon}n^{-p\varepsilon/3}.
\]
The last series is summable. Hence there is an almost surely finite index $N_{\varepsilon,p}$ such that
\[
 V_n\le h_n^{2-\varepsilon}\quad(n\ge N_{\varepsilon,p}),
 \qquad h_n^{-1}V_n\le h_n^{1-\varepsilon}\longrightarrow0.
\]
No independence between matrix sizes is needed. Any finite pathwise multiplicative constant from the resolvent estimate can be included after fixing the outcome; it does not affect
this limit.

In contrast, Proposition~\ref{A:res:A.4} gives an exponential tail at scale $h$ for
\[
 \sum_{j=0}^nh(1+x_j)^{-3/2}|\eta_{j,n}|^2\le C_\omega h
\]
almost surely for all sufficiently large $n$. The final $O(h)$ bound therefore does not come from the scalar Gaussian approximation alone. After these two coefficient comparisons,
Proposition~\ref{A:res:7.5} compares the mean matrix with Brownian diagonal added to $\mathcal B_{n,h,\beta}^B$, again with squared error $O_\omega(h)$. The scalar error is
$o_\omega(h)$ in squared resolvent norm; the first-order form and the deterministic comparison each contribute $O_\omega(h)$. The inequality
$\|D_1+D_2+D_3\|_{\mathfrak S_2}^2\le3\sum_{j=1}^3\|D_j\|_{\mathfrak S_2}^2$ then proves Theorem~\ref{A:res:8.5}.

\subsubsection{The general Jacobi model}\label{A:sec:the-two-matrix-models}
The coefficients are now dependent. After the boundary restriction, we approximate a conditional coefficient block by a Gaussian vector, whose covariance is the inverse Hessian at
the conditional minimizer. Write $I_n^{\mathrm{mid}}$ for the interval where this approximation is used. It lies in the regular region and is separated from both boundaries of the
conditional block by logarithmically many sites; the precise intervals are defined in Subsection~\ref{A:sec:gaussian-approximation-of-the-conditional-gibbs-measure}. These
separations allow a translation-invariant covariance calculation to be used without discarding the finite boundary errors.

The reconstruction is performed under the original exterior distribution. Let $\tau$ denote the retained exterior variables, let $\nu_n$ be their law, and let $\mu_{n,\tau}$ be the
conditional law of the interior coefficients. Suppose $\mathcal Q_n(\tau;\dd x,\dd g,\dd b)$ is the conditional coupling used in the comparison, whose $x$-marginal is exactly
$\mu_{n,\tau}$. The unconditional joint law is
\[
 \nu_n(\dd\tau)\,\mathcal Q_n(\tau;\dd x,\dd g,\dd b).
\]
For every bounded Borel function $F$ of the original coefficient array,
\[
 \begin{split}
 \int F(\tau,x)\,\nu_n(\dd\tau)\mathcal Q_n(\tau;\dd x,\dd g,\dd b)
 &=\int\left(\int F(\tau,x)\,\mu_{n,\tau}(\dd x)\right)\nu_n(\dd\tau)\\
 &=\mathbb E_{\mathrm{original}}F(\tau,x).
 \end{split}
\]
Thus coupling the interior coefficients preserves the original joint coefficient distribution. A restriction to an event on which estimates hold is not part of this sampling law.
If $E$ implies that $\tau$ belongs to an admissible set $T_n$, then for every nonnegative Borel
$H$,
\[
 \mathbb E[\mathbf1_EH]
 \le\int_{T_n}\mathbb E[H\mid\tau]\,\nu_n(\dd\tau).
\]
This inequality permits conditional moment estimates without replacing $\nu_n$ by its conditional law on $T_n$. The coupling kernel must still be defined on the complement of $T_n$
with the original conditional coefficient marginal. Estimates on that complement are handled separately. Equation~\eqref{A:eq:exact-kernel-marginal} supplies the required
conditional marginal in this construction.

Let $\zeta_j$ be the scalar coefficient in the decomposition of Proposition~\ref{A:res:13.1}, and let $\widehat\zeta_j^{\mathrm G}$ be its conditional Gaussian approximation. The
Gaussian approximation is centered at the minimizer of its own conditional block. Its difference from the minimizer used to define $\zeta_j$ is retained in
\eqref{A:eq:14.13}--\eqref{A:eq:14.15}; it is part of the error below. The proof of Theorem~\ref{A:res:14.7} gives the unconditional estimate
\[
 \left\|\sum_{j\in I_n^{\mathrm{mid}}}(1+x_j)^{-7/2}
       |\zeta_j-\widehat\zeta_j^{\mathrm G}|^2\right\|_{L^p}
 \le C_ph^5.
\]
For a nonnegative error $W_n$ with $\|W_n\|_{L^p}\le C_p h_n^5$: choosing $p>3$ gives
\[
 \mathbb P\{W_n>h_n^4\}\le C_p h_n^p,
 \qquad W_n\le h_n^4
\]
for every $n$ above an almost surely finite index. This last implication gives $h_n^{-3}W_n=O(h_n)$, not $o(h_n)$ by itself. A smaller order must use the stronger moment estimate
directly, for example a threshold $h_n^{5-\varepsilon}$ with $0<\varepsilon<1$ and $p\varepsilon>3$, which yields $h_n^{-3}W_n=o(h_n)$.

Fix the parameter in the stationary Gaussian approximation and suppress it from the notation. Let $(Z_j)_{j\in\mathbb Z}$ be centered and stationary, with $c_r=\mathbb E[Z_0Z_r]$
and $\sum_{r\in\mathbb Z}|c_r|<\infty$. For $m\ge1$, exactly $m-|r|$ ordered pairs $(i,j)\in\{1,\ldots,m\}^2$ have difference $r$ when $|r|<m$. Therefore
\[
 \begin{split}
 \frac1m\operatorname{Var}\left(\sum_{j=1}^m Z_j\right)
 &=\frac1m\sum_{i,j=1}^m c_{j-i}\\
 &=\sum_{r\in\mathbb Z}\left(1-\frac{|r|}{m}\right)_+c_r.
 \end{split}
\]
For each fixed $r$ the coefficient tends to one and has absolute value at most one. Given $\varepsilon>0$, choose $R$ with $\sum_{|r|>R}|c_r|<\varepsilon$. The difference from
$\sum_rc_r$ is then bounded by
\[
 \frac1m\sum_{|r|\le R}|r|\,|c_r|+2\varepsilon.
\]
Taking $m\to\infty$ and then $\varepsilon\downarrow0$ proves convergence to $\widehat c(0)=\sum_rc_r$. In particular this sum is nonnegative because it is a limit of variances.

For the stationary field at parameter $s$, this limiting variance is $\widehat C_s(0)$. The Gaussian coordinates remain correlated, and their partial sums determine the independent
increment variance.

Set $b(s)=b^\dagger(s)$, $\tau(s)=-E'(s)$ and $\mathsf H_s=D^2F_{V,s}(a^\dagger(s),b^\dagger(s))$, with $F_{V,s}$ from \eqref{A:eq:local-energy}. Lemma~\ref{A:res:15.1}, with
$\ell=(1,2)$, gives
\[
 \widehat C_s(0)
 =\frac{2}{\beta n b_\star^2}\ell\mathsf H_s^{-1}\ell^*
 =\frac{4b(s)\tau(s)}{\beta n b_\star^2},\qquad
 \widehat C_0(0)=\frac4\beta h^3.
\]
Indeed, $\operatorname{Var}(\zeta_{j,n}^B)=(4/\beta)h^2\operatorname{Var}(B(jh)-B((j-1)h)) =(4/\beta)h^3$. Thus the variance per term of the partial sums agrees with the variance
$(4/\beta)h^3$ of the chosen increments.

Spectral factorization gives a coupling with the Brownian increments having this variance. For the moment keep $s$ fixed. The normalized covariance symbol has an analytic factor
with real coefficients:
\[
 q_s(\mathrm e^{\mathrm i\theta})
   =h^{-3}\widehat C_s(\theta)
   =|g_s(\mathrm e^{\mathrm i\theta})|^2,\qquad g_s(1)>0.
\]
For this calculation write $g_l=g_l(s)$ and $\sigma=g_s(1)$. Suppose the real coefficients of the fixed spectral factor satisfy $|g_l|\le C\rho^l$ for some $0<\rho<1$. Set
\[
 \sigma=\sum_{l\ge0}g_l,\qquad
 p_l=-\sum_{r=l+1}^\infty g_r.
\]
Then $|p_l|\le C\rho^{l+1}/(1-\rho)$, $p_0=g_0-\sigma$, and $p_l-p_{l-1}=g_l$ for $l\ge1$. Consequently, throughout a disk of radius greater than one,
\[
 g(z)=\sum_{l\ge0}g_lz^l
      =\sigma+(1-z)\sum_{l\ge0}p_lz^l.
\]
Let $\xi_j$, $j\in\mathbb Z$, be independent standard real Gaussian variables. The series
\[
 \widehat Z_j=\sum_{l\ge0}g_l\xi_{j-l},\qquad
 V_j=\sum_{l\ge0}p_l\xi_{j-l}
\]
converge in $L^2$ because the coefficient sequences are square summable. Subtracting their partial sums and passing in $L^2$ to the limit gives
\[
 \widehat Z_j=\sigma\xi_j+V_j-V_{j-1}.
\]
For a finite set of indices these equalities hold simultaneously almost surely; intersecting over $j\in\mathbb Z$ gives a single event for the whole sequence. Also
\[
 \mathbb E V_j^2=\sum_{l\ge0}p_l^2\le C_0,
 \qquad
 \sum_{j=r+1}^s \widehat Z_j
 =\sigma\sum_{j=r+1}^s\xi_j+V_s-V_r.
\]
At spacing $h$, set $e_j=h^{1/2}V_j$. Then $\mathbb E e_j^2\le C_0h$ and
\[
 Z_j=h^{3/2}\widehat Z_j=h^{3/2}\sigma\xi_j+h(e_j-e_{j-1}).
\]
This is the precise difference representation used to estimate the latter term through its primitive. When the factor depends on a slowly varying parameter, this identity is first
used with that parameter fixed; the error from changing parameters is the separate covariance remainder already specified in the comparison. At $s=0$, $\sigma=2/\sqrt\beta$. Taking
$\xi_j=h^{-1/2}(B(jh)-B((j-1)h))$ for $j\ge1$ identifies the independent term with $\zeta_{j,n}^B$; the nonpositive indices use independent standard normal variables. The analytic
factor and its uniform coefficient bounds are established in Subsection~\ref{A:sec:spectral-factorization-and-coupling-to-brownian-increments}.

In the actual model, the parameter $s=s_j$ varies across a finite conditional block. Theorem~\ref{A:res:15.3} gives, on $I_n^{\mathrm{mid}}$,
\[
 \zeta_j-\zeta_j^B
 =h(e_j-e_{j-1})+(\zeta_j-\widehat\zeta_j^{\mathrm G})
       +\rho_j+\varepsilon_j^{\mathrm{fr}}.
\]
Here $\zeta_j^B=\zeta_{j,n}^B$, $\rho_j$ includes the variation of $s_j$ and the difference $g_{s_j}(1)-2/\sqrt\beta$, and $\varepsilon_j^{\mathrm{fr}}$ is the error in coupling
the conditional Gaussian vector to the vector constructed from the varying spectral factors. The original first-order coefficient $\eta_j$ of Proposition~\ref{A:res:13.1} is
unchanged.

Extend $e$ constantly beyond each end of $I_n^{\mathrm{mid}}$, as in \eqref{A:eq:15.51}, and call the extension $\widetilde e$. Its difference then vanishes outside that interval,
without introducing boundary jumps. The weighted estimate includes both constant pieces. Theorems~\ref{A:res:14.7} and~\ref{A:res:15.3}, together with the coefficient estimates in
Subsection~\ref{A:sec:quadratic-form-bounds}, give almost surely
\[
 h^{-3}\sum_{j\in I_n^{\mathrm{mid}}}(1+x_j)^{-7/2}
 \bigl(|\zeta_j-\widehat\zeta_j^{\mathrm G}|^2
       +|\rho_j|^2+|\varepsilon_j^{\mathrm{fr}}|^2\bigr)
 \le C_\omega h,
\]
\[
 \sum_{j=0}^Nh(1+x_j)^{-3/2}
       \bigl(|\widetilde e_j|^2+|\eta_j|^2\bigr)\le C_\omega h.
\]
These are the two weighted requirements established above. Once the common form bounds of Theorem~\ref{A:res:16.2} are in place, they give squared resolvent error $O_\omega(h)$ on
the middle interval.

The two ends are estimated directly. For a banded perturbation $W$ of fixed bandwidth on the first $\ell=O(\log n)$ coordinates, with entries bounded by
$C_\omega h^{-1/2}(\log n)^q$, the Dirichlet boundary estimate \eqref{A:eq:3.58d} gives
\[
 \sup_{z\in K}\|R_1(z)WR_0(z)\|_{\mathfrak S_2}^2
 \le C_{\omega,K}h^3(\log n)^{2q+4}=o_\omega(h).
\]
At the other end $x_j\asymp h^{-2}$, and Corollary~\ref{A:res:3.18} gives $\|R_a(z)\boldsymbol e_j\|^2\le C_{\omega,K}h^4$ for the $j$th coordinate vector $\boldsymbol e_j$. The
replacements of the scalar coefficients on these two regions, and the comparisons using Schur complements to include the removed coordinates, are included in \eqref{A:eq:15.65}.
The remaining comparison changes only the deterministic coefficients. It uses the product rule \eqref{A:eq:product-rule} on all edges and retains an $O_\omega(h)$ bound for the
contribution of the resolvent on the final coordinates. Theorem~\ref{A:res:17.5} completes this comparison, and Theorem~\ref{A:res:17.6} gives the bound for the full finite matrix.

It remains to apply the comparison maps. For an isometry $J$ and a finite matrix $D$,
\[
 \|JDJ^*\|_{\mathfrak S_2}^2
 =\operatorname{Tr}(JD^*DJ^*)
 =\operatorname{Tr}(D^*D)=\|D\|_{\mathfrak S_2}^2.
 \eqtag{isometric-hs}
\]
In the Gaussian case this identity applies directly to Theorem~\ref{A:res:8.5}. In the general case apply the contraction $\mathcal J_n$ to the difference for the full matrix in
Theorem~\ref{A:res:17.6}. Since $\mathcal J_nP_{I_n}^*U_{r_n}^*=\iota_h^{(N)}$, the ideal property gives the same conclusion:

\[
 \sup_{z\in K}\|S_n(z)-S_{N,h}^B(z)\|_{\mathfrak S_2}^2
 \le C_{\omega,K}h.
\]
The expression $S_n$ therefore remains the compressed resolvent of the full matrix.

\subsubsection{The continuum approximation}\label{A:sec:passage-to-the-continuum}
The remaining comparison depends only on the common Brownian path. The interval length is $L=(N+1)h\asymp h^{-2}$, so estimates on a fixed interval would not suffice. We extend
finite interval resolvents by zero to $L^2(\mathbb R_+)$. The Dirichlet restriction at $L$ has squared error $O_{B,\beta,K}(L^{-1/2})$ by \eqref{A:eq:4.23}.

On $(0,L)$, the proof in Subsection~\ref{A:sec:finite-difference-approximation-of-the-stochastic-airy-operator} first samples the integrated Brownian potential. The deterministic
potential is changed to $x-h/2$, whose value at the interpolation vertex $jh$ is the center $x_j$. Equation~\eqref{A:eq:4.31} bounds this change of forms by $O_{B,\beta,K}(h)$ in
squared resolvent norm. The continuous piecewise linear Galerkin approximation has the same bound by \eqref{A:eq:4.41}.

There are two further comparisons before reaching the step embedding in the theorem. The piecewise linear basis has mass matrix
\[
 M_h=I_h^*I_h=I-\frac{h^2}{6}D_h^*D_h,
\]
so the resolvent of the restricted form cannot be identified with the inverse of the difference matrix by replacing $M_h$ by $I$. Equation~\eqref{A:eq:4.48} controls the mass and
quadrature terms. Finally, \eqref{A:eq:4.55} compares the piecewise linear map $I_h$ with $\iota_h$. Both squared errors are $O_{B,\beta,K}(h)$ with constants valid up to
$L\asymp h^{-2}$. In particular, the mass matrix comparison is an operator estimate in its own right; equivalence of the two finite inner products alone would not give this rate.

All five differences can be written on $L^2(\mathbb R_+)$. Let $\iota_L$ extend functions on $(0,L)$ by zero. Denote the half-line resolvent by $T_0(z)=R_\beta(z)$, the Dirichlet
resolvent on $(0,L)$ by $R_L(z)$, and the resolvent with the sampled integrated potential on that interval by $R_L^{\mathrm{fr}}(z)$. With the piecewise linear map $I_h$, the step
isometry $\iota_h$, the stiffness matrix $K_h^{\mathrm c}$, and the mass matrix $M_h$ from the subsequent discretization argument, put
\[
 \begin{aligned}
 T_1(z)&=\iota_LR_L(z)\iota_L^*,&
 T_2(z)&=\iota_LR_L^{\mathrm{fr}}(z)\iota_L^*,\\
 T_3(z)&=I_h(K_h^{\mathrm c}-zM_h)^{-1}I_h^*,&
 T_4(z)&=I_h(\mathcal B_{N,h,\beta}^B-z)^{-1}I_h^*,\\
 T_5(z)&=\iota_h(\mathcal B_{N,h,\beta}^B-z)^{-1}\iota_h^*.
 \end{aligned}
\]
The maps in the last two lines include extension by zero. The displayed formula for $T_3$ uses the mass matrix because $I_h$ is not an isometry. For example, testing the weak
equation with $I_hv$ gives $(K_h^{\mathrm c}-zM_h)u=I_h^*f$, and hence exactly the displayed operator.

There is the exact identity
\[
 T_5(z)-T_0(z)=\sum_{j=1}^5\bigl(T_j(z)-T_{j-1}(z)\bigr).
\]
If $D_j=\sup_{z\in K}\|T_j(z)-T_{j-1}(z)\|_{\mathfrak S_2}^2$, then
\[
 \sup_{z\in K}\|T_5(z)-T_0(z)\|_{\mathfrak S_2}^2
 \le5\sum_{j=1}^5D_j.
\]
The five terms are, in order, the interval restriction, sampling of the potential, Galerkin approximation, replacement of the finite matrices, and change of embedding. The
discretization subsections estimate these five differences. In the last step, writing $R_h=(\mathcal B_{N,h,\beta}^B-z)^{-1}$, one uses
\[
 T_5-T_4=(\iota_h-I_h)R_h\iota_h^*
                    +I_hR_h(\iota_h^*-I_h^*),
\]
so both occurrences of the embedding are accounted for. The estimates \eqref{A:eq:4.23}, \eqref{A:eq:4.31}, \eqref{A:eq:4.41}, \eqref{A:eq:4.48}, and \eqref{A:eq:4.55} give
$D_1\le C_{B,\beta,K}L^{-1/2}$ and $D_j\le C_{B,\beta,K}h$ for $2\le j\le5$. Since $L\asymp h^{-2}$, their sum is at most $C_{B,\beta,K}h$. This proves Theorem~\ref{A:res:1.3};
together with \eqref{A:eq:outline-comparison} and the finite matrix comparison, it gives the upper bound.

For the lower bound, $\operatorname{rank}S_n(z)\le N$. The singular value identity \eqref{A:eq:5.4} and the spectral asymptotics \eqref{A:eq:5.3} give
\[
 \inf_{z\in K}\|S_n(z)-R_\beta(z)\|_{\mathfrak S_2}^2
 \ge\inf_{z\in K}\sum_{k>N}s_k(R_\beta(z))^2
 \ge c_{\omega,K}N^{-1/3}\ge c'_{\omega,K}h.
\]
Section~\ref{A:sec:inverse-coupling-and-spectral-consequences} constructs the common event and proves the real parameter comparison, including the contribution of the omitted
coordinates in the general model.

\subsection{Weighted estimates for discrete resolvents}
\label{sec:weighted-airy-resolvents-and-comparison-principles}
\label{A:sec:weighted-resolvent-estimates-and-finite-difference-approximation}

\label{A:sec:A}
We estimate the change in a resolvent caused by a diagonal coefficient and by the discrete derivative of a coefficient. The relevant norms have weights $(1+x)^{-7/2}$ and
$(1+x)^{-3/2}$, respectively. Their exponents will follow from the local length $(1+x)^{-1/2}$ of the Green kernel. The lower bound in
Subsection~\ref{A:sec:rank-lower-bound-and-optimality} does not use these perturbation estimates. Throughout this subsection, the matrix size varies while the form constants,
coefficient bounds and compact spectral set are fixed. The constants may depend on the sampled outcome; the estimates are pathwise and use only their finiteness.

\subsubsection{Weighted resolvent estimates}\label{A:sec:weighted-resolvent-estimates-and-comparison-principles}

All estimates in this subsection are deterministic. The form bound controls $D_hR(z)$ in operator norm. To account for the location of a coefficient error, we prove decay of the
resolvent kernel on the local length scale $(1+x)^{-1/2}$. Applying this decay to the two kernels of $R(z)^*R(z)$ in the diagonal perturbation estimate gives the weight
$(1+x)^{-7/2}$.

Let $H_0,H_1$ be self-adjoint operators associated with closed, lower semibounded Hermitian sesquilinear forms $\mathfrak h_0,\mathfrak h_1$ on a common form domain. Set
$\mathfrak w=\mathfrak h_1-\mathfrak h_0$ and $R_j(z)=(H_j-z)^{-1}$. The form resolvent identity is
\[
 \langle f,(R_1(z)-R_0(z))g\rangle
 =-\mathfrak w[R_1(\bar z)f,R_0(z)g].
\]
Indeed, for $u=R_1(\bar z)f$ and $v=R_0(z)g$, the weak equations give
\[
\mathfrak w[u,v] =\bigl(\mathfrak h_1[u,v]-\bar z\langle u,v\rangle\bigr) -\bigl(\mathfrak
 h_0[u,v]-\bar z\langle u,v\rangle\bigr) =\langle f,R_0(z)g\rangle-\langle R_1(\bar z)f,g\rangle.
\]
For a bounded perturbation $W=H_1-H_0$ this is $R_1-R_0=-R_1WR_0$. We use the form identity when $W$ contains derivatives.

\paragraph{Quadratic forms and spectral bounds}
Let $h=h_n\downarrow0$ and $N_nh\asymp h^{-2}$. We retain $x_j=(j-\tfrac12)h$ and use counting measure on both $\ell_n^2=\mathbb C^{N_n}$ and
$\ell_{n,\mathrm e}^2=\mathbb C^{\{0,\ldots,N_n\}}$: $\langle u,v\rangle_{\ell^2}=\sum_{j=1}^{N_n}u_j\overline{v_j}$. Vertex vectors are extended by the Dirichlet values
$u_0=u_{N_n+1}=0$.

The edge indices are $0,\ldots,N_n$; both boundary differences are included. Define \(D_h:\ell_n^2\to\ell_{n,\mathrm e}^2\) by \((D_hu)_j=(u_{j+1}-u_j)/h\), \(0\le j\le N_n\). Its
adjoint is \((D_h^*v)_j=(v_{j-1}-v_j)/h\), \(1\le j\le N_n\). The Dirichlet Laplacian is \(-\Delta_{h,D}=D_h^*D_h.\) In coordinates,
\((-\Delta_{h,D}u)_j = h^{-2}(2u_j-u_{j-1}-u_{j+1}).\) On $\mathbb C^{N_n}$ define the positive quadratic form
\[
\ell_h[u]
=
\sum_{j=0}^{N_n}|(D_hu)_j|^2
+
\sum_{j=1}^{N_n}(1+x_j)|u_j|^2.
\]
Let $X_h$ denote multiplication by $x_j$. In sums over edges we set $x_0=0$ and, when a boundary expression requires it, $x_{N_n+1}=(N_n+\tfrac12)h$. The edge weight $1+x_j$ is
comparable to $1+jh$ for $0\le j\le N_n$. Define
\[
 (J_-u)_j=u_j,\qquad (J_+u)_j=u_{j+1},\qquad
 J_h=\tfrac12(J_-+J_+).
\]
The definitions include both boundary edges. We will write $\Phi_j=\Phi(x_j)$ once the exponential weight $\Phi$ is introduced.

A restriction to indices $p,\ldots,q$ means taking the corresponding principal matrix and setting the values at $p-1,q+1$ equal to zero in its difference representation. Deleting
$h^{-2}|u_{j+1}-u_j|^2$ instead lowers the Laplacian form and separates the intervals with free boundary values. We use this Neumann bracketing for eigenvalue counts and principal
restrictions for resolvent comparisons.

The interpolation at the nodes $jh$ is normalized to agree with counting measure: let $I_hu$ be continuous and linear on each interval, with $(I_hu)(jh)=h^{-1/2}u_j$ for
$0\le j\le N_n+1$. Its energy, unlike its $L^2$ norm, agrees exactly with the discrete derivative energy:
\[
\int_0^{(N_n+1)h}|(I_hu)'(x)|^2\,\mathrm{d}x
=
\sum_{j=0}^{N_n}|(D_hu)_j|^2.
\]
For $x=jh+th$, $0\le t\le1$, the linear interpolation is
\[
 (I_hu)(x)=h^{-1/2}\bigl((1-t)u_j+tu_{j+1}\bigr),\qquad u_0=u_{N_n+1}=0.
\]
The exact unweighted cell integral is
\[
 \int_{jh}^{(j+1)h}|I_hu|^2\,\dd x
 =\frac{|u_j|^2+\operatorname{Re}(u_j\overline{u_{j+1}})+|u_{j+1}|^2}{3}.
\]
Since $|\operatorname{Re}(u_j\overline{u_{j+1}})|\le(|u_j|^2+|u_{j+1}|^2)/2$, this integral lies between one sixth and one half of $|u_j|^2+|u_{j+1}|^2$. The actual minimum and
maximum of $1+x$ on the cell are $1+jh$ and $1+(j+1)h$. Therefore
\[
 \begin{split}
 \frac{1+jh}{6}(|u_j|^2+|u_{j+1}|^2)
 &\le\int_{jh}^{(j+1)h}(1+x)|I_hu(x)|^2\,\dd x\\
 &\le\frac{1+(j+1)h}{2}(|u_j|^2+|u_{j+1}|^2).
 \end{split}
\]
Summing over $j=0,\ldots,N_n$, and using $x_i=(i-1/2)h$, gives
\[
 \frac13\sum_{i=1}^{N_n}(1+x_i)|u_i|^2
 \le\int_0^{(N_n+1)h}(1+x)|I_hu|^2\,\dd x
 \le\sum_{i=1}^{N_n}(1+x_i+h)|u_i|^2.
\]
For $0<h\le1$ the last expression is at most $2\sum_i(1+x_i)|u_i|^2$. On every cell the derivative is constant, so
\[
 \|(I_hu)'\|_2^2=h^{-2}\sum_{j=0}^{N_n}|u_{j+1}-u_j|^2.
\]
Combining these two estimates proves the claimed equivalence with $\ell_h[u]$, with deterministic constants independent of
$N_n,h$.

Consequently,
\[
 \ell_h[u]\asymp
 \int_0^{(N_n+1)h}|(I_hu)'|^2\,\mathrm dx
 +\int_0^{(N_n+1)h}(1+x)|I_hu|^2\,\mathrm dx.
 \eqtag{3.7}
\]
The unweighted bound $\frac13\|u\|^2\le\|I_hu\|_2^2\le\|u\|^2$ is the mass estimate in Lemma~\ref{W:lem:operator-mesh}.

Let $H_n$ be a real symmetric Jacobi matrix with sesquilinear form
\[
\mathfrak h_n[u,v] ={} \sum_{j=0}^{N_n}c_{j,n}(D_hu)_j\overline{(D_hv)_j} + \sum_{j=1}^{N_n}q_{j,n}u_j\overline{v_j} + \sum_{\nu=1}^{L_0} \gamma_{\eta_n^{(\nu)},\sigma_\nu}[u,v], \quad \sigma_\nu\in\{-,+\},
\eqtag{3.8}
\]
We denote the first term in \eqref{A:eq:3.8} by $\mathfrak t_n[u,v]$. Here $L_0$ is fixed and all coefficient sequences in this representation are real. Set
\[
\gamma_{\eta,\sigma}[u,v]
=
-\langle D_hu,M_\eta J_\sigma v\rangle
-\langle M_\eta J_\sigma u,D_hv\rangle.
\]
Both \(J_-\) and \(J_+\) map \(\ell_n^2\) into \(\ell_{n,\mathrm e}^2\). For an edge vector \(v\), their adjoints satisfy
\[
(J_-^*v)_j=v_j,
\quad
(J_+^*v)_j=v_{j-1},
\quad 1\le j\le N_n.
\]
In particular, $J_\pm^*J_\pm=I$ and $\|J_\pm\|_{\mathrm{op}}=\|J_\pm^*\|_{\mathrm{op}}=1$. We assume
\[
 0\le c_{j,n}\le C_\omega,\qquad
 \sum_{\nu=1}^{L_0}\|\eta_n^{(\nu)}\|_{\ell^\infty}\le C_\omega.
 \eqtag{3.10}
\]
We do not impose a pointwise bound on $q_{j,n}$. For the finite difference matrix driven by Brownian increments, the random potential is placed in the diagonal part of
\eqref{A:eq:3.8}. Its size is controlled as a form; the diagonal part cancels exactly in the exponential conjugation below. The bounds in \eqref{A:eq:3.10} instead control the
off-diagonal entries.

To obtain the form bound from an integrated potential, its local oscillations must be separated from its possibly large absolute value. The following estimate does this by
subtracting the linear interpolation of its averages on the unit intervals. For a real, measurable, locally bounded function $W$, write
$\overline W_p=\int_p^{p+1}W(x)\,\mathrm dx$. Fix $\alpha<1/2$ and assume
\[
 \operatorname{osc}_{[p,p+1]}W+
 |\overline W_{p+1}-\overline W_p|\le A(1+p)^\alpha,
 \qquad p\in\mathbb N_0.
\]
For measurable functions, oscillation may be interpreted as essential oscillation. The following estimate is the deterministic form estimate used to define the quadratic form with
white noise potential in \cite{ramirez2011beta}.

\begin{lemma}\label{A:res:3.1}
For every $\varepsilon>0$ and compactly supported $f\in H_0^1(\mathbb R_+)$,
\[
 \left|-\int_0^\infty W(|f|^2)'\,\mathrm dx\right|
 \le\varepsilon\left(\|f'\|_2^2+\int_0^\infty x|f|^2\,\mathrm dx\right)
       +C_{\varepsilon,A,\alpha}\|f\|_2^2.
\]
\end{lemma}

\begin{proof}
Define $\widetilde W(p+t)=(1-t)\overline W_p+t\overline W_{p+1}$ for $p\in\mathbb N_0$ and $0\le t\le1$. Put $a=\max\{\alpha,0\}<1/2$. For $x=p+t$, $0\le t<1$, the definition of
$\widetilde W$ gives
\[
 |W(x)-\widetilde W(x)|
 \le \operatorname{osc}_{[p,p+1]}W
       +t|\overline W_{p+1}-\overline W_p|,
 \qquad \widetilde W'(x)=\overline W_{p+1}-\overline W_p.
\]
The assumed bounds therefore give, for almost every $x\ge0$,
\[
 |W(x)-\widetilde W(x)|+|\widetilde W'(x)|
 \le C_\alpha A(1+x)^a.
\]
Since $f(0)=0$ and $f$ has compact support, integration by parts gives
\[
 -\int W(|f|^2)'\,\mathrm dx
 =-2\operatorname{Re}\int(W-\widetilde W)f'\bar f\,\mathrm dx
   +\int\widetilde W'|f|^2\,\mathrm dx.
\]
For $0<\varepsilon\le1$, Young's inequality bounds the absolute value by
\[
 \frac{\varepsilon}{2}\|f'\|_2^2+
 C_\alpha(A^2/\varepsilon+A)
       \int(1+x)^{2a}|f(x)|^2\,\mathrm dx.
\]
For $0<a<1/2$, the maximum of $tX^{2a}-\delta X$ over $X\ge0$ is attained at $X=(2at/\delta)^{1/(1-2a)}$ and equals $(1-2a)(2a)^{2a/(1-2a)}\delta^{-2a/(1-2a)}t^{1/(1-2a)}$.
Consequently, for $t,\delta>0$ and $X\ge1$,
\[
 tX^{2a}\le\delta X+
 C_a\delta^{-2a/(1-2a)}t^{1/(1-2a)}.
\]
For $a=0$ the inequality reduces to $t\le\delta X+t$. Take $\delta=\varepsilon/2$ and $t=C_\alpha(A^2/\varepsilon+A)$. This gives the required form bound; for example,
\[
 C_{\varepsilon,A,\alpha}
 \le C_\alpha\varepsilon^{-(1+2a)/(1-2a)}
                  (1+A)^{2/(1-2a)}
 \qquad(0<\varepsilon\le1).
\]
The case $\varepsilon>1$ follows from $\varepsilon=1$.
\end{proof}

Suppose $\mathfrak h_n=\mathfrak h_n^{(0)}+\mathfrak b_n+\mathfrak r_n$, where
\[
 c_0\ell_h[u]\le\mathfrak h_n^{(0)}[u]+\kappa_0\|u\|^2
              \le C_0\ell_h[u],\qquad
 \mathfrak b_n[u]=-\int W_n(|I_hu|^2)'\,\mathrm dx.
\]
Assume that the functions $W_n$ satisfy Lemma~\ref{A:res:3.1} with a common exponent and constant, and that for every
$\varepsilon>0$,
\[
 |\mathfrak r_n[u]|\le\varepsilon\ell_h[u]
                         +C_{\varepsilon,\omega}\|u\|^2.
\]
All constants in these assumptions are independent of $n$.

\begin{proposition}\label{A:res:3.2}
Under these assumptions, there are constants $0<c_\omega\le C_\omega<\infty$ and $\kappa_\omega<\infty$ such that
\[
 c_\omega\ell_h[u]\le\mathfrak h_n[u]+\kappa_\omega\|u\|^2
                       \le C_\omega\ell_h[u].
 \eqtag{3.12}
\]
\end{proposition}

\begin{proof}
By Lemma~\ref{A:res:3.1}, the interpolation identity \eqref{A:eq:3.7}, and $\|I_hu\|_2\le\|u\|$, for every $\delta>0$,
\[
 |\mathfrak b_n[u]|+|\mathfrak r_n[u]|
 \le\delta\ell_h[u]+D_{\delta,\omega}\|u\|^2.
\]
Choose $\delta=c_0/2$ and put $\kappa_\omega=\kappa_0+D_{c_0/2,\omega}$. Then
\[
 \mathfrak h_n[u]+\kappa_\omega\|u\|^2
 \ge(c_0-\delta)\ell_h[u]=\tfrac12c_0\ell_h[u],
\]
whereas
\[
 \mathfrak h_n[u]+\kappa_\omega\|u\|^2
 \le(C_0+\delta)\ell_h[u]+2D_{c_0/2,\omega}\|u\|^2
 \le(C_0+\delta+2D_{c_0/2,\omega})\ell_h[u].
\]
This proves both bounds with constants uniform in $n$.
\end{proof}

We use \eqref{A:eq:3.12} with the same shift and constants for every pair of operators being compared. The lower bound supplies the energy and spectral estimates. The
representation \eqref{A:eq:3.8} and the coefficient bounds \eqref{A:eq:3.10} are additionally needed for exponential localization; the off-diagonal entries may have either sign.

For the remainder of this subsection, assume that \eqref{A:eq:3.8}, \eqref{A:eq:3.10}, and \eqref{A:eq:3.12} hold uniformly for all sufficiently large $n$. The coefficient
estimates for each matrix model establish these bounds.

Let $R_n(z)=(H_n-z)^{-1}$ and let $K\Subset\mathbb C\setminus\mathbb R$ be compact. We may replace $K$ by $K\cup\{\bar z:z\in K\}$, and do so throughout the weighted estimates for
the resolvents of the matrices $H_n$. Write $\delta_K=\operatorname{dist}(K,\mathbb R)>0$.
\begin{proposition}\label{A:res:3.3}
For every \(H_n\) satisfying \eqref{A:eq:3.12},
\[
\sup_{z\in K}
\left(
\|D_hR_n(z)\|_{\mathrm{op}}
+
\|(1+X_h)^{1/2}R_n(z)\|_{\mathrm{op}}
\right)
\le
C_{\omega,K}.
\eqtag{3.13}
\]
\end{proposition}

\begin{proof}
Put $u=R_n(z)f$. Self-adjointness gives
\[
 \|u\|\le\delta_K^{-1}\|f\|.
 \eqtag{3.14}
\]
Taking real parts in $\langle (H_n-z)u,u\rangle=\langle f,u\rangle$ then yields
\[
 \mathfrak h_n[u]
 =\operatorname{Re}\langle f,u\rangle+\operatorname{Re}z\,\|u\|^2
 \le C_K\|f\|^2.
 \eqtag{3.15}
\]
By \eqref{A:eq:3.12},
\[
 \|D_hu\|^2+\|(1+X_h)^{1/2}u\|^2
 =\ell_h[u]\le c_\omega^{-1}
       (\mathfrak h_n[u]+\kappa_\omega\|u\|^2)
 \le C_{\omega,K}\|f\|^2.
\]
In particular, if $r_K=\sup_{z\in K}\operatorname{Re}z$, the constant in this energy estimate can be chosen as
\[
 c_\omega^{-1}\bigl[\delta_K^{-1}
       +(r_K+\kappa_\omega)_+\delta_K^{-2}\bigr].
\]
Each operator norm in \eqref{A:eq:3.13} is at most the square root of this constant. Their sum is at most twice that quantity, uniformly in the number of vertices.
\end{proof}

For Hilbert--Schmidt estimates we also need the number of low energy states. Let $\mathcal L_h=-\Delta_{h,D}+1+X_h$, so that $\langle u,\mathcal L_hu\rangle=\ell_h[u]$. On a unit
interval the kinetic energy contributes at most $C\sqrt E$ states below energy $E$; only $O(E)$ such intervals can contribute. The following bracketing argument makes this count
uniform in the lattice spacing.
\begin{lemma}\label{A:res:3.4}
There exists a deterministic constant \(C\) such that for every \(E\ge1\),
\[
\#\{\mu\in\sigma(\mathcal L_h):\mu\le E\}
\le
CE^{3/2}.
\eqtag{3.17}
\]
Consequently, if \(\mu_{k,h}\) are the eigenvalues of \(\mathcal L_h\), then
\[
\mu_{k,h}\ge ck^{2/3}.
\eqtag{3.18}
\]
\end{lemma}

\begin{proof}
Group the indices into the disjoint sets $J_p=\{j:p\le x_j<p+1\}$, $p=0,1,\ldots$. Delete every squared difference joining two different sets and both external Dirichlet boundary
terms. The resulting quadratic form is smaller and is a direct sum over the $J_p$. Each block has Neumann boundary conditions. On $J_p$ the potential is at least $1+p$, so only
$p\le E-1$ can contribute. For a nonempty block of $m=|J_p|$ sites, the vectors $u_r=\cos(\pi k(r-\tfrac12)/m)$, $1\le r\le m$, diagonalize its Neumann difference Laplacian.
Indeed, with $\theta=\pi k/m$, they satisfy $u_0=u_1$, $u_{m+1}=u_m$ and $2u_r-u_{r-1}-u_{r+1}=4\sin^2(\theta/2)u_r$. The eigenvalues are therefore
\[
 4h^{-2}\sin^2\!\left(\frac{\pi k}{2m}\right),
 \qquad k=0,\ldots,m-1.
\]
Since $mh\le1+h\le2$ and $\sin t\ge2t/\pi$ for $0\le t\le\pi/2$, the $k$th displayed eigenvalue is at least $4k^2/(mh)^2\ge k^2$. Thus the number below $A\ge0$ is at most
$1+\sqrt A$, uniformly in $h,m$. The min--max principle gives
\[
 N_{\mathcal L_h}(E)
 \le C\sum_{0\le p\le E-1}(1+\sqrt{E-p-1})
 \le C\left(E+\int_0^E\sqrt t\,\mathrm dt\right)
 \le CE^{3/2}.
\]
Taking $E=\mu_{k,h}$ gives $k\le N_{\mathcal L_h}(\mu_{k,h})\le C\mu_{k,h}^{3/2}$, which proves \eqref{A:eq:3.18}.
\end{proof}

\begin{proposition}\label{A:res:3.5}
Let \(H_n\) satisfy \eqref{A:eq:3.12}, and let \(\lambda_{k,n}\) be its eigenvalues in increasing order. Then
\(\lambda_{k,n}+\kappa_\omega\ge c_\omega k^{2/3}\). In particular,
\[
\sup_{z\in K}\|R_n(z)\|_{\mathfrak S_2}^2
\le C_{\omega,K},
\eqtag{3.20}
\]
and
\[
\sup_{z\in K}\|D_hR_n(z)\|_{\mathfrak S_2}^2
\le C_{\omega,K}N_n^{1/3}
\le C_{\omega,K}h^{-1}.
\eqtag{derivative-HS}
\]
The derivative in the second estimate includes both Dirichlet boundary edges.
\end{proposition}

\begin{proof}
Let \(\{\psi_{k,n}\}_{k=1}^{N_n}\) be an orthonormal eigenbasis and put \(t_{k,n}=\lambda_{k,n}+\kappa_\omega\). The lower form estimate gives
\(H_n+\kappa_\omega I\ge c_\omega\mathcal L_h\). The min--max principle and Lemma~\ref{A:res:3.4} therefore give
\[
t_{k,n}\ge c_\omega\mu_{k,h}\ge c_\omega k^{2/3}>0.
\]
For \(z\in K\), the nonreal resolvent bound implies
\[
\frac{t_{k,n}}{|\lambda_{k,n}-z|}
\le
1+\frac{|z+\kappa_\omega|}{|\lambda_{k,n}-z|}
\le 1+\delta_K^{-1}\sup_{w\in K}|w+\kappa_\omega|
=:C_{\omega,K}.
\]
Consequently,
\[
|\lambda_{k,n}-z|^{-2}
\le C_{\omega,K}t_{k,n}^{-2}
\le C_{\omega,K}k^{-4/3}.
\]
Expansion in the eigenbasis gives
\[
 \|R_n(z)\|_{\mathfrak S_2}^2
 \le C_{\omega,K}\sum_{k=1}^{N_n}k^{-4/3}
 \le C_{\omega,K}\left(1+\int_1^\infty t^{-4/3}\,\dd t\right)
 =4C_{\omega,K}.
\]
This proves \eqref{A:eq:3.20}. The derivative estimate below uses the different sum $\sum_{k=1}^{N_n}k^{-2/3}\le 1+\int_1^{N_n}t^{-2/3}\,\dd t\le3N_n^{1/3}$.

Applying \eqref{A:eq:3.12} to \(\psi_{k,n}\) also gives
\[
\|D_h\psi_{k,n}\|^2
\le\ell_h[\psi_{k,n}]
\le c_\omega^{-1}t_{k,n}.
\]
Using $t_{k,n}/|\lambda_{k,n}-z|^2\le C_{\omega,K}t_{k,n}^{-1}$, we obtain
\[
\begin{aligned}
\|D_hR_n(z)\|_{\mathfrak S_2}^2& =\sum_{k=1}^{N_n} \frac{\|D_h\psi_{k,n}\|^2}{|\lambda_{k,n}-z|^2} \le C_\omega\sum_{k=1}^{N_n} \frac{t_{k,n}}{|\lambda_{k,n}-z|^2} \le C_{\omega,K}\sum_{k=1}^{N_n}t_{k,n}^{-1}\\
&\le C_{\omega,K}\sum_{k=1}^{N_n}k^{-2/3} \le C_{\omega,K}N_n^{1/3}.
\end{aligned}
\]
Since $N_nh\asymp h^{-2}$, the last bound is at most $C_{\omega,K}h^{-1}$, proving \eqref{A:eq:derivative-HS}.
\end{proof}

\paragraph{Resolvent kernels}
The spectral estimates above do not distinguish the locations of the coefficient errors. To recover this information, first estimate point evaluation in the energy norm. On an
interval of length $\rho$, the one-dimensional trace inequality has coefficients $\rho^{-1}$ and $\rho$ in front of the $L^2$ and derivative terms. Near $x_j$, these are balanced
against the potential by choosing $\rho_j=(1+x_j)^{-1/2}$. Since $x_j\le Ch^{-2}$, this length remains at least a fixed multiple of $h$.
\begin{lemma}\label{A:res:3.6}
For every \(u\in\ell_n^2\) and every vertex \(j\), \(|u_j|^2\le Ch\rho_j\ell_h[u]\).
\end{lemma}

\begin{proof}
Put \(y_j=jh\), so that \(u_j=h^{1/2}(I_hu)(y_j)\). The zero extension \(f\) of \(I_hu\) from \([0,(N_n+1)h]\) belongs to \(H^1(\mathbb R)\), since its boundary values vanish. Let
\(I_j=[y_j-\rho_j,y_j+\rho_j]\). For every \(t\in I_j\), the fundamental theorem of calculus and Cauchy--Schwarz give
\[
 |f(y_j)|^2
 \le 2|f(t)|^2+2|y_j-t|\int_{I_j}|f'(x)|^2\,\mathrm{d}x.
\]
Average this inequality over the interval of length \(2\rho_j\) to get
\[
 |f(y_j)|^2
 \le \rho_j^{-1}\int_{I_j}|f(x)|^2\,\mathrm{d}x
        +2\rho_j\int_{I_j}|f'(x)|^2\,\mathrm{d}x.
\]
On $I_j\cap\operatorname{supp}f$, we have $x\ge0$ and \(|x-x_j|\le\rho_j+h/2\le3/2\) since \(0<h\le1\). Hence
\(1+x\asymp1+x_j\) on this part of \(I_j\), with absolute constants. Using \(\rho_j^{-1}(1+x_j)^{-1}=\rho_j\) gives
\[
 |f(y_j)|^2
 \le C\rho_j\left(
       \int_0^{(N_n+1)h}|f'|^2
       +\int_0^{(N_n+1)h}(1+x)|f|^2\right).
\]
Multiplication by $h$ and \eqref{A:eq:3.7} prove the bound for every vertex. For a principal subinterval $\{p,\ldots,q\}$, extend its vector by zero to $\{1,\ldots,N_n\}$. Its
energy is $\sum_{k=p-1}^{q}|D_hu_k|^2+\sum_{k=p}^{q}(1+x_k)|u_k|^2$. The bound above applies with this energy and the original value of $x_j$.
\end{proof}

Choose \(\kappa_\omega\) large enough so that \(A_n:=H_n+\kappa_\omega I\) satisfies \(A_n\ge c_\omega \mathcal L_h\ge c_\omega I.\) Let \(G_n=A_n^{-1}.\)
\begin{proposition}\label{A:res:3.7}
For every \(j\),
\[
G_n(j,j)
\le
C_\omega h(1+x_j)^{-1/2}.
\eqtag{3.26}
\]
\end{proposition}

\begin{proof}
Set $u=G_ne_j$. Then $A_nu=e_j$, $u_j=G_n(j,j)>0$, and $\langle u,A_nu\rangle=G_n(j,j)$. Lemma~\ref{A:res:3.6} and \eqref{A:eq:3.12} give
\[
 G_n(j,j)^2=|u_j|^2\le Ch(1+x_j)^{-1/2}\ell_h[u]
 \le C_\omega h(1+x_j)^{-1/2}G_n(j,j).
\]
Division by $G_n(j,j)$ proves \eqref{A:eq:3.26}. \end{proof}

The reciprocal of the local length is $\sqrt{1+x}$. Integrating it gives the distance in which the Green function decays exponentially. We use exponential weights in the form
estimate for the positive inverse \cite{agmon1982lectures} (Chapter~1) and in the resolvent conjugation for nonreal spectral parameters \cite{combes1973asymptotic}. The diagonal
potential cancels in both calculations. Set \(\Phi(x) = \int_0^x\sqrt{1+s}\,\mathrm{d}s = \frac23\left((1+x)^{3/2}-1\right)\) and fix an index \(i\). We work first with a bounded
weight: for \(R>0\), let
\[
\psi_{i,R}(j)
=
\min\left\{|\Phi(x_j)-\Phi(x_i)|,R\right\},
\quad
\theta_j=\exp\bigl(a\psi_{i,R}(j)\bigr).
\]
Set $\theta_0=\theta_1$ and $\theta_{N_n+1}=\theta_{N_n}$. For interior edges, $|\psi_{i,R}(j+1)-\psi_{i,R}(j)|\le h\sqrt{1+x_{j+1}}$. The bound $h\sqrt{1+x_{N_n}}\le C$, together
with $|\mathrm e^t-1|\le |t|\mathrm e^{|t|}$, gives
\[
|\theta_{j+1}-\theta_j|
\le
C a h\sqrt{1+x_j}\max(\theta_j,\theta_{j+1}),
\quad
\frac{\max(\theta_j,\theta_{j+1})}
{\min(\theta_j,\theta_{j+1})}
\le \mathrm{e}^{Ca}.
\eqtag{3.31}
\]

\begin{lemma}\label{A:res:3.8}
For $0<a\le1$,
\[
 \left|\mathfrak h_n[\theta u,\theta u]
      -\operatorname{Re}\mathfrak h_n[u,\theta^2u]\right|
 \le C_\omega a^2\ell_h[\theta u],
 \eqtag{3.32}
\]
uniformly in $R$ and the center $i$.
\end{lemma}

\begin{proof}
The diagonal contribution vanishes, since
\[
\sum_j q_j|\theta_ju_j|^2
=
\sum_jq_ju_j\overline{\theta_j^2u_j}.
\]
For the off-diagonal entries, let \(t_{j,n}\) be the total real matrix entry joining vertices \(j\) and \(j+1\). The divergence form difference term contributes \(-c_{j,n}h^{-2}\).
A one-sided first-order perturbation contributes
\(-\eta_jh^{-1}\) for \(J_-\) and \(+\eta_jh^{-1}\) for \(J_+\). Hence
\[
|t_{j,n}|
\le
C_\omega h^{-2}
+
h^{-1}\sum_{\nu=1}^{L_0}|\eta_{j,n}^{(\nu)}|
\le
C_\omega h^{-2}.
\eqtag{3.32b}
\]
For an arbitrary real symmetric tridiagonal matrix, direct expansion of the two quadratic forms gives the edge identity
\[
\mathfrak h_n[\theta u,\theta u]
-
\operatorname{Re}\mathfrak h_n[u,\theta^2u]
=
-\sum_{j=1}^{N_n-1}
 t_{j,n}(\theta_{j+1}-\theta_j)^2
 \operatorname{Re}(u_{j+1}\overline{u_j}),
\eqtag{3.32c}
\]
where \(u_0=u_{N_n+1}=0\). Indeed, the edge contribution on the left is
\[
2t_{j,n}\theta_j\theta_{j+1}
\operatorname{Re}(u_{j+1}\overline{u_j})
-
t_{j,n}(\theta_j^2+\theta_{j+1}^2)
\operatorname{Re}(u_{j+1}\overline{u_j}).
\]
Using \(2|u_ju_{j+1}|\le |u_j|^2+|u_{j+1}|^2\), \eqref{A:eq:3.31}, and \eqref{A:eq:3.32b},
\[
\begin{aligned}
\big|\mathfrak h_n[\theta u,\theta u]
 -\operatorname{Re}\mathfrak h_n[u,\theta^2u]\big|
&\le C_\omega a^2 \sum_j(1+x_j) \max(\theta_j^2,\theta_{j+1}^2)
 \bigl(|u_j|^2+|u_{j+1}|^2\bigr)\\
&\le C_\omega a^2 \sum_j(1+x_j)\theta_j^2|u_j|^2
 \le C_\omega a^2\ell_h[\theta u].
\end{aligned}
\]
Estimate \eqref{A:eq:3.32} follows.
\end{proof}

\begin{proposition}\label{A:res:3.9}
There are constants $a_\omega,C_\omega>0$ such that
\[
 |G_n(j,i)|\le C_\omega h
 \frac{\mathrm e^{-a_\omega|\Phi_j-\Phi_i|}}
 {(1+x_j)^{1/4}(1+x_i)^{1/4}}.
 \eqtag{3.33}
\]
\end{proposition}

\begin{proof}
Let \(u=G_ne_i.\) Then \(A_nu=e_i.\) Since \(\theta_i=1\), \(\mathfrak a_n[u,\theta^2u] = \langle e_i,\theta^2u\rangle = \overline{u_i} = G_n(i,i)\), where
\(\mathfrak a_n[u,v]=\langle u,A_nv\rangle\). Lemma~\ref{A:res:3.8}, applied to $A_n$, gives
\[
 c_\omega\ell_h[\theta u]
 \le G_n(i,i)+C_\omega a^2\ell_h[\theta u].
\]
Choose $a$ so that $C_\omega a^2\le c_\omega/2$. Then
\[
\ell_h[\theta u]
\le
C_\omega G_n(i,i).
\eqtag{3.34}
\]
By the diagonal bound \eqref{A:eq:3.26},
\[
G_n(i,i)
\le
C_\omega h(1+x_i)^{-1/2}.
\eqtag{3.35}
\]
Lemma~\ref{A:res:3.6}, \eqref{A:eq:3.34}, and \eqref{A:eq:3.35} give
\[
 \theta_j^2|G_n(j,i)|^2
 \le Ch(1+x_j)^{-1/2}\ell_h[\theta u]
 \le C_\omega h^2(1+x_j)^{-1/2}(1+x_i)^{-1/2}.
\]
Taking square roots and letting $R\to\infty$ proves \eqref{A:eq:3.33}.
\end{proof}

In a kernel convolution, the exponential factor restricts the effective range of summation to the local decay length. In the coordinate $y=\Phi(x)$, the Jacobian is $(1+x)^{-1/2}$.
This factor accounts for the extra decay in the following sums; the factor $h^{-1}$ converts integration back to counting measure.

\begin{lemma}\label{A:res:3.10}
For every \(c>0\),
\[
\sum_k
(1+x_k)^{-1/2}
\mathrm{e}^{-c|\Phi(x_k)-\Phi(x_i)|}
\le
Ch^{-1}(1+x_i)^{-1}.
\eqtag{3.36}
\]
Moreover, for every \(0<c'<c\),
\[
\sum_k (1+x_k)^{-1/2} \mathrm{e}^{-c\left( |\Phi(x_k)-\Phi(x_i)| + |\Phi(x_k)-\Phi(x_j)| \right)} \le Ch^{-1} \frac{ \mathrm{e}^{-c'|\Phi(x_i)-\Phi(x_j)|} }{ (1+x_i)^{1/2}(1+x_j)^{1/2} }.
\eqtag{3.37}
\]
\end{lemma}

\begin{proof}
Write $C_k=[(k-1)h,kh]$ and $y_i=\Phi(x_i)$. For $x\in C_k$,
\[
 1+x\asymp1+x_k,\qquad
 |\Phi(x)-\Phi(x_k)|\le Ch\sqrt{1+x_k}\le C.
\]
Therefore, for every $b\ge0$,
\[
 h\sum_{k=1}^{N_n}(1+x_k)^{-b}\mathrm e^{-c|\Phi(x_k)-y_i|}
 \le C_{b,c}\int_0^\infty(1+x)^{-b}
                  \mathrm e^{-c|\Phi(x)-y_i|}\,\mathrm dx.
\]
Under the substitution $y=\Phi(x)$, one has $1+x=(1+3y/2)^{2/3}$ and $\mathrm dx=(1+x)^{-1/2}\mathrm dy$. For $p\ge0$,
\[
 (1+y)^{-p}\le(1+y_i)^{-p}(1+|y-y_i|)^p.
\]
Consequently,
\[
 \int_0^\infty(1+y)^{-p}\mathrm e^{-c|y-y_i|}\,\mathrm dy
 \le 2(1+y_i)^{-p}\int_0^\infty(1+t)^p\mathrm e^{-ct}\,\mathrm dt.
\]
Taking $p=(2b+1)/3$ gives
\[
 \sum_{k=1}^{N_n}(1+x_k)^{-b}\mathrm e^{-c|\Phi_k-\Phi_i|}
 \le C_{b,c}h^{-1}(1+x_i)^{-b-1/2}.
\]
For $b=1/2$ this is \eqref{A:eq:3.36}. Setting $b=3/2$ gives the exponent $-2$ used in the diagonal perturbation estimate. The extra power $-1/2$ comes from
$\dd x=(1+x)^{-1/2}\dd\Phi(x)$; the coefficient assumptions are unchanged.

For the second sum, set $u=\Phi_i$, $v=\Phi_j$, $d=|u-v|$, and $S(y)=|y-u|+|y-v|$. The bounds on each $C_k$ give
\[
 h\sum_k(1+x_k)^{-1/2}\mathrm e^{-cS(\Phi_k)}
 \le C\int_0^\infty(1+x)^{-1/2}\mathrm e^{-cS(\Phi(x))}\,\mathrm dx
 \le C\int_0^\infty(1+y)^{-2/3}\mathrm e^{-cS(y)}\,\mathrm dy.
\]
Here $S(y)\ge d$, and
\[
 (1+y)^{-2/3}\le
  \frac{(1+|y-u|)^{1/3}(1+|y-v|)^{1/3}}
 {(1+u)^{1/3}(1+v)^{1/3}}.
\]
For $\delta=(c-c')/2>0$, $(1+s)^{1/3}(1+t)^{1/3}\le C_\delta\mathrm e^{\delta(s+t)}$. Consequently,
\[
\int_0^\infty(1+y)^{-2/3}\mathrm e^{-cS(y)}\,\mathrm dy \le\frac{C_\delta\mathrm e^{-c'd}}
 {(1+u)^{1/3}(1+v)^{1/3}} \int_0^\infty\mathrm e^{-\delta S(y)}\,\mathrm dy
 \le\frac{2C_\delta\delta^{-1}\mathrm e^{-c'd}} {(1+u)^{1/3}(1+v)^{1/3}}.
\]
Here the last integral is at most $\int_{\mathbb R}\mathrm e^{-\delta|y-u|}\,\mathrm dy=2/\delta$. Since $(1+\Phi_i)^{1/3}\asymp(1+x_i)^{1/2}$, multiplication by $h^{-1}$ proves
\eqref{A:eq:3.37}.
\end{proof}

\begin{proposition}\label{A:res:3.11}
For \(A_n=H_n+\kappa_\omega I\),
\[
|A_n^{-2}(i,j)|
\le
C_\omega h
\frac{
\mathrm{e}^{-c|\Phi(x_i)-\Phi(x_j)|}
}{
(1+x_i)^{3/4}(1+x_j)^{3/4}
}.
\eqtag{3.38}
\]
In particular,
\[
\|A_n^{-1}e_j\|^2
=
A_n^{-2}(j,j)
\le
C_\omega h(1+x_j)^{-3/2}.
\eqtag{3.39}
\]
\end{proposition}

\begin{proof}
By \eqref{A:eq:3.33} and \eqref{A:eq:3.37}, for a smaller $c'>0$,
\[
\begin{aligned}
 |A_n^{-2}(i,j)|
 &\le\sum_k|G_n(i,k)G_n(k,j)|\\
 &\le\frac{C_\omega h^2}{(1+x_i)^{1/4}(1+x_j)^{1/4}}
       \sum_k\frac{\mathrm e^{-c(|\Phi_i-\Phi_k|+|\Phi_k-\Phi_j|)}}
                    {(1+x_k)^{1/2}}\\
 &\le\frac{C_\omega h\,\mathrm e^{-c'|\Phi_i-\Phi_j|}}
             {(1+x_i)^{3/4}(1+x_j)^{3/4}}.
\end{aligned}
\]
This proves \eqref{A:eq:3.38}; setting $i=j$ gives \eqref{A:eq:3.39}.
\end{proof}

We transfer the decay of the positive inverse kernel to nonreal $z$ through $(H_n+\kappa_\omega)R_n(z)=I+(z+\kappa_\omega)R_n(z)$. For this identity to preserve the exponential
weight, the conjugated resolvent must be bounded in the energy space. Fix a real sequence $\psi$

satisfying
\[
 |\psi_{j+1}-\psi_j|\le Ch\sqrt{1+x_j},
 \qquad 1\le j<N_n,
 \eqtag{3.40}
\]
and set $\Theta_a=\operatorname{diag}(\mathrm e^{a\psi_j})$. Extend $\psi$ constantly across each boundary edge.

\begin{lemma}\label{A:res:3.12}
For sufficiently small $a>0$, depending only on the form constants and $K$,
\[
 \|\Theta_aR_n(z)\Theta_a^{-1}\|_{\mathrm{op}}
 \le C_{\omega,K},\qquad z\in K,
 \eqtag{3.41}
\]
and
\[
 \|\Theta_a(H_n+\kappa_\omega)R_n(z)\Theta_a^{-1}\|_{\mathrm{op}}
 \le C_{\omega,K}.
 \eqtag{3.42}
\]
The constants are uniform under truncation of $\psi$.
\end{lemma}

\begin{proof}
Put $A_n=H_n+\kappa_\omega I$ and equip $\mathcal E_h=\mathbb C^{N_n}$ with $\|u\|_{\mathcal E_h}^2=\langle u,A_nu\rangle$. Identify its dual through the scalar product for
counting measure. For $f\in\mathbb C^{N_n}$,
\[
 \|f\|_{\mathcal E_h^*}
 =\sup_{u\ne0}\frac{|\langle f,u\rangle|}{\|A_n^{1/2}u\|}
 =\|A_n^{-1/2}f\|.
\]
The equality follows from Cauchy--Schwarz after the substitution $v=A_n^{1/2}u$, with equality for $v=A_n^{-1/2}f$ when $f\ne0$. Thus
\[
 \|R_n(z)\|_{\mathcal E_h^*\to\mathcal E_h}
 =\|A_n^{1/2}R_n(z)A_n^{1/2}\|_{\mathrm{op}}
 =\sup_{\lambda\in\sigma(H_n)}
                   \frac{\lambda+\kappa_\omega}{|\lambda-z|}.
\]
By \eqref{A:eq:3.12}, \(c_\omega\ell_h[u] \le \|u\|_{\mathcal E_h}^2 \le C_\omega\ell_h[u].\) In the form sense,
\(H_n-z = A_n^{1/2} \left[I-(z+\kappa_\omega)A_n^{-1}\right] A_n^{1/2}.\) For \(\lambda\in\sigma(H_n)\),
\[
\frac{\lambda+\kappa_\omega}{|\lambda-z|}
\le
1+
\frac{|z+\kappa_\omega|}{|\lambda-z|}
\le
1+
\frac{|z+\kappa_\omega|}{\operatorname{dist}(K,\mathbb R)}.
\]
\[
\|(H_n-z)^{-1}\|_{\mathcal E_h^*\to\mathcal E_h}
\le
C_{\omega,K}.
\eqtag{3.42b}
\]

For the conjugated sesquilinear form of the difference operator, write \(r_j=\mathrm{e}^{a\psi_j}, \allowbreak\, \delta_j=a(\psi_{j+1}-\psi_j).\) The discrete product rule
\(D_h(fg)=J_h(f)D_hg+D_h(f)J_hg\) and the identities
\[
\begin{aligned}
J_h(r)_j
&=\mathrm{e}^{a(\psi_{j+1}+\psi_j)/2}\cosh(\delta_j/2),
&
D_h(r)_j
&=\frac{2}{h}\mathrm{e}^{a(\psi_{j+1}+\psi_j)/2}\sinh(\delta_j/2),
\\
J_h(r^{-1})_j
&=\mathrm{e}^{-a(\psi_{j+1}+\psi_j)/2}\cosh(\delta_j/2),
&
D_h(r^{-1})_j
&=-\frac{2}{h}\mathrm{e}^{-a(\psi_{j+1}+\psi_j)/2}\sinh(\delta_j/2)
\end{aligned}
\]
give
\[
\begin{aligned}
\mathfrak t_n[r^{-1}u,rv]-\mathfrak t_n[u,v]
={}&\sum_jc_{j,n}\Bigl[
\bigl(\cosh^2(\delta_j/2)-1\bigr)(D_hu)_j\overline{(D_hv)_j}
+h^{-1}\sinh(\delta_j)(D_hu)_j\overline{(J_hv)_j}
\\
&-h^{-1}\sinh(\delta_j)(J_hu)_j\overline{(D_hv)_j}
-4h^{-2}\sinh^2(\delta_j/2)(J_hu)_j\overline{(J_hv)_j}\Bigr].
\end{aligned}
\eqtag{3.42d}
\]
Because $x_{N_n}=O(h^{-2})$, \(|\delta_j| \le Cah\sqrt{1+x_j} \le C a.\) For \(a\le a_0\),
\[
|\cosh^2(\delta_j/2)-1| \le Ca^2h^2(1+x_j),\quad h^{-1}|\sinh\delta_j| \le Ca\sqrt{1+x_j},\quad h^{-2}\sinh^2(\delta_j/2) \le Ca^2(1+x_j).
\]
Since $|x_{j+1}-x_j|\le h$ and $0<h\le1$, the weights $1+x_j$ and $1+x_{j+1}$ differ by at most a fixed factor. Consequently,
\(\sum_j(1+x_j)|(J_hu)_j|^2 \le C\sum_j(1+x_j)|u_j|^2.\) The four terms in \eqref{A:eq:3.42d} are bounded respectively by
\[
\begin{aligned}
&C_\omega a^2\|D_hu\|\|D_hv\|,\qquad C_\omega a\|D_hu\|\|(1+X_h)^{1/2}v\|,\\
&C_\omega a\|(1+X_h)^{1/2}u\|\|D_hv\|,\qquad C_\omega a^2\|(1+X_h)^{1/2}u\| \|(1+X_h)^{1/2}v\|.
\end{aligned}
\]
The first bound uses $h^2(1+x_j)\le C$. Their sum is at most $C_\omega(a+a^2)\|u\|_{\mathcal E_h}\|v\|_{\mathcal E_h}$.

All diagonal terms are invariant under exponential conjugation:
\[
(r_j^{-1}u_j)\overline{r_jv_j}=u_j\overline{v_j}.
\eqtag{3.42h}
\]
Thus conjugation changes only the off-diagonal part.

For a one-sided first-order perturbation, the associated tridiagonal matrix has a diagonal part, which cancels by \eqref{A:eq:3.42h}, and an off-diagonal entry
\(s_j=-\eta_jh^{-1}\) for \(J_-\) and \(s_j=+\eta_jh^{-1}\) for \(J_+\). Hence the change in the form under conjugation is a sum of terms of the form
\(s_j \left[ (\mathrm{e}^{-\delta_j}-1)u_{j+1}\overline{v_j} + (\mathrm{e}^{\delta_j}-1)u_j\overline{v_{j+1}} \right]\). Using $|s_j|\le h^{-1}|\eta_j|$ and
$|\mathrm e^{\pm\delta_j}-1|\le Ca h\sqrt{1+x_j}$ gives
\[
\begin{aligned}
 |\gamma_{\eta,\sigma}[r^{-1}u,rv]-\gamma_{\eta,\sigma}[u,v]|
 &\le Ca\|\eta\|_\infty\sum_j\sqrt{1+x_j}
       (|u_{j+1}v_j|+|u_jv_{j+1}|)\\
 &\le Ca\|\eta\|_\infty\bigl(
       \|(1+X_h)^{1/2}u\|\|v\|
       +\|u\|\|(1+X_h)^{1/2}v\|\bigr)\\
 &\le C_\omega a\|\eta\|_\infty
                    \|u\|_{\mathcal E_h}\|v\|_{\mathcal E_h}.
\end{aligned}
\eqtag{3.42j}
\]
Summing over the fixed number of first-order terms and using \eqref{A:eq:3.10}, the full conjugated form \(\mathfrak h_{n,a}[u,v] = \mathfrak h_n[\Theta_a^{-1}u,\Theta_av]\)
satisfies \(|\mathfrak h_{n,a}[u,v]-\mathfrak h_n[u,v]| \le C_\omega a \|u\|_{\mathcal E_h}\|v\|_{\mathcal E_h}\) for $0<a\le a_0\le1$: the quadratic contribution obeys
$C_\omega a^2\|u\|_{\mathcal E_h}\|v\|_{\mathcal E_h} \le C_\omega a\|u\|_{\mathcal E_h}\|v\|_{\mathcal E_h}$, so the sum is bounded by $2C_\omega a$ times these two norms.

Let $\mathcal C_a:\mathcal E_h\to\mathcal E_h^*$ denote this form difference, with $\|\mathcal C_a\|\le C_1a$ for $a\le1$, and let $C_2$ bound \eqref{A:eq:3.42b} on $K$. Decreasing
the preceding admissible $a_0$ to at most $\min\{1,(2C_1C_2)^{-1}\}$ gives $\|(H_n-z)^{-1}\mathcal C_a\|_{\mathcal E_h\to\mathcal E_h}\le1/2$ for every $n$, $z\in K$ and truncated
weight under consideration. The Neumann series yields
\[
 \Theta_aR_n(z)\Theta_a^{-1}
 =\bigl[I+(H_n-z)^{-1}\mathcal C_a\bigr]^{-1}(H_n-z)^{-1},
 \qquad
 \|\Theta_aR_n(z)\Theta_a^{-1}\|_{\mathcal E_h^*\to\mathcal E_h}
 \le2C_2.
\]
Since $A_n\ge c_\omega I$,
\[
 \|f\|_{\mathcal E_h^*}\le c_\omega^{-1/2}\|f\|,
 \qquad
 \|u\|\le c_\omega^{-1/2}\|u\|_{\mathcal E_h},
\]
which proves \eqref{A:eq:3.41}. Equation~\eqref{A:eq:3.42} follows from $(H_n+\kappa_\omega)R_n(z)=I+(z+\kappa_\omega)R_n(z)$. For $\psi^{(M)}_k=\min\{|\Phi(x_k)-\Phi(x_j)|,M\}$,
the $1$-Lipschitz truncation preserves \eqref{A:eq:3.40}. Taking $M\ge\max_{k\le N_n}|\Phi(x_k)-\Phi(x_j)|$ proves both estimates for the untruncated weight with the same
constants.
\end{proof}

\begin{proposition}\label{A:res:3.13}

For every \(z\in K\),
\[
\sum_k
\mathrm{e}^{2c|\Phi(x_k)-\Phi(x_j)|}
|R_n(z)(k,j)|^2
\le
C_{\omega,K}
h(1+x_j)^{-3/2}.
\eqtag{3.44}
\]
In particular,
\[
\|R_n(z)e_j\|^2
\le
C_{\omega,K}
h(1+x_j)^{-3/2}.
\eqtag{3.45}
\]
\end{proposition}

\begin{proof}
Write $A_n=H_n+\kappa_\omega I$ and $C_z=A_nR_n(z)=I+(z+\kappa_\omega)R_n(z)$. Since $A_n$ commutes with $R_n(z)$, we have $R_n(z)=C_zG_n$. Fix $j$ and let
$\Theta_c=\operatorname{diag}(\mathrm e^{c|\Phi_k-\Phi_j|})$. Choose $c>0$ so that $c\le a_\omega/2$ in Proposition~\ref{A:res:3.9} and the conjugation estimate of
Lemma~\ref{A:res:3.12} applies. The kernel bound gives
\[
\|\Theta_cG_ne_j\|^2 \le C_\omega h^2(1+x_j)^{-1/2} \sum_k(1+x_k)^{-1/2} \mathrm
 e^{-2(a_\omega-c)|\Phi_k-\Phi_j|} \le C_\omega h(1+x_j)^{-3/2},
\]
where the second inequality is \eqref{A:eq:3.36}. Therefore
\[
 \|\Theta_cR_n(z)e_j\|
 \le\|\Theta_cC_z\Theta_c^{-1}\|_{\mathrm{op}}
                          \|\Theta_cG_ne_j\|
 \le C_{\omega,K}h^{1/2}(1+x_j)^{-3/4}.
\]
Lemma~\ref{A:res:3.12} already permits the untruncated weight in this finite-dimensional calculation. Squaring the preceding inequality gives \eqref{A:eq:3.44}; omitting its
exponential factors gives \eqref{A:eq:3.45}.
\end{proof}

\begin{proposition}\label{A:res:3.14}

For \(z\in K\),
\[
\left|
[R_n(\bar z)R_n(z)](i,j)
\right|
\le
C_{\omega,K}h
\frac{
\mathrm{e}^{-c|\Phi(x_i)-\Phi(x_j)|}
}{
(1+x_i)^{3/4}(1+x_j)^{3/4}
}.
\eqtag{3.46}
\]
\end{proposition}

\begin{proof}
We have \([R_n(\bar z)R_n(z)](i,j) = \sum_k \overline{R_n(z)(k,i)} R_n(z)(k,j)\). Using \(|\Phi_i-\Phi_j| \le |\Phi_k-\Phi_i|+|\Phi_k-\Phi_j|,\) we write
\(1 \le \mathrm{e}^{c|\Phi_k-\Phi_i|} \mathrm{e}^{c|\Phi_k-\Phi_j|} \mathrm{e}^{-c|\Phi_i-\Phi_j|}.\) Hence, by Cauchy--Schwarz,
\[
\begin{aligned}
|[R_n(\bar z)R_n(z)](i,j)|
&\le \mathrm{e}^{-c|\Phi_i-\Phi_j|} \times
 \left( \sum_k \mathrm{e}^{2c|\Phi_k-\Phi_i|} |R_n(z)(k,i)|^2 \right)^{1/2}\\
&\quad\times \left( \sum_k \mathrm{e}^{2c|\Phi_k-\Phi_j|} |R_n(z)(k,j)|^2 \right)^{1/2}.
\end{aligned}
\]
By \eqref{A:eq:3.44}, the two factors are bounded by $C_{\omega,K}h^{1/2}(1+x_i)^{-3/4}$ and $C_{\omega,K}h^{1/2}(1+x_j)^{-3/4}$, respectively. Their product, with the exponential
factor retained, gives \eqref{A:eq:3.46}. All estimates use the same compact set $K\cup\overline K$.
\end{proof}

\paragraph{Diagonal and first-order perturbations}
For a multiplier next to one resolvent, the row estimate gives the weight $(1+x)^{-3/2}$. Keeping both resolvents in the product yields a stronger weight, by the off-diagonal decay
of their kernels.

\begin{proposition}\label{A:res:3.15}
For every sequence $q=(q_j)_{j=1}^{N_n}$, with $M_q=\operatorname{diag}(q_j)$,
\[
\|M_qR_n(z)\|_{\mathfrak S_2}^2
\le
C_{\omega,K}
h\sum_j|q_j|^2(1+x_j)^{-3/2}.
\eqtag{3.47}
\]
For every sequence $e=(e_j)_{j=0}^{N_n}$, with $M_e$ acting by multiplication on the edge space,
\[
\|M_eJ_hR_n(z)\|_{\mathfrak S_2}^2
\le
C_{\omega,K}
h\sum_j|e_j|^2(1+x_j)^{-3/2}.
\]
The same estimate holds with \(J_h\) replaced by \(J_-\) or \(J_+\).
\end{proposition}

\begin{proof}
By the definition of $M_q$, \(\|M_qR_n(z)\|_{\mathfrak S_2}^2 = \sum_j|q_j|^2\|e_j^*R_n(z)\|^2\). The adjoint identity gives \(\|e_j^*R_n(z)\|^2 = \|R_n(\bar z)e_j\|^2.\) Applying
\eqref{A:eq:3.45}, \(\|e_j^*R_n(z)\|^2 \le C_{\omega,K}h(1+x_j)^{-3/2}.\) Summing in \(j\) gives \eqref{A:eq:3.47}. For an orthonormal basis $(f_k)$ of the vertex space, extend
$R_n(z)f_k$ by zero at vertices $0,N_n+1$. Since $|(J_hu)_j|^2\le(|u_j|^2+|u_{j+1}|^2)/2$,
\[
 \|M_eJ_hR_n(z)\|_{\mathfrak S_2}^2 \le\frac12\sum_{j=0}^{N_n}|e_j|^2
 \sum_k\bigl(|(R_n(z)f_k)_j|^2+|(R_n(z)f_k)_{j+1}|^2\bigr) \le
 C_{\omega,K}h\sum_{j=0}^{N_n}|e_j|^2(1+x_j)^{-3/2}.
\]
Here \eqref{A:eq:3.45} is applied at $\bar z$ to each nonzero row; the weights at two neighboring vertices differ by at most a fixed factor. For $J_-$ retain only the first term in
the inner sum; for $J_+$ retain only the second. This proves both remaining bounds.
\end{proof}

Let $H_a,H_b$ satisfy \eqref{A:eq:3.8}, \eqref{A:eq:3.10}, and \eqref{A:eq:3.12} with common constants, and put $R_a(z)=(H_a-z)^{-1}$ and $R_b(z)=(H_b-z)^{-1}$. The next estimate
retains both resolvents. Bounding one of them only in operator norm would leave the weight of \eqref{A:eq:3.47}.
\begin{proposition}\label{A:res:3.16}
For every sequence $q=(q_j)_{j=1}^{N_n}$, with $M_q=\operatorname{diag}(q_j)$,
\[
\|R_a(z)M_qR_b(z)\|_{\mathfrak S_2}^2
\le
C_{\omega,K}
h\sum_j |q_j|^2(1+x_j)^{-7/2}.
\eqtag{3.49}
\]
\end{proposition}

\begin{proof}
Let \(K_a=R_a(\bar z)R_a(z), \allowbreak\, K_b=R_b(z)R_b(\bar z).\) Then
\[
\|R_aM_qR_b\|_{\mathfrak S_2}^2
=
\operatorname{Tr}
\left[
R_b(\bar z)M_{\overline q}R_a(\bar z)R_a(z)M_qR_b(z)
\right].
\]
Expanding in coordinates,
\[
\|R_aM_qR_b\|_{\mathfrak S_2}^2
=
\sum_{i,j}
\overline{q_i}q_j
K_a(i,j)K_b(j,i).
\]
By \eqref{A:eq:3.46}, \(|K_a(i,j)K_b(j,i)| \le C_{\omega,K}h^2 (1+x_i)^{-3/2}(1+x_j)^{-3/2} \mathrm{e}^{-c|\Phi_i-\Phi_j|}\). Put
$a_{ij}=(1+x_i)^{-3/2}(1+x_j)^{-3/2} e^{-c|\Phi_i-\Phi_j|}$. Then $a_{ij}=a_{ji}\ge0$, and
\[
 \sum_{i,j}|q_iq_j|a_{ij}
 \le\frac12\sum_{i,j}(|q_i|^2+|q_j|^2)a_{ij}
 =\sum_i|q_i|^2\sum_j a_{ij}.
\]
The coordinate expansion consequently gives
\[
\|R_aM_qR_b\|_{\mathfrak S_2}^2
\le
C_{\omega,K}h^2
\sum_i |q_i|^2(1+x_i)^{-3/2}
S_i,
\eqtag{3.52}
\]
where $S_i=\sum_j(1+x_j)^{-3/2}\mathrm e^{-c|\Phi_i-\Phi_j|}$. Apply the weighted sum proved in Lemma~\ref{A:res:3.10} with $b=3/2$. Its Jacobian contributes another half power, so
that
\[
S_i\le Ch^{-1}(1+x_i)^{-2}.
\eqtag{3.53}
\]
Substitution in \eqref{A:eq:3.52} yields
\[
 \|R_aM_qR_b\|_{\mathfrak S_2}^2
  \le C_{\omega,K}h\sum_i|q_i|^2
                   (1+x_i)^{-3/2}(1+x_i)^{-2},
\]
which is \eqref{A:eq:3.49}.
\end{proof}

A derivative of a coefficient error is handled differently. For a real

edge sequence $e=(e_j)$, define $\Gamma_e=-(D_h^*M_eJ_h+J_h^*M_eD_h)$. Its form is
\[
\langle u,\Gamma_eu\rangle
=
-2\operatorname{Re}
\sum_j e_j(D_hu)_j\overline{(J_hu)_j}.
\]

Indeed, $2\operatorname{Re}((D_hu)_j\overline{(J_hu)_j}) =h^{-1}(|u_{j+1}|^2-|u_j|^2)$, so summation by parts identifies this form with the diagonal multiplier $(e_j-e_{j-1})/h$.
In each term of $\Gamma_e$, we combine one derivative with a resolvent in operator norm; the remaining multiplier is controlled by Proposition~\ref{A:res:3.15}. The same
factorization applies to $J_-$ and $J_+$, although their first-order forms need not be diagonal.

\begin{proposition}\label{A:res:3.17}
For $H_a,H_b$ satisfying the common hypotheses \eqref{A:eq:3.8}, \eqref{A:eq:3.10}, and \eqref{A:eq:3.12},
\[
\|R_a(z)\Gamma_eR_b(z)\|_{\mathfrak S_2}^2
\le
C_{\omega,K}
\sum_j h(1+x_j)^{-3/2}|e_j|^2.
\eqtag{3.55}
\]
The same estimate holds if \(J_h\) is replaced by \(J_-\) or \(J_+\).
\end{proposition}

\begin{proof}
Fix $J\in\{J_h,J_-,J_+\}$. For this choice set $\Gamma_e=-(D_h^*M_eJ+J^*M_eD_h)$. The identity $R_aD_h^*=(D_hR_a(\bar z))^*$ and the ideal property give
\[
 \|R_aD_h^*M_eJR_b\|_{\mathfrak S_2}
 \le\|D_hR_a(\bar z)\|_{\mathrm{op}}
        \|M_eJR_b\|_{\mathfrak S_2}.
\]
For the other term, take the adjoint of the factor containing $M_e$:
\[
 \|R_aJ^*M_eD_hR_b\|_{\mathfrak S_2}
 \le\|M_eJR_a(\bar z)\|_{\mathfrak S_2}
        \|D_hR_b(z)\|_{\mathrm{op}}.
\]
Propositions~\ref{A:res:3.3} and~\ref{A:res:3.15} bound the squared norm of each product by $C_{\omega,K}\sum_jh(1+x_j)^{-3/2}|e_j|^2$. Squaring the sum of the two products gives
$\|R_a\Gamma_eR_b\|_{\mathfrak S_2}^2 \le2\|R_aD_h^*M_eJR_b\|_{\mathfrak S_2}^2 +2\|R_aJ^*M_eD_hR_b\|_{\mathfrak S_2}^2$, which proves \eqref{A:eq:3.55} in all three cases. Both
Dirichlet edges are included in each sum.
\end{proof}

For a diagonal difference $H_a-H_b=h^{-2}M_r$, the resolvent identity and Proposition~\ref{A:res:3.16} show that $\sum_j(1+x_j)^{-7/2}|r_j|^2=O(h^4)$ suffices for squared error
$O(h)$. For $r_j=h(e_j-e_{j-1})$, summation by parts gives $h^{-2}M_r=\Gamma_e$, so the sufficient bound in Proposition~\ref{A:res:3.17} is instead
$\sum_jh(1+x_j)^{-3/2}|e_j|^2=O(h)$. These are estimates on the same resolvent difference, obtained from different forms of its coefficient. They are used separately in the
coupling constructions.

\paragraph{Boundary estimates and changes in the kinetic coefficient}
The boundary points give two different improvements. At the right boundary, the potential is of order $h^{-2}$, and the kernel estimate already gives a small resolvent column. At
the left boundary, it is the Dirichlet condition that supplies the additional small factor.

\begin{corollary}\label{A:res:3.18}
If \(x_j\asymp h^{-2}\), then \(\|R_n(z)e_j\|^2\le C_{\omega,K}h^4\).
\end{corollary}

\begin{proof}
By \eqref{A:eq:3.45}, $\|R_n(z)e_j\|^2\le C_{\omega,K}h(1+x_j)^{-3/2} \le C_{\omega,K}h^4$. Here $1+x_j\ge c h^{-2}$, with $c>0$ independent of $n$.
\end{proof}

At the left Dirichlet boundary we have a different estimate, which

requires only the lower form bound. Let $H_a$ satisfy \eqref{A:eq:3.12} and write $R_a(z)=(H_a-z)^{-1}$.

\begin{proposition}\label{A:res:3.19}
Uniformly for $1\le j\le N_n$ and $z\in K$,
\[
 \|R_a(z)e_j\|^2\le C_{\omega,K}h^2j.
 \eqtag{3.58a}
\]
\end{proposition}

\begin{proof}

Let \(f\in\ell^2\) and put \(u=R_a(\bar z)f\). Since \(u_0=0\),
\[
u_j
=
\sum_{k=0}^{j-1}(u_{k+1}-u_k)
=
h\sum_{k=0}^{j-1}(D_hu)_k.
\]
Cauchy--Schwarz and $\|D_hR_a(z)\|_{\mathrm{op}}\le C_{\omega,K}$ give \(|u_j|^2 \le jh^2\sum_{k=0}^{j-1}|(D_hu)_k|^2 \le C_{\omega,K}jh^2\|f\|^2\). Therefore
\[
 \|R_a(z)e_j\|^2=\|e_j^*R_a(\bar z)\|^2
 =\sup_{\|f\|=1}|(R_a(\bar z)f)_j|^2
 \le C_{\omega,K}jh^2,
\]
which proves \eqref{A:eq:3.58a}.
\end{proof}

Let $H_a,H_b$ satisfy the common form bounds, and let $W_n$ have bandwidth at most a fixed integer $R_0$. Suppose
\[
 \operatorname{supp}W_n\subset\{1,\ldots,L_n+R_0\}^2,\qquad
 L_n=O(\log n),\qquad \max_{i,j}|(W_n)_{ij}|\le\alpha_n.
\]

\begin{corollary}\label{A:res:3.19A}
For the matrices just defined,
\[
 \sup_{z\in K}\|R_a(z)W_nR_b(z)\|_{\mathfrak S_2}^2
 \le C_{\omega,K,R_0}\alpha_n^2h^4(L_n+1)^4.
 \eqtag{3.58d}
\]
\end{corollary}

\begin{proof}

Let \(\mathcal E_n = \{(i,j):(W_n)_{ij}\ne0\}.\) The condition $|i-j|\le R_0$ gives \(|\mathcal E_n|\le C_{R_0}(L_n+1)\). Expanding into matrix units,
\[
R_aW_nR_b
=
\sum_{(i,j)\in\mathcal E_n}
(W_n)_{ij}
(R_ae_i)\otimes(R_b(\bar z)e_j)^*.
\]
Using \(\|\sum_{r=1}^mT_r\|_{\mathfrak S_2}^2 \le m\sum_r\|T_r\|_{\mathfrak S_2}^2\), Proposition~\ref{A:res:3.19}, and
\(i,j\le L_n+R_0\),
\[
\begin{aligned}
\|R_aW_nR_b\|_{\mathfrak S_2}^2
&\le |\mathcal E_n|\sum_{(i,j)\in\mathcal E_n}|(W_n)_{ij}|^2
\|R_ae_i\|^2\|R_b(\bar z)e_j\|^2
\\
&\le C_{\omega,K,R_0}(L_n+1)|\mathcal E_n|\alpha_n^2h^4(L_n+1)^2
\le C_{\omega,K,R_0}\alpha_n^2h^4(L_n+1)^4.
\end{aligned}
\]
Estimate \eqref{A:eq:3.58d} follows.
\end{proof}

In particular, when $h\asymp n^{-1/3}$ and, for a fixed $q$,
\[
 \alpha_n\le C_\omega h^{-1/2}(\log n)^q,
 \eqtag{3.58e}
\]
substitution in \eqref{A:eq:3.58d} gives
\[
 \sup_{z\in K}\|R_a(z)W_nR_b(z)\|_{\mathfrak S_2}^2
 \le C_{\omega,K,R_0}h^3(\log n)^{2q+4}=o_\omega(h).
 \eqtag{3.58f}
\]
These perturbations are also infinitesimally form bounded. Indeed, for every \(u_0=0\),
\[
\sum_{j=1}^{L_n+R_0}|u_j|^2
\le
h^2(L_n+R_0)^2
\sum_{k=0}^{L_n+R_0-1}|D_hu_k|^2,
\eqtag{3.58g}
\]
and
\[
 |\langle W_nu,u\rangle|
 \le\alpha_n\sum_{|i-j|\le R_0}|u_i u_j|
 \le(2R_0+1)\alpha_n\sum_{j=1}^{L_n+R_0}|u_j|^2.
\]
Together with \eqref{A:eq:3.58e}, this gives the relative bound
\[
C_\omega \alpha_n h^2(L_n+1)^2
\le
C_\omega h^{3/2}(\log n)^{q+2}
=o_\omega(1).
\eqtag{3.58h}
\]

Summing \eqref{A:eq:3.45} over $x_j\ge L$ gives the following tail bound.

\begin{proposition}\label{A:res:3.20}
For \(L\ge1\),
\[
\|\mathbf 1_{\{x\ge L\}}R_n(z)\|_{\mathfrak S_2}^2
\le
C_{\omega,K}L^{-1/2}.
\eqtag{3.59}
\]
\end{proposition}

\begin{proof}
Using \eqref{A:eq:3.45},
\[
\|\mathbf 1_{\{x\ge L\}}R_n(z)\|_{\mathfrak S_2}^2 = \sum_{x_j\ge L}\|R_n(\bar z)e_j\|^2 \le C_{\omega,K} \sum_{x_j\ge L} h(1+x_j)^{-3/2}.
\]
For a contributing center $x_j=(j-\tfrac12)h$, let $I_j=[x_j-h/2,x_j+h/2]$. Since $0<h\le1$, $1+x\le\tfrac32(1+x_j)$ on $I_j$. Consequently
\[
 h(1+x_j)^{-3/2}
 \le(3/2)^{3/2}\int_{I_j}(1+x)^{-3/2}\,\dd x.
\]
The intervals have disjoint interiors, and $x_j\ge L$ implies $I_j\subset[L-h/2,\infty)$. Therefore
\[
 \sum_{x_j\ge L}h(1+x_j)^{-3/2}
 \le\frac{2(3/2)^{3/2}}{\sqrt{1+L-h/2}}
 \le 2(3/2)^{3/2}L^{-1/2}.
\]
The same inequality holds for an empty sum. Substitution into the column bound proves \eqref{A:eq:3.59} with a constant independent of the number of retained coordinates.
\end{proof}

The tail can also be bounded in operator norm. The energy estimate gives \(\sum_j(1+x_j)|(R_nf)_j|^2 \le C_{\omega,K}\|f\|^2.\) Hence
\[
\|\mathbf 1_{\{x\ge L\}}R_n(z)f\|^2
\le
L^{-1}
\sum_{x_j\ge L}(1+x_j)|(R_nf)_j|^2
\le
C_{\omega,K}L^{-1}\|f\|^2.
\]
Taking the supremum over \(\|f\|=1\) gives \(\|\mathbf 1_{\{x\ge L\}}R_n(z)\|_{\mathrm{op}}^2 \le C_{\omega,K}L^{-1}\).

For a matrix unit $W=\alpha|e_i\rangle\langle e_j|$,
\[
 R_a(z)WR_b(z)=\alpha(R_a(z)e_i)(R_b(\overline z)e_j)^*,
 \qquad
 \|R_a(z)WR_b(z)\|_{\mathfrak S_2}^2
 =|\alpha|^2\|R_a(z)e_i\|^2\|R_b(\overline z)e_j\|^2.
\]
For boundary indices $i,j$, we use the sharper estimates for $R_a(z)e_i$ and $R_b(\bar z)e_j$ proved above. The logarithmic left block of the general Jacobi model is treated more
sharply by the Schur complement argument of Section~\ref{A:sec:boundary-reduction}.

It remains to treat a change in the second-order coefficient. In the applications the perturbation is $D_h^*M_dD_h+M_r$, with
\[
|d_j|\le Ch^2(1+x_j),\quad
|d_j-d_{j-1}|\le Ch^3,\quad
|r_j|\le Ch^2(1+x_j)^2.
\eqtag{3.62}
\]
At $x_j\asymp h^{-2}$, the bound for $d_j$ in \eqref{A:eq:3.62} is only of order one. Taking its maximum would therefore discard the smallness near the left boundary. The product
rule below keeps $d_j$ as a multiplier and separates its first difference; the remaining Laplacian is then eliminated by the resolvent equation. Let $d=(d_j)_{j=0}^{N_n}$ be a real
or complex edge sequence and set $a_j=d_{j-1}$, $b_j=(d_j-d_{j-1})/h$ for $1\le j\le N_n$. Write $T_h=D_h^*D_h$. For the last identity below, take $H_0=T_h+M_v+Q$, where $Q$ is any
vertex operator, and $R_0(z)=(H_0-z)^{-1}$ with
$z\in\rho(H_0)$.

\begin{lemma}\label{A:res:product-rule}
With this notation,
\[
D_h^*M_dD_h=M_aT_h-M_bJ_-^*D_h.
\eqtag{product-rule}
\]
Moreover, for every matrix \(A\) on the vertex space and every diagonal multiplier \(b\),
\[
\|AM_bJ_-^*\|_{\mathfrak S_2}^2
=\|AM_b\|_{\mathfrak S_2}^2.
\eqtag{restricted-multiplier}
\]
Finally,
\[
T_hR_0(z)=I+zR_0(z)-M_vR_0(z)-QR_0(z).
\eqtag{resolvent-elimination}
\]
\end{lemma}

\begin{proof}
The derivative adjoint gives, for every \(1\le j\le N_n\),
\[
\begin{aligned}
(D_h^*M_dD_hu)_j
&=h^{-2}\bigl[d_{j-1}(u_j-u_{j-1})-d_j(u_{j+1}-u_j)\bigr]\\
&=d_{j-1}h^{-2}(2u_j-u_{j-1}-u_{j+1})
-h^{-2}(d_j-d_{j-1})(u_{j+1}-u_j)\\
&=a_j(T_hu)_j-b_j(J_-^*D_hu)_j.
\end{aligned}
\]
This proves \eqref{A:eq:product-rule}. In particular, the first and last coordinates of its left side are
\[
h^{-2}\bigl[(d_0+d_1)u_1-d_1u_2\bigr],\quad h^{-2}\bigl[(d_{N_n-1}+d_{N_n})u_{N_n}
 -d_{N_n-1}u_{N_n-1}\bigr].
\]
These expressions also follow from the right hand side using $u_0=u_{N_n+1}=0$. Thus the identity uses $d_0$ and $d_{N_n}$ as specified, with no additional boundary condition.

Since \(J_-^*J_-=I\), the finite-dimensional trace identity gives
\[
\|AM_bJ_-^*\|_{\mathfrak S_2}^2 =\operatorname{Tr}\bigl(AM_bJ_-^*J_-M_{\overline b}A^*\bigr)
 =\operatorname{Tr}\bigl(AM_bM_{\overline b}A^*\bigr) =\|AM_b\|_{\mathfrak S_2}^2.
\]
Finally, \((T_h+M_v+Q-z)R_0(z)=I\) gives \eqref{A:eq:resolvent-elimination} by rearrangement.
\end{proof}

The term $Q$ in \eqref{A:eq:resolvent-elimination} must be retained: when the second resolvent contains a random potential or a first-order term, that term also appears after
eliminating $T_h$. For the first-order term in \eqref{A:eq:product-rule}, \eqref{A:eq:restricted-multiplier} and Proposition~\ref{A:res:3.15} give
\[
\|R_a(z)M_bJ_-^*D_hR_b(z)\|_{\mathfrak S_2}^2
\le C_{\omega,K}h\sum_{j=1}^{N_n}|b_j|^2(1+x_j)^{-3/2}.
\]
Indeed, factor off \(D_hR_b\), use Proposition~\ref{A:res:3.3}, and apply Proposition~\ref{A:res:3.15} to the adjoint of
\(R_aM_b\). Sections~\ref{A:sec:the-mean-dumitriu-edelman-coefficients}
and~\ref{A:sec:comparison-of-the-deterministic-coefficients} apply this bound to their smoothly varying coefficients.

After these estimates, a fixed number of comparisons can be summed without changing the order. Let $R_n^{[0]}(z),\ldots,R_n^{[m]}(z)$ act on the same Hilbert space, where $m$ is
independent of $n$, and suppose
\[
 \sup_{z\in K}\|R_n^{[r]}(z)-R_n^{[r-1]}(z)\|_{\mathfrak S_2}^2
 \le\varepsilon_{r,n},\qquad
 \sum_{r=1}^m\varepsilon_{r,n}\le C_\omega h_n.
\]

\begin{lemma}\label{A:res:3.21}
Under these assumptions,
\[
 \sup_{z\in K}\|R_n^{[m]}(z)-R_n^{[0]}(z)\|_{\mathfrak S_2}^2
 \le m C_\omega h_n.
\]
\end{lemma}

\begin{proof}
For fixed $z$ let $E_{r,n}(z)=R_n^{[r]}(z)-R_n^{[r-1]}(z)$, $1\le r\le m$. Cancellation of consecutive terms gives
\[
 R_n^{[m]}(z)-R_n^{[0]}(z)=\sum_{r=1}^m E_{r,n}(z).
\]
These are operators on the same Hilbert space and each difference is Hilbert--Schmidt under the hypotheses. The triangle inequality followed by Cauchy--Schwarz for $m$ real numbers
gives
\[
 \left\|\sum_{r=1}^mE_{r,n}(z)\right\|_{\mathfrak S_2}^2
 \le\left(\sum_{r=1}^m\|E_{r,n}(z)\|_{\mathfrak S_2}\right)^2
 \le m\sum_{r=1}^m\|E_{r,n}(z)\|_{\mathfrak S_2}^2.
\]
Taking the supremum over the prescribed compact spectral set yields
\[
 \sup_{z\in K}\|R_n^{[m]}(z)-R_n^{[0]}(z)\|_{\mathfrak S_2}^2
 \le m\sum_{r=1}^m\sup_{z\in K}\|E_{r,n}(z)\|_{\mathfrak S_2}^2.
\]
If the individual bounds begin at indices $N_r$, the displayed estimate holds with all of them whenever $n\ge\max_{r\le m}N_r$. The number $m$ is fixed, so a finite sum of bounds
of the stated order has the same order. For random constants this deduction is made on the intersection of the finitely many probability-one events.
\end{proof}

\subsection{Bidiagonal inverses and the stochastic Bessel operator}
\label{sec:bidiagonal-inverses-and-the-stochastic-bessel-realization}
\label{B:sec:inverse-representations-and-the-stochastic-bessel-operator}

The entries of $B^{-1}$ contain products of successive ratios $Y_r/X_r$. Taking logarithms turns these products into sums. In the positive inverse $B^{-\mathsf T}B^{-1}$, the
common normalization of these products cancels. We use these identities to compare the finite inverse with the kernel of $\mathcal H_{\beta,a}^{-1}$. The formulas below keep the
factors $n$, $\alpha_n$, and the cell lengths $\Delta_{n,k}$ explicit. For the polynomial model, the change of variables in \eqref{B:eq:1-12}--\eqref{B:eq:1-16} also changes the
embedding; its Jacobian is compared with the normalized cell indicators before the final estimate.

For the limiting operator, we identify the range of the positive product with its differential domain. This identifies the two boundary conditions, while the associated closed form
gives the energy estimates. Restriction of the factor and the logarithmic coefficient errors are treated at the end of this subsection. The bidiagonal identity is proved in
\cite{ramirez2009diffusion}; the limiting factor and the coordinate change are those of \cite{rider2019universality}, with the kernel convention fixed below.

\subsubsection{Integral kernels and embeddings}
\label{B:subsection:kernel-convention-and-normalized-cells}
\label{B:subsection:notation-and-probability-quantifiers}

The inverse formulas involve both matrix entries and integral kernels. Their normalization is fixed by the cell embedding, so the kernel convention and the cell lengths need to be
specified together.

On every complex Hilbert space, scalar products are linear in the first argument and conjugate linear in the second. Thus $\langle f,g\rangle=\int_0^1 f(s)\overline{g(s)}\,\dd s$
on $L^2(0,1)$. For an integral operator $T$ we use the convention
\begin{equation}
 (Tf)(t)=\int_0^1 K_T(s,t)f(s)\,\dd s.
 \label{B:eq:integral-kernel-convention}
\end{equation}
In particular, the first variable of $K_T$ is integrated. Operator and Schatten norms are denoted by $\|\cdot\|_{\rm op}$ and $\|\cdot\|_{\mathfrak S_p}$. For a random variable
$F$, the notation $\|F\|_p=(\mathbb E|F|^p)^{1/p}$ refers to the probability measure specified in the statement. The notation $O_\omega(a_n)$ permits a finite random constant and a
finite random initial index. Almost sure statements are for fixed model and construction parameters.

Let $0=z_{n,0}<\cdots<z_{n,n}=1$ and write
\begin{equation}
 I_{n,k}=[z_{n,k-1},z_{n,k}),\qquad
 \Delta_{n,k}=z_{n,k}-z_{n,k-1},\qquad
 \widetilde{\mathscr J}_ne_k=\Delta_{n,k}^{-1/2}\mathbf1_{I_{n,k}}.
 \label{B:eq:cell-embedding-definition}
\end{equation}
The values at the finitely many cell boundaries do not affect these operators. The cells are disjoint up to null sets, and hence
\[
 \langle\widetilde{\mathscr J}_ne_i,
             \widetilde{\mathscr J}_ne_j\rangle=\delta_{ij},\qquad
 (\widetilde{\mathscr J}_n^*f)_i
 =\Delta_{n,i}^{-1/2}\int_{I_{n,i}}f(s)\,\dd s.
\]
Thus $\widetilde{\mathscr J}_n$ is an isometry. For a matrix $A$ and $T_A=\widetilde{\mathscr J}_nA\widetilde{\mathscr J}_n^*$, substitution of the adjoint formula gives, for
$t\in I_{n,j}$,
\[
 (T_Af)(t)=\sum_{i=1}^n\int_{I_{n,i}}
              \frac{A_{ji}}{\sqrt{\Delta_{n,i}\Delta_{n,j}}}f(s)\,\dd s.
\]
Consequently,
\begin{equation}
 K_{T_A}(s,t)=\frac{A_{ji}}{\sqrt{\Delta_{n,i}\Delta_{n,j}}},
 \qquad (s,t)\in I_{n,i}\times I_{n,j},
 \label{B:eq:2-1}
\end{equation}
and
\begin{equation}
 \|T_A\|_{\mathfrak S_2}^2
 =\sum_{i,j}\int_{I_{n,i}\times I_{n,j}}
       \frac{|A_{ji}|^2}{\Delta_{n,i}\Delta_{n,j}}\,\dd s\,\dd t
 =\sum_{i,j}|A_{ji}|^2=\|A\|_{\mathfrak S_2}^2.
 \label{B:eq:2-2}
\end{equation}
Put $E_n=\widetilde{\mathscr J}_n\widetilde{\mathscr J}_n^*$. Its action is cellwise averaging:
\[
 (E_nf)(t)=\Delta_{n,k}^{-1}\int_{I_{n,k}}f(s)\,\dd s,
 \qquad t\in I_{n,k}.
\]
The map $(v,g)\mapsto\widetilde{\mathscr J}_nv+g$ is unitary from $\mathbb C^n\oplus\ker\widetilde{\mathscr J}_n^*$ onto $L^2(0,1)$: its squared norm is $|v|^2+\|g\|_2^2$, and its
inverse sends $f$ to $(\widetilde{\mathscr J}_n^*f,(I-E_n)f)$. Moreover,
\[
 T_A(\widetilde{\mathscr J}_nv+g)=\widetilde{\mathscr J}_nAv.
\]
Thus $T_A$ is unitarily equivalent to $A\oplus0$. Applying the same identity to $T_A^*T_A$ proves agreement of all nonzero singular values, including multiplicities.

For $k\ge1$ put $u_{n,k}=-\log z_{n,k}$. The two successive logarithmic lengths used in the coefficient formulas are
\begin{equation}
 \begin{aligned}
 h_{n,k}^+&=u_{n,k}-u_{n,k+1}
       =\log(z_{n,k+1}/z_{n,k}),&&1\le k<n,\\
 h_{n,k}^-&=u_{n,k-1}-u_{n,k}
       =\log(z_{n,k}/z_{n,k-1}),&&2\le k\le n.
 \end{aligned}
 \label{B:eq:spectral-cell-conventions}
\end{equation}
The ratio $Y_k/X_k$ uses the increment of length $h_{n,k}^+$; the cell $I_{n,k}$ uses $h_{n,k}^-$. Thus $h_{n,k}^+=h_{n,k+1}^-$. In the operator estimates we write
$h_{n,k}=h_{n,k}^-$. Every logarithmic length in those estimates has a positive left cell boundary.

Conditional expectations below fix the exterior coefficients. The notation $\mathbb E_q$ denotes the conditional Gibbs law with those boundary values. Estimates for a coupling use
its joint law, including any auxiliary variables. Uniformity in the boundary configurations refers to the deterministic set specified in each conditional bound.

\subsubsection{The inverse of a bidiagonal matrix}
\label{B:subsection:bidiagonal-inverses-and-product-variables}

The product structure is already visible in the exact finite inverse. For the matrix $B(X,Y)$ of Section~\ref{B:sec:reorganized-1}, the diagonal entries are $B_{kk}=X_k>0$ and the
subdiagonal entries are $B_{k+1,k}=-Y_k<0$. In the following identities the dimension $n$ is fixed, and empty products equal one.

\begin{lemma}\label{B:result:2-1}
The inverse is lower triangular, and for $i\le j$,
\begin{equation}
 (B^{-1})_{ji}=\frac1{X_i}\prod_{r=i}^{j-1}\frac{Y_r}{X_{r+1}}
             =\frac1{X_j}\prod_{r=i}^{j-1}\frac{Y_r}{X_r}.
 \label{B:eq:2-3}
\end{equation}
\end{lemma}
\begin{proof}
The diagonal entries $X_k$ are positive, so the determinant of the lower bidiagonal matrix $B$ is $\prod_{k=1}^nX_k>0$. To compute its $i$th inverse column, solve $Bv=e_i$, with
$v_0=0$ and the absent subdiagonal term at zero set to zero. The equations are
\[
 X_kv_k-Y_{k-1}v_{k-1}=\mathbf1_{\{k=i\}},\qquad 1\le k\le n.
\]
For $k<i$, induction gives $v_k=0$. The equation at $i$ gives $v_i=X_i^{-1}$. At each $k>i$, $v_k=(Y_{k-1}/X_k)v_{k-1}$, so induction gives
\[
 v_k=\frac{1}{X_i}\prod_{j=i}^{k-1}\frac{Y_j}{X_{j+1}}
     =\frac1{X_k}\prod_{j=i}^{k-1}\frac{Y_j}{X_j},\qquad k\ge i.
\]
Both products are empty and equal to one when $k=i$. The two expressions agree because their numerators are $\prod_{j=i}^{k-1}Y_j$ and their denominators are $\prod_{j=i}^kX_j$.
Thus
\[
 (B^{-1})_{ki}
 =\mathbf1_{\{i\le k\}}\frac1{X_k}
                   \prod_{j=i}^{k-1}\frac{Y_j}{X_j}.
\]
Substitution in the row equation verifies every column, including the first and last rows. This proves the displayed inverse formula with precisely the subdiagonal convention
$-Y_j$ used here; compare the bidiagonal formula in \cite{ramirez2009diffusion}.
\end{proof}

Define
\begin{equation}
 P_{n,k}=\prod_{r=k}^{n-1}\frac{Y_r}{X_r},\qquad P_{n,n}=1.
 \label{B:eq:2-5}
\end{equation}
The ratio $P_{n,i}/P_{n,j}$ is the product in the second expression of \eqref{B:eq:2-3}. Therefore
\begin{equation}
 (B^{-1})_{ji}=\frac{P_{n,i}}{X_jP_{n,j}}\mathbf1_{\{i\le j\}}.
 \label{B:eq:2-6}
\end{equation}

\begin{proposition}\label{B:result:2-2}
For $1\le i,j\le n$,
\begin{equation}
 (BB^{\mathsf T})^{-1}_{ij}
 =P_{n,i}P_{n,j}\sum_{\ell=i\vee j}^{n}
                         \frac1{X_\ell^2P_{n,\ell}^2}.
 \label{B:eq:2-7}
\end{equation}
\end{proposition}
\begin{proof}
The identity $(BB^{\mathsf T})^{-1}=B^{-\mathsf T}B^{-1}$ and \eqref{B:eq:2-6} give
\[
 (B^{-\mathsf T}B^{-1})_{ij}
 =\sum_{\ell=1}^n
    \frac{P_{n,i}P_{n,j}}{X_\ell^2P_{n,\ell}^2}
       \mathbf1_{\{i\le\ell\}}\mathbf1_{\{j\le\ell\}}.
\]
The two indicators restrict the sum to $\ell\ge i\vee j$. For $v\in\mathbb C^n$, this also gives
\[
 \langle(BB^{\mathsf T})^{-1}v,v\rangle
 =\sum_{\ell=1}^n\frac1{X_\ell^2P_{n,\ell}^2}
                  \left|\sum_{i=1}^{\ell}P_{n,i}v_i\right|^2
 =\|B^{-1}v\|_2^2.
\]
Put $S_0=0$ and $S_\ell=\sum_{i=1}^\ell P_{n,i}v_i$. Every weight $(X_\ell^2P_{n,\ell}^2)^{-1}$ is strictly positive. If the quadratic form is zero, each $S_\ell$ is therefore
zero, and
\[
 v_\ell=\frac{S_\ell-S_{\ell-1}}{P_{n,\ell}}=0,
 \qquad 1\le\ell\le n.
\]
The order $B^{-\mathsf T}B^{-1}$ is fixed by this identity.
\end{proof}

The coefficients used in the scaled inverse are
\begin{equation}
 q_{n,k}=\frac1{nX_k^2P_{n,k}^2}.
 \label{B:eq:2-9}
\end{equation}
Since $X_kP_{n,k}=Y_kP_{n,k+1}$ for $k<n$,
\begin{equation}
 q_{n,k}=\frac1{nY_k^2P_{n,k+1}^2}\quad(k<n),\qquad
 q_{n,n}=\frac1{nX_n^2}.
 \label{B:eq:2-10}
\end{equation}

\subsubsection{Embedded inverses and scalar normalization}
\label{B:subsection:embedded-kernels-and-scalar-invariance}

Embedding the inverse introduces a factor determined by the cell length in each variable. 

Let $\alpha_n=1$ for the beta-Laguerre model and $\alpha_n=4/L_V^2$ for the general potential model. With the embedding \eqref{B:eq:cell-embedding-definition}, set
\begin{equation}
 \mathcal R_n^\Delta=\alpha_n\widetilde{\mathscr J}_n
                  (n^2BB^{\mathsf T})^{-1}\widetilde{\mathscr J}_n^*,
 \qquad
 \mathcal K_n^\Delta=\sqrt{\alpha_n}\widetilde{\mathscr J}_n
                  (nB)^{-1}\widetilde{\mathscr J}_n^*.
 \label{B:eq:2-12}
\end{equation}
Define the kernel coefficients
\begin{equation}
 p_{n,k}=\frac{\sqrt{\alpha_n}P_{n,k}}{\sqrt{n\Delta_{n,k}}},\qquad
 r_{n,k}=\frac1{\sqrt{n\Delta_{n,k}}X_kP_{n,k}}.
 \label{B:eq:2-13}
\end{equation}
They satisfy
\begin{equation}
 \Delta_{n,k}r_{n,k}^2=q_{n,k}.
 \label{B:eq:2-19}
\end{equation}

\begin{proposition}\label{B:result:2-3}
The kernels on $I_{n,i}\times I_{n,j}$ are
\begin{align}
 K_n^\Delta(s,t)&=p_{n,i}r_{n,j}\mathbf1_{\{i\le j\}},
       \label{B:eq:2-17}\\
 R_n^\Delta(s,t)&=p_{n,i}p_{n,j}
                        \sum_{\ell=i\vee j}^{n}q_{n,\ell}.
       \label{B:eq:2-14}
\end{align}
Moreover,
\begin{equation}
 (\mathcal K_n^\Delta)^*\mathcal K_n^\Delta=\mathcal R_n^\Delta.
 \label{B:eq:2-16}
\end{equation}
\end{proposition}
\begin{proof}
Equations \eqref{B:eq:2-1} and \eqref{B:eq:2-6} give
\[
 K_n^\Delta(s,t)=
 \frac{\sqrt{\alpha_n}P_{n,i}}
 {n\sqrt{\Delta_{n,i}\Delta_{n,j}}X_jP_{n,j}}
                \mathbf1_{\{i\le j\}},
\]
which is \eqref{B:eq:2-17}. Since $\widetilde{\mathscr J}_n^*\widetilde{\mathscr J}_n=I$,
\[
 (\mathcal K_n^\Delta)^*\mathcal K_n^\Delta
 =\alpha_n\widetilde{\mathscr J}_n
       (nB)^{-\mathsf T}(nB)^{-1}\widetilde{\mathscr J}_n^*
 =\mathcal R_n^\Delta.
\]
In the kernel convention \eqref{B:eq:integral-kernel-convention}, the kernel of $K^*K$ is $\int K(s,r)\overline{K(t,r)}\,\dd r$. For the real coefficients above, this integral
equals
\[
 p_{n,i}p_{n,j}\sum_{\ell=i\vee j}^{n}
               \Delta_{n,\ell}r_{n,\ell}^2,
\]
which proves \eqref{B:eq:2-14} by \eqref{B:eq:2-19}. In particular, integrating the square of this kernel gives
\[
 \|\mathcal R_n^\Delta\|_{\mathfrak S_2}^2
 =\sum_{i,j=1}^n\Delta_{n,i}\Delta_{n,j}
       p_{n,i}^2p_{n,j}^2
          \left(\sum_{\ell\ge i\vee j}q_{n,\ell}\right)^2
 =\frac{\alpha_n^2}{n^4}\|(BB^{\mathsf T})^{-1}\|_{\mathfrak S_2}^2.
\]
For the last equality, substitute $\Delta_{n,i}p_{n,i}^2=\alpha_nP_{n,i}^2/n$ and $q_{n,\ell}=(nX_\ell^2P_{n,\ell}^2)^{-1}$, and use \eqref{B:eq:2-7}. Thus inversion precedes the
isometric embedding; the embedded operator is zero on the orthogonal complement of the cell functions.
\end{proof}

On a diagonal square $I_{n,i}\times I_{n,i}$ the factor has the constant value $p_{n,i}r_{n,i}$. This is required by the nonzero diagonal of $B^{-1}$. To quantify the effect of
imposing $s<t$, let $\mathcal K_n^{<}$ have kernel $K_n^\Delta(s,t)\mathbf1_{\{s<t\}}$. The difference is supported on one half of each diagonal square, so
\[
 \|\mathcal K_n^\Delta-\mathcal K_n^{<}\|_{\mathfrak S_2}^2
 =\frac12\sum_{i=1}^n\Delta_{n,i}^2p_{n,i}^2r_{n,i}^2
 =\frac{\alpha_n}{2n^2}\sum_{i=1}^nX_i^{-2}.
\]
We retain $\mathbf1_{\{i\le j\}}$ throughout. The trace of the positive inverse includes these same diagonal contributions:
\[
 \operatorname{Tr}\mathcal R_n^\Delta
 =\sum_{i=1}^n\Delta_{n,i}p_{n,i}^2\sum_{\ell=i}^nq_{n,\ell}
 =\sum_{\ell=1}^nq_{n,\ell}\sum_{i=1}^{\ell}\Delta_{n,i}p_{n,i}^2
 =\|\mathcal K_n^\Delta\|_{\mathfrak S_2}^2.
\]
All sums are finite, so these identities require no assumption on moments of the reciprocal coefficients.

For any $c\in\mathbb R$, multiplication of $p_{n,k}$ by $e^c$ and of $q_{n,k}$ by $e^{-2c}$ leaves \eqref{B:eq:2-14} unchanged:
\begin{equation}
 (e^cp_{n,i})(e^cp_{n,j})\sum_{\ell\ge i\vee j}e^{-2c}q_{n,\ell}
       =p_{n,i}p_{n,j}\sum_{\ell\ge i\vee j}q_{n,\ell}.
 \label{B:eq:2-20}
\end{equation}
The corresponding change $r_{n,k}\mapsto e^{-c}r_{n,k}$ also leaves the factor kernel unchanged. The product $P_{n,k}$ in \eqref{B:eq:2-5} has the fixed normalization $P_{n,n}=1$;
\eqref{B:eq:2-20} permits a separate scalar normalization when comparing its kernel coefficients.

For the beta-Laguerre embedding, $z_{n,k}=k/n$ and $\Delta_{n,k}=1/n$. Hence $p_{n,k}=P_{n,k}$ and
\begin{equation}
 R_{n,\beta,a}^{\mathrm L}(s,t)
 =P_{n,i}P_{n,j}\sum_{\ell=i\vee j}^{n}
                \frac1{n(X_{n,\ell}^{\mathrm L})^2P_{n,\ell}^2}.
 \label{B:eq:2-21}
\end{equation}
For the general potential, use $z_{n,k}=\theta_V(k/n)$ and the normalization \eqref{B:eq:1-20}. Substitution into \eqref{B:eq:2-13} gives
\begin{equation}
 p_{n,k}^{V}=\frac2{L_V}\frac{P_{n,k}}{\sqrt{n\Delta_{n,k}}},\qquad
 q_{n,k}^{V}=\frac1{nX_k^2P_{n,k}^2}.
 \label{B:eq:2-22}
\end{equation}
The resulting kernel is
\begin{equation}
 R_n^{\Delta,V}(s,t)=p_{n,i}^{V}p_{n,j}^{V}
                    \sum_{\ell=i\vee j}^{n}q_{n,\ell}^{V}.
 \label{B:eq:2-24}
\end{equation}
Here the superscript $\Delta$ records the use of normalized indicators; the embedding $\mathscr J_{n,V}$ in the main theorem is compared with this embedding next.

\subsubsection{Comparison of the two embeddings}
\label{B:subsection:canonical-and-normalized-cell-embeddings}

The general potential changes the spatial variable as well as the coefficients. To compare its original embedding with normalized cell indicators, the Jacobian must be controlled
up to zero. This follows from the profile equation for $\phi$; the notation $L_V$, $\theta_V$ and $U_{\theta_V}$ remains that of \eqref{B:eq:1-12}--\eqref{B:eq:1-24}.

\begin{lemma}\label{B:lem:bessel-coordinate-regularity}
The function $\theta_V$ belongs to $C^2([0,1])$. There are constants $c_V'>0$ and $C_V'<\infty$ such that
\begin{equation}
 c_V'\le\theta_V'(t)\le C_V',\qquad
 |\theta_V''(t)|\le C_V'\quad(0\le t\le1).
 \label{B:eq:bessel-coordinate-regularity}
\end{equation}
\end{lemma}
\begin{proof}
Write $U(x)=V(x^2)$ and choose $c>0$ with $U''\ge c$. For $r\ge0$ put
\[
 J(r)=\frac1{2\pi}\int_0^{2\pi}U(2r\sin(\vartheta/2))\,\dd\vartheta,
 \qquad F(r)=\frac r2J'(r).
\]
For $z=e^{\ii\vartheta}$, $(2\sin(\vartheta/2))^{2m}=(1-z)^m(1-z^{-1})^m$. Expanding both factors and integrating uses
$(2\pi)^{-1}\int_0^{2\pi}e^{\ii(r-s)\vartheta}\,\dd\vartheta =\mathbf1_{\{r=s\}}$ and gives
\[
 \frac1{2\pi}\int_0^{2\pi}(2\sin(\vartheta/2))^{2m}\,\dd\vartheta
 =\sum_{r=0}^m\binom mr^2=\binom{2m}{m}.
\]
The last equality follows by comparing the coefficient of $z^m$ in $(1+z)^m(1+z)^m=(1+z)^{2m}$. Expanding $V$ therefore gives
\[
 J(r)=\sum_{m=1}^{d}g_m\binom{2m}{m}r^{2m},\qquad
 F(r)=\sum_{m=1}^{d}m g_m\binom{2m}{m}r^{2m}.
\]
Differentiation under the finite integral yields
\[
 J''(r)\ge\frac c{2\pi}\int_0^{2\pi}4\sin^2(\vartheta/2)\,\dd\vartheta=2c.
\]
Since $J'(0)=0$, we obtain $J'(r)\ge2cr$ and $F'(r)=(J'(r)+rJ''(r))/2\ge2cr$. Define the polynomial $\Psi(v)=\sum_{m=1}^{d}m g_m\binom{2m}{m}v^m$. Then $F(r)=\Psi(r^2)$, so
$\Psi'(v)\ge c$ for $v>0$ and $\Psi'(0)=2g_1=U''(0)>0$. Since $\Psi(0)=0$ and $\Psi(v)\ge cv$ for $v\ge0$, every $t\in[0,1]$ has a unique solution $v(t)\ge0$ to $\Psi(v(t))=t$. The
inverse function theorem applies at every point of $[0,v(1)]$, including zero, and gives
\[
 v(0)=0,\qquad v'(t)=\Psi'(v(t))^{-1}>0.
\]
Thus $v$ extends smoothly to a neighborhood of $[0,1]$. The nonnegative solution of the profile equation is
$\phi(t)=\sqrt{v(t)}$.

For $t>0$ set $f(t)=\phi(t)/\sqrt t$. The identity
\[
 \frac{v(t)}t=\int_0^1 v'(ut)\,\dd u
\]
extends $f$ to a positive smooth function on $[0,1]$, with $f(0)=(2g_1)^{-1/2}$. Substitution $s=tu$ in $L_V(t)$ gives
\[
 L_V(t)=\sqrt t\,b(t),\qquad
 b(t)=\int_0^1\frac{u^{-1/2}}{f(tu)}\,\dd u.
\]
Put $F_0=1/f$. For every nonnegative integer $r$, differentiation under the integral is justified by
\[
 b^{(r)}(t)=\int_0^1u^{r-1/2}F_0^{(r)}(tu)\,\dd u,\qquad
 |b^{(r)}(t)|\le\frac{\|F_0^{(r)}\|_\infty}{r+1/2}.
\]
Thus $\theta_V(t)=t b(t)^2/L_V^2$ is smooth through zero. Formula \eqref{B:eq:1-16}, first for $t>0$ and then by continuity at zero, gives
\[
 \theta_V'(t)=\frac{2b(t)}{L_V^2f(t)},\qquad
 \theta_V''(t)=\frac2{L_V^2}
               \left(\frac{b'(t)}{f(t)}-\frac{b(t)f'(t)}{f(t)^2}\right).
\]
If $f_-,b_->0$ are the minima and $f_+,b_+$ the maxima on $[0,1]$, these identities imply
\[
 \frac{2b_-}{L_V^2f_+}\le\theta_V'
 \le\frac{2b_+}{L_V^2f_-},\qquad
 \|\theta_V''\|_\infty
 \le\frac2{L_V^2}\left(\frac{\|b'\|_\infty}{f_-}
                   +\frac{b_+\|f'\|_\infty}{f_-^2}\right).
\]
These are \eqref{B:eq:bessel-coordinate-regularity}. Also $b(0)=2/f(0)$, so $\theta_V'(0)=4(2g_1)/L_V^2$, as in \eqref{B:eq:1-18}.
\end{proof}

The two isometries to be compared are $\mathscr J_{n,V}$ in \eqref{B:eq:1-24} and $\widetilde{\mathscr J}_{n,V}$ in \eqref{B:eq:cell-embedding-definition}, on the partition
$z_{n,k}=\theta_V(k/n)$. Their $k$th columns have the same support. The first column function contains the Jacobian of $\theta_V$; the second is constant on that cell. The estimate
below uses orthogonality between different cells, which makes its constant independent of the number of columns.

\begin{lemma}\label{B:result:2-4}
There is a unitary $W_n$ on $L^2(0,1)$ such that
\begin{equation}
 W_n\mathscr J_{n,V}=\widetilde{\mathscr J}_{n,V},\qquad
 \|W_n-I\|_{\rm op}\le C_V/n.
 \label{B:eq:2-25}
\end{equation}
For every $n\times n$ matrix $A$,
\begin{equation}
 \|\mathscr J_{n,V}A\mathscr J_{n,V}^*
       -\widetilde{\mathscr J}_{n,V}A\widetilde{\mathscr J}_{n,V}^*\|_{\mathfrak S_2}
 \le (C_V/n)\|A\|_{\mathfrak S_2}.
 \label{B:eq:2-26}
\end{equation}
\end{lemma}
\begin{proof}
Write $\psi_k=\mathscr J_{n,V}e_k$ and $\widetilde\psi_k=\widetilde{\mathscr J}_{n,V}e_k$. Both are real, positive on $I_{n,k}$, and have norm one. Set $D_k=n\Delta_{n,k}$. For
$t\in[(k-1)/n,k/n]$, Lemma~\ref{B:lem:bessel-coordinate-regularity} gives
\[
 D_k=n\int_{(k-1)/n}^{k/n}\theta_V'(u)\,\dd u,\qquad
 |D_k-\theta_V'(t)|\le\|\theta_V''\|_\infty/n.
\]
For $d=\theta_V'(t)$, both $D_k$ and $d$ are at least $c_V'$, and
\[
 \left|\sqrt{D_k/d}-1\right|
 =\frac{|D_k-d|}{\sqrt d(\sqrt{D_k}+\sqrt d)}
 \le\frac{\|\theta_V''\|_\infty}{2c_V'n}.
\]
Consequently
\[
 \sup_{z\in I_{n,k}}
 \left|\frac{\psi_k(z)}{\widetilde\psi_k(z)}-1\right|
 =\sup_{t\in[(k-1)/n,k/n]}
       \left|\sqrt{D_k/\theta_V'(t)}-1\right|\le C_V/n.
\]
Multiplying by $\widetilde\psi_k$ and integrating proves $\|\psi_k-\widetilde\psi_k\|_2\le C_V/n$. For $v\in\mathbb C^n$, disjoint supports give the exact identity
\[
 \|(\mathscr J_{n,V}-\widetilde{\mathscr J}_{n,V})v\|_2^2
 =\sum_k|v_k|^2\|\psi_k-\widetilde\psi_k\|_2^2
 \le C_V^2n^{-2}|v|^2.
\]
We extend this comparison of the columns to a unitary on $L^2(0,1)$.

Let $c_k=\langle\widetilde\psi_k,\psi_k\rangle\in(0,1]$. If $c_k=1$, use the identity on $L^2(I_{n,k})$. Otherwise put
\[
 s_k=(1-c_k^2)^{1/2},\qquad
 e_k'=(\widetilde\psi_k-c_k\psi_k)/s_k.
\]
The vectors $\psi_k,e_k'$ are orthonormal. Define $W_{n,k}$ on their span by
\[
 W_{n,k}\psi_k=c_k\psi_k+s_ke_k',\qquad
 W_{n,k}e_k'=-s_k\psi_k+c_ke_k',
\]
and let it be the identity on the orthogonal complement. Its matrix on this span is a real rotation. Direct multiplication gives
\[
 (W_{n,k}-I)^*(W_{n,k}-I)=2(1-c_k)I,
 \qquad \|W_{n,k}-I\|_{\rm op}
       =\sqrt{2(1-c_k)}=\|\psi_k-\widetilde\psi_k\|_2.
\]
The orthogonal direct sum $W_n=\bigoplus_kW_{n,k}$ is unitary, sends every $\psi_k$ to $\widetilde\psi_k$, and satisfies \eqref{B:eq:2-25}. For
$T=\mathscr J_{n,V}A\mathscr J_{n,V}^*$,
\[
 W_nTW_n^*-T=(W_n-I)TW_n^*+T(W_n^*-I).
\]
The ideal inequality and the isometry of $\mathscr J_{n,V}$ imply
\[
 \|W_nTW_n^*-T\|_{\mathfrak S_2}
 \le2\|W_n-I\|_{\rm op}\|T\|_{\mathfrak S_2}
 =2\|W_n-I\|_{\rm op}\|A\|_{\mathfrak S_2}.
\]
Increasing $C_V$ proves \eqref{B:eq:2-26}.
\end{proof}

In particular, for $\mathcal A_n=\alpha_n^{-1}n^2B_nB_n^{\mathsf T}$,
\[
 \|\mathscr J_{n,V}\mathcal A_n^{-1}\mathscr J_{n,V}^*
                  -\mathcal R_n^{\Delta,V}\|_{\mathfrak S_2}^2
 \le C_Vn^{-2}\|\mathcal A_n^{-1}\|_{\mathfrak S_2}^2.
\]
The isometry in \eqref{B:eq:2-2} identifies the last norm with $\|\mathcal R_n^{\Delta,V}\|_{\mathfrak S_2}$. Proposition~\ref{B:prop:finite-inverse-upper-bound} first bounds that
norm for the constructed arrays and then uses this inequality. The unitary $W_n$ depends only on $V$ and $n$.

\subsubsection{Coefficients of the stochastic Bessel operator}
\label{B:subsection:brownian-coefficients}

The coefficients of $\mathcal K_{\beta,a}$ are bounded by powers of the spatial variable, using the sublinear growth of the Brownian path in logarithmic time. These bounds apply to
$\mathcal Y$, $p^B$, $q^B$ and $\mathcal K_{\beta,a}$ from \eqref{B:eq:1-29}--\eqref{B:eq:intro-brownian-coefficients}. The identity $(p^B(s))^2q^B(s)=s^{-1}$ also identifies
$(p^B)^2$ as the speed density in the differential expression below.

\begin{lemma}\label{B:lem:bessel-brownian-growth}
For fixed $\beta>0$, almost surely, for every $\eta>0$ there is $C_{\omega,\eta}\ge1$ such that, for $0<s\le1$,
\begin{align}
 |\mathcal Y(s)|&\le\log C_{\omega,\eta}+\eta|\log s|,
       \label{B:eq:bessel-logarithmic-growth}\\
 C_{\omega,\eta}^{-1}s^{a/2+\eta}
 &\le p^B(s)\le C_{\omega,\eta}s^{a/2-\eta},
       \label{B:eq:bessel-p-bounds}\\
 C_{\omega,\eta}^{-2}s^{-a-1+2\eta}
 &\le q^B(s)\le C_{\omega,\eta}^{2}s^{-a-1-2\eta}.
       \label{B:eq:bessel-q-bounds}
\end{align}
\end{lemma}
\begin{proof}
Use the Brownian motion $W$ identified in \eqref{B:eq:1-29}. The reflection principle and the Gaussian tail estimate give, for every $\epsilon>0$ and integer $m\ge1$,
\[
 \mathbb P\!\left(\sup_{0\le u\le m+1}|W(u)|>\epsilon m\right)
 \le4\exp\!\left(-\frac{\epsilon^2m^2}{2(m+1)}\right).
\]
This series is summable in $m$. On the intersection of the events obtained from Borel--Cantelli for positive rational $\epsilon$, $|W(u)|/u\to0$ as $u\to\infty$. On this event set
\[
 \log C_{\omega,\eta}
 =\beta^{-1/2}\sup_{u\ge0}(|W(u)|-\sqrt\beta\eta u)_+.
\]
For sufficiently large $u$, $|W(u)|\le\sqrt\beta\eta u/2$; the remaining supremum is finite by continuity. It is measurable because it equals the supremum over nonnegative rational
$u$. The defining inequality gives $\beta^{-1/2}|W(u)|\le\log C_{\omega,\eta}+\eta u$ for all $u\ge0$. Substitute $u=-\log s$ to prove \eqref{B:eq:bessel-logarithmic-growth}.
Exponentiating $\pm\mathcal Y(s)$ and $\mp2\mathcal Y(s)$ gives \eqref{B:eq:bessel-p-bounds} and \eqref{B:eq:bessel-q-bounds}. The rational intersection was fixed before choosing
$\eta$, so the event works for every positive $\eta$.
\end{proof}

For later deterministic arguments, let $p,q$ be positive continuous functions on $(0,1]$. Write
\begin{equation}
 P(t)=\int_0^t p(s)^2\,\dd s,\qquad
 P(1)<\infty,\qquad
 \mathcal T=\int_0^1q(t)P(t)\,\dd t<\infty.
 \label{B:eq:2-31d}
\end{equation}
The first integral puts $pg$ in $L^1$ for every $g\in L^2$; the second is the integral of the squared triangular kernel. These conditions allow $q$ to be nonintegrable at zero. For
$p=p^B$ and $q=q^B$, the required integral is $\int_0^1q(t)P(t)\,\mathrm dt$. These conditions hold almost surely for $p=p^B,q=q^B$ whenever $a>-1$. Indeed, choose $\eta>0$ with
$2\eta<a+1$ and $4\eta<1$. Then
\begin{equation}
 P(t)\le\frac{C_{\omega,\eta}^2}{a+1-2\eta}t^{a+1-2\eta},\qquad
 \mathcal T\le\frac{C_{\omega,\eta}^4}{a+1-2\eta}
                   \int_0^1t^{-4\eta}\,\dd t<\infty.
 \label{B:eq:2-31e}
\end{equation}

Here $\int_0^1t^{-4\eta}\,\dd t=(1-4\eta)^{-1}$, so the two strict inequalities $a+1-2\eta>0$ and $1-4\eta>0$ supply an explicit finite bound for $\mathcal T$. The bound is
pathwise; no moment bound is claimed for the random constant
$C_{\omega,\eta}$.

\subsubsection{The Green operator}
\label{B:subsection:the-stochastic-bessel-realization}

The triangular factor gives a positive inverse by taking its adjoint product. Before identifying a differential expression, we need to know that this inverse is injective and has
dense range. Under \eqref{B:eq:2-31d}, let $\mathcal K$ have kernel $p(s)\sqrt{q(t)}\mathbf1_{\{s<t\}}$ and write $R=\mathcal K^*\mathcal K$.

\begin{proposition}\label{B:result:2-5}
The operator $R$ is positive, injective, and trace class, with dense range. Its kernel and trace are
\begin{equation}
 R(s,t)=p(s)p(t)\int_{s\vee t}^1q(r)\,\dd r,\qquad
 \operatorname{Tr}R=\|\mathcal K\|_{\mathfrak S_2}^2=\mathcal T.
 \label{B:eq:2-31}
\end{equation}
\end{proposition}
\begin{proof}
Tonelli's theorem gives
\[
 \int_0^1\int_0^t p(s)^2q(t)\,\dd s\,\dd t=\mathcal T.
\]
Thus $\mathcal K$ is Hilbert--Schmidt. For $g,h\in L^2(0,1)$, Cauchy--Schwarz on the square gives
\[
 \int_0^1\int_0^1
 |\mathcal K(s,t)g(s)\overline{h(t)}|\,\dd s\,\dd t
 \le\mathcal T^{1/2}\|g\|_2\|h\|_2.
\]
Fubini therefore identifies the adjoint and gives the pointwise representatives
\begin{equation}
 (\mathcal K g)(t)=\sqrt{q(t)}\int_0^t p(s)g(s)\,\dd s,\qquad
 (\mathcal K^*h)(s)=p(s)\int_s^1\sqrt{q(t)}h(t)\,\dd t.
 \label{B:eq:2-31f}
\end{equation}
The first integral is finite because $pg\in L^1(0,1)$. For the second, $q$ is bounded on every $[s,1]$, with $s>0$.

Let $F_g(t)=\int_0^t p(r)g(r)\,\dd r$. Substitution into the adjoint formula gives
\begin{equation}
 (Rg)(s)=p(s)\int_s^1q(t)F_g(t)\,\dd t.
 \label{B:eq:2-31i}
\end{equation}
For fixed $s>0$, the absolute double integral is at most $\|p\|_2\|g\|_2\int_s^1q(t)\,\dd t$, so its order may be exchanged. This yields the kernel in \eqref{B:eq:2-31}. The kernel
also lies in $L^2((0,1)^2)$: by Cauchy--Schwarz in
$r$,
\[
 |R(s,t)|^2\le d(s)d(t),\qquad
 d(s)=p(s)^2\int_s^1q(r)\,\dd r,\qquad
 \int_0^1d(s)\,\dd s=\mathcal T.
\]
In particular $\|R\|_{\mathfrak S_2}\le\mathcal T$.

For an orthonormal basis $(e_j)$, Parseval in the integrated variable and Tonelli give
\[
 \sum_j\langle Re_j,e_j\rangle
 =\sum_j\|\mathcal K e_j\|_2^2
 =\int_0^1\int_0^1|\mathcal K(s,t)|^2\,\dd s\,\dd t=\mathcal T.
\]
Since $R=\mathcal K^*\mathcal K\ge0$, this proves that $R$ is trace class with the stated trace.

If $\mathcal K g=0$, then $F_g=0$ almost everywhere. The function $F_g$ is absolutely continuous on $[0,1]$, so it vanishes everywhere. Its derivative $pg$ is zero almost
everywhere, and positivity of $p$ implies $g=0$. If $\mathcal K^*h=0$, put $H_h(s)=\int_s^1\sqrt{q(t)}h(t)\,\dd t$. On $[\epsilon,1]$ this function is absolutely continuous,
because $\sqrt q\,h\in L^1$ there. The adjoint formula and $p>0$ give $H_h=0$ almost everywhere, hence everywhere on that interval. Consequently $H_h'=-\sqrt q\,h=0$ almost
everywhere. Taking $\epsilon=1/m$ for all $m\ge2$ proves $h=0$ on $(0,1)$. Finally $Rg=0$ implies $\|\mathcal K g\|_2^2=\langle Rg,g\rangle=0$, so $R$ is injective.
Self-adjointness gives $\overline{\operatorname{ran}R}=(\ker R)^{\perp}=L^2(0,1)$. Also
$\|R\|_{\mathrm{op}}=\|\mathcal K\|_{\mathrm{op}}^2\le\|\mathcal K\|_{\mathfrak S_2}^2=\mathcal T$.
\end{proof}

With $p=p^B,q=q^B$, the proposition proves the representation of $\mathcal H_{\beta,a}^{-1}$ introduced in \eqref{B:eq:1-1}.

For $p^B,q^B$, finiteness of the expected trace imposes a further condition on $a$. For $0<s<t\le1$, the centered Gaussian variable $\mathcal Y(s)-\mathcal Y(t)$ has variance
$\beta^{-1}\log(t/s)$. Its exponential moment is therefore
\[
 \mathbb E e^{2(\mathcal Y(s)-\mathcal Y(t))}
 =\exp\!\left(\frac2\beta\log\frac ts\right)
 =\left(\frac ts\right)^{2/\beta}.
\]
Apply Tonelli's theorem to the nonnegative squared factor kernel and substitute $s=tu$. The resulting equalities hold in
$[0,\infty]$:
\[
 \begin{aligned}
 \mathbb E\operatorname{Tr}\mathcal H_{\beta,a}^{-1}
 &=\mathbb E\|\mathcal K_{\beta,a}\|_{\mathfrak S_2}^2
 =\int_0^1\frac1t\int_0^t
       \left(\frac st\right)^{a-2/\beta}\,\dd s\,\dd t\\
 &=\int_0^1\int_0^1u^{a-2/\beta}\,\dd u\,\dd t
 =\int_0^1u^{a-2/\beta}\,\dd u.
 \end{aligned}
\]
For $v=a-2/\beta\ne-1$, the integral over $[\epsilon,1]$ is $(1-\epsilon^{v+1})/(v+1)$; for $v=-1$ it is $-\log\epsilon$. Letting $\epsilon\downarrow0$ gives
\[
 \mathbb E\operatorname{Tr}\mathcal H_{\beta,a}^{-1}
 =\begin{cases}
 (a+1-2/\beta)^{-1},&a+1>2/\beta,\\
 +\infty,&a+1\le2/\beta.
 \end{cases}
\]
Almost sure finiteness of the trace holds throughout $a>-1$. The tightness arguments below use finite random bounds: if $0\le X_n\le Z<\infty$ almost surely for every $n$, then
\[
 \sup_n\mathbb P(X_n>L)\le\mathbb P(Z>L)\longrightarrow0
 \qquad(L\to\infty).
\]
An $L^p$ bound for a coupled difference requires a moment estimate for that difference; the trace calculation above concerns the limiting operator alone.

\subsubsection{The operator domain}
\label{B:subsection:bessel-operator-domain}

The range of the Green operator determines the correct differential domain. Expressing a function through its weighted derivative reveals the two boundary conditions without
differentiating $p^B$ or $q^B$.

Continue under \eqref{B:eq:2-31d}. For $f\in L^2(0,1)$, use the notation
\begin{equation}
 y_f=f/p,\qquad v_f=q^{-1}y_f'
 \label{B:eq:2-31a}
\end{equation}
whenever $y_f$ has a locally absolutely continuous representative. Let $\mathcal D_{p,q}$ consist of the functions $f\in L^2(0,1)$ for which
\begin{equation}
 \begin{gathered}
 y_f,v_f\in\operatorname{AC}([\epsilon,1])\quad(0<\epsilon<1),
 \qquad p^{-1}v_f'\in L^2(0,1),\\
 y_f(1)=0,\qquad \lim_{s\downarrow0}v_f(s)=0.
 \end{gathered}
 \label{B:eq:bessel-domain-conditions}
\end{equation}
The relation in \eqref{B:eq:2-31a} holds almost everywhere. All differentiation is applied to $y_f$ and $v_f$; continuity of $p,q$ suffices. If the local regularity and the $L^2$
condition in the first line of \eqref{B:eq:bessel-domain-conditions} hold, then
\[
 \int_0^1|v_f'(s)|\,\dd s
 \le\|p\|_2\|p^{-1}v_f'\|_2<\infty.
\]
Thus $v_f(0+)$ exists: for $0<r<s$, $|v_f(s)-v_f(r)|\le\int_r^s|v_f'|$, whose right side tends to zero as $r,s\downarrow0$. Membership in the domain fixes this trace. At the
regular boundary $p(1)>0$, so $y_f(1)=0$ is equivalent to $f(1)=0$.

\begin{theorem}\label{B:result:2-5a}
The inverse $H=R^{-1}$ is positive and self-adjoint, with compact resolvent, and
\begin{equation}
 \operatorname{dom}H=\operatorname{ran}R=\mathcal D_{p,q},\qquad
 Hf=-p^{-1}v_f'.
 \label{B:eq:2-31b}
\end{equation}
In particular these statements apply to $\mathcal H_{\beta,a}$ for fixed $\beta>0$ and $a>-1$, almost surely.
\end{theorem}
\begin{proof}
The compact self-adjoint spectral theorem, together with Proposition~\ref{B:result:2-5}, gives an orthonormal basis $(e_j)$ and positive eigenvalues $\mu_j\downarrow0$, counted
with multiplicity, such that $Re_j=\mu_je_j$ and $\sum_j\mu_j=\mathcal T$. There are infinitely many $\mu_j$, since an injective operator on $L^2(0,1)$ cannot have finite rank. For
$f_j=\langle f,e_j\rangle$, the inverse is
\begin{equation}
 \operatorname{dom}H=\left\{f:\sum_j|f_j|^2\mu_j^{-2}<\infty\right\},
 \qquad Hf=\sum_j\mu_j^{-1}f_je_j.
 \label{B:eq:2-31h}
\end{equation}
This domain equals $\operatorname{ran}R$: $f=Rg$ means $f_j=\mu_jg_j$, and conversely the indicated sum of squares is finite, so the series $g=\sum_j\mu_j^{-1}f_je_j$ converges in
the Hilbert space.

For $f$ in this domain, its spectral truncations $f^{(m)}= \sum_{j\le m}f_je_j$ satisfy
\[
 \|f-f^{(m)}\|_2^2+\|H(f-f^{(m)})\|_2^2
 =\sum_{j>m}(1+\mu_j^{-2})|f_j|^2\longrightarrow0.
\]
In particular the domain is dense. If $h\in\operatorname{dom}H^*$ and $H^*h=k$, testing against $e_j$ gives $k_j=\mu_j^{-1}h_j$. Hence $\sum_j|h_j|^2\mu_j^{-2}=\|k\|_2^2<\infty$.
Conversely, that condition and Cauchy--Schwarz give $\langle Hf,h\rangle=\langle f,Hh\rangle$ for every $f\in\operatorname{dom}H$. Thus $H^*=H$. Formula \eqref{B:eq:2-31h} also
gives $H\ge\mu_1^{-1}I$ and
\[
 \left\|(H+I)^{-1}-\sum_{j=1}^m
       \frac{\mu_j}{1+\mu_j}|e_j\rangle\langle e_j|\right\|_{\rm op}
 =\sup_{j>m}\frac{\mu_j}{1+\mu_j}\longrightarrow0.
\]
This proves compactness of $(H+I)^{-1}$. The spectral formulas also give $HRg=g$ for $g\in L^2$ and $RHf=f$ for $f\in\operatorname{dom}H$; thus $0\in\rho(H)$ and $H^{-1}=R$. For
any $z\in\rho(H)$,
\[
 (H-z)^{-1}e_j=\frac{\mu_j}{1-z\mu_j}e_j,\qquad
 \sup_{j>m}\left|\frac{\mu_j}{1-z\mu_j}\right|\longrightarrow0.
\]
Indeed $|1-z\mu_j|\ge1/2$ for all sufficiently large $j$. The finite spectral sums therefore converge in operator norm and every resolvent is compact.

For $f=Rg$, the function $F_g(t)=\int_0^t p(r)g(r)\,\dd r$ satisfies
\[
 |F_g(t)|\le P(t)^{1/2}\|g\|_{L^2(0,t)},\qquad F_g(0)=0.
\]
It is absolutely continuous on $[0,1]$. On $[\epsilon,1]$, $qF_g$ is continuous and integrable. Differentiating \eqref{B:eq:2-31i} there gives
\begin{equation}
 y_f'=-qF_g,\qquad v_f=-F_g,\qquad v_f'=-pg.
 \label{B:eq:2-31j}
\end{equation}
The local absolute continuity follows from these formulas. Moreover $y_f(1)=0$, $v_f(0+)=0$ because $pg\in L^1$, and $-p^{-1}v_f'=g$. Thus
$\operatorname{ran}R\subset\mathcal D_{p,q}$.

Conversely, take $f\in\mathcal D_{p,q}$ and put $g=-p^{-1}v_f'\in L^2$. Since $p\in L^2$, the derivative $v_f'=-pg$ belongs to $L^1(0,1)$. For $0<\epsilon<s$, absolute continuity
gives $v_f(s)-v_f(\epsilon)=-\int_\epsilon^s pg$. The omitted integral satisfies $|\int_0^\epsilon pg|\le P(\epsilon)^{1/2}\|g\|_2\to0$. Together with $v_f(\epsilon)\to0$, this
gives $v_f(s)=-F_g(s)$. On $[s,1]$ the product $qF_g$ is integrable. Integrating $y_f'=qv_f$ from $s$ to $1$ and using $y_f(1)=0$ now gives
\[
 y_f(s)=\int_s^1q(t)F_g(t)\,\dd t.
\]
Equation \eqref{B:eq:2-31i} implies $f=Rg$, proving the reverse inclusion and the action in \eqref{B:eq:2-31b}. For $p^B,q^B$, the domain identification follows from
\eqref{B:eq:2-31e}.
\end{proof}

\subsubsection{A closed sesquilinear form}
\label{B:subsection:bessel-closed-form}

The same realization can be defined by its energy form. Moving the factor $p$ into the measure leaves a derivative form whose domain expresses the regular boundary condition and
finite energy at the singular boundary.

Let $\mathcal H_p=L^2((0,1),p^2\,\dd s)$, with scalar product $\langle y,z\rangle_{\mathcal H_p}=\int_0^1p^2y\overline z\,\dd s$. Multiplication by $p$ defines a unitary
$M_p:\mathcal H_p\to L^2(0,1)$. Define the complex form domain and the sesquilinear form by
\begin{equation}
 \begin{split}
 \mathcal V=\{y\in\mathcal H_p:\;&y\in\operatorname{AC}([\epsilon,1])\ (0<\epsilon<1),
            \ y(1)=0,\ \int_0^1|y'|^2/q<\infty\},\\
 \mathfrak a(y,z)&=\int_0^1q(s)^{-1}y'(s)\overline{z'(s)}\,\dd s,
                    \qquad y,z\in\mathcal V.
 \end{split}
 \label{B:eq:2-31m}
\end{equation}
The form is linear in $y$, conjugate linear in $z$, and Hermitian. Its quadratic form is $y\mapsto\mathfrak a(y,y)\ge0$; its restriction to real functions is a symmetric bilinear
form. Closedness means that $\mathcal V$ is complete for the norm $(\|y\|_{\mathcal H_p}^2+\mathfrak a(y,y))^{1/2}$.

\begin{proposition}\label{B:result:2-5b}
The sesquilinear form $(\mathfrak a,\mathcal V)$ is densely defined, closed, and coercive on $\mathcal H_p$. The inclusion of $\mathcal V$, with its form norm, into $\mathcal H_p$
is compact. The associated self-adjoint operator is
$M_p^{-1}HM_p$.
\end{proposition}
\begin{proof}
For $y\in\mathcal V$, Cauchy--Schwarz and Tonelli give
\begin{equation}
 \|y\|_{\mathcal H_p}^2
 \le\int_0^1p(s)^2\int_s^1q(t)\,\dd t\,\dd s
             \int_0^1|y'|^2/q
 =\mathcal T\,\mathfrak a(y,y).
 \label{B:eq:2-31n}
\end{equation}
Write $\mathcal H_p=L^2((0,1),p^2\,\dd t)$ and
\[
 \mathcal V=\left\{y\in \mathcal H_p:y\in AC_{\mathrm{loc}}(0,1],\ y(1)=0,
                       \int_0^1\frac{|y'|^2}{q}<\infty\right\}.
\]
For $h\in L^2(0,1)$ define
\begin{equation}
 y_h(t)=-\int_t^1\sqrt{q(s)}h(s)\,\dd s,\qquad
 M_py_h=-\mathcal K^*h,\qquad t>0.
 \label{B:eq:2-31o}
\end{equation}
Local boundedness of $q$ away from zero makes this integral finite and $y_h$ locally absolutely continuous, with $y_h'=\sqrt q\,h$ and $y_h(1)=0$. To prove membership in
$\mathcal V$, Cauchy--Schwarz and Tonelli give
\[
 \begin{split}
 \|y_h\|_{\mathcal H_p}^2
 &\le\|h\|_2^2\int_0^1p(t)^2\int_t^1q(s)\,\dd s\,\dd t\\
 &=\|h\|_2^2\int_0^1q(s)\int_0^sp(t)^2\,\dd t\,\dd s<\infty.
 \end{split}
\]
The last integral is the squared Hilbert--Schmidt norm of the triangular factor already constructed. Hence $y_h\in\mathcal V$ and $q^{-1/2}y_h'=h$.

Conversely, for $y\in\mathcal V$, the fundamental theorem on $[t,1]$ and $y(1)=0$ give $y=y_h$ with $h=q^{-1/2}y'$. Thus $D:y\mapsto q^{-1/2}y'$ is a bijective isometry from the
derivative norm on $\mathcal V$ to $L^2(0,1)$. If $(y_j)$ is Cauchy in that norm, $Dy_j\to h$ in $L^2$, and the preceding bound gives $y_j\to y_h$ also in $\mathcal H_p$. The
derivative norm is therefore complete and controls the $\mathcal H_p$ norm. It is equivalent to the sum of these two norms, which proves closedness of the form.

Under the unitary map $y\mapsto py$ from $\mathcal H_p$ to ordinary $L^2$, the inclusion of $\mathcal V$ into $\mathcal H_p$ has kernel, up to its sign,
\[
 p(t)\sqrt{q(s)}\,\mathbf1_{\{t<s\}}.
\]
Its squared integral is the finite quantity above. This inclusion is Hilbert--Schmidt after identifying $\mathcal V$ with $L^2$ through $D$, and in particular is compact. This
proves the compactness conclusion without an unweighted boundary estimate at zero.

The set $C_c^\infty(0,1)$ is contained in $\mathcal V$ and dense in $\mathcal H_p$: truncation to $[1/m,1-1/m]$ converges in the weighted $L^2$ norm, and on every such interval $p$
and $p^{-1}$ are bounded, so ordinary smooth approximation applies. This proves density of the form domain. If $h_j\in C_c^\infty(0,1)$ and $h_j\to Dy$ in $L^2$, then
\[
 \mathfrak a(y_{h_j}-y,y_{h_j}-y)=\|h_j-Dy\|_2^2\to0,
 \qquad
 \|y_{h_j}-y\|_{\mathcal H_p}^2\le\mathcal T\|h_j-Dy\|_2^2.
\]
Consequently $\{y_h:h\in C_c^\infty(0,1)\}$ is a form core. Each of its elements is constant near zero and vanishes near one.

For $g\in L^2$, set $y=M_p^{-1}Rg$. By \eqref{B:eq:2-31j}, $Dy=-\mathcal K g$. If $z\in\mathcal V$ and $h=Dz$, then $pz=-\mathcal K^*h$, whence
\begin{equation}
 \mathfrak a(y,z)=-\langle\mathcal K g,h\rangle
                =\langle g,pz\rangle.
 \label{B:eq:2-31p}
\end{equation}
The associated operator has domain
\[
 \{y\in\mathcal V:\ \exists b\in\mathcal H_p\text{ such that }
 \mathfrak a(y,z)=\langle b,z\rangle_{\mathcal H_p}
                  \text{ for every }z\in\mathcal V\}
\]
and sends $y$ to $b$. For each $b$, equation~\eqref{B:eq:2-31p} with $g=pb$ gives the solution $y=M_p^{-1}RM_pb$. Coercivity implies its uniqueness. The associated operator is
therefore the inverse of $M_p^{-1}RM_p$, namely
$M_p^{-1}HM_p$.
\end{proof}

The left boundary condition can also be read directly from the form. Suppose $y\in\mathcal V$ and $\mathfrak a(y,z)=\langle b,z\rangle_{\mathcal H_p}$ for all $z\in\mathcal V$.
Testing with $z\in C_c^\infty(0,1)$ gives $v'=-p^2b$ in distributions, where $v=q^{-1}y'$. On compact intervals $v\in L^2$ and $p^2b\in L^1$, so $v$ has a locally absolutely
continuous representative. Moreover,
\[
 \int_0^1|v'|\,\dd s=\int_0^1p^2|b|\,\dd s
 \le\|p\|_2\|b\|_{\mathcal H_p}<\infty.
\]
It therefore has a limit $v(0+)$. For a core function $z=y_h$, integration by parts on $[\epsilon,1]$, followed by $\epsilon\downarrow0$, yields
\[
 \mathfrak a(y,z)=-v(0+)\overline{z(0+)}
                    +\int_0^1p^2b\overline z\,\dd s.
\]
Here $z$ is constant near zero, $z(1)=0$, and its derivative has compact support; the integral converges by Cauchy--Schwarz in $\mathcal H_p$. Taking a nonnegative, nonzero
$h\in C_c^\infty(0,1)$ gives $z(0+)=-\int_0^1\sqrt q\,h\ne0$. The weak equation forces $v(0+)=0$.

The approximation by $C_c^\infty(0,1)$ above was in $\mathcal H_p$. When $\int_0^1q<\infty$, evaluation at zero is continuous in energy:
\[
 |y(0+)|^2\le\left(\int_0^1q\,\dd s\right)\mathfrak a(y,y).
\]
Every compactly supported smooth function has zero trace there, whereas the core function just constructed has nonzero trace. Thus $C_c^\infty(0,1)$ cannot be a form core in this
case. For $q^B$ with $-1<a<0$, choosing $0<2\eta<-a$ in \eqref{B:eq:bessel-q-bounds} gives $\int_0^1q^B\le C_{\omega,\eta}^2/(-a-2\eta)<\infty$. The form core
$\{y_h:h\in C_c^\infty(0,1)\}$ applies throughout the range $a>-1$.

\subsubsection{Eigenfunctions and diffusion coordinates}
\label{B:subsection:bessel-spectrum-and-diffusion}

The equations for $y_f$ and $v_f$ form a first-order system with continuous coefficients on every $[\epsilon,1]$. At the regular boundary the condition $y_f(1)=0$ leaves one scalar
initial value, $v_f(1)$. Uniqueness for this system, with initial data at the regular boundary, proves simplicity.

\begin{corollary}\label{B:cor:bessel-simple-spectrum}
The eigenvalues of $H$ are simple and form a sequence $0<\Lambda_1<\Lambda_2<\cdots$ tending to infinity. Every nonzero eigenfunction satisfies
\begin{equation}
 \int_0^1p(s)f(s)\,\dd s\ne0.
 \label{B:eq:2-31c}
\end{equation}
Each eigenspace has a unique real unit vector for which this integral is positive.
\end{corollary}
\begin{proof}
Let $Hf=\Lambda f$. With $y=f/p$ and $v=q^{-1}y'$, the domain identification gives the differential system
\begin{equation}
 \begin{pmatrix}y\\v\end{pmatrix}'
  =\begin{pmatrix}0&q\\-\Lambda p^2&0\end{pmatrix}
                \begin{pmatrix}y\\v\end{pmatrix},\qquad
 y(1)=0,\quad v(0+)=0.
 \label{B:eq:2-31k}
\end{equation}
Initially the equations hold almost everywhere. Local absolute continuity makes $y,v$ continuous; their displayed derivatives are then continuous on every $[\epsilon,1]$, so $y,v$
are $C^1$ there. Write $A(s)$ for the coefficient matrix and $C_\epsilon=\sup_{[\epsilon,1]}\|A(s)\|_{\rm op}$. If two solutions have the same value at $1$, their difference $Z$
satisfies $Z(s)=-\int_s^1A(t)Z(t)\,\dd t$. Iteration gives for every integer $m$
\[
 |Z(s)|\le\sup_{[\epsilon,1]}|Z|
              \frac{(C_\epsilon(1-s))^m}{m!}.
\]
Letting $m\to\infty$ proves $Z=0$ on $[\epsilon,1]$. Thus a solution with $y(1)=v(1)=0$ is zero on all of $(0,1]$. The linear map $f\mapsto v_f(1)$ is scalar valued and injective
on the eigenspace. The eigenspace therefore has dimension at most one. The spectral description in \eqref{B:eq:2-31h} now gives the simple increasing sequence.

Since $pf\in L^1(0,1)$, integration of the second equation gives
\begin{equation}
 v_f(1)=-\Lambda\int_0^1p(s)f(s)\,\dd s.
 \label{B:eq:2-31l}
\end{equation}
Both $v_f(1)$ and $\Lambda$ are nonzero, proving \eqref{B:eq:2-31c}. The kernel of $R$ is real. The real and imaginary parts of a complex eigenfunction are therefore
eigenfunctions, at least one of them nonzero. A real unit eigenfunction has exactly two possible signs, and \eqref{B:eq:2-31c} selects one of them.
\end{proof}

For $p=p^B$ and $q=q^B$, set $s=e^{-x}$ and define
\begin{equation}
 \mathfrak m_W(x)=e^{-x}p^B(e^{-x})^2
           =e^{-(a+1)x+2W(x)/\sqrt\beta},\qquad
 \mathfrak s_W(x)=e^{-x}q^B(e^{-x})
           =e^{ax-2W(x)/\sqrt\beta}.
 \label{B:eq:2-31q}
\end{equation}
These are the speed and scale densities on the half-line. The map $(\mathscr U_pf)(x)=f(e^{-x})/p^B(e^{-x})$ satisfies
\[
 \int_0^\infty|\mathscr U_pf(x)|^2\mathfrak m_W(x)\,\dd x
 =\int_0^1|f(s)|^2\,\dd s,
\]
and is onto, with inverse $h\mapsto[s\mapsto p^B(s)h(-\log s)]$. For a compactly supported $h$, substitute $t=e^{-y}$ and $r=e^{-u}$ in the Green kernel \eqref{B:eq:2-31}. The
identities $p^B(t)^2\,\dd t=\mathfrak m_W(y)\,\dd y$ after reversing the limits and $q^B(r)\,\dd r=\mathfrak s_W(u)\,\dd u$ give
\begin{equation}
 (\mathscr U_pR\mathscr U_p^{-1}h)(x)
 =\int_0^\infty\left(\int_0^{x\wedge y}\mathfrak s_W(u)\,\dd u\right)
                   h(y)\mathfrak m_W(y)\,\dd y.
 \label{B:eq:2-31r}
\end{equation}
The identity extends to the weighted $L^2$ space by density and the boundedness of $R$.

If $h=\mathscr U_pf$, the chain rule for the locally absolutely continuous function $y_f$ yields
\[
 h'(x)=-e^{-x}y_f'(e^{-x}),\qquad
 \mathfrak s_W(x)^{-1}h'(x)=-v_f(e^{-x}).
\]
Differentiating the second expression gives $(\mathfrak s_W^{-1}h')'(x)=e^{-x}v_f'(e^{-x})$. It follows that
\begin{equation}
 \mathscr U_pH\mathscr U_p^{-1}h
   =-\mathfrak m_W^{-1}(\mathfrak s_W^{-1}h')',\qquad
 h(0)=0,\quad \lim_{x\to\infty}\mathfrak s_W(x)^{-1}h'(x)=0.
 \label{B:eq:bessel-diffusion-action}
\end{equation}
In particular, the complete domain on $L^2((0,\infty),\mathfrak m_W\,\dd x)$ consists of functions $h$ for which, writing $\ell_h=\mathfrak s_W^{-1}h'$, one has
\[
 \begin{gathered}
 h,\ell_h\in\operatorname{AC}([0,L])\quad(L<\infty),\qquad
 \mathfrak m_W^{-1}\ell_h'\in L^2(\mathfrak m_W\,\dd x),\\
 h(0)=0,\qquad\ell_h(\infty)=0.
 \end{gathered}
\]
Indeed, the norm of the differential action transforms as
\[
 \int_0^\infty|\mathfrak m_W^{-1}\ell_h'|^2\mathfrak m_W\,\dd x
 =\int_0^1|p^{-1}v_f'|^2\,\dd s.
\]
The same substitution in the energy gives
\[
 \int_0^\infty\mathfrak s_W^{-1}|h'|^2\,\dd x
 =\int_0^1q^{-1}|y_f'|^2\,\dd s.
\]
Thus the domain and form are transported by the specified unitary. With the choice $b=-W$ in \cite{ramirez2009diffusion}, its speed density, scale density, and inverse kernel equal
\eqref{B:eq:2-31q}--\eqref{B:eq:2-31r}. The process obtained by this sign change is again Brownian motion.

For zero noise, $p(s)=s^{a/2}$ and $q(s)=s^{-a-1}$. Then $P(t)=t^{a+1}/(a+1)$ and $\operatorname{Tr}R=1/(a+1)$. Applying \eqref{B:eq:2-31i} to $g=p$ gives
\begin{equation}
 Rp=\frac{s^{a/2}(1-s)}{a+1},\qquad
 v_{Rp}(s)=-\frac{s^{a+1}}{a+1},\qquad H(Rp)=p.
 \label{B:eq:2-31s}
\end{equation}
For a smooth $f$ in this domain,
\[
 v_f=s^{a/2+1}f'-\frac a2s^{a/2}f,\qquad
 -p^{-1}v_f'=-sf''-f'+\frac{a^2}{4s}f.
\]
When $-1<a<0$, the domain element $Rp$ in \eqref{B:eq:2-31s} diverges at zero. The condition there is the limit of $v_f$ in \eqref{B:eq:bessel-domain-conditions}.

\subsubsection{Dependence on the Brownian path}
\label{B:subsection:bessel-path-measurability}

The pathwise construction must also define random operators in the norm topologies used later. For fixed $\beta>0$ and $a>-1$, this follows by working on the set $\mathscr W$ of
continuous real paths $w$ on $[0,\infty)$ with $w(0)=0$ and $w(u)/u\to0$ at infinity. The following measurability statement will also be used for the spectral projections.

\begin{lemma}\label{B:lem:bessel-path-measurability}
The maps $w\mapsto p^B\in L^2(0,1)$ and $w\mapsto\mathcal K_{\beta,a}\in\mathfrak S_2$ are Borel on $\mathscr W$. The map $w\mapsto R\in\mathfrak S_1$ is Borel there.
\end{lemma}
\begin{proof}
On the space of continuous paths with $w(0)=0$, put $b_m(w)=m^{-1}\sup_{u\in[m,m+1]}|w(u)|$ for $m\ge1$. The ratios $u/m$ lie in $[1,2]$ on this interval, so the defining
sublinearity condition is equivalent to $b_m(w)\to0$. More explicitly,
\[
 \mathscr W=\bigcap_{r=1}^\infty\bigcup_{M=1}^\infty
              \bigcap_{m\ge M}\{w:b_m(w)\le r^{-1}\}.
\]
Each $b_m$ is continuous for local uniform convergence, proving that $\mathscr W$ is Borel. The bound $|w(u)|\le C_{w,\eta}+\sqrt\beta\eta u$ gives
\eqref{B:eq:bessel-p-bounds}--\eqref{B:eq:2-31e} for each of these paths.

For $m\ge2$ let $p_m=\mathbf1_{[1/m,1]}p^B$ and restrict the factor kernel to $s\ge1/m$. Denote the resulting operator by $\mathcal K_m$. Both kernels then depend on $w$ only on
$[0,\log m]$. Put $d=\|w-\widetilde w\|_{C([0,\log m])}$, and assume that both paths have norm at most $A$ there. The identity $e^x-e^y=(x-y)\int_0^1e^{y+t(x-y)}\,\dd t$ gives
\[
 \begin{aligned}
 |p_m(w;s)-p_m(\widetilde w;s)|
 &\le\beta^{-1/2}e^{A/\sqrt\beta}s^{a/2}d\,
                    \mathbf1_{\{s\ge1/m\}},\\
 |K_m(w;s,t)-K_m(\widetilde w;s,t)|
 &\le2\beta^{-1/2}e^{2A/\sqrt\beta}
       t^{-1/2}(s/t)^{a/2}d\,\mathbf1_{\{1/m\le s<t\}}.
 \end{aligned}
\]
The squared deterministic factors have finite integrals on the indicated compact triangle. Integrating these two inequalities yields
\[
 \|p_m(w)-p_m(\widetilde w)\|_2
 +\|\mathcal K_m(w)-\mathcal K_m(\widetilde w)\|_{\mathfrak S_2}
 \le C_{m,A,a,\beta}\|w-\widetilde w\|_{C([0,\log m])}.
\]
Thus these truncated maps are continuous. For every $w\in\mathscr W$, absolute continuity of the defining integrals gives
\[
 \|p_m-p^B\|_2\to0,\qquad
 \|\mathcal K_m-\mathcal K_{\beta,a}\|_{\mathfrak S_2}^2
 =\int_0^1q^B(t)\int_0^{t\wedge m^{-1}}p^B(s)^2\,\dd s\,\dd t\to0.
\]
The integrand in the last limit is bounded by $q^B(t)P(t)$, which is integrable by \eqref{B:eq:2-31e}. Hence these norm limits are Borel. Finally, the product inequality for
Hilbert--Schmidt operators gives
\[
 \|\mathcal K_m^*\mathcal K_m-\mathcal K_{\beta,a}^*\mathcal K_{\beta,a}\|_{\mathfrak S_1}
 \le(\|\mathcal K_m\|_{\mathfrak S_2}
       +\|\mathcal K_{\beta,a}\|_{\mathfrak S_2})
                 \|\mathcal K_m-\mathcal K_{\beta,a}\|_{\mathfrak S_2}\to0.
\]
Each approximating map is continuous in trace norm. The limit map is therefore Borel measurable in trace norm. Extending these maps by zero off $\mathscr W$ gives Borel maps on the
full continuous path space.
\end{proof}

\subsubsection{Restriction to a subinterval}
\label{B:subsection:retained-cells-errors-and-weights}

The scale density cannot be integrated through zero in the parameter range used for the approximation. The comparison is therefore first made on a selected interval, and the
omitted part of the positive inverse is estimated separately.

Assume $a>0$ for the remaining approximation formulas. Choose integers $1\le k_-<k_+\le n$, let $\mathscr I_n^{\mathrm{ret}}=\{k_-+1,\ldots,k_+\}$, and set $l=z_{n,k_-}$,
$u=z_{n,k_+}$. Let $\Pi_{l,u}$ be multiplication by $\mathbf1_{(l,u)}$. All coefficient masses below have finite positive integration limits. In fact, taking $0<2\eta<a$ in
\eqref{B:eq:bessel-q-bounds} gives
\[
 \int_0^{z_{n,1}}q^B(r)\,\dd r
 \ge C_{\omega,\eta}^{-2}\int_0^{z_{n,1}}r^{-a-1+2\eta}\,\dd r=\infty.
\]

For $k\in\mathscr I_n^{\mathrm{ret}}$, define
\begin{equation}
 Q_{n,k}^B=\int_{I_{n,k}}q^B(r)\,\dd r,\qquad
 p_{n,k}^B=p^B(z_{n,k}).
 \label{B:eq:2-32}
\end{equation}
The function $q^B$ is positive and continuous on $[l,u]$, so each $Q_{n,k}^B$ is finite and strictly positive. On the selected cell pairs define
\begin{equation}
 K_n^{B,\mathrm{step},0}(s,t)
   =p_{n,i}^B\sqrt{Q_{n,j}^B/\Delta_{n,j}}\,
                         \mathbf1_{\{i\le j\}},
 \quad (s,t)\in I_{n,i}\times I_{n,j},
 \label{B:eq:bessel-step-factor}
\end{equation}
and set it equal to zero elsewhere. Taking the product of the adjoint with this operator gives the kernel
\begin{equation}
 R_n^{B,\mathrm{step},0}(s,t)
   =p_{n,i}^Bp_{n,j}^B\sum_{\ell=i\vee j}^{k_+}Q_{n,\ell}^B.
 \label{B:eq:2-34}
\end{equation}
This follows by integrating the product of the two factor kernels on each cell in the second kernel variable, whose length cancels the denominator $\Delta_{n,\ell}$ in
\eqref{B:eq:bessel-step-factor}.

For a general factor $\mathcal K$ as in Section~\ref{B:subsection:the-stochastic-bessel-realization}, let $\mathcal K^0=\Pi_{l,u}\mathcal K\Pi_{l,u}$ and
$R^0=(\mathcal K^0)^*\mathcal K^0$. These operators act on $L^2(l,u)$ and are extended by zero on its orthogonal complement. Compression of the full inverse includes the values of
the triangular factor on $(u,1)$. They give the correction calculated below.

\begin{lemma}\label{B:lem:bessel-restricted-factor}
The restricted inverse has kernel
\begin{equation}
 R^0(s,t)=\mathbf1_{(l,u)}(s)\mathbf1_{(l,u)}(t)
                  p(s)p(t)\int_{s\vee t}^{u}q(r)\,\dd r.
 \label{B:eq:bessel-restricted-green-kernel}
\end{equation}
Moreover,
\begin{equation}
 \Pi_{l,u}R\Pi_{l,u}-R^0
   =\left(\int_u^1q(r)\,\dd r\right)
                  |\Pi_{l,u}p\rangle\langle\Pi_{l,u}p|.
 \label{B:eq:bessel-compression-difference}
\end{equation}
\end{lemma}
\begin{proof}
Let $P$ be multiplication by $\mathbf1_{(l,u)}$, and let the factor kernel $K(s,t)=p(s)\sqrt{q(t)}\mathbf1_{\{s<t\}}$ act by integration in its first variable. For $f\in L^2$
supported in $(l,u)$,
\[
 (P\mathcal KPf)(t)
 =\mathbf1_{(l,u)}(t)\sqrt{q(t)}\int_l^t p(s)f(s)\,\dd s.
\]
The triangular factor is Hilbert--Schmidt. Its adjoint product can therefore be computed by Fubini in the scalar pairings, with absolute integrability following from
Cauchy--Schwarz. Its adjoint product has kernel
\[
 R^0(s,t)=\mathbf1_{(l,u)^2}(s,t)\,p(s)p(t)
                     \int_{s\vee t}^u q(v)\,\dd v.
\]
The compressed full inverse instead has kernel
\[
 (PRP)(s,t)=\mathbf1_{(l,u)^2}(s,t)\,p(s)p(t)
                     \int_{s\vee t}^1 q(v)\,\dd v.
\]
Subtracting their integration intervals yields
\[
 PRP-R^0=\left(\int_u^1q(v)\,\dd v\right)|Pp\rangle\langle Pp|.
\]
This is a positive operator of rank at most one. If $c_u=\int_u^1q$ and $v=Pp$, its only possible nonzero eigenvalue is $c_u\|v\|_2^2$, since $(|v\rangle\langle v|)v=\|v\|_2^2v$
and the operator vanishes on $v^\perp$. Consequently
\[
 \|PRP-R^0\|_{\mathfrak S_2}
       =c_u\|Pp\|_2^2
       =\left(\int_u^1q\right)\left(\int_l^u p^2\right).
\]
If either factor vanishes, the same formula gives zero. All identities are equalities of bounded operators because they have been obtained from equal square-integrable kernels. The
value of either kernel on its diagonal is immaterial for Lebesgue integration.
\end{proof}

The boundary conditions of these two operators follow from their kernels. For $g\in L^2(l,u)$ put $F_{l,g}(s)=\int_l^s p(t)g(t)\,\dd t$. If $f=R^0g$, then
\[
 \frac{f(s)}{p(s)}=\int_s^u q(t)F_{l,g}(t)\,\dd t,\qquad
 v_f(s)=-F_{l,g}(s),\qquad y_f(u)=0,\quad v_f(l)=0.
\]
The derivative is $v_f'=-pg$. Conversely, take $f\in L^2(l,u)$ with $y_f,v_f\in\operatorname{AC}([l,u])$, $g=-p^{-1}v_f'\in L^2(l,u)$, $v_f(l)=0$, and $y_f(u)=0$. Integration gives
$v_f(s)=-\int_l^s pg=-F_{l,g}(s)$ and $y_f(s)=-\int_s^u qv_f=\int_s^u qF_{l,g}$. Thus $f=R^0g$. This proves that $R^0$ is the inverse on $L^2(l,u)$ with the stated differential
action and these two boundary conditions.

Put $Q_u=\int_u^1q(t)\,\dd t$. For $f=\Pi_{l,u}R\Pi_{l,u}g$, the rank-one term adds $Q_uF_{l,g}(u)$ to $y_f(s)$ and leaves its derivative unchanged. Consequently
\[
 y_f(s)=\int_s^u q(t)F_{l,g}(t)\,\dd t+Q_uF_{l,g}(u),\qquad
 v_f(l)=0,\quad y_f(u)+Q_uv_f(u)=0.
\]
Conversely, under the same regularity conditions and $v_f(l)=0$, the boundary equation $y_f(u)+Q_uv_f(u)=0$ gives $y_f(u)=Q_uF_{l,g}(u)$. Integration of $y_f'=qv_f$ then gives
exactly the displayed formula. Thus the compression is the inverse with this Robin condition at $u$. For $u=1$, $Q_u=0$ and both inverses agree.

For the finite factor, the same multiplication gives
\begin{equation}
 R_n^{\Delta,0}(s,t)=p_{n,i}p_{n,j}
                       \sum_{\ell=i\vee j}^{k_+}q_{n,\ell},
       \qquad i,j\in\mathscr I_n^{\mathrm{ret}},
 \label{B:eq:bessel-finite-restriction}
\end{equation}
with zero kernel elsewhere. Its difference from the compression of $\mathcal R_n^\Delta$ is the positive operator of rank at most one with kernel
\[
 \mathbf1_{(l,u)}(s)\mathbf1_{(l,u)}(t)
          p_{n,i}p_{n,j}\sum_{\ell>k_+}q_{n,\ell}.
\]
Thus the finite and continuum restrictions both truncate the factor before taking its squared modulus. In \eqref{B:eq:bessel-restricted-green-kernel}, the choice $p=p^B,q=q^B$ is
denoted by $R_{\beta,a}^{0}$ below.

The restricted inverse also has an exact matrix formula. Let $I=\{k_-+1,\ldots,k_+\}$, let $E_I:\mathbb C^{|I|}\to\mathbb C^n$ be the coordinate inclusion, with the indices of $I$
in increasing order, and put $B_I=E_I^*B_nE_I$ and $\widetilde{\mathscr J}_{n,I}= \widetilde{\mathscr J}_nE_I$. Formula \eqref{B:eq:2-3} uses only coefficients between its two
indices, so
\[
 E_I^*B_n^{-1}E_I=B_I^{-1},\qquad
 \mathcal R_n^{\Delta,0}
 =\frac{\alpha_n}{n^2}\widetilde{\mathscr J}_{n,I}
          (B_IB_I^{\mathsf T})^{-1}\widetilde{\mathscr J}_{n,I}^*.
\]
The size in the prefactor remains $n$. The principal submatrix of the full positive matrix has an additional first diagonal term:
\[
 E_I^*B_nB_n^{\mathsf T}E_I
 =B_IB_I^{\mathsf T}+Y_{k_-}^2e_1e_1^{\mathsf T}.
\]
Indeed its first diagonal entry is $X_{k_-+1}^2+Y_{k_-}^2$; all other diagonal and adjacent entries agree. The restriction in \eqref{B:eq:bessel-finite-restriction} uses
$B_IB_I^{\mathsf T}$. To compute the inverse of the principal submatrix instead, write $G_I=(B_IB_I^{\mathsf T})^{-1}$, $v_I=G_Ie_1$, and $t=Y_{k_-}^2$. Then
\[
 (E_I^*B_nB_n^{\mathsf T}E_I)^{-1}
 =G_I-\frac{t}{1+t\langle G_Ie_1,e_1\rangle}\,v_Iv_I^*.
\]
The denominator is positive, since $G_I$ is positive definite. To verify the formula, put $H_I=B_IB_I^{\mathsf T}$ and $c=t/(1+t\langle G_Ie_1,e_1\rangle)$. Since $H_Iv_I=e_1$,
\[
 (H_I+te_1e_1^*)(G_I-cv_Iv_I^*)
 =I+\bigl[t-c-ct\langle G_Ie_1,e_1\rangle\bigr]e_1v_I^*=I.
\]
Thus this inverse is smaller than $G_I$ by a positive rank-one term. The compression of the full inverse, on the other hand, exceeds the inverse obtained from the restricted factor
by \eqref{B:eq:bessel-compression-difference} and its finite counterpart. These two boundary corrections arise at opposite ends of the selected interval. Their embeddings both
retain the original prefactor $\alpha_n/n^2$.

\subsubsection{Logarithmic coefficient errors}
\label{B:subsection:logarithmic-coefficient-errors-on-the-retained-cells}
\label{B:subsection:cumulative-and-local-coefficient-errors}

Multiplying the factors outside the tail sum by a positive scalar and dividing each coefficient in the sum by the square of that scalar leaves the positive inverse kernel
unchanged. Removing that freedom before taking logarithms leaves the two error sequences that enter the stability estimate.

Choose a real scalar $c_n$ and set, on $\mathscr I_n^{\mathrm{ret}}$, $\widehat p_{n,k}=e^{-c_n}p_{n,k}$ and $\widehat q_{n,k}=e^{2c_n}q_{n,k}$. By \eqref{B:eq:2-20} this change
preserves the finite kernels. Define the errors through
\begin{equation}
 \widehat p_{n,k}=p_{n,k}^Be^{\delta_{n,k}},\qquad
 \widehat q_{n,k}=Q_{n,k}^B e^{-2\delta_{n,k}+\rho_{n,k}}.
 \label{B:eq:2-35}
\end{equation}
Positivity of all four coefficients makes the logarithms well-defined.

\begin{proposition}\label{B:result:2-6}
For selected indices $i,j$,
\begin{equation}
 R_n^{\Delta,0}(i,j)
 =p_{n,i}^Bp_{n,j}^Be^{\delta_{n,i}+\delta_{n,j}}
   \sum_{\ell=i\vee j}^{k_+}Q_{n,\ell}^B
                         e^{-2\delta_{n,\ell}+\rho_{n,\ell}}.
 \label{B:eq:2-37}
\end{equation}
\end{proposition}
\begin{proof}
For $i,j,\ell\in\mathscr I_n^{\rm ret}$ with $\ell\ge i\vee j$,
\[
 \begin{aligned}
 p_{n,i}p_{n,j}q_{n,\ell}
 &=e^{c_n}\widehat p_{n,i}\,e^{c_n}\widehat p_{n,j}
                      e^{-2c_n}\widehat q_{n,\ell}\\
 &=p_{n,i}^Bp_{n,j}^BQ_{n,\ell}^B
       e^{\delta_{n,i}+\delta_{n,j}-2\delta_{n,\ell}+\rho_{n,\ell}}.
 \end{aligned}
\]
Summing this equality over $i\vee j\le\ell\le k_+$ in \eqref{B:eq:bessel-finite-restriction} gives \eqref{B:eq:2-37}. Both kernels are zero outside the square of the restricted
kernel, so the identity also holds as an equality of integral operators.
\end{proof}

The definitions yield formulas separating the cumulative ratios from the individual coefficient:
\begin{align}
 \delta_{n,k}
 &=\sum_{r=k}^{n-1}\log\frac{Y_r}{X_r}
   -\mathcal Y(z_{n,k})-\frac a2\log z_{n,k}
   +\frac12\log\frac{\alpha_n}{n\Delta_{n,k}}-c_n,
       \label{B:eq:delta-explicit-coefficients}\\
 \rho_{n,k}
 &=\log\frac{\alpha_n z_{n,k}}
                   {n^2\Delta_{n,k}^2X_k^2}
    -\log\frac{Q_{n,k}^B}{\Delta_{n,k}q^B(z_{n,k})}.
       \label{B:eq:2-38}
\end{align}
The first line follows from \eqref{B:eq:2-5}, \eqref{B:eq:2-13}, and \eqref{B:eq:intro-brownian-coefficients}. For the second, \eqref{B:eq:2-35} gives
\[
 e^{\rho_{n,k}}
 =\frac{q_{n,k}p_{n,k}^2}{Q_{n,k}^B(p_{n,k}^B)^2}
 =\frac{\alpha_n}{n^2\Delta_{n,k}X_k^2Q_{n,k}^B(p_{n,k}^B)^2}.
\]
Use $(p_{n,k}^B)^2q^B(z_{n,k})=z_{n,k}^{-1}$ and take logarithms. The product $P_{n,k}$ and the scalar $c_n$ have canceled from this expression. In the couplings below the
coefficient row and $c_n$ are $\mathcal N_n\vee\mathcal U$-measurable. The mass $Q_{n,k}^B$ also uses the Brownian path between the two sampled positions. It belongs to the error
being estimated; the coefficient samplers use the observation field specified for their row.

For $r\in I_{n,k}$, the definition of $q^B$ gives
\[
 \left|\log\frac{q^B(r)}{q^B(z_{n,k})}\right|
 \le(a+1)h_{n,k}
           +2\sup_{s\in I_{n,k}}|\mathcal Y(s)-\mathcal Y(z_{n,k})|
 =:L_{n,k}.
\]
The ratio lies between $e^{-L_{n,k}}$ and $e^{L_{n,k}}$. Averaging it over the cell and taking its logarithm proves
\begin{equation}
 \left|\log\frac{Q_{n,k}^B}{\Delta_{n,k}q^B(z_{n,k})}\right|
 \le(a+1)h_{n,k}
       +2\sup_{s\in I_{n,k}}|\mathcal Y(s)-\mathcal Y(z_{n,k})|.
 \label{B:eq:bessel-cell-mass-log-bound}
\end{equation}
This is a bound for the second term in \eqref{B:eq:2-38} on each cell.

The deterministic coefficient in \eqref{B:eq:2-38} also has an exact expression in terms of the profile. On the Laguerre partition, $\alpha_n=1$ and $n\Delta_{n,k}=1$, so
\[
 \rho_{n,k}
 =-2\log\frac{X_k}{\sqrt{k/n}}
   -\log\frac{Q_{n,k}^B}{\Delta_{n,k}q^B(z_{n,k})}.
\]
For the general potential, \eqref{B:eq:1-16} implies $\phi(k/n)^2\theta_V'(k/n)^2=\alpha_nz_{n,k}$. Consequently
\[
 \rho_{n,k}
 =-2\log\frac{X_k}{\phi(k/n)}
   +2\log\frac{\theta_V'(k/n)}{n\Delta_{n,k}}
   -\log\frac{Q_{n,k}^B}{\Delta_{n,k}q^B(z_{n,k})}.
\]
If $c\le\theta_V'$ and $\|\theta_V''\|_\infty\le C$, then $n\Delta_{n,k}=n\int_{(k-1)/n}^{k/n}\theta_V'(t)\,\dd t\ge c$ and
\[
 |n\Delta_{n,k}-\theta_V'(k/n)|\le C/n,
 \qquad
 \left|\log\frac{\theta_V'(k/n)}{n\Delta_{n,k}}\right|
 \le\frac{C}{cn}.
\]
The second inequality uses $|\log x-\log y|\le|x-y|/c$ for $x,y\ge c$. Together with the bounded sum of weights proved below, it shows that this deterministic contribution to the
weighted squared error is $O(n^{-2})$. The remaining two contributions are the local coefficient error and the logarithmic error in approximating the integral $Q_{n,k}^B$.

To see the different role of the product error, define for $1\le r<n$
\[
 \xi_{n,r}=\log\frac{Y_r}{X_r}+\frac a2h_{n,r}^+
  -\bigl[\mathcal Y(z_{n,r})-\mathcal Y(z_{n,r+1})\bigr]
  +\frac12\log\frac{\Delta_{n,r+1}}{\Delta_{n,r}}.
\]
Subtracting \eqref{B:eq:delta-explicit-coefficients} at successive indices and using $\log(z_{n,r+1}/z_{n,r})=h_{n,r}^+$ gives $\delta_{n,r}-\delta_{n,r+1}=\xi_{n,r}$ wherever both
errors are considered. Telescoping therefore gives, for retained $i\le\ell$,
\[
 \delta_{n,i}-\delta_{n,\ell}
      =\sum_{r=i}^{\ell-1}\xi_{n,r},\qquad
 \delta_{n,i}=\delta_{n,k_+}+\sum_{r=i}^{k_+-1}\xi_{n,r}.
\]
The Brownian increment in $\xi_{n,r}$ has logarithmic length $h_{n,r}^+=h_{n,r+1}^-$; the mass in $\rho_{n,r}$ uses $h_{n,r}^-$. They belong to successive cells. The exponent of
each summand in \eqref{B:eq:2-37} is exactly
\[
 \delta_{n,i}+\delta_{n,j}-2\delta_{n,\ell}+\rho_{n,\ell}
 =\sum_{r=i}^{\ell-1}\xi_{n,r}
   +\sum_{r=j}^{\ell-1}\xi_{n,r}+\rho_{n,\ell},
 \qquad \ell\ge i\vee j.
\]
Empty sums are zero. This identity removes the value at the upper index $\delta_{n,k_+}$ before any estimate. The scalar normalization used by the coupling is still needed for the
separate maximum and weighted bounds on the individual $\delta_{n,k}$ in the stability theorem.

Small errors at individual indices need not give small cumulative errors. Fix $1<\lambda<\min\{a+1,2\}$. On the uniform partition let $k_+/n\to1$, $\delta_{n,k_+}=0$, and
$\xi_{n,r}=n^{-1}$ on the selected indices. Then $\delta_{n,k}=(k_+-k)/n$. For all sufficiently large $n$ and $n/4\le k\le n/2$, one has $\delta_{n,k}\ge1/3$ and
$h_{n,k}z_{n,k}^{\lambda}\ge4^{-\lambda}/n$. Hence
\[
 \sum_{k\in\mathscr I_n^{\mathrm{ret}}}
       h_{n,k}z_{n,k}^{\lambda}|\delta_{n,k}|^2
 \ge\frac{4^{-\lambda}}{9n}
       \#\{k:n/4\le k\le n/2\}\ge c_\lambda>0.
\]
Here the selected interval contains these indices, as it does for the logarithmic cutoffs used below. Thus bounds for individual logarithmic ratios must be supplemented by
estimates for their partial sums. The coupling below controls the cumulative errors appearing in \eqref{B:eq:delta-explicit-coefficients}.

\subsubsection{Weights in the coefficient estimates}
\label{B:subsection:the-weighted-error-norm}

The kernel determines the weights for the coefficient errors. In logarithmic time these weights discretize an exponentially decaying measure.

For $1<\lambda<\min\{2,a+1\}$, put
\begin{equation}
 w_{n,k}^{(\lambda)}=h_{n,k}z_{n,k}^{\lambda},\qquad
 k\in\mathscr I_n^{\mathrm{ret}}.
 \label{B:eq:2-39}
\end{equation}
These weights discretize the measure $e^{-\lambda u}\,\dd u$ in the logarithmic variable. Direct integration gives
\begin{equation}
 \int_{u_{n,k}}^{u_{n,k-1}}e^{-\lambda u}\,\dd u
 =\frac{z_{n,k}^{\lambda}}{\lambda}(1-e^{-\lambda h_{n,k}}).
 \label{B:eq:spectral-weight-integral}
\end{equation}
For $x\ge0$,
\[
 e^{-x}\le\frac{1-e^{-x}}x=\int_0^1e^{-tx}\,\dd t\le1,
 \qquad 0\le1-\frac{1-e^{-x}}x\le x/2,
\]
where the quotient at zero is one. Consequently, for $h_{n,k}\le h_*$,
\[
 e^{-\lambda h_*}w_{n,k}^{(\lambda)}
 \le\int_{u_{n,k}}^{u_{n,k-1}}e^{-\lambda u}\,\dd u
 \le w_{n,k}^{(\lambda)},
\]
and the relative error is at most $\lambda h_{n,k}/2$.

For either model, write $z_{n,k}=\theta(k/n)$, where $\theta=\operatorname{id}$ or $\theta_V$. If $0<c\le\theta'\le C$, then, for $k\ge2$,
\[
 \frac cn\le\Delta_{n,k}\le\frac Cn,\qquad
 \frac{ck}{n}\le z_{n,k}\le\frac{Ck}{n},\qquad
 \frac{\Delta_{n,k}}{z_{n,k}}
 \le h_{n,k}=\int_{z_{n,k-1}}^{z_{n,k}}\frac{\dd s}{s}
 \le\frac{\Delta_{n,k}}{z_{n,k-1}}.
\]
For $k\ge2$ these inequalities give
\[
 \frac{c}{Ck}\le h_{n,k}\le\frac{2C}{ck},\qquad
 c_\lambda n^{-\lambda}k^{\lambda-1}
 \le w_{n,k}^{(\lambda)}\le C_\lambda n^{-\lambda}k^{\lambda-1}.
\]
They hold also on the first selected cell, whose index exceeds one. We may take $h_*=C/c$. Since the logarithmic intervals have disjoint interiors,
\[
 \sum_{k\in\mathscr I_n^{\mathrm{ret}}}w_{n,k}^{(\lambda)}
 \le e^{\lambda h_*}\int_0^\infty e^{-\lambda u}\,\dd u
 =e^{\lambda h_*}/\lambda.
\]
Thus each later weighted coefficient bound uses the physical cells and their actual logarithmic lengths.

The model estimates used in the stability theorem take the form
\begin{equation}
 \sum_{k\in\mathscr I_n^{\mathrm{ret}}}
       w_{n,k}^{(\lambda)}
           (|\delta_{n,k}|^2+|\rho_{n,k}|^2)\le C_\omega/n.
 \label{B:eq:2-40}
\end{equation}
We prove this estimate separately for the two coefficient laws. Proposition~\ref{B:prop:weighted-kernel-rows} bounds the squared row norm by $Cz^{\lambda-1}$ when
$2\eta<a+1-\lambda$ and $8\eta<2-\lambda$. The latter inequality also bounds the sum in Proposition~\ref{B:prop:kernel-coefficient-variations}. These conditions hold for the choice
\eqref{B:eq:stability-exponents}, completing the connection between \eqref{B:eq:2-40} and the deterministic inverse estimate.

\subsection{Weighted stability and truncation of the inverse at the hard edge}
\label{sec:weighted-hard-edge-stability-and-truncation-estimates}
\label{B:sec:weighted-stability-and-truncation-of-the-inverse}

The representations in Section~\ref{B:sec:inverse-representations-and-the-stochastic-bessel-operator} reduce the comparison of inverses to estimates for their kernel coefficients.
On a selected interval, the comparison is controlled by the logarithmic errors \eqref{B:eq:2-35}. To remove the restriction, we estimate the difference between the full inverse and
the positive product of the restricted triangular factor. The constants in these deterministic estimates are uniform in the partition and the interval of restriction. When $p=p^B$
and $q=q^B$, Lemma~\ref{B:lem:bessel-brownian-growth} bounds these constants by finite random variables.

\subsubsection{Sums on logarithmic partitions}
\label{B:subsection:logarithmic-meshes-and-power-sums}

The coefficient errors are accumulated in logarithmic time, whereas the Hilbert--Schmidt norm uses physical cell lengths. Comparing these lengths converts the weighted coefficient
bounds into kernel estimates.

Relabel the selected cells of Section~\ref{B:subsection:retained-cells-errors-and-weights} consecutively and suppress the matrix size. Thus $0<z_0<z_1<\cdots<z_M\le1$. We keep the
notation $I_i,\Delta_i,h_i$ of \eqref{B:eq:cell-embedding-definition} and \eqref{B:eq:spectral-cell-conventions} for their intervals, lengths, and logarithmic lengths. Fix $h_*>0$
and assume $h_i\le h_*$. The identity $z_{i-1}=z_i e^{-h_i}$ gives
\begin{equation}
 e^{-h_*}z_i h_i\le\Delta_i=z_i(1-e^{-h_i})\le z_i h_i.
 \label{B:eq:logarithmic-cell-length}
\end{equation}
Indeed, $1-e^{-h}=\int_0^h e^{-t}\,\dd t$ lies between $h e^{-h_*}$ and $h$ for $0\le h\le h_*$. In this subsection $R(i,j)$ denotes a kernel value on $I_i\times I_j$, with zero
extension outside the square of the restricted kernel. Its norm is
\begin{equation}
 \|R\|_{\mathfrak S_2}^2
       =\sum_{i,j=1}^M\Delta_i\Delta_j|R(i,j)|^2,
 \label{B:eq:step-kernel-HS-norm}
\end{equation}
as follows by integrating on the cell pairs in \eqref{B:eq:integral-kernel-convention}.

\begin{lemma}\label{B:lem:logarithmic-power-sums}
For $\gamma>0$, uniformly in $i$,
\begin{align}
 \sum_{j\le i}h_jz_j^\gamma
     &\le\frac{e^{\gamma h_*}}\gamma z_i^\gamma,
       \label{B:eq:logarithmic-positive-sum}\\
 \sum_{j\ge i}h_jz_j^{-\gamma}
     &\le\frac{e^{\gamma h_*}}\gamma z_i^{-\gamma}.
       \label{B:eq:logarithmic-negative-sum}
\end{align}
For every $\theta>0$, writing $h_{\max}=\max_jh_j$, one also has
\begin{equation}
 \sum_{j=1}^M h_j^2z_j^\theta
        \le\frac{e^{\theta h_*}}\theta h_{\max}z_M^\theta.
 \label{B:eq:logarithmic-squared-length-sum}
\end{equation}
\end{lemma}
\begin{proof}
For $\gamma>0$, each cell satisfies
\[
 z_j^\gamma-z_{j-1}^\gamma
 =z_j^\gamma(1-e^{-\gamma h_j})
 =z_j^\gamma\int_0^{\gamma h_j}e^{-t}\,\dd t
 \ge\gamma e^{-\gamma h_*}h_jz_j^\gamma.
\]
Summing the differences gives
\begin{align*}
 \sum_{j=1}^i h_jz_j^\gamma
 &\le\frac{e^{\gamma h_*}}\gamma
            \sum_{j=1}^i(z_j^\gamma-z_{j-1}^\gamma)\\
 &=\frac{e^{\gamma h_*}}\gamma(z_i^\gamma-z_0^\gamma)
 \le\frac{e^{\gamma h_*}}\gamma z_i^\gamma.
\end{align*}
For the negative powers,
\[
 z_{j-1}^{-\gamma}-z_j^{-\gamma}
 =z_j^{-\gamma}(e^{\gamma h_j}-1)
 =z_j^{-\gamma}\int_0^{\gamma h_j}e^t\,\dd t
 \ge\gamma h_jz_j^{-\gamma}.
\]
The sum from $i$ to $M$ therefore obeys
\begin{align*}
 \sum_{j=i}^M h_jz_j^{-\gamma}
 &\le\frac1\gamma
           \sum_{j=i}^M(z_{j-1}^{-\gamma}-z_j^{-\gamma})
 =\frac{z_{i-1}^{-\gamma}-z_M^{-\gamma}}\gamma\\
 &\le\frac{z_{i-1}^{-\gamma}}\gamma
 =\frac{e^{\gamma h_i}}\gamma z_i^{-\gamma}
 \le\frac{e^{\gamma h_*}}\gamma z_i^{-\gamma}.
\end{align*}
These are \eqref{B:eq:logarithmic-positive-sum} and \eqref{B:eq:logarithmic-negative-sum}. Finally, $h_j^2\le h_{\max}h_j$ and \eqref{B:eq:logarithmic-positive-sum} give
\[
 \sum_{j=1}^M h_j^2z_j^\theta
 \le h_{\max}\sum_{j=1}^M h_jz_j^\theta
 \le\frac{e^{\theta h_*}}\theta h_{\max}z_M^\theta,
 \qquad \theta>0.
\]
The constants in these inequalities are independent of $z_0,z_M,M$.
\end{proof}

For the partitions $z_{n,k}=k/n$ and $z_{n,k}=\theta_V(k/n)$, the bounds on $\theta_V'$ in Lemma~\ref{B:lem:bessel-coordinate-regularity} give $z_{n,k}\ge c k/n$ and
$\Delta_{n,k}\le C/n$. Consequently, for
$k\ge2$,
\[
 h_{n,k}=\log\left(1+\frac{\Delta_{n,k}}{z_{n,k-1}}\right)
       \le\frac{\Delta_{n,k}}{z_{n,k-1}}\le\frac{C}{c(k-1)}.
\]
On cells with $k>k_-$ this is at most $C/(c k_-)$. On a fixed interval $z_{n,k-1}\ge\varepsilon>0$ it is at most $C/(n\varepsilon)$. Thus the bound by $h_*$ holds when
$k_-\ge C/(ch_*)$, or when $n\varepsilon\ge C/h_*$ on a fixed interval. These conditions are independent of the upper truncation index.

\subsubsection{Interpolation of the kernel coefficients}
\label{B:subsection:interpolation-of-kernel-coefficients}

To compare two kernels, we interpolate their positive coefficients and estimate the derivative of the resulting inverse. The restrictions on the exponents leave room for both the
power singularity and the small losses in the coefficient bounds.

Fix $a>0$ and
\begin{equation}
 1<\lambda<\min\{a+1,2\},\qquad
 0<\eta<\frac1{16}\min\{a+1-\lambda,\,2-\lambda\}.
 \label{B:eq:stability-exponents}
\end{equation}
In particular,
\begin{equation}
 a+1-2\eta>0,\qquad
 a-2\eta>\lambda-1,\qquad 2-\lambda-8\eta>0.
 \label{B:eq:stability-integrability-exponents}
\end{equation}
Let $p_i,Q_i>0$ satisfy, with $D_\eta<\infty$ independent of the partition,
\begin{equation}
 p_i\le D_\eta z_i^{a/2-\eta},\qquad
 Q_i\le D_\eta h_i z_i^{-a-2\eta}.
 \label{B:eq:reference-coefficient-bounds}
\end{equation}
Here $p_i$ is a multiplicative kernel coefficient and $Q_i$ is a mass of the scale measure. The corresponding speed mass on $I_i$ is $\Delta_i p_i^2$.

The coefficients $p^B,q^B$ in \eqref{B:eq:2-32} satisfy \eqref{B:eq:reference-coefficient-bounds}. The first inequality is \eqref{B:eq:bessel-p-bounds}. For the second, put
$c=a+2\eta>0$ and use \eqref{B:eq:bessel-q-bounds}:
\begin{align*}
 Q_i^B&\le C_{\omega,\eta}\int_{z_{i-1}}^{z_i}t^{-c-1}\,\dd t
   = C_{\omega,\eta}z_i^{-c}\frac{e^{ch_i}-1}{c}\\
 &\le C_{\omega,\eta}e^{ch_*}h_i z_i^{-c}.
\end{align*}
The last inequality follows by integrating $e^{ct}$ on $[0,h_i]$. The same $C_{\omega,\eta}$ bounds the integral on every selected interval of this construction.

Let $\delta_i,\rho_i\in\mathbb R$ obey
\begin{equation}
 \max_i(|\delta_i|+|\rho_i|)\le1.
 \label{B:eq:logarithmic-errors-bounded}
\end{equation}
For $0\le\tau\le1$, set
\begin{equation}
 p_i^\tau=p_i e^{\tau\delta_i},\qquad
 Q_i^\tau=Q_i e^{\tau(-2\delta_i+\rho_i)},\qquad
 R^\tau(i,j)=p_i^\tau p_j^\tau
                      \sum_{\ell\ge i\vee j}Q_\ell^\tau.
 \label{B:eq:kernel-coefficient-interpolation}
\end{equation}
These are the logarithmic errors defined in \eqref{B:eq:2-35} when $p_i,Q_i$ are the coefficients defined in \eqref{B:eq:2-32}. Then $R^0$ is \eqref{B:eq:2-34}, and $R^1$ is the
restricted finite inverse \eqref{B:eq:2-37}. In particular, both restrict the triangular factor before taking its squared modulus.

For $g\in L^2(0,1)$, summation over the cell pairs gives
\[
 \langle R^\tau g,g\rangle
 =\sum_{\ell=1}^M Q_\ell^\tau
   \left|\sum_{i\le\ell}p_i^\tau\int_{I_i}g(s)\,\dd s\right|^2\ge0.
\]
This identity includes every diagonal cell pair. The error bound gives
\[
 |\delta_i|\le1,\qquad
 |-2\delta_i+\rho_i|\le2|\delta_i|+|\rho_i|\le2,
\]
and consequently
\[
 p_i^\tau\le eD_\eta z_i^{a/2-\eta},\qquad
 Q_i^\tau\le e^2D_\eta h_i z_i^{-a-2\eta},\qquad 0\le\tau\le1.
\]
All estimates through Section~\ref{B:subsection:weighted-comparison-of-inverses} have constants depending only on $a,\lambda,\eta,h_*,D_\eta$. They are independent of
$M,z_0,z_M,\tau$, and the errors satisfying \eqref{B:eq:logarithmic-errors-bounded}. Define
\begin{equation}
 M_i^\tau=\sum_{j\le i}\Delta_j(p_j^\tau)^2,
 \qquad T_i^\tau=\sum_{\ell\ge i}Q_\ell^\tau.
 \label{B:eq:cumulative-kernel-coefficients}
\end{equation}

\subsubsection{Bounds for the kernel and its rows}
\label{B:subsection:kernel-and-row-estimates}

We bound the two sums defined above: the partial sum of the squared multiplicative coefficients weighted by the cell lengths, and the tail sum of the coefficients inside the
inverse formula. These bounds imply the row estimates needed to sum the coefficient variations.

\begin{lemma}\label{B:lem:cumulative-coefficient-bounds}
Uniformly in $i$ and $\tau$,
\begin{equation}
 M_i^\tau\le C z_i^{a+1-2\eta},\qquad
 T_i^\tau\le C z_i^{-a-2\eta}.
 \label{B:eq:cumulative-coefficient-bounds}
\end{equation}
\end{lemma}
\begin{proof}
By \eqref{B:eq:logarithmic-cell-length} and the bounds on $p_j^\tau$,
\begin{align*}
 M_i^\tau
 &\le e^2D_\eta^2\sum_{j\le i}\Delta_jz_j^{a-2\eta}
 \le e^2D_\eta^2\sum_{j\le i}h_jz_j^{a+1-2\eta}\\
 &\le\frac{e^2D_\eta^2e^{(a+1-2\eta)h_*}}{a+1-2\eta}
                  z_i^{a+1-2\eta}.
\end{align*}
Here \eqref{B:eq:logarithmic-positive-sum} applies with $\gamma=a+1-2\eta>0$. For $T_i^\tau$, use \eqref{B:eq:logarithmic-negative-sum} with $\gamma=a+2\eta>0$:
\[
 T_i^\tau
 \le e^2D_\eta\sum_{\ell\ge i}h_\ell z_\ell^{-a-2\eta}
 \le\frac{e^2D_\eta e^{(a+2\eta)h_*}}{a+2\eta}z_i^{-a-2\eta}.
\]
Both constants are uniform in $i,\tau,M,z_0,z_M$.
\end{proof}

Let
\begin{equation}
 A_i^\tau=\sum_{j=1}^M\Delta_j|R^\tau(i,j)|^2.
 \label{B:eq:kernel-row-square}
\end{equation}
This is the squared $L^2$ norm of the kernel with its first variable fixed in $I_i$.

\begin{proposition}\label{B:prop:weighted-kernel-rows}
For $i\le j$,
\begin{equation}
 0\le R^\tau(i,j)\le C z_i^{a/2-\eta}z_j^{-a/2-3\eta}.
 \label{B:eq:ordered-kernel-bound}
\end{equation}
For every $i$,
\begin{equation}
 A_i^\tau\le C z_i^{\lambda-1}.
 \label{B:eq:weighted-kernel-row-bound}
\end{equation}
\end{proposition}
\begin{proof}
For $i\le j$, substitution of the coefficient bounds gives
\[
 R^\tau(i,j)=p_i^\tau p_j^\tau T_j^\tau
 \le Cz_i^{a/2-\eta}z_j^{a/2-\eta}z_j^{-a-2\eta}
 =Cz_i^{a/2-\eta}z_j^{-a/2-3\eta}.
\]
For $j\le i$, the corresponding expression is
\[
 R^\tau(i,j)=p_i^\tau p_j^\tau T_i^\tau
 \le Cz_j^{a/2-\eta}z_i^{-a/2-3\eta}.
\]
The exact decomposition of the squared row is
\[
 A_i^\tau
 =(p_i^\tau)^2(T_i^\tau)^2M_i^\tau
 +(p_i^\tau)^2\sum_{j>i}\Delta_j(p_j^\tau)^2(T_j^\tau)^2.
\]
The diagonal term is included in the first summand. Its part of the row satisfies
\begin{align}
 \sum_{j\le i}\Delta_j|R^\tau(i,j)|^2
 &\le C z_i^{-a-6\eta}\sum_{j\le i}h_jz_j^{a+1-2\eta}
 \le C z_i^{1-8\eta}
 \le C z_i^{\lambda-1}.
 \label{B:eq:lower-part-row-bound}
\end{align}
Equation~\eqref{B:eq:logarithmic-positive-sum} applies with $a+1-2\eta>0$. Since $1-8\eta-(\lambda-1)=2-\lambda-8\eta>0$ and $z_i\le1$, the last power is bounded as stated. For
$j\ge i$, the identity
\[
 z_i^{a-2\eta}z_j^{1-a-6\eta}
 =z_i^{\lambda-1}
   \left(\frac{z_i}{z_j}\right)^{a+1-\lambda-2\eta}
   z_j^{2-\lambda-8\eta}
\]
has positive exponents $a+1-\lambda-2\eta$ and $2-\lambda-8\eta$. Consequently,
\begin{align}
 \sum_{j\ge i}\Delta_j|R^\tau(i,j)|^2
 &\le C z_i^{a-2\eta}\sum_{j\ge i}h_jz_j^{1-a-6\eta}\notag\\
 &=C z_i^{\lambda-1}\sum_{j\ge i}h_j
       \left(\frac{z_i}{z_j}\right)^{a-2\eta-\lambda+1}
                                   z_j^{2-\lambda-8\eta}\notag\\
 &\le C z_i^{\lambda-1}\sum_{j=1}^M h_jz_j^{2-\lambda-8\eta}
 \le C z_i^{\lambda-1}.
 \label{B:eq:upper-part-row-bound}
\end{align}
In the last step,
\[
 \sum_{j=1}^M h_jz_j^{2-\lambda-8\eta}
 \le\frac{e^{(2-\lambda-8\eta)h_*}}{2-\lambda-8\eta}
                   z_M^{2-\lambda-8\eta}
 \le\frac{e^{(2-\lambda-8\eta)h_*}}{2-\lambda-8\eta}.
\]
The exact row decomposition uses $j>i$ in its second term, which is bounded by the displayed sum over $j\ge i$. Adding the two bounds proves \eqref{B:eq:weighted-kernel-row-bound}.
\end{proof}

The row estimate also gives a bound for the whole kernel:
\[
 \|R^\tau\|_{\mathfrak S_2}^2
 =\sum_i\Delta_i A_i^\tau
 \le C\sum_i h_i z_i^\lambda
 \le C\frac{e^{\lambda h_*}}\lambda.
\]

\subsubsection{Estimates for coefficient variations}
\label{B:subsection:estimates-for-coefficient-variations}

A change in one coefficient of the tail sum contributes an operator of rank at most one, while a change in a multiplicative coefficient acts on a row or column of the kernel.
Treating these contributions separately preserves the weights obtained in the preceding estimates.

Define $v_\ell^\tau=\sum_{i\le\ell}p_i^\tau\mathbf1_{I_i}$ and the positive operator
\begin{equation}
 K_\ell^\tau=Q_\ell^\tau
               |v_\ell^\tau\rangle\langle v_\ell^\tau|.
 \label{B:eq:scale-coefficient-operator}
\end{equation}
Thus its kernel is $Q_\ell^\tau p_i^\tau p_j^\tau\mathbf1_{\{i,j\le\ell\}}$, and $R^\tau=\sum_{\ell=1}^M K_\ell^\tau$. For a sequence $u_i$, let $D_u$ be multiplication by $u_i$ on
$I_i$ and by zero elsewhere.

\begin{lemma}\label{B:lem:scale-coefficient-norm}
For every $\ell$ and $\tau$,
\begin{equation}
 \|K_\ell^\tau\|_{\mathfrak S_2}
        =Q_\ell^\tau M_\ell^\tau
        \le C h_\ell z_\ell^{1-4\eta}.
 \label{B:eq:scale-coefficient-norm}
\end{equation}
\end{lemma}
\begin{proof}
Disjointness of the cells gives
\[
 \|v_\ell^\tau\|_2^2
 =\sum_{i\le\ell}\Delta_i(p_i^\tau)^2=M_\ell^\tau.
\]
Integrating the squared kernel over its cell pairs yields
\begin{align*}
 \|K_\ell^\tau\|_{\mathfrak S_2}^2
 &=(Q_\ell^\tau)^2
   \sum_{i,j\le\ell}\Delta_i\Delta_j(p_i^\tau)^2(p_j^\tau)^2\\
 &=(Q_\ell^\tau)^2
   \left(\sum_{i\le\ell}\Delta_i(p_i^\tau)^2\right)^2
 =(Q_\ell^\tau M_\ell^\tau)^2.
\end{align*}
Taking square roots and applying \eqref{B:eq:reference-coefficient-bounds} and \eqref{B:eq:cumulative-coefficient-bounds}, we obtain
\[
 \|K_\ell^\tau\|_{\mathfrak S_2}
 =Q_\ell^\tau M_\ell^\tau
 \le Ch_\ell z_\ell^{-a-2\eta}z_\ell^{a+1-2\eta}
 =Ch_\ell z_\ell^{1-4\eta}.
\]
\end{proof}

For $\ell,m\in\{1,\ldots,M\}$, integration on the common cells gives
\[
 \langle v_\ell^\tau,v_m^\tau\rangle
 =\sum_{i\le\ell\wedge m}\Delta_i(p_i^\tau)^2
 =M_{\ell\wedge m}^\tau.
\]
Integrating the product of the two kernels gives
\begin{align*}
 \langle K_\ell^\tau,K_m^\tau\rangle_{\mathfrak S_2}
 &=Q_\ell^\tau Q_m^\tau
   \sum_{i,j\le\ell\wedge m}\Delta_i\Delta_j(p_i^\tau)^2(p_j^\tau)^2\\
 &=Q_\ell^\tau Q_m^\tau(M_{\ell\wedge m}^\tau)^2.
\end{align*}
Here $\langle A,B\rangle_{\mathfrak S_2}=\operatorname{Tr}(B^*A)$ is linear in its first argument. Thus, for complex
$f_\ell$,
\[
 \left\|\sum_\ell f_\ell K_\ell^\tau\right\|_{\mathfrak S_2}^2
 =\sum_{\ell,m}f_\ell\overline{f_m}\,
       Q_\ell^\tau Q_m^\tau(M_{\ell\wedge m}^\tau)^2.
\]

\begin{proposition}\label{B:prop:kernel-coefficient-variations}
For arbitrary complex sequences $f_i,u_i$,
\begin{align}
 \left\|\sum_{\ell=1}^M f_\ell K_\ell^\tau\right\|_{\mathfrak S_2}^2
       &\le C\sum_{\ell=1}^M h_\ell z_\ell^\lambda|f_\ell|^2,
         \label{B:eq:scale-variation-bound}\\
 \|D_uR^\tau+R^\tau D_u\|_{\mathfrak S_2}^2
       &\le C\sum_{i=1}^M h_i z_i^\lambda|u_i|^2.
         \label{B:eq:multiplication-variation-bound}
\end{align}
\end{proposition}
\begin{proof}
The sequence $M_i^\tau$ is nonnegative and increasing, so $(M_{\ell\wedge m}^\tau)^2\le M_\ell^\tau M_m^\tau$. Taking absolute values in the preceding double sum gives
\begin{align*}
 \left\|\sum_\ell f_\ell K_\ell^\tau\right\|_{\mathfrak S_2}^2
 &\le\sum_{\ell,m}|f_\ell||f_m|
       Q_\ell^\tau Q_m^\tau M_\ell^\tau M_m^\tau\\
 &=\left(\sum_\ell|f_\ell|Q_\ell^\tau M_\ell^\tau\right)^2
 \le C\left(\sum_\ell|f_\ell|h_\ell z_\ell^{1-4\eta}\right)^2.
\end{align*}
For each summand write
\[
 |f_\ell|h_\ell z_\ell^{1-4\eta}
 =\bigl(\sqrt{h_\ell}z_\ell^{\lambda/2}|f_\ell|\bigr)
  \bigl(\sqrt{h_\ell}z_\ell^{1-4\eta-\lambda/2}\bigr).
\]
Cauchy--Schwarz and \eqref{B:eq:logarithmic-positive-sum} give
\begin{align*}
 \left(\sum_\ell |f_\ell|h_\ell z_\ell^{1-4\eta}\right)^2
 &\le\left(\sum_\ell h_\ell z_\ell^\lambda|f_\ell|^2\right)
       \left(\sum_\ell h_\ell z_\ell^{2-8\eta-\lambda}\right)\\
 &\le\frac{e^{(2-8\eta-\lambda)h_*}}{2-8\eta-\lambda}
                 \sum_\ell h_\ell z_\ell^\lambda|f_\ell|^2.
\end{align*}
The exponent $2-8\eta-\lambda$ is positive, so this proves \eqref{B:eq:scale-variation-bound} uniformly in the ends of the selected interval.

In the convention \eqref{B:eq:integral-kernel-convention}, the first variable is integrated. Hence, for $(s,t)\in I_i\times I_j$,
\[
 K_{D_uR^\tau}(s,t)=u_jR^\tau(i,j),\qquad
 K_{R^\tau D_u}(s,t)=u_iR^\tau(i,j).
\]
Their separate squared norms are
\begin{align*}
 \|R^\tau D_u\|_{\mathfrak S_2}^2
 &=\sum_i\Delta_i|u_i|^2\sum_j\Delta_j|R^\tau(i,j)|^2
 =\sum_i\Delta_i|u_i|^2 A_i^\tau,\\
 \|D_uR^\tau\|_{\mathfrak S_2}^2
 &=\sum_j\Delta_j|u_j|^2\sum_i\Delta_i|R^\tau(i,j)|^2
 =\sum_j\Delta_j|u_j|^2 A_j^\tau.
\end{align*}
The last equality uses $R^\tau(i,j)=R^\tau(j,i)$. These formulas hold for complex $u_i$. Applying the row bound gives
\begin{align*}
 \|D_uR^\tau+R^\tau D_u\|_{\mathfrak S_2}^2
 &\le2\|D_uR^\tau\|_{\mathfrak S_2}^2
      +2\|R^\tau D_u\|_{\mathfrak S_2}^2\\
 &=4\sum_i\Delta_i|u_i|^2 A_i^\tau
 \le C\sum_i\Delta_i z_i^{\lambda-1}|u_i|^2
 \le C\sum_i h_i z_i^\lambda|u_i|^2.
\end{align*}
The final inequality is \eqref{B:eq:logarithmic-cell-length}, and the result is \eqref{B:eq:multiplication-variation-bound}.
\end{proof}

\subsubsection{Weighted comparison of the inverses}
\label{B:subsection:weighted-comparison-of-inverses}

The two types of coefficient variation can now be integrated along the interpolation. The row estimate controls multiplication by the logarithmic error, and the estimate for
operators of rank one controls changes in the coefficients of the tail sum. Under the assumptions of Section~\ref{B:subsection:interpolation-of-kernel-coefficients}, they give one
comparison bound with constants independent of the selected interval.

\begin{theorem}\label{B:thm:weighted-inverse-stability}
The two kernels satisfy
\begin{equation}
 \|R^1-R^0\|_{\mathfrak S_2}^2
 \le C\sum_{i=1}^M h_i z_i^\lambda
                         (|\delta_i|^2+|\rho_i|^2).
 \label{B:eq:weighted-inverse-stability}
\end{equation}
The constant is uniform in the selected interval and the number of cells.
\end{theorem}
\begin{proof}
For the fixed partition, let $r_{ij}(\tau)$ be the value of the interpolated kernel on $I_i\times I_j$. The normalized cell indicators $\Delta_i^{-1/2}\mathbf1_{I_i}$ form an
orthonormal basis of the step-function subspace. With the first variable integrated, the corresponding matrix entry in row $j$, column $i$, is
$\sqrt{\Delta_i\Delta_j}\,r_{ij}(\tau)$. Its squared Hilbert--Schmidt norm is therefore
\[
 \|R^\tau\|_{\mathfrak S_2}^2
 =\sum_{i,j}\Delta_i\Delta_j|r_{ij}(\tau)|^2.
\]
Every $r_{ij}(\tau)$ is a finite sum of products of exponential coefficient functions. In particular it is continuously differentiable on $[0,1]$. For each interior $\tau$, the
square of the norm of the difference quotient minus its entrywise derivative is
\[
 \sum_{i,j}\Delta_i\Delta_j
 \left|\frac{r_{ij}(\tau+t)-r_{ij}(\tau)}t-r_{ij}'(\tau)\right|^2
 \longrightarrow0.
\]
The sum is finite, so its limit is zero. The derivative is continuous, with one-sided derivatives at the endpoints. Thus coefficient differentiation also gives differentiation in
Hilbert--Schmidt norm.

Put $f_\ell=-2\delta_\ell+\rho_\ell$. Differentiating each entry gives
\begin{align}
 \partial_\tau R^\tau(i,j)
 &=(\delta_i+\delta_j)R^\tau(i,j)
     +p_i^\tau p_j^\tau\sum_{\ell\ge i\vee j}f_\ell Q_\ell^\tau,
        \label{B:eq:interpolated-kernel-derivative}\\
 \partial_\tau R^\tau
 &=D_\delta R^\tau+R^\tau D_\delta
                       +\sum_\ell f_\ell K_\ell^\tau.\notag
\end{align}
The maximum assumption implies $p_i^\tau\le e p_i$ and $Q_i^\tau\le e^2Q_i$ for $0\le\tau\le1$. Thus the constants $C_{\mathrm m}$ and $C_{\mathrm s}$ in
\eqref{B:eq:multiplication-variation-bound} and \eqref{B:eq:scale-variation-bound} are uniform in $\tau$ and in the selected interval. Since $|f_i|^2\le8|\delta_i|^2+2|\rho_i|^2$,
those estimates imply
\[
 \begin{aligned}
 \|\partial_\tau R^\tau\|_{\mathfrak S_2}^2
 &\le2C_{\mathrm m}\sum_i h_i z_i^\lambda|\delta_i|^2
       +2C_{\mathrm s}\sum_i h_i z_i^\lambda|f_i|^2\\
 &\le(2C_{\mathrm m}+16C_{\mathrm s})
          \sum_i h_i z_i^\lambda(|\delta_i|^2+|\rho_i|^2).
 \end{aligned}
\]
To apply the weighted estimate for the derivative, write
\[
 \mathfrak d=\sum_i h_i z_i^\lambda
                   (|\delta_i|^2+|\rho_i|^2).
\]
The preceding coefficient estimates give $\|\partial_\tau R^\tau\|_{\mathfrak S_2}^2\le C\mathfrak d$ uniformly for $0\le\tau\le1$. The ordinary fundamental theorem of calculus for
each finite matrix entry gives
\[
 R^1-R^0=\int_0^1\partial_\tau R^\tau\,\dd\tau.
\]
The integrand is continuous in a finite-dimensional subspace of the Hilbert--Schmidt space. The triangle inequality for its integral and scalar Cauchy--Schwarz imply
\[
 \begin{split}
 \|R^1-R^0\|_{\mathfrak S_2}^2
 &\le\left(\int_0^1\|\partial_\tau R^\tau\|_{\mathfrak S_2}\,\dd\tau\right)^2\\
 &\le\int_0^1\|\partial_\tau R^\tau\|_{\mathfrak S_2}^2\,\dd\tau
 \le C\mathfrak d.
 \end{split}
\]
Differentiation and integration take place in a finite-dimensional space. The final constant comes from the weighted coefficient estimates and is uniform in the number and relative
lengths of the cells.

\end{proof}

For the application, fix a Brownian path satisfying Lemma~\ref{B:lem:bessel-brownian-growth}. The coefficients in \eqref{B:eq:2-32} obey \eqref{B:eq:reference-coefficient-bounds}
on every selected interval with one finite constant $D_\eta(\omega)$. Suppose the errors of \eqref{B:eq:2-35} satisfy
\[
 \max_{k\in\mathscr I_n^{\mathrm{ret}}}
          (|\delta_{n,k}|+|\rho_{n,k}|)\longrightarrow0,
 \qquad
 \sum_{k\in\mathscr I_n^{\mathrm{ret}}}h_{n,k}z_{n,k}^\lambda
       (|\delta_{n,k}|^2+|\rho_{n,k}|^2)\le\frac{A_\omega}{n}
\]
for every $n\ge n_0(\omega)$, where $A_\omega<\infty$. Increasing $n_0(\omega)$ makes the displayed maximum at most one for the same values of $n$. The theorem then gives
\[
 \|R_n^{\Delta,0}-R_n^{B,\mathrm{step},0}\|_{\mathfrak S_2}^2
 \le\frac{C(a,\lambda,\eta,h_*,D_\eta(\omega))A_\omega}{n}.
\]
Both operators use the same selected cells and the restricted factors of Section~\ref{B:subsection:retained-cells-errors-and-weights}. The estimate holds on the fixed path and uses
only finiteness of $A_\omega$ and $D_\eta(\omega)$.

The exponent conditions used above are $a-2\eta>\lambda-1$ in \eqref{B:eq:upper-part-row-bound} and $2-8\eta-\lambda>0$ in \eqref{B:eq:scale-variation-bound}. The range of
$\lambda$ in \eqref{B:eq:stability-exponents} permits both with a fixed positive $\eta$ as the lower cutoff tends to zero. The remaining subsubsections estimate the terms removed
by restricting the factors.

\subsubsection{Integration on the full interval}
\label{B:subsection:integration-on-the-full-interval}

We retain $a>0$ and estimate the parts of the positive inverse removed at the two boundary points. The same integral calculations apply to the continuum kernel and to the finite
kernel when each cell is assigned its physical length. We use either Lebesgue measure on $X=(0,1]$, or
\begin{equation}
 X_n=\{z_1,\ldots,z_n\},\qquad
 \mathfrak m_n=\sum_{k=1}^n\Delta_k\delta_{z_k},\qquad
 \frac{c_{\rm m}}n\le\Delta_k\le\frac{C_{\rm m}}n,
 \label{B:eq:full-grid-measure}
\end{equation}
on the full partition of \eqref{B:eq:cell-embedding-definition}, with $z_0=0$ and $z_n=1$. We write $\mathfrak m$ for the measure in use and $X$ for its underlying set. In the
atomic case, the map
\[
 f\longmapsto\sum_{k=1}^n f(z_k)\mathbf1_{I_k}
\]
is an isometry from $L^2(\mathfrak m_n)$ onto the functions constant on each cell, since both squared norms equal $\sum_k\Delta_k|f(z_k)|^2$. An integral kernel on $X_n^2$
therefore has the same Hilbert--Schmidt norm as the corresponding step kernel. For a kernel $S$ on $X_n^2$, this identification gives
\[
 \|S\|_{\mathfrak S_2}^2
 =\sum_{i,j=1}^n\Delta_i\Delta_j|S(z_i,z_j)|^2
 =\sum_{i=1}^n\Delta_i^2|S(z_i,z_i)|^2
   +2\sum_{i<j}\Delta_i\Delta_j|S(z_i,z_j)|^2
\]
when $S$ is symmetric. In particular the atom at $z_n=1$ has mass $\Delta_n$. The identities
\[
 \mathfrak m_n(X_n)=1,\qquad
 \mathfrak m_n((0,z_j])=z_j,\qquad
 \mathfrak m_n((z_j,1])=1-z_j
\]
follow by summing the cell lengths. We keep the diagonal entries and the indicated boundary point conventions throughout the truncation.

For $c>0$ and $0<y\le1$, put
\[
 F_c(y)=
 \begin{cases}
  1,&0<c<1,\\
  1+|\log y|,&c=1,\\
  y^{1-c},&c>1.
 \end{cases}
\]
\begin{lemma}\label{B:lem:full-interval-power-integrals}
For $v>-1$, $c>0$, and $y\in X$,
\begin{equation}
 \int_{s\le y}s^v\,\mathfrak m(\dd s)\le C_v y^{v+1},\qquad
 \int_{t\ge y}t^{-c}\,\mathfrak m(\dd t)\le C_c F_c(y).
 \label{B:eq:full-interval-power-bounds}
\end{equation}
The constants in the atomic case depend only on the indicated exponent and $c_{\rm m},C_{\rm m}$.
\end{lemma}
\begin{proof}
For Lebesgue measure, direct integration gives
\[
 \int_0^y s^v\,\dd s=\frac{y^{v+1}}{v+1},\qquad
 \int_y^1t^{-c}\,\dd t=
 \begin{cases}
  (1-y^{1-c})/(1-c),&c\ne1,\\
  -\log y,&c=1.
 \end{cases}
\]
If $0<c<1$, the second integral is at most $(1-c)^{-1}$. If $c>1$, it is at most $y^{1-c}/(c-1)$. The case $c=1$ is bounded by $F_1(y)$.

For the atomic measure, summing \eqref{B:eq:full-grid-measure} gives
\[
 c_{\rm m}\frac{k}{n}\le z_k\le C_{\rm m}\frac{k}{n},\qquad
 c_{\rm m}\le1\le C_{\rm m}.
\]
For $y=z_j$ and every real $v$, positivity implies $z_k^v\le\max\{c_{\rm m}^v,C_{\rm m}^v\}(k/n)^v$. Hence
\[
 \sum_{k\le j}\Delta_k z_k^v
 \le C_{\rm m}\max\{c_{\rm m}^v,C_{\rm m}^v\}
                 n^{-v-1}\sum_{k=1}^j k^v.
\]
When $v\ge0$, each term is at most $j^v$, so the sum is at most $j^{v+1}$. When $-1<v<0$, monotonicity gives
\[
 \sum_{k=1}^j k^v
 \le1+\int_1^j t^v\,\dd t
 =1+\frac{j^{v+1}-1}{v+1}
 \le\frac{j^{v+1}}{v+1}.
\]
Since $v+1>0$ and $j/n\le y/c_{\rm m}$, these bounds prove the first inequality in \eqref{B:eq:full-interval-power-bounds}.

For the other integral,
\[
 \sum_{k=j}^n\Delta_k z_k^{-c}
 \le C_{\rm m}c_{\rm m}^{-c}n^{c-1}\sum_{k=j}^n k^{-c}.
\]
If $c>1$, then
\[
 \sum_{k=j}^n k^{-c}
 \le j^{-c}+\int_j^\infty t^{-c}\,\dd t
 \le\left(1+\frac1{c-1}\right)j^{1-c}.
\]
Moreover $y\le C_{\rm m}j/n$ and $1-c<0$ imply $(j/n)^{1-c}\le C_{\rm m}^{c-1}y^{1-c}$. If $c=1$, then
\[
 \sum_{k=j}^n k^{-1}\le1+\log(n/j)
 \le1+\log C_{\rm m}+|\log y|
 \le(1+\log C_{\rm m})F_1(y).
\]
If $0<c<1$, then
\[
 \sum_{k=j}^n k^{-c}
 \le1+\int_1^n t^{-c}\,\dd t
 =1+\frac{n^{1-c}-1}{1-c}
 \le\frac{n^{1-c}}{1-c}.
\]
Multiplication by $n^{c-1}$ gives a uniform constant in this case. These estimates include $j=1$ and $j=n$. At $j=n$, the atomic integral is $\Delta_n$ and the bound uses
$F_c(1)=1$.
\end{proof}

\subsubsection{Bounds for the full Green kernel}
\label{B:subsection:bounds-for-the-full-green-kernel}

The truncation argument uses the full inverse. Its factorization gives positivity and identifies its trace with the squared Hilbert--Schmidt norm of the factor. We establish the
required integrability and then estimate the rows of the resulting Green kernel.

Let $p,r$ be finite nonnegative measurable functions on $X$, locally bounded on $(0,1]$ in the continuum case, and set
\begin{equation}
 K(s,t)=p(s)r(t)\mathbf1_{\{s\le t\}},\qquad
 (\mathcal Kf)(t)=\int_X K(s,t)f(s)\,\mathfrak m(\dd s),\qquad
 \mathcal R=\mathcal K^*\mathcal K.
 \label{B:eq:full-triangular-factor}
\end{equation}
This extends the factors of Section~\ref{B:subsection:embedded-kernels-and-scalar-invariance} and Section~\ref{B:subsection:the-stochastic-bessel-realization} to the two measures
above. Assume that for every sufficiently small $\xi>0$ there is a constant $A_\xi<\infty$ such that
\begin{equation}
 K(s,t)\le A_\xi s^{a/2-\xi}t^{-(a+1)/2-\xi},
       \qquad 0<s\le t\le1.
 \label{B:eq:triangular-factor-bound}
\end{equation}
For a family of finite partitions these constants are independent of $n$. Constants below depend only on $a$, the chosen positive exponents, $c_{\rm m},C_{\rm m}$, and the
necessary $A_\xi$. The bound on $K(s,t)=p(s)r(t)$ will be used before each integration. Under $p\mapsto Cp$ and $r\mapsto C^{-1}r$, with $C>0$, the kernel and the constants $A_\xi$
are unchanged.

\begin{proposition}\label{B:prop:full-factor-and-green-kernel}
The factor $\mathcal K$ is Hilbert--Schmidt, and $\mathcal R$ is positive and trace class, with
\begin{equation}
 \operatorname{Tr}\mathcal R=\|\mathcal K\|_{\mathfrak S_2}^2\le C.
 \label{B:eq:full-factor-norm}
\end{equation}
Its kernel is
\begin{equation}
 R(s,t)=p(s)p(t)\int_{u\ge s\vee t}r(u)^2\,\mathfrak m(\dd u).
 \label{B:eq:full-green-kernel}
\end{equation}
For every $\varepsilon>0$,
\begin{equation}
 |R(s,t)|\le C_\varepsilon
       \left(\frac{s\wedge t}{s\vee t}\right)^{a/2}(st)^{-\varepsilon}.
 \label{B:eq:full-green-kernel-bound}
\end{equation}
\end{proposition}
\begin{proof}
Choose a positive $\xi$ for which \eqref{B:eq:triangular-factor-bound} holds, with $2\xi<a+1$ and $4\xi<1$. Tonelli's theorem and Lemma~\ref{B:lem:full-interval-power-integrals}
give
\begin{align*}
 \int_X\int_X|K(s,t)|^2\,\mathfrak m(\dd s)\mathfrak m(\dd t)
 &\le A_\xi^2\int_X t^{-a-1-2\xi}
             \int_{s\le t}s^{a-2\xi}\,\mathfrak m(\dd s)
                                       \mathfrak m(\dd t)\\
 &\le A_\xi^2C_{a-2\xi}\int_X t^{-4\xi}\,\mathfrak m(\dd t)
 \le A_\xi^2C_{a-2\xi}C_{-4\xi}.
\end{align*}
Both exponents to which the first power estimate is applied exceed $-1$. The diagonal atoms are included in the inner integral $s\le t$. Thus $K$ defines a Hilbert--Schmidt
operator. For every $f$, $\langle\mathcal Rf,f\rangle=\|\mathcal Kf\|_2^2\ge0$. For an orthonormal basis $(e_j)$, Parseval in the first kernel variable and Tonelli give
\[
 \sum_j\langle\mathcal Re_j,e_j\rangle
 =\sum_j\|\mathcal Ke_j\|_2^2
 =\int_X\int_X|K(s,t)|^2\,\mathfrak m(\dd s)\mathfrak m(\dd t).
\]
This proves positivity, trace-class membership, and \eqref{B:eq:full-factor-norm}.

For $f,g\in L^2(\mathfrak m)$, the nonnegative triple integral is
\begin{align*}
 &\int_X\int_X\int_X K(s,u)K(t,u)|f(s)||g(t)|\,
          \mathfrak m(\dd s)\mathfrak m(\dd t)\mathfrak m(\dd u)\\
 &\qquad=\int_X(\mathcal K|f|)(u)(\mathcal K|g|)(u)\,
                         \mathfrak m(\dd u)\\
 &\qquad\le\|\mathcal K|f|\|_2\|\mathcal K|g|\|_2
 \le\|\mathcal K\|_{\rm op}^2\|f\|_2\|g\|_2<\infty.
\end{align*}
Fubini therefore applies with the complex factors $f(s)\overline{g(t)}$. The two triangular indicators require $u\ge s\vee t$, and hence
\begin{align*}
 \langle\mathcal Kf,\mathcal Kg\rangle
 &=\int_X\int_X\left(\int_XK(s,u)K(t,u)\,\mathfrak m(\dd u)\right)
            f(s)\overline{g(t)}\,\mathfrak m(\dd s)\mathfrak m(\dd t)\\
 &=\int_X\int_Xp(s)p(t)
       \left(\int_{u\ge s\vee t}r(u)^2\,\mathfrak m(\dd u)\right)
            f(s)\overline{g(t)}\,\mathfrak m(\dd s)\mathfrak m(\dd t).
\end{align*}
The inner integral is finite for $s,t>0$ by the local boundedness assumption, and by finiteness of the sum in the atomic case. Furthermore, Cauchy--Schwarz in $u$ gives
\[
 |R(s,t)|^2\le R(s,s)R(t,t),\qquad
 \int_XR(s,s)\,\mathfrak m(\dd s)=\|\mathcal K\|_{\mathfrak S_2}^2.
\]
Consequently the displayed kernel belongs to $L^2(\mathfrak m\otimes \mathfrak m)$, with squared norm at most $\|\mathcal K\|_{\mathfrak S_2}^4$. The sesquilinear identity
identifies it as the kernel of $\mathcal K^*\mathcal K$.

For $s\le t$, multiply the two bounds on $K(s,u)$ and $K(t,u)$ inside that integral. Since $a+1+2\xi>1$,
\begin{align*}
 R(s,t)
 &\le A_\xi^2s^{a/2-\xi}t^{a/2-\xi}
              \int_{u\ge t}u^{-a-1-2\xi}\,\mathfrak m(\dd u)\\
 &\le A_\xi^2C_{a+1+2\xi}s^{a/2-\xi}t^{-a/2-3\xi}.
\end{align*}
Given $\varepsilon>0$, choose $\xi$ also to satisfy $3\xi\le\varepsilon$. The identity
\[
 s^{a/2-\xi}t^{-a/2-3\xi}
 =\left(\frac st\right)^{a/2}(st)^{-\varepsilon}
                  s^{\varepsilon-\xi}t^{\varepsilon-3\xi}
 \le\left(\frac st\right)^{a/2}(st)^{-\varepsilon}
\]
uses $s,t\le1$ and the nonnegative exponents in the last two factors. For $t\le s$, interchange $s$ and $t$ in the same bound, using the symmetry of \eqref{B:eq:full-green-kernel}.
This proves \eqref{B:eq:full-green-kernel-bound}.
\end{proof}

\subsubsection{Rows near zero and one}
\label{B:subsection:row-bounds-near-zero-and-one}

The exponent at zero can be seen directly for $p(s)=s^{a/2}$, $r(t)=t^{-(a+1)/2}$, $a>0$, with Lebesgue measure. For $s\le t$ their Green kernel is
$R(s,t)=a^{-1}(s/t)^{a/2}(1-t^a)$. Splitting a row at $t=s$ gives
\[
 \int_0^1|R(s,t)|^2\,\dd t
 \le\frac1{a^2}\left[\frac{s}{a+1}
                         +s^a\int_s^1t^{-a}\,\dd t\right].
\]
The last integral is $(1-s^{1-a})/(1-a)$ for $a\ne1$ and $\log(1/s)$ for $a=1$. The resulting upper bounds are $O_a(s^a)$, $O(s(1+|\log s|))$ and $O_a(s)$ when $a<1$, $a=1$ and
$a>1$, respectively. For every $\eta>0$, $s^\eta|\log s|\le(\mathrm e\eta)^{-1}$ on $(0,1]$. Thus the bound $s^{\min\{a,1\}-\eta}$ covers the transition at $a=1$. The next proof
incorporates the exponent losses in the factor bound and applies to both Lebesgue measure and the atomic measure.

Put $\nu_a=\min\{a,1\}$. The next lemma uses only a bound on a symmetric kernel; it also applies to the interpolated Brownian kernels in
Section~\ref{B:sec:approximation-by-brownian-motion-and-the-lower-bound}. Fix $0<\eta<\nu_a$, set $\gamma=\nu_a-\eta$, and choose $0<\varepsilon<\eta/8$. Suppose $S$ is a symmetric
measurable kernel with
\[
 |S(s,t)|\le D_\varepsilon
       \left(\frac{s\wedge t}{s\vee t}\right)^{a/2}(st)^{-\varepsilon}.
\]

\begin{lemma}\label{B:lem:green-kernel-row-bound}
Under this bound,
\begin{equation}
 \int_X|S(s,t)|^2\,\mathfrak m(\dd t)
        \le C_\eta s^{\nu_a-\eta},\qquad
 \|S\|_{\mathfrak S_2}^2\le C_\eta.
 \label{B:eq:full-kernel-row-bound}
\end{equation}
\end{lemma}
\begin{proof}
Recall $\gamma=\nu_a-\eta$ and $0<\varepsilon<\eta/8$. Then
\[
 0<\gamma<1,\qquad
 a-2\varepsilon-\gamma\ge\eta-2\varepsilon>0,\qquad
 1-\gamma-4\varepsilon>1-\nu_a+\eta/2>0.
\]
In particular $a-2\varepsilon>-1$. Splitting the row at $t=s$ counts the diagonal atom in the part $t\le s$. For this part,
\begin{align}
 \int_{t\le s}|S(s,t)|^2\,\mathfrak m(\dd t)
 &\le D_\varepsilon^2s^{-a-2\varepsilon}
                    \int_{t\le s}t^{a-2\varepsilon}\,\mathfrak m(\dd t)
 \le C_\varepsilon s^{1-4\varepsilon}\le C_\eta s^\gamma.
 \label{B:eq:green-kernel-lower-row}
\end{align}
The first power estimate gives the middle inequality. Since $s\le1$ and $1-4\varepsilon>\gamma$, the final inequality follows by comparison of the exponents.

The part $t>s$ is bounded by the following integral over $t\ge s$:
\begin{align*}
 \int_{t\ge s}|S(s,t)|^2\,\mathfrak m(\dd t)
 &\le D_\varepsilon^2s^{a-2\varepsilon}
                \int_{t\ge s}t^{-a-2\varepsilon}\,\mathfrak m(\dd t)\\
 &=D_\varepsilon^2s^\gamma\int_{t\ge s}
       (s/t)^{a-2\varepsilon-\gamma}t^{-\gamma-4\varepsilon}
                                            \,\mathfrak m(\dd t)\\
 &\le D_\varepsilon^2s^\gamma
                  \int_Xt^{-\gamma-4\varepsilon}\,\mathfrak m(\dd t)\\
 &\le D_\varepsilon^2C_{-\gamma-4\varepsilon}s^\gamma.
\end{align*}
Here $s/t\le1$, its exponent is positive, and $-\gamma-4\varepsilon>-1$. Adding the two row estimates gives the first kernel bound. Tonelli and the first power estimate, with
$y=1$, now give
\[
 \|S\|_{\mathfrak S_2}^2
 =\int_X\left(\int_X|S(s,t)|^2\,\mathfrak m(\dd t)\right)
                                     \mathfrak m(\dd s)
 \le C_\eta\int_Xs^\gamma\,\mathfrak m(\dd s)
 \le C_\eta C_\gamma.
\]
\end{proof}

We apply this lemma to $R$ from Proposition~\ref{B:prop:full-factor-and-green-kernel}. For $0<d<b\le1$ and $b\ge1/2$, let $P^-,P^0,P^+$ denote multiplication by the indicators of
$(0,d],(d,b],(b,1]$, respectively. In the atomic case take $d=z_m,b=z_M$ with $1\le m<M\le n$. Write
\begin{equation}
 \tau(t)=\mathfrak m([t,1])=
 \begin{cases}
  1-t,&\mathfrak m\text{ is Lebesgue measure},\\
  1-z_{j-1},&t=z_j\text{ for the atomic measure}.
 \end{cases}
 \label{B:eq:measure-above-a-point}
\end{equation}
In particular, $\tau(z_n)=\Delta_n$ on a finite partition.

\begin{proposition}\label{B:prop:truncation-row-bounds}
For $0<\eta<\nu_a$,
\begin{align}
 \|P^-\mathcal R\|_{\mathfrak S_2}^2
      +\|\mathcal RP^-\|_{\mathfrak S_2}^2
      &\le C_\eta d^{1+\nu_a-\eta},
          \label{B:eq:inverse-rows-near-zero}\\
 \|P^+\mathcal R\|_{\mathfrak S_2}^2
      +\|\mathcal RP^+\|_{\mathfrak S_2}^2
      &\le C_\eta(1-b)^3.
          \label{B:eq:truncation-upper-rows}
\end{align}
\end{proposition}
\begin{proof}
Multiplication by $P^-$ on the left restricts the second kernel variable; multiplication on the right restricts the first. Since $R(s,t)=R(t,s)$,
\begin{align*}
 \|P^-\mathcal R\|_{\mathfrak S_2}^2
       +\|\mathcal RP^-\|_{\mathfrak S_2}^2
 &=2\int_{s\le d}\int_X|R(s,t)|^2\,
                           \mathfrak m(\dd t)\mathfrak m(\dd s)\\
 &\le C_\eta\int_{s\le d}s^{\nu_a-\eta}\,\mathfrak m(\dd s)
 \le C_\eta d^{1+\nu_a-\eta}.
\end{align*}
The first inequality uses Lemma~\ref{B:lem:green-kernel-row-bound} with the bound from Proposition~\ref{B:prop:full-factor-and-green-kernel}. The second uses $\nu_a-\eta>0$ in the
first power estimate. It includes the first sampling point when $d=z_1$.

For the other boundary point, fix $b_*\in(0,1)$ and choose an admissible $\xi>0$ with $2\xi<a+1$. If $s\le t$ and $t\ge b_*$, then
\begin{align*}
 |R(s,t)|
 &\le A_\xi^2s^{a/2-\xi}t^{a/2-\xi}
             \int_{u\ge t}u^{-a-1-2\xi}\,\mathfrak m(\dd u)\\
 &\le A_\xi^2\max\{1,b_*^{a/2-\xi}\}
           b_*^{-a-1-2\xi}s^{a/2-\xi}\mathfrak m([t,1]).
\end{align*}
Thus
\begin{equation}
 |R(s,t)|\le C_{\xi,b_*}s^{a/2-\xi}\tau(t).
 \label{B:eq:green-kernel-near-one}
\end{equation}
Squaring this bound and integrating over $s\le t$ gives
\[
 \int_{s\le t}|R(s,t)|^2\,\mathfrak m(\dd s)
 \le C_{\xi,b_*}\tau(t)^2
            \int_{s\le t}s^{a-2\xi}\,\mathfrak m(\dd s)
 \le C_{\xi,b_*}\tau(t)^2.
\]
For $s>t$, interchange the variables in \eqref{B:eq:green-kernel-near-one}. The bound $t\ge b_*$ then gives $|R(s,t)|\le C_{\xi,b_*}\tau(s)$. We shall use
\begin{equation}
 \int_{s>t}\tau(s)^2\,\mathfrak m(\dd s)\le\tau(t)^3.
 \label{B:eq:upper-tail-length-integral}
\end{equation}
For Lebesgue measure its left side equals $\int_t^1(1-s)^2\,\dd s=(1-t)^3/3$. For the atomic measure, the integral over one physical cell satisfies
\begin{equation}
 3\int_x^y(1-u)^2\,\dd u
 =(y-x)\big[(1-x)^2+(1-x)(1-y)+(1-y)^2\big]
 \ge(y-x)(1-x)^2
 \label{B:eq:quadratic-integral-on-a-cell}
\end{equation}
Indeed, $3\int_x^y(1-u)^2\,\dd u=(1-x)^3-(1-y)^3$; factorization of this difference of cubes gives the displayed identity, and all three terms in brackets are nonnegative. When
$t=z_j$, summing over the cells after $j$ yields
\begin{align*}
 \int_{s>z_j}\tau(s)^2\,\mathfrak m_n(\dd s)
 &=\sum_{i=j+1}^n\Delta_i(1-z_{i-1})^2\\
 &\le3\sum_{i=j+1}^n\int_{z_{i-1}}^{z_i}(1-u)^2\,\dd u
 =(1-z_j)^3\le\tau(z_j)^3.
\end{align*}
The sum is zero at $j=n$. The diagonal atom at $t=z_j$ is already in the integral over $s\le t$ above. Combining the two parts of the row, and using $0\le\tau(t)\le1$, proves
\[
 \int_X|R(s,t)|^2\,\mathfrak m(\dd s)
 \le C_{\xi,b_*}\bigl(\tau(t)^2+\tau(t)^3\bigr)
 \le C_{\xi,b_*}\tau(t)^2\qquad(t\ge b_*).
\]
Take $b_*=1/2$. In the continuum, $\int_{t>b}\tau(t)^2\,\dd t=(1-b)^3/3$. If $b=z_M$, the same cell identity gives
\[
 \int_{t>b}\tau(t)^2\,\mathfrak m_n(\dd t)
 =\sum_{j=M+1}^n\Delta_j(1-z_{j-1})^2
 \le3\int_b^1(1-u)^2\,\dd u=(1-b)^3.
\]
Consequently,
\[
 \|P^+\mathcal R\|_{\mathfrak S_2}^2
       +\|\mathcal RP^+\|_{\mathfrak S_2}^2
 =2\int_{t>b}\int_X|R(s,t)|^2\,
                       \mathfrak m(\dd s)\mathfrak m(\dd t)
 \le C_\eta(1-b)^3.
\]
For $b=z_{n-1}$ the final atomic integral is exactly $\Delta_n^3$, since $\tau(z_n)=\Delta_n$. For $b=1$, both omitted operators vanish. This proves
\eqref{B:eq:truncation-upper-rows} at those boundary points as well.
\end{proof}

\subsubsection{Restriction of the triangular factor}
\label{B:subsection:restriction-of-the-triangular-factor}

Restricting the triangular factor and then forming its adjoint product is not the same as compressing the full inverse. The difference is a positive rank-one term, which must be
retained when the truncation is removed.

With the same projections, write
\begin{equation}
 \mathcal R^0_{d,b}=(P^0\mathcal KP^0)^*(P^0\mathcal KP^0),\qquad
 \mathcal T_{d,b}=P^0\mathcal K^*P^+\mathcal KP^0.
 \label{B:eq:restricted-factor-operators}
\end{equation}
For Lebesgue measure this is the restriction considered in Lemma~\ref{B:lem:bessel-restricted-factor}. The formulation here also includes finite kernels and their diagonal entries.
Let $q_b=\int_{(b,1]}r(u)^2\,\mathfrak m(\dd u)$, which is finite by local boundedness of $r$. All restricted operators are extended by zero on the complement of
$\operatorname{ran}P^0$.

\begin{lemma}\label{B:lem:factor-restriction-correction}
The correction is positive and has rank at most one:
\begin{equation}
 \mathcal T_{d,b}=q_b|P^0p\rangle\langle P^0p|,\qquad
 \|\mathcal T_{d,b}\|_{\mathfrak S_2}^2
       =q_b^2\left(\int_{(d,b]}p(s)^2\,\mathfrak m(\dd s)\right)^2
       \le C(1-b)^2.
 \label{B:eq:factor-restriction-rank-one}
\end{equation}
Moreover,
\begin{equation}
 P^0\mathcal RP^0=\mathcal R^0_{d,b}+\mathcal T_{d,b}.
 \label{B:eq:restriction-compression-identity}
\end{equation}
\end{lemma}
\begin{proof}
The functions $P^0p$ and $P^+r$ belong to $L^2(\mathfrak m)$: in the continuum their supports stay away from zero and the coefficients are locally bounded; on a finite partition
every sum is finite. For $f\in L^2(\mathfrak m)$, the triangular indicator gives
\[
 (P^+\mathcal KP^0f)(u)
 =\mathbf1_{\{u>b\}}r(u)\int_{(d,b]}p(s)f(s)\,\mathfrak m(\dd s).
\]
With the scalar product linear in its first argument, this is $\langle f,P^0p\rangle P^+r$. Thus, for $f,g\in L^2(\mathfrak m)$,
\[
 \langle\mathcal T_{d,b}f,g\rangle
 =\langle P^+\mathcal KP^0f,P^+\mathcal KP^0g\rangle
 =q_b\langle f,P^0p\rangle\overline{\langle g,P^0p\rangle}.
\]
This proves $\mathcal T_{d,b}=q_b|P^0p\rangle\langle P^0p|$ and positivity. For any orthonormal basis $(e_j)$, Parseval gives
\begin{align*}
 \|\mathcal T_{d,b}\|_{\mathfrak S_2}^2
 &=q_b^2\|P^0p\|_2^2\sum_j|\langle e_j,P^0p\rangle|^2\\
 &=q_b^2\|P^0p\|_2^4
 =q_b^2\left(\int_{(d,b]}p(s)^2\,\mathfrak m(\dd s)\right)^2.
\end{align*}
If $q_b=0$ or $P^0p=0$, these formulas give the zero operator. Otherwise its only nonzero eigenvalue is
$q_b\|P^0p\|_2^2$.

To estimate this quantity with the stated constants, write its kernel as
\[
 T_{d,b}(s,t)=\mathbf1_{(d,b]^2}(s,t)
              \int_{u>b}K(s,u)K(t,u)\,\mathfrak m(\dd u).
\]
Choose an admissible $\xi>0$ with $2\xi<a+1$. For $u>b\ge1/2$,
\[
 K(s,u)K(t,u)
 \le A_\xi^2s^{a/2-\xi}t^{a/2-\xi}u^{-a-1-2\xi}
 \le 2^{a+1+2\xi}A_\xi^2s^{a/2-\xi}t^{a/2-\xi}.
\]
For both measures $\mathfrak m((b,1])=1-b$; in the atomic case this is $\sum_{j>M}\Delta_j=1-z_M$. Integration therefore gives
\begin{equation}
 |T_{d,b}(s,t)|\le C_\xi(1-b)
     s^{a/2-\xi}t^{a/2-\xi}\mathbf1_{(d,b]^2}(s,t),
 \label{B:eq:factor-restriction-kernel-bound}
\end{equation}
Squaring and applying the first power estimate yields
\begin{align*}
 \|\mathcal T_{d,b}\|_{\mathfrak S_2}^2
 &\le C_\xi(1-b)^2
       \left(\int_{(d,b]}s^{a-2\xi}\,\mathfrak m(\dd s)\right)^2\\
 &\le C_\xi C_{a-2\xi}^2(1-b)^2b^{2(a+1-2\xi)}
 \le C_\xi(1-b)^2.
\end{align*}
Only the bounds for the products $K(s,u)K(t,u)$ enter this constant. In particular, replacing $p$ by $Cp$ and $r$ by $C^{-1}r$ changes $q_b$ by $C^{-2}$ and $\|P^0p\|_2^2$ by
$C^2$. Their product, and therefore the only possible nonzero eigenvalue of $\mathcal T_{d,b}$, is unchanged.

For $u\le d$ and $s>d$, the indicator $\mathbf1_{\{s\le u\}}$ is zero. Hence $P^-\mathcal KP^0=0$, also with the stated atomic boundary point conventions. Using the orthogonal
decomposition $I=P^-+P^0+P^+$ gives
\begin{align*}
 P^0\mathcal RP^0
 &=\sum_{\sigma\in\{-,0,+\}}
       (P^\sigma\mathcal KP^0)^*(P^\sigma\mathcal KP^0)\\
 &=(P^0\mathcal KP^0)^*(P^0\mathcal KP^0)
   +(P^+\mathcal KP^0)^*(P^+\mathcal KP^0)
 =\mathcal R^0_{d,b}+\mathcal T_{d,b}.
\end{align*}
This proves \eqref{B:eq:restriction-compression-identity}. If $b=1$, $P^+=0$, so the correction vanishes.
\end{proof}

For Lebesgue measure, the power coefficients $p(s)=s^{a/2}$ and $r(s)=s^{-(a+1)/2}$, with $a>0$, give
\[
 q_b=\frac{b^{-a}-1}{a},\qquad
 \|\mathcal T_{d,b}\|_{\mathfrak S_2}^2
 =\left(\frac{b^{-a}-1}{a}\right)^2
  \left(\frac{b^{a+1}-d^{a+1}}{a+1}\right)^2.
\]
For fixed $0<d<1$, differentiation of $b^{-a}$ at $b=1$ gives
\[
 \lim_{b\uparrow1}
 \frac{\|\mathcal T_{d,b}\|_{\mathfrak S_2}^2}{(1-b)^2}
 =\left(\frac{1-d^{a+1}}{a+1}\right)^2>0.
\]
The same limit is $(a+1)^{-2}$ when also $d\downarrow0$. Thus the correction contributes a quadratic boundary term in the full inverse comparison, in addition to the cubic bound
for the omitted rows.

\subsubsection{Removal of the truncation}
\label{B:subsection:removal-of-the-factor-truncation}

There are now two errors to combine: the kernel outside the retained square and the rank-one correction from restricting the factor. Their supports are disjoint, so the squared
Hilbert--Schmidt errors add. The hypotheses remain \eqref{B:eq:full-grid-measure} and \eqref{B:eq:triangular-factor-bound}, with the measures and cutoffs specified above.

\begin{theorem}\label{B:thm:inverse-truncation-bound}
For $0<\eta<\nu_a$,
\begin{equation}
 \|\mathcal R-\mathcal R^0_{d,b}\|_{\mathfrak S_2}^2
 \le C_\eta\big[d^{1+\nu_a-\eta}+(1-b)^3+(1-b)^2\big].
 \label{B:eq:inverse-truncation-bound}
\end{equation}
The constant is independent of $d,b$ and, for finite partitions, of $n$.
\end{theorem}
\begin{proof}
Set $P=P^0$ and $Q=I-P=P^-+P^+$. The restriction identity gives
\[
 \mathcal R-\mathcal R^0_{d,b}
 =(\mathcal R-P\mathcal RP)+\mathcal T_{d,b},\qquad
 \mathcal R-P\mathcal RP=Q\mathcal RQ+Q\mathcal RP+P\mathcal RQ.
\]
The three kernels on the right are supported, respectively, on $(X\setminus(d,b])^2$, $(d,b]\times(X\setminus(d,b])$, and $(X\setminus(d,b])\times(d,b]$, with the first kernel
variable integrated. These sets are disjoint. The correction is supported on $(d,b]^2$. Integration against the product measure therefore gives
\begin{equation}
 \|\mathcal R-\mathcal R^0_{d,b}\|_{\mathfrak S_2}^2
 =\|\mathcal R-P^0\mathcal RP^0\|_{\mathfrak S_2}^2
                          +\|\mathcal T_{d,b}\|_{\mathfrak S_2}^2.
 \label{B:eq:orthogonal-kernel-restrictions}
\end{equation}
Since $(Q\mathcal RP)^*=P\mathcal RQ$, the first term in this identity is exactly
\[
 \|\mathcal R-P\mathcal RP\|_{\mathfrak S_2}^2
 =\|Q\mathcal RQ\|_{\mathfrak S_2}^2
       +2\|Q\mathcal RP\|_{\mathfrak S_2}^2.
\]
Decomposing both copies of $L^2(\mathfrak m)$ into the ranges of $P$ and $Q$ also gives
\begin{align*}
 \|Q\mathcal R\|_{\mathfrak S_2}^2
       +\|\mathcal RQ\|_{\mathfrak S_2}^2
 &=2\|Q\mathcal RQ\|_{\mathfrak S_2}^2
       +2\|Q\mathcal RP\|_{\mathfrak S_2}^2\\
 &\ge\|\mathcal R-P\mathcal RP\|_{\mathfrak S_2}^2.
\end{align*}
Finally $P^-P^+=0$ gives the exact decompositions
\[
 \|Q\mathcal R\|_{\mathfrak S_2}^2
 =\|P^-\mathcal R\|_{\mathfrak S_2}^2+\|P^+\mathcal R\|_{\mathfrak S_2}^2,
 \qquad
 \|\mathcal RQ\|_{\mathfrak S_2}^2
 =\|\mathcal RP^-\|_{\mathfrak S_2}^2+\|\mathcal RP^+\|_{\mathfrak S_2}^2.
\]
Substituting Proposition~\ref{B:prop:truncation-row-bounds} and Lemma~\ref{B:lem:factor-restriction-correction} into \eqref{B:eq:orthogonal-kernel-restrictions} yields
\[
 \|\mathcal R-\mathcal R^0_{d,b}\|_{\mathfrak S_2}^2
 \le C_\eta d^{1+\nu_a-\eta}
          +C_\eta(1-b)^3+C_\eta(1-b)^2.
\]
The constants have the dependence stated before Proposition~\ref{B:prop:full-factor-and-green-kernel}, proving \eqref{B:eq:inverse-truncation-bound}.
\end{proof}

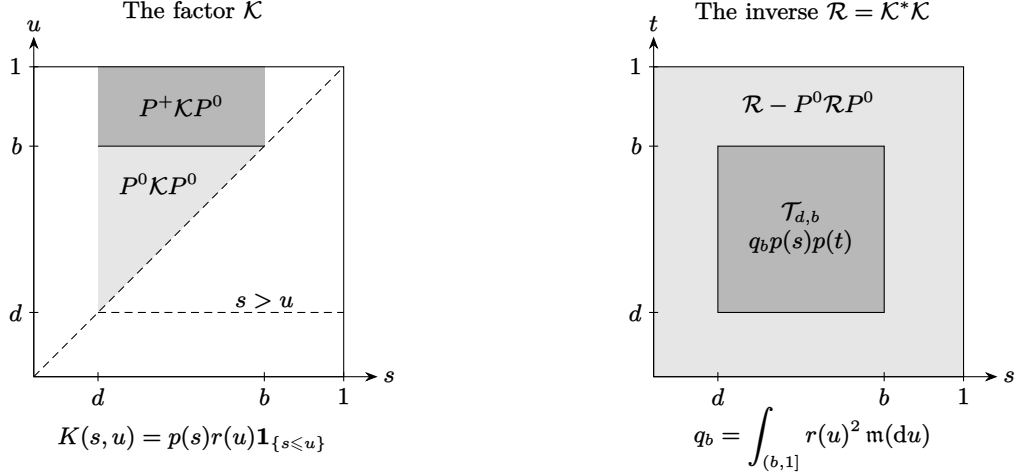
\begin{figure}[H]
\centering
\begin{tikzpicture}[x=1cm,y=1cm,every node/.style={font=\footnotesize,inner sep=2pt},>=Stealth]
 \node at (2.1,4.85) {The factor $\mathcal K$};
 \fill[black!10] (.85,.85)--(.85,3.05)--(3.05,3.05)--cycle;
 \fill[black!28] (.85,3.05) rectangle (3.05,4.1);
 \draw[->] (0,0)--(4.55,0) node[right] {$s$};
 \draw[->] (0,0)--(0,4.5) node[above] {$u$};
 \draw (0,0) rectangle (4.1,4.1);
 \draw[densely dashed] (0,0)--(4.1,4.1);
 \foreach \a/\b in {.85/d,3.05/b,4.1/1}{
   \draw (\a,-.06)--(\a,.06); \node[below] at (\a,-.08) {$\b$};
   \draw (-.06,\a)--(.06,\a); \node[left] at (-.08,\a) {$\b$};
 }
 \draw[densely dashed] (.85,.85)--(4.1,.85);
 \draw (.85,3.05)--(3.05,3.05);
 \node[align=center] at (1.65,2.55) {$P^0\mathcal KP^0$};
 \node at (1.95,3.58) {$P^+\mathcal KP^0$};
 \node at (3.05,1.0) {$s>u$};
 \node at (2.1,-.82) {$K(s,u)=p(s)r(u)\mathbf1_{\{s\le u\}}$};
 \begin{scope}[xshift=8.2cm]
  \node at (2.1,4.85) {The inverse $\mathcal R=\mathcal K^*\mathcal K$};
  \fill[black!10] (0,0) rectangle (4.1,4.1);
  \fill[white] (.85,.85) rectangle (3.05,3.05);
  \fill[black!28] (.85,.85) rectangle (3.05,3.05);
  \draw[->] (0,0)--(4.55,0) node[right] {$s$};
  \draw[->] (0,0)--(0,4.5) node[above] {$t$};
  \draw (0,0) rectangle (4.1,4.1);
  \draw (.85,.85) rectangle (3.05,3.05);
  \foreach \a/\b in {.85/d,3.05/b,4.1/1}{
   \draw (\a,-.06)--(\a,.06); \node[below] at (\a,-.08) {$\b$};
   \draw (-.06,\a)--(.06,\a); \node[left] at (-.08,\a) {$\b$};
  }
  \node[align=center] at (1.95,1.96) {$\mathcal T_{d,b}$\\$q_b p(s)p(t)$};
  \node at (2.05,3.62) {$\mathcal R-P^0\mathcal RP^0$};
  \node at (2.1,-.82) {$q_b=\displaystyle\int_{(b,1]}r(u)^2\,\mathfrak m(\dd u)$};
 \end{scope}
\end{tikzpicture}
\caption{Restriction of the triangular factor and its adjoint product. In the left panel, $s$ is integrated and the resulting function is evaluated at $u$. The darker rectangle gives $P^+\mathcal KP^0=|P^+r\rangle\langle P^0p|$ and hence the correction $\mathcal T_{d,b}=q_b|P^0p\rangle\langle P^0p|$ in \eqref{B:eq:restriction-compression-identity}. In the right panel, this correction occupies $(d,b]^2$, while $\mathcal R-P^0\mathcal RP^0$ is supported outside that square. The orthogonality in \eqref{B:eq:orthogonal-kernel-restrictions} follows from these supports. For the atomic measure, the included boundary points are those in \eqref{B:eq:restricted-factor-operators}.}
\label{B:fig:factor-restriction-regions}
\end{figure}

Figure~\ref{B:fig:factor-restriction-regions} separates the part lost when the triangular factor is restricted from the kernel outside $(d,b]^2$. Their disjoint supports explain
the sum of squared norms in \eqref{B:eq:orthogonal-kernel-restrictions}.

For a family of such factors let $r_n=\lceil(\log(n+2))^A\rceil$, where $A>0$ is fixed, and assume $d_n\le C r_n/n$, $1-b_n\le C r_n/n$. The constants for the partition and the
kernel factors are uniform over this family.

\begin{corollary}\label{B:cor:logarithmic-truncation-error}
For these cutoffs,
\begin{equation}
 \|\mathcal R_n-\mathcal R^0_{n,d_n,b_n}\|_{\mathfrak S_2}^2
                              =o(n^{-1}).
 \label{B:eq:logarithmic-truncation-error}
\end{equation}
\end{corollary}
\begin{proof}
Fix $0<\eta<\nu_a$. Since $r_n=O((\log n)^A)$, the three terms in Theorem~\ref{B:thm:inverse-truncation-bound}, multiplied by $n$, satisfy
\begin{equation}
 \begin{aligned}
 nd_n^{1+\nu_a-\eta}
 &\le Cr_n^{1+\nu_a-\eta}n^{-(\nu_a-\eta)}\longrightarrow0,\\
 n(1-b_n)^3&\le Cr_n^3/n^2\longrightarrow0,\qquad
 n(1-b_n)^2\le Cr_n^2/n\longrightarrow0.
 \end{aligned}
 \label{B:eq:truncation-orders}
\end{equation}
Each limit follows from $(\log n)^q n^{-v}\to0$ for fixed $q,v>0$. Their sum proves \eqref{B:eq:logarithmic-truncation-error}.
\end{proof}

For the partitions $z_{n,k}=k/n$ and $z_{n,k}=\theta_V(k/n)$, the bounds on the cell lengths give
\[
 d_n=z_{n,r_n}=\sum_{k=1}^{r_n}\Delta_{n,k}
       \le C_{\rm m}\frac{r_n}{n},\qquad
 1-b_n=1-z_{n,n-r_n}=\sum_{k=n-r_n+1}^n\Delta_{n,k}
       \le C_{\rm m}\frac{r_n}{n}.
\]
Since $r_n=O((\log n)^A)$, we have $r_n/n\to0$. Thus $1\le r_n<n/2$ and $b_n\ge1/2$ for all sufficiently large $n$, as required by the truncation theorem.

For Brownian or coupled coefficients, fix an outcome on which \eqref{B:eq:triangular-factor-bound} holds for the whole sequence with finite constants $A_\xi(\omega)$. For a fixed
$0<\eta<\nu_a$, the proof selects only finitely many positive losses $\xi$. On this outcome the theorem gives
\[
 n\|\mathcal R_n-\mathcal R^0_{n,d_n,b_n}\|_{\mathfrak S_2}^2
 \le C_{\omega,\eta}
 \left(r_n^{1+\nu_a-\eta}n^{-(\nu_a-\eta)}
                   +\frac{r_n^3}{n^2}+\frac{r_n^2}{n}\right)
 \longrightarrow0.
\]
This application uses finiteness of the factor constants on the chosen outcome. Scalar normalization in \eqref{B:eq:2-20} preserves $K(s,t)$, so it preserves all of these bounds.

\subsection{A criterion for convergence of the inverses at the hard edge}
\label{sec:a-criterion-for-convergence-of-the-inverses-at-the-hard-edge}
The deterministic estimates can now be collected into a criterion for the finite ensembles. It separates the weighted coefficient comparison on the selected interval from the
bounds needed to restore the omitted part of the inverse.

\subsubsection{A criterion for the finite ensembles}
\label{B:subsection:the-criterion-for-convergence-of-the-inverses}

Fix $\beta>0$, $a>0$, and $1<\lambda<\min\{a+1,2\}$. Use the model normalization $\alpha$ and the embeddings of Section~\ref{B:subsection:embedded-kernels-and-scalar-invariance}:
$\alpha=1$ on the partition $k/n$ and $\alpha=4/L_V^2$ on the partition $\theta_V(k/n)$. For a bidiagonal matrix $B_n$ with positive coefficients write
\begin{equation}
 \begin{split}
 \mathcal A_n&=\alpha^{-1}n^2B_nB_n^{\mathsf T},\qquad
 \mathcal R_n=\mathscr J_n\mathcal A_n^{-1}\mathscr J_n^*,\\
 \widetilde{\mathcal R}_n&=\mathcal R_n^\Delta,
 \qquad \widetilde{\mathcal K}_n=\mathcal K_n^\Delta,
 \qquad \widetilde{\mathcal R}_n
           =\widetilde{\mathcal K}_n^*\widetilde{\mathcal K}_n.
 \end{split}
 \label{B:eq:criterion-operators}
\end{equation}
The operators on the second line are those of \eqref{B:eq:2-12}--\eqref{B:eq:2-16} with this fixed normalization. The coefficients $p_{n,k},q_{n,k}$ remain those defined in
\eqref{B:eq:2-13} and \eqref{B:eq:2-9}.

Choose the intervals in \eqref{B:eq:brownian-retained-interval} with a fixed $A>12$, as used by the coefficient constructions below, and put
\[
 P_n^-=\mathbf1_{(0,d_n]},\qquad
 P_n^+=\mathbf1_{(\bar z_n,1]},\qquad
 \widetilde{\mathcal R}_n^0
  =(P_n^0\widetilde{\mathcal K}_nP_n^0)^*
                          (P_n^0\widetilde{\mathcal K}_nP_n^0).
\]
The operators $R_n^{B,\mathrm{step},0}$ and $R_n^{B,0}$ constructed from the Brownian path are the step and restricted kernels in \eqref{B:eq:2-34} and
\eqref{B:eq:bessel-restricted-green-kernel}, respectively, on these same cells. In formulas comparing them with \eqref{B:eq:criterion-operators}, we also write
\begin{equation}
 \mathcal R_n^{B,\mathrm{step},0}=R_n^{B,\mathrm{step},0},\qquad
 \mathcal R_n^{B,0}=R_n^{B,0},\qquad
 \mathcal R=R^B=\mathcal H_{\beta,a}^{-1}.
 \label{B:eq:criterion-brownian-cuts}
\end{equation}

Fix the complete array and Brownian motion on one probability space. The following hypotheses refer to those same random variables and the same deterministic cutoffs. In (i) the
observations determine each finite row. Conditions (ii) and (iii) estimate its difference from operators defined from the complete Brownian path.

\noindent\textup{(i)} For every $n$, $B_n$ has the coefficient law of the model. There is one sigma-field $\mathcal U$ independent of $\sigma\{W(u):u\ge0\}$ such that
\begin{equation}
 B_n\text{ is measurable with respect to }\mathcal N_n\vee\mathcal U,
 \label{B:eq:criterion-information}
\end{equation}
where $\mathcal N_n$ is defined in \eqref{B:eq:brownian-sampled-information}.

\noindent\textup{(ii)} There is a positive finite $\mathcal N_n\vee\mathcal U$-measurable scalar $C_n$. For the errors in \eqref{B:eq:2-35} with $c_n=\log C_n$ and weights
$w_{n,k}=w_{n,k}^{(\lambda)}$ from \eqref{B:eq:2-39}, almost surely,
\begin{equation}
 \mathfrak d_n:=\sum_{k\in\mathscr I_n^{\mathrm{ret}}}w_{n,k}
                 (|\delta_{n,k}|^2+|\rho_{n,k}|^2)=O(n^{-1}),
 \qquad
 \max_{k\in\mathscr I_n^{\mathrm{ret}}}
                 (|\delta_{n,k}|+|\rho_{n,k}|)\longrightarrow0.
 \label{B:eq:criterion-coefficient-errors}
\end{equation}

\noindent\textup{(iii)} For the same arrays, almost surely,
\begin{equation}
 \mathfrak e_n:=\|\widetilde{\mathcal R}_n
                  -\widetilde{\mathcal R}_n^0\|_{\mathfrak S_2}^2
                   =o(n^{-1}).
 \label{B:eq:criterion-truncation-error}
\end{equation}

To check (iii), one can estimate the removed rows and columns and the factor correction separately. Put
$\mathcal T_n=P_n^0\widetilde{\mathcal K}_n^*P_n^+ \widetilde{\mathcal K}_nP_n^0$ and
\[
 \mathfrak f_n=
 \sum_{\sigma\in\{-,+\}}
    \bigl(\|P_n^\sigma\widetilde{\mathcal R}_n\|_{\mathfrak S_2}^2
           +\|\widetilde{\mathcal R}_nP_n^\sigma\|_{\mathfrak S_2}^2\bigr)
                   +\|\mathcal T_n\|_{\mathfrak S_2}^2.
\]

\begin{lemma}\label{B:lem:finite-truncation-error-equivalence}
The exact error satisfies $\mathfrak e_n\le\mathfrak f_n\le2\mathfrak e_n$.
\end{lemma}
\begin{proof}
Suppress $n$ and write $P=P_n^0$, $Q=I-P=P_n^-+P_n^+$, and $\widetilde R=\widetilde{\mathcal R}_n$. Lemma~\ref{B:lem:factor-restriction-correction} gives
$P\widetilde RP=\widetilde R^0+\mathcal T$. For a direct coefficient formula, put $\mathscr I=\mathscr I_n^{\mathrm{ret}}$ and $T_i=\sum_{\ell=i}^nq_{n,\ell}$. On the full
partition, $\widetilde R(i,j)=p_{n,i}p_{n,j}T_{i\vee j}$, whereas the correction on the square of the restricted kernel has coefficient $\sum_{\ell>k_+}q_{n,\ell}$. Thus
\[
 \mathfrak e_n={}
 \sum_{(i,j)\notin\mathscr I^2}\Delta_{n,i}\Delta_{n,j}
           p_{n,i}^2p_{n,j}^2T_{i\vee j}^{\,2} 
 +\left(\sum_{\ell>k_+}q_{n,\ell}\right)^2
        \left(\sum_{i\in\mathscr I}\Delta_{n,i}p_{n,i}^2\right)^2.
\]
The first sum is over $1\le i,j\le n$. In particular it includes the first and last cells. Every $q_{n,\ell}$ is a finite matrix coefficient; no integral of $q^B$ over the first
cell occurs. To compare this exact expression with $\mathfrak f_n$, use the four terms in
\[
 \widetilde R-\widetilde R^0
 =Q\widetilde RQ+Q\widetilde RP+P\widetilde RQ+\mathcal T.
\]
They have disjoint kernel supports. Since $\widetilde R$ is self-adjoint, $\|Q\widetilde RP\|_{\mathfrak S_2}=\|P\widetilde RQ\|_{\mathfrak S_2}$. It follows that
\[
 \mathfrak e_n=\|Q\widetilde RQ\|_{\mathfrak S_2}^2
       +2\|Q\widetilde RP\|_{\mathfrak S_2}^2
       +\|\mathcal T\|_{\mathfrak S_2}^2.
\]
Orthogonality of $P_n^-$ and $P_n^+$ gives
\[
 \sum_{\sigma\in\{-,+\}}\|P_n^\sigma\widetilde R\|_{\mathfrak S_2}^2
 =\|Q\widetilde R\|_{\mathfrak S_2}^2
 =\|Q\widetilde RQ\|_{\mathfrak S_2}^2+\|Q\widetilde RP\|_{\mathfrak S_2}^2.
\]
Taking adjoints proves the same identity for the two column terms. Therefore
\[
 \mathfrak f_n=2\|Q\widetilde RQ\|_{\mathfrak S_2}^2
       +2\|Q\widetilde RP\|_{\mathfrak S_2}^2
       +\|\mathcal T\|_{\mathfrak S_2}^2.
\]
Subtracting the two formulas gives
\[
 \mathfrak f_n-\mathfrak e_n=\|Q\widetilde RQ\|_{\mathfrak S_2}^2\ge0,
 \qquad
 2\mathfrak e_n-\mathfrak f_n
 =2\|Q\widetilde RP\|_{\mathfrak S_2}^2+\|\mathcal T\|_{\mathfrak S_2}^2\ge0.
\]
These are the two required inequalities, including when any of the blocks or the correction vanishes.
\end{proof}

In particular, the bounds of Proposition~\ref{B:prop:truncation-row-bounds} and Lemma~\ref{B:lem:factor-restriction-correction} imply (iii) when the factor bound
\eqref{B:eq:triangular-factor-bound} holds uniformly for the constructed arrays. The three powers tend to zero after multiplication by $n$, as calculated in
\eqref{B:eq:truncation-orders}.

\begin{proposition}\label{B:prop:finite-inverse-upper-bound}
Under \textup{(ii)} and \textup{(iii)}, almost surely,
\begin{equation}
 \begin{gathered}
 \sup_{n\ge1}\|\widetilde{\mathcal R}_n\|_{\mathfrak S_2}<\infty,\\
 \|\mathcal R_n-\widetilde{\mathcal R}_n\|_{\mathfrak S_2}^2=O(n^{-2}),
 \qquad
 \|\mathcal R_n-\mathcal R\|_{\mathfrak S_2}^2=O(n^{-1}).
 \end{gathered}
 \label{B:eq:finite-inverse-upper-bound}
\end{equation}
\end{proposition}
\begin{proof}
Fix the common event in (ii), (iii), and the path growth and sampling estimates for this partition, these cutoffs, and $\lambda$. Choose $\eta$ as in
\eqref{B:eq:stability-exponents}. The coefficients $p^B,q^B$ satisfy \eqref{B:eq:reference-coefficient-bounds} with one finite pathwise constant, and (ii) gives a finite
$n_0(\omega)$ such that the maximum is at most one for $n\ge n_0(\omega)$. Scalar normalization preserves each inverse summand:
\[
 (C_n^{-1}p_{n,i})(C_n^{-1}p_{n,j})(C_n^2q_{n,l})
 =p_{n,i}p_{n,j}q_{n,l}.
\]
Thus Proposition~\ref{B:result:2-6} and Theorem~\ref{B:thm:weighted-inverse-stability} give
\[
 \|\widetilde{\mathcal R}_n^0-
        \mathcal R_n^{B,\mathrm{step},0}\|_{\mathfrak S_2}^2
 \le C_\omega\mathfrak d_n=O_\omega(n^{-1}).
\]
The assumptions on the partition are those already verified for the same selected cells. Theorem~\ref{B:thm:brownian-inverse-approximation} and
Proposition~\ref{B:prop:brownian-factor-truncation} give
\[
 \|\mathcal R_n^{B,\mathrm{step},0}-\mathcal R_n^{B,0}\|_{\mathfrak S_2}^2
 =O_\omega(n^{-1}),\qquad
 \|\mathcal R_n^{B,0}-\mathcal R\|_{\mathfrak S_2}^2=o_\omega(n^{-1}).
\]
The latter estimate uses \eqref{B:eq:truncation-orders}, including the quadratic correction from the columns outside the selected interval. Together with (iii), these estimates
apply to the exact identity
\begin{equation}
 \begin{split}
 \widetilde{\mathcal R}_n-\mathcal R={}&
 (\widetilde{\mathcal R}_n-\widetilde{\mathcal R}_n^0)
 +(\widetilde{\mathcal R}_n^0-\mathcal R_n^{B,\mathrm{step},0})\\
 &+(\mathcal R_n^{B,\mathrm{step},0}-\mathcal R_n^{B,0})
 +(\mathcal R_n^{B,0}-\mathcal R).
 \end{split}
 \label{B:eq:criterion-before-embedding}
\end{equation}
The inequality $\|\sum_{j=1}^4E_j\|_{\mathfrak S_2}^2 \le4\sum_{j=1}^4\|E_j\|_{\mathfrak S_2}^2$ yields
$\|\widetilde{\mathcal R}_n-\mathcal R\|_{\mathfrak S_2}^2 =O_\omega(n^{-1})$. Since $\mathcal R\in\mathfrak S_2$ and every earlier finite inverse exists, this proves
\begin{equation}
 M_\omega^{\mathrm{cell}}
 :=\sup_{n\ge1}\|\widetilde{\mathcal R}_n\|_{\mathfrak S_2}<\infty,
 \qquad \|\mathcal A_n^{-1}\|_{\mathfrak S_2}
 =\|\widetilde{\mathcal R}_n\|_{\mathfrak S_2}.
 \label{B:eq:criterion-full-bound}
\end{equation}
The equality uses the isometry $\widetilde{\mathscr J}_n$. Only after obtaining this bound do we apply Lemma~\ref{B:result:2-4}:
\[
 \|\mathcal R_n-\widetilde{\mathcal R}_n\|_{\mathfrak S_2}
 \le C_Vn^{-1}\|\mathcal A_n^{-1}\|_{\mathfrak S_2}
 \le C_VM_\omega^{\mathrm{cell}}n^{-1}.
\]
For the Laguerre partition the embeddings coincide, so this difference vanishes. The squared embedding error is $O_\omega(n^{-2})$, and the triangle inequality now proves
$\|\mathcal R_n-\mathcal R\|_{\mathfrak S_2}^2=O_\omega(n^{-1})$. Together these estimates prove \eqref{B:eq:finite-inverse-upper-bound}.
\end{proof}

\begin{proposition}\label{B:result:common-operator-criterion}
Under \textup{(i)--(iii)}, almost surely,
\begin{equation}
 0<\liminf_{n\to\infty}n\|\mathcal R_n-\mathcal R\|_{\mathfrak S_2}^2
 \le\sup_{n\ge1}n\|\mathcal R_n-\mathcal R\|_{\mathfrak S_2}^2<\infty.
 \label{B:eq:criterion-inverse-rate}
\end{equation}
On the same event, the embedded resolvents satisfy \eqref{B:eq:inverse-and-resolvent-rates} on every nonempty compact subset of $\rho(\mathcal H_{\beta,a})$.
\end{proposition}
\begin{proof}
Proposition~\ref{B:prop:finite-inverse-upper-bound} gives the upper bound for all sufficiently large $n$. Each earlier inverse is a finite matrix embedded by an isometry, and its
norm is finite. Including these finitely many terms gives the supremum in \eqref{B:eq:criterion-inverse-rate}.

For fixed $n$, the map $B\mapsto\alpha^{-1}n^2BB^{\mathsf T}$ is continuous on the space of real bidiagonal matrices with positive coefficients and takes values in the positive
definite matrices. For positive definite $A,A'$,
\[
 (A')^{-1}-A^{-1}=(A')^{-1}(A-A')A^{-1},\qquad
 \|(A')^{-1}-A^{-1}\|_{\mathfrak S_2}
 \le\|(A')^{-1}\|\,\|A-A'\|_{\mathfrak S_2}\|A^{-1}\|.
\]
The local Neumann series bounds $\|(A')^{-1}\|$ near each $A$, so inversion is continuous. For the deterministic isometry $\mathscr J_n$,
\[
 \|\mathscr J_nT\mathscr J_n^*-\mathscr J_nT'\mathscr J_n^*\|_{\mathfrak S_2}
 =\|T-T'\|_{\mathfrak S_2}.
\]
Composition proves that $\mathcal R_n$ is a Borel $\mathfrak S_2$-valued function of $B_n$. Assumption (i) consequently makes $\mathcal R_n$ measurable with respect to
$\mathcal N_n\vee\mathcal U$. The auxiliary field $\mathcal U$ is independent of the entire Brownian path. Theorem~\ref{B:thm:brownian-information-lower-bound} therefore applies
with $T_n=\mathcal R_n$ and $\mathcal U_n=\mathcal U$. It supplies the positive lower limit on an event for this array. Intersecting with the upper bound event proves
\eqref{B:eq:criterion-inverse-rate} for the same matrices.

The cell masses and scalar normalization used to estimate the error leave $\mathcal R_n=\mathscr J_n\mathcal A_n^{-1}\mathscr J_n^*$ as the function of $B_n$ just considered. The
finite coefficient laws in (i) are therefore the laws in the conclusion. Equations~\eqref{B:eq:embedded-resolvent-identification} and \eqref{B:eq:bessel-resolvent-identification}
identify its bounded transforms with the full embedded resolvents. Fix one outcome on which
\[
 Z=\sup_{n\ge1}n\|\mathcal R_n-\mathcal R\|_{\mathfrak S_2}^2<\infty,
 \qquad \ell=\liminf_{n\to\infty}n\|\mathcal R_n-\mathcal R\|_{\mathfrak S_2}^2>0.
\]
These are the inverse estimates already obtained for this same array. Put $M_* =\|\mathcal R\|+\sqrt Z$ and choose $N_\ell$ with
$n\|\mathcal R_n-\mathcal R\|_{\mathfrak S_2}^2\ge\ell/2$ for $n\ge N_\ell$.

Now choose any nonempty compact $K\subset\rho(\mathcal H_{\beta,a})$. Set
\[
 r_K=\sup_K|z|,\qquad B_K=\sup_K\|(I-z\mathcal R)^{-1}\|<\infty.
\]
For every $n$ with $r_KB_K\sqrt{Z/n}\le1/2$, the factorization
\[
 I-z\mathcal R_n=(I-z\mathcal R)
       [I-z(I-z\mathcal R)^{-1}(\mathcal R_n-\mathcal R)]
\]
gives inverses bounded by $2B_K$ uniformly on $K$. On the finite isometric range this is equivalent to $z\in\rho(\mathcal A_n)$; on its orthogonal complement the operator is the
identity. Let $G_n,G$ be the corresponding embedded and limiting resolvents. The exact identity gives
\[
 G_n(z)-G(z)=(I-z\mathcal R_n)^{-1}
                  (\mathcal R_n-\mathcal R)(I-z\mathcal R)^{-1}.
\]
It implies the upper bound $\|G_n-G\|_{\mathfrak S_2}\le2B_K^2\|\mathcal R_n-\mathcal R\|_{\mathfrak S_2}$. Multiplying back by the two denominators gives the lower bound
\[
 \|G_n(z)-G(z)\|_{\mathfrak S_2}
 \ge(1+r_KM_*)^{-2}\|\mathcal R_n-\mathcal R\|_{\mathfrak S_2}.
\]
Thus for every
\[
 n\ge\max\{N_\ell,1,\lceil4r_K^2B_K^2Z\rceil\}
\]
one has
\[
 \frac{\ell}{2(1+r_KM_*)^4n}
 \le\inf_{z\in K}\|G_n(z)-G(z)\|_{\mathfrak S_2}^2
 \le\sup_{z\in K}\|G_n(z)-G(z)\|_{\mathfrak S_2}^2
 \le\frac{4B_K^4Z}{n}.
\]
The compact set was chosen after the path was fixed. Consequently all eligible compact sets are covered by this deterministic argument on the original event, with their own
constants and indices; no intersection indexed by compact sets is taken.

\end{proof}

\begin{corollary}\label{B:cor:finite-inverse-resolvent-tightness}
Under \textup{(i)--(iii)}, the variables $n\|\mathcal R_n-\mathcal R\|_{\mathfrak S_2}^2$ are tight. For every fixed nonempty compact $K\subset\mathbb C\setminus(0,\infty)$,
\begin{equation}
 n\sup_{z\in K}
    \|\mathscr J_n(\mathcal A_n-z)^{-1}\mathscr J_n^*
                     -(\mathcal H_{\beta,a}-z)^{-1}\|_{\mathfrak S_2}^2
       =O_{\mathbb P}(1).
 \label{B:eq:finite-resolvent-tightness}
\end{equation}
\end{corollary}
\begin{proof}
Set $Z=\sup_{n\ge1}n\|\mathcal R_n-\mathcal R\|_{\mathfrak S_2}^2$. The maps are Borel by the preceding proof and Lemma~\ref{B:lem:bessel-path-measurability}, so $Z$ is measurable.
It is finite almost surely by \eqref{B:eq:criterion-inverse-rate}. For every $L>0$,
\[
 \sup_n\mathbb P(n\|\mathcal R_n-\mathcal R\|_{\mathfrak S_2}^2>L)
 \le\mathbb P(Z>L)\longrightarrow0\qquad(L\to\infty).
\]
The operator norms satisfy the measurable bound
\[
 M_\omega=\max\{\|\mathcal R\|,\sup_n\|\mathcal R_n\|\}
 \le\|\mathcal R\|+\sqrt Z<\infty.
\]
For the fixed $K\subset\mathbb C\setminus(0,\infty)$, positivity gives existence of every transform on $K$ for every $n$. The continuous function $M\mapsto d_K(M)$ is strictly
positive at each finite $M$. Lemma~\ref{B:lem:positive-inverse-resolvent-bounds} therefore gives
\begin{equation}
 \sup_{n\ge1}n\sup_{z\in K}
     \|G_n(z)-G(z)\|_{\mathfrak S_2}^2
 \le d_K(M_\omega)^{-4}Z<\infty\quad\text{almost surely}.
 \label{B:eq:criterion-resolvent-envelope}
\end{equation}
Continuity in $z$ reduces each supremum over $K$ to a countable dense subset, as proved in Proposition~\ref{B:prop:resolvent-tightness}. The right hand side is a finite measurable
random variable. Hence
\[
 \sup_n\mathbb P\left(n\sup_{z\in K}
          \|G_n(z)-G(z)\|_{\mathfrak S_2}^2>L\right)
 \le\mathbb P\bigl(d_K(M_\omega)^{-4}Z>L\bigr)\longrightarrow0.
\]
The identifications \eqref{B:eq:embedded-resolvent-identification}--\eqref{B:eq:bessel-resolvent-identification} turn this into \eqref{B:eq:finite-resolvent-tightness}. The proof
uses the tails of finite random bounds throughout.
\end{proof}

The upper estimate uses the weighted sum and maximum in (ii), together with the full inverse truncation in (iii). The law and observation condition (i) enter when the lower bound
is applied to that same array. The coefficient constructions below verify these conditions separately. Condition (i) remains an assumption on the law and observations, rather than
on rank. The resulting finite random supremum is a pathwise bound, without a bound on its moments.

\subsection{Reduction to the error from sampling the Brownian path}
\label{sec:the-circular-comparison-scheme}
\label{S:sec:proof-strategy}
For the circular ensemble, we compute the leading coefficient of the squared inverse error. The error in coupling the circular ensemble operator to the operator formed from
Brownian samples is smaller than the error from sampling the Brownian path, so the latter determines the coefficient. Its conditional mean and conditional covariance both
contribute.

Fix $\beta>0$. All operators are conjugated by the matrix $\mathcal A$ defined in the introduction, with the boundary vectors $u_0,u_1$ fixed there. Put
\[
 \mathcal H=L^2([0,1);\mathbb R^2),\qquad
 \mathscr K=\mathfrak S_2(\mathcal H),\qquad
 a_\beta=\min\{1,2/\beta\}.
\]
Norms and inner products below are those of the real Hilbert space $\mathscr K$, unless another space is indicated. Conditional on the sampled positions, this error has a Gaussian
first variation and a mean determined by those positions. Their conditional second moment gives the coefficient, while the covariance bound controls fluctuations around it.
Section~\ref{S:sec:proof-of-main-theorems} combines these estimates with the operator expansion below.

\subsubsection{The finite comparison and the boundary point}
\label{S:sec:reduction-to-the-brownian-sampling-error}
Proposition~\ref{S:prop:7.2} constructs the Brownian path $Z$ and all circular walks on one probability space. Each walk has the circular walk law and terminal point $q_\infty$.
The enlargement in Theorem~\ref{S:thm:joint-probability-space} adds the Gaussian bridges used below without changing the joint law of $Z$ and the entire circular array.

Recall the two differences
\[
 F_n^{\mathrm{full}}=\mathtt r\mathrm{Circ}_{\beta,n}
                         -\mathtt r\mathrm{Br}_{\beta,n},
 \qquad
 \varepsilon_{\beta,n}^{\mathrm{samp}}
       =\mathtt r\mathrm{Br}_{\beta,n}-\mathtt r\mathrm{Sine}_\beta.
\]
Their sum is $\Delta_{\beta,n}$. Proposition~\ref{S:prop:circular-brownian-inverse-comparison} gives
\begin{equation}\label{S:eq:s2-finite-full}
 d_n:=\sqrt n\,\|F_n^{\mathrm{full}}\|
       =O_{\beta,\omega}((\log n)^{-501/2})\longrightarrow0.
\end{equation}
In particular,
\begin{equation}\label{S:eq:s2-error-reduction}
 n\bigl|\|\Delta_{\beta,n}\|^2
              -\|\varepsilon_{\beta,n}^{\mathrm{samp}}\|^2\bigr|
 \le 2d_n\sqrt n\,\|\varepsilon_{\beta,n}^{\mathrm{samp}}\|+d_n^2.
\end{equation}
Thus a finite limit for the normalized squared sampling error will also be the limit for the circular error. The mixed term tends to zero because \eqref{S:eq:s2-finite-full} gives
$d_n\to0$.

The expansion into a conditionally Gaussian term and a remainder is made away from the singular boundary. Use the cutoff of \eqref{S:eq:cutoff-choice},
\[
 m_n=\lceil(\log n)^{512}\rceil,\qquad T_n=1-m_n/n,
\]
and restrict to sufficiently large $n$ that $m_n<n$. Let $P_n$ be multiplication by $\mathbf1_{[0,T_n)}$ on $\mathcal H$ and set $\Pi_n D=P_nDP_n$ for $D\in\mathscr K$. If $D$ has
kernel $K_D$, then $\Pi_nD$ has kernel $\mathbf1_{[0,T_n)^2}K_D$. Consequently $\Pi_n$ is an orthogonal projection on $\mathscr K$, and
\begin{equation}\label{S:eq:s2-orthogonal-restriction}
 \|D\|^2=\|\Pi_nD\|^2+\|D-\Pi_nD\|^2.
\end{equation}
The complementary support consists of the square $[T_n,1)^2$ and the two rectangles with exactly one coordinate beyond $T_n$. By Theorem~\ref{S:thm:singular-regular-truncations},
these compressions equal the regular inverses used in the expansion, extended by zero with the same coordinate boundary vectors.

Theorem~\ref{S:thm:full-boundary-tails} proves, for every $0<\epsilon<a_\beta$,
\begin{equation}\label{S:eq:s2-all-tails}
 \mathfrak t_n^{\mathrm{end}}
 :=n\!\!\sum_{K\in\{\mathtt r\mathrm{Circ}_{\beta,n},
                    \mathtt r\mathrm{Br}_{\beta,n},
                    \mathtt r\mathrm{Sine}_\beta\}}
           \|K-\Pi_nK\|^2
 =O_{\beta,\omega,\epsilon}(n^{-a_\beta+\epsilon})\longrightarrow0.
\end{equation}
The inequality $\|A-B\|^2\le2\|A\|^2+2\|B\|^2$ gives
\begin{equation}\label{S:eq:s2-sampling-tail}
 n\|\varepsilon_{\beta,n}^{\mathrm{samp}}
             -\Pi_n\varepsilon_{\beta,n}^{\mathrm{samp}}\|^2
 \le 2\mathfrak t_n^{\mathrm{end}}\longrightarrow0.
\end{equation}
Equations \eqref{S:eq:s2-orthogonal-restriction} and \eqref{S:eq:s2-sampling-tail} reduce the asymptotic for its squared norm to the square of the restricted kernel. The disjoint
supports make the two squared kernel norms additive.

\subsubsection{The conditional expansion}
\label{S:sec:the-brownian-path-between-samples}
Conditional on the Brownian samples at a fixed matrix size and on the terminal point, the path segments between consecutive samples have hyperbolic Brownian bridge distributions.
The resulting expansion keeps the part determined by the samples separate from the centered Gaussian part; this distinction is needed when the squared norm is computed.

At each matrix size, the conditioning sigma-field is
\[
 \mathscr S_n=\sigma\{q_\infty,Z_{\vartheta_{n,j}}:1\le j\le n\}.
\]
This determines the sampled matrices, their boundary columns, and the derivative of the Green kernel. All conditional distributions and estimates below are taken with respect to
this sigma-field. The regular conditional laws and the enlargement preserving the original joint distribution are constructed in
Section~\ref{S:sec:conditioning-fields-and-simultaneous-versions}.

The matrix expansion in Proposition~\ref{S:prop:bridge-affine-expansion} separates the relative logarithm into a term $\mathcal B_n$ determined by the samples, a conditionally
centered Gaussian field $\mathcal Z_n$, and a smaller field. Let $\mathscr D_n$ be the first variation at the sampled coefficient, defined by
\eqref{S:eq:continuous-derivative-kernel}, and put
\[
 G_n=-\mathscr D_n\mathcal Z_n,\qquad
 B_n=-\mathscr D_n\mathcal B_n.
\]
The sign is fixed by sampled inverse minus continuous inverse. The field $\mathcal B_n$ is determined by the sampled positions, whereas $\mathcal Z_n$ is centered under their
conditional law. The conditional second moment includes both $B_n$ and $G_n$, whose weighted sizes are of the same order.

Theorems~\ref{S:thm:gaussian-finite-mode} and~\ref{S:thm:15.1} give the exact identity
\begin{equation}\label{S:eq:s2-operator-expansion}
 \Pi_n\varepsilon_{\beta,n}^{\mathrm{samp}}
       =G_n+B_n+R_n^{\mathrm{br}},\qquad
 G_n\mid\mathscr S_n\sim\mathcal N_{\mathscr K}(0,Q_n),
\end{equation}
where $Q_n$ is positive and trace class, and $Q_n,B_n$ are $\mathscr S_n$-measurable. All three terms are supported on the square of the restricted kernel. The Gaussian term is a
measurable limit of finite sums in $\mathscr K$, and its covariance is their limit in trace norm. The sampled chord determines $B_n$; measurability of the kernel maps gives the
measurability of $B_n$ and $R_n^{\mathrm{br}}$. The remainder satisfies
\begin{equation}\label{S:eq:s2-operator-remainder}
 u_n:=\sqrt n\,\|R_n^{\mathrm{br}}\|
       =O_{\beta,\omega}((\log n)^{-250})\longrightarrow0.
\end{equation}
There are two errors to control: the error in the matrix field and the nonlinear remainder of the kernel. For the sampled coefficient $X$ and the weighted quadratic form
$\mathcal E_n^+$, Theorem~\ref{S:thm:continuous-frechet} gives, for $U$ in its weighted field space and $\eta(U)=\operatorname{ess\,sup}\|U\|_{\mathrm{op}}\le1$,
\[
 n\|\mathsf K_X(U)-\mathsf K_X(0)-\mathscr D_{n,X}U\|_{\mathfrak S_2}^2
 \le C\eta(U)^2\mathcal E_n^+(U).
\]
The weights in $\mathcal E_n^+$ come from integration against the two boundary columns. A uniform bound on the growing column is therefore unnecessary. The field estimates make the
right side tend to zero and bound the linear image of the error in approximating the matrix field. Both contributions are included in \eqref{S:eq:s2-operator-remainder}.

\subsubsection{The conditional second moment}
\label{S:sec:concentration-of-the-gaussian-energy}
Squaring the leading expansion produces a conditional expectation and a centered fluctuation. The expectation includes the squared conditional mean as well as the covariance trace,
so neither part can be discarded in the leading order coefficient.

Set
\begin{equation}\label{S:eq:final-mean-fluctuation}
 \mathfrak C_n=n\operatorname{tr}Q_n+n\|B_n\|^2,\qquad
 \mathfrak s_n=n\|G_n+B_n\|^2-\mathfrak C_n.
\end{equation}
Fix a conditional law given $\mathscr S_n$ for which $Q_n$ is positive trace class and $B_n\in\mathscr K$. In this paragraph all coefficients are deterministic under that law.
Write $Q=Q_n$, $b=B_n$ and $G=G_n$. A spectral representation of a centered Gaussian with this covariance is
\[
 G=\sum_{j\ge1}\sqrt{\lambda_j}\,\xi_j e_j,
 \qquad \lambda_j\ge0,\quad\sum_j\lambda_j<\infty,
\]
where the $\xi_j$ are independent standard real Gaussians. Include an orthonormal basis of $\ker Q$ when expanding $b$, and put $b_j=\langle b,e_j\rangle$. Parseval and monotone
convergence give $\mathbb E\|G\|^2=\sum_j\lambda_j=\operatorname{tr}Q$. The series for $\langle G,b\rangle$ converges in $L^2$, because
\[
 \sum_j\lambda_j b_j^2\le\|Q\|\|b\|^2<\infty.
\]
Likewise $\sum_j\lambda_j(\xi_j^2-1)$ converges in $L^2$, with variance $2\sum_j\lambda_j^2$. For every finite partial sum, independence and $\mathbb E\xi_j^3=0$ make the
covariance between its quadratic and linear parts zero. Passing in $L^2$ to the limit yields
\[
 \begin{split}
 \mathbb E\|G+b\|^2&=\operatorname{tr}Q+\|b\|^2,\\
 \operatorname{Var}(\|G+b\|^2)
 &=2\operatorname{tr}(Q^2)+4\langle Qb,b\rangle.
 \end{split}
\]
Indeed the centered random variable is the sum of the quadratic series and $2\langle G,b\rangle$. The spectral inequalities
\[
 \operatorname{tr}(Q^2)\le\|Q\|\operatorname{tr}Q,
 \qquad \langle Qb,b\rangle\le\|Q\|\|b\|^2
\]
then give exactly the stated conditional variance bound after multiplication by $n^2$.

This calculation is made for each regular conditional distribution with finite $\operatorname{tr}Q+\|b\|^2$. It remains valid when these conditional moments have infinite
unconditional expectation. The squared norm of the conditional mean $b$ contributes to the leading conditional second moment.

Corollary~\ref{S:cor:coefficient-conditional-mean} proves
\begin{equation}\label{S:eq:s2-mean-limit}
 \mathfrak C_n\longrightarrow\PathCoeff(\omega)
   =\frac4\beta\int_0^1\kappa_\beta(r)
                  \frac{\Gamma_{\beta,\omega}(r)}r\,\dd r
   \in(0,\infty)\qquad\text{almost surely}.
\end{equation}
The conditional variance and squared conditional mean produce separate terms in this calculation. Their sum gives the Riemann sum for \eqref{S:eq:s2-mean-limit}. Terms involving
distinct intervals also occur; Theorem~\ref{S:thm:trace-conditional-complete} estimates them before the limit is taken. Lemma~\ref{S:lem:weighted-qv} applies pathwise to weights
depending on the whole limiting path. Proposition~\ref{S:prop:coefficient-positive-measurable} proves finiteness, strict positivity, and measurability of the integral from $Z$
alone.

The preceding calculation gives the following bound for the centered fluctuation:

\begin{equation}\label{S:eq:s2-conditional-variance}
 n^2\operatorname{Var}(\|G_n+B_n\|^2\mid\mathscr S_n)
 \le 2n^2\|Q_n\|_{\mathrm{op}}
                  \bigl(\operatorname{tr}Q_n+2\|B_n\|^2\bigr).
\end{equation}
Proposition~\ref{S:prop:covariance-upper-limit} and \eqref{S:eq:gaussian-trace-bias-mean}, with the field bounds used there, supply
\begin{equation}\label{S:eq:s2-covariance-rate}
 \operatorname{tr}Q_n+\|B_n\|^2=O_{\beta,\omega}(n^{-1}),\qquad
 \|Q_n\|_{\mathrm{op}}
       =O_{\beta,\omega,\epsilon}(n^{-1-a_\beta+\epsilon}).
\end{equation}
For $0<\epsilon<a_\beta$, the right side of \eqref{S:eq:s2-conditional-variance} is $O_{\beta,\omega,\epsilon}(n^{-a_\beta+\epsilon})$. The bound for $\|Q_n\|_{\mathrm{op}}$
therefore makes the conditional variance of the normalized squared error tend to zero. This deduction uses only covariance upper bounds.

We derive the concentration estimate with deterministic constants under each conditional law. With $Q,b,G$ as above set $v=\operatorname{tr}(Q^2)$ and $a=\|Q\|$. For
$0<t<(2a)^{-1}$, the scalar Gaussian integral gives
\[
 \log\mathbb E\exp\{t\lambda(\xi^2-1)\}
 =-t\lambda-\tfrac12\log(1-2t\lambda)
 \le\frac{t^2\lambda^2}{1-2ta}.
\]
The inequality follows by expanding the logarithm and using $2^{k-1}/k\le2^{k-2}$ for $k\ge2$. Independence, followed by passage from finite sums using Fatou's lemma, gives the
same upper bound $t^2v/(1-2ta)$ for $\|G\|^2-\operatorname{tr}Q$. For the negative parameter the bound is $t^2v$, since $-\log(1+x)+x\le x^2/2$ for $x\ge0$. Exponential Markov,
with $t=\sqrt x/(\sqrt v+2a\sqrt x)$ for the upper tail, yields
\[
 \mathbb P\{\|G\|^2-\operatorname{tr}Q>2\sqrt{vx}+2ax\}\le e^{-x}.
\]
Optimizing the lower-tail quadratic bound gives
$\mathbb P\{\operatorname{tr}Q-\|G\|^2>2\sqrt{vx}\}\le e^{-x}$.
When the covariance vanishes, the Gaussian random element is deterministic and both inequalities hold. Since $2\langle G,b\rangle$ is a real Gaussian of variance
$4\langle Qb,b\rangle$,
\[
 \mathbb P\{|2\langle G,b\rangle|>2\sqrt{2\langle Qb,b\rangle x}\}\le2e^{-x}.
\]
A union bound consequently controls the absolute centered noncentral square by
\[
 2\sqrt{vx}+2ax+2\sqrt{2\langle Qb,b\rangle x}
\]
with exceptional probability at most $4e^{-x}$.

Let $0<\kappa<a_\beta$ and localize on the $\mathscr S_n$-measurable event
\[
 E_{L,n}=\{\operatorname{tr}Q_n+\|B_n\|^2\le L/n,
                    \|Q_n\|\le L n^{-1-\kappa}\},\qquad L\in\mathbb N.
\]
On this event $v\le L^2n^{-2-\kappa}$ and $\langle Q_nb,b\rangle\le L^2n^{-2-\kappa}$. Conditional integration gives
\[
 \mathbb P\{E_{L,n},\ |\mathfrak s_n|>
 12L(n^{-\kappa/2}\sqrt x+n^{-\kappa}x)\}\le4e^{-x}.
\]
The indicator of $E_{L,n}$ can be taken outside the conditional probability precisely because this event is measurable with respect to $\mathscr S_n$.

Put $x=8\log n$ for $n\ge2$. For each fixed integer $L$ the exceptional probabilities sum to at most $4\sum_{n\ge2}n^{-8}<\infty$. Intersect the resulting probability-one events
over $L$. The pathwise covariance estimates supply, on this intersection, an integer $L=L(\omega)$ and an index $N_0(\omega)$ for which $E_{L,n}$ holds whenever $n\ge N_0$.
Borel--Cantelli for this fixed integer $L$ then gives
\[
 |\mathfrak s_n|\le C L(\omega)
       (n^{-\kappa/2}\sqrt{\log n}+n^{-\kappa}\log n)
\]
for every $n$ above another finite index. To obtain the stated bound with any $0<\varepsilon<a_\beta/2$, choose $\kappa=a_\beta-\varepsilon$ and use
$\sqrt{\log n}=o(n^{\varepsilon/2})$. This proves $\mathfrak s_n=O_\omega(n^{-a_\beta/2+\varepsilon})$. The argument uses the conditional law at each fixed size and summability
over the sizes. Independence between different sizes is not required.

In particular,

\begin{equation}\label{S:eq:s2-concentration-rate}
 |\mathfrak s_n|
   =O_{\beta,\omega,\epsilon}(n^{-a_\beta/2+\epsilon})
      \longrightarrow0,\qquad 0<\epsilon<a_\beta/2.
\end{equation}
 The same theorem gives a bound simultaneous over all Gaussian truncation orders; only the full sum is needed here.

Combining the conditional second moment limit with the concentration estimate gives
\[
 n\|G_n+B_n\|^2=\mathfrak C_n+\mathfrak s_n
           \longrightarrow\PathCoeff(\omega)<\infty.
\]
In particular $v_n=\sqrt n\|G_n+B_n\|$ is bounded on the chosen path. For $u_n=\sqrt n\|R_n^{\mathrm{br}}\|\to0$, the exact operator expansion gives
\[
 n\left|\|\Pi_n\varepsilon_{\beta,n}^{\mathrm{samp}}\|^2
                        -\|G_n+B_n\|^2\right|
 \le2u_nv_n+u_n^2\longrightarrow0.
\]
The kernel restricted to the square and the kernel on its complement are orthogonal in $\mathscr K$. Therefore
\[
 n\|\varepsilon_{\beta,n}^{\mathrm{samp}}\|^2
 =n\|\Pi_n\varepsilon_{\beta,n}^{\mathrm{samp}}\|^2
  +n\|(I-\Pi_n)\varepsilon_{\beta,n}^{\mathrm{samp}}\|^2
 \longrightarrow\PathCoeff(\omega),
\]
because the last summand is at most $2\mathfrak t_n^{\mathrm{end}}\to0$. This proves the boundedness of $\sqrt n\|\varepsilon_{\beta,n}^{\mathrm{samp}}\|$ required for the final
comparison. Now $d_n=\sqrt n\|F_n^{\mathrm{full}}\|\to0$ implies
\[
 n\left|\|\Delta_{\beta,n}\|^2
                 -\|\varepsilon_{\beta,n}^{\mathrm{samp}}\|^2\right|
 \le2d_n\sqrt n\|\varepsilon_{\beta,n}^{\mathrm{samp}}\|+d_n^2\longrightarrow0.
\]
Hence the same coefficient is obtained for the full circular error. The finite comparison, boundary estimate, Gaussian expansion, covariance estimate, and conditional coefficient
limit used here are stated and proved separately in the later sections.

The conclusion concerns only the original Brownian path and circular array. On an enlargement preserving their joint distribution, the event expressing this limit is the inverse
image of the corresponding event on that original distribution. Probability one on the enlargement therefore gives probability one for the original pair, without keeping the
auxiliary Gaussian bridge variables in the conclusion.

\subsection{Green kernels for canonical systems}
\label{sec:green-kernels-for-canonical-systems}
\label{S:sec:green-kernels-for-canonical-systems}
\label{S:sec:external-aux}

At spectral parameter zero the canonical equation has constant solutions. The Green kernel can therefore be written in terms of the two boundary vectors and the coefficient factor
$X$. Since the formula uses only values of $X$, it applies to both the Brownian and step coefficients. Absolute continuity is imposed on $f$ in canonical coordinates.

We first derive this formula on a regular interval. At the singular boundary, the order of the two kernel variables determines the weighted double integral for square
integrability. The formula also identifies the self-adjoint domain and the boundary conditions corresponding to a restriction of the kernel. Perturbations of $X$ are treated in
Section~\ref{S:sec:stability-and-fr-echet-expansion-of-the-inverse-kernel}.

\subsubsection{Canonical systems on a regular interval}
\label{S:sec:canonical-operators-and-their-green-kernels}

The boundary vectors determine the Green kernel through their antisymmetric pairing. With $J$ as in the introduction, write $[v,w]=v^{\mathsf T}Jw$, a nondegenerate antisymmetric
bilinear form on real vectors, extended complex bilinearly. The normalization $[u_0,u_1]=1$ makes $u_0,u_1$ linearly independent and identifies $[v,u_i]=0$ with
$v\in\mathbb C u_i$.

\begin{lemma}\label{S:lem:3.1}
The normalized boundary vectors satisfy
\begin{equation}\label{S:eq:3.2}
 J=u_0u_1^{\mathsf T}-u_1u_0^{\mathsf T}.
\end{equation}
\end{lemma}
\begin{proof}
Write $v=\alpha u_0+\nu u_1$. Pairing with $u_1$ and $u_0$ gives $\alpha=[v,u_1]$ and $\nu=-[v,u_0]$. Thus
\begin{equation}\label{S:eq:3.4}
 v=[v,u_1]u_0-[v,u_0]u_1.
\end{equation}
Since $J^{\mathsf T}J=I$,
\begin{equation}\label{S:eq:3.5}
 [Jw,u]=w^{\mathsf T}u.
\end{equation}
Substitution of $v=Jw$ in \eqref{S:eq:3.4} proves $Jw=(u_0u_1^{\mathsf T}-u_1u_0^{\mathsf T})w$ for every $w$.
\end{proof}

Let $0<T\le1$, and suppose that $X:[0,T]\to\mathrm{SL}_2(\mathbb R)$ is measurable, with $X,X^{-1}$ essentially bounded. Put $R=R_X=X^{\mathsf T}X$ and
\[
 \mathcal H_R(0,T)=L^2((0,T),R(t)\,\dd t;\mathbb C^2),\qquad
 \langle f,g\rangle_R=\int_0^T f(t)^*R(t)g(t)\,\dd t.
\]
The inner product is linear in its second argument. There are constants $0<c\le C<\infty$ such that $cI\le R(t)\le CI$ almost everywhere. Consequently the weighted and ordinary
$L^2$ norms are equivalent. The map
\[
 \mathcal U_X:\mathcal H_R(0,T)\longrightarrow L^2((0,T);\mathbb C^2),
 \qquad \mathcal U_Xf=Xf,
\]
is unitary: its inverse is multiplication by $X^{-1}$, and $\|Xf\|_2^2=\int_0^T f^*Rf=\|f\|_R^2$. All coefficients and boundary vectors are real. The operators and inverses
constructed below therefore preserve real functions; their restrictions are the real operators used in the probabilistic estimates. A real orthonormal basis is also an orthonormal
basis of the complexification, so the real and complex Hilbert--Schmidt norms of a real kernel agree.

On the regular interval $[0,T]$, define $\mathcal D_X^{u_0,u_1}f=2R^{-1}Jf'$ on the domain
\[
 \mathcal D(\mathcal D_X^{u_0,u_1})
 =\{f\in H^1((0,T);\mathbb C^2):
             [f(0),u_0]=[f(T),u_1]=0\}.
\]
Equivalently, impose the same boundary conditions and require $f$ to be absolutely continuous, $f\in\mathcal H_R$, and $R^{-1}Jf'\in\mathcal H_R$: indeed $f'=-JR(\mathcal D_Xf)/2$,
and $R,R^{-1}$ are bounded. The boundary conditions use the Sobolev traces of $f$. Set
\[
 a_X=Xu_0,\qquad c_X=Xu_1,\qquad
 p=|a_X|^2,\qquad q=|c_X|^2.
\]
When the matrix function is fixed we omit the subscript on $a_X,c_X$. The columns $a,c$ are the conjugated solutions at spectral parameter zero, since the solutions before
conjugation are constant vectors. Define their integrals on the two sides of $t$ by
\[
 P_X(t)=\int_0^t p(s)\,\dd s,\qquad Q_X(t)=\int_t^T q(s)\,\dd s.
\]
The kernel estimates use these $t$-dependent integrals, which have different behavior at the two boundary points.

\begin{theorem}\label{S:thm:3.2}
The regular operator $\mathcal D_X^{u_0,u_1}$ is self-adjoint and $0\in\rho(\mathcal D_X^{u_0,u_1})$. Its conjugated inverse has kernel
\begin{equation}\label{S:eq:3.23}
 K_X(s,t)=\frac12\left[
 a(s)c(t)^{\mathsf T}\mathbf1_{\{s<t\}}
 +c(s)a(t)^{\mathsf T}\mathbf1_{\{s>t\}}\right].
\end{equation}
This inverse is Hilbert--Schmidt.
\end{theorem}
\begin{proof}
Let $g\in\mathcal H_R(0,T)$. Uniform upper and lower bounds on $R$ make this the same vector space as ordinary $L^2$, with equivalent norms. Define
\begin{equation}
 f(s)=\frac12u_0\int_s^T u_1^{\mathsf T}R(t)g(t)\,\dd t
       +\frac12u_1\int_0^s u_0^{\mathsf T}R(t)g(t)\,\dd t.
 \label{S:eq:s3-regular-solution}
\end{equation}
Both integrands lie in $L^2(0,T)$ and hence in $L^1(0,T)$. The function $f$ is absolutely continuous and
\[
 f'(s)=\tfrac12(-u_0u_1^{\mathsf T}+u_1u_0^{\mathsf T})R(s)g(s)
       =-\tfrac12JR(s)g(s).
\]
This derivative is in $L^2$; the function is also in $L^2$, since its two scalar integrals are bounded by constants times $\|g\|_R$. Thus $f\in H^1(0,T)$. At $s=0$ it lies in
$\mathbb C u_0$ and at $s=T$ it lies in $\mathbb C u_1$. The boundary conditions hold, and $2R^{-1}Jf'=g$ almost everywhere.

Conversely, if $f$ belongs to the stated differential domain and $g=2R^{-1}Jf'$, subtract the function defined above. The difference has zero derivative and is therefore a constant
vector $v$. Its left boundary condition gives $v\in\mathbb C u_0$, and its right condition gives $v\in\mathbb C u_1$. These lines intersect only at zero, because $[u_0,u_1]=1$.
Thus the formula is a two-sided inverse on the stated domain.

Conjugation by the unitary $\mathcal U_X$ turns the formula into the stated kernel. For any real vectors $v,w$, $\|vw^{\mathsf T}\|_{\mathrm{HS}}^2=|v|^2|w|^2$. The two triangles
therefore give
\[
 \|\mathsf K_X\|_{\mathfrak S_2}^2
 =\frac14\int_{s<t}p(s)q(t)\,\dd s\,\dd t
  +\frac14\int_{s>t}q(s)p(t)\,\dd s\,\dd t
 =\frac12\int_0^Tq(t)\int_0^tp(s)\,\dd s\,\dd t<\infty.
\]
The symmetry $K_X(t,s)^*=K_X(s,t)$, with Fubini justified by square integrability of the kernel, proves that $\mathsf K_X$ is bounded and self-adjoint. Thus its unitarily conjugate
$B$ on $\mathcal H_R$ is bounded and self-adjoint. The first inverse identity implies $\ker B=\{0\}$, and hence $\overline{\operatorname{Ran}B}=(\ker B^*)^\perp=\mathcal H_R$.

For completeness, the inverse $B^{-1}$ on $\operatorname{Ran}B$ is self-adjoint. If $v$ belongs to its adjoint domain, there is $w\in\mathcal H_R$ such that
\[
 \langle B^{-1}f,v\rangle_R=\langle f,w\rangle_R
                 \quad(f\in\operatorname{Ran}B).
\]
Taking $f=Bh$ gives $\langle h,v\rangle_R=\langle h,Bw\rangle_R$ for all $h$, so $v=Bw$ and $B^{-1}v=w$. The converse follows from self-adjointness of $B$. Therefore the adjoint
domain and action equal those of $B^{-1}$. The two inverse identities identify $B^{-1}$ with $\mathcal D_X^{u_0,u_1}$. Its inverse is bounded on the whole weighted space, so zero
belongs to its resolvent set. This also proves the regular inverse statement used from \cite{valko2017sine} in the present normalization.
\end{proof}

We use the notation
\[
 \mathtt r\mathcal D_X^{-1}
 =\mathcal U_X(\mathcal D_X^{u_0,u_1})^{-1}\mathcal U_X^{-1}
 =\mathsf K_X.
\]
For $T<1$, extend this bounded inverse by zero to $\mathcal H_{\mathbb C}=L^2((0,1);\mathbb C^2)$, the complexification of $\mathcal H$. We will use the same notation for the
singular inverse constructed below.

We shall use two unitary changes of coordinates. First, for a constant $S\in\mathrm{SL}_2(\mathbb R)$ put
\[
 \widetilde u_i=Su_i,\qquad \widetilde X=XS^{-1},\qquad
 \widetilde R=S^{-\mathsf T}RS^{-1}.
\]
For any real $2\times2$ matrix $M$, multiplication gives $M^{\mathsf T}JM=(\det M)J$ and $MJM^{\mathsf T}=(\det M)J$. In particular,
\[
 (Sf)^*\widetilde R(Sg)=f^*Rg,\qquad
 2\widetilde R^{-1}J(Sf)'=2SR^{-1}S^{\mathsf T}JSf'
                        =S(2R^{-1}Jf'),
\]
and $[Sf,Su_i]=[f,u_i]$. The map $f\mapsto Sf$ consequently identifies the differential domains and the weighted operators. The conjugated columns are unchanged:
$\widetilde X\widetilde u_i=Xu_i$. Also $u_0\mapsto\kappa u_0$, $u_1\mapsto\kappa^{-1}u_1$, for a nonzero real constant $\kappa$, leaves every product in $K_X$ unchanged.

Second, if $\widetilde X^{\mathsf T}\widetilde X=X^{\mathsf T}X$, then $O=\widetilde XX^{-1}$ satisfies
\[
 O^{\mathsf T}O=X^{-\mathsf T}
                \widetilde X^{\mathsf T}\widetilde XX^{-1}=I,
 \qquad \det O=1.
\]
Thus multiplication by $O(t)$ is unitary on ordinary $L^2$, and
\[
 K_{\widetilde X}(s,t)=O(s)K_X(s,t)O(t)^{\mathsf T}.
\]
A common multiplication by $O$ preserves Hilbert--Schmidt distances. The following example shows the effect of changing the factor for one kernel. Take $u_0=-e_1$, $u_1=e_2$, and
$X=Y=I$ on $(0,T)$. Replacing only $Y$ by $J$ leaves both Hamiltonians equal to $I$. On $s<t$ the kernel difference is $(e_2e_1^{\mathsf T}+e_1e_2^{\mathsf T})/2$, and the other
triangle has its transpose. Its squared Hilbert--Schmidt norm is $T^2/2$. We therefore use $\mathcal A$ for all three random operators.

\subsubsection{The maximal operator and boundary traces}
\label{S:sec:s3-maximal-domain}

To extend the regular formula to a singular interval, we establish the boundary traces and weighted kernel integrability. In particular, square integrability of each boundary
column separately is stronger than what the Green kernel requires.

Take $T=1$ and allow a singular right boundary. Assume only that $X,X^{-1}$ are measurable and essentially bounded on each $[0,T]$, $T<1$. Keep $R,a,c,p,q$ as above, and abbreviate
$\mathcal H_R=L^2((0,1),R(t)\,\dd t;\mathbb C^2)$. The deterministic integrability assumptions are
\begin{equation}\label{S:eq:singular-integrability-assumptions}
 Q_X:=\int_0^1q(t)\,\dd t<\infty,\qquad
 I_X:=\int_0^1q(t)\int_0^t p(s)\,\dd s\,\dd t<\infty.
\end{equation}
Thus $Q_X=Q_X(0)$ in the notation for the tail primitive. The first condition puts the terminal direction $u_1$ in $\mathcal H_R$. For the kernel norm, integration must follow the
order $s<t$:
\[
 \int_{0<s<t<1}\|a(s)c(t)^{\mathsf T}\|_{\mathrm{HS}}^2\,\dd s\,\dd t
 =\int_0^1q(t)P_X(t)\,\dd t=I_X.
\]
Bounding this by $(\int p)(\int q)$ would exclude cases with $\int p=\infty$ even when $I_X<\infty$. The value of $\int p$ determines the terminal boundary condition; $I_X$
controls the inverse. Theorem~\ref{S:thm:5.1} verifies these assumptions for $X=\mathcal A(Z)$, including $\beta=2$.

Write $\tau_Xf=2R^{-1}Jf'$ and define
\begin{equation}\label{S:eq:singular-maximal-domain}
 \mathscr D_{\max}(X)=
 \{f\in\mathcal H_R:f\in AC_{\mathrm{loc}}([0,1);\mathbb C^2),
                         \ \tau_Xf\in\mathcal H_R\}.
\end{equation}
The maximal operator $T_{X,\max}$ acts by $\tau_X$ on this domain. Let $T_{X,0}$ be the same expression on $C_c^\infty((0,1);\mathbb C^2)$, and let $T_{X,\min}$ be its closure. For
$f\in\mathscr D_{\max}(X)$ put
\[
 \alpha_f(t)=[f(t),u_1],\qquad
 \nu_f(t)=-[f(t),u_0],\qquad
 \mathfrak b_t(f,g)=f(t)^*Jg(t).
\]
The coordinate maps $\alpha_f,\nu_f$ are complex linear. For each $t<1$, $\mathfrak b_t$ is a sesquilinear form on $\mathscr D_{\max}(X)\times\mathscr D_{\max}(X)$, conjugate
linear in its first argument, linear in its second, and skew-Hermitian: $\mathfrak b_t(g,f)=-\overline{\mathfrak b_t(f,g)}$. By \eqref{S:eq:3.4}, $f=\alpha_fu_0+\nu_fu_1$. The
limiting boundary form $\mathfrak b_1$ has the same conventions whenever its limit exists.

\begin{lemma}\label{S:lem:singular-maximal-traces}
The operator $T_{X,0}$ is densely defined and symmetric, and $T_{X,0}^*=T_{X,\max}$. On the maximal domain the terminal limits $\mathfrak b_1(f,g)$ and $\alpha_f(1)$ exist. The
traces $f(0)$ and $\alpha_f(1)$ are continuous in the graph norm. If $\int_0^1p<\infty$, the full coordinate trace $f(1)$ also exists and is continuous in the graph norm.
\end{lemma}
\begin{proof}
Local equivalence of the weighted and ordinary $L^2$ norms gives $C_c^\infty((0,1);\mathbb C^2)$ dense in $\mathcal H_R$: first replace $f$ by $f_j=\mathbf1_{[1/j,1-1/j]}f$, and
then choose a smooth approximation supported in $[1/(2j),1-1/(2j)]$ with $\|f_j-\varphi_j\|_R<1/j$. Indeed,
\[
 \|f-f_j\|_R^2
 =\int_{(0,1/j)\cup(1-1/j,1)}f^*Rf\,\dd t\longrightarrow0.
\]
Local boundedness of $R^{-1}$ ensures $\tau_X\varphi\in\mathcal H_R$ for every test function. Since $J^*=-J$, integration by parts gives
\[
 \langle T_{X,0}\varphi,\psi\rangle_R
 =-2\int\varphi'^*J\psi
 =2\int\varphi^*J\psi'
 =\langle\varphi,T_{X,0}\psi\rangle_R.
\]
Moreover, for $h,v\in\mathcal H_R$,
\[
 \langle T_{X,0}\varphi,h\rangle_R=\langle\varphi,v\rangle_R
 \quad(\varphi\in C_c^\infty)
 \quad\Longleftrightarrow\quad 2Jh'=Rv
 \quad\text{in distributions}.
\]
The right side implies $h\in H^1((0,T);\mathbb C^2)$ for every $T<1$, because $h,Rv\in L^2(0,T)$. The representatives agree on overlaps, so $h\in\mathscr D_{\max}(X)$ and
$\tau_Xh=v$. Conversely, this differential identity holds for every member of the maximal domain. Hence $T_{X,0}^*=T_{X,\max}$, and the maximal operator is closed.

For $f,g\in\mathscr D_{\max}(X)$ set $h=T_{X,\max}f$ and $v=T_{X,\max}g$. The identities $f'=-JRh/2$ and $g'=-JRv/2$ give
\[
 \frac{\dd}{\dd t}(f^*Jg)
 =f'^*Jg+f^*Jg'=\frac12(f^*Rv-h^*Rg).
\]
The absolute integrals of the two terms on the right are bounded by $\|f\|_R\|v\|_R$ and $\|h\|_R\|g\|_R$, respectively. The terminal limit therefore exists. Integration up to $T$
followed by $T\uparrow1$ gives the Lagrange identity
\begin{equation}\label{S:eq:singular-lagrange-identity}
 \langle f,T_{X,\max}g\rangle_R
 -\langle T_{X,\max}f,g\rangle_R
 =2\bigl(\mathfrak b_1(f,g)-\mathfrak b_0(f,g)\bigr).
\end{equation}

Differentiating the two bilinear coordinates gives
\begin{equation}\label{S:eq:singular-coordinate-derivatives}
 f=\alpha_fu_0+\nu_fu_1,\qquad
 \alpha_f'=-\tfrac12u_1^{\mathsf T}Rh,\qquad
 \nu_f'=\tfrac12u_0^{\mathsf T}Rh.
\end{equation}
For example, $\alpha_f'=(-JRh/2)^{\mathsf T}Ju_1=-u_1^{\mathsf T}Rh/2$. Because $u_1\in\mathcal H_R$, Cauchy--Schwarz gives
\[
 \int_T^1|\alpha_f'(t)|\,\dd t
 \le\frac12\left(\int_T^1q(t)\,\dd t\right)^{1/2}
             \|h\|_{\mathcal H_R(T,1)}.
\]
Thus $\alpha_f(1)$ exists and
\begin{equation}\label{S:eq:singular-terminal-trace-bound}
 |\alpha_f(1)-\alpha_f(T)|
 \le\frac12\left(\int_T^1q(t)\,\dd t\right)^{1/2}
             \|T_{X,\max}f\|_{\mathcal H_R(T,1)}.
\end{equation}

For the bound in the graph norm, fix $T_0\in(0,1)$ and $c_0I\le R\le C_0I$ on $[0,T_0]$. Averaging $f(0)=f(t)-\int_0^t f'(s)\,\dd s$ in $t$ yields
\[
 |f(0)| \le T_0^{-1}\int_0^{T_0}|f(t)|\,\dd t +\int_0^{T_0}|f'(t)|\,\dd t \le (c_0T_0)^{-1/2}\|f\|_R +\tfrac12(C_0T_0)^{1/2}\|h\|_R.
\]
In the second term we used $|Rh|\le\sqrt{C_0}|R^{1/2}h|$. Taking $T=0$ in \eqref{S:eq:singular-terminal-trace-bound} now gives
\[
 |\alpha_f(1)|\le |u_1|\,|f(0)|
                                  +\tfrac12Q_X^{1/2}\|h\|_R,
\]
This bounds the trace by the graph norm. If $\int_0^1p<\infty$, then $\nu_f'=u_0^{\mathsf T}Rh/2$ and $|u_0^{\mathsf T}Rh|\le|R^{1/2}u_0|\,|R^{1/2}h|$. Consequently
\[
 \int_T^1|\nu_f'(t)|\,\dd t
 \le\tfrac12\left(\int_T^1p(t)\,\dd t\right)^{1/2}
            \|h\|_{\mathcal H_R(T,1)}<\infty.
\]
The coordinate $\nu_f(t)$ therefore has a terminal limit, and integration of its derivative gives
\[
 |\nu_f(1)-\nu_f(T)|
 \le\frac12\left(\int_T^1p(t)\,\dd t\right)^{1/2}
                          \|h\|_{\mathcal H_R(T,1)}.
\]
Both coordinates then converge, so $f(t)\to\alpha_f(1)u_0+\nu_f(1)u_1=:f(1)$, and the two trace bounds prove its continuity in the graph norm.
\end{proof}

\subsubsection{The inverse at a singular boundary}
\label{S:sec:s3-singular-inverse}

Under the preceding integrability conditions, the Green kernel defines an integral operator on the full interval. We identify its range with the self-adjoint differential domain,
including the condition at the singular boundary.

Under \eqref{S:eq:singular-integrability-assumptions}, use the same Green kernel as in the regular case:
\begin{equation}\label{S:eq:3.56}
 K_X(s,t)=\frac12\left[
 a(s)c(t)^{\mathsf T}\mathbf1_{\{s<t\}}
 +c(s)a(t)^{\mathsf T}\mathbf1_{\{s>t\}}\right],\qquad 0<s,t<1.
\end{equation}
Write $\mathsf K_X$ for its integral operator on $\mathcal H_{\mathbb C}$ and $\mathsf B_X=\mathcal U_X^{-1}\mathsf K_X\mathcal U_X$ on $\mathcal H_R$. Tonelli gives
\begin{equation}\label{S:eq:3.58}
 \int_0^1 p(s)\int_s^1q(t)\,\dd t\,\dd s
 +\int_0^1 q(s)\int_0^sp(t)\,\dd t\,\dd s=2I_X.
\end{equation}
In particular the kernel is square integrable. Let $\mathcal C_X=\{v\in\mathbb C^2:t\mapsto v\text{ belongs to }\mathcal H_R\}$ be the space of square-integrable constant
functions.

\begin{lemma}\label{S:lem:singular-green-operator}
The operators $\mathsf K_X,\mathsf B_X$ are self-adjoint and injective, with dense ranges, and
\begin{equation}\label{S:eq:singular-hs-norm}
 \|\mathsf K_X\|_{\mathfrak S_2}^2
 =\|\mathsf B_X\|_{\mathfrak S_2(\mathcal H_R)}^2=I_X/2.
\end{equation}
Every $f\in\mathscr D_{\max}(X)$ has the unique decomposition
\begin{equation}\label{S:eq:singular-maximal-decomposition}
 f=\mathsf B_X(T_{X,\max}f)+v_f,\qquad v_f\in\mathcal C_X.
\end{equation}
\end{lemma}
\begin{proof}
The rank-one norm formula and \eqref{S:eq:3.58} prove \eqref{S:eq:singular-hs-norm}. For $F,G\in\mathcal H_{\mathbb C}$, Cauchy--Schwarz on the product space gives
\[
 \int_0^1\!\int_0^1 |F(s)|\|K_X(s,t)\|_{\mathrm{HS}}|G(t)|\,\dd t\,\dd s
 \le\|K_X\|_{L^2((0,1)^2)}\|F\|_2\|G\|_2<\infty.
\]
Thus Fubini applies. Using $K_X(t,s)^*=K_X(s,t)$ after interchanging $s,t$ gives $\langle F,\mathsf K_XG\rangle=\langle\mathsf K_XF,G\rangle$. The bounded operator $\mathsf K_X$ is
therefore self-adjoint, and
\[
 \mathsf B_X^*=(\mathcal U_X^{-1}\mathsf K_X\mathcal U_X)^*
             =\mathcal U_X^{-1}\mathsf K_X^*\mathcal U_X=\mathsf B_X.
\]
The unitary $\mathcal U_X$ also preserves the Hilbert--Schmidt norm.

For $g\in\mathcal H_R$, conjugating the kernel gives the representative $v=\mathsf B_Xg$,
\begin{equation}\label{S:eq:singular-green-solution}
 v(t)=\frac12u_0\int_t^1u_1^{\mathsf T}R(s)g(s)\,\dd s
      +\frac12u_1\int_0^tu_0^{\mathsf T}R(s)g(s)\,\dd s.
\end{equation}
The first integral is absolutely convergent since
\[
 \int_t^1|u_1^{\mathsf T}Rg|\,\dd s
 \le\left(\int_t^1q(s)\,\dd s\right)^{1/2}
                         \|g\|_{\mathcal H_R(t,1)}.
\]
The second integral satisfies
\[
 \int_0^t|u_0^{\mathsf T}R(s)g(s)|\,\dd s
 \le\left(\int_0^tp(s)\,\dd s\right)^{1/2}\|g\|_R<\infty,
 \qquad t<1.
\]
Both indefinite integrals are locally absolutely continuous. Differentiating this representation gives
\[
 v'=\tfrac12(-u_0u_1^{\mathsf T}+u_1u_0^{\mathsf T})Rg
    =-\tfrac12JRg.
\]
Since $v=\mathsf B_Xg$ almost everywhere, $\|v\|_R\le(I_X/2)^{1/2}\|g\|_R<\infty$. The displayed derivative identity gives $2R^{-1}Jv'=g\in\mathcal H_R$. Consequently
$v\in\mathscr D_{\max}(X)$, and
\begin{equation}\label{S:eq:singular-green-boundaries}
 T_{X,\max}\mathsf B_Xg=g,\qquad
 [v(0),u_0]=0,\qquad \alpha_v(1)=0.
\end{equation}
The left condition follows by setting $t=0$; the terminal condition follows from $\alpha_v(t)=\frac12\int_t^1u_1^{\mathsf T}Rg$ and the tail bound.

If $\mathsf B_Xg=0$, its locally absolutely continuous representative vanishes identically on each compact interval; applying $T_{X,\max}$ gives $g=0$. Hence $\mathsf B_X$ is
injective. Self-adjointness gives $(\operatorname{Ran}\mathsf B_X)^\perp=\ker\mathsf B_X=\{0\}$. The same conclusions for $\mathsf K_X$ follow by conjugation. For an arbitrary
$f\in\mathscr D_{\max}(X)$,
\[
 (f-\mathsf B_XT_{X,\max}f)'=0.
\]
The difference is therefore constant and lies in $\mathcal H_R$, which proves \eqref{S:eq:singular-maximal-decomposition}. If $f=\mathsf B_Xg+v$ with $v\in\mathcal C_X$, applying
$T_{X,\max}$ forces $g=T_{X,\max}f$ and then determines $v$.
\end{proof}

The condition $\int_0^1q<\infty$ gives $u_1\in\mathcal C_X$. If $v=\alpha u_0+\nu u_1\in\mathcal C_X$ and $\alpha\ne0$, then $u_0=\alpha^{-1}(v-\nu u_1)\in\mathcal H_R$. Hence
\[
 \mathcal C_X=
 \begin{cases}
  \mathbb C u_1,&\displaystyle\int_0^1p=\infty,\\
  \mathbb C^2,&\displaystyle\int_0^1p<\infty.
 \end{cases}
\]
For $S=(u_0\ u_1)$, $p+q=\operatorname{tr}(S^{\mathsf T}RS)=\|XS\|_{\mathrm{HS}}^2$. The inequalities $\|XS\|_{\mathrm{HS}}\le\|S\|_{\mathrm{op}}\|X\|_{\mathrm{HS}}$ and
$\|X\|_{\mathrm{HS}}\le\|S^{-1}\|_{\mathrm{op}}\|XS\|_{\mathrm{HS}}$ give
\[
 \|S^{-1}\|_{\mathrm{op}}^{-2}\operatorname{tr}R
 \le p+q\le\|S\|_{\mathrm{op}}^2\operatorname{tr}R.
\]
Under \eqref{S:eq:singular-integrability-assumptions}, the condition $\int p<\infty$ is therefore equivalent to
$\int\operatorname{tr}R<\infty$.

Define $T_X$ as the restriction of $T_{X,\max}$ to
\begin{equation}\label{S:eq:singular-selfadjoint-domain}
 \mathcal D(T_X)=
 \begin{cases}
 \{f\in\mathscr D_{\max}(X):[f(0),u_0]=0\},
       &\displaystyle\int_0^1p=\infty,\\[1mm]
 \{f\in\mathscr D_{\max}(X):[f(0),u_0]=0,\ \alpha_f(1)=0\},
       &\displaystyle\int_0^1p<\infty.
 \end{cases}
\end{equation}

\begin{theorem}\label{S:thm:singular-inverse-domain}
The operator $T_X$ is self-adjoint, has compact resolvent, and $0\in\rho(T_X)$. Its inverse is $\mathsf B_X$; more precisely,
\begin{equation}\label{S:eq:singular-two-sided-inverse}
 T_X^{-1}=\mathsf B_X,\qquad T_X\mathsf B_X=I_{\mathcal H_R},
 \qquad \mathsf B_XT_X=I_{\mathcal D(T_X)}.
\end{equation}
\end{theorem}
\begin{proof}
Put $B=\mathsf B_X$. Lemma~\ref{S:lem:singular-green-operator} shows that $B$ is bounded, self-adjoint, injective and Hilbert--Schmidt, and that
\[
 T_{X,\max}Bg=g,
 \qquad [Bg(0),u_0]=0,\qquad\alpha_{Bg}(1)=0.
\]
Hence $\operatorname{Ran}B\subset\mathcal D(T_X)$. We prove the reverse inclusion in both cases.

If $\int_0^1p=\infty$, the square-integrable constant vectors form $\mathcal C_X=\mathbb C u_1$. For $f\in\mathcal D(T_X)$ the decomposition of a function in the maximal domain
gives $f=Bg+\nu u_1$ with $g=T_{X,\max}f$. Its left trace satisfies
\[
 0=[f(0),u_0]=[Bg(0),u_0]+\nu[u_1,u_0]=-\nu.
\]
Thus $\nu=0$ and $f=Bg$. No right boundary condition has been imposed in this argument. The right trace $\alpha_f(1)=0$ follows as a consequence for these functions.

If $\int_0^1p<\infty$, write the same decomposition as $f=Bg+\alpha u_0+\nu u_1$. The left condition again gives $\nu=0$. Since $\alpha_{Bg}(1)=0$, the right condition gives
\[
 0=\alpha_f(1)=[\alpha u_0+\nu u_1,u_1]=\alpha.
\]
Therefore $f=Bg$ in this case as well. We have proved $\mathcal D(T_X)=\operatorname{Ran}B$, and the action of $T_X$ on this range is $B^{-1}$.

The range of $B$ is dense because $B=B^*$ and $\ker B=0$. The calculation of the adjoint of the inverse in the regular case applies unchanged: if
$\langle B^{-1}f,v\rangle_R=\langle f,w\rangle_R$ for all $f=Bh$, then $v=Bw$. It follows that $(B^{-1})^*=B^{-1}$ with equality of domains. Thus $T_X$ is self-adjoint. The range
and action identities give $T_XB=I$ on $\mathcal H_R$ and $BT_X=I$ on $\mathcal D(T_X)$. The bound $\|T_X^{-1}\|=\|B\|<\infty$ proves $0\in\rho(T_X)$.

For $z\in\rho(T_X)$, the spectral theorem, or multiplication on $\mathcal D(T_X)$, gives
\begin{equation}
 (T_X-z)^{-1}=B(I-zB)^{-1},\qquad
 (I-zB)^{-1}=I+z(T_X-z)^{-1}.
 \label{S:eq:singular-resolvent-from-inverse}
\end{equation}
The right side of the first identity maps into $\operatorname{Ran}B=\mathcal D(T_X)$, and multiplying on either side verifies the inverse equations. Since $B$ is Hilbert--Schmidt
and the other factor is bounded, every resolvent is Hilbert--Schmidt and therefore compact. The argument establishes directly the domain used in the canonical formulation of
\cite{valko2017sine}.
\end{proof}

The domain also has an explicit graph core:
\[
 \{\mathsf B_Xg:g\in C_c^\infty((0,1);\mathbb C^2)\}.
\]
Indeed, for $f=\mathsf B_Xg$ choose $g_j\in C_c^\infty$ with $g_j\to g$ in $\mathcal H_R$. Then
\[
 \|\mathsf B_Xg_j-f\|_R+\|T_X\mathsf B_Xg_j-T_Xf\|_R
 \le(\|\mathsf B_X\|+1)\|g_j-g\|_R\longrightarrow0.
\]
By \eqref{S:eq:singular-green-solution}, each approximating function is a constant multiple of $u_0$ near zero and a constant multiple of $u_1$ near one. Its constant boundary
pieces extend to the boundary points. The spectral theorem for $\mathsf B_X$ yields an orthonormal basis of eigenvectors with nonzero real eigenvalues tending to zero. Their
reciprocals are the eigenvalues of $T_X$, and
\[
 \sum_{\lambda\in\operatorname{Spec}(T_X)}\lambda^{-2}=I_X/2,
\]
with multiplicities included.

\subsubsection{The minimal domain and the Weyl alternative}
\label{S:sec:s3-minimal-domain}

The Green operator identifies one self-adjoint realization. The minimal and maximal domains explain why a right boundary condition is needed in one case and not in the other.
Throughout this discussion, \eqref{S:eq:singular-integrability-assumptions} remains in force.

\begin{proposition}\label{S:prop:singular-minimal-domain}
One has $T_{X,\min}^*=T_{X,\max}$ and
\begin{equation}\label{S:eq:singular-minimal-domain}
 \mathcal D(T_{X,\min})=
 \begin{cases}
 \{f\in\mathscr D_{\max}(X):f(0)=0\},&\displaystyle\int_0^1p=\infty,\\
 \{f\in\mathscr D_{\max}(X):f(0)=0,\ f(1)=0\},
                                      &\displaystyle\int_0^1p<\infty.
 \end{cases}
\end{equation}
\end{proposition}
\begin{proof}
For $f,g\in\mathscr D_{\max}(X)$, the boundary form at $1$ is
\begin{equation}\label{S:eq:singular-terminal-boundary-form}
 \mathfrak b_1(f,g)=
 \begin{cases}
 0,&\displaystyle\int_0^1p=\infty,\\
 f(1)^*Jg(1),&\displaystyle\int_0^1p<\infty.
 \end{cases}
\end{equation}
In the second case this follows from the coordinate limits in Lemma~\ref{S:lem:singular-maximal-traces}. In the first case write $f=f_0+\lambda u_1$, $g=g_0+\mu u_1$ with
$f_0,g_0\in\operatorname{Ran}\mathsf B_X =\mathcal D(T_X)$, using \eqref{S:eq:singular-maximal-decomposition}. Self-adjointness and the Lagrange identity imply
$\mathfrak b_1(f_0,g_0)=\mathfrak b_0(f_0,g_0)=0$, because $f_0(0),g_0(0)\in\mathbb C u_0$ and $u_0$ is real. The remaining terms are
\[
 \mathfrak b_t(f_0,\mu u_1)=\mu\overline{\alpha_{f_0}(t)},
 \qquad
 \mathfrak b_t(\lambda u_1,g_0)=-\overline\lambda\alpha_{g_0}(t),
 \qquad
 \mathfrak b_t(\lambda u_1,\mu u_1)=0.
\]
The limits vanish by \eqref{S:eq:singular-green-boundaries}, proving \eqref{S:eq:singular-terminal-boundary-form}.

Since $T_{X,0}$ is densely defined and symmetric, it is closable. Lemma~\ref{S:lem:singular-maximal-traces} gives
\[
 T_{X,\min}^*=T_{X,0}^*=T_{X,\max},\qquad
 T_{X,\min}=T_{X,\max}^*.
\]
The inclusion $T_{X,0}\subset T_{X,\max}$ implies $T_{X,\max}^*\subset T_{X,0}^*=T_{X,\max}$. The adjoint identity, tested first on compactly supported functions, therefore says
that $f\in\mathcal D(T_{X,\max}^*)$ exactly when $f\in\mathscr D_{\max}(X)$ and
\[
 \mathfrak b_1(f,g)-\mathfrak b_0(f,g)=0
 \qquad\text{for every }g\in\mathscr D_{\max}(X).
\]
To test the left boundary take $g(t)=\chi_0(t)v$, where $\chi_0=1$ near zero and $\chi_0=0$ near one. Local bounds on $R,R^{-1}$ imply
\[
 \|\chi_0v\|_R^2+\|\tau_X(\chi_0v)\|_R^2
 \le C_{\chi_0,X}|v|^2<\infty,
 \qquad v\in\mathbb C^2.
\]
Thus the adjoint boundary identity implies $f(0)^*Jv=0$ for every $v$, or $f(0)=0$. If $\int p=\infty$, this is sufficient by \eqref{S:eq:singular-terminal-boundary-form}. If
$\int p<\infty$, choose $\chi_1=0$ near zero and $\chi_1=1$ near one. Every constant vector belongs to $\mathcal H_R$, and
\[
 \|\chi_1v\|_R\le\|\chi_1\|_\infty\|v\|_R,\qquad
 \|\tau_X(\chi_1v)\|_R^2
 =4\int_{\operatorname{supp}\chi_1'}
       |\chi_1'|^2v^*J^*R^{-1}Jv\,\dd t<\infty.
\]
Testing with $\chi_1v$ therefore gives $f(1)^*Jv=0$ for every $v$, or $f(1)=0$. Conversely, these vanishing traces make the Lagrange boundary expression zero for every
$g\in\mathscr D_{\max}(X)$, which proves \eqref{S:eq:singular-minimal-domain}.

\end{proof}

For $z\in\mathbb C\setminus\mathbb R$, consider $2Jf'=zRf$. The right boundary is limit circle when the space of solutions square integrable near one has dimension two, and limit
point when that dimension is one. We determine this dimension for each nonreal $z$.

\begin{proposition}\label{S:prop:singular-weyl-alternative}
The right boundary is limit circle if $\int_0^1p<\infty$ and limit point otherwise. The deficiency indices of $T_{X,\min}$ are, respectively, $(2,2)$ and $(1,1)$.
\end{proposition}
\begin{proof}
Fix $z\in\mathbb C\setminus\mathbb R$, and put $B=\mathsf B_X$. Since $B$ is bounded and self-adjoint, $I-zB$ is invertible: for real $\mu\in\sigma(B)$, $1-z\mu\ne0$, and the
continuous modulus has a positive minimum on the compact set $\sigma(B)$. The decomposition of a function in the maximal domain gives, for $f\in\ker(T_{X,\max}-z)$,
\[
 f=B(T_{X,\max}f)+v=zBf+v,\qquad v\in\mathcal C_X.
\]
Consequently $f=(I-zB)^{-1}v$. Conversely, for every $v\in\mathcal C_X$ let
\[
 f_v=(I-zB)^{-1}v=v+zB(I-zB)^{-1}v.
\]
The first term is a square-integrable constant and the second belongs to $\operatorname{Ran}B\subset\mathscr D_{\max}(X)$. Thus $f_v\in\mathscr D_{\max}(X)$ and
\[
 T_{X,\max}f_v=z(I-zB)^{-1}v=zf_v.
\]
The map $v\mapsto f_v$ is injective because $I-zB$ is invertible, and the preceding decomposition proves surjectivity onto the nonreal solution space. It follows that
\[
 \dim\ker(T_{X,\max}-z)=\dim\mathcal C_X
 =\begin{cases}
 1,&\int_0^1p=\infty,\\
 2,&\int_0^1p<\infty.
 \end{cases}
\]
This calculation concerns the whole interval. It also gives the dimension of solutions square integrable near one. Indeed, the coefficient in $f'=-zJRf/2$ is bounded on every
compact subinterval. Its integral equation has a unique solution for each initial value, by successive integration and the estimate $(C|t-s|)^m/m!$ for the $m$th iterated integral.
A solution defined near one therefore extends uniquely to zero and remains square integrable on the added compact interval. Restriction and extension preserve the dimension of the
space of square-integrable solutions.

By the definitions of the limit-point and limit-circle cases, this proves the claimed classification. Finally $T_{X,\min}^*=T_{X,\max}$, so taking $z=\mathrm i$ and $z=-\mathrm i$
gives
\[
 n_+(T_{X,\min})=\dim\ker(T_{X,\max}-\mathrm i),\qquad
 n_-(T_{X,\min})=\dim\ker(T_{X,\max}+\mathrm i).
\]
Both equal one in the first case and two in the second. This supplies the dimension conclusion of the Weyl alternative directly under the integrability hypotheses used here.
\end{proof}

\subsubsection{Regular restrictions and their resolvents}
\label{S:sec:s3-regular-restrictions}

Regular restrictions approximate the singular inverse provided the terminal boundary direction is kept fixed. The exact compression identity below makes this choice explicit and
then transfers convergence from inverses to resolvents.

Continue under \eqref{S:eq:singular-integrability-assumptions}. For $0<T<1$ let $T_{X,T}$ act by $\tau_X$ in $\mathcal H_R(0,T)$ on
\begin{equation}\label{S:eq:singular-regular-domain}
 \mathcal D(T_{X,T})=
 \{f\in H^1((0,T);\mathbb C^2):[f(0),u_0]=[f(T),u_1]=0\}.
\end{equation}
The coordinate vector $u_1$ is fixed as $T$ varies. Let $\iota_T$ be extension by zero from $L^2((0,T);\mathbb C^2)$ to $\mathcal H_{\mathbb C}$, and $P_T=\iota_T\iota_T^*$. Denote
the restricted multiplication map by $\mathcal U_{X,T}$ and put
\[
 \mathsf T_X=\mathcal U_XT_X\mathcal U_X^{-1},\qquad
 \mathsf K_{X,T}=\iota_T\mathcal U_{X,T}T_{X,T}^{-1}
                       \mathcal U_{X,T}^{-1}\iota_T^*.
\]
The inverse of $T_{X,T}$ exists by Theorem~\ref{S:thm:3.2}.

\begin{theorem}\label{S:thm:singular-regular-truncations}
The regular inverse is the exact compression
\begin{equation}\label{S:eq:singular-exact-compression}
 \mathsf K_{X,T}=P_T\mathsf K_XP_T.
\end{equation}
Moreover, as $T\uparrow1$,
\begin{equation}\label{S:eq:singular-truncation-hs-tail}
 \|\mathsf K_X-\mathsf K_{X,T}\|_{\mathfrak S_2}^2
 =\frac12\int_T^1q(t)\int_0^tp(s)\,\dd s\,\dd t\longrightarrow0.
\end{equation}
\end{theorem}
\begin{proof}
For $0<s<T$, the regular Green formula gives, after conjugation,
\[
 (\mathsf K_{X,T}F)(s)
 =\frac12a(s)\int_s^Tc(t)^{\mathsf T}F(t)\,\dd t
  +\frac12c(s)\int_0^sa(t)^{\mathsf T}F(t)\,\dd t.
\]
For $s\ge T$, $(\mathsf K_{X,T}F)(s)=0$. In $P_T\mathsf K_XP_TF$, the function integrated against the full kernel is $P_TF$. For $s<T$, the upper limit in the first integral is
therefore $T$ and the second integral is unchanged. Multiplication by $P_T$ makes the resulting function zero for $s\ge T$. Hence $\mathsf K_{X,T}=P_T\mathsf K_XP_T$ as operators
on the common $L^2$ space.

The complementary kernel is supported on $\{s\ge T\}\cup\{t\ge T\}$. In the triangle $s<t$ this is just the region $T<t<1$, $0<s<t$, where the integral of the squared kernel norm
is
\[
 \frac14\int_T^1q(t)\int_0^tp(s)\,\dd s\,\dd t.
\]
The transposed triangle contributes the same amount. The diagonal has product measure zero, so there is no additional term. Consequently
\[
 \|\mathsf K_X-\mathsf K_{X,T}\|_{\mathfrak S_2}^2
 =\frac12\int_T^1q(t)\int_0^tp(s)\,\dd s\,\dd t.
\] The function $q(t)\int_0^tp$ is integrable by hypothesis, and its integral over $(T,1)$ tends to zero by dominated convergence. This proves the limit for every $T\uparrow1$, not only along a prescribed sequence of cutoffs.
\end{proof}

For the embedded regular resolvent, work on $P_T\mathcal H_{\mathbb C}$, since $\mathsf K_{X,T}$ vanishes on its orthogonal complement. For $z\in\rho(T_{X,T})$, define
\[
 \mathsf R_T(z)=\iota_T\mathcal U_{X,T}(T_{X,T}-z)^{-1}
                  \mathcal U_{X,T}^{-1}\iota_T^*,\qquad
 \mathsf R(z)=(\mathsf T_X-z)^{-1}.
\]
Applying \eqref{S:eq:singular-resolvent-from-inverse} on $P_T\mathcal H_{\mathbb C}$ and on its orthogonal complement separately gives
\[
 \mathsf R_T(z)=(I-z\mathsf K_{X,T})^{-1}\mathsf K_{X,T},
 \qquad
 \mathsf R(z)=(I-z\mathsf K_X)^{-1}\mathsf K_X.
\]
Let $E\subset\rho(T_X)$ be compact, and set
\[
 M_E=\sup_{z\in E}\|(I-z\mathsf K_X)^{-1}\|_{\mathrm{op}},
 \qquad L_E=\sup_{z\in E}|z|.
\]
Compactness bounds $L_E$. To bound $M_E$, fix $z_0\in E$ and factor
\[
 I-z\mathsf K_X=(I-z_0\mathsf K_X)
 \left[I-(z-z_0)(I-z_0\mathsf K_X)^{-1}\mathsf K_X\right].
\]
If $|z-z_0|\|(I-z_0\mathsf K_X)^{-1}\mathsf K_X\|\le1/2$, the Neumann series gives $\|(I-z\mathsf K_X)^{-1}\|\le2\|(I-z_0\mathsf K_X)^{-1}\|$. A finite collection of these
neighborhoods covers $E$, so $M_E<\infty$. By Theorem~\ref{S:thm:singular-regular-truncations}, for every $T$ sufficiently close to one,
\[
 L_EM_E\|\mathsf K_{X,T}-\mathsf K_X\|_{\mathrm{op}}<\frac12.
\]
The factorization
\[
 I-z\mathsf K_{X,T} =(I-z\mathsf K_X)\bigl[I-z(I-z\mathsf K_X)^{-1} (\mathsf K_{X,T}-\mathsf K_X)\bigr]
\]
then proves, uniformly for $z\in E$, invertibility and the bound $\|(I-z\mathsf K_{X,T})^{-1}\|_{\mathrm{op}}\le2M_E$. Restriction to $P_T\mathcal H_{\mathbb C}$ and the regular
inverse formula show that $E\subset\rho(T_{X,T})$.

Multiplying $\mathsf R_T(z)-\mathsf R(z)$ on the left by $I-z\mathsf K_{X,T}$ and on the right by $I-z\mathsf K_X$ gives
\[
 \mathsf K_{X,T}(I-z\mathsf K_X)
 -(I-z\mathsf K_{X,T})\mathsf K_X
 =\mathsf K_{X,T}-\mathsf K_X.
\]
Hence
\[
 \mathsf R_T(z)-\mathsf R(z)
 =(I-z\mathsf K_{X,T})^{-1}
   (\mathsf K_{X,T}-\mathsf K_X)(I-z\mathsf K_X)^{-1}.
\]
The Hilbert--Schmidt ideal inequality yields
\begin{equation}\label{S:eq:singular-embedded-resolvents}
 \sup_{z\in E}\|\mathsf R_T(z)-\mathsf R(z)\|_{\mathfrak S_2}
 \le 2M_E^2\|\mathsf K_{X,T}-\mathsf K_X\|_{\mathfrak S_2}
 \longrightarrow0.
\end{equation}
This proves the resolvent convergence for all sufficiently large cutoffs, uniformly on every compact subset of the limiting resolvent set.

For later use, if $A$ is any Hilbert--Schmidt integral operator with kernel $K$, write $A^T=P_TAP_T$. Multiplication of its kernel by $\mathbf1_{(0,T)^2}$ gives
\[
 \|A-A^T\|_{\mathfrak S_2}^2
 =\int_{\{s\ge T\}\cup\{t\ge T\}}\|K(s,t)\|_{\mathrm{HS}}^2
                                      \,\dd s\,\dd t,
 \qquad
 \langle A^T,A-A^T\rangle_{\mathfrak S_2}=0.
\]
The integration region consists of the square beyond $T$ and the two rectangles with one variable on each side of $T$.

The same restriction identity applies to the finite operators on a regular interval, where the coefficients are piecewise constant.

\begin{lemma}\label{S:lem:regular-interval-restriction}
Let $X$ be piecewise constant on a finite partition of $[0,L]$, with values in $\mathrm{SL}_2(\mathbb R)$, and use the boundary vectors $u_0,u_1$ at $0,L$. The regular operator
$\mathcal D_X$ is self-adjoint and invertible. For $0<T\le L$, $P_T\mathsf K_XP_T$ is the conjugated inverse on $[0,T]$ with the same coordinate boundary vectors, extended by zero
to
$[0,L]$.
\end{lemma}
\begin{proof}
Let the coefficient values be $X_1,\ldots,X_M$. Each is invertible. The finite constants
\[
 C=\max_{k\le M}\|X_k\|_{\mathrm{op}}^2,
 \qquad c=\left(\max_{k\le M}\|X_k^{-1}\|_{\mathrm{op}}^2\right)^{-1}>0
\]
satisfy $cI\le R\le CI$ almost everywhere on $(0,L)$ and on every restricted interval $(0,T)$. Splitting a cell at $T$ changes neither inequality. The regular interval theorem
therefore gives the self-adjoint domain and inverse on both intervals.

For $s,t<T$, the columns are the same $a(s)=X(s)u_0$, $c(t)=X(t)u_1$, and the same coordinate boundary vectors enter the two Green formulas. Their kernels agree pointwise off the
diagonal on $(0,T)^2$. Extension by zero therefore has kernel $\mathbf1_{(0,T)^2}K_X$, which is exactly $P_T\mathsf K_XP_T$.

On $P_TL^2$ this is the inverse of the regular differential operator on $(0,T)$. On $(I-P_T)L^2$ it is zero. The latter zero operator is only an extension convention; it is not the
inverse of a differential operator on the discarded part of the interval. No continuity of $X$ across a partition boundary is needed, because the differential domain is expressed
through the canonical function $f$, not through its possibly discontinuous image $Xf$.
\end{proof}

\subsubsection{Piecewise constant coefficients}
\label{S:sec:s3-piecewise-coefficients}

For sampled coefficients, the Green formula remains exact, but a diagonal cell contains both orderings of the two kernel variables. Keeping its two triangles separate is necessary
when computing the squared Hilbert--Schmidt error.

Let $X$ be constant on a finite partition of a regular interval. Index the intervals in reverse order: if $i>j$, the interval $I_i$ lies to the left of $I_j$ in the coordinate $t$.
Write $\Delta_i=|I_i|$, $X_i=X|_{I_i}$, $a_i=X_iu_0$, and $c_i=X_iu_1$. Theorem~\ref{S:thm:3.2} applies since the finitely many matrices and their inverses are bounded. For
$s\in I_i,t\in I_j$ with $i\ne j$,
\[
 K_X(s,t)=
 \begin{cases}
  a_ic_j^{\mathsf T}/2,&i>j,\\
  c_ia_j^{\mathsf T}/2,&i<j.
 \end{cases}
\]
On $I_i\times I_i$ both orders of $s,t$ occur. The two values are $a_ic_i^{\mathsf T}/2$ and $c_ia_i^{\mathsf T}/2$, each on a triangle of area $\Delta_i^2/2$. Their difference is
exactly
\[
 \tfrac12(a_ic_i^{\mathsf T}-c_ia_i^{\mathsf T})
 =\tfrac12X_iJX_i^{\mathsf T}=J/2.
\]
Thus the kernel takes two different values on each diagonal square. Integrating the two triangles and all rectangles gives
\begin{equation}\label{S:eq:s3-step-kernel-norm}
 \|\mathsf K_X\|_{\mathfrak S_2}^2 =\frac14\sum_{i>j}\Delta_i\Delta_j|a_i|^2|c_j|^2 +\frac14\sum_{i<j}\Delta_i\Delta_j|c_i|^2|a_j|^2 +\frac14\sum_i\Delta_i^2|a_i|^2|c_i|^2.
\end{equation}
The last sum is the contribution of the diagonal squares. Their common edges have two-dimensional measure zero.

The jump $J/2$ is independent of the coefficient matrix. It cancels in a difference of normalized kernels on the same partition, and also in the first variation along
$\mathrm{SL}_2(\mathbb R)$. These two operators admit the finite matrix representations proved in Lemma~\ref{S:lem:exact-mesh-kernel-stability}. The individual inverse is injective
on an infinite-dimensional $L^2$ space and therefore has infinite rank. In particular, replacing its two values on a diagonal square by a single block would change the operator.

We use the intervals $I_{n,k}=[1-(k+1)/n,1-k/n)$, $m\le k<n$, on $[0,1-m/n)$. Define the sums that omit the current interval by
\[
 A_{n,k}=\frac1n\sum_{\ell=m}^{k-1}|c_\ell|^2,\qquad
 C_{n,k}=\frac1n\sum_{\ell=k+1}^{n-1}|a_\ell|^2.
\]
They are the contributions of complete intervals to the two primitives. With $r=1-t=(k+v)/n$, $0<v<1$, the exact remaining fractions are
\begin{equation}\label{S:eq:s3-step-primitives}
 Q_X(1-r)=A_{n,k}+\frac v n|c_k|^2,\qquad
 P_X(1-r)=C_{n,k}+\frac{1-v}n|a_k|^2.
\end{equation}
For example, $\int_t^{1-m/n}q$ contains all intervals with $\ell<k$ and the part of $I_{n,k}$ of length $v/n$ to the right of $t$. The other formula follows from the intervals with
$\ell>k$ and the part of length $(1-v)/n$ to the left. This also explains the names of the reflected weights in Sections
\ref{S:sec:sampling-times-and-estimates-for-the-boundary-solutions} and~\ref{S:sec:stability-and-fr-echet-expansion-of-the-inverse-kernel}: $A$ integrates $|c|^2$ below $r$, while
$C$ integrates $|a|^2$ above $r$.

The fractions in \eqref{S:eq:s3-step-primitives} matter for perturbations that vary inside an interval. Their exact contribution to the weighted norm, including the factor $1/2$
available only for constant fields, is calculated in \eqref{S:eq:continuous-integral-energy}--\eqref{S:eq:diagonal-triangles-integral}.

\subsubsection{A diagonal example}
\label{S:sec:s3-diagonal-example}

The distinction between the two integrability conditions is already visible without randomness. In the reflected coordinate take
\[
 X_\beta^{\mathrm e}(r)=\mathcal A(\ii r^{2/\beta})
 =\begin{pmatrix}r^{-1/\beta}&0\\0&r^{1/\beta}\end{pmatrix},
 \qquad u_0=-e_1,\quad u_1=e_2.
\]
Then, with $\gamma=2/\beta$,
\[
 p_\beta(r)=r^{-\gamma},\qquad q_\beta(r)=r^\gamma,\qquad
 Q_\beta(r)=\frac{r^{1+\gamma}}{1+\gamma},\qquad
 P_\beta(r)=
 \begin{cases}
 (1-r^{1-\gamma})/(1-\gamma),&\gamma\ne1,\\
 \log(1/r),&\gamma=1.
 \end{cases}
\]
These are exact integrals of $q_\beta$ over $(0,r)$ and $p_\beta$ over $(r,1)$. The column $a$ fails to be square integrable when $\gamma\ge1$. Nevertheless the kernel criterion is
finite for every $\gamma>0$, since Tonelli gives
\[
 \int_0^1r^\gamma P_\beta(r)\,\dd r
 =\int_0^1u^{-\gamma}\int_0^u r^\gamma\,\dd r\,\dd u
 =\frac1{2(1+\gamma)}.
\]
Thus the squared Hilbert--Schmidt norm in this example is $1/[4(1+\gamma)]$, on either side of $\gamma=1$. The quadratic form from the introduction is
$\mathfrak h_{\mathrm{aff}}(v)=|v|^2/4+v_2^2$. On the two columns of this example it equals $p_\beta/4$ and $5q_\beta/4$. The corresponding integrand is therefore
\[
 G_\beta^{\mathrm{geo}}(r)
 =\frac18p_\beta(r)Q_\beta(r)+\frac58q_\beta(r)P_\beta(r)
 =\frac{r}{8(1+\gamma)}+\frac58r^\gamma P_\beta(r).
\]
As $r\downarrow0$ it has order $r$ for $\gamma>1$, order $r\log(1/r)$ for $\gamma=1$, and order $r^\gamma$ for $0<\gamma<1$, as follows by substitution of the displayed formula for
$P_\beta$. The integral is finite for every fixed $\gamma>0$. Indeed, Tonelli gives
\[
 \int_0^1r^{\gamma-1}P_\beta(r)\,\dd r
 =\int_0^1u^{-\gamma}\int_0^u r^{\gamma-1}\,\dd r\,\dd u
 =\frac1\gamma,
\]
and consequently
\[
 \int_0^1\frac{G_\beta^{\mathrm{geo}}(r)}r\,\dd r
 =\frac1{8(1+\gamma)}+\frac5{8\gamma}.
\]
This includes $\gamma=1$ without taking a limit at a boundary point. For the boundary columns determined by the Brownian path, Theorem~\ref{S:thm:5.1} and
Proposition~\ref{S:prop:5.2} give the corresponding pathwise estimates.

\subsubsection{Normalization of the operator and finite walk}
\label{S:sec:normalization-of-the-operator-and-finite-walk}

We apply the preceding results to the three coefficient functions from the introduction, always with $X=\mathcal A(z)$. Direct multiplication gives
\[
 \mathcal A(x+\ii y)^{\mathsf T}\mathcal A(x+\ii y)
 =\frac1y\begin{pmatrix}1&-x\\-x&x^2+y^2\end{pmatrix}.
\]
It is positive definite because it is $X^{\mathsf T}X$ with $X$ invertible, and its determinant is one.

For the limiting path, $X_\beta(t)=\mathcal A(Z_{\tau_\beta(t)})$, with $u_0=(-1,0)^{\mathsf T}$ and $u_1=(q_\infty,1)^{\mathsf T}$. Theorem~\ref{S:thm:sine-selfadjoint-domain}
verifies \eqref{S:eq:singular-integrability-assumptions} and shows that $\int p=\infty$ for $\beta\le2$ and $\int p<\infty$ for $\beta>2$. Thus
\eqref{S:eq:singular-selfadjoint-domain} defines $\mathrm{Sine}_\beta$ with no right condition in the first range and $\lim_{t\uparrow1}[f(t),u_1]=0$ in the second. Its conjugated
inverse is $\mathtt r\mathrm{Sine}_\beta=\mathsf K_{X_\beta}$. At $\beta=2$, the divergence is the critical case of Theorem~\ref{S:thm:5.1}. By the spectral identification in
\cite{valko2017sine}, the spectrum has the law of the Sine$_\beta$ point process; the deterministic domain and inverse formulas have been proved above.

For the circular path write $b_0,\ldots,b_{n-1}\in\mathbb H$ and $b_n\in\partial\mathbb H$ for the walk specified in the introduction. Apply $\mathcal A$ only to its interior
points:
\[
 X_{\beta,n}^{\mathrm{circ}}(t)=\mathcal A(b_{\lfloor nt\rfloor}),
 \qquad 0\le t<1.
\]
The final point specifies the boundary vector $u_1=(b_n,1)^{\mathsf T}$. The last interval uses $b_{n-1}$ in the step coefficient. For $z=x+\ii y$ its harmonic measure on
$\mathbb R$ has density
\[
 \nu_z(\dd q)=\frac{y\,\dd q}{\pi((q-x)^2+y^2)}.
\]
Its integral over $\mathbb R$ is $\pi^{-1}[\arctan((q-x)/y)]_{-\infty}^{\infty}=1$. Hence there is no atom at infinity, and $b_n$ is finite almost surely. For fixed $n$ the
finitely many interior heights are positive, and $X_{\beta,n}^{\mathrm{circ}}$ and its inverse are bounded. Theorem~\ref{S:thm:3.2} therefore defines the regular operator
$\mathrm{Circ}_{\beta,n}$ and its conjugated inverse $\mathtt r\mathrm{Circ}_{\beta,n}$. Under the coupling of Proposition~\ref{S:prop:7.2}, $b_n=q_\infty$, so all comparisons use
the same terminal vector.

The finite spectrum is
\begin{equation}\label{S:eq:2.56}
 \operatorname{Spec}(\mathrm{Circ}_{\beta,n})
 \overset{\mathrm d}=\operatorname{Spec}_{\mathrm{per}}(C\beta E_n)
 =\{n(\theta_{j,n}+2\pi m):1\le j\le n,\ m\in\mathbb Z\},
\end{equation}
where the angles have the circular $\beta$-ensemble distribution. This is the circular ensemble representation of \cite{valko2017sine}. The nonzero eigenvalues of the conjugated
inverse are the reciprocals of the direct eigenvalues, with the same multiplicities. Indeed, $Tf=\lambda f$ implies $T^{-1}f=\lambda^{-1}f$; conversely $T^{-1}f=\mu f$ with
$\mu\ne0$ puts $f=\mu^{-1}T^{-1}f$ in $\mathcal D(T)$ and gives $Tf=\mu^{-1}f$. The separate ordering of positive and negative inverse eigenvalues used for comparison is fixed in
Section~\ref{S:sec:spectral-consequences}.

The auxiliary sampled coefficients are
\[
 X_{\beta,n}^{\mathrm{br}}(t)=\mathcal A(Z_{\vartheta_{n,k+1}}),
 \qquad t\in I_{n,k},\quad 0\le k<n.
\]
They define a regular operator with the same two boundary vectors,
\[
 \begin{split}
 \mathrm{Br}_{\beta,n}f&=2(R_{\beta,n}^{\mathrm{br}})^{-1}Jf',
 \qquad R_{\beta,n}^{\mathrm{br}}
          =(X_{\beta,n}^{\mathrm{br}})^{\mathsf T}X_{\beta,n}^{\mathrm{br}},\\
 \mathcal D(\mathrm{Br}_{\beta,n})
 &=\{f\in H^1((0,1);\mathbb C^2):
                       [f(0),u_0]=[f(1),u_1]=0\}.
 \end{split}
\]
For each fixed $n$, the finitely many positive sampled heights make the coefficient and its inverse bounded. Theorem~\ref{S:thm:3.2} gives $\mathtt r\mathrm{Br}_{\beta,n}$ with
columns $a_n=X_{\beta,n}^{\mathrm{br}}u_0$ and $c_n=X_{\beta,n}^{\mathrm{br}}u_1$ in \eqref{S:eq:3.23}. The conjugated domain in each case is
$\{g:X^{-1}g\in\mathcal D(\mathcal D_X)\}$.

On an interval with constant matrix $R_k$, the equation is $f'=-\lambda JR_kf/2$. Its solution is
\[
 f(t)=\exp\left(-\frac\lambda2(t-s)JR_k\right)f(s).
\]
Differentiating the exponential verifies the equation; uniqueness follows from the Volterra argument above. Canonical solutions join continuously at the partition points. For
$g=Xf$, the one-sided values satisfy
\[
 (X^-)^{-1}g^-=(X^+)^{-1}g^+.
\]
This is the transmission condition for $g$ at a jump of $X$; it expresses continuity of the canonical solution $f$.

Fix the cutoffs used below. Put $\eta_\beta=\min\{1/\beta,1/2\}$ and, for all sufficiently large $n$,
\begin{equation}\label{S:eq:cutoff-choice}
 \underline m_n=\lceil(\log n)^6\rceil,\qquad m_n=\lceil(\log n)^{512}\rceil,\qquad \underline T_n=1-\underline m_n/n,\qquad T_n=1-m_n/n.
\end{equation}

The following identities compare our normalization with \cite{valko2020operator} and \cite{valko2017sine}. Our deterministic sampling times for $\mathrm{Br}_{\beta,n}$ replace the
stopping times used for the circular coupling in \cite{valko2020operator}.

\begin{proposition}\label{S:prop:operator-normalization}
The differential and Brownian normalizations above agree with the circular and $\mathrm{Sine}_\beta$ operators of \cite{valko2017sine,valko2020operator}. A reversed finite jump of
index $j$ has beta parameter $\beta j/2$, and the periodically extended spectrum \eqref{S:eq:2.56} is invariant under translation by $2\pi n$.
\end{proposition}
\begin{proof}
The stochastic equations $\dd x=y\,\dd W^1$ and $\dd y=y\,\dd W^2$, with independent standard Brownian motions, give by It\^o's formula
\[
 \dd f(x,y)=yf_x\,\dd W^1+yf_y\,\dd W^2
                     +\tfrac12y^2(f_{xx}+f_{yy})\,\dd\theta.
\]
Thus the generator is $\Delta_{\mathbb H}/2$. For $\sigma^2=4/\beta$, the process $Z_{\sigma^2 v}$ has generator $\sigma^2\Delta_{\mathbb H}/2$; equivalently use the independent
Brownian motions $\sigma^{-1}W^i_{\sigma^2v}$ in its stochastic equations. Taking $v=-\log(1-t)$ gives exactly $\sigma^2v=\tau_\beta(t)$. The boundary lines remain $\infty$ and the
terminal point of that path.

Let $F(z)=\left(\begin{smallmatrix}1&-x\\0&y\end{smallmatrix}\right)$, so $\mathcal A(z)=y^{-1/2}F(z)$ and $\det F=y$.
The Hamiltonian in \cite{valko2017sine} is
\begin{equation}\label{S:eq:hamiltonian-normalization}
 H_{\mathrm{VV}}=\frac{F^{\mathsf T}F}{2\det F}
          =\tfrac12R_X,\qquad
 H_{\mathrm{VV}}^{-1}J\partial_t=2R_X^{-1}J\partial_t.
\end{equation}
The weighted norms satisfy $\|f\|_{H_{\mathrm{VV}}}^2=\|f\|_{R_X}^2/2$. Multiplication by $X/\sqrt2$ identifies the first weighted space with ordinary $L^2$, and multiplication by
$X$ identifies the second. The scalar factors cancel in the conjugated operator and in its inverse. Hence the spectral parameter and the factor $1/2$ in \eqref{S:eq:3.23} are
unchanged. For the representatives used here, $[(-1,0)^{\mathsf T},(q,1)^{\mathsf T}]=1$.

Set $z_{n,j}=b_{n-j}$, $0\le j\le n$. The jump $z_{n,j+1}\mapsto z_{n,j}$ is forward jump number $n-j-1$; its second beta parameter is
\[
 \frac\beta2\bigl(n-(n-j-1)-1\bigr)=\frac{\beta j}{2},
 \qquad 1\le j<n.
\]
Level zero is the harmonic terminal point. For an interior point of $I_{n,k}$, $n-k-1<nt<n-k$, and therefore
\begin{equation}\label{S:eq:subinterval-index-conversion}
 b_{\lfloor nt\rfloor}=b_{n-k-1}=z_{n,k+1},\qquad
 \lceil n(1-t)\rceil=k+1.
\end{equation}
The value on interval $k$ is therefore $z_{n,k+1}$, and its Brownian counterpart is sampled at $\vartheta_{n,k+1}$. Values at the finitely many partition points have no effect on
the $L^2$ coefficients or integral kernels.

For the choice $h_j=4/(\beta j+1)$,
\[
 \frac2{h_j}-\frac12=\frac{\beta j}{2}.
\]
This is the parameter relation in the radial comparison of \cite{valko2020operator}; the estimate with hypothesis $h_j\le1$ applies exactly when $\beta j\ge3$. The stronger local
estimate of Lemma~\ref{S:lem:7.1} has its own smaller fixed upper bound on $h_j$. Both restrictions hold uniformly for $j\ge\underline m_n$ for all sufficiently large $n$, since
$\underline m_n\to\infty$.

Finally let $R_k=\left(\begin{smallmatrix}a&b\\b&c\end{smallmatrix}\right)$ be a constant Hamiltonian on an interval of
length $1/n$. Since $ac-b^2=1$,
\[
 JR_k=\begin{pmatrix}-b&-c\\a&b\end{pmatrix},\qquad
 (JR_k)^2=(b^2-ac)I=-I.
\]
The transfer through this interval is therefore
\[
 T_k(\lambda)=\exp\left(-\frac\lambda{2n}JR_k\right)
 =\cos\left(\frac\lambda{2n}\right)I
       -\sin\left(\frac\lambda{2n}\right)JR_k,
\]
and
\begin{equation}\label{S:eq:aux-period-transfer}
 T_k(\lambda+2\pi n)=-T_k(\lambda).
\end{equation}
Let $\Phi_n(\lambda)$ be the ordered product of all $n$ transfers. An eigenfunction has initial value proportional to $u_0$, and its right boundary equation is
$[\Phi_n(\lambda)u_0,u_1]=0$. Replacing $\lambda$ by $\lambda+2\pi n$ multiplies this equation by $(-1)^n$, without changing its zeros. Equation~\eqref{S:eq:2.56} gives the
circular eigenangle law in one period. Replacing an angular representative by $\theta_{j,n}+2\pi$ leaves that set unchanged.
\end{proof}

\subsubsection{Measurability and conditioning of the sampled path}
\label{S:sec:s3-measurability}

The random inverse is determined measurably by the coefficients just specified. For a step function, its kernel is a finite combination of products of sampled matrix entries and
fixed rectangle or triangle indicators. Thus $\mathtt r\mathrm{Br}_{\beta,n}$ is measurable with respect to $\mathscr S_n=\sigma\{q_\infty,Z_{\vartheta_{n,j}}:1\le j\le n\}$, as
used in Section~\ref{S:sec:proof-strategy}. More generally, let $K_\omega(s,t)$ be a jointly measurable real matrix kernel with finite squared integral almost surely. Its norm is
measurable by Tonelli. Fix an orthonormal basis $(e_j)$ of $L^2((0,1)^2;\mathbb R^{2\times2})$. On the measurable set $\{\|K_\omega\|_2\le m\}$, Cauchy--Schwarz bounds
$|\langle e_j,K_\omega\rangle|$ by $m$. Measurable integration, first for truncated integrands and then by their pointwise limits, therefore makes every coefficient measurable.
Letting $m\to\infty$ and using Parseval gives
\[
 K_\omega=\lim_{N\to\infty}\sum_{j=1}^N
       \langle e_j,K_\omega\rangle e_j
 \quad\text{in }L^2((0,1)^2;\mathbb R^{2\times2}).
\]
The isometry from square-integrable kernels to Hilbert--Schmidt operators proves strong measurability. The limiting inverse is an instance of this argument, on the event
$I_X<\infty$; it may be set to zero on the null complement.

Let $p_h$ be the heat kernel relative to hyperbolic area and let $P_\eta$ be the normalized Poisson kernel from Lemma~\ref{S:lem:7.2}. For a fixed sampling row, choose finitely
many additional times in each interval between consecutive sampled positions. Under conditioning on the terminal point $\eta$, the transition densities are
$p_h(z,w)P_\eta(w)/P_\eta(z)$. Multiplying the densities within one interval makes all intermediate Poisson factors cancel. Dividing by the transition density between its two fixed
ends removes the remaining quotient. The resulting density is
\[
 \frac{\prod_{i=0}^l p_{t_{i+1}-t_i}(z_i,z_{i+1})}
      {p_h(z_0,z_{l+1})},
\]
the ordinary bridge density for that interval. The Markov factorization over the entire finite row groups these normalized densities into a product over its distinct intervals.
Consequently, for bounded cylinder functions $F_1,\ldots,F_M$ of those path segments,
\[
 \mathbb E\left[\prod_{j=1}^MF_j\,\middle|\,\mathscr S_n\right]
 =\prod_{j=1}^M\int F_j\,\dd\mathsf P_j,
\]
where $\mathsf P_j$ is the bridge law with the two sampled ends and the corresponding duration.

To extend the identity to the whole segments, identify each interval with $[0,1]$ by a deterministic linear change of time. Evaluation at rational times generates the Borel
sigma-field of the continuous path space. Cylinder sets with these times and a countable generating family of spatial Borel sets form a countable generating $\pi$-system. The
finite-dimensional identity holds simultaneously on this system outside one conditional null set. Fix all but one cylinder function and apply the monotone-class theorem to the
remaining coordinate; repeat over the finitely many segments. This gives the product identity for bounded Borel functions of the full segments, which is conditional independence
with the specified marginal bridge laws.

The construction is performed for each finite row. A countable intersection covers all row indices $n$, while different rows may share path segments and need not be independent.
The sampling sigma-field and the preserved original joint law remain exactly those used in the stated coupling.

\subsubsection{The columns $Xu_0$ and $Xu_1$}
\label{S:sec:s3-affine-columns}

For the fixed factor $\mathcal A$, the coefficient can also be expressed through the squared lengths of the boundary columns. In the reflected coordinate write
$Z_{\theta_\beta(r)}=x+\ii y$. Then
\[
 a_\omega(r)=(-y^{-1/2},0)^{\mathsf T},\qquad
 c_\omega(r)=((q_\infty-x)y^{-1/2},y^{1/2})^{\mathsf T}.
\]
Thus $p_\omega=|a_\omega|^2=y^{-1}$, $q_\omega=|c_\omega|^2=y+(q_\infty-x)^2/y$, and
\[
 \mathfrak h_{\mathrm{aff}}(a_\omega)=p_\omega/4,\qquad
 \mathfrak h_{\mathrm{aff}}(c_\omega)=q_\omega/4+1/p_\omega.
\]
With the primitives $A_\omega,C_\omega$ of the introduction this gives exactly
\begin{equation}\label{S:eq:s3-affine-scalar-coefficient}
 \Gamma_{\beta,\omega}(r)
 =\frac18A_\omega(r)p_\omega(r)
       +\frac18C_\omega(r)q_\omega(r)
       +\frac{C_\omega(r)}{2p_\omega(r)}.
\end{equation}
The first column is horizontal and $c_{\omega,2}^2=1/p_\omega$. These identities express $\Gamma_{\beta,\omega}$ in terms of the random squared column lengths. Substitution of
$p_\beta,q_\beta$ from the diagonal example gives the deterministic functional evaluated in Section~\ref{S:sec:s13-diagonal-example}.

Retain $\Delta_{\beta,n}$ from the \hyperref[S:sec:proof-strategy]{reduction to the error from sampling the Brownian path} and set
$D_{\beta,n}=\|\Delta_{\beta,n}\|_{\mathfrak S_2}^2$. With $A^T=P_TAP_T$ as above, the truncated differences are
\[
 \begin{split}
 \Delta_{\beta,n}^{T}
   &=\mathtt r\mathrm{Circ}_{\beta,n}^{\,T_n}
                          -\mathtt r\mathrm{Sine}_\beta^{\,T_n},\\
 \varepsilon_{\beta,n}^{\mathrm{fin},T}
   &=\mathtt r\mathrm{Circ}_{\beta,n}^{\,T_n}
                          -\mathtt r\mathrm{Br}_{\beta,n}^{\,T_n},\\
 \varepsilon_{\beta,n}^{\mathrm{samp},T}
   &=\mathtt r\mathrm{Br}_{\beta,n}^{\,T_n}
                          -\mathtt r\mathrm{Sine}_\beta^{\,T_n}.
 \end{split}
\]
These are compressions by the same projection. In particular, $\Delta_{\beta,n}^{T}=\varepsilon_{\beta,n}^{\mathrm{fin},T} +\varepsilon_{\beta,n}^{\mathrm{samp},T}$, with no
additional boundary term. The reduction in Section~\ref{S:sec:proof-strategy} uses this identity and the full kernel tails; the final deduction is in
Section~\ref{S:sec:proof-of-main-theorems}.

\subsection{Stability and Fr\'echet expansion of the inverse kernel}
\label{sec:weighted-stability-and-fr-echet-expansion-of-canonical-kernels}
\label{S:sec:stability-and-fr-echet-expansion-of-the-inverse-kernel}

The inverse is bilinear in its two boundary columns. To estimate a matrix perturbation, we integrate first in the unchanged variable. This gives a weighted bound for the perturbed
column that remains valid when the singular boundary solution has infinite $L^2$ norm.

We first prove the column identities for a general coefficient function. For the sampled coefficients, each integral includes the contribution from the current partition interval.
This contribution is needed even for the linear estimate. All kernel bounds are uniform in the reference matrices and the number of intervals. Equivalence with an unweighted norm
is used only after those data have been fixed.

\subsubsection{Differentiation of the Green kernel}
\label{S:sec:s3-kernel-variation}

Use the original coordinate $t$ on $(0,T)$ for the moment, and the real boundary columns $a=Xu_0$, $c=Xu_1$ of Section~\ref{S:sec:green-kernels-for-canonical-systems}. Initially
both columns lie in $L^2$. The weighted formulas below then remove that restriction, so that they also apply to the singular inverse. For column fields $v,w$, let $\mathcal K(v,w)$
denote the kernel $v(s)w(t)^{\mathsf T}/2$ on $s<t$ and its transposed value on $s>t$. Thus $K_X=\mathcal K(a,c)$, and $\mathcal K$ is bilinear. For column increments
$\alpha,\gamma$, define
\begin{equation}\label{S:eq:3.43}
 L_X[\alpha,\gamma](s,t) =\frac12[\alpha(s)c(t)^{\mathsf T} +a(s)\gamma(t)^{\mathsf T}]\mathbf1_{\{s<t\}} +\frac12[\gamma(s)a(t)^{\mathsf T} +c(s)\alpha(t)^{\mathsf T}]\mathbf1_{\{s>t\}},
\end{equation}
and
\begin{equation}\label{S:eq:3.44}
 Q[\alpha,\gamma](s,t)
 =\frac12[\alpha(s)\gamma(t)^{\mathsf T}\mathbf1_{\{s<t\}}
          +\gamma(s)\alpha(t)^{\mathsf T}\mathbf1_{\{s>t\}}].
\end{equation}
The same symbols denote their integral operators when the kernels are square integrable. Multiplication on $s<t$ gives
\[
 (a+\alpha)(s)(c+\gamma)(t)^{\mathsf T}
 -a(s)c(t)^{\mathsf T}
 =\alpha(s)c(t)^{\mathsf T}+a(s)\gamma(t)^{\mathsf T}
                          +\alpha(s)\gamma(t)^{\mathsf T}.
\]
Transposing gives the other triangle. Hence, exactly,
\begin{equation}\label{S:eq:s3-exact-column-expansion}
 \mathcal K(a+\alpha,c+\gamma)-\mathcal K(a,c)
 =L_X[\alpha,\gamma]+Q[\alpha,\gamma].
\end{equation}

The derivative contains a perturbation of $a$ at the first position and a perturbation of $c$ at the second. Integrating the unchanged column gives different weights for these two
terms. Write
\[
 P_X(t)=\int_0^t p(s)\,\dd s,\qquad
 Q_X(t)=\int_t^T q(s)\,\dd s.
\]
The rank-one inner product identity $\langle vw^{\mathsf T},xy^{\mathsf T}\rangle_{\mathrm{HS}} =\langle v,x\rangle\langle w,y\rangle$ gives
\begin{equation}\label{S:eq:s3-linear-norm-exact}
 \|L_X[\alpha,\gamma]\|_{\mathfrak S_2}^2 ={}\frac12\int_0^T [Q_X(t)|\alpha(t)|^2+P_X(t)|\gamma(t)|^2]\,\dd t +\int_{s<t}\langle\alpha(s),a(s)\rangle \langle c(t),\gamma(t)\rangle\,\dd s\,\dd t.
 \end{equation}
For example, the first term on $s<t$ is $|\alpha(s)|^2|c(t)|^2/4$; reflection across $s=t$ doubles its integral, giving the first weighted term with factor $1/2$. The other square
and the cross product give the remaining terms in \eqref{S:eq:s3-linear-norm-exact}. Pointwise,
\[
 2|\langle\alpha(s),a(s)\rangle
       \langle c(t),\gamma(t)\rangle|
 \le |\alpha(s)|^2|c(t)|^2+|a(s)|^2|\gamma(t)|^2.
\]
The cross integral is therefore absolutely convergent whenever the weighted integral in the next display is finite, and
\begin{equation}\label{S:eq:s3-linear-weighted-bound}
 \|L_X[\alpha,\gamma]\|_{\mathfrak S_2}^2
 \le\int_0^T[Q_X(t)|\alpha(t)|^2+P_X(t)|\gamma(t)|^2]\,\dd t.
\end{equation}
If the weighted integral is finite, Tonelli applies to the two square terms and the preceding inequality makes the cross term absolutely integrable. Thus
\eqref{S:eq:s3-linear-norm-exact} and \eqref{S:eq:s3-linear-weighted-bound} hold under this weighted assumption as well.

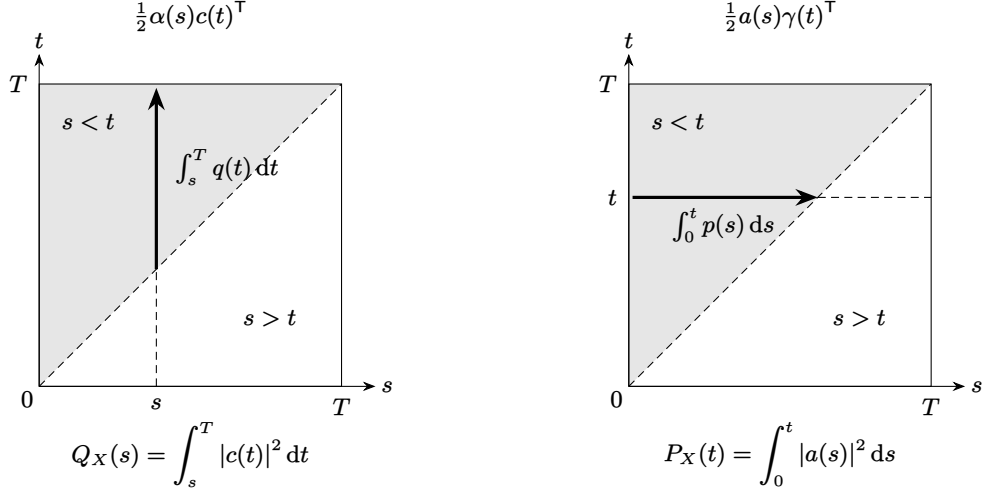
\begin{figure}[H]
\centering
\begin{tikzpicture}[x=1cm,y=1cm,every node/.style={font=\footnotesize,inner sep=2pt},>=Stealth]
 \foreach \xx in {0,7.8}{
  \begin{scope}[xshift=\xx cm]
   \fill[black!10] (0,0)--(0,4)--(4,4)--cycle;
   \draw[->] (0,0)--(4.45,0) node[right] {$s$};
   \draw[->] (0,0)--(0,4.4) node[above] {$t$};
   \draw (0,0) rectangle (4,4);
   \draw[densely dashed] (0,0)--(4,4);
   \node[below left] at (0,0) {$0$};
   \node[below] at (4,-.08) {$T$};
   \node[left] at (-.08,4) {$T$};
   \node at (.65,3.5) {$s<t$};
   \node at (3.05,.9) {$s>t$};
  \end{scope}
 }
 \node at (2,4.9) {$\frac12\alpha(s)c(t)^{\mathsf T}$};
 \draw[very thick,->] (1.55,1.55)--(1.55,3.95);
 \draw[densely dashed] (1.55,0)--(1.55,1.55);
 \node[below] at (1.55,-.08) {$s$};
 \node[anchor=west] at (1.72,2.9) {$\int_s^T q(t)\,\dd t$};
 \node at (2,-.9) {$Q_X(s)=\displaystyle\int_s^T|c(t)|^2\,\dd t$};
 \node at (9.8,4.9) {$\frac12 a(s)\gamma(t)^{\mathsf T}$};
 \draw[very thick,->] (7.85,2.5)--(10.3,2.5);
 \draw[densely dashed] (10.3,2.5)--(11.8,2.5);
 \node[left] at (7.72,2.5) {$t$};
 \node at (9.05,2.13) {$\int_0^t p(s)\,\dd s$};
 \node at (9.8,-.9) {$P_X(t)=\displaystyle\int_0^t|a(s)|^2\,\dd s$};
\end{tikzpicture}
\caption{Integration of the two terms in the first variation of the canonical Green kernel. Here $p=|a|^2$ and $q=|c|^2$. For fixed $s$, integration of $|\alpha(s)|^2|c(t)|^2$ over $s<t<T$ gives $Q_X(s)|\alpha(s)|^2$. For fixed $t$, integration of $|a(s)|^2|\gamma(t)|^2$ over $0<s<t$ gives $P_X(t)|\gamma(t)|^2$. Transposition supplies the triangle $s>t$, so the factor $1/4$ from each squared kernel becomes $1/2$. The cross product is the separate integral in \eqref{S:eq:s3-linear-norm-exact}, bounded there before obtaining \eqref{S:eq:s3-linear-weighted-bound}.}
\label{S:fig:kernel-integration-weights}
\end{figure}

The vertical and horizontal integrations in Figure~\ref{S:fig:kernel-integration-weights} give the two weights in \eqref{S:eq:s3-linear-weighted-bound}. The separate weights allow
the estimate to apply even when a boundary column is not square integrable.

On a regular interval, square-integrability of the boundary columns allows the following bilinear estimate in ordinary
$L^2$.

\begin{lemma}\label{S:lem:3.3}
For $\alpha,\gamma\in L^2((0,T);\mathbb R^2)$,
\begin{equation}\label{S:eq:3.45}
 \|Q[\alpha,\gamma]\|_{\mathfrak S_2}^2
 =\frac12\int_{s<t}|\alpha(s)|^2|\gamma(t)|^2\,\dd s\,\dd t
 \le\frac12\|\alpha\|_2^2\|\gamma\|_2^2.
\end{equation}
\end{lemma}
\begin{proof}
The squared norms on the two triangles are
\[
 \frac14\int_{s<t}|\alpha(s)|^2|\gamma(t)|^2\,\dd s\,\dd t,
 \qquad
 \frac14\int_{s>t}|\gamma(s)|^2|\alpha(t)|^2\,\dd s\,\dd t.
\]
Changing $s,t$ in the second integral proves the equality. Enlarging the integration region to the whole square proves the inequality. If the support of $\alpha$ lies strictly to
the left of the support of $\gamma$, the latter inequality is an equality; the constant $1/2$ cannot be reduced.
\end{proof}

In particular $\mathcal K:L^2\times L^2\to\mathfrak S_2$ is a bounded bilinear map with norm at most $1/\sqrt2$. Suppose $a_\varepsilon,c_\varepsilon$ satisfy
\begin{equation}\label{S:eq:column-variation-l2}
 \|a_\varepsilon-a-\varepsilon\dot a\|_2
 +\|c_\varepsilon-c-\varepsilon\dot c\|_2=o(|\varepsilon|),
 \qquad a,c,\dot a,\dot c\in L^2((0,T);\mathbb R^2).
\end{equation}
Here the norms are on the whole interval. Write $r_a=a_\varepsilon-a-\varepsilon\dot a$ and $r_c=c_\varepsilon-c-\varepsilon\dot c$. The derivative is
\begin{equation}\label{S:eq:3.35}
 \dot K=L_X[\dot a,\dot c].
\end{equation}
Indeed, \eqref{S:eq:s3-exact-column-expansion} gives
\[
 K_\varepsilon-K_X-\varepsilon\dot K
 =\mathcal K(r_a,c)+\mathcal K(a,r_c)
          +\mathcal K(a_\varepsilon-a,c_\varepsilon-c),
\]
so Lemma~\ref{S:lem:3.3} gives
\begin{equation}\label{S:eq:column-variation-hs-remainder}
 \|K_\varepsilon-K_X-\varepsilon\dot K\|_{\mathfrak S_2} \le\frac1{\sqrt2}\bigl(\|r_a\|_2\|c\|_2+\|a\|_2\|r_c\|_2 +\|a_\varepsilon-a\|_2\|c_\varepsilon-c\|_2\bigr) =o(|\varepsilon|).
\end{equation}
For arbitrary increments the quadratic remainder is bounded by $\|\alpha\|_2\|\gamma\|_2/\sqrt2$, which is at most $(\|\alpha\|_2+\|\gamma\|_2)^2/(4\sqrt2)$. Dividing by
$\|\alpha\|_2+\|\gamma\|_2$ proves Fr\'echet differentiability of the column map.

If second-order expansions hold in $L^2$,
\[
 a_\varepsilon=a+\varepsilon\dot a +\tfrac12\varepsilon^2\ddot a+o_{L^2}(\varepsilon^2),\qquad c_\varepsilon=c+\varepsilon\dot c +\tfrac12\varepsilon^2\ddot c+o_{L^2}(\varepsilon^2),
\]
then bilinearity yields
\begin{equation}\label{S:eq:column-second-order-expansion}
 K_\varepsilon-K_X
 =\varepsilon L_X[\dot a,\dot c]
 +\varepsilon^2\left(\tfrac12L_X[\ddot a,\ddot c]
                         +Q[\dot a,\dot c]\right)
 +o_{\mathfrak S_2}(\varepsilon^2).
\end{equation}
Write
\[
 \begin{aligned}
 a_\varepsilon&=a+\varepsilon\dot a+\tfrac12\varepsilon^2\ddot a+r_a,
 &\|r_a\|_2&=o(\varepsilon^2),\\
 c_\varepsilon&=c+\varepsilon\dot c+\tfrac12\varepsilon^2\ddot c+r_c,
 &\|r_c\|_2&=o(\varepsilon^2).
 \end{aligned}
\]
The bilinear map $\mathcal K(v,w)$ has norm at most $\|v\|_2\|w\|_2/\sqrt2$. After subtracting the constant, first-order and displayed second-order terms, the remaining kernel is
exactly
\[
 \begin{split}
 &\mathcal K(r_a,c)+\mathcal K(a,r_c)\\
 &+\tfrac12\varepsilon^3\mathcal K(\dot a,\ddot c)
   +\tfrac12\varepsilon^3\mathcal K(\ddot a,\dot c)
   +\tfrac14\varepsilon^4\mathcal K(\ddot a,\ddot c)\\
 &+\mathcal K(\varepsilon\dot a+\tfrac12\varepsilon^2\ddot a,r_c)
   +\mathcal K(r_a,\varepsilon\dot c+\tfrac12\varepsilon^2\ddot c)
   +\mathcal K(r_a,r_c).
 \end{split}
\]
The first line is $o(\varepsilon^2)$ in Hilbert--Schmidt norm. The second line is $O(|\varepsilon|^3)+O(\varepsilon^4)$. The first two terms of the last line are
$o(|\varepsilon|^3)$ and its last term is $o(\varepsilon^4)$. There are finitely many terms, so their sum is $o(\varepsilon^2)$. In particular the coefficient of $\varepsilon^2$ is
$\tfrac12L_X[\ddot a,\ddot c]+Q[\dot a,\dot c]$; both terms involving the second derivatives of the columns must be retained to obtain the second derivative. The estimate is in the
norm on the entire fixed regular interval and does not replace weighted estimates when a column fails to belong to $L^2$.

All relative perturbations below are written on the left: $Y=\me^UX$. In particular the column derivatives at $U=0$ are $Ua$ and $Uc$. For a right perturbation $Y=X\me^V$,
$X\me^VX^{-1}=\me^{XVX^{-1}}$, term by term in the exponential series. On a fixed neighborhood of the identity where the matrix logarithm is defined, this gives
\begin{equation}\label{S:eq:2.13e}
 U=\operatorname{Ad}_XV=XVX^{-1}.
\end{equation}
The logarithm is real and traceless there: it is real by its convergent power series, while $1=\det(\me^U)=\me^{\operatorname{tr}U}$ implies $\operatorname{tr}U=0$ for real $U$.

For a bounded measurable $U:(0,T)\to\mathfrak{sl}_2(\mathbb R)$, write
\[
 \eta=\operatorname*{ess\,sup}_{t\in(0,T)}\|U(t)\|_{\mathrm{op}},
 \qquad
 \mathfrak e_X(U)=\int_0^T
       [Q_X(t)|U(t)a(t)|^2+P_X(t)|U(t)c(t)|^2]\,\dd t.
\]
When $I_X<\infty$, $\mathfrak e_X$ is a finite nonnegative quadratic form on the real vector space $L^\infty((0,T);\mathfrak{sl}_2(\mathbb R))$. Its associated symmetric bilinear
form is $(\mathfrak e_X(U+V)-\mathfrak e_X(U-V))/4$ on the product of this space with itself. Let $\mathscr D_XU=L_X[Ua,Uc]$. Then
\begin{equation}\label{S:eq:s3-exponential-estimates}
 \begin{split}
 \|\mathscr D_XU\|_{\mathfrak S_2}^2&\le\mathfrak e_X(U),\\
 \|K_{\me^UX}-K_X\|_{\mathfrak S_2}^2&\le C\mathfrak e_X(U),\\
 \|K_{\me^UX}-K_X-\mathscr D_XU\|_{\mathfrak S_2}^2
     &\le C\eta^2\mathfrak e_X(U),\qquad \eta\le1,
 \end{split}
\end{equation}
with an absolute constant. For any vector $v$ the exponential series, with $|U^jv|\le\eta^{j-1}|Uv|$, gives
\begin{equation}\label{S:eq:continuous-exponential-vector-estimates}
 |(\me^U-I)v|\le\me^\eta|Uv|,\qquad
 |(\me^U-I-U)v|\le\tfrac12\eta\me^\eta|Uv|.
\end{equation}
Set $\delta a=(\me^U-I)a$, $\delta c=(\me^U-I)c$, $r_a=(\me^U-I-U)a$, and $r_c=(\me^U-I-U)c$. The exact remainder in \eqref{S:eq:s3-exact-column-expansion} is
\[
 K_{\me^UX}-K_X-\mathscr D_XU=L_X[r_a,r_c]+Q[\delta a,\delta c].
\]
The first term has norm at most $\eta\me^\eta\mathfrak e_X(U)^{1/2}/2$ by \eqref{S:eq:s3-linear-weighted-bound}. For the second,
\[
 \begin{split}
 \|Q[\delta a,\delta c]\|_{\mathfrak S_2}^2
 &=\frac12\int_{s<t}|\delta a(s)|^2|\delta c(t)|^2\,\dd s\,\dd t\\
 &\le\frac{\eta^2\me^{4\eta}}2
       \int_0^T Q_X(s)|U(s)a(s)|^2\,\dd s
 \le\frac{\eta^2\me^{4\eta}}2\mathfrak e_X(U).
 \end{split}
\]
The triangle inequality gives the last bound in \eqref{S:eq:s3-exponential-estimates}, and addition of the linear term gives the middle bound. The quadratic estimate uses
$|\delta c(t)|\le\eta\me^\eta|c(t)|$ and retains the weighted integral in the other variable. On a regular interval, these estimates also prove Fr\'echet differentiability in
$L^\infty$, since
\[
 \mathfrak e_X(U)
 \le2\eta^2\int_0^T q(t)P_X(t)\,\dd t.
\]
Under \eqref{S:eq:singular-integrability-assumptions}, Tonelli gives $\int Q_Xp=\int P_Xq=I_X$. Thus $\mathfrak e_X(U)\le2I_X\eta^2$, and the remainder divided by $\eta$ is at most
$C\sqrt{2I_X}\eta$, which tends to zero. This proves the same differentiability at the singular boundary; separate $L^2$ norms of $a$ and $c$ are not used.

We next identify the perturbed kernel with the inverse on the stated differential domain. For bounded $U$, $Y=\me^UX$ satisfies
\[
 \me^{-\eta}|Xv|\le |Yv|\le\me^\eta|Xv|\qquad(v\in\mathbb R^2).
\]
Thus the integrability or nonintegrability of each constant direction is unchanged; moreover $Q_Y\le\me^{2\eta}Q_X$ and $I_Y\le\me^{4\eta}I_X$.
Theorems~\ref{S:thm:singular-inverse-domain} and~\ref{S:thm:3.2} identify the kernels in \eqref{S:eq:s3-exponential-estimates} with the appropriate inverses, with the same
coordinate boundary lines.

\subsubsection{Changes in the matrix factor and boundary vectors}
\label{S:sec:s5-boundary-variations}

Boundary variations can be included on a regular interval, or when $\int_0^1\operatorname{tr}R<\infty$. Let $X_\varepsilon=\me^{\varepsilon U}X$ with $U$ bounded, and let
$u_{i,\varepsilon}=u_i+\varepsilon\dot u_i+o(\varepsilon)$ be real, normalized boundary vectors. Differentiating their normalization gives $[\dot u_0,u_1]+[u_0,\dot u_1]=0$, and
\[
 \dot a=Ua+X\dot u_0,\qquad \dot c=Uc+X\dot u_1.
\]
To check \eqref{S:eq:column-variation-l2}, write, for example,
\[
 \me^{\varepsilon U}Xu_{0,\varepsilon}-a -\varepsilon(Ua+X\dot u_0) ={}(\me^{\varepsilon U}-I-\varepsilon U)a +\varepsilon(\me^{\varepsilon U}-I)X\dot u_0 +\me^{\varepsilon U}X (u_{0,\varepsilon}-u_0-\varepsilon\dot u_0).
 \]
Put $C_X=(\int\operatorname{tr}R)^{1/2}$ and $\eta=\|U\|_\infty$. For $i=0,1$, set $r_{i,\varepsilon}=u_{i,\varepsilon}-u_i-\varepsilon\dot u_i$. The three terms of the same
expansion with index $i$ satisfy
\[
 \|\me^{\varepsilon U}Xu_{i,\varepsilon}-Xu_i -\varepsilon(UXu_i+X\dot u_i)\|_2 \le C_X\me^{|\varepsilon|\eta} \left(\tfrac12\varepsilon^2\eta^2|u_i| +\varepsilon^2\eta|\dot u_i|+|r_{i,\varepsilon}|\right) =o(|\varepsilon|).
 \]
Thus both column variations satisfy \eqref{S:eq:column-variation-l2}, and Lemma~\ref{S:lem:3.3} gives \eqref{S:eq:3.35}. At a limit-point boundary point, a boundary variation would
require the new constant direction to belong to the weighted space; that direction may have infinite weighted norm.

The multiplication order for relative products follows from the same left convention. If $H_\ell(s)=\me^{sD_\ell}H_\ell(0)$ and $P_{k,j}(s)=H_k(s)\cdots H_{j-1}(s)$, with
$P_{k,k}=I$, the ordinary product rule gives
\[
 \partial_sP_{k,j}
 =\sum_{\ell=k}^{j-1}P_{k,\ell}D_\ell H_\ell\cdots H_{j-1}.
\]
Multiplication on the right by $P_{k,j}^{-1}$ cancels the factors after $D_\ell$, and hence
\[
 (\partial_sP_{k,j})P_{k,j}^{-1}
   =\sum_{\ell=k}^{j-1}\operatorname{Ad}_{P_{k,\ell}}D_\ell.
\]
Likewise, from $X_k=H_kX_{k+1}$, $\widetilde H_k=\me^{\zeta_k}H_k$, and $\widetilde X_k=\me^{U_k}X_k$, one obtains
\[
 \me^{U_k}
 =\widetilde X_kX_k^{-1}
 =\me^{\zeta_k}H_k\me^{U_{k+1}}H_k^{-1}
 =\me^{\zeta_k}\me^{\operatorname{Ad}_{H_k}U_{k+1}}.
\]

\subsubsection{The sampled weights and the space of matrix fields}
\label{S:sec:s5-weighted-spaces}

The column estimates suggest a weighted space for perturbations of a sampled matrix factor. Its weights must include the current interval as well as the intervals strictly before
or after it; otherwise the linear kernel estimate would miss part of a diagonal square.

Fix integers $n>m\ge1$, put $a_*=m/n$, and work in the reflected coordinate $r=1-t$ on $I_*=(a_*,1)$. Write
\[
 I_k=(k/n,(k+1)/n),\qquad r=(k+s)/n,\quad 0<s<1,\quad m\le k<n.
\]
Values at partition points have no effect on the operators. Let $X(r)=X_k\in\mathrm{SL}_2(\mathbb R)$ on $I_k$, and use the normalization $[u_0,u_1]=1$ of
Section~\ref{S:sec:green-kernels-for-canonical-systems}. Set
\[
 a_k=X_ku_0,\qquad c_k=X_ku_1,\qquad p_k=|a_k|^2,\qquad q_k=|c_k|^2.
\]
The corresponding step functions are denoted by $a,c,p,q$. In the application, $m=m_n$ and $X_k=\mathcal A(Z_{\vartheta_{n,k+1}})$. Thus $k$ labels the operator interval, whose
coefficient is sampled at recursion index $k+1$.

Use the weights already introduced in Section~\ref{S:sec:sampling-times-and-estimates-for-the-boundary-solutions}:
\[
 A_{n,k}=\frac1n\sum_{\ell=m}^{k-1}q_\ell,\qquad
 C_{n,k}=\frac1n\sum_{\ell=k+1}^{n-1}p_\ell,\qquad
 A_{n,k}^+=A_{n,k}+\frac{q_k}{n},\qquad
 C_{n,k}^+=C_{n,k}+\frac{p_k}{n}.
\]
We sometimes omit $n$ from these subscripts. All matrix fields in this section are real. For $U(r)=U_k(s)$ with $U_k\in L^2((0,1);\mathbb R^{2\times2})$, define the nonnegative
quadratic forms
\begin{equation}\label{S:eq:continuous-field-energies}
 \begin{aligned}
 \mathcal E_n(U)&=\sum_{k=m}^{n-1}\int_0^1
       \bigl[A_{n,k}|U_k(s)a_k|^2+C_{n,k}|U_k(s)c_k|^2\bigr]\,\dd s,\\
 \mathcal E_n^{\mathrm{diag}}(U)&=\frac1n\sum_{k=m}^{n-1}\int_0^1
       \bigl[q_k|U_k(s)a_k|^2+p_k|U_k(s)c_k|^2\bigr]\,\dd s,\\
 \mathcal E_n^+(U)&=\sum_{k=m}^{n-1}\int_0^1
       \bigl[A_{n,k}^+|U_k(s)a_k|^2+C_{n,k}^+|U_k(s)c_k|^2\bigr]\,\dd s.
 \end{aligned}
\end{equation}
The associated symmetric bilinear form of each expression $\mathcal E$ is $(\mathcal E(U+V)-\mathcal E(U-V))/4$ on the same real vector space. For fields constant on each interval,
the integrals over $s$ reduce to the corresponding finite sums. In either case,
\begin{equation}\label{S:eq:weighted-energy-decomposition}
 \mathcal E_n^+(U)=\mathcal E_n(U)+\mathcal E_n^{\mathrm{diag}}(U).
\end{equation}
The added diagonal term accounts for $I_k\times I_k$. To compute its contribution, write the reflected kernel as $a(r)c(r')^{\mathsf T}/2$ for $r>r'$. Integration in $r'$ when
$a(r)$ is perturbed gives $\int_{a_*}^r q$; integration in $r$ when $c(r')$ is perturbed gives $\int_{r'}^1p$. Thus define \phantomsection\label{S:lem:continuous-primitives}
\[
 A_X(r)=\int_{a_*}^r q(v)\,\dd v,\qquad
 C_X(r)=\int_r^1 p(v)\,\dd v.
\]
Splitting each integral into complete partition intervals and the remaining part of $I_k$ gives
\begin{equation}\label{S:eq:continuous-inside-subinterval-primitives}
 A_X((k+s)/n)=A_{n,k}+\frac{s q_k}{n},\qquad
 C_X((k+s)/n)=C_{n,k}+\frac{(1-s)p_k}{n}.
\end{equation}
On the same domain define the nonnegative quadratic form
\[
 \mathfrak e_{n,X}(U)=\int_{I_*}
       \bigl[A_X(r)|U(r)a(r)|^2+C_X(r)|U(r)c(r)|^2\bigr]\,\dd r.
\]
The substitution $\dd r=\dd s/n$ yields
\begin{equation}\label{S:eq:continuous-integral-energy}
 \begin{aligned}
 \mathfrak e_{n,X}(U)
 &=\frac1n\sum_k\int_0^1\left[
       \left(A_{n,k}+\frac{s q_k}{n}\right)|U_k(s)a_k|^2
      +\left(C_{n,k}+\frac{(1-s)p_k}{n}\right)|U_k(s)c_k|^2
                              \right]\,\dd s\\
 &\le\frac1n\mathcal E_n^+(U).
 \end{aligned}
\end{equation}
For step fields, integration of $s$ and $1-s$ gives the sharper identity
\begin{equation}\label{S:eq:diagonal-triangles-integral}
 \mathfrak e_{n,X}(U)=\frac1n\left(
       \mathcal E_n(U)+\tfrac12\mathcal E_n^{\mathrm{diag}}(U)\right).
\end{equation}
For a field varying inside $I_k$, the factors $s$ and $1-s$ in \eqref{S:eq:continuous-integral-energy} must remain inside the integral.

The contribution from the current interval cannot be dropped. There is no bound of $\mathcal E_n^+$ by $\mathcal E_n$ even for a fixed partition. For example, set $U_k=0$ for
$k\ne m$ and
\[
 S_m=(a_m\ c_m),\qquad U_m=S_m\begin{pmatrix}0&0\\1&0\end{pmatrix}S_m^{-1}.
\]
Then $U_ma_m=c_m$, $U_mc_m=0$, and $\operatorname{tr}U_m=0$. Since $A_{n,m}=0$, one has $\mathcal E_n(U)=0$ but $\mathcal E_n^{\mathrm{diag}}(U)=|c_m|^4/n>0$. The first variation
of the Green kernel is $c_mc_m^{\mathsf T}/2$ on $I_m\times I_m$ and zero off that square. Hence its squared Hilbert--Schmidt norm is $|c_m|^4/(4n^2)$. This also shows why
$\mathcal E_n$ alone cannot control the derivative.

Let $\mathfrak H_{n,X}^+$ be the space of measurable $\mathfrak{sl}_2(\mathbb R)$-valued fields with finite $\mathcal E_n^+$, modulo equality almost everywhere. Its inner product
is
\begin{equation}\label{S:eq:weighted-inner-product}
 \langle U,V\rangle_{n,+}
 =\sum_{k=m}^{n-1}\int_0^1
 \bigl[A_{n,k}^+\langle U_k(s)a_k,V_k(s)a_k\rangle
      +C_{n,k}^+\langle U_k(s)c_k,V_k(s)c_k\rangle\bigr]\,\dd s.
\end{equation}
Write $\|U\|_{n,+}=\mathcal E_n^+(U)^{1/2}$; for step fields this is also the norm denoted by $\|U\|_{\mathcal E_n^+}$. The optional constraint $U_k(s)c_k\in\mathbb Rc_k$ defines a
linear subspace. Set
\[
 D_k=A_{n,k}^+a_ka_k^{\mathsf T}+C_{n,k}^+c_kc_k^{\mathsf T},\qquad
 \mathcal H_0=\bigoplus_{k=m}^{n-1}L^2((0,1);\mathbb R^{2\times2}).
\]
The same definitions apply to any nonempty subset of the indices, with the sums for the weights restricted to that subset in its natural order.

\begin{lemma}\label{S:lem:weighted-hilbert}
The real symmetric bilinear form \eqref{S:eq:weighted-inner-product} is a positive definite inner product on $\mathfrak H_{n,X}^+$ and makes it a separable real Hilbert space. The
step fields and the fields satisfying $U_kc_k\in\mathbb Rc_k$ are closed subspaces. Multiplication by $D_k^{1/2}$ gives an isometry into $\mathcal H_0$, measurable in the reference
columns.
\end{lemma}
\begin{proof}
Since $[a_k,c_k]=1$, the matrix $S_k=(a_k\ c_k)$ has determinant $-1$. Moreover $A_{n,k}^+\ge q_k/n>0$ and $C_{n,k}^+\ge p_k/n>0$. Thus
\[
 D_k=S_k\operatorname{diag}(A_{n,k}^+,C_{n,k}^+)S_k^{\mathsf T}>0,
\]
and multiplication of these matrices gives
\begin{equation}\label{S:eq:weighted-gram-determinant}
 \det D_k=A_{n,k}^+C_{n,k}^+>0,\qquad
 \mathcal E_n^+(U)=\sum_k\int_0^1
                 \operatorname{tr}(U_k(s)D_kU_k(s)^{\mathsf T})\,\dd s.
\end{equation}
In particular $\mathcal E_n^+(U)=0$ implies $U_ka_k=U_kc_k=0$ almost everywhere, and hence $U_k=0$ almost everywhere. For this fixed finite list put
\[
 d_- =\min_k\lambda_{\min}(D_k)>0,\qquad
 d_+ =\max_k\lambda_{\max}(D_k)<\infty.
\]
The inequalities $d_-I\le D_k\le d_+I$ imply
\begin{equation}\label{S:eq:weighted-fixed-frame-equivalence}
 d_-\sum_k\int_0^1\|U_k(s)\|_{\mathrm{HS}}^2\,\dd s
 \le\mathcal E_n^+(U)
 \le d_+\sum_k\int_0^1\|U_k(s)\|_{\mathrm{HS}}^2\,\dd s.
\end{equation}
For step fields the integrals may be omitted. The constants $d_-$ and $d_+$ may depend on $n$ and $X$. Their dependence can be read without choosing eigenvectors:
\[
 \lambda_{\min}(D_k)\ge
 \frac{A_{n,k}^+C_{n,k}^+}{A_{n,k}^+p_k+C_{n,k}^+q_k}>0,
 \qquad \lambda_{\max}(D_k)\le A_{n,k}^+p_k+C_{n,k}^+q_k.
\]
The first inequality follows from $\lambda_{\min}\lambda_{\max} =\det D_k$ and $\lambda_{\max}\le\operatorname{tr}D_k$. These bounds establish completeness for a fixed list. They
are not used as uniform constants in the kernel estimates.

The trace constraint is the kernel of a continuous linear map on $\mathcal H_0$. So is $[U_kc_k,c_k]=0$, since $|[Vc_k,c_k]|\le |c_k|^2\|V\|_{\mathrm{HS}}$. These kernels are
closed. The step fields form a finite-dimensional closed subspace. Completeness and separability now follow from those of $\mathcal H_0$ and
\eqref{S:eq:weighted-fixed-frame-equivalence}.

Define
\begin{equation}\label{S:eq:weighted-space-isometry}
 (\mathcal J_nU)_k(s)=U_k(s)D_k^{1/2}.
\end{equation}
Equation \eqref{S:eq:weighted-gram-determinant} proves that this map is an isometry. Its inverse on its image is right multiplication by $D_k^{-1/2}$, a bounded map for this fixed
finite list, so the image is closed. To check measurability without choosing eigenvectors, for a positive definite $2\times2$ matrix $D$ one has
\[
 D^{1/2}=\frac{D+\sqrt{\det D}\,I}
                   {\sqrt{\operatorname{tr}D+2\sqrt{\det D}}}.
\]
Put $d=\sqrt{\det D}$. Cayley--Hamilton gives $(D+dI)^2=(\operatorname{tr}D+2d)D$. The displayed matrix therefore squares to $D$; it is positive definite because $D+dI>0$ and its
denominator is positive. Both this expression and its inverse are continuous in $D$. Consequently $\mathcal J_n$ and its inverse on its image are measurable whenever the columns
and weights are.
\end{proof}

The scalar weight from Section~\ref{S:sec:sampling-times-and-estimates-for-the-boundary-solutions} therefore has the extension
\begin{equation}\label{S:eq:weighted-scalar-weight}
 \mathfrak G_{n,k}^+=A_{n,k}^+p_k+C_{n,k}^+q_k
                 =\mathfrak G_{n,k}+\frac2n p_kq_k.
\end{equation}
Since $|Va_k|^2\le p_k\|V\|_{\mathrm{op}}^2$ and $|Vc_k|^2\le q_k\|V\|_{\mathrm{op}}^2$,
\[
 \mathcal E_n^+(U)\le\sum_k\mathfrak G_{n,k}^+
                         \int_0^1\|U_k(s)\|_{\mathrm{op}}^2\,\dd s.
\]

\subsubsection{The first variation and the nonlinear remainder}
\label{S:sec:s5-first-variation}

The weighted space now permits a norm estimate for the linear variation and its remainder. Reflection changes the order of the kernel variables but not the Hilbert--Schmidt norm,
so the expansion can be written in the coordinate used for the sampling intervals.

Write $K_X(r,r')$ for the reflected Green kernel and $\mathsf K_X$ for its integral operator, extended by zero outside $I_*^2$. Thus
\[
 K_X(r,r')=\tfrac12a(r)c(r')^{\mathsf T}\quad(r>r'),\qquad
 K_X(r',r)=K_X(r,r')^{\mathsf T}.
\]
Reflection is unitary on the common $L^2$ space. We use the same symbols $L_X[\alpha,\gamma]$ and $Q[\alpha,\gamma]$ for the reflected kernels in
\eqref{S:eq:s3-exact-column-expansion}. On $r>r'$ their values are
\[
 \tfrac12\bigl[\alpha(r)c(r')^{\mathsf T}
                       +a(r)\gamma(r')^{\mathsf T}\bigr],
 \qquad \tfrac12\alpha(r)\gamma(r')^{\mathsf T},
\]
respectively; on $r<r'$ take the transposed values with the arguments interchanged.

Define $\mathscr D_{n,X}U=L_X[Ua,Uc]$, abbreviated to $\mathscr D_nU$. Its kernel is
\begin{equation}\label{S:eq:continuous-derivative-kernel}
 \begin{aligned}
 \dot K_X[U](r,r')={}&\tfrac12\bigl[
       (U(r)a(r))c(r')^{\mathsf T}
          +a(r)(U(r')c(r'))^{\mathsf T}\bigr]\mathbf1_{\{r>r'\}}\\
 &+\tfrac12\bigl[(U(r)c(r))a(r')^{\mathsf T}
          +c(r)(U(r')a(r'))^{\mathsf T}\bigr]\mathbf1_{\{r<r'\}}.
 \end{aligned}
\end{equation}

\begin{proposition}\label{S:prop:6.1}
The map $\mathscr D_{n,X}:\mathfrak H_{n,X}^+\to\mathfrak S_2(\mathcal H)$ is bounded and linear, and
\begin{equation}\label{S:eq:6.21}
 \|\mathscr D_{n,X}U\|_{\mathfrak S_2}^2
 \le\mathfrak e_{n,X}(U)\le\frac1n\mathcal E_n^+(U).
\end{equation}
Its values are self-adjoint operators.
\end{proposition}
\begin{proof}
Under $r=1-t$, the primitives $Q_X(t)$ and $P_X(t)$ in \eqref{S:eq:s3-linear-weighted-bound} become $A_X(r)$ and $C_X(r)$. Hence, for arbitrary real column fields with finite right
side,
\begin{equation}\label{S:eq:triangular-linear-estimate}
 \|L_X[\alpha,\gamma]\|_{\mathfrak S_2}^2
 \le\int_{I_*}\bigl[A_X(r)|\alpha(r)|^2+C_X(r)|\gamma(r)|^2\bigr]\,\dd r.
\end{equation}
Taking $\alpha=Ua$, $\gamma=Uc$ and using \eqref{S:eq:continuous-integral-energy} proves
\begin{equation}\label{S:eq:continuous-frechet-bound}
 \|\mathscr D_{n,X}U\|_{\mathfrak S_2}^2
 \le\mathfrak e_{n,X}(U)\le n^{-1}\mathcal E_n^+(U).
\end{equation}
For later local estimates, the pointwise bounds on the two triangles are
\begin{equation}\label{S:eq:6.23}
 \|\dot K_X[U](r,r')\|_{\mathrm{HS}}^2
 \le\tfrac12|U(r)a(r)|^2q(r')+\tfrac12p(r)|U(r')c(r')|^2
 \quad(r>r'),
\end{equation}
\begin{equation}\label{S:eq:6.24}
 \|\dot K_X[U](r,r')\|_{\mathrm{HS}}^2
 \le\tfrac12|U(r)c(r)|^2p(r')+\tfrac12q(r)|U(r')a(r')|^2
 \quad(r<r').
\end{equation}
They follow from $|vw^{\mathsf T}|_{\mathrm{HS}}=|v||w|$ and $|A+B|^2\le2|A|^2+2|B|^2$. The integrals of these two majorants have sum at most $\mathfrak e_{n,X}(U)<\infty$.
Changing $U$ on a null set changes the kernel only on the union of two null strips, so the operator is well-defined on equivalence classes. For $f,g\in\mathcal H$,
\[
 \iint |f(r)|\,\|\dot K_X[U](r,r')\|_{\mathrm{HS}}\,|g(r')|
 \,\dd r\,\dd r'
 \le\|\dot K_X[U]\|_{L^2}\|f\|_2\|g\|_2<\infty.
\]
Fubini and transposition across the diagonal now prove self-adjointness. The displayed kernel is linear in $U$. The estimate holds for all real matrix fields with finite weighted
norm. Tracelessness ensures that $\me^UX$ remains in $\mathrm{SL}_2(\mathbb R)$ in the nonlinear map.
\end{proof}

For bounded $U$ write
\[
 \eta(U)=\operatorname*{ess\,sup}_{r\in I_*}\|U(r)\|_{\mathrm{op}},
 \qquad Y=\me^UX,\qquad
 \alpha=Ua,\quad\gamma=Uc,
\]
\[
 \rho^a=(\me^U-I-U)a,\qquad \rho^c=(\me^U-I-U)c.
\]
For step fields $\eta(U)$ is the same as $\eta_n(U)=\max_k\|U_k\|_{\mathrm{op}}$. Applying \eqref{S:eq:continuous-exponential-vector-estimates} to the individual columns preserves
$|Ua|$ and $|Uc|$ in the weighted integrals. In particular,
\begin{equation}\label{S:eq:6.15}
 |\rho^a(r)|\le\tfrac12\eta\me^\eta|U(r)a(r)|,
\end{equation}
\begin{equation}\label{S:eq:6.16}
 |\rho^c(r)|\le\tfrac12\eta\me^\eta|U(r)c(r)|.
\end{equation}
These bounds keep the same column weights as the linear estimate.

Using the bilinear map $Q$ defined above, put
\[
 \mathscr Q_n(U)=Q[\alpha,\gamma],\qquad \mathscr N_n(U)=L_X[\rho^a,\rho^c] +Q[\alpha,\rho^c]+Q[\rho^a,\gamma]+Q[\rho^a,\rho^c].
\]
Expanding $(a+\alpha+\rho^a)(c+\gamma+\rho^c)^{\mathsf T}$ on each ordered triangle proves the exact identity
\begin{equation}\label{S:eq:6.19}
 \mathsf K_Y-\mathsf K_X
       =\mathscr D_nU+\mathscr Q_n(U)+\mathscr N_n(U).
\end{equation}
The term $\mathscr Q_n$ is the product of the two linear column increments. The other quadratic terms come from $\rho^a,\rho^c$. For a step field, the two triangular values of
$\mathscr Q_n$ on $I_k\times I_k$ differ by
\[
 \tfrac12U_k(a_kc_k^{\mathsf T}-c_ka_k^{\mathsf T})U_k^{\mathsf T}
                 =\tfrac12(\det U_k)J.
\]
Consequently $\mathscr Q_n$ and $\mathscr N_n$ are kept as triangular kernels separately. Their sum has matching values on the two diagonal triangles, as proved in
Lemma~\ref{S:lem:exact-mesh-kernel-stability}.

\begin{proposition}\label{S:prop:6.2}
For every bounded field $U$,
\begin{equation}\label{S:eq:6.27}
 \|\mathscr Q_n(U)\|_{\mathfrak S_2}^2
       \le\frac{\eta(U)^2}{2}\mathfrak e_{n,X}(U)
       \le\frac{\eta(U)^2}{2n}\mathcal E_n^+(U).
\end{equation}
\end{proposition}
\begin{proof}
The two transposed triangular kernels give the exact identity
\[
 \|\mathscr Q_n(U)\|_{\mathfrak S_2}^2
        =\frac12\int_{r>r'}|U(r)a(r)|^2|U(r')c(r')|^2\,\dd r\,\dd r'.
\]
In the second factor use $|U(r')c(r')|\le\eta|c(r')|$. Thus
\begin{equation}\label{S:eq:6.28}
 \|\mathscr Q_n(U)\|_{\mathfrak S_2}^2
       \le\frac{\eta^2}{2}\int_{I_*}A_X(r)|U(r)a(r)|^2\,\dd r.
\end{equation}
Alternatively, bounding the first factor by $\eta|a(r)|$ gives
\begin{equation}\label{S:eq:6.29}
 \|\mathscr Q_n(U)\|_{\mathfrak S_2}^2
       \le\frac{\eta^2}{2}\int_{I_*}C_X(r)|U(r)c(r)|^2\,\dd r.
\end{equation}
Using $\mathfrak e_{n,X}(U)\le n^{-1}\mathcal E_n^+(U)$ in either bound gives \eqref{S:eq:6.27}. For a step field the exact norm includes
\[
 \frac1{2n^2}\sum_{i>j}|U_ia_i|^2|U_jc_j|^2
       +\frac1{4n^2}\sum_i|U_ia_i|^2|U_ic_i|^2.
\]
Each diagonal triangle has area $1/(2n^2)$, giving the last coefficient. \end{proof}

\begin{proposition}\label{S:prop:6.3}
If $\eta(U)\le1$, then
\[
 \|\mathscr N_n(U)\|_{\mathfrak S_2}^2
       \le C\eta(U)^2\mathfrak e_{n,X}(U)
       \le\frac{C\eta(U)^2}{n}\mathcal E_n^+(U),
\]
with an absolute constant $C$.
\end{proposition}
\begin{proof}
Put $b=\eta\me^\eta/2$. Equations \eqref{S:eq:6.15}--\eqref{S:eq:6.16} and \eqref{S:eq:triangular-linear-estimate} give
\[
 \|L_X[\rho^a,\rho^c]\|_{\mathfrak S_2}
                     \le b\,\mathfrak e_{n,X}(U)^{1/2}.
\]
For the three remaining terms the same ordered integral gives, respectively,
\[
 \begin{aligned}
 \|Q[\alpha,\rho^c]\|_{\mathfrak S_2}^2
 &\le\tfrac12b^2\eta^2\int_{I_*}A_X|Ua|^2,\\
 \|Q[\rho^a,\gamma]\|_{\mathfrak S_2}^2
 &\le\tfrac12b^2\eta^2\int_{I_*}A_X|Ua|^2,\\
 \|Q[\rho^a,\rho^c]\|_{\mathfrak S_2}^2
 &\le\tfrac12b^4\eta^2\int_{I_*}A_X|Ua|^2.
 \end{aligned}
\]
For example, the first inequality uses $|\rho^c(r')|\le b|U(r')c(r')|\le b\eta|c(r')|$. The second uses $|\rho^a(r)|\le b|U(r)a(r)|$ and $|U(r')c(r')|\le\eta|c(r')|$; using both
remainder bounds gives the third. Taking square roots and adding proves
\[
 \|\mathscr N_n(U)\|_{\mathfrak S_2}
 \le\left(b+\sqrt2b\eta+\frac{b^2\eta}{\sqrt2}\right)
                          \mathfrak e_{n,X}(U)^{1/2}.
\]
For $0\le\eta\le1$, the coefficient is at most $C_0\eta$, where $C_0=\me/2+\me/\sqrt2+\me^2/(4\sqrt2)$. Squaring gives $C_0^2\eta^2\mathfrak e_{n,X}(U)$ and then
$C_0^2\eta^2\mathcal E_n^+(U)/n$ by \eqref{S:eq:continuous-integral-energy}. In particular, the constant contains no norm of a boundary column or inverse reference matrix.
\end{proof}

Combining these estimates gives the following result for $Y=\me^UX$, for either step fields or bounded measurable fields.

\begin{theorem}\label{S:thm:6.4}
If $\eta(U)\le1$, then, with an absolute constant $C$,
\begin{equation}\label{S:eq:6.31}
 \|\mathsf K_Y-\mathsf K_X\|_{\mathfrak S_2}^2
                              \le\frac Cn\mathcal E_n^+(U),
\end{equation}
\begin{equation}\label{S:eq:6.32}
 \|\mathsf K_Y-\mathsf K_X-\mathscr D_nU\|_{\mathfrak S_2}^2
                              \le\frac{C\eta(U)^2}{n}\mathcal E_n^+(U).
\end{equation}
\end{theorem}
\begin{proof}
By \eqref{S:eq:6.19} and Propositions~\ref{S:prop:6.2}--\ref{S:prop:6.3},
\[
 \|\mathsf K_Y-\mathsf K_X-\mathscr D_nU\|_{\mathfrak S_2}
 \le\|\mathscr Q_n(U)\|_{\mathfrak S_2}
       +\|\mathscr N_n(U)\|_{\mathfrak S_2}
 \le C\eta\,\mathfrak e_{n,X}(U)^{1/2}.
\]
This and \eqref{S:eq:continuous-integral-energy} give \eqref{S:eq:6.32}. Adding the derivative, whose norm is at most $\mathfrak e_{n,X}(U)^{1/2}$, gives
$\|\mathsf K_Y-\mathsf K_X\|\le(1+C\eta)\mathfrak e_{n,X}(U)^{1/2}$. Squaring and using $\eta\le1$ proves \eqref{S:eq:6.31}.
\end{proof}

Now allow $m$, $X$, and $U=U^{(n)}$ to vary with $n$. The assumptions for the next statement are
\begin{equation}\label{S:eq:6.33}
 \mathcal E_n^+(U^{(n)})=O(1),
\end{equation}
\begin{equation}\label{S:eq:6.34}
 \eta(U^{(n)})\longrightarrow0.
\end{equation}
Set $Y^{(n)}=\me^{U^{(n)}}X^{(n)}$. For random fields the assumptions and conclusions below are interpreted on each path on which these two bounds hold.

\begin{theorem}\label{S:thm:6.5}
Under \eqref{S:eq:6.33}--\eqref{S:eq:6.34},
\begin{equation}\label{S:eq:6.35}
 \sqrt n\,\|\mathsf K_{Y^{(n)}}-\mathsf K_{X^{(n)}}
                         -\mathscr D_nU^{(n)}\|_{\mathfrak S_2}\longrightarrow0,
\end{equation}
\begin{equation}\label{S:eq:6.36}
 \|\mathsf K_{Y^{(n)}}-\mathsf K_{X^{(n)}}\|_{\mathfrak S_2}^2
       =\|\mathscr D_nU^{(n)}\|_{\mathfrak S_2}^2+o(n^{-1}).
\end{equation}
\end{theorem}
\begin{proof}
Write $D_n=\mathscr D_nU^{(n)}$ and $R_n=\mathsf K_{Y^{(n)}}-\mathsf K_{X^{(n)}}-D_n$. For all large $n$,
\[
 \sqrt n\,\|D_n\|_{\mathfrak S_2}\le\mathcal E_n^+(U^{(n)})^{1/2},
 \qquad
 \sqrt n\,\|R_n\|_{\mathfrak S_2}
       \le C\eta(U^{(n)})\mathcal E_n^+(U^{(n)})^{1/2}.
\]
The second quantity tends to zero and the first is bounded. Expanding the square and applying Cauchy--Schwarz gives
\[
 n\bigl|\|D_n+R_n\|_{\mathfrak S_2}^2-\|D_n\|_{\mathfrak S_2}^2\bigr| \le2(\sqrt n\|D_n\|_{\mathfrak S_2})(\sqrt n\|R_n\|_{\mathfrak S_2}) +n\|R_n\|_{\mathfrak S_2}^2\longrightarrow0.
\]
Choose $M<\infty$ and $n_0$ such that $\mathcal E_n^+(U^{(n)})\le M$ for $n\ge n_0$. For these $n$, the last right side is at most $M[2C\eta(U^{(n)})+C^2\eta(U^{(n)})^2]$, which
tends to zero. For random fields fix an outcome on which $M<\infty$ and $\eta(U^{(n)})\to0$. The same inequality applies to each such outcome.
\end{proof}

\subsubsection{Matrix coefficients of the kernel difference}
\label{S:sec:s5-matrix-coefficients}

Lemma~\ref{S:lem:regular-interval-restriction} identifies the compressed kernels with the finite regular inverses. The next lemma gives finite matrix representations for the kernel
difference and its first variation. Work on $I_*=(m/n,1)$ in the reflected coordinate. In the next lemma $W$ is constant on each $I_k$, $Y_k=\me^{W_k}X_k$, and
$\eta=\max_k\|W_k\|_{\mathrm{op}}$. Write $\mathscr D_{n,X}^{\mathrm{step}}$ for the restriction of $\mathscr D_{n,X}$ to step fields.

\begin{lemma}\label{S:lem:exact-mesh-kernel-stability}
The operators $\mathsf K_Y-\mathsf K_X$ and $\mathscr D_{n,X}^{\mathrm{step}}W$ are represented exactly by finite matrices on the space of functions constant on each $I_k$ and
vanish on its orthogonal complement. If $\eta\le1$, they satisfy
\begin{equation}\label{S:eq:inverse-kernel-perturbation-bound}
 \begin{aligned}
 \|\mathscr D_{n,X}^{\mathrm{step}}W\|_{\mathfrak S_2}^2
                  &\le n^{-1}\mathcal E_n^+(W),\\
 \|\mathsf K_Y-\mathsf K_X\|_{\mathfrak S_2}^2
                  &\le Cn^{-1}\mathcal E_n^+(W),\\
 \|\mathsf K_Y-\mathsf K_X-\mathscr D_{n,X}^{\mathrm{step}}W\|_{\mathfrak S_2}^2
                  &\le Cn^{-1}\eta^2\mathcal E_n^+(W).
 \end{aligned}
\end{equation}
\end{lemma}
\begin{proof}
Put $\widetilde a_k=\me^{W_k}a_k$ and $\widetilde c_k=\me^{W_k}c_k$. Since $\operatorname{tr}W_k=0$, $\det\me^{W_k}=1$, and therefore
\[
 \widetilde a_k\widetilde c_k^{\mathsf T}
       -\widetilde c_k\widetilde a_k^{\mathsf T}
       =a_kc_k^{\mathsf T}-c_ka_k^{\mathsf T}=J.
\]
For $k>\ell$ define
\[
 \Delta_{k\ell}=\tfrac12
       (\widetilde a_k\widetilde c_\ell^{\mathsf T}-a_kc_\ell^{\mathsf T}),
 \qquad
 D_{k\ell}(W)=\tfrac12\bigl[(W_ka_k)c_\ell^{\mathsf T}
                                  +a_k(W_\ell c_\ell)^{\mathsf T}\bigr].
\]
For $k<\ell$ take the transposes with the indices interchanged. For $k=\ell$ use the same displayed expressions. The normalization identity proves
$\Delta_{kk}=\Delta_{kk}^{\mathsf T}$. For the first variation, direct multiplication of a $2\times2$ matrix gives
\[
 W_kJ+JW_k^{\mathsf T}=(\operatorname{tr}W_k)J=0,
\]
so $D_{kk}(W)=D_{kk}(W)^{\mathsf T}$ as well. Thus both triangular values on $I_k\times I_k$ coincide for the difference and for its first variation. Equivalently, in block
notation,
\begin{equation}\label{S:eq:6.17}
 \mathscr D_nW(i,j)=\tfrac12\bigl[(W_ia_i)c_j^{\mathsf T} +a_i(W_jc_j)^{\mathsf T}\bigr]\mathbf1_{\{i>j\}} +\tfrac12\bigl[(W_ic_i)a_j^{\mathsf T} +c_i(W_ja_j)^{\mathsf T}\bigr]\mathbf1_{\{i\le j\}}.
\end{equation}
The identities $\Delta_{kk}=\Delta_{kk}^{\mathsf T}$ and $D_{kk}(W)=D_{kk}(W)^{\mathsf T}$ give symmetric diagonal blocks, so each represents both triangular values.

Let $\iota_n:\bigoplus_{k=m}^{n-1}\mathbb R^2\to\mathcal H$ be the isometry taking $(v_k)$ to $\sqrt n\,v_k$ on $I_k$ and to zero elsewhere. Its adjoint is
\[
 (\iota_n^*f)_k=\sqrt n\int_{I_k}f(r)\,\dd r.
\]
An operator with constant kernel block $B_{k\ell}$ acts on $I_k$ as $\sum_\ell B_{k\ell}\int_{I_\ell}f$. Consequently
\begin{equation}\label{S:eq:kernel-difference-matrix}
 \mathsf K_Y-\mathsf K_X
       =\iota_n(\Delta_{k\ell}/n)_{k,\ell}\iota_n^*,\qquad
 \mathscr D_{n,X}^{\mathrm{step}}W
       =\iota_n(D_{k\ell}(W)/n)_{k,\ell}\iota_n^*.
\end{equation}
If $f$ is orthogonal to all step functions, each integral over $I_k$ vanishes, showing that the operator vanishes on the orthogonal complement. For the orthonormal basis
$\sqrt n\,\mathbf1_{I_k}e_\mu$, $\mu=1,2$, the matrix elements are the entries of $B_{k\ell}/n$. Therefore
\[
 \left\|\iota_n(B_{k\ell}/n)\iota_n^*\right\|_{\mathfrak S_2}^2
 =\frac1{n^2}\sum_{k,\ell}\sum_{\mu,\nu=1}^2
                         |(B_{k\ell})_{\mu\nu}|^2
 =\frac1{n^2}\sum_{k,\ell}\|B_{k\ell}\|_{\mathrm{HS}}^2.
\]
Here $B$ is either $\Delta$ or $D(W)$; their diagonal blocks have already been identified with both triangular values.

The first estimate in \eqref{S:eq:inverse-kernel-perturbation-bound} follows from Proposition~\ref{S:prop:6.1}, with the more precise diagonal weight
\eqref{S:eq:diagonal-triangles-integral}. Theorem~\ref{S:thm:6.4} gives the other two. All three bounds include the diagonal triangles. The individual Green inverse has infinite
rank because it is injective on the retained, infinite-dimensional $L^2$ space.
\end{proof}

\subsubsection{Fr\'echet differentiability for measurable fields}
\label{S:sec:s5-measurable-frechet}

The linear estimate requires only a weighted square integral. The nonlinear map requires more: exponentiating a square-integrable matrix field can produce nonintegrable columns. We
therefore use the intersection with $L^\infty$ for the nonlinear map,
\[
 \mathfrak B_{n,X}=\mathfrak H_{n,X}^+
             \cap L^\infty(I_*;\mathfrak{sl}_2(\mathbb R)),\qquad
 \|U\|_{\mathfrak B_{n,X}}=\|U\|_{n,+}+\eta(U).
\]
This is a Banach space. Indeed, a Cauchy sequence converges in the weighted Hilbert norm and in $L^\infty$. By \eqref{S:eq:weighted-fixed-frame-equivalence} the former convergence
implies ordinary $L^2$ convergence for this fixed $X$; the latter also implies $L^2$ convergence since $|I_*|<\infty$. The limits therefore agree and lie in the intersection.

For $U\in\mathfrak B_{n,X}$ write $a_U=\me^Ua$, $c_U=\me^Uc$, and $\mathsf K_X(U)=\mathsf K_{\me^UX}$. Its kernel is
\begin{equation}\label{S:eq:continuous-perturbed-kernel}
 K_X(U;r,r')=\tfrac12\bigl[
 a_U(r)c_U(r')^{\mathsf T}\mathbf1_{\{r>r'\}}
       +c_U(r)a_U(r')^{\mathsf T}\mathbf1_{\{r<r'\}}\bigr].
\end{equation}
In particular $\mathsf K_X(0)=\mathsf K_X$.

\begin{theorem}\label{S:thm:continuous-frechet}
The map $\mathsf K_X:\mathfrak B_{n,X}\to\mathfrak S_2(\mathcal H)$ is Fr\'echet differentiable at zero with derivative $\mathscr D_{n,X}$. For $\eta(U)\le1$,
\begin{equation}\label{S:eq:continuous-frechet-quantitative}
 \begin{aligned}
 \|\mathsf K_X(U)-\mathsf K_X(0)\|_{\mathfrak S_2}^2
                     &\le Cn^{-1}\mathcal E_n^+(U),\\
 \|\mathsf K_X(U)-\mathsf K_X(0)-\mathscr D_{n,X}U\|_{\mathfrak S_2}^2
                     &\le Cn^{-1}\eta(U)^2\mathcal E_n^+(U).
 \end{aligned}
\end{equation}
\end{theorem}
\begin{proof}
To identify \eqref{S:eq:continuous-perturbed-kernel} as an inverse, put $Y(t)=\me^{U(1-t)}X(1-t)$ on $[0,1-a_*]$. Since the reference matrices form a fixed finite list,
\[
 \|Y\|_\infty\le\me^{\eta(U)}\max_k\|X_k\|_{\mathrm{op}},\qquad
 \|Y^{-1}\|_\infty\le\me^{\eta(U)}\max_k\|X_k^{-1}\|_{\mathrm{op}}.
\]
Also $\det Y=1$. The version of Theorem~\ref{S:thm:3.2} for measurable coefficients thus applies with the unchanged boundary vectors $u_0,u_1$. Its Green formula, followed by
reflection and extension by zero, is \eqref{S:eq:continuous-perturbed-kernel}.

Theorem~\ref{S:thm:6.4} applies to this $Y$ and proves both estimates in \eqref{S:eq:continuous-frechet-quantitative}. In particular,

\[
 \frac{\|\mathsf K_X(U)-\mathsf K_X(0)-\mathscr D_{n,X}U\|_{\mathfrak S_2}}
            {\|U\|_{\mathfrak B_{n,X}}}
 \le\frac C{\sqrt n}\,
       \frac{\eta(U)\|U\|_{n,+}}{\|U\|_{n,+}+\eta(U)}
 \le\frac C{\sqrt n}\eta(U)\longrightarrow0
\]
as $U\to0$ in $\mathfrak B_{n,X}$. Proposition~\ref{S:prop:6.1} proves boundedness of the linear map on that Banach space and on all of $\mathfrak H_{n,X}^+$.
\end{proof}

In particular Theorem~\ref{S:thm:6.5} applies to measurable fields: if $\eta(U_n)\to0$ and $\mathcal E_n^+(U_n)=O(1)$, then
\begin{equation}\label{S:eq:continuous-frechet-remainder}
 \sqrt n\,\|\mathsf K_{X_n}(U_n)-\mathsf K_{X_n}(0)
                          -\mathscr D_{n,X_n}U_n\|_{\mathfrak S_2}\longrightarrow0,
\end{equation}
The difference of the squared norms is $o(n^{-1})$, by the same square expansion as in Theorem~\ref{S:thm:6.5}. These estimates are uniform in the reference matrices.

The nonlinear map need not be defined without the boundedness assumption. For example, take $X=I$, $u_0=(-1,0)^{\mathsf T}$, $u_1=(0,1)^{\mathsf T}$ and
\[
 U(r)=\operatorname{diag}(f(r),-f(r)),\qquad f(r)=(1-r)^{-1/4}.
\]
For a fixed $n,m$ the weights are bounded and $f\in L^2(I_*)$, so $U\in\mathfrak H_{n,X}^+$. Choose $a_*<b<c<1$. The proposed squared kernel norm is bounded below by
\[
 \frac12\left(\int_{a_*}^{b}\me^{-2f(r')}\,\dd r'\right)
             \left(\int_c^1\me^{2f(r)}\,\dd r\right)=\infty.
\]
The first integral is strictly positive, and the second diverges: with $v=(1-r)^{-1/4}$ it is a positive multiple of a tail integral of $\me^{2v}v^{-5}$. Hence the exponential map
does not in general define a Hilbert--Schmidt kernel on the whole weighted Hilbert space.

\subsubsection{The second derivative}
\label{S:sec:s5-second-derivative}

At second order there are two contributions: the second derivative of the matrix exponential in each column, and the product of the first variations of the two columns. Their sum
gives the symmetric bilinear map below.

For $U,V\in\mathfrak B_{n,X}$ set $P(U,V)=(UV+VU)/2$ pointwise. Define the symmetric bilinear map
\[
 \mathscr H_{n,X}(U,V)
   =L_X[P(U,V)a,P(U,V)c]+Q[Ua,Vc]+Q[Va,Uc].
\]
On $r>r'$ its kernel is
\begin{equation}\label{S:eq:continuous-second-variation-kernel}
 \begin{aligned}
 \tfrac12\bigl[&(P(U,V)(r)a(r))c(r')^{\mathsf T}
                +a(r)(P(U,V)(r')c(r'))^{\mathsf T}\\
              &+(U(r)a(r))(V(r')c(r'))^{\mathsf T}
                +(V(r)a(r))(U(r')c(r'))^{\mathsf T}\bigr].
 \end{aligned}
\end{equation}
On the other triangle use transposition and interchange the arguments.

\begin{proposition}\label{S:prop:continuous-second-variation}
The map $\mathsf K_X$ is twice continuously Fr\'echet differentiable on $\mathfrak B_{n,X}$, with second derivative $\mathscr H_{n,X}$ at zero. Moreover
\begin{equation}\label{S:eq:continuous-second-variation-bound}
 \|\mathscr H_{n,X}(U,V)\|_{\mathfrak S_2}
 \le\frac C{\sqrt n}
           \bigl[\eta(U)\|V\|_{n,+}+\eta(V)\|U\|_{n,+}\bigr],
\end{equation}
and, for $\eta(U)\le1$,
\begin{equation}\label{S:eq:continuous-cubic-remainder}
 \|\mathsf K_X(U)-\mathsf K_X(0)-\mathscr D_{n,X}U
                   -\tfrac12\mathscr H_{n,X}(U,U)\|_{\mathfrak S_2}
               \le\frac C{\sqrt n}\eta(U)^2\|U\|_{n,+}.
\end{equation}
\end{proposition}
\begin{proof}
For every vector $v$,
\[
 |P(U,V)v|\le\tfrac12\bigl[\eta(U)|Vv|+\eta(V)|Uv|\bigr].
\]
Minkowski's inequality in the weighted integrals and \eqref{S:eq:triangular-linear-estimate} therefore give
\[
 \|L_X[P(U,V)a,P(U,V)c]\|_{\mathfrak S_2}
 \le\tfrac12\bigl[\eta(U)\mathfrak e_{n,X}(V)^{1/2}
                      +\eta(V)\mathfrak e_{n,X}(U)^{1/2}\bigr].
\]
For the two product terms,
\[
 \begin{aligned}
 \|Q[Ua,Vc]\|_{\mathfrak S_2}^2
 &=\frac12\int_{r>r'}|U(r)a(r)|^2|V(r')c(r')|^2\,\dd r\,\dd r'\\
 &\le\frac{\eta(V)^2}{2}\int_{I_*}A_X(r)|U(r)a(r)|^2\,\dd r
 \le\frac{\eta(V)^2}{2}\mathfrak e_{n,X}(U),
 \end{aligned}
\]
Interchanging $U,V$ gives $\|Q[Va,Uc]\|_{\mathfrak S_2}\le\eta(U)\mathfrak e_{n,X}(V)^{1/2}/\sqrt2$. Adding the three norm bounds and applying
\eqref{S:eq:continuous-integral-energy} proves \eqref{S:eq:continuous-second-variation-bound}.

To compute the cubic remainder, retain $\alpha=Ua$, $\gamma=Uc$ and $\rho^a,\rho^c$ from the first-order calculation, and set
\[
 s_a=(\me^U-I-U-U^2/2)a,\qquad
 s_c=(\me^U-I-U-U^2/2)c.
\]
Taylor's formula with integral remainder gives
\[
 (\me^U-I-U-U^2/2)v
       =\frac12\int_0^1(1-t)^2\me^{tU}U^3v\,\dd t,
 \qquad |s_a|\le\tfrac16\me^\eta\eta^2|Ua|,
\]
with the same bound for $s_c$. Since $P(U,U)=U^2$, the exact remainder after subtracting $\tfrac12\mathscr H_{n,X}(U,U)$ is
\begin{equation}\label{S:eq:s5-cubic-exact}
 L_X[s_a,s_c]+Q[\alpha,\rho^c]+Q[\rho^a,\gamma]+Q[\rho^a,\rho^c].
\end{equation}
The first term has norm at most $\me^\eta\eta^2\mathfrak e_{n,X}(U)^{1/2}/6$. With $b=\eta\me^\eta/2$, the last three have the respective norm bounds
\[
 \frac{b\eta}{\sqrt2}\mathfrak e_{n,X}(U)^{1/2},\qquad
 \frac{b\eta}{\sqrt2}\mathfrak e_{n,X}(U)^{1/2},\qquad
 \frac{b^2\eta}{\sqrt2}\mathfrak e_{n,X}(U)^{1/2},
\]
as calculated in Proposition~\ref{S:prop:6.3}. For $\eta\le1$ their sum is at most $C\eta^2\mathfrak e_{n,X}(U)^{1/2}$, proving \eqref{S:eq:continuous-cubic-remainder}.

For the Fr\'echet derivative, fix $v=a$ or $c$ and put $F_\ell(U)=U^\ell v/\ell!$ in $L^2(I_*;\mathbb R^2)$. The reference list is finite, so $v\in L^2$. The derivatives of this
polynomial are
\[
 F_\ell'(U)[V]=\frac1{\ell!}\sum_{i+j=\ell-1}U^iVU^jv,
\]
\[
 F_\ell''(U)[V,W]=\frac1{\ell!}\sum_{i+j+k=\ell-2}
           (U^iVU^jWU^k+U^iWU^jVU^k)v.
\]
On $\eta(U)\le R$, $R\ge1$, their operator norms satisfy
\[
 \begin{aligned}
 \|F_\ell(U)\|_2&\le\|v\|_2R^\ell/\ell!,\\
 \|F_\ell'(U)[V]\|_2
       &\le\|v\|_2\eta(V)R^{\ell-1}/(\ell-1)!,\\
 \|F_\ell''(U)[V,W]\|_2
       &\le\|v\|_2\eta(V)\eta(W)R^{\ell-2}/(\ell-2)!.
 \end{aligned}
\]
The last two sums contain respectively $\ell$ and $\ell(\ell-1)$ words. Terms with $\ell<1$ or $\ell<2$ are zero. Since $\eta(V)\le\|V\|_{\mathfrak B_{n,X}}$, all three series
converge uniformly on bounded Banach balls, the latter two in their operator norms. Write $F=\sum_\ell F_\ell$, $S_1=\sum_\ell F_\ell'$, and $S_2=\sum_\ell F_\ell''$. Uniform
convergence on each ball makes $S_1$ and $S_2$ continuous in their respective operator norms. For a finite sum, the fundamental theorem of calculus along $U+tH$ gives the
identities below. Uniform convergence on a ball containing that segment passes both identities to the limit:
\[
 F(U+H)-F(U)=\int_0^1S_1(U+tH)[H]\dd t,
\]
\[
 S_1(U+H)[V]-S_1(U)[V]=\int_0^1S_2(U+tH)[H,V]\dd t.
\]
For $H\ne0$, subtracting the proposed linear terms gives
\[
 \frac{\|F(U+H)-F(U)-S_1(U)[H]\|_2}{\|H\|_{\mathfrak B_{n,X}}}
 \le\sup_{0\le t\le1}\|S_1(U+tH)-S_1(U)\|_{\mathrm{op}}\longrightarrow0,
\]
\[
 \frac{\|S_1(U+H)-S_1(U)-S_2(U)[H,\cdot]\|_{\mathrm{op}}}
      {\|H\|_{\mathfrak B_{n,X}}}
 \le\sup_{0\le t\le1}\|S_2(U+tH)-S_2(U)\|_{\mathrm{op}}\longrightarrow0.
\]
Both suprema tend to zero by continuity at $U$, so $DF=S_1$ and $D^2F=S_2$ are continuous in operator norm.

Finally, the bilinear kernel map $(v,w)\mapsto Q[v,w]$ has norm at most $1/\sqrt2$ from $L^2\times L^2$ to $\mathfrak S_2$, by Lemma~\ref{S:lem:3.3} after reflection. The product
rule therefore applies to $\mathsf K_X(U)=Q[\me^Ua,\me^Uc]$. At zero the first column derivative is $Uv$, and the second is $(UV+VU)v/2$. The four terms of the product rule are
exactly \eqref{S:eq:continuous-second-variation-kernel}.
\end{proof}

\subsubsection{Subdivision and approximation}
\label{S:sec:s5-subdivision}

The estimates must remain valid when a matrix field is approximated on a finer partition. We compare the weights before and after subdivision to keep the constants uniform in the
relative cell lengths.

Let $\pi$ be a finite subdivision of the intervals $I_k$, with no restriction on the relative lengths of its members. For $Q=(b,b+h)\subset I_k$ write $a_Q=a_k$, $c_Q=c_k$ and put
\[
 \overline A_Q=A_X(b+h),\qquad \overline C_Q=C_X(b).
\]
On the real vector space of matrix fields constant on each $Q\in\pi$, define the nonnegative quadratic form
\[
 \mathfrak E_{\pi,X}(U)=\sum_{Q\in\pi}|Q|
       \bigl[\overline A_Q|U_Qa_Q|^2
                       +\overline C_Q|U_Qc_Q|^2\bigr].
\]
The weights contain the lengths $|Q|$ of the subintervals in the coordinate $r$.

\begin{proposition}\label{S:prop:continuous-refinement}
For every field constant on $\pi$,
\begin{equation}\label{S:eq:continuous-refinement-energy}
 \begin{aligned}
 \mathfrak E_{\pi,X}(U)-\mathfrak e_{n,X}(U)
   &=\frac12\sum_{Q\in\pi}|Q|^2
          \bigl[|c_Q|^2|U_Qa_Q|^2+|a_Q|^2|U_Qc_Q|^2\bigr]\ge0,\\
 \mathfrak E_{\pi,X}(U)&\le\frac1n\mathcal E_n^+(U).
 \end{aligned}
\end{equation}
The kernel estimates hold with the same constants for every such subdivision.
\end{proposition}
\begin{proof}
On $Q=(b,b+h)$ the exact primitives are
\[
 A_X(r)=A_X(b)+(r-b)|c_Q|^2,\qquad
 C_X(r)=C_X(b+h)+(b+h-r)|a_Q|^2.
\]
Consequently
\[
 h\overline A_Q-\int_Q A_X(r)\,\dd r=\tfrac12h^2|c_Q|^2,\qquad h\overline C_Q-\int_Q C_X(r)\,\dd r=\tfrac12h^2|a_Q|^2.
\]
Multiplying by $|U_Qa_Q|^2$, $|U_Qc_Q|^2$ and summing proves the first identity. If $Q\subset I_k$, then $\overline A_Q\le A_{n,k}^+$ and $\overline C_Q\le C_{n,k}^+$. Summing
these inequalities times $|Q|$ gives
\[
 \mathfrak E_{\pi,X}(U)
 \le\int_{I_*}\bigl[A_{n,k(r)}^+|U(r)a(r)|^2
                       +C_{n,k(r)}^+|U(r)c(r)|^2\bigr]\,\dd r
 =n^{-1}\mathcal E_n^+(U).
\]
Here $k(r)$ is the parent interval index; the last equality is exactly $\dd r=\dd s/n$ on each parent interval. Every kernel estimate in
Propositions~\ref{S:prop:6.1}--\ref{S:prop:6.3} and Proposition~\ref{S:prop:continuous-second-variation} was proved with $\mathfrak e_{n,X}$ before replacing it by
$n^{-1}\mathcal E_n^+$. Since $\mathfrak e_{n,X}(U)\le\mathfrak E_{\pi,X}(U)$, substitution gives those bounds with $\mathfrak E_{\pi,X}$ and the same constants. No quotient of the
smallest and largest subdivision lengths enters: each discarded contribution from a single interval is precisely $|Q|^2/2$ times the displayed nonnegative expression. Thus the
bound remains uniform for nonuniform subdivisions.
\end{proof}

Let $\pi_j$ be subdivisions with $|\pi_j|=\max_{Q\in\pi_j}|Q|\to0$, and let $U^{(j)}$ be the interval average of
$U\in\mathfrak H_{n,X}^+$:
\[
 U^{(j)}(r)=\frac1{|Q|}\int_Q U(v)\,\dd v\quad(r\in Q\in\pi_j).
\]
These integrals exist since, for fixed $n,X$, the field is in ordinary $L^2$. Averaging preserves tracelessness and, when present, the linear condition $Uc_k\in\mathbb Rc_k$ on the
parent interval. The averages satisfy
\begin{equation}\label{S:eq:continuous-simple-approximation}
 \|U^{(j)}-U\|_{n,+}\longrightarrow0,\qquad
 \|\mathscr D_{n,X}U^{(j)}-\mathscr D_{n,X}U\|_{\mathfrak S_2}
                                                    \longrightarrow0.
\end{equation}
If $\eta(U)<\infty$, they also satisfy $\eta(U^{(j)})\le\eta(U)$ and $\mathsf K_X(U^{(j)})\to\mathsf K_X(U)$ in
$\mathfrak S_2$.

The matrix $D_k$ is constant on every $Q\subset I_k$. Jensen's inequality in the matrix Hilbert space gives
\[
 |Q|\left\|\left(\frac1{|Q|}\int_Q U\right)D_k^{1/2}\right\|_{\mathrm{HS}}^2
 \le\int_Q\|U(r)D_k^{1/2}\|_{\mathrm{HS}}^2\,\dd r.
\]
Multiplication by $n$ and summation show that interval averaging is a contraction in $\|\cdot\|_{n,+}$. Choose a field $V$ continuous on the closure of each parent interval, with
$\|U-V\|_{n,+}<\varepsilon$. Such fields are dense by ordinary $L^2$ approximation on a finite union of intervals and \eqref{S:eq:weighted-fixed-frame-equivalence}; approximation
can be performed in the fixed finite-dimensional matrix constraint subspace. Writing $V^{(j)}$ for its interval average, one has
\[
 \sup_{r\in I_k}\|V^{(j)}(r)-V(r)\|_{\mathrm{HS}}
 \le\sup_{\substack{r,v\in\overline I_k\\|r-v|\le|\pi_j|}}
                       \|V(r)-V(v)\|_{\mathrm{HS}}\longrightarrow0.
\]
There are only finitely many parent intervals, and their weights are fixed. Hence $\|V^{(j)}-V\|_{n,+}\to0$, and contraction gives
\[
 \|U^{(j)}-U\|_{n,+}
 \le2\|U-V\|_{n,+}+\|V^{(j)}-V\|_{n,+}.
\]
Taking $\limsup_j$ and then $\varepsilon\downarrow0$ proves the first limit in \eqref{S:eq:continuous-simple-approximation}. The second follows from
$\|\mathscr D_{n,X}(U^{(j)}-U)\|_{\mathfrak S_2} \le n^{-1/2}\|U^{(j)}-U\|_{n,+}$. Convexity of the matrix operator norm proves $\eta(U^{(j)})\le\eta(U)$.

For the nonlinear limit, let $A,B$ have norms at most $R$. Differentiating $\me^{(1-t)A}\me^{tB}$ and integrating gives
\[
 \me^A-\me^B=\int_0^1\me^{(1-t)A}(A-B)\me^{tB}\,\dd t,
 \qquad
 \|\me^A-\me^B\|_{\mathrm{HS}}\le\me^R\|A-B\|_{\mathrm{HS}}.
\]
By \eqref{S:eq:weighted-fixed-frame-equivalence}, $U^{(j)}\to U$ in ordinary $L^2(I_*)$. Since $a,c$ are bounded for this fixed reference list, the last inequality with $R=\eta(U)$
gives $a_{U^{(j)}}\to a_U$ and $c_{U^{(j)}}\to c_U$ in $L^2$. Their $L^2$ norms are uniformly bounded. Bilinearity and Lemma~\ref{S:lem:3.3}, after reflection, give
\[
 \|\mathsf K_X(U^{(j)})-\mathsf K_X(U)\|_{\mathfrak S_2} \le\frac1{\sqrt2}\bigl(\|a_{U^{(j)}}-a_U\|_2\|c_{U^{(j)}}\|_2 +\|a_U\|_2\|c_{U^{(j)}}-c_U\|_2\bigr) \longrightarrow0.
\]
Only this approximation for fixed reference matrices uses bounds depending on $n,X,R$.

The quadratic forms for the subdivisions converge to the integral form as well. Equation \eqref{S:eq:continuous-refinement-energy} gives
\[
 \begin{aligned}
 0&\le\mathfrak E_{\pi_j,X}(U^{(j)})-\mathfrak e_{n,X}(U^{(j)})\\
  &\le\frac{|\pi_j|}{2}\int_{I_*}
       \bigl[q|U^{(j)}a|^2+p|U^{(j)}c|^2\bigr]\,\dd r
   =\frac{|\pi_j|}{2}\mathcal E_n^{\mathrm{diag}}(U^{(j)})
   \le\frac{|\pi_j|}{2}\mathcal E_n^+(U).
 \end{aligned}
\]
Cauchy--Schwarz for the symmetric bilinear form associated with $\mathfrak e_{n,X}$ gives
\[
 |\mathfrak e_{n,X}(U^{(j)})-\mathfrak e_{n,X}(U)|
 \le n^{-1}\|U^{(j)}-U\|_{n,+}
                     (\|U^{(j)}\|_{n,+}+\|U\|_{n,+})\longrightarrow0.
\]
Thus $\mathfrak E_{\pi_j,X}(U^{(j)})\to\mathfrak e_{n,X}(U)$, including the triangles within the subdivided intervals.

\subsubsection{The exact Hilbert--Schmidt inner product}
\label{S:sec:s5-bilinear-form}

To compute the leading order coefficient of the squared error, we use the exact identity and retain the cross products in \eqref{S:eq:s3-linear-norm-exact}. With sampled
coefficients they separate into terms on different intervals and terms inside one interval; the latter are still two triangular integrals.

Keep the same reference columns and weights. For $F,G\in\mathfrak H_{n,X}^+$ write $F_k(s)=F((k+s)/n)$ and put
\[
 \alpha_k^F(s)=F_k(s)a_k,\qquad \gamma_k^F(s)=F_k(s)c_k,
\]
with the corresponding notation for $G$. Define
\[
 r_k^a(F)=\int_0^1\langle\alpha_k^F(s),a_k\rangle\,\dd s,
 \qquad
 r_k^c(F)=\int_0^1\langle\gamma_k^F(s),c_k\rangle\,\dd s.
\]
The contribution from the two square terms on distinct partition intervals is
\[
 \mathsf L_n(F,G)=\frac12\sum_k\int_0^1
       \bigl[A_{n,k}\langle\alpha_k^F(s),\alpha_k^G(s)\rangle
             +C_{n,k}\langle\gamma_k^F(s),\gamma_k^G(s)\rangle\bigr]\,\dd s.
\]
The mixed products on those intervals give
\[
 \mathsf O_n(F,G)=\frac1{2n}\sum_{m\le k<\ell<n}
       \bigl[r_k^c(F)r_\ell^a(G)+r_k^c(G)r_\ell^a(F)\bigr].
\]
The contribution from $I_k\times I_k$, with both triangles included, is
\[
 \mathsf S_{n,k}(F,G)=\frac1{2n}\int_{0<t<s<1} \Big\langle\alpha_k^F(s)c_k^{\mathsf T} +a_k\gamma_k^F(t)^{\mathsf T},\qquad \alpha_k^G(s)c_k^{\mathsf T} +a_k\gamma_k^G(t)^{\mathsf T}\Big\rangle_{\mathrm{HS}} \,\dd t\,\dd s.
\]
The equality of the two triangular integrals, proved below, accounts for the factor $1/(2n)$. Here all inner products are real, and
$\langle M,N\rangle_{\mathrm{HS}}=\operatorname{tr}(M^{\mathsf T}N)$. The maps $\mathsf L_n$, $\mathsf O_n$, and $\mathsf S_{n,k}$ are symmetric bilinear forms on
$\mathfrak H_{n,X}^+\times\mathfrak H_{n,X}^+$. Their quadratic forms are obtained by setting $G=F$. The forms $\mathsf L_n$ and $\mathsf S_{n,k}$ are nonnegative on the diagonal;
$\mathsf O_n$ need not be. The fields may vary within each interval.

\begin{proposition}\label{S:prop:local-trace-aux}
For all $F,G\in\mathfrak H_{n,X}^+$,
\begin{equation}\label{S:eq:bridge-local-trace-formula}
 n\langle\mathscr D_nF,\mathscr D_nG\rangle_{\mathfrak S_2}
       =\mathsf L_n(F,G)+\sum_{k=m}^{n-1}\mathsf S_{n,k}(F,G)
                              +\mathsf O_n(F,G).
\end{equation}
\end{proposition}
\begin{proof}
Begin with bounded simple fields. For $k<\ell$ and $x=(k+s)/n$, $y=(\ell+t)/n$, the derivative on $I_k\times I_\ell$ is
\[
 \dot K_X[F](x,y)=\tfrac12
       \bigl[\gamma_k^F(s)a_\ell^{\mathsf T}
                              +c_k\alpha_\ell^F(t)^{\mathsf T}\bigr].
\]
The derivative on the transposed rectangle is its transpose. The change of variables has Jacobian $n^{-2}$, and the normalized inner product multiplies by $n$. Thus the pair of
rectangles contributes
\[
 \begin{aligned}
 \frac1{2n}\int_0^1\!\int_0^1\bigl[&
 p_\ell\langle\gamma_k^F(s),\gamma_k^G(s)\rangle
       +q_k\langle\alpha_\ell^F(t),\alpha_\ell^G(t)\rangle\\
 &+\langle\gamma_k^F(s),c_k\rangle
                         \langle a_\ell,\alpha_\ell^G(t)\rangle\\
 &+\langle\gamma_k^G(s),c_k\rangle
                         \langle a_\ell,\alpha_\ell^F(t)\rangle\bigr]
                                               \,\dd s\,\dd t.
 \end{aligned}
\]
This follows by expanding the inner product and using the rank-one identity in the proof of
Proposition~\ref{S:prop:6.1}. Summing the first term in $\ell>k$ yields $\tfrac12
C_{n,k}\int_0^1\langle\gamma_k^F,\gamma_k^G\rangle$; summing the second in $k<\ell$ yields $\tfrac12
A_{n,\ell}\int_0^1\langle\alpha_\ell^F,\alpha_\ell^G\rangle$. The two cross integrals factor and give $\mathsf
O_n(F,G)$.

On the triangle $0<t<s<1$, the derivative is $[\alpha_k^F(s)c_k^{\mathsf T}+a_k\gamma_k^F(t)^{\mathsf T}]/2$. On $0<s<t<1$, it is
$[\gamma_k^F(s)a_k^{\mathsf T}+c_k\alpha_k^F(t)^{\mathsf T}]/2$. Interchanging $s,t$ in the second integral and transposing both matrices makes its inner product equal to the
first. Each carries the prefactor $1/(4n)$, so their sum is $\mathsf S_{n,k}(F,G)$ as defined above. These squares and the pairs of distinct rectangles partition $I_*^2$ up to a
null set. This proves \eqref{S:eq:bridge-local-trace-formula} for the simple fields.

Cauchy--Schwarz in the weighted Hilbert norm gives
\[
 |\mathsf L_n(F,G)|\le\tfrac12
                 \mathcal E_n(F)^{1/2}\mathcal E_n(G)^{1/2}.
\]
For a diagonal square, applying $|v+w|^2\le2|v|^2+2|w|^2$ as in \eqref{S:eq:6.23} gives
\[
 \begin{aligned}
 0\le\mathsf S_{n,k}(F,F)
 &\le\frac1n\int_0^1
       \bigl[sq_k|\alpha_k^F(s)|^2
                 +(1-s)p_k|\gamma_k^F(s)|^2\bigr]\,\dd s\\
 &\le\frac1n\int_0^1
       \bigl[q_k|\alpha_k^F(s)|^2+p_k|\gamma_k^F(s)|^2\bigr]\,\dd s.
 \end{aligned}
\]
Cauchy--Schwarz on each square, and then in $k$, therefore implies
\[
 \sum_k|\mathsf S_{n,k}(F,G)|
 \le\left(\sum_k\mathsf S_{n,k}(F,F)\right)^{1/2}
       \left(\sum_k\mathsf S_{n,k}(G,G)\right)^{1/2}
 \le\mathcal E_n^{\mathrm{diag}}(F)^{1/2}
                         \mathcal E_n^{\mathrm{diag}}(G)^{1/2}.
\]
For a cross term use
\[
 |r_k^c(F)|\le\sqrt{q_k}
                 \left(\int_0^1|\gamma_k^F(s)|^2\,\dd s\right)^{1/2},
 \qquad
 |r_\ell^a(G)|\le\sqrt{p_\ell}
                 \left(\int_0^1|\alpha_\ell^G(t)|^2\,\dd t\right)^{1/2}.
\]
Cauchy--Schwarz over the ordered pairs gives
\[
 \frac1{2n}\sum_{k<\ell}|r_k^c(F)r_\ell^a(G)| \le\frac12\left(\sum_k C_{n,k} \int_0^1|\gamma_k^F(s)|^2\,\dd s\right)^{1/2} \cdot\left(\sum_\ell A_{n,\ell} \int_0^1|\alpha_\ell^G(t)|^2\,\dd t\right)^{1/2}.
\]
Interchanging $F,G$ bounds the other cross term. Applying the two-dimensional Cauchy--Schwarz inequality to their sum gives
$|\mathsf O_n(F,G)|\le\mathcal E_n(F)^{1/2}\mathcal E_n(G)^{1/2}/2$. These estimates bound $\mathsf L_n$, $\sum_k\mathsf S_{n,k}$, and $\mathsf O_n$ by $C\|F\|_{n,+}\|G\|_{n,+}$,
with respective constants $1/2$, $1$, and $1/2$. Proposition~\ref{S:prop:6.1} gives the same bound, with constant $1$, for
$n\langle\mathscr D_nF,\mathscr D_nG\rangle_{\mathfrak S_2}$. For the interval averages $F_j,G_j$ of \eqref{S:eq:continuous-simple-approximation}, each of these four bilinear
expressions $\mathfrak b$ therefore satisfies
\[
 |\mathfrak b(F_j,G_j)-\mathfrak b(F,G)|
 \le C\bigl(\|F_j-F\|_{n,+}\|G_j\|_{n,+}
           +\|F\|_{n,+}\|G_j-G\|_{n,+}\bigr)\longrightarrow0.
\]
Passing to the limit in the identity for $F_j,G_j$ proves \eqref{S:eq:bridge-local-trace-formula} for arbitrary $F,G$.

\end{proof}

Define the bounded symmetric bilinear form on $\mathfrak H_{n,X}^+\times\mathfrak H_{n,X}^+$ by $\mathfrak T_n(F,G)=n\langle\mathscr D_nF,\mathscr D_nG\rangle_{\mathfrak S_2}$. Its
associated quadratic form is $\mathfrak T_n(F,F)=n\|\mathscr D_nF\|_{\mathfrak S_2}^2\ge0$. Proposition~\ref{S:prop:6.1} implies
\begin{equation}\label{S:eq:continuous-bilinear-bound}
 |\mathfrak T_n(F,G)|\le\|F\|_{n,+}\|G\|_{n,+}.
\end{equation}
For step fields, \eqref{S:eq:kernel-difference-matrix} gives also the exact block formula
\[
 \mathfrak T_n(F,G)=\frac1n\sum_{k,\ell=m}^{n-1}
                 \langle D_{k\ell}(F),D_{k\ell}(G)\rangle_{\mathrm{HS}}.
\]
The equality with the integral formula uses the diagonal symmetry from Lemma~\ref{S:lem:exact-mesh-kernel-stability}. For varying fields we use the integral identity
\eqref{S:eq:bridge-local-trace-formula}.

\subsubsection{Measurability and conditional covariance}
\label{S:sec:s5-measurability-covariance}
\phantomsection\label{S:prop:continuous-measurability}

Suppose the finite reference columns and weights are measurable with respect to a sigma-field $\mathcal G$. On the full space of real matrix fields the map
$\mathcal J_X=\mathcal J_n$ of \eqref{S:eq:weighted-space-isometry} is onto $\mathcal H_0$. The traceless fields have a closed image in that fixed space. Define
\[
 \mathcal T_X=\mathscr D_{n,X}\mathcal J_X^{-1}:\mathcal H_0
                                              \longrightarrow\mathfrak S_2(\mathcal H).
\]
Proposition~\ref{S:prop:6.1} applies to arbitrary real matrix fields, so $\|\mathcal T_X\|\le n^{-1/2}$.

This is a measurable family in the strong operator sense. To prove it, choose a countable dense set in $\mathcal H_0$ consisting of step functions with rational boundary points and
rational matrix entries. For one such deterministic function $v$, the kernel of $\mathcal T_Xv$ is a finite linear combination of deterministic rectangle and triangle indicators
with measurable matrix coefficients. The coefficients involve only $a_k,c_k,D_k^{-1/2}$, the last measurable by Lemma~\ref{S:lem:weighted-hilbert}. Thus $\mathcal T_Xv$ is strongly
measurable as an $L^2$ kernel and as a Hilbert--Schmidt operator. For arbitrary deterministic $v\in\mathcal H_0$, choose $v_j$ from that dense set with $v_j\to v$. The
deterministic bound
\[
 \|\mathcal T_Xv_j-\mathcal T_Xv\|_{\mathfrak S_2}
                          \le n^{-1/2}\|v_j-v\|_{\mathcal H_0}
\]
proves strong operator measurability: for each fixed $v$, the image is a strongly measurable vector.

Suppose $U(\omega,r)$ is jointly measurable and $\mathcal E_n^+(U)<\infty$ almost surely. On the measurable event $\{\mathcal E_n^+(U)<\infty\}$ put $v=\mathcal J_XU$, and set
$v=0$ on its null complement. Fix an orthonormal basis $(e_l)$ of $\mathcal H_0$. The coefficients $c_l(\omega)=\langle v(\omega),e_l\rangle$ are measurable integrals. Parseval
gives
\[
 \left\|v-\sum_{l=1}^Lc_le_l\right\|_{\mathcal H_0}^2
 =\|v\|_{\mathcal H_0}^2-\sum_{l=1}^L|c_l|^2\longrightarrow0.
\]
This constructs a strongly measurable $\mathcal H_0$-valued version. For a fixed countable dense set $(v_l^0)$, define
\[
 k_j(\omega)=\min\{l:\|v(\omega)-v_l^0\|_{\mathcal H_0}<2^{-j}\},
 \qquad v_j(\omega)=v_{k_j(\omega)}^0.
\]
Density gives an index at every outcome, and its measurability follows from
\[
 \{k_j=l\}=\{\|v-v_l^0\|<2^{-j}\}
                \cap\bigcap_{i<l}\{\|v-v_i^0\|\ge2^{-j}\}.
\]
Thus $v_j$ is countably valued and measurable. On $\{k_j=l\}$, $\mathcal T_Xv_j=\mathcal T_Xv_l^0$, whose measurability was proved above. Finally,
\[
 \|\mathcal T_Xv_j-\mathscr D_{n,X}U\|_{\mathfrak S_2}
 \le n^{-1/2}\|v_j-v\|_{\mathcal H_0}<n^{-1/2}2^{-j}
\]
on the event $\{\mathcal E_n^+(U)<\infty\}$. The limit is strongly measurable.

If in addition $\eta(U)<\infty$ almost surely, take deterministic dyadic subdivisions of the parent intervals and their averages $U^{(j)}$. Their finitely many entries are
measurable. The kernels $\mathsf K_X(U^{(j)})$ are measurable finite sums of rectangle and triangle indicators, with almost surely finite $L^2$ norms. The nonlinear approximation
proved above gives $\mathsf K_X(U^{(j)})\to\mathsf K_X(U)$ in $\mathfrak S_2$ almost surely. Thus the nonlinear kernel is strongly measurable as well. The derivative is defined on
the measurable event of finite reference bounds and finite $\mathcal E_n^+(U)$; the nonlinear map additionally requires $\eta(U)<\infty$. Whenever a quantity is defined almost
surely, extend it by zero on its exceptional null set.

The derivative also inherits a support restriction from the perturbation. If $U$ vanishes outside a measurable set $S\subset I_*$, each term of
\eqref{S:eq:continuous-derivative-kernel} vanishes when both arguments lie outside $S$. Hence the derivative kernel is supported in $(S\times I_*)\cup(I_*\times S)$. In particular,
if $U$ is supported in $I_k$, Proposition~\ref{S:prop:6.1} gives
\begin{equation}\label{S:eq:continuous-one-subinterval-bound}
 \|\mathscr D_{n,X}U\|_{\mathfrak S_2}^2
       \le\frac1n\int_0^1\mathcal M_{n,k}^+(U_k(s))\,\dd s,\qquad
 \mathcal M_{n,k}^+(V)=A_{n,k}^+|Va_k|^2+C_{n,k}^+|Vc_k|^2.
\end{equation}
The maps depend only on the Lebesgue equivalence class of the field.

Finally, let $\mathcal G$ admit a regular conditional law, let $Z$ be conditionally centered in $\mathfrak H_{n,X}^+$, and assume
\[
 M=\mathbb E[\mathcal E_n^+(Z)\mid\mathcal G]<\infty
                                  \quad\text{almost surely}.
\]
The reference columns are $\mathcal G$-measurable. Use their isometry $\mathcal J_X$ to put $v=\mathcal J_XZ$ in the fixed separable Hilbert space $\mathcal H_0$. Then
\[
 \|v\|_{\mathcal H_0}^2=\mathcal E_n^+(Z),\qquad
 G_n=\mathcal T_Xv,\qquad\|\mathcal T_X\|\le n^{-1/2}.
\]
Under a regular conditional law given $\mathcal G$, both $\mathcal J_X$ and $\mathcal T_X$ are fixed operators. Centering in the original weighted space is therefore equivalent to
centering of $v$ in this fixed space under that law.

Set $M=\mathbb E[\|v\|_{\mathcal H_0}^2\mid\mathcal G]$ and $E_L=\{M\le L\}$ for positive integers $L$. These events belong to $\mathcal G$ and increase to an event of probability
one. Although $\mathbb E M$ need not be finite,
\[
 \mathbb E[\mathbf1_{E_L}\|v\|_{\mathcal H_0}^2]
 =\mathbb E[\mathbf1_{E_L}M]\le L,
 \qquad
 \mathbb E[\mathbf1_{E_L}\|G_n\|_{\mathscr K}^2]\le L/n.
\]
Thus the localized vectors have Bochner conditional expectations. Conditional bounded linearity gives $\mathbb E[\mathbf1_{E_L}G_n\mid\mathcal G]=0$.

The localized covariance
\[
 Q_n^{(L)}=\mathbb E[\mathbf1_{E_L}(G_n\otimes G_n)\mid\mathcal G]
\]
exists in the trace-class space, because the trace norm of the integrand is $\mathbf1_{E_L}\|G_n\|^2$, which is integrable. If $L'\ge L$, taking the $\mathcal G$-measurable
indicator outside the conditional integral gives
\[
 \mathbf1_{E_L}Q_n^{(L')}=Q_n^{(L)}.
\]
For almost every outcome, $M$ is finite, so these localized covariances stabilize once $L\ge M$. Their limit defines a strongly measurable positive trace-class operator $Q_n$, and
on every $E_L$ it equals the covariance computed under the regular conditional law. The trace and bilinear identities can consequently be established first on each $E_L$ and then
on their union. This justifies those identities without assuming that $G_n$ has an unconditional second moment. Conditional Gaussian characteristic functions are bounded
integrands, so the same localization introduces no extra integrability requirement for them.

For one conditional law satisfying the hypotheses, let $\mu$ be the law of $G_n$ and use the rank-one convention $(v\otimes w)h=\langle h,w\rangle v$. Since
$\|v\otimes v\|_{\mathfrak S_1(\mathscr K)}=\|v\|_{\mathscr K}^2$, the Bochner integral
\[
 Q_n=\int_{\mathscr K}v\otimes v\,\mu(\dd v)
\]
exists in trace norm and is a positive operator. For $A,B\in\mathscr K$ it gives
\begin{equation}\label{S:eq:s5-conditional-covariance}
 \langle Q_nA,B\rangle_{\mathscr K}
       =\mathbb E[\langle G_n,A\rangle_{\mathscr K}
                      \langle G_n,B\rangle_{\mathscr K}\mid\mathcal G].
\end{equation}
For a deterministic orthonormal basis $(e_j)$ of $\mathscr K$, Parseval and monotone convergence give
\begin{equation}\label{S:eq:s5-conditional-trace}
 \operatorname{tr}Q_n =\sum_j\mathbb E[|\langle G_n,e_j\rangle|^2\mid\mathcal G] =\mathbb E[\|G_n\|_{\mathscr K}^2\mid\mathcal G],\qquad n\operatorname{tr}Q_n =\mathbb E[\mathfrak T_n(Z,Z)\mid\mathcal G]\le M.
\end{equation}
If $Z$ is conditionally Gaussian, let $A_1,\ldots,A_d\in\mathscr K$ and $t_1,\ldots,t_d\in\mathbb R$. For
$A=\sum_jt_jA_j$,
\[
 \mathbb E[\exp(\mathrm i\langle G_n,A\rangle_{\mathscr K})\mid\mathcal G]
 =\exp\left(-\tfrac12\langle Q_nA,A\rangle_{\mathscr K}\right).
\]
Indeed $\langle G_n,A\rangle=\langle Z,\mathscr D_{n,X}^*A\rangle_{n,+}$, and the conditional reference columns fix that bounded linear functional. This characteristic function
identifies the conditional Gaussian distribution, including covariance with a nontrivial kernel. The covariance and trace identities hold for every centered field with the stated
conditional second moment.

These covariances are strongly $\mathcal G$-measurable in trace norm. Let $P_N$ project onto $\operatorname{span}\{e_1,\ldots,e_N\}$. The entries of $P_NQ_nP_N$ are measurable by
\eqref{S:eq:s5-conditional-covariance}. The rank-one decomposition
\[
 v\otimes v-w\otimes w=(v-w)\otimes v+w\otimes(v-w)
\]
gives $\|v\otimes v-w\otimes w\|_1\le(\|v\|+\|w\|)\|v-w\|$. Applying it with $w=P_Nv$, then conditional Cauchy--Schwarz, yields
\[
 \|Q_n-P_NQ_nP_N\|_1
 \le2\mathbb E[\|G_n\|^2\mid\mathcal G]^{1/2}
          \mathbb E[\|(I-P_N)G_n\|^2\mid\mathcal G]^{1/2}
 \longrightarrow0.
\]
The final convergence follows from dominated convergence under the conditional law. Thus $Q_n$ is the almost sure limit in trace norm of measurable finite matrices, which proves
the claim.

Let $\widetilde Q_n^+$ be the covariance of $Z$ in the weighted Hilbert space after the reference columns have been fixed. The adjoint $\mathscr D_{n,X}^*$ in the following
identity is taken with respect to $\langle\cdot,\cdot\rangle_{n,+}$. Applying the covariance definition to that adjoint gives, for $A,B\in\mathscr K$,
\[
 \langle Q_nA,B\rangle =\mathbb E[\langle Z,\mathscr D_{n,X}^*A\rangle_{n,+} \langle Z,\mathscr D_{n,X}^*B\rangle_{n,+}\mid\mathcal G] =\langle\widetilde Q_n^+\mathscr D_{n,X}^*A, \mathscr D_{n,X}^*B\rangle_{n,+}.
\]
Consequently
\begin{equation}\label{S:eq:s5-covariance-of-image}
 Q_n=\mathscr D_{n,X}\widetilde Q_n^+\mathscr D_{n,X}^*,\qquad
 \|Q_n\|_{\mathrm{op}}\le n^{-1}\|\widetilde Q_n^+\|_{\mathrm{op}}.
\end{equation}
On the fixed space $\mathcal H_0$, the same identity uses $\mathcal T_X$ and the covariance of $\mathcal J_XZ$, under the same conditional law.

The factor $n^{-1}$ in \eqref{S:eq:s5-covariance-of-image} is the square of the deterministic norm bound for $\mathscr D_{n,X}$. For example, for $\|A\|_{\mathscr K}=1$,
\[
 \langle Q_nA,A\rangle
 \le\|\widetilde Q_n^+\|_{\mathrm{op}}
                   \|\mathscr D_{n,X}^*A\|_{n,+}^2
 \le n^{-1}\|\widetilde Q_n^+\|_{\mathrm{op}}.
\]
For the bridge field of Section~\ref{S:sec:the-error-from-sampling-the-brownian-path}, Proposition~\ref{S:prop:weighted-field-covariance}, specifically
\eqref{S:eq:weighted-field-covariance-bound}, proves under the stated sampling sigma-field that
\[
 \operatorname{tr}\widetilde Q_n^+=O_{\beta,\omega}(1),\qquad
 \|\widetilde Q_n^+\|_{\mathrm{op}}
       \le C_{\beta,\omega,\varepsilon}n^{-a_\beta+\varepsilon},
 \quad a_\beta=\min\{1,2/\beta\}.
\]
Substitution into \eqref{S:eq:s5-conditional-trace} and \eqref{S:eq:s5-covariance-of-image} gives respectively
\[
 \operatorname{tr}Q_n=O_{\beta,\omega}(n^{-1}),\qquad
 \|Q_n\|_{\mathrm{op}}
       \le C_{\beta,\omega,\varepsilon}n^{-1-a_\beta+\varepsilon}.
\]
Proposition~\ref{S:prop:weighted-field-covariance} supplies the required bound on the field covariance. Section~\ref{S:sec:covariance-estimates} also estimates the operator
covariance directly from the bridge modes.

\clearpage

\section{Discrete approximations of the limiting operators}\label{sec:discrete-approximations-of-the-limiting-operators}

Here the limiting Brownian environment is fixed, and only the resolution of its approximation changes. At the soft edge we use a finite difference operator on an expanding
interval; at the hard edge we use logarithmic cells and the limiting inverse kernel; for the circular ensemble we sample a hyperbolic Brownian path while retaining the terminal
boundary line. This section estimates the error from sampling the Brownian environment. The comparison with coefficients having the exact finite ensemble laws is carried out in
Section~\ref{sec:coupling-matrix-coefficients-to-brownian-motion}.

\subsection{Finite difference approximation of the stochastic Airy operator}
\label{sec:finite-difference-approximation-of-the-airy-operator}
Both couplings for the tridiagonal models are compared with the same finite difference operator driven by the Brownian increments. Its approximation error must be controlled on an
expanding interval, with constants that do not deteriorate with the interval length.

\subsubsection{Pathwise estimates for the finite difference approximation}\label{A:sec:finite-difference-approximation-of-the-stochastic-airy-operator}

Fix $\beta>0$ and a nonempty compact set $K\Subset\mathbb C\setminus\mathbb R$. Enlarge $K$ to $K\cup\overline K$ so that the same estimates apply to adjoints, and set
$\delta_K=\operatorname{dist}(K,\mathbb R)>0$. The comparison is made for each $B\in\Omega_B$; the probability estimates below establish that this Borel set has full Wiener
measure, independently of the spacings, boundary points and spectral sets.

Write $L=(N+1)h$. Except in the final half-line estimate, it suffices to assume
\[
 0<h\le1,\qquad 1\le L\le C_+h^{-2},\qquad
 N=L/h-1\le C_+h^{-3},
\]
where $C_+$ is deterministic and fixed. The theorem also assumes $L\ge c_-h^{-2}$ for a fixed $c_->0$. Constants denoted by $C_{\omega,\beta,K}$ may depend on $C_+$ and on a fixed
exponent $0<\alpha<1/2$ in the bound on Brownian oscillations over unit intervals. They are independent of $h,N,L$; dependence on $c_-$ enters only in the final truncation. We
prove
\[
\sup_{z\in K}
\left\|
\iota_h(\mathcal B_{N,h,\beta}^{B}-z)^{-1}\iota_h^*
-
(\mathcal H_\beta-z)^{-1}
\right\|_{\mathfrak S_2}^2
\le
C_{\omega,\beta,K}h.
\]
We first impose a Dirichlet condition at $L$, then sample the Brownian primitive on $(0,L)$. Its distributional derivative is supported at the nodes, so the sampled resolvent is
$H^2$ on each open cell. This permits linear interpolation and comparison with the variational solution in the piecewise linear space. Exact integration of the hat functions then
identifies the difference between the matrix of the restricted quadratic form and $\mathcal B_{N,h,\beta}^B$. The last comparison replaces the linear interpolation map by
$\iota_h$.

The discretization errors and the truncation error impose opposite bounds on the interval length. The finite-dimensional part of the error from restricting the variational equation
has squared norm bounded by $CNh^4$, and the potential quadrature contributes $Ch^4(1+L)^{3/2}$. Both are $O(h)$ when $L\le C_+h^{-2}$. The half-line tail contributes $CL^{-1/2}$,
which is $O(h)$ when $L\ge c_-h^{-2}$. The proof below retains these quantities until the corresponding sums have been evaluated.

\paragraph{Brownian oscillations and partial sums}\label{A:sec:brownian-and-gaussian-primitive-bounds}
We need two controls on the same Brownian path. Oscillation on unit intervals gives common form bounds. A weighted sum of squared oscillations on intervals of length $h$ gives the
approximation rate. Use $\mathbb W$, $\mathcal O_\alpha$, $T_r$, $Q$ and $\Omega_B$ defined in Section~\ref{A:sec:intro-operator-embeddings}. We prove their probability bounds
here.

\begin{lemma}\label{A:res:A.6}
For $R\ge1$ and $\alpha>0$,
\[
 \mathbb W(\mathcal O_\alpha>R)\le C_\alpha\mathrm e^{-c_\alpha R^2},
 \qquad
 \mathbb W(Q>C+R)\le C\mathrm e^{-cR}.
\]
Thus $Q$ and $\mathcal O_\alpha$ have all positive moments and $\mathbb W(\Omega_B)=1$.
\end{lemma}

\begin{proof}
The reflection principle and $\operatorname{osc}_{[0,1]}B+|B(1)|\le3\sup_{t\le1}|B(t)|$ give
\[
 \mathbb W(\operatorname{osc}_{[k,k+1]}B+|B(k+1)-B(k)|>u)
 \le C\mathrm e^{-cu^2}.
\]
A union bound therefore yields
\[
\mathbb W(\mathcal O_\alpha>R) \le C\sum_{k\ge0}\mathrm e^{-cR^2(1+k)^{2\alpha}} \le C\mathrm
 e^{-cR^2/2} \sum_{k\ge0}\mathrm e^{-c(1+k)^{2\alpha}/2} \le C_\alpha\mathrm e^{-c_\alpha R^2}.
\]

For fixed $r$, Brownian scaling and independent increments imply that
\[
 Z_{r,k}=d_r^{-1}\operatorname{osc}_{[kd_r,(k+1)d_r]}(B)^2
\]
are independent and have the law of $\operatorname{osc}_{[0,1]}(B)^2$. Integration of the reflection bound gives, for integers $p\ge1$,
\[
 \mathbb E Z_{r,k}^p
 =p\int_0^\infty u^{p-1}\mathbb P(Z_{r,k}>u)\,\mathrm du
 \le C_1^p p!.
\]
For $Y_{r,k}=Z_{r,k}-\mathbb EZ_{r,k}$, the same bound holds with $C_1$ enlarged. Since $\mathbb EY_{r,k}=0$,
\[
 \mathbb E\mathrm e^{t aY_{r,k}}
 \le1+\sum_{p\ge2}(C_1|t|a)^p
 \le\exp(2C_1^2t^2a^2)
 \quad(C_1|t|a\le1/2).
\]
With $a_{r,k}=d_r(1+kd_r)^{-3/2}$, integral comparison gives
\[
 \sum_k a_{r,k}\le C,\qquad
 \max_k a_{r,k}\le d_r,\qquad
 \sum_k a_{r,k}^2\le Cd_r.
\]
Thus, by independence, for a finite sum and $0\le t\le c/d_r$,
\[
 \log\mathbb E\exp\left(t\sum_k a_{r,k}Y_{r,k}\right)
 \le Ct^2d_r.
\]
The sums defining $T_r$ converge almost surely and in $L^1$, because their terms are nonnegative and $\sum_k a_{r,k}\mathbb EZ_{r,k}<\infty$. Fatou's lemma extends the exponential
bound to the infinite sum. Taking $t=c_2\min\{u,1\}/d_r$, with $c_2$ sufficiently small, in exponential Markov gives
\[
 \mathbb W(T_r-\mathbb ET_r>u)
 \le\exp\{-c_3\min(u^2,u)/d_r\}.
\]
Since $\sup_r\mathbb ET_r\le C$, another union bound gives
\[
 \mathbb W(Q>C+R)
 \le\sum_{r\ge0}\mathrm e^{-cR2^r}
 \le\mathrm e^{-cR}\sum_{r\ge0}\mathrm e^{-cr}
 \le C'\mathrm e^{-cR},\qquad R\ge1.
\]
Here $2^r\ge r+1$. For every $p>0$, the two tail bounds give
\[
 \mathbb E\mathcal O_\alpha^p
 \le C_\alpha p\int_0^\infty t^{p-1}\mathrm e^{-c_\alpha t^2}\,\mathrm dt<\infty,
 \qquad
 \mathbb E(Q-C-1)_+^p
 \le C'p\int_0^\infty t^{p-1}\mathrm e^{-ct}\,\mathrm dt<\infty.
\]
Thus $Q$ and every $\mathcal O_q$ with rational $0<q<1/2$ are finite almost surely. Their countable intersection is $\Omega_B$. For any real $\alpha>0$, choose rational
$0<q<\min\{\alpha,1/2\}$ and use $\mathcal O_\alpha\le\mathcal O_q$. In particular, on $\Omega_B$,
\[
 \operatorname{osc}_{[m,m+1]}B+|B(m+1)-B(m)|
 \le\mathcal O_\alpha(B)(1+m)^\alpha,
 \qquad m\in\mathbb N_0.
\]
\end{proof}

\par\Needspace{5\baselineskip} The later coefficient comparisons require a related estimate for Gaussian or sub-Gaussian
partial sums. The next lemma assumes increment bounds for each array and gives moment and tail estimates uniform in $n$. Its hypotheses allow separate sampling at different sizes;
the example following the proof explains the resulting limitation on almost sure uniformity. Write
\[
 \|X\|_{\psi_2}=\inf\{a>0:\mathbb E\mathrm e^{X^2/a^2}\le2\}.
\]
Then $\mathbb P(|X|>u)\le2\exp(-u^2/\|X\|_{\psi_2}^2)$.

For each \(n\), let \(0<h_n\le1\), let \(N_n\in\mathbb N\cup\{\infty\}\), and let \((\xi_{j,n})_{1\le j\le N_n}\) be a centered jointly Gaussian array. Suppose that a deterministic
constant \(C_0>0\), independent of \(n\), satisfies
\[
\operatorname{Var}\left(\sum_{j=p}^{q}\xi_{j,n}\right)
\le C_0h_n(q-p+1),
\quad 1\le p\le q\le N_n.
\eqtag{A.23-var}
\]
Put \(\xi_{j,n}=0\) for \(j>N_n\) when \(N_n<\infty\), and set \(S_{0,n}=0\) and \(S_{m,n}=\sum_{j=1}^m\xi_{j,n}\). The continuous piecewise linear interpolation is
\[
S_{h_n}(x)=S_{m,n}+\theta\xi_{m+1,n},
\quad m=\lfloor x/h_n\rfloor,
\quad \theta=x/h_n-m\in[0,1).
\eqtag{A.23-interpolation}
\]
In particular, \(S_{h_n}(mh_n)=S_{m,n}\), with no additional factor of $h_n$, and the process is constant for \(x\ge N_nh_n\) in the finite case. For \(\alpha>0\), define
\[
\mathcal O_{\alpha,n}
=\sup_{\ell\in\mathbb N_0}(1+\ell)^{-\alpha}
\left(
\operatorname{osc}_{[\ell,\ell+1]}S_{h_n}
+|S_{h_n}(\ell+1)-S_{h_n}(\ell)|
\right).
\]

Alternatively, the array may be centered but not Gaussian, provided
\[
 \left\|\sum_{j=p}^q\xi_{j,n}\right\|_{\psi_2}
 \le K_0\sqrt{h_n(q-p+1)},\qquad 1\le p\le q\le N_n,
 \eqtag{A.23-subgaussian}
\]
with $K_0$ independent of $n$. The coordinates, blocks, and different matrix sizes may be dependent; their joint distribution is subject to the displayed bounds. These hypotheses
may also hold under a conditional law, with $C_0$ or $K_0$ uniform over the stated conditioning parameters.

\begin{lemma}\label{A:res:A.7}
For every $\alpha>0$ and $1\le p<\infty$,
\[
 \sup_n\mathbb E\mathcal O_{\alpha,n}^{\,p}
 \le C_{\alpha,p,C_0},
 \eqtag{A.23}
\]
and, for $R\ge1$,
\[
 \sup_n\mathbb P(\mathcal O_{\alpha,n}>R)
 \le C_{\alpha,C_0}\mathrm e^{-c_{\alpha,C_0}R^2}.
 \eqtag{A.23-tail}
\]
The same bounds hold for all real block boundary points and under the alternative or conditional hypotheses above, with the corresponding uniform constants.
\end{lemma}

\begin{proof}
Fix \(n\), write \(h=h_n\), and extend increments by zero at the finite boundary point. For a portion of one interval of length \(a\le h\), its interpolated increment has \(L^2\)
norm at most \(\sqrt{C_0}a/\sqrt h\le\sqrt{C_0a}\). Split any
\([s,t]\) into two partial intervals of lengths \(a,b\) and \(k\) complete intervals. By \eqref{A:eq:A.23-var} and the
triangle inequality,
\[
 \|S_h(t)-S_h(s)\|_2
 \le\sqrt{C_0}(\sqrt a+\sqrt{kh}+\sqrt b)
 \le\sqrt{3C_0(t-s)}.
\]
This also holds within one interval, and gives
\[
 \operatorname{Var}(S_h(t)-S_h(s))\le3C_0|t-s|,
 \quad s,t\ge0.
\eqtag{A.23-metric}
\]
On \([\ell,\ell+1]\), let \(M_k\) be the largest absolute increment over its \(2^k\) dyadic subintervals, and put
\(\sigma^2=3C_0\). Gaussian tails and a union bound give, for \(u\ge0\),
\[
 \mathbb P\left(\exists k\ge0:
 M_k>2\sigma2^{-k/2}(\sqrt{k+1}+u)\right)
 \le\sum_{k\ge0}2^{k+1}\mathrm{e}^{-2(\sqrt{k+1}+u)^2}
 \le C\mathrm{e}^{-2u^2}.
\]
For $t\in[\ell,\ell+1)$, put $t_k=\ell+2^{-k}\lfloor2^k(t-\ell)\rfloor$. Then $t_0=\ell$, $t_k\to t$, and each nonzero increment $t_{k+1}-t_k$ is one interval of the $(k+1)$st
dyadic partition. Continuity therefore gives
\[
 |S_h(t)-S_h(\ell)|
 \le\sum_{k\ge0}|S_h(t_{k+1})-S_h(t_k)|
 \le\sum_{k\ge1}M_k.
\]
For $t=\ell+1$, take $t\uparrow\ell+1$ in this inequality; continuity of $S_h$ preserves its right side. On the complement of the preceding exceptional event,
\[
 \sum_{k\ge0}M_k
 \le2\sigma\sum_{k\ge0}2^{-k/2}(\sqrt{k+1}+u)
 \le C\sigma(1+u).
\]
Since the oscillation is at most twice the supremum of $|S_h(t)-S_h(\ell)|$, these bounds give
\[
 \mathbb P\bigl(\operatorname{osc}_{[\ell,\ell+1]}S_h+
 |S_h(\ell+1)-S_h(\ell)|>C\sigma(1+u)\bigr)
 \le C\mathrm e^{-2u^2}.
\]
The constants are uniform in $n,h,N,\ell$. Therefore, for $R\ge1$,
\[
 \mathbb P(\mathcal O_{\alpha,n}>R)
 \le C\sum_{\ell\ge0}\mathrm{e}^{-cR^2(1+\ell)^{2\alpha}}
 \le C_\alpha \mathrm{e}^{-c_\alpha R^2}.
\]
For $R\ge1$ and $\ell\ge0$,
\[
 e^{-cR^2(1+\ell)^{2\alpha}}
 \le e^{-cR^2/2}e^{-c(1+\ell)^{2\alpha}/2}.
\]
The series in the second factor converges: its tail is bounded by $1+\int_0^\infty e^{-cx^{2\alpha}/2}\,\mathrm dx =1+(2/c)^{1/(2\alpha)}\Gamma(1/(2\alpha))/(2\alpha)$. For
$p\ge1$, Tonelli's theorem now gives
\[
 \mathbb E\mathcal O_{\alpha,n}^p
 =p\int_0^\infty R^{p-1}\mathbb P(\mathcal O_{\alpha,n}>R)\,\mathrm dR
 \le1+\frac{pC_\alpha}{2c_\alpha^{p/2}}\Gamma(p/2).
\]
For a real left boundary, if $\ell\in[m,m+1)$, then
\[
 \operatorname{osc}_{[\ell,\ell+1]}S_h
 \le\operatorname{osc}_{[m,m+1]}S_h+
       \operatorname{osc}_{[m+1,m+2]}S_h,
\]
and the difference of the boundary values is no larger than this oscillation. Also $(1+\ell)^{-\alpha}\le(1+m)^{-\alpha}$ and $(2+m)^\alpha\le2^\alpha(1+m)^\alpha$. Hence
\[
 \sup_{\ell\ge0}(1+\ell)^{-\alpha}
 \bigl(\operatorname{osc}_{[\ell,\ell+1]}S_h+
       |S_h(\ell+1)-S_h(\ell)|\bigr)
 \le2(1+2^\alpha)\mathcal O_{\alpha,n}.
\]
Applying the preceding tail and moment estimates to the right side proves the tail and moment bounds for real interval boundaries.

Under \eqref{A:eq:A.23-subgaussian}, homogeneity and the triangle inequality for the Luxemburg norm give $\|S_h(t)-S_h(s)\|_{\psi_2}\le\sqrt3K_0|t-s|^{1/2}$ by the same partition
of $[s,t]$. Put $\sigma=\sqrt3K_0$. The exponential tail in the definition of this norm gives explicitly
\[
 \mathbb P\bigl(M_k>2\sigma2^{-k/2}(\sqrt{k+1}+u)\bigr)
 \le2^{k+1}\mathrm e^{-4(\sqrt{k+1}+u)^2}.
\]
Its sum over $k$ is bounded by $C\mathrm e^{-4u^2}$. The dyadic telescoping and the weighted union bound already proved then give the non-Gaussian conclusion. For a conditional law
satisfying the stated increment bound, each Markov estimate above holds conditionally with the same $\sigma$. The sums over $k$ and $\ell$ are deterministic. Their conditional tail
bound is therefore also $C_\alpha\mathrm e^{-c_\alpha R^2}$, and integration gives \eqref{A:eq:A.23} under that law.
\end{proof}

For the Jacobi coefficients, $\xi_j=h^{-1}\bar\zeta_j$ and hence $S_{mh}=h^{-1}\sum_{j\le m}\bar\zeta_j$. Moreover, whenever \(c_{\alpha,C_0}A^2>1\), \eqref{A:eq:A.23-tail} implies
\[
 \sum_n\mathbb P\left(\mathcal O_{\alpha,n}>
        A\sqrt{\log(n+1)}\right)<\infty.
\]
It does not imply an almost sure bound uniform in \(n\): independent Gaussian arrays with $h_n=1/n$, $N_n=n$, and independent increments of variance $h_n$ have independent
$S_{h_n}(1)\sim N(0,1)$. For every fixed $M$, the independent events $\{|S_{h_n}(1)|>M\}$ have a common positive probability, so the second Borel--Cantelli lemma makes their
supremum infinite almost surely. The bound uniform in the sampling interval in Subsection~\ref{A:sec:finite-difference-approximation-of-the-stochastic-airy-operator} instead uses
samples of a single Brownian path.

\paragraph{The quadratic form and uniform kernel estimates}\label{A:sec:galerkin-approximation-and-the-resolvent-error}

First pass from the dyadic intervals in $Q$ to an arbitrary spacing. For $b\in\Omega_B$, set $b_h(x)=b(h\lfloor x/h\rfloor)$. The same path satisfies, simultaneously for every
\(0<h\le1\),
\[
 \int_0^\infty(1+x)^{-3/2}|b(x)-b_h(x)|^2\,\mathrm{d}x \le ChQ(b),\qquad
 \sum_{j\ge1}(1+jh)^{-3/2}|b(jh)-b((j-1)h)|^2 \le CQ(b).
\]
To prove both inequalities, choose $d_r/2<h\le d_r$. Put $C_j=[jh,(j+1)h]$ and $D_k=[kd_r,(k+1)d_r]$. Each $C_j$ meets at most two $D_k$ in positive length, and continuity gives
\[
 \operatorname{osc}_{C_j}(b)^2
 \le2\sum_{k:\,|C_j\cap D_k|>0}
                    \operatorname{osc}_{D_k}(b)^2.
\]
Each $D_k$ meets at most four $C_j$. On intersecting intervals, $|jh-kd_r|\le1$, so their weights are comparable. Therefore
\[
 \sum_{j\ge0}(1+jh)^{-3/2}\operatorname{osc}_{C_j}(b)^2
 \le C T_r(b)\le C Q(b).
\]
On $C_j$ one has $|b(x)-b_h(x)|\le\operatorname{osc}_{C_j}(b)$. Integration over $C_j$ supplies the factor $h$ in the first inequality; using the two boundary values gives the
second. Restriction to any finite interval only decreases the left sides. The comparison of the two partitions is deterministic once $b$ is fixed. Thus the single condition
$Q(b)<\infty$ gives both estimates for every real $h\in(0,1]$ and every boundary point; no intersection over uncountably many probability one events is involved. Multiplication by
$c_\beta$ gives the corresponding bounds for $e=c_\beta B$.

Fix $B\in\Omega_B$. We prove the closedness and lower bound of the form defined in Section~\ref{A:sec:intro-limiting-operator}. The lower bound gives a shift independent of the
finite boundary point and the sampling interval, as required by the continuum comparisons. Equip $\mathcal D_{\mathrm A}$ with the norm
\[
 \|f\|_{\mathcal E}^2=\|f'\|_2^2+
        \int_0^\infty(1+x)|f(x)|^2\,\mathrm dx.
 \eqtag{airy-energy}
\]
Let $\mathcal L=-\partial_x^2+x+1$ with Dirichlet boundary condition; its quadratic form is $\|f\|_{\mathcal E}^2$. Use $\mathfrak a$, $\mathfrak w_B$ and $\mathfrak q_\beta$ from
\eqref{A:eq:airy-form}. For $\overline B_p=\int_p^{p+1}B(x)\,\mathrm dx$,
\[
|\overline B_{p+1}-\overline B_p| \le\int_0^1|B(p+1+t)-B(p+t)|\,\mathrm dt
 \le\operatorname{osc}_{[p,p+1]}B+ \operatorname{osc}_{[p+1,p+2]}B.
\]
Consequently,
\[
 \operatorname{osc}_{[p,p+1]}B+
 |\overline B_{p+1}-\overline B_p|
 \le C_\alpha\mathcal O_\alpha(B)(1+p)^\alpha.
 \eqtag{4.4}
\]
For $0<\alpha<1/2$, Lemma~\ref{A:res:3.1} gives, on compactly supported functions, the explicit estimate
\[
 |\mathfrak w_B[f]|
 \le\varepsilon\mathfrak a[f]
 +C_\alpha\varepsilon^{-(1+2\alpha)/(1-2\alpha)}
 (1+\mathcal O_\alpha(B))^{2/(1-2\alpha)}\|f\|_2^2,
 \qquad 0<\varepsilon\le1.
 \eqtag{form-relative}
\]
Polarization gives, for every $t>0$,
\[
 |\mathfrak w_B[f,g]|=|\mathfrak w_B[tf,t^{-1}g]|
 \le C_B\bigl(t^2\|f\|_{\mathcal E}^2+t^{-2}\|g\|_{\mathcal E}^2\bigr).
\]
For nonzero $f,g$, take $t^2=\|g\|_{\mathcal E}/\|f\|_{\mathcal E}$. After changing $C_B$, this yields
\[
 |\mathfrak w_B[f,g]|\le C_B\|f\|_{\mathcal E}\|g\|_{\mathcal E}.
\]
If $f=0$ or $g=0$, sesquilinearity gives $\mathfrak w_B[f,g]=0$.

Compactly supported functions are dense in $\mathcal D_{\mathrm A}$ for this norm. Indeed, choose cutoffs $\chi_m$ equal to one on $[0,m]$, zero on $[2m,\infty)$, and with
$|\chi_m'|\le C/m$. Then
\[
 \|(1-\chi_m)f\|_{\mathcal E}^2
 \le C\int_m^\infty\bigl(|f'|^2+(1+x)|f|^2\bigr)\,\mathrm dx
       +Cm^{-2}\int_m^{2m}|f|^2\,\mathrm dx\longrightarrow0.
\]
The product bound therefore extends $\mathfrak w_B$ uniquely to $\mathcal D_{\mathrm A}$. Choose
$\varepsilon=\min\{1,(2c_\beta)^{-1}\}$ in \eqref{A:eq:form-relative}, let $C_{B,\varepsilon}$ denote its last
coefficient, and set $\kappa_\omega=1+c_\beta C_{B,\varepsilon}$. Since $\mathfrak q_\beta=\mathfrak a+c_\beta\mathfrak
w_B$,
\[
 \tfrac12\mathfrak a[f]+\|f\|_2^2
 \le\mathfrak q_\beta[f]+\kappa_\omega\|f\|_2^2
 \le\tfrac32\mathfrak a[f]+(1+2c_\beta C_{B,\varepsilon})\|f\|_2^2.
\]
These inequalities make the shifted form norm equivalent to the complete norm \eqref{A:eq:airy-energy}; hence the form is closed and bounded below. To prove compactness of the
domain embedding, fix $M>0$ and partition $[0,M]$ into $q$ intervals of length $\delta=M/q$. Let $P_{M,q}f$ equal the mean of $f$ on each interval and vanish on $[M,\infty)$. For
such an interval $I$, the fundamental theorem of calculus gives
\[
 \int_I|f-f_I|^2\,\mathrm dx
 =\frac1{2\delta}\int_I\int_I|f(x)-f(y)|^2\,\mathrm dx\,\mathrm dy
 \le\frac{\delta^2}{2}\int_I|f'|^2\,\mathrm dx.
\]
Also
\[
 \int_M^\infty|f|^2\,\mathrm dx\le M^{-1}\|f\|_{\mathcal E}^2.
\]
Hence
\[
  \|f-P_{M,q}f\|_2^2
 \le (\delta^2/2+M^{-1})\|f\|_{\mathcal E}^2.
\]
Since $P_{M,q}$ has finite rank and $\delta^2/2+M^{-1}\to0$ when $M\to\infty$ and $M/q\to0$, the embedding $\mathcal D_{\mathrm A}\hookrightarrow L^2$ is compact. The
representation theorem for closed forms defines $\mathcal H_\beta$. Its shifted inverse maps $L^2$ boundedly into $\mathcal D_{\mathrm A}$ and is therefore compact; the resolvent
identity gives compactness of every $R_\beta(z)$. This is the construction of \cite{ramirez2011beta}.

Retain the lattice and Dirichlet conventions of Subsection~\ref{A:sec:weighted-resolvent-estimates-and-comparison-principles}.

The operator $\mathcal B_{N,h,\beta}^B$ and its diagonal sequence are defined in \eqref{A:eq:brownian-difference}. Its quadratic form is
\[
\mathfrak b_{N,h}[u]
=
\sum_{k=0}^{N}|D_hu_k|^2
+
\sum_{k=1}^{N}x_k|u_k|^2
+
\sum_{k=1}^{N}\omega_{k,h}|u_k|^2.
\eqtag{4.11}
\]
Put $e=c_\beta B$, $e_k=e(kh)$, and $e_h(x)=e(h\lfloor x/h\rfloor)$. For $u_0=u_{N+1}=0$, the product identity
$2\operatorname{Re}((u_{k+1}-u_k)\overline{(u_{k+1}+u_k)/2}) =|u_{k+1}|^2-|u_k|^2$ gives
\[
\begin{aligned}
 \sum_{k=1}^N\omega_{k,h}|u_k|^2
 &=-\sum_{k=0}^N\frac{e_k}{h}
                   (|u_{k+1}|^2-|u_k|^2)\\
 &=-2\operatorname{Re}\sum_{k=0}^N e_k(D_hu)_k\overline{(J_hu)_k}
 =-\int_0^{(N+1)h}e_h(x)(|I_hu(x)|^2)'\,\mathrm dx.
\end{aligned}
\]
The last equality follows by integrating the derivative separately on each interval; its boundary values are $h^{-1}|u_k|^2$. This is the form of $\Gamma_e$ from
Subsection~\ref{A:sec:weighted-resolvent-estimates-and-comparison-principles}.

For $p\in\mathbb N_0$, all samples of $B_h$ on $[p,p+2]$ lie in $[\max\{0,p-1\},p+2]$. The oscillation on this enlarged interval is at most the sum of the oscillations on its at
most three unit intervals. The difference of the two adjacent block averages of $B_h$ is bounded by the same oscillation. Thus $e_h$ satisfies Lemma~\ref{A:res:3.1} with a constant
at most $C_{\alpha,\beta}\mathcal O_\alpha(B)$, uniformly for $0<h\le1$. Hence
\[
 \left|\sum_{k=1}^N\omega_{k,h}|u_k|^2\right|
 \le\varepsilon\ell_h[u]+C_{\varepsilon,\omega,\beta}\|u\|^2.
\]
Proposition~\ref{A:res:3.2} gives the common form bounds uniformly in $N,h$. In \eqref{A:eq:3.8}, take $c_{j,n}=1$, $q_{j,n}=x_j+\omega_{j,h}$, and no first-order term; the
coefficient condition \eqref{A:eq:3.10} is then satisfied.

Set
\[
 L=L_h:=(N+1)h,\qquad Nh\asymp h^{-2}.
 \eqtag{4.19}
\]
The next kernel estimates hold on every interval $(0,L)$, with constants independent of $L$. The relation \eqref{A:eq:4.19} is used only when these estimates are applied to the
approximation. Let $\mathcal H_{\beta,L}^{D}$ be the operator obtained by restricting $\mathfrak q_\beta$ to $H_0^1(0,L)$, with Dirichlet conditions at both boundary points, and
let $R_{\beta,L}^{D}(z)=(\mathcal H_{\beta,L}^{D}-z)^{-1}$ act on $L^2(0,L)$. Write $P_L:L^2(\mathbb R_+)\to L^2(0,L)$ for restriction and $E_L=P_L^*$ for extension by zero. The
operator used in the half-line comparison is $E_LR_{\beta,L}^{D}(z)P_L$. More generally, an interval operator $T$ will be extended as $E_LTP_L$. This operation preserves its
nonzero singular values, since $E_L$ is an isometry and $P_LE_L=I$.

For the continuum estimates, let \(I_L=(0,L)\), where \(1\le L\le\infty\). At a finite right boundary we impose a Dirichlet condition. In both cases we use the form domain
\[
\mathcal D_L
=
\left\{f\in H_0^1(I_L):\int_{I_L}x|f(x)|^2\,\mathrm{d}x<\infty\right\}.
\]
For finite \(L\) the weighted condition is automatic. Set
\[
\|f\|_{\mathcal E_L}^2
:=
\int_{I_L}|f'(x)|^2\,\mathrm{d}x
+
\int_{I_L}(1+x)|f(x)|^2\,\mathrm{d}x.
\]
Let \(E:\mathbb R_+\to\mathbb R\) be locally bounded and satisfy the hypothesis on oscillation over blocks of Lemma~\ref{A:res:3.1}, with fixed exponent \(\alpha<1/2\) and block
constant at most \(A_E\); in the form below we use its restriction to \(I_L\). Let \(v\in L^\infty(I_L;\mathbb R)\), with \(\|v\|_\infty\le V_0\), and define
\[
\mathfrak h_{E,v,L}[f,g]
=
\int_{I_L}f'\overline{g'}\,\mathrm{d}x
+
\int_{I_L}(x+v)f\overline g\,\mathrm{d}x
-
\int_{I_L}E(x)\frac {\mathrm{d}}{\mathrm{d}x}(f\overline g)\,\mathrm{d}x.
\]
When $L=\infty$, the last integral is initially defined on compactly supported functions and extended continuously to $\mathcal D_L$. Indeed, Lemma~\ref{A:res:3.1} gives
infinitesimal form boundedness, and polarization gives the product bound in $\mathcal E_L$ proved above. For finite $L$, extension of both test functions by zero gives exactly the
same block estimate; their zero traces create no boundary term. On the half line the displayed expression is defined by the sum of its interval contributions in the preceding form
estimate. The shifted form norm is uniformly equivalent to $\|\cdot\|_{\mathcal E_L}$, so the form is closed and lower semibounded. Denote its self-adjoint operator by
\(H_{E,v,L}\) and its resolvent by \(R_{E,v,L}(z)=(H_{E,v,L}-z)^{-1}\).

The continuum argument uses the same distance $|\Phi(x)-\Phi(y)|$ as the discrete one. Here the integrated potential, as well as ordinary multiplication, is invariant under
reciprocal exponential weights. For bounded, compactly supported $g$, define the sesquilinear form on $\mathcal D_L\times\mathcal D_L$ by
\[
 \Gamma_g[f,w]
 =-\int_{I_L}g(x)(f(x)\overline{w(x)})'\,\mathrm dx.
\]
It represents the distributional derivative of $g$ and is Hermitian when $g$ is real. The resolvent product is the bounded operator specified by
\[
 \langle R_1(z)\Gamma_gR_2(z)f,w\rangle
 =\Gamma_g[R_2(z)f,R_1(\overline z)w],\qquad f,w\in L^2(I_L).
\]
Estimate~\eqref{A:eq:4.21} extends this product to every $g$ in the stated weighted $L^2$ space.

\begin{proposition}\label{A:res:4.1}
Uniformly in the class above and $z\in K$,

\[
 \|\partial_xR_{E,v,L}(z)\|_{\mathrm{op}}
 +\|(1+x)^{1/2}R_{E,v,L}(z)\|_{\mathrm{op}}\le C.
 \eqtag{4.20}
\]
The resolvent has a continuous kernel, and for some $a>0$,
\[
 |R_{E,v,L}(z;x,y)|\le
 C(1+x)^{-1/4}(1+y)^{-1/4}\mathrm e^{-a|\Phi(x)-\Phi(y)|},
\]
\[
 \int_{I_L}\mathrm e^{2a|\Phi(y)-\Phi(x)|}
 |R_{E,v,L}(z;y,x)|^2\,\mathrm dy\le C(1+x)^{-3/2}.
\]
The constants depend only on $A_E,\alpha,V_0,K$ and are uniform in $L$.
\end{proposition}

\begin{proof}
Use the energy norm $\|\cdot\|_{\mathcal E_L}$ defined above. Lemma~\ref{A:res:3.1}, applied on \(I_L\), and the bound on \(v\) imply that for every \(\varepsilon>0\),
\[
\left|
-\int_{I_L}E(|f|^2)'\,\mathrm{d}x
+\int_{I_L}v|f|^2\,\mathrm{d}x
\right|
\le
\varepsilon\|f\|_{\mathcal E_L}^2
+C_{\varepsilon,A_E,\alpha,V_0}\|f\|_2^2.
\]
When \(L<\infty\), extend \(f\in H_0^1(0,L)\) by zero to \(\mathbb R_+\). The extension belongs to \(H_0^1(\mathbb R_+)\), and the term involving the integrated potential is
unchanged. Lemma~\ref{A:res:3.1} therefore gives the same constant independently of \(L\). Choose \(\kappa\), independent of \(L\), so large that the positive form
\(\mathfrak a[f,g]=\mathfrak h_{E,v,L}[f,g]+\kappa\langle f,g\rangle\) obeys
\[
c\|f\|_{\mathcal E_L}^2
\le
\mathfrak a[f]
\le
C\|f\|_{\mathcal E_L}^2.
\]
Let \(A=H_{E,v,L}+\kappa\). In the form scale,
\[
H_{E,v,L}-z
=
A^{1/2}\bigl[I-(z+\kappa)A^{-1}\bigr]A^{1/2}.
\]
For \(\lambda\in\sigma(H_{E,v,L})\),
\[
\frac{\lambda+\kappa}{|\lambda-z|}
\le
1+\frac{|z+\kappa|}{\operatorname{dist}(K,\mathbb R)}.
\]
Let $C_A$ be a common norm equivalence constant for $\|A^{1/2}f\|_2$ and $\|f\|_{\mathcal E_L}$. By duality the same constant controls $A^{-1/2}:\mathcal E_L^*\to L^2$.
Consequently,
\[
 \|R_{E,v,L}(z)\|_{\mathcal E_L^*\to\mathcal E_L}
 \le C_A^2\sup_{\lambda\in\sigma(H_{E,v,L})}
                   \frac{\lambda+\kappa}{|\lambda-z|}
 \le C_A^2\bigl(1+\delta_K^{-1}\sup_{w\in K}|w+\kappa|\bigr).
\]
Every quantity on the right is independent of $L$. This gives the required inverse bound from $\mathcal E_L^*$ to $\mathcal E_L$. Since \(L^2(I_L)\) embeds continuously into
\(\mathcal E_L^*\),
\[
\|R_{E,v,L}(z)F\|_{\mathcal E_L}
\le C\|F\|_2.
\]
Estimate \eqref{A:eq:4.20} follows.

For the exponential estimate, fix \(y\in I_L\), let \(\psi_{y,R}(x)=\min\{|\Phi(x)-\Phi(y)|,R\}\), and put
\(\theta=\mathrm{e}^{a\psi_{y,R}}\). The function $\psi_{y,R}$ is Lipschitz and $|\psi_{y,R}'(x)|\le(1+x)^{1/2}$ almost
everywhere. For fixed $R$, $\theta^{\pm1}$ preserve the form domain, since they are bounded and $|(\theta^{\pm1})'|\le a(1+x)^{1/2}\theta^{\pm1}$. For compactly supported
\(u,w\in\mathcal D_L\), the multiplication term and the term involving the integrated potential are invariant under exponential conjugation; the identity then extends to the full
form domain by the form bound above:
\[
(\theta^{-1}u)\overline{\theta w}=u\overline w,
\quad
-\int E\frac {\mathrm{d}}{\mathrm{d}x}\bigl((\theta^{-1}u)\overline{\theta w}\bigr)
=-\int E(u\overline w)'.
\]
Only the derivative term changes. Direct differentiation gives
\[
\mathfrak h_{E,v,L}[\theta^{-1}u,\theta w]
-\mathfrak h_{E,v,L}[u,w]
=
 a\int\psi' u'\overline w
-a\int\psi'u\overline{w'}
-a^2\int(\psi')^2u\overline w.
\]
\[
\left|
\mathfrak h_{E,v,L}[\theta^{-1}u,\theta w]
-\mathfrak h_{E,v,L}[u,w]
\right|
\le
C(a+a^2)\|u\|_{\mathcal E_L}\|w\|_{\mathcal E_L},
\]
uniformly in \(R,L,y\). Let \(T_z:\mathcal E_L\to\mathcal E_L^*\) be the form operator associated with \(H_{E,v,L}-z\), and let \(C_{a,R}\) be the preceding form difference. The
preceding resolvent bound gives \(\|T_z^{-1}\|_{\mathcal E_L^*\to\mathcal E_L}\le C_K\). If $C_1$ is the constant in the conjugation bound, choose $0<a\le1$ with $2C_KC_1a\le1/2$.
Then $\|T_z^{-1}C_{a,R}\|_{\mathcal E_L\to\mathcal E_L}\le1/2$ for every $L,y,R$ and $z\in K$. The Neumann series then yields
\[
\|\theta R_{E,v,L}(z)\theta^{-1}\|_{\mathcal E_L^*\to\mathcal E_L}
\le C,
\]
uniformly in the truncation threshold \(R\).

For point evaluation, extend $f$ by zero outside $I_L$ and put $\rho_x=(1+x)^{-1/2}$, $J_x=[x-\rho_x,x+\rho_x]$. As in Lemma~\ref{A:res:3.6}, averaging
$|f(x)|^2\le2|f(t)|^2+2|x-t|\int_{J_x}|f'|^2$ over $t\in J_x$ gives
\[
|f(x)|^2
\le
C\left(\rho_x^{-1}\int_{J_x}|f|^2
+\rho_x\int_{J_x}|f'|^2\right)
\le
C(1+x)^{-1/2}\|f\|_{\mathcal E_L}^2.
\]
Define $\delta_x\in\mathcal E_L^*$ by $\langle\delta_x,v\rangle=\overline{v(x)}$. Then $\|\delta_x\|_{\mathcal E_L^*}\le C(1+x)^{-1/4}$. Also
$\|\delta_x-\delta_y\|_{\mathcal E_L^*}\le|x-y|^{1/2}$, by Cauchy--Schwarz applied to $f(x)-f(y)$. Thus

$R(z;x,y)=(R(z)\delta_y)(x)$ is jointly continuous. For compactly supported $F\in L^2(I_L)$, the identity $F=\int F(y)\delta_y\,\mathrm dy$ in $\mathcal E_L^*$ gives
\[
 (R(z)F)(x)=\int_{I_L}R(z;x,y)F(y)\,\mathrm dy.
\]
The kernel bound below implies $R(z;x,\cdot)\in L^2(I_L)$. Since evaluation composed with $R(z)$ is also continuous on $L^2(I_L)$, density extends this identity to every
$F\in L^2(I_L)$. Taking \(R>|\Phi(x)-\Phi(y)|\) and applying the conjugated $\mathcal E_L^*\to\mathcal E_L$ bound to \(\delta_y\) gives
\[
\mathrm{e}^{a|\Phi(x)-\Phi(y)|}|R_{E,v,L}(z;x,y)| = \left|\bigl(\theta
 R_{E,v,L}(z)\theta^{-1}\delta_y\bigr)(x)\right| \le C\|\delta_x\|_{\mathcal
 E_L^*}\|\delta_y\|_{\mathcal E_L^*},
\]
which proves the pointwise kernel estimate.

Choose \(0<a_1<a\). Squaring the pointwise estimate and integrating in the first variable gives
\[
\int_{I_L}\mathrm{e}^{2a_1|\Phi(t)-\Phi(x)|} |R_{E,v,L}(z;t,x)|^2\,\mathrm{d}t \le C(1+x)^{-1/2}
 \int_0^L(1+t)^{-1/2} \mathrm{e}^{-2(a-a_1)|\Phi(t)-\Phi(x)|}\,\mathrm{d}t.
\]
With \(s=\Phi(t)\), one has \(\mathrm{d}t=(1+t)^{-1/2}\,\mathrm{d}s\) and \(1+t=(1+3s/2)^{2/3}\). Hence the last integral is bounded by
\[
C\int_0^\infty(1+s)^{-2/3}\mathrm{e}^{-c|s-\Phi(x)|}\,\mathrm{d}s
\le C(1+\Phi(x))^{-2/3}
\le C(1+x)^{-1}.
\]
For the middle inequality use $(1+s)^{-2/3}\le(1+\Phi(x))^{-2/3}(1+|s-\Phi(x)|)^{2/3}$ and integrate the polynomial factor against $\mathrm e^{-c|s-\Phi(x)|}$, as in the proof of
Lemma~\ref{A:res:3.10}. Renaming \(a_1\) as \(a\) gives the weighted $L^2$ estimate for $R_{E,v,L}(z;\cdot,x)$. This proves the proposition.
\end{proof}

The kernel estimate supplies the continuum analogue of Propositions~\ref{A:res:3.15} and~\ref{A:res:3.17}. We state the consequences here so that the sampling error and the
interval truncation use exactly the same constants. Let $H_1,H_2$ belong to the preceding class on the same interval and write $R_j(z)=(H_j-z)^{-1}$. Applying the kernel estimate
at $z$ and $\bar z$ gives
\[
[R_{E,v,L}(\bar z)R_{E,v,L}(z)](x,x)
+
[R_{E,v,L}(z)R_{E,v,L}(\bar z)](x,x)
\le C(1+x)^{-3/2}.
\]
It follows that, for every measurable \(g\) in the stated weighted space,
\[
\|M_gR_{E,v,L}(z)\|_{\mathfrak S_2}^2
+
\|R_{E,v,L}(z)M_g\|_{\mathfrak S_2}^2
\le
C\int_{I_L}(1+x)^{-3/2}|g(x)|^2\,\mathrm{d}x.
\]

Initially take \(g\) bounded and compactly supported. In the sense of forms,
\(\Gamma_g=-\partial_x^*M_g-M_g\partial_x\), and hence
\[
R_1(z)\Gamma_gR_2(z)
=-R_1(z)\partial_x^*M_gR_2(z)
-R_1(z)M_g\partial_xR_2(z).
\]
By the ideal property of \(\mathfrak S_2\), \eqref{A:eq:4.20}, and the preceding multiplication estimate,
\[
\|R_1\partial_x^*M_gR_2\|_{\mathfrak S_2} \le \|\partial_xR_1(\bar z)\|_{\mathrm{op}}
 \|M_gR_2(z)\|_{\mathfrak S_2},\qquad \|R_1M_g\partial_xR_2\|_{\mathfrak S_2} \le
 \|R_1(z)M_g\|_{\mathfrak S_2} \|\partial_xR_2(z)\|_{\mathrm{op}}.
\]
Squaring the two bounds and adding gives
\[
 \|R_1(z)\Gamma_gR_2(z)\|_{\mathfrak S_2}^2
 \le C\int_{I_L}(1+x)^{-3/2}|g(x)|^2\,\mathrm dx.
 \eqtag{4.21}
\]
Put $w(x)=(1+x)^{-3/2}$ and let $g\in L^2(I_L,w(x)\,\mathrm dx)$. Define
\[
 g_m=g\mathbf1_{\{x\le m,\ |g(x)|\le m\}},\qquad m\ge1.
\]
Then $\|g_m-g\|_{L^2(w)}\to0$ by dominated convergence. The multiplication estimate applied to differences gives
\[
 \|M_{g_m-g_l}R_2(z)\|_{\mathfrak S_2}
 +\|R_1(z)M_{g_m-g_l}\|_{\mathfrak S_2}
 \le C\|g_m-g_l\|_{L^2(w)}.
\]
For the second product, interpret $R_1(z)M_g$ as the adjoint of $M_{\bar g}R_1(\bar z)$, initially on compactly supported functions. Consequently both products have limits in
$\mathfrak S_2$. The first limit agrees with multiplication by $g$ after the resolvent: along a subsequence $g_m\to g$ almost everywhere, while convergence in $\mathfrak S_2$
implies convergence in operator norm. For each fixed forcing vector these two limits agree almost everywhere, and hence in $L^2$.

Let
\[
 T_m=-[\partial_xR_1(\bar z)]^*M_{g_m}R_2(z)
       -R_1(z)M_{g_m}\partial_xR_2(z).
\]
The derivative factors are bounded uniformly in $L$ and $z\in K$. Thus
\[
 \sup_{z\in K}\|T_m(z)-T_l(z)\|_{\mathfrak S_2}
 \le C\|g_m-g_l\|_{L^2(w)}.
\]
Completeness gives a limit $T_g(z)$, uniformly on $K$. Applying the same estimate to any two approximating sequences shows that this limit does not depend on the truncation. Taking
$l\to\infty$ proves
\[
 \sup_{z\in K}\|T_g(z)\|_{\mathfrak S_2}^2
 \le C\int_{I_L}w(x)|g(x)|^2\,\mathrm dx.
\]
For real $g$ this is the operator denoted by $R_1\Gamma_gR_2$. Indeed, if $u=R_2(z)f$ and $v=R_1(\bar z)k$, the products $gu,gv$ belong to $L^2$, and $u',v'\in L^2$. Each term of
the form pairing $-\int g(u\bar v)'$ is therefore absolutely integrable by Cauchy--Schwarz. Passing to the limit from $g_m$ identifies that pairing with $\langle T_gf,k\rangle$.
This construction defines the resolvent product; it does not require $\Gamma_g$ itself to be bounded on $L^2$.

Finally,
\[
\|\mathbf1_{\{x\ge A_0\}}R_{E,v,L}(z)\|_{\mathfrak S_2}^2
=
\int_{I_L\cap[A_0,\infty)}
[R_{E,v,L}(z)R_{E,v,L}(\bar z)](x,x)\,\mathrm{d}x
\le
C\int_{A_0}^\infty(1+x)^{-3/2}\,\mathrm{d}x.
\]
The integral equals $2(1+A_0)^{-1/2}$. Applying the same bound at $\bar z$ and taking adjoints gives
\[
 \|\mathbf1_{\{x\ge A_0\}}R_{E,v,L}(z)\|_{\mathfrak S_2}^2
 +\|R_{E,v,L}(z)\mathbf1_{\{x\ge A_0\}}\|_{\mathfrak S_2}^2
 \le C A_0^{-1/2},\qquad A_0\ge1.
 \eqtag{4.22}
\]
Integration over all of $I_L$ also gives $\|R_{E,v,L}(z)\|_{\mathfrak S_2}^2\le C$, uniformly in $L$.

\paragraph{Truncation and sampling of the potential}
We first remove the part of the half-line beyond $L$, and then replace $e=c_\beta B$ by $e_h$ on the remaining interval. Proposition~\ref{A:res:4.1} applies to both $E=c_\beta B$,
$v=0$ and $E=e_h$, $v=-h/2$, with common constants. The bounded correction $-h/2$ will make the nodal potential agree with $x_j=(j-\tfrac12)h$.

Use the restriction and extension maps $P_L,E_L$ defined above. To account for the new Dirichlet condition at $L$, use a smooth cutoff in the weak equation. Its kinetic commutator
is supported in a strip of fixed width near $L$, where the column bound is $O(L^{-3/2})$. The omitted half-line has the larger squared mass $O(L^{-1/2})$.

\begin{lemma}\label{A:res:4.2}
Uniformly for \(z\in K\),
\[
\left\|
R_\beta(z)
-
E_LR_{\beta,L}^{D}(z)P_L
\right\|_{\mathfrak S_2}^2
\le
C_{\omega,\beta,K}L^{-1/2}.
\eqtag{4.23}
\]
In particular, under \eqref{A:eq:4.19}, this error is \(O_\omega(h)\).
\end{lemma}

\begin{proof}
For $1\le L\le4$, the uniform Hilbert--Schmidt bounds just proved and the triangle inequality give the result, since $L^{-1/2}\ge1/2$. Suppose $L>4$. Choose
$\chi_L\in C_c^\infty([0,L))$ equal to one on $[0,L-2]$ and zero on $[L-1,\infty)$, with $|\chi_L'|+|\chi_L''|\le C$.

For $u=R_\beta(z)f$, multiplication by $\chi_L$ preserves the domain of the quadratic form with Dirichlet boundary conditions on $(0,L)$. If $v\in H_0^1(0,L)$, extend $\chi_Lv$ by
zero. Subtraction of the two form equations gives
\[
 \mathfrak q_{\beta,L}[\chi_Lu,v]
      -\mathfrak q_\beta[u,E_L(\chi_Lv)]
 =\int_0^L\chi_L'(u\overline{v'}-u'\overline v)\,\mathrm dx
 =-\int_0^L(2\chi_L'u'+\chi_L''u)\overline v\,\mathrm dx.
\]
The boundary terms vanish because $\chi_L'$ is supported inside $(0,L)$. The potential terms cancel before any derivative of the Brownian path is taken. Thus the weak equation
gives
\[
  (\mathcal H_{\beta,L}^D-z)(\chi_Lu)
 =P_L(\chi_Lf-2\chi_L'u'-\chi_L''u).
\]
With $S_L=E_LR_{\beta,L}^D(z)P_L$, this identity becomes
\[
 R_\beta(z)-S_L
 =(1-\chi_L)R_\beta(z)-S_L(1-\chi_L)
       -2S_LM_{\chi_L'}\partial_xR_\beta(z)
       -S_LM_{\chi_L''}R_\beta(z).
\]
By \eqref{A:eq:4.22},
\[
 \|(1-\chi_L)R_\beta(z)\|_{\mathfrak S_2}^2
 \le C\int_{L-2}^\infty(1+x)^{-3/2}\,\mathrm dx
 \le C L^{-1/2}.
\]
On the interval $(0,L)$, each of $1-\chi_L$, $\chi_L'$, $\chi_L''$ is supported in $[L-2,L]$. The multiplication estimate following Proposition~\ref{A:res:4.1} therefore gives
\[
 \|S_L(1-\chi_L)\|_{\mathfrak S_2}^2
 +\|S_LM_{\chi_L'}\|_{\mathfrak S_2}^2
 +\|S_LM_{\chi_L''}\|_{\mathfrak S_2}^2
 \le C\int_{L-2}^L(1+x)^{-3/2}\,\mathrm dx
 \le C L^{-3/2}.
\]
Since $\|\partial_xR_\beta(z)\|_{\mathrm{op}}+ \|R_\beta(z)\|_{\mathrm{op}}\le C$, the ideal property bounds each commutator term by $C L^{-3/2}$ in squared Hilbert--Schmidt norm.
Finally, $\|\sum_{j=1}^4 A_j\|_{\mathfrak S_2}^2\le 4\sum_{j=1}^4\|A_j\|_{\mathfrak S_2}^2$ proves \eqref{A:eq:4.23}.
\end{proof}

On $(0,L)$, replace the integrated potential by $e_h$ and put $p_h(x)=x-h/2$, so that $p_h(jh)=x_j$. The resulting form is
\[
\mathfrak q_{\beta,h,L}^{\mathrm{fr}}[f,v] ={} \int_0^L f'\overline{v'}\,\mathrm{d}x + \int_0^L p_h(x)f\overline v\,\mathrm{d}x - \int_0^L e_h(x) \frac {\mathrm{d}}{\mathrm{d}x}\bigl(f\overline v\bigr)\,\mathrm{d}x, \quad f,v\in H_0^1(0,L).
\]
Let \(\mathcal H_{\beta,h,L}^{\mathrm{fr}}\) be the associated self-adjoint operator and set
\(R_{\beta,h,L}^{\mathrm{fr}}(z) = (\mathcal H_{\beta,h,L}^{\mathrm{fr}}-z)^{-1}\).

The bounds on the oscillations and differences of consecutive unit interval averages of $e_h$, together with $\|p_h-x\|_\infty\le1/2$, give common upper and lower form bounds,
uniformly in $h,L$, for the original and sampled potentials.

Define
\[
 \mathcal E_{\mathrm{int}}(h,L)
 =\int_0^L(1+x)^{-3/2}|e(x)-e_h(x)|^2\,\mathrm dx.
\]
The deterministic partition comparison above gives $\mathcal E_{\mathrm{int}}(h,L)\le C_\beta hQ(B)$ on $\Omega_B$. The following lemma bounds its centered fluctuations.

\begin{lemma}\label{A:res:4.3}
Uniformly in $0<h\le1$ and $L>0$, $\mathbb E\mathcal E_{\mathrm{int}}(h,L)\le C_\beta h$, and for $t>0$,
\[
 \mathbb P\bigl(\mathcal E_{\mathrm{int}}(h,L)
             -\mathbb E\mathcal E_{\mathrm{int}}(h,L)>t\bigr)
 \le2\exp\!\left\{-c(\beta)\min\left(\frac{t^2}{h^3},
                                               \frac{t}{h^2}\right)\right\}.
 \eqtag{4.28}
\]
\end{lemma}

\begin{proof}
Fix \(h,L\), and put \(I_j=(jh,(j+1)h)\) for \(j\ge0\). On \(I_j\cap(0,L)\), use the local coordinate \(s=x-jh\), put
\(l_j=\min\{h,(L-jh)_+\}\), and set \(w_j(s)=(1+jh+s)^{-3/2}\). The Gaussian element
\[
 Y_j(s)=c_\beta w_j(s)^{1/2}
                  (B(jh+s)-B(jh)),\quad 0<s<l_j,
\]
has covariance kernel \(A_j(s,t)=c_\beta^2w_j(s)^{1/2}\min(s,t)w_j(t)^{1/2}\). The elements \(Y_j\) are independent, and for their direct sum \(Y\) with covariance
\(A=\bigoplus_jA_j\),
\[
 \mathcal E_{\mathrm{int}}(h,L)=\|Y\|^2,
 \quad \mathbb E\mathcal E_{\mathrm{int}}(h,L)=\operatorname{Tr}A.
\]
The covariance on a single interval satisfies
\[
\begin{aligned}
 \operatorname{Tr}A_j
 &=c_\beta^2\int_0^{l_j}(1+jh+s)^{-3/2}s\,\mathrm ds
 \le\tfrac12c_\beta^2h^2(1+jh)^{-3/2},\\
 \sup_s\int_0^{l_j}|A_j(s,t)|\,\mathrm dt
 &\le c_\beta^2(1+jh)^{-3/2}
              \sup_{s\le h}(sh-s^2/2)
 \le\tfrac12c_\beta^2h^2(1+jh)^{-3/2},\\
 \|A_j\|_{\mathfrak S_2}^2
 &\le c_\beta^4(1+jh)^{-3}
             \int_0^h\int_0^h\min(s,t)^2\,\mathrm ds\,\mathrm dt
 =\tfrac16c_\beta^4h^4(1+jh)^{-3}.
\end{aligned}
\]
Schur's test and the direct sum decomposition therefore give
\[
 \operatorname{Tr}A \le C_\beta h^2\sum_{j\ge0}(1+jh)^{-3/2} \le C_\beta h,\quad
 \|A\|_{\mathrm{op}} \le C_\beta h^2,\quad \|A\|_{\mathfrak S_2}^2 \le C_\beta h^4
 \sum_{j\ge0}(1+jh)^{-3}\le C_\beta h^3.
\]
Write $V=\operatorname{Tr}(A^2)$ and $\Lambda=\|A\|_{\mathrm{op}}$. If $A=0$, then $Y=0$ almost surely and the conclusion is immediate. Otherwise the spectral representation of a
Gaussian element gives independent standard real normals $Z_k$ such that
\[
 \|Y\|^2-\operatorname{Tr}A
 =\sum_{k\ge1}\lambda_k(Z_k^2-1),\qquad
 \sum_k\lambda_k<\infty,\quad \sum_k\lambda_k^2=V.
\]
The nonnegative series $\sum_k\lambda_kZ_k^2$ is finite almost surely because its expectation is $\operatorname{Tr}A$. The centered series also converges in $L^2$, since the
squared norm of its tail is $2\sum_{k>d}\lambda_k^2$.

For $0<2\theta\Lambda<1$, direct Gaussian integration gives
\[
 \mathbb E\exp\{\theta\lambda(Z^2-1)\}
 =e^{-\theta\lambda}(1-2\theta\lambda)^{-1/2}.
\]
The convergent logarithmic series satisfies, for $0\le u<1$,
\[
 -\log(1-u)-u=\sum_{r\ge2}\frac{u^r}{r}
 \le\frac{u^2}{2(1-u)}.
\]
Multiplying the scalar expectations for a finite sum therefore gives
\[
 \log\mathbb E\exp\left\{\theta\sum_{k=1}^d
               \lambda_k(Z_k^2-1)\right\}
 \le\frac{\theta^2V}{1-2\theta\Lambda}.
\]
Fatou's lemma passes this inequality to the almost sure limit of the series. This is the moment generating function calculation underlying the estimate in
\cite{laurent2000adaptive}.

For $s>0$, set $\theta=\sqrt s/(\sqrt V+2\Lambda\sqrt s)$. Markov's inequality and substitution give
\[
 \begin{split}
 &\log\mathbb P\bigl(\|Y\|^2-\operatorname{Tr}A
                         >2\sqrt{Vs}+2\Lambda s\bigr)\\
 &\qquad\le-\theta(2\sqrt{Vs}+2\Lambda s)
                +\frac{\theta^2V}{1-2\theta\Lambda}=-s.
 \end{split}
\]
For a threshold $t>0$, choose
\[
 s=\frac1{16}\min\{t^2/V,t/\Lambda\}.
\]
Then $2\sqrt{Vs}+2\Lambda s\le t/2+t/8<t$. The already proved bounds $V\le C_\beta h^3$ and $\Lambda\le C_\beta h^2$ imply
\[
 \mathbb P\bigl(\mathcal E_{\mathrm{int}}(h,L)
             -\mathbb E\mathcal E_{\mathrm{int}}(h,L)>t\bigr)
 \le \exp\left[-c_\beta
             \min\{t^2/h^3,t/h^2\}\right].
\]
Together with $\operatorname{Tr}A\le C_\beta h$, this proves \eqref{A:eq:4.28}, uniformly in $L$.
\end{proof}

Fix a deterministic sequence $h_n\asymp n^{-1/3}$ and deterministic lengths $L_{h_n}$. Let $C_0$ satisfy $\mathbb E\mathcal E_{\mathrm{int}}(h,L)\le C_0h$ for all permitted $h,L$.
Taking $t=h_n$ in \eqref{A:eq:4.28} yields
\[
 \mathbb P\bigl(\mathcal E_{\mathrm{int}}(h_n,L_{h_n})
                       >(C_0+1)h_n\bigr)
 \le 2e^{-c_\beta/h_n}.
\]
There are deterministic $c_1>0$ and $n_1$ with $h_n^{-1}\ge c_1n^{1/3}$ for $n\ge n_1$. The last probabilities are summable: grouping the integers between $m^3$ and $(m+1)^3-1$
bounds their sum by a constant times $\sum_{m\ge1}(m+1)^2e^{-c_\beta c_1m}<\infty$. Borel--Cantelli gives a finite random $N_0$ such that
\[
 \mathcal E_{\mathrm{int}}(h_n,L_{h_n})\le (C_0+1)h_n
 \quad(n\ge N_0).
 \eqtag{4.29}
\]
Neither the Gaussian elements nor their errors at different sizes have been assumed independent. This argument is for the fixed deterministic sequence of pairs $(h_n,L_{h_n})$. The
stronger estimate $\mathcal E_{\mathrm{int}}(h,L)\le C_\beta hQ(B)$ established above is a deterministic inequality on $\Omega_B$ and hence applies there to every $0<h\le1$ and
every $L>0$ simultaneously. That simultaneous conclusion comes from the pathwise estimate, not from an uncountable intersection in Borel--Cantelli.

The original and approximating forms differ by
\[
\mathfrak q_{\beta,L}[f,v]
-
\mathfrak q_{\beta,h,L}^{\mathrm{fr}}[f,v]
=
\frac h2\langle f,v\rangle
-
\int_0^L(e-e_h)
\frac {\mathrm{d}}{\mathrm{d}x}(f\overline v)\,\mathrm{d}x.
\]
For finite $L$, both forms have domain $H_0^1(0,L)$ and $e-e_h$ is bounded there. Subtract their weak resolvent equations. As an operator from $L^2(0,L)$ to itself, the exact
difference is
\[
 R_{\beta,L}^{D}(z)-R_{\beta,h,L}^{\mathrm{fr}}(z)
 =-R_{\beta,L}^{D}(z)
       \bigl(\tfrac h2 I+\Gamma_{e-e_h}\bigr)
                       R_{\beta,h,L}^{\mathrm{fr}}(z).
\]
The products containing $\Gamma_{e-e_h}$ are evaluated through the form pairing. Their Hilbert--Schmidt bound is \eqref{A:eq:4.21}. The inequality
$\|U+V\|_{\mathfrak S_2}^2\le 2\|U\|_{\mathfrak S_2}^2+2\|V\|_{\mathfrak S_2}^2$ now gives
\[
\|R_{\beta,L}^D-R_{\beta,h,L}^{\mathrm{fr}}\|_{\mathfrak S_2}^2 \le\frac{h^2}{2}\|R_{\beta,L}^D
 R_{\beta,h,L}^{\mathrm{fr}} \|_{\mathfrak S_2}^2 +2\|R_{\beta,L}^D\Gamma_{e-e_h}
 R_{\beta,h,L}^{\mathrm{fr}}\|_{\mathfrak S_2}^2 \le C_{\omega,\beta,K} \bigl(h^2+\mathcal
 E_{\mathrm{int}}(h,L)\bigr).
\]
Here the first product is bounded by the operator norm of one resolvent and the uniform Hilbert--Schmidt norm of the other. Thus, on $\Omega_B$,
\[
\sup_{z\in K}
\left\|
R_{\beta,L}^{D}(z)
-
R_{\beta,h,L}^{\mathrm{fr}}(z)
\right\|_{\mathfrak S_2}^2
\le
C_{\omega,\beta,K}h.
\eqtag{4.31}
\]

\paragraph{Galerkin approximation}
Sampling the primitive gives the local regularity needed for interpolation. Let $I_j=(jh,(j+1)h)$, $0\le j\le N$, and put $u=R_{\beta,h,L}^{\mathrm{fr}}(z)F$. For a test function
supported inside $I_j$, the term with $e_h$ is the integral of a derivative multiplied by a constant, and therefore vanishes. The weak equation on $I_j$ is
\[
 -u''+(x-h/2-z)u=F.
\]
Its right side belongs to $L^2(I_j)$, so $u\in H^2(I_j)$. To see this directly, integrate $(x-h/2-z)u-F$ twice; the difference from $u$ has zero second distributional derivative
and therefore has the form $a_j+b_jx$ on $I_j$. Integration by parts on the intervals then gives, for every $v\in H_0^1(0,L)$,
\[
 \sum_{j=1}^N\bigl(-u'(jh+)+u'(jh-)
             +(e(jh)-e((j-1)h))u(jh)\bigr)\overline{v(jh)}=0.
\]
The nodal values of $v$ can be chosen independently, so
\[
 u'(jh+)-u'(jh-)
   =(e(jh)-e((j-1)h))u(jh),\qquad 1\le j\le N.
\]
Thus $u$ need not belong to $H^2(0,L)$: its derivative jumps at the nodes. The second derivative in the interpolation estimate is the ordinary $L^2$ derivative on each open cell.
In distributions on $(0,L)$ one has
\[
 u''=\sum_{j=0}^N\mathbf1_{I_j}u''|_{I_j}
       +\sum_{j=1}^N\bigl(e(jh)-e((j-1)h)\bigr)u(jh)\delta_{jh}.
\]
The point masses cancel the sampled potential in the equation for $u$. The sum of $L^2(I_j)$ norms below concerns only the derivatives on the open cells. The kernel estimate bounds
$(1+x)u$ uniformly in $L$; the equation on each cell then controls its second derivative.

\begin{lemma}\label{A:res:4.3A}
For \(z\in K\),
\[
\|(1+x)u\|_{L^2(0,L)}
+
\left(
\sum_{j=0}^{N}\|u''\|_{L^2(I_j)}^2
\right)^{1/2}
\le
C_{\omega,\beta,K}\|F\|_{L^2(0,L)}.
\eqtag{4.32}
\]
\end{lemma}

\begin{proof}
Write $R=R_{\beta,h,L}^{\mathrm{fr}}(z)$, $X=1+x$, and $Y=1+y$. By Proposition~\ref{A:res:4.1},
\[
 |R(x,y)|\le C_{\omega,K}X^{-1/4}Y^{-1/4}
                  \mathrm e^{-c|\Phi(x)-\Phi(y)|}.
\]
Since $1+\Phi(x)\asymp X^{3/2}$ and $1+\Phi(x)\le(1+\Phi(y))(1+|\Phi(x)-\Phi(y)|)$,
\[
 X/Y\le C(1+|\Phi(x)-\Phi(y)|)^{2/3}.
\]
For the row integral of the kernel $XR(x,y)$, put $v=\Phi(y)$. The Jacobian $\mathrm dy=Y^{-1/2}\mathrm dv$ gives
\[
 \int_0^L X|R(x,y)|\,\mathrm dy
 \le C_{\omega,K}\int_0^{\Phi(L)}(X/Y)^{3/4}
                         \mathrm e^{-c|\Phi(x)-v|}\,\mathrm dv
 \le C_{\omega,K}.
\]
Here $Y=(1+3v/2)^{2/3}$, and the last bound follows by integrating $(1+|\Phi(x)-v|)^{1/2}\mathrm e^{-c|\Phi(x)-v|}$ over the real line. For the column integral put $t=\Phi(x)$
instead. Then
\[
 \int_0^L X|R(x,y)|\,\mathrm dx
 \le C_{\omega,K}\int_0^{\Phi(L)}(X/Y)^{1/4}
                         \mathrm e^{-c|t-\Phi(y)|}\,\mathrm dt
 \le C_{\omega,K},
\]
using the integrable factor $(1+|t-\Phi(y)|)^{1/6}$. Both integrals are bounded by their counterparts over the whole real line in the $\Phi$ coordinate. The upper boundary
$\Phi(L)$ therefore contributes no factor. Schur's test gives $\|(1+x)R\|_{\mathrm{op}}\le C_{\omega,K}$ uniformly in $L$. The equation on each $I_j$ now yields
\[
 \left(\sum_{j=0}^N\|u''\|_{L^2(I_j)}^2\right)^{1/2}
 \le\|(x-h/2-z)u\|_2+\|F\|_2
 \le C_{\omega,K}\|F\|_2.
\]
Together these bounds prove \eqref{A:eq:4.32}.
\end{proof}

Define
\[
 V_{h,L}=\{v\in H_0^1(0,L):v|_{I_j}\text{ has degree at most one for }0\le j\le N\}.
\]
The values at $h,2h,\ldots,Nh$ determine $v$, so $\dim V_{h,L}=N$. Since $V_{h,L}\subset H_0^1(0,L)$, the form restricts to this space.

Let $\mathcal I_hu$ be the continuous piecewise linear function that agrees with $u$ at the points $jh$, where $u$ is continuous and belongs to $H^2$ on each interval. The map
\(\mathcal I_h\) acts on functions, whereas \(I_h\) acts on coordinate vectors. If \(u(0)=u(L)=0\), then \(\mathcal I_hu=I_h\bigl((h^{1/2}u(jh))_{j=1}^N\bigr)\). On $I_j$, put
$w=u-\mathcal I_hu$. Then $w$ vanishes at both boundary points and $w'=u'-h^{-1}\int_{I_j}u'$. Applying Poincar\'e's inequality first to the function $w'$ of mean zero and then to
$w$, whose values vanish at both ends, gives
\[
 \|w'\|_{L^2(I_j)}\le\frac h\pi\|u''\|_{L^2(I_j)},\qquad
 \|w\|_{L^2(I_j)}\le\frac h\pi\|w'\|_{L^2(I_j)}
                   \le\frac{h^2}{\pi^2}\|u''\|_{L^2(I_j)}.
\]
Consequently,
\[
 \|u-\mathcal I_hu\|_{\mathcal E_L}^2
 \le\left(\frac{h^2}{\pi^2}+
                  \frac{(1+L)h^4}{\pi^4}\right)
              \sum_{j=0}^N\|u''\|_{L^2(I_j)}^2.
\]
 The weighted $L^2$ error contributes $(1+L)h^4$, which remains of order $h^2$ because $1+L\le Ch^{-2}$. Therefore
\[
\|u-\mathcal I_hu\|_{\mathcal E_L}
\le
Ch
\left(
\sum_j\|u''\|_{L^2(I_j)}^2
\right)^{1/2}.
\eqtag{4.33}
\]
For the interpolation error $w=u-\mathcal I_hu$, all nodal values vanish, including the two boundary points. Therefore, for every $v\in H_0^1(0,L)$,
\[
 -\int_0^Le_h(w\overline v)'\,\mathrm dx
 =\sum_{j=1}^N\bigl(e(jh)-e((j-1)h)\bigr)w(jh)\overline{v(jh)}=0.
\]
The interpolation estimate thus follows from the derivatives inside each interval and the vanishing interpolation error at its two ends. Before using $L\le C_+h^{-2}$, the
preceding inequalities give
\[
 \|R_{\beta,h,L}^{\mathrm{fr}}(z)F
       -\mathcal I_hR_{\beta,h,L}^{\mathrm{fr}}(z)F\|_{\mathcal E_L}
 \le C_{\omega,K}h\sqrt{1+(1+L)h^2}\,\|F\|_2.
\]
In particular, \eqref{A:eq:4.32} and \eqref{A:eq:4.33} give $\inf_{v_h\in V_{h,L}}\|R_{\beta,h,L}^{\mathrm{fr}}(z)F-v_h\|_{\mathcal E_L} \le C_{\omega,K}h\|F\|_2$ in the stated
range of $L$.

Let \(P_{h,L}^{\mathrm{lin}}\) be the orthogonal projection of \(L^2(0,L)\) onto \(V_{h,L}\). For \(z\in K\), write
\[
a_z(u,v)
=
\mathfrak q_{\beta,h,L}^{\mathrm{fr}}[u,v]
-z\langle u,v\rangle,
\quad u,v\in H_0^1(0,L).
\]

For $F\in L^2(0,L)$, the exact solution $u=R_{\beta,h,L}^{\mathrm{fr}}(z)F$ satisfies $a_z(u,v)=\langle F,v\rangle$ for every $v\in H_0^1(0,L)$. We seek $u_h\in V_{h,L}$ for which
this equality holds for every $v_h\in V_{h,L}$. To write these $N$ equations in coordinates, let
\[
 \varphi_{j,h}(x)=h^{-1/2}\max\{1-|x-jh|/h,0\},
 \qquad 1\le j\le N.
\]
These functions are supported on $[(j-1)h,(j+1)h]$ and satisfy $\varphi_{j,h}(kh)=h^{-1/2}\delta_{jk}$. Thus $I_hc=\sum_{j=1}^Nc_j\varphi_{j,h}$, with $c_j=h^{1/2}(I_hc)(jh)$.
Define
\[
 (M_h)_{ij}=\langle\varphi_{j,h},\varphi_{i,h}\rangle_{L^2},
 \qquad
 (K_h^{\mathrm c})_{ij}
 =\mathfrak q_{\beta,h,L}^{\mathrm{fr}}
       [\varphi_{j,h},\varphi_{i,h}].
\]
Both matrices are Hermitian. The mass matrix satisfies $\|I_hc\|_{L^2}^2=c^*M_hc$ and is positive definite. The derivative contribution to $K_h^{\mathrm c}$ is $T_h=D_h^*D_h$.

Testing with $v_h=\varphi_{i,h}$ now gives
\[
 \sum_{j=1}^N\bigl((K_h^{\mathrm c})_{ij}-z(M_h)_{ij}\bigr)c_j
 =\langle F,\varphi_{i,h}\rangle,
 \qquad
 (K_h^{\mathrm c}-zM_h)c=I_h^*F.
\]
Here $(I_h^*F)_i=\int_0^L F(x)\overline{\varphi_{i,h}(x)}\,\mathrm{d}x$ is defined for every $F\in L^2(0,L)$. The matrix is invertible for nonreal $z$: taking the imaginary part of
$c^*(K_h^{\mathrm c}-zM_h)c=0$ gives $(\operatorname{Im}z)c^*M_hc=0$, hence $c=0$. The operator taking $F$ to $u_h$ is consequently
\[
 F\longmapsto I_h(K_h^{\mathrm c}-zM_h)^{-1}I_h^*F.
\]

The $L^2$ orthogonal projection is
\[
 P_{h,L}^{\mathrm{lin}}u=I_hM_h^{-1}I_h^*u,
 \qquad
 \langle u-P_{h,L}^{\mathrm{lin}}u,v_h\rangle=0
 \quad(v_h\in V_{h,L}).
\]
The formula follows by writing the projected function as $I_hc$ and testing its defining orthogonality against each $\varphi_{i,h}$. The solution $u_h\in V_{h,L}$ of the restricted
weak equation satisfies $a_z(u-u_h,v_h)=0$ for every $v_h\in V_{h,L}$, by subtraction of the two weak equations.

The imaginary part of $a_z$ controls the $L^2$ norm; its shifted real part then controls the energy norm. This gives stability on the space $V_{h,L}$ with the same constants as on
the full form domain. Let $u=R_{\beta,h,L}^{\mathrm{fr}}(z)F$, and let $u_h$ solve the finite weak equation above.

\begin{lemma}\label{A:res:4.3B}
Uniformly for $z\in K$,
\[
 \inf_{0\ne v_h\in V_{h,L}}\sup_{0\ne w_h\in V_{h,L}}
 \frac{|a_z(v_h,w_h)|}{\|v_h\|_{\mathcal E_L}\|w_h\|_{\mathcal E_L}}
 \ge c_{\omega,K}>0.
\]
The bound also holds with the two arguments interchanged. The solution of the restricted weak equation is unique and satisfies
\[
 \|u-u_h\|_{\mathcal E_L}
 \le C_{\omega,K}\inf_{v_h\in V_{h,L}}\|u-v_h\|_{\mathcal E_L}.
\]
\end{lemma}

\begin{proof}
Choose $\kappa_\omega\ge0$ so that $b(v,w)=\mathfrak q_{\beta,h,L}^{\mathrm{fr}}[v,w] +\kappa_\omega\langle v,w\rangle$ satisfies
$c_0\|v\|_{\mathcal E_L}^2\le b(v,v)\le C_0\|v\|_{\mathcal E_L}^2$. For $v_h\in V_{h,L}$, the reality of the quadratic form gives
\[
 \delta_K\|v_h\|_2^2\le
 |\operatorname{Im}a_z(v_h,v_h)|\le|a_z(v_h,v_h)|.
\]
Taking real parts and using this inequality,
\[
c_0\|v_h\|_{\mathcal E_L}^2 \le \operatorname{Re}a_z(v_h,v_h)
 +(\operatorname{Re}z+\kappa_\omega)\|v_h\|_2^2
 \le\left(1+\frac{|\operatorname{Re}z+\kappa_\omega|}{\delta_K}\right) |a_z(v_h,v_h)|.
\]
Put $C_K'=1+\delta_K^{-1}\sup_{z\in K}|\operatorname{Re}z+\kappa_\omega|$. Taking $w_h=v_h$ in the supremum gives the lower bound $c_0/C_K'$. It depends only on the common form
constants and $K$, not on the dimension or the interval length. The identity $a_z(v,w)=\overline{a_{\bar z}(w,v)}$ and $\bar z\in K$ give, for every nonzero $w_h$,
\[
 \sup_{v_h\ne0}\frac{|a_z(v_h,w_h)|}{\|v_h\|_{\mathcal E_L}}
 =\sup_{v_h\ne0}\frac{|a_{\bar z}(w_h,v_h)|}{\|v_h\|_{\mathcal E_L}}
 \ge c_{\omega,K}\|w_h\|_{\mathcal E_L}.
\]
Taking the infimum after division by $\|w_h\|_{\mathcal E_L}$ proves the second inf--sup bound. For the homogeneous finite system the preceding inequality gives
$\|v_h\|_{\mathcal E_L}=0$. Its kernel is therefore $\{0\}$, and the square coefficient matrix is invertible.

Cauchy--Schwarz for the positive form $b$ gives
\[
 |a_z(v,w)|\le |b(v,w)|+|z+\kappa_\omega|\|v\|_2\|w\|_2
 \le C_{\omega,K}\|v\|_{\mathcal E_L}\|w\|_{\mathcal E_L}.
\]
Subtracting the exact and finite weak equations yields $a_z(u-u_h,w_h)=0$ for every $w_h\in V_{h,L}$. Hence, for any
$v_h\in V_{h,L}$,
\[
c_{\omega,K}\|u_h-v_h\|_{\mathcal E_L} \le\sup_{w_h\ne0}
 \frac{|a_z(u_h-v_h,w_h)|}{\|w_h\|_{\mathcal E_L}} =\sup_{w_h\ne0}
 \frac{|a_z(u-v_h,w_h)|}{\|w_h\|_{\mathcal E_L}} \le C_{\omega,K}\|u-v_h\|_{\mathcal E_L}.
\]
The triangle inequality gives $\|u-u_h\|_{\mathcal E_L}\le (1+C_{\omega,K}/c_{\omega,K})\|u-v_h\|_{\mathcal E_L}$. Taking the infimum proves the result.
\end{proof}

Write $\mathcal G_{h,L}^{\mathrm{fr}}(z)F=u_h$ for the solution of the restricted weak equation, regarded as a function in $L^2(0,L)$. Applying Lemma~\ref{A:res:4.3B} with
\(v_h=\mathcal I_hu\), followed by \eqref{A:eq:4.32} and \eqref{A:eq:4.33}, gives
\[
\|u-u_h\|_{\mathcal E_L}
\le
C_{\omega,K}h\|F\|.
\eqtag{4.35}
\]

The energy error of order $h$ gains one further factor of $h$ when paired with the adjoint resolvent. This is the Aubin--Nitsche duality argument \cite{ciarlet1978finiteelement}
(Chapter~3); it uses piecewise $H^2$ regularity for the adjoint solution as well. Put \(e=u-u_h\) and let \(w=R_{\beta,h,L}^{\mathrm{fr}}(\bar z)e\). Since \(\bar z\in K\),
Lemma~\ref{A:res:4.3A} and all preceding estimates apply to \(w\) with the same constants. The adjoint variational equation is
\[
a_z(v,w)=\langle v,e\rangle,
\quad v\in H_0^1(0,L),
\]
which follows from \(a_z(v,w)=\overline{a_{\bar z}(w,v)}\). Taking \(v=e\) and using the orthogonality obtained by subtracting the weak equations, with
\(\mathcal I_hw\in V_{h,L}\) gives
\[
\|e\|^2
=
a_z(e,w)
=
a_z(e,w-\mathcal I_hw).
\]
The adjoint solution has zero boundary values and the same cell regularity as $u$. The interpolation estimate therefore applies on the identical partition. The continuity estimate
in Lemma~\ref{A:res:4.3B}, together with \eqref{A:eq:4.32}--\eqref{A:eq:4.33} applied to $w$, yields
\[
\|e\|^2 \le C_{\omega,K}\|e\|_{\mathcal E_L} \|w-\mathcal I_hw\|_{\mathcal E_L} \le
 C_{\omega,K}h\|e\|_{\mathcal E_L}\|e\|.
\]
Combining this with \eqref{A:eq:4.35} proves
\[
\left\|
\mathcal G_{h,L}^{\mathrm{fr}}(z)
-
R_{\beta,h,L}^{\mathrm{fr}}(z)
\right\|_{\mathrm{op}}
\le
C_{\omega,K}h^2.
\eqtag{4.37}
\]

Keeping the factor $\sqrt{1+(1+L)h^2}$ in both interpolation estimates gives, before specializing $L$,
\[
 \|\mathcal G_{h,L}^{\mathrm{fr}}(z)
                   -R_{\beta,h,L}^{\mathrm{fr}}(z)\|_{\mathrm{op}}
 \le C_{\omega,K}\bigl(h^2+(1+L)h^4\bigr).
\]
The two factors come respectively from the primal energy error and the adjoint approximation. For $e\ne0$, divide the duality estimate by $\|e\|_2$. For $e=0$, its left side is
zero.

The exact resolvent has infinite rank, so the operator norm estimate cannot be multiplied by $N$ for the full difference. We first control its part outside $V_{h,L}$, and apply the
rank bound only to the remaining finite-dimensional difference. Let $\mu_k$ and $\psi_k$ be the ordered eigenvalues and an orthonormal eigenbasis of
$\mathcal H_{\beta,h,L}^{\mathrm{fr}}$.

\begin{lemma}\label{A:res:4.3C}
For every \(v\in H_0^1(0,L)\),
\[
\|(I-P_{h,L}^{\mathrm{lin}})v\|_{L^2(0,L)}
\le
Ch\|v'\|_{L^2(0,L)}.
\]
With the fixed shift $\kappa_\omega$,
\[
\mu_k+\kappa_\omega\ge c_\omega k^{2/3},
\quad
\|(I-P_{h,L}^{\mathrm{lin}})\psi_k\|^2
\le
C_\omega\min\{1,h^2(\mu_k+\kappa_\omega)\}.
\]
Consequently, uniformly for \(z\in K\),
\[
\|(I-P_{h,L}^{\mathrm{lin}})
R_{\beta,h,L}^{\mathrm{fr}}(z)\|_{\mathfrak S_2}^2
+
\|R_{\beta,h,L}^{\mathrm{fr}}(z)
(I-P_{h,L}^{\mathrm{lin}})\|_{\mathfrak S_2}^2
\le
C_{\omega,K}h.
\eqtag{4.39}
\]
\end{lemma}

\begin{proof}
Every \(v\in H^1(0,L)\) has a continuous representative. On an interval \(I_j\), put
\[
c_j=\frac1h\int_{I_j}v'(x)\,\mathrm{d}x.
\]
The piecewise linear interpolant satisfies \((\mathcal I_hv)'=c_j\) on \(I_j\). Hence \(w_j:=v-\mathcal I_hv\) vanishes at both boundary points of \(I_j\), and \(w_j'=v'-c_j\).
Since \(c_j\) is the \(L^2(I_j)\)-orthogonal projection of \(v'\) onto the constants,
\[
\|w_j'\|_{L^2(I_j)}^2
=
\|v'\|_{L^2(I_j)}^2
-
\frac1h\left|\int_{I_j}v'(x)\,\mathrm{d}x\right|^2
\le
\|v'\|_{L^2(I_j)}^2.
\]
The Dirichlet Poincar\'e inequality on an interval of length \(h\) gives
\[
\|v-\mathcal I_hv\|_{L^2(I_j)}
\le
\frac h\pi\|w_j'\|_{L^2(I_j)}
\le
\frac h\pi\|v'\|_{L^2(I_j)}.
\]
After summing over the intervals and using the best approximation property of the orthogonal projection,
\[
\|(I-P_{h,L}^{\mathrm{lin}})v\|
\le
\|v-\mathcal I_hv\|
\le
\frac h\pi\|v'\|.
\]

The shifted form comparison gives
\[
c_\omega\|\psi_k\|_{\mathcal E_L}^2
\le
\mathfrak q_{\beta,h,L}^{\mathrm{fr}}[\psi_k]
+\kappa_\omega\|\psi_k\|^2
=
\mu_k+\kappa_\omega.
\]
In particular, \(\|\psi_k'\|^2\le C_\omega(\mu_k+\kappa_\omega)\), and the projection estimate gives the second inequality, together with
$\|(I-P_{h,L}^{\mathrm{lin}})\psi_k\|\le\|\psi_k\|=1$.

The eigenvalue lower bound is uniform in \(L\). Let \(\nu_k(L)\) be the eigenvalues of the deterministic form
\[
\mathfrak e_L[v]
=
\int_0^L|v'|^2\,\mathrm{d}x
+
\int_0^L(1+x)|v|^2\,\mathrm{d}x
\]
on \(H_0^1(0,L)\). The min--max principle and the shifted form comparison give
\[
\mu_k+\kappa_\omega = \min_{\substack{S\subset H_0^1(0,L)\\\dim S=k}} \max_{0\ne v\in S}
 \frac{\mathfrak q_{\beta,h,L}^{\mathrm{fr}}[v] +\kappa_\omega\|v\|^2}{\|v\|^2} \ge c_\omega
 \min_{\substack{S\subset H_0^1(0,L)\\\dim S=k}} \max_{0\ne v\in S} \frac{\mathfrak
 e_L[v]}{\|v\|^2} = c_\omega\nu_k(L).
\]
We prove the bound on $\nu_k(L)$ with a constant independent of $L$. Partition $(0,L)$ into the nonempty intervals
\[
 I_m=(m,\min\{m+1,L\}),\qquad m\in\mathbb N_0,\quad m<L.
\]
Enlarge the form domain from $H_0^1(0,L)$ to $\bigoplus_mH^1(I_m)$, thereby allowing independent boundary values on adjacent intervals, and replace $1+x$ on $I_m$ by its lower
bound $1+m$. The resulting form is smaller. On an interval of length $l_m\le1$, its eigenvalues are
\[
 1+m+(\pi j/l_m)^2,\qquad j=0,1,2,\ldots.
\]
They follow by solving the constant-coefficient equation with zero derivative at both ends: the eigenfunctions are the constant function and $\cos(\pi j(x-m)/l_m)$. Thus a block
with $m>E-1$ contributes no eigenvalue below $E$, and every other block contributes at most $1+\pi^{-1}\sqrt{E-1-m}$. Min--max and this count give, for $E\ge1$,
\[
 \begin{split}
 N_L(E)&:=\#\{k:\nu_k(L)\le E\}\\
 &\le\sum_{0\le m\le E-1}
       \bigl(1+\pi^{-1}\sqrt{E-1-m}\bigr)
 \le C E^{3/2}.
 \end{split}
\]
The count remains valid if the last interval has length less than one, because decreasing $l_m$ only increases its nonconstant Neumann eigenvalues. Since $\nu_k(L)\ge1$ and
$k\le N_L(\nu_k(L))$, we obtain
\[
 \nu_k(L)\ge C^{-2/3}k^{2/3},\qquad
 \mu_k+\kappa_\omega\ge c_\omega k^{2/3}.
\]
For comparison, extension by zero also gives $\nu_k(L)\ge\nu_k(\infty)$, and the half-line values $\nu_k(\infty)=1-a_k$ have the classical asymptotics for the zeros of the Airy
function in \cite{nistAiryZeros}. The bracketing argument above proves the lower bound used here without requiring those asymptotics.

Put \(\widehat\mu_k=\mu_k+\kappa_\omega\). For \(z\in K\),
\[
\frac{\widehat\mu_k}{|\mu_k-z|}
\le
1+\frac{|z+\kappa_\omega|}{|\mu_k-z|}
\le
1+\frac{\sup_{w\in K}|w+\kappa_\omega|}{\delta_K},
\]
and hence
\[
|\mu_k-z|^{-2}
\le
C_{\omega,K}\widehat\mu_k^{-2}.
\]
Choose the summation index $k_h=\lceil h^{-3}\rceil$. For $k\le k_h$ use $h^2\widehat\mu_k^{-1}\le C_\omega h^2k^{-2/3}$; for $k>k_h$ use $\widehat\mu_k^{-2}\le C_\omega k^{-4/3}$.
This split requires only the lower eigenvalue bound, not an upper bound or an identification of $k_h$ with an actual spectral cutoff. Expanding in the eigenbasis gives
\[
\begin{aligned}
\|(I-P_{h,L}^{\mathrm{lin}}) R_{\beta,h,L}^{\mathrm{fr}}(z)\|_{\mathfrak S_2}^2& = \sum_{k\ge1} \frac{\|(I-P_{h,L}^{\mathrm{lin}})\psi_k\|^2} {|\mu_k-z|^2} \le C_{\omega,K} \sum_{k\ge1} \min\{\widehat\mu_k^{-2},h^2\widehat\mu_k^{-1}\}\\
&\le C_{\omega,K} \left( h^2\sum_{k\le k_h}k^{-2/3} + \sum_{k>k_h}k^{-4/3} \right).
\end{aligned}
\]
Since
\[
h^2\sum_{k\le k_h}k^{-2/3}
\le Ch^2k_h^{1/3}\le Ch,
\quad
\sum_{k>k_h}k^{-4/3}
\le Ck_h^{-1/3}\le Ch,
\]
the first term in \eqref{A:eq:4.39} is bounded by \(C_{\omega,K}h\). The other term has the same bound because
\[
\|R_{\beta,h,L}^{\mathrm{fr}}(z)(I-P_{h,L}^{\mathrm{lin}})\|_{\mathfrak S_2}
=
\|(I-P_{h,L}^{\mathrm{lin}})R_{\beta,h,L}^{\mathrm{fr}}(\bar z)\|_{\mathfrak S_2}.
\]
Estimate \eqref{A:eq:4.39} follows.
\end{proof}

The right hand sides $F$ and $P_{h,L}^{\mathrm{lin}}F$ have the same inner products with every test function in $V_{h,L}$. Therefore the solution of the restricted weak equation
depends only on \(P_{h,L}^{\mathrm{lin}}F\), and its range is contained in
\(V_{h,L}\). Consequently,
\[
\mathcal G_{h,L}^{\mathrm{fr}}(z)
=
P_{h,L}^{\mathrm{lin}}
\mathcal G_{h,L}^{\mathrm{fr}}(z)
P_{h,L}^{\mathrm{lin}}.
\]
\[
\mathcal G_{h,L}^{\mathrm{fr}}(z) - P_{h,L}^{\mathrm{lin}}R_{\beta,h,L}^{\mathrm{fr}}(z)
 P_{h,L}^{\mathrm{lin}} = P_{h,L}^{\mathrm{lin}} \bigl(\mathcal G_{h,L}^{\mathrm{fr}}(z)
 -R_{\beta,h,L}^{\mathrm{fr}}(z)\bigr) P_{h,L}^{\mathrm{lin}}.
\]
The operator on the right has rank at most \(\dim V_{h,L}=N\), while \eqref{A:eq:4.37} and
\(\|P_{h,L}^{\mathrm{lin}}\|_{\mathrm{op}}=1\) give operator norm at most \(C_{\omega,K}h^2\). Since \(Nh\asymp
h^{-2}\),
\[
\left\|
\mathcal G_{h,L}^{\mathrm{fr}}(z)
-
P_{h,L}^{\mathrm{lin}}R_{\beta,h,L}^{\mathrm{fr}}(z)
P_{h,L}^{\mathrm{lin}}
\right\|_{\mathfrak S_2}^2
\le
NC_{\omega,K}h^4
\le
C_{\omega,K}h.
\eqtag{4.40}
\]
Finally,
\[
R_{\beta,h,L}^{\mathrm{fr}}(z) - P_{h,L}^{\mathrm{lin}}R_{\beta,h,L}^{\mathrm{fr}}(z)
 P_{h,L}^{\mathrm{lin}} ={} (I-P_{h,L}^{\mathrm{lin}})R_{\beta,h,L}^{\mathrm{fr}}(z) +
 P_{h,L}^{\mathrm{lin}}R_{\beta,h,L}^{\mathrm{fr}}(z) (I-P_{h,L}^{\mathrm{lin}}).
\]
Using the length dependent operator norm bound before \eqref{A:eq:4.39}, the finite-rank and projection contributions together give the more explicit estimate
\[
 \|\mathcal G_{h,L}^{\mathrm{fr}}(z)
                -R_{\beta,h,L}^{\mathrm{fr}}(z)\|_{\mathfrak S_2}^2
 \le C_{\omega,K}\left\{h+N\bigl[h^2+(1+L)h^4\bigr]^2\right\}.
\]
Here $h$ comes from the two infinite-dimensional projection tails; only the second term uses rank at most $N$. In the stated range, $N\le C_+h^{-3}$ and
$h^2+(1+L)h^4\le(2+C_+)h^2$. Thus the Hilbert--Schmidt triangle inequality, \eqref{A:eq:4.39}, and \eqref{A:eq:4.40} give
\[
\sup_{z\in K}
\left\|
\mathcal G_{h,L}^{\mathrm{fr}}(z)
-
R_{\beta,h,L}^{\mathrm{fr}}(z)
\right\|_{\mathfrak S_2}^2
\le
C_{\omega,\beta,K}h.
\eqtag{4.41}
\]

\paragraph{Quadrature and the mass matrix}
The weak equations in the piecewise linear basis involve $K_h^{\mathrm c}-zM_h$, whereas the difference resolvent involves $\mathcal B_{N,h,\beta}^B-zI$. We compute both
discrepancies, in the mass and in the potential, before estimating their effect between inverses. Integrating the hat basis functions gives
\[
 (M_h)_{jj}=2\int_0^1t^2\,\mathrm dt=\frac23,
 \qquad
 (M_h)_{j,j+1}=\int_0^1t(1-t)\,\mathrm dt=\frac16,
 \qquad (M_h)_{ij}=0\quad(|i-j|>1).
\]
Every hat function has two cells of support, including at the Dirichlet boundary points and when $N=1$. Hence
\[
\langle u,v\rangle
-
\langle I_hu,I_hv\rangle_{L^2}
=
\frac{h^2}{6}
\langle D_hu,D_hv\rangle.
\eqtag{4.43}
\]
Hence \(M_h=I-\frac{h^2}{6}T_h\) and \(\frac13 I\le M_h\le I\). Since $(I-M_h)_{jj}=1/3$ for every $h$, $M_h$ does not converge to $I$ in operator norm. We estimate this difference
between resolvents.

For a continuous function $g$, the composite trapezoidal rule is
\[
 Q_h(g)=\sum_{j=0}^{N}\frac h2\bigl(g(jh)+g((j+1)h)\bigr).
\]
If $g=(I_hu)\overline{(I_hv)}$, its boundary values vanish and the normalization of $I_h$ gives $Q_h(g)=\sum_{j=1}^Nu_j\overline{v_j}$. Thus the quadrature rule replaces the Gram
matrix by $I$. Apply the same quadrature rule to the deterministic potential. The kinetic term is already the exact difference energy. The Brownian term is also exact: the
derivative of $e_h$ is supported at the nodes, where the interpolation agrees with the coefficient vector. Thus this term produces no additional quadrature error.

For the potential term, on the interval \(I_j\), write \(p_h(x)=jh+h\left(t-\frac12\right), \allowbreak\, x=jh+ht.\) Writing $y=jh$, the three integrals required for the products
of the functions $1-t$ and $t$ are
\[
\begin{aligned}
\int_0^1\bigl(y+h(t-\tfrac12)\bigr)(1-t)^2\,\mathrm dt& =\frac y3-\frac h{12},\\
\int_0^1\bigl(y+h(t-\tfrac12)\bigr)t(1-t)\,\mathrm dt& =\frac y6,\qquad \int_0^1\bigl(y+h(t-\tfrac12)\bigr)t^2\,\mathrm dt =\frac y3+\frac h{12}.
\end{aligned}
\]
For the linear interpolants of $U_0,U_1$ and of $V_0,V_1$, substitution of these three values and subtraction from the trapezoidal sum gives
\[
\begin{aligned}
&\frac h2\left[ \left(jh-\frac h2\right)U_0\overline{V_0}
 + \left(jh+\frac h2\right)U_1\overline{V_1} \right]
 - \int_{I_j}p_hU\overline V\,\mathrm{d}x\\
&\quad= \frac{h\,jh}{6} (U_1-U_0)\overline{(V_1-V_0)}
 + \frac{h^2}{6} \left(U_1\overline{V_1}-U_0\overline{V_0}\right).
\end{aligned}
\]
The last terms telescope over all $N+1$ cells to $h^2(U(L)\overline{V(L)}-U(0)\overline{V(0)})/6=0$. Summing before taking absolute values preserves this cancellation at the two
Dirichlet boundary points. With $U=I_hu$, $V=I_hv$, their differences between the two boundary values are $U_1-U_0=h^{1/2}(D_hu)_j$ and $V_1-V_0=h^{1/2}(D_hv)_j$. Hence
\[
\sum_{k=1}^{N}x_k u_k\overline{v_k}
-
\int_0^Lp_h
(I_hu)\overline{(I_hv)}\,\mathrm{d}x
=
\sum_{j=0}^{N}
\frac{h^2(jh)}6
(D_hu)_j\overline{(D_hv)_j}.
\eqtag{4.44}
\]
For the random term, integration on each interval gives
\[
 -\int_0^L e_h\bigl((I_hu)\overline{(I_hv)}\bigr)'\,\mathrm dx
 =-\frac1h\sum_{j=0}^N e(jh)
              (u_{j+1}\overline{v_{j+1}}-u_j\overline{v_j}).
\]
Shifting the first sum by one index and using both zero boundary values, we obtain
\[
-
\int_0^Le_h
\frac {\mathrm{d}}{\mathrm{d}x}
\left((I_hu)
\overline{(I_hv)}\right)\,\mathrm{d}x
=
\sum_{k=1}^{N}
\frac{e(kh)-e((k-1)h)}h
u_k\overline{v_k}.
\eqtag{4.45}
\]
The right hand side is the Brownian diagonal potential of \(\mathcal B_{N,h,\beta}^{B}\).

Let $K_h^{\mathrm l}$ be the matrix obtained by quadrature. The derivative term is $T_h$, and \eqref{A:eq:4.44}--\eqref{A:eq:4.45} give its diagonal terms. Thus
$K_h^{\mathrm l}=T_h+X_h+M_{\omega_h}=\mathcal B_{N,h,\beta}^{B}$. For \(z\in K\), set
\(A_h^{\mathrm c}(z)=K_h^{\mathrm c}-zM_h, \allowbreak\, A_h^{\mathrm l}(z)=K_h^{\mathrm l}-zI\). Equations \eqref{A:eq:4.43}--\eqref{A:eq:4.45} give the matrix identity
\[
A_h^{\mathrm l}(z)-A_h^{\mathrm c}(z)
=
D_h^*M_{d_h(z)}D_h,
\quad
 d_{j,h}(z)=\frac{h^2}{6}(jh-z).
\eqtag{4.46}
\]
The mass and potential errors have combined into a single edge coefficient; in particular $d_{0,h}(z)=-h^2z/6$ is present. At the first and last vertices, for $u_0=u_{N+1}=0$,
\[
 \begin{aligned}
 (D_h^*M_dD_hu)_1&=h^{-2}\bigl[(d_0+d_1)u_1-d_1u_2\bigr],\\
 (D_h^*M_dD_hu)_N&=h^{-2}\bigl[(d_{N-1}+d_N)u_N-d_{N-1}u_{N-1}\bigr].
 \end{aligned}
\]
These formulas also agree for $N=1$, with both exterior values zero. Thus \eqref{A:eq:4.46} has no additional boundary correction. In particular,
\[
|d_{j,h}(z)|
\le
C_Kh^2(1+x_j),
\quad
|d_{j,h}(z)-d_{j-1,h}(z)|
\le
Ch^3.
\eqtag{4.47}
\]

\paragraph{Weighted estimates for the inverse of the matrix pencil}
Before applying the identity \eqref{A:eq:4.46}, we need weighted estimates for $(A_h^{\mathrm c}(z))^{-1}$ as well as for the ordinary difference resolvent. The mass matrix enters
the spectral representation and is tridiagonal in the coefficient coordinates. We keep those coordinates for the kernel estimates and use the $M_h$ inner product for
self-adjointness.

Set \(m_h[u,v]:=\langle u,M_hv\rangle, \allowbreak\, k_h^{\mathrm c}[u,v]:=\langle u,K_h^{\mathrm c}v\rangle\), and use $\ell_h$ from
Subsection~\ref{A:sec:weighted-resolvent-estimates-and-comparison-principles} with $N_n=N$. The identity
\(k_h^{\mathrm c}[u,v] = \mathfrak q_{\beta,h,L}^{\mathrm{fr}} [I_hu,I_hv]\) and the finite element norm equivalences
imply that, for each \(B\in\Omega_B\), there is a finite \(\kappa_\omega\), independent of \(h,L,z\), such that
\[
c_\omega\ell_h[u]
\le
k_h^{\mathrm c}[u,u]
+\kappa_\omega m_h[u,u]
\le
C_\omega\ell_h[u].
\eqtag{4.47b}
\]
Indeed, $\|(I_hu)'\|^2=\|D_hu\|^2$, and the weighted cell integration in Subsection~\ref{A:sec:weighted-resolvent-estimates-and-comparison-principles} gives
\[
\int_0^L(1+x)|I_hu|^2\,\mathrm dx
\asymp\sum_{j=1}^N(1+x_j)|u_j|^2.
\]
The term $-\int e_h(|I_hu|^2)'$ is infinitesimally form bounded relative to $\ell_h$. This follows from Lemma~\ref{A:res:3.1}, applied to $e_h$ and $I_hu$. Estimate
\eqref{A:eq:4.47b} follows.

The mass form has the coordinate representation
\[
m_h[u,v] ={} \frac23\sum_{j=1}^{N}u_j\overline{v_j} + \frac16\sum_{j=1}^{N-1} \left( u_{j+1}\overline{v_j} +u_j\overline{v_{j+1}} \right).
\eqtag{4.47c}
\]
Since $M_h=I-h^2T_h/6$ and $0\le h^2T_h\le4I$,
\[
\frac13I\le M_h\le I,
\quad
\|M_h^{-1}\|_{\mathrm{op}}\le3.
\eqtag{4.47d}
\]
We have $\langle u,v\rangle_{M_h}:=m_h[u,v]=\langle I_hu,I_hv\rangle_{L^2}$. The bounds $I/3\le M_h\le I$ give uniform equivalence with the Euclidean inner product. The operator
\(H_h^{\mathrm c}:=M_h^{-1}K_h^{\mathrm c}\) is self-adjoint in this inner product, because
\[
\langle H_h^{\mathrm c}u,v\rangle_{M_h}
=
\langle K_h^{\mathrm c}u,v\rangle
=
\langle u,K_h^{\mathrm c}v\rangle
=
\langle u,H_h^{\mathrm c}v\rangle_{M_h}.
\]
Moreover,
\[
A_h^{\mathrm c}(z)
=
M_h(H_h^{\mathrm c}-zI),
\quad
(A_h^{\mathrm c}(z))^{-1}
=
(H_h^{\mathrm c}-zI)^{-1}M_h^{-1}.
\eqtag{4.47f}
\]
Equivalently, in an orthonormal basis the matrix is $M_h^{-1/2}K_h^{\mathrm c}M_h^{-1/2}$, which is Hermitian in the ordinary Euclidean inner product. The map $I_hM_h^{-1/2}$ is an
isometry onto $V_{h,L}$, since its adjoint times itself is $I$. Consequently the same solution operator can be written as
\[
 \mathcal G_{h,L}^{\mathrm{fr}}(z)
 =(I_hM_h^{-1/2})
   (M_h^{-1/2}K_h^{\mathrm c}M_h^{-1/2}-zI)^{-1}
   (I_hM_h^{-1/2})^*.
\]
The last formula represents the embedded inverse of $K_h^{\mathrm c}-zM_h$ as an isometric image of a self-adjoint resolvent. For the kernel estimates we work in the original
coefficient coordinates, where $K_h^{\mathrm c}-zM_h$ is tridiagonal, and control the mass matrix in the conjugated form. The matrix $M_h^{-1/2}K_h^{\mathrm c}M_h^{-1/2}$ may have
entries outside this band.

For the inverse estimate between the energy space and its dual, let \(a_{h,z}^{\mathrm c}[u,v] := k_h^{\mathrm c}[u,v]-z m_h[u,v]\) and let \(\mathcal E_{h,L}\) be \(\mathbb C^N\)
with norm \(\ell_h^{1/2}\). If \(F\in\mathcal E_{h,L}^*\) and
\(a_{h,z}^{\mathrm c}[u,v]=F(v) \allowbreak\;(v\in\mathcal E_{h,L}),\) then
\[
\|u\|_{\mathcal E_{h,L}}
\le
C_{\omega,K}\|F\|_{\mathcal E_{h,L}^*}.
\eqtag{4.47g}
\]
To prove this, take \(v=u\). The imaginary part gives \(|\operatorname{Im}z|\,m_h[u,u] \le |F(u)| \le \|F\|_{\mathcal E_{h,L}^*} \|u\|_{\mathcal E_{h,L}}\). Using
\eqref{A:eq:4.47b},
\[
c_\omega\|u\|_{\mathcal E_{h,L}}^2 \le k_h^{\mathrm c}[u,u] +\kappa_\omega m_h[u,u] = \operatorname{Re}F(u) +(\operatorname{Re}z+\kappa_\omega)m_h[u,u] \le C_{\omega,K} \|F\|_{\mathcal E_{h,L}^*} \|u\|_{\mathcal E_{h,L}}.
\]
Division by \(\|u\|_{\mathcal E_{h,L}}\) proves \eqref{A:eq:4.47g}. In particular, for \(R_{\mathrm c}(z):=(A_h^{\mathrm c}(z))^{-1},\) we have
\[
\|R_{\mathrm c}(z)\|_{\mathrm{op}}
+
\|D_hR_{\mathrm c}(z)\|_{\mathrm{op}}
+
\|(1+X_h)^{1/2}R_{\mathrm c}(z)\|_{\mathrm{op}}
\le
C_{\omega,K}.
\eqtag{4.47h}
\]
Exponential conjugation of this matrix also requires control of the contribution of the mass form.

With the same $\Phi$ as above, for an index $i$ and $R<\infty$, define $\psi_j^{i,R}=\min\{|\Phi(x_j)-\Phi(x_i)|,R\}$ and
$\Theta_{i,a,R}=\operatorname{diag}(\mathrm e^{a\psi_j^{i,R}})$. Extend the weight constantly across each boundary edge, as in Lemma~\ref{A:res:3.12}. Then
\[
|\psi_{j+1}^{i,R}-\psi_j^{i,R}|
\le
Ch\sqrt{1+x_j}.
\eqtag{4.47i}
\]
The tridiagonal localization identity \eqref{A:eq:3.32c} applies separately to $K_h^{\mathrm c}$ and $M_h$; its proof applies to off-diagonal entries of either sign.

For a real weight at the interpolation points \(\theta_j\), \eqref{A:eq:4.47c} gives
\[
m_h[\theta u,\theta u] - \operatorname{Re}m_h[u,\theta^2u] = -\frac16\sum_{j=1}^{N-1} (\theta_{j+1}-\theta_j)^2 \operatorname{Re}(u_{j+1}\overline{u_j}).
\eqtag{4.47j}
\]
For \(\theta_j=\mathrm{e}^{a\psi_j^{i,R}}\),
\[
\left|
 m_h[\theta u,\theta u]
-
\operatorname{Re}m_h[u,\theta^2u]
\right|
\le
C a^2h^2\ell_h[\theta u].
\eqtag{4.47k}
\]
For reciprocal weights, put \(\delta_j=a(\psi_{j+1}^{i,R}-\psi_j^{i,R}).\) Again from \eqref{A:eq:4.47c},
\[
m_h[\Theta_{i,a,R}^{-1}u,\Theta_{i,a,R}v] -m_h[u,v] = \frac16\sum_{j=1}^{N-1} \Bigl[ (\mathrm{e}^{-\delta_j}-1)u_{j+1}\overline{v_j} + (\mathrm{e}^{\delta_j}-1)u_j\overline{v_{j+1}} \Bigr].
\eqtag{4.47l}
\]
For every edge in this interval, $h^2(1+x_j)\le h^2(1+L)\le1+C_+$. Consequently $|\delta_j|\le C a\sqrt{1+C_+}$ and $|\mathrm e^{\pm\delta_j}-1|\le C_{C_+}a$ for $0<a\le1$.
Substitution in \eqref{A:eq:4.47l}, followed by Cauchy--Schwarz and comparability of adjacent weights, gives
\[
\left|
 m_h[\Theta_{i,a,R}^{-1}u,\Theta_{i,a,R}v]
-m_h[u,v]
\right|
\le
Ca\|u\|_{\mathcal E_{h,L}}
\|v\|_{\mathcal E_{h,L}}.
\eqtag{4.47m}
\]

For $K_h^{\mathrm c}$, the derivative and the Brownian diagonal terms are already covered by the earlier conjugation calculation. The new contribution is the off-diagonal part of
the integrated deterministic potential. Write
\[
K_h^{\mathrm c}=T_h+P_h+W_h,
\eqtag{4.47n}
\]
where $(W_h)_{kk}=\omega_{k,h}$ and $(P_h)_{ij}=\int_0^Lp_h(x)\varphi_{j,h}(x)\varphi_{i,h}(x)\,\mathrm dx$. The latter matrix is tridiagonal. On the common support $I_j$ of
$\varphi_{j,h}$ and $\varphi_{j+1,h}$, the middle integral above gives $(P_h)_{j,j+1}=jh/6$. In particular,
\[
|(P_h)_{j,j+1}|
\le
C(1+x_j).
\eqtag{4.47o}
\]
The diagonal matrix \(W_h\), which contains the diagonal potential formed from scaled Brownian increments, cancels in both the localization and conjugation identities. For the
difference term, with \(r_j=\mathrm{e}^{a\psi_j^{i,R}}\), the product formulas are
\[
\begin{aligned}
D_h(r^{-1}u)_j
={}&
\mathrm{e}^{-a(\psi_j+\psi_{j+1})/2}
\left[
\cosh(\delta_j/2)D_hu_j
-
\frac2h\sinh(\delta_j/2)J_hu_j
\right],
\\
D_h(rv)_j
={}&
\mathrm{e}^{a(\psi_j+\psi_{j+1})/2}
\left[
\cosh(\delta_j/2)D_hv_j
+
\frac2h\sinh(\delta_j/2)J_hv_j
\right].
\end{aligned}
\eqtag{4.47p}
\]
Using
\[
|\cosh^2(\delta_j/2)-1| \le Ca^2h^2(1+x_j),\quad h^{-1}|\sinh\delta_j| \le Ca\sqrt{1+x_j},\quad h^{-2}\sinh^2(\delta_j/2) \le Ca^2(1+x_j),
\]
the four products in \eqref{A:eq:4.47p} have the form of \eqref{A:eq:3.42d} with $c_{j,n}=1$. The bounds for those four terms in Lemma~\ref{A:res:3.12} therefore give
\[
\left|
 \langle D_h(r^{-1}u),D_h(rv)\rangle
-
\langle D_hu,D_hv\rangle
\right|
\le
C(a+a^2)
\|u\|_{\mathcal E_{h,L}}
\|v\|_{\mathcal E_{h,L}}.
\eqtag{4.47q}
\]
For the off-diagonal part of \(P_h\), the reciprocal factor is \(\mathrm{e}^{\pm\delta_j}-1\). Equations \eqref{A:eq:4.47i} and \eqref{A:eq:4.47o}, together with
\(h\sqrt{1+x_j}\le C\), give
\[
\left| \langle \Theta_{i,a,R}^{-1}u, P_h\Theta_{i,a,R}v\rangle - \langle u,P_hv\rangle \right| \le Ca\sum_j(1+x_j) \left(|u_j||v_{j+1}|+|u_{j+1}||v_j|\right) \le Ca\|u\|_{\mathcal E_{h,L}} \|v\|_{\mathcal E_{h,L}}.
\eqtag{4.47r}
\]
For the localization identity, the potential contribution is bounded by
\[
 C\sum_j(1+x_j)(\theta_{j+1}-\theta_j)^2|u_ju_{j+1}|
 \le Ca^2\sum_j h^2(1+x_j)^2\theta_j^2|u_j|^2
 \le Ca^2\ell_h[\theta u],
\]
where adjacent weights are comparable and $h^2(1+x_j)\le C$. Combining \eqref{A:eq:4.47k}--\eqref{A:eq:4.47r}, and including the term \(-z m_h\), yields
\[
\left|
 a_{h,z}^{\mathrm c}
 [\Theta_{i,a,R}^{-1}u,\Theta_{i,a,R}v]
-
 a_{h,z}^{\mathrm c}[u,v]
\right|
\le
C_{\omega,K}a
\|u\|_{\mathcal E_{h,L}}
\|v\|_{\mathcal E_{h,L}}.
\eqtag{4.47s}
\]
Let $C_{i,a,R}(z):\mathcal E_{h,L}\to\mathcal E_{h,L}^*$ be the form difference in \eqref{A:eq:4.47s}. By \eqref{A:eq:4.47g},
\[
 \|R_{\mathrm c}(z)C_{i,a,R}(z)\|_{\mathcal E_{h,L}\to\mathcal E_{h,L}}
 \le C_{\omega,K}a.
\]
Choose one $0<a\le1$ so that this bound is at most $1/2$ for all $z\in K$. Its permitted size depends on $C_+$ through the bound on $\delta_j$, but is independent of $i,R,h,N$.
Then
\[
 \Theta_{i,a,R}R_{\mathrm c}(z)\Theta_{i,a,R}^{-1}
 =\sum_{m=0}^\infty
       (-R_{\mathrm c}(z)C_{i,a,R}(z))^m R_{\mathrm c}(z)
\]
converges as an operator from $\mathcal E_{h,L}^*$ to $\mathcal E_{h,L}$, with norm at most twice the bound in \eqref{A:eq:4.47g}. The constants are independent of $i,R,h$.

For $B_h=K_h^{\mathrm c}+\kappa_\omega M_h$, the same tridiagonal identity gives
\[
\langle\theta u,B_h\theta u\rangle -\operatorname{Re}\langle u,B_h\theta^2u\rangle = -\sum_{j=1}^{N-1}(B_h)_{j,j+1} (\theta_{j+1}-\theta_j)^2 \operatorname{Re}(u_{j+1}\overline{u_j}).
\]
The interval formulas imply \(|(B_h)_{j,j+1}|\le C_\omega h^{-2}\). Together with \eqref{A:eq:4.47i}, this yields
\[
\left|
\langle\theta u,B_h\theta u\rangle
-\operatorname{Re}\langle u,B_h\theta^2u\rangle
\right|
\le
C_\omega a^2\ell_h[\theta u].
\]

Set $G_h=B_h^{-1}$ and $X_j=1+x_j$. The form bounds for $B_h$ and the point-evaluation estimate give
\[
 c_\omega\ell_h[u]\le\langle u,B_hu\rangle,
 \qquad |u_j|^2\le ChX_j^{-1/2}\ell_h[u].
 \eqtag{4.47t}
\]
Fix $i$ and put $u=G_he_i$. Since $B_h$ is positive definite,
\[
 G_h(i,i)=\langle u,B_hu\rangle>0,
 \qquad u_i=G_h(i,i).
\]
Applying the point-evaluation inequality at $i$ yields
\[
 G_h(i,i)^2\le ChX_i^{-1/2}\ell_h[u]
             \le C_\omega hX_i^{-1/2}G_h(i,i).
\]
Division by the positive diagonal entry proves $G_h(i,i)\le C_\omega hX_i^{-1/2}$.

For the truncated weight $\theta_j=\exp(a_0\psi_j^{i,R})$, the identity $B_hu=e_i$ and $\theta_i=1$ give
\[
 \operatorname{Re}\langle u,B_h\theta^2u\rangle=G_h(i,i).
\]
The conjugation estimate just proved therefore implies
\[
 (c_\omega-C_\omega a_0^2)\ell_h[\theta u]\le G_h(i,i).
\]
Choose $a_0>0$ so small that $C_\omega a_0^2\le c_\omega/2$, independently of $i,R,h,N$. A second point-evaluation estimate now gives, for every vertex $j$,
\[
 \theta_j^2|G_h(j,i)|^2
 \le ChX_j^{-1/2}\ell_h[\theta u]
 \le C_\omega h^2X_i^{-1/2}X_j^{-1/2}.
\]
For fixed $h,N$, taking $R$ above $\max_j|\Phi_j-\Phi_i|$ removes the truncation without changing any constant. Hence
\[
 |G_h(j,i)|\le C_\omega hX_i^{-1/4}X_j^{-1/4}
                  e^{-a_0|\Phi_j-\Phi_i|}.
\]
Finally choose $0<a\le a_0/2$, also within the permitted range of the preceding estimate for the conjugated resolvent. The weighted sum in Lemma~\ref{A:res:3.10} gives
\[
 \begin{split}
 \sum_j e^{2a|\Phi_j-\Phi_i|}|G_h(j,i)|^2
 &\le C_\omega h^2X_i^{-1/2}
       \sum_jX_j^{-1/2}e^{-2(a_0-a)|\Phi_j-\Phi_i|}\\
 &\le C_\omega h^2X_i^{-1/2}\,h^{-1}X_i^{-1}
 =C_\omega hX_i^{-3/2}.
 \end{split}
 \eqtag{4.47w}
\]
This proves the normalized column bound for the positive matrix $B_h$, with constants independent of the number of vertices.

To pass from \eqref{A:eq:4.47w} for the shifted positive inverse to $A_h^{\mathrm c}(z)^{-1}$ at nonreal $z$, set
\[
C_{h,z}
:=
R_{\mathrm c}(z)B_h
=
I+(z+\kappa_\omega)R_{\mathrm c}(z)M_h.
\eqtag{4.47x}
\]
The factor $M_h$ in \eqref{A:eq:4.47x} cannot be commuted through $R_{\mathrm c}$. We bound it in its displayed position. The conjugated form estimate \eqref{A:eq:4.47s} gives
\[
\sup_{i,R}
\left\|
\Theta_{i,a,R}R_{\mathrm c}(z)M_h
\Theta_{i,a,R}^{-1}
\right\|_{\mathrm{op}}
\le
C_{\omega,K}.
\eqtag{4.47y}
\]
Indeed, the conjugated mass matrix has diagonal entries $2/3$ and adjacent entries $\mathrm e^{\pm\delta_j}/6$. Schur's test gives
\[
 \|\Theta_{i,a,R}M_h\Theta_{i,a,R}^{-1}\|_{\mathrm{op}}
 \le\frac23+\frac13\mathrm e^{Ca}.
\]
The conjugated inverse bound just proved also holds on $\ell^2$, since $\ell^2\subset\mathcal E_{h,L}^*$ and $\mathcal E_{h,L}\subset\ell^2$ continuously. Therefore
\[
\|\Theta R_{\mathrm c}(z)M_h\Theta^{-1}\|_{\mathrm{op}} \le\|\Theta R_{\mathrm
 c}(z)\Theta^{-1}\|_{\mathrm{op}} \|\Theta M_h\Theta^{-1}\|_{\mathrm{op}} \le C_{\omega,K},\qquad
 \Theta=\Theta_{i,a,R},
\]
which proves \eqref{A:eq:4.47y} for every $f\in\ell^2$. Equation \eqref{A:eq:4.47x} consequently implies
\[
\sup_{i,R}
\|\Theta_{i,a,R}C_{h,z}\Theta_{i,a,R}^{-1}\|_{\mathrm{op}}
\le
C_{\omega,K}.
\eqtag{4.47z}
\]
Since \(R_{\mathrm c}(z)=C_{h,z}G_h,\) \eqref{A:eq:4.47w}, \eqref{A:eq:4.47z}, and Fatou's lemma yield
\[
\sum_j
 \mathrm{e}^{2a|\Phi(x_j)-\Phi(x_i)|}
|R_{\mathrm c}(z)(j,i)|^2
\le
C_{\omega,K}h(1+x_i)^{-3/2}.
\eqtag{4.47aa}
\]
In particular,
\[
\|R_{\mathrm c}(z)e_i\|^2
\le
C_{\omega,K}h(1+x_i)^{-3/2}.
\eqtag{4.47ab}
\]
Write $C=R_{\mathrm c}(z)$. Since both coefficient matrices $K_h^{\mathrm c}$ and $M_h$ are self-adjoint,
\[
 C^*=(K_h^{\mathrm c}-\bar zM_h)^{-1}=R_{\mathrm c}(\bar z).
\]
Thus \eqref{A:eq:4.47aa} at $z$ controls the columns of $C$, and the same estimate at $\bar z$ controls its rows. Neither statement requires $C$ to be normal. For its first Gram
matrix,
\[
 (C^*C)(i,j)=\sum_k\overline{C(k,i)}C(k,j).
\]
The inequality $|\Phi_i-\Phi_j|\le|\Phi_k-\Phi_i|+|\Phi_k-\Phi_j|$ and Cauchy--Schwarz imply
\[
 \begin{split}
 |(C^*C)(i,j)|
 &\le e^{-a|\Phi_i-\Phi_j|}
  \left(\sum_ke^{2a|\Phi_k-\Phi_i|}|C(k,i)|^2\right)^{1/2}\\
 &\qquad\qquad\times
  \left(\sum_ke^{2a|\Phi_k-\Phi_j|}|C(k,j)|^2\right)^{1/2}.
 \end{split}
\]
Substitution of the two column bounds proves
\[
 \left|[R_{\mathrm c}(\bar z)R_{\mathrm c}(z)](i,j)\right|
 \le C_{\omega,K}h\,
 \frac{e^{-a|\Phi_i-\Phi_j|}}
      {(1+x_i)^{3/4}(1+x_j)^{3/4}}.
 \eqtag{4.47ac}
\]
For the opposite product use
\[
 (CC^*)(i,j)=\sum_k C(i,k)\overline{C(j,k)}
\]
and the row estimates, equivalently the column estimates for $R_{\mathrm c}(\bar z)$. The identical Cauchy--Schwarz calculation gives the same right hand side. In particular both
products in the trace expansion for a diagonal perturbation are bounded separately; no interchange of $R_{\mathrm c}(z)$ and $R_{\mathrm c}(\bar z)$ is used.

For the next proposition, each $R_1,R_2$ may be an inverse $(K_h^{\mathrm c}-zM_h)^{-1}$ as above or an ordinary resolvent whose operator satisfies \eqref{A:eq:3.8},
\eqref{A:eq:3.10}, and \eqref{A:eq:3.12}. All constants are common to the pair. Let $q$ be a complex vertex sequence, let $e$ be an edge sequence, and let $J\in\{J_-,J_+,J_h\}$,
with the boundary conventions already fixed.

\begin{proposition}\label{A:res:4.5A}
The estimates \eqref{A:eq:4.47h}, \eqref{A:eq:4.47aa}, and \eqref{A:eq:4.47ac} hold for $R_{\mathrm c}$. Moreover,
\[
 \|M_qR_{\mathrm c}(z)\|_{\mathfrak S_2}^2+
 \|R_{\mathrm c}(z)M_q\|_{\mathfrak S_2}^2
 \le C_{\omega,K}h\sum_j|q_j|^2(1+x_j)^{-3/2},
 \eqtag{4.47ad}
\]
\[
 \|R_1(z)M_qR_2(z)\|_{\mathfrak S_2}^2
 \le C_{\omega,K}h\sum_j|q_j|^2(1+x_j)^{-7/2},
 \eqtag{4.47ae}
\]
\[
 \|R_1(z)(D_h^*M_eJ+J^*M_eD_h)R_2(z)\|_{\mathfrak S_2}^2
 \le C_{\omega,K}h\sum_j|e_j|^2(1+x_j)^{-3/2}.
 \eqtag{4.47af}
\]
Finally, $\|D_hR_{\mathrm c}(z)\|_{\mathfrak S_2}^2\le C_{\omega,K}h^{-1}$.
\end{proposition}

\begin{proof}
By \eqref{A:eq:4.47ab} and $R_{\mathrm c}(z)^*=R_{\mathrm c}(\bar z)$,
\[
 \|M_qR_{\mathrm c}(z)\|_{\mathfrak S_2}^2
 +\|R_{\mathrm c}(z)M_q\|_{\mathfrak S_2}^2
 =\sum_j|q_j|^2\bigl(\|R_{\mathrm c}(\bar z)e_j\|^2
                    +\|R_{\mathrm c}(z)e_j\|^2\bigr)
 \le Ch\sum_j|q_j|^2(1+x_j)^{-3/2}.
\]
This is \eqref{A:eq:4.47ad}. To prove \eqref{A:eq:4.47ae}, write
\[
\|R_1M_qR_2\|_{\mathfrak S_2}^2 = \sum_{i,j}\overline{q_i}q_j K_1(i,j)K_2(j,i),
\]
where \(K_1=R_1(\bar z)R_1(z)\) and \(K_2=R_2(z)R_2(\bar z)\). Equation~\eqref{A:eq:4.47ac} applies to the inverse of the matrix pencil, and Proposition~\ref{A:res:3.14} applies to
an ordinary resolvent. Hence, for $\nu=1,2$, $|K_\nu(i,j)|\le Ch(1+x_i)^{-3/4}(1+x_j)^{-3/4} \mathrm e^{-c|\Phi_i-\Phi_j|}$. Using \(2|q_iq_j|\le |q_i|^2+|q_j|^2\), the estimate
reduces to
\[
Ch^2\sum_i|q_i|^2(1+x_i)^{-3/2}
\sum_j(1+x_j)^{-3/2}
 \mathrm{e}^{-c|\Phi(x_i)-\Phi(x_j)|}.
\]
The case $b=3/2$ in the proof of Lemma~\ref{A:res:3.10} bounds the last sum by $Ch^{-1}(1+x_i)^{-2}$, proving \eqref{A:eq:4.47ae}.

For \eqref{A:eq:4.47af}, split the first-order form into two terms. For example,
\[
\|R_1D_h^*M_eJR_2\|_{\mathfrak S_2}
\le
\|D_hR_1(\bar z)\|_{\mathrm{op}}
\|M_eJR_2\|_{\mathfrak S_2}.
\]
The first factor is bounded by \eqref{A:eq:4.47h}, and the second is bounded by \eqref{A:eq:4.47ad} together with $|(Ju)_j|^2\le C(|u_j|^2+|u_{j+1}|^2)$. For the other term,
\[
 \|R_1J^*M_eD_hR_2\|_{\mathfrak S_2}
 \le\|M_{\bar e}JR_1(\bar z)\|_{\mathfrak S_2}
                         \|D_hR_2\|_{\mathrm{op}}.
\]
This also proves the estimate for complex $e$.

For the derivative Hilbert--Schmidt estimate, use the Hermitian matrix
\[
\widetilde H_h^{\mathrm c}
=
M_h^{-1/2}K_h^{\mathrm c}M_h^{-1/2}.
\]
Then
\[
R_{\mathrm c}(z)
=
M_h^{-1/2}(\widetilde H_h^{\mathrm c}-z)^{-1}M_h^{-1/2}.
\]
By \eqref{A:eq:4.47d}, \(\|M_h^{-1/2}\|_{\mathrm{op}}^2\le3\). Therefore the Schatten ideal property gives
\[
\|D_hR_{\mathrm c}(z)\|_{\mathfrak S_2}^2 \le \|M_h^{-1/2}\|_{\mathrm{op}}^2 \left\|
D_hM_h^{-1/2}(\widetilde H_h^{\mathrm c}-z)^{-1} \right\|_{\mathfrak S_2}^2 \le 3\left\|
 D_hM_h^{-1/2}(\widetilde H_h^{\mathrm c}-z)^{-1} \right\|_{\mathfrak S_2}^2.
\]
Let \(\{\varphi_k\}_{k=1}^{N}\) be an orthonormal eigenbasis of \(\widetilde H_h^{\mathrm c}\), with \(\widetilde H_h^{\mathrm c}\varphi_k=\lambda_k\varphi_k\), and put
\(u_k=M_h^{-1/2}\varphi_k\). Then
\[
m_h[u_k,u_k]
=
\langle u_k,M_hu_k\rangle
=
1,
\quad
k_h^{\mathrm c}[u_k,u_k]
=
\langle u_k,K_h^{\mathrm c}u_k\rangle
=
\lambda_k.
\]
The lower form estimate in \eqref{A:eq:4.47b} consequently yields
\[
c_\omega\ell_h[u_k]
\le
\lambda_k+\kappa_\omega,
\quad
\|D_hM_h^{-1/2}\varphi_k\|^2
=
\|D_hu_k\|^2
\le
c_\omega^{-1}(\lambda_k+\kappa_\omega).
\]
Expanding the Hilbert--Schmidt norm in the eigenbasis \(\{\varphi_k\}\) gives the identity
\[
\left\|
D_hM_h^{-1/2}(\widetilde H_h^{\mathrm c}-z)^{-1}
\right\|_{\mathfrak S_2}^2
=
\sum_{k=1}^{N}
\frac{\|D_hM_h^{-1/2}\varphi_k\|^2}{|\lambda_k-z|^2}.
\]
The generalized min--max formula and \eqref{A:eq:4.47b}--\eqref{A:eq:4.47d} give
\[
\lambda_k+\kappa_\omega = \min_{\substack{S\subset\mathbb C^N\\ \dim S=k}} \max_{0\ne u\in S}
 \frac{k_h^{\mathrm c}[u,u]+\kappa_\omega m_h[u,u]}{m_h[u,u]} \ge c_\omega
 \min_{\substack{S\subset\mathbb C^N\\ \dim S=k}} \max_{0\ne u\in S} \frac{\ell_h[u]}{m_h[u,u]}
 \ge c_\omega \min_{\substack{S\subset\mathbb C^N\\ \dim S=k}} \max_{0\ne u\in S}
 \frac{\ell_h[u]}{\|u\|^2}.
\]
The last minimum is the \(k\)-th eigenvalue of the deterministic difference operator \(\mathcal L_h\). Hence \eqref{A:eq:3.18} implies
\(\lambda_k+\kappa_\omega\ge c_\omega k^{2/3}\). For \(t\ge0\) and \(z\in K\), the function
\[
(t,z)\longmapsto
\frac{t(1+t)}{|t-(z+\kappa_\omega)|^2}
\]
is bounded on \([0,\infty)\times K\): the denominator is bounded away from zero on bounded \(t\)-intervals because
\(K\Subset\mathbb C\setminus\mathbb R\), and the quotient tends to \(1\) as \(t\to\infty\), uniformly for \(z\in K\).
Hence
\[
\frac{\lambda_k+\kappa_\omega}{|\lambda_k-z|^2}
\le
\frac{C_{\omega,K}}{1+\lambda_k+\kappa_\omega}
\le
C_{\omega,K}k^{-2/3}.
\]
Combining the preceding estimates and using \(N\asymp h^{-3}\), which follows from \(Nh\asymp h^{-2}\), we obtain
\[
\|D_hR_{\mathrm c}(z)\|_{\mathfrak S_2}^2
\le C_{\omega,K}\sum_{k=1}^{N}k^{-2/3}
\le C_{\omega,K}N^{1/3}\le C_{\omega,K}h^{-1}.
\eqtag{4.47ag}
\]
This proves \eqref{A:eq:4.47ag}.
\end{proof}

\paragraph{Completion of the continuum comparison}
We can now estimate the quadrature error between the two inverses. The coefficient $d_{j,h}(z)$ in \eqref{A:eq:4.46} is not uniformly small when $jh\asymp h^{-2}$. As in
Lemma~\ref{A:res:product-rule}, we retain its dependence on the coefficient index and replace the Laplacian acting on the difference resolvent by its equation. The weighted
estimates then give the required order on the whole interval.

\begin{proposition}\label{A:res:4.5B}
Uniformly for \(z\in K\),
\[
\left\|
I_h(A_h^{\mathrm l}(z))^{-1}I_h^*
-
\mathcal G_{h,L}^{\mathrm{fr}}(z)
\right\|_{\mathfrak S_2}^2
\le
C_{\omega,\beta,K}h.
\eqtag{4.48}
\]
\end{proposition}

\begin{proof}
For \(F\in L^2(0,L)\), the coefficient vector of \(\mathcal G_{h,L}^{\mathrm{fr}}(z)F\) is \((A_h^{\mathrm c}(z))^{-1}I_h^*F\). Set \(R_{\mathrm l}=(A_h^{\mathrm l}(z))^{-1}\) and
\(R_{\mathrm c}=(A_h^{\mathrm c}(z))^{-1}\). Proposition~\ref{A:res:4.5A} supplies the weighted bounds for individual basis vectors and the mixed diagonal estimates for
\(R_{\mathrm c}\), with the same powers of \(h\) as in Section~\ref{A:sec:weighted-resolvent-estimates-and-comparison-principles}. Also,
\(R_{\mathrm c}(z)^*=R_{\mathrm c}(\overline z)\), because both \(K_h^{\mathrm c}\) and \(M_h\) are self-adjoint. The compact set was made conjugation-invariant at the start of the
subsection, so these adjoint bounds use the same constants.

By \eqref{A:eq:4.46},
\[
R_{\mathrm l}-R_{\mathrm c}
=-R_{\mathrm c}D_h^*M_{d_h(z)}D_hR_{\mathrm l}.
\]
Lemma~\ref{A:res:product-rule} applies to the complex edge sequence \(d_{j,h}(z)=h^2(jh-z)/6\), \(0\le j\le N\), with
\[
a_j=\frac{h^2}{6}\bigl((j-1)h-z\bigr),
\quad b_j=\frac{h^2}{6},
\quad 1\le j\le N.
\]
In particular, \(|a_j|\le C_Kh^2(1+x_j)\) and \(|b_j|\le Ch^2\), including the first and last vertices. The left edge value \(d_{0,h}(z)=-h^2z/6\) is retained in this calculation.
The resolvent after quadrature is that of
\(T_h+X_h+M_{\omega^B}\), so \(T_hR_{\mathrm l}=I+zR_{\mathrm l}-M_{x+\omega^B}R_{\mathrm l}\). Consequently,
\[
R_{\mathrm c}D_h^*M_dD_hR_{\mathrm l} ={} R_{\mathrm c}M_a+zR_{\mathrm c}M_aR_{\mathrm l}
 -R_{\mathrm c}M_{ax}R_{\mathrm l} -R_{\mathrm c}M_{a\omega^B}R_{\mathrm l} -R_{\mathrm
 c}M_bJ_-^*D_hR_{\mathrm l}.
\eqtag{4.51}
\]
For $x\in[(j-1)h,jh]$, $1+x$ and $1+x_j$ are comparable with absolute constants. Therefore the two sums needed below satisfy
\[
 \begin{aligned}
 h\sum_{j=1}^N(1+x_j)^{1/2}
 &\le C\int_0^L(1+x)^{1/2}\,\mathrm dx
   =\tfrac{2C}{3}\bigl((1+L)^{3/2}-1\bigr),\\
 h\sum_{j=1}^N(1+x_j)^{-3/2}
 &\le C\int_0^L(1+x)^{-3/2}\,\mathrm dx
   =2C\bigl(1-(1+L)^{-1/2}\bigr)\le2C.
 \end{aligned}
\]
The first and third terms of \eqref{A:eq:4.51} use the growing first integral. The other terms use the integrable weight in the second. In particular the first sum is $O(h^{-3})$
for $L\le C_+h^{-2}$. The multiplication estimate, applied to the adjoint where necessary, gives
\[
\|R_{\mathrm c}M_a\|_{\mathfrak S_2}^2
\le C_{\omega,K}h\sum_{j=1}^Nh^4(1+x_j)^{1/2}
\le C_{\omega,K}h^4(1+L)^{3/2}\le C_{\omega,K}h.
\]
For the second term, retain both inverses in \eqref{A:eq:4.47ae}. Since $|a_j|\le C_Kh^2(1+x_j)$,
\[
 \|zR_{\mathrm c}M_aR_{\mathrm l}\|_{\mathfrak S_2}^2
 \le C_{\omega,K}h^5\sum_{j=1}^N(1+x_j)^{-3/2}
 \le C_{\omega,K}h^4.
\]
The weighted estimate for the mixed diagonal product yields
\[
\|R_{\mathrm c}M_{ax}R_{\mathrm l}\|_{\mathfrak S_2}^2
 \le C_{\omega,K}h^5\sum_{j=1}^N(1+x_j)^{1/2}
 \le C_{\omega,K}h^4(1+L)^{3/2}\le C_{\omega,K}h.
\]
Write \(\Delta_jB=B(jh)-B((j-1)h)\), so \(\omega_j^B=c_\beta\Delta_jB/h\). The weighted increment bound on $\Omega_B$ gives
\[
\|R_{\mathrm c}M_{a\omega^B}R_{\mathrm l}\|_{\mathfrak S_2}^2 \le C_{\omega,\beta,K}h^3
 \sum_{j=1}^N(1+x_j)^{-3/2}|\Delta_jB|^2 \le C_{\omega,\beta,K}Q(B)h^3=o_\omega(h).
\]
For the remaining term, \eqref{A:eq:restricted-multiplier} gives
\[
\|R_{\mathrm c}M_bJ_-^*\|_{\mathfrak S_2}^2
=\|R_{\mathrm c}M_b\|_{\mathfrak S_2}^2
\le C_{\omega,K}h\sum_{j=1}^Nh^4(1+x_j)^{-3/2}
\le C_{\omega,K}h^4.
\]
Multiplication on the right by \(D_hR_{\mathrm l}\), whose operator norm is bounded by Proposition~\ref{A:res:3.3}, proves
\[
\|R_{\mathrm c}M_bJ_-^*D_hR_{\mathrm l}\|_{\mathfrak S_2}^2
\le C_{\omega,K}h^4.
\]
The inequality
\[
\Bigl\|\sum_{i=1}^5A_i\Bigr\|_{\mathfrak S_2}^2
\le5\sum_{i=1}^5\|A_i\|_{\mathfrak S_2}^2
\]
now gives
\[
\|R_{\mathrm l}-R_{\mathrm c}\|_{\mathfrak S_2}^2
\le C_{\omega,\beta,K}
       \bigl[h^4(1+L)^{3/2}+h^3Q(B)+h^4\bigr]
\le C_{\omega,\beta,K}h.
\]
Finally,
\[
I_hR_{\mathrm l}I_h^*
-\mathcal G_{h,L}^{\mathrm{fr}}(z)
=I_h(R_{\mathrm l}-R_{\mathrm c})I_h^*.
\]
Since \(\|I_h\|_{\mathrm{op}}^2=\|M_h\|_{\mathrm{op}}\le1\), the embedding cannot increase the Hilbert--Schmidt norm. This proves \eqref{A:eq:4.48}.
\end{proof}

Only the embedding remains to be changed. We compare $I_h$ and $\iota_h$ after composing with the difference resolvent. Their ranges are both contained in $L^2(0,L)$, although
$\iota_hu$ vanishes on the last cell $[Nh,L)$. The exact cell calculation gives a factor $h^2$ multiplying the squared derivative norm. On $[jh,(j+1)h)$, $0\le j<N$, the isometry
$\iota_hu$ equals $h^{-1/2}u_{j+1}$, while $I_hu=h^{-1/2}((1-t)u_j+tu_{j+1})$, $t=(x-jh)/h$. On the final interval $[Nh,(N+1)h)$, $\iota_hu=0$ and $u_{N+1}=0$. Thus all $N+1$
intervals give the same formula,
\[
 \|(I_h-\iota_h)u\|_{L^2(0,L)}^2
 =\sum_{j=0}^N|u_{j+1}-u_j|^2\int_0^1(1-t)^2\,\mathrm dt
 =\frac{h^2}{3}\|D_hu\|^2.
\]
Let $R_{N,h}^B(z)=(\mathcal B_{N,h,\beta}^B-z)^{-1}$ and $A_h=I_h-\iota_h$. Proposition~\ref{A:res:3.5} gives
\[
 \|A_hR_{N,h}^B(z)\|_{\mathfrak S_2}^2
 =\frac{h^2}{3}\|D_hR_{N,h}^B(z)\|_{\mathfrak S_2}^2
 \le C_{\omega,K}h^2N^{1/3}\le C_{\omega,K}h.
\]
The cell identity has been applied to every vector of an orthonormal basis and summed; all terms are nonnegative, so it also gives the stated Hilbert--Schmidt identity. This
derivative norm contains the last difference $-u_N/h$ as well as the first difference $u_1/h$. Since
\[
 I_hR_{N,h}^BI_h^*-\iota_hR_{N,h}^B\iota_h^*
 =A_hR_{N,h}^BI_h^*+\iota_hR_{N,h}^BA_h^*,
\]
$\|I_h\|_{\mathrm{op}}\le1$, $\|\iota_h\|_{\mathrm{op}}=1$, and $R_{N,h}^B(z)^*=R_{N,h}^B(\bar z)$, the two terms are bounded by the last display at $z$ and $\bar z$. Hence
\[
 \sup_{z\in K}
 \|I_hR_{N,h}^B(z)I_h^*-\iota_hR_{N,h}^B(z)\iota_h^*\|_{\mathfrak S_2}^2
 \le C_{\omega,K}h.
 \eqtag{4.55}
\]
Combining \eqref{A:eq:4.31}, \eqref{A:eq:4.41}, \eqref{A:eq:4.48}, and \eqref{A:eq:4.55}, we obtain
\[
\sup_{z\in K}
\left\|
\iota_h
(\mathcal B_{N,h,\beta}^{B}-z)^{-1}
\iota_h^*
-
R_{\beta,L}^{D}(z)
\right\|_{\mathfrak S_2}^2
\le
C_{\omega,\beta,K}h.
\eqtag{4.56}
\]
All finite interval operators in \eqref{A:eq:4.56} are extended by zero to $L^2(\mathbb R_+)$. The remaining comparison with the half-line is Lemma~\ref{A:res:4.2}.

\begin{proof}[Proof of Theorem~\ref{A:res:1.3}]
Work on the fixed set $\Omega_B$, and use the same path in every operator. Put $R_{N,h}^B(z)=(\mathcal B_{N,h,\beta}^B-z)^{-1}$. Here $\iota_h$ takes values in $L^2(\mathbb R_+)$,
while $I_h$ takes values in $L^2(0,L)$. Adding and subtracting the successive comparison operators gives the identity on $L^2(\mathbb R_+)$
\[
\begin{aligned}
 &\iota_hR_{N,h}^B(z)\iota_h^*-R_\beta(z)\\
 &\quad=\bigl(\iota_hR_{N,h}^B(z)\iota_h^*
                      -E_LI_hR_{N,h}^B(z)I_h^*P_L\bigr)\\
 &\qquad+E_L\bigl(I_hR_{N,h}^B(z)I_h^*-\mathcal G_{h,L}^{\mathrm{fr}}(z)\bigr)P_L\\
 &\qquad+E_L\bigl(\mathcal G_{h,L}^{\mathrm{fr}}(z)-R_{\beta,h,L}^{\mathrm{fr}}(z)\bigr)P_L\\
 &\qquad+E_L\bigl(R_{\beta,h,L}^{\mathrm{fr}}(z)-R_{\beta,L}^{D}(z)\bigr)P_L\\
 &\qquad+\bigl(E_LR_{\beta,L}^{D}(z)P_L-R_\beta(z)\bigr).
\end{aligned}
\]
The five squared Hilbert--Schmidt errors, uniformly for $z\in K$, are bounded respectively by
\[
\begin{aligned}
 &C_{\omega,\beta,K}h^2N^{1/3},\\
 &C_{\omega,\beta,K}\bigl[h^4(1+L)^{3/2}+h^3Q(B)+h^4\bigr],\\
 &C_{\omega,\beta,K}\bigl\{h+N[h^2+(1+L)h^4]^2\bigr\},\\
 &C_{\omega,\beta,K}\bigl[h^2+\mathcal E_{\mathrm{int}}(h,L)\bigr],\\
 &C_{\omega,\beta,K}L^{-1/2}.
\end{aligned}
\]
The first is the cell identity preceding \eqref{A:eq:4.55}; the second is Proposition~\ref{A:res:4.5B} with the weighted sums retained. The third combines the projection tails in
\eqref{A:eq:4.39} with the rank-$N$ part of the error from the variational restriction. The fourth follows from the form resolvent identity preceding \eqref{A:eq:4.31}, and the
fifth is

Lemma~\ref{A:res:4.2}. Thus each summand compares operators on the same Hilbert space with the stated boundary conditions.

Since $N\le C_+h^{-3}$, $1+L\le(1+C_+)h^{-2}$ and $\mathcal E_{\mathrm{int}}(h,L)\le C_\beta hQ(B)$,
\[
\begin{aligned}
 h^2N^{1/3}&\le C_+^{1/3}h,&
 h^4(1+L)^{3/2}&\le(1+C_+)^{3/2}h,\\
 N[h^2+(1+L)h^4]^2&\le C_+(2+C_+)^2h,&
 L^{-1/2}&\le c_-^{-1/2}h.
\end{aligned}
\]
The remaining terms are at most $C_\beta(1+Q(B))h$ for $h\le1$. Apply $\|\sum_{j=1}^5T_j\|_{\mathfrak S_2}^2\le 5\sum_{j=1}^5\|T_j\|_{\mathfrak S_2}^2$ to the exact difference
above. It follows that
\[
 \sup_{z\in K}\|\iota_hR_{N,h}^B(z)\iota_h^*-R_\beta(z)\|_{\mathfrak S_2}^2
 \le C_{\omega,\beta,K}h.
\]
All continuum form and kernel constants are independent of $L$. The discrete conjugation constants depend on $C_+$ through $h\sqrt{1+L}\le\sqrt{1+C_+}$, and the final tail constant
also depends on $c_-$. These are fixed constants in the hypotheses of the theorem. For the deterministic sequences $h_n,N_n$, the required inequalities hold for all sufficiently
large $n$. Substitution proves \eqref{A:eq:1.4} on $\Omega_B$, without a new exceptional set depending on the sequences or on $K$.
\end{proof}

\subsection{Discrete approximation of the inverse of the stochastic Bessel operator}
\label{sec:discrete-approximation-of-the-stochastic-bessel-inverse}
\label{B:sec:approximation-by-brownian-motion-and-the-lower-bound}

The comparison between the step kernel and the continuous kernel uses one Brownian path at every matrix size. On the selected cells, the step kernel has the exact scale masses
$Q_{n,k}^B$. We compare it first with the kernel having sampled $p^B$ and the original lower integration limit, and then with the inverse obtained from the restricted triangular
factor. The first difference is estimated cell by cell, including both triangles of each diagonal square. The second uses the weighted Brownian oscillations and the row estimates
from Section~\ref{B:sec:weighted-stability-and-truncation-of-the-inverse}. Their squared orders are $O((nr_n)^{-1})$ and $O(n^{-1})$, respectively. The factor restriction estimates
then remove the two intervals adjacent to the boundary. The maximum oscillation is estimated separately for the logarithmic coefficient comparisons.
Section~\ref{B:subsection:conditional-brownian-bridges} treats the conditional bridges needed for the lower bound. The resolvent identities in
Section~\ref{B:subsection:inverse-identities-and-embedded-resolvents} and the criterion in Section~\ref{B:subsection:the-criterion-for-convergence-of-the-inverses} are then applied
to the two finite ensembles.

\subsubsection{Partitions and sampled coefficients}
\label{B:subsection:brownian-data-and-grids}

We use $W,\mathcal Y,p^B,q^B$ and $R^B=\mathcal H_{\beta,a}^{-1}$ from \eqref{B:eq:1-29}--\eqref{B:eq:2-31}, with fixed $\beta>0$. The upper bounds below assume $a>0$; the lower
bound in Theorem~\ref{B:thm:brownian-information-lower-bound} holds for $a>-1$. The representation by an integral operator is the one used in \cite{ramirez2009diffusion}; its form
under the deterministic coordinate change for general potentials is described in \cite{rider2019universality}.

Let $z_{n,k}=\theta(k/n)$, where $\theta$ is the identity or the fixed map $\theta_V$ from \eqref{B:eq:1-15}. The cells $I_{n,k}$, their lengths $\Delta_{n,k}$, the times
$u_{n,k}=-\log z_{n,k}$, and the logarithmic lengths $h_{n,k}=h_{n,k}^-$ are as in Section~\ref{B:subsection:kernel-convention-and-normalized-cells}.
Lemma~\ref{B:lem:bessel-coordinate-regularity} gives constants $0<c<C$ with $c\le\theta'\le C$. Therefore
\begin{equation}
 c\frac{k}{n}\le z_{n,k}\le C\frac{k}{n},\qquad
 \frac cn\le\Delta_{n,k}\le\frac Cn,\qquad
 \frac{c'}k\le h_{n,k}\le\frac{C'}k\quad(k\ge2).
 \label{B:eq:brownian-grid-bounds}
\end{equation}
Indeed, the first two inequalities follow by integrating $\theta'$. For the last, integrate $s^{-1}$ over the cell to obtain
$\Delta_{n,k}/z_{n,k}\le h_{n,k}\le \Delta_{n,k}/z_{n,k-1}$, and use $k/(k-1)\le2$. Also $\Delta_{n,k}=z_{n,k}(1-e^{-h_{n,k}})$, as in \eqref{B:eq:logarithmic-cell-length}.

Fix $A>1$ and put $r_n=\lceil(\log(n+2))^A\rceil$. Choose deterministic integers $k_-(n),k_+(n)$ such that $k_-(n)\asymp r_n$ and $n-k_+(n)\asymp r_n$. All estimates concern
sufficiently large $n$, when $1\le k_-<k_+<n$. Set
\begin{equation}
 d_n=z_{n,k_-},\qquad \bar z_n=z_{n,k_+},\qquad
 P_n^0=\Pi_{d_n,\bar z_n},\qquad
 \mathscr I_n^{\mathrm{ret}}=\{k_-+1,\ldots,k_+\}.
 \label{B:eq:brownian-retained-interval}
\end{equation}
Here $\Pi_{d_n,\bar z_n}$ is the multiplication projection defined in Section~\ref{B:subsection:retained-cells-errors-and-weights}. Equation~\eqref{B:eq:brownian-grid-bounds} gives
$d_n\asymp r_n/n$, $1-\bar z_n\asymp r_n/n$, and $\max_{\mathscr I_n^{\mathrm{ret}}}h_{n,k}\le C/r_n$. We use the coefficients $p_{n,k}^B,Q_{n,k}^B$ and the step kernel
$R_n^{B,\mathrm{step},0}$ from \eqref{B:eq:2-32}--\eqref{B:eq:2-34} on these cells. Denote the inverse obtained from the restricted triangular factor by $R_n^{B,0}=R_{\beta,a}^0$
on $(d_n,\bar z_n)$; its kernel is \eqref{B:eq:bessel-restricted-green-kernel}.

Lemma~\ref{B:lem:bessel-brownian-growth} supplies all the power bounds for $p^B,q^B$ used below. Its event can be fixed for a countable sequence of positive exponents decreasing to
zero. A bound with a smaller exponent implies the bound with a larger one on $(0,1]$. Every scale mass in \eqref{B:eq:2-32} has positive integration limits. The divergence of
$\int_0^d q^B$ for $a>0$ was established in Section~\ref{B:subsection:retained-cells-errors-and-weights}.

\begin{lemma}\label{B:lem:brownian-grid-sums}
For $1<\lambda<2$, $0\le\gamma<1$, and $\zeta>0$,
\begin{align}
 \sum_{k\in\mathscr I_n^{\mathrm{ret}}}z_{n,k}^{\lambda}h_{n,k}^2
       &\le C_\lambda n^{-1},
       \label{B:eq:brownian-weight-sum}\\
 \sum_{k\in\mathscr I_n^{\mathrm{ret}}}h_{n,k}^2z_{n,k}^{2-\gamma}
       &\le C_\gamma n^{-1},\qquad
 \sum_{k\in\mathscr I_n^{\mathrm{ret}}}h_{n,k}z_{n,k}^{\zeta}
       \le C_\zeta.
       \label{B:eq:brownian-power-sums}
\end{align}
The constants are uniform in $n$ and in the selected interval.
\end{lemma}
\begin{proof}
All logarithmic lengths in the sums have index $k\ge2$. The sums of powers used to bound them may start at $1$. For $v\ge0$, monotonicity gives
\[
 \sum_{k=1}^n k^v
 \le\int_1^{n+1}x^v\,\dd x
 =\frac{(n+1)^{v+1}-1}{v+1}
 \le\frac{2^{v+1}}{v+1}n^{v+1}.
\]
For $-1<v<0$, the terms with $k\ge2$ satisfy $k^v\le\int_{k-1}^{k}x^v\,\dd x$, and hence
\[
 \sum_{k=1}^n k^v
 \le1+\int_1^n x^v\,\dd x
 =1+\frac{n^{v+1}-1}{v+1}
 \le\left(1+\frac1{v+1}\right)n^{v+1}.
\]
Use \eqref{B:eq:brownian-grid-bounds} separately in the three sums:
\begin{align*}
 \sum_{k\in\mathscr I_n^{\mathrm{ret}}}z_{n,k}^{\lambda}h_{n,k}^2
 &\le C n^{-\lambda}\sum_{k=1}^n k^{\lambda-2}
 \le C_\lambda n^{-\lambda}n^{\lambda-1}
 =C_\lambda n^{-1},\\
 \sum_{k\in\mathscr I_n^{\mathrm{ret}}}h_{n,k}^2z_{n,k}^{2-\gamma}
 &\le C n^{\gamma-2}\sum_{k=1}^n k^{-\gamma}
 \le C_\gamma n^{\gamma-2}n^{1-\gamma}
 =C_\gamma n^{-1},\\
 \sum_{k\in\mathscr I_n^{\mathrm{ret}}}h_{n,k}z_{n,k}^{\zeta}
 &\le C n^{-\zeta}\sum_{k=1}^n k^{\zeta-1}
 \le C_\zeta n^{-\zeta}n^\zeta=C_\zeta.
\end{align*}
The required exponents are respectively $\lambda-2>-1$, $-\gamma>-1$, and $\zeta-1>-1$. Nonnegativity permits the enlargement of each index set to $\{1,\ldots,n\}$, so the
constants are independent of the ends of the selected interval.
\end{proof}

\subsubsection{Weighted Brownian oscillations}
\label{B:subsection:weighted-brownian-oscillations}

Sampling and quadrature errors are both controlled by the oscillation of the Brownian path inside each logarithmic interval. The normalized oscillations have a common law within a
fixed row; weighted sums, rather than only their maximum, give the required rate.

For $k\ge2$ let $J_{n,k}=[u_{n,k},u_{n,k-1}]$ and define
\begin{equation}
 M_{n,k}=h_{n,k}^{-1/2}\operatorname{osc}_{J_{n,k}}W.
 \label{B:eq:brownian-normalized-oscillation}
\end{equation}
For a fixed $n$, define the continuous processes
\[
 B_{n,k}(t)=h_{n,k}^{-1/2}
       [W(u_{n,k}+h_{n,k}t)-W(u_{n,k})],\qquad 0\le t\le1.
\]
Their covariance within one interval is
\[
 \mathbb E[B_{n,k}(s)B_{n,k}(t)]
 =h_{n,k}^{-1}\min\{h_{n,k}s,h_{n,k}t\}=\min\{s,t\}.
\]
For distinct $k,\ell$, the two Brownian increments have disjoint interiors and their covariance is zero. Every finite family of evaluations is jointly Gaussian; its characteristic
function therefore factors over the interval index. Rational evaluations generate the Borel sigma-field on $C([0,1])$. Applying the monotone class theorem gives independence of the
processes $(B_{n,k})_{2\le k\le n}$ at this size. Thus $M_{n,k}=\operatorname{osc}_{[0,1]}B_{n,k}$ are independent and have a common distribution. Oscillation is a continuous
functional of a path in the uniform norm. The reflection principle gives
\begin{equation}
 \mathbb P(M_{n,k}>x)\le4e^{-x^2/8},\qquad x\ge0.
 \label{B:eq:brownian-oscillation-tail}
\end{equation}
Indeed, $\operatorname{osc}_{[0,1]}B\le2\sup_{t\le1}|B(t)|$ and
\[
 \mathbb P\left(\sup_{t\le1}|B(t)|>y\right)
 \le\mathbb P\left(\sup_{t\le1}B(t)>y\right)
   +\mathbb P\left(\inf_{t\le1}B(t)<-y\right)
 =4\mathbb P(B(1)>y)\le4e^{-y^2/2}.
\]
Take $y=x/2$ to obtain \eqref{B:eq:brownian-oscillation-tail}. The processes at different sizes are all obtained from
$W$.

\begin{lemma}\label{B:lem:weighted-brownian-oscillations}
For fixed $1<\lambda<2$, almost surely,
\begin{equation}
 \sum_{k\in\mathscr I_n^{\mathrm{ret}}}
    z_{n,k}^{\lambda}h_{n,k}^2M_{n,k}^2=O(n^{-1}).
 \label{B:eq:weighted-brownian-oscillations}
\end{equation}
On the same event,
\begin{equation}
 \max_{k\in\mathscr I_n^{\mathrm{ret}}}
       \operatorname{osc}_{J_{n,k}}W
 =O\!\left(\sqrt{\frac{\log(n+2)}{r_n}}\right)=o(1).
 \label{B:eq:maximum-brownian-oscillation}
\end{equation}
\end{lemma}
\begin{proof}
Fix $n$ and sum over $k\in\mathscr I_n^{\mathrm{ret}}$. Put $X_k=M_{n,k}^2$ and $a_k=z_{n,k}^{\lambda}h_{n,k}^2$. For every integer $m\ge1$, the tail formula and
\eqref{B:eq:brownian-oscillation-tail} give
\[
 \mathbb E X_k^m
 =m\int_0^\infty x^{m-1}\mathbb P(X_k>x)\,\dd x
 \le4m\int_0^\infty x^{m-1}e^{-x/8}\,\dd x
 =4\,8^m m!.
\]
In particular $\mu:=\mathbb E X_k\le32$. For $Y_k=X_k-\mu$, Jensen's inequality and the elementary power inequality yield
\[
 \mathbb E|Y_k|^m
 \le2^{m-1}(\mathbb E X_k^m+\mu^m)
 \le2^m\mathbb E X_k^m
 \le4\,16^m m!\le64^m m!.
\]
Take $K=64$. For $|t|Ka_k<1$, Tonelli's theorem applied to the nonnegative exponential series shows that $\mathbb E e^{|t|a_k|Y_k|}<\infty$. The series for $e^{ta_kY_k}$ may
consequently be integrated term by term. Since $\mathbb E Y_k=0$, for $|t|Ka_k\le1/2$,
\[
 \mathbb E e^{ta_kY_k}
 \le1+\sum_{m=2}^\infty(K|t|a_k)^m
 =1+\frac{K^2t^2a_k^2}{1-K|t|a_k}
 \le\exp(2K^2t^2a_k^2).
\]
Let $a_*:=\max_k a_k$ and $v_n:=\sum_k a_k^2$. At all sizes under consideration the retained set is nonempty, so $a_*,v_n>0$. Independence at this fixed size gives
\[
 \log\mathbb E\exp\left(t\sum_k a_kY_k\right)
 =\sum_k\log\mathbb E e^{ta_kY_k}
 \le2K^2t^2v_n,\qquad 0\le t\le(2Ka_*)^{-1}.
\]
For $s>0$ choose $t=\min\{s/(8K^2v_n),(4Ka_*)^{-1}\}$. Exponential Markov gives
\[
 \mathbb P\left(\sum_k a_kY_k>s\right)
 \le e^{-ts+2K^2t^2v_n}
 \le e^{-3ts/4},\qquad
 \frac{3ts}{4}
 =\min\left\{\frac{3s^2}{32K^2v_n},
                  \frac{3s}{16Ka_*}\right\}.
\]
This proves
\begin{equation}
 \mathbb P\left(\sum_k a_kY_k>s\right)
 \le\exp\left[-c\min\left\{
          \frac{s^2}{\sum_k a_k^2},\frac{s}{a_*}\right\}\right].
 \label{B:eq:weighted-oscillation-concentration}
\end{equation}
The deterministic coefficients satisfy
\[
 \sum_k a_k\le C/n,\qquad
 a_k\le C n^{-\lambda}k^{\lambda-2},\qquad
 a_*\le C n^{-\lambda}r_n^{\lambda-2},\qquad
 v_n\le a_*\sum_k a_k\le Ca_*/n.
\]
The first inequality is Lemma~\ref{B:lem:brownian-grid-sums}; the bound for $a_*$ uses $k\ge k_-+1\ge c_-r_n$ and $\lambda<2$. Choose a deterministic $C_0$ with
$\mathbb E\sum_k a_kX_k=\mu\sum_k a_k\le C_0/n$. Taking $s=1/n$ in \eqref{B:eq:weighted-oscillation-concentration}, we obtain
\begin{align*}
 \mathbb P\left(\sum_k a_kX_k>\frac{C_0+1}{n}\right)
 &\le\exp\left[-c\min\left\{\frac1{n^2v_n},\frac1{na_*}\right\}\right]\\
 &\le\exp\left[-\frac{c'}{na_*}\right]
 \le\exp[-c''n^{\lambda-1}r_n^{2-\lambda}].
\end{align*}
Since $r_n\ge1$ and $\lambda>1$, choose $n_0$ such that $c''n^{\lambda-1}\ge3\log(n+2)$ for every $n\ge n_0$. The displayed probabilities are then at most $(n+2)^{-3}$. Their sum
is finite, so the first Borel--Cantelli lemma gives \eqref{B:eq:weighted-brownian-oscillations} along the whole sequence.

For the maximum, choose $C_1$ such that $h_{n,k}\le C_1/r_n$ on every selected interval. The union bound and \eqref{B:eq:brownian-oscillation-tail} give
\begin{align*}
 &\mathbb P\left(\max_k\operatorname{osc}_{J_{n,k}}W
       >L\sqrt{\frac{\log(n+2)}{r_n}}\right)\\
 &\qquad\le\sum_k
 \mathbb P\left(M_{n,k}>L\sqrt{\frac{\log(n+2)}{r_nh_{n,k}}}\right)
 \le4n(n+2)^{-L^2/(8C_1)}.
\end{align*}
Take $L^2/(8C_1)>5$ and apply Borel--Cantelli. Finally,
\[
 \sqrt{\frac{\log(n+2)}{r_n}}
 \le[\log(n+2)]^{(1-A)/2}\longrightarrow0,
\]
since $A>1$. Intersecting the two events proves both limits on a common event. Both probability series were summed under the original law of the single path $W$; independence was
used only within a fixed row when computing the exponential moment.
\end{proof}

\subsubsection{Sampling and quadrature errors}
\label{B:subsection:brownian-sampling-and-quadrature}

There are two different uses of the path inside a cell. Sampling replaces the multiplicative coefficient by its value at one end, whereas quadrature replaces the integrated scale
density by its sampled value times the cell length.

On a selected cell $I_{n,k}$ set
\begin{equation}
 \epsilon_n^p(z)=\log\frac{p^B(z_{n,k})}{p^B(z)},\qquad
 \rho_{n,k}^B=\log\frac{\Delta_{n,k}q^B(z_{n,k})}{Q_{n,k}^B}.
 \label{B:eq:brownian-logarithmic-errors}
\end{equation}
The second expression is the quadrature term appearing in \eqref{B:eq:2-38}. The cell integrals depend on the path throughout each cell; the values $p^B(z_{n,k})$ and
$q^B(z_{n,k})$ use only its sampled values.

\begin{proposition}\label{B:prop:brownian-logarithmic-errors}
For fixed $1<\lambda<2$, almost surely,
\begin{equation}
 \int_{d_n}^{\bar z_n}z^{\lambda-1}|\epsilon_n^p(z)|^2\,\dd z
 +\sum_{k\in\mathscr I_n^{\mathrm{ret}}}
       h_{n,k}z_{n,k}^{\lambda}|\rho_{n,k}^B|^2
 =O(n^{-1}).
 \label{B:eq:brownian-logarithmic-error-bound}
\end{equation}
Moreover, $\|\epsilon_n^p\|_{L^\infty(d_n,\bar z_n)} +\max_{\mathscr I_n^{\mathrm{ret}}}|\rho_{n,k}^B|\to0$ almost surely.
\end{proposition}
\begin{proof}
For $z\in I_{n,k}$ write $z=z_{n,k}e^{-v}$, where $0\le v\le h_{n,k}$; boundary values can be assigned by continuity. The defining functions give the exact identities
\begin{align*}
 \epsilon_n^p(z_{n,k}e^{-v})
 &=\frac a2v-\beta^{-1/2}
       [W(u_{n,k}+v)-W(u_{n,k})],\\
 \log\frac{q^B(z_{n,k}e^{-v})}{q^B(z_{n,k})}
 &=(a+1)v-2\beta^{-1/2}
       [W(u_{n,k}+v)-W(u_{n,k})].
\end{align*}
The change of variables in this integral uses $\dd z=-z_{n,k}e^{-v}\,\dd v$. Hence
\[
 \frac{Q_{n,k}^B}{\Delta_{n,k}q^B(z_{n,k})}
 =\frac1{1-e^{-h_{n,k}}}
   \int_0^{h_{n,k}}
      e^{av-2\beta^{-1/2}[W(u_{n,k}+v)-W(u_{n,k})]}\,\dd v.
\]
Equivalently this is the average of $q^B(z_{n,k}e^{-v})/q^B(z_{n,k})$ with probability density $e^{-v}/(1-e^{-h_{n,k}})$ on $[0,h_{n,k}]$. Put
$L_{n,k}=(a+1)h_{n,k}+2\beta^{-1/2}h_{n,k}^{1/2}M_{n,k}$. The second logarithm lies between $-L_{n,k}$ and $L_{n,k}$. Averaging the corresponding exponential inequalities, and
taking the negative logarithm, gives
\[
 e^{-L_{n,k}}\le
 \frac{Q_{n,k}^B}{\Delta_{n,k}q^B(z_{n,k})}
 \le e^{L_{n,k}},\qquad
 |\rho_{n,k}^B|\le L_{n,k}.
\]
Together with the first exact identity, this yields
\begin{align*}
 |\epsilon_n^p(z)|
 &\le\frac a2h_{n,k}+\beta^{-1/2}h_{n,k}^{1/2}M_{n,k},\\
 |\rho_{n,k}^B|
 &\le(a+1)h_{n,k}+2\beta^{-1/2}h_{n,k}^{1/2}M_{n,k}.
\end{align*}
The weight integrates to
\[
 \int_{I_{n,k}}z^{\lambda-1}\,\dd z
 =\frac{z_{n,k}^{\lambda}}{\lambda}(1-e^{-\lambda h_{n,k}})
 \le h_{n,k}z_{n,k}^{\lambda}.
\]
Squaring the two bounds, using $(x+y)^2\le2x^2+2y^2$, and summing therefore bounds the left side of \eqref{B:eq:brownian-logarithmic-error-bound} by
\[
 C_{a,\beta}\sum_{k\in\mathscr I_n^{\mathrm{ret}}}
 z_{n,k}^{\lambda}(h_{n,k}^3+h_{n,k}^2M_{n,k}^2).
\]
The deterministic part satisfies
\[
 \sum_k z_{n,k}^{\lambda}h_{n,k}^3
 \le\left(\max_k h_{n,k}\right)
         \sum_k z_{n,k}^{\lambda}h_{n,k}^2
 \le\frac{C}{r_n}\frac{C_\lambda}{n}.
\]
Lemma~\ref{B:lem:weighted-brownian-oscillations} gives order $n^{-1}$ for the remaining sum on its event of probability one. On that same event, taking maxima in the two pointwise
bounds gives
\[
 \|\epsilon_n^p\|_{L^\infty(d_n,\bar z_n)}
 +\max_k|\rho_{n,k}^B|
 \le\frac{C_{a}}{r_n}
   +C_\beta\max_k\operatorname{osc}_{J_{n,k}}W
 \le\frac{C_a}{r_n}
   +C_{\beta,\omega}\sqrt{\frac{\log(n+2)}{r_n}}
 \longrightarrow0.
\]
This proves both estimates, with the sign of $\rho_{n,k}^B$ fixed by \eqref{B:eq:brownian-logarithmic-errors}.
\end{proof}

\subsubsection{Approximation of the Green kernel}
\label{B:subsection:brownian-green-kernel-approximation}

To separate the sampling of $p^B$ from the integration of $q^B$, set $p_\tau(z)=p^B(z)e^{\tau\epsilon_n^p(z)}$ on the selected interval, $0\le\tau\le1$, and let
\begin{equation}
 R_{n,\tau}(z,w)=\mathbf1_{(d_n,\bar z_n)^2}(z,w)
              p_\tau(z)p_\tau(w)\int_{z\vee w}^{\bar z_n}q^B(s)\,\dd s.
 \label{B:eq:brownian-kernel-interpolation}
\end{equation}
Then $R_{n,0}=R_n^{B,0}$. Write $R_n^{B,\mathrm{sample}}=R_{n,1}$.

\begin{proposition}\label{B:prop:brownian-sampled-kernel}
For fixed $1<\lambda<\min\{a+1,2\}$, almost surely,
\begin{equation}
 \|R_n^{B,\mathrm{sample}}-R_n^{B,0}\|_{\mathfrak S_2}^2
 \le C_\omega\int_{d_n}^{\bar z_n}
               z^{\lambda-1}|\epsilon_n^p(z)|^2\,\dd z=O(n^{-1}).
 \label{B:eq:brownian-sampled-kernel-bound}
\end{equation}
\end{proposition}
\begin{proof}
Choose $0<\eta<\tfrac1{16}\min\{a+1-\lambda,2-\lambda\}$. For a fixed path in the event of Lemma~\ref{B:lem:bessel-brownian-growth}, the same coefficient bounds hold on every
selected interval. Let $h_*$ bound all $h_{n,k}$ with $k\ge2$, as permitted by \eqref{B:eq:brownian-grid-bounds}. For $z\in I_{n,k}$, $1\le z_{n,k}/z\le e^{h_*}$, so the geometric
mean satisfies
\begin{align*}
 p_\tau(z)
 &=p^B(z)^{1-\tau}p^B(z_{n,k})^\tau\\
 &\le C_{\omega,\eta}z^{a/2-\eta}
       (z_{n,k}/z)^{\tau(a/2-\eta)}
 \le C_{\omega,\eta}e^{|a/2-\eta|h_*}z^{a/2-\eta}.
\end{align*}
The constant is independent of $n,k,\tau$. For $z\le w$,
\begin{align*}
 |R_{n,\tau}(z,w)|
 &\le C_\omega z^{a/2-\eta}w^{a/2-\eta}
          \int_w^1s^{-a-1-2\eta}\,\dd s\\
 &=\frac{C_\omega}{a+2\eta}
     z^{a/2-\eta}w^{a/2-\eta}(w^{-a-2\eta}-1)
 \le C_\omega z^{a/2-\eta}w^{-a/2-3\eta}.
\end{align*}
For $w\le z$ the same kernel formula gives $|R_{n,\tau}(z,w)|\le C_\omega w^{a/2-\eta}z^{-a/2-3\eta}$. The two parts of the squared row can now be integrated explicitly:
\begin{align*}
 \int_{d_n}^{z}|R_{n,\tau}(z,w)|^2\,\dd w
 &\le C_\omega z^{-a-6\eta}\int_0^z w^{a-2\eta}\,\dd w
 =\frac{C_\omega}{a+1-2\eta}z^{1-8\eta}
 \le C_\omega z^{\lambda-1},\\
 \int_z^{\bar z_n}|R_{n,\tau}(z,w)|^2\,\dd w
 &\le C_\omega z^{a-2\eta}\int_z^{\bar z_n}w^{-a-6\eta}\,\dd w\\
 &=C_\omega z^{\lambda-1}\int_z^{\bar z_n}
       (z/w)^{a+1-\lambda-2\eta}w^{1-\lambda-8\eta}\,\dd w\\
 &\le C_\omega z^{\lambda-1}\int_0^1w^{1-\lambda-8\eta}\,\dd w
 =\frac{C_\omega}{2-\lambda-8\eta}z^{\lambda-1}.
\end{align*}
Here $a+1-2\eta>0$, $a+1-\lambda-2\eta>0$, and $2-\lambda-8\eta>0$. Adding the two estimates gives
\begin{equation}
 \int_{d_n}^{\bar z_n}|R_{n,\tau}(z,w)|^2\,\dd w
       \le C_\omega z^{\lambda-1}.
 \label{B:eq:brownian-interpolated-row-bound}
\end{equation}
It also gives $\|R_{n,\tau}\|_{\mathfrak S_2}^2 \le C_\omega\int_0^1z^{\lambda-1}\,\dd z=C_\omega/\lambda$.

For fixed $n$ put $b_n(z,w)=\epsilon_n^p(z)+\epsilon_n^p(w)$ and $L_n=2\|\epsilon_n^p\|_\infty<\infty$. For $\tau,\tau+t\in[0,1]$, $R_{n,\tau+t}=R_{n,\tau}e^{tb_n}$ as kernels. The
scalar integral remainder
\[
 e^x-1-x=x^2\int_0^1(1-s)e^{sx}\,\dd s
\]
implies
\[
 \left\|\frac{R_{n,\tau+t}-R_{n,\tau}}t
                 -b_nR_{n,\tau}\right\|_{\mathfrak S_2}
 \le\frac{|t|L_n^2}{2}e^{|t|L_n}
                   \|R_{n,\tau}\|_{\mathfrak S_2}\longrightarrow0.
\]
In this formula $b_nR_{n,\tau}$ denotes the pointwise product of the kernels. Moreover,
\[
 \|b_nR_{n,\tau}-b_nR_{n,\sigma}\|_{\mathfrak S_2}
 \le L_n^2|\tau-\sigma|e^{L_n|\tau-\sigma|}
                         \|R_{n,\sigma}\|_{\mathfrak S_2}.
\]
Thus the Hilbert--Schmidt derivative exists and is continuous, with $\partial_\tau R_{n,\tau}(z,w)=b_n(z,w)R_{n,\tau}(z,w)$. Symmetry of $R_{n,\tau}$ gives
\begin{align*}
 \|\partial_\tau R_{n,\tau}\|_{\mathfrak S_2}^2
 &\le2\int\!\!\int
 (|\epsilon_n^p(z)|^2+|\epsilon_n^p(w)|^2)
                  |R_{n,\tau}(z,w)|^2\,\dd z\,\dd w\\
 &=4\int_{d_n}^{\bar z_n}|\epsilon_n^p(z)|^2
           \int_{d_n}^{\bar z_n}|R_{n,\tau}(z,w)|^2\,\dd w\,\dd z\\
 &\le C_\omega\int_{d_n}^{\bar z_n}
                    z^{\lambda-1}|\epsilon_n^p(z)|^2\,\dd z.
\end{align*}
Both variables in the first line range over $(d_n,\bar z_n)$. The last constant is independent of $n,\tau$. The fundamental theorem of calculus and Cauchy--Schwarz in $\tau$ yield
\begin{align*}
 \|R_{n,1}-R_{n,0}\|_{\mathfrak S_2}^2
 &=\int\!\!\int
   \left|\int_0^1\partial_\tau R_{n,\tau}(z,w)\,\dd\tau\right|^2
                                             \,\dd z\,\dd w\\
 &\le\int_0^1\|\partial_\tau R_{n,\tau}\|_{\mathfrak S_2}^2\,\dd\tau
 \le C_\omega\int_{d_n}^{\bar z_n}
                     z^{\lambda-1}|\epsilon_n^p(z)|^2\,\dd z.
\end{align*}
Proposition~\ref{B:prop:brownian-logarithmic-errors} proves the stated order on its event of probability one.
\end{proof}

\begin{proposition}\label{B:prop:brownian-integration-limit}
Almost surely,
\begin{equation}
 \|R_n^{B,\mathrm{step},0}-R_n^{B,\mathrm{sample}}\|_{\mathfrak S_2}^2
       \le \frac{C_\omega}{nr_n}=o(n^{-1}).
 \label{B:eq:brownian-integration-limit-bound}
\end{equation}
\end{proposition}
\begin{proof}
Fix $z\in I_{n,i}$, $w\in I_{n,j}$ and put $m=i\vee j$. The consecutive cells from $m$ to $k_+$ partition $[z_{n,m-1},\bar z_n)$, so
\[
 \sum_{\ell=m}^{k_+}Q_{n,\ell}^B
 =\int_{z_{n,m-1}}^{\bar z_n}q^B(s)\,\dd s.
\]
The two kernels use the same factors $p_{n,i}^Bp_{n,j}^B$. Their difference is therefore exactly
\begin{equation}
 E_n(z,w)=p_{n,i}^Bp_{n,j}^B
                         \int_{z_{n,m-1}}^{z\vee w}q^B(s)\,\dd s.
 \label{B:eq:brownian-integration-limit-error}
\end{equation}
For $i<j$, the upper limit is $w$; for $i>j$, it is $z$. For $i=j$, both triangles of the cell square are included by $z\vee w$. In every case
$0\le E_n(z,w)\le p_{n,i}^Bp_{n,j}^BQ_{n,m}^B$. This is a pointwise kernel inequality.

Define $a_i=\Delta_{n,i}(p_{n,i}^B)^2$, $M_m^B=\sum_{k_-<i\le m}a_i$, and $M_{k_-}^B=0$. The pointwise estimate above and the exact identity
\[
 \sum_{i\vee j=m}a_ia_j
 =a_m^2+2a_mM_{m-1}^B
 =(M_m^B)^2-(M_{m-1}^B)^2\le2a_mM_m^B
\]
control the complete cell sum. In particular, the diagonal cell contributes $a_m^2$ once; both cell masses are retained.

Choose $0<\eta<\min\{(a+1)/4,1/16\}$, so $a+1-2\eta>0$ and $0<8\eta<1$. The Brownian power bounds, \eqref{B:eq:logarithmic-cell-length}, and
Lemma~\ref{B:lem:logarithmic-power-sums} give
\begin{align*}
 a_m&\le C_\omega h_{n,m}z_{n,m}^{a+1-2\eta},\\
 M_m^B&\le C_\omega\sum_{k_-<i\le m}
                       h_{n,i}z_{n,i}^{a+1-2\eta}
       \le C_\omega z_{n,m}^{a+1-2\eta},\\
 Q_{n,m}^B&\le C_\omega
       \int_{z_{n,m-1}}^{z_{n,m}}s^{-a-1-2\eta}\,\dd s
       \le C_\omega h_{n,m}z_{n,m}^{-a-2\eta}.
\end{align*}
The last inequality uses $(e^{(a+2\eta)h}-1)/(a+2\eta)\le e^{(a+2\eta)h_*}h$. Integrating the pointwise bound over each cell pair now gives
\begin{align*}
 \|E_n\|_{\mathfrak S_2}^2
 &\le\sum_m(Q_{n,m}^B)^2\sum_{i\vee j=m}a_ia_j
 \le2\sum_m(Q_{n,m}^B)^2a_mM_m^B\\
 &\le C_\omega\sum_m h_{n,m}^3
     z_{n,m}^{-2a-4\eta+(a+1-2\eta)+(a+1-2\eta)}\\
 &=C_\omega\sum_m h_{n,m}^3z_{n,m}^{2-8\eta}
 \le\frac{C_\omega}{r_n}
       \sum_mh_{n,m}^2z_{n,m}^{2-8\eta}
 \le\frac{C_\omega}{nr_n}.
\end{align*}
The final estimate is \eqref{B:eq:brownian-power-sums} with $\gamma=8\eta$. Since $r_n\to\infty$, multiplication by $n$ gives the required little-o bound.
\end{proof}

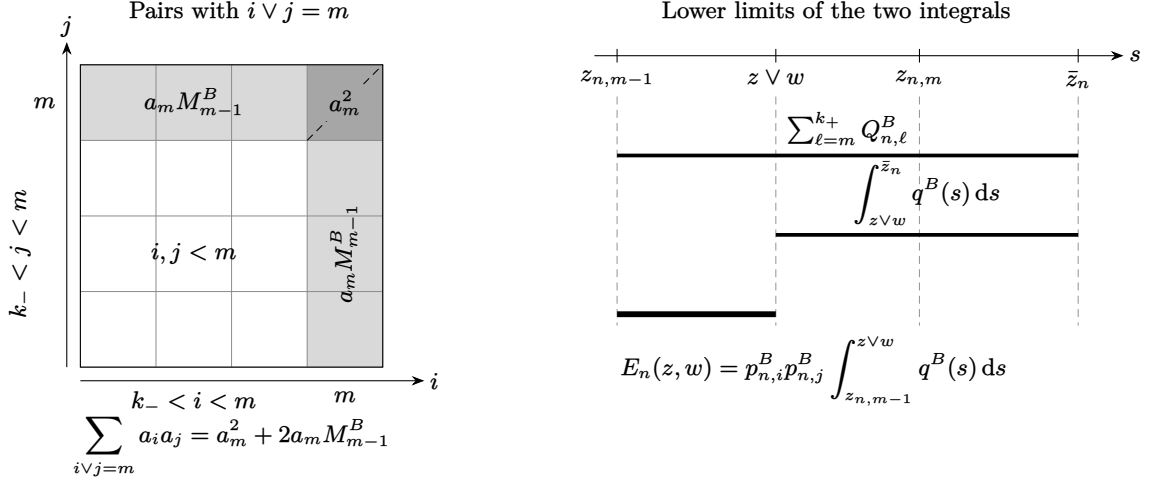
\begin{figure}[H]
\centering
\begin{tikzpicture}[x=1cm,y=1cm,every node/.style={font=\footnotesize,inner sep=2pt},>=Stealth]
 \node at (2.1,4.7) {Pairs with $i\vee j=m$};
 \fill[black!15] (0,3) rectangle (3,4);
 \fill[black!15] (3,0) rectangle (4,3);
 \fill[black!35] (3,3) rectangle (4,4);
 \foreach \x in {0,1,2,3,4}{\draw[black!45,line width=.3pt] (\x,0)--(\x,4);\draw[black!45,line width=.3pt] (0,\x)--(4,\x);}
 \draw (0,0) rectangle (4,4);
 \draw[densely dashed] (3,3)--(4,4);
 \draw[->] (0,-.18)--(4.55,-.18) node[right] {$i$};
 \draw[->] (-.18,0)--(-.18,4.3) node[above] {$j$};
 \node[below] at (3.5,-.25) {$m$};
 \node[below] at (1.5,-.25) {$k_-<i<m$};
 \node[left] at (-.25,3.5) {$m$};
 \node[rotate=90] at (-.8,1.5) {$k_-<j<m$};
 \node at (1.5,1.5) {$i,j<m$};
 \node at (1.5,3.5) {$a_mM_{m-1}^B$};
 \node[rotate=90] at (3.5,1.5) {$a_mM_{m-1}^B$};
 \node[fill=black!35] at (3.5,3.5) {$a_m^2$};
 \node at (2,-1.05) {$\displaystyle\sum_{i\vee j=m}a_ia_j=a_m^2+2a_mM_{m-1}^B$};
 \node at (10,4.7) {Lower limits of the two integrals};
 \draw[->] (6.8,4.12)--(13.8,4.12) node[right] {$s$};
 \foreach \x/\y in {7.1/{z_{n,m-1}},9.2/{z\vee w},11.1/{z_{n,m}},13.2/{\bar z_n}}{
   \draw (\x,4.03)--(\x,4.21);
   \node[below] at (\x,4.0) {$\y$};
   \draw[black!40,densely dashed] (\x,3.58)--(\x,.55);
 }
 \draw[line width=1.4pt] (7.1,2.8)--(13.2,2.8);
 \node[above] at (10.15,2.85) {$\sum_{\ell=m}^{k_+}Q_{n,\ell}^B$};
 \draw[line width=1.4pt] (9.2,1.75)--(13.2,1.75);
 \node[above] at (11.2,1.8) {$\displaystyle\int_{z\vee w}^{\bar z_n}q^B(s)\,\dd s$};
 \draw[line width=2.2pt] (7.1,.7)--(9.2,.7);
 \node[below,align=center] at (9.7,.47)
  {$\displaystyle E_n(z,w)=p_{n,i}^Bp_{n,j}^B\int_{z_{n,m-1}}^{z\vee w}q^B(s)\,\dd s$};
\end{tikzpicture}
\caption{The sum over cell pairs in Proposition~\ref{B:prop:brownian-integration-limit}. For fixed $m$, the shaded row and column contain exactly the pairs $i\vee j=m$. The two off-diagonal parts each have mass $a_mM_{m-1}^B$, and their intersection has mass $a_m^2$. The dashed diagonal divides that cell according to whether $z\vee w=z$ or $z\vee w=w$. On the right, the step kernel integrates from $z_{n,m-1}$, while the sampled kernel integrates from $z\vee w$. Their difference is the displayed integral over part of one cell. Its squared contribution is bounded by $(Q_{n,m}^B)^2(a_m^2+2a_mM_{m-1}^B)$ before summation over $m$.}
\label{B:fig:green-kernel-cell-sum}
\end{figure}

The grouping by $i\vee j$ in Figure~\ref{B:fig:green-kernel-cell-sum} accounts for both copies of an off-diagonal cell and for the single diagonal cell. The second panel shows the
interval responsible for the error in \eqref{B:eq:brownian-integration-limit-error}.

\subsubsection{Restriction of the triangular factor of the limiting inverse}
\label{B:subsection:restriction-of-the-brownian-factor}

The truncation estimates from Section~\ref{B:subsection:restriction-of-the-triangular-factor} extend the comparison on the retained interval to the full inverses when
$\mathcal K=\mathcal K_{\beta,a}$. For $0<d<b\le1$ with $b\ge1/2$, denote the corresponding restricted inverse by $R_{d,b}^{B,0}$.

\begin{proposition}\label{B:prop:brownian-factor-truncation}
For $0<\eta<\nu_a:=\min\{a,1\}$, almost surely,
\begin{equation}
 \|R^B-R_{d,b}^{B,0}\|_{\mathfrak S_2}^2
 \le C_{\omega,\eta}
       \bigl[d^{1+\nu_a-\eta}+(1-b)^3+(1-b)^2\bigr].
 \label{B:eq:brownian-factor-truncation}
\end{equation}
The estimates of Proposition~\ref{B:prop:truncation-row-bounds} and Lemma~\ref{B:lem:factor-restriction-correction} also hold for the Brownian step factors on any consecutive cells
with indices at least two, with constants uniform in the chosen cells and in $n$.
\end{proposition}
\begin{proof}
Fix an outcome in the event of Lemma~\ref{B:lem:bessel-brownian-growth}, and abbreviate $p=p^B$, $q=q^B$, $R=R^B$. Choose
\[
 0<\xi<\frac1{16}\min\{a,\eta,1\},\qquad
 \nu=\nu_a=\min\{a,1\}.
\]
One finite constant, depending on this path and $\xi$, gives
\[
 p(s)\le C s^{a/2-\xi},\qquad
 q(s)\le C s^{-a-1-2\xi},\qquad 0<s\le1.
\]
In particular $p\in L^2(0,1)$. For $s\le t$, integration of the second bound yields
\[
 |R(s,t)|=p(s)p(t)\int_t^1q(u)\,\mathrm du
 \le C s^{a/2-\xi}t^{-a/2-3\xi}.
\]
Let $A(s)=\int_0^1|R(s,t)|^2\,\mathrm dt$. Using symmetry and splitting the integral at $t=s$ gives
\[
 \begin{split}
 A(s)&\le C s^{-a-6\xi}\int_0^s t^{a-2\xi}\,\mathrm dt
        +Cs^{a-2\xi}\int_s^1t^{-a-6\xi}\,\mathrm dt\\
 &\le C s^{1-8\xi}
        +Cs^{a-2\xi}\int_s^1t^{-a-6\xi}\,\mathrm dt.
 \end{split}
\]
If $a+6\xi<1$, the last integral is at most $(1-a-6\xi)^{-1}$. If $a+6\xi>1$, it is at most $s^{1-a-6\xi}/(a+6\xi-1)$. In the equality case it is $|\log s|$; for any $\delta>0$,
$s^\delta|\log s|\le(e\delta)^{-1}$. Our choice $8\xi<\eta$ ensures in all three cases that
\[
 A(s)\le C_{\omega,\eta}s^{\nu-\eta},\qquad
 \int_0^d A(s)\,\mathrm ds
       \le C_{\omega,\eta}d^{1+\nu-\eta}.
\]
In the equality case $a+6\xi=1$, use $\delta=(\eta-2\xi)/2>0$ in the logarithmic estimate; then $a-2\xi-\delta>a-\eta=\nu-\eta$.

The other boundary admits a stronger estimate. For $t\ge1/2$, $p$ and $q$ are bounded on $[1/2,1]$. If $s\le t$, then
\[
 |R(s,t)|\le C p(s)(1-t),
\]
whereas if $s\ge t$, then $|R(s,t)|\le C(1-s)$. Consequently
\[
 A(t)\le C(1-t)^2\int_0^t p(s)^2\,\mathrm ds
            +C\int_t^1(1-s)^2\,\mathrm ds
       \le C(1-t)^2,
 \qquad \int_b^1A(t)\,\mathrm dt\le C(1-b)^3.
\]
All constants so far are independent of $d,b$.

Let $P$ be multiplication by $\mathbf1_{(d,b)}$. The kernel of $R-PRP$ is the original kernel restricted to the complement of $(d,b)^2$. Symmetry and a union bound on these two
coordinate strips give
\[
 \|R-PRP\|_{\mathfrak S_2}^2
 \le2\int_{(0,d]\cup[b,1)}A(s)\,\mathrm ds
 \le C_{\omega,\eta}\bigl(d^{1+\nu-\eta}+(1-b)^3\bigr).
\]
Restriction of the triangular factor is not compression of its adjoint product. Direct kernel subtraction gives
\[
 PRP-R_{d,b}^{B,0}
   =Q_b\,|Pp\rangle\langle Pp|,
 \qquad Q_b=\int_b^1q(t)\,\mathrm dt.
\]
Since $b\ge1/2$, $Q_b\le C(1-b)$, and therefore
\[
 \|PRP-R_{d,b}^{B,0}\|_{\mathfrak S_2}^2
       =Q_b^2\|Pp\|_2^4\le C(1-b)^2.
\]
This rank-one kernel is supported on $(d,b)^2$, disjoint from the kernel of $R-PRP$. The two errors are orthogonal in Hilbert--Schmidt space. Adding their squared norms proves
\eqref{B:eq:brownian-factor-truncation}. This is the continuum estimate of Theorem~\ref{B:thm:inverse-truncation-bound}, here with the boundary contributions calculated directly.

For the discrete estimate choose $2\le L\le U\le n$. For $L\le k\le U$ define $p_k=p^B(z_{n,k})$ and $r_k=(Q_{n,k}^B/\Delta_{n,k})^{1/2}$; set $p_k=r_k=0$ at every other index. In
particular $p_1=r_1=0$. The full measure remains $\mathfrak m_n=\sum_{k=1}^n\Delta_{n,k}\delta_{z_{n,k}}$. For $k\ge2$, \eqref{B:eq:brownian-grid-bounds} gives
\[
 z_{n,k-1}\ge \frac{c(k-1)}n
 \ge\frac{c}{2C}z_{n,k}.
\]
Put $c_1=c/(2C)>0$. For every point of the $k$th cell,
\[
 q^B(s)\le C_{\omega,\xi}s^{-a-1-2\xi}
 \le C_{\omega,\xi}c_1^{-a-1-2\xi}z_{n,k}^{-a-1-2\xi}.
\]
Integration over its physical length and taking square roots give
\[
 Q_{n,k}^B\le C_{\omega,\xi}\Delta_{n,k}z_{n,k}^{-a-1-2\xi},
 \qquad r_k\le C_{\omega,\xi}z_{n,k}^{-(a+1)/2-\xi}.
\]
Thus, at nonzero indices with $i\le j$,
\[
 p_i r_j\le C_{\omega,\xi}
 z_{n,i}^{a/2-\xi}z_{n,j}^{-(a+1)/2-\xi}
 =C_{\omega,\xi}z_{n,j}^{-1/2}
      (z_{n,i}/z_{n,j})^{a/2}(z_{n,i}z_{n,j})^{-\xi}.
\]
At the remaining indices the left side is zero. The full atomic measure satisfies \eqref{B:eq:full-grid-measure} with the original constants $c,C$. The atomic versions of
Proposition~\ref{B:prop:truncation-row-bounds} and Lemma~\ref{B:lem:factor-restriction-correction} therefore apply to this factor extended by zero, with constants independent of
$L,U,n$. For the last atom, the inclusive mass is
\[
 \mathfrak m_n([z_{n,n},1])=\Delta_{n,n}>0.
\]
Consequently the estimates include the whole last cell in the second kernel variable. The restriction is made in the coefficients; the measure keeps all of its original cell
masses.
\end{proof}

\begin{theorem}\label{B:thm:brownian-inverse-approximation}
For the fixed parameters and retained intervals above, almost surely,
\begin{equation}
 \|R_n^{B,\mathrm{step},0}-R_n^{B,0}\|_{\mathfrak S_2}^2=O(n^{-1}),
 \qquad
 \|R_n^{B,\mathrm{step},0}-\mathcal H_{\beta,a}^{-1}\|_{\mathfrak S_2}^2
       =O(n^{-1}).
 \label{B:eq:brownian-inverse-approximation}
\end{equation}
\end{theorem}
\begin{proof}
Fix $1<\lambda<\min\{a+1,2\}$, choose $\eta_0$ as in \eqref{B:eq:stability-exponents}, and fix $0<\eta_1<\nu_a$. Work on the intersection of the events in
Lemmas~\ref{B:lem:bessel-brownian-growth} and~\ref{B:lem:weighted-brownian-oscillations}. For the fixed deterministic partition and cutoffs this intersection has probability one.
Proposition~\ref{B:prop:brownian-logarithmic-errors} gives
\[
 \int_{d_n}^{\bar z_n}z^{\lambda-1}|\epsilon_n^p(z)|^2\,\dd z
 \le C\sum_{k\in\mathscr I_n^{\mathrm{ret}}}
                z_{n,k}^{\lambda}h_{n,k}^2(M_{n,k}^2+h_{n,k})
 =O_\omega(n^{-1}).
\]
Indeed, the random sum is $O_\omega(n^{-1})$ by that lemma, and the deterministic sum is bounded by $\max_k h_{n,k}\sum_k z_{n,k}^{\lambda}h_{n,k}^2\le C/(nr_n)$ by
Lemma~\ref{B:lem:brownian-grid-sums}. Propositions~\ref{B:prop:brownian-sampled-kernel} and \ref{B:prop:brownian-integration-limit} therefore imply
\[
 \begin{split}
 \|R_n^{B,\mathrm{sample}}-R_n^{B,0}\|_{\mathfrak S_2}^2
 &\le C_\omega\int_{d_n}^{\bar z_n}
                   z^{\lambda-1}|\epsilon_n^p(z)|^2\,\dd z
   =O_\omega(n^{-1}),\\
 \|R_n^{B,\mathrm{step},0}-R_n^{B,\mathrm{sample}}\|_{\mathfrak S_2}^2
 &\le C_\omega\sum_{k\in\mathscr I_n^{\mathrm{ret}}}
                        h_{n,k}^3z_{n,k}^{2-8\eta_0}
 \le \frac{C_\omega}{nr_n}.
 \end{split}
\]
For the last inequality use $8\eta_0<1$, $\max_k h_{n,k}\le C/r_n$, and \eqref{B:eq:brownian-power-sums} with $\gamma=8\eta_0$. All kernels are extended by zero to $L^2(0,1)$, so
\[
 \begin{split}
 \|R_n^{B,\mathrm{step},0}-R_n^{B,0}\|_{\mathfrak S_2}^2
 &\le 2\|R_n^{B,\mathrm{step},0}-R_n^{B,\mathrm{sample}}\|_{\mathfrak S_2}^2
       +2\|R_n^{B,\mathrm{sample}}-R_n^{B,0}\|_{\mathfrak S_2}^2\\
 &=O_\omega(n^{-1}).
 \end{split}
\]
Proposition~\ref{B:prop:brownian-factor-truncation} gives
\[
 \|R_n^{B,0}-R^B\|_{\mathfrak S_2}^2
 \le C_{\omega,\eta_1}
       [d_n^{1+\nu_a-\eta_1}+(1-\bar z_n)^3+(1-\bar z_n)^2].
\]
Since $d_n,1-\bar z_n\le Cr_n/n$, equation~\eqref{B:eq:truncation-orders} yields
\[
 n\|R_n^{B,0}-R^B\|_{\mathfrak S_2}^2
 \le C_{\omega,\eta_1}
 \left[r_n^{1+\nu_a-\eta_1}n^{-(\nu_a-\eta_1)}
               +r_n^3n^{-2}+r_n^2n^{-1}\right]\longrightarrow0.
\]
The inequality
\[
 \|R_n^{B,\mathrm{step},0}-R^B\|_{\mathfrak S_2}^2
 \le2\|R_n^{B,\mathrm{step},0}-R_n^{B,0}\|_{\mathfrak S_2}^2
       +2\|R_n^{B,0}-R^B\|_{\mathfrak S_2}^2
\]
and $R^B=\mathcal H_{\beta,a}^{-1}$ prove the second convergence estimate. For a fixed countable family of deterministic partitions and cutoffs, intersect the corresponding events
$\Omega_j$: $\mathbb P((\bigcap_j\Omega_j)^c)\le\sum_j\mathbb P(\Omega_j^c)=0$. This includes both partitions used for the finite ensembles.
\end{proof}

The restrictions $\lambda<a+1$ and $\lambda<2$ are the integrability conditions in the deterministic kernel estimates. The weighted oscillation sum uses $\lambda>1$, and removal of
the interval near zero uses $\nu_a>0$. These conditions explain the range $a>0$ of the upper approximation theorem. The intermediate masses $Q_{n,k}^B$ and the resulting step
kernel may depend on the complete Brownian path inside each cell. The observation condition in the lower bound is checked for the actual finite inverse, as a function of the
coefficient row, in Proposition~\ref{B:result:common-operator-criterion}.

\subsection{Sampling times and estimates for the boundary solutions}
\label{sec:sampling-times-and-boundary-estimates-for-the-circular-ensemble}
\label{S:sec:sampling-times-and-estimates-for-the-boundary-solutions}

The two boundary columns have opposite powers near the singular boundary. We first estimate their growth and decay to establish the existence of the singular inverse. We then
locate the samples and estimate the resulting discrete weights. The required weighted bounds include the local variance of the sampling error.

Fix $\beta>0$. Constants denoted by $C_\beta$ are deterministic and may change from one occurrence to the next. A subscript $\omega$ allows dependence on the Brownian path. We use
the reflected coordinate $r=1-t$ and write
\[
 r_{n,k}=r_k=\frac{k}{n},\qquad t_{n,k}=1-\frac{k}{n},\qquad
 \theta_\beta(r)=\frac4\beta\log\frac1r.
\]
In this coordinate $I_{n,k}=[k/n,(k+1)/n)$; its reflection differs from the interval in the original coordinate only at its boundary points. All kernel identities are understood up
to sets of Lebesgue measure zero. The value assigned to $I_{n,k}$ is sampled at $\vartheta_{n,k+1}$, where
\[
 \vartheta_{n,k}=\sum_{j=k}^{n-1}\frac4{\beta j+1},\qquad 1\le k\le n.
\]
Thus $\vartheta_{n,n}=0$, whereas the Brownian time at $r=k/n$ is $\theta_\beta(k/n)$. These times generally differ.

\subsubsection{Boundary solutions from the Brownian path}
\label{S:sec:s4-brownian-boundary-solutions}

To apply the inverse criterion of Section~\ref{S:sec:green-kernels-for-canonical-systems}, we need bounds on both boundary solutions, including the one determined by the random
terminal point. We obtain these lengths directly from the upper half-plane Brownian motion defining $\mathrm{Sine}_\beta$. In Brownian time its coordinates satisfy
\[
 \dd x_\theta=y_\theta\,\dd W_\theta^{(1)},\qquad
 \dd y_\theta=y_\theta\,\dd W_\theta^{(2)},\qquad
 (x_0,y_0)=(0,1),
\]
where $W^{(1)}$ and $W^{(2)}$ are independent standard Brownian motions. It\^o's formula applied to $\log y$ gives
\[
 y_\theta=\exp(W_\theta^{(2)}-\theta/2),\qquad
 x_\theta=\int_0^\theta y_v\,\dd W_v^{(1)}.
\]
The limit $q_\infty=\lim_{\theta\to\infty}x_\theta$ exists almost surely, as is verified in the next proof. The boundary vectors and columns remain those of
Sections~\ref{S:sec:introduction} and \ref{S:sec:green-kernels-for-canonical-systems}:
\[
 u_0=(-1,0)^{\mathsf T},\qquad u_1=(q_\infty,1)^{\mathsf T},\qquad
 a_\omega(r)=\mathcal A(Z_{\theta_\beta(r)})u_0,\quad
 c_\omega(r)=\mathcal A(Z_{\theta_\beta(r)})u_1.
\]
In particular $[u_0,u_1]=1$.

\begin{lemma}\label{S:lem:5.0}
Almost surely, for every $\delta>0$ there is $C_{\delta,\omega}<\infty$ such that, for all $\theta\ge0$,
\[
 C_{\delta,\omega}^{-1}\me^{-(1/2+\delta)\theta} \le y_\theta\le C_{\delta,\omega}\me^{-(1/2-\delta)\theta},\qquad |q_\infty-x_\theta|\le C_{\delta,\omega}\me^{-(1/2-\delta)\theta}.
\]
\end{lemma}

\begin{proof}
Use the original independent driving Brownian motions $W^1,W^2$, with $Z_0=\ii$. Solving the vertical stochastic equation gives
\[
 y_t=\exp(W_t^2-t/2),\qquad x_t=\int_0^t y_s\,\mathrm dW_s^1.
\]
For each rational $d>0$, the reflection principle gives
\[
 \mathbb P\left(\sup_{0\le s\le m+1}|W_s^2|>dm\right)
       \le4\exp\left(-\frac{d^2m^2}{2(m+1)}\right),\qquad m\ge1.
\]
The sum is finite. Taking a countable intersection over positive rational $d$ proves $W_t^2/t\to0$ almost surely. On this event, for each $d>0$,
\[
 C_d:=\exp\left(\sup_{t\ge0}(|W_t^2|-dt)_+\right)<\infty,
 \qquad
 C_d^{-1}e^{-(1/2+d)t}\le y_t\le C_de^{-(1/2-d)t}.
\]
The supremum may be taken over rational times, so $C_d$ is measurable.

Fix a rational $0<d<1/4$ and condition on the entire path of $W^2$. The integrand in the expression for $x_t$ is then deterministic, while $W^1$ remains Brownian. Its total
conditional quadratic variation is finite, and its tail obeys
\[
 V_m:=\int_m^\infty y_s^2\,\mathrm ds
 \le\frac{C_d^2}{1-2d}e^{-(1-2d)m},\qquad m\ge0.
\]
The conditional Gaussian process $x$ has the law of a Brownian motion evaluated at $\int_0^t y_s^2\,\mathrm ds$. Brownian continuity at the finite terminal time proves that $x_t$
has a finite limit $q_\infty$, conditionally almost surely and therefore almost surely under the original law.

Put $D_m=\sup_{t\ge m}|q_\infty-x_t|$. If $M_t=x_t-x_m$ for $t\ge m$, then $D_m\le2\sup_{t\ge m}|M_t|$. The conditional time change and reflection consequently imply
\[
 \mathbb P(D_m>u\mid W^2)\le4\exp\left(-\frac{u^2}{8V_m}\right).
\]
On $\{C_d\le M\}$, set $u_m=e^{-(1/2-2d)m}$. The tower property gives
\[
 \mathbb P(C_d\le M,\ D_m>u_m)
 \le4\exp\left(-\frac{1-2d}{8M^2}e^{2dm}\right).
\]
This is summable in $m$. Borel--Cantelli, followed by the countable intersections over integers $M$ and rational $d\in(0,1/4)$, gives an event on which $D_m\le u_m$ for every
sufficiently large integer $m$, for each such $d$. Since $C_d<\infty$, the localization covers almost every path.

For $m\le t<m+1$, the bound $|q_\infty-x_t|\le D_m$ and the preceding estimate imply
\[
 |q_\infty-x_t|\le C_{d,\omega}e^{-(1/2-2d)t}.
\]
Continuity on the remaining compact interval absorbs the finitely many earlier times. Given any $\delta>0$, choose a rational $0<d<\min\{\delta/2,1/4\}$. The vertical bounds with
exponent $d$ and the horizontal bound with exponent $2d$ imply all the stated bounds with exponent $\delta$, on the same event. This proves directly the exponential envelopes used
in \cite{valko2017sine}, with the diffusion coefficients used here.
\end{proof}

The explicit columns from Section~\ref{S:sec:green-kernels-for-canonical-systems} give, with $\theta=\theta_\beta(r)$,
\[
 a_\omega(r)=\begin{pmatrix}-y_\theta^{-1/2}\\0\end{pmatrix},\qquad
 c_\omega(r)=\begin{pmatrix}(q_\infty-x_\theta)y_\theta^{-1/2}\\
                                  y_\theta^{1/2}\end{pmatrix},
\]
and hence
\[
 p_\omega(r)=y_\theta^{-1},\qquad
 q_\omega(r)=y_\theta+\frac{(q_\infty-x_\theta)^2}{y_\theta}.
\]
Write $\gamma=2/\beta$ and $a_\beta=\min(1,\gamma)$. The primitives $A_\omega,C_\omega$ retain their definitions from the introduction: $A_\omega(r)=\int_0^r q_\omega(u)\dd u$ and
$C_\omega(r)=\int_r^1p_\omega(u)\dd u$.

\begin{theorem}\label{S:thm:5.1}
Almost surely, for every $\varepsilon>0$,
\begin{equation}\label{S:eq:5.2}
 C_{\varepsilon,\omega}^{-1}r^{-\gamma+\varepsilon} \le p_\omega(r)\le C_{\varepsilon,\omega}r^{-\gamma-\varepsilon},\qquad C_{\varepsilon,\omega}^{-1}r^{\gamma+\varepsilon} \le q_\omega(r)\le C_{\varepsilon,\omega}r^{\gamma-\varepsilon}, \qquad 0<r\le1.
\end{equation}
The boundary point is limit point for $0<\beta\le2$ and limit circle for $\beta>2$; the terminal solution at spectral parameter zero is square integrable in both cases.
\end{theorem}

\begin{proof}
For the bounds on $p_\omega$, use Lemma~\ref{S:lem:5.0} with $\delta=\beta\varepsilon/4$, invert its two inequalities, and substitute $\theta=(4/\beta)\log(1/r)$. The lower bound
on $q_\omega$ follows from $q_\omega\ge y_\theta$ and the same substitution. To prove its upper bound, choose $0<\delta\le\beta\varepsilon/12$ in Lemma~\ref{S:lem:5.0}. Then
\[
 \frac{(q_\infty-x_\theta)^2}{y_\theta}
 \le C_{\delta,\omega}^3
       \me^{-(1-2\delta)\theta+(1/2+\delta)\theta}
 =C_{\delta,\omega}^3\me^{-(1/2-3\delta)\theta}.
\]
At $\theta=(4/\beta)\log(1/r)$,
\[
 \me^{-(1/2-3\delta)\theta}=r^{\gamma-12\delta/\beta}
       \le r^{\gamma-\varepsilon},\qquad
 y_\theta\le C_{\delta,\omega}r^{\gamma-4\delta/\beta}
       \le C_{\delta,\omega}r^{\gamma-\varepsilon},
\]
because $0<r\le1$ and $12\delta/\beta\le\varepsilon$. Adding the bounds for the two terms defining $q_\omega$ proves its claimed upper bound.

Choose $0<\nu<\min(\gamma/2,1/4)$, and, if $\gamma<1$, also $\nu<(1-\gamma)/2$. Apply the preceding bounds with $\varepsilon=\nu$ and integrate to obtain
\[
 \int_0^1q_\omega(r)\dd r\le\frac{C_{\nu,\omega}}{1+\gamma-\nu}<\infty
\]
and
\[
 C_\omega(r)\le C_{\nu,\omega}
 \begin{cases}
 (r^{1-\gamma-\nu}-1)/(\gamma+\nu-1),&\gamma\ge1,\\
 1/(1-\gamma-\nu),&0<\gamma<1.
 \end{cases}
\]
Multiplication by $q_\omega(r)\le C_{\nu,\omega}r^{\gamma-\nu}$ shows that $q_\omega C_\omega$ is bounded by a constant times $r^{1-2\nu}$ when $\gamma\ge1$, and by
$r^{\gamma-\nu}$ when $\gamma<1$. These powers are integrable at zero. Continuity gives the local coefficient bounds away from zero. We have therefore verified both assumptions of
\eqref{S:eq:singular-integrability-assumptions}, in the coordinate $t=1-r$, before using the deterministic classification at the boundary. Multiplying the two bounds for $q_\omega$
and $C_\omega$, with $2\nu<\varepsilon$ when $\gamma\ge1$ and $\nu<\min(\varepsilon,1-\gamma)$ when $\gamma<1$, gives
\begin{equation}\label{S:eq:sine-singular-integrability-bound}
 q_\omega(r)C_\omega(r)
       \le C_{\beta,\omega,\varepsilon}r^{a_\beta-\varepsilon},
 \qquad 0<r<1,\quad 0<\varepsilon<a_\beta.
\end{equation}

For the first constant solution, substitution of $r=\me^{-\beta\theta/4}$ gives
\begin{equation}\label{S:eq:sine-zero-solution-integrability}
 \int_0^1p_\omega(r)\,\dd r
 =\frac\beta4\int_0^\infty
    \exp\left(-W_\theta^{(2)}+\frac{2-\beta}{4}\theta\right)\dd\theta.
\end{equation}
If $\beta<2$, take $\eta<(2-\beta)/4$ in $|W_\theta^{(2)}|\le K_{\eta,\omega}+\eta\theta$. The integrand is bounded below by a positive constant times an increasing exponential, so
its integral is infinite. If $\beta>2$, take $\eta<(\beta-2)/4$; the integrand is then bounded above by an integrable exponential.

At $\beta=2$ the integral is a positive constant times $\int_0^\infty e^{-W_t^{(2)}}\,\mathrm dt$. We verify its divergence directly, including the endpoint case in
\cite{valko2017sine}. Write $W=W^{(2)}$, and define stopping times
\[
 \tau_1=\inf\{t\ge0:W_t=-1\},\qquad
 \tau_{m+1}=\inf\{t\ge\tau_m+1:W_t=-(m+1)\}.
\]
All these times are finite almost surely. To see the required hitting property, for Brownian motion started at $x>y$, reflection gives
$\mathbb P_x(\inf_{s\le T}W_s\le y)=2\mathbb P(W_T\le y-x)\to1$ as $T\to\infty$; symmetry treats $x<y$. Apply this fact after the almost surely finite time $\tau_m+1$ and use
induction. The intervals $[\tau_m,\tau_m+1]$ are disjoint apart from possible endpoints.

Let
\[
 E_m=\left\{\sup_{0\le s\le1}(W_{\tau_m+s}-W_{\tau_m})\le1\right\},
 \qquad p_0=\mathbb P\left(\sup_{0\le s\le1}W_s\le1\right)>0.
\]
The strong Markov property gives $\mathbb P(E_m\mid\mathcal F_{\tau_m})=p_0$. For $j<m$, the event $E_j$ belongs to $\mathcal F_{\tau_m}$ because $\tau_j+1\le\tau_m$. Successive
conditioning therefore gives, for every $M,r\ge1$,
\[
 \mathbb P\left(\bigcap_{m=M}^{M+r-1}E_m^c\right)=(1-p_0)^r.
\]
Letting $r\to\infty$, then taking the union over $M$, proves that infinitely many $E_m$ occur almost surely. On $E_m$, $W_t\le-m+1$ throughout $[\tau_m,\tau_m+1]$, so
\[
 \int_{\tau_m}^{\tau_m+1}e^{-W_t}\,\mathrm dt\ge e^{m-1}.
\]
The disjoint intervals and infinitely many such lower bounds prove $\int_0^\infty e^{-W_t}\,\mathrm dt=\infty$.

At spectral parameter zero every coordinate solution is constant, so write it as $v=\alpha u_0+\tau u_1$. In the unweighted representation its columns are $Xu_0=a_\omega$ and
$Xu_1=c_\omega$, whose squared norms are $p_\omega$ and $q_\omega$. We have proved $c_\omega\in L^2$ for every $\beta>0$.

If $\beta>2$, also $a_\omega\in L^2$, and the triangle inequality gives $Xv\in L^2$ for every $(\alpha,\tau)\in\mathbb C^2$. If $\beta\le2$ and $Xv\in L^2$ with $\alpha\ne0$, then
\[
 a_\omega=\alpha^{-1}(Xv-\tau c_\omega)\in L^2,
\]
contradicting the divergence of $\int_0^1p_\omega$. Hence $\alpha=0$ in this case. Thus
\[
 \mathcal C_X=
 \begin{cases}
  \mathbb C u_1,&0<\beta\le2,\\
  \mathbb C^2,&\beta>2.
 \end{cases}
\]
The deterministic characterization of $\mathcal C_X$ and Proposition~\ref{S:prop:singular-weyl-alternative} give the limit point and limit circle alternatives, respectively. In
both cases the terminal solution $u_1$ is square integrable. No cancellation with $u_1$ can make a nonzero $u_0$ component square integrable in the first case.
\end{proof}

\subsubsection{The singular operator and its inverse}
\label{S:sec:s4-singular-operator}

The preceding bounds on the boundary solutions allow us to apply the deterministic domain criterion. The resulting distinction at $\beta=2$ concerns whether a second boundary
condition is needed; it does not prevent the Green inverse from existing.

Fix a path on the event in the preceding proof. The continuous positive height makes $X(t)=\mathcal A(Z_{\theta_\beta(1-t)})$ and $X(t)^{-1}$ bounded on every compact subinterval
of $[0,1)$. Put $R=X^{\mathsf T}X$ and use the maximal domain \eqref{S:eq:singular-maximal-domain}. Consider $2R^{-1}J\partial_t$ on
\begin{equation}\label{S:eq:sine-explicit-singular-domain}
 \mathcal D(\mathrm{Sine}_\beta)=
 \begin{cases}
 \{f\in\mathscr D_{\max}(X):[f(0),u_0]=0\},&0<\beta\le2,\\[1mm]
 \{f\in\mathscr D_{\max}(X):[f(0),u_0]=0,\ [f(1),u_1]=0\},&\beta>2.
 \end{cases}
\end{equation}
When $\beta>2$, the trace $f(1)$ exists by the trace theorem for the maximal domain in Section~\ref{S:sec:green-kernels-for-canonical-systems}, since both boundary directions are
square integrable. When $\beta\le2$, the operator domain is specified by the condition at zero and the integrability condition at one.

\begin{theorem}\label{S:thm:sine-selfadjoint-domain}
Almost surely, $\mathrm{Sine}_\beta$ on the domain \eqref{S:eq:sine-explicit-singular-domain} is self-adjoint and $0\in\rho(\mathrm{Sine}_\beta)$. Its conjugated inverse is the
Hilbert--Schmidt operator with kernel \eqref{S:eq:3.56}.
\end{theorem}

\begin{proof}
Proposition~\ref{S:prop:operator-normalization} identifies the Hamiltonian as $R/2$ and the canonical boundary lines as $\infty,q_\infty$, consistently with \cite{valko2017sine}.
In the coordinate $t=1-r$, Theorem~\ref{S:thm:5.1} and \eqref{S:eq:sine-singular-integrability-bound} give
\[
 Q_X=\int_0^1q_\omega(r)\,\dd r<\infty,\qquad
 I_X=\int_0^1q_\omega(r)C_\omega(r)\,\dd r
 \le C_{\beta,\omega,\varepsilon}\int_0^1r^{a_\beta-\varepsilon}\,\dd r<\infty,
 \quad 0<\varepsilon<a_\beta.
\]
The same theorem gives $\int_0^1p_\omega=\infty$ exactly when $\beta\le2$. Thus the two cases of \eqref{S:eq:singular-selfadjoint-domain} are precisely
\eqref{S:eq:sine-explicit-singular-domain}, and Theorem~\ref{S:thm:singular-inverse-domain} applies. It gives self-adjointness, $0\in\rho(\mathrm{Sine}_\beta)$, and the inverse
$\mathsf B_X$; conjugation by $\mathcal U_X$ gives the kernel \eqref{S:eq:3.56}. Proposition~\ref{S:prop:singular-minimal-domain} identifies the minimal domain.

Use the notation $\mathsf R_\beta^{\mathrm{sine}}=\mathtt r\mathrm{Sine}_\beta$ for this conjugated inverse. For $0<T<1$, the superscript $T$ denotes the restriction of both kernel
variables to $[0,T]$, extended by zero on the common unweighted space. Theorem~\ref{S:thm:singular-regular-truncations} identifies it with the regular inverse having the same
canonical terminal vector $u_1$. Its exact norm and tail norm are
\begin{equation}\label{S:eq:sine-singular-hs-identities}
 \begin{aligned}
 \|\mathsf R_\beta^{\mathrm{sine}}\|_{\mathfrak S_2}^2
   &=\frac12\int_0^1q_\omega(r)C_\omega(r)\,\dd r,\\
 \|\mathsf R_\beta^{\mathrm{sine}}-
                 \mathsf R_\beta^{\mathrm{sine},1-\delta}\|_{\mathfrak S_2}^2
   &=\frac12\int_0^\delta q_\omega(r)C_\omega(r)\,\dd r
     \le C_{\beta,\omega,\varepsilon}\delta^{1+a_\beta-\varepsilon}.
 \end{aligned}
\end{equation}
Here $0<\varepsilon<a_\beta$ and $0<\delta<1/2$. The two equalities are \eqref{S:eq:singular-hs-norm} and \eqref{S:eq:singular-truncation-hs-tail} after $r=1-t$. The last
inequality follows by integrating \eqref{S:eq:sine-singular-integrability-bound}:
$\int_0^\delta r^{a_\beta-\varepsilon}\dd r =\delta^{1+a_\beta-\varepsilon}/(1+a_\beta-\varepsilon)$. Thus every sequence $T\uparrow1$ gives Hilbert--Schmidt convergence;
\eqref{S:eq:singular-embedded-resolvents} gives the corresponding embedded resolvent convergence on compact subsets of the resolvent set.

\end{proof}

All the preceding estimates hold on a single event for this fixed $\beta$. Take the countable intersection of the events of Lemma~\ref{S:lem:5.0} for positive rational losses, of
path continuity, and, at $\beta=2$, of the recurrence event used in Theorem~\ref{S:thm:5.1}. A smaller rational loss gives every positive loss. The remaining arguments are
deterministic on this event.

The conjugated unbounded operator has domain
\[
 \mathcal D(\mathcal U_X\mathrm{Sine}_\beta\mathcal U_X^{-1})
 =\{g\in\mathcal H:X^{-1}g\in\mathcal D(\mathrm{Sine}_\beta)\}.
\]
The domain condition requires $X^{-1}g$ to be locally absolutely continuous. In the limit-point case, Proposition~\ref{S:prop:singular-minimal-domain} gives
$\lim_{t\uparrow1}[f(t),u_1]=0$ for every function in the maximal domain.

The jointly measurable kernel satisfies the finite norm criterion just verified. Subsection~\ref{S:sec:s3-measurability} therefore proves strong measurability of
$\mathsf R_\beta^{\mathrm{sine}}$ in $\mathfrak S_2$. For self-adjoint $K$, $\|(I-\ii K)^{-1}\|\le1$, and the resolvent identity gives
\[
 \|(I-\ii K)^{-1}K-(I-\ii L)^{-1}L\|_{\mathfrak S_2}
          \le\|K-L\|_{\mathfrak S_2}
\]
for self-adjoint Hilbert--Schmidt $K,L$. Thus the conjugated Sine operator has a strongly measurable resolvent at $\ii$.

We shall write $\mathsf R_{\beta,n}^{\mathrm{br}}=\mathtt r\mathrm{Br}_{\beta,n}$ and $\mathsf R_{\beta,n}^{\mathrm{circ}}=\mathtt r\mathrm{Circ}_{\beta,n}$. The regular operators
and their inverses are defined in Subsection~\ref{S:sec:normalization-of-the-operator-and-finite-walk}. By Lemma~\ref{S:lem:regular-interval-restriction}, restriction to a selected
interval gives the corresponding compressed kernel. The countable intersection over $n$ preserves probability one. The coupling and its common terminal point are constructed in
Section~\ref{S:sec:coupling-the-finite-transfer-increments-to-brownian-increments}.

\subsubsection{Deterministic sampling times}
\label{S:sec:sampling-times-and-weights-of-the-green-kernel}

An error $O(k^{-1})$ in Brownian time is as large as the Brownian duration of the $k$th operator interval, which is itself of order $k^{-1}$. We therefore expand the relative
position of the sampling time within that interval of Brownian time. The next expansion retains the term that survives after division by the interval length.

Set
\[
 b_\beta=\frac12-\frac1\beta,\qquad
 \mathfrak d_\beta=\frac1{2\beta^2}-\frac1{2\beta}+\frac1{12},\qquad
 K_\beta=\frac14+|b_\beta|+3|\mathfrak d_\beta|+\beta^{-3}.
\]
Define $e_{n,k}^{\mathrm{time}}$ by the identity
\begin{equation}\label{S:eq:sampling-times-expansion}
 \frac\beta4\vartheta_{n,k}
 =\log\frac nk+b_\beta\left(\frac1k-\frac1n\right)
  +\mathfrak d_\beta\left(\frac1{k^2}-\frac1{n^2}\right)
  +e_{n,k}^{\mathrm{time}}.
\end{equation}

\begin{lemma}\label{S:lem:4.1}
For $1\le k\le n$,
\begin{equation}\label{S:eq:sampling-times-remainder}
 |e_{n,k}^{\mathrm{time}}|
 \le K_\beta\sum_{j=k}^{n-1}j^{-4}
 \le\frac{16K_\beta}{3}(k^{-3}-n^{-3}).
\end{equation}
\end{lemma}

\begin{proof}
Set $c=\beta^{-1}$, $b=b_\beta$, $d=\mathfrak d_\beta$, and
\[
 F(x):=\log x-\frac b x-\frac d{x^2},  x\ge1.
\]
For every integer $j\ge1$, the alternating series for the logarithm gives
\[
 \left|\log(1+j^{-1})-j^{-1}+\tfrac12j^{-2}-\tfrac13j^{-3}\right|
 \le\tfrac14j^{-4}.
\]
The other two differences satisfy
\[
 \begin{aligned}
 j^{-1}-(j+1)^{-1}&=j^{-2}-j^{-3}+r_{1,j},& |r_{1,j}|&\le j^{-4},\\
 j^{-2}-(j+1)^{-2}&=2j^{-3}+r_{2,j},& |r_{2,j}|&\le3j^{-4}.
 \end{aligned}
\]
For the first line, the remainder is exactly $1/(j^3(j+1))$. For the second, apply Taylor's formula to $(1+x)^{-2}$ on $0\le x\le1$, whose second derivative is at most $6$.
Consequently,
\[
 F(j+1)-F(j)
 =j^{-1}+(-\tfrac12+b)j^{-2}+(\tfrac13-b+2d)j^{-3}+r_j,
   |r_j|\le(\tfrac14+|b|+3|d|)j^{-4}.
\]
The definitions give $-\tfrac12+b=-c$ and $\tfrac13-b+2d=c^2$. The identity
\[
 \frac1{j+c}=\frac1j-\frac c{j^2}+\frac{c^2}{j^3}
               -\frac{c^3}{j^3(j+c)}
\]
therefore implies
\[
 \left|\frac1{j+c}-\bigl(F(j+1)-F(j)\bigr)\right|
 \le K_\beta j^{-4}.
\]
Sum this estimate from $j=k$ to $n-1$. The $F$-terms telescope, and $\sum_{j=k}^{n-1}(j+c)^{-1}=(\beta/4)\vartheta_{n,k}$, proving \eqref{S:eq:sampling-times-expansion} and the
first remainder bound. For $x\in[j,j+1]$, $x\le2j$, so
\[
 \sum_{j=k}^{n-1}j^{-4}\le16\int_k^n x^{-4}\,\dd x
 =\frac{16}{3}(k^{-3}-n^{-3}).
\]
Summation gives \eqref{S:eq:sampling-times-remainder}; both sides are zero when $k=n$.
\end{proof}

The linear and quadratic corrections in \eqref{S:eq:sampling-times-expansion} are bounded by $C_\beta/k$ and $C_\beta/k^2$, respectively. Hence
$\vartheta_{n,k}=(4/\beta)\log(n/k)+O_\beta(k^{-1})$. Replacing $k$ by $k+1$ and retaining the linear correction gives
\begin{equation}\label{S:eq:subinterval-sampling-time-expansion}
 \vartheta_{n,k+1}
 =\frac4\beta\log\frac n{k+1}
  +\frac4\beta b_\beta\left(\frac1{k+1}-\frac1n\right)
  +O_\beta(k^{-2}),\qquad 1\le k<n.
\end{equation}
The constants are uniform in $n,k$.

For $0\le k<n$, put
\[
 \bar r_{n,k}=\frac{k+1}{n},\qquad
 \rho_{n,k}=\exp\left(-\frac\beta4\vartheta_{n,k+1}\right),\qquad
 C_\beta^0=1+\beta^{-1}.
\]
Thus $\rho_{n,k}$ is the position in the $r$ coordinate of the sample assigned to $I_{n,k}$.

\begin{lemma}\label{S:lem:sampling-position}
Uniformly for $0\le k<n$,
\begin{equation}\label{S:eq:sampling-position-expansion}
 \left|\rho_{n,k}-\bar r_{n,k}
       +\frac{b_\beta}{n}(1-\bar r_{n,k})\right|
 \le\frac{C_\beta(1-\bar r_{n,k})}{n(k+1)},
\end{equation}
\begin{equation}\label{S:eq:sampling-position-mesh-distance}
 |n\rho_{n,k}-(k+1)|\le C_\beta^0.
\end{equation}
\end{lemma}

\begin{proof}
Write $l=k+1$ and $D=l^{-1}-n^{-1}\ge0$. Lemma~\ref{S:lem:4.1} gives
\[
 \rho_{n,k}=\frac ln \me^{-v_{n,l}},
 v_{n,l}=b_\beta D+\mathfrak d_\beta(l^{-2}-n^{-2})+e_{n,l}^{\mathrm{time}}.
\]
Since $l^{-2}-n^{-2}\le2D/l$ and $l^{-3}-n^{-3}\le3D/l^2$,
\[
 |v_{n,l}-b_\beta D|\le C_\beta D/l,  |v_{n,l}|\le C_\beta D\le C_\beta/l.
\]
Taylor's formula $|\me^{-v}-1+v|\le\tfrac12\me^{|v|}v^2$ and $D^2\le D/l$ imply $|\me^{-v_{n,l}}-1+b_\beta D|\le C_\beta D/l$. Multiplication by $l/n$ proves
\eqref{S:eq:sampling-position-expansion}, since $(l/n)D=(1-l/n)/n$. The uniform ratio bounds follow from $|v_{n,l}|\le C_\beta$. For the bound on $|n\rho_{n,k}-(k+1)|$, set
$c=\beta^{-1}$ and use monotonicity of $(x+c)^{-1}$:
\[
 \int_l^n\frac{\dd x}{x+c}
 \le\sum_{j=l}^{n-1}\frac1{j+c}
 \le\int_{l-1}^{n-1}\frac{\dd x}{x+c}.
\]
Exponentiating yields
\[
 \frac{n(l-1+c)}{n-1+c}\le n\rho_{n,k}\le\frac{n(l+c)}{n+c}.
\]
The right side is at most $l+c$. The left side minus $l$ is $(c-1)(n-l)/(n-1+c)\ge-1$, since $l\ge1$. Thus $|n\rho_{n,k}-l|\le1+c$, proving
\eqref{S:eq:sampling-position-mesh-distance} with the stated constant. At $l=n$ the defining sum is empty, so $\rho_{n,n-1}=1$.
\end{proof}

The proof also gives $\rho_{n,k}/\bar r_{n,k}=\exp(O_\beta((k+1)^{-1}))$. In particular this ratio and its reciprocal are uniformly bounded, $\rho_{n,n-1}=1$, and the error in
\eqref{S:eq:sampling-position-expansion} is $O_{\beta,\delta}(n^{-2})$ when $k/n\ge\delta>0$.

For $n\ge2$, $1\le k<n$, and $0\le s\le1$, put
\begin{equation}\label{S:eq:sampling-offset-definitions}
 \begin{aligned}
 r_k&:=k/n,&
 \theta_{n,k}(s)&:=\theta_\beta((k+s)/n)=\frac4\beta\log\frac n{k+s},\\
 \theta_{k,+}&:=\theta_\beta((k+1)/n),&
 \tau_{n,k}&:=\theta_\beta(k/n)-\theta_{k,+}
              =\frac4\beta\log(1+k^{-1}),\\
 \widehat\sigma_{n,k}&:=\frac{\vartheta_{n,k+1}-\theta_{k,+}}{\tau_{n,k}},&
 \ell_{n,k}(s)&:=|\theta_{n,k}(s)-\vartheta_{n,k+1}|.
 \end{aligned}
\end{equation}
The denominator is positive. Let $\kappa_{n,k}:=\int_0^1|u-\widehat\sigma_{n,k}|\,\dd u$.

Define $\widehat\sigma_\beta(r)=b_\beta(1-r)$ and $\kappa_\beta(r)=\int_0^1|u-\widehat\sigma_\beta(r)|\,\dd u$.

\begin{proposition}\label{S:prop:sampling-offset}
Uniformly for $1\le k<n$,
\begin{equation}\label{S:eq:sampling-normalized-offset}
 \widehat\sigma_{n,k}=b_\beta(1-r_k)+O_\beta(k^{-1}).
\end{equation}
Moreover,
\begin{equation}\label{S:eq:sampling-variance-two-errors}
 \left|\tau_{n,k}\kappa_{n,k}-\frac4{\beta k}\kappa_\beta(r_k)\right|
 +\left|\int_0^1\ell_{n,k}(s)\,\dd s-\frac4{\beta k}\kappa_\beta(r_k)\right|
 \le C_\beta k^{-2}.
\end{equation}
\end{proposition}

\begin{proof}
The scaled numerator is, by \eqref{S:eq:sampling-times-expansion},
\[
 \frac\beta4(\vartheta_{n,k+1}-\theta_{k,+})
 =b_\beta\left(\frac1{k+1}-\frac1n\right)
  +\mathfrak d_\beta\left(\frac1{(k+1)^2}-\frac1{n^2}\right)
  +O_\beta(k^{-3}).
\]
Put $x=k^{-1}$ and $r=k/n$, so $n^{-1}=rx$. Expanding $(1+x)^{-1}$ and $(1+x)^{-2}$ on $0\le x\le1$ with uniformly bounded Taylor remainders gives the numerator
\[
 b_\beta(1-r)x+\bigl(-b_\beta+\mathfrak d_\beta(1-r^2)\bigr)x^2
       +O_\beta(x^3).
\]
The scaled denominator is $\log(1+x)=x(1-x/2+O(x^2))$, and $\log(1+x)\ge x/2$ on this interval. Thus its reciprocal, after removing $x$, is $1+x/2+O(x^2)$ with a uniform remainder.
Multiplication gives
\begin{equation}\label{S:eq:sampling-normalized-offset-refined}
 \widehat\sigma_{n,k}
 =b_\beta(1-r_k)
  +\frac{-\tfrac12b_\beta(1+r_k)+\mathfrak d_\beta(1-r_k^2)}{k}
  +O_\beta(k^{-2}).
\end{equation}
The coefficients are bounded on $0\le r\le1$, so the remainder in \eqref{S:eq:sampling-normalized-offset} is uniform in
$k<n$.

The map $K(v):=\int_0^1|u-v|\,\dd u$ is globally Lipschitz with constant one, because $||u-v|-|u-w||\le|v-w|$. Direct integration gives
\[
 K(v)=
 \begin{cases}
 \tfrac12-v,&v\le0,\\
 \tfrac12-v+v^2,&0\le v\le1,\\
 v-\tfrac12,&v\ge1.
 \end{cases}
\]
For $0\le v\le1$, $K(v)=(v-1/2)^2+1/4$; outside this interval $K(v)\ge1/2$. For $\beta\le2$, $b_\beta\le0$; for $\beta>2$, $0<b_\beta<1/2$. Substituting $v=b_\beta(1-r)$ gives
\begin{equation}\label{S:eq:sampling-factor-branches}
 \kappa_\beta(r)=
 \begin{cases}
 \tfrac12-b_\beta(1-r),&0<\beta\le2,\\
 \tfrac12-b_\beta(1-r)+b_\beta^2(1-r)^2,&\beta>2,
 \end{cases}
 \qquad 0\le r\le1.
\end{equation}
Furthermore,
\begin{equation}\label{S:eq:sampling-factor-approximation}
 |\kappa_{n,k}-\kappa_\beta(r_k)|\le C_\beta/k,
   \tfrac14\le\kappa_{n,k},\kappa_\beta(r_k)\le C_\beta.
\end{equation}
The first estimate follows from the Lipschitz bound and \eqref{S:eq:sampling-normalized-offset}; the second follows from $K(v)\ge1/4$ and $K(v)\le1/2+|v|$.

To compare $\ell_{n,k}(s)$ with $\tau_{n,k}|1-s-\widehat\sigma_{n,k}|$, apply the remainder formula for linear interpolation to $f(s)=\log(k+s)$:
\[
 0\le\log(k+s)-(1-s)\log k-s\log(k+1)
 \le\frac{s(1-s)}{2k^2}.
\]
The interpolation remainder is bounded by $s(1-s)\sup|f''|/2$, and $f''(s)=-(k+s)^{-2}$. Therefore
\[
 \theta_{n,k}(s)-\vartheta_{n,k+1}
 =\tau_{n,k}(1-s-\widehat\sigma_{n,k})+\epsilon_{n,k}(s),
 |\epsilon_{n,k}(s)|\le\frac{2s(1-s)}{\beta k^2}.
\]
Since $\big||x+y|-|x|\big|\le|y|$, integration and the substitution $u=1-s$ give
\begin{equation}\label{S:eq:time-interpolation-error}
 \left|\int_0^1\ell_{n,k}(s)\,\dd s-\tau_{n,k}\kappa_{n,k}\right|
 \le\frac1{3\beta k^2}.
\end{equation}
Here the constant is $(2/\beta)\int_0^1s(1-s)\dd s=1/(3\beta)$. Finally,
\[
 \left|\tau_{n,k}-\frac4{\beta k}\right|\le\frac2{\beta k^2},
   \tau_{n,k}\le\frac4{\beta k}.
\]
Combine these inequalities with \eqref{S:eq:sampling-factor-approximation}, first for $\tau_{n,k}\kappa_{n,k}$ and then for the integral of $\ell_{n,k}$, to prove
\eqref{S:eq:sampling-variance-two-errors}. At $k=n-1$, both $\vartheta_{n,n}$ and $\theta_\beta(1)$ are zero, so $\widehat\sigma_{n,n-1}=0$.
\end{proof}

At $\beta=2$, $\kappa_2\equiv1/2$, while \eqref{S:eq:sampling-normalized-offset-refined} gives the finite displacement
\begin{equation}\label{S:eq:sampling-critical-finite-offset}
 \widehat\sigma_{n,k}=-\frac{1-r_k^2}{24k}+O(k^{-2})
 , (\beta=2).
\end{equation}
In fact it is strictly negative when $k+1<n$: the strict convexity of $x\mapsto x^{-1}$ gives $(j+1/2)^{-1}<\int_j^{j+1}x^{-1}\,\dd x$, and summing from $j=k+1$ to $n-1$ gives
$\vartheta_{n,k+1}<\theta_\beta((k+1)/n)$. Equality occurs when $k+1=n$.

For later weighted sums, \eqref{S:eq:time-interpolation-error} and \eqref{S:eq:sampling-variance-two-errors} have the following direct consequence. For any finite set
$\mathcal K\subset\{m,\ldots,n-1\}$, $m\ge1$, and any scalar weights
$w_k$,
\begin{equation}\label{S:eq:sampling-weighted-error}
 \begin{aligned}
 &\left|\sum_{k\in\mathcal K}w_k
       \left(\int_0^1\ell_{n,k}(s)\,\dd s-\tau_{n,k}\kappa_{n,k}\right)\right|\\
 & +\left|\sum_{k\in\mathcal K}w_k
       \left(\tau_{n,k}\kappa_{n,k}-\frac4{\beta k}\kappa_\beta(r_k)\right)\right|
 \le C_\beta\sum_{k\in\mathcal K}\frac{|w_k|}{k^2}
 \le\frac{C_\beta}{m}\sum_{k\in\mathcal K}\frac{|w_k|}{k}.
 \end{aligned}
\end{equation}
Indeed, the triangle inequality bounds the two sums by the sum of their absolute summands. The individual errors are at most $C_\beta k^{-2}$, and $k^{-2}\le m^{-1}k^{-1}$ on
$\mathcal K$. Thus \eqref{S:eq:sampling-weighted-error} also holds after arbitrary random values have been substituted for the weights.

\subsubsection{Sampled columns and their integrals}
\label{S:sec:s4-sampled-columns}

The location of the samples transfers the continuous column bounds to the discrete coefficients. Their cumulative integrals then provide the weights for the perturbation estimate.

Since $\theta_\beta(\rho_{n,k})=\vartheta_{n,k+1}$, the sampled columns satisfy the exact identities
\[
 a_{n,k}^{\mathrm{br}}=a_\omega(\rho_{n,k}),\qquad
 c_{n,k}^{\mathrm{br}}=c_\omega(\rho_{n,k}),\qquad 0\le k<n.
\]
We use $p_{n,k}=|a_{n,k}^{\mathrm{br}}|^2$ and $q_{n,k}=|c_{n,k}^{\mathrm{br}}|^2$.

\begin{proposition}\label{S:prop:5.2}
Almost surely, for every $\varepsilon>0$, uniformly in $n$ and $0\le k<n$,
\[
 C_{\varepsilon,\omega}^{-1}\bar r_{n,k}^{-\gamma+\varepsilon} \le p_{n,k}\le C_{\varepsilon,\omega}\bar r_{n,k}^{-\gamma-\varepsilon},\qquad C_{\varepsilon,\omega}^{-1}\bar r_{n,k}^{\gamma+\varepsilon} \le q_{n,k}\le C_{\varepsilon,\omega}\bar r_{n,k}^{\gamma-\varepsilon}.
\]
For $k\ge1$, the same bounds hold with $\bar r_{n,k}$ replaced by $r_k$.
\end{proposition}

\begin{proof}
Lemma~\ref{S:lem:sampling-position} gives a deterministic $L_\beta$ with $|\log(\rho_{n,k}/\bar r_{n,k})|\le L_\beta$. Consequently, for every real $\alpha$,
\[
 \me^{-|\alpha|L_\beta}\bar r_{n,k}^{\alpha}
 \le\rho_{n,k}^{\alpha}
 \le\me^{|\alpha|L_\beta}\bar r_{n,k}^{\alpha}.
\]
Apply Theorem~\ref{S:thm:5.1} at $\rho_{n,k}$ and use this inequality for $\alpha=-\gamma\pm\varepsilon$ and $\gamma\pm\varepsilon$. All constants are independent of $n,k$. For
$k\ge1$, $1\le\bar r_{n,k}/r_k\le2$, which proves the same bounds with $r_k$ in place of $\bar r_{n,k}$. The use of $\bar r_{n,0}=1/n$ covers $k=0$ without evaluating a singular
power at zero.
\end{proof}

To estimate the diagonal cell squares, we use the product of the two column bounds. Fix $0<\sigma<1$ and apply the two upper bounds with loss $\sigma/2$. The exponents $\gamma$
cancel, giving
\[
 p_{n,k}q_{n,k}\le C_{\beta,\omega,\sigma}
       \bar r_{n,k}^{-\sigma},\qquad 0\le k<n.
\]
For $1\le m<n$ and $d\in\{1,2\}$ it follows that
\[
 \frac1n\sum_{k=m}^{n-1}\frac{p_{n,k}q_{n,k}}{k^d}
 \le C_{\beta,\omega,\sigma}n^{-1+\sigma}
           \sum_{k=m}^\infty k^{-d-\sigma}
 \le C_{\beta,\omega,\sigma}n^{-1+\sigma}m^{1-d-\sigma}.
\]
Indeed, the decreasing summand is bounded by its first value plus $\int_m^\infty x^{-d-\sigma}\dd x$, and $m\ge1$ absorbs the first value. These estimates hold for every $\beta>0$.
The bounds for these products control the contributions from the current interval in $\mathcal E_n^+$.

The next estimates hold for every deterministic cutoff sequence with $1\le m_n<n$ and $m_n/n\to0$. They apply, in particular, to the fixed logarithmic cutoff in
\eqref{S:eq:cutoff-choice} once $n$ is sufficiently large. The sums $A_{n,k}$ and $C_{n,k}$ introduced in Section~\ref{S:sec:green-kernels-for-canonical-systems} are
\[
 A_{n,k}=\frac1n\sum_{\ell=m_n}^{k-1}q_{n,\ell},\qquad
 C_{n,k}=\frac1n\sum_{\ell=k+1}^{n-1}p_{n,\ell},\qquad m_n\le k<n.
\]
In particular $A_{n,m_n}=C_{n,n-1}=0$. For a matrix $U\in\mathfrak{sl}_2(\mathbb R)$ put
\[
 \mathcal M_{n,k}(U)
     =A_{n,k}|Ua_{n,k}^{\mathrm{br}}|^2
             +C_{n,k}|Uc_{n,k}^{\mathrm{br}}|^2,
\]
and write
\[
 \mathfrak G_{n,k}=A_{n,k}p_{n,k}+C_{n,k}q_{n,k},\qquad
 \mathfrak G_\omega(r)=A_\omega(r)p_\omega(r)+C_\omega(r)q_\omega(r).
\]
The weights are defined by the random sampled columns; the diagonal example of Section~\ref{S:sec:green-kernels-for-canonical-systems} is used only for comparison.

The combination $\mathfrak G=A p+C q$ is more useful than separate maxima of the columns. Indeed, in the first variation of $a(r)c(r')^{\mathsf T}/2$, perturbing $a(r)$ and
integrating the other variable produces $A(r)|Ua(r)|^2$, while perturbing $c(r)$ produces $C(r)|Uc(r)|^2$. Both are at most the corresponding part of
$\mathfrak G(r)\|U(r)\|_{\mathrm{op}}^2$. Brownian increments on the $k$th interval have time scale $k^{-1}$, so the sums needed later have the form $\sum_k\mathfrak G_{n,k}/k$.
The following estimates apply directly to these sums of products.

\begin{proposition}\label{S:prop:5.3}
For every sufficiently small $\varepsilon>0$, almost surely,
\begin{equation}\label{S:eq:5.20}
 C_{\varepsilon,\omega}^{-1}r^{1+\gamma+\varepsilon}
     \le A_\omega(r)\le C_{\varepsilon,\omega}r^{1+\gamma-\varepsilon},
 \qquad 0<r\le1,
\end{equation}
and
\[
 C_\omega(r)\le C_{\varepsilon,\omega}
 \begin{cases}
 r^{1-\gamma-\varepsilon},&\gamma>1,\\
 r^{-\varepsilon},&\gamma=1,\\
 1,&0<\gamma<1.
 \end{cases}
\]
The same upper estimates hold for $A_{n,k},C_{n,k}$ with $r=r_k$, uniformly over the selected indices.
\end{proposition}

\begin{proof}
We may take $0<\varepsilon<\min(\gamma,1)/2$ and, if $\gamma<1$, also $\varepsilon<(1-\gamma)/2$. Integrating the two bounds on $q_\omega$ gives
\[
 \frac{C_{\varepsilon,\omega}^{-1}}{1+\gamma+\varepsilon}
                 r^{1+\gamma+\varepsilon}
 \le A_\omega(r)\le
 \frac{C_{\varepsilon,\omega}}{1+\gamma-\varepsilon}
                 r^{1+\gamma-\varepsilon}.
\]
For $\alpha>0$,
\[
 \int_r^1u^{-\alpha}\dd u=
 \begin{cases}
 (r^{1-\alpha}-1)/(\alpha-1),&\alpha\ne1,\\
 \log(1/r),&\alpha=1.
 \end{cases}
\]
Apply it with $\alpha=\gamma+\varepsilon$. When $\gamma\ge1$, it is at most $r^{1-\gamma-\varepsilon}/(\gamma+\varepsilon-1)$. When $\gamma<1$, the choice
$\varepsilon<(1-\gamma)/2$ gives $1-\gamma-\varepsilon>(1-\gamma)/2>0$, so
\[
 \int_r^1u^{-\gamma-\varepsilon}\,\dd u
 =\frac{1-r^{1-\gamma-\varepsilon}}{1-\gamma-\varepsilon}
 \le\frac1{1-\gamma-\varepsilon}.
\]

For the discrete $A$ sum, $\gamma-\varepsilon>0$ and Proposition~\ref{S:prop:5.2} give
\[
 A_{n,k}\le\frac Cn\sum_{\ell=m_n}^{k-1}(\ell/n)^{\gamma-\varepsilon}
 \le C\sum_{\ell=m_n}^{k-1}
                  \int_{\ell/n}^{(\ell+1)/n}u^{\gamma-\varepsilon}\dd u
 \le\frac{C}{1+\gamma-\varepsilon}r_k^{1+\gamma-\varepsilon}.
\]
The inequality uses monotonicity on each indicated interval, whose length is exactly $1/n$. Since $u^{-\gamma-\varepsilon}$ is decreasing,
\[
 C_{n,k}\le\frac Cn\sum_{\ell=k+1}^{n-1}(\ell/n)^{-\gamma-\varepsilon}
 \le C\sum_{\ell=k+1}^{n-1}
                  \int_{(\ell-1)/n}^{\ell/n}u^{-\gamma-\varepsilon}\dd u
 \le C\int_{r_k}^1u^{-\gamma-\varepsilon}\dd u.
\]
The same integral calculation proves all discrete upper bounds. The inequalities include empty sums at the two extreme indices.

We also record the lower estimates needed away from a cutoff. If $\gamma>1$ and $0<r\le1/2$, then
\[
 C_\omega(r)\ge C_{\varepsilon,\omega}^{-1}
        \int_r^{2r}u^{-\gamma+\varepsilon}\dd u
 =C_{\varepsilon,\omega}^{-1}r^{1-\gamma+\varepsilon}
        \int_1^2v^{-\gamma+\varepsilon}\dd v.
\]
The last integral is a finite positive constant. There is no such positive lower power bound at $r=1$, where $C_\omega(1)=0$. For $k\ge2m_n$, the indices
$\lceil k/2\rceil\le\ell\le k-1$ belong to the $A$ sum and number $\lfloor k/2\rfloor\ge k/3$. Therefore
\[
 A_{n,k}\ge\frac{k}{3n}
       C_{\varepsilon,\omega}^{-1}(r_k/2)^{\gamma+\varepsilon}
 \ge c_{\varepsilon,\omega}r_k^{1+\gamma+\varepsilon}.
\]
Finally, if $0<\gamma<1$, the integrability established above and monotone convergence give $C_\omega(r)\uparrow\int_0^1p_\omega(u)\dd u\in(0,\infty)$. At $\gamma=1$, this
primitive tends to infinity by \eqref{S:eq:sine-zero-solution-integrability} and its proof at the critical value.
\end{proof}

\begin{proposition}\label{S:prop:5.4}
Almost surely, for every $\varepsilon>0$,
\begin{equation}\label{S:eq:5.30}
 \mathfrak G_\omega(r)\le C_{\varepsilon,\omega}r^{a_\beta-\varepsilon},
 \qquad 0<r\le1.
\end{equation}
Uniformly for $m_n\le k<n$, one also has $\mathfrak G_{n,k}\le C_{\varepsilon,\omega}r_k^{a_\beta-\varepsilon}$.
\end{proposition}

\begin{proof}
Choose $0<\nu<\min(\varepsilon/4,\gamma/4,1/4)$; when $\gamma<1$ also take $\nu<(1-\gamma)/2$. Apply Theorem~\ref{S:thm:5.1} and Proposition~\ref{S:prop:5.3} with
$\varepsilon=\nu$. Their two product bounds are
\[
 p_\omega(r)A_\omega(r)
       \le C r^{-\gamma-\nu}r^{1+\gamma-\nu}=Cr^{1-2\nu},
\]
\[
 q_\omega(r)C_\omega(r)\le C
 \begin{cases}
 r^{\gamma-\nu}r^{1-\gamma-\nu}=r^{1-2\nu},&\gamma\ge1,\\
 r^{\gamma-\nu},&0<\gamma<1.
 \end{cases}
\]
In both cases the exponents on the right are at least $a_\beta-\varepsilon$:
\[
 1-2\nu\ge a_\beta-\varepsilon,\qquad
 \gamma-\nu\ge a_\beta-\varepsilon\quad(\gamma<1).
\]
Their sum is therefore bounded by $Cr^{a_\beta-\varepsilon}$ for $0<r\le1$, which proves \eqref{S:eq:5.30}. The discrete inequalities use the same products,
$p_{n,k}A_{n,k}\le Cr_k^{1-2\nu}$ and $q_{n,k}C_{n,k}\le Cr_k^{1-2\nu}$ or $Cr_k^{\gamma-\nu}$. They also hold when $A_{n,k}=0$ or $C_{n,k}=0$. Interior lower bounds for these
quantities are proved in Proposition~\ref{S:prop:5.6} below.
\end{proof}

\begin{proposition}\label{S:prop:5.5}
Almost surely,
\[
 \int_0^1\frac{\mathfrak G_\omega(r)}r\,\dd r<\infty,\qquad
 \sup_n\sum_{k=m_n}^{n-1}\frac{\mathfrak G_{n,k}}k<\infty,\qquad
 \max_{m_n\le k<n}\frac{\mathfrak G_{n,k}}k\longrightarrow0.
\]
\end{proposition}

\begin{proof}
Choose $0<\zeta<a_\beta\le1$ and apply Proposition~\ref{S:prop:5.4} with $\varepsilon=a_\beta-\zeta$. The integral is at most $C\int_0^1r^{\zeta-1}\dd r=C/\zeta$. For every integer
$N\ge1$, monotonicity gives
\[
 \sum_{k=1}^Nk^{\zeta-1}
       \le1+\int_1^N x^{\zeta-1}\dd x
       \le1+N^\zeta/\zeta.
\]
Consequently
\[
 \sum_{k=m_n}^{n-1}\frac{\mathfrak G_{n,k}}k
       \le Cn^{-\zeta}\sum_{k=1}^nk^{\zeta-1}
       \le C(1+1/\zeta),
\]
and $\max_{k\ge m_n}\mathfrak G_{n,k}/k \le Cn^{-\zeta}m_n^{\zeta-1}\le Cn^{-\zeta}$. If $N_0$ is the starting index of these bounds, then
\[
 \sup_n\sum_{k=m_n}^{n-1}\frac{\mathfrak G_{n,k}}k
 \le\max\left\{C(1+1/\zeta),
       \max_{n<N_0}\sum_{k=m_n}^{n-1}\frac{\mathfrak G_{n,k}}k\right\}<\infty.
\]
The finite maximum is well defined because every sampled height is positive; empty sums and maxima are zero. The estimate $Cn^{-\zeta}$ also proves convergence of the maxima. For
the sum restricted to $k\le n\delta$, the bound for the sum of powers gives
\begin{equation}\label{S:eq:s4-small-index-sum}
 \sum_{\substack{m_n\le k<n\\k/n\le\delta}}
          \frac{\mathfrak G_{n,k}}k
       \le C(n^{-\zeta}+\delta^\zeta/\zeta),\qquad 0<\delta\le1.
\end{equation}
\end{proof}

Let $a_n,c_n,p_n,q_n$ be the step functions with respective values $a_{n,k}^{\mathrm{br}},c_{n,k}^{\mathrm{br}},p_{n,k},q_{n,k}$ on $I_{n,k}$, $m_n\le k<n$. Assign to $r=1$ the
value of the last interval. Put $D_\beta=2+\beta^{-1}$ and $d_n=D_\beta/n$. For a continuous scalar or vector function $f$ on a compact interval $I$, write
\[
 \omega_{f,I}(h)=\sup\{|f(u)-f(v)|:u,v\in I,\ |u-v|\le h\}.
\]
The supremum tends to zero as $h\downarrow0$.

\begin{lemma}\label{S:lem:sampled-exclusive-convergence}
Almost surely, $c_\omega,q_\omega$ extend continuously to zero, and $c_n,q_n$ converge uniformly to $c_\omega,q_\omega$, respectively, on $[m_n/n,1]$. The functions $a_n,p_n$
converge to $a_\omega,p_\omega$ uniformly on $[\delta,1]$ for every $\delta>0$. Moreover,
\[
 \max_{m_n\le k<n}|A_{n,k}-A_\omega(r_k)|\to0,\qquad
 \max_{\substack{m_n\le k<n\\r_k\ge\delta}}
              |C_{n,k}-C_\omega(r_k)|\to0.
\]
\end{lemma}

\begin{proof}
We work on the single event specified after Theorem~\ref{S:thm:sine-selfadjoint-domain}. Lemma~\ref{S:lem:sampling-position} and $|u-(k+1)/n|\le1/n$ for $u\in[k/n,(k+1)/n]$ give
\begin{equation}\label{S:eq:sampled-radius-uniform-distance}
 |\rho_{n,k}-r_k|\le d_n,\qquad
 \sup_{u\in[k/n,(k+1)/n]}|\rho_{n,k}-u|\le d_n,
 \qquad \rho_{n,n-1}=1.
\end{equation}
For $0<\eta<\gamma$, Theorem~\ref{S:thm:5.1} gives $q_\omega(r)=|c_\omega(r)|^2\le K_{\eta,\omega}r^{\gamma-\eta}$. Since $\gamma-\eta>0$, both $c_\omega(r)$ and $q_\omega(r)$ tend
to zero as $r\downarrow0$. Their continuous extensions therefore satisfy, for $f=c_\omega,q_\omega$,
\[
 \sup_{u\in[m_n/n,1]}|f_n(u)-f(u)|
       \le\omega_{f,[0,1]}(d_n)\longrightarrow0.
\]
For fixed $\delta>0$ and large $n$, one has $m_n/n<\delta$ and $d_n<\delta/2$. If $u\ge\delta$, its sample belongs to $[\delta/2,1]$. Thus, for $f=a_\omega,p_\omega$,
\begin{equation}\label{S:eq:sampled-compact-modulus}
 \sup_{u\in[\delta,1]}|f_n(u)-f(u)|
       \le\omega_{f,[\delta/2,1]}(d_n)\longrightarrow0.
\end{equation}
These inequalities hold also at the assigned interval boundary points.

Integrating the step functions on their complete intervals gives
\begin{equation}\label{S:eq:exclusive-primitives-step-identities}
 A_{n,k}=\int_{m_n/n}^{r_k}q_n(u)\dd u,\qquad
 C_{n,k}=\int_{r_k+1/n}^{1}p_n(u)\dd u.
\end{equation}
Subtracting the continuum integrals gives
\begin{equation}\label{S:eq:exclusive-primitive-error-decompositions}
 \begin{aligned}
 A_{n,k}-A_\omega(r_k)
   &=\int_{m_n/n}^{r_k}(q_n-q_\omega)(u)\dd u
                      -\int_0^{m_n/n}q_\omega(u)\dd u,\\
 C_{n,k}-C_\omega(r_k)
   &=\int_{r_k+1/n}^1(p_n-p_\omega)(u)\dd u
                      -\int_{r_k}^{r_k+1/n}p_\omega(u)\dd u.
 \end{aligned}
\end{equation}
The first integral in either line has length at most one. The global modulus for $q_\omega$ and its power bound show that
\begin{equation}\label{S:eq:exclusive-A-uniform-error}
 \max_{m_n\le k<n}|A_{n,k}-A_\omega(r_k)|
 \le\omega_{q_\omega,[0,1]}(d_n)
       +\frac{K_{\eta,\omega}}{1+\gamma-\eta}
                    (m_n/n)^{1+\gamma-\eta}\longrightarrow0.
\end{equation}
For $r_k\ge\delta$, the second line of \eqref{S:eq:exclusive-primitive-error-decompositions} is entirely in $[\delta,1]$. Its omitted interval has length $1/n$, so
\begin{equation}\label{S:eq:exclusive-C-compact-error}
 \max_{\substack{m_n\le k<n\\r_k\ge\delta}}|C_{n,k}-C_\omega(r_k)|
 \le\omega_{p_\omega,[\delta/2,1]}(d_n)
            +\frac1n\sup_{u\in[\delta,1]}p_\omega(u)\longrightarrow0.
\end{equation}
At $k=m_n$, the first error is the omitted integral of $q_\omega$. At $k=n-1$, the second is $-\int_{(n-1)/n}^1p_\omega$. The latter estimate uses $p_\omega$ only on
$[\delta/2,1]$. The preceding deterministic
bounds hold on the chosen event for every admissible cutoff sequence and every $\delta>0$. \end{proof}

Fix $0<r_-<r_+<1$. The constants in the next statement may depend on this compact interval.

\begin{proposition}\label{S:prop:5.6}
Almost surely, for all sufficiently large $n$ and all $r_k\in[r_-,r_+]$,
\[
 0<c_\omega\le p_{n,k},q_{n,k},A_{n,k},C_{n,k},\mathfrak G_{n,k}
             \le C_\omega<\infty.
\]
\end{proposition}

\begin{proof}
Both columns are continuous and nonzero on $[r_-,r_+]$, because $\mathcal A(Z_{\theta_\beta(r)})$ is invertible and $u_0,u_1\ne0$. The extreme value theorem thus gives positive
minima and finite maxima of $p_\omega,q_\omega$ there. Also, for every $r\in[r_-,r_+]$,
\[
 A_\omega(r)\ge\int_{r_-/2}^{r_-}q_\omega(u)\dd u>0,\qquad
 C_\omega(r)\ge\int_{r_+}^{(1+r_+)/2}p_\omega(u)\dd u>0.
\]
Their upper bounds are $A_\omega(r_+)$ and $C_\omega(r_-)$. Let $c_0$ be the minimum of these four positive lower bounds and $C_0$ the maximum of the four upper bounds. By
Lemma~\ref{S:lem:sampled-exclusive-convergence}, for all large $n$ each sampled quantity differs from its limit at $r_k$ by at most $\min(c_0/2,1)$. Thus the four quantities are
between $c_0/2$ and $C_0+1$. Their products show $c_0^2/2\le\mathfrak G_{n,k}\le2(C_0+1)^2$. Taking the smaller lower constant and the larger upper constant proves the statement.
\end{proof}

The preceding estimates justify the weighted Riemann sums directly. For any continuous $f\colon[0,1]\to\mathbb R$, on the same event,
\begin{equation}\label{S:eq:s4-weighted-riemann}
 \sum_{k=m_n}^{n-1}f(r_k)\frac{\mathfrak G_{n,k}}k
       \longrightarrow\int_0^1 f(r)\frac{\mathfrak G_\omega(r)}r\dd r.
\end{equation}
Fix $0<\delta<1$ and put $k_\delta=\lceil n\delta\rceil$. For large $n$, $m_n\le k_\delta<n$. On this range define
\[
 d_n(\delta)=\max_{k_\delta\le k<n}|\mathfrak G_{n,k}-\mathfrak G_\omega(r_k)|,
 \qquad g(r)=f(r)\mathfrak G_\omega(r)/r.
\]
The convergence of $p,q,A,C$ on $[\delta,1]$ gives $d_n(\delta)\to0$, and $g$ is continuous there. Since $1/k=1/(nr_k)$,
\[
 \left|\sum_{k=k_\delta}^{n-1}f(r_k)\frac{\mathfrak G_{n,k}}k
       -\frac1n\sum_{k=k_\delta}^{n-1}g(r_k)\right|
 \le\frac{\|f\|_\infty}{\delta}d_n(\delta).
\]
On each full interval $[k/n,(k+1)/n]$ the integral error is at most $n^{-1}\omega_{g,[\delta,1]}(1/n)$. The remaining interval $[\delta,k_\delta/n]$ has length at most $1/n$.
Consequently,
\[
 \left|\frac1n\sum_{k=k_\delta}^{n-1}g(r_k)-\int_\delta^1g(r)\dd r\right|
 \le\omega_{g,[\delta,1]}(1/n)+n^{-1}\|g\|_{L^\infty([\delta,1])}.
\]
Choose $0<\zeta<a_\beta$. The bound $\mathfrak G_\omega(r)\le Cr^\zeta$ and \eqref{S:eq:s4-small-index-sum} bound the two parts below $\delta$ by
$C\|f\|_\infty(n^{-\zeta}+\delta^\zeta/\zeta)$ and $C\|f\|_\infty\delta^\zeta/\zeta$. It follows that the limsup of the absolute difference in \eqref{S:eq:s4-weighted-riemann} is
at most $2C\|f\|_\infty\delta^\zeta/\zeta$. Sending $\delta\downarrow0$ proves the convergence of the integrals. Since the argument is deterministic after fixing the path, it also
holds for a continuous function $f$ chosen from that path.

\subsubsection{Intervals of Brownian time}
\label{S:sec:s4-brownian-time-intervals}

An interval in the operator variable may intersect several sampling intervals after the change to Brownian time. The signed fractions below specify the portion of each increment
used in the covariance calculation.

For $1\le j<n$, let
\[
 J_{n,j}:=[\vartheta_{n,j+1},\vartheta_{n,j}],
 h_{n,j}:=\vartheta_{n,j}-\vartheta_{n,j+1}=\frac4{\beta j+1}.
\]
The intervals $J_{n,j}$ partition $[0,\vartheta_{n,1}]$ in Brownian time. The intervals $I_{n,k}$ have length $1/n$ in the operator domain, and their images under the time change
may meet several $J_{n,j}$. For $\theta\ge0$, define
\begin{equation}\label{S:eq:interval-relative-time}
 u_{n,j}(\theta):=\frac1{h_{n,j}}\int_0^\theta
                  \mathbf1_{J_{n,j}}(v)\,\dd v,
 \lambda_{n,j}(\theta_0,\theta):=u_{n,j}(\theta)-u_{n,j}(\theta_0),
\end{equation}
and $\lambda_{n,k,j}(s):=\lambda_{n,j}(\vartheta_{n,k+1},\theta_{n,k}(s))$. In particular,
\[
 h_{n,j}\lambda_{n,j}(\theta_0,\theta)
    =\int_{\theta_0}^{\theta}\mathbf1_{J_{n,j}}(v)\,\dd v,
\]
where the integral is oriented when $\theta<\theta_0$. In $\lambda_{n,k,j}(s)$, the point of the operator domain is $r=(k+s)/n$ in the closure of the reflected subinterval, its
Brownian time is $\theta_{n,k}(s)$, and the sampling time is $\vartheta_{n,k+1}$.

For $d\in\mathbb Z$, put
\begin{equation}\label{S:eq:interval-fractions-limit}
 \Lambda_{\beta,d}(r,s):=
 \min\{1,\max\{0,d+1-s-b_\beta(1-r)\}\}
                  -\mathbf1_{\{d\ge1\}},  (r,s)\in[0,1]^2.
\end{equation}

Set
\[
 R_\beta=\left\lceil2C_\beta^0+2|b_\beta|+8\right\rceil,\qquad
 k_\beta=\lceil C_\beta^0\rceil+3.
\]

\begin{theorem}\label{S:thm:interval-fractions}
One has $\lambda_{n,k,j}=0$ for $|j-k|>R_\beta$ and $\Lambda_{\beta,d}=0$ for $|d|>R_\beta$. For every fixed $d$, uniformly over $1\le k,k+d<n$,
\begin{equation}\label{S:eq:interval-fractions-error}
 \sup_{0\le s\le1}
 \left|\lambda_{n,k,k+d}(s)-\Lambda_{\beta,d}(k/n,s)\right|
 \le C_{\beta,d}/k.
\end{equation}
\end{theorem}

\begin{proof}
Choose $C_\beta^0$ from \eqref{S:eq:sampling-position-mesh-distance}. Under the inverse time change $r=\me^{-\beta\theta/4}$, the interval $J_{n,j}$ becomes
$[\rho_{n,j-1},\rho_{n,j}]$, whose boundary points differ from $j/n$ and $(j+1)/n$ by at most $C_\beta^0/n$. The segment joining $(k+s)/n$ and $\rho_{n,k}$ lies inside
\[
 \left[(k-C_\beta^0)/n,\ (k+1+C_\beta^0)/n\right].
\]
A nonzero $\lambda_{n,k,j}(s)$ requires the segment between the two times to meet the interior of $J_{n,j}$. Consequently $k-1-2C_\beta^0\le j\le k+1+2C_\beta^0$. If
$k\ge k_\beta$, then $(k+s)/n>\rho_{n,0}$, so $0\le\theta_{n,k}(s)<\vartheta_{n,1}$. The sample time also lies in $[0,\vartheta_{n,1}]$, which the intervals $J_{n,j}$ partition up
to boundary points. Since $h_{n,j}\lambda_{n,j}(\theta_0,\theta)$ is the signed length traversed inside sampling interval $j$, summing the indicator functions gives
\begin{equation}\label{S:eq:interval-length-partition}
 \sum_{j=1}^{n-1}h_{n,j}\lambda_{n,k,j}(s)
     =\theta_{n,k}(s)-\vartheta_{n,k+1},
 \sum_{j=1}^{n-1}h_{n,j}|\lambda_{n,k,j}(s)|=\ell_{n,k}(s).
\end{equation}
For the second identity, each integral has the sign of $\theta_{n,k}(s)-\vartheta_{n,k+1}$.

Set $j=k+d$ and $c=\beta^{-1}$. Since $\vartheta_{n,j}$ decreases strictly with $j$,
\begin{equation}\label{S:eq:relative-time-at-sample}
 u_{n,k+d}(\vartheta_{n,k+1})=\mathbf1_{\{d\ge1\}}.
\end{equation}
For $d\le0$ the time is at or before the initial point of sampling interval $k+d$; for $d\ge1$ it is at or after its terminal point. This includes the two equalities at $d=0,1$.
Before applying the minimum and maximum, the fraction at the target time is
\[
 z^{\mathrm{time}}_{n,k,d}(s)
 :=\frac{\theta_{n,k}(s)-\vartheta_{n,j+1}}{h_{n,j}}
 =(j+c)\left[
 \log\frac{j+1}{k+s}
 -b_\beta\left(\frac1{j+1}-\frac1n\right)+O_\beta(j^{-2})
 \right].
\]
Here \eqref{S:eq:subinterval-sampling-time-expansion} applies uniformly also when $j=n-1$. For $k\ge2|d|+2$, one has $j\asymp_d k$. Taylor's formula for $\log(1+(d+1-s)/(k+s))$
therefore gives
\[
 (j+c)\log\frac{j+1}{k+s}=d+1-s+O_{\beta,d}(k^{-1}).
\]
Also, directly,
\[
 (j+c)\left(\frac1{j+1}-\frac1n\right)
 =1-\frac kn+\frac{c-1}{j+1}-\frac{d+c}{n}
 =1-\frac kn+O_{\beta,d}(k^{-1}).
\]
The remaining error is $O_{\beta,d}(k^{-1})$. It follows that $z^{\mathrm{time}}_{n,k,d}(s)=d+1-s-b_\beta(1-k/n)+O_{\beta,d}(k^{-1})$, uniformly in $s$. The map
$v\mapsto\min(1,\max(0,v))$ is Lipschitz with constant one, so this estimate and \eqref{S:eq:relative-time-at-sample} prove \eqref{S:eq:interval-fractions-error}. The finitely many
smaller values of $k$ are covered by increasing $C_{\beta,d}$, since both fractions have absolute value at most one. For any $(r,s),(r',s')\in[0,1]^2$, the same Lipschitz bound
gives
\[
 |\Lambda_{\beta,d}(r,s)-\Lambda_{\beta,d}(r',s')|
       \le |b_\beta||r-r'|+|s-s'|.
\]
If $d>1+|b_\beta|$, the minimum and maximum in \eqref{S:eq:interval-fractions-limit} give one, which cancels the indicator. If $d<-1-|b_\beta|$, they give zero and the indicator
vanishes. This proves the stated finite support. At
$r=1$,
\begin{equation}\label{S:eq:fractions-regular-boundary}
\Lambda_{\beta,0}(1,s)=1-s,
 \Lambda_{\beta,d}(1,s)=0, (d\ne0).
\end{equation}

Put $v=1-s-b_\beta(1-r)$. Since $\min(1,\max(0,d))=\mathbf1_{\{d\ge1\}}$ for integer $d$,
\[
 \Lambda_{\beta,d}(r,s)
 =\int_0^v\mathbf1_{(-d,\,1-d)}(x)\,\dd x,
\]
where $\int_0^v=-\int_v^0$ for $v<0$. These intervals partition the real line up to integers. Their integrals all have the sign of $v$, and only finitely many are nonzero on the
bounded range of $v$ considered here. Consequently,
\begin{equation}\label{S:eq:limiting-fractions-sum}
 \sum_{d\in\mathbb Z}\Lambda_{\beta,d}(r,s)=1-s-b_\beta(1-r),
 \sum_{d\in\mathbb Z}|\Lambda_{\beta,d}(r,s)|=|1-s-b_\beta(1-r)|.
\end{equation}
Integration of the absolute identity gives $\int_0^1\sum_d|\Lambda_{\beta,d}(r,s)|\dd s=\kappa_\beta(r)$. The substitution is $u=1-s$.

On $k/n\ge\delta>0$, \eqref{S:eq:interval-fractions-error} gives the uniform error $O_{\beta,\delta,d}(n^{-1})$, including subintervals approaching the regular boundary.
\end{proof}

For $r\in I_{n,k}$ with $k\ge m$, the segment between $\vartheta_{n,k+1}$ and $\theta_\beta(r)$ can traverse a sampling interval $J_{n,j}$ with $j<m$. More precisely, for any
$m\ge k_\beta$, the intervals occurring for $m\le k<n$ belong to
\begin{equation}\label{S:eq:additional-sampling-intervals}
 \{j:\max(1,m-R_\beta)\le j<n\}.
\end{equation}
This adds at most $R_\beta$ indices to the family of sampling intervals. The subintervals $I_{n,k}$ still have indices $m\le k<n$. The sums over sampling intervals include all of
these additional indices.

\phantomsection\label{S:lem:integrated-fractions} Let $m_n\to\infty$, $m_n/n\to0$, and let $F_{n,k}(s)$ be measurable
matrix functions for $m_n\le k<n$, $0\le s\le1$. Suppose that $F\colon(0,1]\times[0,1]\to\mathbb R^{d_1\times d_2}$ is continuous and, for every $\delta>0$,
\[
 e_n(\delta):=\sup_{\substack{m_n\le k<n,\ k/n\ge\delta\\0\le s\le1}}
       \|F_{n,k}(s)-F(k/n,s)\|\longrightarrow0,
 \qquad
 \sup_s\|F_{n,k}(s)\|\le C(k/n)^\zeta
\]
for some $\zeta>0$. We use any fixed matrix norm. For a fixed integer $q$, set
\[
 G_{n,j}^{(q)}=\sum_{k=m_n}^{n-1}\int_0^1
       \lambda_{n,k,j}(s)\lambda_{n,k,j+q}(s)F_{n,k}(s)\,\dd s,
 \qquad 1\le j<n,
\]
with $\lambda_{n,k,i}=0$ when $i\notin\{1,\ldots,n-1\}$. The limiting function is
\begin{equation}\label{S:eq:integrated-fractions-limit}
 G^{(q)}(r)=\sum_{d\in\mathbb Z}\int_0^1
       \Lambda_{\beta,d}(r,s)\Lambda_{\beta,d+q}(r,s)F(r,s)\,\dd s.
\end{equation}
We shall prove uniform convergence on $j/n\ge\delta$, including indices next to $n$.

By Theorem~\ref{S:thm:interval-fractions}, only $|j-k|\le R_\beta$ can occur in the discrete sum. For $m_n>2R_\beta$, such indices satisfy $k\le2j$. Therefore
\[
 G_{n,j}^{(q)}=0\quad(j<m_n-R_\beta),\qquad
 \|G_{n,j}^{(q)}\|\le C(2R_\beta+1)2^\zeta(j/n)^\zeta.
\]
The limit in \eqref{S:eq:integrated-fractions-limit} is a finite sum: only $|d|\le R_\beta$ and $|d+q|\le R_\beta$ occur. It is continuous on each compact subinterval of $(0,1]$.
Approximating any $r>0$ by $k/n$ in the assumed bounds also gives $\sup_s\|F(r,s)\|\le Cr^\zeta$; hence $F$ and $G^{(q)}$ extend continuously to zero with value zero.

Fix $0<\delta<1$ and take $n$ so large that $m_n+R_\beta<\delta n/2$ and $\delta n>2(|q|+R_\beta)$. Put
\[
 M_\delta=\sup_{(r,s)\in[\delta/2,1]\times[0,1]}\|F(r,s)\|,
\]
and let $\omega_{F,\delta}(h)$ be the supremum of $\|F(r,s)-F(r',s)\|$ over this set with $|r-r'|\le h$. For $j/n\ge\delta$, write $k=j-d$. Whenever all the indices are valid,
Theorem~\ref{S:thm:interval-fractions} and $|ab-cd|\le |a-c|+|b-d|$ for $|a|,|b|,|c|,|d|\le1$ yield
\[
 \begin{split}
 &\big\|\lambda_{n,j-d,j}(s)\lambda_{n,j-d,j+q}(s)F_{n,j-d}(s)
       -\Lambda_{\beta,d}(j/n,s)\Lambda_{\beta,d+q}(j/n,s)F(j/n,s)\big\|\\
 &\hspace{8mm}\le e_n(\delta/2)+\omega_{F,\delta}(R_\beta/n)
                      +C_{\beta,q,\delta}M_\delta/n.
 \end{split}
\]
Here we first evaluate both limiting factors at $(j-d)/n$, and then use their Lipschitz bounds to replace this point by
$j/n$.

A missing $k=j-d\ge n$ has $d<0$ and $1-j/n\le R_\beta/n$. By \eqref{S:eq:fractions-regular-boundary}, $|\Lambda_{\beta,d}(j/n,s)|\le |b_\beta|R_\beta/n$. A missing interval
$j+q\ge n$ requires $q>0$ and $1-j/n\le q/n$. At $r=1$ the product $\Lambda_{\beta,d}\Lambda_{\beta,d+q}$ is zero, since $q\ne0$; its absolute value at $j/n$ is thus at most
$C_{\beta,q}/n$. There are no missing indices at the lower boundary for these $n$. Summing at most $2R_\beta+1$ integrals proves
\[
 \sup_{\delta\le j/n<1}\|G_{n,j}^{(q)}-G^{(q)}(j/n)\|
 \le C_{\beta,q,\delta}\bigl[e_n(\delta/2)+\omega_{F,\delta}(R_\beta/n)
                         +(1+M_\delta)/n\bigr]\longrightarrow0.
\]
Since $|\rho_{n,j}-j/n|\le(2+\beta^{-1})/n$, uniform continuity of $G^{(q)}$ permits $j/n$ to be replaced by $\rho_{n,j}$ in this conclusion. For random functions, the argument
applies on every path on which these assumptions hold.

\subsubsection{Integrated coefficients and summation by parts}
\label{S:sec:s4-integrated-coefficients}

The conditional mean and the summation by parts argument use the signed displacement from each sampling time. The local Brownian variance involves the absolute displacement. We
first compute the signed integrated coefficients exactly, and then bound their variation in the index of the sampling interval.

For the sum involving $\mathcal B_n$, integrate $\lambda_{n,k,j}$ before taking absolute values.

\phantomsection\label{S:lem:interval-weight-formula} Let $n\ge2$, $1\le m<n$, and $1\le j<n$. Set
\[
 x_j:=j+\beta^{-1},
 \rho^-_{n,j}:=\rho_{n,j-1},
 \rho^+_{n,j}:=\rho_{n,j},  a:=m/n.
\]
Thus $\rho^-_{n,j}=\rho^+_{n,j}\me^{-1/x_j}$ and $0<\rho^-_{n,j}<\rho^+_{n,j}\le1$. Define
\[
 \omega_{n,j}^{[m]}:=\sum_{k=m}^{n-1}\int_0^1\lambda_{n,k,j}(s)\,\dd s.
\]
Writing $v_+=\max\{v,0\}$, the formula is
\begin{equation}\label{S:eq:interval-weight-integral}
 \omega_{n,j}^{[m]}
 =n x_j\int_{\rho^-_{n,j}}^{\rho^+_{n,j}}
                  \frac{(u-a)_+}{u}\,\dd u-(j-m)_+.
\end{equation}
Equivalently, if $\mathcal I_{n,j}(a)$ denotes the integral in \eqref{S:eq:interval-weight-integral} including the factor $x_j$ but not $n$, then
\begin{equation}\label{S:eq:interval-weight-cutoff-formula}
 \mathcal I_{n,j}(a)=
 \begin{cases}
  x_j(\rho^+_{n,j}-\rho^-_{n,j})-a,
      &0\le a\le\rho^-_{n,j},\\[1mm]
  x_j\bigl(\rho^+_{n,j}-a-a\log(\rho^+_{n,j}/a)\bigr),
      &\rho^-_{n,j}<a<\rho^+_{n,j},\\[1mm]
  0,&\rho^+_{n,j}\le a\le1.
 \end{cases}
\end{equation}
The expressions agree at the two dividing points. If $j\ge m$ and $a\le\rho^-_{n,j}$, then
\begin{equation}\label{S:eq:interval-weight-interior}
 \omega_{n,j}^{[m]}
 =n x_j\rho_{n,j}(1-\me^{-1/x_j})-j.
\end{equation}
In particular, this coefficient is independent of the cutoff $m$.

The target coordinate on sampling interval $j$, expressed in the coordinate $r=1-t$, is
\[
 u_{n,j}(\theta_\beta(r))
 =\min\{1,\max\{0,x_j\log(\rho^+_{n,j}/r)\}\}.
\]
It is one for $r\le\rho^-_{n,j}$, zero for $r\ge\rho^+_{n,j}$, and $x_j\log(\rho^+_{n,j}/r)$ in between. The identity \eqref{S:eq:relative-time-at-sample} says exactly that
$u_{n,j}(\vartheta_{n,k+1})=\mathbf1_{\{k<j\}}$. There are $(j-m)_+$ such indices in $\{m,\ldots,n-1\}$. The substitution $r=(k+s)/n$, with $\dd s=n\,\dd r$, therefore gives
\[
 \omega_{n,j}^{[m]}
 =n\int_a^1 u_{n,j}(\theta_\beta(r))\,\dd r-(j-m)_+.
\]
Here $s=nr-k$, which accounts for the Jacobian $\dd s=n\,\dd r$.

The function $u_{n,j}$ has the nonnegative integral representation
\[
 u_{n,j}(\theta_\beta(r))
 =x_j\int_{\rho^-_{n,j}}^{\rho^+_{n,j}}
           \mathbf1_{\{r<u\}}\,\frac{\dd u}{u}.
\]
For $r\le\rho^-_{n,j}$ this is one because $x_j\log(\rho^+_{n,j}/\rho^-_{n,j})=1$; the other two cases follow by changing the lower integration limit to $r$. Tonelli's theorem
yields
\[
 \int_a^1u_{n,j}(\theta_\beta(r))\,\dd r
 =x_j\int_{\rho^-_{n,j}}^{\rho^+_{n,j}}
               \frac{(u-a)_+}{u}\,\dd u.
\]
This proves \eqref{S:eq:interval-weight-integral}. Integrating $1-a/u$ in each of the three ranges proves \eqref{S:eq:interval-weight-cutoff-formula}, including equality at its
dividing points. If $j\ge m$ and $a\le\rho^-_{n,j}$, subtracting $j-m$ cancels the term $na=m$ and gives \eqref{S:eq:interval-weight-interior}.

\begin{lemma}\label{S:lem:interval-weight-variation}
Fix $\beta>0$ and $\max(k_\beta,2R_\beta+2)\le m<n$. With $C_\beta$ independent of $m,n$,
\begin{equation}\label{S:eq:interval-weight-variation}
 \sup_{m-R_\beta\le j<n}|\omega_{n,j}^{[m]}|\le2R_\beta+1,
 \sum_{j=m-R_\beta}^{n-2}
       |\omega_{n,j+1}^{[m]}-\omega_{n,j}^{[m]}|\le C_\beta.
\end{equation}
The same uniform variation bound holds after extending the weights by zero to $\mathbb Z$.
\end{lemma}

\begin{proof}

Theorem~\ref{S:thm:interval-fractions} gives $\lambda_{n,k,j}=0$ for $|k-j|>R_\beta$, and every fraction has absolute value at most one. At most $2R_\beta+1$ partition subintervals
contribute to one coefficient. Thus $|\omega_{n,j}^{[m]}|\le2R_\beta+1$ throughout the stated range. By \eqref{S:eq:sampling-position-mesh-distance},
\[
 n\rho^-_{n,j}=n\rho_{n,j-1}\ge j-C_\beta^0.
\]
The choice of $R_\beta$ in Theorem~\ref{S:thm:interval-fractions} has $R_\beta>C_\beta^0$. Hence $j\ge m+R_\beta$ implies $\rho^-_{n,j}\ge m/n$, and
\eqref{S:eq:interval-weight-interior} applies. Within the range under consideration, the cutoff can affect that expression only at indices in $[m-R_\beta,m+R_\beta-1]$, a set of at
most $2R_\beta$ indices. We therefore treat these indices separately at the cutoff.

Put $c=\beta^{-1}$, $b=b_\beta$, $x=j+c$, and $y=n\rho_{n,j}$. The expansion of $\rho_{n,j}$ \eqref{S:eq:sampling-position-expansion} gives
\begin{equation}\label{S:eq:interval-weight-position-expansion}
 y=x+\frac12+b\frac{j+1}{n}+e_{n,j}^{\rho},
   |e_{n,j}^{\rho}|\le\frac{C_\beta}{j+1}.
\end{equation}
Here $j+1-b=x+1/2$. This bound is uniform through $j=n-1$, where $\rho_{n,n-1}=1$ and $e_{n,n-1}^{\rho}=0$. For $x\ge1$, Taylor's formula gives
\[
 x(1-\me^{-1/x})=1-\frac1{2x}+r_-(x),
   0\le r_-(x)\le\frac1{6x^2}.
\]
Write $d=b(j+1)/n+e_{n,j}^{\rho}$, which is bounded by a constant depending only on $\beta$. Substitution into \eqref{S:eq:interval-weight-interior} gives the exact rearrangement
\[
 yx(1-\me^{-1/x})-j
 =c+d-\frac1{4x}-\frac{d}{2x}+y r_-(x).
\]
Here $y\le C_\beta x$, $x\ge j$, and $1/n\le1/j$. Since $c+bj/n=\tfrac12-b(1-j/n)$, every remaining term is bounded by $C_\beta/j$. We have obtained
\begin{equation}\label{S:eq:interval-weight-limit}
 \omega_{n,j}^{[m]}=\frac12-b_\beta(1-j/n)+O_\beta(j^{-1}),
 \qquad m+R_\beta\le j<n.
\end{equation}

For the difference of consecutive coefficients, use the recurrence
\[
 \rho_{n,j+1}=\rho_{n,j}\me^{1/(j+1+c)}.
\]
If $j\ge m+R_\beta$, both coefficients have the form \eqref{S:eq:interval-weight-interior}, and therefore
\begin{equation}\label{S:eq:interval-weight-exact-difference}
 \omega_{n,j+1}^{[m]}-\omega_{n,j}^{[m]}
 =yD(x)-1,
 D(x):=(x+1)(\me^{1/(x+1)}-1)-x(1-\me^{-1/x}).
\end{equation}
Taylor expansion through order three gives
\begin{equation}\label{S:eq:interval-weight-taylor}
 D(x)=\frac1x-\frac1{2x^2}+r_D(x),
   |r_D(x)|\le x^{-3}, (x\ge1).
\end{equation}
Expand the two exponentials through their cubic terms. After multiplication by $x+1$ and $x$, respectively, the remainder bounds are $\me/(24(x+1)^3)$ and $1/(24x^3)$. The
displayed part of the difference is
\[
 \frac1{2(x+1)}+\frac1{2x}
       +\frac16\bigl((x+1)^{-2}-x^{-2}\bigr).
\]
Replacing $1/(2(x+1))$ by $1/(2x)-1/(2x^2)$ has error at most $1/(2x^3)$, and the last parenthesis has absolute value at most $2x^{-3}$. The total error is at most
$[\tfrac12+\tfrac13+(\me+1)/24]x^{-3}\le x^{-3}$, which bounds $r_D$ in \eqref{S:eq:interval-weight-taylor}.

Substitute \eqref{S:eq:interval-weight-position-expansion} into \eqref{S:eq:interval-weight-exact-difference}. Subtracting $b/n$ gives
\begin{equation}\label{S:eq:interval-weight-difference-error}
 yD(x)-1-\frac bn =-\frac1{4x^2} +\frac bn\left(\frac{1-c}{x}-\frac{j+1}{2x^2}\right) +e_{n,j}^{\rho}\left(\frac1x-\frac1{2x^2}\right) +y r_D(x).
\end{equation}
Each term is bounded by $C_\beta/j^2$: use $n\ge j+2$, $x\ge j$, $(j+1)/x\le2$, $|e_{n,j}^{\rho}|\le C_\beta/(j+1)$, and $y\le C_\beta x$. Consequently,
\begin{equation}\label{S:eq:interval-weight-difference}
 \left|\omega_{n,j+1}^{[m]}-\omega_{n,j}^{[m]}-\frac{b_\beta}{n}\right|
 \le\frac{C_\beta}{j^2},\qquad m+R_\beta\le j\le n-2.
\end{equation}

Set $J_0=m+R_\beta$. There are at most $2R_\beta$ adjacent pairs whose first index lies between $m-R_\beta$ and $\min(J_0-1,n-2)$. Each has difference at most $2(2R_\beta+1)$. For
the remaining pairs, \eqref{S:eq:interval-weight-difference} gives
\[
 \sum_{j=J_0}^{n-2}|\omega_{n,j+1}^{[m]}-\omega_{n,j}^{[m]}|
 \le |b_\beta|+C_\beta\sum_{j=J_0}^\infty j^{-2}
 \le |b_\beta|+\frac{2C_\beta}{m}.
\]
If $J_0>n-2$, this sum is empty and all pairs were counted in the preceding finite set. Thus the total variation is bounded by $4R_\beta(2R_\beta+1)+|b_\beta|+2C_\beta/m$, proving
\eqref{S:eq:interval-weight-variation}. Extending by zero adds exactly the two jumps $|\omega_{n,m-R_\beta}^{[m]}|+|\omega_{n,n-1}^{[m]}|$. Their sum is at most $2(2R_\beta+1)$,
which proves \eqref{S:eq:interval-weight-zero-extension}.

For $m+R_\beta\le J\le K\le n-1$, summing \eqref{S:eq:interval-weight-difference} over $J\le j<K$ gives
\begin{equation}\label{S:eq:interval-weight-local-variation}
 \sum_{j=J}^{K-1}|\omega_{n,j+1}^{[m]}-\omega_{n,j}^{[m]}|
 \le |b_\beta|\frac{K-J}{n}+\frac{C_\beta}{J}.
\end{equation}
For the zero extension $\bar\omega_{n,j}^{[m]}$ this becomes
\begin{equation}\label{S:eq:interval-weight-zero-extension}
 \sum_{j\in\mathbb Z}
       |\bar\omega_{n,j+1}^{[m]}-\bar\omega_{n,j}^{[m]}|
 \le C_\beta.
\end{equation}

\end{proof}

For comparison, the functions with a fixed index difference satisfy the weaker but sufficient bound
\begin{equation}\label{S:eq:interval-fractions-adjacent-bound}
 \sup_{0\le s\le1}|\lambda_{n,k+1,k+1+d}(s)-\lambda_{n,k,k+d}(s)|
 \le\frac{C_{\beta,d}}k
\end{equation}
whenever both index pairs are valid. Indeed, apply \eqref{S:eq:interval-fractions-error} at $k$ and $k+1$ and use that $\Lambda_{\beta,d}(r,s)$ is $|b_\beta|$-Lipschitz in $r$; the
intermediate difference is at most $|b_\beta|/n$. This proof also holds where the minimum or maximum changes branch. The exact recurrence for the integrated coefficient gives the
sharper, summable difference in \eqref{S:eq:interval-weight-difference}.

The signed coefficients have limit $\tfrac12-b_\beta(1-r)$. The factor $\kappa_\beta(r)=\int_0^1\sum_d|\Lambda_{\beta,d}(r,s)|\,\dd s$ integrates the absolute time difference in
units of the limiting interval length. For $\beta>2$, the difference is $b_\beta^2(1-r)^2$. At $\beta=2$, the leading adjacent difference in \eqref{S:eq:interval-weight-difference}
vanishes, while the discrete remainder can still vary with $j$. The final sampling interval is $j=n-1$, with $\rho^+_{n,n-1}=1$.

\phantomsection\label{S:cor:interval-weight-summation} We shall use the following identities obtained by summation by
parts. Let $m,n$ satisfy Lemma~\ref{S:lem:interval-weight-variation}, and let $m-R_\beta\le A\le B<n$. For vectors $z_A,\ldots,z_B$ in a normed space, put $S_r=\sum_{j=A}^rz_j$.
Then
\begin{equation}\label{S:eq:interval-weight-summation}
 \sum_{j=A}^B\omega_{n,j}^{[m]}z_j =\omega_{n,B}^{[m]}S_B -\sum_{j=A}^{B-1}(\omega_{n,j+1}^{[m]}-\omega_{n,j}^{[m]})S_j,\qquad \left\|\sum_{j=A}^B\omega_{n,j}^{[m]}z_j\right\| \le C_\beta\max_{A\le r\le B}\|S_r\|.
\end{equation}
The constant is independent of the boundary points $A,B$ and of the vectors. More generally, for any real weights $w_j$ and scalars $z_j=\tfrac12(q_{j-1}-q_j)+e_j$, one has the
identity
\begin{equation}\label{S:eq:interval-weight-telescoping}
 \sum_{j=A}^Bw_jz_j =\frac12\left(w_Aq_{A-1}-w_Bq_B +\sum_{j=A}^{B-1}(w_{j+1}-w_j)q_j\right) +\sum_{j=A}^Bw_je_j.
\end{equation}
If $q_j\ge0$, $W=\max_{A\le j\le B}|w_j|$, and $V=\sum_{j=A}^{B-1}|w_{j+1}-w_j|$, it follows that
\begin{equation}\label{S:eq:interval-weight-telescoping-bound}
 \left|\sum_{j=A}^Bw_jz_j\right|
 \le\frac W2(q_{A-1}+q_B)
       +\frac V2\max_{A\le j\le B}q_j
       +W\sum_{j=A}^B|e_j|.
\end{equation}
For $w_j=\omega_{n,j}^{[m]}$, one has $W+V\le C_\beta$.

Indeed, with $S_{A-1}=0$,
\[
 \sum_{j=A}^B\omega_{n,j}^{[m]}(S_j-S_{j-1}) =\sum_{j=A}^B\omega_{n,j}^{[m]}S_j -\sum_{j=A-1}^{B-1}\omega_{n,j+1}^{[m]}S_j =\omega_{n,B}^{[m]}S_B -\sum_{j=A}^{B-1}(\omega_{n,j+1}^{[m]}-\omega_{n,j}^{[m]})S_j.
\]
Taking norms bounds this expression by $(|\omega_{n,B}^{[m]}|+\sum_{j=A}^{B-1}|\omega_{n,j+1}^{[m]}-\omega_{n,j}^{[m]}|) \max_{A\le r\le B}\|S_r\|$.
Lemma~\ref{S:lem:interval-weight-variation} proves the claimed constant. For \eqref{S:eq:interval-weight-telescoping}, collect each coefficient of $q_j$ directly. The two boundary
contributions are $w_Aq_{A-1}$ and $-w_Bq_B$. Taking absolute values gives \eqref{S:eq:interval-weight-telescoping-bound}. The supremum and variation bounds in
Lemma~\ref{S:lem:interval-weight-variation} apply to every subinterval, proving the uniform bound over subintervals. The identities hold for arbitrary vectors and weights,
including nonmonotone sequences $q_j$.

\subsubsection{Truncation of the Green kernel}
\label{S:sec:s4-truncation}

The smaller cutoff is used in the construction of the terminal part of the finite walk. The larger cutoff is used in the kernel expansion. For $k\ge m_n$ the Brownian time length
of an operator interval satisfies
\[
 \tau_\beta(1-k/n)-\tau_\beta(1-(k+1)/n)
 =\frac4\beta\log(1+1/k)\le\frac4{\beta m_n}\longrightarrow0.
\]
Also $m_n/n\to0$: with $x=\log n$, the ratio is bounded by $(x^{512}+1)\me^{-x}$, which tends to zero. Choose $n_0$ such that $1\le\underline m_n\le m_n<n$ for all $n\ge n_0$. The
error estimates proved later in \eqref{S:eq:circular-brownian-retained-rate} and \eqref{S:eq:continuous-bridge-remainder-rate} are, after multiplication of their norms by
$\sqrt n$, bounded by constant multiples of
\[
 \frac{(\log n)^{11/2}}{\sqrt{m_n}}
 \le(\log n)^{-501/2},\qquad
 \frac{(1+\log n)^6}{\sqrt{m_n}}
 \le2^6(\log n)^{-250}\quad(n\ge\me).
\]
These inequalities use only $m_n\ge(\log n)^{512}$. A squared tail bound of the form $C(m_n/n)^{1+\delta}$, with any fixed $\delta>0$, also tends to zero after multiplication by
$n$, since
\[
 n(m_n/n)^{1+\delta}
 \le 2^{1+\delta}(\log n)^{512(1+\delta)}n^{-\delta}
 \longrightarrow0.
\]
Theorem~\ref{S:thm:full-boundary-tails} supplies these estimates for the three inverse kernels. The cutoffs are used only to estimate their differences.

We use the reflected coordinate only through the common unitary map $g(t)\mapsto g(1-r)$. For $0\le\delta\le1$, let $P_\delta^{\mathrm e}$ be multiplication by
$\mathbf1_{[\delta,1)}$ on $L^2((0,1);\mathbb R^2)$, and put
\[
 \Pi_\delta A=P_\delta^{\mathrm e}AP_\delta^{\mathrm e},\qquad
 \mathcal T_\delta A=A-\Pi_\delta A.
\]
The superscript $1-\delta$ for an inverse still refers to its cutoff in the original $t$ coordinate, so $R_X^{1-\delta}=\Pi_\delta R_X$ after reflection. For a Green inverse write
$p=|a|^2$, $q=|c|^2$ and
$C(r)=\int_r^1p(u)\dd u$.

\begin{lemma}\label{S:lem:orthogonal-truncation-masses}
The map $\Pi_\delta$ is an orthogonal projection on $\mathfrak S_2$. For $A,B\in\mathfrak S_2$,
\begin{equation}\label{S:eq:boundary-pair-orthogonality}
 \|A-B\|_{\mathfrak S_2}^2
 =\|\Pi_\delta(A-B)\|_{\mathfrak S_2}^2
       +\|\mathcal T_\delta A-\mathcal T_\delta B\|_{\mathfrak S_2}^2.
\end{equation}
For the Green inverse,
\begin{equation}\label{S:eq:4.1}
 \|R_X-R_X^{1-\delta}\|_{\mathfrak S_2}^2
       =\frac12\int_0^\delta q(r)C(r)\dd r.
\end{equation}
\end{lemma}

\begin{proof}
Identify a Hilbert--Schmidt operator with its square-integrable matrix kernel. Under this isometry, $\Pi_\delta$ is multiplication by the real indicator $\mathbf1_{[\delta,1)^2}$.
It is idempotent and self-adjoint, so it is an orthogonal projection. Its complement is multiplication by $\mathbf1_{\{\min(r,r')<\delta\}}$. Integrating the squared norm on the
two disjoint sets proves \eqref{S:eq:boundary-pair-orthogonality}. The disjoint decomposition
\[
 \{\min(r,r')<\delta\}
 =\{\min(r,r')<d\}\mathbin{\dot\cup}
                \{d\le\min(r,r')<\delta\}
\]
for $0\le d\le\delta$ also gives
\begin{equation}\label{S:eq:boundary-nested-orthogonality}
 \|A-\Pi_\delta A\|_{\mathfrak S_2}^2
 =\|A-\Pi_d A\|_{\mathfrak S_2}^2
       +\|\Pi_d A-\Pi_\delta A\|_{\mathfrak S_2}^2.
\end{equation}
Writing $P=P_\delta^{\mathrm e}$ and $Q=I-P$ gives $\mathcal T_\delta A=PAQ+QAP+QAQ$. The three terms have disjoint kernel supports.

Equation \eqref{S:eq:4.1} is the exact compression identity \eqref{S:eq:singular-truncation-hs-tail} in the reflected coordinate; for regular coefficients it follows from the same
restricted kernel in Lemma~\ref{S:lem:regular-interval-restriction}. It includes the terminal square and both mixed rectangles.

For a coefficient constant on the $n$ intervals, put $p_k=|a_k|^2$ and $q_k=|c_k|^2$. At $r=(k+s)/n$, $0<s<1$,
\[
 C(r)=\frac1n\sum_{l=k+1}^{n-1}p_l+\frac{1-s}{n}p_k.
\]
Integrating $q_kC(r)$ with $\dd r=\dd s/n$ and using $\int_0^1(1-s)\dd s=1/2$ gives
\begin{equation}\label{S:eq:boundary-exact-mesh-tail}
 \|R_X-R_X^{1-m/n}\|_{\mathfrak S_2}^2
 =\frac1{2n^2}\sum_{k=0}^{m-1}q_k\sum_{l=k+1}^{n-1}p_l
       +\frac1{4n^2}\sum_{k=0}^{m-1}p_kq_k,
 \qquad 0\le m\le n.
\end{equation}
The last term includes both triangles of each square $I_{n,k}\times I_{n,k}$; their total area is $1/n^2$ and each has squared kernel norm $p_kq_k/4$. For $0\le u\le m\le n$,
subtraction of the two integrals in \eqref{S:eq:4.1}, or \eqref{S:eq:boundary-nested-orthogonality}, shows that $\|R_X^{1-u/n}-R_X^{1-m/n}\|_{\mathfrak S_2}^2$ is the same sum with
$u\le k<m$. At $m=n$ the compressed operator is zero.
\end{proof}

\begin{proposition}\label{S:prop:4.3}
Almost surely, for every $\varepsilon>0$,
\begin{equation}\label{S:eq:boundary-brownian-uniform-tails}
 \begin{aligned}
 \|\mathsf R_\beta^{\mathrm{sine}}
            -\mathsf R_\beta^{\mathrm{sine},1-\delta}\|_{\mathfrak S_2}^2
   &\le C_{\beta,\omega,\varepsilon}\delta^{1+a_\beta-\varepsilon},
       &&0<\delta\le\tfrac12,\\
 \|\mathsf R_{\beta,n}^{\mathrm{br}}
            -\mathsf R_{\beta,n}^{\mathrm{br},1-m/n}\|_{\mathfrak S_2}^2
   &\le C_{\beta,\omega,\varepsilon}(m/n)^{1+a_\beta-\varepsilon},
       &&1\le m\le n/2.
 \end{aligned}
\end{equation}
One constant works for all $m$ and all sufficiently large $n$.
\end{proposition}

\begin{proof}
The first estimate is \eqref{S:eq:sine-singular-hs-identities} when $\varepsilon<a_\beta$; for $0<\varepsilon_0<\min\{\varepsilon,a_\beta\}$, use
$\delta^{1+a_\beta-\varepsilon_0}\le\delta^{1+a_\beta-\varepsilon}$ on $0<\delta\le1$. For the sampled inverse choose
\[
 0<\nu<\min(\varepsilon/4,a_\beta/4,1/4),\qquad
 \nu<(1-\gamma)/2\quad\hbox{if }\gamma<1.
\]
Write $\bar r_k=(k+1)/n$, $p_k=p_{n,k}$, $q_k=q_{n,k}$. Proposition~\ref{S:prop:5.2} gives, including $k=0$,
\[
 p_k\le C\bar r_k^{-\gamma-\nu},\qquad
 q_k\le C\bar r_k^{\gamma-\nu},\qquad
 p_kq_k\le C\bar r_k^{-2\nu}.
\]
Monotonicity of $u^{-\gamma-\nu}$ gives
\[
 \frac1n\sum_{l=k+1}^{n-1}p_l
 \le\frac Cn\sum_{v=k+2}^n(v/n)^{-\gamma-\nu}
 \le C\int_{\bar r_k}^1u^{-\gamma-\nu}\dd u.
\]
If $\gamma\ge1$, the right side is at most $C\bar r_k^{1-\gamma-\nu}$, including $\bar r_k^{-\nu}$ when $\gamma=1$. If $\gamma<1$, it is bounded by a constant. Multiplication by
$q_k$ therefore gives
\[
 q_k\frac1n\sum_{l=k+1}^{n-1}p_l\le C\bar r_k^{a_\beta-2\nu}.
\]
Set $\delta=m/n$. Since $a_\beta-2\nu>0$, the part of \eqref{S:eq:boundary-exact-mesh-tail} with distinct indices is at most
\[
 \frac Cn\sum_{v=1}^m(v/n)^{a_\beta-2\nu}
    \le\frac{Cm}{n}(m/n)^{a_\beta-2\nu}
    =C\delta^{1+a_\beta-2\nu}.
\]
For the contribution of the diagonal squares, $2\nu<1$ and
\[
 \sum_{v=1}^m v^{-2\nu}
    \le1+\int_1^m x^{-2\nu}\dd x
    \le C_\nu m^{1-2\nu}.
\]
Consequently that contribution is at most
\[
 \frac C{n^2}\sum_{v=1}^m(v/n)^{-2\nu}
       \le\frac Cn\delta^{1-2\nu}\le C\delta^{2-2\nu}.
\]
The last inequality uses $1/n\le\delta$, also when $m=1$. Both $1+a_\beta-2\nu$ and $2-2\nu$ are at least $1+a_\beta-\varepsilon$. This proves the second estimate with a constant
independent of $m,n$.
\end{proof}

The circular path requires a separate estimate near its terminal point. The following bound concerns the auxiliary simultaneous construction used in
Section~\ref{S:sec:coupling-the-finite-transfer-increments-to-brownian-increments}; its comparison with the Brownian samples is made in
Section~\ref{S:sec:comparison-of-the-finite-and-sampled-brownian-operators}.

\subsubsection{A uniform estimate near the singular boundary}
\label{S:sec:a-uniform-estimate-near-the-singular-boundary}

The portion nearest the singular boundary requires a bound that is uniform over matrix size. An auxiliary construction based on stopping the Brownian motion supplies that estimate.
It is used only for this comparison; the sampled operator itself retains its deterministic observation times.

Let $B$, with $B(0)=\ii$, and $(\tau_{n,j}^{\star})$ be the simultaneous stopped Brownian array of~\cite{valko2020operator}, and write $z_{n,j}^{\star}=B(\tau_{n,j}^{\star})$. That
paper constructs the stopping times simultaneously. The superscript distinguishes these stopping times from the deterministic lengths $\tau_{n,k}$ in
\eqref{S:eq:sampling-offset-definitions}. The sampled Brownian operator uses the deterministic times $\vartheta_{n,j}$. Set
\[
 u_n=\lceil(\log n)^6\rceil=\underline m_n,\qquad
 L_n=\log\log n.
\]
All estimates involving these quantities are for sufficiently large $n$, so $L_n\ge1$ and $u_n<n$. We use the following facts from that construction. The same random radii
\[
 Y_k=2\operatorname{arctanh}\sqrt{\xi_k},\qquad
 \xi_k\sim\operatorname{Beta}(1,\beta k/2),
\]
are used in every row when $k<u_n$. At these indices, $\tau_{n,k}^{\star}$ is the first exit after $\tau_{n,k+1}^{\star}$ from the hyperbolic ball of radius $Y_k$ centered at
$B(\tau_{n,k+1}^{\star})$. Moreover, on a single event of probability one, the estimate for the retained stopping times gives
\[
 \tau_{n,u_n}^{\star}\le\frac4\beta\log(n/u_n)
                    +\frac{(\log n)^{9/2}}{u_n}
             \le c_\beta\log n,\qquad c_\beta=4/\beta+1,
\]
for every sufficiently large $n$. These are the only estimates on the stopped array used below.

A lower bound on Brownian oscillations controls the time needed for all terminal exits. The modulus of continuity on this time interval then bounds every position in the terminal
part of a row. Since the radii are shared across rows, the construction gives one event of probability one and a finite $n_0(\omega)$ on which the estimates hold for every
$n\ge n_0(\omega)$.

\begin{lemma}\label{S:lem:auxiliary-boundary-oscillation}
Almost surely, there is a finite integer $N$ such that, for every $n\ge N$,
\[
 \max_{1\le j\le u_n}
 d_{\mathbb H}(z_{n,j}^{\star},z_{n,u_n}^{\star})
                  \le(1+\log\log n)^4.
\]
\end{lemma}

\begin{proof}
The beta tail gives, for $r\ge0$,
\[
 \mathbb P(Y_k\ge r)
 =\mathbb P(\xi_k\ge\tanh^2(r/2))
 =(1-\tanh^2(r/2))^{\beta k/2}
 =\operatorname{sech}^{\beta k}(r/2).
\]
For $x\ge0$, $\tanh x\le x$ and
\[
 (\tanh x-x+x^3/3)'=x^2-\tanh^2x\ge0.
\]
Integration gives $\log\cosh x\ge x^2/2-x^4/12$. Hence $\log\cosh(r/2)\ge r^2/12$ for $0\le r\le1$ and $\mathbb P(Y_k\ge r)\le\exp(-\beta kr^2/12)$ on this range. For all
sufficiently large $k$ we may set $r=\sqrt{36\log(2k)/(\beta k)}\le1$; the resulting probabilities are at most $(2k)^{-3}$. Borel--Cantelli, with the finitely many remaining radii
included in the constant, proves
\begin{equation}\label{S:eq:aux-tail-radius-bound}
 Y_k\le C_\xi\sqrt{\frac{\log(2k)}k},\qquad k\ge1.
\end{equation}
This application of Borel--Cantelli uses only the marginal tail bounds.

In upper half-plane coordinates write $B(t)=x_t+\ii y_t$ and $V_t=\log y_t=W_t^{(2)}-t/2$. For an interval $I$ put $\operatorname{osc}_I V=\sup_{v,w\in I}|V_v-V_w|$. There are
absolute $c,C>0$ such that, for deterministic $s,h,a>0$,
\[
 \mathbb P(\operatorname{osc}_{[s,s+h]}V\le a)
                   \le C\exp(-ch/a^2).                              \tag{**}
\]
If $h/a^2\ge32$, divide the interval into $M=\lfloor h/(16a^2)\rfloor\ge h/(32a^2)$ equal pieces. The increments of $V$ on these pieces are independent normal variables of variance
$h/M\ge16a^2$. Whatever their means, their densities are at most $(2\pi h/M)^{-1/2}$. Each has probability at most $2a/(\sqrt{2\pi}\,4a)<1/4$ to lie in $[-a,a]$. Oscillation at
most $a$ forces all $M$ increments into this interval, so its probability is at most $4^{-M}$. This proves \textup{(**)} in this range; enlarging $C$ covers $h/a^2<32$.

For $l\in\mathbb Z$, $j\in\mathbb N_0$, let
\[
 I_{l,j}=[j2^l,(j+1)2^l],\qquad D_{l,j}=1+|l|+\log(j+2).
\]
Choose a deterministic $b>0$ sufficiently small that $c/b^2\ge4$. By \textup{(**)},
\[
 \mathbb P\left(\operatorname{osc}_{I_{l,j}}V
          \le b\sqrt{2^l/D_{l,j}}\right)
      \le C\me^{-4|l|}(j+2)^{-4}.
\]
The double sum is finite. Also $\mathbb P(\operatorname{osc}_{I_{l,j}}V=0)=0$, since the difference of the boundary values has a nondegenerate normal distribution. Borel--Cantelli
therefore implies
\[
 b_\omega:=\inf_{l\in\mathbb Z,\ j\ge0}
       \sqrt{D_{l,j}/2^l}\,\operatorname{osc}_{I_{l,j}}V>0
 \quad\text{almost surely}:
\]
all but finitely many entries of this infimum are at least $b$, and every exceptional entry is positive. Hence $\operatorname{osc}_{I_{l,j}}V\ge b_\omega\sqrt{2^l/D_{l,j}}$ for
every $l,j$.

For arbitrary $s\ge0,h>0$, choose $l$ with $h/8<2^l\le h/4$ and put $j=\lceil s/2^l\rceil$. Then $I_{l,j}\subset[s,s+h]$, $|l|\le C(1+|\log h|)$ and $j+2\le8s/h+3$. It follows that
\[
 D_{l,j}\le C\{\log(2+s+h)+\log(h+h^{-1})\}.
\]
For example, use $\log(h+h^{-1})\ge|\log h|$ and $8s/h+3\le C(2+s+h)(1+h^{-1})$; the two logarithms on the right also have a fixed positive sum. The hyperbolic metric satisfies
$|\log y_v-\log y_w|\le d_{\mathbb H}(B(v),B(w))$. The triangle inequality through $B(s)$ thus gives
\[
 \operatorname{osc}_{I_{l,j}}V
       \le2\sup_{s\le t\le s+h}d_{\mathbb H}(B(s),B(t)).
\]
Combining the last three estimates proves, on one event, simultaneously for all real $s\ge0,h>0$,
\begin{equation}\label{S:eq:aux-tail-lower-fluctuation}
 \sup_{s\le t\le s+h}d_{\mathbb H}(B(s),B(t))^2
       \ge\frac{c_Bh}{\log(2+s+h)+\log(h+h^{-1})},\qquad c_B>0.
\end{equation}
In particular this bound may be applied at arbitrary random starting times on this event.

We also need an upper bound, for which we give a direct calculation. For an integer $m\ge0$ and $0\le v\le1$, use the isometry
$z\mapsto(z-\operatorname{Re}B(m))/\operatorname{Im}B(m)$, which sends $B(m)$ to $\ii$. The new coordinates are
\[
 \eta_{m,v}=\frac{y_{m+v}}{y_m}
      =\exp(W_{m+v}^{(2)}-W_m^{(2)}-v/2),\qquad
 \xi_{m,v}=\frac{x_{m+v}-x_m}{y_m}
      =\int_0^v\eta_{m,w}\,\dd W_{m+w}^{(1)}.
\]
Set
\[
 V_m^*=\sup_{0\le v\le1}|W_{m+v}^{(2)}-W_m^{(2)}|,\qquad
 Z_m^*=\me^{-V_m^*}\sup_{0\le v\le1}|\xi_{m,v}|.
\]
Given the vertical increments on this interval and the past at $m$, $\xi_{m,v}$ is a centered Gaussian process with covariance $\int_0^{v\wedge w}\eta_{m,a}^2\dd a$. Its law on the
space of continuous paths is that of a real Brownian motion evaluated at these deterministic integrated variances. Their terminal value is at most $\me^{2V_m^*}$. The reflection
principle therefore gives
\[
 \mathbb P(Z_m^*>z\mid\mathcal F_m,
                W_{m+v}^{(2)}-W_m^{(2)},\ 0\le v\le1)
       \le4\me^{-z^2/2},\qquad z>0.
\]
The same principle bounds $\mathbb P(V_m^*>z)\le4\me^{-z^2/2}$. Since $\cosh d_{\mathbb H}(\ii,\xi+\ii)=1+\xi^2/2$, the geodesic from $\ii$ to $\xi+\ii$ has length
$2\operatorname{arsinh}(|\xi|/2)$. Following it by the vertical segment from $\xi+\ii$ to $\xi+\ii\eta$ gives a path of length
\[
 2\operatorname{arsinh}(|\xi|/2)+|\log\eta|
       \le2\log(1+|\xi|)+|\log\eta|.
\]
Since $|\xi_{m,v}|\le\me^{V_m^*}Z_m^*$ and $|\log\eta_{m,v}|\le V_m^*+1/2$, the triangle inequality gives
\[
 \sup_{0\le v\le1}d_{\mathbb H}(B(m),B(m+v))
       \le3V_m^*+2\log(1+Z_m^*)+1/2
       \le3V_m^*+2Z_m^*+1/2.
\]
The two Gaussian tails imply, with absolute constants independent of $m$,
\[
 \mathbb P\left(\sup_{0\le v\le1}
          d_{\mathbb H}(B(m),B(m+v))>R\right)
       \le C\me^{-cR^2},\qquad R\ge2.
\]
Taking $R=A\sqrt{\log(m+2)}$ with $A$ sufficiently large makes the series of these probabilities finite. Borel--Cantelli and continuity on the finitely many exceptional unit
intervals give a finite random constant valid for every integer $m$. Set $D_m=\sup_{m\le u\le m+1}d_{\mathbb H}(B(m),B(u))$. For $s\le t\le s+1$ and $m=\lfloor s\rfloor$, the
triangle inequality through $B(m)$ and, if needed, $B(m+1)$ gives
\[
 d_{\mathbb H}(B(s),B(t))\le2D_m+D_{m+1}
       \le C_B\sqrt{\log(m+3)}.
\]
Taking the supremum over $0\le s\le t\le T$ proves
\begin{equation}\label{S:eq:s4-hyperbolic-unit-modulus}
 \sup_{\substack{0\le s\le t\le T\\t-s\le1}}
       d_{\mathbb H}(B(s),B(t))
       \le C_B\sqrt{\log(3+T)},\qquad T\ge1.
\end{equation}

Fix a path satisfying \eqref{S:eq:aux-tail-radius-bound}, \eqref{S:eq:aux-tail-lower-fluctuation}, \eqref{S:eq:s4-hyperbolic-unit-modulus}, and the stated bound on
$\tau_{n,u_n}^{\star}$. Choose a finite $A\ge1$ with $c_BA/10>C_\xi^2$, and set
\[
 H_{n,k}=AL_n\frac{\log(2k)}k,\qquad 1\le k<u_n.
\]
For $u\ge3$, monotonicity of $\log(2x)/x$ on $x\ge2$ gives
\[
 \sum_{k=1}^{u-1}\frac{\log(2k)}k
 \le\log2+\frac{\log4}{2}+\int_2^u\frac{\log(2x)}x\,\dd x
 \le C(1+\log u)^2.
\]
Increasing $C$ includes $u=1,2$. Since $u_n\le2(\log n)^6$ for large $n$,
\begin{equation}\label{S:eq:aux-tail-time-budget}
 \sum_{k=1}^{u_n-1}H_{n,k}
       \le CAL_n(1+\log u_n)^2
       \le CAL_n^3=o(\log n).
\end{equation}
Also $H_{n,k}\le AL_n$ and $H_{n,k}^{-1}\le u_n/(AL_n\log2)$. It follows, uniformly over $k<u_n$ and $0\le s\le2c_\beta\log n$, that
\[
 \log(2+s+H_{n,k})\le2L_n,\qquad
 \log(H_{n,k}+H_{n,k}^{-1})\le7L_n
\]
for every sufficiently large $n$. For the second bound, its logarithm is at most $\log(AL_n+2(\log n)^6/(AL_n\log2))=6L_n+O_{A}(\log L_n+1)$. The first follows from
$\log(C_\beta\log n+AL_n+2)=L_n+O_{A,\beta}(1)$. The random $A$ is fixed before $n$ tends to infinity.

For $k<u_n$, descending induction gives $\tau_{n,k}^{\star}-\tau_{n,k+1}^{\star}\le H_{n,k}$. If it holds for $k<j<u_n$, then \eqref{S:eq:aux-tail-time-budget} and the bound on the
initial time give
\[
 \tau_{n,k+1}^{\star}\le c_\beta\log n+\sum_{j=k+1}^{u_n-1}H_{n,j}
                   \le2c_\beta\log n.
\]
This includes the first induction step, where the sum is empty. Apply \eqref{S:eq:aux-tail-lower-fluctuation} with $s=\tau_{n,k+1}^{\star}$ and $h=H_{n,k}$. The denominator is at
most $10L_n$, so
\[
 \sup_{\tau_{n,k+1}^{\star}\le t\le\tau_{n,k+1}^{\star}+H_{n,k}}
 d_{\mathbb H}(B(\tau_{n,k+1}^{\star}),B(t))^2
       \ge\frac{c_BA}{10}\frac{\log(2k)}k>Y_k^2.
\]
Since the distance exceeds $Y_k$ before time $s+H_{n,k}$, continuity gives $\tau_{n,k}^{\star}\le s+H_{n,k}$. Summing these bounds gives
\begin{equation}\label{S:eq:aux-tail-total-time}
 S_n:=\tau_{n,1}^{\star}-\tau_{n,u_n}^{\star}\le C_\omega L_n^3,\qquad
 \tau_{n,1}^{\star}\le2c_\beta\log n.
\end{equation}

Divide $[\tau_{n,u_n}^{\star},\tau_{n,1}^{\star}]$ into at most $1+S_n$ intervals of length at most one. Applying \eqref{S:eq:s4-hyperbolic-unit-modulus} to each piece and summing
their lengths in the hyperbolic metric gives
\[
 \begin{split}
 \sup_{\tau_{n,u_n}^{\star}\le t\le\tau_{n,1}^{\star}}
        d_{\mathbb H}(B(\tau_{n,u_n}^{\star}),B(t))
 &\le C_B(1+S_n)\sqrt{\log(3+\tau_{n,1}^{\star})}\\
 &\le C_\omega(1+L_n)^{7/2}\le(1+L_n)^4
 \end{split}
\]
for all sufficiently large $n$. The last inequality holds once $(1+L_n)^{1/2}\ge C_\omega$. Thus the coefficient in the last bound is deterministic. Every $\tau_{n,j}^{\star}$ with
$1\le j\le u_n$ lies in this time interval, which proves the lemma.
\end{proof}

We finally check the use of this estimate under terminal conditioning. Let $\eta=B(\infty)$ and let $\mathsf Q$ be the joint law of $\eta$ and the entire countable array
$(z_{n,j}^{\star})$. For a sufficiently large deterministic $n_0$, the event in the lemma is
\[
 E=\bigcup_{N\ge n_0}\ \bigcap_{n\ge N}
     \left\{\max_{1\le j\le u_n}
       d_{\mathbb H}(z_{n,j}^{\star},z_{n,u_n}^{\star})
                                      \le(1+\log\log n)^4\right\}.
\]
It is a Borel event of the countable array and $\mathsf Q(E)=1$. The boundary and array spaces are standard Borel, so disintegration gives
$\mathsf Q(\dd\eta,\dd z)=\nu_\ii(\dd\eta)\mathsf Q^\eta(\dd z)$. Since the integrand is nonnegative,
\[
 0=\mathsf Q(E^c)=\int\mathsf Q^\eta(E^c)\nu_\ii(\dd\eta)
 \quad\Longrightarrow\quad
 \mathsf Q^\eta(E)=1\quad\hbox{for }\nu_\ii\hbox{-almost every }\eta.
\]
Proposition~\ref{S:prop:7.2} samples this conditional law as a single array, so the same event and starting index apply to all its rows.

For each row let $h_n$ be the isometry fixing $\eta$ and mapping $z_{n,u_n}^{\star}$ to the point at which the auxiliary path is attached. A measurable choice is obtained by
setting
\[
 g_\eta(z)=-\frac1{z-\eta}\quad(\eta\in\mathbb R),\qquad g_\infty(z)=z.
\]
This is an orientation-preserving hyperbolic isometry sending $\eta$ to infinity. If the images of the two points are $x+\ii y$ and $x'+\ii y'$, put
\[
 h_n(z)=g_\eta^{-1}\!\left(\frac{y'}y g_\eta(z)+x'-\frac{y'}y x\right).
\]
The formula is Borel measurable in the points and $\eta$. It fixes $\eta$ and has the required value at the first point. Every orientation-preserving isometry fixing infinity is
$w\mapsto aw+b$ with $a>0,b\in\mathbb R$; the required image forces $a=y'/y$ and $b=x'-(y'/y)x$. This proves uniqueness as well. Distance invariance yields
\[
 d_{\mathbb H}(h_nz_{n,j}^{\star},h_nz_{n,u_n}^{\star})
       =d_{\mathbb H}(z_{n,j}^{\star},z_{n,u_n}^{\star}).
\]
Hence the same bound holds for the attached auxiliary path and implies \eqref{S:eq:auxiliary-tail-oscillation} for any fixed $C_\beta^{\mathrm{tail}}\ge1$. Identification of its
conditional distribution for the circular $\beta$-ensemble is the separate concatenation argument in Proposition~\ref{S:prop:7.2}. The operator interval with index $u_n$ uses the
position with index $u_n+1$. Proposition~\ref{S:prop:finite-extreme-tail} estimates the retained increment between the positions with indices $u_n+1$ and $u_n$.

\subsection{The error from sampling the Brownian path}
\label{sec:brownian-bridge-expansion-and-the-sampling-remainder}
\label{S:sec:the-error-from-sampling-the-brownian-path}

The remaining discretization error comes from replacing the Brownian path by its values at the deterministic times $\vartheta_{n,k+1}$. For fixed $\beta>0$, we work with the
Brownian marginal of Proposition~\ref{S:prop:7.2}. The sampling error supplies the leading order coefficient, while
Section~\ref{S:sec:comparison-of-the-finite-and-sampled-brownian-operators} treats the smaller error in the comparison with the circular $\beta$-ensemble.

For $r\in I_{n,k}$, recall that the sampled boundary columns are $a_{n,k}^{\mathrm{br}}=\mathcal A(Z_{\vartheta_{n,k+1}})u_0$ and
$c_{n,k}^{\mathrm{br}}=\mathcal A(Z_{\vartheta_{n,k+1}})u_1$. In the reflected coordinate, the sampled inverse has kernel
\[
 K_n^{\mathrm{br}}(r,r')=\tfrac12
 a_{n,k}^{\mathrm{br}}(c_{n,l}^{\mathrm{br}})^{\mathsf T},
 \qquad r\in I_{n,k},\quad r'\in I_{n,l},\quad r>r'.
\]
The other triangle is defined by transposition. This formula applies also when $k=l$; the diagonal $r=r'$ is irrelevant for the integral operator. Replacing these sampled columns
by $a_\omega(r)$ and $c_\omega(r')$ gives the inverse of the limiting operator. Their existence and boundary conditions were established in
Theorem~\ref{S:thm:sine-selfadjoint-domain}. We first compress both kernels to $(m_n/n,1)^2$, without changing either boundary vector.

Conditioning on the sampled positions at a fixed matrix size and on the terminal point fixes these columns. In normal coordinates, the path between consecutive observations is
decomposed into the linear interpolation between its two endpoints, a centered Euclidean Brownian bridge, and a remainder. The linear interpolation term contributes to the
conditional mean and must be retained. The derivative of the inverse is bounded from the field space with squared norm $\mathcal E_n^+$ to $\mathfrak S_2$ by $n^{-1/2}$. We shall
prove that the linear interpolation term and the bridge have bounded $\mathcal E_n^+$, while the remainder has vanishing $\mathcal E_n^+$.

Throughout, the matrix factor is $\mathcal A$ and the boundary vectors are $u_0=(-1,0)^{\mathsf T}$ and $u_1=(\eta_\infty,1)^{\mathsf T}$. We construct the Euclidean bridges on an
extension of the existing probability space. The Brownian path and the circular rows keep their original joint law. The expansion of the inverse is proved in
Proposition~\ref{S:prop:continuous-bridge-linearization}.

Use the notation $\vartheta_{n,j}$, $h_{n,j}$, $J_{n,j}$, $\theta_{n,k}(s)$ and $\lambda_{n,k,j}(s)$ from
Section~\ref{S:sec:sampling-times-and-estimates-for-the-boundary-solutions}. The operator intervals have indices $m_n\le k<n$. By \eqref{S:eq:additional-sampling-intervals}, all
intervals between consecutive Brownian observation times traversed from $\vartheta_{n,k+1}$ to $\theta_{n,k}(s)$ have indices in $\{\ell_n,\ldots,n-1\}$, where
\begin{equation}\label{S:eq:sampling-interval-lower-index}
 \ell_n:=m_n-R_\beta\ge m_n/2,\qquad
 \ell_n\sim(\log n)^{512}.
\end{equation}
These inequalities and all constructions on the selected interval below are understood for sufficiently large deterministic $n$. For each $k$ at most $2R_\beta+1$ indices $j$
occur, and $h_{n,j}\le C_\beta/m_n$ throughout this enlarged set.

Write $X_\theta=\mathcal A(Z_\theta)$ and $X_{n,k}^{\mathrm{step}}=X_{\vartheta_{n,k+1}}$. The relative matrix is upper triangular with positive diagonal. Consequently its real
upper triangular logarithm is defined, and we set
\[
 \mathcal U_{n,k}(s)
 :=\log\bigl(X_{\theta_{n,k}(s)}
                    (X_{n,k}^{\mathrm{step}})^{-1}\bigr),
 \qquad r=(k+s)/n,\quad 0\le s\le1.
\]
In particular, $\exp(\mathcal U_{n,k}(s))X_{n,k}^{\mathrm{step}} =X_{\theta_{n,k}(s)}$ and $\operatorname{tr}\mathcal
U_{n,k}(s)=0$.

Put $U_n=\mathcal A(Z_{\vartheta_{n,1}})\cdot\eta_\infty$. Conditionally on the path through $\vartheta_{n,1}$, the terminal point has harmonic measure from $Z_{\vartheta_{n,1}}$.
The isometry represented by $\mathcal A(Z_{\vartheta_{n,1}})$ sends $Z_{\vartheta_{n,1}}$ to $\ii$. Thus $U_n$ has law $\nu_\ii$, independently of this past. Define
\begin{equation}\label{S:eq:sample-chain-sigma-field}
 \mathscr S_n
 :=\sigma(Z_{\vartheta_{n,j}}:1\le j\le n)\vee\sigma(U_n).
\end{equation}
The identity $\eta_\infty=\Re Z_{\vartheta_{n,1}} +\Im Z_{\vartheta_{n,1}}\,U_n$ has a measurable inverse. Hence
\begin{equation}\label{S:eq:row-conditioning-global-terminal}
 \mathscr S_n
 =\sigma(Z_{\vartheta_{n,j}}:1\le j\le n)\vee\sigma(\eta_\infty).
\end{equation}
The sampled matrices, boundary columns and weights are $\mathscr S_n$-measurable. We condition on one row at a time; Lemma~\ref{S:lem:conditional-all-skeletons} identifies the
larger sigma-field obtained by observing all rows.

Conditioning on the terminal point fixes the boundary vector. Once both finite boundary points of a bridge are specified, the Poisson factors in the conditional transition
densities cancel, leaving the usual bridge law. The Euclidean bridges are conditionally independent within each row. All rows retain their dependence through the original Brownian
path.

All regular conditional distributions are taken on Borel path spaces. A filtration is augmented only under the probability law used for the stochastic integral. Write
$\mathcal F_t^Z=\sigma(Z_s:0\le s\le t)$ before augmentation, $L=\tfrac12\Delta_{\mathbb H}$, and let $p_t$ denote its transition density with respect to hyperbolic area. We use
$e_1=y\partial_x$, $e_2=y\partial_y$ and $b_\eta^\mu=e_\mu\log P_\eta$.

\subsubsection{Conditional laws of the short bridges}
We construct the short bridges under the terminal conditional law, which changes the Brownian drift.

\begin{lemma}\label{S:lem:terminal-affine-innovations}
For $\nu_\ii$-almost every $\eta$, the terminal conditional law $\mathbb P^\eta$ satisfies, on every finite horizon,
\begin{equation}\label{S:eq:terminal-path-density}
 \left.\frac{\dd\mathbb P^\eta}{\dd\mathbb P}\right|_{\mathcal F_t^Z}
 =P_\eta(Z_t),
 p_t^\eta(z,z')=p_t(z,z')\frac{P_\eta(z')}{P_\eta(z)}.
\end{equation}
Moreover,
\begin{equation}\label{S:eq:busemann-affine-bounded}
 |b_\eta(z)|=1,
 L^\eta=L+\nabla\log P_\eta\cdot\nabla,
 \widetilde W_t^\eta=W_t-\int_0^t b_\eta(Z_s)\,\dd s
\end{equation}
Here $\widetilde W^\eta$ is a standard two-dimensional Brownian motion under $\mathbb P^\eta$.
\end{lemma}
\begin{proof}
The harmonic-measure disintegration in Lemma~\ref{S:lem:7.2} identifies the terminal conditional density with $P_\eta(Z_t)$ on each finite horizon. We verify the Brownian change of
drift and its consistency in the present coordinates. For a finite terminal point $\eta\in\mathbb R$, the Poisson kernel normalized at $\ii$ is
\[
 P_\eta(x+\ii y)=\frac{(1+\eta^2)y}{(x-\eta)^2+y^2}.
\]
For the boundary point at infinity it is $P_\infty(x+\ii y)=y$. Write $D=(x-\eta)^2+y^2$. Differentiation gives, for finite $\eta$,
\[
 b_\eta^1=-\frac{2y(x-\eta)}D,
 \qquad b_\eta^2=\frac{(x-\eta)^2-y^2}D,
 \qquad (b_\eta^1)^2+(b_\eta^2)^2=1.
\]
The identity for the last sum is $4y^2(x-\eta)^2+((x-\eta)^2-y^2)^2=D^2$. At infinity the vector is $(0,1)$. These are the coefficients in
Proposition~\ref{S:prop:aux-stochastic-normalization}.

Since $LP_\eta=0$, It\^o's formula along the original diffusion gives
\[
 \mathrm d\log P_\eta(Z_t)
       =b_\eta(Z_t)\cdot\mathrm dW_t-\tfrac12\mathrm dt.
\]
Here $L\log P_\eta=-|\nabla\log P_\eta|^2/2=-1/2$ and $P_\eta(\ii)=1$. Thus
\[
 P_\eta(Z_t)
 =\exp\left(\int_0^t b_\eta(Z_s)\cdot\mathrm dW_s
                   -\frac12\int_0^t|b_\eta(Z_s)|^2\,\mathrm ds\right).
\]
For every finite $T$, $\mathbb E\exp(\frac12\int_0^T|b_\eta(Z_s)|^2\,\mathrm ds)=e^{T/2}<\infty$. The exponential is therefore a true martingale on $[0,T]$, and Girsanov's theorem
makes
\[
 \widetilde W_t^\eta=W_t-\int_0^t b_\eta(Z_s)\,\mathrm ds
\]
a standard two-dimensional Brownian motion under this changed law. The density martingale makes the laws on integer horizons consistent. The disintegration identifies these laws
with $\mathbb P^\eta$ for almost every terminal value. Intersecting the full-measure terminal sets for the integer horizons gives a single such set on which every finite-time
conclusion holds. No absolute continuity on the infinite-time path sigma-field is required.

For $s<t$ and bounded Borel $f$, Bayes' formula and the original Markov property give
\[
 \begin{split}
 \mathbb E^\eta[f(Z_t)\mid\mathcal F_s^Z]
 &=\frac{\mathbb E[P_\eta(Z_t)f(Z_t)\mid\mathcal F_s^Z]}
           {P_\eta(Z_s)}\\
 &=\int_{\mathbb H}f(z')p_{t-s}(Z_s,z')
               \frac{P_\eta(z')}{P_\eta(Z_s)}\,\mathrm d\mathrm{vol}(z').
 \end{split}
\]
This proves the transition formula. Finally, for a smooth compactly supported $f$, the diffusion product rule and $LP_\eta=0$ yield
\[
 L^\eta f=P_\eta^{-1}L(P_\eta f)
          =Lf+\nabla\log P_\eta\cdot\nabla f.
\]
These formulas use the original Brownian driving coordinates of the increment comparisons below.

In the original upper half-plane coordinates, this reads
\begin{equation}\label{S:eq:conditioned-original-affine-sde}
 \dd x_t=y_t\,\dd\widetilde W_t^{\eta,1}
                     +y_tb_\eta^1(Z_t)\,\dd t,\qquad
 \dd y_t=y_t\,\dd\widetilde W_t^{\eta,2}
                     +y_tb_\eta^2(Z_t)\,\dd t.
\end{equation}

For the metric $g=y^{-2}(\dd x^2+\dd y^2)$ the nonzero Christoffel symbols are
\[
 \Gamma^x_{xy}=\Gamma^x_{yx}=-y^{-1},\qquad
 \Gamma^y_{xx}=y^{-1},\qquad \Gamma^y_{yy}=-y^{-1}.
\]
If $f=\log y$, then $\partial_xf=0$, $\partial_yf=y^{-1}$, and $\partial_y^2f=-y^{-2}$. Hence
\[
 (\operatorname{Hess}f)_{xx}=-y^{-2},\qquad
 (\operatorname{Hess}f)_{xy}=0,\qquad
 (\operatorname{Hess}f)_{yy}=0,
\]
which is $\operatorname{Hess}f=-g+\dd f\otimes\dd f$. An isometry taking $\eta$ to infinity gives
\begin{equation}\label{S:eq:s9-poisson-hessian}
 |\dd\log P_\eta|=1,\qquad
 \nabla^2\log P_\eta=-g+
       \dd\log P_\eta\otimes\dd\log P_\eta.
\end{equation}
Put $\alpha=\dd\log P_\eta$ and $T_m=\nabla^m\alpha$. Equation~\eqref{S:eq:s9-poisson-hessian} gives $|T_0|=1$ and $\|T_1\|_{\mathrm{op}}=1$. For every $m\ge1$, differentiate
$T_1=-g+\alpha\otimes\alpha$ a further $m$ times. Since $\nabla g=0$, the Leibniz rule and the tensor operator norm give
\[
 \|T_{m+1}\|_{\mathrm{op}}
 \le\sum_{j=0}^m\binom mj
       \|T_j\|_{\mathrm{op}}\|T_{m-j}\|_{\mathrm{op}}.
\]
There are $\binom mj$ choices of the $j$ differentiation slots falling on the first factor; permutations of slots preserve the norm. With $C_0=C_1=1$ and
$C_{m+1}=\sum_{j=0}^m\binom mj C_jC_{m-j}$, induction gives $\sup_{z,\eta}\|\nabla^m\dd\log P_\eta(z)\|_{\mathrm{op}}\le C_m$. In particular $C_2=2$ bounds the third derivative of
$\log P_\eta$. In a fixed normal ball, ordinary derivatives are obtained from covariant derivatives by repeated use of
$\partial_i\alpha_j=(\nabla_i\alpha)_j+\Gamma^k_{ij}\alpha_k$. The Christoffel symbols and their derivatives are bounded in that chart. Induction expresses each coordinate
derivative of $\alpha$ as a finite sum of covariant derivatives $T_j$ multiplied by Christoffel symbols and their coordinate derivatives. The bounds $\|T_j\|\le C_j$ and the chart
bounds therefore bound every such derivative. Isometries preserve these constants when the center changes.
\end{proof}

Let $\mathsf P_h^{z,z'}$ be the duration-$h$ Brownian bridge law on $C([0,1];\mathbb H)$, after rescaling time to $[0,1]$, and put
\[
 \gamma_{n,j}(u)=Z_{\vartheta_{n,j+1}+h_{n,j}u},\qquad 0\le u\le1.
\]

\begin{lemma}\label{S:lem:bridge-disintegration}
The laws $\mathsf P_h^{z,z'}$ form a Borel probability kernel. For each $n$,
\begin{equation}\label{S:eq:row-bridge-product-law}
 \mathcal L((\gamma_{n,j})_{1\le j<n}\mid\mathscr S_n)
 =\bigotimes_{j=1}^{n-1}
   \mathsf P_{h_{n,j}}^{Z_{\vartheta_{n,j+1}},Z_{\vartheta_{n,j}}}.
\end{equation}
The finite bridge law is unchanged by conditioning on the terminal point.

\end{lemma}
\begin{proof}
Use Proposition and Definition~2.1 and Remark~2.2 of \cite{Guneysu2019}. For $0=u_0<u_1<\cdots<u_d<u_{d+1}=1$, $x_0=z$, and $x_{d+1}=z'$, the density of
$(\gamma(u_1),\ldots,\gamma(u_d))$ under $\mathsf P_h^{z,z'}$, relative to hyperbolic volume in the $d$ interior variables, is
\[
 \frac{\prod_{a=0}^d p_{h(u_{a+1}-u_a)}(x_a,x_{a+1})}{p_h(z,z')}.
\]
It is Borel in $(h,z,z')$, since $p$ is smooth and positive. Cylinder sets at rational times generate the Borel sigma-field of $C([0,1];\mathbb H)$. The collection of sets whose
bridge probabilities are Borel is closed under complements and countable disjoint unions; the monotone class theorem therefore proves kernel measurability.

For $p_t^\eta(v,w)=p_t(v,w)P_\eta(w)/P_\eta(v)$, the density ratio is unchanged, because
\[
 \begin{aligned}
 \frac{\prod_{a=0}^d p_{h(u_{a+1}-u_a)}^\eta(x_a,x_{a+1})}{p_h^\eta(z,z')}
 &=\frac{\prod_{a=0}^d p_{h(u_{a+1}-u_a)}(x_a,x_{a+1})}{p_h(z,z')}
   \frac{\prod_{a=0}^d P_\eta(x_{a+1})/P_\eta(x_a)}{P_\eta(z')/P_\eta(z)}
 \\
 &=\frac{\prod_{a=0}^d p_{h(u_{a+1}-u_a)}(x_a,x_{a+1})}{p_h(z,z')}.
 \end{aligned}
\]
Fix a row and list its sampling times in chronological order as $0=t_0<t_1<\cdots<t_{n-1}$. In each interval $[t_{l-1},t_l]$, insert finitely many rational relative observation
times, and let $D_l$ be the product of the transition densities across the resulting subintervals. Under the terminal conditional law the joint density of all sampled and inserted
positions is $\prod_lD_l^\eta$. The marginal density of the sampled positions is $\prod_lp_{t_l-t_{l-1}}^\eta(z_{l-1},z_l)$. Their ratio factorizes as
\[
 \prod_{l=1}^{n-1}
       \frac{D_l^\eta}{p_{t_l-t_{l-1}}^\eta(z_{l-1},z_l)}
 =\prod_{l=1}^{n-1}
       \frac{D_l}{p_{t_l-t_{l-1}}(z_{l-1},z_l)},
\]
where the equality follows by cancellation of the Poisson kernel factors, as proved above. Every denominator is positive. Integrating the inserted coordinates against bounded
cylinder functionals therefore gives
\[
 \mathbb E\left[\prod_{j=1}^{n-1}f_j(\gamma_{n,j})
                         \,\middle|\,\mathscr S_n\right]
 =\prod_{j=1}^{n-1}\int f_j(\gamma)
    \,\mathsf P_{h_{n,j}}^{Z_{\vartheta_{n,j+1}},Z_{\vartheta_{n,j}}}
                                    (\mathrm d\gamma).
\]
To verify this as a conditional identity, multiply both sides by any bounded function of the sampled positions and the terminal point and integrate their joint law. The density
factorization gives equality of these integrals. Equation~\eqref{S:eq:row-conditioning-global-terminal} identifies the resulting conditioning sigma-field with $\mathscr S_n$.

For each path variable choose the countable cylinder algebra generated by rational relative times and a countable basis of the hyperbolic plane. Fix versions on which the displayed
equality holds for all finite products from these algebras. For fixed choices of the other $f_i$, the bounded functions $f_j$ satisfying the equality form a vector space closed
under bounded monotone limits. The monotone class theorem extends the equality to every bounded Borel $f_j$. Repeating this for the finitely many path variables identifies the
conditional law with the product in \eqref{S:eq:row-bridge-product-law}. Thus independence concerns the whole path segments, not only any particular finite set of observations.

For nonnegative functionals and nonnegative $\mathscr S_n$-measurable weights, first truncate all factors and then use monotone convergence. A signed identity follows by
subtraction whenever its conditional absolute integral is finite. The countable cylinder algebras and the countable set of rows permit a single exceptional null set to be removed.
This does not assert that Gaussian bridge laws remain unchanged after conditioning on the union of the observation sigma-fields from all rows.
\end{proof}

Let $C_0([0,1];\mathbb R^2)$ denote continuous paths vanishing at both boundary points.

\begin{lemma}\label{S:lem:pinned-gaussian-kernel}
There is a Borel kernel $\mathsf K_{h,z,z'}$ coupling $\mathsf P_h^{z,z'}$ with a standard two-dimensional Euclidean bridge using the same driving Brownian motion. The same kernel
is used for all parameters in Lemma~\ref{S:lem:short-bridge-expansion}.
\end{lemma}
\begin{proof}
Express $Y_t=\exp_z^{-1}(Z_t)$ in the orthonormal basis $(e_1(z),e_2(z))$. The normal chart is global on $\mathbb H$. Let $a$ be its inverse metric, $\sigma=a^{1/2}$ its positive
square root, and $b_0^i=-\tfrac12g^{jk}\Gamma^i_{jk}$. For $t<h$, the bridge equation is
\begin{equation}\label{S:eq:pinned-kernel-sde}
\dd Y_t=\sigma(Y_t)\,\dd W_t^{\mathrm{pin}}
 +\left[b_0(Y_t)+a(Y_t)\nabla\log
          p_{h-t}(\exp_zY_t,z')\right]\dd t.
\end{equation}
Here $\nabla$ denotes the Euclidean gradient in the first normal coordinate variable. To verify the drift, set $H(t,x)=p_{h-t}(\exp_zx,z')$ and let $L$ act in this chart. Since
$(\partial_t+L)H=0$, the diffusion product rule gives
\[
 H^{-1}\bigl((\partial_t+L)(Hf)-f(\partial_t+L)H\bigr)
 =Lf+\langle a\nabla\log H,\nabla f\rangle.
\]
The density $H(t,Y_t)/p_h(z,z')$ and the bridge formula of Lemma~\ref{S:lem:bridge-disintegration} therefore imply this martingale problem on every compact set and
$t\le h-\epsilon$. Apply it to cutoff coordinate functions and their pairwise products. The martingale part $M$ of $Y$ has $\dd[M^i,M^j]_t=a^{ij}(Y_t)\dd t$. Thus
$W^{\mathrm{pin}}=\int\sigma(Y_t)^{-1}\dd M_t$ has identity quadratic covariation. Its exponential characteristic martingales have bounded modulus on finite time intervals and give
\[
 \mathbb E\left[\me^{\ii\lambda\cdot
       (W_t^{\mathrm{pin}}-W_s^{\mathrm{pin}})}
           \mid\mathcal F_s\right]
 =\me^{-|\lambda|^2(t-s)/2}.
\]
This identifies a two-dimensional Brownian motion and proves weak existence of~\eqref{S:eq:pinned-kernel-sde} with the required bridge law. The semimartingale property at the
terminal time follows from Theorem~2.8 of \cite{Guneysu2019}: the plane is connected and geodesically complete, and the bridge is Markov in its natural filtration by the preceding
density formula.

Normalize time to $[0,1]$. On a compact parameter set with $0<h_-\le h\le h_+$, a fixed ball in normal coordinates, and $0\le u\le u_*:=1-1/r$, the coefficients in
\[
 \dd Y_{hu}=\sqrt h\,\sigma(Y_{hu})\dd\mathcal W_u
 +h\left[b_0(Y_{hu})+a(Y_{hu})\nabla\log
        p_{h(1-u)}(\exp_zY_{hu},z')\right]\dd u
\]
are smooth. After a smooth cutoff in the normal coordinates they are bounded and Lipschitz in the coordinate variables and time, with a common constant $L$. Denote the resulting
coefficients by $b_{\mathfrak p}, \sigma_{\mathfrak p}$, where $\mathfrak p=(h,z,z')$. Starting with $Y^{(0)}=0$, the Picard iterates satisfy
\[
 Y_t^{(j+1)}=\int_0^t b_{\mathfrak p}(v,Y_v^{(j)})\dd v
       +\int_0^t\sigma_{\mathfrak p}(v,Y_v^{(j)})\dd\mathcal W_v.
\]
Cauchy--Schwarz and the $L^2$ martingale maximal inequality give
\[
 d_j(t):=\mathbb E\sup_{v\le t}|Y_v^{(j+1)}-Y_v^{(j)}|^2
 \le C_L\int_0^t d_{j-1}(v)\dd v,
 \qquad d_0(t)\le C_Lt.
\]
Iteration yields $d_j(u_*)\le C_L^{j+1}u_*^{j+1}/(j+1)!$, and hence $\sum_j d_j(u_*)^{1/2}<\infty$. The iterates converge in $L^2$ of the supremum norm to a solution. For two
solutions with the same driver and initial value, let $d(t)=\mathbb E\sup_{u\le t}|Y_u-\widehat Y_u|^2$ under the stopped law. Subtraction gives $d(t)\le C_L\int_0^t d(u)\,\dd u$.
Gronwall's inequality gives $d(t)=0$.

For measurability use continuous Euler approximations on a deterministic partition with maximal interval length $\Delta$. If $t_i\le t<t_{i+1}$, set
\[
 Y_t^\Delta=Y_{t_i}^\Delta
  +b_{\mathfrak p}(t_i,Y_{t_i}^\Delta)(t-t_i)
  +\sigma_{\mathfrak p}(t_i,Y_{t_i}^\Delta)
                         (\mathcal W_t-\mathcal W_{t_i}).
\]
This is a Borel function of $(\mathfrak p,\mathcal W)$. Boundedness of the coefficients implies $\mathbb E|Y_t^\Delta-Y_{t_i}^\Delta|^2\le C_L\Delta$. For
$D(t)=\mathbb E\sup_{v\le t}|Y_v-Y_v^\Delta|^2$, subtract the two integral equations. The drift term is bounded by $t\int_0^t\mathbb E|b(v,Y_v)-b(t_i,Y_{t_i}^\Delta)|^2\dd v$; the
martingale term is bounded by four times the corresponding integral for $\sigma$. Lipschitz continuity gives
\[
 D(t)\le C_L\int_0^tD(v)\dd v+C_L\Delta.
\]
Gronwall proves
\begin{equation}\label{S:eq:estimate-measurable-euler-error}
 \mathbb E\sup_{t\le u_*}|Y_t-Y_t^\Delta|^2\le C_L\Delta.
\end{equation}
For $\Delta_j=2^{-4j}$, the probability that the supremum error exceeds $2^{-j}$ is at most $C_L2^{-2j}$. The approximations therefore converge almost surely for each fixed
parameter. Their uniform Cauchy criterion is Borel in the parameter and path; define their limit there and a fixed continuous path elsewhere. This gives a Borel solution map. For
any weak solution with driver $\mathcal W$, the Euler estimate gives $Y^{\Delta_j}(\mathfrak p,\mathcal W)\to Y$ in probability in the supremum norm. The Borel limit above is a
function of that same driver and has the same limit. Every weak solution therefore has the law of this function of Wiener measure.

Exhaust the parameter space by compact sets, space by balls, and time by $1-1/r$. These are countable exhaustions. Pathwise uniqueness identifies the stopped maps on their overlaps
outside a Wiener null set for each parameter. Patch the maps in a fixed order. Uniqueness in law identifies every such restriction with the corresponding restriction of the
continuous bridge. It follows, by the countable exhaustion, that the patched path has no explosion before one and converges to $z'$ at one. Both properties are Borel conditions on
rational path evaluations. On the exceptional set assign the path $u\mapsto\exp_z(u\exp_z^{-1}z')$. The resulting Borel map has law $\mathsf P_h^{z,z'}$ for every fixed
$(h,z,z')$.

Using the same normalized Brownian motion, define
\begin{equation}\label{S:eq:pinned-euclidean-bridge-map}
 \beta(u)=(1-u)\int_0^u\frac{\dd\mathcal W_v}{1-v},
 \quad u<1,\qquad \beta(1)=0.
\end{equation}
For $u<1$, integration by parts writes the right side as $\mathcal W_u-(1-u)\int_0^u(1-v)^{-2}\mathcal W_v\dd v$, which is a Borel path functional. For $u\le v<1$ its covariance is
\[
 (1-u)(1-v)\int_0^u(1-r)^{-2}\dd r=u(1-v).
\]
Its finite-dimensional laws agree with those of $\mathcal W_u-u\mathcal W_1$, which has continuous paths vanishing at one. Equality on rational times and continuity on $[0,1)$
therefore give the continuous extension of \eqref{S:eq:pinned-euclidean-bridge-map}. In the coupling, $\beta$ is the stochastic integral \eqref{S:eq:pinned-euclidean-bridge-map}
driven by the same $\mathcal W$ as the pinned path. Equality in law with $\mathcal W_u-u\mathcal W_1$ proves that $\beta$ is a standard Euclidean bridge.

The joint image of Wiener measure under these two Borel maps is $\mathsf K_{h,z,z'}$. In particular,
\[
 v=h^{-1/2}\exp_z^{-1}z',\qquad
 \rho_h(u)=h^{-1}\bigl(\exp_z^{-1}\gamma(u)
                         -\sqrt h(uv+\beta(u))\bigr)
\]
is a measurable continuous path, with $\rho_h(0)=\rho_h(1)=0$. It is defined also when $\gamma$ exits a fixed normal ball. This same joint law is used for every localization radius
and moment order below.
\end{proof}

\subsubsection{The joint law of the bridge pairs}
\label{S:sec:joint-bridge-law}

The preceding coupling uses a common driver for fixed boundary points. We now apply it to the segments of the Brownian path already used in the circular coupling. We preserve the
existing hyperbolic path by disintegrating the joint bridge law with respect to that path and sampling the Euclidean bridge from the resulting conditional law. Measurability in the
duration and boundary points permits this construction on all sampling intervals.

For $\mathfrak p=(h,z,z')$, write $\mathsf P_{\mathfrak p}=\mathsf P_h^{z,z'}$ and $\mathsf K_{\mathfrak p}=\mathsf K_{h,z,z'}$. Let $\nu$ be any probability law on the parameter
space and let $(\mathfrak p,\gamma)$ have law
$\nu(\dd\mathfrak p)\mathsf P_{\mathfrak p}(\dd\gamma)$.

\begin{proposition}\label{S:prop:aux-probability-kernels}
There is a Borel probability kernel $\mathsf L(\mathfrak p,\gamma;\dd\beta)$ such that
\[
 \nu(\dd\mathfrak p)\mathsf P_{\mathfrak p}(\dd\gamma)
       \mathsf L(\mathfrak p,\gamma;\dd\beta)
 =\nu(\dd\mathfrak p)\mathsf K_{\mathfrak p}(\dd\gamma,\dd\beta).
\]
It can be sampled using an independent uniform variable without changing the joint law of any previously constructed variables.
\end{proposition}
\begin{proof}
The parameter space and both continuous path spaces are standard Borel. Apply Theorem~8.5 of \cite{kallenberg2021foundations} to
$\nu(\dd\mathfrak p)\mathsf K_{\mathfrak p}(\dd\gamma,\dd\beta)$ conditioned on $(\mathfrak p,\gamma)$. Lemma~\ref{S:lem:pinned-gaussian-kernel} identifies its first marginal,
giving the stated identity; a countable determining algebra makes it an identity of kernels for $\nu$-almost every $\mathfrak p$. By Lemma~4.22 of \cite{kallenberg2021foundations},
there is a Borel map $\Phi$ such that, for every bounded Borel $f$,
\[
 \int_0^1 f(\Phi(\mathfrak p,\gamma,u))\,\dd u
 =\int f(\beta)\,\mathsf L(\mathfrak p,\gamma;\dd\beta).
\]
Let $\lambda$ be Lebesgue probability measure on $[0,1]$ and extend an original probability space $(\Omega,\mathcal A,\mathbb P)$ to $(\Omega\times[0,1],\mathbb P\otimes\lambda)$.
Set $\beta(\omega,u)=\Phi(\mathfrak p(\omega),\gamma(\omega),u)$. Fubini gives, for every bounded Borel function $F$ of the original variables,
\[
 \widetilde{\mathbb E}[F(\omega)f(\beta)]
 =\mathbb E\!\left[F(\omega)\int f(b)\,
              \mathsf L(\mathfrak p(\omega),\gamma(\omega);\dd b)\right].
\]
Taking $f=1$ proves $\widetilde{\mathbb E}F=\mathbb EF$. For countably many samples, replace $[0,1]$ by $[0,1]^{\mathbb N}$ and use distinct coordinate uniforms. For every finite
index set $I$, the same calculation gives
\[
 \widetilde{\mathbb E}\!\left[F\prod_{i\in I}f_i(\beta_i)\right]
 =\mathbb E\!\left[F\prod_{i\in I}\int f_i(b)\,
                       \mathsf L_i(\mathfrak p_i,\gamma_i;\dd b)\right].
\]
These identities for finite collections of auxiliary paths determine the joint distribution of the original variables and the countable family of auxiliary paths and preserve the
original marginal. All disintegrations are fixed on the uncompleted Borel spaces; only then do we take the usual augmentation for the probability law used for stochastic
integration. As explained in Remark~2.2(5) of \cite{Guneysu2019}, completion under the unpinned law before pinning is inappropriate. The later augmentation does not change the
constructed probability measure on Borel sets.
\end{proof}

\begin{theorem}\label{S:thm:joint-probability-space}
For each fixed $\beta>0$, there is a standard Borel probability space carrying the finite arrays and bridge pairs, on which the Brownian path, terminal point and finite coupling
retain their original joint law. For each row,
\begin{equation}\label{S:eq:row-joint-bridge-kernel}
 \mathcal L((\gamma_{n,j},\beta_{n,j})_{1\le j<n}\mid\mathscr S_n)
 =\bigotimes_{j=1}^{n-1}
 \mathsf K_{h_{n,j},Z_{\vartheta_{n,j+1}},Z_{\vartheta_{n,j}}}.
\end{equation}
\end{theorem}
\begin{proof}
Start with the exact coupling with a common terminal point in Proposition~\ref{S:prop:7.2}. The path space $C([0,\infty);\mathbb H)$ and the countable spaces of auxiliary variables
are standard Borel. Assign a fixed value to the terminal point map on the null set where its limit is undefined.

The retained increment variables are functions of this path. For $z=Z_{\vartheta_{n,j+1}}$ and $z'=Z_{\vartheta_{n,j}}$, write $g_{z,\eta_\infty}z'$ in polar coordinates
$(R,\phi)$, using the local isometry of Lemma~\ref{S:lem:7.2}, and set
\begin{equation}\label{S:eq:actual-path-angular-seed}
 R_{n,j}=R,
 U_{n,j}=F_R^{\mathrm{ang}}(\phi),
 F_R^{\mathrm{ang}}(v)=\int_0^v\pi_R(\dd\phi),  0\le v<2\pi.
\end{equation}
The positive angular density gives the conditional probability transform. The preceding boundary points and $\eta_\infty$ generate $\mathcal F_{n,j+1}^{\eta}$, as in
\eqref{S:eq:conditioning-row-filtration}. Outside the null set where polar coordinates are undefined,
\[
 \mathbb P(U_{n,j}\le u\mid R_{n,j},\mathcal F_{n,j+1}^{\eta})=u,
   0\le u\le1.
\]
For a bounded Borel $f_j$, Lemma~\ref{S:lem:7.2} and the conditional probability transform give
\[
 \mathbb E[f_j(R_{n,j},U_{n,j})\mid\mathcal F_{n,j+1}^{\eta}]
 =\mathbb E[f_j(R_{n,j},U_{n,j})\mid\eta_\infty].
\]
The right side depends on $j,n$ and the terminal data, not on the preceding boundary points. Iterating this identity in chronological order $j=n-1,n-2,\ldots,1$ yields
\[
 \mathbb E\!\left[\prod_{j=1}^{n-1}f_j(R_{n,j},U_{n,j})
                      \,\middle|\,\eta_\infty\right]
 =\prod_{j=1}^{n-1}
       \mathbb E[f_j(R_{n,j},U_{n,j})\mid\eta_\infty].
\]
These identities for products of bounded functions prove conditional independence. Since the angular density is positive, $(F_R^{\mathrm{ang}})^{-1}(F_R^{\mathrm{ang}}(\phi))=\phi$
for $0\le\phi<2\pi$. Inverting polar coordinates and then $g_{z,\eta_\infty}$ recovers the original boundary point and the variables in \eqref{S:eq:7.46}.

The circular arrays, including their simultaneous auxiliary construction near the terminal point, are retained from Proposition~\ref{S:prop:7.2}. The extension adds the Euclidean
bridges.

Let $\mathfrak m_{n,j}$ be the actual law of $\mathfrak p=(h_{n,j},Z_{\vartheta_{n,j+1}},Z_{\vartheta_{n,j}})$. Apply Proposition~\ref{S:prop:aux-probability-kernels} with
$\nu=\mathfrak m_{n,j}$. For $\mathfrak m_{n,j}$-almost every
$\mathfrak p$,
\begin{equation}\label{S:eq:bridge-coupling-fibre-disintegration}
 \mathsf P_{\mathfrak p}(\dd\gamma)
 \mathsf L_{n,j}(\mathfrak p,\gamma;\dd\beta)
 =\mathsf K_{\mathfrak p}(\dd\gamma,\dd\beta).
\end{equation}
At the existing segment $\gamma_{n,j}$, use a new independent uniform to sample $\beta_{n,j}$ from $\mathsf L_{n,j}$. For a bounded Borel function $H$ of all the original data,
\[
 \int H(\omega)\prod_{n,j}
       \mathsf L_{n,j}(\mathfrak p_{n,j}(\omega),
                        \gamma_{n,j}(\omega);\dd\beta_{n,j})
       \,\mathbb P(\dd\omega)=\mathbb EH.
\]
The product is a countable product of probability kernels. It preserves the original joint marginal, and the countably many null sets for the parameter law have null union.

Fix a row. For bounded Borel $f_j(\gamma,\beta)$, first integrate the new uniforms, then apply Lemma~\ref{S:lem:bridge-disintegration}:
\[
 \mathbb E\left[\prod_{j=1}^{n-1} f_j(\gamma_{n,j},\beta_{n,j})\mid\mathscr S_n\right] =\prod_{j=1}^{n-1}\int\!\int f_j(\gamma,\beta) \mathsf L_{n,j}(\mathfrak p_{n,j},\gamma;\dd\beta) \mathsf P_{\mathfrak p_{n,j}}(\dd\gamma) =\prod_{j=1}^{n-1}\int f_j\,\dd\mathsf K_{\mathfrak p_{n,j}}.
\]
This is the product identity \eqref{S:eq:row-joint-bridge-kernel}; each second marginal is the Euclidean bridge law.

Define its Gaussian coordinates by
\begin{equation}\label{S:eq:joint-kl-coordinates}
 \xi_{n,j,m,\nu}=\sqrt2\,\pi m\int_0^1\beta_{n,j}^\nu(u)\sin(\pi m u)\,\dd u.
\end{equation}
For $K(u,v)=\min(u,v)-uv$,
\[
 \int_0^1K(u,v)\sin(\pi m v)\,\dd v
      =\frac{\sin(\pi m u)}{(\pi m)^2},
 \mathbb E[\xi_{n,j,m,\nu}\xi_{n,j',m',\nu'}\mid\mathscr S_n]
      =\delta_{jj'}\delta_{mm'}\delta_{\nu\nu'}.
\]
To evaluate the first integral, put $a=\pi m$ and split it at $v=u$. Integration of $v\sin(av)$ and $(1-v)\sin(av)$ gives
\[
 \begin{aligned}
 \int_0^1K(u,v)\sin(av)\dd v
 &=(1-u)\left[-\frac{u\cos(au)}a+\frac{\sin(au)}{a^2}\right]\\
 &\quad+u\left[\frac{(1-u)\cos(au)}a
                         +\frac{\sin(au)}{a^2}\right]
 =\frac{\sin(au)}{a^2}.
 \end{aligned}
\]
Multiplying this by the two normalization factors in \eqref{S:eq:joint-kl-coordinates}, and using $\int_0^1\sin(\pi m u)\sin(\pi m'u)\dd u =\mathbf1_{\{m=m'\}}/2$, proves the
covariance formula. The coordinates are jointly Gaussian as integrals of the independent Gaussian bridges, so zero mixed covariances prove independence. For each fixed component,
Parseval gives the conditional mean squared tail
\[
 \mathbb E\left[\left\|\beta_{n,j}^{\nu}
       -\sqrt2\sum_{m=1}^{M}\frac{\xi_{n,j,m,\nu}}{\pi m}
                      \sin(\pi m\,\cdot)\right\|_{L^2(0,1)}^2
                       \,\middle|\,\mathscr S_n\right]
 =\sum_{m>M}\frac1{(\pi m)^2}\longrightarrow0.
\]
Thus the later Gaussian coordinates are measurable functions of the chosen bridges, and the original joint distribution is retained.
\end{proof}

\subsubsection{Conditional laws and simultaneous versions}
\label{S:sec:conditioning-fields-and-simultaneous-versions}
\label{S:tab:conditioning-conditional}
The increment coupling reveals the increments successively, while the bridge approximation conditions on the whole sampled row. The row determines each increment from its boundary
points; the attached Euclidean bridge remains random under this conditioning. We state the conditional identities for both arguments in a common notation, also used for the
conditional expectations in Section~\ref{S:sec:the-conditional-gaussian-operator}. Subsection~\ref{S:sec:simultaneous-pathwise-estimates} constructs the common event.

Throughout this discussion, $\beta>0$ is fixed. We work on the Borel probability space $(\Omega,\mathcal A,\mathbb P)$ constructed in Theorem~\ref{S:thm:joint-probability-space};
the conditioning sigma-fields are defined before completion.

Write $\mathcal H_\eta=\sigma(\eta_\infty)$ and $\mathcal F_t^Z=\sigma(Z_u:0\le u\le t)$, before completion. For a fixed row, write the increment variables and their distributions
conditional on the terminal point as
\[
 \mathsf s_{n,j}=(R_{n,j},U_{n,j}),\qquad \nu_{n,j}(\dd r,\dd u)=p_{h_{n,j}}^{\mathrm{rad}}(r)\,\dd r\,\dd u, r>0, 0<u<1, 1\le j<n.
\]
Equation~\eqref{S:eq:actual-path-angular-seed} recovers these variables from the sampled positions. The distributions $\nu_{n,j}$ are the same for every terminal value
$\eta_\infty$.

We use two filtrations within a row:
\begin{equation}\label{S:eq:conditional-discrete-fields}
 \begin{aligned}
  \mathcal F_{n,k}^{\eta}
    &=\mathcal H_\eta\vee\sigma(\mathsf s_{n,j}:k\le j<n),
       &&1\le k\le n,\\
  \mathcal G_{n;k,j}
    &=\mathcal H_\eta\vee\sigma(\mathsf s_{n,i}:k\le i\le j),
       &&k-1\le j<n.
 \end{aligned}
\end{equation}
An empty list generates the trivial sigma-field. The first filtration reveals the Brownian increments in decreasing index order, the second in the order $k,k+1,\ldots$ used for
relative products. Both are defined by the original Brownian variables.

\phantomsection\label{S:prop:conditional-conditional-laws} Lemma~\ref{S:lem:7.2} and
Theorem~\ref{S:thm:joint-probability-space} give
\begin{equation}\label{S:eq:conditional-seed-product}
 \mathcal L((\mathsf s_{n,j})_{1\le j<n}\mid\mathcal H_\eta)
          =\bigotimes_{j=1}^{n-1}\nu_{n,j}.
\end{equation}
In coordinates sending $\eta_\infty$ to infinity, the conditional distribution of one increment does not depend on the past. Successive conditioning therefore gives, for bounded
Borel functions $f_j$,
\[
 \mathbb E\left[\prod_{j=1}^{n-1}f_j(\mathsf s_{n,j})
                 \,\middle|\,\mathcal H_\eta\right]
 =\prod_{j=1}^{n-1}\int f_j\,\dd\nu_{n,j}.
\]
Thus conditioning on the terminal point leaves the row increments independent. Each increment is a measurable function of its boundary points and $\eta_\infty$. Conversely,
starting from $Z_0=\ii$, the inverse angular transform and the local isometry recover the boundary points one at a time. These maps are measurable and inverse to one another
outside the null set where polar coordinates are undefined. Together with \eqref{S:eq:row-conditioning-global-terminal}, this gives, modulo null sets,
\begin{equation}\label{S:eq:conditional-skeleton-seeds}
 \mathscr S_n
  =\sigma(\eta_\infty,Z_{\vartheta_{n,1}},\ldots,Z_{\vartheta_{n,n}})
  =\mathcal F_{n,1}^{\eta}.
\end{equation}
In particular $\mathcal F_{n,j+1}^{\eta}\subset \mathcal F_{\vartheta_{n,j+1}}^Z\vee\mathcal H_\eta$ and
$\mathcal F_{n,j}^{\eta}\subset \mathcal F_{\vartheta_{n,j}}^Z\vee\mathcal H_\eta$.

In the coordinates of Theorem~\ref{S:thm:7.4}, write $f_{n,j}(\mathsf s_{n,j})=\zeta_{n,j}^{\mathrm{loc}}$. The moment bound in that theorem allows us to separate the mean from the
centered error:
\begin{equation}\label{S:eq:conditional-drift-kernel}
 b_{n,j}=\int f_{n,j}\,\dd\nu_{n,j},
 \xi_{n,j}=f_{n,j}(\mathsf s_{n,j})-\int f_{n,j}\,\dd\nu_{n,j}.
\end{equation}
The mean $b_{n,j}$ is deterministic in these coordinates. The product law \eqref{S:eq:conditional-seed-product} implies, for $k\le j$ and $p\ge2$,
\begin{equation}\label{S:eq:conditional-increment-moments}
 \mathbb E[\xi_{n,j}\mid\mathcal G_{n;k,j-1}]=0,\qquad
 \mathbb E[\|\xi_{n,j}\|^p\mid\mathcal G_{n;k,j-1}]^{1/p}
 \le C_\beta p^{3/2}j^{-3/2}.
\end{equation}
The same statements hold conditional on $\mathcal F_{n,j+1}^{\eta}$. Indeed, either sigma-field reveals a subset of the other increments, which are independent of $\mathsf s_{n,j}$
conditional on the terminal point. The conditional expectation identity follows first for bounded functions by the product law, and then for moments by truncation and monotone
convergence. In increasing index order, $\operatorname{Ad}_{P_{n;k,j}}$ is $\mathcal G_{n;k,j-1}$-measurable. This is the predictability used in
\eqref{S:eq:s7-reversed-filtration-moments}.

Conditioning on the whole row is different: $\mathscr S_n$ determines the error, so
\begin{equation}\label{S:eq:conditional-skeleton-reveals-finite}
 \mathbb E[\xi_{n,j}\mid\mathscr S_n]=\xi_{n,j}.
\end{equation}

The variables $\xi_{n,j,m,\nu}$, with four indices, are the Gaussian coordinates of the auxiliary Euclidean bridges. Unlike the increment errors $\xi_{n,j}$, they remain random
after conditioning on $\mathscr S_n$. The second marginal of \eqref{S:eq:row-joint-bridge-kernel} and the covariance calculation in \eqref{S:eq:joint-kl-coordinates} give, for
every finitely supported real array $(t_{j,m,\nu})$,
\[
 \mathbb E\left[
  \exp\left(\ii\sum_{j,m,\nu}t_{j,m,\nu}\xi_{n,j,m,\nu}\right)
       \,\middle|\,\mathscr S_n\right]
 =\exp\left(-\frac12\sum_{j,m,\nu}t_{j,m,\nu}^2\right).
\]
Thus the coordinates are independent standard normals conditional on $\mathscr S_n$. In contrast, the errors $\xi_{n,j}$ in \eqref{S:eq:conditional-drift-kernel} are
$\mathscr S_n$-measurable.

In particular, for a fixed row,
\begin{equation}\label{S:eq:s9-conditioning-distinction}
 \mathbb E[\xi_{n,j}\mid\mathscr S_n]=\xi_{n,j},\qquad
 \mathbb E[\xi_{n,j,m,\nu}\mid\mathscr S_n]=0,\qquad
 \mathbb E[\xi_{n,j,m,\nu}^2\mid\mathscr S_n]=1.
\end{equation}
Here $\xi_{n,j}$ is the centered error of a finite increment, whereas $\xi_{n,j,m,\nu}$ is a Gaussian bridge coordinate. Stochastic integrals use $W$ under $\mathbb P$,
$\widetilde W^\eta$ under $\mathbb P^\eta$, and $W^{\mathrm{pin}}$ under the finite bridge law, in the respective augmented natural filtrations.

Observing every sampling row determines the whole Brownian path, as the next lemma shows.

\begin{lemma}\label{S:lem:conditional-all-skeletons}
Up to null sets,
\begin{equation}\label{S:eq:conditional-all-row-information}
 \bigvee_{n\ge2}\mathscr S_n=\sigma(Z_t:t\ge0).
\end{equation}
\end{lemma}
\begin{proof}
Let $\mathcal A=\bigvee_{n\ge2}\mathscr S_n$. Each sampling time is deterministic, and each sampled position is an evaluation of $Z$. The terminal point is the limit of its path,
and the terminal coordinate used in $\mathscr S_n$ is a measurable function of that point and the sampled positions. Consequently
\[
 \mathcal A\subseteq\sigma(Z_t:t\ge0)
\]
up to the completions specified in the statement.

For the opposite inclusion, fix a positive rational $t$ and put $r=e^{-\beta t/4}\in(0,1)$. Choose $j_n=\lfloor nr\rfloor$ for every $n$ with $nr\ge2$ and $n(1-r)>1$. These
inequalities ensure $1\le j_n<n$, and
\[
 |j_n/n-r|\le n^{-1},\qquad j_n\ge nr/2.
\]
The deterministic sampling estimate gives
\[
 \left|\vartheta_{n,j_n}-t\right|
 \le\frac4\beta\left|\log\frac{r}{j_n/n}\right|
                          +\frac{C_\beta}{j_n}
 \le\frac{C_{\beta,t}}n\longrightarrow0.
\]
For the logarithmic bound we used the mean value theorem on $[r/2,r]$. Since $Z_{\vartheta_{n,j_n}}$ is $\mathcal A$-measurable and $Z$ has continuous paths,
$Z_t=\lim_nZ_{\vartheta_{n,j_n}}$ is $\mathcal A$-measurable. A countable intersection gives this conclusion for every positive rational $t$ on one event; $Z_0=\ii$ is
deterministic. For any real $t\ge0$, continuity and rational approximation then give $Z_t$ as a limit of these measurable evaluations. Thus $\sigma(Z_t:t\ge0)\subseteq\mathcal A$,
proving equality. In particular, the path segment on a fixed sampling interval is measurable with respect to $\mathcal A$, so its conditional distribution given $\mathcal A$ is a
point mass.
\end{proof}
Consequently, conditional on all the $\mathscr S_n$, each hyperbolic path segment has a point mass distribution. The bridge distribution used above is conditional on a single
$\mathscr S_n$. Auxiliary variables attached to the path can still be random when the whole path is known.

There are countably many conditional identities used in the proof. We choose versions on one event and localize the moment identities on events where the conditional bounds are at
most a fixed integer, before allowing that integer to increase.

Let $\mathfrak T_\beta$ consist of the nonnegative rational times and the deterministic sampling times $\vartheta_{n,j}$. Let $\mathfrak C_\beta$ be the countable family of
observation maps that generate $\mathcal H_\eta$, the sigma-fields in \eqref{S:eq:conditional-discrete-fields}, $\mathscr S_n$, and the path sigma-fields in
Subsection~\ref{S:sec:conditioning-fields-and-simultaneous-versions} at times $t\in\mathfrak T_\beta$. Where a bridge argument requires it, we also include the observation of its
boundary points and its past at a rational fraction of its duration. These observations take values in standard Borel spaces. A stopped path is continued constantly after its
stopping time, so it takes values in a fixed path space.

\phantomsection\label{S:lem:conditional-countable-versions} For each observation $T_a$, choose a regular conditional
probability kernel and a countable algebra $\mathcal D_a$ that determines probability measures on the space of variables being conditioned. Such an algebra exists because the space
is standard Borel; for a product space, one may use finite cylinders. Let $N_{a,D}$ be the null set on which the conditional identity for $D\in\mathcal D_a$ fails. Then
\[
 N=\bigcup_a\bigcup_{D\in\mathcal D_a}N_{a,D},
   \mathbb P(N)=0.
\]
Outside $N$, the monotone class theorem extends the identity from $\mathcal D_a$ to all Borel sets. We apply this construction to
\eqref{S:eq:conditional-seed-product}--\eqref{S:eq:conditional-skeleton-reveals-finite}, Lemma~\ref{S:lem:bridge-disintegration}, and Theorem~\ref{S:thm:joint-probability-space},
for all rows, subintervals and finite mode cutoffs at once.

For the Gaussian random elements, apply the isometry $\mathcal J_n$ of Lemma~\ref{S:lem:weighted-hilbert} to work in a fixed separable Hilbert space. Their conditional
distributions are determined by a countable dense family of linear functionals, with rational parameters in the characteristic functions; continuity gives the remaining parameters.
Once the conditional kernel $K_a$ has been fixed, its nonnegative moments satisfy
\[
 \int f\,\dd K_a=\lim_{M\to\infty}\int(f\wedge M)\,\dd K_a,
   f\ge0.
\]
For a Hilbert-valued random element $V$ with conditional kernel $K_a$, put $A_M=\{\int\|v\|^2K_a(T_a,\dd v)\le M\}$. Then $A_M\in\sigma(T_a)$ and
\[
 \mathbb E[\mathbf1_{A_M}\|V\|^2]
 =\mathbb E\left[\mathbf1_{A_M}\int\|v\|^2K_a(T_a,\dd v)\right]\le M.
\]
The localized Bochner conditional expectations therefore exist. If the conditional second moment is finite almost surely, the sets $A_M$ increase to an event of probability one.
This is the localization in the \hyperref[S:prop:gaussian-covariance-measurability]{covariance construction}.

The same kernels also define the moments of noninteger order. When an admissible integer $q\ge p$ is available, the conditional Lyapunov inequality gives
$\|Y\|_{L^p(K_a)}\le\|Y\|_{L^q(K_a)}$. The increasing moment exponents used later belong to fixed countable sequences. Stochastic identities are established in the usual
augmentation for the probability law in question and have continuous versions on finite time intervals. We intersect the events on which they hold at rational times with the events
constructed above, and denote a Borel subset of this intersection of probability one by $\Omega_\beta^{\mathrm{ver}}$. These choices require conditional identities only for almost
every value of the conditioning variables.

\subsubsection{Normal coordinates and displacement estimates}
\label{S:sec:s9-normal-coordinates}

In normal coordinates the symmetric square root of the inverse metric is $I+O(|x|^2)$. This quadratic error is the gain used in the pinned coupling. The later comparison of
boundary point increments uses the diffusion matrix for the original driving Brownian motions $W^1,W^2$, which has a linear correction. A pathwise comparison keeps this driver
fixed, so the diffusion matrix is specified separately in each construction.

Let $x\in T_o\mathbb H\simeq\mathbb R^2$, $r=|x|$, and, for $r>0$, set $P_x=xx^{\mathsf T}/r^2$, $P_x^\perp=I-P_x$. The polar metric $\dd r^2+\sinh^2r\,\dd\phi^2$ gives
\begin{equation}\label{S:eq:s9-normal-metric}
 \begin{aligned}
 g(x)&=P_x+\left(\frac{\sinh r}{r}\right)^2P_x^\perp,\qquad
 a(x)=P_x+\left(\frac r{\sinh r}\right)^2P_x^\perp,\\
 \sigma(x)&=a(x)^{1/2}
       =P_x+\frac r{\sinh r}P_x^\perp.
 \end{aligned}
\end{equation}
These expressions extend smoothly to $x=0$. For example,
\[
 g(x)=I+\psi(r^2)(r^2I-xx^{\mathsf T}),\qquad
 \psi(r^2)=\frac{(\sinh r/r)^2-1}{r^2}
          =\frac13+O(r^2).
\]
The power series and their differentiated series converge uniformly on a fixed ball. They give $g-I=O(r^2)$, $\partial g=O(r)$, and hence $\Gamma^i_{jk}=O(r)$ by
$\Gamma^i_{jk}=\tfrac12g^{i\ell} (\partial_jg_{k\ell}+\partial_kg_{j\ell}-\partial_\ell g_{jk})$. The coordinate drift of $\tfrac12\Delta_{\mathbb H}$ is explicit:
\begin{equation}\label{S:eq:s9-normal-drift}
 b_0(x)=\frac12\left(\frac{\coth r}{r}
                         -\frac1{\sinh^2r}\right)x
       =\frac13x+O(r^3).
\end{equation}
For $e(\phi)=(\cos\phi,\sin\phi)$, $\partial_r^2(re)=0$, $\partial_r(re)=e$, and $\partial_\phi^2(re)=-re$. Hence $\tfrac12\Delta(re)=\tfrac12(\coth r-r/\sinh^2r)e=b_0(re)$. In
particular, on a fixed ball,
\[
 \|\sigma(x)-I\|\le C|x|^2,\qquad |b_0(x)|\le C|x|.
\]
The constants are independent of $o$, since an isometry identifies the normal charts and the orthonormal frames.

Consider either the original diffusion or its terminal conditional law. Their generators are $L+v\cdot\nabla$, where respectively $v=0$ or $v=\nabla\log P_\eta$, and $|v|\le1$. For
$Z_0=z$ set $\tau_{r_0}=\inf\{t:d(z,Z_t)\ge r_0\}$.

\begin{lemma}\label{S:lem:global-hyperbolic-displacement}
For $p\ge2$ and $h>0$,
\begin{equation}\label{S:eq:global-displacement-moment}
 \left\|\sup_{0\le t\le h}d(z,Z_t)\right\|_{L^p}
 \le C(\sqrt{ph}+h),
\end{equation}
with an absolute constant under either law, also conditionally on the past at a deterministic starting time.
\end{lemma}
\begin{proof}
Write $r=d(z,\cdot)$. The function $r^2$ is smooth everywhere on the hyperbolic plane, including $z$. In polar coordinates the metric is $\dd r^2+\sinh^2r\,\dd\phi^2$, so
\[
 \tfrac12\Delta(r^2)=1+r\coth r,
 |\nabla(r^2)|^2=4r^2.
\]
At $r=0$ the first expression has limit $2$. Since $\coth r=1+2/(\me^{2r}-1)\le1+r^{-1}$ for $r>0$, and $|v|\le1$ by Lemma~\ref{S:lem:terminal-affine-innovations},
\[
 (\tfrac12\Delta+v\cdot\nabla)(r^2)\le2+3r.
\]
Stop when $r$ reaches $M$. It\^o's formula gives a continuous martingale $N$ such that
\[
 r(Z_{t\wedge\tau_M})^2
 =N_{t\wedge\tau_M}+\int_0^{t\wedge\tau_M}
           (\tfrac12\Delta+v\cdot\nabla)(r^2)(Z_u)\,\dd u,
 \langle N\rangle_{h\wedge\tau_M}=4\int_0^{h\wedge\tau_M}r(Z_u)^2\dd u.
\]
Apply Lemma~\ref{S:lem:aux-martingale-moments} to $N$ at $q=p/2\ge1$. Its stated range gives
\[
 \|N^*\|_{L^{p/2}}
 \le C\sqrt p\|\langle N\rangle^{1/2}\|_{L^{p/2}}.
\]
Put $A_M=\|\sup_{t\le h}r(Z_{t\wedge\tau_M})\|_{L^p}$. Applying this martingale inequality with $q=p/2$, and using monotonicity of norms on a probability space,
$\|\cdot\|_{L^{p/2}}\le\|\cdot\|_{L^p}$, gives
\[
 A_M^2\le C\sqrt{ph}\,A_M+2h+3hA_M.
\]
Set $D=C\sqrt{ph}+3h$. Since $DA_M\le A_M^2/2+D^2/2$, the last inequality gives $A_M^2\le D^2+4h$. Taking square roots and using $p\ge2$ proves $A_M\le C(\sqrt{ph}+h)$ uniformly in
$M$. At $p=2$ this also implies $\mathbb P(\tau_M\le h)\le C_h/M^2$, so the diffusion does not explode. Letting $M$ tend to infinity proves \eqref{S:eq:global-displacement-moment}.

For each fixed $r_0>0$ we also obtain constants $h_0,c,C>0$ such that
\begin{equation}\label{S:eq:fixed-normal-ball-exit}
 \mathbb P(\tau_{r_0}\le h)\le C\me^{-c/h},\qquad 0<h\le h_0.
\end{equation}
To prove this, choose a small fixed $a>0$ with $C\sqrt a\le(4\me)^{-1}$, and set $p=ar_0^2/h$. Decrease $h_0$ until $p\ge2$ and $Ch\le r_0/(4\me)$ for $h\le h_0$. Then
\eqref{S:eq:global-displacement-moment} and Markov's inequality give
\[
 \mathbb P(\tau_{r_0}\le h)
 \le\left(\frac{C(\sqrt{ph}+h)}{r_0}\right)^p
 \le(2\me)^{-p}\le \me^{-c/h}.
\]
At a deterministic time $s$, the Markov property gives the same law started at $Z_s$. Uniformity in that starting point gives
\[
 \mathbb E\!\left[\sup_{0\le t\le h}d(Z_s,Z_{s+t})^p
       \,\middle|\,\mathcal F_s^Z\right]^{1/p}
 \le C(\sqrt{ph}+h).
\]
Under $\mathbb P^\eta$ the identical bound holds, since the drift bound $|v|\le1$ is uniform in $\eta$.
\end{proof}

\subsubsection{Approximation of a short Brownian bridge}
\label{S:sec:s9-short-bridge}

In normal coordinates, a short hyperbolic Brownian bridge is approximated by the linear interpolation between its endpoints and a Euclidean Brownian bridge. The error estimate
below is for the fixed coupling already constructed, including paths that leave the localization ball; this allows the same bridge variables to be used as the parameters vary.

For the fixed coupling $\mathsf K_{h,z,z'}$ put
\[
 v=h^{-1/2}\exp_z^{-1}z',\qquad
 \exp_z^{-1}\gamma(s)=\sqrt h\,(sv+\beta(s))+h\rho_h(s).
\]
This defines $\rho_h$ on every path of the coupling, including paths which leave a specified ball. For moment order $p\ge2$, introduce $L\ge1$ and $R\ge4(1+L+\sqrt p)$.

\begin{lemma}\label{S:lem:short-bridge-expansion}
There are absolute $h_0,c_0,c,C>0$ such that, if $0<h\le h_0$, $|v|\le L$ and $\sqrt h R^2\le c_0$, then for $u\ge1$,
\begin{equation}\label{S:eq:short-bridge-probability}
 \mathsf K_{h,z,z'}\bigl(\|\rho_h\|_\infty>CuR\bigr)
 \le C\me^{-cR^2}+(C\sqrt h)^p+u^{-p}.
\end{equation}
\end{lemma}

The Euclidean and hyperbolic bridges have the same singular pinning drift to first order. Subtracting their equations leaves a stable linear term $-\eta_s/(1-s)$. Its integrating
factor controls the error up to the pinned boundary point.

Let $B$ be a Euclidean Brownian motion and let $F,g$ be progressively measurable on $[0,1)$, with the integrability required in the next display. Matrix norms of $F$ are
Hilbert--Schmidt norms. Define
\[
 \mathcal Y_s=(1-s)\int_0^s\frac{F_u}{1-u}\,\dd B_u
             +(1-s)\int_0^s\frac{g_u}{1-u}\,\dd u,
 \qquad 0\le s<1.
\]
Write $\mathcal Y^F$ for its stochastic part.

\begin{lemma}\label{S:lem:pinned-linear-maximal}
For $p\ge2$,
\begin{equation}\label{S:eq:pinned-stochastic-hardy}
\left\|\sup_{0\le s<1}|\mathcal Y_s|\right\|_{L^p}
\le C\sqrt p\left\|\left(\int_0^1\|F_u\|^2\dd u\right)^{1/2}\right\|_{L^p}
 +\left\|\int_0^1|g_u|\dd u\right\|_{L^p}.
\end{equation}
If $\|F_u\|_{L^p}\le f_0$ for all $u$, then
\[
\|\mathcal Y_s^{F}\|_{L^p}
\le C\sqrt p\,f_0\sqrt{s(1-s)}.
\]

\end{lemma}
\begin{proof}
Put $s=1-\me^{-r}$ and define
\[
\widetilde B_r:=\int_0^{1-\me^{-r}}(1-u)^{-1/2}\,\dd B_u,
\phi_r:=\me^{-r/2}F_{1-\me^{-r}}.
\]
The integrand is deterministic and $[\widetilde B^i,\widetilde B^j]_r =\delta_{ij}\int_0^{1-\me^{-r}}(1-u)^{-1}\dd u=\delta_{ij}r$. Disjoint increments are independent Gaussian
vectors, so $\widetilde B$ is Brownian. The time substitution turns the stochastic part of $\mathcal Y$ into the Ornstein--Uhlenbeck stochastic convolution
\[
Y_r=\int_0^r \me^{-(r-q)}\phi_q\,\dd\widetilde B_q,
  \dd Y_r=-Y_r\,\dd r+\phi_r\,\dd\widetilde B_r.
\]
Moreover, $\int_0^\infty\|\phi_r\|^2\dd r=\int_0^1\|F_u\|^2\dd u$. Applying It\^o's formula to $|Y_r|^2$ up to time $R$ gives
\[
\sup_{r\le R}|Y_r|^2
\le 2\sup_{r\le R}\left|\int_0^r\langle Y_q,\phi_q\,\dd\widetilde B_q\rangle\right|
 +\int_0^R\|\phi_q\|^2\dd q.
\]
Take the $L^{p/2}$-norm. Lemma~\ref{S:lem:aux-martingale-moments} at exponent $p/2\ge1$ and H\"older's inequality yield
\[
\left\|\sup_{r\le R}|Y_r|\right\|_{L^p}^2
\le C\sqrt p\left\|\sup_{r\le R}|Y_r|\right\|_{L^p}
 \left\|\left(\int_0^R\|\phi_q\|^2\dd q\right)^{1/2}\right\|_{L^p}
 +\left\|\left(\int_0^R\|\phi_q\|^2\dd q\right)^{1/2}\right\|_{L^p}^2.
\]
For $m\in\mathbb N$, set $\tau_m=\inf\{r:|Y_r|\ge m\ \text{or}\ \int_0^r\|\phi_q\|^2\dd q\ge m\}$. Apply the preceding calculation to $Y_{r\wedge\tau_m}$. Young's inequality
absorbs half the squared norm on the left and gives
\[
 \mathbb E\sup_{r\le R}|Y_{r\wedge\tau_m}|^p
 \le(C\sqrt p)^p\,
       \mathbb E\left(\int_0^{R\wedge\tau_m}\|\phi_q\|^2\dd q\right)^{p/2}
 \le(C\sqrt p)^p\,
       \mathbb E\left(\int_0^R\|\phi_q\|^2\dd q\right)^{p/2}.
\]
If the last expectation is finite, continuity and local square integrability give, for almost every outcome, an $m_0(\omega,R)$ such that $\tau_m>R$ for $m\ge m_0(\omega,R)$.
Fatou's lemma removes the stopping; monotone convergence as $R\uparrow\infty$ proves the required bound. If the expectation is infinite, the bound is automatic. For the drift part,
$(1-s)/(1-u)\le1$ whenever $u\le s$, and hence
\[
\sup_{s<1}(1-s)\int_0^s\frac{|g_u|}{1-u}\dd u\le\int_0^1|g_u|\dd u.
\]
Adding the stochastic and drift bounds gives \eqref{S:eq:pinned-stochastic-hardy}. For the pointwise estimate, Lemma~\ref{S:lem:aux-martingale-moments} at exponent $p$ and
Minkowski in $L^{p/2}$, gives
\[
\|\mathcal Y_s^F\|_{L^p}
\le C\sqrt p\,(1-s)
\left(\int_0^s\frac{f_0^2}{(1-u)^2}\dd u\right)^{1/2}
=C\sqrt p\,f_0\sqrt{s(1-s)}.
\]
\end{proof}

Use the normal chart at $o$, and abbreviate $p_\tau(\exp_ox,\exp_oy)$ by $p_\tau(x,y)$. Let $\nabla_x$ be its Euclidean coordinate gradient and define
\[
 \mathfrak r_\tau(x,y)
 :=a(x)\nabla_x\log p_\tau(x,y)-\frac{y-x}{\tau}.
\]

\begin{lemma}\label{S:lem:hyperbolic-heat-gradient}
There are absolute $r_0,\tau_0,C>0$ such that, for $|x|,|y|\le r_0$ and $0<\tau\le\tau_0$,
\begin{equation}\label{S:eq:hyperbolic-heat-gradient-remainder}
 |\mathfrak r_\tau(x,y)|
 \le C(1+|x|+|y|)
 +\frac{C}{\tau}(|x|+|y|)^2|x-y|.
\end{equation}
The constants are independent of the center $o$.

\end{lemma}
\begin{proof}
Lemma~\ref{S:lem:heat-kernel-estimates}, specifically \eqref{S:eq:s6-heat-factorization}, gives on the fixed normal ball
\begin{equation}\label{S:eq:strong-local-heat-factorization}
 p_\tau(x,y)=\frac{\me^{-d(x,y)^2/(2\tau)}}{2\pi\tau}
 \bigl(u_0(x,y)+\tau\psi_\tau(x,y)\bigr),
\end{equation}
with
\[
 \sup_{0<\tau\le\tau_0}\sup_{|x|,|y|\le2r_0}
 \bigl(|\psi_\tau(x,y)|+|\nabla_x\psi_\tau(x,y)|\bigr)\le C.
\]
The leading coefficient is computed directly in Lemma~\ref{S:lem:heat-kernel-estimates}:
\[
 u_0(x,y)=\left(\frac{d(x,y)}{\sinh d(x,y)}\right)^{1/2},  u_0(x,x)=1.
\]
The function $(r/\sinh r)^{1/2}$ is a smooth positive function of $r^2$ near zero. Squared distance is smooth in these charts, so $u_0$ is smooth also at $x=y$, and its minimum on
the compact product is strictly positive. After decreasing $\tau_0$, the amplitude
\[
 q_\tau(x,y):=u_0(x,y)+\tau\psi_\tau(x,y)
\]
satisfies $q_\tau\ge c>0$, and therefore
\begin{equation}\label{S:eq:heat-amplitude-gradient}
 \sup_{0<\tau\le\tau_0}\sup_{|x|,|y|\le2r_0}
 \bigl|a(x)\nabla_x\log q_\tau(x,y)\bigr|\le C.
\end{equation}

Let
\[
 v_{x,y}:=\exp_x^{-1}y\in T_x\mathbb H
\]
and write its components in the coordinate basis associated with the normal chart. For a smooth variation of the initial point with velocity $w$, differentiate the energy of the
minimizing geodesic on $[0,1]$. Metric compatibility and integration by parts give
\[
 \frac{\dd}{\dd\epsilon}\frac12\int_0^1
       |\dot\gamma_\epsilon(t)|_g^2\dd t\bigg|_{\epsilon=0}
 =\bigl\langle\partial_\epsilon\gamma_\epsilon,
                   \dot\gamma_0\bigr\rangle_g\big|_0^1
  -\int_0^1\langle\partial_\epsilon\gamma_\epsilon,
                      \nabla_t\dot\gamma_0\rangle_g\dd t
 =-\langle w,v_{x,y}\rangle_g.
\]
The terminal variation is zero and the last integral vanishes by the geodesic equation. Thus
\begin{equation}\label{S:eq:first-variation-distance-coordinate}
 a(x)\nabla_x\frac{d(x,y)^2}{2}=-v_{x,y}.
\end{equation}
The geodesic equation gives
\begin{equation}\label{S:eq:normal-log-coordinate-comparison}
 |v_{x,y}-(y-x)|
 \le C(|x|+|y|)^2|x-y|.
\end{equation}
Choose $r_0$ so small that the coordinate ball of radius $2r_0$ is geodesically convex and the Euclidean and Riemannian norms there are uniformly equivalent. Write the constant
speed minimizing geodesic from $x$ to $y$ in these coordinates as $\gamma:[0,1]\to\mathbb R^2$. Then $\gamma'(0)=v_{x,y}$ and
\[
 y-x-v_{x,y}=\int_0^1(1-s)\gamma''(s)\,\dd s.
\]
The calculation after~\eqref{S:eq:s9-normal-metric} gives $|\Gamma(z)|\le C|z|$. The geodesic equation and uniform metric equivalence give
\[
 |\gamma''(s)|\le C|\gamma(s)|\,|\gamma'(s)|^2,
 |\gamma'(s)|\le C|x-y|.
\]
Moreover,
\[
 |\gamma(s)|\le |x|+\int_0^1|\gamma'(u)|\,\dd u
 \le C(|x|+|y|).
\]
Consequently
\[
 |v_{x,y}-(y-x)|
 \le C(|x|+|y|)|x-y|^2
 \le C(|x|+|y|)^2|x-y|,
\]
which proves \eqref{S:eq:normal-log-coordinate-comparison}.

Taking the logarithmic derivative of \eqref{S:eq:strong-local-heat-factorization}, and using \eqref{S:eq:first-variation-distance-coordinate}, yields
\[
 a(x)\nabla_x\log p_\tau(x,y)
 =\frac{v_{x,y}}{\tau}+a(x)\nabla_x\log q_\tau(x,y).
\]
Subtracting $(y-x)/\tau$ gives $\mathfrak r_\tau(x,y)=\tau^{-1}(v_{x,y}-(y-x)) +a(x)\nabla_x\log q_\tau(x,y)$. The two bounds \eqref{S:eq:normal-log-coordinate-comparison} and
\eqref{S:eq:heat-amplitude-gradient} give \eqref{S:eq:hyperbolic-heat-gradient-remainder}. Finally, hyperbolic isometries identify the calculations in normal
coordinates at different centers, so the constants are independent of $o$. \end{proof}

\begin{proof}[Proof of Lemma~\ref{S:lem:short-bridge-expansion}]
By hyperbolic homogeneity we may apply an isometry sending $z$ to a fixed base point. Use Riemannian normal coordinates centered there and choose the orthonormal frame so that the
differential of $x\mapsto\mathcal A(\exp_zx)$ at the origin is given by the two matrices $\mathsf A_1,\mathsf A_2$ used below. Write $Y_t=\exp_z^{-1}(Z_t)$ and $y'=h^{1/2}v$.

Equations~\eqref{S:eq:s9-normal-metric}--\eqref{S:eq:s9-normal-drift} give $\|\sigma(y)-I\|\le C|y|^2$ and $|b_0(y)|\le C|y|$ on a fixed ball.
Lemma~\ref{S:lem:hyperbolic-heat-gradient} supplies its pinned drift, with $a=\sigma\sigma^{\mathsf T}$. In the following equations $b=b_0$ and $B$ is the driving Brownian motion
in~\eqref{S:eq:pinned-kernel-sde}.

Under the bridge law the coordinate process is the Doob transform
\[
\dd Y_t=\sigma(Y_t)\dd B_t+
\left[b(Y_t)+a(Y_t)\nabla_1\log p_{h-t}(Y_t,y')\right]\dd t,
  t<h.
\]
Set $Y_{hs}=h^{1/2}\xi_s$ and rescale the pinned driver by $B_s=h^{-1/2}W_{hs}^{\mathrm{pin}}$; we use this same rescaled driver in both equations below. Up to the first exit from
the normal ball,
\[
\dd\xi_s
=\sigma(h^{1/2}\xi_s)\dd B_s
 +\frac{v-\xi_s}{1-s}\dd s
 +r_h(s,\xi_s,v)\dd s,
\]
where
\[
|r_h(s,x,v)|
\le C\sqrt h(1+|x|+|v|)
 +Ch(1+|x|+|v|)^2\frac{|x-v|}{1-s}.
\]
Let $\bar\xi_s=sv+\beta(s)$ be the solution of the Euclidean pinned equation driven by the same Brownian motion,
\[
\dd\bar\xi_s=\dd B_s+\frac{v-\bar\xi_s}{1-s}\dd s,
  \bar\xi_0=0.
\]
Put $\eta_s=\xi_s-\bar\xi_s$. For
\[
\tau_R:=\inf\{s<1:|\xi_s|\ge R\}\wedge1,
\]
define the stopped mild extension
\[
\eta_s^{(R)}=(1-s)\int_0^s
\frac{\mathbf1_{\{q<\tau_R\}}[\sigma(h^{1/2}\xi_q)-I]}{1-q}\,\dd B_q
 +(1-s)\int_0^s
\frac{\mathbf1_{\{q<\tau_R\}}r_h(q,\xi_q,v)}{1-q}\,\dd q.
\]
For $s\le\tau_R$ this equals $\eta_s$. Since $|\xi_q|\le R$ before $\tau_R$,
\[
\|\sigma(h^{1/2}\xi_q)-I\|\le ChR^2.
\]
The singular remainder contains $|\xi_s-v|/(1-s)$. We use the square-root decay of $\xi_s-v$ at the pinned boundary point to integrate this term. Define
\[
K_p:=\sup_{0\le s<1}(1-s)^{-1/2}
\|\eta_s^{(R)}\|_{L^p}.
\]
In the next inequalities first take the supremum only over $0\le s\le1-\epsilon$ and write $K_{p,\epsilon}$. This quantity is finite: before $\tau_R$ the coefficients of the scaled
equation are bounded, and the remaining time denominators are bounded on this interval. All constants below are independent of $\epsilon$. The pointwise part of
Lemma~\ref{S:lem:pinned-linear-maximal} bounds the stochastic contribution to $K_{p,\epsilon}$ by $C\sqrt p\,hR^2$. The nonsingular part of the drift contributes at most
$C\sqrt h(1+R+L)$, because
\[
(1-s)\int_0^s\frac{\dd q}{1-q}
=(1-s)\log\frac1{1-s}
\le C\sqrt{1-s}.
\]
For the singular part,
\[
\|\bar\xi_q-v\|_{L^p}
\le (1-q)L+C\sqrt p\sqrt{q(1-q)}
\le C(L+\sqrt p)\sqrt{1-q}.
\]
Moreover, before $\tau_R$, $\eta_q=\eta_q^{(R)}$. Since
\[
(1-s)\int_0^s(1-q)^{-3/2}\dd q
\le2\sqrt{1-s},
\]
the singular drift contributes at most
\[
Ch(1+R+L)^2(L+\sqrt p+K_{p,\epsilon})
\]
to $K_{p,\epsilon}$. Consequently,
\[
K_{p,\epsilon}\le C\left[
\sqrt p\,hR^2+\sqrt h(1+R+L)
+h(1+R+L)^2(L+\sqrt p+K_{p,\epsilon})
\right].
\]
Here $1+R+L\le5R/4$ and $L+\sqrt p\le R/4$. Since $h\le1$,
\[
 h(1+R+L)^2\le\frac{25}{16}\sqrt h(\sqrt hR^2)
                    \le\frac{25}{16}c_0,
\]
while
\[
 \sqrt p\,hR^2\le\frac{c_0}{4}\sqrt hR,
 \qquad h(1+R+L)^2(L+\sqrt p)
       \le\frac{25c_0}{64}\sqrt hR.
\]
Choose $c_0$ so that the coefficient of $K_{p,\epsilon}$ is at most $1/2$. Moving that term to the left proves
\begin{equation}\label{S:eq:s9-pinned-pointwise-error}
 K_{p,\epsilon}\le C\sqrt h\,R,
\end{equation}
with $C$ independent of $\epsilon,h,p,L,R$. Letting $\epsilon\downarrow0$ and using monotone convergence of the suprema gives $K_p\le C\sqrt h\,R$. For the maximal norm, the
singular drift satisfies the following bound, whose right side is independent of $s$. Since $(1-s)/(1-q)^2\le(1-q)^{-1}$ for $q\le s$,
\[
 \sup_{s<1}(1-s)\int_0^s
 \frac{\mathbf1_{\{q<\tau_R\}}|\xi_q-v|}{(1-q)^2}\,\dd q
 \le\int_0^1
 \frac{|\bar\xi_q-v|+|\eta_q^{(R)}|}{1-q}\,\dd q.
\]
Minkowski and the bound $K_p\le C\sqrt hR$ give
\[
 \left\|\int_0^1
 \frac{|\bar\xi_q-v|+|\eta_q^{(R)}|}{1-q}\,\dd q\right\|_{L^p}
 \le C(L+\sqrt p+K_p)\int_0^1(1-q)^{-1/2}\dd q
 \le C(L+\sqrt p+K_p).
\]
The stochastic part is bounded in maximal $L^p$ by $C\sqrt p\,hR^2$ using Lemma~\ref{S:lem:pinned-linear-maximal}; the nonsingular drift is bounded pathwise by $C\sqrt h(1+R+L)$.
After multiplication of the preceding singular majorant by $Ch(1+R+L)^2$, the assumptions $R\ge4(1+L+\sqrt p)$ and $\sqrt hR^2\le c_0$ yield
\[
 \left\|\sup_{s<1}|\eta_s^{(R)}|\right\|_{L^p}
 \le C\sqrt h\,R.
\]
All stochastic estimates may first be stopped when their brackets reach an integer. The coefficient bound controls these brackets on $[0,1-\epsilon]$; Fatou removes this auxiliary
stopping. The drift majorant has finite $L^p$ norm, so its integrals are almost surely finite and the estimates extend to $s<1$.

Since $\sup_s|\bar\xi_s|\le L+\sup_s|\beta(s)|$, a bound on the bridge maximum controls exit from the stopped region. For a Euclidean Brownian motion $W$, the process $W_s-sW_1$
has the bridge law and satisfies $\sup_s|W_s-sW_1|\le2\sup_s|W_s|$. If the latter vector supremum exceeds $x/2$, one of the two scalar coordinate suprema exceeds $x/(2\sqrt2)$.
Reflection and the one-dimensional Gaussian tail give $\mathbb P(\sup_s|\beta(s)|>x)\le8\me^{-x^2/16}$. This estimate depends only on the bridge law. Since $R\ge4(1+L)$,
\[
\mathbb P\left(\sup_s|\bar\xi_s|>R/2\right)
\le C\me^{-cR^2}.
\]
If $\tau_R<1$ while $\sup_s|\bar\xi_s|\le R/2$, then $\sup_s|\eta_s^{(R)}|\ge R/2$. Markov's inequality therefore gives
\[
\mathbb P(\tau_R<1)
\le C\me^{-cR^2}+(C\sqrt h)^p.
\]
On $\{\tau_R=1\}$, the original and stopped equations agree and $\rho_h=h^{-1/2}\eta$. A second application of Markov's inequality yields, for $u\ge1$,
\[
\mathbb P\left(\sup_s|\rho_h(s)|>CuR,\ \tau_R=1\right)
\le u^{-p}.
\]
Combining the last two estimates proves \eqref{S:eq:short-bridge-probability}. Since $R\ge4$ and $c_0<4r_0$, the assumption $\sqrt hR^2\le c_0$ implies $\sqrt hR<r_0$. Hence the
coordinate calculation is valid throughout the bridge on $\{\tau_R=1\}$. The remainder was defined by the full path difference under the fixed kernel of
Lemma~\ref{S:lem:pinned-gaussian-kernel}; this definition also applies on the exit event and is independent of $p,L,R$.
\end{proof}

\subsubsection{Uniform bridge estimates and the matrix expansion}
\label{S:sec:s9-matrix-expansion}

The preceding estimate uses the same coupling as the localization parameters vary within their stated range. We can therefore let the moment order grow and take a union over a row
while keeping the attached Gaussian bridges fixed. The cutoff $\ell_n$ makes every retained duration small enough for that choice.

For $1\le j<n$, put $v_{n,j}=h_{n,j}^{-1/2} \exp_{Z_{\vartheta_{n,j+1}}}^{-1}Z_{\vartheta_{n,j}}$ in the frame $(e_1,e_2)$ at the initial point. Let $\rho_{n,j}$ be the remainder
of the fixed bridge coupling. Only $\ell_n\le j<n$ enters the following maximum.

\begin{corollary}\label{S:cor:bridge-array-bound}
On the common probability space of Theorem~\ref{S:thm:joint-probability-space}, almost surely for all sufficiently large
$n$,
\begin{equation}\label{S:eq:bridge-array-remainder}
\max_{j\ge\ell_n}\left(|v_{n,j}|+\|\beta_{n,j}\|_\infty\right)
\le C_{\beta,\omega}\sqrt{\log n},
\max_{j\ge\ell_n}\|\rho_{n,j}\|_\infty
\le C_{\beta,\omega}(\log n)^6.
\end{equation}

\end{corollary}
\begin{proof}
Condition on $\mathcal F_{\vartheta_{n,j+1}}^Z$ and apply Lemma~\ref{S:lem:global-hyperbolic-displacement} to the restarted original diffusion. Division by $\sqrt{h_{n,j}}$ gives
$\|v_{n,j}\|_p\le C(\sqrt p+\sqrt{h_{n,j}})$. Choose $p=c_1A_0^2\log n$, with fixed $c_1>0$ small enough that $C\sqrt{c_1}\le1/(4\me)$. For all sufficiently large $n$, the
remaining term is at most $A_0\sqrt{\log n}/(4\me)$. Markov's inequality at $A_0\sqrt{\log n}$ therefore gives
\[
 \mathbb P\left(|v_{n,j}|>A_0\sqrt{\log n}\right)
 \le(2\me)^{-p}\le n^{-c_\beta A_0^2}.
\]
The constants are uniform over $\ell_n\le j<n$. Conditionally on the sampled path, $\beta_{n,j}$ has the law of $W_s-sW_1$, $0\le s\le1$, for a two-dimensional Brownian motion $W$.
Pointwise, $\sup_s|W_s-sW_1|\le2\sup_s|W_s|$. The reflection bound in each coordinate gives $\mathbb P(\|\beta_{n,j}\|_\infty>u\mid\mathscr S_n) \le8e^{-u^2/16}$. With
$u=A_0\sqrt{\log n}$, this is $8n^{-A_0^2/16}$, uniformly in the conditioning data. Choose $A_0$ so that $\min\{c_\beta A_0^2,A_0^2/16\}>4$. A union bound over $j$ yields
\[
 \mathbb P\!\left(\max_{j\ge\ell_n}
  (|v_{n,j}|+\|\beta_{n,j}\|_\infty)>2A_0\sqrt{\log n}\right)
 \le n^{1-c_\beta A_0^2}+8n^{1-A_0^2/16}.
\]
The series in $n$ converges, so Borel--Cantelli proves the first estimate in \eqref{S:eq:bridge-array-remainder}.

For the remainder choose
\[
p_n=\left\lfloor A\frac{\log n}{\log\log n}\right\rfloor,
L_n=A_0\sqrt{\log n},
R_n=(\log n)^2.
\]
We apply Lemma~\ref{S:lem:short-bridge-expansion} with the scalar tail parameter $u=(\log n)^4$. The threshold in \eqref{S:eq:short-bridge-probability} is therefore
$CuR_n=C(\log n)^6$. For all large $n$, $R_n\ge4(1+L_n+\sqrt{p_n})$, and, uniformly for $j\ge\ell_n\ge m_n/2$,
\[
\sqrt{h_{n,j}}R_n^2
\le C_\beta(\log n)^{-252}.
\]
Apply Lemma~\ref{S:lem:short-bridge-expansion}, conditionally on the two boundary points of each sampling interval, on the event $|v_{n,j}|\le L_n$. A union bound over at most $n$
sampling intervals and \eqref{S:eq:short-bridge-probability} give
\[
\begin{aligned}
\mathbb P\left(
\max_{j\ge\ell_n}\|\rho_{n,j}\|_\infty>C(\log n)^6
\right)
\le{}& Cn^{1-c_\beta A_0^2}
 +Cn \me^{-c(\log n)^4}\\
&+Cn\left(C(\log n)^{-256}\right)^{p_n}
 +Cn(\log n)^{-4p_n}.
\end{aligned}
\]
Choose $A_0$ so that $c_\beta A_0^2>4$. Next choose $A>2$. For all sufficiently large $n$,
\[
p_n\ge \frac A2\frac{\log n}{\log\log n}
,\qquad
C(\log n)^{-256}\le(\log n)^{-128}.
\]
Consequently,
\[
 n\left(C(\log n)^{-256}\right)^{p_n}
 \le n^{1-64A},
 n(\log n)^{-4p_n}
 \le n^{1-2A}.
\]
For large $n$, the four terms are bounded by $Cn^{-3}$, $Cn^{-3}$, $Cn^{1-64A}$, and $Cn^{1-2A}$, respectively; for the second term use $c(\log n)^4\ge4\log n$. Their sums over $n$
are finite because $A>2$. Borel--Cantelli gives the second estimate in \eqref{S:eq:bridge-array-remainder}, with threshold $C(\log n)^6$. The product identity
\eqref{S:eq:row-joint-bridge-kernel} gives conditional independence within each row. The path segments remain the
original restrictions of $Z$; only the auxiliary Gaussian bridges are added to the probability space. \end{proof}

Apply the matrix map $\mathcal A$ to the coordinate expansion. Since $(y\partial_\mu\mathcal A)\mathcal A^{-1}=\mathsf A_\mu$, the linear terms have the same two matrix
coefficients on every sampling interval. With
\[
E_{12}=\begin{pmatrix}0&1\\0&0\end{pmatrix},\qquad
H_0=\frac12\begin{pmatrix}1&0\\0&-1\end{pmatrix},
\]
It\^o's formula gives
\[
\dd X_\theta=\left(-E_{12}\,\dd W_\theta^1-H_0\,\dd W_\theta^2+\bigl(\tfrac18I+\tfrac12H_0\bigr)\dd\theta\right)X_\theta.
\]
In Stratonovich form the diffusion matrices are $\mathsf A_1=-E_{12}$ and $\mathsf A_2=-H_0$, while the drift matrix is $\mathsf A_0=H_0/2$. Thus the two coordinates of the
displacements from linear interpolation and from the Brownian bridge act through $\mathsf A_1$ and $\mathsf A_2$. The drift and the higher order terms are included in the remainder
estimated below.

For one sampling interval, let $v_{n,j}\in\mathbb R^2$ be the normalized boundary point increment in Lemma~\ref{S:lem:short-bridge-expansion}, using the orthonormal basis
$(e_1,e_2)$ at $Z_{\vartheta_{n,j+1}}$. Let $\beta_{n,j}$ be the corresponding standard Brownian bridge. Conditionally on $\mathscr S_n$, the family $(\beta_{n,j})_{1\le j<n}$ is
independent.

Recall the clipped coordinates and signed fractions from \eqref{S:eq:interval-relative-time} and Theorem~\ref{S:thm:interval-fractions}. Thus
\[
 \lambda_{n,k,j}(s)=u_{n,j}(\theta_{n,k}(s))-u_{n,j}(\vartheta_{n,k+1}),
   u_{n,j}(\vartheta_{n,k+1})=\mathbf1_{\{j\ge k+1\}}.
\]
Equation~\eqref{S:eq:interval-fractions-error} gives the limits of the fractions for each fixed difference $j-k$; they vanish for $|j-k|>R_\beta$. All sums over $j$ in this
subsection run over $\ell_n\le j<n$. The length identity is \eqref{S:eq:interval-length-partition}. In particular, $u_{n,j}(\vartheta_{n,k+1})$ is zero or one for every $j$, also
when the segment is traversed in the reverse direction.

The displacement obtained by linear interpolation and the displacement from the Gaussian bridge are
\begin{align}
 m_{n,k}(s)&:=\sum_j h_{n,j}^{1/2}\lambda_{n,k,j}(s)v_{n,j},\label{S:eq:bridge-chord-field}\\
 \widehat W_{n,k}(s)&:=\sum_j h_{n,j}^{1/2}\left[\beta_{n,j}\!\left(u_{n,j}(\theta_{n,k}(s))\right)-\beta_{n,j}\!\left(u_{n,j}(\vartheta_{n,k+1})\right)\right].\notag
\end{align}
Since the clipped coordinate of every sampling time is zero or one and $\beta_{n,j}(0)=\beta_{n,j}(1)=0$, completely traversed intervals make no contribution to the second sum.
Define the fields with values in the Lie algebra
\[
\mathcal B_{n,k}(s):=\mathsf A_1m_{n,k}^1(s)+\mathsf A_2m_{n,k}^2(s),
\mathcal Z_{n,k}(s):=\mathsf A_1\widehat W_{n,k}^1(s)+\mathsf A_2\widehat W_{n,k}^2(s).
\]
The field $\mathcal B_n$ is $\mathscr S_n$-measurable. Given $\mathscr S_n$, finite collections of values of $\mathcal Z_n$ are centered Gaussian.
Proposition~\ref{S:prop:weighted-field-covariance} identifies this field as a Gaussian random element of the weighted Hilbert space and computes its covariance.

Let $A_1,\ldots,A_N$ be real square matrices, all of the same size, and put $a=\sum_{i=1}^N\|A_i\|_{\mathrm{op}}$. Products below have the displayed order. We use the logarithm
defined by its power series at $I$.
\begin{lemma}\label{S:lem:estimate-finite-products}
If $a\le a_0:=\log(3/2)$, then
\begin{equation}\label{S:eq:estimate-finite-log-product}
 \left\|\log(\me^{A_N}\cdots\me^{A_1})-\sum_{i=1}^NA_i\right\|
 \le Ca^2,\qquad
 \|\log(\me^{A_N}\cdots\me^{A_1})\|\le Ca,
\end{equation}
with an absolute constant $C$.
\end{lemma}
\begin{proof}

Expand each exponential in its absolutely norm convergent series, keeping the indicated product order. The term of degree zero is $I$ and the sum of terms of degree one is
$\sum_iA_i$. Submultiplicativity and multiplication of the scalar majorants give, for $P=\me^{A_N}\cdots \me^{A_1}$,
\[
 \|P-I\|\le \me^a-1=:b,
 \left\|P-I-\sum_iA_i\right\|
 \le \me^a-1-a\le\tfrac12\me^a a^2.
\]
Take $a_0=\log(3/2)$, so $b\le1/2$. The norm convergent series for $\log(I+(P-I))$ has remainder
\[
 \|\log P-(P-I)\|
 \le\sum_{l=2}^\infty b^l/l
 \le\frac{b^2}{2(1-b)}\le b^2.
\]
Since $b\le \me^{a_0}a$, the two bounds imply \eqref{S:eq:estimate-finite-log-product}. If $A_i=L_i+E_i$, subtracting $\sum_iL_i$ gives
\[
 \left\|\log P-\sum_iL_i\right\|\le
       \sum_i\|E_i\|+Ca^2.
\]
The identity $\sum_iA_i-\sum_iL_i=\sum_iE_i$ gives this bound. For the positive diagonal triangular matrices used below this local logarithm agrees with the unique real triangular
logarithm; the empty product gives zero.

\end{proof}

The left logarithmic derivatives of $\mathcal A$ in the orthonormal frame are the constant matrices $\mathsf A_1,\mathsf A_2$. Consequently the displacements from linear
interpolation and from the bridge enter through the fixed matrices $\mathsf A_1,\mathsf A_2$, without an additional random rotation. Products of the few traversed segments change
this linear sum only quadratically. The next estimate retains the decay in the interval index needed to sum these errors in the weighted norm.

\begin{proposition}\label{S:prop:bridge-affine-expansion}

On an event of probability one, for all sufficiently large $n$ and all retained subintervals,
\begin{equation}\label{S:eq:bridge-joint-expansion}
\mathcal U_{n,k}(s)=\mathcal B_{n,k}(s)+\mathcal Z_{n,k}(s)+\mathcal R_{n,k}(s).
\end{equation}
The remainder satisfies
\begin{equation}\label{S:eq:bridge-remainder-energy}
\mathcal E_n^{+}(\mathcal R_n)\le
 C_{\beta,\omega}\frac{(1+\log n)^{12}}{m_n}\longrightarrow0,
\end{equation}

\end{proposition}
\begin{proof}

Work on the event of probability one in Corollary~\ref{S:cor:bridge-array-bound}. For a sampling interval on which the path starts at $z$, identify $T_z\mathbb H$ with
$\mathbb R^2$ using the orthonormal basis $(e_1(z),e_2(z))$, where $e_1=y\partial_x$ and $e_2=y\partial_y$, and put
\[
\Phi_z(w):=\log\!\left(\mathcal A(\exp_z w)\mathcal A(z)^{-1}\right).
\]
Let $z=x_0+\ii y_0$ and $g_z(w)=(w-x_0)/y_0$. The exact composition identity \eqref{S:eq:affine-section-exact-composition} gives
$\mathcal A(z')\mathcal A(z)^{-1}=\mathcal A(g_z z')$. Moreover $g_z$ is an isometry and $Dg_z$ maps $(e_1(z),e_2(z))$ to $(e_1(\ii),e_2(\ii))$. Thus $g_z(\exp_z w)=\exp_\ii w$
when $w$ is written in these bases. Therefore
\begin{equation}\label{S:eq:estimate-normal-affine-homogeneity}
 \Phi_z(w)=\log\mathcal A(\exp_\ii w)
\end{equation}
in these bases. To calculate this function, put $r=|w|$ and
\[
 S(r)=\frac{\sinh r}{r},  S(0)=1,
 d(w)=\cosh r-w_2S(r).
\]
The normal exponential at $\ii$ has coordinates $x=w_1S(r)/d(w)$ and $y=1/d(w)$: in \eqref{S:eq:7.10}, take $\sin\phi=-w_1/r$ and $\cos\phi=w_2/r$ for $r>0$. At the origin the
absolutely convergent power series
\[
 S(|w|)=\sum_{m=0}^\infty\frac{(w_1^2+w_2^2)^m}{(2m+1)!},
 \qquad
 \cosh|w|=\sum_{m=0}^\infty\frac{(w_1^2+w_2^2)^m}{(2m)!}
\]
give $S(0)=d(0)=1$. Thus the coordinate expressions are smooth near $w=0$ and have values $x=0$, $y=1$, as required by $\exp_\ii0=\ii$. Substituting them in the matrix function
$\mathcal A$ gives
\[
 \mathcal A(\exp_\ii w)
 =\begin{pmatrix}
 \sqrt{d(w)}&-w_1S(r)/\sqrt{d(w)}\\
 0&1/\sqrt{d(w)}
 \end{pmatrix}.
\]
Since $|w_2|\le r$, one has $\me^{-r}\le d(w)\le\me^r$. For $d>0$, set
\[
 f(d):=\int_0^1\frac{\dd t}{1+t(d-1)}
       =\begin{cases}\log d/(d-1),&d\ne1,\\1,&d=1.\end{cases}
\]
The exponential formula for a triangular matrix, or its power series, now gives the exact identity
\begin{equation}\label{S:eq:s9-exact-affine-logarithm}
 \Phi_z(w)=\log d(w)\,H_0-w_1S(r)f(d(w))E_{12}.
\end{equation}
Its diagonal entries are $\pm\tfrac12\log d$, and the upper right entry of its exponential is $-w_1S(r)f(d)(\sqrt d-1/\sqrt d)/\log d =-w_1S(r)/\sqrt d$, with the continuous value
when $d=1$. For $0\le j\le3$, differentiating the integral defining $f$ gives
\[
 f^{(j)}(d)=(-1)^j j!\int_0^1
       \frac{t^j\,\dd t}{[1+t(d-1)]^{j+1}}.
\]
These derivatives are uniformly bounded for $\me^{-r_0}\le d\le \me^{r_0}$, for a fixed $r_0>0$. The series for $S(r)$ and $\cosh r$ are convergent power series in $r^2$. Their
derivatives in $w$, through order three, are bounded on $|w|\le r_0$. Consequently \eqref{S:eq:s9-exact-affine-logarithm} has a uniform bound for derivatives of order three there,
independent of $z$. Taylor's formula gives
\[
 \begin{aligned}
 d(w)&=1-w_2+\tfrac12(w_1^2+w_2^2)+O(|w|^3),\\
 \log d(w)&=-w_2+\tfrac12w_1^2+O(|w|^3),\\
 S(r)f(d(w))&=1+\tfrac12w_2+O(|w|^2).
 \end{aligned}
\]
Here, for example, the second line follows by substituting the first into $\log(1+x)=x-x^2/2+O(|x|^3)$; the $w_2^2/2$ terms cancel. Thus the quadratic expansion is
\begin{equation}\label{S:eq:s9-quadratic-affine-logarithm}
 \Phi_z(w)=\mathsf A_1w_1+\mathsf A_2w_2
       +\tfrac12w_1^2H_0-\tfrac12w_1w_2E_{12}+O(|w|^3).
\end{equation}
In particular $D\Phi_z(0)e_\mu=\mathsf A_\mu$, and the form needed below is
\begin{equation}\label{S:eq:estimate-normal-affine-taylor}
 \Phi_z(w)=\mathsf A_1w_1+\mathsf A_2w_2+\mathcal Q_z(w),
 \qquad \|\mathcal Q_z(w)\|\le C|w|^2,\quad |w|\le r_0.
\end{equation}
Every remainder here is bounded by Taylor's integral remainder on this fixed ball. For the $j$-th conditioned bridge, Corollary~\ref{S:cor:bridge-array-bound} and the expansion of
a bridge on a short interval give
\[
w_{n,j}(u)=h_{n,j}^{1/2}\bigl(uv_{n,j}+\beta_{n,j}(u)\bigr)
+h_{n,j}\rho_{n,j}(u).
\]
Uniformly in $u$ and $j\ge\ell_n$,
\[
|w_{n,j}(u)|
\le C_{\beta,\omega}\left(
\sqrt{\frac{\log n}{j}}+\frac{(\log n)^6}{j}
\right)<r_0
\]
for all large $n$. Taylor's formula therefore gives
\[
\left\|
\Phi_z(w_{n,j}(u))
-\sum_{\mu=1}^2\mathsf A_\mu h_{n,j}^{1/2}
\bigl(uv_{n,j}^\mu+\beta_{n,j}^\mu(u)\bigr)
\right\|
\le C_{\beta,\omega}\frac{(\log n)^6}{j}.
\]
The term $h_{n,j}\rho_{n,j}$ is bounded by $C(\log n)^6/j$, while the quadratic remainder of the expansion is $O((\log n)/j)$.

For fixed $k,s$, the segment from $\vartheta_{n,k+1}$ to $\theta_{n,k}(s)$ passes through at most $2R_\beta+1$ intervals $J_{n,j}$. On $J_{n,j}$ write
$\Phi_{n,j}(u)=\Phi_{z_{n,j}}(w_{n,j}(u))$, where $z_{n,j}=Z_{\vartheta_{n,j+1}}$. The matrix relating the positions at relative times $u_0$ and $u_1$ is
\[
\mathcal A(Z_{n,j}(u_1))\mathcal A(Z_{n,j}(u_0))^{-1}
=\me^{\Phi_{n,j}(u_1)}\me^{-\Phi_{n,j}(u_0)}.
\]
Apply Lemma~\ref{S:lem:estimate-finite-products} to $A_1=-\Phi_{n,j}(u_0)$ and $A_2=\Phi_{n,j}(u_1)$. Their total norm is $O_{\beta,\omega}(\sqrt{(\log n)/j})$, uniformly in the
boundary points, and is, for all sufficiently large $n$, at most $a_0$. Together with \eqref{S:eq:estimate-normal-affine-taylor}, this gives
\[
 \log\left(\me^{\Phi_{n,j}(u_1)}\me^{-\Phi_{n,j}(u_0)}\right) ={}\sum_{\mu=1}^2\mathsf A_\mu h_{n,j}^{1/2} \Bigl[(u_1-u_0)v_{n,j}^\mu +\beta_{n,j}^\mu(u_1)-\beta_{n,j}^\mu(u_0)\Bigr] +\mathcal R_{n,j}(u_0,u_1),
 \]
where $\|\mathcal R_{n,j}(u_0,u_1)\|\le C_{\beta,\omega}(\log n)^6/j$. The quadratic Taylor errors and the entire error in the logarithm of the product from
\eqref{S:eq:estimate-finite-log-product} are $O_\omega((\log n)/j)$; the term containing $\rho_{n,j}(u_1)-\rho_{n,j}(u_0)$ is $O_\omega((\log n)^6/j)$.

Compose these segment transfers in their traversal order. Their matrix products telescope to
\[
 \mathcal A(Z_{\theta_{n,k}(s)})
 \mathcal A(Z_{\vartheta_{n,k+1}})^{-1},
\]
also when Brownian time decreases. There are at most $2R_\beta+1$ nonidentity factors and their indices satisfy $j\asymp_\beta k$. Writing their logarithms as $L_i+E_i$, the
preceding estimates give
\[
 \sum_i\|L_i\|\le C_{\beta,\omega}\sqrt{(\log n)/k},
 \sum_i\|E_i\|\le C_{\beta,\omega}(\log n)^6/k.
\]
Their total norm tends to zero uniformly for $k\ge m_n$. The square of that total is bounded by $C_{\beta,\omega}[(\log n)/k+(\log n)^{12}/k^2]$, which is at most
$C_{\beta,\omega}(\log n)^6/k$ for all sufficiently large $n$, because $k\ge(\log n)^{512}$. Applying Lemma~\ref{S:lem:estimate-finite-products} a second time bounds the full
logarithm minus $\sum_iL_i$ by the same quantity. The linear sum is exactly $\mathcal B_{n,k}(s)+\mathcal Z_{n,k}(s)$: the signed fractions are the differences between sampled
values of the clipped coordinates and each bridge difference has the same two boundary points. This proves \eqref{S:eq:bridge-joint-expansion} with
\[
\sup_{0\le s\le1}\|\mathcal R_{n,k}(s)\|
\le C_{\beta,\omega}\frac{(\log n)^6}{k}.
\]
Consequently,
\[
\mathcal E_n(\mathcal R_n)
\le C_{\beta,\omega}(\log n)^{12}
\sum_{k=m_n}^{n-1}\frac{\mathfrak G_{n,k}}{k^2}
\le C_{\beta,\omega}\frac{(\log n)^{12}}{m_n}
\sum_{k=m_n}^{n-1}\frac{\mathfrak G_{n,k}}k.
\]
For the diagonal part, the definition of $\mathcal E_n^{\mathrm{diag}}$ and the inequality $|Rv|\le\|R\|_{\mathrm{op}}|v|$ give
\[
 \mathcal E_n^{\mathrm{diag}}(\mathcal R_n)
 \le\frac{C_{\beta,\omega}(\log n)^{12}}n
       \sum_{k=m_n}^{n-1}\frac{p_{n,k}q_{n,k}}{k^2}
 \le C_{\beta,\omega,\sigma}(\log n)^{12}
               n^{-1+\sigma}m_n^{-1-\sigma}.
\]
Here $0<\sigma<1$ is fixed, and the last inequality is the product sum in Section~\ref{S:sec:s4-sampled-columns} with $d=2$. Since $n^{-1+\sigma}m_n^{-\sigma}\le1$, this is at most
$C_{\beta,\omega,\sigma}(\log n)^{12}/m_n$. Together with Proposition~\ref{S:prop:5.5}, the off-diagonal estimate proves \eqref{S:eq:bridge-remainder-energy} for the full quadratic
form. All terms in this estimate use the same restricted operator cells; the additional sampling intervals in Brownian time with $\ell_n\le j<m_n$ have already been included in the
pointwise bound through $j\asymp_\beta k$. Finally, the linear part is bounded by $C_\omega\sqrt{(\log n)/k}$, whereas the remainder is bounded by $C_\omega(\log n)^6/k$. Since
$k\ge m_n\ge(\log n)^{512}$, the latter is smaller for all sufficiently large $n$. Combining the linear and remainder bounds gives
\begin{equation}\label{S:eq:bridge-uniform-smallness}
 \sup_{m_n\le k<n}\sup_{0\le s\le1}\|\mathcal U_{n,k}(s)\|
 \le C_{\beta,\omega}\sqrt{\frac{\log n}{m_n}}\longrightarrow0.
\end{equation}
\end{proof}

Write the boundary points of $J_{n,j}$ as
\[
 s_{n,j}=\vartheta_{n,j+1},
 t_{n,j}=\vartheta_{n,j},  h_{n,j}=t_{n,j}-s_{n,j},
 \varrho_{n,j}=\me^{-\beta s_{n,j}/4}.
\]
The normal increment $\Delta_{n,j}=\sqrt{h_{n,j}}v_{n,j}$ is $\exp_{Z_{s_{n,j}}}^{-1}Z_{t_{n,j}}$, expressed in the frame $(e_1,e_2)$ at $Z_{s_{n,j}}$. In the following comparison,
$B=W$ is the Brownian motion driving the original diffusion in the upper half-plane and $\Delta B_{n,j}=B_{t_{n,j}}-B_{s_{n,j}}$. These increments are independent Gaussians within
each row under the original law. Under the law conditioned on the terminal point $\eta$, we use $\widetilde W^\eta$ from \eqref{S:eq:conditioned-original-affine-sde}. The pinned
driver of \eqref{S:eq:pinned-kernel-sde} is used only for the preceding bridge comparison.

\subsubsection{Normal increments and the driving Brownian motion}
\label{S:sec:s9-normal-increments}
The bridge approximation and the observed increments use different descriptions of the diffusion matrix. To compare the increments pathwise, the original driving Brownian motion
must be retained. The calculation below identifies its coefficients in normal coordinates.

Restart the original diffusion at a deterministic point $z=x_0+\ii y_0$. In the frame $(e_1(z),e_2(z))$, write $Y_t=\exp_z^{-1}Z_t$. The Brownian motion $B$ in this subsubsection
is the original driving Brownian motion; under the law conditioned on $\eta$, it is replaced by $\widetilde W^\eta$. The drift $b_0$ and the inverse metric $a$ are those computed
in \eqref{S:eq:s9-normal-metric}--\eqref{S:eq:s9-normal-drift}.
\begin{lemma}\label{S:lem:affine-driver-normal-coefficients}
On a fixed ball $|Y|\le r_0$, the two equations are
\[
 \dd Y=\sigma_{\mathrm{aff}}(Y)\dd B+b_0(Y)\dd t,
 \qquad
 \dd Y=\sigma_{\mathrm{aff}}(Y)\dd\widetilde W^\eta
                              +b_{z,\eta}(Y)\dd t.
\]
Uniformly in $z,\eta$,
\begin{equation}\label{S:eq:affine-normal-local-coefficients}
 \|\sigma_{\mathrm{aff}}(Y)\|\le C,\qquad
 \|\sigma_{\mathrm{aff}}(Y)-I\|\le C|Y|,\qquad
 |b_0(Y)|+|b_{z,\eta}(Y)|\le C.
\end{equation}
\end{lemma}
\begin{proof}
The isometry $g_z(w)=(w-x_0)/y_0$ satisfies $Dg_z\,e_\mu(w)=e_\mu(g_zw)$. Applying it preserves the original driving coordinates themselves. In normal coordinates at $\ii$, write
$\exp_\ii Y=(x(Y),y(Y))$. The formula used in \eqref{S:eq:s9-exact-affine-logarithm} gives
\[
 \begin{aligned}
 x(Y)&=Y_1+Y_1Y_2+O(|Y|^3),\\
 y(Y)&=1+Y_2+\tfrac12(Y_2^2-Y_1^2)+O(|Y|^3),\\
 D\exp_\ii(Y)&=
 \begin{pmatrix}1+Y_2&Y_1\\-Y_1&1+Y_2\end{pmatrix}
       +O(|Y|^2).
 \end{aligned}
\]
Taylor's integral formula bounds the derivatives of the remainders on a smaller fixed ball. The diffusion matrix in the equation driven by $B$ is consequently
\begin{equation}\label{S:eq:s9-affine-driver-rotation}
 \sigma_{\mathrm{aff}}(Y)
 =y(Y)(D\exp_\ii(Y))^{-1}
 =I+Y_1J+O(|Y|^2).
\end{equation}
For the inverse in this computation, write $D\exp_\ii=I+Y_2I-Y_1J+E$, $\|E\|\le C|Y|^2$, and use $(I+A)^{-1}=I-A+A^2(I+A)^{-1}$ with $\|A\|\le1/2$.

It\^o's formula gives drift $\tfrac12\Delta_{\mathbb H}\log_\ii=b_0$. The explicit expression \eqref{S:eq:s9-normal-drift} bounds it by $C|Y|$ on this ball. Under terminal
conditioning the additional drift is
\[
 \sigma_{\mathrm{aff}}(Y)b_{g_z\eta}(\exp_\ii Y),
\]
where $b_\xi$ is the two-dimensional drift obtained from the Poisson kernel of Proposition~\ref{S:prop:aux-stochastic-normalization}. Its Euclidean norm is one, independently of
$\xi$. Thus $|b_{z,\eta}(Y)|\le |b_0(Y)|+\|\sigma_{\mathrm{aff}}(Y)\|\,|b_{g_z\eta}(\exp_\ii Y)|\le C$. Finally \eqref{S:eq:s9-affine-driver-rotation} generally has a nonzero
linear term, whereas $a^{1/2}=I+O(|Y|^2)$. The increment comparison uses $\sigma_{\mathrm{aff}}$, with the original driver; the symmetric square root represents the same generator
with a different driver.
\end{proof}

\begin{lemma}\label{S:lem:qv-normal-increments}
There are deterministic $h_0,c_0,C>0$ such that, whenever $p\ge2$, $h=h_{n,j}\le h_0$, and $ph\le c_0$,
\begin{equation}\label{S:eq:normal-increment-conditional-moment}
 \mathbb E\left[|\Delta_{n,j}-\Delta B_{n,j}|^p
                  \mid\mathcal F_{s_{n,j}}^Z\right]^{1/p}
 \le Cph.
\end{equation}
Under $\mathbb P^\eta$, with its natural augmented filtration,
\begin{equation}\label{S:eq:terminal-normal-increment-moment}
 \mathbb E^\eta\left[|\Delta_{n,j}-\Delta\widetilde W^\eta_{n,j}|^p
                  \mid\mathcal F_{s_{n,j}}^Z\right]^{1/p}
 \le Cph.
\end{equation}
The constant is independent of the terminal point and the starting state. Equivalently, this estimate holds under the joint law conditioned on
$\mathcal F_{s_{n,j}}^Z\vee\sigma(\eta_\infty)$.

\end{lemma}
\begin{proof}
Restart at time $s=s_{n,j}$ and at the observed point $z=Z_s$. Write $D_t=B_{s+t}-B_s$ in the original law and
$D_t=\widetilde W_{s+t}^\eta-\widetilde W_s^\eta$ in the conditioned law. In both cases $D$ is Brownian under the
conditional law given the past. Lemma~\ref{S:lem:affine-driver-normal-coefficients} gives, before $\tau=\inf\{t:|Y_t|\ge
r_0\}$,
\[
 \dd Y_t=\sigma_{\mathrm{aff}}(Y_t)\dd D_t+b_*(Y_t)\dd t,
   |b_*|\le C.
\]
All expectations in this proof are under this conditional law. The global identity $|Y_t|=d(z,Z_{s+t})$ and Lemma~\ref{S:lem:global-hyperbolic-displacement} imply
\begin{equation}\label{S:eq:global-normal-and-driver-moments}
 \|\sup_{t\le h}|Y_t|\|_{L^{2p}}
 +\|\sup_{t\le h}|D_t|\|_{L^{2p}}
 \le C(\sqrt{ph}+h).
\end{equation}
These estimates apply to the unstopped variables.

Define the stopped error integral
\[
 E_h^{\mathrm{stop}}
 =\int_0^h\mathbf1_{\{t<\tau\}}
          (\sigma_{\mathrm{aff}}(Y_t)-I)\dd D_t
  +\int_0^h\mathbf1_{\{t<\tau\}}b_*(Y_t)\dd t.
\]
It equals $Y_h-D_h$ on $\{\tau>h\}$. The inequality for continuous martingales with its $\sqrt p$ constant, Minkowski in $L^{p/2}$, and the global displacement bound at time $t$
give
\begin{equation}\label{S:eq:stopped-normal-error-moment}
 \begin{aligned}
 \|E_h^{\mathrm{stop}}\|_{L^p}
 &\le C\sqrt p\left(\int_0^h
             \|\mathbf1_{\{t<\tau\}}Y_t\|_{L^p}^2\dd t\right)^{1/2}+Ch\\
 &\le C\sqrt p\left(\int_0^h(pt+t^2)\dd t\right)^{1/2}+Ch
 \le Cph.
 \end{aligned}
\end{equation}
Indeed,
\[
 \sqrt p\left(\int_0^h(pt+t^2)\,\dd t\right)^{1/2}
 \le\frac{ph}{\sqrt2}+\frac{\sqrt p\,h^{3/2}}{\sqrt3}.
\]
Since $0<h\le h_0\le1$ and $p\ge2$, one has $\sqrt p\,h^{3/2}/(ph)=\sqrt{h/p}\le1$ and $h/(ph)=1/p\le1$. These inequalities give the last bound in
\eqref{S:eq:stopped-normal-error-moment}.

Let $E_h=Y_h-D_h$ without stopping. H\"older's inequality, \eqref{S:eq:fixed-normal-ball-exit}, and \eqref{S:eq:global-normal-and-driver-moments} imply
\begin{equation}\label{S:eq:normal-exit-hoelder}
 \|\mathbf1_{\{\tau\le h\}}E_h\|_{L^p} \le\mathbb P(\tau\le h)^{1/(2p)}\|E_h\|_{L^{2p}} \le C \me^{-c/(ph)}(\sqrt{ph}+h) \le C'ph.
\end{equation}
For $x=ph$ and $u=c/x$, $x^{-1/2}\me^{-c/x}=c^{-1/2}\sqrt u\,\me^{-u}\le(2c\me)^{-1/2}$; also $h\le x$. These bounds give the last inequality.

The identity $E_h=\mathbf1_{\{\tau>h\}}E_h^{\mathrm{stop}} +\mathbf1_{\{\tau\le h\}}E_h$ implies
\[
 \|E_h\|_{L^p}
 \le\|E_h^{\mathrm{stop}}\|_{L^p}
       +\|\mathbf1_{\{\tau\le h\}}E_h\|_{L^p}\le Cph.
\]
The preceding estimates were made under either of the two conditional laws specified at the start of the proof. Since $Y_h=\Delta_{n,j}$ and the corresponding choices of $D_h$ are
$\Delta B_{n,j}$ and $\Delta\widetilde W^\eta_{n,j}$, they give
\[
 \begin{aligned}
 \mathbb E[|\Delta_{n,j}-\Delta B_{n,j}|^p
                        \mid\mathcal F_s^Z]&\le(Cph)^p,\\
 \mathbb E^\eta[|\Delta_{n,j}-\Delta\widetilde W^\eta_{n,j}|^p
                        \mid\mathcal F_s^Z]&\le(Cph)^p.
 \end{aligned}
\]
In the second line the drift, displacement and exit constants are uniform in $\eta$. Taking $p$th roots proves the two claimed inequalities. For the joint conditioning statement,
let $\nu_\infty$ be the law of $\eta_\infty$. For bounded nonnegative $F$ measurable with respect to $\mathcal F_s^Z$ and bounded nonnegative $g$, disintegration and the second
inequality give
\[
 \begin{aligned}
 &\mathbb E\!\big[Fg(\eta_\infty)
      |\Delta_{n,j}-\Delta\widetilde W^{\eta_\infty}_{n,j}|^p\big]\\
 &\quad=\int g(\eta)\mathbb E^\eta\!\big[F
      |\Delta_{n,j}-\Delta\widetilde W^\eta_{n,j}|^p\big]
                                      \,\nu_\infty(\dd\eta)\\
 &\quad\le(Cph)^p\mathbb E[Fg(\eta_\infty)].
 \end{aligned}
\]
The products $Fg(\eta_\infty)$ generate the joined sigma-field. A monotone class argument extends the inequality to every bounded nonnegative function measurable with respect to
that sigma-field, which is the required conditional moment bound given the past and terminal
point. \end{proof}

\begin{corollary}\label{S:cor:normal-increment-full-array}
On the probability space of Theorem~\ref{S:thm:joint-probability-space} there is an event $\Omega_{\mathrm{incr}}$ of probability one such that, for all sufficiently large $n$, the
following bounds hold for every sampling interval
$\ell_n\le j<n$,
\begin{equation}\label{S:eq:qv-normal-comparison}
 \Delta_{n,j}=\Delta B_{n,j}+R_{n,j},
 |R_{n,j}|\le C_\beta h_{n,j}\log n,
 |\Delta B_{n,j}|\le C_\beta\sqrt{h_{n,j}\log n},
\end{equation}
and, writing $\widetilde W=\widetilde W^{\eta_\infty}$,
\begin{equation}\label{S:eq:terminal-normal-full-array}
 |\Delta_{n,j}-\Delta\widetilde W_{n,j}|
       \le C_\beta h_{n,j}\log n,
 |\Delta\widetilde W_{n,j}|\le C_\beta\sqrt{h_{n,j}\log n}.
\end{equation}
The constants may be deterministic; only the starting index is random.
\end{corollary}
\begin{proof}
Take $p_n=8\log n$ for all sufficiently large $n$. Since $h_{n,j}\le C_\beta/\ell_n$ and $\ell_n\ge m_n/2$, both $h_{n,j}\le h_0$ and $p_nh_{n,j}\le c_0$ hold uniformly on the
array over the enlarged index range. Markov's inequality applied at $\me Cp_nh_{n,j}$ gives probability at most $\me^{-p_n}=n^{-8}$ for each error in
Lemma~\ref{S:lem:qv-normal-increments}. The Brownian increment moment $\|D_h\|_{L^p}\le C\sqrt{ph}$ gives the same bound at $\me C\sqrt{p_nh}$ for its driver. In the conditioned
calculation integrate first over the past and then over the terminal point. Let $E_n$ be the event that all four bounds in \eqref{S:eq:qv-normal-comparison} and
\eqref{S:eq:terminal-normal-full-array} hold throughout row $n$. A union bound over at most $n$ sampling intervals gives $\mathbb P(E_n^c)\le4n^{-7}$. Hence the sum over $n$ is
finite and Borel--Cantelli shows that $E_n$ occurs for all sufficiently large $n$ almost surely. For a sufficiently large deterministic $N_0(\beta)$,
\[
 \mathbb P\!\left(\bigcup_{n\ge N}E_n^c\right)
 \le C_\beta N^{-6},  N\ge N_0(\beta).
\]
The pathwise identity $W_t-\widetilde W_t=\int_0^t b_{\eta_\infty}(Z_s)\dd s$ also gives $|\Delta B_{n,j}-\Delta\widetilde W_{n,j}|\le h_{n,j}$. It verifies consistency of the two
comparisons, without claiming that $B$ has independent centered increments after conditioning on $\eta_\infty$. All events concerning the original increments have the same
probabilities on this enlarged space. For the resulting event $E$, $0\le\mathbb P(E\mid\eta_\infty)\le1$ and $\mathbb E\mathbb P(E\mid\eta_\infty)=\mathbb P(E)=1$. Hence
$\mathbb P(E\mid\eta_\infty)=1$ almost surely.
\end{proof}

\subsubsection{Conditional localization}
\label{S:sec:conditional-tail-bounds-and-localization}

The constants in the conditional estimates below are determined by the sampled positions and are finite almost surely. We first restrict to events where these constants are bounded
by a fixed positive integer. These events are measurable with respect to the current row. We use the following form of Borel--Cantelli.

For each $n$, let $\mathcal I_n$ be a finite set. Suppose that $Y_{n,i}\ge0$ is measurable, $u_{n,i}>0$ is $\mathcal H_{n,i}$-measurable, and $A_{n,i}(M)\in\mathcal H_{n,i}$ for
each integer $M\ge1$. Assume that deterministic numbers $b_{n,i}(M)\ge0$ satisfy
\begin{equation}\label{S:eq:conditional-localized-tail-hypothesis}
 \mathbf1_{A_{n,i}(M)}
 \mathbb P(Y_{n,i}>u_{n,i}\mid\mathcal H_{n,i})
 \le b_{n,i}(M),\qquad
 \sum_n\sum_{i\in\mathcal I_n}b_{n,i}(M)<\infty.
\end{equation}

\begin{lemma}\label{S:lem:conditional-localized-borel-cantelli}
Almost surely, for every integer $M$, only finitely many rows contain an event $A_{n,i}(M)\cap\{Y_{n,i}>u_{n,i}\}$.
\end{lemma}
\begin{proof}
Since $A_{n,i}(M)\in\mathcal H_{n,i}$, the tower property in Theorem~8.1 of \cite{kallenberg2021foundations} gives
\begin{equation}\label{S:eq:conditional-localization-tower}
 \mathbb P(A_{n,i}(M)\cap\{Y_{n,i}>u_{n,i}\}) =\mathbb E\left[\mathbf1_{A_{n,i}(M)} \mathbb P(Y_{n,i}>u_{n,i}\mid\mathcal H_{n,i})\right] \le b_{n,i}(M).
\end{equation}
A union bound over $i\in\mathcal I_n$ and the assumed summability allow the first Borel--Cantelli lemma
\cite{kallenberg2021foundations} (Theorem~4.18) to apply for each $M$. Take the countable intersection over $M\in\mathbb
N$.
\end{proof}

On this single event, suppose that finite integers $M_\omega,N_\omega$ satisfy $\omega\in A_{n,i}(M_\omega)$ for every $n\ge N_\omega$ and every $i\in\mathcal I_n$. The bound for
the integer $M_\omega$ then gives $Y_{n,i}\le u_{n,i}$ throughout every sufficiently large row.

The same conclusion can be obtained under a terminal conditional law. Indeed, if $\mathcal H_\eta\subset\mathcal H_{n,i}$, then
\[
 \mathbb E[\mathbf1_{A_{n,i}(M)} \mathbf1_{\{Y_{n,i}>u_{n,i}\}}\mid\mathcal H_\eta] =\mathbb E\left[\mathbf1_{A_{n,i}(M)} \mathbb P(Y_{n,i}>u_{n,i}\mid\mathcal H_{n,i}) \mid\mathcal H_\eta\right] \le b_{n,i}(M).
\]
The union bound also applies under almost every terminal conditional law. Since $n,i,M$ range over countable sets, the versions can be fixed on one event. In the Gaussian
applications, each localization event belongs to $\mathscr S_n$, so the tower identity \eqref{S:eq:conditional-localization-tower} applies.

\subsubsection{Gaussian series and quadratic forms}
\label{S:sec:s9-gaussian-quadratic-forms}

The local bridge has a sine series with Gaussian coefficients. Applying the kernel derivative requires convergence of this expansion in the weighted field norm. The estimates for
quadratic forms that follow will be used both for weighted sums of increments with deterministic coefficients and for conditional estimates of squared weighted norms of fields.

Let $\mathscr K$ be a separable real Hilbert space and $\mathcal G$ admit a regular conditional law. Let $v_i$ be strongly $\mathcal G$-measurable, and let $g_i$ be conditionally
independent standard normals. Assume that, almost surely,
\[
 V:=\sum_{i\ge1}\|v_i\|_{\mathscr K}^2<\infty.
\]
Write $S_N=\sum_{i=1}^Ng_iv_i$ and $(v\otimes w)x=\langle x,w\rangle v$.

\begin{lemma}\label{S:lem:conditional-gaussian-series}

Conditionally on $\mathcal G$, $S_N$ converges in $L^2(\mathscr K)$ and almost surely to a strongly measurable $S$ with law $\mathcal N_{\mathscr K}(0,Q)$. The covariance operators
satisfy, with convergence in $\mathfrak S_1(\mathscr K)$,
\begin{equation}\label{S:eq:conditional-series-covariance}
 Q_N:=\sum_{i=1}^Nv_i\otimes v_i\longrightarrow
 Q:=\sum_{i\ge1}v_i\otimes v_i,
 \qquad 0\le Q_N\le Q_{N+1}\le Q.
\end{equation}

\end{lemma}
\begin{proof}

Fix a conditional law for which the hypotheses hold. In this proof the vectors $v_i$ are deterministic and $V$ is a finite number. Independence and centering give, for $K>N$,
\[
 \mathbb E\|S_K-S_N\|^2=\sum_{N<i\le K}\|v_i\|^2.
\]
Thus the partial sums are Cauchy in $L^2(\mathscr K)$. For fixed $N$, $\|S_r-S_N\|$, $r\ge N$, is a nonnegative submartingale by conditional Jensen's inequality. Doob's $L^2$
maximal inequality, Theorem~9.17 of \cite{kallenberg2021foundations}, gives

\[
 \mathbb E\max_{N\le r\le K}\|S_r-S_N\|^2
 \le4\sum_{i>N}\|v_i\|^2.
\]
Choose integers $N_l$ with $\sum_{i>N_l}\|v_i\|^2\le2^{-4l}$ and let $K\to\infty$ by monotone convergence. Markov gives
\[
 \mathbb P\!\left(\sup_{r\ge N_l}\|S_r-S_{N_l}\|>2^{-l}\right)
 \le4\,2^{-2l}.
\]
Borel--Cantelli and completeness prove almost sure convergence of the entire sequence. This limit agrees with the $L^2$ limit. Under the original law the event of convergence has
conditional probability one and hence probability one. Defining $S$ to be that limit there and zero elsewhere makes it strongly measurable. Equivalently, all unconditional $L^2$
arguments can be localized on the $\mathcal G$-measurable sets $\{V\le L\}$, $L\in\mathbb N$.

For a rank-one operator, $(v\otimes w)^*(v\otimes w)=\|v\|^2(w\otimes w)$. Its only possible nonzero singular value is $\|v\|\|w\|$, so $\|v\otimes w\|_{\mathfrak S_1}=\|v\|\|w\|$.
Therefore the series in \eqref{S:eq:conditional-series-covariance} is absolutely convergent in trace norm. Its partial sums and tails are positive, and continuity of the trace
gives
\[
 \operatorname{tr}Q=\sum_i\|v_i\|^2,
   \|Q-Q_N\|_{\mathfrak S_1}
       =\operatorname{tr}(Q-Q_N)=\sum_{i>N}\|v_i\|^2.
\]
Passing to the $L^2$ limit gives the tail identity
\begin{equation}\label{S:eq:conditional-series-tail}
 \mathbb E[\|S-S_N\|_{\mathscr K}^2\mid\mathcal G]
 =\|Q-Q_N\|_{\mathfrak S_1(\mathscr K)}
 =\sum_{i>N}\|v_i\|_{\mathscr K}^2.
\end{equation}
  Throughout this fixed conditional law, $\mathbb E\|S\|^2=\operatorname{tr}Q$. Put $r_N=\sum_{i>N}\|v_i\|^2$. Under this conditional law, $\|S-S_N\|_{L^2}=\sqrt{r_N}$ and
  $\|S_N\|_{L^2},\|S\|_{L^2}\le\sqrt V$. Hence
\[
 \mathbb E\left|\langle S_N,x\rangle\langle S_N,y\rangle
              -\langle S,x\rangle\langle S,y\rangle\right|
 \le2\|x\|\|y\|\sqrt{Vr_N}\longrightarrow0.
\]
Passing the finite covariance identity to the limit gives
\[
 \mathbb E[\langle S,x\rangle\langle S,y\rangle]
 =\langle Qx,y\rangle.
\]
For any finite family $x_1,\ldots,x_r$ and real $t_1,\ldots,t_r$, independence gives
\[
 \mathbb E\exp\!\left(\mathrm i\sum_{a=1}^rt_a\langle S_N,x_a\rangle\right)
 =\exp\!\left(-\frac12\sum_{a,b=1}^rt_at_b\langle Q_Nx_a,x_b\rangle\right).
\]
For $x=\sum_at_ax_a$, the two limiting errors satisfy
\[
 \left|\mathbb E\me^{\mathrm i\langle S_N,x\rangle}
              -\mathbb E\me^{\mathrm i\langle S,x\rangle}\right|
 \le\|x\|\sqrt{r_N},\qquad
 \left|\me^{-\langle Q_Nx,x\rangle/2}-\me^{-\langle Qx,x\rangle/2}\right|
 \le\tfrac12\|x\|^2r_N.
\]
Thus the conditional characteristic identity holds for each $x$. Choose a countable orthonormal basis and let $D_0$ be its rational linear span. Remove the union of the exceptional
sets for $x\in D_0$. For a remaining conditioning value $\omega$, let $K_\omega$ be the conditional law of $S$. Its second moment is $V(\omega)$, and
\[
 \left|\int\me^{\mathrm i\langle v,x\rangle}K_\omega(\dd v)
             -\int\me^{\mathrm i\langle v,y\rangle}K_\omega(\dd v)\right|
 \le\sqrt{V(\omega)}\,\|x-y\|,
\]
while
\[
 \left|\me^{-\langle Qx,x\rangle/2}-\me^{-\langle Qy,y\rangle/2}\right|
 \le\tfrac12\|Q\|\,\|x-y\|\,\|x+y\|.
\]
Approximating $x$ by elements of $D_0$ yields
\[
 \int\me^{\mathrm i\langle v,x\rangle}K_\omega(\dd v)
 =\me^{-\langle Q(\omega)x,x\rangle/2},\qquad x\in\mathscr K.
\]
The basis coordinates generate the Borel sigma-field, so these characteristic functions identify $K_\omega=\mathcal N_{\mathscr K}(0,Q(\omega))$. Finally, under the fixed
conditional law,
\[
 S-S_N=\lim_{K\to\infty}\sum_{i=N+1}^Kg_iv_i
\]
is measurable with respect to the completed sigma-field of $(g_i)_{i>N}$. Its independence of $(g_i)_{i\le N}$, and hence of $S_N$, follows from the conditional product law.

Finally $Q_N$ is strongly $\mathcal G$-measurable in the separable trace-class space, and its limit in trace norm is measurable. Conditional expectations in this lemma mean
integration under the regular conditional law. On $A_L=\{V\le L\}\in\mathcal G$, Tonelli gives
\[
 \mathbb E[\mathbf1_{A_L}\|S\|^2]
 =\mathbb E[\mathbf1_{A_L}V]\le L,
 \qquad \|S\otimes S\|_{\mathfrak S_1}=\|S\|^2.
\]
Thus $\mathbf1_{A_L}S$ and $\mathbf1_{A_L}(S\otimes S)$ have Bochner integrals. The sets $A_L$ increase to an event of probability one. All estimates above hold for any fixed
enumeration, including successive finite groups of modes.

\end{proof}

Let $G$ be a centered Gaussian element of a separable real Hilbert space, with positive trace-class covariance $Q$, and let $b$ be a fixed vector. For this calculation, put
\[
 T=\operatorname{tr}Q,\qquad \Lambda=\|Q\|,\qquad
 V=\operatorname{tr}(Q^2),\qquad v_b=\langle Qb,b\rangle,
 \qquad Z_b=\|G+b\|^2-T-\|b\|^2.
\]
The case $b=0$ of the following calculation is the weighted chi-square estimate of equations~(4.1)--(4.2) in Lemma~1 of \cite{laurent2000adaptive}. We prove the noncentered
estimate and its infinite-dimensional version. The maximal estimate for partial sums is proved in Section~\ref{S:sec:s11-gaussian-quadratic-forms}.

\begin{lemma}\label{S:lem:covariance-gaussian-concentration}
For every $x>0$,
\begin{equation}\label{S:eq:covariance-gaussian-explicit-tails}
 \begin{aligned}
 \mathbb P\bigl(|\|G\|^2-T|>2\sqrt{Vx}+2\Lambda x\bigr)
       &\le2\me^{-x},\\
 \mathbb P\bigl(|\langle G,b\rangle|>\sqrt{2v_bx}\bigr)
       &\le2\me^{-x},\\
 \mathbb P\bigl(|Z_b|>2\sqrt{(V+2v_b)x}+2\Lambda x\bigr)
       &\le2\me^{-x}.
 \end{aligned}
\end{equation}
These inequalities hold under any regular conditional Gaussian distribution under which $Q$ and $b$ are fixed.
\end{lemma}
\begin{proof}
First work in a fixed finite-dimensional covariance eigenbasis. If $Z$ is standard normal, completing the square in its density gives, for $1-2\theta\lambda>0$,
\[
 \mathbb E\exp\{\theta[\lambda(Z^2-1)+2\sqrt\lambda\,bZ]\}
 =\frac{e^{-\theta\lambda}}{\sqrt{1-2\theta\lambda}}
       \exp\left(\frac{2\theta^2\lambda b^2}{1-2\theta\lambda}\right).
\]
Indeed, the coefficient of $Z^2$ in the exponent of the Gaussian integral is $-(1-2\theta\lambda)/2$; translating its center contributes the last exponential. For $|u|<1$,
\[
 -\log(1-u)-u=\sum_{r\ge2}\frac{u^r}{r}
 \le\sum_{r\ge2}\frac{|u|^r}{r}
 \le\frac{u^2}{2(1-|u|)}.
\]
Writing $G_i=\sqrt{\lambda_i}Z_i$ and $b_i=\langle b,e_i\rangle$, independence and the preceding scalar formulas therefore imply
\begin{equation}\label{S:eq:covariance-gaussian-mgf}
 \begin{split}
 \log\mathbb E e^{\theta Z_b}
 &=\sum_i\left[-\theta\lambda_i
       -\tfrac12\log(1-2\theta\lambda_i)
       +\frac{2\theta^2\lambda_i b_i^2}{1-2\theta\lambda_i}\right]\\
 &\le\frac{\theta^2(V+2v_b)}{1-2|\theta|\Lambda},
                  \qquad 2|\theta|\Lambda<1.
 \end{split}
\end{equation}
This proves directly the moment generating function estimate used in \cite{bechar2009bernstein}.

Put $A=V+2v_b$. When $A>0$, choose
\[
 \theta_x=\frac{\sqrt x}{\sqrt A+2\Lambda\sqrt x},
 \qquad t_x=2\sqrt{Ax}+2\Lambda x.
\]
One has $2\theta_x\Lambda<1$ and
\[
 -\theta_xt_x+\frac{\theta_x^2A}{1-2\theta_x\Lambda}=-x.
\]
Exponential Markov inequality applied to $Z_b$ and $-Z_b$ gives $\mathbb P(|Z_b|>t_x)\le2e^{-x}$. Taking $b=0$ proves the first line of
\eqref{S:eq:covariance-gaussian-explicit-tails}. The scalar variable $\langle G,b\rangle$ is normal with variance $v_b$, so its moment generating function is $e^{\theta^2v_b/2}$.
Optimizing at $\theta=\sqrt{2x/v_b}$ proves the second line when $v_b>0$; when $v_b=0$, that scalar is zero almost surely.

For a positive trace-class covariance on a separable Hilbert space, use an orthonormal eigenbasis on the closure of $\operatorname{Ran}Q$. The component of $G$ in $\ker Q$ is zero
almost surely. The projections $G^{[d]}$ satisfy
\[
 \mathbb E\|G-G^{[d]}\|^2=\sum_{i>d}\lambda_i\longrightarrow0,
 \qquad G^{[d]}\longrightarrow G\quad\text{almost surely}.
\]
The almost sure statement also follows directly from the convergent coordinate expansion of the Hilbert-valued random element $G$. The constants in the finite-dimensional estimates
are bounded by
\[
 V_d\le V\le\Lambda T,
 \qquad (v_b)_d\le v_b\le\Lambda\|b\|^2,
 \qquad \Lambda_d\le\Lambda.
\]
Define $Z_b^{[d]}=\|G^{[d]}\|^2-T_d+2\langle G^{[d]},b\rangle$. It converges almost surely to $Z_b$, including when $b$ has a component in $\ker Q$, since that component cancels
from the centered expression. Fatou's lemma for $e^{\theta Z_b^{[d]}}$ passes the uniform moment generating function bound to $Z_b$. The same choice of $\theta_x$ proves the
infinite-dimensional tails. If $A=0$, then $V=0$, hence $Q=0$ and $Z_b=0$ almost surely.

Under a regular conditional Gaussian law all these calculations apply with $Q$ and $b$ fixed. The quantities $T,\Lambda,V,v_b$ are measurable functions of the covariance and shift.
Thus the resulting inequalities are conditional inequalities with measurable thresholds. No unconditional moment of $T$ or $\|b\|$ is needed. When a Bochner conditional expectation
is used, restrict first to $\{T+\|b\|^2\le M\}$, and then let the integer $M$ increase.
\end{proof}

For a fixed threshold $t>0$, $V\le T\Lambda$ and the choice $x=\frac1{16}\min\{t^2/(T\Lambda),t/\Lambda\}$ give $2\sqrt{Vx}+2\Lambda x\le t/2+t/8<t$. Hence
\begin{equation}\label{S:eq:bridge-gaussian-concentration}
 \mathbb P\bigl(|\|G\|^2-T|\ge t\bigr)
 \le2\exp\left[-\frac1{16}
       \min\left\{\frac{t^2}{T\Lambda},\frac t\Lambda\right\}\right],
 \qquad t>0,
\end{equation}
when $T\Lambda>0$; for $Q=0$ the probability is zero. The moment identities are
\begin{equation}\label{S:eq:covariance-gaussian-moments}
 \operatorname{Var}(\|G\|^2)=2V,
 \qquad \operatorname{Var}(\langle G,b\rangle)=v_b,
 \qquad \operatorname{Cov}(\|G\|^2,\langle G,b\rangle)=0.
\end{equation}
For the moment identities, let $g$ be standard normal. Then symmetry gives $\mathbb E(g^3-g)=0$, and integration by parts gives $\mathbb Eg^4=3\mathbb Eg^2=3$. Therefore
$\mathbb E(g^2-1)^2=2$. Independence makes the finite centered quadratic sums orthogonal in $L^2$, with total variance $2\sum_i\lambda_i^2$. Their $L^2$ limit is $\|G\|^2-T$, by
the almost sure convergence above. The linear sums converge in $L^2$ with variance $\sum_i\lambda_i b_i^2$; every mixed covariance is zero. Passing to these two $L^2$ limits proves
\eqref{S:eq:covariance-gaussian-moments}, and in particular
\begin{equation}\label{S:eq:gaussian-noncentral-variance}
 \operatorname{Var}(\|G+b\|^2)=2V+4v_b.
\end{equation}
The quadratic and linear terms have zero covariance; they may still be dependent. Under a conditional law these are identities of conditional moments. On the original space,
restricting to $\{T+\|b\|^2\le L\}$ makes all the relevant second and fourth moments integrable, by these same formulas; letting $L$ increase specifies the identities even without
unconditional integrability.

\subparagraph{Finite arrays} For estimates of the finite arrays below and in the next subsubsection, probabilities are under the original law on the space of
Theorem~\ref{S:thm:joint-probability-space}, and $B=W$ is the original driving Brownian motion. We first estimate sums with deterministic weights under this law. A pathwise
extension will then allow random weights.

For all sufficiently large $n$, put
\[
 \mathcal J_n^{\mathrm{at}}:=\{\ell_n,\ldots,n-1\},
 \ell_n=m_n-R_\beta\ge m_n/2,
 h_n^*:=\max_{j\in\mathcal J_n^{\mathrm{at}}}h_{n,j}\le C_\beta/\ell_n.
\]
We disregard the finitely many $n$ for which this set is empty or $\ell_n<1$. For every integer $q$, define
\begin{equation}\label{S:eq:qv-pair-index-set}
 \mathcal I_n(q):=\{j\in\mathcal J_n^{\mathrm{at}}:j+q\in\mathcal J_n^{\mathrm{at}}\}.
\end{equation}
Thus $\mathcal I_n(0)=\mathcal J_n^{\mathrm{at}}$. The coordinates $s_{n,j},t_{n,j},h_{n,j},\varrho_{n,j}$ are the ones used in the preceding increment lemmas, and we abbreviate
$g_{n,j}=\Delta B_{n,j}$. For a fixed $n$, the vectors $g_{n,j}$ are independent and have laws $\mathcal N(0,h_{n,j}I_2)$. Different rows are constructed from the same Brownian
path. In this subsection, matrix norms are Hilbert--Schmidt norms unless marked otherwise, and $F:M=\operatorname{tr}(F^{\mathsf T}M)$.

We first treat a finite Gaussian array with deterministic coefficients.

Let $g_j$ be independent centered Gaussian vectors in $\mathbb R^2$ with covariance $h_jI_2$, $h_j>0$, over a finite index set, and let $F_j$ be deterministic matrices. Define
\[
 Z_0=\sum_j F_j:(g_jg_j^{\mathsf T}-h_jI_2),
 A_j=\tfrac12(F_j+F_j^{\mathsf T}),
\]
Set $v_0=\sum_jh_j^2\|A_j\|^2$ and $b_0=\max_jh_j\|A_j\|_{\mathrm{op}}$. For $q\ne0$ and any set $\mathcal J$ on which $j,j+q$ exist, put
\[
 Z_q=\sum_{j\in\mathcal J}g_j^{\mathsf T}F_jg_{j+q},
 v_q=\sum_{j\in\mathcal J}h_jh_{j+q}\|F_j\|^2,
 b_q=\max_{j\in\mathcal J}\sqrt{h_jh_{j+q}}\|F_j\|_{\mathrm{op}}.
\]

\begin{lemma}\label{S:lem:qv-gaussian-arrays}
For $u>0$, with an absolute $c>0$,
\begin{equation}\label{S:eq:qv-diagonal-bernstein}
 \mathbb P(|Z_0|>u)
 \le2\exp\left[-c\min\left(\frac{u^2}{v_0},\frac{u}{b_0}\right)\right].
\end{equation}
and
\begin{equation}\label{S:eq:qv-shifted-bernstein}
 \mathbb P(|Z_q|>u)
 \le4\exp\left[-c\min\left(\frac{u^2}{v_q},\frac{u}{b_q}\right)\right].
\end{equation}
When $v_q=0$ (including $q=0$), the corresponding sum is zero almost surely.
\end{lemma}

\begin{proof}
For a real symmetric matrix $H$ and a standard Gaussian vector $\xi$, orthogonal diagonalization gives
\[
 \xi^{\mathsf T}H\xi-\operatorname{tr}H
       =\sum_i\lambda_i(Z_i^2-1),
\]
where the $Z_i$ are independent standard normals and $\lambda_i$ may have either sign. Gaussian integration and the logarithmic series used in the preceding lemma give
\[
 \log\mathbb E e^{\theta(\xi^{\mathsf T}H\xi-\operatorname{tr}H)}
 \le\frac{\theta^2\operatorname{tr}(H^2)}
             {1-2|\theta|\|H\|_{\mathrm{op}}},
 \qquad 2|\theta|\|H\|_{\mathrm{op}}<1.
\]
Only $|\lambda_i|\le\|H\|_{\mathrm{op}}$ is used in this bound, so positive definiteness is unnecessary. The same optimization as before yields
\[
 \mathbb P\left(\left|\xi^{\mathsf T}H\xi-\operatorname{tr}H\right|
       >2\sqrt{\operatorname{tr}(H^2)x}+2\|H\|_{\mathrm{op}}x\right)
 \le2e^{-x}.
\]
This is the matrix estimate of \cite{bechar2009bernstein}.

Write $g_j=\sqrt{h_j}\xi_j$. Since the skew-symmetric part of $F_j$ pairs to zero with $g_jg_j^{\mathsf T}$ and has trace zero,
\[
 Z_0=\sum_jh_j(\xi_j^{\mathsf T}A_j\xi_j-\operatorname{tr}A_j).
\]
Concatenate the independent vectors $\xi_j$. The matrix in this quadratic form is block diagonal with blocks $h_jA_j$, and hence has squared Frobenius norm $v_0$ and operator norm
$b_0$. If $v_0>0$, take $x=\frac1{16}\min\{u^2/v_0,u/b_0\}$. The threshold is at most $u/2+u/8<u$, which proves \eqref{S:eq:qv-diagonal-bernstein}.

For $q\ne0$, set $d=|q|$ and write each index as $j=r+md$, with $0\le r<d$. Divide the pairs $(j,j+q)$ into two classes according to the parity of $m$. Within a fixed residue class
a pair joins the two successive integers $m,m+1$ when $q>0$, or $m,m-1$ when $q<0$. Pairs with the same parity of $m$ have no common endpoint. Distinct residues also have no common
endpoint. Thus all Gaussian vectors used by distinct pairs within either class are independent.

For one pair put $C_j=\sqrt{h_jh_{j+q}}F_j$. Its bilinear expression is the quadratic form of the concatenated vector $(\xi_j,\xi_{j+q})$ with symmetric matrix
\[
 H_j=\begin{pmatrix}0&C_j/2\\ C_j^{\mathsf T}/2&0\end{pmatrix},
 \quad \operatorname{tr}H_j=0,
 \quad \operatorname{tr}(H_j^2)=\tfrac12\|C_j\|^2,
 \quad \|H_j\|_{\mathrm{op}}=\tfrac12\|C_j\|_{\mathrm{op}}.
\]
The last equality follows by squaring $H_j$, which gives diagonal blocks $C_jC_j^{\mathsf T}/4$ and $C_j^{\mathsf T}C_j/4$. Each parity class is therefore a block diagonal
quadratic form with squared Frobenius norm at most $v_q/2$ and operator norm at most $b_q/2$.

If $v_q>0$, take $x=\frac1{64}\min\{u^2/v_q,u/b_q\}$. For either class the preceding tail threshold is at most
\[
 2\sqrt{(v_q/2)x}+b_qx
       \le\frac{u}{\sqrt{32}}+\frac u{64}<\frac u2.
\]
If $|Z_q|>u$, at least one class exceeds $u/2$ in absolute value. A union bound, which does not require independence between the two classes, gives
\eqref{S:eq:qv-shifted-bernstein}. When $v_0=0$ all $A_j$ are zero and $Z_0=0$; when $v_q=0$ for $q\ne0$, all contributing $F_j$ are zero and $Z_q=0$.
\end{proof}

For $T>0$, collect the intervals contained in $[0,T]$ in the index sets
\begin{equation}\label{S:eq:qv-compact-index-sets}
 \mathcal J_n(T):=\{j\in\mathcal J_n^{\mathrm{at}}:t_{n,j}\le T\},\qquad \mathcal J_n(T,q):=\{j\in\mathcal I_n(q):t_{n,j}\le T,\ t_{n,j+q}\le T\}.
\end{equation}
The second set equals the first when $q=0$. Each index in $\mathcal J_n(T,q)$ represents two sampling intervals contained in $[0,T]$.

\subsubsection{Weighted quadratic variation}
\label{S:sec:s9-weighted-variation}

The matrix $m_{n,k}(s)m_{n,k}(s)^{\mathsf T}$ contains products of the observed boundary point increments. Its kernel weights depend on the full path and terminal point, including
Brownian increments later than the interval being tested.

We construct an event on which the matrix-valued increment measures converge weakly against every continuous function. The bounds for total variation near zero then extend
convergence to weights satisfying \eqref{S:eq:qv-random-weight-hypotheses}.

Consider an array $F_{n,j}\in\mathbb R^{2\times2}$, $j\in\mathcal J_n^{\mathrm{at}}$, and a locally continuous $F_\omega:(0,1]\to\mathbb R^{2\times2}$. We require, pathwise for
every $0<\delta<1$,
\begin{equation}\label{S:eq:qv-random-weight-hypotheses}
 \max_{\substack{j\in\mathcal J_n^{\mathrm{at}}\\\varrho_{n,j}\ge\delta}}
      \|F_{n,j}-F_\omega(\varrho_{n,j})\|\longrightarrow0,
 \|F_{n,j}\|\le C_\omega\varrho_{n,j}^{\zeta}
\end{equation}
for some $\zeta>0$, with the second inequality holding for all sufficiently large $n$.

\begin{lemma}\label{S:lem:weighted-qv}
There is an event of probability one $\Omega_{\mathrm{qv}}$ on which every array satisfying \eqref{S:eq:qv-random-weight-hypotheses} obeys
\begin{equation}\label{S:eq:weighted-qv-compact}
 \sum_{j\in\mathcal J_n^{\mathrm{at}}}F_{n,j}:
       (\Delta_{n,j}\Delta_{n,j}^{\mathsf T}-h_{n,j}I_2)
       \longrightarrow0,
\end{equation}
and, for every fixed nonzero integer $q$,
\begin{equation}\label{S:eq:weighted-cross-variation}
 \sum_{j\in\mathcal I_n(q)}F_{n,j}:
       \Delta_{n,j}\Delta_{n,j+q}^{\mathsf T}\longrightarrow0.
\end{equation}
\end{lemma}
\begin{proof}
Fix an integer $T\ge1$ and an offset $q\in\mathbb Z$.

\phantomsection\label{S:lem:qv-compact-measures} There is an event $\Omega_{\mathrm{qv}}^{\mathrm{cmp}}$ of probability
one with the following properties, simultaneously for all integers $T\ge1$ and $q\in\mathbb Z$. Put $r_T=\me^{-\beta T/4}$ and define the finite matrix-valued measures on $[r_T,1]$
by
\begin{equation}\label{S:eq:qv-compact-matrix-measures}
 \mu_{n,T}^{(0)}:=\sum_{j\in\mathcal J_n(T)} (g_{n,j}g_{n,j}^{\mathsf T}-h_{n,j}I_2)\,\delta_{\varrho_{n,j}},\qquad \mu_{n,T}^{(q)}:=\sum_{j\in\mathcal J_n(T,q)} g_{n,j}g_{n,j+q}^{\mathsf T}\,\delta_{\varrho_{n,j}} , (q\ne0).
\end{equation}
We prove that, on one event of probability one, these measures converge weakly to zero and their total variations are bounded by a finite $C_\omega$ for all $n\ge n_0(\omega)$. The
event is common to all continuous matrix-valued functions.

On $\Omega_{\mathrm{qv}}^{\mathrm{cmp}}\cap\Omega_{\mathrm{incr}}$, replacing $\mathcal J_n(T,q)$ by $\{j\in\mathcal I_n(q):s_{n,j}<T\}$ changes a sum with uniformly bounded
coefficients by $O_{\beta,T,q,\omega}((1+\log n)/n)$. This estimate holds for both $g_{n,j}$ and the normal increments $\Delta_{n,j}$, with centering $h_{n,j}I_2$ when $q=0$.

If $s_{n,j}\le T+1$, then $\varrho_{n,j}\ge \me^{-\beta(T+1)/4}$. The explicit estimate $|n\varrho_{n,j}-(j+1)|\le C_\beta$ in Lemma~\ref{S:lem:sampling-position} gives
$j\ge c_{\beta,T}n$ and hence $h_{n,j}\le C_{\beta,T}/n$ for all large $n$. In particular
\begin{equation}\label{S:eq:qv-compact-mesh-bound}
 \max_{j\in\mathcal J_n(T)}h_{n,j}\le C_{\beta,T}/n,
   \sum_{j\in\mathcal J_n(T)}h_{n,j}\le T.
\end{equation}
The last time $t_{n,\ell_n}$ covered by the retained intervals tends to infinity, since Lemma~\ref{S:lem:4.1} gives $t_{n,\ell_n}=(4/\beta)\log(n/\ell_n)+O_\beta(\ell_n^{-1})$ and
$\ell_n/n\to0$. Consequently this range covers $[0,T]$ except for at most one interval crossing $T$.

For a deterministic function $f$ with $M=\sup_{[r_T,1]}\|f\|>0$, apply Lemma~\ref{S:lem:qv-gaussian-arrays} to $F_j=f(\varrho_{n,j})$. For the diagonal sum its parameters obey
$v_0\le C_{\beta,T}M^2/n$ and $b_0\le C_{\beta,T}M/n$. For a shifted sum on $\mathcal J_n(T,q)$, both durations are at most $C_{\beta,T}/n$, and
\[
 \sum_{j\in\mathcal J_n(T,q)}h_{n,j}h_{n,j+q}
 \le \left(\max_{\ell\in\mathcal J_n(T)}h_{n,\ell}\right)
                 \sum_{j\in\mathcal J_n(T)}h_{n,j}
 \le C_{\beta,T}/n.
\]
Thus for every $u>0$, including $q=0$ with the first definition in \eqref{S:eq:qv-compact-matrix-measures},
\begin{equation}\label{S:eq:qv-compact-test-tail}
 \mathbb P\left(\left|\int f:\dd\mu_{n,T}^{(q)}\right|>u\right)
 \le4\exp\left[-c_{\beta,T}n
          \min\left(\frac{u^2}{M^2},\frac{u}{M}\right)\right].
\end{equation}
The zero function contributes zero. For each fixed function and positive threshold, the probability is summable in $n$. Borel--Cantelli gives the almost sure limit.

For each $T$, take the countable family of piecewise linear matrix functions whose vertices have rational entries and whose knots, in the coordinate $(r-r_T)/(1-r_T)$, lie on a
finite uniform rational partition. This family is uniformly dense: linear interpolation on a partition with intervals of length $1/N$ has error bounded by the modulus of continuity
at $(1-r_T)/N$, and the finitely many vertex values can then be approximated by rational matrices. Include the constant matrix $I_2$. Apply Borel--Cantelli to these functions, all
positive rational thresholds, and all integers $T,q$. This constructs one event on which every function in these countable families has limiting integral zero.

The choice of the constant function $I_2$ in the diagonal family gives
\[
 \sum_{j\in\mathcal J_n(T)}|g_{n,j}|^2
 =2\sum_{j\in\mathcal J_n(T)}h_{n,j}+o(1)\le2T+o(1).
\]
For the diagonal measure, the sum of the norms of its point masses is at most $\sum_j|g_{n,j}|^2+2T$, and is therefore bounded. For $q\ne0$, Cauchy--Schwarz gives
\[
 \|\mu_{n,T}^{(q)}\|_{\mathrm{TV}}
 \le\left(\sum_{j\in\mathcal J_n(T,q)}|g_{n,j}|^2\right)^{1/2}
      \left(\sum_{j\in\mathcal J_n(T,q)}|g_{n,j+q}|^2\right)^{1/2}
 \le\sum_{j\in\mathcal J_n(T)}|g_{n,j}|^2.
\]
The shifted indices in the second sum are distinct members of $\mathcal J_n(T)$. Density and these bounds extend convergence to every continuous function: for any approximating
member $f_m$,
\[
 \left|\int f:\dd\mu_{n,T}^{(q)}\right|
 \le\left|\int f_m:\dd\mu_{n,T}^{(q)}\right|
       +\|f-f_m\|_\infty\|\mu_{n,T}^{(q)}\|_{\mathrm{TV}}.
\]
Taking $n\to\infty$ and then $m\to\infty$ proves weak convergence for every continuous function on the chosen event. This extension is made after fixing the path. In particular, if
$f_\omega$ is a random continuous function and $F_{n,j}$ converges to $f_\omega(\varrho_{n,j})$ uniformly on this compact interval, the additional error is bounded by
\[
 \sup_j\|F_{n,j}-f_\omega(\varrho_{n,j})\|
                  \|\mu_{n,T}^{(q)}\|_{\mathrm{TV}}\longrightarrow0.
\]
The probability estimate was applied only to the countable family of deterministic functions. The last inequality is evaluated separately for each continuous path $f_\omega$ and
the sampled increments, including when they are dependent.

The sampling intervals are consecutive and disjoint in their interiors. At most one has $s_{n,j}<T<t_{n,j}$. If the interval for $j$ is contained in $[0,T]$ but its partner $j+q$
is not, $j$ is within $|q|$ indices of that crossing. Thus the difference of the two sets in the statement has at most $|q|+1$ indices. Their starts and the partner starts are
within $(|q|+1)h_n^*$ of $T$, and hence are at most $T+1$ for large $n$. All their durations are $O_{\beta,T,q}(n^{-1})$ by the first paragraph. On $\Omega_{\mathrm{incr}}$,
Corollary~\ref{S:cor:normal-increment-full-array} gives
\[
 |g_{n,j}|+|\Delta_{n,j}|
 \le C_\beta\bigl(\sqrt{h_{n,j}\log n}+h_{n,j}\log n\bigr)
 \le C_{\beta,T,q}\sqrt{\log n/n}
\]
for these indices, for all sufficiently large $n$. Each product or centered square has norm $O_{\beta,T,q}((1+\log n)/n)$. Multiplying by a bounded coefficient and summing the
bounded number of indices proves the estimate for changing the set of retained intervals. If $T$ is itself a sampling time, no interval has $T$ in its interior; the same count and
bound remain valid.

\phantomsection\label{S:lem:qv-weighted-masses} There is an event $\Omega_{\mathrm{qv}}^{\mathrm{tail}}$ of probability
one such that the following statements hold for every positive rational $c$ and every integer $q$. Define, for real
$T\ge0$,
\[
 H_n(c,T):=\sum_{\substack{j\in\mathcal J_n^{\mathrm{at}}\\s_{n,j}\ge T}}
                 \me^{-cs_{n,j}}h_{n,j},
 S_n(c,T):=\sum_{\substack{j\in\mathcal J_n^{\mathrm{at}}\\s_{n,j}\ge T}}
                 \me^{-cs_{n,j}}|g_{n,j}|^2.
\]
There are a deterministic $C_{\beta,c}<\infty$ and an almost surely finite $N_{\beta,c}$ such that, for every $n\ge N_{\beta,c}$ and every real $T\ge0$,
\begin{equation}\label{S:eq:qv-tail-uniform}
 H_n(c,T)+S_n(c,T)\le C_{\beta,c}\me^{-cT}.
\end{equation}
The same random starting index works for every $T$.

On $\Omega_{\mathrm{qv}}^{\mathrm{tail}}\cap\Omega_{\mathrm{incr}}$, put
\[
 \epsilon_n^{\mathrm{qv}}:=(\log n)\sqrt{h_n^*}+(\log n)^2h_n^*
 =O_\beta((\log n)^{-255})\longrightarrow0.
\]
For each fixed $c,q$ and all sufficiently large $n$, uniformly in $T\ge0$,
\begin{equation}\label{S:eq:qv-normal-weighted-replacement}
 \begin{split}
 &\sum_{\substack{j\in\mathcal J_n^{\mathrm{at}}\\s_{n,j}\ge T}}\me^{-cs_{n,j}}
   \|\Delta_{n,j}\Delta_{n,j}^{\mathsf T}-g_{n,j}g_{n,j}^{\mathsf T}\|
       \le C_{\beta,c}\epsilon_n^{\mathrm{qv}}\me^{-cT},\\
 &\sum_{\substack{j\in\mathcal I_n(q)\\s_{n,j}\ge T}}\me^{-cs_{n,j}}
   \|\Delta_{n,j}\Delta_{n,j+q}^{\mathsf T}
                    -g_{n,j}g_{n,j+q}^{\mathsf T}\|
       \le C_{\beta,c,q}\epsilon_n^{\mathrm{qv}}\me^{-cT}
       , (q\ne0).
 \end{split}
\end{equation}
The weighted absolute masses of the normal centered squares and shifted products satisfy
\begin{equation}\label{S:eq:qv-normal-weighted-total-mass}
 \begin{split}
 &\sum_{\substack{j\in\mathcal J_n^{\mathrm{at}}\\s_{n,j}\ge T}}\me^{-cs_{n,j}}
   \|\Delta_{n,j}\Delta_{n,j}^{\mathsf T}-h_{n,j}I_2\|
       \le C_{\beta,c}\me^{-cT},\\
 &\sum_{\substack{j\in\mathcal I_n(q)\\s_{n,j}\ge T}}\me^{-cs_{n,j}}
   |\Delta_{n,j}|\,|\Delta_{n,j+q}|
       \le C_{\beta,c,q}\me^{-cT}, (q\ne0).
 \end{split}
\end{equation}
For the driver increments the two inequalities are
\[
 \sum_{\substack{j\in\mathcal J_n^{\mathrm{at}}\\s_{n,j}\ge T}}
 e^{-cs_{n,j}}\|g_{n,j}g_{n,j}^{\mathsf T}-h_{n,j}I_2\|
 \le S_n(c,T)+2H_n(c,T),
\]
\[
 \sum_{\substack{j\in\mathcal I_n(q)\\s_{n,j}\ge T}}
e^{-cs_{n,j}}|g_{n,j}|\,|g_{n,j+q}|
 \le e^{cd_{n,q}/2}\sqrt{S_n(c,T)S_n(c,T_q)}.
\]
Here $d_{n,q}=|q|h_n^*$ and $T_q=\max(0,T-d_{n,q})$. The first follows from the triangle inequality and $\|gg^{\mathsf T}\|=|g|^2$; the second is proved in
\eqref{S:eq:qv-shifted-mass-bound} below. The uniform estimate \eqref{S:eq:qv-tail-uniform} bounds each right side by $C_{\beta,c,q}e^{-cT}$ when $d_{n,q}\le1$.

Take $n$ so large that $h_n^*\le1$. For an interval $[s,s+h]$,
\[
 \me^{-cs}h\le \me^{ch_n^*}\int_s^{s+h}\me^{-cu}\,\dd u.
\]
Disjointness gives $H_n(c,T)\le(\me^c/c)\me^{-cT}$. With
\[
 V_n(c,T):=\sum_{s_{n,j}\ge T}\me^{-2cs_{n,j}}h_{n,j}^2,
   b_n(c,T):=\max_{s_{n,j}\ge T}\me^{-cs_{n,j}}h_{n,j},
\]
we have $V_n(c,T)\le C_c h_n^*\me^{-2cT}$ and $b_n(c,T)\le h_n^*\me^{-cT}$. Apply \eqref{S:eq:qv-diagonal-bernstein} to $F_j=\me^{-cs_{n,j}}I_2$ on this deterministic tail set and
take threshold $\me^{-cT}$. The mean of $S_n(c,T)$ is $2H_n(c,T)$, so
\begin{equation}\label{S:eq:qv-tail-row-probability}
 \mathbb P\left(S_n(c,T)>(2\me^c/c+1)\me^{-cT}\right)
 \le2\me^{-d_c/h_n^*}\le2\me^{-d_{\beta,c}\ell_n}.
\end{equation}
For $T>s_{n,\ell_n}$ the sums are zero. The number of integers $T\ge0$ with a nonempty sum is at most $C_\beta(1+\log n)$, by the bound on the last sampled time. The union over
these integers has probability at most
\[
 C_\beta(1+\log n)\me^{-d_{\beta,c}\ell_n}.
\]
This is summable: $\ell_n\ge\tfrac12(\log n)^{512}$ and the displayed expression is at most $n^{-3}$ for all sufficiently large $n$. Borel--Cantelli therefore proves the estimate
simultaneously for every integer $T$ after one random starting index. For real $T$, monotonicity gives $S_n(c,T)\le S_n(c,\lfloor T\rfloor)\le C_c\me^c\me^{-cT}$. The deterministic
bound for $H_n$ already holds for real $T$. Intersect over positive rational $c$ to obtain \eqref{S:eq:qv-tail-uniform}.

For a fixed $q$ and $j\in\mathcal I_n(q)$, the intervening $|q|$ sampling intervals give
\begin{equation}\label{S:eq:qv-shift-time-distance}
 |s_{n,j+q}-s_{n,j}|\le |q|h_n^*=:d_{n,q}\longrightarrow0.
\end{equation}
Let $T_q=\max(0,T-d_{n,q})$. Cauchy--Schwarz, the injectivity of $j\mapsto j+q$, and comparison of the two exponential weights imply
\begin{equation}\label{S:eq:qv-shifted-mass-bound}
 \sum_{\substack{j\in\mathcal I_n(q)\\s_{n,j}\ge T}}
     \me^{-cs_{n,j}}|g_{n,j}|\,|g_{n,j+q}|
 \le \me^{cd_{n,q}/2}\bigl(S_n(c,T)S_n(c,T_q)\bigr)^{1/2}
 \le C_{\beta,c,q}\me^{-cT}
\end{equation}
for $d_{n,q}\le1$. For $T_q=\max(0,T-d_{n,q})$ one has $0\le T-T_q\le d_{n,q}$. Thus the exponential factor after substituting \eqref{S:eq:qv-tail-uniform} is bounded by
\[
 \me^{cd_{n,q}/2}\me^{-c(T+T_q)/2}
 =\me^{-cT}\me^{c(d_{n,q}+T-T_q)/2}
 \le\me^c\me^{-cT},
\]
whenever $d_{n,q}\le1$. This includes $T<d_{n,q}$, for which $T_q=0$.

Write $R_{n,j}=\Delta_{n,j}-g_{n,j}$. On $\Omega_{\mathrm{incr}}$, its bound in \eqref{S:eq:qv-normal-comparison} gives
\begin{equation}\label{S:eq:qv-weighted-error-square}
 D_n(c,T):=\sum_{s_{n,j}\ge T}\me^{-cs_{n,j}}|R_{n,j}|^2
 \le C_\beta(\log n)^2h_n^*H_n(c,T)
 \le C_{\beta,c}(\log n)^2h_n^*\me^{-cT}.
\end{equation}
For the diagonal tensors, use
\[
 \Delta_j\Delta_j^{\mathsf T}-g_jg_j^{\mathsf T}
 =R_jg_j^{\mathsf T}+g_jR_j^{\mathsf T}+R_jR_j^{\mathsf T}.
\]
The weighted sum of its norms is at most $2\sqrt{S_n(c,T)D_n(c,T)}+D_n(c,T)$, which is the first bound in \eqref{S:eq:qv-normal-weighted-replacement}.

For a shifted product the exact difference is
\[
 R_jg_{j+q}^{\mathsf T}+g_jR_{j+q}^{\mathsf T}
                         +R_jR_{j+q}^{\mathsf T}.
\]
The three weighted absolute sums are at most, respectively,
\[
 \me^{cd_{n,q}/2}\sqrt{D_n(c,T)S_n(c,T_q)},
 \me^{cd_{n,q}/2}\sqrt{S_n(c,T)D_n(c,T_q)},
 \me^{cd_{n,q}/2}\sqrt{D_n(c,T)D_n(c,T_q)}.
\]
Equations \eqref{S:eq:qv-tail-uniform} and \eqref{S:eq:qv-weighted-error-square} give the second bound in \eqref{S:eq:qv-normal-weighted-replacement}. Since
$h_n^*\le C_\beta(\log n)^{-512}$, the stated rate for $\epsilon_n^{\mathrm{qv}}$ follows. Finally the mass of the centered square of the driver is bounded by $S_n(c,T)+2H_n(c,T)$,
and the driver cross mass is bounded by \eqref{S:eq:qv-shifted-mass-bound}. Add the corresponding normal replacement bounds to prove \eqref{S:eq:qv-normal-weighted-total-mass}. The
bounds are for the absolute weighted sums.

Take
\[
 \Omega_{\mathrm{qv}}
 =\Omega_{\mathrm{incr}}\cap\Omega_{\mathrm{qv}}^{\mathrm{cmp}}
                           \cap\Omega_{\mathrm{qv}}^{\mathrm{tail}}.
\]
Only the countable families of functions, integer cutoffs, integer offsets, and positive rational exponential weights used in the preceding compact and tail estimates enter this
event. Fix a path on it and an array satisfying \eqref{S:eq:qv-random-weight-hypotheses}. Choose a rational $c$ with $0<c<\beta\zeta/4$. The assumed power bound implies
\begin{equation}\label{S:eq:qv-random-weight-exponential}
 \|F_{n,j}\|\le C_\omega \me^{-cs_{n,j}}.
\end{equation}
Consequently \eqref{S:eq:qv-normal-weighted-replacement}, with $T=0$, shows that replacing all normal products by their driver products changes either entire weighted sum by at
most $C_{\beta,c,q,\omega}\epsilon_n^{\mathrm{qv}}\to0$.

Fix an integer $T\ge1$. Uniform convergence of the weights on $[r_T,1]$, bounded total variation, and the \hyperref[S:lem:qv-compact-measures]{measure convergence on compact
intervals} give
\[
 \begin{split}
 \left|\sum_{j\in\mathcal J_n(T,q)}
       F_{n,j}:\mu_{n,T}^{(q)}(\{\varrho_{n,j}\})\right|
 &\le \left|\int_{[r_T,1]}F_\omega:\dd\mu_{n,T}^{(q)}\right|\\
 &\quad+\max_{\varrho_{n,j}\ge r_T}
       \|F_{n,j}-F_\omega(\varrho_{n,j})\|
       \|\mu_{n,T}^{(q)}\|_{\mathrm{TV}}\longrightarrow0.
 \end{split}
\]
The first term tends to zero for the fixed continuous function $F_\omega$ on the chosen path. This pathwise estimate applies to the weights and increments with their joint
distribution. Use $\mathcal I_n(0)=\mathcal J_n^{\mathrm{at}}$ for the diagonal sum. The earlier cutoff estimate shows that replacing $\mathcal J_n(T,q)$ by
$\{j\in\mathcal I_n(q):s_{n,j}<T\}$ changes this sum by at most $C_{\beta,T,q,\omega}(1+\log n)/n$: at most $C_q$ intervals cross the cutoff, their weights are bounded, and each
measure atom has norm at most $C_{\beta,T,q,\omega}(1+\log n)/n$. Thus the entire compact part tends to zero.

For the remaining intervals with $s_{n,j}\ge T$, combine \eqref{S:eq:qv-random-weight-exponential} with the driver version of \eqref{S:eq:qv-normal-weighted-total-mass}. The limsup
of the absolute tail is at most $C_{\beta,c,q,\omega}\me^{-cT}$. The full sum is the compact part plus this tail and the controlled normal replacement error. For either sum $Z_n$
in the statement, the compact part tends to zero and the replacement error is $o(1)$. Hence $\limsup_{n\to\infty}|Z_n|\le C_{\beta,c,q,\omega}\me^{-cT}$ for every integer $T\ge1$.
Letting $T\to\infty$ gives
$Z_n\to0$.

For every $0<\zeta'<\zeta$, choose rational $\beta\zeta'/4<c<\beta\zeta/4$. Equation \eqref{S:eq:qv-normal-weighted-total-mass}, with $T=(4/\beta)\log(1/\delta)$, bounds the
absolute weighted tail with $\varrho_{n,j}\le\delta$ by $C_{\omega,\zeta',q}\delta^{\zeta'}$, uniformly for all large $n$. This supplies uniform integrability at the singular
boundary, including any fixed number of additional intervals with $j=\ell_n+O(1)$ or
$j=m_n+O(1)$.

The estimates hold pathwise on $\Omega_{\mathrm{qv}}$ for every array satisfying the hypotheses. They give convergence for the random weights, without a rate unless further
regularity is imposed.

The estimate includes the earlier formulation with $\|F_{n,j}\|\le C_{\omega,\varepsilon}\varrho_{n,j}^{a-\varepsilon}$ for all sufficiently small $\varepsilon>0$: fix any one such
$\varepsilon<a$ and use $\zeta=a-\varepsilon$. It also includes continuous compactly supported functions. Local convergence expressed on $j/n\ge\delta$ is equivalent after
shrinking $\delta$, because $|\varrho_{n,j}-j/n|\le C_\beta/n$.

The event depends only on the increment arrays and applies to every weight array satisfying
\eqref{S:eq:qv-random-weight-hypotheses}, for all fixed offsets. \end{proof}

\subsubsection{The local variance}
\label{S:sec:s9-local-variance}

For a fraction $\lambda$ of a sampling interval of length $h$, the square of the linear interpolation term contributes $h\lambda^2$ after replacing the boundary point increment
square by its quadratic variation. The bridge variance supplies the missing part:
\[
 h\lambda^2+h(|\lambda|-\lambda^2)=h|\lambda|.
\]
The contribution from linear interpolation enters this covariance identity after taking the weighted quadratic variation limit. When several intervals are traversed, the square of
the interpolation term also contains products of distinct increments. Their weighted sums are handled by Lemma~\ref{S:lem:weighted-qv}.

Theorem~\ref{S:thm:interval-fractions} gives the limits of $\lambda_{n,k,j}$ and their support $|j-k|\le R_\beta$. Keep $m_n\le k<n$ for partition indices and
$j\in\mathcal J_n^{\mathrm{at}}$ for sampling indices. The following calculation verifies \eqref{S:eq:qv-random-weight-hypotheses} for the resulting weights.

To compute the local variance, write
\[
 \Sigma_{n,k}(s):=\mathbb E[\widehat W_{n,k}(s)
                 \widehat W_{n,k}(s)^{\mathsf T}\mid\mathscr S_n],
   \ell_{n,k}(s)=|\theta_{n,k}(s)-\vartheta_{n,k+1}|.
\]

\phantomsection\label{S:lem:qv-induced-weights} Suppose that $F_{n,k}(s)$ are measurable matrix fields on the partition
subintervals, that $F_\omega$ is continuous on $(0,1]\times[0,1]$, and that, pathwise, the following bounds hold for some $\zeta>0$:
\begin{equation}\label{S:eq:qv-mesh-weight-hypotheses}
 \begin{split}
 e_n(\delta)&:=\sup_{\substack{m_n\le k<n,\ k/n\ge\delta\\0\le s\le1}}
       \|F_{n,k}(s)-F_\omega(k/n,s)\|\longrightarrow0
       , (0<\delta<1),\\
 \sup_s\|F_{n,k}(s)\|&\le C_\omega(k/n)^\zeta.
 \end{split}
\end{equation}
For an integer $q$, define on $\mathcal J_n^{\mathrm{at}}$
\begin{equation}\label{S:eq:qv-induced-weight-definition}
 G_{n,j}^{(q)}:=\sum_{k=m_n}^{n-1}\int_0^1
       \lambda_{n,k,j}(s)\lambda_{n,k,j+q}(s)F_{n,k}(s)\,\dd s,
\end{equation}
with fractions whose interval index lies outside $\mathcal J_n^{\mathrm{at}}$ set to zero. The support theorem shows that this convention changes no nonzero fraction from a
selected subinterval. These weights satisfy \eqref{S:eq:qv-random-weight-hypotheses}, with locally continuous limit
\begin{equation}\label{S:eq:qv-induced-weight-explicit-limit}
 G_\omega^{(q)}(r)=\sum_{d\in\mathbb Z}\int_0^1
   \Lambda_{\beta,d}(r,s)\Lambda_{\beta,d+q}(r,s)F_\omega(r,s)\,\dd s,
 \Lambda_{\beta,d}(r,s)=\operatorname{clip}_{[0,1]}
       (d+1-s-b_\beta(1-r))-\mathbf1_{\{d\ge1\}}.
\end{equation}
Both the weights and their limits vanish if $|q|>2R_\beta$. For a fixed $\delta>0$, let $M_\delta$ be the supremum of $\|F_\omega\|$ on $[\delta/2,1]\times[0,1]$, and let
$\omega_\delta(u)$ be its modulus of continuity in the first coordinate on that set. For all large $n$,
\begin{equation}\label{S:eq:qv-induced-weight-modulus}
 \max_{j/n\ge\delta}
 \|G_{n,j}^{(q)}-G_\omega^{(q)}(\varrho_{n,j})\|
 \le C_{\beta,q,\delta}
       \left[e_n(\delta/2)+\frac{M_\delta}{n}
                       +\omega_\delta(C_{\beta,q}/n)\right].
\end{equation}
Moreover $\|G_{n,j}^{(q)}\|\le C_{\beta,q,\zeta}C_\omega(j/n)^\zeta$ throughout the enlarged index range.

Only $k$ with $|k-j|\le R_\beta$ can occur, and both fractions can be nonzero only if $|q|\le2R_\beta$. There are at most $2R_\beta+1$ terms. Since $j\ge m_n-R_\beta\to\infty$ and
$k\ge m_n$, their radii $j/n$ and $k/n$ are uniformly comparable. This gives the power bound, including the additional indices below $m_n$.

For convergence on a compact interval write $k=j-d$. For indices in the retained range, the estimate \eqref{S:eq:interval-fractions-error} bounds the error in each fraction by
$C_{\beta,q,\delta}/n$. The product error is bounded by the sum of the two errors, because all factors have magnitude at most one. Replacing $F_{n,j-d}(s)$ by $F_\omega((j-d)/n,s)$
gives an error at most $e_n(\delta/2)$. Replacing its first argument by $j/n$ gives an error at most $\omega_\delta(R_\beta/n)$. The limiting functions are uniformly Lipschitz in
$r$, so changing their first argument gives an error at most $C_\beta M_\delta/n$ after multiplication by the field.

For large $n$ and $j/n\ge\delta$, a missing partition index can only be $j-d\ge n$. Then $d<0$ and $1-j/n\le R_\beta/n$. The identities $\Lambda_{\beta,0}(1,s)=1-s$ and
$\Lambda_{\beta,d}(1,s)=0$ for $d\ne0$ show that the missing factor is $O_\beta(n^{-1})$ uniformly in $s$. If $j+q\ge n$, then $q\ne0$, $1-j/n\le |q|/n$, and every product
$\Lambda_{\beta,d}(1,s)\Lambda_{\beta,d+q}(1,s)$ vanishes. Its omitted contribution is again $O_{\beta,q}(M_\delta/n)$. The case $j+q<\ell_n$ cannot occur on this compact range for
large $n$. These observations verify the last partition subintervals and both signs of the shift.

Summing the finite number of terms proves the bound at $j/n$. Since $|\varrho_{n,j}-j/n|\le C_\beta/n$, evaluation at $\varrho_{n,j}$ introduces an additional error at most
$C_{\beta,q}(M_\delta/n+\omega_\delta(C_{\beta,q}/n))$. This proves \eqref{S:eq:qv-induced-weight-modulus}. The finite sum of continuous integrands in
\eqref{S:eq:qv-induced-weight-explicit-limit} is continuous, and Lemma~\ref{S:lem:sampling-position} converts the power bound to one in $\varrho_{n,j}$. This is the quantitative
form of the induced weight limit in \eqref{S:eq:integrated-fractions-limit} needed here.

\begin{proposition}\label{S:prop:bridge-and-interpolation-second-moment}

On $\Omega_{\mathrm{qv}}$, every array satisfying \eqref{S:eq:qv-mesh-weight-hypotheses} obeys
\begin{equation}\label{S:eq:qv-second-moment-error}
 \mathcal V_n(F):=\sum_{k=m_n}^{n-1}\int_0^1 F_{n,k}(s):
    \left[m_{n,k}(s)m_{n,k}(s)^{\mathsf T}
               +\Sigma_{n,k}(s)-\ell_{n,k}(s)I_2\right]\dd s
       \longrightarrow0.
\end{equation}
The covariance in this identity is conditional on $\mathscr S_n$.

\end{proposition}
\begin{proof}
For two operator positions $(k,s)$ and $(l,t)$, put $u_j=u_{n,j}(\theta_{n,k}(s))$ and $\widetilde u_j=u_{n,j}(\theta_{n,l}(t))$. The coordinates of a Euclidean bridge have
covariance $K(u,v)=u\wedge v-uv$ and mean zero. Each sampled argument $u_{n,j}(\vartheta_{n,k+1})$ is zero or one, where the bridge vanishes. Conditional independence over $j$
therefore gives the full covariance identity for the two positions
\[
 \mathbb E[\widehat W_{n,k}^{\mu}(s)\widehat W_{n,l}^{\nu}(t)
                      \mid\mathscr S_n]
 =\mathbf1_{\{\mu=\nu\}}\sum_{j\in\mathcal J_n^{\mathrm{at}}}
            h_{n,j}K(u_j,\widetilde u_j).
\]
A common Brownian observation interval can contribute to two distinct operator cells, so we do not assume independence of those cell fields. The linear interpolation term
$m_{n,k}(s)$ is $\mathscr S_n$-measurable. Hence the conditional mean of $m_{n,k}(s)+\widehat W_{n,k}(s)$ is $m_{n,k}(s)$, and its conditional second moment at a single position is
\[
 \mathbb E[(m_{n,k}(s)+\widehat W_{n,k}(s))
            (m_{n,k}(s)+\widehat W_{n,k}(s))^{\mathsf T}
                \mid\mathscr S_n]
 =m_{n,k}(s)m_{n,k}(s)^{\mathsf T}+\Sigma_{n,k}(s).
\]
To evaluate $\Sigma_{n,k}$, put $v_j=u_{n,j}(\vartheta_{n,k+1})$. By \eqref{S:eq:relative-time-at-sample}, $v_j\in\{0,1\}$ and $\lambda_{n,k,j}(s)=u_j-v_j$. The Euclidean bridge
vanishes at $v_j$. Its covariance at $u_j$ is
\[
 u_j(1-u_j)I_2=
       (|\lambda_{n,k,j}(s)|-\lambda_{n,k,j}(s)^2)I_2,
\]
whether $v_j=0$ or $v_j=1$. Conditional independence of the bridges on different intervals, proved in Theorem~\ref{S:thm:joint-probability-space}, therefore gives
\begin{equation}\label{S:eq:qv-exact-bridge-covariance}
 \Sigma_{n,k}(s)=\sum_{j\in\mathcal J_n^{\mathrm{at}}}h_{n,j}
       (|\lambda_{n,k,j}(s)|-\lambda_{n,k,j}(s)^2)I_2.
\end{equation}
Completely traversed sampling intervals have $|\lambda|=1$ and covariance zero. The signed fractions occur in the linear interpolation term
$m_{n,k}(s)=\sum_{j\in\mathcal J_n^{\mathrm{at}}}\lambda_{n,k,j}(s)\Delta_{n,j}$. The exact length partition \eqref{S:eq:interval-length-partition} gives
$\ell_{n,k}(s)=\sum_{j\in\mathcal J_n^{\mathrm{at}}}h_{n,j}|\lambda_{n,k,j}(s)|$ once $m_n\ge k_\beta$. Subtracting this identity and expanding the square of the linear
interpolation term yields
\begin{equation}\label{S:eq:qv-interval-second-moment}
 \begin{split}
 &m_{n,k}(s)m_{n,k}(s)^{\mathsf T}+\Sigma_{n,k}(s)-\ell_{n,k}(s)I_2\\
 &=\sum_{j\in\mathcal J_n^{\mathrm{at}}}\lambda_{n,k,j}(s)^2
             (\Delta_{n,j}\Delta_{n,j}^{\mathsf T}-h_{n,j}I_2)\\
 & +\sum_{\substack{i,j\in\mathcal J_n^{\mathrm{at}}\\i\ne j}}
         \lambda_{n,k,i}(s)\lambda_{n,k,j}(s)
                                   \Delta_{n,i}\Delta_{n,j}^{\mathsf T}.
 \end{split}
\end{equation}
All sums are finite; the coefficient assumptions make all integrals in $s$ absolutely finite for this fixed row. Thus integration and collection by $q=j-i$ are legitimate and give
\begin{equation}\label{S:eq:qv-sum-over-increments}
 \mathcal V_n(F) =\sum_{j\in\mathcal J_n^{\mathrm{at}}}G_{n,j}^{(0)}: (\Delta_{n,j}\Delta_{n,j}^{\mathsf T}-h_{n,j}I_2) +\sum_{\substack{q\in\mathbb Z\\0<|q|\le2R_\beta}} \sum_{j\in\mathcal I_n(q)}G_{n,j}^{(q)}: \Delta_{n,j}\Delta_{n,j+q}^{\mathsf T}.
 \end{equation}
The cross sum runs over both positive and negative offsets, with each ordered pair counted once. The identity therefore holds also for nonsymmetric matrices. Indices of the
sampling intervals begin at $\ell_n$; indices of the partition intervals begin at
$m_n$.

The \hyperref[S:lem:qv-induced-weights]{weight estimates} verify the hypotheses of Lemma~\ref{S:lem:weighted-qv} for $G^{(0)}$ and each of the finitely many nonzero offsets. Their
sums tend to zero on $\Omega_{\mathrm{qv}}$, proving \eqref{S:eq:qv-second-moment-error}. The index set includes the additional intervals below $m_n$, whose contributions satisfy
the same boundary point bound.

With $\tau_{n,k},\kappa_{n,k}$ as in Proposition~\ref{S:prop:sampling-offset},
\begin{equation}\label{S:eq:integrated-time-difference}
 \left|\int_0^1\ell_{n,k}(s)\,\dd s-\tau_{n,k}\kappa_{n,k}\right|
       \le\frac{1}{3\beta k^2},
 \tau_{n,k}=\frac4\beta\log(1+1/k),
 |\kappa_{n,k}-\kappa_\beta(k/n)|\le C_\beta/k.
\end{equation}
More generally, for fields depending on $s$, the weighted integral has the limit
\begin{equation}\label{S:eq:qv-weighted-second-moment-limit}
 \begin{split}
 &\sum_{k=m_n}^{n-1}\int_0^1F_{n,k}(s):
       [m_{n,k}(s)m_{n,k}(s)^{\mathsf T}+\Sigma_{n,k}(s)]\,\dd s\\
 &\longrightarrow
 \frac4\beta\int_0^1\!\int_0^1
       \frac{|1-s-b_\beta(1-r)|}{r}\,
                         \operatorname{tr}F_\omega(r,s)\,\dd s\,\dd r.
 \end{split}
\end{equation}
The integral is absolutely finite. These limits concern, for each fixed path, the covariance already conditioned on $\mathscr S_n$. If the weights are $\mathscr S_n$-measurable,
the left side of \eqref{S:eq:qv-weighted-second-moment-limit} is the conditional weighted second moment of the sum of the linear interpolation term and the centered bridge. For
other random weights, the formula gives only the displayed pathwise pairing, not a conditional second moment. After conditioning on $\mathscr S_n$, we retain the conditional means
and cross covariances of the driver increments.

The first two estimates in \eqref{S:eq:integrated-time-difference} are precisely the deterministic estimate \eqref{S:eq:time-interpolation-error}; the sampling error is
\eqref{S:eq:sampling-factor-approximation}. For the limit involving the field, use the following pointwise consequence of the same estimates:
\begin{equation}\label{S:eq:qv-time-difference-approximation}
 \sup_{0\le s\le1}\left|
 \ell_{n,k}(s)-\frac4{\beta k}|1-s-b_\beta(1-k/n)|\right|
 \le C_\beta k^{-2}.
\end{equation}
Linear interpolation in \eqref{S:eq:time-interpolation-error} gives $\theta_{n,k}(s)-\vartheta_{n,k+1} =\tau_{n,k}(1-s-\widehat\sigma_{n,k})+O_\beta(k^{-2})$ uniformly in $s$.
Substitute $\tau_{n,k}=4/(\beta k)+O_\beta(k^{-2})$ and $\widehat\sigma_{n,k}=b_\beta(1-k/n)+O_\beta(k^{-1})$, then use the reverse triangle inequality. Multiplication by the
weight and integration show that the error is at most
\begin{equation}\label{S:eq:qv-time-error-sum}
 C_\beta\sum_{k=m_n}^{n-1}k^{-2}\int_0^1\|F_{n,k}(s)\|\,\dd s
 \le\frac{C_{\beta,\zeta,\omega}}{m_n}
       \sum_{k=1}^{n-1}\frac{(k/n)^\zeta}{k}
 \le\frac{C_{\beta,\zeta,\omega}}{m_n}\longrightarrow0.
\end{equation}
For $0<\zeta\le1$, monotonicity gives $\sum_{k=1}^{n-1}k^{\zeta-1}\le\int_0^{n-1}x^{\zeta-1}\dd x\le n^\zeta/\zeta$. For $\zeta\ge1$, each summand is at most $n^{\zeta-1}$, so the
sum is at most $n^\zeta$.

On each compact rectangle $[\delta,1]\times[0,1]$, uniform convergence of $F_{n,k}$ and continuity of $F_\omega$ imply convergence of the resulting Riemann sum to the double
integral in \eqref{S:eq:qv-weighted-second-moment-limit} restricted to $r\ge\delta$. The limit inherits the bound $\|F_\omega(r,s)\|\le C_\omega r^\zeta$: approximate a fixed $r>0$
by selected points in a compact interval and pass to the limit uniformly in $s$. Since $|1-s-b_\beta(1-r)|$ is bounded, the integral over $0<r<\delta$ is at most
$C_\omega\int_0^\delta r^{\zeta-1}\dd r=C_\omega\delta^\zeta/\zeta$ in absolute value. The same comparison of sums of powers with integrals bounds the discrete tail by
$C_{\zeta,\omega}\delta^\zeta$, uniformly in $n$. Letting $\delta\downarrow0$ proves convergence of the whole deterministic Riemann sum and absolute finiteness of its limit.
Combine this with \eqref{S:eq:qv-second-moment-error} and \eqref{S:eq:qv-time-error-sum} to prove \eqref{S:eq:qv-weighted-second-moment-limit}.

\end{proof}

We apply the preceding result to the following two random weights.

Use the sampled columns and the sums $A_{n,k}$ and $C_{n,k}$ from \eqref{S:eq:exclusive-primitives-step-identities}, and write $\mathfrak G_{n,k}=A_{n,k}p_{n,k}+C_{n,k}q_{n,k}$ as
before. Define
\begin{equation}\label{S:eq:qv-matrix-H-definition}
 \mathbf H_{n,k}^{\mu\nu}
 :=\frac12\left[A_{n,k}\langle\mathsf A_\mu a_{n,k},\mathsf A_\nu a_{n,k}\rangle
       +C_{n,k}\langle\mathsf A_\mu c_{n,k},\mathsf A_\nu c_{n,k}\rangle\right],
   \Gamma_{n,k}:=\operatorname{tr}\mathbf H_{n,k}.
\end{equation}

\begin{corollary}\label{S:cor:qv-sampled-green-weights}
The fields $\mathfrak G_{n,k}I_2$ and $\mathbf H_{n,k}$ satisfy \eqref{S:eq:qv-mesh-weight-hypotheses} for every fixed $0<\zeta<a_\beta$. Almost surely,
\begin{equation}\label{S:eq:qv-scalar-weighted-limit}
 \sum_{k=m_n}^{n-1}\mathfrak G_{n,k}\int_0^1
        (|m_{n,k}(s)|^2+\operatorname{tr}\Sigma_{n,k}(s))\,\dd s
 =\frac8\beta\sum_{k=m_n}^{n-1}\frac{\kappa_\beta(k/n)\mathfrak G_{n,k}}{k}
       +o_\omega(1)=O_\omega(1),
\end{equation}
and
\begin{equation}\label{S:eq:qv-matrix-weighted-limit}
 \sum_{k=m_n}^{n-1}\int_0^1\mathbf H_{n,k}:
        [m_{n,k}(s)m_{n,k}(s)^{\mathsf T}+\Sigma_{n,k}(s)]\,\dd s
 =\frac4\beta\sum_{k=m_n}^{n-1}\frac{\kappa_\beta(k/n)\Gamma_{n,k}}{k}
       +o_\omega(1).
\end{equation}
\end{corollary}

\begin{proof}
For $v,w\in\mathbb R^2$ and $\mu,\nu\in\{1,2\}$,
\[
 \bigl|\langle\mathsf A_\mu v,\mathsf A_\nu v\rangle
       -\langle\mathsf A_\mu w,\mathsf A_\nu w\rangle\bigr|
 \le\|\mathsf A_\mu\|\|\mathsf A_\nu\|(|v|+|w|)|v-w|.
\]
Lemma~\ref{S:lem:sampled-exclusive-convergence} bounds all columns and primitives on $r_k\ge\delta>0$ and makes their differences from the continuous values tend uniformly to zero.
Apply the displayed inequality and $|xy-x'y'|\le|x-x'||y|+|x'||y-y'|$ to each entry of \eqref{S:eq:qv-matrix-H-definition}. This gives the continuous compact limits required in
\eqref{S:eq:qv-mesh-weight-hypotheses}, also at $r=1$, where $C_\omega(1)=0$.

For $z\in\mathbb R^2$, the definition of $\mathbf H_{n,k}$ gives
\[
 z^{\mathsf T}\mathbf H_{n,k}z
 =\tfrac12A_{n,k}\left|\sum_\mu z_\mu\mathsf A_\mu a_{n,k}\right|^2
  +\tfrac12C_{n,k}\left|\sum_\mu z_\mu\mathsf A_\mu c_{n,k}\right|^2\ge0.
\]
Since $\sum_\mu|\mathsf A_\mu v|^2=(v_1^2+5v_2^2)/4$,
\[
 \|\mathbf H_{n,k}\|\le\Gamma_{n,k}\le\tfrac58\mathfrak G_{n,k},\qquad
 \|\mathfrak G_{n,k}I_2\|\le\sqrt2\mathfrak G_{n,k}
 \le C_{\omega,\zeta}(k/n)^\zeta,
 \quad 0<\zeta<a_\beta,
\]
where the last inequality is Proposition~\ref{S:prop:5.4}. Thus both fields satisfy all the hypotheses of Proposition~\ref{S:prop:bridge-and-interpolation-second-moment}. For
either $F_{n,k}=\mathfrak G_{n,k}I_2$ or $F_{n,k}=\mathbf H_{n,k}$, it gives
\[
 \sum_{k=m_n}^{n-1}\int_0^1 F_{n,k}:
 \bigl[m_{n,k}(s)m_{n,k}(s)^{\mathsf T}+\Sigma_{n,k}(s)-\ell_{n,k}(s)I_2\bigr]\,\dd s
 =o_\omega(1).
\]
The deterministic time estimate \eqref{S:eq:sampling-variance-two-errors} and Proposition~\ref{S:prop:5.5} yield
\[
 \begin{aligned}
 &\left|\sum_{k=m_n}^{n-1}\operatorname{tr}F_{n,k}
   \left(\int_0^1\ell_{n,k}(s)\,\dd s-\frac{4\kappa_\beta(k/n)}{\beta k}\right)\right|\\
 &\qquad\le C_\beta\sum_{k=m_n}^{n-1}\frac{\mathfrak G_{n,k}}{k^2}
 \le\frac{C_\beta}{m_n}\sum_{k=m_n}^{n-1}\frac{\mathfrak G_{n,k}}k
 \le\frac{C_{\beta,\omega}}{m_n}\longrightarrow0.
 \end{aligned}
\]
For the scalar field, $F:(mm^{\mathsf T}+\Sigma)=\mathfrak G(|m|^2+\operatorname{tr}\Sigma)$ and $\operatorname{tr}F=2\mathfrak G$; for the matrix field,
$\operatorname{tr}F=\Gamma$. Substitution gives \eqref{S:eq:qv-scalar-weighted-limit} and \eqref{S:eq:qv-matrix-weighted-limit}. Finally,
\[
 0\le\sum_{k=m_n}^{n-1}\frac{\kappa_\beta(k/n)\mathfrak G_{n,k}}k
 \le\|\kappa_\beta\|_\infty\sup_n\sum_{k=m_n}^{n-1}\frac{\mathfrak G_{n,k}}k<\infty,
\]
which proves the stated boundedness.
\end{proof}

\subsubsection{The Gaussian field in the weighted Hilbert space}
\label{S:sec:s9-weighted-gaussian-field}

We realize the bridge field as a Gaussian random element in the domain of the kernel derivative. Its sine expansion converges in $\mathfrak H_n^+$ and yields a finite conditional
second moment and a measurable covariance. The same norm, including its diagonal contribution, also controls the nonlinear remainder.

We apply Section~\ref{S:sec:stability-and-fr-echet-expansion-of-the-inverse-kernel} with $m=m_n$ and the sampled columns on the partition intervals $a_{n,k},c_{n,k}$. For
measurable fields $U_{n,k}(s)$, the three quadratic forms are defined in \eqref{S:eq:continuous-field-energies}; in particular,
\[
 \mathcal E_n^{+}(U)
 =\sum_{k=m_n}^{n-1}\int_0^1\mathcal M_{n,k}^{+}(U_{n,k}(s))\dd s,
 \mathcal M_{n,k}^{+}(V)
 =A_{n,k}^{+}|Va_{n,k}|^2+C_{n,k}^{+}|Vc_{n,k}|^2.
\]
The exact trace decomposition uses the off-diagonal quadratic form; the derivative bound uses the full quadratic form. Write $\mathscr D_n=\mathscr D_{n,X^{\mathrm{step}}}$ for the
first derivative defined by the kernel \eqref{S:eq:continuous-derivative-kernel}. Theorem~\ref{S:thm:continuous-frechet} proves \eqref{S:eq:continuous-frechet-bound} and
\eqref{S:eq:continuous-frechet-remainder} with absolute constants. Proposition~\ref{S:prop:continuous-refinement} justifies arbitrary finite subdivisions with their lengths in $r$,
and the \hyperref[S:prop:continuous-measurability]{measurability argument for the kernel maps} proves that these kernels are strongly measurable in $\mathfrak S_2$.

Conditional on $\mathscr S_n$, the bridge fields in distinct sampling intervals are independent and Gaussian. Since
\[
\mathcal E_n(\mathcal B_n)+
\mathbb E[\mathcal E_n(\mathcal Z_n)\mid\mathscr S_n]
\le C\sum_{k=m_n}^{n-1}\mathfrak G_{n,k}\int_0^1
\left(|m_{n,k}(s)|^2+
\mathbb E[|\widehat W_{n,k}(s)|^2\mid\mathscr S_n]\right)\dd s,
\]
Corollary~\ref{S:cor:qv-sampled-green-weights} gives compact convergence and the required power bound near the singular boundary for $\mathfrak G_{n,k}$. Its scalar identity
\eqref{S:eq:qv-scalar-weighted-limit} therefore gives
\[
\mathcal E_n(\mathcal B_n)+
\mathbb E[\mathcal E_n(\mathcal Z_n)\mid\mathscr S_n]
\le C_{\beta,\omega}
\sum_{k=m_n}^{n-1}\frac{\mathfrak G_{n,k}}k+o_\omega(1).
\]
The missing diagonal part is smaller. On the event of probability one in Corollary~\ref{S:cor:bridge-array-bound}, at most $O_\beta(1)$ indices $j$ contribute to a fixed
subinterval, and hence
\[
\sup_{0\le s\le1}\bigl(|m_{n,k}(s)|+|\widehat W_{n,k}(s)|\bigr)
\le C_{\beta,\omega}\sqrt{\frac{\log n}{k}}.
\]
For $U=\sum_\mu\mathsf A_\mu y_\mu$, the bounds on $\mathsf A_1,\mathsf A_2$ give
\[
|c_{n,k}|^2|Ua_{n,k}|^2+|a_{n,k}|^2|Uc_{n,k}|^2
\le C|a_{n,k}|^2|c_{n,k}|^2|y|^2.
\]
The logarithmic product bound \eqref{S:eq:s7-column-product-bound}, which follows from the Brownian estimate \eqref{S:eq:s7-continuous-bounds}, gives
$p_{n,k}q_{n,k}\le C_{\beta,\omega}(\log n)^{5/2}$. Hence

\begin{equation}\label{S:eq:bridge-diagonal-square-energy}
\mathcal E_n^{\mathrm{diag}}(\mathcal B_n)
+\mathcal E_n^{\mathrm{diag}}(\mathcal Z_n)
\le \frac{C_{\beta,\omega}(\log n)^{7/2}}n
\sum_{k=m_n}^{n-1}\frac1k
\le \frac{C_{\beta,\omega}(\log n)^{9/2}}n
\longrightarrow0.
\end{equation}
The displayed bound shows that the diagonal contribution tends to zero for $\mathcal B_n$ and $\mathcal Z_n$. We use the full quadratic form as the squared norm of
$\mathfrak H_n^+$, as in Lemma~\ref{S:lem:weighted-hilbert}. Its reference columns, matrices $D_{n,k}$, and isometry $\mathcal J_n$ are $\mathscr S_n$-measurable, with partition
indices $m_n\le k<n$. Conditional on $\mathscr S_n$, this is a fixed separable Hilbert space. The form $\mathcal E_n$ remains a separate term in the trace identity.

The diagonal estimate \eqref{S:eq:bridge-diagonal-square-energy} was a pathwise estimate; to control a conditional trace we also need its conditional mean version. The bounds on
$\mathsf A_1,\mathsf A_2$ and the Brownian bridge covariance give
\[
 \mathbb E[|\widehat W_{n,k}(s)|^2\mid\mathscr S_n]\le C_\beta k^{-1},
 \int_0^1|m_{n,k}(s)|^2\dd s\le C_{\beta,\omega}\frac{\log n}{k}.
\]
Only a bounded number of indices $j$ contribute to each subinterval, so
\[
 \mathcal E_n^{\mathrm{diag}}(\mathcal B_n)
 +\mathbb E[\mathcal E_n^{\mathrm{diag}}(\mathcal Z_n)\mid\mathscr S_n]
 \le\frac{C_{\beta,\omega}(1+\log n)}n
       \sum_{k=m_n}^{n-1}\frac{p_{n,k}q_{n,k}}k
\]
The conditional covariance bound is deterministic; the sampled columns and the linear interpolation term are fixed under the conditional integral. The constant depending on the
path is used only after that integral is evaluated. Applying the logarithmic product bound to the resulting sum gives

\begin{equation}\label{S:eq:weighted-bridge-diagonal-mean}
 \mathcal E_n^{\mathrm{diag}}(\mathcal B_n)
 +\mathbb E[\mathcal E_n^{\mathrm{diag}}(\mathcal Z_n)\mid\mathscr S_n]
 \le \frac{C_{\beta,\omega}(\log n)^{7/2}}n
          \sum_{k=m_n}^{n-1}\frac1k
 \le\frac{C_{\beta,\omega}(\log n)^{9/2}}n\longrightarrow0.
\end{equation}
Adding this to the preceding bound with $\mathcal E_n$ in place of $\mathcal E_n^+$ gives
\begin{equation}\label{S:eq:weighted-bridge-mean-bound}
 \mathcal E_n^{+}(\mathcal B_n)
 +\mathbb E[\mathcal E_n^{+}(\mathcal Z_n)\mid\mathscr S_n]
 \le C_{\beta,\omega}.
\end{equation}

Take $n$ large enough that $m_n<n$ and $\ell_n=m_n-R_\beta\ge1$. Sampling indices range over $\mathcal J_n^{\mathrm{at}}=\{\ell_n,\ldots,n-1\}$, with zero coefficients for indices
that do not contribute; partition indices range over $\{m_n,\ldots,n-1\}$. For $m\ge1$, define the scalar function and the $\mathscr S_n$-measurable field
\begin{equation}\label{S:eq:gaussian-mode-profiles}
 \begin{aligned}
 q_{n,k,j,m}(s)
 &:=\frac{\sqrt2}{\pi m}
   \left[\sin\!\bigl(\pi m u_{n,j}(\theta_{n,k}(s))\bigr)
       -\sin\!\bigl(\pi m u_{n,j}(\vartheta_{n,k+1})\bigr)\right],\\
 \mathcal V_{n,j,m,\nu;k}(s)
 &:=\sqrt{h_{n,j}}q_{n,k,j,m}(s)\mathsf A_\nu,
   \nu\in\{1,2\}.
 \end{aligned}
\end{equation}
The second sine is zero, because the clipped coordinate of a sampling time is zero or one. The coefficient vanishes unless $|k-j|\le R_\beta$. Use the normal variables
$\xi_{n,j,m,\nu}$ extracted from the attached bridges in \eqref{S:eq:joint-kl-coordinates}. Put
\begin{equation}\label{S:eq:gaussian-field-truncation}
 \mathcal Z_n^{(M)}:=
 \sum_{j\in\mathcal J_n^{\mathrm{at}}}\sum_{m=1}^M\sum_{\nu=1}^2
 \xi_{n,j,m,\nu}\mathcal V_{n,j,m,\nu},
   M\in\mathbb N.
\end{equation}
The order of summation is increasing $m$, with the finitely many $(j,\nu)$ in a fixed deterministic order at each mode. Set
\begin{equation}\label{S:eq:gaussian-sine-coefficient-bound}
 e_{n,j}:=h_{n,j}
       \sum_{\substack{m_n\le k<n\\|k-j|\le R_\beta}}
             \mathfrak G_{n,k}^{+},
 e_n^{\mathrm{KL}}:=\sum_{j\in\mathcal J_n^{\mathrm{at}}}e_{n,j}.
\end{equation}
These variables are $\mathscr S_n$-measurable and finite for each retained finite matrix list. They bound $\|\mathcal
V_{n,j,m,\nu}\|_{n,+}^2$.

With $\otimes_{+}$ denoting the rank-one operator for the weighted inner product, define
\[
 \widetilde Q_n^{+,(M)}
 :=\sum_{j\in\mathcal J_n^{\mathrm{at}}}\sum_{m=1}^M\sum_{\nu=1}^2
       \mathcal V_{n,j,m,\nu}\otimes_{+}\mathcal V_{n,j,m,\nu}.
\]

\begin{proposition}\label{S:prop:weighted-field-covariance}

Conditionally on $\mathscr S_n$, $\mathcal Z_n^{(M)}$ converges to $\mathcal Z_n$ in $L^2(\mathfrak H_n^+)$ and almost
surely. The limit is centered Gaussian, and its covariance is the increasing limit in trace norm of $\widetilde
Q_n^{+,(M)}$.

\end{proposition}
\begin{proof}
The constant $C$ is deterministic. Set
\[
 \mathfrak t_n^{\mathcal E}:=\operatorname{tr}\widetilde Q_n^{+}
   =\mathbb E[\mathcal E_n^{+}(\mathcal Z_n)\mid\mathscr S_n],
 \mathfrak l_n^{\mathcal E}:=\|\widetilde Q_n^{+}\|_{\mathrm{op}}.
\]
Fix the sampled Brownian path and use its regular conditional law. The normal coordinates $\xi_{n,j,m,\nu}$ are independent by Theorem~\ref{S:thm:joint-probability-space}. Since
the second sine in \eqref{S:eq:gaussian-mode-profiles} is zero,
\[
 |q_{n,k,j,m}(s)|\le\frac{\sqrt2}{\pi m},
 \mathcal M_{n,k}^{+}(\mathsf A_\nu)
 \le\|\mathsf A_\nu\|_{\mathrm{op}}^2\mathfrak G_{n,k}^{+}.
\]
Integration over $s$, with the normalized subinterval coordinate used in the full quadratic form, gives
\begin{equation}\label{S:eq:weighted-single-mode-bound}
 \|\mathcal V_{n,j,m,\nu}\|_{n,+}^2
 \le\frac{C e_{n,j}}{m^2}.
\end{equation}
Consequently,
\begin{equation}\label{S:eq:weighted-field-interval-trace}
 \sum_{m\ge1}\sum_{\nu=1}^2
       \|\mathcal V_{n,j,m,\nu}\|_{n,+}^2
 \le C e_{n,j},
 \sum_{m>M}\sum_{\nu=1}^2
       \|\mathcal V_{n,j,m,\nu}\|_{n,+}^2
 \le C e_{n,j}/M,
\end{equation}
where $\sum_{m>M}m^{-2}\le\int_M^\infty x^{-2}\dd x=M^{-1}$. For each fixed $n$ these sums are finite without any uniform bound for the boundary solutions.

Apply Lemma~\ref{S:lem:conditional-gaussian-series} to the vectors $\mathcal J_n\mathcal V_{n,j,m,\nu}$ in the fixed space
$\mathcal H_{0,n}=\bigoplus_{k=m_n}^{n-1}L^2((0,1);\mathbb R^{2\times2})$. The image of the weighted field space is closed, so the Gaussian series limit lies in that image, and the
inverse isometry gives convergence of $\mathcal Z_n^{(M)}$ in the specified conditional norms and almost surely.

To identify the series limit with the bridge field, we compute its conditional mean squared error. For a scalar bridge $\beta$ and its coefficients in
\eqref{S:eq:joint-kl-coordinates}, direct integration of the bridge covariance gives
\[
 \mathbb E[\beta(u)\xi_m]=\frac{\sqrt2}{\pi m}\sin(\pi m u),
 \mathbb E[\xi_m\xi_l]=\mathbf1_{\{m=l\}}.
\]
The sine expansion of the covariance then gives, for every deterministic $u\in[0,1]$,
\begin{equation}\label{S:eq:scalar-bridge-mode-tail}
 \mathbb E\left|\beta(u)-\sqrt2\sum_{m=1}^M
                   \frac{\xi_m}{\pi m}\sin(\pi m u)\right|^2
 =\frac2{\pi^2}\sum_{m>M}\frac{\sin^2(\pi m u)}{m^2}
 \le\frac2{\pi^2M}.
\end{equation}
Indeed, the integral calculation in Theorem~\ref{S:thm:joint-probability-space} identifies the sine coefficients of $v\mapsto K(u,v)$. Completeness of the sine basis gives equality
in $L^2(0,1)$. The series
\[
 \frac2{\pi^2}\sum_{m\ge1}
    \frac{\sin(\pi m u)\sin(\pi m v)}{m^2}
\]
converges absolutely and uniformly on $[0,1]^2$, since every term is bounded by $2/(\pi^2m^2)$. Its sum is continuous. Equality almost everywhere in $v$ therefore implies equality
everywhere, for each $u$. This justifies evaluation at the relative times $u_{n,j}(\theta_{n,k}(s))$, including zero and one.

After conditioning, every argument $u_{n,j}(\theta_{n,k}(s))$ is deterministic. Apply \eqref{S:eq:scalar-bridge-mode-tail} at that argument and at $u_{n,j}(\vartheta_{n,k+1})$,
where both the bridge and every sine function vanish. Independence between different $(j,\nu)$ eliminates mixed products of approximation errors. Tonelli, first for the nonnegative
squared weighted norm and then for the mode sum, gives exactly
\[
 \mathbb E[\mathcal E_n^{+}(\mathcal Z_n-\mathcal Z_n^{(M)})
                     \mid\mathscr S_n]
 =\sum_{j,m>M,\nu}\|\mathcal V_{n,j,m,\nu}\|_{n,+}^2.
\]
In particular $\mathbb E[\mathcal E_n^+(\mathcal Z_n)\mid\mathscr S_n]<\infty$, and the original field equals the series limit. Lemma~\ref{S:lem:conditional-gaussian-series}
identifies the limit as Gaussian and proves
\begin{equation}\label{S:eq:weighted-field-mode-tail}
 \mathbb E[\mathcal E_n^{+}(\mathcal Z_n-\mathcal Z_n^{(M)}) \mid\mathscr S_n] =\|\widetilde Q_n^{+}-\widetilde Q_n^{+,(M)}\|_{\mathfrak S_1(\mathfrak H_n^{+})} =\sum_{j\in\mathcal J_n^{\mathrm{at}}}\sum_{m>M}\sum_{\nu=1}^2 \|\mathcal V_{n,j,m,\nu}\|_{n,+}^2 \le\frac{C e_n^{\mathrm{KL}}}{M}.
\end{equation}
  After applying $\mathcal J_n$, the finite covariances are sums of $(\mathcal J_n\mathcal V)\otimes(\mathcal J_n\mathcal V)$. Their convergence in trace norm proves strong
  measurability of the limit by Lemma~\ref{S:lem:conditional-gaussian-series}. A countable intersection over $n$ gives these statements simultaneously for all rows.

For the uniform bounds, finite interaction range and $j\asymp_\beta k$ whenever the corresponding scalar function is nonzero give
$e_n^{\mathrm{KL}}\le C_\beta\sum_k\mathfrak G_{n,k}^+/k$. Including the current interval in both weights gives the identity
$\mathfrak G_{n,k}^+=\mathfrak G_{n,k}+2p_{n,k}q_{n,k}/n$. Consequently the earlier product estimate yields
\[
 e_n^{\mathrm{KL}}
 \le C_\beta\sum_{k=m_n}^{n-1}\frac{\mathfrak G_{n,k}}k
       +C_{\beta,\omega,\sigma}n^{-1+\sigma}m_n^{-\sigma}
 \le C_{\beta,\omega,\sigma},\qquad 0<\sigma<1,
\]
by Proposition~\ref{S:prop:5.5}. Using instead the sharper logarithmic bound for $p_{n,k}q_{n,k}$ gives the estimate
\begin{equation}\label{S:eq:weighted-sine-coefficient-sum}
 e_n^{\mathrm{KL}} \le C_\beta\sum_{k=m_n}^{n-1}\frac{\mathfrak G_{n,k}^{+}}{k} \le C_\beta\sum_k\frac{\mathfrak G_{n,k}}k +\frac{C_{\beta,\omega}(\log n)^{5/2}}n\sum_k\frac1k \le C_{\beta,\omega}.
\end{equation}
This uses Proposition~\ref{S:prop:5.5} and the bound for $|a|^2|c|^2$. For the bias and trace, use the estimate for the conditional mean \eqref{S:eq:weighted-bridge-mean-bound},
derived with the sampled columns fixed.

For $j\in\mathcal J_n^{\mathrm{at}}$, define
\begin{equation}\label{S:eq:s9-bridge-supports}
 E_{n,j}:=\{r\in(m_n/n,1):
       \vartheta_{n,j+1}<\theta_\beta(r)<\vartheta_{n,j}\}.
\end{equation}
These sets are pairwise disjoint. Since a sine mode vanishes at zero and one, \eqref{S:eq:gaussian-mode-profiles} gives $\mathcal V_{n,j,m,\nu}(r)=0$ for almost every
$r\notin E_{n,j}$. Even within a single operator interval, the nonzero mode fields from different sampling intervals occupy disjoint portions.

Let $T_j:\ell^2(\mathbb N\times\{1,2\})\to\mathfrak H_n^+$ send coefficients $(v_{m,\nu})$ to their sum of mode fields. For finite sequences and each $r$, Cauchy--Schwarz and
\[
 \frac2{\pi^2}\sum_{m\ge1}\frac{\sin^2(\pi m u)}{m^2}
       =u(1-u)\le\frac14
\]
give
\[
 \left\|\sum_{m,\nu}v_{m,\nu}
       \mathcal V_{n,j,m,\nu}(r)\right\|_{\mathrm{op}}^2
 \le C h_{n,j}\sum_{m,\nu}|v_{m,\nu}|^2.
\]
Integration with the full weights proves $\|T_j\|^2\le Ce_{n,j}$ and extends $T_j$ to $\ell^2$. The weighted inner product is an integral of pointwise products. Thus the
disjointness in \eqref{S:eq:s9-bridge-supports} proves
\begin{equation}\label{S:eq:s9-orthogonal-interval-modes}
 T_j^*T_l=0\quad(j\ne l),\qquad
 \left\|\sum_jT_jv_j\right\|_{n,+}^2
       =\sum_j\|T_jv_j\|_{n,+}^2.
\end{equation}
Define $T(v_j)_j=\sum_jT_jv_j$ on the finite Hilbert direct sum of the coefficient spaces. Orthogonality gives
\[
 \|T\|^2
 =\sup_{\sum_j\|v_j\|^2=1}\sum_j\|T_jv_j\|_{n,+}^2
 =\max_j\|T_j\|^2.
\]
For $U,V\in\mathfrak H_n^+$, the covariance identity is
\[
 \langle\widetilde Q_n^+U,V\rangle_{n,+}
 =\sum_{j,m,\nu}\langle\mathcal V_{n,j,m,\nu},U\rangle_{n,+}
                    \langle\mathcal V_{n,j,m,\nu},V\rangle_{n,+}
 =\langle T^*U,T^*V\rangle.
\]
Thus $\widetilde Q_n^+=TT^*$, and

\begin{equation}\label{S:eq:s9-field-covariance-direct-sum}
 \|\widetilde Q_n^+\|_{\mathrm{op}}
 =\|T\|^2\le C\max_j e_{n,j}
 \le C_\beta\max_{m_n\le k<n}\frac{\mathfrak G_{n,k}^+}{k}.
\end{equation}
For the last inequality, at most $2R_\beta+1$ terms occur in $e_{n,j}$, with $h_{n,j}\le C_\beta/k$ whenever $|k-j|\le R_\beta$. The orthogonality is in the weighted space of
fields. The kernel derivative sends these disjoint field supports to intersecting strips, so its images $\mathscr D_nT_jv_j$ can have nonzero Hilbert--Schmidt inner products. For
any bounded linear map $D$ on these comparison spaces,
\[
 \|D(\widetilde Q_n^{+,(M)}-\widetilde Q_n^+)D^*\|_{\mathfrak S_1}
 \le\|D\|^2\,
       \|\widetilde Q_n^{+,(M)}-\widetilde Q_n^+\|_{\mathfrak S_1}
 \le \frac{C\|D\|^2e_n^{\mathrm{KL}}}{M}\longrightarrow0.
\]
Taking $D=I$ gives convergence of the field covariances; taking $D=\mathscr D_n$ gives convergence of the kernel covariances. The trace changes by at most this trace norm and each
bilinear matrix element by at most this norm times the norms of the two vectors. These inequalities justify all of the finite covariance identities used above in the limit. Choose
$0<\eta<\min(a_\beta,\varepsilon)$. The off-diagonal part of this maximum is at most
$C_{\beta,\omega,\eta}n^{-a_\beta+\eta}\max_{k\ge m_n}k^{a_\beta-1-\eta} \le\allowbreak C_{\beta,\omega,\eta}n^{-a_\beta+\eta}$.

For the diagonal part choose $0<\sigma<\min(1,\varepsilon)$. The product estimate in Section~\ref{S:sec:s4-sampled-columns} gives
\[
 \max_{m_n\le k<n}\frac{2p_{n,k}q_{n,k}}{nk}
 \le C_{\beta,\omega,\sigma}n^{-1+\sigma}m_n^{-1-\sigma}
 \le C_{\beta,\omega,\sigma}n^{-a_\beta+\varepsilon},
\]
since $a_\beta\le1$ and $m_n\ge1$. This proves \eqref{S:eq:weighted-field-covariance-bound} without using the sharper logarithmic column maximum.

For the remaining identities fix a row and its regular conditional law given $\mathscr S_n$. Write $D=\mathscr D_n$, $b=\mathcal B_n$, $Q_M=\widetilde Q_n^{+,(M)}$, and
$Q=\widetilde Q_n^+$. Under this law the reference columns, $D$, $b$, and every mode field are deterministic. The derivative estimate gives
\[
 \|D\|_{\mathfrak H_n^+\to\mathfrak S_2}\le n^{-1/2}.
\]
For finite mode sums and any $F,G$ in the real Hilbert space of self-adjoint Hilbert--Schmidt operators, independence and centering of the Gaussian coefficients give
\[
 \begin{split}
 \mathbb E[D\mathcal Z_n^{(M)}\mid\mathscr S_n]&=0,\\
 \mathbb E[\langle D\mathcal Z_n^{(M)},F\rangle
           \langle D\mathcal Z_n^{(M)},G\rangle\mid\mathscr S_n]
       &=\langle DQ_MD^*F,G\rangle,\\
 \mathbb E[\|D\mathcal Z_n^{(M)}\|_{\mathfrak S_2}^2\mid\mathscr S_n]
       &=\sum_j\sum_{m=1}^{M}\sum_{\nu=1}^{2}\|D\mathcal V_{n,j,m,\nu}\|_{\mathfrak S_2}^2
        =\operatorname{tr}(DQ_MD^*).
 \end{split}
\]
The convergence of the full field in conditional $L^2$ implies
\[
 \mathbb E[\|D(\mathcal Z_n-\mathcal Z_n^{(M)})\|_{\mathfrak S_2}^2
                 \mid\mathscr S_n]
 \le\frac1n\operatorname{tr}(Q-Q_M)\longrightarrow0.
\]
Meanwhile the ideal inequality gives
\[
 \|DQD^*-DQ_MD^*\|_{\mathfrak S_1}
 \le\|D\|^2\|Q-Q_M\|_{\mathfrak S_1}\longrightarrow0.
\]
Conditional Cauchy--Schwarz passes each bilinear covariance to the limit. For example, the difference of the two products is the sum of a tail paired with one full variable and a
finite variable paired with a tail; its conditional absolute expectation is bounded by the product of their conditional $L^2$ norms. Continuity of the trace passes the last
identity to the limit. Thus the conditional mean and covariance of $D(b+\mathcal Z_n)$ are $Db$ and $DQD^*$, and
\[
 \begin{split}
 \mathbb E[\|D(b+\mathcal Z_n)\|_{\mathfrak S_2}^2\mid\mathscr S_n]
   &=\|Db\|_{\mathfrak S_2}^2+\operatorname{tr}(DQD^*),\\
 \operatorname{tr}(DQD^*)
   &=\sum_{j,m,\nu}\|D\mathcal V_{n,j,m,\nu}\|_{\mathfrak S_2}^2
    \le\frac1n\operatorname{tr}Q
     =\frac1n\mathfrak t_n^{\mathcal E}.
 \end{split}
\]
The mixed term vanishes in conditional expectation, but need not vanish on each path.

To interpret these statements on the original probability space, apply the isometry $\mathcal J_n$ into $\mathcal H_{0,n}$, which is measurable with respect to the sampled Brownian
positions and the terminal point at this size. The map $D\mathcal J_n^*$ acts on this fixed space and is zero on the orthogonal complement of the isometric image. Its covariance
products are the same as above. On
\[
 A_L=\{\operatorname{tr}Q+\|b\|_{n,+}^2\le L\}\in\mathscr S_n
\]
the conditional second moment is integrable after averaging, and $\mathbb E[\mathbf1_{A_L}(\|\mathcal Z_n\|_{n,+}^2+\|b\|_{n,+}^2)]\le L$. All Bochner conditional means and
covariance integrals used here are therefore defined on $A_L$. These events increase to an event of probability one. This localization justifies the identities without presuming an
unconditional moment bound for the random reference columns.

\end{proof}

Almost surely $e_n^{\mathrm{KL}}=O_{\beta,\omega}(1)$, and, for every $\varepsilon>0$ and all sufficiently large $n$,
\begin{equation}\label{S:eq:weighted-field-covariance-bound}
 \mathcal E_n^{+}(\mathcal B_n)+\mathfrak t_n^{\mathcal E}
 \le C_{\beta,\omega},
 \mathfrak l_n^{\mathcal E}
 \le C_\beta\max_{m_n\le k<n}\frac{\mathfrak G_{n,k}^{+}}{k}
 \le C_{\beta,\omega,\varepsilon}n^{-a_\beta+\varepsilon}.
\end{equation}
All scalar quantities in these statements are $\mathscr S_n$-measurable. The conjugated covariances $\mathcal J_n\widetilde Q_n^{+,(M)}\mathcal J_n^*$ and
$\mathcal J_n\widetilde Q_n^+\mathcal J_n^*$ are strongly measurable, trace-class operators on $\mathcal H_{0,n}$. They vanish on the orthogonal complement of
$\mathcal J_n\mathfrak H_n^+$.

Fix $0<\varepsilon_0<a_\beta$ and define the $\mathscr S_n$-measurable events
\[
 \Omega_{M,n}^{\mathcal E}:=
 \left\{\mathcal E_n^{+}(\mathcal B_n)+\mathfrak t_n^{\mathcal E}\le M,
  \mathfrak l_n^{\mathcal E}\le Mn^{-a_\beta+\varepsilon_0}\right\},
   M\in\mathbb N.
\]
By \eqref{S:eq:weighted-field-covariance-bound}, almost surely one finite integer $M_\omega$ contains these quantities for all sufficiently large $n$.

Apply Lemma~\ref{S:lem:covariance-gaussian-concentration} to $\mathcal J_n\mathcal Z_n\in\mathcal H_{0,n}$, with zero shift, under the conditional law given $\mathscr S_n$. Its
hypotheses follow from Proposition~\ref{S:prop:weighted-field-covariance}. In particular
\[
 \operatorname{tr}((\widetilde Q_n^+)^2)
 \le\mathfrak t_n^{\mathcal E}\mathfrak l_n^{\mathcal E}.
\]
For $t>0$, the lemma therefore gives
\begin{equation}\label{S:eq:s9-field-energy-concentration}
 \mathbb P\left(
 \left|\mathcal E_n^+(\mathcal Z_n)-\mathfrak t_n^{\mathcal E}\right|>t
      \,\middle|\,\mathscr S_n\right)
 \le2\exp\left[-c\min\left\{
 \frac{t^2}{\mathfrak t_n^{\mathcal E}\mathfrak l_n^{\mathcal E}},
 \frac{t}{\mathfrak l_n^{\mathcal E}}\right\}\right].
\end{equation}
If the covariance is zero, the field vanishes under that conditional law and the probability is zero. For a fixed integer $M\ge1$, on $\Omega_{M,n}^{\mathcal E}$ both denominators
are bounded by $M^2n^{-a_\beta+\varepsilon_0}$ and $Mn^{-a_\beta+\varepsilon_0}$, respectively. Taking $t=M$ shows
\[
 \mathbf1_{\Omega_{M,n}^{\mathcal E}}
 \mathbb P\left(\mathcal E_n^+(\mathcal Z_n)>2M
                    \,\middle|\,\mathscr S_n\right)
 \le2\exp(-c n^{a_\beta-\varepsilon_0}).
\]
The indicator belongs to $\mathscr S_n$. The tower property thus gives
\[
 \mathbb P\left(\Omega_{M,n}^{\mathcal E}
       \cap\{\mathcal E_n^+(\mathcal Z_n)>2M\}\right)
 \le2\exp(-c n^{a_\beta-\varepsilon_0}).
\]
The bound is summable, since $c n^{a_\beta-\varepsilon_0}\ge3\log n$ for all sufficiently large $n$. Lemma~\ref{S:lem:conditional-localized-borel-cantelli} applies with one index
per row, the localization $\Omega_{M,n}^{\mathcal E}$, and threshold $2M$. It gives one event on which the failure occurs only finitely often for every integer $M$. The bounds in
\eqref{S:eq:weighted-field-covariance-bound} then supply integers $M_\omega$ and $n_0(\omega)$ such that the corresponding localization holds for all $n\ge n_0(\omega)$. On that
event,
\begin{equation}\label{S:eq:bridge-weighted-norm-bound}
 \mathcal E_n^+(\mathcal Z_n)+\mathcal E_n^+(\mathcal B_n)
 \le3M_\omega=:C_{\beta,\omega}
\end{equation}
for all sufficiently large $n$. The Hilbert space triangle inequality and \eqref{S:eq:bridge-remainder-energy} now imply
$\mathcal E_n^+(\mathcal U_n)=O_\omega(1)$.

\subsubsection{Expansion of the inverse kernel}
\label{S:sec:s9-inverse-expansion}

The relative logarithm tends to zero in the essential supremum norm, and its weighted norm remains bounded. The weighted norm of the remainder tends to zero. These two estimates
bound the nonlinear terms without requiring uniform bounds on the boundary columns.

Use the sampled matrix function $X^{\mathrm{step}}$ as reference, and let $\mathscr D_n$ be the derivative in \eqref{S:eq:continuous-derivative-kernel}. Write
$a_n(r)=a_{n,k}^{\mathrm{br}}$ and $c_n(r)=c_{n,k}^{\mathrm{br}}$ on $I_{n,k}$. On the retained triangle $r>r'$ the first variation is the explicit kernel
\[
 (\mathscr D_nU)(r,r')=\tfrac12\bigl[
       (U(r)a_n(r))c_n(r')^{\mathsf T}
       +a_n(r)(U(r')c_n(r'))^{\mathsf T}\bigr].
\]
Its value on the other triangle is the transpose with interchanged arguments, and it vanishes outside the square of the restricted kernel. For a bounded matrix field $U$, put
\[
 E_U(r)=\me^{U(r)}-I,\qquad
 F_U(r)=\me^{U(r)}-I-U(r)
       =\int_0^1(1-v)\me^{vU(r)}U(r)^2\dd v.
\]
We define the nonlinear remainder by subtracting the constant and linear terms:
\[
 Q_n^{\mathrm{nl,br}}
 =\mathsf K_{X^{\mathrm{step}}}(\mathcal U_n)
      -\mathsf K_{X^{\mathrm{step}}}(0)-\mathscr D_n\mathcal U_n,
 \qquad
 R_n^{\mathrm{br}}=-\mathscr D_n\mathcal R_n-Q_n^{\mathrm{nl,br}}.
\]
In particular, with $U=\mathcal U_n$ its kernel on $r>r'$ is exactly
\[
 Q_n^{\mathrm{nl,br}}(r,r')=\tfrac12\bigl[ (F_U(r)a_n(r))c_n(r')^{\mathsf T} +a_n(r)(F_U(r')c_n(r'))^{\mathsf T}\bigr] +\tfrac12(E_U(r)a_n(r))(E_U(r')c_n(r'))^{\mathsf T}.
 \]
This follows by multiplying the two perturbed columns and subtracting the constant and linear terms. The expansion contains two terms in which one column is multiplied by $F_U$,
together with the product of the column increments. These formulas also apply to the two ordered triangles inside a single operator cell. The sampling error is the difference
between the inverse formed from the sampled coefficients and the perturbed inverse, which accounts for both minus signs in $R_n^{\mathrm{br}}$.

\begin{proposition}\label{S:prop:continuous-bridge-linearization}
Almost surely, for all sufficiently large $n$,
\begin{equation}\label{S:eq:bridge-operator-linearization}
 \varepsilon_{\beta,n}^{\mathrm{samp},T}
 =-\mathscr D_n\mathcal Z_n-\mathscr D_n\mathcal B_n+R_n^{\mathrm{br}},
 \qquad
 \sqrt n\,\|R_n^{\mathrm{br}}\|_{\mathfrak S_2}
       \le C_{\beta,\omega}(\log n)^{-250}.
\end{equation}
\end{proposition}
\begin{proof}
The singular inverse and its regular restriction use the same coordinate boundary vectors. The Green formulas in Theorems~\ref{S:thm:singular-inverse-domain}
and~\ref{S:thm:singular-regular-truncations} give
\[
 \mathsf R_\beta^{\mathrm{sine},T_n}
       =\mathsf K_{X^{\mathrm{step}}}(\mathcal U_n),
 \qquad
 \mathsf R_{\beta,n}^{\mathrm{br},T_n}
       =\mathsf K_{X^{\mathrm{step}}}(0).
\]
The pointwise identity of the coefficients is $X_{\theta_{n,k}(s)}=\me^{\mathcal U_{n,k}(s)}X_{n,k}^{\mathrm{step}}$. It holds on every selected interval, and all kernels are
extended by zero outside the square of the restricted kernel. Substituting $\mathcal U_n=\mathcal B_n+\mathcal Z_n+\mathcal R_n$ into the definitions of $Q_n^{\mathrm{nl,br}}$ and
$R_n^{\mathrm{br}}$ proves the first identity, including its two minus signs.

Equation~\eqref{S:eq:bridge-weighted-norm-bound} and \eqref{S:eq:bridge-remainder-energy} give
\[
 \mathcal E_n^+(\mathcal U_n)^{1/2}
 \le\mathcal E_n^+(\mathcal B_n)^{1/2}
   +\mathcal E_n^+(\mathcal Z_n)^{1/2}
   +\mathcal E_n^+(\mathcal R_n)^{1/2}=O_\omega(1).
\]
Also $\mathcal U_n$ is real and traceless, by its definition as the real logarithm of a positive diagonal triangular matrix of determinant one. Its essential supremum satisfies
$\eta(\mathcal U_n)\le C_{\beta,\omega}\sqrt{\log n/m_n}\to0$. For this fixed row the reference matrices and their inverses are bounded. Thus
$\mathcal U_n\in\mathfrak B_{n,X^{\mathrm{step}}}$, and $\eta(\mathcal U_n)\le1$ for all $n\ge n_0(\omega)$, since this quantity tends to zero. The derivative bound retains the
integrals within each partition interval. Indeed, for $r=(k+s)/n$ the exact primitives of the reference columns are $A_X(r)=A_{n,k}+s q_{n,k}/n$ and
$C_X(r)=C_{n,k}+(1-s)p_{n,k}/n$. The triangular kernel estimate gives
\[
 \begin{split}
 \|\mathscr D_nU\|_{\mathfrak S_2}^2
 \le\frac1n\sum_{k=m_n}^{n-1}\int_0^1\bigl[&
 (A_{n,k}+s q_{n,k}/n)|U_{n,k}(s)a_{n,k}|^2\\
 &+(C_{n,k}+(1-s)p_{n,k}/n)|U_{n,k}(s)c_{n,k}|^2\bigr]\dd s
 \le\frac1n\mathcal E_n^+(U).
 \end{split}
\]
Thus the normalized coordinate has measure $\dd s/n$, and the current cell is included even at $k=m_n$ or $k=n-1$. The factors $s$ and $1-s$ cannot be replaced by their averages
for a field depending on $s$. Put $U=\mathcal U_n$ and $\eta=\eta(U)\le1$. The matrix exponential commutes with $U$. Its integral formulas consequently give, for every vector $v$,
\[
 \begin{split}
 |E_Uv|&=\left|\int_0^1e^{tU}Uv\,\mathrm dt\right|
                   \le e^\eta|Uv|,\\
 |F_Uv|&=\left|\int_0^1(1-t)e^{tU}U^2v\,\mathrm dt\right|
                   \le\tfrac12\eta e^\eta|Uv|,\\
 |E_Uv|&\le\eta e^\eta|v|.
 \end{split}
\]
The first two terms of the exact remainder kernel form the linear kernel with the field $F_U$ in place of $U$. The preceding derivative estimate is an estimate of the explicit
triangular kernel and holds for arbitrary matrix fields, even though $F_U$ need not be traceless. It gives
\[
 \|\mathscr D_nF_U\|_{\mathfrak S_2}
 \le\frac1{\sqrt n}\mathcal E_n^+(F_U)^{1/2}
 \le\frac{\eta e^\eta}{2\sqrt n}\mathcal E_n^+(U)^{1/2}.
\]
Let $B_U$ denote the remaining bilinear kernel. On the triangle $r>r'$ it is
\[
 \tfrac12(E_U(r)a_n(r))(E_U(r')c_n(r'))^{\mathsf T},
\]
and on the other triangle it is the transposed reflection of this expression. The rank-one matrix norm and the equality of the two triangular integrals give
\[
 \begin{split}
 \|B_U\|_{\mathfrak S_2}^2
 &=\frac12\int_{r>r'}|E_U(r)a_n(r)|^2
                           |E_U(r')c_n(r')|^2\,\mathrm dr'\,\mathrm dr\\
 &\le\frac{\eta^2e^{4\eta}}2
    \int |U(r)a_n(r)|^2
           \left(\int_{m_n/n}^{r}|c_n(r')|^2\,\mathrm dr'\right)\mathrm dr\\
 &\le\frac{\eta^2e^{4\eta}}{2n}\mathcal E_n^+(U).
 \end{split}
\]
For the last step use $r=(k+s)/n$ and the exact primitive $A_X(r)=A_{n,k}+sq_{n,k}/n\le A_{n,k}^+$. The integral includes $r,r'$ in the same partition cell; only the measure-zero
diagonal is absent. Thus the integral retains the region where both variables belong to the same partition interval, with their stated ordering.

The triangle inequality, $Q_n^{\mathrm{nl,br}}=\mathscr D_nF_U+B_U$, and $\eta\le1$ now yield
\[
 \sqrt n\,\|Q_n^{\mathrm{nl,br}}\|_{\mathfrak S_2}
 \le\left(\tfrac12e+\tfrac1{\sqrt2}e^2\right)
       \eta(\mathcal U_n)\mathcal E_n^+(\mathcal U_n)^{1/2}.
\] Together with the
derivative bound for $\mathcal R_n$, this yields
\begin{equation}\label{S:eq:continuous-bridge-remainder-rate}
 \sqrt n\,\|R_n^{\mathrm{br}}\|_{\mathfrak S_2} \le\mathcal E_n^+(\mathcal R_n)^{1/2} +C\eta(\mathcal U_n)\mathcal E_n^+(\mathcal U_n)^{1/2} \le C_{\beta,\omega}\frac{(1+\log n)^6}{\sqrt{m_n}} \le C_{\beta,\omega}(\log n)^{-250}.
\end{equation}
For the second inequality, the term involving $\mathcal R_n$ is bounded by \eqref{S:eq:bridge-remainder-energy}, and the other term is at most $C_{\beta,\omega}\sqrt{\log n/m_n}$.
For $n\ge\me$, $(1+\log n)^6\le2^6(\log n)^6$ and $\sqrt{m_n}\ge(\log n)^{256}$, which proves the last inequality.

The full difference estimate of the same theorem also gives
\[
 n\|\varepsilon_{\beta,n}^{\mathrm{samp},T}\|_{\mathfrak S_2}^2
 \le C\mathcal E_n^+(\mathcal U_n)\le C_{\beta,\omega}.
\]
The kernel maps are strongly measurable by Section~\ref{S:sec:stability-and-fr-echet-expansion-of-the-inverse-kernel}. The reference derivative is determined by $\mathscr S_n$,
while the remainders are measurable on the full probability
space. \end{proof}

We now remove the cutoff using only the Brownian estimates of Proposition~\ref{S:prop:4.3}. Put $\delta_n=m_n/n$, let $P=P_{\delta_n}^{\mathrm e}$, and write
$\Delta_n=\mathsf R_{\beta,n}^{\mathrm{br}}-\mathsf R_\beta^{\mathrm{sine}}$. Then $P\Delta_nP=\varepsilon_{\beta,n}^{\mathrm{samp},T}$ and
\[
 H_n:=\Delta_n-P\Delta_nP
    =\mathsf T_{\beta,n}^{\mathrm{br}}-\mathsf T_{\beta,n}^{\mathrm{sine}}.
\]
The kernel of $P\Delta_nP$ is supported on $(\delta_n,1)^2$ and that of $H_n$ on its complement. Orthogonality in the space of Hilbert--Schmidt kernels, as proved in
Lemma~\ref{S:lem:orthogonal-truncation-masses}, gives
\[
 n\|\Delta_n\|_{\mathfrak S_2}^2
 =n\|\varepsilon_{\beta,n}^{\mathrm{samp},T}\|_{\mathfrak S_2}^2
       +n\|H_n\|_{\mathfrak S_2}^2.
\]
Fix $0<\varepsilon<a_\beta$. The two Brownian tail estimates imply
\[
 n\|H_n\|_{\mathfrak S_2}^2
 \le 2n\bigl(\|\mathsf T_{\beta,n}^{\mathrm{br}}\|_{\mathfrak S_2}^2
                  +\|\mathsf T_{\beta,n}^{\mathrm{sine}}\|_{\mathfrak S_2}^2\bigr)
 \le C_{\beta,\omega,\varepsilon}
               n^{-a_\beta+\varepsilon}m_n^{1+a_\beta-\varepsilon}
 \longrightarrow0.
\]
For $0<\beta\le2$ the tail exponent is $2-\varepsilon$; for $\beta>2$ it is $1+2/\beta-\varepsilon$. Fix the value of $\beta$ and put $a=a_\beta>0$. Choosing $\varepsilon=a/2$
makes the exponent $1+a/2>1$ in both cases. The preceding estimate therefore gives
\[
 \sqrt n\,\|H_n\|_{\mathfrak S_2}
 \le C_{\beta,\omega}n^{-a/4}m_n^{1/2+a/4}.
\]
For $n\ge e$, $m_n=\lceil(\log n)^{512}\rceil\le2(\log n)^{512}$. Consequently
\[
 (\log n)^{250}\sqrt n\,\|H_n\|_{\mathfrak S_2}
 \le C_{\beta,\omega}
     n^{-a/4}(\log n)^{506+128a}.
\]
For this fixed $a$, choose a deterministic $N_a$ such that $(506+128a)\log\log n\le(a/8)\log n$ for every $n\ge N_a$. Such an index exists because $\log\log n/\log n\to0$. The last
right hand side is then at most $C_{\beta,\omega}n^{-a/8}\to0$. Taking the maximum with the finite random index supplied by the tail estimate proves
\[
 \sqrt n\,\|H_n\|_{\mathfrak S_2}=o_\omega((\log n)^{-250}).
\]
The constants and the starting index may depend on $\beta$; this argument makes no uniform claim as $a_\beta\downarrow0$.

Extend $\mathscr D_n\mathcal B_n$ and $\mathscr D_n\mathcal Z_n$ by zero as before and set $\widehat R_n^{\mathrm{br}}=R_n^{\mathrm{br}}+H_n$. The identity on the full interval is
\[
 \Delta_n=-\mathscr D_n(\mathcal B_n+\mathcal Z_n)
                     +\widehat R_n^{\mathrm{br}},\qquad
 \sqrt n\|\widehat R_n^{\mathrm{br}}\|_{\mathfrak S_2}
       \le C_{\beta,\omega}(\log n)^{-250}.
\]
Consequently $n\|\Delta_n\|_{\mathfrak S_2}^2\le C_{\beta,\omega}$. Moreover, writing $L_n=\mathscr D_n(\mathcal B_n+\mathcal Z_n)$, the weighted norm bound gives
$\sqrt n\|L_n\|_{\mathfrak S_2}=O_\omega(1)$, and Cauchy--Schwarz gives
\[
 \left|n\|\Delta_n\|_{\mathfrak S_2}^2-n\|L_n\|_{\mathfrak S_2}^2\right|
 \le2(\sqrt n\|L_n\|_{\mathfrak S_2})
              (\sqrt n\|\widehat R_n^{\mathrm{br}}\|_{\mathfrak S_2})
       +n\|\widehat R_n^{\mathrm{br}}\|_{\mathfrak S_2}^2
\longrightarrow0.
\]
This proves the order of the sampling error and the bound for the remainder on the full interval in Corollary~\ref{S:cor:gaussian-full-interval}. The squared norm of $L_n$ is the
sum of the squared norms of the linear interpolation and Gaussian terms and twice their Hilbert--Schmidt inner product on each path. Their leading order coefficient is calculated
in Section~\ref{S:sec:the-limiting-coefficient-of-the-error}.

\clearpage

\section{Coupling matrix coefficients to Brownian motion}\label{sec:coupling-matrix-coefficients-to-brownian-motion}

We construct the coefficient vectors jointly with Brownian motion, preserving their prescribed joint distribution at every matrix size. Independence makes this direct for the
Gaussian and beta-Laguerre models. Polynomial potentials require conditional constructions, while the circular model requires control of products of coupled increments.

The constructions of the tridiagonal coefficient arrays in Sections~\ref{sec:combined-a-gaussian}--\ref{sec:combined-a-gaussian-comparison} proceed from independent coefficients to
conditional Jacobi laws. At the hard edge, Section~\ref{sec:the-hard-edge-coefficient-construction} identifies the errors needed by the inverse formula;
Section~\ref{B:sec:the-beta-laguerre-coupling} treats beta-Laguerre, and
Sections~\ref{B:sec:general-potentials-and-convex-reference-measures}--\ref{B:sec:normalization-and-coefficient-estimates} treat dependent coefficients. The circular comparison
passes from increments to transfer matrix products and then inverses in Sections~\ref{sec:combined-s-increments}--\ref{sec:combined-s-finite-comparison}. Parameters and Brownian
probability spaces remain those of the respective main theorems.

Section~\ref{sec:finite-gaussian-estimates-for-the-wigner-comparison} proves the estimates for the GOE and GUE coefficient process, the Gram matrices, and the expected traces used
in the comparison for Wigner matrices. That proof combines these finite size estimates with qualitative convergence of the Gaussian ensembles to the stochastic Airy operator.

\subsection{The Gaussian $\beta$-ensemble}
\label{sec:combined-a-gaussian}
\label{A:sec:the-gaussian-beta-ensemble}
We decompose the centered tridiagonal matrix into a diagonal multiplier and the operator $\Gamma_{\eta_n}^-$ defined below. Each entry of the diagonal multiplier depends on a
single independent coefficient pair. The weighted estimates in Propositions~\ref{A:res:3.16} and~\ref{A:res:3.17} then control the two terms in the comparison with Brownian
increments.

Throughout the construction for the Gaussian $\beta$-ensemble, $\beta>0$ is fixed and $h=n^{-1/3}$. We retain the coordinates $x_k=(k-\tfrac12)h$, the difference operator $D_h$,
and the embeddings of Section~\ref{A:sec:weighted-resolvent-estimates-and-finite-difference-approximation}. Every restriction uses the original coordinates and the same value of
$h$. For resolvent estimates, fix a compact set $K\Subset\mathbb C\setminus\mathbb R$; replacing it by $K\cup\overline K$, we may assume that it is invariant under complex
conjugation. Set $\delta_K=\operatorname{dist}(K,\mathbb R)>0$.

The comparison of the finite matrices is stated in Theorem~\ref{A:res:8.5}. We compare the mean coefficients with the coefficients of the deterministic finite difference operator
on the initial half of the matrix. Decay of the resolvent kernel controls the remaining coordinates.

\subsubsection{Coefficient decomposition}\label{A:sec:dumitriu-edelman-decomposition}
The quadratic form gives the decomposition below. Auxiliary boundary coefficients make the identity uniform across the matrix; they are not additional random entries.

We use counting measure on the vertex space $\ell^2(\{1,\ldots,n\})$ and the edge space $\ell^2(\{0,\ldots,n\})$. Let $a_{k,n},b_{k,n}$ be the independent coefficients of
$T_{n,\beta}^{\mathrm{DE}}$ specified in Subsection~\ref{A:sec:intro-models-scaling}. For $\nu>0$, set $m_\nu=\mathbb E\chi_\nu$ and $X_\nu=\chi_\nu-m_\nu$, with
$\chi_0=m_0=X_0=0$.

We use the auxiliary coefficients
\[
 b_{0,n}=b_{n,n}=1.
 \eqtag{6.2}
\]
They are not entries of $T_{n,\beta}^{\mathrm{DE}}$. For $u_0=u_{n+1}=0$, the matrix quadratic form is
\[
 \mathfrak a_{n,\beta}^{\mathrm{DE}}[u]
 =h^{-2}\sum_{k=1}^n(2-a_{k,n})|u_k|^2
      -2h^{-2}\operatorname{Re}\sum_{k=1}^{n-1}b_{k,n}u_{k+1}\overline{u_k}.
 \eqtag{6.4}
\]
For any choice of the two auxiliary coefficients,
\[
 \sum_{k=0}^n b_{k,n}|u_{k+1}-u_k|^2
 =\sum_{k=1}^n(b_{k-1,n}+b_{k,n})|u_k|^2
       -2\operatorname{Re}\sum_{k=1}^{n-1}b_{k,n}u_{k+1}\overline{u_k}.
\]
It follows that
\[
 \mathfrak a_{n,\beta}^{\mathrm{DE}}[u]
 =\sum_{k=0}^n b_{k,n}|(D_hu)_k|^2
    +h^{-2}\sum_{k=1}^n(2-a_{k,n}-b_{k-1,n}-b_{k,n})|u_k|^2.
 \eqtag{6.5}
\]
In particular, the contributions of $b_{0,n}$ and $b_{n,n}$ cancel between the two sums. Choosing them to be one does not change the matrix.

Write
\[
 \overline b_{k,n}=\frac{m_{\beta(n-k)}}{\sqrt{\beta n}},\qquad
 \delta_{k,n}=b_{k,n}-\overline b_{k,n}\quad(1\le k<n),
\]
and set $\overline b_{0,n}=\overline b_{n,n}=1$ and $\delta_{0,n}=\delta_{n,n}=0$. Let $\overline H_{n,\beta}$ have form
\[
 \overline{\mathfrak h}_{n,\beta}[u]
 =\sum_{k=0}^n\overline b_{k,n}|(D_hu)_k|^2
    +h^{-2}\sum_{k=1}^n(2-\overline b_{k-1,n}-\overline b_{k,n})|u_k|^2.
 \eqtag{6.11}
\]
Thus $\overline H_{n,\beta}=\mathbb E\mathcal A_{n,\beta}^{\mathrm{DE}}$. In the divergence form representation \eqref{A:eq:6.5}, the scalar fluctuation at vertex $k$ contains both
$\delta_{k-1,n}$ and $\delta_{k,n}$. To give each scalar coefficient its own independent pair, we write on each edge $u_{k+1}=u_k+h(D_hu)_k$ in the off-diagonal fluctuation. The
part without a difference then involves only $(a_{k,n},\delta_{k,n})$; the other part is estimated by the first-order resolvent bound.

The scalar coefficient and the coefficient of the first-order form are
\[
 \zeta_{k,n}=-a_{k,n}-2\delta_{k,n}\quad(1\le k\le n),
 \eqtag{6.14}
\]
\[
 \eta_{k,n}=h^{-1}\delta_{k,n}\quad(0\le k\le n),\qquad
 \Gamma_{\eta_n}^{-}
 =-D_h^*M_{\eta_n}J_- -J_-^*M_{\eta_n}D_h.
\]
Here $M_{\zeta_n}$ acts on the vertex space and $M_{\eta_n}$ on the edge space. The map $J_-$ and the adjoint of $D_h$ use counting measure, as in
Subsection~\ref{A:sec:weighted-resolvent-estimates-and-comparison-principles}, with $N_n=n$. For vectors $u,v$ extended by zero at vertices $0,n+1$, the Hermitian sesquilinear form
of $\Gamma_{\eta_n}^{-}$ is
\[
 \langle\Gamma_{\eta_n}^{-}u,v\rangle
 =-\sum_{k=0}^n\eta_{k,n}
       \bigl((D_hu)_k\overline{v_k}
                   +u_k\overline{(D_hv)_k}\bigr).
\]

\begin{proposition}\label{A:res:6.1}
For every $n$,
\[
 \mathcal A_{n,\beta}^{\mathrm{DE}}
 =\overline H_{n,\beta}+h^{-2}M_{\zeta_n}+\Gamma_{\eta_n}^{-}.
\]
\end{proposition}
\begin{proof}
Subtracting \eqref{A:eq:6.11} from \eqref{A:eq:6.5} gives
\[
 \mathfrak r_{n,\beta}[u] :=\mathfrak a_{n,\beta}^{\mathrm{DE}}[u] -\overline{\mathfrak
 h}_{n,\beta}[u] =h^{-2}\sum_{k=0}^n\delta_{k,n}|u_{k+1}-u_k|^2
 -h^{-2}\sum_{k=1}^n(a_{k,n}+\delta_{k-1,n}+\delta_{k,n})|u_k|^2.
 \eqtag{6.12}
\]
Expansion of the squared differences cancels the terms $\delta_{k-1,n}+\delta_{k,n}$, leaving
\[
 \mathfrak r_{n,\beta}[u]
 =-h^{-2}\sum_{k=1}^na_{k,n}|u_k|^2
       -2h^{-2}\operatorname{Re}\sum_{k=0}^n\delta_{k,n}u_{k+1}\overline{u_k}.
 \eqtag{6.13}
\]
At each edge,
\[
 -2\delta_{k,n}\operatorname{Re}(u_{k+1}\overline{u_k})
 =-2\delta_{k,n}|u_k|^2
      -2\delta_{k,n}\operatorname{Re}((u_{k+1}-u_k)\overline{u_k}).
\]
By the definition of $\eta_n$ and the adjoints of $D_h$ and $J_-$,
\[
 \langle u,\Gamma_{\eta_n}^{-}u\rangle
 =-2\operatorname{Re}\sum_{k=0}^n\eta_{k,n}(D_hu)_k\overline{u_k}
 =-2h^{-2}\operatorname{Re}\sum_{k=0}^n
            \delta_{k,n}(u_{k+1}-u_k)\overline{u_k}.
 \eqtag{6.18}
\]
Substitution in \eqref{A:eq:6.13} gives $\mathfrak r_{n,\beta}[u] =h^{-2}\sum_k\zeta_{k,n}|u_k|^2+\langle u,\Gamma_{\eta_n}^{-}u\rangle$. All operators are Hermitian, so equality
of quadratic forms implies the matrix identity by polarization. The conventions $\delta_{0,n}=\delta_{n,n}=0$ include both boundary edges. The diagonal entry $2h^{-2}\delta_{k,n}$
of $\Gamma_{\eta_n}^{-}$ cancels the term $-2h^{-2}\delta_{k,n}$ in $h^{-2}M_{\zeta_n}$. Its off-diagonal entry at $(k,k+1)$ is $-h^{-2}\delta_{k,n}$, as required by the original
matrix.
\end{proof}

The variables $\zeta_{k,n}$ are independent over $k$, since each depends only on the independent pair $(a_{k,n},\delta_{k,n})$. The fields $\zeta_n$ and $\eta_n$ share $\delta_n$;
the resolvent comparison therefore estimates their contributions without assuming independence.

The leading variances determine the normalization of the coupling. By Lemma~\ref{A:res:A.1}, for every fixed $L$, uniformly over $x_k\le L$ as $n\to\infty$,
\[
 \operatorname{Var}(a_{k,n}) =\frac2{\beta n},\quad \operatorname{Var}(\delta_{k,n})
 =\frac1{\beta n}\left(\frac12+O_\beta(n^{-1})\right),\quad \operatorname{Var}(\zeta_{k,n})
 =\frac4{\beta n}+O_\beta(n^{-2}) =\frac4\beta h^3+O_\beta(h^6).
\]
Indeed $k\le Lh^{-1}+1$ implies $\beta(n-k)\ge\beta n/2$ for large $n$, and $a_{k,n}$ is independent of $\delta_{k,n}$. The corresponding coefficients constructed from Brownian
increments are
\[
 \zeta_{k,n}^{B}=\frac2{\sqrt\beta}h\bigl(B(kh)-B((k-1)h)\bigr),
 \qquad
 h^{-2}\zeta_{k,n}^{B}
    =\frac2{\sqrt\beta}\frac{B(kh)-B((k-1)h)}h.
\]
Their variances are $4h^3/\beta$ and $4/(\beta h)$, respectively.

This variance match also gives the required distributional potential. For this identification only, fix $f\in C_c^\infty((0,\infty))$ and put $u_k=\sqrt h f(x_k)$. For sufficiently
large $n$, its support lies inside the matrix interval. With $g_h(x)=|f(x_k)|^2$ on $[(k-1)h,kh)$ and zero outside these intervals,
\[
 h^{-2}\sum_{k=1}^n\zeta_{k,n}^B|u_k|^2
 =\frac2{\sqrt\beta}\int_0^\infty g_h(x)\,\mathrm dB(x).
\]
The functions $g_h$ and $g=|f|^2$ are supported in one fixed compact set, and $|g_h(x)-g(x)|\le h\|g'\|_\infty/2$. The It\^o isometry gives
\[
 \mathbb E\left|h^{-2}\sum_{k=1}^n\zeta_{k,n}^B|u_k|^2
       -\frac2{\sqrt\beta}\int_0^\infty|f(x)|^2\,\mathrm dB(x)\right|^2
 =\frac4\beta\int_0^\infty|g_h-g|^2\,\mathrm dx
 \le C_{f,\beta}h^2.
\]
Integration by parts identifies the limiting integral with $-(2/\sqrt\beta)\int B(x)(|f(x)|^2)'\,\mathrm dx$, as in \cite{ramirez2011beta}. Its variance is $(4/\beta)\int|f|^4>0$
for $f\ne0$. This calculation with a fixed test function identifies the noise coefficient. Operator convergence will follow from the weighted comparison below. At the final vertex,
$\zeta_{n,n}=-a_{n,n}$ has variance $2/(\beta n)$, half the interior variance $4/(\beta n)$. Proposition~\ref{A:res:A.5} bounds this variance difference in the weighted error at
the far boundary point.

\subsubsection{Moments of the chi distribution and quantile coupling}\label{A:sec:chi-moments-and-gaussian-quantile-coupling}
To compare the mean coefficients, we need an expansion of $m_\nu$. To compare their fluctuations, we need uniform centered moments and a coupling of $X_\nu$ to a normal variable.
These two estimates must include small $\nu$: the degrees of freedom $\beta(n-k)$ do not tend to infinity near the last matrix coordinate.

For $\nu>0$, the chi law has density
\[
 p_\nu(t)=\frac{2^{1-\nu/2}}{\Gamma(\nu/2)}
             t^{\nu-1}\mathrm e^{-t^2/2}\mathbf1_{\{t>0\}}.
\]
The substitution $y=t^2/2$ gives
\[
 \mathbb E\chi_\nu^q
 =2^{q/2}\frac{\Gamma((\nu+q)/2)}{\Gamma(\nu/2)}
 \quad(q\ge0).
\]

\par\Needspace{5\baselineskip}
\begin{lemma}\label{A:res:A.1}
As $\nu\to\infty$,
\[
 m_\nu=\sqrt\nu\left(1-\frac1{4\nu}+\frac1{32\nu^2}
                         +O(\nu^{-3})\right),\qquad
 \operatorname{Var}(\chi_\nu)=\frac12-\frac1{8\nu}+O(\nu^{-2}).
\]
\end{lemma}
\begin{proof}
Apply the classical Gamma ratio expansion \cite{nistGammaAsymptotics} (\S5.11(iii), equations~(5.11.13) and~(5.11.15)) with $x=\nu/2$, $a=1/2$ and $b=0$:
\[
 \frac{\Gamma(x+1/2)}{\Gamma(x)}
 =x^{1/2}\left(1-\frac1{8x}+\frac1{128x^2}+O(x^{-3})\right).
\]
Since $m_\nu=\sqrt2\,\Gamma((\nu+1)/2)/\Gamma(\nu/2)$, this gives the expansion of $m_\nu$. The second follows by squaring and using $\operatorname{Var}(\chi_\nu)=\nu-m_\nu^2$.
\end{proof}

The expansion above applies as $\nu\to\infty$. Near the last matrix coordinate, we use a concentration estimate uniform in $\nu$. For a real random variable $Y$, use the norm
$\|Y\|_{\psi_1}=\inf\{a>0:\mathbb E\exp(|Y|/a)\le2\}$.

\begin{lemma}\label{A:res:A.2}
There are absolute constants $c,C>0$ such that, for all $p\ge1$ and $t\ge0$,
\[
 \sup_{\nu\ge0}\|X_\nu\|_{L^p}\le C\sqrt p,
 \eqtag{A.7}
\]
\[
 \sup_{\nu\ge0}\mathbb P(|X_\nu|>t)\le2\mathrm e^{-ct^2}.
 \eqtag{A.8}
\]
Also $\sup_{\nu\ge0}\|X_\nu^2-\mathbb EX_\nu^2\|_{\psi_1}<\infty$.
\end{lemma}
\begin{proof}
For $\nu>0$ and $Y=\chi_\nu^2$, integration of the Gamma density gives $\mathbb E\mathrm e^{tY}=(1-2t)^{-\nu/2}$ for $t<1/2$. For $y>\nu$, minimize $-ty-(\nu/2)\log(1-2t)$ over
$0<t<1/2$; for $0<y<\nu$, minimize it over $t<0$. In both cases the minimizer is $t=(1-\nu/y)/2$, and hence
\[
 \mathbb P(Y\ge y)\le\mathrm e^{-I_\nu(y)}\quad(y>\nu),\quad \mathbb P(Y\le y)\le\mathrm
 e^{-I_\nu(y)}\quad(0<y<\nu),\quad I_\nu(y)=\frac{y-\nu-\nu\log(y/\nu)}2.
\]
For $u=t/\sqrt\nu\ge0$, direct substitution gives
\[
 I_\nu((\sqrt\nu+t)^2)
 =\nu\left(u+\frac{u^2}2-\log(1+u)\right)\ge\frac{t^2}2.
\]
For $u=t/\sqrt\nu\in(0,1)$, substitution of $y=\nu(1-u)^2$ into $I_\nu(y)$ gives
\[
 I_\nu((\sqrt\nu-t)^2)
 =\nu\left(-u+\frac{u^2}2-\log(1-u)\right)\ge\frac{t^2}2.
\]
Here we used $\log(1+u)\le u$ and $\log(1-u)\le-u$. The lower event is empty when $t\ge\sqrt\nu$. Therefore
\[
 \mathbb P(|\chi_\nu-\sqrt\nu|>t)\le2\mathrm e^{-t^2/2},\qquad
 |m_\nu-\sqrt\nu|
 \le\int_0^\infty2\mathrm e^{-t^2/2}\,\mathrm dt=:C_0.
\]
For $t\ge2C_0$, it follows that $\mathbb P(|X_\nu|>t)\le2\mathrm e^{-t^2/8}$. Choose $c\le1/8$ also satisfying $4cC_0^2\le\log2$. For $t<2C_0$, the right side of \eqref{A:eq:A.8}
is then at least one. This proves \eqref{A:eq:A.8} for every $\nu>0$. At $\nu=0$, $X_0=0$ almost surely, and $\mathbb P(|X_0|>t)=0$ for $t\ge0$. Finally,
\[
 \mathbb E|X_\nu|^p
 \le2p\int_0^\infty t^{p-1}\mathrm e^{-ct^2}\,\mathrm dt
 =2c^{-p/2}\Gamma(1+p/2).
\]
Stirling's upper bound for $p\ge2$, and monotonicity of $L^p$ norms for $1\le p\le2$, imply \eqref{A:eq:A.7}.

Put $Z_\nu=X_\nu^2-\mathbb EX_\nu^2$. Integrating \eqref{A:eq:A.8} at the even integer moments gives
\[
 \mathbb E|X_\nu|^{2m}\le2c^{-m}m!,\qquad
 \mathbb E|Z_\nu|^m
 \le2^{m-1}\bigl(\mathbb E|X_\nu|^{2m}
                       +(\mathbb EX_\nu^2)^m\bigr)
 \le C_1^m m!\quad(m\ge1).
\]
Thus $\mathbb E\exp(|Z_\nu|/(2C_1)) \le1+\sum_{m\ge1}2^{-m}=2$. Since $\mathbb EZ_\nu=0$, the same power series gives
\[
 \mathbb E\mathrm e^{tZ_\nu}
 \le1+\sum_{m\ge2}(C_1|t|)^m
 \le\exp(2C_1^2t^2),\qquad |t|\le(2C_1)^{-1}.
\]
The constants do not depend on $\nu$, including $\nu=0$.
\end{proof}

We also record a consequence for the squared variables. Since $\sup_\nu\mathbb EX_\nu^2\le C_2$, for $t\ge2C_2$ the event $|X_\nu^2-\mathbb EX_\nu^2|>t$ can occur only if
$X_\nu^2>t$. Equation~\eqref{A:eq:A.8}, followed by an increase of the constant for $0\le t<2C_2$, gives
\[
 \mathbb P\bigl(|X_\nu^2-\mathbb EX_\nu^2|>t\bigr)
 \le C\exp\{-c\min(t,t^2)\},\qquad \nu\ge0,\ t\ge0.
 \eqtag{A.9}
\]

Let $F_\nu$ be the distribution function of $X_\nu$, and let $\Phi_{\mathrm N}$ be the standard normal distribution function. Define
\[
 F_\nu^{-1}(u)=\inf\{x:F_\nu(x)\ge u\}\quad(0<u<1),\qquad
 Q_\nu(z)=F_\nu^{-1}(\Phi_{\mathrm N}(z)).
\]
Thus $Q_0=0$. If $Z\sim N(0,1)$, then $\Phi_{\mathrm N}(Z)$ is uniform on $(0,1)$ and $Q_\nu(Z)$ has the law of $X_\nu$. For $1\le p<\infty$ and probability measures $\mu,\lambda$
on $\mathbb R$ with finite $p$th moments, write $W_p(\mu,\lambda)=\inf\{\|U-V\|_{L^p}:\operatorname{Law}(U)=\mu,\operatorname{Law}(V)=\lambda\}$, where the infimum is over joint
laws. The coupling fixes the Gaussian variable determined by the Brownian increment, so we need a single coupling that satisfies all the moment estimates below.

\par\Needspace{5\baselineskip}
\begin{lemma}\label{A:res:A.3}
For $Z\sim N(0,1)$ and every $1\le p<\infty$,
\[
 W_p\bigl(\operatorname{Law}(X_\nu),N(0,1/2)\bigr)
 =\left\|Q_\nu(Z)-\frac Z{\sqrt2}\right\|_{L^p}
 \le\frac{C_p}{\sqrt{1+\nu}},\qquad \nu\ge0.
 \eqtag{A.12}
\]
\end{lemma}
\begin{proof}
Fix $1\le p<\infty$. A Stein kernel for a centered real random variable $Y$ is a measurable function $\tau$ with $\tau(Y)\in L^1$ satisfying
$\mathbb E[Yf(Y)]=\mathbb E[\tau(Y)f'(Y)]$ for $f\in C_c^\infty(\mathbb R)$. Proposition~3.4 of \cite{ledoux2015stein}, applied in dimension one, gives
\[
 W_p\bigl(\operatorname{Law}(Y),N(0,1)\bigr)
 \le\|Z_0\|_{L^p}\,\|\tau(Y)-1\|_{L^p}.
 \eqtag{A.12a}
\]
Here $Z_0\sim N(0,1)$, $\tau(Y)\in L^p$, and $Y\in L^p$ when $p\ge2$. The displayed scalar Stein identity is the one-dimensional instance of the vector identity required there. The
monotone quantile coupling is optimal for the convex cost $|x-y|^p$; see \cite{santambrogio2015optimal}, Section~2.2. Consequently, if $F_\mu,F_\lambda$ are the distribution
functions of $\mu,\lambda$ and $U$ is uniform on $(0,1)$,
\[
 W_p(\mu,\lambda)=\|F_\mu^{-1}(U)-F_\lambda^{-1}(U)\|_{L^p}.
\]
This formula includes $p=1$, so the same quantile coupling can be used for every finite moment order.

We apply \eqref{A:eq:A.12a} to the Gamma variable underlying the chi law. Suppose $\nu\ge2$, put $r=\nu/2$, and set
\[
 G=\frac{\chi_\nu^2}{2}\sim\operatorname{Gamma}(r,1),\qquad
 Y=\frac{G-r}{\sqrt r}=\frac{\chi_\nu^2-\nu}{\sqrt{2\nu}}.
\]
The concentration estimate proved in Lemma~\ref{A:res:A.2} implies $\sup_{\nu>0}\|\chi_\nu-\sqrt\nu\|_{L^q}\le C_q$ for every finite $q\ge1$. The identity
\[
 Y=\sqrt2(\chi_\nu-\sqrt\nu)
       +\frac{(\chi_\nu-\sqrt\nu)^2}{\sqrt{2\nu}}
\]
therefore gives, for $\nu=2r\ge2$,
\[
 \|Y\|_{L^q}
 \le\sqrt2\|\chi_\nu-\sqrt\nu\|_{L^q}
       +\frac{\|\chi_\nu-\sqrt\nu\|_{L^{2q}}^2}{\sqrt{2\nu}}
 \le C_q.
 \eqtag{A.12b}
\]
The Gamma density $g_r(x)=x^{r-1}\mathrm e^{-x}/\Gamma(r)$ satisfies $(xg_r(x))'=(r-x)g_r(x)$. Its boundary term $xg_r(x)\psi(x)$ vanishes at zero, since
$xg_r(x)=x^r\mathrm e^{-x}/\Gamma(r)$ and $r>0$, and at infinity for every smooth $\psi$ of polynomial growth. Integration by parts, with $\psi(x)=f((x-r)/\sqrt r)$, gives
\[
 \mathbb E[(G-r)\psi(G)]=\mathbb E[G\psi'(G)],\qquad
 \mathbb E[Yf(Y)]=\mathbb E\left[\left(1+\frac Y{\sqrt r}\right)f'(Y)\right].
\]
For smooth $f$ whose value and derivative have polynomial growth, apply the compactly supported identity to $f\theta_R$, where $\theta_R$ is smooth, $0\le\theta_R\le1$,
$\theta_R=1$ on $[-R,R]$, $\theta_R=0$ outside $[-2R,2R]$, and $|\theta_R'|\le C/R$. The extra derivative term is bounded by
\[
 \left|\mathbb E\!\left[\left(1+\frac Y{\sqrt r}\right)
                      f(Y)\theta_R'(Y)\right]\right|
 \le\frac C R\mathbb E\bigl[(1+|Y|)|f(Y)|
                         \mathbf1_{\{|Y|\ge R\}}\bigr]\longrightarrow0.
\]
All other terms converge by dominated convergence and \eqref{A:eq:A.12b}. In particular, $\mathbb EY=0$, $\tau_r(y)=1+y/\sqrt r$ is a Stein kernel, and \eqref{A:eq:A.12b} gives all
the moment hypotheses above. Equations~\eqref{A:eq:A.12a}--\eqref{A:eq:A.12b} imply
\[
 W_p\bigl(\operatorname{Law}(Y),N(0,1)\bigr)
 \le\frac{\|Z_0\|_{L^p}\|Y\|_{L^p}}{\sqrt r}
 \le\frac{C_p}{\sqrt\nu}.
 \eqtag{A.12c}
\]
Rationalizing twice gives, for every $u\ge-1$,
\[
 \sqrt{1+u}-1-\frac u2=-\frac{u^2}{2(1+\sqrt{1+u})^2}.
\]
Since $\chi_\nu=\sqrt\nu(1+Y/\sqrt r)^{1/2}$, it follows that
\[
 R_\nu:=\chi_\nu-\sqrt\nu-\frac Y{\sqrt2}
 =-\frac{Y^2}{\sqrt\nu(1+\sqrt{1+Y/\sqrt r})^2},\qquad
 \|R_\nu\|_{L^p}\le\nu^{-1/2}\|Y\|_{L^{2p}}^2.
\]
Since $\mathbb EY=0$, $\mathbb ER_\nu=m_\nu-\sqrt\nu$ and $X_\nu=Y/\sqrt2+R_\nu-\mathbb ER_\nu$. Let $F_Y$ denote the distribution function of $Y$, and put
$\widetilde Y=F_Y^{-1}(\Phi_{\mathrm N}(Z))$. The map $y\mapsto\sqrt{\nu+\sqrt{2\nu}\,y}-m_\nu$ is continuous and strictly increasing on the support $[-\sqrt r,\infty)$ of $Y$. It
therefore sends its quantiles to those of $X_\nu$:
\[
 F_\nu^{-1}(u)=\sqrt{\nu+\sqrt{2\nu}\,F_Y^{-1}(u)}-m_\nu,\qquad 0<u<1.
\]
With $\widetilde R_\nu=\sqrt{\nu+\sqrt{2\nu}\,\widetilde Y}-\sqrt\nu-\widetilde Y/\sqrt2$, we obtain
\[
 Q_\nu(Z)=\frac{\widetilde Y}{\sqrt2}
                  +\widetilde R_\nu-\mathbb E\widetilde R_\nu,
 \qquad
 \|\widetilde R_\nu\|_{L^p}=\|R_\nu\|_{L^p}.
\]
The same $Z$ therefore satisfies
\[
 \left\|Q_\nu(Z)-\frac Z{\sqrt2}\right\|_{L^p}
 \le\frac1{\sqrt2}\|\widetilde Y-Z\|_{L^p}
                  +2\|\widetilde R_\nu\|_{L^p}
 \le C_p\nu^{-1/2}.
\]
For $0\le\nu\le2$, Lemma~\ref{A:res:A.2} bounds $\|Q_\nu(Z)\|_{L^p}+\|Z/\sqrt2\|_{L^p}$ by $C_p$. Since $(1+\nu)^{-1/2}\ge1/\sqrt3$ on this interval, the triangle inequality proves
\eqref{A:eq:A.12} there as well.
\end{proof}

Lemma~\ref{A:res:A.1} and $\operatorname{Var}(X_\nu)\le\mathbb E\chi_\nu^2=\nu$ on bounded intervals give
\[
 \left|\operatorname{Var}(X_\nu)-\frac12\right|\le\frac C{1+\nu},
 \qquad \nu\ge0.
\]

\subsubsection{Coupling with Brownian motion}\label{A:sec:coupling-the-dumitriu-edelman-fluctuations-to-brownian-increments}
The Gaussian approximation to $\zeta_{k,n}$ contains two independent normal variables with equal coefficients. We construct them as an orthogonal transformation of a normalized
Brownian increment and an auxiliary independent normal variable. Their sum recovers the Brownian increment with the sign and normalization in \eqref{A:eq:8.18}.

We first check that the quantile map is jointly measurable in its parameter and Gaussian argument. This makes the matrix coefficients measurable functions of $\beta$, as required
in Subsection~\ref{A:sec:moment-estimates-for-the-dumitriu-edelman-model}.
\begin{lemma}\label{A:res:8.2}
The map $(\nu,z)\mapsto Q_\nu(z)$ is Borel on $[0,\infty)\times\mathbb R$. Consequently the coupling in Lemma~\ref{A:res:A.3} is measurable in $\nu$ and uses the same Gaussian
variable for every moment order.
\end{lemma}
\begin{proof}
For $\nu>0$, the chi density $p_\nu(t)$ and its mean $m_\nu$ are continuous in their parameters. Hence
\[
 F_\nu(x)=\int_0^\infty
       \mathbf1_{\{t\le x+m_\nu\}}p_\nu(t)\,\mathrm dt,
 \qquad F_0(x)=\mathbf1_{\{x\ge0\}},
\]
is jointly Borel by integration of a nonnegative Borel function. For $0<u<1$, right continuity of $F_\nu$ gives $F_\nu^{-1}(u)\le x$ if and only if $u\le F_\nu(x)$. Consequently,
for every real $a$,
\[
 \{(\nu,z):Q_\nu(z)<a\}
 =\bigcup_{q\in\mathbb Q,\ q<a}
       \{(\nu,z):F_\nu(q)\ge\Phi_{\mathrm N}(z)\}.
\]
The right side is Borel. This proves measurability of the quantile map, including $\nu=0$. If $Z$ is standard normal, then $\Phi_{\mathrm N}(Z)$ is uniform on $(0,1)$, and
\[
 \mathbb P\{Q_\nu(Z)\le x\}=F_\nu(x).
\]
Thus the same map realizes the centered chi law for every moment estimate in Lemma~\ref{A:res:A.3}.
\end{proof}

Start with a standard Brownian motion $B$ and a countable independent family $(\rho_{k,n})_{n\ge1,\,1\le k\le n}$ of standard normal variables, independent of $B$. For
$h=n^{-1/3}$, put
\[
 \varepsilon_{k,n}=h^{-1/2}\bigl(B(kh)-B((k-1)h)\bigr),\qquad
 G_{k,n}^{(0)}=-\frac{\varepsilon_{k,n}+\rho_{k,n}}{\sqrt2},\qquad
 G_{k,n}^{(1)}=-\frac{\varepsilon_{k,n}-\rho_{k,n}}{\sqrt2}.
\]
For fixed $n$, the pairs $(\varepsilon_{k,n},\rho_{k,n})$ are independent over $k$. Within each pair, for $s,t\in\mathbb
R$,
\[
 \mathbb E\exp\{\mathrm i(sG_{k,n}^{(0)}+tG_{k,n}^{(1)})\}
 =\exp\left\{-\frac{(s+t)^2+(s-t)^2}{4}\right\}
 =\mathrm e^{-s^2/2}\mathrm e^{-t^2/2}.
\]
Thus $G_{k,n}^{(0)}$ and $G_{k,n}^{(1)}$ are independent standard normals, and the pairs are independent over $k$. They satisfy the exact identity
\[
 -\frac{G_{k,n}^{(0)}+G_{k,n}^{(1)}}{\sqrt2}=\varepsilon_{k,n}.
 \eqtag{8.18}
\]
Define
\[
 a_{k,n}=\sqrt{\frac2{\beta n}}\,G_{k,n}^{(0)},\qquad 1\le k\le n,
 \eqtag{8.19}
\]
and set $\nu_{k,n}=\beta(n-k)$ and $X_{k,n}=Q_{\nu_{k,n}}(G_{k,n}^{(1)})$ for $1\le k\le n$. In particular $X_{n,n}=0$. Lemma~\ref{A:res:A.3} gives
\[
 \left\|X_{k,n}-\frac{G_{k,n}^{(1)}}{\sqrt2}\right\|_{L^p}
 \le\frac{C_{p,\beta}}{\sqrt{1+n-k}},\qquad 1\le k\le n.
 \eqtag{8.21}
\]
Here $1+\beta(n-k)\ge\min\{1,\beta\}(1+n-k)$, so the constant can be chosen uniformly for $\beta$ in each compact subinterval of $(0,\infty)$. Set
\[
 \delta_{k,n}=\frac{X_{k,n}}{\sqrt{\beta n}}\quad(1\le k\le n),
 \qquad \delta_{0,n}=0.
 \eqtag{8.23}
\]
For $k<n$ the coefficients $b_{k,n}=\overline b_{k,n}+\delta_{k,n}$ obey
\[
 b_{k,n}
 =\frac{m_{\nu_{k,n}}+Q_{\nu_{k,n}}(G_{k,n}^{(1)})}{\sqrt{\beta n}}
 \ \stackrel{\mathrm d}=\ \frac{\chi_{\beta(n-k)}}{\sqrt{\beta n}}.
\]
They are positive almost surely. The characteristic function calculation and the separate quantile transformation of each $G_{k,n}^{(1)}$ give the independent coefficient law of
Subsection~\ref{A:sec:dumitriu-edelman-decomposition} for $(a_{k,n})_{k\le n},(b_{k,n})_{k<n}$. For $\nu=\beta(n-k)>0$, the change of variables $t=\sqrt{\beta n}\,b$ in the chi
density gives
\[
 f_{k,n}(b)=\frac{2(\beta n/2)^{\nu/2}}{\Gamma(\nu/2)}
                b^{\nu-1}e^{-\beta n b^2/2}\mathbf1_{\{b>0\}}.
\]
Multiplying these densities and the $n$ normal densities verifies the full coefficient law, including its normalization:
\[
 p_{n,\beta}(a,b)=c_{n,\beta}e^{-\beta n\operatorname{Tr}T(a,b)^2/4}
       \prod_{k=1}^{n-1}b_k^{\beta(n-k)-1}\mathbf1_{\{b_k>0\}},
\]
\[
 c_{n,\beta}=\left(\frac{\beta n}{4\pi}\right)^{n/2}
       \prod_{k=1}^{n-1}
       \frac{2(\beta n/2)^{\beta(n-k)/2}}{\Gamma(\beta(n-k)/2)}.
\]
Here $\operatorname{Tr}T^2=\sum_ka_k^2+2\sum_{k<n}b_k^2$. For $n=1$ the product is empty. The auxiliary values remain $b_{0,n}=b_{n,n}=1$; $X_{n,n}=0$ specifies the centered
fluctuation at this auxiliary boundary point. These auxiliary values do not occur in the density. The density agrees with \eqref{A:eq:jacobi-law} for $V(t)=t^2/2$; throughout this
construction the off-diagonal coefficients are centered at their means $\overline b_{k,n}$. The Brownian path $B$ couples the arrays for different $n$. Independence over $k$ holds
at each fixed size. For $k<n$, the recovered diagonal and off-diagonal coefficients within a pair are dependent conditional on $B$. The product coefficient law just verified is
unconditional. Indeed, conditioning on $B$ fixes $\varepsilon_{k,n}$ and leaves only $\rho_{k,n}$ random, so
\[
 \operatorname{Cov}\!\left((G_{k,n}^{(0)},G_{k,n}^{(1)})\mid B\right)
   =\frac12\begin{pmatrix}1&-1\\-1&1\end{pmatrix}.
\]
The vector $(1,1)$ is in its kernel, since the sum of the two coordinates is fixed by $B$. Their conditional mean is $-\varepsilon_{k,n}(1,1)/\sqrt2$ and has covariance
$\frac12\left(\begin{smallmatrix}1&1\\1&1\end{smallmatrix}\right)$; the sum of these two covariance matrices is the
unconditional identity covariance.

Using the definition of $\zeta_n$ and \eqref{A:eq:8.18}, we have
\[
 \zeta_{k,n}=-\sqrt{\frac2{\beta n}}G_{k,n}^{(0)}
                    -\frac2{\sqrt{\beta n}}X_{k,n},\qquad
 \zeta_{k,n}^{B}=\frac2{\sqrt{\beta n}}\varepsilon_{k,n}
       =-\sqrt{\frac2{\beta n}}(G_{k,n}^{(0)}+G_{k,n}^{(1)}).
\]
Thus
\[
 \Delta\zeta_{k,n}:=\zeta_{k,n}-\zeta_{k,n}^{B}
 =-\frac2{\sqrt{\beta n}}
          \left(X_{k,n}-\frac{G_{k,n}^{(1)}}{\sqrt2}\right).
\]
The diagonal Gaussian terms cancel. Equation~\eqref{A:eq:8.21} implies
\[
 \|\Delta\zeta_{k,n}\|_{L^p}
 \le C_{p,\beta}\frac{h^{3/2}}{\sqrt{1+n-k}},\qquad 1\le k\le n.
 \eqtag{8.28}
\]
At the last coordinate the identity reads $\Delta\zeta_{n,n}=\sqrt{2/(\beta n)}G_{n,n}^{(1)}$; its variance is $2/(\beta n)$, so this case is included in the estimate. Write
$q_{k,n}=h^{-2}\Delta\zeta_{k,n}$. The following maximum estimates bound $|q_{k,n}|$ by a vanishing multiple of $1+x_k$ and make $\|\eta_n\|_\infty$ tend to zero. These bounds are
used to compare the quadratic forms. The weighted sums in Subsection~\ref{A:sec:weighted-errors-of-the-gaussian-coefficients} give the resolvent error.

\begin{lemma}\label{A:res:8.2A}
Almost surely, for all sufficiently large $n$,
\[
 \max_{k\le n/2}|\Delta\zeta_{k,n}|\le h^{11/4},\qquad
 \max_{k>n/2}|\Delta\zeta_{k,n}|\le h^{5/4},
 \eqtag{8.28a}
\]
\[
 \|\eta_n\|_\infty\le h^{1/4},
 \eqtag{8.28b}
\]
and consequently
\[
 |q_{k,n}|\le C h^{3/4}(1+x_k),\qquad 1\le k\le n.
 \eqtag{8.28c}
\]
\end{lemma}
\begin{proof}
Fix a finite $p\ge1$. For $k\le n/2$, \eqref{A:eq:8.28} gives $\|\Delta\zeta_{k,n}\|_{L^p}\le C_{p,\beta}h^3$. The union bound and Markov's inequality yield
\[
 \mathbb P\left(\max_{k\le n/2}|\Delta\zeta_{k,n}|>h^{11/4}\right)
 \le C_{p,\beta}n(h^3/h^{11/4})^p
 =C_{p,\beta}n^{1-p/12}.
 \eqtag{8.28d}
\]
For $k>n/2$ use $\|\Delta\zeta_{k,n}\|_{L^p}\le C_{p,\beta}h^{3/2}$ to obtain
\[
 \mathbb P\left(\max_{k>n/2}|\Delta\zeta_{k,n}|>h^{5/4}\right)
 \le C_{p,\beta}n(h^{3/2}/h^{5/4})^p
 =C_{p,\beta}n^{1-p/12}.
 \eqtag{8.28e}
\]
Finally $\eta_{k,n}=h^{1/2}X_{k,n}/\sqrt\beta$, so Lemma~\ref{A:res:A.2} gives
\[
 \mathbb P(\|\eta_n\|_\infty>h^{1/4})
 \le C_{p,\beta}n(h^{1/2}/h^{1/4})^p
 =C_{p,\beta}n^{1-p/12}.
 \eqtag{8.28f}
\]
All constants in these three probability bounds are uniform for $\beta$ in a fixed compact subinterval of $(0,\infty)$. For the fixed $\beta$ under consideration, taking $p=36$
makes the bounds summable in $n$. Borel--Cantelli on the joint space carrying $B$ and the auxiliary normals proves \eqref{A:eq:8.28a}--\eqref{A:eq:8.28b}, with the deterministic
thresholds displayed there and a random starting size. For $k\le n/2$, multiplication by $h^{-2}$ gives $|q_{k,n}|\le h^{3/4}$. For $k>n/2$ and $n\ge2$,
$x_k=(k-\tfrac12)h\ge(n-1)h/2\ge h^{-2}/4$, and hence
\[
 \frac{|q_{k,n}|}{1+x_k}\le C h^{-2}h^{5/4}h^2
 =C h^{5/4}\le C h^{3/4}.
\]
This proves \eqref{A:eq:8.28c}.
\end{proof}

\subsubsection{Weighted coefficient estimates}\label{A:sec:weighted-errors-of-the-gaussian-coefficients}
The perturbation estimates require different weighted sums for the first-order and diagonal errors. With the coefficients just constructed, define
\[
 \mathcal E_{\eta,n}^{\mathrm{DE}}
 =\sum_{k=1}^n h(1+x_k)^{-3/2}|\eta_{k,n}|^2,\qquad
 \mathcal E_{\mathrm{scal},n}^{\mathrm{DE}}
 =h^{-3}\sum_{k=1}^n(1+x_k)^{-7/2}|\Delta\zeta_{k,n}|^2.
\]

Once the common form bounds are established, Proposition~\ref{A:res:3.17} bounds the first-order contribution by $\mathcal E_{\eta,n}^{\mathrm{DE}}$, and
Proposition~\ref{A:res:3.16} bounds the scalar contribution by $\mathcal E_{\mathrm{scal},n}^{\mathrm{DE}}$. The first-order sum will be $O(h)$ almost surely. The scalar sum has
$L^p$ size $O(h^2)$; higher moments will make it $o(h)$ almost surely.

All expectations here are under the joint coupling law, before conditioning on $B$. For fixed $n$, the $X_{k,n}$ are independent under this law. The constants in the moment and
probability estimates below are uniform, at every matrix size, for $\beta$ in a fixed compact interval $I_\beta\Subset(0,\infty)$: Lemma~\ref{A:res:A.2} is uniform in the degrees
of freedom, and $1+\beta(n-k)\ge\min\{1,\inf I_\beta\}(1+n-k)$ in Lemma~\ref{A:res:A.3}. The proof uses the marginal moment bounds at each size and then sums their failure
probabilities.

For $x\in[(k-1)h,kh]$ and $0<h\le1$, the ratio $(1+x)/(1+x_k)$ lies between $1/2$ and $2$. Hence, for every real $a$,
\[
 h\sum_{k=1}^{M}(1+x_k)^{-a}
 \le 2^{|a|}\int_0^{Mh}(1+x)^{-a}\,\mathrm dx.
\]
The integral equals $\bigl((1+Mh)^{1-a}-1\bigr)/(1-a)$ when $a\ne1$. Thus, for any fixed $C_0>0$,
\[
 h\sum_{k\ge1}(1+x_k)^{-a}\le C_a\quad(a>1),\qquad
 h\sum_{k=1}^{M}(1+x_k)^{-a}\le C_{a,C_0}h^{-2(1-a)}
 \quad(a<1,\ Mh\le C_0h^{-2}).
\]
In particular, the second bound applies to negative $a$, when the weight grows along the interval.

\begin{proposition}\label{A:res:A.4}
For every $p\ge1$,
\[
 \|\mathcal E_{\eta,n}^{\mathrm{DE}}\|_{L^p}\le C_{p,\beta}h.
 \eqtag{A.14}
\]
There are $A_0(\beta),c(\beta)>0$ such that $\mathbb P(\mathcal E_{\eta,n}^{\mathrm{DE}}>Ah)\le2\exp(-c(\beta)A/h)$ for $A\ge A_0(\beta)$. The moment constant and the two tail
constants can be chosen uniformly for $\beta\in I_\beta$.
\end{proposition}
\begin{proof}
Since $\eta_{k,n}=h^{1/2}X_{k,n}/\sqrt\beta$, put
\[
 w_k=(1+x_k)^{-3/2},\qquad c_k=\frac{h^2}{\beta}w_k,
 \qquad \mathcal E_{\eta,n}^{\mathrm{DE}}=\sum_{k=1}^n c_kX_{k,n}^2.
\]
The moment bound in Lemma~\ref{A:res:A.2} and the integral comparison give
\[
 \mathbb E\mathcal E_{\eta,n}^{\mathrm{DE}}\le \frac C\beta h,
 \qquad \max_k c_k\le \frac{h^2}\beta,\qquad
  \sum_{k=1}^n c_k^2
 \le \frac{h^4}{\beta^2}\sum_{k=1}^n(1+x_k)^{-3}\le \frac C{\beta^2}h^3.
\]
Let $S_n=\mathcal E_{\eta,n}^{\mathrm{DE}}-\mathbb E\mathcal E_{\eta,n}^{\mathrm{DE}}$. Apply the exponential estimate at the end of the proof of Lemma~\ref{A:res:A.2} to the
independent centered variables $X_{k,n}^2-\mathbb EX_{k,n}^2$. Its constants are independent of $\beta(n-k)$, including the value zero at $k=n$. Consequently,
\[
 \log\mathbb E\exp(tS_n)
 =\sum_{k=1}^n\log\mathbb E
      \exp\{tc_k(X_{k,n}^2-\mathbb EX_{k,n}^2)\}
 \le C\beta^{-2}t^2h^3,\qquad |t|\le c\beta h^{-2}.
\]
For $u>0$, choose $t=\min\{\beta^2u/(4Ch^3),c\beta h^{-2}/2\}$. Then $C\beta^{-2}th^3\le u/4$. Exponential Markov at $t$ and $-t$ gives
\[
 \mathbb P(|S_n|>u)
 \le2\exp(-tu+C\beta^{-2}t^2h^3)
 \le2\exp\{-c\min(\beta^2u^2/h^3,\beta u/h^2)\}.
\]
Choose $A_0(\beta)=(2C+2)/\beta$, with $C$ from the mean bound. For $A\ge A_0(\beta)$, the event $\{\mathcal E_{\eta,n}^{\mathrm{DE}}>Ah\}$ implies $S_n>Ah/2$, and $\beta A\ge2$.
Hence
\[
 \min\{\beta^2(Ah/2)^2/h^3,\beta(Ah/2)/h^2\}
 =\frac{\beta A}{2h}.
\]
This proves the tail bound with $c(\beta)=c\beta/2$. On $I_\beta$, use the largest threshold and the smallest exponential constant. Finally, for $p\ge1$ and $h\le1$,
\[
 \mathbb E\left[(\mathcal E_{\eta,n}^{\mathrm{DE}}/h)^p\right] =p\int_0^\infty A^{p-1} \mathbb
 P(\mathcal E_{\eta,n}^{\mathrm{DE}}>Ah)\,\mathrm dA \le A_0(\beta)^p+2p\int_{A_0(\beta)}^\infty A^{p-1}\mathrm
 e^{-c(\beta) A/h}\,\mathrm dA \le C_{p,\beta}.
\]
Taking the $p$-th root proves \eqref{A:eq:A.14}.
\end{proof}

\begin{proposition}\label{A:res:A.5}
For every $p\ge1$,
\[
 \|\mathcal E_{\mathrm{scal},n}^{\mathrm{DE}}\|_{L^p}
 \le C_{p,\beta}h^2.
 \eqtag{A.18}
\]
For each $\epsilon>0$, almost surely, for all sufficiently large $n$,
\[
 \mathcal E_{\mathrm{scal},n}^{\mathrm{DE}}
 \le C_{\omega,\beta,\epsilon}h^{2-\epsilon}.
 \eqtag{A.19}
\]
The constant in \eqref{A:eq:A.18} is uniform for $\beta\in I_\beta$.
\end{proposition}
\begin{proof}
The identity for $\Delta\zeta_{k,n}$ in Subsection~\ref{A:sec:coupling-the-dumitriu-edelman-fluctuations-to-brownian-increments}, followed by the quantile estimate, gives
\eqref{A:eq:8.28} at every coordinate. In particular,
\[
 \bigl\||\Delta\zeta_{k,n}|^2\bigr\|_{L^p}
 =\|\Delta\zeta_{k,n}\|_{L^{2p}}^2
 \le C_{p,\beta}\frac{h^3}{1+n-k}.
\]
Minkowski's inequality therefore gives
\[
 \|\mathcal E_{\mathrm{scal},n}^{\mathrm{DE}}\|_{L^p} \le h^{-3}\sum_{k=1}^n(1+x_k)^{-7/2}
 \bigl\||\Delta\zeta_{k,n}|^2\bigr\|_{L^p} \le C_{p,\beta}\sum_{k=1}^n
 \frac{(1+x_k)^{-7/2}}{1+n-k}.
\]
With $J=\lfloor n/2\rfloor$, the first part of this sum satisfies
\[
 \sum_{k=1}^J\frac{(1+x_k)^{-7/2}}{1+n-k}
 \le\frac2n\sum_{k=1}^J(1+x_k)^{-7/2}
 \le C n^{-1}h^{-1}=Ch^2.
\]
For $k>J$, one has $1+x_k\ge c h^{-2}$. Setting $m=n-k$ gives
\[
 \sum_{k=J+1}^n\frac{(1+x_k)^{-7/2}}{1+n-k}
 \le Ch^7\sum_{m=0}^{n-J-1}\frac1{1+m}
 \le Ch^7(1+\log n)\le Ch^2.
\]
Here $h^5(1+\log n)=n^{-5/3}(1+\log n)$ is bounded for $n\ge1$. At $k=n$, the exact identity $\Delta\zeta_{n,n}=\sqrt{2/\beta}\,h^{3/2}G_{n,n}^{(1)}$ gives
\[
 \begin{aligned}
 h^{-3}(1+x_n)^{-7/2}|\Delta\zeta_{n,n}|^2
 &=\frac2\beta(1+x_n)^{-7/2}|G_{n,n}^{(1)}|^2,\\
 \left\|\frac2\beta(1+x_n)^{-7/2}|G_{n,n}^{(1)}|^2\right\|_{L^p}
 &\le C_{p,\beta}h^7.
 \end{aligned}
\]
Thus the difference in variances at the last coordinate is included in the same sum. The two parts prove \eqref{A:eq:A.18}; all constants are uniform on $I_\beta$ by the quantile
estimate stated above. For later use, for every $\epsilon>0$ and $p>3/\epsilon$, Markov's inequality also gives
\[
 \mathbb P\left(
 h^{-3}(1+x_n)^{-7/2}|\Delta\zeta_{n,n}|^2>h^{7-\epsilon}\right)
 \le C_{p,\beta}h^{p\epsilon}
 =C_{p,\beta}n^{-p\epsilon/3}.
\]
Summability gives an event of probability one and a finite $n_0(\omega,\epsilon)$ such that the contribution from the last coordinate is at most $h^{7-\epsilon}$ for
$n\ge n_0(\omega,\epsilon)$. For fixed $\epsilon>0$, choose $p>3/\epsilon$. Markov's inequality yields
\[
 \mathbb P(\mathcal E_{\mathrm{scal},n}^{\mathrm{DE}}>h^{2-\epsilon})
 \le h^{-p(2-\epsilon)}
          \mathbb E[(\mathcal E_{\mathrm{scal},n}^{\mathrm{DE}})^p]
 \le C_{p,\beta}h^{p\epsilon}=C_{p,\beta}n^{-p\epsilon/3}.
\]
Summability and Borel--Cantelli prove \eqref{A:eq:A.19}. Intersect both families of events over positive rational
$\epsilon$.
For any real $\epsilon>0$, choose rational $0<\epsilon'<\epsilon$ and use $h^{2-\epsilon'}\le h^{2-\epsilon}$ and $h^{7-\epsilon'}\le h^{7-\epsilon}$ for $h\le1$. This gives one
event for all $\epsilon>0$, with $\beta$ fixed.
\end{proof}

Fix $A_\eta$ above the threshold in Proposition~\ref{A:res:A.4}, and let $E_{n,\mathrm{DE}}$ be the event on which
\[
 \max_{k\le n/2}|\Delta\zeta_{k,n}|\le h^{11/4},\quad \max_{k>n/2}|\Delta\zeta_{k,n}|\le
 h^{5/4},\quad \|\eta_n\|_\infty\le h^{1/4},\quad \mathcal E_{\eta,n}^{\mathrm{DE}}\le A_\eta
 h,\quad \mathcal E_{\mathrm{scal},n}^{\mathrm{DE}}\le h^{3/2}.
\]
An empty maximum is taken to be zero. The maximum estimates of Lemma~\ref{A:res:8.2A} and Propositions~\ref{A:res:A.4}--\ref{A:res:A.5} give
\[
 \mathbb P(E_{n,\mathrm{DE}}^c) \le C_\beta n^{-2}+2\mathrm e^{-c(\beta)/h} +h^{-27}\mathbb
 E[(\mathcal E_{\mathrm{scal},n}^{\mathrm{DE}})^{18}] \le C_\beta n^{-2}+2\mathrm
 e^{-c(\beta)/h}+C_\beta h^9.
\]
Here we used moment order $36$ for the maxima and order $18$ for the scalar sum. In the latter case,
$h^{-27}\mathbb E[(\mathcal E_{\mathrm{scal},n}^{\mathrm{DE}})^{18}] \le C_\beta h^{-27}h^{36}=C_\beta h^9$. All three terms are summable because $h=n^{-1/3}$. Set
\[
 \Omega_{\mathrm{DE}}
 =\{B\in\Omega_B\}\cap\bigcup_{m\ge1}\bigcap_{n\ge m}E_{n,\mathrm{DE}}.
\]
This is a measurable event of probability one on the joint space carrying $B$ and the auxiliary normals. The comparisons involving $\zeta_n$ and $\eta_n$ below are made on
$\Omega_{\mathrm{DE}}$. These events are used only to estimate errors; we do not condition or resample the matrices on them. Their marginal laws therefore remain the exact
Dumitriu--Edelman laws.

\subsubsection{Expected tridiagonal coefficients}\label{A:sec:the-mean-dumitriu-edelman-coefficients}

For $B\in\Omega_B$, define
\[
 \overline H_{n,\beta}^{B}=\overline H_{n,\beta}+h^{-2}M_{\zeta_n^B}.
\]
We compare this operator with $\mathcal B_{n,h,\beta}^B$, which has the same diagonal potential generated by the Brownian increments. The coefficients $\overline b_{k,n}$ are close
to $\sqrt{1-k/n}$, whose variation over the whole matrix is of order one. On $I_0=\{1,\ldots,J_n\}$, where $J_n=\lfloor n/2\rfloor$, we compare the coefficient of $|(D_hu)_k|^2$
with $1$ and the deterministic potential with $x_k$ by using Section~\ref{A:sec:weighted-resolvent-estimates-and-comparison-principles}. On the remaining coordinates,
$x_k\asymp h^{-2}$ supplies a lower bound for the shifted operator. A subscript $J_n$ always denotes a principal block of the original $n$-dimensional operator, with the original
coefficients, coordinates and spacing
$h=n^{-1/3}$.

Fix a compact interval $I_\beta=[\beta_-,\beta_+]\subset(0,\infty)$ when uniformity in $\beta$ is needed. In this subsubsection the deterministic estimates hold for
$n\ge N(I_\beta)$, where $N(I_\beta)$ is deterministic. For $B\in\Omega_B$, the same threshold is used in all estimates involving Brownian increments; their constants may depend on
$B$, $I_\beta$ and $K$.

\begin{lemma}\label{A:res:7.1}
For every $\nu\ge1$, one has $c\sqrt\nu\le m_\nu\le\sqrt\nu$, with an absolute $c>0$.
\end{lemma}
\begin{proof}
Jensen's inequality gives $m_\nu\le(\mathbb E\chi_\nu^2)^{1/2}=\sqrt\nu$. For $Y=\chi_\nu^2$, the Gamma moments are $\mathbb EY=\nu$ and $\mathbb EY^2=\nu(\nu+2)$. Writing
$A=\{Y\ge\nu/2\}$ and using Cauchy--Schwarz,
\[
 \frac\nu2\le\mathbb E[Y\mathbf1_A]
 \le\sqrt{\nu(\nu+2)\mathbb P(A)},\qquad
 \mathbb P(A)\ge\frac{\nu}{4(\nu+2)}\ge\frac1{12}.
\]
Therefore $m_\nu\ge\sqrt{\nu/2}\,\mathbb P(A) \ge\sqrt\nu/(12\sqrt2)$.
\end{proof}

For $1\le k\le J_n$, put $s_k=k/n$. Then $\nu_k=\beta(n-k)=\beta n(1-s_k)\ge\beta n/2$. Lemma~\ref{A:res:A.1}, with a remainder uniform on this interval, gives
\[
\begin{aligned}
\overline b_{k,n}
&=
\sqrt{1-s_k}
\left(
1
-
\frac{1}{4\beta n(1-s_k)}
+
\frac{1}{32\beta^2n^2(1-s_k)^2}
+
O\!\left(n^{-3}(1-s_k)^{-3}\right)
\right)
\\
&=
\sqrt{1-s_k}
-
\frac{1}{4\beta n}(1-s_k)^{-1/2}
+
O\!\left(n^{-2}(1-s_k)^{-3/2}\right).
\end{aligned}
\eqtag{7.11}
\]
Define $b^0_{k,n}=\sqrt{1-k/n}$ for $0\le k\le J_n$, and $c_{k,n}=\overline b_{k,n}-b^0_{k,n}$ for $1\le k\le J_n$. Equation~\eqref{A:eq:7.11} gives
\[
 c_{k,n}=-\frac1{4\beta n}(1-k/n)^{-1/2}
       +O_\beta\bigl(n^{-2}(1-k/n)^{-3/2}\bigr).
 \eqtag{7.14}
\]
This expansion is used only on the initial block. At the auxiliary left edge, $c_{0,n}=\overline b_{0,n}-b^0_{0,n}=0$. It is not the value obtained by extending the chi mean
formula to $k=0$.

Put $\phi(s)=\sqrt{1-s}$ and $d^0_{k,n}=b^0_{k,n}-1$. On $[0,1/2]$, all derivatives of $\phi$ are bounded.
\begin{lemma}\label{A:res:7.2}
For \(0\le k\le J_n\),
\[
|d^0_{k,n}|
\le
C h^2(1+x_k).
\eqtag{7.19}
\]
For \(1\le k\le J_n\),
\[
|d^0_{k,n}-d^0_{k-1,n}|
\le
C h^3.
\eqtag{7.20}
\]
\end{lemma}

\begin{proof}
For $0\le k\le J_n$,
\[
 d^0_{k,n}=-\frac{k/n}{1+\sqrt{1-k/n}},\qquad
 |d^0_{k,n}|\le\frac kn=kh^3\le Ch^2(1+x_k).
\]
For $1\le k\le J_n$, rationalizing the difference gives
\[
 d^0_{k,n}-d^0_{k-1,n}
 =-\frac{n^{-1}}{\sqrt{1-k/n}+\sqrt{1-(k-1)/n}},
\]
whose absolute value is at most $Cn^{-1}=Ch^3$.
\end{proof}

Define the diagonal error from the deterministic coefficients
\[
r^0_{k,n}
:=
h^{-2}
\left(
2-b^0_{k-1,n}-b^0_{k,n}
\right)
-
x_k.
\eqtag{7.21}
\]

\begin{lemma}\label{A:res:7.3}
For \(1\le k\le J_n\),
\[
|r^0_{k,n}|
\le
C h^2(1+x_k)^2.
\eqtag{7.22}
\]
\end{lemma}

\begin{proof}
Set $s=(k-1/2)/n=h^2x_k$ and $\delta=1/(2n)=h^3/2$. Both $s-\delta$ and $s+\delta$ lie in $[0,1/2]$. Twice applying the fundamental theorem of calculus gives
\[
 \phi(s+\delta)+\phi(s-\delta)-2\phi(s)
 =\int_0^\delta(\delta-t)
          \bigl(\phi''(s+t)+\phi''(s-t)\bigr)\,\mathrm dt.
\]
As $\phi''(s)=-1/[4(1-s)^{3/2}]$, this integral is nonpositive and has absolute value at most $C\delta^2$. Hence
\[
 2-b^0_{k-1,n}-b^0_{k,n}=2(1-\phi(s))+O(\delta^2).
 \eqtag{7.23}
\]
Moreover, direct rationalization gives the exact identity
\[
 2(1-\phi(s))=s+\frac{s^2}{(1+\sqrt{1-s})^2}.
 \eqtag{7.24}
\]
Put $E_{s,\delta}=-\int_0^\delta(\delta-t) [\phi''(s+t)+\phi''(s-t)]\,\mathrm dt$. Then $0\le E_{s,\delta}\le C\delta^2$, and the definition of $r^0_{k,n}$ gives
\[
 r^0_{k,n}=h^{-2}\left\{
 \frac{s^2}{(1+\sqrt{1-s})^2}+E_{s,\delta}\right\}.
\]
Here $h^{-2}s=x_k$. Consequently
\[
 0\le r^0_{k,n}\le C h^{-2}(s^2+\delta^2)
 =C(h^2x_k^2+h^4)\le C h^2(1+x_k)^2.
\]
For $k=1$ the boundary point $s-\delta$ is zero; the bound on $\phi''$ still holds on the entire interval of integration. This proves \eqref{A:eq:7.22}.
\end{proof}

For the coefficients $c_{k,n}$ above, the first difference at $k=1$ uses the value $c_{0,n}=0$ and is estimated separately below.

Define
\[
 s_{k,n}^{\chi}=-h^{-2}(c_{k-1,n}+c_{k,n}),\qquad 1\le k\le J_n.
 \eqtag{7.28}
\]
The constants in the following bounds are uniform for $\beta$ in a fixed compact subinterval of $(0,\infty)$.

\begin{lemma}\label{A:res:7.4}
One has $\max_{0\le k\le J_n}|c_{k,n}|\le C_\beta h^3$ and $\max_{1\le k\le J_n}|s_{k,n}^{\chi}|\le C_\beta h$. In addition,
\[
 |c_{k,n}-c_{k-1,n}|\le C_\beta h^6\quad(2\le k\le J_n),\qquad
 c_{1,n}-c_{0,n}=-\frac{h^3}{4\beta}+O_\beta(h^6).
\]
\end{lemma}

\begin{proof}
Define \(F_\beta(s)=-(4\beta)^{-1}(1-s)^{-1/2}\), \(0\le s\le1/2\). Equation~\eqref{A:eq:7.14} gives, for \(1\le k\le
J_n\),
\[
c_{k,n}=n^{-1}F_\beta(k/n)+\varepsilon_{k,n},
\quad |\varepsilon_{k,n}|\le C_\beta n^{-2}.
\]
Both \(F_\beta\) and \(F_\beta'\) are bounded on this interval. Thus $|c_{k,n}|\le C_\beta n^{-1}=C_\beta h^3$ for $1\le k\le J_n$. The value $c_{0,n}=0$ satisfies this bound as
well. For \(2\le k\le J_n\), the mean-value theorem and the triangle inequality give
\[
|c_{k,n}-c_{k-1,n}| \le n^{-1}|F_\beta(k/n)-F_\beta((k-1)/n)|
 +|\varepsilon_{k,n}|+|\varepsilon_{k-1,n}| \le n^{-2}\|F_\beta'\|_\infty+2C_\beta n^{-2} \le
 C_\beta h^6.
\]
  At the first edge, \(F_\beta(1/n)=F_\beta(0)+O_\beta(n^{-1})\) and \(c_{0,n}=0\), so
\[
c_{1,n}-c_{0,n}
=n^{-1}F_\beta(0)+O_\beta(n^{-2})
=-\frac{h^3}{4\beta}+O_\beta(h^6).
\]
Finally, for every \(1\le k\le J_n\),
\[
|s_{k,n}^{\chi}|
\le h^{-2}(|c_{k-1,n}|+|c_{k,n}|)
\le C_\beta h.
\]
For \(\beta\) in a compact positive interval, the expansion for the ratio of Gamma functions is uniform because
\(\beta(n-k)\ge\beta_-n/2\); the bounds on \(F_\beta\), \(F_\beta'\), and the remainder are uniform as well.
\end{proof}

The mean off-diagonal coefficient tends to zero near $k=n$. To obtain a uniform lower bound for the full difference energy, we use $x_k\asymp h^{-2}$ in that region: the potential
controls the difference energy by $|u_{k+1}-u_k|^2\le2(|u_k|^2+|u_{k+1}|^2)$. We now establish this form bound uniformly for all principal intervals.

Let \(I=\{p,\ldots,q\}\subset\{1,\ldots,n\}\) be a consecutive interval. For a vector supported on \(I\), put
\(u_{p-1}=u_{q+1}=0\) and define
\[
\ell_{h,I}[u]
=
\sum_{k=p-1}^{q}|(D_hu)_k|^2
+
\sum_{k=p}^{q}(1+x_k)|u_k|^2.
\]
Write \(\overline V_{k,n} = h^{-2}\bigl(2-\overline b_{k-1,n}-\overline b_{k,n}\bigr)\) and let
\[
\overline{\mathfrak h}_{n,\beta}^{I}[u]
=
\sum_{k=p-1}^{q}\overline b_{k,n}|(D_hu)_k|^2
+
\sum_{k=p}^{q}\overline V_{k,n}|u_k|^2
\]
be the Dirichlet restriction of the deterministic form with mean off-diagonal coefficients. On the initial block
\(I_0=\{1,\ldots,J_n\}\), also set \(V^0_{k,n} = h^{-2}\bigl(2-b^0_{k-1,n}-b^0_{k,n}\bigr)\) and
\[
\mathfrak h_{n}^{0,I_0}[u]
=
\sum_{k=0}^{J_n}b^0_{k,n}|(D_hu)_k|^2
+
\sum_{k=1}^{J_n}V^0_{k,n}|u_k|^2.
\]

The Brownian diagonal multiplier is
\[
 \omega_{k,n}^{B}=h^{-2}\zeta_{k,n}^{B}
       =\frac2{\sqrt\beta}\frac{B(kh)-B((k-1)h)}h.
 \eqtag{7.29f}
\]
The weights on $I$ use the original coordinates $x_k$.

\begin{proposition}\label{A:res:7.4A}
For $n\ge N(I_\beta)$ and $\beta\in I_\beta$, uniformly in every consecutive interval $I$,
\[
 c(\beta)\ell_{h,I}[u]
 \le\overline{\mathfrak h}_{n,\beta}^{I}[u]+\kappa_\beta\|u\|^2
 \le C_\beta\ell_{h,I}[u].
 \eqtag{7.29e}
\]
The constants may be chosen uniformly over $\beta\in I_\beta$. The same form bounds hold for $\mathfrak h_n^{0,I_0}$ and its principal restrictions.

For every $B\in\Omega_B$, let $\mathfrak h_I^B$ denote any of these quadratic forms after addition of $\sum_{k\in I}\omega_{k,n}^{B}|u_k|^2$. There is one shift $\kappa_B$, valid
for all these forms, such that
\[
 \begin{gathered}
 1\le\kappa_B\le C_{I_\beta}(1+\mathcal O_{1/4}(B))^4,\\
 c_*\ell_{h,I}[u]\le\mathfrak h_I^B[u]+\kappa_B\|u\|^2
       \le C_{I_\beta}(1+\kappa_B)\ell_{h,I}[u].
 \end{gathered}
\]
Here $c_*>0$ depends only on $I_\beta$, and the threshold $N(I_\beta)$ is independent of $B$.
\end{proposition}

\begin{proof}
Jensen's inequality and \(\mathbb E\chi_\nu^2=\nu\) give, for \(1\le k\le n-1\),
\[
0\le \overline b_{k,n}
=\frac{\mathbb E\chi_{\beta(n-k)}}{\sqrt{\beta n}}
\le
\sqrt{1-\frac{k}{n}}
\le1.
\eqtag{7.29g}
\]
The auxiliary boundary values \(\overline b_{0,n}=\overline b_{n,n}=1\) also lie in this interval. Consequently
\(\overline V_{k,n}\ge0\). For \(1\le k\le n-1\), the identity and bound \(1-\sqrt{1-t}=\frac{t}{1+\sqrt{1-t}}\ge\frac
t2, \allowbreak\, 0\le t\le1\), and \eqref{A:eq:7.29g} imply
\[
\overline V_{k,n} = h^{-2}\bigl[(1-\overline b_{k-1,n}) +(1-\overline b_{k,n})\bigr] \ge h^{-2}\left(\frac{k-1}{2n}+\frac{k}{2n}\right) = \left(k-\frac12\right)h = x_k.
\eqtag{7.29h}
\]
At the last vertex,
\[
\overline V_{n,n}
=
h^{-2}(1-\overline b_{n-1,n})
\ge
h^{-2}(1-n^{-1/2})
\ge
\frac13x_n
\eqtag{7.29i}
\]
for all sufficiently large \(n\). Hence, after decreasing a deterministic constant,
\[
\overline V_{k,n}\ge c x_k,
\quad 1\le k\le n.
\eqtag{7.29j}
\]

For the upper bound, uniformly for \(0\le j\le J_n\), with the convention at \(j=0\), the expansion \eqref{A:eq:7.11} gives
\[
\left|
\overline b_{j,n}-\sqrt{1-\frac jn}
\right|
\le
\frac{C_\beta}{n}.
\eqtag{7.29k}
\]
Since \(1-\sqrt{1-t}\le t\), for \(k\le J_n\) we obtain
\[
\overline V_{k,n}
\le
C_\beta h^{-2}\left(\frac{k}{n}+\frac1n\right)
\le
C_\beta(1+x_k).
\eqtag{7.29l}
\]
For $k>J_n$, the estimate \(0\le\overline V_{k,n}\le2h^{-2}\) is sufficient, because \(x_k\ge c h^{-2}\). Therefore
\[
0\le\overline V_{k,n}\le C_\beta(1+x_k)
\quad(1\le k\le n).
\eqtag{7.29m}
\]

For $1\le k\le J_n$, choosing $N(I_\beta)\ge2/\beta_-$ ensures $\beta(n-k)\ge1$. Lemma~\ref{A:res:7.1} therefore applies throughout this block; at $k=0$ use $\overline b_{0,n}=1$.
Thus, for $0\le k\le J_n$,
\[
\overline b_{k,n}
\ge
c\sqrt{1-\frac{k}{n}}
\ge
c(\beta)>0.
\eqtag{7.29n}
\]
Let \(u\) be supported on \(I\). Split the derivative edges in \(\{p-1,\ldots,q\}\) according to whether \(k\le J_n\) or
\(k>J_n\). Equation \eqref{A:eq:7.29n} gives
\[
\sum_{\substack{p-1\le k\le q\\ k\le J_n}}
|(D_hu)_k|^2
\le
C_\beta
\sum_{k=p-1}^{q}\overline b_{k,n}|(D_hu)_k|^2.
\eqtag{7.29o}
\]
For the remaining edges, each nonzero boundary point has index at least \(J_n+1\), so \(x_j\ge c h^{-2}\). Since every vertex occurs in at most two edges,
\[
\sum_{\substack{p-1\le k\le q\\ k>J_n}} |(D_hu)_k|^2 \le 2h^{-2} \sum_{\substack{p-1\le k\le q\\ k>J_n}} \bigl(|u_k|^2+|u_{k+1}|^2\bigr) \le C\sum_{j=p}^{q}x_j|u_j|^2 \le C\sum_{j=p}^{q}\overline V_{j,n}|u_j|^2,
\eqtag{7.29p}
\]
where the last inequality uses \eqref{A:eq:7.29j}. Equations \eqref{A:eq:7.29j}, \eqref{A:eq:7.29o}, and \eqref{A:eq:7.29p} show that
\(\ell_{h,I}[u] \le C_\beta\bigl( \overline{\mathfrak h}_{n,\beta}^{I}[u]+\|u\|^2 \bigr)\). After increasing \(\kappa_\beta\), this gives the lower estimate in \eqref{A:eq:7.29e}.
For the upper estimate, $0\le\overline b_{k,n}\le1$ and \eqref{A:eq:7.29m} give
\[
 \overline{\mathfrak h}_{n,\beta}^{I}[u]+\kappa_\beta\|u\|^2
 \le\sum_{k=p-1}^{q}|D_hu_k|^2
      +(C_\beta+\kappa_\beta)\sum_{k=p}^{q}(1+x_k)|u_k|^2
 \le C_\beta'\ell_{h,I}[u].
\]

For the quadratic form with coefficients $b^0_{k,n}=\sqrt{1-k/n}$ on \(I_0\), we have $b^0_{k,n}\ge2^{-1/2}$ for $0\le k\le J_n$, and the calculation in
\eqref{A:eq:7.29h}--\eqref{A:eq:7.29m} gives \(c x_k\le V^0_{k,n}\le C(1+x_k), \allowbreak\, 1\le k\le J_n.\) The positive lower bound on the coefficient of $|D_hu_j|^2$ therefore
controls $\|D_hu\|^2$ directly, proving the upper and lower form bounds.

For the Brownian term, put $e_h(x)=(2/\sqrt\beta)B(h\lfloor x/h\rfloor)$, and let $f=I_hu$ for the zero extension of a vector supported on $I$. The function $f$ vanishes at
$(p-1)h$ and $(q+1)h$. Integration on the intervals where $e_h$ is constant yields
\[
 -\int_0^\infty e_h(x)(|f(x)|^2)'\,\mathrm dx =-\frac2{\sqrt\beta h}\sum_{k=0}^n
 B(kh)(|u_{k+1}|^2-|u_k|^2) =\sum_{k\in I}\omega_{k,n}^{B}|u_k|^2.
\]
The oscillation estimate for $e_h$ in Subsection~\ref{A:sec:galerkin-approximation-and-the-resolvent-error} gives
\[
 \operatorname{osc}_{[m,m+1]}e_h+
 \left|\int_{m+1}^{m+2}e_h-\int_m^{m+1}e_h\right|
 \le C_\beta\mathcal O_{1/4}(B)(1+m)^{1/4}.
\]
This bound is uniform in $h\le1$ and $m\ge0$. Lemma~\ref{A:res:3.1}, the interpolation equivalence \eqref{A:eq:3.7}, and zero extension now give, for $0<\varepsilon\le1$,
\[
 \left|\sum_{k\in I}\omega_{k,n}^{B}|u_k|^2\right|
 \le\varepsilon\ell_{h,I}[u]
       +C_{\varepsilon,\beta}(1+\mathcal O_{1/4}(B))^4\|u\|^2.
 \eqtag{7.29r}
\]
There is no boundary term, since $f$ vanishes at both ends of its support. Choose common deterministic constants $c_0,C_0,\kappa_0$, with $\kappa_0\ge1$, for the forms already
considered, uniformly over $\beta\in I_\beta$, so that
$c_0\ell_{h,I}[u]\le\mathfrak h_I[u]+\kappa_0\|u\|^2\le C_0\ell_{h,I}[u]$.
Take $\varepsilon=c_0/2$ in \eqref{A:eq:7.29r}, decreasing $c_0$ first if necessary so that $c_0\le1$, and put
\[
 D_B=C_{c_0/2,I_\beta}(1+\mathcal O_{1/4}(B))^4,
 \qquad \kappa_B=1+\kappa_0+D_B.
\]
The same $D_B$ bounds the added quadratic form for every interval and every $\beta\in I_\beta$. Consequently,
\[
 \begin{aligned}
 \mathfrak h_I^B[u]+\kappa_B\|u\|^2
   &\ge\tfrac12c_0\ell_{h,I}[u]+\|u\|^2,\\
 \mathfrak h_I^B[u]+\kappa_B\|u\|^2
   &\le(C_0+c_0/2)\ell_{h,I}[u]+(1+2D_B)\|u\|^2\\
   &\le C_{I_\beta}(1+\kappa_B)\ell_{h,I}[u].
 \end{aligned}
\]
Thus $c_*=c_0/2$ is independent of $B$. All estimates used above hold once $n\ge N(I_\beta)$, and none requires an additional event or a size threshold depending on the path. This
proves the proposition.
\end{proof}

In \eqref{A:eq:3.8}, the coefficients of the squared differences for the operators just considered are $\overline b_{k,n}$ or $b^0_{k,n}$. Both lie in $[0,1]$, and there is no
additional first-order coefficient. The finite difference operator driven by the Brownian increments has coefficient $1$ and satisfies the same form bounds by
Subsection~\ref{A:sec:galerkin-approximation-and-the-resolvent-error}. Increase the common $\kappa_B$ by a deterministic multiple if needed to include that operator.

We check that the estimates remain uniform for principal restrictions. If $E_I$ is extension by zero from $\ell^2(I)$ to $\mathbb C^n$, then
\[
 \ell_{h,I}[u]=\ell_h[E_Iu].
\]
The min--max principle and Lemma~\ref{A:res:3.4} imply that the $r$th eigenvalue $\mu_{r,I}$ of this restricted positive form satisfies $\mu_{r,I}\ge\mu_{r,h}\ge cr^{2/3}$ for
$1\le r\le |I|$. Thus the proof of Proposition~\ref{A:res:3.5} gives
\[
 \sup_{z\in K}\|R_I(z)\|_{\mathfrak S_2}^2\le C_{\omega,\beta,K},
 \qquad
 \sup_{z\in K}\|D_hR_I(z)\|_{\mathfrak S_2}^2
       \le C_{\omega,\beta,K}|I|^{1/3}\le C_{\omega,\beta,K}h^{-1}.
\]
Here $R_I$ denotes the resolvent of any of these restrictions, and $D_h$ includes both boundary edges of $I$.

For the kernel estimates, the point evaluation bound of Lemma~\ref{A:res:3.6} already applies to $E_Iu$ with the original $x_j$. In the exponential form estimates, extend the
weight constantly across each boundary edge of $I$. The edge identities in Lemmas~\ref{A:res:3.8} and~\ref{A:res:3.12} then sum only over edges inside $I$, with the same constants,
since $h\sqrt{1+x_j}\le C$ throughout $\{1,\ldots,n\}$. All sums in Lemma~\ref{A:res:3.10} have nonnegative terms and only decrease when restricted to $I$. The proofs of
Propositions~\ref{A:res:3.13}--\ref{A:res:3.16} consequently apply on $I$ and give, in particular,
\[
 \sup_{z\in K}\|R_I(z)e_j\|^2
       \le C_{\omega,\beta,K}h(1+x_j)^{-3/2},
 \qquad j\in I.
\]
The constants are independent of the boundary points of $I$.

Let $J=J_n=\lfloor n/2\rfloor$ and work on $\ell^2(\{1,\ldots,J\})$, with $u_0=u_{J+1}=0$. Let $\overline H_{J,\beta}^{B}$ and $\mathcal B_{J,h,\beta}^{B}$ be the corresponding
principal restrictions, retaining the original coefficients in the boundary terms at $J$. Since $\overline b_{k,n}-1=d^0_{k,n}+c_{k,n}$, their exact difference is
\[
 \begin{aligned}
 \overline H_{J,\beta}^{B}-\mathcal B_{J,h,\beta}^{B}
 &=\Delta_0+\Delta_\chi,\\
 \Delta_0&=D_h^*M_{d^0}D_h+M_{r^0},\qquad
 \Delta_\chi=D_h^*M_cD_h+M_{s^\chi}.
 \end{aligned}
\]
Here $r^0$ and $s^\chi$ are defined in \eqref{A:eq:7.21} and \eqref{A:eq:7.28}. The diagonal entries of $\Delta_0$ and $\Delta_\chi$ are $-x_k$ and $0$, respectively, because
\[
 \begin{aligned}
 h^{-2}(d^0_{k-1,n}+d^0_{k,n})+r^0_{k,n}&=-x_k,\\
 h^{-2}(c_{k-1,n}+c_{k,n})+s^\chi_{k,n}&=0.
 \end{aligned}
\]
Their off-diagonal entries sum to $-h^{-2}(\overline b_{k,n}-1)$, and the Brownian diagonal potential cancels. These identities include $k=1,J$ with the stated boundary
coefficients. The three operators in this comparison satisfy the common form and coefficient bounds verified above.

Let \(R_0(z) = (\mathcal B_{J,h,\beta}^{B}-z)^{-1},\) and let \(R_1(z) = (\mathcal B_{J,h,\beta}^{B}+\Delta_0-z)^{-1}.\) We prove
\[
\sup_{z\in K}
\|R_1(z)-R_0(z)\|_{\mathfrak S_2}^{2}
\le
C_{\omega,\beta,K}h.
\eqtag{7.33}
\]
By Lemmas~\ref{A:res:7.2} and \ref{A:res:7.3}, \(|d^0_{k,n}|\le Ch^2(1+x_k),\) \(|d^0_{k,n}-d^0_{k-1,n}|\le Ch^3,\) and
\(|r^0_{k,n}|\le Ch^2(1+x_k)^2.\) Applying Lemma~\ref{A:res:product-rule} on the edges $0,\ldots,J$ gives
\[
D_h^*M_{d^0}D_h
=
M_aT_h
-
M_bJ_-^*D_h,
\eqtag{7.34}
\]
where \(T_h=D_h^*D_h,\) \(a_k=d^0_{k-1,n},\) and \(b_k=h^{-1}(d^0_{k,n}-d^0_{k-1,n}).\) By Lemma~\ref{A:res:7.2},
\(|a_k|\le Ch^2(1+x_k),\) and
\[
|b_k|\le Ch^2.
\eqtag{7.36}
\]
The finite difference operator driven by the Brownian increments is \(\mathcal B_{J,h,\beta}^{B} = T_h+X_h+h^{-2}M_{\zeta^B}.\) Let \(V_k^{B} = x_k+h^{-2}\zeta_{k,n}^B.\) Then
\[
T_hR_0(z)
=
I+zR_0(z)-M_{V^B}R_0(z).
\eqtag{7.37}
\]
The resolvent identity gives
\[
R_1-R_0
=
-R_1(D_h^*M_{d^0}D_h+M_{r^0})R_0.
\eqtag{7.38}
\]
Using \eqref{A:eq:7.34} and \eqref{A:eq:7.37},
\[
R_1D_h^*M_{d^0}D_hR_0 ={} R_1M_a+zR_1M_aR_0-R_1M_{a x}R_0
 -R_1M_{a\omega^B}R_0-R_1M_bJ_-^*D_hR_0,
\eqtag{7.39}
\]
where $\omega_k^B=(2/\sqrt\beta)\,[B(kh)-B((k-1)h)]/h$. The deterministic partition comparison in Subsection~\ref{A:sec:galerkin-approximation-and-the-resolvent-error} gives, for
every $B\in\Omega_B$ and every $0<h\le1$,
\[
 \sum_{k=1}^{J}(1+x_k)^{-3/2}|B(kh)-B((k-1)h)|^2
 \le C Q(B).
 \eqtag{7.40}
\]
Its weights $(1+kh)^{-3/2}$ and the weights here are comparable, since $x_k=(k-\tfrac12)h$; restricting the sum to $k\le J$ only decreases it. We use this bound for the Brownian
term in \eqref{A:eq:7.39}.

Proposition~\ref{A:res:3.15}, applied to the adjoint of $R_1M_a$ in \eqref{A:eq:7.39}, gives the weight $(1+x_k)^{-3/2}$. The growth of the interval is included in the following
sum. Since $|a_k|\le Ch^2(1+x_k)$,
\[
h\sum_{k\le J}|a_k|^2(1+x_k)^{-3/2} \le Ch\sum_{k\le J}h^4(1+x_k)^2(1+x_k)^{-3/2} = Ch^5\sum_{k\le J}(1+x_k)^{1/2}.
\]
Since \(x_J\asymp h^{-2}\),
\[
h\sum_{k\le J}(1+x_k)^{1/2}
\le
C\int_0^{Ch^{-2}}(1+x)^{1/2}\,\mathrm{d}x
\le
Ch^{-3}.
\]
Dividing the preceding estimate by \(h\) gives \(\sum_{k\le J}(1+x_k)^{1/2}\le Ch^{-4}\), and therefore
\[
h\sum_{k\le J}|a_k|^2(1+x_k)^{-3/2}
\le
Ch.
\eqtag{7.41}
\]
By the weighted estimate for a multiplier and a resolvent,
\[
\|R_1M_a\|_{\mathfrak S_2}^{2}\le C_{\omega,K}h.
\eqtag{7.42}
\]
The same estimate, with one additional bounded resolvent factor, gives
\[
\|R_1M_aR_0\|_{\mathfrak S_2}^{2}\le C_{\omega,K}h.
\eqtag{7.43}
\]
For the \(a x\) term, \(|a_kx_k|\le Ch^2(1+x_k)^2.\) The diagonal estimate \eqref{A:eq:3.49} gives
\[
\|R_1M_{a x}R_0\|_{\mathfrak S_2}^{2} \le C_{\omega,K} h\sum_{k\le J} h^4(1+x_k)^4(1+x_k)^{-7/2} = C_{\omega,K} h^5\sum_{k\le J}(1+x_k)^{1/2} \le C_{\omega,K}h.
\eqtag{7.44}
\]
For the Brownian part, the diagonal perturbation estimate gives
\[
 \begin{aligned}
 \|R_1M_{a\omega^B}R_0\|_{\mathfrak S_2}^2
 &\le C_{\omega,\beta,K}h^5
            \sum_{k\le J}(1+x_k)^{-3/2}|\omega_k^B|^2\\
 &\le C_{\omega,\beta,K}h^3
            \sum_{k\le J}(1+x_k)^{-3/2}|B(kh)-B((k-1)h)|^2\\
 &\le C_{\omega,\beta,K}Q(B)h^3=o_\omega(h).
 \end{aligned}
 \eqtag{7.45}
\]
The last inequality is \eqref{A:eq:7.40}. Since $Q(B)<\infty$ is fixed, $h^3Q(B)=o(h)$ as $n\to\infty$.

For the first-order remainder, \eqref{A:eq:7.36} gives \(|b_k|\le Ch^2\). Equation~\eqref{A:eq:restricted-multiplier} identifies $\|R_1M_bJ_-^*\|_{\mathfrak S_2}$ with
$\|R_1M_b\|_{\mathfrak S_2}$. The sequence $b_k=h^{-1}(d^0_{k,n}-d^0_{k-1,n})$ is real, so $(R_1(z)M_b)^*=M_bR_1(\overline z)$. Applying Proposition~\ref{A:res:3.15} to this
adjoint gives
\[
\|R_1M_bJ_-^*\|_{\mathfrak S_2}^{2}
\le
C_{\omega,K}
h\sum_k h^4(1+x_k)^{-3/2}
\le
C_{\omega,K}h^4.
\eqtag{7.46}
\]
Using \(\|D_hR_0\|_{\mathrm{op}}\le C_{\omega,K},\) we get
\[
\|R_1M_bJ_-^*D_hR_0\|_{\mathfrak S_2}^{2}
\le
C_{\omega,K}h^4.
\eqtag{7.47}
\]
Finally, for the diagonal remainder \(r^0\), \(|r^0_k|\le Ch^2(1+x_k)^2.\) Thus
\[
\|R_1M_{r^0}R_0\|_{\mathfrak S_2}^{2} \le C_{\omega,K} h\sum_k h^4(1+x_k)^4(1+x_k)^{-7/2} \le C_{\omega,K}h.
\eqtag{7.48}
\]
The right side of \eqref{A:eq:7.38}, after \eqref{A:eq:7.39}, has six terms. Their squared Hilbert--Schmidt norms are bounded above by $C_{\omega,\beta,K}(h,h,h,h^3Q(B),h^4,h)$,
respectively. Using $\|\sum_{j=1}^6T_j\|_{\mathfrak S_2}^2 \le6\sum_{j=1}^6\|T_j\|_{\mathfrak S_2}^2$ proves \eqref{A:eq:7.33}.

Set $R_2(z)=(\overline H_{J,\beta}^{B}-z)^{-1}$. Then
\[
R_1(z)=(\mathcal B_{J,h,\beta}^{B}+\Delta_0-z)^{-1},
\quad
\overline H_{J,\beta}^{B}
=\mathcal B_{J,h,\beta}^{B}+\Delta_0+\Delta_\chi,
\]
where \(\Delta_\chi=D_h^*M_cD_h+M_{s^\chi}\). The exact resolvent identity is therefore
\[
R_2-R_1=-R_2D_h^*M_cD_hR_1-R_2M_{s^\chi}R_1.
\]
All off-diagonal coefficients in this identity are indexed by \(0\le k\le J\). In particular, \(c_{0,n}=0\), whereas the coefficient \(c_{J,n}\) at the cut edge retains its
original value. Lemma~\ref{A:res:7.4} gives
\(\|c\|_{\ell^\infty(\{0,\ldots,J\})}\le C_\beta h^3\) and \(\|s^\chi\|_{\ell^\infty(\{1,\ldots,J\})}\le C_\beta h\).
For the term $R_2D_h^*M_cD_hR_1$, we use the operator norm of the first differentiated resolvent and the Hilbert--Schmidt norm of the second. This uses only the bound on $c$,
including $c_{0,n}=0$.

Proposition~\ref{A:res:7.4A} has already established common form bounds for the operators with square-root and mean off-diagonal coefficients with the same Brownian potential.
Propositions~\ref{A:res:3.3} and~\ref{A:res:3.5}, applied on
\(K\cup\overline K\), consequently give
\[
\sup_{z\in K}\|D_hR_2(\overline z)\|_{\mathrm{op}}^2
\le C_{\omega,\beta,K},
\quad
\sup_{z\in K}\|D_hR_1(z)\|_{\mathfrak S_2}^2
\le C_{\omega,\beta,K}J^{1/3}
\le C_{\omega,\beta,K}h^{-1}.
\]
Since \(R_2(z)D_h^*=(D_hR_2(\overline z))^*\), the ideal property of the Hilbert--Schmidt norm gives
\[
\begin{aligned}
\|R_2D_h^*M_cD_hR_1\|_{\mathfrak S_2}^2& =\|(D_hR_2(\overline z))^*M_cD_hR_1(z)\|_{\mathfrak S_2}^2\\
&\le \|D_hR_2(\overline z)\|_{\mathrm{op}}^2 \|c\|_\infty^2\|D_hR_1(z)\|_{\mathfrak S_2}^2 \le C_{\omega,\beta,K}h^6h^{-1} =C_{\omega,\beta,K}h^5.
\end{aligned}
\eqtag{chi-kinetic}
\]
For the diagonal correction, the diagonal estimate \eqref{A:eq:3.49} gives
\[
\begin{aligned}
\|R_2M_{s^\chi}R_1\|_{\mathfrak S_2}^2
&\le C_{\omega,\beta,K}h
\sum_{k=1}^{J}|s_{k,n}^{\chi}|^2(1+x_k)^{-7/2}\\
&\le C_{\omega,\beta,K}h^2
\left(h\sum_{k=1}^{J}(1+x_k)^{-7/2}\right)
\le C_{\omega,\beta,K}h^2,
\end{aligned}
\eqtag{chi-diagonal}
\]
since the last Riemann sum is bounded by a constant independent of \(J\) and \(h\). Combining \eqref{A:eq:chi-kinetic} and \eqref{A:eq:chi-diagonal} with
\(\|A+B\|_{\mathfrak S_2}^2\le2\|A\|_{\mathfrak S_2}^2 +2\|B\|_{\mathfrak S_2}^2\), we obtain
\[
\sup_{z\in K}
\|R_2(z)-R_1(z)\|_{\mathfrak S_2}^{2}
\le C_{\omega,\beta,K}(h^5+h^2)
\le C_{\omega,\beta,K}h^2
=o_\omega(h).
\eqtag{7.57}
\]
Combining \eqref{A:eq:7.33} and \eqref{A:eq:7.57}, we obtain the interior comparison:
\[
\sup_{z\in K}
\left\|
(\overline H_{J,\beta}^{B}-z)^{-1}
-
(\mathcal B_{J,h,\beta}^{B}-z)^{-1}
\right\|_{\mathfrak S_2}^{2}
\le
C_{\omega,\beta,K}h.
\eqtag{7.58}
\]

To compare the full matrices, remove the two off-diagonal entries joining $J=\lfloor n/2\rfloor$ to $J+1$. For either full operator $H$, let $H^\oplus$ be the direct sum of its
principal restrictions to $\{1,\ldots,J\}$ and $\{J+1,\ldots,n\}$. Then
\[
 H-H^\oplus=t_J(e_Je_{J+1}^*+e_{J+1}e_J^*),\qquad
 t_J=\begin{cases}
 -h^{-2}\overline b_{J,n},&H=\overline H_{n,\beta}^B,\\
 -h^{-2},&H=\mathcal B_{n,h,\beta}^B.
 \end{cases}
\]
The perturbation has rank two, zero diagonal, and $|t_J|\le h^{-2}$. The principal blocks retain both Dirichlet boundary energies with the original coefficients.

Write $R=(H-z)^{-1}$ and $R^\oplus=(H^\oplus-z)^{-1}$. Proposition~\ref{A:res:7.4A} and the form bounds for the Brownian potential of
Section~\ref{A:sec:finite-difference-approximation-of-the-stochastic-airy-operator} apply to both principal blocks. Their coordinates at the cut satisfy $x_J,x_{J+1}\asymp h^{-2}$.
Proposition~\ref{A:res:3.13}, applied to the full operator and separately to each block, gives
\[
 \|R(z)e_i\|^2+\|R^\oplus(z)e_i\|^2
 \le C_{\omega,\beta,K}h^4,\qquad i=J,J+1.
\]
The resolvent identity and the norm of a rank-one operator imply
\[
 \|R-R^\oplus\|_{\mathfrak S_2}^2 \le2|t_J|^2\bigl( \|Re_J\|^2\|R^\oplus(\overline z)e_{J+1}\|^2
 +\|Re_{J+1}\|^2\|R^\oplus(\overline z)e_J\|^2\bigr) \le C_{\omega,\beta,K}h^{-4}h^8
 =C_{\omega,\beta,K}h^4.
\]
With $\overline H_{n,\beta}^{B,\oplus} =\overline H_{J,\beta}^B\oplus\overline H_{\mathrm{tail},\beta}^B$ and
$\mathcal B_{n,h,\beta}^{B,\oplus} =\mathcal B_{J,h,\beta}^B\oplus\mathcal B_{\mathrm{tail},h,\beta}^B$, this proves
\[
 \| (\overline H_{n,\beta}^B-z)^{-1}
       -(\overline H_{n,\beta}^{B,\oplus}-z)^{-1}\|_{\mathfrak S_2}^2
 \le C_{\omega,\beta,K}h^4,
 \eqtag{7.60}
\]
\[
 \| (\mathcal B_{n,h,\beta}^B-z)^{-1}
       -(\mathcal B_{n,h,\beta}^{B,\oplus}-z)^{-1}\|_{\mathfrak S_2}^2
 \le C_{\omega,\beta,K}h^4.
 \eqtag{7.61}
\]
For either operator $H_{\mathrm{tail}}$ on $\{J+1,\ldots,n\}$, the common lower form bound gives
\[
 H_{\mathrm{tail}}+\kappa_B I\ge c_*(1+x_{J+1})I.
\]
Write $r_K=\sup_{z\in K}|z|$. For every eigenvalue $\lambda$ of $H_{\mathrm{tail}}$ and $z\in K$,
\[
 \frac{\lambda+\kappa_B}{|\lambda-z|}
 \le1+\frac{|z+\kappa_B|}{|\lambda-z|}
 \le1+\frac{r_K+\kappa_B}{\delta_K}.
\]
Since $1+x_{J+1}\ge c h^{-2}$, the spectral theorem yields
\[
 \begin{aligned}
 \sup_{z\in K}\|(H_{\mathrm{tail}}-z)^{-1}\|_{\mathrm{op}}
 &\le\frac{1+(r_K+\kappa_B)/\delta_K}{c_*(1+x_{J+1})}
       \le C_{I_\beta,K}(1+\kappa_B)h^2,\\
 \sup_{z\in K}\|(H_{\mathrm{tail}}-z)^{-1}\|_{\mathfrak S_2}^2
 &\le(n-J)C_{I_\beta,K}(1+\kappa_B)^2h^4
       \le C_{I_\beta,K}(1+\kappa_B)^2h.
 \end{aligned}
\]
The use of the shifted eigenvalue in the numerator makes these inequalities valid for every $n\ge N(I_\beta)$, with no further threshold depending on $\kappa_B$. Thus
\[
 \| (\overline H_{\mathrm{tail},\beta}^B-z)^{-1}
       -(\mathcal B_{\mathrm{tail},h,\beta}^B-z)^{-1}\|_{\mathfrak S_2}^2
 \le C_{\omega,\beta,K}h.
 \eqtag{7.64}
\]

\begin{proposition}\label{A:res:7.5}
For every $B\in\Omega_B$, every $\beta\in I_\beta$ and every $n\ge N(I_\beta)$,
\[
 \sup_{z\in K}\| (\overline H_{n,\beta}^B-z)^{-1}
        -(\mathcal B_{n,h,\beta}^B-z)^{-1}\|_{\mathfrak S_2}^2
 \le C_{\omega,\beta,K}h.
 \eqtag{7.65}
\]
The constant may be chosen uniformly over $\beta\in I_\beta$ and is independent of $n$.
\end{proposition}
\begin{proof}
Fix $B\in\Omega_B$, $\beta\in I_\beta$, and $n\ge N(I_\beta)$. Write
\[
 \overline R=(\overline H_{n,\beta}^{B}-z)^{-1},\qquad
 R^B=(\mathcal B_{n,h,\beta}^{B}-z)^{-1},
\]
and let $\overline R^{\oplus}$ and $R^{B,\oplus}$ be the inverses after cutting the edge between $J=\lfloor n/2\rfloor$ and $J+1$, as above. All four operators act on the same
Hilbert space with counting measure. We decompose their difference as
\[
 \overline R-R^B
 =(\overline R-\overline R^{\oplus})
 +(\overline R^{\oplus}-R^{B,\oplus})
 +(R^{B,\oplus}-R^B).
\]
For the middle term, the orthogonal decomposition of the coordinate space gives
\[
 \begin{split}
 \|\overline R^{\oplus}-R^{B,\oplus}\|_{\mathfrak S_2}^{2}
 ={}&\|(\overline H_{J,\beta}^{B}-z)^{-1}
       -(\mathcal B_{J,h,\beta}^{B}-z)^{-1}\|_{\mathfrak S_2}^{2}\\
 &+\|(\overline H_{\mathrm{tail},\beta}^{B}-z)^{-1}
       -(\mathcal B_{\mathrm{tail},h,\beta}^{B}-z)^{-1}\|_{\mathfrak S_2}^{2}.
 \end{split}
\]
Indeed, the squared Hilbert--Schmidt norm is additive over these two orthogonal blocks. Taking the supremum in $z$ and using \eqref{A:eq:7.58} and \eqref{A:eq:7.64} bounds this
term by $C_{B,I_\beta,K}h$. Equations \eqref{A:eq:7.60} and \eqref{A:eq:7.61} bound each of the other two squared norms by $C_{B,I_\beta,K}h^4$. Thus
\[
 \begin{split}
 \sup_{z\in K}\|\overline R-R^B\|_{\mathfrak S_2}^{2}
 &\le3\sup_K\|\overline R-\overline R^{\oplus}\|_{\mathfrak S_2}^{2}
     +3\sup_K\|\overline R^{\oplus}-R^{B,\oplus}\|_{\mathfrak S_2}^{2}\\
 &\quad+3\sup_K\|R^{B,\oplus}-R^B\|_{\mathfrak S_2}^{2}
 \le C_{B,I_\beta,K}(h+2h^4)\le3C_{B,I_\beta,K}h.
 \end{split}
\]
All the estimates have the deterministic size threshold $N(I_\beta)$ and constants uniform over $\beta\in I_\beta$, so their maximum has the same properties.
\end{proof}

\subsubsection{Resolvent comparison}\label{A:sec:gaussian-resolvent-comparison}
The maximum bounds for the coupled coefficients give common upper and lower quadratic form bounds. Once those are established, the weighted square sums estimate the difference of
the resolvents.

For the resolvent comparison, write
\[
 \begin{aligned}
 H_{n,\mathrm{DE}}^{[0]}&=\mathcal A_{n,\beta}^{\mathrm{DE}},\qquad
 H_{n,\mathrm{DE}}^{[1]}=\overline H_{n,\beta}+h^{-2}M_{\zeta_n},\\
 H_{n,\mathrm{DE}}^{[2]}&=\overline H_{n,\beta}^B,
 \qquad R_r(z)=(H_{n,\mathrm{DE}}^{[r]}-z)^{-1}\quad(r=0,1,2).
 \end{aligned}
\]
Proposition~\ref{A:res:6.1} and the definition of $q_n$ give
\[
 H_{n,\mathrm{DE}}^{[0]}-H_{n,\mathrm{DE}}^{[1]}=\Gamma_{\eta_n}^{-},
 \qquad H_{n,\mathrm{DE}}^{[1]}-H_{n,\mathrm{DE}}^{[2]}=M_{q_n}.
\]

\begin{proposition}\label{A:res:8.2B}
On $\Omega_{\mathrm{DE}}$, the three operators $H_{n,\mathrm{DE}}^{[r]}$ and $\mathcal B_{n,h,\beta}^B$ satisfy \eqref{A:eq:3.8}, \eqref{A:eq:3.10}, and \eqref{A:eq:3.12} with
common constants, for all sufficiently large $n$. The same holds for their principal restrictions with inherited coordinates.
\end{proposition}
\begin{proof}
Fix a compact interval $I_\beta=[\beta_-,\beta_+]\subset(0,\infty)$. Proposition~\ref{A:res:7.4A} supplies a deterministic $N_0=N_0(I_\beta)$, a deterministic lower coefficient
$c_*>0$, and a shift $\kappa_B\ge1$ such that
\[
 c_*\ell_{h,I}[u]\le
 \mathfrak h_{n,\mathrm{DE}}^{[2],I}[u]+\kappa_B\|u\|^2
 \le C_{I_\beta}(1+\kappa_B)\ell_{h,I}[u],
 \qquad
 \kappa_B\le C_{I_\beta}(1+\mathcal O_{1/4}(B))^4
\]
for every $B\in\Omega_B$, $\beta\in I_\beta$, $n\ge N_0$, and principal interval $I$. The deterministic part of the form of $\mathcal B_{n,h,\beta}^B$ is $\ell_h[u]-\|u\|^2$.
Applying \eqref{A:eq:7.29r} to its Brownian term, and increasing $\kappa_B$ and decreasing $c_*$ if necessary, gives the same bounds for this matrix and its principal restrictions.

At a fixed size, assume the three maximum bounds \eqref{A:eq:8.28a}--\eqref{A:eq:8.28b}. Equation~\eqref{A:eq:8.28c} and Cauchy--Schwarz give, for every vector supported on $I$,
\[
 \begin{aligned}
 |\langle M_{q_n}u,u\rangle|
 &\le C h^{3/4}\sum_{k\in I}(1+x_k)|u_k|^2,\\
 |\langle\Gamma_{\eta_n}^{-}u,u\rangle|
 &\le2\|\eta_n\|_\infty\|D_hu\|\,\|J_-u\|
 \le h^{1/4}\bigl(\|D_hu\|^2+\|u\|^2\bigr).
 \end{aligned}
\]
Thus the sum of these two absolute values is at most $(C h^{3/4}+h^{1/4})\ell_{h,I}[u]$. Choose a deterministic $N_1=N_1(I_\beta)\ge N_0$ for which $C h^{3/4}+h^{1/4}\le c_*/2$
when $n\ge N_1$. The exact identities defining $H_{n,\mathrm{DE}}^{[0]}$ and $H_{n,\mathrm{DE}}^{[1]}$ then imply, for $r=0,1$,
\[
 \frac{c_*}{2}\ell_{h,I}[u]
 \le\mathfrak h_{n,\mathrm{DE}}^{[r],I}[u]+\kappa_B\|u\|^2
 \le\bigl(C_{I_\beta}(1+\kappa_B)+c_*/2\bigr)\ell_{h,I}[u].
\]
This proves \eqref{A:eq:3.12} at every such size on the stated maximum events. For each fixed $\beta$ and $\omega\in\Omega_{\mathrm{DE}}$, they occur at every size
$n\ge n_0(\omega,\beta)$ for some finite $n_0(\omega,\beta)$.

For \eqref{A:eq:3.8}, take $c_{k,n}=\overline b_{k,n}$ for the three matrices $H_{n,\mathrm{DE}}^{[r]}$ and place their scalar terms in the diagonal coefficient. Jensen's
inequality and the auxiliary boundary point conventions give $0\le\overline b_{k,n}\le1$. Only $H_{n,\mathrm{DE}}^{[0]}$ has the additional term $\Gamma_{\eta_n}^{-}$, with
$\|\eta_n\|_\infty\le h^{1/4}\le1$. For $\mathcal B_{n,h,\beta}^B$ the kinetic coefficient is one. These choices verify \eqref{A:eq:3.10} with common deterministic coefficient
bounds.

For completeness, let $P_I$ be coordinate restriction and $\widetilde u=P_I^*u$ the zero extension of a vector on $I=\{p,\ldots,q\}$. For each matrix $H$ above,
\[
 \langle P_IHP_I^*u,u\rangle=\langle H\widetilde u,\widetilde u\rangle,\qquad
 \ell_h[\widetilde u]
 =\sum_{k=p-1}^{q}|D_h\widetilde u_k|^2
       +\sum_{k=p}^{q}(1+x_k)|u_k|^2
 =\ell_{h,I}[u].
\]
The boundary edges $p-1$ and $q$ retain their original coefficients. Hence the bounds hold with the original coordinates and the same constants for every principal restriction. For
a general principal index set, decompose it into its maximal consecutive intervals. Tridiagonality makes the restricted operator their direct sum, with the derivative and energy
form taken separately on each interval. Adding the form inequalities proves the same bounds, independently of the number of intervals.
\end{proof}

\begin{proposition}\label{A:res:8.4}
On $\Omega_{\mathrm{DE}}$, uniformly for all sufficiently large $n$,
\[
 \sup_{z\in K}\|R_0(z)-R_1(z)\|_{\mathfrak S_2}^2
 \le C_{\omega,\beta,K}h,
 \eqtag{8.37}
\]
\[
 \sup_{z\in K}\|R_1(z)-R_2(z)\|_{\mathfrak S_2}^2
 \le C_{\omega,\beta,K}h.
 \eqtag{8.38}
\]
\end{proposition}
\begin{proof}
The first weighted error satisfies
\[
 \sum_{k=0}^n h(1+x_k)^{-3/2}|\eta_{k,n}|^2
 =\mathcal E_{\eta,n}^{\mathrm{DE}}\le A_\eta h,
 \eqtag{8.29}
\]
since $\eta_{0,n}=0$. The resolvent identity is $R_0-R_1=-R_0\Gamma_{\eta_n}^{-}R_1$. Proposition~\ref{A:res:8.2B} verifies all the hypotheses of Proposition~\ref{A:res:3.17}.
Applying its first-order perturbation estimate gives
\[
 \sup_{z\in K}\|R_0(z)-R_1(z)\|_{\mathfrak S_2}^2
 \le C_{\omega,\beta,K}\mathcal E_{\eta,n}^{\mathrm{DE}}
 \le C_{\omega,\beta,K}h.
 \eqtag{8.32}
\]
For the second difference, the definition of $E_{n,\mathrm{DE}}$ gives
\[
 \sum_{k=1}^n(1+x_k)^{-7/2}|\Delta\zeta_{k,n}|^2
 =h^3\mathcal E_{\mathrm{scal},n}^{\mathrm{DE}}
 \le h^{9/2}\le h^4.
 \eqtag{8.33}
\]
Using $R_1-R_2=-R_1h^{-2}M_{\Delta\zeta_n}R_2$ and Proposition~\ref{A:res:3.16}, we obtain
\[
 \sup_{z\in K}\|R_1(z)-R_2(z)\|_{\mathfrak S_2}^2 \le C_{\omega,\beta,K}h^{-3}
 \sum_{k=1}^n(1+x_k)^{-7/2}|\Delta\zeta_{k,n}|^2 \le C_{\omega,\beta,K}h^{3/2}\le
 C_{\omega,\beta,K}h.
 \eqtag{8.36}
\]
Both bounds hold on $\Omega_{\mathrm{DE}}$. The scalar comparison is $O_\omega(h^{3/2})=o_\omega(h)$ on this event; only the weaker bound $O_\omega(h)$ is stated. The first-order
term is bounded at order $h$.

More precisely, Proposition~\ref{A:res:A.5} and the same diagonal estimate give, for each $0<\epsilon<1$, on a probability one subset of $\Omega_{\mathrm{DE}}$,
\[
 \sup_{z\in K}\|R_1(z)-R_2(z)\|_{\mathfrak S_2}^2
 \le C_{\omega,\beta,K}\mathcal E_{\mathrm{scal},n}^{\mathrm{DE}}
 \le C_{\omega,\beta,K,\epsilon}h^{2-\epsilon}=o(h).
\]
Take the intersection of these events over rational $\epsilon\in(0,1)$. For any real $\epsilon\in(0,1)$, choose a rational $0<\epsilon'<\epsilon$ and use
$h^{2-\epsilon'}\le h^{2-\epsilon}$ for $h\le1$. Thus the refinement holds for all $\epsilon\in(0,1)$ on one probability one subset of $\Omega_{\mathrm{DE}}$. The estimates
defining $\Omega_{\mathrm{DE}}$ alone already give the bound $O_\omega(h^{3/2})$ in \eqref{A:eq:8.36}.

\end{proof}

Combining Proposition~\ref{A:res:8.4} with the comparison of the mean coefficients in Proposition~\ref{A:res:7.5} gives the estimate for the full matrices.

\begin{theorem}\label{A:res:8.5}
On this coupling, almost surely, for all sufficiently large $n$,

\[
 \sup_{z\in K}\| (\mathcal A_{n,\beta}^{\mathrm{DE}}-z)^{-1}
                   -(\mathcal B_{n,h,\beta}^B-z)^{-1}\|_{\mathfrak S_2}^2
 \le C_{\omega,\beta,K}h.
\]
\end{theorem}
\begin{proof}
By Proposition~\ref{A:res:7.5},
\[
 \sup_{z\in K}\|R_2(z)-(\mathcal B_{n,h,\beta}^B-z)^{-1}\|_{\mathfrak S_2}^2
 \le C_{\omega,\beta,K}h.
 \eqtag{8.39}
\]
Decompose the desired difference as $(R_0-R_1)+(R_1-R_2)+(R_2-(\mathcal B_{n,h,\beta}^B-z)^{-1})$. The bounds just proved give, on the same event,
\[
 \sup_{z\in K}\|R_0(z)-(\mathcal B_{n,h,\beta}^B-z)^{-1}\|_{\mathfrak S_2}^2
 \le C_{\omega,\beta,K}(h+h^{3/2}+h)\le C_{\omega,\beta,K}h.
\]
Here we used $\|T_1+T_2+T_3\|_{\mathfrak S_2}^2\le3\sum_{j=1}^3\|T_j\|_{\mathfrak S_2}^2$ under the joint coupling law. The event $\Omega_{\mathrm{DE}}$ does not depend on $K$; the
weighted resolvent estimates apply to each compact $K$ on this same event.

The proof of Proposition~\ref{A:res:A.5} separately estimates the summand with $k=n$, using its exact Gaussian distribution. It has $L^p$ size $O(h^7)$ and is
$O_\omega(h^{7-\epsilon})$ almost surely for each $\epsilon>0$. This summand is already included in $\mathcal E_{\mathrm{scal},n}^{\mathrm{DE}}$. That sharper estimate holds on the
probability one event constructed in the proof of Proposition~\ref{A:res:A.5}; the bound on $\Omega_{\mathrm{DE}}$ above uses only
$\mathcal E_{\mathrm{scal},n}^{\mathrm{DE}}\le h^{3/2}$.

\end{proof}

The isometry $\iota_h$ preserves the Hilbert--Schmidt norm of a finite matrix, as in \eqref{A:eq:isometric-hs}. Combining this theorem with the estimate for the finite difference
approximation driven by the Brownian increments therefore gives the upper bound for the Gaussian ensemble at nonreal spectral parameters. The common spectral event, the lower
bound, and the passage to real inverses are completed in Section~\ref{A:sec:inverse-coupling-and-spectral-consequences}.

\subsubsection{Moment estimates for the Dumitriu--Edelman model}\label{A:sec:moment-estimates-for-the-dumitriu-edelman-model}
To obtain moment bounds at a fixed matrix size from the preceding comparison, we must also control the random constant and the starting size. We return to the coefficient bounds
and keep track of the dependence of the resolvent estimates on the shift in the form estimate for the Brownian potential.

Fix $1\le p<\infty$ and a compact interval $I_\beta=[\beta_-,\beta_+]\Subset(0,\infty)$. The coefficient bounds in Propositions~\ref{A:res:A.4}--\ref{A:res:A.5} are uniform for
$\beta\in I_\beta$. We use Lemma~\ref{A:res:analytic-constants} to retain the dependence of the resolvent estimates on the common form shift. Its proof in
Section~\ref{A:sec:parameter-dependence-and-further-directions} uses the deterministic form lower bound and the bound on the off-diagonal entries, and does not use the moment
conclusion proved here.

For each $n$ and $\beta$, let $F_{n,\beta}$ be the event consisting of the three maximum bounds in \eqref{A:eq:8.28a}--\eqref{A:eq:8.28b}, with the deterministic thresholds
displayed there and an empty maximum interpreted as zero. For every $q\ge1$, the proofs of \eqref{A:eq:8.28d}--\eqref{A:eq:8.28f} give
\[
 \sup_{\beta\in I_\beta}\mathbb P(F_{n,\beta}^c)
 \le C_{q,I_\beta}n^{1-q/12}.
\]
Given $A>0$, choosing $q>12(A+1)$ therefore gives
\[
 \sup_{\beta\in I_\beta}\mathbb P(F_{n,\beta}^c)
 \le C_{A,I_\beta}n^{-A}.
\]
On $F_{n,\beta}$ the three maxima give the form bounds. The two weighted errors will be estimated in expectation.

The proof of Proposition~\ref{A:res:7.4A} permits a common choice, uniform over $\beta\in I_\beta$,
\[
 \kappa_B=C_{I_\beta}(1+\mathcal O_{1/4}(B))^4\ge1.
\]
Increase the deterministic prefactor to cover the finite difference operator driven by the Brownian increments as well, and set $\kappa_B=1$ on the Wiener null set $\Omega_B^c$.
This is a measurable function of $B$. In the estimates below we work on $F_{n,\beta}\cap\{B\in\Omega_B\}$; removing the second event changes no expectation. There are deterministic
$n_0=n_0(I_\beta)$ and $c_*=c_*(I_\beta)>0$ such that, for $n\ge n_0$, Proposition~\ref{A:res:8.2B} gives $H+\kappa_B I\ge c_*\mathcal L_h$ on $F_{n,\beta}$ for every comparison
matrix $H$ and each principal restriction. The scaled off-diagonal entries are bounded by a deterministic constant. For the original matrix,
\[
 h^2|\mathcal A_{n,\beta}^{\mathrm{DE}}(k,k+1)|
 =b_{k,n}
 \le\overline b_{k,n}+h|\eta_{k,n}|
 \le1+h^{5/4}\le2.
\]
For the other comparison matrices use $\overline b_{k,n}\le1$, $b^0_{k,n}\le1$, or the unit kinetic coefficient. Their dimensions are at most $n=h^{-3}$, and their inherited
coordinates satisfy $1+x_k\le2h^{-2}$. Thus all the hypotheses of Lemma~\ref{A:res:analytic-constants} hold. Write $P_B=1+\kappa_B$ and $C=C_{I_\beta,K}$ in the estimates below.
For any of these resolvents $R$, and any pair $R_a,R_b$, that lemma gives
\[
 \begin{aligned}
 \|D_hR\|_{\mathrm{op}}^2&\le CP_B^2,\qquad
 \|D_hR\|_{\mathfrak S_2}^2\le CP_B^2h^{-1},\\
 \|RM_v\|_{\mathfrak S_2}^2&\le CP_B^4h\sum_k(1+x_k)^{-3/2}|v_k|^2,\\
 \|R_aM_vR_b\|_{\mathfrak S_2}^2&\le CP_B^{11}h\sum_k(1+x_k)^{-7/2}|v_k|^2.
 \end{aligned}
\]
The sums run over the relevant principal interval, and the estimates are uniform for $z\in K$. Combining the derivative operator norm with the multiplication estimate gives the
first-order constant $CP_B^6$. Lemma~\ref{A:res:A.6} gives a Gaussian tail for $\mathcal O_{1/4}(B)$, so
\[
 \|P_B^{12}\|_{L^r}\le C_{r,I_\beta}\qquad(1\le r<\infty).
\]
These moments also hold under the joint coupling law, whose Brownian marginal is Wiener measure.

Set
\[
 D_n=\sup_{z\in K}
  \|(\mathcal A_{n,\beta}^{\mathrm{DE}}-z)^{-1}
          -(\mathcal B_{n,h,\beta}^B-z)^{-1}\|_{\mathfrak S_2}^2,
 \qquad
 \mathcal Z_{\omega,h}
 =h^5\sum_{k=1}^n(1+x_k)^{-3/2}|\omega_{k,n}^B|^2.
\]
For the full resolvents $R_r$ of Subsection~\ref{A:sec:gaussian-resolvent-comparison}, the two exact resolvent identities in Proposition~\ref{A:res:8.4} and the preceding bounds
give, on $F_{n,\beta}$,
\[
 \begin{aligned}
 \sup_{z\in K}\|R_0-R_1\|_{\mathfrak S_2}^2
 &\le CP_B^6\mathcal E_{\eta,n}^{\mathrm{DE}},\\
 \sup_{z\in K}\|R_1-R_2\|_{\mathfrak S_2}^2
 &\le CP_B^{11}h\,h^{-4}
       \sum_{k=1}^n(1+x_k)^{-7/2}|\Delta\zeta_{k,n}|^2
 =CP_B^{11}\mathcal E_{\mathrm{scal},n}^{\mathrm{DE}}.
 \end{aligned}
\]

For the comparison with the mean coefficients, put $J=\lfloor n/2\rfloor$ and, on this initial block, write
\[
 \begin{aligned}
 \widehat R_0&=(\mathcal B_{J,h,\beta}^B-z)^{-1},&
 \widehat R_1&=(\mathcal B_{J,h,\beta}^B+\Delta_0-z)^{-1},\\
 \widehat R_2&=(\overline H_{J,\beta}^B-z)^{-1},
 \end{aligned}
\]
using $\Delta_0,a,b,r^0,c,s^\chi$ from Subsection~\ref{A:sec:the-mean-dumitriu-edelman-coefficients}. The coefficient estimates and integral comparison give
\[
 \begin{gathered}
 |a_k|\le Ch^2(1+x_k),\quad |b_k|\le Ch^2,\quad
 |r_k^0|\le Ch^2(1+x_k)^2,\\
 h^5\sum_{k\le J}(1+x_k)^{1/2}\le Ch.
 \end{gathered}
\]
Apply the displayed resolvent estimates to the six terms in \eqref{A:eq:7.38}--\eqref{A:eq:7.39}. Uniformly on $K$,
\[
 \begin{aligned}
 \|\widehat R_1M_a\|_{\mathfrak S_2}^2
 +\|z\widehat R_1M_a\widehat R_0\|_{\mathfrak S_2}^2
 &\le CP_B^4h^5\sum_{k\le J}(1+x_k)^{1/2}\le CP_B^4h,\\
 \|\widehat R_1M_{ax}\widehat R_0\|_{\mathfrak S_2}^2
 +\|\widehat R_1M_{r^0}\widehat R_0\|_{\mathfrak S_2}^2
 &\le CP_B^{11}h^5\sum_{k\le J}(1+x_k)^{1/2}\le CP_B^{11}h,\\
 \|\widehat R_1M_{a\omega^B}\widehat R_0\|_{\mathfrak S_2}^2
 &\le CP_B^{11}h^5\sum_{k\le J}(1+x_k)^{-3/2}|\omega_{k,n}^B|^2
 \le CP_B^{11}\mathcal Z_{\omega,h},\\
 \|\widehat R_1M_bJ_-^*D_h\widehat R_0\|_{\mathfrak S_2}^2
 &\le CP_B^6h^5\sum_{k\le J}(1+x_k)^{-3/2}\le CP_B^6h^4.
 \end{aligned}
\]
For the first line use also $\|\widehat R_0\|_{\mathrm{op}}\le\delta_K^{-1}$; for the last line use the factorization through $D_h\widehat R_0$ and
\eqref{A:eq:restricted-multiplier}. The remaining two terms satisfy
\[
 \begin{aligned}
 \|\widehat R_2D_h^*M_cD_h\widehat R_1\|_{\mathfrak S_2}^2
 &\le \|D_h\widehat R_2(\overline z)\|_{\mathrm{op}}^2
       \|c\|_\infty^2\|D_h\widehat R_1(z)\|_{\mathfrak S_2}^2
 \le CP_B^4h^5,\\
 \|\widehat R_2M_{s^\chi}\widehat R_1\|_{\mathfrak S_2}^2
 &\le CP_B^{11}h\sum_{k\le J}(1+x_k)^{-7/2}|s_{k,n}^\chi|^2
 \le CP_B^{11}h^2.
 \end{aligned}
\]
Here $\|c\|_\infty\le Ch^3$ includes the auxiliary edge $k=0$, where $c_{0,n}=0$, and the edge $k=J$.

To restore the full interval, use the two cuts in \eqref{A:eq:7.60}--\eqref{A:eq:7.61}. At their two adjacent vertices, the column bound is $CP_B^4h(1+x_k)^{-3/2}\le CP_B^4h^4$.
The rank-two resolvent identity therefore bounds each squared cut error by
\[
 C h^{-4}(P_B^4h^4)^2=CP_B^8h^4.
\]
Write $H_{\mathrm{tail}}$ for either principal restriction to $\{J+1,\ldots,n\}$. Its form bound gives $H_{\mathrm{tail}}+\kappa_B I\ge c_*(1+x_{J+1})I$. For each eigenvalue
$\lambda$ of $H_{\mathrm{tail}}$,
\[
 \frac{\lambda+\kappa_B}{|\lambda-z|}
 \le1+\frac{\kappa_B+|z|}{\delta_K}\le C_KP_B,\qquad
 |\lambda-z|^{-1}\le CP_Bh^2.
\]
Since there are at most $n=h^{-3}$ such eigenvalues,
$\|(H_{\mathrm{tail}}-z)^{-1}\|_{\mathfrak S_2}^2\le CP_B^2h$.
This bound holds for every $B\in\Omega_B$ and $n\ge n_0$, even when $\kappa_B$ is large. There is a fixed number of resolvent differences in these comparisons. The squared triangle
inequality and $P_B\ge1$, $h\le1$ therefore prove
\[
 \mathbf1_{F_{n,\beta}}D_n
 \le C_{I_\beta,K}\mathbf1_{F_{n,\beta}}P_B^{12}
       \bigl(h+\mathcal E_{\eta,n}^{\mathrm{DE}}
                 +\mathcal E_{\mathrm{scal},n}^{\mathrm{DE}}
                 +\mathcal Z_{\omega,h}\bigr).
\]

Each $\omega_{k,n}^B$ is a centered Gaussian of variance $4/(\beta h)$. Minkowski's inequality and the integral comparison give, for $r\ge1$,
\[
 \|\mathcal Z_{\omega,h}\|_{L^r} \le h^5\sum_{k=1}^n(1+x_k)^{-3/2} \|\omega_{k,n}^B\|_{L^{2r}}^2
 \le C_{r,I_\beta}h^4\sum_{k=1}^n(1+x_k)^{-3/2} \le C_{r,I_\beta}h^3.
\]
The coefficient errors and $P_B$ may be dependent, since all use the same Brownian motion. Propositions~\ref{A:res:A.4}--\ref{A:res:A.5} and H\"older's inequality control their
products under this joint law and give
\[
 \begin{aligned}
 \|\mathbf1_{F_{n,\beta}}D_n\|_{L^p}
 &\le C_{I_\beta,K}\|P_B^{12}\|_{L^{2p}}
  \bigl(h+\|\mathcal E_{\eta,n}^{\mathrm{DE}}\|_{L^{2p}}
          +\|\mathcal E_{\mathrm{scal},n}^{\mathrm{DE}}\|_{L^{2p}}
          +\|\mathcal Z_{\omega,h}\|_{L^{2p}}\bigr)\\
 &\le C_{p,I_\beta,K}(h+h^2+h^3)
 \le C_{p,I_\beta,K}h.
 \end{aligned}
\]
On $F_{n,\beta}^c$, we apply the nonreal resolvent bound to the same self-adjoint matrices defined by the coupling. Their laws remain unchanged, and their dimension gives
\[
 D_n\le2\sup_{z\in K}\|(\mathcal A_{n,\beta}^{\mathrm{DE}}-z)^{-1}
                      \|_{\mathfrak S_2}^2
       +2\sup_{z\in K}\|(\mathcal B_{n,h,\beta}^B-z)^{-1}
                      \|_{\mathfrak S_2}^2
 \le\frac{4n}{\delta_K^2}.
\]
Thus, choosing $A>4p/3$ gives
\[
 \|\mathbf1_{F_{n,\beta}^c}D_n\|_{L^p}
 \le\frac{4n}{\delta_K^2}\mathbb P(F_{n,\beta}^c)^{1/p}
 \le C_{p,A,I_\beta,K}n^{1-A/p}
 \le C_{p,A,I_\beta,K}h.
\]
For $n<n_0$ the same deterministic bound gives
$D_n\le4n/\delta_K^2\le(4n_0^{4/3}/\delta_K^2)h$.
Combining this with the two event estimates proves
\[
 \sup_{\beta\in I_\beta}
 \left\|\sup_{z\in K}
  \|(\mathcal A_{n,\beta}^{\mathrm{DE}}-z)^{-1}
       -(\mathcal B_{n,h,\beta}^B-z)^{-1}\|_{\mathfrak S_2}^2
       \right\|_{L^p}
 \le C_{p,I_\beta,K}h.
\]

For the remaining discretization error, $L=(n+1)h=h^{-2}+h$ satisfies $h^{-2}\le L\le2h^{-2}$. Thus the spacing and interval hypotheses of Lemma~\ref{A:res:discretization-moments}
hold. Its proof integrates the estimates for sampling the potential and restricting the variational equation of
Section~\ref{A:sec:finite-difference-approximation-of-the-stochastic-airy-operator}, under the full Wiener law. The Brownian marginal of our coupling is this same law, so the lemma
gives
\[
 \sup_{\beta\in I_\beta}
 \left\|\sup_{z\in K}\|\mathcal R_n^B(z)-R_\beta(z)\|_{\mathfrak S_2}^2
       \right\|_{L^p}\le C_{p,I_\beta,K}h.
\]
The isometry identity \eqref{A:eq:isometric-hs}, the squared triangle inequality in $\mathfrak S_2$, and Minkowski's inequality in $L^p$ therefore give
\[
 \sup_{\beta\in I_\beta}
 \left\|\sup_{z\in K}\|\mathcal R_n(z)-R_\beta(z)\|_{\mathfrak S_2}^2
       \right\|_{L^p}\le C_{p,I_\beta,K}h.
\]
The supremum in $z$ is measurable, since the resolvents are continuous on $K$ and a fixed countable dense subset suffices. Lemma~\ref{A:res:8.2} makes the finite coefficient
construction Borel in $\beta$ as well. Uniformity in $\beta$ throughout this argument means the displayed supremum of the separate $L^p$ norms under their specified coupling laws.

\subsection{The general Jacobi model: deterministic reduction}
\label{sec:combined-a-deterministic}
\label{A:sec:the-general-jacobi-model-deterministic-reduction}
For the Gaussian model, independence allowed us to couple each pair of coefficients separately. For a polynomial potential, the coefficient law has finite-range interactions. We
therefore fix the exterior of an interval and expand its exact conditional energy at its minimizer. The local profile from the introduction serves a different purpose: it
identifies the deterministic edge and its slope. We must compare these two minimizers before replacing one by the other.

Uniform convexity controls the fluctuations and the inverse Hessian; finite range gives decay of the response to exterior data. These estimates compare the conditional minimizer
with the local profile. The boundary resolvent estimates and the coefficient decomposition follow in Subsections~\ref{A:sec:boundary-reduction}
and~\ref{A:sec:decomposition-of-the-general-jacobi-coefficients}. The Gaussian approximation itself is constructed in
Section~\ref{A:sec:the-general-jacobi-model-gaussian-approximation-and-comparison}.

\subsubsection{Conditional laws and minimizers}\label{A:sec:log-concave-measures-and-conditional-minimizers}
The polynomial trace makes the dependence between Jacobi coefficients local. Writing it as a sum over closed paths identifies which exterior coordinates enter an interval's
conditional energy and hence its minimizer.

The coefficient law \eqref{A:eq:jacobi-law} has energy
\[
 \mathscr H_n(a,b)
 =\operatorname{Tr}V(T(a,b))
       -2\sum_{k=1}^{n-1}\left(1-\frac kn-\frac1{\beta n}\right)\log b_k,
 \qquad
 \mathrm d\mu_{n,V,\beta}\ \propto\
       \exp\left(-\frac{\beta n}{2}\mathscr H_n\right)\mathrm da\,\mathrm db.
 \eqtag{11.1}
\]

Let \(V(x)=\sum_{\ell=0}^{d}v_\ell x^\ell\) with \(d=\deg V\). Since \(T(a,b)\) is tridiagonal, the trace of \(T^\ell\) is a sum over closed nearest-neighbor paths of length
\(\ell\):
\[
\operatorname{Tr}T^\ell
=
\sum_{i=1}^{n}
\sum_{\substack{\gamma_0=i,\ \gamma_\ell=i\\
|\gamma_{r+1}-\gamma_r|\le1}}
\prod_{r=0}^{\ell-1}
T_{\gamma_r,\gamma_{r+1}} .
\]
Each factor is a diagonal coefficient \(a_j\) or an off-diagonal coefficient \(b_j\). A path of length \(\ell\) visits only sites within distance $\ell$ of its starting point. Thus
differentiating $\operatorname{Tr}V(T)$ with respect to a coefficient at index $k$ involves only coefficients at distance at most $R_V$ from $k$, and $R_V=d$ is sufficient.

For a coefficient $\xi=a_k$ or $b_k$, set $\jmath(a_k)=k$ and $\jmath(b_k)=k+\tfrac12$. Distances between coefficients refer to these site labels. In the interlaced order
$a_1,b_1,a_2,b_2,\ldots$, they are equivalent, up to a factor two, to distances between coordinate indices. Fix $R=R_V=\deg V$; the path expansion gives
\[
 \partial_\eta\partial_\xi\mathscr H_n=0
  \quad\text{if }|\jmath(\eta)-\jmath(\xi)|>R.
\]
There are at most two coefficients per unit of site distance. Thus the constants in the estimates for band matrices below do not depend on the length of the coefficient interval.

Fix $0<s_0<s_1<1$, retain the terminal cut $s_\beta$ of the

introduction, and put $N_n=n-s_\beta$. A coefficient interval $J$ consists of free coordinates in $a_1,b_1,\ldots,b_{N_n-1},a_{N_n}$; subintervals and intervals of diagonal
coordinates are also allowed. Put $w_{k,n}=1-k/n-1/(\beta n)$. Every free off-diagonal coordinate has $w_{k,n}>0$, since
\[
 w_{k,n}\ge\frac{s_\beta+1-\beta^{-1}}n>0,
 \qquad k\le N_n-1.
\]
This restriction is essential for the conditional minimization. For $\beta<1$, $w_{n-1,n}<0$ and the full energy tends to $-\infty$ as $b_{n-1}\downarrow0$ with the other
coordinates fixed. Uniform convexity of $V$ therefore does not make the unrestricted coefficient energy uniformly convex for every $\beta$.

Let $\partial J$ consist of the actual fixed coefficient variables at site distance at most $R_V$ from $J$, and let $\mathscr F_{J^c}$ be the sigma-field generated by all fixed
coordinates. Given their values $\tau$, the conditional law is
\[
 \mathrm{d}\mu_{n,J}^{\tau}(y)
 =Z_{n,J}(\tau)^{-1}\mathrm{e}^{-\Phi_{n,J}^{\tau}(y)}
   \mathbf1_{\mathcal D_J}(y)\,\mathrm{d}y,
 \quad
 \Phi_{n,J}^{\tau}(y)=\frac{\beta n}{2}\mathscr H_n(y,\tau)
       +c_n(\tau),
\eqtag{cond-density}
\]
where $\mathcal D_J$ is the product of real lines and positive half-lines. Choose $c_n(\tau)$ to cancel the terms depending only on $\tau$. The conditional density then depends on
the exterior only through $\partial J$. Terms containing both a free and a fixed coefficient remain in the conditional energy and contribute to its stationary equations and
Hessian. For example, if $V(t)=t^2/2+g t^4/4$ with $g>0$, the trace contains $g b_{k-1}^2b_k^2$. Fixing $b_{k-1}=\tau_{b,k-1}$ leaves the term $g\tau_{b,k-1}^2b_k^2$ in the
conditional energy of $b_k$. Consequently, conditioning the coefficient law is not the same operation as deleting rows and columns of the Jacobi matrix. For a deterministic
\(M\ge1\), set
\[
 \mathcal T_{n,J}(M)
 =\{\tau:\max_{i\in\partial J}|\tau_i|\le M\}.
\eqtag{good-data}
\]
The set $\mathcal T_{n,J}(M)$ restricts only the exterior coordinates. For each such $\tau$, $y$ still ranges over the full domain $\mathcal D_J$. In particular, this is not a
conditional law restricted to a neighborhood of its minimizer. The bounds proved uniformly on $\mathcal T_{n,J}(M)$ will later be integrated against the original exterior
distribution in Lemma~\ref{A:res:disintegration}.

We first record the convexity estimates used for these conditional laws. Lemma~\ref{A:res:A.8} bounds the tail probabilities of the Euclidean distance from the minimizer by those
of a Gaussian norm. The concentration inequality \eqref{A:eq:convex-concentration} bounds each coordinate about its mean, with constants independent of the dimension.
Lemma~\ref{A:res:mean} estimates the difference between that mean and the minimizer.

Let $\mathcal D\subset\mathbb R^d$ be a convex set with nonempty interior, and let
\[
 \mu(\mathrm dy)=Z^{-1}\mathrm e^{-\Phi(y)}\mathbf1_{\mathcal D}(y)\,\mathrm dy,
 \qquad \nabla^2\Phi\ge cnI_d.
\]
We assume that $\Phi$ is twice continuously differentiable in the interior and has a finite minimum at $y^\circ\in\overline{\mathcal D}$. At a boundary minimum we use its
continuous value and one-sided directional derivatives along admissible rays. These assumptions include the product domains used below.

\begin{lemma}\label{A:res:A.8}
If \(Y\sim\mu\) and \(G\sim N(0,(cn)^{-1}I_d)\), then \(\|Y-y^\circ\|_2\) is stochastically dominated by \(\|G\|_2\). Consequently, for every finite \(p\),
\[
\bigl\|\,\|Y-y^\circ\|_2\,\bigr\|_{L^p(\mu)}
\le C_{p,c}\sqrt{\frac{d}{n}},
\]
and there are constants \(C,c_1>0\) such that
\[
\mathbb P\left(
\|Y-y^\circ\|_2>C\sqrt{\frac dn}+t
\right)
\le 2\mathrm{e}^{-c_1nt^2},
\quad t\ge0.
\]
\end{lemma}

\begin{proof}
Put $\rho=cn$ and translate the minimizing point to zero. For a unit vector $\theta$, convexity makes the admissible radii an interval $(0,R_\theta)$, with $R_\theta\in[0,\infty]$.
A direction which meets the interior has a nonnegative one-sided derivative of $r\mapsto\Phi(r\theta)$ at zero. Strong convexity along that ray implies
\[
 \frac{\mathrm d}{\mathrm dr}\left(\Phi(r\theta)-\Phi(0)-\frac{\rho r^2}{2}\right)\ge0,
 \qquad 0<r<R_\theta.
\]
The same conclusion follows by first starting at a positive radius and then taking its limit to zero. Define
\[
 w_\theta(r)=\mathbf1_{\{r<R_\theta\}}
 \exp\left\{-\Phi(r\theta)+\Phi(0)+\frac{\rho r^2}{2}\right\}.
\]
For almost every direction this is nonnegative, bounded by one, and nonincreasing. Directions which do not meet the interior contribute zero to polar integration. Thus
\[
 W(r)=\int_{\mathbb S^{d-1}}w_\theta(r)\,\mathrm d\theta
\]
is nonincreasing. Polar coordinates show that the density of $\|Y\|_2$ is proportional to $r^{d-1}e^{-\rho r^2/2}W(r)$. For $d=1$ the sphere integral is the sum over the two
directions, and the same formula holds.

Let $R$ have the radial law of $N(0,\rho^{-1}I_d)$, and let $R'$ be an independent copy. The normalizing constant $\mathbb EW(R)$ is finite and strictly positive. For any bounded
nondecreasing function $f$,
\[
 \begin{split}
 2\operatorname{Cov}(W(R),f(R))
 &=\mathbb E[(W(R)-W(R'))(f(R)-f(R'))]\\
 &\le0.
 \end{split}
\]
Consequently $\mathbb Ef(\|Y\|_2)=\mathbb E[W(R)f(R)]/\mathbb EW(R)\le\mathbb Ef(R)$. Taking $f=\mathbf1_{(r,\infty)}$ proves the stochastic domination, including when the minimum
lies on the boundary.

For completeness, if $Z\sim N(0,I_d)$ and $X=\|Z\|_2^2$, independence and the one-dimensional Gaussian integral give
\[
 \mathbb Ee^{tX}=(1-2t)^{-d/2},\qquad
 \log\mathbb Ee^{t(X-d)}
 \le\frac{dt^2}{1-2t},\qquad 0<t<\tfrac12.
\]
Here $-\log(1-u)-u=\sum_{k\ge2}u^k/k\le u^2/[2(1-u)]$. Markov's inequality with $t=\sqrt{s}/(\sqrt d+2\sqrt{s})$ yields
\[
 \mathbb P\{X>d+2\sqrt{ds}+2s\}\le e^{-s},\qquad s>0.
\]
Since $(\sqrt d+\sqrt{2s})^2\ge d+2\sqrt{ds}+2s$, domination gives
\[
 \mathbb P\left\{\|Y-y^\circ\|_2>
       \sqrt{d/\rho}+t\right\}\le e^{-\rho t^2/2},\qquad t\ge0.
\]
Finally, for $p>0$, integrate this bound for
$U=(\|Y-y^\circ\|_2-\sqrt{d/\rho})_+$:
\[
 \mathbb EU^p\le p\int_0^\infty t^{p-1}e^{-\rho t^2/2}\,\mathrm dt
 =2^{p/2}\Gamma(1+p/2)\rho^{-p/2}.
\]
The inequality $(a+b)^p\le\max\{1,2^{p-1}\}(a^p+b^p)$ and $d\ge1$ give the stated moment bound. All constants are independent of $d$, $n$, and the location of the minimum.
\end{proof}

For $d\asymp n$, this gives a bound of order one on the Euclidean displacement, sufficient to bound the coefficients. The sharper coordinate estimate about the minimizer follows by
combining concentration about the mean with Lemma~\ref{A:res:mean}.

Let $\mu(\mathrm dy)=Z^{-1}e^{-\Phi(y)}\mathbf1_D(y)\,\mathrm dy$ on a convex domain $D\subset\mathbb R^d$ with nonempty interior,

and suppose $\nabla^2\Phi\ge\rho I$, $\rho>0$. The functional inequalities used below are
\[
 \operatorname{Var}_\mu f\le\rho^{-1}\int|\nabla f|^2\,\mathrm d\mu,
 \qquad
 \operatorname{Ent}_\mu(f^2)\le2\rho^{-1}\int|\nabla f|^2\,\mathrm d\mu.
 \eqtag{convex-functional-inequalities}
\]
Here $\operatorname{Ent}_\mu g=\int g\log g\,\mathrm d\mu -(\int g\,\mathrm d\mu)\log\int g\,\mathrm d\mu$ for $g\ge0$. For the stated convex support, extend $F=\Phi-\rho|x|^2/2$
lower semicontinuously to $\overline D$ and by $+\infty$ outside. Caffarelli's contraction theorem \cite{caffarelli2000monotonicity} (Theorem~11, pp.~559--561) includes these
convex support constraints and gives a one-Lipschitz transport from $N(0,\rho^{-1}I)$ to the normalized measure $\mu$. The Gaussian logarithmic Sobolev constant is $2/\rho$ by the
Bakry--\'Emery criterion and scaling \cite{bakry1985diffusions} (Proposition~4 and Corollary~2, pp.~194--199). Composition with the contraction preserves this constant by
Corollary~8 of \cite{caffarelli2000monotonicity} (p.~562); linearizing at $f=1$ gives the variance constant $1/\rho$. This proves \eqref{A:eq:convex-functional-inequalities} for
bounded smooth functions. Value truncations, cutoffs and closure in the form norm extend it to the Sobolev and Lipschitz functions used here, by lower semicontinuity of entropy and
convergence of the gradient norms.

For an $L_f$-Lipschitz function $f$, put $\psi(t)=\log\mathbb E_\mu\mathrm e^{t(f-\mathbb E_\mu f)}$. Applying the entropy inequality to $\mathrm e^{tf/2}$ gives
\[
 t\psi'(t)-\psi(t)\le\frac{t^2L_f^2}{2\rho},\qquad
 \psi(t)\le\frac{t^2L_f^2}{2\rho}.
\]
The second inequality follows by integrating $(\psi(t)/t)'$ from zero; bounded truncations justify the calculation before taking limits. Optimizing exponential Markov and applying
the same bound to $-f$ yields
\[
 \mu\{|f-\mathbb E_\mu f|>u\}
 \le2\exp\{-\rho u^2/(2L_f^2)\},\qquad
 \|f-\mathbb E_\mu f\|_{L^p(\mu)}\le C_pL_f/\sqrt\rho.
 \eqtag{convex-concentration}
\]
For a linear functional $f(y)=v\cdot y$, $L_f=\|v\|_2$; the corresponding $\psi_2$ norm is at most $C\|v\|_2/\sqrt\rho$. These estimates are centered at $\mathbb E_\mu f$.

To compare conditional minimizers with different exterior data, we will invert an averaged Hessian. Its polynomial part has finite range. Small free off-diagonal coefficients can
make the positive logarithmic Hessian large, so we keep that diagonal term separate. Exponential conjugation, as in \cite{combes1973asymptotic}, leaves this diagonal part unchanged
and only needs a bound on the polynomial part.

Let $K=P+D$ be a real symmetric matrix on a finite graph, with graph distance $\mathrm d$. Suppose
\[
 D\ge0\text{ is diagonal},\qquad K\ge cnI,\qquad
 P_{ij}=0\quad(\mathrm d(i,j)>R_0),\qquad
 \sup_i\sum_j|P_{ij}|\le Cn.
\]
The estimate below is independent of the graph size and of the entries of $D$.

\begin{lemma}\label{A:res:A.9}
There are $a,C_1>0$, depending only on $c,C,R_0$, such that
\[
 |(K^{-1})_{ij}|\le\frac{C_1}{n}\mathrm e^{-a\mathrm d(i,j)}.
\]
\end{lemma}

\begin{proof}
Set $\overline R=\max(1,R_0)$, $S=\max(c,C)$, and $a=\overline R^{-1}\log(1+c/(2S))$. Fix $i$ and let $W_i=\operatorname{diag}(\mathrm e^{-a\mathrm d(i,k)})_k$ on its connected
component. Entries between different components vanish, so it suffices to work there. The diagonal part commutes with $W_i$, and
\[
 W_iKW_i^{-1}=K+E_i,\qquad
 (E_i)_{k\ell}=P_{k\ell}
   \bigl(\mathrm e^{-a(\mathrm d(i,k)-\mathrm d(i,\ell))}-1\bigr).
\]
On the support of $P_{k\ell}$ the two distances differ by at most $R_0$. Each absolute row sum of $E_i$ is at most $Cn(\mathrm e^{aR_0}-1)\le cn/2$. Symmetry of $P$ gives the same
column bound, although $E_i$ need not be symmetric. Schur's test yields $\|E_i\|_{\mathrm{op}}\le cn/2$, hence
\[
 \|(I+K^{-1}E_i)^{-1}\|_{\mathrm{op}}
 \le\sum_{m\ge0}\|K^{-1}E_i\|_{\mathrm{op}}^m\le2.
\]
Therefore $\|W_iK^{-1}W_i^{-1}\|_{\mathrm{op}}\le2/(cn)$. Its $(i,j)$ entry equals $\mathrm e^{a\mathrm d(i,j)}(K^{-1})_{ij}$, proving the claim with $C_1=2/c$.
\end{proof}

With upper and lower spectral bounds, a geometric series gives exponential decay of the entries of the inverse of a band matrix. This estimate is the case of positive band matrices
of Theorem~2.4 of \cite{demko1984decay}. Let \(A\) be a real symmetric positive definite \(d\times d\) matrix. Assume \(0<\lambda\le\Lambda\), \(\lambda I\le A\le\Lambda I\), and
\(A_{ij}=0\) whenever \(|i-j|>R\), where

\(R\in\mathbb N_0\) is an integer bandwidth.
\begin{lemma}\label{A:res:A.10}
There exist constants \(C,c>0\), depending only on \(R\) and \(\Lambda/\lambda\), such that
\[
|(A^{-1})_{ij}|
\le
\frac{C}{\lambda}\mathrm{e}^{-c|i-j|}.
\eqtag{A.30}
\]
In particular, if \(cnI\le A\le CnI,\) then \(|(A^{-1})_{ij}| \le \frac{C}{n}\mathrm{e}^{-c|i-j|}.\)
\end{lemma}

\begin{proof}
If $R=0$, the matrix is diagonal and $|(A^{-1})_{ij}|\le\lambda^{-1}\mathbf1_{\{i=j\}}$. If $\Lambda=\lambda$, the spectral theorem gives $A=\lambda I$ and the same conclusion.
Hence assume $R\ge1$ and $\Lambda>\lambda$, and put
\[
 B=I-A/\Lambda,\qquad q=1-\lambda/\Lambda\in(0,1).
\]
The spectral bounds imply $0\le B\le qI$. Thus the geometric series converges in operator norm and
\[
 A^{-1}=\Lambda^{-1}\sum_{m=0}^{\infty}B^m,
 \qquad \|B^m\|_{\mathrm{op}}\le q^m.
\]
The matrix $B$ has bandwidth at most $R$. In the expansion of $(B^m)_{ij}$, a nonzero product uses $m$ successive index differences, each at most $R$ in absolute value. Therefore
$(B^m)_{ij}=0$ if $|i-j|>mR$. With $m_0=\lceil |i-j|/R\rceil$, including $m_0=0$ for $i=j$, we obtain
\[
 |(A^{-1})_{ij}|
 \le\Lambda^{-1}\sum_{m=m_0}^{\infty}q^m
 =\lambda^{-1}q^{m_0}
 \le\lambda^{-1}\exp\left\{-\frac{-\log q}{R}|i-j|\right\}.
\]
This proves \eqref{A:eq:A.30} with constants depending only on $R$ and $\Lambda/\lambda$. For $cnI\le A\le CnI$, apply this argument to $A/n$ and use $A^{-1}=n^{-1}(A/n)^{-1}$.
\end{proof}

The stationary equations express the displacement of a trial profile from the minimizer through an inverse Hessian. The gradient at the profile may be small in the interior and
larger near the boundary. Exponential decay of the inverse Hessian then makes the effect of the boundary errors small far from the boundary.

Let \(I=\{\ell,\ldots,m\}\subset\mathbb Z\) be nonempty, let \(\mathcal D_I\subset\mathbb R^I\) be open and convex, and let \(\Phi_I\in C^2(\mathcal D_I)\) be strongly convex.
Suppose that its minimizer \(y^\star\) belongs to \(\mathcal D_I\), and fix \(y^\dagger\in\mathcal D_I\). All derivatives below are with respect to the free coordinates in \(I\);
any conditional exterior data are held fixed. In particular, \(\nabla\Phi_I(y^\star)=0\). Write
\[
 d_\partial(i)=\min\{i-\ell,m-i\},\quad
 S_\partial^{(w)}=\{j\in I:d_\partial(j)\le w\},
\]
where the integer \(w\ge0\) is independent of \(n\) and \(I\). The two neighborhoods of the boundary may overlap. Assume a decomposition
\[
 g:=\nabla\Phi_I(y^\dagger) =g_0+g_\partial,\quad \|g_0\|_\infty\le\varepsilon_n,\qquad
 \operatorname{supp}g_\partial \subset S_\partial^{(w)},\quad B_\partial:=\|g_\partial\|_\infty.
\eqtag{boundary-residual-data}
\]
Here $\varepsilon_n\ge0$ and $B_\partial\ge0$ bound the two parts separately. A bound on \(g\) only outside
\(S_\partial^{(w)}\) fits this setup by taking \(g_0=g\mathbf1_{I\setminus S_\partial^{(w)}}\) and
\(g_\partial=g\mathbf1_{S_\partial^{(w)}}\); norms on empty sets are zero.

Assume that the Hessian has bandwidth $R$ on the segment joining $y^\dagger$ to $y^\star$. We allow either $\lambda I\le\nabla^2\Phi_I\le\Lambda I$ on this segment, or the more
general condition
\[
 K:=\int_0^1\nabla^2\Phi_I(y^\dagger+t(y^\star-y^\dagger))\,\mathrm dt
     =P+D\ge\lambda I,
 \qquad D\ge0\text{ diagonal},\quad
 \sup_i\sum_j|P_{ij}|\le L,
\]
where $P$ is symmetric and has bandwidth $R$. The first case is included by taking $P=K$, $D=0$, and $L=(2R+1)\Lambda$. For separate boundary estimates we use a specified
decomposition $g=g_0+g_{\mathrm L}+g_{\mathrm R}$, where
\[
 \begin{gathered}
 \operatorname{supp}g_{\mathrm L}\subset\{j:j-\ell\le w\},\qquad
 \operatorname{supp}g_{\mathrm R}\subset\{j:m-j\le w\},\qquad
 \|g_0\|_\infty\le\varepsilon_n,\\
 B_{\mathrm L}=\|g_{\mathrm L}\|_\infty,\qquad
 B_{\mathrm R}=\|g_{\mathrm R}\|_\infty.
 \end{gathered}
\]

\begin{proposition}\label{A:res:A.11}
Under these hypotheses,
\[
 |y_i^\star-y_i^\dagger|
 \le\frac C\lambda\bigl(\varepsilon_n+B_\partial\mathrm e^{-ad_\partial(i)}\bigr),
 \qquad i\in I.
\eqtag{A.34}
\]
For the separate boundary decomposition,
\[
 |y_i^\star-y_i^\dagger|
 \le\frac C\lambda\bigl(\varepsilon_n+
 B_{\mathrm L}\mathrm e^{-a(i-\ell)}+B_{\mathrm R}\mathrm e^{-a(m-i)}\bigr).
\eqtag{boundary-residual-two-sides}
\]
Here $C,a>0$ depend only on $R,w,L/\lambda$.
\end{proposition}

\begin{proof}
Set \(\delta=y^\star-y^\dagger\). Convexity of \(\mathcal D_I\) keeps the entire segment inside the domain. The fundamental theorem of calculus, applied to the gradient in the free
coordinates, gives
\[
 0=\nabla\Phi_I(y^\star)=g+K\delta,
 \quad \delta=-K^{-1}g.
\]
Apply Lemma~\ref{A:res:A.9} with scale one, lower bound $\lambda$, and row bound $L$. Its proof permits the choice
\[
 a=\max(1,R)^{-1}\log\left(1+\frac{\lambda}{2\max(\lambda,L)}\right)
\]
and gives
\[
 |(K^{-1})_{ij}|\le\frac2\lambda\mathrm e^{-a|i-j|}.
\eqtag{boundary-residual-kernel}
\]

By \eqref{A:eq:boundary-residual-data} and \eqref{A:eq:boundary-residual-kernel}, the part involving $g_0$ satisfies
\[
 |(K^{-1}g_0)_i|
 \le\frac{2\varepsilon_n}{\lambda}
       \sum_{j\in I}\mathrm{e}^{-a|i-j|}
 \le\frac{2\varepsilon_n}{\lambda}
       \left(1+\frac{2}{\mathrm{e}^a-1}\right).
\]
There are at most \(2(w+1)\) distinct sites in \(S_\partial^{(w)}\), including when the neighborhoods of the boundary overlap. For each such \(j\), one of the boundary points
\(b\in\{\ell,m\}\) satisfies \(|j-b|\le w\). The triangle inequality gives \(|i-j|\ge|i-b|-w\ge d_\partial(i)-w\). Consequently
\[
 |(K^{-1}g_\partial)_i|
 \le\frac{4(w+1)\mathrm{e}^{aw}}{\lambda}
        B_\partial \mathrm{e}^{-a d_\partial(i)}.
\]
Combining the last two displays proves \eqref{A:eq:A.34}, for example with
\(C=2(1+2/(\mathrm{e}^a-1))+4(w+1)\mathrm{e}^{aw}\). This constant depends only on the quantities stated in the
proposition. For $K^{-1}g_{\mathrm L}$, the sum has at most $w+1$ terms, and $|i-j|\ge(i-\ell)-w$ on its support. For $K^{-1}g_{\mathrm R}$, use \(|i-j|\ge(m-i)-w\). Summing these
two contributions proves \eqref{A:eq:boundary-residual-two-sides}. Linearity gives the same bound when the two supports overlap.
\end{proof}

The dependence on \(B_\partial\) cannot be dropped. For instance, \(\Phi_I(y)=\frac12\|y\|^2-M y_\ell\) and
\(y^\dagger=0\) have \(\varepsilon_n=0\), \(B_\partial=M\), and \(y^\star=M e_\ell\) for every \(M>0\). Thus fixed
support alone does not give a uniform displacement even at the boundary point. More generally, suppose $\varepsilon_n/\lambda\le C_0/n$ and $B_\partial/\lambda\le C_b n^p$, with
$p\ge0$ and uniform $R,w,L/\lambda$. Then \eqref{A:eq:A.34} gives $|y_i^\star-y_i^\dagger|\le C/n$ whenever $d_\partial(i)\ge(p+1)a^{-1}\log n$. It is the ratio of each gradient
error to the Hessian scale that enters this estimate. Multiplying an energy by a positive constant changes both quantities by the same factor and leaves the displacement bound
unchanged.

We apply the preceding estimates to the exact density \eqref{A:eq:cond-density}. The factor $\beta n/2$ belongs to its exponent: a lower Hessian bound for $\mathscr H_n$
independent of $n$ becomes a lower bound proportional to $n$ for $\Phi_{n,J}^\tau$. Upper bounds also require control of the coefficients; they do not follow from convexity. The
minimizer is the same for both normalizations.

\begin{lemma}\label{A:res:12.1A}
The conditional law of the initial coefficients and every further conditional law just described satisfy
\[
 D^2\Phi_{n,J}^{\tau}(y)\ge c_{V,\beta}nI
 \quad(y\in\mathcal D_J)
\eqtag{12.15a}
\]
for every exterior value in the domain with real diagonal and positive off-diagonal coefficients. Their unique minimizers \(y_J^\circ(\tau)\) satisfy, for some integer \(D_V\) and
constants independent of \(n,J,\tau\),
\[
 \|y_J^\circ(\tau)\|_\infty
 \le C_V(1+\|\tau\|_{\infty,\partial J})^{D_V}.
\eqtag{mode-growth}
\]
For \(\tau\in\mathcal T_{n,J}(M)\), every free off-diagonal coefficient with \(k/n\le s_1\) is at least \(b_{M,s_1}>0\) at the minimizer, for all sufficiently large \(n\).
\end{lemma}

Write $d_J=|J|\le2n$. Lemma~\ref{A:res:A.8}, applied to $\mu_{n,J}^\tau$, gives
\[
 \begin{split}
 \bigl\|\,\|Y-y_J^\circ(\tau)\|_2\,\bigr\|_{L^p(\mu_{n,J}^{\tau})}
 &\le C_{p,V,\beta}\sqrt{d_J/n},\\
 \mu_{n,J}^{\tau}\{\|Y-y_J^\circ(\tau)\|_2>
              C\sqrt{d_J/n}+t\}
 &\le 2\mathrm{e}^{-cnt^2},\quad t\ge0,
 \end{split}
\eqtag{conditional-radius}
\]
Consequently, there is a deterministic \(M_1=M_1(M,V,\beta)\) such that
\[
 \sup_{\tau\in\mathcal T_{n,J}(M)}
 \mu_{n,J}^{\tau}\{\max_{i\in J}|Y_i|>M_1\}
 \le C\mathrm{e}^{-cn}.
\eqtag{conditional-upper}
\]
The constants are uniform for \(\beta\) in a fixed compact subinterval of \((0,\infty)\), when one terminal block removes all the corresponding negative logarithmic coefficients.

\begin{proof}

Write a change of the free coefficients as $v=(v^A,v^B)$. Let $S(v)$ have diagonal entries $v_i^A$ and entries $v_i^B$ in positions $(i,i+1)$ and $(i+1,i)$, with zero entries for
the coordinates held fixed. Then $T(y+tv,\tau)=T(y,\tau)+tS(v)$. We first justify the trace Hessian bound; scalar convexity of $V$ is sufficient here, without operator convexity.
Diagonalize $T=U\operatorname{diag}(\lambda_1,\ldots,\lambda_n)U^*$ and put $\widehat S=U^*S(v)U$. For a monomial $V(t)=t^m$, $m\ge2$, cyclicity of the trace gives
\[
 D^2\operatorname{Tr}(T^m)[S,S]
 =m\sum_{\ell=0}^{m-2}\operatorname{Tr}(T^\ell S T^{m-2-\ell}S).
\]
The coefficient of $|\widehat S_{ij}|^2$ in this expression is $m\sum_{\ell=0}^{m-2}\lambda_i^\ell\lambda_j^{m-2-\ell}$. It equals the divided difference of $V'$ for unequal
eigenvalues and $V''(\lambda_i)$ when they coincide. Summing over the monomials yields
\[
 D^2\operatorname{Tr}V(T)[S,S]
 =\sum_{i,j}V'^{[1]}(\lambda_i,\lambda_j)|\widehat S_{ij}|^2,
 \eqtag{finite-trace-hessian}
\]
where $V'^{[1]}(x,y)=\int_0^1V''((1-t)x+ty)\,\mathrm dt\ge c_V$. Since $S=S(v)$ is supported on the free matrix entries,
\[
 D^2\operatorname{Tr}V(T)[S(v),S(v)]
 \ge c_V\operatorname{Tr}S(v)^2
 =c_V\left(\sum_i(v_i^A)^2+2\sum_i(v_i^B)^2\right).
 \eqtag{finite-trace-convexity}
\]
Each retained logarithm has nonnegative Hessian. Multiplication by \(\beta n/2\) proves \eqref{A:eq:12.15a}, without any restriction on \(\tau\). The polynomial part tends to
infinity when the free coefficients tend to infinity. For every free off-diagonal coefficient, the term $-2w_{k,n}\log b_k$ tends to infinity as $b_k\downarrow0$ because
$w_{k,n}>0$. A minimizing sequence therefore stays bounded and away from zero in every positive coordinate. It has a convergent subsequence in the domain, and strict convexity
makes its minimizer unique. If only diagonal coordinates are free, there is no positivity boundary to consider.

Apply Theorem~12 of \cite{krishnapur2016universality} to the potential $V/2$. The energies have the same minimizer, since
\[
 \frac12\mathscr H_n
 =\operatorname{Tr}(V/2)(T)-\sum_k w_{k,n}\log b_k,
 \qquad 0<w_{k,n}\le1.
\]
In that theorem take $y_2=y_3=1$ and $y_1=1+\|\tau\|_{\infty,\partial J}+|a_{\mathrm t}|$, where $a_{\mathrm t}$ minimizes $V$. All fixed coefficients entering a differentiated
trace term lie within $R_V$ indices of $J$ and satisfy its boundary bounds. There is no negative logarithmic coefficient among the free variables, so its lower cutoff imposes no
restriction. With $d=\deg V$, the theorem gives
\[
 \|y_J^\circ(\tau)\|_\infty
 \le |a_{\mathrm t}|+C_V(1+y_1)^{d^2}
 \le C'_V(1+\|\tau\|_{\infty,\partial J})^{d^2}.
\]
Thus $D_V=d^2$ is admissible in \eqref{A:eq:mode-growth}. The conditional polynomial contains the original terms involving coefficients on both sides of each cut.

On \(\mathcal T_{n,J}(M)\), the coefficients of the conditional minimizer and the exterior coefficients within $R_V$ indices of the block lie in a fixed bounded cube. For a free
edge \(k\), fix all other coordinates at their minimizing values and write \(g_k(t)=\operatorname{Tr}V(T^{(k,0)}+tS_k)\). To see that $g_k$ is even, conjugate by the diagonal
matrix which is $+1$ on indices at most $k$ and $-1$ on the remaining indices. This changes only the sign of the entry at edge $k$ and leaves the trace of every polynomial
unchanged. Let $I_k=\{k-R_V,\ldots,k+R_V\}\cap\{1,\ldots,n\}$ and let $\widehat T_k(t)$ be the principal matrix of $T^{(k,0)}+tS_k$ on $I_k$. Every differentiated monomial in
$g_k''(t)$ contains two entries of $S_k$; its path has length at most $R_V=\deg V$ and lies in $I_k$. Therefore
\[
 g_k''(t)=D^2\operatorname{Tr}V(\widehat T_k(t))[S_k,S_k].
\]
The free coefficients and the exterior coefficients in this interval are bounded on the stated exterior event, so
$\sup_{0\le t\le b_k^\circ}\|\widehat T_k(t)\|_{\mathrm{op}}\le C_M$. Formula~\eqref{A:eq:finite-trace-hessian} and $\|S_k\|_{\mathfrak S_2}^2=2$ give
$0\le g_k''(t)\le2\sup_{|x|\le C_M}|V''(x)|$. Since $g_k'(0)=0$,
\[
 0<g_k'(b_k^\circ)=\int_0^{b_k^\circ}g_k''(t)\,\mathrm dt
 \le C_M b_k^\circ.
\]
The stationary equation gives
\[
 \frac{2w_{k,n}}{b_k^\circ}=g_k'(b_k^\circ),
 \quad
 (b_k^\circ)^2\ge\frac{2w_{k,n}}{C_M}
 \ge\frac{1-s_1}{C_M}
 \quad(k/n\le s_1).
\eqtag{regular-mode-lower}
\]
The last inequality holds for \(n\ge2/[\beta(1-s_1)]\). The bound on \(g_k''\) uses the bounded exterior data; its constant is not uniform over arbitrary exterior values.

Lemma~\ref{A:res:A.8}, in dimension \(d_J\), gives \eqref{A:eq:conditional-radius}. Since \(d_J\le2n\), choose a fixed radius larger than \(C\sqrt2+1\) and use
\eqref{A:eq:mode-growth} to obtain \eqref{A:eq:conditional-upper}. The curvature constants can be chosen using \(\beta_-\) for parameters in \([\beta_-,\beta_+]\), and the
remaining constants are uniform on the same interval.

\end{proof}

\subsubsection{The local minimizer and edge scaling}\label{A:sec:the-local-minimizer-and-edge-scaling}

The conditional minimizer depends on the exterior data and solves a finite system of stationary equations. The functions $a^\dagger(s),b^\dagger(s)$ from the introduction instead
minimize the two variable energy $F_{V,s}$ in \eqref{A:eq:local-energy}. They describe the interior coefficients when the coefficient of the logarithmic term is held constant. We
first establish their regularity and the resulting edge expansion, then estimate the gradient of the finite energy evaluated at these coefficients.

Write $h=h_{n,V}$. The off-diagonal parameter must match the coefficient of $\log b_k$ in \eqref{A:eq:11.1}. We place the diagonal parameter halfway between the two adjacent
off-diagonal parameters:
\[
 s^a_{n,k}=\frac{k-\tfrac12}{n}+\frac1{\beta n},\qquad
 s^b_{n,k}=\frac kn+\frac1{\beta n},\qquad
 a^\dagger_{n,k}=a^\dagger(s^a_{n,k}),\quad
 b^\dagger_{n,k}=b^\dagger(s^b_{n,k}).
 \eqtag{local-parameters}
\]
The last two expressions are used only for parameters in $[0,1)$. The finite coefficient law supplies the term $1/(\beta n)$. Its contribution after edge scaling will be the scalar
shift $h/\beta$ in Proposition~\ref{A:res:11.6}. This choice of the diagonal parameter cancels the first-order variation in the sum of the two adjacent $b^\dagger$ coefficients.

The Gaussian case gives an explicit comparison of the mean and the mode. For $V(t)=t^2/2$, one has $a^\dagger(s)=0$ and $b^\dagger(s)=\sqrt{1-s}$. Hence $b^\dagger_{n,k}$ is the
mode of the scaled chi density when $\beta(n-k)>1$. The chi mean expansion in Lemma~\ref{A:res:A.1} gives, uniformly for $1\le k\le s_0n$, with $s_0<1$,
\[
 \overline b_{k,n}-b^\dagger_{n,k}
 =\frac{1}{4\beta n\sqrt{1-k/n}}+O_{\beta,s_0}(n^{-2}).
\]
The conditional Gaussian approximation is centered at the mode, with the mean displacement given above.

The Fourier average $W_V$ in \eqref{A:eq:11.4} is the zeroth Laurent coefficient of $V(a+b(z+z^{-1}))$. Its derivatives determine the minimizer of \eqref{A:eq:local-energy}.

\begin{lemma}\label{A:res:11.1}
For $0\le s<1$, the Hessian of $F_{V,s}$ is bounded below by $c_V\operatorname{diag}(1,2)$. In particular, $F_{V,s}$ is strictly convex.
\end{lemma}

\begin{proof}
Apply the trace convexity bound in \cite{krishnapur2016universality} to Jacobi matrices with constant coefficients. In the normalization of \eqref{A:eq:11.4}, its two contributions
are
\[
 D^2W_V(a,b)[(\alpha,\gamma),(\alpha,\gamma)]
 =\frac1{2\pi}\int_0^{2\pi}V''(a+2b\cos\theta)
       (\alpha+2\gamma\cos\theta)^2\,\mathrm d\theta
 \ge c_V(\alpha^2+2\gamma^2),
 \eqtag{11.6}
\]
\[
 D^2[-2(1-s)\log b][(\alpha,\gamma),(\alpha,\gamma)]
 =\frac{2(1-s)}{b^2}\gamma^2\ge0.
 \eqtag{11.7}
\]
Their sum is the claimed Hessian bound.
\end{proof}

\par\Needspace{5\baselineskip}
\begin{lemma}\label{A:res:11.2}
For every \(s\in[0,1)\), there exists a unique minimizer \((a^\dagger(s),b^\dagger(s)) \in\mathbb R\times(0,\infty)\) of
\(F_{V,s}\).
\end{lemma}

\begin{proof}
Let $a_{\mathrm t}$ be the unique minimizer of $V$. Existence follows from $V''\ge c_V>0$, which implies quadratic growth at infinity; uniqueness follows from strict convexity.
Taylor's formula with integral remainder gives
\[
 V(x)\ge V(a_{\mathrm t})+\frac{c_V}{2}(x-a_{\mathrm t})^2.
\]
Averaging at $x=a+2b\cos\theta$, and using the averages of $\cos\theta$ and $\cos^2\theta$, yields
\[
 W_V(a,b)\ge V(a_{\mathrm t})+
 \frac{c_V}{2}(a-a_{\mathrm t})^2+c_Vb^2.
\]
For fixed $s<1$, it follows that
\[
 F_{V,s}(a,b)\ge V(a_{\mathrm t})+
 \frac{c_V}{2}(a-a_{\mathrm t})^2+
 c_Vb^2-2(1-s)\log b.
\]
The last two terms tend to infinity both as $b\downarrow0$ and as $b\to\infty$, and have a finite lower bound on $(0,\infty)$. Thus every nonempty sublevel set of $F_{V,s}$ has $a$
bounded and $b$ bounded above and away from zero. It is closed in a compact rectangle contained in $\mathbb R\times(0,\infty)$, and hence compact. The value at $(a_{\mathrm t},1)$
is finite, so a minimizing sequence lies in one such sublevel set and has a limit there. Continuity gives a minimizer with positive second coordinate. The strict convexity in
Lemma~\ref{A:res:11.1} excludes two distinct minimizers.
\end{proof}

Differentiating $F_{V,s}$ at its minimizer gives the stationary equations
\[
\frac{1}{2\pi}
\int_0^{2\pi}
V'(a^\dagger(s)+2b^\dagger(s)\cos\theta)\,\mathrm{d}\theta
=0,
\]
and
\[
\frac{1}{2\pi}
\int_0^{2\pi}
2\cos\theta\,
V'(a^\dagger(s)+2b^\dagger(s)\cos\theta)\,\mathrm{d}\theta
=
\frac{2(1-s)}{b^\dagger(s)}.
\]

\begin{proposition}\label{A:res:11.3}
For every \(s_0<1\), \(s\mapsto(a^\dagger(s),b^\dagger(s))\) is real analytic on \([0,s_0]\). Moreover, for each
\(m\ge0\), there exists \(C_{m,s_0,V}<\infty\) such that \(\sup_{0\le s\le s_0} \left( |(a^\dagger)^{(m)}(s)| +
|(b^\dagger)^{(m)}(s)| \right) \le C_{m,s_0,V}\). Also, \(0<c_{s_0,V}\le b^\dagger(s)\le C_{s_0,V}<\infty \allowbreak\;
(0\le s\le s_0)\).
\end{proposition}

\begin{proof}
Define, for $s<1$ and $b>0$,
\[
 G(s,a,b)=
 \begin{pmatrix}\partial_aW_V(a,b)\\
 \partial_bW_V(a,b)-2(1-s)/b\end{pmatrix}.
\]
Since $V$ is a polynomial, $W_V$ is a polynomial in $a,b$, and $G$ is real analytic on this open set. The stationary equations give $G(s,a^\dagger(s),b^\dagger(s))=0$. Its
derivative in $(a,b)$ is
\[
 D_{(a,b)}G(s,a^\dagger,b^\dagger)
 =D^2F_{V,s}(a^\dagger,b^\dagger)
 \ge c_V\operatorname{diag}(1,2),
\]
so this derivative is invertible. The analytic implicit function theorem gives an analytic solution near every fixed $s_*<1$. The argument for existence and uniqueness in the
preceding lemma also applies to negative $s$, because $1-s>0$. The solution near $s_*$ has $b>0$ and is a stationary point of a strictly convex function, so it agrees with its
unique minimizer. These local solutions consequently agree on overlaps and define the required analytic function on a neighborhood of every compact interval
$[0,s_0]\subset(-\infty,1)$.

Differentiating the stationary equations gives
\[
 \begin{pmatrix}(a^\dagger)'(s)\\(b^\dagger)'(s)\end{pmatrix}
 =-\bigl[D^2F_{V,s}(a^\dagger(s),b^\dagger(s))\bigr]^{-1}
       \begin{pmatrix}0\\2/b^\dagger(s)\end{pmatrix}.
\]
Every derivative of an analytic local solution is continuous. A finite collection of the preceding neighborhoods covers $[0,s_0]$, so each fixed derivative is bounded on that
interval. Finally, the positive continuous function $b^\dagger$ attains its minimum and maximum there; positivity at every point makes the minimum strictly positive. This proves
all three uniform bounds.
\end{proof}

By the Fourier representation \eqref{A:eq:11.4}, the Jacobi matrix with constant coefficients $a^\dagger(s),b^\dagger(s)$ has spectrum $[a^\dagger(s)-2b^\dagger(s),E(s)]$. The
constants in \eqref{A:eq:edge-data} give $E(s)=E_\star-\tau_\star s+O(s^2)$.

The normalization can be checked on a fixed $f\in C_c^2((0,\infty))$. Set $u_k=h^{1/2}f(kh)$. Riemann sums and the fundamental theorem of calculus give
\[
 \sum_k|u_k|^2\longrightarrow\|f\|_2^2,\qquad
 h^{-2}\sum_k|u_{k+1}-u_k|^2\longrightarrow\|f'\|_2^2.
\]
On the support of $f$, $k=O(h^{-1})$ and $E_\star-E(k/n)=\tau_\star k/n+O(k^2/n^2)$. Consequently
\[
 \frac{h^{-2}}{b_\star}\sum_k(E_\star-E(k/n))|u_k|^2
 =\frac{\tau_\star}{b_\star nh^3}
       \sum_k kh|u_k|^2+O(h^{-4}n^{-2}).
\]
Choosing $b_\star nh^3=\tau_\star$ makes the first term converge to $\int x|f(x)|^2\,\mathrm dx$, and the remainder is $O(h^2)$. On the same compact support,
$b^\dagger_{n,k}/b_\star=1+O(h^2)$, so
\[
 \sum_k\frac{b^\dagger_{n,k}}{b_\star}|D_hu_k|^2
 =\sum_k|D_hu_k|^2+O(h^2)\longrightarrow\|f'\|_2^2.
\]
These two calculations fix the deterministic normalization. The variance of the limiting noise is identified by the covariance calculation in
Section~\ref{A:sec:the-general-jacobi-model-gaussian-approximation-and-comparison}. We now retain the exact vertex coordinates and finite logarithmic correction needed for the
quantitative comparison.

At the two boundary point indices, \eqref{A:eq:local-parameters} gives $s^b_{n,0}=1/(\beta n)$ and $s^b_{n,n}=1+1/(\beta n)$; no local coefficient is evaluated at a parameter
greater than one. Lemma~\ref{A:res:profile-boundary} treats the limit as $s\uparrow1$.

Since $n^{-1}=(b_\star/\tau_\star)h^3$ and $x_k=(k-\tfrac12)h$, the diagonal parameter is
\[
 s^a_{n,k}=\frac{b_\star}{\tau_\star}h^2x_k
             +\frac{b_\star}{\beta\tau_\star}h^3.
 \eqtag{11.30}
\]
The adjacent edge parameters satisfy
\[
 s^b_{n,k}=s^a_{n,k}+\frac1{2n},\qquad
 s^b_{n,k-1}=s^a_{n,k}-\frac1{2n}.
 \eqtag{11.33}
\]
Equivalently, replace $x_k$ in \eqref{A:eq:11.30} by $x_k+h/2$ or $x_k-h/2$. For each fixed $s_0<1$, we take $n$ sufficiently large that all adjacent parameters used below lie in a
fixed larger compact subinterval of $[0,1)$. Let
\(a_1=(a^\dagger)'(0), \allowbreak\, b_1=(b^\dagger)'(0).\) Then \(a_1+2b_1=-\tau_\star.\) For \(s\) small,
\(a^\dagger(s) = a_\star+a_1s+\frac12a_2s^2+O(s^3),\) and \(b^\dagger(s) = b_\star+b_1s+\frac12b_2s^2+O(s^3).\) The
local edge function \(E(s)=a^\dagger(s)+2b^\dagger(s)\) therefore satisfies \(E(s) = E_\star-\tau_\star s+\frac12E_2s^2+O(s^3),\) where \(E_2=a_2+2b_2.\) We shall use the following
uniform version of this expansion.

\begin{lemma}\label{A:res:11.4}\label{A:res:11.5}
Fix \(s_0<1\). There exists \(C=C(V,\beta,s_0)\) such that, whenever \(s^a_{n,k}\le s_0,\) we have
\[
|a^\dagger_{n,k}-a_\star|
+
|b^\dagger_{n,k}-b_\star|
\le
C h^2(1+x_k).
\eqtag{11.38}
\]
Moreover,
\[
|b^\dagger_{n,k}-b^\dagger_{n,k-1}|
\le
C h^3.
\eqtag{11.39}
\]
Finally,
\[
\left|
b^\dagger_{n,k}-b_\star
-
b_1\frac{b_\star}{\tau_\star}h^2x_k
\right|
\le
C\left(h^3+h^4(1+x_k)^2\right).
\eqtag{11.40}
\]
The analogous estimate holds for \(a^\dagger_{n,k}\).
\end{lemma}

\begin{proof}
Choose $s_1$ with $s_0<s_1<1$ and take $n$ large enough that $s_0+(2n)^{-1}\le s_1$. All the vertex and adjacent edge parameters below then lie in $[0,s_1]$. By \eqref{A:eq:11.30},
\(s^a_{n,k} = \frac{b_\star}{\tau_\star}h^2x_k + O(h^3).\) Thus
\[
s^a_{n,k}\le C h^2(1+x_k).
\eqtag{11.41}
\]
Since \(a^\dagger\) and \(b^\dagger\) have bounded first derivatives on \([0,s_1]\),
\(|a^\dagger(s)-a_\star|+|b^\dagger(s)-b_\star| \le Cs.\) Using \eqref{A:eq:11.41} gives \eqref{A:eq:11.38}. For the
difference, \(s^b_{n,k}-s^b_{n,k-1}=n^{-1}=O(h^3)\). The bounded derivative of \(b^\dagger\) then gives
\(|b^\dagger_{n,k}-b^\dagger_{n,k-1}|\le Ch^3\), which is \eqref{A:eq:11.39}. For \eqref{A:eq:11.40}, Taylor expand:
\(b^\dagger(s)=b_\star+b_1s+O(s^2).\) Substitute \(s=s^b_{n,k} = \frac{b_\star}{\tau_\star}h^2x_k + O(h^3).\) Then
\(b^\dagger_{n,k} = b_\star + b_1\frac{b_\star}{\tau_\star}h^2x_k + O(h^3) + O(s^2)\). Since \(s^2\le C h^4(1+x_k)^2,\)
we obtain \eqref{A:eq:11.40}. For the diagonal coefficient, use
\[
 a^\dagger(s)-a_\star-a_1s
 =\int_0^s(s-t)(a^\dagger)''(t)\,\mathrm dt.
\]
At $s=s^a_{n,k}$, \eqref{A:eq:11.30} and boundedness of the second derivative give
\[
 \left|a^\dagger_{n,k}-a_\star
       -a_1\frac{b_\star}{\tau_\star}h^2x_k\right|
 \le C\{h^3+h^4(1+x_k)^2\}.
\]
\end{proof}

The coefficient of $|D_hu_k|^2$ in the scaled deterministic KRV operator is \(c^\dagger_{n,k} = \frac{b^\dagger_{n,k}}{b_\star}.\) The coefficient of $|D_hu_k|^2$ in
$-\Delta_{h,D}+X_h$ is $1$. Define the difference between these coefficients \(d^\dagger_{n,k} = c^\dagger_{n,k}-1 = \frac{b^\dagger_{n,k}-b_\star}{b_\star}\).

If $s^b_{n,k}\le s_0$, then $s^a_{n,k}\le s_0$ by \eqref{A:eq:11.33}. Apply Lemma~\ref{A:res:11.4} and divide its bounds \eqref{A:eq:11.38}--\eqref{A:eq:11.39} by $b_\star>0$. The
vertex and predecessor edge parameters remain in the enlarged compact interval used in that lemma. With $C=C(V,\beta,s_0)$, this gives
\[
|d^\dagger_{n,k}|
\le
C h^2(1+x_k),
\eqtag{11.44}
\]
and
\[
|d^\dagger_{n,k}-d^\dagger_{n,k-1}|
\le
C h^3.
\eqtag{11.45}
\]

The scaled deterministic potential at vertex \(k\) is
\[
V^\dagger_{n,k}
=
\frac{h^{-2}}{b_\star}
\left[
E_\star
-
a^\dagger(s^a_{n,k})
-
b^\dagger(s^b_{n,k-1})
-
b^\dagger(s^b_{n,k})
\right].
 \eqtag{local-potential}
\]
The center identity \eqref{A:eq:11.33} makes the linear terms in the two adjacent $b^\dagger$ values cancel. Their sum therefore differs from $2b^\dagger(s^a_{n,k})$ by
$O(n^{-2})$. The remaining linear term comes from $E_\star-E(s^a_{n,k})$ and includes $h/\beta$ after scaling. The following estimate keeps the dependence of the quadratic
remainder on the coefficient index on the growing interval; only on bounded intervals does it reduce to
$V_{n,k}^\dagger=x_k+O(h)$. \par\Needspace{5\baselineskip}

\begin{proposition}\label{A:res:11.6}
Fix \(s_0<1\). There exists \(C=C(V,\beta,s_0)\) such that, whenever \(s^a_{n,k}\le s_0\),
\[
\left|
V^\dagger_{n,k}
-
x_k
-
\frac{h}{\beta}
\right|
\le
C h^2(1+x_k)^2.
\eqtag{11.47}
\]
In particular,
\[
|V^\dagger_{n,k}-x_k|
\le
C\left(h+h^2(1+x_k)^2\right).
\eqtag{11.48}
\]
\end{proposition}

\begin{proof}
Let \(s_k=s^a_{n,k}.\) By \eqref{A:eq:11.33}, \(s^b_{n,k}=s_k+\frac{1}{2n}, \allowbreak\, s^b_{n,k-1}=s_k-\frac{1}{2n}\). Choose $s_1\in(s_0,1)$ and $n$ large enough that the
adjacent parameters lie in $[0,s_1]$. Put $\epsilon=(2n)^{-1}$. The integral remainder gives
\[
 b^\dagger(s_k+\epsilon)+b^\dagger(s_k-\epsilon)-2b^\dagger(s_k)
 =\int_0^\epsilon(\epsilon-t)
       \big((b^\dagger)''(s_k+t)+(b^\dagger)''(s_k-t)\big)\,\mathrm dt.
\]
Its absolute value is at most $C\epsilon^2=Cn^{-2}$. Hence
\[
 a^\dagger(s^a_{n,k})+b^\dagger(s^b_{n,k-1})+b^\dagger(s^b_{n,k})
 =a^\dagger(s_k)+2b^\dagger(s_k)+O(n^{-2}).
\]
The first two terms are \(E(s_k)=a^\dagger(s_k)+2b^\dagger(s_k).\) By the edge expansion, \(E(s_k)=E_\star-\tau_\star s_k+O(s_k^2).\) Therefore,
\(E_\star - a^\dagger(s^a_{n,k}) - b^\dagger(s^b_{n,k-1}) - b^\dagger(s^b_{n,k}) = \tau_\star s_k+O(s_k^2)+O(n^{-2})\). Multiply by \(h^{-2}/b_\star\):
\[
V^\dagger_{n,k}
=
\frac{\tau_\star}{b_\star}h^{-2}s_k
+
O\left(h^{-2}s_k^2\right)
+
O\left(h^{-2}n^{-2}\right).
\eqtag{11.53}
\]
The coordinate identity \eqref{A:eq:11.30} gives \(s_k = \frac{b_\star}{\tau_\star}h^2x_k + \frac{b_\star}{\beta\tau_\star}h^3\). Hence
\[
\frac{\tau_\star}{b_\star}h^{-2}s_k
=
x_k+\frac{h}{\beta}.
\eqtag{11.54}
\]
Next, \(s_k\le C h^2(1+x_k),\) so
\[
h^{-2}s_k^2
\le
C h^2(1+x_k)^2.
\eqtag{11.55}
\]
Finally, $n^{-1}\le C h^3$, $0<h\le1$ and $x_k\ge0$ imply
\[
 h^{-2}n^{-2}\le C^2h^4\le C^2h^2(1+x_k)^2.
\] Combining \eqref{A:eq:11.53}--\eqref{A:eq:11.55} proves \eqref{A:eq:11.47}. For \eqref{A:eq:11.48}, the triangle inequality gives
$|V^\dagger_{n,k}-x_k| \le|V^\dagger_{n,k}-x_k-h/\beta|+h/\beta \le C(h+h^2(1+x_k)^2)$.
\end{proof}

On a Dirichlet interval whose vertex and adjacent edge parameters belong to $[0,s_0]$, $s_0<1$, define $\mathcal A_{n,V}^{\dagger}$ by the quadratic form
\[
\mathfrak a_{n,V}^{\dagger}[u] ={} \sum_{k} c^\dagger_{n,k}|D_hu_k|^2 + \sum_{k} V^\dagger_{n,k}|u_k|^2,
\]

where \(c^\dagger_{n,k}=\frac{b^\dagger_{n,k}}{b_\star},\) and $V^\dagger_{n,k}$ is given by \eqref{A:eq:local-potential}. The deterministic finite difference operator with linear
potential is \(\mathcal A_h^{\mathrm{det}} = -\Delta_{h,D}+X_h,\) with form
\[
\mathfrak a_h^{\mathrm{det}}[u]
=
\sum_k|D_hu_k|^2+\sum_k x_k|u_k|^2.
\]
Using the previous estimates, \(\mathcal A_{n,V}^{\dagger} - \mathcal A_h^{\mathrm{det}}\) has the form
\(D_h^*M_{d^\dagger_n}D_h + M_{r^\dagger_n} + \frac{h}{\beta}I,\) where \(d^\dagger_{n,k}=c^\dagger_{n,k}-1,\) and
\(r^\dagger_{n,k} = V^\dagger_{n,k} - x_k - \frac{h}{\beta}.\) The estimates are \(|d^\dagger_{n,k}| \le C h^2(1+x_k),\)
\(|d^\dagger_{n,k}-d^\dagger_{n,k-1}| \le C h^3,\) and \(|r^\dagger_{n,k}| \le C h^2(1+x_k)^2.\) These are the bounds
required by the product rule and weighted resolvent estimates of Subsection~\ref{A:sec:weighted-resolvent-estimates-and-comparison-principles}. They are not uniform smallness
estimates on an interval of length $O(h^{-2})$: the weights retain the growth of the coefficient errors. The scalar shift
\(\frac{h}{\beta}I\) does not affect the stated estimate. Indeed, for any operator \(H\) satisfying the upper and lower
form bounds, \((H+h\beta^{-1}-z)^{-1}-(H-z)^{-1} = -\frac{h}{\beta} (H+h\beta^{-1}-z)^{-1}(H-z)^{-1}\).
\[
\left\|
(H+h\beta^{-1}-z)^{-1}-(H-z)^{-1}
\right\|_{\mathfrak S_2}^2
\le
C_{\omega,K}h^2
\| (H-z)^{-1}\|_{\mathfrak S_2}^2
\le
C_{\omega,K}h^2,
\]
which is \(o(h)\). Thus the scalar shift is below the full \(O(h)\) squared Hilbert--Schmidt scale.

We next check how well this profile solves the finite stationary equations. Away from the boundary points, the same finite set of paths appears in the finite and constant
coefficient derivatives. The local parameters vary by $O(n^{-1})$ over those paths, so the gradient evaluated at the local profile is $O(n^{-1})$ away from the boundary points. We
use the following identities for the derivatives of the trace. Write \(V'(x)=\sum_{\ell=0}^{d-1}c_\ell x^\ell\), where \(d=\deg V\).

For $|u-v|\le1$, let $\mathcal N_d(u,v)$ be the coefficients visited by paths of length at most $d-1$ from $u$ to $v$. In comparing two Jacobi matrices we assume that their
coefficients on this set have modulus at most $M$. The finite matrices have the same index set; absent paths at a boundary are treated separately. Derivatives in $a_k$ and $b_k$
below are taken for $1\le k\le n$ and $1\le k<n$, respectively.

\begin{lemma}\label{A:res:11.6A}
For the finite Jacobi matrix $T$ and the Jacobi operator $J(a,b)$ with constant coefficients, respectively,
\[
 \begin{aligned}
 \partial_{a_k}\operatorname{Tr}V(T)&=[V'(T)]_{kk},&
 \partial_{b_k}\operatorname{Tr}V(T)&=2[V'(T)]_{k,k+1},\\
 \partial_aW_V(a,b)&=[V'(J(a,b))]_{00},&
 \partial_bW_V(a,b)&=2[V'(J(a,b))]_{01}.
 \end{aligned}
\]
Under the local bound just specified,
\[
 |[V'(T)]_{uv}-[V'(\widetilde T)]_{uv}|
 \le C_{V,M}\max_{\xi\in\mathcal N_d(u,v)}|y_\xi-\widetilde y_\xi|.
\]
\end{lemma}

\begin{proof}
For a matrix direction \(\dot T\), cyclicity of the trace gives, for every \(\ell\ge1\),
\[
\frac{\mathrm{d}}{\mathrm{d}t}\operatorname{Tr}(T+t\dot T)^\ell\bigg|_{t=0}
=\sum_{r=0}^{\ell-1}\operatorname{Tr}
\bigl(T^r\dot T T^{\ell-1-r}\bigr)
=\ell\operatorname{Tr}(T^{\ell-1}\dot T).
\]
Summing over the monomials of \(V\) yields \(D\operatorname{Tr}V(T)[\dot T]=\operatorname{Tr}(V'(T)\dot T)\). Taking
\(\dot T=E_{kk}\) proves the first identity, while taking \(\dot T=E_{k,k+1}+E_{k+1,k}\) proves the second because
\(V'(T)\) is symmetric.

Under the Fourier transform, \(J(a,b)\) is multiplication by \(a+2b\cos\theta\). Therefore
\[
[V'(J(a,b))]_{00}
=\frac1{2\pi}\int_0^{2\pi}V'(a+2b\cos\theta)\,\mathrm{d}\theta
=\partial_aW_V(a,b),
\]
and, since the integrand is real and even,
\[
2[V'(J(a,b))]_{01}
=\frac1{2\pi}\int_0^{2\pi}2\cos\theta\,
V'(a+2b\cos\theta)\,\mathrm{d}\theta
=\partial_bW_V(a,b).
\]

Finally, the matrix entry of a power has the path expansion
\[
(T^\ell)_{uv}
=\sum_{\substack{\gamma_0=u,\ \gamma_\ell=v\\
|\gamma_{r+1}-\gamma_r|\le1}}
\prod_{r=0}^{\ell-1}T_{\gamma_r,\gamma_{r+1}}.
\]
Every path has at most \(\ell\le d-1\) factors and there are at most \(3^\ell\) such paths. For two products whose factors have modulus at most \(M\), telescoping gives
\[
\left|\prod_{r=1}^{\ell}x_r-\prod_{r=1}^{\ell}\widetilde x_r\right|
\le \ell M^{\ell-1}\max_r|x_r-\widetilde x_r|.
\]
Applying this estimate to every path and then summing over \(\ell=0,\ldots,d-1\) proves the local Lipschitz bound.
\end{proof}

Use the local minimizing values only at indices whose parameters are less than one. Complete the other finitely many terminal coordinates by any bounded real diagonal and positive
off-diagonal values, for example $a_{\mathrm t}$ and $b_\star$. All derivatives considered below depend only on neighborhoods where the local values are defined, so they are
independent of this completion. Let \(T_n^\dagger=T(a_n^\dagger,b_n^\dagger)\) be the resulting finite Jacobi matrix formed from the local minimizer. The next proposition compares
the finite derivatives of
\(\operatorname{Tr}V(T_n^\dagger)\) with the exact constant coefficient derivatives in Lemma~\ref{A:res:11.6A}.

Fix \(s_0<1\), and put \(d=\deg V\). Define
\[
\mathcal I_{n}^{a}(s_0)
=\{k:d\le k\le n-d+1,\ s^a_{n,k}\le s_0\},
\quad
\mathcal I_{n}^{b}(s_0)
=\{k:d\le k\le n-d,\ s^b_{n,k}\le s_0\}.
\]

\begin{proposition}\label{A:res:11.7}
For all sufficiently large $n$,
\[
\max_{k\in\mathcal I_n^a(s_0)}
\left|\partial_{a_k}\mathscr H_n(\mathbf y_n^\dagger)\right|
+
\max_{k\in\mathcal I_n^b(s_0)}
\left|\partial_{b_k}\mathscr H_n(\mathbf y_n^\dagger)\right|
\le \frac{C_{V,\beta,s_0}}{n}.
\eqtag{11.63}
\]
\end{proposition}

\begin{proof}
Choose $s_1\in(s_0,1)$. Proposition~\ref{A:res:11.3} bounds the local minimizer and its first derivative on $[0,s_1]$. Every coefficient in a path of length at most $\deg V-1$
through an admissible index has parameter in this interval and differs from its value at the starting parameter by $O(n^{-1})$. The finite path expansion and local Lipschitz bound
in Lemma~\ref{A:res:11.6A} therefore give
\[
 \left|[V'(T_n^\dagger)]_{kk}
       -[V'(J(a^\dagger(s),b^\dagger(s)))]_{00}\right|
 \le Cn^{-1},\qquad s=s^a_{n,k},
\]
and
\[
 \left|2[V'(T_n^\dagger)]_{k,k+1}
       -2[V'(J(a^\dagger(s),b^\dagger(s)))]_{01}\right|
 \le Cn^{-1},\qquad s=s^b_{n,k}.
\]
The admissible index sets keep these paths inside the finite matrix, so the constants are uniform in $k$.

The local Euler--Lagrange equations imply
\[
 [V'(J(a^\dagger(s),b^\dagger(s)))]_{00}=0,
 \qquad
 2[V'(J(a^\dagger(s),b^\dagger(s)))]_{01}
 =\frac{2(1-s)}{b^\dagger(s)}.
\]
The logarithmic energy is independent of $a_k$. At $s=s^b_{n,k}=k/n+1/(\beta n)$, its derivative in $b_k$ is
\[
 -\frac{2(1-k/n-1/(\beta n))}{b^\dagger(s^b_{n,k})}
 =-\frac{2(1-s)}{b^\dagger(s)}.
\]
The constant terms thus cancel in both coordinates, leaving the two $O(n^{-1})$ errors above. This proves \eqref{A:eq:11.63}.
\end{proof}

Let $\Lambda$ be a consecutive block of coefficient variables and

let $\partial_R\Lambda$ be its actual exterior coefficient variables within site distance $R$. This is the set denoted by $\partial J$ when $J=\Lambda$ in
\eqref{A:eq:cond-density}. Let $\mathscr H_{n,\Lambda}^\tau$ be that conditional energy in the normalization of \eqref{A:eq:11.1}. The energy includes the terms involving both free
and fixed coordinates, and its minimizer is $y_\Lambda^\circ(\tau)$. The distance to $\partial\Lambda$ means distance to the ends of the free coefficient interval in these site
coordinates. A physical boundary point of the finite matrix is not an extra conditioned coefficient. Its effect is the absence of paths crossing that boundary point; this changes
the gradient terms supported near the boundary.

Fix \(s_0<1\), \(M_0<\infty\), and \(m_0>0\). Assume that every macroscopic parameter attached to
\(\Lambda\cup\partial_R\Lambda\) is at most \(s_0\), and that the conditioned values satisfy \(|\tau_{a,k}|\le M_0\) and
\(m_0\le\tau_{b,k}\le M_0\). Let \(\mathbf y_{n,\Lambda}^{\star,\tau}\) be the minimizer of $\mathscr
H_{n,\Lambda}^{\tau}$; it is the same vector as $y_\Lambda^\circ(\tau)$, with the dependence on $n$ and the comparison
with the local profile made explicit.

\begin{lemma}\label{A:res:11.7A}
There are constants \(0<m<M<\infty\), depending only on \(V,\beta,s_0,M_0,m_0\), such that
\[
|a_k^{\star,\tau}|\le M,
\quad
m\le b_k^{\star,\tau}\le M
\]
for every free coordinate. The same bounds hold, after changing \(m,M\), on the whole segment joining \(\mathbf y_n^\dagger|_\Lambda\) to \(\mathbf y_{n,\Lambda}^{\star,\tau}\).
Along this segment the conditional Hessian is banded with fixed bandwidth and satisfies
\[
cI\le D^2\mathscr H_{n,\Lambda}^{\tau}\le CI.
\]
\end{lemma}

\begin{proof}
The conditional domain is a product of real lines and positive half-lines, and its logarithmic weights are positive. Existence and uniqueness, and the uniform upper bound under the
stated exterior hypotheses, follow from Lemma~\ref{A:res:12.1A}: for large $n$, every interval with parameters at most $s_0<1$ lies in the retained initial segment. Choose
$s_1\in(s_0,1)$ and take $n$ large enough that every parameter used here is at most $s_1$. Then $1-s^b_{n,k}\ge1-s_1$. The fixed exterior bound $M_0$ gives a uniform bound on the
minimizer for every number of free coordinates.

A lower bound for each free off-diagonal minimizer follows directly from its Euler--Lagrange equation. By Lemma~\ref{A:res:11.6A},
\[
2[V'(T(\mathbf y_{n,\Lambda}^{\star,\tau},\tau))]_{k,k+1}
=\frac{2(1-s^b_{n,k})}{b_k^{\star,\tau}}.
\]
All coefficients entering the matrix entry on the left have modulus at most \(M+M_0\). Its path expansion therefore gives a uniform upper bound
\(|2[V'(T)]_{k,k+1}|\le C_{V,M,M_0}\). The right side is positive and
\(1-s^b_{n,k}\ge1-s_1\), hence
\[
b_k^{\star,\tau}\ge
\frac{2(1-s_1)}{C_{V,M,M_0}}=:m>0.
\]
Proposition~\ref{A:res:11.3} gives the upper and lower bounds with the same parameter dependence for the local minimizer. Every convex combination of the two coefficient vectors
consequently has \(|a_k|\le M'\) and \(m'\le b_k\le M'\) for fixed \(m',M'\).

For the lower Hessian bound, extend a free coefficient direction by zero in the exterior and apply \eqref{A:eq:finite-trace-convexity}. This gives
\[
 D^2\operatorname{Tr}V(T)[\dot y,\dot y]
 \ge c_V\left(\sum_k\dot a_k^2+2\sum_k\dot b_k^2\right).
\]
Every retained logarithmic term has nonnegative second derivative, so this proves the lower bound. For the upper bound, write $V(t)=\sum_{m=0}^dv_mt^m$ and let $S_i=\partial_iT$
for a free coefficient coordinate. Each $S_i$ has Hilbert--Schmidt norm at most $\sqrt2$. Set the coefficients outside the free block and its $R_V$-neighborhood to zero, and denote
the resulting matrix by $\widehat T$. Every path of length at most $d$ that enters a free derivative stays in this neighborhood. Thus the free second derivatives of the two
polynomial traces agree. The retained coefficients are bounded by a fixed $M'$, including the conditioned exterior values in that neighborhood, and
$\|\widehat T\|_{\mathrm{op}}\le3M'$. Differentiation gives
\[
 \partial_{ij}^2\operatorname{Tr}V(T)
 =\sum_{m=2}^dmv_m\sum_{r=0}^{m-2}
       \operatorname{Tr}(\widehat T^rS_i\widehat T^{m-2-r}S_j).
\]
The trace inequality $|\operatorname{Tr}(AB)|\le \|A\|_{\mathfrak S_2}\|B\|_{\mathfrak S_2}$ implies
\[
 |\partial_{ij}^2\operatorname{Tr}V(T)|
 \le2\sum_{m=2}^dm(m-1)|v_m|(3M')^{m-2}=:C_{V,M'}.
\]
The path expansion makes this entry zero when the two coefficient sites are separated by more than $d$. There are at most $C_d$ nonzero entries in each row, including the two
coefficient types. The polynomial Hessian therefore has operator norm at most $C_dC_{V,M'}$ by Schur's test. The logarithmic Hessian is diagonal and has entries
$2(1-s^b_{n,k})/b_k^2\le2/(m')^2$. Adding these bounds proves the uniform upper bound along the segment; the same path restriction gives its fixed bandwidth.
\end{proof}

For boundary data \(\tau\), put
\[
\Delta_\partial(\tau)
=1+\max_{\xi\in\partial_R\Lambda}
|\tau_\xi-y_\xi^\dagger|,
\]
with the maximum equal to zero when $\partial_R\Lambda$ is empty. The maximum is over actual exterior coefficients only. The additional one covers the finitely many missing path
terms at a physical boundary point, even when all available exterior values equal the local profile. The gradient of the unscaled energy evaluated at this profile has an
$O(n^{-1})$ contribution and an $O(\Delta_\partial(\tau))$ contribution supported near the boundary. Multiplication by $\beta n/2$ changes these bounds to $O(1)$ and
$O(n\Delta_\partial(\tau))$, but also changes the Hessian scale from one to $n$. The resulting displacement estimate is identical.

\begin{proposition}\label{A:res:11.8} Under the assumptions of
Lemma~\ref{A:res:11.7A}, there exist constants \(C,c>0\), independent of \(n\), \(\Lambda\), and its length, such that every free coefficient \(\xi\in\Lambda\) satisfies
\[
|y_{n,\Lambda}^{\star,\tau}(\xi)-y_n^\dagger(\xi)|
\le \frac{C}{n}
+C\Delta_\partial(\tau)
 \mathrm{e}^{-c\operatorname{dist}(\xi,\partial\Lambda)}.
\eqtag{11.66}
\]
In particular, when the boundary data range in a fixed compact set and \(\operatorname{dist}(\xi,\partial\Lambda)\ge A\log n\), with \(A\) large enough,
\[
|y_{n,\Lambda}^{\star,\tau}(\xi)-y_n^\dagger(\xi)|
\le \frac{C}{n}.
\eqtag{11.67}
\]
\end{proposition}

\begin{proof}
Set \(F=\nabla\mathscr H_{n,\Lambda}^{\tau}(\mathbf y_n^\dagger|_\Lambda)\). Let \(\mathcal B_\Lambda\) be the free coefficients whose distance from a conditioned boundary or a
boundary point of the finite matrix is at most \(R\). For
\(\xi\notin\mathcal B_\Lambda\), the conditional derivative agrees with the corresponding derivative of the full finite
energy at the local minimizer. Proposition~\ref{A:res:11.7} therefore gives \(|F_\xi|\le C/n\).

For \(\xi\in\mathcal B_\Lambda\), compare the conditional coefficients in the set $\{\eta:|\jmath(\eta)-\jmath(\xi)|\le R\}$ with the local minimizing coefficients.
Lemma~\ref{A:res:11.6A} gives a Lipschitz bound for the derivative of the polynomial trace. The logarithmic derivative is local to the free off-diagonal coordinate and is uniformly
bounded because the minimizing off-diagonal coefficient is bounded away from zero. To bound paths excluded by a finite boundary point, write $V(t)=\sum_{m=0}^dv_mt^m$ and bound all
coefficients in the fixed neighborhood by $M\ge1$. For a diagonal coefficient the derivative of $\operatorname{Tr}T^m$ is $m(T^{m-1})_{kk}$; for an off-diagonal coefficient it is
$2m(T^{m-1})_{k,k+1}$. Each entry of $T^{m-1}$ has at most $3^{m-1}$ path terms, each of absolute value at most $M^{m-1}$. The absolute sum of excluded derivative terms is
consequently at most
\[
 2\sum_{m=1}^d m|v_m|(3M)^{m-1}=:C_{V,M}.
\]
This is bounded by $C_{V,M}\Delta_\partial(\tau)$ because $\Delta_\partial(\tau)\ge1$. Consequently
\[
|F_\xi|\le \frac{C}{n}+C\Delta_\partial(\tau),
\quad \xi\in\mathcal B_\Lambda.
\]
Decompose \(F=F_{\mathrm{bulk}}+F_\partial\), where
\[
\|F_{\mathrm{bulk}}\|_{\ell^\infty}\le \frac{C}{n},
\quad
\operatorname{supp}F_\partial\subset\mathcal B_\Lambda,
\quad
\|F_\partial\|_{\ell^\infty}
\le C\Delta_\partial(\tau).
\]

Let \(\delta=\mathbf y_{n,\Lambda}^{\star,\tau} -\mathbf y_n^\dagger|_\Lambda\). The mean-value theorem for the gradient gives
\[
0=F+\widetilde K\delta,
\quad
\widetilde K
=\int_0^1D^2\mathscr H_{n,\Lambda}^{\tau}
(\mathbf y_n^\dagger|_\Lambda+t\delta)\,\mathrm{d}t.
\]
Lemma~\ref{A:res:11.7A} gives $cI\le\widetilde K\le CI$ and a fixed bandwidth. Apply Proposition~\ref{A:res:A.11} with this Hessian scale, the bounds
$\|F_{\mathrm{bulk}}\|_\infty\le C/n$ and $\|F_\partial\|_\infty\le C\Delta_\partial(\tau)$, and boundary width at most $R$ in site distance. Passing to the interlaced coefficient
order changes these constants only by fixed factors. Since $\delta=-\widetilde K^{-1}F$, that proposition gives
\[
 |\delta_\xi|\le\frac Cn+C\Delta_\partial(\tau)
                 \mathrm e^{-c\operatorname{dist}(\xi,\partial\Lambda)}.
\]
This proves \eqref{A:eq:11.66}. On a fixed compact set of exterior values, $\Delta_\partial(\tau)\le C$. Taking $A\ge1/c$ makes the second term $O(n^{-1})$ at distance $A\log n$
and proves \eqref{A:eq:11.67}. The exponential decay of boundary influence thus determines the logarithmic distance.
\end{proof}

Under these compact parameter and exterior hypotheses, the local profile can replace the conditional minimizer at distance $A\log n$ from the boundary with a coordinate error
$O(n^{-1})$. We next estimate the fluctuations about the conditional minimizer and the effect of the two matrix boundaries.

\subsubsection{Comparison with principal submatrices}\label{A:sec:boundary-reduction}

The two deleted blocks require different arguments. At the right cut, the off-diagonal coefficient is small and the terminal spectrum stays a fixed unscaled distance below
$E_\star$. Its Schur complement can therefore be estimated perturbatively. At the left cut, the scaled off-diagonal coefficient has order $h^{-2}$ and is not a small perturbation.
There we need a spectral estimate for the original leading block and a sharper bound for the boundary column of the retained resolvent.

We first obtain these bounds from the coefficient law, before constructing the Brownian coupling. The form shift at this stage may grow as $(1+\log n)^2$; this is sufficient for an
$o(h)$ boundary error. The common shift independent of $n$ is proved later in Theorem~\ref{A:res:16.2}. The comparison controls every block of the full resolvent, with the cross
terms at both cuts retained.

Use the matrix scaling \eqref{A:eq:edge-scaling} and the cuts $r_n,s_\beta,I_n$ fixed in the introduction. Set $N_n=n-s_\beta$ and $r=r_n$. Set the auxiliary boundary values
$b_0=b_n=b_\star$. The boundary conditions $u_0=u_{n+1}=0$ then give the quadratic form
\[
\mathfrak a_n[u] ={} \sum_{k=0}^{n} \frac{b_k}{b_\star}|D_hu_k|^2 + \frac{h^{-2}}{b_\star} \sum_{k=1}^{n} \left(E_\star-a_k-b_{k-1}-b_k\right)|u_k|^2 .
\eqtag{12.1}
\]
Indeed,
\[
\sum_{k=0}^{n}b_k|u_{k+1}-u_k|^2
=
\sum_{k=1}^{n}(b_{k-1}+b_k)|u_k|^2
-
2\operatorname{Re}
\sum_{k=1}^{n-1}b_ku_{k+1}\overline{u_k}.
\]
Equation \eqref{A:eq:12.1} is
\[
\frac{h^{-2}}{b_\star}
\left[
E_\star\sum_{k=1}^{n}|u_k|^2
-
\langle u,T(a,b)u\rangle
\right].
\]

For each free edge $k\le N_n-1$, the logarithmic term in \eqref{A:eq:11.1} has second derivative
\[
 \frac{2(1-k/n-1/(\beta n))}{b_k^2}
 \ge\frac{2(s_\beta+1-\beta^{-1})}{n b_k^2}>0.
\]
Every edge with a negative logarithmic coefficient has $k>n-\beta^{-1}$. A term involving that edge uses only sites at distance at most $R_V$, and therefore lies to the right of
$N_n$, since $s_\beta>\beta^{-1}+R_V$. Such coordinates are all fixed in the conditional law of the selected interval. The deleted matrix block has vertices \(\{N_n+1,\ldots,n\}\).
For the conditional law of the initial coefficients, the free coordinates are
\(a_1,\ldots,a_{N_n},b_1,\ldots,b_{N_n-1}\), and the fixed terminal coordinates are
\(a_{N_n+1},\ldots,a_n,b_{N_n},\ldots,b_{n-1}\). Only fixed coordinates within $R_V$ indices of the initial segment
affect its conditional density. The bound $R_V$ identifies the terms that contain coefficients on both sides of the cut. These terms remain in the conditional energy. Since
\(s_\beta\) is fixed, \(N_n=n-O_{\beta,V}(1).\) The coordinate of the terminal region is \(x_{N_n} = \left(N_n-\frac12\right)h \asymp nh \asymp h^{-2}.\)

The deleted left block is $\mathcal L_n=\{1,\ldots,r\}$, and the selected block is $I_n$. Its displacement is $\delta_n=rh$, so
\[
\delta_n^2\le Ch^2(\log n)^2=o(h).
\eqtag{12.11}
\]
This estimate will be used after the interior interval is reindexed from its left boundary.

Also $\mathrm e^{-cr_n}\le n^{-cA_0}$, so choosing $A_0$ sufficiently large supplies the polynomial bounds used below.

Retain the conditional laws and the constants $0<s_0<s_1<1$ from Subsection~\ref{A:sec:log-concave-measures-and-conditional-minimizers}. The estimates below concern the original
exterior distribution and the full conditional domains $\mathcal D_J$; bounded exterior data do not restrict the free coefficients.

Let $Y_n$ denote the full coefficient vector. There are constants \(M_0\ge1\) and \(b_{\min}>0\) such that the coefficient event
\[
 \mathcal E_n
 =\left\{\max_k(|a_k|\vee b_k)\le M_0,\quad
             \min_{1\le k\le\lfloor s_1n\rfloor}b_k\ge b_{\min}\right\}
\eqtag{raw-event}
\]
has an exponentially small complement. The maximum runs over the original diagonal and off-diagonal coordinates.

\begin{lemma}\label{A:res:raw}
For all sufficiently large \(n\),
\[
 \mathbb P(\mathcal E_n^c)\le C\mathrm{e}^{-cn},
 \quad
 \sup_n\mathbb E\left[\max_k(|a_k|\vee b_k)^p\right]<\infty
 \quad(1\le p<\infty).
\eqtag{raw-probability}
\]
If \(\nu_{n,J}\) is the distribution of the original exterior data for any coefficient interval \(J\), then
\[
 \nu_{n,J}(\mathcal T_{n,J}(M_0)^c)\le C\mathrm{e}^{-cn}.
\eqtag{exterior-probability}
\]
The constants can be chosen uniformly on a compact positive \(\beta\)-interval. The event $\mathcal E_n$ is defined solely by the displayed coefficient bounds.
\end{lemma}

\begin{proof}
We prove the upper bound, including the case \(\beta<1\); it is the argument for bounding the coefficient vector of \cite{krishnapur2016universality} in the present normalization.
Write \(\nu_k=\beta(n-k)-1\). There are only a bounded number of negative \(\nu_k\), uniformly on any compact positive \(\beta\)-interval. Reveal which of their off-diagonal
coefficients are below one and condition on the values of those off-diagonal coefficients. The other negative- \(\nu_k\) coordinates are then restricted to \([1,\infty)\). This is
a convex product domain. The fixed off-diagonal coefficients are at most one, whereas, on the restricted coordinates, \(\partial_b^2(-\nu_k\log b)=\nu_k/b^2\ge-1\). The polynomial
Hessian is at least \(\beta_-c_VnI/2\), so the conditional density has Hessian at least \(cnI\) for all sufficiently large \(n\).

The theorem on conditional minimizers just used applies also to these restricted negative logarithms: its cutoff parameter is one, its fixed data are bounded by one, and
\(|w_{k,n}|\le1\). It bounds every coordinate of the conditional minimizer by a fixed \(C_V\). Lemma~\ref{A:res:A.8} applies on the convex product domain, including a minimizer on
its boundary. Let \(\mathscr R_n\) be the sigma-field generated by the terminal data just revealed. In dimension at most \(2n\), it gives
\[
 \mathbb P\{\max_k(|a_k|\vee b_k)>C_V+C_1+t
       \mid\mathscr R_n\}
 \le2\mathrm{e}^{-cnt^2},\quad t\ge0,
\eqtag{field-upper-tail}
\]
after increasing \(C_1\) to include the fixed off-diagonal coefficients. The constants are independent of the subset revealed and of its values. Averaging this inequality proves
the upper bound in \eqref{A:eq:raw-probability}, and integrating its tail proves all the stated moments for large \(n\). Fix one of the remaining sizes $n$ and
$I_\beta=[\beta_-,\beta_+]\Subset(0,\infty)$. Uniform convexity implies $V(t)\ge c_Vt^2/4-C_V$, hence
\[
 \exp\{-\beta n\operatorname{Tr}V(T)/2\}
 \le C_{n,I_\beta,V}
       \exp\{-c_{n,I_\beta,V}(\sum_k a_k^2+\sum_k b_k^2)\}.
\]
On $0<b_k\le1$ the power in the density is bounded above by $b_k^{\beta_-(n-k)-1}$; on $b_k>1$ it is bounded above by $b_k^{\beta_+(n-k)-1}$. The first exponent is greater than
$-1$, and the Gaussian bound makes the second part integrable after multiplication by any fixed polynomial in the coefficients. These bounds give a common integrable majorant on
$I_\beta$. The normalizing integral is continuous and positive there by dominated convergence, and so has a positive minimum. Every fixed coefficient moment is consequently bounded
uniformly on $I_\beta$. Taking the maximum over the finitely many remaining sizes completes the proof of the moment bound. This auxiliary decomposition is used only for the bound
on the coefficients. All subsequent conditional laws are those of the original Gibbs measure.

For the lower bound, fix \(k\le s_1n\) and condition on all other coordinates. On the set, measurable with respect to the exterior coordinates, where they are bounded by \(M_0\),
the minimizer with respect to that coordinate belongs to
\([b_*,C_*]\), by Lemma~\ref{A:res:12.1A} and \eqref{A:eq:regular-mode-lower}. Its conditional density has curvature at
least \(cn\). The one-dimensional case of Lemma~\ref{A:res:A.8} gives
\[
 \mathbb P\{b_k<b_*/2\mid\mathscr F_{\{b_k\}^c}\}
 \mathbf1_{\{\max_{i\ne b_k}|Y_{n,i}|\le M_0\}}
 \le 2\mathrm{e}^{-cn}.
\]
Let $M_n=\max_k(|a_k|\vee b_k)$ and $E_k=\{\max_{i\ne b_k}|Y_{n,i}|\le M_0\}$. For each regular edge, integration of the conditional bound gives
$\mathbb P(\{b_k<b_*/2\}\cap E_k)\le2e^{-cn}$. Moreover, $\{M_n\le M_0\}\subset E_k$. Thus
\[
 \mathbb P\!\left(\{M_n>M_0\}\cup
       \{\min_{k\le s_1n}b_k<b_*/2\}\right)
 \le \mathbb P(M_n>M_0)+2n e^{-cn}.
\]
Choose $M_0$ above the fixed threshold in \eqref{A:eq:field-upper-tail}. That estimate bounds the first term by $Ce^{-c_1n}$. Since $n e^{-cn}\le[2/(ce)]e^{-cn/2}$ for $n>0$, the
last display is at most $C'e^{-c'n}$ for $c'=\min\{c_1,c/2\}$. Putting $b_{\min}=b_*/2$ proves \eqref{A:eq:raw-probability}. The inclusion
\(\mathcal E_n\subset\{\tau_J\in\mathcal T_{n,J}(M_0)\}\) proves \eqref{A:eq:exterior-probability}. Here the conditioning is on the exterior coordinates.
\end{proof}

The next lemma integrates the conditional estimates over the exterior law.

Suppose \(\mu^\tau\) is an exact conditional law, its exterior marginal is \(\nu\), and a measurable set of exterior data \(T\) satisfies \(\nu(T^c)\le\epsilon\).

\begin{lemma}\label{A:res:disintegration}
For any event \(F\) and any random variable \(W\) defined under these laws,
\[
 \mathbb P(F)\le\sup_{\tau\in T}\mu^\tau(F)+\epsilon,
\eqtag{integrated-tail}
\]
and, if \(\sup_{\tau\in T}\|W\|_{L^p(\mu^\tau)}\le a\), then
\[
 \|W\|_{L^p}\le a+\|W\|_{L^{2p}}\epsilon^{1/(2p)}.
\eqtag{integrated-moment}
\]
\end{lemma}

For the coefficient fluctuations \(W=Y_{n,i}-y_{J,i}^\circ(\tau_J)\), the coarse norm on the right is bounded uniformly in \(n,J,i\). More generally, a polynomial bound in \(n\)
for that norm is sufficient when \(\epsilon\le C\mathrm{e}^{-cn}\).

\begin{proof}
The conditional probability identity gives
\[
 \mathbb P(F)=\int_T\mu^\tau(F)\,\nu(\mathrm d\tau)
             +\int_{T^c}\mu^\tau(F)\,\nu(\mathrm d\tau)
 \le\sup_{\tau\in T}\mu^\tau(F)+\nu(T^c)
 \le\sup_{\tau\in T}\mu^\tau(F)+\epsilon.
\]
This is \eqref{A:eq:integrated-tail}. On \(T\), integration gives \(\mathbb E[|W|^p\mathbf1_T]\le a^p\). On its complement, H\"older gives
\(\mathbb E[|W|^p\mathbf1_{T^c}] \le\|W\|_{L^{2p}}^p\epsilon^{1/2}\). Taking a \(p\)-th root proves \eqref{A:eq:integrated-moment}. Put $Z_n=\max_i|Y_{n,i}|$. Since the exterior
coordinates are coordinates of the same vector, \eqref{A:eq:mode-growth} gives
\[
 |Y_{n,i}-y_{J,i}^\circ(\tau_J)|
 \le Z_n+C_V(1+Z_n)^{D_V}
 \le C_V'(1+Z_n)^{\max\{1,D_V\}}.
\]
Lemma~\ref{A:res:raw} bounds the $2p\max\{1,D_V\}$ moment of $1+Z_n$, uniformly in $n,J,i$. This proves the coarse norm bound. More generally, write a polynomial
$P=\sum_\alpha p_\alpha x^\alpha$ of degree at most $d$ and set $\|P\|_1=\sum_\alpha|p_\alpha|$. For its arguments formed from coefficients and their conditional minimizers,
\[
 |P|\le C_{V,d}\|P\|_1(1+Z_n)^{d\max\{1,D_V\}}.
\]
Thus its moment bound is uniform when $\|P\|_1$ is uniform, and polynomial in $n$ when $\|P\|_1$ is polynomial in $n$.
\end{proof}

Fix \(M\ge1\) and \(1\le p<\infty\). For every coefficient interval \(J\) and every \(\tau\in\mathcal T_{n,J}(M)\), all expectations below are under the full law
\(\mu_{n,J}^{\tau}\) in \eqref{A:eq:cond-density}.

\begin{lemma}\label{A:res:mean}
If all coefficients within $R_V$ indices of $i$ lie in $\{k/n\le s_0\}$, then
\[
 \|Y_i-y_{J,i}^\circ\|_{L^p(\mu_{n,J}^{\tau})}
 \le C_{p,M}n^{-1/2},
 \quad
 |\mathbb E_\tau Y_i-y_{J,i}^\circ|\le C_Mn^{-1}.
\eqtag{mean-regular}
\]
For every free coordinate, including the nonregular region,
\[
 \|Y_i-y_{J,i}^\circ\|_{L^p(\mu_{n,J}^{\tau})}
 \le C_{p,M}\sqrt{\frac{\log(n+1)}{n}}.
\eqtag{mean-rough}
\]
\end{lemma}

The constants are deterministic and uniform in the indicated exterior values. For a regular coordinate and \(t\ge0\),
\[
 \mu_{n,J}^{\tau}\{|Y_i-y_{J,i}^\circ|>C_M/n+t\}
 \le2\mathrm{e}^{-cnt^2}.
\eqtag{mode-tail}
\]
For every \(q>0\),
\[
 \sup_{\tau\in\mathcal T_{n,J}(M)}
 \mu_{n,J}^{\tau}\left\{
   \max_{i\in J}|Y_i-y_{J,i}^\circ|
       >C_{q,M}\sqrt{\frac{\log(n+1)}n}\right\}
 \le C_{q,M}n^{-q}.
\eqtag{mode-maximum}
\]
For uniformity in $\beta\in I_\beta=[\beta_-,\beta_+]\Subset(0,\infty)$, use the same final $s_*$ omitted indices for every $\beta$, where $s_*$ is given in \eqref{A:eq:B.45}. The
constants above are uniform with this convention. Lemma~\ref{A:res:disintegration} then integrates the conditional estimates over the exterior distribution: the moment rates are
unchanged, and the probability bounds acquire at most $C\mathrm e^{-cn}$.

\begin{proof}
Fix the ambient law \(\mu_{n,J}^{\tau}\), with \(\tau\in\mathcal T_{n,J}(M)\), and write \(y^\circ=y_J^\circ(\tau)\). Set \(L=\lceil A\log(n+1)\rceil\). Around a free coordinate
\(i\) choose a block \(I_i\subset J\) of at most \(CL\) coordinates. Each boundary of \(I_i\) inside \(J\), where the exterior data can differ from \(y^\circ\), is at distance at
least \(L\) from \(i\). Near a fixed boundary point of \(J\) the block is one-sided; if \(J\) is shorter than this, use
\(I_i=J\). In the regular case choose \(I_i\subset\{k/n\le s_1\}\), which is possible for large \(n\) since \(s_0<s_1\).

Condition further on \(\mathscr F_{I_i^c}\). The conditional density is the exact block density; denote its minimizer by
\(y^{i,\circ}\). Let \(\mathcal B_i\in\mathscr F_{I_i^c}\) be the event that the additional, random exterior coordinates
in \(\partial I_i\cap J\) are bounded by \(M_1\), where \(M_1\) is chosen in \eqref{A:eq:conditional-upper}. The already fixed coordinates are bounded by \(M\). Thus
\[
 \mu_{n,J}^{\tau}(\mathcal B_i^c)\le C\mathrm{e}^{-cn}.
\eqtag{local-exterior-tail}
\]
This follows from \eqref{A:eq:conditional-upper}, since only the additional exterior coordinates enter $\mathcal B_i$. On $\mathcal B_i$, the conditional minimizers \(y^\circ\) and
\(y^{i,\circ}\) and the coefficients within $R_V$ indices of them are bounded by a deterministic constant, by \eqref{A:eq:mode-growth}.

Evaluate the block gradient at \(y^\circ|_{I_i}\). It vanishes except within distance \(R_V\) of a boundary of \(I_i\) inside \(J\), and there its entries are at most \(Cn\). The
averaged Hessian along the segment from \(y^\circ|_{I_i}\) to \(y^{i,\circ}\) has the form \(K_i=P_i+D_i\), where \(D_i\ge0\) is diagonal, \(P_i\) has fixed range and row sums at
most \(Cn\), and \(K_i\ge cnI\). The diagonal matrix $D_i$ is the Hessian of the logarithmic part of the potential. Lemma~\ref{A:res:A.9} and the gradient identity imply, on
\(\mathcal B_i\),
\[
 |y_i^{i,\circ}-y_i^\circ|
 \le C \mathrm{e}^{-cL}\le n^{-4},
\eqtag{local-mode-comparison}
\]
after fixing \(A\) sufficiently large. This boundary contribution is zero when \(I_i=J\).

The block has at most \(CL\) coordinates. Applying Lemma~\ref{A:res:A.8} to its conditional law gives
\[
 \bigl\|\,\|Y_{I_i}-y^{i,\circ}\|_2\,\bigr\|_{L^p(
            \,\cdot\mid\mathscr F_{I_i^c})}
 \le C_p\sqrt{L/n}.
\eqtag{local-radius}
\]
On \(\mathcal B_i\), combine this with \eqref{A:eq:local-mode-comparison}. On \(\mathcal B_i^c\), use Lemma~\ref{A:res:A.8} under $\mu_{n,J}^{\tau}$
\(\|Y_i-y_i^\circ\|_{L^{2p}(\mu_{n,J}^{\tau})}\le C_p\) and H\"older. Precisely,
\[
 \|Y_i-y_i^\circ\|_{L^p(\mu_{n,J}^{\tau})} \le C_p\sqrt{L/n}+n^{-4} +C_p\mu_{n,J}^{\tau}(\mathcal
 B_i^c)^{1/(2p)} \le C_{p,M}\sqrt{\log(n+1)/n}.
\eqtag{rough-disintegration}
\]
This proves the rough moment bound.

For the sharper estimate, work under the exact conditional law $\mu_i$ on $I_i$, on the exterior event $\mathcal B_i$. All its free off-diagonal minimizing values are bounded below
by a deterministic positive constant. Choose $a_0>0$, uniformly over these exterior values, so that they are at least $8a_0$. Reducing $a_0$ if necessary, every coefficient is
bounded on the cube of radius $2a_0$ around $y^{i,\circ}$. On that cube, the path expansion and the logarithmic derivatives give
\[
 cnI\le\nabla^2\Phi_i\le CnI,\qquad
 |\partial_{jkl}^3\Phi_i|\le Cn.
\]
Each derivative has a fixed finite range. These local derivative bounds are evaluated directly for the original potential on the indicated intervals.

Put $m=\mathbb E_{\mu_i}Y$ and $U=Y-m$. Jensen's inequality and Lemma~\ref{A:res:A.8} give
\[
 \|m-y^{i,\circ}\|_2\le\mathbb E_{\mu_i}\|Y-y^{i,\circ}\|_2
 \le C\sqrt{d_i/n}\le C\sqrt{\log(n+1)/n}.
\]
Thus $m$ lies in the cube of radius $a_0/2$ for large $n$. Define $\mathcal C_i=\{\|Y-y^{i,\circ}\|_\infty\le a_0\}$. The Gaussian exponential moment in Lemma~\ref{A:res:A.8} gives
\[
 \mu_i(\mathcal C_i^c)
 \le2^{d_i/2}\mathrm e^{-cna_0^2/4}
 \le C\mathrm e^{-c'n},\qquad d_i\le C\log(n+1).
\eqtag{local-extension-tail}
\]
Choose $\vartheta\in C_c^\infty((-2,2))$ with $0\le\vartheta\le1$ and $\vartheta=1$ on $[-1,1]$, and set
\[
 \chi(y)=\prod_{j\in I_i}\vartheta((y_j-y_j^{i,\circ})/a_0).
\]
Its support is compactly contained in the positive coordinate domain. For each free coordinate $k$, integration by parts against the exact density is therefore legitimate and gives
\[
 \mathbb E_{\mu_i}[\chi\,\partial_k\Phi_i]
 =\mathbb E_{\mu_i}\partial_k\chi.
\]
The derivative of $\chi$ is bounded by a fixed constant and vanishes on $\mathcal C_i$, so the right side is $O(\mathrm e^{-c'n})$. Also $\mathbb E\chi\ge1/2$ for large $n$.

The variance estimate in Subsection~\ref{A:sec:log-concave-measures-and-conditional-minimizers} gives $\mathbb E|U_j|^2\le C/n$ for each free coordinate. On the support of $\chi$,
Taylor's formula at $m$ gives
\[
 \partial_k\Phi_i(Y)=\partial_k\Phi_i(m)
       +\sum_j(\nabla^2\Phi_i(m))_{kj}U_j+R_k(Y),\qquad
 |R_k(Y)|\le Cn\sum_{j\in N_k}|U_j|^2,
\]
where $N_k$ has uniformly bounded size. The segment from $m$ to $Y$ remains in the cube on which the derivative bounds hold. Hence $\mathbb E[\chi|R_k|]\le C$. Since
$\mathbb EU_j=0$,
\[
 |\mathbb E(\chi U_j)|
 =|\mathbb E((1-\chi)U_j)|
 \le(\mathbb EU_j^2)^{1/2}\mu_i(\mathcal C_i^c)^{1/2}
 \le Cn^{-1/2}\mathrm e^{-c'n/2}.
\]
The Hessian row sums are at most $Cn$. Substitution in the integration by parts identity now yields
\[
 |\partial_k\Phi_i(m)|\le C
\]
uniformly in $k$ and in $d_i$. Finally, on the segment from $y^{i,\circ}$ to $m$,
\[
 \nabla\Phi_i(m)=K_i(m-y^{i,\circ}),\qquad
 K_i=\int_0^1\nabla^2\Phi_i(y^{i,\circ}+t(m-y^{i,\circ}))\,\mathrm dt.
\]
Lemma~\ref{A:res:A.9} gives
\[
 |m_j-y_j^{i,\circ}|
 \le\sum_k|(K_i^{-1})_{jk}|\,|\partial_k\Phi_i(m)|
 \le\frac Cn\sum_k\mathrm e^{-c|j-k|}\le\frac Cn.
\]
Concentration about $m$ therefore proves, under $\mu_i$,
\[
 |\mathbb E(Y_i-y_i^{i,\circ}\mid\mathscr F_{I_i^c})| \le C/n,\qquad
 \|Y_i-y_i^{i,\circ}\|_{L^p(\,\cdot\mid\mathscr F_{I_i^c})} \le C_pn^{-1/2}.
\eqtag{local-sharp}
\]
The integration by parts identity uses the cutoff under the original coefficient law.

Use \eqref{A:eq:local-mode-comparison} on \(\mathcal B_i\), and the radial estimate under $\mu_{n,J}^{\tau}$ and \eqref{A:eq:local-exterior-tail} on its complement. The tower
property gives
\[
 \begin{split}
 |\mathbb E_\tau Y_i-y_i^\circ|
 &\le C/n+n^{-4}+C\mu_{n,J}^{\tau}(\mathcal B_i^c)^{1/2},\\
 \|Y_i-y_i^\circ\|_{L^p(\mu_{n,J}^{\tau})}
 &\le C_pn^{-1/2}+n^{-4}
       +C_p\mu_{n,J}^{\tau}(\mathcal B_i^c)^{1/(2p)}.
 \end{split}
\eqtag{sharp-disintegration}
\]
These inequalities prove \eqref{A:eq:mean-regular}.

Finally, uniform conditional convexity and the logarithmic Sobolev inequality on the convex domain give
\[
 \mu_{n,J}^{\tau}\{|Y_i-\mathbb E_\tau Y_i|>t\}
 \le2\mathrm{e}^{-cnt^2}.
\]
The bound $|\mathbb E_\tau Y_i-y_i^\circ|\le C/n$ at a regular coordinate now gives \eqref{A:eq:mode-tail}. For all coordinates, \eqref{A:eq:mean-rough} with \(p=1\) bounds the
displacement of the mean from the minimizer by
\(C\sqrt{\log(n+1)/n}\). A union bound over at most \(2n\) coordinates, with \(t=C_q\sqrt{\log(n+1)/n}\), proves
\eqref{A:eq:mode-maximum}. The corresponding unconditional tail and moment bounds follow from \eqref{A:eq:exterior-probability} and Lemma~\ref{A:res:disintegration}, using the
coarse moment bounds established there. The constants are uniform for $\beta$ in a compact subinterval of $(0,\infty)$ by Lemma~\ref{A:res:12.1A}.
\end{proof}

For any deterministic family \(\mathfrak I_n\) of at most \(Cn^3\) coefficient intervals of the above types, let
\[
 \mathcal D_{n,Q}:=\mathcal E_n\cap \bigcap_{J\in\mathfrak I_n} \left\{\max_{i\in
 J}|Y_{n,i}-y_{J,i}^\circ(\tau_J)| \le C_Q\sqrt{\frac{\log(n+1)}n}\right\}.
\eqtag{finite-coefficient-event}
\]
Fix this family once so that it includes the initial segment, reduced and transport intervals, every logarithmic block used above, and the terminal windows of diagonal coordinates.
Such a deterministic choice has at most \(Cn^3\) members; in fact all coefficient intervals and all intervals of diagonal coordinates already number at most \(Cn^2\). Apply
\eqref{A:eq:mode-maximum} with exponent \(Q+4\), integrate over bounded exterior data, and use a union bound. It gives
\[
 \mathbb P(\mathcal D_{n,Q}^c)\le C_Qn^{-Q}.
\eqtag{finite-coefficient-probability}
\]
Set \(\Omega_{\mathrm{bd}}=\liminf_{n\to\infty}\mathcal D_{n,3}\). By Borel--Cantelli this event has probability one under any coupling of the laws of the finite matrices,
regardless of dependence across \(n\). We use
\(\Omega_{\mathrm{bd}}\) for pathwise conclusions only. Conditional laws remain those in \eqref{A:eq:cond-density}; to
integrate their estimates, apply

\eqref{A:eq:mean-regular}--\eqref{A:eq:mode-maximum} on \(\mathcal T_{n,J}(M_0)\) and then \eqref{A:eq:integrated-tail} or \eqref{A:eq:integrated-moment}.

We derive the terminal spectral gap directly from the KRV Euler--Lagrange equations.

Write $(a_s,b_s)=(a^\dagger(s),b^\dagger(s))$ and $\varepsilon_s=1-s$. The minimum $a_{\mathrm t}$ of $V$ satisfies
$V'(a_{\mathrm t})=0$.
\begin{lemma}\label{A:res:T1}
For every \(s<1\),
\[
0<b_s^2\le \frac{1-s}{c_V},
\quad
|a_s-a_{\mathrm t}|<2b_s.
\eqtag{T2}
\]
Consequently,
\[
a_s\longrightarrow a_{\mathrm t},
\quad
b_s\longrightarrow0,
\quad
E(s):=a_s+2b_s\longrightarrow a_{\mathrm t}
\eqtag{T3}
\]
as \(s\uparrow1\). Moreover,
\[
d_{\mathrm t}:=E_\star-a_{\mathrm t}>0.
\eqtag{T4}
\]
\end{lemma}

\begin{proof}
The Euler--Lagrange equations are
\[
\frac1{2\pi}\int_0^{2\pi}
V'(a_s+2b_s\cos\theta)\,\mathrm{d}\theta=0
\eqtag{T5}
\]
and
\[
\partial_bW_V(a_s,b_s)
=
\frac{2(1-s)}{b_s}.
\eqtag{T6}
\]
For fixed \(a\), the function \(b\mapsto W_V(a,b)\) is even, so its derivative vanishes at \(b=0\). Also
\[
\partial_b^2W_V(a,b)
=
\frac1{2\pi}\int_0^{2\pi}
4\cos^2\theta\,V''(a+2b\cos\theta)\,\mathrm{d}\theta
\ge 2c_V.
\]
Hence \(\partial_bW_V(a_s,b_s)\ge2c_Vb_s.\) Together with \eqref{A:eq:T6}, this gives \(b_s^2\le(1-s)/c_V\).

Because \(V'\) is strictly increasing and vanishes at \(a_{\mathrm t}\), \eqref{A:eq:T5} rules out both \(a_s-2b_s\ge a_{\mathrm t}\) and \(a_s+2b_s\le a_{\mathrm t}.\) Indeed, in
the first case the integrand in \eqref{A:eq:T5} is nonnegative and strictly positive on a set of positive measure; in the second it is nonpositive and strictly negative on a set of
positive measure. Therefore \(a_s-2b_s<a_{\mathrm t}<a_s+2b_s,\) which proves the second estimate in \eqref{A:eq:T2}. It follows that \(a_{\mathrm t}<E(s)<a_{\mathrm t}+4b_s,\) and
\eqref{A:eq:T3} follows. At \(s=0\), the logarithmic term forces \(b_\star=b_0>0\), so the strict left inequality gives \(E_\star>a_{\mathrm t}\). This is \eqref{A:eq:T4}.
\end{proof}

To estimate the stationarity error near $s=1$, we use the squared off-diagonal coefficient.

Put $q(s)=(b^\dagger(s))^2$ for $s<1$ and $b^\dagger(1)=0$.

\begin{lemma}\label{A:res:profile-boundary}
The functions \(a^\dagger,q\) extend real analytically through \(s=1\), with \(a^\dagger(1)=a_{\mathrm t}\) and
\(q(1)=0\). Their derivatives of every fixed order are bounded on \([0,1]\). In particular, for \(s,t\in[0,1]\),
\[
 |a^\dagger(s)-a^\dagger(t)|+|q(s)-q(t)|\le C|s-t|,
 \quad
 |b^\dagger(s)-b^\dagger(t)|\le C|s-t|^{1/2}.
\eqtag{profile-increments}
\]
\end{lemma}

At a free coefficient, suppose every coefficient within $R_V$ indices is evaluated at the local minimizer with parameter in $[0,1)$, and its site is at distance at least $R_V$ from
each original matrix boundary point. Then the unscaled energy has errors
\[
 \left|\partial_{a_k}\mathscr H_n(y_n^\dagger)\right|\le Cn^{-1},
 \quad
 \left|\partial_{b_k}\mathscr H_n(y_n^\dagger)\right|
       \le Cn^{-1}b_{n,k}^\dagger.
\eqtag{profile-global-residual}
\]
These constants do not deteriorate as the macroscopic parameters approach one. Moreover,
\(E(s)=a^\dagger(s)+2b^\dagger(s)\) is strictly decreasing on \([0,1)\). For each fixed \(s_0\in(0,1)\),
\[
 E_\star-E(s)\ge c_{s_0}s\quad(0\le s\le s_0),
 \quad
 E_\star-E(s)\ge E_\star-E(s_0)>0\quad(s_0\le s\le1).
\eqtag{profile-edge-gap}
\]

\begin{proof}
The polynomial identity for the local trace is
\[
 W_V(a,b)=P(a,b^2),
 \quad
 P(a,q)=\sum_{m=0}^{\lfloor\deg V/2\rfloor}
             \frac{V^{(2m)}(a)}{(m!)^2}q^m.
\eqtag{profile-polynomial}
\]
Indeed, Taylor expansion is finite, odd powers of \(\cos\theta\) have zero average, and
\(\frac1{2\pi}\int_0^{2\pi}(2\cos\theta)^{2m}\,\mathrm{d}\theta =(2m)!/(m!)^2\). With \(\varepsilon=1-s\), the local
stationary equations become
\[
 P_a(a,q)=0,\quad qP_q(a,q)=\varepsilon.
\eqtag{profile-q-equations}
\]
Their Jacobian with respect to \((a,q)\) at \((a_{\mathrm t},0,0)\) is
\[
 \begin{pmatrix}
 V''(a_{\mathrm t})&V'''(a_{\mathrm t})\\
 0&V''(a_{\mathrm t})
 \end{pmatrix}.
\]
It is invertible. The analytic implicit function theorem and Lemma~\ref{A:res:T1} identify its branch with the local minimizer for \(\varepsilon>0\) small. More explicitly,
\[
 q(1-\varepsilon)
   =\frac{\varepsilon}{V''(a_{\mathrm t})}+O(\varepsilon^2),
 \quad
 a^\dagger(1-\varepsilon)
   =a_{\mathrm t}
    -\frac{V'''(a_{\mathrm t})}{V''(a_{\mathrm t})^2}
       \varepsilon+O(\varepsilon^2).
\eqtag{profile-terminal-expansion}
\]
Together with Proposition~\ref{A:res:11.3}, this gives the required regularity and derivative bounds on the entire compact interval. The square-root inequality
\(|\sqrt u-\sqrt v|\le\sqrt{|u-v|}\), for \(u,v\ge0\), proves \eqref{A:eq:profile-increments}. Only $a$ and $q$ are continued for $s>1$; we do not take a real square root of the
continuation of $q$.

For the error, a path returning to its starting vertex crosses each edge an even number of times. Thus \([V'(T)]_{kk}\) is a polynomial in the finitely many diagonals \(a_i\) and
squared off-diagonal coefficients \(q_i\) within $R_V$ indices of $k$. A path from \(k\) to \(k+1\) crosses edge \(k\) an odd number of times and every other edge an even number of
times. Consequently
\[
 [V'(T)]_{k,k+1}=b_k\Psi_k((a_i),(q_i)),
\eqtag{profile-open-path}
\]
where \(\Psi_k\) is a fixed finite-range polynomial. On bounded coefficient sets these polynomials are uniformly Lipschitz in \((a_i,q_i)\), independently of \(n,k\). For a
constant array they are respectively \(P_a(a,q)\) and
\(P_q(a,q)\): this follows either from the path expansion or from \(2[V'(J(a,b))]_{01}=\partial_bW_V=2bP_q\). Every
neighboring macroscopic parameter differs from the parameter at the specified vertex or edge by at most \(C_V/n\). Equation~\eqref{A:eq:profile-increments} therefore bounds the
diagonal error by \(C/n\). At edge \(k\), keep the factor
\(b_k=b^\dagger(s^b_{n,k})\) in \eqref{A:eq:profile-open-path} exactly, and replace only the polynomial \(\Psi_k\) by
its constant value. Since \(w_{k,n}=1-s^b_{n,k}\), the constant terms cancel:
\[
 2b_kP_q(a^\dagger(s^b_{n,k}),q(s^b_{n,k}))
       -\frac{2w_{k,n}}{b_k}=0.
\]
The remaining error is at most \(Cb_k/n\). This proves \eqref{A:eq:profile-global-residual} without dividing an error by a vanishing off-diagonal coefficient.

The edge monotonicity also follows from these stationary equations. Write \(a=a^\dagger(s)\), \(b=b^\dagger(s)\), and let angle brackets denote normalized integration over
\([0,2\pi]\). Set \(A=\langle V''(a+2b\cos\theta)\rangle\) and
\(C=\langle2\cos\theta\,V''(a+2b\cos\theta)\rangle\). Since the integral of a periodic derivative is zero,
\[
 0=\langle\cos\theta\,V'(a+2b\cos\theta)
       -2b\sin^2\theta\,V''(a+2b\cos\theta)\rangle.
\]
It follows that $\partial_b W_V/b=4\langle\sin^2\theta V''\rangle$. The $(b,b)$ entry of the Hessian at the minimizer is consequently
\[
 \partial_b^2 W_V+2(1-s)/b^2
 =4\langle\cos^2\theta V''\rangle+\partial_b W_V/b=4A.
\]
Thus $\mathsf H_s=\left(\begin{smallmatrix}A&C\\C&4A\end{smallmatrix}\right)$, and
\[
 2A\pm C=\langle2(1\pm\cos\theta)V''\rangle\ge2c_V.
\]
Differentiating both stationary equations with respect to $s$ gives
\[
 \mathsf H_s\begin{pmatrix}a'\\b'\end{pmatrix}
 =\begin{pmatrix}0\\-2/b\end{pmatrix},\qquad
 a'=\frac{2C}{b(4A^2-C^2)},\quad
 b'=-\frac{2A}{b(4A^2-C^2)}.
 \eqtag{local-hessian-system}
\]
Since $E'=a'+2b'$, cancellation of $2A-C$ gives
\[
 -E'(s)=\frac{2}{b(2A+C)}>0.
\eqtag{profile-edge-derivative}
\]
Its continuous positive minimum on \([0,s_0]\), followed by integration and the limit at the boundary of Lemma~\ref{A:res:T1}, proves \eqref{A:eq:profile-edge-gap}.
\end{proof}

On the entire retained initial segment, the logarithmic diagonal of the Hessian may be unbounded. We keep it separate from the bounded finite-range part when estimating the inverse
Hessian.

Let the free coordinates be \(a_{\ell+1},\ldots,a_m\) and \(b_{\ell+1},\ldots,b_{m-1}\), where \(0\le\ell<m\le N_n\). Fix the actual exterior values in \(\mathcal T_{n,J}(M_0)\),
and let \(y_J^\circ\) be the exact conditional minimizer. At the left or right boundary of the conditional interval, let \(\Delta_{\mathrm L}\) or \(\Delta_{\mathrm R}\) be the
maximum coefficient discrepancy between these exterior values and the local minimizer in the set of exterior coefficients within $R_V$ indices of the relevant boundary point, with
an empty maximum equal to zero. Add one to
\(\Delta_{\mathrm L}\) when \(\ell<R_V\), to account for missing full line paths within $R_V$ indices of the original
matrix boundary point; in particular \(\Delta_{\mathrm L}=1\) for the full initial segment \(\ell=0\).

\begin{proposition}\label{A:res:profile-global}
For a free coefficient with site label \(k\),
\[
 |y_J^\circ(k)-y_n^\dagger(k)|
 \le \frac Cn+C\Delta_{\mathrm L}\mathrm{e}^{-c(k-\ell)}
                 +C\Delta_{\mathrm R}\mathrm{e}^{-c(m-k)}.
\eqtag{global-mode-profile}
\]
For edge coordinates a shift by one half in \(k\) changes only the constant. The constants are independent of the interval length and of its distance from the boundary point $s=1$.
They are uniform on compact positive \(\beta\)-intervals when the final $s_*$ indices are omitted as in \eqref{A:eq:B.45}.
\end{proposition}

\begin{proof}
The terminal deletion ensures that every coefficient within $R_V$ indices of the selected block still has $1-s^b_{n,k}>0$. Thus all values of the local minimizer used here are
defined and all free logarithmic weights are positive. Lemma~\ref{A:res:12.1A} bounds the conditional minimizer on a deterministic coefficient cube. The local minimizer and every
segment joining it to the minimizer are bounded on a possibly larger cube. Along that segment the averaged unscaled Hessian has the following form, with \(D\) diagonal:
\[
 K=P+D\ge cI,\quad D\ge0,\quad
 P_{ij}=0\ (|i-j|>R),\quad \sup_i\sum_j|P_{ij}|\le C.
\eqtag{global-mode-hessian}
\]
The polynomial trace gives the finite-range part \(P\), including its positive lower convexity, and the free logarithms give \(D\). Lemma~\ref{A:res:A.9}, with scale parameter one,
yields \(|(K^{-1})_{ij}|\le C\mathrm{e}^{-c|i-j|}\), uniformly in \(D\).

Let \(F=\nabla\mathscr H_{n,J}^\tau(y_n^\dagger|_J)\). Outside the fixed neighborhoods of the boundary, Lemma~\ref{A:res:profile-boundary} gives \(|F_i|\le C/n\). In those
neighborhoods the logarithmic term at a free coordinate is unchanged when the exterior is replaced by the local minimizer; only the polynomial terms involving exterior coefficients
change. Lemma~\ref{A:res:11.6A} bounds those changes by \(C\Delta_{\mathrm L}\) and
\(C\Delta_{\mathrm R}\), respectively. Missing paths at the original boundary point have uniformly bounded weights. We can
therefore decompose \(F=F_0+F_{\mathrm L}+F_{\mathrm R}\), where \(\|F_0\|_\infty\le C/n\), and the other two vectors have supports of fixed size at the boundary points and norms
at most \(C\Delta_{\mathrm L}\) and \(C\Delta_{\mathrm R}\). The mean-value identity is \(y_J^\circ-y_n^\dagger|_J=-K^{-1}F\). Let $\nu(i)$ denote the site of coefficient $i$.
There are at most two coefficient coordinates at each site. Increasing a fixed integer $w=w(V)$ if needed, the supports of $F_{\mathrm L}$ and $F_{\mathrm R}$ lie within $w$ sites
of $\ell$ and $m$, respectively. For a coordinate $i$ at site $k$, the three contributions obey
\[
 |(K^{-1}F_0)_i|
 \le \frac{C}{n}\sum_j e^{-c|k-\nu(j)|}
 \le \frac{2C}{n}\left(1+\frac{2}{e^c-1}\right),
\]
\[
 |(K^{-1}F_{\mathrm L})_i|
 \le C\Delta_{\mathrm L}\sum_{|\nu(j)-\ell|\le w}
                 e^{-c|k-\nu(j)|}
 \le C_w\Delta_{\mathrm L}e^{-c(k-\ell)}.
\]
For the right contribution, $|k-\nu(j)|\ge m-k-w$ on its support, whose cardinality is at most $2(2w+1)$. Thus
\[
 |(K^{-1}F_{\mathrm R})_i|
 \le2(2w+1)C\Delta_{\mathrm R}e^{cw}e^{-c(m-k)}.
\]
Substitution in the mean-value identity proves \eqref{A:eq:global-mode-profile}. The constants just displayed depend on $w,c,C$, not on the interval length. For an edge coordinate
at $k+\delta$, $|\delta|\le1/2$, the two factors satisfy
\[
 e^{-c(k+\delta-\ell)}\le e^{c/2}e^{-c(k-\ell)},\qquad
 e^{-c(m-k-\delta)}\le e^{c/2}e^{-c(m-k)}.
\]
These inequalities prove the bounds with the half site convention.
\end{proof}

Fix an off-diagonal coordinate $b_k$, and retain the logarithmic weight $w_{k,n}$ defined above. Let
\(S_k=|e_k\rangle\langle e_{k+1}|+|e_{k+1}\rangle\langle e_k|\), and let \(T^{(k,0)}\) be the Jacobi matrix obtained by
setting only the \(k\)-th off-diagonal coefficient equal to zero. Changing the sign of the single coefficient \(b_k\) is unitary conjugation by a diagonal sign matrix. Hence, with
all other coefficients fixed,
\(g_k(t):=\operatorname{Tr}V\bigl(T^{(k,0)}+tS_k\bigr)\) is even. By the trace Hessian identity
\eqref{A:eq:finite-trace-hessian},
\[
g_k''(t)
=
D^2\operatorname{Tr}V(T^{(k,0)}+tS_k)[S_k,S_k]
\ge
c_V\operatorname{Tr}S_k^2
=2c_V.
\eqtag{T8}
\]

Let \(Y_{n,\tau}^\circ\) be any conditional minimizer for which the coordinate \(b_k\) is retained and \(w_{k,n}\!>0\).

Let $J_b$ index the free off-diagonal coordinates. For the actual coefficient $b_k$, write $\mathscr F_{n,k}^{(b)}$ for the sigma-field generated by all other coefficients; in this
conditional law $\beta(n-k)-1$ may be negative.

\begin{lemma}\label{A:res:T2}
Every free off-diagonal minimizer satisfies
\[
 (b_k^\circ)^2\le w_{k,n}/c_V.
\eqtag{T9}
\]
For an actual edge $k=n-m$, its exact conditional law satisfies
\[
 \mathbb P\left(\sqrt n\,b_k>C_{V,\beta}(\sqrt m+t)
       \,\middle|\,\mathscr F_{n,k}^{(b)}\right)
 \le2\mathrm e^{-c_{V,\beta}t^2},\qquad t\ge0.
\eqtag{T11}
\]
\end{lemma}

In particular, the minimizing values obey
\[
 \max_{\substack{n-L\le k\le N_n\\k\in J_b}} b_k^\circ
 \le C\sqrt{(1+L)/n},\qquad L\ge1,
\eqtag{T10}
\]
with an empty maximum equal to zero. This does not estimate a fixed cut coefficient. For the actual coefficients, every deterministic sequence $1\le L_n=O(\log n)$ satisfies
\[
 \max_{n-L_n\le k\le n-1} b_k
 \le C_\omega\sqrt{\frac{L_n+\log(n+1)}n}
\eqtag{T12}
\]
almost surely for all sufficiently large $n$. The logarithm remains necessary when $L_n$ is bounded.

\begin{proof}
At a free minimizer the Euler equation is \(g_k'(b_k^\circ)=2w_{k,n}/b_k^\circ\). Since $g_k'(0)=0$, \eqref{A:eq:T8} gives $g_k'(b)=\int_0^b g_k''(t)\,\mathrm dt\ge2c_Vb$. At
$b=b_k^\circ$, stationarity therefore implies $2c_V(b_k^\circ)^2\le2w_{k,n}\le2(n-k)/n$, which proves \eqref{A:eq:T9} and \eqref{A:eq:T10}. No Euler equation is imposed on a fixed
edge.

For the exact law, condition on every coefficient except \(b_k\). Its density is proportional to
\[
 b^{\beta m-1}\exp\{-\beta n g_k(b)/2\}\mathbf1_{b>0}.
\]
The exponent \(\beta m-1\) is greater than \(-1\), also when it is negative. Since \(g_k'(b)-2c_Vb\ge0\), the ratio of this density to
\(r_{m,n}(b)\propto b^{\beta m-1}\mathrm{e}^{-\beta nc_Vb^2/2}\) is nonincreasing. Let $w$ denote the density ratio and let $X,Y$ be independent with law $r_{m,n}$. Monotonicity
gives
\[
 2\operatorname{Cov}\bigl(w(X),\mathbf1_{\{X>t\}}\bigr)
 =\mathbb E\bigl[(w(X)-w(Y))
       (\mathbf1_{\{X>t\}}-\mathbf1_{\{Y>t\}})\bigr]\le0.
\]
Therefore
\[
 \mathbb P(b_k>t\mid\mathscr F_{n,k}^{(b)})
 =\frac{\mathbb E_{r_{m,n}}[w(b)\mathbf1_{\{b>t\}}]}
 {\mathbb E_{r_{m,n}}w(b)}
 \le\mathbb P_{r_{m,n}}(b>t).
\]
Thus $r_{m,n}$ dominates every conditional law under consideration.

We give the tail argument also for noninteger \(\beta m\). If \(X\) has the gamma law with shape \(\nu/2\) and scale two, then for \(0<\lambda<1/2\),
\[
 \log\mathbb E \mathrm{e}^{\lambda(X-\nu)}
 =\frac\nu2[-\log(1-2\lambda)-2\lambda]
 \le\frac{\nu\lambda^2}{1-2\lambda}.
\]
For \(u>0\), Markov's inequality with \(\lambda=\sqrt u/(\sqrt\nu+2\sqrt u)\) gives \(\mathbb P(X>\nu+2\sqrt{\nu u}+2u)\le \mathrm{e}^{-u}\). Since \((\sqrt\nu+\sqrt{2u})^2\) is at
least this threshold,
\[
 \mathbb P\left(
 b_{n-m}>\frac{\sqrt{\beta m}+\sqrt{2u}}
                  {\sqrt{\beta c_Vn}}
 \,\middle|\,\mathscr F_{n,n-m}^{(b)}\right)\le \mathrm{e}^{-u}.
\eqtag{terminal-chi-tail}
\]
The argument applies to every real \(\nu=\beta m>0\) and proves \eqref{A:eq:T11}.

Integrate \eqref{A:eq:terminal-chi-tail} before taking a union over edges. For every \(Q>0\) one can choose a deterministic \(C_Q\) so that
\[
 \mathcal B_{n,Q} =\left\{b_{n-m}\le C_Q \sqrt{\frac{m+\log(n+1)}n},\quad1\le m\le
 n-1\right\},\qquad \mathbb P(\mathcal B_{n,Q}^{\,c}) \le C_Qn^{-Q}.
\eqtag{terminal-edge-event}
\]
Indeed take \(u=(Q+3)\log(n+1)\) for each edge and sum at most \(n\) probabilities. The constants are uniform for
\(\beta\in[\beta_-,\beta_+]\Subset(0,\infty)\). Borel--Cantelli, with \(Q>1\), proves \eqref{A:eq:T12} on any coupling
in which each finite matrix has its original coefficient law. Squaring \eqref{A:eq:T12} gives
\[
 \max_{n-L_n\le k\le n-1}b_k^2
 \le C_\omega\frac{L_n+\log(n+1)}n.
\]
Its right side is $O_\omega(L_n/n)$ when $L_n\asymp\log n$ and $O_\omega(\log n/n)$ when $L_n$ is bounded.
\end{proof}

\begin{lemma}\label{A:res:diagonal-path} Suppose that all coefficients within $R_V$ indices of $k$ have modulus at most
\(M\). With absent edges at the original matrix boundary points set to zero,
\[
 \big|[V'(T)]_{kk}-V'(a_k)\big|
       \le C_{V,M}(b_{k-1}^2+b_k^2).
\eqtag{diagonal-path-bound}
\]
If \(a_k\) is a free diagonal of a conditional minimizer, then
\[
 |a_k^\circ-a_{\mathrm t}|
       \le C_{V,M}(\widetilde b_{k-1}^{\,2}
                         +\widetilde b_k^{\,2}),
\eqtag{diagonal-localization}
\]
where $\widetilde b_i$ is the minimizing value on a free edge and the boundary value on a fixed edge.
\end{lemma}

\begin{proof}
Write $d=\deg V$, $V'(t)=\sum_{m=0}^{d-1}v_mt^m$, and $M_1=\max\{1,M\}$. For $1\le m\le d-1$, let $\mathcal P_m(k)$ consist of paths $\pi=(i_0,\ldots,i_m)$ in $\{1,\ldots,n\}$ with
$i_0=i_m=k$ and $|i_{r+1}-i_r|\le1$. Matrix multiplication gives
\[
 [T^m]_{kk}=\sum_{\pi\in\mathcal P_m(k)}
               \prod_{r=0}^{m-1}T_{i_r,i_{r+1}}.
\]
The constant path contributes $a_k^m$. For every other path, let $r$ be the first index with $i_{r+1}\ne k$. If $i_{r+1}=k-1$, return to $k$ requires a later crossing of the same
edge; its product therefore contains $b_{k-1}^2$. A first step to $k+1$ gives a factor $b_k^2$. There are at most $3^m$ paths, and the other $m-2$ factors have modulus at most
$M_1$. Thus, for $2\le m\le d-1$,
\[
 |[T^m]_{kk}-a_k^m|
 \le 3^mM_1^{m-2}(b_{k-1}^2+b_k^2).
\]
The difference vanishes for $m=0,1$. Summing with coefficients $v_m$ proves \eqref{A:eq:diagonal-path-bound}, with $C_{V,M}=\sum_{m=2}^{d-1}|v_m|3^mM_1^{m-2}$. A path crossing an
absent boundary point edge has zero weight under the stated convention.

At a free diagonal of a conditional minimizer, the logarithmic terms are independent of $a_k$ and stationarity gives $[V'(T)]_{kk}=0$. Since $V'(a_{\mathrm t})=0$ and $V''\ge c_V$,
\[
 V'(a_k^\circ)
 =(a_k^\circ-a_{\mathrm t})\int_0^1
 V''(a_{\mathrm t}+t(a_k^\circ-a_{\mathrm t}))\,\mathrm dt,
 \qquad |V'(a_k^\circ)|\ge c_V|a_k^\circ-a_{\mathrm t}|.
\]
Apply the first estimate with the free coefficients at their minimizing values and the boundary coefficients unchanged. Dividing by $c_V$ proves \eqref{A:eq:diagonal-localization}.
\end{proof}

Choose \(L_n=\lceil A_{\mathrm t}\log(n+1)\rceil\), with a fixed \(A_{\mathrm t}>0\). For the actual matrix let
\(I_n^{\mathrm{diag}}=\{n-L_n+1,\ldots,n\}\), and let \(\widehat a^\circ\) be the joint conditional minimizer over just
these diagonals, with every off-diagonal and every exterior diagonal fixed at its actual value. This is a different conditional problem from the initial segment or reduced
minimizer \(Y_{n,\tau}^\circ\); in particular, the latter's fixed terminal diagonals are never varied.

For each fixed \(R<\infty\), on bounded exterior data and on \(\mathcal B_{n,Q}\), the initial segment and reduced conditional minimizers satisfy
\[
 \max_{N_n-R\le k\le N_n}|a_k^\circ-a_{\mathrm t}|
       \le C_{R,Q}\frac{\log(n+1)}{n}.
\eqtag{T15}
\]
Their last free off-diagonal coefficients satisfy \eqref{A:eq:T10}; a fixed terminal cut edge instead satisfies \eqref{A:eq:terminal-edge-event}. For the actual matrix, there are
finite size events of complementary probability at most \(C_{R,Q}n^{-Q}\) on which
\[
 \max_{n-R\le k\le n}|a_k-a_{\mathrm t}|
 \le C_{R,Q}\left(\sqrt{\frac{\log(n+1)}n}
                         +\frac{\log(n+1)}n\right).
\eqtag{T16}
\]

The estimates in the next lemma and the two quantitative bounds above hold on finite size events whose complements have probability at most $C_{R,Q}n^{-Q}$. The constants are
uniform for $\beta\in I_\beta$ with the final $s_*$ indices omitted as in \eqref{A:eq:B.45}.

For an integer $R\ge1$ and $n\ge R$, let $T_R$ be the principal matrix of $T_n$ on $\{n-R+1,\ldots,n\}$.

\begin{lemma}\label{A:res:T3}
For every fixed $R$, almost surely, for all sufficiently large \(n\),
\[
 \lambda_{\max}(T_R)
       \le a_{\mathrm t}+o_\omega(1)
       \le E_\star-\frac{d_{\mathrm t}}2,
\eqtag{T17}
\]
and its scaled Dirichlet operator satisfies
\[
 \frac{h^{-2}}{b_\star}(E_\star I-T_R)
       \ge\frac{d_{\mathrm t}}{2b_\star}h^{-2}I.
\eqtag{T18}
\]
For the last fixed number of free sites of either conditional minimizer, using actual fixed off-diagonal coefficients where applicable,
\[
 E_\star-a_k^\circ-\widetilde b_{k-1}-\widetilde b_k
       \ge\frac{d_{\mathrm t}}2
\eqtag{T19}
\]
for all sufficiently large \(n\), under the same bounds on exterior data and the coefficient events specified above.
\end{lemma}

\begin{proof}
On \(\mathcal T_{n,J}(M_0)\), Lemma~\ref{A:res:12.1A} gives
\[
 \max_i|Y_{n,\tau,i}^\circ|\le C_{M_0,V}.
\eqtag{T20}
\]
Apply Lemma~\ref{A:res:12.1A} with only the coordinates in $I_n^{\mathrm{diag}}$ free. Its bound \eqref{A:eq:mode-growth} gives $\|\widehat a^\circ\|_\infty\le C_V(1+M_0)^{D_V}$ on
the stated exterior event.

For the minimizer for the initial or reduced block, apply the exact equations only at its free diagonals:
\[
 [V'(T(Y_{n,\tau}^\circ,\tau))]_{kk}=0.
\eqtag{T21}
\]
For \(N_n-R\le k\le N_n\), every incident free off-diagonal coefficient has square at most \(C_R/n\) by \eqref{A:eq:T9}. The possible incident fixed off-diagonal coefficient at the
right cut is \(b_{N_n}\), whose square is at most
\(C_Q\log(n+1)/n\) on \(\mathcal B_{n,Q}\). Equation~\eqref{A:eq:diagonal-localization} proves \eqref{A:eq:T15}. The
fixed diagonals to the right enter as bounded coefficients in the path polynomials. Equation~\eqref{A:eq:T21} is imposed only on the free diagonal coordinates.

For the actual matrix, let \(\mathscr F_n^{\mathcal W}\) be the sigma-field generated by all off-diagonal coefficients and all diagonals outside \(I_n^{\mathrm{diag}}\). Let
\(\mathcal U_{n,Q}\) be its measurable event on which the coefficients outside $I_n^{\mathrm{diag}}$ and within $R_V$ indices of it are bounded by \(M_0\) and the off-diagonal
coefficients incident to \(I_n^{\mathrm{diag}}\) obey the bounds in \eqref{A:eq:terminal-edge-event}. The complement has probability at most \(C\mathrm{e}^{-cn}+C_Qn^{-Q}\), by
Lemma~\ref{A:res:raw} and \eqref{A:eq:terminal-edge-event}. The free diagonal vector is unrestricted on this exterior event. On \(\mathcal U_{n,Q}\), compactness and
Lemma~\ref{A:res:diagonal-path}, applied to \(\widehat a^\circ\), give
\[
 \max_{k\in I_n^{\mathrm{diag}}}|\widehat a_k^\circ-a_{\mathrm t}|
       \le C_Q\frac{L_n+\log(n+1)}n
       \le C_{Q,A_{\mathrm t}}\frac{\log(n+1)}n.
\eqtag{terminal-diagonal-mode}
\]
Even the first diagonal of this window is covered: an excursion into the fixed left part contains its small incident off-diagonal coefficient squared.

Conditional on this exterior sigma-field, the diagonal density is strongly log-concave on \(\mathbb R^{L_n}\), with Hessian lower bound \(\kappa nI\), where \(\kappa=\beta c_V/2\).
Let $G\sim N(0,(\kappa n)^{-1}I_{L_n})$. Lemma~\ref{A:res:A.8} compares the distance from $\widehat a^\circ$ with $\|G\|_2$. The Gaussian exponential moment is
\[
 \mathbb E e^{\kappa n\|G\|_2^2/4}
 =\prod_{j=1}^{L_n}\mathbb E e^{N_j^2/4}=2^{L_n/2},
 \qquad N_j\sim N(0,1)\text{ independent}.
\]
Applying Markov's inequality at threshold $e^{\kappa nt^2/4}$ gives
\[
 \mathbb P\left(\|a_{I_n^{\mathrm{diag}}}-\widehat a^\circ\|_2>t
                     \,\middle|\,\mathscr F_n^{\mathcal W}\right)
       \le 2^{L_n/2}\mathrm{e}^{-\kappa nt^2/4}.
\eqtag{terminal-diagonal-radius}
\]
Choose $t=C_{Q,A_{\mathrm t}}\sqrt{\log(n+1)/n}$. Since $L_n\le C_{A_{\mathrm t}}\log(n+1)$, a sufficiently large constant makes the right side at most $n^{-Q-2}$. Integrate over
$\mathcal U_{n,Q+2}$ and add its exceptional probability by Lemma~\ref{A:res:disintegration}. Together with \eqref{A:eq:terminal-diagonal-mode}, this proves \eqref{A:eq:T16}. The
last $R+1$ sites are in this window for large $n$.

Fix \(R\) and \(Q\) for the finite size estimate. Intersect the bounded coefficient event under the original law,
\(\mathcal B_{n,Q+2}\), and the preceding actual diagonal event, with constants enlarged once. Call the result
\(\mathcal T_{n,R,Q}^{\mathrm{term}}\). Then
\[
 \mathbb P((\mathcal T_{n,R,Q}^{\mathrm{term}})^c)
       \le C_{R,Q}n^{-Q},
\eqtag{terminal-good-event}
\]
and \eqref{A:eq:T15}, \eqref{A:eq:T16} and the bounds on adjacent off-diagonal coefficients hold there with deterministic constants. Gershgorin's inequality gives
\[
 \lambda_{\max}(T_R)
 \le\max_{n-R<k\le n}\left(a_k+
 \mathbf1_{\{k>1\}}b_{k-1}+\mathbf1_{\{k<n\}}b_k\right)
 \le a_{\mathrm t}+C_{R,Q}\sqrt{\log(n+1)/n},
\]
where an edge omitted from $T_R$ contributes only a nonnegative term to this upper bound. Since \(d_{\mathrm t}>0\), this proves \eqref{A:eq:T17}--\eqref{A:eq:T18} for a
deterministic sufficiently large \(n\) on \(\mathcal T_{n,R,Q}^{\mathrm{term}}\). At any of the last fixed number of free diagonal sites, the preceding bounds also give
\[
 |a_k^\circ-a_{\mathrm t}|+\widetilde b_{k-1}+\widetilde b_k
 \le C_{R,Q}\left(\frac{\log(n+1)}n+
                      \sqrt{\frac{\log(n+1)}n}\right)=:\epsilon_{n,R,Q}.
\]
Here a fixed edge uses \eqref{A:eq:terminal-edge-event}, and a free minimizing edge uses \eqref{A:eq:T10}. Choose $n$ so large that $\epsilon_{n,R,Q}\le d_{\mathrm t}/2$. Then
$E_\star-a_k^\circ-\widetilde b_{k-1}-\widetilde b_k \ge d_{\mathrm t}-\epsilon_{n,R,Q}\ge d_{\mathrm t}/2$, which proves \eqref{A:eq:T19}. Taking \(Q>1\) and applying
Borel--Cantelli proves the almost sure conclusions. On a compact \(\beta\)-interval use \(\beta_-\) in \(\kappa\), the same $s_*$ omitted indices, and the uniform constants of the
preceding lemmas. These estimates are uniform in the stated parameters; the finite array coupling is constructed for each fixed $\beta$.
\end{proof}

For the remaining pathwise arguments we also work on the event \(\bigcap_{R\in\mathbb N}\liminf_n\mathcal T_{n,R,3}^{\mathrm{term}}\). Its intersection with the coefficient event
\(\Omega_{\mathrm{bd}}\) still has probability one and is again denoted by \(\Omega_{\mathrm{bd}}\). The probability estimates at each size are
\eqref{A:eq:finite-coefficient-probability} and \eqref{A:eq:terminal-good-event}.

For an edge coefficient set \(\alpha_{k,n}:=\frac{h^{-2}}{b_\star}b_k.\)
\begin{lemma}\label{A:res:12.1}
On the coefficient event \(\Omega_{\mathrm{bd}}\) defined above, the following statements hold for all sufficiently large \(n\). The first two are coefficient estimates; the
resolvent bound on the final coordinate block additionally uses Lemma~\ref{A:res:T3}. At the logarithmic left cut,
\[
0<c_{\mathrm L}\le \frac{b_k}{b_\star}\le C_{\mathrm L}<\infty,
\quad
0\le k\le r_n+2R_V+4,
\eqtag{12.13}
\]
and consequently \(c_{\mathrm L}h^{-2}\le \alpha_{r_n,n}\le C_{\mathrm L}h^{-2}\). At the terminal cut,
\[
b_{N_n}\le C_\omega n^{-1/2}(\log n)^C
\le C_\omega h^{3/2}(\log n)^C,
\quad
|\alpha_{N_n,n}|
\le C_\omega h^{-1/2}(\log n)^C.
\eqtag{12.14}
\]
If \(A_R\) is the Dirichlet restriction to the fixed terminal block, then
\[
A_R+\kappa_\omega\ge c_\omega h^{-2}I,
\quad
\sup_{z\in K}\|(A_R-z)^{-1}\|_{\mathrm{op}}
\le C_{\omega,K}h^2.
\eqtag{12.15}
\]
\end{lemma}

\begin{proof}

Since $r_n=O(\log n)$, we have $r_n+2R_V+4\le s_1n$ for all sufficiently large $n$. On $\mathcal E_n$, Lemma~\ref{A:res:raw} therefore gives
\[
 b_{\min}\le b_k\le M_0\quad(1\le k\le r_n+2R_V+4),
 \qquad b_0=b_\star.
\]
Dividing by $b_\star>0$ proves \eqref{A:eq:12.13}. The exceptional probabilities are summable.

For $k=n-m$ in the fixed terminal block, Lemma~\ref{A:res:T2} gives
\[
 \mathbb P\!\left\{\sqrt n\,b_k>C_{V,\beta}(\sqrt m+t)
       \,\middle|\,\mathscr F_{n,k}^{(b)}\right\}
 \le2e^{-c_{V,\beta}t^2}.
\]
Take $t=(3\log(n+1)/c_{V,\beta})^{1/2}$ and integrate the conditional bound. The block has a fixed number of indices and bounded $m$, so
\[
 \mathbb P\!\left\{\max_{N_n\le k<n}b_k
       >C_{V,\beta}n^{-1/2}\sqrt{\log(n+1)}\right\}
 \le C_{V,\beta}(n+1)^{-3}.
\]
Borel--Cantelli, $n^{-1/2}\asymp h^{3/2}$, and $\alpha_{N_n,n}=h^{-2}b_{N_n}/b_\star$ prove \eqref{A:eq:12.14}.

Lemma~\ref{A:res:T3} and its event \eqref{A:eq:terminal-good-event} give, almost surely for all sufficiently large $n$,
\[
 \max_k\left(|(T_R)_{kk}-a_{\mathrm t}|
                +\sum_{\ell\ne k}|(T_R)_{k\ell}|\right)
 \le d_{\mathrm t}/2.
\]
An absent edge contributes zero to the matrix and can only decrease the row sum. Gershgorin's inequality consequently gives
\[
 \lambda_{\max}(T_R)\le a_{\mathrm t}+d_{\mathrm t}/2
       =E_\star-d_{\mathrm t}/2,
 \qquad
 A_R=\frac{h^{-2}}{b_\star}(E_\star I-T_R)
       \ge\frac{d_{\mathrm t}}{2b_\star}h^{-2}I.
\]
For a compact $K$, choose $n$ large enough that $\sup_{z\in K}|z|+|\kappa_\omega|\le d_{\mathrm t}h^{-2}/(4b_\star)$. Then
\[
 A_R+\kappa_\omega I\ge\frac{d_{\mathrm t}}{4b_\star}h^{-2}I,
 \qquad
 \sup_{z\in K}\|(A_R-z)^{-1}\|_{\mathrm{op}}
 \le\frac{4b_\star}{d_{\mathrm t}}h^2.
\]
This is \eqref{A:eq:12.15}.

\end{proof}

The left Schur complement involves $\alpha-\alpha^2e_C^*(A_C-z)^{-1}e_C$. The two terms can both have order $h^{-2}$, so bounding them separately would lose their cancellation. We
will interpret their difference through a solution with boundary value one. Before doing so, we need coercivity of the retained quadratic form with Dirichlet boundary conditions,
obtained from Lemmas~\ref{A:res:12.1A} and~\ref{A:res:mean}.

Let \(Y_{n,\tau}^\circ\) be the conditional minimizer and let \(H_{n,\tau}^\circ\) be the scaled Jacobi operator obtained by substituting its coefficients into \eqref{A:eq:12.1}.
We first prove its form estimate on coefficient events with deterministic constants. Near the left cut the error from the local minimizer is small only as a form, whereas in the
terminal region a macroscopic potential gap replaces a uniform lower bound on the off-diagonal coefficients.

Use original coefficient indices \(k\) up to the cut, and reindex the reduced vertices by \(k=r+j\), \(1\le j\le M=N_n-r\). Thus \(x_j=(j-\tfrac12)h\) and the macroscopic parameter
at vertex \(j\) is
\[
 s_{j,n}:=s^a_{n,r+j}
   =\frac{b_\star}{\tau_\star}h^2
           \left(x_j+rh+\frac h\beta\right).
\eqtag{recentered-profile-coordinate}
\]
We retain the displacement \(rh\) separately: although it tends to zero, the displacement in the rescaled coefficient index \(r/n\) is larger than \(O(n^{-1})\). Let
\(J_0=\max\{j\le M:s_{j,n}\le s_0\}\), so \(J_0\asymp n\) and
\(x_{J_0}\asymp h^{-2}\).

Use the form $\ell_h$ on these $M$ reindexed vertices, with $u_0=u_{M+1}=0$. In the proof, $R_\partial=R_\partial(V,\beta)$ denotes a fixed number of terminal sites.

\begin{lemma}\label{A:res:12.1B}
For every $Q>0$ there are events $\mathcal F_{n,Q}$, with $\mathbb P(\mathcal F_{n,Q}^c)\le C_Qn^{-Q}$, such that, for all sufficiently large $n$, the reduced form at the
conditional minimizer satisfies on $\mathcal F_{n,Q}$
\[
 c_0\ell_h[u]-C_0\|u\|^2
 \le\mathfrak h_{n,\tau}^\circ[u]
 \le C_0\ell_h[u]+C_0\|u\|^2.
\eqtag{12.15d}
\]
The constants $c_0,C_0>0$ are deterministic.
\end{lemma}

The same bounds hold for principal Dirichlet restrictions. For $\omega\in\Omega_{\mathrm{bd}}$, they hold for all $n\ge n_0(\omega)$. For boundary values $\xi_-,\xi_+$, the proof
gives \eqref{A:eq:affine-center-control} after writing $u=v+e$ with an extension $e$ of those values. Constants and finite size probabilities are uniform for $\beta\in I_\beta$
when the final $s_*$ indices of \eqref{A:eq:B.45} are omitted.

\begin{proof}
Throughout, fixed off-diagonal coefficients in the conditional energy are the original coefficients. For estimates of the row potentials at the right boundary we also use the
original cut off-diagonal coefficient. We change auxiliary boundary coefficients only through identities that leave the complete Dirichlet forms unchanged.

First consider the conditional initial segment before removing the left block. Its exterior is bounded on \(\mathcal E_n\), so Proposition~\ref{A:res:profile-global}, including the
original discrepancy in the left boundary data, gives for every fixed \(|\nu|\le R_V+2\)
\[
 |Y_{n,\tau}^{\mathrm{pre},\circ}(r+\nu)
             -Y_n^\dagger(r+\nu)|
       \le Cn^{-1}+C\mathrm{e}^{-cr}.
\eqtag{12.15e0}
\]
The right boundary of the fixed regular block is at a positive macroscopic distance from $r+\nu$. Thus its site distance $d_{n,R}$ satisfies $d_{n,R}\ge c_Rn$ for some fixed
$c_R>0$ and all sufficiently large $n$. Since
$\sup_{t\ge0}te^{-cc_Rt}=1/(\mathrm e cc_R)$,
\[
 e^{-cd_{n,R}}\le e^{-cc_Rn}\le\frac1{\mathrm e cc_Rn}.
\]
For the left boundary, choose $A_0\ge1/c$. Then $r\ge A_0\log n$ gives $e^{-cr}\le n^{-cA_0}\le n^{-1}$. The finite event \(\mathcal D_{n,Q}\) of
\eqref{A:eq:finite-coefficient-event} includes the initial segment fluctuation estimate, hence
\[
 \max_{|\nu|\le R_V+2}
 |Y_n(r+\nu)-Y_{n,\tau}^{\mathrm{pre},\circ}(r+\nu)|
       \le C_Q\sqrt{\frac{\log(n+1)}n}.
\eqtag{12.15e1}
\]
Only the fixed coefficients on the left of the new cut are needed here; the larger fixed neighborhood includes all coefficients entering those equations. It follows that the actual
new left exterior differs from the local minimizer by at most \(C_Q\varepsilon_n\), where \(\varepsilon_n=\sqrt{\log(n+1)/n}\). For the right exterior and the local minimizer, only
boundedness is used. Apply Proposition~\ref{A:res:profile-global} to the reduced conditional law. In reindexed coordinates, with \(Y_n^\dagger(j)\) now denoting the local minimizer
at original site \(r+j\), it gives
\[
 |Y_{n,\tau}^\circ(j)-Y_n^\dagger(j)|
 \le Cn^{-1}+C_Q\varepsilon_n \mathrm{e}^{-c(j-1)}
                         +C\mathrm{e}^{-c(M-j)}.
\eqtag{12.15e}
\]
For an edge the distance from the corresponding vertex differs by at most $1/2$. For either boundary distance $d$ this changes the bound by at most the fixed factor $e^{c/2}$,
since $e^{-c(d-1/2)}=e^{c/2}e^{-cd}$. The right boundary distance is $M-j$. The constants in \eqref{A:eq:profile-global-residual} and \eqref{A:eq:global-mode-hessian} are uniform
up to that boundary.

Choose
\[
 L_n^{\mathrm L}=\lceil A_{\mathrm L}\log(n+1)\rceil
\eqtag{12.15e2}
\]
with \(A_{\mathrm L}\) large enough that the left exponential in \eqref{A:eq:12.15e} is \(O(n^{-1})\) after this index. On the regular middle region the distant right exponential
is also \(O(n^{-1})\), so
\[
 |Y_{n,\tau}^\circ(j)-Y_n^\dagger(j)|\le C_Qn^{-1},
       \quad L_n^{\mathrm L}\le j\le J_0+2.
\eqtag{12.15f}
\]
Enlarging the regular cutoff slightly accommodates the two extra sites. For the local minimizer, use the adjacent minimizing off-diagonal coefficients also at its auxiliary
boundary edges. On this regular region, smoothness and the midpoint choice of the diagonal parameter give
\(b^\dagger(s^b_{n,r+j-1})+b^\dagger(s^b_{n,r+j}) =2b^\dagger(s_{j,n})+O(n^{-2})\). The edge gap in \eqref{A:eq:profile-edge-gap} and the exact coordinate
\eqref{A:eq:recentered-profile-coordinate} therefore imply
\[
 \frac{B_j^\dagger}{b_\star}\ge c_0,
 \quad
 V_{j,n}^\dagger\ge c_0x_j-C_0,
 \quad
 |V_{j,n}^\dagger|\le C_0(1+x_j).
\eqtag{12.15g}
\]
For the upper bound, \(rh\le1\) for all sufficiently large \(n\) and \(E_\star-E(s)\le Cs\) on the regular interval. The coefficient error \(O(n^{-1})\) changes the scaled
potential by \(O(h^{-2}n^{-1})=O(h)\).

On the initial logarithmic strip one has only
\[
 \max_{j\le L_n^{\mathrm L}+2}
 |Y_{n,\tau}^\circ(j)-Y_n^\dagger(j)|
       \le C_Q\varepsilon_n
       \le C_Qh^{3/2}\sqrt{\log(n+1)}.
\eqtag{12.15f1}
\]
Define a comparison Jacobi matrix by replacing the diagonal coefficients at the conditional minimizer and internal off-diagonal coefficients with their local minimizing values for
indices at most \(L_n^{\mathrm L}\), leaving all later minimizing coefficients unchanged. Denote its scaled form by \(\mathfrak q_n\). For its left auxiliary coefficient choose the
local minimizing value \(b^\dagger(s^b_{n,r})\), and for its right auxiliary coefficient choose the actual
\(b_{N_n}\). The transition at \(L_n^{\mathrm L}+1\) has only an \(O(n^{-1})\) coefficient error, by
\eqref{A:eq:12.15f}. Thus every regular row of this comparison satisfies \eqref{A:eq:12.15g} up to an \(O(h)\) potential error. As operators, the difference from the matrix at the
conditional minimizer is supported on
\(\{1,\ldots,L_n^{\mathrm L}+1\}\), is tridiagonal, and has entries bounded by
\[
 C_Qh^{-2}\varepsilon_n
       \le C_Qh^{-1/2}\sqrt{\log(n+1)}.
\eqtag{12.15f2}
\]
Changing the left auxiliary coefficient has made no operator contribution: its two occurrences cancel by Lemma~\ref{A:res:auxiliary-boundary}. For \(u_0=0\), telescoping gives
\(u_j=h\sum_{i=0}^{j-1}D_hu_i\), and hence
\(\sum_{j\le L+1}|u_j|^2\le C h^2(L+1)^2\sum_i|D_hu_i|^2\). It follows that
\[
 |\mathfrak h_{n,\tau}^\circ[u]-\mathfrak q_n[u]|
 \le C_Q\varepsilon_n(L_n^{\mathrm L}+2)^2\|D_hu\|^2
 \le C_Qh^{3/2}(L_n^{\mathrm L}+2)^2
               \sqrt{\log(n+1)}\,\|D_hu\|^2
 =o(1)\|D_hu\|^2.
\eqtag{12.15f3}
\]
This is a relative form bound with deterministic constants on the specified finite size event.

For the far rows, \eqref{A:eq:profile-edge-gap} gives
\[
 E_\star-E(s_{j,n})\ge d_0:=E_\star-E(s_0)>0,
       \quad j>J_0.
\eqtag{12.15h}
\]
The global square-root estimate \eqref{A:eq:profile-increments} gives
\[
 \left|b^\dagger(s^b_{n,r+j-1})+
          b^\dagger(s^b_{n,r+j})-2b^\dagger(s_{j,n})\right|
       \le Cn^{-1/2}.
\eqtag{far-profile-row-error}
\]
This unscaled error tends to zero, so the macroscopic gap persists near the terminal point. Choose a fixed \(R_\partial\) large enough that the sum of the three errors at the right
boundary in \eqref{A:eq:12.15e} is at most \(d_0/4\) when \(j\le M-R_\partial\). This choice depends only on the fixed coefficient bounds. The other errors tend to zero, so with
the actual off-diagonal coefficients in the representation of the boundary terms,
\[
 V_{j,n}^\circ\ge c h^{-2},
       \quad J_0<j\le M-R_\partial.
\eqtag{12.15i}
\]
For the last \(R_\partial\) rows, Lemma~\ref{A:res:T3} instead controls the exact free diagonals and the minimizing or fixed incident off-diagonal coefficients. Thus
\[
 \widehat V_{j,n}^\circ
 :=\frac{h^{-2}}{b_\star}
   (E_\star-A_j^\circ-\widehat B_{j-1}^\circ-\widehat B_j^\circ)
 \ge\frac{d_{\mathrm t}}{2b_\star}h^{-2},
       \quad M-R_\partial<j\le M,
\eqtag{12.15j}
\]
where \(\widehat B_M^\circ=b_{N_n}\) is the original fixed right edge. At interior edges the hatted coefficients agree with the unhatted ones. The comparison form \(\mathfrak q_n\)
coincides with the minimizer form in the terminal region.

Its coefficients of $|D_hu_j|^2$ are nonnegative and bounded above everywhere, and are bounded below in the regular region. Its row potential is bounded below by \(c(1+x_j)-C\)
everywhere and in absolute value by \(C(1+x_j)\): use \eqref{A:eq:12.15g} near the edge, and bounded coefficients with \(1+x_j\asymp h^{-2}\) in the terminal region. Every
derivative edge outside the regular region, including the last Dirichlet edge, satisfies
\[
 |D_hu_j|^2
 \le2h^{-2}(|u_j|^2+|u_{j+1}|^2)
  \le C\big((1+x_j)|u_j|^2+(1+x_{j+1})|u_{j+1}|^2\big).
\eqtag{12.15k}
\]
At the right boundary the absent value is zero and its term is omitted; at the junction of the two intervals both coordinates already have \(1+x\asymp h^{-2}\). Sum this estimate
to obtain \(\ell_h[u]\le C(\mathfrak q_n[u]+\kappa\|u\|^2)\)

for a fixed \(\kappa\). Conversely the upper coefficient and potential bounds give \(\mathfrak q_n[u]+\kappa\|u\|^2\le C\ell_h[u]\). Write the two bounds just proved as
$c_1\ell_h[u]\le\mathfrak q_n[u]+\kappa\|u\|^2 \le C_1\ell_h[u]$, with $c_1>0$. Let $b_n^{\rm form}$ be the coefficient of $\|D_hu\|^2$ in \eqref{A:eq:12.15f3}; it tends to zero at
a deterministic rate on this event. Choose $n_0$ so that $b_n^{\rm form}\le c_1/2$ for $n\ge n_0$. Then
\[
 \frac{c_1}{2}\ell_h[u]
 \le \mathfrak h_{n,\tau}^\circ[u]+\kappa\|u\|^2
 \le (C_1+c_1/2)\ell_h[u].
\]
This is \eqref{A:eq:12.15d}.

All estimates used the finite coefficient event and a fixed terminal event. For example, with \(R_*\) larger than
\(s_\beta+R_\partial+2R_V+4\), take
\[
 \mathcal F_{n,Q}
  =\mathcal D_{n,Q}\cap\mathcal T_{n,R_*,Q}^{\mathrm{term}}.
 \quad
 \mathbb P(\mathcal F_{n,Q}^c)\le C_Qn^{-Q}.
\eqtag{confinement-good-event}
\]
The choices of \(R_\partial\), \(A_0\), and \(A_{\mathrm L}\) may be made uniformly on a compact positive
\(\beta\)-interval; \(Q\) changes only the tail constants and the sufficiently large deterministic threshold. By
Borel--Cantelli these bounds hold almost surely for all sufficiently large matrix sizes.

Finally, a principal Dirichlet restriction is obtained by zero extension, so both its quadratic form and $\ell_h$ are obtained by restriction. The cut edge contributions remain in
the quadratic form with Dirichlet boundary conditions. For nonzero boundary conditions, fix boundary values \(\xi_-\) and \(\xi_+\), and keep the original cut off-diagonal
coefficients \(b_r,b_{N_n}\). The resulting quadratic functional on $\mathbb C^M$ is
\[
 \mathfrak q_\xi[u] ={} \mathfrak h_{n,\tau}^\circ[u] -\frac{2h^{-2}}{b_\star}\operatorname{Re}
 (b_r\overline{\xi_-}u_1+ b_{N_n}\overline{\xi_+}u_M) +\frac{h^{-2}}{b_\star}
 (b_r|\xi_-|^2+b_{N_n}|\xi_+|^2).
\eqtag{affine-center-form}
\]
Its homogeneous quadratic part is $\mathfrak h_{n,\tau}^\circ[u]$. For an interior vector \(e\) chosen as an extension of these boundary data, put \(u=v+e\); \(v\) has zero
boundary values. The bounded positive shifted form from \eqref{A:eq:12.15d} gives

\(|\mathfrak h_{n,\tau}^\circ(e,v)| \le C\ell_h[e]^{1/2}\ell_h[v]^{1/2}\). Also $|v_1|,|v_M|\le h\|D_hv\|$. Put
$\Xi=(|\xi_-|^2+|\xi_+|^2)^{1/2}$. Expansion gives
\[
 \mathfrak q_\xi[v+e]-\mathfrak h_{n,\tau}^\circ[v]
 =2\operatorname{Re}\mathfrak h_{n,\tau}^\circ[e,v]
  -\frac{2h^{-2}}{b_\star}\operatorname{Re}
     (b_r\overline{\xi_-}v_1+b_{N_n}\overline{\xi_+}v_M)
  +\mathfrak q_\xi[e].
\]
The bounded cut coefficients imply $|\mathfrak q_\xi[e]|\le C(\ell_h[e]+h^{-2}\Xi^2)$. For each $\varepsilon>0$, the two terms involving $v$ have absolute sum at most
\[
 C\ell_h[e]^{1/2}\ell_h[v]^{1/2}+Ch^{-1}\Xi\ell_h[v]^{1/2}
 \le\varepsilon\ell_h[v]
      +C_\varepsilon(\ell_h[e]+h^{-2}\Xi^2).
\]
Taking $\varepsilon=c_0/2$ in the lower bound \eqref{A:eq:12.15d}, and using its upper bound, yields
\[
 \begin{split}
 \tfrac12c_0\ell_h[v]-C\|v\|^2-C_\xi(e)
  &\le\mathfrak q_\xi[v+e]
  \le C\ell_h[v]+C\|v\|^2+C_\xi(e),\\
 C_\xi(e)&=C\big(\ell_h[e]
                  +h^{-2}(|\xi_-|^2+|\xi_+|^2)\big).
 \end{split}
\eqtag{affine-center-control}
\]
Here $\ell_h[e]$ is evaluated on the interior vector extended by zero at both boundary points. For a principal subinterval, we use its original cut coefficients and evaluate
$\ell_h$ with the same zero boundary values. The homogeneous quadratic part therefore has the same coercivity constants, while the additive term accounts for the fixed boundary
values. The cancellation of the auxiliary boundary terms applies to homogeneous Dirichlet forms; the terms involving the nonzero boundary values are retained.
\end{proof}

Write $H_n^\circ=H_{n,\tau}^\circ$ for the matrix at the conditional minimizer in Lemma~\ref{A:res:12.1B}. Center the conditional coefficients at $Y_{n,\tau}^\circ$. In interior
coordinates write \(\alpha_j=A_j-A_j^\circ\) and
\(\delta_j=B_j-B_j^\circ\). Use the same auxiliary off-diagonal coefficients in the two Dirichlet representations, so
\(\delta_0=\delta_M=0\); the original exterior coefficients remain fixed in the conditional energy. Define
\[
\zeta_j=-\frac{\alpha_j+2\delta_j}{b_\star},
\quad
\eta_j=\frac{\delta_j}{b_\star h}.
\eqtag{12.15l}
\]
The corresponding matrix decomposition is proved in Proposition~\ref{A:res:13.1}. Set
\[
\bar\zeta_j=\zeta_j-\mathbb E_\tau\zeta_j,
\quad
S_m^\zeta=h^{-1}\sum_{j=1}^m\bar\zeta_j.
\]
Let $S_n^\zeta$ be the piecewise linear interpolation of $S_m^\zeta$ at $x=mh$, extended constantly beyond the last vertex. The following bound for the interpolated centered
partial sums holds under every conditional coefficient law, for every exterior value.

\begin{lemma}\label{A:res:12.1C}
For every \(\alpha\in(0,1/2)\) and \(1\le p<\infty\), under the exact conditional law with its actual exterior values fixed,
\[
\left\|
\sup_{q\ge0}(1+q)^{-\alpha}
\left(
\operatorname{osc}_{[q,q+1]}S_n^\zeta
+|S_n^\zeta(q+1)-S_n^\zeta(q)|
\right)
\right\|_{L^p}
\le C_{p,\alpha,V,\beta},
\eqtag{12.15n}
\]
\end{lemma}

For the mean estimates, require \(\tau\in\mathcal T_{n,J}(M_0)\). With deterministic constants, the scalar means satisfy
\[
|\mathbb E_\tau\zeta_j|
\le
\begin{cases}
C n^{-1},&j\le J_0,\\
C n^{-1/2}\sqrt{\log n},&j>J_0,
\end{cases}
\eqtag{12.15o}
\]
and, almost surely, for all sufficiently large \(n\),
\[
\|\eta_n\|_{\ell^\infty}
\le C h^{1/2}\sqrt{\log n}.
\eqtag{12.15p}
\]

\begin{proof}
By Lemma~\ref{A:res:12.1A}, the conditional law has logarithmic Sobolev constant at most \(C/n\). Every increment
\(S_s^\zeta-S_r^\zeta\) is a centered linear functional of the coefficient vector. Write $\ell_j$ for the coefficient
row defining $\zeta_j$. Each coordinate belongs to at most $D$ row supports, and $\|\ell_j\|_2\le C$, with $D,C$ fixed. Cauchy--Schwarz at each coordinate gives
\[
 \left\|\sum_{j=r+1}^s\ell_j\right\|_2^2
 \le D\sum_{j=r+1}^s\|\ell_j\|_2^2\le DC^2(s-r).
\] Therefore
\[
\|S_s^\zeta-S_r^\zeta\|_{\psi_2}
\le
Cn^{-1/2}h^{-1}(s-r)^{1/2}
\le C(h(s-r))^{1/2}.
\eqtag{12.15q}
\]
Choose $s_2$ with $s_0<s_2<s_1$. For large $n$, every coefficient within $R_V$ indices of a row $j\le J_0$ has parameter at most $s_2$. Apply Lemma~\ref{A:res:mean} with this
regular cutoff to the two coefficients in $\zeta_j$. This proves \eqref{A:eq:12.15o} for $j\le J_0$. On $j>J_0$, the rough bound in the same lemma gives the factor $\sqrt{\log n}$.
After multiplication by \(h^{-2}\), the far scalar mean is at most \(C h^{-1/2}\sqrt{\log n}\); since \(1+x_j\asymp h^{-2}\) there, its relative coefficient is
\(O_\omega(h^{3/2}\sqrt{\log n})=o_\omega(1)\).

Apply the conditional subgaussian version of Lemma~\ref{A:res:A.7} to \(\xi_{j,n}=h^{-1}\bar\zeta_j\), using \eqref{A:eq:12.15q}. The weighted block tail \eqref{A:eq:A.23-tail},
integrated under the conditional law, proves \eqref{A:eq:12.15n} from this increment estimate.

The finite size maximal estimate \eqref{A:eq:mode-maximum}, divided by \(b_\star h\), gives the bound on \(\eta_n\) with probability at least \(1-C_Qn^{-Q}\) on good exterior data.
Add \(C\mathrm{e}^{-cn}\) by \eqref{A:eq:integrated-tail} and use Borel--Cantelli to obtain \eqref{A:eq:12.15p}.
\end{proof}

The centered partial sums enter the matrix as a discrete derivative. We record the deterministic estimate before using their probability bounds. The estimate is invariant under
addition of a constant to the primitive and uses its local oscillations. It will also apply to the samples of a single Brownian path.

Let $E=(E_j)_{j=0}^{M}$ be real, and use the operator $\Gamma_E$ of Proposition~\ref{A:res:3.17}, with $J=J_h$. For complex test vectors with \(u_0=u_{M+1}=0\), the identity
\(2\operatorname{Re}((u_{j+1}-u_j) \overline{(u_{j+1}+u_j)/2})=|u_{j+1}|^2-|u_j|^2\) gives the equality
\[
 \langle u,\Gamma_Eu\rangle
 =-h^{-1}\sum_{j=0}^{M}E_j(|u_{j+1}|^2-|u_j|^2)
 =\sum_{j=1}^{M}\frac{E_j-E_{j-1}}h|u_j|^2.
\eqtag{16.28}
\]
Polarization gives the sesquilinear identity. Adding a constant to $E$ changes neither side, since $u_0=u_{M+1}=0$.

Let \(0<h\le1\) and \(M\ge1\). Extend \(E_j\) constantly for \(j>M\), and let \(S_E\) be its continuous piecewise linear interpolation at \(x=jh\). Assume
$\mathcal O_\alpha(S_E)\le A$ for a fixed $0<\alpha<1/2$, with $\mathcal O_\alpha$ as in Subsection~\ref{A:sec:finite-difference-approximation-of-the-stochastic-airy-operator}.

\begin{lemma}\label{A:res:16.1} For every \(\varepsilon>0\), with constants independent of
\(h,M,E\),
\[
 |\langle u,\Gamma_Eu\rangle|
 \le\varepsilon\ell_h[u]
   +C_{\varepsilon,\alpha}(1+A)^{2/(1-2\alpha)}\|u\|^2.
\eqtag{16.30}
\]
In particular the exponent of \(1+A\) is four when \(\alpha=1/4\).
\end{lemma}

\begin{proof}
Put \(k_h=\lceil h^{-1}\rceil\), so \(1\le k_hh\le2\), and introduce the moving average and its fluctuation,
\[
 \overline E_j=k_h^{-1}\sum_{t=j}^{j+k_h-1}E_t,
 \quad F_j=E_j-\overline E_j,
 \quad j\ge0.
\eqtag{moving-average}
\]
For $p=\lfloor jh\rfloor$ and $j\le t\le j+k_h-1$, the interval from $jh$ to $th$ meets at most the three unit intervals with indices $p,p+1,p+2$. Hence
\[
 |E_j-\overline E_j|
 \le k_h^{-1}\sum_{t=j}^{j+k_h-1}|E_j-E_t|
 \le A\sum_{q=p}^{p+2}(1+q)^\alpha
 \le C_\alpha A(1+jh)^\alpha.
\]
Applying this bound to the oscillation between $(j-1)h$ and $(j+k_h-1)h$, and using $k_hh\ge1$, gives
\[
 |F_j|\le C_\alpha A(1+jh)^\alpha,\quad
 \left|\frac{\overline E_j-\overline E_{j-1}}h\right|
 =\frac{|E_{j+k_h-1}-E_{j-1}|}{k_hh}
 \le C_\alpha A(1+jh)^\alpha\quad(j\ge1).
\eqtag{moving-average-bounds}
\]
These estimates include the terminal extension. With \(f_j=|u_j|^2\), \(f_0=f_{M+1}=0\), summation by parts is
\[
 \sum_{j=1}^{M}\frac{E_j-E_{j-1}}h f_j =\sum_{j=1}^{M} \frac{\overline E_j-\overline E_{j-1}}h
 f_j -\sum_{j=0}^{M}\frac{F_j}{h}(f_{j+1}-f_j).
\eqtag{moving-average-identity}
\]
For \(X_j=1+x_j\), neighboring vertex and edge weights are comparable uniformly for \(h\le1\). Since
\(|J_hu_j|^2\le(|u_j|^2+|u_{j+1}|^2)/2\), Cauchy--Schwarz gives, for every \(\delta>0\),
\[
 \begin{split}
 |\langle u,\Gamma_Eu\rangle|
 &\le C_\alpha A\sum_jX_j^\alpha|u_j|^2
  +2C_\alpha A\sum_{j=0}^{M}(1+jh)^\alpha
                     |D_hu_j|\,|J_hu_j|\\
 &\le\delta\|D_hu\|^2
   +C_{\delta,\alpha}(A+A^2)
                         \sum_jX_j^{2\alpha}|u_j|^2.
 \end{split}
\eqtag{moving-average-energy}
\]
Here \(X_j\ge1\) and \(\alpha>0\) permit \(X_j^\alpha\le X_j^{2\alpha}\). For \(0<\theta<1\), maximization of
\(tX^\theta-aX\) over \(X\ge0\) gives
\[
 tX^\theta\le aX
 +(1-\theta)\theta^{\theta/(1-\theta)}
       a^{-\theta/(1-\theta)}t^{1/(1-\theta)}
 \quad(t,a>0).
\]
Take \(\theta=2\alpha\), \(\delta=\varepsilon/2\), and \(a=\varepsilon/2\). Since \(A+A^2\le2(1+A)^2\), this proves \eqref{A:eq:16.30}.

\end{proof}

For boundary removal we use the following finite size form bounds, with a logarithmic shift. They follow from the tails in \eqref{A:eq:12.15n}. The shift uniform in $n$ is obtained
later in Theorem~\ref{A:res:16.2}, after the Brownian coupling has been constructed.

\begin{lemma}\label{A:res:core-finite} For every \(Q>0\) there are
coefficient events \(\mathcal V_{n,Q}\), with \(\mathbb P(\mathcal V_{n,Q}^c)\le C_Qn^{-Q}\), such that the following holds for all sufficiently large \(n\). On
\(\mathcal V_{n,Q}\), the actual reduced quadratic form with Dirichlet boundary conditions \(\mathfrak a_{\mathrm{red}}\) satisfies \eqref{A:eq:12.16} below with deterministic
\(c,C>0\) and
\(\kappa_{n,Q}\le C_Q(1+\log(n+1))^2\). The constants and probabilities are uniform for $\beta\in I_\beta$ with the final $s_*$ indices omitted as in
\eqref{A:eq:B.45}.
\end{lemma}

\begin{proof}
Let \(\mathcal O_n=\mathcal O_{1/4}(S_n^\zeta)\) be the weighted oscillation of $S_n^\zeta$ over intervals $[p,p+1]$, with the continuous interpolation and constant terminal
extension of Lemma~\ref{A:res:A.7}. The conditional increment estimate \eqref{A:eq:12.15q} and the subgaussian version of that lemma give
\[
 \mathbb P(\mathcal O_n>t\mid\tau_{\mathrm{bd}})
 \le C \mathrm{e}^{-ct^2},\quad t\ge1.
\]
This centered estimate is uniform in the exterior values; the estimates for the mean are used only when the exterior coefficients satisfy the stated bounds. Define
\[
 \mathcal V_{n,Q}
  =\mathcal F_{n,Q+2}\cap
       \{\mathcal O_n\le C_Q\sqrt{1+\log(n+1)}\}.
\eqtag{core-good-event}
\]
Here \(\mathcal F_{n,Q+2}\) is the coefficient and terminal event in \eqref{A:eq:confinement-good-event}. Integration of the conditional tail and a sufficiently large \(C_Q\) give
the required probability bound. On this event the form bound \eqref{A:eq:12.15d}, the mean bounds \eqref{A:eq:12.15o}, and the finite size bound
\(\|\eta_n\|_\infty\le C_Qh^{1/2}\sqrt{1+\log(n+1)}\) all hold with deterministic constants.

Apply Lemma~\ref{A:res:16.1} to $E_j=S_j^\zeta$ with $\alpha=1/4$. The identity $S_j^\zeta-S_{j-1}^\zeta=h^{-1}\bar\zeta_j$ gives
\[
 \left|h^{-2}\sum_j\bar\zeta_j|u_j|^2\right|
 \le\varepsilon\ell_h[u]
       +C_\varepsilon(1+\mathcal O_n)^4\|u\|^2.
\]
The regular scalar mean is at most $Ch^{-2}n^{-1}=Ch$. On the far interval, \eqref{A:eq:12.15o} and $1+x_j\asymp h^{-2}$ instead give
\[
 \frac{h^{-2}|\mathbb E_\tau\zeta_j|}{1+x_j}
 \le Ch^{3/2}\sqrt{\log(n+1)}=o(1).
\]
Adding the mean and decreasing the relative coefficient in the first estimate proves, for every fixed $\varepsilon>0$ and all sufficiently large $n$ on $\mathcal V_{n,Q}$,
\[
 \left|h^{-2}\sum_j\zeta_j|u_j|^2\right|
 \le\varepsilon\ell_h[u]
       +C_{\varepsilon,Q}(1+\mathcal O_n)^4\|u\|^2.
 \eqtag{12.15r}
\]
Lemma~\ref{A:res:16.1} applies to the fixed coefficient sequence on this event. Since $\mathcal O_n\le C_Q\sqrt{1+\log(n+1)}$ on the specified event, its additive constant is
$O_Q((1+\log(n+1))^2)$.

The first-order form obeys, on the same event,
\[
 |\langle u,\Gamma_{\eta_n}^{-}u\rangle|
 \le2\|\eta_n\|_\infty\|D_hu\|\,\|u\|
 \le\varepsilon\ell_h[u]+C_{\varepsilon,Q}\|u\|^2.
\eqtag{12.15s}
\]
Choose the two relative bounds smaller than \(c_0/4\) in \eqref{A:eq:12.15d} and use the decomposition \eqref{A:eq:12.15m}. This proves
\[
 c\ell_h[u]-\kappa_{n,Q}\|u\|^2 \le\mathfrak a_{\mathrm{red}}[u] \le
 C\ell_h[u]+\kappa_{n,Q}\|u\|^2,\qquad 0\le\kappa_{n,Q} \le C_Q(1+\log(n+1))^2.
\eqtag{12.16}
\]
The same inequalities hold for principal Dirichlet restrictions by zero extension. All conditional estimates used above have common constants on compact positive
\(\beta\)-intervals, which proves uniformity in $\beta$.
\end{proof}

Write $A_C=\widetilde H_n^{\mathrm{red}}$ before reindexing and use $\mathcal L_h$ from Subsection~\ref{A:sec:weighted-resolvent-estimates-and-comparison-principles}. On
\(\mathcal V_{n,Q}\), \(\mathsf H=A_C+(\kappa_{n,Q}+1)I\ge c\mathcal L_h\). For \(K\Subset\mathbb C\setminus\mathbb R\), functional calculus gives
\[
 \begin{split}
 \|\mathcal L_h^{1/2}(A_C-z)^{-1}\mathcal L_h^{1/2}\|_{\mathrm{op}}
 &\le c^{-1}\sup_{\lambda\in\sigma(A_C)}
       \frac{\lambda+\kappa_{n,Q}+1}{|\lambda-z|}\\
 &\le C_K(1+\kappa_{n,Q})
 \le C_{Q,K}(1+\log(n+1))^2.
 \end{split}
\eqtag{core-form-resolvent}
\]
Indeed, factor the left side through \(\mathsf H^{1/2}\) and use \(\|\mathcal L_h^{1/2}\mathsf H^{-1/2}\|\le c^{-1/2}\);
\((\lambda+\kappa_{n,Q}+1)/(\lambda-z) =1+(z+\kappa_{n,Q}+1)/(\lambda-z)\). Borel--Cantelli gives the displayed
logarithmic bound for all sufficiently large $n$. Uniform pathwise form bounds for the coupled operators are proved in Theorem~\ref{A:res:16.2}.

The terminal block is removed by its Schur complement. Decompose
\[
\mathbb C^n=\mathcal H_P\oplus\mathcal H_R,
\quad
\mathcal H_P=\ell^2(\{1,\ldots,N_n\}),
\quad
\mathcal H_R=\ell^2(\{N_n+1,\ldots,n\}).
\]
With respect to this decomposition,
\[
\mathcal A_n-z
=
\begin{pmatrix}
A_P-z&-\alpha_R|e_P\rangle\langle e_R|\\
-\alpha_R|e_R\rangle\langle e_P|&A_R-z
\end{pmatrix},
\eqtag{12.17}
\]
where \(e_P=e_{N_n}\), \(e_R=e_{N_n+1}\), and \(\alpha_R=\alpha_{N_n,n}\). Put \(B_R(z)=(A_R-z)^{-1}, \allowbreak\, m_R(z)=e_R^*B_R(z)e_R.\) By Lemma~\ref{A:res:12.1} and the fixed
dimension of \(A_R\),
\[
\|B_R(z)\|_{\mathrm{op}}\le C_{\omega,K}h^2,
\quad
\|B_R(z)\|_{\mathfrak S_2}^2\le C_{\omega,K}h^4.
\eqtag{12.18}
\]
Eliminating the terminal coordinates changes the retained equation by the rank-one matrix
\(\Sigma_R(z)=\sigma_R(z)|e_P\rangle\langle e_P|, \allowbreak\, \sigma_R(z)=\alpha_R^2m_R(z)\), and
\eqref{A:eq:12.14}--\eqref{A:eq:12.18} give
\[
|\sigma_R(z)|
\le
C_{\omega,K}h(\log n)^C=o(1).
\eqtag{12.19}
\]
Let \(R_P(z)=(A_P-z)^{-1}\). Since \(z\) stays a fixed distance from the real axis,
\(\|R_P(z)\|_{\mathrm{op}}\le\delta_K^{-1}\). Hence, for all sufficiently large \(n\), \(S_P(z)^{-1} :=
(A_P-z-\Sigma_R(z))^{-1} = R_P(z)\bigl(I-\Sigma_R(z)R_P(z)\bigr)^{-1}\) exists and has operator norm bounded by
\(2\delta_K^{-1}\). The full block inverse is
\[
 (\mathcal A_n-z)^{-1}
 =\begin{pmatrix}
 S_P^{-1}&\alpha_R S_P^{-1}e_Pe_R^*B_R\\
 \alpha_RB_Re_Re_P^*S_P^{-1}&
 B_R+\alpha_R^2B_Re_R(e_P^*S_P^{-1}e_P)e_R^*B_R
 \end{pmatrix}.
\]
All terms are evaluated at $z$. Also
\[
 S_P^{-1}-R_P=\sigma_R S_P^{-1}e_Pe_P^*R_P.
\]
An outer product $uv^*$ has squared Hilbert--Schmidt norm $\|u\|^2\|v\|^2$. For $K$ enlarged by its complex conjugate, the uniform operator bounds therefore give
\[
 \begin{split}
 \|S_P^{-1}-R_P\|_{\mathfrak S_2}^2
 &\le C_K|\sigma_R|^2
     \le C_{\omega,K}h^2(\log n)^C,\\
 \|\alpha_RB_Re_Re_P^*S_P^{-1}\|_{\mathfrak S_2}^2
 &\le C_K\alpha_R^2\|B_R(z)e_R\|^2
     \le C_{\omega,K}h^3(\log n)^C,\\
 \|\alpha_R^2B_Re_R(e_P^*S_P^{-1}e_P)e_R^*B_R\|_{\mathfrak S_2}^2
 &\le C_K\alpha_R^4\|B_R(z)e_R\|^2\|B_R(\bar z)e_R\|^2
     \le C_{\omega,K}h^6(\log n)^C.
 \end{split}
\]
The other off-diagonal block has the same estimate. The four matrix blocks are orthogonal in Hilbert--Schmidt inner product. Combining their bounds, and using
$\|B_R\|_{\mathfrak S_2}^2\le C_{\omega,K}h^4$, proves
\[
 \sup_{z\in K}\left\|(\mathcal A_n-z)^{-1}
       -P_P^*(A_P-z)^{-1}P_P\right\|_{\mathfrak S_2}^2
 \le C_{\omega,K}\bigl(h^2(\log n)^C+h^3(\log n)^C+h^4\bigr)
 =o_\omega(h).
\eqtag{12.24}
\]

At the left cut, the off-diagonal entry has size $\alpha\asymp h^{-2}$. We compare the original leftmost coefficients with the half-line Jacobi matrix $J^V$ constructed in
\cite{krishnapur2016universality}. The comparison retains the dependence of its coefficients on their position and gives random coefficient error of scale \(n^{-1/2}\). The
estimates below use the spectral gap of its leading principal matrix after a rank-one modification.

\begin{lemma}\label{A:res:left-deterministic} There is a bounded
half-line Jacobi matrix \(J^V\), with diagonal \(a_j^V\) and positive off-diagonal \(b_j^V\), and constants
\(c,C,\gamma>0\), such that
\[
 |a_j^V-a_\star|+|b_j^V-b_\star|\le C\mathrm{e}^{-cj}, \quad \inf_j b_j^V>0,\qquad [V'(J^V)]_{jj}=0, \quad
 [V'(J^V)]_{j,j+1}=1/b_j^V,
\eqtag{left-profile-data}
\]
where \(a_\star=a^\dagger(0)\). For all sufficiently large integers \(r\), its original leading minor satisfies
\[
 E_\star I_r-J^V[1,r]-b_\star|e_r\rangle\langle e_r|
       \ge \gamma r^{-2}I_r.
\eqtag{left-profile-rankone}
\]
\end{lemma}

\begin{proof}

Apply Lemma~30 of \cite{krishnapur2016universality} with potential $V/2$. Its local energy is half of ours, so it constructs the same bounded half-line Jacobi matrix $J^V$ and
gives
\[
 |a_j^V-a_\star|+|b_j^V-b_\star|\le Ce^{-cj}.
\]
Its local stationarity equations are
\[
 [V'(J^V)]_{jj}=0,
 \qquad 2[V'(J^V)]_{j,j+1}-2/b_j^V=0.
\]
Every $b_j^V$ is positive. Choose $j_0$ so that $b_j^V\ge b_\star/2$ for $j\ge j_0$; then $\inf_jb_j^V\ge\min\{b_1^V,\ldots,b_{j_0}^V,b_\star/2\}>0$. Proposition~33 of the same
reference gives
\[
 \lambda_{\max}(J^V[1,m])
 \le E_\star-b_\star j_{0,0}^{\,2}m^{-2}+o(m^{-2}),
\eqtag{left-profile-dirichlet}
\]
where $j_{0,0}$ is the first positive zero of $J_0$.

Set $\widehat J=(J^V-a_\star I)/b_\star$ and split $\widehat J[1,4r]$ after vertex $r$. Write
\[
 q=\widehat J_{r,r+1},\qquad U=\widehat J[1,r],\qquad
 W=\widehat J[r+1,4r],\qquad m=3r.
\]
The coefficient estimate gives $|q-1|+\|W-D_m\|_{\mathrm{op}}\le Ce^{-cr}$, where $D_m$ has zero diagonal and unit adjacent entries. The row-sum bound proves the operator estimate
without a factor depending on $m$.

We first justify the comparison used at this cut. For any $t>0$, take normalized real top eigenvectors $u,v$ of $U+qt\,e_re_r^*$ and $W+(q/t)e_1e_1^*$. Their boundary coordinates
$u_r,v_1$ are nonzero: the tridiagonal eigenvalue equation would otherwise force every coordinate to vanish. Put $a=v_1$, $b=t u_r$. The identity
\[
 \begin{split}
 &\langle (U+qt\,e_re_r^*)au,au\rangle+
   \langle(W+(q/t)e_1e_1^*)bv,bv\rangle\\
 &\qquad-\langle\widehat J[1,4r](au,bv),(au,bv)\rangle
 =q\left|\sqrt t\,a u_r-\frac{b v_1}{\sqrt t}\right|^2=0
 \end{split}
\]
shows that the Rayleigh quotient of $(au,bv)$ is a weighted average of the two top eigenvalues. Therefore
\[
 \lambda_{\max}(\widehat J[1,4r])
 \ge\min\{\lambda_{\max}(U+qt\,e_re_r^*),
           \lambda_{\max}(W+(q/t)e_1e_1^*)\}.
 \eqtag{left-rankone-minmax}
\]
This proves the needed direction of the rank-one comparison directly.

Choose $t=q^{-1}$. The second matrix is $W+q^2e_1e_1^*$ and differs from $D_m+e_1e_1^*$ by $O(e^{-cr})$ in operator norm. Put $\theta=\pi/(2m+1)$ and
$v_j=\cos((j-\tfrac12)\theta)$, $1\le j\le m$. The fictitious values satisfy $v_0=v_1$ and $v_{m+1}=0$. The identity $v_{j-1}+v_{j+1}=2\cos\theta\,v_j$ consequently gives an
eigenvector of $D_m+e_1e_1^*$ with eigenvalue $2\cos\theta$. Its Rayleigh quotient, and then the perturbation bound, imply
\[
 \lambda_{\max}(W+q^2e_1e_1^*)
 \ge2\cos\frac{\pi}{6r+1}-Ce^{-cr}
 =2-\frac{\pi^2}{36r^2}+O(r^{-3})-Ce^{-cr}.
\]
On the other hand, \eqref{A:eq:left-profile-dirichlet} at size $4r$ gives
\[
 \lambda_{\max}(\widehat J[1,4r])
 \le2-\frac{j_{0,0}^{\,2}}{16r^2}+o(r^{-2}).
\]
The constants satisfy $j_{0,0}>2\pi/3$. This inequality can be checked without a numerical root approximation: for $t=x^2/4\le5/4$, the alternating series for $J_0(x)$ gives
\[
 J_0(x)\ge1-t+\frac{t^2}{4}-\frac{t^3}{36}
 \ge\frac{199}{2304}>0.
\]
The cubic polynomial is decreasing on $[0,5/4]$, and the omitted alternating terms have decreasing absolute values. Thus $j_{0,0}>\sqrt5>2\pi/3$. It follows that the second block
has top eigenvalue strictly larger than that of $\widehat J[1,4r]$ for every $r$ above a deterministic index. The minimum in \eqref{A:eq:left-rankone-minmax} must then be attained
by the first block. For a further deterministic index,
\[
 \lambda_{\max}(U+e_re_r^*)
 \le\lambda_{\max}(\widehat J[1,4r])
 \le2-\frac{j_{0,0}^{\,2}}{32r^2}.
\]
Multiplying by $b_\star$ and adding $a_\star$ proves
\[
 E_\star I_r-J^V[1,r]-b_\star e_re_r^*
 \ge\frac{b_\star j_{0,0}^{\,2}}{32r^2}I_r.
\]
This proves \eqref{A:eq:left-profile-rankone}, with a constant and starting index depending only on $V$.

\end{proof}

The rank-one shift in \eqref{A:eq:left-profile-rankone} is needed below: a lower bound for $A_L$ alone does not control the difference $\alpha-\alpha^2e_r^*(A_L-z)^{-1}e_r$. We
transfer that shifted spectral bound to the actual leading block. Put \(r=\lceil A_0\log n\rceil\), \(A_L=(h^{-2}/b_\star)(E_\star I_r-T_n[1,r])\), and
\(\alpha=(h^{-2}/b_\star)b_r\).

\begin{proposition}\label{A:res:left-spectral}
For every \(Q>0\) there are coefficient events \(\mathcal L_{n,Q}\) with \(\mathbb P(\mathcal L_{n,Q}^c)\le C_Qn^{-Q}\) such that, for all sufficiently large \(n\),
\[
 c h^{-2}\le\alpha\le C h^{-2},\quad
 Q_L:=A_L-\alpha e_re_r^*
       \ge c h^{-2}r^{-2}I_r.
\eqtag{left-spectral-gap}
\]
The constants are deterministic and can be chosen uniformly for \(\beta\) in a compact positive interval. The event $\mathcal L_{n,Q}$ is defined in terms of the coefficients.
\end{proposition}

\begin{proof}
Let \(\ell=2r+2R_V+4\), and condition on all coefficients outside \(a_1,\ldots,a_\ell,b_1,\ldots,b_{\ell-1}\), retaining the actual terms involving coefficients on both sides of
the boundary. On the event of bounded exterior data in Lemma~\ref{A:res:raw}, its conditional minimizer \(y^\circ\) is bounded, and its free off-diagonal coefficients have a
positive lower bound because \(\ell/n\to0\); use Lemma~\ref{A:res:12.1A}. As a comparison vector, restrict
\((a_j^V,b_j^V)\) to these free coordinates and keep the actual fixed exterior unchanged.

At an interior free coordinate the polynomial stationary equations in \eqref{A:eq:left-profile-data} cancel the logarithmic derivative with constant weight. The remaining
logarithmic term is \(2(k/n+1/(\beta n))/b_k^V=O(\ell/n)\). There is no boundary error at the original left boundary, because \(J^V\) has that same boundary point. Changing the
exterior on the right creates a bounded error supported within a fixed distance of \(\ell\). Thus, for the unscaled gradient at this comparison vector,
\[
 g=g_{\mathrm{bulk}}+g_{\mathrm R},\quad
 \|g_{\mathrm{bulk}}\|_\infty\le C\ell/n,\quad
 \|g_{\mathrm R}\|_\infty\le C,
\]
with \(g_{\mathrm R}\) supported at the right conditional boundary. The averaged Hessian between this vector and
\(y^\circ\) is \(P+D\ge cI\), where \(P\) is bounded and finite range and \(D\ge0\) is the logarithmic diagonal.
Lemma~\ref{A:res:A.9} and the gradient identity give, for \(j\le r+1\),
\[
 |y_j^\circ-y_j^V|
    \le C\ell/n+C\mathrm{e}^{-c(\ell-j)}
    \le Cr/n+C\mathrm{e}^{-cr}.
\eqtag{left-mode-comparison}
\]
This estimate applies to both coefficient types with that site label.

Lemma~\ref{A:res:mean}, integrated by Lemma~\ref{A:res:disintegration}, bounds the maximum deviation of the actual free coefficients from this mode by \(C_Q\sqrt{\log(n+1)/n}\),
with complementary probability at most \(C_Qn^{-Q}\). The event where the exterior data are not bounded adds only \(C\mathrm{e}^{-cn}\). In particular, if \(J_\ell\) denotes this
free initial segment, define
\[
 \mathcal L_{n,Q}
 =\left\{\tau_{J_\ell}\in\mathcal T_{n,J_\ell}(M_0),\quad
       \max_{i\in J_\ell}|Y_{n,i}-y_i^\circ(\tau_{J_\ell})|
          \le C_Q\sqrt{\frac{\log(n+1)}n}\right\}.
\]
On this event, tridiagonality and \eqref{A:eq:left-profile-data} imply
\[
 \|T_n[1,r]-J^V[1,r]\|_{\mathrm{op}}+|b_r-b_\star| \le C_Q\left( \sqrt{\frac{\log(n+1)}n}+\frac
 rn+\mathrm{e}^{-cr}\right) =o(r^{-2}).
\eqtag{left-random-comparison}
\]
For the matrix difference $E=T_n[1,r]-J^V[1,r]$, symmetry and tridiagonality give
\[
 \|E\|_{\mathrm{op}}\le\max_i\sum_j|E_{ij}|\le3\max_{i,j}|E_{ij}|.
\]
Since $r=\lceil A_0\log n\rceil$, multiplying the coefficient bound by $r^2$ gives a quantity tending to zero for fixed $Q,A_0$. Let $\epsilon_n$ denote the left side of
\eqref{A:eq:left-random-comparison}. The rank-one norm identity gives
\[
 \left\|T_n[1,r]+b_re_re_r^*-J^V[1,r]-b_\star e_re_r^*\right\|_{\mathrm{op}}
 \le\epsilon_n=o(r^{-2}).
\]
Thus \eqref{A:eq:left-profile-rankone} implies $Q_L\ge(h^{-2}/b_\star)(\gamma r^{-2}-\epsilon_n)I_r \ge(\gamma/(2b_\star))h^{-2}r^{-2}I_r$ for all sufficiently large $n$. Also
$b_r\to b_\star>0$, so $\alpha\asymp h^{-2}$. These are \eqref{A:eq:left-spectral-gap}. The deterministic matrix \(J^V\) is independent of \(\beta\), and every conditional estimate
above is uniform on compact positive \(\beta\)-intervals. This proves the uniformity in $\beta$.
\end{proof}

\begin{lemma}\label{A:res:left-resolvent} On the events in
Proposition~\ref{A:res:left-spectral}, let \(Z\subset\mathbb C\) be any fixed bounded set, possibly meeting the real axis. For all sufficiently large \(n\), every \(A_L-z\),
\(z\in Z\), is invertible. With \(e=e_r\),
\(B(z)=(A_L-z)^{-1}\), and \(\kappa(z)=\alpha-\alpha^2e^*B(z)e\),
\[
 \begin{gathered}
 \|B(z)\|_{\mathrm{op}}\le C_Zh^2r^2,\quad
 \|B(z)\|_{\mathfrak S_2}^2\le C_Zh^4r^5,\quad
 \|B(z)e\|^2\le C_Zh^4r^2,\\
 \operatorname{Re}\kappa(z)\ge c\alpha/r^2,
 \quad |\kappa(z)|\le2\alpha.
 \end{gathered}
\eqtag{left-resolvent-bounds}
\]
For real \(a\in Z\), one has the sharper upper bound \(0<\kappa(a)\le\alpha\).
\end{lemma}

\begin{proof}
Set \(g=c h^{-2}r^{-2}\) from \eqref{A:eq:left-spectral-gap}. Then \(Q_L\ge gI\) and \(A_L=Q_L+\alpha ee^*\ge gI\). Since \(g\to\infty\), uniformly on \(Z\),
\[
 B(z)=(I-zB(0))^{-1}B(0),\quad
 \|(I-zB(0))^{-1}\|_{\mathrm{op}}\le2.
\]
This proves the operator bound; dimension \(r\) gives the Hilbert--Schmidt bound. Put \(q=e^*Q_L^{-1}e>0\). The rank-one inverse identity gives
\[
 m(0)=e^*B(0)e=\frac{q}{1+\alpha q},
 \quad \kappa(0)=\frac{\alpha}{1+\alpha q},
 \quad 0<q\le g^{-1}.
\eqtag{left-sherman-morrison}
\]
In particular \(c\alpha/r^2\le\kappa(0)\le\alpha\), and \(m(0)\le\alpha^{-1}\). Positivity of \(B(0)\) gives
\[
 \|B(0)e\|^2\le\|B(0)\|_{\mathrm{op}}m(0)
       \le(g\alpha)^{-1}\le Ch^4r^2.
\]
The displayed factorization of \(B(z)\) transfers this bound for the boundary resolvent uniformly to \(Z\). The resolvent identity then gives
\[
 |\kappa(z)-\kappa(0)|
 \le\alpha^2|z|\,\|B(\bar z)e\|\,\|B(0)e\|
 \le C_Zr^2=o(\alpha/r^2),
\eqtag{left-schur-coefficient-continuation}
\]
because \(h^2r^4\to0\). This proves the bounds for the Schur complement coefficient at complex parameters. For real $a\in Z$, put $q_a=e^*(Q_L-a)^{-1}e$. Since $Q_L-a\ge(g/2)I$,
one has $0<q_a\le2/g$, and
\[
 e^*B(a)e=\frac{q_a}{1+\alpha q_a},\qquad
 \kappa(a)=\frac{\alpha}{1+\alpha q_a},\qquad
 \frac{\alpha}{1+2\alpha/g}\le\kappa(a)\le\alpha.
\]
This proves the estimates for real spectral parameters. These constants are uniform for $a$ in the bounded set $Z\cap\mathbb R$.
\end{proof}

For the principal matrix on indices $1,\ldots,n-s_\beta$ \(A_P\), split
\[
\mathcal H_P=\mathcal H_L\oplus\mathcal H_C,
\quad
\mathcal H_L=\ell^2(\{1,\ldots,r\}),
\quad
\mathcal H_C=\ell^2(\{r+1,\ldots,N_n\}).
\]
Let \(A_L\) and \(A_C\) be the two Dirichlet restrictions, put \(B_L(z)=(A_L-z)^{-1}\), \(R_C(z)=(A_C-z)^{-1}\), and
\(\alpha=\alpha_{r,n}\). If \(e_L=e_r\) and \(e_C=e_{r+1}\), define
\[
 \begin{aligned}
 p_L(z)&=\alpha B_L(z)e_L,
 &\kappa_L(z)&=\alpha-\alpha^2 e_L^*B_L(z)e_L,\\
 p_C(z)&=\alpha R_C(z)e_C,
 &d_C(z)&=\alpha-\alpha^2 e_C^*R_C(z)e_C.
 \end{aligned}
\]
The vectors $p_L(z)$ and $p_C(z)$ solve the restricted equations with boundary source $\alpha e_L$ and $\alpha e_C$, respectively. These are solutions with boundary value one on
the other side of the cut. Estimating $p_C=\alpha R_Ce_C$ only by $\|R_C\|_{\mathrm{op}}$ would give a squared bound of order $h^{-4}$. Instead, a boundary extension supported on
an interval of fixed physical length has energy $O(h^{-1})$. Subtracting that extension leaves a homogeneous Dirichlet problem, to which the finite size form estimate applies.

Work on $\mathcal V_{n,Q}$.

\begin{lemma}\label{A:res:core-poisson}
For \(K\Subset\mathbb C\setminus\mathbb R\), uniformly on \(K\),
\[
 \begin{aligned}
\|p_C(z)\|^2& \le C_{Q,K}h^{-1}(1+\log(n+1))^6,\qquad |d_C(z)| \le C_{Q,K}h^{-1}(1+\log(n+1))^4,\\
\|R_C(z)e_C\|^2& \le C_{Q,K}h^3(1+\log(n+1))^6.
\end{aligned}
\eqtag{core-poisson-bounds}
\]
The same bound as for \(\|p_C(z)\|^2\) holds for $\ell_h[p_C(z)]$, evaluated with the boundary values \(p_{C,0}=1\),
\(p_{C,M+1}=0\). It is not an estimate for the energy of the zero extension of its interior vector.
\end{lemma}

\begin{proof}
Reindex the reduced block by \(1,\ldots,M\), so \(e_C=e_1\). Keep the original cut off-diagonal coefficient \(b_r\). For interior vectors \(u,v\in\mathbb C^M\) and left boundary
values \(\xi,\eta\in\mathbb C\), define the sesquilinear form on $(\mathbb C^M\times\mathbb C)^2$ by
\[
 \mathfrak a_z^{\mathrm{aff}}((u,\xi),(v,\eta)) ={} \mathfrak a_{\mathrm{red}}[u,v]-z\langle
 u,v\rangle -\alpha(u_1\overline\eta+\xi\overline{v_1}) +\alpha\xi\overline\eta.
\eqtag{core-affine-form}
\]
At $z=0$ this form is Hermitian. Its quadratic form follows by expanding $\alpha|u_1-\xi|^2$; the term $\alpha|u_1|^2$ is already included in $\mathfrak a_{\mathrm{red}}[u]$. For
fixed $\xi$, the diagonal of \eqref{A:eq:core-affine-form} is a quadratic functional of $u$.

Take the cutoff \(\chi_j=(1-jh)_+\), with \(\chi_0=1\) and \(\chi_{M+1}=0\), and let \(v_0=v_{M+1}=0\) for a test vector.

The cutoff satisfies
\[
 \sum_{j=0}^{M}|D_h\chi_j|^2+
       \sum_{j=1}^{M}(1+x_j)|\chi_j|^2\le Ch^{-1},
\eqtag{core-affine-cutoff}
\]
where its derivative uses the boundary value \(\chi_0=1\). Indeed, $|D_h\chi_j|,|\chi_j|\le1$, and the support contains at most $h^{-1}+1$ edges with $x_j\le1$.

Put \(F_z(v)=\mathfrak a_z^{\mathrm{aff}}((\chi,1),(v,0))\). We claim
\[
 \begin{split}
 |F_z(v)|
   &\le C_{Q,K}h^{-1/2}(1+\log(n+1))\ell_h[v]^{1/2},\\
 |\mathfrak a_z^{\mathrm{aff}}((\chi,1),(\chi,1))|
   &\le C_{Q,K}h^{-1}(1+\log(n+1)).
 \end{split}
\eqtag{core-affine-forcing}
\]
Here and in what follows a cutoff denotes its interior vector when used in a matrix expression. On its support \(jh\le1\), the operator evaluated at the local minimizer has bounded
scaled potential: this follows by Taylor expansion at \(s=0\), with \(s=(r+j+\beta^{-1})/n\), retaining the shift \(rh\). The comparison with the conditional minimizer
\eqref{A:eq:12.15e}, including the actual exterior coefficient at \(j=0\), bounds the difference between its potential and the potential at the local minimizer by
\[
 C h+C_Qh^{-2}\varepsilon_n \mathrm{e}^{-cj}
       +Ch^{-2}\mathrm{e}^{-c(M-j)},\quad
 \varepsilon_n=\sqrt{\log(n+1)/n}.
\]
Its coefficients of $|D_hu_j|^2$ are bounded above. Thus the difference term and the bounded potential term at the local minimizer in the first estimate of
\eqref{A:eq:core-affine-forcing} are at most \(Ch^{-1/2}\ell_h[v]^{1/2}\), by \eqref{A:eq:core-affine-cutoff}. The localized left error is at most
\[
 C_Qh^{-2}\varepsilon_n\sum_{j\ge1}\mathrm{e}^{-cj}|v_j|
 \le C_Qh^{-1}\varepsilon_n\|D_hv\|,
\]
since \(|v_j|\le h\sqrt j\|D_hv\|\). This is smaller than \(C_Qh^{-1/2}\ell_h[v]^{1/2}\) for large \(n\). For the cutoff paired with itself, its contribution is
\(O(h^{-2}\varepsilon_n)=O(h^{-1/2}\sqrt{\log(n+1)})\), again within the second claimed bound. The right error is exponentially small on the cutoff support.

For the centered scalar fluctuation, the partial sums start at $S_0=0$ and obey \(\max_{jh\le1+h}|S_j|\le C\mathcal O_n\). Summation by parts gives
\[
 \begin{split}
 h^{-2}\sum_j\bar\zeta_j\chi_j\overline{v_j}
  &=-h^{-1}\sum_{j=0}^{M}S_j
       (\chi_{j+1}\overline{v_{j+1}}-\chi_j\overline{v_j}),\\
 \left|h^{-2}\sum_j\bar\zeta_j\chi_j\overline{v_j}\right|
  &\le C\mathcal O_n
       (\|D_h\chi\|\,\|v\|+\|\chi\|\,\|D_hv\|).
 \end{split}
\]
With $v=\chi$, the boundary products are $S_0|\chi_0|^2=0$ and $S_M|\chi_{M+1}|^2=0$. Hence
\[
 \left|h^{-2}\sum_j\bar\zeta_j|\chi_j|^2\right|
 \le2\max_{jh\le1+h}|S_j|\,\|D_h\chi\|\,\|\chi\|
 \le C\mathcal O_n h^{-1}.
\]  The scalar means are regular on this
support and contribute at most \(Ch\|\chi\|\,\|v\|\), or \(Ch\|\chi\|^2\) in the second estimate. Finally, the first-order form contributes at most
\[
 \|\eta_n\|_\infty
              (\|D_h\chi\|\,\|v\|+\|\chi\|\,\|D_hv\|),
\]
and the corresponding bound with \(v=\chi\). The decomposition \eqref{A:eq:12.15m} also holds on the enlarged space: the actual and minimizing cut coefficients coincide and
$\delta_0=\delta_M=0$, so their difference has no term involving the boundary values. The spectral term is bounded by \(|z|\|\chi\|\,\|v\|\). Using
\(\mathcal O_n\le C_Q\sqrt{1+\log(n+1)}\) proves \eqref{A:eq:core-affine-forcing}.

The equation \((A_C-z)p_C=\alpha e_1\) says \(\mathfrak a_z^{\mathrm{aff}}((p_C,1),(w,0))=0\) for every test vector
\(w\) with zero boundary values. Consequently \(v=p_C-\chi\) has zero boundary values and solves \((A_C-z)v=-F_z\) in
the form sense. The estimate for $|F_z(w)|$ in \eqref{A:eq:core-affine-forcing} bounds the dual norm of $F_z$ relative to $\ell_h$. Apply \eqref{A:eq:core-form-resolvent} to obtain
\[
 \ell_h[v]^{1/2}
    \le C_{Q,K}h^{-1/2}(1+\log(n+1))^3.
\eqtag{core-poisson-energy}
\]
Since $p_C=v+\chi$, the triangle inequality in the energy norm and the bounds for $v$ and $\chi$ give the stated estimates for $p_C$, with boundary values $(1,0)$.

The same solution controls $d_C$ without estimating its two defining terms separately. Evaluating the sesquilinear form at this solution gives
\[
 d_C(z) =\alpha(1-(p_C)_1) =\mathfrak a_z^{\mathrm{aff}}((p_C,1),(p_C,1)) =\mathfrak
 a_z^{\mathrm{aff}}((\chi,1),(\chi,1)) +\overline{F_{\bar z}(v)}.
\eqtag{core-schur-coefficient-affine}
\]
For the second equality substitute \((A_C-z)p_C=\alpha e_1\) into \eqref{A:eq:core-affine-form}; for the last use stationarity against \((v,0)\) and
\(\overline{\mathfrak a_{\bar z}^{\mathrm{aff}}(x,y)} =\mathfrak a_z^{\mathrm{aff}}(y,x)\). Enlarge \(K\) by its conjugate if needed. Equations \eqref{A:eq:core-affine-forcing} and
\eqref{A:eq:core-poisson-energy} prove the stated bound on \(d_C\). Finally \(\alpha\asymp h^{-2}\), by the cut coefficient estimates, so
\(\|R_C(z)e_C\|^2=\alpha^{-2}\|p_C(z)\|^2\) gives the last bound. All estimates are uniform on the stated coefficient events.
\end{proof}

\begin{lemma}\label{A:res:12.2}
Uniformly for \(z\in K\), almost surely for all sufficiently large \(n\),
\[
 \frac{c\alpha}{r^2}\le\operatorname{Re}\kappa_L(z),
 \quad |\kappa_L(z)|\le2\alpha,
\eqtag{12.27}
\]
\[
 \|B_L(z)\|_{\mathrm{op}}\le C_Kh^2r^2,
 \quad \|B_L(z)\|_{\mathfrak S_2}^2\le C_Kh^4r^5,
 \quad \|B_L(z)e_L\|^2\le C_Kh^4r^2,
\eqtag{12.28}
\]
and
\[
 |d_C(z)|\le C_{\omega,K}h^{-1}(\log n)^C,
 \quad \|R_C(z)e_C\|^2\le C_{\omega,K}h^3(\log n)^C.
\eqtag{12.29}
\]
\end{lemma}

In particular, for the diagonal matrix entries \(m_L(z)=e_L^*B_L(z)e_L\) and \(m_C(z)=e_C^*R_C(z)e_C\),
\[
 m_L(z)=\frac{\alpha-\kappa_L(z)}{\alpha^2},
 \quad m_C(z)=\frac{\alpha-d_C(z)}{\alpha^2}.
\eqtag{12.30}
\]
More precisely, these bounds hold with deterministic constants on \(\mathcal L_{n,Q}\cap\mathcal V_{n,Q}\) for all sufficiently large \(n\); the exponents for the reduced block can
be taken to be four and six, as in \eqref{A:eq:core-poisson-bounds}.

\begin{proof}
Lemma~\ref{A:res:left-resolvent}, applied with the rank-one spectral gap of the original left block, gives the left bounds for the boundary resolvent.
Lemma~\ref{A:res:core-poisson} gives the bounds for the reduced block, and the two scalar identities follow from the definitions. For $Q>1$,
\[
 \sum_n\mathbb P\bigl((\mathcal L_{n,Q}\cap\mathcal V_{n,Q})^c\bigr)
 \le\sum_n\bigl(\mathbb P(\mathcal L_{n,Q}^c)+\mathbb P(\mathcal V_{n,Q}^c)\bigr)
 \le C_Q\sum_n n^{-Q}<\infty.
\]
Borel--Cantelli gives a finite $n_0(\omega)$ such that both events occur for every $n\ge n_0(\omega)$, almost surely. Their probabilities depend only on the specified finite
coefficient laws.
\end{proof}

The coupling at the cut is $-\alpha e_Le_C^*$, with $\alpha\asymp h^{-2}$. We use the preceding bounds on $(A_L-z)^{-1}e_L$ and $(A_C-z)^{-1}e_C$, together with the denominator
$1-\alpha^2m_L(z)m_C(z)$ in the Schur complement formula.

With respect to \(\mathcal H_L\oplus\mathcal H_C\),
\[
A_P-z
=
\begin{pmatrix}
A_L-z&-\alpha|e_L\rangle\langle e_C|\\
-\alpha|e_C\rangle\langle e_L|&A_C-z
\end{pmatrix}.
\]
Set \(D_n(z)=1-\alpha^2m_L(z)m_C(z).\) The Schur complement formula for the rank-one coupling gives, with
\(x_z=B_L(z)e_L\) and \(y_z=R_C(z)e_C\),
\[
\begin{aligned}
P_C(A_P-z)^{-1}P_C^*& =R_C(z)+\frac{\alpha^2m_L(z)}{D_n(z)} |y_z\rangle\langle y_{\bar z}|,\qquad P_L(A_P-z)^{-1}P_C^* =\frac{\alpha}{D_n(z)} |x_z\rangle\langle y_{\bar z}|,\\
P_L(A_P-z)^{-1}P_L^*& =B_L(z)+\frac{\alpha^2m_C(z)}{D_n(z)} |x_z\rangle\langle x_{\bar z}|.
\end{aligned}
\eqtag{12.35}
\]
The denominator is controlled by the boundary estimates. From \eqref{A:eq:12.30},
\[
D_n(z)
=
\frac{\kappa_L(z)+d_C(z)}{\alpha}
-
\frac{\kappa_L(z)d_C(z)}{\alpha^2}.
\eqtag{12.36}
\]
The lower bound for the denominator has order $r^{-2}$, so the error must be $o(r^{-2})$. The exact identity \eqref{A:eq:12.36} gives
\[
 \left|D_n(z)-\frac{\kappa_L(z)}\alpha\right|
 =\left|\frac{d_C(z)}\alpha
             \left(1-\frac{\kappa_L(z)}\alpha\right)\right|
 \le\frac{3|d_C(z)|}\alpha
 \le C_{\omega,K}h(\log n)^C.
\]
Since $\operatorname{Re}\kappa_L/\alpha\ge c/r^2$, it follows that $\operatorname{Re}D_n\ge c/r^2-C_{\omega,K}h(\log n)^C$. The last error is $o(r^{-2})$, because $r=O(\log n)$.
Hence
\[
 \inf_{z\in K}|D_n(z)|\ge\frac{c}{r^2}
\eqtag{left-schur-denominator}
\]
for all sufficiently large \(n\) on the finite size events above. We estimate the three blocks in \eqref{A:eq:12.35}. The error after compression to the interior satisfies
\[
\sup_{z\in K}
\left\|
P_C(A_P-z)^{-1}P_C^*-R_C(z)
\right\|_{\mathfrak S_2}^2
\le
C_{\omega,K}h^2r^4(\log n)^C.
\eqtag{12.38}
\]
Indeed, \(|m_L|+|m_C|\le C/\alpha\), by \eqref{A:eq:12.27}--\eqref{A:eq:12.30}. Each rank-one coefficient in \eqref{A:eq:12.35} is therefore at most \(C\alpha r^2\). For vectors
\(u,v\), \(\||u\rangle\langle v|\|_{\mathfrak S_2}^2 =\|u\|^2\|v\|^2\). The correction to the reduced resolvent has factors $R_C(z)e_C$ and $R_C(\bar z)e_C$, each with squared norm
$O(h^3(\log n)^C)$. The cross terms have one factor $B_L(z)e_L$, with squared norm $O(h^4r^2)$, and one factor $R_C(\bar z)e_C$, with squared norm $O(h^3(\log n)^C)$, or their
adjoints. Thus the two off-diagonal resolvent blocks satisfy
\[
\sup_{z\in K}
\|P_L(A_P-z)^{-1}P_C^*\|_{\mathfrak S_2}^2
\le
C_{\omega,K}h^3r^6(\log n)^C.
\eqtag{12.39}
\]
Finally, the factors $B_L(z)e_L$ and $B_L(\bar z)e_L$ give
\[
\sup_{z\in K}
\left\|
P_L(A_P-z)^{-1}P_L^*-B_L(z)
\right\|_{\mathfrak S_2}^2
\le
C_{\omega,K}h^4r^8(\log n)^C,
\eqtag{12.40}
\]
and \(\|B_L(z)\|_{\mathfrak S_2}^2\le C h^4r^5\). Since \(r=O(\log n)\), every right hand side in \eqref{A:eq:12.38}--\eqref{A:eq:12.40}, as well as \(h^4r^5\), is \(o(h)\). We
have proved:

\par\Needspace{9\baselineskip}
\begin{proposition}\label{A:res:12.3}
Let \(A_P\) be the principal matrix on indices $1,\ldots,n-s_\beta$ and \(A_C\) its Dirichlet restriction to \(I_n\). Then
\[
\sup_{z\in K}
\left\|
(A_P-z)^{-1}-P_C^*(A_C-z)^{-1}P_C
\right\|_{\mathfrak S_2}^2
=o_\omega(h).
\]
\end{proposition}

\begin{proof}
Let $R=(A_P-z)^{-1}$ and $E=R-P_C^*R_CP_C$. Orthogonality of matrix blocks gives
\[
 \|E\|_{\mathfrak S_2}^2
 =\|P_CEP_C^*\|_{\mathfrak S_2}^2
   +\|P_LRP_C^*\|_{\mathfrak S_2}^2
   +\|P_CRP_L^*\|_{\mathfrak S_2}^2
   +\|P_LRP_L^*\|_{\mathfrak S_2}^2.
\]
The first three terms are bounded by \eqref{A:eq:12.38}--\eqref{A:eq:12.39}. For the last, \eqref{A:eq:12.40} and $(a+b)^2\le2a^2+2b^2$ give
$C_{\omega,K}(h^4r^8(\log n)^C+h^4r^5)$. Every bound is $o(h)$ because $r=O(\log n)$ and $h\asymp n^{-1/3}$. Taking the supremum over $K$ proves the uniform error bound.
\end{proof}

\begin{proposition}\label{A:res:12.4}
With \(A_P\) as above,
\[
\sup_{z\in K}
\left\|
(\mathcal A_{n,V,\beta}-z)^{-1}
-P_P^*(A_P-z)^{-1}P_P
\right\|_{\mathfrak S_2}^2
=o_\omega(h).
\]
\end{proposition}

\begin{proof}
Set $R=(\mathcal A_{n,V,\beta}-z)^{-1}$, $R_P=(A_P-z)^{-1}$, and $E=R-P_P^*R_PP_P$. Let $P_R$ restrict to the deleted terminal coordinates. Orthogonality of the domain and range
blocks gives
\[
 \|E\|_{\mathfrak S_2}^2
 =\|P_PEP_P^*\|_{\mathfrak S_2}^2
  +\|P_PRP_R^*\|_{\mathfrak S_2}^2
  +\|P_RRP_P^*\|_{\mathfrak S_2}^2
  +\|P_RRP_R^*\|_{\mathfrak S_2}^2.
\]
The four estimates preceding \eqref{A:eq:12.24} therefore yield
\[
 h^{-1}\sup_{z\in K}\|E\|_{\mathfrak S_2}^2
 \le C_{\omega,K}\bigl(h(\log n)^C+h^2(\log n)^C+h^3\bigr)
 \longrightarrow0,
\]
since $h\asymp n^{-1/3}$.
\end{proof}
Combining Propositions~\ref{A:res:12.3} and \ref{A:res:12.4} gives
\[
\sup_{z\in K}
\left\|
(\mathcal A_{n,V,\beta}-z)^{-1}
-P_{I_n}^*(H_n^{\mathrm{red}}-z)^{-1}P_{I_n}
\right\|_{\mathfrak S_2}^2
=o_\omega(h),
\eqtag{12.43}
\]
where \(H_n^{\mathrm{red}}=A_C\).

The bounds \eqref{A:eq:12.24} and \eqref{A:eq:12.38}--\eqref{A:eq:12.40} also hold on the coefficient events used to prove them. For each $Q>0$, take
\[
 \mathcal G_{\partial,n,Q}
   =\mathcal L_{n,Q+2}\cap\mathcal V_{n,Q+2},\quad
 \mathbb P(\mathcal G_{\partial,n,Q}^{\,c})\le C_Qn^{-Q}.
\eqtag{boundary-good-event}
\]
The event \(\mathcal V_{n,Q+2}\) includes the coefficient deviations, the bounds on the final coordinate block, and the bound on $S_n^\zeta$. Combining \eqref{A:eq:12.24},
\eqref{A:eq:core-poisson-bounds}, and \eqref{A:eq:12.38}--\eqref{A:eq:12.40}, we obtain on this event
\[
 \begin{split}
 &\sup_{z\in K}
 \left\|(\mathcal A_{n,V,\beta}-z)^{-1}
       -P_{I_n}^*(A_C-z)^{-1}P_{I_n}\right\|_{\mathfrak S_2}^2\\
 &\quad\le C_{Q,K}\left[
 h^2(1+\log(n+1))^{16}
 +h^3(1+\log(n+1))^{12}
 +h^4(1+\log(n+1))^{8}\right]=o(h).
 \end{split}
\eqtag{boundary-finite-rate}
\]
For example, the correction to the reduced resolvent is bounded by \(Ch^2r^4(1+\log(n+1))^{12}\), the off-diagonal resolvent block by \(Ch^3r^6(1+\log(n+1))^6\), and the left
correction by \(Ch^4r^8\). The uncorrected left block contributes \(Ch^4r^5\). The Schur complement correction from the terminal block is \(O(h\log(n+1))\) on these same events,
which is covered by the displayed bound. The constants and probabilities are uniform on compact positive \(\beta\)-intervals with the common terminal deletion. Set
\[
 \Omega_{\partial}=\Omega_{\mathrm{bd}}\cap
       \liminf_{n\to\infty}\mathcal G_{\partial,n,3}.
\]
Its probability is one, since the complements in \eqref{A:eq:boundary-good-event} are summable. For every outcome in
\(\Omega_{\partial}\), the boundary estimates hold for all sufficiently large \(n\) on every compact \(K\Subset\mathbb
C\setminus\mathbb R\). The events are defined by coefficients and do not depend on \(K\); only the constants and the threshold in \(n\) depend on \(K\).

Reindex the principal block on $I_n$ by the unitary $(U_ru)_j=u_{r+j}$, $1\le j\le M_n$. This is an exact change of coordinates:
\[
  (U_rH_n^{\mathrm{red}}U_r^*-z)^{-1}
 =U_r(H_n^{\mathrm{red}}-z)^{-1}U_r^*.
\]
The original coordinate is $x_{r+j}=x_j+rh$. In the later comparison with the linear potential $X_h$, this deterministic displacement contributes the scalar $\delta_n I$, where
$\delta_n=rh$.

For a real scalar $\delta$ and an operator $H$ with the uniform form bounds of Section~\ref{A:sec:weighted-resolvent-estimates-and-comparison-principles}, the resolvent identity
and Proposition~\ref{A:res:3.5} give

\[
 \begin{split}
 (H+\delta I-z)^{-1}-(H-z)^{-1}
   &=-\delta(H+\delta I-z)^{-1}(H-z)^{-1},\\
 \|(H+\delta I-z)^{-1}-(H-z)^{-1}\|_{\mathfrak S_2}^2
   &\le\delta^2\delta_K^{-2}\|(H-z)^{-1}\|_{\mathfrak S_2}^2
    \le C_{\omega,K}\delta^2.
 \end{split}
\]
In particular,
\[
 \sup_{z\in K}\|(H+\delta_n I-z)^{-1}-(H-z)^{-1}\|_{\mathfrak S_2}^2
 \le C_{\omega,K}h^2(\log n)^2=o_\omega(h).
\eqtag{12.47}
\]
This scalar potential comparison leaves the Brownian path and comparison maps fixed.

Let \(\mathcal P_n=\{a_1,\ldots,a_{N_n},b_1,\ldots,b_{N_n-1}\}\) and \(\mathcal T_n=\{a_{N_n+1},\ldots,a_n,b_{N_n},\ldots,b_{n-1}\}\). Condition on $\mathcal T_n=\tau$. All
coordinates in $\mathcal P_n$ remain free. For $k<N_n$, the choice of $s_\beta$ gives $1-k/n-1/(\beta n)>0$. Denote the conditional negative log density by
\[
\Phi_{n,\tau}^{\mathrm{pre}}(y)
=
\frac{\beta n}{2}
\mathscr H_n(y,\tau)
+
c_{n,\tau}.
\]
Here $y\in\mathbb R^{N_n}\times(0,\infty)^{N_n-1}$ and $c_{n,\tau}$ is independent of $y$.

\begin{proposition}\label{A:res:12.5}
For every admissible terminal value $\tau$ and every free coefficient vector $v$,
\[
 D^2\Phi_{n,\tau}^{\mathrm{pre}}(y)[v,v]
 \ge\frac{\beta c_Vn}{2}\|v\|_2^2.
\]
\end{proposition}

\begin{proof}

For a free direction $v=(v^a,v^b)$, the trace Hessian bound \eqref{A:eq:finite-trace-convexity} gives
\[
 \begin{aligned}
 D^2\Phi_{n,\tau}^{\mathrm{pre}}(y)[v,v]
 &=\frac{\beta n}{2}D^2\operatorname{Tr}V(T(y,\tau))[S(v),S(v)]
   +\sum_{k\text{ free}}\frac{\beta(n-k)-1}{b_k^2}(v_k^b)^2\\
 &\ge\frac{\beta c_Vn}{2}
          \left(\sum_k(v_k^a)^2+2\sum_k(v_k^b)^2\right)
 \ge\frac{\beta c_Vn}{2}\|v\|_2^2.
 \end{aligned}
\]
Here $S(v)$ is zero at the fixed coordinates and $\beta(n-k)-1>0$ at every free off-diagonal coordinate. Neither inequality restricts the exterior values.

\end{proof}

Let
\[
 y_{n,\tau}^{\star}=\operatorname*{argmin}_y\Phi_{n,\tau}^{\mathrm{pre}}(y),
 \qquad K_{n,\tau}=D^2\Phi_{n,\tau}^{\mathrm{pre}}(y_{n,\tau}^{\star}).
\]
The minimizer is unique, and $K_{n,\tau}\ge cnI$. On $\mathcal T_{n,\mathcal P_n}(M_0)$, Lemma~\ref{A:res:12.1A} bounds its coordinates. Separating the polynomial and logarithmic
parts of the Hessian gives the following decomposition, with
\(D_{n,\tau}\) diagonal:
\[
 K_{n,\tau}=P_{n,\tau}+D_{n,\tau},\qquad
 D_{n,\tau}\ge0,\qquad
 \sup_i\sum_j |(P_{n,\tau})_{ij}|\le Cn.
\]
Since $P_{n,\tau}$ has fixed range, Lemma~\ref{A:res:A.9} yields
\[
 |(K_{n,\tau}^{-1})_{ij}|\le Cn^{-1}\mathrm e^{-c|i-j|}.
\]
This is the estimate used for the covariance comparison in Section~\ref{A:sec:application-to-the-conditional-jacobi-coefficients}. An upper bound $K_{n,\tau}\le CnI$ is needed only
on regular coefficient blocks, where \eqref{A:eq:regular-mode-lower} controls the logarithmic Hessian; we do not use such an upper bound for the full initial segment.

Proposition~\ref{A:res:11.7} bounds the stationarity error of $Y_n^\dagger$ by $C/n$ in the regular interior. Conditioning on the first $r_n$ coefficients introduces a new boundary
at $r_n+1$. We first estimate the minimizer with only the terminal coefficients fixed, and then the minimizer with both boundary blocks fixed.

Condition first only on the terminal variables and let \(Y_{n,\tau}^{\mathrm{pre},\circ}\) be the minimizer of the remaining energy of the initial coefficients, which still
contains the left block \(\{1,\ldots,r_n\}\). Fix \(s_0<1\). Restrict this minimizer to the regular coefficient block ending at \(\lfloor s_0n\rfloor\); it is the conditional
minimizer on that block with boundary values inherited from the minimizer for the initial coefficients. On \(\mathcal T_{n,\mathcal P_n}(M_0)\), Lemma~\ref{A:res:12.1A} bounds this
minimizer for the initial coefficients and hence the inherited boundary data of the regular block, uniformly in \(n\). The boundedness event \eqref{A:eq:raw-probability} supplies
the required control of the actual terminal data. The regular lower bound follows from \eqref{A:eq:regular-mode-lower}. For every fixed \(|\nu|\le 2R_V+4\), the coefficient with
original label \(r_n+\nu\) is at distance \(r_n+O(1)\) from the original boundary and at distance of order \(n\) from the right boundary of this regular block.
Proposition~\ref{A:res:11.8} therefore gives
\[
\left|
Y_{n,\tau}^{\mathrm{pre},\circ}(r_n+\nu)
-Y_n^\dagger(r_n+\nu)
\right|
\le
C n^{-1}+C \mathrm{e}^{-c r_n}.
\eqtag{12.57}
\]
Choose \(A_0\) in \(r_n=\lceil A_0\log n\rceil\) so that
\[
\mathrm{e}^{-c r_n}\le n^{-10}.
\eqtag{12.58}
\]

Apply conditional concentration for the strongly log-concave initial segment law to the finitely many coordinates adjacent to the cut. On $\Omega_{\mathrm{bd}}$ this gives
\[
\max_{|\nu|\le 2R_V+4}
\left|
Y_{n,\tau}(r_n+\nu)
-Y_{n,\tau}^{\mathrm{pre},\circ}(r_n+\nu)
\right|
\le
C_\omega n^{-1/2}\sqrt{\log n}.
\eqtag{12.59}
\]
Indeed, \eqref{A:eq:mean-regular} first bounds the displacement of the conditional mean from the minimizer for the initial coefficients by \(C/n\). The tail bound about the mean
obtained from the logarithmic Sobolev inequality then gives \eqref{A:eq:mode-tail}. Integrating on \(\mathcal T_{n,\mathcal P_n}(M_0)\), adding \(C\mathrm{e}^{-cn}\) for its
complement, and choosing a summable polynomial tail proves \eqref{A:eq:12.59} by Borel--Cantelli. Combining \eqref{A:eq:12.57}--\eqref{A:eq:12.59}, the exterior data seen by the
reindexed interior interval at its new left boundary differ from the local minimizer by \(O_\omega\bigl(n^{-1/2}\sqrt{\log n}\bigr),\) a larger scale than
\(O(n^{-1})\).

Let \(Y_{n,\tau}^{\circ}\) be the conditional minimizer on the reindexed interior interval after the exterior cut variables have been fixed at their actual values. The mean-value
theorem for the unscaled conditional energy on the interior interval gives
\[
Y_{n,\tau}^{\circ}-Y_n^\dagger
=-\mathscr K_{n,\tau}^{-1}
\bigl(g_{\mathrm{bulk}}+g_{\mathrm L}+g_{\mathrm R}\bigr).
\eqtag{12.59a}
\]
Here \(\mathscr K_{n,\tau}\) is the averaged unscaled Hessian. By \eqref{A:eq:global-mode-hessian}, it is a bounded finite-range polynomial part plus a nonnegative logarithmic
diagonal, and is bounded below by \(cI\). Lemma~\ref{A:res:A.9}, uniformly in that diagonal, therefore gives
\[
|\mathscr K_{n,\tau}^{-1}(i,j)|
\le C_\omega \mathrm{e}^{-c|i-j|}.
\eqtag{12.59b}
\]
The three errors in \eqref{A:eq:12.59a} have the bounds
\[
\|g_{\mathrm{bulk}}\|_{\ell^\infty}\le Cn^{-1},
\eqtag{12.59c}
\]
\[
\operatorname{supp}g_{\mathrm L}
\subset\{1,\ldots,2R_V+4\},
\quad
\|g_{\mathrm L}\|_{\ell^\infty}
\le C_\omega n^{-1/2}\sqrt{\log n},
\eqtag{12.59d}
\]
and
\[
\operatorname{supp}g_{\mathrm R}
\subset\{M_n-2R_V-3,\ldots,M_n\},
\quad
\|g_{\mathrm R}\|_{\ell^\infty}\le C_\omega.
\eqtag{12.59e}
\]
The first estimate is \eqref{A:eq:profile-global-residual}, whose proof uses the squared off-diagonal coefficient \(q=(b^\dagger)^2\) up to $s=1$. The second follows from
\eqref{A:eq:12.59} and the local Lipschitz bounds for $\operatorname{Tr}V(T)$: only coefficients at most $R_V$ indices apart can occur in the same term. The third is the error from
the conditioned terminal coefficients,

which may be of order one because it is supported at $x\asymp h^{-2}$, where the resolvent has a smaller norm.

Convolving \eqref{A:eq:12.59b} with \eqref{A:eq:12.59c}--\eqref{A:eq:12.59e} gives the global minimizer estimate
\[
|Y_{n,\tau}^{\circ}(j)-Y_n^\dagger(j)| \le{} C_\omega n^{-1} +C_\omega n^{-1/2}\sqrt{\log n}\,\mathrm{e}^{-c(j-1)} +C_\omega \mathrm{e}^{-c(M_n-j)}.
\]
Near the new left boundary, the prefactor of the exponential term gives a larger error than \(O(n^{-1})\). The latter rate applies only after separation from that boundary.

Choose \(L_n^{\mathrm L}=\lceil A_{\mathrm L}\log n\rceil, \allowbreak\, J_n^{\mathrm L}=L_n^{\mathrm L}+3R_V+10\), with
\(A_{\mathrm L}\) large enough that \(n^{-1/2}\sqrt{\log n}\,\mathrm{e}^{-cL_n^{\mathrm L}}\le n^{-4}.\) Then, on every
fixed interval of macroscopic parameters bounded away from $s=1$, \(|Y_{n,\tau}^{\circ}(j)-Y_n^\dagger(j)| \le C_\omega n^{-1} \allowbreak\; (j\ge J_n^{\mathrm L}+1)\). On the left
boundary strip \(1\le j\le J_n^{\mathrm L}+R_V\), the weaker bound \(|Y_{n,\tau}^{\circ}(j)-Y_n^\dagger(j)| \le C_\omega n^{-1/2}\sqrt{\log n}\) is used together with the estimate
\eqref{A:eq:3.58g} for vectors vanishing at the left boundary. Near the right boundary, the resolvent tail contributes at \(O(h)\) scale.
Sections~\ref{A:sec:decomposition-of-the-general-jacobi-coefficients}--\ref{A:sec:comparison-of-the-deterministic-coefficients} use these separate estimates for the two boundaries.

The boundary comparison is \eqref{A:eq:12.43}, with its finite size version \eqref{A:eq:boundary-finite-rate}. The scalar displacement is estimated by \eqref{A:eq:12.47};
conditional convexity is Proposition~\ref{A:res:12.5}.

\subsubsection{Decomposition of the general Jacobi coefficients}\label{A:sec:decomposition-of-the-general-jacobi-coefficients}

After the boundary restriction, the retained Jacobi block must be written in the form used by the weighted perturbation estimates. Its cut coefficients remain part of the
conditional energy, even though they are not entries of the block. We retain $r=r_n$, $N=N_n$, $M=M_n=N-r$ and write $\kappa(j)=r+j$. The maps $D_h,J_-,J_+,I_h$ and the coordinates
$x_j$ follow Section~\ref{A:sec:weighted-resolvent-estimates-and-comparison-principles}, with $M$ in place of $N_n$. The first-order perturbation estimate is applied only after
Theorem~\ref{A:res:16.2} establishes the common form bounds.

The reindexed principal Jacobi block has diagonal entries \(A_j=a_{r+j}\), \(1\le j\le M\), and off-diagonal entries
\(B_j=b_{r+j}\), \(1\le j\le M-1\). Its actual exterior cut coefficients are \(\widehat B_0=b_r\) and \(\widehat
B_M=b_N\). They are not entries of this principal block, although they remain part of the exterior data of the coefficient energy. To represent the quadratic form with Dirichlet
boundary conditions, choose any real numbers
\(g_0,g_M\) and set \(B_j^g=B_j\) for \(1\le j\le M-1\), \(B_0^g=g_0\), and \(B_M^g=g_M\). The numbers $g_0,g_M$ are
auxiliary boundary coefficients in this representation; the original random variables $\widehat B_0,\widehat B_M$ remain fixed.

Let \(M\ge1\), \(h,b_\star>0\), and impose \(u_0=u_{M+1}=0\). For fixed retained coefficients \(A,B\), the form
\[
\mathfrak q_g[u]
=\sum_{j=0}^{M}\frac{B_j^g}{b_\star}|D_hu_j|^2
+\frac{h^{-2}}{b_\star}\sum_{j=1}^{M}
 (E_\star-A_j-B_{j-1}^g-B_j^g)|u_j|^2
\eqtag{auxiliary-boundary-form}
\] will be used with arbitrary auxiliary values $g_0,g_M$.

\begin{lemma}\label{A:res:auxiliary-boundary}
The form is independent of $g$, and
\[
\mathfrak q_g[u]
=\frac{h^{-2}}{b_\star}
 \left\{\sum_{j=1}^{M}(E_\star-A_j)|u_j|^2
 -2\operatorname{Re}\sum_{j=1}^{M-1}B_ju_{j+1}\overline{u_j}\right\}.
\eqtag{auxiliary-boundary-matrix}
\]
Thus it is the form of \(h^{-2}b_\star^{-1}(E_\star I-T(A,B))\), where \(T(A,B)\) contains only the retained coefficients. In particular, bounds on the complete form are
independent of the chosen auxiliary coefficients.
\end{lemma}

\begin{proof}
Expand the sum of squared differences. The zero boundary values give
\[
\sum_{j=0}^{M}B_j^g|u_{j+1}-u_j|^2
=\sum_{j=1}^{M}(B_{j-1}^g+B_j^g)|u_j|^2
 -2\operatorname{Re}\sum_{j=1}^{M-1}B_ju_{j+1}\overline{u_j}.
\]
The first sum cancels the two off-diagonal coefficient terms in the potential in \eqref{A:eq:auxiliary-boundary-form}, proving \eqref{A:eq:auxiliary-boundary-matrix}. Equivalently,
for two choices \(g,\widetilde g\), their difference is
\[
\mathfrak q_g[u]-\mathfrak q_{\widetilde g}[u] ={} \frac{g_0-\widetilde g_0}{b_\star}
 \bigl(|D_hu_0|^2-h^{-2}|u_1|^2\bigr) +\frac{g_M-\widetilde g_M}{b_\star}
 \bigl(|D_hu_M|^2-h^{-2}|u_M|^2\bigr)=0.
\eqtag{auxiliary-boundary-cancellation}
\]
Both identities also hold for \(M=1\): the interior edge sum is empty and both boundary contributions cancel at the single vertex. Polarization proves equality of the sesquilinear
forms and hence of their finite self-adjoint matrices.
\end{proof}

From now on we use \(g_0=g_M=b_\star\) and write \(B_0=B_M=b_\star\) only in the representation of the quadratic form with Dirichlet boundary conditions. The lemma, also applied
with \(g=(\widehat B_0,\widehat B_M)\), identifies this matrix as
\[
\widetilde H_n^{\mathrm{red}}=U_rH_n^{\mathrm{red}}U_r^*.
\eqtag{auxiliary-boundary-operator}
\]
Consequently, for every \(z\in\rho(\widetilde H_n^{\mathrm{red}})\),
\[
\left\|(\widetilde H_n^{\mathrm{red}}-z)^{-1}
 -U_r(H_n^{\mathrm{red}}-z)^{-1}U_r^*\right\|_{\mathfrak S_2}^2=0.
\eqtag{13.7}
\]
Equation~\eqref{A:eq:auxiliary-boundary-operator} only reindexes the principal block. Removing the exterior blocks still requires the resolvent comparison in
Section~\ref{A:sec:boundary-reduction}; the scalar shift by $rhI$ belongs to the later comparison of deterministic coefficients.

The Dirichlet condition in the lemma is essential. If instead \(u_0=\xi_0\) and \(u_{M+1}=\xi_{M+1}\) are imposed, the same calculation gives
\[
\mathfrak q_g[u]-\mathfrak q_{\widetilde g}[u] =\frac{h^{-2}}{b_\star}\bigl\{ (g_0-\widetilde
 g_0) (|\xi_0|^2-2\operatorname{Re}(u_1\overline{\xi_0})) +{} (g_M-\widetilde g_M)
 (|\xi_{M+1}|^2-2\operatorname{Re}(\xi_{M+1}\overline{u_M}))\bigr\}.
\]
Thus energies with nonzero boundary values and the Schur complements in Section~\ref{A:sec:boundary-reduction} must retain the original cut coefficients. The cancellation above
applies to homogeneous Dirichlet principal blocks.

Let \(Y_n=(A_1,\ldots,A_M,B_1,\ldots,B_{M-1})\) be the unchanged retained coefficient vector. Let \(\tau_{\mathrm{bd}}\) consist of \emph{all} remaining original Jacobi
coefficients, including the removed left block, the terminal variables, and the two original cut coefficients \(\widehat B_0=b_r\), \(\widehat B_M=b_N\). The auxiliary boundary
values are used in the extension; $Y_n$ and $\tau_{\mathrm{bd}}$ retain their original coordinates. Put \(\mathcal D_n^{\mathrm{int}}=\mathbb R^M\times(0,\infty)^{M-1}\).

\begin{corollary}\label{A:res:conditioning}
For exterior data with strictly positive off-diagonal coefficients, a version of the conditional law of \(Y_n\) is
\[
\mathrm{d}\mu_n^{\tau_{\mathrm{bd}}}(y) =\frac{1}{Z_n(\tau_{\mathrm{bd}})}
 \exp\{-\Phi_n^{\tau_{\mathrm{bd}}}(y)\} \mathbf1_{\mathcal
 D_n^{\mathrm{int}}}(y)\,\mathrm{d}y,\qquad \Phi_n^{\tau_{\mathrm{bd}}}(y) =\frac{\beta
 n}{2}\mathscr H_n(y,\tau_{\mathrm{bd}}),
\eqtag{conditional-density}
\]
up to a \(y\)-independent additive constant in \(\Phi\). Here \(\mathscr H_n\) is the full energy in \eqref{A:eq:11.1}, with all exterior coefficients held at their original
values. The law and its unique conditional minimizer are independent of $g_0,g_M$.
\end{corollary}

\begin{proof}
Split the original coefficient coordinates into \((y,t)\). The joint Lebesgue density is proportional to \(\exp\{-\beta n\mathscr H_n(y,t)/2\}\) on the original positive
coefficient domain. Fixing \(t\) and dividing by its integral over
\(y\) gives \eqref{A:eq:conditional-density}. The trace \(\operatorname{Tr}V(T(y,t))\) retains every interaction
crossing the boundary of the free block. In particular, terms involving coefficients on both sides of either cut retain the original exterior coefficients. Changing an auxiliary
boundary value in \eqref{A:eq:auxiliary-boundary-form} changes no coordinate of \((y,t)\); it therefore leaves the conditional energy and all its derivatives unchanged.

The free off-diagonal coefficient indices satisfy \(r+1\le k\le N-1\), so \(\beta(n-k)-1\ge\beta(s_\beta+1)-1>0\). Their logarithmic barriers are convex and diverge at zero.
Conditional strong convexity in Proposition~\ref{A:res:12.5}, or its restriction to these coordinates, gives \(\nabla_y^2\Phi_n^{\tau_{\mathrm{bd}}}\ge c_{V,\beta}nI\). Coercivity
and the barriers give a unique minimizer in \(\mathcal D_n^{\mathrm{int}}\) for each such exterior configuration.

Finally, let \(\nu_n(\mathrm{d}t)\) be the exact exterior marginal and let \(F\) be a bounded Borel function of the full coefficient vector. Disintegration gives
\[
\int\!\int F(y,t)\,\mu_n^t(\mathrm{d}y)\,\nu_n(\mathrm{d}t)
=\int F(a,b)\,\mu_{n,V,\beta}(\mathrm{d}a\,\mathrm{d}b).
\eqtag{disintegration}
\]
Neither factor involves $g_0,g_M$; the exterior variables $b_r,b_N$ remain those of the original ensemble.
\end{proof}

Write the conditional minimizer as \(Y_n^\circ=(A_1^\circ,\ldots,A_M^\circ, B_1^\circ,\ldots,B_{M-1}^\circ)\). Extend it in the quadratic form with Dirichlet boundary conditions by
the same auxiliary boundary values
\(B_0^\circ=B_M^\circ=b_\star\). The alternative extension using original exterior data is \(\widehat
B_j^\circ=B_j^\circ\) for \(1\le j\le M-1\), \(\widehat B_0^\circ=\widehat B_0\), and \(\widehat B_M^\circ=\widehat B_M\); the hats at the boundary points do not denote additional
minimizer coordinates. Use the centered coordinates $\alpha_j,\delta_j$ of \eqref{A:eq:12.15l}. Then \(\delta_0=\delta_M=0\) exactly because the two forms use the same auxiliary
boundary coefficients.

Lemma~\ref{A:res:13.2} bounds the displacement of the conditional mean from $Y_n^\circ$. The Gaussian approximation in
Subsection~\ref{A:sec:gaussian-approximation-of-the-conditional-gibbs-measure} is centered at $Y_n^\circ$ and uses the Hessian there. For internal edges and vertices put
$a_j^\dagger=a^\dagger(s^a_{n,r+j})$ and $b_j^\dagger=b^\dagger(s^b_{n,r+j})$; at the two auxiliary edges use the adjacent local profile values. Fix a deterministic cutoff
$s_0\in(0,1)$. The refined estimate for the influence of the boundary from Section~\ref{A:sec:boundary-reduction} gives
\[
|A_j^\circ-a_j^\dagger| + |B_j^\circ-b_j^\dagger| \le C_{\omega,s_0} n^{-1} +C_\omega n^{-1/2}\sqrt{\log n}\,\mathrm{e}^{-c(j-1)} +C_\omega \mathrm{e}^{-c(M-j)}
\eqtag{13.17}
\]
whenever every coefficient within $R_V$ indices lies in the regular region, \(s^b_{n,\kappa(j)+R_V}\le s_0.\) The first exponential is measured from the left conditional boundary
of the reindexed interval. In particular, this estimate need not be \(O(n^{-1})\) at the first bounded number of reindexed sites. Choose
\[
L_n^{\mathrm L}=\lceil A_{\mathrm L}\log n\rceil,
\quad
J_n^{\mathrm L}=L_n^{\mathrm L}+3R_V+10,
\eqtag{13.17b}
\]
with \(A_{\mathrm L}\) large enough that \(n^{-1/2}\sqrt{\log n}\,\mathrm{e}^{-cL_n^{\mathrm L}} \le n^{-4}.\) The left boundary contribution is then bounded by $C_\omega n^{-4}$
for \(j\ge J_n^{\mathrm L}+1\). The left boundary strip \(I_n^{\mathrm L}=\{1,\ldots,J_n^{\mathrm L}\}\) is retained exactly through Gaussian transport and local covariance
approximation; it will be replaced separately using the bound on $R(z)e_j$ near the left boundary. The terminal region is likewise retained and is estimated by the resolvent bound
on $\{x\ge c h^{-2}\}$.

Define the conditional minimizer operator \(H_n^\circ\) by the quadratic form
\[
\mathfrak h_n^\circ[u] ={} \sum_{j=0}^{M} \frac{B_j^\circ}{b_\star}|D_hu_j|^2 + \frac{h^{-2}}{b_\star} \sum_{j=1}^{M} \left( E_\star-A_j^\circ-B_{j-1}^\circ-B_j^\circ \right)|u_j|^2.
\eqtag{13.18}
\]
The same form with \(A_j,B_j\) in place of \(A_j^\circ,B_j^\circ\) represents the unchanged principal operator
\(\widetilde H_n^{\mathrm{red}}\):
\[
\mathfrak h_n^{\mathrm{int}}[u] ={} \sum_{j=0}^{M} \frac{B_j}{b_\star}|D_hu_j|^2 + \frac{h^{-2}}{b_\star} \sum_{j=1}^{M} \left( E_\star-A_j-B_{j-1}-B_j \right)|u_j|^2.
\eqtag{13.19}
\]
Subtracting \eqref{A:eq:13.18} from \eqref{A:eq:13.19}, we obtain
\[
\mathfrak h_n^{\mathrm{int}}[u] - \mathfrak h_n^\circ[u] ={} \sum_{j=0}^{M} \frac{\delta_j}{b_\star}|D_hu_j|^2 - \frac{h^{-2}}{b_\star} \sum_{j=1}^{M} \left( \alpha_j+\delta_{j-1}+\delta_j \right)|u_j|^2.
\eqtag{13.20}
\]
Expanding the fluctuation in the difference term in \eqref{A:eq:13.20} cancels the adjacent diagonal terms and leaves the fluctuating off-diagonal product
$-2h^{-2}\delta_j\operatorname{Re}(u_{j+1}\overline{u_j})/b_\star$, as in \eqref{A:eq:6.13}. For each off-diagonal fluctuation,
\[
-\frac{2h^{-2}}{b_\star}\delta_j \operatorname{Re}(u_{j+1}\overline{u_j}) ={}
 -\frac{2h^{-2}}{b_\star}\delta_j|u_j|^2 -2\frac{\delta_j}{b_\star h}
 \operatorname{Re}((D_hu)_j\overline{u_j}).
\]
Adding $-h^{-2}\alpha_j|u_j|^2/b_\star$ to the first term gives $h^{-2}\zeta_j|u_j|^2$. The second term has coefficient
$\eta_j=\delta_j/(b_\star h)$.

Use $\zeta_j,\eta_j$ from \eqref{A:eq:12.15l} and the first-order operator $\Gamma_\eta^-$ of Subsection~\ref{A:sec:dumitriu-edelman-decomposition}. Its quadratic form is
\[
\langle u,\Gamma_{\eta}^{-}u\rangle = -2\operatorname{Re} \sum_{j=0}^{M} \eta_j(D_hu)_j\overline{u_j} = -\frac{2h^{-2}}{b_\star} \operatorname{Re} \sum_{j=0}^{M} \delta_j(u_{j+1}-u_j)\overline{u_j}.
\eqtag{13.26}
\]
The sum of the diagonal and first-order terms is \(h^{-2}\sum_{j=1}^{M}\zeta_j|u_j|^2 + \langle u,\Gamma_{\eta}^{-}u\rangle\). Using \eqref{A:eq:12.15l} and \eqref{A:eq:13.26},
\[
\begin{aligned}
h^{-2}\sum_{j=1}^{M}\zeta_j|u_j|^2+\langle u,\Gamma_{\eta}^{-}u\rangle
&=-\frac{h^{-2}}{b_\star}\sum_{j=1}^{M}(\alpha_j+2\delta_j)|u_j|^2
-\frac{2h^{-2}}{b_\star}\operatorname{Re}\sum_{j=0}^{M}
\delta_j(u_{j+1}-u_j)\overline{u_j}
\\
&=-\frac{h^{-2}}{b_\star}\sum_{j=1}^{M}\alpha_j|u_j|^2
-\frac{2h^{-2}}{b_\star}\operatorname{Re}\sum_{j=0}^{M}\delta_j u_{j+1}\overline{u_j}.
\end{aligned}
\eqtag{13.27}
\]
On the other hand, expanding the first term of \eqref{A:eq:13.20},
\[
\sum_{j=0}^{M}
\frac{\delta_j}{b_\star}|D_hu_j|^2
=
\frac{h^{-2}}{b_\star}
\sum_{j=0}^{M}
\delta_j|u_{j+1}-u_j|^2,
\]
and using \(\delta_0=\delta_M=0\), we get
\[
\frac{h^{-2}}{b_\star} \sum_{j=0}^{M} \delta_j|u_{j+1}-u_j|^2 ={} \frac{h^{-2}}{b_\star} \sum_{j=1}^{M} (\delta_{j-1}+\delta_j)|u_j|^2 - \frac{2h^{-2}}{b_\star} \operatorname{Re} \sum_{j=0}^{M} \delta_j u_{j+1}\overline{u_j}.
\eqtag{13.28}
\]
Substituting \eqref{A:eq:13.28} into \eqref{A:eq:13.20} cancels the diagonal \(\delta_{j-1}+\delta_j\) terms and gives the right hand side of \eqref{A:eq:13.27}. This proves the
following identity.

\begin{proposition}\label{A:res:13.1}
For every $n$, the reduced matrix has the exact decomposition
\[
 \widetilde H_n^{\mathrm{red}}
 =H_n^\circ+h^{-2}M_{\zeta_n}+\Gamma_{\eta_n}^{-}.
 \eqtag{12.15m}
\]
The boundary contributions in the two forms agree when their auxiliary coefficients agree.
\end{proposition}

\par\smallskip Define the operator retaining the diagonal coefficient fluctuations after removing the first-order term
\(H_n^{\mathrm{diag}} := H_n^\circ+h^{-2}M_{\zeta_n}.\) Then
\[
\widetilde H_n^{\mathrm{red}}
=
H_n^{\mathrm{diag}}+\Gamma_{\eta_n}^{-}.
\eqtag{13.32}
\]

Let
\[
\mathcal R_n
=
\left\{j:\ s^b_{n,\kappa(j)+R_V}\le s_0\right\}
\]
be the regular coefficient region corresponding to Lemma~\ref{A:res:mean}.

\begin{lemma}\label{A:res:13.2}
For every \(1\le p<\infty\), under the exact conditional law with \(\tau_{\mathrm{bd}}\in\mathcal T_{n,I_n}(M_0)\),
\[
\sup_{j\in\mathcal R_n}
\left(
\|A_j-A_j^\circ\|_{L^p}
+\|B_j-B_j^\circ\|_{L^p}
\right)
\le C_pn^{-1/2},
\]
and
\[
\sup_{j\in\mathcal R_n}
\left(
|\mathbb E(A_j-A_j^\circ)|
+|\mathbb E(B_j-B_j^\circ)|
\right)
\le C n^{-1}.
\]
\end{lemma}

For all retained coordinates,
\[
\|A_j-A_j^\circ\|_{L^p}
+\|B_j-B_j^\circ\|_{L^p}
\le C_p\sqrt{\frac{\log n}{n}},
\]
and almost surely, for all sufficiently large \(n\),
\[
\max_j\bigl(|A_j-A_j^\circ|+|B_j-B_j^\circ|\bigr)
\le C n^{-1/2}\sqrt{\log n}.
\]

\begin{proof}
For $1\le j\le M$, the reindexing identifies $A_j$ with the original coefficient $a_{r+j}$; for $1\le j<M$, it identifies $B_j$ with $b_{r+j}$. Corollary~\ref{A:res:conditioning}
identifies their joint conditional law and minimizer with those of Lemma~\ref{A:res:mean}. For $j\in\mathcal R_n$, all coefficients in the required interaction neighborhood have
parameter at most $s_0$. Consequently,
\[
 \|A_j-A_j^\circ\|_p+\|B_j-B_j^\circ\|_p\le C_pn^{-1/2},
 \qquad
 |\mathbb E_\tau(A_j-A_j^\circ)|
 +|\mathbb E_\tau(B_j-B_j^\circ)|\le Cn^{-1}.
\]
The rough estimate in that lemma gives $C_p\sqrt{\log(n+1)/n}$ at every remaining free coordinate. The deterministic auxiliary boundary points have zero fluctuation.

To obtain the almost-sure bound on the maximum, fix $Q>2$. The maximum estimate \eqref{A:eq:mode-maximum}, applied under the exact conditional law, gives
\[
 \mathbb P_\tau\!\left(
 \max_j(|A_j-A_j^\circ|+|B_j-B_j^\circ|)
       >C_Q\sqrt{\log(n+1)/n}\right)\le C_Qn^{-Q}
\]
on the bounded set of exterior values, after enlarging $C_Q$ to combine the two coordinate types. Integrating over that set and using \eqref{A:eq:integrated-tail} bounds the
unconditional probability by $C_Qn^{-Q}+Ce^{-cn}$. Its sum over $n$ is finite. Borel--Cantelli therefore proves the stated maximum bound on the original joint coupling. The
equivalence $\log(n+1)\asymp\log n$ for $n\ge2$ gives the displayed normalization.
\end{proof}

Define the weighted squared fluctuation \(\mathcal E_{\eta,n} := \sum_{j=0}^{M} h(1+x_j)^{-3/2}|\eta_j|^2\). Since
\(\eta_j=\frac{\delta_j}{b_\star h},\) we have
\[
\mathcal E_{\eta,n}
=
\frac{1}{b_\star^2h}
\sum_{j=0}^{M}
(1+x_j)^{-3/2}|\delta_j|^2.
\]
Because \(\delta_0=\delta_M=0\), the boundary terms vanish.

\begin{proposition}\label{A:res:13.3}
Almost surely, for all sufficiently large \(n\),
\[
\mathcal E_{\eta,n}
\le
C_{\omega,V,\beta}h.
\eqtag{13.48}
\]
More precisely, uniformly for \(\tau_{\mathrm{bd}}\in\mathcal T_{n,I_n}(M_0)\),
\[
\mathbb E_{\tau_{\mathrm{bd}}}\mathcal E_{\eta,n}
\le
C_{V,\beta}h,
\eqtag{13.49}
\]
and for \(A\ge A_0\),
\[
\mathbb P_{\tau_{\mathrm{bd}}}
\left(
\mathcal E_{\eta,n}>Ah
\right)
\le
2\exp(-cA h^{-1}).
\eqtag{13.50}
\]
\end{proposition}

\begin{proof}
First compute the conditional expectation. On \(\mathcal R_n\), Lemma~\ref{A:res:13.2} gives \(\mathbb E_{\tau_{\mathrm{bd}}}|\delta_j|^2\le C n^{-1}.\) Since
\(\sum_j(1+x_j)^{-3/2}\le Ch^{-1}\), the regular contribution is at most
\[
\frac{C}{h n}
\sum_{j\in\mathcal R_n}(1+x_j)^{-3/2}
\le C h.
\]
On the complement of \(\mathcal R_n\), the rough part of Lemma~\ref{A:res:13.2} gives \(\mathbb E|\delta_j|^2\le C(\log n)/n\). Moreover \(x_j\ge c h^{-2}\) there, and hence
\((1+x_j)^{-3/2}\le Ch^3\). Since there are at most \(Ch^{-3}\) such sites, their contribution is bounded by
\[
\frac{C}{h}\,
\frac{\log n}{n}
\sum_{j\notin\mathcal R_n}(1+x_j)^{-3/2}
\le C h^2\log n=o(h).
\]
Estimate \eqref{A:eq:13.49} follows. Set \(F(Y_n)=\mathcal E_{\eta,n}^{1/2}\). This is a Lipschitz function of the coefficient vector \(Y_n\). Indeed,
\[
\mathcal E_{\eta,n}
=
\sum_j
\frac{(1+x_j)^{-3/2}}{b_\star^2h}
(B_j-B_j^\circ)^2.
\]
Hence the largest coefficient in the quadratic form defining \(\mathcal E_{\eta,n}\) is at most \(C h^{-1}\). Therefore
\(\operatorname{Lip}(F)\le C h^{-1/2}.\) The conditional measure satisfies a logarithmic Sobolev inequality with
constant \(C/n\), because \(\nabla^2\Phi_n^{\tau_{\mathrm{bd}}}\ge cnI.\) Thus
\[
\mathbb P_{\tau_{\mathrm{bd}}}
\left(
F-\mathbb E_{\tau_{\mathrm{bd}}}F>t
\right)
\le
\exp\left(
-\frac{c n t^2}{\operatorname{Lip}(F)^2}
\right)
\le
\exp(-c n h\,t^2).
\]
By Jensen,
\[
\mathbb E_{\tau_{\mathrm{bd}}}F
\le
\left(
\mathbb E_{\tau_{\mathrm{bd}}}\mathcal E_{\eta,n}
\right)^{1/2}
\le
C_{V,\beta}h^{1/2}.
\]
Choose $A_0\ge4C_{V,\beta}^2$. For $A\ge A_0$, the event $F^2>Ah$ implies $F-\mathbb EF>\sqrt{Ah}/2$. With $t=\sqrt{Ah}/2$, $nh t^2=(A/4)nh^2\ge cAh^{-1}$. Therefore
\(\mathbb P_{\tau_{\mathrm{bd}}} \left( \mathcal E_{\eta,n}>Ah \right) \le 2\mathrm{e}^{-cAh^{-1}}\). Estimate \eqref{A:eq:13.50} follows. By \eqref{A:eq:integrated-tail}, the
unconditional probability is bounded by \(2\mathrm{e}^{-cAh^{-1}}+C\mathrm{e}^{-cn}\). Since \(h^{-1}\asymp n^{1/3}\), this is summable for any sufficiently large fixed \(A\);
Borel--Cantelli proves \eqref{A:eq:13.48}. For every fixed $p\ge1$, tail integration gives
\[
 \mathbb E_\tau\mathcal E_{\eta,n}^p
 \le(A_0h)^p+2ph^p\int_{A_0}^\infty A^{p-1}\mathrm e^{-cA/h}\,\mathrm dA
 \le C_ph^p.
\]
For arbitrary exterior data, writing $Z_n=\max_i|Y_{n,i}|$ and using \eqref{A:eq:mode-growth} gives
\[
 \mathcal E_{\eta,n}\le Ch^{-2}(1+Z_n)^{2D_V+2}.
\]
Thus $\|\mathcal E_{\eta,n}\|_{2p}\le C_ph^{-2}$ by Lemma~\ref{A:res:raw}. Equation~\eqref{A:eq:integrated-moment} bounds the contribution of unbounded exterior data by
$C_ph^{-2}\mathrm e^{-c n/(2p)}$. It follows that the unconditional $L^p$ norm is also at most $C_ph$.
\end{proof}

The coefficient estimates above will give the first-order resolvent comparison after the common form bounds have been established in Theorem~\ref{A:res:16.2}. We first construct
the scalar Gaussian approximation, keeping the exact conditional law.

To express $\zeta_n$ as a linear image, set \(\Delta Y_n = (\alpha_1,\ldots,\alpha_M,\delta_1,\ldots,\delta_{M-1}).\) Define
\(\mathcal L_{\mathrm{diag}}:\mathbb R^{2M-1}\to\mathbb R^M\) by
\[
(\mathcal L_{\mathrm{diag}}\Delta Y_n)_j
=
-\frac{1}{b_\star}
\left(
\alpha_j+2\delta_j
\right),
\quad
1\le j\le M,
\]
with the convention \(\delta_M=0\). Then \(\zeta_n=\mathcal L_{\mathrm{diag}}\Delta Y_n.\) The first-order coefficient vector is
\[
\eta_n=\mathcal L_{\mathrm{g}}\Delta Y_n,
\quad
(\mathcal L_{\mathrm{g}}\Delta Y_n)_j
=
\frac{\delta_j}{b_\star h}.
\]
Let $K_n=\nabla^2\Phi_n^{\tau_{\mathrm{bd}}}(Y_n^\circ)$ and, conditionally on the exterior, let $G_n\sim N(0,K_n^{-1})$. For its image under the linear coefficient map
$\zeta_n^G=\mathcal L_{\mathrm{diag}}G_n$,
\[
 \operatorname{Cov}(\zeta_n^G\mid\tau_{\mathrm{bd}})
 =\mathcal L_{\mathrm{diag}}K_n^{-1}\mathcal L_{\mathrm{diag}}^*.
\]
In contrast, the exact scalar coefficients satisfy
\[
 \operatorname{Cov}(\zeta_n\mid\tau_{\mathrm{bd}})
 =\mathcal L_{\mathrm{diag}}
       \operatorname{Cov}(Y_n\mid\tau_{\mathrm{bd}})
          \mathcal L_{\mathrm{diag}}^*.
\]
Section~\ref{A:sec:gaussian-approximation-of-the-conditional-gibbs-measure} couples these two coefficient laws.
Section~\ref{A:sec:application-to-the-conditional-jacobi-coefficients} estimates the Gaussian covariance displayed above.

On a regular coefficient interval, Lemma~\ref{A:res:13.2} gives
\[
 \mathbb E_{\tau_{\mathrm{bd}}}|\zeta_j|^2
 \le\frac{2}{b_\star^2}
       \bigl(\mathbb E|\alpha_j|^2+4\mathbb E|\delta_j|^2\bigr)
 \le Cn^{-1}\le Ch^3.
\]
Consequently $\|h^{-2}\zeta_j\|_2\le Ch^{-1/2}$. The Brownian variance is computed in Lemma~\ref{A:res:15.1}. On the complementary region the corresponding bound is
$Ch^{-1/2}\sqrt{\log(n+1)}$, and the diagonal resolvent estimate uses the weight $(1+x_j)^{-7/2}$.

The conditional minimizer is uniformly close to the local minimizer only at sites separated from both conditioning boundaries.

The gradient error at the new left boundary is bounded by $C_\omega n^{-1/2}\sqrt{\log(n+1)}$, whereas the ordinary error in the stationarity equations is \(O(n^{-1})\). Equation
\eqref{A:eq:13.17} keeps these two scales separate.

Choose a fixed number \(0<s_{\mathrm{cp}}<1\) and set \(J_n^{\mathrm{cp}} = \max\{j\le M:s^b_{n,\kappa(j)+R_V}\le s_{\mathrm{cp}}\}\). Recall the length \(J_n^{\mathrm L}\) of the
initial interval in \eqref{A:eq:13.17b}, and define the interval \(I_n^{\mathrm{cp,mid}} = \{J_n^{\mathrm L}+1,\ldots,J_n^{\mathrm{cp}}\}.\) For each index in this block, the
coefficients at distance at most $R_V$ lie strictly between the two conditioned boundaries. Since
\(s_{\mathrm{cp}}>0\) is fixed, \(J_n^{\mathrm{cp}}\asymp n, \allowbreak\, x_{J_n^{\mathrm{cp}}}\asymp h^{-2}.\)
Choosing \(A_{\mathrm L}\) sufficiently large, equation \eqref{A:eq:13.17} gives
\[
\max_{J_n^{\mathrm L}+1\le j\le J_n^{\mathrm{cp}}+1}
\left(
|A_j^\circ-a_j^\dagger|
+
|B_j^\circ-b_j^\dagger|
\right)
\le
C_\omega n^{-1}
=
C_\omega h^3.
\eqtag{13.67}
\]
On the left boundary strip itself we use the boundary estimate
\[
\max_{1\le j\le J_n^{\mathrm L}+R_V}
\left(
|A_j^\circ-a_j^\dagger|
+
|B_j^\circ-b_j^\dagger|
\right)
\le
C_\omega n^{-1/2}\sqrt{\log n}.
\eqtag{13.67a}
\]

The bounds \eqref{A:eq:13.67}--\eqref{A:eq:13.67a} are used in Subsection~\ref{A:sec:comparison-of-the-deterministic-coefficients}. There the corresponding matrix difference is
split into entries meeting the left strip and the remaining regular entries. The left strip entries are $O_\omega(h^{-1/2}\sqrt{\log(n+1)})$ and are controlled by the estimate for
vectors vanishing at the left boundary. The remaining entries are $O_\omega(h)$ with fixed bandwidth. Equations~\eqref{A:eq:17.23}--\eqref{A:eq:17.27} give squared resolvent error
$o(h)$ on the regular interval, also with the common Brownian potential added. The terminal region is treated separately by its resolvent tail. The coefficient estimates here
follow directly from the conditional distributions and the potential bounds established above.

We henceforth suppress the reindexing unitaries and compare the diagonal coefficients in
\(H_n^{\mathrm{diag}}=H_n^\circ+h^{-2}M_{\zeta_n}\).

\subsection{The general Jacobi model: Gaussian approximation and comparison}
\label{sec:combined-a-gaussian-comparison}
\label{A:sec:the-general-jacobi-model-gaussian-approximation-and-comparison}
We first couple the conditional coefficient vector to a Gaussian centered at the minimizer of its conditional energy, with covariance equal to the inverse Hessian there. The
required error is $O(n^{-1})$ in each coordinate, uniformly in the length of the interval. The coordinates of the resulting Gaussian vector are correlated; its comparison with
Brownian increments is a separate argument.

Throughout this comparison for the general Jacobi model, $1\le p<\infty$. Conditional estimates have deterministic constants, uniform over the specified bounded exterior data.
Almost sure statements refer to the joint construction of those data, the coefficients, and the Brownian motion at all matrix sizes.

\subsubsection{Transport from a Gaussian measure}\label{A:sec:localized-gaussian-transport}
The Langevin argument requires a potential defined on the whole space with bounded Hessian. Starting from the strongly convex conditional potential, we continue the logarithmic
terms through the positivity boundary and modify the polynomial growth at infinity. Both changes preserve the potential and its Hessian near the conditional minimizer, and hence
the Gaussian covariance. A probability kernel then recovers the original law, including its mass outside this neighborhood.

We estimate each coordinate of the transport velocity through an exponentially weighted sum of derivatives of the nonquadratic part of the potential. Summability of the weights
makes the error bound independent of the number of coefficients.

\par\medskip\noindent Let \(\mathfrak G_d\) be a graph on \(\{1,\ldots,d\}\) of uniformly bounded degree, and let
\(\mathrm d(\cdot,\cdot)\) be its graph distance. Suppose that \(\Phi^\sharp\in C^\infty(\mathbb R^d)\) has minimizer
\(0\), satisfies
\[
\nabla^2\Phi^\sharp(u)\ge c n I,
\eqtag{A.35a}
\]
and is finite-range in the sense that, for one fixed \(R_0<\infty\),
\(\partial_{i_1}\cdots\partial_{i_m}\Phi^\sharp(u)=0\) whenever \(\max_{r,s}\mathrm d(i_r,i_s)>R_0\). Assume also that,
on every fixed coordinate cube \(Q_R=\{u:\|u\|_{\ell^\infty}\le R\},\) the local derivatives obey
\(\sup_{u\in Q_R} |\partial_{i_1}\cdots\partial_{i_m}\Phi^\sharp(u)| \le C_{R,m}n\) for every fixed \(m\ge1\).

\begin{lemma}\label{A:res:A.11A}
For every fixed \(R>0\), there is a \(C^\infty\) finite-range potential \(\Phi^{\mathrm{ext}}\) such that
\[
\Phi^{\mathrm{ext}}(u)=\Phi^\sharp(u)
\quad (u\in Q_R),
\eqtag{A.35d}
\]
\[
 c_1nI\le \nabla^2\Phi^{\mathrm{ext}}(u)\le C_1nI
\quad (u\in\mathbb R^d),
\eqtag{A.35e}
\]
and, for every fixed \(m\ge2\),
\[
|\partial_{i_1}\cdots\partial_{i_m}
\Phi^{\mathrm{ext}}(u)|
\le C_m n.
\eqtag{A.35f}
\]
\end{lemma}

The constants depend only on the constants above, on \(R\), and on the graph and distance bound specified above. In particular,
\[
\nabla\Phi^{\mathrm{ext}}(0)=0,
\quad
\nabla^2\Phi^{\mathrm{ext}}(0)=\nabla^2\Phi^\sharp(0).
\eqtag{A.35g}
\]

\begin{proof}
Composing the potential with a bounded smooth map in each coordinate bounds its derivatives and introduces the Hessian term $\chi_R''\partial_k\Phi^\sharp$, whose sign is
uncontrolled. We add one diagonal term to dominate it and another to restore the curvature where the coordinate map becomes constant.

Choose an odd function \(\chi_R\in C^\infty(\mathbb R)\) satisfying \(\chi_R(x)=x\allowbreak\; (|x|\le R), \allowbreak\, |\chi_R(x)|\le R+2,\)
\(0\le \chi_R'(x)\le1, \allowbreak\, \chi_R'(x)=0\allowbreak\; (|x|\ge R+2).\) For example, take a smooth decreasing function $\eta$ with $\eta=1$ on $(-\infty,0]$, $\eta=0$ on
$[2,\infty)$, and $0\le\eta\le1$, and set $\chi_R(x)=\int_0^x\eta(|t|-R)\,\mathrm dt$. The derivatives of $\eta$ vanish at the transition boundary points. Thus $\chi_R$ is smooth,
and all its derivatives are bounded. Put \(q_R''(x)=1-\chi_R'(x)^2, \allowbreak\, q_R(0)=q_R'(0)=0.\) Since
\(\chi_R'=1\) on \([-R,R]\), the smooth convex function \(q_R\) and its first derivative vanish there. Moreover,
\[
\chi_R'(x)^2+q_R''(x)=1
\eqtag{A.35k}
\]
and \(q_R''=1\) outside a compact interval, so \(q_R\) has quadratic growth. Choose another smooth convex function
\(\rho_R\), also vanishing together with its first derivative on \([-R,R]\), such that
\[
\rho_R''(x)\ge |\chi_R''(x)|.
\eqtag{A.35l}
\]
For the preceding choice, $\chi_R''$ has a fixed sign on each transition interval and vanishes to all orders at its boundary points. Hence $|\chi_R''|$ is smooth. We may take
$\rho_R''=|\chi_R''|$ and $\rho_R(0)=\rho_R'(0)=0$.

Let \(\chi_R(u)=(\chi_R(u_1),\ldots,\chi_R(u_d))\), and define
\[
\Phi^{\mathrm{ext}}(u)
=
\Phi^\sharp(\chi_R(u))
+
\lambda n\sum_{k=1}^d
\bigl(\rho_R(u_k)+q_R(u_k)\bigr),
\]
where \(\lambda\) will be fixed below. The construction is coordinatewise, so mixed derivatives at sufficiently separated indices remain zero. It also gives \eqref{A:eq:A.35d}.

Write \(D_\chi(u)=\operatorname{diag} \bigl(\chi_R'(u_1),\ldots,\chi_R'(u_d)\bigr).\) A direct differentiation gives
\[
\begin{aligned}
\nabla^2\Phi^{\mathrm{ext}}(u)
={}&D_\chi(u)\nabla^2\Phi^\sharp(\chi_R(u))D_\chi(u)
+\operatorname{diag}\left(\partial_k\Phi^\sharp(\chi_R(u))\chi_R''(u_k)\right)_{k=1}^d
\\
&+\lambda n\operatorname{diag}\left(\rho_R''(u_k)+q_R''(u_k)\right)_{k=1}^d.
\end{aligned}
\eqtag{A.35n}
\]
Because every coordinate of \(\chi_R(u)\) lies in \([-R-2,R+2]\), the local derivative assumption gives
\[
|\partial_k\Phi^\sharp(\chi_R(u))|
+
\|\nabla^2\Phi^\sharp(\chi_R(u))\|_{\mathrm{op}}
\le C_R n.
\eqtag{A.35o}
\]
For $v\in\mathbb R^d$, the three terms in \eqref{A:eq:A.35n} give
\[
\langle v,\nabla^2\Phi^{\mathrm{ext}}(u)v\rangle \ge{} cn\sum_k\chi_R'(u_k)^2|v_k|^2
 -C_Rn\sum_k|\chi_R''(u_k)|\,|v_k|^2 +\lambda n\sum_k\bigl(\rho_R''(u_k)+q_R''(u_k)\bigr)|v_k|^2.
\]
Choose $\lambda\ge\max\{2C_R,c\}$. By \eqref{A:eq:A.35k} and \eqref{A:eq:A.35l}, the right side is at least $cn\sum_k(\chi_R'(u_k)^2+q_R''(u_k))|v_k|^2=cn\|v\|^2$. For the upper
bound, each row of $\nabla^2\Phi^\sharp$ has at most $N_{R_0}$ nonzero entries, where $N_{R_0}$ depends only on the fixed range and the degree bound. Its row and column sums are at
most $N_{R_0}C_{R+2,2}n$. The Schur test bounds the first term of \eqref{A:eq:A.35n} by $Cn\|v\|^2$; boundedness of $\chi_R''$, $\rho_R''$, and $q_R''$ gives the same bound for the
two diagonal terms. This proves \eqref{A:eq:A.35e}.

For indices $i_1,\ldots,i_m$, let $\Pi_m$ be the set of partitions of $\{1,\ldots,m\}$. In a partition $\pi$, a block $A$ contributes only when all $i_r$, $r\in A$, are equal;
denote that common index by $i_A$. Repeated differentiation gives
\[
 \partial_{i_1}\cdots\partial_{i_m}
       (\Phi^\sharp\circ\chi_R)(u)
 =\sum_{\substack{\pi\in\Pi_m:\\
              i_r=i_s\text{ for }r,s\in A,\ A\in\pi}}
   \left(\prod_{A\in\pi}\partial_{i_A}\right)
       \Phi^\sharp(\chi_R(u))
   \prod_{A\in\pi}\chi_R^{(|A|)}(u_{i_A}).
\]
Every derivative of $\Phi^\sharp$ in this sum is evaluated on $Q_{R+2}$ and has order at most $m$. Thus its absolute value is at most
\[
 n\,|\Pi_m|\,
   \max_{1\le j\le m}C_{R+2,j}
   \left(1+\max_{1\le j\le m}\|\chi_R^{(j)}\|_\infty\right)^m.
\]
For $m\ge2$ the derivative of each penalty is zero unless all $i_r$ agree; in that case it is $\lambda n(\rho_R^{(m)}+q_R^{(m)})(u_{i_1})$, which is bounded by $C_mn$. This proves
\eqref{A:eq:A.35f} independently of $d$. Equality on $Q_R$ gives \eqref{A:eq:A.35g}.
\end{proof}

The extension is centered at its minimizer, whereas concentration controls fluctuations about the mean. To pass between these centers without a dimension factor, we expand the
gradient at the mean. Its linear term then has zero expectation, and finite range bounds the expected remainder in each coordinate. The inverse Hessian converts that bound into an
$O(n^{-1})$ displacement.

Let \(\mu(\mathrm{d}y)=Z^{-1}\mathrm{e}^{-\Phi(y)}\,\mathrm{d}y\) be a probability measure on \(\mathbb R^d\), and let
\(y^\circ\) be the minimizer of \(\Phi\). Suppose that, uniformly in \(d\) and \(n\),
\[
cnI_d\le\nabla^2\Phi(y)\le CnI_d.
\]
Let $\mathrm d$ be the distance on a graph $\mathfrak G_d$. Assume $\partial_{ij}^2\Phi=0$ for $\mathrm d(i,j)>R_0$, with $R_0$ independent of $d,n$, and assume that the numbers of
vertices at each fixed distance are uniformly bounded:
\[
\sup_{d,i,r}\#\{j:\mathrm d(i,j)=r\}<\infty,
\]
and that
\[
|\partial_{ijk}^3\Phi(y)|\le Cn,
\]
with the derivative equal to zero unless \(i,j,k\) lie in one fixed graph neighborhood.

\begin{proposition}\label{A:res:A.11B}
If \(Y\sim\mu\), then, for every finite \(p\),
\[
\sup_i|\mathbb E(Y_i-y_i^\circ)|\le \frac{C}{n},
\quad
\sup_i\|Y_i-y_i^\circ\|_{L^p(\mu)}\le C_pn^{-1/2}.
\]
The constants are independent of the dimension.
\end{proposition}

\begin{proof}
Put \(U=Y-y^\circ\), \(m=\mathbb EU\), and \(X=U-m\). The Brascamp--Lieb inequality~\cite{brascamp1976extensions} gives
\[
\operatorname{Cov}(X)=\operatorname{Cov}(Y)\le (cn)^{-1}I_d,
\quad
\|X_i\|_{L^p}\le C_pn^{-1/2},
\]
where the second estimate follows from the Bakry--\'Emery logarithmic Sobolev inequality~\cite{bakry1985diffusions}. These estimates concern fluctuations about the mean and do not
involve \(m\).

The lower Hessian bound gives $\Phi(y)\ge\Phi(y^\circ)+cn|y-y^\circ|^2/2$; the upper bound gives $|\nabla\Phi(y)|\le Cn|y-y^\circ|$. Choose a smooth cutoff
$\chi_R(y)=\chi((y-y^\circ)/R)$, equal to one on the ball of radius $R$, supported on the ball of radius $2R$, with $|\partial_i\chi_R|\le C/R$. Integrating a compactly supported
derivative gives
\[
 0=\int\partial_i(\chi_Re^{-\Phi})\,\mathrm dy,
 \qquad \mathbb E[\chi_R\partial_i\Phi]
       =\mathbb E[\partial_i\chi_R],
 \qquad |\mathbb E[\partial_i\chi_R]|\le C/R.
\]
The quadratic density bound makes $|\nabla\Phi|$ integrable. Dominated convergence as $R\to\infty$ therefore yields
\[
\mathbb E\,\partial_i\Phi(y^\circ+m+X)=0.
\]
Taylor's formula at the mean \(y^\circ+m\) gives
\[
\partial_i\Phi(y^\circ+m+X)
=
\partial_i\Phi(y^\circ+m)
+
\sum_j\partial_{ij}^2\Phi(y^\circ+m)X_j
+
\mathcal R_i(X),
\]
where
\[
\mathcal R_i(X)
=
\sum_{j,k}X_jX_k
\int_0^1(1-t)\,
\partial_{ijk}^3\Phi(y^\circ+m+tX)\,\mathrm{d}t.
\]
The linear term has zero expectation because \(\mathbb EX=0\). For each fixed \(i\), finite range leaves only a bounded number of pairs \((j,k)\), and Brascamp--Lieb gives
\[
\mathbb E|X_jX_k|
\le
\operatorname{Var}(X_j)^{1/2}
\operatorname{Var}(X_k)^{1/2}
\le \frac{C}{n}.
\]
\[
\sup_i|\mathbb E\mathcal R_i(X)|\le C.
\]
Taking expectations in the Taylor formula yields
\[
\sup_i|\partial_i\Phi(y^\circ+m)|\le C.
\]

Since \(\nabla\Phi(y^\circ)=0\), the mean-value theorem gives
\[
\nabla\Phi(y^\circ+m)=\overline K m,
\quad
\overline K=\int_0^1\nabla^2\Phi(y^\circ+tm)\,\mathrm{d}t.
\]
The matrix $\overline K$ also vanishes when $\mathrm d(i,j)>R_0$ and satisfies \(cnI\le\overline K\le CnI\). The bound on the number of nonzero entries in each row and the upper
operator bound imply \(\sup_i\sum_j|\overline K_{ij}|\le C'n\). Lemma~\ref{A:res:A.9}, with zero additional diagonal part, therefore gives
\[
|(\overline K^{-1})_{ij}|\le \frac{C}{n}\mathrm{e}^{-c\mathrm d(i,j)}.
\]
Summing over distances, using the uniform bound on their multiplicities, gives \(\sup_i\sum_j|(\overline K^{-1})_{ij}|\le C/n\). Hence
\[
\sup_i|m_i|
\le
\sup_i\sum_j|(\overline K^{-1})_{ij}|
\sup_j|\partial_j\Phi(y^\circ+m)|
\le \frac{C}{n}.
\]
Finally, \(U_i=X_i+m_i\), so the centered \(L^p\)-estimate and the bound on \(m_i\) give \(\|U_i\|_{L^p}\le C_pn^{-1/2}\), completing both claims.
\end{proof}

For the comparison of the two densities we need an event containing a fixed coordinate neighborhood of the minimizer. The next consequence of concentration keeps the mean
displacement explicit, so that it applies to both the original law and its extension.

Let \(\mathcal D\subset\mathbb R^d\) be convex and let
\[
\mu(\mathrm{d}y)=Z^{-1}\mathrm{e}^{-\Phi(y)}\mathbf 1_{\mathcal D}(y)\,\mathrm{d}y,
\]
where \(\Phi\in C^2(\mathcal D)\), \(\nabla^2\Phi\ge \rho n I_d\), and \(y^\circ\) is the minimizer of \(\Phi\) on
\(\mathcal D\). Let $Y\sim\mu$ and put $b_n=\|\mathbb EY-y^\circ\|_{\ell^\infty}$.

\begin{lemma}\label{A:res:A.11C}
For every coordinate \(i\) and every \(t>0\),
\[
\mathbb P\bigl(|Y_i-\mathbb EY_i|\ge t\bigr)
\le 2\exp\{-\rho nt^2/2\}.
\]
\end{lemma}

Consequently, for \(r>b_n\),
\[
\mathbb P(Y_i-y_i^\circ\le-r)
\le \exp\{-\rho n(r-b_n)^2/2\},
\quad
\mathbb P(|Y_i-y_i^\circ|\ge r)
\le2\exp\{-\rho n(r-b_n)^2/2\}.
\]
Let \(J\subset\{1,\ldots,d\}\), let \(r_j\ge r_*>0\) for \(j\in J\), and let \(R>0\). If
\[
\mathcal C=
\{y:y_j-y_j^\circ\ge-r_j\ (j\in J),
\ \|y-y^\circ\|_{\ell^\infty}\le R\}
\]
and \(b_n\le\frac12\min\{r_*,R\}\), then
\[
\mu(\mathcal C^c)
\le
3d\exp\left\{-\frac{\rho n}{8}\min\{r_*^2,R^2\}\right\}.
\]

\begin{proof}

For $Z_i=Y_i-\mathbb EY_i$, the concentration estimate \eqref{A:eq:convex-concentration}, with curvature $\rho n$ and the $1$-Lipschitz coordinate map, gives
\[
 \mathbb P(\pm Z_i\ge t)\le e^{-\rho nt^2/2},
 \qquad
 \mathbb P(|Z_i|\ge t)\le2e^{-\rho nt^2/2}.
\]
Since $|\mathbb EY_i-y_i^\circ|\le b_n$, for $r>b_n$,
\[
 \{Y_i-y_i^\circ\le-r\}\subset\{Z_i\le-(r-b_n)\},
 \qquad
 \{|Y_i-y_i^\circ|\ge r\}\subset\{|Z_i|\ge r-b_n\}.
\]
These inclusions prove the two bounds centered at the minimizer. Moreover,
\[
 \begin{aligned}
 \mu(\mathcal C^c)
 &\le\sum_{j\in J}e^{-\rho n(r_j-b_n)^2/2}
       +2\sum_{i=1}^d e^{-\rho n(R-b_n)^2/2}\\
 &\le |J|e^{-\rho nr_*^2/8}+2d e^{-\rho nR^2/8}
 \le3d e^{-\rho n\min\{r_*^2,R^2\}/8}.
 \end{aligned}
\]
We used $r_j-b_n\ge r_*/2$, $R-b_n\ge R/2$, and $|J|\le d$.

\end{proof}

When two potentials agree on a set, normalizing their restrictions to that set gives the same probability law, even if their full partition functions differ. We couple a common
submeasure on the diagonal and couple the two nonnegative measures obtained by subtracting that submeasure from the original laws. Both marginals are therefore retained. The
conditional construction recovers the original coefficients while keeping the Gaussian variable fixed.

\par\Needspace{7\baselineskip}
\begin{lemma}\label{A:res:A.11D}
Let \(\mu\) and \(\nu\) be probability measures on a standard Borel space and let \(\mathcal C\) be a Borel set. Assume
\(\mu(\mathcal C),\nu(\mathcal C)>0\) and
\[
\mu(\,\cdot\mid\mathcal C)=\nu(\,\cdot\mid\mathcal C).
\]
Then there is a coupling \((X,Y)\) of \((\mu,\nu)\) such that
\[
\mathbb P(X\ne Y)
\le
\max\{\mu(\mathcal C^c),\nu(\mathcal C^c)\}.
\]
If the state space is contained in \(\mathbb R^d\), then for every \(1\le p<q<\infty\),
\[
\|X_i-Y_i\|_{L^p}
\le
\bigl(\|X_i\|_{L^q}+\|Y_i\|_{L^q}\bigr)
\mathbb P(X\ne Y)^{1/p-1/q}.
\]
\end{lemma}

\begin{proof}
Set \(a=\mu(\mathcal C)\), \(b=\nu(\mathcal C)\), \(m=\min\{a,b\}\), and \(\rho=\mu(\,\cdot\mid\mathcal C)=\nu(\,\cdot\mid\mathcal C)\). The measure \(m\rho\) is dominated by both
\(\mu\) and \(\nu\). For $m<1$, put $\mu_0=(\mu-m\rho)/(1-m)$ and define the coupling measure
\[
\pi
=
m\,\Delta_\#\rho
+
\frac{(\mu-m\rho)\otimes(\nu-m\rho)}{1-m},
\]
where \(\Delta(x)=(x,x)\). Its first marginal is \(\mu\), its second marginal is \(\nu\), and the diagonal component has mass \(m\). Thus
\[
\pi\{X\ne Y\}
\le1-m
=\max\{1-a,1-b\}.
\]
If $m=1$, then $a=b=1$ and $\mu=\nu=\rho$; $\pi=\Delta_\#\rho$ has the required marginals and $\pi\{X\ne Y\}=0$. The construction can also be made while keeping the second random
variable unchanged. For $m<1$ and a given $Y\sim\nu$, put
\[
r(y)=\begin{cases}(m/b)\mathbf 1_{\mathcal C}(y),&b>0,\\0,&b=0.\end{cases}
\]
Given $Y=y$, set $X=y$ with probability $r(y)$; otherwise sample $X$ from $\mu_0$. Since $r(y)\nu(\mathrm{d}y)=m\rho(\mathrm{d}y)$ and
$(1-r(y))\nu(\mathrm{d}y)=\nu(\mathrm{d}y)-m\rho(\mathrm{d}y)$, the resulting joint law is $\pi$. This formula also shows that if $\mu$, $\nu$, and their common conditional law
depend measurably on an external parameter, then the conditional sampling construction is a measurable coupling kernel. For the moment estimate, put \(E=\{X\ne Y\}\). H\"older's
inequality gives
\[
\|(X_i-Y_i)\mathbf 1_E\|_{L^p}
\le
\|X_i-Y_i\|_{L^q}\mathbb P(E)^{1/p-1/q},
\]
and Minkowski's inequality bounds the first factor by \(\|X_i\|_{L^q}+\|Y_i\|_{L^q}\).
\end{proof}

Distances between coefficient coordinates refer to their matrix indices and types: there are only a fixed number of types, such as $a_j$ and $b_j$, at each index $j$. Consequently
the number of coordinates at distance $r$ from a given one is bounded independently of $r,n$. Hence
$\sup_i\sum_j\mathrm e^{-a\mathrm d(i,j)}\le C\sum_{r\ge0}\mathrm e^{-ar}<\infty$ for every $a>0$.

Let \(\Phi\) denote the globally uniformly convex extension on \(\mathbb R^d\), where \(d\asymp n\). We assume that its minimizer is the origin and that
\(cnI\le \nabla^2\Phi(u)\le CnI.\) Let \(K=\nabla^2\Phi(0), \allowbreak\, \gamma=N(0,K^{-1}),\) and write \(\Phi(u)=\frac12\langle u,Ku\rangle+\Psi(u).\) The coordinate graph is a
one-dimensional strip of fixed width. In particular, there is a deterministic constant \(C_{\mathcal G}\) such that
\(\sup_{n}\sup_{i}\sup_{r\ge0} \#\{j:\mathrm d(i,j)=r\} \le C_{\mathcal G}\), and all Hessians are finite-range:
\((\nabla^2\Phi(u))_{ij}=0\) if \(\mathrm d(i,j)>R_0\). For the interpolating potentials
\[
\Phi_t(u)
=
\frac12\langle u,Ku\rangle+t\Psi(u),
\quad 0\le t\le1,
\eqtag{A.37b}
\]
we therefore have, uniformly in \(t\),
\[
cnI\le\nabla^2\Phi_t(u)\le CnI,
\eqtag{A.37c}
\]
and $\partial_{ij}^2\Phi_t=0$ beyond the same fixed graph distance. We also assume the local bound on third derivatives furnished by Lemma~\ref{A:res:A.11A}:
\[
|\partial_{ijk}^3\Phi_t(u)|\le Cn,
\eqtag{A.37d}
\]
with the derivative equal to zero unless \(i,j,k\) lie in one fixed graph neighborhood.

The potential varies affinely with $t$; its exponential defines the interpolating density after normalization. Both endpoint potentials have gradient zero at the origin and the
same Hessian there. Let
\(\mu_t(\mathrm{d}u)=Z_t^{-1}\mathrm{e}^{-\Phi_t(u)}\,\mathrm{d}u.\) For every \(t\in[0,1]\), the potential \(\Phi_t\)
has minimizer zero and satisfies the hypotheses of Proposition~\ref{A:res:A.11B} with constants independent of the parameter. Consequently, for every finite \(p\),
\[
\sup_{t\in[0,1]}\sup_i
\left|\int u_i\,\mathrm{d}\mu_t(u)\right|
\le \frac{C}{n},
\quad
\sup_{t\in[0,1]}\sup_i
\|u_i\|_{L^p(\mu_t)}
\le C_pn^{-1/2}.
\eqtag{A.37e}
\]

\begin{lemma}\label{A:res:A.12}
For every \(p<\infty\),
\[
\sup_{t\in[0,1]}\sup_i
\|\partial_i\Psi\|_{L^p(\mu_t)}
\le C_p.
\eqtag{A.38}
\]
\end{lemma}

\begin{proof}
The choice $K=\nabla^2\Phi(0)$ gives $\nabla\Psi(0)=0$ and $\nabla^2\Psi(0)=0$. Thus the constant and linear terms in the expansion of $\partial_i\Psi$ vanish, and
\[
 \partial_i\Psi(u)
 =\sum_{j,k}\int_0^1(1-s)
       \partial_{ijk}^3\Phi(su)u_ju_k\,\mathrm ds.
\]
Only indices $j,k$ in a fixed neighborhood $N_i$ contribute. Its cardinality is bounded independently of $d$, so H\"older's inequality and \eqref{A:eq:A.37e} give
\[
 \|\partial_i\Psi\|_{L^p(\mu_t)}
 \le Cn\sum_{j,k\in N_i}
        \|u_j\|_{L^{2p}(\mu_t)}\|u_k\|_{L^{2p}(\mu_t)}
 \le C_p.
\]
The constants are uniform in $i,t,n,d$. No sum over all coordinates occurs in this estimate.
\end{proof}

To turn this derivative bound into a transport displacement, we solve a Poisson equation for the Langevin generator. The spectral gap alone would estimate all coordinates together.
We retain the decay with index distance of the derivative flow, as in the Helffer--Sj\"ostrand representation \cite{helffer1994correlation}, and integrate its decay in time. These
two operations give a summable weight and a factor $n^{-1}$, respectively.

\par\medskip\noindent Fix \(t\in[0,1]\), and let \(L_t=\Delta-\nabla\Phi_t\cdot\nabla\) be the reversible Langevin
generator in \(L^2(\mu_t)\). Let \(1\le p<\infty\),

and let \(f\in W^{1,p}(\mu_t)\cap L^2(\mu_t)\) satisfy \(\int f\,\mathrm{d}\mu_t=0.\) \par\Needspace{7\baselineskip}
\begin{lemma}\label{A:res:A.13}
The centered weak solution of \(L_t\varphi=f, \allowbreak\, \int\varphi\,\mathrm{d}\mu_t=0,\) satisfies, for some deterministic \(a,c_0>0\),
\[
\|\partial_i\varphi\|_{L^p(\mu_t)}
\le
\frac{1}{c_0n}
\sum_{j=1}^d
 \mathrm{e}^{-a\mathrm d(i,j)}
\|\partial_jf\|_{L^p(\mu_t)}.
\eqtag{A.39}
\]
\end{lemma}

The constants are uniform in $n,d,t$. The bound on the number of coordinates at each distance gives
\[
\sup_i\sum_j \mathrm{e}^{-a\mathrm d(i,j)}
\le
C_{\mathcal G}\sum_{r\ge0}\mathrm{e}^{-ar}
<\infty,
\]
and consequently
\[
\sup_i\|\partial_i\varphi\|_{L^p(\mu_t)}
\le
C_pn^{-1}
\sup_j\|\partial_jf\|_{L^p(\mu_t)}.
\eqtag{A.39a}
\]

\begin{proof}

Fix $t$. The parameter $s\ge0$ below is semigroup time and is distinct from the interpolation parameter $t\in[0,1]$. Consider the Langevin equation
\[
 \mathrm dX_s^x=-\nabla\Phi_t(X_s^x)\,\mathrm ds+\sqrt2\,\mathrm dW_s,
 \qquad X_0^x=x.
\]
By \eqref{A:eq:A.37c}, the drift is globally Lipschitz. The solution defines a differentiable stochastic flow with Jacobian
\[
 J_s^x=D_xX_s^x,\qquad
 \partial_sJ_s^x=-H_s^xJ_s^x,\qquad
 J_0^x=I,\qquad H_s^x=\nabla^2\Phi_t(X_s^x).
\]
The semigroup $P_s^tf(x)=\mathbb E[f(X_s^x)]$ is reversible with invariant measure $\mu_t$.

Fix a coordinate \(i\). Let \(W_i=\operatorname{diag} \bigl(\mathrm{e}^{a\mathrm d(i,1)},\ldots,\mathrm{e}^{a\mathrm d(i,d)}\bigr)\), where $0<a\le1$ will be chosen below. If
$\mathrm d(k,\ell)\le R_0$, then $|\mathrm d(i,k)-\mathrm d(i,\ell)|\le R_0$. Thus
\[
 \left|\mathrm e^{a(\mathrm d(i,k)-\mathrm d(i,\ell))}-1\right|
 \le C_{R_0}a.
\]
For $E_{i,s}^x=W_iH_s^xW_i^{-1}-H_s^x$, this gives
\[
 (E_{i,s}^x)_{k\ell}=(H_s^x)_{k\ell}
       \left[\mathrm e^{a(\mathrm d(i,k)-\mathrm d(i,\ell))}-1\right],
 \qquad \|E_{i,s}^x\|_{\mathrm{op}}\le C_0an.
\]
Indeed, each row and column has a bounded number of nonzero entries, and $|(H_s^x)_{k\ell}|\le Cn$ by \eqref{A:eq:A.37c}; the last bound follows from the Schur test.

The Jacobian equation retains the time ordering of the Hessians $H_s^x$. Set \(Z_s=W_iJ_s^xW_i^{-1}.\) By the Jacobian equation,
\(\frac {\mathrm{d}}{\mathrm{d}s}Z_s = -\bigl(H_s^x+E_{i,s}^x\bigr)Z_s.\) For every \(v\in\mathbb R^d\),
\[
\frac12\frac {\mathrm{d}}{\mathrm{d}s}\|Z_sv\|^2 = -\operatorname{Re} \langle Z_sv,(H_s^x+E_{i,s}^x)Z_sv\rangle \le -(c-C_0a)n\|Z_sv\|^2.
\]
Choose \(a\le c/(2C_0)\), and put \(c_0=c/2\). Gronwall's inequality gives
\[
\|W_iJ_s^xW_i^{-1}\|_{\mathrm{op}}
\le \mathrm{e}^{-c_0ns}.
\eqtag{A.40}
\]
Since \(W_ie_i=e_i\), taking the \((j,i)\)-entry in \eqref{A:eq:A.40} yields the following bound on the derivative of the flow
\[
|(J_s^x)_{ji}|
\le
 \mathrm{e}^{-c_0ns}\mathrm{e}^{-a\mathrm d(i,j)}.
\eqtag{A.40a}
\]
For smooth \(f\), differentiation under the expectation gives
\[
\partial_iP_s^tf(x)
=
\mathbb E
\left[
\sum_{j=1}^d
(J_s^x)_{ji}\,\partial_jf(X_s^x)
\right].
\]
Using \eqref{A:eq:A.40a},
\[
|\partial_iP_s^tf(x)|
\le
 \mathrm{e}^{-c_0ns}
\sum_j \mathrm{e}^{-a\mathrm d(i,j)}
 P_s^t(|\partial_jf|)(x).
\]
The Markov property, Jensen's inequality, and invariance of \(\mu_t\) give
\[
\|P_s^t|g|\|_{L^p(\mu_t)}^p
\le
\int P_s^t(|g|^p)\,\mathrm{d}\mu_t
=
\|g\|_{L^p(\mu_t)}^p.
\]
Taking \(L^p(\mu_t)\)-norms in the preceding bound on each coordinate derivative and using Minkowski's inequality gives
\[
\|\partial_iP_s^tf\|_{L^p(\mu_t)}
\le
 \mathrm{e}^{-c_0ns}
\sum_j \mathrm{e}^{-a\mathrm d(i,j)}
\|\partial_jf\|_{L^p(\mu_t)}.
\eqtag{A.41}
\]
For the approximation step, choose smooth cutoffs $\chi_R$ equal to one on $\{|x|\le R\}$ and satisfying $|\nabla\chi_R|\le C/R$. The product rule gives
\[
 \|\partial_j(\chi_Rf)-\partial_jf\|_{L^p(\mu_t)}
 \le\|\mathbf1_{\{|x|>R\}}\partial_jf\|_{L^p(\mu_t)}
       +\frac C R\|\mathbf1_{\{|x|>R\}}f\|_{L^p(\mu_t)}\longrightarrow0.
\]
Also $\chi_Rf\to f$ in $L^p(\mu_t)$ and in $L^2(\mu_t)$ when the latter norm is required. For fixed $R$, mollification of $\chi_Rf$ converges in these norms and in each derivative
norm: on its compact support the smooth positive density is bounded above and below. A diagonal choice gives smooth compactly supported approximants. The $L^p$ contraction and
\eqref{A:eq:A.41}, applied to differences of these approximants, pass the derivative estimate to the limit.

The Brascamp--Lieb inequality following from \eqref{A:eq:A.37c} gives the Poincar\'e inequality
\(\operatorname{Var}_{\mu_t}(g) \le \frac1{cn} \int|\nabla g|^2\,\mathrm{d}\mu_t\). Since $f$ is centered, invariance
and the energy identity imply
\[
 \frac{\mathrm d}{\mathrm ds}\|P_s^tf\|_{L^2(\mu_t)}^2
 =-2\|\nabla P_s^tf\|_{L^2(\mu_t)}^2
 \le-2cn\|P_s^tf\|_{L^2(\mu_t)}^2.
\]
Integration gives $\|P_s^tf\|_2\le e^{-cns}\|f\|_2$; approximation in $L^2$ covers every centered $f$ in the statement. Define $\varphi_S=-\int_0^SP_s^tf\,\mathrm ds$ and
$\varphi=-\int_0^\infty P_s^tf\,\mathrm ds$ as $L^2(\mu_t)$ integrals. The semigroup property gives, for $\epsilon>0$,
\[
 \frac{P_\epsilon^t\varphi_S-\varphi_S}{\epsilon}
 =\frac1\epsilon\int_0^\epsilon P_s^tf\,\mathrm ds
    -\frac1\epsilon\int_S^{S+\epsilon}P_s^tf\,\mathrm ds
 \longrightarrow f-P_S^tf
 \quad\text{in }L^2(\mu_t).
\]
Thus $\varphi_S\in\operatorname{Dom}L_t$ and $L_t\varphi_S=f-P_S^tf$. Moreover,
\[
 \|\varphi-\varphi_S\|_2\le\frac{e^{-cnS}}{cn}\|f\|_2,
 \qquad \|L_t\varphi_S-f\|_2\le e^{-cnS}\|f\|_2.
\]
Closedness of the generator now gives $\varphi\in\operatorname{Dom}L_t$ and $L_t\varphi=f$. Invariance gives $\int\varphi\,\mathrm d\mu_t=0$. The centered solution is unique: if
\(L_tg=0\), then \(0=-\langle g,L_tg\rangle_{L^2(\mu_t)} =\int|\nabla g|^2\,\mathrm{d}\mu_t\), so \(g\) is constant and its centered version vanishes.

Differentiating $L_t\varphi=f$ gives
\[
 -L_t(\partial_i\varphi)
 +\sum_j(\partial_{ij}^2\Phi_t)\partial_j\varphi
 =-\partial_i f.
\]
Equation~\eqref{A:eq:A.41} gives, for $T>S$,
\[
 \|\partial_i\varphi_T-\partial_i\varphi_S\|_{L^p(\mu_t)}
 \le\frac{e^{-c_0nS}}{c_0n}\sum_j
       e^{-a\mathrm d(i,j)}\|\partial_jf\|_{L^p(\mu_t)}.
\]
Let $v_i$ be the resulting $L^p(\mu_t)$ limit. On compact sets the density is bounded below, so both convergences also hold in $L^1(\mathrm dx)$. For
$\eta\in C_c^\infty(\mathbb R^d)$,
\[
 \int\varphi\,\partial_i\eta\,\mathrm dx
 =\lim_{S\to\infty}\int\varphi_S\partial_i\eta\,\mathrm dx
 =-\lim_{S\to\infty}\int\partial_i\varphi_S\eta\,\mathrm dx
 =-\int v_i\eta\,\mathrm dx.
\]
Hence $v_i=\partial_i\varphi$ weakly, including when $p=1$. The differentiated equation above is understood in distributions. Applying \eqref{A:eq:A.41} to the integral yields
\[
\begin{aligned}
\|\partial_i\varphi\|_{L^p(\mu_t)}
&\le \int_0^\infty\|\partial_iP_s^tf\|_{L^p(\mu_t)}\,\mathrm{d}s
\\
&\le \left(\int_0^\infty \mathrm{e}^{-c_0ns}\,\mathrm{d}s\right)
\sum_j \mathrm{e}^{-a\mathrm d(i,j)}\|\partial_jf\|_{L^p(\mu_t)}
=\frac1{c_0n}\sum_j \mathrm{e}^{-a\mathrm d(i,j)}\|\partial_jf\|_{L^p(\mu_t)}.
\end{aligned}
\]
Equation~\eqref{A:eq:A.39} follows, while the bound on the number of coordinates at each distance gives \eqref{A:eq:A.39a}.
\end{proof}

It remains to choose the velocity so that its flow transports the interpolating densities exactly. If $\rho_t=Z_t^{-1}e^{-\Phi_t}$ and $f_t=\Psi-\mathbb E_{\mu_t}\Psi$, then
\[
 \partial_t\rho_t=-f_t\rho_t,\qquad
 \operatorname{div}(\rho_t\nabla\varphi_t)
       =\rho_tL_t\varphi_t.
\]
Thus the equation $L_t\varphi_t=f_t$ has precisely the sign required by $\partial_t\rho_t+\operatorname{div}(\rho_tv_t)=0$ with $v_t=\nabla\varphi_t$. We construct the transport
from this continuity equation by the flow argument in \S4 of \cite{moser1965volume}, applied on Euclidean space. The resulting flow carries $\rho_0(x)\,\mathrm dx$ to
$\rho_t(x)\,\mathrm dx$. The Jacobian identity is proved below. Related transports for interpolated potentials appear in \cite{bekerman2015transport}.

\begin{theorem}\label{A:res:A.14}
Under the assumptions above, there exists a coupling of \(U\sim Z^{-1}\mathrm{e}^{-\Phi(u)}\,\mathrm{d}u\) and \(G\sim N(0,K^{-1})\) such that, for every \(p<\infty\),
\[
\sup_i\|U_i-G_i\|_{L^p}
\le C_pn^{-1}.
\eqtag{A.42}
\]
\end{theorem}

\begin{proof}
For the interpolation \eqref{A:eq:A.37b}, set $f_t=\Psi-\mathbb E_{\mu_t}\Psi$, let $L_t\varphi_t=f_t$ be the centered solution, and put $v_t=\nabla\varphi_t$.
Lemmas~\ref{A:res:A.12} and~\ref{A:res:A.13} give
\[
 \sup_{0\le t\le1}\sup_i\|v_{t,i}\|_{L^p(\mu_t)}\le C_pn^{-1}.
\]
We prove $|v_t(x)|+\|D_xv_t(x)\|\le C(1+|x|)$ locally uniformly in $t$. These bounds give a global flow on $0\le t\le1$; its Jacobian will identify the transported density.

Fix $t$ temporarily. Write $X_s^x$ for the Langevin process in the proof of Lemma~\ref{A:res:A.13} and $J_s^x=D_xX_s^x$. Strong convexity and It\^o's formula give
\[
 \mathbb E|X_s^x|^2\le e^{-2cns}|x|^2+\frac d{cn},
 \qquad \|J_s^x\|_{\mathrm{op}}\le e^{-cns}.
\]
Here $d\le Cn$. Since $\nabla\Psi(0)=0$ and $\|\nabla^2\Psi\|_{\mathrm{op}}\le Cn$, the semigroup representation is
\[
 v_t(x)=-\int_0^\infty
       \mathbb E\bigl[(J_s^x)^*\nabla\Psi(X_s^x)\bigr]\,\mathrm ds,
 \qquad |v_t(x)|\le C(1+|x|).
\]
All constants in these estimates are uniform in $t$.

The coordinate bounds on third derivatives also give the operator bound
\[
 \|\nabla^3\Phi_t(x)[v,w]\|_2\le Cn\|v\|_2\|w\|_2.
\]
Indeed, the $i$th component is bounded by $Cn(\sum_{j\in N_i}|v_j|^2)^{1/2} (\sum_{k\in N_i}|w_k|^2)^{1/2}$, where $N_i$ is a fixed radius neighborhood. Its cardinality and the
number of such neighborhoods containing any coordinate are bounded. Square, sum in $i$, and use $\sum_{k\in N_i}|w_k|^2\le\|w\|_2^2$ to obtain the derivative bound. For
$Q_s^x[v,w]=D_x^2X_s^x[v,w]$, differentiation of the flow equation gives
\[
 \partial_sQ_s^x[v,w]
 =-\nabla^2\Phi_t(X_s^x)Q_s^x[v,w]
  -\nabla^3\Phi_t(X_s^x)[J_s^xv,J_s^xw],\qquad Q_0^x=0.
\]
The evolution generated by $-\nabla^2\Phi_t(X_s^x)$ has norm at most $e^{-cn(s-r)}$ between times $r$ and $s$. Variation of constants yields
\[
 \|Q_s^x\|_{\mathrm{op}}
 \le Cn\int_0^s e^{-cn(s-r)}e^{-2cnr}\,\mathrm dr
 \le C e^{-cns}.
\]
Consequently differentiation of the velocity integral is justified uniformly on compact sets by the bound
\[
 \|D_xv_t(x)\|_{\mathrm{op}} \le Cn\int_0^\infty
 \left[e^{-cns}\bigl(e^{-cns}|x|+C\bigr)+e^{-2cns}\right]\,\mathrm ds \le C(1+|x|).
\]
The two integrands arise from $(D_xJ_s^x)^*\nabla\Psi(X_s^x)$ and $(J_s^x)^*\nabla^2\Psi(X_s^x)J_s^x$, respectively. Here the differentiations on a finite interval $[0,S]$ can be
justified directly. For a unit vector $e$ and $\epsilon\ne0$, put $D_s^\epsilon=(X_s^{x+\epsilon e}-X_s^x)/\epsilon$. The difference of the two stochastic equations gives
\[
 D_s^\epsilon=e-\int_0^s H_r^\epsilon D_r^\epsilon\,\mathrm dr,
 \qquad
 H_r^\epsilon=\int_0^1\nabla^2\Phi_t
 (X_r^x+u(X_r^{x+\epsilon e}-X_r^x))\,\mathrm du.
\]
The Hessian bound gives $\sup_{s\le S}|X_s^{x+\epsilon e}-X_s^x| \le e^{CnS}|\epsilon|$ and $\sup_{s\le S}|D_s^\epsilon|\le e^{CnS}$. Its Lipschitz bound, which follows from the
bound on $\nabla^3\Phi_t$, gives $\sup_{s\le S}\|H_s^\epsilon-\nabla^2\Phi_t(X_s^x)\| \le C_{n,S}|\epsilon|$. Subtracting the equation for $J_s^xe$ and applying the integral
inequality therefore yields
\[
 \sup_{s\le S}|D_s^\epsilon-J_s^xe|
 \le C_{n,S}|\epsilon|\longrightarrow0.
\]
Taking difference quotients in the equation for $J_s^x$ gives the displayed equation for $Q_s^x$; continuity and boundedness of $\nabla^3\Phi_t$ on the ranges of these paths
justify the limit in its forcing term. This proves the finite time derivatives.

To check the parameter dependence, drive $X^{t,x}$ and $X^{u,y}$ by the same Brownian motion. Since $\nabla\Phi_t-\nabla\Phi_u=(t-u)\nabla\Psi$, their integral equations imply
\[
 \sup_{s\le S}|X_s^{t,x}-X_s^{u,y}|
 \le e^{CnS}\left(|x-y|+Cn|t-u|
                 \int_0^S|X_r^{u,y}|\,\mathrm dr\right).
\]
On bounded sets of $x,y$, the right side tends to zero in $L^2$ as $(t,x)\to(u,y)$: the second moment bound above and the finite time integral equation bound the integral
uniformly. The equations for $J$ and $Q$ then give convergence of these variations on $[0,S]$. For $Q$, first restrict the paths to a fixed ball, use uniform continuity of the
third derivatives there, and then let the ball increase; the displayed uniform bounds on $J,Q$ and the second moments of $X$ justify this last limit. Finally, for $|x|\le R$ the
omitted parts of the velocity and its derivative satisfy
\[
 \int_S^\infty Cn\bigl[e^{-cns}(1+e^{-cns}R)
                         +e^{-2cns}\bigr] \,\mathrm ds
 \le C_R e^{-cnS}.
\]
Finite time continuity and this uniform tail bound prove continuity of $v_t(x)$ and $D_xv_t(x)$ in $(t,x)$, and justify both infinite time differentiations.

Let $T_t(x)$ solve $\dot T_t(x)=v_t(T_t(x))$, $T_0(x)=x$. Local existence and uniqueness follow from the local derivative bound. Gronwall gives $1+|T_t(x)|\le e^{Ct}(1+|x|)$;
solving backwards from any time in $[0,1]$ gives the same bound. Thus $T_t$ is a global $C^1$ diffeomorphism. Its Jacobian determinant is positive and satisfies
\[
 \det D_xT_t(x)
 =\exp\left\{\int_0^t\operatorname{div}v_r(T_r(x))\,\mathrm dr\right\}.
\]
Write $\rho_t=Z_t^{-1}e^{-\Phi_t}$. Differentiation under its normalizing integral is legitimate by the common quadratic lower bound and the growth of at most quadratic order of
$\Psi$. The Poisson equation gives
\[
 \partial_t\rho_t=-f_t\rho_t,
 \qquad \operatorname{div}(\rho_tv_t)=\rho_tL_t\varphi_t=f_t\rho_t.
\]
The weak Poisson identity is pointwise here: $v_t$ is $C^1$ and $\rho_t$ is positive and smooth. The chain rule therefore implies
\[
 \frac{\mathrm d}{\mathrm dt}
 \left[\rho_t(T_t(x))\det D_xT_t(x)\right]=0.
\]
At $t=0$ this expression equals $\rho_0(x)$. Change of variables now proves $(T_t)_\#\mu_0=\mu_t$.

Take $G\sim\mu_0=N(0,K^{-1})$ and $U=T_1(G)$. By Minkowski,
\[
 \|U_i-G_i\|_{L^p}
 \le\int_0^1\|v_{t,i}(T_t(G))\|_{L^p}\,\mathrm dt
 =\int_0^1\|v_{t,i}\|_{L^p(\mu_t)}\,\mathrm dt
 \le C_pn^{-1}.
\]
This proves \eqref{A:eq:A.42}. The map $T_1$ is independent of $p$.
\end{proof}

For conditional coefficient laws on a positive domain, we combine the transport with Lemma~\ref{A:res:A.11D}. Let
\(\mathcal D=\mathbb R^{d_A}\times(0,\infty)^{d_B}\), with \(d_A+d_B\asymp n\), and let
\(\mu(\mathrm{d}y)=Z^{-1}\mathrm{e}^{-\Phi(y)}\mathbf 1_{\mathcal D}(y)\,\mathrm{d}y\). Assume that \(\Phi\) is
finite-range and \(\nabla^2\Phi\ge cnI\) on \(\mathcal D\), and let \(y^\circ\) be its minimizer. Suppose every positive coordinate satisfies \(y_j^\circ\ge4a\) for one fixed
\(a>0\), and suppose
\[
\sup_i|\mathbb E(Y_i-y_i^\circ)|\le \frac{C_0}{n},
\quad Y\sim\mu.
\eqtag{A.42d}
\]
Suppose further that a globally defined \(C^3\), finite-range extension \(\widetilde\Phi\) agrees with \(\Phi\) whenever every positive coordinate is at least \(a\) and
\(\|y-y^\circ\|_{\ell^\infty}\le R\), and that
\[
cnI\le\nabla^2\widetilde\Phi(y)\le CnI,
\quad
|\partial_{ijk}^3\widetilde\Phi(y)|\le Cn,
\]
where mixed derivatives vanish unless all their indices lie within a fixed graph distance of one another. Assume that the numbers of vertices at each distance are uniformly
bounded, as in Proposition~\ref{A:res:A.11B}, and that
\(\nabla\widetilde\Phi(y^\circ)=0\) and \(K=\nabla^2\Phi(y^\circ)=\nabla^2\widetilde\Phi(y^\circ)\).

\begin{proposition}\label{A:res:A.14A}
There is a coupling of the centered variable \(U=Y-y^\circ\) and \(G\sim N(0,K^{-1})\) such that, for every finite
\(p\),
\[
\sup_i\|U_i-G_i\|_{L^p}\le C_pn^{-1}.
\eqtag{A.42e}
\]
\end{proposition}

\begin{proof}
Let \(\widetilde\mu\propto \mathrm{e}^{-\widetilde\Phi}\), and let \(J_+\) index the positive coordinates. Put
\[
 \mathcal C=\{y:y_j\ge a\ (j\in J_+),
                    \ \|y-y^\circ\|_\infty\le R\}.
\]
For \(Y\sim\mu\), \(\widetilde Y\sim\widetilde\mu\), hypothesis \eqref{A:eq:A.42d}, Proposition~\ref{A:res:A.11B}, and the logarithmic Sobolev bound on the convex domain give, for
every fixed \(q<\infty\),
\[
 \|\mathbb EY-y^\circ\|_\infty+ \|\mathbb E\widetilde Y-y^\circ\|_\infty \le C/n,\qquad
 \sup_i\bigl(\|Y_i-y_i^\circ\|_q+ \|\widetilde Y_i-y_i^\circ\|_q\bigr) \le C_qn^{-1/2}.
\]
Put $r_*=\min\{3a,R\}>0$. For either variable $Z=Y$ or $Z=\widetilde Y$, the event $Z\notin\mathcal C$ implies
$|Z_i-y_i^\circ|\ge r_*$ for some $i$, since $y_j^\circ\ge4a$ in every positive coordinate. For $n$ such that $C/n\le
r_*/2$,
\[
 \mathbb P(Z\notin\mathcal C)
 \le\sum_{i=1}^d\mathbb P(|Z_i-\mathbb EZ_i|\ge r_*/2)
 \le2d\exp(-cnr_*^2)\le C\exp(-c'n).
\]
Here $d\le Cn$ and $r_*$ is fixed. Adding the two bounds gives
\[
 \mu(\mathcal C^c)+\widetilde\mu(\mathcal C^c)\le C\mathrm{e}^{-cn}.
\]
The unnormalized densities agree on \(\mathcal C\), so their conditional laws there coincide. Lemma~\ref{A:res:A.11D}, used with \(q=2p\) and coordinates centered at \(y^\circ\),
provides a coupling with
\[
 \mathbb P(Y\ne\widetilde Y)\le C\mathrm{e}^{-cn},\quad
 \sup_i\|Y_i-\widetilde Y_i\|_p
       \le C_pn^{-1/2}\mathrm{e}^{-cn/(2p)}.
\]
Theorem~\ref{A:res:A.14} couples \(\widetilde Y-y^\circ\) to \(G\sim N(0,K^{-1})\) with coordinate error \(C_p/n\). Let $T(G)$ denote this extended variable and let
$\kappa(\widetilde y,\mathrm dy)$ be the conditional kernel in Lemma~\ref{A:res:A.11D}, conditional on the extended variable. The probability measure
\[
 \gamma(\mathrm dg)\,
 \kappa\bigl(y^\circ+T(g),\mathrm dy\bigr)
\]
has marginals $N(0,K^{-1})$ and $\mu$. On this space put $U=Y-y^\circ$. Then
\[
 \|U_i-G_i\|_p
 \le\|Y_i-\widetilde Y_i\|_p
       +\|\widetilde Y_i-y_i^\circ-G_i\|_p
 \le C_p n^{-1/2}\mathrm e^{-cn/(2p)}+C_pn^{-1}
 \le C_p'n^{-1}.
\]
This is \eqref{A:eq:A.42e}. The transport map and the kernel in Lemma~\ref{A:res:A.11D} are fixed before choosing the moment order. Thus a single coupling satisfies all of the
displayed $L^p$ estimates.
\end{proof}

In the application to the Jacobi coefficients, Lemma~\ref{A:res:mean} supplies the displacement of the conditional mean. We construct the extension and identify its Hessian in
Subsection~\ref{A:sec:gaussian-approximation-of-the-conditional-gibbs-measure}. That subsection also applies the local scalar map and converts the coordinate error into the
weighted sum required by the resolvent estimate.

\subsubsection{Partial sums under the conditional Gibbs law}\label{A:sec:primitive-bounds-for-conditional-gibbs-measures}
The form comparison needs control of cumulative fluctuations, not just individual coordinates. Strong convexity supplies concentration for these partial sums, while locality and
bounded overlap keep the constants independent of the number of coefficients.

Let \(\mu_n^\tau\) be a probability measure on a convex domain in \(\mathbb R^{d_n}\), possibly conditional on boundary data \(\tau\), with potential \(\Phi_n^\tau\) satisfying
\(\nabla^2\Phi_n^\tau\ge c n I\) uniformly in \(n\) and in the allowed boundary data. Let \(\ell_j\in\mathbb R^{d_n}\) be a family of local row vectors satisfying, for fixed
constants
\(R,C<\infty\),
\[
\operatorname{supp}\ell_j\subset\{i:d_{\mathcal G}(i,j)\le R\},
\quad
\|\ell_j\|_2\le C,
\eqtag{A.23b}
\]
and assume bounded overlap of their supports. For \(U\sim\mu_n^\tau\), put
\[
\bar\zeta_j
=
\ell_j\cdot U-\mathbb E_{\mu_n^\tau}(\ell_j\cdot U),
\quad
S_m=h^{-1}\sum_{j=1}^{m}\bar\zeta_j,
\quad
h^3\asymp n^{-1}.
\]
Extend \(S_m\) piecewise linearly in the coordinate \(x=mh\), and keep it constant past the finite boundary point.

\begin{lemma}\label{A:res:A.7A}
For every \(\alpha>0\) and every \(1\le p<\infty\),
\[
\sup_n\sup_\tau
\left\|
\sup_{q\ge0}(1+q)^{-\alpha}
\left(
\operatorname{osc}_{[q,q+1]}S
+
|S(q+1)-S(q)|
\right)
\right\|_{L^p(\mu_n^\tau)}
<\infty.
\eqtag{A.23c}
\]
The same conclusion holds without centering each \(\zeta_j\), provided
\[
\sup_j|\mathbb E\zeta_j|\le Cn^{-1}.
\eqtag{A.23d}
\]
\end{lemma}

\begin{proof}
Apply the concentration calculation preceding \eqref{A:eq:convex-concentration} to $a\cdot U$. Its gradient has norm $\|a\|_2$, so, with $Z=a\cdot(U-\mathbb EU)$,
\[
 \log\mathbb E e^{tZ}\le\frac{Ct^2\|a\|_2^2}{n},\qquad
 \mathbb P(|Z|>u)\le2e^{-cn u^2/\|a\|_2^2}.
\]
For $a\ne0$, tail integration gives
\[
 \mathbb E\exp\!\left(\frac{nZ^2}{C_1\|a\|_2^2}\right)
 =1+\int_1^\infty
       \mathbb P\!\left(|Z|>\|a\|_2\sqrt{C_1\log v/n}\right)\,\mathrm dv
 \le1+2\int_1^\infty v^{-cC_1}\,\mathrm dv\le2
\]
when $C_1$ is sufficiently large. For $a=0$, $Z=0$ and $\|Z\|_{\psi_2}=0$. Thus $\|Z\|_{\psi_2}\le Cn^{-1/2}\|a\|_2$. If at most $m_0$ row supports contain a given coordinate, then
\[
 \left\|\sum_{j=r+1}^s\ell_j\right\|_2^2
 \le m_0\sum_{j=r+1}^s\|\ell_j\|_2^2\le C(s-r).
\]
Consequently,
\[
\|S_s-S_r\|_{\psi_2}
\le
C n^{-1/2}h^{-1}(s-r)^{1/2}
\le
C\bigl(h(s-r)\bigr)^{1/2}.
\eqtag{A.23e}
\]
Apply the subgaussian part of Lemma~\ref{A:res:A.7} with \(\xi_{j,n}=h^{-1}\bar\zeta_j\). Estimate \eqref{A:eq:A.23e} is \eqref{A:eq:A.23-subgaussian} for these increments. The
interpolation in \eqref{A:eq:A.23-interpolation} agrees with the partial sum defined above, including its constant extension at the finite boundary point. The conditional version
of \eqref{A:eq:A.23-tail} gives the weighted block tail with constants uniform in \(n\) and in the allowed boundary data. Its integrated form proves \eqref{A:eq:A.23c}.

For the uncentered sums, subtract the centered partial sums from the uncentered ones. On the interval $[(j-1)h,jh)$ this deterministic difference has slope
\(h^{-2}\mathbb E\zeta_j\); past the finite boundary point its slope is zero. Under \eqref{A:eq:A.23d}, its absolute slope is at most \(C/(nh^2)\le C'h\), since
\(h^3\asymp n^{-1}\). Thus its oscillation plus difference of the boundary values on every interval $[p,p+1]$ is at most \(2C'h\). This bound also covers partial intervals and
blocks crossing the finite boundary point, and the weight \((1+q)^{-\alpha}\) can only decrease it. Adding this deterministic contribution proves the bound for the uncentered sums.
\end{proof}

\subsubsection{Gaussian approximation of the conditional Gibbs measure}\label{A:sec:gaussian-approximation-of-the-conditional-gibbs-measure}

We now apply the transport estimate to the exact conditional coefficient law. The Gaussian is centered at its conditional minimizer and has inverse Hessian covariance, as in
\cite{krishnapur2016universality}. We require a coupling whose error is $C_pn^{-1}$ in each coordinate, on an interval containing order $n$ variables. We verify the extension and
transport hypotheses before making this comparison.

The interval used for transport is larger than the interval on which Gaussian coefficients replace the scalar potential. The extra sites separate the latter interval from both
conditional boundaries. The coupling itself is defined on the entire transport interval; outside $I_n^{\mathrm{mid}}$, the comparison operator keeps the original scalar
coefficients.

Fix once and for all \(0<s_{\mathrm{reg}}<s_0<1\), where \(s_0\) is the regular cutoff used in Lemma~\ref{A:res:mean}. We use the same integer $R_V$ as in the introduction. Retain
the initial interval from \eqref{A:eq:13.17b},
\(L_n^{\mathrm L}=\lceil A_{\mathrm L}\log n\rceil, \allowbreak\, J_n^{\mathrm L}=L_n^{\mathrm L}+3R_V+10\), and choose
independently \(L_n^{\mathrm R}=\lceil A_{\mathrm R}\log n\rceil, \allowbreak\, g_n^{\mathrm R}=L_n^{\mathrm R}+3R_V+10\). The constants \(A_{\mathrm L},A_{\mathrm R}\) will be
fixed large enough below. If \(s_{j,n}\) denotes the macroscopic coordinate of the \(j\)-th site after reindexing, define
\[
J_n^{\mathrm{reg}}
=
\max\left\{
J_n^{\mathrm L}+1\le j\le M:
 s_{j+g_n^{\mathrm R},n}\le s_{\mathrm{reg}}
\right\},
\]
and \(J_n^{\mathrm{tr}} = J_n^{\mathrm{reg}}+g_n^{\mathrm R}.\) For all sufficiently large \(n\), \(J_n^{\mathrm L}<J_n^{\mathrm{reg}}<J_n^{\mathrm{tr}}<M.\) We write
\[
I_n^{\mathrm L}
=
\{1,\ldots,J_n^{\mathrm L}\},
\quad
I_n^{\mathrm{mid}}
=
\{J_n^{\mathrm L}+1,\ldots,J_n^{\mathrm{reg}}\},
\]
\[
I_n^{\mathrm{tr}}
=
\{1,\ldots,J_n^{\mathrm{tr}}\},
\quad
I_n^{\mathrm{far}}
=
\{J_n^{\mathrm{reg}}+1,\ldots,M\}.
\]
The characteristic functions of the first two regions are
\[
\chi_{j,n}^{\mathrm L}=\mathbf 1_{\{j\in I_n^{\mathrm L}\}},
\quad
\chi_{j,n}^{\mathrm{mid}}=\mathbf 1_{\{j\in I_n^{\mathrm{mid}}\}}.
\]

For the vertex coordinate $s_{j,n}$ in \eqref{A:eq:recentered-profile-coordinate}, put \(\lambda_n=rh+h/\beta\). Since
\(h^3=\tau_\star/(b_\star n)\),
\[
 s_{j,n} =\frac{r+j-\tfrac12+\beta^{-1}}{n} =\frac{b_\star}{\tau_\star}h^2(x_j+\lambda_n),\quad
 s_{j,n}-\frac{b_\star}{\tau_\star}h^2x_j =\frac{b_\star}{\tau_\star}h^3(r+\beta^{-1}), \quad
 x_j=(j-\tfrac12)h.
\eqtag{exact-profile-shift}
\]
Since $\lambda_n=O(h\log(n+1))$, the last difference is $O(h^3\log(n+1))$. We also have
\[
 0\le s_{j,n}\le C h^2(1+x_j),\quad
 s_{j+1,n}-s_{j,n}=\frac1n=\frac{b_\star}{\tau_\star}h^3,
 \quad 1\le j\le J_n^{\mathrm{tr}},
\eqtag{profile-size-increment}
\]
for all sufficiently large \(n\), uniformly on every compact positive \(\beta\)-interval. Here we used
\(0\le\lambda_n\le1\), which holds for the same range of \(n\).

The cutoff definition gives
\[
 s_{\mathrm{reg}}-\frac1n<s_{J_n^{\mathrm{tr}},n}
 \le s_{\mathrm{reg}},\quad
 x_{J_n^{\mathrm{reg}}}
 =\frac{\tau_\star s_{\mathrm{reg}}}{b_\star}h^{-2}
    -\lambda_n-g_n^{\mathrm R}h+O(h).
\]
Here maximality is used only after the cutoff lies strictly inside \(\{1,\ldots,M\}\), which follows from
\(s_{\mathrm{reg}}<1\). Both logarithmic intervals have length in the coordinate $x$ \(O(h\log n)=o(1)\). Hence
deterministic constants \(0<c_{\mathrm{reg}}<C_{\mathrm{reg}}\) give \(c_{\mathrm{reg}}h^{-2}\le x_{J_n^{\mathrm{reg}}} \le x_{J_n^{\mathrm{tr}}}\le C_{\mathrm{reg}}h^{-2}\) for
all sufficiently large \(n\). Also $J_n^{\mathrm L}h=O(h\log n)$. All minimizers and covariance symbols below are evaluated at $s_{j,n}$.

Let \(\partial_{\mathrm L}I_n^{\mathrm{tr}}\) and \(\partial_{\mathrm R}I_n^{\mathrm{tr}}\) denote the two boundaries of the graph of free coefficient coordinates. Every diagonal
coefficient indexed by \(I_n^{\mathrm{mid}}\) satisfies
\[
\operatorname{dist}
\bigl(i,\partial_{\mathrm L}I_n^{\mathrm{tr}}\bigr)
\ge L_n^{\mathrm L}+R_V+4,
\eqtag{L10}
\]
and
\[
\operatorname{dist}
\bigl(i,\partial_{\mathrm R}I_n^{\mathrm{tr}}\bigr)
\ge L_n^{\mathrm R}+R_V+4.
\eqtag{L11}
\]
These are the two boundary separations used in the approximation of the inverse Hessian. The translation-invariant covariance approximation is restricted to \(I_n^{\mathrm{mid}}\).

Let $\mathscr F_n^{\mathrm{out}}$ be generated by the external boundary data of Section~\ref{A:sec:boundary-reduction} and by the reduced coordinates
\[
 (A_j)_{J_n^{\mathrm{tr}}<j\le M},\qquad
 (B_j)_{J_n^{\mathrm{tr}}<j<M}.
\]
Thus $B_{J_n^{\mathrm{tr}}}$ is free; it is not included in the exterior sigma-field. Conditional on $\mathscr F_n^{\mathrm{out}}$, the variables to be transported are
\(Y_n^{\mathrm{tr}} = (A_1,\ldots,A_{J_n^{\mathrm{tr}}}, B_1,\ldots,B_{J_n^{\mathrm{tr}}})\) on the open convex domain
\[
\mathcal D_n^+
=
\mathbb R^{J_n^{\mathrm{tr}}}
\times
(0,\infty)^{J_n^{\mathrm{tr}}}.
\eqtag{14.2}
\]
We denote the vector of conditioned variables by \(\tau\). The exact conditional density is
\(\mathrm{d}\mu_{n,+}^{\tau}(y) = \frac1{Z_{n,+}^{\tau}} \exp\{-\Phi_{n,+}^{\tau}(y)\} \mathbf 1_{\mathcal
D_n^+}(y)\,\mathrm{d}y\). The conditional potential is the restriction of \(\frac{\beta n}{2}\mathscr H_n\) to the transport variables, with all exterior variables fixed at their
actual values \(\tau\). In particular, \(\tau\) retains
\(b_r\) and \(b_N\); the auxiliary boundary values \(B_0=B_M=b_\star\) do not replace either coefficient in this energy,
by Corollary~\ref{A:res:conditioning}. The cutoff is defined using vertex parameters. For the edge with the same site index the parameter is larger by $1/(2n)$. All transported
edges, and all coordinates within $R_V$ indices of them, therefore have parameter at most $s_{\mathrm{reg}}+C_{R_V}/n<s_0$ for large $n$. In particular, their logarithmic
coefficients have the deterministic lower bound
\[
1-\frac{k}{n}-\frac1{\beta n}
\ge
(1-s_{\mathrm{reg}})/2>0
\eqtag{14.4}
\]
for every transported off-diagonal coefficient coordinate. The polynomial part is uniformly convex. Therefore
\(\nabla^2\Phi_{n,+}^{\tau}(y) \ge c n I \allowbreak\; (y\in\mathcal D_n^+)\) for every exterior value in the domain
with real diagonal and positive off-diagonal coefficients.

For the conditional constructions here and in \hyperref[A:sec:localized-inverse-hessians-and-weighted-gaussian-couplings]{the comparison of inverse Hessians}, \(\tau\) is the
vector of actual coefficients fixed by \(\mathscr F_n^{\mathrm{out}}\). Use
\[
 \mathcal T_n^{\mathrm{out}}(M_0)
 =\{\tau:\|\tau\|_\infty\le M_0\}.
\eqtag{good-exterior}
\]
This set also controls the exterior data of the full reduced interval; it bounds more than just the exterior coefficients within $R_V$ indices of the interval on which transport is
applied. It is measurable in the conditioning sigma-field, and its complement has probability at most \(C\mathrm{e}^{-cn}\) by Lemma~\ref{A:res:raw}. In all conditional statements
in this transport construction, fix \(\tau\in\mathcal T_n^{\mathrm{out}}(M_0)\). Their constants are deterministic, uniform in these data and in \(n\). Lemma~\ref{A:res:12.1A},
applied to the conditional minimizer, bounds its coordinate magnitudes and bounds its off-diagonal coordinates away from zero. We compare the two conditional minimizers separately
on the event \(\mathcal D_{n,Q}\). To recover an unconditional probability, use \eqref{A:eq:integrated-tail}; the probability of the complementary set of exterior data is at most
\(C\mathrm{e}^{-cn}\). This convention is also used in the conditional covariance construction below. It never denotes
conditioning the Gibbs law on \(\mathcal D_{n,Q}\) or on \(\Omega_{\mathrm{bd}}\).

Let
\[
\overline Y_n^\circ
=
(\overline A_1^\circ,\ldots,\overline A_{J_n^{\mathrm{tr}}}^\circ,
 \overline B_1^\circ,\ldots,\overline B_{J_n^{\mathrm{tr}}}^\circ)
\]
be the unique minimizer of \(\Phi_{n,+}^{\tau}\) on \(\mathcal D_n^+\). We reserve \(Y_n^\circ = (A_1^\circ,\ldots,A_M^\circ,B_1^\circ,\ldots,B_{M-1}^\circ)\) for the conditional
minimizer on the full reduced interval introduced in Section~\ref{A:sec:decomposition-of-the-general-jacobi-coefficients}.

The Gaussian approximation must be centered at $\overline Y_n^\circ$: this is the stationary point of the density being transported. The scalar decomposition, however, uses
$Y_n^\circ$. The two energies have the same left boundary data and differ only through terms meeting the right boundary of $I_n^{\mathrm{tr}}$. We keep their difference separate
from the transport error. The next estimate shows how the additional sites at the right boundary make that difference negligible on $I_n^{\mathrm{mid}}$.

\par\Needspace{7\baselineskip}
\begin{proposition}\label{A:res:14.1}
On the finite coefficient event \(\mathcal D_{n,Q}\), there are deterministic \(c>0\) and \(C_Q<\infty\) such that
\[
\left|
\overline Y_{n,i}^\circ-Y_{n,i}^\circ
\right|
\le
C_Q n^{-1/2}\sqrt{\log n}\,
 \mathrm{e}^{-c\,\operatorname{dist}(i,\partial_{\mathrm R}I_n^{\mathrm{tr}})}
\eqtag{14.7}
\]
for every transported coefficient coordinate \(i\).
\end{proposition}  By choosing \(A_{\mathrm R}\) sufficiently large,
\[
\max_{j\in I_n^{\mathrm{mid}}}
\left(
|\overline A_j^\circ-A_j^\circ|
+
|\overline B_j^\circ-B_j^\circ|
\right)
\le
n^{-4}
\eqtag{14.8}
\]
almost surely for all sufficiently large \(n\).

Independently of this comparison of conditional minimizers, for every exterior value in \(\mathcal T_{n,I_n^{\mathrm{tr}}}(M_0)\), there is a deterministic number
\(b_{\mathrm{reg}}>0\) such that
\[
\min_{1\le j\le J_n^{\mathrm{tr}}}
\overline B_j^\circ
\ge
4b_{\mathrm{reg}}
\eqtag{14.10}
\]
for every sufficiently large \(n\), uniformly in those data.

\begin{proof}
Put $z=Y_n^\circ|_{I_n^{\mathrm{tr}}}$, with both coefficient types retained, and $\delta=\overline Y_n^\circ-z$. Let $\Phi^0$ denote the shorter conditional potential whose right
exterior values are taken from the full minimizer. Then $\nabla\Phi^0(z)=0$. The actual potential differs only through terms meeting the right boundary. Its gradient at $z$ is
\[
 g_i=\partial_i\Phi_{n,+}^{\tau}(z)-\partial_i\Phi^0(z),
 \qquad g_i=0\quad\text{if }\operatorname{dist}
 (i,\partial_{\mathrm R}I_n^{\mathrm{tr}})>R_V.
\]
On $\mathcal D_{n,Q}$, the changed exterior coordinates differ from the full minimizer by at most $C_Qn^{-1/2}\sqrt{\log n}$. The event $\mathcal D_{n,Q}$ is contained in
$\mathcal E_n$, which bounds the original coefficients by the fixed number $M_0$. Lemma~\ref{A:res:12.1A} therefore confines both minimizers to a fixed cube independently of $Q$.
The derivatives of each polynomial interaction are bounded there by $Cn$. The logarithmic derivative depends only on its own free coordinate and cancels in $g_i$. The fundamental
theorem of calculus in the fixed exterior coordinates consequently gives
\[
 \|g\|_\infty\le C_Q\sqrt{n\log n}.
\]
There are at most a fixed number of nonzero entries of $g$. The free gradient equation is
\[
 0=g+\overline K\delta,\qquad
 \overline K=\int_0^1\nabla^2\Phi_{n,+}^{\tau}(z+t\delta)\,\mathrm dt.
\]
By Lemma~\ref{A:res:12.1A}, this segment lies in a bounded coefficient cube independent of $Q$. The polynomial part of $\overline K$ has fixed range and row sums at most $Cn$; the
logarithmic part is nonnegative and diagonal. The lower bound is $cnI$. Lemma~\ref{A:res:A.9} therefore yields $|\overline K^{-1}_{ij}|\le Cn^{-1}e^{-c\mathrm d(i,j)}$, with $c,C$
independent of $Q$. Let $S_{\mathrm R}=\operatorname{supp}g$, with $|S_{\mathrm R}|\le C$. Its elements are within $R_V$ of the right conditional boundary, so
\[
 |\delta_i|
 \le Cn^{-1}\|g\|_\infty\sum_{j\in S_{\mathrm R}}\mathrm e^{-c\mathrm d(i,j)}
 \le C_Qn^{-1/2}\sqrt{\log n}\,
       \mathrm e^{-c\operatorname{dist}(i,\partial_{\mathrm R}I_n^{\mathrm{tr}})}.
\]
The factor $\mathrm e^{cR_V}$ is included in $C_Q$. This proves \eqref{A:eq:14.7}. Choose $cA_{\mathrm R}>5$. Then $n^{-1/2}\sqrt{\log n}\,e^{-cL_n^{\mathrm R}}=o(n^{-4})$, and
Borel--Cantelli for $\mathcal D_{n,Q}$, with $Q>1$, gives \eqref{A:eq:14.8}.

To prove \eqref{A:eq:14.10} without this comparison, fix the other coordinates at the shorter minimizer and write $g_k(t)=\operatorname{Tr}V(T^{(k,0)}+tS_k)$. This is even, and its
second derivative is bounded by $C$ on the bounded coefficient cube. Stationarity at $b_k^\circ$ gives
\[
 \frac{2(1-k/n-1/(\beta n))}{b_k^\circ}
 =g_k'(b_k^\circ)=\int_0^{b_k^\circ}g_k''(t)\,\mathrm dt
 \le Cb_k^\circ.
\]
Together with \eqref{A:eq:14.4}, this implies $(b_k^\circ)^2\ge c>0$. Taking $4b_{\mathrm{reg}}\le\sqrt c$ proves the claim, including the first and last free edges. The comparison
of conditional minimizers in \cite{krishnapur2016universality} uses this Hessian estimate. The displayed gradient term accounts for the error in this comparison.

There are two uses of the boundary comparison. The factor $n^{-1/2}\sqrt{\log n}$ in \eqref{A:eq:14.7} is available on $\mathcal D_{n,Q}$, which involves the free variables. For
conditional moment estimates one instead fixes $\tau\in\mathcal T_n^{\mathrm{out}}(M_0)$. The same gradient calculation then gives $\|g\|_\infty\le Cn$, because both the fixed
exterior coordinates and the full minimizer are bounded. Consequently
\[
 \max_{j\in I_n^{\mathrm{mid}}}
 (|\overline A_j^\circ-A_j^\circ|+|\overline B_j^\circ-B_j^\circ|)
 \le C e^{-c g_n^{\mathrm R}}.
\]
Choosing the deterministic logarithmic length so that $c g_n^{\mathrm R}\ge5\log n$ makes this at most $n^{-4}$ for all large $n$, uniformly in the fixed admissible exterior. Thus
the center error can be included in a conditional $L^p$ estimate without conditioning on $\mathcal D_{n,Q}$.

\end{proof}

For conditional expectations at a fixed size, we use the bound just proved before restricting to the event $\mathcal D_{n,Q}$, which also involves free coefficients. On
\(\mathcal T_n^{\mathrm{out}}(M_0)\), both conditional minimizers and all changed exterior values are bounded deterministically. The same gradient comparison has a gradient error
of size \(Cn\) at the right boundary, so Lemma~\ref{A:res:A.9} gives
\[
 |\overline Y_{n,i}^\circ-Y_{n,i}^\circ|
 \le C \mathrm{e}^{-c\operatorname{dist}(i,\partial_{\mathrm R}I_n^{\mathrm{tr}})}.
\eqtag{finite-center-comparison}
\]
Increasing \(A_{\mathrm R}\) if necessary makes \eqref{A:eq:14.8} hold uniformly on this set of exterior data, for every sufficiently large \(n\). This bound applies under each
original conditional law with exterior data in the specified set.

Set \(U_n^+ = Y_n^{\mathrm{tr}}-\overline Y_n^\circ.\) Its support is the translated positive domain
\[
\mathcal U_n^+
=
\left\{
(u^A,u^B):
 u_j^A\in\mathbb R,
\quad
u_j^B> -\overline B_j^\circ,
\quad
1\le j\le J_n^{\mathrm{tr}}
\right\}.
\]
For \(j\in I_n^{\mathrm{mid}}\), center the diagonal coefficients at the conditional minimizer on the transport interval
\[
\zeta_j^+
=
-\frac1{b_\star}
\left(
U_{n,j}^{+,A}+2U_{n,j}^{+,B}
\right).
\eqtag{14.12}
\]
The sequence of diagonal coefficients from Section~\ref{A:sec:decomposition-of-the-general-jacobi-coefficients} satisfies
\[
\zeta_j
=
\zeta_j^++m_{j,n},
\eqtag{14.13}
\]
where
\[
m_{j,n}
=
-\frac1{b_\star}
\left[
(\overline A_j^\circ-A_j^\circ)
+2(\overline B_j^\circ-B_j^\circ)
\right].
\eqtag{14.14}
\]
By Proposition~\ref{A:res:14.1},
\[
\max_{j\in I_n^{\mathrm{mid}}}|m_{j,n}|
\le
C n^{-4}.
\eqtag{14.15}
\]

There are two obstructions to a potential on the whole space with bounded Hessian. The logarithm is undefined across the positive coordinate boundary, while the polynomial trace
has unbounded derivatives at infinity. We first continue $-\log$ convexly, using the deterministic lower bound on $\overline B_j^\circ$. We then apply the coordinatewise extension
of Lemma~\ref{A:res:A.11A}. Neither operation changes the potential near $\overline Y_n^\circ$.

\begin{lemma}\label{A:res:14.2}
For every \(a>0\), there exists \(\ell_a\in C^\infty(\mathbb R)\) such that
\[
\ell_a(x)=-\log x
\quad
(x\ge a),
\eqtag{14.16}
\] \(\ell_a''(x)>0 \allowbreak\; (x\in\mathbb R),\) and, for every \(m\ge2\),
\(\sup_{x\in\mathbb R} |\ell_a^{(m)}(x)| \le C_{a,m}\). Moreover, \(\ell_a(x)\ge c_a x^2-C_a \allowbreak\; (x\le0)\).
\end{lemma}

\begin{proof}

Choose $\chi\in C^\infty(\mathbb R;[0,1])$ equal to zero on $(-\infty,a/2]$ and one on $[a,\infty)$. Define
\[
 r_a(x)=\begin{cases}
 a^{-2},&x\le a/2,\\
 (1-\chi(x))a^{-2}+\chi(x)x^{-2},&x>a/2,
 \end{cases}
 \qquad
 \ell_a(x)=-\log a-\frac{x-a}{a}
                 +\int_a^x(x-s)r_a(s)\,\mathrm ds.
\]
Smoothness of $\chi$ at its constant regions gives $r_a\in C^\infty(\mathbb R)$, with
\[
 r_a>0,\qquad r_a(x)=x^{-2}\ (x\ge a),\qquad
 \sup_{x\in\mathbb R}|r_a^{(k)}(x)|<\infty\quad(k\ge0).
\]
Direct differentiation yields
\[
 \ell_a(a)=-\log a,\qquad\ell_a'(a)=-a^{-1},\qquad
 \ell_a''=r_a>0,\qquad\ell_a^{(m)}=r_a^{(m-2)}\quad(m\ge2).
\]
On $[a,\infty)$ the function $\ell_a+\log$ has zero second derivative and vanishing value and first derivative at $a$; it is therefore zero. This proves \eqref{A:eq:14.16} and all
derivative bounds. For $x\le0$, $r_a(x)=a^{-2}$, and
\[
 \begin{aligned}
 \ell_a(x)&=\ell_a(0)+\ell_a'(0)x+\frac{x^2}{2a^2}\\
 &\ge\frac{x^2}{4a^2}+\ell_a(0)-a^2|\ell_a'(0)|^2
 \ge\frac{x^2}{4a^2}-C_a.
 \end{aligned}
\]
The middle inequality is $|\ell_a'(0)x|\le x^2/(4a^2)+a^2|\ell_a'(0)|^2$.

\end{proof}

Use the lower bound $b_{\mathrm{reg}}$ from \eqref{A:eq:14.10}. In the exact conditional potential, replace every transported term \(-\log B_j\) by
\(\ell_{b_{\mathrm{reg}}}(B_j)\). Denote the preliminary potential on all of \(\mathbb R^{2J_n^{\mathrm{tr}}}\) by \(\Phi_{n,\log}^{\tau}.\) The polynomial trace is defined for
arbitrary real off-diagonal coordinates. Its Hessian lower bound holds for all real values of these coordinates, and the continued logarithms have nonnegative coefficients by
\eqref{A:eq:14.4}. Hence $\nabla^2\Phi_{n,\log}^\tau\ge cnI$ on the whole space. Near $\overline Y_n^\circ$ all logarithms are unchanged, so its gradient there is zero and this
remains the unique minimizer. Finite range is preserved. On each fixed coordinate cube about the minimizer, the coefficients and the relevant exterior values are bounded; the path
expansion and Lemma~\ref{A:res:14.2} therefore bound every local derivative by $C_mn$, uniformly in the interval length.

Let
\[
\widehat\Phi_{n,\log}^{\tau}(u)
=
\Phi_{n,\log}^{\tau}(\overline Y_n^\circ+u)
-
\Phi_{n,\log}^{\tau}(\overline Y_n^\circ).
\]
Fix a deterministic \(R_{\mathrm{ext}}>1\). Applying the explicit finite-range continuation of Lemma~\ref{A:res:A.11A} to \(\widehat\Phi_{n,\log}^{\tau}\), we obtain a globally
uniformly convex extension
\(\widehat\Phi_{n,\mathrm{ext}}^{\tau}\) satisfying
\[
\widehat\Phi_{n,\mathrm{ext}}^{\tau}(u)
=
\widehat\Phi_{n,\log}^{\tau}(u)
\quad
(\|u\|_{\ell^\infty}\le R_{\mathrm{ext}}),
\]
\[
cnI
\le
\nabla^2\widehat\Phi_{n,\mathrm{ext}}^{\tau}(u)
\le
C nI
\quad
(u\in\mathbb R^{2J_n^{\mathrm{tr}}}),
\]
and
\[
\bigl(\nabla^2\widehat\Phi_{n,\mathrm{ext}}^{\tau}(u)\bigr)_{ik}=0
\quad
(\operatorname{dist}(i,k)>R_1),
\]
for one fixed \(R_1\). Moreover, for every relevant \(m\ge2\), \(|\partial_{i_1}\cdots\partial_{i_m} \widehat\Phi_{n,\mathrm{ext}}^{\tau}(u)| \le C_mn\), and the derivative
vanishes unless the indices lie in one fixed graph neighborhood.

Return to uncentered variables by setting
\[
\Phi_{n,\mathrm{ext}}^{\tau}(y)
=
\widehat\Phi_{n,\mathrm{ext}}^{\tau}(y-\overline Y_n^\circ)
+
\Phi_{n,\log}^{\tau}(\overline Y_n^\circ).
\]
Define
\[
\mathcal S_n
=\bigcap_{j=1}^{J_n^{\mathrm{tr}}}\{B_j\ge b_{\mathrm{reg}}\},
\]
and the coordinate neighborhood
\[
\mathcal Q_n
=
\left\{
\|y-\overline Y_n^\circ\|_{\ell^\infty}
\le R_{\mathrm{ext}}
\right\}.
\]
By construction, on \(\mathcal S_n\cap\mathcal Q_n\),
\[
\Phi_{n,\mathrm{ext}}^{\tau}
=
\Phi_{n,+}^{\tau}.
\eqtag{14.24}
\]
Since \(\overline B_j^\circ\ge4b_{\mathrm{reg}}\) and the extension agrees with the original potential in a neighborhood of the minimizer, \(\overline Y_n^\circ\) is the unique
minimizer of the final extended potential, and
\[
K_{n,\mathrm{tr}}^{\tau}
:=
\nabla^2\Phi_{n,+}^{\tau}(\overline Y_n^\circ)
=
\nabla^2\Phi_{n,\mathrm{ext}}^{\tau}(\overline Y_n^\circ).
\]

Let \(\widetilde\mu_n^\tau\) be the centered globally extended law on the full Euclidean space:
\[
\mathrm{d}\widetilde\mu_n^\tau(u)
=
\widetilde Z_n^{-1}
\exp\left\{
-\widehat\Phi_{n,\mathrm{ext}}^{\tau}(u)
\right\}\,\mathrm{d}u.
\]
Write $\mu_{n,+}^{\tau,\circ}$ for the law of $U_n^+$ defined above. In centered coordinates, the two potentials agree on
\[
\mathcal C_n^\circ
=
\left\{
\min_{1\le j\le J_n^{\mathrm{tr}}}(\overline B_j^\circ+u_j^B)\ge b_{\mathrm{reg}},
\quad
\|u\|_{\ell^\infty}\le R_{\mathrm{ext}}
\right\}.
\eqtag{14.26c}
\]
Subtract the common value of the potentials at $\overline Y_n^\circ$ when writing their centered densities. Their unnormalized densities then agree on $\mathcal C_n^\circ$.
Normalizing their restrictions to this set gives the same conditional law, even when their full partition functions differ.

\begin{proposition}\label{A:res:14.3}
For every \(p<\infty\), conditionally on the exterior data, there exists a coupling of
\(U_n^+\sim\mu_{n,+}^{\tau,\circ}\) and \(\widetilde U_n\sim\widetilde\mu_n^\tau\) such that
\[
\mathbb P_\tau(U_n^+\ne\widetilde U_n)
\le C \mathrm{e}^{-cn},
\eqtag{14.27}
\]
and
\[
\sup_i
\|U_{n,i}^+-\widetilde U_{n,i}\|_{L^p(\mathbb P_\tau)}
\le C_p\mathrm{e}^{-c_pn}.
\eqtag{14.28}
\]
\end{proposition}

\begin{proof}
The interval $I_n^{\mathrm{tr}}$ lies in a fixed regular region. Lemma~\ref{A:res:mean} for the original conditional law and Proposition~\ref{A:res:A.11B} for the extended law
give, for every finite $q$,
\[
 \|\mathbb E_\tau U_n^+\|_\infty
       +\|\mathbb E_\tau\widetilde U_n\|_\infty\le Cn^{-1},
 \qquad
 \sup_i\bigl(\|U_{n,i}^+\|_{L^q(\mathbb P_\tau)}
       +\|\widetilde U_{n,i}\|_{L^q(\mathbb P_\tau)}\bigr)
 \le C_qn^{-1/2}.
\]
Both potentials have Hessian at least $cnI$ on their convex domains. For large $n$, leaving $\mathcal C_n^\circ$ requires a coordinate deviation from its mean of at least
$\min\{2b_{\mathrm{reg}},R_{\mathrm{ext}}/2\}$. Lemma~\ref{A:res:A.11C} and a union bound over at most $Cn$ coordinates yield
\[
 \mu_{n,+}^{\tau,\circ}((\mathcal C_n^\circ)^c)
 +\widetilde\mu_n^\tau((\mathcal C_n^\circ)^c)
 \le Cn\exp\{-cn\min(b_{\mathrm{reg}}^2,R_{\mathrm{ext}}^2)\}
 \le Ce^{-cn}.
\]
By \eqref{A:eq:14.24}, the unnormalized densities agree on $\mathcal C_n^\circ$, and hence the two laws have the same conditional law there. Lemma~\ref{A:res:A.11D} supplies a
coupling kernel that retains $\widetilde U_n$ and satisfies
\[
 \mathbb P_\tau(U_n^+\ne\widetilde U_n)\le Ce^{-cn}.
\]
The normalized densities and the indicator of $\mathcal C_n^\circ(\tau)$ are jointly Borel in the exterior parameter and coefficient vector. Their integrals, restricted masses, and
common conditional law are therefore Borel. The explicit kernel in Lemma~\ref{A:res:A.11D}, with the identity kernel when its common mass is one, is consequently a Borel
probability kernel.

The same kernel works for every $p$. Apply the moment bound in that lemma with $q=2p$ to obtain
\[
 \sup_i\|U_{n,i}^+-\widetilde U_{n,i}\|_{L^p(\mathbb P_\tau)}
 \le C_pn^{-1/2}
       \mathbb P_\tau(U_n^+\ne\widetilde U_n)^{1/(2p)}
 \le C_pe^{-c_pn}.
\]
These are \eqref{A:eq:14.27} and \eqref{A:eq:14.28}.
\end{proof}

Let

\[
G_n^{\mathrm{tr}}
=
(g_1^A,\ldots,g_{J_n^{\mathrm{tr}}}^A,
 g_1^B,\ldots,g_{J_n^{\mathrm{tr}}}^B)
\sim
N\left(0,(K_{n,\mathrm{tr}}^{\tau})^{-1}\right).
\]
We set \(g_0^B=0\). The centered Gaussian diagonal coefficients used on the interior block are
\(\widehat\zeta_j^{\mathrm G} = -\frac1{b_\star} \left(g_j^A+2g_j^B\right), \allowbreak\, j\in I_n^{\mathrm{mid}}\).
Define the diagonal coefficient sequence
\[
\zeta_j^{\mathrm G}
=
\begin{cases}
\zeta_j,&j\in I_n^{\mathrm L},\\
\widehat\zeta_j^{\mathrm G},&j\in I_n^{\mathrm{mid}},\\
\zeta_j,&j\in I_n^{\mathrm{far}},
\end{cases}
\]
and set
\[
 H_n^{\mathrm G}=H_n^\circ+h^{-2}M_{\zeta_n^{\mathrm G}}.
 \eqtag{gaussian-comparison-operator}
\]
In $H_n^{\mathrm G}$, the scalar coefficients in the middle interval have the Gaussian distribution constructed above. The kinetic coefficients and the scalar coefficients outside
that interval are unchanged.

Put
\[
\Psi_{n,\mathrm{ext}}^{\tau}(u)
=
\widehat\Phi_{n,\mathrm{ext}}^{\tau}(u)
-
\frac12\langle u,K_{n,\mathrm{tr}}^{\tau}u\rangle.
\]
For \(0\le t\le1\), define
\[
\Phi_t(u)
=
\frac12\langle u,K_{n,\mathrm{tr}}^{\tau}u\rangle
+t\Psi_{n,\mathrm{ext}}^{\tau}(u)
\eqtag{14.34b}
\]
and \(\mathrm{d}\mu_t(u)=Z_t^{-1}\mathrm{e}^{-\Phi_t(u)}\,\mathrm{d}u.\) Then
\(\mu_0=N(0,(K_{n,\mathrm{tr}}^{\tau})^{-1}), \allowbreak\, \mu_1=\widetilde\mu_n^\tau\). The Hessian bounds for the
extension and its quadratic approximation give
\[
cnI
\le
\nabla^2\Phi_t(u)
\le
C nI
\quad
(0\le t\le1),
\eqtag{14.34d}
\]
and
\[
(\nabla^2\Phi_t(u))_{ik}=0
\quad
(\operatorname{dist}(i,k)>R_1).
\eqtag{14.34e}
\]
All local derivatives of order at least two are bounded by \(C_mn\), uniformly in \(t,u,n\).

\begin{lemma}\label{A:res:14.4}
For every \(p<\infty\), \(\sup_{0\le t\le1}\sup_i \|u_i\|_{L^p(\mu_t)} \le C_pn^{-1/2}\).
\end{lemma}

\begin{proof}
Write $K=K_{n,\mathrm{tr}}^\tau$ and $\widehat\Phi=\widehat\Phi_{n,\mathrm{ext}}^\tau$. Their interpolation satisfies
\[
 \nabla\Phi_t(0)=0,\qquad
 \nabla^2\Phi_t(u)=(1-t)K+t\nabla^2\widehat\Phi(u),\qquad
 cnI\le\nabla^2\Phi_t(u)\le CnI.
\]
Both summands have range at most $R_1$. For derivatives of order at least three, the quadratic part vanishes and
$\partial_{i_1}\cdots\partial_{i_m}\Phi_t =t\partial_{i_1}\cdots\partial_{i_m}\widehat\Phi$. Thus all hypotheses of Proposition~\ref{A:res:A.11B} hold with constants independent of
$t\in[0,1]$. With $m_{t,i}=\int u_i\,\mathrm d\mu_t$, that proposition gives
\[
 |m_{t,i}|\le Cn^{-1},\qquad
 \|u_i-m_{t,i}\|_{L^p(\mu_t)}\le C_pn^{-1/2}.
\]
The triangle inequality yields $\|u_i\|_{L^p(\mu_t)}\le C_pn^{-1/2}+Cn^{-1}\le C_p'n^{-1/2}$, uniformly in $i,t$ and the admissible exterior data.
\end{proof}

\begin{lemma}\label{A:res:14.5}
For every \(p<\infty\), \(\sup_{0\le t\le1}\sup_i \|\partial_i\Psi_{n,\mathrm{ext}}^{\tau}\|_{L^p(\mu_t)} \le C_p\).
\end{lemma}

\begin{proof}
Write $\Psi=\Psi_{n,\mathrm{ext}}^\tau$. The extension preserves the Hessian at its minimizer, so $\nabla\Psi(0)=0$ and $\nabla^2\Psi(0)=0$. In the centered coordinate $u$,
Taylor's formula is
\[
 \partial_i\Psi(u)=\sum_{j,k}\int_0^1(1-r)
   \partial_{ijk}^3\widehat\Phi_{n,\mathrm{ext}}^\tau(ru)
   u_ju_k\,\mathrm dr.
\]
In the uncentered potential, the same derivative is evaluated at $\overline Y_n^\circ+ru$. Each nonzero third derivative is bounded by $Cn$, and only a bounded number of pairs
$(j,k)$ occur for each $i$. Lemma~\ref{A:res:14.4} therefore gives
\[
 \|\partial_i\Psi\|_{L^p(\mu_t)}
 \le Cn\sum_{j,k:\,\mathrm d(i,j),\mathrm d(i,k)\le R_1}
        \|u_j\|_{L^{2p}(\mu_t)}\|u_k\|_{L^{2p}(\mu_t)}
 \le C_p.
\]
This verifies the local remainder bound of Lemma~\ref{A:res:A.12} for the conditional coefficients.
\end{proof}

Let \(L_t=\Delta-\nabla\Phi_t\cdot\nabla\) be the Langevin generator reversible with respect to \(\mu_t\). Lemma~\ref{A:res:A.13} applies because
\eqref{A:eq:14.34d}--\eqref{A:eq:14.34e}, Lemma~\ref{A:res:14.4}, and Lemma~\ref{A:res:14.5} verify its hypotheses on uniform convexity, finite range, and the nonquadratic
remainder. Thus, if \(f\) is centered and \(L_t\varphi=f\), then
\[
\|\partial_i\varphi\|_{L^p(\mu_t)}
\le
\frac{C}{n}
\sum_j \mathrm{e}^{-a\operatorname{dist}(i,j)}
\|\partial_jf\|_{L^p(\mu_t)}.
\]
Applying this estimate to \(f=\Psi_{n,\mathrm{ext}}^{\tau} -\mathbb E_{\mu_t}\Psi_{n,\mathrm{ext}}^{\tau}\), Theorem~\ref{A:res:A.14} constructs a transport from
$N(0,(K_{n,\mathrm{tr}}^\tau)^{-1})$ to $\widetilde\mu_n^\tau$ satisfying
\[
\sup_i
\|\widetilde U_{n,i}-G_{n,i}^{\mathrm{tr}}\|_{L^p(\mathbb P_\tau)}
\le
C_pn^{-1}.
\eqtag{14.34i}
\]

Apply the transport map, followed by the kernel of Lemma~\ref{A:res:A.11D}, to obtain the following coordinate estimate under the exact conditional law. The covariance is the
inverse of the conditional Hessian, as in \cite{krishnapur2016universality}, where averages of Gaussian probability measures are compared in total variation on intervals of length
at most $n^{1/2-\epsilon}$. Here the transported interval has order $n$ sites, and the conclusion is an $L^p$ coupling bound for each coordinate. This does not give convergence in
total variation of the full transported vector to its Gaussian approximation.

\par\Needspace{7\baselineskip}
\begin{proposition}\label{A:res:14.6}
For each admissible exterior value $\tau$, there is a coupling of $U_n^+\sim\mu_{n,+}^{\tau,\circ}$ and $G_n^{\mathrm{tr}}\sim N(0,(K_{n,\mathrm{tr}}^\tau)^{-1})$ such that, for
every $p<\infty$,
\[
\sup_i
\|U_{n,i}^+-G_{n,i}^{\mathrm{tr}}\|_{L^p(\mathbb P_\tau)}
\le
C_pn^{-1}.
\eqtag{14.35}
\]
\end{proposition}

\begin{proof}
Write $\gamma_\tau=N(0,(K_{n,\mathrm{tr}}^\tau)^{-1})$, and let $T_\tau$ be the transport in \eqref{A:eq:14.34i}. Let $\kappa_\tau(\widetilde u,\mathrm du)$ be the kernel of
Proposition~\ref{A:res:14.3}, conditional on the extended centered variable. Construct $(G_n^{\mathrm{tr}},U_n^+)$ with joint law
\[
 \gamma_\tau(\mathrm dg)\,
 \kappa_\tau(T_\tau(g),\mathrm du),
 \qquad \widetilde U_n=T_\tau(G_n^{\mathrm{tr}}).
\]
Integration in $g$ shows that $\widetilde U_n$ has the extended law. Integration of $\kappa_\tau$ against that law gives the exact centered conditional law of $U_n^+$. The Gaussian
marginal is unchanged. For each coordinate, Proposition~\ref{A:res:14.3} and \eqref{A:eq:14.34i} give
\[
 \|U_{n,i}^+-G_{n,i}^{\mathrm{tr}}\|_{L^p(\mathbb P_\tau)}
 \le\|U_{n,i}^+-\widetilde U_{n,i}\|_{L^p(\mathbb P_\tau)}
        +\|\widetilde U_{n,i}-G_{n,i}^{\mathrm{tr}}\|_{L^p(\mathbb P_\tau)}
 \le C_pe^{-c_pn}+C_pn^{-1}.
\]
Since $e^{-c_pn}\le C_pn^{-1}$, this proves \eqref{A:eq:14.35}. The map $T_\tau$ and the kernel $\kappa_\tau$ do not depend on $p$. The kernel acts on the given $G_n^{\mathrm{tr}}$
with its stated conditional Gaussian distribution, preserving its joint distribution with the other variables.

This preservation also holds when $G_n^{\mathrm{tr}}$ has already been coupled to a random element $Y$ used in the coupling with the Brownian path. Let
$\lambda_\tau(\mathrm dy,\mathrm dg)$ be their joint conditional law. Define the extension by
\[
 \lambda_\tau(\mathrm dy,\mathrm dg)
       \kappa_\tau(T_\tau(g),\mathrm du).
\]
Integration in $u$ leaves $\lambda_\tau$ unchanged. Integration in $y$ gives $\gamma_\tau(\mathrm dg)\kappa_\tau(T_\tau(g),\mathrm du)$, so $u$ has the exact centered conditional
law proved above. In particular, the $L^p$ estimate for $u-g$ holds under this extension, even though $u$ and $Y$ need not be independent given $\tau$.

\end{proof}

For exterior data outside \(\mathcal T_n^{\mathrm{out}}(M_0)\), complete this coefficient coupling by taking the exact centered conditional vector and \(G_n^{\mathrm{tr}}\)
independently. The conditional minimizer and its positive Hessian are defined for every exterior coefficient configuration. The minimizer and transport are Borel by
Lemma~\ref{A:res:B.7B}; the explicit extension uses fixed cutoff functions and uniform constants on the bounded exterior set. Parameter integration makes the density integrals
Borel, and Lemma~\ref{A:res:B.7D} gives a Borel kernel with the original conditional law as marginal. Use this kernel on $\mathcal T_n^{\mathrm{out}}(M_0)$ and the product of the
two marginals on its measurable complement. Lemma~\ref{A:res:A.8}, in dimension at most \(2n\), and the Gaussian covariance bound \((K_{n,\mathrm{tr}}^\tau)^{-1}\le Cn^{-1}I\) give
\[
 \sup_{\tau}\sup_i
 \|U_{n,i}^+-G_{n,i}^{\mathrm{tr}}\|_{L^{2p}(\mathbb P_\tau)}
 \le C_p
\]
for any coupling of the two marginals. Together with \eqref{A:eq:integrated-moment}, this gives the unconditional
\(O(n^{-1})\) coordinate bound. To check the full coefficient law, let $\nu_n^{\mathrm{out}}$ be the actual exterior
marginal. For measurable sets $A$ of exterior values and $D$ of transported coefficients, the construction gives
\[
 \mathbb P\{\tau\in A,\ Y_n^{\mathrm{tr}}\in D\}
   =\int_A\mu_{n,+}^\tau(D)\,\nu_n^{\mathrm{out}}(\mathrm d\tau).
\]
This is the original disintegration on both the specified exterior set and its complement. It determines the joint distribution of the transported coefficients and the original
scalar coefficients outside the middle interval.

For \(j\in I_n^{\mathrm{mid}}\), equations \eqref{A:eq:14.12}--\eqref{A:eq:14.15} and Proposition~\ref{A:res:14.6} give
\[
\zeta_j-\widehat\zeta_j^{\mathrm G} ={} -\frac1{b_\star} \left[ (U_{n,j}^{+,A}-g_j^A) +2(U_{n,j}^{+,B}-g_j^B) \right] +m_{j,n}.
\]
For every \(p<\infty\),
\[
\sup_{j\in I_n^{\mathrm{mid}}}
\|\zeta_j-\widehat\zeta_j^{\mathrm G}\|_{L^p(\mathbb P_\tau)}
\le
C_pn^{-1}
=
C_ph^3.
\]
Define
\[
\mathcal T_n^{\mathrm{mid}}
=
\sum_{j\in I_n^{\mathrm{mid}}}
(1+x_j)^{-7/2}
|\zeta_j-\widehat\zeta_j^{\mathrm G}|^2.
\]
Under $\mathbb P_\tau$, Minkowski's inequality gives, for every $p<\infty$,
\[
\|\mathcal T_n^{\mathrm{mid}}\|_{L^p} \le \sum_{j\in I_n^{\mathrm{mid}}} (1+x_j)^{-7/2} \|\zeta_j-\widehat\zeta_j^{\mathrm G}\|_{L^{2p}}^2 \le C_ph^6 \sum_{j\ge1}(1+x_j)^{-7/2} \le C_ph^5.
\]
The same order holds after integrating the exterior data. Indeed, the polynomial minimizer bound \eqref{A:eq:mode-growth}, the coefficient moments of Lemma~\ref{A:res:raw}, and
Lemma~\ref{A:res:A.8} give uniformly bounded unconditional $4p$-moments of each scalar difference before using its smallness. H\"older's inequality on the exceptional exterior set
then gives
\[
 \bigl\|\mathbf1_{\{\tau\notin\mathcal T_n^{\mathrm{out}}(M_0)\}}
             \mathcal T_n^{\mathrm{mid}}\bigr\|_{L^p}
 \le C_ph^{-1}e^{-c_pn}.
\]
Thus $\|\mathcal T_n^{\mathrm{mid}}\|_{L^p}\le C_ph^5$ under the unconditional coupling. On bounded exterior data, Markov at threshold $h^4$ gives conditional probability at most
$C_ph^p$; integrating adds at most $Ce^{-cn}$. Since $h\asymp n^{-1/3}$, $\sum_n(h_n^p+e^{-cn})<\infty$ for $p>3$. Borel--Cantelli gives
\[
\mathcal T_n^{\mathrm{mid}}
\le
C_\omega h^4
\eqtag{14.39}
\]
almost surely for all sufficiently large \(n\).

For the form comparison we also need a bound on the largest coordinate error.

Put \(t_{j,n}=\zeta_j-\widehat\zeta_j^{\mathrm G}\) on \(I_n^{\mathrm{mid}}\).

\begin{lemma}\label{A:res:transport-max} For every fixed \(0<\varepsilon<1/3\)
and every finite \(p\ge1\), on the coupling constructed from the common Brownian path,
\[
 \mathbb P\!\left(\max_{j\in I_n^{\mathrm{mid}}}|t_{j,n}|
                  >n^{-1+\varepsilon}\right)
 \le C_{p,\varepsilon}n^{1-p\varepsilon}+C\mathrm{e}^{-cn}.
\eqtag{transport-max-tail}
\]
\end{lemma}

In particular, with \(\varepsilon=1/6\), almost surely for all sufficiently large \(n\),
\[
 \max_{j\in I_n^{\mathrm{mid}}}|t_{j,n}|\le n^{-5/6},
 \quad
 h^{-2}\max_{j\in I_n^{\mathrm{mid}}}|t_{j,n}|
       \le C_V n^{-1/6}=O(h^{1/2})=o(1).
\eqtag{transport-max}
\]
For each $Q>0$, the event in \eqref{A:eq:transport-max} has complementary probability at most \(C_Qn^{-Q}\), after choosing a single sufficiently large fixed \(p\).

\begin{proof}
For admissible exterior data, the coordinate estimate preceding \eqref{A:eq:14.39}, including the change of conditional minimizer, gives
\(\mathbb E_\tau|t_{j,n}|^p\le C_p n^{-p}\). There are at most \(n\) middle coordinates. Markov's inequality and a union bound therefore give
\[
 \mathbb P_\tau\!\left(\max_j|t_{j,n}|>n^{-1+\varepsilon}\right)
 \le\sum_j n^{p(1-\varepsilon)}\mathbb E_\tau|t_{j,n}|^p
 \le C_p n^{1-p\varepsilon}.
\]
Integrate on \(\mathcal T_n^{\mathrm{out}}(M_0)\) and add its complementary probability using Lemma~\ref{A:res:disintegration}. This proves \eqref{A:eq:transport-max-tail} for the
extension obtained by taking the product of the marginals as well. Choose \(p\varepsilon>2\) and apply Borel--Cantelli. Since \(h^{-2}\asymp n^{2/3}\), the scaled maximum is
\(O(n^{-1/3+\varepsilon})\); taking \(\varepsilon=1/6\) proves \eqref{A:eq:transport-max}. Taking
\(p/6>Q+1\) gives the required probability bound, with deterministic constants uniform on compact positive \(\beta\)-intervals. The
conditional sampling in Section~\ref{A:sec:application-to-the-conditional-jacobi-coefficients} preserves the conditional joint coefficient law, so these probability estimates also
hold on that final coupling.
\end{proof}

The difference of the two operators is
\[
H_n^{\mathrm G}-H_n^{\mathrm{diag}}
=
h^{-2}M_{\chi^{\mathrm{mid}}
(\widehat\zeta^{\mathrm G}-\zeta)}.
\eqtag{14.40}
\]
By \eqref{A:eq:transport-max}, this multiplier has operator norm tending to zero, so a common form bound for one of the two operators transfers to the other. The starting form
bound is proved in Theorem~\ref{A:res:16.2} from the common Brownian path. For the resolvent error we use the sharper weighted estimate \eqref{A:eq:14.39}.

Theorem~\ref{A:res:14.7} uses the common form bounds established in Theorem~\ref{A:res:16.2} from the coefficient estimates, before any of the resolvent conclusions are used.

The same diagonal estimate gives a moment bound on events with bounded resolvent constants. Let $\mathcal A_{n,L}$ be an event on which the form constants and coefficient bounds
for the pair are bounded by a fixed $L$. With $R_n^{\mathrm{diag}}(z)=(H_n^{\mathrm{diag}}-z)^{-1}$ and $R_n^{\mathrm G}(z)=(H_n^{\mathrm G}-z)^{-1}$, the diagonal estimate and the
unconditional moment bound just proved give
\[
 \left\|\mathbf1_{\mathcal A_{n,L}}\sup_{z\in K}
   \|R_n^{\mathrm{diag}}(z)-R_n^{\mathrm G}(z)\|_{\mathfrak S_2}^2
       \right\|_{L^p}
 \le C_{L,K}h^{-3}\|\mathcal T_n^{\mathrm{mid}}\|_{L^p}
 \le C_{p,L,V,\beta,K}h^2.
\]
This is an estimate of the error multiplied by an event indicator under the original coupling law. Subsection~\ref{A:sec:uniform-moment-estimates} integrates the random form
constants to remove this restriction.

Let
\[
\mathcal L_{\mathrm{diag}}^{\mathrm{mid}}:
\mathbb R^{2J_n^{\mathrm{tr}}}
\longrightarrow
\mathbb R^{I_n^{\mathrm{mid}}}
\]
be the local scalar map \((\mathcal L_{\mathrm{diag}}^{\mathrm{mid}}u)_j = -\frac1{b_\star}(u_j^A+2u_j^B), \allowbreak\, j\in I_n^{\mathrm{mid}}\). The covariance used in
Section~\ref{A:sec:application-to-the-conditional-jacobi-coefficients} is
\[
C_{n,\mathrm{mid}}^{\tau,+}
=
\mathcal L_{\mathrm{diag}}^{\mathrm{mid}}
(K_{n,\mathrm{tr}}^\tau)^{-1}
(\mathcal L_{\mathrm{diag}}^{\mathrm{mid}})^*.
\eqtag{14.45}
\]
On $\mathcal T_n^{\mathrm{out}}(M_0)$, the conditional minimizer lies in a fixed bounded set on which the off-diagonal coordinates are bounded away from zero. Hence
$cnI\le K_{n,\mathrm{tr}}^\tau\le CnI$, with deterministic constants, and its bandwidth is fixed. The banded inverse estimate and the two entries in each row of
$\mathcal L_{\mathrm{diag}}^{\mathrm{mid}}$ give
\[
 |C_{n,\mathrm{mid}}^{\tau,+}(i,j)|
 \le Cn^{-1}e^{-c|i-j|}\le C'h^3e^{-c|i-j|}.
 \eqtag{14.46}
\]
Here $C,C',c$ are deterministic and uniform over the stated exterior set. As the covariance of the specified Gaussian linear image, this matrix is nonnegative. By
\eqref{A:eq:L10}--\eqref{A:eq:L11}, both boundaries of the finite conditional Hessian are at logarithmic graph distance from every diagonal coefficient used in \eqref{A:eq:14.45}.
These distances enter the two boundary terms in Proposition~\ref{A:res:A.15}.

We summarize the coefficient estimate and its resolvent consequence. The operator $H_n^{\mathrm G}$ is the one defined in \eqref{A:eq:gaussian-comparison-operator}; the original
coefficients on $I_n^{\mathrm L}\cup I_n^{\mathrm{far}}$ have not been changed.

\begin{theorem}\label{A:res:14.7}
The coefficient coupling above satisfies, almost surely for all sufficiently large $n$ and every compact
$K\Subset\mathbb C\setminus\mathbb R$,
\[
\sum_{j\in I_n^{\mathrm{mid}}}
|\zeta_j-\widehat\zeta_j^{\mathrm G}|^2
(1+x_j)^{-7/2}
\le
C_\omega h^4,
\eqtag{14.47}
\]
\[
\sup_{z\in K}
\left\|
(H_n^{\mathrm{diag}}-z)^{-1}
-
(H_n^{\mathrm G}-z)^{-1}
\right\|_{\mathfrak S_2}^2
\le
C_{\omega,V,\beta,K}h.
\eqtag{14.49}
\]
The middle Gaussian covariance is \eqref{A:eq:14.45}, with the exponential estimate \eqref{A:eq:14.46}.
\end{theorem}

Theorems~\ref{A:res:14.7}, \ref{A:res:15.4}, \ref{A:res:15.5}, and \ref{A:res:15.6} are proved together after Theorem~\ref{A:res:16.2}, which supplies the common form bounds for
their resolvent comparisons.

\subsubsection{Inverse Hessians and comparison of Gaussian vectors}\label{A:sec:localized-inverse-hessians-and-weighted-gaussian-couplings}
The preceding transport leaves a Gaussian vector whose covariance is the inverse of a finite conditional Hessian. We compare this inverse with those of the constant coefficient
Hessians. The weights $(1+x_j)^{-7/2}$ in the expected squared error require estimates that retain decay with index distance. The transport estimate relates this Gaussian
distribution to the conditional Gibbs law.

We form an approximate inverse from the inverses of constant coefficient Hessians and estimate the errors due to the varying parameter and the two boundaries. To compare the
Gaussian vectors, we also need a positive approximating covariance. The factorization in Subsection~\ref{A:sec:spectral-factorization-and-coupling-to-brownian-increments}
constructs this covariance.

Let $\Lambda_n$ be the coefficient vertices over the sites $1,\ldots,N_n$, with a fixed finite number of types at each site. Write $p=(j,\sigma)$, $\nu(p)=j$, and
$\sigma(p)=\sigma$. In the application to Jacobi matrices $a_j$ and $b_j$ have the same integer index $j$. The graph distance $d(p,q)$ is comparable to $|\nu(p)-\nu(q)|$, up to
constants depending only on this strip. Put
\[
 d_{\mathrm L}(p)=d(p,\partial_{\mathrm L}\Lambda_n),\qquad
 d_{\mathrm R}(p)=d(p,\partial_{\mathrm R}\Lambda_n).
\]
Let $\mathscr K_n$ be real symmetric, with $cI\le\mathscr K_n\le CI$ and $\mathscr K_n(p,q)=0$ for $d(p,q)>R_0$. Let $\mathscr K_s$, $0\le s\le s_0$, be translation-invariant
finite-range operators on the full strip, with Hermitian symbols $\mathsf K_s(\theta)$. Assume $cI\le\mathsf K_s(\theta)\le CI$, uniform analyticity in a fixed complex strip, and
uniform $C^1$ dependence on $s$. All range and spectral constants are independent of $n$.

Take $s_j\in[0,s_0]$ with $|s_{j+1}-s_j|\le Ch^3$, and use the integer center convention $s_{pq}=s_{\lfloor(\nu(p)+\nu(q))/2\rfloor}$. Suppose
\[
 |\mathscr K_n(p,q)-\mathscr K_{s_{pq}}(p-q)|
 \le C\bigl[h^3+\mathrm e^{-cd_{\mathrm L}(p)}
                   +\mathrm e^{-cd_{\mathrm R}(p)}\bigr]
                                        \mathbf1_{\{d(p,q)\le R_0\}}.
 \eqtag{A.49}
\]
The two exponential terms allow errors of order one in the boundary rows. For the Hessian of the conditional potential for the Jacobi coefficients, these include the displacement
of its minimizer from the local profile. They cannot be discarded before taking the inverse, since inversion couples the boundary to every interior site.

Here $\mathscr K_s(p-q)$ means the $(\sigma(p),\sigma(q))$ entry of the matrix coefficient at lag $\nu(p)-\nu(q)$. Our Fourier convention is
$\mathsf K_s(\theta)=\sum_{r\in\mathbb Z}\mathscr K_s(r)e^{-\mathrm i r\theta}$; thus the inverse coefficients have the opposite sign in the exponent. Define
\[
 q_s(p-q)=\frac1{2\pi}\int_{-\pi}^{\pi}
   \mathrm e^{\mathrm i(\nu(p)-\nu(q))\theta}
   \mathsf K_s(\theta)^{-1}_{\sigma(p)\sigma(q)}\,\mathrm d\theta,
 \qquad Q_n(p,q)=q_{s_{pq}}(p-q).
\]
The center convention makes $Q_n$ Hermitian; it is real symmetric for the real local Hessians used below. The eigenvalues of $Q_n$ may have either sign.

\begin{proposition}\label{A:res:A.15}
There are \(c_1,C_1>0\) such that for all \(p,q\in\Lambda_n\),
\[
|\mathscr K_n^{-1}(p,q)-Q_n(p,q)| \le{} C_1\mathrm{e}^{-c_1d(p,q)} \Bigl[ h^3 +\mathrm{e}^{-c_1\min\{d_{\mathrm L}(p),d_{\mathrm L}(q)\}} +\mathrm{e}^{-c_1\min\{d_{\mathrm R}(p),d_{\mathrm R}(q)\}} \Bigr].
\eqtag{C1}
\]
\end{proposition}

In particular, if \(S_n\) is separated from the two boundaries by distances \(L_n^{\mathrm L}\) and \(L_n^{\mathrm R}\), then for \(p,q\in S_n\),
\[
|\mathscr K_n^{-1}(p,q)-Q_n(p,q)|
\le
C_1\left(
 h^3+\mathrm{e}^{-c_1L_n^{\mathrm L}}+\mathrm{e}^{-c_1L_n^{\mathrm R}}
\right)\mathrm{e}^{-c_1d(p,q)}.
\eqtag{C2}
\]
If \(K_n=\frac{\beta n}{2}\mathscr K_n, \allowbreak\, n^{-1}\asymp h^3\), and both exponential boundary terms are at most \(h^3\), then
\[
\left|
K_n^{-1}(p,q)-\frac{2}{\beta n}Q_n(p,q)
\right|
\le Ch^6\mathrm{e}^{-cd(p,q)},
\quad p,q\in S_n.
\eqtag{C3}
\]

\begin{proof}
We first prove the estimates for the Fourier coefficients. Uniform analyticity gives, for real $\theta$ and sufficiently small $|\eta|$,
\[
\|\mathsf K_s(\theta+\mathrm i\eta)-\mathsf K_s(\theta)\|
\le C|\eta|.
\]
Choose $\eta_0>0$ so that the right side is at most $c/2$ when $|\eta|\le\eta_0$. Since $\|\mathsf K_s(\theta)^{-1}\|\le c^{-1}$, the inverse expansion gives
\[
\sup_{s,\theta,|\eta|\le\eta_0}
\|\mathsf K_s(\theta+\mathrm i\eta)^{-1}\|\le 2c^{-1}.
\]
Differentiating the identity $\mathsf K_s\mathsf K_s^{-1}=I$ yields
\[
\partial_s\mathsf K_s^{-1}
=-\mathsf K_s^{-1}(\partial_s\mathsf K_s)\mathsf K_s^{-1}.
\]
The finite range and the uniform bounds on the derivatives of the coefficients bound this derivative throughout the same smaller strip. In the integral defining $q_s(p-q)$, move
the contour from the real axis to $\operatorname{Im}\theta=\eta_0\operatorname{sgn}(\nu(p)-\nu(q))$. Periodicity cancels the two vertical sides. Put $\ell=\nu(p)-\nu(q)$ and
$\eta=\eta_0\operatorname{sgn}\ell$. For either the inverse symbol $F_s=\mathsf K_s^{-1}$ or its derivative $F_s=-\mathsf K_s^{-1}(\partial_s\mathsf K_s)\mathsf K_s^{-1}$,
\[
 \left\|\frac1{2\pi}\int_{-\pi}^{\pi}
     \mathrm e^{\mathrm i\ell(\theta+\mathrm i\eta)}
                 F_s(\theta+\mathrm i\eta)\,\mathrm d\theta\right\|
 \le \mathrm e^{-\eta_0|\ell|}
       \sup_{\theta}\|F_s(\theta+\mathrm i\eta)\|
 \le C\mathrm e^{-\eta_0|\ell|}.
\]
The bounded number of coordinate types at each site then gives
\[
|q_s(p-q)|+|\partial_sq_s(p-q)|
\le C\mathrm{e}^{-cd(p,q)},
\quad
|\partial_s\mathscr K_s(p-q)|
\le C\mathbf 1_{\{d(p,q)\le R_0\}}.
\eqtag{C4}
\]
Set $\mathcal R_n=\mathscr K_nQ_n-I$. To estimate $\mathcal R_n$, evaluate both factors at $s_{pq}$ and subtract the sum over $r\notin\Lambda_n$ from the sum over
$r\in\mathfrak G_\infty$. At a fixed parameter their product is exactly the identity: Fourier inversion gives
\[
\sum_{r\in\mathfrak G_\infty}
\mathscr K_{s_{pq}}(p-r)q_{s_{pq}}(r-q)=\delta_{pq}.
\eqtag{C5}
\]
For $s=s_{pq}$, the entries of $\mathcal R_n=\mathscr K_nQ_n-I$ satisfy
\[
\begin{aligned}
\mathcal R_n(p,q)={}&\sum_{r\in\Lambda_n}
 [\mathscr K_n(p,r)-\mathscr K_{s_{pr}}(p-r)]q_{s_{rq}}(r-q)\\
&+\sum_{r\in\Lambda_n}
 [\mathscr K_{s_{pr}}(p-r)-\mathscr K_s(p-r)]q_{s_{rq}}(r-q)\\
&+\sum_{r\in\Lambda_n}\mathscr K_s(p-r)
 [q_{s_{rq}}(r-q)-q_s(r-q)]
 -\sum_{r\notin\Lambda_n}\mathscr K_s(p-r)q_s(r-q).
\end{aligned}
\]
The four sums account, respectively, for the finite Hessian, the parameter in its coefficients, the parameter in the inverse coefficients, and the missing indices. Each contains
only a bounded number of indices, all with $d(p,r)\le R_0$. For the first sum, \eqref{A:eq:A.49} and \eqref{A:eq:C4} give
$C[h^3+e^{-cd_{\mathrm L}(p)}+e^{-cd_{\mathrm R}(p)}]e^{-c'd(p,q)}$. For the second, the two center indices differ by at most $C(1+d(r,q))$, so $|s_{pr}-s|\le Ch^3(1+d(p,q))$. For
$0<c'<c$ and $x\ge0$, $(1+x)\mathrm e^{-(c-c')x}\le1+[\mathrm e(c-c')]^{-1}$. The second sum is therefore bounded by
$Ch^3(1+d(p,q))\mathrm e^{-cd(p,q)}\le C_{c,c'}h^3\mathrm e^{-c'd(p,q)}$. For the third sum, the mean-value formula and \eqref{A:eq:C4} give
\[
 |q_{s_{rq}}(r-q)-q_s(r-q)|
 \le |s_{rq}-s|\sup_\theta|\partial_\theta q_\theta(r-q)|
 \le Ch^3e^{-cd(r,q)}.
\]
On its support $d(p,r)\le R_0$, hence $d(r,q)\ge d(p,q)-R_0$. The number and size of the nonzero coefficients $\mathscr K_s(p-r)$ are uniformly bounded. Their sum is consequently
at most $C h^3e^{cR_0}e^{-cd(p,q)}$. In the last sum, a nonzero term requires $d(p,\partial\Lambda_n)\le R_0$. The finite number of such terms is bounded by
$C[e^{-cd_{\mathrm L}(p)}+e^{-cd_{\mathrm R}(p)}]e^{-c'd(p,q)}$. Thus
\[
|\mathcal R_n(p,q)|
\le
C\left[
 h^3+\mathrm{e}^{-cd_{\mathrm L}(p)}+\mathrm{e}^{-cd_{\mathrm R}(p)}
\right]\mathrm{e}^{-cd(p,q)}.
\eqtag{C7}
\]

The entries of $\mathcal R_n$ can be of order one in boundary rows. We control its contribution using exponential decay of the actual finite inverse, which follows from coercivity.
The uniform bounds $cI\le\mathscr K_n\le CI$ and the fixed bandwidth allow a direct application of Theorem~2.4 of \cite{demko1984decay}, giving
\[
|\mathscr K_n^{-1}(p,r)|\le C\mathrm{e}^{-cd(p,r)}.
\eqtag{C8}
\]
From \(\mathscr K_nQ_n=I+\mathcal R_n\),
\[
\mathscr K_n^{-1}-Q_n=-\mathscr K_n^{-1}\mathcal R_n.
\eqtag{C9}
\]
For the $h^3$ term, split the exponential using $d(p,r)+d(r,q)\ge d(p,q)$ and sum the remaining exponential over $r$. This gives $Ch^3e^{-c_1d(p,q)}$. For a boundary
$\sigma\in\{\mathrm L,\mathrm R\}$, we need to keep its distance as well:

\[
\sum_r \mathrm{e}^{-a d(p,r)}\mathrm{e}^{-b d_\sigma(r)}\mathrm{e}^{-a d(r,q)} \le C\mathrm{e}^{-c d(p,q)}\mathrm{e}^{-c\min\{d_\sigma(p),d_\sigma(q)\}}.
\eqtag{C10}
\]
To prove \eqref{A:eq:C10}, put $S=d(p,r)+d(r,q)$, $B=d_\sigma(r)$, $D=d(p,q)$, and $m=\min\{d_\sigma(p),d_\sigma(q)\}$. The triangle inequality gives $S\ge D$. Since distance to a
set is $1$-Lipschitz, $S+B\ge m$. In particular,
\[
 d(p,r)+d_\sigma(r)+d(r,q)
 \ge\tfrac12\bigl(d(p,q)+
       \min\{d_\sigma(p),d_\sigma(q)\}\bigr).
\eqtag{C11}
\]
For $\lambda=\min\{a,b\}$, the summand in \eqref{A:eq:C10} is at most $e^{-\lambda(S+B)}\le e^{-\lambda(D+m)/4}e^{-\lambda S/2}$. The sum of the last factor is at most
$\sum_r e^{-\lambda d(p,r)/2}\le C_\lambda$, by the bounded number of coordinate types at each distance. This proves \eqref{A:eq:C10}. Equations \eqref{A:eq:C7}--\eqref{A:eq:C11}
prove \eqref{A:eq:C1}. Restriction to \(S_n\) gives \eqref{A:eq:C2}, and multiplication by \(2/(\beta n)\asymp h^3\) gives \eqref{A:eq:C3}.
\end{proof}

The hypotheses for the Hessian of the conditional potential for the Jacobi coefficients are checked in Subsection~\ref{A:sec:application-to-the-conditional-jacobi-coefficients}.

To obtain a coupling from the covariance estimate, we multiply by the weights before taking the trace norm. The estimate avoids a factor equal to the number of coordinates and
requires no lower bound on the smallest covariance eigenvalue. This is useful for Gaussian vectors projected to fewer coordinates.

Let $C,D$ be real positive semidefinite $M\times M$ covariance matrices and $W$ a positive diagonal matrix. Write
$\|x\|_W^2=\langle x,Wx\rangle$. \par\Needspace{8\baselineskip}
\begin{lemma}\label{A:res:A.16}
There is a jointly Gaussian coupling of $X\sim N(0,C)$ and $Y\sim N(0,D)$ such that
\[
\mathbb E\|X-Y\|_W^2
\le
\left\|W^{1/2}(C-D)W^{1/2}\right\|_{\mathfrak S_1}.
\eqtag{A.52}
\]
For every $p<\infty$ the same coupling satisfies
\[
\bigl\|\|X-Y\|_W^2\bigr\|_{L^p}
\le
C_p\left\|W^{1/2}(C-D)W^{1/2}\right\|_{\mathfrak S_1}.
\eqtag{A.53}
\]
The coupling can be constructed starting from either Gaussian vector already defined on the probability space. If $C,D,W$ depend Borel measurably on a parameter, the corresponding
conditional Gaussian kernel can also be chosen Borel measurably.
\end{lemma}

\begin{proof}
Set $\widetilde C=W^{1/2}CW^{1/2}$ and $\widetilde D=W^{1/2}DW^{1/2}$, and write $S=\widetilde C^{1/2}$ and $T=\widetilde D^{1/2}$. For a standard real Gaussian vector $Z$, define
\[
 X=W^{-1/2}SZ,\qquad Y=W^{-1/2}TZ.
\]
Then $\mathbb EXX^*=C$ and $\mathbb EYY^*=D$. No commutation between a covariance matrix and $W$ is required. Moreover,
\[
 \mathbb E\|X-Y\|_W^2=\operatorname{Tr}(S-T)^2.
\]
We prove the square-root inequality needed to estimate this trace. Put $H=S-T$ and let $J=\operatorname{sgn}H$, with value zero on $\ker H$. The spectral theorem gives $\|J\|\le1$
and $JH=HJ=|H|$. Since $\widetilde C-\widetilde D=SH+HT$, cyclicity of the finite trace gives
\[
 \operatorname{Tr}J(\widetilde C-\widetilde D)
 =\operatorname{Tr}|H|(S+T).
\]
Choose an orthonormal eigenbasis $e_j$ of $H$, with eigenvalues $h_j$. Positivity of $S,T$ implies
\[
 \langle (S+T)e_j,e_j\rangle
 \ge \big|\langle (S-T)e_j,e_j\rangle\big|=|h_j|.
\]
Consequently
\[
 \operatorname{Tr}H^2
 \le\sum_j |h_j|\langle(S+T)e_j,e_j\rangle
 =\operatorname{Tr}J(\widetilde C-\widetilde D)
 \le\|\widetilde C-\widetilde D\|_{\mathfrak S_1}.
\]
For the last inequality, diagonalize the self-adjoint matrix $\widetilde C-\widetilde D$ and use $|\langle Je,e\rangle|\le1$ on each unit eigenvector. This is the
finite-dimensional square-root inequality of \cite{powers1970free}, proved here with the required constant. It proves \eqref{A:eq:A.52}.

The vector $W^{1/2}(X-Y)$ is centered Gaussian. Diagonalizing its covariance shows that its squared norm has the distribution of $\sum_j\lambda_j Z_j^2$, where $\lambda_j\ge0$ and
$\sum_j\lambda_j=\mathbb E\|X-Y\|_W^2$. For $p\ge1$, Minkowski's inequality gives
\[
 \left\|\sum_j\lambda_j Z_j^2\right\|_{L^p}
 \le\|Z_1^2\|_{L^p}\sum_j\lambda_j.
\]
For $0<p<1$, Jensen's inequality gives the same conclusion with constant one. Thus the single coupling above proves \eqref{A:eq:A.53} for every positive finite $p$.

Its cross covariance is
\[
 K=\mathbb EXY^*=W^{-1/2}STW^{-1/2}.
\]
If $v\in\ker D$, then $v^*Y=0$ almost surely, since its second moment is $v^*Dv=0$. Hence $Kv=0$, so $KD^\dagger D=K$. The Gaussian vector
\[
 V=X-KD^\dagger Y
\]
has covariance $R=C-KD^\dagger K^*\ge0$ and cross covariance with $Y$ equal to $K-KD^\dagger D=0$. Its joint characteristic function with $Y$ therefore factors, proving that $V$
and $Y$ are independent, also when either covariance is singular. It follows that, for an already defined $Y_0\sim N(0,D)$ and an independent standard Gaussian $Z'$, the vector
\[
 X_0=KD^\dagger Y_0+R^{1/2}Z'
 \eqtag{A.52a}
\]
has exactly the same joint law with $Y_0$ as $X$ has with $Y$. In particular, both norm estimates remain valid. Reversing the roles of the two coordinates gives
\[
 Y_0=K^*C^\dagger X_0+(D-K^*C^\dagger K)^{1/2}Z'
\]
when $X_0$ is prescribed instead.

Finally, square roots of nonnegative matrices are continuous, inversion of the positive matrix $W$ is continuous, and
\[
 D^\dagger=\lim_{m\to\infty}D(D+m^{-1}I)^{-2}
\]
is Borel measurable. The mean and covariance in the displayed conditional formulas are therefore Borel functions of the parameters and of the prescribed vector. Integrating a
bounded Borel function of these formulas against the fixed standard Gaussian law defines a Borel probability kernel. This construction does not depend on $p$.
\end{proof}

The covariance comparison determines a linear image of the coefficient vector, whereas the transport acts on the full vector. We therefore sample the remaining Gaussian coordinates
from their conditional distribution given this image. The following formula preserves the prescribed image, including when its covariance is singular.

Let $\Sigma\ge0$ be a $d\times d$ covariance matrix, $L:\mathbb R^d\to\mathbb R^m$, and $C=L\Sigma L^*$. Fix $X\sim N(0,C)$ and an independent standard Gaussian vector
$Z\in\mathbb R^d$.

\begin{lemma}\label{A:res:A.16A}
The vector
\[
G=\Sigma L^*C^\dagger X+
\left(\Sigma-\Sigma L^*C^\dagger L\Sigma\right)^{1/2}Z
\eqtag{A.53a}
\]
has law $N(0,\Sigma)$ and satisfies $LG=X$ almost surely.
\end{lemma}

Borel dependence of $\Sigma$ and $L$ makes this a Borel conditional Gaussian kernel.

\begin{proof}
Put $B=\Sigma^{1/2}L^*$, so that $C=B^*B$, and define
\[
 P=BC^\dagger B^*.
\]
The identities $C^\dagger C C^\dagger=C^\dagger$ and $(C^\dagger)^*=C^\dagger$ give $P^2=P=P^*$. Moreover, $\ker C=\ker B$, since $v^*Cv=\|Bv\|^2$, and hence $BC^\dagger C=B$. Thus
$P$ is the orthogonal projection onto $\operatorname{ran}B$. In particular,
\[
 R:=\Sigma-\Sigma L^*C^\dagger L\Sigma
   =\Sigma^{1/2}(I-P)\Sigma^{1/2}\ge0,
\]
which justifies its square root in \eqref{A:eq:A.53a}.

The two terms defining $G$ are independent centered Gaussian vectors. Since $C^\dagger C C^\dagger=C^\dagger$, their covariances sum to
\[
 \Sigma L^*C^\dagger C C^\dagger L\Sigma+R
 =\Sigma L^*C^\dagger L\Sigma+R=\Sigma.
\]
Therefore $G\sim N(0,\Sigma)$.

The orthogonal projection onto $\operatorname{ran}C$ is $CC^\dagger$. The complementary component of $X$ has zero covariance, and therefore $CC^\dagger X=X$ almost surely. Also
\[
 LRL^*=C-CC^\dagger C=0.
\]
For every vector $v$, this gives $\|R^{1/2}L^*v\|^2=v^*LRL^*v=0$, so $LR^{1/2}=0$. Applying $L$ to the formula for $G$ now gives
\[
 LG=CC^\dagger X+LR^{1/2}Z=X
\]
almost surely. These arguments include $C=0$. The square-root and pseudoinverse formulas in Lemma~\ref{A:res:A.16} show that the displayed sampling map is Borel in
$(\Sigma,L,X,Z)$; integration in the independent $Z$ gives the claimed conditional kernel.
\end{proof}

\begin{corollary}\label{A:res:A.17}
Let \(C,D\ge0\) be finite covariance matrices. Suppose \(|C(i,j)-D(i,j)| \le Ch^6\mathrm{e}^{-c|i-j|}.\) Let
\(W_j=(1+x_j)^{-7/2}.\) Then there is a coupling of \(X\sim N(0,C)\) and \(Y\sim N(0,D)\) such that
\[
\left\|
\sum_jW_j|X_j-Y_j|^2
\right\|_{L^p}
\le
C_ph^5.
\eqtag{A.54}
\]
\end{corollary}

\begin{proof}
For a finite matrix, the triangle inequality applied to the matrix units gives $\|M\|_{\mathfrak S_1}\le\sum_{i,j}|M_{ij}|$. Using $2\sqrt{W_iW_j}\le W_i+W_j$ and symmetry in
$i,j$, we obtain
\[
 \begin{aligned}
 \|W^{1/2}(C-D)W^{1/2}\|_{\mathfrak S_1}
 &\le Ch^6\sum_{i,j}\sqrt{W_iW_j}e^{-c|i-j|}\\
 &\le Ch^6\sum_iW_i\sum_je^{-c|i-j|}
 \le Ch^6\sum_iW_i\le Ch^5.
 \end{aligned}
\]
Lemma~\ref{A:res:A.16} now proves \eqref{A:eq:A.54}.
\end{proof}

\subsubsection{Spectral factorization and coupling to Brownian increments}\label{A:sec:spectral-factorization-and-coupling-to-brownian-increments}

We now construct a positive covariance from the local scalar symbols. The factorization is used twice: multiplying a factor by its adjoint makes the finite covariance positive, and
evaluating the factor at $1$ identifies the part of the field that accumulates in partial sums. The remaining part will be a discrete derivative. The classical spectral
factorization is given in Theorem~7.13 of \cite{wiener1957prediction}. Proposition~\ref{A:res:A.18} proves the scalar analytic factorization, including the uniform annulus and
parameter estimates needed here. In this spectral factorization argument $q_s(z)$ is a scalar covariance symbol, distinct from the matrix-valued inverse coefficients $q_s(p-q)$
above. For $\rho>0$, write $\mathcal U_\rho=\{z\in\mathbb C:\mathrm e^{-\rho}<|z|<\mathrm e^\rho\}$ and let
\(H^\infty(\mathcal U_\rho)\) have its supremum norm.

For an absolutely summable real sequence $(g_m)$ and independent standard normals $(\xi_j)_{j\in\mathbb Z}$, the Gaussian series $Y_j=\sum_mg_m\xi_{j-m}$ satisfies
\[
 \mathbb E[Y_iY_j]=\sum_mg_mg_{m+j-i},\qquad
 \sum_r\mathbb E[Y_0Y_r]e^{\mathrm i r\theta}
   =\left|\sum_mg_me^{\mathrm i m\theta}\right|^2.
\]
Absolute summability justifies both rearrangements. In particular,
\[
 \sum_{r\in\mathbb Z}\mathbb E[Y_0Y_r]
   =\left(\sum_m g_m\right)^2,
 \qquad \mathbb E Y_0^2=\sum_m g_m^2.
\]
The positive square root of this covariance sum is the coefficient multiplying the Brownian increments. We construct $g_s$ with real exponentially decaying coefficients and with
$g_s(1)>0$, fixing the sign needed to couple to the increments of the driving Brownian motion.

Suppose that \(s\mapsto q_s\) belongs to \(C^1([0,s_0];H^\infty(\mathcal U_\rho))\), with
\[
 \sup_{0\le s\le s_0}
 \bigl(\|q_s\|_{H^\infty(\mathcal U_\rho)}
       +\|\partial_sq_s\|_{H^\infty(\mathcal U_\rho)}\bigr)
 \le M.
\eqtag{factor-analytic-norm}
\]
At the parameter boundary points the derivative is one-sided. Suppose also that the circle values are real and even, and
\(q_s(\mathrm{e}^{\mathrm{i}\theta})\ge c_0>0\) for every \(s,\theta\).

\begin{proposition}\label{A:res:A.18}
There is a nonvanishing analytic factor $g_s$ on a common disk of radius greater than one, with real coefficients, such that $q_s(e^{\mathrm i\theta})=|g_s(e^{\mathrm i\theta})|^2$
and $g_s(1)=\sqrt{q_s(1)}>0$. It is $C^1$ in $s$ and its coefficients satisfy
\[
 |g_m(s)|+|\partial_sg_m(s)|\le Ce^{-c|m|},
 \qquad m\in\mathbb Z,\quad 0\le s\le s_0.
\eqtag{factor-coefficient-bound}
\]
Here $g_m=0$ for $m<0$, and the radii and constants depend only on $\rho,M,c_0$.
\end{proposition}

\begin{proof}
Cauchy's formula bounds \(\partial_zq_s\) uniformly on \(\mathcal U_{\rho/2}\), because the closure of this annulus has positive distance from the boundary of \(\mathcal U_\rho\).
For \(z=t\mathrm{e}^{\mathrm{i}\theta}\) in a still smaller annulus,
\[
 |q_s(z)-q_s(\mathrm{e}^{\mathrm{i}\theta})|\le C|t-1|.
\]
Choose \(0<\rho_1<\rho/2\) so that the right side is at most \(c_0/2\) whenever
\(\mathrm{e}^{-\rho_1}<t<\mathrm{e}^{\rho_1}\). Then \(\operatorname{Re}q_s(z)\ge c_0/2\) throughout this annulus. The
principal logarithm \(L_s(z)=\log q_s(z)\) is consequently single-valued and holomorphic there, with its circle value equal to the real logarithm. Both \(L_s\) and
\(\partial_sL_s=(\partial_sq_s)/q_s\) are uniformly bounded. The identity
\[
 \log(q_s+v)-\log q_s
 =\int_0^1\frac{v}{q_s+tv}\,\mathrm{d}t
\]
for sufficiently small norm \(v\) proves differentiability and continuity of the derivative in this norm, with the same logarithm branch for all parameters.

Let \(c_m(s)\) be the Laurent coefficients of \(L_s\). Integrating on \(|z|=\mathrm{e}^{\rho_1/2}\) for \(m\ge0\) and on
\(|z|=\mathrm{e}^{-\rho_1/2}\) for \(m<0\), and doing the same for \(\partial_sL_s\), gives
\[
 |c_m(s)|+|\partial_sc_m(s)|\le C \mathrm{e}^{-\rho_1|m|/2}.
\eqtag{factor-log-coefficients}
\]
The logarithm is real and even on the circle. Its sine Fourier coefficients vanish, and its cosine coefficients are real, so $c_{-m}(s)=c_m(s)\in\mathbb R$.

Define
\[
 F_s(z)=\frac12c_0(s)+\sum_{m\ge1}c_m(s)z^m,
 \quad g_s(z)=\exp F_s(z).
\eqtag{canonical-factor}
\]
By \eqref{A:eq:factor-log-coefficients}, the series for \(F_s\) and \(\partial_sF_s\) converge uniformly in \(s\) on the closed disk \(|z|\le \mathrm{e}^{\rho_1/4}\). They are
bounded there; hence so are \(g_s\) and
\(\partial_sg_s=g_s\partial_sF_s\). The exponential never vanishes. Cauchy's formula on the circle of radius
\(\mathrm{e}^{\rho_1/4}\) gives \eqref{A:eq:factor-coefficient-bound}, with negative coefficients zero. The real
coefficients of \(F_s\) give real coefficients of the expansion of \(g_s\). On the unit circle the absolutely convergent series satisfy
\[
 2\operatorname{Re}F_s(\mathrm{e}^{\mathrm{i}\theta})
 =c_0(s)+\sum_{m\ge1}c_m(s)
                  (\mathrm{e}^{\mathrm{i}m\theta}+\mathrm{e}^{-\mathrm{i}m\theta})
 =\log q_s(\mathrm{e}^{\mathrm{i}\theta}).
\]
Exponentiation proves the modulus identity. Finally \(F_s(1)\) is real, so \(g_s(1)=\mathrm{e}^{F_s(1)}>0\), and the modulus identity at \(1\) gives its stated value.
\end{proof}

For a fixed parameter the Fourier coefficients of $q_s$ form a covariance sequence. To construct a positive covariance with varying parameters, we evaluate the factor row by row
and multiply by its adjoint.

Let $s_j\in[0,s_0]$ be deterministic, with $|s_{j+1}-s_j|\le Lh^3$, and use the coefficients of Proposition~\ref{A:res:A.18}. Only parameter differences enter the covariance
estimate; the size of $s_j$ will enter the Brownian remainder estimate later. For a finite row interval define
\[
(\mathcal G_nv)_j
=
h^{3/2}\sum_{m\in\mathbb Z}g_m(s_j)v_{j-m},
\quad
C_n^{\mathrm{fac}}=\mathcal G_n\mathcal G_n^*.
\]

Define the matrix
\[
C_n^{\mathrm{mid}}(i,j)
=
\frac{h^3}{2\pi}
\int_{-\pi}^{\pi}
\mathrm{e}^{\mathrm{i}(i-j)\theta}q_{s_{ij}}(\mathrm{e}^{\mathrm{i}\theta})\,\mathrm{d}\theta,
\quad
s_{ij}=s_{\lfloor(i+j)/2\rfloor}.
\]
The matrix \(C_n^{\mathrm{mid}}\) need not be positive.

\begin{proposition}\label{A:res:A.18A}
\(C_n^{\mathrm{fac}}\ge0\), and \(C_n^{\mathrm{fac}}(i,j) = h^3\sum_m g_m(s_i)g_{m+j-i}(s_j).\) Moreover,
\[
|C_n^{\mathrm{fac}}(i,j)-C_n^{\mathrm{mid}}(i,j)|
\le
Ch^6\mathrm{e}^{-c|i-j|}.
\eqtag{A.57d}
\]
\end{proposition}

Consequently, if an exact positive covariance \(C_n\) satisfies
\[
|C_n(i,j)-C_n^{\mathrm{mid}}(i,j)|
\le
Ch^6\mathrm{e}^{-c|i-j|},
\eqtag{A.57e}
\]
then
\[
|C_n(i,j)-C_n^{\mathrm{fac}}(i,j)|
\le
Ch^6\mathrm{e}^{-c|i-j|}.
\eqtag{A.57f}
\]

\begin{proof}
The exponential coefficient bound makes every row of $\mathcal G_n$ square summable. Hence, for every finite vector $v$,
$\langle v,C_n^{\mathrm{fac}}v\rangle=\|\mathcal G_n^*v\|_2^2\ge0$. Multiplication of the two rows gives the stated covariance formula. For $r=j-i$, Fourier inversion gives
$C_n^{\mathrm{mid}}(i,j)=h^3\sum_mg_m(s_{ij})g_{m+r}(s_{ij})$. Subtract the summands as
\[
\begin{aligned}
 &g_m(s_i)g_{m+r}(s_j)-g_m(s_{ij})g_{m+r}(s_{ij})\\
 &\qquad=[g_m(s_i)-g_m(s_{ij})]g_{m+r}(s_j)
       +g_m(s_{ij})[g_{m+r}(s_j)-g_{m+r}(s_{ij})].
\end{aligned}
\]
The center distances give $|s_i-s_{ij}|+|s_j-s_{ij}|\le Ch^3(1+|r|)$. The mean-value theorem and the exponential bounds therefore yield
\[
|C_n^{\mathrm{fac}}(i,j)-C_n^{\mathrm{mid}}(i,j)| \le Ch^6(1+|r|) \sum_m \mathrm{e}^{-c|m|}\mathrm{e}^{-c|m+r|} \le Ch^6(1+|r|)^2\mathrm{e}^{-c|r|} \le Ch^6\mathrm{e}^{-c'|r|}.
\]
Here the convolution is explicitly bounded by
\[
 \sum_{m\in\mathbb Z}e^{-c|m|}e^{-c|m+r|}
 =e^{-c|r|}\left(|r|+1+\frac{2}{e^{2c}-1}\right).
\]
For $0<c'<c$ put $d=c-c'$. Differentiating $(1+t)^2e^{-dt}$ shows that its maximum on $[0,\infty)$ occurs at $t=0$ or $t=2/d-1$, when the latter is nonnegative. In particular,
\[
 (1+t)^2e^{-ct}\le C_d e^{-c't},\qquad
 C_d=\max\{1,4e^d/(e^2d^2)\}.
\]
This proves \eqref{A:eq:A.57d}. Combining it with \eqref{A:eq:A.57e} gives
\[
 |C_n(i,j)-C_n^{\mathrm{fac}}(i,j)|
 \le|C_n(i,j)-C_n^{\mathrm{mid}}(i,j)|
       +|C_n^{\mathrm{mid}}(i,j)-C_n^{\mathrm{fac}}(i,j)|
 \le C h^6e^{-c'|i-j|},
\]
after replacing the two exponential rates by their minimum.
\end{proof}

For the almost sure estimate below, take $h=h_n\asymp n^{-1/3}$ as in the matrix models, with the coupled vectors for all $n$ defined on one probability space. Their dependence
across sizes is arbitrary.

\begin{corollary}\label{A:res:A.18B}
Under the hypotheses of Proposition~\ref{A:res:A.18A}, let \(X\sim N(0,C_n)\) and \(Y=\mathcal G_n\xi\), where the coordinates of \(\xi=(\xi_k)_{k\in\mathbb Z}\) are independent
standard normal variables. If \eqref{A:eq:A.57e} holds and \(W_j=(1+x_j)^{-7/2}\), then there is a jointly Gaussian coupling of \(X\) and \(Y\) such that
\(\left\| \sum_jW_j|X_j-Y_j|^2 \right\|_{L^p} \le C_ph^5\) for every \(p<\infty\). In particular, almost surely, for all sufficiently large \(n\),
\(\sum_jW_j|X_j-Y_j|^2 \le C_\omega h^4\).
\end{corollary}

\begin{proof}
We define each coordinate of $\mathcal G_n\xi$ by its Gaussian series. The exponential coefficient bounds give convergence in every finite $L^p$ and absolute convergence almost
surely, since $\mathbb E\sum_m|g_m(s_j)\xi_{j-m}|<\infty$. Its covariance is $C_n^{\mathrm{fac}}$. Both covariances are positive semidefinite. Apply Corollary~\ref{A:res:A.17} to
\eqref{A:eq:A.57f}. If $E_n=\sum_jW_j|X_j-Y_j|^2$, the resulting moment bound gives
\[
 \mathbb P(E_n>h_n^4)\le h_n^{-4p}\mathbb E E_n^p
       \le C_ph_n^p.
\]
For $p>3$ these probabilities are summable. Under any joint distribution with the specified coupled pair at each size, Borel--Cantelli gives a finite $n_0(\omega)$ beyond which the
bound holds. Lemma~\ref{A:res:A.16} supplies the conditional law of $X$ given $Y$. We may therefore take $\xi$ to be the normalized increments of one Brownian motion and sample $X$
conditionally, keeping those increments fixed.
\end{proof}

The independent part is obtained by subtracting the value of the factor at $1$. Indeed, if $g_s(z)-g_s(1)=(1-z)p_s(z)$ and $E_j=\sum_mp_m(s)\xi_{j-m}$, then, at a fixed $s$,
\[
 \sum_mg_m(s)\xi_{j-m}=g_s(1)\xi_j+E_j-E_{j-1},\qquad
 \sum_{j=\ell}^{r}(E_j-E_{j-1})=E_r-E_{\ell-1}.
\]
Thus summing the correlated correction leaves only its two boundary values. This is the Beveridge--Nelson decomposition \cite{beveridge1981decomposition}; the algebraic identity
and its finite filter antecedents are recorded in \cite{phillips1992asymptotics}. In the operator comparison, summation by parts preserves the cancellation in this derivative.
Proposition~\ref{A:res:A.20} accounts for the variation of $s_j$.

Under the hypotheses of Proposition~\ref{A:res:A.18}, suppose additionally that
\[
 |q_s(1)-4/\beta|\le C_0s,\quad 0\le s\le s_0.
\eqtag{A.58}
\]
For the factor $g_s$ define
\[
 p_s(z)=\frac{g_s(z)-g_s(1)}{1-z}.
\eqtag{A.59}
\]

\begin{corollary}\label{A:res:A.19} The singularity at $z=1$ is removable.  On a common disk of radius
greater than one, \(p_s\) is analytic and \(C^1\) in \(s\) in the supremum norm. Its power series coefficients are real and, when extended by zero to negative indices, satisfy
\(|p_m(s)|+|\partial_sp_m(s)|\le C\mathrm{e}^{-c|m|}\). Moreover, with
\(a_s=g_s(1)-2/\sqrt\beta\),
\[
 g_s(z)=\frac2{\sqrt\beta}+(1-z)p_s(z)+a_s,
 \quad |a_s|\le Cs,
 \quad |\partial_sa_s|\le C.
\eqtag{factor-zero-mode}
\]
\end{corollary}

The last remainder is independent of \(z\). The constants are uniform for \(\beta\) in compact positive intervals whenever the preceding hypotheses hold uniformly.

\begin{proof}
For $m\ge0$, absolute convergence at $z=1$ gives
\[
 p_m(s)=-\sum_{r>m}g_r(s),\qquad
 \partial_sp_m(s)=-\sum_{r>m}\partial_sg_r(s).
\]
Indeed, multiplication by $1-z$ gives coefficient $g_m(s)$ for $m\ge1$ and constant coefficient $g_0(s)-g_s(1)$. The exponential bounds on $g_r$ justify the derivative and give
$|p_m|+|\partial_sp_m|\le Ce^{-cm}$.

For the analytic extension at $z=1$, choose radii \(1<R_2<R_1\) inside the disk supplied by Proposition~\ref{A:res:A.18}. Cauchy's formula bounds \(\partial_zg_s\) and
\(\partial_s\partial_zg_s\) uniformly on
\(|z|\le R_2\). On that disk, convexity of the disk and the fundamental theorem of calculus along the segment from \(1\)
to \(z\) give the formula
\[
 p_s(z)=-\int_0^1
             (\partial_zg_s)(1+t(z-1))\,\mathrm{d}t.
\]
The uniform bound on $\partial_s\partial_zg_s$ allows differentiation under this integral:
\[
 \partial_sp_s(z)=-\int_0^1
      (\partial_s\partial_zg_s)(1+t(z-1))\,\mathrm dt.
\]
Both integrands are uniformly bounded on $|z|\le R_2$, and the second is continuous in $s$ there. These formulas extend $p_s$ through $z=1$ with uniform bounds on $p_s$ and
$\partial_sp_s$. The coefficients are real by the preceding series formula. At zero frequency,
\[
 |a_s|=
 \frac{|q_s(1)-4/\beta|}{g_s(1)+2/\sqrt\beta}\le Cs,
 \quad
 \partial_sa_s=\frac{\partial_sq_s(1)}{2\sqrt{q_s(1)}}.
\]
The denominators are bounded below by the preceding positivity hypotheses. This proves \eqref{A:eq:factor-zero-mode}. Finally the coefficients of the expansion, extended by zero to
negative Laurent indices, obey \([(1-z)p_s]_m=p_m(s)-p_{m-1}(s)\) for all \(m\in\mathbb Z\).

For a fixed parameter $s$, let $Z_j=\sum_{m\ge0}g_m(s)\xi_{j-m}$. Absolute summability of the coefficients just proved permits summation of the covariance series:
\[
 \sum_{r\in\mathbb Z}\operatorname{Cov}(Z_0,Z_r)
 =\sum_{m,\ell\ge0}g_m(s)g_\ell(s)
 =g_s(1)^2=q_s(1).
\]
The variance at a single site is instead $\sum_mg_m(s)^2$. The identity $q_0(1)=4/\beta$ gives the coefficient $2/\sqrt\beta$ multiplying the Brownian increments. The remaining
correlations are represented by $(1-z)p_s(z)$.

\end{proof}

Variation of $s$ produces two errors. The difference $g_{s_j}(1)-2/\sqrt\beta$ is controlled by $s_j$. Comparing the shifted coefficient in $e_{j-1}$ with its value at $s_j$
introduces an error controlled by $s_j-s_{j-1}$. We therefore impose bounds on both the parameter and its increment. Let \(h=h_n\asymp n^{-1/3}\), \(0<h\le1\), and let
\(I_n=\{a_n,\ldots,b_n\}\subset\{1,\ldots,\lfloor Cn\rfloor\}\) be any nonempty deterministic interval. We evaluate its
weights at \(x_j=(j-\tfrac12)h\), without recentering at \(a_n\), and \(X_j=1+x_j\). Let the deterministic sequence
\(s_j\) be defined also at \(a_n-1\) and satisfy
\[
 s_j\in[0,s_0]\quad(a_n-1\le j\le b_n),\qquad |s_j-s_{j-1}|\le Lh^3,\quad 0\le s_j\le
 Lh^2X_j\quad(j\in I_n).
\eqtag{factor-parameter-bounds}
\]
The size bound is imposed only at the stated nonnegative indices. Let \(g_s,p_s\) be the factors in Proposition~\ref{A:res:A.18} and Corollary~\ref{A:res:A.19}. On one two-sided
Brownian motion put
\(\xi_j=h^{-1/2}(B(jh)-B((j-1)h))\), \(j\in\mathbb Z\). For each fixed \(n\) these variables are independent standard
real Gaussians. Define
\[
 \zeta_j^{\mathrm{fac}}=h^{3/2}\sum_mg_m(s_j)\xi_{j-m},\quad
 \zeta_j^B=\frac2{\sqrt\beta}h^{3/2}\xi_j,\quad
 e_j=h^{1/2}\sum_mp_m(s_j)\xi_{j-m}.
\]
The last definition includes \(a_n-1\). The sums converge in every finite \(L^p\) and absolutely almost surely: the deterministic exponential upper bound has summable expectation
after multiplication by \(|\xi_{j-m}|\). This is simultaneous over the countably many indices and matrix sizes.

\begin{proposition}\label{A:res:A.20} Under
\eqref{A:eq:factor-parameter-bounds}, there is a centered Gaussian remainder \(\rho\) such that, for every \(j\in I_n\),
\(\zeta_j^{\mathrm{fac}}=\zeta_j^B+h(e_j-e_{j-1})+\rho_j\). It satisfies
\[
 \mathbb E\sum_{j\in I_n}X_j^{-7/2}|\rho_j|^2\le Ch^6,
\eqtag{A.63}
\]
and, for every \(1\le p<\infty\),
\[
 \left\|\sum_{j\in I_n}X_j^{-7/2}|\rho_j|^2\right\|_{L^p}
       \le C_ph^6,\quad
 \left\|\sum_{j\in I_n}hX_j^{-3/2}|e_j|^2\right\|_{L^p}
       \le C_ph.
\eqtag{factor-moment-bounds}
\]
\end{proposition}

On the given coupling, almost surely for all sufficiently large \(n\),
\[
 \sum_{j\in I_n}X_j^{-7/2}|\rho_j|^2\le C_\omega h^4,
\eqtag{A.64}
\]
and
\[
 \sum_{j\in I_n}hX_j^{-3/2}|e_j|^2\le C_\omega h.
\eqtag{A.65}
\]
The finite size constants are uniform over all such coordinate sequences, index intervals, and additional compact parameter families satisfying the stated bounds. The almost sure
estimates also hold simultaneously along any fixed countable sequence of parameters, allowing dependence between matrix sizes.

\begin{proof}
Write \(a_s=g_s(1)-2/\sqrt\beta\), as in \eqref{A:eq:factor-zero-mode}. Define
\[
 \rho_j^{(0)} =h^{3/2}a_{s_j}\xi_j,\qquad \rho_j^{(1)} =h^{3/2}\sum_m
 \bigl(p_{m-1}(s_{j-1})-p_{m-1}(s_j)\bigr)\xi_{j-m},\qquad \rho_j =\rho_j^{(0)}+\rho_j^{(1)}.
\eqtag{factor-remainder-formula}
\]
Reindexing the absolutely convergent series gives
\[
 h(e_j-e_{j-1})=h^{3/2}\sum_m
       \bigl(p_m(s_j)-p_{m-1}(s_{j-1})\bigr)\xi_{j-m}.
\]
Add \(\rho_j^{(1)}\) and use \([(1-z)p_{s_j}]_m=p_m(s_j)-p_{m-1}(s_j)\), then add \(\rho_j^{(0)}\) and the Brownian term. This proves the identity, including \(j=a_n\), where the
stipulated predecessor \(s_{a_n-1}\) is needed. We estimate the variances of the two terms in their joint distribution and then apply the Gaussian quadratic form bound to each row.

The two bounds in \eqref{A:eq:factor-parameter-bounds} give
\[
 \begin{aligned}
  \mathbb E|\rho_j^{(0)}|^2
  &=h^3|a_{s_j}|^2\le Ch^3s_j^2\le Ch^7X_j^2,\\
  \mathbb E|\rho_j^{(1)}|^2
  &=h^3\sum_m|p_{m-1}(s_{j-1})-p_{m-1}(s_j)|^2\le Ch^9.
 \end{aligned}
\eqtag{factor-remainder-variances}
\]
For every \(q>1\), monotonicity of \((1+x)^{-q}\) gives \(\sum_{j\ge1}X_j^{-q}\le C_qh^{-1}\); this follows by comparison with its integral on \([0,\infty)\), with at most one
initial interval term. Therefore
\[
 \mathbb E\sum_{j\in I_n}X_j^{-7/2}|\rho_j|^2
 \le Ch^7\sum_{j\ge1}X_j^{-3/2}
       +Ch^9\sum_{j\ge1}X_j^{-7/2}
 \le Ch^6.
\]
For a finite centered real Gaussian vector \(Z\) and a nonnegative weight matrix \(W\), diagonalize
\(Q=W^{1/2}\operatorname{Cov}(Z)W^{1/2}\). If its eigenvalues are \(\lambda_k\ge0\), then \(Z^{\mathsf T}WZ\) has the
law of \(\sum_k\lambda_kN_k^2\), with independent standard Gaussians \(N_k\). Minkowski's inequality gives
\[
 \|Z^{\mathsf T}WZ\|_{L^p}
 \le\|N_1^2\|_{L^p}\operatorname{Tr}Q.
\eqtag{gaussian-weighted-quadratic-moment}
\]
Apply this to \(\rho\), whose coordinates need not be independent. This proves its moment bound in \eqref{A:eq:factor-moment-bounds}. Markov gives failure probability \(C_ph^{2p}\)
at threshold \(h^4\). Since \(h\asymp n^{-1/3}\), choosing \(p>3/2\) makes the sum over \(n\) finite and proves \eqref{A:eq:A.64}.

Put $\mathcal E_e=\sum_{j\in I_n}hX_j^{-3/2}|e_j|^2$. To prove $\mathcal E_e\le Ah$ for every $n\ge n_0(\omega)$ on an event of probability one, we apply the concentration bound
for Gaussian quadratic forms to its covariance. Let $P_h(j,k)=p_{j-k}(s_j)$, with rows \(j\in I_n\) and columns \(k\in\mathbb Z\). The exponential coefficient bound gives

\[
 \sup_j\sum_k|P_h(j,k)|\le C,\quad
 \sup_k\sum_{j\in I_n}|P_h(j,k)|\le C.
\]
Cauchy--Schwarz with these two absolute row and column sums shows \(\|P_h\|_{\mathrm{op}}\le C\). Thus
\(\operatorname{Cov}(e)=hP_hP_h^{\mathsf T}\) has norm at most \(Ch\), and each diagonal entry is at most \(Ch\). Let
\(W_e=\operatorname{diag}(hX_j^{-3/2})\) and \(Q_e=W_\mathrm{e}^{1/2}\operatorname{Cov}(e)W_\mathrm{e}^{1/2}\). It
follows that
\[
 \operatorname{Tr}Q_e\le C_1h,\quad
 \|Q_e\|_{\mathrm{op}}\le C_2h^2,\quad
 \operatorname{Tr}(Q_\mathrm{e}^2)\le C_1C_2h^3.
\eqtag{factor-primitive-covariance}
\]
The trace and \eqref{A:eq:gaussian-weighted-quadratic-moment} prove the second estimate in \eqref{A:eq:factor-moment-bounds}. For the almost sure bound apply equation~(4.1) in
Lemma~1 of \cite{laurent2000adaptive} to the eigenvalues of $Q_e$. With $A\ge4C_1$, put $s=(A-2C_1)/(4C_2h)$. The threshold in that inequality is at most
\[
 \operatorname{Tr}Q_e+2\|Q_e\|_{\mathfrak S_2}\sqrt s
                         +2\|Q_e\|_{\mathrm{op}}s
 \le h\{A/2+\sqrt{C_1(A-2C_1)}\}\le Ah.
\]
Consequently
\[
 \mathbb P(\mathcal E_e>Ah)
 \le\exp\!\left(-\frac{A-2C_1}{4C_2h}\right)
 \le\mathrm e^{-A/(8C_2h)}.
 \eqtag{factor-primitive-tail}
\]
This is summable in \(n\) and proves \eqref{A:eq:A.65}. The bounds are uniform over the stated coordinate sequences and index intervals. Borel--Cantelli therefore applies on the
common probability space carrying the Brownian motion, also along any fixed countable parameter sequence.

For integers $a\le b$ in the parameter interval, summing the exact identity in the proposition gives
\[
 \sum_{j=a}^b(\zeta_j^{\mathrm{fac}}-\zeta_j^B-\rho_j)
       =h(e_b-e_{a-1}).
\]
Since $\mathbb E e_j^2\le Ch$, Cauchy--Schwarz bounds the variance of the right side by $4Ch^3$, independently of $b-a$. The boundary point $e_{a-1}$ is the value supplied by the
same convolution. Restriction to the middle interval retains this value, as does the constant extension used later in summation by parts.

\end{proof}

\subsubsection{Application to the conditional Jacobi coefficients}\label{A:sec:application-to-the-conditional-jacobi-coefficients}

The preceding estimates now give the covariance comparison needed for the conditional Jacobi law. For fixed exterior data $\tau$, the covariance $C_{n,\mathrm{mid}}^{\tau,+}$ in
\eqref{A:eq:14.45} is obtained from the inverse Hessian of the original conditional potential at its minimizer by the linear coefficient map: the extension in
Subsection~\ref{A:sec:gaussian-approximation-of-the-conditional-gibbs-measure} leaves that Hessian unchanged. Comparing it with $C_n^{\mathrm{fac}}=\mathcal G_n\mathcal G_n^*$ then
separates $\mathcal G_n\xi$ into Brownian increments and a discrete derivative. We first compare the finite Hessians, prove positivity, and compute the variance of the increments.

Throughout the conditional estimates, $\tau$ belongs to $\mathcal T_n^{\mathrm{out}}(M_0)$. Constants are deterministic and uniform over that set. The completed construction
retains the original conditional law on the complementary exterior configurations.

Use $H_n^{\mathrm G}$ from \eqref{A:eq:gaussian-comparison-operator} and the coefficient $\zeta_{j,n}^B=h^2\omega_{j,h}$ constructed from the Brownian increment from
\eqref{A:eq:brownian-difference}. Replacing the diagonal coefficients on the middle interval gives
\[
\zeta_{j,n}^{B,\mathrm{mid}}
=
\begin{cases}
\zeta_j,&j\in I_n^{\mathrm L},\\
\zeta_{j,n}^{B},&j\in I_n^{\mathrm{mid}},\\
\zeta_j,&j\in I_n^{\mathrm{far}},
\end{cases}
\quad
H_n^{B,\mathrm{mid}}
=
H_n^\circ+h^{-2}M_{\zeta_n^{B,\mathrm{mid}}}.
\]

On $I_n^{\mathrm{mid}}$,
\[
\widehat\zeta_n^{\mathrm G}
=
\mathcal L_{\mathrm{diag}}^{\mathrm{mid}}G_n^{\mathrm{tr}},
\quad
G_n^{\mathrm{tr}}
\sim N(0,(K_{n,\mathrm{tr}}^\tau)^{-1}),
\]
and \(|C_{n,\mathrm{mid}}^{\tau,+}(i,j)| \le C h^3\mathrm{e}^{-c|i-j|}.\) In this application to the conditional Jacobi coefficients, scalar indices $i,j$ and unqualified sums run
over $I_n^{\mathrm{mid}}=\{J_n^{\mathrm L}+1,\ldots,J_n^{\mathrm{reg}}\}$. Write $C_n^\tau=C_{n,\mathrm{mid}}^{\tau,+}$ and $X_j=1+x_j$. Every coefficient used in the linear map on
the middle interval is separated from both boundaries of the conditional Hessian by the logarithmic distances in \eqref{A:eq:L10}--\eqref{A:eq:L11}. We keep the Hessian on the
larger transport interval when inverting it. For these indices and their immediate predecessor use \(s_j=s_{j,n}\) from \eqref{A:eq:exact-profile-shift}. Whenever a sequence
indexed by all integers is convenient, use the explicit extension
\[
 s_j=\min\left\{s_{\mathrm{reg}},\,
         \max\left\{0,\frac{r+j-\tfrac12+\beta^{-1}}n\right\}\right\},
 \quad j\in\mathbb Z.
\eqtag{profile-extension}
\]
It agrees with the exact coordinate sequence on \(I_n^{\mathrm{mid}}\) and at \(J_n^{\mathrm L}\), stays in
\([0,s_{\mathrm{reg}}]\), and has \(|s_{j+1}-s_j|\le1/n\), because projection onto a real interval is \(1\)-Lipschitz.
At the indices where the formula for the projection changes we use finite differences. The bound \(s_j\le Ch^2(1+x_j)\) holds at every nonnegative site used in the weighted sums,
by \eqref{A:eq:profile-size-increment} and monotonicity of the projection. The Gaussian driving increments remain indexed by all of \(\mathbb Z\); they are never truncated at the
scalar cutoffs.

Set $s_0=s_{\mathrm{reg}}$ and abbreviate $a(s)=a^\dagger(s)$, $b(s)=b^\dagger(s)$ for $s\in[0,s_0]$. Use the local energy $F_{V,s}$ from \eqref{A:eq:local-energy}. Let
\(\mathsf H_s = D^2_{a,b}F_{V,s}(a(s),b(s))\) be the \(2\times2\) Hessian of the local energy at the minimizer. The term $\operatorname{Tr}V(T)$ couples neighboring coefficients.
In a constant environment, the Hessian of the full energy is a translation-invariant finite-range operator $\mathscr K_s$ on sequences $u=(u_j^A,u_j^B)_{j\in\mathbb Z}$. Let
$\mathsf K_s(\theta)$ be its $2\times2$ Hermitian Fourier symbol. Constant perturbations of the two coefficient types give
\[
  \mathsf K_s(0)=D^2_{a,b}F_{V,s}(a(s),b(s))=\mathsf H_s.
\]
One can verify this without using a non-square-summable test vector: put a constant vector on $m$ consecutive sites and divide the Hessian form by $m$. Only a bounded number of
boundary interactions differ from the constant environment. Their contribution tends to zero as $m\to\infty$, leaving the second derivative of the local energy. Uniform convexity
of $V$ and positivity of the logarithmic curvature for $s\le s_0<1$ imply
\[
 cI_2\le\mathsf K_s(\theta)\le CI_2,
 \qquad s\in[0,s_0],\quad\theta\in[-\pi,\pi],
\]
with constants uniform in $n$ and depending on $V,\beta,s_0$. The scalar map $\mathcal L_{\mathrm{diag}}$ corresponds to $-b_\star^{-1}\ell$, where $\ell=(1,2)$. Its Gaussian
covariance spectral density is
\[
\widehat C_s(\theta)
=
\frac{2}{\beta n\,b_\star^2}
\ell\,\mathsf K_s(\theta)^{-1}\ell^*.
\eqtag{15.3}
\]
The factor $2/(\beta n)$ is the reciprocal of the prefactor in the conditional energy $(\beta n/2)\mathscr H_n$. The comparison with the original coefficient law is
Proposition~\ref{A:res:14.6}; Equation~\eqref{A:eq:15.3} specifies the covariance used in this Gaussian approximation. Finite range makes $\mathsf K_s(\theta)$ a trigonometric
polynomial. Its uniform positivity implies that $\widehat C_s(\theta)$ is real analytic and strictly positive.

The coefficient multiplying the Brownian noise is determined by the variance per unit length of the partial sums. For a centered stationary sequence $Z_j$ with exponentially
decaying covariance $c_r=\mathbb E[Z_0Z_r]$, one has
\[
 \frac1m\operatorname{Var}\!\left(\sum_{j=1}^m Z_j\right)
 =\sum_{|r|<m}\left(1-\frac{|r|}{m}\right)c_r
 \longrightarrow\sum_{r\in\mathbb Z}c_r.
\]
The limiting sum is the covariance symbol at frequency zero. For the Jacobi coefficients it is given by the following identity of \cite{krishnapur2016universality}. To compare
normalizations, use potential $V/2$ in that proposition: its local energy is $F_{V,s}/2$, so its inverse Hessian is $2\mathsf H_s^{-1}$. Set $\tau(s)=-E'(s)$, with
$\tau(0)=\tau_\star$.

\par\Needspace{8\baselineskip}
\begin{lemma}\label{A:res:15.1}
For every \(s\in[0,s_0]\),
\[
\ell\,\mathsf H_s^{-1}\ell^*
=
2b(s)\tau(s).
\eqtag{15.4}
\]
\end{lemma}

Consequently,
\[
\widehat C_s(0)
=
\frac{4}{\beta n b_\star^2}b(s)\tau(s).
\eqtag{15.5}
\]
In particular, at the edge,
\[
\widehat C_0(0)
=
\frac4\beta h^3,
\eqtag{15.6}
\]
because \(h^3=\frac{\tau_\star}{b_\star n}.\) Moreover, in the range used for the soft edge scaling \(s=O(h^2(1+x))\),
\[
\widehat C_s(0)
=
\frac4\beta h^3
\left[
1+O\bigl(h^2(1+x)\bigr)
\right].
\eqtag{15.7}
\]

\begin{proof}
Write $W=W_V$ for the Fourier average in \eqref{A:eq:11.4}, and evaluate its derivatives at $(a(s),b(s))$. The stationary equations are
\[
 W_a=0,\qquad W_b=\frac{2(1-s)}b.
\]
We first derive an identity between its second derivatives. Periodicity gives
\[
 \begin{aligned}
 0&=\frac1{2\pi}\int_0^{2\pi}
       \frac{\mathrm d}{\mathrm d\theta}
       \{\sin\theta\,V'(a+2b\cos\theta)\}\,\mathrm d\theta\\
  &=\frac1{2\pi}\int_0^{2\pi}
       \{\cos\theta V'(a+2b\cos\theta)
          -2b\sin^2\theta V''(a+2b\cos\theta)\}\,\mathrm d\theta.
 \end{aligned}
\]
Since $W_b=2\langle\cos\theta V'\rangle$, $W_{aa}=\langle V''\rangle$, and $W_{bb}=4\langle\cos^2\theta V''\rangle$, where the brackets denote the normalized integral, this implies
\[
 W_{bb}+b^{-1}W_b=4W_{aa}.
\]
The logarithmic term in $F_{V,s}$ contributes $2(1-s)/b^2=b^{-1}W_b$ to its second $b$ derivative and nothing to the other entries. Thus, with $A=W_{aa}$ and $C=W_{ab}$,
\[
 \mathsf H_s=\begin{pmatrix}A&C\\C&4A\end{pmatrix},
 \qquad \Delta=4A^2-C^2>0.
\]
Strict positivity follows from Lemma~\ref{A:res:11.1}. Differentiating the two stationary equations in $s$, justified by Proposition~\ref{A:res:11.3}, yields
\[
 \mathsf H_s\binom{a'}{b'}=\binom{0}{-2/b}.
\]
Consequently $a'=2C/(b\Delta)$ and $b'=-2A/(b\Delta)$, and direct matrix inversion gives
\[
 \mathsf H_s^{-1}
 =\frac1\Delta\begin{pmatrix}4A&-C\\-C&A\end{pmatrix}
 =-\frac b2\begin{pmatrix}4b'&a'\\a'&b'\end{pmatrix}.
\]
Multiplication by $\ell=(1,2)$ on both sides now gives
\[
 \ell\mathsf H_s^{-1}\ell^*=-2b(a'+2b')=2b\tau(s),
\]
because $E(s)=a(s)+2b(s)$ and $\tau(s)=-E'(s)$. This proves \eqref{A:eq:15.4} with the normalization of the present local energy.

Substitution in \eqref{A:eq:15.3} proves \eqref{A:eq:15.5}. At $s=0$, the equality $h^3=\tau_\star/(b_\star n)$ gives \eqref{A:eq:15.6}. Finally, the derivative bounds on the fixed
compact parameter interval imply
\[
 |b(s)\tau(s)-b_\star\tau_\star|
 \le s\sup_{0\le t\le s_0}|(b\tau)'(t)|.
\]
Since $b_\star\tau_\star>0$ and $s=O(h^2(1+x))$ in the stated range, division by $b_\star\tau_\star$ proves \eqref{A:eq:15.7}.
\end{proof}

The conditional Hessian acts on the full coefficient block, including the boundary neighborhoods removed in the covariance comparison. Write
\(K_{n,\mathrm{tr}}^\tau = \frac{\beta n}{2}\,\mathscr K_{n,\mathrm{tr}}^\tau\), where
\(\mathscr K_{n,\mathrm{tr}}^\tau\) is the unscaled Hessian. Its coordinate graph is a one-dimensional strip of fixed
width. We denote its graph distance by \(d(\cdot,\cdot)\), and set
\[
d_{\mathrm L}(p)=d(p,\partial_{\mathrm L}I_n^{\mathrm{tr}}),
\quad
d_{\mathrm R}(p)=d(p,\partial_{\mathrm R}I_n^{\mathrm{tr}}).
\]
For the Fourier calculation write each coefficient index as $p=(j,\sigma)$, where $j\in\mathbb Z$ and $\sigma\in\{A,B\}$. Set $\nu(p)=j$ and $\sigma(p)=\sigma$. Thus
$\nu(a_j)=\nu(b_j)=j$. The integer index $\nu$ groups the two coefficient types into one vector; it is distinct from the geometric position $\jmath(b_j)=j+1/2$ used in
Subsection~\ref{A:sec:the-local-minimizer-and-edge-scaling}. The two distance conventions differ by at most a fixed constant. For two coefficient indices put
\[
s_{pq}=s_{\lfloor(\nu(p)+\nu(q))/2\rfloor}.
\]
Since $|s_{j+1}-s_j|\le Ch^3$, changing from the geometric center to this integer center changes the local parameter by $O(h^3)$, which is included in the error below. For the
finite Hessian use the minimizer of the conditional potential on the transport interval, \(\overline Y_n^\circ\), on \(\mathcal T_n^{\mathrm{out}}(M_0)\). Lemma~\ref{A:res:12.1A}
bounds its free coefficients and bounds the regular off-diagonal ones away from zero. Proposition~\ref{A:res:11.3} gives the same bounds for the local profile. The segment joining
them therefore stays in a fixed bounded set with a positive lower bound on the off-diagonal coordinates. Evaluating the gradient of the unscaled conditional potential at the local
minimizer gives an \(O(n^{-1})\) error away from the two boundaries, by Proposition~\ref{A:res:11.7}, and an \(O(1)\) error within $R_V$ indices of either boundary. The averaged
Hessian is uniformly positive and vanishes outside a fixed band. Lemma~\ref{A:res:A.9}, or its special case with a bounded Hessian in the regular region, therefore gives
\[
 |\overline Y_{n,p}^\circ-Y_{n,p}^\dagger|
 \le C\left[n^{-1}+\mathrm{e}^{-cd_{\mathrm L}(p)}+\mathrm{e}^{-cd_{\mathrm R}(p)}\right].
\eqtag{finite-center}
\]
The error estimate is applied on a compact regular interval, with deterministic constants uniform over the admissible exterior data. To compare Hessians, put $Y^\dagger$ at the
local coefficients and use the fundamental theorem of calculus on the segment to $\overline Y_n^\circ$. For $d(p,q)\le R_2$,
\[
 |\partial_{pq}^2\mathscr H_n(\overline Y_n^\circ)
       -\partial_{pq}^2\mathscr H_n(Y^\dagger)|
 \le C\sum_{r:\,d(r,p)\le C R_V}
           |\overline Y_{n,r}^\circ-Y_{n,r}^\dagger|.
\]
There are a bounded number of summands. On the regular cube the third derivatives of the unscaled polynomial and logarithmic terms are bounded; fixed exterior coordinates only
change the constant. Next freeze the coefficients in $\partial_{pq}^2\mathscr H_n(Y^\dagger)$ at $s_{pq}$. The path expansion involves only a bounded number of nearby coefficients.
Their parameters differ from $s_{pq}$ by $O(n^{-1})$, so the polynomial Hessian changes by $O(n^{-1})$. On a free edge, the logarithmic Hessian is $2(1-s^b_{n,k})/b_k^2$. Its
derivatives in both arguments are bounded on the same regular cube, giving the same error for this diagonal term. Paths missing from the finite matrix or exterior values occur only
in the fixed width boundary rows, where their contribution is bounded. For $d(p,q)>R_2$ both entries are zero. Together with \eqref{A:eq:finite-center}, these calculations give

\[
\left| \mathscr K_{n,\mathrm{tr}}^\tau(p,q) - \bigl[\mathscr K_{s_{pq}}(\nu(p)-\nu(q))\bigr]_{\sigma(p)\sigma(q)} \right| \le C \left[ h^3+\mathrm{e}^{-c d_{\mathrm L}(p)}+\mathrm{e}^{-c d_{\mathrm R}(p)} \right] \mathbf 1_{\{d(p,q)\le R_2\}}.
\eqtag{15.11}
\]
Here \(R_2\) is fixed. Because \(p\) and \(q\) are at uniformly bounded distance whenever the left hand side is nonzero, one may interchange their boundary distances after changing
the constant. The two exponential terms bound the contributions of the boundary rows; the \(O(h^3)\) approximation of the local Hessian is used only at points separated from both
boundaries.

Apply Proposition~\ref{A:res:A.15} with $\Lambda_n=I_n^{\mathrm{tr}}$ and $\mathscr K_n=\mathscr K_{n,\mathrm{tr}}^\tau$, and use its inverse coefficients $q_s$ and approximate
inverse $Q_n$. There are two coefficient types at each site. Lemma~\ref{A:res:12.1A} puts all minimizing coefficients in a fixed cube, with every regular off-diagonal coefficient
at least $b_{\mathrm{reg}}>0$. The unscaled Hessian is the trace Hessian plus the diagonal entries $2(1-s^b_{n,k})/b_k^2$ on free edges. The trace Hessian lower bound and the
bounded finite path sums give
\[
 c_VI\le\mathscr K_{n,\mathrm{tr}}^\tau
 \le\bigl(C_{V,M_0}+2b_{\mathrm{reg}}^{-2}\bigr)I.
\]
Indeed, each row of $D^2\operatorname{Tr}V$ has a fixed number of entries bounded by $C_{V,M_0}$, and $0<1-s^b_{n,k}\le1$ by \eqref{A:eq:14.4}. The same bounds hold for the
translation-invariant local Hessians. The scaled Hessian is $K_{n,\mathrm{tr}}^\tau=(\beta n/2)\mathscr K_{n,\mathrm{tr}}^\tau$. The interaction range is $R_2$, and
\eqref{A:eq:15.11} verifies \eqref{A:eq:A.49}. The parameter increments satisfy $|s_{j+1}-s_j|\le1/n\le Ch^3$ by \eqref{A:eq:profile-extension}. The local coefficients and their
first derivatives are bounded on $[0,s_{\mathrm{reg}}]$ by Proposition~\ref{A:res:11.3}; their finite Fourier sums are analytic on a common strip. Proposition~\ref{A:res:A.15} now
gives \eqref{A:eq:C1}--\eqref{A:eq:C3}, uniformly over these exterior configurations.

Every coefficient coordinate entering \(\mathcal L_{\mathrm{diag}}^{\mathrm{mid}}\) is at graph distance at least
\(L_n^{\mathrm L}+R_V+4\) and \(L_n^{\mathrm R}+R_V+4\) from the two boundaries. Choose \(A_{\mathrm L},A_{\mathrm R}\)
so large that \(\mathrm{e}^{-cL_n^{\mathrm L}}+\mathrm{e}^{-cL_n^{\mathrm R}} \le h^3.\) Then, for all coefficient coordinates occurring in the middle scalar map,
\(\left| (\mathscr K_{n,\mathrm{tr}}^\tau)^{-1}(p,q)-Q_n(p,q) \right| \le C h^3\mathrm{e}^{-c d(p,q)}\). Multiplication by \(2/(\beta n)\asymp h^3\) gives
\[
\left|
(K_{n,\mathrm{tr}}^\tau)^{-1}(p,q)
-
\frac{2}{\beta n}Q_n(p,q)
\right|
\le
C h^6\mathrm{e}^{-c d(p,q)}.
\eqtag{15.17}
\]

For middle scalar sites define the matrix
\[
C_n^{\mathrm{aux}}(i,j)
=
\frac1{2\pi}\int_{-\pi}^{\pi}
 \mathrm{e}^{\mathrm i(i-j)\theta}\widehat C_{s_{ij}}(\theta)\,\mathrm{d}\theta,
\quad
s_{ij}=s_{\lfloor(i+j)/2\rfloor}.
\]
The linear functional at $i$ involves only the two coefficient types at integer index $i$. Hence the four inverse Hessian entries in its covariance with the functional at $j$ use
the same center $s_{ij}$. More explicitly,
\[
 C_n^\tau(i,j)=b_\star^{-2}
 \sum_{\sigma,\sigma'\in\{A,B\}}\ell_\sigma\ell_{\sigma'}
 (K_{n,\mathrm{tr}}^\tau)^{-1}((i,\sigma),(j,\sigma')).
\]
Replacing the inverse entries by $(2/(\beta n))Q_n$ gives exactly $C_n^{\mathrm{aux}}(i,j)$ by \eqref{A:eq:15.3}. The four errors in \eqref{A:eq:15.17} preserve exponential decay,
so
\[
|C_n^\tau(i,j)-C_n^{\mathrm{aux}}(i,j)|
\le
C h^6\mathrm{e}^{-c|i-j|},
\quad i,j\in I_n^{\mathrm{mid}}.
\eqtag{15.18}
\]
The matrix $C_n^{\mathrm{aux}}$ need not be positive. We use it to estimate the difference between $C_n^\tau$ and the positive covariance $\mathcal G_n\mathcal G_n^*$.

We verify the hypotheses of Proposition~\ref{A:res:A.18}, including its parameter norm. Write the finite-range Hessian symbol as
\[
 \mathsf K_s(z)=\sum_{|d|\le R}K_d(s)z^d,
 \quad K_d(s)\in\mathbb R^{2\times2},
 \quad K_{-d}(s)=K_d(s)^{\mathsf T}.
\]
In this Laurent convention $K_d(s)=\mathscr K_s(-d)$, so evaluation at $z=e^{\mathrm i\theta}$ agrees with the Fourier symbol used in Proposition~\ref{A:res:A.15}. On the fixed
interval $[0,s_0]$, the functions $a(s),b(s)$ are \(C^1\),
\(b(s)\ge b_0>0\), and their derivatives are bounded. These bounds follow from Proposition~\ref{A:res:11.3}. Each entry
of the Hessian of the trace term is a polynomial in \(a,b\). The additional logarithmic coefficient \(2(1-s)/b(s)^2\) has a bounded continuous derivative as well. Thus every
\(K_d\) is \(C^1\) with a common bound, including the one-sided derivatives at the two parameter boundary points.

The circle bound \(\mathsf K_s(\mathrm{e}^{\mathrm{i}\theta})\ge cI\) and finite range imply, for
\(z=t\mathrm{e}^{\mathrm{i}\theta}\) with \(\mathrm{e}^{-\rho}\le t\le \mathrm{e}^\rho\),
\[
 \|\mathsf K_s(z)-\mathsf K_s(\mathrm{e}^{\mathrm{i}\theta})\|
 \le \sum_{|d|\le R}\|K_d(s)\|\,|t^d-1|
 \le C\rho.
\]
Choose \(\rho>0\) so small that this is at most \(c/2\). The Neumann series gives $\mathsf K_s(z)^{-1}$ on this annulus, uniformly in $s$, with norm at most $2/c$, and
\[
 \partial_s\mathsf K_s(z)^{-1}
 =-\mathsf K_s(z)^{-1}(\partial_s\mathsf K_s(z))
                    \mathsf K_s(z)^{-1}.
\eqtag{symbol-analytic-derivative}
\]
The inverse identity, applied also to difference quotients, shows \(C^1\) dependence in the supremum norm for bounded holomorphic functions on any fixed smaller annulus, with
bounds uniform in $s$.

For a complex argument \(z\), define
\[
 q_s(z)=\frac{2}{\beta\tau_\star b_\star}
             \ell\,\mathsf K_s(z)^{-1}\ell^{\mathsf T}.
\eqtag{normalized-symbol}
\]
Then \(\widehat C_s(\theta)=h^3q_s(\mathrm{e}^{\mathrm{i}\theta})\), by \(h^3=\tau_\star/(b_\star n)\). Both \(q_s\) and
\(\partial_sq_s\) are uniformly bounded holomorphic functions on the fixed smaller annulus, and \(s\mapsto q_s\) is
\(C^1\) in that norm. On the circle, positivity and both spectral bounds for \(\mathsf K_s\) give \(0<c_1\le
q_s(\mathrm{e}^{\mathrm{i}\theta})\le C_1\). Moreover, real coefficients give \(\mathsf K_s(\mathrm{e}^{-\mathrm{i}\theta})=\overline{\mathsf K_s(\mathrm{e}^{\mathrm{i}\theta})}\).
The scalar quadratic form of the Hermitian inverse is real, so \(q_s(\mathrm{e}^{-\mathrm{i}\theta})=q_s(\mathrm{e}^{\mathrm{i}\theta})\). This proves evenness of the symbol.
Constants in these analytic and positivity estimates are uniform for \(\beta\) in a compact positive interval. Formula \eqref{A:eq:normalized-symbol} also makes the dependence on
\(\beta\) continuous.

Apply Proposition~\ref{A:res:A.18}. It gives real coefficients of the expansion \(g_m(s)\), extended by zero for \(m<0\), such that
\[
 |g_m(s)|+|\partial_sg_m(s)|\le C\mathrm{e}^{-c|m|},
\eqtag{15.21}
\]
and
\[
 q_s(\mathrm{e}^{\mathrm{i}\theta})=|g_s(\mathrm{e}^{\mathrm{i}\theta})|^2.
\eqtag{15.22}
\]
The factor is analytic on a disk of radius greater than one, is nonvanishing there, and has \(g_s(1)>0\) by its defining exponential. This fixes the phase of the factor.
Lemma~\ref{A:res:15.1} gives \(q_0(1)=4/\beta\); bounded
\(\partial_sq_s(1)\) gives \(|q_s(1)-4/\beta|\le Cs\). Hence
\[
 g_s(1)=\sqrt{q_s(1)}=\frac2{\sqrt\beta}+O(s).
\eqtag{15.23}
\]
Here $q_s(1)$ denotes angular frequency zero, whereas $(2\pi)^{-1}\int q_s(e^{\mathrm i\theta})\,\mathrm d\theta$ is the normalized variance at a single site. These need not agree.
For the quadratic potential $V(t)=t^2/2$, the Hessian at the constant coefficients has symbol $\operatorname{diag}(1,4)$ and $b_\star=\tau_\star=1$. In this case
$q_s\equiv4/\beta$, $g_s\equiv2/\sqrt\beta$, and $p_s=0$. Thus the covariance decomposition adds no derivative in this case; the separate term $\Gamma_{\eta_n}^-$ in the
coefficient decomposition is still present. Corollary~\ref{A:res:A.19} bounds the derivatives of $p_s$ for general $V$. Under the convolution convention in \eqref{A:eq:15.27},
multiplication by \(z\) shifts the driving sequence from $\xi_j$ to $\xi_{j-1}$. Use the function defined in \eqref{A:eq:A.59},
\[
p_s(z) = \frac{g_s(z)-g_s(1)}{1-z}.
\]
The numerator vanishes at $z=1$. Corollary~\ref{A:res:A.19} gives its analytic extension on a common disk of radius greater than one; in particular,
\[
p_s(z)=\sum_{m\in\mathbb Z}p_m(s)z^m,
\quad
|p_m(s)|+|\partial_sp_m(s)|
\le
C\mathrm{e}^{-c|m|},
\]
and
\[
g_s(z)=g_s(1)+(1-z)p_s(z).
\eqtag{15.26}
\]

The variance per site of long sums of the stationary Gaussian coefficients was computed above from the Hessian symbol. For varying $s_j$, we use the exact remainder
\eqref{A:eq:factor-remainder-formula} of Proposition~\ref{A:res:A.20}.

Define the slowly varying factor operator
\[
\mathcal G_n:\ell^2(\mathbb Z)\longrightarrow\mathbb R^{I_n^{\mathrm{mid}}},
\quad
(\mathcal G_nv)_j
=
h^{3/2}\sum_{m\in\mathbb Z}g_m(s_j)v_{j-m}.
\eqtag{15.27}
\]
The covariance of the Gaussian vector constructed with this operator is
\[
C_n^{\mathrm{fac}}
:=
\mathcal G_n\mathcal G_n^*,
\quad
C_n^{\mathrm{fac}}(i,j)
=
h^3\sum_{m\in\mathbb Z}
g_m(s_i)g_{m+j-i}(s_j).
\]
Positivity follows from $C_n^{\mathrm{fac}}=\mathcal G_n\mathcal G_n^*$, as in Proposition~\ref{A:res:A.18A}.

Equation~\eqref{A:eq:15.22} identifies $C_n^{\mathrm{aux}}$ with the matrix $C_n^{\mathrm{mid}}$ in Proposition~\ref{A:res:A.18A}, with the parameter at the site
$\lfloor(i+j)/2\rfloor$ in entry $(i,j)$. The parameter increments and coefficient derivative bounds have just been verified. Applying that proposition gives
\[
 |C_n^{\mathrm{fac}}(i,j)-C_n^{\mathrm{aux}}(i,j)|
 \le Ch^6\mathrm e^{-c|i-j|}.
 \eqtag{15.30}
\]
The differences in \eqref{A:eq:15.18} and \eqref{A:eq:15.30} both have order $h^6$: the variation of the parameter has order $h^3$, and the covariance normalization contributes
another factor $h^3$. Both comparisons use $C_n^{\mathrm{aux}}$ as the intermediate matrix. Combining \eqref{A:eq:15.18} and \eqref{A:eq:15.30},
\[
|C_n^\tau(i,j)-C_n^{\mathrm{fac}}(i,j)|
\le
C h^6\mathrm{e}^{-c|i-j|}.
\eqtag{15.31}
\]

For \(\tau\in\mathcal T_n^{\mathrm{out}}(M_0)\), let \(\zeta_n^{\mathrm{fac}}\) be centered Gaussian with covariance
\(C_n^{\mathrm{fac}}=\mathcal G_n\mathcal G_n^*\), and let \((\widehat\zeta_j^{\mathrm G})_{j\in I_n^{\mathrm{mid}}}\)
be the Gaussian vector with conditional covariance $C_n^\tau$.

\par\Needspace{8\baselineskip}
\begin{proposition}\label{A:res:15.2}
These vectors admit a jointly Gaussian coupling such that, for every \(p<\infty\),
\[
\left\|
\sum_{j\in I_n^{\mathrm{mid}}}
(1+x_j)^{-7/2}
|\widehat\zeta_j^{\mathrm G}-\zeta_j^{\mathrm{fac}}|^2
\right\|_{L^p(\mathbb P_\tau)}
\le
C_ph^5.
\eqtag{15.32}
\]
\end{proposition}

Consequently, almost surely, for all sufficiently large \(n\),
\[
\sum_{j\in I_n^{\mathrm{mid}}}
(1+x_j)^{-7/2}
|\widehat\zeta_j^{\mathrm G}-\zeta_j^{\mathrm{fac}}|^2
\le
C h^4.
\eqtag{15.33}
\]

\begin{proof}
Let \(W=\operatorname{diag}\bigl((1+x_j)^{-7/2}\bigr)_{j\in I_n^{\mathrm{mid}}}.\) Both \(C_n^\tau\) and
\(C_n^{\mathrm{fac}}\) are positive semidefinite. Apply Lemma~\ref{A:res:A.16} with weight \(W\). It gives
\[
\mathbb E_\tau
\sum_j(1+x_j)^{-7/2}
|\widehat\zeta_j^{\mathrm G}-\zeta_j^{\mathrm{fac}}|^2
\le
\left\|
W^{1/2}(C_n^\tau-C_n^{\mathrm{fac}})W^{1/2}
\right\|_{\mathfrak S_1}.
\eqtag{15.34}
\]
Use the matrix unit trace bound and weighted exponential summation proved in Corollary~\ref{A:res:A.17}. With \eqref{A:eq:15.31}, they give
\[
\left\| W^{1/2}(C_n^\tau-C_n^{\mathrm{fac}})W^{1/2} \right\|_{\mathfrak S_1} \le C h^6 \sum_{i,j}W_i^{1/2}W_j^{1/2}\mathrm{e}^{-c|i-j|} \le C h^6\sum_iW_i \le C h^5.
\eqtag{15.35}
\]
Conditional on the exterior, the difference is Gaussian. Applying \eqref{A:eq:gaussian-weighted-quadratic-moment} to its weighted squared norm, and Jensen's inequality when
$0<p<1$, gives
\[
 \left\|\sum_jW_j
   |\widehat\zeta_j^{\mathrm G}-\zeta_j^{\mathrm{fac}}|^2
       \right\|_{L^p(\mathbb P_\tau)}
 \le C_p\mathbb E_\tau\sum_jW_j
   |\widehat\zeta_j^{\mathrm G}-\zeta_j^{\mathrm{fac}}|^2
 \le C_ph^5.
\]
This proves \eqref{A:eq:15.32} uniformly on $\mathcal T_n^{\mathrm{out}}(M_0)$. For exterior data outside that set, retain the given Gaussian vector $\zeta_n^{\mathrm{fac}}$ and
sample $\widehat\zeta^{\mathrm G}$

independently of it, conditional on the exterior data, with covariance $C_n^\tau$ and mean zero. Both covariances have diagonal entries at most $Ch^3$, the first by the
unconditional Hessian lower bound. For this conditionally independent pair, the estimate for Gaussian quadratic forms gives
\[
 \left\|\sum_jW_j|\widehat\zeta_j^{\mathrm G}
                 -\zeta_j^{\mathrm{fac}}|^2\right\|_{L^p(\mathbb P_\tau)}
 \le C_p\sum_jW_j\bigl(C_n^\tau(j,j)+C_n^{\mathrm{fac}}(j,j)\bigr)
 \le C_ph^3\sum_jW_j\le C_ph^2.
\]
If $\nu_n$ is the actual exterior marginal, disintegration gives
\[
 \mathbb E\left(\sum_jW_j
   |\widehat\zeta_j^{\mathrm G}-\zeta_j^{\mathrm{fac}}|^2\right)^p
 \le C_ph^{5p}+C_ph^{2p}\nu_n(\mathcal T_n^{\mathrm{out}}(M_0)^c)
 \le C_ph^{5p}.
\]
Here the last inequality uses the exponentially small exterior probability. The vector $\widehat\zeta^{\mathrm G}$ need not be independent of the exterior. In particular, Markov's
inequality on the admissible exterior configurations and \eqref{A:eq:integrated-tail} give
\[
 \mathbb P\left\{\sum_j(1+x_j)^{-7/2}
       |\widehat\zeta_j^{\mathrm G}-\zeta_j^{\mathrm{fac}}|^2>h^4\right\}
 \le C_ph^p+C\mathrm{e}^{-cn}.
\]
For \(p>3\), the right side is summable in \(n\), and Borel--Cantelli proves \eqref{A:eq:15.33}.

The weighted estimate here is for this particular Gaussian coupling. It does not give a coordinatewise difference of order $h^3$. Writing $W_j=(1+x_j)^{-7/2}$, it gives
$W_j\mathbb E_\tau|\varepsilon_j^{\mathrm{fr}}|^2\le Ch^5$; the two marginal variances separately give $\mathbb E_\tau|\varepsilon_j^{\mathrm{fr}}|^2\le Ch^3$. Both bounds are
needed to verify the form hypotheses below. The resolvent comparison then uses the weighted sum itself.

\end{proof}

For the common form bound we need a maximum estimate for the potential relative to the confining weight. We combine the weighted covariance estimate with a rough coordinate
variance bound, using each for the range of indices where it is stronger. On the conditional Gaussian coupling of Proposition~\ref{A:res:15.2}, put \(X_j=1+x_j\) and
\(\varepsilon_j^{\mathrm{fr}} =\widehat\zeta_j^{\mathrm G}-\zeta_j^{\mathrm{fac}}\).

\par\Needspace{8\baselineskip}
\begin{lemma}\label{A:res:freezing-max} For admissible exterior data,
\[
 \operatorname{Var}_\tau\!\left(
       \frac{h^{-2}\varepsilon_j^{\mathrm{fr}}}{X_j}\right)
 \le C\min\{hX_j^{3/2},\,h^{-1}X_j^{-2}\}
 \le Ch^{1/7}.
\eqtag{freezing-variance}
\]
\end{lemma}

Consequently, for each \(Q>0\) there is a deterministic \(C_Q\) such that
\[
 \mathbb P\!\left(
   \max_{j\in I_n^{\mathrm{mid}}}
       \frac{h^{-2}|\varepsilon_j^{\mathrm{fr}}|}{X_j}
       > C_Qh^{1/14}\sqrt{\log(n+1)}\right)
 \le C_Qn^{-Q}+C\mathrm{e}^{-cn}.
\eqtag{freezing-max-tail}
\]
In particular the maximum is \(O_\omega(h^{1/14}\sqrt{\log(n+1)})=o_\omega(1)\) for all sufficiently large \(n\) on the final common coupling. All constants may be chosen uniformly
on compact positive \(\beta\)-intervals.

\begin{proof}
Both conditional marginals have coordinate variances bounded by \(Ch^3\): for the Gaussian block this follows from the Hessian lower bound and the bounded local scalar map, and for
the factor block from the exponentially summable coefficients of \(\mathcal G_n\). The inequality \(|a-b|^2\le2|a|^2+2|b|^2\) gives
\(\operatorname{Var}_\tau\varepsilon_j^{\mathrm{fr}}\le Ch^3\) for any coupling of the marginals. On the selected
coupling, the nonnegative summands in \eqref{A:eq:15.34}--\eqref{A:eq:15.35} also give
\[
 X_j^{-7/2}\operatorname{Var}_\tau
       \varepsilon_j^{\mathrm{fr}}
 \le\mathbb E_\tau\sum_iX_i^{-7/2}
                   |\varepsilon_i^{\mathrm{fr}}|^2
 \le Ch^5.
\]
Multiply the two variance bounds by \(h^{-4}X_j^{-2}\). The increasing bound $ChX_j^{3/2}$ and the decreasing bound $Ch^{-1}X_j^{-2}$ meet at $X_j=h^{-4/7}$. For
$1\le X_j\le h^{-4/7}$, the first bound is at most
\(h(h^{-4/7})^{3/2}=h^{1/7}\); for \(X_j\ge h^{-4/7}\), the second is at most \(h^{-1}(h^{-4/7})^{-2}=h^{1/7}\). This
proves \eqref{A:eq:freezing-variance}.

Conditionally on the exterior data the difference is a centered Gaussian vector. The union bound gives, for every
\(u>0\),
\[
 \mathbb P_\tau\!\left(
   \max_j\frac{h^{-2}|\varepsilon_j^{\mathrm{fr}}|}{X_j}>u\right)
 \le 2n\exp\!\left(-\frac{u^2}{2Ch^{1/7}}\right).
\eqtag{freezing-union}
\]
Take \(u=C_Qh^{1/14}\sqrt{\log(n+1)}\) with \(C_Q^2/(2C)>Q+1\). Integrating over admissible exterior data and adding
\(C\mathrm{e}^{-cn}\) proves \eqref{A:eq:freezing-max-tail}. The Gaussian tail bound is conditional on the exterior; averaging over the exterior need not give a Gaussian distribution. With \(Q>1\), Borel--Cantelli gives the bound for all sufficiently large matrix sizes, almost surely. Lemma~\ref{A:res:A.16} preserves the joint Gaussian distribution conditional on the exterior
while keeping the factor vector fixed. The estimate therefore also holds when that vector is generated by the common Brownian path.
\end{proof}

We construct the coefficient vectors conditionally on a common Brownian motion. Sample the exterior vector $\tau$ from its original marginal, independently of one two-sided
Brownian motion $B$. Conditional on $\tau$, the normalized increments of $B$ are still independent standard normals. Define $\zeta^{\mathrm{fac}}$ from these increments as below,
and use Lemma~\ref{A:res:A.16} conditional on this vector. Let $L$ map the full coefficient vector to the middle scalar sequence. Then
\[
 \Sigma=(K_{n,\mathrm{tr}}^\tau)^{-1},\qquad
 C_n^\tau=L\Sigma L^*,\qquad
 \widehat\zeta^{\mathrm G}\sim N(0,C_n^\tau).
\]
Lemma~\ref{A:res:A.16A} supplies a Gaussian vector \(G\sim N(0,\Sigma)\) such that, almost surely,
\[
 LG=\widehat\zeta^{\mathrm G}.
\]
On $\mathcal T_n^{\mathrm{out}}(M_0)$, transport $G$ as in Section~\ref{A:sec:gaussian-approximation-of-the-conditional-gibbs-measure}, then apply Lemma~\ref{A:res:A.11D} to
recover the conditional law on the positive coefficient domain. On the complement, sample that law independently of the Gaussian variables, as specified after
Proposition~\ref{A:res:14.6}. Neither operation changes $\zeta^{\mathrm{fac}}$. If $\mu_{n,+}^\tau$ denotes the conditional law of $Y_n^{\mathrm{tr}}$ and $\nu_n$ the original
exterior marginal, the completed coefficient vector satisfies, for every bounded Borel $F$,
\[
 \mathbb E F(Y_n^{\mathrm{tr}},\tau)
 =\int\!\int F(y,\tau)\,\mu_{n,+}^\tau(\mathrm dy)\,\nu_n(\mathrm d\tau).
\]
This is the original joint coefficient law, including its dependence across the cut. Using these kernels and a countable family of auxiliary variables, we sample the coefficient
arrays for all sizes over the same $B$, as detailed in Subsection~\ref{A:sec:measurable-couplings-of-the-finite-models}. The scalar pair is jointly Gaussian conditionally on
$\tau$. Averaging this conditional distribution over $\tau$ need not give a Gaussian distribution. The recovered coefficient vector has the Gibbs law.

For the two-sided Brownian motion just chosen, put $\xi_j=h^{-1/2}(B(jh)-B((j-1)h))$, $j\in\mathbb Z$. For each $n$ these are independent $N(0,1)$ variables, and
\[
\zeta_j^{\mathrm{fac}}
=
(\mathcal G_n\xi)_j
=
h^{3/2}\sum_{m\in\mathbb Z}g_m(s_j)\xi_{j-m}
\]
has covariance $C_n^{\mathrm{fac}}$. These coordinatewise convergent series define $\mathcal G_n\xi$. Only indices $j\ge1$ enter the finite difference operator driven by the
Brownian increments. The independent negative time half of $B$ supplies the remaining increments with $j\le0$ for the convolution. By exponential coefficient decay, their
contribution to a row $j\ge1$ has variance at most $Ch^3\sum_{m\ge j}e^{-2cm}\le Ch^3e^{-2cj}$. We keep these increments; deleting them would change the factor covariance near the
left boundary. Define \(\zeta_j^B = \frac{2}{\sqrt\beta}h^{3/2}\xi_j = \frac{2}{\sqrt\beta} h\bigl(B(jh)-B((j-1)h)\bigr)\). Proposition~\ref{A:res:15.2} gives
\eqref{A:eq:15.32}--\eqref{A:eq:15.33} for the covariance error $\varepsilon^{\mathrm{fr}}$. The index interval here is \(\{J_n^{\mathrm L}+1,\ldots,J_n^{\mathrm{reg}}\}\), with
the additional parameter \(s_{J_n^{\mathrm L}}\) to define $e_{J_n^{\mathrm L}}$. Equations \eqref{A:eq:profile-size-increment} and \eqref{A:eq:profile-extension} verify both
hypotheses on the coordinate sequence in Proposition~\ref{A:res:A.20}, including this predecessor. Proposition~\ref{A:res:A.20} applies with these indices. Its sequences are
\[
 \begin{aligned}
e_j& =h^{1/2}\sum_mp_m(s_j)\xi_{j-m},\qquad \rho_j^{(0)} =h^{3/2}\bigl(g_{s_j}(1)-2/\sqrt\beta\bigr)\xi_j,\\
\rho_j^{(1)}& =h^{3/2}\sum_m \bigl[p_{m-1}(s_{j-1})-p_{m-1}(s_j)\bigr]\xi_{j-m}, \qquad \rho_j=\rho_j^{(0)}+\rho_j^{(1)}.
\end{aligned}
\]
The absolute convergence and reindexing were proved in \eqref{A:eq:factor-remainder-formula}. The definitions give the decomposition
\[
\zeta_j^{\mathrm{fac}}
=
\zeta_j^B
+
h(e_j-e_{j-1})
+
\rho_j.
\eqtag{15.37}
\]

Proposition~\ref{A:res:A.20} applies on this interval with the inherited weights. Its two coordinate variance bounds are $\mathbb E|\rho_j^{(0)}|^2\le Ch^7X_j^2$ and
$\mathbb E|\rho_j^{(1)}|^2\le Ch^9$. The first uses the size of $s_j$; the second uses its increment. After weighted summation, \eqref{A:eq:factor-moment-bounds} gives
\[
 \left\|\sum_jX_j^{-7/2}|\rho_j|^2\right\|_{L^p}\le C_ph^6.
\]
Markov at threshold $h^4$ has probability at most $C_ph^{2p}$, summable for $p>3/2$. Hence
\[
 \sum_j(1+x_j)^{-7/2}|\rho_j|^2\le C_\omega h^4
 \eqtag{15.41}
\]
for all sufficiently large $n$, almost surely.

For $\mathcal E_e=\sum_jhX_j^{-3/2}|e_j|^2$, use \eqref{A:eq:factor-primitive-covariance}: the weighted covariance has trace $O(h)$ and operator norm $O(h^2)$. The exponential
estimate \eqref{A:eq:factor-primitive-tail} gives $\mathbb P(\mathcal E_e>Ah)\le e^{-cA/h}$ for a sufficiently large fixed $A$. Borel--Cantelli therefore yields
\[
 \sum_jh(1+x_j)^{-3/2}|e_j|^2\le C_\omega h.
 \eqtag{15.45}
\]
Both conclusions hold on the same probability space carrying the Brownian motion used by the covariance coupling; they require no independence between errors or between sizes.

The first-order form estimate will be applied after extending $e$ constantly outside $I_n^{\mathrm{mid}}$. We specify this extension below, after collecting the Gaussian estimates.

There exists a coupling of the Gaussian vector \((\widehat\zeta_j^{\mathrm G})_{j\in I_n^{\mathrm{mid}}}\) with a Brownian motion \(B\), a Gaussian sequence \(e=(e_j)\), the
remainder \(\rho=(\rho_j)\) from spectral factorization, and the remainder \(\varepsilon^{\mathrm{fr}}=(\varepsilon_j^{\mathrm{fr}})\) from the covariance approximation. On the
admissible set of exterior data, use the construction above; on its complement, take the product of the marginals and define
$\varepsilon^{\mathrm{fr}}= \widehat\zeta^{\mathrm G}-\zeta^{\mathrm{fac}}$. Write
$\zeta_j^B=(2/\sqrt\beta)h(B(jh)-B((j-1)h))$.

\begin{theorem}\label{A:res:15.3}
For every \(j\in I_n^{\mathrm{mid}}\),
\[
\widehat\zeta_j^{\mathrm G}
=
\zeta_j^B
+
h(e_j-e_{j-1})
+
\rho_j
+
\varepsilon_j^{\mathrm{fr}}.
\eqtag{15.47}
\]
Almost surely, the weighted estimates \eqref{A:eq:15.48}--\eqref{A:eq:15.50} hold for all sufficiently large $n$.
\end{theorem}

The estimates are
\[
\sum_{j\in I_n^{\mathrm{mid}}}
(1+x_j)^{-7/2}|\varepsilon_j^{\mathrm{fr}}|^2
\le
C_\omega h^4,
\eqtag{15.48}
\]
\[
\sum_{j\in I_n^{\mathrm{mid}}}
h(1+x_j)^{-3/2}|e_j|^2
\le
C_\omega h,
\eqtag{15.49}
\]
and
\[
\sum_{j\in I_n^{\mathrm{mid}}}
(1+x_j)^{-7/2}|\rho_j|^2
\le
C_\omega h^4.
\eqtag{15.50}
\]

\begin{proof}
For every middle index, \eqref{A:eq:15.37} and the definition of $\varepsilon_j^{\mathrm{fr}}$ give
\[
 \widehat\zeta_j^{\mathrm G}-\zeta_j^B
 =(\zeta_j^{\mathrm{fac}}-\zeta_j^B)+\varepsilon_j^{\mathrm{fr}}
 =h(e_j-e_{j-1})+\rho_j+\varepsilon_j^{\mathrm{fr}}.
\]
The definition of $e$ includes the predecessor at the first middle site, so this proves \eqref{A:eq:15.47} at every stated index. Proposition~\ref{A:res:15.2} gives
\eqref{A:eq:15.48}; \eqref{A:eq:15.45} and \eqref{A:eq:15.41} give \eqref{A:eq:15.49} and \eqref{A:eq:15.50}. Intersecting the finitely many probability one events gives all three
estimates for the same coupling and all sufficiently large sizes.
\end{proof}

Adding the transport error from Subsection~\ref{A:sec:gaussian-approximation-of-the-conditional-gibbs-measure} to \eqref{A:eq:15.47} now expresses the original scalar coefficient,
for $j\in I_n^{\mathrm{mid}}$, as
\[
\zeta_j={} \zeta_j^B+h(e_j-e_{j-1})+\rho_j +(\zeta_j-\widehat\zeta_j^{\mathrm G})
 +(\widehat\zeta_j^{\mathrm G}-\zeta_j^{\mathrm{fac}}).
\]
The last two terms are bounded by Theorem~\ref{A:res:14.7} and Proposition~\ref{A:res:15.2}, respectively. The term $\Gamma_{\eta_n}^-$ in Proposition~\ref{A:res:13.1} is
unchanged.

We restrict the derivative to the middle interval by extending its primitive constantly at both ends, using the predecessor value at the left cut. This preserves the first middle
difference and makes the differences vanish beyond the last middle site. Put \(j_n^-=J_n^{\mathrm L}+1, \allowbreak\, j_n^+=J_n^{\mathrm{reg}}.\) The convolution also defines
\(e_j\) at \(j=j_n^--1\). Set
\[
\widetilde e_j
=
\begin{cases}
e_{j_n^--1},&0\le j<j_n^--1,\\
e_j,&j_n^--1\le j\le j_n^+,\\
e_{j_n^+},&j_n^+<j\le M,
\end{cases}
\]
and
\[
\widetilde\rho_j
=
\chi_{j,n}^{\mathrm{mid}}\rho_j,
\quad
\widetilde\varepsilon_j^{\mathrm{fr}}
=
\chi_{j,n}^{\mathrm{mid}}\varepsilon_j^{\mathrm{fr}}.
\]
The extended sequence satisfies \(\widetilde e_j-\widetilde e_{j-1}=0 \allowbreak\; (j\in I_n^{\mathrm L}\cup I_n^{\mathrm{far}})\).

The vector being paired with the form vanishes at $0$ and $M+1$, so its boundary contributions are zero for arbitrary values of $\widetilde e_0$ and $\widetilde e_M$.

For $u_0=u_{M+1}=0$, summation over the full interval gives
\[
\sum_{j=1}^M\frac{\widetilde e_j-\widetilde e_{j-1}}h|u_j|^2 =-\sum_{j=0}^M\widetilde e_j
 \frac{|u_{j+1}|^2-|u_j|^2}h =-2\operatorname{Re}\sum_{j=0}^M \widetilde e_j(D_hu)_j
 \overline{\frac{u_{j+1}+u_j}{2}}.
\]
Here the second equality uses $|v|^2-|w|^2=2\operatorname{Re}((v-w)\overline{(v+w)/2})$. Consequently the diagonal multiplier on the left is
\[
\Gamma_{\widetilde e}
=-\bigl(D_h^*M_{\widetilde e}J_h
        +J_h^*M_{\widetilde e}D_h\bigr),\quad
(J_hu)_j=\frac{u_{j+1}+u_j}{2}.
\]
This identity uses the center map, which differs from $J_-$ in the off-diagonal decomposition: $\Gamma_{\widetilde e}$ is a diagonal operator represented by a first-order form,
whereas $\Gamma_{\eta}^{-}$ generally retains off-diagonal entries. Both are estimated by the first-order bound of
Subsection~\ref{A:sec:weighted-resolvent-estimates-and-comparison-principles}. The remaining constant portions of $\widetilde e$ are included in the weighted norm below. Put
$a=j_n^-$ and $b=j_n^+$. The extra constant pieces contribute at most
\[
 C ha\,|e_{a-1}|^2+C h\,|e_b|^2,
\]
since $a=O(\log n)$ and $\sum_{j>b}h(1+x_j)^{-3/2}\le C(1+x_b)^{-1/2}\le Ch$. Both variables have variance at most $Ch$. For each $Q>0$ the Gaussian tail bound gives
\[
 \mathbb P\!\left(|e_{a-1}|^2+|e_b|^2>C_Qh\log(n+1)\right)
 \le C_Qn^{-Q}.
\]
On the complementary event the two contributions are bounded by $C_Qh^2(\log(n+1))^2+C_Qh^2\log(n+1)=o(h)$. Combine this with \eqref{A:eq:15.45} and Borel--Cantelli, taking $Q>1$,
to obtain
\[
 \sum_{j=0}^Mh(1+x_j)^{-3/2}|\widetilde e_j|^2\le C_\omega h.
\eqtag{15.51}
\]
Since the two Gaussian variances are at most $Ch$, their squared variables have $L^p$ norms at most $C_ph$. Multiplying by the deterministic weights above bounds the extra pieces
by $C_ph^2\log(n+1)$ and $C_ph^2$, respectively. Thus the whole weighted sum has $L^p$ norm at most $C_ph$. The edge $j=0$ causes no change, since the convention $x_0=0$ gives
$1+x_0=1$. By the definitions of \(\zeta^{B,\mathrm{mid}}\), \(\widetilde e\), \(\widetilde\rho\), and
\(\widetilde\varepsilon^{\mathrm{fr}}\), equation~\eqref{A:eq:15.47} becomes the identity on the full reduced interval
\[
\zeta_j^{\mathrm G}
=
\zeta_j^{B,\mathrm{mid}}
+h(\widetilde e_j-\widetilde e_{j-1})
+\widetilde\rho_j
+\widetilde\varepsilon_j^{\mathrm{fr}},
\quad 1\le j\le M.
\eqtag{15.52}
\]

The form estimates also require coordinate maxima. These follow from the Gaussian variances and a union bound at each matrix size.

\par\Needspace{8\baselineskip}
\begin{lemma}\label{A:res:factor-max} Under the coupling just constructed,
for every \(Q>0\), outside an event of probability at most \(C_Qn^{-Q}\),
\[
 \begin{aligned}
\|\widetilde e_n\|_\infty& \le C_Qh^{1/2}\sqrt{\log(n+1)},\\
\max_{j\in I_n^{\mathrm{mid}}} \frac{h^{-2}|\rho_j|}{1+x_j}& \le C_Qh^{3/2}\sqrt{\log(n+1)},\qquad \max_{1\le j\le M}|\zeta_j^B| \le C_Qh^{3/2}\sqrt{\log(n+1)}.
\end{aligned}
\eqtag{factor-maxima}
\]
\end{lemma}

On intersecting with the coefficient events of Section~\ref{A:sec:boundary-reduction}, with the same probability order,
\[
 \|\eta_n\|_\infty\le C_Qh^{1/2}\sqrt{\log(n+1)},
 \quad
 \max_{1\le j\le M}|\zeta_j|
       \le C_Qh^{3/2}\sqrt{\log(n+1)}.
\eqtag{exact-maxima}
\]
All these estimates hold almost surely for all sufficiently large \(n\) on the final coupling, and their finite size constants are uniform on compact positive \(\beta\)-intervals.

\begin{proof}
The convolution formula gives \(\mathbb E|e_j|^2=h\sum_m|p_m(s_j)|^2\le Ch\). This includes the one additional left value used in the constant extension. That extension creates no
new values, so a union bound uses at most \(n+2\) Gaussian coordinates. For the factor remainder, the exact recentering identity \eqref{A:eq:recentered-profile-coordinate} gives
\[
 s_j\le Ch^2(1+x_j),\quad |s_j-s_{j-1}|\le Ch^3,
\]
for large \(n\), since \(rh\to0\). The coordinate variance estimates \eqref{A:eq:factor-remainder-variances} therefore imply
\[
 \mathbb E|\rho_j|^2
 \le 2\mathbb E|\rho_j^{(0)}|^2
       +2\mathbb E|\rho_j^{(1)}|^2
 \le C\bigl(h^7(1+x_j)^2+h^9\bigr).
\]
After dividing by \(h^4(1+x_j)^2\), the variance is at most \(Ch^3\). The Brownian scalar coordinates have variance $4h^3/\beta$. Every variable in the three lists
\[
 h^{-1/2}e_j,\qquad h^{-7/2}(1+x_j)^{-1}\rho_j,\qquad
 h^{-3/2}\zeta_j^B
\]
has variance at most a fixed $C_0$. There are at most $3(n+2)$ entries. For $t>0$, the union bound for their one coordinate Gaussian tails gives
\[
 \mathbb P\!\left(\max_{Z\text{ in these lists}}|Z|>t\right)
 \le6(n+2)\exp\{-t^2/(2C_0)\}.
\]
Choose $t=C_Q\sqrt{\log(n+1)}$ with $C_Q^2/(2C_0)>Q+2$ to obtain \eqref{A:eq:factor-maxima}.

For the exact coefficients, use the finite size bound on the maximum after centering at the minimizer from \eqref{A:eq:mode-maximum}, including its rough terminal bias, and
integrate the bad exterior probability as in Lemma~\ref{A:res:disintegration}. The fixed local scalar map gives \(\zeta\), and division of the off-diagonal coefficient fluctuation
by \(b_\star h\) gives \(\eta\). The auxiliary boundary coefficients are identical in the actual and centered Dirichlet representations, so \(\delta_0=\delta_M=0\). This proves
\eqref{A:eq:exact-maxima}. Finally choose \(Q>1\) and apply Borel--Cantelli to the finite collection of failure events.
\end{proof}

Return to $H_n^{B,\mathrm{mid}}$ defined at the start of \hyperref[A:sec:application-to-the-conditional-jacobi-coefficients]{the application to conditional distributions of Jacobi
coefficients}. Multiplying the coefficient identity \eqref{A:eq:15.52} by $h^{-2}$ turns its difference term into $h^{-1}(\widetilde e_j-\widetilde e_{j-1})$. The summation by
parts identity above identifies this multiplier with $\Gamma_{\widetilde e}$, so
\[
H_n^{\mathrm G}
=
H_n^{B,\mathrm{mid}}
+
\Gamma_{\widetilde e}
+
h^{-2}M_{\widetilde\rho+\widetilde\varepsilon^{\mathrm{fr}}}.
\eqtag{15.54}
\]
Set \(\widetilde H_n = H_n^{B,\mathrm{mid}}+\Gamma_{\widetilde e}.\) The coefficient estimates imply the following comparison once the common form bounds of
Theorem~\ref{A:res:16.2} are established.

\begin{theorem}\label{A:res:15.4}
For every compact \(K\Subset\mathbb C\setminus\mathbb R,\) almost surely, for all sufficiently large \(n\),
\[
\sup_{z\in K}
\left\|
(H_n^{\mathrm G}-z)^{-1}
-
(H_n^{B,\mathrm{mid}}-z)^{-1}
\right\|_{\mathfrak S_2}^2
\le
C_{\omega,V,\beta,K}h.
\]
\end{theorem}

It remains to replace the scalar coefficients outside the middle interval. Near the left boundary the estimate for the resolvent columns controls this replacement; on the far
interval the weight in the diagonal multiplier estimate gives the required bound. Set $H_n^{B,\circ}=H_n^\circ+h^{-2}M_{\zeta_n^B}$. The difference
\[
H_n^{B,\mathrm{mid}}-H_n^{B,\circ}
=
h^{-2}M_{\mathbf 1_{I_n^{\mathrm L}\cup I_n^{\mathrm{far}}}
(\zeta_n-\zeta_n^B)}
\eqtag{15.57a}
\]
is estimated separately on the two regions.

The common form proof starts from $H_n^{B,\circ}$ and the coefficient maxima below. The proofs of Theorems~\ref{A:res:15.5}--\ref{A:res:15.6} are placed after that proof.

We also need two bounds that hold on the entire reduced interval. For the exact KRV sequence of diagonal coefficients, Lemma~\ref{A:res:factor-max}, including its rough terminal
mean bound, gives \(\max_{1\le j\le M} \left( |A_j-A_j^\circ|+|B_j-B_j^\circ| \right) \le C_\omega h^{3/2}\sqrt{\log n}\) for all sufficiently large \(n\). Since
\(\zeta_j = -b_\star^{-1} \left( A_j-A_j^\circ+2(B_j-B_j^\circ) \right),\) we obtain
\[
\max_j|\zeta_j|
\le
C_\omega h^{3/2}\sqrt{\log n}.
\eqtag{15.58}
\]
For the Brownian increment sequence, the variables \(\zeta_j^B\) are centered Gaussian with variance \(4h^3/\beta\). Since \(M=O(h^{-3})\), the Gaussian union bound gives
\(\mathbb P \left( \max_{1\le j\le M}|\zeta_j^B| >A h^{3/2}\sqrt{\log(h^{-1})} \right) \le C h^{cA^2-3}\). Choosing \(A\) large and applying Borel--Cantelli yields
\[
\max_j|\zeta_j^B|
\le
C_{\omega,\beta}h^{3/2}\sqrt{\log n}
\eqtag{15.59}
\]
for all sufficiently large \(n\). Therefore
\[
\max_{1\le j\le M}
\bigl(|\zeta_j|+|\zeta_j^B|\bigr)
\le
C_{\omega,\beta}h^{3/2}\sqrt{\log n}.
\eqtag{15.60}
\]

Define the coefficient sequence
\[
\zeta_{j,n}^{B,\mathrm{reg}}
=
\begin{cases}
\zeta_{j,n}^B,&1\le j\le J_n^{\mathrm{reg}},\\
\zeta_j,&j\in I_n^{\mathrm{far}},
\end{cases}
\]
and \(H_n^{B,\mathrm{reg}} = H_n^\circ+h^{-2}M_{\zeta_n^{B,\mathrm{reg}}}.\) Thus \(H_n^{B,\mathrm{mid}}\) and
\(H_n^{B,\mathrm{reg}}\) differ only on \(I_n^{\mathrm L}\). By \eqref{A:eq:15.58}--\eqref{A:eq:15.59}, \(\max_{j\in
I_n^{\mathrm L}} |\zeta_j-\zeta_j^B| \le C_\omega h^{3/2}\sqrt{\log n}\) for all sufficiently large \(n\). Hence
\[
W_n^{\mathrm L,sc}
:=
H_n^{B,\mathrm{mid}}-H_n^{B,\mathrm{reg}}
=
h^{-2}M_{\mathbf 1_{I_n^{\mathrm L}}(\zeta_n-\zeta_n^B)}
\]
is a diagonal matrix supported on \(J_n^{\mathrm L}=O(\log n)\) sites and satisfies \(\max_{i,j}|(W_n^{\mathrm L,sc})_{ij}| \le C_\omega h^{-1/2}\sqrt{\log n}.\) The support and
coefficient bound will be used both in the common form proof and in the following resolvent comparison.

\begin{theorem}\label{A:res:15.5}
Under the coupling of the Jacobi coefficients and Brownian motion constructed above,
\[
\sup_{z\in K}
\left\|
(H_n^{B,\mathrm{mid}}-z)^{-1}
-
(H_n^{B,\mathrm{reg}}-z)^{-1}
\right\|_{\mathfrak S_2}^2
=
o_\omega(h).
\]
\end{theorem}

On the remaining interval $I_n^{\mathrm{far}}$, we have $x\asymp h^{-2}$. Therefore
\[
\sum_{j\in I_n^{\mathrm{far}}}(1+x_j)^{-7/2} \le C h^{-1} \int_{c_{\mathrm{reg}}h^{-2}}^{C h^{-2}} (1+x)^{-7/2}\,\mathrm{d}x \le C h^{-1}\cdot h^5 = C h^4.
\eqtag{15.61}
\]
Define the diagonal error on \(I_n^{\mathrm{far}}\)
\[
\mathcal K_{7,n}
:=
h^{-3}
\sum_{j\in I_n^{\mathrm{far}}}
|\zeta_j-\zeta_j^B|^2
(1+x_j)^{-7/2}.
\]
Using \eqref{A:eq:15.60} and \eqref{A:eq:15.61},
\[
\mathcal K_{7,n} \le C_\omega h^{-3} \bigl(h^3\log n\bigr) \sum_{j\in I_n^{\mathrm{far}}}(1+x_j)^{-7/2} \le C_\omega h^4\log n = o_\omega(h).
\eqtag{15.63}
\]
The bound uses $|\zeta_j-\zeta_j^B|^2\le 2|\zeta_j|^2+2|\zeta_j^B|^2$ and the separate maxima. Summable bounds for the two marginal failure probabilities give both maxima on one
event of probability one.

\par\Needspace{9\baselineskip}
\begin{theorem}\label{A:res:15.6}
Let
\[
H_n^{B,\mathrm{reg}}
=
H_n^\circ+h^{-2}M_{\zeta_n^{B,\mathrm{reg}}},
\quad
H_n^{B,\circ}
=
H_n^\circ+h^{-2}M_{\zeta_n^B}.
\]
Then, for every compact \(K\Subset\mathbb C\setminus\mathbb R,\) almost surely, for all sufficiently large \(n\),
\[
\sup_{z\in K}
\left\|
(H_n^{B,\mathrm{reg}}-z)^{-1}
-
(H_n^{B,\circ}-z)^{-1}
\right\|_{\mathfrak S_2}^2
\le
C_{\omega,V,\beta,K}h^4\log n
=
o_\omega(h).
\]
\end{theorem}

We will obtain the full coefficient comparison
\[
\sup_{z\in K}
\left\|
(\mathcal A_{n,V,\beta}-z)^{-1}
-
P_{I_n}^*U_{r_n}^*(H_n^{B,\circ}-z)^{-1}U_{r_n}P_{I_n}
\right\|_{\mathfrak S_2}^2
\le
C_{\omega,V,\beta,K}h.
\eqtag{15.65}
\]
The common form proof and the resolvent deductions below establish this estimate in that order. Subsection~\ref{A:sec:comparison-of-the-deterministic-coefficients} then keeps the
Brownian potential fixed and compares $H_n^{B,\circ}$ with $\mathcal B_{M,h,\beta}^B$.

\subsubsection{Quadratic form bounds}\label{A:sec:quadratic-form-bounds}

We use the coefficient estimates to prove the common form bounds in Theorem~\ref{A:res:16.2}, which are needed for the weighted resolvent estimates. We start with the form at a
minimizer, add the potential formed from the Brownian increments using Lemma~\ref{A:res:16.1}, and then estimate the remaining coefficient changes relative to $\ell_h$.

The finite size argument of Lemma~\ref{A:res:core-finite} allowed a logarithmic shift. It cannot supply the present uniform conclusion by itself: uniform tails for separate
coefficient arrays do not bound their supremum over all sizes. Lemma~\ref{A:res:brownian-mesh} bounds the oscillations of all piecewise linear interpolations of one Brownian path
by the same random constant. Lemma~\ref{A:res:16.1} then gives a common form shift; the maxima of the coupling errors give relative coefficients tending to zero. The removed matrix
blocks are not included in these interpolations. Their contributions to the full resolvent have already been estimated using Schur complements.

On the reindexed Dirichlet interval $\{1,\ldots,M\}$, we use $\ell_h$ from Subsection~\ref{A:sec:weighted-resolvent-estimates-and-comparison-principles}. We prove the form
comparison \eqref{A:eq:3.12} uniformly over the operators used below: for some $0<c_\omega\le C_\omega<\infty$ and $\kappa_\omega<\infty$, and all sufficiently large $n$,
\[
c_\omega\ell_h[u]
\le
\mathfrak h_n[u]+\kappa_\omega\|u\|^2
\le
C_\omega\ell_h[u]
\eqtag{16.2}
\]
for all \(u\). The shift and comparison constants are common to all these operators.

The KRV operator evaluated at the local minimizer has the form
\[
\mathfrak h_n^{\dagger}[u]
=
\sum_{j=0}^{M}c_{j,n}^{\dagger}|D_hu_j|^2
+
\sum_{j=1}^{M}V_{j,n}^{\dagger}|u_j|^2,
\]
where \(c_{j,n}^{\dagger} = \frac{b_{j,n}^{\dagger}}{b_\star},\) and
\[
V_{j,n}^{\dagger}
=
\frac{h^{-2}}{b_\star}
\left[
E_\star
-
a_{j,n}^{\dagger}
-
b_{j-1,n}^{\dagger}
-
b_{j,n}^{\dagger}
\right].
\]
The coordinates are those of Section~\ref{A:sec:the-local-minimizer-and-edge-scaling}, recentered on the interior interval. The potential estimates use the auxiliary boundary
coefficients \(b_{0,n}^\dagger,b_{M,n}^\dagger\) at the two boundary points. Lemma~\ref{A:res:auxiliary-boundary} shows that this convention leaves the matrix evaluated at the
local minimizer unchanged. We prove
\[
c_0\ell_h[u]
\le
\mathfrak h_n^{\dagger}[u]+\kappa_0\|u\|^2
\le
C_0\ell_h[u].
\eqtag{16.4}
\]
Lemma~\ref{A:res:profile-boundary} gives the edge monotonicity directly from the local stationary equations. Equivalently, in the Hessian notation of
Section~\ref{A:sec:application-to-the-conditional-jacobi-coefficients},
\(-E'(s)=(2b^\dagger(s))^{-1}\ell\mathsf H_s^{-1}\ell^*>0\), where \(\ell=(1,2)\). For a fixed \(s_0\in(0,1)\), it gives
\[
 E_\star-E(s)\ge c_{s_0}s,\quad 0\le s\le s_0,
\eqtag{16.6}
\]
and
\[
 E_\star-E(s)\ge E_\star-E(s_0)=:d_{s_0}>0,
                  \quad s_0\le s\le1.
\eqtag{16.7}
\]
The exact recentered coordinate is \(s_j=(b_\star/\tau_\star)h^2(x_j+rh+h/\beta)\), by \eqref{A:eq:recentered-profile-coordinate}. In the regular region, the midpoint choice of the
diagonal parameter and bounded second derivatives give
\[
 V_{j,n}^\dagger
 =\frac{h^{-2}}{b_\star}(E_\star-E(s_j))+O(h^{-2}n^{-2}).
\]
Thus \(V_{j,n}^\dagger\ge c x_j-C\) there, while \(|V_{j,n}^\dagger|\le C(1+x_j)\), since \(rh\le1\) for all sufficiently large \(n\). The more precise expansion near the edge has
linear part \(x_j+rh+h/\beta\); the displacement
\(rh\) is kept separate from the \(O(h)\) remainder. In the terminal region, the regularity estimate for the squared
off-diagonal coefficients gives
\[
 V_{j,n}^\dagger
 =\frac{h^{-2}}{b_\star}
       \big(E_\star-E(s_j)+O(n^{-1/2})\big),
\]
by \eqref{A:eq:far-profile-row-error}. This is sufficient: \(d_{s_0}+O(n^{-1/2})\ge d_{s_0}/2\) for large \(n\). Also
\(1+x_j\asymp h^{-2}\) in that region, and the local minimizing coefficients are bounded. Combining both regions proves
\[
 V_{j,n}^\dagger\ge c(1+x_j)-C,
\eqtag{16.13}
\]
and
\[
 |V_{j,n}^\dagger|\le C(1+x_j).
\eqtag{16.14}
\]
These potential statements use the choice of auxiliary coefficients at the boundary points specified above. Choosing both auxiliary boundary coefficients equal to $b_\star$
preserves the complete form, but need not preserve the displayed pointwise potential inequalities.

Uniform ellipticity of the difference coefficient is needed only in the regular region. There $b^\dagger$ is bounded below, so $c_{j,n}^\dagger\ge c>0$. For every remaining edge,
each nonzero boundary point has $1+x\asymp h^{-2}$, after enlarging the regular interval by one site if necessary. Hence
\[
 |D_hu_j|^2\le2h^{-2}(|u_j|^2+|u_{j+1}|^2)
 \le C\bigl((1+x_j)|u_j|^2+(1+x_{j+1})|u_{j+1}|^2\bigr).
\]
At the terminal Dirichlet edge the absent boundary point is zero. Each vertex occurs in at most two edges, so summing uses the confining potential to control the difference energy
even where $c_{j,n}^\dagger$ tends to zero. Combining the lower off-diagonal coefficient bound in the regular region with the potential control in the terminal region gives
\[
\sum_{j=0}^{M}|D_hu_j|^2
\le
C
\left[
\sum_{j=0}^{M}c_{j,n}^{\dagger}|D_hu_j|^2
+
\sum_{j=1}^{M}(1+x_j)|u_j|^2
\right].
\eqtag{16.18}
\]
Equations \eqref{A:eq:16.13} and \eqref{A:eq:16.18} give \(\ell_h[u] \le C(\mathfrak h_n^\dagger[u]+\kappa\|u\|^2)\). Conversely, \(0\le c_{j,n}^{\dagger}\le C\) and
\eqref{A:eq:16.14} imply \(\mathfrak h_n^\dagger[u]+\kappa\|u\|^2 \le C\ell_h[u]\). Together these estimates are \eqref{A:eq:16.4}.

For the form at the conditional minimizer in \eqref{A:eq:13.18}, write
\[
c_{j,n}^{\circ}=\frac{B_j^\circ}{b_\star},
\quad
V_{j,n}^{\circ}
=
\frac{h^{-2}}{b_\star}
\left[E_\star-A_j^\circ-B_{j-1}^\circ-B_j^\circ\right].
\]
The normalized boundary values are \(B_0^\circ=B_M^\circ=b_\star\). To apply the confinement estimate of Lemma~\ref{A:res:12.1B}, instead write the same matrix at the conditional
minimizer with \(\widehat B_0^\circ=b_r\),
\(\widehat B_M^\circ=b_N\), and \(\widehat B_j^\circ=B_j^\circ\) for \(1\le j\le M-1\). Let
\[
\widehat V_{j,n}^\circ
=\frac{h^{-2}}{b_\star}
 (E_\star-A_j^\circ-\widehat B_{j-1}^\circ-\widehat B_j^\circ).
\]
The potential coefficients in the two representations satisfy
\[
V_{j,n}^\circ=\widehat V_{j,n}^\circ
 -\frac{h^{-2}}{b_\star}
 \bigl[(b_\star-\widehat B_0^\circ)\mathbf1_{\{j=1\}}
                +(b_\star-\widehat B_M^\circ)\mathbf1_{\{j=M\}}\bigr].
\eqtag{auxiliary-boundary-potential}
\]
For \(M\ge2\), at the right boundary, for example,
\[
|D_hu_M|^2+V_{M,n}^\circ|u_M|^2
=\frac{\widehat B_M^\circ}{b_\star}|D_hu_M|^2
 +\widehat V_{M,n}^\circ|u_M|^2.
\eqtag{auxiliary-boundary-boundary}
\]
At the left boundary, when $M\ge2$,
\[
 |D_hu_0|^2+V_{1,n}^\circ|u_1|^2
 =\frac{\widehat B_0^\circ}{b_\star}|D_hu_0|^2
   +\widehat V_{1,n}^\circ|u_1|^2,
\]
since $D_hu_0=u_1/h$ and $V_{1,n}^\circ-\widehat V_{1,n}^\circ =h^{-2}(\widehat B_0^\circ/b_\star-1)$. For \(M=1\) both boundary corrections are paired together as in
\eqref{A:eq:auxiliary-boundary-cancellation}. Thus \(\widehat V_{M,n}^\circ\ge c h^{-2}\) controls the boundary term at the last site in the normalized representation even when
\(V_{M,n}^\circ\) itself is negative. A positive lower bound for \(\widehat V^\circ\) need not give a positive lower
bound for \(V^\circ\).

Lemma~\ref{A:res:12.1B} concerns the complete form with original exterior coefficients. By Corollary~\ref{A:res:conditioning}, it is the same conditional minimizer as here, and
Lemma~\ref{A:res:auxiliary-boundary} identifies its complete form exactly with \(\mathfrak h_n^\circ\). Hence, on $\Omega_{\mathrm{bd}}$, the same confinement estimate gives
\[
c_\omega\ell_h[u]
\le
\mathfrak h_n^\circ[u]+\kappa_\omega\|u\|^2
\le
C_\omega\ell_h[u]
\eqtag{16.22}
\]
for every \(u\), for all sufficiently large \(n\).

For the interval $I_n^{\mathrm{cp,mid}}$ use the coefficient estimate \eqref{A:eq:13.67}. On $I_n^{\mathrm{cp,mid}}$
\(|c_{j,n}^\circ-c_{j,n}^\dagger| \le C_\omega h^3, \allowbreak\, |V_{j,n}^\circ-V_{j,n}^\dagger| \le C_\omega h\). On
the logarithmic left boundary strip we use
\[
\max_{j\le J_n^{\mathrm L}+R_V}
\left(
|A_j^\circ-a_j^\dagger|
+
|B_j^\circ-b_j^\dagger|
\right)
\le
C_\omega h^{3/2}\sqrt{\log n}.
\]
After edge scaling, the corresponding banded matrix has entries at most \(C_\omega h^{-1/2}\sqrt{\log n}\), but its relative bound as a quadratic form, obtained from the Hardy
inequality, is \(C_\omega h^{3/2}(J_n^{\mathrm L}+1)^2\sqrt{\log n}=o_\omega(1).\) Near the right boundary the displacement need not be small;
Section~\ref{A:sec:comparison-of-the-deterministic-coefficients} retains its $O(h)$ resolvent tail.

Lemma~\ref{A:res:16.1} bounds the form of $\Gamma_E$ through oscillations of the interpolated sequence $E$. The corresponding continuum bound is proved in

\cite{ramirez2011beta}. We have already used it with the centered coefficient partial sums. We now apply the same deterministic estimate to samples of a single Brownian path, for
which the block constant is independent of $h$.

On the set \(\Omega_B\) defined in Subsection~\ref{A:sec:finite-difference-approximation-of-the-stochastic-airy-operator}, for every \(0<\alpha<1/2\),
\[
 A_{B,\alpha}:=
 \sup_{k\in\mathbb N_0}(1+k)^{-\alpha}
 \bigl(\operatorname{osc}_{[k,k+1]}B
                  +|B(k+1)-B(k)|\bigr)<\infty.
\eqtag{brownian-block}
\]

\begin{lemma}\label{A:res:brownian-mesh} For each fixed \(\alpha>0\),
\(\mathbb P(A_{B,\alpha}>R)\le C_\alpha \mathrm{e}^{-c_\alpha R^2}\) for \(R\ge1\). The piecewise linear interpolations
of \(E_j^B=2\beta^{-1/2}B(jh)\), kept constant beyond \(Mh\), have \(\mathcal O_\alpha\le C_{\alpha,\beta}A_{B,\alpha}\) simultaneously for every \(0<h\le1\) and \(M\ge1\).
\end{lemma}

\begin{proof}
The tail bound for $A_{B,\alpha}=\mathcal O_\alpha(B)$ is Lemma~\ref{A:res:A.6}. Let $S_{h,M}$ be the continuous linear interpolation of $c_\beta B(jh)$, kept constant beyond $Mh$.
For $x\in[k,k+1]$, its value is a convex combination of sample values at times in $[\max\{0,k-1\},k+2]$. If $Mh<k$, it is constant throughout this block. In either case,
\[
 \operatorname{osc}_{[k,k+1]}S_{h,M}
 \le c_\beta\sum_{r=\max\{0,k-1\}}^{k+1}
                         \operatorname{osc}_{[r,r+1]}B,
 \qquad
 |S_{h,M}(k+1)-S_{h,M}(k)|
 \le\operatorname{osc}_{[k,k+1]}S_{h,M}.
\]
For each summand, $(1+r)^\alpha\le 2^\alpha(1+k)^\alpha$. There are at most three summands. Multiplication by $(1+k)^{-\alpha}$ and the definition of $A_{B,\alpha}$ therefore give
\[
 \mathcal O_\alpha(S_{h,M})
 \le 6\,2^\alpha c_\beta A_{B,\alpha}.
\]
This deterministic inequality holds for every $0<h\le1$ and $M\ge1$ on each path for which $A_{B,\alpha}<\infty$.
\end{proof}

Since $h^{-2}\zeta_j^B=h^{-1}(E_j^B-E_{j-1}^B)$, Lemma~\ref{A:res:16.1} gives, using \(\alpha=1/4\),
\[
 \left|\sum_{j=1}^{M}h^{-2}\zeta_j^B|u_j|^2\right|
 \le\varepsilon\ell_h[u]
   +C_{\varepsilon,\beta}(1+A_{B,1/4})^4\|u\|^2,
\eqtag{16.33}
\]
with the same random constant for all \(n\), all terminal lengths, and \(\beta\) in a fixed compact positive interval.

The remaining first-order terms are \(\Gamma_{\eta_n}^{-}\), from the off-diagonal coefficient fluctuations, and
\(\Gamma_{\widetilde e_n}\), with the averaging map $J_h=(J_-+J_+)/2$. Lemma~\ref{A:res:factor-max}, with
Borel--Cantelli, gives
\[
 \|\eta_n\|_\infty+\|\widetilde e_n\|_\infty
 \le C_\omega h^{1/2}\sqrt{\log(n+1)}=o_\omega(1).
\eqtag{16.36}
\]
For either evaluation map at a cell boundary \(J_-\), \(J_+\), and for \(J_h=(J_-+J_+)/2\), its operator norm is at most one. Hence
\[
 |\langle u,\Gamma_v^{\sigma}u\rangle| \le2\|v\|_\infty\|D_hu\|\,\|u\|
 \le\|v\|_\infty(\|D_hu\|^2+\|u\|^2) \le\|v\|_\infty\ell_h[u], \quad
 \sigma\in\{-,+,\mathrm{mid}\}.
\eqtag{16.38}
\]
The maximal bounds therefore give vanishing relative form coefficients for both perturbations. By \eqref{A:eq:16.28}, the center term is a real diagonal multiplier; only the
one-sided term has off-diagonal entries.

For the weighted resolvent estimates we must also retain the tridiagonal representation \eqref{A:eq:3.8} and the coefficient bound \eqref{A:eq:3.10}. All scalar replacements,
including $\Gamma_{\widetilde e}$ by \eqref{A:eq:16.28}, are real diagonal multipliers. Reciprocal exponential conjugation leaves these diagonal multipliers unchanged, and their
contributions therefore cancel identically.

The only additional off-diagonal term is $\Gamma_\eta^-$, whose coefficient is bounded by \eqref{A:eq:16.36}. Its off-diagonal matrix entry is $-\eta_j/h$. Use the truncated
exponential weights $\theta_j=r_j=\mathrm e^{a\psi_j}$ of Lemma~\ref{A:res:3.12}. Substitution in the tridiagonal localization identity \eqref{A:eq:3.32c} and the weight estimate
\eqref{A:eq:3.31} gives
\[
 \left|\langle\theta u,\Gamma_\eta^-\theta u\rangle
 -\operatorname{Re}\langle u,\Gamma_\eta^-\theta^2u\rangle\right|
 \le C_\omega a^2h\|\eta\|_\infty
           \sum_j(1+x_j)\theta_j^2|u_j|^2.
\]
For reciprocal weights the calculation \eqref{A:eq:3.42j} gives
\[
 |\gamma_{\eta,-}[r^{-1}u,rv]-\gamma_{\eta,-}[u,v]|
 \le C_\omega a\|\eta\|_\infty\ell_h[u]^{1/2}\ell_h[v]^{1/2}.
\]
Thus the hypotheses for localization are preserved as well as coercivity.

Let \(q=(q_j)\) be a diagonal multiplier. The elementary estimate below will be applied to each diagonal remainder. Suppose
\[
|q_j|
\le
\varepsilon_n(1+x_j)+C_{\omega},
\quad
\varepsilon_n\to0.
\eqtag{16.40}
\]
The corresponding diagonal form obeys
\[
\left|
\sum_jq_j|u_j|^2
\right|
\le
\varepsilon_n
\sum_j(1+x_j)|u_j|^2
+
C_\omega\|u\|^2.
\]
For large \(n\), \(| \langle u,M_qu\rangle | \le \varepsilon\ell_h[u] + C_{\varepsilon,\omega}\|u\|^2\). The concrete multipliers and their relative bounds are estimated in the
proof of Theorem~\ref{A:res:16.2}. The boundary terms of the full matrix are estimated separately by the Schur complement formulas in Subsection~\ref{A:sec:boundary-reduction}.

Fix $V$ and $\beta$. Let $\mathbb W$ be the law of the two-sided Brownian path, let $\nu_n$ be the exact exterior marginal at size $n$, and let $\mu_n^\tau$ be the original
conditional distribution of the remaining coefficients given exterior data $\tau$. In this paragraph write $Q_n(B,\tau,\mathrm dz)$ for the Borel kernel for the joint distribution
of the coefficients and auxiliary Gaussian variables supplied by Lemma~\ref{A:res:B.8}; the coordinate $z$ includes the interior coefficients and any auxiliary Gaussian variables.
Its projection onto the interior coefficients satisfies
\[
 \int\mathbb W(\mathrm dB)\,
        Q_n(B,\tau,\{z:\operatorname{int}(z)\in A\})
       =\mu_n^\tau(A)
\]
for every Borel set $A$ and each admissible exterior parameter. On the complement of the bounded exterior set, the specified product coupling has this same marginal identity. Thus
this identity uses the original conditional law, not that law conditioned on a coefficient event.

Use the countable product of these kernels with the single Brownian coordinate:
\[
 \mathbb W(\mathrm dB)\prod_{n\ge1}
       \bigl[\nu_n(\mathrm d\tau_n)Q_n(B,\tau_n,\mathrm dz_n)\bigr].
\]
Equivalently, the Borel samplers of Lemma~\ref{A:res:B.8} realize this product using independent uniform random variables, independent of $B$. This construction leaves the Brownian
marginal equal to $\mathbb W$ and makes each exterior variable independent of $B$. Let $C_n(\tau,z)$ denote the full coefficient vector obtained by placing the exterior and
interior coordinates back in their original positions. For a bounded Borel function $f$ of that vector, integration of the other factors in the product gives
\[
 \begin{aligned}
 \mathbb E f(C_n(\tau_n,z_n))
 &=\int\nu_n(\mathrm d\tau)\int\mathbb W(\mathrm dB)
                  \int Q_n(B,\tau,\mathrm dz)f(C_n(\tau,z))\\
 &=\int\nu_n(\mathrm d\tau)\int\mu_n^\tau(\mathrm dy)
                                      f(C_n(\tau,y)).
 \end{aligned}
\]
The last expression is exactly the original finite coefficient expectation by disintegration. Hence every finite matrix has its prescribed law, while all sizes use the same
Brownian path. The construction is fixed before choosing a moment order or an exponent in the probability bounds. Coefficient vectors at different matrix sizes may be dependent
through the common Brownian path.

Fix a deterministic \(C_*\) large enough for the following bounds, whose exceptional events have probability at most a constant times \(n^{-3}\).

Let \(\mathcal E_n^{\mathrm{KRV}}\) be the intersection of \(\mathcal F_{n,3}\) from \eqref{A:eq:confinement-good-event}, \(\mathcal T_n^{\mathrm{out}}(M_0)\), and
\[
 \begin{gathered}
 \max_{I_n^{\mathrm{mid}}}|\zeta_j-\widehat\zeta_j^{\mathrm G}|
       \le n^{-5/6},\quad
 \|\eta_n\|_\infty+\|\widetilde e_n\|_\infty
       \le C_*h^{1/2}\sqrt{\log(n+1)},\\
 \max_{1\le j\le M}(|\zeta_j|+|\zeta_j^B|)
       \le C_*h^{3/2}\sqrt{\log(n+1)},\\
 \max_{I_n^{\mathrm{mid}}}
       \frac{h^{-2}|\rho_j|}{1+x_j}
       \le C_*h^{3/2}\sqrt{\log(n+1)},\quad
 \max_{I_n^{\mathrm{mid}}}
       \frac{h^{-2}|\varepsilon_j^{\mathrm{fr}}|}{1+x_j}
       \le C_*h^{1/14}\sqrt{\log(n+1)}.
 \end{gathered}
\]
Use \(p=30\) in Lemma~\ref{A:res:transport-max} and exponent \(4\) in the probability bounds of Lemmas~\ref{A:res:freezing-max} and \ref{A:res:factor-max}. Together with
\eqref{A:eq:confinement-good-event} and disintegration over the exterior data, these bounds give
\[
 \mathbb P((\mathcal E_n^{\mathrm{KRV}})^c)\le Cn^{-3}+C\mathrm{e}^{-cn},
 \quad
 \Omega_{\mathrm{KRV}}:=\{B|_{\mathbb R_+}\in\Omega_B\}
                \cap\liminf_{n\to\infty}\mathcal E_n^{\mathrm{KRV}},
 \quad \mathbb P(\Omega_{\mathrm{KRV}})=1.
\eqtag{common-good-event}
\]
Here \(\liminf_n E_n=\bigcup_m\bigcap_{n\ge m}E_n\). The comparison with the conditional minimizer in the regular region also holds on \(\mathcal F_{n,3}\), by
\eqref{A:eq:global-mode-profile}.

On this event the bounds for the conditional minimizer \eqref{A:eq:13.67}--\eqref{A:eq:13.67a} and the exact maxima \eqref{A:eq:exact-maxima} hold for all sufficiently large \(n\).
The multiplier bounds are, with \(X_j=1+x_j\),
\[
 \max_{j\in I_n^{\mathrm{mid}}}
       |\zeta_j-\widehat\zeta_j^{\mathrm G}|
 \le n^{-5/6},\quad
 h^{-2}\max_{j\in I_n^{\mathrm{mid}}}|\zeta_j-\widehat\zeta_j^{\mathrm G}|
 \le C_Vn^{-1/6},
\eqtag{16.51}
\]
and
\[
 \max_{j\in I_n^{\mathrm{mid}}} \frac{h^{-2}|\varepsilon_j^{\mathrm{fr}}|}{X_j} \le C_\omega
 h^{1/14}\sqrt{\log(n+1)},\qquad \max_{j\in I_n^{\mathrm{mid}}} \frac{h^{-2}|\rho_j|}{X_j} \le
 C_\omega h^{3/2}\sqrt{\log(n+1)}.
\eqtag{16.52}
\]
The cutoff extensions of \(\rho\) and \(\varepsilon^{\mathrm{fr}}\) vanish outside the middle interval, so these bounds also control the corresponding multipliers on the full
interval.

The proof starts from the form bound for the Brownian potential \eqref{A:eq:16.33} and uses the preceding coordinate estimates to compare the other forms. The boundary blocks
satisfy the Schur complement estimates of Subsection~\ref{A:sec:boundary-reduction}.

Assume throughout the following theorem that $V$ is uniformly convex and the right edge is regular. Let $\mathcal H_n$ be the finite family of operators in
\eqref{A:eq:common-assembly}, together with $H_n^{B,\dagger}$ and $\mathcal B_{M,h,\beta}^B$.

\begin{theorem}\label{A:res:16.2}
On $\Omega_{\mathrm{KRV}}$, the form bound \eqref{A:eq:16.2} and the coefficient bounds \eqref{A:eq:3.8}, \eqref{A:eq:3.10} hold, for all sufficiently large $n$, uniformly for
$H\in\mathcal H_n$, their convex interpolations, and their principal Dirichlet restrictions with inherited coordinates.
\end{theorem}

The constants are independent of $n$ and of the member of this family. Hence the weighted resolvent estimates of
Section~\ref{A:sec:weighted-resolvent-estimates-and-comparison-principles} apply to each comparison in the argument for the general Jacobi model and to the first-order comparison
deferred from Section~\ref{A:sec:decomposition-of-the-general-jacobi-coefficients}.

\begin{proof}
On the coefficient events in \eqref{A:eq:common-good-event}, Lemma~\ref{A:res:12.1B} gives the comparison for the form at the conditional minimizer with fixed constants for all
sufficiently large \(n\). The proof of \eqref{A:eq:16.4} gives the same conclusion for the form evaluated at the local minimizer. For the quadratic form of the deterministic
difference operator with linear potential, \(\langle u,(-\Delta_{h,D}+X_h)u\rangle =\ell_h[u]-\|u\|^2\). The coefficient bounds are uniform on these events. Choose deterministic
constants $c_*>0,C_*,\kappa_*$ such that the three forms without the Brownian potential, uniformly over the coefficients in these events, are bounded below by
\(c_*\ell_h-\kappa_*\|\cdot\|^2\) and above by \(C_*\ell_h+\kappa_*\|\cdot\|^2\), for all sufficiently large \(n\). Take \(\varepsilon=c_*/4\) in \eqref{A:eq:16.33} and put
\[
 \kappa_B=\kappa_*+
 C_{c_*/4,\beta}(1+A_{B,1/4})^4.
\eqtag{common-brownian-shift}
\]
The three operators with the same Brownian term \(H_n^{B,\circ}\), \(H_n^{B,\dagger}\) and \(\mathcal B_{M,h,\beta}^B\) then satisfy
\[
 \frac{3c_*}{4}\ell_h[u]-\kappa_B\|u\|^2
 \le\mathfrak h_n^{B}[u]
 \le(C_*+c_*/4)\ell_h[u]+\kappa_B\|u\|^2.
\eqtag{common-brownian-form}
\]
The same constant controlling Brownian oscillations over unit intervals determines $\kappa_B$ for every matrix size. For each remaining perturbation with quadratic form $r_n$ we
establish a bound $|r_n[u]|\le a_n\ell_h[u]$ with $a_n\to0$. Adding such a perturbation to the lower bound in \eqref{A:eq:common-brownian-form} changes its coefficient from
$3c_*/4$ to $3c_*/4-a_n$. Thus the coefficient maxima will establish coercivity before any of the weighted resolvent estimates is applied.

For a diagonal multiplier \(q\),
\[
 |\langle u,M_qu\rangle|
 \le\left(\max_j\frac{|q_j|}{1+x_j}\right)\ell_h[u].
\]
Define the real diagonal multipliers
\[
 \begin{aligned}
Q_n^{\mathrm L}& =h^{-2}M_{\mathbf1_{I_n^{\mathrm L}}(\zeta-\zeta^B)},\qquad Q_n^{\mathrm{far}} =h^{-2}M_{\mathbf1_{I_n^{\mathrm{far}}}(\zeta-\zeta^B)},\\
Q_n^\rho& =h^{-2}M_{\widetilde\rho},\qquad Q_n^{\mathrm{fr}} =h^{-2}M_{\widetilde\varepsilon^{\mathrm{fr}}},\qquad Q_n^{\mathrm{tr}} =h^{-2}M_{\chi^{\mathrm{mid}}(\zeta-\widehat\zeta^{\mathrm G})}.
\end{aligned}
\]
By \eqref{A:eq:16.51}, the last multiplier has operator norm at most \(C_Vn^{-1/6}\). The relative form coefficients of
\(Q_n^{\mathrm{fr}}\) and \(Q_n^\rho\) are at most, respectively, \(C_\omega h^{1/14}\sqrt{\log(n+1)}\) and \(C_\omega
h^{3/2}\sqrt{\log(n+1)}\), by \eqref{A:eq:16.52}. Since \(1+x_j\ge ch^{-2}\) on the far interval, the exact and Brownian maxima give the same latter relative coefficient for
\(Q_n^{\mathrm{far}}\).

The left multiplier has the growing maximum $h^{-1/2}\sqrt{\log(n+1)}$, so we also use its support on the left strip. There the Dirichlet condition bounds each value by the
discrete derivative. For \(1\le j\le J_n^{\mathrm L}\),
\[
 |u_j|^2
 =h^2\left|\sum_{k=0}^{j-1}D_hu_k\right|^2
 \le jh^2\sum_{k=0}^{j-1}|D_hu_k|^2.
\]
Since the entries of \(Q_n^{\mathrm L}\) are at most \(C_\omega h^{-1/2}\sqrt{\log(n+1)}\), summing in \(j\) yields
\[
 |\langle u,Q_n^{\mathrm L}u\rangle|
 \le C_\omega h^{3/2}(J_n^{\mathrm L}+1)^2
                  \sqrt{\log(n+1)}\,\ell_h[u]
 =o_\omega(1)\ell_h[u].
\eqtag{16.55b}
\]
The constant extension of \(\widetilde e\) creates no derivative at either scalar cutoff: its adjacent values agree outside the middle block. Its maximum is bounded by
Lemma~\ref{A:res:factor-max}. Equation~\eqref{A:eq:16.38} therefore bounds the relative coefficients of both \(\Gamma_{\widetilde e}\) and \(\Gamma_\eta^-\) by
\(C_\omega h^{1/2}\sqrt{\log(n+1)}\).

Consequently the sum of the absolute relative coefficients of all seven perturbations just considered is bounded by
\[
 a_n=C_\omega\bigl[ n^{-1/6} +(h^{1/14}+h^{1/2}+h^{3/2})\sqrt{\log(n+1)} +h^{3/2}(J_n^{\mathrm
 L}+1)^2\sqrt{\log(n+1)}\bigr] \longrightarrow0.
\eqtag{common-relative-error}
\]

The defining coefficient identities give
\[
 \begin{aligned}
H_n^{B,\mathrm{reg}}& =H_n^{B,\circ}+Q_n^{\mathrm{far}},\qquad H_n^{B,\mathrm{mid}} =H_n^{B,\circ}+Q_n^{\mathrm{far}}+Q_n^{\mathrm L},\qquad \widetilde H_n =H_n^{B,\mathrm{mid}}+\Gamma_{\widetilde e},\\
H_n^{\mathrm G}& =\widetilde H_n+Q_n^\rho+Q_n^{\mathrm{fr}},\qquad H_n^{\mathrm{diag}} =H_n^{\mathrm G}+Q_n^{\mathrm{tr}},\qquad \widetilde H_n^{\mathrm{red}} =H_n^{\mathrm{diag}}+\Gamma_\eta^-.
\end{aligned}
\eqtag{common-assembly}
\]
For the seven perturbations $Q_n^{\mathrm{far}},Q_n^{\mathrm L},\Gamma_{\widetilde e}, Q_n^\rho,Q_n^{\mathrm{fr}},Q_n^{\mathrm{tr}},\Gamma_\eta^-$, write $r_{\ell,n}$ for their
quadratic forms. The preceding estimates have the form
\[
 |r_{\ell,n}[u]|\le a_{\ell,n}\ell_h[u],\quad
 \sum_{\ell=1}^7a_{\ell,n}\le a_n\longrightarrow0.
\]
Every form in \eqref{A:eq:common-assembly} is obtained from $\mathfrak h_n^{B,\circ}$ by adding a subset of these terms. More generally, for $0\le t_\ell\le1$,
\[
 \sup_{t\in[0,1]^7}\left|\sum_{\ell=1}^7t_\ell r_{\ell,n}[u]\right|
 \le a_n\ell_h[u].
\]
Choose the random threshold so that $a_n\le c_*/4$. Together with \eqref{A:eq:common-brownian-form}, this proves for each of the stated operators
\[
 \frac{c_*}{2}\ell_h[u]
 \le\mathfrak h_n[u]+\kappa_B\|u\|^2
 \le(C_*+c_*/2+2\kappa_B)\ell_h[u].
\eqtag{16.63}
\]
Convex combinations with the forms of $H_n^{B,\dagger}$ and $\mathcal B_{M,h,\beta}^B$ obey the same inequalities by \eqref{A:eq:common-brownian-form}. For a principal Dirichlet
restriction, extend its vector by zero in the full interval and keep the inherited \(x_j\) in the comparison form. The zero extension includes the energies of both cut edges,
proving the restriction bound. The localized comparison replacements in Section~\ref{A:sec:comparison-of-the-deterministic-coefficients} additionally use
\eqref{A:eq:13.67}--\eqref{A:eq:13.67a}: the regular difference is \(O(h)\ell_h\), and the difference on the left, which is banded, has the Hardy bound in \eqref{A:eq:3.58h}. Hence
those localized interpolations have the same common form bounds after one enlargement of the constants.

For a scalar recentering by $\sigma rh$, $\sigma\in\{-1,1\}$,
\[
 |\sigma rh\|u\|^2|\le rh\,\ell_h[u].
\]
Choose $n$ so that $rh\le c_*/4$. Subtracting or adding this term to \eqref{A:eq:16.63} gives a lower bound $(c_*/4)\ell_h[u]$ and an upper bound $(C_*+3c_*/4+2\kappa_B)\ell_h[u]$
with the same shift $\kappa_B$. These inequalities follow directly from the bounds for the kinetic and potential parts of the forms. This recentered form inequality is used only
for the reduced block. The full matrix is related to it by the Schur complement estimates of Section~\ref{A:sec:boundary-reduction}.

The coefficients in the difference terms at the conditional and local minimizers are nonnegative and bounded above uniformly on the coefficient events; the difference operator with
linear potential has coefficient one in front of $\sum_j|D_hu_j|^2$. All the \(Q\)'s are diagonal. By \eqref{A:eq:16.28}, \(\Gamma_{\widetilde e}\) is also diagonal, while
\(\Gamma_\eta^-\) is a one-sided first-order term with \(\|\eta\|_\infty=o_\omega(1)\). Thus every listed operator satisfies \eqref{A:eq:3.8} and \eqref{A:eq:3.10}, and so do their
convex interpolations. Real diagonal terms cancel in exponential conjugation and in the error in the weighted quadratic form; the one-sided term has the two bounds displayed above
\eqref{A:eq:16.40}. Together with \eqref{A:eq:16.63}, these are all the hypotheses used by the weighted estimates in
Section~\ref{A:sec:weighted-resolvent-estimates-and-comparison-principles}. Their constants can therefore be chosen uniformly for each fixed

\(K\Subset\mathbb C\setminus\mathbb R\).
\end{proof}

For the resolvent bounds, fix a sufficiently large deterministic \(A\) and let \(\mathcal W_n\) be the event
\[
 \begin{gathered}
\mathcal E_{\eta,n}\le Ah,\quad \mathcal T_n^{\mathrm{mid}}\le h^4,\qquad \sum_{j\in I_n^{\mathrm{mid}}}(1+x_j)^{-7/2} |\varepsilon_j^{\mathrm{fr}}|^2\le h^4,\quad \sum_{j\in I_n^{\mathrm{mid}}}(1+x_j)^{-7/2}|\rho_j|^2\le Ah^4,\\
\sum_{j=0}^{M}h(1+x_j)^{-3/2}|\widetilde e_j|^2\le Ah.
\end{gathered}
\]
For the first bound use the exponential tail \eqref{A:eq:13.50}. For the next three sums use their conditional moment bounds from
Subsections~\ref{A:sec:gaussian-approximation-of-the-conditional-gibbs-measure}--\ref{A:sec:application-to-the-conditional-jacobi-coefficients}. The first two have $L^p$ size
$O(h^5)$ and the remainder sum has size $O(h^6)$. Markov's inequality at the displayed $h^4$ thresholds therefore gives $O(h^p)$ and $O(h^{2p})$, respectively. Taking $p=12$ gives
bounds at most $Cn^{-4}$ and $Cn^{-8}$, since $h\asymp n^{-1/3}$. The last sum uses the tail estimate for Gaussian quadratic forms of Proposition~\ref{A:res:A.20}; the constant
extensions contribute at most \(C(h^2\log^2(n+1)+h^2\log(n+1))=o(h)\) on the maximal event of Lemma~\ref{A:res:factor-max}. Thus, after disintegration,
\[
 \mathbb P(\mathcal W_n^c)\le Cn^{-3}+C\mathrm{e}^{-c/h}+C\mathrm{e}^{-cn},\quad
 \Omega_{V,\beta}:=\Omega_{\partial}\cap\Omega_{\mathrm{KRV}}
                      \cap\liminf_n\mathcal W_n,
 \quad \mathbb P(\Omega_{V,\beta})=1.
\]
The common form estimate and all weighted comparisons now hold on this specified event on the joint probability space. In particular, the Gaussian covariance comparison contributes
\(h^{-3}\sum_j(1+x_j)^{-7/2}|\varepsilon_j^{\mathrm{fr}}|^2\le h\); the weaker \(h^{1/14}\) maximum was used only for
coercivity.

\paragraph*{Resolvent consequences of the coefficient estimates}
We now apply the common form bounds of Theorem~\ref{A:res:16.2} to the coefficient estimates.

Write $R_n^{\mathrm{int}}(z)=(\widetilde H_n^{\mathrm{red}}-z)^{-1}$ and $R_n^{\mathrm{diag}}(z)=(H_n^{\mathrm{diag}}-z)^{-1}$. The matrix identity gives
\[
 R_n^{\mathrm{int}}-R_n^{\mathrm{diag}}
 =-R_n^{\mathrm{int}}\Gamma_{\eta_n}^-R_n^{\mathrm{diag}}.
\]
The coefficient of $|D_hu_j|^2$ can be chosen as $B_j^\circ/b_\star$, which is nonnegative and bounded; all scalar terms are diagonal. The remaining coefficient $\eta_j$ is
uniformly bounded by \eqref{A:eq:12.15p}. Thus the tridiagonal and coefficient hypotheses, as well as the common form hypothesis, of the first-order estimate in
Section~\ref{A:sec:weighted-resolvent-estimates-and-comparison-principles} are satisfied. For its two terms the Hilbert--Schmidt ideal property gives
\[
 \begin{split}
 \|R_n^{\mathrm{int}}D_h^*M_\eta J_-R_n^{\mathrm{diag}}\|_{\mathfrak S_2}
 &\le\|D_hR_n^{\mathrm{int}}(\bar z)\|_{\mathrm{op}}
        \|M_\eta J_-R_n^{\mathrm{diag}}(z)\|_{\mathfrak S_2},\\
 \|R_n^{\mathrm{int}}J_-^*M_\eta D_hR_n^{\mathrm{diag}}\|_{\mathfrak S_2}
 &\le\|M_\eta J_-R_n^{\mathrm{int}}(\bar z)\|_{\mathfrak S_2}
        \|D_hR_n^{\mathrm{diag}}(z)\|_{\mathrm{op}}.
 \end{split}
\]
The derivative factors are uniformly bounded. Squaring the multiplication factors and summing the resolvent columns bounds each by $C_{\omega,K}\sum_j h(1+x_j)^{-3/2}|\eta_j|^2$.
Proposition~\ref{A:res:13.3} therefore yields, on this coupling,
\[
 \sup_{z\in K}\|R_n^{\mathrm{int}}(z)-R_n^{\mathrm{diag}}(z)\|_{\mathfrak S_2}^2
  \le C_{\omega,V,\beta,K}h.
\]

Combining the resolvent estimate for the full matrix \eqref{A:eq:12.43}, the identification \eqref{A:eq:auxiliary-boundary-operator}, Proposition~\ref{A:res:13.1},
Proposition~\ref{A:res:13.3}, and the first-order resolvent estimate of Section~\ref{A:sec:weighted-resolvent-estimates-and-comparison-principles} gives

\[
\sup_{z\in K}
\left\|
(\mathcal A_{n,V,\beta}-z)^{-1}
-
P_{I_n}^*U_r^*(H_n^{\mathrm{diag}}-z)^{-1}U_rP_{I_n}
\right\|_{\mathfrak S_2}^2
\le
C_{\omega,V,\beta,K}h.
\]
\begin{proof}[Proof of Theorem~\ref{A:res:14.7}]
The left side of \eqref{A:eq:14.47} is $\mathcal T_n^{\mathrm{mid}}$. The unconditional estimate $\|\mathcal T_n^{\mathrm{mid}}\|_{L^p} \le C_ph^5$, including the exceptional
exterior data, was proved before \eqref{A:eq:14.39}. Markov's inequality at threshold $h^4$ has summable right side $C_ph^p$ for $p>3$. Thus \eqref{A:eq:14.39} proves
\eqref{A:eq:14.47}. The covariance and its decay are \eqref{A:eq:14.45}--\eqref{A:eq:14.46}.

For the operator estimate, use the common form bounds of Theorem~\ref{A:res:16.2}, proved from the coefficient estimates alone. Write
$R_{\mathrm d}(z)=(H_n^{\mathrm{diag}}-z)^{-1}$ and $R_{\mathrm G}(z)=(H_n^{\mathrm G}-z)^{-1}$. Their exact difference is
\[
 R_{\mathrm d}(z)-R_{\mathrm G}(z)
 =-R_{\mathrm d}(z)h^{-2}
   M_{\mathbf1_{I_n^{\mathrm{mid}}}(\zeta-\widehat\zeta^{\mathrm G})}
   R_{\mathrm G}(z).
\]
The diagonal perturbation estimate \eqref{A:eq:3.49} therefore gives
\[
 \sup_{z\in K}\|R_{\mathrm d}(z)-R_{\mathrm G}(z)\|_{\mathfrak S_2}^2
  \le C_{\omega,K}h^{-3}\mathcal T_n^{\mathrm{mid}}
 \le C_{\omega,V,\beta,K}h.
\]
The bounds for all finite moments give a sharper rate. For each $0<\epsilon<1$, choose $p>3/\epsilon$. Then
\[
 \mathbb P(\mathcal T_n^{\mathrm{mid}}>h^{5-\epsilon})
 \le C_ph^{p\epsilon},\qquad \sum_nh^{p\epsilon}<\infty.
\]
On the resulting probability one event,
\[
 \sup_{z\in K}\|R_{\mathrm d}(z)-R_{\mathrm G}(z)\|_{\mathfrak S_2}^2
 \le C_{\omega,V,\beta,K,\epsilon}h^{2-\epsilon}=o(h).
\]
The deterministic form estimates on each path supply the random resolvent constant; no expectation involving that constant is taken. In particular this proves \eqref{A:eq:14.49}.
The coefficient event is independent of $K$; on that event the common form bounds allow every compact $K\Subset\mathbb C\setminus\mathbb R$. The form bounds were established from
the coefficient estimates.
\end{proof}

\begin{proof}[Proof of Theorem~\ref{A:res:15.4}]
Theorem~\ref{A:res:16.2} supplies common form constants for the three operators in \eqref{A:eq:15.54}. For the diagonal part, \eqref{A:eq:15.48} and \eqref{A:eq:15.50} give
\[
\sum_{j=1}^{M}
(1+x_j)^{-7/2}
|\widetilde\rho_j+\widetilde\varepsilon_j^{\mathrm{fr}}|^2
\le
C_\omega h^4.
\]
The diagonal estimate \eqref{A:eq:3.49} therefore yields
\[
\sup_{z\in K}
\left\|
(H_n^{\mathrm G}-z)^{-1}
-
(\widetilde H_n-z)^{-1}
\right\|_{\mathfrak S_2}^2
\le
C_{\omega,V,\beta,K}h.
\]
For the first-order part, \eqref{A:eq:15.51} and the weighted first-order perturbation estimate give
\[
\sup_{z\in K}
\left\|
(\widetilde H_n-z)^{-1}
-
(H_n^{B,\mathrm{mid}}-z)^{-1}
\right\|_{\mathfrak S_2}^2
\le
C_{\omega,V,\beta,K}h.
\]

Put $R_G(z)=(H_n^{\mathrm G}-z)^{-1}$, $\widetilde R(z)=(\widetilde H_n-z)^{-1}$, and $R_{B,\mathrm{mid}}(z)=(H_n^{B,\mathrm{mid}}-z)^{-1}$. The two estimates just proved give
\[
 \sup_{z\in K}\|R_G(z)-R_{B,\mathrm{mid}}(z)\|_{\mathfrak S_2}^2
 \le2\sup_{z\in K}\|R_G(z)-\widetilde R(z)\|_{\mathfrak S_2}^2
   +2\sup_{z\in K}\|\widetilde R(z)-R_{B,\mathrm{mid}}(z)\|_{\mathfrak S_2}^2
 \le C_{\omega,K}h.
\]
The common form bounds used in those estimates are supplied by Theorem~\ref{A:res:16.2} from the coefficient estimates.

The squared Hilbert--Schmidt norm of the resolvent difference arising from the covariance and spectral factorization errors, divided by $h$, tends to zero almost surely. Let
$W_j=(1+x_j)^{-7/2}$. Propositions~\ref{A:res:A.20} and~\ref{A:res:15.2} give, after disintegration,
\[
 \left\|\sum_{I_n^{\mathrm{mid}}}W_j|\rho_j|^2\right\|_{L^p}\le C_ph^6,
 \qquad
 \left\|\sum_{I_n^{\mathrm{mid}}}W_j|\varepsilon_j^{\mathrm{fr}}|^2
          \right\|_{L^p}\le C_ph^5.
\]
For $p>3/\epsilon$, Markov's inequality at thresholds $h^{6-\epsilon}$ and $h^{5-\epsilon}$ gives summable probabilities $C_ph^{p\epsilon}$. Zero extension does not change either
sum. Thus
\[
 \sup_{z\in K}\|R_G(z)-\widetilde R(z)\|_{\mathfrak S_2}^2
 \le C_{\omega,K,\epsilon}(h^{3-\epsilon}+h^{2-\epsilon})=o(h),
 \qquad 0<\epsilon<1.
\]
The squared Hilbert--Schmidt norm of $\widetilde R-R_{B,\mathrm{mid}}$ is controlled by the weighted squared norm of the discrete primitive, which is bounded by a constant times
$h$. The maximum estimate of Lemma~\ref{A:res:freezing-max} is used for the form bound only; substituting it in the weighted error sum would lose this improvement.
\end{proof}

\begin{proof}[Proof of Theorem~\ref{A:res:15.5}]
Write $R_0(z)=(H_n^{B,\mathrm{mid}}-z)^{-1}$ and $R_1(z)=(H_n^{B,\mathrm{reg}}-z)^{-1}$. Their operator difference is $W_n^{\mathrm L,sc}$. Its support has length
$J_n^{\mathrm L}=O(\log n)$ and its entries are bounded by $C_\omega h^{-1/2}\sqrt{\log n}$. The form hypotheses of Corollary~\ref{A:res:3.19A} follow from
Theorem~\ref{A:res:16.2}. The resolvent identity and that corollary give
\[
 \sup_{z\in K}\|R_0(z)-R_1(z)\|_{\mathfrak S_2}^2
 \le C_{\omega,K}h^4(J_n^{\mathrm L}+1)^4
                  (h^{-1/2}\sqrt{\log n})^2
 \le C_{\omega,K}h^3(\log n)^5.
\]
Its ratio to $h$ is bounded by $C_{\omega,K}h^2(\log n)^5\to0$.
\end{proof}

\begin{proof}[Proof of Theorem~\ref{A:res:15.6}]
The operator difference is \(H_n^{B,\mathrm{reg}}-H_n^{B,\circ} = h^{-2}M_{\mathbf 1_{I_n^{\mathrm{far}}}(\zeta_n-\zeta_n^B)}\). Both operators satisfy the upper and lower form
bounds established in Section~\ref{A:sec:quadratic-form-bounds}. The resolvent identity and the diagonal estimate \eqref{A:eq:3.49} give
\[
\left\| (H_n^{B,\mathrm{reg}}-z)^{-1} - (H_n^{B,\circ}-z)^{-1} \right\|_{\mathfrak S_2}^2 \le C_{\omega,K}h^{-3} \sum_{j\in I_n^{\mathrm{far}}} |\zeta_j-\zeta_j^B|^2 (1+x_j)^{-7/2} = C_{\omega,K}\mathcal K_{7,n}.
\]
Equation~\eqref{A:eq:15.63} bounds this by $C_{\omega,K}h^4\log n$; division by $h$ gives $C_{\omega,K}h^3\log n\to0$.
\end{proof}

We now write the coefficient comparison on the original finite space. Let $E_n=P_{I_n}^*U_{r_n}^*$, so that $E_n^*E_n=I_M$, and put
\[
 \begin{gathered}
 A_0=\widetilde H_n^{\mathrm{red}},\quad A_1=H_n^{\mathrm{diag}},\quad
 A_2=H_n^{\mathrm G},\quad A_3=H_n^{\mathrm G}-Q_n^{\mathrm{fr}},\quad
 A_4=\widetilde H_n,\\
 A_5=H_n^{B,\mathrm{mid}},\quad A_6=H_n^{B,\mathrm{reg}},\quad
 A_7=H_n^{B,\circ},\qquad R_j(z)=(A_j-z)^{-1}.
 \end{gathered}
\]
The definitions in \eqref{A:eq:common-assembly} show that
\[
 \begin{aligned}
 A_0-A_1&=\Gamma_\eta^-,& A_1-A_2&=Q_n^{\mathrm{tr}},&
 A_2-A_3&=Q_n^{\mathrm{fr}},& A_3-A_4&=Q_n^\rho,\\
 A_4-A_5&=\Gamma_{\widetilde e},& A_5-A_6&=Q_n^{\mathrm L},&
 A_6-A_7&=Q_n^{\mathrm{far}}.
 \end{aligned}
\]
In particular $A_3$ has the same form bounds: it is obtained from $H_n^{B,\circ}$ by selecting a subset of the seven perturbations estimated in the proof of
Theorem~\ref{A:res:16.2}.

Define the full boundary error by
\[
 D_{\partial,n}(z)=(\mathcal A_{n,V,\beta}-z)^{-1}-E_nR_0(z)E_n^*.
\]
Equation~\eqref{A:eq:12.43} estimates this entire difference, including the deleted blocks and both off-diagonal blocks, by
$\sup_{z\in K}\|D_{\partial,n}(z)\|_{\mathfrak S_2}^2=d_{\partial,n}=o_\omega(h)$, after increasing the bound if necessary. The telescoping identity is
\[
 (\mathcal A_{n,V,\beta}-z)^{-1}-E_nR_7(z)E_n^*
 =D_{\partial,n}(z)+\sum_{j=0}^6 E_n(R_j(z)-R_{j+1}(z))E_n^*.
\]
For every finite matrix $T$, the isometry gives
\[
 \|E_nTE_n^*\|_{\mathfrak S_2}^2
 =\operatorname{Tr}(E_nT^*TE_n^*)=\operatorname{Tr}(T^*T).
\]
There are eight terms in this identity. Cauchy--Schwarz in the Hilbert space $\mathfrak S_2$ therefore bounds its squared norm, uniformly on $K$, by
\[
 8d_{\partial,n}+8\sum_{j=0}^6\sup_{z\in K}
                      \|R_j(z)-R_{j+1}(z)\|_{\mathfrak S_2}^2.
\]
Fix $0<\epsilon<1$. On the intersection of the coefficient events already constructed, the seven summands, in order, are bounded by one finite random constant times
\[
 h,\quad h^{2-\epsilon},\quad h^{2-\epsilon},\quad h^{3-\epsilon},
 \quad h,\quad h^3(\log n)^5,\quad h^4\log n.
\]
The two first-order forms give the two terms of order $h$. The diagonal bounds give the transport, covariance and factorization terms with the powers displayed. The last two
estimates are Theorems~\ref{A:res:15.5} and~\ref{A:res:15.6}. Their sum is $O_\omega(h)$, since $2-\epsilon>1$, $3-\epsilon>1$ and $h\asymp n^{-1/3}$. Taking the maximum of the
finitely many starting indices proves \eqref{A:eq:15.65}. The estimates hold on the same joint coupling; combining them requires no independence between the errors.

\subsubsection{Comparison of the deterministic coefficients}\label{A:sec:comparison-of-the-deterministic-coefficients}

We now keep the Brownian potential fixed and compare the coefficients at the conditional minimizer with the local profile, then the local profile with its linear edge expansion.
The conditional minimizer depends on the exterior coefficient values. Fixing these values and the Brownian path makes the comparison deterministic.

Neither replacement is uniformly small over the whole interval. Near the left cut we use the logarithmic support and the bound for the resolvent columns at the left boundary. In a
fixed regular macroscopic region the minimizing coefficients are $O(n^{-1})$ apart. Near the right boundary we use an $O(h)$ bound for the whole resolvent tail, beyond the region
of the regular Taylor expansion. We obtain
\[
\sup_{z\in K}
\left\|
(H_n^{B,\circ}-z)^{-1}
-(\mathcal B_{M,h,\beta}^B-z)^{-1}
\right\|_{\mathfrak S_2}^2
\le C_{\omega,V,\beta,K}h.
\]
Both operators act on the same retained coordinate space and use the same Brownian diagonal.

Use $H_n^\circ,H_n^\dagger$ and their coefficient notation from Subsection~\ref{A:sec:quadratic-form-bounds}. Write $\omega^B=\omega_h$ for the coefficients defined by the Brownian
increments in \eqref{A:eq:brownian-difference}. Thus $H_n^{B,\circ}=H_n^\circ+M_{\omega^B}$ and $H_n^{B,\dagger}=H_n^\dagger+M_{\omega^B}$.

Proposition~\ref{A:res:3.5}, applied with the common form constants of Theorem~\ref{A:res:16.2}, supplies \eqref{A:eq:derivative-HS} for each resolvent used below.

The exact conditional minimizer displacement has two different boundary scales. With the reindexed interior coordinates starting at one,
Section~\ref{A:sec:decomposition-of-the-general-jacobi-coefficients} gives
\[
|Y_{n,j}^\circ-Y_{n,j}^\dagger| \le{} C_\omega n^{-1} + C_\omega n^{-1/2}\sqrt{\log n}\,\mathrm{e}^{-c(j-1)} + C_\omega \mathrm{e}^{-c(M-j)}.
\]
We use the two boundary terms on their respective intervals below.

Use $L_n^{\mathrm L},J_n^{\mathrm L}$ from \eqref{A:eq:13.17b}, and fix the terminal cutoff
\[
J=J_n^{\mathrm{cp}}
=
\max\left\{j:\ s^b_{n,\kappa(j)+R_V}\le s_{\mathrm{cp}}\right\},
\quad 0<s_{\mathrm{cp}}<1.
\]
These choices give \(J\asymp n, \allowbreak\, x_J\asymp h^{-2}, \allowbreak\, M-J\asymp n,\) and, after taking
\(A_{\mathrm L}\) large enough,
\[
\max_{J_n^{\mathrm L}+1\le j\le J+1}
\left(
|A_j^\circ-a_j^\dagger|
+
|B_j^\circ-b_j^\dagger|
\right)
\le
C_\omega n^{-1}
=
C_\omega h^3.
\eqtag{17.15}
\]
On the left boundary strip we have
\[
\max_{1\le j\le J_n^{\mathrm L}+R_V}
\left(
|A_j^\circ-a_j^\dagger|
+
|B_j^\circ-b_j^\dagger|
\right)
\le
C_\omega n^{-1/2}\sqrt{\log n}
=
C_\omega h^{3/2}\sqrt{\log n}.
\eqtag{17.16}
\]

Let \(H_{\circ,L}^{B}\) and \(H_{\dagger,L}^{B}\) be the Dirichlet restrictions of \(H_n^{B,\circ}\) and
\(H_n^{B,\dagger}\) to \(\{1,\ldots,J\}\), and let \(H_{\circ,R}^{B}\) and \(H_{\dagger,R}^{B}\) be their restrictions
to \(\{J+1,\ldots,M\}\), with the global rescaled coordinate retained. Write
\[
H_{\circ}^{B,\oplus}
=H_{\circ,L}^{B}\oplus H_{\circ,R}^{B},
\quad
H_{\dagger}^{B,\oplus}
=H_{\dagger,L}^{B}\oplus H_{\dagger,R}^{B}.
\]

For either operator $H$ above, let $H^\oplus$ denote the direct sum of its two principal blocks and put $t_J=H_{J,J+1}$. The difference is exactly
\[
 H-H^\oplus=t_J(e_Je_{J+1}^*+e_{J+1}e_J^*),\qquad |t_J|\le C_\omega h^{-2}.
\]
Its diagonal is zero, so both principal blocks have the original Dirichlet boundary energies. All four principal blocks inherit the form bounds of Theorem~\ref{A:res:16.2} by zero
extension. Let $R(z)=(H-z)^{-1}$ and $R^\oplus(z)=(H^\oplus-z)^{-1}$. The resolvent identity gives
\[
 R(z)-R^\oplus(z)
 =-t_J\left[(R(z)e_J)(R^\oplus(\bar z)e_{J+1})^*
           +(R(z)e_{J+1})(R^\oplus(\bar z)e_J)^*\right].
\]
Since $x_J,x_{J+1}\asymp h^{-2}$, each squared vector norm on the right is at most $C_{\omega,K}h^4$. A rank-one operator $uv^*$ has squared Hilbert--Schmidt norm $\|u\|^2\|v\|^2$;
hence
\[
 \|R(z)-R^\oplus(z)\|_{\mathfrak S_2}^2
 \le2|t_J|^2\left[
  \|R(z)e_J\|^2\|R^\oplus(\bar z)e_{J+1}\|^2
 +\|R(z)e_{J+1}\|^2\|R^\oplus(\bar z)e_J\|^2\right]
 \le C_{\omega,K}h^4.
\]
Applied to the two operators this proves
\[
\sup_{z\in K}
\left\|
(H_n^{B,\circ}-z)^{-1}
-
(H_{\circ}^{B,\oplus}-z)^{-1}
\right\|_{\mathfrak S_2}^2
\le C_{\omega,K}h^4,
\eqtag{17.18}
\]
and
\[
\sup_{z\in K}
\left\|
(H_n^{B,\dagger}-z)^{-1}
-
(H_{\dagger}^{B,\oplus}-z)^{-1}
\right\|_{\mathfrak S_2}^2
\le C_{\omega,K}h^4.
\eqtag{17.19}
\]

In $H_{\circ,L}^{B}-H_{\dagger,L}^{B}$ the Brownian terms cancel. We take this difference in matrix coordinates; the auxiliary boundary coefficients cancel by
Lemma~\ref{A:res:auxiliary-boundary}. Split the finite-range difference as
\(H_{\circ,L}^{B}-H_{\dagger,L}^{B} = W_n^{\mathrm L}+W_n^{\mathrm{mid}},\) where \(W_n^{\mathrm L}\) contains every
matrix entry meeting \(\{1,\ldots,J_n^{\mathrm L}\}\), and \(W_n^{\mathrm{mid}}\) contains the remaining entries. The matrix $W_n^{\mathrm L}$ vanishes for $|i-j|>R_V$ and is
supported in \(\{1,\ldots,J_n^{\mathrm L}+R_V\}\). By \eqref{A:eq:17.16}, after multiplication by the scaling factor \(h^{-2}\),
\(\max_{i,j}|(W_n^{\mathrm L})_{ij}| \le C_\omega h^{-1/2}\sqrt{\log n}.\) For \(0\le t\le1\), the interpolation \(H_{\dagger,L}^{B}+tW_n^{\mathrm L}\) satisfies the upper and
lower form bounds uniformly: the discrete Hardy estimate \eqref{A:eq:3.58g} gives
\[
|\langle u,W_n^{\mathrm L}u\rangle| \le C_\omega h^{-1/2}\sqrt{\log n} \sum_{j\le J_n^{\mathrm L}+R_V}|u_j|^2 \le C_\omega h^{3/2}(\log n)^{5/2} \sum_{k<J_n^{\mathrm L}+R_V}|D_hu_k|^2 =o_\omega(1)\,\ell_h[u].
\]
Corollary~\ref{A:res:3.19A} then gives, for any two resolvents along this interpolation,
\[
\sup_{z\in K}
\|R_a(z)W_n^{\mathrm L}R_b(z)\|_{\mathfrak S_2}^2
\le
C_{\omega,K}h^3(\log n)^5
=o_\omega(h).
\eqtag{17.23}
\]
For the remaining matrix $S=W_n^{\mathrm{mid}}$, every entry is $O_\omega(h)$ by \eqref{A:eq:17.15}, and $S$ is tridiagonal. In particular
$\|S\|_{\mathrm{op}}\le3\max_{i,j}|S_{ij}|=O_\omega(h)$, so its interpolation preserves the form bounds. We also record its kinetic and potential parts separately. The difference
across the splitting index is treated separately in the formulas below. Define
\[
 d_{j,n}^c=-h^2S_{j,j+1}\quad(1\le j<J),\qquad d_{0,n}^c=d_{J,n}^c=0,
 \qquad r_{j,n}^c=S_{jj}-h^{-2}(d_{j-1,n}^c+d_{j,n}^c).
\]
The off-diagonal entries of $D_h^*M_{d_n^c}D_h$ are $-h^{-2}d_{j,n}^c=S_{j,j+1}$, and the definition of $r_{j,n}^c$ matches the diagonal. This proves the exact decomposition
\[
W_n^{\mathrm{mid}}
=
D_h^*M_{d_n^c}D_h+M_{r_n^c},
\quad
\|d_n^c\|_\infty\le C_\omega h^3,
\quad
\|r_n^c\|_\infty\le C_\omega h,
\eqtag{17.24}
\]
with both coefficients zero on the left boundary strip. The values $d_{0,n}^c=d_{J,n}^c=0$ are part of this exact representation of $S$; they do not alter either original operator.
Let $R_a,R_b$ be the resolvents before and after this middle perturbation. Proposition~\ref{A:res:3.5} and the bound for $D_hR_a(z)$ in operator norm give
\[
\|R_aD_h^*M_{d_n^c}D_hR_b\|_{\mathfrak S_2}^2 \le C_{\omega,K}\|d_n^c\|_\infty^2 \|D_hR_b\|_{\mathfrak S_2}^2 \le C_{\omega,K}h^6h^{-1} = C_{\omega,K}h^5.
\eqtag{17.25}
\]
For the diagonal part, Proposition~\ref{A:res:3.16} gives
\[
\|R_aM_{r_n^c}R_b\|_{\mathfrak S_2}^2 \le C_{\omega,K}h \sum_{j\le J}|r_{j,n}^c|^2(1+x_j)^{-7/2} \le C_{\omega,K}h^3 \sum_{j\le J}(1+x_j)^{-7/2} \le C_{\omega,K}h^2.
\eqtag{17.26}
\]
Combining \eqref{A:eq:17.23}--\eqref{A:eq:17.26} by applying the resolvent identity twice proves
\[
\sup_{z\in K}
\left\|
(H_{\circ,L}^{B}-z)^{-1}
-
(H_{\dagger,L}^{B}-z)^{-1}
\right\|_{\mathfrak S_2}^2
\le
C_{\omega,K}h^2+o_\omega(h).
\eqtag{17.27}
\]

The regular block contributes $o(h)$. We estimate the right block as a whole, at order $h$. Its support lies in \(x\ge x_J\asymp h^{-2}\), so Proposition~\ref{A:res:3.20} gives
\[
\sup_{z\in K}
\left(
\|(H_{\circ,R}^{B}-z)^{-1}\|_{\mathfrak S_2}^2
+
\|(H_{\dagger,R}^{B}-z)^{-1}\|_{\mathfrak S_2}^2
\right)
\le
C_{\omega,K}x_J^{-1/2}
\le C_{\omega,K}h.
\]
\[
\sup_{z\in K}
\left\|
(H_{\circ,R}^{B}-z)^{-1}
-
(H_{\dagger,R}^{B}-z)^{-1}
\right\|_{\mathfrak S_2}^2
\le C_{\omega,K}h.
\eqtag{17.27b}
\]
The left and right terms act on orthogonal summands. Combining \eqref{A:eq:17.18}, \eqref{A:eq:17.19}, \eqref{A:eq:17.27}, and \eqref{A:eq:17.27b} proves the following proposition.

\par\Needspace{9\baselineskip}
\begin{proposition}\label{A:res:17.2}
For every compact \(K\Subset\mathbb C\setminus\mathbb R\), almost surely, for all sufficiently large \(n\),
\[
\sup_{z\in K}
\left\|
(H_n^{B,\circ}-z)^{-1}
-
(H_n^{B,\dagger}-z)^{-1}
\right\|_{\mathfrak S_2}^2
\le
C_{\omega,V,\beta,K}h.
\]
\end{proposition}

\begin{proof}
Cut both matrices at $J=J_n^{\mathrm{cp}}$ and add and subtract the resolvents of the resulting direct sums. The two cut errors are \eqref{A:eq:17.18}--\eqref{A:eq:17.19}. On the
direct sums the squared Hilbert--Schmidt norm is the sum of the squared norms on the two orthogonal blocks, bounded by \eqref{A:eq:17.27} and \eqref{A:eq:17.27b}. If $R_1,R_0$ are
the two full reduced resolvents and $R_1^\oplus,R_0^\oplus$ their cut resolvents, then
\[
 \|R_1-R_0\|_{\mathfrak S_2}^2
 \le3\bigl(\|R_1-R_1^\oplus\|_{\mathfrak S_2}^2
       +\|R_1^\oplus-R_0^\oplus\|_{\mathfrak S_2}^2
       +\|R_0^\oplus-R_0\|_{\mathfrak S_2}^2\bigr)
 \le C_{\omega,K}h
\]
uniformly for $z\in K$, by the four cited bounds.
\end{proof}

The estimates for the local minimizer from Section~\ref{A:sec:the-local-minimizer-and-edge-scaling} must be written in the recentered interior coordinate. The original reduced
interior block begins at index \(r_n+1\). Its displacement in the coordinate $x$ is \(\delta_n=r_nh.\) Since \(r_n=\lceil A_0\log n\rceil\), we have \(\delta_n=A_0h\log n+O(h)\)
and
\(\delta_n^2=o(h)\). In the original coordinate, the deterministic linear potential is \(x_j+\delta_n.\) After
recentering the lattice at the new Dirichlet boundary, Section~\ref{A:sec:the-local-minimizer-and-edge-scaling} gives the expansion
\(V_{j,n}^{\dagger} = x_j+\delta_n+\frac{h}{\beta} + r_{j,n}^{\dagger}\). On every fixed regular edge region \(s\le s_0<1\), \(|r_{j,n}^{\dagger}| \le C h^2(1+x_j)^2.\) The
deterministic off-diagonal coefficient difference is \(d_{j,n}^{\dagger} = c_{j,n}^{\dagger}-1.\) On the same regular region, \(|d_{j,n}^{\dagger}| \le C h^2(1+x_j),\) and
\(|d_{j,n}^{\dagger}-d_{j-1,n}^{\dagger}| \le C h^3.\) Therefore, in the regular region,
\(H_n^\dagger = D_h^*D_h+X_h + \left(\delta_n+\frac{h}{\beta}\right)I + D_h^*M_{d_n^\dagger}D_h + M_{r_n^\dagger}\). Adding the Brownian potential \(M_{\omega^B}\), we get
\[
H_n^{B,\dagger}
=
\mathcal B_{M,h,\beta}^B
+
\left(\delta_n+\frac{h}{\beta}\right)I
+
D_h^*M_{d_n^\dagger}D_h
+
M_{r_n^\dagger}
\]
in the regular region. The scalar shift \(\lambda_n:=\delta_n+\frac{h}{\beta}\) has size
\[
|\lambda_n|\le C h\log n.
\eqtag{17.37}
\]
In particular, $\lambda_n^2=O(h^2(\log n)^2)=o(h)$. The resolvent identity leading to \eqref{A:eq:12.47}, with $\delta=\lambda_n$, gives the same squared bound. This comparison
concerns the deterministic potential; it does not translate the Brownian path. In the regular Taylor expansion, \(s_{j,n}=(b_\star/\tau_\star)h^2(x_j+\lambda_n)\) gives
\(h^{-2}s_{j,n}^2\le Ch^2(1+x_j)^2\), since \(|\lambda_n|\le1\) for all sufficiently large \(n\). The only linear
correction is the displayed scalar \(\lambda_n\); bounded derivatives of the local minimizer and
\(s_{j+1,n}-s_{j,n}=1/n\) give the same bound \(Ch^3\) on differences of neighboring off-diagonal coefficients. The
bounds are uniform on compact positive $\beta$-intervals.

Choose a fixed number \(s_0\in(0,1).\) Let \(J=J_n(s_0)\) be the largest reindexed vertex such that the corresponding original macroscopic coordinate is at most \(s_0\). Then
\(J\asymp n, \allowbreak\, x_J\asymp Jh\asymp nh\asymp h^{-2}.\) Work first on the interior block \(\{1,\ldots,J\}\) with Dirichlet boundary conditions
\(u_0=0, \allowbreak\, u_{J+1}=0.\) Let \(H_J^{B,\dagger}\) be the principal restriction of \(H_n^{B,\dagger}\) to this block. In the representation of its quadratic form as a
difference term plus a potential, retain the local minimizing value at the cut edge in both terms, so the slowly varying estimates below use the same auxiliary boundary
coefficients. This is the same Dirichlet matrix for any paired choice of auxiliary boundary coefficients, by Lemma~\ref{A:res:auxiliary-boundary}. For the finite difference matrix
driven by the Brownian increments in \eqref{A:eq:brownian-difference}, put $\mathcal B_{J,h,\beta}^{B,\lambda} =\mathcal B_{J,h,\beta}^{B}+\lambda_n I$. Then
\(H_J^{B,\dagger} - \mathcal B_{J,h,\beta}^{B,\lambda} = D_h^*M_{d_n^\dagger}D_h+M_{r_n^\dagger}\). The coefficients satisfy
\[
|d_{j,n}^{\dagger}|\le C h^2(1+x_j),
\eqtag{17.43}
\]
\[
|d_{j,n}^{\dagger}-d_{j-1,n}^{\dagger}|\le C h^3,
\eqtag{17.44}
\]
and $|r_{j,n}^\dagger|\le Ch^2(1+x_j)^2$. The bound on $d^\dagger$ is only $O(1)$ at $x\asymp h^{-2}$. Using it uniformly in a product of two differentiated resolvents would give
only $O(h^{-1})$ from \eqref{A:eq:derivative-HS}. Instead we retain its dependence on $x_j$, apply the discrete product rule, and eliminate one second difference by the resolvent
equation.

The diagonal perturbation estimate uses the weighted sum
\[
 h^4\int_0^{c h^{-2}}(1+x)^4(1+x)^{-7/2}\,\mathrm{d}x
 =h^4\int_0^{c h^{-2}}(1+x)^{1/2}\,\mathrm{d}x
 =O(h).
\]
For the off-diagonal difference we use Lemma~\ref{A:res:product-rule} and the equation $(\mathcal B_{J,h,\beta}^{B,\lambda}-z)R_0(z)=I$.

\begin{proposition}\label{A:res:17.3}
For every compact \(K\Subset\mathbb C\setminus\mathbb R,\) almost surely, for all sufficiently large \(n\),
\[
\sup_{z\in K}
\left\|
(H_J^{B,\dagger}-z)^{-1}
-
(\mathcal B_{J,h,\beta}^{B,\lambda}-z)^{-1}
\right\|_{\mathfrak S_2}^2
\le
C_{\omega,V,\beta,K}h.
\eqtag{17.46}
\]
\end{proposition}

\begin{proof}
Write \(R_1(z)=(H_J^{B,\dagger}-z)^{-1},\) and \(R_0(z)=(\mathcal B_{J,h,\beta}^{B,\lambda}-z)^{-1}.\) By the resolvent identity,
\[
R_1(z)-R_0(z)
=
-R_1(z)
\left(
D_h^*M_{d^\dagger}D_h+M_{r^\dagger}
\right)
R_0(z).
\eqtag{17.47}
\]
Let \(T_h=D_h^*D_h\). Lemma~\ref{A:res:product-rule}, on the edge space indexed by \(0,\ldots,J\), gives
\[
D_h^*M_{d^\dagger}D_h
=
M_aT_h
-
M_bJ_-^*D_h,
\eqtag{17.48}
\]
where
\[
a_j=d_{j-1,n}^\dagger,
\eqtag{17.49}
\]
and
\[
b_j=h^{-1}(d_{j,n}^\dagger-d_{j-1,n}^\dagger).
\eqtag{17.50}
\]
By \eqref{A:eq:17.43}--\eqref{A:eq:17.44},
\[
|a_j|\le C h^2(1+x_j),
\eqtag{17.51}
\]
and
\[
|b_j|\le C h^2.
\eqtag{17.52}
\]
The shifted finite difference operator driven by the Brownian increments satisfies \(\mathcal B_{J,h,\beta}^{B,\lambda} = T_h+X_h+M_{\omega^B}+\lambda_n I.\) Therefore,
\[
T_hR_0(z)
=
I+(z-\lambda_n)R_0(z)
-
M_{x+\omega^B}R_0(z).
\eqtag{17.53}
\]
Substituting \eqref{A:eq:17.48} and \eqref{A:eq:17.53} into the resolvent identity yields
\[
R_1D_h^*M_{d^\dagger}D_hR_0 ={} R_1M_a+(z-\lambda_n)R_1M_aR_0-R_1M_{a x}R_0
 -R_1M_{a\omega^B}R_0-R_1M_bJ_-^*D_hR_0.
\eqtag{17.54}
\]
Although the Brownian diagonal cancelled in $H_J^{B,\dagger}- \mathcal B_{J,h,\beta}^{B,\lambda}$, it reappears in \eqref{A:eq:17.54} when $T_hR_0$ is eliminated. Its weighted
increment estimate is therefore still needed. The first term has only one resolvent, so we begin with the multiplier estimate
\[
\|R_1M_a\|_{\mathfrak S_2}^2
=
\|M_aR_1(\bar z)\|_{\mathfrak S_2}^2
\le
C_{\omega,K}
h\sum_{j\le J}|a_j|^2(1+x_j)^{-3/2}.
\eqtag{17.55}
\]
Using \eqref{A:eq:17.51},
\[
h\sum_{j\le J}|a_j|^2(1+x_j)^{-3/2} \le C h\sum_{j\le J} h^4(1+x_j)^2(1+x_j)^{-3/2} = C h^5\sum_{j\le J}(1+x_j)^{1/2}.
\eqtag{17.56}
\]
Since \(x_J\asymp h^{-2},\) we have
\[
h\sum_{j\le J}(1+x_j)^{1/2}
\le
C\int_0^{Ch^{-2}}(1+x)^{1/2}\,\mathrm{d}x
\le
C h^{-3}.
\eqtag{17.57}
\]
Dividing \eqref{A:eq:17.57} by \(h\) gives \(\sum_{j\le J}(1+x_j)^{1/2}\le Ch^{-4}\), and hence
\[
\|R_1M_a\|_{\mathfrak S_2}^2
\le
C_{\omega,K}h.
\eqtag{17.58}
\]
Since \(|z-\lambda_n|\le C_K\) for all large \(n\), and \(\|R_0(z)\|_{\mathrm{op}}\le C_K,\) we also get
\[
\|(z-\lambda_n)R_1M_aR_0\|_{\mathfrak S_2}^2
\le
C_{\omega,K}h.
\eqtag{17.59}
\]
Retaining the second resolvent gives the sharper estimate
\[
 \|(z-\lambda_n)R_1M_aR_0\|_{\mathfrak S_2}^2
 \le C_{\omega,K}h\sum_{j\le J}|a_j|^2(1+x_j)^{-7/2}
 \le C_{\omega,K}h^4,
\]
since $h\sum_j(1+x_j)^{-3/2}\le C$. The order-$h$ term in \eqref{A:eq:17.58} comes from the part $R_1M_a$ with only one resolvent. Next, \(|a_jx_j| \le C h^2(1+x_j)^2.\) By the
diagonal estimate \eqref{A:eq:3.49},
\[
\|R_1M_{a x}R_0\|_{\mathfrak S_2}^2 \le C_{\omega,K} h\sum_{j\le J} h^4(1+x_j)^4(1+x_j)^{-7/2} = C_{\omega,K} h^5\sum_{j\le J}(1+x_j)^{1/2} \le C_{\omega,K}h.
\eqtag{17.60}
\]
The maximum estimate \eqref{A:eq:15.59}, multiplied by $h^{-2}$, also gives
\[
  \max_{1\le j\le M}|\omega_j^B|
  \le C_{\omega,\beta}h^{-1/2}\sqrt{\log(h^{-1})}
\eqtag{17.61}
\]
for all sufficiently large $n$. For the weighted product we estimate the Brownian increments separately, using their moments instead of this maximum. Write
$\Delta_jB=B(jh)-B((j-1)h)$. Since $\omega_j^B=2\Delta_jB/(\sqrt\beta h)$, \eqref{A:eq:17.51} gives
\[
 |a_j\omega_j^B|^2
 \le C_\beta h^2(1+x_j)^2|\Delta_jB|^2.
\eqtag{17.62}
\]
The diagonal estimate \eqref{A:eq:3.49} consequently yields
\[
 \|R_1M_{a\omega^B}R_0\|_{\mathfrak S_2}^2
 \le C_{\omega,\beta,K}h^3
       \sum_{j\le J}(1+x_j)^{-3/2}|\Delta_jB|^2.
\eqtag{17.63}
\]
The sum is at most $CQ(B)$ by the common Brownian increment estimate of Subsection~\ref{A:sec:finite-difference-approximation-of-the-stochastic-airy-operator}. Thus
\[
 \|R_1M_{a\omega^B}R_0\|_{\mathfrak S_2}^2
 \le C_{\omega,\beta,K}Q(B)h^3=o_\omega(h).
\eqtag{17.64}
\]
For the first-order remainder, \(|b_j|\le Ch^2\). By \eqref{A:eq:restricted-multiplier}, \(\|R_1M_bJ_-^*\|_{\mathfrak S_2}^2 =\|R_1M_b\|_{\mathfrak S_2}^2\). The weighted
multiplication bound, applied at \(\overline z\) to the adjoint, therefore gives
\[
\|R_1M_bJ_-^*\|_{\mathfrak S_2}^2
\le
C_{\omega,K}
h\sum_{j\le J}
h^4(1+x_j)^{-3/2}
\le
C_{\omega,K}h^4.
\eqtag{17.65}
\]
Using \(\|D_hR_0(z)\|_{\mathrm{op}}\le C_{\omega,K},\) we get
\[
\|R_1M_bJ_-^*D_hR_0\|_{\mathfrak S_2}^2
\le
C_{\omega,K}h^4.
\eqtag{17.66}
\]
Finally, the diagonal remainder satisfies \(|r_j^\dagger|\le Ch^2(1+x_j)^2.\) Therefore,
\[
\|R_1M_{r^\dagger}R_0\|_{\mathfrak S_2}^2 \le C_{\omega,K} h\sum_{j\le J} h^4(1+x_j)^4(1+x_j)^{-7/2} \le C_{\omega,K}h.
\eqtag{17.67}
\]
Combining \eqref{A:eq:17.47}--\eqref{A:eq:17.67}, we obtain \(\sup_{z\in K} \|R_1(z)-R_0(z)\|_{\mathfrak S_2}^2 \le C_{\omega,V,\beta,K}h\). Estimate \eqref{A:eq:17.46} follows.

Before substituting $J\asymp h^{-3}$, put $L=(J+1)h$. The integral in \eqref{A:eq:17.57} is bounded by $C(1+L)^{3/2}$. The estimates for the six terms of
\eqref{A:eq:17.47}--\eqref{A:eq:17.54} therefore sum to
\[
 C_{\omega,V,\beta,K}
       \{h^4(1+L)^{3/2}+h^4+h^3Q(B)\}.
\]
The term with $r^\dagger$ has the same $h^4(1+L)^{3/2}$ bound as the term with $ax$. Since $L\le C_{s_0,V}h^{-2}$, this expression is $O(h)$ with a constant independent of $J$ and
$n$. Thus the estimate holds on the full regular macroscopic interval, whose length grows in the $x$ coordinate.
\end{proof}

Apply the resolvent identity to $\mathcal B_{J,h,\beta}^{B,\lambda}=\mathcal B_{J,h,\beta}^{B}+\lambda_n I$. Its Hilbert--Schmidt norm is uniformly bounded by
Theorem~\ref{A:res:16.2} and Proposition~\ref{A:res:3.5}, so
\[
 \sup_{z\in K}\|(\mathcal B_{J,h,\beta}^{B,\lambda}-z)^{-1}
                  -(\mathcal B_{J,h,\beta}^B-z)^{-1}\|_{\mathfrak S_2}^2
 \le C_{\omega,K}\lambda_n^2
 \le C_{\omega,K}h^2(\log n)^2=o_\omega(h).
 \eqtag{17.70}
\]
Combining Proposition~\ref{A:res:17.3} with \eqref{A:eq:17.70}, we obtain the interior comparison without the scalar shift:
\[
\sup_{z\in K}
\left\|
(H_J^{B,\dagger}-z)^{-1}
-
(\mathcal B_{J,h,\beta}^{B}-z)^{-1}
\right\|_{\mathfrak S_2}^2
\le
C_{\omega,V,\beta,K}h.
\eqtag{17.71}
\]

The estimates above were proved only on the regular region \(s\le s_0\). We pass from \(\{1,\ldots,J\}\) to the full reduced interior interval \(\{1,\ldots,M\}\) by controlling the
remaining tail. Let \(H_n^{B,\dagger,\oplus} = H_J^{B,\dagger}\oplus H_{\mathrm{tail}}^{B,\dagger}\) be the operator obtained from \(H_n^{B,\dagger}\) by cutting the edge between
\(J\) and \(J+1\). Let \(\mathcal B_{M,h,\beta}^{B,\oplus} = \mathcal B_{J,h,\beta}^{B}\oplus \mathcal B_{\mathrm{tail},h,\beta}^{B}\) be the corresponding cut finite difference
approximation of the stochastic Airy operator. For both operators, use the rank-two identity preceding \eqref{A:eq:17.18}. Their cut entries have modulus at most $C_\omega h^{-2}$,
and the four squared column norms at the cut are at most $C_{\omega,K}h^4$. Thus the two conclusions are
\[
\sup_{z\in K}
\left\|
(H_n^{B,\dagger}-z)^{-1}
-
(H_n^{B,\dagger,\oplus}-z)^{-1}
\right\|_{\mathfrak S_2}^2
\le
C_{\omega,K}h^4.
\]
Likewise,
\[
\sup_{z\in K}
\left\|
(\mathcal B_{M,h,\beta}^{B}-z)^{-1}
-
(\mathcal B_{M,h,\beta}^{B,\oplus}-z)^{-1}
\right\|_{\mathfrak S_2}^2
\le
C_{\omega,K}h^4.
\]
The Hilbert--Schmidt tail estimate gives, for the tail blocks, \(\|\mathbf 1_{\{x\ge x_J\}}R(z)\|_{\mathfrak S_2}^2 \le C_{\omega,K}x_J^{-1/2}\). Since \(x_J\asymp h^{-2},\) we get
\(x_J^{-1/2}\asymp h.\) Thus
\[
\|(H_{\mathrm{tail}}^{B,\dagger}-z)^{-1}\|_{\mathfrak S_2}^2
+
\|(\mathcal B_{\mathrm{tail},h,\beta}^{B}-z)^{-1}\|_{\mathfrak S_2}^2
\le
C_{\omega,K}h.
\]
\[
\left\| (H_{\mathrm{tail}}^{B,\dagger}-z)^{-1} - (\mathcal B_{\mathrm{tail},h,\beta}^{B}-z)^{-1} \right\|_{\mathfrak S_2}^2 \le 2\|(H_{\mathrm{tail}}^{B,\dagger}-z)^{-1}\|_{\mathfrak S_2}^2 + 2\|(\mathcal B_{\mathrm{tail},h,\beta}^{B}-z)^{-1}\|_{\mathfrak S_2}^2 \le C_{\omega,K}h.
\eqtag{17.79}
\]
The interior and tail estimates now give the full comparison.

\begin{proposition}\label{A:res:17.4}
For every compact \(K\Subset\mathbb C\setminus\mathbb R,\) almost surely, for all sufficiently large \(n\),
\[
\sup_{z\in K}
\left\|
(H_n^{B,\dagger}-z)^{-1}
-
(\mathcal B_{M,h,\beta}^{B}-z)^{-1}
\right\|_{\mathfrak S_2}^2
\le
C_{\omega,V,\beta,K}h.
\]
\end{proposition}

\begin{proof}
The difference of the direct sum resolvents is itself the direct sum of the two block differences. Its squared Hilbert--Schmidt norm is the sum of their squared norms, bounded by
$C_{\omega,K}h$ by \eqref{A:eq:17.71} and \eqref{A:eq:17.79}. Each of the two cut errors is $O_\omega(h^4)$ by the rank-two calculation. Writing $R_\dagger,R_B$ for the full
resolvents and $R_\dagger^\oplus,R_B^\oplus$ for the cut resolvents, we obtain
\[
 \|R_\dagger-R_B\|_{\mathfrak S_2}^2
 \le3\bigl(\|R_\dagger-R_\dagger^\oplus\|_{\mathfrak S_2}^2
        +\|R_\dagger^\oplus-R_B^\oplus\|_{\mathfrak S_2}^2
        +\|R_B^\oplus-R_B\|_{\mathfrak S_2}^2\bigr)
 \le C_{\omega,K}(h+h^4)\le C'_{\omega,K}h.
\]
All four estimates are uniform on $K$.
\end{proof}

Combining Propositions~\ref{A:res:17.2} and~\ref{A:res:17.4} gives the next comparison.
\begin{theorem}\label{A:res:17.5}
For every compact \(K\Subset\mathbb C\setminus\mathbb R,\) almost surely, for all sufficiently large \(n\),
\[
\sup_{z\in K}
\left\|
(H_n^{B,\circ}-z)^{-1}
-
(\mathcal B_{M,h,\beta}^{B}-z)^{-1}
\right\|_{\mathfrak S_2}^2
\le
C_{\omega,V,\beta,K}h.
\]
\end{theorem}

\begin{proof}
Put $R_\circ(z)=(H_n^{B,\circ}-z)^{-1}$, $R_\dagger(z)=(H_n^{B,\dagger}-z)^{-1}$, and $R_B(z)=(\mathcal B_{M,h,\beta}^B-z)^{-1}$. Propositions~\ref{A:res:17.2} and~\ref{A:res:17.4}
give
\[
 \sup_{z\in K}\|R_\circ(z)-R_B(z)\|_{\mathfrak S_2}^2
 \le2\sup_{z\in K}\|R_\circ(z)-R_\dagger(z)\|_{\mathfrak S_2}^2
    +2\sup_{z\in K}\|R_\dagger(z)-R_B(z)\|_{\mathfrak S_2}^2
 \le C_{\omega,K}h.
\]
\end{proof}

The probabilistic comparison \eqref{A:eq:15.65} and Theorem~\ref{A:res:17.5} have the same reduced coordinate space. The ideal property of the Hilbert--Schmidt norm and the
triangle inequality therefore give the following bound.

\begin{theorem}\label{A:res:17.6}
Let \(V\) be a uniformly convex even degree polynomial with regular right edge, with \(h=h_{n,V}\) as defined in the introduction. On the coupling constructed in
Sections~\ref{A:sec:boundary-reduction}--\ref{A:sec:application-to-the-conditional-jacobi-coefficients}, for every compact \(K\Subset\mathbb C\setminus\mathbb R\), almost surely
for all sufficiently large \(n\),
\[
\sup_{z\in K}
\left\|
(\mathcal A_{n,V,\beta}-z)^{-1}
-
P_{I_n}^*U_{r_n}^*(\mathcal B_{M,h,\beta}^{B}-z)^{-1}U_{r_n}P_{I_n}
\right\|_{\mathfrak S_2}^2
\le
C_{\omega,V,\beta,K}h.
\]
\end{theorem}

\begin{proof}
Let $E_n=P_{I_n}^*U_{r_n}^*:\mathbb C^M\to\mathbb C^n$. Since $E_n^*E_n=I_M$, \eqref{A:eq:isometric-hs}, with $J=E_n$, gives $\|E_nTE_n^*\|_{\mathfrak S_2}=\|T\|_{\mathfrak S_2}$.
Thus extension by zero to the removed coordinates introduces no dimension dependent factor. Write the difference in the theorem as the sum of the difference in \eqref{A:eq:15.65}
and
\[
 E_n\left[(H_n^{B,\circ}-z)^{-1}
          -(\mathcal B_{M,h,\beta}^B-z)^{-1}\right]E_n^*.
\]
Denote the first term by $D_n(z)$; \eqref{A:eq:15.65} gives $\sup_K\|D_n(z)\|_{\mathfrak S_2}^2\le C_{\omega,K}h$, including the removed blocks. By Theorem~\ref{A:res:17.5} and the
isometry identity, the full difference has squared norm at most
\[
 2\sup_{z\in K}\|D_n(z)\|_{\mathfrak S_2}^2
 +2\sup_{z\in K}\|(H_n^{B,\circ}-z)^{-1}
              -(\mathcal B_{M,h,\beta}^B-z)^{-1}\|_{\mathfrak S_2}^2
 \le C_{\omega,K}h.
\]

\end{proof}

\subsubsection{Measurable couplings of the finite models}\label{A:sec:measurable-couplings-of-the-finite-models}

The conditional couplings must fit together into one construction for all matrix sizes. This requires the samplers to be Borel in the Brownian path and exterior coefficients and to
preserve the Gaussian marginal. Once these properties are established, the coupling is fixed independently of the later choices of moment orders and probability thresholds.

Fix $I_\beta=[\beta_-,\beta_+]\Subset(0,\infty)$ and $K\Subset\mathbb C\setminus\mathbb R$. Here the parameter interval concerns measurability of the kernels. The moment bounds and
the distinct probability quantifiers for varying parameters are treated in Section~\ref{A:sec:parameter-dependence-and-further-directions}.

Let \(\mathfrak B=\{b\in C([0,\infty)):b(0)=0\}\), with the compact-open topology and Wiener measure \(\mathbb W\). Put
\[
 \mathcal X_K=C\bigl(K;\mathfrak S_2(L^2(\mathbb R_+))\bigr),
 \quad \mathfrak X_n=\mathbb R^n\times(0,\infty)^{n-1}.
\]
The space $\mathcal X_K$ has the supremum norm, and all spaces carry their Borel sigma fields. Write $\mu_{n,\beta}$ for the finite coefficient law. We first check measurability of
the operator-valued error itself. For the limiting resolvent, the solutions of the weak equations on the fixed finite-dimensional subspaces give Borel approximations in this norm;
these approximations depend measurably on the coefficients of the forms.

\begin{lemma}\label{A:res:B.6}
The map \((\beta,B)\mapsto R_\beta(B,\cdot)\) has a Borel version from \(I_\beta\times\mathfrak B\) to \(\mathcal X_K\). For both finite models, the embedded finite resolvent is
Borel in \((\beta,x)\). Their squared \(\mathfrak S_2\) errors, with the supremum over \(K\), are therefore Borel.
\end{lemma}

\begin{proof}
For a fixed coordinate interval, matrix inversion is continuous off the real spectrum, uniformly on $K$. The embedding is also continuous in $h$. For its $j$th column, put
$I_j(h)=[(j-1)h,jh)$; then
\[
 \left\|h^{-1/2}\mathbf1_{I_j(h)}
       -h_0^{-1/2}\mathbf1_{I_j(h_0)}\right\|_2^2
 =2-\frac{2|I_j(h)\cap I_j(h_0)|}{\sqrt{hh_0}}
 \longrightarrow0.
\]
There are finitely many columns, so the embeddings converge in operator norm. Composition with the continuous matrix resolvent proves continuity into $\mathcal X_K$. An
integer-valued cutoff partitions $I_\beta$ into countably many Borel sets, on each of which this argument applies. This covers the map $\mathcal J_n$ of the introduction, including
the empty interval convention. For $\beta\in I_\beta$, the common terminal cutoff in \eqref{A:eq:B.45} fixes the selected interval.

For the continuum, let \(V_r\) be the space of piecewise linear functions with zero boundary values on \([0,2^r]\) with intervals of length \(2^{-r}\), extended by zero. These
spaces are nested and their union is dense in $\mathcal D_{\mathrm A}$ in its form norm. To see density, first approximate a function in the form domain by a compactly supported
smooth function, using a smooth cutoff followed by mollification. For such a function $f$, its piecewise linear interpolants on intervals of length $2^{-r}$ are supported in a
fixed compact interval and have

$H^1$ error at most $C_f2^{-r}$. The weight $1+x$ is bounded on that support, so the same convergence holds in the form norm. Choose the piecewise linear basis described in
Subsection~\ref{A:sec:finite-difference-approximation-of-the-stochastic-airy-operator}. Its Gram matrix consists of deterministic \(L^2\) inner products. The form entries depend
continuously on \((\beta,B)\), since the noise entries are integrals of \(B\) against fixed compactly supported piecewise polynomials. The finite weak equation again has matrix
\(K-zM\), with these form and Gram matrices. Hence the embedded resolvents of the restricted forms \(G_r(\beta,B,\cdot)\) are Borel
into \(\mathcal X_K\).

Fix \((\beta,B)\in I_\beta\times\Omega_B\), and write \(A=\mathcal H_\beta(B)+\kappa_B>0\). If \(P_r\) is the orthogonal projection onto \(A^{1/2}V_r\), then the inverse of the
restricted form at \(-\kappa_B\) is
\[
 S_r=A^{-1/2}P_rA^{-1/2},\quad S=A^{-1}.
\]
Indeed, for \(F\in L^2\), its finite weak solution \(u_r\in V_r\) is characterized by
\[
 \langle A^{1/2}u_r,A^{1/2}v_r\rangle
 =\langle F,v_r\rangle
 =\langle A^{-1/2}F,A^{1/2}v_r\rangle,
 \qquad v_r\in V_r.
\]
Thus \(A^{1/2}u_r=P_rA^{-1/2}F\), which proves the formula. The form norm is equivalent to $\|A^{1/2}f\|_2$, so density gives $P_r\to I$ strongly. The form bound and min--max give
$\lambda_k(A)\ge c_B k^{2/3}$, so $\|A^{-1/2}\|_{\mathfrak S_4}^4=\sum_k\lambda_k(A)^{-2}<\infty$. If $F$ is a finite-rank approximation to $A^{-1/2}$ in $\mathfrak S_4$, then
\[
 \|(I-P_r)A^{-1/2}\|_{\mathfrak S_4}
 \le\|A^{-1/2}-F\|_{\mathfrak S_4}
     +\|(I-P_r)F\|_{\mathfrak S_4}.
\]
Write $F=\sum_{j=1}^m a_ju_jv_j^*$ with $u_j,v_j$ unit vectors. Strong convergence gives
\[
 \|(I-P_r)F\|_{\mathfrak S_4}
 \le\sum_{j=1}^m|a_j|\|(I-P_r)u_j\|\longrightarrow0.
\] Thus the Schatten ideal
inequality gives
\[
 \|S_r-S\|_{\mathfrak S_2}
 \le\|A^{-1/2}\|_{\mathfrak S_4}
      \|(P_r-I)A^{-1/2}\|_{\mathfrak S_4}\longrightarrow0.
\]
Let $c_z=z+\kappa_B$. The positive operator associated with the restricted form and $A$ have real spectra, so the spectral theorem gives
\[
 \|(I-c_zS_r)^{-1}\|_{\mathrm{op}}
 +\|(I-c_zS)^{-1}\|_{\mathrm{op}}
 \le2\left(1+\frac{|z+\kappa_B|}{|\operatorname{Im}z|}\right).
\]
On the nullspace of $S_r$ the inverse is the identity. Since $G_r(z)=S_r(I-c_zS_r)^{-1}$, the exact identity
\[
 G_r(z)-R_\beta(z)
 =(I-c_zS_r)^{-1}(S_r-S)(I-c_zS)^{-1}
\]
proves uniform convergence on $K$ in $\mathfrak S_2$. Thus $G_r(\beta,B,\cdot)$ converges in $\mathcal X_K$ for each $(\beta,B)\in I_\beta\times\Omega_B$. The shift $\kappa_B$ is
used only to prove convergence on a fixed path: the maps defined by the finite weak equations whose limit is taken are the unshifted Borel maps already constructed. Set
$R_\beta(B,\cdot)=0$ for $B\notin\Omega_B$. Since $\Omega_B$ is Borel, the resulting $\mathcal X_K$-valued map is Borel.
\end{proof}

\begin{definition}\label{A:res:B.7}
Let $Q_n(\beta,B,\mathrm dx)$ be a jointly Borel Markov kernel from $I_\beta\times\mathfrak B$ to $\mathfrak X_n$ satisfying
\[
 \int_{\mathfrak B}Q_n(\beta,B,A)\,\mathbb W(\mathrm{d}B)
       =\mu_{n,\beta}(A)
 \quad(A\in\mathcal B(\mathfrak X_n)).
\]
Its joint law is \(\Pi_{n,\beta}(\mathrm{d}B,\mathrm{d}x)=\mathbb W(\mathrm{d}B)Q_n(\beta,B,\mathrm{d}x)\).
\end{definition}

\begin{lemma}\label{A:res:B.7A}
Let \(\Theta\) be standard Borel and \(p:\Theta\times\mathbb R^d\to[0,\infty)\) be Borel, with \(0<Z_\theta=\int p(\theta,x)\,\mathrm{d}x<\infty\). Then
\(Z_\theta^{-1}p(\theta,x)\,\mathrm{d}x\) is a Borel probability kernel and admits a Borel sampler \(T(\theta,u)\), \(u\in(0,1)^d\), from independent uniform coordinates.
Additional conditioning variables may be included in \(\theta\).
\end{lemma}

\begin{proof}
For $0\le k\le d$, define the unnormalized marginal densities
\[
 m_k(\theta,x_1,\ldots,x_k)
 =\int_{\mathbb R^{d-k}}p(\theta,x_1,\ldots,x_d)\,
                         \mathrm dx_{k+1}\cdots\mathrm dx_d.
\]
Here $m_0=Z_\theta$ and $m_d=p$. Parameter integration of a nonnegative Borel function is Borel; hence every $m_k$ is Borel. In particular, for each Borel $A\subset\mathbb R^d$,
the function $\theta\mapsto Z_\theta^{-1}\int_Ap(\theta,x)\,\mathrm dx$ is Borel, proving the kernel statement.

We construct the sampler successively. If $0<m_{k-1}(\theta,x_{<k})<\infty$, set
\[
 f_k(\theta,x_{<k};t)
 =\frac{m_k(\theta,x_{<k},t)}{m_{k-1}(\theta,x_{<k})}.
\]
Tonelli gives $\int_{\mathbb R}f_k\,\mathrm dt=1$ at such a prefix. At all other prefixes use the standard normal density instead. This defines a jointly Borel probability density
in $t$ for every prefix. Write
\[
 F_k(\theta,x_{<k};t)=\int_{-\infty}^t f_k(\theta,x_{<k};s)\,\mathrm ds.
\]
It is jointly Borel, continuous and nondecreasing in $t$, with limits zero and one at the two infinities. Its generalized inverse
\[
 Q_k(\theta,x_{<k};u)
 =\inf\{t\in\mathbb R:F_k(\theta,x_{<k};t)\ge u\},\qquad 0<u<1,
\]
is finite. It is Borel: for each real $a$, the event $\{Q_k<a\}$ is the union over rational $r<a$ of $\{F_k(r)\ge u\}$. Continuity of the distribution function also gives
$\mathbb P(Q_k(U)\le t)=F_k(t)$ for a uniform $U$.

For independent uniform coordinates $u_1,\ldots,u_d$, define recursively
\[
 x_1=Q_1(\theta;u_1),\qquad
 x_k=Q_k(\theta,x_1,\ldots,x_{k-1};u_k),\qquad
 T(\theta,u)=(x_1,\ldots,x_d).
\]
This map is Borel by induction. The first $k$ coordinates have density $m_k/Z_\theta$: the induction step follows by multiplying the density $m_{k-1}/Z_\theta$ by $f_k$. The
exceptional prefixes where $m_{k-1}=0$ or $m_{k-1}=\infty$ have zero probability under that marginal, since $m_{k-1}$ is integrable with total mass $Z_\theta<\infty$. The arbitrary
normal choice there does not affect the law. At $k=d$ the density is $p/Z_\theta$, as required. All additional conditioning variables can be included in $\theta$ throughout this
construction.
\end{proof}

Let $D\subset\mathbb R^d$ be open. Assume that $\Phi_\theta$ is jointly Borel, continuous in $x\in D$, has compact sublevel sets in $D$, and has a unique minimizer $m_\theta$. For
the potentials on the whole space $\Phi_{\theta,t}(x)=\langle x,K_\theta x\rangle/2+t\Psi_\theta(x)$, assume in addition the centering and uniform derivative hypotheses of
Theorem~\ref{A:res:A.14}.

\begin{lemma}\label{A:res:B.7B}
The minimizer $m_\theta$ and the transport $\mathcal T_\theta(x)$ at time one are Borel in their arguments. The latter satisfies
\[
 (\mathcal T_\theta)_\#N(0,K_\theta^{-1})
 =Z_{\theta,1}^{-1}e^{-\Phi_{\theta,1}(x)}\,\mathrm dx.
\]
\end{lemma}

\begin{proof}
Enumerate \(D\cap\mathbb Q^d\) as \((q_j)\). The infimum \(a_\theta=\inf_j\Phi_\theta(q_j)\) is Borel. For each \(r\), choose the first \(q_j\) with
\(\Phi_\theta(q_j)\le a_\theta+r^{-1}\). Writing the chosen point as $q_{j_r(\theta)}$, its selection sets are
\[
 \{j_r=j\}=\{\Phi_\theta(q_j)\le a_\theta+r^{-1}\}
 \cap\bigcap_{i<j}\{\Phi_\theta(q_i)>a_\theta+r^{-1}\}.
\]
Thus each selected point is Borel. Density of the rational points and continuity give $a_\theta=\Phi_\theta(m_\theta)$. For fixed $\theta$, all selected points lie in the compact
subset
\[
 C_\theta=\{x\in D:\Phi_\theta(x)\le\Phi_\theta(m_\theta)+1\}.
\]
Every subsequence has a convergent further subsequence in $C_\theta$. If its limit is $x$, continuity and the choice of the approximants give
\[
 0\le\Phi_\theta(x)-\Phi_\theta(m_\theta)
   =\lim_{k\to\infty}
       [\Phi_\theta(q_{j_{r_k}(\theta)})-a_\theta]=0.
\]
Uniqueness implies $x=m_\theta$. If the full sequence did not converge to $m_\theta$, some subsequence would stay at distance at least $\varepsilon>0$ from it and would have a
limit with that property, a contradiction. Hence $q_{j_r(\theta)}\to m_\theta$, and the coordinatewise limit is Borel.

For the transport, put \(f_{\theta,t}=\Psi_\theta-\int\Psi_\theta\,\mathrm{d}\mu_{\theta,t}\). The integral and derivatives in the coefficient variables are Borel, by parameter
integration and difference quotients. With \(W\) an auxiliary Brownian motion, the equations
\[
 X_s =x+\sqrt2W_s-\int_0^s\nabla\Phi_{\theta,t}(X_r)\,\mathrm{d}r,\qquad J_s
 =I-\int_0^s\nabla^2\Phi_{\theta,t}(X_r)J_r\,\mathrm{d}r
\]
have Borel solutions in $(\theta,t,x,W)$. For the first equation, start with $X_s^{(0)}=x+\sqrt2W_s$ and set
\[
 X_s^{(m+1)}=x+\sqrt2W_s
          -\int_0^s\nabla\Phi_{\theta,t}(X_r^{(m)})\,\mathrm dr.
\]
If the drift Lipschitz constant is $L$, then on $[0,T]$
\[
 \sup_{s\le T}|X_s^{(m+1)}-X_s^{(m)}|
 \le LT\sup_{s\le T}|X_s^{(m)}-X_s^{(m-1)}|.
\]
Choose $LT\le1/2$. Each iterate is Borel, so its uniform limit is Borel. The Jacobian equation is linear with a bounded Borel matrix coefficient and obeys the same iteration
estimate. Successive time intervals of length $T$ define both solutions on every finite interval. Lemma~\ref{A:res:A.13} gives
\[
 v_{\theta,t}(x)
 =-\lim_{N\to\infty}\int_0^N
       \mathbb E[J_s^*\nabla f_{\theta,t}(X_s)]\,\mathrm{d}s.
\]
Each finite integral is Borel. The bounds \(\|J_s\|\le \mathrm{e}^{-cns}\), \(\mathbb E|X_s|\le \mathrm{e}^{-cns}|x|+C\), and \(|\nabla\Psi_\theta(x)|\le Cn|x|\) give an integrable
upper bound and
\(|v_{\theta,t}(x)|\le C(1+|x|)\). The second variation calculation in Theorem~\ref{A:res:A.14} gives
$\|D_xv_{\theta,t}(x)\|\le C(1+|x|)$, locally uniformly in the parameters considered here. For the final ODE, fix $R$ and multiply $v_{\theta,t}$ by a smooth cutoff which is one on
$|x|\le R$ and zero for $|x|\ge2R$. The resulting field is globally Lipschitz with a common constant on each bounded parameter set. Successive Picard iterates on sufficiently short
time intervals are Borel and converge uniformly. Their limit is therefore Borel. The growth estimate $1+|T_t(x)|\le e^{Ct}(1+|x|)$ shows that for each $x$ these cutoff solutions
agree with the full solution once $R$ is sufficiently large. Taking the limit through integer $R$ proves joint Borel measurability. The Jacobian identity in
Theorem~\ref{A:res:A.14} gives the stated push-forward law.
\end{proof}

Let \(q_{\theta,s}\) be jointly Borel on the unit circle and satisfy the hypotheses of Proposition~\ref{A:res:A.18}, with one annulus and uniform \(C^1([0,s_0];H^\infty)\) bounds
in \(\theta\). Define
\[
 \ell_m(\theta,s)=\frac1{2\pi}\int_{-\pi}^{\pi}
 \mathrm{e}^{-\mathrm{i}m\vartheta}\log q_{\theta,s}(\mathrm{e}^{\mathrm{i}\vartheta})\,\mathrm{d}\vartheta,
\]
\[
 g_{\theta,s}(z)=\exp\left\{\tfrac12\ell_0(\theta,s)
                    +\sum_{m\ge1}\ell_m(\theta,s)z^m\right\}.
 \eqtag{B.40}
\]

\par\Needspace{6\baselineskip}
\begin{lemma}\label{A:res:B.7C}
The factor is Borel in $\theta$ and $C^1$ in $s$. It satisfies \(|g_{\theta,s}|^2=q_{\theta,s}\) on the circle,
\(g_{\theta,s}(1)>0\), and
\[
 |g_m(\theta,s)|+|\partial_sg_m(\theta,s)|\le C\mathrm{e}^{-c|m|}.
 \eqtag{B.40a}
\]
Negative coefficients are zero. The coefficients of \(p_{\theta,s}(z)=[g_{\theta,s}(z)-g_{\theta,s}(1)]/(1-z)\) have the same measurability and decay.
\end{lemma}

\begin{proof}
The real logarithm on the circle and parameter integration make \(\ell_m\) Borel. Proposition~\ref{A:res:A.18} supplies a common smaller annulus without zeros and
\(|\ell_m|+|\partial_s\ell_m|\le C\mathrm{e}^{-c|m|}\). Choose $1<r<e^c$ and decrease $r$ if necessary to stay inside the common analytic disk. The series for the exponent in
\eqref{A:eq:B.40}, and its $s$-derivative, are uniformly bounded on $|z|\le r$ by the geometric sum $C\sum_{m\ge0}(re^{-c})^m$. Exponentiation gives
$\sup_{|z|=r}(|g_{\theta,s}(z)|+|\partial_sg_{\theta,s}(z)|)\le C_r$. Cauchy's formula gives
\[
 g_m(\theta,s)=\frac{1}{2\pi r^m}\int_{-\pi}^{\pi}
           g_{\theta,s}(re^{\mathrm i\vartheta})
                          e^{-\mathrm im\vartheta}\,\mathrm d\vartheta,
 \qquad
 |g_m|+|\partial_sg_m|\le C_r r^{-m}.
\]
The partial sums are Borel in $\theta$, their uniform limits remain Borel, and integration in this formula preserves measurability. Differentiation in $s$ is justified by the
uniform derivative bound. This proves \eqref{A:eq:B.40a}. The exponent in \eqref{A:eq:B.40} is real at $z=1$, so its exponential satisfies \(g_{\theta,s}(1)>0\); the factor
identity is that of Proposition~\ref{A:res:A.18}. Corollary~\ref{A:res:A.19} applies to \(p_{\theta,s}\).

For the Jacobi model the symbol bounds are \eqref{A:eq:symbol-analytic-derivative}--\eqref{A:eq:normalized-symbol}. The coordinates $s_j$ and $s_{j-1}$ satisfy
\eqref{A:eq:profile-size-increment} and \eqref{A:eq:profile-extension}. These are the hypotheses of Proposition~\ref{A:res:A.20}, uniformly on \(I_\beta\); the exterior
coefficients do not enter this deterministic symbol.
\end{proof}

\begin{lemma}\label{A:res:B.7C1}
The Gaussian constructions of Lemmas~\ref{A:res:A.16} and \ref{A:res:A.16A} give Borel conditional kernels for Borel covariance matrices and positive diagonal weights, including
singular covariances. Either Gaussian vector may be sampled first.
\end{lemma}

\begin{proof}
Let \(K=\mathbb E[XY^*]\) be the cross covariance in Lemma~\ref{A:res:A.16}. The conditional sampler given \(Y=y\) is
\[
 X=KD^\dagger y+(C-KD^\dagger K^*)^{1/2}Z,
 \quad Z\sim N(0,I).
\]
Given the linear image \(X=x\), it is
\[
 G=\Sigma L^*C^\dagger x+
       (\Sigma-\Sigma L^*C^\dagger L\Sigma)^{1/2}Z,
 \quad C=L\Sigma L^*.
\]
The matrices $C-KD^\dagger K^*$ and $\Sigma-\Sigma L^*C^\dagger L\Sigma$ are nonnegative by the cited lemmas. The formulas keep the given vector $Y$ or $X$ fixed and sample its
remaining conditional coordinates. Lemmas~\ref{A:res:A.16} and~\ref{A:res:A.16A} prove the required marginal identities.

Square roots are continuous on nonnegative matrices. Their Moore--Penrose inverses are Borel; for such a matrix $D$, the formula
\[
 D^\dagger=\lim_{\varepsilon\downarrow0}(D^2+\varepsilon I)^{-1}D
\]
proves this directly by spectral calculus. Thus both samplers are Borel even when the rank changes. A Gaussian vector lies in the range of its covariance almost surely, so only
those values are needed for the specified conditional law.
\end{proof}

Let \(\mu_\theta,\nu_\theta\) be Borel probability kernels on a standard Borel space \(E\), and let \(\mathcal C\subset\Theta\times E\) be Borel. Suppose their normalized
restrictions to \(\mathcal C_\theta\) agree whenever both masses are positive.

\begin{lemma}\label{A:res:B.7D}
Given \(Y_\theta\sim\nu_\theta\), there is a Borel conditional sampler for \(X_\theta\sim\mu_\theta\) with
\[
 \mathbb P(X_\theta\ne Y_\theta)
 \le\max\{\mu_\theta(\mathcal C_\theta^c),
                  \nu_\theta(\mathcal C_\theta^c)\}.
 \eqtag{B.43}
\]
\end{lemma}

\begin{proof}
Put \(a=\mu_\theta(\mathcal C_\theta)\), \(b=\nu_\theta(\mathcal C_\theta)\), and \(m=\min(a,b)\). For \(ab>0\), let
\(\rho\) be the common normalized restriction; otherwise use a fixed probability measure. All are Borel in \(\theta\).
If \(m=1\), set \(X=Y\). If \(m<1\), put
\[
 \mu^{\rm r}=\frac{\mu_\theta-m\rho}{1-m},\quad
 r(y)=\begin{cases}(m/b)\mathbf1_{\mathcal C_\theta}(y),&b>0,\\
 0,&b=0,\end{cases}
\]
and use the kernel \(r(y)\delta_y+(1-r(y))\mu^{\rm r}\). Since \(r(y)\nu_\theta(\mathrm{d}y)=m\rho(\mathrm{d}y)\), integration gives \(\mu_\theta\) as the first marginal and
diagonal mass at least \(m\). This proves \eqref{A:eq:B.43}, including the cases of zero mass. The displayed probability kernel has a jointly Borel uniform sampler by Lemma~4.22 of
\cite{kallenberg2021foundations}, as used in Lemma~\ref{A:res:B.7A}; include $(\theta,y)$ in the kernel parameter. More generally, let $E_1,E_2,E_3$ be standard Borel spaces, let
$\nu$ be a probability measure on $E_2$, and let $K_j$ be probability kernels from $E_2$ to $E_j$ for $j=1,3$. The measure
$\pi(\mathrm dx,\mathrm dy,\mathrm dz)=\nu(\mathrm dy)K_1(y,\mathrm dx)K_3(y,\mathrm dz)$ satisfies
\[
 \pi(A\times B\times E_3)=\int_B K_1(y,A)\,\nu(\mathrm dy),
 \qquad
 \pi(E_1\times B\times C)=\int_B K_3(y,C)\,\nu(\mathrm dy).
\]
\end{proof}

For later use, uniform convexity bounds the unnormalized Jacobi density, at each fixed \(n\), by the integrable function
\[
 C_n \mathrm{e}^{-cn(\sum a_i^2+\sum b_i^2)}
 \prod_{k=1}^{n-1}\left[
 b_k^{\beta_-(n-k)-1}\mathbf1_{\{b_k\le1\}}
 +b_k^{\beta_+(n-k)-1}\mathbf1_{\{b_k>1\}}\right].
 \eqtag{B.44}
\]
Indeed \(\operatorname{Tr}V(T)\ge c(\sum a_i^2+2\sum b_i^2)-Cn\), and all exponents at zero exceed \(-1\). Dominated convergence makes \(\beta\mapsto\mu_{n,\beta}\) continuous in
total variation. Marginal densities and the ratios defining conditional densities are Borel. The Gaussian and chi densities give the same measurability and continuity for the
Dumitriu--Edelman law.

\begin{lemma}\label{A:res:B.8}
The couplings in the main proof admit kernels as in Definition~\ref{A:res:B.7}. They also admit Borel kernels retaining all auxiliary Gaussian and transport coordinates and one
common two-sided Brownian path. For any countable family of sizes and parameters, these kernels can be sampled over that same path.
\end{lemma}

\begin{proof}
Write \(B=(B^+,B^-)\in\mathfrak B^2\) for the two independent halves and \(\mathbb W_\pm=\mathbb W\otimes\mathbb W\). We identify them with $B(t)=B^+(t)$ for $t\ge0$ and
$B(t)=-B^-(-t)$ for $t<0$. The Borel set
\[
 \mathfrak B_{\rm temp}
 =\left\{B:\sup_{t\in\mathbb Q}
                \frac{|B(t)|}{1+|t|}<\infty\right\}
\]
has \(\mathbb W_\pm\)-measure one. For example, the reflection bound gives \(\mathbb P(\sup_{0\le t\le k+1}|B^\pm(t)|>k) \le4\mathrm{e}^{-k^2/(2(k+1))}\), summable in \(k\ge1\). On
this set, the increment sequence for each fixed lattice spacing has at most linear growth. The exponential bound on the coefficients in \eqref{A:eq:B.40a} makes all convolutions
with the spectral factor absolutely convergent. Their partial sums are Borel in the path and parameters, hence so are their limits. Off this set, sample the matrix independently
with its original law and set the unused auxiliary variables to zero.

For the Dumitriu--Edelman law use the Gaussian rotations and the chi quantiles of Section~\ref{A:sec:coupling-the-dumitriu-edelman-fluctuations-to-brownian-increments}. The
rotations are continuous in
\(\beta\), and the quantiles are Borel by Lemma~\ref{A:res:B.7A}; a chi variable of degree zero is set to zero. The
covariance and marginal calculations in that section therefore define a Borel probability kernel with the joint law just computed.

For the general model at one fixed $\beta$, use the terminal cutoff $s_\beta$ of the introduction. For a common cutoff on $I_\beta$, use
\[
 s_*:=\lceil\beta_-^{-1}\rceil-1+2R_V+10.
 \eqtag{B.45}
\]
Use fixed logarithmic separation lengths. Choose $n_0$ so that, for every $\beta\in I_\beta$ and $n\ge n_0$, all interval inequalities
$J_n^{\mathrm L}<J_n^{\mathrm{reg}}<J_n^{\mathrm{tr}}<M_n$ and regular parameter bounds used above hold. Such a deterministic $n_0$ exists by the fixed positive parameter margins
and the uniform relation $h\asymp n^{-1/3}$. For $n<n_0$, sample the matrix independently with its original law. This also covers an empty selected interval. For $n\ge n_0$, let
$\tau$ consist of the exact complement of the free transport coordinates in \eqref{A:eq:14.2}, as specified by $\mathscr F_n^{\mathrm{out}}$. In particular the original
coefficients $b_r,b_N$ remain exterior data, whereas $B_{J_n^{\mathrm{tr}}}$ is a free transported edge. Its law
\(\nu_{n,\beta}(\mathrm{d}\tau)\) and the conditional law \(\mu_{\beta,\tau}^+(\mathrm{d}y)\) are Borel by
\eqref{A:eq:B.44}. Sample \(\tau\) independently of \(B\). With \(\theta=(\beta,\tau)\), set
\[
 \mu_\theta^0(A)=\mu_{\beta,\tau}^+(m_\theta+A),\quad
 K_\theta=D^2\Phi_\theta(m_\theta),\quad
 \Sigma_\theta=K_\theta^{-1},\quad
 C_\theta=L_\theta\Sigma_\theta L_\theta^*.
 \eqtag{centered-kernel}
\]
Here \(m_\theta\) is the conditional minimizer and \(L_\theta\) is the scalar coefficient map. The barriers on the positive coordinates and strict convexity give compact sublevels
and a unique interior minimizer for each exterior value with positive off-diagonal coefficients. Lemma~\ref{A:res:B.7B} therefore applies on the positive domain. All matrices in
the display are Borel.

Let \(T=\mathcal T_n^{\rm out}(M_0)\), and form \(Y^{\rm fac}=\mathcal G_\beta\xi(B)\) on every exterior configuration. On \(T\), use Lemma~\ref{A:res:B.7C1} to sample
\(X\sim N(0,C_\theta)\) conditionally on this factor vector, then
\(G\sim N(0,\Sigma_\theta)\) conditionally on \(L_\theta G=X\). The extension in Lemma~\ref{A:res:A.11A} uses fixed
cutoff functions and Borel minimizers and derivatives, so its potential is jointly Borel. Lemma~\ref{A:res:B.7B} gives
\(\widetilde U=\mathcal T_\theta(G)\sim\widetilde\mu_\theta\). Apply Lemma~\ref{A:res:B.7D} to the common restricted law
of \(\widetilde\mu_\theta\) and \(\mu_\theta^0\), retaining \(\widetilde U\), to obtain \(U^0\sim\mu_\theta^0\).

On \(T^c\), retain the same \(Y^{\rm fac}\), sample \(G\sim N(0,\Sigma_\theta)\) independently of it, set \(X=L_\theta G\), and independently sample \(U^0\sim\mu_\theta^0\). The
coefficient block is \(m_\theta+U^0\) in both cases. Thus for every bounded Borel \(F\),
\[
 \mathbb E F(\tau,m_\theta+U^0) =\int\nu_{n,\beta}(\mathrm{d}\tau) \int
 F(\tau,y)\,\mu_{\beta,\tau}^+(\mathrm{d}y) =\int F(x)\,\mu_{n,\beta}^V(\mathrm{d}x).
 \eqtag{exact-kernel-marginal}
\]

Denote this kernel by \(\overline Q_n(\beta,B^+,B^-,\mathrm{d}x,\mathrm{d}\mathsf a)\), where \(\mathsf a\in\mathsf A_n\) lists the auxiliary coordinates in a fixed standard Borel
product space, and define
\[
 \begin{split}
 \overline\Pi_{n,\beta}
   &=\mathbb W_\pm(\mathrm{d}B^+,\mathrm{d}B^-)
                   \overline Q_n(\beta,B^+,B^-,\mathrm{d}x,\mathrm{d}\mathsf a),\\
 Q_n(\beta,B^+,A)
   &=\int_{\mathfrak B}\overline Q_n
       (\beta,B^+,B^-,A\times\mathsf A_n)\,\mathbb W(\mathrm{d}B^-).
 \end{split}
\]
The projection of $\overline\Pi_{n,\beta}$ onto $(B^+,x)$ is $\Pi_{n,\beta}$, and \eqref{A:eq:exact-kernel-marginal} proves the required finite coefficient law. Integrating $B^-$
out of the joint construction gives the marginal kernel
$Q_n$.

Indeed, for a countable list $(n_k,\beta_k)$, use the same two-sided path and conditionally independent auxiliary samplers. For cylinder sets $E,A_1,\ldots,A_m$, the resulting law
is
\[
 \int_E\prod_{k=1}^m
 \overline Q_{n_k}(\beta_k,B^+,B^-,A_k)\,
                 \mathbb W_\pm(\mathrm dB^+,\mathrm dB^-).
\]
Applying the Borel samplers to countably many independent uniform variables gives the displayed cylinder laws. Integrating out the other factors leaves $\overline\Pi_{n_k,\beta_k}$
for every $k$. Thus each coefficient and auxiliary joint law used in the estimates is preserved. The matrices need not be independent without conditioning on the path. No sampler
here depends on a moment order or a failure threshold; the summable estimates can all be applied to this one construction.

\end{proof}

\par\Needspace{8\baselineskip}
\subsection{Construction of the coefficient coupling at the hard edge}
\label{sec:the-hard-edge-coefficient-construction}
The hard edge comparison uses logarithmic ratios rather than additive fluctuations of the tridiagonal coefficients. Before constructing their coupling, we return to the inverse
formula to identify the weighted errors. The lower bound will use the Brownian observations available to the approximation, rather than the rank of the approximating operator.

\subsubsection{Method of proof}
\label{B:subsection:intro-green-kernel-comparison}
\label{B:subsection:intro-logarithmic-comparison}

The inverse determines the coefficient errors that must be controlled. We describe that reduction here; the inverse formula and its weighted stability estimate are proved in
Sections~\ref{B:sec:inverse-representations-and-the-stochastic-bessel-operator} and~\ref{B:sec:weighted-stability-and-truncation-of-the-inverse}. Fix $\beta>0$ and $a>0$. Retain
the model normalization $\mathcal A_n=\alpha_n^{-1}n^2B_nB_n^{\mathsf T}$, where $\alpha_n=1$ for beta-Laguerre and $\alpha_n=4/L_V^2$ for the polynomial potential. The respective
partitions are $z_{n,k}=k/n$ and $z_{n,k}=\theta_V(k/n)$, with cells $I_{n,k}$ and lengths $\Delta_{n,k}$.

With $P_{n,k},p_{n,k},q_{n,k}$ as defined in Lemma~\ref{B:result:2-1} and Propositions~\ref{B:result:2-2}--\ref{B:result:2-3}, the embedded inverse has cell values
\begin{equation}
 R_n^\Delta(i,j)=p_{n,i}p_{n,j}\sum_{\ell=i\vee j}^n q_{n,\ell}.
 \label{B:eq:1-53}
\end{equation}
Here $p_{n,k}=\sqrt{\alpha_n}P_{n,k}/\sqrt{n\Delta_{n,k}}$ and $q_{n,k}=(nX_{n,k}^2P_{n,k}^2)^{-1}$; in particular the lengths of both cells appear when the squared kernel is
integrated. The embedding $\mathscr J_{n,V}$ is compared with this normalized cell embedding only after the inverse has been bounded.

\localheading{The two logarithmic errors}
Retain cells $\mathscr I_n^{\rm ret}=\{k_-+1,\ldots,k_+\}$ with $1\le k_-<k_+\le n$. Put $p_{n,k}^B=p^B(z_{n,k})$, $Q_{n,k}^B=\int_{I_{n,k}}q^B(s)\,\dd s$, and
$h_{n,k}=\log(z_{n,k}/z_{n,k-1})$. These integrals have positive lower boundary points. The normalization and errors in \eqref{B:eq:2-35} satisfy
\[
 e^{-c_n}p_{n,k}=p_{n,k}^Be^{\delta_{n,k}},\qquad
 e^{2c_n}q_{n,k}=Q_{n,k}^Be^{-2\delta_{n,k}+\rho_{n,k}}.
\]
The exact cumulative and local formulas are derived in Section~\ref{B:subsection:cumulative-and-local-coefficient-errors}. The error $\delta_{n,k}$ contains the sum of logarithmic
ratios from $k$ to the upper index, whereas $\rho_{n,k}$ contains one coefficient logarithm and the logarithmic error in approximating $Q_{n,k}^B$. The scalar $c_n$ cancels from
every summand of \eqref{B:eq:1-53}. The finite coefficients and $c_n$ use only $\mathcal N_n\vee\mathcal U$; the integrals $Q_{n,k}^B$ are comparison quantities, not additional
observations used to construct the coefficients.

\localheading{The weighted inverse estimate}
Choose $1<\lambda<\min\{a+1,2\}$ and $0<\eta<\frac1{16}\min\{a+1-\lambda,2-\lambda\}$. On selected cells $k\ge2$, $h_{n,k}z_{n,k}^{\lambda}\asymp n^{-\lambda}k^{\lambda-1}$. The
interpolation of Theorem~\ref{B:thm:weighted-inverse-stability} has exponent $\delta_{n,i}+\delta_{n,j}-2\delta_{n,\ell}+\rho_{n,\ell}$. Its derivative is
\eqref{B:eq:interpolated-kernel-derivative}, with $f_\ell=-2\delta_{n,\ell}+\rho_{n,\ell}$. The estimate for the kernel rows controls the two multiplicative factors outside the
sum, and the rank-one estimate controls the terms containing the cell integrals. To apply that theorem, use $p_i=p_{n,i}^B$, $Q_i=Q_{n,i}^B$ and the physical cell lengths
$\Delta_{n,i}$. The bounds for $p^B,q^B$ give the constant $D_\eta$ in that theorem; the bound on the partition intervals gives $h_*$. The strict exponent restrictions above ensure
$a-2\eta>\lambda-1$ and $2-8\eta-\lambda>0$. If $\max_{\mathscr I_n^{\rm ret}}(|\delta_{n,k}|+|\rho_{n,k}|)\le1$, all interpolation constants are uniform in the selected interval.
The theorem therefore gives
\begin{equation}
 \|R_n^{\Delta,0}-R_n^{B,\mathrm{step},0}\|_{\mathfrak S_2}^2
 \le\int_0^1\|\partial_tR_n^t\|_{\mathfrak S_2}^2\,\dd t
 \le C_\omega\sum_{k\in\mathscr I_n^{\rm ret}}
 h_{n,k}z_{n,k}^{\lambda}(|\delta_{n,k}|^2+|\rho_{n,k}|^2).
 \label{B:eq:1-55}
\end{equation}
Both endpoints of the interpolation are products of restricted triangular factors, as in \eqref{B:eq:bessel-finite-restriction} and \eqref{B:eq:2-34}. They are not compressions of
the full inverses. Thus the two coefficient constructions must establish
\begin{equation}
 \sum_{k\in\mathscr I_n^{\rm ret}}h_{n,k}z_{n,k}^{\lambda}
       (|\delta_{n,k}|^2+|\rho_{n,k}|^2)=O_\omega(n^{-1}),\qquad
 \max_{k\in\mathscr I_n^{\rm ret}}(|\delta_{n,k}|+|\rho_{n,k}|)\longrightarrow0.
 \label{B:eq:1-54}
\end{equation}
The maximum estimate makes the exponential factors uniformly bounded; with this bound in place, the weighted square sum controls the inverse error. The proof therefore uses both
estimates.

\localheading{Independent coefficients in the beta-Laguerre model}
For beta-Laguerre, $c_n=0$. At each size the Brownian increments on the intervals of lengths $h_k^+=\log((k+1)/k)$ are independent. The quantile construction in
Section~\ref{B:subsection:laguerre-scalar-distributions} uses them for the ratios $Y_{n,k}/X_{n,k}$; the integral over physical cell $k$ instead uses the preceding time interval of
length $h_{k-1}^+$. Lemma~\ref{B:lem:laguerre-means-variances} and Proposition~\ref{B:prop:laguerre-quantile-bound} separate the centered ratio error of size $O_{L^p}(k^{-3/2})$
from its mean of size $O(k^{-2})$. Their cumulative and weighted sums are proved in Proposition~\ref{B:prop:laguerre-product-errors}, with
Lemma~\ref{B:lem:laguerre-initial-products} treating the initial indices. Proposition~\ref{B:prop:laguerre-retained-errors} includes the local coefficient errors and the errors in
the cell integrals and verifies both parts of \eqref{B:eq:1-54}.

A ratio does not determine the coefficient pair. The representation by independent beta and gamma variables in Proposition~\ref{B:prop:laguerre-scalar-reconstruction} supplies the
independent radial variable. Together with the initial pairs and unpaired last diagonal in \eqref{B:eq:laguerre-final-array}, Proposition~\ref{B:prop:laguerre-array-law} proves the
prescribed joint distribution of the coefficient vector and its $\mathcal N_n\vee\mathcal U$-measurability on the same probability space as the Brownian path.

\localheading{Dependent coefficients and partial sums}
For a polynomial potential, condition on the exterior of $J_n=\{r_n,\ldots,R_n\}$ and denote the boundary data by $q$. The minimizer and conditional Hessian are treated in
Sections~\ref{B:sec:general-potentials-and-convex-reference-measures} and~\ref{B:sec:conditional-estimates-on-the-interior-interval}, following \cite{krishnapur2016universality}
and \cite{rider2019universality}. Let $\nu_q$ be the exact law of $\eta_k=\log(y_k/m_{y,k}^q)-\log(x_k/m_{x,k}^q)$ and $\Gamma_q=N(0,K_q)$, where
$K_q=(n\beta)^{-1}L_qC_qL_q^{\mathsf T}$ and $(L_qv)_k=v_{y,k}/m_{y,k}^q-v_{x,k}/m_{x,k}^q$. Here $C_q$ is the inverse of the conditional Hessian, before any coordinate
restriction.

Section~\ref{B:sec:gaussian-approximation-of-the-coefficient-field} constructs this Gaussian distribution jointly with the Brownian increments. The decomposition
\eqref{B:eq:gaussian-increment-decomposition} includes the discrete difference $\Xi_{n,k}-\Xi_{n,k+1}$. Lemma~\ref{B:lem:gaussian-partial-sum-bounds} estimates the cumulative error
with the correction at the upper index $\Xi_{n,M_n+1}$ included. Proposition~\ref{B:prop:coefficient-scalar-normalization} absorbs that correction into the common scalar
normalization, which is kept fixed when the interval is restricted.

The exact logarithmic law must be approximated in partial sum norm. Theorem~\ref{B:thm:logarithmic-ratio-coupling} constructs a coupling of $\eta\sim\nu_q$ and $u\sim\Gamma_q$
satisfying
\begin{equation}
 \left\|\frac1{n^2}\sum_{j=r_n}^{R_n}j
       \left|\sum_{k=j}^{R_n}(\eta_k-u_k)\right|^2\right\|_{L^p}
 \le Cp\ell_n^2(p+\ell_n)n^{-2},\qquad
 2\le p\le C_0\ell_n,\quad \ell_n=\log(e+n).
 \label{B:eq:intro-ratio-estimate}
\end{equation}
The random variable inside the $L^p$ norm is $|\mathsf T_2(\eta-u)|^2$, where $\mathsf T_2=n^{-1}\operatorname{diag}(\sqrt j)S_n$. It is the nonnegative quadratic form associated
with $\mathsf T_2^{\mathsf T}\mathsf T_2$, evaluated at $\eta-u$. The heat kernel and Lyapunov calculations in Lemmas~\ref{B:lem:partial-sum-heat-kernel} and
\ref{B:lem:partial-sum-lyapunov-bound} give the factor $\ell_n/n^2$. The Poisson and diffusion estimates retain that dependence through
Propositions~\ref{B:prop:conditional-diffusion-matrix}, \ref{B:prop:diffusion-matrix-deviation}, and \ref{B:prop:diffusion-coupling-moments}. These constructions are related to
\cite{courtadeFathiPananjady2019stein,fathiMikulincer2020stability}. The proof of Theorem~\ref{B:thm:logarithmic-ratio-coupling} also includes the conditional mean and the return
from the smooth extension to the original conditional law.

Theorem~\ref{B:thm:summable-ratio-cost} chooses one Borel coupling using the combined cost \eqref{B:eq:ratio-selected-cost}. Corollary~\ref{B:cor:ratio-coupling-errors} integrates
this cost on admissible exteriors and proves the almost sure weighted and maximum estimates. Summability gives these estimates for the dependent sequence of matrix sizes. The
argument bounds the cost of the selected plan by that of the competitor used to estimate the infimum.

\localheading{Recovery of the coefficient law}
The Gaussian vector is already coupled to the Brownian path. Section~\ref{B:subsection:joint-logarithmic-reconstruction} disintegrates the selected ratio coupling over that vector
and then samples the remaining logarithmic coordinates jointly, conditional on the ratios. Its change of variables and conditional integration prove the original coefficient
marginal. Theorem~\ref{B:thm:simultaneous-coefficient-construction} applies those kernels at all sizes with the same Brownian path. Row $n$ still uses only
$\mathcal N_n\vee\mathcal U$.

For $0<\beta<1$, the reference measure in Section~\ref{B:sec:general-potentials-and-convex-reference-measures} changes only the finitely many negative density exponents outside
$J_n$. Their factors cancel from the interior conditional normalizer. The reconstruction uses the original exterior marginal, with the likelihood moments restricted to
$p(1-\beta)<1$. It therefore recovers \eqref{B:eq:1-10} throughout the stated range.

The exact normalization identities in Section~\ref{B:sec:normalization-and-coefficient-estimates} separate the deterministic profile error, the coupling error for the logarithmic
ratios, the Gaussian partial sum, the local coefficient error, and the quadrature error for the integrals of $q^B$. Theorem~\ref{B:thm:normalized-coefficient-energy} and
Proposition~\ref{B:prop:normalized-coefficient-maxima} then verify \eqref{B:eq:1-54} for the same arrays and normalization.

\localheading{Boundary points and the final upper bound}
The retained comparison does not integrate $q^B$ across the singular first cell. Fix $\nu_a=\min\{a,1\}$ and $0<\eta<\nu_a$. The finite and continuum truncation bounds have the
three squared contributions $d_n^{1+\nu_a-\eta}$, $\varepsilon_n^3$ and $\varepsilon_n^2$. For $d_n,\varepsilon_n\le C(\log n)^A/n$, multiplication by $n$ gives
\[
 \begin{aligned}
 n d_n^{1+\nu_a-\eta}
 &\le C n^{-(\nu_a-\eta)}(\log n)^{A(1+\nu_a-\eta)},\\
 n\varepsilon_n^3&\le Cn^{-2}(\log n)^{3A},\qquad
 n\varepsilon_n^2\le Cn^{-1}(\log n)^{2A}.
 \end{aligned}
\]
All three right sides tend to zero because $\nu_a-\eta>0$. The finite bounds are proved in Proposition~\ref{B:prop:laguerre-truncation-error} and
Theorem~\ref{B:thm:finite-inverse-truncation}; Proposition~\ref{B:prop:brownian-factor-truncation} gives the continuum bound. These bounds retain the exterior factors and the
positive correction from the discarded columns of the triangular factor.

For the final comparison, write $R_n^\Delta$ for the finite inverse embedded with normalized cell indicators, $R_n^{\Delta,0}$ for the inverse obtained from its restricted factor,
$R_n^{B,\mathrm{step},0}$ for the inverse obtained from the sampled coefficients on those cells, $R^{B,0}$ for the restricted continuum inverse, and $R=\mathcal H_{\beta,a}^{-1}$
for the full limit. All five operators act on the same $L^2$ space, with zero extensions for the restrictions. Their exact difference is
\[
 \begin{aligned}
 R_n^\Delta-R={}&(R_n^\Delta-R_n^{\Delta,0})
 +(R_n^{\Delta,0}-R_n^{B,\mathrm{step},0})\\
 &+(R_n^{B,\mathrm{step},0}-R^{B,0})+(R^{B,0}-R).
 \end{aligned}
\]
On the common event supplied by the cited estimates, the first and fourth squared norms are $o(n^{-1})$, by the preceding calculation. The coefficient estimate and the weighted
stability theorem bound the second by $C_\omega/n$. Theorem~\ref{B:thm:brownian-inverse-approximation} bounds the third by $C_\omega/n$. The squared norm of a sum of four
Hilbert--Schmidt operators is at most four times the sum of their squared norms. Therefore
\[
 \|R_n^\Delta-R\|_{\mathfrak S_2}^2\le C_\omega/n,
 \qquad
 \|R_n^\Delta\|_{\mathfrak S_2}
 \le\|R\|_{\mathfrak S_2}+C_\omega n^{-1/2}.
\]
The second estimate is obtained before changing the embedding.

Let $R_n$ use the embedding of the theorem. Lemma~\ref{B:result:2-4}, together with isometric invariance of the finite inverse norm, now gives
\[
 E_n:=R_n-R_n^\Delta,\qquad
 \|E_n\|_{\mathfrak S_2}\le(C/n)\|R_n^\Delta\|_{\mathfrak S_2}
                         \le C_\omega/n.
\]
For the Laguerre embedding this difference is zero. In the general case the same bound proves the full squared upper order by the triangle inequality. More precisely, setting
$D_n=R_n^\Delta-R$ and expanding the squared norm gives
\[
 n\left|\|D_n+E_n\|_{\mathfrak S_2}^2-\|D_n\|_{\mathfrak S_2}^2\right|
 \le n\bigl(2\|D_n\|_{\mathfrak S_2}\|E_n\|_{\mathfrak S_2}
                       +\|E_n\|_{\mathfrak S_2}^2\bigr)
 \le C_\omega(n^{-1/2}+n^{-1})\longrightarrow0.
\]
Thus the embedding error cannot affect the squared order $n^{-1}$ or a normalized leading coefficient, should the latter be established.
Proposition~\ref{B:prop:finite-inverse-upper-bound} uses this comparison to combine the coefficient, sampling and boundary estimates.

\localheading{The lower bound}
The lower bound is stated and proved in Theorem~\ref{B:thm:brownian-information-lower-bound}. It applies to each fixed sequence measurable from the row observations and auxiliary
data independent of the whole Brownian path. On a rectangle away from zero, the proof separates the conditional Gaussian bridge term from the nonlinear remainder. After
localization, the shifted small ball estimate and the remainder probability are summable. The arrays constructed above satisfy the required measurability, so the theorem applies to
them. It does not apply to approximations using further observations of the Brownian path. Proposition~\ref{B:result:common-operator-criterion} states the resulting upper and lower
conclusions. Lemma~\ref{B:lem:compact-inverse-resolvent-bounds} and the proofs in Section~\ref{B:subsection:resolvent-rates-and-tightness} then give the resolvent consequences,
including positive real gaps. The spectral consequences are proved in Section~\ref{B:sec:spectral-consequences-and-further-questions}.

\subsubsection{Previous results}
\label{B:subsection:comparison-with-distributional-convergence-rates}

Dumitriu and Edelman \cite{dumitriu2002matrix} gave the beta-Laguerre bidiagonal model, and Edelman and Sutton \cite{edelman2007random} described its limit as a stochastic
differential operator. Ram\'irez and Rider \cite{ramirez2009diffusion,ramirez2011erratum} proved the hard edge limit for every $\beta>0$ and $a>-1$. Their coupling places the
inverse factors on a probability space on which every sequence has a subsequence converging almost surely in Hilbert--Schmidt norm.

Rider and Waters \cite{rider2019universality} proved convergence to the stochastic Bessel operator for the models with a polynomial potential. Their result assumes that
$x\mapsto V(x^2)$ is uniformly convex, $\beta\ge1$, and $a\ge0$. Their coupling for $\beta\ge1$ establishes almost sure Hilbert--Schmidt convergence of the embedded inverse factors
along subsequences on a common probability space. Theorem~\ref{B:result:1-3} supplies the squared rate $n^{-1}$ using Brownian observations at the model's sampling times for fixed
$\beta>0$ and $a>0$. In particular, its conditional construction includes $0<\beta<1$; the boundary point $a=0$ remains outside the parameter range of our upper bound.

Convergence of the factors also gives trace norm convergence of their positive products. Indeed, for $K_n\to K$ in
$\mathfrak S_2$,
\[
 \begin{aligned}
 K_n^*K_n-K^*K&=K_n^*(K_n-K)+(K_n^*-K^*)K,\\
 \|K_n^*K_n-K^*K\|_{\mathfrak S_1}
 &\le(\|K_n\|_{\mathfrak S_2}+\|K\|_{\mathfrak S_2})
                                      \|K_n-K\|_{\mathfrak S_2}\longrightarrow0.
 \end{aligned}
\]
The inequality for the trace norm of a product of Hilbert--Schmidt operators is given in Theorem~2.8 of \cite{simon2005trace}. Our quantitative estimate is for the Hilbert--Schmidt
distance between the positive inverses on one coupling of the whole sequence. The weights in the explicit inverse kernel determine that estimate.

For comparison with another operator coupling, Valk\'o and Vir\'ag \cite{valko2020operator} obtain a bound $C_\omega(\log n)^6/n$ for the squared Hilbert--Schmidt inverse error of
the circular beta ensemble, valid almost surely for every $n\ge n_0(\omega)$, using stopping times of hyperbolic Brownian motion. The sampling times and the independence
assumptions in the present lower bound are specified in \eqref{B:eq:1-30} and \eqref{B:eq:1-32}. The lower bound applies to couplings satisfying those assumptions.

Hard edge universality can also be formulated through correlation kernels. Vanlessen \cite{vanlessen2007strong} proves the unitary Laguerre type results by asymptotics of
orthogonal polynomials. Deift, Gioev, Kriecherbauer, and Vanlessen \cite{deift2007laguerre} treat the orthogonal and symplectic Laguerre type ensembles under their assumptions on a
fixed weight. These results concern rescaled correlation kernels, gap probabilities, and eigenvalue distributions.

To compare distributional rates with our normalization, let $x_{j,n}$ be the eigenvalues of $nB_{n,\beta,a}^{\mathrm L}(B_{n,\beta,a}^{\mathrm L})^{\mathsf T}$. Their Laguerre
density is proportional to
\[
 \prod_{i<j}|x_i-x_j|^\beta
       \prod_{j=1}^n x_j^\nu e^{-\beta x_j/2},
 \qquad \nu=\frac\beta2(a+1)-1,
\]
as follows from \cite{dumitriu2002matrix,ramirez2009diffusion}. The exponent $\nu$ in this density is the parameter used in the following distributional comparisons. Since
$\mathcal A_{n,\beta,a}^{\mathrm L}=n(nBB^{\mathsf T})$,
\begin{equation}
 \lambda_{j,n}^{\mathrm L}=nx_{j,n},\qquad
 s=4nx_{j,n}=4\lambda_{j,n}^{\mathrm L}.
 \label{B:eq:literature-parameter-conversion}
\end{equation}
Bornemann \cite{bornemann2016hardedge} (Section~2, equation~(5)) proves an $O(n^{-2})$ error for the distribution of the smallest eigenvalue at $\beta=2$, for every fixed $\nu>-1$,
after adjusting the hard edge scale. Forrester and Trinh \cite{forresterTrinh2019hardedge} obtain an $O(n^{-2})$ error for the gap probability for fixed $\beta>0$ and
$\nu\in\mathbb Z_{\ge0}$ at the threshold
\begin{equation}
 x=\frac{s}{4(n+\nu/\beta)},\qquad
 \lambda^{\mathrm L}=\frac{ns}{4(n+\nu/\beta)}.
 \label{B:eq:literature-refined-threshold}
\end{equation}
At $\beta=2$, the difference from Bornemann's threshold is
\[
 \frac{s}{4n}\left[\frac1{1+\nu/(2n)}
                      -\left(1-\frac\nu{2n}\right)\right]
 =\frac{s\nu^2}{16n^3(1+\nu/(2n))}.
\]
Thus the relative difference is $O(n^{-2})$ for fixed $\nu$ and $s>0$. At $\beta=2$, the Fredholm determinant of the deterministic Bessel correlation kernel gives the gap
probability. The determinant $\det{}_2(I-z\mathcal H_{\beta,a}^{-1})$ in Corollary~\ref{B:cor:intro-determinant-rates} is a random entire function on the coupling space. The
estimates above and our pathwise operator estimate address these respective quantities.

\subsubsection{Other operator limits}
\label{B:subsection:comparison-with-other-operator-limits}
\label{B:subsection:proof-strategy-and-organization}

The other parts of this work concern soft edge and circular limits. For the operator limit at the soft edge in Section~\ref{sec:soft-edge-ensembles-and-quantitative-airy-limits},
the squared Hilbert--Schmidt resolvent error has order $n^{-1/3}$ for the Gaussian and polynomial Jacobi models. The reciprocal spectral tail supplies the lower bound for every
approximation of rank at most $n$. The proof estimates an additive Brownian potential on an expanding half-line.

For the limit of the canonical operators in Section~\ref{sec:circular-ensembles-canonical-systems-and-the-sine-coefficient}, $n$ times the squared inverse error converges almost
surely to a specified positive functional of the Brownian path. The comparison of the circular inverse with the inverse obtained by sampling the Brownian path has norm
$o(n^{-1/2})$, so the coefficient is determined by the error from sampling the Brownian path. For the coupling of the bidiagonal matrices to the limiting inverse the two
comparisons contribute bounds at order $n^{-1/2}$ in norm. Their inner product is retained in Section~\ref{B:sec:an-exact-leading-error-coefficient}; a leading order coefficient
requires its asymptotics together with those of the two squared terms.

The result for Wigner matrices in Section~\ref{sec:wigner-matrices-and-the-airy-limit-for-a-deterministic-starting-vector} gives convergence in distribution of the nonreal
resolvents, locally uniformly in operator norm, for Lanczos coordinates generated by a deterministic unit vector. Its comparison is with the GOE and GUE. The recurrence
coefficients are functions of a growing matrix of polynomial moments; the entry comparison supplies their local limit and the required spectral localization.

Section~\ref{B:sec:inverse-representations-and-the-stochastic-bessel-operator} proves the inverse formulas and identifies the limiting operator.
Sections~\ref{B:sec:weighted-stability-and-truncation-of-the-inverse} and~\ref{B:sec:approximation-by-brownian-motion-and-the-lower-bound} give the weighted estimates and Brownian
approximation. Section~\ref{B:subsection:brownian-information-lower-bound} states and proves the lower bound. The beta-Laguerre construction and its theorem are proved in
Section~\ref{B:sec:the-beta-laguerre-coupling}. For general potentials,
Sections~\ref{B:sec:general-potentials-and-convex-reference-measures}--\ref{B:sec:transport-of-logarithmic-ratios-and-construction-of-the-coupling} establish the conditional
estimates and the approximation by Gaussian vectors. Section~\ref{B:sec:construction-of-the-coupling-for-all-matrix-sizes} constructs the coefficient arrays jointly from the
probability kernels. Sections~\ref{B:sec:normalization-and-coefficient-estimates} and~\ref{B:sec:removal-of-the-truncation-and-proof-of-the-general-potential-theorem} verify the
two coefficient bounds, remove the boundary point restrictions, and prove Theorem~\ref{B:result:1-3}. Section~\ref{B:sec:spectral-consequences-and-further-questions} contains the
spectral consequences and the questions about the leading error.

\subsection{The beta-Laguerre coupling}
\label{B:sec:the-beta-laguerre-coupling}

The inverse formula \eqref{B:eq:2-7} depends on cumulative logarithms of $Y_k/X_k$. We couple each logarithmic ratio to a Brownian increment and then recover the two coefficients
using an independent gamma variable. This gives the exact finite law and the approximation on one probability space. The fixed difference between the two chi degrees reduces the
cubic term in the scalar density expansion; the centered increment error is consequently $O_{L^p}(k^{-3/2})$ for each fixed $p$. Its deterministic mean is $O(k^{-2})$ and is summed
separately.

Section~\ref{B:subsection:laguerre-array-construction} uses one Brownian path for all sizes and one independent sequence for the initial coefficients. The remaining gamma variables
include the unpaired last diagonal coefficient. Weighted sums of the cumulative errors and the boundary estimates verify Proposition~\ref{B:result:common-operator-criterion} for
these same arrays. Throughout this subsection, $\beta>0$ and $a>0$ are fixed; constants may depend on them.

\subsubsection{The scalar coefficient distributions}
\label{B:subsection:laguerre-scalar-distributions}

The independent chi coefficients are those of the bidiagonal model \cite{dumitriu2002matrix}; the product representation of its inverse is used in \cite{ramirez2009diffusion}. For
the scalar calculation, let $\mathsf Y_k\sim\chi_{\beta k}$ and $\mathsf X_k\sim\chi_{\beta(k+a)}$ be independent, and set
\begin{equation}
 \alpha_k=\frac{\beta k}{2},\qquad
 \gamma_k=\frac{\beta(k+a)}2,\qquad S_k=\alpha_k+\gamma_k,
 \qquad Z_k=\log\frac{\mathsf Y_k}{\mathsf X_k}.
 \label{B:eq:laguerre-shapes}
\end{equation}
Write $T_k=\mathsf X_k^2+\mathsf Y_k^2$, $Q_k=\mathsf Y_k^2/T_k$, and $q(z)=e^{2z}/(1+e^{2z})$. The normalized logarithmic coefficients are
\begin{equation}
 L_k^X=\log\frac{\mathsf X_k}{\sqrt{\beta k}},\qquad
 L_k^Y=\log\frac{\mathsf Y_k}{\sqrt{\beta k}}.
 \label{B:eq:laguerre-scalar-logarithms}
\end{equation}
Gamma distributions are parametrized by shape and scale, in that order.

\begin{proposition}\label{B:prop:laguerre-scalar-reconstruction}
The variables $T_k$ and $Q_k$ are independent, with distributions $\Gamma(S_k,2)$ and $\operatorname{Beta}(\alpha_k,\gamma_k)$. The density of $Z_k$ is
\begin{equation}
 f_k(z)=\frac{2}{\mathrm B(\alpha_k,\gamma_k)}
       \frac{e^{2\alpha_kz}}{(1+e^{2z})^{S_k}},\qquad z\in\mathbb R.
 \label{B:eq:laguerre-ratio-density}
\end{equation}
Conversely, independent $Z_k$ and $T_k$ with these laws give independent chi variables through
\begin{equation}
 \mathsf Y_k=\sqrt{T_kq(Z_k)},\qquad
 \mathsf X_k=\sqrt{T_k(1-q(Z_k))}.
 \label{B:eq:laguerre-reconstruction}
\end{equation}
\end{proposition}
\begin{proof}
Put $\alpha=\alpha_k$, $\gamma=\gamma_k$, and $S=S_k$. The chi density and the substitution $u=x^2$ give the density $u^{\alpha-1}e^{-u/2}/(2^\alpha\Gamma(\alpha))$ for
$U=\mathsf Y_k^2$. Independence gives the joint density
\[
 \frac{u^{\alpha-1}v^{\gamma-1}e^{-(u+v)/2}}
 {2^S\Gamma(\alpha)\Gamma(\gamma)},\qquad u,v>0,
\]
for $U$ and $V=\mathsf X_k^2$. The change of variables is a bijection between $(0,\infty)\times(0,1)$ and $(0,\infty)^2$, with
\[
 (u,v)=(tq,t(1-q)),\qquad
 t=u+v,\quad q=\frac{u}{u+v},\qquad
 \left|\det\begin{pmatrix}q&t\\1-q&-t\end{pmatrix}\right|=t.
\]
Its density therefore factors as
\begin{equation}
 \frac{t^{S_k-1}e^{-t/2}}{2^{S_k}\Gamma(S_k)}
 \frac{q^{\alpha_k-1}(1-q)^{\gamma_k-1}}
      {\mathrm B(\alpha_k,\gamma_k)},\qquad t>0,\quad0<q<1.
 \label{B:eq:laguerre-density-factorization}
\end{equation}
Integration in $t$ and in $q$ gives one for each factor, proving both marginal laws and their independence. The inverse of $z=\tfrac12\log(q/(1-q))$ satisfies
\[
 q(z)=\frac{e^{2z}}{1+e^{2z}},\qquad
 q'(z)=2q(z)(1-q(z)),\qquad
 q(z)^{\alpha-1}(1-q(z))^{\gamma-1}q'(z)
 =\frac{2e^{2\alpha z}}{(1+e^{2z})^S}.
\]
The one-dimensional change of variables proves \eqref{B:eq:laguerre-ratio-density}.

Conversely, for independent $T_k,Z_k$ with the stated laws, the preceding bijection makes $(T_k,q(Z_k))$ have density \eqref{B:eq:laguerre-density-factorization}. Dividing that
density by the Jacobian $t$ gives the displayed product density of $U,V$. Taking square roots gives, on $x,y>0$,
\[
 \frac{2^{1-\gamma}}{\Gamma(\gamma)}x^{2\gamma-1}e^{-x^2/2}
 \frac{2^{1-\alpha}}{\Gamma(\alpha)}y^{2\alpha-1}e^{-y^2/2}.
\]
These are the two normalized chi densities in \eqref{B:eq:laguerre-reconstruction}. All shapes are positive for every $k\ge1$, so these changes of variables include the initial
indices.
\end{proof}

Put $\varrho_k=\tfrac12\log(T_k/(2S_k))$ and $c_k=\tfrac12\log(1+a/(2k))$. The identities $q(z)=e^z/(2\cosh z)$ and $1-q(z)=e^{-z}/(2\cosh z)$ give
\begin{equation}
 \begin{split}
 L_k^Y&=\varrho_k+c_k+\frac{Z_k}{2}-\frac12\log\cosh Z_k,\\
 L_k^X&=\varrho_k+c_k-\frac{Z_k}{2}-\frac12\log\cosh Z_k.
 \end{split}
 \label{B:eq:laguerre-radial-logarithms}
\end{equation}
Indeed,
\[
 \frac{\mathsf Y_k^2}{\beta k}
 =\frac{T_k}{2S_k}\frac{S_k}{\beta k}
        \frac{e^{Z_k}}{\cosh Z_k},\qquad
 \frac{\mathsf X_k^2}{\beta k}
 =\frac{T_k}{2S_k}\frac{S_k}{\beta k}
        \frac{e^{-Z_k}}{\cosh Z_k},\qquad
 \frac{S_k}{\beta k}=1+\frac{a}{2k}.
\]
Taking half the logarithm gives \eqref{B:eq:laguerre-radial-logarithms}. The variable $\varrho_k$ is a function of $T_k$ and is independent of $Z_k$. Subtraction of the two
identities gives $L_k^Y-L_k^X=Z_k$. The nonlinear term has a global integral remainder. Since $(\log\cosh)'(0)=\log\cosh(0)=0$ and $(\log\cosh)''(z)=1-\tanh^2z$,
\[
 \log\cosh z=\frac{z^2}{2}
       -z^2\int_0^1(1-s)\tanh^2(sz)\,\dd s,
 \qquad 0\le\frac{z^2}{2}-\log\cosh z\le\frac{z^4}{12}.
\]
The last bound uses $|\tanh(sz)|\le s|z|$ and $\int_0^1(1-s)s^2\,\dd s=1/12$. Thus the common quadratic term in the two logarithms is $-Z_k^2/4$, with a remainder between zero and
$Z_k^4/24$. It cancels exactly in their difference, before any moment estimate or bound for the absolute value is taken.

\subsubsection{Means and variances}
\label{B:subsection:laguerre-means-and-variances}

The Brownian time increments prescribe both a drift and a variance. The exact logarithmic chi moments show how close these are to the mean and variance of the ratio; the two
discrepancies are kept separate in the cumulative error.

Write $\psi=(\log\Gamma)'$ and $\psi_j=\psi^{(j)}$. The Gaussian increments used below have mean $d_k$ and standard deviation $b_k$, where
\begin{equation}
 h_k^+=\log\frac{k+1}{k},\qquad
 d_k=-\frac a2h_k^+,\qquad b_k=\sqrt{h_k^+/\beta}.
 \label{B:eq:laguerre-gaussian-parameters}
\end{equation}
Also let $\mu_k=\mathbb EZ_k$ and $\sigma_k^2=\operatorname{Var}Z_k$.

\begin{lemma}\label{B:lem:laguerre-means-variances}
For all $k\ge1$,
\begin{equation}
 \mu_k=\tfrac12[\psi(\alpha_k)-\psi(\gamma_k)],\qquad
 \sigma_k^2=\tfrac14[\psi_1(\alpha_k)+\psi_1(\gamma_k)].
 \label{B:eq:laguerre-ratio-moments}
\end{equation}
There are positive constants $c,C$ such that
\begin{equation}
 \frac ck\le\sigma_k^2\le\frac Ck,\qquad
 |\mu_k-d_k|\le Ck^{-2},\qquad
 |\sigma_k-b_k|\le Ck^{-3/2},\qquad
 \mathbb E[(L_k^X)^2+(L_k^Y)^2]=\frac1{\beta k}+O(k^{-2}).
 \label{B:eq:laguerre-moment-estimates}
\end{equation}
\end{lemma}
\begin{proof}
For $\nu>0$, integration of the chi density with $u=x^2/2$ gives
\[
 \frac{2^{1-\nu/2}}{\Gamma(\nu/2)}
 \int_0^\infty x^{\nu+t-1}e^{-x^2/2}\,\dd x
 =\frac{2^{t/2}}{\Gamma(\nu/2)}
 \int_0^\infty u^{(\nu+t)/2-1}e^{-u}\,\dd u.
\]
The integral is finite exactly when $t>-\nu$. Apply this identity to the two independent chi variables to obtain
\begin{equation}
 \mathbb E\chi_\nu^t=2^{t/2}
       \frac{\Gamma((\nu+t)/2)}{\Gamma(\nu/2)}\quad(t>-\nu),
 \qquad
 \mathbb Ee^{tZ_k}=
 \frac{\Gamma(\alpha_k+t/2)\Gamma(\gamma_k-t/2)}
      {\Gamma(\alpha_k)\Gamma(\gamma_k)}
 \quad(-2\alpha_k<t<2\gamma_k).
 \label{B:eq:laguerre-Mellin}
\end{equation}
Differentiation is justified on every compact subinterval of these open domains. More precisely, for every integer $r\ge0$ and $\epsilon>0$,
\[
 |\log x|^r\le r!\epsilon^{-r}(x^\epsilon+x^{-\epsilon})
 \quad(x>0),
\]
by the exponential series on each of $x\ge1$ and $x\le1$. Choose $\epsilon$ smaller than the distance of the exponents to the integrability boundary points. This inequality bounds
the differentiated integrands by two integrable powers at zero and by those powers times the Gaussian density at infinity. Dominated convergence therefore permits every fixed
number of derivatives.

The first and second derivatives at zero of the logarithm of \eqref{B:eq:laguerre-Mellin} give \eqref{B:eq:laguerre-ratio-moments}. For the individual logarithms they give
\begin{equation}
 \begin{aligned}
 \mathbb EL_k^Y&=\tfrac12[\psi(\alpha_k)-\log\alpha_k],&
 \operatorname{Var}L_k^Y&=\tfrac14\psi_1(\alpha_k),\\
 \mathbb EL_k^X&=\tfrac12[\psi(\gamma_k)-\log\alpha_k],&
 \operatorname{Var}L_k^X&=\tfrac14\psi_1(\gamma_k).
 \end{aligned}
 \label{B:eq:laguerre-local-moments}
\end{equation}
For $r\ge2$, the same differentiation gives the cumulant $2^{-r}[\psi_{r-1}(\alpha_k)+(-1)^r\psi_{r-1}(\gamma_k)]$.

We use the series and integral formulas \cite{nistGammaFunctions} (equations~(5.15.1) and~(5.9.15))
\begin{align}
 (-1)^{j+1}\psi_j(x)&=j!\sum_{l=0}^\infty(x+l)^{-j-1},\qquad j\ge1,
       \label{B:eq:laguerre-polygamma-bounds}\\
 \psi(x)&=\log x-\frac1{2x}
       -2\int_0^\infty\frac{t}{(t^2+x^2)(e^{2\pi t}-1)}\,\dd t.
       \label{B:eq:laguerre-Binet}
\end{align}
For $x$ in a compact subset of $(0,\infty)$, the differentiated terms of the $j=1$ series are bounded by a convergent series of powers of $x_0+l$. This gives the formulas for
$j\ge2$ by termwise differentiation. Since $t\mapsto(x+t)^{-2}$ is decreasing,
\[
 \frac1x=\int_0^\infty(x+t)^{-2}\,\dd t
 \le\sum_{l=0}^\infty(x+l)^{-2}
 \le\frac1{x^2}+\int_0^\infty(x+t)^{-2}\,\dd t
 =\frac1{x^2}+\frac1x.
\]
Substitution of $\alpha_k=\beta k/2$ and $\gamma_k=\beta(k+a)/2$ proves $c/k\le\sigma_k^2\le C/k$ for all $k\ge1$.

For $r>0$ put $I_r=\int_0^\infty t^r/(e^{2\pi t}-1)\,\dd t$. The positive geometric series and Tonelli give
\[
 I_r=\sum_{m=1}^\infty\int_0^\infty t^r e^{-2\pi mt}\,\dd t
 =\Gamma(r+1)\sum_{m=1}^\infty(2\pi m)^{-r-1},
 \qquad I_1=\frac1{24},\quad I_3=\frac1{240}.
\]
The identity $(t^2+x^2)^{-1}-x^{-2}=-t^2/[x^2(t^2+x^2)]$ gives
\[
 \left|\psi(x)-\log x+\frac1{2x}+\frac{2I_1}{x^2}\right|
 \le\frac{2I_3}{x^4}.
\]
In particular,
\begin{equation}
 \psi(x)=\log x-\frac1{2x}-\frac1{12x^2}+O(x^{-4}).
 \label{B:eq:6-52}
\end{equation}
For $x$ in a compact subset of $(0,\infty)$, the derivative of the integrand in \eqref{B:eq:laguerre-Binet} has an integrable majorant. Differentiation under the integral gives
\[
 \psi_1(x)=\frac1x+\frac1{2x^2}
 +4x\int_0^\infty
       \frac{t\,\dd t}{(t^2+x^2)^2(e^{2\pi t}-1)}.
\]
For $y\ge0$, $0\le1-(1+y)^{-2}\le2y$. Hence
\[
 \left|\psi_1(x)-\frac1x-\frac1{2x^2}-\frac{4I_1}{x^3}\right|
 \le\frac{8I_3}{x^5},
\]
which proves
\begin{equation}
 \psi_1(x)=x^{-1}+\tfrac12x^{-2}+\tfrac16x^{-3}+O(x^{-5}).
 \label{B:eq:6-53}
\end{equation}
These remainder bounds hold for every $x>0$.

Write $c_a=\gamma_k-\alpha_k=\beta a/2$. On $0\le s\le c_a$,
\[
 \psi_1(\alpha_k+s)
 =\alpha_k^{-1}+(\tfrac12-s)\alpha_k^{-2}
       +O(\alpha_k^{-3}).
\]
Integration in $\mu_k=-\tfrac12\int_0^{c_a}\psi_1(\alpha_k+s)\,\dd s$ and substitution into the variance formula give
\[
 \mu_k=-\frac a{2k}
       +\left(\frac{a^2}{4}-\frac{a}{2\beta}\right)k^{-2}
       +O(k^{-3}),\qquad
 \sigma_k^2=\frac1{\beta k}
       +\left(\frac1{\beta^2}-\frac{a}{2\beta}\right)k^{-2}
       +O(k^{-3}).
\]
The identity $(1+t)^{-1}=1-t+t^2/(1+t)$, integrated over $[0,1/k]$, gives
\[
 h_k^+=k^{-1}-\tfrac12 k^{-2}+O(k^{-3}),\qquad
 |h_k^+-k^{-1}|\le\frac1{2k^2},\qquad
 h_k^+\ge\frac1{k+1}.
\]
For the deterministic centering, Taylor's integral formula is
\[
 h_k^+=\frac1k-\frac1{k^2}\int_0^1\frac{1-t}{(1+t/k)^2}\,\dd t,
 \qquad
 \mu_k-d_k=-\frac12\int_0^{c_a}
       [\psi_1(\alpha_k+s)-\alpha_k^{-1}]\,\dd s
       -\frac{a}{2k^2}\int_0^1\frac{1-t}{(1+t/k)^2}\,\dd t.
\]
The remainders above and $0\le s\le c_a$ give
\[
 \mu_k-d_k=
 \left(\frac{a(a-1)}4-\frac{a}{2\beta}\right)k^{-2}+O(k^{-3}),
 \qquad
 \sigma_k^2-b_k^2=
 \left(\frac1{\beta^2}+\frac{1-a}{2\beta}\right)k^{-2}+O(k^{-3}).
\]
The choice of Brownian increment therefore retains both the mean and variance corrections. The displayed integrals bound the corresponding remainders with constants depending only
on $\beta,a$. Since $\sigma_k+b_k\ge c k^{-1/2}$,
\[
 |\sigma_k-b_k|
 =\frac{|\sigma_k^2-b_k^2|}{\sigma_k+b_k}
 \le Ck^{-3/2}.
\]
Finally, \eqref{B:eq:laguerre-local-moments} gives
\[
 \mathbb EL_k^Y=-\frac1{2\beta k}+O(k^{-2}),\qquad
 \mathbb EL_k^X=\frac{a-\beta^{-1}}{2k}+O(k^{-2}),
\]
and
\[
 \mathbb E[(L_k^X)^2+(L_k^Y)^2]
 =\tfrac14[\psi_1(\alpha_k)+\psi_1(\gamma_k)]
       +(\mathbb EL_k^X)^2+(\mathbb EL_k^Y)^2
 =\frac1{\beta k}+O(k^{-2}).
\]
Each estimate holds for all large $k$ with constants depending on $\beta,a$. The finitely many remaining indices have finite logarithmic moments by \eqref{B:eq:laguerre-Mellin};
increasing the constants proves the estimates for every
$k\ge1$.
\end{proof}

\subsubsection{Logarithmic moments and sums of independent variables}
\label{B:subsection:logarithmic-moments-and-independent-sums}

Moment estimates are needed beyond the Gaussian approximation region, and they must also control sums of independent ratios. The logarithmic moment generating functions provide
both the tail bounds and the dependence on the moment order.

\begin{lemma}\label{B:lem:laguerre-logarithmic-moments}
For $H=Z_k-\mu_k$ or $H=L_k^\sigma-\mathbb EL_k^\sigma$, $\sigma\in\{X,Y\}$,
\begin{equation}
 \log\mathbb Ee^{tH}\le Ct^2/k\quad(|t|\le ck),\qquad
 \mathbb P(|H|>u)\le2e^{-ck\min(u^2,u)}.
 \label{B:eq:laguerre-log-tails}
\end{equation}
For $p\ge1$ and $\xi_k=(Z_k-\mu_k)/\sigma_k$,
\begin{equation}
 \|Z_k-\mu_k\|_p+\|L_k^X\|_p+\|L_k^Y\|_p
 \le C(\sqrt{p/k}+p/k),\qquad
 \|\xi_k\|_p\le C(\sqrt p+p/\sqrt k).
 \label{B:eq:laguerre-log-moments}
\end{equation}
\end{lemma}
\begin{proof}
For $\nu\in\{\alpha_k,\gamma_k\}$, the centered logarithm of the corresponding chi variable has logarithmic moment generating function
\[
 \ell_\nu(t)=\log\Gamma(\nu+t/2)-\log\Gamma(\nu)
                      -\frac t2\psi(\nu).
\]
The corresponding function for $Z_k-\mu_k$ is
\[
 \ell_Z(t)=\log\Gamma(\alpha_k+t/2)+\log\Gamma(\gamma_k-t/2)
             -\log\Gamma(\alpha_k)-\log\Gamma(\gamma_k)-t\mu_k.
\]
Thus
\[
 \ell_\nu''(t)=\tfrac14\psi_1(\nu+t/2),\qquad
 \ell_Z''(t)=\tfrac14[\psi_1(\alpha_k+t/2)+\psi_1(\gamma_k-t/2)].
\]
For $|t|\le c_\beta k$, with $c_\beta=\beta/2$, every argument on the right is at least $\beta k/4$. The bound $\psi_1(x)\le x^{-1}+x^{-2}$ gives $0\le\ell''(t)\le C/k$ for all
three functions. Centering gives $\ell(0)=\ell'(0)=0$, and Taylor's integral formula, valid for both signs of $t$, yields
\[
 \ell(t)=t^2\int_0^1(1-s)\ell''(st)\,\dd s\le Ct^2/k.
\]

For $u>0$, set $t=k\min\{u/(4C),c_\beta/2\}$. Exponential Markov gives $\mathbb P(H>u)\le\exp(-tu+Ct^2/k)$. When
$u\le2Cc_\beta$,
\[
 tu-Ct^2/k=\frac{3ku^2}{16C}.
\]
When $u>2Cc_\beta$,
\[
 tu-Ct^2/k
 =\frac{c_\beta ku}{2}-\frac{Cc_\beta^2k}{4}
 \ge\frac{c_\beta ku}{4}.
\]
The bound for $\ell(-t)$ gives the same two estimates for $\mathbb P(-H>u)$. Adding the probabilities and reducing the constant proves the second inequality in
\eqref{B:eq:laguerre-log-tails}.

Tonelli applied to $|H|^p=\int_0^{|H|}p u^{p-1}\,\dd u$ gives, for $p\ge1$,
\begin{align*}
 \mathbb E|H|^p
 &\le2p\int_0^\infty u^{p-1}(e^{-cku^2}+e^{-cku})\,\dd u\\
 &=p(ck)^{-p/2}\Gamma(p/2)+2p(ck)^{-p}\Gamma(p).
\end{align*}
Here the substitutions are $v=cku^2$ and $v=cku$, respectively. For $v\ge1/2$,
\[
 \Gamma(v+1)=\int_0^\infty x^ve^{-x/2}e^{-x/2}\,\dd x
 \le2(2v/e)^v,
\]
since $x^ve^{-x/2}$ has its maximum at $x=2v$. Taking $p$th roots in the preceding moment inequality gives $\|H\|_p\le C(\sqrt{p/k}+p/k)$. The mean estimates in
Lemma~\ref{B:lem:laguerre-means-variances} then give
\[
 \|L_k^\sigma\|_p
 \le\|L_k^\sigma-\mathbb EL_k^\sigma\|_p+C/k
 \le C(\sqrt{p/k}+p/k),\qquad \sigma\in\{X,Y\}.
\]
Divide the centered estimate for $Z_k$ by $\sigma_k\ge c k^{-1/2}$ to obtain $\|\xi_k\|_p\le C(\sqrt p+p/\sqrt k)$. All constants in these inequalities are independent of $p,k$.
For the reconstruction remainder, use the integral formula following \eqref{B:eq:laguerre-radial-logarithms} and the mean bound $|\mu_k|\le C/k$. For each fixed $p\ge1$ this gives
\[
 \left\|\frac{Z_k^2}{4}-\frac12\log\cosh Z_k\right\|_p
 \le\frac1{24}\|Z_k\|_{4p}^4\le C_pk^{-2},
 \qquad \|Z_k^2\|_p\le C_pk^{-1}.
\]
The quadratic term in each logarithm has order $k^{-1}$, while the integral remainder has order $O_{L^p}(k^{-2})$. Centering that remainder changes its $L^p$ bound by at most a
factor two, since $|\mathbb ER|\le\|R\|_p$. The exact deterministic means remain those in \eqref{B:eq:laguerre-local-moments}.
\end{proof}

\begin{lemma}\label{B:lem:independent-sum-moments}
If $U_i$ are finitely many independent centered real random variables, then, for $p\ge2$,
\begin{equation}
 \left\|\sum_iU_i\right\|_p
 \le C\sqrt p\left(\sum_i\|U_i\|_p^2\right)^{1/2}.
 \label{B:eq:laguerre-independent-sum}
\end{equation}
\end{lemma}
\begin{proof}
It suffices to consider the case in which all the displayed $L^p$ norms are finite. First let $a_1,\ldots,a_m$ be deterministic real numbers and let $\epsilon_i$ be independent
uniform signs. Write $A^2=\sum_i a_i^2$. Since $(\log\cosh t)''=\operatorname{sech}^2t\le1$ and both $\log\cosh0$ and its first derivative vanish, $\cosh t\le e^{t^2/2}$.
Independence therefore gives
\[
 \mathbb E\exp\!\left(t\sum_i a_i\epsilon_i\right)
 =\prod_i\cosh(ta_i)\le e^{t^2A^2/2}.
\]
If $A>0$, applying Markov's inequality with $t=u/A^2$ to the sum and its negative yields
\[
 \mathbb P\!\left(\left|\sum_i a_i\epsilon_i\right|>u\right)
 \le2e^{-u^2/(2A^2)},\qquad u>0.
\]
For $A=0$ the sum is zero. Integration of this tail gives, for $p\ge2$,
\[
 \mathbb E\left|\sum_i a_i\epsilon_i\right|^p
 \le 2p\int_0^\infty u^{p-1}e^{-u^2/(2A^2)}\,\mathrm du
 =p(2A^2)^{p/2}\Gamma(p/2)
 \le(C\sqrt p\,A)^p.
\]
To verify an absolute constant for every real $p\ge2$, let $q=\lceil p/2\rceil$. Then $p\le2q$ and $q\le p$, and the same integral at exponent $2q$ gives
\[
 \mathbb E\left|\sum_i a_i\epsilon_i\right|^{2q}
 \le 2(2A^2)^q q!\le2(2qA^2)^q.
\]
Monotonicity of the probability norms and $2^{1/(2q)}\sqrt{2q}\le2\sqrt q$ therefore imply
\[
 \left\|\sum_i a_i\epsilon_i\right\|_p
 \le\left\|\sum_i a_i\epsilon_i\right\|_{2q}
 \le2\sqrt q\,A\le2\sqrt p\,A.
\]

Now take independent copies $U_i'$ of the entire family, independent also of the signs, and put $D_i=U_i-U_i'$. Conditional Jensen, applied after fixing $(U_i)$, gives
\[
 \mathbb E\left|\sum_i U_i\right|^p
 =\mathbb E\left|\mathbb E\!\left[\sum_iD_i\mid(U_i)\right]\right|^p
 \le\mathbb E\left|\sum_iD_i\right|^p.
\]
Each $D_i$ has a symmetric distribution, and the $D_i$ are independent. Thus $(D_i)$ and $(\epsilon_iD_i)$ have the same joint law. Conditional on $(D_i)$, the deterministic sign
estimate applies, so
\[
 \left\|\sum_iU_i\right\|_p
 \le C\sqrt p\,\left\|\left(\sum_iD_i^2\right)^{1/2}\right\|_p.
\]
Because $p/2\ge1$, the triangle inequality in $L^{p/2}$ implies
\[
 \left\|\left(\sum_iD_i^2\right)^{1/2}\right\|_p^2
 =\left\|\sum_iD_i^2\right\|_{p/2}
 \le\sum_i\|D_i\|_p^2
 \le4\sum_i\|U_i\|_p^2.
\]
Absorbing the factor two after taking square roots proves \eqref{B:eq:laguerre-independent-sum}. Only independence within this finite family is used.
\end{proof}

\subsubsection{A common sequence of initial coefficients}
\label{B:subsection:laguerre-initial-sequence}

The first coefficients cannot be treated by a large-degree approximation uniformly in size. A single independent sequence is used there instead. Bounds holding for all its index
intervals make this choice compatible with the later pathwise kernel estimates.

Let $(\widehat{\mathsf X}_k,\widehat{\mathsf Y}_k)_{k\ge1}$ be mutually independent chi variables with degrees $\beta(k+a)$ and $\beta k$. We will use this same sequence for the
first $m_n$ coefficient pairs in every row of the construction in Section~\ref{B:subsection:laguerre-array-construction}. The estimates below hold simultaneously for every interval
of indices $i\le k<j$. They also hold after changing finitely many degrees to fixed positive values, with constants allowed to depend on those degrees.

\begin{lemma}\label{B:lem:laguerre-initial-products}
Almost surely, for every $\eta>0$ and all $1\le i\le j$,
\begin{equation}
 C_{\omega,\eta}^{-1}j^{-\eta}(i/j)^{a/2}
 \le\prod_{k=i}^{j-1}\frac{\widehat{\mathsf Y}_k}{\widehat{\mathsf X}_k}
 \le C_{\omega,\eta}j^\eta(i/j)^{a/2}.
 \label{B:eq:laguerre-initial-product-bounds}
\end{equation}
There are positive finite random constants $c_\omega,C_\omega$ such that $c_\omega\sqrt k\le\widehat{\mathsf X}_k,\widehat{\mathsf Y}_k \le C_\omega\sqrt k$ for every $k\ge1$.
\end{lemma}
\begin{proof}
Put $\widehat Z_k=\log(\widehat{\mathsf Y}_k/\widehat{\mathsf X}_k)$, $\widehat\mu_k=\mathbb E\widehat Z_k$, and $S_j=\sum_{k=1}^j(\widehat Z_k-\widehat\mu_k)$, with $S_0=0$. Let
$J$ be the finite set of indices at which the degrees have changed. Write $\widehat\nu_k^X,\widehat\nu_k^Y>0$ for the degrees at such an index. The Mellin formula
\eqref{B:eq:laguerre-Mellin} gives
\[
 \ell_k(t):=\log\mathbb E e^{t(\widehat Z_k-\widehat\mu_k)}
 =\log\frac{\Gamma((\widehat\nu_k^Y+t)/2)
                  \Gamma((\widehat\nu_k^X-t)/2)}
                 {\Gamma(\widehat\nu_k^Y/2)\Gamma(\widehat\nu_k^X/2)}
     -t\widehat\mu_k,
 \qquad -\widehat\nu_k^Y<t<\widehat\nu_k^X.
\]
Choose $t_0>0$ smaller than half of every degree with index in $J$, and smaller than the constant $c$ in \eqref{B:eq:laguerre-log-tails}. When $J$ is empty only the latter
restriction is needed. The function $\ell_k''$ is bounded on $[-t_0,t_0]$ for every $k\in J$. Since $\ell_k(0)=\ell_k'(0)=0$, Taylor's integral identity gives
\[
 \ell_k(t)=t^2\int_0^1(1-s)\ell_k''(st)\,\dd s
 \le \frac{t^2}{2}\sup_{|u|\le t_0}|\ell_k''(u)|
 \qquad(k\in J,\ |t|\le t_0).
\]
Choose $C$ at least as large as the constant at the unchanged indices and as $\tfrac{k}{2}\sup_{|u|\le t_0}|\ell_k''(u)|$ for every $k\in J$. Then $\ell_k(t)\le Ct^2/k$ for all
$k\ge1$ and $|t|\le t_0$. For the means, the identity $\sum_{k=i}^{j-1}d_k=-\tfrac a2\log(j/i)$ gives
\[
 \left|\sum_{k=i}^{j-1}\widehat\mu_k+\frac a2\log(j/i)\right|
 \le \sum_{k=i}^{j-1}|\mu_k-d_k|
      +\sum_{k\in J}|\widehat\mu_k-\mu_k|
 \le C\sum_{k=1}^{\infty}k^{-2}
      +\sum_{k\in J}|\widehat\mu_k-\mu_k|<\infty.
\]
We have proved, with one constant for all indices and product intervals,
\begin{equation}
 \log\mathbb Ee^{t(\widehat Z_k-\widehat\mu_k)}\le Ct^2/k
 \quad(|t|\le t_0),\qquad
 \left|\sum_{k=i}^{j-1}\widehat\mu_k+\frac a2\log(j/i)\right|\le C.
 \label{B:eq:laguerre-initial-mgf}
\end{equation}

Let $\widehat{\mathcal F}_j$ be generated by the first $j$ coefficient pairs, and fix $0<t\le t_0$. For $M_j=\exp(tS_j-Ct^2\sum_{k=1}^j k^{-1})$, independence gives
\[
 \mathbb E[M_j\mid\widehat{\mathcal F}_{j-1}]
 =M_{j-1}e^{-Ct^2/j}
       \mathbb E e^{t(\widehat Z_j-\widehat\mu_j)}
 \le M_{j-1},\qquad M_0=1.
\]
Thus $M_j$ is a nonnegative supermartingale. Fix $v>0$, $L\ge0$ and an integer $N\ge2$, and let $\tau=\inf\{1\le j\le N:S_j>v\log N+L\}$, with $\inf\varnothing=\infty$. The stopped
sum is finite:
\[
 M_{\tau\wedge N}
 =1+\sum_{j=1}^N\mathbf1_{\{\tau\ge j\}}(M_j-M_{j-1}),
 \qquad
 \mathbb E M_{\tau\wedge N}\le1.
\]
Indeed, $\{\tau\ge j\}\in\widehat{\mathcal F}_{j-1}$, so each summand has nonpositive expectation. On $\{\tau\le N\}$ the bound for the harmonic sum
$\sum_{k\le N}k^{-1}\le1+\log N$ yields
\[
 M_{\tau\wedge N}
 \ge\exp\bigl(t(v\log N+L)-Ct^2(1+\log N)\bigr).
\]
Nonnegativity on the complementary event now gives
\[
 \mathbb P\!\left(\max_{j\le N}S_j>v\log N+L\right)
 \le e^{Ct^2-tL}N^{-(tv-Ct^2)}.
\]
Applying the same finite stopping calculation with $-t$ in \eqref{B:eq:laguerre-initial-mgf} bounds the lower deviation of $S_j$. Taking their union at $N=2^l$ proves
\begin{equation}
 \mathbb P\left(\max_{j\le2^l}|S_j|>v\log(2^l)+L\right)
 \le2e^{Ct^2-tL}2^{-l(tv-Ct^2)}.
 \label{B:eq:laguerre-initial-maximal-bound}
\end{equation}

Choose $t=\min\{t_0,v/(2C)\}$ and put $\kappa_v=tv-Ct^2\ge tv/2>0$. Define the nonnegative, measurable variable
\[
 H_v=\left[\sup_{l\ge1}
       \left(\max_{0\le j\le2^l}|S_j|-v\log(2^l)\right)\right]_+.
\]
The union bound and the convergent geometric series give, for $L\ge0$,
\[
 \mathbb P(H_v>L)
 \le2e^{Ct^2-tL}\sum_{l=1}^{\infty}2^{-l\kappa_v}
 =\frac{2e^{Ct^2}}{2^{\kappa_v}-1}e^{-tL}.
\]
Letting $L$ increase proves $H_v<\infty$ almost surely. For $j\ge1$ choose $l$ with $j+1\le2^l<2(j+1)$ and set $C_{\omega,v}=H_v+v\log2$. Then $|S_j|\le C_{\omega,v}+v\log(j+1)$;
the same bound holds for $j=0$. It follows that
\[
 |S_{j-1}-S_{i-1}|
 \le2C_{\omega,v}+v\log i+v\log j
 \le2C_{\omega,v}+2v\log j\qquad(1\le i\le j).
\]
With $v=\eta/2$, add the mean bound in \eqref{B:eq:laguerre-initial-mgf} to obtain
\[
 \left|\log\prod_{k=i}^{j-1}
          \frac{\widehat{\mathsf Y}_k}{\widehat{\mathsf X}_k}
          +\frac a2\log(j/i)\right|
 \le2C_{\omega,\eta/2}+C+\eta\log j.
\]
Exponentiation proves \eqref{B:eq:laguerre-initial-product-bounds}, with the constant $\exp(2C_{\omega,\eta/2}+C)\ge1$. For $i=j$ the logarithmic sum is zero and the product is
one.

For $U\sim\chi_\nu^2$, integration of the gamma density gives
\[
 \mathbb E e^{tU}
 =\frac1{2^{\nu/2}\Gamma(\nu/2)}
       \int_0^\infty u^{\nu/2-1}e^{-(1/2-t)u}\,\dd u
 =(1-2t)^{-\nu/2},\qquad t<1/2.
\]
Exponential Markov with $t=1/4$ gives $\mathbb P(U>4\nu)\le e^{-\nu}\mathbb Ee^{U/4}$. With $t=-1$, the event $U<\nu/4$ implies $e^{-U}>e^{-\nu/4}$. Thus
\[
 \mathbb P(U>4\nu)\le e^{-\nu}2^{\nu/2},\qquad
 \mathbb P(U<\nu/4)\le e^{\nu/4}3^{-\nu/2}.
\]
The exponents $1-\tfrac12\log2$ and $\tfrac12\log3-\tfrac14$ are positive. Away from $J$, both degrees lie between $c_1k$ and $C_1k$, for fixed $0<c_1<C_1<\infty$. Consequently the
sum over $k$ of the probabilities that either coefficient lies outside $[\tfrac12\sqrt{c_1k},\,2\sqrt{C_1k}]$ is finite. The first Borel--Cantelli lemma supplies a finite random
index $K_\omega$ after which both coefficients belong to this interval. Enlarge it to contain $J$ and define
\[
 \begin{split}
 c_\omega&=\min\left\{\frac{\sqrt{c_1}}2,\,
        \min_{1\le k\le K_\omega}
          \frac{\widehat{\mathsf X}_k\wedge\widehat{\mathsf Y}_k}{\sqrt k}
            \right\}>0,\\
 C_\omega&=\max\left\{2\sqrt{C_1},\,
        \max_{1\le k\le K_\omega}
          \frac{\widehat{\mathsf X}_k\vee\widehat{\mathsf Y}_k}{\sqrt k}
            \right\}<\infty.
 \end{split}
\]
Positivity and finiteness follow from the positive degrees and the finite number of variables in these extrema.

Finally, intersect the coordinate event with the events $H_{1/(2m)}<\infty$, $m\ge1$. The intersection has probability one and gives the product estimate for $\eta=1/m$ for every
$m$. For any $\eta>0$, choose $m$ with $1/m\le\eta$. Since $j\ge1$, the estimates with $j^{1/m}$ and $j^{-1/m}$ imply the same bounds with $j^\eta$ and $j^{-\eta}$.
\end{proof}

\subsubsection{Gaussian approximation of the scalar density}
\label{B:subsection:laguerre-density-comparison}

We compare the standardized logarithmic ratio density with the Gaussian density and then use this comparison to estimate the quantile map. Denote the distribution function and
density of $\xi_k=(Z_k-\mu_k)/\sigma_k$ by $F_k$ and $g_k$, and the standard Gaussian distribution function and density by $\Phi$ and $\varphi$.

\begin{lemma}\label{B:lem:laguerre-density-comparison}
There are $c,C>0$ and $k_0$ such that, for $k\ge k_0$ and $|x|\le ck^{1/12}$,
\begin{equation}
 \left|\log\frac{g_k(x)}{\varphi(x)}\right|
       \le\frac Ck(1+|x|^5),\qquad
 |F_k(x)-\Phi(x)|\le\frac Ck(1+|x|^5)\varphi(x).
 \label{B:eq:laguerre-density-CDF}
\end{equation}
\end{lemma}
\begin{proof}
Write the density in \eqref{B:eq:laguerre-ratio-density} as a constant times $e^{-V_k}$, with
\begin{equation}
 V_k(z)=S_k\log(1+e^{2z})-2\alpha_kz,\qquad
 z_k^*=\tfrac12\log(\alpha_k/\gamma_k),\qquad
 \varkappa_k=\frac{4\alpha_k\gamma_k}{S_k},\qquad
 s_k=\varkappa_k^{-1/2}.
 \label{B:eq:laguerre-mode-curvature}
\end{equation}
The identity $q'=2q(1-q)$ gives
\begin{align*}
 V_k'&=2(S_kq-\alpha_k),&
 V_k''&=4S_kq(1-q),\\
 V_k'''&=8S_kq(1-q)(1-2q),&
 V_k^{(4)}&=16S_kq(1-q)(1-6q+6q^2),\\
 V_k^{(5)}&=32S_kq(1-q)(1-14q+36q^2-24q^3).
\end{align*}
For example, differentiating a polynomial $P(q(z))$ applies $2q(1-q)\partial_q$ to $P$. This verifies the successive displayed polynomials. The equation $q(z)=\alpha_k/S_k$ has the
unique solution $z=z_k^*$, and $V_k''>0$ on $\mathbb R$. The curvature at this point is $\varkappa_k$, so $s_k\asymp k^{-1/2}$. Since $0<q<1$ and $S_k=O(k)$, the displayed
derivatives of orders two through five are bounded by $Ck$ on $\mathbb R$. At the mode one has
\begin{equation}
 V_k'''(z_k^*)=
 \frac{8\alpha_k\gamma_k(\gamma_k-\alpha_k)}{S_k^2}=O(1),
 \label{B:eq:laguerre-cubic-cancellation}
\end{equation}
Here $\gamma_k-\alpha_k=\beta a/2$ is fixed. This bound is used before multiplying by the cube of $s_k$.

Set $W_k(y)=V_k(z_k^*+s_ky)-V_k(z_k^*)$. Taylor's integral formula gives
\begin{equation}
 W_k(y)=\frac{y^2}{2}+c_{3,k}y^3+c_{4,k}y^4+\mathcal E_k(y),
 \qquad c_{3,k}=\frac{s_k^3V_k'''(z_k^*)}{6},\quad
 c_{4,k}=\frac{s_k^4V_k^{(4)}(z_k^*)}{24},
 \label{B:eq:laguerre-mode-Taylor}
\end{equation}
with
\[
 \mathcal E_k(y)=\frac{(s_ky)^5}{4!}
       \int_0^1(1-t)^4 V_k^{(5)}(z_k^*+ts_ky)\,\dd t,
 \qquad |\mathcal E_k(y)|\le Ck^{-3/2}|y|^5.
\]
The two coefficients can be evaluated at $q=\alpha_k/S_k$:
\[
 c_{3,k}=\frac{\gamma_k-\alpha_k}{6\sqrt{\alpha_k\gamma_k S_k}},
 \qquad
 c_{4,k}=\frac{S_k^2-6\alpha_k\gamma_k}{24\alpha_k\gamma_k S_k}.
\]
Thus $|c_{3,k}|\le Ck^{-3/2}$ and $|c_{4,k}|\le Ck^{-1}$.

We estimate the normalizing integral on the whole real line. Fix $\delta>0$. The points $z_k^*$ stay in a bounded interval, and $q(1-q)$ has a positive minimum on its closed
$\delta$-neighborhood. For all large $k$,
\[
 c_1k\le V_k''(z)\le C_1k\quad(|z-z_k^*|\le\delta),
 \qquad c\le k s_k^2\le C.
\]
Integration from $V_k'(z_k^*)=0$ gives $W_k(y)\ge c_2y^2$ for $|y|\le y_{0,k}:=\delta/s_k$. At the two boundary points it gives
\[
 V_k(z_k^*\pm\delta)-V_k(z_k^*)\ge c_1k\delta^2/2,
 \quad V_k'(z_k^*+\delta)\ge c_1k\delta,
 \quad V_k'(z_k^*-\delta)\le-c_1k\delta.
\]
The graph of the convex function $V_k$ lies above its tangent at each boundary point. The preceding bounds on these tangents give
\[
 W_k(y)\ge c_3k+c_4\sqrt k(|y|-y_{0,k})
                    \quad(|y|\ge y_{0,k}).
\]
For fixed $b>0$, take $L_k=bk^{1/12}\le y_{0,k}$. Integrating the quadratic and linear lower bounds proves
\begin{equation}
 \int_{|y|>L_k}e^{-W_k(y)}\,\dd y
 \le2\int_{L_k}^{\infty}e^{-c_2y^2}\,\dd y
       +\frac{2e^{-c_3k}}{c_4\sqrt k}
 \le C_b e^{-c_b k^{1/6}}.
 \label{B:eq:laguerre-density-tail}
\end{equation}
The Gaussian integral in this display is at most $e^{-c_2L_k^2}/(2c_2L_k)$, since $y/L_k\ge1$ on its domain. The two linear tails each integrate to $e^{-c_3k}/(c_4\sqrt k)$.

On $[-L_k,L_k]$ put $E_k(y)=W_k(y)-y^2/2$. Its three terms give
\[
 \sup_{|y|\le L_k}|E_k(y)|
 \le C_b(k^{-5/4}+k^{-2/3}+k^{-13/12})\le C_b k^{-2/3}.
\]
For $r\ge0$, the substitution $v=y^2/2$ gives
\[
 \int_{\mathbb R}|y|^re^{-y^2/2}\,\dd y
 =2^{(r+1)/2}\Gamma((r+1)/2).
\]
In particular, squaring the three terms in $E_k$ and integrating gives
\begin{align*}
 \int_{-L_k}^{L_k}e^{-y^2/2}|E_k(y)|^2\,\dd y
 &\le C\int_{\mathbb R}e^{-y^2/2}
       (k^{-3}y^6+k^{-2}y^8+k^{-3}y^{10})\,\dd y\le Ck^{-2},\\
 \int_{-L_k}^{L_k}e^{-y^2/2}|\mathcal E_k(y)|\,\dd y
 &\le Ck^{-3/2}\int_{\mathbb R}|y|^5e^{-y^2/2}\,\dd y
 \le Ck^{-3/2}.
\end{align*}
Taylor's formula also gives $|e^{-v}-1+v|\le e^{|v|}v^2/2$. Hence
\begin{align*}
 \int_{-L_k}^{L_k}e^{-W_k(y)}\,\dd y
 &=\int_{-L_k}^{L_k}e^{-y^2/2}\,\dd y
      -\int_{-L_k}^{L_k}e^{-y^2/2}E_k(y)\,\dd y+O(k^{-2})\\
 &=\sqrt{2\pi}(1-3c_{4,k})+O(k^{-3/2}).
\end{align*}
The cubic integral is zero by oddness. The full fourth moment is $3\sqrt{2\pi}$; for each fixed $r$ the omitted integral of $|y|^re^{-y^2/2}$ is at most $C_r e^{-cL_k^2}$.
Combining with \eqref{B:eq:laguerre-density-tail} gives
\[
 I_k:=\int_{\mathbb R}e^{-W_k(y)}\,\dd y
 =\sqrt{2\pi}(1-3c_{4,k})+O(k^{-3/2}),\qquad
 |\log(I_k/\sqrt{2\pi})|\le C/k.
\]
For the logarithmic bound use $|\log(1+v)|\le2|v|$ when $|v|\le1/2$, which applies for all large $k$. The density $\widetilde g_k=I_k^{-1}e^{-W_k}$ of $(Z_k-z_k^*)/s_k$ therefore
satisfies
\begin{equation}
 \left|\log\frac{\widetilde g_k(y)}{\varphi(y)}\right|
 \le Ck^{-1}(1+|y|^4)+Ck^{-3/2}(|y|^3+|y|^5)
       \quad(|y|\le L_k).
 \label{B:eq:laguerre-density-at-mode}
\end{equation}

To recenter at the mean and normalize by the variance, we use Equations~\eqref{B:eq:6-52}--\eqref{B:eq:6-53}, which give
\[
 \mu_k-z_k^*
 =-\frac{\gamma_k-\alpha_k}{4\alpha_k\gamma_k}+O(k^{-3}),\qquad
 \sigma_k^2-s_k^2
 =\tfrac18(\alpha_k^{-2}+\gamma_k^{-2})+O(k^{-3}),\qquad
 s_k^2=\tfrac14(\alpha_k^{-1}+\gamma_k^{-1}).
\]
Put $A_k=(\mu_k-z_k^*)/s_k$ and $B_k=\sigma_k/s_k$. Then $A_k=O(k^{-3/2})$, $B_k=1+O(k^{-1})$, and
\[
 g_k(x)=B_k\widetilde g_k(A_k+B_kx),\qquad
 \log\frac{g_k(x)}{\varphi(x)}
 =\log B_k+\log\frac{\widetilde g_k(A_k+B_kx)}{\varphi(A_k+B_kx)}
        -\frac{(A_k+B_kx)^2-x^2}{2}.
\]
The square difference is $A_k^2+2A_kB_kx+(B_k^2-1)x^2$, of absolute value at most $C(k^{-3}+k^{-3/2}|x|+k^{-1}x^2)$. Choose $0<b_1<b/4$. For large $k$, $1/2\le B_k\le2$ and
$|A_k|\le b_1k^{1/12}/4$, so $|x|\le b_1k^{1/12}$ implies $|A_k+B_kx|\le L_k$. Substitution into \eqref{B:eq:laguerre-density-at-mode} proves the logarithmic density bound on this
interval. Its right side is at most $Ck^{-7/12}$. Using $|e^v-1|\le e^{|v|}|v|$ gives
\begin{equation}
 |g_k(x)-\varphi(x)|\le Ck^{-1}(1+|x|^5)\varphi(x).
 \label{B:eq:laguerre-density-difference}
\end{equation}
For $|x|>b_1k^{1/12}$ the change of variable $x\mapsto A_k+B_kx$ gives $|A_k+B_kx|>b_1k^{1/12}/4$. Applying \eqref{B:eq:laguerre-density-tail} with $b=b_1/4$ and using the positive
lower bound on $I_k$ shows
\[
 \int_{|x|>b_1k^{1/12}}g_k(x)\,\dd x
 +\int_{|x|>b_1k^{1/12}}\varphi(x)\,\dd x
 \le Ce^{-c_5 k^{1/6}}
\]
for a fixed $c_5>0$.

Let $b_1k^{1/12}=L$. For $0\le x\le L$, normalization of the two densities gives
\[
 F_k(x)-\Phi(x)=-\int_x^L(g_k(t)-\varphi(t))\,\dd t
                     -\int_L^\infty(g_k(t)-\varphi(t))\,\dd t.
\]
The two Gaussian tail integrals needed to bound its first term are
\[
 \int_x^\infty\varphi(t)\,\dd t
 =\varphi(x)\int_0^\infty e^{-xv-v^2/2}\,\dd v
 \le\sqrt{\pi/2}\,\varphi(x),\qquad
 \int_x^\infty t^5\varphi(t)\,\dd t
 =(x^4+4x^2+8)\varphi(x).
\]
The second identity follows by three integrations by parts using $\varphi'=-t\varphi$. Thus
\[
 |F_k(x)-\Phi(x)|
 \le Ck^{-1}(1+x^4)\varphi(x)+Ce^{-c_5k^{1/6}}.
\]
For $-L\le x\le0$, use instead
\[
 F_k(x)-\Phi(x)=\int_{-L}^x(g_k(t)-\varphi(t))\,\dd t
                      +\int_{-\infty}^{-L}(g_k(t)-\varphi(t))\,\dd t.
\]
The substitution $u=-t$ bounds the first integral by $Ck^{-1}\int_{|x|}^\infty(1+u^5)\varphi(u)\,\dd u$. The Gaussian integrals evaluated above, now with lower boundary $|x|$, give
for $-L\le x\le0$
\[
 |F_k(x)-\Phi(x)|
 \le\frac Ck\bigl(1+|x|^4+4|x|^2+8\bigr)\varphi(|x|)
                         +Ce^{-c_5k^{1/6}}
 \le\frac {C'}k(1+|x|^4)\varphi(x)+Ce^{-c_5k^{1/6}}.
\]
Choose $0<b_2<b_1$ with $d:=c_5-b_2^2/2>0$. On $|x|\le b_2k^{1/12}$,
\[
 \frac{e^{-c_5k^{1/6}}}{\varphi(x)}
 \le\sqrt{2\pi}e^{-d k^{1/6}}
 \le\frac{\sqrt{2\pi}\,6!}{d^6 k},
\]
where the last bound is a term of the exponential series. Both signs of $x$ now give the stated bound for the distribution function, after decreasing the interval constant to
$b_2$.
\end{proof}

\subsubsection{The quantile coupling}
\label{B:subsection:laguerre-quantile-coupling}

The density $f_k$ is strictly positive on $\mathbb R$. Thus $F_k$ is continuous and strictly increasing, and the map
\begin{equation}
 \tau_k(g)=F_k^{-1}(\Phi(g))
 \label{B:eq:laguerre-quantile-map}
\end{equation}
is continuous. For $G\sim N(0,1)$ it satisfies $\mathbb P(\tau_k(G)\le x)=\mathbb P(\Phi(G)\le F_k(x))=F_k(x)$.

\begin{proposition}\label{B:prop:laguerre-quantile-bound}
For suitable $c,C>0$ and all sufficiently large $k$,
\begin{equation}
 |\tau_k(g)-g|\le Ck^{-1}(1+|g|^5)
       \quad(|g|\le ck^{1/12}).
 \label{B:eq:laguerre-quantile-pointwise}
\end{equation}
For each fixed $p\ge1$, $\|\tau_k(G)-G\|_p\le C_p/k$ for every $k\ge1$. Moreover, with constants independent of $p,k$,
\begin{equation}
 \|\tau_k(G)-G\|_p\le Cp^{5/2}/k,
       \qquad 2\le p\le c_0 k^{1/12}.
 \label{B:eq:laguerre-quantile-growing}
\end{equation}
\end{proposition}
\begin{proof}
Let $c_d,C_d$ be interval and error constants in Lemma~\ref{B:lem:laguerre-density-comparison}. Fix $0<c<c_d/4$ and $D\ge4C_d$. For $|g|\le ck^{1/12}$ put $r=Dk^{-1}(1+|g|^5)$.
Uniformly over this interval,
\[
 r\le C_D k^{-7/12},\qquad |g|r\le C_D k^{-1/2},\qquad
 [g-r,g+r]\subset[-c_dk^{1/12},c_dk^{1/12}]
\]
for all sufficiently large $k$. If $|y-g|\le r$, then
\[
 \frac{\varphi(y)}{\varphi(g)}
 =e^{-(y^2-g^2)/2}\ge e^{-|g|r-r^2/2}\ge\frac12,
 \qquad
 \frac{g_k(y)}{\varphi(y)}
 \ge\exp[-C_d k^{-1}(1+|y|^5)]\ge\frac12.
\]
The starting index in these bounds depends on the fixed $D$. Integration over each interval of length $r$ gives
\[
 F_k(g+r)-F_k(g)\ge r\varphi(g)/4,\qquad
 F_k(g)-F_k(g-r)\ge r\varphi(g)/4.
\]
Since $|F_k(g)-\Phi(g)|\le C_d k^{-1}(1+|g|^5)\varphi(g)$, our choice of $D$ implies $F_k(g-r)\le\Phi(g)\le F_k(g+r)$. Strict monotonicity of $F_k$ then puts $\tau_k(g)$ in
$[g-r,g+r]$, proving \eqref{B:eq:laguerre-quantile-pointwise}.

Set $A_k=\{|G|\le ck^{1/12}\}$. The Gaussian moment identity
\[
 \mathbb E|G|^r=2^{r/2}\Gamma((r+1)/2)/\sqrt\pi
\]
and the gamma bound proved above give $\|G\|_r\le C\sqrt r$ for $r\ge2$; for $1\le r\le2$, use $\|G\|_r\le\|G\|_2=1$. Consequently
\[
 \|1+|G|^5\|_p\le1+\|G\|_{5p}^5\le Cp^{5/2},\qquad
 \|(\tau_k(G)-G)\mathbf1_{A_k}\|_p\le Cp^{5/2}/k.
\]
For the complement, H\"older's inequality gives $\|X\mathbf1_{A_k^c}\|_p\le\|X\|_{2p}\mathbb P(A_k^c)^{1/(2p)}$. The exact marginal law is $\tau_k(G)\stackrel{\mathrm{law}}=\xi_k$;
both terms in the difference are functions of this same $G$. The Gaussian tail and \eqref{B:eq:laguerre-log-moments} yield
\begin{equation}
 \begin{split}
 \|(\tau_k(G)-G)\mathbf1_{A_k^c}\|_p
 &\le(\|\xi_k\|_{2p}+\|G\|_{2p})\mathbb P(A_k^c)^{1/(2p)}\\
 &\le C(\sqrt p+p/\sqrt k)e^{-c_1 k^{1/6}/p}.
 \end{split}
 \label{B:eq:laguerre-quantile-tail-moment}
\end{equation}
For fixed $p$, the exponential series implies $e^{-c_1 k^{1/6}/p}\le6!(p/c_1)^6/k$. The prefactor in \eqref{B:eq:laguerre-quantile-tail-moment} is bounded for this fixed $p$ and
$k\ge1$, so the tail contribution is
$O_p(k^{-1})$.

Choose a fixed $0<c_0\le1$. If $2\le p\le c_0k^{1/12}$, then $p/\sqrt k\le\sqrt p$ and
\[
 e^{-c_1k^{1/6}/p}
 \le e^{-(c_1/c_0)k^{1/12}}
 \le\frac{12!(c_0/c_1)^{12}}{k}.
\]
The tail contribution is therefore at most $C\sqrt p/k$ throughout this range. Adding the contribution on $A_k$ proves \eqref{B:eq:laguerre-quantile-growing} for all sufficiently
large $k$. For the remaining finite set $k<K$, the permitted $p$ is bounded by $c_0K^{1/12}$. Put $P=\max\{2,c_0K^{1/12}\}$. The exact law and \eqref{B:eq:laguerre-log-moments}
give
\[
 \max_{1\le k<K}k\,\|\tau_k(G)-G\|_P<\infty.
\]
Monotonicity of $L^p$ norms on a probability space covers every permitted $p$ at these indices after increasing one constant. For a fixed $p\ge1$, let $K_p$ be the threshold for
the large-$k$ estimate and define
\[
 C_p^{\rm init}=\max_{1\le k<K_p}
       k\bigl(\|\xi_k\|_p+\|G\|_p\bigr)<\infty.
\]
The maximum over an empty set is zero. The exact marginal $\tau_k(G)\stackrel{\rm law}=\xi_k$ gives $\|\tau_k(G)-G\|_p\le C_p^{\rm init}/k$ at each remaining index. Taking the
larger of this and the large-$k$ constant proves the bound for all $k$.

The same map $\tau_k$ has been used for every moment order. Its exact standardized marginal gives $\mathbb E\tau_k(G)=0$ and $\mathbb E\tau_k(G)^2=1$. In particular,
\[
 \mathbb E[\tau_k(G)-G]=0,\qquad
 \mathbb E|\tau_k(G)-G|^2\le Ck^{-2},\qquad
 \mathbb E|\sigma_k(\tau_k(G)-G)|^2\le Ck^{-3}.
\]
Here the final inequality uses $\sigma_k^2\le C/k$ from Lemma~\ref{B:lem:laguerre-means-variances}. To compare with the Gaussian parameters in
\eqref{B:eq:laguerre-gaussian-parameters}, write
\[
 \mu_k+\sigma_k\tau_k(G)-d_k-b_kG
 =\mu_k-d_k+\sigma_k(\tau_k(G)-G)+(\sigma_k-b_k)G.
\]
The last two terms have mean zero. Their norm is bounded by
\[
 \|\sigma_k(\tau_k(G)-G)+(\sigma_k-b_k)G\|_p
 \le\sigma_k\|\tau_k(G)-G\|_p+|\sigma_k-b_k|\|G\|_p
 \le C p^{5/2}k^{-3/2}
\]
for $2\le p\le c_0k^{1/12}$, and by $C_pk^{-3/2}$ for each fixed $p\ge1$. The remaining deterministic term satisfies $|\mu_k-d_k|\le Ck^{-2}$ by
Lemma~\ref{B:lem:laguerre-means-variances}. There is also an exact second moment identity for this coupling. Both $G$ and $\tau_k(G)$ are centered with variance one, so
\[
 2\mathbb E[G(\tau_k(G)-G)]
       =-\mathbb E|\tau_k(G)-G|^2.
\]
Expanding the square gives
\[
 \begin{split}
 &\mathbb E|\mu_k+\sigma_k\tau_k(G)-d_k-b_kG|^2\\
 &\quad=(\mu_k-d_k)^2+(\sigma_k-b_k)^2
              +\sigma_kb_k\mathbb E|\tau_k(G)-G|^2=O(k^{-3}).
 \end{split}
\]
The two random summands in the centered error are functions of the same $G$, and the identity accounts for their covariance without assuming independence.
\end{proof}

\subsubsection{Construction for all matrix sizes}
\label{B:subsection:laguerre-array-construction}

We apply the scalar quantile maps to increments of one Brownian path. The increments used at a fixed size are independent. Different sizes use the same path, which determines their
dependence. The initial coefficient sequence and the remaining Gamma variables give the finite coefficient distribution.

Fix $A>12$ and let $r_n=\lceil(\log n)^A\rceil$. The bound $r_n/n\le((\log n)^A+1)/n\to0$ permits a deterministic choice of $n_0\ge2$ such that, for $n\ge n_0$,
\begin{equation}
 m_n=r_n<M_n=n-r_n<n,\qquad
 \mathscr I_n^{\mathrm{ret}}=\{m_n+1,\ldots,M_n\}.
 \label{B:eq:laguerre-cuts}
\end{equation}
These are the intervals in \eqref{B:eq:brownian-retained-interval} on the partition $z_{n,k}=k/n$. In particular, the two definitions using $\log n$ and $\log(n+2)$ give comparable
cutoffs. The Brownian times are $u_{n,k}=\log(n/k)$, $1\le k\le n$, and $u_{n,n}=0$. Set
\begin{equation}
 G_{n,k}=\frac{W(u_{n,k})-W(u_{n,k+1})}{\sqrt{h_k^+}},
       \qquad 1\le k<n.
 \label{B:eq:laguerre-normal-increments}
\end{equation}
For a fixed $n$, the intervals $[u_{n,k+1},u_{n,k}]$, $1\le k<n$, have disjoint interiors and lengths $h_k^+$. The Brownian covariance gives
\[
 \mathbb E G_{n,k}=0,\qquad
 \mathbb E[G_{n,k}G_{n,\ell}]
 =\frac{|[u_{n,k+1},u_{n,k}]\cap[u_{n,\ell+1},u_{n,\ell}]|}
        {\sqrt{h_k^+h_\ell^+}}
 =\mathbf1_{\{k=\ell\}}.
\]
Here the bars denote interval length. Hence, for any real $t_1,\ldots,t_{n-1}$,
\[
 \mathbb E\exp\left(\ii\sum_{k=1}^{n-1}t_kG_{n,k}\right)
 =\exp\left(-\frac12\sum_{k=1}^{n-1}t_k^2\right).
\]
The vector consists of independent standard Gaussians. Conversely, the increments reconstruct every observation:
\[
 W(u_{n,j})=\sum_{k=j}^{n-1}\sqrt{h_k^+}\,G_{n,k},\qquad
 \mathcal N_n=\sigma(G_{n,k}:1\le k<n).
\]
The last sum is empty at $j=n$.

On a product extension independent of $W$, take the unchanged sequence from Lemma~\ref{B:lem:laguerre-initial-products}, mutually independent variables
\begin{equation}
 T_{n,k}\sim\Gamma(S_k,2)\quad(n\ge n_0,\ m_n<k<n),\qquad
 V_n\sim\Gamma(\gamma_n,2)\quad(n\ge n_0),
 \label{B:eq:laguerre-auxiliary-variables}
\end{equation}
and independent arrays with the exact chi laws for $n<n_0$. More explicitly, take the countable product of copies of $(0,\infty)$ carrying the initial chi variables, the displayed
gamma variables, and the finitely many early rows, with their marginal distributions. The probability space is the product of this space with the space carrying $W$. If
$U_1,\ldots,U_q$ are distinct auxiliary coordinates with marginal laws $\nu_1,\ldots,\nu_q$, then, for Borel sets $A_1,\ldots,A_q$ and $E\in\sigma\{W(t):t\ge0\}$,
\[
 \mathbb P\left(E\cap\bigcap_{j=1}^q\{U_j\in A_j\}\right)
 =\mathbb P(E)\prod_{j=1}^q\nu_j(A_j).
\]
Finite cylinder sets generate the auxiliary sigma-field, denoted by $\mathcal U$. The monotone class theorem therefore gives $\mathcal U\perp\sigma\{W(t):t\ge0\}$. This specifies
the joint law of all auxiliary variables, including their independence across sizes.

For $m_n<k<n$, define
\begin{equation}
 Z_{n,k}=\mu_k+\sigma_k\tau_k(G_{n,k}),\qquad Q_{n,k}=q(Z_{n,k}).
 \label{B:eq:laguerre-coupled-ratios}
\end{equation}
The quantile map $\tau_k$ is the one fixed in \eqref{B:eq:laguerre-quantile-map}. For $n\ge n_0$ put
\begin{equation}
 (X_{n,k}^{\mathrm L},Y_{n,k}^{\mathrm L})=
 \begin{cases}
 (\widehat{\mathsf X}_k,\widehat{\mathsf Y}_k)/\sqrt{n\beta},
                              &1\le k\le m_n,\\[1mm]
 (\sqrt{T_{n,k}(1-Q_{n,k})},\sqrt{T_{n,k}Q_{n,k}})/\sqrt{n\beta},
                              &m_n<k<n,
 \end{cases}
 \qquad X_{n,n}^{\mathrm L}=\sqrt{V_n/(n\beta)}.
 \label{B:eq:laguerre-final-array}
\end{equation}
Use the separately sampled rows for $n<n_0$. The gamma variable $V_n$ supplies the last diagonal coefficient. There is no off-diagonal coefficient with index $n$. The second branch
of \eqref{B:eq:laguerre-final-array} is used at every $m_n<k<n$, including the indices beyond $M_n$; the latter cutoff specifies the interval used in the inverse comparison.

Let $B_n^{\mathrm L}$ be the bidiagonal matrix defined in Section~\ref{B:sec:inverse-representations-and-the-stochastic-bessel-operator} from these coefficients. We retain the
products $P_{n,j}$, the matrices $\mathcal A_n^{\mathrm L}$, and their embeddings; here $\alpha_n=1$ and $\mathscr J_n=\widetilde{\mathscr J}_n=\mathscr E_n$. For a fixed $k$, the
identities
\[
 \sqrt{n\beta}\,X_{n,k}^{\mathrm L}=\widehat{\mathsf X}_k,\qquad
 \sqrt{n\beta}\,Y_{n,k}^{\mathrm L}=\widehat{\mathsf Y}_k
\]
hold for every sufficiently large $n$, since $m_n\to\infty$. The product estimate from Lemma~\ref{B:lem:laguerre-initial-products} therefore applies with one random constant to
\[
 \prod_{k=i}^{j-1}\frac{Y_{n,k}^{\mathrm L}}{X_{n,k}^{\mathrm L}}
 =\prod_{k=i}^{j-1}\frac{\widehat{\mathsf Y}_k}{\widehat{\mathsf X}_k},
 \qquad n\ge n_0,\quad 1\le i\le j\le m_n+1.
\]
Its coordinate bounds give, on the same event and for $1\le k\le m_n$,
\[
 c_\omega\sqrt{\frac{k}{n\beta}}
 \le X_{n,k}^{\mathrm L},Y_{n,k}^{\mathrm L}
 \le C_\omega\sqrt{\frac{k}{n\beta}}.
\]

\begin{proposition}\label{B:prop:laguerre-array-law}
For every $n$ the coefficients in the $n$th row are independent, with
\begin{equation}
 \sqrt{n\beta}X_{n,k}^{\mathrm L}\sim\chi_{\beta(k+a)}\quad(k\le n),
 \qquad
 \sqrt{n\beta}Y_{n,k}^{\mathrm L}\sim\chi_{\beta k}\quad(k<n).
 \label{B:eq:laguerre-exact-row-law}
\end{equation}
The row and its embedded inverse are $\mathcal N_n\vee\mathcal U$-measurable, where $\mathcal N_n$ is defined in \eqref{B:eq:brownian-sampled-information}. The Brownian bridges
between its consecutive sampling times are jointly independent of this sigma-field.
\end{proposition}
\begin{proof}
Fix $n\ge n_0$. By \eqref{B:eq:laguerre-quantile-map}, $\tau_k(G_{n,k})$ has the law of $(Z_k-\mu_k)/\sigma_k$. Thus $Z_{n,k}$ has density $f_k$. It is independent of $T_{n,k}$,
because it is a function of $W$ and $T_{n,k}$ is an auxiliary coordinate. For bounded Borel functions $F,G$ and $m_n<k<n$, this gives
\[
 \begin{split}
 &\mathbb E\!\left[
 F(\sqrt{n\beta}\,X_{n,k}^{\mathrm L})
 G(\sqrt{n\beta}\,Y_{n,k}^{\mathrm L})\right]\\
 &\quad=\int_{\mathbb R}\int_0^\infty
 F\bigl(\sqrt{t(1-q(z))}\bigr)G\bigl(\sqrt{tq(z)}\bigr)
 \frac{t^{S_k-1}e^{-t/2}}{2^{S_k}\Gamma(S_k)}\,\dd t\, f_k(z)\,\dd z\\
 &\quad=\mathbb E F(\chi_{\beta(k+a)})\,\mathbb E G(\chi_{\beta k}).
 \end{split}
\]
The last equality is the change of variables \eqref{B:eq:laguerre-density-factorization} in Proposition~\ref{B:prop:laguerre-scalar-reconstruction}. It proves independence of the
two reconstructed coefficients within each pair.

The pairs $(G_{n,k},T_{n,k})$, $m_n<k<n$, are mutually independent. The product construction also makes them independent of the initial pairs with $k\le m_n$ and of $V_n$. Since
$V_n\sim\Gamma(\gamma_n,2)$ and $2\gamma_n=\beta(n+a)$, its square root has the required chi law for the last coefficient. Taking the product of these integrals gives, for bounded
Borel functions $f_k,g_k$,
\[
 \mathbb E\left[\prod_{k=1}^n f_k(\sqrt{n\beta}X_{n,k}^{\mathrm L})
                    \prod_{k=1}^{n-1}g_k(\sqrt{n\beta}Y_{n,k}^{\mathrm L})\right]
 =\prod_{k=1}^n\mathbb Ef_k(\chi_{\beta(k+a)})
        \prod_{k=1}^{n-1}\mathbb Eg_k(\chi_{\beta k}).
\]
Rectangles generate the finite coefficient space, so this identity is its full joint law. The separately sampled rows satisfy the same identity when $n<n_0$, including $n=1$, where
the second product is empty.

For the conditional law, write $m=m_n$ and keep $n\ge n_0$ fixed. For bounded Borel functions $F_k$ of each pair and a bounded Borel function $H$ of the last coefficient,
integration over the auxiliary coordinates gives
\[
 \begin{split}
 &\mathbb E\!\left[
   \prod_{k<n}F_k(\sqrt{n\beta}X_{n,k}^{\mathrm L},
                 \sqrt{n\beta}Y_{n,k}^{\mathrm L})
       H(\sqrt{n\beta}X_{n,n}^{\mathrm L})\,\middle|\,\mathcal N_n\right]\\
 &=\left(\prod_{k\le m}
       \mathbb EF_k(\widehat{\mathsf X}_k,\widehat{\mathsf Y}_k)\right)
       \left(\int_0^\infty H(\sqrt v)\,\Gamma(\gamma_n,2)(\dd v)\right)\\
 &\quad\times\prod_{m<k<n}\int_0^\infty
       F_k\!\left(\sqrt{t(1-Q_{n,k})},\sqrt{tQ_{n,k}}\right)
                         \Gamma(S_k,2)(\dd t).
 \end{split}
\]
The equality follows first after multiplication by an indicator in $\mathcal N_n$, using the product law of the auxiliary coordinates, and then for all bounded
$\mathcal N_n$-measurable functions by bounded convergence. Given $\mathcal N_n$, the ratio in each reconstructed pair is fixed and its two entries share the same radial variable.
Averaging this conditional law over the independent $G_{n,k}$ gives the chi product law proved above. The conditional formula is unchanged when $\mathcal N_n$ on the left is
replaced by the sigma-field of the whole Brownian path, since $\mathcal U$ is independent of that path and the displayed integrand uses only $Q_{n,k}$.

Each map in \eqref{B:eq:laguerre-normal-increments}--\eqref{B:eq:laguerre-final-array} is Borel: evaluations of a continuous path, the continuous quantile map, the function $q$,
and positive square roots are Borel. The row uses $W$ only through $\mathcal N_n$, and the remaining arguments of these maps belong to $\mathcal U$. All coefficients are positive
and finite on one event of probability one, obtained by intersecting the corresponding events over the countable array. Hence
\[
 \det B_n^{\mathrm L}=\prod_{k=1}^nX_{n,k}^{\mathrm L}>0,\qquad
 (\mathcal A_n^{\mathrm L})^{-1}
 =n^{-2}(B_n^{\mathrm L})^{-\mathsf T}(B_n^{\mathrm L})^{-1}.
\]
Finite inversion is continuous on this set. For any matrix $D\in\mathbb C^{n\times n}$, the kernel convention of \eqref{B:eq:integral-kernel-convention} gives
\[
 K_{\mathscr E_nD\mathscr E_n^*}(s,t)
 =n\sum_{i,j=1}^nD_{ji}\mathbf1_{I_{n,i}}(s)\mathbf1_{I_{n,j}}(t),
 \qquad
 \|\mathscr E_nD\mathscr E_n^*\|_{\mathfrak S_2}^2
 =\sum_{i,j=1}^n|D_{ij}|^2.
\]
Indeed, each rectangle formed by two cells has area $n^{-2}$, which cancels the squared kernel factor $n^2$. The embedding is therefore continuous in Hilbert--Schmidt norm. This
proves measurability of the inverse.

To verify bridge independence for the complete row, fix $1\le k<n$, put $v=u_{n,k+1}$, $w=u_{n,k}$, $h=w-v=h_k^+$, and define
\[
 \mathbb B_{n,k}^+(t)
 =W(v+t)-W(v)-\frac{t}{h}[W(w)-W(v)],
 \qquad 0\le t\le h.
\]
The covariance calculation in Lemma~\ref{B:lem:conditional-bridge-law} gives, for $s\notin(v,w)$,
\[
 \operatorname{Cov}(\mathbb B_{n,k}^+(t),W(s))=0,\qquad
 \mathbb E[\mathbb B_{n,k}^+(t)\mathbb B_{n,k}^+(t')]
 =t\wedge t'-tt'/h.
\]
Different intervals have zero cross covariance. Joint Gaussianity and continuity give independence of the bridge paths from one another and from $\mathcal N_n$. Since $\mathcal U$
is independent of the entire Brownian path, bounded functions $F,g,H$ of the bridges, $\mathcal N_n$, and $\mathcal U$, respectively, satisfy
\[
 \mathbb E[FgH]=\mathbb E[Fg]\mathbb EH
 =\mathbb EF\,\mathbb E[gH].
\]
The monotone class theorem extends this to $\mathcal N_n\vee\mathcal U$. Thus, writing $\mathsf{Br}_h$ for the centered Brownian bridge law on $C([0,h])$,
\[
 \mathbb E[F(\mathbb B_n^+)\mid\mathcal N_n\vee\mathcal U]
 =\int F(b_1,\ldots,b_{n-1})
           \bigotimes_{k=1}^{n-1}\mathsf{Br}_{h_k^+}(\dd b_k).
\]
The same argument holds for $n<n_0$, whose coefficients are sampled from $\mathcal U$. For $n=1$ the bridge vector is empty and its law is a point mass.

These are the conditional laws used in Theorem~\ref{B:thm:brownian-information-lower-bound}.
\end{proof}

The construction fixes the dependence between rows as well. For two sizes $n,n'\ge2$,
\[
 \operatorname{Cov}(G_{n,k},G_{n',\ell})
 =\frac{|[u_{n,k+1},u_{n,k}]\cap[u_{n',\ell+1},u_{n',\ell}]|}
        {\sqrt{h_k^+h_\ell^+}}.
\]
The initial coefficient pairs are also shared through \eqref{B:eq:laguerre-final-array}. The probability space, cutoff sequence and quantile maps have been fixed independently of
the moment orders and error thresholds used in the estimates that follow. Bridge independence always refers to one row's observation field. Indeed,
$\{\log(n/k):n\ge2,\ 1\le k\le n\}$ is dense in $[0,\infty)$, since the positive rational numbers $n/k\ge1$ are dense in $[1,\infty)$. Continuity of $W$ therefore gives
\[
 \bigvee_{n\ge2}\mathcal N_n=\sigma\{W(u):u\ge0\}.
\]
Conditioning on all these observations determines every bridge. The lower bound therefore uses the individual observation fields $\mathcal N_n\vee\mathcal U$ for the already fixed
sequence of finite inverses. All subsequent Borel--Cantelli arguments use summable probability bounds and require no independence between rows.

\subsubsection{Sums of logarithmic ratio errors}
\label{B:subsection:laguerre-product-errors}

The products in the inverse kernel involve sums from the current index to the last coefficient. We estimate these sums first, then the individual logarithmic coefficients appearing
in the scale error. Write $\ell_n=\log(e+n)$ throughout the remainder of this subsection. For the beta-Laguerre normalization, $p_{n,j}=P_{n,j}$. With $c_n=0$, the error in
\eqref{B:eq:delta-explicit-coefficients} is
\begin{equation}
 \delta_{n,j}=\log P_{n,j}-\frac a2\log(j/n)
                 -\beta^{-1/2}W(\log(n/j)),\qquad1\le j\le n.
 \label{B:eq:laguerre-delta}
\end{equation}
We use this formula to extend the error to all indices; at $j=n$ it is zero. For fixed $1<\lambda<2$, put
\begin{equation}
 \mathcal P_{n,\lambda}=\sum_{j=1}^{n-1}
                 h_j^+(j/n)^\lambda|\delta_{n,j}|^2.
 \label{B:eq:laguerre-product-energy}
\end{equation}

\begin{proposition}\label{B:prop:laguerre-product-errors}
For the constructed arrays and every fixed $1<\lambda<2$, almost surely,
\begin{equation}
 n\mathcal P_{n,\lambda}\longrightarrow0,\qquad
 \max_{m_n<j\le n}|\delta_{n,j}|\longrightarrow0.
 \label{B:eq:laguerre-product-error-limits}
\end{equation}
\end{proposition}
\begin{proof}
For $m_n<k<n$ use the decomposition
\begin{equation}
 \begin{split}
 e_{n,k}&=Z_{n,k}-d_k-b_kG_{n,k}
            =\widetilde e_{n,k}+\mu_k-d_k,\\
 \widetilde e_{n,k}
   &=\sigma_k[\tau_k(G_{n,k})-G_{n,k}]
                         +(\sigma_k-b_k)G_{n,k}.
 \end{split}
 \label{B:eq:laguerre-increment-error}
\end{equation}
For fixed $n$, each $\widetilde e_{n,k}$ is a function of $G_{n,k}$. Proposition~\ref{B:prop:laguerre-array-law} gives independence across $k$, and the exact quantile law gives
\[
 \mathbb E\widetilde e_{n,k}
 =\sigma_k(\mathbb E\xi_k-\mathbb EG_{n,k})
   +(\sigma_k-b_k)\mathbb EG_{n,k}=0.
\]
Lemma~\ref{B:lem:laguerre-means-variances} and Proposition~\ref{B:prop:laguerre-quantile-bound} therefore imply
\begin{equation}
 \|\widetilde e_{n,k}\|_p\le C_pk^{-3/2},\qquad
 |\mu_k-d_k|\le Ck^{-2}.
 \label{B:eq:laguerre-increment-moments}
\end{equation}
The identities $\sum_{k=j}^{n-1}\log((k+1)/k)=\log(n/j)$, $u_{n,n}=0$, and $W(0)=0$ give
\begin{equation}
 \sum_{k=j}^{n-1}d_k=\frac a2\log(j/n),\qquad
 \sum_{k=j}^{n-1}b_kG_{n,k}=\beta^{-1/2}W(\log(n/j)),\qquad
 \delta_{n,j}=\sum_{k=j}^{n-1}e_{n,k}\quad(j>m_n).
 \label{B:eq:laguerre-regular-telescope}
\end{equation}
The three sums are empty when $j=n$. For $q>1$ and $j\ge1$, monotonicity yields
\[
 \sum_{k=j}^\infty k^{-q}
 \le j^{-q}+\int_j^\infty x^{-q}\,\dd x
 \le\left(1+\frac1{q-1}\right)j^{1-q}.
\]
For every fixed $p\ge2$, Lemma~\ref{B:lem:independent-sum-moments} and the triangle inequality now give
\[
 \|\delta_{n,j}\|_p
 \le C\sqrt p\left(\sum_{k=j}^{n-1}\|\widetilde e_{n,k}\|_p^2\right)^{1/2}
       +\sum_{k=j}^{n-1}|\mu_k-d_k|
 \le C_p/j \qquad(j>m_n).
\]
The cumulative errors themselves are correlated. For $i,j>m_n$, independence of the centered increments and the second moment identity from the quantile proof give
\[
 \begin{split}
 \operatorname{Cov}(\delta_{n,i},\delta_{n,j})
 &=\sum_{k=i\vee j}^{n-1}
   \left[(\sigma_k-b_k)^2
         +\sigma_kb_k\mathbb E|\tau_k(G)-G|^2\right]\\
 &\le C\sum_{k=i\vee j}^{\infty}k^{-3}
 \le\frac{C}{(i\vee j)^2},\qquad G\sim N(0,1).
 \end{split}
\]
Their means are $\sum_{k=j}^{n-1}(\mu_k-d_k)$ and have absolute value at most $C/j$. Thus the quadratic estimate below uses the triangle inequality for the correlated squared sums,
while independence is used only for the increments. At $i=n$ or $j=n$ the covariance sum is empty, consistently with $\delta_{n,n}=0$. For $p\ge1$, apply the preceding $L^p$ bound
at order $2p$ to each square. Since $h_j^+\le j^{-1}$, the triangle inequality in $L^p$ gives
\begin{equation}
 \begin{split}
 \left\|\sum_{j=m_n+1}^{n-1}h_j^+(j/n)^\lambda
                            |\delta_{n,j}|^2\right\|_p
 &\le C_p n^{-\lambda}\sum_{j=m_n+1}^\infty j^{\lambda-3}\\
 &\le C_{p,\lambda}n^{-\lambda}m_n^{\lambda-2}.
 \end{split}
 \label{B:eq:laguerre-regular-energy-moments}
\end{equation}
Here
\[
 \sum_{j>m}j^{\lambda-3}
 \le\int_m^\infty x^{\lambda-3}\,\dd x
 =\frac{m^{\lambda-2}}{2-\lambda},
 \qquad 1<\lambda<2.
\]
Fix an integer $p$ with $p(\lambda-1)>1$. For every $\epsilon>0$, Markov's inequality gives
\[
 \sum_{n\ge n_0}\mathbb P\left(
 n\sum_{j=m_n+1}^{n-1}h_j^+(j/n)^\lambda|\delta_{n,j}|^2>\epsilon\right) 
 \le C_{p,\lambda}\epsilon^{-p}
 \sum_{n\ge n_0}n^{-p(\lambda-1)}m_n^{-p(2-\lambda)}<\infty.
\]
Borel--Cantelli and intersection over $\epsilon=1/l$, $l\ge1$, prove that $n$ times this sum tends to zero.

For the maximum choose $L>0$ with $L^2/2>3$. The normal tail bound and a union over the $n-1$ increments give
\[
 \sum_{n\ge n_0}\mathbb P\left(\max_{k<n}|G_{n,k}|>L\sqrt{\ell_n}\right)
 \le\sum_{n\ge n_0}2n(e+n)^{-L^2/2}<\infty.
\]
On the resulting event, all increments lie in the stated interval for every $n\ge n_0(\omega)$. The quantile estimate applies at every $k>m_n$, since
\[
 \frac{L\sqrt{\ell_n}}{(m_n+1)^{1/12}}
 \le C_L(\log n)^{1/2-A/12}\longrightarrow0.
\]
Substitute the pointwise quantile bound, $\sigma_k\le Ck^{-1/2}$, and $|\sigma_k-b_k|\le Ck^{-3/2}$ into \eqref{B:eq:laguerre-increment-error}. Summing the resulting powers gives
\begin{equation}
 |e_{n,k}|\le C(1+\ell_n^{5/2})k^{-3/2}+Ck^{-2},\qquad
 |\delta_{n,j}|\le
 \frac{C(1+\ell_n^{5/2})}{\sqrt j}+\frac Cj
        \quad(m_n<j\le n).
 \label{B:eq:laguerre-regular-maximum}
\end{equation}
Indeed, the two tail sums are bounded by $3j^{-1/2}$ and $2j^{-1}$. The maximum in this display is at most
\[
 C\bigl[(\log n)^{-A/2}+(\log n)^{(5-A)/2}+(\log n)^{-A}\bigr]
 \longrightarrow0.
\]

For $1\le j\le m_n$, divide $P_{n,j}$ by $P_{n,m_n+1}$ and sum $d_k$ through $m_n$. This gives the exact identity
\begin{equation}
 \delta_{n,j}=\sum_{k=j}^{m_n}(\widehat Z_k-d_k)
 -\beta^{-1/2}[W(u_{n,j})-W(u_{n,m_n+1})]
 +\delta_{n,m_n+1}.
 \label{B:eq:laguerre-initial-error-decomposition}
\end{equation}
Apply the bound for the initial coefficient products with exponent one. Its logarithm gives, uniformly in these $j$,
\[
 \left|\sum_{k=j}^{m_n}(\widehat Z_k-d_k)\right|
 =\left|\log\prod_{k=j}^{m_n}
     \frac{\widehat{\mathsf Y}_k}{\widehat{\mathsf X}_k}
       -\frac a2\log\frac{j}{m_n+1}\right|
 \le\log C_{\omega,1}+\log(m_n+1).
\]
The growth bound for the Brownian path gives $\sup_{0\le u\le\log n}|W(u)|\le C_\omega(1+\ell_n)$. The two times in \eqref{B:eq:laguerre-initial-error-decomposition} belong to this
interval, and $\delta_{n,m_n+1}\to0$ by \eqref{B:eq:laguerre-regular-maximum}. Thus
\begin{equation}
 \begin{split}
 \max_{j\le m_n}|\delta_{n,j}|&\le C_\omega(1+\ell_n),\\
 n\sum_{j=1}^{m_n}h_j^+(j/n)^\lambda|\delta_{n,j}|^2
 &\le C_{\omega,\lambda}n^{1-\lambda}
                        m_n^\lambda(1+\ell_n)^2\longrightarrow0.
 \end{split}
 \label{B:eq:laguerre-initial-error-energy}
\end{equation}
To evaluate the sum, use $\sum_{j\le m}j^{\lambda-1}\le m^\lambda$ and $h_j^+\le j^{-1}$. For large $n$, $m_n\le2(\log n)^A$ and $\ell_n\le2\log n$; hence its final bound is at
most
\[
 C_{\omega,\lambda}n^{-(\lambda-1)}(\log n)^{A\lambda+2}
 \longrightarrow0.
\]
For every $b,c>0$, this last limit follows by setting $x=\log n$ and choosing an integer $N>b$: the exponential series implies $x^b e^{-cx}\le N!c^{-N}x^{b-N}\to0$. Add the sums
over $j\le m_n$ and $j>m_n$ to prove \eqref{B:eq:laguerre-product-error-limits}. All probability estimates above use the original coupled rows; summability is along the full
sequence of matrix sizes.
\end{proof}

\subsubsection{Squared logarithmic coefficients}
\label{B:subsection:laguerre-logarithmic-coefficients}

The cumulative ratio error is only one part of the inverse comparison. The scale coefficient also contains individual logarithmic coefficients, so their weighted square sum must be
controlled separately under the joint distribution of the coefficient vector.

For the constructed coefficient vectors define
\begin{equation}
 \begin{gathered}
 L_{n,k}^X=\log\frac{X_{n,k}^{\mathrm L}}{\sqrt{k/n}}\quad(k\le n),
 \qquad L_{n,k}^Y=\log\frac{Y_{n,k}^{\mathrm L}}{\sqrt{k/n}}\quad(k<n),\\
 \mathcal M_{n,\lambda}=\sum_{k=1}^{n-1}h_k^+(k/n)^\lambda
                   [(L_{n,k}^X)^2+(L_{n,k}^Y)^2].
 \end{gathered}
 \label{B:eq:laguerre-logarithmic-energy}
\end{equation}

\begin{proposition}\label{B:prop:laguerre-logarithmic-energy}
For every fixed $1<\lambda<2$, almost surely,
\begin{equation}
 n\mathcal M_{n,\lambda}\longrightarrow\frac1{\beta(\lambda-1)}.
 \label{B:eq:laguerre-local-energy-limit}
\end{equation}
\end{proposition}
\begin{proof}
Put $A_{n,k}=(L_{n,k}^X)^2+(L_{n,k}^Y)^2$ for $k<n$. Proposition~\ref{B:prop:laguerre-array-law} gives independence of these variables within each row and the scalar marginal laws
for all indices. For every fixed $p\ge1$,
\[
 \|A_{n,k}\|_p
 \le\|L_{n,k}^X\|_{2p}^2+\|L_{n,k}^Y\|_{2p}^2\le C_p/k,
 \qquad
 \|A_{n,k}-\mathbb EA_{n,k}\|_p\le2\|A_{n,k}\|_p.
\]
The mean and variance calculation in Lemma~\ref{B:lem:laguerre-means-variances} consequently gives
\begin{equation}
 \mathbb EA_{n,k}=\frac1{\beta k}+O(k^{-2}),\qquad
 \|A_{n,k}\|_p+\|A_{n,k}-\mathbb EA_{n,k}\|_p\le C_p/k.
 \label{B:eq:laguerre-local-inputs}
\end{equation}
Writing $h_k^+=k^{-1}+\epsilon_k$ with $|\epsilon_k|\le(2k^2)^{-1}$, we have $|h_k^+\mathbb EA_{n,k}-1/(\beta k^2)|\le Ck^{-3}$. Multiplying by $(k/n)^\lambda$ and summing proves
\begin{equation}
 \mathbb E\mathcal M_{n,\lambda}
 =\frac{n^{-\lambda}}\beta\sum_{k=1}^{n-1}k^{\lambda-2}
       +O(n^{-\lambda}).
 \label{B:eq:laguerre-local-mean}
\end{equation}
The remainder is bounded by $Cn^{-\lambda}\sum_{k\ge1}k^{\lambda-3}\le C_\lambda n^{-\lambda}$. For $v=\lambda-2\in(-1,0)$, monotonicity gives
\[
 \frac{n^{v+1}-1}{v+1}
 \le\sum_{k=1}^{n-1}k^v
 \le1+\frac{(n-1)^{v+1}-1}{v+1}.
\]
Since $|n^{v+1}-(n-1)^{v+1}|\le v+1$ for $n\ge2$, these bounds imply
\[
 n\mathbb E\mathcal M_{n,\lambda}
 =\frac1{\beta(\lambda-1)}+O_\lambda(n^{1-\lambda})
 \longrightarrow\frac1{\beta(\lambda-1)}.
\]

For fixed $p\ge2$, the centered variables $U_k=h_k^+(k/n)^\lambda(A_{n,k}-\mathbb EA_{n,k})$ satisfy $\|U_k\|_p\le C_pn^{-\lambda}k^{\lambda-2}$.
Lemma~\ref{B:lem:independent-sum-moments} therefore gives
\begin{equation}
 \begin{split}
 \|\mathcal M_{n,\lambda}-\mathbb E\mathcal M_{n,\lambda}\|_p
 &\le C_p n^{-\lambda}
                  \left(\sum_{k<n}k^{2\lambda-4}\right)^{1/2}\\
 &\le C_{p,\lambda}
 \begin{cases}
 n^{-\lambda},&1<\lambda<3/2,\\
 n^{-3/2}\sqrt{\ell_n},&\lambda=3/2,\\
 n^{-3/2},&3/2<\lambda<2.
 \end{cases}
 \end{split}
 \label{B:eq:laguerre-local-fluctuation}
\end{equation}
The three summations used here are
\[
 \sum_{k<n}k^{2\lambda-4}\le
 \begin{cases}
  1+(3-2\lambda)^{-1},&1<\lambda<3/2,\\
  1+\log n,&\lambda=3/2,\\
  1+n^{2\lambda-3}/(2\lambda-3),&3/2<\lambda<2.
 \end{cases}
\]
Choose $0<\delta<\min\{\lambda-1,1/2\}$ and a fixed $p\ge2$ with $p\delta>1$. For $\lambda<3/2$, $n^{-\lambda}\le n^{-1-\delta}$; for $\lambda>3/2$, $n^{-3/2}\le n^{-1-\delta}$. At
$\lambda=3/2$, $\sup_{n\ge2}n^{-1/2+\delta}\sqrt{\ell_n}<\infty$ by the comparison of logarithmic and polynomial powers in the preceding proof. Hence
\[
 \begin{split}
 \sum_{n\ge n_0}\mathbb P\left(
 n|\mathcal M_{n,\lambda}-\mathbb E\mathcal M_{n,\lambda}|>\epsilon\right)
 &\le\epsilon^{-p}\sum_{n\ge n_0}n^p
       \|\mathcal M_{n,\lambda}-\mathbb E\mathcal M_{n,\lambda}\|_p^p\\
 &\le C_{p,\lambda,\delta}\epsilon^{-p}
        \sum_{n\ge n_0}n^{-p\delta}<\infty.
 \end{split}
\]
Intersect the Borel--Cantelli events for $\epsilon=1/l$, $l\ge1$. On this event the centered sum multiplied by $n$ tends to zero. Adding its deterministic mean proves
\eqref{B:eq:laguerre-local-energy-limit}. The argument uses independence only within each fixed row.
\end{proof}

\begin{lemma}\label{B:lem:laguerre-local-maximum}
Almost surely,
\begin{equation}
 \max_{m_n<k\le n}|L_{n,k}^X|
 +\max_{m_n<k<n}|L_{n,k}^Y|
 =O\!\left(\sqrt{\frac{\ell_n}{m_n+1}}\right)=o(1).
 \label{B:eq:laguerre-local-maximum}
\end{equation}
\end{lemma}
\begin{proof}
Fix $D>0$ and set $t_{n,k}=D\sqrt{\ell_n/k}$. The scalar mean bound applies because the coefficient marginals are preserved, and gives
\[
 \sup_{k>m_n}t_{n,k}\le D\sqrt{\ell_n/(m_n+1)}\longrightarrow0,
 \qquad
 \frac{|\mathbb EL_{n,k}^\sigma|}{t_{n,k}}
 \le\frac{C}{D\sqrt{k\ell_n}}
 \le\frac{C}{D\sqrt{(m_n+1)\ell_n}}\longrightarrow0.
\]
These inequalities apply to $X$ through $k=n$ and to $Y$ through $k=n-1$. In particular $\sqrt{n\beta}X_{n,n}^{\mathrm L}$ has law $\chi_{\beta(n+a)}$, so the last coefficient has
the same mean and tail estimates. For large $n$, $t_{n,k}\le1$ and $|\mathbb EL_{n,k}^\sigma|\le t_{n,k}/2$ uniformly. Thus
\[
 \begin{split}
 \mathbb P(|L_{n,k}^\sigma|>t_{n,k})
 &\le\mathbb P\left(
 |L_{n,k}^\sigma-\mathbb EL_{n,k}^\sigma|>t_{n,k}/2\right)\\
 &\le2\exp(-ckt_{n,k}^2/4)=2(e+n)^{-cD^2/4}.
 \end{split}
\]
There are fewer than $2n$ coordinates in the two ranges. Choose $D$ with $cD^2/4>3$; the sum of their union probabilities is at most $4\sum_n n(e+n)^{-cD^2/4}<\infty$.
Borel--Cantelli gives one event of probability one and a finite $n_0(\omega)$ such that both maxima are bounded by $D\sqrt{\ell_n/(m_n+1)}$ for $n\ge n_0(\omega)$. Their sum is
bounded by twice this quantity, which tends to zero since $A>12$. The coordinatewise bound obtained before taking maxima also yields $|L_{n,n}^X|^2\le D^2\ell_n/n$. Hence a last
diagonal term added to the weighted logarithmic energy would satisfy
\[
 n h_n^+|L_{n,n}^X|^2\le D^2\ell_n/n\longrightarrow0.
\]
There is no corresponding $L_{n,n}^Y$. In the inverse calculation the last diagonal is retained through $q_{n,n}^{\mathrm L}$, regardless of the empty logarithmic product
$P_{n,n}=1$.
\end{proof}

\subsubsection{The errors on the selected cells}
\label{B:subsection:laguerre-retained-errors}

The estimates just obtained can now be expressed in the two error variables of the stability criterion. The ratio sums determine $\delta_{n,k}$, while the scale error combines the
diagonal logarithm with the quadrature error $\rho_{n,k}^B$.

Use $h_k^-=\log(k/(k-1))$ for $k\ge2$, and the integrals $Q_{n,k}^B$ already defined in \eqref{B:eq:2-32}. The quadrature error $\rho_{n,k}^B$ is the one in
\eqref{B:eq:brownian-logarithmic-errors}: $Q_{n,k}^B=n^{-1}q^B(k/n)e^{-\rho_{n,k}^B}$. The finite coefficients $q_{n,k}^{\mathrm L}$ are those of \eqref{B:eq:2-9}. In the scalar
normalization $c_n=0$, the second error from \eqref{B:eq:2-35} takes the explicit form
\begin{equation}
 \rho_{n,k}=-2L_{n,k}^X+\rho_{n,k}^B,
                  \qquad k\in\mathscr I_n^{\mathrm{ret}}.
 \label{B:eq:laguerre-scale-error}
\end{equation}

\begin{proposition}\label{B:prop:laguerre-retained-errors}
For every fixed $1<\lambda<\min\{a+1,2\}$, the constructed arrays satisfy
\begin{equation}
 \sum_{k\in\mathscr I_n^{\mathrm{ret}}}h_k^-(k/n)^\lambda
             (|\delta_{n,k}|^2+|\rho_{n,k}|^2)=O_\omega(n^{-1}),
 \qquad
 \max_{k\in\mathscr I_n^{\mathrm{ret}}}
             (|\delta_{n,k}|+|\rho_{n,k}|)\longrightarrow0
 \label{B:eq:laguerre-defect-bounds}
\end{equation}
almost surely, with the exact identities in \eqref{B:eq:2-35} and $c_n=0$.
\end{proposition}
\begin{proof}
The finite coefficient definition and the identity for the kernel coefficients $[p^B(z)]^2q^B(z)=z^{-1}$ give
\begin{equation}
 P_{n,k}^2q_{n,k}^{\mathrm L}=\frac1{n(X_{n,k}^{\mathrm L})^2},
 \qquad
 [p^B(k/n)]^2 Q_{n,k}^B=k^{-1}e^{-\rho_{n,k}^B}.
 \label{B:eq:laguerre-scale-cancellation}
\end{equation}
The second equality uses $Q_{n,k}^B=n^{-1}q^B(k/n)e^{-\rho_{n,k}^B}$. The definitions of the logarithmic errors also give $P_{n,k}=p^B(k/n)e^{\delta_{n,k}}$ and
$n(X_{n,k}^{\mathrm L})^2=k e^{2L_{n,k}^X}$. Consequently
\[
 \frac{P_{n,k}^2q_{n,k}^{\mathrm L}}
 {[p^B(k/n)]^2Q_{n,k}^B}
 =e^{-2L_{n,k}^X+\rho_{n,k}^B},
 \qquad
 \log\frac{q_{n,k}^{\mathrm L}}{Q_{n,k}^B}
 =-2\delta_{n,k}-2L_{n,k}^X+\rho_{n,k}^B.
\]
Thus
\begin{equation}
 P_{n,k}=p^B(k/n)e^{\delta_{n,k}},\qquad
 q_{n,k}^{\mathrm L}=Q_{n,k}^B
                     e^{-2\delta_{n,k}-2L_{n,k}^X+\rho_{n,k}^B}.
 \label{B:eq:laguerre-exact-errors}
\end{equation}
This proves \eqref{B:eq:laguerre-scale-error} and verifies both identities in \eqref{B:eq:2-35} with $c_n=0$.

For the cell integral in the Brownian comparison, set $z=k/n$, $h=h_k^-$ and $u=u_{n,k}$. The substitution $y=ze^{-t}$ in its defining integral gives
\[
 e^{-\rho_{n,k}^B}
 =\frac{Q_{n,k}^B}{n^{-1}q^B(z)}
 =\frac1{1-e^{-h}}\int_0^h
       \exp\!\left(at-2\beta^{-1/2}[W(u+t)-W(u)]\right)\,\dd t.
\]
Here $1-e^{-h}=1/k$. The increment on this cell is
\[
 W(u+t)-W(u)=\frac{t}{\sqrt{h_k^-}}G_{n,k-1}
                  +\mathbb B_{n,k-1}^+(t),\qquad 0\le t\le h_k^-.
\]
Indeed $h_k^-=h_{k-1}^+$ and the bridge in the proof of the coefficient distribution has boundary points $u_{n,k}$ and $u_{n,k-1}$. The ratio $Z_{n,k}$, in contrast, uses $G_{n,k}$
on the following interval. At the first selected index $k=m_n+1$ the cell integral still uses $G_{n,m_n}$, although the coefficient pair with index $m_n$ belongs to the common
initial sequence. All these Brownian observations exist for the original row. Put $M=(a+1)h+2\beta^{-1/2}\operatorname{osc}_{[u,u+h]}W$. For $0\le t\le h$,
$|(a+1)t-2\beta^{-1/2}(W(u+t)-W(u))|\le M$, and therefore
\[
 e^{-M}\le
 \int_0^h e^{(a+1)t-2\beta^{-1/2}[W(u+t)-W(u)]}
                         \frac{e^{-t}\,\dd t}{1-e^{-h}}
 =e^{-\rho_{n,k}^B}\le e^M.
\]
Taking logarithms proves $|\rho_{n,k}^B|\le(a+1)h+2\beta^{-1/2}\operatorname{osc}_{[u,u+h]}W$.

The two logarithmic lengths satisfy
\begin{equation}
 h_k^+\le h_k^-\le3h_k^+.
 \label{B:eq:laguerre-cell-lengths}
\end{equation}
Indeed, $k/(k-1)\ge(k+1)/k$ gives the lower inequality, and
\[
 \frac{h_k^-}{h_k^+}
 \le\frac{(k-1)^{-1}}{(k+1)^{-1}}
 =\frac{k+1}{k-1}\le3 \qquad(k\ge2).
\]
Every selected index lies between $m_n+1$ and $M_n<n$, so all its terms occur in the two sums already estimated. In particular,
\[
 \sum_{\mathscr I_n^{\mathrm{ret}}}h_k^-(k/n)^\lambda|\delta_{n,k}|^2
 \le3\mathcal P_{n,\lambda},
 \qquad
 \sum_{\mathscr I_n^{\mathrm{ret}}}h_k^-(k/n)^\lambda|L_{n,k}^X|^2
 \le3\mathcal M_{n,\lambda}.
\]
Using $|-2x+y|^2\le8|x|^2+2|y|^2$, the complete weighted error is therefore at most
\[
 3\mathcal P_{n,\lambda}+24\mathcal M_{n,\lambda}
 +2\sum_{\mathscr I_n^{\mathrm{ret}}}
           h_k^-(k/n)^\lambda|\rho_{n,k}^B|^2
 =O_\omega(n^{-1}).
\]
Here the first term is $o_\omega(n^{-1})$, the second has the limit in \eqref{B:eq:laguerre-local-energy-limit}, and the third is bounded by
Proposition~\ref{B:prop:brownian-logarithmic-errors}.

For the maximum, $h_k^-\le1/(k-1)\le1/m_n$ and the cellwise Brownian logarithm estimate give
\[
 \begin{split}
 \max_{\mathscr I_n^{\mathrm{ret}}}|\rho_{n,k}^B|
 &\le\frac{a+1}{m_n}
 +2\beta^{-1/2}\max_{\mathscr I_n^{\mathrm{ret}}}
             \operatorname{osc}_{[u_{n,k},u_{n,k-1}]}W\\
 &=O_\omega\bigl(m_n^{-1}+\sqrt{\ell_n/m_n}\bigr)\longrightarrow0.
 \end{split}
\]
Finally,
\[
 \max_{\mathscr I_n^{\mathrm{ret}}}(|\delta_{n,k}|+|\rho_{n,k}|)
 \le\max_{j>m_n}|\delta_{n,j}|
   +2\max_{k>m_n}|L_{n,k}^X|
   +\max_{\mathscr I_n^{\mathrm{ret}}}|\rho_{n,k}^B|\longrightarrow0.
\]
Intersect the events on which the bounds for coefficient products, logarithms of individual coefficients, and Brownian oscillations hold for this $\lambda$. All identities and both
estimates hold on that intersection. The cell integrals occur in the Brownian comparison; the coefficient vector and the constant $C_n=1$ retain their $\mathcal N_n\vee\mathcal U$
measurability.
\end{proof}

\subsubsection{Bounds for the triangular factor}
\label{B:subsection:laguerre-full-factor-bounds}

To remove the restriction to the selected interval, we need a bound on the entire triangular factor. The exact product formula supplies it, including the diagonal cells and the
unpaired last coefficient.

Use the operators $\mathcal K_n^{\mathrm L},\mathcal R_n^{\mathrm L}$ from \eqref{B:eq:2-12} with $\alpha_n=1$ and $\widetilde{\mathscr J}_n=\mathscr E_n$. The formula
\eqref{B:eq:2-6} for the entries of the inverse matrix and the embedding identity \eqref{B:eq:2-1} give the kernel
\begin{equation}
 K_n^{\mathrm L}(i,j)=
 \frac{P_{n,i}}{X_{n,j}^{\mathrm L}P_{n,j}}\mathbf1_{\{i\le j\}},
 \qquad z_{n,k}=k/n.
 \label{B:eq:laguerre-factor-kernel}
\end{equation}
In fact the factor $n^{-1}$ in $(nB_n^{\mathrm L})^{-1}$ cancels the factor $(\Delta_{n,i}\Delta_{n,j})^{-1/2}=n$ in the kernel. The diagonal cells, including $j=n$, are included.

\begin{proposition}\label{B:prop:laguerre-factor-bound}
Almost surely, for every $\eta>0$,
\begin{equation}
 |K_n^{\mathrm L}(i,j)|\le C_{\omega,\eta}
 z_{n,j}^{-1/2}(z_{n,i}/z_{n,j})^{a/2}(z_{n,i}z_{n,j})^{-\eta},
 \qquad 1\le i\le j\le n,
 \label{B:eq:laguerre-factor-bound}
\end{equation}
for all sufficiently large $n$.
\end{proposition}
\begin{proof}
Write $m=m_n$ and $z_k=k/n$. On the event of Lemma~\ref{B:lem:laguerre-local-maximum}, choose $n_0(\omega)$ so that $\max_{k>m}|L_{n,k}^X|\le1$ for every $n\ge n_0(\omega)$. Hence
$(X_{n,k}^{\mathrm L})^{-1}=z_k^{-1/2}e^{-L_{n,k}^X} \le e z_k^{-1/2}$ for $k>m$, including $k=n$. For $k\le m$, the bound for the common initial sequence gives
\[
 (X_{n,k}^{\mathrm L})^{-1}
 =\frac{\sqrt{n\beta}}{\widehat{\mathsf X}_k}
 \le\frac{\sqrt\beta}{c_\omega}\sqrt{n/k}.
\]
These two ranges prove
\begin{equation}
 (X_{n,k}^{\mathrm L})^{-1}\le C_\omega z_k^{-1/2}
             \quad(1\le k\le n).
 \label{B:eq:laguerre-diagonal-lower-bound}
\end{equation}

Fix $\eta>0$. If $m<i\le j\le n$, the exact product identity is
\begin{equation}
 \frac{P_{n,i}}{P_{n,j}}=(i/j)^{a/2}
  \exp\{\mathcal Y(z_i)-\mathcal Y(z_j)+\delta_{n,i}-\delta_{n,j}\}.
 \label{B:eq:laguerre-regular-products}
\end{equation}
For $\epsilon>0$, the sublinear growth estimate for the Brownian path gives
\[
 |\mathcal Y(z_i)|+|\mathcal Y(z_j)|
 \le 2C_{\omega,\epsilon}
       +\epsilon\log(z_i^{-1}z_j^{-1}).
\]
Increase $n_0(\omega)$ so that $\max_{k>m}|\delta_{n,k}|\le1$ for the same values of $n$. It follows that
\[
 \frac{P_{n,i}}{P_{n,j}}
 \le e^{2+2C_{\omega,\epsilon}}(i/j)^{a/2}(z_iz_j)^{-\epsilon}.
\]
Take $\epsilon\le\eta$ and use $0<z_iz_j\le1$ to obtain the required ratio bound.

If $1\le i\le j\le m$, all coefficients in the ratio come from the common initial sequence. For any fixed $\epsilon>0$,
\[
 \frac{P_{n,i}}{P_{n,j}}
 \le C_{\omega,\epsilon}j^\epsilon(i/j)^{a/2},
 \qquad
 j^\epsilon(z_iz_j)^\eta
 \le n^{-2\eta}m^{\epsilon+2\eta}
 \le C n^{-2\eta}(\log n)^{A(\epsilon+2\eta)}\longrightarrow0.
\]
The limit follows from the comparison of logarithmic and polynomial powers proved above. Thus $j^\epsilon\le(z_iz_j)^{-\eta}$ for all indices in this range once $n$ is large
enough.

For $i\le m<j$, split the product after index $m$:
\begin{align}
 \frac{P_{n,i}}{P_{n,j}}
 &=\left(\prod_{k=i}^{m}
          \frac{\widehat{\mathsf Y}_k}{\widehat{\mathsf X}_k}\right)
                         \frac{P_{n,m+1}}{P_{n,j}}\notag\\
 &\le C_{\omega,\epsilon}
 (i/(m+1))^{a/2}(m+1)^\epsilon
 ((m+1)/j)^{a/2}(z_{m+1}z_j)^{-\epsilon}\notag\\
 &=C_{\omega,\epsilon}(i/j)^{a/2}n^{2\epsilon}j^{-\epsilon}.
 \label{B:eq:laguerre-products-across-initial-interval}
\end{align}
The first estimate applies the lemma for the initial products with upper index $m+1$; the second is the bound already obtained for reconstructed indices. The factors
$(m+1)^\epsilon$ and $z_{m+1}^{-\epsilon}$ multiply to $n^\epsilon$. Set $\epsilon=\eta/4$. Since $i\le m$, $j\le n$, and $\eta-\epsilon>0$,
\[
 n^{2\epsilon}j^{-\epsilon}(z_iz_j)^\eta
 =n^{2\epsilon-2\eta}i^\eta j^{\eta-\epsilon}
 \le n^{\epsilon-\eta}m^\eta
 \le Cn^{-3\eta/4}(\log n)^{A\eta}\longrightarrow0.
\]
This proves the ratio bound in the third range, including $j=m+1$. For $i=j$, the ratio is the empty product and equals one. Multiplying each ratio bound by
\eqref{B:eq:laguerre-diagonal-lower-bound} proves \eqref{B:eq:laguerre-factor-bound}.

For each rational $d>0$, intersect the events for the Brownian bounds and the initial coefficient products, and also include the two maximum events. This is a countable
intersection of events of probability one. For any $\eta>0$ choose rational $0<d\le\eta$. Then $(z_iz_j)^{-d}\le(z_iz_j)^{-\eta}$, so the bound at $d$ implies the bound at $\eta$.
If $N_d(\omega)$ and $C_d(\omega)$ are the starting index and constant at $d$, define
\[
 C_d^*(\omega)=\max\left\{C_d(\omega),
 \max_{1\le n<N_d(\omega)}\max_{1\le i\le j\le n}
 \frac{|K_n^{\mathrm L}(i,j)|}
 {z_j^{-1/2}(z_i/z_j)^{a/2}(z_iz_j)^{-d}}\right\}.
\]
Every denominator is positive. Every numerator is finite because all chi coefficients are positive and finite on a common countable intersection of events of probability one. Thus
$C_d^*<\infty$, and it bounds every size. The finite maxima are measurable.
\end{proof}

\subsubsection{Bounds for the inverse and its restrictions}
\label{B:subsection:laguerre-inverse-truncation}

The full factor bound supplies the hypotheses of the deterministic truncation estimates for the coupled rows. Those estimates bound the omitted kernel and the rank-one correction
separately, giving a comparison of the full positive inverses.

Let $\nu_a=\min\{a,1\}$. For the fixed cuts in \eqref{B:eq:laguerre-cuts}, the projections already used in Section~\ref{B:subsection:the-criterion-for-convergence-of-the-inverses}
are $P_n^-=\mathbf1_{(0,m_n/n]}$, $P_n^0=\mathbf1_{(m_n/n,M_n/n]}$, and $P_n^+=\mathbf1_{(M_n/n,1]}$. Set $d_n=m_n/n$ and $\varepsilon_n=1-M_n/n=r_n/n$. Use the restricted factor
and its correction in the form
\begin{equation}
 \mathcal R_n^{\mathrm L,0}
    =(P_n^0\mathcal K_n^{\mathrm L}P_n^0)^*
                       (P_n^0\mathcal K_n^{\mathrm L}P_n^0),\qquad
 \mathcal T_n^{\mathrm L}
    =P_n^0(\mathcal K_n^{\mathrm L})^*P_n^+\mathcal K_n^{\mathrm L}P_n^0.
 \label{B:eq:laguerre-cut-operators}
\end{equation}
These are the operators in the criterion, specialized to this model.

\begin{lemma}\label{B:lem:laguerre-inverse-kernel-bounds}
Almost surely, $\sup_n\|\mathcal R_n^{\mathrm L}\|_{\mathfrak S_2}<\infty$. For every $0<\eta<\nu_a$, there is a finite $n_0(\omega,\eta)$ such that, for $n\ge n_0(\omega,\eta)$,
\begin{equation}
 |R_n^{\mathrm L}(i,j)|\le C_{\omega,\eta}
       (z_{n,i\wedge j}/z_{n,i\vee j})^{a/2}
       (z_{n,i}z_{n,j})^{-\eta}.
 \label{B:eq:laguerre-full-inverse-bound}
\end{equation}
For $i\le j$ and $j\ge n/2$, $|R_n^{\mathrm L}(i,j)|\le C_{\omega,\eta}z_{n,i}^{a/2-\eta}(1-z_{n,j-1})$.
\end{lemma}
\begin{proof}
Use the atomic measure $\mathfrak m_n=n^{-1}\sum_{k=1}^n\delta_{k/n}$. Its cell masses are all $1/n$, its cumulative mass through $z_k$ is $k/n$, and its inclusive mass above $z_j$
is $(n-j+1)/n$. Take $p_i=P_{n,i}$ and $r_j=(X_{n,j}^{\mathrm L}P_{n,j})^{-1}$. For each sufficiently small $\xi>0$, the finite factor bound reads
\[
 |p_ir_j|\le C_{\omega,\xi}z_{n,j}^{-1/2}
                (z_{n,i}/z_{n,j})^{a/2}(z_{n,i}z_{n,j})^{-\xi},
 \qquad i\le j.
\]
This is \eqref{B:eq:triangular-factor-bound} with constants uniform in $n$. The positive product has kernel
\[
 R_n^{\mathrm L}(i,j)
 =\frac1n\sum_{l\ge i\vee j}K_n^{\mathrm L}(i,l)K_n^{\mathrm L}(j,l)
 =P_{n,i}P_{n,j}\sum_{l\ge i\vee j}q_{n,l}^{\mathrm L},
 \qquad q_{n,l}^{\mathrm L}=\frac1{n(X_{n,l}^{\mathrm L}P_{n,l})^2}.
\]
The summation includes $l=i\vee j$ and $l=n$. Thus it is the full embedded inverse from \eqref{B:eq:2-16}.

Fix $0<\eta<\nu_a$ and choose $0<\epsilon<\eta/8$. Apply Proposition~\ref{B:prop:full-factor-and-green-kernel} with the coefficient exponent $\xi$ required in its proof for this
$\epsilon$. The preceding factor bound holds for every such $\xi>0$, and gives
\[
 |R_n^{\mathrm L}(i,j)|
 \le C_{\omega,\epsilon}
 (z_{n,i\wedge j}/z_{n,i\vee j})^{a/2}
 (z_{n,i}z_{n,j})^{-\epsilon}.
\]
Because $\epsilon<\eta$, this proves \eqref{B:eq:laguerre-full-inverse-bound}. With this same $\epsilon$, Lemma~\ref{B:lem:green-kernel-row-bound} gives
\[
 \frac1n\sum_{j=1}^n|R_n^{\mathrm L}(i,j)|^2
 \le C_{\omega,\eta}z_{n,i}^{\nu_a-\eta}.
\]
The required inequality between its exponents is $8\epsilon<\eta$, which holds by their choice. Since $\nu_a-\eta>0$ and
$z_{n,i}\le1$,
\[
 \|\mathcal R_n^{\mathrm L}\|_{\mathfrak S_2}^2
 =\frac1n\sum_i\frac1n\sum_j|R_n^{\mathrm L}(i,j)|^2
 \le\frac{C_{\omega,\eta}}n\sum_i z_{n,i}^{\nu_a-\eta}
 \le C_{\omega,\eta}.
\]
All earlier sizes have finite inverse norms and are included in the supremum by increasing the random constant.

For the estimate near one, choose $0<\xi<\min\{\eta,a/4,1/8\}$, so $a-2\xi>0$. If $i\le j$ and $j\ge n/2$, multiplication of the two factor bounds inside the sum gives
\[
 \begin{split}
 |R_n^{\mathrm L}(i,j)|
 &\le C_{\omega,\xi}z_{n,i}^{a/2-\xi}z_{n,j}^{a/2-\xi}
       \frac1n\sum_{l=j}^n z_{n,l}^{-a-1-2\xi}\\
 &\le C_{\omega,\xi}z_{n,i}^{a/2-\xi}\frac{n-j+1}{n}
 \le C_{\omega,\eta}z_{n,i}^{a/2-\eta}(1-z_{n,j-1}).
 \end{split}
\]
In the second line all $z_{n,l}$ and $z_{n,j}$ lie in $[1/2,1]$, so their powers have fixed bounds. At $j=n$ the sum contains exactly one term, and $1-z_{n,n-1}=1/n$. For the next
estimate use the bound with exponent $\xi$; since $0<\xi<\eta$ and $0<z\le1$, one has $z^{a/2-\xi}\le z^{a/2-\eta}$. For $i>j\ge n/2$, symmetry of the positive inverse and the
already proved estimate with indices $(j,i)$ give
\[
 |R_n^{\mathrm L}(i,j)|=|R_n^{\mathrm L}(j,i)|
 \le C_{\omega,\xi}z_{n,j}^{a/2-\xi}(1-z_{n,i-1})
 \le C_{\omega,\xi}(1-z_{n,j-1}).
\]
Here $a/2-\xi>0$, $z_{n,j}\le1$, and $z_{n,i-1}\ge z_{n,j-1}$. Consequently the entire row satisfies
\[
 \frac1n\sum_i|R_n^{\mathrm L}(i,j)|^2
 \le C_{\omega,\xi}(1-z_{n,j-1})^2
       \left(\frac1n\sum_{i\le j}z_{n,i}^{a-2\xi}
                          +\frac{n-j}{n}\right)
 \le C_{\omega,\xi}(1-z_{n,j-1})^2.
\]
This includes the full final row and its diagonal cell.
\end{proof}

\begin{proposition}\label{B:prop:laguerre-truncation-error}
For $0<\eta<\nu_a$, almost surely,
\begin{equation}
 \begin{split}
 \|P_n^-\mathcal R_n^{\mathrm L}\|_{\mathfrak S_2}^2
   +\|\mathcal R_n^{\mathrm L}P_n^-\|_{\mathfrak S_2}^2
       &\le C_{\omega,\eta}d_n^{1+\nu_a-\eta},\\
 \|P_n^+\mathcal R_n^{\mathrm L}\|_{\mathfrak S_2}^2
   +\|\mathcal R_n^{\mathrm L}P_n^+\|_{\mathfrak S_2}^2
       &\le C_{\omega,\eta}\varepsilon_n^3,\\
 \|\mathcal T_n^{\mathrm L}\|_{\mathfrak S_2}^2
       &\le C_{\omega,\eta}\varepsilon_n^2.
 \end{split}
 \label{B:eq:laguerre-truncation-bounds}
\end{equation}
In particular, $\|\mathcal R_n^{\mathrm L}-\mathcal R_n^{\mathrm L,0}\|_{\mathfrak S_2}^2 =o_\omega(n^{-1})$.
\end{proposition}
\begin{proof}
Work on the event of the preceding factor and inverse bounds, and choose $0<\eta<\nu_a$. The cuts are sampling points, and $M_n/n\ge1/2$ for all sufficiently large $n$. Symmetry
and the row estimate near zero give
\[
 \begin{split}
 \|P_n^-\mathcal R_n^{\mathrm L}\|_{\mathfrak S_2}^2
 +\|\mathcal R_n^{\mathrm L}P_n^-\|_{\mathfrak S_2}^2
 &\le\frac{2C_{\omega,\eta}}n
           \sum_{i\le m_n}(i/n)^{\nu_a-\eta}\\
 &\le2C_{\omega,\eta}(m_n/n)^{1+\nu_a-\eta}.
 \end{split}
\]
The last inequality uses $i/n\le m_n/n$ and $\nu_a-\eta>0$. The bound for rows near one in the preceding proof gives
\[
 \begin{split}
 \|P_n^+\mathcal R_n^{\mathrm L}\|_{\mathfrak S_2}^2
 +\|\mathcal R_n^{\mathrm L}P_n^+\|_{\mathfrak S_2}^2
 &\le\frac{2C_\omega}{n^3}\sum_{l=1}^{r_n}l^2\\
 &=\frac{C_\omega r_n(r_n+1)(2r_n+1)}{3n^3}
 \le2C_\omega\varepsilon_n^3.
 \end{split}
\]
Here $l=n-j+1$; thus the index $l=1$ includes the last row.

Write $P=P_n^0$, $K=\mathcal K_n^{\mathrm L}$, and $u_n=\sum_{i=m_n+1}^{M_n}P_{n,i}\mathbf1_{I_{n,i}}$. Integration on the cells in the second kernel variable gives
\[
 (P_n^+KPf)(t)
 =\frac1{X_{n,l}^{\mathrm L}P_{n,l}}
       \sum_{i=m_n+1}^{M_n}P_{n,i}\int_{I_{n,i}}f(s)\,\dd s,
 \qquad t\in I_{n,l},\quad l>M_n.
\]
Its positive product is the rank-one operator with vector $u_n$ and coefficient $\sum_{l>M_n}q_{n,l}^{\mathrm L}$. Hence
\begin{equation}
 \begin{split}
 T_n^{\mathrm L}(i,j)&=P_{n,i}P_{n,j}
                     \sum_{l=M_n+1}^n q_{n,l}^{\mathrm L},
                      \qquad m_n<i,j\le M_n,\\
 \|\mathcal T_n^{\mathrm L}\|_{\mathfrak S_2}^2
  &=\left(\sum_{l=M_n+1}^n q_{n,l}^{\mathrm L}\right)^2
       \left(\frac1n\sum_{i=m_n+1}^{M_n}P_{n,i}^2\right)^2.
 \end{split}
 \label{B:eq:laguerre-rank-one-tail}
\end{equation}
This formula includes the case where the retained vector or the omitted sum vanishes. For a uniform estimate, take $0<\xi<\min\{a/4,1/8\}$. For $l>M_n$ the factor bounds give
\[
 |K_n^{\mathrm L}(i,l)K_n^{\mathrm L}(j,l)|
 \le C_{\omega,\xi}z_{n,i}^{a/2-\xi}z_{n,j}^{a/2-\xi},
 \qquad m_n<i,j\le M_n,
\]
since $z_{n,l}\ge1/2$. Summing the $r_n$ cells in the second kernel variable with their mass $1/n$ therefore yields
\[
 \begin{split}
 |T_n^{\mathrm L}(i,j)|
 &\le C_{\omega,\xi}\varepsilon_n
             z_{n,i}^{a/2-\xi}z_{n,j}^{a/2-\xi},\\
 \|\mathcal T_n^{\mathrm L}\|_{\mathfrak S_2}^2
 &\le C_{\omega,\xi}\varepsilon_n^2
       \left(\frac1n\sum_{i=m_n+1}^{M_n}z_{n,i}^{a-2\xi}\right)^2
 \le C_{\omega,\xi}\varepsilon_n^2.
 \end{split}
\]
The exponent $a-2\xi$ is positive, so the last average is at most one. In particular the last summand in \eqref{B:eq:laguerre-rank-one-tail} is the rank-one operator
$q_{n,n}^{\mathrm L}u_n\otimes u_n$, where $P_{n,n}=1$ and $q_{n,n}^{\mathrm L}=1/[n(X_{n,n}^{\mathrm L})^2]$. The factor estimate at $l=n$ gives
\[
 (q_{n,n}^{\mathrm L})^2\|u_n\|_2^4
 =\frac1{n^2}
   \left(\frac1n\sum_{i=m_n+1}^{M_n}
                     |K_n^{\mathrm L}(i,n)|^2\right)^2
 \le\frac{C_\omega}{n^2}.
\]
Thus $\delta_{n,n}=0$ does not remove the unpaired last diagonal from the inverse. Its contribution is included in the quadratic upper end correction. All three inequalities in
\eqref{B:eq:laguerre-truncation-bounds} follow.

Let $Q=I-P$ and $R=\mathcal R_n^{\mathrm L}$. Triangularity gives $P_n^-KP=0$, and the decomposition $I=P+P_n^-+P_n^+$ gives
$PRP=\mathcal R_n^{\mathrm L,0}+\mathcal T_n^{\mathrm L}$. The four terms in
\[
 R-\mathcal R_n^{\mathrm L,0}
 =QRQ+QRP+PRQ+\mathcal T_n^{\mathrm L}
\]
have pairwise disjoint kernel supports. Since $R=R^*$,
\[
 \begin{split}
 \|R-\mathcal R_n^{\mathrm L,0}\|_{\mathfrak S_2}^2
 &=\|QRQ\|_{\mathfrak S_2}^2+2\|QRP\|_{\mathfrak S_2}^2
                         +\|\mathcal T_n^{\mathrm L}\|_{\mathfrak S_2}^2,\\
 \|QR\|_{\mathfrak S_2}^2+\|RQ\|_{\mathfrak S_2}^2
 &=2\|QRQ\|_{\mathfrak S_2}^2+2\|QRP\|_{\mathfrak S_2}^2.
 \end{split}
\]
Also $Q=P_n^-+P_n^+$ with orthogonal ranges. The second line is therefore the sum of the four row and column norms already bounded. Adding the correction estimate proves
\begin{equation}
 \|\mathcal R_n^{\mathrm L}-\mathcal R_n^{\mathrm L,0}\|_{\mathfrak S_2}^2
 \le C_{\omega,\eta}
       [d_n^{1+\nu_a-\eta}+\varepsilon_n^3+\varepsilon_n^2].
 \label{B:eq:laguerre-truncation-restoration}
\end{equation}
Because $d_n=\varepsilon_n=r_n/n$,
\[
 nd_n^{1+\nu_a-\eta}=r_n^{1+\nu_a-\eta}n^{-(\nu_a-\eta)},
 \qquad n\varepsilon_n^3=r_n^3/n^2,
 \qquad n\varepsilon_n^2=r_n^2/n.
\]
For large $n$, $r_n\le2(\log n)^A$. The three expressions are bounded respectively by constant multiples of
\[
 (\log n)^{A(1+\nu_a-\eta)}n^{-(\nu_a-\eta)},
 \qquad (\log n)^{3A}n^{-2},
 \qquad (\log n)^{2A}n^{-1}.
\]
Each tends to zero by the comparison of logarithmic and polynomial powers, since $\nu_a-\eta>0$. This proves the required $o_\omega(n^{-1})$ bound.
\end{proof}

\subsubsection{Proof of the beta-Laguerre theorem}
\label{B:subsection:proof-of-theorem-1-1}

\begin{proof}[Proof of Theorem~\ref{B:result:1-1}]
Fix $1<\lambda<\min\{a+1,2\}$, which is possible because $a>0$. Use the arrays of Proposition~\ref{B:prop:laguerre-array-law}. In
Proposition~\ref{B:result:common-operator-criterion}, take
\[
 \theta=\operatorname{id},\qquad \alpha_n=1,\qquad
 \mathscr J_n=\widetilde{\mathscr J}_n=\mathscr E_n,
 \qquad \mathcal A_n^{\mathrm L}=n^2B_n^{\mathrm L}(B_n^{\mathrm L})^{\mathsf T}.
\]
Proposition~\ref{B:prop:laguerre-array-law} verifies condition (i): the coefficient distribution is exactly \eqref{B:eq:intro-laguerre-coefficients}, and each matrix is
$\mathcal N_n\vee\mathcal U$-measurable with $\mathcal U$ independent of the whole Brownian path. Proposition~\ref{B:prop:laguerre-retained-errors} verifies condition (ii) with
$C_n=1$ and the logarithmic errors \eqref{B:eq:laguerre-exact-errors}. Proposition~\ref{B:prop:laguerre-truncation-error} verifies condition (iii) for the same selected cells and
the same full inverse. The criterion therefore proves
\[
 0<\liminf_{n\to\infty}
 n\|\mathcal R_n^{\mathrm L}-\mathcal H_{\beta,a}^{-1}\|_{\mathfrak S_2}^2
 \le\sup_{n\ge1}
 n\|\mathcal R_n^{\mathrm L}-\mathcal H_{\beta,a}^{-1}\|_{\mathfrak S_2}^2
 <\infty
\]
almost surely for this sequence, as required.
\end{proof}

The limit in \eqref{B:eq:laguerre-local-energy-limit} concerns $\mathcal M_{n,\lambda}$. To describe what it implies for the operator error, set
\[
 F_n=\mathcal R_n^{\mathrm L,0}-\mathcal R_n^{B,\mathrm{step},0},
 \qquad G_n=\mathcal R_n^{B,\mathrm{step},0}-\mathcal R_n^{B,0}.
\]
The coefficient comparison and the estimate for sampling the path give $\|F_n\|_{\mathfrak S_2}+\|G_n\|_{\mathfrak S_2}=O_\omega(n^{-1/2})$. The two truncation estimates give
$\|\mathcal R_n^{\mathrm L}-\mathcal H_{\beta,a}^{-1}-F_n-G_n\|_{\mathfrak S_2} =o_\omega(n^{-1/2})$. Hence Cauchy--Schwarz yields
\[
 n\|\mathcal R_n^{\mathrm L}-\mathcal H_{\beta,a}^{-1}\|_{\mathfrak S_2}^2
 =n\|F_n\|_{\mathfrak S_2}^2+n\|G_n\|_{\mathfrak S_2}^2
 +2n\operatorname{Re}\langle F_n,G_n\rangle_{\mathfrak S_2}+o_\omega(1).
\]
Determining a leading order coefficient requires the joint asymptotics of these two differences, including their inner product. The limit of $\mathcal M_{n,\lambda}$ alone does not
determine that coefficient.

\subsection{General potentials and convex reference measures}
\label{B:sec:general-potentials-and-convex-reference-measures}

For a general potential, the first few logarithmic powers can prevent uniform convexity of the original coefficient Hamiltonian. The parameters $V$, $\beta>0$ and $a>0$ remain
those of Theorem~\ref{B:result:1-3}. Replacing the finitely many negative powers $\beta(k+a)-1$ and $\beta k-1$ by zero gives a convex reference measure. Estimates must then be
transferred back through the density ratio. Its integrability is determined by the powers at zero, so the admissible moment range is retained in every application of H\"older's
inequality.

On an interval beyond the changed coordinates the two conditional laws agree. Their minimizers are compared with the profile and first correction of Rider and Waters
\cite{rider2019universality} by evaluating the stationary equations. Exponential decay of the inverse Hessian bounds the effect of the fixed exterior coordinates. The final
subsubsections construct exterior vectors with their original joint marginals and give the conditional integral identities used to complete each coefficient row. The Gaussian
approximation developed below concerns the interior conditional law. Its covariance is the inverse Hessian on the whole conditional interval. Transport is applied to the
logarithmic ratios; the remaining coordinates are then sampled from the original density conditional on those ratios. Throughout the construction the joint exterior marginal is
unchanged, including the dependence between its initial and final coordinates. This preserves the finite coefficient law and permits the subsequent comparison with a single
Brownian path.

\subsubsection{The polynomial in the coefficient density}
\label{B:subsection:the-polynomial-hamiltonian}

We first establish locality and convexity of the polynomial part of the density. These properties are unchanged when the first few logarithmic terms are modified.

Use the bidiagonal matrix $B=B(x,y)$ and polynomial $V(\lambda)=\sum_{m=1}^d g_m\lambda^m$ defined in Section~\ref{B:sec:reorganized-1}. A constant in $V$ has no effect on the
normalized coefficient law. Set
\begin{equation}
 \begin{gathered}
 \mathscr S_n=\{(x,k):1\le k\le n\}\cup\{(y,k):1\le k<n\},\qquad
 \iota(\sigma,k)=k,\qquad \kappa=n\beta,\\
 \mathcal P_n(z)=\operatorname{Tr}V(B(z)B(z)^{\mathsf T}),\qquad
 c_{x,k}=\beta(k+a)-1,\quad c_{y,k}=\beta k-1.
 \end{gathered}
 \label{B:eq:general-coefficient-notation}
\end{equation}
Here $z=(x_1,\ldots,x_n;y_1,\ldots,y_{n-1})$. The original probability measure is
\begin{equation}
 \mathbb P_n(\dd z)=Z_n^{-1}\exp[-\kappa\mathcal P_n(z)]
                  \prod_{\alpha\in\mathscr S_n}z_\alpha^{c_\alpha}\,\dd z
                  \quad\text{on }(0,\infty)^{\mathscr S_n}.
 \label{B:eq:original-coefficient-density}
\end{equation}
Its Hamiltonian on the open orthant is $H_n=\mathcal P_n-\kappa^{-1}\sum_\alpha c_\alpha\log z_\alpha$. This is the coefficient density in \cite{rider2019universality}. Choose
$c_*>0$ so that $U(u)=V(u^2)$ satisfies $U''(u)\ge c_*$ for all real $u$.

The normalization below is valid for $a>-1$. $V(x^2)\ge c_*x^2/2$ by evenness and \eqref{B:eq:1-8}. Applied to the singular values of $B$, this gives
\[
 \operatorname{Tr}V(BB^{\mathsf T})
 \ge\frac{c_*}{2}\operatorname{Tr}BB^{\mathsf T}
 =\frac{c_*}{2}\left(\sum_kx_k^2+\sum_ky_k^2\right).
\]
Every exponent in the density exceeds $-1$, and
\[
 \int_0^\infty t^{\nu-1}e^{-ct^2/2}\,\dd t
 =2^{\nu/2-1}c^{-\nu/2}\Gamma(\nu/2)<\infty
 \quad(c,\nu>0).
\]
The product of these integrals bounds the normalizer from above; integration over $[1,2]^{2n-1}$ makes it strictly positive.

When $V(\lambda)=\lambda/2$, the exponential in \eqref{B:eq:1-10} splits into its individual coordinate factors. Normalizing a factor of degree $\nu$ gives
\[
 \frac{(n\beta)^{\nu/2}}{2^{\nu/2-1}\Gamma(\nu/2)}
 t^{\nu-1}e^{-n\beta t^2/2}\,\dd t,
\]
which is the density of $(n\beta)^{-1/2}\chi_\nu$. Thus \eqref{B:eq:1-10} gives exactly \eqref{B:eq:intro-laguerre-coefficients} in the linear case.

\begin{lemma}\label{B:lem:coefficient-polynomial-convexity}
The polynomial $\mathcal P_n$ is separately even in every coordinate, and $\nabla^2\mathcal P_n\ge c_*I$. There is an integer $R_*$ depending only on $d$ such that every derivative
involving a coordinate at site $k$ depends only on sites in $[k-R_*,k+R_*]$. On each fixed coefficient box, its derivatives of fixed order have bounds independent of $n$.
\end{lemma}
\begin{proof}
For a real matrix $C$ put $\mathcal A(C)=\left(\begin{smallmatrix}0&C\\C^{\mathsf T}&0\end{smallmatrix}\right)$. The
eigenvalues of this symmetric matrix are the positive and negative singular values of $C$. Consequently
\[
 \mathcal P_n(z)=\tfrac12\operatorname{Tr}U(\mathcal A(B(z))),\qquad
 \|\mathcal A(C)\|_{\rm HS}^2=2\|C\|_{\rm HS}^2.
\]
For a monomial and a real symmetric $A$, the product rule gives
\[
 D(A^m)[E]=\sum_{j=0}^{m-1}A^jEA^{m-1-j},\qquad
 D\operatorname{Tr}U(A)[E]=\operatorname{Tr}(U'(A)E).
\]
In an eigenbasis of $A$, the coefficient of $E_{ij}$ in the first sum is $\sum_{l=0}^{m-1}\lambda_i^l\lambda_j^{m-1-l}$. For unequal eigenvalues summing the finite geometric series
gives $(\lambda_i^m-\lambda_j^m)/(\lambda_i-\lambda_j)$; at equal eigenvalues it gives $m\lambda_i^{m-1}$. Linear combination over the monomials proves
$D F(A)[E]_{ij}=F^{[1]}(\lambda_i,\lambda_j)E_{ij}$ for a polynomial $F$. For unequal arguments, $F^{[1]}(s,t)=(F(s)-F(t))/(s-t)$, and for equal arguments it equals $F'(s)$.
Differentiating $\operatorname{Tr}U(A)$ twice therefore yields
\[
 D^2\operatorname{Tr}U(A)[E,E]
 =\sum_{i,j}(U')^{[1]}(\lambda_i,\lambda_j)|E_{ij}|^2
 \ge c_*\|E\|_{\rm HS}^2.
\]
Indeed $(U')^{[1]}(s,t)=\int_0^1U''((1-v)s+vt)\,\dd v\ge c_*$. Applying this to $\mathcal A(B)$ proves the Hessian bound because $z\mapsto B(z)$ preserves the Euclidean norm of the
coefficients.

The nonzero positions of $B$ form the path with successive vertices $\mathrm R_1,\mathrm C_1,\mathrm R_2,\mathrm C_2,\ldots,\mathrm R_n,\mathrm C_n$. Remove the edge corresponding
to a specified coefficient. Give the vertices in one component sign $-1$ and those in the other sign $1$. The resulting diagonal matrices $D_{\mathrm R},D_{\mathrm C}$ change only
that coefficient in $D_{\mathrm R}BD_{\mathrm C}$. Moreover
\[
 (D_{\mathrm R}BD_{\mathrm C})(D_{\mathrm R}BD_{\mathrm C})^{\mathsf T}
       =D_{\mathrm R}BB^{\mathsf T}D_{\mathrm R}.
\]
Invariance of the trace proves separate evenness. In particular, with $z^{\alpha,s}$ denoting replacement of the $\alpha$ coordinate by $s$,
\begin{equation}
 \partial_\alpha\mathcal P_n(z)
 =z_\alpha\int_0^1\partial_{\alpha\alpha}^2
                   \mathcal P_n(z^{\alpha,tz_\alpha})\,\dd t
 \ge c_*z_\alpha\quad(z_\alpha\ge0).
 \label{B:eq:polynomial-coordinate-derivative}
\end{equation}

To check locality and uniformity, expand $\operatorname{Tr}(BB^{\mathsf T})^m=\tfrac12\operatorname{Tr}\mathcal A(B)^{2m}$. Every term is a product of the edge labels of a closed
walk of length $2m$ on this path. A term surviving differentiation in a specified edge contains that edge. Its starting vertex is within $2m$ steps of that edge, and there are at
most $2^{2m}$ walks from any starting vertex. An $h$th derivative has at most $(2m)^h$ assignments of its labels to the factors. Thus the number and coefficients of the monomials
in a fixed derivative are bounded in terms of $d,h,V$. Each such walk stays in a fixed neighborhood of every edge it visits. Taking, for example, $R_*=2d+2$ gives the required
interaction range. A surviving monomial of an $h$th derivative has degree $2m-h$ and is bounded by $(1+M)^{2d}$ when $\|z\|_\infty\le M$. Summing the bounded number of terms proves
the uniform derivative bounds.
\end{proof}

\subsubsection{Reference measures and their minimizers}
\label{B:subsection:reference-measures-and-minimizers}

The reference measure replaces the negative coordinate powers by zero. Its Hamiltonian is defined also at zero when the corresponding power vanishes, so we minimize it on the
closed positive orthant.

Put $c_\alpha^+=\max\{c_\alpha,0\}$ and $\alpha_\alpha^\sharp=c_\alpha^+/\kappa$. Define
\begin{equation}
 \overline{\mathbb P}_n(\dd z)=\overline Z_n^{-1}
       \exp[-\kappa\mathcal P_n(z)]\prod_\alpha z_\alpha^{c_\alpha^+}\,\dd z.
 \label{B:eq:reference-coefficient-density}
\end{equation}
Both normalizing constants are positive and finite. In fact, $\mathcal P_n(z)\ge\mathcal P_n(0)+c_*|z|^2/2$ by the preceding lemma, and each original exponent is greater than $-1$.
The integral of this Gaussian bound times the coordinate powers factors into finite one-dimensional integrals. Positivity follows by integrating on a compact cube in the open
orthant.

For $t\ge0$ use the extended convex function
\[
 \mathfrak b_t(s)=
 \begin{cases}
 -t\log s,&t>0,\ s>0,\\
 0,&t=0,\ s\ge0,\\
 +\infty,&\text{otherwise}.
 \end{cases}
\]
The reference Hamiltonian on the closed orthant is
\begin{equation}
 H_n^\sharp(z)=\mathcal P_n(z)+\sum_\alpha
                  \mathfrak b_{\alpha_\alpha^\sharp}(z_\alpha).
 \label{B:eq:reference-hamiltonian}
\end{equation}
For an interval of sites $I$, let $\mathscr S(I)=\{\alpha\in\mathscr S_n: \iota(\alpha)\in I\}$ and let $\partial_{R_*}I$ contain the existing sites outside $I$ at distance at most
$R_*$ from $I$. A boundary configuration $q$ specifies the coefficients in $\mathscr S(\partial_{R_*}I)$. After subtracting terms independent of the interior coefficients, write
\begin{equation}
 \mathcal P_{n,I}^q(u)=\mathcal P_n(u,e)-\mathcal P_n(0,e),\qquad
 H_{n,I}^{\sharp,q}(u)=\mathcal P_{n,I}^q(u)
       +\sum_{\alpha\in\mathscr S(I)}\mathfrak b_{\alpha_\alpha^\sharp}(u_\alpha),
 \label{B:eq:conditional-reference-hamiltonian}
\end{equation}
where $e$ is any extension of $q$ to the exterior. Locality makes the first difference independent of the remaining entries of $e$. Gradients omit the zero logarithmic terms and
extend smoothly in those coordinates to their boundary value zero.

\begin{lemma}\label{B:lem:reference-minimizer-coordinates}
For every finite nonnegative $q$, $H_{n,I}^{\sharp,q}$ has a unique minimizer $m^{q,I}$. Its smooth gradient is zero. Its $\alpha$ coordinate is zero exactly when $c_\alpha^+=0$.
If $\|q\|_\infty\le M$ and $c_\alpha^+>0$, then
\begin{equation}
 \sqrt{\alpha_\alpha^\sharp/C_M}\le m_\alpha^{q,I}
          \le\sqrt{\alpha_\alpha^\sharp/c_*}.
 \label{B:eq:reference-minimizer-coordinate-bounds}
\end{equation}
The upper bound holds for arbitrary finite nonnegative $q$.
\end{lemma}
\begin{proof}
The polynomial in \eqref{B:eq:conditional-reference-hamiltonian} is separately even in its interior coordinates and has Hessian at least $c_*I$. Thus
$\mathcal P_{n,I}^q(u)\ge c_*|u|^2/2$. The quadratic bound dominates the finite sum of logarithms as $|u|\to\infty$. For each $t\ge0$, the function $t\log s-c_*s^2/4$ is bounded
above on $s>0$ when $t>0$; for $t=0$ its supremum is zero. Therefore
\[
 H_{n,I}^{\sharp,q}(u)\ge\frac{c_*}{4}|u|^2-C_{n,I}.
\]
The extended Hamiltonian is lower semicontinuous and finite at the vector with all coordinates one. The sublevel set through that vector is compact and nonempty, so it contains a
minimizer. For distinct feasible $u,v$, strong convexity gives
\[
 H_{n,I}^{\sharp,q}((u+v)/2)
 \le\tfrac12 H_{n,I}^{\sharp,q}(u)+\tfrac12 H_{n,I}^{\sharp,q}(v)
       -\tfrac{c_*}{8}|u-v|^2.
\]
Two distinct minimizers would contradict this inequality.

A positive logarithmic coefficient forces the minimizing coordinate to be positive. For a zero coefficient, \eqref{B:eq:polynomial-coordinate-derivative} gives
\[
 \mathcal P_{n,I}^q(u^{\alpha,s})-
 \mathcal P_{n,I}^q(u^{\alpha,0})\ge c_*s^2/2\quad(s\ge0),
\]
so its minimizing value is zero. Stationarity holds at each positive coordinate; evenness gives a zero derivative at each zero coordinate. At a positive coordinate it follows that
\begin{equation}
 \alpha_\alpha^\sharp
 =(m_\alpha^{q,I})^2\int_0^1\partial_{\alpha\alpha}^2
       \mathcal P_{n,I}^q((m^{q,I})^{\alpha,t m_\alpha^{q,I}})\,\dd t.
 \label{B:eq:reference-minimizer-stationarity}
\end{equation}
The integral is at least $c_*$. Since $\alpha_{\sigma,k}^\sharp\le(k+a)/n$, this gives a uniform upper bound for every minimizing coordinate, independent of $q$. If $q$ is bounded
by $M$, all coefficients in the integral are in a fixed box. Lemma~\ref{B:lem:coefficient-polynomial-convexity} bounds the integral from above by $C_M$, proving the lower estimate.
\end{proof}

Let $z_n^\sharp$ be the full reference minimizer. The preceding proof gives
\begin{equation}
 0\le(z_n^\sharp)_{\sigma,k}\le C\sqrt{k/n}.
 \label{B:eq:full-reference-minimizer-upper}
\end{equation}
For $\beta\ge1$ the two probability measures and their closed Hamiltonians coincide, and we also write $z_n^o$ for the minimizer. For $\beta=1$ and $n\ge2$, $(z_n^\sharp)_{y,1}=0$
and all other existing coordinates are positive. The probability measures in \eqref{B:eq:original-coefficient-density} and \eqref{B:eq:reference-coefficient-density} have densities
on the open orthant.

\subsubsection{Constant coefficients and the deterministic profile}
\label{B:subsection:constant-coefficients-and-profile}

The finite minimizer varies with position. A constant coefficient calculation identifies the local profile against which it will be compared, together with the derivative bounds
needed near zero.

For constant coefficients define
\begin{equation}
 \mathcal W_V(x,y)=\frac1{2\pi}\int_0^{2\pi}
        V(x^2+y^2-2xy\cos t)\,\dd t
 =\sum_{m=1}^d g_m\sum_{j=0}^m\binom mj^2x^{2j}y^{2m-2j}.
 \label{B:eq:constant-coefficient-polynomial}
\end{equation}
For the equality, expand $(x-ye^{\ii t})^m(x-ye^{-\ii t})^m$; integration keeps exactly the terms with equal powers of $e^{\ii t}$ and $e^{-\ii t}$. Use $\phi$, $L_V(t)$, and
$\theta_V$ from \eqref{B:eq:1-12}--\eqref{B:eq:1-16}. In this subsection write $\ell(t)=L_V(t)$ and $\theta=\theta_V$.

\begin{lemma}\label{B:lem:constant-coefficient-minimizer}
For every $0<t\le1$, the unique minimizer of $\mathcal W_V(x,y)-t\log x-t\log y$ on $(0,\infty)^2$ is $(\phi(t),\phi(t))$. On $0<t\le1$, for every fixed integer $j\ge0$,
\begin{equation}
 |\phi^{(j)}(t)|\le C_j t^{1/2-j},\qquad
 |(\ell^{-1})^{(j)}(t)|\le C_j t^{-1/2-j}.
 \label{B:eq:profile-derivative-bounds}
\end{equation}
\end{lemma}
\begin{proof}
Choose $N>4d+4$, let $S_N$ be the cyclic shift, and put $C=xI-yS_N$. The eigenvalues of $CC^{\mathsf T}$ are the values of $x^2+y^2-2xy\cos t$ at the $N$th roots of unity. Every
nonconstant Fourier mode in the integrand of \eqref{B:eq:constant-coefficient-polynomial} has order at most $d<N$, so $N^{-1}\operatorname{Tr}V(CC^{\mathsf T})=\mathcal W_V(x,y)$
exactly. For a direction $(h,k)$, $N^{-1}\|hI-kS_N\|_{\rm HS}^2=h^2+k^2$. The dilation calculation in Lemma~\ref{B:lem:coefficient-polynomial-convexity} therefore gives
$D^2\mathcal W_V\ge c_*I$. Its gradient at $(0,0)$ vanishes. Quadratic growth at infinity and the logarithmic barriers at the coordinate axes give existence of a minimizer. Strict
convexity gives uniqueness, and symmetry in $x,y$ places it on the diagonal.

Set
\begin{equation}
 A_m=m\binom{2m}{m},\qquad F(r)=\sum_{m=1}^d A_mg_mr^{2m}.
 \label{B:eq:profile-polynomial-and-constants}
\end{equation}
Vandermonde's identity and symmetry give $\partial_x\mathcal W_V(r,r)=F(r)/r$. Thus the critical equation is $F(r)=t$, which is the defining equation for $\phi$. The proof of
Lemma~\ref{B:lem:bessel-coordinate-regularity} shows that
\[
 \phi(t)=\sqrt t\,v(t),\qquad \ell(t)=\sqrt t\,w(t),
 \qquad v,w\in C^\infty([0,1]),\qquad \min v,\min w>0.
\]
For $b\in\{1/2,-1/2\}$ and smooth $v_0$ on $[0,1]$, Leibniz' rule is
\[
 \frac{\mathrm d^j}{\mathrm dt^j}(t^bv_0(t))
 =\sum_{l=0}^j\binom jl(b)_l t^{b-l}v_0^{(j-l)}(t),
 \qquad (b)_l=b(b-1)\cdots(b-l+1),\quad (b)_0=1.
\]
All derivatives of $v$ and $1/w$ are bounded. Since $t^{b-l}\le t^{b-j}$ for $0<t\le1$ and $l\le j$, this formula proves \eqref{B:eq:profile-derivative-bounds}. Here $\ell^{-1}$
denotes the reciprocal function $1/\ell$. It also gives $\phi(t)\asymp\sqrt t$, $\ell(t)\asymp\sqrt t$, and $\phi'(t)\asymp t^{-1/2}$, the last using $F'(\phi)\phi'=1$ and
$F'(r)\asymp r$ on the compact interval $[0,\phi(1)]$.
\end{proof}

The coordinate regularity and estimates for the partition intervals established in Lemmas~\ref{B:lem:bessel-coordinate-regularity} and \ref{B:lem:brownian-grid-sums} apply to this
same $\theta$. In particular, $z_{n,k}=\theta(k/n)$ and $\Delta_{n,k}$ retain their previous definitions. For example, the estimate for the cell average follows directly from
\[
 n\Delta_{n,k}-\theta'(k/n)
 =-n\int_{(k-1)/n}^{k/n}\int_t^{k/n}\theta''(s)\,\dd s\,\dd t,
 \qquad |n\Delta_{n,k}-\theta'(k/n)|\le\|\theta''\|_\infty/(2n).
\]

\subsubsection{Derivatives at constant coefficients}
\label{B:subsection:closed-walks-and-force-coefficients}

The first correction to that profile depends on how the polynomial gradient responds to neighboring coefficients. The following finite-range calculation records those derivatives
before they are substituted into the stationary equations.

Write $J=BB^{\mathsf T}$, so $J_{kk}=x_k^2+y_{k-1}^2$ with $y_0=0$ and $J_{k,k+1}=-x_ky_k$. For an index whose interaction neighborhood is inside the finite path, put
$F_{i,m}^x=\partial_{x_i}\operatorname{Tr}J^m$ and $F_{i,m}^y=\partial_{y_i}\operatorname{Tr}J^m$. Define
\begin{equation}
 B_m=\frac{2m^2-2m+1}{2m-1}A_m,\qquad
 C_m=\frac{2m^2-2m}{2m-1}A_m,\qquad D_m=-C_m/2.
 \label{B:eq:polynomial-gradient-constants}
\end{equation}
At the constant configuration $x_j=y_j=r>0$, let
\[
 \alpha_{m,h}^{\sigma\tau}
 =r^{2-2m}\partial_{z_{\tau,i+h}}F_{i,m}^\sigma.
\]
Homogeneity makes these coefficients independent of $r$, and locality makes them zero outside a fixed finite range.

\begin{lemma}\label{B:lem:constant-gradient-coefficients}
At constant coefficients, $F_{i,m}^x=F_{i,m}^y=A_mr^{2m-1}$. The sums of the derivative coefficients are
\begin{equation}
 \sum_h\begin{pmatrix}\alpha_{m,h}^{xx}&\alpha_{m,h}^{xy}\\
                       \alpha_{m,h}^{yx}&\alpha_{m,h}^{yy}\end{pmatrix}
 =\begin{pmatrix}B_m&C_m\\C_m&B_m\end{pmatrix},\qquad
 \sum_hh\begin{pmatrix}\alpha_{m,h}^{xx}&\alpha_{m,h}^{xy}\\
                       \alpha_{m,h}^{yx}&\alpha_{m,h}^{yy}\end{pmatrix}
 =\begin{pmatrix}0&D_m\\-D_m&0\end{pmatrix}.
 \label{B:eq:gradient-coefficient-sums}
\end{equation}
\end{lemma}
\begin{proof}
The closed walks defining a derivative are contained in a fixed neighborhood of its index. They can therefore be computed on the cycle used in the preceding proof, with its length
increased if necessary. At constant coefficients all derivatives of each species agree by translation. Differentiating the constant coefficient trace and dividing by the cycle
length identifies differentiation in one coordinate with differentiation of $\mathcal W_m(x,y)=\sum_j\binom mj^2x^{2j}y^{2m-2j}$. Differentiating once more, with all coordinates of
one species changed simultaneously, identifies the four unweighted sums with the entries of its Hessian.

Put $S_j=\sum_{l=0}^m l^j\binom ml^2$, $j=0,1,2$. Vandermonde and the substitution $l\mapsto m-l$ give $S_0=\binom{2m}{m}$ and $S_1=mS_0/2$. For $m\ge2$,
\[
 S_2-S_1=m(m-1)\sum_l\binom{m-2}{l-2}\binom ml
 =m(m-1)\binom{2m-2}{m-2}
 =\frac{m(m-1)^2}{2(2m-1)}S_0.
\]
Hence $S_2=m^3S_0/[2(2m-1)]$. For $m=1$ these identities follow from $S_0=2$ and $S_1=S_2=1$. At $x=y=r$ we consequently have
\[
 \partial_x\mathcal W_m=2S_1r^{2m-1},\quad
 \partial_x^2\mathcal W_m=(4S_2-2S_1)r^{2m-2},\quad
 \partial_x\partial_y\mathcal W_m=(4mS_1-4S_2)r^{2m-2}.
\]
These are the value and constants in \eqref{B:eq:polynomial-gradient-constants}.

For the first moments, place the vertices $\mathrm R_j,\mathrm C_j$ of the dilation at $2j,2j+1$. Reflection about $2i+1/2$ fixes the $x_i$ edge and sends
\[
 x_{i+h}\mapsto x_{i-h},\qquad y_{i+h}\mapsto y_{i-h-1}.
\]
The trace of the even power of the dilation is invariant under this permutation; the coefficient signs on the reflected edges agree. It follows, on a finite path containing the
whole differentiated neighborhood, that $\alpha_{m,h}^{xx}=\alpha_{m,-h}^{xx}$ and $\alpha_{m,h}^{xy}=\alpha_{m,-h-1}^{xy}$. The neighborhoods agree with those of the original
matrix by locality. Thus
\[
 \sum_hh\alpha_{m,h}^{xx}=-\sum_hh\alpha_{m,h}^{xx}=0,
 \qquad
 2\sum_hh\alpha_{m,h}^{xy}=-\sum_h\alpha_{m,h}^{xy}=-C_m.
\]
Reflection about $2i+3/2$, the center of the $-y_i$ edge, sends $y_{i+h}\mapsto y_{i-h}$ and $x_{i+h}\mapsto x_{i+1-h}$. Consequently $\alpha_{m,h}^{yy}=\alpha_{m,-h}^{yy}$ and
$\alpha_{m,h}^{yx}=\alpha_{m,1-h}^{yx}$. For the two remaining entries, changing variables to $l=-h$ and $l=1-h$, respectively, gives
\[
 \sum_hh\alpha_{m,h}^{yy}=-\sum_l l\alpha_{m,l}^{yy},\qquad
 \sum_hh\alpha_{m,h}^{yx}=\sum_l(1-l)\alpha_{m,l}^{yx}
                       =C_m-\sum_l l\alpha_{m,l}^{yx}.
\]
Their first moments are consequently zero and $C_m/2$. Together with the two preceding computations these are all four entries of the second matrix in
\eqref{B:eq:gradient-coefficient-sums}.
\end{proof}

\subsubsection{The first correction and the critical equations}
\label{B:subsection:first-correction-and-critical-equations}

To reduce the gradient of the finite energy at the leading profile, we introduce two correction functions. They account for the different logarithmic powers and the variation of
the profile from site to site.

Define $g_x=x^{(1)}$ and $g_y=y^{(1)}$ by
\begin{equation}
 \begin{split}
 g_x(t)+g_y(t)&=(a-2/\beta)\phi'(t),\\
 g_x(t)-g_y(t)&=(a+\tfrac12)\ell(t)^{-1}-\tfrac12\phi'(t),\qquad
 f_{n,\sigma,k}=\phi(k/n)+n^{-1}g_\sigma(k/n).
 \end{split}
 \label{B:eq:first-correction-definition}
\end{equation}
These are the functions in \cite{rider2019universality}. Lemma~\ref{B:lem:constant-coefficient-minimizer} implies
\begin{equation}
 |g_\sigma^{(j)}(t)|\le C_jt^{-1/2-j},\qquad
 |f_{n,\sigma,k}/\phi(k/n)-1|\le C/k.
 \label{B:eq:first-correction-derivative-bounds}
\end{equation}
In particular $f_{n,\sigma,k}>0$ for $k$ above a fixed threshold.

\begin{lemma}\label{B:lem:first-correction-residual}
There are fixed integers $K_0,R_*$ and a constant $C$ such that, when $K_0+R_*<i\le n-R_*$, the local gradient evaluated at $f_n$ satisfies
\begin{equation}
 |\partial_{x_i}H_n(f_n)|+|\partial_{y_i}H_n(f_n)|
            \le C n^{-2}(i/n)^{-3/2}=C/(\sqrt n\,i^{3/2}).
 \label{B:eq:first-correction-gradient-residual}
\end{equation}
Only the coordinates in the interaction neighborhood are used in this formula.
\end{lemma}
\begin{proof}
Fix $t=i/n$, put $r=\phi(t)$, and abbreviate $g_x=g_x(t)$, $g_y=g_y(t)$ and $r'=\phi'(t)$. Choose $K_0>2R_*$ sufficiently large that all coefficients used below are positive and
their ratio to $\phi(k/n)$ lies in $[1/2,3/2]$. For $|h|\le R_*$, the derivative bounds for $\phi$ and $g_\sigma$, with $t+h/n\asymp t$, give
\[
 f_{n,\sigma,i+h}
 =r+n^{-1}(hr'+g_\sigma)+\epsilon_{\sigma,h},
 \qquad |\epsilon_{\sigma,h}|\le Cn^{-2}t^{-3/2}.
\]
Indeed, the second-order remainder of $\phi(t+h/n)$ is bounded by $Cn^{-2}t^{-3/2}$, and the first-order difference of $n^{-1}g_\sigma(t+h/n)$ from $n^{-1}g_\sigma(t)$ has the same
bound.

Define the three scalar coefficients
\[
 B=\sum_{m=1}^d g_mB_mr^{2m-2},\qquad
 C=\sum_{m=1}^d g_mC_mr^{2m-2},\qquad
 D=\sum_{m=1}^d g_mD_mr^{2m-2}=-C/2.
\]
Lemma~\ref{B:lem:constant-gradient-coefficients}, including its first moments in the site offset, gives
\[
 \begin{aligned}
 \partial_{x_i}\mathcal P_n(f_n)
 &=\frac{F(r)}r+n^{-1}(Bg_x+Cg_y+Dr')+O(n^{-2}t^{-3/2}),\\
 \partial_{y_i}\mathcal P_n(f_n)
 &=\frac{F(r)}r+n^{-1}(Cg_x+Bg_y-Dr')+O(n^{-2}t^{-3/2}).
 \end{aligned}
\]
To check these errors uniformly down to $t=i/n$, a differentiated polynomial monomial of degree $2m$ has degree $2m-1$. Its linear response to the $\epsilon_{\sigma,h}$ is bounded
by $Cn^{-2}t^{-3/2}$, because its first derivatives are bounded on the fixed coefficient box. For $m\ge2$, its second derivatives are bounded on the relevant segment by
$Cr^{2m-3}$. The coordinate change is at most $Cn^{-1}t^{-1/2}$, so the quadratic error is bounded by
\[
 C r^{2m-3}n^{-2}t^{-1}=Cn^{-2}t^{m-5/2}
 \le Cn^{-2}t^{-3/2}.
\]
For $m=1$ the differentiated monomial is linear and this quadratic error vanishes. There are only a fixed number of local monomials, so these bounds can be summed.

Set $b_x=a-\beta^{-1}$ and $b_y=-\beta^{-1}$. The logarithmic terms of the two stationary equations have coefficients $t+b_\sigma/n$. The exact identity
\[
 \frac1{r+u}=\frac1r-\frac{u}{r^2}+\frac{u^2}{r^2(r+u)}
\]
with $u=g_\sigma/n$, together with $|u|\le r/2$, $r\asymp t^{1/2}$ and $|g_\sigma|\le Ct^{-1/2}$, yields
\[
 \frac{t+b_\sigma/n}{r+g_\sigma/n}
 =\frac tr+\frac1n\left(\frac{b_\sigma}r-\frac{tg_\sigma}{r^2}\right)
       +O(n^{-2}t^{-3/2}).
\]
For example, the two omitted terms are bounded by $Ct u^2/r^3$ and $C|u|/(nr^2)$, both of the stated order.

Since $F(r)=t$, the zeroth-order terms cancel. Put $H=B+t/r^2$. The coefficients at order $n^{-1}$ are
\[
 H g_x+Cg_y+Dr'-\frac{a-\beta^{-1}}r,
 \qquad Cg_x+Hg_y-Dr'+\frac{\beta^{-1}}r.
\]
Their cancellation follows from the following two identities. The constants in \eqref{B:eq:polynomial-gradient-constants} satisfy
\[
 B_m+C_m=(2m-1)A_m,\qquad B_m-C_m=\frac{A_m}{2m-1}.
\]
Consequently
\[
 H+C=\frac{F'(r)}r=\frac1{rr'},\qquad
 H-C=\sum_{m=1}^d\frac{2mA_mg_m}{2m-1}r^{2m-2}=\frac{\ell(t)}r.
\]
For the last equality, the change of variables $s=F(u)$ gives
\[
 \ell(t)=\int_0^t\frac{\mathrm ds}{\phi(s)}
 =\int_0^r\frac{F'(u)}u\,\mathrm du
 =\sum_{m=1}^d\frac{2mA_mg_m}{2m-1}r^{2m-1}.
\]
The sum of the two order-$n^{-1}$ coefficients vanishes because
\[
 (H+C)(g_x+g_y)=\frac{a-2/\beta}{r}.
\]
Their difference vanishes because $2D=-C$ and
\[
 \begin{aligned}
 (H-C)(g_x-g_y)+2Dr'
 &=\frac{\ell(t)}r\left(\frac{a+1/2}{\ell(t)}-\frac{r'}2\right)-Cr'\\
 &=\frac{a+1/2}r-\frac{H+C}{2}r'=\frac ar.
 \end{aligned}
\]
Thus both gradient components are bounded by $Cn^{-2}t^{-3/2}=C/(\sqrt n\,i^{3/2})$. The prescribed distance from the two ends ensures that the local differentiated polynomial has
not been altered by a boundary restriction or a changed logarithmic power. The argument uses the original fixed $\beta>0$ throughout.
\end{proof}

\subsubsection{Products of deterministic coefficients}
\label{B:subsection:products-of-deterministic-coefficients}

The deterministic centering of the inverse involves products of the corrected coefficients. Taking logarithms makes their cumulative contribution explicit in the spatial coordinate
used by the limiting inverse kernel.

\begin{lemma}\label{B:lem:deterministic-coefficient-products}
Uniformly for $K_0<i<j\le n$, after increasing $K_0$ if necessary,
\begin{equation}
 \sum_{k=i}^{j-1}\log\frac{f_{n,y,k}}{f_{n,x,k}}
 =-\left(\frac a2+\frac14\right)\log\frac{\theta(j/n)}{\theta(i/n)}
       +\frac12\log\frac{\phi(j/n)}{\phi(i/n)}+O(i^{-1}).
 \label{B:eq:deterministic-product-centering}
\end{equation}
\end{lemma}
\begin{proof}
Set $v=(g_y-g_x)/\phi$. Lemma~\ref{B:lem:constant-coefficient-minimizer} and \eqref{B:eq:first-correction-derivative-bounds} give $|v(t)|\le C/t$ and $|v'(t)|\le C/t^2$. Increase
$K_0$ until $|g_\sigma(k/n)/(n\phi(k/n))|\le1/2$ for $k>K_0$. The bound $|\log(1+x)-x|\le x^2$ for $|x|\le1/2$ yields
\[
 \left|\log\frac{f_{n,y,k}}{f_{n,x,k}}-\frac1n v(k/n)\right|
 \le \frac{C}{k^2}.
\]
For each cell,
\[
 \left|\frac1n v(k/n)-\int_{k/n}^{(k+1)/n}v(t)\,\dd t\right|
 \le \int_{k/n}^{(k+1)/n}\int_{k/n}^{t}|v'(s)|\,\dd s\,\dd t
 \le \frac{C}{k^2}.
\]
Thus the error after summation from $i$ to $j-1$ is at most $C\sum_{k=i}^\infty k^{-2}\le C'/i$, including $j=n$. Finally, \eqref{B:eq:first-correction-definition} gives
\[
 v(t)=-\left(\frac a2+\frac14\right)\frac{\theta'(t)}{\theta(t)}
                +\frac12\frac{\phi'(t)}{\phi(t)}.
\]
Integration over $[i/n,j/n]$ proves \eqref{B:eq:deterministic-product-centering}.
\end{proof}

\subsubsection{Comparison of the two coefficient measures}
\label{B:subsection:reference-measures-and-density-comparison-at-the-singular-coordinates}

We transfer the convexity estimates to the original measure by estimating the density ratio. The exceptional sets below distinguish coordinates where the reference logarithmic
coefficient vanishes from those where the two measures differ.

Let $\mathscr S_- =\{\alpha:c_\alpha<0\}$ and $\mathscr S_0=\{\alpha:c_\alpha\le0\}$. These sets contain existing coordinates only. For sufficiently large $n$ they are fixed finite
sets, and
\[
 \mathscr S_-=
 \{(y,k):k<\beta^{-1}\}\cup\{(x,k):k+a<\beta^{-1}\}.
\]
The reference Hamiltonian has zero logarithmic coefficient on $\mathscr S_0$, whereas the measures differ only at $\mathscr S_-$. Increase a fixed $n_0$ so that these sets and
their interaction neighborhoods have stabilized and $\kappa\ge1$ for $n\ge n_0$.

\begin{lemma}\label{B:lem:reference-positive-moments}
For every fixed $s\ge0$, $\sup_{n\ge n_0}\sup_{\alpha\in\mathscr S_n} \overline{\mathbb E}_n Z_\alpha^s<\infty$. If $c_\alpha^+=0$, then
$\overline{\mathbb E}_n(\sqrt\kappa Z_\alpha)^2\le c_*^{-1}$.
\end{lemma}
\begin{proof}
Integrate the derivative of $z_\alpha^{s+1}$ times the unnormalized reference density. At zero the boundary factor has exponent $s+1+c_\alpha^+>0$. At infinity it vanishes by the
quadratic lower bound for $\mathcal P_n$. Polynomial derivatives are integrable against that bound; for each fixed $n$ the other coordinate powers also have finite Gaussian
integrals. Integration in one coordinate followed by Fubini therefore gives
\begin{equation}
 \begin{split}
 \kappa\overline{\mathbb E}_n
 [Z_\alpha^{s+1}\partial_\alpha\mathcal P_n(Z)]
   &=(s+1+c_\alpha^+)\overline{\mathbb E}_n Z_\alpha^s,\\
 \overline{\mathbb E}_n Z_\alpha^{s+2}
   &\le\frac{s+1+c_\alpha^+}{\kappa c_*}
                       \overline{\mathbb E}_n Z_\alpha^s.
 \end{split}
 \label{B:eq:reference-moment-recursion}
\end{equation}
The second inequality uses \eqref{B:eq:polynomial-coordinate-derivative}. Since $c_\alpha^+\le\beta(n+a)$, its coefficient is at most $c_*^{-1}[(s+1)/\beta+1+a]$. At even orders
the recursion gives, for each integer $l\ge1$,
\[
 \overline{\mathbb E}_n Z_\alpha^{2l}
 \le\prod_{j=0}^{l-1}
       \frac{2j+1+c_\alpha^+}{\kappa c_*}
 \le\prod_{j=0}^{l-1}c_*^{-1}\bigl[(2j+1)/\beta+1+a\bigr].
\]
For $0<s\le2l$, Jensen's inequality applied to $x^{s/(2l)}$ bounds $\overline{\mathbb E}_n Z_\alpha^s$ by the $s/(2l)$ power of this product. The order $s=0$ has expectation one.
At a zero exponent the same inequality with $s=0$ proves the second moment bound.
\end{proof}

For $\alpha\in\mathscr S_-$ put $\rho_\alpha=-c_\alpha\in(0,1)$ and set
\begin{equation}
 w_n(z)=\prod_{\alpha\in\mathscr S_-}
              (\sqrt\kappa z_\alpha)^{-\rho_\alpha},\qquad
 L_n(z)=\frac{w_n(z)}{\overline{\mathbb E}_n w_n},\qquad
 \rho_*=(1-\beta)_+.
 \label{B:eq:initial-density-weight}
\end{equation}
Empty products equal one. For $0<\beta<1$ and $n\ge n_0$, $\rho_*=\max_{\mathscr S_-}\rho_\alpha$, attained at $(y,1)$.

\begin{theorem}\label{B:thm:initial-density-comparison}
For each finite $p>1$ with $p\rho_*<1$, uniformly in $n\ge n_0$,
\begin{equation}
 \frac{\dd\mathbb P_n}{\dd\overline{\mathbb P}_n}=L_n,\qquad
 \overline{\mathbb E}_n w_n\ge c>0,\qquad
 \overline{\mathbb E}_n w_n^p\le C_p,\qquad
 \|L_n\|_{L^p(\overline{\mathbb P}_n)}\le C_p.
 \label{B:eq:initial-density-moment-bounds}
\end{equation}
For $0<\beta<1$, $\overline{\mathbb E}_n L_n^p=\infty$ when $p\rho_*\ge1$. For $\beta\ge1$, $w_n=L_n=1$.
\end{theorem}
\begin{proof}
It suffices to treat $0<\beta<1$. Condition on all coordinates outside $\mathscr S_-$, denoting them by $v$, and put $x_\alpha=\sqrt\kappa z_\alpha$ on this fixed set. Its
conditional reference density is proportional to $\exp[-F_v(x)]$ on $\mathbb R_+^{m_-}$, where $m_-=|\mathscr S_-|$ and
\[
 F_v(x)=\kappa[\mathcal P_n(v,x/\sqrt\kappa)-\mathcal P_n(v,0)]
                 \ge c_*|x|^2/2.
\]
Separate evenness gives zero gradient in these coordinates at zero. All monomials not meeting $\mathscr S_-$ cancel. The remaining terms involve a fixed set $\mathscr N$ of
neighboring coordinates of $v$. For $0\le x_\alpha\le1$, the integral Taylor formula and the local polynomial derivative bound give
\[
 F_v(x)=\int_0^1(1-t)
 D^2_{\mathscr S_-}\mathcal P_n(v,tx/\sqrt\kappa)[x,x]\,\dd t
 \le C(1+\|v_{\mathscr N}\|)^{2d}|x|^2.
\]
Choose $h_v=h_0(1+\|v_{\mathscr N}\|)^{-d}$ with fixed small $0<h_0\le1$. On $[0,h_v]^{m_-}$ the last expression is bounded by a fixed constant. Hence
\[
 \int_{\mathbb R_+^{m_-}}e^{-F_v(x)}\,\dd x\ge c h_v^{m_-}.
\]
For $p\rho_\alpha<1$ the numerator in the conditional moment satisfies
\[
 \begin{split}
 \int_{\mathbb R_+^{m_-}}\prod_\alpha x_\alpha^{-p\rho_\alpha}
                  e^{-F_v(x)}\,\dd x
 &\le\prod_\alpha\int_0^\infty
                  x^{-p\rho_\alpha}e^{-c_*x^2/2}\,\dd x\\
 &=\prod_\alpha\frac12
       (2/c_*)^{(1-p\rho_\alpha)/2}
       \Gamma((1-p\rho_\alpha)/2)<\infty.
 \end{split}
\]
It follows that $\overline{\mathbb E}_n[w_n^p\mid v] \le C_p(1+\|v_{\mathscr N}\|)^{dm_-}$. Write $b=dm_-$ and $N_*=|\mathscr N|$, which is bounded independently of $n$. For $b>0$
the finite-dimensional inequality
\[
 (1+\|v_{\mathscr N}\|)^b
 \le C_{b,N_*}\left(1+\sum_{\gamma\in\mathscr N}|v_\gamma|^b\right)
\]
and Lemma~\ref{B:lem:reference-positive-moments} give a uniform expectation bound. Thus $\overline{\mathbb E}_n w_n^p\le C_p$. The conditioning used here leaves the full exterior
marginal of $\overline{\mathbb P}_n$ unchanged.

The same lemma gives $\overline{\mathbb E}_n(\sqrt\kappa Z_\alpha)^2\le c_*^{-1}$ on $\mathscr S_-$. Choose $K\ge\max\{1,(2m_-/c_*)^{1/2}\}$. The union bound implies
\[
 \overline{\mathbb P}_n\{\max_{\mathscr S_-}\sqrt\kappa Z_\alpha\le K\}
       \ge\tfrac12,\qquad
 \overline{\mathbb E}_n w_n\ge\tfrac12K^{-\sum_\alpha\rho_\alpha}>0.
\]
Put $b_- =\sum_{\alpha\in\mathscr S_-}\rho_\alpha$. Dividing the unnormalized densities and then integrating gives
\[
 \prod_{\alpha\in\mathscr S_-}z_\alpha^{-\rho_\alpha}
       =\kappa^{b_-/2}w_n(z),\qquad
 \frac{Z_n}{\overline Z_n}
       =\kappa^{b_-/2}\overline{\mathbb E}_nw_n.
\]
Consequently
\[
 \frac{\dd\mathbb P_n}{\dd\overline{\mathbb P}_n}
 =\frac{\overline Z_n}{Z_n}\kappa^{b_-/2}w_n
 =\frac{w_n}{\overline{\mathbb E}_nw_n},\qquad
 \|L_n\|_p\le c^{-1}(\overline{\mathbb E}_nw_n^p)^{1/p}\le C_p.
\]
Both factors involving $\kappa$ have canceled. The finite positive normalizer and the established lower bound also justify this formula when taking moments.

Finally suppose $p\rho_\alpha\ge1$. Restrict every other coordinate to a fixed compact interval in $(0,\infty)$ and let this coordinate range in $(0,1)$. At fixed $n$ its exponent
in the reference density is zero, and all remaining density factors are bounded below by a positive constant on this set. The integral of $z_\alpha^{-p\rho_\alpha}$ diverges. This
proves sharpness.
\end{proof}

We shall use the fixed exponents
\begin{equation}
 \begin{gathered}
 p_\beta=\begin{cases}(2-\beta)/(2(1-\beta)),&0<\beta<1,\\
                         \infty,&\beta\ge1,\end{cases}\qquad
 \gamma_\beta=1-p_\beta^{-1}
       =\begin{cases}\beta/(2-\beta),&0<\beta<1,\\1,&\beta\ge1,
          \end{cases}\\
 \mathbb P_n(A)\le C_\beta\overline{\mathbb P}_n(A)^{\gamma_\beta},\qquad
 \mathbb E_n F\le C_\beta
          (\overline{\mathbb E}_n F^{1/\gamma_\beta})^{\gamma_\beta}
          \quad(F\ge0).
 \end{gathered}
 \label{B:eq:coefficient-measure-comparison}
\end{equation}
For $\beta<1$, $p_\beta\rho_*=(2-\beta)/2<1$, so H\"older's inequality proves both bounds, the second by monotone truncation if needed. For $\beta\ge1$ they hold with equality and
$C_\beta=1$. In particular, an exponential probability bound remains exponential, and a bound $Cn^{-D}$ becomes $C_\beta n^{-D\gamma_\beta}$. The exponent is fixed before choosing
a probability threshold. For events $A_n$ with $\overline{\mathbb P}_n(A_n)\le Cn^{-D}$,
\[
 \sum_{n\ge N}\mathbb P_n(A_n)
 \le C_\beta\sum_{n\ge N}n^{-D\gamma_\beta}<\infty
 \quad\text{provided }D\gamma_\beta>1.
\]
When $A_n$ is a union of coordinate events, its whole reference probability is estimated before applying this inequality. For a nonnegative random variable $F$, truncation and the
same H\"older inequality give, for every fixed
$s>0$,
\[
 \|F\|_{L^s(\mathbb P_n)}
 \le C_\beta^{1/s}\|F\|_{L^{s/\gamma_\beta}
                                     (\overline{\mathbb P}_n)}.
\]
Thus a reference moment of order $s$ alone does not give an original moment of that order when $0<\beta<1$. For the interior estimates below, we use the common conditional
probability kernel and thus avoid this change of moment order.

\subsubsection{Inverse Hessians and dependence on boundary data}
\label{B:subsection:reference-minimizers-and-conditional-centering}

A change in the exterior coordinates enters the stationary equations near the boundary of the interval. Exponential decay of the inverse Hessian makes its effect small in the
interior. The estimate below applies both to Hessians of $H_{n,I}^{\sharp,q}$ and to their averages along line segments, as needed for comparing minimizers. Boundary dependence in
this setting also appears in \cite{rider2019universality}.

\begin{lemma}\label{B:lem:inverse-hessian-decay}
Suppose a real symmetric matrix $A$, indexed by a subset of $\mathscr S_n$, satisfies $A\ge c_*I$, $A_{\alpha\gamma}=0$ for $|\iota(\alpha)-\iota(\gamma)|>R_*$, and
$|A_{\alpha\gamma}|\le C_M$ for $\alpha\ne\gamma$. Then
\begin{equation}
 |(A^{-1})_{\alpha\gamma}|
        \le (2/c_*)e^{-\varepsilon_M|\iota(\alpha)-\iota(\gamma)|}
 \label{B:eq:inverse-hessian-entry-bound}
\end{equation}
for a positive $\varepsilon_M$ depending only on $c_*,C_M,R_*$.
\end{lemma}
\begin{proof}
Fix $\gamma$ and let $D_{\alpha\alpha}= e^{\varepsilon|\iota(\alpha)-\iota(\gamma)|}$. The diagonal of $E=DAD^{-1}-A$ vanishes. A nonzero off-diagonal entry has absolute value at
most $C_M(e^{\varepsilon R_*}-1) \le C_M\varepsilon R_*e^{R_*}$ for $0<\varepsilon\le1$. Each row and column has at most $2(2R_*+1)$ potentially nonzero entries. Thus
\[
 \|E\|\le\bigl(\max_\alpha\sum_\delta|E_{\alpha\delta}|
                   \max_\delta\sum_\alpha|E_{\alpha\delta}|\bigr)^{1/2}
 \le C_M'\varepsilon.
\]
Choose $\varepsilon=\min\{1,c_*/(2C_M')\}$ when $C_M'>0$ and $\varepsilon=1$ otherwise. Then $A+E=A(I+A^{-1}E)$ and $\|A^{-1}E\|\le1/2$. The geometric series converges in operator
norm, and its finite partial sums multiplied on either side by $I+A^{-1}E$ have remainder of norm at most $2^{-N-1}$. Hence
\[
 (A+E)^{-1}=\sum_{j=0}^\infty(-A^{-1}E)^jA^{-1},\qquad
 \|(A+E)^{-1}\|\le\frac{\|A^{-1}\|}{1-\|A^{-1}E\|}\le2/c_*.
\]
Since $(A+E)^{-1}=DA^{-1}D^{-1}$ and $D_{\gamma\gamma}=1$,
\[
 e^{\varepsilon|\iota(\alpha)-\iota(\gamma)|}
 |(A^{-1})_{\alpha\gamma}|
 =|\langle e_\alpha,(A+E)^{-1}e_\gamma\rangle|\le2/c_*.
\]
Only the off-diagonal entries enter $E$; the proof requires no upper bound on the diagonal of $A$.
\end{proof}

\begin{proposition}\label{B:prop:conditional-minimizer-dependence}
For $\|q\|_\infty,\|q'\|_\infty\le M$,
\begin{equation}
 |m_\alpha^{q,I}-m_\alpha^{q',I}|
 \le C_M\sum_{\gamma\in\mathscr S(\partial_{R_*}I)}
 e^{-c_M|\iota(\alpha)-\iota(\gamma)|}|q_\gamma-q'_\gamma|.
 \label{B:eq:conditional-minimizer-boundary-dependence}
\end{equation}
The minimizer is a continuous function of finite nonnegative boundary data. Its restriction to a subinterval is the minimizer with the induced boundary data.
\end{proposition}
\begin{proof}
Let $v$ be a bounded feasible vector, with strictly positive coordinates at every positive logarithmic coefficient, and put $m=m^{q,I}$. Along the segment joining $v$ to $m$ the
averaged Hessian is
\begin{equation}
 A=\int_0^1\nabla^2\mathcal P_{n,I}^q(v+t(m-v))\,\dd t
   +\operatorname{diag}\left(
    \begin{cases}\alpha_\alpha^\sharp/(v_\alpha m_\alpha),
                       &\alpha_\alpha^\sharp>0,\\
                  0,&\alpha_\alpha^\sharp=0.
    \end{cases}\right).
 \label{B:eq:averaged-conditional-hessian}
\end{equation}
For $v_\alpha\ne m_\alpha$, integration gives
\[
 \int_0^1\frac{\dd t}{((1-t)v_\alpha+tm_\alpha)^2}
 =\frac{v_\alpha^{-1}-m_\alpha^{-1}}{m_\alpha-v_\alpha}
 =\frac1{v_\alpha m_\alpha}.
\]
For $v_\alpha=m_\alpha>0$, the integrand is constant $v_\alpha^{-2}$, with the same value for its integral. The smooth extension is polynomial at zero logarithmic coefficients. The
mean-value identity and stationarity give
\begin{equation}
 m-v=-A^{-1}\nabla H_{n,I}^{\sharp,q}(v),\qquad
 |m_\alpha-v_\alpha|\le C_M\sum_{\delta\in\mathscr S(I)}
       e^{-c_M|\iota(\alpha)-\iota(\delta)|}
           |\partial_\delta H_{n,I}^{\sharp,q}(v)|.
 \label{B:eq:minimizer-residual-estimate}
\end{equation}
The Hessian lower bound, finite range, and off-diagonal bounds follow from Lemma~\ref{B:lem:coefficient-polynomial-convexity}; hence the preceding lemma applies to $A$.

Take $v=m^{q',I}$ and subtract the stationary equation at $q'$. The logarithmic terms do not depend on the boundary values, so they cancel. For each $\delta$, the resulting
expression for $\partial_\delta H_{n,I}^{\sharp,q}(m^{q',I})$ is
\[
 \int_0^1\sum_{\gamma\in\mathscr S(\partial_{R_*}I)}
 \partial_{u_\delta q_\gamma}^2
 \mathcal P_{n,I}^{q'+t(q-q')}(m^{q',I})(q_\gamma-q'_\gamma)\,\dd t.
\]
The integrand is zero unless $|\iota(\delta)-\iota(\gamma)|\le R_*$ and otherwise is bounded by $C_M|q_\gamma-q'_\gamma|$. For each $\gamma$ at most $2(2R_*+1)$ indices $\delta$
occur. For each contributing pair $(\delta,\gamma)$, the triangle inequality gives
\[
 e^{-c_M|\iota(\alpha)-\iota(\delta)|}
 \le e^{c_MR_*}e^{-c_M|\iota(\alpha)-\iota(\gamma)|}.
\]
Substituting this expression for $\partial_\delta H_{n,I}^{\sharp,q}(m^{q',I})$ in \eqref{B:eq:minimizer-residual-estimate} and summing over the at most $2(2R_*+1)$ choices of
$\delta$ gives
\[
 |m_\alpha^{q,I}-m_\alpha^{q',I}|
 \le C_M^2\,2(2R_*+1)e^{c_MR_*}
 \sum_{\gamma\in\mathscr S(\partial_{R_*}I)}
 e^{-c_M|\iota(\alpha)-\iota(\gamma)|}|q_\gamma-q'_\gamma|.
\]
Enlarging the constant proves \eqref{B:eq:conditional-minimizer-boundary-dependence}. It proves continuity on every bounded set of boundary data, hence continuity on their
finite-dimensional nonnegative domain. To identify the minimizer on a smaller interval, let $I'$ be a subinterval of $I$. Let $q''$ be the boundary values induced by $(m,q)$ on
$\partial_{R_*}I'$, and denote by $m[u]$ the vector obtained by replacing $m|_{I'}$ by a feasible vector $u$. The definition of the conditional Hamiltonian gives the exact energy
difference
\[
 H_{n,I}^{\sharp,q}(m[u])-H_{n,I}^{\sharp,q}(m)
   =H_{n,I'}^{\sharp,q''}(u)
                 -H_{n,I'}^{\sharp,q''}(m|_{I'}).
\]
Terms independent of $u$ cancel, and all terms involving an interior coordinate occur in both differences with the same boundary values. Minimality of $m$ makes the left side
nonnegative for every $u$. Hence $m|_{I'}$ minimizes the conditional energy on $I'$; its strict convexity makes this minimizer unique.
\end{proof}

For nonnegative exponent arrays bounded above by $C(k+1)$ before division by $\kappa$, the coordinate stationarity equation bounds each minimizing coordinate by $C\sqrt{(k+1)/n}$.
The polynomial Hessian, its finite range, and the boundary derivative bound are unchanged, so \eqref{B:eq:conditional-minimizer-boundary-dependence} holds for these arrays. In
particular, if two such arrays agree with $c_\alpha$ beyond a fixed site $K$ and differ only at the first $K$ sites, their full minimizers $z^\sharp,\widetilde z^\sharp$ satisfy
\begin{equation}
 |z^\sharp_{\sigma,k}-\widetilde z^\sharp_{\sigma,k}|
      \le C n^{-1/2}e^{-c(k-K)}\quad(k>K).
 \label{B:eq:finite-exponent-change-minimizer}
\end{equation}
To see this, restrict to $[K+1,n]$. The two conditional Hamiltonians there have the same exponents, and the differing boundary coordinates are $O(n^{-1/2})$ by
\eqref{B:eq:reference-minimizer-coordinate-bounds}. There are only finitely many such coordinates. Applying \eqref{B:eq:conditional-minimizer-boundary-dependence} gives the
estimate. For $K$ sufficiently large, both minimizing coordinates are bounded below by $c\sqrt{k/n}$. The mean-value formula for the logarithm thus also gives, for $K<j\le v<n$,
\[
 \sum_{k=j}^v\left|\log\frac{z^\sharp_{y,k}}{z^\sharp_{x,k}}
                  -\log\frac{\widetilde z^\sharp_{y,k}}{\widetilde z^\sharp_{x,k}}\right|
 \le C\sum_{k=j}^v k^{-1/2}e^{-c(k-K)}
 \le Cj^{-1/2}e^{-c(j-K)}.
\]

\subsubsection{Approximation of the full minimizer}
\label{B:subsection:approximation-of-the-full-minimizer}

The next estimate extends the deterministic approximation of \cite{rider2019universality} to the convex reference Hamiltonian \eqref{B:eq:reference-hamiltonian}. Its proof keeps
the original $\beta$ in \eqref{B:eq:first-correction-definition}.

\begin{theorem}\label{B:thm:reference-minimizer-approximation}
There are constants $C,c,K$ and $n_0$ such that every existing coordinate satisfies, for $n\ge n_0$,
\begin{equation}
 |(z_n^\sharp)_{\sigma,k}-f_{n,\sigma,k}|
 \le \frac C{\sqrt n\,k^{3/2}}
       +\frac C{\sqrt n}e^{-c(k-K)_+}+Ce^{-c(n-k)}.
 \label{B:eq:reference-minimizer-absolute-error}
\end{equation}
\end{theorem}
\begin{proof}
Increase $K_0$ from Lemma~\ref{B:lem:first-correction-residual} so that it contains $\mathscr S_0$ and its interaction neighborhood, $K_0>2/\beta+2R_*$, and
$1/2\le f_{n,\sigma,k}/\phi(k/n)\le3/2$ for $k>K_0$. Use the feasible comparison vector
\begin{equation}
 v_{n,\sigma,k}=
 \begin{cases}
 0,&c_{\sigma,k}\le0,\\
 \sqrt{k/n},&c_{\sigma,k}>0,\ k\le K_0,\\
 f_{n,\sigma,k},&k>K_0.
 \end{cases}
 \label{B:eq:feasible-reference-comparison}
\end{equation}
It is bounded uniformly, is positive at every positive logarithmic coefficient, and agrees with $z_n^\sharp$ at its zero coordinates. Also
$|v_{n,\sigma,k}-f_{n,\sigma,k}| \le Cn^{-1/2}\mathbf1_{\{k\le K_0\}}$.

For $K_0+R_*<k\le n-R_*$ the entire interaction neighborhood equals $f_n$ and the exponents agree with the original ones. Thus \eqref{B:eq:first-correction-gradient-residual}
applies. For $k\le K_0+R_*$ all comparison coordinates in this neighborhood are $O(n^{-1/2})$. Equation~\eqref{B:eq:polynomial-coordinate-derivative} and the derivative bound on a
fixed box give a polynomial derivative of order $n^{-1/2}$. A positive logarithmic coefficient is $O(n^{-1})$ and the corresponding $v_{n,\sigma,k}$ is at least $cn^{-1/2}$, so its
logarithmic derivative has the same bound. At a zero coordinate the smooth derivative is zero. At the last $R_*$ sites all comparison coordinates are bounded above and below by
positive constants, giving a bounded derivative. Consequently, for $K=K_0+R_*$,
\begin{equation}
 |\partial_{\sigma,k}H_n^\sharp(v_n)|
 \le C\left[n^{-1/2}\mathbf1_{\{k\le K\}}
       +\frac{\mathbf1_{\{K<k\le n-R_*\}}}{\sqrt n\,k^{3/2}}
       +\mathbf1_{\{k>n-R_*\}}\right].
 \label{B:eq:reference-comparison-residual}
\end{equation}
Use \eqref{B:eq:minimizer-residual-estimate} on the full interval. The two finite boundary sums obey
\[
 \sum_{j=1}^K e^{-c|k-j|}\le K e^{-c(k-K)_+},\qquad
 \sum_{j=n-R_*+1}^{n}e^{-c|k-j|}
 \le R_*e^{cR_*}e^{-c(n-k)}.
\]
Multiplying these sums by the respective derivative bounds $Cn^{-1/2}$ and $C$ gives $Cn^{-1/2}e^{-c(k-K)_+}$ and $Ce^{-c(n-k)}$. For the middle contribution, split the sum at
$j=k/2$:
\[
 \begin{split}
 \sum_{j\ge k/2}e^{-c|k-j|}j^{-3/2}&\le Ck^{-3/2}
                      \sum_{l\in\mathbb Z}e^{-c|l|},\\
 \sum_{j<k/2}e^{-c|k-j|}j^{-3/2}&\le e^{-ck/2}
                       \sum_{j\ge1}j^{-3/2}\le Ck^{-3/2}.
 \end{split}
\]
There are at most two coordinates at each site. Adding the finite initial difference $v_n-f_n$ proves \eqref{B:eq:reference-minimizer-absolute-error}.
\end{proof}

Division by $\phi(k/n)\asymp\sqrt{k/n}$ gives
\begin{equation}
 \frac{|(z_n^\sharp)_{\sigma,k}-f_{n,\sigma,k}|}{\phi(k/n)}
 \le Ck^{-2}+Ck^{-1/2}e^{-c(k-K)_+}
                 +C\sqrt{n/k}\,e^{-c(n-k)}.
 \label{B:eq:reference-minimizer-relative-error}
\end{equation}
On $r_n\le k\le n-r_n$ with $r_n=\lceil(\log n)^A\rceil$ and $A>1$, the last two terms are bounded by $Ce^{-c'r_n}$ after adjusting constants. For the last term, use
$\sqrt{n/k}\le\sqrt2$ when $k\ge n/2$ and $\sqrt n e^{-cn/2}\le C e^{-cn/4}$ otherwise. Therefore
\begin{equation}
 \begin{split}
 |(z_n^\sharp)_{\sigma,k}/\phi(k/n)-1|
   +|\log((z_n^\sharp)_{\sigma,k}/\phi(k/n))|
       \le Ck^{-1}+Ce^{-c'r_n}.
 \end{split}
 \label{B:eq:reference-minimizer-interior-ratio}
\end{equation}
Choose $n_0$ so that the displayed relative error is at most $1/2$ for $n\ge n_0$. The ratio then lies in $[1/2,3/2]$, where $|\log x|\le2|x-1|$. For $V(\lambda)=\lambda/2$ the
reference minimizer is explicitly
\[
 (z_n^\sharp)_{x,k}=\sqrt{(k+a-\beta^{-1})_+/n},\qquad
 (z_n^\sharp)_{y,k}=\sqrt{(k-\beta^{-1})_+/n}.
\]
For fixed $b$ and $k>2|b|$, the integral remainder of the expansion for $\sqrt{k+b}$ gives $|\sqrt{(k+b)/n}-\sqrt{k/n}-b/(2\sqrt{nk})| \le C_b/(\sqrt n\,k^{3/2})$, consistent with
the estimate above.

\subsubsection{Conditional distributions on the interior interval}
\label{B:subsection:conditional-distributions-on-the-interior-interval}

Beyond the exceptional coordinates, the original and reference measures have the same conditional density. This is the conditional law used for the interior construction; the
exterior coordinates keep their original joint distribution.

If $\mathscr S(I)$ is disjoint from $\mathscr S_-$, all its exponents are nonnegative. For finite nonnegative boundary data define
\begin{equation}
 \begin{split}
 Z_{n,I}(q)&=\int_{(0,\infty)^{\mathscr S(I)}}
        e^{-\kappa\mathcal P_{n,I}^q(u)}\prod_{\alpha\in\mathscr S(I)}
                            u_\alpha^{c_\alpha}\,\dd u,\\
 \mu_{n,I}^q(\dd u)&=Z_{n,I}(q)^{-1}
        e^{-\kappa\mathcal P_{n,I}^q(u)}\prod_{\alpha\in\mathscr S(I)}
                            u_\alpha^{c_\alpha}\,\dd u.
 \end{split}
 \label{B:eq:initial-coefficients-common-kernel}
\end{equation}

\begin{proposition}\label{B:prop:unchanged-conditional-distributions}
Let $E$ be the full exterior of $I$. The original and reference conditional laws on $I$ agree, and
\begin{equation}
 \frac{\dd\mathbb P_n^{\rm ext}}{\dd\overline{\mathbb P}_n^{\rm ext}}(e)=L_n(e),
 \qquad
 \mathcal L_{\mathbb P_n}(Z_I\mid E=e)
   =\mathcal L_{\overline{\mathbb P}_n}(Z_I\mid E=e)
   =\mu_{n,I}^{q(e)}
 \label{B:eq:exterior-density-and-conditional-law}
\end{equation}
for almost every $e$. Formula~\eqref{B:eq:initial-coefficients-common-kernel} is a probability kernel continuous in total variation on all finite nonnegative boundary
configurations.
\end{proposition}
\begin{proof}
The density $L_n$ depends only on $E$. Disintegrating the reference measure and using the density identity, for every nonnegative Borel $F$,
\begin{equation}
 \mathbb E_nF(E,Z_I)=\int L_n(e)\overline{\mathbb P}_n^{\rm ext}(\dd e)
                         \int F(e,u)\mu_{n,I}^{q(e)}(\dd u).
 \label{B:eq:exterior-disintegration-identity}
\end{equation}
Taking $F$ independent of $u$ identifies the exterior marginal. The displayed identity then identifies the conditional law; the exterior null sets agree since $L_n$ is positive and
finite almost surely. Conditioning only on the boundary $Q=q(E)$ gives a different density ratio for its marginal. Let $\overline\lambda^Q_n$ be the reference law of $Q$ and choose
the conditional expectation $\ell_n^Q(q)=\overline{\mathbb E}_n[L_n(E)\mid Q=q]$. For bounded Borel $b$ and $f$,
\[
 \mathbb E_n[b(Q)f(Z_I)]
 =\int b(q)\ell_n^Q(q)
       \left(\int f(u)\,\mu_{n,I}^{q}(\dd u)\right)
                         \overline\lambda^Q_n(\dd q).
\]
Indeed, apply \eqref{B:eq:exterior-disintegration-identity} to $b(q(e))f(u)$ and then take the reference conditional expectation given $Q$. In particular
$\dd\lambda^Q_n/\dd\overline\lambda^Q_n =\ell_n^Q$, whereas the full exterior ratio is $L_n(e)$. The kernel on $I$ is still $\mu_{n,I}^q$. In general $L_n(e)$ cannot be evaluated
from the boundary coordinates alone.

For its pointwise definition, the polynomial bound $\mathcal P_{n,I}^q(u)\ge c_*|u|^2/2$ gives $0<Z_{n,I}(q)<\infty$ for every finite nonnegative $q$. On a compact set of $q$ the
unnormalized densities are dominated by $e^{-\kappa c_*|u|^2/2}\prod u_\alpha^{c_\alpha}$, an integrable function. They are continuous in $q$ at each $u$. Dominated convergence
therefore gives continuity of the normalizer and of the densities in $L^1$. The normalizer is bounded below on a compact set of $q$: integration over $[1,2]^{\mathscr S(I)}$ has a
uniform positive lower bound there. Write $f_q$ for the unnormalized density. For $q_j\to q$ the preceding domination gives $\|f_{q_j}-f_q\|_1\to0$ and
$Z_{n,I}(q_j)\to Z_{n,I}(q)>0$. In particular,
\[
 \left\|\frac{f_{q_j}}{Z_{n,I}(q_j)}-\frac{f_q}{Z_{n,I}(q)}\right\|_1
 \le\frac{\|f_{q_j}-f_q\|_1+|Z_{n,I}(q_j)-Z_{n,I}(q)|}{Z_{n,I}(q_j)}
 \longrightarrow0.
\]
Half this $L^1$ distance is the total variation distance. This is a continuity statement at fixed $n,I$ on the complete nonnegative boundary domain. For $h\ge0$ and $c\ge0$, the
substitution $u=e^{-s}$ gives
\[
 \int_0^1u^c|\log u|^h\,\dd u
 =\int_0^\infty e^{-(c+1)s}s^h\,\dd s
 =\frac{\Gamma(h+1)}{(c+1)^{h+1}}<\infty.
\]
On $u\ge1$, $(\log u)^h\le C_h(1+u^{h+1})$, so $u^c(\log u)^he^{-\kappa c_*u^2/2}$ is integrable. Multiplying the majorant for the product density by $\sum_\alpha|\log u_\alpha|^h$
and applying Tonelli therefore gives a finite integral independent of $q$ in each fixed compact set. Division by the positive lower bound for $Z_{n,I}(q)$ proves the locally
uniform logarithmic moment bound.
\end{proof}

\subsubsection{Conditional minimizers and their derivatives}
\label{B:subsection:conditional-minimizers-and-derivatives}

The conditional Gaussian approximation is made on an interval separated from both matrix boundaries. We sum the coefficient errors on a smaller interval, where the influence of the
conditional boundary data has decayed.

Fix $A>20$ and henceforth use
\begin{equation}
 \begin{gathered}
 r_n=\lceil(\log n)^A\rceil,\quad R_n=n-r_n,\quad
 J_n=[r_n,R_n]\cap\mathbb Z,\quad b_n=\lceil r_n/4\rceil,\\
 M_n=n-r_n-b_n,\qquad s_n=r_n-R_*-1,\qquad
 J_n^-=[1,s_n]\cap\mathbb Z,\qquad\vartheta_n=r_n^{-1/8}.
 \end{gathered}
 \label{B:eq:general-potential-intervals}
\end{equation}
Increase $n_0$ so that these intervals exist, $s_n\ge1$ contains the exceptional coordinates, and $r_n>2/\beta+2R_*+2$ for $n\ge n_0$. The right boundary of $J_n^-$ is the site set
$[r_n-R_*,r_n-1]$. For $n<n_0$ we take the exterior of $J_n$ to be the whole coefficient row; its law will be sampled directly. Define
\begin{equation}
 \mathcal Q_n^\sharp=
 \left\{q>0:\max_{\gamma=(\sigma,k)\in\mathscr S(\partial_{R_*}J_n)}
       \frac{|q_\gamma-(z_n^\sharp)_\gamma|}{\phi(k/n)}
                   \le\vartheta_n\right\}.
 \label{B:eq:admissible-boundary-configurations}
\end{equation}
On $J_n$ the two conditional Hamiltonians coincide. Denote their minimizer by $m^q$.

\begin{proposition}\label{B:prop:conditional-minimizer-centering}
Uniformly for $q\in\mathcal Q_n^\sharp$ and $k\in J_n$,
\begin{equation}
 |m_{\sigma,k}^q-(z_n^\sharp)_{\sigma,k}|
 \le C\vartheta_n\left[\phi(r_n/n)e^{-c(k-r_n)}
                     +\phi(R_n/n)e^{-c(R_n-k)}\right].
 \label{B:eq:conditional-minimizer-localization}
\end{equation}
In particular,
\begin{equation}
 \begin{split}
 |m_{\sigma,k}^q/\phi(k/n)-1|
 +|\log(m_{\sigma,k}^q/\phi(k/n))|
 \le C\bigl[k^{-1}+e^{-c'r_n}
       +\vartheta_n(e^{-c'(k-r_n)}+e^{-c'(R_n-k)})\bigr].
 \end{split}
 \label{B:eq:conditional-minimizer-relative-centering}
\end{equation}
\end{proposition}
\begin{proof}
By the restriction property in Proposition~\ref{B:prop:conditional-minimizer-dependence}, $z_n^\sharp|_{J_n}$ is the conditional minimizer for its own boundary values.
Equation~\eqref{B:eq:conditional-minimizer-boundary-dependence} and the admissibility bound give
\[
 \begin{split}
 |m_{\sigma,k}^q-(z_n^\sharp)_{\sigma,k}|
 &\le C\sum_{\alpha\in\mathscr S(\partial_{R_*}J_n)}
       e^{-c|k-\iota(\alpha)|}
                    |q_\alpha-(z_n^\sharp)_\alpha|\\
 &\le C\vartheta_n\sum_{\alpha\in\mathscr S(\partial_{R_*}J_n)}
                  e^{-c|k-\iota(\alpha)|}\phi(\iota(\alpha)/n)\\
 &\le C\vartheta_n\left[
        \phi(r_n/n)e^{-c(k-r_n)}
                 +\phi(R_n/n)e^{-c(R_n-k)}\right].
 \end{split}
\]
For the last inequality, each boundary layer has a fixed number of sites, its indices differ from $r_n$ or $R_n$ by at most $R_*$, and $\phi(t)\asymp\sqrt t$ makes the profile
ratios on each layer uniformly bounded. This proves \eqref{B:eq:conditional-minimizer-localization}.

Divide by $\phi(k/n)\asymp\sqrt{k/n}$ and use \eqref{B:eq:reference-minimizer-interior-ratio}. The left boundary ratio is at most $C\sqrt{r_n/k}\le C$. For the right one,
$\sqrt{R_n/k}\le\sqrt{1+R_n-k}$, and
\[
 (1+j)^{1/2}e^{-cj}
 \le\left(\sup_{u\ge0}(1+u)^{1/2}e^{-cu/2}\right)e^{-cj/2}.
\]
This proves the relative estimate without its logarithm. Uniformly over $q\in\mathcal Q_n^\sharp$ and $k\in J_n$, its right side is at most $C(r_n^{-1}+e^{-c'r_n}+2\vartheta_n)$,
which tends to zero. For large $n$ the ratio lies in $[1/2,3/2]$ and
\[
 \left|\log\frac{m_{\sigma,k}^q}{\phi(k/n)}\right|
 =\left|\int_1^{m_{\sigma,k}^q/\phi(k/n)}\frac{\dd s}{s}\right|
 \le2\left|\frac{m_{\sigma,k}^q}{\phi(k/n)}-1\right|.
\]
This proves the logarithmic part with the same uniformity. On $r_n+b_n\le k\le R_n-b_n$, the boundary terms are bounded by $Ce^{-c''r_n}$.
\end{proof}

For $q\in\mathcal Q_n^\sharp$ let $\mathcal V_{n,q}=\{u:m_\alpha^q/2\le u_\alpha\le2m_\alpha^q\}$. The constants in the next estimate are uniform on this set and over
$q\in\mathcal Q_n^\sharp$.

\begin{lemma}\label{B:lem:conditional-scaled-derivatives}
On $\mathcal V_{n,q}$,
\begin{equation}
 cI\le\nabla^2H_{n,J_n}^{q}(u)\le CI,\qquad
 m_{\sigma,k}^q\asymp\sqrt{k/n},\qquad
 \sqrt\kappa m_{\sigma,k}^q\asymp\sqrt k.
 \label{B:eq:conditional-hessian-and-scales}
\end{equation}
For $h=3,4,5$, every nonzero mixed coordinate derivative of $\kappa H_{n,J_n}^{q}(m^q+x/\sqrt\kappa)$ is bounded by $C_h k_*^{1-h/2}$, where $k_*$ is any site among its derivative
indices.
\end{lemma}
\begin{proof}
Proposition~\ref{B:prop:conditional-minimizer-centering} gives the coordinate bounds. The polynomial Hessian is bounded below by $c_*I$ and has uniformly bounded entries and a
fixed range on the relevant box. Its operator norm is therefore bounded by its largest absolute row sum. The logarithmic Hessian is diagonal, nonnegative, and has entries at most
$C(k/n)/\phi(k/n)^2\le C$. This proves the Hessian bounds.

Every coefficient in the neighborhood of a nonzero mixed derivative, including the fixed boundary coordinates, is comparable to $\phi(k_*/n)$. A monomial of degree $2m$ has an
$h$th derivative bounded by $C_h\phi(k_*/n)^{2m-h}$ when $2m\ge h$, and is zero otherwise. Since $\phi$ is bounded above, this is at most $C_h\phi(k_*/n)^{2-h}$. For a logarithm
the derivative is zero unless all indices agree; its absolute value is $(h-1)!\alpha_\alpha^\sharp/u_\alpha^h \le C_h(k_*/n)\phi(k_*/n)^{-h}$. For derivative coordinates
$\alpha_1,\ldots,\alpha_h$, the chain rule is exactly
\[
 \partial_{x_{\alpha_1}}\cdots\partial_{x_{\alpha_h}}
 [\kappa H_{n,J_n}^{q}(m^q+x/\sqrt\kappa)]
 =\kappa^{1-h/2}
 (\partial_{\alpha_1}\cdots\partial_{\alpha_h}H_{n,J_n}^{q})
          (m^q+x/\sqrt\kappa).
\]
The coordinates lie in $\mathcal V_{n,q}$ as specified above. A nonzero local derivative has all its sites within a fixed distance of $k_*$, so their positive scales are
comparable. Both estimates are consequently at most
\begin{equation}
 C_h\kappa^{1-h/2}\phi(k_*/n)^{2-h}
                  \le C_{h,\beta}k_*^{1-h/2}.
 \label{B:eq:conditional-higher-derivative-bound}
\end{equation}
The number of monomials and the allowed site offsets are bounded in terms of $d,h$, so summation preserves the bound.
\end{proof}

\subsubsection{A logarithmic Sobolev inequality}
\label{B:subsection:reference-logarithmic-sobolev-inequality}

Uniform convexity gives concentration under the conditional reference laws. The logarithmic Sobolev estimate is recorded with its boundary dependence explicit, since it will be
used for both coordinate fluctuations and weighted logarithmic sums.

Write $\overline\mu_{n,I}^{q}$ for the conditional reference law on any nonempty interval $I$. This notation agrees with $\mu_{n,I}^{q}$ when $I$ misses $\mathscr S_-$. Write
$\overline{\mathbb E}_q$ for expectation under this reference measure. For a probability measure $\mu$ and a real function $f$, set
$\operatorname{Ent}_\mu f^2=\int f^2\log f^2\,\dd\mu -(\int f^2\,\dd\mu)\log\int f^2\,\dd\mu$, with $0\log0=0$.

\begin{lemma}\label{B:lem:reference-logarithmic-sobolev}
For every finite nonnegative boundary configuration $q$,
\begin{equation}
 \operatorname{Ent}_{\overline\mu_{n,I}^{q}}f^2
   \le\frac{2}{\kappa c_*}
                     \int|\nabla f|^2\,\dd\overline\mu_{n,I}^{q}.
 \label{B:eq:reference-logarithmic-sobolev}
\end{equation}
The inequality holds for real smooth functions and extends to their closure in the norm $\bigl(\int(|f|^2+|\nabla f|^2)\,\dd\overline\mu_{n,I}^{q}\bigr)^{1/2}$. For every real
one-Lipschitz function $F$,
\begin{equation}
 \log\overline{\mathbb E}_q e^{t(F-\overline{\mathbb E}_qF)}\le\frac{t^2}{2\kappa c_*},\qquad
 \overline\mu_{n,I}^{q}\{|F-\overline{\mathbb E}_qF|>u\}
                    \le2e^{-\kappa c_*u^2/2}.
 \label{B:eq:reference-lipschitz-concentration}
\end{equation}
\end{lemma}
\begin{proof}
For the fixed finite boundary configuration $q$, write the extended Hamiltonian as $h(z)=c_*|z|^2/2+F_1(z)$, where $F_1$ is proper, lower semicontinuous, and convex, including its
value $+\infty$ outside the coefficient domain. By the contraction theorem \cite{caffarelli2000monotonicity} (Theorem~11, pp.~559--561), including its convex support constraints,
there is a one-Lipschitz transport from $\gamma=N(0,(\kappa c_*)^{-1}I)$ to $\overline\mu_{n,I}^{q}$. The Gaussian logarithmic Sobolev inequality has constant $2/(\kappa c_*)$ by
scaling \cite{bakry1985diffusions} (Proposition~4, (13), Corollary~2, and the Ornstein--Uhlenbeck example, pp.~195 and 199--200). Writing $T$ for this transport, the chain rule for
Lipschitz maps gives, for bounded smooth $f$,
\[
 \begin{split}
 \operatorname{Ent}_{\overline\mu_{n,I}^{q}}f^2
 &=\operatorname{Ent}_\gamma(f\circ T)^2\\
 &\le\frac{2}{\kappa c_*}\int|\nabla(f\circ T)|^2\,\dd\gamma\\
 &\le\frac{2}{\kappa c_*}\int|\nabla f|^2\circ T\,\dd\gamma
 =\frac{2}{\kappa c_*}\int|\nabla f|^2\,\dd\overline\mu_{n,I}^{q}.
 \end{split}
\]
This is the contraction argument of Corollary~8 in \cite{caffarelli2000monotonicity} (p.~562). Smooth truncation in the values and smooth cutoffs extend it to every smooth function
of finite form norm. If $f_j\to f$ in that norm, then $f_j^2\to f^2$ in $L^1$; lower semicontinuity of entropy and convergence of the gradient norms give the stated extension to
the closure.

For bounded smooth one-Lipschitz $F$, apply this inequality to $e^{tF/2}$. With $\Lambda(t)=\log\overline{\mathbb E}_q e^{t(F-\overline{\mathbb E}_qF)}$, it gives
$t\Lambda'(t)-\Lambda(t)\le t^2/(2\kappa c_*)$. Integrating $(\Lambda(t)/t)'$ from zero proves $\Lambda(t)\le t^2/(2\kappa c_*)$ for $t>0$; apply the same argument to $-F$ for
$t<0$. A lower bound $F_1(x)\ge a+b\cdot x$, with fixed $a\in\mathbb R$ and a fixed vector $b$, gives an integrable Gaussian upper bound for the density. Hence every Lipschitz $F$
has all exponential moments, and Lipschitz extension, truncation, and mollification pass the bound to the stated class by dominated convergence. Exponential Markov with
$t=\kappa c_*u$, applied to $F$ and $-F$, gives the two tails in \eqref{B:eq:reference-lipschitz-concentration}.
\end{proof}

\subsubsection{Concentration about the conditional minimizer}
\label{B:subsection:concentration-about-conditional-minimizers}

Concentration is naturally centered at the mean, whereas the later Gaussian approximation is centered at the conditional minimizer. The first estimate below compares these centers;
the remaining estimates then control fluctuations about the minimizer.

\begin{proposition}\label{B:prop:reference-coordinate-concentration}
For each $M<\infty$, uniformly over $n,I$ and $\|q\|_\infty\le M$,
\begin{equation}
 \begin{split}
 \max_\alpha|\overline{\mathbb E}_q Z_\alpha-m_\alpha^{q,I}|&\le C_Mn^{-1/2},\\
 \overline\mu_{n,I}^{q}\{|Z_\alpha-m_\alpha^{q,I}|>t\}
       &\le2e^{-c_Mnt^2}\quad(t\ge C_Mn^{-1/2}),\\
 \|Z_\alpha-m_\alpha^{q,I}\|_{L^p(\overline\mu_{n,I}^{q})}
       &\le C_M\sqrt{p/n}\quad(p\ge1).
 \end{split}
 \label{B:eq:reference-coordinate-concentration}
\end{equation}
\end{proposition}
\begin{proof}
We first obtain a conditional energy estimate for an arbitrary interval $I'$ and finite nonnegative boundary configuration $Q$. Write $h=H_{n,I'}^{\sharp,Q}$ and $m=m^{Q,I'}$.
Integrating $\partial_\alpha[(z_\alpha-m_\alpha)e^{-\kappa h(z)}]$ gives
\begin{equation}
 1=\kappa\overline{\mathbb E}_Q[(Z_\alpha-m_\alpha)\partial_\alpha h(Z)],
 \qquad \overline{\mathbb E}_Q|Z-m|^2\le\frac{|\mathscr S(I')|}{\kappa c_*}.
 \label{B:eq:conditional-energy-identity}
\end{equation}
The boundary terms in this integration by parts vanish as follows. If $c_\alpha^+>0$, the boundary factor is a bounded multiple of $(z_\alpha-m_\alpha)z_\alpha^{c_\alpha^+}$ near
zero, and its derivative has integrable lowest power $z_\alpha^{c_\alpha^+-1}$. If $c_\alpha^+=0$, Lemma~\ref{B:lem:reference-minimizer-coordinates} gives $m_\alpha=0$, so the
factor has an additional $z_\alpha$ and vanishes at the face. The boundary term at infinity vanishes by quadratic growth; polynomial derivatives and the indicated powers are
integrable. Integrate first over a finite rectangle with positive lower bounds and then let its sides tend to the boundary and infinity. Fubini is justified by these same bounds.
Sum the first identity over $\alpha$ and use
\[
 (z-m)\cdot\nabla h(z)
 =(z-m)\cdot(\nabla h(z)-\nabla h(m))\ge c_*|z-m|^2
\]
to obtain the second estimate.

Fix now $I,q$ with $\|q\|_\infty\le M$, and set $m=m^{q,I}$. There is a constant $C_*$ bounding every conditional minimizer coordinate, even for unbounded finite boundary
configurations, by Lemma~\ref{B:lem:reference-minimizer-coordinates}. Equation~\eqref{B:eq:conditional-energy-identity} implies
$\overline{\mathbb E}_q|Z-m|\le\sqrt{2/(\beta c_*)}$, since $|\mathscr S(I)|\le2n$. Choose a fixed $M_2\ge M$ greater than $C_*+\sqrt{2/(\beta c_*)}+1$. The concentration
inequality for the one-Lipschitz function $|Z-m|$ gives
\begin{equation}
 \overline\mu_{n,I}^{q}\{|Z-m|>M_2-C_*\}\le2e^{-cn}.
 \label{B:eq:bounded-random-boundary-probability}
\end{equation}

Let $L>2R_*$ be an integer to be chosen below. For $k\in I$ set $I_k=I\cap[k-L,k+L]$. Condition on all coordinates outside $I_k$. Let $Q_k$ be its boundary data, comprising the
fixed entries of $q$ and the random entries of $Z$. Replace only the latter by the entries of $m$ to obtain $Q_k^m$. The restriction property gives $m^{Q_k^m,I_k}=m|_{I_k}$. Every
differing boundary coordinate is at distance at least $L$ from $k$; if $I_k$ is clipped by an end of $I$, the boundary beyond that end is fixed and agrees in $Q_k,Q_k^m$.
Conditionally, \eqref{B:eq:conditional-energy-identity} gives $\mathbb E[|Z_\alpha-m_\alpha^{Q_k,I_k}|\mid Q_k] \le C\sqrt{L/n}$ for $\iota(\alpha)=k$. On $\|Q_k\|_\infty\le M_2$,
the boundary estimate gives
\[
 |m_\alpha^{Q_k,I_k}-m_\alpha|
 \le C_{M_2}e^{-c_{M_2}L}
       \sum_{\gamma\in\mathscr S(\partial_{R_*}I_k)\cap\mathscr S(I)}
                         |Z_\gamma-m_\gamma|.
\]
There are at most $4R_*$ summands. On the complement, the difference of these two minimizing coordinates is at most $2C_*$. That complement is contained in the event in
\eqref{B:eq:bounded-random-boundary-probability}. With $A_{n,I,q}=\max_\alpha\overline{\mathbb E}_q|Z_\alpha-m_\alpha|$, which is finite by
\eqref{B:eq:conditional-energy-identity}, we obtain
\[
 A_{n,I,q}\le C\sqrt{L/n}
       +4R_*C_{M_2}e^{-c_{M_2}L}A_{n,I,q}+Ce^{-cn}.
\]
Choose a fixed integer $L>2R_*$ for which $4R_*C_{M_2}e^{-c_{M_2}L}\le1/2$. The inequality becomes
\[
 A_{n,I,q}\le2C\sqrt{L/n}+2Ce^{-cn}\le C_Mn^{-1/2},
\]
since $\sup_{n\ge1}\sqrt n e^{-cn}<\infty$. This choice is independent of $n,I,q$. If $I_k=I$, then $Q_k=Q_k^m=q$ and the two minimizers coincide. The conditional energy bound
directly gives $\overline{\mathbb E}_q|Z_\alpha-m_\alpha| \le\sqrt{|\mathscr S(I)|/(\kappa c_*)} \le C\sqrt{L/n}$, since $I\subset[k-L,k+L]$. There is no term involving the
boundary difference in this case. Finally $|\overline{\mathbb E}_qZ_\alpha-m_\alpha| \le\overline{\mathbb E}_q|Z_\alpha-m_\alpha|\le A_{n,I,q}$ proves the stated mean displacement.

For a coordinate function, \eqref{B:eq:reference-lipschitz-concentration} and $|\overline{\mathbb E}_qZ_\alpha-m_\alpha|\le C_Mn^{-1/2}$ give
\[
 \overline\mu_{n,I}^q\{|Z_\alpha-m_\alpha|>t\}
 \le 2\exp[-\kappa c_*(t-C_Mn^{-1/2})^2/2].
\]
For $t\ge2C_Mn^{-1/2}$ this proves the required tail bound. If $t_0=2C_Mn^{-1/2}$, integration yields
\[
 \overline{\mathbb E}_q|Z_\alpha-m_\alpha|^p
 \le t_0^p+2p\int_0^\infty t^{p-1}e^{-c_Mnt^2}\,\dd t
 =t_0^p+2(c_Mn)^{-p/2}\Gamma(p/2+1).
\]
The bound $\Gamma(p/2+1)^{1/p}\le C\sqrt p$ proves the moment estimate, with an adjusted constant for $1\le p\le2$.
\end{proof}

\subsubsection{Probability estimates at the boundaries}
\label{B:subsection:concentration-boundary-configurations-and-terminal-estimates}

We now estimate the probability that the boundary values are admissible under the full laws, and separately bound the coordinates near the right matrix boundary.

Let $Q_n$ denote the boundary coordinates of $J_n$ under either full law, and put $\mathcal G_n=\{Q_n\in\mathcal Q_n^\sharp\}$. Let $\ell_n=\log(\mathrm e+n)$ and
$\mathscr T_n=\{(\sigma,k)\in\mathscr S_n:k\ge M_n-R_*\}$. All constants below are for the fixed parameters of \eqref{B:eq:general-potential-intervals}.

\begin{theorem}\label{B:thm:boundary-and-terminal-probabilities}
For all sufficiently large $n$,
\begin{equation}
 \overline{\mathbb P}_n(\mathcal G_n^c)\le Ce^{-cr_n^{3/4}},\qquad
 \mathbb P_n(\mathcal G_n^c)\le C_\beta e^{-c\gamma_\beta r_n^{3/4}}.
 \label{B:eq:boundary-configuration-probability}
\end{equation}
For $L$ above a fixed positive constant,
\begin{equation}
 \begin{split}
 \overline{\mathbb P}_n\{\max_{\alpha\in\mathscr T_n}
 |Z_\alpha-(z_n^\sharp)_\alpha|>L\sqrt{\ell_n/n}\}
       &\le Cr_n e^{-cL^2\ell_n},\\
 \mathbb P_n\{\max_{\alpha\in\mathscr T_n}
 |Z_\alpha-(z_n^\sharp)_\alpha|>L\sqrt{\ell_n/n}\}
       &\le C_\beta r_n^{\gamma_\beta}e^{-c\gamma_\beta L^2\ell_n}.
 \end{split}
 \label{B:eq:terminal-coordinate-probability}
\end{equation}
The original boundary failures are summable; the original right boundary failures are summable for fixed sufficiently large $L$.
\end{theorem}
\begin{proof}
The boundary of $J_n$ has at most $4R_*$ coordinates, each with $k\ge r_n-R_*$. For the full reference law the exterior is empty, so
Proposition~\ref{B:prop:reference-coordinate-concentration} applies with $m=z_n^\sharp$. The threshold for the relative error is $\vartheta_n\phi(k/n)$ and
\[
 \sqrt n\,\vartheta_n\phi(k/n)\ge cr_n^{3/8}\longrightarrow\infty,
 \qquad n\vartheta_n^2\phi(k/n)^2\ge cr_n^{3/4}.
\]
Apply Proposition~\ref{B:prop:reference-coordinate-concentration} to the full reference law and sum over these coordinates. This proves the first bound in
\eqref{B:eq:boundary-configuration-probability}. Equation~\eqref{B:eq:coefficient-measure-comparison} gives the second. Since $r_n^{3/4}/\log n\to\infty$, choose $n_0$ so that the
displayed upper bound is at most $Cn^{-2}$ for $n\ge n_0$.

The set $\mathscr T_n$ has cardinality at most $Cr_n$. Apply the same concentration bound at threshold $L\sqrt{\ell_n/n}$ and take a union bound. This proves the bound under the
reference measure in \eqref{B:eq:terminal-coordinate-probability}. For this union event $A_{n,L}$, H\"older gives
\[
 \mathbb P_n(A_{n,L})
 =\overline{\mathbb E}_n[L_n\mathbf1_{A_{n,L}}]
 \le C_\beta\bigl(Cr_n e^{-cL^2\ell_n}\bigr)^{\gamma_\beta}
 \le C_\beta r_n^{\gamma_\beta}e^{-c\gamma_\beta L^2\ell_n}.
\]
The likelihood exponent is applied to the whole union probability. Choose $L$ so that $c\gamma_\beta L^2>4$. Its right side is at most $C(\log n)^{A\gamma_\beta}n^{-4}$, a summable
sequence.
\end{proof}

For later products we record the logarithmic consequences. For $k\ge n/2$, the coordinate bounds for the full minimizer give
\begin{equation}
 c\le (z_n^\sharp)_{\sigma,k}\le C.
 \label{B:eq:terminal-minimizer-bounds}
\end{equation}
For any jointly defined arrays with the original law at each size, the preceding theorem and Borel--Cantelli give a finite $n_0(\omega)$ such that, for every $n\ge n_0(\omega)$,
almost surely,
\begin{equation}
 \max_{\alpha\in\mathscr T_n}
       \left|\log\frac{Z_\alpha}{(z_n^\sharp)_\alpha}\right|
                  \le C L\sqrt{\ell_n/n}.
 \label{B:eq:terminal-logarithmic-deviation}
\end{equation}
In fact, for large $n$ the coordinate deviation is at most $c/2$, and the mean-value formula bounds the logarithmic deviation by $2|Z_\alpha-(z_n^\sharp)_\alpha|/c$. Since at most
$Cr_n$ indices are involved,
\begin{equation}
 \sup_{M_n-R_*\le i<j\le n}
 \left|\sum_{k=i}^{j-1}\left[
 \log\frac{Y_{n,k}}{(z_n^\sharp)_{y,k}}
 -\log\frac{X_{n,k}}{(z_n^\sharp)_{x,k}}\right]\right|
            \le Cr_n\sqrt{\ell_n/n}\longrightarrow0.
 \label{B:eq:terminal-product-fluctuation}
\end{equation}
The deterministic product is uniformly bounded above and below. Indeed \eqref{B:eq:first-correction-definition} and \eqref{B:eq:reference-minimizer-absolute-error}, together with
\eqref{B:eq:terminal-minimizer-bounds}, give
\begin{equation}
 \left|\log\frac{(z_n^\sharp)_{y,k}}{(z_n^\sharp)_{x,k}}\right|
 \le C/n+Ce^{-c(n-k)}+Cn^{-1/2}e^{-c(k-K)}\quad(n/2\le k<n).
 \label{B:eq:terminal-deterministic-product-bound}
\end{equation}
The sum of the first term over this interval near the right boundary is $O(r_n/n)$; the other two sums are bounded geometric series. This proves the product bound. The coordinate
estimates include $X_{n,n}$, whereas the products end at $k=j-1<n$.

\subsubsection{Comparison with a product distribution near the first coordinate}
\label{B:subsection:initial-product-distribution-comparison}

Near the first coordinate, we split off the quadratic part of the potential to obtain a product of chi laws. The remaining interaction appears in an explicit density ratio.

Retain $c_0=2V'(0)>0$ and define $\bar\nu_{x,k}=\max\{\beta(k+a),1\}$ and $\bar\nu_{y,k}=\max\{\beta k,1\}$. Let $\bar\nu_n^-$ be the product probability law on $J_n^-$ under which
\begin{equation}
 \sqrt{\kappa c_0}\,X_k\sim\chi_{\bar\nu_{x,k}},\qquad
 \sqrt{\kappa c_0}\,Y_k\sim\chi_{\bar\nu_{y,k}}
                  \quad(1\le k\le s_n)
 \label{B:eq:initial-product-distribution}
\end{equation}
independently. Its exponents are $\bar\nu_{\sigma,k}-1=c_{\sigma,k}^+$. Splitting off the quadratic term $g_1\operatorname{Tr}BB^{\mathsf T}$ in the polynomial gives the exact
conditional reference density
\begin{equation}
 \overline\mu_{n,-}^{q^-}(\dd u)
 =(Z_{n,q^-}^-)^{-1}e^{-\kappa\mathcal R_{n,q^-}(u)}\bar\nu_n^-(\dd u).
 \label{B:eq:initial-conditional-density}
\end{equation}
Here $\mathcal R_{n,q^-}$ is the sum of the original monomials of degree at least four that contain an interior coefficient, with the boundary entries set equal to $q^-$ and terms
independent of $u$ subtracted.

Fix $C_0<\infty$ and restrict to $\|q^-\|_\infty\le C_0\sqrt{r_n/n}$. For a sufficiently large fixed $L$ set
\begin{equation}
 \delta_n=Cs_n e^{-cL^2s_n^2},\qquad
 \epsilon_n=C_L\sum_{m=2}^d\frac{s_n^{2m+1}}{n^{m-1}}.
 \label{B:eq:initial-comparison-errors}
\end{equation}
The sum is zero for a linear potential. Both sequences tend to zero.

\begin{proposition}\label{B:prop:initial-conditional-comparison}
Uniformly over the stated boundary data and every Borel set $F$ in the initial coordinates,
\begin{equation}
 \overline\mu_{n,-}^{q^-}(F)
 \le\frac{e^{2\epsilon_n}}{1-\delta_n}\bar\nu_n^-(F)+\delta_n
 \le3\bar\nu_n^-(F)+\delta_n
 \label{B:eq:initial-probability-comparison}
\end{equation}
for all sufficiently large $n$.
\end{proposition}
\begin{proof}
Let $\mathcal T_{n,L}=\{\max_{\mathscr S(J_n^-)}Z_\alpha \le Ls_n/\sqrt n\}$. The conditional minimizer is at most $C\sqrt{s_n/n}$ by
Lemma~\ref{B:lem:reference-minimizer-coordinates}. For large $n$, the distance to the threshold is at least $Ls_n/(2\sqrt n)$. The boundary configurations belong to a fixed bounded
box. Proposition~\ref{B:prop:reference-coordinate-concentration} and the union bound over at most $2s_n$ coordinates yield
$\overline\mu_{n,-}^{q^-}(\mathcal T_{n,L}^c)\le\delta_n$. For the product law, $\bar\nu_{\sigma,k}\le C_{\beta,a}s_n$ and $\mathbb E e^{\chi_\nu^2/4}=2^{\nu/2}$. Therefore
\[
 \mathbb P\{\chi_{\bar\nu_{\sigma,k}}>L\sqrt{\beta c_0}\,s_n\}
 \le\exp[(\bar\nu_{\sigma,k}/2)\log2-\beta c_0L^2s_n^2/4]
 \le C e^{-cL^2s_n^2}.
\]
Summation gives $\bar\nu_n^-(\mathcal T_{n,L}^c)\le\delta_n$ after adjusting its constants.

On $\mathcal T_{n,L}$, every coordinate in a monomial of $\mathcal R_{n,q^-}$, including each fixed boundary factor, is at most $C_Ls_n/\sqrt n$, since $r_n=s_n+R_*+1$. The walk
expansion from Lemma~\ref{B:lem:coefficient-polynomial-convexity} bounds the number of terms of degree $2m$ meeting $J_n^-$ by $C_ms_n$: their starting vertices lie within a fixed
distance of these sites and each walk has at most $2^{2m}$ choices. Thus
\[
 \sup_{\mathcal T_{n,L}}|\kappa\mathcal R_{n,q^-}|
 \le C_L\sum_{m=2}^d n s_n(s_n/\sqrt n)^{2m}=\epsilon_n.
\]
The product law gives $Z_{n,q^-}^-\ge e^{-\epsilon_n}(1-\delta_n)$. The already proved conditional cutoff probability also gives
\[
 1-\delta_n\le\overline\mu_{n,-}^{q^-}(\mathcal T_{n,L})
   =(Z_{n,q^-}^-)^{-1}\int_{\mathcal T_{n,L}}
               e^{-\kappa\mathcal R_{n,q^-}}\,\dd\bar\nu_n^-
   \le e^{\epsilon_n}/Z_{n,q^-}^-.
\]
The two normalizer bounds are
\[
 e^{-\epsilon_n}(1-\delta_n)\le Z_{n,q^-}^-
                         \le e^{\epsilon_n}/(1-\delta_n).
\]
For every Borel $F$, splitting at $\mathcal T_{n,L}$ therefore gives
\[
 \overline\mu_{n,-}^{q^-}(F)
 \le\frac{e^{\epsilon_n}}{Z_{n,q^-}^-}
           \bar\nu_n^-(F\cap\mathcal T_{n,L})
       +\overline\mu_{n,-}^{q^-}(\mathcal T_{n,L}^c)
 \le\frac{e^{2\epsilon_n}}{1-\delta_n}\bar\nu_n^-(F)+\delta_n.
\]
The proof uses the absolute bound on $\kappa\mathcal R_{n,q^-}$ inside $\mathcal T_{n,L}$ and its separately established probability bound outside that set. It imposes no sign
condition on the coefficients $g_m$ of the terms of higher degree. The coefficient in the second inequality tends to one, since every fixed power of $s_n=O((\log n)^A)$ is
dominated by a positive power of $n$.
\end{proof}

\subsubsection{Tightness of the initial coefficients}
\label{B:subsection:tightness-of-initial-coefficients}

The inverse kernel contains both individual initial coefficients and their cumulative ratios. We use the product comparison to control them jointly.

For $m\ge1$ put $\varepsilon_m=2^{-m-4}$. For $n\ge n_0$ define
\begin{equation}
 \begin{split}
 D_{n;i,j}&=\left(\frac{z_{n,i}}{z_{n,j}}\right)^{a/2+1/4}
                \left(\frac{\phi(j/n)}{\phi(i/n)}\right)^{1/2},\\
 \mathcal E_{n,m}^-&=1+
 \max_{1\le k\le s_n}z_{n,k}^{\varepsilon_m}
 \left(\frac{X_k}{\phi(k/n)}+\frac{\phi(k/n)}{X_k}
       +\frac{Y_k}{\phi(k/n)}+\frac{\phi(k/n)}{Y_k}\right)\\
 &\quad+\max_{1\le i\le j\le s_n+1}(z_{n,i}z_{n,j})^{\varepsilon_m}
 \left(\frac{\prod_{k=i}^{j-1}Y_k/X_k}{D_{n;i,j}}
             +\frac{D_{n;i,j}}{\prod_{k=i}^{j-1}Y_k/X_k}\right).
 \end{split}
 \label{B:eq:initial-coordinate-functionals}
\end{equation}
The products end at $s_n$ and thus depend only on $J_n^-$. Empty products are one. For $n<n_0$ set these functionals equal to one.

\begin{lemma}\label{B:lem:initial-product-tightness}
There are jointly defined random vectors with laws $\bar\nu_n^-$ for which $\sup_n\mathcal E_{n,m}^-<\infty$ almost surely for every $m$. In particular, their laws are uniformly
tight.
\end{lemma}
\begin{proof}
Take one independent sequence of chi pairs with degrees $\bar\nu_{x,k},\bar\nu_{y,k}$ and use its first $s_n$ pairs, divided by $\sqrt{\kappa c_0}$, at size $n$. Only finitely many
of these degrees differ from $\beta(k+a),\beta k$. Lemma~\ref{B:lem:laguerre-initial-products}, including the case of finitely many changed degrees, gives, on one event,
\[
 c_\omega\sqrt k\le\chi_{x,k},\chi_{y,k}\le C_\omega\sqrt k,
 \qquad
 C_{\omega,\eta}^{-1}(ij)^{-\eta}
 \le\frac{\prod_{k=i}^{j-1}\chi_{y,k}/\chi_{x,k}}{(i/j)^{a/2}}
 \le C_{\omega,\eta}(ij)^\eta.
\]
The product estimate follows from the slightly stronger bound with $j^\eta$ in that lemma, since $i\ge1$. The functions $\theta(t)/t$ and $\phi(t)/\sqrt t$ are bounded above and
below by positive constants. Hence
\begin{equation}
 \sqrt{\kappa c_0}\,\phi(k/n)\asymp\sqrt k,\qquad
 D_{n;i,j}\asymp(i/j)^{a/2}.
 \label{B:eq:initial-profile-product-comparison}
\end{equation}
The first maximum in \eqref{B:eq:initial-coordinate-functionals} is therefore uniformly bounded by a finite random constant. For the second, choose $\eta=\varepsilon_m/4$. Its two
summands are at most
\[
 C_{\omega,m}n^{-2\varepsilon_m}
                  (s_n+1)^{2(\varepsilon_m+\eta)}.
\]
For large $n$, its logarithm is at most
\[
 -2\varepsilon_m\log n
   +2A(\varepsilon_m+\eta)\log\log n+C_m\longrightarrow-\infty.
\]
The finitely many other factors are bounded. Intersect the countably many events for the chosen $\eta=\varepsilon_m/4$. On this event put $H_m=\sup_n\mathcal E_{n,m}^-$; it is
measurable and finite, since it is a countable supremum. For every $L>0$,
\[
 \sup_n\bar\nu_n^-(\mathcal E_{n,m}^->L)
 \le\mathbb P(H_m>L)\longrightarrow0.
\]
This proves uniform tightness from the constructed common sequence.
\end{proof}

\begin{proposition}\label{B:prop:original-initial-tightness}
For every fixed $m$,
\begin{equation}
 \lim_{L\to\infty}\sup_n\overline{\mathbb P}_n(\mathcal E_{n,m}^->L)=0,
 \qquad
 \lim_{L\to\infty}\sup_n\mathbb P_n(\mathcal E_{n,m}^->L)=0.
 \label{B:eq:initial-coordinate-tightness}
\end{equation}
\end{proposition}
\begin{proof}
On $\mathcal G_n$, the boundary data $q^-$ for $J_n^-$ satisfy $\|q^-\|_\infty\le C\sqrt{r_n/n}$ by \eqref{B:eq:full-reference-minimizer-upper} and
\eqref{B:eq:admissible-boundary-configurations}. The event $\mathcal G_n$ is measurable with respect to the exterior of $J_n^-$, because both boundary layers of $J_n$ lie in that
exterior. Condition the reference full law on this exterior and apply Proposition~\ref{B:prop:initial-conditional-comparison}. After integration,
\begin{equation}
 \overline{\mathbb P}_n(\mathcal E_{n,m}^->L)
 \le\overline{\mathbb P}_n(\mathcal G_n^c)
             +3\bar\nu_n^-(\mathcal E_{n,m}^->L)+\delta_n.
 \label{B:eq:initial-tightness-comparison}
\end{equation}
Given $\epsilon>0$, choose $N$ so that the first and last terms have sum at most $\epsilon$ for all $n\ge N$. Then
\[
 \sup_{n\ge N}\overline{\mathbb P}_n(\mathcal E_{n,m}^->L)
 \le\epsilon+3\sup_n\bar\nu_n^-(\mathcal E_{n,m}^->L).
\]
Lemma~\ref{B:lem:initial-product-tightness} makes the second term tend to zero as $L\to\infty$. For each of the finitely many $n<N$, all coefficients are positive and finite almost
surely, so the functional is finite and its tail tends to zero. Letting $L\to\infty$ and then $\epsilon\downarrow0$ proves tightness under the reference measure. The uniform
measure comparison \eqref{B:eq:coefficient-measure-comparison} gives
\[
 \sup_n\mathbb P_n(\mathcal E_{n,m}^->L)
 \le C_\beta\left[
     \sup_n\overline{\mathbb P}_n(\mathcal E_{n,m}^->L)
                   \right]^{\gamma_\beta}\longrightarrow0.
\]
Here $\gamma_\beta>0$, so the established reference tightness proves tightness under the original measure.
\end{proof}

\subsubsection{The joint distribution of the exterior coordinates}
\label{B:subsection:joint-realization-of-exterior-coordinates}

Bounds in distribution for each size do not by themselves give one almost surely bounded array. The exterior coordinates are therefore realized jointly across sizes, preserving
their joint distribution at each matrix size and keeping the entire exterior sequence independent of the Brownian path.

Let $\lambda_n^{\rm ext}$ be the original joint law of all coordinates outside $J_n$, using the whole row for $n<n_0$. Each $\mathcal E_{n,m}^-$ is a Borel function of these
exterior coordinates.

\begin{theorem}\label{B:thm:original-exterior-realization}
There are jointly defined exterior vectors $E_n$ with laws $\lambda_n^{\rm ext}$, independent of $W$ as a sequence, for which
\begin{equation}
 \sup_n\mathcal E_{n,m}^-(E_n)<\infty\quad\text{almost surely for every }m.
 \label{B:eq:joint-exterior-initial-bounds}
\end{equation}
The boundary configurations belong to $\mathcal Q_n^\sharp$ for all sufficiently large $n$, almost surely.
\end{theorem}
\begin{proof}
Write $F_{n,m}=\mathcal E_{n,m}^-$. Uniform tightness permits constants $A_m\ge1$ with $\sup_n\lambda_n^{\rm ext}(F_{n,m}>A_m)\le2^{-m}$. Set $F_n^\sharp=\sup_mF_{n,m}/A_m$. For
$L\ge1$ the union bound gives
\[
 \sup_n\lambda_n^{\rm ext}(F_n^\sharp>L)
 \le\sum_{m\ge1}\sup_n\lambda_n^{\rm ext}(F_{n,m}>LA_m)
                    \longrightarrow0.
\]
For an explicit uniform limit, given $\epsilon>0$, choose $m_0$ with $\sum_{m>m_0}2^{-m}<\epsilon/2$. Uniform tightness of the first $m_0$ functionals gives $L$ for which the sum
of their tail bounds is below $\epsilon/2$. This proves the displayed limit. For fixed $n$, continuity from above gives $\lambda_n^{\rm ext}(F_n^\sharp=\infty)=0$. Set $F_n^\sharp$
equal to one on this Borel null set for the subsequent disintegration. All inequalities involving its original supremum remain valid outside that set.

Let $H_n$ be the distribution function of $F_n^\sharp$, and let $h_n(u)=\inf\{t:H_n(t)\ge u\}$ for $0<u<1$. Uniform tightness implies that, for each such $u$, there is an $L(u)$
with $H_n(L(u))>u$ for every $n$. Consequently
\begin{equation}
 \sup_n h_n(u)\le L(u)<\infty.
 \label{B:eq:uniform-quantile-bound}
\end{equation}
The generalized inverse satisfies $h_n(u)\le t$ exactly when $u\le H_n(t)$, by right continuity and monotonicity of $H_n$. It is a Borel function, and for a uniform $U_0$,
\[
 \mathbb P(h_n(U_0)\le t)=H_n(t).
\]
This includes atoms of the scalar distribution. The exterior spaces are finite-dimensional Borel spaces. Disintegration with respect to $F_n^\sharp$ gives a probability kernel
$K_n(t,\dd e)$ with
\[
 \lambda_n^{\rm ext}(\dd e)=\int K_n(t,\dd e)\,\dd H_n(t),\qquad
 K_n(t,\{e:F_n^\sharp(e)=t\})=1
                          \quad\text{for }\dd H_n\text{-almost every }t.
\]
We use the disintegration and randomization theorems for Borel probability kernels \cite{kallenberg2021foundations} (Lemma~4.22 and Theorem~8.5). On a space with independent
uniform variables $U_0,U_1,\ldots$, use $U_n$ to sample $E_n$ from $K_n(h_n(U_0),\dd e)$. For any bounded Borel $f$,
\[
 \mathbb Ef(E_n)=\int_0^1\int f(e)K_n(h_n(u),\dd e)\,\dd u
             =\int f(e)\lambda_n^{\rm ext}(\dd e),
\]
because $h_n(U_0)$ has distribution $H_n$. The support property and the null set convention give
\[
 \mathbb P(F_n^\sharp(E_n)=h_n(U_0))=1,
 \qquad
 \mathbb P\!\left(\bigcap_{n\ge1}
       \{F_n^\sharp(E_n)=h_n(U_0)\}\right)=1.
\]
Intersect also with the complements of the original null sets on which $F_n^\sharp$ was changed; each $E_n$ has its original marginal. On the resulting event
\[
 F_{n,m}(E_n)\le A_mF_n^\sharp(E_n)=A_mh_n(U_0),\qquad
 \sup_nF_{n,m}(E_n)\le A_m\sup_nh_n(U_0)<\infty.
\]
Take the product of the space carrying the entire sequence $(U_0,U_1,\ldots)$ with Wiener space. For bounded Borel functions $F$ of the exterior sequence and $G$ of the Brownian
path, Fubini gives
\[
 \mathbb E[F((E_n)_{n\ge1})G(W)]
   =\mathbb E F((E_n)_{n\ge1})\,\mathbb E G(W).
\]
Thus the complete exterior sequence is independent of $W$. The shared $U_0$ determines the scalar bounds across sizes; each $U_n$ supplies the remaining conditional exterior
coordinates. The probability of each boundary failure remains its original marginal probability. The summability in Theorem~\ref{B:thm:boundary-and-terminal-probabilities} and
Borel--Cantelli give a finite $n_0(\omega)$ beyond which all the stated bounds hold.
\end{proof}

For every $\eta>0$, choose $m$ with $\varepsilon_m\le\eta$. The definition of $\mathcal E_{n,m}^-$ and \eqref{B:eq:joint-exterior-initial-bounds} give
\begin{equation}
 \begin{gathered}
 C_{\omega,\eta}^{-1}\phi(k/n)z_{n,k}^{\eta}
       \le X_{n,k},Y_{n,k}
       \le C_{\omega,\eta}\phi(k/n)z_{n,k}^{-\eta}\quad(k\le s_n),\\
 C_{\omega,\eta}^{-1}D_{n;i,j}(z_{n,i}z_{n,j})^\eta
       \le\prod_{k=i}^{j-1}Y_{n,k}/X_{n,k}
       \le C_{\omega,\eta}D_{n;i,j}(z_{n,i}z_{n,j})^{-\eta}
                         \quad(1\le i\le j\le s_n+1).
 \end{gathered}
 \label{B:eq:original-initial-coordinate-bounds}
\end{equation}
Indeed, each nonnegative term of the functional is bounded by its supremum, and $z^{\varepsilon_m}\ge z^\eta$ for
$0<z\le1$.

\subsubsection{Extension by conditional probability kernels}
\label{B:subsection:extension-by-conditional-probability-kernels}

The exterior construction preserves all correlations between the initial coordinates, the boundary data, and the exterior coordinates in each row. Write $c_n(e)$ for the map
selecting the boundary configuration of $J_n$ from $e$.

\par\begingroup
\dimen0=\pagegoal\advance\dimen0 by -\pagetotal \ifdim\dimen0<8\baselineskip\newpage\fi
\endgroup
\begin{lemma}\label{B:lem:conditional-filling-and-costs}
If $\mathcal L(Z_{J_n}\mid E_n=e)=\mu_{n,J_n}^{c_n(e)}$, the full coefficient vector has law $\mathbb P_n$. Its deviations near the right boundary satisfy the almost sure
conclusions of Theorem~\ref{B:thm:boundary-and-terminal-probabilities}. If Borel plans $\pi_{n,q}(\dd g,\dd u)$ have $u$-marginal $\mu_{n,J_n}^q$ and specify the conditional law of
$(G_n,Z_{J_n})$ given $E_n$, then every nonnegative Borel cost satisfies
\begin{equation}
 \begin{split}
 &\mathbb E[\mathbf1_{\mathcal G_n}C_n(c_n(E_n),G_n,Z_{J_n})]\\
 &\quad=\int_{c_n(e)\in\mathcal Q_n^\sharp}\lambda_n^{\rm ext}(\dd e)
           \int C_n(c_n(e),g,u)\pi_{n,c_n(e)}(\dd g,\dd u)\\
 &\quad\le\sup_{q\in\mathcal Q_n^\sharp}\int C_n(q,g,u)\pi_{n,q}(\dd g,\dd u).
 \end{split}
 \label{B:eq:initial-coefficients-conditional-risk}
\end{equation}
\end{lemma}
\begin{proof}
For a nonnegative Borel function $F$ of the full coefficient vector, the specified conditional law and \eqref{B:eq:exterior-disintegration-identity} give
\begin{equation}
 \mathbb EF(E_n,Z_{J_n})
 =\int\lambda_n^{\rm ext}(\dd e)\int F(e,u)\mu_{n,J_n}^{c_n(e)}(\dd u)
 =\int F\,\dd\mathbb P_n.
 \label{B:eq:original-row-filling}
\end{equation}
Thus the failure probabilities near the right boundary are those in \eqref{B:eq:terminal-coordinate-probability}; their sum is finite for the fixed choice of $L$. Borel--Cantelli
gives the bound at the right boundary for every $n\ge n_0(\omega)$, including when the set of indices near the right boundary contains some filled interior coordinates. The
finitely many rows below $n_0$ are sampled with their original full laws. For the cost identity, first take a bounded nonnegative Borel $C_n$. Since $\mathbf1_{\mathcal G_n}$
depends only on $E_n$, the specified conditional law gives
\[
 \mathbb E[\mathbf1_{\mathcal G_n}C_n\mid E_n=e]
 =\mathbf1_{\{c_n(e)\in\mathcal Q_n^\sharp\}}
    \int C_n(c_n(e),g,u)\,\pi_{n,c_n(e)}(\dd g,\dd u).
\]
Integration against $\lambda_n^{\rm ext}$ proves the equality in \eqref{B:eq:initial-coefficients-conditional-risk}. This probability measure has total mass one, so the integral is
bounded by the stated supremum. For arbitrary nonnegative costs, apply the identity to $C_n\wedge l$ and let $l\to\infty$ by monotone convergence.
\end{proof}

If the last supremum in \eqref{B:eq:initial-coefficients-conditional-risk} is bounded by deterministic $a_n\ge0$ with $\sum_na_n<\infty$, then
\[
 \mathbb E\sum_n\mathbf1_{\mathcal G_n}C_n
   =\sum_n\mathbb E[\mathbf1_{\mathcal G_n}C_n]
   \le\sum_na_n<\infty.
\]
The nonnegative series is therefore finite almost surely. Boundary failures occur only finitely often. Intersecting with the events that each $C_n$ is finite, when this additional
hypothesis holds, also makes the sum over those exceptional indices finite.

Finally suppose a joint law $\Xi_n(\dd w,\dd e,\dd g)$ of a Brownian path, the original exterior, and a Gaussian vector has already been constructed. If $K_{n,q}(g,\dd u)$ is a
Borel probability kernel, then for every bounded Borel $F$ of the first three variables,
\[
 \int F(w,e,g)\,\Xi_n(\dd w,\dd e,\dd g)
                     K_{n,c_n(e)}(g,\dd u)
       =\int F\,\dd\Xi_n.
\]
Thus adjoining these coordinates preserves the existing joint law. When its conditional $(g,u)$ distribution is $\pi_{n,c_n(e)}$, \eqref{B:eq:original-row-filling} and
\eqref{B:eq:initial-coefficients-conditional-risk} apply to that extension. We distinguish the conditional distribution used for sampling from the marginal identity needed to
recover the coefficient law. Under the extension, the conditional law of $u$ given $(w,e,g)$ is $K_{n,c_n(e)}(g,\dd u)$. If the conditional law of $g$ given $e$ is
$\gamma_{n,c_n(e)}$, the required coefficient marginal is the integral identity
\[
 \int\gamma_{n,q}(\dd g)K_{n,q}(g,\dd u)=\mu_{n,J_n}^{q}(\dd u).
\]
This identity averages over $g$ at fixed exterior. Conditional on the Brownian path, the coefficient law need not equal
$\mu_{n,J_n}^{q}$.
Sections~\ref{B:sec:transport-of-logarithmic-ratios-and-construction-of-the-coupling}--\ref{B:sec:construction-of-the-coupling-for-all-matrix-sizes} construct the required Gaussian
vectors and sampling kernels. The estimates above provide their original exterior distribution and the bounds \eqref{B:eq:original-initial-coordinate-bounds}.

\subsection{Conditional estimates on the interior interval}
\label{B:sec:conditional-estimates-on-the-interior-interval}

The conditional law on $J_n$ from Section~\ref{B:sec:general-potentials-and-convex-reference-measures} now provides the setting for the logarithmic estimates. Moments and weighted
square sums are first controlled under that original law. Expanding the means also requires covariance estimates, which are obtained after a smooth extension in the coordinates
$x=\sqrt\kappa(Z-m^q)$. Agreement of the potentials on a specified set then transfers the estimates back; expectations under the two distributions are kept distinct throughout.

An individual logarithmic mean contains a term of order $k^{-1}$, whose sum would exceed the error allowed in the inverse comparison. Reflection at constant coefficients exchanges
the $x$ and $y$ coordinates and makes their leading logarithmic means equal. Comparison of the finite inverse Hessian with the constant coefficient inverse leaves an error of order
$k^{-2}$, together with a term exponentially decreasing from each conditional boundary. The sum of the means therefore includes a contribution from each boundary.

\subsubsection{The conditional distribution and its Hessian}
\label{B:subsection:conditional-distribution-and-hessian}

The interior estimates are conditional throughout. Fixing the boundary values determines the minimizer and Hessian; the constants must remain uniform over the admissible boundary
configurations when those values are later sampled.

Fix $V$, $a>0$, $\beta>0$, and $A_0>20$. Use the intervals of Section~\ref{B:sec:general-potentials-and-convex-reference-measures}, writing
\begin{equation}
 \begin{gathered}
 r=r_n=\lceil(\log n)^{A_0}\rceil,\quad R=R_n=n-r,\quad
 J=J_n=\{r,\ldots,R\},\quad b=\lceil r/4\rceil,\\
 I=\{r+b,\ldots,R-b\},\qquad
 I^{\rm ret}=\{r+b+1,\ldots,R-b\}.
 \end{gathered}
 \label{B:eq:conditional-intervals}
\end{equation}
All matrices in this subsection act on $\mathscr S(J)=\{x,y\}\times J$. We write $k_\alpha=k$ for $\alpha=(\sigma,k)$ and
$d(\alpha,\gamma)=|k_\alpha-k_\gamma|+ \mathbf1_{\{\sigma_\alpha\ne\sigma_\gamma\}}$. The interaction radius $R_*$ may be increased by a fixed amount. Take $n\ge n_0$, with $n_0$
large enough that $r<n/4$, $r\ge8/\beta$, and every modified density exponent lies outside the interaction neighborhood of $J$.

For $q\in\mathcal Q_n^\sharp$ from \eqref{B:eq:admissible-boundary-configurations}, abbreviate the conditional kernel \eqref{B:eq:initial-coefficients-common-kernel} by $\mu_q$,
its Hamiltonian by $H^q$, and its minimizer by $m=m^q$. Retain $\kappa=n\beta$ and the density exponents $c_\alpha$ of \eqref{B:eq:general-coefficient-notation}. Thus
\begin{equation}
 H^q(Z)=\mathcal P_{n,J}^q(Z)-\kappa^{-1}
             \sum_{\alpha\in\mathscr S(J)}c_\alpha\log Z_\alpha,
 \qquad c_{\sigma,k}\ge\beta k/2>1.
 \label{B:eq:conditional-hamiltonian}
\end{equation}
Expectation under $\mu_q$ is denoted by $\mathbb E_q$. Define
\begin{equation}
 A=\nabla^2H^q(m),\quad C=A^{-1},\quad M=\operatorname{diag}(m),
 \quad d_\alpha=\sqrt\kappa\,m_\alpha,\quad t_k=k/n.
 \label{B:eq:conditional-matrix-data}
\end{equation}
Constants below may depend on $V,a,\beta,A_0$ and on a subsequently fixed extension parameter, but are uniform in $n$ and in $q\in\mathcal Q_n^\sharp$. Moment constants are
independent of the moment order in the stated ranges.

\begin{proposition}\label{B:prop:interior-conditional-geometry}
There are positive constants $a_-,a_+,d_-,d_+,c$ such that
\begin{equation}
 a_-I\le A\le a_+I,\qquad
 d_-\sqrt{k_\alpha}\le d_\alpha\le d_+\sqrt{k_\alpha},\qquad
 |C_{\alpha\gamma}|\le C_0e^{-cd(\alpha,\gamma)}.
 \label{B:eq:interior-hessian-bounds}
\end{equation}
With $\vartheta_n=r^{-1/8}$ and the first correction $f_{n,\sigma,k}$ of \eqref{B:eq:first-correction-definition},
\begin{equation}
 \begin{split}
 |m_{\sigma,k}-f_{n,\sigma,k}|/\phi(t_k)
 &\le C_0k^{-2}+C_0e^{-cr}
       +C_0\vartheta_n(e^{-c(k-r)}+e^{-c(R-k)}),\\
 |m_{\sigma,k}/\phi(t_k)-1|+|\log(m_{\sigma,k}/\phi(t_k))|
 &\le C_0k^{-1}+C_0e^{-cr}
       +C_0\vartheta_n(e^{-c(k-r)}+e^{-c(R-k)}).
 \end{split}
 \label{B:eq:interior-minimizer-bounds}
\end{equation}
\end{proposition}
\begin{proof}
Write $z^\sharp=z_n^\sharp$ and $t_k=k/n$. The error splits as
\[
 \frac{|m_{\sigma,k}-f_{n,\sigma,k}|}{\phi(t_k)}
 \le\frac{|m_{\sigma,k}-z^\sharp_{\sigma,k}|}{\phi(t_k)}
    +\frac{|z^\sharp_{\sigma,k}-f_{n,\sigma,k}|}{\phi(t_k)}.
\]
The first summand is bounded by
\[
 C\vartheta_n\left[
 \frac{\phi(r/n)}{\phi(t_k)}e^{-c(k-r)}
 +\frac{\phi(R/n)}{\phi(t_k)}e^{-c(R-k)}\right].
\]
Theorem~\ref{B:thm:reference-minimizer-approximation} and Proposition~\ref{B:prop:conditional-minimizer-centering} give the absolute error and the two boundary contributions. Since
$\phi(t)\asymp\sqrt t$, their division by $\phi(t_k)$ uses
\[
 \frac{\phi(r/n)}{\phi(k/n)}\le C_0,\qquad
 \frac{\phi(R/n)}{\phi(k/n)}\le C_0\sqrt{R/k}
                              \le C_0\sqrt{1+R-k}.
\]
For $j=R-k\ge0$,
\[
 (1+j)^{1/2}e^{-cj}
 \le\left(\sup_{u\ge0}(1+u)^{1/2}e^{-cu/2}\right)e^{-cj/2}.
\]
For the error in \eqref{B:eq:reference-minimizer-relative-error}, the term at the upper index satisfies
\[
 \sqrt{n/k}\,e^{-c(n-k)}
 \le
 \begin{cases}
 \sqrt2\,e^{-cr},&k\ge n/2,\\
 \sqrt n\,e^{-cn/2}\le C_ce^{-cn/4},&k<n/2.
 \end{cases}
\]
Since $r<n/4$, both lines are bounded by $C e^{-c'r}$. The finite initial term is bounded by $C_0e^{-c'r}$ on $J$. This proves the first line of
\eqref{B:eq:interior-minimizer-bounds}. The relative first correction is $O(k^{-1})$. The resulting ratio $m_{\sigma,k}/\phi(t_k)$ lies in $[1/2,3/2]$ for large $n$;
$|\log u|\le2|u-1|$ there proves the second line.

The Hessian lower bound follows from Lemma~\ref{B:lem:coefficient-polynomial-convexity}. Its logarithmic part is the diagonal $c_\alpha/(\kappa m_\alpha^2)$, which is bounded above
since $m_\alpha^2\asymp k_\alpha/n$. The polynomial entries are bounded on this coefficient range and have finite bandwidth. Their row and column sums are bounded uniformly,
proving the upper bound. The bounds for $d_\alpha$ follow as well. For the inverse, set $Q=I-A/(2a_+)$. Then $\|Q\|\le1-a_-/(2a_+)=:\rho<1$ and
\[
 C=(2a_+)^{-1}\sum_{j=0}^\infty Q^j,\qquad
 (Q^j)_{\alpha\gamma}=0
       \quad\hbox{if }jR_*<d(\alpha,\gamma).
\]
With $j_0=\lceil d(\alpha,\gamma)/R_*\rceil$, we have
\[
 |C_{\alpha\gamma}|
 \le\frac1{2a_+}\sum_{j=j_0}^{\infty}\rho^j
 =\frac{\rho^{j_0}}{2a_+(1-\rho)}
 \le a_-^{-1}e^{-c_1d(\alpha,\gamma)},\qquad
 c_1=-R_*^{-1}\log\rho>0.
\]
The series converges in operator norm. Each power has the stated range because a product of $j$ entries can join indices at distance at most $jR_*$. All constants therefore apply
to every matrix size and every permitted boundary configuration.
\end{proof}

We also retain the following two consequences of Section~\ref{B:sec:general-potentials-and-convex-reference-measures}:
\begin{equation}
 \operatorname{Ent}_{\mu_q}(F^2)
       \le\frac2{\kappa c_*}\mathbb E_q|\nabla F|^2,
 \qquad \|Z_\alpha-m_\alpha\|_p\le C_0\sqrt{p/\kappa}
       \quad(p\ge2).
 \label{B:eq:interior-concentration}
\end{equation}
They are Lemma~\ref{B:lem:reference-logarithmic-sobolev} and Proposition~\ref{B:prop:reference-coordinate-concentration}, applied using
Proposition~\ref{B:prop:unchanged-conditional-distributions}. On $|Z_\alpha/m_\alpha-1|\le1/2$, every nonzero coordinate derivative of $\kappa H^q(m+x/\sqrt\kappa)$ of order
$h=3,4,5$ is bounded by $C_hk_*^{1-h/2}$, where $k_*$ is any site in its fixed interaction neighborhood, by Lemma~\ref{B:lem:conditional-scaled-derivatives}. These constants
include the fixed boundary coordinates.

\subsubsection{Positive, reciprocal, and logarithmic moments}
\label{B:subsection:interior-logarithmic-moments}

Logarithms and reciprocal coefficients are sensitive to small values, so concentration alone is insufficient. The following moment bounds address that issue before the logarithmic
ratios are compared with their linear Gaussian approximation.

Put
\begin{equation}
 \zeta_\alpha=\log(Z_\alpha/m_\alpha),\qquad
 \eta_k=\zeta_{y,k}-\zeta_{x,k},\qquad
 (L_qv)_k=v_{y,k}/m_{y,k}-v_{x,k}/m_{x,k},\qquad
 K_q=\kappa^{-1}L_qCL_q^{\mathsf T}.
 \label{B:eq:conditional-logarithms}
\end{equation}
In particular, $K_q$ is the covariance of $L_qG$ for $G\sim N(0,\kappa^{-1}C)$.

\begin{lemma}\label{B:lem:interior-reciprocal-logarithmic-moments}
There is $c_{\rm mom}>0$ such that
\begin{equation}
 \|Z_\alpha/m_\alpha\|_s+\|m_\alpha/Z_\alpha\|_s\le C_0
       \quad(1\le s\le c_{\rm mom}k_\alpha).
 \label{B:eq:interior-reciprocal-moments}
\end{equation}
After decreasing $c_{\rm mom}$, the same range gives
\begin{equation}
 \begin{aligned}
 \|\zeta_\alpha\|_p&\le C_0\sqrt{p/k_\alpha}
                 &&(2\le p\le c_{\rm mom}k_\alpha),\\
 \|\eta_k\|_p&\le C_0\sqrt{p/k}
                 &&(2\le p\le c_{\rm mom}k).
 \end{aligned}
 \label{B:eq:interior-log-moments}
\end{equation}
\end{lemma}
\begin{proof}
For $s\le ck_\alpha$, \eqref{B:eq:interior-concentration} and $m_\alpha\asymp\phi(t_{k_\alpha})$ give
\[
 \|Z_\alpha\|_s\le m_\alpha+C_0\sqrt{s/\kappa}
        \le C_0\phi(t_{k_\alpha}).
\]
The same bound holds for each of the finitely many coordinates in an interaction neighborhood of $\alpha$, with fixed boundary values in place of random coordinates where
necessary. Choose an integer $D_V\ge\max\{2,2\deg V\}$. A monomial of $\partial_\alpha\mathcal P_{n,J}^q$ has degree between one and $2\deg V-1$, by separate evenness. H\"older,
with at most $D_V$ factors, therefore gives
\begin{equation}
 \|\partial_\alpha\mathcal P_{n,J}^q(Z)\|_s
 \le C_0\sum_{m=1}^{\deg V}\phi(t_{k_\alpha})^{2m-1}
 \le C_0\phi(t_{k_\alpha})
 \label{B:eq:conditional-polynomial-force-moment}
\end{equation}
when $D_Vs\le ck_\alpha$. The number of monomials is bounded independently of $n$.

For an integer $s\ge1$ with $s-1<c_\alpha$, integrate $\partial_\alpha[Z_\alpha^{-(s-1)}e^{-\kappa H^q(Z)}]$ in $Z_\alpha$. This yields
\begin{equation}
 (c_\alpha-s+1)\mathbb E_qZ_\alpha^{-s}
 =\kappa\mathbb E_q[Z_\alpha^{-(s-1)}
                         \partial_\alpha\mathcal P_{n,J}^q(Z)].
 \label{B:eq:conditional-negative-moment-identity}
\end{equation}
Here $n,q$ are fixed during integration. The bound $\mathcal P_{n,J}^q(Z)\ge c_*|Z|^2/2$ gives
\[
 \int Z_\alpha^{-s}e^{-\kappa\mathcal P_{n,J}^q(Z)}
                 \prod_\gamma Z_\gamma^{c_\gamma}\,\dd Z
 \le
 \left(\int_0^\infty u^{c_\alpha-s}e^{-\kappa c_*u^2/2}\,\dd u\right)
 \prod_{\gamma\ne\alpha}
 \left(\int_0^\infty u^{c_\gamma}e^{-\kappa c_*u^2/2}\,\dd u\right)
 <\infty .
\]
The condition is $c_\alpha-s>-1$. The boundary expression at $Z_\alpha=\varepsilon$ is bounded, after integrating the other coordinates, by a constant times
$\varepsilon^{c_\alpha-s+1}$. At $Z_\alpha=T$ it is bounded by a polynomial in $T$ times $e^{-\kappa c_*T^2/2}$. Both limits vanish. Every monomial in the differentiated integrand
has an integrable power at zero and Gaussian decay at infinity. Fubini and dominated convergence therefore justify integration first over $\varepsilon<Z_\alpha<T$ and then over all
coordinates. Choose $c_{\rm mom}$ so that $s\le c_\alpha/2$. For $s>1$, H\"older applied to \eqref{B:eq:conditional-negative-moment-identity} gives
\[
 (c_\alpha-s+1)\|Z_\alpha^{-1}\|_s^s
 \le\kappa\|Z_\alpha^{-1}\|_s^{s-1}
                   \|\partial_\alpha\mathcal P_{n,J}^q\|_s.
\]
For $s=1$, the identity reads $c_\alpha\mathbb E_q Z_\alpha^{-1} =\kappa\mathbb E_q\partial_\alpha\mathcal P_{n,J}^q$. For $s>1$ the inverse moment just proved finite is also
positive, so division gives
\[
 m_\alpha\|Z_\alpha^{-1}\|_s
 \le \frac{\kappa m_\alpha}
              {c_\alpha-s+1}
       \|\partial_\alpha\mathcal P_{n,J}^q\|_s
 \le
 \frac{C\kappa\phi(t_{k_\alpha})^2}{c_\alpha}\le C.
\]
The first inequality holds at $s=1$ as well. Fix a constant for the integer moment range $c_{\rm int}>0$ satisfying the moment bound for the derivative of the polynomial part of
the potential and the positive lower bound for the denominator. Decrease $c_{\rm mom}$ to at most $c_{\rm int}/4$ and increase $n_0$ until $c_{\rm int}r\ge2$. For
$1\le s\le c_{\rm mom}k_\alpha$, the integer $j=\lceil s\rceil$ satisfies $j\le c_{\rm int}k_\alpha$; monotonicity gives $\|m_\alpha/Z_\alpha\|_s\le\|m_\alpha/Z_\alpha\|_j\le C$.

Integration of $t^{-1}$ gives $0\le\log u\le u-1$ for $u\ge1$ and $0\le-\log u\le(1-u)/u$ for $0<u\le1$. Hence $|\log u|\le|u-1|(1+u^{-1})$. After another fixed decrease of
$c_{\rm mom}$, the reciprocal estimate is available at order $2p$ throughout the stated logarithmic range. Thus
\[
 \|\zeta_\alpha\|_p
 \le\|(Z_\alpha-m_\alpha)/m_\alpha\|_{2p}
                 (1+\|m_\alpha/Z_\alpha\|_{2p})
 \le C_0\sqrt{p/k_\alpha}.
\]
The triangle inequality proves the estimate for $\eta_k$. All decreases of the range constant are by fixed factors depending only on $V,a,\beta$.
\end{proof}

Choose henceforth $c_{\rm rc}>0$ with $8c_{\rm rc}\le c_{\rm mom}$. For each fixed $\delta>0$, coordinate concentration gives
\begin{equation}
 \sup_q\mu_q\{\max_\alpha|Z_\alpha-m_\alpha|/m_\alpha>\delta\}
       \le C_\delta n e^{-c_\delta r}.
 \label{B:eq:conditional-relative-deviation}
\end{equation}
Indeed the squared threshold in coordinate $\alpha$, multiplied by $\kappa$, is at least $c\delta^2k_\alpha$; sum the resulting tail bounds over at most $2n$ coordinates.

\subsubsection{Sums of squared logarithms}
\label{B:subsection:conditional-squared-logarithms}

For deterministic $w_\alpha\ge0$, define $\mathcal M(w)=\sum_\alpha w_\alpha\zeta_\alpha^2$.

\begin{proposition}\label{B:prop:conditional-logarithmic-energy}
For $2\le p\le c_{\rm rc}r$,
\begin{equation}
 \begin{split}
 \mathbb E_q\mathcal M(w)&\le C_0\sum_\alpha w_\alpha/k_\alpha,\\
 \|\mathcal M(w)\|_p&\le C_0p\sum_\alpha w_\alpha/k_\alpha,\\
 \|\mathcal M(w)-\mathbb E_q\mathcal M(w)\|_p
 &\le C_0p^{3/2}
       \left(\sum_\alpha w_\alpha^2/k_\alpha^2\right)^{1/2}.
 \end{split}
 \label{B:eq:conditional-energy-moments}
\end{equation}
\end{proposition}
\begin{proof}
For bounded smooth $g$ with $\mathbb E_qg=0$, apply \eqref{B:eq:interior-concentration} to $1+t g$. Taylor expansion on the bounded range of $g$ gives
\[
 \operatorname{Ent}_{\mu_q}(1+t g)^2
   =2t^2\mathbb E_qg^2+O_g(|t|^3),\qquad
 \mathbb E_q|\nabla(1+t g)|^2=t^2\mathbb E_q|\nabla g|^2.
\]
Divide by $2t^2$ and let $t\to0$. Approximation in the form norm then yields $\operatorname{Var}_{\mu_q}(g)\le(\kappa c_*)^{-1} \mathbb E_q|\nabla g|^2$ on the form domain. We
first derive the $L^p$ consequence of the Poincar\'e inequality in \eqref{B:eq:interior-concentration}. For a smooth $F$, let $N_p=\|F-\mathbb E_qF\|_p$ and $G_p=\||\nabla F|\|_p$.
For $p>2$, applying Poincar\'e to $|F-\mathbb E_qF|^{p/2}$, followed by H\"older, gives
\[
 N_p^p-N_{p/2}^p
 \le\frac{p^2}{4\kappa c_*}
        \mathbb E_q[|F-\mathbb E_qF|^{p-2}|\nabla F|^2]
 \le\frac{p^2}{4\kappa c_*}N_p^{p-2}G_p^2.
\]
For $0\le u\le1$ and $p\ge2$, $1-u^p\ge1-u^2$. Division by $N_p^{p-2}$, with the zero case treated directly, gives $N_p^2\le N_{p/2}^2+p^2G_p^2/(4\kappa c_*)$. Let $j_*$ be the
least integer for which $p/2^{j_*}\le2$. Monotonicity of $L^s$ norms gives $G_{p/2^j}\le G_p$ and $N_{p/2^{j_*}}\le N_2$. Iterating the last inequality therefore gives
\[
 N_p^2
 \le N_2^2+\frac{G_p^2}{4\kappa c_*}
                      \sum_{j=0}^{j_*-1}(p/2^j)^2
 \le \frac{G_p^2}{\kappa c_*}
                \left(1+\frac{p^2}{3}\right).
\]
This includes $p=2$ by Poincar\'e, and yields
\begin{equation}
 \|F-\mathbb E_qF\|_p
       \le C_0p\kappa^{-1/2}\||\nabla F|\|_p.
 \label{B:eq:conditional-poincare-moments}
\end{equation}
For the regularization $g_\varepsilon=(u^2+\varepsilon)^{p/4}$, where $u=F-\mathbb E_qF$ and $0<\varepsilon\le1$,
\[
 g_\varepsilon^2\le C_p(1+|u|^p),\qquad
 |\nabla g_\varepsilon|^2
 =\frac{p^2}{4}u^2(u^2+\varepsilon)^{p/2-2}|\nabla F|^2
 \le C_p(1+|u|^{p-2})|\nabla F|^2.
\]
For $2<p<4$ the last inequality uses the negative exponent; for $p\ge4$ use $u^2+\varepsilon\le u^2+1$. The last majorant is integrable by H\"older when $u,\nabla F\in L^p$.
Dominated convergence in Poincar\'e therefore gives the inequality for $|u|^{p/2}$, including the set $u=0$. One first applies this to bounded smooth truncations if necessary and
then uses the same $L^p$ majorants.

Since $2p\le2c_{\rm rc}r\le c_{\rm mom}k_\alpha$,
\[
 \mathbb E_q\mathcal M(w)
 =\sum_\alpha w_\alpha\|\zeta_\alpha\|_2^2
 \le C_0\sum_\alpha\frac{w_\alpha}{k_\alpha},\qquad
 \|\mathcal M(w)\|_p
 \le\sum_\alpha w_\alpha\|\zeta_\alpha\|_{2p}^2
 \le C_0p\sum_\alpha\frac{w_\alpha}{k_\alpha}.
\]
These estimates also give the required function moments before the centered inequality is applied. For that inequality,
\[
 \partial_\alpha\mathcal M=2w_\alpha\zeta_\alpha/Z_\alpha,
 \qquad
 \|\zeta_\alpha/Z_\alpha\|_p
 \le\|\zeta_\alpha\|_{2p}\|Z_\alpha^{-1}\|_{2p}
 \le C_0\sqrt{p\kappa}/k_\alpha.
\]
Minkowski in $L^{p/2}$ therefore yields
\[
 \||\nabla\mathcal M|\|_p^2
 \le4\sum_\alpha w_\alpha^2
                 \|\zeta_\alpha/Z_\alpha\|_p^2
 \le C_0p\kappa\sum_\alpha w_\alpha^2/k_\alpha^2.
\]
Substitution into \eqref{B:eq:conditional-poincare-moments} proves the claim. For its application to the unbounded logarithm, replace $\zeta$ by
$\tau_N(\zeta)=\int_0^\zeta\chi(s/N)\,\dd s$, where $0\le\chi\le1$ is smooth, equals one near zero and has compact support. We have $|\tau_N(\zeta)|\le|\zeta|$ and
$|\tau_N'(\zeta)|\le1$. For $\mathcal M_N=\sum_\alpha w_\alpha\tau_N(\zeta_\alpha)^2$,
\[
 |\mathcal M_N|\le\mathcal M,\qquad
 |\partial_\alpha\mathcal M_N|
 \le2w_\alpha|\zeta_\alpha|/Z_\alpha .
\]
Each gradient vanishes near $Z_\alpha=0$ and near $Z_\alpha=\infty$ for fixed $N$; it is bounded in that coordinate. Smooth cutoffs place $\mathcal M_N$ in the form domain.
Dominated convergence, using the $L^p$ majorants above and the finite number of coordinates at fixed $n$, gives $\mathcal M_N\to\mathcal M$ and
$\nabla\mathcal M_N\to\nabla\mathcal M$ in $L^p$. Also $|\mathbb E_q(\mathcal M_N-\mathcal M)| \le\|\mathcal M_N-\mathcal M\|_p\to0$. The centered inequality consequently passes to
$\mathcal M$.
\end{proof}

Fix $1<\lambda<2$ and weights $0\le w_\alpha\le C_wn^{-\lambda}k_\alpha^{\lambda-1}$. There are at most two coordinates at each site, and elementary integral comparison gives
\begin{equation}
 \begin{split}
 \sum_\alpha w_\alpha/k_\alpha&\le C_\lambda/n,\\
 a_{n,\lambda}:=
 n^{-\lambda}\left(\sum_{k=r}^Rk^{2\lambda-4}\right)^{1/2}
 &\le C_\lambda
 \begin{cases}
 n^{-\lambda}r^{\lambda-3/2},&1<\lambda<3/2,\\
 n^{-3/2}\sqrt{1+\log(n/r)},&\lambda=3/2,\\
 n^{-3/2},&3/2<\lambda<2.
 \end{cases}
 \end{split}
 \label{B:eq:conditional-weight-sums}
\end{equation}
For example, when $2\lambda-4<-1$, the sum is at most $r^{2\lambda-4}+\int_r^\infty t^{2\lambda-4}\,\dd t \le C_\lambda r^{2\lambda-3}$. At exponent $-1$ it is at most
$1+\log(R/r)$, and above $-1$ it is at most $C_\lambda n^{2\lambda-3}$.

\begin{corollary}\label{B:cor:conditional-logarithmic-probabilities}
There is a deterministic $K_\lambda$ such that for every $D>0$,
\begin{equation}
 \begin{split}
 \sup_q\mu_q\{\mathcal M(w)>K_\lambda/n\}&\le C_Dn^{-D},\\
 \sup_q\mu_q\{\max_\alpha\sqrt{k_\alpha}|\zeta_\alpha|
                              >C_D\sqrt{\log(\me+n)}\}&\le C_Dn^{-D}.
 \end{split}
 \label{B:eq:conditional-energy-and-maximum-probabilities}
\end{equation}
The starting index may depend on $D$; $K_\lambda$ does not.
\end{corollary}
\begin{proof}
Put $\delta_\lambda=\tfrac12\min\{\lambda-1,1/2\}>0$. Equation~\eqref{B:eq:conditional-weight-sums} implies $a_{n,\lambda}\le C_\lambda n^{-1-\delta_\lambda}$ for large $n$. For
every fixed $p\ge2$ and $\varepsilon>0$, Markov's inequality gives
\begin{equation}
 \sup_q\mu_q\{|\mathcal M(w)-\mathbb E_q\mathcal M(w)|
                                      >\varepsilon/n\}
 \le(C_\lambda p^{3/2}/\varepsilon)^p n^{-p\delta_\lambda}.
 \label{B:eq:conditional-centered-energy-probability}
\end{equation}
Choose an integer $p\ge2$ with $p\delta_\lambda>D+1$. This number is fixed before $n$ tends to infinity, so $p\le c_{\rm rc}r_n$ for all sufficiently large $n$. Fix
$K_\lambda=C_\lambda+1$, using the constant in the mean bound. Then
\[
 \{\mathcal M(w)>K_\lambda/n\}
 \subseteq
 \{|\mathcal M(w)-\mathbb E_q\mathcal M(w)|>1/n\},
 \qquad
 \sup_q\mu_q\{\mathcal M(w)>K_\lambda/n\}
 \le (C_\lambda p^{3/2})^p n^{-p\delta_\lambda}.
\]
The constant $K_\lambda$ is independent of the chosen decay exponent $D$; the displayed prefactor and the starting index can depend on $D$. For the second take
$p_n=\lceil(D+4)\log(\me+n)\rceil$. From \eqref{B:eq:interior-log-moments}, at threshold $eC_0\sqrt{p_n}$ the probability for each $\sqrt{k_\alpha}|\zeta_\alpha|$ is at most
$e^{-p_n}$. Summing over the coordinates gives $2ne^{-p_n}\le C_Dn^{-D}$. Since $p_n/r\to0$, choose $n_0$ so that all fixed multiples of $p_n$ used here lie in the permitted moment
range for
$n\ge n_0$.
\end{proof}

\subsubsection{Exponential sums and band matrices}
\label{B:subsection:exponential-sums-and-band-matrices}

The covariance and mean comparisons below use the following exponential summation bounds. Site indices below are positive integers. For $h\ge0$ and $\nu>0$,
\begin{equation}
 (k/l)^h\le(1+|k-l|)^h,\qquad
 \sum_{l\ge1}(k/l)^he^{-\nu|k-l|}\le C_{h,\nu}.
 \label{B:eq:exponential-site-sums}
\end{equation}
Indeed, for $k>l$ the first bound follows from $k/l=1+(k-l)/l$; for $k\le l$ its left side is at most one. The sum is bounded by $1+2\sum_{j\ge1}(1+j)^he^{-\nu j}<\infty$. Its
finiteness follows, for instance, by comparison with
$e^\nu\int_0^\infty(2+t)^he^{-\nu t}\,\dd t$.

\begin{lemma}\label{B:lem:exponential-convolution-bounds}
For $k,p\in J$ and fixed $h\ge0$, the sum
\begin{equation}
 \begin{split}
 &\sum_{j\in J}(1+|k-j|+|j-p|)^he^{-\nu(|k-j|+|j-p|)}\\
 &\hspace{8mm}\times
 \left[\frac{1+|j-k|}{k}+e^{-\nu(j-r)}+e^{-\nu(R-j)}\right]
 \le C_{h,\nu}e^{-\nu|k-p|/4}
     [k^{-1}+e^{-\nu(k-r)/4}+e^{-\nu(R-k)/4}].
 \end{split}
 \label{B:eq:boundary-exponential-convolution}
\end{equation}
The same estimate without the square brackets is bounded by $C_{h,\nu}e^{-\nu|k-p|/2}$.
\end{lemma}
\begin{proof}
Set $S_j=|k-j|+|j-p|$. Then $S_j\ge|k-p|$ and $S_j\ge|k-j|$. The boundedness of $(1+t)^{h+1}e^{-\nu t/4}$ gives
\[
 (1+S_j)^h(1+|j-k|)e^{-\nu S_j}
 \le C_{h,\nu}e^{-\nu|k-p|/4}e^{-\nu|k-j|/2}.
\]
Summation proves the $k^{-1}$ term. For the left boundary use also $S_j+j-r\ge k-r$. Then
\[
 \tfrac34S_j+j-r
 \ge\tfrac14|k-p|+\tfrac14(k-r)+\tfrac14|k-j|.
\]
Set $A_{h,\nu}=\sup_{s\ge0}(1+s)^he^{-\nu s/4}<\infty$. The left boundary summand is consequently at most
\[
 A_{h,\nu}e^{-\nu|k-p|/4}e^{-\nu(k-r)/4}e^{-\nu|k-j|/4}.
\]
Summing the last factor gives at most $1+2e^{-\nu/4}/(1-e^{-\nu/4})$, proving the left boundary bound. At the right boundary, $S_j+R-j\ge R-k$, and hence
\[
 \tfrac34 S_j+R-j
 \ge\tfrac14|k-p|+\tfrac14(R-k)+\tfrac14|k-j|.
\]
It follows that the right boundary summand is at most $C_{h,\nu}e^{-\nu|k-p|/4}e^{-\nu(R-k)/4} e^{-\nu|k-j|/4}$. Its sum has the required bound because
\[
 \sum_{j\in J}e^{-\nu|k-j|/4}
 \le 1+\frac{2e^{-\nu/4}}{1-e^{-\nu/4}}.
\]
Without the brackets, put $C'_{h,\nu}=\sup_{s\ge0}(1+s)^he^{-\nu s/4}<\infty$. The inequalities $S_j\ge|k-p|$ and $S_j\ge|k-j|$ give
\[
 (1+S_j)^he^{-\nu S_j}
 \le C'_{h,\nu}e^{-\nu|k-p|/2}e^{-\nu|k-j|/4}.
\]
The same geometric sum proves the bound without the brackets. For a shift $|h_0|\le R_*$, the inequalities $|k-(j+h_0)|\ge|k-j|-R_*$ and $1+|j+h_0-k|\le(1+R_*)(1+|j-k|)$ show
explicitly that a fixed interaction neighborhood changes only these constants.
\end{proof}

Let $P_t$ be a symmetric Markov semigroup, strongly continuous on $L^2(\nu)$, with its contraction extensions to $L^p(\nu)$. The integrals below are strong integrals for
$1\le p<\infty$ and weak-star integrals for $p=\infty$. For a component vector $v=(v_\alpha)$ and a site $i$, use the norm
$\sup_\alpha e^{\eta|k_\alpha-i|}\|v_\alpha\|_{L^p(\nu)}$; for two component indices use $\sup_{\alpha,\gamma}e^{\eta(|k_\alpha-i|+|k_\gamma-j|)} \|v_{\alpha\gamma}\|_{L^p(\nu)}$.

\begin{lemma}\label{B:lem:weighted-band-semigroups}
For some $\eta_0>0$, the operators
\begin{equation}
 \mathcal K_1=\int_0^\infty e^{-tA}P_t\,\dd t,\qquad 
 \mathcal K_2=\int_0^\infty(e^{-tA}\otimes e^{-tA})P_t\,\dd t
 \label{B:eq:weighted-free-inverses}
\end{equation}
have bounds $K_1,K_2$ in the respective component norms, uniform in $n,q,p,i,j$ for $0\le\eta\le\eta_0$.
\end{lemma}
\begin{proof}
Fix a coordinate $\gamma$ and let $D_{\alpha\alpha}=e^{\eta d(\alpha,\gamma)}$. If $A_{\alpha\rho}\ne0$, then $d(\alpha,\rho)\le R_*$, and
\[
 |(DAD^{-1}-A)_{\alpha\rho}|
 \le |A_{\alpha\rho}|(e^{\eta R_*}-1).
\]
There are at most $C_{R_*}$ entries in each row and column. Schur's test bounds the operator norm of the difference by $C_0\eta$ for $0<\eta\le1$. Choose $\eta_1>0$ with
$C_0\eta_1\le a_-/2$. For $v(t)=e^{-tDAD^{-1}}v(0)$,
\[
 \frac\dd{\dd t}|v(t)|^2
 =-2\operatorname{Re}\langle DAD^{-1}v(t),v(t)\rangle
 \le-a_-|v(t)|^2.
\]
Hence $\|De^{-tA}D^{-1}\|\le e^{-a_-t/2}$ and
\begin{equation}
 |(e^{-tA})_{\alpha\gamma}|
       \le e^{-a_-t/2-\eta_1d(\alpha,\gamma)}.
 \label{B:eq:band-matrix-heat-kernel}
\end{equation}
Take $\eta_0=\eta_1/2$. For $0\le\eta\le\eta_0$,
\[
 e^{\eta d(\alpha,\gamma)}e^{-\eta_1d(\alpha,\gamma)}
 \le e^{-\eta_1d(\alpha,\gamma)/2}.
\]
For the weight $w_\alpha=e^{\eta|k_\alpha-i|}$, the triangle inequality gives $w_\alpha/w_\gamma\le e^{\eta d(\alpha,\gamma)}$. Thus
\[
 \sup_\alpha\sum_\gamma
 \frac{w_\alpha}{w_\gamma}|(e^{-tA})_{\alpha\gamma}|
 \le e^{-a_-t/2}
       \sup_\alpha\sum_\gamma e^{-\eta_1d(\alpha,\gamma)/2}
 \le C_0e^{-a_-t/2}.
\]
For the weighted component norm $\|\cdot\|_{p,w}$ in the statement, contractivity of $P_t$ and the triangle inequality yield
\[
 \|e^{-tA}P_tv\|_{p,w}
 \le C_0e^{-a_-t/2}\|v\|_{p,w},\qquad
 \left\|\int_T^\infty e^{-tA}P_tv\,\dd t\right\|_{p,w}
 \le\frac{2C_0}{a_-}e^{-a_-T/2}\|v\|_{p,w}.
\]
For two indices put $w_{\alpha\theta}=w_\alpha\widetilde w_\theta$, with the two centers as in the statement. The weighted row sum is
\[
 \sum_{\gamma,\rho}\frac{w_{\alpha\theta}}{w_{\gamma\rho}}
 |(e^{-tA})_{\alpha\gamma}(e^{-tA})_{\theta\rho}|
 =\left(\sum_\gamma\frac{w_\alpha}{w_\gamma}
                         |(e^{-tA})_{\alpha\gamma}|\right)
  \left(\sum_\rho\frac{\widetilde w_\theta}{\widetilde w_\rho}
                         |(e^{-tA})_{\theta\rho}|\right)
 \le C_0^2e^{-a_-t}.
\]
Contractivity of $P_t$ and integration from $T$ to infinity give $C_0^2a_-^{-1}e^{-a_-T}$ for the operator tail. In particular the constants in the statement may be taken as
$K_1=2C_0/a_-$ and $K_2=C_0^2/a_-$. For $p<\infty$ the semigroup is strongly continuous and these are Bochner integral bounds. For $p=\infty$, define each integral by pairing its
component with an arbitrary element of $L^1(\nu)$. The absolute value of that pairing is bounded by the corresponding integrable majorant times the $L^1$ norm; duality gives the
weak-star integral and the displayed tail estimate.
\end{proof}

\subsubsection{A smooth uniformly convex potential}
\label{B:subsection:smooth-conditional-potential}

A smooth uniformly convex potential on the whole space permits the covariance and mean expansions below. We choose it to agree with the conditional potential on the stated set. The
comparison with the original distribution is then controlled by the probabilities of its complement.

Let $x=\sqrt\kappa(Z-m)$ and
\begin{equation}
 V_*(x)=\kappa[H^q(m+x/\sqrt\kappa)-H^q(m)],
 \qquad x_\alpha>-d_\alpha.
 \label{B:eq:original-normalized-potential}
\end{equation}
Expand the polynomial in $H^q$ into its local monomials and keep each logarithm as a separate term depending on one coordinate. For a term depending on the coordinate set $B$, let
$u_B$ be its remainder after subtracting the Taylor polynomial of degree two at zero. The sets $B$ have bounded size and diameter, and only a bounded number meet a given
coordinate. Since $\nabla V_*(0)=0$,
\[
 V_*(x)=\tfrac12x^{\mathsf T}Ax+\sum_Bu_B(x_B).
\]
Set $k_B=\min_{\alpha\in B}k_\alpha$ and $X_B=\sum_{\alpha\in B}|x_\alpha|$. For $|x_\alpha|\le3\delta d_\alpha$ and
$3\delta<1/2$,
\begin{equation}
 |D^ju_B|\le C_j k_B^{-1/2}X_B^{3-j}\ (0\le j\le2),
 \qquad |D^hu_B|\le C_hk_B^{1-h/2}\ (3\le h\le5).
 \label{B:eq:local-potential-remainders}
\end{equation}
The last bounds follow from the scaled derivatives established above. For the first ones apply Taylor's integral formula to $D^ju_B$; its first $3-j$ coefficients of the expansion
vanish.

Fix an even smooth $0\le\chi\le1$ which equals one on $[-1,1]$ and zero off $(-2,2)$. Define
\begin{equation}
 U(x)=\sum_Bu_B(x_B)\chi_B(x_B),\qquad
 \chi_B=\prod_{\alpha\in B}\chi(x_\alpha/(\delta d_\alpha)),
 \qquad V_0(x)=\tfrac12x^{\mathsf T}Ax+U(x).
 \label{B:eq:extended-conditional-potential}
\end{equation}
Each summand is zero outside its cutoff support. Extend also the logarithms by
\begin{equation}
 \begin{split}
 l_\alpha(t)&=t/d_\alpha+
 \chi(t/(\delta d_\alpha))[\log(1+t/d_\alpha)-t/d_\alpha],\\
 f_k(x)&=\sqrt k\,[l_{y,k}(x_{y,k})-l_{x,k}(x_{x,k})],
 \qquad L=f'(0),\qquad B=LCL^{\mathsf T}.
 \end{split}
 \label{B:eq:extended-logarithmic-functions}
\end{equation}
The cutoff correction in $l_\alpha$ is zero for $|t|\ge2\delta d_\alpha$. Set
\begin{equation}
 E_\delta=\{x:|x_\alpha|\le\delta d_\alpha
                              \text{ for every }\alpha\}.
 \label{B:eq:conditional-agreement-event}
\end{equation}

\begin{proposition}\label{B:prop:smooth-conditional-extension}
There is $\delta_0>0$ such that, for $0<\delta\le\delta_0$, $V_0$ and $l_\alpha$ are smooth on their full Euclidean spaces and
\begin{equation}
 a_-I/2\le V_0''\le2a_+I,\qquad \|U''\|_{\rm op}\le C_0\delta.
 \label{B:eq:extended-hessian-bounds}
\end{equation}
They agree with $V_*$ and the original logarithms on $E_\delta$. Their local derivatives satisfy
\begin{equation}
 \begin{split}
 |\partial_\alpha U|&\le C_\delta k_\alpha^{-1/2}
                   \Big(\sum_{d(\alpha,\rho)\le R_*}|x_\rho|\Big)^2,\\
 |\partial_{\alpha\gamma}^2U|&\le C_\delta k_\alpha^{-1/2}
                   \sum_{d(\alpha,\rho)\le R_*}|x_\rho|,\\
 |D^hU|&\le C_{h,\delta} k_B^{1-h/2}\quad(3\le h\le5).
 \end{split}
 \label{B:eq:extended-local-derivatives}
\end{equation}
All derivatives vanish outside their fixed interaction range. Moreover $l_\alpha$ is increasing and $\|f'-L\|_{\rm op}\le C_0\delta$, $cI\le B\le CI$.
\end{proposition}
\begin{proof}
For a coordinate multi-index $\nu$ of length $h$ the exact product rule is
\[
 D^\nu(u_B\chi_B)=
       \sum_{\omega\le\nu}\binom\nu\omega
                   D^\omega u_B\,D^{\nu-\omega}\chi_B.
\]
On its support $X_B\le C\delta\sqrt{k_B}$ and $|D^j\chi_B|\le C_j(\delta\sqrt{k_B})^{-j}$. For $h=2$ and $|\omega|=j\le2$, each term is bounded by
\[
 C k_B^{-1/2}X_B^{3-j}
          (\delta\sqrt{k_B})^{-(2-j)}\le C\delta.
\]
For $3\le h\le5$ and $j\le2$ the same expression with $2$ replaced by $h$ is bounded by $C_{h,\delta}k_B^{1-h/2}$. For $j\ge3$ it is bounded by
$C_jk_B^{1-j/2}(\delta\sqrt{k_B})^{-(h-j)}$, which has the same order. If a cutoff is differentiated, a coordinate in $B$ has size at least $\delta d_\rho\ge c\delta\sqrt{k_B}$.
Thus every term of a first or second derivative is also bounded by $C_\delta k_B^{-1/2}X_B^2$ or $C_\delta k_B^{-1/2}X_B$, respectively. Summing over the bounded number of sets $B$
gives \eqref{B:eq:extended-local-derivatives}. The row and column sums of $U''$ are at most $C_{\rm cut}\delta$, so Schur's test proves its operator bound. For a vector $v$, the
absolute row and column bounds give
\[
 |v^{\mathsf T}U''(x)v|\le C_{\rm cut}\delta|v|^2,\qquad
 (a_--C_{\rm cut}\delta)|v|^2
 \le v^{\mathsf T}V_0''(x)v
 \le(a_++C_{\rm cut}\delta)|v|^2.
\]
Every local remainder and its first two derivatives vanish at zero, so $V_0(0)=0$ and $\nabla V_0(0)=0$. For every $n$ for which $\nabla^2 V_0\ge(a_-/2)I$, integration along the
segment from $0$ to $x$ gives $V_0(x)\ge a_-|x|^2/4$. The cutoff vanishes before the boundary of the box where the original logarithm is evaluated; the zero extensions are
consequently smooth.

Writing $u=t/d_\alpha$, we have
\[
 g(u)=\log(1+u)-u=-\int_0^u\frac{s}{1+s}\,\dd s,\qquad
 g'(u)=-\frac{u}{1+u}.
\]
For $|u|\le2\delta<1/2$, these give $|g(u)|\le u^2$ and $|g'(u)|\le2|u|$. The first derivative is
\[
 d_\alpha l_\alpha'(t)
 =1+\delta^{-1}\chi'(u/\delta)g(u)
                        +\chi(u/\delta)g'(u).
\]
On the support of $\chi'$ the second term is bounded by $4\|\chi'\|_\infty\delta$, and the third is bounded by $4\delta$. Outside the support of the cutoff correction the
derivative equals $d_\alpha^{-1}$. Higher derivatives follow from the finite product rule, using $g^{(h)}(u)=(-1)^{h-1}(h-1)!(1+u)^{-h}$ for $h\ge2$. Consequently
\[
 |d_\alpha l_\alpha'(t)-1|
 \le C\delta,\qquad
 |l_\alpha^{(h)}(t)|\le C_{h,\delta}d_\alpha^{-h}
       \quad(2\le h\le4).
\]
Choose $C\delta\le1/2$. The derivative is positive everywhere and $l_\alpha$ is linear for large $|t|$, hence a bijection. The only nonzero entries of row $k$ of $L$ are
\begin{equation}
 L_{k,xk}=-\sqrt k/d_{x,k},\qquad
 L_{k,yk}=\sqrt k/d_{y,k},\qquad
 B=\operatorname{diag}(\sqrt k)K_q\operatorname{diag}(\sqrt k).
 \label{B:eq:normalized-logarithmic-covariance}
\end{equation}
The row supports are disjoint and $|L^{\mathsf T}v|^2=\sum_k k(d_{x,k}^{-2}+d_{y,k}^{-2})|v_k|^2$. For constants $c_L,C_L>0$, the last identity implies
$c_L|v|^2\le|L^{\mathsf T}v|^2\le C_L|v|^2$. Since $a_+^{-1}I\le C\le a_-^{-1}I$,
\[
 \frac{c_L}{a_+}|v|^2
 \le v^{\mathsf T}Bv=(L^{\mathsf T}v)^{\mathsf T}C(L^{\mathsf T}v)
 \le\frac{C_L}{a_-}|v|^2.
\]
The two entries in row $k$ of $f'(x)-L$ have absolute values bounded by $C\delta\sqrt{k}/d_{\sigma,k}\le C'\delta$. Their disjoint supports give $\|f'-L\|_{\rm op}\le C_0\delta$.
Also $B_{ij}$ is a sum of four entries of $C$, multiplied by uniformly bounded coefficients; this proves its exponential decay. Finally
$\sqrt{k}\,d_{\sigma,k}^{-h}\le Ck^{(1-h)/2}$ gives the derivative bounds for $D^hf_k$, $2\le h\le4$.

All smallness constants can be fixed independently of $n,q$. In addition to the two conditions above choose $\delta_0$ so that
\begin{equation}
 C_{\rm cut}\delta_0\le a_-/2,\quad
 K_1C_1\delta_0\le\tfrac12,\quad
 K_2C_2\delta_0\le\tfrac12,\quad
 C_{\rm mean}\delta_0\le\tfrac12,\quad 3\delta_0<\tfrac12.
 \label{B:eq:conditional-extension-smallness}
\end{equation}
Here $C_1,C_2$ bound multiplication by $U''$ on the weighted component spaces of Lemma~\ref{B:lem:weighted-band-semigroups}, with one or two indices. They depend only on the fixed
range and $C_{\rm cut}$. The constant $C_{\rm mean}$ is the product of a bound for the row sums of $C$ and a bound for the row sums of $\delta^{-1}U''$, increased by the uniform
coordinate standard deviation implied by the lower bound $a_-/2$. These choices imply \eqref{B:eq:extended-hessian-bounds}; the upper bound uses $a_-\le a_+$.
\end{proof}

Fix one such $\delta$ and let $\widetilde\mu_q(\dd x)=\widetilde Z_q^{-1}e^{-V_0(x)}\dd x$. Brackets $\langle\cdot\rangle$ denote its expectation. The definitions are Borel in $q$:
the original coefficients and $m^q$ are continuous, the cutoffs are fixed smooth functions, and every normalizing integral is a Borel integral of a positive density.

\subsubsection{Moments under the extended distribution}
\label{B:subsection:extended-coordinate-moments}

The extension gives concentration of centered coordinates. A separate estimate on their means is needed to control the coordinates themselves and to start the expansion about the
minimizer.

\begin{lemma}\label{B:lem:extended-coordinate-moments}
Writing $\epsilon_\alpha=k_\alpha^{-1/2}$, we have
\begin{equation}
 |\langle x_\alpha\rangle|\le C_\delta\epsilon_\alpha,
 \qquad \|x_\alpha\|_{L^p(\widetilde\mu_q)}\le C_\delta\sqrt p
          \quad(p\ge2).
 \label{B:eq:extended-coordinate-moments}
\end{equation}
\end{lemma}
\begin{proof}
The Bakry--\'Emery inequality \cite{bakry1985diffusions}, with $V_0''\ge a_-I/2$, gives $\|x_\alpha-\langle x_\alpha\rangle\|_p\le C_0\sqrt p$. This follows also by the Herbst
calculation in the proof of Lemma~\ref{B:lem:reference-logarithmic-sobolev}. Put $N=\nabla U$ and $\bar x=\langle x\rangle$. Choose a smooth cutoff $\chi_T$ equal to one on
$|x|\le T$ and zero on $|x|\ge2T$, with $|\nabla\chi_T|\le C/T$. Integration of $\partial_\alpha(\chi_Te^{-V_0})$ gives
\[
 \langle\chi_T(Ax+N)_\alpha\rangle
             =\langle\partial_\alpha\chi_T\rangle.
\]
The right side tends to zero. The bound $|\nabla V_0(x)|\le2a_+|x|$ and $e^{-V_0(x)}\le e^{-a_-|x|^2/4}$ give an integrable majorant for the left side at fixed dimension. Thus
$A\bar x+\langle N\rangle=0$. Since $N(0)=0$,
\[
 |N_\alpha(x)|
 \le\int_0^1\sum_{\rho:d(\alpha,\rho)\le R_*}
              |U''_{\alpha\rho}(tx)|\,|x_\rho|\,\dd t
 \le C_0\delta\sum_{d(\alpha,\rho)\le R_*}|x_\rho|.
\]
The centered second moments give $\langle|x_\rho|\rangle\le|\bar x_\rho|+C_0$. Let $s_C=\sup_\alpha\sum_\gamma|C_{\alpha\gamma}|$ and let $d_*$ bound the number of coordinates in
any interaction neighborhood. Both are independent of the dimension. Multiplying the mean identity by $C$ therefore gives, with
$m_*=\|\bar x\|_\infty$,
\[
 m_*\le s_C C_0d_*\delta(m_*+C_0)
          \le C_{\rm mean}\delta(m_*+1).
\]
Here the constant in \eqref{B:eq:conditional-extension-smallness} has been chosen large enough for the second inequality. The quantity $m_*$ is finite in each dimension by the
density bound. Since $C_{\rm mean}\delta\le1/2$, subtraction yields $m_*\le1$. In particular $\langle x_\rho^2\rangle\le2C_0^2+2$ uniformly in $\rho,n,q$. The quadratic bound in
\eqref{B:eq:extended-local-derivatives} now implies $\langle|N_\alpha|\rangle\le C_\delta k_\alpha^{-1/2}$. Finally,
\[
 |\bar x_\alpha|
 \le C_\delta\sum_\gamma e^{-cd(\alpha,\gamma)}k_\gamma^{-1/2}
 \le C_\delta k_\alpha^{-1/2}
             \sum_{j\ge1}(1+|j-k_\alpha|)^{1/2}e^{-c|j-k_\alpha|}.
\]
Equation~\eqref{B:eq:exponential-site-sums} proves the mean bound; adding it to the centered moment estimate proves the lemma.
\end{proof}

\subsubsection{Exponential decay of covariances}
\label{B:subsection:conditional-covariance-representation}

The mean and covariance expansions require information about the distance between coordinates, not only a global variance bound. The Helffer--Sj\"ostrand representation
\cite{helffer1994correlation} gives the component estimate below for the present conditional distributions. The exponential weighting is related to the estimates discussed in
Section~3.1 of \cite{menz2014brascamp}.

\begin{lemma}\label{B:lem:extended-covariance-decay}
For real $F,G\in H^1(\widetilde\mu_q)$,
\begin{equation}
 |\operatorname{Cov}_{\widetilde\mu_q}(F,G)|
 \le C_\delta\sum_{\alpha,\gamma}e^{-cd(\alpha,\gamma)}
               \|\partial_\alpha F\|_2\|\partial_\gamma G\|_2.
 \label{B:eq:extended-covariance-decay}
\end{equation}
\end{lemma}
\begin{proof}
Use the real Hilbert space $L^2(\widetilde\mu_q;\mathbb R)$ and let $H^1(\widetilde\mu_q;\mathbb R)$ be the completion of real $C_c^\infty$ functions for the norm
$(\|h\|_2^2+\int|\nabla h|^2\,\dd\widetilde\mu_q)^{1/2}$. The closed symmetric bilinear form
\[
 (u,v)\longmapsto\int\nabla u\cdot\nabla v\,\dd\widetilde\mu_q,
 \qquad u,v\in H^1(\widetilde\mu_q;\mathbb R),
\]
has quadratic form $u\mapsto\int|\nabla u|^2\,\dd\widetilde\mu_q$. Write $-\mathscr L$ for its associated nonnegative self-adjoint operator. On smooth functions in its domain,
$\mathscr L=\Delta-\nabla V_0\cdot\nabla$. The smooth density is positive, $V_0(x)\ge a_-|x|^2/4$, and $|\nabla V_0(x)|\le2a_+|x|$. Smooth cutoffs $\chi_R(x)=\chi_0(x/R)$ tending
to one satisfy $\int|\nabla\chi_R|^2\,\dd\widetilde\mu_q\to0$. Thus the form contains the constant function with zero energy. Its semigroup $P_t$ preserves constants and the
measure, and contracts $L^p$ for $1\le p\le\infty$. The Poincar\'e inequality from uniform convexity gives $\|P_tg\|_2\le e^{-a_-t/2}\|g\|_2$ when $\langle g\rangle=0$. Put
$g=G-\langle G\rangle$ and, for $T<\infty$, $h_T=\int_0^T P_tg\,\dd t$. For $\epsilon>0$,
\[
 \frac{P_\epsilon h_T-h_T}{\epsilon}
 =\frac1\epsilon\int_T^{T+\epsilon}P_tg\,\dd t
       -\frac1\epsilon\int_0^\epsilon P_tg\,\dd t
 \longrightarrow P_Tg-g \quad\hbox{in }L^2 .
\]
Thus $h_T\in\operatorname{dom}\mathscr L$ and $-\mathscr Lh_T=g-P_Tg$. Moreover
\[
 \left\|\int_T^\infty P_tg\,\dd t\right\|_2
 \le(2/a_-)e^{-a_-T/2}\|g\|_2,\qquad
 \|P_Tg\|_2\le e^{-a_-T/2}\|g\|_2.
\]
Closedness of the generator gives $h=\int_0^\infty P_tg\,\dd t\in\operatorname{dom}\mathscr L$ and $-\mathscr Lh=g$. It has mean zero because each $P_tg$ does. The form identity
and Poincar\'e give uniqueness among mean zero solutions.

For compactly supported smooth $v$, integration by parts gives $\langle f,-\mathscr Lg\rangle=\langle\nabla f,\nabla g\rangle$. Using
$\partial_\alpha\mathscr Lv=\mathscr L\partial_\alpha v -\sum_\gamma V_{0,\alpha\gamma}''\partial_\gamma v$, we obtain
\[
 \begin{split}
 \|\mathscr Lv\|_2^2
 &=-\sum_\alpha\langle\partial_\alpha v,
                                 \partial_\alpha\mathscr Lv\rangle\\
 &=-\sum_\alpha\langle\partial_\alpha v,
                                 \mathscr L\partial_\alpha v\rangle
   +\sum_{\alpha,\gamma}\int V_{0,\alpha\gamma}''
               \partial_\alpha v\,\partial_\gamma v\,\dd\widetilde\mu_q.
 \end{split}
\]
Another integration by parts in the first sum yields
\begin{equation}
 \|\mathscr Lv\|_2^2
 =\sum_{\alpha,\gamma}\|\partial_{\alpha\gamma}v\|_2^2
   +\int(\nabla v)^{\mathsf T}V_0''\nabla v\,
                                          \dd\widetilde\mu_q.
 \label{B:eq:conditional-bochner-identity}
\end{equation}
To extend these identities to $h$, keep the finite dimension fixed during approximation. On every bounded ball the density and its reciprocal are bounded, so
$h\in H^1_{\rm loc}(\dd x)$ and
\[
 \Delta h=\mathscr Lh+\nabla V_0\cdot\nabla h\in L^2_{\rm loc}(\dd x).
\]
For a smooth compactly supported $\theta$, the distributional product rule gives $\Delta(\theta h)=\theta\Delta h+ 2\nabla\theta\cdot\nabla h+h\Delta\theta\in L^2(\dd x)$.
Plancherel's identity gives
\[
 \sum_{\alpha,\gamma}
       \|\partial_{\alpha\gamma}(\theta h)\|_{L^2(\dd x)}^2
 =\int|\xi|^4|\widehat{\theta h}(\xi)|^2\,\dd\xi
 =\|\Delta(\theta h)\|_{L^2(\dd x)}^2<\infty .
\]
It follows that $h\in H^2_{\rm loc}$. Smooth approximation of $\theta h$ converges in $H^2$, and hence in the local generator graph norm. For the cutoff $\chi_R$ used below,
\[
 \mathscr L(\chi_Rv)-\chi_R\mathscr Lv
 =2\nabla\chi_R\cdot\nabla v+v\mathscr L\chi_R,
 \qquad |\mathscr L\chi_R|\le C_n\mathbf1_{\{R\le|x|\le2R\}}.
\]
The first error has norm at most $C R^{-1}\|\mathbf1_{\{|x|\ge R\}}\nabla h\|_2$. The second has norm at most $C_n\|\mathbf1_{\{|x|\ge R\}}h\|_2$. Both tend to zero; also
$(1-\chi_R)\mathscr Lh\to0$ in $L^2$. Thus $\chi_Rh\to h$ in the generator graph norm. Local $H^2$ mollification supplies $v_j\in C_c^\infty$ with $v_j\to h$ and
$\mathscr Lv_j\to\mathscr Lh$ in $L^2$. For their differences, the form identity and the Bochner identity give
\[
 \|\nabla(v_i-v_j)\|_2^2
 \le\|v_i-v_j\|_2\|\mathscr L(v_i-v_j)\|_2,\qquad
 \|D^2(v_i-v_j)\|_2^2
 \le\|\mathscr L(v_i-v_j)\|_2^2.
\]
The limits agree with the distributional derivatives of $h$. Since $V_0''$ is bounded, these estimates put $\nabla h$ in the vector form domain. The commutator formula passes to
the weak equation
\begin{equation}
 (-\mathscr L+V_0'')\nabla h=\nabla G,
 \qquad
 \operatorname{Cov}(F,G)=\langle\nabla F,\nabla h\rangle_{L^2}.
 \label{B:eq:conditional-covariance-representation}
\end{equation}
For the second identity use $-\mathscr Lh=G-\langle G\rangle$ and the definition of the closed form.

The vector operator $-\mathscr L+V_0''$ has form lower bound $a_-/2$. The constant matrix $A$ commutes with $P_t$, so $(-\mathscr L+A)^{-1}=\mathcal K_1$ from
\eqref{B:eq:weighted-free-inverses}. In the weighted component norm, multiplication by $U''$ has norm at most $C_1\delta$. Hence
\begin{equation}
 (-\mathscr L+V_0'')^{-1}
 =\sum_{j=0}^\infty(-\mathcal K_1U'')^j\mathcal K_1,
 \qquad
 \|[(-\mathscr L+V_0'')^{-1}]_{\alpha\gamma}\|_{2\to2}
       \le2K_1e^{-c d(\alpha,\gamma)}.
 \label{B:eq:conditional-vector-inverse}
\end{equation}
For $v$ in a weighted component space, the tail after order $m$ has norm at most
\[
 K_1\sum_{j=m+1}^{\infty}(K_1C_1\delta)^j\|v\|
 \le2K_1\,2^{-m-1}\|v\|.
\]
The sum $u$ satisfies $u=\mathcal K_1(v-U''u)$. In the fixed finite component dimension this is also an equality in the direct sum of the $L^2$ spaces, and the definition of
$\mathcal K_1$ implies $(-\mathscr L+A+U'')u=v$. Coercivity identifies it with the form inverse. For $v$ supported at $\gamma$, center the weight at $k_\gamma$ to obtain
$\|u_\alpha\|_2\le2K_1e^{-c|k_\alpha-k_\gamma|}\|v_\gamma\|_2$. The species contribution to $d$ is at most one and changes the constant by a fixed factor.

Apply this estimate to each component of $\nabla G$ and sum over its coordinate index. The covariance identity gives
\[
 |\operatorname{Cov}(F,G)|
 \le\sum_\alpha\|\partial_\alpha F\|_2
                         \|(\nabla h)_\alpha\|_2
 \le C_\delta\sum_{\alpha,\gamma}e^{-cd(\alpha,\gamma)}
                \|\partial_\alpha F\|_2\|\partial_\gamma G\|_2.
\]
Approximation in $H^1(\widetilde\mu_q)$ passes this estimate to the stated class: covariance is continuous in the two $L^2$ arguments, and each of the finitely many derivative
norms converges. For local polynomials, multiplying by $\chi_R$ gives this approximation directly; the derivative error is bounded by
$|(1-\chi_R)\nabla F|+C R^{-1}\mathbf1_{\{|x|\ge R\}}|F|$, whose $L^2$ norm tends to zero by the coordinate moments.
\end{proof}

\subsubsection{Covariance of the coefficient fluctuations}
\label{B:subsection:coefficient-fluctuation-covariance}

To compare the covariance with the inverse Hessian, we first estimate the signed cubic contributions. With $N=\nabla U$, $\bar x=\langle x\rangle$ and
$\Sigma_{\alpha\gamma}=\operatorname{Cov}_{\widetilde\mu_q}(x_\alpha,x_\gamma)$, the following signed third-moment estimate supplies this improvement.

\begin{lemma}\label{B:lem:local-third-moment-estimate}
If $b,c$ belong to a fixed interaction neighborhood of site $k$, then
\begin{equation}
 |\operatorname{Cov}(x_bx_c,x_j)|
      \le C_\delta k^{-1/2}e^{-c\operatorname{dist}(j,\{b,c\})}.
 \label{B:eq:local-third-covariance}
\end{equation}
For three coordinates $b,c,j$ in that neighborhood, $|\langle x_bx_cx_j\rangle|\le C_\delta k^{-1/2}$.
\end{lemma}
\begin{proof}
The integration by parts identity
\begin{equation}
 \langle(Ax+N)_aF\rangle=\langle\partial_aF\rangle
 \label{B:eq:extended-integration-by-parts}
\end{equation}
is valid for the polynomial functions used here. A cutoff proves it first on bounded sets; the Gaussian density bound makes every boundary term vanish. Apply it to
$F=x_bx_c-\langle x_bx_c\rangle$ and multiply by row $j$ of $C$. The resulting exact formula is
\begin{equation}
 \operatorname{Cov}(x_bx_c,x_j)
 =C_{jb}\bar x_c+C_{jc}\bar x_b
            -\sum_aC_{ja}\operatorname{Cov}(x_bx_c,N_a).
 \label{B:eq:third-moment-score-identity}
\end{equation}
The derivative is $\partial_a(x_bx_c)=\mathbf1_{\{a=b\}}x_c+ \mathbf1_{\{a=c\}}x_b$. At $b=c$ this is $2x_b$ in the single component $a=b$. Hence
\[
 \sum_a\|\partial_a(x_bx_c)\|_2
 \le\|x_b\|_2+\|x_c\|_2\le C_\delta.
\] By
\eqref{B:eq:extended-local-derivatives} and the fourth coordinate moments, each nonzero component of $\nabla N_a$ has $L^2$ norm at most $C_\delta\epsilon_a$.
Lemma~\ref{B:lem:extended-covariance-decay} therefore gives
\[
 |\operatorname{Cov}(x_bx_c,N_a)|
 \le C_\delta\epsilon_a e^{-c\operatorname{dist}(a,\{b,c\})}.
\]
Use $\epsilon_a\le C k^{-1/2}(1+|k_a-k|)^{1/2}$ and \eqref{B:eq:interior-hessian-bounds} in the last sum of \eqref{B:eq:third-moment-score-identity}. For either $\rho=b$ or
$\rho=c$, the contributing sum is bounded by
\[
 C_\delta k^{-1/2}
 \sum_a(1+|k_a-k|)^{1/2}
                  e^{-c d(j,a)-c d(a,\rho)}
 \le C_\delta k^{-1/2}e^{-c'd(j,\rho)}.
\]
This is the estimate for the exponential convolution, with the center shifted a bounded distance to $\rho$. Summing the two choices of $\rho$ gives the required distance to
$\{b,c\}$. The first two terms of \eqref{B:eq:third-moment-score-identity} are bounded separately by $C_\delta k^{-1/2}e^{-cd(j,b)}$ and $C_\delta k^{-1/2}e^{-cd(j,c)}$ using the
mean estimate. Finally $\langle x_bx_cx_j\rangle=\operatorname{Cov}(x_bx_c,x_j) +\langle x_bx_c\rangle\bar x_j$. Its second term has the same order by Cauchy--Schwarz and the
coordinate moment bounds.
\end{proof}

\begin{proposition}\label{B:prop:extended-coefficient-covariance}
The covariance error satisfies
\begin{equation}
 |\Sigma_{\alpha\gamma}-C_{\alpha\gamma}|
       \le C_\delta(k_\alpha k_\gamma)^{-1/2}
                                             e^{-cd(\alpha,\gamma)}.
 \label{B:eq:extended-coefficient-covariance}
\end{equation}
\end{proposition}
\begin{proof}
Let $T^0_{\alpha bc}=\kappa^{-1/2} \partial^3_{\alpha bc}H^q(m)$. Taylor's formula gives
\begin{equation}
 \begin{split}
 N_\alpha(x)&=\tfrac12\sum_{b,c}T^0_{\alpha bc}x_bx_c+R_\alpha(x),\\
 \partial_jR_\alpha(x)&=
 \sum_{b,c}x_bx_c\int_0^1(1-t)
                         \partial^4_{\alpha jbc}U(tx)\,\dd t.
 \end{split}
 \label{B:eq:quadratic-force-expansion}
\end{equation}
There are a bounded number of terms, all in the neighborhood of $\alpha$, and $|T^0_{\alpha bc}|\le C\epsilon_\alpha$. Moreover
$\|\partial_jR_\alpha\|_2\le C_\delta\epsilon_\alpha^2$ for each nonzero component, by the fourth derivative bound for $U$. Applying Lemma~\ref{B:lem:local-third-moment-estimate}
to the first line and Lemma~\ref{B:lem:extended-covariance-decay} to its remainder proves
\begin{equation}
 |\operatorname{Cov}(N_\alpha,x_j)|
       \le C_\delta\epsilon_\alpha^2e^{-cd(\alpha,j)}.
 \label{B:eq:force-coordinate-covariance}
\end{equation}
More explicitly, the finite quadratic sum is bounded by
\[
 \frac12\sum_{b,c}|T^0_{\alpha bc}|\,
                 |\operatorname{Cov}(x_bx_c,x_j)|
 \le C_\delta\epsilon_\alpha^2e^{-c'd(\alpha,j)}.
\]
For the remainder, the covariance estimate and $\partial_\rho x_j=\mathbf1_{\{\rho=j\}}$ give
\[
 |\operatorname{Cov}(R_\alpha,x_j)|
 \le C_\delta\sum_{\rho:d(\alpha,\rho)\le R_*}
       e^{-cd(\rho,j)}\|\partial_\rho R_\alpha\|_2
 \le C_\delta\epsilon_\alpha^2e^{-c'd(\alpha,j)}.
\]
The bounded number and displacement of the derivative indices account for the change in the exponential constant.

Apply \eqref{B:eq:extended-integration-by-parts} to $F=x_j-\bar x_j$. In matrix form this is $A\Sigma+\operatorname{Cov}(N,x)=I$. Thus
\[
 \Sigma_{\alpha j}-C_{\alpha j}
       =-\sum_a C_{\alpha a}\operatorname{Cov}(N_a,x_j).
\]
The inverse entry bound, \eqref{B:eq:force-coordinate-covariance}, and $k_\alpha/k_a\le1+|k_\alpha-k_a|$ give
$|\Sigma_{\alpha j}-C_{\alpha j}| \le C_\delta k_\alpha^{-1}e^{-c'd(\alpha,j)}$. Both matrices are symmetric. Applying the estimate to the entry $(j,\alpha)$ yields the bound with
$k_j^{-1}$. Consequently
\[
 |\Sigma_{\alpha j}-C_{\alpha j}|
 \le C_\delta\min\{k_\alpha^{-1},k_j^{-1}\}
                         e^{-c'd(\alpha,j)}
 \le C_\delta(k_\alpha k_j)^{-1/2}e^{-c'd(\alpha,j)}.
\]
This is \eqref{B:eq:extended-coefficient-covariance}.
\end{proof}

\subsubsection{Covariance of the logarithmic ratios}
\label{B:subsection:logarithmic-ratio-covariance}

The logarithmic ratios are nonlinear functions of the coefficient fluctuations. Their covariance comparison combines the linear covariance estimate with the signed third-moment
bound; estimating only absolute cubic moments would lose the stated order.

\begin{proposition}\label{B:prop:extended-logarithmic-covariance}
For the functions $f_k$ in \eqref{B:eq:extended-logarithmic-functions},
\begin{equation}
 |\operatorname{Cov}_{\widetilde\mu_q}(f_i,f_j)-B_{ij}|
       \le C_\delta(ij)^{-1/2}e^{-c|i-j|}.
 \label{B:eq:extended-logarithmic-covariance}
\end{equation}
\end{proposition}
\begin{proof}
Write $f_i=v_i+q_i+r_i$, where
\[
 v_i=L_ix,\qquad
 q_i=\frac{\sqrt i}{2}
       \left(\frac{x_{x,i}^2}{d_{x,i}^2}
                         -\frac{x_{y,i}^2}{d_{y,i}^2}\right).
\]
For $\alpha\in\{(x,i),(y,i)\}$, Taylor's formula applied to $\partial_\alpha f_i$ gives
\[
 \partial_\alpha r_i(x)
 =\sum_{b,c\in\{(x,i),(y,i)\}}x_bx_c
       \int_0^1(1-t)\partial_{\alpha bc}^3f_i(tx)\,\dd t.
\]
The other derivative components are zero. The global third derivative bound therefore yields $|\nabla r_i|\le C_\delta i^{-1}(|x_{x,i}|+|x_{y,i}|)^2$. Consequently
\[
 \|\nabla v_i\|_2\le C_0,\qquad
 \|\nabla q_i\|_2\le C_0i^{-1/2},\qquad
 \|\nabla r_i\|_2\le C_\delta i^{-1}.
\]
Expanding the covariance into its nine terms gives
\[
 \begin{split}
 \operatorname{Cov}(f_i,f_j)-B_{ij}
 ={}&L_i(\Sigma-C)L_j^{\mathsf T}
       +\operatorname{Cov}(v_i,q_j)+\operatorname{Cov}(q_i,v_j)\\
 &+\operatorname{Cov}(q_i,q_j)
       +\operatorname{Cov}(v_i,r_j)+\operatorname{Cov}(r_i,v_j)\\
 &+\operatorname{Cov}(q_i,r_j)+\operatorname{Cov}(r_i,q_j)
       +\operatorname{Cov}(r_i,r_j).
 \end{split}
\]
Proposition~\ref{B:prop:extended-coefficient-covariance} bounds the first term by $C_\delta(ij)^{-1/2}e^{-c|i-j|}$. For example, one of the mixed terms is the finite sum
\[
 \operatorname{Cov}(v_i,q_j)
 =\frac{\sqrt j}{2}
   \sum_{\alpha\in\{(x,i),(y,i)\}}L_{i\alpha}
   \left[
    \frac{\operatorname{Cov}(x_\alpha,x_{x,j}^2)}{d_{x,j}^2}
   -\frac{\operatorname{Cov}(x_\alpha,x_{y,j}^2)}{d_{y,j}^2}
   \right].
\]
Both covariances involve repeated quadratic indices covered by the preceding third moment estimate. Each coefficient in $q_j$ is $O(j^{-1/2})$. Equation
\eqref{B:eq:local-third-covariance} therefore gives
\[
 |\operatorname{Cov}(v_i,q_j)|\le C_\delta j^{-1}e^{-c|i-j|},
 \qquad
 |\operatorname{Cov}(q_i,v_j)|\le C_\delta i^{-1}e^{-c|i-j|}.
\]
All three gradients at site $i$ are supported on $\{(x,i),(y,i)\}$, and the analogous statement holds at site $j$. Applying Lemma~\ref{B:lem:extended-covariance-decay} to each
pair, with the three gradient bounds above, gives
\[
 \begin{aligned}
 |\operatorname{Cov}(q_i,q_j)|
   &\le C_\delta(ij)^{-1/2}e^{-c|i-j|},\\
 |\operatorname{Cov}(v_i,r_j)|
   &\le C_\delta j^{-1}e^{-c|i-j|},\qquad
 |\operatorname{Cov}(r_i,v_j)|
   \le C_\delta i^{-1}e^{-c|i-j|},\\
 |\operatorname{Cov}(q_i,r_j)|
   &\le C_\delta i^{-1/2}j^{-1}e^{-c|i-j|},\qquad
 |\operatorname{Cov}(r_i,q_j)|
   \le C_\delta i^{-1}j^{-1/2}e^{-c|i-j|},\\
 |\operatorname{Cov}(r_i,r_j)|
   &\le C_\delta i^{-1}j^{-1}e^{-c|i-j|}.
 \end{aligned}
\]
The last three prefactors are at most $(ij)^{-1/2}$ because $i,j\ge1$. For the other powers use
\[
 (i^{-1}+j^{-1})e^{-c|i-j|}
 =(ij)^{-1/2}(\sqrt{i/j}+\sqrt{j/i})e^{-c|i-j|}
 \le C_c(ij)^{-1/2}e^{-c|i-j|/2}.
\]
This proves the displayed bound.
\end{proof}

\subsubsection{Expansion of the means}
\label{B:subsection:conditional-mean-expansion}

The mean correction is determined by the third derivatives of the conditional energy contracted with its inverse Hessian. Keeping that contraction explicit allows the two
logarithmic means to be compared before their difference is summed.

Define the correction in the original coefficient coordinates by
\begin{equation}
 b^q_\alpha=\frac1{2\kappa}\sum_{\rho,\tau}
               \partial^3_{\alpha\rho\tau}H^q(m)C_{\rho\tau},
 \qquad s^q=-Cb^q,\qquad \mu^q=M^{-1}s^q.
 \label{B:eq:conditional-laplace-correction}
\end{equation}
Every sum in $b^q_\alpha$ is over its fixed interaction neighborhood in $\mathscr S(J)$.

\begin{proposition}\label{B:prop:extended-mean-expansion}
Uniformly in $\alpha$,
\begin{equation}
 \begin{split}
 \langle x_\alpha\rangle
       &=\sqrt\kappa\,s^q_\alpha+O_\delta(k_\alpha^{-3/2}),\\
 \langle l_\alpha(x_\alpha)\rangle
       &=s^q_\alpha/m_\alpha
             -C_{\alpha\alpha}/(2\kappa m_\alpha^2)
                       +O_\delta(k_\alpha^{-2}).
 \end{split}
 \label{B:eq:extended-mean-expansion}
\end{equation}
\end{proposition}
\begin{proof}
For coordinates $b,c,j$ in the neighborhood of site $k$, the preceding results imply
\begin{equation}
 |\langle x_bx_c\rangle-C_{bc}|\le C_\delta/k,
 \qquad |\langle x_bx_cx_j\rangle|\le C_\delta k^{-1/2},
 \qquad
 \left\langle\Big(\sum_{b\text{ in the neighborhood}}|x_b|\Big)^4
                         \right\rangle\le C_\delta.
 \label{B:eq:local-raw-moments}
\end{equation}
The first bound is $\Sigma_{bc}-C_{bc}+\bar x_b\bar x_c$; the second is Lemma~\ref{B:lem:local-third-moment-estimate}; the third follows from the bounded number of coordinates and
their fourth moments.

In addition to $T^0$ in \eqref{B:eq:quadratic-force-expansion}, put $T^1_{\alpha bcj}=\kappa^{-1}\partial^4_{\alpha bcj}H^q(m)$. All derivatives at zero agree for the two
potentials. Taylor's formula for $N_\alpha$ through degree three is
\begin{equation}
 \begin{split}
 N_\alpha(x)&=\tfrac12\sum_{b,c}T^0_{\alpha bc}x_bx_c
   +\tfrac16\sum_{b,c,j}T^1_{\alpha bcj}x_bx_cx_j+R_\alpha^{[4]}(x),\\
 R_\alpha^{[4]}(x)&=\frac1{3!}\int_0^1(1-t)^3
                          D^4N_\alpha(tx)[x,x,x,x]\,\dd t.
 \end{split}
 \label{B:eq:cubic-force-expansion}
\end{equation}
The coefficients $T^0,T^1$ have sizes $C\epsilon_\alpha,C\epsilon_\alpha^2$, respectively. The global fifth derivative bound for $U$ gives
$|R_\alpha^{[4]}(x)|\le C_\delta\epsilon_\alpha^3 (\sum_{d(\alpha,b)\le R_*}|x_b|)^4$. Taking expectations in \eqref{B:eq:cubic-force-expansion} and applying
\eqref{B:eq:local-raw-moments} yields
\[
 \begin{split}
 \langle N_\alpha\rangle-\tfrac12\sum_{b,c}T^0_{\alpha bc}C_{bc}
 ={}&\tfrac12\sum_{b,c}T^0_{\alpha bc}
                       (\langle x_bx_c\rangle-C_{bc})\\
 &+\tfrac16\sum_{b,c,j}T^1_{\alpha bcj}\langle x_bx_cx_j\rangle
                       +\langle R_\alpha^{[4]}\rangle
 =O_\delta(\epsilon_\alpha^3).
 \end{split}
\]
The three orders are, respectively, $\epsilon_\alpha\epsilon_\alpha^2$, $\epsilon_\alpha^2\epsilon_\alpha$ and $\epsilon_\alpha^3$. The deterministic vector in the leading
contraction is exactly $\sqrt\kappa\,b^q$, since
\[
 \tfrac12\sum_{b,c}T^0_{\alpha bc}C_{bc}
 =\frac1{2\sqrt\kappa}
       \sum_{b,c}\partial^3_{\alpha bc}H^q(m)C_{bc}
 =\sqrt\kappa\,b^q_\alpha .
\]
Write its expectation error as $e_\alpha$, with $|e_\alpha|\le C_\delta k_\alpha^{-3/2}$. Then $A\bar x+\sqrt\kappa b^q+e=0$ and $s^q=-Cb^q$ give $\bar x=\sqrt\kappa s^q-Ce$. In
components,
\[
 \bar x_\alpha=\sqrt\kappa\,s^q_\alpha
       +O_\delta\left(\sum_b e^{-cd(\alpha,b)}k_b^{-3/2}\right)
 =\sqrt\kappa\,s^q_\alpha+O_\delta(k_\alpha^{-3/2})
\]
by \eqref{B:eq:exponential-site-sums}.

The derivatives of $l_\alpha$ at zero are those of $\log(1+x_\alpha/d_\alpha)$, and its fourth derivative is bounded by $C_\delta d_\alpha^{-4}$ globally. Hence
\[
 l_\alpha(x_\alpha)=\frac{x_\alpha}{d_\alpha}
    -\frac{x_\alpha^2}{2d_\alpha^2}
    +\frac{x_\alpha^3}{3d_\alpha^3}+r_\alpha(x_\alpha),
 \qquad |r_\alpha(t)|\le C_\delta |t|^4/d_\alpha^4.
\]
The linear mean error is $O(k_\alpha^{-3/2}/d_\alpha)=O(k_\alpha^{-2})$. The error in the quadratic term is $O(k_\alpha^{-1}/d_\alpha^2)$. The cubic expectation is
$O(k_\alpha^{-1/2}/d_\alpha^3)$ and the fourth remainder is $O(d_\alpha^{-4})$. All three are also $O(k_\alpha^{-2})$, proving the second estimate. In particular, the four
contributions to this error bound are
\[
 C_\delta\left(
 \frac{k_\alpha^{-3/2}}{d_\alpha}
 +\frac{k_\alpha^{-1}}{d_\alpha^2}
 +\frac{k_\alpha^{-1/2}}{d_\alpha^3}
 +\frac1{d_\alpha^4}\right)
 \le C'_\delta k_\alpha^{-2}.
\]
Their constants depend on the fixed extension parameter and the local derivative bounds.
\end{proof}

\subsubsection{Hessians at constant coefficients}
\label{B:subsection:constant-coefficient-hessians}

We now compute the leading corrections to the two logarithmic means at constant coefficients. A periodic realization removes the finite boundary while retaining the local Hessian
and third derivatives.

For $0<t\le1$ let $\mathbb H_t$ be the Hessian at the constant configuration $x_j=y_j=\phi(t)$ of the coefficient Hamiltonian with logarithmic coefficient $t$ at both species, with
each coordinate direction multiplied by $\phi(t)$. More precisely, on a periodic cycle with $N>8R_*+8$ sites set
\[
 H_{N,t}^{\rm per}(z)=\operatorname{Tr}V(B_N^{\rm per}(z)
              B_N^{\rm per}(z)^{\mathsf T})
                  -t\sum_{\sigma,j}\log z_{\sigma,j}.
\]
For $|h|,|l|\le R_*$, using indices modulo $N$, put
\[
 \begin{split}
 [H_t(h)]_{\sigma\tau}
 &=\phi(t)^2\partial_{\sigma,0}\partial_{\tau,-h}
                      H_{N,t}^{\rm per}(\phi(t)\mathbf1),\\
 \mathbb T_t[(\sigma,0),(\tau,h),(\upsilon,l)]
 &=\phi(t)^3\partial_{\sigma,0}\partial_{\tau,h}
       \partial_{\upsilon,l}H_{N,t}^{\rm per}(\phi(t)\mathbf1).
 \end{split}
\]
The periodic bidiagonal matrix has diagonal $x_j$, subdiagonal $-y_j$, and entry $-y_{N-1}$ in the periodic corner. The expansion as a sum over closed walks in
Section~\ref{B:sec:general-potentials-and-convex-reference-measures} shows that these derivatives are independent of $N$ in this range. Extend them by translation and by zero
beyond the interaction range. Write
\begin{equation}
 (\mathbb H_tv)_k=\sum_{p\in\mathbb Z}H_t(k-p)v_p,
 \qquad G_t(k-p)=(\mathbb H_t^{-1})_{kp}.
 \label{B:eq:constant-hessian-definition}
\end{equation}
These are operators on $\ell^2(\mathbb Z;\mathbb R^2)$.

\begin{lemma}\label{B:lem:constant-hessian-bounds}
For constants independent of $0<t\le1$,
\begin{equation}
 ctI\le\mathbb H_t\le CtI,\quad
 \|H_t(h)\|+|\mathbb T_t(0,h,l)|\le Ct,
 \quad
 \|\partial_tH_t(h)\|+|\partial_t\mathbb T_t(0,h,l)|\le C,
 \label{B:eq:constant-hessian-coefficients}
\end{equation}
and
\begin{equation}
 \|G_t(h)\|\le Ct^{-1}e^{-c|h|},\qquad
 \|\partial_tG_t(h)\|\le Ct^{-2}e^{-c|h|}.
 \label{B:eq:constant-inverse-hessian-bounds}
\end{equation}
\end{lemma}
\begin{proof}
A polynomial monomial of degree $2m$, differentiated two or three times and then multiplied by the corresponding factors $\phi(t)$, has degree $2m$. Its value is a fixed constant
times $\phi(t)^{2m}$. The logarithmic contribution is $t$ on the Hessian diagonal and $-2t$ when all three tensor indices coincide. Since $\phi(t)=\sqrt t\,v_V(t)$ with
$v_V=\phi/\sqrt t$ smooth and positive by the proof of Lemma~\ref{B:lem:bessel-coordinate-regularity}, $|\phi(t)^{2m}|\le Ct$ and $|\partial_t\phi(t)^{2m}|\le C$. There are a
bounded number of terms. These observations prove the entry bounds in \eqref{B:eq:constant-hessian-coefficients}. Convexity of the polynomial part gives
$\mathbb H_t\ge c_*\phi(t)^2I\ge ctI$. Bounded row sums prove the upper operator bound, initially on finitely supported vectors and then on $\ell^2$. Choose constants $c_2,C_2>0$
with $c_2I\le t^{-1}\mathbb H_t\le C_2I$, independent of $t$, and set $Q_t=I-\mathbb H_t/(2C_2t)$. Then $\|Q_t\|\le1-c_2/(2C_2)=:\rho_2<1$. On $\ell^2$,
\[
 \mathbb H_t^{-1}=(2C_2t)^{-1}\sum_{j=0}^{\infty}Q_t^j,\qquad
 (Q_t^j)_{0h}=0\quad\hbox{when }jR_*<|h|.
\]
The sum of its remaining geometric tail is at most $C t^{-1}\rho_2^{\lceil|h|/R_*\rceil}$, proving the bound for each entry of the inverse. This construction also proves bounded
invertibility on the full line.

For a fixed $t>0$ and $s$ close to $t$, the inverse identity is
\[
 \mathbb H_s^{-1}-\mathbb H_t^{-1}
 =-\mathbb H_s^{-1}(\mathbb H_s-\mathbb H_t)\mathbb H_t^{-1}.
\]
Finite range and the entry derivative bounds imply operator norm differentiability of $\mathbb H_t$. The displayed identity gives continuity of its inverse near $t$ and, on
division by $s-t$, the derivative $-\mathbb H_t^{-1}\mathbb H_t'\mathbb H_t^{-1}$. In entries this is
\[
 \partial_tG_t(h)=-\sum_{u,v\in\mathbb Z}
                   G_t(u)\,\partial_tH_t(v-u)\,G_t(h-v).
\]
Only $|v-u|\le R_*$ contributes, and every term has norm at most $C t^{-2}e^{-c|u|-c|h-v|}$. Therefore
\[
 \sum_{\substack{u,v\in\mathbb Z\\|v-u|\le R_*}}
       C t^{-2}e^{-c|u|-c|h-v|}
 \le C't^{-2}\sum_{u\in\mathbb Z}e^{-c|u|-c|h-u|}
 \le C''t^{-2}e^{-c'|h|}.
\]
For fixed $u$ there are at most $2R_*+1$ choices of $v$, and $|h-v|\ge|h-u|-R_*$. These two facts give the factor $(2R_*+1)e^{cR_*}$ in the first inequality. For the second,
$|u|+|h-u|\ge|h|$ gives
\[
 \sum_{u\in\mathbb Z}e^{-c|u|-c|h-u|}
 \le e^{-c|h|/2}\sum_{u\in\mathbb Z}e^{-c|u|/2}
 =\left(1+\frac2{e^{c/2}-1}\right)e^{-c|h|/2}.
\]
The bound is independent of $t,h$, apart from the displayed factors. It proves absolute convergence of the differentiated sum and the $t^{-2}$ estimate in
\eqref{B:eq:constant-inverse-hessian-bounds}.
\end{proof}

\subsubsection{Comparison of inverse Hessians}
\label{B:subsection:comparison-of-conditional-inverse-hessians}

We compare the constant coefficient Hessian with the finite conditional Hessian in relative coordinates, which account for the varying coefficient scale. Decay of the inverse gives
separate estimates for slow variation and for the two boundaries.

In the coordinates scaled by $M$, the Hessian, its inverse, and the third derivatives are
\begin{equation}
 \mathsf A=MAM,\qquad \mathsf C=M^{-1}CM^{-1},\qquad
 \mathsf T_{\alpha\rho\tau}
       =m_\alpha m_\rho m_\tau\partial^3_{\alpha\rho\tau}H^q(m).
 \label{B:eq:conditional-relative-matrices}
\end{equation}
For $u,v\in J$, subscripts $u,v$ on these matrices denote the $2\times2$ blocks. Write $B_\nu(u)=e^{-\nu(u-r)}+e^{-\nu(R-u)}$. The minimizer estimates imply the following
coefficient bound, for some fixed
$\nu>0$:
\begin{equation}
 \|\mathsf A_{uv}-H_{t_k}(u-v)\|
 \le C_0\left[t_u\{u^{-1}+B_\nu(u)\}
                            +\frac{|u-k|}{n}\right]
           \quad(|u-v|\le R_*).
 \label{B:eq:conditional-hessian-coefficient-comparison}
\end{equation}
To verify it first set $k=u$. Each relative polynomial Hessian entry is a fixed sum of monomials of degree $2m$ in the minimizer and its boundary coordinates. In a fixed
neighborhood of $u$, those coordinates equal $\phi(t_u)$ times $1+O(u^{-1}+B_\nu(u))$ by \eqref{B:eq:interior-minimizer-bounds}; the term $e^{-cr}$ is at most $C/n$ after
increasing $n_0$. For a monomial with $2m$ factors, put $a_j=\phi(t_u)(1+\epsilon_j)$ and $\widehat a_j=\phi(t_u)$. The exact telescoping identity
\[
 \prod_{j=1}^{2m}a_j-\prod_{j=1}^{2m}\widehat a_j
 =\sum_{j=1}^{2m}(a_j-\widehat a_j)
               \prod_{l<j}a_l\prod_{l>j}\widehat a_l
\]
bounds its difference by $C\phi(t_u)^{2m}(u^{-1}+B_\nu(u))$. Here the factors $1+\epsilon_j$ are uniformly bounded, and $2m\le2\deg V$ is fixed. For an entry involving boundary
coordinates, $B_\nu(u)$ is bounded below by a fixed positive number; their comparability with $\phi(t_u)$ gives the same inequality. The relative logarithmic Hessian is exactly
$c_{\sigma,u}/\kappa$, whose difference from $t_u$ is $O(n^{-1})$. Since $\phi(t_u)^{2m}\le Ct_u$, this proves the bound at $k=u$. For the change of parameter, write
\[
 H_{t_u}(u-v)-H_{t_k}(u-v)
 =\int_{t_k}^{t_u}\partial_tH_t(u-v)\,\dd t .
\]
Its norm is at most $C|u-k|/n$.

For three derivative indices at $p,p+h,p+l$ within the interaction range and in $J$, multiplication by their three minimizing coordinates restores the degree of each differentiated
polynomial monomial to $2m$. The displayed product identity bounds those terms. The logarithmic third derivative after this multiplication equals $-2c_{\sigma,p}/\kappa$ when all
indices coincide and equals zero otherwise; its constant coefficient value is $-2t_p$. Thus
\[
 \begin{split}
 &|\mathsf T_{(\sigma,p),(\tau,p+h),(\upsilon,p+l)}
       -\mathbb T_{t_k}[(\sigma,0),(\tau,h),(\upsilon,l)]|\\
 &\qquad\le
 C_0\left[t_p\{p^{-1}+B_\nu(p)\}+\frac{|p-k|}{n}\right].
 \end{split}
\]
The last term follows by integrating $\partial_t\mathbb T_t$. Near a conditional boundary, $B_\nu(p)$ is bounded below by a fixed positive constant, which also bounds terms
involving its boundary coordinates.

\begin{proposition}\label{B:prop:conditional-inverse-hessian-comparison}
There is $c>0$ such that, with $\varepsilon_{n,k}=k^{-1}+e^{-c(k-r)}+e^{-c(R-k)}$,
\begin{equation}
 \begin{split}
 \|\mathsf C_{kp}\|&\le
                  \frac{C_0}{\sqrt{t_kt_p}}e^{-c|k-p|},\\
 \|\mathsf C_{kp}-G_{t_k}(k-p)\|&\le
        \frac{C_0}{\sqrt{t_kt_p}}e^{-c|k-p|}\varepsilon_{n,k}.
 \end{split}
 \label{B:eq:conditional-inverse-row-comparison}
\end{equation}
For fixed $|h|,|l|\le R_*$ and $p+h,p+l\in J$, if $|p-k|\le k/2$,
\begin{equation}
 \|\mathsf C_{p+h,p+l}-G_{t_k}(h-l)\|
 \le\frac{C_0}{t_k}\left[
      \frac{1+|p-k|}{k}+e^{-c(p-r)}+e^{-c(R-p)}\right].
 \label{B:eq:conditional-local-inverse-comparison}
\end{equation}
For all $p,k$ the left side of the last inequality is at most $C_0(t_p^{-1}+t_k^{-1})$.
\end{proposition}
\begin{proof}
The first bound follows from $m_{\sigma,k}\asymp\sqrt{t_k}$ and \eqref{B:eq:interior-hessian-bounds}. To compare inverses, restrict the kernel $G_{t_k}(u-v)$ to $u,v\in J$ and call
the resulting matrix $F$. Multiplying by $\mathsf A$, and using the inverse identity on the whole integer lattice, gives, entry by entry,
\begin{equation}
 \begin{split}
 \mathsf C_{uv}-G_{t_k}(u-v)
 ={}&\sum_{a,b\in J}\mathsf C_{ua}
         [H_{t_k}(a-b)-\mathsf A_{ab}]G_{t_k}(b-v)\\
 &+\sum_{\substack{a\in J,\ b\notin J}}
         \mathsf C_{ua}H_{t_k}(a-b)G_{t_k}(b-v).
 \end{split}
 \label{B:eq:conditional-inverse-residual}
\end{equation}
To compute $(I-\mathsf AF)_{av}$ for $a,v\in J$, split the identity on the whole integer lattice as
\[
 \begin{split}
 \delta_{av}I_2
 &=\sum_{b\in J}H_{t_k}(a-b)G_{t_k}(b-v)
   +\sum_{b\notin J}H_{t_k}(a-b)G_{t_k}(b-v),\\
 (I-\mathsf AF)_{av}
 &=\sum_{b\in J}[H_{t_k}(a-b)-\mathsf A_{ab}]G_{t_k}(b-v)
   +\sum_{b\notin J}H_{t_k}(a-b)G_{t_k}(b-v).
 \end{split}
\]
Multiplication on the left by $\mathsf C=\mathsf A^{-1}$ gives $\mathsf C(I-\mathsf AF)=\mathsf C-F$, which proves \eqref{B:eq:conditional-inverse-residual}. All terms with
$|a-b|>R_*$ vanish, and the remaining sums over the whole integer lattice converge absolutely.

Take $u=k,v=p$. Divide the bound for a term of the first sum by $(t_kt_p)^{-1/2}$. Equations \eqref{B:eq:conditional-hessian-coefficient-comparison} and
\eqref{B:eq:constant-inverse-hessian-bounds} leave the coefficient
\[
 C_0\sqrt{p/a}
 \left[\frac1k+\frac akB_\nu(a)+\frac{|a-k|}{k}\right]
                    e^{-\nu(|k-a|+|b-p|)}.
\]
This calculation uses $(t_a/t_k)a^{-1}=k^{-1}$ in the first local term. The ratios $p/a$ and $a/k$ are bounded by fixed powers of $1+|k-a|+|a-p|$. The fixed displacement
$|a-b|\le R_*$ changes only a constant in the exponential. Summing over the bounded number of $b$ and applying \eqref{B:eq:boundary-exponential-convolution} bounds this sum by
\[
 C_0(t_kt_p)^{-1/2}e^{-c|k-p|}
        [k^{-1}+e^{-c(k-r)}+e^{-c(R-k)}].
\]
In the second sum of \eqref{B:eq:conditional-inverse-residual}, $a$ lies within $R_*$ of a boundary and $b$ within $R_*$ of $a$. The factors $t_k$ from $H_{t_k}$ and $t_k^{-1}$
from $G_{t_k}$ cancel. After division by $(t_kt_p)^{-1/2}$ its terms are at most $C_0\sqrt{p/a}\,e^{-\nu(|k-a|+|b-p|)}$. For a left boundary term, $|k-a|+|b-p|\ge |k-p|-R_*$ and
$|k-a|\ge k-r-R_*$. Put $S=|k-a|+|b-p|$. Since $a\ge1$ and $|a-b|\le R_*$, $p/a\le1+|p-a|\le1+S+R_*$. Hence
\[
 \sqrt{p/a}\,e^{-\nu S}
 \le C_{\nu,R_*}e^{-3\nu S/4},\qquad
 C_{\nu,R_*}=\sup_{s\ge0}\sqrt{1+s+R_*}\,e^{-\nu s/4}<\infty.
\]
The inequalities $S\ge|k-p|-R_*$ and $S\ge k-r-R_*$ give $e^{-3\nu S/4}\le e^{\nu R_*/2} e^{-\nu|k-p|/4}e^{-\nu(k-r)/4}$. This proves the bound with $c=\nu/4$. For a right boundary
term, $R-a\le R_*$ and $|k-a|\ge R-k-R_*$. Together with $|k-a|+|b-p|\ge|k-p|-R_*$, the preceding estimate for $\sqrt{p/a}$ and the inequalities $S\ge|k-p|-R_*$, $S\ge R-k-R_*$
give
\[
 \sqrt{p/a}\,e^{-\nu S}
 \le C_{\nu,R_*}e^{\nu R_*/2}
                 e^{-\nu|k-p|/4}e^{-\nu(R-k)/4}.
\]
At each boundary point there are at most a fixed multiple of $R_*^2$ pairs $(a,b)$, independently of $|J|$. This proves \eqref{B:eq:conditional-inverse-row-comparison}.

For the local comparison put $u=p+h,v=p+l$ and decompose
\[
 \mathsf C_{uv}-G_{t_k}(h-l)
 =[\mathsf C_{uv}-G_{t_u}(h-l)]
                      +[G_{t_u}(h-l)-G_{t_k}(h-l)].
\]
The first term is bounded by $C_0t_p^{-1}\varepsilon_{n,p}$, with a change in constant for the fixed offsets. When $|p-k|\le k/2$, all parameters between $t_u$ and $t_k$ are
comparable to $t_k$ for large $n$. The derivative bound in \eqref{B:eq:constant-inverse-hessian-bounds} bounds the second term by
\[
 C_0|t_u-t_k|t_k^{-2}
       \le C_0t_k^{-1}(1+|p-k|)/k.
\]
This proves \eqref{B:eq:conditional-local-inverse-comparison}. For the unrestricted bound estimate the two inverses separately by $C_0/t_p$ and $C_0/t_k$.
\end{proof}

\subsubsection{Cancellation between the two logarithmic means}
\label{B:subsection:mean-logarithmic-ratio-cancellation}

We compare the leading corrections to the two logarithmic means before taking absolute values. The constant coefficient symmetry gives their cancellation, and the finite Hessian
estimates bound the resulting error.

Retain $\mu^q=M^{-1}s^q$ and define
\begin{equation}
 \mathfrak b=Mb^q,\qquad
 \mathfrak b_\alpha=\frac1{2\kappa}
                 \sum_{\rho,\tau}\mathsf T_{\alpha\rho\tau}
                                      \mathsf C_{\rho\tau},
 \qquad \mu^q=-\mathsf C\mathfrak b.
 \label{B:eq:relative-third-derivative-contraction}
\end{equation}
For the constant coefficients set
\begin{equation}
 \begin{split}
 (\mathfrak b_t^{\rm fr})_\sigma
 &=\frac1{2\kappa}\sum_{\tau,\upsilon}\sum_{h,l\in\mathbb Z}
   \mathbb T_t[(\sigma,0),(\tau,h),(\upsilon,l)]
                                      [G_t(h-l)]_{\tau\upsilon},\\
 \mu_t^{\rm fr}&=-\sum_{h\in\mathbb Z}G_t(h)\mathfrak b_t^{\rm fr}.
 \end{split}
 \label{B:eq:constant-mean-contraction}
\end{equation}
The first sum is finite and the second converges absolutely by \eqref{B:eq:constant-inverse-hessian-bounds}. Finally put
\begin{equation}
 a_k^{\rm lead}=\mu^q_{y,k}-\mu^q_{x,k}
       -\frac1{2\kappa}
                  (\mathsf C_{yk,yk}-\mathsf C_{xk,xk}).
 \label{B:eq:leading-logarithmic-ratio-mean}
\end{equation}

\begin{proposition}\label{B:prop:logarithmic-mean-cancellation}
For the same type of constants as above,
\begin{equation}
 |s^q_{\sigma,k}|\le\frac{C_0}{\sqrt\kappa\sqrt k},\qquad
 |a_k^{\rm lead}|\le\frac{C_0}{k}\varepsilon_{n,k},\qquad
 \sum_{k=j}^v|a_k^{\rm lead}|\le C_0/j
                   \quad(r\le j\le v\le R).
 \label{B:eq:leading-logarithmic-mean-bounds}
\end{equation}
\end{proposition}
\begin{proof}
The local tensor $\mathsf T$ has size $Ct_p$ at site $p$, and its contracted inverse blocks have size $C/t_p$. Thus
\begin{equation}
 |\mathfrak b_{\sigma,p}|+|\mathfrak b_{t_k}^{\rm fr}|\le C_0/\kappa.
 \label{B:eq:relative-contraction-bound}
\end{equation}
For a contraction whose three indices lie in $J$, subtract its constant coefficient value using
\[
 \mathsf T_{\alpha\rho\tau}\mathsf C_{\rho\tau}
 -\mathbb T_{\alpha\rho\tau}^{\,t_k}G_{\rho\tau}^{\,t_k}
 =
 (\mathsf T_{\alpha\rho\tau}-\mathbb T_{\alpha\rho\tau}^{\,t_k})
                                      \mathsf C_{\rho\tau}
 +\mathbb T_{\alpha\rho\tau}^{\,t_k}
                       (\mathsf C_{\rho\tau}-G_{\rho\tau}^{\,t_k}).
\]
Here the superscript $t_k$ denotes the already defined arrays at constant coefficients, translated to the site $p$ of $\alpha$. For $|p-k|\le k/2$, the difference of third
derivatives is bounded by
\[
 C_0t_k\left[(1+|p-k|)/k+B_c(p)\right]
\]
by the calculation of third derivatives following \eqref{B:eq:conditional-hessian-coefficient-comparison}. The inverse difference has the reciprocal factor $t_k^{-1}$ by
\eqref{B:eq:conditional-local-inverse-comparison}. The unchanged tensor and inverse are bounded by $Ct_k$ and $C/t_k$, respectively. There are a fixed number of contraction terms.
Therefore
\begin{equation}
 |\mathfrak b_{\sigma,p}-(\mathfrak b_{t_k}^{\rm fr})_\sigma|
 \le\frac{C_0}{\kappa}
           \left[(1+|p-k|)/k+B_c(p)\right]
           \quad(|p-k|\le k/2).
 \label{B:eq:relative-contraction-comparison}
\end{equation}
If an index within the interaction range lies outside $J$, each omitted term evaluated at constant coefficients has size $C_0/\kappa$, and $B_c(p)\ge e^{-cR_*}$. This proves the
same inequality there. Outside $|p-k|\le k/2$ use \eqref{B:eq:relative-contraction-bound}.

At site $k$ the difference of the relative means is exactly
\[
 \begin{split}
 \mu^q_k-\mu_{t_k}^{\rm fr}
 ={}&-\sum_{p\in J}[\mathsf C_{kp}-G_{t_k}(k-p)]\mathfrak b_p\\
 &-\sum_{p\in J}G_{t_k}(k-p)
                        [\mathfrak b_p-\mathfrak b_{t_k}^{\rm fr}]
       +\sum_{p\notin J}G_{t_k}(k-p)\mathfrak b_{t_k}^{\rm fr}.
 \end{split}
\]
For the first sum, \eqref{B:eq:conditional-inverse-row-comparison} and \eqref{B:eq:relative-contraction-bound} give
\[
 \frac{C_0\varepsilon_{n,k}}{\kappa t_k}
              \sum_{p\in J}\sqrt{k/p}\,e^{-c|k-p|}
       \le\frac{C_0\varepsilon_{n,k}}{\kappa t_k}.
\]
For the second sum with $|p-k|\le k/2$ use \eqref{B:eq:relative-contraction-comparison} and Lemma~\ref{B:lem:exponential-convolution-bounds}. This gives the same bound. Its
remaining part is at most
\[
 \frac{C_0}{\kappa t_k}\sum_{|p-k|>k/2}e^{-c|p-k|}
       \le\frac{C_0e^{-c'k}}{\kappa t_k}
       \le\frac{C_0}{\kappa t_k k}.
\]
The omitted exterior sum is at most $C_0(\kappa t_k)^{-1}[e^{-c(k-r)}+e^{-c(R-k)}]$. Thus
\begin{equation}
 |\mu^q_k-\mu_{t_k}^{\rm fr}|
       \le\frac{C_0}{\kappa t_k}\varepsilon_{n,k},
 \qquad |\mu^q_{\sigma,k}|+\kappa^{-1}\mathsf C_{\sigma k,\sigma k}
       \le C_0/k.
 \label{B:eq:relative-mean-comparison}
\end{equation}
The second inequality follows directly from the inverse row bound, \eqref{B:eq:relative-contraction-bound}, and $\kappa t_k=\beta k$. Multiplying its mean part by
$m_{\sigma,k}\asymp\sqrt{k/n}$ proves the bound on $s^q$.

To prove equality of the constant coefficient contractions, we use a reflection of the bipartite dilation path. Place its row vertex $\mathrm R_j$ at $2j$ and its column vertex
$\mathrm C_j$ at $2j+1$. Reflection of vertices by $v\mapsto2-v$ exchanges the coefficient locations $x_j$ and $-y_{-j}$. Conjugating additionally by the diagonal signs $(-1)^v$
changes the sign of every edge. The induced coefficient map is consequently
\[
 (\mathcal Jx)_j=y_{-j},\qquad (\mathcal Jy)_j=x_{-j}.
\]
On the periodic cycle these are orthogonal conjugations of the dilation. They preserve the traces of its even powers and the sum of its equal logarithmic terms. With
$\jmath(\sigma,j)=(\bar\sigma,-j)$, where $\bar x=y$ and $\bar y=x$, differentiation at the fixed constant configuration gives
\[
 \mathbb T_t[\jmath\alpha,\jmath\rho,\jmath\tau]
        =\mathbb T_t[\alpha,\rho,\tau],\qquad
 (\mathbb H_t)_{\jmath\alpha,\jmath\rho}
        =(\mathbb H_t)_{\alpha,\rho}.
\]
The same permutation commutes with the bounded inverse. For every specified local derivative one may take the cycle longer than the interaction neighborhood. The resulting
identities therefore hold for the tensors on the whole integer
lattice. Let $P=\left(\begin{smallmatrix}0&1\\1&0\end{smallmatrix}\right)$. Differentiation at the constant
configuration gives
\[
 H_t(-h)=PH_t(h)P,\quad G_t(-h)=PG_t(h)P,
 \quad [G_t(0)]_{xx}=[G_t(0)]_{yy}.
\]
In the finite sum defining the first contraction, make the bijection $(\tau,\upsilon,h,l)\mapsto (\bar\tau,\bar\upsilon,-h,-l)$. It gives
\[
 \begin{split}
 (\mathfrak b_t^{\rm fr})_{\bar\sigma}
 &=\frac1{2\kappa}\sum_{\tau,\upsilon,h,l}
   \mathbb T_t[(\bar\sigma,0),(\bar\tau,-h),(\bar\upsilon,-l)]
                       [G_t(l-h)]_{\bar\tau\bar\upsilon}\\
 &=\frac1{2\kappa}\sum_{\tau,\upsilon,h,l}
   \mathbb T_t[(\sigma,0),(\tau,h),(\upsilon,l)]
                       [G_t(h-l)]_{\tau\upsilon}
 =(\mathfrak b_t^{\rm fr})_\sigma .
 \end{split}
\]
Thus $P\mathfrak b_t^{\rm fr}=\mathfrak b_t^{\rm fr}$. Since the inverse row is absolutely summable,
\[
 P\mu_t^{\rm fr}
 =-\sum_hPG_t(h)P\,P\mathfrak b_t^{\rm fr}
 =-\sum_hG_t(-h)\mathfrak b_t^{\rm fr}
 =\mu_t^{\rm fr}.
\]
In particular its two coordinates are equal.

Substituting these equal quantities into the finite expression gives
\[
 \begin{split}
 a_k^{\rm lead}
 ={}&(\mu^q_{y,k}-(\mu_{t_k}^{\rm fr})_y)
      -(\mu^q_{x,k}-(\mu_{t_k}^{\rm fr})_x)\\
 &-\frac1{2\kappa}
   \left[(\mathsf C_{yk,yk}-[G_{t_k}(0)]_{yy})
        -(\mathsf C_{xk,xk}-[G_{t_k}(0)]_{xx})\right].
 \end{split}
\]
The first two differences are bounded by $C(\kappa t_k)^{-1}\varepsilon_{n,k}$ from \eqref{B:eq:relative-mean-comparison}. Each inverse difference is bounded by
$Ct_k^{-1}\varepsilon_{n,k}$ from \eqref{B:eq:conditional-inverse-row-comparison} at $p=k$. Since $\kappa t_k=\beta k$, their sum is at most $C_0k^{-1}\varepsilon_{n,k}$. Finally,
\[
 \sum_{k=j}^vk^{-2}\le C/j,\qquad
 \sum_{k=j}^v\frac{B_c(k)}k
 \le\frac1j\left(\sum_{h=0}^\infty e^{-ch}
                      +\sum_{h=0}^\infty e^{-ch}\right)\le C_c/j.
\]
These inequalities prove the final summation bound.
\end{proof}

Combining \eqref{B:eq:extended-mean-expansion} with \eqref{B:eq:leading-logarithmic-ratio-mean} gives, with
$b_k:=\langle f_k\rangle/\sqrt k$,
\begin{equation}
 b_k=a_k^{\rm lead}+O_\delta(k^{-2}),\qquad
 |b_k|\le C_\delta k^{-2}+C_\delta k^{-1}B_c(k),\qquad
 \sum_{k=j}^v|b_k|\le C_\delta/j.
 \label{B:eq:extended-logarithmic-mean-sums}
\end{equation}
Here $b_k$ denotes a scalar expectation; the coefficient vector $b^q$ remains defined by \eqref{B:eq:conditional-laplace-correction}. The inverse matrices and all tensor
contractions in these estimates are on $J$, while observations may subsequently be restricted to $I$.

\subsubsection{Comparison with the original conditional distribution}
\label{B:subsection:comparison-of-conditional-distributions}

Using the local corrections for the extended potential, we construct coefficients with their original conditional distribution. Agreement of the potentials on the localization set
gives a common submeasure. We couple this submeasure on the diagonal and estimate the remaining mass separately.

Let $\mu_q^*$ be the law of $\sqrt\kappa(Z-m)$ for $Z\sim\mu_q$. It has density proportional to $e^{-V_*}$ on $x_\alpha>-d_\alpha$. Use the event $E_\delta$ in
\eqref{B:eq:conditional-agreement-event}. For fixed $\delta$ set
$q_n=\min\{1,C_\delta n e^{-c_\delta r}\}$.

\begin{proposition}\label{B:prop:conditional-density-coupling}
There is a Borel family of couplings of $X^*\sim\mu_q^*$ and $\widetilde X\sim\widetilde\mu_q$ such that
\begin{equation}
 \Pr_q\{X^*\ne\widetilde X\text{ or }X^*\notin E_\delta\}\le q_n.
 \label{B:eq:conditional-equality-coupling}
\end{equation}
For $2\le p\le c_{\rm rc}r$, their logarithmic ratios obey
\begin{equation}
 \|\eta_k-\widetilde\eta_k\|_p
       \le C_\delta\sqrt{p/k}\,q_n^{1/(2p)},
 \qquad
 \widetilde\eta_k=f_k(\widetilde X)/\sqrt k.
 \label{B:eq:conditional-coupling-moments}
\end{equation}
\end{proposition}
\begin{proof}
Equation~\eqref{B:eq:conditional-relative-deviation} bounds $\mu_q^*(E_\delta^c)$ by $C_\delta ne^{-c_\delta r}$. For the extension, \eqref{B:eq:extended-coordinate-moments} and
Lipschitz concentration give
\[
 \widetilde\mu_q\{|x_\alpha|>\delta d_\alpha\}
 \le2e^{-c_\delta k_\alpha}
\]
for large $n$: its mean is $O(k_\alpha^{-1/2})$ and $\delta d_\alpha\ge c\delta\sqrt{k_\alpha}$. The complement of $E_\delta$ is the union of the events
$\{|x_\alpha|>\delta d_\alpha\}$. There are at most $2n$ coordinates and $k_\alpha\ge r$, so
\[
 \widetilde\mu_q(E_\delta^c)
 \le2\sum_{\alpha\in\mathscr S(J)}e^{-c_\delta k_\alpha}
 \le4n e^{-c_\delta r}.
\]

Let $a_q=\mu_q^*(E_\delta)$, $\widetilde a_q=\widetilde\mu_q(E_\delta)$, and $\lambda_q=\min\{a_q,\widetilde a_q\}$. Put $Z_q^*=\int_{\{x_\alpha>-d_\alpha\}}e^{-V_*(x)}\,\dd x$ and
$Z_{q,E}=\int_{E_\delta}e^{-V_0(x)}\,\dd x>0$. The box $E_\delta$ lies in the original domain and its two potentials agree. Hence
\[
 a_q=Z_{q,E}/Z_q^*,\qquad
 \widetilde a_q=Z_{q,E}/\widetilde Z_q,\qquad
 \nu_q(\dd x)=Z_{q,E}^{-1}
                   \mathbf1_{E_\delta}(x)e^{-V_0(x)}\,\dd x.
\]
These identities prove that $\nu_q$ is the common conditional probability on $E_\delta$. Set $\rho_q=\mu_q^*-\lambda_q\nu_q$ and $\widetilde\rho_q=\widetilde\mu_q-\lambda_q\nu_q$.
Both are positive measures of mass $1-\lambda_q$. An explicit coupling is
\begin{equation}
 \pi_q=\lambda_q\,\operatorname{diag}_{\#}\nu_q
           +\frac{\rho_q\otimes\widetilde\rho_q}{1-\lambda_q}
                 \quad(\lambda_q<1),
 \label{B:eq:explicit-conditional-density-coupling}
\end{equation}
where $\operatorname{diag}(x)=(x,x)$; if $\lambda_q=1$ keep only the diagonal term. For a Borel set $D$ and $\lambda_q<1$, its first marginal is
\[
 \pi_q(D\times\mathbb R^{\mathscr S(J)})
 =\lambda_q\nu_q(D)
       +\frac{\rho_q(D)(1-\lambda_q)}{1-\lambda_q}
 =\mu_q^*(D).
\]
Its second marginal is $\lambda_q\nu_q(D)+\widetilde\rho_q(D)=\widetilde\mu_q(D)$. For $\lambda_q=1$ the two marginals equal $\nu_q$. The diagonal term is supported on
$x=y\in E_\delta$, so
\[
 \pi_q\{x\ne y\text{ or }x\notin E_\delta\}
 \le1-\lambda_q
 =\max\{\mu_q^*(E_\delta^c),\widetilde\mu_q(E_\delta^c)\}
 \le q_n .
\]
The constants in $q_n$ may be enlarged to cover the finitely many initial indices. All densities and event integrals are Borel functions of $q$, as are $\rho_q/(1-\lambda_q)$ and
$\widetilde\rho_q/(1-\lambda_q)$ on $\{q:\lambda_q<1\}$. When $\lambda_q=1$, these two kernels are unused and may be assigned a fixed probability measure. Thus
\eqref{B:eq:explicit-conditional-density-coupling} is a Borel family. It can be used without replacing a previously sampled extended coordinate. For $\lambda_q<1$ define
\[
 a_q^{\rm acc}(y)=\frac{\lambda_q}{\widetilde a_q}
                         \mathbf1_{E_\delta}(y),\qquad
 K_q^{\rm orig}(y,\dd x)=a_q^{\rm acc}(y)\delta_y(\dd x)
        +(1-a_q^{\rm acc}(y))\frac{\rho_q(\dd x)}{1-\lambda_q}.
\]
Here $0\le a_q^{\rm acc}\le1$ and $\widetilde\mu_q(\dd y)a_q^{\rm acc}(y)=\lambda_q\nu_q(\dd y)$. The other part has $y$-marginal $\widetilde\rho_q$, so
\[
 \widetilde\mu_q(\dd y)K_q^{\rm orig}(y,\dd x)
 =\lambda_q\operatorname{diag}_{\#}\nu_q(\dd x,\dd y)
       +\frac{\rho_q(\dd x)\widetilde\rho_q(\dd y)}{1-\lambda_q}.
\]
For $\lambda_q=1$ use the identity kernel. The original coordinate therefore has law $\mu_q^*$, while every joint law involving the previously sampled coordinate $y$ is preserved
by integrating the kernel to one. On $E_\delta$, the acceptance factor is the ratio of the common normalized mass. It may differ from one because agreement of the potentials on
this set does not imply agreement of their normalizing constants.

On the diagonal event the logarithms agree. The original ratios have the moment bound \eqref{B:eq:interior-log-moments}. For the extended ones, $l_\alpha(0)=0$ and
$|l_\alpha'|\le C_\delta/d_\alpha$ give $\|\widetilde\eta_k\|_{2p}\le C_\delta\sqrt{p/k}$. On the remaining event $F$ use H\"older:
\[
 \|\eta_k-\widetilde\eta_k\|_p
 =\|(\eta_k-\widetilde\eta_k)\mathbf1_F\|_p
 \le(\|\eta_k\|_{2p}+\|\widetilde\eta_k\|_{2p})
                                        \Pr_q(F)^{1/(2p)}.
\]
This proves \eqref{B:eq:conditional-coupling-moments} for one coupling and every stated order. The choice $8c_{\rm rc}\le c_{\rm mom}$ permits all the moments here and the fixed
multiples used later. For each fixed $C_1$, choose $n_0(C_1)$ so that $8C_1\log(\me+n)\le c_{\rm mom}r$ for every $n\ge n_0(C_1)$.
\end{proof}

\begin{corollary}\label{B:cor:original-conditional-mean-expansion}
Under the original conditional distribution,
\begin{equation}
 \begin{split}
 \mathbb E_q(Z_\alpha-m_\alpha)
     &=s^q_\alpha+O_\delta(\kappa^{-1/2}k_\alpha^{-3/2}),\\
 \mathbb E_q\zeta_\alpha
     &=s^q_\alpha/m_\alpha-C_{\alpha\alpha}/(2\kappa m_\alpha^2)
                                      +O_\delta(k_\alpha^{-2}),\\
 \mathbb E_q\eta_k&=a_k^{\rm lead}+O_\delta(k^{-2}),\qquad
 \sum_{k=j}^v\bigl(|\mathbb E_q\eta_k|+|b_k|\bigr)\le C_\delta/j.
 \end{split}
 \label{B:eq:original-conditional-mean-expansion}
\end{equation}
\end{corollary}
\begin{proof}
The coupling of Proposition~\ref{B:prop:conditional-density-coupling} also compares coordinates. Their second moments are bounded by \eqref{B:eq:interior-concentration} after
scaling, and by \eqref{B:eq:extended-coordinate-moments}. Therefore
\[
 |\mathbb E_qX_\alpha^*-\langle x_\alpha\rangle|
       \le C_\delta q_n^{1/2},\qquad
 |\mathbb E_q\zeta_\alpha-\langle l_\alpha\rangle|
       \le C_\delta k_\alpha^{-1/2}q_n^{1/2}.
\]
For each fixed $D>0$,
\[
 q_n^{1/2}\le C_\delta^{1/2}n^{1/2}e^{-c_\delta r/2}
                  \le C_{D,\delta}n^{-D}
\]
for $n\ge n_0(D,\delta)$: choose this index so that $(D+1/2)\log n-c_\delta r/2$ is below the logarithm of the required constant. With $D=2$ and $k_\alpha\le n$,
$n^{-2}\le k_\alpha^{-3/2}$ and $k_\alpha^{-1/2}n^{-2}\le k_\alpha^{-2}$. These are the two required error orders. The coordinate comparison is divided by $\sqrt\kappa$ only after
this estimate in fluctuation coordinates. Substituting \eqref{B:eq:extended-mean-expansion} and dividing the coordinate identity by $\sqrt\kappa$ proves the first two formulas.
Subtraction between the species proves the third. For the partial sums, the preceding identities give
\[
 \sum_{k=j}^v|\mathbb E_q\eta_k|
 \le\sum_{k=j}^v|a_k^{\rm lead}|
       +C_\delta\sum_{k=j}^vk^{-2}\le C_\delta/j.
\]
The same bound for $\sum_{k=j}^v|b_k|$ is \eqref{B:eq:extended-logarithmic-mean-sums}; adding the two estimates proves the bound for the sum of the means. Both means refer to the
separately specified conditional laws, and the estimate compares those laws using the coupling above.
\end{proof}

\subsubsection{Estimates for the coupled coefficient arrays}
\label{B:subsection:conditional-estimates-for-coupled-arrays}

Integrating the conditional bounds against the original exterior distribution gives a failure probability at each size. Their sum is finite, so Borel--Cantelli applies to the
chosen joint construction of all sizes.

Fix at each size one deterministic weight array from Corollary~\ref{B:cor:conditional-logarithmic-probabilities}. Let $E_n$ have the original exterior law and suppose
$\mathcal L(Z_{J_n}\mid E_n=e)=\mu_{n,J_n}^{c_n(e)}$. For admissible $q$ set $b_n(q)=\mathbb E_q\mathcal M_n(w)$ and set $b_n(q)=0$ otherwise. This is Borel: the conditional
density is a Borel kernel, the positive minimizer is continuous in $q$, and integration of the nonnegative Borel energy is measurable.

\begin{corollary}\label{B:cor:conditional-estimates-under-realization}
For these arrays, the following limits hold almost surely, and the displayed inequality holds for every $n\ge n_0(\omega)$ for a finite random index $n_0(\omega)$:
\begin{equation}
 \begin{gathered}
 n[\mathcal M_n(w)-b_n(c_n(E_n))]\longrightarrow0,
 \qquad \limsup_n n\mathcal M_n(w)\le C_\lambda,\\
 \max_{\alpha\in\mathscr S(J_n)}|\zeta_{n,\alpha}|
      \le C\sqrt{\log(\me+n)/r_n}\longrightarrow0
                     \quad(n\ge n_0(\omega)).
 \end{gathered}
 \label{B:eq:conditional-all-size-estimates}
\end{equation}
\end{corollary}
\begin{proof}
Let $\mathcal G_n=\{c_n(E_n)\in\mathcal Q_n^\sharp\}$. Theorem~\ref{B:thm:boundary-and-terminal-probabilities} gives $\sum_n\Pr(\mathcal G_n^c)<\infty$ for the original exterior
laws. For every $\varepsilon>0$, conditional integration and \eqref{B:eq:conditional-centered-energy-probability} give
\[
 \begin{split}
 \Pr\{n|\mathcal M_n-b_n(c_n(E_n))|>\varepsilon\}
 &\le\Pr(\mathcal G_n^c)
   +\int_{c_n(e)\in\mathcal Q_n^\sharp}
       \mu_{n,J_n}^{c_n(e)}\{|\mathcal M_n-\mathbb E_{c_n(e)}
                              \mathcal M_n|>\varepsilon/n\}
                                          \,\mathcal L(E_n)(\dd e)\\
 &\le\Pr(\mathcal G_n^c)
             +(C_\lambda p^{3/2}/\varepsilon)^p n^{-p\delta_\lambda}.
 \end{split}
\]
Choose a fixed integer $p$ with $p\delta_\lambda>2$. It is in the allowed moment range after a deterministic starting index. For each integer $j\ge1$ the displayed probabilities
with $\varepsilon=1/j$ are summable. Their Borel--Cantelli events have a countable intersection on which $n|\mathcal M_n-b_n(c_n(E_n))|\to0$. On this probability-one event,
$\mathcal G_n$ holds for every $n\ge n_0(\omega)$, and
\[
 n\mathcal M_n
 \le n|\mathcal M_n-b_n(c_n(E_n))|
                       +n b_n(c_n(E_n))
 \le o(1)+C_\lambda.
\]

For the maximum, conditional integration with $D=3$ gives
\[
 \Pr\left\{\max_\alpha\sqrt{k_\alpha}|\zeta_{n,\alpha}|
                              >C_3\sqrt{\log(\me+n)}\right\}
 \le\Pr(\mathcal G_n^c)+C_3n^{-3}.
\]
The series is finite. Its Borel--Cantelli event and $k_\alpha\ge r_n$ give $\max_\alpha|\zeta_{n,\alpha}| \le C_3\sqrt{\log(\me+n)/r_n}\to0$. We integrate throughout against the
original exterior marginal, using the equality of the original and reference interior kernels proved in Section~\ref{B:sec:general-potentials-and-convex-reference-measures}. These
integrations require no further change of measure. The argument applies to the given dependence between different sizes. For any fixed countable collection of weights, intersection
of the corresponding events of probability one gives all the estimates for that collection.
\end{proof}

\subsection{Gaussian approximation of the coefficient vector}
\label{B:sec:gaussian-approximation-of-the-coefficient-field}

The coefficient products in the inverse involve partial sums of logarithmic ratios. We therefore couple the Gaussian vector with covariance proportional to the inverse Hessian to
the Brownian increments, with estimates for the partial sums of the coupling error. The variance at frequency zero determines the Brownian time change. The Fourier factorization
separates this variance from a discrete difference, and a covariance correction then gives the required conditional Gaussian distribution. Summing the discrete difference leaves
two boundary values; the other three errors have partial sum variances bounded by $Cj^{-2}$. The estimates for the weighted squared sums include their cross covariances.

\subsubsection{The conditional Gaussian distribution}
\label{B:subsection:conditional-gaussian-distribution}

The conditional mean is $m^q+s^q$. The linear approximations to the logarithmic ratios are formed from $G-m^q-s^q$, so their covariance can be coupled to the Brownian increments
while $s^q$ contributes to the coefficient error.

Use the parameters and matrices of Section~\ref{B:sec:conditional-estimates-on-the-interior-interval}. In particular, $V,a,\beta,A_0$ are fixed, $A_0>20$, $A=\nabla^2H^q(m^q)$,
$C=A^{-1}$ and $M=\operatorname{diag}(m^q)$. Write $m_n=r+b$, $M_n=R-b$, so $I=\{m_n,\ldots,M_n\}$. The matrices $A,C$ act on $\mathscr S(J)$, and $K_q$ is indexed by $J$. The
relative matrices $\mathsf A,\mathsf C$ are defined in \eqref{B:eq:conditional-relative-matrices}. With the correction $s^q$ from \eqref{B:eq:conditional-laplace-correction}, set
\begin{equation}
 \Omega_q=\kappa^{-1}C,
 \qquad \gamma_{n,q}=N(m^q+s^q,\Omega_q),
 \qquad v_-=(-1,1)^{\mathsf T}.
 \label{B:eq:gaussian-coefficient-law}
\end{equation}
Thus $K_q=L_q\Omega_qL_q^{\mathsf T}$, and its entries are
\begin{equation}
 K_q(k,l)=\kappa^{-1}v_-^{\mathsf T}\mathsf C_{kl}v_-.
 \label{B:eq:covariance-log-ratio}
\end{equation}
For a vector $G$ with law $\gamma_{n,q}$, the centered linear observable is $L_q(G-m^q-s^q)$. Its characteristic function is
\[
 \mathbb E_q e^{\ii v^{\mathsf T}L_q(G-m^q-s^q)}
 =\exp\{-\tfrac12v^{\mathsf T}K_qv\},\qquad v\in\mathbb R^J.
\]
The Gaussian distribution is defined on the whole coefficient space. The logarithms of the positive coefficients under $\mu_q$ have the separate estimates in
Section~\ref{B:sec:conditional-estimates-on-the-interior-interval}. Here $L_q$ is the linear map defined there.

Unless stated otherwise, the estimates below are uniform in $q\in\mathcal Q_n^\sharp$ and $n\ge n_0$. Constants denoted by $c,C,C_0$ may depend on $V,a,\beta,A_0$; their value may
change from one line to the next. We increase the deterministic $n_0$ when necessary. The formulas defining the Gaussian vectors will also be given for every finite nonnegative
boundary configuration $q$. For deterministic $q$, $\Pr_q$, $\operatorname{Cov}_q$ and $\operatorname{Var}_q$ refer to the joint law of the Gaussian vectors defined below.
Expectations without a subscript are over the Brownian increments and independent normal variables in the construction.

We shall use the bounds
\begin{equation}
 cI\le A\le CI,\qquad m_{\sigma,k}\asymp\sqrt{k/n},\qquad
 \|\mathsf C_{kl}\|\le
        \frac{C}{\sqrt{t_kt_l}}e^{-c|k-l|},\quad k,l\in J,
 \label{B:eq:covariance-finite-decay}
\end{equation}
proved in Propositions~\ref{B:prop:interior-conditional-geometry} and \ref{B:prop:conditional-inverse-hessian-comparison}. In particular the inverse in this formula is taken before
restricting its entries to $I$. To identify the resulting boundary term, partition the coordinate set into $\mathscr S(I)$ and its complement $O$ in $\mathscr S(J)$. Block
elimination in $AC=I$ gives
\[
 C_{II}=(A_{II}-A_{IO}A_{OO}^{-1}A_{OI})^{-1},\qquad
 C_{II}-A_{II}^{-1}
 =A_{II}^{-1}A_{IO}
 (A_{OO}-A_{OI}A_{II}^{-1}A_{IO})^{-1}A_{OI}A_{II}^{-1}\ge0.
\]
Both Schur complements are positive since $A>0$. The notation $II$ in this display denotes coefficient coordinates, not the scalar ratio indices. Applying the ratio map on those
coordinates gives $K_I^G$; replacing $C_{II}$ by $A_{II}^{-1}$ would change the conditional Gaussian distribution.

\subsubsection{Fourier representation at constant coefficients}
\label{B:subsection:gaussian-fourier-representation}

At constant coefficients, translation invariance turns the Hessian into a matrix Fourier multiplier. This representation is useful both for its inverse and for identifying the
scalar covariance seen by the logarithmic ratios.

Use $\mathbb H_t,H_t(h),G_t(h)$ from \eqref{B:eq:constant-hessian-definition}. For a finitely supported vector $v$ define $\widehat v(\vartheta)=\sum_kv_ke^{\ii k\vartheta}$. Then
\begin{equation}
 \widehat H_t(\vartheta)=\sum_{|h|\le R_*}H_t(h)e^{\ii h\vartheta},
 \qquad
 G_t(h)=\frac1{2\pi}\int_{-\pi}^{\pi}
       e^{-\ii h\vartheta}\widehat H_t(\vartheta)^{-1}\,\dd\vartheta.
 \label{B:eq:covariance-fourier}
\end{equation}
Indeed, for a finitely supported vector,
\[
 \widehat{\mathbb H_tv}(\vartheta)
 =\sum_{k,h}H_t(h)v_{k-h}e^{\ii k\vartheta}
 =\widehat H_t(\vartheta)\widehat v(\vartheta),\qquad
 \frac1{2\pi}\int_{-\pi}^{\pi}|\widehat v(\vartheta)|^2\,\dd\vartheta
 =\sum_k|v_k|^2.
\]
The second identity follows by integrating each exponential $e^{\ii(k-l)\vartheta}$. Completion gives the unitary Fourier map on $\ell^2(\mathbb Z;\mathbb C^2)$. The inverse
multiplier in \eqref{B:eq:covariance-fourier} is bounded by the positive lower bound proved below, and its two products with $\widehat H_t$ equal $I_2$. Thus
$(\mathbb H_t^{-1})_{kl}=G_t(k-l)$. Since $H_t(-h)=H_t(h)^{\mathsf T}$ and the operator is real symmetric, $G_t(-h)=G_t(h)^{\mathsf T}$. Define
\begin{equation}
 \Gamma_t(h)=\beta^{-1}v_-^{\mathsf T}G_t(h)v_-,\qquad
 \widehat\Gamma_t(\vartheta)
      =\beta^{-1}v_-^{\mathsf T}\widehat H_t(\vartheta)^{-1}v_-.
 \label{B:eq:covariance-scalar-symbol}
\end{equation}
The scalar coefficients are real and even.

\begin{lemma}\label{B:lem:analytic-covariance-symbol}
For some $\rho>0$, the inverse symbol is analytic on $|\operatorname{Im}\vartheta|<\rho$ and satisfies
\begin{equation}
 \|\widehat H_t^{-1}\|\le Ct^{-1},\qquad
 \|\partial_t\widehat H_t^{-1}\|\le Ct^{-2}.
 \label{B:eq:analytic-inverse-symbol}
\end{equation}
On the real axis, $ctI_2\le\widehat H_t\le CtI_2$ and $c/t\le\widehat\Gamma_t\le C/t$. Moreover,
\begin{equation}
 \|G_t(h)\|+|\Gamma_t(h)|\le Ct^{-1}e^{-c|h|},\qquad
 \|\partial_tG_t(h)\|+|\partial_t\Gamma_t(h)|
                  \le Ct^{-2}e^{-c|h|}.
 \label{B:eq:covariance-frozen-bounds}
\end{equation}
\end{lemma}
\begin{proof}
Lemma~\ref{B:lem:constant-hessian-bounds} gives $ctI\le\mathbb H_t\le CtI$ and
\[
 \|H_t(h)\|\le Ct\mathbf1_{\{|h|\le R_*\}},\qquad
 \|\partial_tH_t(h)\|\le C\mathbf1_{\{|h|\le R_*\}}.
\]
To obtain the symbol bounds directly, fix $w\in\mathbb C^2$ and put $v_k=e^{-\ii k\vartheta}w$ for $|k|\le N$, and zero otherwise. For $2N+1>R_*$ the exact quadratic form is
\[
 \langle v,\mathbb H_tv\rangle
 =\sum_{|h|\le R_*}(2N+1-|h|)e^{\ii h\vartheta}
                              w^*H_t(h)w.
\]
There are $2N+1-|h|$ pairs with $k-l=h$ in the selected interval. Consequently
\[
 \left|(2N+1)^{-1}\langle v,\mathbb H_tv\rangle
                  -w^*\widehat H_t(\vartheta)w\right|
 \le\frac{Ct|w|^2}{2N+1}\sum_{|h|\le R_*}|h|\longrightarrow0.
\]
Also $|v|^2=(2N+1)|w|^2$. Dividing the lower and upper form bounds by $2N+1$ proves $ctI_2\le\widehat H_t(\vartheta)\le CtI_2$. Inverting this inequality and using $|v_-|^2=2$
gives
\[
 \frac{2}{\beta Ct}\le
 \beta^{-1}v_-^{\mathsf T}\widehat H_t(\vartheta)^{-1}v_-
 \le\frac{2}{\beta ct}.
\]

For real $\vartheta$ and $|y|\le1$,
\[
 \begin{split}
 \|\widehat H_t(\vartheta+\ii y)-\widehat H_t(\vartheta)\|
 &\le Ct\sum_{|h|\le R_*}|e^{-hy}-1|\\
 &\le Ct|y|\sum_{|h|\le R_*}|h|e^{|h|}
 \le C_1t|y|.
 \end{split}
\]
Choose $\rho>0$ such that $C_1\rho/c\le1/2$. If $E$ denotes this difference, then
\[
 (\widehat H_t(\vartheta)+E)^{-1}
 =\sum_{j=0}^\infty(-\widehat H_t(\vartheta)^{-1}E)^j
                                     \widehat H_t(\vartheta)^{-1}.
\]
The series is locally uniformly convergent and has norm at most $2/(ct)$. Differentiating the identity for the inverse gives
\[
 \partial_t\widehat H_t^{-1}
       =-\widehat H_t^{-1}(\partial_t\widehat H_t)\widehat H_t^{-1}.
\]
The finite Fourier sum for $\partial_t\widehat H_t$ is bounded on the strip, proving \eqref{B:eq:analytic-inverse-symbol}.

For $h\ne0$ put $y_h=-(\rho/2)\operatorname{sgn}(h)$. Cauchy's theorem on the rectangle with horizontal sides at heights zero and $y_h$ gives
\[
 G_t(h)=\frac1{2\pi}\int_{-\pi}^{\pi}
 e^{-\ii h(\vartheta+\ii y_h)}
 \widehat H_t(\vartheta+\ii y_h)^{-1}\,\dd\vartheta.
\]
The vertical integrals cancel: their matrix factors are periodic and $e^{-\ii h(\pi+\ii y)}=e^{-\ii h(-\pi+\ii y)}$ for integer $h$. On the displaced contour the exponential has
modulus $e^{hy_h}=e^{-\rho|h|/2}$. Hence $|G_t(h)|\le Ct^{-1}e^{-\rho|h|/2}$. For its derivative, the bound $Ct^{-2}$ is uniform on each compact $t$-interval in $(0,1]$ and on this
contour. Difference quotients may therefore be integrated, giving the identical formula with $\partial_t\widehat H_t^{-1}$ and the bound $Ct^{-2}e^{-\rho|h|/2}$. At $h=0$ both
estimates follow on the real contour. Contracting each formula with $v_-$ and dividing by $\beta$ proves \eqref{B:eq:covariance-frozen-bounds}.
\end{proof}

\subsubsection{The variance at frequency zero}
\label{B:subsection:variance-at-frequency-zero}

The following identity determines the Brownian time change. It is the circulant Hessian calculation, written in relative coordinates, from \cite{rider2019universality}. We use
$\phi,L_V(t),\theta_V$ from \eqref{B:eq:1-12}--\eqref{B:eq:1-15}, and the notation $\ell(t)=L_V(t)$ and $\theta=\theta_V$ used in
Section~\ref{B:sec:general-potentials-and-convex-reference-measures}.

\begin{lemma}\label{B:lem:zero-frequency-variance}
For $0<t\le1$,
\begin{equation}
 \widehat\Gamma_t(0)=\sum_{h\in\mathbb Z}\Gamma_t(h)
       =\frac{2}{\beta\phi(t)\ell(t)}
       =\frac1\beta\frac{\theta'(t)}{\theta(t)}.
 \label{B:eq:gaussian-zero-frequency-identity}
\end{equation}
\end{lemma}
\begin{proof}
Write $\mathcal W_V$ for the polynomial in \eqref{B:eq:constant-coefficient-polynomial}, and use the constants $A_m,B_m,C_m$ from \eqref{B:eq:polynomial-gradient-constants}. A
constant perturbation on a periodic cycle sums all Hessian offsets. Consequently
\begin{equation}
 \widehat H_t(0)
   =\phi(t)^2\nabla^2\mathcal W_V(\phi(t),\phi(t))+tI_2.
 \label{B:eq:covariance-zero-hessian}
\end{equation}
Indeed the polynomial trace is $N\mathcal W_V(x,y)$ on such a configuration, and its second derivative in a constant vector is $N$ times the sum of a row of blocks. The logarithms
contribute $tI_2$ after multiplication of each direction by $\phi(t)$.

The two diagonal entries of the polynomial Hessian agree. Since $B_m-C_m=A_m/(2m-1)$, the eigenvalue of \eqref{B:eq:covariance-zero-hessian} corresponding to $v_-$ is
\[
 t+\sum_{m=1}^{\deg V}\frac{g_m A_m}{2m-1}\phi(t)^{2m}
   =\sum_{m=1}^{\deg V}\frac{2m g_m A_m}{2m-1}\phi(t)^{2m}
   =\phi(t)\ell(t).
\]
For the first equality we used $t=F(\phi(t))$; for the last, the change of variable $u=F(x)$ gives
\[
 \ell(t)=\int_0^{\phi(t)}\frac{F'(x)}x\,\dd x
       =\sum_{m=1}^{\deg V}\frac{2m g_m A_m}{2m-1}
                                             \phi(t)^{2m-1}.
\]
The integral converges at zero since $F'(x)/x$ is bounded there. The vector $v_+=(1,1)^{\mathsf T}$ is the other eigenvector. Its eigenvalue is
\[
 t+\sum_mg_m(B_m+C_m)\phi(t)^{2m}
 =\phi(t)F'(\phi(t)).
\]
Both eigenvalues are positive. With $P_\pm=v_\pm v_\pm^{\mathsf T}/2$, the inverse at frequency zero is
\[
 \widehat H_t(0)^{-1}
 =\frac{P_+}{\phi(t)F'(\phi(t))}
       +\frac{P_-}{\phi(t)\ell(t)}.
\]
Since $P_+v_-=0$ and $P_-v_-=v_-$, its contraction equals $2/(\beta\phi(t)\ell(t))$. Differentiating $\ell$ and $\theta$ gives
\[
 \ell'(t)=\phi(t)^{-1},\qquad
 \frac{\theta'(t)}{\theta(t)}
   =\frac{2\ell'(t)}{\ell(t)}=\frac{2}{\phi(t)\ell(t)}.
\]
Finally $\sum_h|\Gamma_t(h)|<\infty$ by \eqref{B:eq:covariance-frozen-bounds}. Its Fourier series converges uniformly, and evaluation at zero gives
$\sum_h\Gamma_t(h)=\widehat\Gamma_t(0)$.
\end{proof}

For example, if $V(x)=x/2+gx^2$ with $g\ge0$, then
\[
 t=\phi^2+12g\phi^4,\qquad \ell=2\phi+16g\phi^3,\qquad
 \widehat H_t(0)=\phi^2
 \begin{pmatrix}2+32g\phi^2&16g\phi^2\\16g\phi^2&2+32g\phi^2\end{pmatrix}.
\]
Thus its eigenvalue in direction $v_-$ is $\phi\ell$. For $g=0$, $\Gamma_t(h)=(\beta t)^{-1}\mathbf1_{\{h=0\}}$.

\subsubsection{Comparison of covariance entries}
\label{B:subsection:comparison-gaussian-covariance-entries}

The constant coefficient covariance is only a local approximation to the conditional covariance. The estimate below separates the changes in the profile from the two conditional
boundary contributions and gives their decay as the indices separate.

For $k,l\in I$ put $t_{kl}=(k+l)/(2n)$. Choose a positive decay constant $c_1$ smaller than the constants in \eqref{B:eq:conditional-inverse-row-comparison} and
\eqref{B:eq:covariance-frozen-bounds}, and set
$\epsilon_n=e^{-c_1b}+e^{-c_1r}$.

\begin{theorem}\label{B:thm:gaussian-covariance-comparison}
Uniformly for $k,l\in I$,
\begin{equation}
 \begin{split}
 K_q(k,l)&=n^{-1}\Gamma_{t_{kl}}(k-l)+\mathcal E_{n,q}(k,l),\\
 |\mathcal E_{n,q}(k,l)|
 &\le\frac{C}{nt_{kl}}e^{-c|k-l|}
         \left[\frac{1+|k-l|}{1+k\wedge l}+\epsilon_n\right].
 \end{split}
 \label{B:eq:covariance-midpoint}
\end{equation}
The corresponding estimate for $\mathsf C_{kl}-G_{t_{kl}}(k-l)$ has prefactor $C/t_{kl}$.
\end{theorem}
\begin{proof}
Proposition~\ref{B:prop:conditional-inverse-hessian-comparison} gives
\begin{equation}
 \|\mathsf C_{kl}-G_{t_k}(k-l)\|
 \le \frac{C}{\sqrt{t_kt_l}}e^{-c|k-l|}
          [k^{-1}+e^{-c(k-r)}+e^{-c(R-k)}].
 \label{B:eq:covariance-row}
\end{equation}
This estimate includes the second sum in \eqref{B:eq:conditional-inverse-residual}, over $a\in J,b\notin J$. For $k\in I$ the two distances in the last bracket are at least $b$.

Set $u=k\wedge l$, $v=k\vee l$, and $d=v-u$. If $v\le2u$, $t_k,t_l,t_{kl}$ and every point between $t_k$ and $t_{kl}$ lie in $[u/n,2u/n]$. Hence
\[
 \begin{split}
 \|G_{t_k}(k-l)-G_{t_{kl}}(k-l)\|
 &\le\left|\int_{t_{kl}}^{t_k}
                     Ct^{-2}e^{-cd}\,\dd t\right|\\
 &\le \frac{Cd}{n t_{kl}^2}e^{-cd}
 \le \frac{C}{t_{kl}}e^{-cd}\frac{d}{1+u}.
 \end{split}
\]
Also $(t_kt_l)^{-1/2}\le C/t_{kl}$ and $k^{-1}\le2/(1+u)$. Adding this bound to \eqref{B:eq:covariance-row} proves the matrix estimate in this case.

If $v>2u$, then $d>u$. The separate inverse bounds give
\[
 \|\mathsf C_{kl}-G_{t_{kl}}(k-l)\|
 \le\frac{C}{t_{kl}}e^{-cd}
               \left(1+\frac{t_{kl}}{\sqrt{t_kt_l}}\right).
\]
Put $x=d/u>1$. Then
\[
 1+\frac{t_{kl}}{\sqrt{t_kt_l}}
 =1+\frac{2+x}{2\sqrt{1+x}}
 \le 2+\frac x2\le\frac52x,\qquad
 \frac{1+d}{1+u}\ge\frac{d}{2u}=\frac x2.
\]
Thus this factor is at most $5(1+d)/(1+u)$, proving the estimate on the remaining pairs. Finally,
\[
 \begin{split}
 K_q(k,l)-n^{-1}\Gamma_{t_{kl}}(k-l)
 &=\frac1{n\beta}v_-^{\mathsf T}
       [\mathsf C_{kl}-G_{t_{kl}}(k-l)]v_-,\\
 \left|v_-^{\mathsf T}
       [\mathsf C_{kl}-G_{t_{kl}}(k-l)]v_-\right|
 &\le2\|\mathsf C_{kl}-G_{t_{kl}}(k-l)\|.
 \end{split}
\]
This proves the scalar prefactor and error bound. The matrix $\mathsf C$ remains the inverse on $J$; both sums over indices beyond the conditional boundaries are already included
in \eqref{B:eq:covariance-row}.
\end{proof}

\subsubsection{Covariances of partial sums}
\label{B:subsection:covariances-gaussian-partial-sums}

Let $\eta^G$ be a centered Gaussian vector with covariance $K_q$. The following estimates depend only on $K_q$. In particular, for real $c_k$ supported on $\{j,\ldots,M_n\}$,
\[
 \operatorname{Var}_q\left(\sum_kc_k\eta_k^G\right)
       =\sum_{k,l}c_kc_lK_q(k,l).
\]

\begin{corollary}\label{B:cor:gaussian-partial-sum-variance}
For $j\in I$,
\begin{equation}
 \sum_{k,l=j}^{M_n}|\mathcal E_{n,q}(k,l)|
 \le C/j+C\epsilon_n\log(2+n/j),
 \label{B:eq:covariance-error-partial-sums}
\end{equation}
and
\begin{equation}
 \operatorname{Var}_q\!\left(\sum_{k=j}^{M_n}\eta_k^G\right)
 =\frac1\beta\log\frac{z_{n,M_n+1}}{z_{n,j}}+O(j^{-1}),
 \qquad z_{n,k}=\theta(k/n).
 \label{B:eq:gaussian-partial-sum-variance}
\end{equation}
\end{corollary}
\begin{proof}
Assign a pair $(k,l)$ to $u=k\wedge l$ and $d=|k-l|$. There are at most two pairs for each $(u,d)$, and $nt_{kl}\ge u$. Thus
\[
 \begin{split}
 \sum_{k,l=j}^{M_n}|\mathcal E_{n,q}(k,l)|
 &\le C\sum_{u=j}^{M_n}\sum_{d\ge0}e^{-cd}
                          [(1+d)u^{-2}+\epsilon_nu^{-1}]\\
 &\le C\sum_{u=j}^{M_n}(u^{-2}+\epsilon_nu^{-1})
 \le C/j+C\epsilon_n\log(2+n/j).
 \end{split}
\]
For comparable $k,l$, the mean value estimate from the preceding proof, contracted with $v_-$, gives
\[
 n^{-1}|\Gamma_{t_{kl}}(k-l)-\Gamma_{t_k}(k-l)|
                         \le Cd\,u^{-2}e^{-cd}.
\]
For $v>2u$, the separate bounds give at most $Cu^{-1}e^{-cd}$; since $d>u$, this is at most $Cd\,u^{-2}e^{-cd}$. With $q=e^{-c}\in(0,1)$, the two possible orderings of each pair
give
\[
 \frac1n\sum_{k,l=j}^{M_n}
 |\Gamma_{t_{kl}}(k-l)-\Gamma_{t_k}(k-l)|
 \le2C\sum_{u=j}^\infty u^{-2}\sum_{d=0}^\infty d q^d
 \le\frac{4Cq}{(1-q)^2j}.
\]
We used $\sum_{u=j}^\infty u^{-2}\le2/j$ and $\sum_{d\ge0}dq^d=q/(1-q)^2$.

For a fixed $k\in[j,M_n]$, the omitted part of the full $l$-sum is bounded by
\[
 \frac1n\sum_{l\notin[j,M_n]}|\Gamma_{t_k}(k-l)|
 \le \frac Ck[e^{-c(k-j)}+e^{-c(M_n-k)}].
\]
Summing its left and right terms separately gives
\[
 \sum_{k=j}^{M_n}\frac{e^{-c(k-j)}}k
 \le\frac1{j(1-e^{-c})},\qquad
 \sum_{k=j}^{M_n}\frac{e^{-c(M_n-k)}}k
 \le\frac1{j(1-e^{-c})}.
\]
The full sum in $l$ equals $\widehat\Gamma_{t_k}(0)$. Thus, with $f(t)=\widehat\Gamma_t(0)$,
\[
 \operatorname{Var}_q\left(\sum_{k=j}^{M_n}\eta_k^G\right)
 =\frac1n\sum_{k=j}^{M_n}f(t_k)
       +O(j^{-1})+O\bigl(\epsilon_n\log(2+n/j)\bigr).
\]
The derivative bound is $|f'(t)|\le Ct^{-2}$. On each cell,
\[
 \frac1n f(k/n)-\int_{k/n}^{(k+1)/n}f(t)\,\dd t
 =-\int_{k/n}^{(k+1)/n}\int_{k/n}^{t}f'(s)\,\dd s\,\dd t.
\]
The integration triangle has area $1/(2n^2)$, so the absolute value is at most $C/k^2$. The sum of these errors is at most $C(j^{-2}+\int_j^\infty x^{-2}\,\dd x)\le C/j$. The
boundary point of the last cell is $(M_n+1)/n$, and therefore
\[
 \int_{j/n}^{(M_n+1)/n}f(t)\,\dd t
 =\frac1\beta\left[\log\theta((M_n+1)/n)-\log\theta(j/n)\right].
\]
Since $b\ge r/4$, we have $\epsilon_n\le2e^{-c_1r/4}$. The relation $r/\log n\to\infty$ gives $\epsilon_n\log n\le C/n\le C/j$ after a deterministic starting index. Combining these
estimates proves \eqref{B:eq:gaussian-partial-sum-variance}.
\end{proof}

If $0\le w_{n,j}\le C_0n^{-\lambda}j^{\lambda-1}$ with $1<\lambda<2$, \eqref{B:eq:covariance-error-partial-sums} also implies
\begin{equation}
 \sum_{j\in I}w_{n,j}\sum_{k,l=j}^{M_n}|\mathcal E_{n,q}(k,l)|\le C/n.
 \label{B:eq:covariance-weighted-error}
\end{equation}
Indeed its two contributions are bounded by $Cn^{-\lambda}\sum_{j\le n}j^{\lambda-2}\le C/n$ and $C\epsilon_n\log n\,n^{-\lambda}\sum_{j\le n}j^{\lambda-1} \le C\epsilon_n\log n$,
respectively.

\subsubsection{Square roots of the scalar Fourier multiplier}
\label{B:subsection:scalar-fourier-square-root}

The square root of the scalar covariance multiplier provides a Gaussian convolution. Subtracting its value at zero frequency leaves a discrete difference, which is the part whose
cumulative sums will reduce to boundary values.

On the real axis set $\widehat a_t(\vartheta)=\widehat\Gamma_t(\vartheta)^{1/2}$, using the positive square root, and put
\begin{equation}
 \begin{split}
 a_t(h)&=\frac1{2\pi}\int_{-\pi}^{\pi}
                   \widehat a_t(\vartheta)e^{-\ii h\vartheta}\,\dd\vartheta,
       \qquad a_t^\circ=\widehat a_t(0),\\
 \widehat b_t(\vartheta)&=
       \frac{\widehat a_t(\vartheta)-a_t^\circ}{1-e^{-\ii\vartheta}},
       \qquad
 b_t(h)=\frac1{2\pi}\int_{-\pi}^{\pi}
                   \widehat b_t(\vartheta)e^{-\ii h\vartheta}\,\dd\vartheta.
 \end{split}
 \label{B:eq:gaussian-fourier-factor-definitions}
\end{equation}
The value of the quotient at zero is given by continuous extension.

\begin{lemma}\label{B:lem:scalar-fourier-factors}
These multipliers extend analytically to a fixed strip. Their Fourier coefficients are real and satisfy
\begin{equation}
 |a_t(h)|+|b_t(h)|\le Ct^{-1/2}e^{-c|h|},\qquad
 |\partial_ta_t(h)|+|\partial_tb_t(h)|\le Ct^{-3/2}e^{-c|h|}.
 \label{B:eq:coupling-root-decay}
\end{equation}
They also satisfy
\begin{equation}
 a_t^\circ=\frac1{\sqrt\beta}
                       \sqrt{\frac{\theta'(t)}{\theta(t)}},\qquad
 a_t(h)=a_t^\circ\mathbf1_{\{h=0\}}+b_t(h)-b_t(h+1).
 \label{B:eq:coupling-coefficient-difference}
\end{equation}
\end{lemma}
\begin{proof}
For real $\vartheta$ and small $y$, absolute convergence of the Fourier series gives
\[
 |\widehat\Gamma_t(\vartheta+\ii y)-\widehat\Gamma_t(\vartheta)|
 \le \sum_h|\Gamma_t(h)|\,|e^{-hy}-1|
 \le \frac{C|y|}{t}\sum_h|h|e^{-c|h|/2}\le C_1|y|/t.
\]
Reduce the strip until this is at most half the lower bound on the real axis. Its real part is then positive, so the principal square root is analytic, with
$c t^{-1/2}\le|\widehat a_t|\le Ct^{-1/2}$. The identity $\partial_t\widehat a_t=(\partial_t\widehat\Gamma_t)/(2\widehat a_t)$ gives $|\partial_t\widehat a_t|\le Ct^{-3/2}$ there.

The only zeros of $1-e^{-\ii\vartheta}$ are $2\pi\mathbb Z$. Near zero,
\[
 \widehat a_t(\vartheta)-\widehat a_t(0)
       =\vartheta\int_0^1\partial_\vartheta
                                  \widehat a_t(s\vartheta)\,\dd s.
\]
Choose a fixed circle radius $r_0>0$ smaller than the distance between two nested strips. Cauchy's formula on $|\zeta-\vartheta|=r_0$ gives
\[
 |\partial_\vartheta\widehat a_t(\vartheta)|
 \le r_0^{-1}\sup_{|\zeta-\vartheta|=r_0}|\widehat a_t(\zeta)|
 \le Ct^{-1/2},\qquad
 |\partial_\vartheta\partial_t\widehat a_t(\vartheta)|
 \le Ct^{-3/2}.
\]
Also $1-e^{-\ii\vartheta}=\ii\vartheta\int_0^1e^{-\ii s\vartheta}\,\dd s$, whose quotient by $\vartheta$ is bounded away from zero in a fixed disk. The quotient defining
$\widehat b_t$ thus has a removable singularity and satisfies $|\widehat b_t|\le Ct^{-1/2}$ there. Applying the displayed derivative estimate to
$\partial_t\widehat a_t(\vartheta)-\partial_t\widehat a_t(0)$ gives $|\partial_t\widehat b_t|\le Ct^{-3/2}$. The principal square root is periodic because
$\widehat\Gamma_t(\vartheta+2\pi)=\widehat\Gamma_t(\vartheta)$ and both values lie in the same right half-plane. Hence these estimates hold around every zero $2\pi\mathbb Z$ of the
denominator. On the complement of their disks in a closed smaller strip over one period, the denominator has a positive minimum.

For any one of $f_t=\widehat a_t,\widehat b_t,\partial_t\widehat a_t, \partial_t\widehat b_t$, and $h\ne0$, displacement by $y_h=-\rho_0\operatorname{sgn}(h)$ within that smaller
strip gives
\[
 \frac1{2\pi}\int_{-\pi}^{\pi}f_t(\vartheta)e^{-\ii h\vartheta}\,\dd\vartheta
 =\frac{e^{-\rho_0|h|}}{2\pi}\int_{-\pi}^{\pi}
                f_t(\vartheta+\ii y_h)e^{-\ii h\vartheta}\,\dd\vartheta.
\]
Periodicity cancels the two vertical integrals. The preceding bounds, and the undisplaced integrals at $h=0$, prove \eqref{B:eq:coupling-root-decay}.

On the real axis $\widehat a_t$ is real and even, so $a_t(h)$ is real and even. Moreover $\widehat b_t(-\vartheta)=\overline{\widehat b_t(\vartheta)}$, which gives real $b_t(h)$.
Evenness also gives $\partial_\vartheta\widehat a_t(0)=0$, hence $\widehat b_t(0)=0$. The absolutely convergent Fourier series satisfies
\[
 (1-e^{-\ii\vartheta})\sum_hb_t(h)e^{\ii h\vartheta}
       =\sum_h[b_t(h)-b_t(h+1)]e^{\ii h\vartheta}.
\]
Comparing coefficients in $\widehat a_t=a_t^\circ+(1-e^{-\ii\vartheta})\widehat b_t$ proves the second identity in \eqref{B:eq:coupling-coefficient-difference}. The first follows
from the positive square root of \eqref{B:eq:gaussian-zero-frequency-identity}. In addition, uniform convergence at $\vartheta=0$ gives $\sum_hb_t(h)=0$. For later covariance
computations, Fourier multiplication and $a_t(-h)=a_t(h)$ give
\[
 \sum_h a_t(h)a_t(h-d)=\sum_h a_t(h)a_t(d-h)=\Gamma_t(d),
 \qquad d\in\mathbb Z.
\]
Both convolutions converge absolutely by \eqref{B:eq:coupling-root-decay}.
\end{proof}

\subsubsection{Brownian increments and a Gaussian convolution}
\label{B:subsection:brownian-increments-gaussian-convolution}

The convolution can be driven by normalized increments of the common Brownian path. Its ratio increments use forward logarithmic lengths; these must be distinguished from the
backward lengths of the operator cells.

Use the sampling positions $z_{n,k}$ and logarithmic times $u_{n,k}=-\log z_{n,k}$ of Section~\ref{B:sec:approximation-by-brownian-motion-and-the-lower-bound}. For $1\le k<n$ set
\begin{equation}
 h^+_{n,k}=u_{n,k}-u_{n,k+1}
       =\log\frac{z_{n,k+1}}{z_{n,k}},\qquad
 h^-_{n,k}=h_{n,k}\quad(k\ge2).
 \label{B:eq:coupling-clocks}
\end{equation}
Thus the first length is used for the ratio at site $k$; the second is the previously defined length of the logarithmic image of cell $k$.

Let $W$ be a standard Brownian motion. Independently of $W$, take standard normal variables $\xi_{n,l}^{\perp}$ for $l\notin I$. Define
\begin{equation}
 \xi_{n,l}=\begin{cases}
 [W(u_{n,l})-W(u_{n,l+1})]/\sqrt{h^+_{n,l}},&l\in I,\\
 \xi_{n,l}^{\perp},&l\notin I.
 \end{cases}
 \label{B:eq:coupling-white-noise}
\end{equation}
Here $l$ ranges over $\mathbb Z$. For each fixed $n$ and every finite $F\subset\mathbb Z$, the sampling intervals used for the Brownian increments for $l\in F\cap I$ have disjoint
interiors. Their normalized increments have covariance $I_{F\cap I}$. The auxiliary variables are mutually independent and independent of the Brownian path. Thus, for real $s_l$,
\[
 \mathbb E\exp\!\left(\ii\sum_{l\in F}s_l\xi_{n,l}\right)
      =\exp\!\left(-\tfrac12\sum_{l\in F}s_l^2\right).
\]
This proves independence and the standard normal laws for each row. Put
\begin{equation}
 \eta_{n,k}^{\rm fr}=n^{-1/2}\sum_{h\in\mathbb Z}
                         a_{t_k}(h)\xi_{n,k-h},\qquad k\in I.
 \label{B:eq:coupling-frozen-field}
\end{equation}
For each $n,k$,
\[
 \mathbb E\sum_h|a_{t_k}(h)\xi_{n,k-h}|
   =\sqrt{2/\pi}\sum_h|a_{t_k}(h)|<\infty.
\]
Thus the series converges absolutely almost surely. Its tails converge in every finite $L^p$ by Minkowski's inequality and \eqref{B:eq:coupling-root-decay}. These conclusions hold
simultaneously for all the countably many indices used below. The covariance is
\begin{equation}
 K_n^{\rm fr}(k,l)
    =\frac1n\sum_{h\in\mathbb Z}a_{t_k}(h)a_{t_l}(h+l-k).
 \label{B:eq:coupling-frozen-covariance}
\end{equation}
For finite sums, the only surviving products have equal noise indices: $k-h=l-g$ means $g=h+l-k$. The $L^2$ tail estimate
\[
 \left\|n^{-1/2}\sum_{|h|>H}a_{t_k}(h)\xi_{n,k-h}\right\|_2^2
   =\frac1n\sum_{|h|>H}|a_{t_k}(h)|^2
   \le\frac Ck e^{-cH}
\]
permits passage to the covariance limit. More generally, any finite linear combination of the convolution coordinates is an $L^2$ limit of centered Gaussian variables. Their
characteristic functions converge to the centered Gaussian characteristic function with the limiting variance. The vector is therefore jointly Gaussian, with
\eqref{B:eq:coupling-frozen-covariance}.

\begin{lemma}\label{B:lem:gaussian-convolution-covariance}
For $k,l\in I$,
\begin{equation}
 |K_n^{\rm fr}(k,l)-n^{-1}\Gamma_{t_{kl}}(k-l)|
 \le\frac{C}{nt_{kl}}e^{-c|k-l|}
                           \frac{1+|k-l|}{1+k\wedge l}.
 \label{B:eq:coupling-frozen-error}
\end{equation}
\end{lemma}
\begin{proof}
Write $d=k-l$, $u=k\wedge l$, $v=k\vee l$, and $t_0=t_{kl}$. If $v\le2u$, the parameters are comparable and
\[
 |a_{t_k}(h)-a_{t_0}(h)|
 \le C|d|n^{-1}t_0^{-3/2}e^{-c|h|}.
\]
The convolution with both factors evaluated at the same parameter, computed in Lemma~\ref{B:lem:scalar-fourier-factors} gives $\sum_h a_{t_0}(h)a_{t_0}(h-d)=\Gamma_{t_0}(d)$. For
the other factor, the derivative estimate is
\[
 |a_{t_l}(h-d)-a_{t_0}(h-d)|
 \le \left|\int_{t_0}^{t_l} Cs^{-3/2}e^{-c|h-d|}\,\dd s\right|
 \le C|d|n^{-1}t_0^{-3/2}e^{-c|h-d|}.
\]
Subtract the two products in \eqref{B:eq:coupling-frozen-covariance} by writing
\[
 a_{t_k}(h)a_{t_l}(h-d)-a_{t_0}(h)a_{t_0}(h-d)
 =[a_{t_k}(h)-a_{t_0}(h)]a_{t_l}(h-d)
  +a_{t_0}(h)[a_{t_l}(h-d)-a_{t_0}(h-d)].
\]
Since $|h|+|h-d|\ge|d|$,
\[
 \sum_h e^{-c(|h|+|h-d|)}
 \le e^{-c|d|/2}\sum_he^{-c|h|/2}\le C_ce^{-c|d|/2}.
\]
After multiplication by $1/n$, each of the two product differences has sum at most $C|d|n^{-2}t_0^{-2}e^{-c'|d|}$. Since $nt_0=(k+l)/2\ge u\ge(1+u)/2$, this is bounded by
\[
 \frac C{nt_0}e^{-c'|d|}
                  \frac{1+|d|}{1+u}.
\]
When $d=0$ the two parameters equal $t_0$ and the covariance error is exactly zero. For $v>2u$, direct summation of the coefficient bounds gives
\[
 |K_n^{\rm fr}(k,l)|+n^{-1}|\Gamma_{t_0}(d)|
 \le \frac C{nt_0}e^{-c'|d|}
          \left(1+\frac{t_0}{\sqrt{t_kt_l}}\right).
\]
The final factor is at most $5(1+|d|)/(1+u)$ by the calculation for noncomparable indices in Theorem~\ref{B:thm:gaussian-covariance-comparison}. This proves
\eqref{B:eq:coupling-frozen-error} for every pair.
\end{proof}

\subsubsection{Square roots of covariance matrices}
\label{B:subsection:square-roots-covariance-matrices}

To correct the covariance while retaining exponential decay of the matrix entries, we need an estimate for each entry of the difference of two matrix square roots. Let $P,Q$ be
real symmetric matrices indexed by a finite interval $\mathcal I\subset\mathbb N$. Suppose, with fixed positive constants,
\begin{equation}
 \begin{gathered}
 c_0I\le P,Q\le C_0I,\qquad
 |P_{kl}|+|Q_{kl}|\le C_1e^{-c_1|k-l|},\\
 |(P-Q)_{kl}|\le C_1e^{-c_1|k-l|}
                           \frac{1+|k-l|}{1+k\wedge l}.
 \end{gathered}
 \label{B:eq:matrix-square-root-hypotheses}
\end{equation}

\begin{lemma}\label{B:lem:localized-matrix-square-root}
Under \eqref{B:eq:matrix-square-root-hypotheses},
\begin{equation}
 |(P^{1/2}-Q^{1/2})_{kl}|
       \le\frac{C_2e^{-c_2|k-l|}}{\sqrt{(1+k)(1+l)}}.
 \label{B:eq:coupling-root-comparison}
\end{equation}
The constants depend only on those in \eqref{B:eq:matrix-square-root-hypotheses}.
\end{lemma}
\begin{proof}
Fix a column $l_0$ and set $D_\rho=\operatorname{diag}(e^{\rho|k-l_0|})$. For $0<\rho\le c_1/2$,
\[
 |(D_\rho P D_\rho^{-1}-P)_{kj}|
 \le C_1\rho|k-j|e^{-(c_1-\rho)|k-j|}.
\]
Its absolute row and column sums are at most $C_1\rho\sum_{h\in\mathbb Z}|h|e^{-c_1|h|/2}=C_3\rho$. The Schur test therefore gives operator norm at most $C_3\rho$. Choose $\rho$
with $C_3\rho\le c_0/2$. For $s\ge0$, a Neumann series about $P+sI$ gives
\[
 \|(D_\rho P D_\rho^{-1}+sI)^{-1}\|
 \le\frac1{c_0+s-C_3\rho}\le\frac2{c_0+s}.
\]
Taking its $(k,l_0)$ entry, and then taking $l_0=l$, yields
\begin{equation}
 |(P+sI)^{-1}_{kl}|+|(Q+sI)^{-1}_{kl}|
                  \le\frac4{c_0+s}e^{-\rho|k-l|}.
 \label{B:eq:coupling-resolvent-decay}
\end{equation}
This estimates the inverse entries by conjugating with $D_\rho$; compare the exponential conjugation method of \cite{combes1973asymptotic}. For band matrices, exponential inverse
decay is given by Theorem~2.4 of \cite{demko1984decay}.

For $x>0$, the substitution $s=xu^2$ gives
\[
 \frac1\pi\int_0^\infty s^{-1/2}\frac{x}{x+s}\,\dd s
 =\frac{2\sqrt x}{\pi}\int_0^\infty\frac{\dd u}{1+u^2}
 =\sqrt x.
\]
The integral converges at zero and infinity. Applying the scalar identity to each eigenvalue gives the corresponding matrix identity. At each $s>0$, the ordered inverse identity
yields
\[
 \begin{split}
 P(P+sI)^{-1}-Q(Q+sI)^{-1}
 &=s[(Q+sI)^{-1}-(P+sI)^{-1}]\\
 &=s(P+sI)^{-1}(P-Q)(Q+sI)^{-1}.
 \end{split}
\]
Indeed the last product expands to $(Q+sI)^{-1}-(P+sI)^{-1}$ after writing $P-Q=(P+sI)-(Q+sI)$. Subtraction of the two convergent integrals gives
\begin{equation}
 P^{1/2}-Q^{1/2}
 =\frac1\pi\int_0^\infty s^{1/2}
             (P+sI)^{-1}(P-Q)(Q+sI)^{-1}\,\dd s.
 \label{B:eq:coupling-root-integral}
\end{equation}
The individual integrals defining the square roots converge in matrix norm, and the difference integrand is bounded by $s^{1/2}(c_0+s)^{-2}\|P-Q\|$. This also proves absolute
convergence of \eqref{B:eq:coupling-root-integral}.

For $u=1+p\wedge q$ and $d=|p-q|$,
\[
 \frac{1+d}{u}\sqrt{(1+p)(1+q)}
       =(1+d)\sqrt{1+d/u}\le(1+d)^{3/2}.
\]
Thus the last hypothesis implies $|(P-Q)_{pq}|\le C e^{-c_1|p-q|/2}/\sqrt{(1+p)(1+q)}$. Also
\[
 \sqrt{\frac{1+k}{1+p}}\le\sqrt{1+|k-p|},\qquad
 \sqrt{\frac{1+l}{1+q}}\le\sqrt{1+|l-q|}.
\]
For $a>0$ and $b\ge0$, differentiation gives the finite constant
\[
 M_{a,b}:=\sup_{r\ge0}(1+r)^b e^{-ar}
 \le\max\{1,(b/a)^b e^{a-b}\},
\]
where the right side is read as $1$ for $b=0$ when necessary. Use $a=\rho/2$, $b=1/2$ for each square-root factor, and use $a=c_1/2$, $b=3/2$ in the preceding estimate for $P-Q$.
Thus
\[
 \sqrt{1+|k-p|}e^{-\rho|k-p|}
 \le M_{\rho/2,1/2}e^{-\rho|k-p|/2},
\]
and the identical scalar inequality with $|l-q|$ treats the second factor. With $S=|k-p|+|p-q|+|q-l|$ and some fixed $\alpha>0$, the $(k,l)$ entry of the integrand in
\eqref{B:eq:coupling-root-integral}, without $s^{1/2}$, is bounded by
\[
 \frac{C}{(c_0+s)^2\sqrt{(1+k)(1+l)}}\sum_{p,q\in\mathcal I}e^{-\alpha S}.
\]
Since $S\ge|k-l|$, its sum is at most
\[
 e^{-\alpha|k-l|/2}
   \sum_p e^{-\alpha|k-p|/2}\sum_qe^{-\alpha|q-l|/2}
 \le C_\alpha e^{-\alpha|k-l|/2}.
\]
The last integral is evaluated by $s=c_0u^2$, followed by $u=\tan\theta$:
\[
 \begin{split}
 \int_0^\infty\frac{s^{1/2}}{(c_0+s)^2}\,\dd s
 &=\frac2{\sqrt{c_0}}\int_0^\infty\frac{u^2}{(1+u^2)^2}\,\dd u\\
 &=\frac2{\sqrt{c_0}}\int_0^{\pi/2}\sin^2\theta\,\dd\theta
 =\frac{\pi}{2\sqrt{c_0}}.
 \end{split}
\]
Combining this value with \eqref{B:eq:coupling-root-integral} and the bound for the sum over both indices gives $C_2e^{-c_2|k-l|}/\sqrt{(1+k)(1+l)}$. All sums were bounded by sums
over $\mathbb Z$, and the resolvent constants depend only on the stated hypotheses. The same constants therefore apply to every finite interval $\mathcal I$.
\end{proof}

\subsubsection{An exact Gaussian vector on the observation interval}
\label{B:subsection:exact-gaussian-vector-observation-interval}

The convolution has the local covariance computed above. We compare it with the conditional covariance on the observation interval. Rescaling by the local variances brings both
matrices into the uniform range where the square-root comparison applies.

Set $\mathscr D=\operatorname{diag}_{k\in I}(k^{-1/2})$, $K_I^G=(K_q(k,l))_{k,l\in I}$, and
\begin{equation}
 P_q=\mathscr D^{-1}K_I^G\mathscr D^{-1},\qquad
 Q=\mathscr D^{-1}K_n^{\rm fr}\mathscr D^{-1}.
 \label{B:eq:coupling-normalized-covariances}
\end{equation}
The matrix $Q$ depends on $n,V,\beta$, whereas $P_q$ also depends on the boundary configuration.

\begin{lemma}\label{B:lem:normalized-gaussian-covariances}
For all sufficiently large $n$, $P_q,Q$ satisfy \eqref{B:eq:matrix-square-root-hypotheses}, uniformly for $q\in\mathcal Q_n^\sharp$. In addition,
\begin{equation}
 \|P_q-Q\|\le C/m_n.
 \label{B:eq:coupling-normalized-error}
\end{equation}
\end{lemma}
\begin{proof}
Let $F=L_{q,I}$ be the restriction of $L_q$ to rows in $I$. For $v\in\mathbb R^I$, the disjoint supports of these rows give
\[
 |F^{\mathsf T}\mathscr D^{-1}v|^2
   =\sum_{k\in I}k(m_{x,k}^{-2}+m_{y,k}^{-2})v_k^2
   \asymp n|v|^2.
\]
The inverse Hessian obeys $a_+^{-1}I\le C\le a_-^{-1}I$. Consequently
\[
 \frac{|F^{\mathsf T}\mathscr D^{-1}v|^2}{n\beta a_+}
 \le v^{\mathsf T}P_qv
 =\frac1{n\beta}(F^{\mathsf T}\mathscr D^{-1}v)^{\mathsf T}
                       C(F^{\mathsf T}\mathscr D^{-1}v)
 \le\frac{|F^{\mathsf T}\mathscr D^{-1}v|^2}{n\beta a_-}.
\]
The preceding identity and $m_{\sigma,k}^2\asymp k/n$ give $c|v|^2\le v^{\mathsf T}P_qv\le C|v|^2$ with constants uniform for admissible $q$. Equations
\eqref{B:eq:covariance-midpoint} and \eqref{B:eq:coupling-frozen-error}, together with $\sqrt{kl}/(nt_{kl})\le1$, imply
\[
 |(P_q-Q)_{kl}|\le Ce^{-c|k-l|}
       \left[\frac{1+|k-l|}{1+k\wedge l}+\epsilon_n\right].
\]
For large $n$, $\epsilon_n\le n^{-1}$. Since $k,l\le n$,
\[
 \epsilon_n\le\frac1n\le\frac2{1+k\wedge l}
 \le2\frac{1+|k-l|}{1+k\wedge l}.
\]
The sum of the two terms in brackets is therefore at most three times its first term. For a fixed row $k$, split the sum at $l=k/2$. On the first part, $1+\min(k,l)\ge k/2$, so
\[
 \sum_{\substack{l\in I\\l\ge k/2}}|(P_q-Q)_{kl}|
 \le\frac Ck\sum_{d\in\mathbb Z}(1+|d|)e^{-c|d|}
 \le C/k.
\]
On the second part, $1+\min(k,l)\ge1+m_n$ and $k-l>k/2$. Hence
\[
 \sum_{\substack{l\in I\\l<k/2}}|(P_q-Q)_{kl}|
 \le\frac C{m_n}\sum_{d>k/2}(1+d)e^{-cd}
 \le\frac C{m_n}e^{-c'k}.
\]
Indeed, for an integer $D\ge0$ and $q=e^{-c}$,
\[
 \sum_{d=D}^\infty(1+d)q^d
 =q^D\left(\frac{D+1}{1-q}+\frac q{(1-q)^2}\right).
\]
Taking $D=\lfloor k/2\rfloor+1$ and using $\sup_{x\ge0}(1+x)e^{-cx/4}<\infty$ proves the second row bound with, for example, $c'=c/4$. The first row bound is at most $C/m_n$
because $k\ge m_n$. Thus all row sums are at most $C/m_n$. Since $(P_q-Q)_{kl}=(P_q-Q)_{lk}$, its $l$th absolute column sum is its $l$th absolute row sum. Hence
\[
 \|P_q-Q\|\le
 \left(\max_k\sum_l|(P_q-Q)_{kl}|
       \max_l\sum_k|(P_q-Q)_{kl}|\right)^{1/2}\le C/m_n.
\]
This proves \eqref{B:eq:coupling-normalized-error}.

Choose $q$ equal to the boundary restriction of $z_n^\sharp$. For large $n$ these coordinates are positive and this $q$ belongs to $\mathcal Q_n^\sharp$. The preceding bounds imply
$\lambda_{\min}(Q)\ge\lambda_{\min}(P_q)-C/m_n\ge c/2$, and also a uniform upper bound. This proves positivity of the fixed matrix $Q$. Finally,
\eqref{B:eq:covariance-finite-decay} gives $|(P_q)_{kl}|\le Ce^{-c|k-l|}$; the corresponding bound for $Q$ follows either from the difference estimate or directly from
\eqref{B:eq:coupling-frozen-covariance}.
\end{proof}

For the field in \eqref{B:eq:coupling-frozen-field}, define
\begin{equation}
 g_n=Q^{-1/2}\mathscr D^{-1}\eta_n^{\rm fr},\qquad
 \eta_{n,I}^G=\mathscr D P_q^{1/2}g_n,
 \qquad e_{n,k}^{(1)}=\eta_{n,k}^G-\eta_{n,k}^{\rm fr}.
 \label{B:eq:coupling-exact-comparison}
\end{equation}
These definitions use the previously constructed convolution.

\begin{proposition}\label{B:prop:exact-gaussian-covariance-coupling}
For every deterministic finite nonnegative $q$, $\eta_{n,I}^G$ is centered Gaussian with covariance $K_I^G$. For admissible $q$,
\begin{equation}
 |\operatorname{Cov}_q(e_{n,k}^{(1)},e_{n,l}^{(1)})|
                       \le C(kl)^{-3/2}e^{-c|k-l|}.
 \label{B:eq:gaussian-comparison-error-covariance}
\end{equation}
\end{proposition}
\begin{proof}
The matrix $Q$ is positive after the deterministic starting index in Lemma~\ref{B:lem:normalized-gaussian-covariances}. Since the convolution is jointly Gaussian,
\[
 \operatorname{Cov}(g_n)
 =Q^{-1/2}\mathscr D^{-1}K_n^{\rm fr}\mathscr D^{-1}Q^{-1/2}=I,
 \qquad
 \mathbb E e^{\ii v^{\mathsf T}g_n}=e^{-|v|^2/2}.
\]
The last expression factors over the coordinates, proving their independence. Thus
\[
 \operatorname{Cov}_q(\eta_{n,I}^G)
 =\mathscr D P_q^{1/2}I P_q^{1/2}\mathscr D
 =\mathscr D P_q\mathscr D=K_I^G.
\]
For every finite nonnegative $q$, the minimizer is positive on $J$ and $C^q$ is positive definite. The map $F$ is onto: for $v\in\mathbb R^I$ take $u_{y,k}=m_{y,k}v_k$, $u_{x,k}=0$
on $I$ and set the other coordinates to zero, obtaining $Fu=v$. It follows that $v^{\mathsf T}K_I^Gv=(F^{\mathsf T}v)^{\mathsf T}\Omega_q(F^{\mathsf T}v)>0$ for $v\ne0$. Hence
$P_q^{1/2}$ is defined on the full stated boundary space. Only the uniform error estimates require admissibility of $q$.

For admissible $q$, write $e_{n,k}^{(1)}=\sum_{p\in I}D_{kp}g_{n,p}$, where $D=\mathscr D(P_q^{1/2}-Q^{1/2})$. Lemma~\ref{B:lem:localized-matrix-square-root} gives
\[
 |D_{kp}|\le\frac{Ce^{-c|k-p|}}{(1+k)\sqrt{1+p}}
 \le Ck^{-3/2}\sqrt{1+|k-p|}\,e^{-c|k-p|}
 \le Ck^{-3/2}e^{-c'|k-p|}.
\]
Consequently
\[
 |\operatorname{Cov}_q(e_{n,k}^{(1)},e_{n,l}^{(1)})|
 \le C(kl)^{-3/2}\sum_p e^{-c'(|k-p|+|l-p|)}
 \le C(kl)^{-3/2}e^{-c''|k-l|}.
\]
This proves \eqref{B:eq:gaussian-comparison-error-covariance}. The joint law with the original convolution is retained by the same $g_n$. In particular,
\[
 \eta_n^{\rm fr}=\mathscr D Q^{1/2}g_n,\qquad
 \operatorname{Cov}_q(\eta_{n,I}^G,\eta_n^{\rm fr})
       =\mathscr D P_q^{1/2}Q^{1/2}\mathscr D.
\]
The factors in this covariance identity occur in the displayed order. The components of $g_n$ are functions of the original Brownian increments and auxiliary normals in
\eqref{B:eq:coupling-white-noise}. In particular, for this joint construction the complete error covariance is
\[
 \operatorname{Cov}_q(e_n^{(1)})
 =\mathscr D(P_q^{1/2}-Q^{1/2})^2\mathscr D.
\]
\end{proof}

\subsubsection{Decomposition using the Brownian increments}
\label{B:subsection:gaussian-brownian-increment-decomposition}

The factorization now separates the coupled Gaussian increment into its Brownian part, a discrete difference, and the remaining covariance and profile errors. The Gaussian sequence
in the discrete difference is defined at one additional site so that summation gives the required boundary term.

For $m_n\le k\le M_n+1$ set
\begin{equation}
 \Xi_{n,k}=n^{-1/2}\sum_{h\in\mathbb Z}b_{t_k}(h)\xi_{n,k-h},
 \qquad
 \eta_{n,k}^0=n^{-1/2}a_{t_k}^\circ\xi_{n,k}.
 \label{B:eq:coupling-primitive}
\end{equation}
For each of these series,
\[
 \mathbb E\sum_h n^{-1/2}|b_{t_k}(h)\xi_{n,k-h}|
 \le Ck^{-1/2}\sum_h e^{-c|h|}<\infty,\qquad
 \left\|n^{-1/2}\sum_{|h|>H}b_{t_k}(h)\xi_{n,k-h}\right\|_2^2
 \le Ck^{-1}e^{-cH}.
\]
Tonelli proves absolute almost sure convergence, and the Gaussian moment formula gives convergence in every finite $L^p$. The event of convergence may be intersected over all sizes
and indices. The index $M_n+1$ lies in $J$ for large $n$; its noise coordinate in \eqref{B:eq:coupling-white-noise} is auxiliary because $M_n+1\notin I$. For $k\in I$ define
\begin{equation}
 \begin{split}
 B_{n,k}&=\beta^{-1/2}[W(u_{n,k})-W(u_{n,k+1})],\\
 e_{n,k}^{(2)}&=n^{-1/2}\sum_h
        [b_{t_{k+1}}(h+1)-b_{t_k}(h+1)]\xi_{n,k-h},\\
 e_{n,k}^{(3)}&=
       \left(n^{-1/2}a_{t_k}^\circ-\sqrt{h^+_{n,k}/\beta}\right)\xi_{n,k}.
 \end{split}
 \label{B:eq:gaussian-increment-errors}
\end{equation}

\begin{proposition}\label{B:prop:gaussian-difference-decomposition}
The coupled vector satisfies, for $k\in I$,
\begin{equation}
 \eta_{n,k}^G=B_{n,k}+\Xi_{n,k}-\Xi_{n,k+1}
                         +\sum_{\ell=1}^3e_{n,k}^{(\ell)}.
 \label{B:eq:gaussian-increment-decomposition}
\end{equation}
For $\ell=1,2,3$,
\begin{equation}
 |\operatorname{Cov}_q(e_{n,k}^{(\ell)},e_{n,l}^{(\ell)})|
                       \le C(kl)^{-3/2}e^{-c|k-l|}.
 \label{B:eq:gaussian-error-covariance}
\end{equation}
Also, for $m_n\le j,l\le M_n+1$,
\begin{equation}
 |\operatorname{Cov}(\Xi_{n,j},\Xi_{n,l})|
                       \le C(jl)^{-1/2}e^{-c|j-l|}.
 \label{B:eq:coupling-primitive-covariance}
\end{equation}
\end{proposition}
\begin{proof}
The coefficient identity in \eqref{B:eq:coupling-coefficient-difference} and absolute convergence give
\[
 \eta_{n,k}^{\rm fr}
 =\eta_{n,k}^0+\Xi_{n,k}
             -n^{-1/2}\sum_hb_{t_k}(h+1)\xi_{n,k-h}.
\]
Reindexing the series for $\Xi_{n,k+1}$ yields
\[
 \Xi_{n,k+1}=n^{-1/2}\sum_hb_{t_{k+1}}(h+1)\xi_{n,k-h}.
\]
Their difference is exactly $e_{n,k}^{(2)}$, with the sign in \eqref{B:eq:gaussian-increment-errors}. For adjacent parameters,
\[
 n^{-1/2}|b_{t_{k+1}}(h+1)-b_{t_k}(h+1)|
 \le C n^{-3/2}t_k^{-3/2}e^{-c|h+1|}
 \le C k^{-3/2}e^{-c|h|}.
\]
In terms of the noise index $p$, its coefficient is
\[
 d^{(2)}_{kp}=n^{-1/2}
       [b_{t_{k+1}}(k-p+1)-b_{t_k}(k-p+1)],\qquad
 |d^{(2)}_{kp}|\le Ck^{-3/2}e^{-c|k-p|}.
\]
The covariance is $\sum_p d^{(2)}_{kp}d^{(2)}_{lp}$. Its absolute value is at most $C(kl)^{-3/2}\sum_p e^{-c(|k-p|+|l-p|)} \le C(kl)^{-3/2}e^{-c'|k-l|}$. The bound follows from the
convolution sum in Lemma~\ref{B:lem:gaussian-convolution-covariance} and proves \eqref{B:eq:gaussian-error-covariance} for $\ell=2$.

For $f(t)=\widehat\Gamma_t(0)$, Lemma~\ref{B:lem:zero-frequency-variance} gives
\begin{equation}
 h^+_{n,k}/\beta=\int_{k/n}^{(k+1)/n}f(t)\,\dd t.
 \label{B:eq:gaussian-brownian-variance-identity}
\end{equation}
The inequalities $c/t\le f(t)\le C/t$ imply
\[
 \frac c{k+1}\le\int_{k/n}^{(k+1)/n}f(t)\,\dd t\le\frac Ck,
 \qquad \frac ck\le n^{-1}f(t_k)\le\frac Ck.
\]
The exact difference is
\[
 n^{-1}f(t_k)-h^+_{n,k}/\beta
 =-\int_{k/n}^{(k+1)/n}\int_{k/n}^t f'(s)\,\dd s\,\dd t,
\]
whose absolute value is at most $C/k^2$. Both positive square roots are therefore bounded below by $c k^{-1/2}$. Hence
\[
 \left|n^{-1/2}a_{t_k}^\circ-\sqrt{h^+_{n,k}/\beta}\right|
 =\frac{|n^{-1}f(t_k)-h^+_{n,k}/\beta|}
          {\sqrt{n^{-1}f(t_k)}+\sqrt{h^+_{n,k}/\beta}}
 \le C k^{-3/2}.
\]
The variables $\xi_{n,k}$ for $k\in I$ are independent, so $\operatorname{Cov}(e_{n,k}^{(3)},e_{n,l}^{(3)})=0$ for $k\ne l$, and the variance is at most $Ck^{-3}$. This proves the
third error bound. The first is Proposition~\ref{B:prop:exact-gaussian-covariance-coupling}. Since $\eta_{n,k}^0=B_{n,k}+e_{n,k}^{(3)}$, the defining equation for $e^{(1)}$ proves
\eqref{B:eq:gaussian-increment-decomposition}.

Finally, independence of the noise coordinates gives
\[
 \operatorname{Cov}(\Xi_{n,j},\Xi_{n,l})
   =\frac1n\sum_h b_{t_j}(h)b_{t_l}(h+l-j).
\]
The coefficient bounds and the exponential convolution estimate in Lemma~\ref{B:lem:gaussian-convolution-covariance} give \eqref{B:eq:coupling-primitive-covariance}.
\end{proof}

All components in this decomposition are linear functions of the Gaussian variables already chosen. For example, their joint law has
\begin{equation}
 \operatorname{Cov}(\eta_{n,k}^0,\Xi_{n,l})
       =n^{-1}a_{t_k}^\circ b_{t_l}(l-k).
 \label{B:eq:coupling-cross-covariance}
\end{equation}
This follows by selecting the coefficient of $\xi_{n,k}$ in \eqref{B:eq:coupling-primitive}.

\subsubsection{The Gaussian vector on the full conditional interval}
\label{B:subsection:full-conditional-gaussian-vector}

The ratio variables on the observation interval do not determine the full coefficient vector. Gaussian regression fills the remaining directions while keeping those ratios
unchanged and restoring the exact full conditional covariance.

Let $F=L_{q,I}$ and $\Omega=\Omega_q$, retaining every coefficient in $\mathscr S(J)$ in $\Omega$. With $K_I^G=F\Omega F^{\mathsf T}$, define
\begin{equation}
 \begin{split}
 \Omega^\perp&=\Omega-\Omega F^{\mathsf T}(K_I^G)^{-1}F\Omega,\\
 \zeta_n^q&=\Omega F^{\mathsf T}(K_I^G)^{-1}\eta_{n,I}^G
                           +(\Omega^\perp)^{1/2}\omega_n,\\
 G_n^q&=m^q+s^q+\zeta_n^q,
 \end{split}
 \label{B:eq:gaussian-regression-formula}
\end{equation}
where $\omega_n$ is a standard normal vector on $\mathscr S(J)$ independent of all the previous Gaussian variables. These are the usual conditional Gaussian mean and covariance,
written on the full coefficient space.

\begin{proposition}\label{B:prop:full-gaussian-regression}
For every finite nonnegative $q$, \eqref{B:eq:gaussian-regression-formula} defines a vector with distribution $\gamma_{n,q}$, and
\begin{equation}
 \operatorname{Cov}_q(\zeta_n^q)=\Omega_q,\qquad
 L_{q,I}\zeta_n^q=\eta_{n,I}^G\quad\hbox{almost surely}.
 \label{B:eq:coupling-lift-verification}
\end{equation}
The construction is Borel measurable in $q$, $\eta_{n,I}^G$, and $\omega_n$.
\end{proposition}
\begin{proof}
Positive logarithmic coefficients on $J$ and Lemma~\ref{B:lem:reference-minimizer-coordinates} give $m^q>0$. The Hessian is positive definite, and $F$ is onto because its rows are
nonzero and have disjoint supports. Thus $K_I^G>0$ for every finite nonnegative $q$. We check directly that $\Omega^\perp$ is nonnegative. Set $T=\Omega^{1/2}F^{\mathsf T}$ and
\[
 P=T(T^{\mathsf T}T)^{-1}T^{\mathsf T}
   =\Omega^{1/2}F^{\mathsf T}(K_I^G)^{-1}F\Omega^{1/2}.
\]
The onto property of $F$ makes $T^{\mathsf T}T=K_I^G$ positive definite. Multiplication gives $P^{\mathsf T}=P$ and $P^2=P$, so $P$ is the orthogonal projection onto
$\operatorname{Ran}T$. Consequently
\[
 \Omega^\perp=\Omega^{1/2}(I-P)\Omega^{1/2}\ge0.
\]
The two summands in \eqref{B:eq:gaussian-regression-formula} are independent centered Gaussian vectors, because $\omega_n$ is independent of $\eta_{n,I}^G$. Their covariance sum is
\[
 \begin{aligned}
 &\Omega F^{\mathsf T}(K_I^G)^{-1}
          K_I^G(K_I^G)^{-1}F\Omega+\Omega^\perp\\
 &\qquad=\Omega F^{\mathsf T}(K_I^G)^{-1}F\Omega
       +\Omega-\Omega F^{\mathsf T}(K_I^G)^{-1}F\Omega=\Omega.
 \end{aligned}
\]
This proves $\zeta_n^q\sim N(0,\Omega)$, including when $\Omega^\perp$ is singular, and adding $m^q+s^q$ gives the required law of $G_n^q$. These formulas also agree with Gaussian
regression \cite{eaton2007multivariate}. For the exact image identity,
\[
 F\Omega^\perp F^{\mathsf T}
 =K_I^G-K_I^G(K_I^G)^{-1}K_I^G=0,
\]
which implies $F(\Omega^\perp)^{1/2}=0$. Applying $F$ to the regression formula therefore gives $F\zeta_n^q=\eta_{n,I}^G$, proving \eqref{B:eq:coupling-lift-verification}.
Conditionally on the Gaussian variables already constructed, its characteristic function at $v\in\mathbb R^{\mathscr S(J)}$ is
\[
 \exp\left\{\ii v^{\mathsf T}\Omega F^{\mathsf T}(K_I^G)^{-1}
                   \eta_{n,I}^G-\tfrac12v^{\mathsf T}\Omega^\perp v\right\}.
\]
Integrating this conditional characteristic function against the centered law of $\eta_{n,I}^G$ gives $\mathbb E_qe^{\ii v^{\mathsf T}\zeta_n^q} =\exp(-v^{\mathsf T}\Omega v/2)$.
The complete joint covariance is
\[
 \operatorname{Cov}_q
 \begin{pmatrix}\zeta_n^q\\ \eta_{n,I}^G\end{pmatrix}
 =
 \begin{pmatrix}
  \Omega&\Omega F^{\mathsf T}\\ F\Omega&K_I^G
 \end{pmatrix},
\]
because $\operatorname{Cov}_q(\zeta_n^q,\eta_{n,I}^G) =\Omega F^{\mathsf T}(K_I^G)^{-1}K_I^G$. The conditional characteristic function above also specifies the joint law with every
previously chosen Brownian increment and convolution coordinate. More explicitly, let $H$ be any finite centered vector of those variables and put
$C_H=\operatorname{Cov}_q(\eta_{n,I}^G,H)$. The independent regression normal $\omega_n$ is independent of the entire vector $(\eta_{n,I}^G,H)$. Consequently
\[
 \operatorname{Cov}_q(\zeta_n^q,H)
       =\Omega F^{\mathsf T}(K_I^G)^{-1}C_H.
\]
Together with $\operatorname{Cov}_q(\zeta_n^q)=\Omega$ and the unchanged covariance of $H$, this determines the joint Gaussian distribution. At fixed boundary data,
$(\Omega^\perp)^{1/2}\omega_n$ is independent of the previously constructed Gaussian variables.

Continuity of $m^q$ in the boundary data is Proposition~\ref{B:prop:conditional-minimizer-dependence}. The Hessian and the third derivatives at this positive minimizer are
continuous, as are $C^q,s^q$. Inversion and the positive square root are continuous on the positive definite cone. For positive semidefinite matrices of norm at most $M_0$, choose
polynomials $p_j$ converging uniformly to $\sqrt t$ on $[0,M_0]$. The spectral theorem gives $\|p_j(S)-S^{1/2}\|\le\sup_{[0,M_0]}|p_j(t)-\sqrt t|$. The square root is therefore a
locally uniform limit of continuous matrix maps, which proves its continuity also at singular matrices. This applies to $\Omega^\perp$. All finite matrix maps above are therefore
Borel. The infinite series are measurable limits of finite sums; on the measurable null set of divergence give them value zero. This proves measurability in the parameters on the
entire stated domain.
\end{proof}

For later changes of coordinates, put $\Sigma=M^{-1}\Omega M^{-1}$ and $F_{\rm rel}=FM$. The conditional covariance in these coordinates is
$\Sigma^\perp=\Sigma-\Sigma F_{\rm rel}^{\mathsf T}(K_I^G)^{-1} F_{\rm rel}\Sigma$. Direct cancellation of $M^{-1}M$ gives
\begin{equation}
 M\Sigma F_{\rm rel}^{\mathsf T}(K_I^G)^{-1}
     =\Omega F^{\mathsf T}(K_I^G)^{-1},\qquad
 M\Sigma^\perp M=\Omega^\perp.
 \label{B:eq:coupling-lift-coordinate-identity}
\end{equation}
Conditionally on $\eta_{n,I}^G$, multiplication of a relative Gaussian vector by $M$ gives the characteristic function
\[
 \exp\left\{\ii v^{\mathsf T}M\Sigma F_{\rm rel}^{\mathsf T}
       (K_I^G)^{-1}\eta_{n,I}^G
           -\tfrac12v^{\mathsf T}M\Sigma^\perp Mv\right\}.
\]
Substitution of \eqref{B:eq:coupling-lift-coordinate-identity} gives exactly the characteristic function in the proof of Proposition~\ref{B:prop:full-gaussian-regression}. This
comparison uses the covariance $M\Sigma^\perp M$ of the independent Gaussian term after multiplication by $M$; the square root itself is taken in the specified coefficient
coordinates.

\subsubsection{Partial sums of the coupling error}
\label{B:subsection:partial-sums-coupling-error}

The inverse uses cumulative logarithmic ratios. Summing the increment decomposition therefore gives the relevant error directly: the first differences sum to the difference of the
two endpoint values, and the other errors remain partial sums.

For $j\in I$ define
\begin{equation}
 \begin{split}
 C_{n,j}^{(\ell)}&=\sum_{k=j}^{M_n}e_{n,k}^{(\ell)},\quad\ell=1,2,3,
             \qquad \chi_n^\nabla=\Xi_{n,M_n+1},\\
 U_n^G&=M^{-1}\zeta_n^q,\qquad
 \eta_{n,k}^G=(U_n^G)_{y,k}-(U_n^G)_{x,k}\quad(k\in J),\\
 V_{n,j}^G&=\chi_n^\nabla+
                    \sum_{k=j}^{M_n}(\eta_{n,k}^G-B_{n,k})
          =\Xi_{n,j}+\sum_{\ell=1}^3 C_{n,j}^{(\ell)}.
 \end{split}
 \label{B:eq:coupling-decomposition}
\end{equation}
The notation for $\eta^G$ on $I$ is consistent by \eqref{B:eq:coupling-lift-verification}. The last equality follows by summing \eqref{B:eq:gaussian-increment-decomposition}; the
first differences sum to $\Xi_{n,j}-\Xi_{n,M_n+1}$.

\begin{lemma}\label{B:lem:gaussian-partial-sum-bounds}
Uniformly for admissible $q$ and $j\in I$,
\begin{equation}
 \operatorname{Var}_q C_{n,j}^{(\ell)}\le C/j^2,\qquad
 \operatorname{Var}_q V_{n,j}^G\le C/j,\qquad
 \operatorname{Var}\chi_n^\nabla\le C/n.
 \label{B:eq:coupling-variances}
\end{equation}
For $k,l\in J$ and both species,
\begin{equation}
 |\operatorname{Cov}_q((U_n^G)_{\sigma,k},(U_n^G)_{\tau,l})|
                         \le C(kl)^{-1/2}e^{-c|k-l|}.
 \label{B:eq:coupling-local-covariance}
\end{equation}
\end{lemma}
\begin{proof}
For any $k,l\ge1$, $k/l\le1+|k-l|$. Applying \eqref{B:eq:gaussian-error-covariance} and summing first over $l$,
\[
 \begin{split}
 \operatorname{Var}_q C_{n,j}^{(\ell)}
 &\le C\sum_{k,l=j}^{M_n}(kl)^{-3/2}e^{-c|k-l|}\\
 &\le C\sum_{k=j}^{M_n}k^{-3}
           \sum_{l=j}^{M_n}(1+|k-l|)^{3/2}e^{-c|k-l|}
 \le C\sum_{k=j}^{\infty}k^{-3}\le C/j^2.
 \end{split}
\]
The last sum is at most $j^{-3}+\int_j^\infty x^{-3}\,\dd x$. All four variables in the last line of \eqref{B:eq:coupling-decomposition} are centered and may be correlated.
Pointwise Cauchy--Schwarz gives
\[
 \operatorname{Var}_q V_{n,j}^G
 =\mathbb E_q\left|\Xi_{n,j}+\sum_{\ell=1}^3C_{n,j}^{(\ell)}\right|^2
 \le4\operatorname{Var}\Xi_{n,j}
         +4\sum_{\ell=1}^3\operatorname{Var}_qC_{n,j}^{(\ell)}
 \le C/j+C/j^2\le C'/j.
\]
Here \eqref{B:eq:coupling-primitive-covariance} bounds the first variance. As $M_n=R-b=n-r-b$ and $(r+b)/n\to0$, choose $n_0$ such that $M_n+1\ge n/2$ for $n\ge n_0$. Thus
\[
 \operatorname{Var}\chi_n^\nabla
 =\operatorname{Var}\Xi_{n,M_n+1}
 \le C/(M_n+1)\le 2C/n.
\]
Finally, the full covariance identity gives
\[
 \operatorname{Cov}_q(U_n^G)
 =M^{-1}\Omega_qM^{-1}
 =\frac1{n\beta}\mathsf C,\qquad
 \frac1{n\beta\sqrt{t_kt_l}}=\frac1{\beta\sqrt{kl}}.
\]
Substitution of \eqref{B:eq:covariance-finite-decay} proves \eqref{B:eq:coupling-local-covariance}, including both species and every coordinate in $J$.
\end{proof}

\subsubsection{Weighted Gaussian quadratic forms}
\label{B:subsection:weighted-gaussian-quadratic-forms}

The weighted square sum depends on covariances between distinct partial sums and between the components of the Gaussian approximation error. We estimate these covariances to obtain
the operator and trace bounds used in concentration of the quadratic form.

Fix $1<\lambda<2$ and $C_0<\infty$. Let deterministic weights satisfy
\begin{equation}
 0\le w_{n,j}\le C_0n^{-\lambda}j^{\lambda-1},\qquad
 0\le w_{n,\alpha}\le C_0n^{-\lambda}k_\alpha^{\lambda-1},\qquad
 \tau_{n,\lambda}=n^{-\lambda}m_n^{\lambda-2}.
 \label{B:eq:coupling-weights}
\end{equation}
The first array is indexed by $I$, the second by $\mathscr S(J)$. For example $h^-_{n,j}z_{n,j}^{\lambda}$ has the first bound, by the estimates for the partition intervals in
Section~\ref{B:sec:approximation-by-brownian-motion-and-the-lower-bound}. Note that $n\tau_{n,\lambda}\to0$. Let $\Sigma^V,\Sigma^\Xi,\Sigma^U$ denote the covariance matrices of
$(\sqrt{w_{n,j}}V_{n,j}^G)_j$, $(\sqrt{w_{n,j}}\Xi_{n,j})_j$, and $(\sqrt{w_{n,\alpha}}(U_n^G)_\alpha)_\alpha$, respectively. Write $\Sigma^{(\ell)}$ for the covariance of
$(\sqrt{w_{n,j}}C_{n,j}^{(\ell)})_j$.

\begin{proposition}\label{B:prop:weighted-gaussian-covariances}
For each $\Sigma\in\{\Sigma^V,\Sigma^\Xi,\Sigma^U\}$,
\begin{equation}
 \operatorname{Tr}\Sigma\le T/n,\qquad
 \|\Sigma\|\le C\tau_{n,\lambda}.
 \label{B:eq:coupling-covariance-energy}
\end{equation}
For $\ell=1,2,3$,
\begin{equation}
 \operatorname{Tr}\Sigma^{(\ell)}\le C\tau_{n,\lambda},\qquad
 \|\Sigma^{(\ell)}\|\le C\tau_{n,\lambda}.
 \label{B:eq:gaussian-small-error-energy}
\end{equation}
The constants are uniform in admissible $q$ and in the weights satisfying \eqref{B:eq:coupling-weights}.
\end{proposition}
\begin{proof}
Equation \eqref{B:eq:coupling-primitive-covariance} yields
\[
 |\Sigma^\Xi_{jl}|\le C\sqrt{w_{n,j}/j}\sqrt{w_{n,l}/l}
                                                    e^{-c|j-l|}.
\]
Put $a_j=\sqrt{w_{n,j}/j}$. On coordinates of positive weight the matrix $E_{jl}=\Sigma^\Xi_{jl}/(a_ja_l)$ has absolute row and column sums bounded by
$C\sum_{d\in\mathbb Z}e^{-c|d|}$; set its remaining entries to zero. Thus
\[
 \|\Sigma^\Xi\|=\|\operatorname{diag}(a)E\operatorname{diag}(a)\|
       \le C\max_j a_j^2\le C n^{-\lambda}m_n^{\lambda-2}.
\]
The trace is at most
\[
 C\sum_{j\in I}w_{n,j}/j
 \le Cn^{-\lambda}\sum_{j=1}^n j^{\lambda-2}
 \le Cn^{-\lambda}\left(1+\int_1^nx^{\lambda-2}\,\dd x\right)
 \le C_\lambda/n.
\]
For the full coefficient vector put $a_\alpha=(w_{n,\alpha}/k_\alpha)^{1/2}$ on coordinates with positive weight and zero elsewhere. Equation \eqref{B:eq:coupling-local-covariance}
gives
\[
 |\Sigma^U_{\alpha\gamma}|
 \le Ca_\alpha a_\gamma e^{-c|k_\alpha-k_\gamma|}.
\]
After division by $a_\alpha a_\gamma$ on their nonzero rows, the absolute row and column sums are at most $2C\sum_{d\in\mathbb Z}e^{-c|d|}$, because there are two species. Zero
weights give identically zero rows and columns. Hence
\[
 \|\Sigma^U\|\le C\max_\alpha\frac{w_{n,\alpha}}{k_\alpha}
 \le Cn^{-\lambda}r^{\lambda-2},\qquad
 \operatorname{Tr}\Sigma^U
 \le2Cn^{-\lambda}\sum_{k=r}^R k^{\lambda-2}\le C_\lambda/n.
\]
The ratio $r^{\lambda-2}/m_n^{\lambda-2}=(m_n/r)^{2-\lambda}$ is bounded, proving the required operator bound by
$C\tau_{n,\lambda}$.

For $\Sigma^{(\ell)}$, Lemma~\ref{B:lem:gaussian-partial-sum-bounds} gives
\[
 \operatorname{Tr}\Sigma^{(\ell)}
 \le C\sum_{j\in I}w_{n,j}/j^2
 \le Cn^{-\lambda}\sum_{j=m_n}^\infty j^{\lambda-3}
 \le C_\lambda n^{-\lambda}m_n^{\lambda-2}.
\]
Its operator norm is at most its trace because it is positive semidefinite. For any deterministic vector $v$, set $X_0=\sum_jv_j\sqrt{w_{n,j}}\Xi_{n,j}$ and
$X_\ell=\sum_jv_j\sqrt{w_{n,j}}C_{n,j}^{(\ell)}$. These variables are centered, so
\[
 v^{\mathsf T}\Sigma^Vv
 =\mathbb E|X_0+X_1+X_2+X_3|^2
 \le4\mathbb E\sum_{\ell=0}^3|X_\ell|^2.
\]
Equivalently,
\begin{equation}
 \Sigma^V\le4(\Sigma^\Xi+\Sigma^{(1)}+\Sigma^{(2)}+\Sigma^{(3)}).
 \label{B:eq:coupling-cross-control}
\end{equation}
For positive semidefinite matrices, an inequality of quadratic forms passes to the trace and to the largest eigenvalue. Thus
\[
 \operatorname{Tr}\Sigma^V
 \le4\operatorname{Tr}\Sigma^\Xi
       +4\sum_{\ell=1}^3\operatorname{Tr}\Sigma^{(\ell)}
 \le C/n+C\tau_{n,\lambda},
 \qquad
 \|\Sigma^V\|\le C\tau_{n,\lambda}.
\]
Since $n\tau_{n,\lambda}=n^{1-\lambda}m_n^{\lambda-2}\le1$ and tends to zero, this proves \eqref{B:eq:coupling-covariance-energy}. The inequality between quadratic forms
\eqref{B:eq:coupling-cross-control} includes all cross covariances among the four terms of $V^G$.
\end{proof}

We turn the covariance bounds for these quadratic forms into probability estimates with the following Gaussian inequality. For a finite centered real Gaussian vector $X$ with
covariance $\Sigma$, write $S=\operatorname{Tr}\Sigma$ and choose $L>0$ with $\|\Sigma\|\le L$.

\begin{lemma}\label{B:lem:gaussian-quadratic-form-tail}
For every $v>2S$,
\begin{equation}
 \Pr\{|X|^2>v\}\le\exp[-(v-2S)/(4L)].
 \label{B:eq:coupling-mgf-tail}
\end{equation}
\end{lemma}
\begin{proof}
By an orthogonal diagonalization, $|X|^2$ has the same distribution as $\sum_i\lambda_i g_i^2$, where $g_i$ are independent standard normals and $0\le\lambda_i\le L$,
$\sum_i\lambda_i=S$. A zero eigenvalue contributes zero and needs no inverse covariance. For $0\le\theta<(2L)^{-1}$, the one-dimensional Gaussian integral gives
\[
 \mathbb E e^{\theta\lambda_i g_i^2}
 =\frac1{\sqrt{2\pi}}\int_{\mathbb R}
       e^{-(1-2\theta\lambda_i)x^2/2}\,\mathrm dx
 =(1-2\theta\lambda_i)^{-1/2}.
\]
Take $\theta=(4L)^{-1}$. Since $0\le\lambda_i/(2L)\le1/2$ and
\[
 -\log(1-u)=\int_0^u\frac{\mathrm dt}{1-t}\le2u
       \qquad(0\le u\le1/2),
\]
independence yields
\[
 \log\mathbb E e^{|X|^2/(4L)}
 =-\frac12\sum_i\log\left(1-\frac{\lambda_i}{2L}\right)
 \le\frac1{2L}\sum_i\lambda_i=\frac S{2L}.
\]
Exponential Markov now gives
\[
 \Pr\{|X|^2>v\}
 \le e^{-v/(4L)}\mathbb E e^{|X|^2/(4L)}
 \le\exp\!\left[-\frac{v-2S}{4L}\right],
\]
which is \eqref{B:eq:coupling-mgf-tail}. If $\Sigma=0$, all coordinates vanish almost surely and the conclusion holds directly. The constants use only $S$ and the stated upper
bound $L$, not the dimension.
\end{proof}

For a fixed $K>2T$, Proposition~\ref{B:prop:weighted-gaussian-covariances} therefore gives
\begin{equation}
 \Pr_q\!\left\{\sum_{j\in I}w_{n,j}|V_{n,j}^G|^2>K/n\right\}
 \le \exp[-c n^{\lambda-1}m_n^{2-\lambda}].
 \label{B:eq:coupling-energy-tail}
\end{equation}
For each of the three covariance matrices, $S\le T/n$ and $L=C\tau_{n,\lambda}$. At $v=K/n$ with $K>2T$ the exponent is at least
\[
 \frac{K-2T}{4Cn\tau_{n,\lambda}}
 =\frac{K-2T}{4C}n^{\lambda-1}m_n^{2-\lambda}.
\]
Choosing common values of $T,C,K$ therefore gives \eqref{B:eq:coupling-energy-tail} also for the sums of $w_{n,j}|\Xi_{n,j}|^2$ and $w_{n,\alpha}|(U_n^G)_\alpha|^2$. For each fixed
$\varepsilon>0$, once
$2C\tau_{n,\lambda}\le\varepsilon/(2n)$,
\begin{equation}
 \Pr_q\!\left\{\sum_jw_{n,j}|C_{n,j}^{(\ell)}|^2>
                                         \varepsilon/n\right\}
 \le \exp[-c\varepsilon n^{\lambda-1}m_n^{2-\lambda}].
 \label{B:eq:coupling-small-tail}
\end{equation}
Indeed, $v-2S\ge\varepsilon/(2n)$ at the stated starting index. Taking $L=C\tau_{n,\lambda}$ in \eqref{B:eq:coupling-mgf-tail} gives the exponent
$\varepsilon/(8Cn\tau_{n,\lambda}) =(\varepsilon/(8C))n^{\lambda-1}m_n^{2-\lambda}$. The starting index may depend on $\varepsilon$; the Gaussian construction remains fixed.

\subsubsection{Maximum estimates}
\label{B:subsection:gaussian-maximum-estimates}

The weighted bounds determine the rate, while the nonlinear logarithmic comparison also needs uniform smallness. The same variance estimates give maximum bounds after a union over
the finite row.

\begin{lemma}\label{B:lem:coupled-gaussian-maxima}
For every fixed $D>0$, there is $K_D<\infty$ such that, uniformly for admissible $q$, the following inequalities fail with probability at most $C_Dn^{-D}$:
\begin{equation}
 \begin{split}
 &\max_{j\in I}\sqrt j(|\Xi_{n,j}|+|V_{n,j}^G|)
   +\max_{\alpha\in\mathscr S(J)}\sqrt{k_\alpha}|(U_n^G)_\alpha|
   +\sqrt n|\chi_n^\nabla|
                      \le K_D\sqrt{\log(\me+n)},\\
 &\max_{j\in I}j\sum_{\ell=1}^3|C_{n,j}^{(\ell)}|
                      \le K_D\sqrt{\log(\me+n)}.
 \end{split}
 \label{B:eq:coupling-maxima}
\end{equation}
\end{lemma}
\begin{proof}
Every individual rescaled coordinate on the left is centered Gaussian with variance at most a deterministic $\sigma_0^2$, by \eqref{B:eq:coupling-primitive-covariance},
\eqref{B:eq:coupling-variances} and \eqref{B:eq:coupling-local-covariance}. For a centered Gaussian $Y$ with variance at most $\sigma_0^2$,
$\mathbb E e^{tY}\le e^{t^2\sigma_0^2/2}$. Optimizing exponential Markov's inequality at $t=x/\sigma_0^2$ gives $\Pr(|Y|>x)\le2e^{-x^2/(2\sigma_0^2)}$. The variables to which this
scalar bound is applied are $\sqrt j\,\Xi_{n,j}$, $\sqrt j\,V_{n,j}^G$, $\sqrt{k_\alpha}(U_n^G)_\alpha$, $\sqrt n\,\chi_n^\nabla$, and $jC_{n,j}^{(\ell)}$. Their number is
$2|I|+2|J|+1+3|I|\le8n+1$. Choose $\sigma_0>0$ above all their standard deviations and choose $L_D$ so that $L_D^2/(2\sigma_0^2)>D+3$. Their union event at threshold
$L_D\sqrt{\log(\me+n)}$ has probability at most
\[
 2(8n+1)(\me+n)^{-L_D^2/(2\sigma_0^2)}\le C_Dn^{-D}.
\]
On the complement, the first maximum is at most $2L_D\sqrt{\log(\me+n)}$, the maximum over the full coefficient vector is at most $L_D\sqrt{\log(\me+n)}$, and the scalar term has
the same bound. The second line is at most $3L_D\sqrt{\log(\me+n)}$. Thus $K_D=4L_D$ gives both bounds. The union bound uses only the individual variance estimates and applies to
their given joint Gaussian distribution.
\end{proof}

\subsubsection{Joint construction for all matrix sizes}
\label{B:subsection:gaussian-construction-all-sizes}

The conditional Gaussian formulas can now be evaluated at the random exterior configuration. The construction keeps the exterior sequence and auxiliary randomness independent of
the Brownian path and records exactly which observations are used at each size.

Let $(E_n)_{n\ge n_0}$ have the joint original exterior law constructed in Theorem~\ref{B:thm:original-exterior-realization}, and put $q_n=c_n(E_n)$ and
$\mathscr E=\sigma(E_n:n\ge n_0)$. Take a Brownian motion $W$ independent of $\mathscr E$ and all the auxiliary standard normal variables in \eqref{B:eq:coupling-white-noise} and
\eqref{B:eq:gaussian-regression-formula} independent of both. Let $\mathcal U$ be the sigma-field generated by $\mathscr E$ and these auxiliary variables. Recall $\mathcal N_n$
from \eqref{B:eq:brownian-sampled-information}.

\begin{theorem}\label{B:thm:joint-gaussian-brownian-construction}
The preceding formulas, with $q=q_n$, give
\begin{equation}
 \mathcal L(G_n^{q_n}\mid\mathscr E)=\gamma_{n,q_n},\qquad
 G_n^{q_n},\Xi_n,\chi_n^\nabla,\eta_n^G
                  \text{ are }\mathcal N_n\vee\mathcal U\text{-measurable}.
 \label{B:eq:coupling-information}
\end{equation}
Here $\mathcal U$ is independent of the entire Brownian path. Conditionally on $\mathscr E$, all fields at a fixed size are jointly Gaussian and satisfy
\eqref{B:eq:gaussian-increment-decomposition} and \eqref{B:eq:coupling-decomposition}.
\end{theorem}
\begin{proof}
Let $\mathbf P_{\rm ext}$ be the law of the entire exterior sequence from Theorem~\ref{B:thm:original-exterior-realization}. Take its product with the countable product of standard
normal laws needed for all $\xi_{n,l}^\perp$ and all coordinates of $\omega_n$, and with Wiener measure. Denote the first two factors together by $\mathbf P_{\rm aux}$. For bounded
Borel functions $f,g$ of the indicated variables, Fubini gives
\[
 \mathbb E[f(W)g(E,\xi^\perp,\omega)]
  =\left(\int f\,\dd\mathbf P_W\right)
     \left(\int g\,\dd\mathbf P_{\rm aux}\right).
\]
This proves independence of $\mathcal U$ from the entire Brownian path while preserving the dependence between exterior vectors at different sizes. For fixed $n$, finite sets
$F\subset\mathbb Z$, and real coefficients $s_l$ and vector $v$, the conditional characteristic function is
\[
 \mathbb E\!\left[
 e^{\ii\sum_{l\in F}s_l\xi_{n,l}+\ii v^{\mathsf T}\omega_n}
 \,\middle|\,\mathscr E\right]
 =\exp\!\left(-\tfrac12\sum_{l\in F}s_l^2-\tfrac12|v|^2\right).
\]
Here the sampling intervals for the Brownian increments and their lengths are deterministic. Consequently the calculation in
Proposition~\ref{B:prop:exact-gaussian-covariance-coupling} applies conditionally with $q=q_n$. The regression in Proposition~\ref{B:prop:full-gaussian-regression} then gives, for
any $v\in\mathbb R^{\mathscr S(J)}$,
\[
 \mathbb E[e^{\ii v^{\mathsf T}G_n^{q_n}}\mid\mathscr E]
 =\exp\{\ii v^{\mathsf T}(m^{q_n}+s^{q_n})
                         -\tfrac12v^{\mathsf T}\Omega_{q_n}v\}.
\]
This identifies the conditional distribution for every finite nonnegative value of $q_n$. For any finite linear combination of the centered fields, truncation of the series gives a
linear combination of the same independent normals, with coefficients determined by $q_n$. The $L^2$ convergence proved above gives convergence of the conditional variances and
characteristic functions. Thus the complete tuple at fixed $n$ is conditionally jointly Gaussian. Adding the specified mean to $G_n^{q_n}$ preserves the conditional Gaussian
distribution. The two decomposition identities are algebraic on the common event of absolute convergence.

The parameter $q_n$ is $\mathcal U$-measurable. In \eqref{B:eq:coupling-white-noise} the Brownian variables are evaluated only at the nodes of row $n$; all remaining variables are
$\mathcal U$-measurable. Finite sums and the Borel matrix maps thus have the required measurability. The absolute convergence events have probability one for each row; their
countable intersection permits using the measurable limits simultaneously for all rows. This proves \eqref{B:eq:coupling-information}.

For consecutive sampling times $u<v$, set $B_{u,v}(t)=W(t)-W(u)-(t-u)(W(v)-W(u))/(v-u)$. For every sampling time $s$ at size $n$, Brownian covariance gives
\[
 \operatorname{Cov}(B_{u,v}(t),W(s))
 =\min(t,s)-\min(u,s)
       -\frac{t-u}{v-u}[\min(v,s)-\min(u,s)]=0.
\]
For $s\le u$ all differences are zero; for $s\ge v$ the two terms are $t-u$ and $t-u$. These cases include $s=u,v$. Gaussian characteristic functions give independence from the
vector of sampled Brownian values. Dense evaluations generate the path sigma-field of each continuous bridge, as in Lemma~\ref{B:lem:conditional-bridge-law}. Finally, let $f,g,h$
be bounded functions of the bridges, the sampled Brownian values, and the auxiliary variables, respectively. Then
\[
 \mathbb E[f(B)g(W_{\rm nodes})h(E,\xi^\perp,\omega)]
 =\mathbb E[f(B)g(W_{\rm nodes})]\,\mathbb Eh
 =\mathbb Ef(B)\,\mathbb E[g(W_{\rm nodes})h(E,\xi^\perp,\omega)].
\]
A monotone class argument proves independence of the bridge tuple from $\mathcal N_n\vee\mathcal U$. These distributional identities concern the chosen sampling times at size $n$;
the joint dependence between sizes is preserved.
\end{proof}

\par\begingroup
\dimen0=\pagegoal\advance\dimen0 by -\pagetotal \ifdim\dimen0<6\baselineskip\newpage\fi
\endgroup
\begin{corollary}\label{B:cor:gaussian-coupling-almost-sure-bounds}
For any fixed deterministic weight arrays satisfying \eqref{B:eq:coupling-weights}, the following limits hold almost surely. The inequality holds for every $n\ge n_0(\omega)$ for a
finite random index $n_0(\omega)$:
\begin{equation}
 \begin{split}
 n\sum_{j\in I}w_{n,j}|C_{n,j}^{(\ell)}|^2&\longrightarrow0
                                         \quad(\ell=1,2,3),\\
 \sum_{j\in I}w_{n,j}(|V_{n,j}^G|^2+|\Xi_{n,j}|^2)
   +\sum_{\alpha\in\mathscr S(J)}w_{n,\alpha}|(U_n^G)_\alpha|^2
                                &\le K/n\quad(n\ge n_0(\omega)),\\
 \max_{j\in I}|V_{n,j}^G|
  +\max_{\alpha\in\mathscr S(J)}|(U_n^G)_\alpha|
  +|\chi_n^\nabla|&\longrightarrow0.
 \end{split}
 \label{B:eq:coupling-alln}
\end{equation}
The construction is independent of the weights and of $\lambda$. The conclusions may be intersected for a fixed countable collection of weight arrays and exponents.
\end{corollary}
\begin{proof}
Let $\mathcal G_n=\{q_n\in\mathcal Q_n^\sharp\}$. Theorem~\ref{B:thm:boundary-and-terminal-probabilities} gives $\sum_n\Pr(\mathcal G_n^c)<\infty$ for the original exterior law.
For any $\varepsilon>0$, conditional integration of \eqref{B:eq:coupling-small-tail} gives
\[
 \begin{split}
 \Pr\!\left\{n\sum_jw_{n,j}|C_{n,j}^{(\ell)}|^2>\varepsilon\right\}
 &\le\Pr(\mathcal G_n^c)
   +\mathbb E\left[\mathbf1_{\mathcal G_n}
       \Pr\!\left(n\sum_jw_{n,j}|C_{n,j}^{(\ell)}|^2>
                                       \varepsilon\mid\mathscr E\right)\right]\\
 &\le\Pr(\mathcal G_n^c)
           +e^{-c\varepsilon n^{\lambda-1}m_n^{2-\lambda}}.
 \end{split}
\]
For fixed $\lambda>1$, the exponent divided by $\log n$ tends to infinity: $m_n\ge1$ and $2-\lambda>0$, while $n^{\lambda-1}/\log n\to\infty$. Thus the second term is summable.
Borel--Cantelli at $\varepsilon=1/h$, $h\in\mathbb N$, proves the first limit for all three errors.

Write $Q_n^V,Q_n^\Xi,Q_n^U$ for the three weighted squared sums in the second line of \eqref{B:eq:coupling-alln}, and fix thresholds $K_V,K_\Xi,K_U$ as in
\eqref{B:eq:coupling-energy-tail}. Conditional integration against the original exterior law gives
\[
 \Pr\!\left(\bigcup_{F\in\{V,\Xi,U\}}\{Q_n^F>K_F/n\}\right)
 \le\Pr(\mathcal G_n^c)+3e^{-c n^{\lambda-1}m_n^{2-\lambda}}.
\]
The series is finite. Borel--Cantelli gives a finite $n_0(\omega)$ such that $Q_n^V+Q_n^\Xi+Q_n^U\le(K_V+K_\Xi+K_U)/n$ for every $n\ge n_0(\omega)$. This proves the second line
with $K=K_V+K_\Xi+K_U$.

Apply Lemma~\ref{B:lem:coupled-gaussian-maxima} with $D=3$ under the same conditional law. If $\mathscr B_n$ is its failure event, then
\[
 \Pr(\mathscr B_n)\le
 \Pr(\mathcal G_n^c)+
 \mathbb E[\mathbf1_{\mathcal G_n}\Pr(\mathscr B_n\mid\mathscr E)]
 \le\Pr(\mathcal G_n^c)+C_3n^{-3}.
\]
These probabilities are summable. On the resulting event, the last line of \eqref{B:eq:coupling-alln} is bounded, for every $n\ge n_1(\omega)$, by
\[
 C\sqrt{\log(\me+n)/m_n}
   +C\sqrt{\log(\me+n)/r}
   +C\sqrt{\log(\me+n)/n},
\]
which tends to zero since $r\asymp m_n\asymp(\log n)^{A_0}$. For the first limit, take the intersection over $\varepsilon=1/h$, $h\in\mathbb N$, and $\ell=1,2,3$. Intersect this
event with those for the three fixed thresholds and the maximum bound. Summability and the first Borel--Cantelli lemma give probability one for each of these events. For a fixed
countable family of weights and exponents, a further countable intersection gives the estimates for all those weights and exponents. The Gaussian maps were fixed before
\eqref{B:eq:coupling-weights} and depend only on the coefficient data, Brownian samples, and auxiliary normals.
\end{proof}

\subsubsection{The scalar normalization of the inverse kernel}
\label{B:subsection:scalar-normalization-gaussian-sums}

The correction at the upper index $\chi_n^\nabla$ in \eqref{B:eq:coupling-decomposition} can be absorbed into the scalar normalization of the inverse kernel. Use the coefficient
errors and scalar normalization from Section~\ref{B:sec:inverse-representations-and-the-stochastic-bessel-operator}, so, on the selected cells,
\[
 C_n^{-1}p_{n,k}=p^B(z_{n,k})e^{\delta_{n,k}},\qquad
 C_n^2q_{n,k}=Q_{n,k}^B e^{-2\delta_{n,k}+\rho_{n,k}}.
\]
For a finite real scalar $\chi$, put
\begin{equation}
 C_n^{[\chi]}=e^{-\chi}C_n,\qquad
 \delta_{n,k}^{[\chi]}=\delta_{n,k}+\chi,\qquad
 \rho_{n,k}^{[\chi]}=\rho_{n,k}.
 \label{B:eq:gaussian-scalar-normalization}
\end{equation}
Multiplication by $e^\chi$ in the first identity and by $e^{-2\chi}$ in the second gives
\[
 (C_n^{[\chi]})^{-1}p_{n,k}=p^B(z_{n,k})e^{\delta_{n,k}^{[\chi]}},\qquad
 (C_n^{[\chi]})^2q_{n,k}
       =Q_{n,k}^Be^{-2\delta_{n,k}^{[\chi]}+\rho_{n,k}^{[\chi]}}.
\]
For every summand in the inverse kernel the exponential factor is unchanged, since
\[
 (\delta_{n,i}+\chi)+(\delta_{n,j}+\chi)
                -2(\delta_{n,l}+\chi)+\rho_{n,l}
   =\delta_{n,i}+\delta_{n,j}-2\delta_{n,l}+\rho_{n,l}.
\]
In particular, writing $\widehat p_{n,k}=C_n^{-1}p_{n,k}$ and $\widehat q_{n,k}=C_n^2q_{n,k}$, each summand of the positive inverse obeys
\[
 (e^\chi\widehat p_{n,i})(e^\chi\widehat p_{n,j})
       (e^{-2\chi}\widehat q_{n,l})
       =\widehat p_{n,i}\widehat p_{n,j}\widehat q_{n,l}.
\]
Thus the scalar choice leaves the finite inverse, its law, and every comparison embedding unchanged. For $\chi=\chi_n^\nabla$ the exact partial sum identity is
\[
 \chi_n^\nabla+\sum_{k=j}^{M_n}(\eta_{n,k}^G-B_{n,k})
 =\Xi_{n,M_n+1}+\Xi_{n,j}-\Xi_{n,M_n+1}
                       +\sum_{\ell=1}^3C_{n,j}^{(\ell)}
 =V_{n,j}^G.
\]
This removes the same value at the upper index for every $j$. The scalar $e^{-\chi_n^\nabla}$ is positive and $\mathcal N_n\vee\mathcal U$-measurable by
Theorem~\ref{B:thm:joint-gaussian-brownian-construction}. Corollary~\ref{B:cor:gaussian-coupling-almost-sure-bounds} gives $\chi_n^\nabla\to0$, hence $e^{-\chi_n^\nabla}\to1$.
Section~\ref{B:sec:normalization-and-coefficient-estimates} applies this normalization to the subsequently constructed finite coefficients.

\subsection{Transport of logarithmic ratios and construction of the coupling}
\label{B:sec:transport-of-logarithmic-ratios-and-construction-of-the-coupling}

The inverse formula involves products of coefficient ratios. After division by the corresponding products of conditional minimizer coordinates, these products have logarithm
\[
 \log\prod_{k=j}^{R}
       \frac{y_km_{x,k}^q}{x_km_{y,k}^q}
   =\sum_{k=j}^{R}\eta_k,\qquad j\in J.
\]
We therefore compare the original conditional distribution of $\eta$ with the centered Gaussian distribution constructed in
Section~\ref{B:sec:gaussian-approximation-of-the-coefficient-field} in the quadratic form defined by these partial sums. We begin with the change of variables and the two marginal
densities, then prove the coupling estimate. Throughout, the distributions are conditioned on the same fixed exterior data $q$.

Theorem~\ref{B:thm:logarithmic-ratio-coupling} states the moment estimate. Its proof in Section~\ref{B:subsection:proof-of-ratio-coupling} uses Poisson equations under the extended
conditional distribution, the discrete heat kernel estimate, and stationary diffusions. The subsequent measurable selection and conditional sampling retain the Gaussian vector
already coupled to the Brownian observations and recover the original coefficient distribution.

\subsubsection{The conditional logarithmic distribution}
\label{B:subsection:conditional-logarithmic-distribution}

The ratio comparison is made under one fixed conditional distribution. The interval choices ensure that its logarithmic powers are unchanged when the original measure is replaced
by the reference measure and that the estimates remain uniform over the admissible boundary configurations.

Fix $V$, $a>0$, $\beta>0$, and $A_0>20$. Retain the intervals $J=J_n=\{r_n,\ldots,R_n\}$, $I$, and $I^{\rm ret}$ from \eqref{B:eq:conditional-intervals}. Write $r=r_n$, $R=R_n$,
$m_n=r+\lceil r/4\rceil$, $M_n=R-\lceil r/4\rceil$, and $d=|J|$. The starting size is chosen so that
\begin{equation}
 r<n/4,\qquad r\ge8/\beta,\qquad n\beta\ge1,
 \quad\mathscr S(J)\text{ is disjoint from the modified initial coordinates}.
 \label{B:eq:ratio-parameter-range}
\end{equation}
The conditional distribution $\mu_q$, Hamiltonian $H^q$, minimizer $m^q$, and matrices $A,C$ are those of
\eqref{B:eq:conditional-hamiltonian}--\eqref{B:eq:conditional-matrix-data}. Their definitions extend to every finite nonnegative boundary configuration by
Proposition~\ref{B:prop:unchanged-conditional-distributions}. Unless specified otherwise, estimates uniform in $q$ refer to $q\in\mathcal Q_n^\sharp$. All constants may depend on
the fixed parameters above. The polynomial part is denoted by $\mathcal P_{n,J}^q$ as in Section~\ref{B:sec:general-potentials-and-convex-reference-measures}.

The logarithmic ratios $\eta_k$ are defined in \eqref{B:eq:conditional-logarithms}. Introduce the additional coordinates
\begin{equation}
 s_k=\frac12\left(\log\frac{x_k}{m_{x,k}^q}
                      +\log\frac{y_k}{m_{y,k}^q}\right),\qquad
 x_k=m_{x,k}^qe^{s_k-\eta_k/2},\quad
 y_k=m_{y,k}^qe^{s_k+\eta_k/2}.
 \label{B:eq:ratio-coordinate-map}
\end{equation}
We denote the inverse map by $\mathfrak Z_q(s,\eta)$. Every $m_{\sigma,k}^q$ is strictly positive. The two displayed exponentials therefore give a smooth bijection from
$\mathbb R^d\times\mathbb R^d$ onto $(0,\infty)^{2d}$, whose smooth inverse is the pair of logarithmic formulas. In site order its Jacobian is block diagonal, with
\[
 D_{(s_k,\eta_k)}(x_k,y_k)
 =\begin{pmatrix}x_k&-x_k/2\\y_k&y_k/2\end{pmatrix},
 \qquad
 |\det D\mathfrak Z_q(s,\eta)|
 =\prod_{k\in J}x_ky_k
 =\prod_{k\in J}m_{x,k}^qm_{y,k}^q\,e^{2s_k}.
\]
The absolute determinant has the same value if coordinates are grouped by species. Multiplication of the original density by this determinant increases each coordinate exponent by
one. At site $k$ the resulting power is
\[
 \begin{split}
 x_k^{c_{x,k}+1}y_k^{c_{y,k}+1}
 &=(m_{x,k}^q)^{\beta(k+a)}(m_{y,k}^q)^{\beta k}\\
 &\quad\times
  \exp\{\beta(k+a)(s_k-\eta_k/2)+\beta k(s_k+\eta_k/2)\}.
 \end{split}
\]
The exponent equals $\beta(2k+a)s_k-\beta a\eta_k/2$. Set
\[
 a_q^0=\prod_{k\in J}(m_{x,k}^q)^{\beta(k+a)}
                         (m_{y,k}^q)^{\beta k}>0.
\]
Dividing the transformed density by this factor gives
\begin{equation}
 g_q(s,\eta)=\exp\left\{-\kappa\mathcal P_{n,J}^q(\mathfrak Z_q(s,\eta))
   +\sum_{k\in J}\beta(2k+a)s_k-\frac{\beta a}{2}\sum_{k\in J}\eta_k\right\}.
 \label{B:eq:ratio-joint-density}
\end{equation}
Set
\begin{equation}
 I_q(\eta)=\int_{\mathbb R^d}g_q(s,\eta)\,\dd s,\qquad
 \mathcal Z_q^{\rm rat}=\int_{\mathbb R^d}I_q(\eta)\,\dd\eta,
 \qquad \nu_q=\mathcal L_{\mu_q}(\eta).
 \label{B:eq:ratio-normalizing-integrals}
\end{equation}
Writing the coefficient normalizer as $\mathcal Z_q^{\rm coef}=\int_{(0,\infty)^{2d}}e^{-\kappa H^q(z)}\,\dd z$, the change of variables gives
\[
 e^{-\kappa H^q(\mathfrak Z_q(s,\eta))}
       |\det D\mathfrak Z_q(s,\eta)|=a_q^0g_q(s,\eta),
 \qquad \mathcal Z_q^{\rm coef}=a_q^0\mathcal Z_q^{\rm rat}.
\]
In particular the joint probability density of $(s,\eta)$ is $g_q/\mathcal Z_q^{\rm rat}$. These identities hold for each fixed value of $q$ and its associated conditional density.

\begin{proposition}\label{B:prop:logarithmic-ratio-density}
For every $q$, the function $\eta\mapsto I_q(\eta)$ is finite, strictly positive, and smooth. The density of the logarithmic ratios and the conditional distribution of $s$ are
\begin{equation}
 f_q(\eta)=\frac{I_q(\eta)}{\mathcal Z_q^{\rm rat}},\qquad
 \mathfrak r_q(\dd s\mid\eta)=\frac{g_q(s,\eta)}{I_q(\eta)}\,\dd s.
 \label{B:eq:ratio-marginal-and-conditional}
\end{equation}
These formulas define Borel probability kernels for every value of their parameters.
\end{proposition}
\begin{proof}
Fix $n$ and a compact set $\mathcal Q$ of finite nonnegative boundary data. Separate evenness gives $\nabla\mathcal P_{n,J}^q(0)=0$. Integrating the Hessian on $[0,z]$ therefore
gives
\[
 \begin{split}
 \mathcal P_{n,J}^q(z)-\mathcal P_{n,J}^q(0)
 &=\int_0^1(1-t)\,
      z^{\mathsf T}\nabla^2\mathcal P_{n,J}^q(tz)z\,\dd t\\
 &\ge c_*|z|^2/2.
 \end{split}
\]
The values of $\mathcal P_{n,J}^q(0)$ are bounded on $\mathcal Q$. On $\mathcal Q\times[1,2]^{2d}$ the polynomial has a finite maximum, so integration on that box bounds
$\mathcal Z_q^{\rm coef}$ below by one positive constant. Its normalized density is consequently bounded by
\begin{equation}
 C\exp(-\kappa c_*|z|^2/2)\prod_{\alpha\in\mathscr S(J)}z_\alpha^{c_\alpha}.
 \label{B:eq:ratio-coordinate-majorant}
\end{equation}
For any $a>0$ and $c>-1$,
\[
 \int_0^\infty e^{-az^2}z^c\,\dd z
   =\frac12 a^{-(c+1)/2}\Gamma((c+1)/2).
\]
The substitution is $v=az^2$, with $\dd z=\tfrac12a^{-1/2}v^{-1/2}\dd v$. Each $c_\alpha>1$, so the product of these integrals is finite. The integral on the box is positive. Thus
$0<\mathcal Z_q^{\rm coef}<\infty$, and the normalizer identity above proves $0<\mathcal Z_q^{\rm rat}<\infty$.

Fix also a compact set of $\eta$ values. Continuity and strict positivity of $m^q$ give positive lower and finite upper bounds for all its coordinates on $\mathcal Q$. On the
resulting compact parameter set,
\[
 |\mathfrak Z_q(s,\eta)|^2
 =\sum_{k\in J}e^{2s_k}
     \bigl[(m_{x,k}^q)^2e^{-\eta_k}
                  +(m_{y,k}^q)^2e^{\eta_k}\bigr]
 \ge c\sum_{k\in J}e^{2s_k}.
\]
The term $-\beta a\sum_k\eta_k/2$ is bounded there. Hence
\begin{equation}
 g_q(s,\eta)\le C\prod_{k\in J}
            \exp\{-ce^{2s_k}+\beta(2k+a)s_k\}.
 \label{B:eq:ratio-density-majorant}
\end{equation}
The constants in this bound may depend on the compact sets and on $n$.

Differentiation in $\eta_k$ acts on a polynomial in $(x,y)$ by $\tfrac12(y_k\partial_{y_k}-x_k\partial_{x_k})$. In particular,
\[
 \partial_{\eta_k}\log g_q
  =-\frac{\kappa}{2}
      \bigl(y_k\partial_{y_k}\mathcal P_{n,J}^q
                    -x_k\partial_{x_k}\mathcal P_{n,J}^q\bigr)
                      (\mathfrak Z_q(s,\eta))-\frac{\beta a}{2}.
\]
Every further $\eta$ derivative is a polynomial in the coordinates of $\mathfrak Z_q$, with a possible constant term. For a labelled list $i_1,\ldots,i_h$, the product rule is
\[
 \partial_{i_1}\cdots\partial_{i_h}g_q
  =g_q\sum_{\pi\in\mathcal P(\{1,\ldots,h\})}
                   \prod_{D\in\pi}\partial_{i_D}\log g_q,
 \qquad
 \partial_{i_D}=\prod_{\ell\in D}\partial_{i_\ell}.
\]
Repeated coordinate indices are allowed in this formula. Each product on the right is a finite sum of monomials $\prod_ke^{b_ks_k}$ with $b_k\ge0$, times coefficients bounded on
the chosen compact parameter set. For $A,B>0$,
\[
 \int_{\mathbb R}e^{-Ae^{2s}+Bs}\,\dd s
   =\frac12A^{-B/2}\int_0^\infty e^{-v}v^{B/2-1}\,\dd v
   =\frac12A^{-B/2}\Gamma(B/2),
 \qquad v=Ae^{2s}.
\]
Here $B=\beta(2k+a)+b_k>0$. A finite sum of products of these integrable functions dominates every derivative of a fixed order. The mean-value formula for difference quotients and
dominated convergence therefore give, for each such list,
\[
 \partial_{i_1}\cdots\partial_{i_h}I_q(\eta)
   =\int_{\mathbb R^d}
       \partial_{i_1}\cdots\partial_{i_h}g_q(s,\eta)\,\dd s.
\]
The same majorant proves continuity of each derivative, so $I_q$ is smooth in $\eta$. Its value is positive at every $\eta$, because $g_q(s,\eta)>0$ for every $s$.

All expressions are Borel in their parameters, including $m^q$. Integrating their nonnegative densities over Borel sets proves measurability of the kernels. The denominators are
finite and positive at every parameter. Their masses are one by the definitions of $I_q$ and $\mathcal Z_q^{\rm rat}$. For any nonnegative Borel function $\Psi$ on the coefficient
space, Tonelli and the change of variables yield
\[
 \begin{split}
 &\int_{\mathbb R^d}f_q(\eta)
      \int_{\mathbb R^d}\Psi(\mathfrak Z_q(s,\eta))
                   \,\mathfrak r_q(\dd s\mid\eta)\,\dd\eta\\
 &\qquad=\frac1{\mathcal Z_q^{\rm rat}}
       \int_{\mathbb R^{2d}}\Psi(\mathfrak Z_q(s,\eta))
                                     g_q(s,\eta)\,\dd s\,\dd\eta
       =\int\Psi(z)\,\mu_q(\dd z).
 \end{split}
\]
Define $\eta_q(z)=(\log(y_k/m_{y,k}^q) -\log(x_k/m_{x,k}^q))_{k\in J}$. For bounded Borel $F$ on $\mathbb R^d$ and bounded Borel $\Psi$ on the coefficient space, the same change of
variables gives
\[
 \int F(\eta_q(z))\Psi(z)\,\mu_q(\dd z) 
 =\int F(\eta)
       \left[\int\Psi(\mathfrak Z_q(s,\eta))
                       \,\mathfrak r_q(\dd s\mid\eta)\right]
                                  \nu_q(\dd\eta).
\]
The quantity in brackets is a Borel function of $\eta$. Taking $F$ to be the indicator of any Borel set is the defining identity for a conditional distribution given $\eta_q$.
Since $0<I_q(\eta)<\infty$ at every $\eta$, the displayed normalized density defines this version at every $\eta$, including points in a $\nu_q$-null set. The same integrable
majorant gives every fixed logarithmic moment. For $b>0$, $c>-1$, and $p\ge0$, splitting at one gives
\[
 \int_0^\infty z^c|\log z|^p e^{-bz^2}\,\dd z
 \le \frac{\Gamma(p+1)}{(c+1)^{p+1}}
 +\frac12 b^{-(c+p+1)/2}\Gamma((c+p+1)/2).
\]
On $(0,1)$ this follows from $z=e^{-t}$ and $e^{-bz^2}\le1$; on $[1,\infty)$ use $\log z\le z$ and then $t=bz^2$. For the product majorant in
\eqref{B:eq:ratio-coordinate-majorant}, Tonelli applies this estimate to the coordinate containing the logarithmic factor and the same estimate with $p=0$ to every other
coordinate. A finite sum over coordinates bounds any fixed power of the Euclidean norm of their logarithms. On each compact boundary set at fixed $n$, the normalizer has a positive
lower bound and the minimizer coordinates are bounded above and below; the translated logarithmic ratios therefore have locally uniform moments of every fixed order.

\end{proof}

\begin{lemma}\label{B:lem:logarithmic-marginal-continuity}
At fixed $n$, $q\mapsto\nu_q$ is continuous in total variation and has locally uniformly integrable moments of every finite order.
\end{lemma}
\begin{proof}
Proposition~\ref{B:prop:logarithmic-ratio-density} gives, on each compact set $\mathcal Q$ of boundary parameters at fixed $n$, a common integrable bound for the coefficient
densities and their logarithmic moments of every fixed order. If $q_j\to q$ in this set, dominated convergence yields convergence of the normalizers and
\[
 \|p_{q_j}^{\rm coef}-p_q^{\rm coef}\|_1\longrightarrow0.
\]
Pushforward by the fixed map $(x,y)\mapsto(\log y_k-\log x_k)_{k\in J}$ contracts total variation. Denote its density by $\widehat f_q$ and put
$c(q)=(\log m_{y,k}^q-\log m_{x,k}^q)_{k\in J}$. Then $f_q(v)=\widehat f_q(v+c(q))$, $c(q_j)\to c(q)$, and
\[
 \|f_{q_j}-f_q\|_1
 \le\|p_{q_j}^{\rm coef}-p_q^{\rm coef}\|_1
 +\|\widehat f_q(\,\cdot+c(q_j))-
                     \widehat f_q(\,\cdot+c(q))\|_1\longrightarrow0.
\]
The last limit uses continuity of translation on $L^1$, obtained by approximating $\widehat f_q$ by compactly supported continuous functions. For every fixed $p>0$, the logarithmic
moment bound at order $p+1$ and boundedness of $c(q)$ on $\mathcal Q$ give
\[
 \sup_{q\in\mathcal Q}\int_{\{|\eta|>L\}}|\eta|^p\,\nu_q(\dd\eta)
 \le L^{-1}\sup_{q\in\mathcal Q}\int|\eta|^{p+1}\,\nu_q(\dd\eta)
 \longrightarrow0.
\]
Truncation and total variation convergence prove convergence of each fixed moment. These constants are local in $q$ at fixed $n$; uniform large-$n$ estimates retain the
admissibility assumptions of Section~\ref{B:sec:conditional-estimates-on-the-interior-interval}.
\end{proof}

\subsubsection{The Gaussian marginal}
\label{B:subsection:gaussian-logarithmic-marginal}

The comparison concerns centered linear ratios of the Gaussian coordinates. Their logarithms are not used, since the Gaussian vector can have negative coordinates. The mean of the
original logarithmic ratios is treated separately.

Retain $L_q$, $K_q$, and the corrected full Gaussian distribution $\gamma_{n,q}$ from \eqref{B:eq:gaussian-coefficient-law}. Write $\Gamma_q=N(0,K_q)$. If $G\sim\gamma_{n,q}$ and
$v=G-m^q-s^q$, then $v\sim N(0,\kappa^{-1}C)$ and
\[
 u=L_qv,\qquad
 \mathbb E e^{\ii t^{\mathsf T}u}
   =\exp\{-\tfrac1{2\kappa}t^{\mathsf T}L_qCL_q^{\mathsf T}t\},
 \qquad t\in\mathbb R^d.
\]
Thus $u$ has law $\Gamma_q$. This formula uses the centered linear displacement $v$. The full Gaussian distribution is defined on $\mathbb R^{2d}$, whereas the original
coefficients are positive and their logarithms have law $\nu_q$ after taking differences. The means are
\[
 \int u_k\,\Gamma_q(\dd u)=0,\qquad
 \int\eta_k\,\nu_q(\dd\eta)
    =a_k^{\rm lead}+O_\delta(k^{-2})
       \quad(q\in\mathcal Q_n^\sharp),
\]
by Corollary~\ref{B:cor:original-conditional-mean-expansion}. The estimate below compares these two marginals with their specified means.

The two nonzero entries in row $k$ of $L_q$ are $-1/m_{x,k}^q$ and $1/m_{y,k}^q$. Hence
\[
 |L_q^{\mathsf T}t|^2
   =\sum_{k\in J}\bigl[(m_{x,k}^q)^{-2}
                       +(m_{y,k}^q)^{-2}\bigr]t_k^2,
 \qquad
 t^{\mathsf T}K_qt
   =\kappa^{-1}|C^{1/2}L_q^{\mathsf T}t|^2>0
       \quad(t\ne0).
\]
This proves surjectivity of $L_q$ and positive definiteness of $K_q$ for every finite nonnegative $q$. Continuity of $m^q$ and $A$ gives continuity of $K_q$.

Let $T_q=D\mathfrak Z_q(0,0)$ and partition $\mathsf B_q=T_q^{\mathsf T}AT_q$ according to the two coordinate groups. The next calculation uses linear coordinates
$(\widehat s,u)=T_q^{-1}v$ for the centered Gaussian displacement. The symbol $s$ in $\mathfrak Z_q(s,\eta)$ continues to denote the exact nonlinear logarithmic coordinate.

\begin{lemma}\label{B:lem:gaussian-ratio-schur-complement}
The covariance of the Gaussian linear observables is
\begin{equation}
 K_q=\kappa^{-1}
 (\mathsf B_{\eta\eta}-\mathsf B_{\eta s}
           \mathsf B_{ss}^{-1}\mathsf B_{s\eta})^{-1}.
 \label{B:eq:ratio-gaussian-schur-complement}
\end{equation}
Moreover $q\mapsto\Gamma_q$ is continuous in total variation, with locally uniformly integrable moments of every finite order.
\end{lemma}
\begin{proof}
Suppress $q$ in the minimizer and put $M_x=\operatorname{diag}_{k\in J}(m_{x,k})$ and $M_y=\operatorname{diag}_{k\in J}(m_{y,k})$. In species order, direct multiplication gives
\[
 T_q=
 \begin{pmatrix}M_x&-M_x/2\\M_y&M_y/2\end{pmatrix},
 \qquad
 T_q^{-1}=
 \begin{pmatrix}M_x^{-1}/2&M_y^{-1}/2\\-M_x^{-1}&M_y^{-1}\end{pmatrix}.
\]
In particular, the second block row is $L_q$ and $\mathsf B_q^{-1}=T_q^{-1}CT_q^{-\mathsf T}$. Its lower right block is
$L_qCL_q^{\mathsf T}$.

Both $A$ and $\mathsf B_q$ are positive definite. Put
\[
 S_q=\mathsf B_{\eta\eta}
        -\mathsf B_{\eta s}\mathsf B_{ss}^{-1}\mathsf B_{s\eta}.
\]
The principal block $\mathsf B_{ss}$ is positive definite. For $u\ne0$, setting $\widehat s=-\mathsf B_{ss}^{-1}\mathsf B_{s\eta}u$ shows
$u^{\mathsf T}S_qu =(\widehat s,u)^{\mathsf T}\mathsf B_q(\widehat s,u)>0$. The block factorization is
\[
 \mathsf B_q=
 \begin{pmatrix}I&0\\\mathsf B_{\eta s}\mathsf B_{ss}^{-1}&I\end{pmatrix}
 \begin{pmatrix}\mathsf B_{ss}&0\\0&S_q\end{pmatrix}
 \begin{pmatrix}I&\mathsf B_{ss}^{-1}\mathsf B_{s\eta}\\0&I\end{pmatrix}.
\]
Inverting the three factors in reverse order gives lower right block $S_q^{-1}$. Comparison with $T_q^{-1}CT_q^{-\mathsf T}$ therefore yields $L_qCL_q^{\mathsf T}=S_q^{-1}$, which
proves \eqref{B:eq:ratio-gaussian-schur-complement}.

The same identity can be checked directly from the density. The centered vector $v$ has density
\[
 \left(\frac{\kappa}{2\pi}\right)^d(\det A)^{1/2}
                         e^{-\kappa v^{\mathsf T}Av/2}.
\]
Since $(\det A)^{1/2}|\det T_q|=(\det\mathsf B_q)^{1/2}$, the density of $(\widehat s,u)$ is
\[
 \left(\frac{\kappa}{2\pi}\right)^d
 (\det\mathsf B_q)^{1/2}
 \exp\{-\tfrac{\kappa}{2}
           (\widehat s,u)^{\mathsf T}\mathsf B_q(\widehat s,u)\}.
\]
Completing the square gives
\[
 \begin{split}
 (\widehat s,u)^{\mathsf T}\mathsf B_q(\widehat s,u)
 &=\bigl(\widehat s+\mathsf B_{ss}^{-1}\mathsf B_{s\eta}u\bigr)^{\mathsf T}
       \mathsf B_{ss}
       \bigl(\widehat s+\mathsf B_{ss}^{-1}\mathsf B_{s\eta}u\bigr)
       +u^{\mathsf T}S_qu.
 \end{split}
\]
The integral of the first exponential in $\widehat s$ is $(2\pi/\kappa)^{d/2}(\det\mathsf B_{ss})^{-1/2}$. The factorization gives $\det\mathsf B_q=\det\mathsf B_{ss}\det S_q$.
Thus the marginal density of $u$ is exactly
\[
 \frac{(\det K_q)^{-1/2}}{(2\pi)^{d/2}}
             e^{-u^{\mathsf T}K_q^{-1}u/2}
 =\left(\frac{\kappa}{2\pi}\right)^{d/2}
            (\det S_q)^{1/2}e^{-\kappa u^{\mathsf T}S_qu/2}.
\]
The conditional Gaussian distribution of $\widehat s$ given $u$ has mean $-\mathsf B_{ss}^{-1}\mathsf B_{s\eta}u$ and covariance $\kappa^{-1}\mathsf B_{ss}^{-1}$. This conditional
Gaussian distribution and the exact density $\mathfrak r_q(\dd s\mid\eta)$ describe the two distributions used in the subsequent reconstruction.

At fixed $n$, continuity and positivity of $K_q$ imply that on a compact set $\mathcal Q$ of boundary data there are $0<c_{\mathcal Q}\le C_{\mathcal Q}<\infty$ with
$c_{\mathcal Q}I\le K_q\le C_{\mathcal Q}I$. The Gaussian density is consequently bounded by
\[
 (2\pi c_{\mathcal Q})^{-d/2}
            \exp\{-|u|^2/(2C_{\mathcal Q})\}.
\]
For $q_j\to q$ these densities converge pointwise, since their determinants and inverses converge. Dominated convergence gives convergence in $L^1$. Multiplication by $|u|^{p+1}$
preserves an integrable majorant for each fixed $p>0$, and hence
\[
 \sup_{q\in\mathcal Q}
 \int_{\{|u|>L\}}|u|^p\,\Gamma_q(\dd u)
 \le L^{-1}\sup_{q\in\mathcal Q}
                \int|u|^{p+1}\,\Gamma_q(\dd u)\longrightarrow0.
\]
This proves continuity and local uniform integrability.
\end{proof}

\subsubsection{Weighted partial sums}
\label{B:subsection:weighted-ratio-partial-sums}

The errors enter the inverse through cumulative sums. We therefore use a cost that measures these sums with the weights from the kernel comparison.

Fix $1<\lambda<\min\{a+1,2\}$. Use the backward cell lengths $h_{n,j}=\log(z_{n,j}/z_{n,j-1})$ on $z_{n,j}=\theta_V(j/n)$. The weights of the coefficient comparison are
\begin{equation}
 w_{n,j}=h_{n,j}z_{n,j}^{\lambda}\asymp n^{-\lambda}j^{\lambda-1},\qquad
 (S_ne)_j=\sum_{k=j}^{R}e_k,\qquad
 Q_n(e)=\sum_{j\in J}w_{n,j}|(S_ne)_j|^2.
 \label{B:eq:ratio-partial-sum-cost}
\end{equation}
In the ordered coordinates $r,\ldots,R$, the matrix of $S_n$ has entries $(S_n)_{jk}=\mathbf1_{\{j\le k\}}$. If $y=S_ne$ and $y_{R+1}=0$, then
\[
 e_j=y_j-y_{j+1}\quad(r\le j\le R),\qquad
 (S_n^{-1})_{jk}=\mathbf1_{\{j=k\}}-\mathbf1_{\{k=j+1\}}.
\]
Both boundary point conventions enter this inverse formula. Since every $w_{n,j}>0$, it follows that $\sqrt{Q_n}$ is a norm. Its matrix and its entries are
\[
 Q_n(e)=e^{\mathsf T}S_n^{\mathsf T}
                 \operatorname{diag}(w_{n,j})S_ne,\qquad
 [S_n^{\mathsf T}\operatorname{diag}(w_{n,j})S_n]_{kl}
                 =\sum_{j=r}^{k\wedge l}w_{n,j}.
\]
These formulas include the products of distinct coordinates. They also give a finite-dimensional lower bound:
\[
 |e|^2=\sum_{j=r}^R|y_j-y_{j+1}|^2
       \le4\sum_{j=r}^R|y_j|^2
       \le\frac{4}{\min_{j\in J}w_{n,j}}Q_n(e).
\]
For example, if $E_k$ is the $k$th coordinate vector and $k<R$, then $S_n(E_k-E_{k+1})$ equals $-E_{k+1}$, so $Q_n(E_k-E_{k+1})=w_{n,k+1}$. A single coordinate instead satisfies
$Q_n(E_k)=\sum_{j=r}^kw_{n,j}$. This is the effect of summing the logarithms before estimating the coefficient products.

All sums in this subsection run over the full interval $J$. Define
\begin{equation}
 \mathsf T_2=n^{-1}\operatorname{diag}_{j\in J}(\sqrt j)S_n,
 \qquad \mathsf Q_2(e)=|\mathsf T_2e|^2,
 \qquad \ell_n=\log(\me+n).
 \label{B:eq:ratio-auxiliary-quadratic-form}
\end{equation}
Explicitly, $\mathsf Q_2(e)=n^{-2}\sum_{j=r}^Rj|(S_ne)_j|^2$. The ratio of the two weights is bounded by
\[
 \frac{w_{n,j}}{j/n^2}
   \le C n^{2-\lambda}j^{\lambda-2}
   \le C(n/r)^{2-\lambda}.
\]
The last inequality uses $j\ge r$ and $\lambda-2<0$. Multiplication by $|(S_ne)_j|^2$ and summation, and then retention of a single term in $\mathsf Q_2$, give
\begin{equation}
 Q_n(e)\le C(n/r)^{2-\lambda}\mathsf Q_2(e),\qquad
 \sqrt j\,|(S_ne)_j|\le n|\mathsf T_2e|.
 \label{B:eq:ratio-norm-comparisons}
\end{equation}
For an arbitrary coupling $\pi$ of $\nu_q$ and $\Gamma_q$, the moments proved above imply that the integrals of these quadratic forms are finite. Indeed,
\[
 Q_n(\eta-u)\le2Q_n(\eta)+2Q_n(u)
          \le2\|S_n^{\mathsf T}\operatorname{diag}(w_{n,j})S_n\|
                           (|\eta|^2+|u|^2).
\]
For every $p\ge1$, the inequalities $|v-w|^2\le2(|v|^2+|w|^2)$ and $(a+b)^p\le2^{p-1}(a^p+b^p)$ for $a,b\ge0$ give
\[
 \begin{split}
 Q_n(\eta-u)^p
 &\le2^{2p-1}
   \|S_n^{\mathsf T}\operatorname{diag}(w_{n,j})S_n\|^p
                     (|\eta|^{2p}+|u|^{2p}),\\
 \mathsf Q_2(\eta-u)^p
 &\le2^{2p-1}\|\mathsf T_2\|^{2p}
                     (|\eta|^{2p}+|u|^{2p}).
 \end{split}
\]
Integrating the right sides uses only the two marginals. Their moments are finite at each fixed $n,q$. The next theorem gives the required dependence on $n$ and on the moment order
as the interval grows. Since $u$ is centered, convexity also gives
\[
 \mathsf Q_2(\mathbb E_{\nu_q}\eta)
       =\mathsf Q_2(\mathbb E_\pi(\eta-u))
       \le\mathbb E_\pi\mathsf Q_2(\eta-u).
\]
Thus the discrepancy of the two means is included in the cost.

The argument below applies to any nonnegative weights bounded above by the right hand scale in \eqref{B:eq:ratio-partial-sum-cost}; the specific positive weights are used when
choosing optimal couplings.

\begin{theorem}\label{B:thm:logarithmic-ratio-coupling}
For every fixed $C_0<\infty$, sufficiently large $n$, and $q\in\mathcal Q_n^\sharp$, there is a coupling of $\eta\sim\nu_q$ and $u\sim\Gamma_q$ satisfying, simultaneously for
$2\le p\le C_0\ell_n$,
\begin{equation}
 \|\mathsf Q_2(\eta-u)\|_{L^p}
       \le Cp\ell_n^2(p+\ell_n)n^{-2}.
 \label{B:eq:ratio-coupling-moment-bound}
\end{equation}
The constant is uniform in $q,n,p$; the starting size may depend on $C_0$.
\end{theorem}
The coupling in the theorem has both specified marginals and is fixed for all $p$ in the displayed range. For its error $e=\eta-u$, the deterministic comparison yields
\[
 \|Q_n(e)\|_p
   \le C(n/r)^{2-\lambda}\|\mathsf Q_2(e)\|_p
   \le Cp\ell_n^2(p+\ell_n)n^{-\lambda}r^{\lambda-2}.
\]
For each $j$, and also for their pointwise maximum,
\[
 \left\|\max_{j\in J}\sqrt j\,|(S_ne)_j|\right\|_{2p}
   \le n\|\mathsf Q_2(e)^{1/2}\|_{2p}
   =n\|\mathsf Q_2(e)\|_p^{1/2}
   \le C\sqrt p\,\ell_n\sqrt{p+\ell_n}.
\]
In particular, we obtain
\begin{equation}
 \begin{split}
 \|Q_n(\eta-u)\|_{L^p}
   &\le Cp\ell_n^2(p+\ell_n)n^{-\lambda}r^{\lambda-2},\\
 \max_{j\in J}\sqrt j\,\|(S_n(\eta-u))_j\|_{L^{2p}}
   &\le C\sqrt p\,\ell_n\sqrt{p+\ell_n}.
 \end{split}
 \label{B:eq:ratio-coupling-consequences}
\end{equation}
For a fixed $p\ge2$ the first bound gives
\[
 \|nQ_n(\eta-u)\|_p
       \le C_p\ell_n^3n^{1-\lambda}r^{\lambda-2}
       \longrightarrow0.
\]
For $\varepsilon>0$, Markov's inequality on the same coupling therefore gives, uniformly over admissible $q$,
\[
 \Pr_q\{nQ_n(\eta-u)>\varepsilon\}
 \le(C_p/\varepsilon)^p
       \ell_n^{3p}n^{-p(\lambda-1)}r^{-p(2-\lambda)}.
\]
Choose a fixed $p$ with $p(\lambda-1)>1$; this order lies below $C_0\ell_n$ for all sufficiently large $n$ when $C_0>0$. For $0<\epsilon<p(\lambda-1)-1$, one has
$\ell_n^{3p}\le C_{p,\epsilon}n^\epsilon$ and $r\ge1$. The displayed probability bound is consequently summable in $n$. The maximum estimate likewise gives
\[
 \left\|\max_{j\in J}|(S_n(\eta-u))_j|\right\|_{2p}
       \le C_p\ell_n^{3/2}r^{-1/2}\longrightarrow0,
\]
since $r=\lceil(\log n)^{A_0}\rceil$ and $A_0>20$. These consequences follow from the stated coupling estimate, which we prove in
Section~\ref{B:subsection:proof-of-ratio-coupling}. We then select the conditional plans measurably and construct the joint coupling over all sizes. The sums defining both
quadratic forms still end at $R$. Restriction to the observation interval $I$ is treated in Section~\ref{B:subsection:restricted-ratio-partial-sums}.

The diffusion comparison uses a matrix obtained from Poisson solutions. The Poisson construction of Stein kernels under a spectral gap is developed in Theorem~2.4 of
\cite{courtadeFathiPananjady2019stein}. The following estimates describe how the entries depend on the coordinate indices under the conditional law used here.

\subsubsection{The Poisson equations}
\label{B:subsection:ratio-poisson-equations}

The transport estimate is obtained from integration by parts for the logarithmic observables. Solutions of the Poisson equations for the extended conditional measure provide a
diffusion matrix with the spatial decay required to estimate the partial sums.

Use $V_0$, $\widetilde\mu_q$, $l_\alpha$, $f$, and $L=f'(0)$ from Proposition~\ref{B:prop:smooth-conditional-extension} and \eqref{B:eq:extended-logarithmic-functions}. Expectation
under $\widetilde\mu_q$ is denoted by brackets. The covariance $B=LCL^{\mathsf T}$ in \eqref{B:eq:normalized-logarithmic-covariance} satisfies $b_-I\le B\le b_+I$ for fixed
positive constants. We use the self-adjoint generator from the proof of Lemma~\ref{B:lem:extended-covariance-decay}. We continue to work with real functions. Its closed symmetric
bilinear form is
\[
 \mathcal E_q(u,v)=\langle\nabla u\cdot\nabla v\rangle,
 \qquad \operatorname{dom}\mathcal E_q=H^1(\widetilde\mu_q),
\]
where $H^1(\widetilde\mu_q)$ is the completion of real $C_c^\infty(\mathbb R^{\mathscr S(J)})$ functions for the norm $(\|u\|_2^2+\langle|\nabla u|^2\rangle)^{1/2}$. An element
$u\in H^1(\widetilde\mu_q)$ belongs to $\operatorname{dom}\mathscr L$ precisely when there is $v\in L^2(\widetilde\mu_q)$ such that
\[
 \mathcal E_q(u,\psi)=-\langle v\psi\rangle
          \quad(\psi\in H^1(\widetilde\mu_q));
 \qquad\text{then }\mathscr Lu=v.
\]
On smooth functions in this domain, $\mathscr L=\Delta-\nabla V_0\cdot\nabla$. The graph approximation established in that proof makes $C_c^\infty$ a core. The constant function
belongs to the form domain with zero energy, so $\mathscr L1=0$. Its semigroup $P_t=e^{t\mathscr L}$ consequently satisfies
\[
 P_t1=1,\qquad \langle P_tv\rangle=\langle v\rangle,
 \qquad \|P_tv\|_p\le\|v\|_p\quad(1\le p\le\infty).
\]
For mean zero $v$, the Poincar\'e inequality gives $\|P_tv\|_2\le e^{-ct}\|v\|_2$. The bounds from Section~\ref{B:sec:conditional-estimates-on-the-interior-interval} are
\begin{equation}
 \begin{gathered}
 cI\le V_0''\le CI,\quad V_0''=A+U'',\quad
 |U''_{\alpha\gamma}|\le C\delta\mathbf1_{\{d(\alpha,\gamma)\le R_*\}},\\
 |\partial^h U|\le C_\delta k_*^{1-h/2}\quad(3\le h\le5),\qquad
 |\partial^h f_i|\le C_\delta i^{(1-h)/2}\quad(2\le h\le4).
 \end{gathered}
 \label{B:eq:ratio-extension-inputs}
\end{equation}
A nonzero derivative of $U$ has all its indices in a fixed interaction neighborhood; a derivative of $f_i$ uses only the two coordinates at $i$. The parameter $\delta$ is fixed
below, within the permitted range of \eqref{B:eq:conditional-extension-smallness}.

\begin{lemma}\label{B:lem:local-poisson-gradients}
There are smooth functions $h_i$ satisfying $-\mathscr Lh_i=f_i-\langle f_i\rangle$. Their gradients $g_{i,\alpha}=\partial_\alpha h_i$ obey
\begin{equation}
 |g_{i,\alpha}(x)|\le Ce^{-c|i-k_\alpha|},\qquad
 |g_{i,\alpha}(x)-(CL_i^{\mathsf T})_\alpha|
           \le C\delta e^{-c|i-k_\alpha|}.
 \label{B:eq:ratio-poisson-gradient-bounds}
\end{equation}
Both estimates hold on the whole fluctuation space.
\end{lemma}
\begin{proof}
Put $F_i=f_i-\langle f_i\rangle$. Since $f_i(0)=0$, $\|\nabla f_i\|_\infty\le C$, and $V_0(x)\ge c|x|^2/2$, we have $F_i\in L^2(\widetilde\mu_q)$ and $\langle F_i\rangle=0$. The
spectral gap therefore defines
\begin{equation}
 h_i=\int_0^\infty P_tF_i\,\dd t,
 \qquad -\mathscr Lh_i=F_i,\qquad \langle h_i\rangle=0.
 \label{B:eq:ratio-poisson-integral}
\end{equation}
Indeed, for $h_{i,T}=\int_0^TP_tF_i\,\dd t$,
\[
 \mathscr Lh_{i,T}=P_TF_i-F_i,\qquad
 \|h_i-h_{i,T}\|_2\le c^{-1}e^{-cT}\|F_i\|_2.
\]
Closedness of $\mathscr L$ proves that $h_i$ belongs to its domain. The identity $\mathcal E_q(u,u)\ge c\|u\|_2^2$ on mean zero functions proves uniqueness and gives
$\|h_i\|_2\le c^{-1}\|F_i\|_2$.

We justify smoothness at each fixed $n,q,\delta$ before proving estimates uniform in these parameters. Consider the additive noise equation
\[
 X_t^x=x-\int_0^t\nabla V_0(X_s^x)\,\dd s+\sqrt2\,\mathsf B_t.
\]
The drift is globally Lipschitz. Picard differences on $[0,T]$ are bounded by $A_TC^jt^j/j!$, so their sum converges uniformly to the unique global solution. The mean-value formula
and $V_0''\ge cI$ give
\[
 |X_t^x-X_t^y|\le e^{-ct}|x-y|,\qquad
 \sup_{s\le T}|X_s^x|
 \le (|x|+\sqrt2\sup_{s\le T}|\mathsf B_s|)e^{CT}.
\]
For a labelled set $E$ of initial directions, let $J_t^E=D_x^{|E|}X_t^x[v_a:a\in E]$. The first variation satisfies $\dot J_t=-V_0''(X_t^x)J_t$, $J_0=I$, and $\|J_t\|\le e^{-ct}$.
At orders $|E|\ge2$ its equation is
\[
 \dot J_t^E=-V_0''(X_t^x)J_t^E
 -\sum_{\substack{\pi\in\mathcal P(E)\\|\pi|\ge2}}
 D^{|\pi|}\nabla V_0(X_t^x)[J_t^D:D\in\pi],\qquad J_0^E=0.
\]
All positive order derivatives of the drift are bounded at this fixed dimension. Variation of constants and
\[
 \int_0^t e^{-c(t-s)}e^{-|\pi|cs}\,\dd s
 \le \frac{e^{-ct}}{(|\pi|-1)c}
\]
prove inductively $|J_t^E|\le C_Ee^{-ct}\prod_{a\in E}|v_a|$. To justify these flow derivatives, observe that the difference quotients obey the same integral equations with an
error $\varepsilon_h$ tending uniformly to zero on $[0,T]$, by the mean-value formula for the next drift derivative. Gronwall bounds the difference from the displayed variation by
$Te^{CT}\|\varepsilon_h\|_{\infty,[0,T]}$.

For $T_t\psi(x)=\mathbb E\psi(X_t^x)$, this differentiation and It\^o's formula imply $\partial_tT_t\psi=\mathscr LT_t\psi$ on $C_c^\infty$. Integration by parts with
$\chi_R(x)=\chi_0(x/R)$ gives
\[
 \langle\chi_R\mathscr LT_t\psi\rangle
 =-\langle\nabla\chi_R\cdot\nabla T_t\psi\rangle\longrightarrow0.
\]
The derivative bounds and Gaussian tails justify time integration, so $\langle T_t\psi\rangle=\langle\psi\rangle$. Density in $C_0$ and uniqueness of finite Borel measures extend
invariance to Borel functions. Jensen and path continuity then make $T_t$ a strongly continuous contraction semigroup on $L^2$. Its generator agrees with $\mathscr L$ on the smooth
core established in Lemma~\ref{B:lem:extended-covariance-decay}; consequently $T_t=P_t$. One may verify this semigroup identity on the dense range of
$(\lambda-\mathscr L)C_c^\infty$: the two resolvents both send $(\lambda-\mathscr L)\psi$ to $\psi$ and have norm at most $\lambda^{-1}$.

Couple $X^x$ to $X^Y$ with $Y\sim\widetilde\mu_q$ independent of the Brownian motion. Invariance, contraction, and the flow derivative bounds give
\[
 |P_tF_i(x)|\le Ce^{-ct}(|x|+\mathbb E|Y|),\qquad
 \|D^hP_tF_i\|_\infty\le C_{n,q,\delta,h}e^{-ct}\quad(h\ge1).
\]
Thus \eqref{B:eq:ratio-poisson-integral} converges pointwise, defines a smooth function of at most linear growth, and may be differentiated at every order. Every positive order
derivative of $h_i$ is bounded. The Poisson equation holds everywhere, and differentiation gives
\begin{equation}
 (-\mathscr L+V_0'')g_i=\nabla f_i.
 \label{B:eq:ratio-differentiated-poisson}
\end{equation}
Each positive order derivative $v$ of $h_i$ belongs to $\operatorname{dom}\mathscr L$. In fact $v,\nabla v,\nabla^2v$ are bounded, $\mathscr Lv$ has at most linear growth, and
\[
 \mathscr L(\chi_Rv)-\chi_R\mathscr Lv
 =2\nabla\chi_R\cdot\nabla v+v\mathscr L\chi_R\longrightarrow0
 \quad\text{in }L^2,\qquad \chi_Rv\longrightarrow v\quad\text{in }L^2.
\]
Here $|\mathscr L\chi_R|\le C_n\mathbf1_{\{R\le|x|\le2R\}}$; Gaussian tails control both errors.

For the uniform estimates, equip vector functions with $\|v\|_{i,\gamma,\infty} =\max_\alpha e^{\gamma|k_\alpha-i|}\|v_\alpha\|_\infty$. Lemma~\ref{B:lem:weighted-band-semigroups}
gives, for small fixed
$\gamma>0$,
\[
 \mathcal K_1=(-\mathscr L+A)^{-1}
 =\int_0^\infty e^{-tA}P_t\,\dd t,\qquad
 \|\mathcal K_1\|_{i,\gamma,\infty}\le K_1.
\]
The integral is first defined on vector-valued $L^2$ and, by pairing with $L^1$, on $L^\infty$. Both bounds follow from
\[
 C\int_0^\infty e^{-ct}\,\dd t\,
 \sup_\alpha\sum_\rho e^{-c_1d(\alpha,\rho)/2}<\infty.
\]
The inverse identity follows by integrating the product semigroup up to $T$: its image under $-\mathscr L+A$ is $I-e^{-TA}P_T$, and the latter error tends to zero in operator norm.
Finite range gives $\|U''\|_{i,\gamma,\infty}\le C_1\delta$. Choose $K_1C_1\delta\le1/2$. Then
\begin{equation}
 g_i=\sum_{l=0}^\infty(-\mathcal K_1U'')^l\mathcal K_1\nabla f_i.
 \label{B:eq:ratio-poisson-neumann-series}
\end{equation}
Indeed the series solves \eqref{B:eq:ratio-differentiated-poisson} in the generator domain, and the homogeneous equation has only the zero solution, since
\[
 \sum_\alpha\langle|\nabla v_\alpha|^2\rangle
 +\langle v^{\mathsf T}V_0''v\rangle\ge c\sum_\alpha\|v_\alpha\|_2^2.
\]
The right hand side $\nabla f_i$ is bounded and supported at site $i$; hence $\|g_i\|_{i,\gamma,\infty}\le C$.

For the constant vector $g_i^0=CL_i^{\mathsf T}$,
\[
 (-\mathscr L+V_0'')(g_i-g_i^0)
 =\nabla f_i-L_i^{\mathsf T}-U''g_i^0.
\]
The first term is supported at site $i$ and is bounded by $C\delta$. Exponential decay of $C$ and finite range of $U''$ give $|(U''g_i^0)_\alpha|\le C\delta e^{-c_2|k_\alpha-i|}$.
The same inverse series, with $\gamma<c_2$, yields $\|g_i-g_i^0\|_{i,\gamma,\infty}\le C\delta$. The functions are continuous and the density is positive, so these essential
supremum bounds hold at every point. This proves both bounds in \eqref{B:eq:ratio-poisson-gradient-bounds}.
\end{proof}

\par\begingroup
\dimen0=\pagegoal\advance\dimen0 by -\pagetotal \ifdim\dimen0<6\baselineskip\vfil\break\fi
\endgroup
\subsubsection{Derivatives of the Poisson solutions}
\label{B:subsection:ratio-poisson-spatial-derivatives}

The derivative of the diffusion matrix involves second derivatives of the Poisson solutions. The following estimate keeps their decay in both coordinate indices, which is needed
when these derivatives are summed.

\begin{lemma}\label{B:lem:local-poisson-second-derivatives}
The second derivatives $r_{i,\alpha\gamma}=\partial_\gamma g_{i,\alpha}$ satisfy
\begin{equation}
 |r_{i,\alpha\gamma}(x)|
 \le C k_\gamma^{-1/2}
       e^{-c(|i-k_\gamma|+|k_\alpha-k_\gamma|)}.
 \label{B:eq:ratio-poisson-second-bound}
\end{equation}
\end{lemma}
\begin{proof}
For a smooth scalar $v$, the product rule gives
\[
 \begin{split}
 \partial_\gamma\mathscr Lv
 &=\partial_\gamma\Delta v
   -\sum_\rho\partial_\gamma
                      [(\partial_\rho V_0)\partial_\rho v]\\
 &=\mathscr L\partial_\gamma v
       -\sum_\rho V''_{0,\gamma\rho}\partial_\rho v.
 \end{split}
\]
Differentiating the $\alpha$ component of \eqref{B:eq:ratio-differentiated-poisson} therefore gives
\begin{equation}
 \begin{split}
 -\mathscr L r_{i,\alpha\gamma}
 +\sum_\rho V''_{0,\alpha\rho}r_{i,\rho\gamma}
 +\sum_\rho V''_{0,\gamma\rho}r_{i,\alpha\rho}
 =\partial^2_{\alpha\gamma}f_i
        -\sum_\rho V'''_{0,\alpha\rho\gamma}g_{i,\rho}.
 \end{split}
 \label{B:eq:ratio-second-poisson-equation}
\end{equation}
The first Hessian term differentiates $V_0''g_i$ through $g_i$; the second is the commutator with $\mathscr L$. The remaining derivative of $V_0''$ is moved to the right. The flow
argument above gives bounded derivatives of $h_i$ of every positive order at fixed dimension. The cutoff approximation in the generator graph norm therefore applies to each
component of $r_i$.

For a matrix $v=(v_{\alpha\gamma})$ of component functions, the constant coefficient operator is $-\mathscr Lv+Av+vA$. Its inverse on matrix-valued $L^2$ is
\[
 \mathcal K_2v
 =\int_0^\infty e^{-tA}(P_tv)e^{-tA}\,\dd t
 =\int_0^\infty(e^{-tA}\otimes e^{-tA})P_tv\,\dd t.
\]
Here $P_t$ acts entry by entry and the tensor formula uses ordered pairs of indices. In fact, for the integral $v_T$ up to time $T$,
\[
 -\mathscr Lv_T+Av_T+v_TA
 =v-e^{-TA}(P_Tv)e^{-TA}.
\]
The norm of the last term is at most $e^{-2a_-T}\|v\|_2$. The integral converges in $L^2$, so closedness proves the inverse identity.

For a fixed pair $(\rho,\sigma)$ use
\[
 \|v\|_{\rho,\sigma,\gamma_0,\infty}
 =\max_{\alpha,\theta}
 e^{\gamma_0(d(\alpha,\rho)+d(\theta,\sigma))}
                            \|v_{\alpha\theta}\|_\infty.
\]
The heat kernel estimate of Lemma~\ref{B:lem:weighted-band-semigroups} and the triangle inequality for $d$ show, for sufficiently small $\gamma_0>0$, that the absolute weighted row
sum at time $t$ is bounded by
\[
 C^2e^{-2ct}
 \left(\sup_\alpha\sum_a e^{-c_1d(\alpha,a)/2}\right)^2.
\]
Its integral is a constant $K_2$ independent of $n,q,\rho,\sigma$. The same bound defines the weak-star integral on $L^\infty$. The remaining multiplication operator is
$\mathcal U_2v=U''v+vU''$. For its first term,
\[
 e^{\gamma_0(d(\alpha,\rho)+d(\theta,\sigma))}
                      \|(U''v)_{\alpha\theta}\|_\infty 
 \le C\delta
   \sum_{d(\alpha,a)\le R_*}e^{\gamma_0d(\alpha,a)}
                   \|v\|_{\rho,\sigma,\gamma_0,\infty}.
\]
For the other term the corresponding bound is
\[
 e^{\gamma_0(d(\alpha,\rho)+d(\theta,\sigma))}
        \|(vU'')_{\alpha\theta}\|_\infty
 \le C\delta\sum_{d(\theta,a)\le R_*}e^{\gamma_0d(\theta,a)}
                 \|v\|_{\rho,\sigma,\gamma_0,\infty}.
\]
Both sums have bounded cardinality. Adding them gives $\|\mathcal U_2\|\le C_2\delta$. Choose the fixed $\delta$ so that $K_2C_2\delta\le1/2$, as permitted by
\eqref{B:eq:conditional-extension-smallness}. Thus
\[
 \mathcal V=\sum_{l=0}^\infty(-\mathcal K_2\mathcal U_2)^l\mathcal K_2,
 \qquad \|\mathcal V\|_{\rho,\sigma,\gamma_0,\infty}\le2K_2.
\]
The series value $v=\mathcal V F$ satisfies $v=\mathcal K_2(F-\mathcal U_2v)$ and hence belongs to the constant coefficient operator domain. It solves
$-\mathscr Lv+V_0''v+vV_0''=F$. If $v$ solves the homogeneous equation in this domain, summing its scalar products with the components of $v$ gives
\[
 \begin{split}
 0={}&\sum_{\alpha,\gamma}\langle|\nabla v_{\alpha\gamma}|^2\rangle
  +\left\langle\operatorname{Tr}(v^{\mathsf T}V_0''v)
                        +\operatorname{Tr}(v^{\mathsf T}vV_0'')\right\rangle\\
 &\ge2c\sum_{\alpha,\gamma}\|v_{\alpha\gamma}\|_2^2.
 \end{split}
\]
For the last inequality the first trace is the sum of the quadratic forms of the columns of $v$, and the second is the sum of those of its rows. This proves uniqueness without any
symmetry assumption on $v$. It identifies the inverse with the smooth matrix $r_i$ in \eqref{B:eq:ratio-second-poisson-equation}. Applying the weighted bound to a matrix supported
in component $(\rho,\sigma)$ proves
\begin{equation}
 \|\mathcal V_{\alpha\gamma,\rho\sigma}\|_{L^\infty\to L^\infty}
       \le Ce^{-c(d(\alpha,\rho)+d(\gamma,\sigma))}.
 \label{B:eq:ratio-two-component-inverse}
\end{equation}

Write the right hand side of \eqref{B:eq:ratio-second-poisson-equation} as $F_{i,\rho\sigma}$. Since $f_i$ uses only the two coordinates at site $i$,
\[
 |\partial^2_{\rho\sigma}f_i|
 \le Ci^{-1/2}\mathbf1_{\{k_\rho=i=k_\sigma\}}.
\]
In its other term, $V_0'''=U'''$. A nonzero derivative $U'''_{\rho a\sigma}$ requires its three coordinates to belong to a fixed interaction neighborhood. There are a bounded
number of possible $a$. The derivative estimate and the first Poisson bound therefore give
\[
 \left|\sum_aU'''_{\rho a\sigma}g_{i,a}\right|
 \le C_\delta k_\rho^{-1/2}e^{-c|i-k_\rho|}
                         \mathbf1_{\{d(\rho,\sigma)\le R_*\}}.
\]
This formula includes repeated coordinate indices. Substitute the bounds for the two terms on the right hand side into \eqref{B:eq:ratio-two-component-inverse}. For the first, the
at most four pairs of species at site $i$ contribute $Ci^{-1/2}e^{-c(|k_\alpha-i|+|k_\gamma-i|)}$. For the second, put $v=k_\rho$. The index $\sigma$ is within $R_*$ of $v$, and
its possible sites and species have bounded cardinality, at most $2(2R_*+1)$. For each such coordinate,
\[
 |k_\gamma-k_\sigma|\ge|k_\gamma-v|-R_*,\qquad
 e^{-c|k_\gamma-k_\sigma|}
 \le e^{cR_*}e^{-c|k_\gamma-v|}.
\]
Summing over $\sigma$ therefore contributes at most $2(2R_*+1)e^{cR_*}$. Consequently
\[
 |r_{i,\alpha\gamma}|
 \le C\sum_{v\in J}v^{-1/2}
           e^{-c(|i-v|+|k_\alpha-v|+|k_\gamma-v|)}.
\]
The first contribution is included by the term $v=i$ in this sum.

Set $D_v=|i-v|+|k_\alpha-v|+|k_\gamma-v|$ and $D=|i-k_\gamma|+|k_\alpha-k_\gamma|$. The triangle inequality gives $D\le2D_v$, and positivity of the site indices gives
\[
 v^{-1/2}\le k_\gamma^{-1/2}(1+|v-k_\gamma|)^{1/2}.
\]
Splitting the exponential into two equal factors yields
\[
 \begin{split}
 \sum_{v\in J}v^{-1/2}e^{-cD_v}
 &\le k_\gamma^{-1/2}e^{-cD/4}
       \sum_{v\in J}(1+|v-k_\gamma|)^{1/2}
                                    e^{-c|v-k_\gamma|/2}\\
 &\le\frac{2k_\gamma^{-1/2}e^{-cD/4}}
                         {(1-e^{-c/2})^2}.
 \end{split}
\]
For the last inequality, extend the sum of distances to all integers and use $(1+h)^{1/2}\le1+h$ together with $\sum_{h\ge0}(1+h)q^h=(1-q)^{-2}$ for $0<q<1$. This bounds both
boundary point cases with the same constant and proves \eqref{B:eq:ratio-poisson-second-bound}. Continuity extends these bounds to every point of the fluctuation space.
\end{proof}

\subsubsection{The diffusion matrix}
\label{B:subsection:ratio-diffusion-matrix}

The Poisson gradients yield the matrix in the integration by parts identity. Its symmetric part must remain uniformly positive, and its mean must approximate the target covariance;
the two properties are verified together below.

Put $Y=f-\langle f\rangle$ and define
\begin{equation}
 \bar\sigma_{ij}(x)=\sum_\alpha g_{i,\alpha}(x)\partial_\alpha f_j(x),
 \qquad \bar\tau=\frac12(\bar\sigma+\bar\sigma^{\mathsf T}).
 \label{B:eq:ratio-diffusion-matrix-definition}
\end{equation}

\begin{lemma}\label{B:lem:local-diffusion-matrix}
After a fixed choice of $\delta>0$, these matrices satisfy
\begin{equation}
 \frac12B\le\bar\tau(x)\le\frac32B,\qquad
 |\partial_\alpha\bar\tau_{ij}(x)|
 \le Ck_\alpha^{-1/2}e^{-c(|i-k_\alpha|+|j-k_\alpha|)}.
 \label{B:eq:ratio-diffusion-matrix-bounds}
\end{equation}
Furthermore,
\begin{equation}
 \langle\bar\tau_{ij}\rangle=\operatorname{Cov}(f_i,f_j),\qquad
 |\langle\bar\tau_{ij}\rangle-B_{ij}|
                  \le C(ij)^{-1/2}e^{-c|i-j|}.
 \label{B:eq:ratio-diffusion-matrix-mean}
\end{equation}
\end{lemma}
\begin{proof}
Write $g_i^0=CL_i^{\mathsf T}$. Since $B_{ij}=L_jCL_i^{\mathsf T}$, subtraction gives
\[
 \bar\sigma_{ij}-B_{ij}
 =\sum_\alpha(g_{i,\alpha}-g_{i,\alpha}^0)\partial_\alpha f_j
       +\sum_\alpha g_{i,\alpha}^0(\partial_\alpha f_j-L_{j,\alpha}).
\]
Each sum has only the two coordinates of site $j$. The bounds $|\nabla f_j|\le C$, $|\nabla f_j-L_j^{\mathsf T}|\le C\delta$ and \eqref{B:eq:ratio-poisson-gradient-bounds} imply
$|\bar\sigma_{ij}-B_{ij}|\le C\delta e^{-c|i-j|}$. The transpose has the same bound, because $B$ is symmetric. Thus
\[
 \sup_i\sum_j|\bar\tau_{ij}-B_{ij}|
 +\sup_j\sum_i|\bar\tau_{ij}-B_{ij}|
 \le 2C\delta\sum_{h\in\mathbb Z}e^{-c|h|}
 =2C\delta\frac{1+e^{-c}}{1-e^{-c}}.
\]
Schur's test gives $\|\bar\tau-B\|\le C_{\rm pos}\delta$. Choose $\delta$ to meet the convergence conditions for both Neumann series in the preceding subsubsections and
$C_{\rm pos}\delta\le b_-/2$. Then, for every $v$,
\[
 |v^{\mathsf T}(\bar\tau-B)v|
 \le\frac{b_-}{2}|v|^2\le\frac12v^{\mathsf T}Bv,
 \qquad \frac12B\le\bar\tau\le\frac32B.
\]
All constants used to select $\delta$ are uniform in the admissible boundary data and the dimension. Derivative constants below may depend on this fixed choice.

Differentiation gives the four contributions
\begin{equation}
 2\partial_\gamma\bar\tau_{ij}
 =\sum_\alpha\bigl[
 r_{i,\alpha\gamma}\partial_\alpha f_j
 +g_{i,\alpha}\partial^2_{\alpha\gamma}f_j
 +r_{j,\alpha\gamma}\partial_\alpha f_i
 +g_{j,\alpha}\partial^2_{\alpha\gamma}f_i\bigr].
 \label{B:eq:ratio-diffusion-derivative-expansion}
\end{equation}
Put $v=k_\gamma$. The first and third contributions obey
\[
 \begin{split}
 \sum_{k_\alpha=j}|r_{i,\alpha\gamma}\partial_\alpha f_j|
 &\le Cv^{-1/2}e^{-c(|i-v|+|j-v|)},\\
 \sum_{k_\alpha=i}|r_{j,\alpha\gamma}\partial_\alpha f_i|
 &\le Cv^{-1/2}e^{-c(|j-v|+|i-v|)}.
 \end{split}
\]
The second contribution is zero unless $v=j$; on that event it is at most $Cj^{-1/2}e^{-c|i-j|}$. The fourth is zero unless $v=i$ and is at most $Ci^{-1/2}e^{-c|j-i|}$. Adding
these four bounds proves \eqref{B:eq:ratio-diffusion-matrix-bounds}, including $i=j$ and repeated coordinate indices.

To take the mean, choose $0\le\chi_N\le1$ equal to one on the ball of radius $N$, supported in the ball of radius $2N$, and satisfying $|\nabla\chi_N|\le C/N$. Test the weak
Poisson equation with $\chi_Nf_j$. The term where the derivative meets the cutoff satisfies
\[
 |\langle f_j\nabla h_i\cdot\nabla\chi_N\rangle|
 \le \frac C N\|f_j\mathbf1_{\{|x|\ge N\}}\|_2
                 \|\nabla h_i\|_2\longrightarrow0.
\]
For the term without a derivative on the cutoff, Cauchy--Schwarz gives
\[
 \left|\langle(1-\chi_N)(f_i-\langle f_i\rangle)f_j\rangle\right|
 \le\|(f_i-\langle f_i\rangle)\mathbf1_{\{|x|\ge N\}}\|_2\|f_j\|_2
 \longrightarrow0.
\]
Also, at this fixed dimension,
\[
 \left|\langle(1-\chi_N)\nabla h_i\cdot\nabla f_j\rangle\right|
 \le\|\nabla h_i\|_\infty\|\nabla f_j\|_\infty
                \widetilde\mu_q\{|x|\ge N\}\longrightarrow0.
\]
Together with the estimate for the term containing the gradient of the cutoff, these limits give
\[
 \langle(f_i-\langle f_i\rangle)f_j\rangle
 =\sum_\alpha\langle g_{i,\alpha}\partial_\alpha f_j\rangle
 =\langle\bar\sigma_{ij}\rangle.
\]
The left side is symmetric in $i,j$. Taking the symmetric part proves $\langle\bar\tau_{ij}\rangle=\operatorname{Cov}(f_i,f_j)$.
Proposition~\ref{B:prop:extended-logarithmic-covariance} supplies the remaining mean estimate.
\end{proof}

\begin{proposition}\label{B:prop:conditional-diffusion-matrix}
The law of $Y$ has a strictly positive smooth density. There is a smooth symmetric matrix $\tau$ with $B/2\le\tau\le3B/2$ such that
\begin{equation}
 \mathbb E[Y\cdot\nabla\psi(Y)]
   =\mathbb E[\tau(Y):\nabla^2\psi(Y)]\qquad
                         (\psi\in C_c^\infty(\mathbb R^d)).
 \label{B:eq:ratio-scalar-integration-by-parts}
\end{equation}
It can be chosen as $\tau(y)=\mathbb E[\bar\tau(x)\mid Y=y]$.
\end{proposition}
\begin{proof}
Let $\Phi\in C_c^\infty(\mathbb R^d;\mathbb R^d)$. The composition $\Phi_i(Y)$ is bounded, and its gradient is $\sum_j(\partial_j\Phi_i)(Y)\nabla f_j$. At fixed dimension this
gradient is bounded. Testing the Poisson equation with $\chi_N\Phi_i(Y)$ is thus justified. Its cutoff term is bounded by $CN^{-1}\|\Phi_i\|_\infty\|\nabla h_i\|_2$ and tends to
zero. For the left side, $Y_i\in L^1$ and $\mathbb E[|Y_i|\mathbf1_{\{|x|\ge N\}}]\to0$, while $|\Phi_i(Y)|\le\|\Phi_i\|_\infty$. For the gradient term, each
$g_{i,\alpha}(\partial_j\Phi_i)(Y)\partial_\alpha f_j$ is bounded at this fixed dimension. Its omitted integral is at most its supremum times $\widetilde\mu_q\{|x|\ge N\}\to0$.
Thus all terms without a derivative on the cutoff converge. Summing in $i$ gives
\begin{equation}
 \mathbb E\sum_iY_i\Phi_i(Y)
  =\mathbb E\sum_{i,j}\bar\sigma_{ij}(x)\partial_j\Phi_i(Y).
 \label{B:eq:ratio-vector-stein-identity}
\end{equation}
Conditional expectation therefore gives the identity for vector fields with $\sigma(y)=\mathbb E[\bar\sigma(x)\mid Y=y]$. In particular $\sigma$ is a Stein kernel. For
$\Phi=\nabla\psi$, the antisymmetric part contributes zero, since
\[
 \sum_{i,j}(\bar\sigma_{ij}-\bar\sigma_{ji})\partial_{ij}\psi(Y)
 =\sum_{i,j}\bar\sigma_{ij}
        [\partial_{ij}\psi(Y)-\partial_{ji}\psi(Y)]=0.
\]
This yields the scalar identity with the symmetric conditional matrix. The identity for vector fields continues to use
$\sigma$.

To obtain a smooth version of the conditional matrix, introduce $v_k=(l_{x,k}(x_{x,k})+l_{y,k}(x_{y,k}))/2$ and use $(f,v)$ as coordinates. The inverse is
\begin{equation}
 x_{x,k}=l_{x,k}^{-1}(v_k-f_k/(2\sqrt k)),\qquad
 x_{y,k}=l_{y,k}^{-1}(v_k+f_k/(2\sqrt k)).
 \label{B:eq:ratio-extended-inverse-coordinates}
\end{equation}
Indeed $f_k=\sqrt k(l_{y,k}-l_{x,k})$. The forward derivative at one site is
\[
 \begin{pmatrix}
 -\sqrt k\,l'_{x,k}&\sqrt k\,l'_{y,k}\\
 l'_{x,k}/2&l'_{y,k}/2
 \end{pmatrix},\qquad
 \det=-\sqrt k\,l'_{x,k}l'_{y,k}.
\]
Each $l_\alpha$ is strictly increasing, onto $\mathbb R$, and has derivative between two positive multiples of $d_\alpha^{-1}$. Its inverse is smooth, with
$(l_\alpha^{-1})'=1/(l_\alpha'\circ l_\alpha^{-1})$; differentiating this identity repeatedly bounds every derivative of the inverse of positive order at fixed $n,q,\delta$.
Outside bounded intervals, each map has the form $x\mapsto a_\alpha x+b_\alpha$ with $a_\alpha>0$. Thus the coordinate change is a global smooth diffeomorphism and its absolute
inverse Jacobian is
\[
 J_q(f,v)=\prod_{k\in J}
 \frac1{\sqrt k\,l'_{x,k}(x_{x,k}(f,v))l'_{y,k}(x_{y,k}(f,v))}>0.
\]
The bounds $|l_\alpha(t)|\le C|t|/d_\alpha$, with $l_\alpha(0)=0$, show that $|x(f,v)|^2\ge c_{n,q,\delta}|v|^2$. Conversely $|x(f,v)|\le C_{n,q,\delta}(1+|f|+|v|)$. Uniform
convexity of $V_0$ gives, on each compact set of $f$, $V_0(x(f,v))\ge c_{n,q,\delta}|v|^2-C_{n,q,\delta}$. Write $F(f,v)=V_0(x(f,v))$. On each compact set of $f$, the inverse
coordinate map has bounded derivatives of positive order and $|x(f,v)|\le C(1+|v|)$. Since $\nabla V_0$ has at most linear growth and all its positive order derivatives are
bounded, the chain rule gives
\[
 |\partial_f^\alpha F(f,v)|\le C_\alpha(1+|v|)\quad(|\alpha|\ge1).
\]
The positive lower bounds on $l_\alpha'$, and the bounded derivatives of $l_\alpha$, its inverse, and the Poisson gradients, give bounded $f$-derivatives of $J_q$ and
$\bar\tau(x(f,v))J_q$. The product rule and $F(f,v)\ge c|v|^2-C$ therefore imply
\[
 |\partial_f^\alpha(e^{-F}J_q)|+
 |\partial_f^\alpha(\bar\tau(x(f,v))e^{-F}J_q)|
 \le C_\alpha(1+|v|)^{|\alpha|}e^{-c|v|^2}.
\]
This bound is integrable in $v$ and uniform on the chosen compact set of $f$, including $|\alpha|=0$. It justifies differentiation under both integrals below at every fixed order.

Define
\[
 D_q(f)=\int e^{-V_0(x(f,v))}J_q(f,v)\,\dd v,\qquad
 N_q(f)=\int\bar\tau(x(f,v))e^{-V_0(x(f,v))}J_q(f,v)\,\dd v.
\]
The first integral is finite and positive, and both functions are smooth by the displayed differentiated bounds. With $Z_0=\int e^{-V_0(x)}\,\dd x$, change of variables and Tonelli
give $\int D_q(f)\,\dd f=Z_0$. Thus $Y$ has the positive smooth density $D_q(y+\langle f\rangle)/Z_0$. For any bounded Borel scalar $\varphi$, the same change of variables yields
\[
 \mathbb E[\bar\tau(x)\varphi(Y)]
 =Z_0^{-1}\int N_q(f)\varphi(f-\langle f\rangle)\,\dd f
 =\mathbb E\!\left[
 \frac{N_q(Y+\langle f\rangle)}{D_q(Y+\langle f\rangle)}\varphi(Y)\right].
\]
This verifies the conditional version $\tau(y)=N_q(y+\langle f\rangle)/D_q(y+\langle f\rangle)$. For every vector $v$, integration of the pointwise inequalities for
$v^{\mathsf T}\bar\tau v$ and division by $D_q(f)>0$ give $\tfrac12v^{\mathsf T}Bv\le v^{\mathsf T}\tau(y)v\le\tfrac32v^{\mathsf T}Bv$ for every $y$. The scalar identity already
proved now holds for this smooth version. In addition, the Lipschitz map $f$ and the Gaussian tails of $\widetilde\mu_q$ give $\mathbb E|Y|^p<\infty$ for every finite $p$.
\end{proof}

\subsubsection{A discrete heat kernel estimate}
\label{B:subsection:partial-sum-heat-kernel}

The inverse of the matrix defining the partial sum cost has a discrete difference structure. Its heat kernel estimate will control the matrix norm used in the coupling.

Let $D_J=\operatorname{diag}_{k\in J}(k^{-1/2})$ and set
\begin{equation}
 T=\mathsf T_2D_J,
 \qquad \mathsf M=(T^{\mathsf T}T)^{-1}.
 \label{B:eq:partial-sum-matrix}
\end{equation}
The inverse $D=S_n^{-1}$ satisfies $(Dv)_k=v_k-v_{k+1}$ for $k<R$ and $(Dv)_R=v_R$. Therefore
\begin{equation}
 \mathsf M=n^2\operatorname{diag}(\sqrt k)D
               \operatorname{diag}(1/k)D^{\mathsf T}
               \operatorname{diag}(\sqrt k).
 \label{B:eq:partial-sum-inverse-formula}
\end{equation}

\begin{lemma}\label{B:lem:partial-sum-heat-kernel}
The entries of $P_t^{\mathsf M}=e^{-t\mathsf M}$ are nonnegative, and its row sums are at most one. For constants uniform in $n$,
\begin{equation}
 \begin{split}
 P_t^{\mathsf M}(i,i)&\le\min\{1,C/(n\sqrt t)\}\quad(t>0),\\
 P_t^{\mathsf M}(i,i)&\le(C/n)e^{-c(t-1)}\quad(t\ge1).
 \end{split}
 \label{B:eq:partial-sum-heat-bound}
\end{equation}
\end{lemma}
\begin{proof}
Direct multiplication in \eqref{B:eq:partial-sum-inverse-formula} gives
\[
 n^{-2}\mathsf M_{kk}=1+\frac{k}{k+1}\ (k<R),\qquad
 n^{-2}\mathsf M_{RR}=1,\qquad
 n^{-2}\mathsf M_{k,k+1}=-a_k=-\sqrt{\frac{k}{k+1}}.
\]
All other off-diagonal entries vanish. For an interior site put $x=k/(k+1)$ and $y=(k-1)/k$. Its row sum divided by $n^2$ is
\[
 s_k=1+x-\sqrt x-\sqrt y
 \ge1+x-\frac{1+x}{2}-\frac{1+y}{2}
 =\frac1{2k(k+1)}.
\]
At the left boundary $s_r=1+x-\sqrt x=(\sqrt x-1/2)^2+3/4$; at the right boundary $s_R=1-\sqrt{(R-1)/R}\ge0$. For a one site interval the single entry is $1$. Expanding the edge
squares proves
\begin{equation}
 n^{-2}v^{\mathsf T}\mathsf Mv
 =\sum_{k=r}^{R-1}a_k(v_{k+1}-v_k)^2+\sum_{k=r}^Rs_kv_k^2
 \ge c\left(v_r^2+\sum_{k=r}^{R-1}|v_{k+1}-v_k|^2\right).
 \label{B:eq:partial-sum-quadratic-form}
\end{equation}
Here $a_k\ge1/\sqrt2$, and the left boundary term has coefficient at least $3/4$; these lower bounds also give a fixed $c$ in the one site case.

Choose $b\ge\max_k\mathsf M_{kk}$. The finite matrix $bI-\mathsf M$ has nonnegative entries, so its exponential series proves $P_t^{\mathsf M}=e^{-bt}e^{t(bI-\mathsf M)}\ge0$
entrywise. Moreover
\[
 \mathbf1-P_t^{\mathsf M}\mathbf1
 =\int_0^tP_s^{\mathsf M}\mathsf M\mathbf1\,\dd s\ge0.
\]
Its rows sum to at most one. Symmetry gives the same column bound, whence $\|P_tv\|_\infty\le\|v\|_\infty$ and
$\sum_i|(P_tv)_i|\le\sum_j|v_j|\sum_iP_t(i,j)\le\|v\|_1$.

Writing $v_k=v_r+\sum_{j=r}^{k-1}(v_{j+1}-v_j)$ and using Cauchy--Schwarz gives
\[
 \|v\|_2^2\le\frac{d(d+1)}2
 \left(v_r^2+\sum_{j=r}^{R-1}|v_{j+1}-v_j|^2\right)
 \le d^2\left(v_r^2+\sum_j|v_{j+1}-v_j|^2\right),\qquad d=|J|.
\]
Since $d\le n$, \eqref{B:eq:partial-sum-quadratic-form} yields $\mathsf M\ge cI$ uniformly.

We derive the Nash inequality using an extension that preserves this boundary condition. For a finitely supported sequence $w$ on $\mathbb Z$, put $E(w)=\sum_j|w_{j+1}-w_j|^2$.
Parseval, the bound $|\widehat w|\le\|w\|_1$, and $4\sin^2(\vartheta/2)\ge4b^2/\pi^2$ on $b\le|\vartheta|\le\pi$ give, for $0<b\le1$,
\[
 \|w\|_2^2
 \le\frac b\pi\|w\|_1^2+\frac{\pi^2}{4b^2}E(w).
\]
If $w=0$ the desired inequality holds. Otherwise $E(w)>0$. When $E(w)\le\|w\|_1^2$, take $b=(E(w)/\|w\|_1^2)^{1/3}$ and cube the bound. When $E(w)>\|w\|_1^2$, use
$\|w\|_2^6\le\|w\|_1^6<E(w)\|w\|_1^4$. In both cases $\|w\|_2^6\le CE(w)\|w\|_1^4$. Extend $v$ by
\[
 w_k=v_k\ (r\le k\le R),\qquad
 w_k=v_{2R+1-k}\ (R+1\le k\le2R+1-r),
\]
and set $w_k=0$ elsewhere. Each coordinate occurs twice and the central edge has equal boundary values. Consequently
\[
 \|w\|_2^2=2\|v\|_2^2,\quad \|w\|_1=2\|v\|_1,\quad
 E(w)=2v_r^2+2\sum_{k=r}^{R-1}|v_{k+1}-v_k|^2.
\]
Substitution of these three identities and \eqref{B:eq:partial-sum-quadratic-form} proves
\begin{equation}
 \|v\|_2^6\le Cn^{-2}(v^{\mathsf T}\mathsf Mv)\|v\|_1^4.
 \label{B:eq:partial-sum-nash-inequality}
\end{equation}
For $v(t)=P_t^{\mathsf M}e_i$, put $y(t)=\|v(t)\|_2^2$. The matrix exponential is invertible, so $y(t)>0$. The derivative is $y'(t)=-2v(t)^{\mathsf T}\mathsf Mv(t)$. The $\ell^1$
contraction and the Nash inequality give $y'\le-cn^2y^3$. Hence
\[
 \frac{\dd}{\dd t}y(t)^{-2}=-2y(t)^{-3}y'(t)\ge2cn^2,
 \qquad y(t)\le(1+2cn^2t)^{-1/2}.
\]
The semigroup identity and symmetry give $P_{2t}^{\mathsf M}(i,i)=\|P_t^{\mathsf M}e_i\|_2^2=y(t)$. Replacing $2t$ by $t$ proves both $P_t(i,i)\le1$ and $P_t(i,i)\le C/(n\sqrt t)$.
For $t\ge1$ the lower spectral bound gives
\[
 P_t^{\mathsf M}(i,i)
 =\langle P_{1/2}^{\mathsf M}e_i,
          P_{t-1}^{\mathsf M}P_{1/2}^{\mathsf M}e_i\rangle
 \le e^{-c(t-1)}P_1^{\mathsf M}(i,i)
 \le\frac Cn e^{-c(t-1)}.
\]
This proves both semigroup bounds with constants independent of $r,R,n$.
\end{proof}

\subsubsection{The matrix Lyapunov equation}
\label{B:subsection:ratio-lyapunov-equation}

The heat semigroup represents the solution of the matrix Lyapunov equation as an integral. The resulting quadratic form controls the difference between square roots of covariance
matrices in terms of the error in the diffusion matrix.

On real $d\times d$ matrices use $\langle E,F\rangle_F=\operatorname{Tr}(E^{\mathsf T}F)$ and define
\begin{equation}
 \mathcal R_{\mathsf M}(E)=\int_0^\infty
                 e^{-t\mathsf M}Ee^{-t\mathsf M}\,\dd t,
 \qquad \|E\|_{\mathscr H_{\mathsf M}}^2
                  =\langle E,\mathcal R_{\mathsf M}(E)\rangle_F.
 \label{B:eq:ratio-lyapunov-norm}
\end{equation}
The lower spectral bound gives $\|e^{-t\mathsf M}Ee^{-t\mathsf M}\|_F\le e^{-2ct}\|E\|_F$, so the integral converges in Frobenius norm. Its finite time version obeys
\[
 \mathsf M\int_0^T e^{-t\mathsf M}Ee^{-t\mathsf M}\,\dd t
 +\left(\int_0^T e^{-t\mathsf M}Ee^{-t\mathsf M}\,\dd t\right)\mathsf M
 =E-e^{-T\mathsf M}Ee^{-T\mathsf M}.
\]
The last term tends to zero. In an orthonormal eigenbasis of $\mathsf M$, $\mathcal R_{\mathsf M}(E)_{ij}=E_{ij}/(m_i+m_j)$, and therefore
\[
 \langle E,\mathcal R_{\mathsf M}(F)\rangle_F
 =\sum_{i,j}\frac{E_{ij}F_{ij}}{m_i+m_j},\qquad
 \|E\|_{\mathscr H_{\mathsf M}}^2
 =\sum_{i,j}\frac{|E_{ij}|^2}{m_i+m_j}.
\]
Every denominator is positive. These formulas prove self-adjointness, strict positivity, and uniqueness of the inverse of $E\mapsto\mathsf ME+E\mathsf M$ on all real matrices.

\begin{lemma}\label{B:lem:partial-sum-lyapunov-bound}
In the basis of matrix units,
\begin{equation}
 0\le(\mathcal R_{\mathsf M})_{ij,kl}\le C\ell_n/n^2,
 \qquad
 \|E\|_{\mathscr H_{\mathsf M}}^2
          \le\frac{C\ell_n}{n^2}\left(\sum_{i,j}|E_{ij}|\right)^2.
 \label{B:eq:ratio-lyapunov-entry-bound}
\end{equation}
\end{lemma}
\begin{proof}
For the matrix unit $E^{kl}$, direct multiplication gives $(e^{-t\mathsf M}E^{kl}e^{-t\mathsf M})_{ij}=P_t^{\mathsf M}(i,k)P_t^{\mathsf M}(j,l)$. Thus
\[
 (\mathcal R_{\mathsf M})_{ij,kl}
 =\int_0^\infty P_t^{\mathsf M}(i,k)P_t^{\mathsf M}(j,l)\,\dd t\ge0.
\]
Its diagonal entries can be bounded by splitting the time integral at $n^{-2}$ and $1$:
\[
 \begin{split}
 (\mathcal R_{\mathsf M})_{ij,ij}
 &\le\int_0^{n^{-2}}1\,\dd t
 +\frac{C}{n^2}\int_{n^{-2}}^1\frac{\dd t}{t}
 +\frac{C}{n^2}\int_1^\infty e^{-2c(t-1)}\,\dd t\\
 &=\frac1{n^2}+\frac{2C\log n}{n^2}+\frac{C}{2cn^2}
 \le\frac{C\ell_n}{n^2}.
 \end{split}
\]
When $n=1$ the middle integral is empty. The positivity as an operator on Frobenius space, proved above, also gives Cauchy--Schwarz:
\[
 |\langle E^{ij},\mathcal R_{\mathsf M}(E^{kl})\rangle_F|^2
 \le\langle E^{ij},\mathcal R_{\mathsf M}(E^{ij})\rangle_F
       \langle E^{kl},\mathcal R_{\mathsf M}(E^{kl})\rangle_F.
\]
Together with entrywise nonnegativity this proves the first bound. For an arbitrary real matrix $E$, expansion in matrix units gives
\[
 \|E\|_{\mathscr H_{\mathsf M}}^2
 =\sum_{i,j,k,l}E_{ij}E_{kl}(\mathcal R_{\mathsf M})_{ij,kl}
 \le\frac{C\ell_n}{n^2}\sum_{i,j,k,l}|E_{ij}||E_{kl}|
 =\frac{C\ell_n}{n^2}\left(\sum_{i,j}|E_{ij}|\right)^2.
\]
All cross terms have been included.
\end{proof}

\begin{proposition}\label{B:prop:diffusion-matrix-deviation}
For $p\ge1$, the diffusion matrix of Proposition~\ref{B:prop:conditional-diffusion-matrix} satisfies
\begin{equation}
 \left\|\|\tau(Y)-B\|_{\mathscr H_{\mathsf M}}\right\|_{L^{2p}}
       \le \frac Cn\left(\sqrt p\,\ell_n+\ell_n^{3/2}\right).
 \label{B:eq:ratio-diffusion-deviation}
\end{equation}
\end{proposition}
\begin{proof}
For a fixed coordinate $\alpha$, sum its derivative bound over the two observation indices to obtain
\[
 \sum_{i,j}|\partial_\alpha\bar\tau_{ij}|
 \le Ck_\alpha^{-1/2}
       \left(\sum_{h\in\mathbb Z}e^{-c|h|}\right)^2
 \le Ck_\alpha^{-1/2}.
\]
The Lyapunov entry bound then yields
\begin{equation}
 \|\partial_\alpha\bar\tau\|_{\mathscr H_{\mathsf M}}^2
     \le\frac{C\ell_n}{n^2k_\alpha},\qquad
 \sup_x\sum_\alpha\|\partial_\alpha\bar\tau(x)\|_{\mathscr H_{\mathsf M}}^2
     \le C\ell_n^2/n^2=:L_n^2.
 \label{B:eq:ratio-diffusion-lipschitz-bound}
\end{equation}
For its second inequality we used both species and $\sum_{\alpha\in\mathscr S(J)}k_\alpha^{-1}=2\sum_{k=r}^Rk^{-1} \le2(1+\log(R/r))\le C\ell_n$.

Let $(E^a)_{a=1}^{d^2}$ be an orthonormal basis for $\mathscr H_{\mathsf M}$, and put $F_a(x)=\langle\bar\tau(x),E^a\rangle_{\mathscr H_{\mathsf M}}$. Each $F_a$ and its gradient
are bounded at fixed dimension. Poincar\'e followed by Parseval gives
\[
 \begin{split}
 \mathbb E\|\bar\tau-\langle\bar\tau\rangle\|_{\mathscr H_{\mathsf M}}^2
 &=\sum_a\operatorname{Var}(F_a)
 \le C\sum_{a,\alpha}\mathbb E|\partial_\alpha F_a|^2\\
 &=C\mathbb E\sum_\alpha
       \|\partial_\alpha\bar\tau\|_{\mathscr H_{\mathsf M}}^2
 \le CL_n^2.
 \end{split}
\]
This summation introduces no additional power of $d$. For $F(x)=\|\bar\tau(x)-\langle\bar\tau\rangle\|_{\mathscr H_{\mathsf M}}$, the reverse triangle inequality and integration
along a line give
\[
 |F(x)-F(y)|
 \le\int_0^1\left\|\sum_\alpha(x_\alpha-y_\alpha)
 \partial_\alpha\bar\tau(y+t(x-y))\right\|_{\mathscr H_{\mathsf M}}\dd t
 \le L_n|x-y|.
\]
The concentration consequence of the logarithmic Sobolev inequality in Section~\ref{B:sec:conditional-estimates-on-the-interior-interval} applies to this Lipschitz function. At a
zero of the norm one may first use $(\|\bar\tau-\langle\bar\tau\rangle\|_{\mathscr H_{\mathsf M}}^2+ \epsilon^2)^{1/2}$, whose Lipschitz constant is also $L_n$, and let
$\epsilon\downarrow0$. Hence $\Pr(|F-\mathbb EF|>u)\le2e^{-cu^2/L_n^2}$. For real $p\ge1$, Tonelli and the substitution $v=cu^2/L_n^2$ give
\[
 \mathbb E|F-\mathbb EF|^{2p}
 \le4p\int_0^\infty u^{2p-1}e^{-cu^2/L_n^2}\,\dd u
 =2\Gamma(p+1)(L_n^2/c)^p.
\]
To bound the gamma factor directly, use $v^pe^{-v/2}\le(2p/e)^p$ and integrate the remaining $e^{-v/2}$. This gives $\Gamma(p+1)\le2(2p/e)^p$ and
$\|F-\mathbb EF\|_{2p}\le C\sqrt p L_n$. Together with $\mathbb EF\le(\mathbb EF^2)^{1/2}\le CL_n$, this proves $\|F\|_{2p}\le C\sqrt p\ell_n/n$.

For the mean, the inequality $i/j\le1+|i-j|$ gives
\[
 \sum_{i,j}|\langle\bar\tau_{ij}\rangle-B_{ij}|
 \le C\sum_i i^{-1}\sum_j(1+|i-j|)^{1/2}e^{-c|i-j|}
 \le C\ell_n.
\]
The last geometric sum is uniformly finite: replace its polynomial by $1+|i-j|$ and sum the geometric series and its derivative. Thus the Lyapunov entry estimate gives
$\|\langle\bar\tau\rangle-B\|_{\mathscr H_{\mathsf M}} \le C\ell_n^{3/2}/n$. Finally the conditional version proved above implies
\[
 \mathbb E\|\tau(Y)-B\|_{\mathscr H_{\mathsf M}}^{2p}
 =\mathbb E\|\mathbb E[\bar\tau-B\mid Y]\|_{\mathscr H_{\mathsf M}}^{2p}
 \le\mathbb E\|\bar\tau-B\|_{\mathscr H_{\mathsf M}}^{2p}.
\]
Conditional Jensen applies because the norm to the power $2p$ is convex. The triangle inequality with the fluctuation and mean estimates proves the stated bound for every $p\ge1$.
\end{proof}

\subsubsection{Factorization of the covariance matrix}
\label{B:subsection:gaussian-covariance-factors}

Two centered Gaussian vectors can be coupled using square roots of their covariance matrices and an orthogonal matrix. Minimizing over this matrix gives the following bound in
terms of the inverse of the Lyapunov operator.

For positive definite real matrices $P,Q$, put
\[
 d_{\rm B}(P,Q)^2=\min_{O^{\mathsf T}O=I}
                  \|P^{1/2}-Q^{1/2}O\|_F^2.
\]
This is the squared Bures--Wasserstein distance \cite{bhatiaJainLim2019bures}. For real positive definite matrices, the polar factor in that result is orthogonal. For a positive
matrix $Q$, write $\mathcal L_QX=QX+XQ$.

\begin{lemma}\label{B:lem:gaussian-factor-distance}
If $P\ge Q/2$ and $E=P-Q$, then
\begin{equation}
 d_{\rm B}(P,Q)^2\le\langle E,\mathcal L_Q^{-1}E\rangle_F.
 \label{B:eq:gaussian-factor-lyapunov-bound}
\end{equation}
A minimizing orthogonal matrix is $O=F(F^{\mathsf T}F)^{-1/2}$, where $F=Q^{1/2}P^{1/2}$.
\end{lemma}
\begin{proof}
First identify the minimizing orthogonal matrix. Put $F=Q^{1/2}P^{1/2}$ and take a singular value decomposition $F=UDV^{\mathsf T}$, with $D=\operatorname{diag}(d_1,\ldots,d_d)$
and $d_i>0$. For every orthogonal $O$,
\[
 \|P^{1/2}-Q^{1/2}O\|_F^2
 =\operatorname{Tr}P+\operatorname{Tr}Q-2\operatorname{Tr}(F^{\mathsf T}O).
\]
Writing $Z=U^{\mathsf T}OV$, which is orthogonal, gives $\operatorname{Tr}(F^{\mathsf T}O)=\sum_i d_iZ_{ii}\le\sum_i d_i$, since $Z_{ii}\le1$. Equality holds for
$O=UV^{\mathsf T}=F(F^{\mathsf T}F)^{-1/2}$. This proves the formula for the minimizer without a choice of eigenvector signs.

For the bound, set $A_s=Q+sE$ for $0\le s\le1$, where $E=P-Q$. These matrices are positive definite and satisfy $A_s\ge Q/2$. Define the symmetric matrix
\[
 X_s=\mathcal L_{A_s}^{-1}E,\qquad A_sX_s+X_sA_s=E.
\]
In a diagonal basis for $A_s$, the inverse divides the $(i,j)$ entry by the positive number $a_i(s)+a_j(s)$. Thus $X_s$ exists, is symmetric, and depends continuously on $s$. Solve
the linear matrix equation $\dot F_s=X_sF_s$, $F_0=Q^{1/2}$. The product $C_s=F_sF_s^{\mathsf T}$ obeys $\dot C_s=X_sC_s+C_sX_s$. The matrix $A_s$ obeys the same equation because
$\dot A_s=E=X_sA_s+A_sX_s$, and $C_0=A_0=Q$. Uniqueness for this finite-dimensional linear equation gives $F_sF_s^{\mathsf T}=A_s$. Consequently $F_1=P^{1/2}O_1$ for an orthogonal
$O_1$. Right multiplication by $O_1^{\mathsf T}$ preserves the Frobenius norm, so
\[
 d_{\rm B}(P,Q)\le\|F_1-Q^{1/2}\|_F
 \le\int_0^1\|\dot F_s\|_F\,\dd s.
\]
The squared speed is computed directly:
\[
 \|\dot F_s\|_F^2
 =\operatorname{Tr}(X_sA_sX_s)
 =\tfrac12\operatorname{Tr}\bigl(X_s(A_sX_s+X_sA_s)\bigr)
 =\tfrac12\langle E,\mathcal L_{A_s}^{-1}E\rangle_F.
\]
For an arbitrary real matrix $Z$, positivity of $A_s-Q/2$ implies
\[
 \langle Z,(\mathcal L_{A_s}-\tfrac12\mathcal L_Q)Z\rangle_F
 =\operatorname{Tr}(Z^{\mathsf T}(A_s-Q/2)Z)
  +\operatorname{Tr}(Z^{\mathsf T}Z(A_s-Q/2))\ge0.
\]
Inversion reverses this inequality: conjugating by $\mathcal L_Q^{-1/2}$ reduces it to the scalar bound $t^{-1}\le2$ for the eigenvalues $t\ge1/2$ of a positive operator. Hence
$\mathcal L_{A_s}^{-1}\le2\mathcal L_Q^{-1}$. Cauchy--Schwarz in $s$ now yields
\[
 d_{\rm B}(P,Q)^2
 \le\int_0^1\|\dot F_s\|_F^2\,\dd s
 \le\int_0^1\langle E,\mathcal L_Q^{-1}E\rangle_F\,\dd s
 =\langle E,\mathcal L_Q^{-1}E\rangle_F.
\]
This proves the estimate with the stated constant.
\end{proof}

The following identity will be used after a change of coordinates. For any invertible $T$ and real matrix $E$,
\begin{equation}
 \langle TET^{\mathsf T},\mathcal L_{TT^{\mathsf T}}^{-1}
                         (TET^{\mathsf T})\rangle_F
 =\langle E,\mathcal L_{(T^{\mathsf T}T)^{-1}}^{-1}E\rangle_F.
 \label{B:eq:lyapunov-congruence-identity}
\end{equation}
Indeed, let $T=U\operatorname{diag}(s_i)V^{\mathsf T}$, where every $s_i>0$, and write $E'=V^{\mathsf T}EV$. Then $TT^{\mathsf T}=U\operatorname{diag}(s_i^2)U^{\mathsf T}$, while
$U^{\mathsf T}TET^{\mathsf T}U$ has entry $s_is_jE'_{ij}$. The Lyapunov inverse in this basis divides that entry by $s_i^2+s_j^2$. Its Frobenius pairing is therefore
\[
 \sum_{i,j}\frac{s_i^2s_j^2|E'_{ij}|^2}{s_i^2+s_j^2}
 =\sum_{i,j}\frac{|E'_{ij}|^2}{s_i^{-2}+s_j^{-2}}.
\]
Since $(T^{\mathsf T}T)^{-1}=V\operatorname{diag}(s_i^{-2})V^{\mathsf T}$, the last sum is the right side of \eqref{B:eq:lyapunov-congruence-identity}. This calculation holds for
every real $E$, including nonsymmetric matrices.

\subsubsection{Invariant diffusions}
\label{B:subsection:invariance-of-ratio-diffusion}

We prove invariance of the distribution used in the coupling. Let $a$ be a smooth symmetric matrix field on $\mathbb R^d$ with $cI\le a(x)\le CI$, and let $\nu$ be a probability
measure with finite moments of every order. Suppose that
\begin{equation}
 \int[-x\cdot\nabla\psi+a:\nabla^2\psi]\,\dd\nu=0
                       \qquad(\psi\in C_c^\infty(\mathbb R^d)).
 \label{B:eq:invariant-diffusion-test-identity}
\end{equation}

\begin{lemma}\label{B:lem:invariant-diffusion}
The measure $\nu$ is invariant for
\begin{equation}
 \dd X_t=-X_t\,\dd t+\sqrt{2a(X_t)}\,\dd B_t.
 \label{B:eq:stationary-comparison-diffusion}
\end{equation}
\end{lemma}
\begin{proof}
Write $\sigma(x)=\sqrt{2a(x)}$ and $\mathcal A=a:\nabla^2-x\cdot\nabla$. Smoothness of the positive square root on the positive cone makes $\sigma$ smooth and locally Lipschitz.
Localize $-x$ and $\sigma$ with a smooth cutoff on a ball of radius $R+1$, equal to one on the ball of radius $R$. The resulting coefficients are bounded and globally Lipschitz.
Picard iteration satisfies
\[
 \mathbb E\sup_{s\le t}|X_s^{(m+1)}-X_s^{(m)}|^2
 \le K_{R,T}\frac{(C_{R,T}t)^m}{m!}\quad(0\le t\le T),
\]
by the It\^o isometry, the $L^2$ martingale maximal inequality, and iteration of the time integral. The square roots are summable, so the iterates converge in
$L^2(\Omega;C([0,T];\mathbb R^d))$. The same estimates and Gronwall give pathwise uniqueness and compatibility of the solutions stopped on nested balls. Thus a unique solution
exists up to its explosion time.

For $\tau_R=\inf\{t:|X_t|\ge R\}$, It\^o's formula gives
\[
 \mathbb E|X_{t\wedge\tau_R}|^2
 =\mathbb E|X_0|^2+
 \mathbb E\int_0^{t\wedge\tau_R}
        (2\operatorname{Tr}a(X_s)-2|X_s|^2)\,\dd s
 \le\mathbb E|X_0|^2+2dCt.
\]
At exit the norm is at least $R$, including paths started outside the ball. Thus $\Pr(\tau_R\le t)\le R^{-2}(\mathbb E|X_0|^2+2dCt)\to0$, which proves nonexplosion. The local
diffusion martingale problem gives a continuous local martingale $M_t=X_t-X_0+\int_0^tX_s\,\dd s$ with $\dd\langle M^i,M^j\rangle_t=2a_{ij}(X_t)\,\dd t$. These identities are
obtained using the coordinate functions and their products, multiplied by compactly supported cutoffs before $\tau_R$. Since $a$ is uniformly positive,
\[
 B_t=\int_0^t(2a(X_s))^{-1/2}\,\dd M_s,
 \qquad \langle B^i,B^j\rangle_t=t\delta_{ij}.
\]
L\'evy's characterization identifies the driving Brownian motion. Localized pathwise uniqueness and the preceding nonexplosion estimate therefore give uniqueness in law for the
martingale problem, for every initial point.

We verify the domain and semigroup assumptions of the criterion for invariant measures in Section~1 of \cite{echeverria1982criterion} (p.~2). For two initial points in a fixed
ball, the stopped difference equation, the local Lipschitz constant, and the maximal inequality give
\[
 \mathbb E\sup_{s\le t\wedge\tau_R^x\wedge\tau_R^y}
 |X_s^x-X_s^y|^2\le C_{R,t}|x-y|^2.
\]
The exit estimate is uniform for initial points in a compact set. First taking $y\to x$ and then $R\to\infty$ proves continuity of $P_tf$ for bounded continuous $f$. The mild
equation and It\^o isometry give
\[
 X_t^x=e^{-t}x+\int_0^te^{-(t-s)}\sigma(X_s^x)\,\dd B_s,
 \qquad
 \mathbb E|X_t^x-e^{-t}x|^2\le dC(1-e^{-2t}).
\]
For $|x|>e^tK$,
\[
 \Pr(|X_t^x|\le K)
 \le\frac{dC(1-e^{-2t})}{(e^{-t}|x|-K)^2}.
\]
Hence $P_t$ preserves $C_0(\mathbb R^d)$, first for compactly supported continuous functions and then by uniform approximation. For $|x|\le L$,
$\mathbb E|X_t^x-x|^2\le2(1-e^{-t})^2L^2+2dC(1-e^{-2t})$. If $f$ has compact support, uniform continuity and Markov's inequality give
\[
 \sup_{|x|\le L}|P_tf(x)-f(x)|
 \le\omega_f(\epsilon)+2\|f\|_\infty\epsilon^{-2}
 [2(1-e^{-t})^2L^2+2dC(1-e^{-2t})].
\]
Outside a fixed larger ball the preceding probability of reaching the support is $O(t)$ uniformly for $0<t\le1$. Letting $t\downarrow0$ and then $\epsilon\downarrow0$ proves strong
continuity on $C_0$.

Use the unital space $C_0(\mathbb R^d)+\mathbb R$, the continuous functions with a finite limit at infinity. The conservative semigroup extends by $P_t(f+c)=P_tf+c$ and is strongly
continuous on this space. For $f\in C_c^2$, $\mathcal Af\in C_0$, and the stopped It\^o formula, followed by bounded convergence, gives
\[
 P_tf-f=\int_0^tP_s\mathcal Af\,\dd s,
 \qquad
 \left\|\frac{P_tf-f}{t}-\mathcal Af\right\|_\infty
 \le\sup_{s\le t}\|P_s\mathcal Af-\mathcal Af\|_\infty\longrightarrow0.
\]
Thus $\mathscr M=C_c^2+\mathbb R$ lies in the generator domain on the unital space; constants have generator zero. It is uniformly dense there, by cutoff and convolution for the
$C_0$ part. If $f_j=c_j+g_j$ with $g_j\in C_c^2$ and $\Theta$ is smooth and convex, then $\Theta(f_1,\ldots,f_p)-\Theta(c_1,\ldots,c_p)$ is $C_c^2$. The chain rule gives
\[
 \mathcal A\Theta(f_1,\ldots,f_p)
 -\sum_j\Theta_j(f_1,\ldots,f_p)\mathcal Af_j
 =\sum_{i,j}\Theta_{ij}(f_1,\ldots,f_p)
                  \nabla f_i^{\mathsf T}a\nabla f_j\ge0.
\]
The last inequality is the trace pairing of two positive semidefinite matrices. At a minimum of $f\in\mathscr M$ the gradient is zero and the Hessian is positive semidefinite, so
$\mathcal Af\ge0$ there.

For $f\in C_c^2$, mollification gives $f_\epsilon\in C_c^\infty$ on a fixed larger compact set with $\|f_\epsilon-f\|_{C^2}\to0$. Boundedness of $a$ and of $x$ on this set implies
$\|\mathcal Af_\epsilon-\mathcal Af\|_\infty\to0$. The assumed integral identity therefore extends to every $f\in\mathscr M$. The criterion now applies: the state space is locally
compact and separable, the class of functions is dense and satisfies the minimum and convex composition conditions, its martingale problem is well posed, and
$\int\mathcal Af\,\dd\nu=0$ on that class. It gives $\int P_tf\,\dd\nu=\int f\,\dd\nu$, proving invariance. A diffusion defined by a positive Stein kernel is also treated in
\cite{fathiMikulincer2020stability}.
\end{proof}

\subsubsection{A moment estimate for the coupling}
\label{B:subsection:diffusion-coupling-moment-estimate}

We apply the preceding transport estimate after the invertible change of variables $T$. This gives one coupling for which the weighted partial sums of the errors satisfy the stated
bounds for all moment orders.

Let $Y,B,\tau$ be as in Proposition~\ref{B:prop:conditional-diffusion-matrix}. More generally, the following proof uses a centered vector with a smooth positive density and finite
moments of every order, together with a smooth $\tau$ satisfying \eqref{B:eq:ratio-scalar-integration-by-parts} and $B/2\le\tau\le3B/2$. For an arbitrary invertible $T$, put
$\mathsf M=(T^{\mathsf T}T)^{-1}$ and use the norm \eqref{B:eq:ratio-lyapunov-norm}.

\begin{proposition}\label{B:prop:diffusion-coupling-moments}
There is one coupling of $Y$ and $Z\sim N(0,B)$ such that, for every real $p\ge1$,
\begin{equation}
 \|T(Y-Z)\|_{L^{2p}}^2
 \le (2p-1)\lambda_{\min}(B)^{-1}
       \left\|\|\tau(Y)-B\|_{\mathscr H_{\mathsf M}}^2\right\|_{L^p}.
 \label{B:eq:diffusion-coupling-moment-bound}
\end{equation}
\end{proposition}
\begin{proof}
Put $\nu=\mathcal L(TY)$, $C_1=TBT^{\mathsf T}$, and $a(x)=T\tau(T^{-1}x)T^{\mathsf T}$. For $\psi(y)=\varphi(Ty)$, $y\cdot\nabla\psi(y)=(Ty)\cdot\nabla\varphi(Ty)$ and
$\tau(y):\nabla^2\psi(y)=a(Ty):\nabla^2\varphi(Ty)$. Thus \eqref{B:eq:ratio-scalar-integration-by-parts} gives the integral identity for $\nu,a$. The bounds
\[
 \tfrac12C_1\le a(x)\le\tfrac32C_1,\qquad
 \tfrac12\lambda_{\min}(C_1)I\le a(x)
 \le\tfrac32\lambda_{\max}(C_1)I
\]
verify the hypotheses of Lemma~\ref{B:lem:invariant-diffusion} at this fixed invertible $T$.

For each $x$, choose the smooth orthogonal polar factor $O(x)$ of $C_1^{1/2}a(x)^{1/2}$. Use one Brownian motion in
\begin{equation}
 \begin{split}
 \dd X_t&=-X_t\,\dd t+\sqrt2\,a(X_t)^{1/2}\,\dd B_t,\\
 \dd Z_t^T&=-Z_t^T\,\dd t+\sqrt2\,C_1^{1/2}O(X_t)\,\dd B_t.
 \end{split}
 \label{B:eq:coupled-comparison-diffusions}
\end{equation}
Start the two coordinates with independent laws $\nu,N(0,C_1)$, independent of the Brownian motion. The two diffusion coefficients are smooth and locally Lipschitz, and both are
bounded in Frobenius norm. For the joint quadratic Lyapunov function the generator gives
\[
 \mathcal A_{\rm joint}(|x|^2+|z|^2)
 =-2(|x|^2+|z|^2)+2\operatorname{Tr}a(x)+2\operatorname{Tr}C_1
 \le C_T-2(|x|^2+|z|^2).
\]
The stopped argument in the preceding lemma gives a unique nonexplosive joint solution. Its first marginal remains $\nu$. The process $\widetilde B_t=\int_0^tO(X_s)\,\dd B_s$ is a
continuous martingale in the common filtration, with
\[
 \dd\langle\widetilde B^i,\widetilde B^j\rangle_t
 =(O(X_t)O(X_t)^{\mathsf T})_{ij}\,\dd t=\delta_{ij}\,\dd t.
\]
L\'evy's characterization makes it Brownian in that filtration, including independence of its future increments from the initial sigma-field. The second equation has the mild form
$Z_t^T=e^{-t}Z_0^T+\sqrt2\int_0^te^{-(t-s)}C_1^{1/2}\,\dd\widetilde B_s$. The Brownian increments of $\widetilde B$ are independent of $\mathcal F_0$, and the integrand in this
mild formula is deterministic. For $u\in\mathbb R^d$ its conditional characteristic function is
\[
 \mathbb E[e^{\ii u^{\mathsf T}Z_t^T}\mid\mathcal F_0]
 =\exp\!\left(\ii e^{-t}u^{\mathsf T}Z_0^T
           -\tfrac12(1-e^{-2t})u^{\mathsf T}C_1u\right).
\]
Taking expectations and using $Z_0^T\sim N(0,C_1)$ gives
\[
 \mathbb E e^{\ii u^{\mathsf T}Z_t^T}
 =e^{-(1-e^{-2t})u^{\mathsf T}C_1u/2}
      e^{-e^{-2t}u^{\mathsf T}C_1u/2}
 =e^{-u^{\mathsf T}C_1u/2}.
\]
This proves $Z_t^T\sim N(0,C_1)$. Both marginals are therefore fixed in time, although the two coordinates need not be independent at positive times.

Set $D_t=X_t-Z_t^T$, $H_t=a(X_t)^{1/2}-C_1^{1/2}O(X_t)$, and $A_*(X_t)=\|H_t\|_F^2$. Then $\dd D_t=-D_t\,\dd t+\sqrt2H_t\,\dd B_t$. For a finite horizon $t_0$, the mild equations
express $D_t$, $X_t$, and $Z_t^T$ as their exponentially damped initial values plus $\sqrt2e^{-t}M_t$, where each corresponding martingale has the form
$M_t=\int_0^te^sC_s\,\dd B_s$ with deterministically bounded $\|C_s\|_F$ at this fixed $T$. Hence each coordinate has quadratic variation at most a finite $V_*$ on $[0,t_0]$. The
exponential supermartingale bound gives
\[
 \mathbb P\{\sup_{t\le t_0}|M_t|>u\}
 \le2d\exp\{-u^2/(2dV_*)\}.
\]
If $V_*=0$, the martingale vanishes. Integrating this tail and using the finite moments of the initial marginals proves $\mathbb E\sup_{t\le t_0}(|X_t|+|Z_t^T|)^q<\infty$ for every
finite $q>0$. These bounds justify removal of stopping in It\^o's formula below by dominated convergence. Their constants may depend on the fixed dimension and time horizon; the
final estimate follows from the uniform bound for the quadratic form defined by the matrix Lyapunov equation.

For $p\ge1$, use $F_\epsilon(d)=(|d|^2+\epsilon)^p$. Its derivatives are
\[
 \nabla F_\epsilon(d)=2p(|d|^2+\epsilon)^{p-1}d,\qquad
 \nabla^2F_\epsilon(d)
 =2p(|d|^2+\epsilon)^{p-1}I
  +4p(p-1)(|d|^2+\epsilon)^{p-2}dd^{\mathsf T}.
\]
The second term is zero for $p=1$. For $p>1$, letting $\epsilon\downarrow0$ in the stopped It\^o formula is justified by the finite time moments and
$(|d|^2+\epsilon)^{p-2}|d|^2\le(|d|^2+\epsilon)^{p-1}$. The drift contribution converges to $-2p|D_t|^{2p}$. The two diffusion contributions satisfy
\[
 2p|D_t|^{2p-2}\|H_t\|_F^2
 +4p(p-1)|D_t|^{2p-4}|H_t^{\mathsf T}D_t|^2 
 \le2p(2p-1)|D_t|^{2p-2}A_*(X_t),
\]
where $|H_t^{\mathsf T}D_t|^2\le\|H_t\|_F^2|D_t|^2$. For $p>1$ the expression at $D_t=0$ is its continuous limit; for $p=1$ only the first term remains.

Let $y(t)=\mathbb E|D_t|^{2p}$ and $A_p=\|A_*(X_0)\|_{L^p}$. Stationarity of $X_t$ and H\"older give, for almost every
$t$,
\[
 y'(t)\le-2py(t)+2p(2p-1)y(t)^{1-1/p}A_p.
\]
For $p=1$ this is a linear inequality. For $p>1$, put $z_\epsilon(t)=(y(t)+\epsilon)^{1/p}$. Division by $p(y+\epsilon)^{1-1/p}$ gives
$z_\epsilon'\le-2z_\epsilon+2\epsilon^{1/p}+2(2p-1)A_p$. Multiplication by $e^{2t}$, integration, and then $\epsilon\downarrow0$ yield
\begin{equation}
 \|D_t\|_{L^{2p}}^2
 \le e^{-2t}\|D_0\|_{L^{2p}}^2
 +(2p-1)(1-e^{-2t})\|A_*(X_0)\|_{L^p}.
 \label{B:eq:diffusion-distance-evolution}
\end{equation}
The estimate follows by taking expectation in It\^o's formula for the joint process $(D_t,X_t)$ and solving the resulting differential inequality.

The joint laws at integer times are tight: for compact sets $K_1,K_2$,
\[
 \Pr((X_m,Z_m^T)\notin K_1\times K_2)
 \le\nu(K_1^c)+N(0,C_1)(K_2^c).
\]
Choose one weakly convergent subsequence, independently of $p$. Its limit has the same two marginals. For each real $p\ge1$ and $K<\infty$, $\min\{|x-z|^{2p},K\}$ is bounded and
continuous. Passing to the limit in its expectation and then letting $K\uparrow\infty$ gives
\[
 \left(\int|x-z|^{2p}\,\dd\pi(x,z)\right)^{1/p}
 \le(2p-1)A_p.
\]
This holds for every $p$ under the same limit law: the subsequence was chosen solely from tightness, and the argument with bounded continuous test functions applies to any
specified $p$.

The minimizing polar factor and Lemma~\ref{B:lem:gaussian-factor-distance} give, with $E=\tau(Y)-B$,
\[
 A_*(TY)\le
 \langle TET^{\mathsf T},\mathcal L_{C_1}^{-1}(TET^{\mathsf T})\rangle_F.
\]
Since $C_1\ge\lambda_{\min}(B)TT^{\mathsf T}$, the two trace terms in the Lyapunov quadratic form imply $\mathcal L_{C_1}\ge\lambda_{\min}(B)\mathcal L_{TT^{\mathsf T}}$. The
variational formula for inverse quadratic forms and \eqref{B:eq:lyapunov-congruence-identity} therefore yield
\[
 A_*(TY)\le\lambda_{\min}(B)^{-1}
   \langle E,\mathcal L_{(T^{\mathsf T}T)^{-1}}^{-1}E\rangle_F
 =\lambda_{\min}(B)^{-1}\|\tau(Y)-B\|_{\mathscr H_{\mathsf M}}^2.
\]
Apply $T^{-1}$ to both coordinates under $\pi$ to obtain the required coupling of $Y$ and $Z\sim N(0,B)$ and its moment inequality.
\end{proof}

\subsubsection{Proof of Theorem~\ref{B:thm:logarithmic-ratio-coupling}}
\label{B:subsection:proof-of-ratio-coupling}

\begin{proof}[Proof of Theorem~\ref{B:thm:logarithmic-ratio-coupling}]
Choose the fixed $\delta$ for which the estimates on the Poisson solutions and the diffusion matrix hold, and put $T=\mathsf T_2D_J$. The map $T$ is invertible because both its
diagonal factors and $S_n$ are invertible. The smooth positive density of $Y=f-\langle f\rangle$, its finite moments, and $B/2\le\tau\le3B/2$ verify the hypotheses of
Proposition~\ref{B:prop:diffusion-coupling-moments}. Together with Proposition~\ref{B:prop:diffusion-matrix-deviation} and $\lambda_{\min}(B)\ge b_->0$, it supplies one coupling of
$Y$ and $Z\sim N(0,B)$ for which
\begin{equation}
 \|T(Y-Z)\|_{2p}^2
 \le Cp\left[n^{-1}(\sqrt p\,\ell_n+\ell_n^{3/2})\right]^2
 \le Cp\ell_n^2(p+\ell_n)/n^2\quad(p\ge1).
 \label{B:eq:ratio-centered-coupling}
\end{equation}
The last inequality follows from $(\sqrt p\,\ell_n+\ell_n^{3/2})^2\le2p\ell_n^2+2\ell_n^3$. The coupling was chosen before the order $p$.

The normalized linear map in the extension obeys $D_JL=\kappa^{-1/2}L_q$. Consequently $D_JBD_J=\kappa^{-1}L_qCL_q^{\mathsf T}=K_q$, so $u=D_JZ$ has exactly the law $\Gamma_q$. Put
$v=\widetilde\eta=D_J(Y+\langle f\rangle)$. Its marginal is the extended logarithmic distribution. For $b_k=\langle f_k\rangle/\sqrt k$, the signed mean estimate gives
$|\sum_{k=j}^Rb_k|\le C/j$. Hence
\begin{equation}
 |\mathsf T_2b|^2\le\frac C{n^2}\sum_{j=r}^Rj^{-1}
                         \le C\ell_n/n^2.
 \label{B:eq:ratio-partial-sum-mean}
\end{equation}
Here $\sum_{j=r}^Rj^{-1}\le1+\log(R/r)\le C\ell_n$. Thus $\mathsf T_2(v-u)=T(Y-Z)+\mathsf T_2b$ with a deterministic mean term of squared order $\ell_n/n^2$.

Let $\chi_q(\dd\eta,\dd v)$ be the coupling of the original and extended ratios from Proposition~\ref{B:prop:conditional-density-coupling}, and let $\widetilde\nu_q$ be its second
marginal. Disintegrate this probability as $K_q^{\rm ex}(v,\dd\eta)\widetilde\nu_q(\dd v)$. Regular conditional probabilities exist on these Euclidean spaces; assign an arbitrary
fixed probability to the kernel on its unused null set. If $\pi_q^{\rm ext}(\dd v,\dd u)$ is the coupling just constructed, use
\begin{equation}
 \pi_q^{\rm ext}(\dd v,\dd u)K_q^{\rm ex}(v,\dd\eta)
 \label{B:eq:ratio-exact-coupling-composition}
\end{equation}
For Borel sets $A,D,C$, respectively in the $\eta,v,u$ spaces, its marginals are verified by
\[
 \begin{split}
 \int K_q^{\rm ex}(v,A)\,\pi_q^{\rm ext}(\dd v,\dd u)&=\nu_q(A),\\
 \int\mathbf1_C(u)K_q^{\rm ex}(v,\mathbb R^d)\,
           \pi_q^{\rm ext}(\dd v,\dd u)&=\Gamma_q(C),\\
 \int\mathbf1_D(v)K_q^{\rm ex}(v,A)\,
           \pi_q^{\rm ext}(\dd v,\dd u)&=\chi_q(A\times D).
 \end{split}
\]
The full $(v,u)$ marginal is also unchanged, since the conditional kernel has mass one. Therefore $\Pr(\eta\ne v)\le\varepsilon_n=Cne^{-cr}$. At fixed $n,q$ this defines a single
probability before any moment order is chosen. The following subsubsections make the choice of coupling measurable in the boundary parameter.

For the original ratios the moment bound at order $4p$ is $\|\eta_k\|_{4p}\le C\sqrt{p/k}$ in the stated range. For the extended ratios, $l_\alpha(0)=0$ and
$|l_\alpha'|\le C/d_\alpha$ give the same bound from the coordinate moments. Since $\sum_{k=j}^Rk^{-1/2}\le2\sqrt R\le2\sqrt n$, Minkowski gives
\[
 \begin{split}
 \left\|\sum_{k=j}^R\eta_k\right\|_{4p}
 +\left\|\sum_{k=j}^Rv_k\right\|_{4p}&\le C\sqrt{pn},\\
 \|\mathsf T_2\eta\|_{4p}^2
 &=\left\|n^{-2}\sum_{j=r}^Rj
            \left|\sum_{k=j}^R\eta_k\right|^2\right\|_{2p}
 \le Cpn\,n^{-2}\sum_{j=r}^Rj\le Cpn.
 \end{split}
\]
Replacing $\eta$ by $v$ gives the same inequality. On the event $\{\eta=v\}$ their difference is zero. H\"older with exponents two and two thus gives
\begin{equation}
 \|\mathsf T_2(\eta-v)\|_{2p}
                    \le C\sqrt{pn}\,(Cne^{-cr})^{1/(4p)}.
 \label{B:eq:ratio-mismatch-moment}
\end{equation}
For $2\le p\le C_0\ell_n$, choose $n$ with $4C_0\ell_n\le c_{\rm rc}r$ and $\varepsilon_n<1$. The negative numerator $\log\varepsilon_n$ is divided by at most $4C_0\ell_n$, and
$p\le C_0\ell_n$. Thus the logarithm of the bound is at most
\[
 C+\tfrac12\log(C_0\ell_n n)
       +\frac{\log C+\log n-cr}{4C_0\ell_n}.
\]
Here $r\asymp(\log n)^{A_0}$, so $r/(\ell_n\log n)\to\infty$; the fixed assumption $A_0>20$ suffices. For large $n$ this expression is at most $-10\log n$, uniformly over the
permitted orders. Hence $\|\mathsf T_2(\eta-v)\|_{2p}\le n^{-10}$.

Finally $\mathsf T_2(\eta-u)=\mathsf T_2(\eta-v)+T(Y-Z)+\mathsf T_2b$. Taking the $L^{2p}$ norm and squaring gives
\[
 \begin{split}
 \|\mathsf Q_2(\eta-u)\|_p
 &=\|\mathsf T_2(\eta-u)\|_{2p}^2\\
 &\le3n^{-20}
   +Cp\ell_n^2(p+\ell_n)n^{-2}+C\ell_n n^{-2}
 \le Cp\ell_n^2(p+\ell_n)n^{-2}.
 \end{split}
\]
This proves \eqref{B:eq:ratio-coupling-moment-bound}, and \eqref{B:eq:ratio-norm-comparisons} gives its two stated consequences. The constants are uniform on the admissible
boundary set; only the starting size required for the moment orders and for the bound on the probability of unequal coordinates depends on $C_0$.
\end{proof}

\subsubsection{Measurable optimal couplings}
\label{B:subsection:measurable-optimal-ratio-couplings}

The preceding theorem provides a coupling for each admissible boundary. To use it with the original random exterior, we choose couplings measurably on the full boundary space. The
cost below contains powers of the weighted quadratic form and high powers of the individual partial sums. Fix the weights and $\lambda$ from \eqref{B:eq:ratio-partial-sum-cost}.
Choose an integer $p_0\ge2$ with $p_0(\lambda-1)>4$ and set
\begin{equation}
 P_n=\max\{p_0,\lceil\ell_n\rceil\},\qquad
 \widehat p_n=2\lceil4\ell_n+p_0\rceil.
 \label{B:eq:ratio-cost-exponents}
\end{equation}
For a deterministic $D\ge1$, define
\begin{equation}
 \mathfrak c_n(e)=
 \sum_{p=p_0}^{P_n}\left(\frac{nQ_n(e)}{Dp^4}\right)^p
 +\sum_{j\in J}\left(
        \frac{\sqrt j\,|(S_ne)_j|}{D\widehat p_n^2}\right)^{\widehat p_n}.
 \label{B:eq:ratio-selected-cost}
\end{equation}
The constant $D$ will be fixed in Section~\ref{B:subsection:optimal-ratio-cost-bound}; the construction in this subsubsection works for every deterministic choice. The exponents in
the first sum are integers and $\widehat p_n$ is even, so the cost is a nonnegative continuous polynomial in $e$. Write $W_n=\operatorname{diag}(w_{n,j})$. Positivity of every
weight and invertibility of $S_n$ imply $c_n^Q=\lambda_{\min}(S_n^{\mathsf T}W_nS_n)>0$. Consequently
\[
 \mathfrak c_n(e)\ge
 \left(\frac{nc_n^Q|e|^2}{Dp_0^4}\right)^{p_0}
 \longrightarrow\infty\qquad(|e|\to\infty).
\]
This is coercivity in the difference variable $e$. We minimize its integral over $\Pi(\Gamma_q,\nu_q)$ and write the minimum as $\mathfrak I_n(q)$.

\begin{theorem}\label{B:thm:measurable-ratio-optimizer}
At each fixed $n$, there is a Borel family of optimal plans $q\mapsto\pi_{n,q}^{\rm rat}\in\Pi(\Gamma_q,\nu_q)$ on the full space of finite nonnegative boundary configurations. Its
value is
\begin{equation}
 \int\mathfrak c_n(\eta-u)\,\pi_{n,q}^{\rm rat}(\dd u,\dd\eta)
 =\mathfrak I_n(q)
 =\min_{\pi\in\Pi(\Gamma_q,\nu_q)}
                       \int\mathfrak c_n(\eta-u)\,\dd\pi<\infty.
 \label{B:eq:ratio-optimal-plan}
\end{equation}
The function $\mathfrak I_n$ is continuous.
\end{theorem}
\begin{proof}
Fix $n$ and put $h=\max\{2P_n,\widehat p_n\}$. The parameter space is a closed Euclidean orthant. The preceding continuity lemmas give weak continuity of $q\mapsto(\Gamma_q,\nu_q)$
and locally uniformly integrable moments of every finite order. The continuous nonnegative cost satisfies
\[
 \mathfrak c_n(\eta-u)\le C_n(1+|u|^h+|\eta|^h).
\]
Thus the product coupling has finite cost at every $q$. The set $\Pi(\Gamma_q,\nu_q)$ is tight, because any plan in it obeys
\[
 \pi((K_1\times K_2)^c)\le\Gamma_q(K_1^c)+\nu_q(K_2^c).
\]
It is weakly closed: the two coordinate projections are continuous, so weak limits preserve the two marginals. Prokhorov compactness and lower semicontinuity of the nonnegative
cost integral therefore give a minimizer. The set of minimizers is nonempty and compact.

We verify continuity of the value. Let $q_j\to q$. For a subsequence along which the lower limit of $\mathfrak I_n(q_j)$ is attained, choose minimizing plans $\pi_j$. Tightness of
the two marginal families makes the plans tight. Any weakly convergent further subsequence has a limit in $\Pi(\Gamma_q,\nu_q)$, and lower semicontinuity gives
\[
 \mathfrak I_n(q)\le\int\mathfrak c_n(\eta-u)\,\dd\pi
 \le\liminf_j\mathfrak I_n(q_j).
\]
For the reverse inequality fix a minimizer $\pi$ at $q$. The marginal continuity lemmas give total variation convergence. Couple each $\Gamma_{q_j}$ to $\Gamma_q$, and each
$\nu_{q_j}$ to $\nu_q$, so that the unequal-coordinate probabilities tend to zero. Such couplings retain the pointwise minimum of the two densities on the diagonal and
independently couple the normalized nonnegative differences between each density and that minimum. The unequal-coordinate probability is their total variation distance.
Disintegrate each of these couplings over its limiting coordinate and apply the two resulting kernels independently to $(U,\Xi)\sim\pi$. The pair $(U_j,\Xi_j)$ then has marginals
$\Gamma_{q_j},\nu_{q_j}$, and
\[
 \Pr\{(U_j,\Xi_j)\ne(U,\Xi)\}
 \le d_{\rm TV}(\Gamma_{q_j},\Gamma_q)+d_{\rm TV}(\nu_{q_j},\nu_q)\longrightarrow0.
\]
Its cost converges in probability to $\mathfrak c_n(\Xi-U)$. The marginal moment bounds of order $2h$ imply uniform integrability of these costs, since
\[
 \sup_j\mathbb E[\mathfrak c_n(\Xi_j-U_j)^2]
 \le C_n\sup_j\mathbb E(1+|U_j|^{2h}+|\Xi_j|^{2h})<\infty.
\]
Expectations therefore converge, and these pairs are competitors at $q_j$. Hence $\limsup_j\mathfrak I_n(q_j)\le\mathfrak I_n(q)$. This proves continuity, including at the boundary
of the parameter orthant.

Finally apply the theorem on measurable selection of optimal couplings in \cite{bogachevMalofeev2020parameter}. Its parameter and transport spaces are Polish, both marginal maps
are Borel, the cost is fixed, continuous and nonnegative, and its infimum is finite and attained as just proved. For completeness the minimizing graph is Borel: the marginal
constraints form a closed graph, the cost integral is lower semicontinuous, and the value is continuous; equality to the value may be written as the cost being at most that value
on the feasible graph. The theorem supplies a Borel selector of its nonempty compact sections on the entire orthant. Its marginal identities and minimizing property are exactly
\eqref{B:eq:ratio-optimal-plan}.
\end{proof}

\subsubsection{Conditional probability kernels}
\label{B:subsection:ratio-conditional-probability-kernels}

The selected plans specify joint laws, but the construction needs to sample the ratio vector given the Gaussian vector. A parameterwise disintegration provides that conditional
kernel while preserving measurability in the boundary configuration.

\begin{lemma}\label{B:lem:parameterized-ratio-disintegration}
There is a Borel probability kernel $K_n^{\rm rat}(q,u;\dd\eta)$ satisfying, for every $q$,
\begin{equation}
 \Gamma_q(\dd u)K_n^{\rm rat}(q,u;\dd\eta)
                      =\pi_{n,q}^{\rm rat}(\dd u,\dd\eta).
 \label{B:eq:ratio-parameterized-disintegration}
\end{equation}
\end{lemma}
\begin{proof}
For $m\ge0$, partition $\mathbb R^d$ into the half-open dyadic cubes of side $2^{-m}$, and let $C_m(u)$ be the cube containing $u$. The partitions are countable, nested, and
generate the Borel sigma-field. Define a probability measure
\[
 K_m(q,u;A)=
 \frac{\pi_{n,q}^{\rm rat}(C_m(u)\times A)}{\Gamma_q(C_m(u))}.
\]
The denominator is positive: $\Gamma_q$ is a nondegenerate Gaussian measure for every finite boundary configuration. For each cube, both numerator and denominator are Borel
functions of $q$. Summing over the countably many cubes shows that $(q,u)\mapsto K_m(q,u;A)$ is Borel for every Borel $A$. In particular $K_m$ is a probability kernel, and is a
Borel map into the Polish space $\mathcal P(\mathbb R^d)$ with its weak topology.

Fix $q$ temporarily and work under $\pi_{n,q}^{\rm rat}$, with coordinates $(U,\Xi)$. For every bounded Borel $f$,
\[
 \int f(\eta)K_m(q,U;\dd\eta)
 =\mathbb E_q[f(\Xi)\mid\sigma(C_m(U))].
\]
Choose a countable uniformly dense family in $C_c(\mathbb R^d)$ on each compact cube, together with continuous cutoffs $0\le\chi_R\le1$ increasing to one, equal to one on
$[-R,R]^d$ and supported in $[-R-1,R+1]^d$. Upward martingale convergence, applied to this countable family, gives convergence of all the displayed integrals for $\Gamma_q$-almost
every $u$. The limits for $\chi_R$ increase to one, by conditional monotone convergence. On this same set of $u$, the sequence $K_m(q,u;\cdot)$ is tight: for a given
$\varepsilon>0$, choose $R$ so that its limiting $\chi_R$ integral exceeds $1-\varepsilon/2$, and then choose $m_0$ so that it exceeds $1-\varepsilon$ for $m\ge m_0$. Enlarge the
compact cube to control the finitely many measures with $m<m_0$.

Every subsequential weak limit has the prescribed integrals against the chosen countable family. Uniform approximation on compact sets shows that this family determines a
probability measure, so any two subsequential limits agree. Tightness therefore gives a weak limit of the entire sequence $K_m(q,u;\cdot)$ for $\Gamma_q$-almost every $u$.

To make the construction joint in the parameter, fix a complete metric $d_{\mathcal P}$ inducing weak convergence on $\mathcal P(\mathbb R^d)$. Its Cauchy criterion shows that
\[
 G=\bigcap_{r\ge1}\ \bigcup_{M\ge1}\ \bigcap_{k,l\ge M}
 \{(q,u):d_{\mathcal P}(K_k(q,u),K_l(q,u))<r^{-1}\}
\]
is Borel. On $G$ take the limit, and on $G^c$ take $\delta_0$. Limits of Borel maps into a metric space are Borel on their convergence set; thus the result is a jointly Borel
probability kernel $K_n^{\rm rat}$. The preceding argument proves $\Gamma_q(G_q)=1$ for every $q$, without requiring the exceptional sets $G_q^c$ to coincide.

For a bounded continuous $f$ and a Borel set $D$ in the first coordinate, the martingale identity and bounded convergence give
\[
 \int_D\!\int f(\eta)K_n^{\rm rat}(q,u;\dd\eta)\,\Gamma_q(\dd u)
 =\int_{D\times\mathbb R^d}f(\eta)\,\pi_{n,q}^{\rm rat}(\dd u,\dd\eta).
\]
A monotone class argument extends this equality to bounded Borel $f$. Taking indicators proves the required disintegration on rectangles, and hence on the product Borel
sigma-field, separately for every boundary configuration.
\end{proof}

\begin{lemma}\label{B:lem:ratio-probability-kernel-sampling}
There is a jointly Borel map $\mathfrak A_n(q,u,v)$ such that
\begin{equation}
 \int_0^1\mathbf1_A(\mathfrak A_n(q,u,v))\,\dd v
                        =K_n^{\rm rat}(q,u;A)
 \label{B:eq:ratio-sampler-kernel}
\end{equation}
for every $(q,u)$ and Borel set $A\subset\mathbb R^d$.
\end{lemma}
\begin{proof}
Write $\theta=(q,u)$ and $\mu_\theta=K_n^{\rm rat}(q,u;\cdot)$. Enumerate the half-open unit cubes of $\mathbb R^d$ as $(C_j)_{j\ge1}$. Partition $[0,1)$ into consecutive half-open
intervals of lengths $\mu_\theta(C_j)$, ignoring intervals of length zero. Inside the interval assigned to $C_j$, assign consecutive subintervals to its $2^d$ dyadic children, in a
fixed order, with lengths equal to their $\mu_\theta$ masses. Continue recursively. The subintervals at each level cover $[0,1)$; they are nested with their parent intervals, and
the subinterval assigned to every cube $C$ has length $\mu_\theta(C)$.

Every endpoint is a finite or countable sum of values $\mu_\theta(C)$ and is therefore Borel in $\theta$. For $0\le v<1$, let $C_m(\theta,v)$ be the unique cube of side $2^{-m}$
whose assigned interval contains $v$, and let $c_m(\theta,v)$ be its center. The sets on which a given cube is selected form a countable Borel partition of the product of the
parameter space and the unit interval, so each $c_m$ is jointly Borel. The closed cubes are nested inside one closed unit cube and have diameters $\sqrt d\,2^{-m}$. Consequently
the centers converge to a finite point, and
\[
 |c_m(\theta,v)-c_l(\theta,v)|\le\sqrt d\,2^{-m}
 \qquad(l\ge m).
\]
Define $\mathfrak A_n(\theta,v)$ to be their limit for $v<1$, and set it to zero for $v=1$. This is a jointly Borel map. The calculation below also covers limits on the faces of
the selected half-open cubes.

Indeed, the distribution of $c_m(\theta,V)$ for a uniform $V$ assigns mass $\mu_\theta(C)$ to the center $c_C$ of each level-$m$ cube. For a bounded continuous $f$, bounded
convergence in $v$ and then in $x$ gives
\[
 \begin{split}
 \int_0^1f(\mathfrak A_n(\theta,v))\,\dd v
 &=\lim_{m\to\infty}\sum_C\mu_\theta(C)f(c_C)\\
 &=\lim_{m\to\infty}\int f(c_{C_m(x)})\,\mu_\theta(\dd x)
 =\int f(x)\,\mu_\theta(\dd x).
 \end{split}
\]
In the last limit $|c_{C_m(x)}-x|\le\sqrt d\,2^{-m}$ for every $x$. Equality of integrals against bounded continuous functions identifies the two probability measures. It therefore
holds for every Borel indicator, which is \eqref{B:eq:ratio-sampler-kernel}. The kernel identity holds for every $\theta$.
\end{proof}

\subsubsection{Reconstruction of the coefficient distribution}
\label{B:subsection:joint-logarithmic-reconstruction}

To recover the positive coefficient vector, we sample the remaining logarithmic scale variables from their exact conditional density. The resulting vector has the full conditional
coefficient law.

For the conditional density $\mathfrak r_q$ in \eqref{B:eq:ratio-marginal-and-conditional}, order the sites of $J$ and set
\begin{equation}
 I_q^{(l)}(s_1,\ldots,s_l;\eta)
    =\int g_q(s,\eta)\,\dd s_{l+1}\cdots\dd s_d,
 \qquad I_q^{(0)}=I_q.
 \label{B:eq:ratio-radial-marginals}
\end{equation}
For every fixed finite set of displayed parameters, these integrals are finite and positive by \eqref{B:eq:ratio-density-majorant}. Define their conditional distribution functions
by
\begin{equation}
 F_{q,\eta,l}(t\mid s_{<l})
  =\frac{\displaystyle\int_{-\infty}^t
                  I_q^{(l)}(s_{<l},u;\eta)\,\dd u}
         {I_q^{(l-1)}(s_{<l};\eta)}.
 \label{B:eq:ratio-radial-conditional-cdf}
\end{equation}

\begin{lemma}\label{B:lem:joint-logarithmic-scale-sampling}
The successive inverse distribution functions in \eqref{B:eq:ratio-radial-conditional-cdf} define a Borel map $\mathfrak R_n(q,\eta,U^{\rm rad})$ with distribution
$\mathfrak r_q(\cdot\mid\eta)$ when $U^{\rm rad}$ is a vector of $d$ independent uniform variables.
\end{lemma}
\begin{proof}
Order $J$ as in \eqref{B:eq:ratio-radial-marginals}. The positive integrand $g_q$ is continuous in its coefficient variables and Borel in every parameter. Its bound
\eqref{B:eq:ratio-density-majorant} remains valid after any finite list of $s$ coordinates is fixed on a compact set. Integrating the remaining exponential majorants gives finite
bounds by the gamma integral used in Proposition~\ref{B:prop:logarithmic-ratio-density}. Hence every $I_q^{(l)}$ is finite and strictly positive at every finite value of their
arguments, Borel jointly, and continuous in each displayed real coordinate. Tonelli gives, for $1\le l\le d$,
\[
 \int_{\mathbb R}I_q^{(l)}(s_{<l},t;\eta)\,\dd t
                         =I_q^{(l-1)}(s_{<l};\eta).
\]
The denominator is thus precisely the mass in the numerator's density. Consequently $F_{q,\eta,l}$ is continuous, has strictly positive derivative
\[
 \partial_t F_{q,\eta,l}(t\mid s_{<l})
   =\frac{I_q^{(l)}(s_{<l},t;\eta)}
          {I_q^{(l-1)}(s_{<l};\eta)}>0,
\]
and tends to zero and one at the two ends of the real line.

For $0<v<1$, write its inverse as
\[
 H_l(q,\eta,s_{<l},v)=
 \inf\{t\in\mathbb Q:F_{q,\eta,l}(t\mid s_{<l})\ge v\}.
\]
The set is nonempty and bounded below by the two limits of $F$; continuity and strict monotonicity identify the infimum with the unique real solution of $F(t)=v$. For every real
$a$,
\[
 \{H_l<a\}\cap\{0<v<1\}
  =\bigcup_{\substack{t\in\mathbb Q\\t<a}}
       \{0<v<1,\ v\le F_{q,\eta,l}(t\mid s_{<l})\}.
\]
These are Borel sets, proving joint measurability of the inverse. Assign it the value zero for $v=0,1$. For a uniform variable $V_l$, the equality
$\{H_l(V_l)\le t\}=\{V_l\le F_{q,\eta,l}(t\mid s_{<l})\}$ away from the null boundary points proves the conditional distribution identity.

Define $s_1=H_1(q,\eta,U_1)$ and then $s_l=H_l(q,\eta,s_{<l},U_l)$, with the uniforms independent. Induction and composition give a Borel vector map. Conditional on its first $l-1$
coordinates, the variable $U_l$ remains uniform; the last calculation gives density $I_q^{(l)}(s_{\le l};\eta)/I_q^{(l-1)}(s_{<l};\eta)$. Iterated integration against a bounded
Borel function therefore gives the joint density
\[
 \prod_{j=1}^{l}
 \frac{I_q^{(j)}(s_{\le j};\eta)}
      {I_q^{(j-1)}(s_{<j};\eta)}
          =\frac{I_q^{(l)}(s_{\le l};\eta)}{I_q(\eta)}.
\]
Every cancellation is between positive finite factors. At $l=d$, $I_q^{(d)}=g_q$, so this is exactly the conditional density in \eqref{B:eq:ratio-marginal-and-conditional}. The
union of the $2d$ uniform boundary point events has measure zero. Their assigned values make the map finite and defined everywhere without affecting its law.
\end{proof}

For $G\sim\gamma_{n,q}$, let $u=L_q(G-m^q-s^q)$ and use independent uniforms, independent of $G$, to set
\begin{equation}
 \eta=\mathfrak A_n(q,u,U^{\rm rat}),\qquad
 s=\mathfrak R_n(q,\eta,U^{\rm rad}),\qquad
 Z=\mathfrak Z_q(s,\eta).
 \label{B:eq:ratio-coefficient-reconstruction}
\end{equation}

\begin{proposition}\label{B:prop:conditional-coefficient-reconstruction}
The construction \eqref{B:eq:ratio-coefficient-reconstruction} has $Z\sim\mu_q$ and $\mathcal L(u,\eta)=\pi_{n,q}^{\rm rat}$. All maps are Borel, and the logarithmic ratios of $Z$
are exactly $\eta$.
\end{proposition}
\begin{proof}
Fix any finite nonnegative $q$. The preceding lemmas define finite Borel maps on their stated domains; $\mathfrak Z_q$ takes values in the positive coefficient orthant. If
$u_q(g)=L_q(g-m^q-s^q)$, integration over the new uniforms gives, for bounded Borel $F$,
\[
 \mathbb EF(G,Z)=\int\gamma_{n,q}(\dd g)
   \int K_n^{\rm rat}(q,u_q(g);\dd\eta)
   \int\mathfrak r_q(\dd s\mid\eta)
                      F(g,\mathfrak Z_q(s,\eta)).
\]
The absolute integral is at most $\|F\|_\infty$, because all three measures are probabilities. This justifies the successive integrations also for signed or complex-valued
functions. Taking $F(g,z)=F_0(g)$ integrates the last two kernels to one at each $g$, so the full Gaussian vector keeps its law.

For a bounded Borel $H(u,\eta)$, the distribution formula for the Gaussian linear image and \eqref{B:eq:ratio-parameterized-disintegration} give
\[
 \mathbb EH(u,\eta)
 =\int\Gamma_q(\dd u)\int K_n^{\rm rat}(q,u;\dd\eta)H(u,\eta)
 =\int H\,\dd\pi_{n,q}^{\rm rat}.
\]
In particular the marginal of $\eta$ is $\nu_q$. For a bounded Borel function $F_1$ of the coefficients, substitution of the two densities in
\eqref{B:eq:ratio-marginal-and-conditional} then gives
\[
 \mathbb EF_1(Z)
 =\int\frac{I_q(\eta)}{\mathcal Z_q^{\rm rat}}\,\dd\eta
        \int\frac{g_q(s,\eta)}{I_q(\eta)}
                             F_1(\mathfrak Z_q(s,\eta))\,\dd s
 =\frac1{\mathcal Z_q^{\rm rat}}
        \int g_q(s,\eta)F_1(\mathfrak Z_q(s,\eta))\,\dd s\,\dd\eta.
\]
To check the normalizer in the last change of variables, recall $a_q^0=\prod_{k\in J}(m_{x,k}^q)^{\beta(k+a)}(m_{y,k}^q)^{\beta k}$ from \eqref{B:eq:ratio-joint-density} and its
derivation. That calculation gives
\[
 e^{-\kappa\mathcal P_{n,J}^q(z)}
      \prod_\alpha z_\alpha^{c_\alpha}\,\dd z
  =a_q^0g_q(s,\eta)\,\dd s\,\dd\eta,\qquad
 \mathcal Z_q^{\rm coef}=a_q^0\mathcal Z_q^{\rm rat}.
\]
Thus $\mathbb EF_1(Z)=\int F_1\,\dd\mu_q$ with the original conditional normalization.

For every value of the sampling variables,
\[
 \log\frac{Z_{x,k}}{m_{x,k}^q}=s_k-\eta_k/2,\qquad
 \log\frac{Z_{y,k}}{m_{y,k}^q}=s_k+\eta_k/2.
\]
Their difference is $\eta_k$. This verifies the exact coordinate identity as well as positivity and finiteness of every reconstructed coefficient.
\end{proof}

\subsubsection{The bound for the optimal cost}
\label{B:subsection:optimal-ratio-cost-bound}

We now bound the expected cost of the selected coupling. The moment estimate gives a uniform choice of the constant in the cost, including the increasing moment orders used to
control the maximum of the partial sums.

To fix $D$ in \eqref{B:eq:ratio-selected-cost}, use the moment estimate as follows. Apply Theorem~\ref{B:thm:logarithmic-ratio-coupling} with $C_0=16$. For all sufficiently large
$n$,
\begin{equation}
 P_n\le2\ell_n,\qquad 8\ell_n\le\widehat p_n\le12\ell_n,
 \qquad64\ell_n\le c_{\rm rc}r_n.
 \label{B:eq:ratio-growing-cost-range}
\end{equation}
The last inequality permits the logarithmic moments used when composing the exact distributions. By \eqref{B:eq:ratio-coupling-consequences} at $2p=\widehat p_n$, there is a
deterministic $C_*$ such that $\sqrt j\,\|(S_ne)_j\|_{\widehat p_n}\le C_*\widehat p_n^2$ for the coupling constructed in that theorem. Fix
\begin{equation}
 D\ge\max\{1,e^2C_*\}
 \label{B:eq:ratio-cost-constant}
\end{equation}
and use the family selected in Theorem~\ref{B:thm:measurable-ratio-optimizer} with this cost.

\begin{theorem}\label{B:thm:summable-ratio-cost}
The minimum in \eqref{B:eq:ratio-optimal-plan} satisfies
\begin{equation}
 \sup_{q\in\mathcal Q_n^\sharp}\mathfrak I_n(q)
 \le a_n:=2[Cn^{1-\lambda}r_n^{\lambda-2}\ell_n^3]^{p_0}
                      +ne^{-2\widehat p_n}.
 \label{B:eq:ratio-optimal-cost-bound}
\end{equation}
With $\nu=p_0(\lambda-1)>4$,
\begin{equation}
 a_n\le Cn^{-(\nu+1)/2}+n^{-15},\qquad
 \sum_{n\ge n_0}\sup_{q\in\mathcal Q_n^\sharp}\mathfrak I_n(q)<\infty.
 \label{B:eq:ratio-optimal-cost-summability}
\end{equation}
\end{theorem}
\begin{proof}
Fix the deterministic $D$ in \eqref{B:eq:ratio-cost-constant}. For sufficiently large $n$, $\ell_n\ge p_0+1$, so $P_n=\lceil\ell_n\rceil\le2\ell_n$ and
\[
 8\ell_n\le\widehat p_n
       \le8\ell_n+2p_0+2\le12\ell_n.
\]
Also $r_n/\ell_n\to\infty$, which gives the last condition in \eqref{B:eq:ratio-growing-cost-range}. All cost moments therefore lie in the range of
Theorem~\ref{B:thm:logarithmic-ratio-coupling} with $C_0=16$. Its single coupling can be used for every term of the cost. This coupling is a competitor at the fixed boundary $q$;
optimality will remove any need to select these competitors measurably.

Write $e=\eta-u$ and $s_n=n^{1-\lambda}r_n^{\lambda-2}$. For each integer $p_0\le p\le P_n$, the weighted consequence of that theorem gives
\[
 \mathbb E\left(\frac{nQ_n(e)}{Dp^4}\right)^p
 =\left(\frac{\|nQ_n(e)\|_p}{Dp^4}\right)^p
 \le\left[\frac{Cs_n\ell_n^2(p+\ell_n)}{Dp^3}\right]^p
 \le[Cs_n\ell_n^3]^p.
\]
Here $p+\ell_n\le3\ell_n$, $D\ge1$, and $p\ge p_0\ge2$. The constant in the final base is independent of $q,n,p$. Since $r_n^{\lambda-2}\le1$ and $\lambda>1$, the base
$b_n=Cs_n\ell_n^3$ tends to zero. Enlarge the deterministic starting size until $b_n\le1/2$. The finite sum is then bounded by
\[
 \sum_{p=p_0}^{P_n}b_n^p
   =b_n^{p_0}\frac{1-b_n^{P_n-p_0+1}}{1-b_n}
   \le 2b_n^{p_0}.
\]
The formula includes $P_n=p_0$.

For the other cost terms use the same coupling at moment order $p=\widehat p_n/2$. Since $4\ell_n\le p\le6\ell_n$, the second consequence in
\eqref{B:eq:ratio-coupling-consequences} implies
\[
 \sqrt j\,\|(S_ne)_j\|_{\widehat p_n}
 \le C\sqrt{\widehat p_n/2}\,\ell_n
                         \sqrt{\widehat p_n/2+\ell_n}
 \le C_*\widehat p_n^2.
\]
The constant $C_*$ can be chosen simultaneously for all $j,q,n$ in the indicated ranges. For each $j$, taking the $\widehat p_n$th power therefore bounds the expected value of the
corresponding summand by $(C_*/D)^{\widehat p_n}\le e^{-2\widehat p_n}$. Summing gives
\[
 \sum_{j\in J}\mathbb E
 \left(\frac{\sqrt j\,|(S_ne)_j|}{D\widehat p_n^2}\right)^{\widehat p_n}
 \le |J|e^{-2\widehat p_n}\le ne^{-2\widehat p_n}.
\]
Adding the two estimates proves that this coupling has expected cost at most $a_n$. Its first and second marginals are exactly $\Gamma_q$ and $\nu_q$ after ordering the coordinates
as $(u,\eta)$. The minimum is bounded by this competitor for every admissible $q$, which proves \eqref{B:eq:ratio-optimal-cost-bound} with a common constant.

Set $\nu=p_0(\lambda-1)>4$. The first term satisfies
\[
 2b_n^{p_0}
 \le Cn^{-\nu}\ell_n^{3p_0}
 \le Cn^{-\nu+ (\nu-1)/2}=Cn^{-(\nu+1)/2}.
\]
Indeed $\log\ell_n/\log n\to0$ gives $3p_0\log\ell_n\le(\nu-1)\log n/2$ for all large $n$. For the second term, $\widehat p_n\ge8\ell_n$ implies
\[
 ne^{-2\widehat p_n}\le n(\me+n)^{-16}\le n^{-15}.
\]
Both power series converge because $(\nu+1)/2>1$. At a fixed earlier size the admissible boundary set is bounded, its closure lies in the full nonnegative parameter space, and
$\mathfrak I_n$ is continuous there. Its supremum is finite. Thus any finitely many earlier sizes can be included by increasing the constant in
\eqref{B:eq:ratio-optimal-cost-summability}. Each nonnegative summand is bounded by the total cost of the selected measurable plan. For every integer $p_0\le p\le P_n$, this gives
\[
 \|nQ_n(\eta-u)\|_{L^p(\pi_{n,q}^{\rm rat})}
 \le Dp^4\mathfrak I_n(q)^{1/p}\le Dp^4a_n^{1/p}
                         \quad(q\in\mathcal Q_n^\sharp).
\]
Similarly, writing $\mathfrak m_n(e)=\max_{j\in J}\sqrt j\,|(S_ne)_j|$, the sum of the high powers of the individual partial sums yields
\[
 \|\mathfrak m_n(\eta-u)\|_{L^{\widehat p_n}(\pi_{n,q}^{\rm rat})}
 \le D\widehat p_n^2\mathfrak I_n(q)^{1/\widehat p_n}.
\]
The selected coupling minimizes the expected cost. The sharper moment estimates for the coupling constructed from the stationary diffusions therefore bound the expectation of each
nonnegative summand. The almost sure estimates below use these two moment inequalities.
\end{proof}

\subsubsection{Coupling conditional on the Gaussian vector}
\label{B:subsection:ratio-coupling-prescribed-gaussian}

The selected plans are now applied to the Gaussian vectors already coupled with Brownian motion. Additional independent uniforms sample first the ratios and then the remaining
scale variables; neither step changes the existing Brownian or exterior laws.

Let $(W,(E_n)_n,(G_n)_n)$ have the preceding joint distribution. At size $n$, $\lambda_n^{\rm ext}$ denotes the original marginal law of the complete exterior vector $E_n$. Suppose
$\mathcal L(G_n\mid E_n=e)=\gamma_{n,c_n(e)}$ on the full conditional interval. Put $q_n=c_n(E_n)$ and $u_n=L_{q_n}(G_n-m^{q_n}-s^{q_n})$. The boundary probabilities satisfy
$\sum_n\Pr(q_n\notin\mathcal Q_n^\sharp)<\infty$ by Theorem~\ref{B:thm:boundary-and-terminal-probabilities}. Adjoin independent uniforms, independent of the entire preceding space,
and use the maps of \eqref{B:eq:ratio-coefficient-reconstruction} to define
\begin{equation}
 \eta_n=\mathfrak A_n(q_n,u_n,U_n^{\rm rat}),\qquad
 s_n^{\rm rad}=\mathfrak R_n(q_n,\eta_n,U_n^{\rm rad}),\qquad
 Z_n=\mathfrak Z_{q_n}(s_n^{\rm rad},\eta_n).
 \label{B:eq:ratio-realized-reconstruction}
\end{equation}

\begin{proposition}\label{B:prop:ratio-coupling-realization}
The full coefficient vector $(E_n,Z_n)$ has the original joint distribution, and the preceding joint law of $(W,(E_n)_n,(G_n)_n)$ is preserved. For each $n$,
\begin{equation}
 \mathcal L(u_n,\eta_n\mid E_n=e)=\pi_{n,c_n(e)}^{\rm rat},\qquad
 \mathcal L(Z_n\mid E_n=e)=\mu_{n,c_n(e)}.
 \label{B:eq:ratio-realized-conditional-laws}
\end{equation}
If $E_n$ and $G_n$ are measurable from the Brownian nodes at size $n$ and an auxiliary sigma-field independent of $W$, the reconstructed coefficient vector has the same
measurability after adjoining the uniforms.
\end{proposition}
\begin{proof}
Let $(\Omega_0,\mathcal F_0,\mathbb P_0)$ carry the entire preceding array and Brownian path. Take its product with
\[
 \prod_{n\ge n_0}\bigl([0,1]\times[0,1]^{d_n}\bigr),
                       \qquad d_n=|J_n|,
\]
equipped with the product of Lebesgue probabilities. The new coordinates are $U_n^{\rm rat}$ and $U_n^{\rm rad}$. The projection onto $\Omega_0$ has law $\mathbb P_0$: integrating
any bounded $\mathcal F_0$-measurable function over the extra coordinates returns that function's original expectation. Thus the complete previous joint law, including its
dependence between sizes and its dependence on $W$, is preserved.

The new uniforms are independent of $\mathcal F_0$. For a bounded Borel function $H$ of $E_n,G_n,\eta_n,s_n^{\rm rad}$, their sampling identities give
\[
 \mathbb E[H(E_n,G_n,\eta_n,s_n^{\rm rad})\mid\mathcal F_0] 
 =\int K_n^{\rm rat}(q_n,u_n;\dd\eta)
          \int\mathfrak r_{q_n}(\dd s\mid\eta)
                                    H(E_n,G_n,\eta,s).
\]
The right side is Borel in the preceding random variables and bounded by $\|H\|_\infty$. Given $E_n=e$, the conditional Gaussian distribution gives
\[
 \mathbb E[F(u_n)\mid E_n=e]
   =\int F(L_{c_n(e)}(g-m^{c_n(e)}-s^{c_n(e)}))
                                      \,\gamma_{n,c_n(e)}(\dd g)
   =\int F(u)\,\Gamma_{c_n(e)}(\dd u).
\]
Apply the tower property to the previous conditional formula. The disintegration identity yields the conditional law of $(u_n,\eta_n)$ in
\eqref{B:eq:ratio-realized-conditional-laws}; Proposition~\ref{B:prop:conditional-coefficient-reconstruction} gives its other conditional law. The kernels on the right are defined
for every finite nonnegative configuration. They specify versions of these conditional probabilities; the conditional identities for the random exterior hold almost everywhere
under its marginal.

For a nonnegative Borel function $F$ on the full coefficient space, Tonelli and the identity for the original conditional distribution give
\[
 \mathbb EF(E_n,Z_n)
  =\int\lambda_n^{\rm ext}(\dd e)
             \int F(e,z)\,\mu_{n,c_n(e)}(\dd z)
  =\int F\,\dd\mathbb P_{n,V,\beta,a}.
\]
This uses the original exterior marginal and proves the prescribed joint distribution of the coefficient vector. Conditioning on all the previously constructed variables instead
gives the first displayed kernel, which retains the dependence on $u_n$. In particular, the reconstruction leaves the previously chosen full Gaussian vector unchanged.

To verify measurability and the independence of the auxiliary randomness, let $\mathcal U_0$ be the old auxiliary sigma-field and let $\mathcal V$ be generated by all the new
uniforms. Put $\mathcal U=\mathcal U_0\vee\mathcal V$. For bounded measurable functions $F_W,F_0,F_V$ of the Brownian path, the old auxiliary variables and the new uniforms,
respectively, independence gives
\[
 \mathbb E[F_WF_0F_V]
   =\mathbb E[F_WF_0]\mathbb E[F_V]
   =\mathbb E[F_W]\mathbb E[F_0]\mathbb E[F_V].
\]
A monotone class argument extends this identity from such products to every bounded $\mathcal U$-measurable function. Hence $\mathcal U$ is independent of the whole Brownian path.
At size $n$, the variables $q_n,G_n,u_n$ are measurable with respect to $\mathcal N_n\vee\mathcal U_0$ by hypothesis. Composing their Borel maps with the two samplers and
$\mathfrak Z_{q_n}$ proves $\mathcal N_n\vee\mathcal U$-measurability of the reconstructed coefficient vector. All uniforms are indexed by a countable set, so the same construction
serves every size simultaneously.
\end{proof}

Write $e_n=\eta_n-u_n$ and $\mathcal G_n=\{q_n\in\mathcal Q_n^\sharp\}$.

\begin{corollary}\label{B:cor:ratio-coupling-errors}\label{B:prop:simultaneous-ratio-error-estimates}
For the arrays in \eqref{B:eq:ratio-realized-reconstruction}, almost surely,
\begin{equation}
 \sum_{n\ge n_0}\mathfrak c_n(e_n)<\infty,\qquad
 \sum_{n\ge n_0}[nQ_n(e_n)]^{p_0}<\infty,\qquad nQ_n(e_n)\longrightarrow0.
 \label{B:eq:ratio-realized-energy-series}
\end{equation}
\manualstatementlabel{B:eq:simultaneous-ratio-error-series}{\theequation} Moreover, there is an almost-surely finite $n_0(\omega)$ such that, for every $n\ge n_0(\omega)$,
\begin{equation}
 \max_{j\in J_n}\sqrt j\,|(S_ne_n)_j|\le D\widehat p_n^2,\qquad
 \max_{j\in J_n}|(S_ne_n)_j|
          \le C\ell_n^2/\sqrt{r_n}\longrightarrow0.
 \label{B:eq:ratio-realized-partial-sum-maximum}
\end{equation}
\manualstatementlabel{B:eq:simultaneous-ratio-error-maximum}{\theequation}
\end{corollary}
\begin{proof}
Conditional integration against \eqref{B:eq:ratio-realized-conditional-laws} gives, for every nonnegative Borel $F$,
\begin{equation}
 \mathbb E[\mathbf1_{\mathcal G_n}F(q_n,u_n,\eta_n)]
 =\int_{c_n(e)\in\mathcal Q_n^\sharp}\lambda_n^{\rm ext}(\dd e)
       \int F(c_n(e),u,\eta)\,\pi_{n,c_n(e)}^{\rm rat}(\dd u,\dd\eta).
 \label{B:eq:ratio-cost-conditional-integration}
\end{equation}
Tonelli and Theorem~\ref{B:thm:summable-ratio-cost} imply
\[
 \mathbb E\sum_n\mathbf1_{\mathcal G_n}\mathfrak c_n(e_n)
 =\sum_n\mathbb E[\mathbf1_{\mathcal G_n}\mathfrak I_n(q_n)]
 \le\sum_na_n<\infty.
\]
The localized series is therefore finite almost surely. Borel--Cantelli makes $\mathcal G_n^c$ occur only finitely often; every sampled coordinate and every polynomial cost is
finite. Consequently $\sum_n\mathfrak c_n(e_n)<\infty$ almost surely. Keeping the $p_0$ term of the cost yields
\begin{equation}
 \sum_n[nQ_n(e_n)]^{p_0}
 \le(Dp_0^4)^{p_0}\sum_n\mathfrak c_n(e_n)<\infty.
 \label{B:eq:ratio-weighted-energy-series}
\end{equation}
In particular $nQ_n(e_n)\to0$. Choose $n_0(\omega)$ such that $\mathfrak c_n(e_n)\le1$ for $n\ge n_0(\omega)$. The terms involving individual partial sums then give
\[
 \max_{j\in J_n}\sqrt j\,|(S_ne_n)_j|\le D\widehat p_n^2,\qquad
 \max_{j\in J_n}|(S_ne_n)_j|
 \le144D\ell_n^2r_n^{-1/2}\longrightarrow0.
\]
Here $\widehat p_n\le12\ell_n$ for all $n$ beyond a deterministic $n_0$ and $r_n\asymp\ell_n^{A_0}$ with $A_0>20$. The first Borel--Cantelli lemma gives all these estimates on the
same event, for the joint construction across matrix sizes. Finally each error event forces its corresponding cost term to exceed the indicated threshold. Markov's inequality
applied to \eqref{B:eq:ratio-cost-conditional-integration} gives
\begin{equation}
 \begin{split}
 \Pr(\mathcal G_n\cap\{nQ_n(e_n)>\varepsilon\})
 &\le(Dp_0^4/\varepsilon)^{p_0}a_n,\\
 \Pr\left(\mathcal G_n\cap
   \{\max_j\sqrt j\,|(S_ne_n)_j|>eD\widehat p_n^2\}\right)
 &\le e^{-\widehat p_n}a_n.
 \end{split}
 \label{B:eq:ratio-realized-error-probabilities}
\end{equation}
\end{proof}

\subsubsection{Restriction of the partial sums}
\label{B:subsection:restricted-ratio-partial-sums}

Restriction to the selected interval leaves a constant contribution from the discarded upper indices. Keeping it explicitly preserves the full partial sums. With
$v_{n,j}=(S_ne_n)_j$ and $v_{n,R_n+1}=0$, this contribution is
$\gamma_n^{\rm buf}=v_{n,M_n+1}$.

\begin{corollary}\label{B:cor:restricted-ratio-partial-sums}
For the selected interval,
\begin{equation}
 \gamma_n^{\rm buf}+\sum_{k=j}^{M_n}e_{n,k}=v_{n,j},\qquad
 \sum_{k=i}^le_{n,k}=v_{n,i}-v_{n,l+1}.
 \label{B:eq:restricted-ratio-sum-identities}
\end{equation}
In particular,
\begin{equation}
 \sum_{j\in I^{\rm ret}}w_{n,j}
       \left|\gamma_n^{\rm buf}+\sum_{k=j}^{M_n}e_{n,k}\right|^2
       \le Q_n(e_n)=o_\omega(n^{-1}),
 \label{B:eq:restricted-ratio-energy}
\end{equation}
and the maximum of these partial sums tends to zero almost surely.
\end{corollary}
\begin{proof}
Extend the full partial sums by $v_{n,R_n+1}=0$. For $j\in I^{\rm ret}\subset\{m_n,\ldots,M_n\}$, split the finite sum at $M_n$ to obtain
\[
 v_{n,j}=\sum_{k=j}^{R_n}e_{n,k}
   =\sum_{k=j}^{M_n}e_{n,k}+\sum_{k=M_n+1}^{R_n}e_{n,k}
   =\sum_{k=j}^{M_n}e_{n,k}+\gamma_n^{\rm buf}.
\]
For $r_n\le i\le l\le R_n$, subtraction gives
\[
 v_{n,i}-v_{n,l+1}
   =\sum_{k=i}^{R_n}e_{n,k}-\sum_{k=l+1}^{R_n}e_{n,k}
   =\sum_{k=i}^{l}e_{n,k}.
\]
The second sum is empty when $l=R_n$. These are \eqref{B:eq:restricted-ratio-sum-identities}, including the terminal index. Nonnegativity of the weights gives
\[
 \sum_{j\in I^{\rm ret}}w_{n,j}
   \left|\gamma_n^{\rm buf}+\sum_{k=j}^{M_n}e_{n,k}\right|^2
   =\sum_{j\in I^{\rm ret}}w_{n,j}|v_{n,j}|^2
   \le\sum_{j\in J_n}w_{n,j}|v_{n,j}|^2=Q_n(e_n).
\]
Corollary~\ref{B:cor:ratio-coupling-errors} proves the required rate. It also bounds every shorter sum simultaneously by $2\max_{j\in J_n\cup\{R_n+1\}}|v_{n,j}|$, which tends to
zero. The maximum of the adjusted retained sums is at most
$\max_{j\in J_n}|v_{n,j}|$.

The scalar is itself one of the full partial sums. Since $M_n+1\asymp n$ and this index belongs to $J_n$,
\[
 |\gamma_n^{\rm buf}|=|v_{n,M_n+1}|
       \le D\widehat p_n^2/\sqrt{M_n+1}
       \le C\ell_n^2/\sqrt n\longrightarrow0
\]
for every $n\ge n_0(\omega)$, on an event of probability one. It is measurable from the same Brownian samples and auxiliary variables as $e_n$. With $\gamma=\gamma_n^{\rm buf}$,
replace the positive normalizing scalar by $C_n'=e^{-\gamma}C_n$. The two coefficient identities then give
\[
 (C_n')^{-1}p_{n,k}=p^B(z_{n,k})e^{\delta_{n,k}+\gamma},\qquad
 (C_n')^2q_{n,k}=Q_{n,k}^B
                         e^{-2(\delta_{n,k}+\gamma)+\rho_{n,k}}.
\]
The local error $\rho_{n,k}$ is unchanged, and each kernel summand has the same exponential factor because
\[
 (\delta_{n,i}+\gamma)+(\delta_{n,j}+\gamma)
            -2(\delta_{n,l}+\gamma)+\rho_{n,l}
       =\delta_{n,i}+\delta_{n,j}-2\delta_{n,l}+\rho_{n,l}.
\]
Equivalently, $(e^\gamma p_{n,i})(e^\gamma p_{n,j})(e^{-2\gamma}q_{n,l}) =p_{n,i}p_{n,j}q_{n,l}$ after the corresponding normalization. The ratio $C_n'/C_n=e^{-\gamma_n^{\rm buf}}$
is positive, measurable, and tends to one. This is the adjustment used in Section~\ref{B:sec:normalization-and-coefficient-estimates}.
\end{proof}

Section~\ref{B:sec:construction-of-the-coupling-for-all-matrix-sizes} places these kernels on one probability space for all matrix sizes. The coefficient and truncation estimates
are then combined with \eqref{B:eq:restricted-ratio-energy} in the proof of the general potential theorem.

\subsection{Construction of the coupling for all matrix sizes}
\label{B:sec:construction-of-the-coupling-for-all-matrix-sizes}

The conditional sampling maps of Section~\ref{B:sec:transport-of-logarithmic-ratios-and-construction-of-the-coupling} act on the Gaussian vectors from
Section~\ref{B:sec:gaussian-approximation-of-the-coefficient-field}. The exterior coordinates retain the joint distribution in Theorem~\ref{B:thm:original-exterior-realization}. We
realize these variables on one probability space and verify the full coefficient law at each size. We then establish the estimates needed for the normalization in
Section~\ref{B:sec:normalization-and-coefficient-estimates}.

\subsubsection{Parameters and the probability space}
\label{B:subsection:coupling-product-probability-space}

The separate conditional constructions must now be realized together. Their parameters are fixed before the rows are sampled, so that the resulting probability space and coupling
do not depend on which error moment is subsequently estimated.

Fix $V$, $a>0$, and $\beta>0$ as in Theorem~\ref{B:result:1-3}. Fix $1<\lambda<\min\{a+1,2\}$, $A_0>30$, and an integer $p_0\ge2$ with $p_0(\lambda-1)>4$. Use the intervals in
\eqref{B:eq:conditional-intervals}, with this value of $A_0$, and the notation $m_n=r_n+b_n$, $M_n=R_n-b_n$ of Section~\ref{B:subsection:conditional-gaussian-distribution}. Thus
$I_n=\{m_n,\ldots,M_n\}$ and $I_n^{\rm ret}=I_n\setminus\{m_n\}$. Write $d_n=|J_n|$. The nodes and the sigma-field $\mathcal N_n$ are those of
\eqref{B:eq:brownian-sampled-information}; the weights $w_{n,j}$, partial sum matrix $S_n$, and quadratic form $Q_n$ are those of \eqref{B:eq:ratio-partial-sum-cost}. The two
logarithmic interval lengths are denoted by $h^+_{n,k}$ and $h^-_{n,k}$ as in \eqref{B:eq:coupling-clocks}.

Fix the extension parameter used in Theorem~\ref{B:thm:logarithmic-ratio-coupling}, and then fix $D$ by \eqref{B:eq:ratio-cost-constant}. The selected plans $\pi_{n,q}^{\rm rat}$
use that value of $D$. Increase a deterministic $n_0$ so that, for $n\ge n_0$,
\begin{equation}
 2\le r_n<m_n<M_n<R_n<n,\qquad
 r_n<n/4,\quad r_n\ge8/\beta,\quad n\beta\ge1,
 \qquad64\log(e+n)\le c_{\rm rc}r_n,
 \label{B:eq:simultaneous-parameter-range}
\end{equation}
and so that the modified initial coordinates and their interaction neighborhood lie to the left of $J_n$. Such a choice exists since $r_n/n\to0$ and $\log(e+n)/r_n\to0$. The same
choice ensures that $n$ exceeds the thresholds required for the Gaussian covariance comparison and for \eqref{B:eq:ratio-growing-cost-range}. Indeed, with $x=\log n$,
\[
 \frac{r_n}{n}\le(x^{A_0}+1)e^{-x}\longrightarrow0,
 \qquad \frac{\log(e+n)}{r_n}\le
 \frac{x+\log(1+e)}{x^{A_0}}\longrightarrow0.
\]
All these choices are deterministic and precede the sampling. For $n<n_0$ the full coefficient vector is sampled directly from $\mathbb P_{n,V,\beta,a}$.

Let $\mathsf E_n$ be the space of exterior coordinates of $J_n$, in their original order, and let $\lambda_n^{\rm ext}$ be their original marginal. Retain the joint measure
$\Lambda^{\rm ext}$ on $\prod_{n\ge n_0}\mathsf E_n$ from Theorem~\ref{B:thm:original-exterior-realization}. Write $q_n=c_n(E_n)$ and $\mathscr E=\sigma(E_n:n\ge n_0)$. The
conditional coefficient kernel $\mu_{n,q}$, minimizer $m^q$, matrices $A_q,C_q,\Omega_q,L_q,K_q$, and Gaussian distributions $\gamma_{n,q},\Gamma_q$ retain their definitions from
\eqref{B:eq:conditional-matrix-data}, \eqref{B:eq:gaussian-coefficient-law}, and Section~\ref{B:subsection:gaussian-logarithmic-marginal}. The mean correction is $s^q$ from
\eqref{B:eq:conditional-laplace-correction}. All these parameter maps are defined for every finite nonnegative $q$. Their continuity and the positivity of the required matrices
follow from Proposition~\ref{B:prop:conditional-minimizer-dependence} and Lemma~\ref{B:lem:gaussian-ratio-schur-complement}.

Let $\mathbb W$ denote Wiener measure on the continuous real paths starting at zero, with the topology of uniform convergence on compact intervals. Let $\mathsf N$ denote the
standard normal distribution on $\mathbb R$. Take the product probability measure
\begin{equation}
 \mathbb P_0=\mathbb W\otimes\Lambda^{\rm ext}
 \otimes\bigotimes_{n\ge n_0}
       \left(\mathsf N^{\otimes(\mathbb Z\setminus I_n)}
                    \otimes N(0,I_{2d_n})\right)
 \otimes\bigotimes_{1\le n<n_0}\mathbb P_{n,V,\beta,a}.
 \label{B:eq:simultaneous-base-measure}
\end{equation}
Its coordinates are denoted by $W$, $(E_n)_n$, $(\xi^\perp_{n,j})_{j\notin I_n}$, $\omega_n$, and the initial rows. Let $\mathcal S$ be the collection of all these coordinates and
put $\mathscr B=\sigma(\mathcal S)$. Enlarge the space by
\begin{equation}
 \mathbb P=\mathbb P_0\otimes
       \bigotimes_{n\ge n_0}
              \left(\mathrm{Leb}_{(0,1)}\otimes
                                  \mathrm{Leb}_{(0,1)}^{\otimes d_n}\right).
 \label{B:eq:simultaneous-sampling-measure}
\end{equation}
The new coordinates are $U_n^{\rm rat}$ and $U_n^{\rm rad}$. Let $\mathcal V$ be the collection of all coordinates except $W$, and put $\mathcal U=\sigma(\mathcal V)$. Fubini's
theorem gives, for bounded Borel $f$ and $g$,
\[
 \mathbb E[f(W)g(\mathcal V)]
       =\int f\,\dd\mathbb W\int g\,\dd\mathcal L(\mathcal V).
\]
Hence $\mathcal U$ is independent of the entire path $W$, and $\mathscr E\subset\mathcal U$. All coordinate spaces and their countable products are standard Borel. The probability
kernels and sampling maps used below are the explicit versions proved in Section~\ref{B:subsection:ratio-conditional-probability-kernels}; see also \cite{kallenberg2021foundations}
(Lemma~4.22 and Theorem~8.5) for kernel randomization and disintegration.

Put $\mathcal G_n=\{q_n\in\mathcal Q_n^\sharp\}$ and retain $\gamma_\beta$ from \eqref{B:eq:coefficient-measure-comparison}. The marginal of $E_n$ under
\eqref{B:eq:simultaneous-base-measure} is $\lambda_n^{\rm ext}$. Thus Theorem~\ref{B:thm:boundary-and-terminal-probabilities} gives
\begin{equation}
 \mathbb P(\mathcal G_n^c)
       \le C_\beta e^{-c\gamma_\beta r_n^{3/4}},\qquad
 \sum_{n\ge n_0}\mathbb P(\mathcal G_n^c)<\infty.
 \label{B:eq:simultaneous-boundary-probabilities}
\end{equation}
Indeed $r_n^{3/4}/\log n\to\infty$, so the first bound is at most $n^{-3}$ for all sufficiently large $n$.

\subsubsection{Gaussian vectors constructed from Brownian increments}
\label{B:subsection:simultaneous-gaussian-inputs}

The Gaussian vectors use the sampled exterior configuration as their boundary parameter. They depend on Brownian motion through the prescribed increments, which give the required
coefficient law and measurability for the lower bound.

On \eqref{B:eq:simultaneous-base-measure}, apply the formulas \eqref{B:eq:coupling-white-noise}, \eqref{B:eq:coupling-frozen-field}, \eqref{B:eq:coupling-primitive},
\eqref{B:eq:coupling-exact-comparison}, and \eqref{B:eq:gaussian-regression-formula}, with $q=q_n$. We denote the resulting full coefficient vector by $G_n=G_n^{q_n}$ and set
\begin{equation}
 u_n^G=L_{q_n}(G_n-m^{q_n}-s^{q_n}).
 \label{B:eq:simultaneous-gaussian-ratios}
\end{equation}
Thus $u_n^G$ is the full vector denoted by $\eta_n^G$ in Section~\ref{B:sec:gaussian-approximation-of-the-coefficient-field}. The variables $\Xi_n$ and
$\chi_n^\nabla=\Xi_{n,M_n+1}$ use the definitions in \eqref{B:eq:coupling-primitive} and \eqref{B:eq:coupling-decomposition}.

\begin{lemma}\label{B:lem:simultaneous-gaussian-laws}
These maps are Borel and use the Brownian path only through $\mathcal N_n$. For every $n\ge n_0$,
\begin{equation}
 \mathcal L(G_n\mid\mathscr E)=\gamma_{n,q_n},\qquad
 \mathcal L(u_n^G\mid\mathscr E)=\Gamma_{q_n},\qquad
 (u_n^G)|_{I_n}=\eta_{n,I_n}^G.
 \label{B:eq:simultaneous-gaussian-laws}
\end{equation}
Conditionally on $\mathscr E$, any finite collection of the Gaussian fields, including fields at different sizes, is jointly Gaussian.
\end{lemma}
\begin{proof}
For fixed $n$, the intervals $[u_{n,j+1},u_{n,j}]$, $j\in I_n$, have disjoint interiors. The corresponding normalized Brownian increments are independent standard normals. They are
independent of the auxiliary normals in \eqref{B:eq:simultaneous-base-measure} and of $\mathscr E$. For each series in \eqref{B:eq:coupling-frozen-field} and
\eqref{B:eq:coupling-primitive}, the Fourier bounds give
\[
 \mathbb E\sum_{h\in\mathbb Z}
       n^{-1/2}\bigl(|a_{k/n}(h)|+|b_{k/n}(h)|\bigr)|\xi_{n,k-h}|
 \le C\sqrt{n/k}\,n^{-1/2}\sum_{h\in\mathbb Z}e^{-c|h|}<\infty.
\]
Tonelli therefore proves absolute convergence at every required $(n,k)$ almost surely. The countable Cauchy criterion makes the convergence set Borel; use zero on its complement.
The series also converge in $L^2$, by square summability of their coefficients.

Condition on the complete exterior sequence. The covariance identities of Proposition~\ref{B:prop:exact-gaussian-covariance-coupling} and
Proposition~\ref{B:prop:full-gaussian-regression} now hold with the fixed parameter $q_n$. The matrix $\Omega^\perp$ in \eqref{B:eq:gaussian-regression-formula}, evaluated at
$q=q_n$, is positive semidefinite. Equation~\eqref{B:eq:coupling-lift-verification} gives covariance $\Omega_{q_n}$ for $\zeta_n^{q_n}$ and its exact linear image on $I_n$. For
every $v\in\mathbb R^{2d_n}$, the conditional characteristic function is consequently
\[
 \mathbb E[e^{\ii v^{\mathsf T}G_n}\mid\mathscr E]
 =\exp\left\{\ii v^{\mathsf T}(m^{q_n}+s^{q_n})
                      -\tfrac12v^{\mathsf T}\Omega_{q_n}v\right\}.
\]
For $t\in\mathbb R^{d_n}$ this identity gives
\[
 \mathbb E[e^{\ii t^{\mathsf T}u_n^G}\mid\mathscr E]
 =e^{-t^{\mathsf T}L_{q_n}\Omega_{q_n}L_{q_n}^{\mathsf T}t/2}
 =e^{-t^{\mathsf T}K_{q_n}t/2}.
\]
This proves the second identity in \eqref{B:eq:simultaneous-gaussian-laws}; the regression identity gives its restriction to $I_n$.

For several sizes, fix the exterior sequence and any real linear combination $H$ of the centered fields. Truncate all convolution series to $|h|\le N$ and denote the resulting
variable by $H_N$. The finite vector of Brownian evaluations and auxiliary normals is Gaussian, so $H_N$ is centered Gaussian. Square summability gives
$\mathbb E[|H_N-H|^2\mid\mathscr E]\to0$ under this fixed conditional law. If $v_N=\mathbb E[H_N^2\mid\mathscr E]$, then $v_N\to v=\mathbb E[H^2\mid\mathscr E]$ and
\[
 \left|\mathbb E[e^{\ii tH}-e^{\ii tH_N}\mid\mathscr E]\right|
 \le |t|\mathbb E[|H-H_N|^2\mid\mathscr E]^{1/2},
 \qquad \mathbb E[e^{\ii tH}\mid\mathscr E]=e^{-t^2v/2}.
\]
Thus every such finite collection is conditionally jointly Gaussian. This is $L^2$ convergence for each fixed exterior sequence on the stated full probability event, with its
conditional means fixed. For two Brownian increments on intervals $[a,b]$ and $[c,d]$, their covariance is $|[a,b]\cap[c,d]|$. This also specifies the dependence between fields at
different matrix sizes.

Matrix inversion on positive definite matrices and the positive semidefinite square root are continuous in finite dimension. Together with the Borel series limits, this proves
measurability of the maps. In \eqref{B:eq:coupling-white-noise} the Brownian coordinates occur at the nodes of size $n$; all remaining variables are $\mathcal U$-measurable.
\end{proof}

For later use, on $\mathcal G_n$ the full covariance obeys
\begin{equation}
 |\operatorname{Cov}(u_{n,k}^G,u_{n,l}^G\mid\mathscr E)|
                \le \frac{Ce^{-c|k-l|}}{\sqrt{kl}},\qquad k,l\in J_n.
 \label{B:eq:simultaneous-ratio-covariance}
\end{equation}
Indeed, with $\epsilon_x=-1$ and $\epsilon_y=1$, its entry is
\[
 \frac1{n\beta}\sum_{\sigma,\tau\in\{x,y\}}
 \frac{\epsilon_\sigma\epsilon_\tau
                (C_{q_n})_{\sigma k,\tau l}}
      {m^{q_n}_{\sigma,k}m^{q_n}_{\tau,l}}.
\]
The bounds $m^q_{\sigma,k}\asymp\sqrt{k/n}$ and $|(C_q)_{\sigma k,\tau l}|\le Ce^{-c|k-l|}$ from \eqref{B:eq:interior-hessian-bounds} give
\eqref{B:eq:simultaneous-ratio-covariance} by summing the four terms.

\subsubsection{Conditional sampling of the coefficients}
\label{B:subsection:simultaneous-coefficient-sampling}

The positive coefficients are recovered in two stages. Sampling the ratios first preserves their coupling with the Gaussian variables; sampling the remaining scale variables from
their conditional law then restores the full coefficient distribution.

Use the previously fixed maps $\mathfrak A_n$, $\mathfrak R_n$, and $\mathfrak Z_q$ from \eqref{B:eq:ratio-sampler-kernel}, Lemma~\ref{B:lem:joint-logarithmic-scale-sampling}, and
\eqref{B:eq:ratio-coordinate-map}. On \eqref{B:eq:simultaneous-sampling-measure} set
\begin{equation}
 \begin{split}
 \eta_n&=\mathfrak A_n(q_n,u_n^G,U_n^{\rm rat}),\\
 s_n^{\rm rad}&=\mathfrak R_n(q_n,\eta_n,U_n^{\rm rad}),\\
 Z_n&=\mathfrak Z_{q_n}(s_n^{\rm rad},\eta_n),\qquad
 \mathbf Z_n=(E_n,Z_n).
 \end{split}
 \label{B:eq:simultaneous-coefficient-maps}
\end{equation}
The last expression places the coordinates of the two vectors in the original coefficient order. In particular $X_{n,n}$ is an exterior coordinate. The variables $s_n^{\rm rad}$
are the logarithmic coordinates $s$ of \eqref{B:eq:ratio-coordinate-map}; their conditional distribution is the joint kernel $\mathfrak r_q$. All their sampled values are finite,
and all coordinates of $\mathbf Z_n$ are positive.

\begin{lemma}\label{B:lem:random-boundary-null-sections}
Let $D_n$ be a Borel set of pairs $(q,u)$ with $\Gamma_q\{u:(q,u)\in D_n\}=0$ for every $q$. Then
\begin{equation}
 \mathbb P\{(q_n,u_n^G)\in D_n\text{ for some }n\ge n_0\}=0.
 \label{B:eq:random-boundary-null-sections}
\end{equation}
\end{lemma}
\begin{proof}
The indicator of $D_n$ is Borel, so integration against the Borel Gaussian kernel is measurable. Equation~\eqref{B:eq:simultaneous-gaussian-laws} and the assumed null sections
give, for every nonnegative Borel $H$,
\[
 \mathbb E[H((E_m)_m)\mathbf1_{D_n}(q_n,u_n^G)]
 =\int H(\mathbf e)\,
 \Gamma_{c_n(e_n)}\{u:(c_n(e_n),u)\in D_n\}
       \,\Lambda^{\rm ext}(\dd\mathbf e)=0.
\]
Taking $H=1$ and using countable subadditivity yields
\[
 \mathbb P\left(\bigcup_{n\ge n_0}
          \{(q_n,u_n^G)\in D_n\}\right)=0,
\]
which is \eqref{B:eq:random-boundary-null-sections}.
\end{proof}

This lemma applies to the Borel exceptional set in the parameterwise disintegration of Lemma~\ref{B:lem:parameterized-ratio-disintegration}. The sampling maps are defined there by
their fixed defaults. Consequently their parameterwise identities may be evaluated at $(q_n,u_n^G)$ simultaneously for all sizes.

For each $(q,u)$ define the probability kernel
\begin{equation}
 \mathsf R_n(q,u;\dd\eta,\dd s)
       =K_n^{\rm rat}(q,u;\dd\eta)\mathfrak r_q(\dd s\mid\eta).
 \label{B:eq:simultaneous-sampling-kernel}
\end{equation}
It is Borel: for a bounded Borel function of $(\eta,s)$, integrate first in $s$, then in $\eta$, and use the measurability of the two kernels proved in
Section~\ref{B:sec:transport-of-logarithmic-ratios-and-construction-of-the-coupling}. Tonelli gives total mass one at every $(q,u)$.

\subsubsection{Conditional joint distributions}
\label{B:subsection:simultaneous-conditional-distributions}

The conditional distribution below verifies that the joint coefficient distribution is preserved at every matrix size. The conditional product distribution holds when all the
previously constructed variables are fixed. Conditioning only on the exterior sequence does not make the coefficient vectors independent across matrix sizes, since they share
Brownian observations.

\begin{proposition}\label{B:prop:simultaneous-conditional-product}
For every finite set $F\subset\{n_0,n_0+1,\ldots\}$ and every nonnegative Borel function $\Psi$,
\begin{equation}
 \mathbb E\left[\Psi\bigl(\mathcal S,
                    (\eta_n,s_n^{\rm rad})_{n\in F}\bigr)\mid\mathscr B\right] 
 =\int\Psi\bigl(\mathcal S,(\eta_n,s_n)_{n\in F}\bigr)
              \prod_{n\in F}\mathsf R_n(q_n,u_n^G;\dd\eta_n,\dd s_n).
 \label{B:eq:reconstruction-conditional-product}
\end{equation}
The marginal distribution of $\mathcal S$ is $\mathbb P_0$.
\end{proposition}
\begin{proof}
For Borel sets $A,C\subset\mathbb R^{d_n}$ and fixed $(q,u)$, the laws of the two sampling maps give
\[
 \begin{split}
 &\int_0^1\!\int_{(0,1)^{d_n}}
 \mathbf1_A(\mathfrak A_n(q,u,v))
 \mathbf1_C\bigl(\mathfrak R_n(q,\mathfrak A_n(q,u,v),w)\bigr)
                                                  \,\dd w\,\dd v\\
 &=\int_A K_n^{\rm rat}(q,u;\dd\eta)\mathfrak r_q(C\mid\eta).
 \end{split}
\]
Independence of the uniform coordinates at different sizes now proves \eqref{B:eq:reconstruction-conditional-product} when
\[
 \Psi=\mathbf1_H(\mathcal S)\prod_{n\in F}
      \mathbf1_{A_n}(\eta_n)\mathbf1_{C_n}(s_n).
\]
More precisely, multiplying either side by $\mathbf1_{H'}(\mathcal S)$ and integrating under \eqref{B:eq:simultaneous-sampling-measure} gives identical integrals by Fubini. This is
the defining property of the conditional expectation.

The rectangles just used form a generating pi-system. For bounded functions, the class on which the equality holds is a vector space, contains constants, and is closed under
bounded monotone limits. The monotone class theorem therefore gives all bounded Borel functions. Truncation and monotone convergence give all nonnegative functions. For a bounded
$\mathscr B$-measurable multiplier $H$ and a bounded Borel function $\Psi$, the resulting identity is
\[
 \mathbb E[H(\mathcal S)\Psi(\mathcal S,(\eta_n,s_n^{\rm rad})_{n\in F})]
 =\int H(s)\!\int\Psi(s,(\eta_n,t_n)_{n\in F})
       \prod_{n\in F}\mathsf R_n(q_n(s),u_n^G(s);\dd\eta_n,\dd t_n)
       \,\mathbb P_0(\dd s).
\]
It identifies the conditional distribution on every finite cylinder. Two probability measures on the countable product which agree on these cylinders agree on its Borel sigma-field
by the pi-lambda theorem. The conditional law of all sampled pairs is therefore the product of these kernels. Taking $\Psi=1$ gives $\mathbb EH(\mathcal S)=\int H\,\dd\mathbb P_0$.
Under this conditioning the Gaussian vector is fixed: $\mathcal L(G_n\mid\mathscr B)=\delta_{G_n}$. Its Gaussian conditional law in \eqref{B:eq:simultaneous-gaussian-laws} is
obtained on conditioning on the smaller field $\mathscr E$. For a finite set of sizes $F$, denote the conditional joint Gaussian distribution of $(G_n)_{n\in F}$ given $\mathscr E$
by $\gamma^{\rm joint}_{F,\mathscr E}$. Taking the conditional expectation of the product formula gives, for bounded Borel $f_n$,
\[
 \mathbb E\!\left[\prod_{n\in F}f_n(Z_n)\,\middle|\,\mathscr E\right]
 =\int\gamma^{\rm joint}_{F,\mathscr E}(\dd(g_n)_{n\in F})
        \prod_{n\in F}\left[\int
        f_n(\mathfrak Z_{q_n}(s,\eta))
        \mathsf R_n(q_n,u_{q_n}(g_n);\dd\eta,\dd s)\right],
\]
where $u_q(g)=L_q(g-m^q-s^q)$. Because the Brownian observation intervals overlap across sizes, the conditional Gaussian measure in this identity need not be a product measure.
Conditional independence of the new sampling variables given $\mathscr B$ therefore does not imply conditional independence of the completed rows given $\mathscr E$.
\end{proof}

\begin{theorem}\label{B:thm:simultaneous-coefficient-construction}
The positive arrays in \eqref{B:eq:simultaneous-coefficient-maps}, together with the initial rows, satisfy
\begin{equation}
 \mathcal L(\mathbf Z_n)=\mathbb P_{n,V,\beta,a},\qquad
 \mathbf Z_n\text{ is }\mathcal N_n\vee\mathcal U\text{-measurable}
                         \quad(n\ge1).
 \label{B:eq:simultaneous-exact-arrays}
\end{equation}
For $n\ge n_0$ their conditional distributions are
\begin{equation}
 \mathcal L(u_n^G,\eta_n\mid\mathscr E)=\pi_{n,q_n}^{\rm rat},
 \qquad \mathcal L(Z_n\mid\mathscr E)=\mu_{n,q_n}.
 \label{B:eq:simultaneous-selected-laws}
\end{equation}
\end{theorem}
\begin{proof}
Apply \eqref{B:eq:reconstruction-conditional-product} at one size and then condition on $\mathscr E\subset\mathscr B$. By \eqref{B:eq:simultaneous-gaussian-laws}, for every
nonnegative Borel $\Phi$,
\begin{equation}
 \begin{split}
 &\mathbb E[\Phi(G_n,\eta_n,s_n^{\rm rad})\mid\mathscr E]\\
 &=\int\gamma_{n,q_n}(\dd g)
   \int K_n^{\rm rat}\bigl(q_n,L_{q_n}(g-m^{q_n}-s^{q_n});\dd\eta\bigr)
   \int\mathfrak r_{q_n}(\dd s\mid\eta)\Phi(g,\eta,s).
 \end{split}
 \label{B:eq:reconstruction-conditional-row}
\end{equation}
The image of $\gamma_{n,q}$ under $g\mapsto L_q(g-m^q-s^q)$ is $\Gamma_q$. Equation~\eqref{B:eq:ratio-parameterized-disintegration} therefore identifies the $(u_n^G,\eta_n)$
marginal as the selected plan $\pi_{n,q_n}^{\rm rat}$.

Its second marginal is $\nu_q$. The density formulas \eqref{B:eq:ratio-marginal-and-conditional} imply, for a nonnegative Borel $F$,
\[
 \begin{split}
 &\int\nu_q(\dd\eta)\int\mathfrak r_q(\dd s\mid\eta)
                               F(\mathfrak Z_q(s,\eta))\\
 &=\int\frac{I_q(\eta)}{\mathcal Z_q^{\rm rat}}\,\dd\eta
            \int\frac{g_q(s,\eta)}{I_q(\eta)}
                               F(\mathfrak Z_q(s,\eta))\,\dd s\\
 &=\frac1{\mathcal Z_q^{\rm rat}}
       \int g_q(s,\eta)F(\mathfrak Z_q(s,\eta))\,\dd s\,\dd\eta
   =\int F(z)\,\mu_{n,q}(\dd z).
 \end{split}
\]
For the final equality, retain the normalizers from Proposition~\ref{B:prop:logarithmic-ratio-density}:
\[
 e^{-\kappa H^q(\mathfrak Z_q(s,\eta))}
     |\det D\mathfrak Z_q(s,\eta)|=a_q^0g_q(s,\eta),
 \qquad \mathcal Z_q^{\rm coef}=a_q^0\mathcal Z_q^{\rm rat}.
\]
Substitution in the preceding integral gives $(\mathcal Z_q^{\rm coef})^{-1}\int e^{-\kappa H^q(z)}F(z)\,\dd z$. For a nonnegative $F$ this follows by applying the change of
variables to $F\wedge L$ and then using monotone convergence. Thus it is the original conditional law for every finite boundary configuration. This proves the second conditional
law and also the identity
\begin{equation}
 \eta_{n,k}=\log\frac{(Z_n)_{y,k}}{m^{q_n}_{y,k}}
                -\log\frac{(Z_n)_{x,k}}{m^{q_n}_{x,k}},\qquad k\in J_n.
 \label{B:eq:reconstructed-logarithmic-ratios}
\end{equation}
The kernels on the right of \eqref{B:eq:simultaneous-selected-laws} depend on the exterior only through $E_n$. The tower property gives the same conditional laws given $E_n$.

For any nonnegative Borel function $F(e,z)$ on a full coefficient row,
\begin{equation}
 \begin{split}
 \mathbb EF(E_n,Z_n)
 &=\int\Lambda^{\rm ext}(\dd\mathbf e)
        \int\mu_{n,c_n(e_n)}(\dd z)F(e_n,z)\\
 &=\int\lambda_n^{\rm ext}(\dd e)
        \int\mu_{n,c_n(e)}(\dd z)F(e,z)
   =\int F\,\dd\mathbb P_{n,V,\beta,a}.
 \end{split}
 \label{B:eq:reconstructed-full-row-law}
\end{equation}
The last equality is the original exterior disintegration in \eqref{B:eq:exterior-disintegration-identity}. For example, if $A$ is a Borel set of initial exterior coordinates and
$C$ is a Borel set of exterior coordinates near the right boundary, the same identity applied to $F(e,z)=\mathbf1_A(e)\mathbf1_C(e)\Phi(z)$ gives
\[
 \mathbb E[\mathbf1_A(E_n)\mathbf1_C(E_n)\Phi(Z_n)]
 =\int_{A\cap C}\!\left(\int\Phi(z)\,\mu_{n,c_n(e)}(\dd z)\right)
                           \lambda_n^{\rm ext}(\dd e).
\]
The exterior marginal specifies their joint distribution, including that of $X_{n,n}$. At $n<n_0$ the equality of laws is part of \eqref{B:eq:simultaneous-base-measure}.

The variables $q_n$ and the sampling uniforms belong to $\mathcal U$. Lemma~\ref{B:lem:simultaneous-gaussian-laws} gives the required $\mathcal N_n\vee\mathcal U$-measurability of
$u_n^G$. Composing the Borel maps in \eqref{B:eq:simultaneous-coefficient-maps} proves \eqref{B:eq:simultaneous-exact-arrays}. The initial rows belong to $\mathcal U$ as well. A
countable family of rational rectangles determines each of the conditional distributions; intersecting the corresponding full measure sets over $n$ gives all these conditional
distributions on one probability space.
\end{proof}

The first conditional law in \eqref{B:eq:simultaneous-selected-laws} gives the following integration rule for any nonnegative Borel $H$:
\begin{equation}
 \mathbb E[\mathbf1_{\mathcal G_n}H(q_n,u_n^G,\eta_n)]
 =\int_{c_n(e)\in\mathcal Q_n^\sharp}\lambda_n^{\rm ext}(\dd e)
       \int H(c_n(e),u,\eta)\,\pi_{n,c_n(e)}^{\rm rat}(\dd u,\dd\eta).
 \label{B:eq:selected-ratio-conditional-integration}
\end{equation}
Here $\mathcal G_n\in\mathscr E$, so the indicator remains outside the conditional integral. Equation~\eqref{B:eq:reconstruction-conditional-row} specifies the conditional
distribution when the full Gaussian vector is included among the observations.

\subsubsection{Brownian bridges and the sampled information}
\label{B:subsection:simultaneous-brownian-information}

For each fixed size, the construction uses Brownian samples but not the bridges between them. The independence statement below is made for that size's observation sigma-field; it
does not imply independence after conditioning on the observations from all sizes.

For each $n\ge2$, let $\mathfrak b_{n,k}$ be the Brownian bridge on $[u_{n,k},u_{n,k-1}]$, $2\le k\le n$, obtained by subtracting linear interpolation between its two Brownian
values, as in Lemma~\ref{B:lem:conditional-bridge-law}. The bridge family is empty at $n=1$.

\begin{lemma}\label{B:lem:simultaneous-bridge-independence}
For each fixed $n$, these bridges are mutually independent and independent of $\mathcal N_n\vee\mathcal U$. In particular, their conditional distribution given the constructed row
is their original product distribution.
\end{lemma}
\begin{proof}
At a fixed size $n$, Lemma~\ref{B:lem:conditional-bridge-law}, applied to the same deterministic times, identifies the Brownian bridges as independent centered Gaussian paths with
covariance $\min(s,t)-st/(b-a)$ on an interval $[a,b]$. Its covariance calculation applies to every interval of this partition: for $v\le a$ or $v\ge b$,
\[
 \operatorname{Cov}\left(W(a+s)-W(a)
       -\frac{s}{b-a}(W(b)-W(a)),W(v)\right)=0.
\]
Different intervals have zero cross covariance. Joint Gaussianity and continuity therefore give independence of the full bridge paths from $\mathcal N_n$.

The auxiliary field $\mathcal U$ is independent of the entire Brownian path. For bounded variables $F,G,H$ measurable respectively with respect to the bridges, $\mathcal N_n$, and
$\mathcal U$,
\[
 \mathbb E[FGH]=\mathbb E[FG]\mathbb EH
 =\mathbb EF\,\mathbb EG\,\mathbb EH=\mathbb EF\,\mathbb E[GH].
\]
The monotone class theorem extends the identity to $\mathcal N_n\vee\mathcal U$. Since $\sigma(\mathbf Z_n)\subset\mathcal N_n\vee\mathcal U$, conditioning once more proves the
conditional independence for the constructed row.

These conditioning statements concern one size. Indeed, $\{-\log\theta(k/n):1\le k\le n,\ n\ge2\}$ is dense in $[0,\infty)$, and Brownian continuity gives, up to completion,
\[
 \bigvee_{n\ge2}\mathcal N_n=\sigma(W(t):t\ge0).
\]
Every bridge is determined by this larger field.
\end{proof}

\subsubsection{Partial sums of the error in the logarithmic ratios}
\label{B:subsection:simultaneous-ratio-error-estimates}

The transport cost was chosen to control partial sums of the logarithmic ratio errors. After sampling, its conditional estimates yield pathwise bounds for those sums, with the
failures of the boundary restrictions accounted for separately.

Use the errors and partial sums of Section~\ref{B:subsection:restricted-ratio-partial-sums}, now with $u_n=u_n^G$:
\begin{equation}
 e_{n,k}=\eta_{n,k}-u_{n,k}^G,\qquad
 v_{n,j}=(S_ne_n)_j,\qquad v_{n,R_n+1}=0.
 \label{B:eq:simultaneous-ratio-error}
\end{equation}
The cost $\mathfrak c_n$ and exponents $P_n,\widehat p_n$ use the definitions \eqref{B:eq:ratio-selected-cost} and \eqref{B:eq:ratio-cost-exponents}. Write $\ell_n=\log(e+n)$.

Corollary~\ref{B:cor:ratio-coupling-errors} applies with $u_n=u_n^G$. Indeed, \eqref{B:eq:simultaneous-selected-laws} supplies the same conditional coupling
$\pi_{n,q_n}^{\rm rat}$, and the boundary event $\mathcal G_n$ has summable failure probabilities. Each reconstructed coordinate and its cost is finite, including the finitely many
inadmissible sizes. Thus \eqref{B:eq:ratio-realized-energy-series} and \eqref{B:eq:ratio-realized-partial-sum-maximum} hold for these arrays. The estimate therefore holds on the
prescribed probability space for the whole sequence, with its given dependence between sizes.

For later probability bounds, the conditional integration identity \eqref{B:eq:selected-ratio-conditional-integration} gives
\begin{equation}
 \mathbb E\sum_{n\ge n_0}\mathbf1_{\mathcal G_n}\mathfrak c_n(e_n)
 \le\sum_{n\ge n_0}a_n<\infty.
 \label{B:eq:simultaneous-integrated-cost}
\end{equation}
Adding $\mathbb P(\mathcal G_n^c)$ to the two estimates in \eqref{B:eq:ratio-realized-error-probabilities} gives, for every fixed $\varepsilon>0$,
\begin{equation}
 \begin{split}
 \mathbb P\{nQ_n(e_n)>\varepsilon\}
 &\le\mathbb P(\mathcal G_n^c)
                +(Dp_0^4/\varepsilon)^{p_0}a_n,\\
 \mathbb P\{\max_j\sqrt j\,|v_{n,j}|>eD\widehat p_n^2\}
 &\le\mathbb P(\mathcal G_n^c)+e^{-\widehat p_n}a_n.
 \end{split}
 \label{B:eq:simultaneous-ratio-error-probabilities}
\end{equation}
There is also a bound for the whole tail of the cost series. For $L>0$ and $N\ge n_0$,
\[
 \mathbb P\!\left\{\sum_{n\ge N}\mathfrak c_n(e_n)>L\right\}
 \le\sum_{n\ge N}\mathbb P(\mathcal G_n^c)
                         +L^{-1}\sum_{n\ge N}a_n.
\]
Outside the union of the boundary failure events the series equals its localized version; Markov's inequality and Tonelli bound that version by the second term. Both tails tend to
zero. The argument uses moment bounds only on admissible boundary data and allows the prescribed dependence between sizes.

\subsubsection{Logarithms of the reconstructed coefficients}
\label{B:subsection:simultaneous-local-logarithms}

The ratio estimate must be accompanied by control of the individual reconstructed coefficients. Their logarithms are centered using the conditional mean correction computed from
the original coefficient distribution.

Use the coordinate logarithms $\zeta_{n,\alpha}=\log((Z_n)_\alpha/m^{q_n}_\alpha)$ from \eqref{B:eq:conditional-logarithms}. For a fixed deterministic weight array
$0\le a_{n,\alpha}\le C_0n^{-\lambda}k_\alpha^{\lambda-1}$, let $\mathcal M_n$ be the corresponding energy of \eqref{B:eq:conditional-energy-moments}. Define
\begin{equation}
 b_n(q)=\mathbf1_{\{q\in\mathcal Q_n^\sharp\}}
    \int\sum_{\alpha\in\mathscr S(J_n)}a_{n,\alpha}
                   |\log(z_\alpha/m^q_\alpha)|^2\,\mu_{n,q}(\dd z).
 \label{B:eq:simultaneous-local-energy-mean}
\end{equation}
The integral is finite for each fixed $q$: the majorant in \eqref{B:eq:ratio-coordinate-majorant} is a Gaussian times positive coordinate powers, and
$\int_0^1 t^c|\log t|^2\,\dd t<\infty$ for $c>-1$. The constants $m^q_\alpha$ are positive. The integral is Borel: truncate the nonnegative integrand at an integer, integrate
against the Borel kernel, and let the integer increase. On admissible data, $0\le b_n(q)\le C/n$ by \eqref{B:eq:conditional-energy-moments}.

\begin{proposition}\label{B:prop:simultaneous-local-logarithms}
For the reconstructed coefficients, almost surely,
\begin{equation}
 n[\mathcal M_n-b_n(q_n)]\longrightarrow0,\qquad
 \limsup_n n\mathcal M_n\le C_{\beta,\lambda,C_0}.
 \label{B:eq:simultaneous-local-energy-limit}
\end{equation}
There is a deterministic $C<\infty$ and, almost surely, a finite $n_0(\omega)$ such that for every $n\ge n_0(\omega)$,
\begin{equation}
 \max_{\alpha\in\mathscr S(J_n)}\sqrt{k_\alpha}\,|\zeta_{n,\alpha}|
                                  \le C\sqrt{\ell_n}.
 \label{B:eq:simultaneous-local-maximum}
\end{equation}
\end{proposition}
\begin{proof}
The conditional coefficient law is \eqref{B:eq:simultaneous-selected-laws}. Let $\delta_\lambda=\tfrac12\min\{\lambda-1,1/2\}$ and choose a fixed integer $p$ with
$p\delta_\lambda>2$. For any $\varepsilon>0$, conditioning on $\mathscr E$ and applying \eqref{B:eq:conditional-centered-energy-probability} give
\[
 \begin{split}
 &\mathbb P\{n|\mathcal M_n-b_n(q_n)|>\varepsilon\}\\
 &\quad\le\mathbb P(\mathcal G_n^c)
 +\int_{c_n(e)\in\mathcal Q_n^\sharp}\lambda_n^{\rm ext}(\dd e)
       \mu_{n,c_n(e)}\{n|\mathcal M_n-\mathbb E_{c_n(e)}\mathcal M_n|
                                                          >\varepsilon\}\\
 &\quad\le\mathbb P(\mathcal G_n^c)
                   +(C_\lambda p^{3/2}/\varepsilon)^p n^{-p\delta_\lambda}.
 \end{split}
\]
The integer $p$ is fixed. Since $r_n\to\infty$, it satisfies $p\le c_{\rm rc}r_n$ for all sufficiently large $n$, as required in
\eqref{B:eq:conditional-centered-energy-probability}. The sum of the last display over $n$ is finite. For each $h\ge1$, discard the null limsup event at $\varepsilon=1/h$ and
intersect the remaining events over $h$. On that intersection, for each $h$ the quantity $n|\mathcal M_n-b_n(q_n)|$ is at most $1/h$ for every $n\ge n_0(\omega,h)$, proving its
convergence to zero. The deterministic mean bound then gives
\[
 \limsup_n n\mathcal M_n
 \le\limsup_n nb_n(q_n)
       +\limsup_n n|\mathcal M_n-b_n(q_n)|\le C.
\]
All integrations used the original conditional coefficient law in \eqref{B:eq:simultaneous-selected-laws}.

The maximum estimate in \eqref{B:eq:conditional-energy-and-maximum-probabilities}, with exponent $3$ in the probability bound, gives
\[
 \mathbb P\{\max_\alpha\sqrt{k_\alpha}|\zeta_{n,\alpha}|
                                      >C\sqrt{\ell_n}\}
 \le\mathbb P(\mathcal G_n^c)+C n^{-3}.
\]
Borel--Cantelli proves \eqref{B:eq:simultaneous-local-maximum}. In particular $\max_\alpha|\zeta_{n,\alpha}| \le C\sqrt{\ell_n/r_n}\to0$ almost surely.
\end{proof}

\subsubsection{Gaussian partial sums on the same space}
\label{B:subsection:simultaneous-gaussian-estimates}

Retain $B_{n,k}$ and $V_{n,j}^G$ from \eqref{B:eq:gaussian-increment-errors} and \eqref{B:eq:coupling-decomposition}. With the notation \eqref{B:eq:simultaneous-gaussian-ratios},
the latter identity reads
\[
 V_{n,j}^G=\Xi_{n,M_n+1}
                  +\sum_{k=j}^{M_n}(u_{n,k}^G-B_{n,k}),\qquad j\in I_n.
\]
It holds for the fields on \eqref{B:eq:simultaneous-base-measure}. Proposition~\ref{B:prop:simultaneous-conditional-product} preserves that whole base distribution when the
coefficients are sampled.

\begin{proposition}\label{B:prop:simultaneous-gaussian-estimates}
There are deterministic $K,C<\infty$ and, almost surely, a finite $n_0(\omega)$ such that for every $n\ge n_0(\omega)$,
\begin{equation}
 \sum_{j\in I_n}w_{n,j}|V_{n,j}^G|^2\le K/n,\qquad
 \max_{j\in I_n}\sqrt j\,|V_{n,j}^G|\le C\sqrt{\ell_n}.
 \label{B:eq:simultaneous-gaussian-estimates}
\end{equation}
\end{proposition}
\begin{proof}
Conditionally on $\mathscr E$, the vector $(\sqrt{w_{n,j}}V_{n,j}^G)_{j\in I_n}$ is centered Gaussian by Lemma~\ref{B:lem:simultaneous-gaussian-laws}. On $\mathcal G_n$,
Proposition~\ref{B:prop:weighted-gaussian-covariances} bounds its covariance $\Sigma_n$ by
\[
 \operatorname{Tr}\Sigma_n\le T/n,\qquad
 \|\Sigma_n\|_{\rm op}\le C n^{-\lambda}m_n^{\lambda-2}.
\]
Fix $K>2T$. Apply Lemma~\ref{B:lem:gaussian-quadratic-form-tail} under the conditional law with $S=\operatorname{Tr}\Sigma_n\le T/n$, $L=Cn^{-\lambda}m_n^{\lambda-2}$, and $v=K/n$.
Its exponent satisfies
\[
 \frac{v-2S}{4L}
 \ge \frac{K-2T}{4C}\,n^{\lambda-1}m_n^{2-\lambda}.
\]
Thus
\begin{equation}
 \mathbb P\left(\sum_{j\in I_n}w_{n,j}|V_{n,j}^G|^2>K/n
                                  \,\middle|\,\mathscr E\right)
 \le e^{-c n^{\lambda-1}m_n^{2-\lambda}}
                       \quad\text{on }\mathcal G_n.
 \label{B:eq:simultaneous-gaussian-energy-probability}
\end{equation}
Integration adds at most $\mathbb P(\mathcal G_n^c)$. Since $\lambda>1$, $m_n\ge1$, and $2-\lambda>0$,
\[
 \frac{n^{\lambda-1}m_n^{2-\lambda}}{\log n}
            \ge\frac{n^{\lambda-1}}{\log n}\longrightarrow\infty.
\]
Choose $n_0$ so that the exponential bound is at most $n^{-3}$ for $n\ge n_0$. Together with \eqref{B:eq:simultaneous-boundary-probabilities}, summability gives the first
inequality for all $n\ge n_1(\omega)$ by Borel--Cantelli.

Lemma~\ref{B:lem:gaussian-partial-sum-bounds} gives $\operatorname{Var}(V_{n,j}^G\mid\mathscr E)\le C_1/j$ on $\mathcal G_n$. The scalar Gaussian bound and a union bound give
\[
 \mathbf1_{\mathcal G_n}
 \mathbb P\{\max_{j\in I_n}\sqrt j|V_{n,j}^G|>L\sqrt{\ell_n}
                                 \mid\mathscr E\}
 \le2n(e+n)^{-L^2/(2C_1)}.
\]
Choose $L$ with $L^2/(2C_1)>4$, integrate, and add the boundary failure probability. The resulting probabilities are summable, so the second bound holds for every
$n\ge n_2(\omega)$ on an event of probability one.
\end{proof}

\subsubsection{Initial and final coefficients}
\label{B:subsection:simultaneous-exterior-estimates}

The conditional interior construction leaves the initial and final coefficients to be controlled by their original distributions. The earlier boundary estimates apply to these
coordinates because the sampling has preserved the joint distribution of all coefficients at each matrix size.

The functionals $\mathcal E^-_{n,h}$ defined in \eqref{B:eq:initial-coordinate-functionals} depend on the initial coefficients, whose indices lie outside $J_n$. Let $\mathscr T_n$
be the set of coordinates near the right boundary in \eqref{B:eq:terminal-coordinate-probability}, and let $z_n^\sharp$ be the reference minimizer from
Section~\ref{B:sec:general-potentials-and-convex-reference-measures}.

\begin{proposition}\label{B:prop:simultaneous-exterior-estimates}
The reconstructed rows satisfy, almost surely,
\begin{equation}
 \sup_{n\ge n_0}\mathcal E^-_{n,h}(\mathbf Z_n)<\infty
                              \quad(h\in\mathbb N).
 \label{B:eq:simultaneous-initial-bounds}
\end{equation}
There is a deterministic $C<\infty$ and, almost surely, a finite $n_0(\omega)$ such that for every $n\ge n_0(\omega)$,
\begin{equation}
 \max_{\alpha\in\mathscr T_n}
     |(\mathbf Z_n)_\alpha-(z_n^\sharp)_\alpha|
                                \le C\sqrt{\ell_n/n}.
 \label{B:eq:simultaneous-terminal-bounds}
\end{equation}
The initial product bounds \eqref{B:eq:original-initial-coordinate-bounds} and the bounds for products near the right boundary
\eqref{B:eq:terminal-product-fluctuation}--\eqref{B:eq:terminal-deterministic-product-bound} hold for these rows.
\end{proposition}
\begin{proof}
The reconstruction leaves the exterior coordinates equal to $E_n$. Thus $\mathcal E^-_{n,h}(\mathbf Z_n)=\mathcal E^-_{n,h}(E_n)$. The joint law of the exterior sequence is
$\Lambda^{\rm ext}$, so Theorem~\ref{B:thm:original-exterior-realization} proves \eqref{B:eq:simultaneous-initial-bounds} for all $h$. For any $\eta>0$, choose the index $h$ in
\eqref{B:eq:initial-coordinate-functionals} with $\varepsilon_h\le\eta$. On $(0,1]$, $z^\eta\le z^{\varepsilon_h}$. Put $H_h=\sup_n\mathcal E^-_{n,h}(E_n)<\infty$ on this event.
For either species and $k\le s_n$, its defining terms imply
\[
 \frac{Z_{n,\sigma,k}}{\phi(k/n)}\le H_hz_{n,k}^{-\varepsilon_h},
 \qquad
 \frac{Z_{n,\sigma,k}}{\phi(k/n)}\ge H_h^{-1}z_{n,k}^{\varepsilon_h}.
\]
For $A_{i,j}=P_{n,i}/(D_{n;i,j}P_{n,j})$, the other two terms give $H_h^{-1}(z_{n,i}z_{n,j})^{\varepsilon_h}\le A_{i,j} \le H_h(z_{n,i}z_{n,j})^{-\varepsilon_h}$. Since
$\varepsilon_h\le\eta$ and all nodes are at most one, these inequalities imply \eqref{B:eq:original-initial-coordinate-bounds} with $C_{\omega,\eta}=H_h$.

The equality of the joint coefficient distributions in \eqref{B:eq:reconstructed-full-row-law} and \eqref{B:eq:terminal-coordinate-probability} give
\[
 \mathbb P\left\{\max_{\alpha\in\mathscr T_n}
      |(\mathbf Z_n)_\alpha-(z_n^\sharp)_\alpha|>L\sqrt{\ell_n/n}\right\} 
 \le C_\beta r_n^{\gamma_\beta}
                               e^{-c\gamma_\beta L^2\ell_n}.
\]
Choose $L$ so that $c\gamma_\beta L^2>4$. Since $r_n^{\gamma_\beta}=O((\log n)^{A_0\gamma_\beta})$, the right side is summable. Borel--Cantelli proves
\eqref{B:eq:simultaneous-terminal-bounds}.

On this set of indices near the right boundary the reference coordinates are bounded below by a fixed positive constant, as proved in
Section~\ref{B:sec:general-potentials-and-convex-reference-measures}. For large $n$, $|\log x|\le2|x-1|$ on the relevant ratios gives
\[
 \max_{\alpha\in\mathscr T_n}
   \left|\log\frac{(\mathbf Z_n)_\alpha}{(z_n^\sharp)_\alpha}\right|
                              \le C\sqrt{\ell_n/n}.
\]
In fact, if $z\ge c_0>0$ and $|x-z|\le c_0/2$, then $|\log(x/z)|=|\int_z^x t^{-1}\,\dd t|\le2|x-z|/c_0$. Apply this with the fixed lower bound on $(z_n^\sharp)_\alpha$. Every
product between sites near the right boundary contains at most $n-M_n+R_*\le Cr_n$ ratios. Consequently, uniformly for $M_n-R_*\le i<j\le n$,
\[
 \left|\sum_{k=i}^{j-1}
     \left[\log\frac{Y_{n,k}}{(z_n^\sharp)_{y,k}}
                  -\log\frac{X_{n,k}}{(z_n^\sharp)_{x,k}}\right]\right|
                       \le Cr_n\sqrt{\ell_n/n}\longrightarrow0.
\]
Equation~\eqref{B:eq:terminal-deterministic-product-bound} gives
\[
 \begin{split}
 \sum_{k=M_n-R_*}^{n-1}
 \left|\log\frac{(z_n^\sharp)_{y,k}}{(z_n^\sharp)_{x,k}}\right|
 &\le C\frac{n-M_n+R_*}{n}
      +C\sum_{h=1}^\infty e^{-ch}
      +\frac C{\sqrt n}\sum_{k=M_n-R_*}^\infty e^{-c(k-K)}\\
 &\le C\frac{r_n}{n}
      +\frac C{e^c-1}
      +\frac{Ce^{-c(M_n-R_*-K)}}{\sqrt n(1-e^{-c})}
 \le C.
 \end{split}
\]
Together with the preceding random logarithmic bound this yields
\[
 e^{-C-o(1)}
 \le \prod_{k=i}^{j-1}\frac{Y_{n,k}}{X_{n,k}}
 \le e^{C+o(1)}
 \qquad(M_n-R_*\le i<j\le n),
\]
uniformly in the two indices. The set $\mathscr T_n$ includes the last diagonal coordinate; the product itself stops at
$j-1\le n-1$.
\end{proof}

\subsubsection{Embedded inverses and the common Brownian motion}
\label{B:subsection:simultaneous-embedded-inverses}

The coefficients now define the embedded inverses on the common probability space. They are measurable with respect to the observations specified in the lower bound. The
coefficient estimates give the upper bound for the inverse difference.

Use $\mathcal A_{n,V,\beta,a}$ and $\mathcal R_{n,V,\beta,a}$ from \eqref{B:eq:1-20} and \eqref{B:eq:1-26}, evaluated on the rows just constructed.

\begin{corollary}\label{B:cor:simultaneous-inverse-information}
The operators $\mathcal R_{n,V,\beta,a}$ are Borel $\mathfrak S_2(L^2(0,1))$-valued random variables measurable with respect to $\mathcal N_n\vee\mathcal U$. Almost surely,
\begin{equation}
 \liminf_{n\to\infty}
 n\|\mathcal R_{n,V,\beta,a}-\mathcal H_{\beta,a}^{-1}\|_{\mathfrak S_2}^2>0.
 \label{B:eq:simultaneous-inverse-lower-bound}
\end{equation}
\end{corollary}
\begin{proof}
The coefficient construction gives positive finite $X_{n,k}$, so
\[
 \det\mathcal A_{n,V,\beta,a}
 =(c_Vn^2)^n\prod_{k=1}^nX_{n,k}^2>0.
\]
Inversion is continuous on the invertible matrices, by $A^{-1}-(A')^{-1}=A^{-1}(A'-A)(A')^{-1}$. The deterministic embedding is an isometry in Hilbert--Schmidt norm:
$\|\mathscr J D\mathscr J^*\|_{\mathfrak S_2}= \|D\|_{\mathfrak S_2}$ when $\mathscr J^*\mathscr J=I$. Thus \eqref{B:eq:simultaneous-exact-arrays} makes each embedded inverse
$\mathcal N_n\vee\mathcal U$-measurable. For resolvents the same argument applies jointly in $(A,z)$ on $\det(A-zI)\ne0$. The auxiliary field $\mathcal U$ is independent of the
whole path $W$, and the limiting inverse is constructed from that path. After this entire sequence has been fixed, Theorem~\ref{B:thm:brownian-information-lower-bound} gives
\eqref{B:eq:simultaneous-inverse-lower-bound}.
\end{proof}

For any countable family of deterministic weight arrays, the conclusions of Propositions~\ref{B:prop:simultaneous-local-logarithms} and \ref{B:prop:simultaneous-gaussian-estimates}
hold on the intersection of their full measure events. The coupling of the logarithmic ratios remains the one fixed by the parameters in
Section~\ref{B:subsection:coupling-product-probability-space}. Section~\ref{B:sec:normalization-and-coefficient-estimates} uses \eqref{B:eq:ratio-realized-energy-series},
\eqref{B:eq:simultaneous-local-energy-limit}, and \eqref{B:eq:simultaneous-gaussian-estimates} to estimate the coefficient errors after scalar normalization. The full truncation
estimate is proved in Section~\ref{B:sec:removal-of-the-truncation-and-proof-of-the-general-potential-theorem}.

The beta-Laguerre construction of Section~\ref{B:sec:the-beta-laguerre-coupling} can use the same $W$. Take the product of the present space with the independent chi and gamma
variables specified there, and apply its Borel coefficient maps at the nodes $\log(n/k)$. The product extension preserves the present joint law. Integrating out the present
auxiliary variables leaves Wiener measure and the independent chi and gamma variables, so the joint coefficient distributions and estimates for the beta-Laguerre model are those
proved in Section~\ref{B:sec:the-beta-laguerre-coupling}. The enlarged auxiliary sigma-field remains independent of $W$. Intersecting the two full measure events gives both
constructions on this space.

\subsection{Normalization and coefficient estimates}
\label{B:sec:normalization-and-coefficient-estimates}

We compare the kernel coefficients of the matrices constructed in Section~\ref{B:sec:construction-of-the-coupling-for-all-matrix-sizes} with the coefficients constructed from the
Brownian path in Section~\ref{B:sec:inverse-representations-and-the-stochastic-bessel-operator}. The scalar in \eqref{B:eq:coefficient-normalization-scalar} removes the values at
the upper index of the two partial sums. The remaining coefficient errors are \eqref{B:eq:normalized-p-error} and \eqref{B:eq:normalized-scale-error}. Their maximum controls the
exponential factors, and their weighted square sum then gives the rate through Theorem~\ref{B:thm:weighted-inverse-stability}.

\subsubsection{Notation and the deterministic partition}
\label{B:subsection:normalization-notation-and-grid}

Under the coupling already constructed, we decompose the errors in the two factors of the inverse kernel into deterministic profile errors, Gaussian errors and sampling errors.

Fix the model and construction parameters of Section~\ref{B:subsection:coupling-product-probability-space}. In particular, $\beta>0$, $a>0$, $1<\lambda<\min\{a+1,2\}$, and $A_0>30$
are fixed. We use the arrays of Theorem~\ref{B:thm:simultaneous-coefficient-construction} on its probability space. The intervals $J_n,I_n,I_n^{\rm ret}$, indices
$r_n,b_n,R_n,m_n,M_n$, boundary configurations $q_n$, and sigma-fields $\mathcal N_n,\mathcal U,\mathscr E$ retain their definitions from that section. Thus
$I_n=\{m_n,\ldots,M_n\}$, $I_n^{\rm ret}=\{m_n+1,\ldots,M_n\}$, and $m_n\asymp r_n$. The estimates below concern $n\ge n_0$, where the deterministic starting index may be
increased.

The coordinate map is $\theta=\theta_V$. Use $z_{n,k}$, $\Delta_{n,k}$, and $u_{n,k}$ from \eqref{B:eq:cell-embedding-definition} and \eqref{B:eq:brownian-sampled-information}, the
logarithmic lengths $h_{n,k}^{\pm}$ from \eqref{B:eq:coupling-clocks}, and the weights $w_{n,k}=h_{n,k}^-z_{n,k}^{\lambda}$ from \eqref{B:eq:ratio-partial-sum-cost}.
Equation~\eqref{B:eq:brownian-grid-bounds} implies
\begin{equation}
 h_{n,k}^{\pm}\asymp k^{-1},\qquad
 w_{n,k}\asymp n^{-\lambda}k^{\lambda-1}\quad(k\in J_n),\qquad
 \max_{k\in I_n}h_{n,k}^-\le C/m_n.
 \label{B:eq:normalization-grid-bounds}
\end{equation}
The products $P_{n,k}$ and coefficients $p_{n,k},q_{n,k}$ are those of \eqref{B:eq:2-5}, \eqref{B:eq:2-9}, and \eqref{B:eq:2-13}, evaluated at the constructed array and at
$\alpha_n=\alpha=4/L_V^2$. Their inverse kernel is \eqref{B:eq:2-24}. The Brownian functions $\mathcal Y,p^B,q^B$, cell integrals $Q_{n,k}^B$, and quadrature errors $\rho_{n,k}^B$
use the definitions \eqref{B:eq:1-29}--\eqref{B:eq:2-32} and \eqref{B:eq:brownian-logarithmic-errors}. We use these cell quantities on $I_n$ by taking $k_-=m_n-1$, $k_+=M_n$ in
\eqref{B:eq:brownian-retained-interval}.

Define
\begin{equation}
 \epsilon_{n,k}^{\Delta}
    =\log\frac{n\Delta_{n,k}}{\theta'(k/n)},\qquad
 \mathfrak f(t)=\log\bigl(\theta(t)^{a/2}\sqrt{\theta'(t)}\bigr)
       \quad(0<t\le1).
 \label{B:eq:normalization-reference-function}
\end{equation}
The identity \eqref{B:eq:1-16} for $\theta'$ gives
\begin{align}
 \theta'(t)^2\phi(t)^2&=\alpha\theta(t),\notag\\
 \mathfrak f(t)
   &=\left(\frac a2+\frac14\right)\log\theta(t)
          -\frac12\log\phi(t)+\frac12\log\frac2{L_V},
       \label{B:eq:normalization-profile-identity}\\
 \mathfrak f'(t)
   &=\left(\frac a2+\frac14\right)\frac{\theta'(t)}{\theta(t)}
          -\frac12\frac{\phi'(t)}{\phi(t)}.\notag
\end{align}

\begin{lemma}\label{B:lem:normalization-grid-error}
There is a deterministic $C<\infty$ such that
\begin{equation}
 |n\Delta_{n,k}-\theta'(k/n)|+|\epsilon_{n,k}^{\Delta}|\le C/n
     \quad(1\le k\le n),\qquad
 |\mathfrak f'(t)|\le C/t,\quad |\mathfrak f''(t)|\le C/t^2.
 \label{B:eq:normalization-grid-error}
\end{equation}
The first estimate holds for all sufficiently large $n$; the last two hold for $0<t\le1$.
\end{lemma}
\begin{proof}
Let $c_\theta=\min_{[0,1]}\theta'>0$ and $C_\theta=\|\theta''\|_\infty$, which are finite by Lemma~\ref{B:lem:bessel-coordinate-regularity}. Then
\[
 \begin{split}
 |n\Delta_{n,k}-\theta'(k/n)|
 &\le n\int_{(k-1)/n}^{k/n}
                   |\theta'(s)-\theta'(k/n)|\,\dd s\\
 &\le nC_\theta\int_{(k-1)/n}^{k/n}(k/n-s)\,\dd s
   =\frac{C_\theta}{2n}.
 \end{split}
\]
Consequently $|n\Delta_{n,k}/\theta'(k/n)-1| \le C_\theta/(2nc_\theta)$. For large $n$ this is at most $1/2$, and $|\log(1+x)|\le2|x|$ on $[-1/2,1/2]$ proves the logarithmic
estimate.

The bounds for $\phi$ established in Lemma~\ref{B:lem:constant-coefficient-minimizer} give $\phi(t)\ge c\sqrt t$, $|\phi'(t)|\le Ct^{-1/2}$, and $|\phi''(t)|\le Ct^{-3/2}$. Also
$\theta(t)\ge c_\theta t$, while $\theta'$ and $\theta''$ are bounded. The expression for $\mathfrak f'$ in \eqref{B:eq:normalization-profile-identity} therefore has absolute value
at most $C/t$. Differentiating it gives
\[
 \mathfrak f''(t)
 =\left(\frac a2+\frac14\right)
       \left(\frac{\theta''(t)}{\theta(t)}
                     -\frac{\theta'(t)^2}{\theta(t)^2}\right)
 -\frac12\left(\frac{\phi''(t)}{\phi(t)}
                     -\frac{\phi'(t)^2}{\phi(t)^2}\right).
\]
The coefficient bounds give, for $0<t\le1$,
\[
 \left|\frac{\theta''}{\theta}\right|\le\frac{C}{t},\qquad
 \left|\frac{\theta'}{\theta}\right|^2\le\frac{C}{t^2},\qquad
 \left|\frac{\phi''}{\phi}\right|
       +\left|\frac{\phi'}{\phi}\right|^2\le\frac{C}{t^2}.
\]
Since $t^{-1}\le t^{-2}$ on this interval, substitution proves the second derivative estimate. These constants depend only on the deterministic profile. In particular the preceding
cell calculation also covers $k=1$ and $k=n$, using the entire cell in each case.
\end{proof}

\subsubsection{Ratios of conditional minimizer coefficients}
\label{B:subsection:normalization-minimizer-ratios}

The conditional minimizer contributes a deterministic approximation to each logarithmic ratio. Keeping this contribution separate makes it possible to sum the profile errors before
combining them with random fluctuations.

For a finite nonnegative boundary configuration $q$, let $m^q$ be the conditional minimizer on $J_n$ used in \eqref{B:eq:conditional-matrix-data}. It is positive by
Proposition~\ref{B:prop:conditional-minimizer-dependence}. Put
\begin{equation}
 \varepsilon_{n,k}^q
 =\log\frac{m_{y,k}^q}{m_{x,k}^q}
       -\left[\mathfrak f(k/n)-\mathfrak f((k+1)/n)\right],
                    \qquad k\in J_n.
 \label{B:eq:normalization-ratio-remainder}
\end{equation}
The functions $f_{n,\sigma,k}$ and $g_\sigma$ used below are the first correction in \eqref{B:eq:first-correction-definition}; see \cite{rider2019universality} for these profile
functions. Their derivative bounds and relative error estimates are \eqref{B:eq:first-correction-derivative-bounds}.

\begin{lemma}\label{B:lem:normalization-minimizer-ratio}
For some fixed $c,C>0$, every sufficiently large $n$ and $q\in\mathcal Q_n^\sharp$ satisfy
\begin{equation}
 |\varepsilon_{n,k}^q|
 \le C\left[k^{-2}+e^{-c(n-k)}
       +\vartheta_n\{e^{-c(k-r_n)}+e^{-c(R_n-k)}\}\right],
           \qquad k\in J_n,
 \label{B:eq:normalization-site-remainder}
\end{equation}
where $\vartheta_n=r_n^{-1/8}$ is the boundary tolerance in \eqref{B:eq:admissible-boundary-configurations}.
\end{lemma}
\begin{proof}
Equation~\eqref{B:eq:reference-minimizer-relative-error} implies
\[
 \frac{|(z_n^\sharp)_{\sigma,k}-f_{n,\sigma,k}|}{\phi(k/n)}
 \le Ck^{-2}+Ce^{-c_1(n-k)}.
\]
Indeed, $k^{3/2}e^{-c(k-K)}$ is bounded for $k\ge K$, and $n/k\le1+n-k$ shows that $\sqrt{n/k}e^{-c(n-k)}\le C e^{-c(n-k)/2}$. Next divide
\eqref{B:eq:conditional-minimizer-localization} by $\phi(k/n)$. For its left term, $\phi(r_n/n)/\phi(k/n)\le C\sqrt{r_n/k}\le C$. For its right term,
$\phi(R_n/n)/\phi(k/n)\le C\sqrt{1+R_n-k}$. For
$d=R_n-k\ge0$,
\[
 \sqrt{1+d}\,e^{-c_1d}
 \le A_{c_1}e^{-c_1d/2},\qquad
 A_{c_1}=\sup_{u\ge0}\sqrt{1+u}\,e^{-c_1u/2}<\infty.
\]
Using this bound and taking $c$ no larger than half of each preceding exponential rate gives
\begin{equation}
 \frac{|m_{\sigma,k}^q-f_{n,\sigma,k}|}{\phi(k/n)}
 \le C\left[k^{-2}+e^{-c(n-k)}
       +\vartheta_n\{e^{-c(k-r_n)}+e^{-c(R_n-k)}\}\right].
 \label{B:eq:normalization-minimizer-comparison}
\end{equation}
All constants are uniform on $\mathcal Q_n^\sharp$. The right hand side is uniformly $o(1)$ on $J_n$, and $|f_{n,\sigma,k}/\phi(k/n)-1|\le C/k$. Thus, after increasing $n_0$, both
$m_{\sigma,k}^q/\phi(k/n)$ and $f_{n,\sigma,k}/\phi(k/n)$ lie in $[1/2,3/2]$. Integration of the derivative of $\log$ on this interval yields
\begin{equation}
 \left|\log\frac{m_{\sigma,k}^q}{f_{n,\sigma,k}}\right|
 \le2\frac{|m_{\sigma,k}^q-f_{n,\sigma,k}|}{\phi(k/n)}.
 \label{B:eq:normalization-minimizer-log-difference}
\end{equation}

Set $t=k/n$ and $u_\sigma=g_\sigma(t)/(n\phi(t))$, so that $|u_\sigma|\le C/k\le1/2$. The exact identity
\[
 \log(1+u)-u=-u^2\int_0^1\frac{s}{1+su}\,\dd s
\]
gives $|\log(1+u)-u|\le u^2$ on this range. Therefore
\begin{equation}
 \log\frac{f_{n,y,k}}{f_{n,x,k}}
 =\frac{g_y(t)-g_x(t)}{n\phi(t)}+r_{n,k},\qquad
 |r_{n,k}|\le Ck^{-2}.
 \label{B:eq:normalization-first-correction-ratio}
\end{equation}
By \eqref{B:eq:first-correction-definition} and $\theta'/\theta=2/(\phi\ell)$,
\[
 \frac{g_y-g_x}{\phi}
 =-\frac{a+1/2}{\phi\ell}+\frac{\phi'}{2\phi}
 =-\left(\frac a2+\frac14\right)\frac{\theta'}{\theta}
                                    +\frac{\phi'}{2\phi}
 =-\mathfrak f'.
\]
For the adjacent interval the integral Taylor formula gives
\begin{equation}
 \mathfrak f(t)-\mathfrak f(t+n^{-1})
 =-n^{-1}\mathfrak f'(t)
       -\int_0^{1/n}(n^{-1}-s)\mathfrak f''(t+s)\,\dd s.
 \label{B:eq:normalization-adjacent-integral}
\end{equation}
The integral remainder in this identity has absolute value at most
\[
 \int_0^{1/n}(n^{-1}-s)\frac{C}{(t+s)^2}\,\dd s
 \le\frac{C}{2n^2t^2}=\frac{C}{2k^2}.
\]
The exact subtraction is
\[
 \varepsilon_{n,k}^q
 =\log\frac{m_{y,k}^q}{f_{n,y,k}}
  -\log\frac{m_{x,k}^q}{f_{n,x,k}}+r_{n,k}
  +\int_0^{1/n}(n^{-1}-s)\mathfrak f''(t+s)\,\dd s.
\]
Apply \eqref{B:eq:normalization-minimizer-log-difference} to each logarithm and \eqref{B:eq:normalization-minimizer-comparison} to each species. The last two terms are bounded by
$Ck^{-2}$, proving \eqref{B:eq:normalization-site-remainder}. The interval ends at $(k+1)/n\le1$, since $k\le R_n<n$ throughout this calculation.
\end{proof}

\subsubsection{Summed deterministic errors}
\label{B:subsection:normalization-deterministic-errors}

Products in the inverse kernel become cumulative sums after taking logarithms. The last coefficient has a different exponent from the interior coefficients; the formulas below
separate its contribution before comparing the two kernels.

Set $\jmath_n=M_n+1$ and define, for $j\in I_n$,
\begin{equation}
 \begin{split}
 D_{n,j}^q&=\sum_{k=j}^{M_n}\varepsilon_{n,k}^q,\\
 d_{n,j}^{(1),q}&=D_{n,j}^q-\tfrac12\epsilon_{n,j}^{\Delta},\\
 d_{n,j}^{(2),q}&=-2\epsilon_{n,j}^{\Delta}
                    -2\log\frac{m_{x,j}^q}{\phi(j/n)}.
 \end{split}
 \label{B:eq:normalization-deterministic-errors}
\end{equation}
In particular, the defining sum telescopes exactly as
\begin{equation}
 \sum_{k=j}^{M_n}\log\frac{m_{y,k}^q}{m_{x,k}^q}
 =\mathfrak f(j/n)-\mathfrak f(\jmath_n/n)+D_{n,j}^q.
 \label{B:eq:normalization-minimizer-product}
\end{equation}
Fix the decay rate $c$ in Lemma~\ref{B:lem:normalization-minimizer-ratio}, and put
\begin{equation}
 \epsilon_n=e^{-cr_n}+\vartheta_n e^{-cb_n}.
 \label{B:eq:normalization-boundary-error}
\end{equation}

\begin{proposition}\label{B:prop:normalization-deterministic-bounds}
Uniformly for $q\in\mathcal Q_n^\sharp$ and $j\in I_n$,
\begin{equation}
 |D_{n,j}^q|+|d_{n,j}^{(1),q}|+|d_{n,j}^{(2),q}|
       +\sum_{\sigma\in\{x,y\}}
                    \left|\log\frac{m_{\sigma,j}^q}{\phi(j/n)}\right|
 \le C(j^{-1}+\epsilon_n).
 \label{B:eq:normalization-deterministic-bounds}
\end{equation}
For every fixed $D>0$, $n^D\epsilon_n\to0$.
\end{proposition}
\begin{proof}
The four sums in \eqref{B:eq:normalization-site-remainder} are bounded by
\[
 \begin{split}
 \sum_{k=j}^{M_n}k^{-2}&\le j^{-2}+\int_j^\infty t^{-2}\,\dd t
                                      \le2/j,\\
 \sum_{k=j}^{M_n}e^{-c(n-k)}
     &\le\frac{e^{-c(n-M_n)}}{1-e^{-c}},\\
 \sum_{k=j}^{M_n}e^{-c(k-r_n)}
     &\le\frac{e^{-c(j-r_n)}}{1-e^{-c}},\qquad
 \sum_{k=j}^{M_n}e^{-c(R_n-k)}
     \le\frac{e^{-c(R_n-M_n)}}{1-e^{-c}}.
 \end{split}
\]
For instance, the exact left boundary sum is
\[
 \sum_{k=j}^{M_n}e^{-c(k-r_n)}
 =e^{-c(j-r_n)}\frac{1-e^{-c(M_n-j+1)}}{1-e^{-c}};
\]
reversing the order of the indices gives the two sums decreasing from the other boundary point. Here $j-r_n\ge b_n$, $R_n-M_n=b_n$, and $n-M_n=r_n+b_n$. Thus all their boundary
contributions are at most $C\epsilon_n$, uniformly in the lower summation index $j$, and $|D_{n,j}^q|\le C/j+C\epsilon_n$. To estimate the minimizer logarithms, write
\[
 \log\frac{m_{\sigma,j}^q}{\phi(j/n)}
 =\log\frac{m_{\sigma,j}^q}{f_{n,\sigma,j}}
       +\log\left(1+\frac{g_\sigma(j/n)}{n\phi(j/n)}\right).
\]
The first term is bounded by \eqref{B:eq:normalization-minimizer-comparison} and \eqref{B:eq:normalization-minimizer-log-difference}; the second is at most $C/j$. On $I_n$ the
boundary terms in the first are at most $C\epsilon_n$. Lemma~\ref{B:lem:normalization-grid-error}, with $n^{-1}\le j^{-1}$, now bounds both $d^{(1),q}$ and $d^{(2),q}$.

Choose a deterministic $n_0$ such that $b_n\ge r_n/4$ and $\vartheta_n\le1$ for $n\ge n_0$; so $\epsilon_n\le2e^{-cr_n/4}$. Since $r_n\ge(\log n)^{A_0}$ and $A_0>1$,
\[
 n^D\epsilon_n\le2\exp\{D\log n-c(\log n)^{A_0}/4\}
          \longrightarrow0.
\]
\end{proof}

For later use denote the deterministic weighted error by
\begin{equation}
 E_n^{\rm det}(q)=\sum_{k\in I_n}w_{n,k}
          \left(|d_{n,k}^{(1),q}|^2+|d_{n,k}^{(2),q}|^2\right).
 \label{B:eq:normalization-deterministic-energy-definition}
\end{equation}

\begin{corollary}\label{B:cor:normalization-deterministic-energy}
One has
\begin{equation}
 \sup_{q\in\mathcal Q_n^\sharp} E_n^{\rm det}(q)
 \le C_\lambda\bigl(n^{-\lambda}m_n^{\lambda-2}+\epsilon_n^2\bigr)
          =o(n^{-1}).
 \label{B:eq:normalization-deterministic-energy}
\end{equation}
\end{corollary}
\begin{proof}
Using \eqref{B:eq:normalization-grid-bounds} and Proposition~\ref{B:prop:normalization-deterministic-bounds}, the left side is at most
\[
 C n^{-\lambda}\sum_{k=m_n}^{M_n}k^{\lambda-3}
       +C\epsilon_n^2n^{-\lambda}
                              \sum_{k=m_n}^{M_n}k^{\lambda-1}.
\]
Since $\lambda-3<-1$,
\[
 \sum_{k=m_n}^{M_n}k^{\lambda-3}
 \le m_n^{\lambda-3}+\int_{m_n}^\infty t^{\lambda-3}\,\dd t
 \le\left(1+\frac1{2-\lambda}\right)m_n^{\lambda-2}.
\]
Since $\lambda-1>0$, $\sum_{k=m_n}^{M_n}k^{\lambda-1} \le\int_1^{n+1}t^{\lambda-1}\,\dd t\le C_\lambda n^\lambda$. The square inequality used here is
$(k^{-1}+\epsilon_n)^2\le2k^{-2}+2\epsilon_n^2$. The inequality follows from $2|ab|\le a^2+b^2$ and therefore applies to either sign of $d^{(1),q}$ and $d^{(2),q}$. After
multiplication by $n$,
\[
 0\le n\sup_q E_n^{\rm det}(q)
 \le C_\lambda n^{1-\lambda}m_n^{\lambda-2}
       +C_\lambda n\epsilon_n^2\longrightarrow0.
\]
Indeed, $m_n\ge1$ and $\lambda<2$ bound the first term by $C_\lambda n^{1-\lambda}$, while $n\epsilon_n^2\le4\exp\{\log n-cr_n/2\}\to0$. The supremum here and in the calculation
above is over $\mathcal Q_n^\sharp$.
\end{proof}

\subsubsection{The scalar normalization}
\label{B:subsection:normalization-scalar}

A common scalar normalization leaves the inverse unchanged. Choosing it to remove the upper boundary values of the cumulative errors simplifies the coefficient comparison without
altering the reconstructed matrix or its probability law.

Retain the errors $e_{n,k}$ and their partial sums $v_{n,j}$ from \eqref{B:eq:simultaneous-ratio-error}, the Gaussian field $V_{n,j}^G$ from \eqref{B:eq:coupling-decomposition},
and $\chi_n^\nabla=\Xi_{n,\jmath_n}$ from the same formula. The Brownian increments are $B_{n,k}$ from \eqref{B:eq:coupling-decomposition}. Thus
$B_{n,k}=\mathcal Y(z_{n,k})-\mathcal Y(z_{n,k+1})$. Choose
\begin{equation}
 C_n^{\rm fin}
 =\sqrt\alpha\,P_{n,\jmath_n}
       \exp\!\left\{-\mathfrak f(\jmath_n/n)
                 -\mathcal Y(z_{n,\jmath_n})
                 -\chi_n^\nabla-v_{n,\jmath_n}\right\}.
 \label{B:eq:coefficient-normalization-scalar}
\end{equation}
Use the errors $\delta_{n,k},\rho_{n,k}$ of \eqref{B:eq:2-35} with $c_n=\log C_n^{\rm fin}$, on the cells $k\in I_n$. In particular, the coefficients after scalar normalization are
$\widehat p_{n,k}=(C_n^{\rm fin})^{-1}p_{n,k}$ and $\widehat q_{n,k}=(C_n^{\rm fin})^2q_{n,k}$.

\begin{proposition}\label{B:prop:coefficient-scalar-normalization}
The scalar $C_n^{\rm fin}$ is positive, finite, and $\mathcal N_n\vee\mathcal U$-measurable. For every $j\in I_n$,
\begin{equation}
 \delta_{n,j}=d_{n,j}^{(1),q_n}+v_{n,j}+V_{n,j}^G.
 \label{B:eq:normalized-p-error}
\end{equation}
These conclusions hold for every sampled finite boundary configuration.
\end{proposition}
\begin{proof}
Theorem~\ref{B:thm:simultaneous-coefficient-construction} gives positive finite coefficients and makes the full coefficient vector $\mathcal N_n\vee\mathcal U$-measurable. Thus
$P_{n,\jmath_n}$ is positive, finite, and measurable in that field. The maps $q\mapsto m^q$ are continuous and positive, so the reconstructed logarithmic ratios in
\eqref{B:eq:reconstructed-logarithmic-ratios} are measurable in that field. Lemma~\ref{B:lem:simultaneous-gaussian-laws} gives the same measurability for $u_n^G$ and
$\chi_n^\nabla$. The finite sums defining $v_n$ in \eqref{B:eq:simultaneous-ratio-error} consequently preserve this measurability. Also $\jmath_n\in J_n$ and
$\mathcal Y(z_{n,\jmath_n})$ is one of the sampled Brownian values divided by $\sqrt\beta$. Every term in the exponent in \eqref{B:eq:coefficient-normalization-scalar} is therefore
finite and measurable, so $C_n^{\rm fin}$ is finite, positive, and measurable.

The reconstruction identity \eqref{B:eq:reconstructed-logarithmic-ratios} gives, for $k\in J_n$,
\begin{equation}
 \log\frac{Y_{n,k}}{X_{n,k}}
 =\log\frac{m_{y,k}^{q_n}}{m_{x,k}^{q_n}}+\eta_{n,k}
 =\log\frac{m_{y,k}^{q_n}}{m_{x,k}^{q_n}}
                                      +u_{n,k}^G+e_{n,k}.
 \label{B:eq:normalization-random-ratio}
\end{equation}
The definitions of $v_n$ and $V_n^G$ yield the exact finite sums
\begin{align}
 \sum_{k=j}^{M_n}e_{n,k}
    &=v_{n,j}-v_{n,\jmath_n},\notag\\
 \sum_{k=j}^{M_n}u_{n,k}^G
    &=\mathcal Y(z_{n,j})-\mathcal Y(z_{n,\jmath_n})
                                    +V_{n,j}^G-\chi_n^\nabla.
 \label{B:eq:normalization-random-sums}
\end{align}
The first equality subtracts two sums with the common upper index $R_n$. For the second, subtract the definition $V_{n,j}^G=\chi_n^\nabla+\sum_{k=j}^{M_n}(u_{n,k}^G-B_{n,k})$ and
sum the Brownian increments.

Since $P_{n,j}=P_{n,\jmath_n}\prod_{k=j}^{M_n}(Y_{n,k}/X_{n,k})$, equations \eqref{B:eq:normalization-minimizer-product} and \eqref{B:eq:normalization-random-sums} imply
\begin{equation}
 \begin{split}
 \log P_{n,j}={}&\log P_{n,\jmath_n}
       +\mathfrak f(j/n)-\mathfrak f(\jmath_n/n)+D_{n,j}^{q_n}\\
 &+\mathcal Y(z_{n,j})-\mathcal Y(z_{n,\jmath_n})
       +V_{n,j}^G-\chi_n^\nabla+v_{n,j}-v_{n,\jmath_n}.
 \end{split}
 \label{B:eq:normalization-product-identity}
\end{equation}
Furthermore, \eqref{B:eq:normalization-reference-function} gives
\begin{equation}
 \mathfrak f(j/n)-\tfrac12\log(n\Delta_{n,j})
       =\tfrac a2\log z_{n,j}-\tfrac12\epsilon_{n,j}^\Delta.
 \label{B:eq:normalization-embedding-identity}
\end{equation}
The embedding formula \eqref{B:eq:2-13} is $\log p_{n,j}=\tfrac12\log\alpha+\log P_{n,j} -\tfrac12\log(n\Delta_{n,j})$. First subtract only the scalar:
\[
 \log\frac{p_{n,j}}{C_n^{\rm fin}}
 =\mathfrak f(j/n)-\tfrac12\log(n\Delta_{n,j})
  +D_{n,j}^{q_n}+\mathcal Y(z_{n,j})+v_{n,j}+V_{n,j}^G.
\]
Here the terms at $\jmath_n$ cancel with the signs in \eqref{B:eq:coefficient-normalization-scalar}. Since $\log p^B(z_{n,j})=\tfrac a2\log z_{n,j}+\mathcal Y(z_{n,j})$,
\eqref{B:eq:normalization-embedding-identity} now gives
\[
 \log\frac{p_{n,j}}{C_n^{\rm fin}p^B(z_{n,j})}
       =D_{n,j}^{q_n}-\tfrac12\epsilon_{n,j}^\Delta
                                  +v_{n,j}+V_{n,j}^G,
\]
which is \eqref{B:eq:normalized-p-error} by \eqref{B:eq:normalization-deterministic-errors} and \eqref{B:eq:2-35}.
\end{proof}

For this same scalar, every finite inverse summand satisfies
\[
 \widehat p_{n,i}\widehat p_{n,j}\widehat q_{n,l}
 =(C_n^{\rm fin})^{-2}p_{n,i}p_{n,j}
                    (C_n^{\rm fin})^2q_{n,l}
 =p_{n,i}p_{n,j}q_{n,l}.
\]
Thus the full inverse and every inverse obtained by restricting its triangular factor keep their original normalization, as in \eqref{B:eq:2-20}. The scalar depends only on the
allowed observations and the constructed coordinates; the cell integrals $Q_{n,k}^B$ enter the comparison errors below.

\subsubsection{The scale coefficients}
\label{B:subsection:normalization-scale-coefficients}

The product defining the scale coefficient is unchanged by the scalar normalization. The resulting logarithmic identity separates the local coefficient fluctuation, the
deterministic profile error and the quadrature error $\rho_{n,k}^B$.

The logarithmic coefficient $\zeta_{n,x,k}= \log(X_{n,k}/m_{x,k}^{q_n})$ retains its definition from Section~\ref{B:subsection:simultaneous-local-logarithms}.

\begin{proposition}\label{B:prop:normalized-scale-error}
For every $k\in I_n$ the scale error in \eqref{B:eq:2-35} is
\begin{equation}
 \rho_{n,k}=d_{n,k}^{(2),q_n}-2\zeta_{n,x,k}+\rho_{n,k}^B.
 \label{B:eq:normalized-scale-error}
\end{equation}
\end{proposition}
\begin{proof}
By the coefficient identities of Section~\ref{B:sec:inverse-representations-and-the-stochastic-bessel-operator},
\[
 \widehat q_{n,k}\widehat p_{n,k}^{\,2}
 =q_{n,k}p_{n,k}^{\,2}=\frac{\alpha}{n^2\Delta_{n,k}X_{n,k}^2}.
\]
On the other hand, the definitions of $\delta,\rho$ give this product as $Q_{n,k}^B p^B(z_{n,k})^2e^{\rho_{n,k}}$. All factors are positive. Dividing and using
$p^B(z)^2q^B(z)=z^{-1}$ gives
\[
 e^{\rho_{n,k}}
 =\frac{\alpha z_{n,k}}{n^2\Delta_{n,k}^2X_{n,k}^2}
       \frac{\Delta_{n,k}q^B(z_{n,k})}{Q_{n,k}^B}.
\]
Taking the logarithm, with $\rho_{n,k}^B=\log(\Delta_{n,k}q^B(z_{n,k})/Q_{n,k}^B)$, and using \eqref{B:eq:normalization-profile-identity} give
\begin{equation}
 \begin{split}
 \rho_{n,k}
 &=\log\frac{\alpha z_{n,k}}
                 {n^2\Delta_{n,k}^2X_{n,k}^2}+\rho_{n,k}^B\\
 &=-2\log\frac{n\Delta_{n,k}}{\theta'(k/n)}
           -2\log\frac{X_{n,k}}{\phi(k/n)}+\rho_{n,k}^B\\
 &=-2\epsilon_{n,k}^\Delta
           -2\log\frac{m_{x,k}^{q_n}}{\phi(k/n)}
           -2\zeta_{n,x,k}+\rho_{n,k}^B.
 \end{split}
 \label{B:eq:normalized-scale-calculation}
\end{equation}
The first two terms on the last line equal $d_{n,k}^{(2),q_n}$. This proves the formula for the normalized scale. In particular, the defining scale identity in \eqref{B:eq:2-35}
also reads
\begin{equation}
 \widehat q_{n,k}
 =\Delta_{n,k}q^B(z_{n,k})
       \left(\frac{\theta'(k/n)\phi(k/n)}
                       {n\Delta_{n,k}X_{n,k}}\right)^2
                     e^{-2\delta_{n,k}},
 \label{B:eq:normalized-scale-node-formula}
\end{equation}
since $Q_{n,k}^B e^{\rho_{n,k}^B}=\Delta_{n,k}q^B(z_{n,k})$.
\end{proof}

\subsubsection{Approximation of the cell integrals}
\label{B:subsection:normalization-brownian-quadrature}

The remaining error in sampling the Brownian path comes from replacing an integral of the scale density by its sampled value on each interval. Its maximum must tend to zero, but
the rate is obtained from its weighted square sum.

The set of cells used in this subsubsection is $I_n$; its smallest left boundary is $z_{n,m_n-1}>0$. The normalized oscillations $M_{n,k}$ are those of
\eqref{B:eq:brownian-normalized-oscillation}, with $h_{n,k}=h_{n,k}^-$. Put $\ell_n=\log(e+n)$ as in Section~\ref{B:subsection:simultaneous-ratio-error-estimates}.

\begin{lemma}\label{B:lem:normalization-brownian-quadrature}
There are deterministic $K_B,C_B<\infty$ such that, almost surely for all sufficiently large $n$,
\begin{equation}
 \sum_{k\in I_n}w_{n,k}|\rho_{n,k}^B|^2\le K_B/n,\qquad
 \max_{k\in I_n}|\rho_{n,k}^B|
       \le C_B\left(m_n^{-1}+\sqrt{\ell_n/m_n}\right).
 \label{B:eq:normalization-brownian-quadrature}
\end{equation}
\end{lemma}
\begin{proof}
Equation~\eqref{B:eq:bessel-cell-mass-log-bound} gives
\begin{equation}
 |\rho_{n,k}^B|
 \le(a+1)h_{n,k}^-+\frac2{\sqrt\beta}
                           (h_{n,k}^-)^{1/2}M_{n,k}.
 \label{B:eq:normalization-cell-oscillation}
\end{equation}
For a fixed $n$, put $a_k=z_{n,k}^{\lambda}(h_{n,k}^-)^2$, $a_* =\max_{I_n}a_k$, and $X_k=M_{n,k}^2$. The bounds on the partition intervals give
\begin{equation}
 \sum_{k\in I_n}a_k\le C/n,\qquad
 a_*\le Cn^{-\lambda}m_n^{\lambda-2},\qquad
 \sum_{k\in I_n}a_k^2\le Ca_*/n.
 \label{B:eq:normalization-oscillation-weights}
\end{equation}
For the first inequality, use $n^{-\lambda}\sum_{k=1}^n k^{\lambda-2}\le C_\lambda/n$. The second follows because $\lambda-2<0$ and $k\ge m_n$. The third is
$\sum a_k^2\le a_*\sum a_k$.

The moment calculation in Lemma~\ref{B:lem:weighted-brownian-oscillations} applies to these independent increment paths. In its notation, for a numerical constant $K_0$ and
$0\le t\le(2K_0a_*)^{-1}$,
\[
 \log\mathbb E\exp\left\{t\sum_{k\in I_n}
                   a_k(X_k-\mathbb EX_k)\right\}
       \le2K_0^2t^2\sum_{k\in I_n}a_k^2.
\]
Also \eqref{B:eq:brownian-oscillation-tail} implies $\mathbb EX_k\le32$, so the weighted mean is at most $C_0/n$. Apply \eqref{B:eq:weighted-oscillation-concentration} with
$s=C_0/n$, increasing $C_0$ to be positive if necessary. If $\sum_ka_k^2\le C_2a_*/n$, its two exponents obey
\[
 \frac{s^2}{\sum_ka_k^2}\ge\frac{C_0^2}{C_2na_*},
 \qquad \frac{s}{a_*}=\frac{C_0}{na_*},
 \qquad \frac1{na_*}\ge c n^{\lambda-1}m_n^{2-\lambda}.
\]
All weights are positive on this nonempty interval. After increasing the deterministic threshold we obtain
\begin{equation}
 \mathbb P\left\{\sum_{k\in I_n}a_kM_{n,k}^2>K_1/n\right\}
 \le \exp\{-c n^{\lambda-1}m_n^{2-\lambda}\}.
 \label{B:eq:normalization-quadrature-tail}
\end{equation}
The parameter choice $k_-=m_n-1$, $k_+=M_n$ is admissible in Section~\ref{B:subsection:brownian-data-and-grids}. It includes the cell with index $m_n$ in every displayed sum.

Squaring \eqref{B:eq:normalization-cell-oscillation} yields
\[
 \sum_{k\in I_n}w_{n,k}|\rho_{n,k}^B|^2
 \le 2(a+1)^2\sum_{k\in I_n}z_{n,k}^{\lambda}(h_{n,k}^-)^3
                 +\frac8\beta\sum_{k\in I_n}a_kM_{n,k}^2.
\]
The first sum is at most $C/(nm_n)$, since it is bounded by $(\max h_{n,k}^-)\sum a_k$. The probabilities in \eqref{B:eq:normalization-quadrature-tail} are summable: with
$v=\lambda-1>0$, their upper bound is $e^{-cn^v}$, which is at most $n^{-3}$ whenever $cn^v\ge3\log n$. This condition holds for $n\ge n_0(c,v)$, and Borel--Cantelli gives the
first bound for every $n\ge n_1(\omega)$.

For the maximum, \eqref{B:eq:brownian-oscillation-tail} and the union bound give
\[
 \mathbb P\left\{\max_{k\in I_n}M_{n,k}>L\sqrt{\ell_n}\right\}
 \le4n e^{-L^2\ell_n/8}.
\]
Choose a deterministic $L$ with $L^2/8>5$. This sequence is summable. Equation~\eqref{B:eq:normalization-cell-oscillation}, $\max_{I_n}h_{n,k}^-\le C/m_n$, and Borel--Cantelli give
the second bound for every $n\ge n_2(\omega)$.
\end{proof}

\subsubsection{Weighted bounds for the random terms}
\label{B:subsection:normalization-random-terms}

We estimate the weighted square sums of the random contributions to the logarithmic coefficients. The weights are those required by the deterministic kernel estimate, and the
resulting bounds give the convergence rate.

Use the following nonnegative quantities on the constructed probability space:
\begin{equation}
 \begin{aligned}
 E_n^{\rm rat}&=\sum_{k\in I_n}w_{n,k}|v_{n,k}|^2,&
 E_n^{\rm G}&=\sum_{k\in I_n}w_{n,k}|V_{n,k}^G|^2,\\
 E_n^{\rm loc}&=\sum_{k\in I_n}w_{n,k}
                  (|\zeta_{n,x,k}|^2+|\zeta_{n,y,k}|^2),&
 E_n^{\rm B}&=\sum_{k\in I_n}w_{n,k}|\rho_{n,k}^B|^2.
 \end{aligned}
 \label{B:eq:normalization-random-energies}
\end{equation}

\begin{lemma}\label{B:lem:normalization-random-energies}
There are deterministic $K_G,K_{\rm loc}<\infty$ such that, almost surely,
\begin{equation}
 nE_n^{\rm rat}\longrightarrow0,\qquad
 \limsup_n nE_n^{\rm G}\le K_G,\qquad
 \limsup_n nE_n^{\rm loc}\le K_{\rm loc},\qquad
 \limsup_n nE_n^{\rm B}\le K_B.
 \label{B:eq:normalization-random-energy-bounds}
\end{equation}
\end{lemma}
\begin{proof}
By the definition \eqref{B:eq:ratio-partial-sum-cost}, the weights in $Q_n(e_n)$ agree with $w_{n,k}$ and its sum runs over all of $J_n$. Thus $0\le E_n^{\rm rat}\le Q_n(e_n)$, and
Corollary~\ref{B:cor:ratio-coupling-errors} gives the bound for $E_n^{\rm rat}$.

For the local logarithms choose the deterministic coefficient weights
\[
 a_{n,(\sigma,k)}=w_{n,k}\mathbf1_{I_n}(k)
            \quad((\sigma,k)\in\mathscr S(J_n)).
\]
Equation~\eqref{B:eq:normalization-grid-bounds} gives the required bound $a_{n,(\sigma,k)}\le C_0n^{-\lambda}k^{\lambda-1}$. Their weighted sum $\mathcal M_n$ in
Proposition~\ref{B:prop:simultaneous-local-logarithms} is exactly $E_n^{\rm loc}$. For an admissible $q$, the mean estimate \eqref{B:eq:conditional-energy-moments} gives
\[
 \mathbb E_q E_n^{\rm loc}
 \le C\sum_{\sigma\in\{x,y\}}\sum_{k\in I_n}\frac{w_{n,k}}k
 \le C n^{-\lambda}\sum_{k=1}^n k^{\lambda-2}
 \le K_{\rm loc}/n.
\]
The conditional law of the reconstructed coefficients given $\mathscr E$ is \eqref{B:eq:simultaneous-selected-laws}. Proposition~\ref{B:prop:simultaneous-local-logarithms}, applied
with these weights, gives $n(E_n^{\rm loc}-\mathbb E_{q_n}E_n^{\rm loc})\to0$ on the rows with admissible boundary data. Equation~\eqref{B:eq:simultaneous-boundary-probabilities}
gives $n_0(\omega)$ such that every row with $n\ge n_0(\omega)$ has admissible boundary data. This gives the bound for the local error, increasing $K_{\rm loc}$ if necessary.

The field $V_n^G$ in \eqref{B:eq:normalization-random-energies} is precisely the one in \eqref{B:eq:simultaneous-gaussian-estimates}.
Proposition~\ref{B:prop:simultaneous-gaussian-estimates} gives $E_n^{\rm G}\le K_G/n$ for every $n\ge n_1(\omega)$. The remaining quadrature estimate is
Lemma~\ref{B:lem:normalization-brownian-quadrature}. For these particular weights let $\Omega_{\rm rat}$, $\Omega_{\rm loc}$, $\Omega_{\rm G}$, and $\Omega_{\rm B}$ be the four
events just obtained, and let $\Omega_{\rm ext}=\bigcup_{m\ge1}\bigcap_{n\ge m}\{q_n\in\mathcal Q_n^\sharp\}$. Their intersection has probability one by the union bound applied to
the five null complements. On this intersection, every energy in \eqref{B:eq:normalization-random-energies} is evaluated using the original coefficient array and its coupled
Gaussian field.
\end{proof}

\subsubsection{The weighted coefficient estimate}
\label{B:subsection:normalization-weighted-estimate}

The coefficient comparison follows by combining the deterministic and random contributions with the same kernel weights. The estimate is stated for both logarithmic coefficients
together, since both enter the inverse.

\begin{theorem}\label{B:thm:normalized-coefficient-energy}
For the scalar \eqref{B:eq:coefficient-normalization-scalar} there are a deterministic $K<\infty$ and an almost-surely finite $n_0(\omega)$ such that, for every $n\ge n_0(\omega)$,
\begin{equation}
 \sum_{k\in I_n}w_{n,k}(|\delta_{n,k}|^2+|\rho_{n,k}|^2)\le K/n.
 \label{B:eq:normalized-coefficient-energy}
\end{equation}
\end{theorem}
\begin{proof}
On the event of Lemma~\ref{B:lem:normalization-random-energies}, the boundary configuration is admissible for every $n\ge n_0(\omega)$, and
Corollary~\ref{B:cor:normalization-deterministic-energy} gives $nE_n^{\rm det}(q_n)\to0$. Apply $|x+y+z|^2\le3(|x|^2+|y|^2+|z|^2)$ to \eqref{B:eq:normalized-p-error} and
\eqref{B:eq:normalized-scale-error}. After multiplication by $w_{n,k}$ and summation,
\begin{align}
 \sum_{k\in I_n}w_{n,k}|\delta_{n,k}|^2
 &\le3\sum_{k\in I_n}w_{n,k}|d_{n,k}^{(1),q_n}|^2
                      +3E_n^{\rm rat}+3E_n^{\rm G},\notag\\
 \sum_{k\in I_n}w_{n,k}|\rho_{n,k}|^2
 &\le3\sum_{k\in I_n}w_{n,k}|d_{n,k}^{(2),q_n}|^2
                      +12E_n^{\rm loc}+3E_n^{\rm B}.
 \label{B:eq:normalized-energy-calculation}
\end{align}
Adding the two lines gives the pathwise bound
\[
 \sum_{k\in I_n}w_{n,k}(|\delta_{n,k}|^2+|\rho_{n,k}|^2)
 \le3E_n^{\rm det}(q_n)+3E_n^{\rm rat}+3E_n^{\rm G}
                         +12E_n^{\rm loc}+3E_n^{\rm B}.
\]
Multiplication by $n$ and the preceding limits bound its upper limit by $K_*=3K_G+12K_{\rm loc}+3K_B$. For each path in
the stated event, the definition of upper limit supplies an integer after which the left side multiplied by $n$ is at
most $K_*+1$. Thus $K=K_*+1$ is deterministic, while the starting index may depend on the path. Taking the least such
index is measurable: its sublevel set at $m$ is the intersection of the displayed inequalities over the integers $n\ge
m$.
\end{proof}

\begin{corollary}\label{B:cor:normalized-coefficient-remainders}
Almost surely,
\begin{equation}
 n\sum_{k\in I_n}w_{n,k}
 \left(|\delta_{n,k}-V_{n,k}^G|^2
       +|\rho_{n,k}+2\zeta_{n,x,k}-\rho_{n,k}^B|^2\right)
                  \longrightarrow0.
 \label{B:eq:normalized-coefficient-remainders}
\end{equation}
\end{corollary}
\begin{proof}
The identities \eqref{B:eq:normalized-p-error} and \eqref{B:eq:normalized-scale-error} identify the two errors inside the sum as $d_{n,k}^{(1),q_n}+v_{n,k}$ and
$d_{n,k}^{(2),q_n}$. For each $k$, the elementary identity $2|d|^2+2|v|^2-|d+v|^2=|d-v|^2\ge0$ gives
\[
 |d_{n,k}^{(1),q_n}+v_{n,k}|^2+|d_{n,k}^{(2),q_n}|^2
 \le2\bigl(|d_{n,k}^{(1),q_n}|^2+|d_{n,k}^{(2),q_n}|^2\bigr)
                                      +2|v_{n,k}|^2.
\]
Multiply by the nonnegative $nw_{n,k}$ and sum. The displayed quantity is at most $2nE_n^{\rm det}(q_n)+2nE_n^{\rm rat}$, whose two terms tend to zero by
Corollary~\ref{B:cor:normalization-deterministic-energy} and Lemma~\ref{B:lem:normalization-random-energies}. This also identifies the weighted sum of squared coefficient errors up
to $o_\omega(n^{-1})$. Use $\langle f,g\rangle_w=\sum_{k\in I_n}w_{n,k}f_kg_k$ for real sequences, and the corresponding Hilbert direct sum norm for pairs. Set
$U_n=(V_n^G,-2\zeta_{n,x}+\rho_n^B)$ and $R_n=(d_n^{(1),q_n}+v_n,d_n^{(2),q_n})$. The bounds for the weighted squared norms give $\|U_n\|_w=O_\omega(n^{-1/2})$ and the preceding
calculation gives $\|R_n\|_w=o_\omega(n^{-1/2})$. Therefore
\[
 n\left|\|U_n+R_n\|_w^2-\|U_n\|_w^2\right|
 \le n\|R_n\|_w(2\|U_n\|_w+\|R_n\|_w)\longrightarrow0.
\]
Expanding the second component of $U_n$ yields
\[
 \sum_{k\in I_n}w_{n,k}(|\delta_{n,k}|^2+|\rho_{n,k}|^2)
 =E_n^{\rm G}+4\|\zeta_{n,x}\|_w^2+E_n^{\rm B}
       -4\langle\zeta_{n,x},\rho_n^B\rangle_w
                                      +o_\omega(n^{-1}).
\]
The inner product in this identity involves the constructed coefficients and the quadrature errors $\rho_{n,k}^B$ on the same probability space. Its asymptotics, together with
those of the two squared terms, are needed to determine a leading order coefficient for the inverse error.
\end{proof}

\subsubsection{Uniform coefficient estimates}
\label{B:subsection:normalization-uniform-estimates}

The weighted estimate gives the rate, but the nonlinear stability argument also requires uniform smallness of the logarithmic errors. The following bounds provide that separate
requirement on the selected interval.

\begin{proposition}\label{B:prop:normalized-coefficient-maxima}
There is a deterministic $C<\infty$ and, almost surely, a finite $n_0(\omega)$ such that for every $n\ge n_0(\omega)$,
\begin{align}
 \max_{k\in I_n}|\delta_{n,k}|
 &\le C\left(m_n^{-1}+\epsilon_n
                  +\ell_n^2m_n^{-1/2}+\sqrt{\ell_n/m_n}\right),\notag\\
 \max_{k\in I_n}|\rho_{n,k}|
 &\le C\left(m_n^{-1}+\epsilon_n+\sqrt{\ell_n/m_n}\right).
 \label{B:eq:normalized-coefficient-maxima}
\end{align}
In particular both maxima tend to zero.
\end{proposition}
\begin{proof}
On the event of probability one where $q_n\in\mathcal Q_n^\sharp$ for $n\ge n_0(\omega)$, Proposition~\ref{B:prop:normalization-deterministic-bounds} gives
\[
 \max_{k\in I_n}\left(|d_{n,k}^{(1),q_n}|
                          +|d_{n,k}^{(2),q_n}|\right)
       \le C(m_n^{-1}+\epsilon_n).
\]
The definition \eqref{B:eq:ratio-cost-exponents} implies $\widehat p_n\le C\ell_n$. Thus \eqref{B:eq:ratio-realized-partial-sum-maximum} and
\eqref{B:eq:simultaneous-gaussian-estimates} imply
\[
 \max_{k\in I_n}|v_{n,k}|\le C\ell_n^2/\sqrt{m_n},\qquad
 \max_{k\in I_n}|V_{n,k}^G|\le C\sqrt{\ell_n/m_n}.
\]
Equation~\eqref{B:eq:simultaneous-local-maximum} gives $\max_{k\in I_n}|\zeta_{n,x,k}|\le C\sqrt{\ell_n/m_n}$. Lemma~\ref{B:lem:normalization-brownian-quadrature} bounds
$\max_{I_n}|\rho_{n,k}^B|$. The two pointwise identities imply the separate inequalities
\[
 \begin{split}
 \max_{I_n}|\delta_{n,k}|&\le
   \max_{I_n}|d_{n,k}^{(1),q_n}|+
   \max_{I_n}|v_{n,k}|+\max_{I_n}|V_{n,k}^G|,\\
 \max_{I_n}|\rho_{n,k}|&\le
   \max_{I_n}|d_{n,k}^{(2),q_n}|+
   2\max_{I_n}|\zeta_{n,x,k}|+\max_{I_n}|\rho_{n,k}^B|.
 \end{split}
\]
Combining the four estimates above proves \eqref{B:eq:normalized-coefficient-maxima} on their common event.

Since $m_n\asymp(\log n)^{A_0}$,
\[
 \ell_n^2m_n^{-1/2}=O((\log n)^{2-A_0/2})\to0,\qquad
 \sqrt{\ell_n/m_n}=O((\log n)^{(1-A_0)/2})\to0.
\]
Also $m_n^{-1}\to0$ and $\epsilon_n\to0$. Thus both maxima tend to zero and, in particular, $\max_{I_n}(|\delta_{n,k}|+|\rho_{n,k}|)\le1$ for every $n\ge n_1(\omega)$.
\end{proof}

\subsubsection{The inverse on the selected interval}
\label{B:subsection:the-inverse-comparison-on-the-retained-interval}

The coefficient bounds now fit the deterministic comparison of embedded inverses. The finite inverse and the inverse constructed from the Brownian path are restricted to the same
cells, so the only errors at this stage are those already estimated in their kernel coefficients.

Use the embedding $\widetilde{\mathscr J}_n$ from \eqref{B:eq:cell-embedding-definition} and the operators $\widetilde{\mathcal K}_n,\widetilde{\mathcal R}_n$ from
\eqref{B:eq:criterion-operators}. In the restricted factors of Section~\ref{B:subsection:the-criterion-for-convergence-of-the-inverses}, take $k_-=m_n$ and $k_+=M_n$. Thus
$P_n^0=\Pi_{z_{n,m_n},z_{n,M_n}}$ and
\begin{equation}
 \widetilde{\mathcal R}_n^0
   =(P_n^0\widetilde{\mathcal K}_nP_n^0)^*
                             (P_n^0\widetilde{\mathcal K}_nP_n^0).
 \label{B:eq:normalization-restricted-inverse}
\end{equation}
On these same cells the comparison operator is $\mathcal R_n^{B,\mathrm{step},0}$ from \eqref{B:eq:criterion-brownian-cuts}.

\begin{corollary}\label{B:cor:normalized-restricted-inverse}
Almost surely, for all sufficiently large $n$,
\begin{equation}
 \|\widetilde{\mathcal R}_n^0
             -\mathcal R_n^{B,\mathrm{step},0}\|_{\mathfrak S_2}^2
                  \le C_\omega/n.
 \label{B:eq:normalization-restricted-inverse-bound}
\end{equation}
\end{corollary}
\begin{proof}
For $i,j\in I_n^{\rm ret}$, Proposition~\ref{B:result:2-6} gives
\begin{equation}
 \widetilde{\mathcal R}_n^0(i,j)
 =p^B(z_{n,i})p^B(z_{n,j})e^{\delta_{n,i}+\delta_{n,j}}
       \sum_{l=i\vee j}^{M_n}Q_{n,l}^B
                              e^{-2\delta_{n,l}+\rho_{n,l}}.
 \label{B:eq:normalization-restricted-kernel}
\end{equation}
Setting $\delta_{n,k}=\rho_{n,k}=0$ in this expression gives $\mathcal R_n^{B,\mathrm{step},0}$. Both kernels vanish outside the square of the selected interval. Their squared
Hilbert--Schmidt difference is consequently
\[
 \sum_{i,j\in I_n^{\rm ret}}\Delta_{n,i}\Delta_{n,j}
 \left|\widetilde{\mathcal R}_n^0(i,j)
                     -\mathcal R_n^{B,\mathrm{step},0}(i,j)\right|^2.
\]
Choose $\eta$ in the interval \eqref{B:eq:stability-exponents}. Lemma~\ref{B:lem:bessel-brownian-growth} gives $p^B(z_{n,k})\le C_{\omega,\eta}z_{n,k}^{a/2-\eta}$. With
$v=a+2\eta>0$,
\[
 Q_{n,k}^B
 \le C_{\omega,\eta}\int_{z_{n,k-1}}^{z_{n,k}}s^{-v-1}\,\dd s
 =C_{\omega,\eta}z_{n,k}^{-v}
                         \frac{e^{v h_{n,k}^-}-1}{v}.
\]
For $m_n\ge C$, the bound $h_{n,k}^-\le C/m_n\le1$ gives $(e^{vh}-1)/v=\int_0^h e^{vs}\,\dd s\le e^v h$. These are exactly the two bounds \eqref{B:eq:reference-coefficient-bounds}.
The cell identity $\Delta_{n,k}=z_{n,k}(1-e^{-h_{n,k}^-})$ verifies the hypothesis on the partition of Section~\ref{B:subsection:interpolation-of-kernel-coefficients}, and
Proposition~\ref{B:prop:normalized-coefficient-maxima} gives \eqref{B:eq:logarithmic-errors-bounded} on $I_n^{\rm ret}$. For the interpolation parameter $0\le t\le1$, the exponent
in one summand of \eqref{B:eq:normalization-restricted-kernel} is $t(\delta_{n,i}+\delta_{n,j}-2\delta_{n,l}+\rho_{n,l})$. The maximum assumption gives
\[
 |\delta_{n,i}+\delta_{n,j}-2\delta_{n,l}+\rho_{n,l}|
 \le |\delta_{n,i}|+|\delta_{n,j}|+
                      2|\delta_{n,l}|+|\rho_{n,l}|\le4.
\]
The constants in the derivative estimate of Theorem~\ref{B:thm:weighted-inverse-stability} are consequently uniform in $t$ and $n$ after fixing the constant in the growth bound for
this path. In particular, with $H_{ijl}=\delta_{n,i}+\delta_{n,j}-2\delta_{n,l}+\rho_{n,l}$, the exact difference of the two cell kernels is
\[
 p^B(z_{n,i})p^B(z_{n,j})
 \sum_{l=i\vee j}^{M_n}Q_{n,l}^B H_{ijl}
                                  \int_0^1e^{tH_{ijl}}\,\dd t.
\]
This follows from $e^H-1=H\int_0^1e^{tH}\,\dd t$ and uses the same finite coefficients throughout $0\le t\le1$. The maximum bound controls the exponential factor; the weighted
square sum then controls the occurrences of $\delta$ and $\rho$ in the kernel estimate. Both assumptions are used in this application. Applying that theorem bounds the last double
sum by
\[
 C_{\omega,\eta}\sum_{k\in I_n^{\rm ret}}
                h_{n,k}^-z_{n,k}^{\lambda}
                     (|\delta_{n,k}|^2+|\rho_{n,k}|^2)
 \le C_{\omega,\eta}K/n,
\]
where the last step is Theorem~\ref{B:thm:normalized-coefficient-energy}.
\end{proof}

Proposition~\ref{B:prop:coefficient-scalar-normalization}, Theorem~\ref{B:thm:normalized-coefficient-energy}, and Proposition~\ref{B:prop:normalized-coefficient-maxima} verify the
coefficient hypothesis \eqref{B:eq:criterion-coefficient-errors} of the inverse criterion. Section~\ref{B:sec:removal-of-the-truncation-and-proof-of-the-general-potential-theorem}
estimates the difference between the full inverse and \eqref{B:eq:normalization-restricted-inverse}. The comparison above concerns the inverses formed from the two restricted
triangular factors on $I_n^{\rm ret}$. Passing to the full inverse, or to its compression, also requires the positive rank-one contribution from the discarded columns. We estimate
this contribution in Section~\ref{B:subsection:finite-factor-restriction}.

\subsection{Coupling of the circular and Brownian increments}
\label{sec:combined-s-increments}
\label{S:sec:coupling-the-finite-transfer-increments-to-brownian-increments}
Fix $\beta>0$. We construct the circular walks jointly with the Brownian path, using a common terminal point and preserving the joint distribution of each finite walk. Conditional
on that terminal point, each transition has its original radial law and an angular density depending on the radius. We couple the radii by their quantiles and use the same uniform
quantile for the two conditional angles.

The estimates must also distinguish a centered increment from its mean. The radial discrepancy has size $h^{3/2}$ over a time interval of length $h$; averaging the relative
logarithm over the angle improves its mean to order $h^2$. This extra factor $h^{1/2}$ is needed when the means are summed in
Section~\ref{S:sec:products-of-the-coupled-transfer-matrices}. The construction below is exact before any localization. We localize only the estimates and prove that the localized
and original increments agree at every selected index for all sufficiently large matrix sizes.

\subsubsection{Radial distributions and their matrix representation}
The one-step coupling begins with two radial distributions. Their parameters must be matched to the Brownian time increments before the conditional angles are coupled; the two
cutoffs serve different roles in that comparison.

For $k\ge1$, put
\[
 \gamma_k=\frac{\beta k}{2},\qquad
 h_k=\frac4{\beta k+1}=\frac4{2\gamma_k+1}.
\]
Thus $h_k=\vartheta_{n,k}-\vartheta_{n,k+1}$ for $1\le k<n$. The two cutoffs remain those of~\eqref{S:eq:cutoff-choice}: the coupling at deterministic times is used for
$k\ge\underline m_n$, and the operator perturbation estimates are subsequently restricted to $k\ge m_n$. On either range the largest $h_k$ tends to zero. Throughout this
subsection, $\|\cdot\|$ denotes the matrix operator norm.

Let $R_h$ have the law of the hyperbolic distance traveled in time $h$ by Brownian motion with generator $\frac12\Delta_{\mathbb H}$. Its distribution function and density are
denoted by $F_h$ and $p_h$. Isometry invariance makes this law independent of the starting point. For $\Xi\sim\operatorname{Beta}(1,\gamma)$, set
$Y_\gamma=2\operatorname{arctanh}\sqrt\Xi$ and denote its distribution function and density by $G_\gamma$ and $g_\gamma$. Integration of the beta density gives
\begin{equation}\label{S:eq:7.4}
 1-G_\gamma(r)
 =\mathbb P\bigl(\Xi>\tanh^2(r/2)\bigr)
 =\operatorname{sech}^{2\gamma}(r/2),\qquad r\ge0.
\end{equation}
Differentiation therefore yields
\[
 g_\gamma(r)=\gamma\sinh(r/2)\operatorname{sech}^{2\gamma+1}(r/2),
 \qquad r>0.
\]
For $\gamma=2/h-1/2$, the two radial laws have the same leading Gaussian exponent $-r^2/(2h)$. The quantile map below sends the Brownian radial law exactly to the circular radial
law. At this stage the angle is unconditioned; we next identify the change in its law conditional on the terminal point.

Use the matrices
\[
 H_0=\frac12\begin{pmatrix}1&0\\0&-1\end{pmatrix},\qquad
 K_\phi=\begin{pmatrix}
 \cos(\phi/2)&\sin(\phi/2)\\-\sin(\phi/2)&\cos(\phi/2)
 \end{pmatrix},\qquad
 T(r,\phi)=K_\phi\me^{rH_0}K_\phi^{-1}.
\]
Since $K_\phi\cdot\ii=\ii$, the point $z(r,\phi)=T(r,\phi)\cdot\ii$ is at distance $r$ from $\ii$. To compute it, write $c=\cos(\phi/2)$ and $s=\sin(\phi/2)$. Then
\[
 K_\phi\cdot(\ii\me^r)
 =\frac{s+\ii c\me^r}{c-\ii s\me^r}
 =\frac{cs(1-\me^{2r})+\ii\me^r}{c^2+s^2\me^{2r}}.
\]
Multiplying numerator and denominator by $\me^{-r}$ gives
\begin{equation}\label{S:eq:7.10}
 x(r,\phi)=\frac{-\sin\phi\sinh r}{\cosh r-\cos\phi\sinh r},\qquad
 y(r,\phi)=\frac1{\cosh r-\cos\phi\sinh r}.
\end{equation}
In particular, $\me^{-r}\le y(r,\phi)\le\me^r$. The point $z(r,\phi)$ determines the coefficient through the matrix $\mathcal A$ defined in
Section~\ref{S:sec:green-kernels-for-canonical-systems}:
\[
 \mathsf A(r,\phi)=\mathcal A(z(r,\phi))
 =\begin{pmatrix}
 y(r,\phi)^{-1/2}&-x(r,\phi)y(r,\phi)^{-1/2}\\
 0&y(r,\phi)^{1/2}
 \end{pmatrix}.
\]
The actions of $T$ and $\mathsf A$ have opposite directions: $T\cdot\ii=z$, whereas $\mathsf A\cdot z=\ii$. Thus $\mathsf A(r,\phi)T(r,\phi)$ fixes $\ii$ and belongs to
$\mathrm{SO}_2(\mathbb R)$. Equivalently,
\begin{equation}\label{S:eq:s6-transvection-inverse}
 \mathsf A(r,\phi)=O(r,\phi)T(r,\phi)^{-1},
 \qquad O(r,\phi)\in\mathrm{SO}_2(\mathbb R).
\end{equation}
Here a real matrix with determinant one $\left(\begin{smallmatrix}a&b\\c&d\end{smallmatrix}\right)$ fixes $\ii$ exactly
when $a=d$, $b=-c$, and $a^2+b^2=1$.

For $z_i=x_i+\ii y_i$, let $w=\mathcal A(z_0)\cdot z_1=(x_1-x_0)/y_0+\ii y_1/y_0$. Multiplication gives
\[
 \mathcal A(w)\mathcal A(z_0)
 =\begin{pmatrix}
 y_1^{-1/2}&-x_0y_1^{-1/2}-(x_1-x_0)y_1^{-1/2}\\
 0&y_1^{1/2}
 \end{pmatrix}.
\]
Consequently,
\begin{equation}\label{S:eq:affine-section-exact-composition}
 \mathcal A(z_1)=\mathcal A(w)\mathcal A(z_0).
\end{equation}
All transfer recursions below use this order of multiplication. If both coefficient functions are multiplied on the left by the same measurable $O(t)\in\mathrm{SO}_2$, their kernel
difference is multiplied by $O(s)$ and $O(t)^{\mathsf T}$ in its two variables. Its squared matrix norm, and hence its Hilbert--Schmidt norm, are unchanged. For the finite matrices
this transformation is
\[
 \widetilde H_{n,k}^{\diamond}
 =O_{n,k}H_{n,k}^{\diamond}O_{n,k+1}^{-1}.
\]
We use this identity only for a common multiplier. Section~\ref{S:sec:comparison-of-the-finite-and-sampled-brownian-operators} treats the different multipliers arising from the two
positions.

\subsubsection{Conditioning on the terminal point}
\label{S:sec:s6-terminal-conditioning}
The inverse kernel contains the right boundary vector. Coupling only the interior positions would therefore leave a boundary error. We condition both walks on the same terminal
point and calculate the resulting transition law before choosing the coupling.

For $z=x+\ii y$ and $\eta\in\mathbb R$, harmonic measure and its Poisson density are
\[
 \nu_z(\dd\eta)=\frac1\pi\frac{y\,\dd\eta}{(\eta-x)^2+y^2},
 \qquad
 P_\eta(z)=\frac{\dd\nu_z}{\dd\nu_\ii}(\eta)
 =\frac{(1+\eta^2)y}{(\eta-x)^2+y^2}.
\]
Set $P_\infty(z)=y$. All these positive harmonic functions have value one at $\ii$. The formula for $\nu_z$ follows by sending $z$ to the center of the disk: the terminal angle is
uniform by rotation invariance, its real boundary coordinate is $\cot(\phi/2)$, and $|\dd\phi/\dd\eta|=2/(1+\eta^2)$. Translation by $x$ and scaling by $y$ give the displayed
density.

We shall use the unique orientation-preserving isometry $g_{z,\eta}$ sending $z$ to $\ii$ and $\eta$ to $\infty$. For a concrete Borel choice, put
\[
 G_\eta=\frac1{\sqrt{1+\eta^2}}
       \begin{pmatrix}\eta&1\\-1&\eta\end{pmatrix}
 \quad(\eta\in\mathbb R),\qquad G_\infty=I,
 \qquad g_{z,\eta}=\mathcal A(G_\eta\cdot z)G_\eta.
\]
These formulas are interpreted as M\"obius transformations. Both $G_\eta\cdot\ii=\ii$ and $G_\eta\cdot\eta=\infty$, and $\mathcal A$ fixes $\infty$. Two such isometries differ by
an isometry fixing $\ii$ and $\infty$, which is the identity. The height formula for a real matrix with determinant one acting by a M\"obius transformation gives
\[
 \Im(G_\eta\cdot z)=P_\eta(z),\qquad
 \Im(g_{z,\eta}\cdot z')=
       \frac{P_\eta(z')}{P_\eta(z)}.
\]
The formula for $\eta=\infty$ follows directly from $g_{z,\infty}=\mathcal A(z)$.

Let $K_1,\ldots,K_N$ be probability transition kernels on $\mathbb H$ which are invariant under orientation-preserving isometries. Let $Z_0=\ii$ and let the chain use these
kernels; after its last step, choose $\eta$ according to $\nu_{Z_N}$. Define
\begin{equation}\label{S:eq:7.41}
 K_j^\eta(z,\dd z')=K_j(z,\dd z')\frac{P_\eta(z')}{P_\eta(z)},
 \qquad \varkappa_j(\dd w)=\Im(w)K_j(\ii,\dd w).
\end{equation}

\begin{lemma}\label{S:lem:7.2}
The conditional chain given $\eta$ has transitions $K_j^\eta$. Its local variables $g_{Z_{j-1},\eta}\cdot Z_j$ are independent with laws $\varkappa_j$, which do not depend on
$\eta$.
\end{lemma}
\begin{proof}
In polar coordinates an invariant kernel based at $\ii$ has a radial probability measure and uniform angular measure. By \eqref{S:eq:7.10}, substitution $t=\tan(\phi/2)$ on
$(-\pi,\pi)$ gives
\begin{equation}\label{S:eq:7.42}
 \frac1{2\pi}\int_0^{2\pi}y(r,\phi)\,\dd\phi
 =\frac1\pi\int_{-\infty}^\infty
     \frac{\dd t}{\me^{-r}+\me^r t^2}=1.
\end{equation}
The final integral is $\pi/\sqrt{\me^{-r}\me^r}=\pi$ before division by $\pi$. The statement also holds for an atom at $r=0$. Isometry invariance and the height ratio above now
imply
\[
 \int K_j(z,\dd z')P_\eta(z')=P_\eta(z).
\]
In particular, $\eta$ has distribution $\nu_\ii$, since successive integration gives $\mathbb E P_\eta(Z_N)=1$.

The joint law of the chain and its terminal point can be written
\[
 \nu_\ii(\dd\eta)\,
 \delta_\ii(\dd z_0)\prod_{j=1}^N K_j(z_{j-1},\dd z_j)
 P_\eta(z_N).
\]
Because $P_\eta(z_0)=1$,
\[
 P_\eta(z_N)=\prod_{j=1}^N
       \frac{P_\eta(z_j)}{P_\eta(z_{j-1})}.
\]
The preceding normalization makes the resulting product exactly a probability kernel in $\eta$. Testing against bounded Borel functions proves that it is a version of the
conditional law.

Given $Z_{j-1}=z$ and $\eta$, push $K_j^\eta(z,\cdot)$ forward under $g_{z,\eta}$. Isometry invariance sends $K_j(z,\cdot)$ to $K_j(\ii,\cdot)$, and the density factor becomes
$\Im(w)$. Thus for a bounded Borel function $f$,
\[
 \mathbb E\bigl[f(g_{Z_{j-1},\eta}\cdot Z_j)
       \mid\eta,Z_0,\ldots,Z_{j-1}\bigr]
 =\int f(w)\,\varkappa_j(\dd w).
\]
For bounded Borel $f_1,\ldots,f_N$, put $c_j=\int f_j(w)\Im(w)K_j(\ii,\dd w)$. Conditioning successively at steps $N,N-1,\ldots,1$ gives
\[
 \mathbb E\!\left[\prod_{j=1}^N f_j(g_{Z_{j-1},\eta}Z_j)\,\middle|\,\eta\right]
 =\prod_{j=1}^N c_j.
\]
This product law is independent of $\eta$.

For Brownian motion at deterministic times, its Markov property and the terminal harmonic measure give $\mathbb P(\eta_\infty\in\dd\eta\mid Z_s,\,s\le t) =\nu_{Z_t}(\dd\eta)$. The
same computation therefore applies to any finite list of sampling times. This identifies the conditional finite-dimensional laws with the Poisson Doob transform. A countable
generating family of cylinder events identifies the path law on each finite horizon; intersecting over integer horizons gives a single conditional version for the path.
\end{proof}

Equation~\eqref{S:eq:7.42} shows that the radial distribution is unchanged by terminal conditioning. Given the radius, the angular law is
\begin{equation}\label{S:eq:7.43}
 \pi_r(\dd\phi)=\frac{y(r,\phi)}{2\pi}\,\dd\phi.
\end{equation}
Write $F_r(\phi)=\pi_r([0,\phi])$ and $\Phi_r=F_r^{-1}$ on $(0,1)$. We use $F_r$ for the angular distribution function and $F_h$ for the Brownian radial distribution function. Set
$E_{12}=\left(\begin{smallmatrix}0&1\\0&0\end{smallmatrix}\right)$ and $\mathsf V_0(\phi)=\sin\phi E_{12}-\cos\phi H_0$.

\begin{lemma}\label{S:lem:7.3}
The angular quantile is continuous for $r\ge0$, $0<u<1$, and smooth for $r>0$. Moreover,
\begin{equation}\label{S:eq:7.44}
 |\Phi_r(u)-\Phi_\rho(u)|\le |r-\rho|
 \le2\pi\me^{r\vee\rho}|r-\rho|,
\end{equation}
and
\begin{equation}\label{S:eq:7.45}
 \int\mathsf V_0(\phi)\,\pi_r(\dd\phi)=-\tanh(r/2)H_0.
\end{equation}
\end{lemma}
\begin{proof}
On each bounded interval of $r$, the derivative $\partial_\phi F_r=y(r,\phi)/(2\pi)$ is bounded away from zero. Differentiation under the angular integral is allowed because its
integrand and derivatives are bounded there. In particular,
\[
 \partial_r y(r,\phi)
 =\frac{\cos\phi\cosh r-\sinh r}
        {(\cosh r-\cos\phi\sinh r)^2}
 =\partial_\phi\bigl(y(r,\phi)\sin\phi\bigr).
\]
Integration from zero to $\phi$ gives $\partial_rF_r(\phi)=y(r,\phi)\sin\phi/(2\pi)$. The implicit function theorem, also applied after extending $r$ slightly below zero, proves
the stated regularity and
\[
 \partial_r\Phi_r(u)=-\sin\Phi_r(u).
\]
Integration in $r$ proves \eqref{S:eq:7.44} for the continuous real branch of $\Phi_r$.

The sine integral in~\eqref{S:eq:7.45} is zero by the change $\phi\mapsto2\pi-\phi$. For $r>0$, use~\eqref{S:eq:7.42} to get
\[
 1=\int(\cosh r-\cos\phi\sinh r)\,\pi_r(\dd\phi)
 =\cosh r-\sinh r\int\cos\phi\,\pi_r(\dd\phi).
\]
Thus the last integral is $(\cosh r-1)/\sinh r=\tanh(r/2)$. At $r=0$, $y(0,\phi)=1$, and hence
\[
 \int\sin\phi\,\pi_0(\dd\phi)=\frac1{2\pi}\int_0^{2\pi}\sin\phi\,\dd\phi=0,
 \qquad
 \int\cos\phi\,\pi_0(\dd\phi)=\frac1{2\pi}\int_0^{2\pi}\cos\phi\,\dd\phi=0.
\]
Thus the left side of \eqref{S:eq:7.45} is zero, equal to $-\tanh(0)H_0$.
\end{proof}
The explicit quantile formula below allows us to integrate the relative logarithm over the angular variable and compute its mean.

Write $\rho_r=\tanh(r/2)$ and choose $\Phi_r(u)$ continuously in $[0,2\pi]$, with $\Phi_r(1/2)=\pi$.

\begin{lemma}\label{S:lem:conditioned-angular-explicit}
For $r\ge0$ and $0<u<1$,
\begin{equation}\label{S:eq:conditioned-angular-explicit}
 \me^{\ii\Phi_r(u)}
 =\frac{\me^{2\pi\ii u}+\rho_r}{1+\rho_r\me^{2\pi\ii u}}.
\end{equation}
Equivalently, away from the pole and with the continuous branch,
\begin{equation}\label{S:eq:conditioned-angular-half-angle}
 \tan\frac{\Phi_r(u)}2=\me^{-r}\tan(\pi u).
\end{equation}
\end{lemma}
\begin{proof}
Let $\theta\in[0,2\pi]$ and define $\phi$ by $\tan(\phi/2)=\me^{-r}\tan(\theta/2)$, with the increasing branch fixing $0$, $\pi$, and $2\pi$. If $t=\tan(\phi/2)$, then
\[
 \frac{\dd\theta}{\dd\phi}
 =\frac{\me^r(1+t^2)}{1+\me^{2r}t^2},\qquad
 \cosh r-\cos\phi\sinh r
 =\frac{\me^{-r}+\me^r t^2}{1+t^2}.
\]
Hence $\dd\theta/\dd\phi=y(r,\phi)$ on each side of $\pi$. Both sides extend continuously across $\pi$ and the chosen function has $\theta(0)=0$. Integration gives
$\theta/(2\pi)=F_r(\phi)$, proving~\eqref{S:eq:conditioned-angular-half-angle}. To obtain the complex formula, put $E=\me^{\ii\theta}$ and $c=\me^{-r}$. Away from $\theta=\pi$,
$\ii\tan(\theta/2)=(E-1)/(E+1)$, so
\[
 \me^{\ii\phi}
 =\frac{1+\ii c\tan(\theta/2)}{1-\ii c\tan(\theta/2)}
 =\frac{(1+c)E+(1-c)}{(1-c)E+(1+c)}
 =\frac{E+\rho_r}{1+\rho_rE},\qquad
 \rho_r=\frac{1-c}{1+c}=\tanh(r/2).
\]
Continuity gives the same identity at $\theta=\pi$. Taking $\theta=2\pi u$ proves \eqref{S:eq:conditioned-angular-explicit}. The denominator has modulus at least $1-\rho_r>0$; thus
the expression is smooth in $(r,\me^{2\pi\ii u})$ on every bounded interval of $r$, including $r=0$.
\end{proof}

\subsubsection{The radial quantile comparison}
\label{S:sec:s6-radial-quantiles}
The preceding conditioning leaves the radial laws $F_h$ and $G_\gamma$ unchanged. We compare these laws at equal tail probability. For a uniform quantile estimate, we compare the
densities near zero and the logarithms of the tails farther out. The density comparison uses the scale $\sqrt h$.

We use the radial heat kernel formula recorded in \cite{valko2020operator}, where it is attributed to~\cite{KarpelevichTutubalinShur1959}:
\begin{equation}\label{S:eq:7.14}
 p_h(r)=\frac{\me^{-h/8}\sinh r}{\sqrt{\pi h^3}}
 \int_r^\infty\frac{s\me^{-s^2/(2h)}}{\sqrt{\cosh s-\cosh r}}\,\dd s.
\end{equation}
The density is with respect to $\dd r$. We give the scalar estimates needed for the quantile map, following the comparison of radial laws in~\cite{valko2020operator}. For
$0\le r\le s$,
\begin{equation}\label{S:eq:7.15}
\frac1{\sqrt2}\sqrt{s^2-r^2} \le \sqrt{\cosh s-\cosh r} \le \frac1{\sqrt2}\sqrt{s^2-r^2} \exp\!\left(\frac{r^2+s^2}{24}\right).
\end{equation}
For $v\ge0$,
\[
 1\le\frac{\sinh v}{v}\le \me^{v^2/6},
\]
with continuous values at zero. To prove the upper bound, put $g(v)=(1+v^2/3)\sinh v-v\cosh v$. Then $g'(v)=\frac v3(v\cosh v-\sinh v)\ge0$, because
$(v\cosh v-\sinh v)'=v\sinh v\ge0$. Thus $\coth v-v^{-1}\le v/3$; integration from zero proves the logarithm of the bound. Use
\[
 \cosh s-\cosh r
 =2\sinh\frac{s+r}{2}\sinh\frac{s-r}{2}
\]
to obtain \eqref{S:eq:7.15}, including coincident arguments by continuity. For $a>0$ the substitution $v=\sqrt{s^2-r^2}$ gives
\begin{equation}\label{S:eq:estimate-radial-gaussian-integrals}
 \int_r^\infty\frac{s \me^{-as^2}}{\sqrt{s^2-r^2}}\dd s =\frac{\sqrt\pi}{2\sqrt a}\me^{-ar^2},\qquad \int_r^\infty s\sqrt{s^2-r^2}\me^{-as^2}\dd s =\frac{\sqrt\pi}{4a^{3/2}}\me^{-ar^2}.
\end{equation}
For the lower density bound, use $a=1/(2h)+1/24$ and the additional factor $\sqrt2\me^{-r^2/24}$; for the upper bound, use $a=1/(2h)$ and $\sqrt2$. Substituting these quantities in
\eqref{S:eq:7.14} gives
\begin{equation}\label{S:eq:7.16}
p_h(r) \ge \left(1+\frac h{12}\right)^{-1/2} \frac1h \exp\!\left(-\frac{r^2}{2h}-\frac{r^2}{12}-\frac h8\right) \sinh r,
\end{equation}
and
\begin{equation}\label{S:eq:7.17}
p_h(r) \le \frac1h \exp\!\left(-\frac{r^2}{2h}-\frac h8\right) \sinh r.
\end{equation}
Integrating \eqref{S:eq:7.16} and using $\sinh u\ge u$ yields
\begin{equation}\label{S:eq:7.18}
\mathbb P(R_h>r) \ge a_h\exp\!\left(-\frac{r^2}{2h}-\frac{r^2}{12}\right),
\end{equation}
where
\[
a_h := \left(1+\frac h{12}\right)^{-1/2} \left(1+\frac h6\right)^{-1}\me^{-h/8}.
\]
Thus, for $0<h\le1$,
\begin{equation}\label{S:eq:7.20}
0\le-\log a_h\le Ch.
\end{equation}

For $\gamma(h)=2/h-1/2$ and $0<h\le1$, the comparison is
\begin{equation}\label{S:eq:7.21}
F_h(r)\ge G_{\gamma(h)}(r),   r\ge0.
\end{equation}
Tonelli applied to \eqref{S:eq:7.14} gives
\[
 \mathbb P(R_h>r)=\frac{2\me^{-h/8}}{\sqrt{\pi h^3}}
       \int_r^\infty s \me^{-s^2/(2h)}
                     \sqrt{\cosh s-\cosh r}\,\dd s,
\]
because the inner integral of $\sinh u/\sqrt{\cosh s-\cosh u}$ from $r$ to $s$ equals $2\sqrt{\cosh s-\cosh r}$. Substitute the upper bound in \eqref{S:eq:7.15} and the second
integral in \eqref{S:eq:estimate-radial-gaussian-integrals} with $a=1/(2h)-1/24>0$. The result is
\begin{equation}\label{S:eq:7.22}
\mathbb P(R_h>r) \le \left(1-\frac h{12}\right)^{-3/2} \exp\!\left(-\frac{r^2}{2h}+\frac{r^2}{12}-\frac h8\right),   0<h\le1.
\end{equation}
The circular tail is
\begin{equation}\label{S:eq:7.23}
\mathbb P(Y_{\gamma(h)}>r) = \exp\!\left[-\left(\frac4h-1\right) \log\cosh\!\left(\frac r2\right)\right].
\end{equation}
If $r\ge h$, then $\log\cosh x\le x^2/2$ and
\[
-\frac32\log\left(1-\frac h{12}\right)-\frac h8
 =\frac32\sum_{j\ge2}\frac{(h/12)^j}{j}
 \le\frac{3(h/12)^2}{4(1-h/12)}
 \le\frac{h^2}{24}\le\frac{r^2}{24},
\]
so the right side of \eqref{S:eq:7.22} is bounded by \eqref{S:eq:7.23}.

If $0\le r<h$, the lower density in \eqref{S:eq:7.16} dominates $g_{\gamma(h)}(r)$. After canceling $\sinh(r/2)$, this is equivalent to
\begin{equation}\label{S:eq:7.24}
\me^{-r^2/(2h)-r^2/12} \ge \me^{h/8}\sqrt{1+\frac h{12}} \left(1-\frac h4\right) \exp\!\left[-\left(\frac4h+1\right) \log\cosh\!\left(\frac r2\right)\right].
\end{equation}
For $0<h\le1$, the logarithm of the prefactor is at most $h/8+h/24-h/4=-h/12$, by $\log(1+x)\le x$ and $\log(1-x)\le-x$. Moreover, $0\le\tanh x\le x$ for $x\ge0$, because
$(\tanh x)'=\operatorname{sech}^2x\le1$. The derivative of $\tanh x-x+x^3/3$ is $x^2-\tanh^2x\ge0$, so integration gives
\[
 \frac{x^2}{2}-\frac{x^4}{12}\le\log\cosh x\le\frac{x^2}{2}.
\]
For $r<h\le1$, the lower bound gives
\[
\left(\frac4h+1\right) \left(\frac{r^2}{8}-\frac{r^4}{192}\right) \ge \frac{r^2}{2h}+\frac{r^2}{12}.
\]
Indeed, the difference between the left and right sides of the last inequality is $r^2/24-(4/h+1)r^4/192\ge r^2(8-4h-h^2)/192\ge0$. This proves \eqref{S:eq:7.24}. Integrating from
$0$ to $r$ proves \eqref{S:eq:7.21} in the remaining range. Thus $Y_{\gamma(h)}$ stochastically dominates $R_h$.

For $0<h<4$, define the quantile map
\[
\mathcal Q_h(r) := G_{\gamma(h)}^{-1}(F_h(r)),   \gamma(h)=\frac2h-\frac12.
\]
Since both radial densities are positive on $(0,\infty)$, the distribution functions are continuous and strictly increasing there. If $R_h$ has distribution function $F_h$, then
\begin{equation}\label{S:eq:7.26}
Y_h:=\mathcal Q_h(R_h)
\end{equation}
has distribution function $G_{\gamma(h)}$. When $0<h\le1$, \eqref{S:eq:7.21} gives $Y_h\ge R_h$ almost surely.

In particular, for $r>0$,
\[
 \mathcal Q_h(r)=2\operatorname{arcosh}
 \left((1-F_h(r))^{-1/(2\gamma(h))}\right).
\]
The integral formula for $p_h$ makes $F_h(r)$ jointly Borel in $(h,r)$; hence so is $\mathcal Q_h$. The probability integral transform proves~\eqref{S:eq:7.26} for every $0<h<4$.
Only the stochastic order and the estimates require $h$ to be small.

\begin{lemma}\label{S:lem:7.1}
There exist deterministic constants $h_0,c,C>0$ such that, for $0<h\le h_0$, the coupling \eqref{S:eq:7.26} satisfies
\begin{equation}\label{S:eq:7.28}
0\le Y_h-R_h \le C\bigl(h^{3/2}+Y_h^3\bigr)
\end{equation}
almost surely, and
\begin{equation}\label{S:eq:7.30}
\mathbb E\exp\!\left(c h^{-1}|Y_h-R_h|^{2/3}\right) \le C.
\end{equation}

\end{lemma}
\begin{proof}
Set $z:=R_h$ and $y:=Y_h$. Then $0\le z\le y$ and
\begin{equation}\label{S:eq:7.31}
F_h(z)=G_{\gamma(h)}(y).
\end{equation}
Take $h_0\le1$. Suppose first that $0\le y\le\sqrt h$ and define $q_h(r)=r\me^{-r^2/(2h)}/h$. For $0\le r\le\sqrt h$, the bounds $1\le\sinh r/r\le\me^{h/6}$ and
$\me^{-h/12}\le\me^{-r^2/12}\le1$, substituted in \eqref{S:eq:7.16}--\eqref{S:eq:7.17}, give
\begin{equation}\label{S:eq:7.32}
(1-Ch)q_h(r)\le p_h(r)\le(1+Ch)q_h(r).
\end{equation}
For the circular density, $2\gamma(h)+1=4/h$ gives
\[
\frac{g_{\gamma(h)}(r)}{q_h(r)} = \frac{(2-h/2)\sinh(r/2)}{r} \exp\!\left[ \frac{r^2}{2h}-\frac4h\log\cosh\!\left(\frac r2\right) \right].
\]
The inequalities above give $0\le r^2/(2h)-(4/h)\log\cosh(r/2)\le r^4/(48h)\le h/48$. Also $1\le2\sinh(r/2)/r\le\me^{h/24}$ on this range. Therefore
\begin{equation}\label{S:eq:7.33}
(1-Ch)q_h(r)\le g_{\gamma(h)}(r)\le(1+Ch)q_h(r).
\end{equation}
By \eqref{S:eq:7.31}, $F_h(y)-F_h(z)=F_h(y)-G_{\gamma(h)}(y).$ The stochastic order and \eqref{S:eq:7.32}--\eqref{S:eq:7.33} imply
\begin{equation}\label{S:eq:7.34}
0\le F_h(y)-G_{\gamma(h)}(y) \le Ch\int_0^y q_h(r)\,\dd r.
\end{equation}
Also, $\me^{-1/2}\le \me^{-r^2/(2h)}\le1$ on $[0,\sqrt h]$, so
\[
\int_z^yq_h(r)\,\dd r \ge c\frac{y^2-z^2}{h},   \int_0^yq_h(r)\,\dd r \le C\frac{y^2}{h}.
\]
After reducing $h_0$ so that $1-Ch_0\ge1/2$ in \eqref{S:eq:7.32}, combining these inequalities with \eqref{S:eq:7.34} yields $y^2-z^2\le Chy^2$. For $y>0$, division by $y+z\ge y$
gives the following bound; for $y=0$ it holds since $z=0$. Thus $0\le y-z\le Chy\le Ch^{3/2}.$

Suppose that $\sqrt h\le y\le1$. The quantile identity gives $\mathbb P(R_h>z)=\mathbb P(Y_h>y).$ Define
\[
L_h(y) :=-\log\mathbb P(Y_h>y) = \left(\frac4h-1\right)\log\cosh\!\left(\frac y2\right).
\]
The lower tail estimate \eqref{S:eq:7.18} gives $L_h(y) \le \frac{z^2}{2h}+\frac{z^2}{12}-\log a_h.$ Using \eqref{S:eq:7.20},
\begin{equation}\label{S:eq:7.36}
z^2 \ge \frac{2h}{1+h/6}\bigl(L_h(y)-Ch\bigr).
\end{equation}
For $0\le y\le1$, $\log\cosh(y/2)\ge y^2/8-y^4/192$. Substitution into~\eqref{S:eq:7.36} gives
\[
 z^2\ge\frac{1-h/4}{1+h/6}\left(y^2-\frac{y^4}{24}\right)-Ch^2
 \ge y^2-C(hy^2+y^4+h^2).
\]
Since $y+z\ge y$, $hy\le y^3$, and $h^2/y\le h^{3/2}$,
\[
 0\le y-z\le C\left(hy+y^3+\frac{h^2}{y}\right)
 \le C(h^{3/2}+y^3).
\]
If $y>1$, then $0\le y-z\le y\le y^3$. The three ranges give \eqref{S:eq:7.28}.

For $0\le y\le1$, the preceding polynomial lower bound gives $\log\cosh(y/2)\ge y^2/12$. For $y\ge1$, convexity and the value zero at the origin give
$\log\cosh(y/2)\ge y\log\cosh(1/2)$. Since $4/h-1\ge3/h$ for $h\le1$, formula~\eqref{S:eq:7.4} yields
\begin{equation}\label{S:eq:7.39}
\mathbb P(Y_h>y)
\le
\begin{cases}
\exp(-cy^2/h),&0\le y\le1,\\
C\exp(-cy/h),&y\ge1.
\end{cases}
\end{equation}
For every fixed $q\ge1$, integration of the tail probabilities gives
\[
 \begin{aligned}
 \mathbb EY_h^q
 &\le q\int_0^1 y^{q-1}\me^{-cy^2/h}\dd y
       +Cq\int_1^\infty y^{q-1}\me^{-cy/h}\dd y\\
 &\le (h/c)^{q/2}\Gamma(q/2+1)
       +C(h/c)^q\Gamma(q+1).
 \end{aligned}
\]
Since $h\le1$, this proves
\begin{equation}\label{S:eq:7.40}
 \mathbb EY_h^q\le C_q h^{q/2}.
\end{equation}
The constant in this last estimate may depend on $q$.

For \eqref{S:eq:7.30}, write $D=Y_h-R_h$ and choose $C_0$ such that $D^{2/3}\le C_0(h+Y_h^2)$ on $Y_h\le1$. For $0<a<c$, integration of the tail probabilities of the bounded
variable $Y_h\wedge1$ gives
\[
 \mathbb E \me^{a(Y_h\wedge1)^2/h}
 =1+\int_0^1\frac{2ay}{h}\me^{ay^2/h}\mathbb P(Y_h>y)\dd y
 \le1+\int_0^\infty\frac{2ay}{h}\me^{-(c-a)y^2/h}\dd y
 =1+\frac{a}{c-a}.
\]
Choose $d>0$ with $dC_0<c$. The contribution of $Y_h\le1$ to $\mathbb E \me^{dD^{2/3}/h}$ is therefore bounded uniformly in $h$. On $Y_h>1$, use $D\le Y_h$ and $Y_h^{2/3}\le Y_h$.
The second line of \eqref{S:eq:7.39} gives, for $d<c$,
\[
 \mathbb E[\me^{dY_h/h}\mathbf1_{\{Y_h>1\}}] =\me^{d/h}\mathbb P(Y_h>1) +\int_1^\infty\frac d h \me^{dy/h}\mathbb P(Y_h>y)\dd y \le C\left(1+\frac d{c-d}\right)\me^{-(c-d)/h}.
\]
Reducing $d$ to satisfy both restrictions proves \eqref{S:eq:7.30}.
\end{proof}

The exponential estimate also specifies the dependence on the moment order. Put $D=Y_h-R_h$. Markov's inequality gives $\mathbb P(D>t)\le C\exp(-ct^{2/3}/h)$. Thus, for every
$p\ge1$,
\[
 \mathbb E D^p
 \le Cp\int_0^\infty t^{p-1}\me^{-ct^{2/3}/h}\dd t
 =C(h/c)^{3p/2}\Gamma(3p/2+1).
\]
For $s\ge1$, splitting $\me^{-t}$ into two factors and maximizing $t^s\me^{-t/2}$ gives $\Gamma(s+1)\le2(2s/\me)^s$. Taking $p$th roots therefore yields
\begin{equation}\label{S:eq:7.29}
 \|Y_h-R_h\|_{L^p}\le Cp^{3/2}h^{3/2},\qquad
 \mathbb E|Y_h-R_h|^p\le C_p h^{3p/2}.
\end{equation}
The product argument uses $p$ of order $\log n/\log\log n$. The polynomial dependence on $p$ in the displayed estimate then yields the probability bound required for all starting
indices in a row.
\subsubsection{Construction of all finite walks}
\label{S:sec:s6-exact-walk-construction}
The quantile comparison is useful only when $k$ is large. For the indices below the smaller cutoff we use an exact circular continuation, whose distances from its starting position
are controlled simultaneously in $n$. The equivariance below allows us to change its starting position without changing its circular transition laws.

For a fixed $\eta$, let $\Psi_\eta(z,w)=g_{z,\eta}^{-1}\cdot w$. Let $K_{n,k}$ be the circular transition with index $k$, with radial parameter $\gamma_k=\beta k/2$, and put
\[
 \varkappa_{n,k}(\dd w)=\Im(w)K_{n,k}(\ii,\dd w).
\]
This is a probability measure by~\eqref{S:eq:7.42}. Fix $1\le m\le n$. Given independent $W_k$ with laws $\varkappa_{n,k}$ for $1\le k<m$, and a starting point $z_m$, define
$z_k=\Psi_\eta(z_{k+1},W_k)$ in the order
$k=m-1,\ldots,1$.

\begin{lemma}\label{S:lem:conditioned-tail-splicing}
This recursion has transition kernels $K_{n,k}^\eta$. For a second starting point $z_m'$, the recursion with the same $W_k$ satisfies
\begin{equation}\label{S:eq:tail-splicing-equivariance}
 z_k'=hz_k\quad(1\le k\le m),\qquad
 h=g_{z_m',\eta}^{-1}g_{z_m,\eta}.
\end{equation}
\end{lemma}
\begin{proof}
Lemma~\ref{S:lem:7.2} identifies the pushforward of $\varkappa_{n,k}$ under $w\mapsto\Psi_\eta(z,w)$ with $K_{n,k}^\eta(z,\cdot)$. Since $z_{k+1}$ is a function only of $W_j$ with
$j>k$, it is independent of $W_k$. For every bounded Borel $f$,
\[
 \mathbb E\bigl[f(z_k)\mid W_{k+1},\ldots,W_{m-1}\bigr]
 =\int f(\Psi_\eta(z_{k+1},w))\,\varkappa_{n,k}(\dd w)
 =\int f(z')K_{n,k}^\eta(z_{k+1},\dd z').
\]
Iterating this conditional identity gives the product of the kernels $K_{n,k}^\eta$.

The isometry $h$ fixes $\eta$ and sends $z_m$ to $z_m'$. It is unique: for $\eta=\infty$, the stabilizer consists of $z\mapsto az+b$, $a>0$, $b\in\mathbb R$, and the image of one
point determines $a$ and $b$. Conjugation gives uniqueness for other boundary points. Both $g_{hz,\eta}h$ and $g_{z,\eta}$ have the defining images $z\mapsto\ii$ and
$\eta\mapsto\infty$. Their equality therefore implies
\[
 \Psi_\eta(hz,w)=g_{hz,\eta}^{-1}\cdot w
 =h g_{z,\eta}^{-1}\cdot w=h\Psi_\eta(z,w).
\]
Starting at $k=m$ and applying this equality at each subsequent step proves~\eqref{S:eq:tail-splicing-equivariance}. In particular, all distances from the starting point are
preserved. For random $z_m,z_m'$ independent of the local variables, conditioning on $(z_m,z_m')$ leaves the product law of the $W_k$ unchanged. The conditional transition identity
and the equality $z_k'=hz_k$ therefore hold under this conditional law; integration over $(z_m,z_m')$ gives both conclusions.
\end{proof}

We now use the Brownian path $Z$ and its terminal point $\eta_\infty$. For every $n$ and $1\le k<n$, set
\[
 z_{n,k}^{\rm br}=Z_{\vartheta_{n,k}},\qquad
 w_{n,k}^{\rm br}=g_{z_{n,k+1}^{\rm br},\eta_\infty}
                      \cdot z_{n,k}^{\rm br}.
\]
The local variable has polar coordinates $(R_{n,k},\phi_{n,k})$. Take the angle in $[0,2\pi)$, and put $U_{n,k}=F_{R_{n,k}}(\phi_{n,k})$. The radius is nonzero almost surely and
the angular law has a positive density; arbitrary choices at radius zero and at the angular boundary point change no distribution. These definitions can be made simultaneously
since there are only countably many $(n,k)$. For a Borel set $A\subset(0,\infty)$ and $0\le u\le1$, Lemma~\ref{S:lem:7.2} gives
\[
 \mathbb P(R_{n,k}\in A,\ U_{n,k}\le u\mid\eta_\infty)
 =\int_A p_{h_k}(r)\left(\int_0^{\Phi_r(u)}
               \frac{y(r,\phi)}{2\pi}\,\dd\phi\right)\dd r
 =u\int_A p_{h_k}(r)\,\dd r.
\]
Together with that lemma's product identity this proves
\begin{equation}\label{S:eq:s6-row-product-law}
 \mathcal L\bigl((R_{n,k},U_{n,k})_{1\le k<n}\mid\eta_\infty\bigr)
 =\bigotimes_{k=1}^{n-1}
       \bigl(p_{h_k}(r)\,\dd r\otimes\dd u\bigr).
\end{equation}
These variables are measurable functions of the original Brownian path. Their product law is conditional on the terminal point, for each fixed row.

Choose a deterministic $n_0=n_0(\beta)$ such that $1\le\underline m_n<n$ and $h_{\underline m_n}\le h_0$ for $n\ge n_0$. For $\underline m_n\le k<n$, put
\begin{equation}\label{S:eq:7.46}
 \begin{gathered}
 \Phi_{n,k}^{\rm br}=\Phi_{R_{n,k}}(U_{n,k}),\qquad
 Y_{n,k}=\mathcal Q_{h_k}(R_{n,k}),\qquad
 \Phi_{n,k}^{\rm circ}=\Phi_{Y_{n,k}}(U_{n,k}),\\
 w_{n,k}^{\rm circ}=z(Y_{n,k},\Phi_{n,k}^{\rm circ}).
 \end{gathered}
\end{equation}
The transfers to be used in the product estimates are
\begin{equation}\label{S:eq:7.47}
 \widehat H_{n,k}^{\rm br}=\mathsf A(R_{n,k},\Phi_{n,k}^{\rm br}),
 \qquad
 \widehat H_{n,k}^{\rm circ}=\mathsf A(Y_{n,k},\Phi_{n,k}^{\rm circ}).
\end{equation}
Both are upper triangular with positive diagonal. Their variables have the two required radial laws and, at their respective radii, the required angular laws. The same uniform
quantile is used for the two angles; the angles themselves are usually different.

For the remaining increments take the simultaneous stopped array of~\cite{valko2020operator}. Write its terminal point as $\eta^\star$ and its finite positions as $z_{n,k}^\star$,
$1\le k\le n$. We retain the array's joint law so that the oscillation bound holds simultaneously for all $n\ge n_0(\omega)$. The array space
\[
 \mathsf A=\prod_{n\ge1}\mathbb H^n
\]
is a standard Borel space. Let $\mathsf Q_\eta$ be a regular conditional law of this array given $\eta^\star=\eta$. Fix the deterministic constant $C_\beta^{\rm tail}$ allowed by
Lemma~\ref{S:lem:auxiliary-boundary-oscillation} and set $M_n=C_\beta^{\rm tail}(1+\log\log n)^4$ for $n\ge n_0$, increasing $n_0$ to at least $3$. That lemma proves that the Borel
event
\[
 E=\bigcup_{N\ge n_0}\ \bigcap_{n\ge N}
 \left\{\max_{1\le k\le\underline m_n}
 d_{\mathbb H}(z_{n,k}^\star,z_{n,\underline m_n}^\star)\le M_n\right\}
\]
has probability one under the joint law of the array constructed at stopping times. It follows that
\[
 0=\int\mathsf Q_\eta(E^c)\,\nu_\ii(\dd\eta),\qquad
 \mathsf Q_\eta(E)=1\quad\text{for }\nu_\ii\text{-almost every }\eta.
\]
Intersect this full measure set of $\eta$ with those on which the conditional law of each row is the one in Lemma~\ref{S:lem:7.2}. There are only countably many rows and a
countable generating family of their Borel sets, so one exceptional null set suffices. The kernel can be defined arbitrarily on that set.

Lemma~4.22 of \cite{kallenberg2021foundations} provides a Borel map $F$ such that, for every bounded Borel function $f$,
\[
 \int_0^1 f(F(\eta,v))\,\dd v=\int f(a)\,\mathsf Q_\eta(\dd a).
\]
Thus a single uniform variable samples the entire conditional array with its joint law $\mathsf Q_\eta$.

Adjoin an independent uniform $V$ and use $F(\eta_\infty,V)$ as the auxiliary array. Conditional on the entire main path, its expectation against $f$ is $\mathsf Q_{\eta_\infty}f$.
Consequently, for a bounded function $g$ of that path,
\[
 \mathbb E[f(F(\eta_\infty,V))g(Z)\mid\eta_\infty]
 =\mathsf Q_{\eta_\infty}f\,
       \mathbb E[g(Z)\mid\eta_\infty].
\]
This proves conditional independence of the Brownian path and auxiliary array given the terminal point while preserving the joint law of the auxiliary rows. In particular,
\begin{equation}\label{S:eq:auxiliary-tail-oscillation}
 \max_{1\le k\le\underline m_n}
 d_{\mathbb H}(z_{n,k}^\star,z_{n,\underline m_n}^\star)\le M_n
\end{equation}
holds almost surely for all sufficiently large $n$, by Lemma~\ref{S:lem:auxiliary-boundary-oscillation}.

Extract from this auxiliary array $W_{n,k}^\star=g_{z_{n,k+1}^\star,\eta_\infty}\cdot z_{n,k}^\star$ for $1\le k<\underline m_n$. For $n\ge n_0$ define $z_{n,n}^{\rm circ}=\ii$
and, in decreasing order of $k$, set
\[
 z_{n,k}^{\rm circ}=
 \begin{cases}
 \Psi_{\eta_\infty}(z_{n,k+1}^{\rm circ},w_{n,k}^{\rm circ}),
          &\underline m_n\le k<n,\\
 \Psi_{\eta_\infty}(z_{n,k+1}^{\rm circ},W_{n,k}^\star),
          &1\le k<\underline m_n.
 \end{cases}
\]
Set $z_{n,0}^{\rm circ}=\eta_\infty$. For the finitely many $n<n_0$, take the entire auxiliary row and the same terminal point. This also defines the case $n=1$, in which there is
no interior jump. Finally, the coefficient on $I_{n,k}$ is $\mathcal A(z_{n,k+1}^{\rm circ})$, in the original upper half-plane coordinates.

\begin{proposition}\label{S:prop:7.2}
Each constructed row has the joint law of the circular walk and its terminal point specified in Section~\ref{S:sec:introduction}, with $z_{n,0}^{\rm circ}=\eta_\infty$. Almost
surely, for all sufficiently large $n$,
\begin{equation}\label{S:eq:spliced-tail-oscillation}
 \max_{1\le k\le\underline m_n}
 d_{\mathbb H}(z_{n,k}^{\rm circ},z_{n,\underline m_n}^{\rm circ})
 \le M_n.
\end{equation}
\end{proposition}
\begin{proof}
Condition on $\eta_\infty=\eta$ outside the exceptional set chosen above. Equation~\eqref{S:eq:s6-row-product-law} and the two quantile maps show that $w_{n,k}^{\rm circ}$ have
laws $\varkappa_{n,k}$ and are independent within the part under consideration of the row. More explicitly, for bounded Borel $f$ and a selected index, their joint law with the
prescribed Brownian local variable is
\[
 \mathbb E[f(w_{n,k}^{\rm br},w_{n,k}^{\rm circ})\mid\eta]
 =\int_0^\infty\!\int_0^1
 f\bigl(z(r,\Phi_r(u)),z(\mathcal Q_{h_k}(r),
                      \Phi_{\mathcal Q_{h_k}(r)}(u))\bigr)
                                    \,\dd u\,p_{h_k}(r)\,\dd r.
\]
Taking $f$ to depend on the second variable alone, change the radial distribution by $\mathcal Q_{h_k}$ and use the angular law $\pi_y$ at the new radius $y$. The result is
$g_{\gamma_k}(y)\,\dd y\,\pi_y(\dd\phi)$, precisely $\varkappa_{n,k}$. Taking $f$ to depend only on the first variable leaves the Brownian law unchanged. The same calculation for a
product of test functions over the selected indices proves the product distribution for each fixed row. The two local variables at a fixed index retain the dependence of their
coupling. The extracted $W_{n,k}^\star$ have these same local laws by Lemma~\ref{S:lem:7.2}. They are independent within their part of the row, and the conditional independence of
the auxiliary array makes them independent of the retained variables. Thus Lemma~\ref{S:lem:conditioned-tail-splicing} gives the conditional path law
\[
 \delta_\ii(\dd z_n)
 K_{n,n-1}^\eta(z_n,\dd z_{n-1})\cdots
 K_{n,1}^\eta(z_2,\dd z_1).
\]
Integration with $\nu_\ii(\dd\eta)$ gives precisely
\begin{equation}\label{S:eq:s6-unconditioned-walk-law}
 \nu_\ii(\dd\eta)\delta_\ii(\dd z_n) \prod_{k=n-1}^{1}K_{n,k}^\eta(z_{k+1},\dd z_k) = \delta_\ii(\dd z_n)\prod_{k=n-1}^{1}K_{n,k}(z_{k+1},\dd z_k) \nu_{z_1}(\dd\eta).
\end{equation}
The product of the Poisson factors is $P_\eta(z_1)$, and $P_\eta(z_1)\nu_\ii(\dd\eta)=\nu_{z_1}(\dd\eta)$. The products here are ordered integrations in decreasing $k$. The right
side is the required unconditional distribution of the finite circular walk: at each step the radial beta variable and the uniform local direction are chosen independently of the
preceding steps, and the last boundary point has harmonic measure from $z_1$. The construction of the small rows has the same law by definition. Applying the coefficient map on
$I_{n,k}$ therefore gives exactly the circular operator law of~\cite{valko2020operator}, with the normalization checked in Proposition~\ref{S:prop:operator-normalization}.

For the estimate, define the Borel isometry
\[
 h_n=g_{z_{n,\underline m_n}^{\rm circ},\eta_\infty}^{-1}
              g_{z_{n,\underline m_n}^\star,\eta_\infty}.
\]
Both paths below the cutoff use the same $W_{n,k}^\star$. Equation~\eqref{S:eq:tail-splicing-equivariance} gives $z_{n,k}^{\rm circ}=h_nz_{n,k}^\star$ for every
$1\le k\le\underline m_n$. Preservation of hyperbolic distance turns~\eqref{S:eq:auxiliary-tail-oscillation} into \eqref{S:eq:spliced-tail-oscillation}, on the same event and with
the same lower bound on the matrix size. This is an extension of the original Brownian probability space: conditional on its whole path, the auxiliary array is sampled from
$\mathsf Q_{\eta_\infty}$, and the new row is a Borel function of that path and the auxiliary sample. Because this conditional kernel has total mass one, the enlargement preserves
the full Brownian path law. The pathwise isometry transfers the auxiliary oscillation estimate without changing this joint construction.
\end{proof}

After the common coordinate change $G_{\eta_\infty}$, the constructed paths have the transfers~\eqref{S:eq:7.47}. For either path, set
$\widehat z_{n,k}^{\diamond}=G_{\eta_\infty}\cdot z_{n,k}^{\diamond}$ and $\widehat X_{n,k}^{\diamond}=\mathcal A(\widehat z_{n,k}^{\diamond})$. By the formula for $g_{z,\eta}$,
\[
 w_{n,k}^{\diamond}
 =\mathcal A(\widehat z_{n,k+1}^{\diamond})
                         \cdot\widehat z_{n,k}^{\diamond}.
\]
Equation~\eqref{S:eq:affine-section-exact-composition} therefore gives
\begin{equation}\label{S:eq:s6-conditioned-transfer-recursion}
 \widehat X_{n,k}^{\diamond}
 =\widehat H_{n,k}^{\diamond}\widehat X_{n,k+1}^{\diamond},
 \qquad \widehat X_{n,n}^{\diamond}=I.
\end{equation}
The terminal line is now $\mathbb R e_1$, and these matrices are used in Section~\ref{S:sec:products-of-the-coupled-transfer-matrices}.
Section~\ref{S:sec:comparison-of-the-finite-and-sampled-brownian-operators} returns to the original upper half-plane coordinates, keeping the two orthogonal factors depending on
the two positions.
\subsubsection{Logarithms of the relative matrices}\label{S:sec:s6-relative-logarithms}
The relative matrix depends on the radius through both $r$ and $\Phi_r(u)$. Its derivative therefore includes the angular derivative $\partial_r\Phi_r(u)$.

Put $\theta=2\pi u$ and define
\[
 v(r,\theta)=\cosh r+\cos\theta\sinh r,\qquad
 b(r,\theta)=\sin\theta\sinh r,\qquad
 \mathfrak F(r,u)=\mathsf A(r,\Phi_r(u)).
\]
Substitution in the half angle formula yields
\[
 y(r,\Phi_r(u))=v(r,\theta),\qquad
 x(r,\Phi_r(u))=-b(r,\theta),
\]
and therefore
\begin{equation}\label{S:eq:s6-affine-quantile-matrix}
 \mathfrak F(r,u)=
 \begin{pmatrix}
 v(r,\theta)^{-1/2}&b(r,\theta)v(r,\theta)^{-1/2}\\
 0&v(r,\theta)^{1/2}
 \end{pmatrix}.
\end{equation}
For example, putting $t=\tan(\theta/2)$ gives $y=\me^r(1+\me^{-2r}t^2)/(1+t^2)$ and $x=-2t\sinh r/(1+t^2)$, which are precisely these formulas. Continuity covers $\theta=\pi$.

We use the real logarithm of upper triangular matrices with determinant one and positive diagonal entries. Its explicit formula is
\begin{equation}\label{S:eq:s6-triangular-logarithm}
 \log\begin{pmatrix}\me^a&d\\0&\me^{-a}\end{pmatrix}
 =\begin{pmatrix}a&d\,a/\sinh a\\0&-a\end{pmatrix},
 \qquad \left.\frac a{\sinh a}\right|_{a=0}=1.
\end{equation}
Indeed, the square of the matrix on the right is $a^2I$, so its exponential is the scalar hyperbolic cosine times the identity plus $\sinh a/a$ times that matrix. This proves the
formula also at $a=0$ by continuity. It specifies the logarithm globally and agrees with the logarithm defined by its power series near the identity.

For $0\le r\le y$ define
\[
 \mathcal R(r,y,u)=
 \log\bigl(\mathfrak F(y,u)\mathfrak F(r,u)^{-1}\bigr)
 -(y-r)\mathsf V_0(\Phi_r(u)).
\]

\begin{lemma}\label{S:lem:localized-affine-relative}
There are absolute $r_{\rm loc}\in(0,1]$ and $C<\infty$ such that, for $0\le r\le y\le r_{\rm loc}$,
\begin{equation}\label{S:eq:7.49}
 \|\mathcal R(r,y,u)\|
 \le C\bigl((y-r)r^2+(y-r)^2\bigr)
 \le C\bigl((y-r)r+(y-r)^2\bigr)
\end{equation}
uniformly in $u\in(0,1)$. Moreover,
\begin{equation}\label{S:eq:localized-transfer-log-bound}
 \sup_{0<u<1}\|\log\mathfrak F(r,u)\|\le Cr,
 \qquad 0\le r\le r_{\rm loc}.
\end{equation}
\end{lemma}
\begin{proof}
We keep $\theta$ fixed when differentiating \eqref{S:eq:s6-affine-quantile-matrix}. Write $v=v(r,\theta)$ and $b=b(r,\theta)$. Then
\[
 v''=v,\quad v^2-(v')^2=\sin^2\theta,\quad
 |v'|\le v,\quad \me^{-r}\le v\le\me^r,\quad
 b'=\sin\theta\cosh r,\quad b''=b.
\]
For $0\le r\le1$, these identities bound $v$, $v^{-1}$, $b$, and their first two derivatives by absolute constants. The derivative formulas are
\[
 (v^{-1/2})''=\frac34v^{-5/2}(v')^2-\frac12v^{-1/2},\qquad
 (v^{1/2})''=\frac12v^{1/2}-\frac14v^{-3/2}(v')^2,
\]
and $(bv^{-1/2})''=b''v^{-1/2}+2b'(v^{-1/2})'+b(v^{-1/2})''$. Thus
\[
 \sup_{0\le r\le1,\,0<u<1}
 \bigl(\|\mathfrak F(r,u)^{-1}\|
       +\|\partial_r\mathfrak F(r,u)\|
       +\|\partial_r^2\mathfrak F(r,u)\|\bigr)\le C.
\]

Multiplication by the inverse in~\eqref{S:eq:s6-affine-quantile-matrix} gives the exact logarithmic derivative
\[
 L(r,u):=\partial_r\mathfrak F(r,u)\mathfrak F(r,u)^{-1}
 =-\frac{v'}vH_0+\frac{\sin\theta\cosh r}{v}E_{12}.
\]
The half angle formula also gives $\cos\Phi_r(u)=v'/v$ and $\sin\Phi_r(u)=\sin\theta/v$. Consequently,
\begin{equation}\label{S:eq:localized-log-derivative}
 L(r,u)-\mathsf V_0(\Phi_r(u))
 =\frac{\sin\theta(\cosh r-1)}{v(r,\theta)}E_{12},\qquad
 \|L(r,u)-\mathsf V_0(\Phi_r(u))\|\le Cr^2.
\end{equation}
This derivative includes the dependence of $\Phi_r(u)$ on $r$. To compare it with differentiation at fixed angle, write
\[
 \mathsf A(r,\phi)=
 \begin{pmatrix}
 D^{1/2}&\sin\phi\sinh r\,D^{-1/2}\\0&D^{-1/2}
 \end{pmatrix},\qquad D=\cosh r-\cos\phi\sinh r,
\]
Taylor expansion at $r=0$ gives
\begin{equation}\label{S:eq:affine-radial-first-derivative}
 \mathsf A(r,\phi)=I+r\mathsf V_0(\phi)+E(r,\phi),
 \qquad\sup_\phi\|E(r,\phi)\|\le Cr^2\quad(0\le r\le1).
\end{equation}
The remainder follows from the displayed formulas for the second derivatives, using $D\ge\me^{-1}$, $D''=D$, and
$|D'|\le D$.

Put $\delta=y-r$. Taylor's formula in integral form gives
\[
 \mathfrak F(y,u)\mathfrak F(r,u)^{-1}
 =I+\delta L(r,u)
 +\delta^2\int_0^1(1-t)\partial_r^2\mathfrak F(r+t\delta,u)
                       \mathfrak F(r,u)^{-1}\,\dd t.
\]
The last term has norm at most $C\delta^2$, and the whole difference from $I$ has norm at most $C\delta$. Choose $r_{\rm loc}$ so that this norm is at most $1/2$ whenever
$y\le r_{\rm loc}$. For $\|Q\|\le1/2$ the absolutely convergent logarithm series gives
\[
 \|\log(I+Q)-Q\|
 \le\sum_{j\ge2}\frac{\|Q\|^j}{j}
 \le\frac{\|Q\|^2}{2(1-\|Q\|)}\le\|Q\|^2.
\]
Hence $\|\log(\mathfrak F(y,u)\mathfrak F(r,u)^{-1})-\delta L(r,u)\| \le C\delta^2$. Substitution of \eqref{S:eq:localized-log-derivative} gives \eqref{S:eq:7.49}. Finally,
$\mathfrak F(0,u)=I$ and the bound for the first derivative gives $\|\mathfrak F(r,u)-I\|\le Cr$. Reducing $r_{\rm loc}$ once more and using $\|\log(I+Q)\|\le2\|Q\|$ proves
\eqref{S:eq:localized-transfer-log-bound}.
\end{proof}

Integrating over the uniform angular quantile before estimating the radial difference gains a factor $h^{1/2}$ in the mean. The pointwise estimate controls the centered part of the
increment. Angular integration also makes the off-diagonal mean vanish exactly. Multiplication gives
\[
 \mathfrak F(y,u)\mathfrak F(r,u)^{-1}
 =\begin{pmatrix}
 \sqrt{v(r,\theta)/v(y,\theta)}&
 \displaystyle\frac{\sin\theta(\sinh y-\sinh r)}
                      {\sqrt{v(y,\theta)v(r,\theta)}}\\
 0&\sqrt{v(y,\theta)/v(r,\theta)}
 \end{pmatrix}.
\]
Its logarithmic off-diagonal entry is odd under $\theta\mapsto2\pi-\theta$, by~\eqref{S:eq:s6-triangular-logarithm}. For the diagonal entries, use
\[
 v(r,\theta)=\cosh^2(r/2)
        |1+\tanh(r/2)\me^{\ii\theta}|^2.
\]
For fixed $r<\infty$ the power series for $\log(1+\tanh(r/2)\me^{\ii\theta})$ converges absolutely and uniformly in $\theta$. Each of its nonconstant Fourier terms has integral
zero. Hence
\[
 \frac1{2\pi}\int_0^{2\pi}\log v(r,\theta)\,\dd\theta
 =2\log\cosh(r/2).
\]
Together these two computations prove
\begin{equation}\label{S:eq:s6-exact-angular-log-mean}
 \int_0^1\log\bigl(\mathfrak F(y,u)\mathfrak F(r,u)^{-1}\bigr)\,\dd u
 =-2\bigl(\log\cosh(y/2)-\log\cosh(r/2)\bigr)H_0.
\end{equation}
The integrand is bounded for each fixed finite pair of radii, so all entries are integrable. We will use the identity only for $0\le r\le y\le r_{\rm loc}$.

\subsubsection{Conditional moments of one increment}
\label{S:sec:s6-increment-moments}
The logarithms below are those of the exact coupled transfers. We estimate them on a fixed radial ball, where the Taylor bounds are uniform. The circular walk itself keeps its
exact law, and the discarded part remains in the decomposition until its probability has been summed over the rows.

On the range under consideration define
\[
 \begin{gathered}
 \zeta_{n,k}^{\rm full}
 =\log\bigl(\widehat H_{n,k}^{\rm circ}
                     (\widehat H_{n,k}^{\rm br})^{-1}\bigr),
 \qquad \delta_{n,k}=Y_{n,k}-R_{n,k},\\
 \chi_{n,k}=\mathbf1_{\{Y_{n,k}\le r_{\rm loc}\}},\qquad
 \zeta_{n,k}^{\rm loc}=\chi_{n,k}\zeta_{n,k}^{\rm full},\qquad
 \zeta_{n,k}^{\rm bad}=(1-\chi_{n,k})\zeta_{n,k}^{\rm full}.
 \end{gathered}
\]
The logarithm is the one in~\eqref{S:eq:s6-triangular-logarithm}. Choose $k_0(\beta)$ so that $h_k\le h_0$ for $k\ge k_0(\beta)$. Then $0\le R_{n,k}\le Y_{n,k}$, and $\chi_{n,k}=1$
puts both radii in the range of Lemma~\ref{S:lem:localized-affine-relative}. The original Brownian variables at all indices of the row define
\begin{equation}\label{S:eq:conditioning-row-filtration}
 \mathcal F_{n,k}^{\eta}
 =\sigma\bigl(\eta_\infty,(R_{n,j},U_{n,j}):k\le j<n\bigr),
 \qquad\mathcal F_{n,n}^{\eta}=\sigma(\eta_\infty).
\end{equation}
The filtration increases as $k$ decreases, in increasing Brownian time. It is defined by the original Brownian variables at every index. Put
\[
 b_{n,k}=\mathbb E[\zeta_{n,k}^{\rm loc}\mid\mathcal F_{n,k+1}^{\eta}],
 \qquad \xi_{n,k}=\zeta_{n,k}^{\rm loc}-b_{n,k}.
\]
All these conditional expectations exist, since the localized logarithm is bounded. The exact decomposition is
\begin{equation}\label{S:eq:7.50}
 \zeta_{n,k}^{\rm full}=b_{n,k}+\xi_{n,k}+\zeta_{n,k}^{\rm bad}.
\end{equation}
Equation~\eqref{S:eq:s6-row-product-law} shows that the conditional law of the $k$th pair, given $\mathcal F_{n,k+1}^{\eta}$, is the law of $(R_{n,k},U_{n,k})$ conditional on the
terminal point. In particular, $b_{n,k}$ is a deterministic function of $(\beta,k)$.

\begin{theorem}\label{S:thm:7.4}
For every retained $k\ge k_0(\beta)$ and $p\ge2$,
\begin{equation}\label{S:eq:7.51}
 \|b_{n,k}\|\le C_\beta h_k^2\le C_\beta k^{-2},
\end{equation}
\begin{equation}\label{S:eq:localized-one-step-moment}
 \mathbb E[\|\xi_{n,k}\|^p\mid\mathcal F_{n,k+1}^{\eta}]^{1/p}
 \le C_\beta p^{3/2}k^{-3/2}.
\end{equation}
The variables $\xi_{n,k}$ are independent and centered conditionally on $\eta_\infty$, within each row.
\end{theorem}
\begin{proof}
Use the fixed conditional law of $(R_{n,k},U_{n,k})$ and abbreviate $h=h_k$, $R=R_{n,k}$, $Y=Y_{n,k}$, and $\delta=Y-R$. The radial estimate~\eqref{S:eq:7.29} gives
\begin{equation}\label{S:eq:localized-radial-moment}
 \mathbb E[\delta^p\mid\mathcal F_{n,k+1}^{\eta}]^{1/p}
 \le Cp^{3/2}h^{3/2},\qquad p\ge1.
\end{equation}
On $\{\chi_{n,k}=1\}$, Lemma~\ref{S:lem:localized-affine-relative} and boundedness of $\mathsf V_0$ give
\begin{equation}\label{S:eq:localized-zeta-delta-bound}
 \|\zeta_{n,k}^{\rm loc}\|\le C\delta.
\end{equation}
Conditional Jensen and Minkowski therefore imply
\[
 \|\xi_{n,k}\|_{L^p(\mathbb P(\cdot\mid\mathcal F_{n,k+1}^{\eta}))}
 \le2\|\zeta_{n,k}^{\rm loc}\|_{L^p(\mathbb P(\cdot\mid
                                          \mathcal F_{n,k+1}^{\eta}))}
 \le Cp^{3/2}h^{3/2}.
\]
Since $h_k\asymp_\beta k^{-1}$, this is \eqref{S:eq:localized-one-step-moment}.

To estimate the mean, condition first on $R$. Both $Y=\mathcal Q_h(R)$ and $\chi_{n,k}$ are then fixed, and the quantile variable is uniform. The exact
mean~\eqref{S:eq:s6-exact-angular-log-mean} gives
\begin{equation}\label{S:eq:s6-exact-local-drift}
 b_{n,k}=-2H_0\,
 \mathbb E\left[\chi_{n,k}
     \bigl(\log\cosh(Y/2)-\log\cosh(R/2)\bigr)\right].
\end{equation}
The off-diagonal drift is zero. Since $(\log\cosh(t/2))'=\tfrac12\tanh(t/2)\le t/4$ for $t\ge0$, for $0\le R\le Y$,
\[
 0\le\log\cosh(Y/2)-\log\cosh(R/2)
 =\frac12\int_R^Y\tanh(t/2)\,\dd t
 \le\frac{Y^2-R^2}{8}=\frac{\delta(2R+\delta)}8.
\]
By~\eqref{S:eq:7.40}, $\|R\|_2\le\|Y\|_2\le Ch^{1/2}$; by~\eqref{S:eq:7.29}, $\|\delta\|_2\le Ch^{3/2}$. Thus
\[
 \mathbb E(\delta R)\le\|\delta\|_2\|R\|_2\le Ch^2,
 \qquad \mathbb E\delta^2\le Ch^3.
\]
Since $\|H_0\|=1/2$, substitution in \eqref{S:eq:s6-exact-local-drift} proves~\eqref{S:eq:7.51}. The Brownian transfer matrix and its error depend on the same increment, so their
products must be estimated under their joint law. Write $\mathbb E_k$ for conditional expectation given $\mathcal F_{n,k+1}^{\eta}$. On $\{\chi_{n,k}=1\}$ the bound for the
Brownian logarithm is $C R$. Hence
\[
 \mathbb E_k\!\bigl[\chi_{n,k}
       \|\log\widehat H_{n,k}^{\rm br}\|\,\|\xi_{n,k}\|\bigr]
 \le C\bigl(\mathbb E_k[R\delta]
                       +\|b_{n,k}\|\mathbb E_k R\bigr)\le Ch^2.
\]
Cauchy--Schwarz bounds the first term by $(\mathbb E_kR^2)^{1/2}(\mathbb E_k\delta^2)^{1/2}$. Also $\zeta_{n,k}^{\rm loc}$ is uniformly bounded on its support, so the exponential
series and $\mathbb E_k\delta^2\le Ch^3$ give
\[
 \bigl\|\mathbb E_k[\exp(\zeta_{n,k}^{\rm loc})-I]
                          -b_{n,k}\bigr\|
 \le\mathbb E_k\|\exp(\zeta_{n,k}^{\rm loc})-I
                              -\zeta_{n,k}^{\rm loc}\|\le Ch^3.
\]
The expectation of the relative transfer matrix is consequently $I+b_{n,k}+O(h^3)$. The product estimate below includes both $b_{n,k}$ and the quadratic contribution from the
exponential.

Finally, all localized variables are measurable functions of their own pair $(R_{n,k},U_{n,k})$. The product law \eqref{S:eq:s6-row-product-law} and deterministic centering prove
conditional independence and zero mean.
\end{proof}

The circular radial tail gives
\begin{equation}\label{S:eq:7.55}
 \mathbb P(\chi_{n,k}=0\mid\mathcal F_{n,k+1}^{\eta})
 =\operatorname{sech}^{\beta k}(r_{\rm loc}/2)
 \le\me^{-c_\beta k},\qquad
 c_\beta=\beta\log\cosh(r_{\rm loc}/2)>0.
\end{equation}
The second moment bound is the case $p=2$ of the theorem:
\begin{equation}\label{S:eq:7.53}
 \mathbb E[\|\xi_{n,k}\|^2\mid\mathcal F_{n,k+1}^{\eta}]
 \le C_\beta k^{-3}.
\end{equation}
For the exponential estimate, subadditivity of $t^{2/3}$ and \eqref{S:eq:localized-zeta-delta-bound} imply
\[
 k\|\xi_{n,k}\|^{2/3}
 \le Ck\delta_{n,k}^{2/3}+Ck\|b_{n,k}\|^{2/3}
 \le C_\beta h_k^{-1}\delta_{n,k}^{2/3}+C_\beta k^{-1/3}.
\]
Use~\eqref{S:eq:7.30} with a sufficiently small exponential parameter to conclude that
\begin{equation}\label{S:eq:7.54}
 \mathbb E\left[\exp(c_\beta k\|\xi_{n,k}\|^{2/3})
                    \mid\mathcal F_{n,k+1}^{\eta}\right]\le C_\beta.
\end{equation}
Here and below $c_\beta$ may decrease and $C_\beta$ may increase. The bounds apply to $\xi_{n,k}$, the centered localized logarithm.

Let $C_0$ be the constant in \eqref{S:eq:localized-transfer-log-bound}. On $\{\chi_{n,k}=1\}$,
\[
 \|\log\widehat H_{n,k}^{\rm br}\|
 +\|\log\widehat H_{n,k}^{\rm circ}\|
 \le C_0(R_{n,k}+Y_{n,k})\le2C_0Y_{n,k}.
\]
There is an absolute $c_0>0$ such that $\log\cosh(y/2)\ge c_0y^2$ for $0\le y\le r_{\rm loc}$: the quotient by $y^2$ is continuous with value $1/8$ at zero and is positive on this
compact interval. Hence
\[
 \mathbb P(Y_{n,k}>y\mid\mathcal F_{n,k+1}^{\eta})
 \le\exp(-\beta c_0k y^2)
 \quad(0\le y\le r_{\rm loc}).
\]
Applying this at $y=C_\beta\sqrt{x/k}/(2C_0)$ gives, for $1\le x\le c_\beta k$,
\begin{equation}\label{S:eq:localized-transfer-moderate-tail}
 \mathbb P\left(\chi_{n,k}=1,\
 \|\log\widehat H_{n,k}^{\rm br}\|
 +\|\log\widehat H_{n,k}^{\rm circ}\|>C_\beta\sqrt{x/k}
       \,\middle|\,\mathcal F_{n,k+1}^{\eta}\right)
 \le C_\beta\me^{-c_\beta x}.
\end{equation}
Choose $c_\beta$ small enough that the indicated $y$ is at most $r_{\rm loc}$ on this range. If the threshold exceeds $2C_0r_{\rm loc}$, the localized event is empty.
\subsubsection{Simultaneous bounds}
We now return from the localized logarithms to the exact transfers. The failure probability is exponentially small in $k$, uniformly in the terminal point. Its sum over indices
starting at $\underline m_n$ is summable in $n$, so Borel--Cantelli applies to the rows constructed from the same Brownian path.

\begin{corollary}\label{S:cor:one-step-simultaneous-localization}
There is an event $\Omega_{\rm loc}$ of probability one on which, for all sufficiently large $n$ and all $\underline m_n\le k<n$, $Y_{n,k}\le r_{\rm loc}$ and
\begin{equation}\label{S:eq:simultaneous-localized-exactness}
 \zeta_{n,k}^{\rm bad}=0,\qquad
 \zeta_{n,k}^{\rm full}=b_{n,k}+\xi_{n,k}.
\end{equation}
\end{corollary}
\begin{proof}
Let $B_n=\{\max_{\underline m_n\le k<n}Y_{n,k}>r_{\rm loc}\}$, for $n$ large enough that every retained $k$ is at least $k_0(\beta)$. The tower property applied
to~\eqref{S:eq:7.55}, followed by a union bound, gives
\[
 \mathbb P(B_n\mid\eta_\infty)
 \le\sum_{k=\underline m_n}^{n-1}\me^{-c_\beta k}
 \le\frac{\me^{-c_\beta\underline m_n}}{1-\me^{-c_\beta}}.
\]
Since $\underline m_n\ge(\log n)^6$, for sufficiently large $n$ this is at most $C_\beta n^{-3}$. Integration over $\eta_\infty$ and Borel--Cantelli imply
$\mathbb P(\limsup_n B_n)=0$. Take $\Omega_{\rm loc}$ to be its complement, intersected with the events on which the defining identities hold. On this event $\chi_{n,k}=1$
throughout every sufficiently large row, so \eqref{S:eq:7.50} gives~\eqref{S:eq:simultaneous-localized-exactness}.

\end{proof}

\begin{corollary}\label{S:cor:increment-logarithm-bounds}
There is $\Omega_{\rm env}\subset\Omega_{\rm loc}$ of probability one on which, for all sufficiently large $n$,
\begin{equation}\label{S:eq:7.56}
 \sup_{\underline m_n\le k<n}
 \frac{k^{3/2}}{(\log n)^2}
 \bigl(\|\xi_{n,k}\|+\|\zeta_{n,k}^{\rm full}\|\bigr)\le C_\beta,
\end{equation}
\begin{equation}\label{S:eq:7.57}
 \sup_{\underline m_n\le k<n}\sqrt{\frac{k}{\log n}}
 \bigl(\|\log\widehat H_{n,k}^{\rm br}\|
       +\|\log\widehat H_{n,k}^{\rm circ}\|+R_{n,k}+Y_{n,k}\bigr)
 \le C_\beta.
\end{equation}
The constant is deterministic. The same event has probability one under almost every conditional law given
$\eta_\infty$.
\end{corollary}
\begin{proof}
From~\eqref{S:eq:7.54}, conditional Markov and the tower property give
\[
 \mathbb P\bigl(k^{3/2}\|\xi_{n,k}\|>u\mid\eta_\infty\bigr)
 \le C_\beta\me^{-c_\beta u^{2/3}}.
\]
At $u=(\log n)^2$, the union over a row has probability at most $C_\beta n\exp(-c_\beta(\log n)^{4/3})$. This is summable, since $c_\beta(\log n)^{4/3}\ge4\log n$ for
$\log n\ge(4/c_\beta)^3$. Hence, almost surely,
\[
 \sup_{\underline m_n\le k<n}k^{3/2}\|\xi_{n,k}\|
 \le(\log n)^2
\]
for all sufficiently large $n$. On $\Omega_{\rm loc}$, \eqref{S:eq:simultaneous-localized-exactness} and~\eqref{S:eq:7.51} imply
\[
 k^{3/2}\bigl(\|\xi_{n,k}\|+\|\zeta_{n,k}^{\rm full}\|\bigr)
 \le2(\log n)^2+C_\beta k^{-1/2}.
\]
This proves~\eqref{S:eq:7.56} with a deterministic constant.

Choose $A>0$ deterministically so that $\beta c_0A^2>4$, with $c_0$ from the radial bound preceding \eqref{S:eq:localized-transfer-moderate-tail}. Uniformly on the range under
consideration,
\[
 A\sqrt{\frac{\log n}{k}}
 \le A(\log n)^{-5/2}\longrightarrow0.
\]
Thus this threshold is at most $r_{\rm loc}$ for large $n$, and
\[
 \mathbb P\left(\exists\,\underline m_n\le k<n:
  Y_{n,k}>A\sqrt{\log n/k}\,\middle|\,\eta_\infty\right)
 \le n^{1-\beta c_0A^2}.
\]
The right side is summable. Integrating and applying Borel--Cantelli gives the radius bound throughout every sufficiently large row. There $R_{n,k}\le Y_{n,k}$ and
\eqref{S:eq:localized-transfer-log-bound} show that the parenthesis in~\eqref{S:eq:7.57} is at most $(2C_0+2)Y_{n,k}$. This proves~\eqref{S:eq:7.57}.

Define $\Omega_{\rm env}$ by intersecting $\Omega_{\rm loc}$ with the complements of these two limsup events. This is a Borel event with probability one. A sufficient starting
matrix size is the maximum of the three almost surely finite starting indices. Finally,
\[
 0=\mathbb P(\Omega_{\rm env}^c)
   =\int\mathbb P(\Omega_{\rm env}^c\mid\eta_\infty=\eta)
                          \,\nu_\ii(\dd\eta)
\]
shows that $\Omega_{\rm env}$ has probability one under almost every conditional law. \end{proof}

\subsubsection{The heat kernel on a fixed ball}
\label{S:sec:the-heat-kernel-and-short-brownian-bridges}
The formula for the radial heat kernel also yields the local estimate needed for the Brownian bridges in Section~\ref{S:sec:the-error-from-sampling-the-brownian-path}.

Let $p_h(z,w)=k_h(d(z,w))$ denote the transition density with respect to hyperbolic area. In polar coordinates the area element is $\sinh r\,\dd r\,\dd\phi$, so
\begin{equation}\label{S:eq:aux-radial-spatial-density}
 p_h(r)=2\pi\sinh r\,k_h(r).
\end{equation}
Here $p_h(r)$ is a density with respect to $\dd r$, whereas $p_h(z,w)$ is the transition density from $z$ to $w$ with respect to hyperbolic area. Fix a closed normal coordinate
ball $B$ of finite radius, and set
\[
 u_0(z,w)=\left(\frac{d(z,w)}{\sinh d(z,w)}\right)^{1/2},
 \qquad u_0(z,z)=1.
\]
Derivatives in $z,w$ are taken in this fixed normal chart. Constants may depend on the radius and are uniform in the center.

\begin{lemma}\label{S:lem:heat-kernel-estimates}
For $0<h\le1$ and $z,w\in B$,
\begin{equation}\label{S:eq:s6-heat-factorization}
 p_h(z,w)=\frac{\me^{-d(z,w)^2/(2h)}}{2\pi h}
 \bigl(u_0(z,w)+h\psi_h(z,w)\bigr),
 \qquad \sup_{0<h\le1}\|\psi_h\|_{C^1(B\times B)}<\infty.
\end{equation}
\end{lemma}
\begin{proof}
We derive the estimate from the radial density \eqref{S:eq:7.14}. Dividing by the polar area factor $2\pi\sinh r$ gives, for $r>0$,
\[
 k_h(r)=\frac{e^{-h/8}}{2\pi\sqrt{\pi h^3}}
 \int_r^\infty\frac{s e^{-s^2/(2h)}}{\sqrt{\cosh s-\cosh r}}\,\dd s.
\]
Put $x=r^2$, and define the smooth positive functions on $x,y\ge0$ by
\[
 F(x)=\frac{\sinh\sqrt x}{\sqrt x},\quad F(0)=1,\qquad
 A(x,y)=\int_0^1F(x+ty)\,\dd t.
\]
The power series for $\sinh$ shows smoothness at zero. Integration of the derivative of $\cosh\sqrt x$ gives
\[
 \cosh\sqrt{x+y}-\cosh\sqrt x=\tfrac y2 A(x,y).
\]
In the heat integral make the substitution $s=\sqrt{r^2+v^2}$, followed by $v=\sqrt h\,t$. Since $s\,\dd s=v\,\dd v$, the singular factor at the lower integration limit cancels.
The result is the exact formula
\[
 k_h(r)=\frac{e^{-r^2/(2h)}}{2\pi h}\,e^{-h/8}J_h(r^2),\qquad
 J_h(x)=\sqrt{\frac2\pi}\int_0^\infty e^{-t^2/2}A(x,ht^2)^{-1/2}\,\dd t.
\]
This expression also defines the continuous value at $r=0$.

We show $J_h(x)=F(x)^{-1/2}+O(h)$ in $C^1$ on any fixed interval $0\le x\le D^2$. The inequality for $\coth$ proved above gives
\[
 0\le\frac{F'(x)}{F(x)}
 =\frac{\coth\sqrt x-x^{-1/2}}{2\sqrt x}\le\frac16
 \quad(x>0),
\]
with value $1/6$ at zero. Thus $A\ge1$, $0\le\partial_xA\le A/6$, and
\[
 |A(x,y)^{-1/2}|\le1,\qquad
 |\partial_x A(x,y)^{-1/2}|\le\frac1{12}
 \quad(x,y\ge0).
\]
For $0\le y\le1$, smoothness and positivity on the compact rectangle $[0,D^2]\times[0,1]$ bound both $\partial_y A^{-1/2}$ and $\partial_x\partial_y A^{-1/2}$. Integrating these
derivatives from zero to $y$ gives
\[
 \sup_{0\le x\le D^2}
 \left(|A(x,y)^{-1/2}-F(x)^{-1/2}|
 +|\partial_x A(x,y)^{-1/2}-(F^{-1/2})'(x)|\right)
 \le C_Dy.
\]
For $y\ge1$, the same inequality follows from the two global bounds just proved, after increasing $C_D$. Differentiation under the Gaussian integral is justified by those bounds.
Substitution of $y=ht^2$ and $\sqrt{2/\pi}\int_0^\infty t^2e^{-t^2/2}\,\dd t=1$ therefore yields
\[
 \|J_h-F^{-1/2}\|_{C^1([0,D^2])}\le C_Dh.
\]
Since $|e^{-h/8}-1|\le h/8$ and $J_h$ has a bounded $C^1$ norm, the same estimate holds for $e^{-h/8}J_h-F^{-1/2}$.

On the fixed normal coordinate ball, the function $x(z,w)=d(z,w)^2$ is smooth, including the diagonal, and its first coordinate derivatives are bounded. Take $D$ larger than all
distances in this ball and define
\[
 \psi_h(z,w)=h^{-1}\left[e^{-h/8}J_h(x(z,w))-F(x(z,w))^{-1/2}\right].
\]
The preceding estimate and the chain rule bound its $C^1(B\times B)$ norm uniformly for $0<h\le1$. Moreover
\begin{equation}\label{S:eq:aux-leading-heat-coefficient}
 F(d(z,w)^2)^{-1/2}
 =u_0(z,w)=\left(\frac{d(z,w)}{\sinh d(z,w)}\right)^{1/2},\qquad u_0(z,z)=1.
\end{equation}
Substitution in the exact heat formula proves \eqref{S:eq:s6-heat-factorization}. Normal charts centered at different points are related by hyperbolic isometries, so the constants
depend on the fixed radius, not on the center.
\end{proof}

The coefficient $u_0$ is a smooth positive function of $d^2$. On a fixed ball its minimum is positive. After reducing $h_0$ if necessary,
\[
 u_0+h\psi_h\ge\tfrac12\min_{B\times B}u_0>0,\qquad
 \sup_{0<h\le h_0}\bigl|\nabla_z\log(u_0+h\psi_h)\bigr|<\infty.
\]
Its radial form satisfies
\[
 r u_0'(r)+\tfrac12(r\coth r-1)u_0(r)=0,
\]
as follows by differentiating $\log u_0(r)=\tfrac12(\log r-\log\sinh r)$. This is the transport equation with the normalization used in
\eqref{S:eq:strong-local-heat-factorization}.

\subsection{Products of the coupled transfer matrices}
\label{sec:combined-s-products}
\label{S:sec:products-of-the-coupled-transfer-matrices}
Fix $\beta>0$ and use the coupling of Proposition~\ref{S:prop:7.2}. Each increment error in the relative product is conjugated by the preceding Brownian product. We first write
this relative product exactly. The conditional independence proved in the preceding subsection then lets us choose a filtration in which the conjugating matrices are predictable.

The distinction between means and centered increments is already visible without conjugation. Absolute summation of the increment bounds gives
$\sum_{j\ge k}j^{-3/2}\asymp k^{-1/2}$, whereas the square root of $\sum_{j\ge k}j^{-3}$ and the drift sum $\sum_{j\ge k}j^{-2}$ are both of order $k^{-1}$. The product estimate
gives order $k^{-1}$, which is summable with the weights in the Green kernel estimate in the next subsection. Throughout, $\|\cdot\|$ denotes the Euclidean matrix operator norm.

\subsubsection{The relative product and its filtration}
\label{S:sec:s7-relative-product}
The common terminal point permits an upper triangular representation of the transfers.

The transfers fixing $\infty$ belong to the group
\[
 \mathsf B_\infty=
 \left\{\begin{pmatrix}\me^u&v\\0&\me^{-u}\end{pmatrix}:u,v\in\mathbb R\right\}.
\]
Its Lie algebra is
\[
 \mathfrak b_\infty=
 \left\{\begin{pmatrix}u&v\\0&-u\end{pmatrix}:u,v\in\mathbb R\right\}
 =\operatorname{span}\{H_0,E_{12}\}.
\]
For $D=\left(\begin{smallmatrix}u&v\\0&-u\end{smallmatrix}\right)$ we use $|D|_{\mathfrak b}=(u^2+v^2)^{1/2}$ when
taking the norm of an adjoint map. The second column and the Frobenius norm give
\[
 |D|_{\mathfrak b}\le\|D\|\le\sqrt2\,|D|_{\mathfrak b}.
\]
Put $\eta=\eta_\infty$ and use the fixed rotation $G_\eta$ of Section~\ref{S:sec:s6-terminal-conditioning}. Put $w_t=G_\eta\cdot Z_t$. Conditional on $\eta$, this is a fixed change
of coordinates sending the common terminal point to $\infty$. All matrices in this subsubsection are in these coordinates; the return to the original upper half-plane coordinates
is made in Section~\ref{S:sec:s8-coordinate-change}. Let $\mathsf F_{n,n}^{\diamond}=I$ and
\[
 \mathsf F_{n,j}^{\diamond}
 =\widehat H_{n,j}^{\diamond}\mathsf F_{n,j+1}^{\diamond},
 \qquad\diamond\in\{\mathrm{br},\mathrm{circ}\}.
\]
For the Brownian path, $\mathsf F_{n,j}^{\rm br}=\mathcal A(w_{\vartheta_{n,j}})$. The matrix used on the operator interval $I_{n,k}$ is
$X_{n,k}^{\diamond}=\mathsf F_{n,k+1}^{\diamond}$. In particular the two matrices on $I_{n,n-1}$ are both $I$. The intervening Brownian product is
\[
 P_{n;k,j}=\widehat H_{n,k}^{\rm br}\cdots\widehat H_{n,j-1}^{\rm br}
 =\mathsf F_{n,k}^{\rm br}(\mathsf F_{n,j}^{\rm br})^{-1},\qquad
 P_{n;k,k}=I.
\]
For $\underline m_n\le k\le m<n$ put
\[
 \mathcal R_{n;k,m}
 =\left(\widehat H_{n,k}^{\rm circ}\cdots\widehat H_{n,m}^{\rm circ}\right)
  \left(\widehat H_{n,k}^{\rm br}\cdots\widehat H_{n,m}^{\rm br}\right)^{-1},
 \qquad\mathcal R_{n;k,k-1}=I.
\]
These are products of the original transfers, before any localization. By the definition of the real triangular logarithm,
$\widehat H_{n,j}^{\rm circ}= \exp(\zeta_{n,j}^{\rm full})\widehat H_{n,j}^{\rm br}$. It follows that
\begin{equation}\label{S:eq:s7-exact-relative-product}
 \mathcal R_{n;k,m} =\prod_{j=k}^m \exp\{\operatorname{Ad}_{P_{n;k,j}}\zeta_{n,j}^{\rm full}\},\qquad \mathcal R_{n;k,n-1} =\mathsf F_{n,k}^{\rm circ}(\mathsf F_{n,k}^{\rm br})^{-1}.
\end{equation}
Here too the factors are ordered increasingly from left to right. To check their order, let $C_{m-1}$ and $B_{m-1}=P_{n;k,m}$ be the circular and Brownian products through $m-1$.
Then
\[
 \mathcal R_{n;k,m} =C_{m-1}\exp(\zeta_{n,m}^{\rm full})\widehat H_{n,m}^{\rm br} (B_{m-1}\widehat H_{n,m}^{\rm br})^{-1} =\mathcal R_{n;k,m-1} B_{m-1}\exp(\zeta_{n,m}^{\rm full})B_{m-1}^{-1}.
\]
Conjugation commutes with the exponential series. Starting with $B_{k-1}=I$ proves the first identity by induction; the second is the definition of the full recursion product.

The order in this identity matters probabilistically. The conjugating product $P_{n;k,j}$ uses indices smaller than $j$, whereas chronological Brownian time runs in the other
direction. Conditional independence within a row permits us to reveal the variables in the order used by the product. For fixed $n,k$ consider the increasing filtration
\begin{equation}\label{S:eq:s7-filtration}
 \mathcal G_{n;k,j}
 =\sigma\bigl(\eta,(R_{n,i},U_{n,i}):k\le i\le j\bigr),\qquad
 \mathcal G_{n;k,k-1}=\sigma(\eta).
\end{equation}
Lemma~\ref{S:lem:7.2} and the construction in Section~\ref{S:sec:coupling-the-finite-transfer-increments-to-brownian-increments} give a product law for $(R_{n,i},U_{n,i})$ within
this row, conditional on $\eta$. The localized variable $\zeta_{n,j}^{\rm loc}$ depends only on $(R_{n,j},U_{n,j})$, and its mean $b_{n,j}$ is deterministic. Hence conditional
independence gives, for every bounded Borel function $f$,
\[
 \mathbb E[f(\zeta_{n,j}^{\rm loc})\mid\mathcal G_{n;k,j-1}]
 =\mathbb E[f(\zeta_{n,j}^{\rm loc})\mid\eta].
\]
One verifies this identity first by testing against products of bounded functions of $\eta$ and the preceding variables, and then uses the monotone class theorem. Truncation
extends it to each integrable function used below. In particular, with $\xi_{n,j}=\zeta_{n,j}^{\rm loc}-b_{n,j}$, Theorem~\ref{S:thm:7.4} yields
\begin{equation}\label{S:eq:s7-reversed-filtration-moments}
 \mathbb E[\xi_{n,j}\mid\mathcal G_{n;k,j-1}]=0,\qquad \mathbb E[\|\xi_{n,j}\|^p\mid\mathcal G_{n;k,j-1}]^{1/p} \le C_\beta p^{3/2}j^{-3/2},\qquad \|b_{n,j}\|\le C_\beta j^{-2}.
\end{equation}
Furthermore $P_{n;k,j}$ uses only the variables with indices $k,\ldots,j-1$. Thus $\operatorname{Ad}_{P_{n;k,j}}$ is $\mathcal G_{n;k,j-1}$-measurable. These conditional identities
hold before imposing bounds on the increments or their products. Such an event may depend on variables not yet revealed and cannot be included in the conditioning without
justification. For example, let $B_j$ be predictable linear maps with a deterministic bound $\|B_j\|\le L$, and use the Frobenius inner product on the Lie algebra. For $i<j$,
conditional expectation at $\mathcal G_{n;k,j-1}$ gives
\[
 \mathbb E\!\left[\langle B_i\xi_{n,i},B_j\xi_{n,j}\rangle
                   \mid\eta\right]=0,\qquad
 \mathbb E\!\left[\left\|\sum_{j=k}^m B_j\xi_{n,j}
                      \right\|_{\mathrm{HS}}^2\mid\eta\right]
 \le C_\beta L^2\sum_{j=k}^m j^{-3}\le C_\beta L^2k^{-2}.
\]
The norm comparison in dimension two and the second moment in \eqref{S:eq:s7-reversed-filtration-moments} justify integrability. The cancellation uses distinct, successively
revealed increments; the two coordinates of one increment are still treated jointly. The predictable mean sum has the separate bound $\|\sum_{j=k}^m B_jb_{n,j}\|\le C_\beta L/k$.
The stopped adjoints used below satisfy the required deterministic bound before either estimate is applied.

\begin{figure}[H]
\centering
\begin{tikzpicture}[x=1cm,y=1cm,every node/.style={font=\footnotesize,inner sep=4pt},>=Stealth]
 \node at (3.1,3.85) {Indices $k,\ldots,j-1$};
 \node at (10.6,3.85) {Index $j$};
 \node[draw,minimum width=6.0cm,minimum height=.85cm] (past) at (3.1,2.95)
   {$\eta,\ (R_{n,i},U_{n,i})_{k\le i<j}$};
 \node[draw,minimum width=5.1cm,minimum height=.85cm] (now) at (10.6,2.95)
   {$(R_{n,j},U_{n,j})$};
 \node[draw,minimum width=6.0cm,minimum height=1.05cm,align=center] (pp) at (3.1,1.25)
   {$P_{n;k,j}=\widehat H_{n,k}^{\rm br}\cdots\widehat H_{n,j-1}^{\rm br}$\\$\mathcal G_{n;k,j-1}$-measurable};
 \node[draw,minimum width=5.1cm,minimum height=1.05cm,align=center] (xi) at (10.6,1.25)
   {$\xi_{n,j}=\zeta_{n,j}^{\rm loc}-b_{n,j}$\\$\mathbb E[\xi_{n,j}\mid\mathcal G_{n;k,j-1}]=0$};
 \draw[->] (past)--(pp); \draw[->] (now)--(xi);
 \draw[densely dashed] (6.85,3.45)--(6.85,.55);
 \draw[->] (1,-.05)--(12.8,-.05) node[right] {$j$ increases};
 \node[align=center] at (6.8,-1.0)
  {$\mathcal R_{n;k,j}=\mathcal R_{n;k,j-1}\,\bigl(P_{n;k,j}\exp(\zeta_{n,j}^{\rm full})P_{n;k,j}^{-1}\bigr)$};
 \node at (6.8,-1.8)
  {$\displaystyle\mathcal R_{n;k,m}=\prod_{j=k}^{m}\exp\{\operatorname{Ad}_{P_{n;k,j}}\zeta_{n,j}^{\rm full}\}$};
\end{tikzpicture}
\caption{The order used in \eqref{S:eq:s7-exact-relative-product} and \eqref{S:eq:s7-filtration}, with $n,k$ fixed. The variables to the left generate $\mathcal G_{n;k,j-1}$ and determine $P_{n;k,j}$. Conditional on $\eta$, the pairs are independent as their index varies. This gives the displayed conditional mean of $\xi_{n,j}$. The recurrence below the boxes adds the next conjugated exponential on the right; the product is ordered by increasing $j$ from left to right. The centered estimates are subsequently applied to bounded predictable maps, as in \eqref{S:eq:s7-reversed-filtration-moments} and the calculation following it.}
\label{S:fig:relative-product-conditioning}
\end{figure}
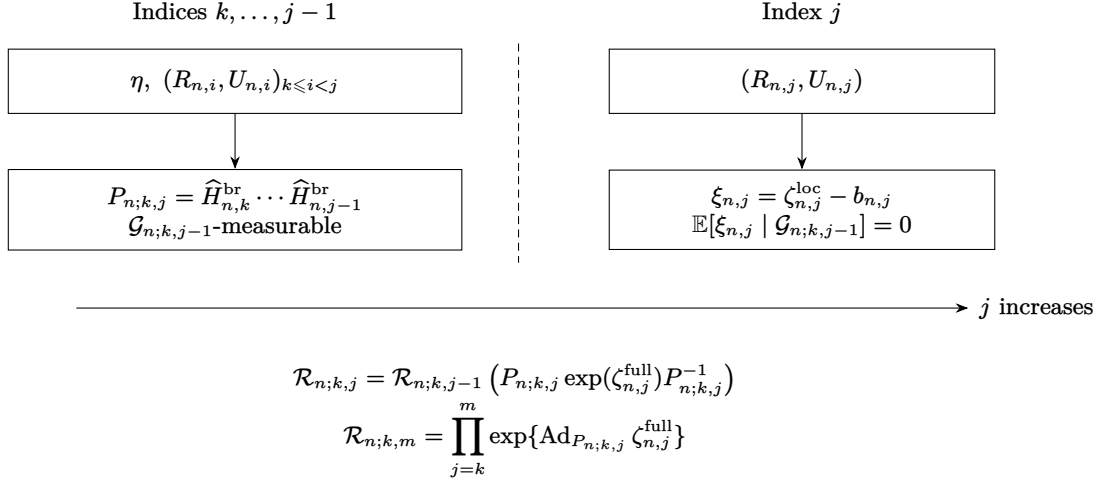

Figure~\ref{S:fig:relative-product-conditioning} records the order in which the transfers are multiplied and the variables are revealed. The product preceding index $j$ is known
under $\mathcal G_{n;k,j-1}$, which is the conditioning used to center the next increment.

The following estimate is proved after the Brownian adjoint bound and the upper triangular product lemma.

\begin{theorem}\label{S:thm:s7-product-comparison}
Almost surely, for all sufficiently large $n$,
\begin{equation}\label{S:eq:s7-all-prefix-rate}
 \sup_{\underline m_n\le k<n}
 \frac{k}{(\log n)^{11/2}}
 \max_{k\le m<n}\|\mathcal R_{n;k,m}-I\|
 \le C_{\beta,\omega}.
\end{equation}
\end{theorem}

\subsubsection{Brownian motion conditioned on its limit}
\label{S:sec:brownian-motion-under-terminal-conditioning}

Conditioning on the terminal point changes the drift of the driving processes. We compute the Brownian motions under the conditional measure and the normalization of the terminal
column used in the conditional estimates of Sections~\ref{S:sec:the-error-from-sampling-the-brownian-path} and \ref{S:sec:the-limiting-coefficient-of-the-error}. In this
subsubsection $Z_t=x_t+\ii y_t$ denotes the Brownian motion in its original upper half-plane coordinates. Let $\mathcal F_t^{\rm aff}$ be its natural filtration, or equivalently
that of its two driving Brownian motions $W^1,W^2$. Thus
\[
 y_t=\exp(W_t^2-t/2),\qquad
 x_t=\int_0^t y_s\,\dd W_s^1,\qquad
 A_t=\int_0^t y_s^2\,\dd s.
\]
The superscripts $1,2$ index the Brownian components. For $\eta\in\partial\mathbb H$ let $\mathbb P^\eta$ be the terminal conditional law from Lemma~\ref{S:lem:7.2}. On each finite
horizon its density with respect to $\mathbb P$ is $P_\eta(Z_T)$. We use the notation
\[
 b_\eta=(y\partial_x\log P_\eta,y\partial_y\log P_\eta),\qquad
 \widetilde W_t=W_t-\int_0^t b_\eta(Z_s)\,\dd s.
\]
Conditional laws and the displayed stochastic processes are first specified on the Borel path space. Each filtration is completed under the probability measure with respect to
which the martingale property is proved.

\begin{proposition}\label{S:prop:aux-stochastic-normalization}
Under $\mathbb P$, $A_\infty<\infty$ almost surely and, on an extension, $x_t=\mathcal B_{A_t}$ for a standard Brownian motion $\mathcal B$. Under $\mathbb P^\eta$, $\widetilde W$
is a two-dimensional Brownian motion on every finite horizon, and $|b_\eta|=1$.
\end{proposition}

\begin{proof}
The bound for the imaginary part in Lemma~\ref{S:lem:5.0} gives $A_\infty=\int_0^\infty y_s^2\,\dd s<\infty$. The continuous local martingale $x=\int y\,\dd W^1$ has bracket $A$;
the Dambis--Dubins--Schwarz theorem allowing finite quadratic variation \cite{kallenberg2021foundations} (Theorem~19.4) gives the claimed Brownian representation on an extension.
For the conditional version used below, take $C(u)=\inf\{t:A_t>u\}$ and an independent Brownian motion $B'$:
\begin{equation}\label{S:eq:s7-brownian-continuation}
 \mathcal B_u=
 \begin{cases}
 x_{C(u)},&0\le u<A_\infty,\\
 q_\infty+B'_{u-A_\infty},&u\ge A_\infty.
 \end{cases}
\end{equation}
Given $\mathcal V=\sigma(W_s^2:s\ge0)$, the integrand $y$ is deterministic, and $x$ is centered Gaussian with covariance $A_{s\wedge t}$. The continuation above is therefore
centered Gaussian with covariance $u\wedge v$ and continuous paths. Its conditional law is Wiener measure for almost every vertical path, so $\mathcal B$ is independent of
$\mathcal V$ and $q_\infty=\mathcal B_{A_\infty}$. For finite $\eta$, write $d=x-\eta$ and $D=d^2+y^2$. Differentiating the Poisson kernel gives
\begin{equation}\label{S:eq:s7-poisson-gradient}
 P_\eta(x+\ii y)=\frac{(1+\eta^2)y}{D},\qquad
 b_\eta=\left(-\frac{2yd}{D},\frac{d^2-y^2}{D}\right),\qquad
 |b_\eta|^2=\frac{4y^2d^2+(d^2-y^2)^2}{D^2}=1.
\end{equation}
For $\eta=\infty$, $P_\infty=y$ and $b_\infty=(0,1)$. Since $LP_\eta=0$ and $P_\eta(\ii)=1$, It\^o's formula yields $P_\eta(Z_t)=\exp\{\int_0^t b_\eta(Z_s)\cdot\dd W_s-t/2\}$. On
each finite horizon $T$, Novikov's condition is immediate from $\mathbb E\exp\{\tfrac12\int_0^T|b_\eta|^2\dd s\}=e^{T/2}$. Girsanov's theorem therefore identifies
$\widetilde W=W-\int b_\eta(Z_s)\dd s$ as a two-dimensional Brownian motion under $\mathbb P^\eta$; see Theorem~19.20, Theorem~19.24 and Corollary~19.26 of
\cite{kallenberg2021foundations}. Taking integer $T$ fixes one version on all finite horizons.
\end{proof}

The terminal column has the exact normalization
\begin{equation}\label{S:eq:s7-terminal-column-poisson}
 q_t^\eta:=|\mathcal A(x_t+\ii y_t)(\eta,1)^{\mathsf T}|^2
 =\frac{(x_t-\eta)^2+y_t^2}{y_t}
 =(1+\eta^2)P_\eta(Z_t)^{-1}.
\end{equation}
Under $\mathbb P^\eta$ the density process satisfies $\dd L=L b_\eta\cdot\dd\widetilde W+L\,\dd t$. Another application of It\^o's formula gives
\begin{equation}\label{S:eq:s7-terminal-column-martingale}
 \dd q_t^\eta=-q_t^\eta b_\eta(Z_t)\cdot\dd\widetilde W_t,\qquad
 \dd\log q_t^\eta=-b_\eta(Z_t)\cdot\dd\widetilde W_t-\frac12\,\dd t.
\end{equation}
Since $|b_\eta|=1$, L\'evy's characterization, Theorem~19.3 of \cite{kallenberg2021foundations}, shows that $-\int_0^t b_\eta(Z_s)\cdot\dd\widetilde W_s$ is Brownian motion in the
conditioned filtration. Solving the logarithmic equation gives
\[
 q_t^\eta=(1+\eta^2)
 \exp\!\left(-\int_0^t b_\eta(Z_s)\cdot\dd\widetilde W_s-t/2\right).
\]
Independent Brownian increments and the Gaussian exponential moment therefore imply, for $0\le s\le t<\infty$,
\[
 \mathbb E^\eta[q_t^\eta\mid\mathcal F_s^Z]=q_s^\eta,
 \qquad
 \mathbb E^\eta[(q_t^\eta)^2]=(1+\eta^2)^2\me^t.
\]
Thus $q^\eta$ is a square-integrable martingale on every finite interval.

The expectation here is conditional on the terminal value. Its integrable normalization under the joint law is $q_t^{\eta_\infty}/(1+\eta_\infty^2)$: conditioning on $\eta_\infty$
gives expectation one. By contrast, Tonelli and the Cauchy law of $\eta_\infty$ give $\mathbb E q_t^{\eta_\infty}=\int(1+\eta^2)\nu_\ii(\dd\eta) =\infty$. The stochastic identities
for the terminal column remain valid when its joint first moment is infinite and are used in~\eqref{S:eq:terminal-q-sde}. The Brownian motion $\widetilde W$ under $\mathbb P^\eta$
is defined in Lemma~\ref{S:lem:terminal-affine-innovations}; a bridge conditioned on two finite positions has a different conditional law.

\subsubsection{A fixed change of coordinates}
\label{S:sec:s7-terminal-coordinates}
Sending the terminal point to infinity makes every transfer upper triangular. The boundary columns are then simple functions of the height and one real coordinate. We calculate
these functions explicitly, keeping the reciprocal normalization of the two boundary vectors.

Write $Z_t=x_t^{\rm aff}+\ii y_t^{\rm aff}$ for the original path from now on. Let $\eta=\eta_\infty$, which is finite almost surely, and recall the rotation and the two rescaled
boundary vectors
\[
 \kappa_\eta=\sqrt{1+\eta^2},\qquad
 G_\eta=\kappa_\eta^{-1}
 \begin{pmatrix}\eta&1\\-1&\eta\end{pmatrix},\qquad
 u_0^{\rm ad}=(-\eta,1)^{\mathsf T},\quad u_1^{\rm ad}=e_1.
\]
For the original vectors $u_0=(-1,0)^{\mathsf T}$ and $u_1=(\eta,1)^{\mathsf T}$, matrix multiplication gives
\[
 G_\eta u_0=\kappa_\eta^{-1}u_0^{\rm ad},\qquad
 G_\eta u_1=\kappa_\eta u_1^{\rm ad},\qquad
 [u_0^{\rm ad},u_1^{\rm ad}]=1.
\]
The M\"obius map represented by $G_\eta$ fixes $\ii$, sends $\eta$ to $\infty$, and sends $\infty$ to $q_0=-\eta$. Put
\[
 w_t=G_\eta\cdot Z_t=x_t+\ii y_t,\qquad
 z_t=\frac{q_0-x_t}{y_t},\qquad d_t=y_t^{-1/2}.
\]
Here $z_t$ is a real-valued coordinate. The stochastic equations and distributions of $w$ are first considered under $\mathbb P^\eta$, where this change of coordinates is fixed.

Multiplication by $\mathcal A(w_t)$ gives
\begin{equation}\label{S:eq:adapted-boundary-columns}
 a_t^{\rm ad}=\mathcal A(w_t)u_0^{\rm ad}
 =y_t^{1/2}(z_t,1)^{\mathsf T},\qquad
 c_t^{\rm ad}=\mathcal A(w_t)u_1^{\rm ad}=y_t^{-1/2}e_1.
\end{equation}
In particular,
\begin{equation}\label{S:eq:adapted-boundary-norms}
 |a_t^{\rm ad}|^2=y_t(1+z_t^2),\qquad
 |c_t^{\rm ad}|^2=y_t^{-1},\qquad
 s_t:=\frac{|c_t^{\rm ad}|}{|a_t^{\rm ad}|}
 =\frac1{y_t\sqrt{1+z_t^2}}.
\end{equation}
For a point $\zeta\in\mathbb H$ define
\[
 O_\eta(\zeta)=\mathcal A(G_\eta\cdot\zeta)
               G_\eta\mathcal A(\zeta)^{-1}.
\]
The three factors send $\ii$ to $\zeta$, then to $G_\eta\cdot\zeta$, and finally to $\ii$. If a real matrix with
determinant one $\left(\begin{smallmatrix}a&b\\c&d\end{smallmatrix}\right)$ fixes $\ii$, then $a=d$, $b=-c$, and
$a^2+b^2=1$. It follows that $O_\eta(\zeta)\in\mathrm{SO}_2$. Thus the columns in the two coordinate systems satisfy
\begin{equation}\label{S:eq:adapted-reciprocal-normalization}
 a^{\rm ad}=\kappa_\eta O_\eta a^{\rm aff},\qquad
 c^{\rm ad}=\kappa_\eta^{-1}O_\eta c^{\rm aff}.
\end{equation}
The reciprocal rescaling of the vectors leaves each product $a(r)c(s)^{\mathsf T}$ unchanged before the orthogonal conjugation. For either an integral or its corresponding sampled
sum, the corresponding scalar quantities therefore obey
\[
 p^{\rm ad}=\kappa_\eta^2p^{\rm aff}, q^{\rm ad}=\kappa_\eta^{-2}q^{\rm aff},\qquad A^{\rm ad}=\kappa_\eta^{-2}A^{\rm aff}, C^{\rm ad}=\kappa_\eta^2C^{\rm aff}.
\]
The factors cancel in $Ap$, $Cq$, and $pq$. These identities compare the column lengths. For the inverse kernels, the orthogonal factors at the two positions are estimated by
Lemma~\ref{S:lem:direct-qr} in Section~\ref{S:sec:comparison-of-the-finite-and-sampled-brownian-operators}.

\begin{lemma}\label{S:lem:s7-conditioned-law}
There is a two-dimensional Brownian motion $B$, independent of $\eta$ under the joint law, such that
\begin{equation}\label{S:eq:s7-conditioned-sde}
 \dd x_t=y_t\,\dd B_t^1,\qquad
 \dd y_t=y_t\,\dd B_t^2+y_t\,\dd t,\qquad w_0=\ii.
\end{equation}
Under the joint law, $z$ is stationary with density $\pi^{-1}(1+z^2)^{-1}$.
\end{lemma}

\begin{proof}
Under $\mathbb P^\eta$, the rotation sends the terminal harmonic function to $P_\infty(w)=\operatorname{Im}w$. Its transformed generator is
\[
 y^{-1}L(yf)=\tfrac12y^2(f_{xx}+f_{yy})+yf_y.
\]
Write $L^\infty$ for this generator. Besides $L^\infty x=0$ and $L^\infty y=y$, direct differentiation gives
\[
 L^\infty(x^2)-2xL^\infty x =L^\infty(y^2)-2yL^\infty y=y^2,\qquad L^\infty(xy)-xL^\infty y-yL^\infty x=0.
\]
Thus the coordinate martingale parts have quadratic covariations $y_t^2\delta_{ij}\,\dd t$, and
\[
 B_t^1=\int_0^t y_u^{-1}\,\dd x_u,\qquad
 B_t^2=\int_0^t y_u^{-1}\,\dd y_u-t
\]
are continuous local martingales with identity quadratic covariance under that law. L\'evy's characterization of Brownian motion, Theorem~19.3 of \cite{kallenberg2021foundations},
identifies the pair as a two-dimensional Brownian motion. For measurability, stop at $\tau_R=\inf\{t:|x_t|\ge R\ \text{or}\ y_t\notin(R^{-1},R)\}\wedge R$. For a dyadic partition
$(t_i)$, the left endpoint sums
\[
 \sum_i y_{t_i\wedge\tau_R}^{-1}
  (x_{t_{i+1}\wedge t\wedge\tau_R}-x_{t_i\wedge t\wedge\tau_R}),
 \qquad
 \sum_i y_{t_i\wedge\tau_R}^{-1}
  (y_{t_{i+1}\wedge t\wedge\tau_R}-y_{t_i\wedge t\wedge\tau_R})
\]
use only path evaluations, and hence are Borel. On $[0,\tau_R]$ the integrands and the diffusion and drift coefficients are bounded. Continuity of $y^{-1}$, the stochastic
isometry, and $|\int_0^t f_s\,\dd s|^2\le R\int_0^R|f_s|^2\,\dd s$ show convergence of these sums in probability, uniformly in $t$ on bounded intervals. Under the joint law, choose
a deterministic subsequence whose successive uniform errors exceed $2^{-l}$ with probability at most $2^{-l}$, simultaneously for $R\le l$. Borel--Cantelli makes their limits Borel
outside a null set; set the limits to zero on that set. As $R\uparrow\infty$ the stopped limits agree on overlaps and give versions of the two displayed integrals on the original
joint path space. The equation has the explicit global solution
\begin{equation}\label{S:eq:s7-conditioned-explicit-solution}
 y_t=\exp(B_t^2+t/2),\qquad x_t=\int_0^t y_u\,\dd B_u^1.
\end{equation}
It is pathwise unique and has finite coordinates on every bounded time interval. The conditional law of $B$ is Wiener measure for almost every $\eta$. For bounded Borel functions
$f$ of the Brownian path and $g$ of the terminal point,
\[
 \mathbb E[f(B)g(\eta)]
 =\int g(v)\,\mathbb E^v[f(B)]\,\nu_\ii(\dd v)
 =\left(\int f\,\dd\mathsf W_2\right)\mathbb E g(\eta),
\]
where $\mathsf W_2$ is two-dimensional Wiener measure. This proves independence and the Wiener distribution under the joint law.

It\^o's formula applied to~\eqref{S:eq:s7-conditioned-sde} gives
\[
 \dd\log y_t=\dd B_t^2+\tfrac12\,\dd t,\qquad
 \dd(y_t^{-1})=-y_t^{-1}\,\dd B_t^2.
\]
The cross variation of $x$ and $y^{-1}$ is zero. The product rule for $(q_0-x)y^{-1}$ consequently yields
\begin{equation}\label{S:eq:s7-conditioned-coordinates}
 \dd\log d_t=-\tfrac12\,\dd B_t^2-\tfrac14\,\dd t,\qquad
 \dd z_t=-\dd B_t^1-z_t\,\dd B_t^2.
\end{equation}
The coefficients in the equation for $z$ are globally Lipschitz. In particular it has a time-homogeneous transition semigroup, denoted here by $Q_t$.

We compute its one-dimensional marginals from the conditional distribution of the terminal point given the original path up to time $t$. For $Z_t=X+\ii Y$ and $D=(\eta-X)^2+Y^2$,
the M\"obius formula gives
\[
 x_t=\frac{(\eta X+1)(\eta-X)-\eta Y^2}{D},\qquad
 y_t=\frac{(1+\eta^2)Y}{D}.
\]
Subtracting the first numerator from $-\eta D$ gives $-(1+\eta^2)(\eta-X)$. Therefore
\begin{equation}\label{S:eq:s7-coordinate-identification}
 z_t=-\frac{\eta-x_t^{\rm aff}}{y_t^{\rm aff}}.
\end{equation}
By the Markov property and the harmonic measure of the terminal point,
\[
 \mathbb P(\eta\in\dd v\mid\mathcal F_t^{\rm aff})
 =\frac{Y}{\pi((v-X)^2+Y^2)}\,\dd v.
\]
This is the harmonic measure computed in Section~\ref{S:sec:s6-terminal-conditioning}. Changing variables $u=-(v-X)/Y$ proves, for every bounded Borel $f$,
\[
 \mathbb E[f(z_t)\mid\mathcal F_t^{\rm aff}]
 =\int_{\mathbb R}\frac{f(u)}{\pi(1+u^2)}\,\dd u.
\]
Write $\pi_C(\dd u)=\dd u/(\pi(1+u^2))$. At time zero, $z_0=-\eta$ has law $\pi_C$ and is independent of $B$. The equation for $z$ and the equality of all its marginals imply
$\pi_C Q_t=\pi_C$. For $0\le t_1<\cdots<t_l$, $h\ge0$, and bounded Borel $f_1,\ldots,f_l$, the Markov property gives
\[
 \mathbb E\prod_{i=1}^l f_i(z_{h+t_i})
 =\int\pi_C(\dd u_1)f_1(u_1)
       \prod_{i=2}^l Q_{t_i-t_{i-1}}(u_{i-1},\dd u_i)f_i(u_i).
\]
The right side is independent of $h$, proving stationarity under the joint law. Under $\mathbb P^\eta$, the initial
value is the fixed number $-\eta$. \end{proof}

\subsubsection{Uniform bounds for transfer matrices along the Brownian path}

We estimate the adjoint maps appearing in the relative product. In the coordinates with terminal point at infinity, these adjoint maps can be expressed by a height ratio and the
real coordinate $z_t$.

For $0\le s\le t$ write
\[
 P_{s,t}=\mathcal A(w_t)\mathcal A(w_s)^{-1}
 =\begin{pmatrix}r_{s,t}&v_{s,t}\\0&r_{s,t}^{-1}\end{pmatrix}.
\]
The inverse matrix is $\mathcal A(w_s)^{-1}=
\left(\begin{smallmatrix}\sqrt{y_s}&x_s/\sqrt{y_s}\\0&1/\sqrt{y_s}\end{smallmatrix}\right)$. Multiplication therefore
gives
\begin{equation}\label{S:eq:s7-relative-transfer-coordinates}
 r_{s,t}^2=\frac{y_s}{y_t} =\exp\left\{-(B_t^2-B_s^2)-\frac{t-s}{2}\right\},\qquad r_{s,t}v_{s,t}=\frac{x_s-x_t}{y_t} =z_t-r_{s,t}^2z_s.
\end{equation}
Since $P_{s,t}^{-1}=\left(\begin{smallmatrix}r_{s,t}^{-1}&-v_{s,t}\\0&r_{s,t}\end{smallmatrix}\right)$, a further
multiplication yields
\begin{equation}\label{S:eq:s7-adjoint-formula}
 \operatorname{Ad}_{P_{s,t}}
 \begin{pmatrix}u&v\\0&-u\end{pmatrix}
 =\begin{pmatrix}
 u&r_{s,t}^2v-2u(z_t-r_{s,t}^2z_s)\\0&-u
 \end{pmatrix}.
\end{equation}
In particular the conjugation preserves $\mathfrak b_\infty$.

\begin{lemma}\label{S:lem:s7-height-ratio}
There is an absolute $C$ such that, for $T,R\ge2$,
\begin{equation}\label{S:eq:s7-height-tail}
 \mathbb P\!\left(\sup_{0\le s\le t\le T}r_{s,t}^2>R\right)
 \le\frac{C(T+1)}R.
\end{equation}
The same estimate holds conditionally on $\eta$.
\end{lemma}

\begin{proof}
Put $A_t=-B_t^2-t/2$ in this proof only, so $r_{s,t}^2=\exp(A_t-A_s)$. For $c\in\mathbb R$ let $\tau_c=\inf\{t\ge0:A_t=c\}$. The process $\me^{A_t}$ is a martingale and
$A_t\to-\infty$ almost surely. If $a,b>0$, the time $\tau_a\wedge\tau_{-b}$ is finite and the stopped martingale takes values in $[\me^{-b},\me^a]$. Optional sampling, first at a
bounded time and then by dominated convergence, gives
\[
 1=\me^a\mathbb P(\tau_a<\tau_{-b})
   +\me^{-b}\mathbb P(\tau_{-b}<\tau_a),\qquad
 \mathbb P(\tau_a<\tau_{-b})
 =\frac{1-\me^{-b}}{\me^a-\me^{-b}}.
\]
The events on the left increase to $\{\tau_a<\infty\}$ as $b\to\infty$. Thus
\begin{equation}\label{S:eq:s7-drifted-maximum}
 \mathbb P\!\left(\sup_{t\ge0}A_t>a\right)=\me^{-a},\qquad a>0.
\end{equation}
The hitting calculation gives $\mathbb P(\sup_tA_t\ge a)=\me^{-a}$. Also $\mathbb P(\sup_tA_t>a)=\lim_{m\to\infty}\me^{-a-1/m}=\me^{-a}$, so the maximum has no atom at $a>0$.

For integers $m\ge0$ define
\[
 \begin{aligned}
 V_m&=\sup_{m\le s\le m+1}(A_{m+1}-A_s),\\
 W_m&=\sup_{m\le s\le t\le m+1}(A_t-A_s),\\
 M_m&=\sup_{t\ge m+1}(A_t-A_{m+1}).
 \end{aligned}
\]
If $K_m=\sup_{0\le u\le1}|B_{m+u}^2-B_m^2|$, then $V_m,W_m\le2K_m+1/2$. The reflection principle and the exponential bound for a normal random variable give
$\mathbb P(K_m>x)\le4\me^{-x^2/2}$ for $x\ge0$. For every fixed $a>0$, integration of the tail yields
\begin{equation}\label{S:eq:s7-unit-brownian-exponential}
 \mathbb E\me^{aK_m}
 =1+a\int_0^\infty\me^{ax}\mathbb P(K_m>x)\,\dd x
 \le1+4a\int_0^\infty\me^{ax-x^2/2}\,\dd x<\infty.
\end{equation}
The bound is independent of $m$. In particular $\mathbb E\me^{V_m}+\mathbb E\me^{W_m}\le C$. The variable $M_m$ is independent of the increments on $[m,m+1]$ and has the
tail~\eqref{S:eq:s7-drifted-maximum}. Hence, for $a>0$,
\[
 \mathbb P(V_m+M_m>a)
 =\mathbb E\min\{1,\me^{V_m-a}\}\le C\me^{-a},\qquad
 \mathbb P(W_m>a)\le C\me^{-a}.
\]
If $s\in[m,m+1]$ and $t\ge s$, then $A_t-A_s\le W_m$ when $t\le m+1$, and $A_t-A_s\le V_m+M_m$ otherwise. Every pair $0\le s\le t\le T$ belongs to one of these cases for some
$0\le m\le\lfloor T\rfloor$. A union bound with $a=\log R$ proves~\eqref{S:eq:s7-height-tail}. Since $B$ is independent of $\eta$, this proof also gives the conditional estimate.
\end{proof}

\begin{lemma}\label{S:lem:s7-cauchy-maximum}
There is an absolute $C$ such that, under the joint law, for every integer $m\ge0$ and $R\ge2$,
\begin{equation}\label{S:eq:s7-cauchy-tail}
 \mathbb P\!\left(\sup_{0\le u\le1}|z_{m+u}|>R\right)\le\frac CR.
\end{equation}
\end{lemma}

\begin{proof}
Let $\Delta B_u^i=B_{m+u}^i-B_m^i$ for $0\le u\le1$, and set
\[
 V_u=\exp(\Delta B_u^2+u/2),\qquad
 J_u=\int_0^u V_v\,\dd\Delta B_v^1.
\]
By~\eqref{S:eq:s7-conditioned-explicit-solution}, $y_{m+u}=y_mV_u$ and $x_{m+u}=x_m+y_mJ_u$. Consequently the exact solution over this interval is
\begin{equation}\label{S:eq:s7-cauchy-interval-solution}
 z_{m+u}=V_u^{-1}(z_m-J_u).
\end{equation}
Put $D=\sup_{u\le1}V_u^{-1}$ and $J^*=\sup_{u\le1}|J_u|$. The future increments $\Delta B$ are independent of $\mathcal H_m=\sigma(\eta,B_v:0\le v\le m)$, and $z_m$ is
$\mathcal H_m$-measurable. Thus the pair $(D,J^*)$ is independent of $z_m$. Conditional on the path $\Delta B^2$, the integrand $V$ is deterministic and $J$ is a square-integrable
martingale in $\Delta B^1$ with terminal second moment $\int_0^1 V_v^2\,\dd v$. The $L^2$ maximal inequality gives
\[
 \mathbb E[J^*\mid\Delta B^2]
 \le2\left(\int_0^1 V_v^2\,\dd v\right)^{1/2}.
\]
This is Doob's $L^2$ maximal inequality \cite{kallenberg2021foundations} (Theorem~9.17), applied under the fixed conditional law. Rational times and continuity give the supremum
over $[0,1]$.

For $K=\sup_{u\le1}|\Delta B_u^2|$, one has $D\le\me^K$ and $(\int_0^1V_v^2\,\dd v)^{1/2}\le\me^{K+1/2}$. The exponential estimate~\eqref{S:eq:s7-unit-brownian-exponential}
therefore implies
\[
 \mathbb ED\le\mathbb E\me^K\le C,\qquad
 \mathbb E(DJ^*)\le2\me^{1/2}\mathbb E\me^{2K}\le C.
\]
Stationarity from Lemma~\ref{S:lem:s7-conditioned-law} gives $\mathbb P(|z_m|>x)\le2/(\pi x)$ for every $x>0$. By~\eqref{S:eq:s7-cauchy-interval-solution} and the independence just
proved,
\[
 \mathbb P\!\left(\sup_{u\le1}|z_{m+u}|>R\right) \le\mathbb P(D|z_m|>R/2)+\mathbb P(DJ^*>R/2) \le\frac4{\pi R}\mathbb ED+\frac2R\mathbb E(DJ^*) \le\frac CR.
\]
The integrations use $\mathbb ED$ and $\mathbb E(DJ^*)$, both finite by the displayed estimates. The Cauchy law for $z_m$ is under the joint law.
\end{proof}

For the products $P_{n;k,j}$ introduced in Section~\ref{S:sec:s7-relative-product}, set
\[
 \Lambda_n=\sup_{\underline m_n\le k\le j\le n}
 \|\operatorname{Ad}_{P_{n;k,j}}\|_{\mathfrak b_\infty\to\mathfrak b_\infty}.
\]
Empty index ranges are ignored; all asymptotic statements concern $n$ sufficiently large that $\underline m_n<n$.

\begin{lemma}\label{S:lem:s7-adjoint}
Almost surely, for all sufficiently large $n$,
\begin{equation}\label{S:eq:s7-adjoint-bound}
 \Lambda_n\le C_{\beta,\omega}(\log n)^3.
\end{equation}
\end{lemma}

\begin{proof}
Take $T_l=2^l$ and $R_l=T_l^{5/4}$, for $l\ge1$. Lemma~\ref{S:lem:s7-height-ratio} bounds the probability that the first supremum below exceeds $R_l$ by $CT_l^{-1/4}$. A union
bound over the at most $T_l+1$ unit intervals, followed by Lemma~\ref{S:lem:s7-cauchy-maximum}, gives the same bound for $\sup_{t\le T_l}|z_t|$. These probabilities are summable in
$l$. Borel--Cantelli, monotonicity of the suprema, and the choice $T\le T_l<2T$ imply that, almost surely, for every
$T\ge2$,
\begin{equation}\label{S:eq:s7-continuous-bounds}
 \sup_{0\le s\le t\le T}r_{s,t}^2\le C_\omega T^{5/4},\qquad
 \sup_{0\le t\le T}(1+|z_t|)\le C_\omega T^{5/4}.
\end{equation}
The finitely many initial dyadic intervals are included because
\[
 \sup_{0\le s\le t\le 2^{l_0}}r_{s,t}^2
       +\sup_{0\le t\le2^{l_0}}(1+|z_t|)<\infty
\]
for every finite $l_0$, by continuity and positivity of $y$. Continuity also identifies these suprema with suprema over rational times, so the event is Borel.

In the coordinates $(u,v)$ of $\mathfrak b_\infty$, the matrix of the adjoint in~\eqref{S:eq:s7-adjoint-formula} is
\[
 \begin{pmatrix}1&0\\-2(z_t-r_{s,t}^2z_s)&r_{s,t}^2\end{pmatrix}.
\]
Its operator norm is at most its Frobenius norm, hence at most $1+2|z_t|+2r_{s,t}^2|z_s|+r_{s,t}^2$. Thus
\begin{equation}\label{S:eq:s7-continuous-adjoint}
 \sup_{0\le s\le t\le T}\|\operatorname{Ad}_{P_{s,t}}\|
 \le C_\omega T^{5/2},\qquad T\ge2.
\end{equation}
Multiplying the matrices $\mathcal A$ gives
\[
 P_{n;k,j}=\mathcal A(w_{\vartheta_{n,k}})
                    \mathcal A(w_{\vartheta_{n,j}})^{-1}
 =P_{\vartheta_{n,j},\vartheta_{n,k}},
\]
including $k=j$. Moreover,
\[
 \vartheta_{n,1}\le\frac4\beta\sum_{j=1}^{n-1}\frac1j
 \le\frac4\beta(1+\log n).
\]
Apply~\eqref{S:eq:s7-continuous-adjoint} with $T=\max\{2,(4/\beta)(1+\log n)\}$. It gives the stronger bound $\Lambda_n\le C_{\beta,\omega}(\log n)^{5/2}$ for $n\ge3$, and in
particular~\eqref{S:eq:s7-adjoint-bound}. The continuous bound holds at all the sampling times on the same event.
\end{proof}

The sampled columns on $I_{n,k}$ use the position with index $k+1$. Set $z_{n,k}=z_{\vartheta_{n,k+1}}$ and recall
\[
 \mathfrak q_{n,k}=
 \frac{q_\infty-x^{\rm aff}_{\vartheta_{n,k+1}}}
      {y^{\rm aff}_{\vartheta_{n,k+1}}}=-z_{n,k}.
\]
Equations~\eqref{S:eq:s7-continuous-bounds} and the bound on $\vartheta_{n,1}$ yield
\begin{equation}\label{S:eq:affine-terminal-coordinate-bound}
 \sup_{\underline m_n\le k<n}(1+|\mathfrak q_{n,k}|)
 \le C_{\beta,\omega}(\log n)^{5/4}.
\end{equation}
Multiplying the two squared column norms in \eqref{S:eq:adapted-boundary-norms} and using the reciprocal factors in~\eqref{S:eq:adapted-reciprocal-normalization} gives
\begin{equation}\label{S:eq:s7-column-product-bound}
 |a_{n,k}^{\rm aff}|^2|c_{n,k}^{\rm aff}|^2
 =|a_{n,k}^{\rm ad}|^2|c_{n,k}^{\rm ad}|^2
 =1+z_{n,k}^2\le C_{\beta,\omega}(\log n)^{5/2}.
\end{equation}

\subsubsection{Maximal inequalities for martingales}
\label{S:sec:martingale-moments-with-a-growing-exponent}

We need the dependence on $p$ in the following two martingale inequalities. The proofs also apply under a regular conditional law whenever the stated martingale and bracket
hypotheses hold under that law.

Let $M$ be a real continuous local martingale with $M_0=0$, and let $T<\infty$ be deterministic.

\begin{lemma}\label{S:lem:aux-martingale-moments}
There is an absolute $C$ such that, for every $p\ge1$,
\begin{equation}\label{S:eq:aux-bdg-growing-p}
 \left\|\sup_{t\le T}|M_t|\right\|_{L^p}
 \le C\sqrt p\,\|[M]_T^{1/2}\|_{L^p}.
\end{equation}
\end{lemma}

\begin{proof}
Write $M^*=\sup_{t\le T}|M_t|$ and $Q=[M]_T^{1/2}$. We first stop the martingale when either $|M|$ or its bracket exceeds a fixed bound. All moments used below are then finite; the
bounds we prove will not depend on this stopping level.

A continuous local martingale $N$ starting from zero satisfies
\[
 \Pr\{\sup_{s\le t}|N_s|\ge a,\ [N]_t\le v\}\le2e^{-a^2/(2v)}
 \qquad(a,v>0).
\]
To see this, stop at the first crossing of $a$ and at bracket level $v$. The exponential supermartingale $\exp\{\theta N_s-\theta^2[N]_s/2\}$ bounds the positive crossing
probability by $e^{-\theta a+\theta^2v/2}$. Choose $\theta=a/v$ and apply the same bound to $-N$. Optional sampling of the nonnegative exponential supermartingale gives the same
inequality conditionally after a stopping time, with the increment bracket. No independence of the increment from the past is required.

For $\lambda>0$, let $\tau$ be the first time $|M_t|=\lambda$. On $\{M^*>2\lambda,Q\le\delta\lambda\}$ one has $\tau\le T$, and after $\tau$ an increment of absolute value at least
$\lambda$ occurs with increment bracket at most $\delta^2\lambda^2$. The conditional inequality therefore gives
\[
 \Pr(M^*>2\lambda,Q\le\delta\lambda)
 \le\alpha_\delta\Pr(M^*\ge\lambda),\qquad
 \alpha_\delta=2e^{-1/(2\delta^2)}.
\]
Consequently $\Pr(M^*>2\lambda)\le\Pr(Q>\delta\lambda)+\alpha_\delta\Pr(M^*\ge\lambda)$. Multiply by $p\lambda^{p-1}$ and integrate over $\lambda>0$. The tail-integral formula
gives
\[
 (2^{-p}-\alpha_\delta)\mathbb E(M^*)^p
 \le\delta^{-p}\mathbb E Q^p.
\]
The use of $\ge$ instead of $>$ in one tail has no effect on its integral. Choose $\delta=[2(p+2)\log2]^{-1/2}$. Then $\alpha_\delta=2^{-p-1}$, and hence
\[
 \|M^*\|_p\le2^{1+1/p}\sqrt{2(p+2)\log2}\,\|Q\|_p
 \le C\sqrt p\,\|Q\|_p\qquad(p\ge1).
\]
Finally increase the stopping levels. Continuity makes the stopping times tend past each fixed finite horizon almost surely. Fatou's lemma applies to the maxima, and the stopped
brackets are at most $[M]_T$. This proves the estimate whenever its right side is finite, and the inequality is immediate in the extended sense when that side is infinite. The
argument uses only martingale and bracket identities, so it also applies under each conditional probability law having those properties.
\end{proof}

For a $d$-dimensional continuous local martingale starting at zero, apply the lemma to its coordinates and use $\sup_t|M_t|\le\sum_{i=1}^d\sup_t|M_t^i|$. Since
$[M^i]_T^{1/2}\le(\sum_j[M^j]_T)^{1/2}$, this gives
\[
 \left\|\sup_{t\le T}|M_t|\right\|_p
 \le Cd\sqrt p\left\|\left(\sum_{i=1}^d[M^i]_T\right)^{1/2}\right\|_p.
\]
The dimension is at most four in the matrix applications. The case $1\le p\le2$ is part of the lemma, and therefore also covers the $p/2$ moment used in the bridge estimates.

For the discrete inequality let $d_k,\ldots,d_N$ be real $L^p$ martingale differences for a filtration $(\mathcal G_j)$, where $p\ge2$, and put $M_m=\sum_{j=k}^m d_j$, $M_{k-1}=0$.

\begin{lemma}\label{S:lem:s7-scalar-maximal}
For the martingale just defined,
\begin{equation}\label{S:eq:s7-scalar-maximal}
 \left\|\max_{k\le m\le N}|M_m|\right\|_{L^p}
 \le2\sqrt{p-1}\left(\sum_{j=k}^N\|d_j\|_{L^p}^2\right)^{1/2}.
\end{equation}
\end{lemma}

\begin{proof}
We first prove a quadratic inequality for the $L^p$ norm. Let $X,Y\in L^p$ be real, $p\ge2$, and suppose $\mathbb E[Y\mid\mathcal A]=0$ with $X$ measurable with respect to
$\mathcal A$. For $\epsilon>0$ put
\[
 F_\epsilon(t)=\mathbb E\bigl((X+tY)^2+\epsilon^2\bigr)^{p/2},\qquad
 \phi_\epsilon(t)=F_\epsilon(t)^{2/p},\quad0\le t\le1.
\]
Differentiation is justified by H\"older and the integrable bound $C_p(|X|+|Y|+\epsilon)^p$. With $Z=X+tY$,
\[
 \begin{split}
 F_\epsilon'(t)&=p\mathbb E[Z(Z^2+\epsilon^2)^{p/2-1}Y],\\
 F_\epsilon''(t)&\le p(p-1)\mathbb E[(Z^2+\epsilon^2)^{(p-2)/2}Y^2]
 \le p(p-1)F_\epsilon(t)^{1-2/p}\|Y\|_p^2.
 \end{split}
\]
For $p=2$ the last exponent is zero and the inequality is an equality. Since $2/p-1\le0$, the square of $F_\epsilon'$ contributes a nonpositive term to $\phi_\epsilon''$. It
follows that $\phi_\epsilon''(t)\le2(p-1)\|Y\|_p^2$. Also $F_\epsilon'(0)=0$: its coefficient of $Y$ is $\mathcal A$-measurable and belongs to $L^{p/(p-1)}$, so conditional
centering applies, first after truncation and then by H\"older. Taylor's integral formula gives
\[
 \phi_\epsilon(1)\le\phi_\epsilon(0)+(p-1)\|Y\|_p^2.
\]
Letting $\epsilon\downarrow0$ proves $\|X+Y\|_p^2\le\|X\|_p^2+(p-1)\|Y\|_p^2$. Apply this with $X=M_{j-1}$, $Y=d_j$, and $\mathcal A=\mathcal G_{j-1}$, and sum in $j$:
\[
 \|M_N\|_p^2\le(p-1)\sum_{j=k}^N\|d_j\|_p^2.
\]

For completeness let $M^*=\max_{k\le m\le N}|M_m|$. Optional sampling for the submartingale $|M_m|$, at its first crossing of $\lambda$, gives
\[
 \lambda\Pr(M^*>\lambda)\le\mathbb E[|M_N|\mathbf1_{\{M^*>\lambda\}}].
\]
Integrate after multiplication by $p\lambda^{p-2}$ and use Tonelli. Since this is a finite sequence of $L^p$ variables, the maximum is in $L^p$, and H\"older yields
\[
 \mathbb E(M^*)^p\le\frac p{p-1}\mathbb E[|M_N|(M^*)^{p-1}]
 \le\frac p{p-1}\|M_N\|_p\|M^*\|_p^{p-1}.
\]
If the maximum has zero norm there is nothing to prove; otherwise divide by its $(p-1)$st power. Combining the two estimates gives
\[
 \|M^*\|_p\le\frac p{\sqrt{p-1}}
      \left(\sum_{j=k}^N\|d_j\|_p^2\right)^{1/2}
 \le2\sqrt{p-1}\left(\sum_{j=k}^N\|d_j\|_p^2\right)^{1/2},
\]
since $p/(p-1)\le2$ for $p\ge2$. This is the required inequality.
\end{proof}

\subsubsection{Products of upper triangular matrices}
\label{S:sec:s7-triangular-products}
Upper triangular structure reduces the product estimate to two scalar quantities: the diagonal is the exponential of a sum, and the remaining entry satisfies an exact recursion. To
apply the martingale inequality, we stop the exponential multiplier using the preceding diagonal sum. This preserves predictability and requires no exponential moments of the
increments.

Fix integers $1\le k\le N$, a real $p\ge4$, and a filtration $(\mathcal G_j)_{j=k-1}^N$. Let $\xi_j\in\mathfrak b_\infty$ be $\mathcal G_j$-measurable, and let
$b_j\in\mathfrak b_\infty$ be $\mathcal G_{j-1}$-measurable. Assume that, for a deterministic $K\ge1$,
\begin{equation}\label{S:eq:s7-product-moment-assumptions}
 \mathbb E[\xi_j\mid\mathcal G_{j-1}]=0,\qquad
 \mathbb E[\|\xi_j\|^p\mid\mathcal G_{j-1}]^{1/p}
 \le Kp^{3/2}j^{-3/2},\qquad
 \|b_j\|\le Kj^{-2}.
\end{equation}
Let $A_j:\mathfrak b_\infty\to\mathfrak b_\infty$ be $\mathcal G_{j-1}$-measurable linear maps, with $\|A_j\|\le L$ for a deterministic $L\ge1$. Set
\[
 \Pi_{k,m}=\prod_{j=k}^m\exp\{A_j(\xi_j+b_j)\},\qquad
 q_{p,L,k}=\frac{p^2L}{k}+\frac{p^3L^2}{k^2}.
\]
The factors are ordered increasingly from left to right, with $\Pi_{k,k-1}=I$. The proof uses the moment of order $p$ in \eqref{S:eq:s7-product-moment-assumptions}.

\begin{lemma}\label{S:lem:s7-product}
There are $C_K<\infty$ and an absolute $C_0<\infty$ such that, for $2C_Kq_{p,L,k}\le r\le1/16$,
\begin{equation}\label{S:eq:s7-product-tail}
 \mathbb P\!\left(\max_{k\le m\le N}\|\Pi_{k,m}-I\|>C_0r\right)
 \le3\left(\frac{C_Kq_{p,L,k}}r\right)^{p/2}.
\end{equation}
The constants are independent of $N,p,L,k$.
\end{lemma}

\begin{proof}
Write
\[
 A_j(\xi_j+b_j)=
 \begin{pmatrix}\alpha_j&\upsilon_j\\0&-\alpha_j\end{pmatrix},\qquad
 \alpha_j=\alpha_j^0+\alpha_j^b,\quad
 \upsilon_j=\upsilon_j^0+\upsilon_j^b,
\]
where superscripts $0$ and $b$ denote the coordinates of $A_j\xi_j$ and $A_jb_j$, respectively. Predictability of $A_j$ and conditional centering imply
\[
 \mathbb E[\alpha_j^0\mid\mathcal G_{j-1}]
 =\mathbb E[\upsilon_j^0\mid\mathcal G_{j-1}]=0.
\]
The two coordinates are fixed linear functionals on $\mathfrak b_\infty$. The norm comparison above and \eqref{S:eq:s7-product-moment-assumptions} give
\begin{equation}\label{S:eq:s7-product-coordinate-bounds}
 \|\alpha_j^0\|_p+\|\upsilon_j^0\|_p \le c_K Lp^{3/2}j^{-3/2},\qquad |\alpha_j^b|+|\upsilon_j^b|\le c_K Lj^{-2}.
\end{equation}
Here $c_K$ and the subsequent constants denoted by that symbol may be enlarged a finite number of times and depend only on $K$. The constant $C_K$ in the conclusion will be chosen
to dominate these bounds before the probability estimate is applied. The conditional covariance of the two centered coordinates need not vanish. With the superscripts just defined,
\[
 \mathbb E[\alpha_j\upsilon_j\mid\mathcal G_{j-1}]
 =\mathbb E[\alpha_j^0\upsilon_j^0\mid\mathcal G_{j-1}]
                         +\alpha_j^b\upsilon_j^b.
\]
The other two terms vanish by predictability and centering. Conditional H\"older bounds the first term by $c_KL^2p^3j^{-3}$, and the second by $c_KL^2j^{-4}$. The map $A_j$ can mix
the two coordinates; neither independence nor zero covariance at the same index is needed.

Put $S_m=\sum_{j=k}^m\alpha_j$ and $S_{k-1}=0$. For $a>1$, comparison of a decreasing function with its integral gives
\[
 \sum_{j=k}^\infty j^{-a}
 \le k^{-a}+\int_k^\infty x^{-a}\,\dd x
 \le\frac a{a-1}k^{1-a}.
\]
Apply Lemma~\ref{S:lem:s7-scalar-maximal} to $\alpha_j^0$ and sum the predictable parts absolutely. Using $a=3$ and $a=2$, respectively, yields
\begin{equation}\label{S:eq:s7-diagonal-max}
 \|S^*\|_p \le2\sqrt{p-1} \left(\sum_{j=k}^N\|\alpha_j^0\|_p^2\right)^{1/2} +\sum_{j=k}^N\|\alpha_j^b\|_\infty \le c_K\frac{p^2L}{k},\qquad S^*=\max_{k-1\le m\le N}|S_m|.
\end{equation}
In particular $|\alpha_j|=|S_j-S_{j-1}|\le2S^*$ for every $j$.

For $a,v\in\mathbb R$ set
\[
 N(v)=\begin{pmatrix}1&v\\0&1\end{pmatrix},\qquad
 A(a)=\begin{pmatrix}\me^a&0\\0&\me^{-a}\end{pmatrix},\qquad
 \phi(a)=\int_0^1\me^{2ta}\,\dd t.
\]
For $D=\left(\begin{smallmatrix}\alpha&\upsilon\\0&-\alpha\end{smallmatrix}\right)$, one has $D^2=\alpha^2I$. Summing
the even and odd powers gives
\[
 \exp D=
 \begin{pmatrix}
 \me^\alpha&\upsilon\,\sinh\alpha/\alpha\\0&\me^{-\alpha}
 \end{pmatrix}
 =N(\phi(\alpha)\upsilon)A(\alpha),
\]
with the quotient defined continuously at zero. The equality follows also from $\phi(a)=\me^a\sinh(a)/a$. Moreover,
\begin{equation}\label{S:eq:s7-triangular-exponential-remainder}
 |\phi(a)|\le\me^{2|a|},\qquad
 |\phi(a)-1|
 \le\int_0^1 2t|a|\me^{2t|a|}\,\dd t
 \le|a|\me^{2|a|}.
\end{equation}
The identity $A(a)N(v)=N(\me^{2a}v)A(a)$ now gives, by induction,
\begin{equation}\label{S:eq:s7-product-shear}
 \Pi_{k,m}=N(Z_m)A(S_m),\qquad
 Z_m=\sum_{j=k}^m\me^{2S_{j-1}}\phi(\alpha_j)\upsilon_j.
\end{equation}
Multiplication by the $m$th factor gives the recurrences
\[
 Z_m=Z_{m-1}+\me^{2S_{m-1}}\phi(\alpha_m)\upsilon_m,
 \qquad S_m=S_{m-1}+\alpha_m.
\]

To bound the predictable multiplier $\me^{2S_{j-1}}$, stop when $S$ first leaves the interval in the following definition. Define
\[
 I_j=\mathbf1_{\{\max_{k-1\le i<j}|S_i|\le1/4\}},\qquad
 M_m=\sum_{j=k}^m I_j\me^{2S_{j-1}}\upsilon_j^0.
\]
Both $I_j$ and $S_{j-1}$ are $\mathcal G_{j-1}$-measurable, and $|I_j\me^{2S_{j-1}}|\le\me^{1/2}$. Thus the summands in $M$ are martingale differences and
\[
 \|M^*\|_p \le2\sqrt{p-1}\,\me^{1/2} \left(\sum_{j=k}^N\|\upsilon_j^0\|_p^2\right)^{1/2} \le c_K\frac{p^2L}{k},\qquad M^*=\max_{k\le m\le N}|M_m|.
\]
For the remaining product of coordinates put $H=\sum_{j=k}^N|\alpha_j\upsilon_j|$. Since $p/2\ge1$, Minkowski's and H\"older's inequalities give
\begin{equation}\label{S:eq:s7-product-quadratic-bound}
 \|H\|_{p/2} \le\sum_{j=k}^N\|\alpha_j\upsilon_j\|_{p/2} \le\sum_{j=k}^N\|\alpha_j\|_p\|\upsilon_j\|_p \le c_Kp^3L^2\sum_{j=k}^\infty j^{-3} \le c_K\frac{p^3L^2}{k^2}.
\end{equation}
In the penultimate step the predictable coordinates are included: $j^{-2}\le p^{3/2}j^{-3/2}$ for $j\ge1$ and $p\ge4$. We also have $\sum_j|\upsilon_j^b|\le c_K L/k$. The order
$p/2$ in \eqref{S:eq:s7-product-quadratic-bound} is essential for this use of the available $p$th moments: for each $j$,
$\mathbb E|\alpha_j\upsilon_j|^{p/2} \le(\mathbb E|\alpha_j|^p)^{1/2} (\mathbb E|\upsilon_j|^p)^{1/2}$. The corresponding contribution in the exact shear recursion is
\[
 (\phi(\alpha_j)-1)\upsilon_j
 =2\alpha_j\upsilon_j
       \int_0^1\!\int_0^t\me^{2s\alpha_j}\,\dd s\,\dd t.
\]
On $|\alpha_j|\le2r$ this is bounded by $\me^{4r}|\alpha_j\upsilon_j|$. We take absolute values at this quadratic term, after the linear coordinates have been centered and their
predictable multiplier has been stopped.

Choose $C_K$ at least as large as all constants in the last three moment bounds and large enough that $\me^{1/8}\sum_j|\upsilon_j^b|\le C_K L/k$. On the event
\[
 E_r=\{S^*\le r,\ M^*\le r,\ H\le r\},\qquad r\le1/16,
\]
all $I_j$ equal one, $|\alpha_j|\le2r$, and $\me^{2|S_{j-1}|}\le\me^{2r}$. The exact difference from $M_m$ is
\[
 Z_m-M_m
 =\sum_{j=k}^m\me^{2S_{j-1}}\upsilon_j^b
  +\sum_{j=k}^m\me^{2S_{j-1}}
                  (\phi(\alpha_j)-1)\upsilon_j.
\]
Using~\eqref{S:eq:s7-triangular-exponential-remainder}, $2C_Kq_{p,L,k}\le r$, and $\me^{6r}\le\me^{3/8}<2$, gives
\[
 \max_m|Z_m|\le r+C_KL/k+\me^{6r}H\le4r.
\]
The diagonal and off-diagonal entries of $N(Z)A(S)-I$ satisfy
\[
 \|N(Z)A(S)-I\|
 \le\|N(Z)A(S)-I\|_{\rm HS}
 \le\me^{|S|}\sqrt{2S^2+Z^2}.
\]
Consequently $E_r$ implies $\max_m\|\Pi_{k,m}-I\|\le8r$. We may therefore take $C_0=8$. The complement of $E_r$ has probability at most
\[
 2\left(\frac{C_Kp^2L}{kr}\right)^p
 +\left(\frac{C_Kp^3L^2}{k^2r}\right)^{p/2}.
\]
If $a=C_Kq_{p,L,k}/r\le1/2$, this is bounded by $2a^p+a^{p/2}\le3a^{p/2}$. This proves \eqref{S:eq:s7-product-tail} with constants fixed independently of $N,p,L,k$.
\end{proof}

\subsubsection{Proof of the product estimate}
\label{S:sec:s7-proof-product-estimate}
\begin{proof}[Proof of Theorem~\ref{S:thm:s7-product-comparison}]
Fix an integer $M\ge1$ and let $L_{n,M}=M(\log n)^3$ for $n\ge3$. For each retained $k$ define
\[
 A_{n;k,j}^{(M)}
 =\mathbf1_{\{\max_{k\le i\le j}
          \|\operatorname{Ad}_{P_{n;k,i}}\|\le L_{n,M}\}}
       \operatorname{Ad}_{P_{n;k,j}},\qquad k\le j<n.
\]
For every $i\le j$, $P_{n;k,i}$ depends only on variables with indices less than $j$. The whole indicator is therefore $\mathcal G_{n;k,j-1}$-measurable. These maps have norm at
most $L_{n,M}$ on the entire probability space. For each stopped centered summand the tower property gives
\[
 \mathbb E[A_{n;k,j}^{(M)}\xi_{n,j}
                   \mid\mathcal G_{n;k,j-1}]=0.
\]
This identity holds on the whole space. Conditioning instead on $\{\Lambda_n\le L_{n,M}\}$ would generally change the law of $\xi_{n,j}$, since that event includes products using
the current and later increments. The stopped products avoid that conditioning. Set
\[
 \Pi_{n;k,m}^{(M)}
 =\prod_{j=k}^m
       \exp\{A_{n;k,j}^{(M)}(b_{n,j}+\xi_{n,j})\}.
\]
Under the law conditional on $\eta$, equation~\eqref{S:eq:s7-reversed-filtration-moments} verifies the moment assumptions of Lemma~\ref{S:lem:s7-product}. The stopped adjoints are
predictable and bounded by the deterministic number $L_{n,M}$, so that lemma applies to the displayed products.

For sufficiently large $n$, put
\[
 p_n=2\left\lceil\frac{64\log n}{\log\log n}\right\rceil,\qquad
 \varepsilon_{n,k}=\frac{(\log n)^{11/2}}k.
\]
Then $p_n\ge4$, $p_n/2\ge64\log n/\log\log n$, and $p_n\le130\log n/\log\log n$. For $k\ge\underline m_n\ge(\log n)^6$, direct substitution gives
\[
 \frac{p_n^2L_{n,M}}k
 \le\frac{C_M}{(\log n)(\log\log n)^2},\qquad
 \frac{p_n^3L_{n,M}^2}{k^2}
 \le\frac{C_M}{(\log n)^3(\log\log n)^3}.
\]
Dividing by $\varepsilon_{n,k}$ before bounding $k$ yields
\begin{equation}\label{S:eq:s7-product-threshold-ratio}
 \begin{aligned}
 \frac{q_{p_n,L_{n,M},k}}{\varepsilon_{n,k}}
 &=\frac{p_n^2L_{n,M}}{(\log n)^{11/2}}
   +\frac{p_n^3L_{n,M}^2}{k(\log n)^{11/2}}\\
 &\le C_M\left(
   \frac1{(\log n)^{1/2}(\log\log n)^2}
   +\frac1{(\log n)^{5/2}(\log\log n)^3}\right).
 \end{aligned}
\end{equation}
Also $\varepsilon_{n,k}\le(\log n)^{-1/2}$. Thus there is a deterministic $n_0(M,\beta)$ such that, for every $n\ge n_0(M,\beta)$ and all retained $k$,
\[
 2C_\beta q_{p_n,L_{n,M},k}\le\varepsilon_{n,k}\le1/16,\qquad
 \frac{C_\beta q_{p_n,L_{n,M},k}}{\varepsilon_{n,k}}
 \le(\log n)^{-1/4}.
\]
Here $C_\beta$ includes the constant in the product lemma. Consequently
\[
 \mathbb P\!\left( \max_{k\le m<n}\|\Pi_{n;k,m}^{(M)}-I\|>C_0\varepsilon_{n,k} \,\middle|\,\eta\right) \le3(\log n)^{-p_n/8} =3\exp\{-(p_n/8)\log\log n\}\le3n^{-16}.
 \]
The estimate is uniform in the conditioning value outside its null set. All indices $M,n,k$ and chosen exponents $p_n$ range over countable sets, so those conditional versions can
be fixed simultaneously outside one null set of terminal values. A union bound over the at most $n$ starting indices gives
\begin{equation}\label{S:eq:s7-stopped-row-probability}
 \mathbb P\!\left(
 \exists\,\underline m_n\le k<n:
 \max_{k\le m<n}\|\Pi_{n;k,m}^{(M)}-I\|>C_0\varepsilon_{n,k}
 \,\middle|\,\eta\right)\le3n^{-15}.
\end{equation}
The maximum over terminal indices $m$ is already included in the product estimate.

Let $L_n^{\rm fail}=\{\exists j\in[\underline m_n,n): \chi_{n,j}=0\}$. The exact radial tail from Theorem~\ref{S:thm:7.4} gives, with $c_\beta=\beta\log\cosh(r_{\rm loc}/2)>0$,
\[
 \mathbb P(L_n^{\rm fail}\mid\eta) \le\sum_{j=\underline m_n}^{n-1} \operatorname{sech}^{\beta j}(r_{\rm loc}/2) \le\frac{\exp(-c_\beta\underline m_n)}{1-\exp(-c_\beta)} \le C_\beta\exp(-c_\beta\underline m_n).
\]
All products in a row use this same set of increments. On $(L_n^{\rm fail})^c\cap\{\Lambda_n\le L_{n,M}\}$, $A_{n;k,j}^{(M)}=\operatorname{Ad}_{P_{n;k,j}}$ and
$b_{n,j}+\xi_{n,j}=\zeta_{n,j}^{\rm full}$ for every retained $k\le j$. Therefore~\eqref{S:eq:s7-exact-relative-product} gives $\Pi_{n;k,m}^{(M)}=\mathcal R_{n;k,m}$ simultaneously
for all $k\le m$. Integrating the conditional probability bound proves
\begin{equation}\label{S:eq:s7-exact-row-probability}
 \begin{aligned}
 &\mathbb P\!\left(\Lambda_n\le L_{n,M},
       \exists\,\underline m_n\le k<n:
       \max_{k\le m<n}\|\mathcal R_{n;k,m}-I\|
                       >C_0\varepsilon_{n,k}\right)\\
 &\hspace{2em}\le3n^{-15}+C_\beta\exp(-c_\beta\underline m_n).
 \end{aligned}
\end{equation}
This is summable in $n$: for large $n$, $c_\beta(\log n)^6\ge3\log n$.

Let $E_{n,M}$ denote the Borel event in \eqref{S:eq:s7-exact-row-probability}. Its probabilities are summable, so
\[
 \mathbb P\!\left(\bigcap_{M\in\mathbb N}
                   (\limsup_n E_{n,M})^c\right)=1.
\]
Intersect this event with that of Lemma~\ref{S:lem:s7-adjoint}. For each resulting path, choose an integer $M_\omega$ exceeding the constant in \eqref{S:eq:s7-adjoint-bound}. There
is $n_0(\omega)$ such that $\Lambda_n\le M_\omega(\log n)^3$ and $E_{n,M_\omega}$ does not occur for any $n\ge n_0(\omega)$. Consequently
\[
 \max_{k\le m<n}\|\mathcal R_{n;k,m}-I\|
 \le C_0\frac{(\log n)^{11/2}}k
 \qquad(\underline m_n\le k<n),
\]
which is \eqref{S:eq:s7-all-prefix-rate}.

\end{proof}

At $m=n-1$, the last identity in~\eqref{S:eq:s7-exact-relative-product} gives the estimate for the recursion matrices. On $I_{n,k}$ with $k\le n-2$, apply it at index $k+1$; on
$I_{n,n-1}$ the difference is zero. Thus, almost surely for all sufficiently large $n$,
\begin{equation}\label{S:eq:s7-relative-frame-rate}
 \|X_{n,k}^{\rm circ}(X_{n,k}^{\rm br})^{-1}-I\|
 \le C_{\beta,\omega}\frac{(\log n)^{11/2}}{k+1},
 \qquad\underline m_n\le k<n.
\end{equation}
In particular,
\begin{equation}\label{S:eq:s7-uniform-smallness}
 \sup_{\underline m_n\le k<n}
 \|X_{n,k}^{\rm circ}(X_{n,k}^{\rm br})^{-1}-I\|
 \le C_{\beta,\omega}(\log n)^{-1/2}\longrightarrow0.
\end{equation}
These estimates are in the common coordinates with terminal point at $\infty$. The column lengths change by the fixed reciprocal rescaling, so the same comparison holds in both
coordinate systems. Indeed, if the supremum in \eqref{S:eq:s7-uniform-smallness} is $\epsilon_n<1$, then for every $v\in\mathbb R^2$,
\[
 |(X_{n,k}^{\rm circ}-X_{n,k}^{\rm br})v|
 \le\epsilon_n|X_{n,k}^{\rm br}v|,
\]
and the triangle inequality gives
\[
 (1-\epsilon_n)|X_{n,k}^{\rm br}v|
 \le|X_{n,k}^{\rm circ}v|
 \le(1+\epsilon_n)|X_{n,k}^{\rm br}v|.
\]
Apply this to $u_0^{\rm ad}$ and $u_1^{\rm ad}$. For the $a$-columns, \eqref{S:eq:adapted-reciprocal-normalization} divides both sides by $\kappa_\eta$; for the $c$-columns it
multiplies both sides by $\kappa_\eta$. Thus the same inequalities hold for the lengths of the columns before the change of coordinates. Their squares and
Proposition~\ref{S:prop:5.2} give the estimates in Proposition~\ref{S:prop:finite-boundary-retained-intervals}.

\subsection{Comparison of the finite operator with the operator obtained by sampling the Brownian path}
\label{sec:combined-s-finite-comparison}
\label{S:sec:comparison-of-the-finite-and-sampled-brownian-operators}

The transfer product estimate was proved after sending the terminal point to infinity. To use it in the inverse kernels, we first recover the original matrix representation, then
sum the interval errors with the kernel weights and finally remove the boundary restriction. Throughout, $\beta>0$ is fixed and the coupling is that of
Proposition~\ref{S:prop:7.2}. The inverse operators use the matrix function $\mathcal A$ and the common boundary vectors $u_0,u_1$ from
Section~\ref{S:sec:green-kernels-for-canonical-systems}. With the cutoffs in~\eqref{S:eq:cutoff-choice}, write
\[
 u_n=\underline m_n,\qquad d_n=u_n/n,\qquad \delta_n=m_n/n.
\]
Thus $\underline T_n=1-d_n$ and $T_n=1-\delta_n$. All estimates involving these cutoffs are for $n$ large enough that
$1\le u_n\le m_n\le n/2$.

For the two step functions, distinguish the positions used on the partition intervals from the recursion indices:
\[
 \zeta_{n,k}^{\rm circ}=z_{n,k+1}^{\rm circ},\qquad
 \zeta_{n,k}^{\rm br}=Z_{\vartheta_{n,k+1}},\qquad
 X_{n,k}^{\diamond}=\mathcal A(\zeta_{n,k}^{\diamond}),
 \quad \diamond\in\{\mathrm{circ},\mathrm{br}\}.
\]
Here $0\le k<n$, and $X_{n,k}^{\diamond}$ is the value on $I_{n,k}$. A superscript $\mathrm{ad}$ denotes coordinates sending the terminal point to infinity; unadorned matrices use
the original upper half-plane coordinates. Reflection $r=1-t$ is the same unitary change of variables for all operators.

\subsubsection{The inverse comparison}
\label{S:sec:s8-inverse-comparison}
The sampling error determines the leading order coefficient only after the difference between the finite inverse and the inverse built from the sampled path has been shown to be
smaller. The product estimate controls this difference on the selected interval. The two boundary estimates control the portions of the kernel outside that interval.

Write
\[
 F_n=\mathsf R_{\beta,n}^{\mathrm{circ},T_n}
           -\mathsf R_{\beta,n}^{\mathrm{br},T_n},
 \qquad
 F_n^{\mathrm{full}}=\mathsf R_{\beta,n}^{\mathrm{circ}}
           -\mathsf R_{\beta,n}^{\mathrm{br}}.
\]
These are the retained and full finite errors. The product comparison applies only above the smaller cutoff; the omitted part requires the boundary estimates below.

\begin{proposition}\label{S:prop:circular-brownian-inverse-comparison}
Almost surely, for all sufficiently large $n$,
\begin{equation}\label{S:eq:circular-brownian-full-rate}
 \sqrt n\,\|F_n^{\mathrm{full}}\|_{\mathfrak S_2}
 \le C_{\beta,\omega}(\log n)^{-501/2}\longrightarrow0.
\end{equation}
\end{proposition}

The weighted inverse estimate determines the larger cutoff $m_n$. The smaller cutoff $u_n$ makes the relative products uniformly close to the identity; discarding the additional
intervals from $u_n$ to $m_n$ makes their summed squared contribution small. Both cutoffs are powers of $\log n$, so the boundary estimates can absorb the extra removal.

For every $\varepsilon>0$, the proof gives the following two contributions:
\[
 n\|F_n^{\mathrm{full}}\|_{\mathfrak S_2}^2
 \le C_{\beta,\omega,\varepsilon}
 \left(\frac{(\log n)^{11}}{m_n}+n^{-a_\beta+\varepsilon}\right),
 \qquad a_\beta=\min\{1,2/\beta\}.
\]
The first term comes from the weighted sum on the retained intervals; the second comes from the omitted kernel, including the rectangles with only one coordinate beyond the cutoff.
We prove the coordinate and boundary estimates before completing this calculation.

\subsubsection{Returning to the original coordinates}
\label{S:sec:s8-coordinate-change}
The two positions give different orthogonal factors under the fixed rotation $G_\eta$. To compare these two factors, multiply the transformed relative matrix by the orthogonal
factor associated with the Brownian position. The upper triangular factor in its QR decomposition is exactly the original relative matrix. The following bound is uniform in the
orthogonal factor and in the terminal point.

Let $R\in\mathrm{SL}_2(\mathbb R)$ and $O\in\mathrm{SO}_2$. Consider a factorization
\[
 RO=O'T,\qquad O'\in\mathrm{SO}_2,
 \qquad T\ \hbox{upper triangular with positive diagonal}.
\]
Put $d=\|R-I\|_{\mathrm{op}}$.

\begin{lemma}\label{S:lem:direct-qr}
The factorization above exists and is unique. If $d\le1/64$, then
\begin{equation}\label{S:eq:direct-qr-bound}
 \|T-I\|_{\mathrm{op}}\le9d,
 \qquad \|\log T\|_{\mathrm{op}}\le18d.
\end{equation}
The logarithm is real, upper triangular and traceless.
\end{lemma}

\begin{proof}
Set $H=O^{\mathsf T}R^{\mathsf T}RO=(h_{ij})$. This is a positive definite matrix of determinant one. If the
factorization exists, then $H=T^{\mathsf T}T$. Writing $T=\left(\begin{smallmatrix}a&b\\0&c\end{smallmatrix}\right)$,
with $a,c>0$, the first column and the mixed entry give $a=\sqrt{h_{11}}$ and $b=h_{12}/\sqrt{h_{11}}$. The determinant gives $c=a^{-1}$. Thus the only possible factor is
\begin{equation}\label{S:eq:direct-cholesky-formula}
 T=\begin{pmatrix}
 \sqrt{h_{11}}&h_{12}/\sqrt{h_{11}}\\
 0&1/\sqrt{h_{11}}
 \end{pmatrix}.
\end{equation}
Its remaining squared column length is $(h_{12}^2+1)/h_{11}=h_{22}$, since $h_{11}h_{22}-h_{12}^2=1$. Hence $T^{\mathsf T}T=H$. The matrix $O'=ROT^{-1}$ then satisfies
\[
 (O')^{\mathsf T}O'=T^{-\mathsf T}HT^{-1}=I,
 \qquad \det O'=1.
\]

Write $R=I+E$. Orthogonal invariance of the operator norm gives
\[
 \|H-I\|_{\mathrm{op}}
 =\|E+E^{\mathsf T}+E^{\mathsf T}E\|_{\mathrm{op}}
 \le2d+d^2\le3d=:\epsilon<\tfrac14.
\]
In particular $3/4\le h_{11}\le5/4$ and $|h_{12}|\le\epsilon$. For $x\in[3/4,5/4]$,
\[
 |\sqrt x-1|=\frac{|x-1|}{\sqrt x+1}\le|x-1|,
 \qquad
 |x^{-1/2}-1|=\frac{|x-1|}{\sqrt x(1+\sqrt x)}\le|x-1|.
\]
The second denominator is at least $\sqrt{3}/2+3/4>1$. Also $|h_{12}|/\sqrt{h_{11}}\le2\epsilon$. Formula \eqref{S:eq:direct-cholesky-formula} therefore gives
\[
 \|T-I\|_{\mathrm{op}}^2
 \le\|T-I\|_{\mathrm{HS}(\mathbb R^2)}^2
 \le\epsilon^2+4\epsilon^2+\epsilon^2
 \le54d^2\le81d^2.
\]
Since $\|T-I\|\le9/64<1/2$,
\[
 \log T=\sum_{j=1}^{\infty}\frac{(-1)^{j+1}}j(T-I)^j,
 \qquad
 \|\log T\|\le\frac{\|T-I\|}{1-\|T-I\|}\le2\|T-I\|.
\]
The series is absolutely convergent and consists of real upper triangular matrices. Its diagonal entries are $\log a$ and $\log(a^{-1})$, so its trace is zero. The scalar identity
for the power series $\exp(\log(1+z))=1+z$ applies to $T-I$ by absolute convergence, which proves $\exp(\log T)=T$. This also identifies the series with the real triangular
logarithm in~\eqref{S:eq:s6-triangular-logarithm}. The orthogonal factors can also be compared without identifying them. From $O'=RO T^{-1}$ and the preceding bound,
\[
 \|O'-O\|
 \le\|(R-I)O+O(I-T)\|\,\|T^{-1}\|
 \le\frac{10d}{1-9d}\le12d.
\]
The last inequality uses $d\le1/64$. The constant in the estimate is uniform in $O$, the terminal value, and the individual matrices $\mathcal A(z),\mathcal A(z^{\prime})$.
\end{proof}

Let $\eta=q_\infty$ and recall the rotation $G_\eta$ and the orthogonal matrix $O_\eta$ from \eqref{S:eq:adapted-reciprocal-normalization}. The terminal point is finite almost
surely. For $z=x+\ii y$, put $D=(\eta-x)^2+y^2$. The formula for the fractional linear action gives
\[
 \operatorname{Re}(G_\eta z)
 =\frac{(\eta x+1)(\eta-x)-\eta y^2}{D},
 \qquad
 \operatorname{Im}(G_\eta z)=\frac{(1+\eta^2)y}{D}.
\]
Substitution in the definition of $\mathcal A$ yields the explicit orthogonal matrix
\begin{equation}\label{S:eq:s8-orthogonal-factor}
 O_\eta(z)=\mathcal A(G_\eta z)G_\eta\mathcal A(z)^{-1}
 =\frac1{\sqrt D}
   \begin{pmatrix}\eta-x&y\\-y&\eta-x\end{pmatrix}.
\end{equation}
Its columns are orthonormal and its determinant is one. For two points $z,z'$, define
\[
 R^{\mathrm{ad}}=\mathcal A(G_\eta z')\mathcal A(G_\eta z)^{-1},
 \qquad
 T^{\mathrm{aff}}=\mathcal A(z')\mathcal A(z)^{-1}.
\]
Then
\begin{equation}\label{S:eq:direct-two-section-identity}
 R^{\mathrm{ad}}O_\eta(z)=O_\eta(z')T^{\mathrm{aff}}.
\end{equation}
Indeed, substitution of~\eqref{S:eq:s8-orthogonal-factor} makes both sides equal to $\mathcal A(G_\eta z')G_\eta\mathcal A(z)^{-1}$. The matrix $T^{\mathrm{aff}}$ is upper
triangular with positive diagonal, so Lemma~\ref{S:lem:direct-qr} applies with no dependence of its constants on $\eta,z,z'$. The calculation allows different factors $O_\eta(z)$
and
$O_\eta(z')$.

Apply this identity with $z=\zeta_{n,k}^{\rm br}$ and $z'=\zeta_{n,k}^{\rm circ}$, and write $R_{n,k}^{\mathrm{ad}}$ for the resulting relative matrix.
Theorem~\ref{S:thm:s7-product-comparison}, with the interval indexing in \eqref{S:eq:s7-relative-frame-rate}, gives
\[
 d_{n,k}:=\|R_{n,k}^{\mathrm{ad}}-I\|_{\mathrm{op}}
 \le C_{\beta,\omega}\frac{(\log n)^{11/2}}{k+1},
 \qquad u_n\le k<n.
\]
Its maximum is at most $C_{\beta,\omega}(\log n)^{-1/2}$. Consequently Lemma~\ref{S:lem:direct-qr} applies simultaneously to all these indices for sufficiently large $n$. Define
\[
 V_{n,k}=\log\bigl(X_{n,k}^{\rm circ}(X_{n,k}^{\rm br})^{-1}\bigr)
\]
using the real triangular logarithm. We have
\begin{equation}\label{S:eq:circular-brownian-affine-log}
 X_{n,k}^{\rm circ}=\me^{V_{n,k}}X_{n,k}^{\rm br},\qquad
 \|V_{n,k}\|_{\mathrm{op}}
 \le C_{\beta,\omega}\frac{(\log n)^{11/2}}{k+1},\qquad
 V_{n,n-1}=0.
\end{equation}
The last equality holds because both positions on the final interval are $\ii$. Each $V_{n,k}$ is real, upper triangular and traceless. The first inequality in
Lemma~\ref{S:lem:direct-qr} gives, at each index,
\[
 \|X_{n,k}^{\rm circ}(X_{n,k}^{\rm br})^{-1}-I\|_{\mathrm{op}}
 \le9d_{n,k}\le C_{\beta,\omega}\frac{(\log n)^{11/2}}{k+1}.
\]
Taking the maximum and using $u_n\ge(\log n)^6$ proves
\begin{equation}\label{S:eq:s8-affine-uniform-comparison}
 \epsilon_n:=\max_{u_n\le k<n}
    \|X_{n,k}^{\rm circ}(X_{n,k}^{\rm br})^{-1}-I\|_{\mathrm{op}}
 \le C_{\beta,\omega}(\log n)^{-1/2}\longrightarrow0.
\end{equation}
Thus the estimates hold for the original matrices that define the inverse kernels.

\subsubsection{Boundary columns above the cutoff}
\label{S:sec:s8-boundary-purpose}
Uniform smallness of the relative matrix allows us to deduce the bounds for $a_{n,k}^{\rm circ}$ and $c_{n,k}^{\rm circ}$ from those for $a_{n,k}^{\rm br}$ and $c_{n,k}^{\rm br}$.
The tail integral contains the squared length of one column multiplied by an integral of the squared length of the other column. These are the products we estimate when summing the
squared norm of the removed kernel.

\label{S:sec:s8-boundary-columns}

Use the boundary columns $a_{n,k}^{\diamond}=X_{n,k}^{\diamond}u_0$ and $c_{n,k}^{\diamond}=X_{n,k}^{\diamond}u_1$. Write
\[
 p_{n,k}^{\diamond}=|a_{n,k}^{\diamond}|^2,
 \qquad q_{n,k}^{\diamond}=|c_{n,k}^{\diamond}|^2,
 \qquad r_k=k/n,
\]
and retain $\gamma=2/\beta$ and $a_\beta=\min(1,\gamma)$ from Section~\ref{S:sec:sampling-times-and-estimates-for-the-boundary-solutions}. For the circular operator, sum the
squared lengths of the first boundary column over the intervals beginning with interval $k$, weighting each term by the cell length:
\[
 C_{n,k}^{\mathrm{circ},+}
       =\frac1n\sum_{\ell=k}^{n-1}p_{n,\ell}^{\mathrm{circ}}.
\]

\begin{proposition}\label{S:prop:finite-boundary-retained-intervals}
Almost surely, for every $\varepsilon>0$ and all sufficiently large $n$, uniformly for $u_n\le k<n$,
\[
 C_{\beta,\omega,\varepsilon}^{-1}r_k^{-\gamma+\varepsilon} \le p_{n,k}^{\mathrm{circ}} \le C_{\beta,\omega,\varepsilon}r_k^{-\gamma-\varepsilon},\qquad C_{\beta,\omega,\varepsilon}^{-1}r_k^{\gamma+\varepsilon} \le q_{n,k}^{\mathrm{circ}} \le C_{\beta,\omega,\varepsilon}r_k^{\gamma-\varepsilon}.
\]
\end{proposition}

\begin{proof}
For any $v\in\mathbb R^2$, the definition of $\epsilon_n$ gives
\[
 |X_{n,k}^{\rm circ}v-X_{n,k}^{\rm br}v|
       \le\epsilon_n|X_{n,k}^{\rm br}v|.
\]
For sufficiently large $n$, $\epsilon_n\le1/2$. Applying the two triangle inequalities and squaring therefore gives
\[
 (1-\epsilon_n)^2p_{n,k}^{\rm br} \le p_{n,k}^{\rm circ}\le(1+\epsilon_n)^2p_{n,k}^{\rm br},\qquad (1-\epsilon_n)^2q_{n,k}^{\rm br} \le q_{n,k}^{\rm circ}\le(1+\epsilon_n)^2q_{n,k}^{\rm br}.
\]
The boundary vectors are the same in these identities. By Proposition~\ref{S:prop:5.2}, the Brownian quantities have the stated bounds at $\bar r_{n,k}=(k+1)/n$. For $k\ge1$ and
any real exponent $\alpha$,
\[
 2^{-|\alpha|}r_k^\alpha\le\bar r_{n,k}^\alpha
                      \le2^{|\alpha|}r_k^\alpha,
 \qquad 1\le\bar r_{n,k}/r_k\le2.
\]
All exponents in the bounds for $p_{n,k}^{\rm br}$ and $q_{n,k}^{\rm br}$ have absolute value at most $\gamma+\varepsilon$, so these comparisons use at most the factor
$2^{\gamma+\varepsilon}$. The factors $(1-\epsilon_n)^2$ and $(1+\epsilon_n)^2$ lie between $1/4$ and $9/4$. Multiplying the Brownian bounds by these factors gives the stated
bounds for every retained $k$, with a common lower bound on the matrix size.
\end{proof}

Given $\varepsilon>0$, choose $\nu$ so that
\begin{equation}\label{S:eq:s8-choice-of-loss}
 0<\nu<\min(\varepsilon/4,a_\beta/4,1/4),
 \qquad \nu<(1-\gamma)/2\quad\hbox{when }\gamma<1.
\end{equation}
In the following calculations set $p_k=p_{n,k}^{\rm circ}$ and $q_k=q_{n,k}^{\rm circ}$. The preceding proposition gives $p_k\le Cr_k^{-\gamma-\nu}$ and $q_k\le Cr_k^{\gamma-\nu}$.
Since $x^{-\gamma-\nu}$ is decreasing,
\[
 \frac1n\sum_{\ell=k+1}^{n-1}p_\ell
 \le \frac Cn\sum_{\ell=k+1}^{n-1}(\ell/n)^{-\gamma-\nu}
 \le C\int_{r_k}^{1}x^{-\gamma-\nu}\dd x.
\]
For $\gamma>1$ this integral is $(r_k^{1-\gamma-\nu}-1)/(\gamma+\nu-1)$; for $\gamma=1$ it is $(r_k^{-\nu}-1)/\nu$; for $\gamma<1$ it is at most $(1-\gamma-\nu)^{-1}$. Thus
\begin{equation}\label{S:eq:s8-circular-primitive}
 \frac1n\sum_{\ell=k+1}^{n-1}p_\ell
 \le C_{\beta,\omega,\nu}
 \begin{cases}
 r_k^{1-\gamma-\nu},&\gamma>1,\\
 r_k^{-\nu},&\gamma=1,\\
 1,&0<\gamma<1.
 \end{cases}
\end{equation}
Multiplication by $q_k$ gives a bound by $Cr_k^{1-2\nu}$ when $\gamma\ge1$, and by $Cr_k^{\gamma-\nu}$ when $\gamma<1$. Both are at most $Cr_k^{a_\beta-2\nu}$. Also
$p_kq_k\le Cr_k^{-2\nu}$. The term included at $k$ satisfies
\[
 \frac{p_kq_k}{n}
 \le \frac Ck r_k^{1-2\nu}
 \le Cr_k^{a_\beta-2\nu}.
\]
Consequently
\begin{equation}\label{S:eq:s8-circular-weight-bounds}
 q_k C_{n,k}^{\mathrm{circ},+}\le Cr_k^{a_\beta-2\nu},
 \qquad p_kq_k\le Cr_k^{-2\nu},
 \qquad u_n\le k<n.
\end{equation}
Let $u_n\le m\le n/2$ and $\delta=m/n$. We claim that
\begin{equation}\label{S:eq:finite-retained-boundary-tail}
 \|\mathsf R_{\beta,n}^{\mathrm{circ},\underline T_n}
       -\mathsf R_{\beta,n}^{\mathrm{circ},1-\delta}\|_{\mathfrak S_2}^2
 \le C_{\beta,\omega,\varepsilon}\delta^{1+a_\beta-\varepsilon}.
\end{equation}
Lemma~\ref{S:lem:orthogonal-truncation-masses} gives
\begin{equation}\label{S:eq:s8-intermediate-kernel-sum}
 \|\mathsf R_{\beta,n}^{\mathrm{circ},\underline T_n} -\mathsf R_{\beta,n}^{\mathrm{circ},1-m/n}\|_{\mathfrak S_2}^2 =\frac1{2n^2}\sum_{k=u_n}^{m-1}q_k \sum_{\ell=k+1}^{n-1}p_\ell +\frac1{4n^2}\sum_{k=u_n}^{m-1}p_kq_k.
\end{equation}
In the reflected coordinates its support is $\{d_n\le\min(r,r')<\delta\}$. For $k<\ell$, the two transposed rectangles contribute $q_kp_\ell/(2n^2)$. On the square of index $k$,
the two triangles have total area $1/n^2$, and their squared kernel norm is $p_kq_k/4$. This proves the coefficients in the displayed formula, including the case $m=u_n$, when both
sums are empty.

Set $b=a_\beta-2\nu>0$. The first sum in \eqref{S:eq:s8-intermediate-kernel-sum} is at most
\[
 \frac Cn\sum_{k=1}^{m-1}(k/n)^b
 \le C\int_0^{m/n}x^b\dd x
 =\frac C{1+b}\delta^{1+b}.
\]
For the inequality, integrate each increasing power on $[k/n,(k+1)/n]$. For the diagonal sum, $0<2\nu<1$ and
\[
 \sum_{k=1}^{m-1}k^{-2\nu}
 \le1+\int_1^m x^{-2\nu}\dd x
 \le C_\nu m^{1-2\nu}.
\]
It is therefore bounded by
\[
 Cn^{-2+2\nu}m^{1-2\nu}
 =\frac Cn\delta^{1-2\nu}
 \le C\delta^{2-2\nu}.
\]
Here $n^{-1}\le\delta$. Since $a_\beta\le1$ and $2\nu<\varepsilon$, both exponents are at least $1+a_\beta-\varepsilon$. This proves~\eqref{S:eq:finite-retained-boundary-tail}. The
constants and the starting index were chosen before $m$, so the estimate holds simultaneously for all the stated cutoffs.

\subsubsection{The remaining intervals near the boundary}
\label{S:sec:s8-remaining-intervals}

Below $u_n$ we estimate the auxiliary circular continuation by its oscillation relative to its value at the cutoff. This compares the terminal columns with the column at the
cutoff. Multiplication by the measure of the removed intervals then bounds their contribution to the kernel.

Let $M_n=C_\beta^{\mathrm{tail}}(1+\log\log n)^4$ be the deterministic sequence in~\eqref{S:eq:auxiliary-tail-oscillation}. We use the oscillation estimate for the coupled circular
path from Proposition~\ref{S:prop:7.2}; the construction and its uniformity over $n$ were established in Lemma~\ref{S:lem:auxiliary-boundary-oscillation} and
Section~\ref{S:sec:coupling-the-finite-transfer-increments-to-brownian-increments}.

\begin{proposition}\label{S:prop:finite-extreme-tail}
For $0<\varepsilon<a_\beta$, almost surely for all sufficiently large $n$,
\begin{equation}\label{S:eq:finite-extreme-tail}
 \|\mathsf R_{\beta,n}^{\mathrm{circ}}
       -\mathsf R_{\beta,n}^{\mathrm{circ},\underline T_n}\|_{\mathfrak S_2}^2
 \le C_{\beta,\omega,\varepsilon}
           \me^{2M_n}d_n^{1+a_\beta-\varepsilon}.
\end{equation}
\end{proposition}

\begin{proof}
Write $u=u_n$ and $\delta=d_n$. The intervals with $0\le k<u$ use precisely the recursion positions $z_{n,1}^{\rm circ},\ldots, z_{n,u}^{\rm circ}$. Interval $u$ uses the position
$z_{n,u+1}^{\rm circ}$. The distance between the positions with indices $u$ and $u+1$ is the circular radius $Y_{n,u}$. By Corollary~\ref{S:cor:increment-logarithm-bounds},
\[
 Y_{n,u}\le C_\beta\sqrt{\frac{\log n}{u}}
       \le C_\beta(\log n)^{-5/2}\longrightarrow0.
\]
The triangle inequality and~\eqref{S:eq:spliced-tail-oscillation} give
\begin{equation}\label{S:eq:boundary-tail-mesh-anchor}
 \max_{0\le k<u}
 d_{\mathbb H}(\zeta_{n,k}^{\rm circ},\zeta_{n,u}^{\rm circ})
 \le M_n+Y_{n,u}\le M_n+1
\end{equation}
for all sufficiently large $n$. This includes the increment between indices $u$ and $u+1$; all points in the maximum are finite positions of the walk.

For arbitrary $z=x+\ii y$ and $w=x'+\ii y'$ in $\mathbb H$,
\[
 D_{z,w}:=\mathcal A(z)\mathcal A(w)^{-1}
 =\begin{pmatrix}
 \sqrt{y'/y}&(x'-x)/\sqrt{yy'}\\0&\sqrt{y/y'}
 \end{pmatrix}.
\]
This matrix has determinant one, and
\[
 \operatorname{tr}(D_{z,w}^{\mathsf T}D_{z,w})
 =\frac{y'}y+\frac y{y'}+\frac{(x-x')^2}{yy'}
 =2\cosh d_{\mathbb H}(z,w).
\]
The last equality follows by applying the isometry $\zeta\mapsto(\zeta-\operatorname{Re}w)/\operatorname{Im}w$. If the resulting point is $X+\ii Y$, its disk radius is
$\rho=|(X+\ii Y-\ii)/(X+\ii Y+\ii)|$. The identity $d=2\operatorname{arctanh}\rho$ gives $\cosh d=(1+\rho^2)/(1-\rho^2)=(X^2+Y^2+1)/(2Y)$, which yields the displayed formula. The
eigenvalues of $D_{z,w}^{\mathsf T}D_{z,w}$ are therefore
\[
 \cosh d\pm\sqrt{\cosh^2d-1}=\me^{\pm d},
 \qquad d=d_{\mathbb H}(z,w).
\]
Applying their bounds to $\mathcal A(w)v$ proves, for every $v\in\mathbb R^2$,
\begin{equation}\label{S:eq:boundary-busemann-column-inequality}
 \me^{-d_{\mathbb H}(z,w)}|\mathcal A(w)v|^2
 \le|\mathcal A(z)v|^2
 \le\me^{d_{\mathbb H}(z,w)}|\mathcal A(w)v|^2.
\end{equation}
With $p_k=p_{n,k}^{\rm circ}$ and $q_k=q_{n,k}^{\rm circ}$, apply this to $u_0,u_1$ and then use \eqref{S:eq:boundary-tail-mesh-anchor}. It gives
\begin{equation}\label{S:eq:extreme-tail-column-comparison}
 p_k\le\me^{M_n+1}p_u,\qquad
 q_k\le\me^{M_n+1}q_u,\qquad 0\le k<u.
\end{equation}
The inequality holds for every $v$ and hence applies pathwise to the random boundary vector $u_1$.

The exact formula for the Green kernel is
\begin{equation}\label{S:eq:extreme-tail-green-sum}
 \|\mathsf R_{\beta,n}^{\mathrm{circ}} -\mathsf R_{\beta,n}^{\mathrm{circ},\underline T_n}\|_{\mathfrak S_2}^2 =\frac1{2n^2}\sum_{k=0}^{u-1}q_k \sum_{\ell=k+1}^{n-1}p_\ell +\frac1{4n^2}\sum_{k=0}^{u-1}p_kq_k.
\end{equation}
Split the inner sum into $\ell<u$ and $\ell\ge u$. For the second part, there are $u$ possible values of $k$, so
\[
 \frac1{2n^2}\sum_{k=0}^{u-1}q_k\sum_{\ell=u}^{n-1}p_\ell
 \le\frac\delta2\me^{M_n+1}q_u C_{n,u}^{\mathrm{circ},+}.
\]
For the first part, there are $u(u-1)/2$ pairs $0\le k<\ell<u$. Together with the $u$ diagonal squares, their contribution is at most
\[
 \frac{u(u-1)+u}{4n^2}\me^{2(M_n+1)}p_uq_u
 =\frac{\delta^2}4\me^{2(M_n+1)}p_uq_u.
\]
Thus
\begin{equation}\label{S:eq:extreme-tail-cutoff-bound}
 \|\mathsf R_{\beta,n}^{\mathrm{circ}} -\mathsf R_{\beta,n}^{\mathrm{circ},\underline T_n}\|_{\mathfrak S_2}^2 \le\frac\delta2\me^{M_n+1}q_u C_{n,u}^{\mathrm{circ},+} +\frac{\delta^2}4\me^{2(M_n+1)}p_uq_u.
\end{equation}
The first term includes all rectangles with exactly one index below $u$. The second includes the entire square with both indices below $u$, including index $k=0$.

Choose $\nu$ as in~\eqref{S:eq:s8-choice-of-loss}. At $k=u$, equation~\eqref{S:eq:s8-circular-weight-bounds} gives $q_u C_{n,u}^{\mathrm{circ},+}\le C\delta^{a_\beta-2\nu}$ and
$p_uq_u\le C\delta^{-2\nu}$. The included summand $p_u/n$ is covered by that equation; explicitly,
\[
 \frac{p_uq_u}{n}\le\frac C{u}\delta^{1-2\nu}
       \le C\delta^{a_\beta-2\nu}.
\]
Since $2\nu<\varepsilon$ and $0<\delta\le1$, these bounds imply
\begin{equation}\label{S:eq:extreme-tail-cutoff-weights}
 q_u C_{n,u}^{\mathrm{circ},+}
       \le C_{\beta,\omega,\varepsilon}\delta^{a_\beta-\varepsilon},
 \qquad
 p_uq_u\le C_{\beta,\omega,\varepsilon}\delta^{-\varepsilon}.
\end{equation}
Substitute these estimates into~\eqref{S:eq:extreme-tail-cutoff-bound}. Because $M_n\ge0$, $a_\beta\le1$ and $\delta\le1$,
\[
 \me^{M_n+1}\delta^{1+a_\beta-\varepsilon}
       +\me^{2(M_n+1)}\delta^{2-\varepsilon}
 \le(\me+\me^2)\me^{2M_n}\delta^{1+a_\beta-\varepsilon}.
\]
\end{proof}

The bound in~\eqref{S:eq:finite-extreme-tail}, multiplied by $n$, tends to zero. Indeed, for
$0<\varepsilon<a_\beta$,
\[
 \log\!\left(n\me^{2M_n}d_n^{1+a_\beta-\varepsilon}\right)
 =-(a_\beta-\varepsilon)\log n+2M_n
                 +(1+a_\beta-\varepsilon)\log u_n.
\]
Here $\log u_n=O(\log\log n)$ and $M_n=O_\beta((1+\log\log n)^4)=o(\log n)$. Choose $n_0(\beta,\varepsilon)$ so that the last two terms are at most $(a_\beta-\varepsilon)\log n/2$
for $n\ge n_0(\beta,\varepsilon)$. The expression inside the logarithm is therefore at most $n^{-(a_\beta-\varepsilon)/2}$ for all sufficiently large $n$.

\subsubsection{Removal of the cutoff}
\label{S:sec:s8-removal-of-cutoff}
The two circular tail estimates apply on disjoint parts of the kernel square. Their squared norms therefore add exactly. We combine them with the Brownian tails on one event, which
will also be used after adjoining the auxiliary Gaussian variables. This step concerns the original operator laws and requires no independence between an interior error and a
boundary error.

For this subsection write $\mathsf R_{\beta,n}^{\mathrm{sine}}=\mathsf R_\beta^{\mathrm{sine}}$. Define
\begin{equation}\label{S:eq:boundary-tail-operator-definitions}
 \mathsf T_{\beta,n}^{\diamond} =\mathsf R_{\beta,n}^{\diamond} -\mathsf R_{\beta,n}^{\diamond,T_n}, \diamond\in\{\mathrm{circ},\mathrm{br},\mathrm{sine}\},\qquad \mathfrak t_n^{\mathrm{end}} =n\sum_{\diamond\in\{\mathrm{circ},\mathrm{br},\mathrm{sine}\}} \|\mathsf T_{\beta,n}^{\diamond}\|_{\mathfrak S_2}^2.
\end{equation}
By the restriction theorem of Section~\ref{S:sec:green-kernels-for-canonical-systems}, each truncated kernel is the regular Green inverse on $[0,T_n]$ with the same coordinate
boundary vectors, extended by zero.

Let $\Omega_\beta^{\mathrm{end}}$ be the intersection of the event on which the Brownian path is continuous and all inverse operators are defined and the events in
Theorem~\ref{S:thm:5.1}, Propositions~\ref{S:prop:5.2}, \ref{S:prop:5.4} and~\ref{S:prop:4.3}, Corollary~\ref{S:cor:increment-logarithm-bounds},
Theorem~\ref{S:thm:s7-product-comparison}, the bound \eqref{S:eq:s7-column-product-bound}, and \eqref{S:eq:spliced-tail-oscillation}. For every exponent loss in these statements
take the intersection over positive rational values. The event has probability one for the fixed $\beta$. It can be chosen Borel on the probability space before completion. In
fact, if the bound is $f_{n,k}\le C_\omega g_{n,k}$ for all $n\ge n_0(\omega)$ for finitely many indices $k\in K_n$ per row, its event is
\[
 \bigcup_{M,N\in\mathbb N}\ \bigcap_{n\ge N}\ \bigcap_{k\in K_n}
       \{f_{n,k}\le M g_{n,k}\}.
\]
The functions here are measurable functions of the coupled positions and increments. For bounds along a continuous path, a countable dense set of times or of positive $r$ suffices,
by continuity. The strong measurability of the inverse operators was proved in Sections~\ref{S:sec:green-kernels-for-canonical-systems} and
\ref{S:sec:sampling-times-and-estimates-for-the-boundary-solutions}. These observations justify the countable construction of the event. On it, any specified real exponent loss
follows from a smaller positive rational one. The estimates already proved in this subsection then hold on this same event.

\begin{theorem}\label{S:thm:full-boundary-tails}
On $\Omega_\beta^{\mathrm{end}}$, for every $\varepsilon>0$ and all sufficiently large $n$,
\begin{equation}\label{S:eq:boundary-tail-polynomial-rate}
 \mathfrak t_n^{\mathrm{end}}
       \le C_{\beta,\omega,\varepsilon}n^{-a_\beta+\varepsilon}.
\end{equation}
In particular, $\mathfrak t_n^{\mathrm{end}}\longrightarrow0$.
\end{theorem}

\begin{proof}
Fix $0<\nu<a_\beta$ and work on the event specified above. The removed part of the circular kernel is the sum of its restrictions to
\[
 \{\min(r,r')<d_n\},\qquad
 \{d_n\le\min(r,r')<\delta_n\}.
\]
Their supports are disjoint up to null sets. Therefore
\begin{equation}\label{S:eq:boundary-circular-two-tail-identity}
 \|\mathsf T_{\beta,n}^{\mathrm{circ}}\|_{\mathfrak S_2}^2 ={}\|\mathsf R_{\beta,n}^{\mathrm{circ}} -\mathsf R_{\beta,n}^{\mathrm{circ},\underline T_n}\|_{\mathfrak S_2}^2 +\|\mathsf R_{\beta,n}^{\mathrm{circ},\underline T_n} -\mathsf R_{\beta,n}^{\mathrm{circ},T_n}\|_{\mathfrak S_2}^2.
 \end{equation}
Proposition~\ref{S:prop:finite-extreme-tail} bounds the first term by $C\me^{2M_n}d_n^{1+a_\beta-\nu}$, and \eqref{S:eq:finite-retained-boundary-tail}, with $m=m_n$, bounds the
second by $C\delta_n^{1+a_\beta-\nu}$. By Proposition~\ref{S:prop:4.3}, each of the two Brownian tail norms squared is at most $C\delta_n^{1+a_\beta-\nu}$. Adding these three
estimates gives
\begin{equation}\label{S:eq:boundary-tail-deterministic-majorant}
 \begin{aligned}
 \mathfrak t_n^{\mathrm{end}}
     &\le C_{\beta,\omega,\nu} b_{\beta,n}^{\mathrm{end}}(\nu),\\
 b_{\beta,n}^{\mathrm{end}}(\nu)
     &:=n\left[\delta_n^{1+a_\beta-\nu}
                   +\me^{2M_n}d_n^{1+a_\beta-\nu}\right].
 \end{aligned}
\end{equation}
The disjoint kernel regions in this sum include all mixed rectangles.

Put $s=1+a_\beta-\nu>0$. The ceilings give $m_n\le2(\log n)^{512}$ and $u_n\le2(\log n)^6$ for large $n$. It follows that
\begin{equation}\label{S:eq:boundary-majorant-explicit-growth}
 b_{\beta,n}^{\mathrm{end}}(\nu) \le C_{\beta,\nu}n^{-a_\beta+\nu} \bigl[(\log n)^{512s} +\me^{2C_\beta^{\mathrm{tail}}(1+\log\log n)^4} (\log n)^{6s}\bigr].
 \end{equation}
Since
\[
 \frac{\log\log n+(1+\log\log n)^4}{\log n}\longrightarrow0,
\]
the bracket in \eqref{S:eq:boundary-majorant-explicit-growth} is $n^{o(1)}$. For any $\varepsilon>0$, choose rational $0<\nu<\min(a_\beta/2,\varepsilon/2)$. The bracket is at most
$n^{\varepsilon-\nu}$ for every $n\ge n_0(\omega,\varepsilon,\nu)$, so \eqref{S:eq:boundary-tail-deterministic-majorant} implies
\[
 \mathfrak t_n^{\mathrm{end}}
 \le C_{\beta,\omega,\nu}n^{-a_\beta+\nu}
                          n^{\varepsilon-\nu}
 =C_{\beta,\omega,\nu}n^{-a_\beta+\varepsilon}.
\]
Taking $\varepsilon<a_\beta$ also proves convergence to zero.

\end{proof}

In particular,
\begin{equation}\label{S:eq:4.5}
 n\|\mathsf R_{\beta,n}^{\mathrm{circ}} -\mathsf R_{\beta,n}^{\mathrm{circ},T_n}\|_{\mathfrak S_2}^2 +n\|\mathsf R_{\beta,n}^{\mathrm{br}} -\mathsf R_{\beta,n}^{\mathrm{br},T_n}\|_{\mathfrak S_2}^2 +n\|\mathsf R_\beta^{\mathrm{sine}} -\mathsf R_\beta^{\mathrm{sine},T_n}\|_{\mathfrak S_2}^2 \longrightarrow0.
 \end{equation}
These conclusions persist under any extension preserving the original joint law. More precisely, if $\pi$ maps an extended probability space to the original one and
$\pi_*\widetilde{\mathbb P}=\mathbb P$, then
\[
 \widetilde{\mathbb P}(\pi^{-1}\Omega_\beta^{\mathrm{end}})
 =\mathbb P(\Omega_\beta^{\mathrm{end}})=1.
\]
The estimates therefore remain valid on the extension of Theorem~\ref{S:thm:joint-probability-space}, after composition of the original variables with $\pi$.

For $\alpha,\chi\in\{\mathrm{circ},\mathrm{br},\mathrm{sine}\}$, put
\[
 D_n^{\alpha,\chi}
 =\mathsf R_{\beta,n}^{\alpha}-\mathsf R_{\beta,n}^{\chi},
 \qquad
 D_n^{\alpha,\chi,T}
 =\mathsf R_{\beta,n}^{\alpha,T_n}
       -\mathsf R_{\beta,n}^{\chi,T_n}.
\]
After the common reflection, $D_n^{\alpha,\chi,T}=\Pi_{\delta_n}D_n^{\alpha,\chi}$.

\begin{corollary}\label{S:cor:4.4}
On $\Omega_\beta^{\mathrm{end}}$,
\begin{equation}\label{S:eq:omitted-kernel-difference}
 0\le n\bigl(\|D_n^{\alpha,\chi}\|_{\mathfrak S_2}^2 -\|D_n^{\alpha,\chi,T}\|_{\mathfrak S_2}^2\bigr) =n\|\mathsf T_{\beta,n}^{\alpha} -\mathsf T_{\beta,n}^{\chi}\|_{\mathfrak S_2}^2\le2\mathfrak t_n^{\mathrm{end}}.
\end{equation}
\end{corollary}

\begin{proof}
Fix an outcome in $\Omega_\beta^{\mathrm{end}}$ and apply the common reflection to every operator. This unitary change of variables preserves the Hilbert--Schmidt norm. Write
$K_n(r,s)$ for the matrix-valued kernel of $D_n^{\alpha,\chi}$ and set
\[
 A_n=[\delta_n,1]^2,\qquad
 B_n=(0,1)^2\setminus A_n.
\]
Boundary lines have two-dimensional Lebesgue measure zero. The exact restriction identity gives the kernels
\[
 K_n^T=\mathbf1_{A_n}K_n,\qquad
 K_n-K_n^T=\mathbf1_{B_n}K_n.
\]
These kernels belong to $L^2((0,1)^2;\mathbb C^{2\times2})$, since both full inverses are Hilbert--Schmidt. Their inner product vanishes:
\[
 \int_{(0,1)^2}\operatorname{Tr}
       \bigl((\mathbf1_{B_n}K_n)^*(\mathbf1_{A_n}K_n)\bigr)
                      \,\mathrm dr\,\mathrm ds=0.
\]
Consequently the integral formula for the Hilbert--Schmidt norm yields
\[
 \|D_n^{\alpha,\chi}\|_{\mathfrak S_2}^2
 =\|D_n^{\alpha,\chi,T}\|_{\mathfrak S_2}^2
  +\|D_n^{\alpha,\chi}-D_n^{\alpha,\chi,T}\|_{\mathfrak S_2}^2.
\]
Subtracting the definitions of the truncated differences gives
\[
 D_n^{\alpha,\chi}-D_n^{\alpha,\chi,T}
       =\mathsf T_{\beta,n}^{\alpha}-\mathsf T_{\beta,n}^{\chi}.
\]
This proves the equality and the nonnegativity in \eqref{S:eq:omitted-kernel-difference}. In particular $B_n$ includes both rectangles with exactly one variable below $\delta_n$,
as well as the square with both variables below it.

If $\alpha=\chi$, all three differences vanish. Otherwise the Hilbert-space inequality $\|U-V\|^2\le2\|U\|^2+2\|V\|^2$ gives
\[
 n\|\mathsf T_{\beta,n}^{\alpha}-\mathsf T_{\beta,n}^{\chi}\|_{\mathfrak S_2}^2
 \le2n\bigl(\|\mathsf T_{\beta,n}^{\alpha}\|_{\mathfrak S_2}^2
              +\|\mathsf T_{\beta,n}^{\chi}\|_{\mathfrak S_2}^2\bigr)
 \le2\mathfrak t_n^{\mathrm{end}}.
\]
The last inequality uses two distinct nonnegative summands of the defining sum for $\mathfrak t_n^{\mathrm{end}}$. Every step is deterministic on the chosen outcome; no
independence of the two tails is used.
\end{proof}

By Theorem~\ref{S:thm:full-boundary-tails}, the right side of \eqref{S:eq:omitted-kernel-difference} is at most $C_{\beta,\omega,\varepsilon}n^{-a_\beta+\varepsilon}$. Taking
$(\alpha,\chi)=(\mathrm{circ},\mathrm{sine})$ and $(\mathrm{br},\mathrm{sine})$, respectively, gives
\[
 n\bigl(\|\Delta_{\beta,n}\|_{\mathfrak S_2}^2 -\|\Delta_{\beta,n}^{T}\|_{\mathfrak S_2}^2\bigr) \longrightarrow0,\qquad n\bigl(\|\mathsf R_{\beta,n}^{\mathrm{br}} -\mathsf R_\beta^{\mathrm{sine}}\|_{\mathfrak S_2}^2 -\|\varepsilon_{\beta,n}^{\mathrm{samp},T}\|_{\mathfrak S_2}^2\bigr) \longrightarrow0.
\]
The same bound also gives, for every $\varepsilon>0$,
\begin{equation}\label{S:eq:s8-tail-difference-norm}
 \sqrt n\,\|D_n^{\alpha,\chi}-D_n^{\alpha,\chi,T}\|_{\mathfrak S_2}
 \le\sqrt{2\mathfrak t_n^{\mathrm{end}}}
 \le C_{\beta,\omega,\varepsilon}n^{-a_\beta/2+\varepsilon}
\end{equation}
for sufficiently large $n$, by using the loss $2\varepsilon$ in the theorem. The exact orthogonal decomposition gives these bounds before estimating the full normalized difference.

\subsubsection{Proof of the inverse comparison}
\label{S:sec:s8-proof-inverse-comparison}
\begin{proof}[Proof of Proposition~\ref{S:prop:circular-brownian-inverse-comparison}]
Let $V_n$ be the restriction of \eqref{S:eq:circular-brownian-affine-log} to $m_n\le k<n$. All weights and columns without an ensemble superscript in this proof are computed from
the sampled Brownian matrix function on those intervals. In particular, $p_{n,k}=|a_{n,k}^{\rm br}|^2$ and $q_{n,k}=|c_{n,k}^{\rm br}|^2$. The weights of
Section~\ref{S:sec:stability-and-fr-echet-expansion-of-the-inverse-kernel} satisfy
\[
 \mathfrak G_{n,k}^{+}
 =A_{n,k}^{+}p_{n,k}+C_{n,k}^{+}q_{n,k}
 =\mathfrak G_{n,k}+\frac2n p_{n,k}q_{n,k}.
\]
Proposition~\ref{S:prop:5.4}, applied with $\varepsilon=a_\beta/2$, gives $\mathfrak G_{n,k}\le Cr_k^{a_\beta/2}\le C$. Equation~\eqref{S:eq:s7-column-product-bound} gives
$p_{n,k}q_{n,k}\le C(\log n)^{5/2}$. Hence
\begin{equation}\label{S:eq:s8-full-weight-bound}
 \max_{m_n\le k<n}\mathfrak G_{n,k}^{+}
 \le C_{\beta,\omega}
           \left(1+\frac{(\log n)^{5/2}}n\right)
 \le C_{\beta,\omega}.
\end{equation}
The second term controls the contributions from the diagonal squares; it is present also when $A_{n,k}=0$ or
$C_{n,k}=0$.

By~\eqref{S:eq:circular-brownian-affine-log},
\[
 \eta_n:=\max_{m_n\le k<n}\|V_{n,k}\|_{\mathrm{op}}
 \le C_{\beta,\omega}\frac{(\log n)^{11/2}}{m_n+1}
 \longrightarrow0.
\]
Substituting the bound for each $\|V_{n,k}\|$ in the full quadratic form and using \eqref{S:eq:s8-full-weight-bound} gives
\begin{equation}\label{S:eq:circular-brownian-weighted-energy}
 \begin{aligned}
 \mathcal E_n^+(V_n)
 &=\sum_{k=m_n}^{n-1}
    \bigl[A_{n,k}^{+}|V_{n,k}a_{n,k}^{\rm br}|^2
                  +C_{n,k}^{+}|V_{n,k}c_{n,k}^{\rm br}|^2\bigr]\\
 &\le\sum_{k=m_n}^{n-1}\mathfrak G_{n,k}^{+}
                           \|V_{n,k}\|_{\mathrm{op}}^2\\
 &\le C_{\beta,\omega}(\log n)^{11}
                         \sum_{k=m_n}^{n-1}(k+1)^{-2}
 \le C_{\beta,\omega}\frac{(\log n)^{11}}{m_n}.
 \end{aligned}
\end{equation}
For the last inequality, integrate $x^{-2}$ on each $[k,k+1]$:
\[
 \sum_{k=m_n}^{\infty}(k+1)^{-2}
       \le\int_{m_n}^{\infty}x^{-2}\dd x=m_n^{-1}.
\]

On the retained intervals the coefficient functions satisfy $X^{\rm circ}=\me^{V_n}X^{\rm br}$, with $\operatorname{tr}V_{n,k}=0$ and the same boundary vectors $u_0,u_1$. They are
constant on the same partition, and $\eta_n\le1$ for large $n$. Lemma~\ref{S:lem:exact-mesh-kernel-stability} therefore gives

\[
 n\|F_n\|_{\mathfrak S_2}^2
 \le C\mathcal E_n^+(V_n)
 \le C_{\beta,\omega}\frac{(\log n)^{11}}{m_n}.
\]
Its use of $\mathcal E_n^+$ includes both triangles of each diagonal square. The third estimate of the same lemma also applies to these particular columns:
\[
 n\|F_n-\mathscr D_{n,X^{\rm br}}^{\rm step}V_n
                                      \|_{\mathfrak S_2}^2
 \le C\eta_n^2\mathcal E_n^+(V_n)
 \le C_{\beta,\omega}\eta_n^2
                              \frac{(\log n)^{11}}{m_n}.
\]
Here $X^{\rm br}$ denotes the matrix-valued coefficient function obtained by sampling the Brownian path on the retained intervals. Both estimates concern the same coefficient
function $X^{\rm br}$ and the same $V_n$. Taking square roots and using $m_n\ge(\log n)^{512}$ yields
\begin{equation}\label{S:eq:circular-brownian-retained-rate}
 \sqrt n\,\|F_n\|_{\mathfrak S_2}
 \le C_{\beta,\omega}\frac{(\log n)^{11/2}}{\sqrt{m_n}}
 \le C_{\beta,\omega}(\log n)^{-501/2},
\end{equation}
since $11/2-512/2=-501/2$.

The full error decomposes exactly as
\[
 F_n^{\mathrm{full}}
 =F_n+\mathsf T_{\beta,n}^{\mathrm{circ}}
              -\mathsf T_{\beta,n}^{\mathrm{br}}.
\]
The two summands on the right have disjoint kernel supports after grouping the two tails together. Corollary~\ref{S:cor:4.4} therefore gives
\[
 n\|F_n^{\mathrm{full}}\|_{\mathfrak S_2}^2 =n\|F_n\|_{\mathfrak S_2}^2 +n\|\mathsf T_{\beta,n}^{\mathrm{circ}} -\mathsf T_{\beta,n}^{\mathrm{br}}\|_{\mathfrak S_2}^2 \le n\|F_n\|_{\mathfrak S_2}^2+2\mathfrak t_n^{\mathrm{end}}.
\]
Theorem~\ref{S:thm:full-boundary-tails} gives, for $0<\nu<a_\beta$,
\begin{equation}\label{S:eq:circular-brownian-full-squared}
 n\|F_n^{\mathrm{full}}\|_{\mathfrak S_2}^2
 \le C_{\beta,\omega,\nu}
       \left(\frac{(\log n)^{11}}{m_n}+n^{-a_\beta+\nu}\right).
\end{equation}
Choose $\nu=a_\beta/2$ and take square roots. The additional term is at most $C_{\beta,\omega}n^{-a_\beta/4}$. Moreover,
\[
 \log\!\left(n^{-a_\beta/4}(\log n)^{501/2}\right)
 =-\frac{a_\beta}{4}\log n+\frac{501}{2}\log\log n
 \longrightarrow-\infty.
\]
Hence there is a deterministic $N_\beta$ such that $n^{-a_\beta/4}\le(\log n)^{-501/2}$ for $n\ge N_\beta$. Taking the square root of \eqref{S:eq:circular-brownian-full-squared}
and using $\sqrt{x+y}\le\sqrt x+\sqrt y$ gives
\[
 \sqrt n\,\|F_n^{\mathrm{full}}\|_{\mathfrak S_2}
 \le C_{\beta,\omega}\bigl((\log n)^{-501/2}+n^{-a_\beta/4}\bigr)
 \le2C_{\beta,\omega}(\log n)^{-501/2}
\]
for all sufficiently large $n$. This is \eqref{S:eq:circular-brownian-full-rate}. This bound holds pathwise on the original coupling and, by the preceding argument using $\pi$, on
every extension preserving its joint law.
\end{proof}

The full comparison also transfers an expansion of the Brownian sampling error. Suppose that on the retained interval this error has a decomposition
\[
 \varepsilon_{\beta,n}^{\mathrm{samp},T}=L_n+E_n,
 \qquad \sqrt n\,\|E_n\|_{\mathfrak S_2}
                    =O_\omega((\log n)^{-250}),
\]
where both terms are extended by zero outside the square of the retained interval. Then the exact identity
\[
 \Delta_{\beta,n}=L_n+
 \bigl(E_n+F_n+\mathsf T_{\beta,n}^{\mathrm{circ}}
                       -\mathsf T_{\beta,n}^{\mathrm{sine}}\bigr)
\]
gives an expansion on the full interval. By \eqref{S:eq:circular-brownian-retained-rate} and \eqref{S:eq:s8-tail-difference-norm}, with loss $a_\beta/4$, the parenthesized
remainder has normalized norm at most
\[
 \sqrt n\,\|E_n\|_{\mathfrak S_2}
   +C_{\beta,\omega}(\log n)^{-501/2}
   +C_{\beta,\omega}n^{-a_\beta/4}
 =O_\omega((\log n)^{-250}).
\]
Section~\ref{S:sec:the-error-from-sampling-the-brownian-path} constructs the Brownian expansion, and Section~\ref{S:sec:the-conditional-gaussian-operator} expresses its leading
term as a conditionally Gaussian random operator.

\subsection{The GOE and GUE in the comparison for Wigner matrices}
\label{sec:finite-gaussian-estimates-for-the-wigner-comparison}
\label{W:sec:gaussian-ensembles-and-airy-operator}

We prove the finite Gaussian estimates needed for entry comparison and state separately the qualitative operator limit for the Gaussian ensembles. The exact coefficient law gives
the process limit and the characteristic function estimates in Theorem~\ref{W:thm:coefficients-gaussian-direct}; the reflected process is treated jointly in
Theorem~\ref{W:thm:tightness-gaussian}. These estimates identify the local potential after the comparison of coefficient observables. Lemma~\ref{W:lem:gaussian-gram-conditioning}
controls orthogonalization at degree $O(N^{1/3})$. The lemma gives the tightness of the Gram matrix and its inverse needed for the comparison.

Theorem~\ref{W:thm:external-gaussian-sao} records the qualitative limits of the operator and of every fixed eigenpair with the present scaling and embedding. It gives the
qualitative localization statement \eqref{W:eq:gauss-sao-qualitative-localization}. The expectation comparison uses the stronger finite size bound in
Proposition~\ref{W:prop:gaussian-localization-expectation}, proved from the quadratic form bound for the Gaussian tridiagonal matrix and eigenvalue count in
Section~\ref{W:sec:coercivity-and-eigenvalue-counting}, independently of that operator limit. This expectation bound is uniform in the parameter of the cutoff in the eigenvalues of
the Gram matrix. Tightness of that matrix and its inverse separately controls the error when the cutoff is removed.

The comparison for Wigner matrices uses three estimates for the independent normal and chi variables. Their weighted sums determine the coefficient process; their polynomial
expansions control the Gram matrix as its dimension increases; and the quadratic form bound gives localization for the full finite matrix. The first calculation includes the two
adjacent off-diagonal coefficients, whose shared contribution determines the covariance.

The rate for the tridiagonal models in Section~\ref{sec:soft-edge-ensembles-and-quantitative-airy-limits} concerns a squared Hilbert--Schmidt coupling error. The Gaussian estimates
below use the qualitative operator limit and the explicit finite-dimensional coefficient distributions. The additional bounded entry ensemble $H_N^{\rm b}$ appears only in
Section~\ref{W:sec:polynomial-coefficients-for-comparison}, where we choose a common function $\chi_N$ for it and for the GOE or GUE. Uniformity in the starting vector means
uniform probabilities or expectations over deterministic unit vectors. This does not give a single event for all vectors chosen after the matrix is observed. Gaussian invariance
identifies the relevant finite-dimensional laws for each fixed vector; it does not identify the law for a random eigenvector used as the initial Lanczos vector. The finite size
probability estimates hold at each matrix size, and the qualitative operator theorem specifies the coupling between sizes.

\subsubsection{The Gaussian ensemble and its tridiagonal law}
\label{W:sec:goe-gue-tridiagonal-law}

Unitary or orthogonal invariance identifies the Lanczos coefficient law from any deterministic unit starting vector. The matrix normalization below fixes the variances of both the
diagonal entries and the independent chi coefficients.

The GOE or GUE matrix $H_N^{\mathrm G}$ has density
\begin{equation}
d\mathbb P({H_N^{\mathrm G}})
=
Z_{N,\beta}^{-1}
\exp\left(
-\frac{\beta N}{4}
\operatorname{Tr}{H_N^{\mathrm G}}^2
\right)d{H_N^{\mathrm G}}.
\label{W:eq:Gaussian-matrix-density}
\end{equation}

For \(\beta=1\), this is the GOE in the normalization of \eqref{W:eq:Gaussian-matrix-density}:
\begin{equation}
({H_N^{\mathrm G}})_{ij}
\sim
N(0,N^{-1}),
\qquad i<j,
\label{W:eq:GOE-offdiagonal}
\end{equation}
and
\begin{equation}
({H_N^{\mathrm G}})_{ii}
\sim
N(0,2N^{-1}).
\label{W:eq:GOE-diagonal}
\end{equation}

For \(\beta=2\), this is the GUE in the normalization of \eqref{W:eq:Gaussian-matrix-density}:
\begin{equation}
({H_N^{\mathrm G}})_{ij}
=
N^{-1/2}Z_{ij},
\qquad i<j,
\label{W:eq:GUE-offdiagonal}
\end{equation}
where
\[
\mathbb EZ_{ij}=0,
\qquad
\mathbb E|Z_{ij}|^2=1,
\qquad
\mathbb EZ_{ij}^2=0,
\]
and
\begin{equation}
({H_N^{\mathrm G}})_{ii}
\sim
N(0,N^{-1}).
\label{W:eq:GUE-diagonal}
\end{equation}

The normalization follows directly from
\[
 \operatorname{Tr}H^2=\sum_{i=1}^N h_{ii}^2+
                         2\sum_{i<j}|h_{ij}|^2.
\]
Thus the density factors over the independent upper triangular entries. In the complex case we may write $Z_{ij}=(X_{ij}+\ii Y_{ij})/\sqrt2$, where $X_{ij},Y_{ij}$ are independent
standard real Gaussian variables. The semicircle law has support $[-2,2]$ in this normalization.

Fix a deterministic unit vector \(w_N\). Choose a deterministic orthogonal or unitary matrix \(V_N\) with
\(V_Ne_1=w_N\). Gaussian invariance gives \(V_N^*{H_N^{\mathrm G}}V_N\stackrel{\mathrm d}={H_N^{\mathrm G}},\) so the
full Lanczos coefficient process from \(w_N\) has the same law as the one from \(e_1\). We write the recursion at this
starting vector to specify the indexing and the choice of positive off-diagonal coefficients: \(q_1=e_1.\) Set \(q_0=0,
\qquad b_{0,N}^G=0,\) and recursively define
\begin{equation}
a_{k,N}^G
=
\langle {H_N^{\mathrm G}}q_k,q_k\rangle,
\label{W:eq:Gaussian-Lanczos-a}
\end{equation}
\begin{equation}
r_k
=
{H_N^{\mathrm G}}q_k-a_{k,N}^Gq_k-b_{k-1,N}^Gq_{k-1},
\label{W:eq:Gaussian-Lanczos-residual}
\end{equation}
\begin{equation}
b_{k,N}^G
=
\|r_k\|,
\qquad
q_{k+1}
=
\frac{r_k}{b_{k,N}^G},
\label{W:eq:Gaussian-Lanczos-b}
\end{equation}
whenever \(b_{k,N}^G>0\).

The resulting tridiagonal matrix is
\begin{equation}
T_N^G
=
\begin{pmatrix}
a_{1,N}^G&b_{1,N}^G\\
b_{1,N}^G&a_{2,N}^G&b_{2,N}^G\\
&b_{2,N}^G&a_{3,N}^G&\ddots\\
&&\ddots&\ddots
\end{pmatrix}.
\label{W:eq:Gaussian-tridiagonal}
\end{equation}

The scaled tridiagonal matrix obtained from the GOE or GUE is
\begin{equation}
\mathcal A_N^G
=
m_N^2(2I-T_N^G)
=
N^{2/3}(2I-T_N^G).
\label{W:eq:Gaussian-edge-operator}
\end{equation}

The following is the tridiagonalization result of \cite{dumitriu2002matrix}, in the normalization \eqref{W:eq:Gaussian-matrix-density}.

\begin{theorem}[Dumitriu--Edelman]
\label{W:thm:Gaussian-Lanczos-law}
The Lanczos coefficients of $H_N^{\mathrm G}$ from any deterministic unit vector have the independent distributions
\begin{align}
 a_{k,N}^{G}&=\frac{g_k}{\sqrt{\beta N}},\qquad
        g_k\sim N(0,2),\quad 1\le k\le N,
        \label{W:eq:Gaussian-a-law}\\
 b_{k,N}^{G}&=\frac{\chi_{\beta(N-k)}}{\sqrt{\beta N}},
        \qquad 1\le k<N.\label{W:eq:Gaussian-b-law}
\end{align}
The recurrence reaches dimension $N$ almost surely.
\end{theorem}

\begin{proof}
First take the deterministic initial vector $e_1$. We work over $\mathbb F_\beta=\mathbb R$ for $\beta=1$ and $\mathbb F_\beta=\mathbb C$ for $\beta=2$. For $1\le m\le N$, let
$\mu_{m,N}$ be the probability law on $m\times m$ self-adjoint matrices with density proportional to
\[
 \exp\left(-\frac{\beta N}{4}\operatorname{Tr} C^2\right).
\]
The scale in this definition is $N$, also when $m<N$. This law is invariant under conjugation by every orthogonal or unitary $m\times m$ matrix: conjugation preserves the trace and
is an orthogonal transformation of the real vector space of self-adjoint matrices, and hence preserves its Lebesgue measure.

If $C\sim\mu_{m,N}$ and $m\ge2$, write
\[
 C=\begin{pmatrix} a&r^*\\ r&D\end{pmatrix},\qquad
 \operatorname{Tr}C^2=a^2+2r^*r+\operatorname{Tr}D^2.
\]
Lebesgue measure factors over these three coordinates. Thus $a,r,D$ are independent, $a\sim N(0,2/(\beta N))$, and $D\sim\mu_{m-1,N}$. The vector $r$ has density proportional to
$\exp(-\beta N|r|^2/2)$ on the real vector space of dimension $\beta(m-1)$. In the complex case the density is with respect to Lebesgue measure on its real and imaginary
coordinates. Consequently
\[
 b=|r|\ \stackrel{\mathrm d}=
       \frac{\chi_{\beta(m-1)}}{\sqrt{\beta N}},\qquad
 \mathbb P(b=0)=0.
\]
Choose a Borel orthogonal or unitary matrix $V(r)$ whose first column is $r/|r|$ when $r\ne0$. Such a choice is obtained by adjoining the coordinate vectors to $r/|r|$ in their
fixed order. Project each new vector onto the orthogonal complement of the vectors already selected, discard zero projections, and normalize every nonzero projection. The
conditions determining the retained indices are Borel; set $V(0)=I$ on the null event $r=0$. Then
\[
 \begin{pmatrix}1&0\\0&V(r)^*\end{pmatrix}
 C
 \begin{pmatrix}1&0\\0&V(r)\end{pmatrix}
 =\begin{pmatrix}a&b e_1^*\\b e_1&D'\end{pmatrix},
 \qquad D'=V(r)^*DV(r).
\]
Conditional on $(a,r)$, invariance and independence of $D$ give, for every bounded Borel function $f$,
\[
 \mathbb E[f(D')\mid a,r]
 =\int f(V(r)^*DV(r))\,\mu_{m-1,N}(\mathrm dD)
 =\int f(D)\,\mu_{m-1,N}(\mathrm dD).
\]
Thus $D'$ has law $\mu_{m-1,N}$ and is independent of $(a,r)$, in particular of $(a,b)$. The change of basis depends on the first column and leaves the conditional law of the
remaining block unchanged.

Apply this construction first with $m=N$. After $k-1$ steps, the first $k-1$ coordinates are tridiagonal, their only coupling to the remaining coordinates is through the first
coordinate of that block, and the remaining block has law $\mu_{N-k+1,N}$ independently of all coefficients already produced. A change of basis that is the identity on the first
$k$ coordinates preserves that coupling. The preceding calculation therefore supplies the next pair
\[
 a_k\sim N(0,2/(\beta N)),\qquad
 b_k\sim\frac{\chi_{\beta(N-k)}}{\sqrt{\beta N}},
\]
independently of every earlier pair, and leaves a block of law $\mu_{N-k,N}$. Induction gives these laws for $1\le k<N$ and the final independent diagonal variable
$a_N\sim N(0,2/(\beta N))$.

Each $b_k$ is strictly positive almost surely, and only finitely many such events are intersected at this fixed $N$. At step $k$, subtract from the current first column its
components in the preceding two basis vectors, then divide the resulting nonzero vector by its norm to obtain the new basis vector. Hence these are exactly the Lanczos vectors in
\eqref{W:eq:Gaussian-Lanczos-a}--\eqref{W:eq:Gaussian-Lanczos-b}, with positive off-diagonal coefficients. The constructed basis has $N$ vectors, so the recurrence does not break
down before the last coordinate. Writing $g_k=\sqrt{\beta N}\,a_k$ proves \eqref{W:eq:Gaussian-a-law} and \eqref{W:eq:Gaussian-b-law}, with mutual independence of the complete
family.

Finally fix any deterministic unit vector $w$, and choose a deterministic orthogonal or unitary $V$ with $Ve_1=w$. Applying the Lanczos recurrence to $V^*H_N^{\mathrm G}V$ from
$e_1$ gives the same coefficients as applying it to $H_N^{\mathrm G}$ from $w$: induction shows that multiplication by $V$ sends each $q_k$ and $r_k$ of the first recurrence to the
corresponding $q_k$ and $r_k$ of the second. The matrix $V^*H_N^{\mathrm G}V$ has the same Gaussian law. The coefficient distributions and absence of breakdown therefore hold for
every deterministic starting vector with exactly the same probabilities.
\end{proof}

\subsubsection{Moments and tails of chi random variables}
\label{W:sec:moments-and-tails-of-chi-variables}

The process comparison needs both an expansion of the chi mean and bounds uniform in its degree parameter. The same estimates will control deterministic drift, centered
fluctuations and the largest coefficient error.

For every real $\nu>0$, define $\chi_\nu$ by the density below. When $\nu$ is a positive integer, it is the Euclidean norm of a standard Gaussian vector in $\mathbb R^\nu$. The
density is
\begin{equation}
f_\nu(x)
=
\frac{2^{1-\nu/2}}{\Gamma(\nu/2)}
x^{\nu-1}e^{-x^2/2},
\qquad
x>0.
\label{W:eq:chi-density}
\end{equation}

Set \(\mu_\nu:=\mathbb E\chi_\nu.\)

\Needspace{6\baselineskip}

\begin{lemma}
\label{W:lem:chi-estimates}
As \(\nu\to\infty\),
\begin{equation}
\mu_\nu
=
\sqrt{\nu}
\left(
1-\frac1{4\nu}
+O(\nu^{-2})
\right),
\label{W:eq:chi-mean-expansion}
\end{equation}
and
\begin{equation}
\operatorname{Var}(\chi_\nu)
=
\frac12+O(\nu^{-1}).
\label{W:eq:chi-variance-expansion}
\end{equation}

For every \(\nu\ge1\) and every \(\lambda\in\mathbb R\),
\begin{equation}
\mathbb E
\exp\left(
\lambda(\chi_\nu-\mu_\nu)
\right)
\le
\exp\left(\frac{\lambda^2}{2}\right).
\label{W:eq:chi-mgf}
\end{equation}
Consequently,
\begin{equation}
\mathbb P
\left(
|\chi_\nu-\mu_\nu|>t
\right)
\le
2e^{-t^2/2},
\qquad
t\ge0,
\label{W:eq:chi-concentration}
\end{equation}
and, for every fixed \(0<p<\infty\),
\begin{equation}
\sup_{\nu\ge1}
\mathbb E
|\chi_\nu-\mu_\nu|^p
<\infty.
\label{W:eq:chi-uniform-moments}
\end{equation}
\end{lemma}

\begin{proof}
For $t>-\nu$, the change of variables $y=x^2/2$ in \eqref{W:eq:chi-density} gives
\[
 \mathbb E\chi_\nu^t
 =\frac{2^{t/2}\Gamma((\nu+t)/2)}{\Gamma(\nu/2)},
 \qquad \mathbb E\chi_\nu^2=\nu.
\]
In particular,
\begin{equation}
 \mu_\nu=\sqrt2\,\frac{\Gamma((\nu+1)/2)}{\Gamma(\nu/2)}.
 \label{W:eq:chi-exact-mean}
\end{equation}
The Gamma ratio expansion in \cite{nistGammaAsymptotics}, equations~(5.11.13) and~(5.11.15), applies with the fixed parameters $a=1/2$, $b=0$ and positive real argument
$z=\nu/2\to\infty$. Its first two coefficients are $1$ and $G_1(1/2,0)=-1/8$. Thus
\[
 \mu_\nu=\sqrt\nu\bigl(1-(4\nu)^{-1}+O(\nu^{-2})\bigr),
 \qquad \nu-\mu_\nu^2=\tfrac12+O(\nu^{-1}),
\]
which proves \eqref{W:eq:chi-mean-expansion}--\eqref{W:eq:chi-variance-expansion}.

For the concentration bound, fix a real $\nu\ge1$ and consider the tilted density
\[
 p_\lambda(x)=\frac{x^{\nu-1}e^{-x^2/2+\lambda x}}{Z(\lambda)},
 \qquad Z(\lambda)=\int_0^\infty x^{\nu-1}e^{-x^2/2+\lambda x}\,\dd x,
 \qquad x>0.
\]
Write $\mathbb E_\lambda$ for integration against this density and $\mu_\lambda=\mathbb E_\lambda X$. On every compact interval of $\lambda$, the integrands obtained by any fixed
number of differentiations of $Z$ are bounded at infinity by $C x^{\nu+r-1}e^{-x^2/4}$ and are integrable at zero. Differentiation under the integral therefore gives
\[
 (\log Z)'(\lambda)=\mu_\lambda,
 \qquad (\log Z)''(\lambda)=\operatorname{Var}_\lambda X.
\]
For $V_\lambda(x)=x^2/2-(\nu-1)\log x-\lambda x$, one has $p_\lambda'=-V_\lambda'p_\lambda$ and $V_\lambda''(x)=1+(\nu-1)/x^2\ge1$. If $\nu>1$, integration by parts on
$[\delta,R]$, followed by $\delta\downarrow0$ and $R\to\infty$, gives
\[
 \mathbb E_\lambda[(X-\mu_\lambda)V_\lambda'(X)]=1.
\]
Indeed $(x-\mu_\lambda)p_\lambda(x)$ tends to zero at both ends, while the only possibly singular integrand at zero is bounded by $C x^{\nu-2}$, which is integrable for $\nu>1$. If
$\nu=1$, the density has a finite value at zero and the same calculation instead gives
\[
 \mathbb E_\lambda[(X-\mu_\lambda)V_\lambda'(X)]
 =1-\mu_\lambda p_\lambda(0)\le1.
\]
The pointwise inequality
\[
 (x-\mu_\lambda)\bigl(V_\lambda'(x)-V_\lambda'(\mu_\lambda)\bigr)
 \ge (x-\mu_\lambda)^2
\]
follows by integrating $V_\lambda''\ge1$ between $x$ and $\mu_\lambda>0$. Since $\mathbb E_\lambda(X-\mu_\lambda)=0$, both cases imply $\operatorname{Var}_\lambda X\le1$. For
$F(\lambda)=\log Z(\lambda)-\log Z(0)-\lambda\mu_\nu$ we consequently have $F(0)=F'(0)=0$ and $F''\le1$. The identity
\[
 F(\lambda)=\lambda^2\int_0^1(1-s)F''(s\lambda)\,\dd s\le\lambda^2/2
\]
holds for either sign of $\lambda$ and proves \eqref{W:eq:chi-mgf}. In particular, $\operatorname{Var}(\chi_\nu)\le1$, so
\[
 \mu_\nu^2=\nu-\operatorname{Var}(\chi_\nu)\ge\nu-1.
\]
Exponential Markov inequalities with $\lambda=t$ and $\lambda=-t$ give the two one-sided bounds $e^{-t^2/2}$ and hence \eqref{W:eq:chi-concentration}. For every $p>0$, Tonelli's
theorem applied to the nonnegative tail integral gives
\[
 \mathbb E|\chi_\nu-\mu_\nu|^p
 \le2p\int_0^\infty t^{p-1}e^{-t^2/2}\,\dd t
 =p\,2^{p/2}\Gamma(p/2).
\]
This is uniform in $\nu\ge1$ and proves \eqref{W:eq:chi-uniform-moments}.
\end{proof}

\subsubsection{Convergence of the coefficient process for the GOE and GUE}
\label{W:sec:convergence-of-gaussian-coefficient-process}

For the GOE and GUE, the diagonal and off-diagonal fluctuations in the tridiagonal representation are independent. We first sum them in the potential $Y_N^G$, whose successive
increments are independent. The difference between $\widetilde Y_N^G$ and $Y_N^G$ is a single boundary coefficient multiplied by $m_N$. The scaled maximum estimate below shows that
this difference tends to zero locally.

Fix $T>0$. For estimates on $[0,T]$, take $N\ge N_T$ so that $\lfloor m_NT\rfloor+4<N$. The sampled values below are defined for $0\le k<N$, and the indicated extra nodes cover
every interpolation on this interval. Set $Y_{N,0}^G=0$ and
\begin{equation}
Y_{N,k}^G
=
\sum_{r=1}^{k}
\left[
-m_Na_{r,N}^G
+
2m_N(1-b_{r,N}^G)
\right].
\label{W:eq:Gaussian-grid-potential}
\end{equation}

Let \(t_{k,N}:=\frac{k}{m_N}.\) The linearly interpolated process \(Y_N^G\) is defined by \(Y_N^G(t_{k,N})=Y_{N,k}^G\) and, for \(t=\frac{k+\theta}{m_N}, \qquad 0\le\theta<1,\) by
\begin{equation}
Y_N^G(t)
=
(1-\theta)Y_{N,k}^G
+
\theta Y_{N,k+1}^G.
\label{W:eq:Gaussian-linear-interpolation}
\end{equation}

We also use the step process
\begin{equation}
Y_N^{G,\mathrm{st}}(t)
=
Y_{N,\lfloor m_Nt\rfloor}^G.
\label{W:eq:Gaussian-step-potential}
\end{equation}

Decompose
\begin{equation}
Y_N^G(t)
=
D_N(t)+M_N^a(t)+M_N^b(t),
\label{W:eq:Gaussian-potential-decomposition}
\end{equation}
where the three processes are linearly interpolated from the sampled values
\begin{equation}
D_{N,k}
=
2m_N
\sum_{r=1}^{k}
(1-\mathbb Eb_{r,N}^G),
\label{W:eq:Gaussian-drift-grid}
\end{equation}
\begin{equation}
M_{N,k}^a
=
-m_N
\sum_{r=1}^{k}
a_{r,N}^G,
\label{W:eq:Gaussian-a-martingale}
\end{equation}
and
\begin{equation}
M_{N,k}^b
=
-2m_N
\sum_{r=1}^{k}
\left(
b_{r,N}^G-\mathbb Eb_{r,N}^G
\right).
\label{W:eq:Gaussian-b-martingale}
\end{equation}

The processes \(M_N^a\) and \(M_N^b\) are independent.

The decomposition \eqref{W:eq:Gaussian-potential-decomposition} is the coefficient decomposition used in \cite{ramirez2011beta}. We state the joint limit with the present
normalization.

\begin{theorem}[Convergence of the coefficient process for the GOE and GUE]
\label{W:thm:Gaussian-potential-limit}
For every fixed $T<\infty$, the deterministic drift satisfies
\begin{equation}
 \sup_{0\le t\le T}\left|D_N(t)-\frac{t^2}{2}\right|
 =O_{\beta,T}(m_N^{-1}),
 \label{W:eq:Gaussian-drift-convergence}
\end{equation}
and the centered processes satisfy
\begin{equation}
 (M_N^a,M_N^b)\Longrightarrow
 \sqrt{\frac2\beta}\,(B_1,B_2)
 \quad\text{in }C([0,T])^2,
 \label{W:eq:Gaussian-joint-martingale-limit}
\end{equation}
where $B_1,B_2$ are independent standard Brownian motions. Consequently,
\begin{equation}
 Y_N^G\Longrightarrow Y_\beta
 \quad\text{in }C([0,T]),
 \label{W:eq:Gaussian-potential-limit}
\end{equation}
with
\begin{equation}
 Y_\beta(t)=\frac{t^2}{2}+\frac2{\sqrt\beta}B(t).
 \label{W:eq:Airy-integrated-potential}
\end{equation}
The drift and the two centered processes converge jointly.
\end{theorem}

\begin{proof}
We compute the drift and prove finite-dimensional convergence and tightness directly for the two independent arrays. Uniformly for $k\le\lfloor m_NT\rfloor+1$,
Lemma~\ref{W:lem:chi-estimates} gives
\begin{equation}
 1-\mathbb Eb_{k,N}^G
 =\frac{k}{2N}+\frac1{4\beta N}
  +O_{\beta,T}\left(\frac{k^2+k+1}{N^2}\right)
 =\frac{k}{2N}+O_{\beta,T}(N^{-1}).
 \label{W:eq:Gaussian-b-mean-coarse}
\end{equation}
Thus, for $K=\lfloor m_Nt\rfloor$,
\[
 D_{N,K}=\frac{m_NK(K+1)}{2N}
          +O_{\beta,T}(m_NK/N)
        =\frac{(K/m_N)^2}{2}+O_{\beta,T}(m_N^{-1}).
\]
Linear interpolation changes the quadratic comparison by at most $1/(8m_N^2)$, proving \eqref{W:eq:Gaussian-drift-convergence}.

Set $\zeta_{k,N}^a=-m_Na_{k,N}^G$ and $\zeta_{k,N}^b=-2m_N(b_{k,N}^G-\mathbb Eb_{k,N}^G)$. Their means vanish, and their variance and fourth moment bounds are
\[
 \begin{aligned}
 m_N\mathbb E(\zeta_{k,N}^a)^2&=2/\beta,&
 m_N\mathbb E(\zeta_{k,N}^b)^2&=2/\beta+O_{\beta,T}(N^{-1}),\\
 \sup_{k\le m_NT+1}\mathbb E|\zeta_{k,N}^a|^4
 +\sup_{k\le m_NT+1}\mathbb E|\zeta_{k,N}^b|^4
 &\le C_{\beta,T}m_N^{-2}.&&
 \end{aligned}
\]

Fix $0\le t_1,\ldots,t_J\le T$ and real numbers $u_{j,a},u_{j,b}$. The corresponding linear combination of the centered step processes is
\[
 S_N=\sum_{k\le m_NT}
       (c_{k,a}\zeta_{k,N}^a+c_{k,b}\zeta_{k,N}^b),\qquad
 c_{k,\tau}=\sum_{j=1}^J u_{j,\tau}
                         \mathbf1_{\{k\le\lfloor m_Nt_j\rfloor\}}.
\]
The coefficients are bounded by $C_u=\sum_{j,\tau}|u_{j,\tau}|$. Independence and the variance bounds give
\[
 \mathbb ES_N^2 =\sum_{\tau=a,b}\sum_{i,j=1}^J u_{i,\tau}u_{j,\tau} \sum_{k\le m_N(t_i\wedge t_j)}\mathbb E(\zeta_{k,N}^\tau)^2 \longrightarrow\frac2\beta\sum_{\tau=a,b} \sum_{i,j=1}^J u_{i,\tau}u_{j,\tau}(t_i\wedge t_j).
\]
For centered independent summands $X_l$, Taylor's formula and telescoping the products of their characteristic functions give
\[
 \left|\mathbb Ee^{\ii\sum_lX_l}
        -e^{-\frac12\sum_l\mathbb EX_l^2}\right|
 \le\sum_l\left|\mathbb Ee^{\ii X_l}-e^{-\mathbb EX_l^2/2}\right|
 \le\frac16\sum_l\mathbb E|X_l|^3
       +\frac18\sum_l(\mathbb EX_l^2)^2.
\]
For the summands of $S_N$,
\[
 \sum_l\mathbb E|X_l|^3
       \le C_u\sum_{k,\tau}(\mathbb E|\zeta_{k,N}^\tau|^4)^{3/4}
       \le C_{u,T}m_N^{-1/2},\qquad
 \sum_l(\mathbb EX_l^2)^2\le C_{u,T}m_N^{-1}.
\]
Consequently,
\[
 |\mathbb Ee^{\ii S_N}-e^{-\mathbb ES_N^2/2}|
       \le C_{u,T}(m_N^{-1/2}+m_N^{-1})\longrightarrow0.
\]
The limiting covariance computed above is that of $\sqrt{2/\beta}(B_1,B_2)$. This proves joint finite-dimensional convergence.

For tightness of the linear interpolants, set $h=m_N^{-1}$. An increment from $s$ to $t$ is $\sum_k\alpha_k\zeta_{k,N}$, where $0\le\alpha_k\le1$ and $\sum_k\alpha_k=(t-s)/h$.
Independence and centering give
\begin{align*}
 \mathbb E\left|\sum_k\alpha_k\zeta_{k,N}\right|^4
 &=\sum_k\alpha_k^4\mathbb E\zeta_{k,N}^4
    +6\sum_{k<l}\alpha_k^2\alpha_l^2
           \mathbb E\zeta_{k,N}^2\mathbb E\zeta_{l,N}^2\\
 &\le Ch^2\left(\sum_k\alpha_k^2\right)^2
 \le Ch^2\min\left\{\frac{t-s}{h},\frac{(t-s)^2}{h^2}\right\}^2
 \le C|t-s|^2.
\end{align*}
Lemma~\ref{W:lem:tightness-dyadic}, with $p=4$ and $\alpha=1$, now gives
\[
 \sup_N\mathbb P\bigl(\omega_T(M_N^\tau,\delta)>u\bigr)
       \le C_{\beta,T}u^{-4}\delta,\qquad \tau=a,b.
\]
Since both processes vanish at zero, the pair is tight in $C([0,T])^2$. Its joint finite-dimensional limits identify the unique limit as $\sqrt{2/\beta}(B_1,B_2)$ with independent
Brownian motions. For the step version, the interpolation error is bounded by the largest increment and
\[
 \mathbb P\left(\max_{k\le m_NT+1}|\zeta_{k,N}|>\varepsilon\right)
 \le C\varepsilon^{-4}m_N^{-1}\longrightarrow0.
\]
Together with $\|D_N-t^2/2\|_\infty\to0$, the joint process convergence and the continuous mapping theorem give
\[
 (D_N,M_N^a,M_N^b)
 \Rightarrow(t^2/2,\sqrt{2/\beta}B_1,\sqrt{2/\beta}B_2).
\]
Applying addition to \eqref{W:eq:Gaussian-potential-decomposition} gives $t^2/2+2B/\sqrt\beta$, where $B=(B_1+B_2)/\sqrt2$ is standard Brownian motion.
\end{proof}

\subsubsection{Weighted sums, increments, and maxima}
\label{W:sec:weighted-gaussian-coefficient-statistics}

The comparison of characteristic functions requires uniform estimates for deterministic weighted sums. We derive these estimates directly from the independent coefficient law,
keeping the weight norms explicit. The increment and maximum bounds at the end of this subsubsection are consequences of the same estimates.

Fix \(T<\infty\), and set \(K_N(T) = \lfloor m_NT\rfloor.\) For vectors
\[
 c=(c_1,\ldots,c_{K_N(T)}),\qquad d=(d_1,\ldots,d_{K_N(T)}),
\]
define
\begin{equation}
Z_N^G(c,d)
=
-m_N\sum_{k=1}^{K_N(T)}c_ka_{k,N}^G
+
2m_N\sum_{k=1}^{K_N(T)}
d_k(1-b_{k,N}^G).
\label{W:eq:Gaussian-weighted-statistic}
\end{equation}

Introduce the normalized discrete norms
\begin{equation}
\|c\|_{p,m_N}
=
\left(
\frac1{m_N}
\sum_{k=1}^{K_N(T)}
|c_k|^p
\right)^{1/p},
\qquad
1\le p<\infty,
\label{W:eq:normalized-lp-norm}
\end{equation}
and \(\|c\|_{\infty,m_N} = \max_{k\le K_N(T)}|c_k|.\) We also write
\begin{equation}
\langle c,d\rangle_{m_N}
=
\frac1{m_N}
\sum_{k=1}^{K_N(T)}c_kd_k.
\label{W:eq:normalized-discrete-inner-product}
\end{equation}

Let
\begin{equation}
\overline Z_N^G(c,d)
=
Z_N^G(c,d)
-
\mathbb EZ_N^G(c,d).
\label{W:eq:centered-Gaussian-weighted-statistic}
\end{equation}

\Needspace{6\baselineskip}

\begin{proposition}
\label{W:prop:weighted-Gaussian-estimates}
For fixed $T$ and all sufficiently large $N$, the following bounds hold for every pair of real coefficient vectors:
\begin{align}
 |\mathbb EZ_N^G(c,d)|&\le C_{\beta,T}\|d\|_{1,m_N},
       \label{W:eq:Gaussian-weighted-mean}\\
 \mathbb Ee^{\lambda\overline Z_N^G(c,d)}
 &\le\exp\!\left[C_\beta\lambda^2
           (\|c\|_{2,m_N}^2+\|d\|_{2,m_N}^2)\right],
       \quad \lambda\in\mathbb R,
       \label{W:eq:Gaussian-weighted-mgf}\\
 \|\overline Z_N^G(c,d)\|_{L^p}
 &\le C_\beta\sqrt p(\|c\|_{2,m_N}+\|d\|_{2,m_N}),\quad p\ge1.
       \label{W:eq:Gaussian-weighted-Lp}
\end{align}
\end{proposition}

\begin{proof}
Since \(\mathbb Ea_{k,N}^G=0,\) we have
\[
\mathbb EZ_N^G(c,d)
=
2m_N
\sum_{k=1}^{K_N(T)}
d_k(1-\mathbb Eb_{k,N}^G).
\]
By \eqref{W:eq:Gaussian-b-mean-coarse},
\[
|1-\mathbb Eb_{k,N}^G|
\le
C_T\frac{k+1}{N}
\le
C_T\frac{m_N}{N}
=
\frac{C_T}{m_N^2}.
\]
Therefore
\[
 |\mathbb EZ_N^G(c,d)| \le \frac{C_T}{m_N} \sum_{k=1}^{K_N(T)}|d_k| = C_T\|d\|_{1,m_N}.
\]
This proves \eqref{W:eq:Gaussian-weighted-mean}.

For the centered statistic, write
\begin{equation}
\overline Z_N^G(c,d)
=
\sum_{k=1}^{K_N(T)}
\left[
c_kA_{k,N}+d_kB_{k,N}
\right],
\label{W:eq:Gaussian-weighted-centered-decomposition}
\end{equation}
where \(A_{k,N} = -m_Na_{k,N}^G\) and \(B_{k,N} = -2m_N \left( b_{k,N}^G-\mathbb Eb_{k,N}^G \right).\) The complete family \((A_{k,N},B_{k,N})_{k}\) is independent.

Since \(A_{k,N}\) is Gaussian with variance \(\frac{2}{\beta m_N},\) we have
\begin{equation}
\mathbb E e^{\lambda c_kA_{k,N}}
=
\exp\left(
\frac{\lambda^2c_k^2}{\beta m_N}
\right).
\label{W:eq:Gaussian-a-mgf-weighted}
\end{equation}

By \eqref{W:eq:chi-mgf},
\begin{equation}
\mathbb E e^{\lambda d_kB_{k,N}}
\le
\exp\left(
\frac{2\lambda^2d_k^2}{\beta m_N}
\right).
\label{W:eq:Gaussian-b-mgf-weighted}
\end{equation}
Multiplying the moment generating functions over \(k\) gives
\[
\mathbb E
e^{\lambda\overline Z_N^G(c,d)}
\le
\exp\left[
\frac{C_\beta\lambda^2}{m_N}
\sum_k
(c_k^2+d_k^2)
\right],
\]
which is \eqref{W:eq:Gaussian-weighted-mgf}.

If $\mathbb Ee^{\lambda X}\le e^{v^2\lambda^2/2}$, Chernoff's bound gives $\mathbb P(|X|>t)\le2e^{-t^2/(2v^2)}$. Therefore
\[
\mathbb E|X|^p
=p\int_0^\infty t^{p-1}\mathbb P(|X|>t)\,dt
\le p(2v^2)^{p/2}\Gamma(p/2),
\]
and Stirling's bound for $\Gamma(p/2)$ gives $\|X\|_{L^p}\le C\sqrt p\,v$. Applying this to \eqref{W:eq:Gaussian-weighted-mgf} proves \eqref{W:eq:Gaussian-weighted-Lp}.
\end{proof}

For several coefficient pairs write $Z_j=Z_N^G(c^{(j)},d^{(j)})$ and $\overline Z_j=Z_j-\mathbb EZ_j$. We denote their joint cumulant by
$\kappa_r(\overline Z_1,\ldots,\overline Z_r)$.

\begin{lemma}
\label{W:lem:gaussian-weighted-cumulants}
For fixed $T$, the covariance is
\begin{align}
 \operatorname{Cov}(Z_1,Z_2)
 &=\frac2\beta\bigl(\langle c^{(1)},c^{(2)}\rangle_{m_N}
             +\langle d^{(1)},d^{(2)}\rangle_{m_N}\bigr)+\mathcal R_{N,2},
             \label{W:eq:Gaussian-weighted-covariance}\\
 |\mathcal R_{N,2}|&\le C_{\beta,T}N^{-1}
                  \|d^{(1)}\|_{2,m_N}\|d^{(2)}\|_{2,m_N}.
             \label{W:eq:Gaussian-covariance-error}
\end{align}
For every fixed $r\ge3$,
\begin{equation}
 |\kappa_r(\overline Z_1,\ldots,\overline Z_r)|
 \le C_{r,\beta}m_N^{-r/2}\sum_{k\le K_N(T)}
                          \prod_{j=1}^r|d_k^{(j)}|.
 \label{W:eq:Gaussian-higher-joint-cumulant}
\end{equation}
\end{lemma}
\begin{proof}
For the covariance, independence of the \(a\)- and \(b\)-families gives
\begin{equation}
 \operatorname{Cov} ={} m_N^2 \frac{2}{\beta N} \sum_k c_k^{(1)}c_k^{(2)} + \frac{4m_N^2}{\beta N} \sum_k d_k^{(1)}d_k^{(2)} \operatorname{Var} \left( \chi_{\beta(N-k)} \right).
 \label{W:eq:Gaussian-covariance-exact}
\end{equation}
Since \(\frac{m_N^2}{N} = \frac1{m_N},\) the first term is \(\frac2\beta \langle c^{(1)},c^{(2)}\rangle_{m_N}.\) Using
\(\operatorname{Var}(\chi_{\beta(N-k)}) = \frac12+O(N^{-1})\) uniformly for \(k\le m_NT\), the second term equals
\(\frac2\beta \langle d^{(1)},d^{(2)}\rangle_{m_N} + \mathcal R_{N,2}.\) By Cauchy--Schwarz,
\[
|\mathcal R_{N,2}|
\le
C_{\beta}
\frac{m_N^2}{N^2}
\sum_k
|d_k^{(1)}d_k^{(2)}|
\le
C_{\beta,T}N^{-1}
\|d^{(1)}\|_{2,m_N}
\|d^{(2)}\|_{2,m_N}.
\]
This proves \eqref{W:eq:Gaussian-weighted-covariance}.

Finally, fix \(r\ge3\). The \(a\)-family is Gaussian, so every cumulant of order at least three involving only the
\(a\)-part vanishes. Mixed cumulants between the independent \(a\)- and \(b\)-families also vanish.

For independent variables, a joint cumulant vanishes unless all arguments arise from the same independent block. Hence
\begin{equation}
 \kappa_r \left( \overline Z_N^G(c^{(1)},d^{(1)}), \ldots, \overline Z_N^G(c^{(r)},d^{(r)}) \right) = \left( -\frac{2m_N}{\sqrt{\beta N}} \right)^r \sum_{k=1}^{K_N(T)} \left( \prod_{j=1}^{r}d_k^{(j)} \right) \kappa_r \left( \chi_{\beta(N-k)} \right).
 \label{W:eq:Gaussian-joint-cumulant-exact}
\end{equation}
For $X_\nu=\chi_\nu-\mathbb E\chi_\nu$, the finite partition formula gives
\[
 \kappa_r(X_\nu)=\sum_{\pi\in\Pi([r])}
   (-1)^{|\pi|-1}(|\pi|-1)!\prod_{B\in\pi}\mathbb EX_\nu^{|B|}.
\]
By H\"older's inequality, $\prod_{B\in\pi}|\mathbb EX_\nu^{|B|}|\le\mathbb E|X_\nu|^r$. Consequently
\[
 |\kappa_r(\chi_\nu)|=|\kappa_r(X_\nu)|
 \le\left(\sum_{\pi\in\Pi([r])}(|\pi|-1)!\right)
             \mathbb E|X_\nu|^r\le C_r
\]
by \eqref{W:eq:chi-uniform-moments}. Since \(\frac{m_N}{\sqrt N} = m_N^{-1/2},\) equation \eqref{W:eq:Gaussian-joint-cumulant-exact} gives
\eqref{W:eq:Gaussian-higher-joint-cumulant}. Equivalently, dividing the sum by $m_N$ gives
\begin{equation}
 |\kappa_r(\overline Z_1,\ldots,\overline Z_r)|
 \le C_{r,\beta}m_N^{1-r/2}\frac1{m_N}
                \sum_{k\le K_N(T)}\prod_{j=1}^r|d_k^{(j)}|.
 \label{W:eq:Gaussian-higher-cumulant-normalized}
\end{equation}
\end{proof}

\subparagraph{Increments and maxima}
\label{W:sec:gaussian-increments-and-coefficient-convergence}

We retain the weights and statistics of Section~\ref{W:sec:weighted-gaussian-coefficient-statistics}. Let \(I=(p,q] \subset \{1,\ldots,K_N(T)\}, \qquad |I|=q-p,\) and define
\begin{equation}
\Delta_IY_N^G
=
Y_{N,q}^G-Y_{N,p}^G.
\label{W:eq:Gaussian-interval-increment}
\end{equation}

This is the weighted statistic \(\Delta_IY_N^G = Z_N^G(\mathbf 1_I,\mathbf 1_I).\)

\Needspace{6\baselineskip}

\begin{corollary}
\label{W:cor:Gaussian-interval-moments}
For every fixed \(T<\infty\) and \(p_0\ge2\), there exists \(C_{\beta,T,p_0}<\infty\) such that
\begin{equation}
\mathbb E
|\Delta_IY_N^G|^{p_0}
\le
C_{\beta,T,p_0}
\left(
\frac{|I|+1}{m_N}
\right)^{p_0/2}
\label{W:eq:Gaussian-interval-p-moment}
\end{equation}
uniformly over all intervals \(I\).

In particular,
\begin{equation}
\mathbb E
|\Delta_IY_N^G|^4
\le
C_{\beta,T}
\left(
\frac{|I|+1}{m_N}
\right)^2.
\label{W:eq:Gaussian-interval-fourth-moment}
\end{equation}
\end{corollary}

\begin{proof}
For the interval indicator, \(\|\mathbf 1_I\|_{2,m_N} = \left( \frac{|I|}{m_N} \right)^{1/2},\) and \(\|\mathbf 1_I\|_{1,m_N} = \frac{|I|}{m_N}.\)
Proposition~\ref{W:prop:weighted-Gaussian-estimates} gives
\[
\left\|
\Delta_IY_N^G-\mathbb E\Delta_IY_N^G
\right\|_{L^{p_0}}
\le
C_{\beta}\sqrt{p_0}
\left(
\frac{|I|}{m_N}
\right)^{1/2},
\]
and \(|\mathbb E\Delta_IY_N^G| \le C_T\frac{|I|}{m_N}.\) Since \(|I|\le m_NT,\) the deterministic term satisfies
\[
\left(
\frac{|I|}{m_N}
\right)^{p_0}
\le
C_T
\left(
\frac{|I|+1}{m_N}
\right)^{p_0/2}.
\]
The triangle inequality in \(L^{p_0}\) proves \eqref{W:eq:Gaussian-interval-p-moment}.
\end{proof}

For $0\le s<t\le T$, define
\[
 \alpha_k=h_N^{-1}|(s,t]\cap((k-1)h_N,kh_N]|.
\]
Then $0\le\alpha_k\le1$, $\sum_k\alpha_k=(t-s)/h_N$, and $Y_N^G(t)-Y_N^G(s)=Z_N^G(\alpha,\alpha)$, with the last possible index $\lceil m_NT\rceil$. The weighted estimates apply on
$[0,T+1]$. They give
\begin{align*}
 \|\overline Z_N^G(\alpha,\alpha)\|_{L^4}^4
 &\le C_\beta\left(h_N\sum_k\alpha_k^2\right)^2
       \le C_\beta|t-s|^2,\\
 |\mathbb EZ_N^G(\alpha,\alpha)|^4
 &\le C_{\beta,T}\left(h_N\sum_k\alpha_k\right)^4
       =C_{\beta,T}|t-s|^4.
\end{align*}
In particular, since $|t-s|\le T$,
\begin{equation}
 \mathbb E|Y_N^G(t)-Y_N^G(s)|^4
 \le C_{\beta,T}|t-s|^2
 \le C_{\beta,T}(|t-s|+m_N^{-1})^2.
 \label{W:eq:Gaussian-continuous-fourth-moment}
\end{equation}

The next estimate includes the last node used when interpolating on $[0,T]$. All coefficient bounds in its proof may be taken on $[0,T+1]$.

\begin{proposition}
\label{W:prop:Gaussian-mesh-jumps}
For every fixed \(T<\infty\) and \(0<p<\infty\),
\begin{equation}
\left\|
\max_{k\le K_N(T)+1}
\left|
Y_{N,k}^G-Y_{N,k-1}^G
\right|
\right\|_{L^p}
\longrightarrow0.
\label{W:eq:Gaussian-mesh-jumps}
\end{equation}
More precisely,
\begin{equation}
\left\|
\max_{k\le K_N(T)+1}
\left|
Y_{N,k}^G-Y_{N,k-1}^G
\right|
\right\|_{L^p}
\le
C_{\beta,T,p}
\sqrt{
\frac{\log(2+m_N)}{m_N}
}
+
\frac{C_T}{m_N}.
\label{W:eq:Gaussian-mesh-jump-rate}
\end{equation}
Consequently,
\begin{equation}
\sup_{0\le t\le T}
|Y_N^G(t)-Y_N^{G,\mathrm{st}}(t)|
\longrightarrow0
\label{W:eq:Gaussian-step-linear-equivalence}
\end{equation}
in every fixed \(L^p\).
\end{proposition}

\begin{proof}
The one step increment is \(Y_{N,k}^G-Y_{N,k-1}^G = A_{k,N}+B_{k,N} + 2m_N(1-\mathbb Eb_{k,N}^G),\) where \(A_{k,N}\) and \(B_{k,N}\) are defined in
\eqref{W:eq:Gaussian-weighted-centered-decomposition}.

Both centered variables are sub-Gaussian with scale \(Cm_N^{-1/2}.\) Indeed, this follows from \eqref{W:eq:Gaussian-a-mgf-weighted} and \eqref{W:eq:Gaussian-b-mgf-weighted} with
\(c_k=d_k=1\).

For $M$ centered variables with sub-Gaussian scale at most $\sigma$, write $Z=\max_j|X_j|$ and $a=\sigma\sqrt{2\log(2M)}$. The union bound implies, for $u\ge0$,
\[
 \mathbb P(Z>a+u)\le 2M e^{-(a+u)^2/(2\sigma^2)}
                      \le e^{-u^2/(2\sigma^2)}.
\]
Tail integration gives $\|(Z-a)_+\|_{L^p}\le C\sigma\sqrt p$. The triangle inequality (and monotonicity from $L^1$ when $p<1$) therefore gives
\begin{equation}
\left\|
\max_{1\le j\le M}|X_j|
\right\|_{L^p}
\le
C\sigma
\sqrt{p+\log(2M)}.
\label{W:eq:subGaussian-maximal-bound}
\end{equation}
Applying this with \(M=K_N(T)+1=O_T(m_N)\) and \(\sigma=Cm_N^{-1/2}\) gives
\[
\left\|
\max_{k\le K_N(T)+1}
(|A_{k,N}|+|B_{k,N}|)
\right\|_{L^p}
\le
C_{\beta,T,p}
\sqrt{
\frac{\log(2+m_N)}{m_N}
}.
\]

By \eqref{W:eq:Gaussian-b-mean-coarse}, \(\max_{k\le K_N(T)+1} 2m_N|1-\mathbb Eb_{k,N}^G| \le \frac{C_T}{m_N}.\) This proves \eqref{W:eq:Gaussian-mesh-jump-rate}.

On every partition interval, the difference between the step path and its linear interpolation is bounded by the corresponding increment between sampling times. Hence
\eqref{W:eq:Gaussian-step-linear-equivalence} follows.
\end{proof}

To pass from $Y_N^G$ to the exact quadratic form potential, we need the scaled maximum in Proposition~\ref{W:prop:Gaussian-offdiagonal-coefficients}. The criterion for convergence
of the embedded resolvents requires only unscaled convergence of the off-diagonal coefficients. The quantitative Fourier estimate is proved separately in
Theorem~\ref{W:thm:coefficients-gaussian-direct}.

\Needspace{6\baselineskip}

\begin{proposition}
\label{W:prop:Gaussian-offdiagonal-coefficients}
For every fixed \(T<\infty\) and \(0<p<\infty\),
\begin{equation}
\left\|
m_N
\max_{1\le k\le K_N(T)+1}
|b_{k,N}^G-1|
\right\|_{L^p}
\longrightarrow0.
\label{W:eq:Gaussian-offdiagonal-coefficients-Lp}
\end{equation}
More precisely,
\begin{equation}
\left\|
m_N
\max_{1\le k\le K_N(T)+1}
|b_{k,N}^G-1|
\right\|_{L^p}
\le
C_{\beta,p}
\sqrt{
\frac{p+\log(2+m_NT)}{m_N}
}
+
\frac{C_T}{m_N}.
\label{W:eq:Gaussian-offdiagonal-coefficients-rate}
\end{equation}
\end{proposition}

\begin{proof}
Write
\[
b_{k,N}^G-1
=
\left(
b_{k,N}^G-\mathbb Eb_{k,N}^G
\right)
+
\left(
\mathbb Eb_{k,N}^G-1
\right).
\]

For the deterministic term, \eqref{W:eq:Gaussian-b-mean-coarse} gives
\begin{equation}
 m_N \max_{k\le K_N(T)+1} |\mathbb Eb_{k,N}^G-1| \le C_Tm_N \left( \frac{m_N}{N}+\frac1N \right) \le \frac{C_T}{m_N}.
 \label{W:eq:Gaussian-offdiagonal-coefficients-drift}
\end{equation}

For the centered term,
\[
m_N
\left(
b_{k,N}^G-\mathbb Eb_{k,N}^G
\right)
=
\frac{m_N}{\sqrt{\beta N}}
\left(
\chi_{\beta(N-k)}
-
\mu_{\beta(N-k)}
\right).
\]
By \eqref{W:eq:chi-mgf}, each of these variables is sub-Gaussian with scale \(C_\beta\frac{m_N}{\sqrt N} = C_\beta m_N^{-1/2}.\) Applying the maximal estimate
\eqref{W:eq:subGaussian-maximal-bound} over \(K_N(T)+1\le m_NT+1\) indices yields
\[
\left\|
m_N
\max_{k\le K_N(T)+1}
|b_{k,N}^G-\mathbb Eb_{k,N}^G|
\right\|_{L^p}
\le
C_{\beta,p}
\sqrt{
\frac{p+\log(2+m_NT)}{m_N}
}.
\]
Combining this with \eqref{W:eq:Gaussian-offdiagonal-coefficients-drift} proves \eqref{W:eq:Gaussian-offdiagonal-coefficients-rate}.
\end{proof}

We now apply the telescoping identity \eqref{W:eq:integrated-potential-difference} to the Gaussian coefficients.

Introduce the auxiliary off-diagonal coefficient at zero
\begin{equation}
b_{0,N}^{G,\sharp}=1,
\qquad
b_{k,N}^{G,\sharp}=b_{k,N}^G,
\quad k\ge1.
\label{W:eq:Gaussian-sharp-offdiagonal-coefficients}
\end{equation}

Define \(\widetilde Y_{N,0}^G=0\) and
\begin{equation}
\widetilde Y_{N,k}^G
=
\sum_{r=1}^{k}
m_N
\left(
2-a_{r,N}^G
-b_{r-1,N}^{G,\sharp}
-b_{r,N}^G
\right).
\label{W:eq:Gaussian-exact-form-potential}
\end{equation}

For the independent sum representation, use the field \(Y_{N,k}^G = \sum_{r=1}^{k} m_N \left( 2-a_{r,N}^G-2b_{r,N}^G \right).\) Their difference is
\begin{equation}
 \widetilde Y_{N,k}^G-Y_{N,k}^G = m_N \sum_{r=1}^{k} \left( b_{r,N}^G-b_{r-1,N}^{G,\sharp} \right) = m_N \left( b_{k,N}^G-b_{0,N}^{G,\sharp} \right) = m_N(b_{k,N}^G-1).
 \label{W:eq:Gaussian-potential-boundary-identity}
\end{equation}

Let \(\widetilde Y_N^G\) denote the linear interpolation of the sampled values \(\widetilde Y_{N,k}^G.\)

By \eqref{W:eq:Gaussian-potential-boundary-identity}, Proposition~\ref{W:prop:Gaussian-offdiagonal-coefficients}, and Theorem~\ref{W:thm:Gaussian-potential-limit}, the difference
of the two integrated potential processes converges to zero in probability, uniformly on compact intervals, and the processes have the same limit $t^2/2+2B(t)/\sqrt\beta$.

\subsubsection{Characteristic functions and increments of the integrated potential}
\label{W:sec:exact-integrated-potential-quantitative-fourier-and-increment-bounds}
The exact integrated potential has dependent adjacent increments: they share an off-diagonal coefficient. We apply the independent sum estimates to the underlying diagonal and chi
variables, after collecting their weights. The characteristic function calculation retains the factor one half multiplying the last off-diagonal term in each partial sum.

Use ${m_N}=N^{1/3}$, ${h_N}={m_N}^{-1}$ and the coefficient classes of Section~\ref{W:sec:boundary-terms-in-the-integrated-potential}. The process limit is established by
\eqref{W:eq:Gaussian-potential-boundary-identity}. For entry comparison we also need a uniform characteristic function error and a fourth moment increment bound. We keep the
boundary weights in both calculations.

The previous process and maximal estimates give
\begin{equation}
 \widetilde Y_N^{\mathrm{lin}}(\cdot;H_N^{\mathrm G},w)
       \Longrightarrow Y_\beta(t)=\frac{t^2}{2}+\frac2{\sqrt\beta}B(t)
       \quad\text{in }C([0,T]).
 \label{W:eq:coefficients-gaussian-process}
\end{equation}
For the quantitative estimate, define
\begin{equation}
\begin{split}
 \mathcal E_N(\boldsymbol t,\boldsymbol\lambda,\boldsymbol\bullet)
 :=\bigg|&\mathbb E\exp\!\left(\ii\sum_{j=1}^{J_0}\lambda_j
                \widetilde Y_N^{\bullet_j}(t_j;H_N^{\mathrm G},w)\right)\\
 &-\exp\!\left(\ii\sum_{j=1}^{J_0}\lambda_jt_j^2/2
       -\frac2\beta\sum_{i,j=1}^{J_0}\lambda_i\lambda_j(t_i\wedge t_j)\right)
       \bigg|.
\end{split}
 \label{W:eq:coefficients-gaussian-fourier-error}
\end{equation}
The coefficients are those of \eqref{W:eq:coefficients-coefficient-class}, so this definition includes both step and linear evaluations.

\begin{theorem}
\label{W:thm:coefficients-gaussian-direct}
For fixed $T,J_0,\Lambda$ and all sufficiently large $N$,
\begin{equation}
 \sup_{\substack{t_j\in[0,T],\ \sum_j|\lambda_j|\le\Lambda\\
                 \bullet_j\in\{\mathrm{st},\mathrm{lin}\}}}
       \mathcal E_N(\boldsymbol t,\boldsymbol\lambda,\boldsymbol\bullet)
       \le C_{\beta,J_0,T,\Lambda}m_N^{-1/2}.
 \label{W:eq:gaussian-characteristic-error-bound}
\end{equation}
Moreover,
\begin{equation}
 \mathbb E|\widetilde Y_N^{\mathrm{lin}}(t;H_N^{\mathrm G})
                -\widetilde Y_N^{\mathrm{lin}}(s;H_N^{\mathrm G})|^4
       \le C_{\beta,T}|t-s|^2,\qquad s,t\in[0,T].
 \label{W:eq:coefficients-gaussian-increment-fourth}
\end{equation}
Both estimates are uniform over deterministic unit starting vectors.
\end{theorem}

\begin{proof}
The process limit \eqref{W:eq:coefficients-gaussian-process} follows from Theorem~\ref{W:thm:Gaussian-potential-limit}, \eqref{W:eq:Gaussian-potential-boundary-identity}, and
Proposition~\ref{W:prop:Gaussian-offdiagonal-coefficients}. The independence of the choice of starting vector is part of Theorem~\ref{W:thm:Gaussian-Lanczos-law}. We prove the two
quantitative bounds.

Set $U_r=-m_Na_r$, $V_r=-2m_N(b_r-\mathbb Eb_r)$ and $V_0=0$. The centered integrated potential increment is
\[
 \Delta\widetilde Y_{N,r}^G-\mathbb E\Delta\widetilde Y_{N,r}^G
       =U_r+\tfrac12(V_{r-1}+V_r).
\]
For $q_\rho=0$, summation and a change of index yield
\[
 \sum_{r<\rho}q_r\bigl(U_r+\tfrac12(V_{r-1}+V_r)\bigr)
 =\sum_{r<\rho}\left(q_rU_r+\frac{q_r+q_{r+1}}2V_r\right).
\]
This expresses the statistic as a sum of independent variables with weights $q$ and $d=A_\rho q$.

Put $\delta_r=1-\mathbb Eb_r$ and $\delta_0=0$. The elementary chi estimates of Lemma~\ref{W:lem:chi-estimates} give, uniformly for $r<\rho$,
\begin{equation}
 \delta_r=\frac{r}{2N}+O_{\beta,T}(N^{-1}),\qquad
 \operatorname{Var}(a_r)=\frac2{\beta N},\qquad
 \operatorname{Var}(b_r)=\frac1{2\beta N}+O_{\beta,T}(N^{-2}).
 \label{W:eq:coefficients-gaussian-input-moments}
\end{equation}
The mean is \eqref{W:eq:Gaussian-b-mean-coarse}; the variances follow from \eqref{W:eq:Gaussian-a-law} and \eqref{W:eq:chi-variance-expansion}.

By \eqref{W:eq:integrated-potential-weighted-sum},
\begin{equation}
 \mathbb E\widetilde Y_{N,k}^G
 =2{m_N}\sum_{r<k}\delta_r+{m_N}\delta_k
 =\frac{{m_N}}{N}\left(\frac{k(k-1)}2+\frac{k}2\right)
       +O_{\beta,T}({m_N} k/N)
 =\frac{(k{h_N})^2}{2}+O_{\beta,T}({h_N}).
 \label{W:eq:coefficients-gaussian-direct-drift}
\end{equation}
This also holds at $k=0$, where both expressions being compared are zero. The error is uniform through the additional sampling interval. Linear interpolation adds at most
${h_N}^2/8$ to the error for the quadratic function, and the step version differs from $t^2/2$ by at most $C_T {h_N}$.

For two evaluation coefficient pairs $(q^{(i)},d^{(i)})$, independence and \eqref{W:eq:coefficients-gaussian-input-moments} give
\begin{equation}
 \operatorname{Cov}(P^{\bullet_i}(t_i),P^{\bullet_j}(t_j))
 =\frac2\beta\left(\langle q^{(i)},q^{(j)}\rangle_{m_N}
                  +\langle d^{(i)},d^{(j)}\rangle_{m_N}\right)
   +O_{\beta,T}(N^{-1}).
 \label{W:eq:coefficients-gaussian-direct-covariance}
\end{equation}
The error from the variance approximation is bounded by $C{m_N}^2N^{-2}\sum_{r<\rho}|d_r^{(i)}d_r^{(j)}|=O_T(N^{-1})$. Let $e_r(t)=\mathbf1_{\{r\le\lfloor m_Nt\rfloor\}}$. The
evaluation weights obey $\|q^{(j)}-e(t_j)\|_{\ell^1}\le1$ and $\|d^{(j)}-e(t_j)\|_{\ell^1}\le2$, with all their entries in $[0,1]$. For either $u=q$ or $u=d$,
\begin{align*}
 \big|\langle u^{(i)},u^{(j)}\rangle_{m_N}-(t_i\wedge t_j)\big|
 &\le h_N\big(\|u^{(i)}-e(t_i)\|_{\ell^1}
                    +\|u^{(j)}-e(t_j)\|_{\ell^1}\big)\\
 &\quad+|h_N\lfloor m_N(t_i\wedge t_j)\rfloor-(t_i\wedge t_j)|
 \le5h_N.
\end{align*}
Thus the covariance is $4(t_i\wedge t_j)/\beta+O_{\beta,T}(h_N)$, uniformly over the observation times and both interpolation conventions.

We give a quantitative characteristic function proof. For $q,d$ in \eqref{W:eq:coefficients-coefficient-class}, the centered statistic is the sum of the independent variables
$X_r^a=-{m_N} q_r a_r$ and $X_r^b=-2{m_N} d_r(b_r-\mathbb Eb_r)$. Gaussian moments and the chi concentration bound imply
\begin{equation}
 \sum_r\bigl(\mathbb E|X_r^a|^3+\mathbb E|X_r^b|^3\bigr)
       \le C_{\beta,T,\Lambda}{h_N}^{1/2},\qquad
 \sum_r\bigl(\operatorname{Var}(X_r^a)^2+
             \operatorname{Var}(X_r^b)^2\bigr)
       \le C_{\beta,T,\Lambda}{h_N}.
 \label{W:eq:coefficients-gaussian-lyapunov}
\end{equation}
To verify the first estimate, each absolute third moment is bounded by $C{h_N}^{3/2}(|q_r|^3+|d_r|^3)$ and there are $O_T({m_N})$ bounded coefficients. For the variances,
\[
 \operatorname{Var}(X_r^a)\le C h_Nq_r^2,\qquad
 \operatorname{Var}(X_r^b)\le C h_Nd_r^2.
\]
Thus the sum of their squares is at most $C h_N^2\sum_r(q_r^4+d_r^4)\le C_{\beta,T,\Lambda}h_N$, since the coefficients are uniformly bounded and their number is $O_T(m_N)$. Write
$\mu_N$ and $v_N$ for the mean and variance of the weighted sum, and put $\mu=\sum_j\lambda_jt_j^2/2$ and $v=(4/\beta)\sum_{i,j}\lambda_i\lambda_j(t_i\wedge t_j)$. The product
identity used in the proof of Theorem~\ref{W:thm:Gaussian-potential-limit} gives
\[
 |\mathbb Ee^{\ii(\mathscr Z-\mu_N)}-e^{-v_N/2}|
 \le\frac16\sum_{r,\tau}\mathbb E|X_r^\tau|^3
             +\frac18\sum_{r,\tau}\operatorname{Var}(X_r^\tau)^2.
\]
The drift and covariance estimates imply $|\mu_N-\mu|+|v_N-v|\le C_{\beta,J_0,T,\Lambda}h_N$. Since $v_N,v\ge0$,
\[
 |e^{\ii\mu_N-v_N/2}-e^{\ii\mu-v/2}|
       \le |\mu_N-\mu|+\tfrac12|v_N-v|.
\]
Adding the last two bounds and using \eqref{W:eq:coefficients-gaussian-lyapunov} proves \eqref{W:eq:gaussian-characteristic-error-bound}.

It remains to prove the increment bound. For $q=q^{\mathrm{lin},t}-q^{\mathrm{lin},s}$ and $d=A_{\rho}q$, the same independence and the chi exponential bound give
$\mathbb E\exp(u(\mathscr Z(q,d)-\mathbb E\mathscr Z(q,d))) \le\exp(C_\beta u^2{h_N}\sum_r(q_r^2+d_r^2))$ for every real $u$. Chernoff's bound and integration of the resulting
Gaussian tail imply
\[
 \mathbb E|\mathscr Z(q,d)-\mathbb E\mathscr Z(q,d)|^4
 \le C_\beta\bigl({h_N}\sum_r(q_r^2+d_r^2)\bigr)^2
 \le C_\beta|t-s|^2
\]
by \eqref{W:eq:coefficients-increment-coefficients} and Lemma~\ref{W:lem:integrated-potential-weight-bounds}. In addition,
$|\mathbb EP_k-\mathbb EP_{k-1}| ={m_N}|\delta_{k-1}+\delta_k|\le C_{\beta,T}{h_N}$. The mean linear path is consequently uniformly Lipschitz. The inequality
$|x+y|^4\le8(|x|^4+|y|^4)$ proves \eqref{W:eq:coefficients-gaussian-increment-fourth}, completing the quantitative estimates.
\end{proof}

\subsubsection{Joint limits at the two spectral edges}
\label{W:sec:joint-limits-at-the-two-spectral-edges}
The two edges use the same diagonal and off-diagonal variables at finite $N$. They are therefore treated jointly. Reflection reverses the diagonal contribution and preserves the
off-diagonal contribution. Their equal limiting variances make the cross covariance vanish, giving independent components in the joint Gaussian limit.

Let $X_N^\sigma$ be the two integrated potentials of Section~\ref{W:sec:reflection-and-the-modulus-of-continuity}.

\Needspace{6\baselineskip}

\begin{theorem}
\label{W:thm:tightness-gaussian}
For $H_N^{\mathrm G}$,
\begin{equation}
 (X_N^+,X_N^-)\Longrightarrow
 \left(t^2/2+2B_+(t)/\sqrt\beta,
       t^2/2+2B_-(t)/\sqrt\beta\right)
       \quad\text{in }C([0,T])^2,
 \label{W:eq:tightness-gaussian-pair}
\end{equation}
where $B_+,B_-$ are independent standard Brownian motions. For each sign,
\begin{equation}
 \mathbb P(\omega_T(X_N^\sigma,\delta)>u)
       \le C_{\beta,T}u^{-4}\delta,\qquad 0<\delta\le T,
 \label{W:eq:tightness-gaussian-modulus}
\end{equation}
for all sufficiently large $N$.
\end{theorem}
\begin{proof}
By the reflection identity of Lemma~\ref{W:lem:tightness-reflection}, the two integrated potentials use the same off-diagonal coefficients and opposite diagonal coefficients. Let
$R_N$ be the linear interpolation of $R_{N,0}=0$ and $R_{N,k}=m_N(b_{k,N}^G-1)$ for $k\ge1$. The exact telescoping identity \eqref{W:eq:Gaussian-potential-boundary-identity} gives,
on $[0,T]$,
\begin{equation}
 X_N^+=D_N+M_N^a+M_N^b+R_N,
 \qquad
 X_N^-=D_N-M_N^a+M_N^b+R_N.
 \label{W:eq:gaussian-reflected-process-decomposition}
\end{equation}
Proposition~\ref{W:prop:Gaussian-offdiagonal-coefficients}, applied with an additional sampling interval, gives $\|R_N\|_{C([0,T])}\to0$ in probability. The joint coefficient
process limit \eqref{W:eq:Gaussian-joint-martingale-limit} and the continuous mapping theorem therefore prove \eqref{W:eq:tightness-gaussian-pair}, with
\[
 B_+=(B_1+B_2)/\sqrt2,
 \qquad B_-=(-B_1+B_2)/\sqrt2.
\]
This pair is jointly Gaussian, and $\mathbb E[B_+(s)B_-(t)]=(-\min(s,t)+\min(s,t))/2=0$. Thus $B_+$ and $B_-$ are independent as processes, although the finite-$N$ integrated
potentials share all their Jacobi coefficients.

The same cancellation is visible in the exact finite-$N$ weights. For $q^s,q^t$ representing evaluations of the integrated potential and $d=A_\rho q$,
\eqref{W:eq:coefficients-gaussian-input-moments} gives
\begin{equation}
 \operatorname{Cov}(X_N^+(s),X_N^-(t))
 =-\frac2\beta\langle q^s,q^t\rangle_{m_N}
  +\frac2\beta\langle A_\rho q^s,A_\rho q^t\rangle_{m_N}
  +O_{\beta,T}(N^{-1})=O_{\beta,T}(h_N).
 \label{W:eq:tightness-gaussian-cross}
\end{equation}
Both inner products equal $\min(s,t)+O_T(h_N)$, including the interpolated weights at the boundary sites.

Each marginal $X_N^\sigma$ has the law of the integrated potential in Theorem~\ref{W:thm:coefficients-gaussian-direct}, since $-H_N^{\mathrm G}$ has the same distribution as
$H_N^{\mathrm G}$. Its fourth moment bound \eqref{W:eq:coefficients-gaussian-increment-fourth}, followed by Lemma~\ref{W:lem:tightness-dyadic} with $p=4$, $\alpha=1$, proves
\eqref{W:eq:tightness-gaussian-modulus}.
\end{proof}

\begin{lemma}
\label{W:lem:gaussian-offdiagonal-maximum}
For $1\le K\le N/2$ and $s\ge0$,
\begin{equation}
 \mathbb P\left(\max_{k\le K}|b_k-1|>
 C_\beta\left\{\frac K N+
       \sqrt{\frac{\log(K+1)+s+1}{N}}\right\}\right)
       \le C_\beta e^{-s}.
 \label{W:eq:tightness-gaussian-offdiagonals}
\end{equation}
\end{lemma}
\begin{proof}
For $K\le N/2$, Lemma~\ref{W:lem:chi-estimates} gives $|\mathbb Eb_k-1|\le k/N+C_\beta/N$. Its concentration bound and a union bound imply
\[
 \mathbb P\left(\max_{k\le K}|b_k-1|>
                   K/N+C_\beta/N+v\right)
 \le 2K\exp(-c_\beta Nv^2),\qquad v\ge0.
\]
Take $v=C_\beta\sqrt{(\log(K+1)+s+1)/N}$ and enlarge $C_\beta$ to obtain \eqref{W:eq:tightness-gaussian-offdiagonals}.
\end{proof}

\subsubsection{Gram matrices for the GOE and GUE}
\label{W:sec:bounds-on-gaussian-gram-matrices}

Coefficient convergence does not by itself control orthogonalization at a growing degree. The estimates below control the eigenvalues of the Gram matrix on the scale used by the
entry comparison.

Fix $T>0$ and a nonnegative integer $r_0$, with the symmetry parameter $\beta\in\{1,2\}$ fixed as above. We consider integers $1\le\rho\le TN^{1/3}+r_0$. Choose
$N_0=N_0(\beta,T,r_0)$ so that, for $N\ge N_0$,
\[
 \rho\le(T+1)N^{1/3}\le N/2,\qquad
 \frac{\rho^3}{N}\le(T+1)^3,\qquad \frac{\beta(N-\rho)}2\ge1.
\]
Every estimate below is uniform over these $N,\rho$. The constants $C_{p,T}$ may also depend on the fixed $\beta$. The additive constant $r_0$ only affects $N_0$. Starting vectors
belong to $\mathbb F_\beta^N$ and are fixed before the matrix is sampled.

We use the monic Chebyshev polynomials $\psi_{-1}=0$, $\psi_0=1$ and $\psi_{n+1}(x)=x\psi_n(x)-\psi_{n-1}(x)$. Thus $\psi_n(x)=U_n(x/2)$, consistently with
\eqref{W:eq:chebyshev-three-term}. For a deterministic unit vector $w$, the Gram matrix is
\[
 G_\rho(H_N^{\mathrm G})_{ij}
 =\bigl(\psi_{i-1}(H_N^{\mathrm G})w\bigr)^*\psi_{j-1}(H_N^{\mathrm G})w,
 \qquad 1\le i,j\le\rho.
\]
It is a Borel function of the matrix entries. The vector $w$ is fixed in this notation. Matrix norms without a subscript are operator norms for the Euclidean norm, and
$\|C\|_{\mathrm{HS}}^2=\operatorname{Tr}(C^*C)$. We prove the following uniform tightness statement.

\begin{lemma}
\label{W:lem:gaussian-gram-conditioning}
For $1\le\rho_N\le TN^{1/3}+r_0$, the two sequences
\begin{equation}
 \lambda_{\max}(G_{\rho_N}(H_N^{\mathrm G})),\qquad
       \lambda_{\min}(G_{\rho_N}(H_N^{\mathrm G}))^{-1}
 \label{W:eq:gram-polynomial-matrix-tight-random-variables}
\end{equation}
are tight, uniformly over deterministic starting vectors.
\end{lemma}

The lemma controls the largest eigenvalue of the Gram matrix and the reciprocal of the smallest. At the present degree, the polynomial expansion has scale $\rho^3/N=O(1)$. We bound
an upper triangular factor in Hilbert--Schmidt norm and use its determinant to control its smallest singular value. The proof is completed in
Section~\ref{W:sec:gaussian-gram-singular-values}.

Apply Lanczos iteration to $H_N^{\mathrm G}$ from $w$. By Theorem~\ref{W:thm:Gaussian-Lanczos-law}, on an event of probability one all $b_{j,N}^G$ are positive and the iteration
produces an orthogonal or unitary matrix $U=(q_1,\ldots,q_N)$ with $q_1=w$ and $H_N^{\mathrm G}U=UJ_N^G$. The real Jacobi matrix $J_N^G$ has the independent coefficient law in that
theorem. All algebraic identities below are on this event; the finitely many coefficients are assigned $a_{j,N}^G=0$ and $b_{j,N}^G=1$ on its null complement, so that the
subsequent formulas define finite measurable matrices everywhere. Define
\begin{equation}
 \mathsf C_\rho
 =\bigl(\psi_0(J_N^G)e_1,\ldots,
             \psi_{\rho-1}(J_N^G)e_1\bigr)_{\{1,\ldots,\rho\}\times\{1,\ldots,\rho\}}.
 \label{W:eq:gram-connection-matrix}
\end{equation}
For every integer $n\ge0$, tridiagonality gives $\operatorname{supp}((J_N^G)^ne_1)\subseteq\{1,\ldots,n+1\}$. If $n<N$, the unique path of length $n$ from coordinate $1$ to $n+1$
uses the successive subdiagonal entries. Since $\psi_n$ is monic,
\[
 \bigl(\psi_n(J_N^G)e_1\bigr)_{n+1}=\prod_{j=1}^n b_{j,N}^G.
\]
The empty product is one. Hence $\mathsf C_\rho$ is upper triangular with positive diagonal. For $Q_\rho=(q_1,\ldots,q_\rho)$, polynomial functional calculus and this support
identity yield
\begin{equation}
 \bigl(\psi_0(H_N^{\mathrm G})w,\ldots,\psi_{\rho-1}(H_N^{\mathrm G})w\bigr)
       =Q_\rho\mathsf C_\rho.
 \label{W:eq:gram-connection-factorization}
\end{equation}
Taking adjoints and using $Q_\rho^*Q_\rho=I_\rho$ gives
\begin{equation}
 G_\rho(H_N^G)=\mathsf C_\rho^*\mathsf C_\rho.
 \label{W:eq:gram-Gram-connection-factorization}
\end{equation}
Thus the eigenvalues of the Gram matrix are the squares of the singular values of $\mathsf C_\rho$. The diagonal change of coordinates below expresses its normalized columns as
polynomials in independent normal and chi-square variables.

Let $J_0$ be the free Jacobi operator on $\ell^2(\mathbb N)$, $J_0e_j=e_{j-1}+e_{j+1}$, with $e_0=0$. If $S e_j=e_{j+1}$ is the unilateral shift, then $J_0=S+S^*$; in particular
$J_0$ is bounded and self-adjoint on all of $\ell^2(\mathbb N)$ and $\|J_0\|\le2$. Put $\mathsf U_r=\psi_r(J_0)$ for $r\ge0$ and $\mathsf U_{-1}=0$. These are bounded polynomial
functions of $J_0$.

\Needspace{6\baselineskip}

\begin{lemma}
\label{W:lem:free-chebyshev-kernel}
For $r\ge0$ and $x,y\ge1$,
\begin{equation}
 (\mathsf U_r)_{xy}
 =\mathbf1_{\{r\ge|x-y|,\ r\equiv x-y\ ({\rm mod}\ 2)\}}
 -\mathbf1_{\{r\ge x+y,\ r\equiv x+y\ ({\rm mod}\ 2)\}}.
 \label{W:eq:gram-free-kernel-formula}
\end{equation}
In particular,
\begin{equation}
 |(\mathsf U_r)_{xy}|\le1,
 \qquad
 (\mathsf U_r)_{xy}=0\quad\hbox{if }|x-y|>r,
 \label{W:eq:free-chebyshev-support}
\end{equation}
and
\begin{equation}
 \mathsf U_re_1=e_{r+1}.
 \label{W:eq:free-chebyshev-first-column}
\end{equation}
\end{lemma}

\begin{proof}
For $k\in\mathbb Z$, set $F_r(k)=\mathbf1_{\{|k|\le r,\ k\equiv r\ ({\rm mod}\ 2)\}}$ and $F_{-1}=0$. The generating polynomial is
\[
 \sum_{k\in\mathbb Z}F_r(k)z^k=z^{-r}+z^{-r+2}+\cdots+z^r.
\]
Multiplication by $z+z^{-1}$ counts each interior monomial twice. Subtracting the corresponding polynomial of order $r-1$ leaves that of order $r+1$. Coefficient comparison
therefore gives $F_{r+1}(k)=F_r(k-1)+F_r(k+1)-F_{r-1}(k)$, including $r=0$.

Define $K_r(x,y)=F_r(x-y)-F_r(x+y)$ for $x,y\ge1$ and set $K_r(0,y)=0$. This boundary value agrees with the formula because $F_r$ is even. For each fixed $r$, the matrix $K_r$ has
at most $2r+1$ nonzero entries in each row and column, and its entries have absolute value at most one. For finitely supported $v$, Cauchy--Schwarz gives
\[
 \sum_x\left|\sum_y K_r(x,y)v_y\right|^2
 \le(2r+1)\sum_{x,y}|K_r(x,y)|\,|v_y|^2
 \le(2r+1)^2\|v\|_2^2.
\]
Thus $K_r$ extends to a bounded operator. The recurrence for $F_r$ yields $K_{r+1}=J_0K_r-K_{r-1}$, first on finitely supported vectors and hence everywhere. Since $K_0=I$ and
$K_{-1}=0$, induction identifies $K_r=\psi_r(J_0)$. This proves \eqref{W:eq:gram-free-kernel-formula}. The second indicator in that formula is bounded by the first: $x+y\ge|x-y|$
and both have the same parity. Consequently the kernel is either zero or one and has the support in \eqref{W:eq:free-chebyshev-support}. For $y=1$ the two indicators differ exactly
when $x=r+1$, which proves \eqref{W:eq:free-chebyshev-first-column}.
\end{proof}

For bounded operators $A,B$ on the same Hilbert space, the same recurrence yields
\begin{equation}
 \psi_n(A)-\psi_n(B)
 =\sum_{r=0}^{n-1}\psi_{n-1-r}(B)(A-B)\psi_r(A).
 \label{W:eq:gram-Chebyshev-Duhamel}
\end{equation}
Indeed, denote the sum by $D_n$ and set $D_0=0$. The terms with $r\le n-1$ in $BD_n-D_{n-1}$ combine by $B\psi_j(B)-\psi_{j-1}(B)=\psi_{j+1}(B)$, giving
\[
 D_{n+1}=BD_n+(A-B)\psi_n(A)-D_{n-1}.
\]
If $D_j=\psi_j(A)-\psi_j(B)$ for $j=n,n-1$, substitution gives
\[
 D_{n+1}=A\psi_n(A)-\psi_{n-1}(A)
             -B\psi_n(B)+\psi_{n-1}(B)
          =\psi_{n+1}(A)-\psi_{n+1}(B).
\]
Since $D_1=A-B$, induction proves \eqref{W:eq:gram-Chebyshev-Duhamel}.

For $1\le k\le {\rho}$, set
\begin{equation}
 d_{1,N}:=1,\qquad
 d_{k,N}:=\prod_{j=1}^{k-1}b_{j,N}^G,
 \qquad
 \mathsf D_{\rho}:=\operatorname{diag}(d_{1,N},\ldots,d_{{\rho},N}).
 \label{W:eq:gaussian-diagonal-products}
\end{equation}
Let $J_{N,{\rho}}^G$ be the leading ${\rho}\times {\rho}$ principal Jacobi block.

\Needspace{6\baselineskip}

Let $J_{0,\rho}$ denote the leading $\rho\times\rho$ block of $J_0$, and put
\begin{equation}
 \widetilde J_\rho=\mathsf D_\rho^{-1}J_{N,\rho}^G\mathsf D_\rho,
 \qquad \widetilde{\mathsf C}_\rho=\mathsf D_\rho^{-1}\mathsf C_\rho.
 \label{W:eq:gaussian-conjugated-matrices}
\end{equation}

\begin{lemma}
\label{W:lem:gaussian-diagonal-conjugation}
The conjugated matrices satisfy
\begin{equation}
 \widetilde J_\rho=J_{0,\rho}+\widetilde D_\rho,\qquad
 \widetilde D_\rho=\sum_{j=1}^\rho a_{j,N}^G E_{jj}
       +\sum_{j=1}^{\rho-1}\bigl((b_{j,N}^G)^2-1\bigr)E_{j,j+1}.
 \label{W:eq:gaussian-conjugated-perturbation}
\end{equation}
Moreover,
\begin{equation}
 \widetilde{\mathsf C}_\rho
   =\bigl(\psi_0(\widetilde J_\rho)e_1,\ldots,
                     \psi_{\rho-1}(\widetilde J_\rho)e_1\bigr),\qquad
 \mathsf C_\rho=\mathsf D_\rho\widetilde{\mathsf C}_\rho.
 \label{W:eq:gaussian-polynomial-conjugation}
\end{equation}
The matrix $\widetilde{\mathsf C}_\rho$ is upper triangular with diagonal one.
\end{lemma}

\begin{proof}
All entries of $\mathsf D_\rho$ are finite and strictly positive on the event of probability one already specified. From $d_{j+1,N}=b_{j,N}^Gd_{j,N}$ we obtain, for $1\le j<\rho$,
\[
 (\mathsf D_\rho^{-1}J_{N,\rho}^G\mathsf D_\rho)_{j+1,j}
       =\frac{b_{j,N}^Gd_{j,N}}{d_{j+1,N}}=1,\qquad
 (\mathsf D_\rho^{-1}J_{N,\rho}^G\mathsf D_\rho)_{j,j+1}
       =\frac{b_{j,N}^Gd_{j+1,N}}{d_{j,N}}=(b_{j,N}^G)^2.
\]
The diagonal is unchanged, proving \eqref{W:eq:gaussian-conjugated-perturbation}. For $0\le n<\rho$, polynomial similarity and $\mathsf D_\rho e_1=e_1$ give
\[
 \psi_n(\widetilde J_\rho)e_1
   =\mathsf D_\rho^{-1}\psi_n(J_{N,\rho}^G)e_1.
\]
For each power of degree at most $n$, a path starting at coordinate one has at most $n$ nearest-neighbour steps. It cannot reach coordinate $\rho+1$, since that requires $\rho$
steps. Hence $\psi_n(J_{N,\rho}^G)e_1$ is precisely the restriction of $\psi_n(J_N^G)e_1$ to the first $\rho$ coordinates, with no discarded coordinates. Applying this to every
column proves \eqref{W:eq:gaussian-polynomial-conjugation}. The columns are upper triangular by the same support bound. Their diagonal entries are
\[
 d_{n+1,N}^{-1}\bigl(\psi_n(J_N^G)e_1\bigr)_{n+1}
 =d_{n+1,N}^{-1}\prod_{j=1}^n b_{j,N}^G=1.
\]
These identities hold pathwise for the two matrices constructed from the same Jacobi coefficients.
\end{proof}

Put
\begin{equation}
 \delta_j:=\mathbb E(b_{j,N}^G)^2-1=-\frac{j}{N},
 \qquad
 \zeta_j:=(b_{j,N}^G)^2-\mathbb E(b_{j,N}^G)^2.
 \label{W:eq:gram-square-offdiagonal-split}
\end{equation}

\Needspace{6\baselineskip}

\begin{lemma}
\label{W:lem:gram-Gaussian-Jacobi-cumulants}
Let $\nu_j=\beta(N-j)$. Uniformly for $j\le TN^{1/3}+r_0$,
\begin{align}
 &\kappa_2(a_{j,N}^G)=\frac{2}{\beta N},
 \qquad
 \kappa_r(a_{j,N}^G)=0\quad(r\ne2),
 \label{W:eq:gram-a-cumulants}\\
 &\kappa_r(\zeta_j)
 =\frac{2^{r-1}(r-1)!\,\nu_j}{(\beta N)^r},
 \qquad r\ge2,
 \label{W:eq:gram-square-offdiagonal-cumulants-exact}\\
 &|\kappa_r(\zeta_j)|
 \le (r-1)!C_\beta^rN^{1-r},
 \qquad r\ge2.
 \label{W:eq:gram-square-offdiagonal-cumulants-bound}
\end{align}
The family consisting of $a_{j,N}^G$, $1\le j\le N$, and $\zeta_j$, $1\le j<N$, is independent.
\end{lemma}

\begin{proof}
For a real variable with moment-generating function finite near zero, use $\kappa_r(X)=\left.\frac{\dd^r}{\dd t^r}\log\mathbb E e^{tX}\right|_{t=0}$ for $r\ge1$. The Gaussian
coefficient law gives
\[
 \log\mathbb E e^{t a_{j,N}^G}=\frac{t^2}{\beta N},\qquad t\in\mathbb R,
\]
so its only nonzero cumulant is that of order two. For $1\le j<N$, write $X_j=\beta N(b_{j,N}^G)^2\sim\chi^2_{\nu_j}$. Its density is
\[
 \frac{x^{\nu_j/2-1}e^{-x/2}}{2^{\nu_j/2}\Gamma(\nu_j/2)}\,\mathbf1_{\{x>0\}}.
\]
For real $t<\beta N/2$, integration with the positive rate $1/2-t/(\beta N)$ gives
\[
 \mathbb E e^{t\zeta_j}
 =e^{-\nu_jt/(\beta N)}
      \left(1-\frac{2t}{\beta N}\right)^{-\nu_j/2}.
\]
On every compact subinterval of $(-\infty,\beta N/2)$, each differentiated integrand is bounded by a constant times $(1+x^r)x^{\nu_j/2-1}e^{-c x}$ with $c>0$. This is integrable at
both zero and infinity. Differentiation under the integral and the Taylor series at zero are therefore valid at every finite order. In particular, for $|t|<\beta N/2$,
\[
 \log\mathbb E e^{t\zeta_j}
 =\sum_{r=2}^{\infty}\frac{\nu_j2^{r-1}}{r(\beta N)^r}\,t^r.
\]
Multiplying each coefficient of the expansion by $r!$ proves \eqref{W:eq:gram-square-offdiagonal-cumulants-exact}. Since $\nu_j\le\beta N$,
\[
 |\kappa_r(\zeta_j)|
 \le(r-1)!\,2^{r-1}\beta^{1-r}N^{1-r}
 \le(r-1)!C_\beta^rN^{1-r},\qquad r\ge2.
\]
One constant $C_\beta$ works for all $r$, including orders growing with $\rho$. Finally, the coefficient law in Theorem~\ref{W:thm:Gaussian-Lanczos-law} is a product law; applying
the separate functions $b\mapsto b^2-\mathbb E b^2$ to its off-diagonal coordinates preserves independence. This independence holds under the unconditional law, separately for each
$N$.
\end{proof}

For the diagonal products we also need $\ell_j=\log b_{j,N}^G$.
\begin{lemma}
\label{W:lem:gaussian-logarithmic-coefficients}
Uniformly for $j\le TN^{1/3}+r_0$ and all sufficiently large $N$,
\begin{align}
 |\mathbb E\ell_j|&\le C_{\beta,T}\frac{j+1}{N},
       \label{W:eq:gram-logb-mean}\\
 |\kappa_r(\ell_j)|&\le(r-1)!C_{\beta,T}^rN^{1-r},\qquad r\ge2.
       \label{W:eq:gram-logb-cumulants}
\end{align}
\end{lemma}
\begin{proof}
Put $x_j=\beta(N-j)/2$. The chi density gives the Mellin transform
\[
 \mathbb E e^{t\ell_j}
 =\left(\frac2{\beta N}\right)^{t/2}
          \frac{\Gamma(x_j+t/2)}{\Gamma(x_j)},\qquad t>-2x_j.
\]
For each finite derivative order, differentiation at zero is justified by the integrable bounds $x^{\nu_j/2-1}(1+|\log x|^r)$ on $(0,1)$ and
$x^{3\nu_j/2-1}(1+|\log x|^r)e^{-x^2/2}$ on $(1,\infty)$, after restricting $|t|\le\nu_j/2$ in the chi-density integral. Thus
\begin{equation}
 \begin{aligned}
 \mathbb E\ell_j
   &=\tfrac12\bigl((\log\Gamma)'(x_j)+\log2-\log(\beta N)\bigr),\\
 \kappa_r(\ell_j)&=2^{-r}(\log\Gamma)^{(r)}(x_j),\qquad r\ge2.
 \end{aligned}
 \label{W:eq:gram-logchi-exact-cumulant}
\end{equation}
Formula (5.7.6) of \cite{nistGammaFunctions} gives, for $x>0$,
\[
 (\log\Gamma)'(x)
 =\lim_{m\to\infty}\left(\log m-\sum_{k=0}^{m-1}\frac1{x+k}\right).
\]
Indeed the difference between the partial harmonic sum and $\log m$ tends to Euler's constant, as used in that formula. The derivative series is locally uniformly convergent on
$(0,\infty)$: on every compact interval $x\ge a>0$, its terms of order $r\ge2$ are bounded by $(r-1)!(a+k)^{-r}$, a summable sequence for each $r$. Termwise differentiation
therefore yields
\[
 (\log\Gamma)^{(r)}(x)=(-1)^r(r-1)!\sum_{k=0}^{\infty}(x+k)^{-r},\qquad r\ge2.
\]

For the first derivative, comparison with the integral of the decreasing function $(x+t)^{-1}$ gives
\[
 0\le\sum_{k=0}^{m-1}\frac1{x+k}
            -\log\frac{x+m}{x}\le\frac1x.
\]
After $m\to\infty$ this implies $|(\log\Gamma)'(x)-\log x|\le x^{-1}$. Since $2x_j/(\beta N)=1-j/N$ and $j/N\le1/2$,
\[
 |\mathbb E\ell_j|
 \le\tfrac12\left(-\log(1-j/N)+x_j^{-1}\right)
 \le\frac jN+\frac1{2x_j}
 \le C_{\beta,T}\frac{j+1}{N}.
\]
For the higher derivatives, the same integral comparison gives
\[
 \sum_{k=0}^{\infty}(x+k)^{-r}
 \le x^{-r}+\frac{x^{1-r}}{r-1},\qquad r\ge2.
\]
Our choice of $N_0$ gives $x_j\ge\beta N/4$ and $x_j\ge1$. Consequently
\[
 |\kappa_r(\ell_j)|
 \le2^{-r}(r-1)!\left(x_j^{-r}+\frac{x_j^{1-r}}{r-1}\right)
 \le(r-1)!C_\beta^rN^{1-r},\qquad r\ge2,
\]
with a constant independent of $r,j,N$. This proves the mean bound and the bounds on cumulants of order at least two.
\end{proof}

\subparagraph{Expansion of the polynomial vectors}
\label{W:sec:gaussian-polynomial-perturbations}

Introduce a scalar parameter $t$ in the finite matrix $J_{0,\rho}+t\widetilde D_\rho$. Its polynomial vector $\psi_n(J_{0,\rho}+t\widetilde D_\rho)e_1$ has degree at most $n$ in
$t$. Applying \eqref{W:eq:gram-Chebyshev-Duhamel} and equating the coefficient of $t^q$ gives, for $q\ge1$,
\[
 [t^q]\psi_n(J_{0,\rho}+t\widetilde D_\rho)e_1
 =\sum_{j=0}^{n-1}\psi_{n-1-j}(J_{0,\rho})\widetilde D_\rho
       [t^{q-1}]\psi_j(J_{0,\rho}+t\widetilde D_\rho)e_1.
\]
Starting from the constant coefficient $\psi_n(J_{0,\rho})e_1=e_{n+1}$ for $n<\rho$, induction on $q$ and $n$ proves the following finite identity at $t=1$:
\begin{equation}
 \psi_n(\widetilde J_\rho)e_1
 =e_{n+1}+\sum_{q=1}^{n}\widetilde{\mathsf V}_{n,q},\qquad 0\le n<\rho,
 \label{W:eq:gaussian-polynomial-expansion}
\end{equation}
where
\begin{equation}
 \widetilde{\mathsf V}_{n,q}
 :=\sum_{\substack{s_0,\ldots,s_q\ge0\\s_0+\cdots+s_q+q=n}}
 \mathsf U_{s_0}\widetilde D_\rho
 \mathsf U_{s_1}\widetilde D_\rho\cdots
 \widetilde D_\rho\mathsf U_{s_q}e_1.
 \label{W:eq:gaussian-perturbation-order}
\end{equation}
In the induction, $s_0=n-1-j$ and the remaining product with $q-1$ factors $\widetilde D_\rho$ has total polynomial degree $j-(q-1)$; this gives exactly the constraint in
\eqref{W:eq:gaussian-perturbation-order}. The variables $s_i$ denote polynomial degrees and are distinct from the fixed integer $r_0$ in the upper bound on $\rho$. Let
$\widetilde{\mathsf C}_\rho^{(q)}$ have column $n+1$ equal to $\widetilde{\mathsf V}_{n,q}$, with zero columns when $n<q$. Then
\begin{equation}
 \widetilde{\mathsf C}_\rho-I_\rho
 =\sum_{q=1}^{\rho-1}\widetilde{\mathsf C}_\rho^{(q)}.
 \label{W:eq:gaussian-coefficient-expansion}
\end{equation}
This is a finite polynomial identity and requires no smallness assumption on $\widetilde D_\rho$. For example,
\[
 \begin{aligned}
 \widetilde{\mathsf V}_{n,1}
 &=\sum_{s_0+s_1=n-1}\mathsf U_{s_0}\widetilde D_\rho\mathsf U_{s_1}e_1,\\
 \widetilde{\mathsf V}_{n,2}
 &=\sum_{s_0+s_1+s_2=n-2}
       \mathsf U_{s_0}\widetilde D_\rho\mathsf U_{s_1}\widetilde D_\rho\mathsf U_{s_2}e_1.
 \end{aligned}
\]
Here $\widetilde D_\rho$ is extended by zero to $\ell^2(\mathbb N)$. To verify the use of the half-line kernels, expand every polynomial into powers. A term in the resulting matrix
product contains at most $s_0+\cdots+s_q+q=n$ nearest-neighbour steps starting at coordinate one. Reaching coordinate $\rho+1$ would require at least $\rho$ steps, so every path
lies in $\{1,\ldots,\rho\}$. Each finite-matrix product therefore agrees with its half-line expression in \eqref{W:eq:gaussian-perturbation-order}. The case $n=0$ gives the column
$e_1$ with no factors $\widetilde D_\rho$.

For the index sums write $E_j^{(0)}=E_{jj}$ and $E_j^{(1)}=E_{j,j+1}$. The coefficient of $E_j^{(0)}$ is $a_j$, and that of $E_j^{(1)}$ is $\delta_j+\zeta_j$. We keep the
dependence on $q$ in the counting estimate below.

\Needspace{6\baselineskip}

Fix integers $p,q\ge1$ and put $M=2pq$. The identity
\[
 \|\widetilde{\mathsf C}_\rho^{(q)}\|_{\rm HS}^{2p}
 =\sum_{\substack{0\le n_1,\ldots,n_p<\rho\\1\le x_1,\ldots,x_p\le\rho}}
       \prod_{a=1}^p
       (\widetilde{\mathsf V}_{n_a,q})_{x_a}
       \overline{(\widetilde{\mathsf V}_{n_a,q})_{x_a}}
\]
contains $2p$ products, with $q$ factors from $\widetilde D_\rho$ in each. Index the positions of these factors by $\{1,\ldots,p\}\times\{1,2\}\times\{1,\ldots,q\}$. At a diagonal
position use $a_j$; at an upper off-diagonal position expand $\delta_j+\zeta_j$. Let $D$ be the positions assigned to $\delta_j$. The moment--cumulant formula partitions the
remaining positions into blocks. More precisely, let $A$ be a set of centered occurrences, and write $\xi_{j,0}=a_j$ and $\xi_{j,1}=\zeta_j$. Independence and the
cumulant-generating functions give, for $t_v$ near zero,
\[
 \mathbb E\exp\left(\sum_{v\in A}t_v\xi_{j_v,\theta_v}\right)
 =\prod_{(j,\theta)}\exp\left\{
     \sum_{r\ge2}\frac{\kappa_r(\xi_{j,\theta})}{r!}
      \left(\sum_{v\in A:(j_v,\theta_v)=(j,\theta)}t_v\right)^r\right\}.
\]
Only coefficient indices occurring in $A$ appear in the product. For example, $\sum_{v\in A}|t_v|<\beta N/4$ lies inside every required chi-square convergence domain. Taking the
coefficient of $\prod_{v\in A}t_v$ gives the moment--cumulant formula: each block consists of occurrences with the same label and type, and a block of size $r$ contributes
$\kappa_r$, because its $r!$ orderings cancel the denominator $r!$. Centered variables have no singleton blocks. Different blocks may have the same coefficient index; these
contributions must also be included. All polynomial sums and expectations at fixed $N,p,q$ are finite, and the normal and chi-square variables have every required moment, so the
expansion and expectation can be interchanged.

\begin{lemma}
\label{W:lem:gaussian-perturbation-index-count}
For fixed choices of diagonal or upper off-diagonal matrix units, drift positions $D$, and a partition $\pi$ of the remaining positions, the number of nonzero index assignments is
at most
\begin{equation}
 \frac{C_p^{M}}{(q!)^{2p}}\rho^{M+|D|+|\pi|}
   =\frac{C_p^{M}}{(q!)^{2p}}
           \rho^{2|D|+\sum_{B\in\pi}(|B|+1)}.
 \label{W:eq:gaussian-perturbation-index-count}
\end{equation}
Here the occurrences in each block of $\pi$ have the same coefficient index and type.
\end{lemma}
The order of a block in this expansion can grow with $\rho$. We therefore use the exact chi-square cumulants in \eqref{W:eq:gram-square-offdiagonal-cumulants-exact}, including
their factor $(r-1)!$ and the power $N^{1-r}$. A fixed order cumulant estimate for a linear statistic of the chi variables would not give the required dependence on $r$. Singleton
centered blocks vanish before absolute values are taken; distinct surviving blocks may have the same coefficient index. These are the partitions counted in the following proof.

\begin{proof}
Fix one of the $p$ factors $(\widetilde{\mathsf V}_{n,q})_x\overline{(\widetilde{\mathsf V}_{n,q})_x}$ in the expansion of $\|\widetilde{\mathsf C}_\rho^{(q)}\|_{\rm HS}^{2p}$. The
two products have a common column degree $n<{\rho}$ and the same row index $x\le\rho$. In either product, after the matrix units have been expanded, the deterministic kernel has
the form
\begin{equation}
 (\mathsf U_{s_0})_{x,\alpha_1}
 \prod_{h=1}^{q-1}
 (\mathsf U_{s_h})_{\beta_h,\alpha_{h+1}}
 (\mathsf U_{s_q})_{\beta_q,1},
 \label{W:eq:gaussian-perturbation-product}
\end{equation}
where $(\alpha_h,\beta_h)=(j_h,j_h)$ for a diagonal matrix unit and $(\alpha_h,\beta_h)=(j_h,j_h+1)$ for an upper off-diagonal matrix unit. The polynomial degrees satisfy
\(s_0+\cdots+s_q=n-q.\)

In Figure~\ref{W:fig:gaussian-paired-products}, a thin segment denotes an entry of $\mathsf U_r$ and a thick segment denotes a matrix unit. The two conjugate factors have the same
indices $(n,x)$. Their polynomial degrees are summed separately; the partition $\pi$ specifies which coefficient indices agree. \begingroup \setlength{\intextsep}{1.2\baselineskip}
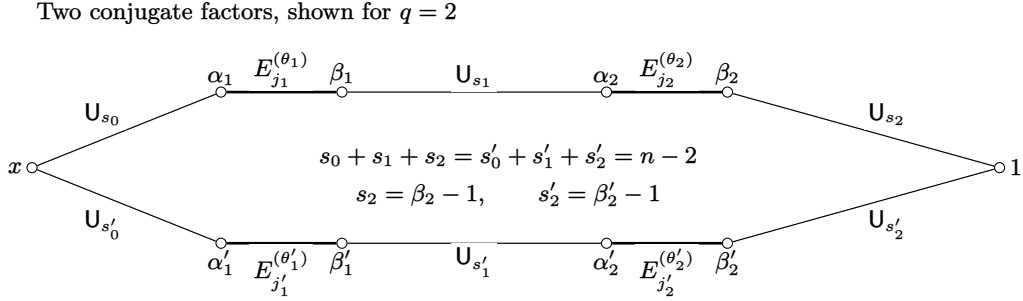
\begin{figure}[H]
\centering
\begin{tikzpicture}[x=1cm,y=1cm,every node/.style={mg text}]
 \node[anchor=west] at (0,2.05) {Two conjugate factors, shown for $q=2$};
 \node[mg site,label=left:{$x$}] (x) at (0,0) {};
 \node[mg site,label=right:{$1$}] (firstcoord) at (12.8,0) {};
 \foreach \yy/\ss in {1/u,-1/v}{
  \node[mg site] (a1\ss) at (2.5,\yy) {};
  \node[mg site] (b1\ss) at (4.1,\yy) {};
  \node[mg site] (a2\ss) at (7.6,\yy) {};
  \node[mg site] (b2\ss) at (9.2,\yy) {};
 }
 \draw[mg kernel] (x)--node[mg label,above left] {$\mathsf U_{s_0}$}(a1u);
 \draw[mg insertion] (a1u)--node[mg label,above] {$E_{j_1}^{(\theta_1)}$}(b1u);
 \draw[mg kernel] (b1u)--node[mg label,above] {$\mathsf U_{s_1}$}(a2u);
 \draw[mg insertion] (a2u)--node[mg label,above] {$E_{j_2}^{(\theta_2)}$}(b2u);
 \draw[mg kernel] (b2u)--node[mg label,above right] {$\mathsf U_{s_2}$}(firstcoord);
 \draw[mg kernel] (x)--node[mg label,below left] {$\mathsf U_{s'_0}$}(a1v);
 \draw[mg insertion] (a1v)--node[mg label,below] {$E_{j'_1}^{(\theta'_1)}$}(b1v);
 \draw[mg kernel] (b1v)--node[mg label,below] {$\mathsf U_{s'_1}$}(a2v);
 \draw[mg insertion] (a2v)--node[mg label,below] {$E_{j'_2}^{(\theta'_2)}$}(b2v);
 \draw[mg kernel] (b2v)--node[mg label,below right] {$\mathsf U_{s'_2}$}(firstcoord);
 \node[above=1pt] at (a1u) {$\alpha_1$};
 \node[above=1pt] at (b1u) {$\beta_1$};
 \node[above=1pt] at (a2u) {$\alpha_2$};
 \node[above=1pt] at (b2u) {$\beta_2$};
 \node[below=1pt] at (a1v) {$\alpha'_1$};
 \node[below=1pt] at (b1v) {$\beta'_1$};
 \node[below=1pt] at (a2v) {$\alpha'_2$};
 \node[below=1pt] at (b2v) {$\beta'_2$};
 \node at (6.3,.15) {$s_0+s_1+s_2=s'_0+s'_1+s'_2=n-2$};
 \node at (6.3,-.35) {$s_2=\beta_2-1,\qquad s'_2=\beta'_2-1$};
\end{tikzpicture}
\caption{The product in \eqref{W:eq:gaussian-perturbation-product} and its conjugate for $q=2$. Open circles denote matrix indices. Thin segments carry entries of $\mathsf U_r$, and thick segments carry $E_j^{(\theta)}$, with $\theta=0$ for a diagonal matrix unit and $\theta=1$ for an upper off-diagonal matrix unit. The factors act from right to left. The last factor imposes \eqref{W:eq:gaussian-last-gap}; for fixed coefficient indices and $n$, the remaining polynomial degrees have at most $\rho^{q-1}/(q-1)!$ choices.}
\label{W:fig:gaussian-paired-products}
\end{figure}
\endgroup

By \eqref{W:eq:free-chebyshev-first-column}, the last factor in \eqref{W:eq:gaussian-perturbation-product} is nonzero only when
\begin{equation}
 s_q=\beta_q-1.
 \label{W:eq:gaussian-last-gap}
\end{equation}
Thus, once the coefficient indices and $n$ are fixed, the number of possible $(s_0,\ldots,s_{q-1})$ is at most
\(\binom{n-s_q-1}{q-1} \le \frac{{\rho}^{q-1}}{(q-1)!},\) with the binomial coefficient interpreted as zero when its
upper index is smaller than $q-1$. For the conjugate product its last polynomial degree is $s'_q=\beta'_q-1$, and its nonnegative remaining polynomial degrees sum to $n-q-s'_q$.
Their number is $\binom{n-s'_q-1}{q-1}\le\rho^{q-1}/(q-1)!$. There are at most $\rho$ values of $n$ and $\rho$ of $x$, hence at most $\rho^2$ common pairs $(n,x)$. Therefore each
pair of conjugate factors contributes at most
\begin{equation}
 {\rho}^2\left(\frac{{\rho}^{q-1}}{(q-1)!}\right)^2
 =q^2\frac{{\rho}^{2q}}{(q!)^2}
 \le 4^q\frac{{\rho}^{2q}}{(q!)^2}
 \label{W:eq:gaussian-paired-index-count}
\end{equation}
possible choices of row, column and polynomial degrees. For $q=1$ the binomial coefficient is either zero or one, so this estimate includes that case. The other free kernels only
restrict the choices. Coinciding coefficient indices in different conjugate pairs can impose further constraints; discarding those constraints gives an upper bound. For fixed
coefficient indices, the Cartesian product of the $p$ sets of choices has cardinality at most
\[
 \left(4^q\frac{\rho^{2q}}{(q!)^2}\right)^p
 =2^M\frac{\rho^M}{(q!)^{2p}},
 \qquad M=2pq.
\]
Each drift position has one coefficient index, and each centered cumulant block has one common coefficient index. There are therefore at most $\rho^{|D|+|\pi|}$ choices of these
indices. Indices belonging to different blocks need not be distinct. There are at most $\rho$ choices for each index, even for an upper off-diagonal entry, whose actual range is
smaller. Thus the total number for the fixed choices of matrix units is bounded by
\[
 \frac{C_p^M}{(q!)^{2p}}\rho^{M+|D|+|\pi|}
 =\frac{C_p^M}{(q!)^{2p}}
       \rho^{\,2|D|+\sum_{B\in\pi}(|B|+1)},
 \qquad M=|D|+\sum_{B\in\pi}|B|.
\]
This proves \eqref{W:eq:gaussian-perturbation-index-count}. Independence sets a centered block to zero unless all its occurrences have the same coefficient index and the same value
of $\theta$; for diagonal Gaussian coefficients, $\kappa_r(a_j)=0$ whenever $r\ne2$.
\end{proof}

\Needspace{6\baselineskip}

Fix $p\in\mathbb N$, $T>0$, and $\rho^3\le C_TN$. Define
\begin{equation}
 w_2=C_T,\qquad
 w_r=(r-1)!C_T^r\rho^{-2(r-2)},\quad r\ge3.
 \label{W:eq:gram-effective-block-weights}
\end{equation}
Let $\mathcal P_{\ge2}([\ell])$ denote the partitions of $\{1,\ldots,\ell\}$ whose blocks have size at least two, and let $S_\ell$ be the sum on the left below.
\Needspace{8\baselineskip}

\begin{lemma}
\label{W:lem:gaussian-cumulant-partition-bound}
For $0\le\ell\le2p\rho$,
\begin{equation}
 \sum_{\pi\in\mathcal P_{\ge2}([\ell])}\prod_{B\in\pi}w_{|B|}
       \le C_{p,T}^{\ell}\ell^{\ell/2}.
 \label{W:eq:gram-partition-sum-bound}
\end{equation}
For $\ell=0$, both sides are interpreted as one.
\end{lemma}

\begin{proof}
For $\ell=0$ the empty partition contributes one, and for $\ell=1$ the sum is zero. Suppose $\ell\ge2$. A partition with $m_r$ blocks of size $r$, $2\le r\le\ell$, has
$\sum r m_r=\ell$. The number of such partitions is $\ell!/\prod_{r=2}^{\ell}(r!)^{m_r}m_r!$: order the $\ell$ occurrences, divide by the orders within the blocks, and then by
permutations of blocks of the same size. Consequently
\begin{equation}
 \frac{S_\ell}{\ell!}
 =\sum_{\sum r m_r=\ell}\prod_{r=2}^{\ell}
        \frac{(w_r/r!)^{m_r}}{m_r!}
 =[z^\ell]\exp H_\rho(z),\qquad
 H_\rho(z)=\sum_{r=2}^{\ell}\frac{w_r}{r!}z^r.
 \label{W:eq:gram-partition-EGF}
\end{equation}
Choose $\alpha=\min\{1,(2C_T\sqrt{2p})^{-1}\}$ and put $R=\alpha\sqrt\ell$. Then, for $\ell\le2p\rho$,
\[
 u:=\frac{C_TR}{\rho^2}
 \le C_T\alpha\sqrt{2p}\,\rho^{-3/2}\le\frac12.
\]
The factorial in $w_r$ leaves the denominator $r$:
\[
 \sum_{r=3}^{\ell}\frac{w_rR^r}{r!}
 =\rho^4\sum_{r=3}^{\ell}\frac{u^r}{r}
 \le\frac{\rho^4u^3}{1-u}
 \le\frac{2C_T^3R^3}{\rho^2}
 \le2C_T^3\alpha^3(2p)^{3/2}\rho^{-1/2}.
\]
Thus $H_\rho(R)\le C_TR^2/2+C_{p,T}\le C_{p,T}\ell$ uniformly in $\rho,\ell$ in the stated range. Every coefficient of $\exp H_\rho$ is nonnegative, so its coefficient of degree
$\ell$ is at most $R^{-\ell}\exp H_\rho(R)$. It follows that
\[
 S_\ell\le\ell!\alpha^{-\ell}\ell^{-\ell/2}e^{C_{p,T}\ell}
       \le C_{p,T}^{\ell}\ell^{\ell/2},
\]
using $\ell!\le\ell^\ell$. No constant depends on the maximum block size $\ell$.
\end{proof}

\Needspace{6\baselineskip}

\begin{lemma}
\label{W:lem:gaussian-perturbation-moments}
Fix $p\in\mathbb N$ and $T<\infty$. Uniformly for ${\rho}\le TN^{1/3}+r_0$ and $1\le q<{\rho}$,
\begin{equation}
 \mathbb E
 \|\widetilde{\mathsf C}_{\rho}^{(q)}\|_{\mathrm{HS}}^{2p}
 \le \frac{C_{p,T}^{2pq}}{(q!)^p}.
 \label{W:eq:gaussian-perturbation-moment-bound}
\end{equation}
\end{lemma}

\begin{proof}
Put $M=2pq$. For a fixed choice of the indices $\theta$, write $\xi_{j,0}=a_j$ and $\xi_{j,1}=\zeta_j$. After choosing the drift positions $D$, independence and the
moment--cumulant formula give
\[
 \mathbb E\prod_{v\notin D}\xi_{j_v,\theta_v}
 =\sum_{\pi\in\mathcal P_{\ge2}(D^c)}
       \prod_{B\in\pi}
       \left[\mathbf1_{\{(j_v,\theta_v)\text{ is constant on }B\}}
                  \kappa_{|B|}(\xi_{j_B,\theta_B})\right].
\]
If $D^c=\varnothing$, the partition sum is the empty product, equal to one. Singleton blocks vanish because every $\xi$ is centered. A block of diagonal Gaussian variables
contributes only when its size is two. Set $\delta_* =\max_{j\le\rho}|\delta_j|$ and $k_{r,\theta}=\max_{j\le\rho}|\kappa_r(\xi_{j,\theta})|$. For $d=|D|$, the index count of
Lemma~\ref{W:lem:gaussian-perturbation-index-count} bounds the absolute contribution of a fixed centered partition $\pi$ by
\begin{equation}
 \frac{C_p^M}{(q!)^{2p}}
 {\rho}^{2d+\sum_{B\in\pi}(|B|+1)}
 \delta_*^d\prod_{B\in\pi}k_{|B|,\theta_B}.
 \label{W:eq:gram-fixed-partition-contribution}
\end{equation}
Here $\theta_B$ is the common value of $\theta$ in the block $B$.

Figure~\ref{W:fig:moment-cumulant-blocks} illustrates one term of this partition sum. Each diamond specifies a cumulant of the coefficient occurrences joined to it. A drift
position contributes its deterministic coefficient separately. A Gaussian diagonal cumulant has order two; a cumulant of a squared off-diagonal coefficient can have any order at
least two, with the explicit value in \eqref{W:eq:gram-square-offdiagonal-cumulants-exact}. \begingroup
\setlength{\intextsep}{1.2\baselineskip}
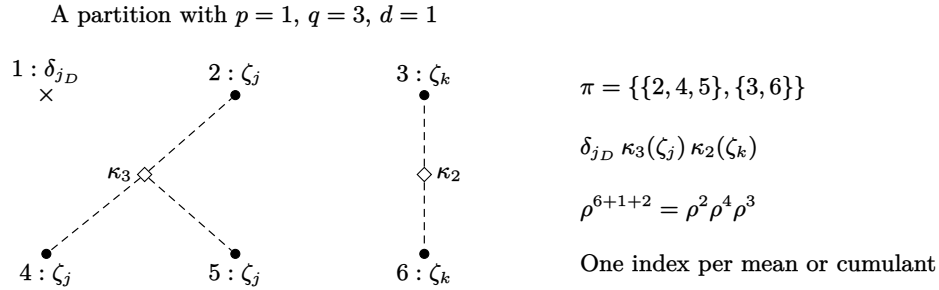
\begin{figure}[H]
\centering
\begin{tikzpicture}[x=1cm,y=1cm,every node/.style={mg text}]
 \node[anchor=west] at (0,1.05) {A partition with $p=1$, $q=3$, $d=1$};
 \node (v1) at (0,0) {$\times$};
 \node[above=3pt] at (v1) {$1:\delta_{j_D}$};
 \node[mg occurrence,label=above:{$2:\zeta_j$}] (v2) at (2.5,0) {};
 \node[mg occurrence,label=above:{$3:\zeta_k$}] (v3) at (5.0,0) {};
 \node[mg occurrence,label=below:{$4:\zeta_j$}] (v4) at (0,-2.1) {};
 \node[mg occurrence,label=below:{$5:\zeta_j$}] (v5) at (2.5,-2.1) {};
 \node[mg occurrence,label=below:{$6:\zeta_k$}] (v6) at (5.0,-2.1) {};
 \node[mg block,label=left:{$\kappa_3$}] (b3) at (1.3,-1.05) {};
 \node[mg block,label=right:{$\kappa_2$}] (b2) at (5.0,-1.05) {};
 \foreach \v in {v2,v4,v5}{\draw[mg incidence] (b3)--(\v);}
 \foreach \v in {v3,v6}{\draw[mg incidence] (b2)--(\v);}
 \node[anchor=west] at (7.0,.1) {$\pi=\{\{2,4,5\},\{3,6\}\}$};
 \node[anchor=west] at (7.0,-.7) {$\delta_{j_D}\,\kappa_3(\zeta_j)\,\kappa_2(\zeta_k)$};
 \node[anchor=west] at (7.0,-1.5) {$\rho^{6+1+2}=\rho^2\rho^4\rho^3$};
 \node[anchor=west] at (7.0,-2.25) {One index per mean or cumulant};
\end{tikzpicture}
\caption{A term of \eqref{W:eq:gram-fixed-partition-contribution} containing one deterministic mean, one third cumulant, and one second cumulant. The two rows show the occurrences in the conjugate products. The cross denotes the deterministic mean, filled dots denote centered variables, and each diamond joins the occurrences belonging to one cumulant. The dashed lines specify their common coefficient index. Formula~\eqref{W:eq:gram-square-offdiagonal-cumulants-exact} gives the third cumulant of $\zeta_j$.}
\label{W:fig:moment-cumulant-blocks}
\end{figure}
\endgroup For the term in Figure~\ref{W:fig:moment-cumulant-blocks},
\[
 [\rho^2|\delta_{j_D}|]\,
 [\rho^4|\kappa_3(\zeta_j)|]\,
 [\rho^3|\kappa_2(\zeta_k)|]
 \le C_{T,\beta}\frac{\rho^{10}}{N^4}
 \le C_{T,\beta}\rho^{-2}.
\]
Here $\rho^3/N\le C_T$, and the drift and cumulant bounds are \eqref{W:eq:gram-square-offdiagonal-cumulants-bound} and $|\delta_j|=j/N$. For a general partition, the corresponding
bounds are
\begin{equation}
 {\rho}^2\max_{j\le {\rho}}|\delta_j|
 \le \frac{{\rho}^3}{N}\le C_T.
 \label{W:eq:gram-drift-effective-weight}
\end{equation}
A Gaussian $a$-block has size two and effective weight at most $C{\rho}^3/N\le C_T$. By \eqref{W:eq:gram-square-offdiagonal-cumulants-bound}, a $\zeta$-block of size $r\ge2$ has
effective weight
\begin{equation}
 {\rho}^{r+1}|\kappa_r(\zeta_j)| \le (r-1)!C_\beta^r{\rho}^{r+1}N^{1-r} \le (r-1)!C_{T,\beta}^r{\rho}^{-2(r-2)}.
 \label{W:eq:gram-high-offdiagonal-block-gain}
\end{equation}
The power of $\rho$ in the last bound follows from the exact identity
\[
 \rho^{r+1}N^{1-r}
   =\left(\frac{\rho^3}{N}\right)^{r-1}\rho^{-2(r-2)}.
\]
Choose the constant in \eqref{W:eq:gram-effective-block-weights} large enough to dominate the size-two bound and $C_\beta^r((T+1)^3)^{r-1}$ at every $r\ge3$. This is possible with
one fixed constant raised to the power $r$. There are at most $2^M$ choices of the indices $\theta$. For each type assignment, the number of sets of $d$ drift positions is at most
$\binom Md$; counting choices that are incompatible with diagonal types only enlarges this bound. For each such set, $\ell=M-d\le2p\rho$ centered occurrences remain. The preceding
lemma applies without any fixed-order restriction and bounds their total partition weight by $C_{p,T}^{\ell}\ell^{\ell/2}$. Consequently, after absorbing the type count into
$C_{p,T}^M$,
\begin{equation}
 \begin{aligned}
 \mathbb E\|\widetilde{\mathsf C}_{\rho}^{(q)}\|_{\mathrm{HS}}^{2p}
 &\le\frac{C_{p,T}^{M}}{(q!)^{2p}}
      \sum_{d=0}^{M}\binom Md C_{p,T}^{d}(M-d)^{(M-d)/2}\\
 &\le\frac{C_{p,T}^{M}M^{M/2}}{(q!)^{2p}}.
 \end{aligned}
 \label{W:eq:gaussian-partition-moment-bound}
\end{equation}
In the second line we used $\ell^{\ell/2}\le M^{M/2}$ for $0\le\ell\le M$, with $0^0=1$, and $\sum_{d=0}^M\binom Md C^d=(1+C)^M$. Thus the constants still depend only on the fixed
parameters and $p$. Finally, $q!\ge(q/e)^q$, obtained by integrating $\log t$ on $[1,q]$, gives
\begin{equation}
 M^{M/2}=(2pq)^{pq}\le(2pe)^{pq}(q!)^p.
 \label{W:eq:gaussian-factorial-bound}
\end{equation}
Combining the last two displays and enlarging the fixed constant gives $C_{p,T}^{2pq}/(q!)^p$, as required. Centering eliminates the singleton terms before absolute values are
taken; the remaining terms are bounded by the displayed products of cumulants.
\end{proof}

\Needspace{6\baselineskip}

\begin{proposition}
\label{W:prop:gram-connection-HS}
For every fixed $0<p<\infty$ and $T<\infty$,
\begin{equation}
 \sup_{{\rho}\le TN^{1/3}+r_0}
 \left\|\|\widetilde{\mathsf C}_{\rho}-I_{\rho}\|_{\mathrm{HS}}\right\|_{L^{2p}}
 \le C_{p,T}.
 \label{W:eq:gram-connection-HS-bound}
\end{equation}
\end{proposition}

\begin{proof}
First let $p\ge1$ be an integer. The triangle inequality in the Hilbert--Schmidt norm and Minkowski's inequality in $L^{2p}$ apply to the finite identity
\eqref{W:eq:gaussian-coefficient-expansion}. They give
\[
 \left\|\|\widetilde{\mathsf C}_{\rho}-I_\rho\|_{\mathrm{HS}}\right\|_{L^{2p}}
 \le\sum_{q=1}^{\rho-1}
       \left(\mathbb E\|\widetilde{\mathsf C}_{\rho}^{(q)}\|_{\mathrm{HS}}^{2p}\right)^{1/(2p)}
 \le\sum_{q=1}^{\infty}\frac{C_{p,T}^{q}}{\sqrt{q!}}.
\]
The ratio of consecutive terms in the last numerical series is $C_{p,T}/\sqrt{q+1}$, which is at most $1/2$ for all sufficiently large $q$, independently of $N,\rho$. The finite
initial part and geometric tail therefore give a finite constant $C_{p,T}$. This bounds each $L^{2p}$ norm before the supremum over $\rho$ is taken. The infinite sum is a
deterministic upper bound for the finite sum above. For arbitrary $p>0$, choose an integer $m\ge\max\{1,p\}$. On a probability space, H\"older's inequality yields
$\|X\|_{L^{2p}}\le\|X\|_{L^{2m}}$. Applying the integer case proves the same moment bound, including $2p<1$. For $\rho=1$ the matrix difference is zero.
\end{proof}

\subsubsection{Bounds for the Gram matrix and its inverse}
\label{W:sec:gaussian-gram-singular-values}

We first show that the diagonal products converge to one. The normalized upper triangular factor has bounded Hilbert--Schmidt distance from the identity. Combining this bound with
its determinant gives the upper and lower bounds on singular values needed to control the Gram matrix and its inverse. Put
\begin{equation}
 L_k:=\log d_{k,N}=\sum_{j=1}^{k-1}\ell_j.
 \label{W:eq:gram-partial-log-products}
\end{equation}

\Needspace{6\baselineskip}

\begin{lemma}
\label{W:lem:gaussian-logarithmic-products}
For every fixed $0<p<\infty$ and ${\rho}\le TN^{1/3}+r_0$,
\begin{equation}
 \left\|\max_{k\le {\rho}}|L_k|\right\|_{L^p}
 \le C_{p,T}N^{-1/3}.
 \label{W:eq:gram-max-partial-log}
\end{equation}
Consequently,
\begin{equation}
 \|\mathsf D_{\rho}-I_{\rho}\|+\|\mathsf D_{\rho}^{-1}-I_{\rho}\|
 \longrightarrow0
 \qquad\hbox{in probability}.
 \label{W:eq:gaussian-diagonal-convergence}
\end{equation}
\end{lemma}

\begin{proof}
For $\rho=1$, $L_1=0$, so the quantity to be bounded is zero. Otherwise write $L_k=S_k+\sum_{j<k}\mathbb E\ell_j$, where $S_k=\sum_{j<k}(\ell_j-\mathbb E\ell_j)$. The deterministic
part obeys
\[
 \max_{k\le\rho}\left|\sum_{j<k}\mathbb E\ell_j\right|
 \le\frac{C_\beta}{N}\sum_{j=1}^{\rho-1}(j+1)
 \le C_\beta\frac{\rho^2}{N}\le C_{T,\beta}N^{-1/3}.
\]
The variables $\ell_j$ are independent functions of the independent $b_{j,N}^G$. Therefore $(S_k)_{1\le k\le\rho}$ is a centered martingale for
$\mathcal F_k=\sigma(\ell_1,\ldots,\ell_{k-1})$. Its log moment-generating function is the sum of the corresponding centered log moment-generating functions. For every fixed
integer $r\ge2$,
\[
 |\kappa_r(S_\rho)|
 \le(\rho-1)(r-1)!C_\beta^rN^{1-r}
 \le C_{r,\beta}(\rho/N)^{r/2},
\]
because $\rho N^{1-r}/(\rho/N)^{r/2}=(\rho N)^{1-r/2}\le1$. Choose an integer $m\ge1$ with $2m\ge p$. The scalar moment--cumulant formula, obtained by differentiating the
exponential of the cumulant-generating function, has no singleton blocks. Since $S_\rho$ is real,
\[
 \mathbb E|S_\rho|^{2m}
 \le\sum_{\pi\in\mathcal P_{\ge2}([2m])}
        \prod_{B\in\pi}|\kappa_{|B|}(S_\rho)|
 \le C_{m,\beta}(\rho/N)^m.
\]
The partition set is finite for this fixed $m$; the product of powers is $(\rho/N)^m$ in every term.

To bound the maximum, put $S^*=\max_{k\le\rho}|S_k|$. All $S_k$ have a finite $2m$th moment, so $\mathbb E(S^*)^{2m}\le\sum_k\mathbb E|S_k|^{2m}<\infty$. Let $\tau$ be the first
$k$ such that $|S_k|>t$, when that set is nonempty. Conditional Jensen's inequality gives $\mathbb E[|S_\rho|\mid\mathcal F_k]\ge|S_k|$. Multiplying by $\mathbf1_{\{\tau=k\}}$ and
summing over $k$ yields
\[
 t\,\mathbb P(S^*>t)
 \le\mathbb E\bigl[|S_\rho|\mathbf1_{\{S^*>t\}}\bigr].
\]
Multiply by $2m t^{2m-2}$ and integrate. Tonelli's theorem applies to the nonnegative integrand, and then H\"older's inequality gives
\[
 \mathbb E(S^*)^{2m}
 \le\frac{2m}{2m-1}\mathbb E\bigl[|S_\rho|(S^*)^{2m-1}\bigr]
 \le\frac{2m}{2m-1}\|S_\rho\|_{L^{2m}}\|S^*\|_{L^{2m}}^{2m-1}.
\]
If $\|S^*\|_{L^{2m}}>0$, divide by its $(2m-1)$st power. If it is zero, the required bound holds with left side zero. Thus, also for $p<1$,
\[
 \|S^*\|_{L^p}\le\|S^*\|_{L^{2m}}
 \le\frac{2m}{2m-1}\|S_\rho\|_{L^{2m}}
 \le C_{p,\beta}(\rho/N)^{1/2}\le C_{p,T}N^{-1/3}.
\]
Combining this with the deterministic part, using the $L^{2m}$ norm first when $p<1$, proves \eqref{W:eq:gram-max-partial-log}. Finally,
\[
 \|\mathsf D_\rho-I_\rho\|+\|\mathsf D_\rho^{-1}-I_\rho\|
 \le2\bigl(e^{\max_{k\le\rho}|L_k|}-1\bigr).
\]
For every $\eta>0$, Markov's inequality bounds the probability that the left side exceeds $\eta$ by $C_{p,T}N^{-p/3}/[\log(1+\eta/2)]^p$. The right side tends to zero uniformly
over the stated range of $\rho$, proving convergence in probability.
\end{proof}

\Needspace{6\baselineskip}

\begin{proof}[Proof of Lemma~\ref{W:lem:gaussian-gram-conditioning}]
Write $K_\rho=\widetilde{\mathsf C}_\rho-I_\rho$ and let $\widetilde s_1\ge\cdots\ge\widetilde s_\rho>0$ be the singular values of $\widetilde{\mathsf C}_\rho$. The matrix $K_\rho$
is strictly upper triangular. Hence
\[
 \operatorname{Tr}K_\rho=0,\qquad
 \sum_{j=1}^\rho\widetilde s_j^2
 =\operatorname{Tr}\bigl((I_\rho+K_\rho)^*(I_\rho+K_\rho)\bigr)
 =\rho+\|K_\rho\|_{\mathrm{HS}}^2,
 \qquad \prod_{j=1}^\rho\widetilde s_j
       =|\det\widetilde{\mathsf C}_\rho|=1.
\]
The norm triangle inequality and Proposition~\ref{W:prop:gram-connection-HS} imply
\begin{equation}
 \widetilde s_1\le1+\|K_\rho\|_{\mathrm{HS}}=O_{\mathbb P}(1).
 \label{W:eq:gram-largest-singular-tight}
\end{equation}
For $x>0$, the function $(x^2-1)/2-\log x$ has its minimum zero at $x=1$. Applying this scalar inequality to the first $\rho-1$ singular values gives
\begin{equation}
 \begin{aligned}
 -\log\widetilde s_\rho
 &=\sum_{j<\rho}\log\widetilde s_j
 \le\frac12\sum_{j<\rho}(\widetilde s_j^2-1)\\
 &=\frac12\bigl(\|K_\rho\|_{\mathrm{HS}}^2-(\widetilde s_\rho^2-1)\bigr)
 \le\frac12\bigl(\|K_\rho\|_{\mathrm{HS}}^2+1\bigr).
 \end{aligned}
 \label{W:eq:gram-smallest-singular-bound}
\end{equation}
For $\rho=1$ the sum is empty and the same inequalities hold. In particular, $\widetilde s_\rho^{-1}\le\exp((\|K_\rho\|_{\mathrm{HS}}^2+1)/2)$ pathwise.

Put $\ell_\rho=\max_{k\le\rho}|L_k|$. The exact factorization $\mathsf C_\rho=\mathsf D_\rho\widetilde{\mathsf C}_\rho$ implies, for every vector $v$,
\[
 \|\mathsf C_\rho v\|
 \le e^{\ell_\rho}\widetilde s_1\|v\|,
 \qquad
 \|\mathsf C_\rho v\|
 \ge e^{-\ell_\rho}\widetilde s_\rho\|v\|.
\]
Thus, on the single event $\{\|K_\rho\|_{\mathrm{HS}}\le M,\ \ell_\rho\le\delta\}$,
\[
 \lambda_{\max}(G_\rho)\le e^{2\delta}(1+M)^2,
 \qquad
 \lambda_{\min}(G_\rho)^{-1}\le e^{2\delta+M^2+1}.
\]
For any fixed $p>0$, the union bound and the two proved moment estimates show that the probability of its complement is at most
\[
 \mathbb P(\|K_\rho\|_{\mathrm{HS}}>M)+\mathbb P(\ell_\rho>\delta)
 \le C_{p,T}M^{-2p}+C_{p,T}\delta^{-p}N^{-p/3}.
\]
The two probabilities are added by the union bound under the joint distribution of $K_\rho$ and $\ell_\rho$. Given $\varepsilon>0$, set $\delta=1$, choose $M$ so that the first
term is at most $\varepsilon/2$, and then choose $N_1$ so that the second is at most $\varepsilon/2$ for every $N\ge N_1$. The two deterministic bounds on the event then prove
tightness uniformly in $1\le\rho\le TN^{1/3}+r_0$. The finitely many permitted sizes below $N_1$ have positive Gram matrices almost surely, so enlarging the threshold handles them
as well. This argument controls probabilities of inverse singular values; it does not integrate the exponential bound above.

For each deterministic unit $w\in\mathbb F_\beta^N$, the Lanczos coefficient law is the same product law used throughout the proof. Equivalently, choosing a deterministic
orthogonal or unitary $V$ with $Ve_1=w$ gives
\[
 G_\rho(H_N^{\mathrm G};w)
   =G_\rho(V^*H_N^{\mathrm G}V;e_1)
   \ \stackrel{\mathrm d}=\ G_\rho(H_N^{\mathrm G};e_1).
\]
The constants and the preceding probability estimates are therefore independent of the deterministic $w$. This is the uniform bound in the lemma.
\end{proof}

For later choices of cutoffs, let $\varepsilon_N\in(0,1)$ be any deterministic sequence tending to zero and define
\begin{equation}
 A_N=\exp\sqrt{\min\{\log(N+2),|\log\varepsilon_N|\}}.
 \label{W:eq:gaussian-gram-threshold}
\end{equation}
Both quantities inside the minimum diverge. Moreover,
\[
 \frac{\log A_N}{\log N}\longrightarrow0,\qquad
 \frac{\log A_N}{|\log\varepsilon_N|}\longrightarrow0.
\]
Indeed, the first ratio is at most $\sqrt{\log(N+2)}/\log N$, and the second is at most $|\log\varepsilon_N|^{-1/2}$. Thus $A_N=N^{o(1)}$ and, for every fixed $c>0$,
\[
 A_N^c\varepsilon_N
 =\exp\bigl(c\log A_N-|\log\varepsilon_N|\bigr)\longrightarrow0.
\]
Set $Z_{N,\rho,w}=\max\{\lambda_{\max}(G_\rho(H_N^{\mathrm G};w)), \lambda_{\min}(G_\rho(H_N^{\mathrm G};w))^{-1}\}$. The preceding proof gives uniform tightness in $\rho,w$. For
any fixed $M$, choose $N_M$ with $A_N\ge M$ for $N\ge N_M$, and hence
\[
 \limsup_{N\to\infty}\sup_{\rho,w}\mathbb P(Z_{N,\rho,w}>A_N)
 \le\limsup_{N\to\infty}\sup_{\rho,w}\mathbb P(Z_{N,\rho,w}>M).
\]
Letting $M\to\infty$ proves, in particular, for every allowed deterministic sequence $\rho_N$,
\begin{equation}
 \mathbb P\left(A_N^{-1}I_{\rho_N}\le G_{\rho_N}(H_N^{\mathrm G})
                              \le A_NI_{\rho_N}\right)\longrightarrow1,
 \label{W:eq:gaussian-gram-spectral-event}
\end{equation}
uniformly over deterministic starting vectors. Here each probability is taken for a vector fixed before sampling; there is no probability of a supremum over matrix-dependent
starting vectors.

\subsubsection{Polynomial coefficient bounds for Gaussian and bounded entries}
\label{W:sec:polynomial-coefficients-for-comparison}
The conditioning estimate above concerns the GOE and GUE. We now choose one deterministic coefficient threshold for both the Gaussian and bounded entry matrices. Here an isotropic
local law is available because both ensembles have uniformly bounded standardized moments of every order. The original truncated Wigner entries are handled by the conditional
resolvent estimates in Section~\ref{W:sec:truncation-and-resolvent-estimates}.

Let $H_N^{\rm b}$ be the matrix with bounded standardized entries constructed in Lemma~\ref{W:lem:cubic-defect}; its off-diagonal third moments match those specified there, and its
standardized diagonal entries are independent symmetric signs. Fix $c_1m_N\le\rho\le C_1m_N$, use $s=4\max\{\rho^{-1},\sqrt d\}$, and let
$\mathcal C=\{z=\zeta+\zeta^{-1}:|\zeta|=e^s\}$. This ellipse uses the logarithmic window fixed earlier: $d\asymp m_N^{-2}(\log N)^\gamma$ with $1<\gamma<2$. In particular, for
large $N$,
\[
 s\asymp m_N^{-1}(\log N)^{\gamma/2},\qquad
 \rho s=O((\log N)^{\gamma/2})=o(\log N),\qquad
 e^{2\rho s}(1+|\log s|)=N^{o(1)}.
\]
The last expression bounds the mass of the contour measure in \eqref{W:eq:comparison-measure}. The two restrictions on $\gamma$ ensure that the window contains the confinement
scale and that the polynomials grow subpolynomially on the ellipse. The cutoff $\Gamma_\rho$, coefficient norms $Q_*$, and variation of the off-diagonal coefficients $B$ are
defined in Section~\ref{W:sec:coefficient-localization}. The estimates below apply to each of these two laws.

\Needspace{6\baselineskip} For $Y_N=H_N^{\rm b}$ or $H_N^{\mathrm G}$, use the same spectral modification as in the definition of $\chi_N$ and put
\[
 F_\rho^{[d]}=[p_0(h_d(Y_N))w,\ldots,p_{\rho-1}(h_d(Y_N))w],\qquad
 G_\rho^{[d]}=(F_\rho^{[d]})^*F_\rho^{[d]},\qquad
 \Delta_N=\max_{i,j<\rho}|(G_\rho^{[d]}-I)_{ij}|.
\]
Thus $\Delta_N$ and $\chi_N$ are computed from the same polynomial matrix. On the confinement event $h_d(Y_N)=Y_N$, it is the unmodified Gram matrix.

\begin{lemma}
\label{W:lem:comparison-ensemble-coefficients}
Uniformly over deterministic unit starting vectors,
\begin{equation}
 \Delta_N\prec\rho^{-1},\qquad
 \mathbf1_{\{\Gamma_\rho>0\}}(Q_*+B)
       \le C_\kappa(1+\rho\Delta_N+\rho^2\Delta_N^2).
 \label{W:eq:comparison-ensemble-regularity}
\end{equation}
Consequently one can choose the same deterministic $\mathfrak r_N\to\infty$, $\mathfrak r_N=N^{o(1)}$, for these two ensembles so that for every fixed $D$,
\begin{equation}
 \sup_{\|w\|=1}\mathbb E D_N^{\rm coef}(Y_N,w)\le N^{-D}
 \label{W:eq:comparison-ensemble-cutoff-error}
\end{equation}
for all sufficiently large $N$.
\end{lemma}
\begin{proof}
The matrix $H_N^{\rm b}$ has centered independent entries on and above the diagonal, variance $N^{-1}$ in every matrix coordinate, and uniformly bounded standardized moments of
every order. Its variance row sums are exactly one because of the chosen diagonal signs. Thus the variance conditions (2.22) and the uniform moment condition (2.23) of
\cite{bloemendal2014isotropic} hold, so its Theorem~2.12 applies. Here a complex off-diagonal entry is one independent block. With $R_Y(z)=(zI-Y)^{-1}$, let
$g(z)=\int(z-x)^{-1}\,\dd\mu_{\rm sc}(x)$ and $\upsilon(z)=|\operatorname{Im}g(z)|$. The cited estimate is
\[
 |w^*R_Y(z)w-g(z)|\prec
 \sqrt{\frac{\upsilon(z)}{N|\operatorname{Im}z|}}
             +\frac{1}{N|\operatorname{Im}z|}.
\]
For $\zeta=e^{s+\ii\theta}$ one has $z=2\cosh s\cos\theta+2\ii\sinh s\sin\theta$ and $g(z)=\zeta^{-1}$. Therefore, when
$\sin\theta\ne0$,
\[
 \frac{\upsilon(z)}{|\operatorname{Im}z|}
 =\frac{e^{-s}|\sin\theta|}{2\sinh s|\sin\theta|}
 =\frac{e^{-s}}{2\sinh s}\le\frac1{2s}.
\]
On the part of $\mathcal C$ where $|\operatorname{Im}z|\ge cs^2$, the two local law terms are bounded by
\[
 C(Ns)^{-1/2}+C(Ns^2)^{-1}
 \le C\sqrt{\rho/N}+C\rho^2/N\le C m_N^{-1}.
\]
Here $s\ge4/\rho$ and $\rho\asymp m_N$. The condition $|\operatorname{Im}z|\ge cN^{-2/3}$ places these points in the domain $\eta\ge N^{-1+\omega}$ for any fixed $0<\omega<1/3$,
for large $N$. The lower half of the contour follows by complex conjugation.

For every fixed $\varepsilon,D>0$, the local law has the uniform probability bound
\[
 \sup_{\|w\|=1}\mathbb P\left(
 |w^*R_Y(z)w-g(z)|
 >N^\varepsilon\left[
 \sqrt{\frac{\upsilon(z)}{N|\operatorname{Im}z|}}
 +\frac1{N|\operatorname{Im}z|}
 \right]\right)\le N^{-D}.
\]
Here each vector inside the probability is deterministic. At the tips, the real spectral gap has logarithmic size and does not meet the $N^{-2/3+\omega}$ condition of the cited
estimate outside the spectrum. We use the following comparison with positive imaginary part.

The upper tail estimate \cite{erdosxu2023small}, applied also to $-Y_N$, gives, for every fixed $D$, a constant $C_D$ such that
\[
 \mathbb P\bigl(\|Y_N\|>2+C_DN^{-2/3}(\log N)^{2/3}\bigr)\le N^{-D}
\]
for all sufficiently large $N$. Indeed, the centered entries of $H_N^{\rm b}$ have off-diagonal variance $N^{-1}$, common diagonal variance $N^{-1}$, and uniformly bounded
standardized moments of every order. In the complex case $\mathbb E(H_{ij}^{\rm b})^2=0$ because the two real coordinates have equal variances and zero covariance. These are the
moment assumptions of the cited theorem. Its allowed right tail range is $1\le x\le M(\log N)^{2/3}$. Choose $A_D$ so that $2\beta A_D^{3/2}/3>D+2$ and set $M=A_D$. Applying the
theorem to $Y_N$ and $-Y_N$ at $x=A_D(\log N)^{2/3}$ gives
\[
 \mathbb P\bigl(\|Y_N\|>2+A_DN^{-2/3}(\log N)^{2/3}\bigr)
 \le 2C_{A_D}(\log N)^{-\beta/2}
                   N^{-2\beta A_D^{3/2}/3}
 \le N^{-D}
\]
for large $N$. Since $\gamma>1$,
\[
 \frac{A_DN^{-2/3}(\log N)^{2/3}}d
 \le C A_D(\log N)^{2/3-\gamma}\longrightarrow0.
\]
The confinement interval is consequently contained in $[-2-d/8,2+d/8]$ for large $N$. At a real tip write $z=E+\mathrm i y$, where
$\operatorname{dist}(E,\operatorname{spec}Y)\ge c s^2$ and $|y|\le c's^2$. Take $h=c's^2$ with a fixed suitable constant. At $z_h=E+\mathrm i h$ the local law and Ward's identity
give
\[
 \|R_Y(z_h)w\|^2
 =\frac{-\operatorname{Im}(w^*R_Y(z_h)w)}h
 \prec \frac{Cs+m_N^{-1}N^{o(1)}}{c's^2}
 \le m_NN^{o(1)}.
\]
For $\mu_w^Y$ the positive spectral measure, the comparison is explicitly
\[
 \int\frac{\mu_w^Y(\mathrm d\lambda)}{(E-\lambda)^2+y^2}
 \le C\int\frac{\mu_w^Y(\mathrm d\lambda)}{(E-\lambda)^2+h^2}.
\]
Lemma~\ref{W:lem:vertical-resolvent-comparison} now gives
\[
 |w^*(R_Y(z)-R_Y(z_h))w|\prec
 Cs^2m_NN^{o(1)}=m_N^{-1}N^{o(1)}.
\]
The formula $g(z)=(z-\sqrt{z^2-4})/2$ gives $|g'(z)|\le C/s$ along this vertical segment, and hence $|g(z)-g(z_h)|\le Cs=m_N^{-1}N^{o(1)}$. To make these bounds simultaneous, take
a finite covering set on the contour of spacing $N^{-10}$. It has at most $CN^{10}$ points. On the confinement event the contour has spectral distance at least $cs^2\ge cN^{-2/3}$;
its matrix elements and column energies have Lipschitz constants at most $CN^{4/3}$ and $CN^2$, respectively, by the same lemma. The errors of the finite discretization are
$O(N^{-8})$, smaller than the required bounds. For a failure exponent $D$, use the pointwise local law with exponent $D+12$ before taking the union. These estimates have constants
independent of the deterministic unit vector $w$. Confinement also implies $h_d(Y)=Y$.

The real Gaussian ensemble requires a variance normalization before the cited local law applies. Its variance row sum is $1+(2/\beta-1)/N$. Put $c_N=(1+(2/\beta-1)/N)^{-1/2}$ and
$Z_N=c_N{H_N^{\mathrm G}}$. The cited local law applies to $Z_N$. Since $R_{{H_N^{\mathrm G}}}(z)=c_NR_{Z_N}(c_Nz)$ and $c_N=1+O(N^{-1})$, the deterministic change in
$c_Ng(c_Nz)-g(z)$ is $O(N^{-1}s^{-1})={m_N}^{-2}N^{o(1)}$, smaller than the required error. Indeed $|c_Nz-z|\le C/N=o(s^2)$, so the scaled contour has the same bounds on its
distance from the spectrum and on the semicircle resolvent, after a fixed enlargement of the constants. Gaussian confinement follows from \cite{erdosxu2023small}, applied to
${H_N^{\mathrm G}}$ and $-{H_N^{\mathrm G}}$. This rescaling makes every GOE variance row sum equal to one.

The resulting estimate is entrywise. Even when $\Delta_N\prec\rho^{-1}$, the bound $\|G_\rho^{[d]}-I\|\le\rho\Delta_N$ does not give convergence to zero. On the support of
$\Gamma_\rho$, the eigenvalues of the Gram matrix lie in $[\kappa,\kappa^{-1}]$. Hence the inverse Gram matrices below have operator norm at most $\kappa^{-1}$.

On the confinement event, Cauchy's formula yields
\[
 (G_\rho^{[d]}(Y,w)-I)_{ij}
   =\frac{1}{2\pi\mathrm i}\int_{\mathcal C}p_i(z)p_j(z)
                         (w^*R_Y(z)w-g(z))\,\mathrm dz.
\]
Its integral weight is $N^{o(1)}$ by \eqref{W:eq:comparison-measure}. The exceptional events for the bounded and Gaussian ensembles have arbitrary polynomially small probabilities.
This proves $\Delta_N\prec\rho^{-1}$ for the spectrally modified Gram matrix too.

For the deterministic bound, abbreviate $G=G_\rho^{[d]}$. On $\{\Gamma_\rho>0\}$, let $z_j=(G_{ij})_{i<j}$ and $\mathfrak h_j=G_{jj}-z_j^{\mathsf T}G_j^{-1}z_j$. The variational
formula for the Schur complement gives $\kappa\le\mathfrak h_j\le\kappa^{-1}$ on this event. Since $\|z_j\|_2^2\le j\Delta_N^2$ and
$\|G_j^{-1}\|\le\kappa^{-1}$,
\[
 |\mathfrak h_j-1|
 \le|G_{jj}-1|+\|G_j^{-1}\|\|z_j\|_2^2
 \le\Delta_N+\kappa^{-1}\rho\Delta_N^2.
\]
The identity $G_j^{-1}z_j=z_j-(G_j-I)G_j^{-1}z_j$ gives coordinatewise
\[
 |(G_j^{-1}z_j)_i|
 \le\Delta_N+\|(G_j-I)_{i\cdot}\|_2\|G_j^{-1}\|\|z_j\|_2
 \le\Delta_N+\kappa^{-1}\rho\Delta_N^2.
\]
Put $\epsilon=\Delta_N+\kappa^{-1}\rho\Delta_N^2$. The polynomial
\[
 q_j=\mathfrak h_j^{-1/2}
           \left(p_j-\sum_{i<j}(G_j^{-1}z_j)_ip_i\right)
\]
has leading coefficient at most $\kappa^{-1/2}$ and each preceding coefficient has absolute value at most $\kappa^{-1/2}\epsilon$. Consequently both a row and a column of its
coefficient matrix satisfy
\[
 \sum_i|C_{ji}|\le\kappa^{-1/2}(1+\rho\epsilon),\qquad
 \sum_j|C_{ji}|\le\kappa^{-1/2}(1+\rho\epsilon).
\]
Using $\mathfrak h_0=1$ and $b_j^2=\mathfrak h_j/\mathfrak h_{j-1}$,
\[
 |b_j-1|=
 \frac{|\mathfrak h_j-\mathfrak h_{j-1}|}
      {\mathfrak h_{j-1}(\sqrt{\mathfrak h_j/\mathfrak h_{j-1}}+1)}
 \le2\kappa^{-1}\epsilon.
\]
It follows that
\[
 B=1+\sum_{j=1}^{\rho-3}|b_{j+1}-b_j|
       \le1+2\sum_{j=1}^{\rho-2}|b_j-1|
       \le1+4\kappa^{-1}\rho\epsilon.
\]
Combining these estimates proves $Q_*+B\le C_\kappa(1+\rho\Delta_N+\rho^2\Delta_N^2)$ on the cutoff support.

The bounds $\mathcal Q\le Q_*+2$ and $\mathcal B\le B+1$ imply on $\{\Gamma_\rho>0\}$ that
\[
 \mathcal Q+\mathcal B\le C_\kappa
              (1+\rho\Delta_N+(\rho\Delta_N)^2).
\]
For fixed $\varepsilon>0$ and large $N$, the event that the left side exceeds $N^\varepsilon$ is contained in $\{\rho\Delta_N>N^{\varepsilon/3}\}$. Since $\Delta_N\prec\rho^{-1}$
uniformly in $w$, for every fixed $D>0$,
\[
 \sup_{\|w\|=1}\mathbb P\{\Gamma_\rho>0,
                          \mathcal Q+\mathcal B>N^\varepsilon\}
       \le N^{-D}
\]
for all sufficiently large $N$, under both matrix laws. Choose increasing integers $N_k\ge\mathrm e^{k^2}$ so that this statement, with $\varepsilon=1/k$ and failure exponent $k$,
holds for every $N\ge N_k$ under both ensemble laws. For $k(N)=\max\{k:N\ge N_k\}$ set $\mathfrak r_N=4N^{1/k(N)}$, with any fixed positive choice for the finitely many smaller
$N$. Then $k(N)\to\infty$, $\mathfrak r_N=N^{o(1)}$, and $\mathfrak r_N\ge4\mathrm e^{k(N)}\to\infty$. By the definition of $\chi_N$,
$\{D_N^{\rm coef}>0\}\subset \{\Gamma_\rho>0,\mathcal Q+\mathcal B>\mathfrak r_N\}$. Hence for $k=k(N)$,
\[
 \sup_{\|w\|=1}\mathbb ED_N^{\rm coef}(Y_N,w)
 \le\sup_{\|w\|=1}\mathbb P\{\Gamma_\rho>0,
                             \mathcal Q+\mathcal B>N^{1/k}\}
 \le N^{-k}.
\]
For every fixed $D$, $k(N)\ge D$ whenever $N\ge N_{\lceil D\rceil}$, which proves \eqref{W:eq:comparison-ensemble-cutoff-error}.
\end{proof}
\subsubsection{The operator limit for the Gaussian ensembles}
\label{W:sec:gaussian-operator-limit}

For the Gaussian matrices, the operator convergence results in \cite{ramirez2011beta} and \cite{krishnapur2016universality} also imply convergence of every fixed ordered eigenpair.
We state this result with our embedding and scaling, then prove the moment bounds needed to estimate expectations of spectral traces. The coupling preserves each finite Gaussian
marginal and specifies the joint law across sizes. The almost sure convergence is with respect to this coupling. The finite size probability and expectation estimates below do not
depend on a choice of such a joint law.

For the operator limit, use the isometry $\mathcal E_N$ defined in \eqref{W:eq:piecewise-constant-embedding} and set
\begin{align}
 \mathcal A_N^G&=N^{2/3}(2I-T_N^G),
       \label{W:eq:gauss-sao-scaled-operator}\\
 \mathfrak r_N^G(z)&=\mathcal E_N(\mathcal A_N^G-z)^{-1}\mathcal E_N^*,
       \qquad z\notin\mathbb R.
       \label{W:eq:gauss-sao-embedded-resolvent}
\end{align}
The limiting operator is the stochastic Airy operator with Dirichlet boundary condition, whose form is
\begin{equation}
 \mathfrak q_\beta[f]=\int_0^\infty(|f'|^2+x|f|^2)\,\dd x
       -\frac2{\sqrt\beta}\int_0^\infty B(x)(|f(x)|^2)'\,\dd x.
 \label{W:eq:gauss-sao-limit-form}
\end{equation}
Its form domain and closure are described in Section~\ref{W:sec:convergence-of-tridiagonal-operators}. Almost surely it has compact resolvent and simple eigenvalues
$\Lambda_1<\Lambda_2<\cdots$. Write $\lambda_{1,N}^G\le\cdots\le\lambda_{N,N}^G$ for the eigenvalues of $\mathcal A_N^G$, and choose normalized eigenvectors $v_{j,N}^G$ and
$\phi_j$. The allowed phases form $\mathbb S_{\mathbb F_1}=\{-1,1\}$ and $\mathbb S_{\mathbb F_2}=\{\omega\in\mathbb C:|\omega|=1\}$.

\begin{theorem}[Ram\'irez--Rider--Vir\'ag; Krishnapur--Rider--Vir\'ag]
\label{W:thm:external-gaussian-sao}
There is a coupling of the Lanczos tridiagonalizations of the GOE or GUE with a standard Brownian motion such that, almost surely, for every compact
$K\Subset\mathbb C\setminus\mathbb R$,
\begin{equation}
 \sup_{z\in K}\|\mathfrak r_N^G(z)-(\mathcal H_\beta-z)^{-1}\|
       \longrightarrow0.
 \label{W:eq:gauss-sao-local-uniform-resolvent}
\end{equation}
On the same coupling, every fixed ordered eigenpair satisfies
\begin{equation}
 \lambda_{j,N}^G\longrightarrow\Lambda_j,\qquad
 \inf_{\omega\in\mathbb S_{\mathbb F_\beta}}
       \|\mathcal E_Nv_{j,N}^G-\omega\phi_j\|_2\longrightarrow0.
 \label{W:eq:gauss-sao-eigenpair-convergence}
\end{equation}
\end{theorem}

\begin{proof}
By Theorem~\ref{W:thm:Gaussian-Lanczos-law},
\begin{equation}
 T_N^G\stackrel{\mathrm d}=H_N^{(\beta)},\qquad
 H_N^{(\beta)}:=\frac1{\sqrt{\beta N}}
 \begin{pmatrix}
 g_1&\chi_{\beta(N-1)}&&\\
 \chi_{\beta(N-1)}&g_2&\ddots&\\
 &\ddots&\ddots&\chi_\beta\\
 &&\chi_\beta&g_N
 \end{pmatrix}.
 \label{W:eq:gauss-sao-law-identification}
\end{equation}
This is the matrix in \cite{krishnapur2016universality}, and it equals $N^{-1/2}$ times the Gaussian tridiagonal matrix in \cite{ramirez2011beta}. Thus $N^{2/3}(2I-H_N^{(\beta)})$
is exactly the scaled tridiagonal operator used in those papers.

The constants in Remark~4 of \cite{krishnapur2016universality} are fixed by $V(x)=x^2/4$, whose limiting recurrence coefficients are $a^\dagger(x)=0$ and $b^\dagger(x)=\sqrt{1-x}$.
Consequently,
\begin{equation}
 \tau=-(a^\dagger)'(0)-2(b^\dagger)'(0)=1,
 \qquad \mathcal E_{\mathrm{KRV}}=2,
 \qquad \gamma_{\mathrm{KRV}}=\vartheta_{\mathrm{KRV}}=1.
 \label{W:eq:gauss-sao-normalization-constants}
\end{equation}
The coordinate embedding in the cited theorem is therefore exactly
\begin{equation}
 e_{k,N}^{\mathrm{KRV}}
 =N^{1/6}\mathbf1_{[(k-1)N^{-1/3},kN^{-1/3})}
 =\mathcal E_Ne_k.
 \label{W:eq:gauss-sao-embedding-identification}
\end{equation}
The associated stochastic Airy operator has Dirichlet boundary condition and the closed quadratic form \eqref{W:eq:gauss-sao-limit-form}: on the compactly supported form core,
$\langle B',|f|^2\rangle=-\int B(|f|^2)'$, and the closure is the one constructed in \cite{ramirez2011beta}. Theorem~1 of \cite{krishnapur2016universality} applies to the uniformly
convex polynomial $V(x)=x^2/4$ and implies norm resolvent convergence and \eqref{W:eq:gauss-sao-eigenpair-convergence} on one probability space. The simple limiting eigenvalues
determine each normalized eigenfunction up to the phase displayed there.

The finite resolvents are extended by zero on $(\operatorname{Ran}\mathcal E_N)^\perp$, as in that operator convergence theorem. Its convention $(z-A)^{-1}$ differs from ours by a
minus sign. For local uniformity, fix $z_0=\mathrm i$ and set $R_N(z)=\mathfrak r_N^G(z)$ and $R(z)=(\mathcal H_\beta-z)^{-1}$. Since $\mathcal E_N^*\mathcal E_N=I$, the usual
resolvent identity also holds for these zero extensions and gives
\begin{equation}
 R_N(z)-R(z)
 =\bigl[I+(z-z_0)R_N(z)\bigr]
   \bigl[R_N(z_0)-R(z_0)\bigr]
   \bigl[I+(z-z_0)R(z)\bigr].
 \label{W:eq:gauss-sao-uniform-resolvent-identity}
\end{equation}
For $d_K=\inf_{z\in K}|\operatorname{Im}z|>0$ and $M_K=\sup_{z\in K}|z-z_0|$, it follows that
\[
 \sup_{z\in K}\|R_N(z)-R(z)\|
 \le(1+M_K/d_K)^2\|R_N(z_0)-R(z_0)\|.
\]
The cited norm resolvent convergence at $z_0$ therefore gives locally uniform convergence on the same coupling.
\end{proof}

The preceding eigenpair limit already implies
\begin{equation}
 \lim_{L\to\infty}\limsup_{N\to\infty}\sup_{\|w_N\|=1}
 \mathbb P\bigl(\|Q_{N,L}\Pi_N^G(E)\|>\eta\bigr)=0,
 \quad
 \Pi_N^G(E):=\mathbf1_{(-\infty,E]}(\mathcal A_N^G),
 \label{W:eq:gauss-sao-qualitative-localization}
\end{equation}
for fixed $E\in\mathbb R$ and $\eta>0$. Here $Q_{N,L}$ is the coordinate tail projection of the introduction. Indeed, on
$\{\lambda_{r+1,N}^G>E\}$,
\[
 \|Q_{N,L}\Pi_N^G(E)\|^2
 \le\sum_{j=1}^r\|Q_{N,L}v_{j,N}^G\|_2^2.
\]
Put $L_N=h_N\lfloor m_NL\rfloor$. For fixed $j$ choose a phase $\omega_{j,N}$ attaining the infimum in \eqref{W:eq:gauss-sao-eigenpair-convergence}. Then
\[
 \left|\|Q_{N,L}v_{j,N}^G\|_2-
               \|\mathbf1_{[L,\infty)}\phi_j\|_2\right|
 \le\|\mathcal E_Nv_{j,N}^G-\omega_{j,N}\phi_j\|_2
       +\left(\int_{L_N}^{L}|\phi_j(x)|^2\,\dd x\right)^{1/2}
 \longrightarrow0.
\]
The second term tends to zero by absolute continuity of the integral. The eigenvalue convergence and this estimate imply
\[
 \limsup_N\mathbb P\bigl(\|Q_{N,L}\Pi_N^G(E)\|>\eta\bigr)
 \le \mathbb P(\Lambda_{r+1}\le E)
 +\mathbb P\left(\sum_{j=1}^r
       \|\mathbf1_{[L,\infty)}\phi_j\|_2^2\ge\eta^2\right).
\]
First let $L\to\infty$, using $\phi_j\in L^2$, and then let $r\to\infty$, using $\Lambda_r\to\infty$. The distribution is independent of the deterministic starting vector, so the
bound is uniform in $w_N$. The expectation estimates used in entry comparison are proved in Sections~\ref{W:sec:coercivity-and-eigenvalue-counting}
and~\ref{W:sec:expectation-of-the-localization-statistic}.

\subsubsection{Coercivity and eigenvalue counting}
\label{W:sec:coercivity-and-eigenvalue-counting}
\label{W:subsec:gram-Gaussian-global-Weyl}

We bound the expectation of $\operatorname{Tr}\theta_E(\mathcal A_N^{\mathrm G})$ by deriving a coercive form estimate and then bounding the eigenvalue count. The
finite-dimensional estimates include the contribution of the exceptional event. We use oscillations of the integrated coefficients on intervals of length two, as in the form
estimates of \cite{ramirez2011beta}, and keep the dependence on the random bound explicit for the subsequent expectations.

Write \(m_N=N^{1/3},\qquad x_k=\frac{k}{{m_N}},\qquad \mathcal L_{N}^{\mathrm{sp}}=\frac{N+1}{{m_N}}.\) Set
\(b_{0,N}^{G,\sharp}=1,\qquad b_{N,N}^{G,\sharp}=b_{N,N}^G=0,\) and continue the exact integrated potential for the
Gaussian matrix \(\widetilde Y_{N,k}^G\) from \eqref{W:eq:Gaussian-exact-form-potential} through \(k=N\). Put
\begin{equation}
 \overline Y_{N,k}^G:=\mathbb E\widetilde Y_{N,k}^G,
 \qquad
 \widetilde{\mathcal M}_{N,k}^G
 :=\widetilde Y_{N,k}^G-\overline Y_{N,k}^G.
 \label{W:eq:gram-Gaussian-centered-exact-potential}
\end{equation}
For a sequence of sampled values $z_k$, its step interpolation has value $z_k$ on $[x_k,x_{k+1})$, $0\le k\le N$. We give it the constant values $z_0$ to the left of zero and $z_N$
to the right of $\mathcal L_N^{\mathrm{sp}}$. This also defines it on every interval $[j,j+2]$ occurring below, including those beyond the last sampling point.

\Needspace{6\baselineskip}

\begin{lemma}
\label{W:lem:gram-Gaussian-local-noise}
There are constants \(c,C>0\) and random variables \(\mathfrak Z_N\ge0\) such that
\begin{equation}
 \sup_N\mathbb E\exp(c\mathfrak Z_N^2)<\infty
 \label{W:eq:gram-Gaussian-noise-exp-moment}
\end{equation}
and, simultaneously for all integers \(0\le j\le\lceil \mathcal L_{N}^{\mathrm{sp}}\rceil\),
\begin{equation}
 \sup_{x,y\in[j,j+2]}
 \left|
 \widetilde{\mathcal M}_{N}^{G,\mathrm{st}}(x)
 -\widetilde{\mathcal M}_{N}^{G,\mathrm{st}}(y)
 \right|
 \le
 \mathfrak Z_N\sqrt{\log(2+j)}.
 \label{W:eq:gram-Gaussian-local-noise}
\end{equation}
\end{lemma}

\begin{proof}
The telescoping identity \eqref{W:eq:Gaussian-potential-boundary-identity}, with $b_{N,N}^G=0$, remains valid through $k=N$. Centering it gives
\begin{equation}
 \widetilde{\mathcal M}_{N,k}^G
 =\mathcal M_{N,k}^G+m_N(b_{k,N}^G-\mathbb Eb_{k,N}^G),
 \label{W:eq:gram-Gaussian-exact-noise-decomposition}
\end{equation}
where
\begin{equation}
 \mathcal M_{N,k}^G
 =-m_N\sum_{r=1}^ka_{r,N}^G
   -2m_N\sum_{r=1}^k(b_{r,N}^G-\mathbb Eb_{r,N}^G).
 \label{W:eq:gram-Gaussian-independent-martingale}
\end{equation}
At $k=0,N$ the centered off-diagonal coefficient is zero. Extend $\mathcal M_{N,k}^G$ constantly for integer indices outside $[0,N]$, and put the centered off-diagonal coefficient
equal to zero there. These conventions give the constant extensions of the step process in the statement.

Let $\mathcal F_{N,k}=\sigma(a_{r,N}^G,b_{r,N}^G:1\le r\le k)$, with deterministic boundary point coefficients included. The coefficient law in
Theorem~\ref{W:thm:Gaussian-Lanczos-law} makes successive pairs independent. For $r<N$ their chi parameter is $\beta(N-r)\ge1$. The Gaussian moment generating function and
\eqref{W:eq:chi-mgf} give, for every real $t$,
\begin{equation}
 \mathbb E\exp\!\left(t[-m_Na_{r,N}^G
                 -2m_N(b_{r,N}^G-\mathbb Eb_{r,N}^G)]\right)
 \le\exp\!\left(\frac{3t^2}{\beta m_N}\right).
 \label{W:eq:gram-Gaussian-increment-mgf}
\end{equation}
At $r=N$ the second summand is zero, so the same bound holds. Separately,
$\mathbb E\exp(t m_N(b_{r,N}^G-\mathbb Eb_{r,N}^G))\le\exp(t^2/(2\beta m_N))$.

Fix an integer $j\ge0$. Put $k_j=\min\{N,\lfloor m_Nj\rfloor\}$ and $Q_j=\min\{N-k_j,\lceil2m_N\rceil+2\}\le5m_N$. All sampled values of the step process on $[j,j+2]$ have indices
between $k_j$ and $k_j+Q_j$. Set $S_q=\mathcal M_{N,k_j+q}^G-\mathcal M_{N,k_j}^G$, $0\le q\le Q_j$. For either sign and $t>0$, conditional expectation in
\eqref{W:eq:gram-Gaussian-increment-mgf} shows that $V_q=\exp(\pm tS_q-3t^2q/(\beta m_N))$ is a nonnegative supermartingale. Stop it at the first $q\le Q_j$ with $\pm S_q>u$, or at
$Q_j$ if there is no such $q$. For the shifted filtration $\mathcal F_{N,k_j+q}$, the event $\{\tau\ge q\}$ is measurable at time $q-1$. The bounded stopping time consequently
satisfies
\[
 \mathbb E V_\tau
 =1+\sum_{q=1}^{Q_j}\mathbb E\bigl[\mathbf1_{\{\tau\ge q\}}(V_q-V_{q-1})\bigr]\le1.
\]
On the crossing event, $V_\tau\ge\exp(tu-15t^2/\beta)$. Summing the two signs and minimizing in $t$ yields
\begin{equation}
 \mathbb P\!\left(\max_{0\le q\le Q_j}|S_q|>u\right)
 \le2\inf_{t>0}e^{-tu+15t^2/\beta}
 =2e^{-\beta u^2/60}.
 \label{W:eq:gram-Gaussian-block-maximal}
\end{equation}
If $Q_j=0$, the maximum is zero and the bound holds directly. For the remaining coefficients, a union bound and \eqref{W:eq:chi-concentration} give, for $u\ge1$,
\begin{equation}
 \mathbb P\!\left(\max_{k_j\le r\le k_j+Q_j}
        m_N|b_{r,N}^G-\mathbb Eb_{r,N}^G|>u\right)
 \le C m_N e^{-\beta m_Nu^2/2}\le C'e^{-\beta u^2/4}.
 \label{W:eq:gaussian-offdiagonal-fluctuation-maximum}
\end{equation}
For the last inequality, factor out $e^{-\beta m_Nu^2/4}$ and use $\sup_{m\ge1}m e^{-\beta m/4}<\infty$. The oscillation of the process in
\eqref{W:eq:gram-Gaussian-exact-noise-decomposition} is bounded by twice each of these two maxima. Rescaling $u$ gives
\[
 \mathbb P\!\left(\sup_{x,y\in[j,j+2]}
 |\widetilde{\mathcal M}_{N}^{G,\mathrm{st}}(x)
  -\widetilde{\mathcal M}_{N}^{G,\mathrm{st}}(y)|>u\right)
 \le C e^{-cu^2},\qquad u\ge u_0,
\]
with constants independent of $N,j$.

Define $\mathfrak Z_N$ to be the maximum, over the finitely many integers in the statement, of this oscillation divided by $\sqrt{\log(2+j)}$. Each oscillation is a maximum of
finitely many differences of sampled values, so $\mathfrak Z_N$ is measurable. The union bound over the intervals gives
\[
 \mathbb P(\mathfrak Z_N>u)
 \le C\sum_{j\ge0}(j+2)^{-cu^2}
 \le C e^{-c'u^2},\qquad u\ge u_1.
\]
Indeed, for $r\ge2$,
$\sum_{j\ge0}(j+2)^{-r}\le2^{-r}+2^{1-r}/(r-1)$.
Choose $0<c_1<c'$. Tonelli's theorem gives
\[
 \mathbb E e^{c_1\mathfrak Z_N^2}
 =1+\int_0^\infty2c_1u e^{c_1u^2}\mathbb P(\mathfrak Z_N>u)\,\dd u
 \le C_1,
\]
where the interval $[0,u_1]$ is bounded directly. This proves the exponential moment uniformly in $N$; the defining maximum gives all the oscillation bounds on the same event. The
argument makes no assumption on the joint law across different matrix sizes.
\end{proof}

\Needspace{6\baselineskip}

On $H_0^1(0,\mathcal L_N^{\mathrm{sp}};\mathbb C)$ the centered potential defines the Hermitian sesquilinear form
\[
 (f,g)\longmapsto-\int_0^{\mathcal L_N^{\mathrm{sp}}}
       \widetilde{\mathcal M}_{N}^{G,\mathrm{st}}(x)(f\overline g)'(x)\,\dd x,
\]
with the first variable linear. The next estimate concerns its quadratic form, obtained by taking $g=f$. For the real ensemble the same formula is restricted to real functions. For
each finite $N$ the coefficient is a bounded step function, so the integral is defined on the stated domain.

\begin{lemma}
\label{W:lem:gram-Gaussian-form-noise}
For every \(\varepsilon>0\), there is \(C_\varepsilon<\infty\) such that, for every \(f\in H_0^1(0,\mathcal L_{N}^{\mathrm{sp}})\),
\begin{equation}
 \left| -\int_0^{\mathcal L_{N}^{\mathrm{sp}}} \widetilde{\mathcal M}_{N}^{G,\mathrm{st}}(x) (|f(x)|^2)'\,\dd x \right| \le{} \varepsilon\int_0^{\mathcal L_{N}^{\mathrm{sp}}}|f'|^2 +\varepsilon\int_0^{\mathcal L_{N}^{\mathrm{sp}}}x|f|^2 +C_\varepsilon(1+\mathfrak Z_N^4)\|f\|_2^2.
 \label{W:eq:gram-Gaussian-form-noise}
\end{equation}
\end{lemma}

\begin{proof}
Fix a realization for which the oscillation bounds of Lemma~\ref{W:lem:gram-Gaussian-local-noise} hold, and write $M=\widetilde{\mathcal M}_{N}^{G,\mathrm{st}}$ and
$\ell=\mathcal L_N^{\mathrm{sp}}$ within this proof. The step function $M$ is bounded for this finite matrix; its constant extension is locally integrable. Set
$\overline M(x)=\int_0^1M(x+t)\,\dd t$. Writing this as the difference of two primitives of $M$ proves absolute continuity and $\overline M'(x)=M(x+1)-M(x)$ almost everywhere. For
$j=\lfloor x\rfloor$, the interval $[x,x+1]$ lies in $[j,j+2]$. Averaging the oscillation bound over $t$ therefore gives
\begin{equation}
 |M(x)-\overline M(x)|+|\overline M'(x)|
 \le C\mathfrak Z_N\sqrt{\log(2+x)},\qquad 0\le x\le\ell,
 \label{W:eq:gram-Gaussian-average-bounds}
\end{equation}
with the derivative statement understood almost everywhere.

For complex as well as real $f\in H_0^1(0,\ell)$, the function $|f|^2$ belongs to $W^{1,1}(0,\ell)$, has zero boundary values and satisfies
$(|f|^2)'=2\operatorname{Re}(f'\overline f)$ almost everywhere. The integration by parts formula for absolutely continuous functions thus gives the exact identity
\[
 -\int_0^\ell M(|f|^2)'\,\dd x
 =-\int_0^\ell(M-\overline M)(|f|^2)'\,\dd x
   +\int_0^\ell\overline M'|f|^2\,\dd x.
\]
All integrals are finite by Cauchy--Schwarz and boundedness on $[0,\ell]$. There is no boundary term because $|f(0)|^2=|f(\ell)|^2=0$. The absolute value of the right side is at
most
\begin{align*}
 &2C\int_0^\ell\mathfrak Z_N\sqrt{\log(2+x)}\,|f'||f|\,\dd x
   +C\int_0^\ell\mathfrak Z_N\sqrt{\log(2+x)}\,|f|^2\,\dd x\\
 &\qquad\le\frac{\varepsilon}{2}\|f'\|_2^2
   +C_\varepsilon\int_0^\ell
       \bigl[\mathfrak Z_N^2\log(2+x)+\mathfrak Z_N\sqrt{\log(2+x)}\bigr]|f|^2\,\dd x.
\end{align*}
For $A\ge0$ and $\delta>0$, maximizing $A\log(2+x)-\delta x$ over $x\ge0$ places the maximum at zero or at $x=A/\delta-2$. Hence
\begin{equation}
 A\log(2+x)\le\delta x+C_\delta(1+A^2),\qquad A,x\ge0.
 \label{W:eq:gram-log-absorption}
\end{equation}
More explicitly the maximum is at most $A\log2+A\log_+(A/(2\delta))$, which is bounded by $C_\delta(1+A^2)$. Since $\mathfrak Z_N\sqrt{\log(2+x)}\le1+\mathfrak Z_N^2\log(2+x)$,
choose $A=\mathfrak Z_N^2$ and then $\delta=\varepsilon/(2C_\varepsilon)$ in \eqref{W:eq:gram-log-absorption}. This bounds the last integral by
\[
 \varepsilon\int_0^\ell x|f|^2\,\dd x
   +C'_\varepsilon(1+\mathfrak Z_N^4)\|f\|_2^2.
\]
The constants are independent of $N$, the realization and $f$. Since the calculation applies deterministically to every $f$ on the one event specifying $\mathfrak Z_N$, the
conclusion holds simultaneously on that event for all functions in the stated form domain.
\end{proof}

\Needspace{6\baselineskip}

Put \(K_N=\lfloor N/2\rfloor\) and define
\begin{equation}
 \mathfrak h_N[v]
 :={m_N}^2\sum_{k=0}^{K_N}|v_{k+1}-v_k|^2
 +\sum_{k=1}^{N}x_k|v_k|^2,
 \qquad v_0=v_{N+1}=0.
 \label{W:eq:gram-discrete-Airy-comparison-form}
\end{equation}

The form $\mathfrak h_N$ is defined on all of $\mathbb F_\beta^N$ with counting-measure norm $\|v\|_2^2=\sum_{k=1}^N|v_k|^2$. Its associated Hermitian matrix is characterized by
polarization of \eqref{W:eq:gram-discrete-Airy-comparison-form}; no additional operator domain is involved.

\begin{proposition}
\label{W:prop:gram-Gaussian-global-coercivity}
There are constants \(c_0,C_0,c>0\) and events \(\Omega_N^{\mathrm{coer}}\) such that
\begin{equation}
 \mathbb P((\Omega_N^{\mathrm{coer}})^c)\le Ce^{-cN}
 \label{W:eq:gram-Gaussian-coercive-event}
\end{equation}
and, on \(\Omega_N^{\mathrm{coer}}\),
\begin{equation}
 \langle v,\mathcal A_N^Gv\rangle
 \ge c_0\mathfrak h_N[v]
 -C_0(1+\mathfrak Z_N^4)\|v\|_2^2
 \label{W:eq:gram-Gaussian-global-coercivity}
\end{equation}
for every \(v\in\mathbb F_\beta^N\).
\end{proposition}

\begin{proof}
Write $\bar b_0=1$, $\bar b_k=b_{k,N}^G$ for $1\le k<N$, and $\bar b_N=0$. These auxiliary coefficients enter only the form identity; the actual Lanczos matrix has no off-diagonal
entry to the left of its first coordinate. Give $v$ the exterior values $v_0=v_{N+1}=0$. Its linear interpolant $f=\mathcal I_Nv$ has values $f(x_k)=\sqrt{m_N}v_k$ and belongs to
$H_0^1(0,\mathcal L_N^{\mathrm{sp}})$. Summing over the cells of the step potential gives
\begin{align*}
 -\int_0^{\mathcal L_N^{\mathrm{sp}}}\widetilde Y_N^{G,\mathrm{st}}(|f|^2)'\,\dd x
 &=\sum_{k=0}^N\widetilde Y_{N,k}^G(|f(x_k)|^2-|f(x_{k+1})|^2)\\
 &=\sum_{k=1}^N(\widetilde Y_{N,k}^G-\widetilde Y_{N,k-1}^G)|f(x_k)|^2\\
 &=m_N^2\sum_{k=1}^N(2-a_{k,N}^G-\bar b_{k-1}-\bar b_k)|v_k|^2.
\end{align*}
Expanding the kinetic sum adds $m_N^2(\bar b_{k-1}+\bar b_k)|v_k|^2$ at coordinate $k$ and the off-diagonal term $-2m_N^2\bar b_k\operatorname{Re}(v_{k+1}\overline v_k)$. The
diagonal coefficients $\bar b_{k-1}+\bar b_k$ cancel, including the two boundary coordinates. Thus
\begin{align}
 \langle v,\mathcal A_N^Gv\rangle
 &=m_N^2\sum_{k=0}^N\bar b_k|v_{k+1}-v_k|^2
   -\int_0^{\mathcal L_N^{\mathrm{sp}}}\widetilde Y_N^{G,\mathrm{st}}(|f|^2)'\,\dd x.
 \label{W:eq:gram-Gaussian-full-form}
\end{align}

The mean of the diagonal coefficient is zero and $\mathbb E(b_{k,N}^G)^2=1-k/N$ for $1\le k<N$. Jensen's inequality, with the deterministic values at $k=0,N$, gives
$\mathbb E\bar b_k\le\sqrt{1-k/N}$ throughout $0\le k\le N$. Using $1-\sqrt{1-y}\ge y/2$ on $[0,1]$, we obtain
\begin{equation}
 2-\mathbb E\bar b_{k-1}-\mathbb E\bar b_k\ge\frac{2k-1}{2N},\qquad 1\le k\le N.
 \label{W:eq:gram-Gaussian-mean-increment}
\end{equation}
Since $m_N^3=N$, the mean potential is at least $\sum_{k=1}^Nx_k|v_k|^2-(2m_N)^{-1}\|v\|_2^2$. In particular,
\begin{equation}
 m_N^2\sum_{k=1}^N(2-\mathbb E\bar b_{k-1}-\mathbb E\bar b_k)|v_k|^2
 \ge\frac12\sum_{k=1}^Nx_k|v_k|^2-\frac1{2m_N}\|v\|_2^2.
 \label{W:eq:gram-Gaussian-mean-potential-lower}
\end{equation}

Define $\Omega_N^{\mathrm{coer}}=\{\min_{0\le k\le K_N}\bar b_k\ge1/4\}$. The variance bound proved in Lemma~\ref{W:lem:chi-estimates} implies $\mu_\nu\ge\sqrt{\nu-1}$ for
$\nu\ge1$. For $N\ge4$ and $1\le k\le K_N$,
\[
 \mathbb Eb_{k,N}^G\ge\sqrt{\frac{\beta(N-k)-1}{\beta N}}
 \ge\sqrt{\frac12-\frac1{\beta N}}\ge\frac12.
\]
The one-sided exponential bound from the same lemma gives
\[
 \mathbb P(b_{k,N}^G<1/4)
 \le\mathbb P\!\left(\chi_{\beta(N-k)}-\mu_{\beta(N-k)}
                           <-\tfrac14\sqrt{\beta N}\right)
 \le e^{-\beta N/32}.
\]
There are at most $N$ such indices and $\bar b_0=1$. Thus $\mathbb P((\Omega_N^{\mathrm{coer}})^c)\le N e^{-\beta N/32}\le C e^{-\beta N/64}$ for $N\ge4$; enlarging $C$ covers
$N<4$. All other $\bar b_k$ are nonnegative, so on this event
\begin{equation}
 m_N^2\sum_{k=0}^N\bar b_k|v_{k+1}-v_k|^2
 \ge\frac14m_N^2\sum_{k=0}^{K_N}|v_{k+1}-v_k|^2.
 \label{W:eq:gram-Gaussian-kinetic-lower}
\end{equation}

For $h=m_N^{-1}$ and $f_k=h^{-1/2}v_k$, convexity along each linear interpolant gives $|f(x_k+th)|^2\le(1-t)|f_k|^2+t|f_{k+1}|^2$, $0\le t\le1$. The total coefficient of $|f_k|^2$
in the corresponding upper bound for $\int|f|^2$ is $h$. In the bound for $\int x|f|^2$, its two adjacent cells contribute
\[
 h\left(x_{k-1}/2+h/3+x_k/2+h/6\right)=h x_k.
\]
Both exterior nodal values are zero. Consequently
\begin{equation}
 \|f\|_2^2\le\|v\|_2^2,
 \qquad \int_0^{\mathcal L_N^{\mathrm{sp}}}x|f|^2\,\dd x
 \le\sum_{k=1}^Nx_k|v_k|^2.
 \label{W:eq:gram-Gaussian-weighted-mass-comparison}
\end{equation}
In the rest of the proof write $G_1=m_N^2\sum_{k=0}^{K_N}|v_{k+1}-v_k|^2$ and $X=\sum_{k=1}^Nx_k|v_k|^2$. The exact cell derivative formula gives
$\|f'\|_2^2=m_N^2\sum_{k=0}^N|v_{k+1}-v_k|^2$. For $k\ge K_N+1$ one has $x_k>m_N^2/2$, and each $|v_k|^2$ occurs at most twice after applying $|a-b|^2\le2|a|^2+2|b|^2$. Therefore
\begin{equation}
 m_N^2\sum_{k=K_N+1}^N|v_{k+1}-v_k|^2
 \le4m_N^2\sum_{k=K_N+1}^N|v_k|^2\le8X.
 \label{W:eq:gram-Gaussian-tail-gradient-absorption}
\end{equation}
Apply Lemma~\ref{W:lem:gram-Gaussian-form-noise} with the deterministic choice $\varepsilon=1/72$.
Equations~\eqref{W:eq:gram-Gaussian-weighted-mass-comparison}--\eqref{W:eq:gram-Gaussian-tail-gradient-absorption} bound the absolute value of the centered potential term by
\[
 \tfrac1{72}(G_1+8X)+\tfrac1{72}X+C(1+\mathfrak Z_N^4)\|v\|_2^2
 =\tfrac1{72}G_1+\tfrac18X+C(1+\mathfrak Z_N^4)\|v\|_2^2.
\]
Subtracting this from the mean and kinetic bounds in \eqref{W:eq:gram-Gaussian-full-form}, and using $m_N\ge1$, gives
\begin{align*}
 \langle v,\mathcal A_N^Gv\rangle
 &\ge(1/4-1/72)G_1+(1/2-1/8)X
          -C_0(1+\mathfrak Z_N^4)\|v\|_2^2\\
 &\ge\tfrac18(G_1+X)-C_0(1+\mathfrak Z_N^4)\|v\|_2^2.
\end{align*}
The choice $c_0=1/8$ is deterministic, and $C_0$ is independent of $N$ and $v$. The form estimate holds for all $f$ on its fixed event, so this final inequality holds
simultaneously for all $v$ on $\Omega_N^{\mathrm{coer}}$.
\end{proof}

\Needspace{6\baselineskip}

\begin{lemma}
\label{W:lem:gram-discrete-Airy-Weyl}
Let \(H_N^{\mathrm{det}}\) be the operator associated with \(\mathfrak h_N\). There is a numerical \(C\) such that
\begin{equation}
 \operatorname{Tr}\mathbf1_{(-\infty,E]}(H_N^{\mathrm{det}})
 \le C(1+E_+)^{3/2}
 \label{W:eq:gram-discrete-Airy-Weyl}
\end{equation}
uniformly in \(N\) and \(E\in\mathbb R\).
\end{lemma}

\begin{proof}
The form $\mathfrak h_N$ is nonnegative, so the count is zero for $E<0$. If $E\ge0$ and $N<64$, it is at most $N\le64(1+E)^{3/2}$. Assume henceforth $N\ge64$. For $E\ge m_N^2/8$,
the dimension bound gives
$\operatorname{Tr}\mathbf1_{(-\infty,E]}(H_N^{\mathrm{det}})\le N=m_N^3\le8^{3/2}E^{3/2}$.
For $0\le E<m_N^2/8$, set $J=\lceil2m_N(E+1)\rceil$. Then
\[
 J\le N/4+2N^{1/3}+1<\lfloor N/2\rfloor=K_N.
\]
In particular, all edges between sites $1,\ldots,J$ occur in $\mathfrak h_N$. Drop the left boundary edge, the edge from $J$ to $J+1$, the potential on the first $J$ sites, and all
kinetic terms on the remaining sites. The resulting form is
\[
 m_N^2\sum_{k=1}^{J-1}|v_{k+1}-v_k|^2+\sum_{k=J+1}^Nx_k|v_k|^2.
\]
It defines a matrix $D_N\le H_N^{\mathrm{det}}$, equal to the direct sum of $m_N^2\Delta_J^{\mathrm N}$ and the diagonal matrix with entries $x_k$, $k>J$. Since
$x_k>J/m_N\ge2(E+1)>E$ on the latter block, only the Neumann block can contribute to the count at $E$. The min--max principle gives $\lambda_j(D_N)\le\lambda_j(H_N^{\mathrm{det}})$
in increasing order, and therefore the count for $H_N^{\mathrm{det}}$ is at most that for this Neumann block.

For $0\le l<J$, the vector $u_k=\cos(\pi l(k-1/2)/J)$, $1\le k\le J$, has the extensions $u_0=u_1$ and $u_{J+1}=u_J$. The identity
\[
 2u_k-u_{k-1}-u_{k+1}=4\sin^2\!\left(\frac{\pi l}{2J}\right)u_k
\]
proves that its Neumann eigenvalue is $4m_N^2\sin^2(\pi l/(2J))$. These $J$ eigenvalues are distinct and their nonzero eigenvectors give the full spectrum. The inequality
$\sin t\ge2t/\pi$ on $[0,\pi/2]$ shows that an eigenvalue at most $E$ must have $l\le J\sqrt E/(2m_N)$. The number of eligible integers, including $l=0$, is at most
\[
 1+\frac{J\sqrt E}{2m_N}
 \le1+\bigl(E+1+(2m_N)^{-1}\bigr)\sqrt E
 \le C(1+E)^{3/2}.
\]
Combining the three cases proves the bound uniformly in $N,E$.
\end{proof}

\Needspace{6\baselineskip}

Write $\mathcal N_N^G(E)=\operatorname{Tr} \mathbf1_{(-\infty,E]}(\mathcal A_N^G)$ for the eigenvalue count, including multiplicities.

\begin{theorem}
\label{W:thm:gram-Gaussian-global-count}
For every $p>0$,
\begin{equation}
 \sup_N\|\mathcal N_N^G(E)\|_{L^p}
       \le C_{p,\beta}(1+E_+)^{3/2},\qquad E\in\mathbb R.
 \label{W:eq:gram-Gaussian-global-count}
\end{equation}
\end{theorem}

\begin{proof}
Put $C_N=C_0(1+\mathfrak Z_N^4)$, with $C_0$ from Proposition~\ref{W:prop:gram-Gaussian-global-coercivity}. For each outcome in $\Omega_N^{\mathrm{coer}}$, its form inequality and
the finite-dimensional min--max principle give
\[
 \lambda_j(\mathcal A_N^G)\ge c_0\lambda_j(H_N^{\mathrm{det}})-C_N,\qquad 1\le j\le N.
\]
Consequently, Lemma~\ref{W:lem:gram-discrete-Airy-Weyl} yields the pointwise bound
\[
 \mathbf1_{\Omega_N^{\mathrm{coer}}}\mathcal N_N^G(E)
 \le C\mathbf1_{\Omega_N^{\mathrm{coer}}}
          \bigl(1+((E+C_N)/c_0)_+\bigr)^{3/2}
 \le C\mathbf1_{\Omega_N^{\mathrm{coer}}}(1+E_++\mathfrak Z_N^4)^{3/2}.
\]
On the complementary event use the deterministic count $\mathcal N_N^G(E)\le N$. For every $p>0$, these two disjoint events give
\begin{align*}
 \mathbb E[\mathcal N_N^G(E)^p]
 &\le C_p\mathbb E\!\left[
       \mathbf1_{\Omega_N^{\mathrm{coer}}}(1+E_++\mathfrak Z_N^4)^{3p/2}\right]
       +N^p\mathbb P((\Omega_N^{\mathrm{coer}})^c)\\
 &\le C_p(1+E_+)^{3p/2}\mathbb E(1+\mathfrak Z_N^4)^{3p/2}
       +C N^p e^{-cN}.
\end{align*}
For fixed $p$, $(1+z^4)^{3p/2}\le C_{p,c_1}e^{c_1z^2}$ on $z\ge0$, where $c_1$ is any exponent allowed in \eqref{W:eq:gram-Gaussian-noise-exp-moment}. Also
$\sup_{N\ge1}N^p e^{-cN}<\infty$. The last display is therefore at most $C_{p,\beta}(1+E_+)^{3p/2}$; taking its $p$th root proves the moment estimate. The pointwise bound and its
probability integral apply for every $p>0$, with $\mathfrak Z_N$ and $\Omega_N^{\mathrm{coer}}$ defined on the same probability space.
\end{proof}

\begin{corollary}
\label{W:cor:gaussian-spectral-count-tails}
For $p,c>0$ and $b\ge0$,
\begin{align}
 \sup_N\|\operatorname{Tr}\mathbf1_{[b,b+1)}(\mathcal A_N^G)\|_{L^p}
    &\le C_{p,\beta}(1+b)^{3/2},
       \label{W:eq:gram-Gaussian-unit-count}\\
 \sup_N\left\|\sum_{\lambda\in\sigma(\mathcal A_N^G),\,\lambda>b}
                  e^{-c(\lambda-b)}\right\|_{L^p}
    &\le C_{p,c,\beta}(1+b)^{3/2}.
       \label{W:eq:gram-Gaussian-exponential-tail}
\end{align}
Moreover,
\begin{equation}
 \lim_{M\to\infty}\limsup_{N\to\infty}
       \mathbb P(\inf\sigma(\mathcal A_N^G)<-M)=0.
 \label{W:eq:gram-Gaussian-bottom-tightness}
\end{equation}
\end{corollary}
\begin{proof}
The interval count satisfies
$\operatorname{Tr}\mathbf1_{[b,b+1)}(\mathcal A_N^G)\le\mathcal N_N^G(b+1)$.
Theorem~\ref{W:thm:gram-Gaussian-global-count} and $2+b\le2(1+b)$ give \eqref{W:eq:gram-Gaussian-unit-count}. For $p\ge1$, first restrict the weighted sum to the finitely many
intervals $[b+r,b+r+1)$, $0\le r\le R$. For every eigenvalue in that interval, $e^{-c(\lambda-b)}\le e^{-cr}$, and all multiplicities are included in the interval count.
Minkowski's inequality for this finite sum therefore gives
\begin{align*}
 \left\|\sum_{r=0}^R\ \sum_{\substack{\lambda>b\\b+r\le\lambda<b+r+1}}
                         e^{-c(\lambda-b)}\right\|_{L^p}
 &\le\sum_{r=0}^R e^{-cr}
       \|\operatorname{Tr}\mathbf1_{[b+r,b+r+1)}(\mathcal A_N^G)\|_{L^p}\\
 &\le C_{p,\beta}(1+b)^{3/2}\sum_{r=0}^R e^{-cr}(1+r)^{3/2}.
\end{align*}
Here $1+b+r\le(1+b)(1+r)$. The scalar series converges by the ratio test. The partial sums on the left are nonnegative and increase to the full weighted sum, so monotone
convergence applied to their $p$th powers proves \eqref{W:eq:gram-Gaussian-exponential-tail}. For $0<p<1$, Jensen's inequality applied to $x^p$ gives $\|X\|_{L^p}\le\|X\|_{L^1}$
for nonnegative $X$, so the established $p=1$ bound suffices.

On $\Omega_N^{\mathrm{coer}}$, coercivity implies $\lambda_{\min}(\mathcal A_N^G)\ge-C_0(1+\mathfrak Z_N^4)$. For $M\ge2C_0$, it follows that
\begin{align*}
 \mathbb P(\lambda_{\min}(\mathcal A_N^G)<-M)
 &\le\mathbb P((\Omega_N^{\mathrm{coer}})^c)
       +\mathbb P\!\left(\mathfrak Z_N^2>\sqrt{M/C_0-1}\right)\\
 &\le C e^{-cN}+\mathbb E e^{c_1\mathfrak Z_N^2}
                         e^{-c_1\sqrt{M/C_0-1}}\\
 &\le C e^{-cN}+C e^{-c'\sqrt M}.
\end{align*}
The last two constants are independent of $N,M$. Taking first $\limsup_N$ and then $M\to\infty$ proves \eqref{W:eq:gram-Gaussian-bottom-tightness}.
\end{proof}

The same estimate controls the two extreme eigenvalues of the original Gaussian matrix. Since $\mathcal A_N^G=m_N^2(2I-T_N^G)$ and $T_N^G$ is unitarily equivalent to
$H_N^{\mathrm G}$,
\[
 \mathbb P\bigl(\lambda_{\max}(H_N^{\mathrm G})>2+d/2\bigr)
 =\mathbb P\bigl(\lambda_{\min}(\mathcal A_N^G)<-m_N^2d/2\bigr)
 \le C e^{-cN}+C e^{-c'\sqrt{m_N^2d}},\qquad d>0.
\]
For $m_N^2d/2<2C_0$, enlarge $C$ so that this bound follows from the bound by one. The Gaussian distributions of $H_N^{\mathrm G}$ and $-H_N^{\mathrm G}$ agree. Applying the same
inequality to $-H_N^{\mathrm G}$ and adding the two probabilities gives
\[
 \mathbb P\bigl(\|H_N^{\mathrm G}\|>2+d/2\bigr)
 \le C e^{-cN}+C e^{-c'\sqrt{m_N^2d}}.
\]
This uses equality of the two marginal laws, without assuming independence. In particular, whenever $m_N^2d_N\to\infty$, this probability tends to zero. The estimate is independent
of the deterministic starting vector because both extreme eigenvalues belong to the original matrix. It supplies the Gaussian spectral confinement used below directly from the
coercive form estimate.

\subsubsection{Expected traces and spectral localization for the GOE and GUE}
\label{W:sec:smoothed-traces-and-gaussian-localization}

Regularizing a polynomial Gram matrix by adding a positive scalar multiple of the identity defines a contraction whose range lies in the first Lanczos coordinates, including at a
breakdown. We use this contraction to control the mass in the tridiagonal coordinates of the whole spectral interval $(-\infty,E]$; the counting moments just proved bound the high
energy remainder in expectation. We then use a smooth cutoff in the eigenvalues of the Gram matrix for entry comparison.

Throughout the localization argument below, $w\in\mathbb F_\beta^N$ is deterministic with $\|w\|=1$, ${m_N}=N^{1/3}$, ${h_N}={m_N}^{-1}$, and
${\rho}=\min\{N,\lfloor {m_N}L\rfloor\}$, with $L\ge1$. Let $U=U^{(w)}(H)$ be the orthogonal or unitary matrix whose columns are the Lanczos vectors of
Section~\ref{W:sec:lanczos-tridiagonalization}, and put
\[
 B={m_N}^2(2I-H),\qquad A=U^*BU,\qquad
 P_{\rho}^{\rm c}=U\Bigl(\sum_{k\le {\rho}}e_ke_k^*\Bigr)U^*,
 \qquad Q_{\rho}^{\rm c}=I-P_{\rho}^{\rm c}.
\]
Thus $P_\rho^{\rm c}$ is the projection onto the first $\rho$ Lanczos vectors. We define the smooth statistic directly in the original matrix coordinates.

\Needspace{6\baselineskip}

Set $X=[\psi_0(H)w,\ldots,\psi_{\rho-1}(H)w]$, $G=X^*X$, and, for $\alpha>0$, define
\begin{equation}
 K_{\rho,\alpha}=X(G+\alpha I_\rho)^{-1}X^*,\qquad
 Q_{\rho,\alpha}=I-K_{\rho,\alpha}.
 \label{W:eq:regularized-polynomial-projection}
\end{equation}

\begin{lemma}
\label{W:lem:regularized-projection-order}
The matrices $K_{\rho,\alpha}$ and $Q_{\rho,\alpha}$ are smooth functions of every Hermitian $H$ and satisfy
\begin{equation}
 0\le K_{\rho,\alpha}\le P_\rho^{\rm c}\le I,
 \qquad Q_\rho^{\rm c}\le Q_{\rho,\alpha}\le I.
 \label{W:eq:resolvent-order}
\end{equation}
If $G>0$ and $g=\lambda_{\min}(G)$, then
\begin{equation}
 0\le P_\rho^{\rm c}-K_{\rho,\alpha}
       \le\frac{\alpha}{g+\alpha}P_\rho^{\rm c}.
 \label{W:eq:regularized-projection-error}
\end{equation}
\end{lemma}
\begin{proof}
The singular value decomposition $X=V\Sigma W^*$ shows that the nonzero eigenvalues of $K_{{\rho},\alpha}$ are $s_j^2/(s_j^2+\alpha)\in[0,1)$ and that its range is
$\operatorname{Ran}X$. The monic polynomials $\psi_j$ give $\operatorname{Ran}X=\operatorname{span}\{w,Hw,\ldots,H^{{\rho}-1}w\}$. If no breakdown precedes ${\rho}$, this is the
span of the first Lanczos vectors. If the first breakdown is at $r<{\rho}$, its $r$-dimensional Krylov subspace is invariant under $H$; all further powers of $H$ applied to $w$
stay in it. The span of the first Lanczos vectors includes this subspace and additional vectors constructed on the orthogonal complement. Thus
$\operatorname{Ran}X\subset\operatorname{Ran}P_{\rho}^{\rm c}$ in both cases. If $v=v_1+v_2$, with $v_1=P_\rho^{\rm c}v$ and $v_2=(I-P_\rho^{\rm c})v$, then $K_{\rho,\alpha}v_2=0$
and
\[
 0\le v^*K_{\rho,\alpha}v
       =v_1^*K_{\rho,\alpha}v_1\le\|v_1\|^2
       =v^*P_\rho^{\rm c}v.
\]
This proves \eqref{W:eq:resolvent-order}. When $G>0$, its rank is $\rho$ and $\operatorname{Ran}X=\operatorname{Ran}P_\rho^{\rm c}$. On the left singular vector corresponding to
$s_j>0$, the difference $P_\rho^{\rm c}-K_{\rho,\alpha}$ has eigenvalue $\alpha/(s_j^2+\alpha)\le\alpha/(g+\alpha)$; it vanishes on the orthogonal complement. This proves
\eqref{W:eq:regularized-projection-error}. The inverse $G+\alpha I_{\rho}$ is positive definite at every $H$; matrix inversion and polynomial functional calculus prove smoothness.
\end{proof}

\Needspace{6\baselineskip}

Fix $E\in\mathbb R$ and define
\begin{equation}
 \begin{gathered}
 f_E(x)=(1+e^{x-E})^{-1},\qquad F_E(H)=f_E(B),\\
 C_{N,L,E,\alpha}^{(w)}(H)
       =\operatorname{Tr}[F_E(H)Q_{\rho,\alpha}(H)F_E(H)],\qquad
 Z_{N,L,E,\alpha}^{(w)}=\frac{C_{N,L,E,\alpha}^{(w)}}{1+C_{N,L,E,\alpha}^{(w)}}.
 \end{gathered}
 \label{W:eq:regularized-localization-statistic}
\end{equation}
In the Lanczos coordinates put $\Pi_N(E)=\mathbf1_{(-\infty,E]}(A)$ and $Q_{N,L}=I-\sum_{k\le\rho}e_ke_k^*$.

\begin{proposition}
\label{W:prop:regularized-localization-bound}
The statistic $Z=Z_{N,L,E,\alpha}^{(w)}$ is smooth on the Hermitian matrices, takes values in $[0,1)$, and satisfies
\begin{equation}
 \|Q_{N,L}\Pi_N(E)\|^2\le4C_{N,L,E,\alpha}^{(w)}(H).
 \label{W:eq:regularized-spectral-bound}
\end{equation}
For any random Hermitian $H$ and $\eta>0$,
\begin{equation}
 \mathbb P(\|Q_{N,L}\Pi_N(E)\|>\eta)
 \le\frac{4+\eta^2}{\eta^2}\mathbb EZ_{N,L,E,\alpha}^{(w)}(H).
 \label{W:eq:regularized-localization-probability}
\end{equation}
\end{proposition}
\begin{proof}
The trace defining $C$ is nonnegative, and at most $N$, since both $F_E$ and $Q_{{\rho},\alpha}$ are positive contractions. This proves boundedness of $Z$ and excludes poles in its
definition. Smoothness follows from Lemma~\ref{W:lem:regularized-projection-order} and Lemma~\ref{W:lem:smooth-functional-calculus}, applied on each compact spectral interval to
the smooth function $f_E$. Let $\Pi=\mathbf1_{(-\infty,E]}(B)$. The scalar inequality $f_E(x)\ge1/2$ for every $x\le E$ gives $\|f_E(B)^{-1}\Pi\|\le2$. Hence
\[
 \|Q_{\rho}^{\rm c}\Pi\|
 \le2\|Q_{\rho}^{\rm c}F_E\|
 \le2\|Q_{\rho}^{\rm c}F_E\|_{\rm HS}
 \le2\{\operatorname{Tr}(F_EQ_{{\rho},\alpha}F_E)\}^{1/2}.
\]
For the last inequality, congruence by $F_E$ preserves the order: $F_EQ_\rho^{\rm c}F_E\le F_EQ_{\rho,\alpha}F_E$. Also $(Q_\rho^{\rm c})^2=Q_\rho^{\rm c}$, so the trace of the
left side is $\|Q_\rho^{\rm c}F_E\|_{\rm HS}^2$. Unitary conjugation proves \eqref{W:eq:regularized-spectral-bound}. Pointwise,
\[
 \mathbf1_{\{\|Q_{N,L}\Pi_N(E)\|>\eta\}}
       \le\frac{4+\eta^2}{\eta^2}Z.
\]
Indeed the indicated event forces $C>\eta^2/4$ and hence $Z>\eta^2/(4+\eta^2)$. Taking expectations proves the result.
\end{proof}

By \eqref{W:eq:resolvent-order}, the regularized complement dominates the projection onto the remaining Lanczos coordinates, including the subsequent invariant blocks. For the GOE
and GUE the off-diagonal Lanczos coefficients are positive almost surely, and \eqref{W:eq:regularized-projection-error} bounds the difference between the two projections.

\Needspace{6\baselineskip}

Suppose a finite self-adjoint matrix $A_N$ in the discrete coordinate space has, on an event $\Omega_N$, the form bound
\begin{equation}
 \langle v,A_Nv\rangle\ge c_0\mathfrak h_N[v]-C_N\|v\|^2,
 \quad C_N\ge0,\qquad
 \operatorname{Tr}\mathbf1_{(-\infty,u]}(H_N^{\rm det})
 \le C(1+u_+)^{3/2},
 \label{W:eq:resolvent-confinement-input}
\end{equation}
where $\mathfrak h_N[v]\ge\sum_k k{h_N}|v_k|^2$ and $H_N^{\rm det}$ is its associated operator.
\begin{lemma}
\label{W:lem:resolvent-weighted-trace}
For every fixed $E$ there is $C_E$ independent of $N,C_N,L$ such that, on $\Omega_N$,
\begin{align}
 \operatorname{Tr} f_E(A_N)^2&\le C_E(1+C_N)^{3/2},
 \label{W:eq:resolvent-trace-zero}\\
 \sum_{\lambda\in\sigma(A_N)}(\lambda+C_N)f_E(\lambda)^2
 &\le C_E(1+C_N)^{5/2},
 \label{W:eq:resolvent-trace-one}\\
 \operatorname{Tr}[f_E(A_N)Q_{N,L}f_E(A_N)]
 &\le\frac{C_E}{L}(1+C_N)^{5/2}.
 \label{W:eq:resolvent-coercive-localization}
\end{align}
Multiplicity is included in both sums.
\end{lemma}
\begin{proof}
Fix an outcome in $\Omega_N$. All the estimates in this proof are deterministic on that outcome. The form inequality and the finite-dimensional min--max principle give, in
increasing eigenvalue order,
\[
 \lambda_j(A_N)\ge c_0\lambda_j(H_N^{\mathrm{det}})-C_N\ge-C_N,
 \qquad
 \mathcal N_{A_N}(u)\le C\bigl[1+((u+C_N)/c_0)_+\bigr]^{3/2}.
\]
The constants below may depend on $c_0$, the numerical constant in the deterministic counting bound, and $E$, but not on $N$, the value of $C_N$, or $L$. Let $a=1+\max\{E,0\}$. For
$q=0,1$, with the factor of degree zero interpreted as one, the part of the spectral sum with $\lambda\le a$ is bounded by
\[
 \sum_{\lambda\le a}(\lambda+C_N)^q f_E(\lambda)^2
 \le(a+C_N)^q\mathcal N_{A_N}(a)
 \le C_E(1+C_N)^{q+3/2}.
\]
For the part in $(a+j,a+j+1]$, $j\ge0$, use $\lambda+C_N\ge0$, $f_E(\lambda)^2\le e^{-2(\lambda-E)}\le e^{-2(a+j-E)}$, and the cumulative count at $a+j+1$. Thus
\begin{align*}
 \sum_{a+j<\lambda\le a+j+1}(\lambda+C_N)^q f_E(\lambda)^2
 &\le e^{-2(a+j-E)}(a+j+1+C_N)^q\mathcal N_{A_N}(a+j+1)\\
 &\le C_E(1+C_N)^{q+3/2}e^{-2j}(1+j)^{q+3/2}.
\end{align*}
The eigenvalues are counted with multiplicity in both lines. Summing the nonnegative terms first over $j\le R$ and then letting $R\to\infty$ bounds the full sum by
\[
 C_E(1+C_N)^{q+3/2}
       \left(1+\sum_{j\ge0}e^{-2j}(1+j)^{q+3/2}\right).
\]
The scalar series is finite by the ratio test. This proves \eqref{W:eq:resolvent-trace-zero}--\eqref{W:eq:resolvent-trace-one}, including eigenvalues close to the lower bound
$-C_N$.

For $\rho=\min\{N,\lfloor m_NL\rfloor\}$, $Q_{N,L}$ is zero if $\rho=N$. Otherwise $k h_N>L$ for every $k>\rho$. For a unit eigenvector $u_j$ of $A_N$ with eigenvalue $\lambda_j$,
it follows that
\[
 L\|Q_{N,L}u_j\|_2^2
 \le\sum_{k>\rho} k h_N|(u_j)_k|^2
 \le\mathfrak h_N[u_j]
 \le c_0^{-1}(\lambda_j+C_N).
\]
Expanding the trace in an orthonormal eigenbasis gives the exact finite sum
\[
 \operatorname{Tr}[f_E(A_N)Q_{N,L}f_E(A_N)]
 =\sum_{j=1}^N f_E(\lambda_j)^2\|Q_{N,L}u_j\|_2^2
 \le\frac1{c_0L}\sum_{j=1}^N(\lambda_j+C_N)f_E(\lambda_j)^2.
\]
The already proved degree-one estimate yields \eqref{W:eq:resolvent-coercive-localization}. The calculation applies Cauchy--Schwarz to the factors in their displayed order and uses
$Q_{N,L}^2=Q_{N,L}=Q_{N,L}^*$.
\end{proof}
For the Gaussian matrix put $C_N=C_0(1+\mathfrak Z_N^4)$ and $T_N=\operatorname{Tr}f_E(\mathcal A_N^G)^2$. The preceding lemma applies on $\Omega_N^{\mathrm{coer}}$, while the
scalar bound $0\le f_E\le1$ gives $0\le T_N\le N$ on every outcome. For any fixed $p>1$,
\begin{align*}
 \mathbb ET_N^p
 &=\mathbb E[\mathbf1_{\Omega_N^{\mathrm{coer}}}T_N^p]
   +\mathbb E[\mathbf1_{(\Omega_N^{\mathrm{coer}})^c}T_N^p]\\
 &\le C_E^p\mathbb E[\mathbf1_{\Omega_N^{\mathrm{coer}}}(1+C_N)^{3p/2}]
   +N^p\mathbb P((\Omega_N^{\mathrm{coer}})^c)\\
 &\le C_{p,E,\beta}+C N^p e^{-cN}\le C'_{p,E,\beta}.
\end{align*}
Here the moment of $C_N$ follows from \eqref{W:eq:gram-Gaussian-noise-exp-moment}, without any independence between $C_N$ and the event. On $\{T_N>M\}$, $T_N\le M^{1-p}T_N^p$.
Hence
\[
 \sup_N\mathbb E[T_N\mathbf1_{\{T_N>M\}}]\le C'_{p,E,\beta}M^{1-p}\longrightarrow0.
\]
This is the uniform integrability used below. For a statistic bounded by one, its integral over the complementary coercivity event is bounded by
$\mathbb P((\Omega_N^{\mathrm{coer}})^c)$, whereas the unnormalized trace required the factor $N^p$ before using the exponentially small probability.

\Needspace{6\baselineskip}

\begin{theorem}
\label{W:thm:gaussian-resolvent-convergence}
For the Gaussian ensemble of Section~\ref{W:sec:goe-gue-tridiagonal-law}, every fixed $E\in\mathbb R$, $L\ge1$, and deterministic $\alpha_N\downarrow0$,
\begin{equation}
 \limsup_{N\to\infty}\sup_{\|w\|=1}
 \mathbb E Z_{N,L,E,\alpha_N}^{(w)}({H_N^{\mathrm G}})
 \le\frac{C_{\beta,E}}L.
 \label{W:eq:gaussian-resolvent-convergence}
\end{equation}
Thus the iterated limit of these expectations is zero. The constants here depend on the fixed energy $E$.
\end{theorem}
\begin{proof}
Let $V$ be orthogonal for $\beta=1$ and unitary for $\beta=2$, and set $H'=VHV^*$ and $w'=Vw$. Functional calculus gives
\begin{align}
 X_\rho^{(w')}(H')&=VX_\rho^{(w)}(H),&
 G_\rho^{(w')}(H')&=G_\rho^{(w)}(H),\notag\\
 B(H')&=VB(H)V^*,&F_E(H')&=VF_E(H)V^*.
 \label{W:eq:resolvent-polynomial-matrix-equivariance}
\end{align}
Substitution in the definitions of the regularized projections gives
\begin{align}
 K_{\rho,\alpha}^{(w')}(H')&=VK_{\rho,\alpha}^{(w)}(H)V^*,\notag\\
 Q_{\rho,\alpha}^{(w')}(H')&=VQ_{\rho,\alpha}^{(w)}(H)V^*.
 \label{W:eq:regularized-projection-equivariance}
\end{align}
Taking traces in \eqref{W:eq:regularized-localization-statistic}, we obtain
\begin{equation}
 C_{N,L,E,\alpha}^{(Vw)}(VHV^*)=C_{N,L,E,\alpha}^{(w)}(H),\qquad
 Z_{N,L,E,\alpha}^{(Vw)}(VHV^*)=Z_{N,L,E,\alpha}^{(w)}(H).
 \label{W:eq:localization-statistic-equivariance}
\end{equation}
These identities hold for every Hermitian $H$.

For each deterministic unit vector $w$, choose an orthogonal or unitary matrix $V$ with $Ve_1=w$. Invariance of the GOE or GUE gives
$V H_N^{\mathrm G}V^*\stackrel{\rm d}=H_N^{\mathrm G}$. Together with \eqref{W:eq:localization-statistic-equivariance}, this implies
\begin{equation}
 \mathbb E Z_{N,L,E,\alpha}^{(w)}(H_N^{\mathrm G}) =\mathbb E Z_{N,L,E,\alpha}^{(w)}(V H_N^{\mathrm G}V^*) =\mathbb E Z_{N,L,E,\alpha}^{(e_1)}(H_N^{\mathrm G}).
 \label{W:eq:resolvent-starting-vector-expectation}
\end{equation}
Thus it suffices to prove the expectation bound for $w=e_1$.

The coefficient distribution for the GOE and GUE is given by Theorem~\ref{W:thm:Gaussian-Lanczos-law}, with the independent normal and chi variables of the tridiagonal model
\cite{dumitriu2002matrix}. Proposition~\ref{W:prop:gram-Gaussian-global-coercivity} and Lemma~\ref{W:lem:gram-discrete-Airy-Weyl} supply \eqref{W:eq:resolvent-confinement-input}
with $C_N=C_0(1+\mathfrak Z_N^4)$, $\mathbb P(\Omega_N^c)\le Ce^{-cN}$, and bounded fixed moments of $C_N$. In particular, $\sup_N\mathbb E(1+C_N)^{5/2}<\infty$ by the exponential
moment of $\mathfrak Z_N^2$. Lemma~\ref{W:lem:resolvent-weighted-trace} gives
\[
 \mathbb E\bigl[\mathbf1_{\Omega_N}
 \operatorname{Tr}(f_E(A_N^G)Q_{N,L}f_E(A_N^G))\bigr]
 \le C_{\beta,E}/L.
\]
The Lanczos algorithm applied to the GOE or GUE has no breakdown almost surely. Therefore \eqref{W:eq:regularized-projection-error} and positivity give, in its original
coordinates,
\[
 C_{N,L,E,\alpha_N}^{(w)}({H_N^{\mathrm G}})
 \le \operatorname{Tr}(f_E(A_N^G)Q_{N,L}f_E(A_N^G))
 +\frac{\alpha_N}{g_N+\alpha_N}\operatorname{Tr} f_E(A_N^G)^2,
 \qquad g_N=\lambda_{\min}(G_{\rho}^{(w)}({H_N^{\mathrm G}})).
\]
Put $T_N=\operatorname{Tr}f_E(A_N^G)^2$ and $V_N=\alpha_N T_N/(g_N+\alpha_N)$. For fixed $L$, both $g_N^{-1}$ and $T_N$ are tight by Lemma~\ref{W:lem:gaussian-gram-conditioning},
\eqref{W:eq:resolvent-trace-zero}, and $\mathbb P(\Omega_N^c)\to0$. For every $M>0$,
\[
 \mathbb E\min\{V_N,1\}
 \le\alpha_NM^2+\mathbb P(g_N^{-1}>M)+\mathbb P(T_N>M).
\]
The inequalities $Z\le\min\{C,1\}$ and $\min\{u+v,1\}\le u+\min\{v,1\}$ imply
\begin{align*}
 \mathbb EZ
 &\le C_{\beta,E}/L+\mathbb E\min\{V_N,1\}
                         +\mathbb P(\Omega_N^c)\\
 &\le C_{\beta,E}/L+\alpha_NM^2
       +\mathbb P(g_N^{-1}>M)+\mathbb P(T_N>M)
       +\mathbb P(\Omega_N^c).
\end{align*}
At fixed $L$ and $M$, $\alpha_NM^2+\mathbb P(\Omega_N^c)\to0$. Tightness gives
\[
 \lim_{M\to\infty}\limsup_{N\to\infty}
   \bigl[\mathbb P(g_N^{-1}>M)+\mathbb P(T_N>M)\bigr]=0.
\]
Taking these limits in the preceding bound and using \eqref{W:eq:resolvent-starting-vector-expectation}, we obtain
\[
 \limsup_{N\to\infty}\sup_{\|w\|=1}
 \mathbb E Z_{N,L,E,\alpha_N}^{(w)}(H_N^{\mathrm G})
 \le C_{\beta,E}/L.
\]
Finally let $L\to\infty$. All estimates used to pass to the $N$ limit were taken at fixed $L$.
\end{proof}

For a deterministic sequence $\varepsilon_N=N^{-o(1)}\downarrow0$, set
\[
 \alpha_N=\exp\bigl(-\sqrt{|\log\varepsilon_N|}\bigr).
\]
Then $\alpha_N\downarrow0$, $\alpha_N=N^{-o(1)}$, and $\log(\alpha_N^{-1})/|\log\varepsilon_N| =|\log\varepsilon_N|^{-1/2}\to0$. The entry derivatives of the regularized
localization statistics are estimated in Section~\ref{W:sec:taylor-remainders-and-bounded-entries}.

\Needspace{6\baselineskip}

For a fixed deterministic sequence $w_N$, choose $\alpha_N(L)\downarrow0$ for each fixed $L$. Define
\begin{equation}
 e_N(L,E)=\left|\mathbb EZ_{N,L,E,\alpha_N(L)}^{(w_N)}(H_N)
          -\mathbb EZ_{N,L,E,\alpha_N(L)}^{(w_N)}(H_N^{\mathrm G})\right|.
 \label{W:eq:resolvent-comparison-error}
\end{equation}
The condition needed for localization is
\begin{equation}
 \lim_{L\to\infty}\limsup_N e_N(L,E)=0.
 \label{W:eq:resolvent-comparison-input}
\end{equation}

\begin{theorem}
\label{W:thm:localization-from-expectation-comparison}
For every fixed $L,E,\eta$,
\begin{equation}
 \limsup_N\mathbb P(\|Q_{N,L}\Pi_N(E)\|>\eta)
 \le\frac{4+\eta^2}{\eta^2}
          \left(\frac{C_{\beta,E}}L+\limsup_N e_N(L,E)\right).
 \label{W:eq:resolvent-transfer-bound}
\end{equation}
Consequently, if \eqref{W:eq:resolvent-comparison-input} holds, then
\begin{equation}
 \lim_{L\to\infty}\limsup_N
 \mathbb P(\|Q_{N,L}\mathbf1_{(-\infty,E]}(\mathcal A_N)\|>\eta)=0.
 \label{W:eq:resolvent-full-localization}
\end{equation}
The same limits hold uniformly over deterministic unit starting vectors when the comparison errors are uniform.
\end{theorem}
\begin{proof}
Write $c_\eta=(4+\eta^2)/\eta^2$. For every $N,L$,
\begin{align*}
 \mathbb P\bigl(\|Q_{N,L}\Pi_N(E)\|>\eta\bigr)
 &\le c_\eta\mathbb E Z_{N,L,E,\alpha_N(L)}^{(w_N)}(H_N)\\
 &\le c_\eta\left(
       \mathbb E Z_{N,L,E,\alpha_N(L)}^{(w_N)}(H_N^{\rm G})
                          +e_N(L,E)\right).
\end{align*}
The first inequality is \eqref{W:eq:regularized-localization-probability}; the second is the definition of $e_N$. Taking $\limsup_N$ and using
Theorem~\ref{W:thm:gaussian-resolvent-convergence} proves \eqref{W:eq:resolvent-transfer-bound}. The two terms on its right tend to zero as $L\to\infty$ by
\eqref{W:eq:resolvent-comparison-input}. For the uniform version, take the supremum over deterministic unit vectors in the two finite-$N$ inequalities before taking either limit;
the Gaussian expectation is independent of that vector.
\end{proof}

\subparagraph{The statistic used in comparison}
\label{W:sec:expectation-of-the-localization-statistic}
\Needspace{6\baselineskip}

\begin{proposition}
\label{W:prop:gaussian-localization-expectation}
For $H_N^{\mathrm G}$ and fixed $L\ge1,E,\kappa$,
\begin{equation}
 \limsup_N\sup_{\|w\|=1}\mathbb E\Phi_E(H_N^{\mathrm G})
       \le C_{\beta,E}/L.
 \label{W:eq:gaussian-gaussian}
\end{equation}
\end{proposition}
\begin{proof}
We use the finite-matrix definitions in \eqref{W:eq:moment-coordinates-J}, \eqref{W:eq:moment-coordinates-projection} and \eqref{W:eq:gaussian-localization}. To distinguish this
statistic from the additive regularization in \eqref{W:eq:regularized-polynomial-projection}, write its polynomial matrix as
\[
 F_r(H,w)=[p_0(h_d(H))w,\ldots,p_{r-1}(h_d(H))w],\qquad
 \mathsf G_r=F_r^*F_r,\qquad r=\rho.
\]
Its contraction is $P_{r,\kappa}^{[d]}=F_ra_\kappa(\mathsf G_r)F_r^*$, where $a_\kappa(t)=\chi_\kappa(t)/t$ with smooth zero extension. An SVD of $F_r$ shows that its eigenvalues
on the range of $F_r$ are $s_j^2a_\kappa(s_j^2)=\chi_\kappa(s_j^2)\in[0,1]$; zero singular values contribute zero. Thus it is a positive contraction for every matrix. The cutoff is
$\Gamma=\det\gamma_\kappa(\mathsf G_{\rho_T})\in[0,1]$, with $\rho_T\ge r$. If $\Gamma>0$, the eigenvalues of $\mathsf G_{\rho_T}$ lie in $(\kappa,\kappa^{-1})$. Since
$\mathsf G_r$ is its leading principal block, interlacing puts its eigenvalues in the same interval, where $a_\kappa(t)=1/t$. Consequently $F_r$ has rank $r$ and
$P_{r,\kappa}^{[d]}$ is exactly the orthogonal projection onto $\operatorname{Ran}F_r$ on this event.

Let $\Omega_d=\{h_d(H)=H\}$. On $\Omega_d\cap\{\Gamma>0\}$, this range is the Krylov space generated by $w,Hw,\ldots,H^{r-1}w$ because the polynomials have successive degrees and
nonzero leading coefficients. Their linear independence excludes breakdown before the $r$th vector. Thus the contraction equals the projection onto the first $r$ Lanczos vectors
there. For $A=m_N^2(2I-H)$ and $\theta_E=f_E^2$, the definitions give
\[
 C_E=\operatorname{Tr}\!\left[\theta_E(A)^{1/2}
             (I-P_{r,\kappa}^{[d]})\theta_E(A)^{1/2}\right]\ge0,
 \qquad \Phi_E=\Gamma\frac{C_E}{1+C_E},\qquad 0\le\Phi_E\le1.
\]
These are Borel functions of the original matrix and $w$. Indeed, on each compact set of Hermitian matrices, uniform polynomial approximation of the continuous scalar functions on
a containing spectral interval implies continuity of their matrix functional calculus. Matrix multiplication, the finite trace and the nonvanishing denominator $1+C_E$ preserve
continuity. This argument includes repeated eigenvalues. On $\{\Gamma=0\}$, $\Phi_E$ is zero exactly, without needing invertibility of any Gram matrix.

For each fixed deterministic $w$, use the recurrence coefficients obtained by Lanczos iteration for the GOE or GUE from this initial vector, as given in
Theorem~\ref{W:thm:Gaussian-Lanczos-law} and set $C_N=C_0(1+\mathfrak Z_N^4)$ and $\Omega_N^{d}=\{\|H_N^{\mathrm G}\|\le2+d_N/2\}\subset\Omega_d$. The previous projection identity
and unitary conjugation, followed by Lemma~\ref{W:lem:resolvent-weighted-trace}, imply
\begin{align*}
 \mathbf1_{\Omega_N^{\mathrm{coer}}\cap\Omega_N^d\cap\{\Gamma>0\}}\Phi_E
 &\le\mathbf1_{\Omega_N^{\mathrm{coer}}\cap\Omega_N^d\cap\{\Gamma>0\}}
       \operatorname{Tr}[f_E(\mathcal A_N^G)Q_{N,L}f_E(\mathcal A_N^G)]\\
 &\le\mathbf1_{\Omega_N^{\mathrm{coer}}\cap\Omega_N^d\cap\{\Gamma>0\}}
       \frac{C_E^{\rm const}}{L}(1+C_N)^{5/2}.
\end{align*}
The constant $C_E^{\rm const}$ is the deterministic constant of that lemma. On the same intersection with $\{\Gamma=0\}$ the left side is zero. Integrating on these events and
bounding $\Phi_E$ by one elsewhere therefore gives
\begin{align*}
 \mathbb E\Phi_E(H_N^{\mathrm G},w)
 &\le\frac{C_E^{\rm const}}L\mathbb E(1+C_N)^{5/2}
       +\mathbb P((\Omega_N^{\mathrm{coer}})^c)+\mathbb P((\Omega_N^d)^c)\\
 &\le\frac{C_{\beta,E}}L+C e^{-cN}+C e^{-c'\sqrt{m_N^2d_N}}.
\end{align*}
The moment is uniform in $N$ by \eqref{W:eq:gram-Gaussian-noise-exp-moment}. The last probability is bounded by the Gaussian norm estimate obtained after
Corollary~\ref{W:cor:gaussian-spectral-count-tails}. At the scale in \eqref{W:eq:moment-coordinates-scales}, the maximum dimension of the Gram matrices obeys
$c_*m_N\le\rho_*\le C_*m_N$ at fixed parameters of the functions under consideration, and
\[
 d_N=A_0\rho_*^{-2}(\log(N+2))^\gamma,\qquad
 m_N^2d_N\ge A_0C_*^{-2}(\log(N+2))^\gamma\longrightarrow\infty,
 \qquad 1<\gamma<2.
\]
Thus the displayed probability bound tends to zero. Its stretched logarithmic decay is used only for the bounded statistic $\Phi_E$ here. The part of this expectation on
$\{\Gamma=0\}$ is zero because the integrand contains $\Gamma$.

For completeness, uniformity in the deterministic vector follows directly in the original coordinates. If $V$ is orthogonal for $\beta=1$ or unitary for $\beta=2$, then
\begin{align*}
 F_r(VHV^*,Vw)&=VF_r(H,w),\\
 \mathsf G_r(VHV^*,Vw)&=\mathsf G_r(H,w),\\
 P_{r,\kappa}^{[d]}(VHV^*,Vw)&=VP_{r,\kappa}^{[d]}(H,w)V^*.
\end{align*}
The same identity for the Gram matrix holds at dimension $\rho_T$, leaving $\Gamma$ unchanged, while functional calculus conjugates $\theta_E(A)$ by $V$. Cyclicity of the finite
trace therefore gives $\Phi_E(VHV^*,Vw)=\Phi_E(H,w)$. For any two deterministic unit vectors choose $Vw=w'$ and use orthogonal or unitary invariance of the GOE or GUE to obtain
equality of the two expectations. The supremum over such vectors is consequently the expectation for any one of them; a common coercivity event for all starting vectors is not
needed. Taking $\limsup_N$ at fixed $L,E,\kappa$ proves the expectation bound. The constant $C_{\beta,E}$ is independent of $\kappa$ within its admissible range: the cutoff in the
eigenvalues of the Gram matrix is at most one and enforces exact projection wherever it is positive.
\end{proof}

The independence of $C_{\beta,E}$ from $\kappa$ comes from the support of $\Gamma$. On $\{\Gamma>0\}$, the regularized function of the Gram matrix equals its inverse; on
$\{\Gamma=0\}$, the statistic $\Phi_E$ vanishes. Thus the inverse of the original Gram matrix is evaluated only where $\Gamma>0$; the statistic is defined to be zero elsewhere.

For an original or truncated matrix $H_N$ and a measurable $0\le U\le1$, put
\begin{align*}
 e_{1,N}&=|\mathbb E[U\Gamma](H_N)-\mathbb E[U\Gamma](H_N^{\mathrm G})|,\\
 e_{2,N}&=|\mathbb E[U\Gamma Z_E](H_N)-\mathbb E[U\Gamma Z_E](H_N^{\mathrm G})|,\\
 d_N&=\mathbb E[1-U](H_N^{\mathrm G}),\qquad
 g_N=\mathbb E[1-\Gamma](H_N^{\mathrm G}),\qquad
 p_{d,N}=\mathbb P_{H_N}(\Omega_d^c).
\end{align*}
These symbols refer to the respective matrix laws: the same bounded measurable functions $U$, $\Gamma$ and $Z_E$ are evaluated under each law, with the deterministic starting
vector fixed. Their products are evaluated under their joint distribution.
\begin{corollary}
\label{W:cor:localization-comparison-errors}
For every $\eta>0$,
\begin{equation}
 \mathbb P_{H_N}(\|Q_{N,L}\Pi_N(E)\|>\eta) \le{}p_{d,N}+d_N+g_N+e_{1,N} +\frac{4+\eta^2}{\eta^2} \bigl(e_{2,N}+\mathbb E\Phi_E(H_N^{\mathrm G})\bigr).
 \label{W:eq:localization-comparison-error-bound}
\end{equation}
\end{corollary}
\begin{proof}
Write $c_\eta=(4+\eta^2)/\eta^2$ and let $\mathcal B=\{\|Q_{N,L}\Pi_N(E)\|>\eta\}$. On $\Omega_d\cap\{\Gamma>0\}$ the exact Krylov projection and $\theta_E(x)\ge1/4$ for $x\le E$
give $\|Q_{N,L}\Pi_N(E)\|^2\le4C_E$, as in \eqref{W:eq:localization-trace-bound}. Thus $\mathbf1_{\mathcal B}\le c_\eta Z_E$ on this intersection. On $\Omega_d$, including the case
$\Gamma=0$, multiplication by $U\Gamma$ and $0\le U\Gamma\le1$ give
\[
 \mathbf1_{\mathcal B}
 =(1-U\Gamma)\mathbf1_{\mathcal B}+U\Gamma\mathbf1_{\mathcal B}
 \le1-U\Gamma+c_\eta U\Gamma Z_E.
\]
Adding $\mathbf1_{\Omega_d^c}$ proves the pointwise inequality on the whole matrix space. All its terms are bounded measurable functions, so taking expectations under the law of
$H_N$ gives
\[
 \mathbb P_{H_N}(\mathcal B)
 \le p_{d,N}+\mathbb E_{H_N}(1-U\Gamma)
                 +c_\eta\mathbb E_{H_N}(U\Gamma Z_E).
\]
The definition of $e_{1,N}$ and the identity $1-U\Gamma=(1-U)+U(1-\Gamma)$ imply
\begin{align*}
 \mathbb E_{H_N}(1-U\Gamma)
 &\le e_{1,N}+\mathbb E_{H_N^{\mathrm G}}[(1-U)+U(1-\Gamma)]\\
 &\le e_{1,N}+d_N+g_N.
\end{align*}
This uses $U\le1$ pointwise and requires no independence between $U$ and $\Gamma$. Likewise, $e_{2,N}$ and $0\le U\le1$ give
\[
 \mathbb E_{H_N}(U\Gamma Z_E)
 \le e_{2,N}+\mathbb E_{H_N^{\mathrm G}}(U\Phi_E)
 \le e_{2,N}+\mathbb E\Phi_E(H_N^{\mathrm G}).
\]
Substitution proves \eqref{W:eq:localization-comparison-error-bound}. The two expectation comparisons are between separate laws and use no coupling of $H_N$ with $H_N^{\mathrm G}$.
\end{proof}
In the later application, $N$ tends to infinity with $L,E,\kappa$ fixed. Gaussian confinement makes the error caused by modifying the eigenvalues vanish, and
Lemma~\ref{W:lem:comparison-ensemble-coefficients} controls the common coefficient cutoff. Lemma~\ref{W:lem:gaussian-gram-conditioning} gives
$\lim_{\kappa\downarrow0}\limsup_N g_N=0$ at each fixed $L$; this is not a polynomial failure estimate at fixed $\kappa$. Proposition~\ref{W:prop:gaussian-localization-expectation}
supplies the remaining bound $C_{\beta,E}/L$ uniformly in $\kappa$. The limits are therefore taken in the order $N$, $\kappa\downarrow0$, $L\to\infty$. The probabilities $p_{d,N}$
and errors $e_{1,N},e_{2,N}$ for the original or truncated matrix are estimated in the separate confinement and comparison arguments that follow. Those estimates require finite
size information beyond the operator limit for the Gaussian ensembles.

\section{Sharp convergence rates and the leading order coefficient}\label{sec:sharp-convergence-rates-and-the-limiting-error-coefficient}

The coefficient estimates and discrete approximations now give the quantitative theorems. Proposition~\ref{prop:shared-resolvent-transfer} supplies the inverse identities; the
arguments below establish uniform bounds on the resolvent factors on the constructed probability spaces.

The lower bounds under rank and observation constraints are stated and proved below. Their quantifiers differ: the bound under a rank constraint holds on one event of probability
one for every admissible approximation, whereas the bound under an observation constraint applies to each fixed sequence measurable from the specified Brownian observations and
independent auxiliary data. In the latter case, the event of probability one may depend on that sequence. The coefficient of the best approximation subject to the rank constraint
need not equal the coefficient for a prescribed matrix law or embedding. For the constructed coupling of the bidiagonal matrices, both error contributions can have the same order,
so their mixed inner product is retained.

For the circular ensemble, the comparison of the circular inverse with the inverse obtained by sampling the Brownian path is of smaller order than the sampling error. The latter is
reduced to a Gaussian term with a nonzero conditional mean. The conditional second moment and an almost sure concentration estimate together identify its leading order coefficient.
Bounds on the mixed terms then give the same coefficient for the circular inverse. All rate statements use their stated fixed parameters; the later pathwise regularity in $\beta$
requires a separate argument.

\subsection{Lower bound under a rank constraint}
\label{sec:the-airy-rank-obstruction}
The lower bound does not depend on how the matrix coefficients are coupled. It follows from the eigenvalues of the limiting operator and applies to every approximation of the rank
specified in Proposition~\ref{A:res:5.1}.

\subsubsection{The lower bound}\label{A:sec:rank-lower-bound-and-optimality}

We first derive the eigenvalue asymptotic needed in the proof. Summing the reciprocal spectral tail then gives the lower bound of Proposition~\ref{A:res:5.1}. The final calculation
explains how the fixed range of an embedding contributes to the same approximation error.

\begin{lemma}\label{A:res:weyl}
Almost surely, with eigenvalues in increasing order,
\[
 N_\beta(E):=\#\{k:\Lambda_k\le E\}
 \sim\frac{2}{3\pi}E^{3/2},\qquad
 \Lambda_k\sim\left(\frac{3\pi k}{2}\right)^{2/3}.
\]
\end{lemma}
\begin{proof}
We first prove the counting asymptotic for $\mathcal L=-\partial_x^2+1+x$ with Dirichlet boundary condition at zero. Partition the half-line into $I_j=(j,j+1)$, $j\ge0$. On $I_j$
the potential lies between $1+j$ and $2+j$. Imposing a Dirichlet condition at every integer restricts the form domain. Conversely, allowing independent $H^1(I_j)$ functions on the
intervals enlarges the domain and gives Neumann boundary conditions. The min--max principle therefore gives, for $E\ge2$,
\[
 \sum_{0\le j\le E-2}
   \left\lfloor\frac{\sqrt{E-2-j}}{\pi}\right\rfloor
 \le N_{\mathcal L}(E)
 \le\sum_{0\le j\le E-1}
   \left(1+\left\lfloor\frac{\sqrt{E-1-j}}{\pi}\right\rfloor\right).
\]
Indeed, the Dirichlet eigenvalues of the Laplacian on a unit interval are $\pi^2l^2$, $l\ge1$, whereas its Neumann eigenvalues have the same formula with $l\ge0$. Intervals with
$1+j>E$ make no contribution to the upper bound. The error from the integer parts and the Neumann zero modes is $O(E)$. For $c=1,2$, monotonicity of $t\mapsto(E-c-t)_+^{1/2}$ gives
\[
 \sum_{j\ge0}(E-c-j)_+^{1/2}
 =\int_0^\infty(E-c-t)_+^{1/2}\,\mathrm dt+O(E^{1/2})
 =\frac23(E-c)_+^{3/2}+O(E^{1/2}).
\]
Consequently $N_{\mathcal L}(E)=2E^{3/2}/(3\pi)+O(E)$.

Fix a path in $\Omega_B$. The relative form estimate in Lemma~\ref{A:res:3.1}, applied to $2B/\sqrt\beta$, gives, for every rational $0<\varepsilon<1$, a finite $C_\varepsilon$
such that
\[
 (1-\varepsilon)\mathcal L-C_\varepsilon I
 \le\mathcal H_\beta
 \le(1+\varepsilon)\mathcal L+C_\varepsilon I
\]
on the common form domain. If $\nu_k$ are the increasing eigenvalues of $\mathcal L$, min--max gives
$(1-\varepsilon)\nu_k-C_\varepsilon\le\Lambda_k\le(1+\varepsilon)\nu_k+C_\varepsilon$.
Thus
\[
 N_{\mathcal L}\!\left(\frac{E-C_\varepsilon}{1+\varepsilon}\right)
 \le N_\beta(E)
 \le N_{\mathcal L}\!\left(\frac{E+C_\varepsilon}{1-\varepsilon}\right).
\]
Put $c_0=2/(3\pi)$. Divide by $E^{3/2}$ and let $E\to\infty$ with $\varepsilon$ fixed. This yields
\[
 c_0(1+\varepsilon)^{-3/2}
 \le\liminf_{E\to\infty}\frac{N_\beta(E)}{E^{3/2}}
 \le\limsup_{E\to\infty}\frac{N_\beta(E)}{E^{3/2}}
 \le c_0(1-\varepsilon)^{-3/2}.
\]
Letting rational $\varepsilon\downarrow0$ proves the counting asymptotic on this same path.

For its inversion, fix $0<\delta<c_0$. There is $E_\delta<\infty$ such that $(c_0-\delta)E^{3/2}\le N_\beta(E)\le(c_0+\delta)E^{3/2}$ for every $E\ge E_\delta$. For all $k$ for
which both quantities below exceed $E_\delta$, the lower energy has at most $k-1$ eigenvalues and the upper energy has at least $k$ eigenvalues. Hence
\[
 \left(\frac{k-1}{c_0+\delta}\right)^{2/3}
 <\Lambda_k\le\left(\frac{k}{c_0-\delta}\right)^{2/3}.
\]
After division by $k^{2/3}$, let $k\to\infty$ and then $\delta\downarrow0$. This proves $\Lambda_k\sim c_0^{-2/3}k^{2/3}$. This calculation uses no estimate for eigenvalue gaps or
inverse-gap moments.
\end{proof}

Lemma~\ref{A:res:weyl} gives, for a sufficiently large finite random $\kappa_\omega$, $c_\omega k^{2/3}\le\Lambda_k+\kappa_\omega\le C_\omega k^{2/3}$ for every $k$. Consequently,
uniformly for $z$ in a fixed compact subset of $\mathbb C\setminus\mathbb R$,
\[
\sum_{k>m}|\Lambda_k-z|^{-2}\asymp_{\omega,K}m^{-1/3}.
\eqtag{5.3}
\]
Choose $k_0$ so that $\Lambda_k\ge2\sup_{z\in K}|z|$ and $c_\omega k^{2/3}\le\Lambda_k\le C_\omega k^{2/3}$ for $k\ge k_0$. Then $|\Lambda_k-z|\asymp\Lambda_k$ uniformly on $K$,
while
\[
 3(m+1)^{-1/3}\le\sum_{k>m}k^{-4/3}\le3m^{-1/3}.
\]
This proves \eqref{A:eq:5.3} with constants uniform in $z$. The same proof applies to a fixed real $a\in\rho(\mathcal H_\beta)$.

To pass from the eigenvalue asymptotic to the approximation error, we use the following identity. For a Hilbert--Schmidt operator $T$ with singular values $s_1\ge s_2\ge\cdots$,
\[
\inf_{\operatorname{rank}A\le m}\|A-T\|_{\mathfrak S_2}^2
=\sum_{k>m}s_k(T)^2.
\eqtag{5.4}
\]
For completeness, let $A$ be any bounded operator of rank at most $m$, and let $Q$ be the orthogonal projection onto its range. Since $(I-Q)A=0$, orthogonality of
$\operatorname{Ran}Q$ and $\ker Q$ gives
\[
 \|T-A\|_{\mathfrak S_2}^2
 =\|(I-Q)T\|_{\mathfrak S_2}^2+\|QT-A\|_{\mathfrak S_2}^2
 \ge\|T\|_{\mathfrak S_2}^2-\|QT\|_{\mathfrak S_2}^2.
\]
Choose orthonormal left singular vectors $u_j$ for the nonzero singular values of $T$, and put $a_j=\|Qu_j\|^2$. Then $0\le a_j\le1$ and
\[
 \sum_j a_j\le\operatorname{rank}Q\le m,
 \qquad
 \|QT\|_{\mathfrak S_2}^2=\sum_j s_j(T)^2a_j.
\]
The first inequality follows by expanding $Q$ in an orthonormal basis of its finite-dimensional range and applying Bessel's inequality to each basis vector. If $T$ has at least
$m\ge1$ nonzero singular values, their ordering gives
\[
 \begin{aligned}
 \sum_j s_j^2a_j
 &\le\sum_{j\le m}s_j^2a_j+s_m^2\left(m-\sum_{j\le m}a_j\right)\\
 &\le\sum_{j\le m}s_j^2.
 \end{aligned}
\]
If $T$ has fewer than $m$ nonzero singular values, the same upper bound is simply $\|T\|_{\mathfrak S_2}^2$, with zeros appended. We have proved
$\|T-A\|_{\mathfrak S_2}^2\ge\sum_{j>m}s_j(T)^2$ without any self-adjointness or positivity assumption on $A$. Conversely, if $Tv_j=s_ju_j$ is a singular value decomposition, the
operator $A_m f=\sum_{j\le m}s_j\langle f,v_j\rangle u_j$ has rank at most $m$ and error squared exactly $\sum_{j>m}s_j^2$. For $m=0$ take $A_0=0$. This proves \eqref{A:eq:5.4}.

The singular values of $R_\beta(z)$ are the decreasing rearrangement of $|\Lambda_k-z|^{-1}$. Uniformly for $z\in K$, this rearrangement changes only a finite initial set. To see
this, choose $k_0\ge2$ so that $\Lambda_k>\sup_{z\in K}\operatorname{Re}z$ for $k\ge k_0$. For each $z\in K$, the function $\lambda\mapsto|\lambda-z|$ is increasing on this
interval. Moreover,
\[
 \min_{k<k_0}\inf_{z\in K}|\Lambda_k-z|^{-1}>0,
 \qquad
 \sup_{k\ge M}\sup_{z\in K}|\Lambda_k-z|^{-1}\longrightarrow0
 \quad(M\to\infty).
\]
Thus, for some finite $M$, the first $M$ singular values are precisely those with indices at most $M$, with a possible rearrangement among them. For $m\ge M$, the singular value
tail equals the sum in \eqref{A:eq:5.3}. Hence, for every $z\in K$ and every operator $A$ of rank at most $m$,
\[
 \|A-R_\beta(z)\|_{\mathfrak S_2}^2\ge c_{\omega,K}m^{-1/3}.
\]
The estimate applies separately to every $z$ and every permitted operator $A$, including an arbitrary choice $A=A(z)$.

More precisely, the Weyl asymptotic and the finite rearrangement just proved give, uniformly for $z\in K$,
\[
 |\Lambda_k-z|^{-2}=\left(\frac2{3\pi}\right)^{4/3}
                  k^{-4/3}(1+o_{B,K}(1)).
\]
The integral bounds above imply $\sum_{k>m}k^{-4/3}=3m^{-1/3}+O(m^{-4/3})$. Hence
\[
 \inf_{\operatorname{rank}A\le m}
       \|A-R_\beta(z)\|_{\mathfrak S_2}^2
 =3\left(\frac2{3\pi}\right)^{4/3}
                  m^{-1/3}(1+o_{B,K}(1)).
 \eqtag{rank-tail-asymptotic}
\]
The spectral truncation of $R_\beta(z)$ attains this minimum, but need not have the coefficient law or the range of a finite matrix in the theorem. The constructed coupling attains
the same order; its leading order coefficient is not identified by this calculation.

\begin{proposition}[Lower bound]\label{A:res:5.1}
Fix $\beta>0$ and a nonempty compact $K\Subset\mathbb C\setminus\mathbb R$. Let $\mathscr R_n(z)$ be any embedded resolvent of rank at most $m_n$, where $m_n=n-O(\log n)$. On an
event of probability one determined by the limiting operator, there is $c_{\omega,K}>0$ such that, for all sufficiently large $n$,
\[
\inf_{z\in K}
\|\mathscr R_n(z)-R_\beta(z)\|_{\mathfrak S_2}^2
\ge c_{\omega,K}n^{-1/3}.
\eqtag{5.6}
\]
The same event and constant apply to every such approximation. The approximating operators may depend on the whole limiting path, additional randomness, and $z$.
\end{proposition}

\begin{proof}
For $z\in K$, the rank bound and \eqref{A:eq:5.4} give
\[
 \|\mathscr R_n(z)-R_\beta(z)\|_{\mathfrak S_2}^2
 \ge\sum_{j>m_n}s_j(R_\beta(z))^2
 \ge c_{\omega,K}m_n^{-1/3}.
\]
The last inequality is \eqref{A:eq:5.3}, with its constant uniform on $K$. Since $m_n=n-O(\log n)$, choose a deterministic $n_0$ such that $m_n\le2n$ for $n\ge n_0$; and hence
$m_n^{-1/3}\ge2^{-1/3}n^{-1/3}$. Absorb this fixed factor into the constant. The event consists of the Weyl and resolvent tail estimates for the limiting operator. Once a path in
that event is fixed, \eqref{A:eq:5.4} holds for every bounded operator of the permitted rank, without any condition on its dependence on the path, additional randomness, or $z$.
Thus the lower bound holds on one event for every permitted approximation.
\end{proof}

For the embedding there is a second, purely geometric comparison. Fix its range projection $P_n$ and suppose $\mathscr R_n=P_n\mathscr R_nP_n$. In the Hilbert--Schmidt inner
product, $P_nR_\beta P_n$ is the orthogonal projection of $R_\beta$ onto operators supported on that range. The resulting error decomposition separates the range restriction from
the remaining matrix comparison. Set $A=\mathscr R_n-P_nR_\beta P_n$ and $D=P_nR_\beta P_n-R_\beta$. Then $A=P_nAP_n$ and $P_nDP_n=0$. Cyclicity of the trace for products of
Hilbert--Schmidt operators gives
\[
 \operatorname{Tr}(AD^*)
 =\operatorname{Tr}(AP_nD^*P_n)=0.
\]
Expanding the square therefore yields the exact orthogonal decomposition
\[
 \|\mathscr R_n-R_\beta\|_{\mathfrak S_2}^2
 =\|\mathscr R_n-P_nR_\beta P_n\|_{\mathfrak S_2}^2
   +\|P_nR_\beta P_n-R_\beta\|_{\mathfrak S_2}^2.
 \eqtag{23.1}
\]

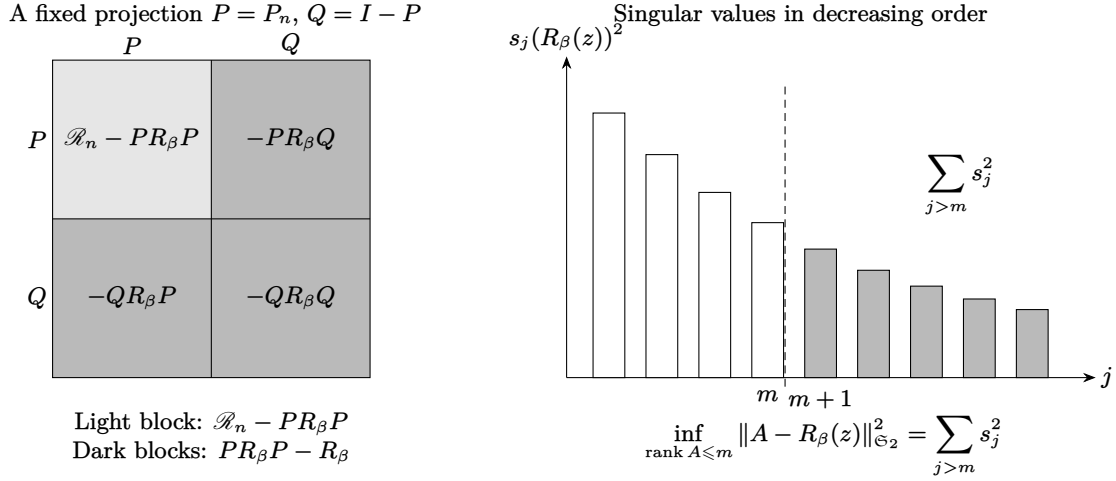
\begin{figure}[H]
\centering
\begin{tikzpicture}[x=1cm,y=1cm,every node/.style={font=\footnotesize,inner sep=2pt},>=Stealth]
 \node at (2.15,4.8) {A fixed projection $P=P_n$, $Q=I-P$};
 \fill[black!10] (0,2.1) rectangle (2.1,4.2);
 \fill[black!27] (2.1,0) rectangle (4.2,4.2);
 \fill[black!27] (0,0) rectangle (2.1,2.1);
 \draw (0,0) rectangle (4.2,4.2);
 \draw (2.1,0)--(2.1,4.2); \draw (0,2.1)--(4.2,2.1);
 \node at (1.05,3.15) {$\mathscr R_n-PR_\beta P$};
 \node at (3.15,3.15) {$-PR_\beta Q$};
 \node at (1.05,1.05) {$-QR_\beta P$};
 \node at (3.15,1.05) {$-QR_\beta Q$};
 \node[above] at (1.05,4.2) {$P$}; \node[above] at (3.15,4.2) {$Q$};
 \node[left] at (0,3.15) {$P$}; \node[left] at (0,1.05) {$Q$};
 \node[align=center] at (2.1,-.8) {Light block: $\mathscr R_n-PR_\beta P$\\Dark blocks: $PR_\beta P-R_\beta$};
 \node at (9.9,4.8) {Singular values in decreasing order};
 \draw[->] (6.8,0)--(13.8,0) node[right] {$j$};
 \draw[->] (6.8,0)--(6.8,4.25) node[above] {$s_j(R_\beta(z))^2$};
 \foreach \x/\h in {7.15/3.5,7.85/2.95,8.55/2.45,9.25/2.05}{
  \draw[fill=white] (\x,0) rectangle +(0.42,\h);
 }
 \foreach \x/\h in {9.95/1.7,10.65/1.42,11.35/1.21,12.05/1.04,12.75/.9}{
  \draw[fill=black!27] (\x,0) rectangle +(0.42,\h);
 }
 \draw[densely dashed] (9.69,-.12)--(9.69,3.85);
 \node[below] at (9.46,-.1) {$m$};
 \node[below] at (10.16,-.1) {$m+1$};
 \node at (12.0,2.55) {$\displaystyle\sum_{j>m}s_j^2$};
 \node at (10.2,-.9) {$\displaystyle\inf_{\operatorname{rank}A\le m}\|A-R_\beta(z)\|_{\mathfrak S_2}^2=\sum_{j>m}s_j^2$};
\end{tikzpicture}
\caption{Two decompositions used to prove the lower bound under a rank constraint. On the left, $\mathscr R_n=P\mathscr R_nP$, and the four blocks of $\mathscr R_n-R_\beta$ are orthogonal in the Hilbert--Schmidt inner product. The light block gives the first summand of \eqref{A:eq:23.1}; the other three give the second. On the right, the singular values are drawn schematically. Their tail is the minimum error among all operators of rank at most $m$, as in \eqref{A:eq:5.4}. The Weyl estimate gives a tail of order $m^{-1/3}$, uniformly on the compact spectral set in Proposition~\ref{A:res:5.1}.}
\label{A:fig:resolvent-approximation-components}
\end{figure}

Figure~\ref{A:fig:resolvent-approximation-components} distinguishes the error for a prescribed range from the minimum over all ranges of a given dimension. The block decomposition
explains the orthogonality in \eqref{A:eq:23.1}, while the singular value tail gives the lower bound in \eqref{A:eq:5.4}.

For the piecewise constant embedding, put $h=h_n$, $N=m_n$, and $P_n=\iota_h\iota_h^*$. On each interval $J_j=[(j-1)h,jh)$, $P_nf$ is the average $f_{J_j}$, and it vanishes on
$[Nh,\infty)$. For $f\in\mathcal D_{\mathrm A}$,
\[
\int_{J_j}|f-f_{J_j}|^2 =\frac1{2h}\int_{J_j}\int_{J_j}|f(x)-f(y)|^2\,\mathrm dx\,\mathrm dy
 \le\frac{h^2}{2}\int_{J_j}|f'(t)|^2\,\mathrm dt.
\]
Here $|f(x)-f(y)|^2\le|x-y|\int_{J_j}|f'|^2$ was used in the last step. Since $Nh\asymp h^{-2}$,
\[
\|(I-P_n)f\|_2^2 \le\frac{h^2}{2}\|f'\|_2^2+ (Nh)^{-1}\int_{Nh}^\infty x|f(x)|^2\,\mathrm dx \le
 Ch^2\left(\|f'\|_2^2+\int_0^\infty(1+x)|f|^2\,\mathrm dx\right).
\]
This is the energy approximation estimate needed for the same spectral summation as in Lemma~\ref{A:res:4.3C}. Here the projection is onto step functions and also removes the
half-line beyond $Nh$; the preceding calculation accounts for both effects. For a normalized eigenfunction $\psi_k$ of $\mathcal H_\beta$, the form comparison and
$\|I-P_n\|_{\mathrm{op}}\le1$ give
\[
 \|(I-P_n)\psi_k\|_2^2\le
 C_\omega\min\{1,h^2(\Lambda_k+\kappa_\omega)\}.
\]
Also $(\Lambda_k+\kappa_\omega)/|\Lambda_k-z|\le C_{\omega,K}$ and $\Lambda_k+\kappa_\omega\ge c_\omega k^{2/3}$. Put $k_h=\lceil h^{-3}\rceil$. Expansion in the eigenbasis gives
\[
\begin{aligned}
 \|(I-P_n)R_\beta(z)\|_{\mathfrak S_2}^2
 &\le C_{\omega,K}\sum_{k\ge1}
       \min\{k^{-4/3},h^2k^{-2/3}\}\\
 &\le C_{\omega,K}\left(h^2\sum_{k\le k_h}k^{-2/3}
                          +\sum_{k>k_h}k^{-4/3}\right)
 \le C_{\omega,K}h.
\end{aligned}
\]
The estimate at $\bar z$ bounds $R_\beta(z)(I-P_n)$ by the same constant. Since
\[
 R_\beta-P_nR_\beta P_n
 =(I-P_n)R_\beta+P_nR_\beta(I-P_n),
\]
we obtain $\sup_{z\in K}\|R_\beta-P_nR_\beta P_n\|_{\mathfrak S_2}^2 \le C_{\omega,K}h$. Its lower bound is $c_{\omega,K}h$ by \eqref{A:eq:5.4}, because
$\operatorname{rank}(P_nR_\beta P_n)\le N$ and $N\asymp h^{-3}$. Thus the fixed range alone already causes squared error of order $h$. The orthogonal decomposition
\eqref{A:eq:23.1} does not say that the other summand is $o(h)$; it does not identify a leading order coefficient for the coupled matrix error.

The same proof at the deterministic point zero gives the lower bound in Theorem~\ref{A:res:1.1}. The passage from the nonreal resolvent estimate to the inverse, together with the
remaining spectral consequences, is carried out in Subsections~\ref{A:sec:atomlessness-real-resolvents-and-inverse-coupling} and~\ref{A:sec:spectral-consequences}.

\subsection{Resolvents and convergence of inverses for the tridiagonal models}
\label{sec:airy-resolvents-atomlessness-and-inverse-coupling}
\label{A:sec:inverse-coupling-and-spectral-consequences}

The coefficient and finite difference comparisons first give the resolvent theorem, Theorem~\ref{A:res:1.2}, at nonreal parameters. Recovering the inverse at zero requires an
additional step: Theorem~\ref{A:res:20.1} excludes each fixed real number from the limiting spectrum. Theorem~\ref{A:res:20.3} then transfers the comparison to real parameters for
the two isometric pairs, while Proposition~\ref{A:res:real-terminal-removal} treats the full general Jacobi matrix. Corollary~\ref{A:res:20.4} at zero yields
Theorem~\ref{A:res:1.1}. The spectral conclusions also use the lower bound of Proposition~\ref{A:res:common-lower}.

\subsubsection{Proof of the resolvent theorem}\label{A:sec:proof-of-the-resolvent-theorem}

\begin{proof}[Proof of Theorem~\ref{A:res:1.2}]
Fix the model and $\beta$. Lemma~\ref{A:res:B.8} places the couplings of Theorems~\ref{A:res:8.5} and~\ref{A:res:17.6} on one probability space. In the general case use $s_\beta$,
so the comparison map is exactly $\mathcal J_n$ of \eqref{A:eq:1.1}. Neither this construction nor its Brownian path depends on the spectral parameter. For $m\ge2$, set
\[
 K_m=\{z\in\mathbb C:|z|\le m,\ |\operatorname{Im}z|\ge m^{-1}\}.
\]
Intersect the events for these countably many sets with $\Omega_{\mathrm{DE}}$ or $\Omega_{V,\beta}$, respectively. The resulting event has probability one. Every compact subset of
$\mathbb C\setminus\mathbb R$ is contained in some $K_m$, so the estimates hold on one event for all such compact sets. Only their constants and starting sizes vary with the set.
In particular, no uncountable intersection over spectral parameters is taken.

For the Gaussian ensemble, put $h=n^{-1/3}$ and $J=\iota_h^{(n)}$. The difference in the theorem is exactly
\[
 J(\mathcal A_{n,\beta}^{\mathrm{DE}}-z)^{-1}J^*-R_\beta(z) =J\bigl[(\mathcal
 A_{n,\beta}^{\mathrm{DE}}-z)^{-1} -(\mathcal B_{n,h,\beta}^B-z)^{-1}\bigr]J^* +J(\mathcal
 B_{n,h,\beta}^B-z)^{-1}J^*-R_\beta(z).
\]
Equation~\eqref{A:eq:isometric-hs} applies because $J^*J=I$. Theorem~\ref{A:res:8.5} bounds the first summand, and Theorem~\ref{A:res:1.3} bounds the second. Thus
\[
 \begin{split}
 &\sup_{z\in K}
 \|J(\mathcal A_{n,\beta}^{\mathrm{DE}}-z)^{-1}J^*-R_\beta(z)\|_{\mathfrak S_2}^2\\
 &\quad\le 2\sup_{z\in K}
 \|(\mathcal A_{n,\beta}^{\mathrm{DE}}-z)^{-1}
                -(\mathcal B_{n,h,\beta}^B-z)^{-1}\|_{\mathfrak S_2}^2
       +2\sup_{z\in K}
 \|J(\mathcal B_{n,h,\beta}^B-z)^{-1}J^*-R_\beta(z)\|_{\mathfrak S_2}^2
 \le C_{\omega,\beta,K}h.
 \end{split}
\]
For the general Jacobi matrix, let $E_n=P_{I_n}^*U_{r_n}^*$ and $I=\iota_h^{(M_n)}$. Then $E_n^*E_n=I_{M_n}$ and $\mathcal J_n=IE_n^*$. Therefore
\[
 \mathcal J_n(\mathcal A_{n,V,\beta}-z)^{-1}\mathcal J_n^* -I(\mathcal
 B_{M_n,h,\beta}^B-z)^{-1}I^* =IE_n^*\left[(\mathcal A_{n,V,\beta}-z)^{-1} -E_n(\mathcal
 B_{M_n,h,\beta}^B-z)^{-1}E_n^*\right]E_nI^*.
\]
The outer maps have norm at most one. Theorem~\ref{A:res:17.6} bounds the squared Hilbert--Schmidt norm by $C_{\omega,V,\beta,K}h$. Since $M_n=n-O(\log n)$ and
$h^3=\tau_\star/(b_\star n)$, we have $M_nh\asymp h^{-2}$. Theorem~\ref{A:res:1.3} and the squared triangle inequality now give
\[
 \sup_{z\in K}\|\mathcal J_n(\mathcal A_{n,V,\beta}-z)^{-1}
                    \mathcal J_n^*-R_\beta(z)\|_{\mathfrak S_2}^2
 \le C_{\omega,V,\beta,K}h.
\]

For the lower bound, use only the rank of the embedded finite resolvent: it is at most $n$ in the Gaussian case and at most $M_n=n-O(\log n)$ in the general case. The approximation
inequality \eqref{A:eq:5.4} therefore gives, in either case,
\[
 \|\mathscr I_n(\mathscr A_n-z)^{-1}\mathscr I_n^*-R_\beta(z)\|_{\mathfrak S_2}^2
 \ge\sum_{j>m_n}s_j(R_\beta(z))^2
 \ge c_{\omega,K}m_n^{-1/3}
 \ge c'_{\omega,K}h_n,
\]
where $m_n=n$ or $M_n$. Proposition~\ref{A:res:5.1} makes the bound uniform on $K$. Since $h_n\asymp n^{-1/3}$, the theorem follows.
\end{proof}

\subsubsection{Atomlessness and resolvents at real parameters}\label{A:sec:atomlessness-real-resolvents-and-inverse-coupling}

To define the limiting inverse at zero, we must first exclude zero from the spectrum. More generally, we prove that every fixed deterministic real number almost surely avoids the
spectrum. The argument varies one independent Gaussian coordinate of the potential in a nonnegative direction; each eigenvalue then moves strictly. The simplicity and strictness
needed for this argument follow from the weak differential equation, without an assumption on eigenvalue densities.

Fix $B\in\Omega_B$. If $u$ is an eigenfunction, its weak equation is
\[
 -u''+xu+c_\beta B'u=\lambda u.
 \eqtag{20.7}
\]
Here $B'u$ denotes $(Bu)'-Bu'$ as a distribution. This agrees with the form definition, since testing it against a smooth compactly supported function $f$ gives
$-\int B(u\overline f)'$.

The weak equation supplies the regularity needed at the boundary point. On $(0,R)$, $u\in H^1(0,R)$ and $B$ is bounded. Subtracting $c_\beta Bu$ from $u'$ removes the
distributional derivative of $B$. Thus the function $v=u'-c_\beta Bu$ belongs to $L^2(0,R)$ and satisfies, in distributions,
\[
 v'=(x-\lambda)u-c_\beta Bu'\in L^2(0,R).
\]
Consequently $v\in H^1(0,R)$. The continuous representatives of $u$ and $v$ satisfy
\[
 \begin{cases}
 u'=v+c_\beta Bu,\\[1mm]
 v'=(x-\lambda-c_\beta^2B^2)u-c_\beta Bv.
 \end{cases}
 \eqtag{20.11}
\]
In particular both have traces at zero, and the right sides are continuous. The integral equation for this system has a unique solution given $(u(0),v(0))$: on $[0,R]$ the
difference of two solutions with the same initial data is bounded by a constant times the integral of its norm, and Gronwall's inequality makes it zero. The Dirichlet condition
fixes $u(0)=0$. If also $v(0)=0$, the solution vanishes. For two eigenfunctions at the same $\lambda$, a linear combination cancels their values of $v(0)$ and hence vanishes. Thus
each eigenspace has dimension one. If an eigenfunction vanishes on an open interval, then $u'=v=0$ there; the same uniqueness argument gives $u=0$ everywhere.

These conclusions also hold for $H(t)=\mathcal H_\beta+tM_q$, with $q\in L^\infty(\mathbb R_+;\mathbb R)$. The second equation then has the additional term $tqu$. Its coefficients
are locally integrable, which is sufficient for the same integral equation and uniqueness proof. Since $M_q$ is bounded and self-adjoint, $H(t)$ has domain
$\operatorname{Dom}\mathcal H_\beta$, is bounded below, and has compact resolvent. For example, at $z$ with $|\operatorname{Im}z|>|t|\|q\|_\infty$, its resolvent is the product of
the compact resolvent of $\mathcal H_\beta$ and a convergent Neumann series. The min--max principle gives
\[
 |\lambda_k(t)-\lambda_k(s)|\le\|q\|_\infty|t-s|
 \qquad(s,t\in\mathbb R).
\]

We give the local perturbation argument on the common operator domain. Fix $t_0\in\mathbb R$ and a simple eigenvalue $\lambda_k(t_0)$. Let $\Gamma$ be a positively oriented circle
enclosing this eigenvalue and no other spectral point of $H(t_0)$, and write $R_0(z)=(z-H(t_0))^{-1}$. Choose $\delta>0$ so that
\[
 \delta\|q\|_\infty\sup_{z\in\Gamma}\|R_0(z)\|<\frac12;
\]
if $q=0$, any sufficiently small $\delta$ will do. For complex $t$ with $|t-t_0|<\delta$,
\[
 (z-H(t))^{-1}
 =R_0(z)\sum_{l\ge0}\bigl((t-t_0)M_qR_0(z)\bigr)^l,
 \qquad z\in\Gamma.
\]
The series converges uniformly on the contour. It also converges as a map into $\operatorname{Dom}H(t_0)$ equipped with its graph norm, because $H(t_0)R_0(z)=zR_0(z)-I$ and the
contour is compact. On every smaller $t$-disk its derivative series converges in this graph norm as well: the additional factors $l$ are summable against a geometric sequence.

The contour integral
\[
 P(t)=\frac1{2\pi\mathrm i}\oint_\Gamma(z-H(t))^{-1}\,\mathrm dz
\]
is therefore analytic with values in bounded maps from $L^2$ to $\operatorname{Dom}H(t_0)$. For real $t$ it is the orthogonal spectral projection on the enclosed eigenvalues.
Shrink $\delta$ so that $\|P(t)-P(t_0)\|<1$. If $u\in\operatorname{Ran}P(t_0)$ and $P(t)u=0$, then $\|u\|\le\|P(t)-P(t_0)\|\|u\|$, so $u=0$. Interchanging the two projections gives
the opposite rank inequality. Thus $P(t)$ has rank one for these real $t$.

Let $\psi_0$ be a unit vector in $\operatorname{Ran}P(t_0)$. Since $\|P(t)\psi_0\|\ge1-\|P(t)-P(t_0)\|>0$, define
\[
 \psi_k(t)=\frac{P(t)\psi_0}{\|P(t)\psi_0\|}.
\]
This vector is differentiable in the common graph norm for real $t$. Its eigenvalue is the $k$-th ordered eigenvalue: the earlier min--max bound gives continuity of every ordered
eigenvalue, and the contour contains only one eigenvalue, with multiplicity one, near $t_0$. Thus $H(t)\psi_k(t)=\lambda_k(t)\psi_k(t)$ may be differentiated in $L^2$. Pairing the
differentiated equation with $\psi_k(t)$ in the second argument gives
\[
 \langle M_q\psi_k,\psi_k\rangle
 +\langle H(t)\psi_k',\psi_k\rangle
 =\lambda_k'(t)+\lambda_k(t)\langle\psi_k',\psi_k\rangle.
\]
Self-adjointness and the eigenvalue equation cancel the last terms, proving
\[
 \lambda_k'(t)=\int_0^\infty q(x)|\psi_k(t,x)|^2\,\mathrm dx.
\]
If $q\ge0$ and $q$ is positive almost everywhere on an open interval, this integral cannot be zero. Otherwise $\psi_k$ vanishes almost everywhere on that interval, hence everywhere
there by continuity, and the previously proved uniqueness of the differential system makes it zero on the whole half-line. This contradicts its unit norm. Since the local argument
applies at every real $t_0$, the ordered eigenvalue is differentiable everywhere with a positive derivative. The mean value theorem on every compact interval then proves that it is
strictly increasing.

\begin{theorem}\label{A:res:20.1}
For every deterministic $a\in\mathbb R$, $\mathbb P\{a\in\sigma(\mathcal H_\beta)\}=0$.
\end{theorem}

\begin{proof}
Choose a nonnegative $\phi\in C_c^\infty(0,\infty)$ with $\|\phi\|_2=1$. It is strictly positive on some open interval. Decompose the noise into its component in this direction and
the orthogonal Gaussian remainder:
\[
 \xi=\int_0^\infty\phi\,\mathrm dB,
 \qquad h_\phi(x)=\int_0^x\phi(s)\,\mathrm ds,
 \qquad B^\perp(x)=B(x)-\xi h_\phi(x).
\]
Then $\xi$ is standard normal. The integral also equals $-\int B\phi'$, so it is a Borel function of the continuous path. For every $x$, the covariance is
\[
 \mathbb E[\xi B^\perp(x)]
 =\int_0^x\phi(s)\,\mathrm ds-h_\phi(x)=0.
\]
Each finite collection of values of $B^\perp$, together with $\xi$, is jointly Gaussian. Hence $\xi$ is independent of the values of $B^\perp$ at rational times, and therefore of
the entire continuous path $B^\perp$.

Let $R$ contain the support of $\phi$. The function $h_\phi$ is Lipschitz and constant on $[R,\infty)$. Thus $\mathcal O_\alpha(h_\phi)<\infty$. At any dyadic spacing $d\le1$, only
at most $R/d+2$ intervals have nonzero oscillation, and each oscillation is at most $d\|\phi\|_\infty$. Consequently
\[
 \sum_{j\ge0}(1+jd)^{-3/2}
     \operatorname{osc}_{[jd,(j+1)d]}(h_\phi)^2
 \le(R/d+2)d^2\|\phi\|_\infty^2\le C_\phi.
\]
This proves $Q(h_\phi)<\infty$. The defining inequalities for $\Omega_B$ give, for every real $t$,
\[
 \mathcal O_\alpha(b+th_\phi) \le\mathcal O_\alpha(b)+|t|\mathcal O_\alpha(h_\phi),\qquad
 Q(b+th_\phi) \le2Q(b)+2t^2Q(h_\phi).
\]
Applying them with $t$ and $-t$ proves that $b\in\Omega_B$ if and only if $b+th_\phi\in\Omega_B$. In particular $B^\perp\in\Omega_B$ almost surely. Its law need not be Wiener
measure; membership in $\Omega_B$ is the pathwise property needed to define the form.

For each such path $b$, let $H^\perp$ be the operator defined by its integrated potential $c_\beta b$. Integration by parts, first for compactly supported functions in the form
domain, gives
\[
 -c_\beta t\int_0^\infty h_\phi(x)(|f|^2)'(x)\,\mathrm dx
 =c_\beta t\int_0^\infty\phi(x)|f(x)|^2\,\mathrm dx.
\]
For $f_m\to f$ in $\mathcal D_{\mathrm A}$, the functions $h_\phi$ and $\phi$ are bounded, and
\[
 \|(|f_m|^2)'-(|f|^2)'\|_1
 \le2\|f_m'-f'\|_2\|f_m\|_2
       +2\|f'\|_2\|f_m-f\|_2\longrightarrow0,
\]
\[
  \||f_m|^2-|f|^2\|_1
 \le(\|f_m\|_2+\|f\|_2)\|f_m-f\|_2\longrightarrow0.
\]
Approximating $f$ by compactly supported functions in the form domain therefore proves the same integral identity for $f$. The entire family on this path is therefore
\[
 H_b(t)=H^\perp+c_\beta tM_\phi,
 \qquad t\in\mathbb R.
\]
Its ordered eigenvalues satisfy
\[
 \frac{\mathrm d}{\mathrm dt}\Lambda_k(b,t)
 =c_\beta\int_0^\infty\phi(x)|\psi_k(b,t,x)|^2\,\mathrm dx>0.
\]
Thus $\{t:\Lambda_k(b,t)=a\}$ has at most one element. Strict monotonicity makes the inverse image of a fixed value contain at most one point, and the remaining scalar Gaussian
assigns zero probability to such a set.

We also verify the measurability needed for conditioning. Choose a countable set of smooth compactly supported real functions dense in $\mathcal D_{\mathrm A}$ for its
deterministic energy norm. To justify restricting min--max to these spans, fix an orthonormal basis $f_1,\ldots,f_k$ of a trial subspace and choose $f_{j,m}$ from the dense set
with $f_{j,m}\to f_j$ in the energy norm. Write $G_m=(\langle f_{i,m},f_{j,m}\rangle)_{i,j}$ and $Q_m=(\mathfrak q_{b,t}[f_{i,m},f_{j,m}])_{i,j}$. Continuity of the form gives
$G_m\to I$ and $Q_m\to Q$ in matrix norm. For large $m$, $G_m\ge I/2$, so the approximating vectors are independent. Moreover
\[
 \sup_{|c|=1}\left|
   \frac{c^*Q_mc}{c^*G_mc}-c^*Qc\right|
 \le2\|Q_m-Q\|+2\|Q\|\|G_m-I\|\longrightarrow0.
\]
Thus the largest Rayleigh quotient of each approximating span converges to that of the original span, which proves the stated restriction of the infimum. On each such span the
supremum of the Rayleigh quotient can be taken over nonzero rational coefficient vectors. Each quotient is Borel in $(b,t)$, since its noise integral has compact support and is
continuous in the local uniform topology on paths. If $\mathcal D_0$ is the chosen countable dense set, the formula is
\[
 \Lambda_k(b,t)=
 \inf_{\substack{f_1,\ldots,f_k\in\mathcal D_0\\
                       f_1,\ldots,f_k\ \mathrm{independent}}}
 \ \sup_{q\in\mathbb Q^k\setminus\{0\}}
 \frac{\mathfrak q_{b,t}[\sum_{j=1}^kq_jf_j]}
      {\|\sum_{j=1}^kq_jf_j\|_2^2}.
\]
Here $\mathfrak q_{b,t}$ is the form of $H_b(t)$. Each denominator is positive. Thus a countable supremum followed by a countable infimum of Borel functions proves measurability of
$\Lambda_k$. If $\nu$ denotes the law of $B^\perp$, independence and Tonelli's theorem now give
\[
 \mathbb P\{\Lambda_k(B^\perp,\xi)=a\}
 =\int\nu(\mathrm db)\int_{\mathbb R}
       \mathbf1_{\{\Lambda_k(b,t)=a\}}
       \frac{\mathrm e^{-t^2/2}}{\sqrt{2\pi}}\,\mathrm dt=0.
\]
The exceptional paths outside $\Omega_B$ have $\nu$-measure zero. Finally the spectrum is the countable set of eigenvalues, so a countable union over $k$ proves the theorem.
\end{proof}

Taking $a=0$ gives $0\in\rho(\mathcal H_\beta)$ almost surely. A countable union gives the same conclusion for any deterministic countable set of real numbers. This does not
exclude the random eigenvalues themselves. Below, after fixing a path, a real parameter may be chosen anywhere in that path's resolvent set.

To change the spectral parameter, we need the finite inverse to be embedded isometrically. Use the pairs
\[
 (A_n,J_n,m_n)=
 \begin{cases}
 (\mathcal A_{n,\beta}^{\mathrm{DE}},\ \iota_h^{(n)},\ n),
                  &\text{Gaussian ensemble},\\
 (\widetilde H_n^{\mathrm{red}},\ \iota_h^{(M_n)},\ M_n),
                  &\text{general Jacobi ensemble}.
 \end{cases}
\]
Here $\widetilde H_n^{\mathrm{red}}=U_{r_n}H_n^{\mathrm{red}}U_{r_n}^*$, $J_n^*J_n=I_{m_n}$, and $h=h_n$ is the corresponding spacing. For the reduced matrix and its isometric
embedding, Theorem~\ref{A:res:1.2} and \eqref{A:eq:12.43}, conjugated by $\mathcal J_n$, give the same nonreal estimate as for the full compressed resolvent. Indeed,
$\mathcal J_nP_{I_n}^*=\iota_h^{(M_n)}U_{r_n}$ and $\|\mathcal J_n\|\le1$, so the additional squared error is $o(h)$.

Fix a path on the event used in Subsection~\ref{A:sec:proof-of-the-resolvent-theorem}. For $a\in\rho(\mathcal H_\beta)\cap\mathbb R$, set
\[
 S(a)=(\mathcal H_\beta-a-\mathrm i)^{-1},\qquad
 S_n(a)=J_n(A_n-a-\mathrm i)^{-1}J_n^*.
\]
The established nonreal comparison gives
\[
 \|S_n(a)-S(a)\|_{\mathfrak S_2}^2\le C_{\omega,a}h_n.
 \eqtag{20.25}
\]
For a nonempty compact real set $L\subset\rho(\mathcal H_\beta)$, the bound is uniform over $a\in L$: the translate $L+\mathrm i$ lies in some $K_m$, even when $L$ is chosen after
the path is fixed. There is nothing to check when the set is empty. For a scalar eigenvalue,
\[
 \frac{(\lambda-a-\mathrm i)^{-1}}
      {1+\mathrm i(\lambda-a-\mathrm i)^{-1}}=(\lambda-a)^{-1}.
\]
The denominator is uniformly invertible on the limiting spectrum. We use the same rational function for the embedded operators. Its value at zero must be zero, since the finite
embedding has an orthogonal complement on which its resolvent image vanishes.

\begin{proposition}\label{A:res:20.2}
For either pair $(A_n,J_n)$, every compact real $L\subset\rho(\mathcal H_\beta)$ satisfies, for all sufficiently large
$n$,
\[
 L\subset\rho(A_n),\qquad
 \sup_{a\in L}\|J_n(A_n-a)^{-1}J_n^*
                 -(\mathcal H_\beta-a)^{-1}\|_{\mathfrak S_2}^2
 \le C_{\omega,L}h_n.
\]
\end{proposition}

\begin{proof}
Put $d=\operatorname{dist}(L,\sigma(\mathcal H_\beta))>0$ and $b=(1+d^{-2})^{1/2}$. For $Q(a)=I+\mathrm iS(a)$, spectral calculus gives
\[
 \|Q(a)^{-1}\|_{\mathrm{op}}
 =\sup_{\lambda\in\sigma(\mathcal H_\beta)}
     \frac{|\lambda-a-\mathrm i|}{|\lambda-a|}\le b.
\]
The spectral value $0$ of the compact operator $S(a)$ contributes $1$ to this bound. Set $Q_n(a)=I+\mathrm iS_n(a)$. Uniformly for $a\in L$,
\[
 Q_n(a)=Q(a)\bigl[I+\mathrm iQ(a)^{-1}(S_n(a)-S(a))\bigr].
\]
For $n$ so large that $b\sup_{a\in L}\|S_n(a)-S(a)\|\le1/2$, the second factor is invertible by a Neumann series and
\[
 \sup_{a\in L}\|Q_n(a)^{-1}\|_{\mathrm{op}}\le2b.
 \eqtag{20.31}
\]
On $\operatorname{Ran}J_n$, the operator $Q_n(a)$ is unitarily equivalent to $(A_n-a)(A_n-a-\mathrm i)^{-1}$; on its orthogonal complement it is the identity. Thus invertibility of
$Q_n(a)$ proves $a\in\rho(A_n)$.

Let $f(w)=w/(1+\mathrm iw)$. Since $f(0)=0$,
\[
 f(S_n(a))=J_n(A_n-a)^{-1}J_n^*,\qquad
 f(S(a))=(\mathcal H_\beta-a)^{-1}.
\]
The identity $f(S)=\mathrm i^{-1}(I-Q^{-1})$ and the inverse difference formula give
\[
 f(S_n)-f(S)=Q_n^{-1}(S_n-S)Q^{-1}.
\]
There is also a uniform first-order expansion with the limiting outer factors. Put $\Delta_n(a)=S_n(a)-S(a)$. Applying the inverse difference identity once more, in the same order,
gives
\[
 f(S_n)-f(S)-Q^{-1}\Delta_nQ^{-1}
 =-\mathrm iQ_n^{-1}\Delta_nQ^{-1}\Delta_nQ^{-1}.
\]
Consequently its Hilbert--Schmidt norm is at most $2b^3\|\Delta_n\|_{\mathrm{op}}\|\Delta_n\|_{\mathfrak S_2}$. The nonreal estimate bounds this by $C_{\omega,L}h_n$, uniformly on
$L$. The real resolvent therefore has the same error order, with a smaller nonlinear remainder. Its leading linear error is $Q(a)^{-1}\Delta_n(a)Q(a)^{-1}$, obtained by multiplying
$\Delta_n(a)$ on both sides by the limiting factors.

The product formula specifies the order of the factors for possibly noncommuting $S_n$ and $S$. The ideal property gives
\[
 \sup_{a\in L}\|J_n(A_n-a)^{-1}J_n^* -(\mathcal H_\beta-a)^{-1}\|_{\mathfrak S_2}^2
 \le4b^4\sup_{a\in L}\|S_n(a)-S(a)\|_{\mathfrak S_2}^2 \le C_{\omega,L}h_n.
 \eqtag{20.35}
\]
For the operator norm bound, the spectral theorem gives $\sup_{a\in L}\|(\mathcal H_\beta-a)^{-1}\|\le d^{-1}$. Hence
\[
 \sup_{a\in L}\|J_n(A_n-a)^{-1}J_n^*\|_{\mathrm{op}}
 \le d^{-1}+C_{\omega,L}^{1/2}h_n^{1/2}.
\]
Because $J_n$ is an isometry for the two pairs under consideration, $\|J_n(A_n-a)^{-1}J_n^*\|_{\mathrm{op}} =\|(A_n-a)^{-1}\|_{\mathrm{op}}$. This proves the uniform bound for all
sufficiently large $n$.
\end{proof}

For real parameters, the proof of the lower bound uses the same singular value identity. The inverse may have negative eigenvalues, so the absolute values of its reciprocal
eigenvalues must be put in decreasing order. The following calculation gives the estimate uniformly on the chosen compact set. Write $A=1+\sup_{a\in L}|a|$. By
Lemma~\ref{A:res:weyl}, there are $c,C>0$ and $k_0\ge2$ such that $ck^{2/3}\le\Lambda_k\le Ck^{2/3}$ and $\Lambda_k>2A$ for $k\ge k_0$. Hence, uniformly for $a\in L$,
\[
 \sum_{k>m}|\Lambda_k-a|^{-2}
 \asymp_{\omega,L}\sum_{k>m}k^{-4/3}\asymp m^{-1/3}.
 \eqtag{20.40}
\]
The last comparison follows by integrating $x^{-4/3}$ over $[m,\infty)$ and $[m+1,\infty)$. To justify the rearrangement, put $C_0=\max_{j<k_0}(|\Lambda_j|+A)$, increasing it to a
positive constant if necessary. For all sufficiently large $k$, $\Lambda_k-A>C_0$. Thus every reciprocal from the first $k_0-1$ eigenvalues precedes this tail in decreasing
absolute value, uniformly in $a\in L$. Within the tail the ordering is already decreasing, because $\Lambda_k>a$. Only a fixed finite initial segment can therefore be permuted.

The finite-rank inequality \eqref{A:eq:5.4} now gives
\[
 \inf_{a\in L}\|J_n(A_n-a)^{-1}J_n^* -(\mathcal H_\beta-a)^{-1}\|_{\mathfrak S_2}^2 \ge\inf_{a\in
 L}\sum_{j>m_n}s_j((\mathcal H_\beta-a)^{-1})^2 \ge c_{\omega,L}m_n^{-1/3}\ge c'_{\omega,L}h_n.
 \eqtag{20.43}
\]
The proof in Section~\ref{A:sec:rank-lower-bound-and-optimality} applies to any finite-rank operator; positivity of the inverse is not needed.

\begin{theorem}\label{A:res:20.3}
For either isometric pair above and every compact real $L\subset\rho(\mathcal H_\beta)$, almost surely, for all sufficiently large $n$ and all $a\in L$,
\[
 c_{\omega,L}h_n\le
 \|J_n(A_n-a)^{-1}J_n^*-(\mathcal H_\beta-a)^{-1}\|_{\mathfrak S_2}^2
 \le C_{\omega,L}h_n.
\]
\end{theorem}

\begin{proof}
Proposition~\ref{A:res:20.2} gives $L\subset\rho(A_n)$ for all sufficiently large $n$. For such $n$ and $a\in L$,
\[
 c_{\omega,L}h_n
 \le\sum_{j>m_n}s_j((\mathcal H_\beta-a)^{-1})^2
 \le\|J_n(A_n-a)^{-1}J_n^*-(\mathcal H_\beta-a)^{-1}\|_{\mathfrak S_2}^2
 \le C_{\omega,L}h_n.
\]
The first two inequalities are \eqref{A:eq:20.43}; the last is Proposition~\ref{A:res:20.2}. Their constants and starting size are uniform on $L$.
\end{proof}

The preceding result concerns the retained Dirichlet matrix. To prove the introduction's inverse theorem, we must still invert the full general Jacobi matrix and then compress.
These operations cannot be interchanged. For a coordinate restriction $P$ and a matrix $R$,
\[
 (PRP^*)^2=PR(P^*P)RP^*,
\]
which is generally different from $PR^2P^*$. We consequently make the following comparisons at real parameters on the full finite space.

Use the notation of Section~\ref{A:sec:boundary-reduction}: $A_C=H_n^{\mathrm{red}}$ on $\{r_n+1,\ldots,N_n\}$, $A_L$ on $\{1,\ldots,r_n\}$, and $A_P$ on $\{1,\ldots,N_n\}$, where
$N_n=n-s_\beta$. Write $r=r_n$, $\ell_n=1+\log(n+1)$, $\alpha=(h^{-2}/b_\star)b_r$, and let $e_L,e_C$ be the boundary point vectors adjacent to this cut. Denote the right boundary
of $A_C$ by $e_P$. For a compact real $L\subset\rho(\mathcal H_\beta)$, set
\[
 \begin{gathered}
 R_C(a)=(A_C-a)^{-1},\quad S_C(a)=(A_C-a-\mathrm i)^{-1},
 \quad B_L(a)=(A_L-a)^{-1},\\
 \kappa_L(a)=\alpha-\alpha^2e_L^*B_L(a)e_L,
 \qquad d_C(a)=\alpha-\alpha^2e_C^*R_C(a)e_C.
 \end{gathered}
\]
The inverse $B_L(a)$ exists uniformly on bounded real sets for large $n$, by Lemma~\ref{A:res:left-resolvent}. All pathwise statements that follow are made on the event already
fixed, intersected with the boundary coefficient events for the countably many bounded spectral sets.

\par\Needspace{6\baselineskip}
\begin{proposition}\label{A:res:real-boundary}
Uniformly for $a\in L$, for all sufficiently large $n$,
\[
 \|R_C(a)\|_{\mathrm{op}}\le C_{\omega,L},\qquad
 \|R_C(a)e_C\|^2\le C_{\omega,L}h^3\ell_n^6,\qquad
 \|R_C(a)e_P\|^2\le C_{\omega,L}h^4,
 \eqtag{20.46}
\]
and
\[
 |d_C(a)|\le C_{\omega,L}h^{-1}\ell_n^6.
 \eqtag{20.47}
\]
The left resolvent satisfies \eqref{A:eq:20.48}--\eqref{A:eq:20.49} below.
\end{proposition}

\begin{proof}
Proposition~\ref{A:res:20.2} gives the bound on $R_C(a)$, since reindexing and its isometric embedding preserve the operator norm. The identities
\[
 R_C=S_C(I+\mathrm iS_C)^{-1},\qquad
 (I+\mathrm iS_C)^{-1}=I-\mathrm iR_C
\]
and commutation of these functions of $A_C$ imply $\|R_Ce\|\le(1+\|R_C\|)\|S_Ce\|$ for every vector $e$. Lemma~\ref{A:res:core-poisson}, on the bounded nonreal set $L+\mathrm i$,
gives $\|S_C(a)e_C\|^2\le Ch^3\ell_n^6$. At the other boundary point $x_{e_P}\asymp h^{-2}$, so Corollary~\ref{A:res:3.18} gives $\|S_C(a)e_P\|^2\le Ch^4$. The common form
hypotheses for this application are supplied by Theorem~\ref{A:res:16.2}. These are the remaining bounds in \eqref{A:eq:20.46}.

The nonreal estimate is $|d_C(a+\mathrm i)|\le Ch^{-1}\ell_n^4$. Since $R_C-S_C=-\mathrm iR_CS_C$ and $R_C(a)$ is self-adjoint,
\[
 |d_C(a)-d_C(a+\mathrm i)| \le\alpha^2\|R_C(a)e_C\|\,\|S_C(a)e_C\| \le
 C_{\omega,L}h^{-4}h^3\ell_n^6 =C_{\omega,L}h^{-1}\ell_n^6.
\]
This proves \eqref{A:eq:20.47}.

For the left block, write $Q_L=A_L-\alpha e_Le_L^*$ and $G_a=(Q_L-a)^{-1}$. Lemma~\ref{A:res:left-resolvent} gives $Q_L\ge c h^{-2}r^{-2}I$, while $\alpha\asymp h^{-2}$. Since
$h^{-2}r^{-2}\to\infty$, choose $n_0=n_0(\omega,L)$ so that $c h^{-2}r^{-2}\ge2\sup_{a\in L}|a|$ for $n\ge n_0$. Then, uniformly on $L$,
\[
 Q_L-a\ge\frac c2h^{-2}r^{-2}I,
 \qquad 0<G_a\le C h^2r^2I.
\]
Put $q_a=e_L^*G_ae_L>0$. Multiplication by $Q_L-a+\alpha e_Le_L^*$ verifies the exact formula
\[
 B_L(a)=G_a-\frac{\alpha}{1+\alpha q_a}G_ae_Le_L^*G_a,
 \qquad B_L(a)e_L=\frac{G_ae_L}{1+\alpha q_a}.
\]
Both matrices being inverted are positive, and the displayed correction is positive semidefinite. Hence $0<B_L(a)\le G_a$, so $\|B_L(a)\|\le Ch^2r^2$ and
$\|B_L(a)\|_{\mathfrak S_2}^2\le r\|B_L(a)\|^2\le Ch^4r^5$.
For the boundary column, the spectral theorem gives $G_a^2\le\|G_a\|G_a$. Consequently
\[
 \|B_L(a)e_L\|^2
 \le\|G_a\|\frac{q_a}{(1+\alpha q_a)^2}
 \le\frac{\|G_a\|}{4\alpha}\le Ch^4r^2.
\]
Here $x/(1+x)^2\le1/4$ for $x\ge0$; this additional denominator is essential for the column estimate. We have obtained
\[
 \|B_L(a)\|_{\mathrm{op}}\le Ch^2r^2,\qquad
 \|B_L(a)\|_{\mathfrak S_2}^2\le Ch^4r^5,\qquad
 \|B_L(a)e_L\|^2\le Ch^4r^2.
 \eqtag{20.48}
\]
The scalar boundary value is $e_L^*B_L(a)e_L=q_a/(1+\alpha q_a)$, and therefore $\kappa_L(a)=\alpha/(1+\alpha q_a)$. The bound $q_a\le Ch^2r^2$ and $r\ge1$ give
$1\le1+\alpha q_a\le1+Cr^2\le C'r^2$. Thus
\[
 \frac{c\alpha}{r^2}\le\kappa_L(a)\le\alpha.
 \eqtag{20.49}
\]
All constants and the starting index are uniform for $a\in L$. Positivity was used only for the left block; the reduced inverse $R_C(a)$ need not be positive.
\end{proof}

Let $P_L,P_C$ be the coordinate restrictions of $\mathcal H_P=\mathcal H_L\oplus\mathcal H_C$ onto its two summands. In this decomposition the off-diagonal block is
$-\alpha e_Le_C^*$. Set
\[
 m_L=e_L^*B_L(a)e_L,\quad m_C=e_C^*R_C(a)e_C,
 \quad D_n(a)=1-\alpha^2m_Lm_C.
\]

\begin{proposition}\label{A:res:real-left-removal}
For all sufficiently large $n$, $L\subset\rho(A_P)$ and
\[
 \sup_{a\in L}\|(A_P-a)^{-1}-P_C^*R_C(a)P_C\|_{\mathfrak S_2}^2
 =o_\omega(h).
 \eqtag{20.50}
\]
Also $\sup_{a\in L}\|(A_P-a)^{-1}\|_{\mathrm{op}}\le C_{\omega,L}$.
\end{proposition}

\begin{proof}
Put $u=\kappa_L(a)/\alpha$ and $v=d_C(a)/\alpha$. Then
\[
 D_n(a)=u+(1-u)v,\qquad
 cr^{-2}\le u\le1,\qquad |v|\le C_{\omega,L}h\ell_n^6.
\]
In particular, uniformly on $L$,
\[
 D_n(a)\ge cr^{-2}-C_{\omega,L}h\ell_n^6\ge\frac{c}{2r^2}>0
\]
for large $n$, because $hr^2\ell_n^6\to0$. This proves that the Schur complement is invertible. The estimate does not imply $R_C(a)\ge0$; its self-adjointness and the bound on
$d_C(a)$ are enough.

Both diagonal blocks are invertible, and the denominator just estimated is nonzero. We may therefore use the block inverse identity obtained from the Schur complement
\eqref{A:eq:12.35} at the real parameter $a$. Put $x=B_L(a)e_L$ and $y=R_C(a)e_C$. The identity reads
\[
 R_P(a):=(A_P-a)^{-1}
 =\begin{pmatrix}
 B_L+\dfrac{\alpha^2m_C}{D_n}xx^*
      &\dfrac{\alpha}{D_n}xy^*\\[2mm]
 \dfrac{\alpha}{D_n}yx^*
      &R_C+\dfrac{\alpha^2m_L}{D_n}yy^*
 \end{pmatrix}.
\]
The adjoints in this formula are valid because the parameter is real. Also $0\le m_L\le\alpha^{-1}$ and $|m_C|=|1-v|/\alpha\le C\alpha^{-1}$. Each coefficient of a rank-one term is
therefore bounded by $C\alpha r^2$. Using $\|xy^*\|_{\mathfrak S_2}^2=\|x\|^2\|y\|^2$ gives
\[
 \begin{aligned}
 \|P_CR_PP_C^*-R_C\|_{\mathfrak S_2}^2
    &\le Ch^{-4}r^4(h^3\ell_n^6)^2
       =Ch^2r^4\ell_n^{12},\\
 \|P_LR_PP_C^*\|_{\mathfrak S_2}^2
    &\le Ch^{-4}r^4(h^4r^2)(h^3\ell_n^6)
       =Ch^3r^6\ell_n^6,\\
 \|P_LR_PP_L^*-B_L\|_{\mathfrak S_2}^2
    &\le Ch^{-4}r^4(h^4r^2)^2=Ch^4r^8.
 \end{aligned}
\]
Blocks with distinct domain or range summands are orthogonal in Hilbert--Schmidt inner product. Thus
\[
 \begin{split}
 \|R_P-P_C^*R_CP_C\|_{\mathfrak S_2}^2
 &=\|P_CR_PP_C^*-R_C\|_{\mathfrak S_2}^2
    +2\|P_LR_PP_C^*\|_{\mathfrak S_2}^2
    +\|P_LR_PP_L^*\|_{\mathfrak S_2}^2\\
 &\le C_{\omega,L}\bigl(h^2r^4\ell_n^{12}
      +h^3r^6\ell_n^6+h^4r^8+h^4r^5\bigr)=o(h).
 \end{split}
\]
This proves \eqref{A:eq:20.50}. The same block formula in operator norm gives the more precise estimate
\[
 \|R_P-P_C^*R_CP_C\|_{\mathrm{op}}
 \le C_{\omega,L}\bigl(hr^2\ell_n^6
       +h^{3/2}r^3\ell_n^3+h^2r^4+h^2r^2\bigr)=o(1).
\]
Together with \eqref{A:eq:20.46}, it proves the operator bound.
\end{proof}

With the left block restored, the real inverse is uniformly bounded. Restoring the block on the last coordinates changes the retained equation by a term tending to zero in operator
norm, so the inverse remains bounded.

Let $A_R$ be the principal matrix on the last $s_\beta$ coordinates, $e_R$ its left boundary, and $\alpha_R=(h^{-2}/b_\star)b_{N_n}$. We also write $e_P$ for the right boundary of
$A_P$, namely the inclusion of the vector already denoted $e_P$ in $\mathcal H_C$. Let $P_P,P_R$ be the restrictions from $\mathbb C^n$ to the initial and terminal summands. When
used on $\mathbb C^n$, $P_C$ means the restriction to the reduced interior, that is, its preceding restriction composed with $P_P$.

\par\Needspace{6\baselineskip}
\begin{proposition}\label{A:res:real-terminal-removal}
For all sufficiently large $n$, $L\subset\rho(\mathcal A_{n,V,\beta})$ and
\[
 \sup_{a\in L}
 \|(\mathcal A_{n,V,\beta}-a)^{-1}-P_P^*R_P(a)P_P\|_{\mathfrak S_2}^2
 =o_\omega(h).
\]
The same estimate holds with $P_P^*R_P(a)P_P$ replaced by $P_C^*R_C(a)P_C$.
\end{proposition}

\begin{proof}
By Lemma~\ref{A:res:12.1}, $A_R+\kappa_\omega\ge c_\omega h^{-2}I$. Since $L$ is bounded, $A_R-a\ge(c_\omega/2)h^{-2}I$ for all large $n$, uniformly on $L$. For
$B_R(a)=(A_R-a)^{-1}$, the fixed dimension of this block gives
\[
 \|B_R(a)\|_{\mathrm{op}}\le Ch^2,
 \qquad \|B_R(a)\|_{\mathfrak S_2}^2\le Ch^4,
 \qquad \|B_R(a)e_R\|^2\le Ch^4.
\]
Choose a fixed exponent $p_0$ large enough for the coefficient bound in \eqref{A:eq:12.14}, so that $|\alpha_R|\le C_\omega h^{-1/2}\ell_n^{p_0}$ for large $n$. If
$m_R=e_R^*B_Re_R$ and $\sigma_R=\alpha_R^2m_R$, then
\[
 0\le\sigma_R\le C_{\omega,L}h\ell_n^{2p_0}=o(1).
\]
By Proposition~\ref{A:res:real-left-removal}, $|\sigma_Re_P^*R_Pe_P|\le1/2$ for large $n$. Consequently
\[
 S_P=(A_P-a-\sigma_Re_Pe_P^*)^{-1}
 =R_P+\frac{\sigma_R}{1-\sigma_Re_P^*R_Pe_P}
                   R_Pe_Pe_P^*R_P
\]
exists, is self-adjoint, and satisfies $\|S_P\|_{\mathrm{op}}\le C$. Eliminating the terminal coordinates as in the previous proof now gives the full inverse
\[
 (\mathcal A_{n,V,\beta}-a)^{-1}
 =\begin{pmatrix}
 S_P&\alpha_RS_Pe_Pe_R^*B_R\\
 \alpha_RB_Re_Re_P^*S_P
   &B_R+\alpha_R^2(e_P^*S_Pe_P)B_Re_Re_R^*B_R
 \end{pmatrix}.
\]
The four block estimates are
\[
 \begin{aligned}
\|S_P-R_P\|_{\mathfrak S_2}^2& \le C\sigma_R^2\|R_Pe_P\|^4 \le Ch^2\ell_n^{4p_0},\\
2\|\alpha_RB_Re_Re_P^*S_P\|_{\mathfrak S_2}^2& =2\alpha_R^2\|B_Re_R\|^2\|S_Pe_P\|^2 \le Ch^3\ell_n^{2p_0},\\
\|B_R+\alpha_R^2(e_P^*S_Pe_P)B_Re_Re_R^*B_R\|_{\mathfrak S_2}^2& \le 2\|B_R\|_{\mathfrak S_2}^2 +2\alpha_R^4|e_P^*S_Pe_P|^2\|B_Re_R\|^4 \le C\bigl(h^4+h^6\ell_n^{4p_0}\bigr).
\end{aligned}
\]
Their sum is $o(h)$, since $h\asymp n^{-1/3}$ and $p_0$ is fixed. Orthogonality of the blocks proves the first error bound. For the second, add the inverse difference from
Proposition~\ref{A:res:real-left-removal}, extended by zero on the terminal coordinates, and use $\|D+E\|_{\mathfrak S_2}^2\le2\|D\|_{\mathfrak S_2}^2+2\|E\|_{\mathfrak S_2}^2$.
\end{proof}

For large $n$, $\|\mathcal J_n\|=1$. Conjugating the last comparison by $\mathcal J_n$ gives
\[
 \sup_{a\in L}\bigl\| \mathcal J_n(\mathcal A_{n,V,\beta}-a)^{-1}\mathcal J_n^*
 -\iota_h^{(M_n)}U_{r_n}R_C(a)U_{r_n}^*(\iota_h^{(M_n)})^* \bigr\|_{\mathfrak S_2}^2=o_\omega(h).
 \eqtag{20.52}
\]
Theorem~\ref{A:res:20.3} and the squared triangle inequality prove the $O(h)$ upper bound for the full compressed inverse. Its lower bound follows directly from \eqref{A:eq:5.4}:
its rank is at most $M_n$, so its squared error is at least the tail in \eqref{A:eq:20.40}.

We now return to the pairs $(\mathscr A_n,\mathscr I_n)$ of the introduction, including the full general Jacobi matrix with $\mathscr I_n=\mathcal J_n$.

\begin{corollary}\label{A:res:20.4}
For every deterministic $a\in\mathbb R$, almost surely, $a\in\rho(\mathscr A_n)$ for all sufficiently large $n$ and
\[
 \|\mathscr I_n(\mathscr A_n-a)^{-1}\mathscr I_n^*
              -(\mathcal H_\beta-a)^{-1}\|_{\mathfrak S_2}^2
 =\Theta_{\omega,a}(n^{-1/3}).
\]
The same conclusion holds for the isometrically embedded reduced Jacobi operator.
\end{corollary}

\begin{proof}
Theorem~\ref{A:res:20.1} puts $a$ outside the limiting spectrum almost surely. Theorem~\ref{A:res:20.3}, with $L=\{a\}$, proves the claim for the two isometric pairs. For the
inverse of the full Jacobi matrix, let $R_n^C(a)$ be the embedded reduced inverse in \eqref{A:eq:20.52}. That equation gives
\[
 \|\mathcal J_n(\mathcal A_{n,V,\beta}-a)^{-1}\mathcal J_n^*
                   -(\mathcal H_\beta-a)^{-1}\|_{\mathfrak S_2}^2
 \le2\|R_n^C(a)-(\mathcal H_\beta-a)^{-1}\|_{\mathfrak S_2}^2
       +o_\omega(h_n)
 \le C_{\omega,a}h_n.
\]
The rank is at most $M_n$, so \eqref{A:eq:5.4} and \eqref{A:eq:20.40} give the lower bound $c_{\omega,a}h_n$. Proposition~\ref{A:res:real-terminal-removal} supplies invertibility
of the full matrix. Finally, $h_n\asymp n^{-1/3}$. In particular the removal of coordinates does not change a leading squared error, should one be identified for the reduced
inverse. Indeed, write $E_n=R_n^C(a)-(\mathcal H_\beta-a)^{-1}$ and let $F_n$ be the difference between the compressed full inverse and $R_n^C(a)$. The established bounds are
$\|E_n\|_{\mathfrak S_2}=O_\omega(\sqrt{h_n})$ and $\|F_n\|_{\mathfrak S_2}=o_\omega(\sqrt{h_n})$. Hence
\[
 \frac{\bigl|\|E_n+F_n\|_{\mathfrak S_2}^2
                     -\|E_n\|_{\mathfrak S_2}^2\bigr|}{h_n}
 \le 2\frac{\|E_n\|_{\mathfrak S_2}}{\sqrt{h_n}}
        \frac{\|F_n\|_{\mathfrak S_2}}{\sqrt{h_n}}
       +\frac{\|F_n\|_{\mathfrak S_2}^2}{h_n}\longrightarrow0.
\]
This does not identify that leading error with the leading coefficient for the best approximation under a rank constraint in \eqref{A:eq:rank-tail-asymptotic}.

\end{proof}

\begin{proof}[Proof of Theorem~\ref{A:res:1.1}]
Take $a=0$ in Corollary~\ref{A:res:20.4}. On the event $0\in\rho(\mathcal H_\beta)$, Lemma~\ref{A:res:weyl} gives
\[
 \|\mathcal H_\beta^{-1}\|_{\mathfrak S_2}^2
 =\sum_{k\ge1}|\Lambda_k|^{-2}<\infty.
\]
The corollary gives a finite $n_0(\omega)$ such that the full finite matrix is invertible and both error bounds hold for every $n\ge n_0(\omega)$. Since $h_n\asymp n^{-1/3}$, those
bounds have order $h_n$.
\end{proof}

The inverse comparison gives a matching of reciprocal eigenvalues. A reciprocal near zero can come from a very large eigenvalue of either sign. We establish the matching statement
here and recover the increasing labels in Subsection~\ref{A:sec:spectral-consequences} after proving a common lower bound.

For either isometric pair $(A_n,J_n)$, write $\lambda_{1,n}\le\cdots\le\lambda_{m_n,n}$ for the finite eigenvalues and put $r_{j,n}(a)=(\lambda_{j,n}-a)^{-1}$,
$r_k(a)=(\Lambda_k-a)^{-1}$. Define
\[
 \mathfrak d_{2,n}(a)^2 :=\inf_{\substack{J\subset\{1,\ldots,m_n\}\\ \pi:J\hookrightarrow\mathbb
 N}} \biggl\{ \sum_{j\in J}|r_{j,n}(a)-r_{\pi(j)}(a)|^2 +\sum_{j\notin J}|r_{j,n}(a)|^2
 +\sum_{k\notin\pi(J)}|r_k(a)|^2\biggr\}.
 \eqtag{20.59}
\]
Unmatched entries are paired with zero. Zeros can be adjoined to either list without changing the infimum. They do not represent additional finite eigenvalues of the original
operators. The admissible partial injections form a countable set. Since the eigenvalues are measurable, each cost is measurable as a sum of nonnegative measurable functions, and
the infimum defines a measurable distance.

We prove the matching inequality directly, including the signs and the zero entries. First let $U,V$ be positive Hilbert--Schmidt operators, with decreasing eigenvalue lists
$u_i,v_i\ge0$, padded by zeros when finite. Choose orthonormal eigenvectors $e_i$ for the positive eigenvalues of $V$ and put $d_i=\langle Ue_i,e_i\rangle\ge0$. The calculation
using the orthogonal projection onto a finite dimensional range in \eqref{A:eq:5.4} gives
\[
 \sum_{i=1}^k d_i\le\sum_{i=1}^k u_i.
\]
Indeed, in an eigenbasis of $U$ the coefficients of this partial trace lie in $[0,1]$ and their sum is at most $k$. For a finite truncation $v_1,\ldots,v_M$, set $v_{M+1}=0$ in the
next calculation. Summation by parts gives
\[
 \begin{aligned}
 \sum_{i=1}^M v_id_i
 &=\sum_{k=1}^M(v_k-v_{k+1})\sum_{i=1}^k d_i\\
 &\le\sum_{k=1}^M(v_k-v_{k+1})\sum_{i=1}^k u_i
 =\sum_{i=1}^M u_iv_i.
 \end{aligned}
\]
The spectral truncations $V_M$ converge to $V$ in Hilbert--Schmidt norm, and
$|\operatorname{Tr}(U(V-V_M))|\le\|U\|_{\mathfrak S_2}\|V-V_M\|_{\mathfrak S_2}\to0$.
It follows that $\operatorname{Tr}(UV)\le\sum_i u_iv_i$.

For self-adjoint Hilbert--Schmidt operators $S,T$, write $S=S_+-S_-$ and $T=T_+-T_-$. All four parts are positive and Hilbert--Schmidt. Traces of products of two positive parts are
nonnegative, by the same eigenbasis expansion. If $s_i^\pm,t_i^\pm$ denote their decreasing nonnegative eigenvalues, then
\[
 \operatorname{Tr}(ST)
 \le\operatorname{Tr}(S_+T_+)+\operatorname{Tr}(S_-T_-)
 \le\sum_i s_i^+t_i^++\sum_i s_i^-t_i^-.
\]
Expanding the squared Hilbert--Schmidt distance therefore gives
\[
 \sum_i(s_i^+-t_i^+)^2+\sum_i(s_i^--t_i^-)^2
 \le\|S-T\|_{\mathfrak S_2}^2.
\]
This pairs positive eigenvalues with positive eigenvalues, negative eigenvalues with negative eigenvalues, and any remaining entry with zero. For a negative pair the difference of
the signed eigenvalues has the same square as the difference of their absolute values.

Apply this inequality to $S=J_n(A_n-a)^{-1}J_n^*$ and $T=(\mathcal H_\beta-a)^{-1}$. There are only $m_n$ nonzero entries on the finite side. Discarding every pair with a zero
entry on that side produces an admissible partial injection in \eqref{A:eq:20.59}; all remaining unmatched contributions are precisely the squared entries paired with zero.
Consequently
\[
 \mathfrak d_{2,n}(a)^2
 \le\|J_n(A_n-a)^{-1}J_n^*-(\mathcal H_\beta-a)^{-1}\|_{\mathfrak S_2}^2
 \le C_{\omega,a}h_n.
\]
This proves the upper bound without imposing an increasing ordering on the direct, possibly signed, spectra. Conversely, $|\pi(J)|\le m_n$ in every matching. Removing at most $m_n$
terms from a square-summable sequence leaves at least the sum of squares after its $m_n$ largest absolute values. Hence
\[
 \sum_{k\notin\pi(J)}|r_k(a)|^2
 \ge\sum_{j>m_n}s_j((\mathcal H_\beta-a)^{-1})^2
 \ge c_{\omega,a}m_n^{-1/3}.
\]
Thus $\mathfrak d_{2,n}(a)^2=\Theta_{\omega,a}(n^{-1/3})$.

For the matching statement for the full matrix, we retain every reciprocal eigenvalue by placing the omitted coordinate subspace in a separate orthogonal summand. Set
\[
 \widetilde P_n=U_{r_n}P_{I_n},\qquad
 Q_n=I-\widetilde P_n^*\widetilde P_n,\qquad
 \widehat{\mathcal H}=L^2(\mathbb R_+)\oplus\ell^2(\mathbb N).
\]
Map the standard basis vectors of $\operatorname{Ran}Q_n$, in increasing order of their original indices, to the first $n-M_n$ standard basis vectors of $\ell^2(\mathbb N)$. This
defines a deterministic isometry $V_n$ on the removed subspace. Put
\[
 \widehat{\mathcal J}_nu
   =\bigl(\iota_h^{(M_n)}\widetilde P_nu,\ V_nQ_nu\bigr).
\]
The two coordinate projections are complementary and orthogonal, so
\[
 \|\widehat{\mathcal J}_nu\|^2
 =\|\widetilde P_nu\|^2+\|Q_nu\|^2=\|u\|^2.
\]
In particular this is an isometry of the entire finite space. The limiting compact operator is $(\mathcal H_\beta-a)^{-1}\oplus0$. The zero summand represents the vanishing inverse
on the omitted directions; it is not the resolvent of a zero operator on that summand. In particular, no densely defined self-adjoint operator on the entire enlarged space is being
assigned this inverse.

Let $K_n(a)=\widetilde P_n^*(\widetilde H_n^{\mathrm{red}}-a)^{-1} \widetilde P_n$. Proposition~\ref{A:res:real-terminal-removal} and isometric invariance of the Hilbert--Schmidt
norm give
\[
 \bigl\|\widehat{\mathcal J}_n
       [(\mathcal A_{n,V,\beta}-a)^{-1}-K_n(a)]
                         \widehat{\mathcal J}_n^*\bigr\|_{\mathfrak S_2}^2
 =o_\omega(h).
 \eqtag{20.63e}
\]
Since $Q_n\widetilde P_n^*=0$ and $\widetilde P_n\widetilde P_n^*=I$, direct multiplication gives
\[
 \widehat{\mathcal J}_nK_n(a)\widehat{\mathcal J}_n^*
 =\bigl[\iota_h^{(M_n)}(\widetilde H_n^{\mathrm{red}}-a)^{-1}
                    (\iota_h^{(M_n)})^*\bigr]\oplus0.
\]
Theorem~\ref{A:res:20.3} therefore implies
\[
 \bigl\|\widehat{\mathcal J}_nK_n(a)\widehat{\mathcal J}_n^*
          -((\mathcal H_\beta-a)^{-1}\oplus0)\bigr\|_{\mathfrak S_2}^2
 \le C_{\omega,V,\beta,a}h.
 \eqtag{20.63g}
\]
Adding the last two differences gives
\[
 \bigl\|\widehat{\mathcal J}_n(\mathcal A_{n,V,\beta}-a)^{-1}
           \widehat{\mathcal J}_n^*
          -((\mathcal H_\beta-a)^{-1}\oplus0)\bigr\|_{\mathfrak S_2}^2
 \le C_{\omega,V,\beta,a}h.
 \eqtag{20.63h}
\]
The embedded inverse has rank $n$. The nonzero singular values of the limiting direct sum are unchanged by its zero summand, so its squared approximation error is at least
$\sum_{j>n}s_j((\mathcal H_\beta-a)^{-1})^2\ge c_{\omega,a}h$. All these comparisons are uniform on compact real subsets of the limiting resolvent set.

Write $\mathfrak d_{2,n}^{\mathrm{full}}(a)$ for \eqref{A:eq:20.59} with all $n$ eigenvalues of the full general Jacobi matrix on the finite side.

\par\Needspace{6\baselineskip}
\begin{theorem}\label{A:res:20.6}
For every deterministic $a\in\mathbb R$, almost surely, for all sufficiently large $n$,
\[
 \bigl\|\widehat{\mathcal J}_n(\mathcal A_{n,V,\beta}-a)^{-1}
            \widehat{\mathcal J}_n^*
            -((\mathcal H_\beta-a)^{-1}\oplus0)\bigr\|_{\mathfrak S_2}^2
 =\Theta_{\omega,a}(n^{-1/3}),
\]
and $\mathfrak d_{2,n}^{\mathrm{full}}(a)^2 =\Theta_{\omega,a}(n^{-1/3})$.
\end{theorem}

\begin{proof}
Theorem~\ref{A:res:20.1} permits the fixed parameter $a$. Invertibility and the inverse estimates were proved above. Isometry preserves each nonzero eigenvalue of the finite
inverse, including its multiplicity. The matching inequality above and \eqref{A:eq:20.63h} give
\[
 \mathfrak d_{2,n}^{\mathrm{full}}(a)^2
 \le\bigl\|\widehat{\mathcal J}_n(\mathcal A_{n,V,\beta}-a)^{-1}
            \widehat{\mathcal J}_n^*
             -((\mathcal H_\beta-a)^{-1}\oplus0)\bigr\|_{\mathfrak S_2}^2
 \le C_{\omega,a}h_n.
\]
Every matching leaves all but at most $n$ limiting entries unmatched. Consequently,
\[
 \mathfrak d_{2,n}^{\mathrm{full}}(a)^2
 \ge\sum_{j>n}s_j((\mathcal H_\beta-a)^{-1})^2
 \ge c_{\omega,a}n^{-1/3}\asymp c_{\omega,a}h_n.
\]
These two inequalities also bound the displayed operator error from below. The construction uses $V_n$ only through
$V_n^*V_n=I$.
\end{proof}

\subsection{Lower bound under an observation constraint and resolvent estimates}
\label{sec:conditional-brownian-bridges-the-information-lower-bound-and-resolvents}
The lower bound for approximation of the limiting inverse is stated in Theorem~\ref{B:thm:brownian-information-lower-bound}. Its proof conditions on the Brownian observations and
estimates the remaining bridges, which are not determined by those observations or by the independent auxiliary variables. We then transfer the inverse estimates to resolvents.

\subsubsection{Conditional Brownian bridges}
\label{B:subsection:conditional-brownian-bridges}

Fix $\beta>0$ and $a>-1$. On a rectangle $z<w$ separated from zero and one, the lower limit in the Green kernel is $w$. Brownian fluctuations in the $z$ variable therefore enter
through the factor $p^B(z)$. We condition on the path outside finitely many intervals in that variable and estimate the independent bridges left inside them. The linear term has
covariance trace of order $n^{-1}$ and operator norm of order $n^{-2}$. These two bounds give a summable probability that a measurable approximation is too close to the inverse.

The inverse $R^B$ is Hilbert--Schmidt by Proposition~\ref{B:result:2-5}. For each $n$ put
\begin{equation}
 \mathcal N_n=\sigma\{W(u_{n,k}):1\le k\le n\},\qquad
 \mathcal F_n=\mathcal N_n\vee\mathcal U_n,
 \label{B:eq:brownian-sampled-information}
\end{equation}
where $\mathcal U_n$ is independent of the entire Brownian path. Let $(T_n)$ be a sequence of Borel $\mathfrak S_2(L^2(0,1))$-valued random variables with $T_n$ measurable with
respect to $\mathcal F_n$. The fields $\mathcal U_n$ and the approximations may be dependent for different $n$. The conditional probability bound below is uniform in the value of
the approximating kernel. We apply it to this sequence after enlarging $\mathcal F_n$ by the exterior path.

Choose fixed intervals $I_-=[x_0,x_1]$ and $I_+=[x_2,x_3]$ with $0<x_0<x_1<x_2<x_3<1$. Write
\begin{equation}
 \mathscr K_n=\{k\ge2:I_{n,k}\subset I_-\},\qquad
 A_n=\bigcup_{k\in\mathscr K_n}I_{n,k},\qquad
 \mathcal D_n=A_n\times I_+,\qquad \mathcal D=I_-\times I_+.
 \label{B:eq:bridge-rectangle}
\end{equation}
Choose $n_0\ge2$ so large that $C/n<x_0$ and $2C/n\le |I_-|/2$ for $n\ge n_0$, using the constants in \eqref{B:eq:brownian-grid-bounds}. The first cell then lies below $I_-$. At
most two other cells meet $I_-$ without being contained in it. Since their total length is at most $2C/n$,
\[
 |I_-|/2\le |A_n|\le |I_-|,\qquad
 \frac{n|I_-|}{2C}\le |\mathscr K_n|\le\frac{n|I_-|}{c}.
\]
For a selected cell, integration of $s^{-1}$ gives
\[
 \frac{c}{nx_1}\le\frac{\Delta_{n,k}}{x_1}
 \le h_{n,k}\le\frac{\Delta_{n,k}}{x_0}\le\frac{C}{nx_0}.
\]
Thus the choices
\[
 s_- =\frac{|I_-|c^2}{2Cx_1^2},\qquad
 s_+ =\frac{|I_-|C^2}{cx_0^2},\qquad H_0=\frac C{x_0}
\]
give
\begin{equation}
 \frac{s_-}{n}\le\sum_{k\in\mathscr K_n}h_{n,k}^2\le\frac{s_+}{n},
 \qquad \max_{k\in\mathscr K_n}h_{n,k}\le H_0/n,
 \qquad n\ge n_0.
 \label{B:eq:bridge-grid-sums}
\end{equation}

For $k\in\mathscr K_n$ put $a_{n,k}=u_{n,k}$, $b_{n,k}=u_{n,k-1}$, and define on $J_{n,k}$
\begin{equation}
 \mathbb B_{n,k}(u)=W(u)-W(a_{n,k})
  -\frac{u-a_{n,k}}{h_{n,k}}
                     [W(b_{n,k})-W(a_{n,k})].
 \label{B:eq:selected-brownian-bridge}
\end{equation}
Enlarge the conditioning field to
\begin{equation}
 \mathcal G_n=\mathcal F_n\vee
 \sigma\{W(u):u\notin\bigcup_{k\in\mathscr K_n}(a_{n,k},b_{n,k})\}.
 \label{B:eq:exterior-brownian-information}
\end{equation}
Every sampling time $u_{n,j}$ lies outside the open intervals removed in \eqref{B:eq:exterior-brownian-information}. In particular, $\mathcal N_n$ is contained in the sigma-field
of the exterior path. The boundary points of the selected intervals are retained in that field.

\begin{lemma}\label{B:lem:conditional-bridge-law}
The processes $(\mathbb B_{n,k})_{k\in\mathscr K_n}$ are mutually independent and independent of $\mathcal G_n$. Their covariance on a cell $[a,b]$ of length $h$ is
\begin{equation}
 C_h(u,v)=u\wedge v-a-\frac{(u-a)(v-a)}h.
 \label{B:eq:brownian-bridge-covariance}
\end{equation}
\end{lemma}
\begin{proof}
Fix one interval $[a,b]$ of length $h$ and write $X(u)=W(u)-W(a)$ and $V=W(b)-W(a)$. For $a\le u,v\le b$,
\[
 \mathbb E X(u)X(v)=u\wedge v-a,\qquad
 \mathbb E X(u)V=u-a,\qquad \mathbb EV^2=h.
\]
Expansion of the covariance of $X(u)-(u-a)V/h$ and $X(v)-(v-a)V/h$ gives
\[
 u\wedge v-a-\frac{(u-a)(v-a)}h
 -\frac{(u-a)(v-a)}h+\frac{(u-a)(v-a)}h
 =C_h(u,v).
\]
The bridge is continuous and vanishes at $a$ and $b$. Its covariance with an exterior evaluation is
\[
 \operatorname{Cov}(\mathbb B(u),W(t))
 =\begin{cases}
  0,&t\le a,\\
  (u-a)-(u-a)h/h=0,&t\ge b.
 \end{cases}
\]
The boundary points are included in these equalities. Each evaluation of a bridge on a distinct interval is a linear combination of exterior evaluations for the present interval,
so the covariances between distinct bridges are also zero.

Let $U$ be a finite real linear combination of bridge evaluations and $V_0$ a finite real linear combination of exterior evaluations. The joint vector is centered Gaussian and
$\mathbb E[UV_0]=0$. Hence
\[
 \mathbb E e^{\ii U+\ii V_0}
 =\exp\!\left[-\tfrac12\mathbb EU^2-\tfrac12\mathbb EV_0^2\right]
 =\mathbb E e^{\ii U}\,\mathbb E e^{\ii V_0}.
\]
Taking separate linear combinations on each selected interval gives factorization over all bridges as well. On every interval choose a countable dense set containing its boundary
points, and do the same on the closed exterior set intersected with $[0,m]$, for each positive integer $m$. Continuity expresses every other evaluation as a limit of evaluations on
these sets. Cylinder events on the countable sets generate the respective path sigma-fields. The monotone class theorem therefore extends the finite-dimensional factorization to
these sigma-fields.

Let $f_k$ be bounded Borel functions of the selected bridges, let $g$ be a bounded function of the exterior path, and let $h_0$ be bounded and $\mathcal U_n$-measurable.
Independence of $\mathcal U_n$ from the entire Brownian path and the factorization just proved give
\[
 \mathbb E\!\left[h_0g\prod_{k\in\mathscr K_n}f_k(\mathbb B_{n,k})\right]
 =\mathbb Eh_0\,\mathbb Eg
       \prod_{k\in\mathscr K_n}\mathbb Ef_k(\mathbb B_{n,k}).
\]
Products $h_0g$ generate the bounded measurable functions for the joined field $\mathcal G_n$. Another monotone class argument yields
\[
 \mathbb E\!\left[\prod_{k\in\mathscr K_n}f_k(\mathbb B_{n,k})
                        \,\middle|\,\mathcal G_n\right]
 =\prod_{k\in\mathscr K_n}\mathbb Ef_k(\mathbb B_{n,k})
 \quad\text{almost surely}.
\]
This proves both independence properties and the stated conditional bridge law.
\end{proof}

Define the continuous path $\widetilde W_n$ by $\widetilde W_n=W$ outside the selected intervals and, on $[a_{n,k},b_{n,k}]$, by
\[
 \widetilde W_n(u)
 =\frac{b_{n,k}-u}{h_{n,k}}W(a_{n,k})
  +\frac{u-a_{n,k}}{h_{n,k}}W(b_{n,k}).
\]
Its evaluations are $\mathcal G_n$-measurable. Since a continuous path is determined by countably many evaluations, it is also measurable with values in continuous path space. On
each selected interval, $W-\widetilde W_n=\mathbb B_{n,k}$. Set
\[
 \xi_n(z)=\beta^{-1/2}\mathbb B_{n,k}(-\log z)
        \quad(z\in I_{n,k},\ k\in\mathscr K_n),
 \qquad \xi_n(z)=0\quad(z\notin A_n).
\]
On the fixed rectangle $\mathcal D$ define
\begin{equation}
 \widetilde R_n(z,w)=z^{a/2}w^{a/2}
    e^{\beta^{-1/2}\widetilde W_n(-\log z)+\mathcal Y(w)}
                       \int_w^1q^B(s)\,\dd s.
 \label{B:eq:conditional-green-kernel}
\end{equation}
For $w\in I_+$ and $s\in[w,1]$, one has $0\le-\log s\le-\log w\le-\log x_2<-\log x_1$. All these times precede the selected logarithmic intervals. The integral and $\mathcal Y(w)$
are consequently determined by $\mathcal G_n$. The positive continuous function $\widetilde R_n$ on $\mathcal D$ is $\mathcal G_n$-measurable. Substituting
$\mathcal Y(z)=\beta^{-1/2}\widetilde W_n(-\log z)+\xi_n(z)$ in \eqref{B:eq:1-2} gives
\begin{equation}
 R^B(z,w)=\widetilde R_n(z,w)e^{\xi_n(z)},\qquad(z,w)\in\mathcal D_n.
 \label{B:eq:bridge-exponential-factorization}
\end{equation}

\begin{figure}[H]
\centering
\begin{tikzpicture}[x=1cm,y=1cm,every node/.style={font=\footnotesize,inner sep=2pt},>=Stealth]
 \node at (2.15,4.85) {$\mathcal D=I_-\times I_+$};
 \draw[->] (0,0)--(4.65,0) node[right] {$z$};
 \draw[->] (0,0)--(0,4.55) node[above] {$w$};
 \draw (0,0) rectangle (4.2,4.2); \draw[densely dashed] (0,0)--(4.2,4.2);
 \fill[black!8] (.8,2.7) rectangle (1.9,3.8);
 \foreach \x/\y in {.88/1.12,1.18/1.42,1.5/1.74}{\fill[black!28] (\x,2.7) rectangle (\y,3.8);}
 \draw (.8,2.7) rectangle (1.9,3.8);
 \draw[densely dashed] (.8,0)--(.8,2.7); \draw[densely dashed] (1.9,0)--(1.9,2.7);
 \draw[densely dashed] (0,2.7)--(.8,2.7); \draw[densely dashed] (0,3.8)--(.8,3.8);
 \node[below] at (.8,-.08) {$x_0$}; \node[below] at (1.9,-.08) {$x_1$};
 \node[left] at (-.08,2.7) {$x_2$}; \node[left] at (-.08,3.8) {$x_3$};
 \node[below] at (4.2,-.08) {$1$}; \node[left] at (-.08,4.2) {$1$};
 \node at (3.1,3.9) {$z<w$};
 \node at (2.1,-.82) {Dark strips: $\mathcal D_n=A_n\times I_+$};
 \node at (10.0,4.85) {One interval in $u=-\log z$};
 \fill[black!8] (6.45,.35) rectangle (7.55,4.15);
 \fill[black!8] (12.1,.35) rectangle (13.5,4.15);
 \draw[->] (6.45,.35)--(13.85,.35) node[right] {$u$};
 \draw[densely dashed] (7.55,.35)--(7.55,4.15);\draw[densely dashed] (12.1,.35)--(12.1,4.15);
 \draw[thick] (6.45,1.75)--(6.75,1.92)--(6.95,1.68)--(7.2,2.02)--(7.55,1.9)
 --(7.82,2.19)--(8.03,2.0)--(8.25,2.61)--(8.5,2.32)--(8.8,2.91)
 --(9.0,2.75)--(9.28,3.3)--(9.6,3.06)--(9.8,3.42)--(10.06,3.18)
 --(10.25,3.55)--(10.54,3.33)--(10.8,3.58)--(11.03,3.14)--(11.32,3.5)
 --(11.58,3.14)--(11.83,3.39)--(12.1,3.15)--(12.38,3.36)--(12.62,3.1)--(12.9,3.3)--(13.2,3.12)--(13.5,3.4);
 \draw[dashed,line width=.8pt] (7.55,1.9)--(12.1,3.15);
 \fill (7.55,1.9) circle (1.7pt); \fill (12.1,3.15) circle (1.7pt);
 \node[below] at (7.55,.25) {$a_{n,k}$}; \node[below] at (12.1,.25) {$b_{n,k}$};
 \node at (10.2,3.94) {$W$}; \node at (10.65,2.45) {$\widetilde W_n$};
 \draw[<->] (9.28,2.375)--(9.28,3.3);
 \node[anchor=west] at (9.48,2.83) {$\mathbb B_{n,k}$};
 \node[align=center] at (10,-.82) {$a_{n,k}=u_{n,k},\quad b_{n,k}=u_{n,k-1}$\\$W-\widetilde W_n=\mathbb B_{n,k}$ on $[a_{n,k},b_{n,k}]$};
\end{tikzpicture}
\caption{The conditioning used for the lower bound under an observation constraint. The left panel shows selected cells in $I_-$ crossed with $I_+$, so $z<w$ throughout. In the right panel, the path is schematic; the exterior values and the values at $a_{n,k},b_{n,k}$ determine the linear interpolant $\widetilde W_n$. Subtracting it gives the bridge of \eqref{B:eq:selected-brownian-bridge}. These bridges are independent of $\mathcal G_n$ by Lemma~\ref{B:lem:conditional-bridge-law}. For $s\in[w,1]$, $-\log s\le-\log x_2<-\log x_1$, so the integral in the Green kernel is determined by $\mathcal G_n$. Consequently \eqref{B:eq:bridge-exponential-factorization} isolates the bridge in the factor $e^{\xi_n(z)}$.}
\label{B:fig:conditional-bridge-region}
\end{figure}
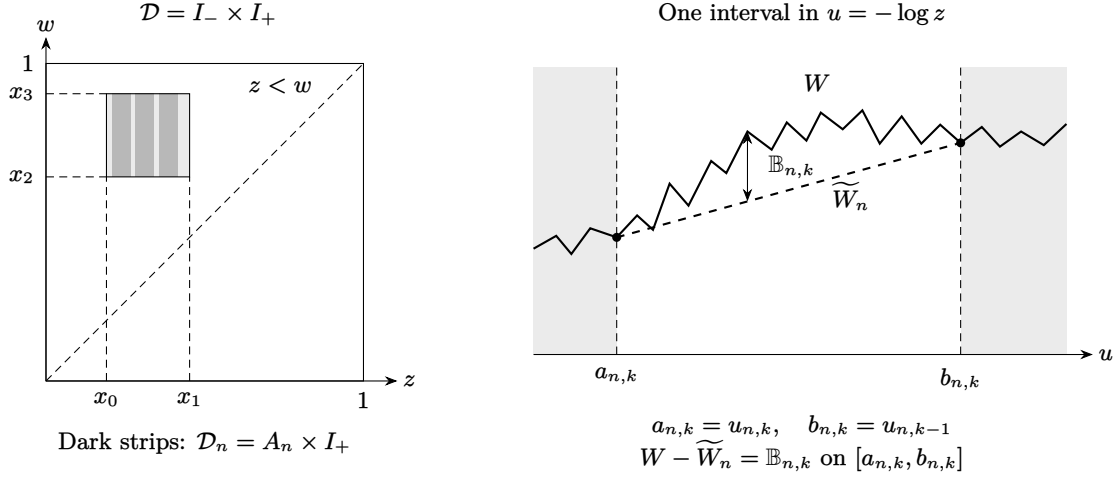

The separated rectangle and a single selected interval of Brownian time are shown in Figure~\ref{B:fig:conditional-bridge-region}. They make explicit why the fluctuation in $z$
factors out of the integral beginning at $w$.

For each integer $M\ge1$, let
\begin{equation}
 \Omega_{n,M}=\{M^{-1}\le\widetilde R_n(z,w)\le M
                          \text{ for all }(z,w)\in\mathcal D\}.
 \label{B:eq:bounded-conditional-kernel}
\end{equation}
Fix a countable dense subset $\mathcal D_*$ of $\mathcal D$. Continuity gives the measurable description
\[
 \Omega_{n,M}
 =\bigcap_{(z,w)\in\mathcal D_*}
       \{M^{-1}\le\widetilde R_n(z,w)\le M\}\in\mathcal G_n.
\]
The events increase with $M$ and satisfy
\begin{equation}
 \mathbb P\left(\bigcup_{M\ge1}\bigcap_{n\ge n_0}\Omega_{n,M}\right)=1.
 \label{B:eq:conditional-kernel-covering}
\end{equation}
To prove it, set $K_W=\sup_{0\le u\le-\log x_0}|W(u)|$. Linear interpolation does not increase this supremum. The positive constants
\[
 c_* =\min_{z\in I_-}z^{a/2}\min_{w\in I_+}w^{a/2}
                              \int_{x_3}^1s^{-a-1}\,\dd s,
 \qquad
 C_* =\max_{z\in I_-}z^{a/2}\max_{w\in I_+}w^{a/2}
                              \int_{x_2}^1s^{-a-1}\,\dd s
\]
are finite and strictly positive for every $a>-1$. All the Brownian values in \eqref{B:eq:conditional-green-kernel} have absolute value at most $K_W$. The two exterior exponential
factors contribute $e^{\pm2K_W/\sqrt\beta}$, and the factor inside the integral contributes another $e^{\pm2K_W/\sqrt\beta}$. Thus
\[
 c_*e^{-4K_W/\sqrt\beta}\le\widetilde R_n(z,w)
       \le C_*e^{4K_W/\sqrt\beta}
 \qquad ((z,w)\in\mathcal D,\ n\ge n_0).
\]
The finite random integer
\[
 \left\lceil\max\{1,c_*^{-1}e^{4K_W/\sqrt\beta},
                         C_*e^{4K_W/\sqrt\beta}\}\right\rceil
\]
satisfies all the inequalities defining \eqref{B:eq:conditional-kernel-covering}. This proves the covering property using only boundedness of a continuous path on a compact
interval.

\subsubsection{The covariance of the linear term}
\label{B:subsection:bridge-covariance-estimates}

The conditional fluctuation is separated into a linear Gaussian term and a nonlinear remainder. Two covariance quantities are needed for the linear term: its trace measures the
total squared error, and its operator norm controls how concentrated that error can be in one direction.

Regard $\mathcal H_n=L^2(\mathcal D_n;\mathbb C)$ as a real Hilbert space, with inner product $\operatorname{Re}\int f\overline g$. Its norm is the usual $L^2$ norm. On
$\mathcal D_n$ put
\begin{equation}
 Z_n(z,w)=\widetilde R_n(z,w)\xi_n(z),\qquad
 N_n(z,w)=\widetilde R_n(z,w)(e^{\xi_n(z)}-1-\xi_n(z)).
 \label{B:eq:bridge-linear-term-and-remainder}
\end{equation}
Conditionally on $\mathcal G_n$, multiplication of the independent bridges by $\widetilde R_n$ gives a centered Gaussian element $Z_n$ of this real Hilbert space. Denote its
covariance by $Q_n$, so that
\[
 \langle Q_nF,G\rangle_{\mathcal H_n}
 =\mathbb E[\langle Z_n,F\rangle_{\mathcal H_n}
                 \langle Z_n,G\rangle_{\mathcal H_n}
                         \mid\mathcal G_n].
\]
Although $Z_n$ is real-valued, the Hilbert space includes arbitrary complex kernels with their usual $L^2$ distance. Define the deterministic constants
\begin{equation}
 s_M=\frac{x_0|I_+|s_-}{6\beta M^2},\qquad
 S_M=\frac{x_1|I_+|M^2s_+}{6\beta},\qquad
 L_M=\frac{x_1|I_+|M^2H_0^2}{\beta\pi^2}.
 \label{B:eq:conditional-covariance-constants}
\end{equation}

\begin{proposition}\label{B:prop:conditional-covariance-bounds}
For $n\ge n_0$, on $\Omega_{n,M}$,
\begin{equation}
 \frac{s_M}{n}\le\operatorname{Tr}Q_n\le\frac{S_M}{n},
 \qquad \|Q_n\|_{\mathrm{op}}\le\frac{L_M}{n^2}.
 \label{B:eq:conditional-covariance-bounds}
\end{equation}
\end{proposition}
\begin{proof}
On an interval $[a,b]$ of length $h$, the covariance diagonal is $C_h(u,u)=(u-a)(b-u)/h$. Integration gives
\begin{equation}
 \int_a^b C_h(u,u)\,\dd u
 =\frac1h\int_0^h v(h-v)\,\dd v=\frac{h^2}{6}.
 \label{B:eq:bridge-covariance-trace}
\end{equation}
For a fixed exterior path, $\widetilde R_n$ is bounded on $\mathcal D$. Each bridge is an $L^2$-valued centered Gaussian variable: Riemann sums give Gaussian linear functionals for
continuous functions. For $f,g\in L^2(a,b)$, Cauchy--Schwarz gives
\[
 \mathbb E\left|\int_a^b\mathbb B(u)(f(u)-g(u))\,\dd u\right|^2
 \le\mathbb E\|\mathbb B\|_{L^2(a,b)}^2\|f-g\|_2^2
 =\frac{h^2}{6}\|f-g\|_2^2.
\]
Approximation by continuous functions therefore gives the Gaussian linear functional for every $L^2$ function. Multiplication and the change of variables defining $Z_n$ are bounded
linear maps on this finite collection of $L^2$ spaces. Conditional independence from Lemma~\ref{B:lem:conditional-bridge-law} therefore gives its centered Gaussian distribution.

Conditional Tonelli and $z=e^{-u}$ yield
\[
 \mathbb E[\|Z_n\|_{\mathcal H_n}^2\mid\mathcal G_n]
 =\frac1\beta\sum_{k\in\mathscr K_n}\int_{J_{n,k}}
 e^{-u}C_{h_{n,k}}(u,u)
           \int_{I_+}|\widetilde R_n(e^{-u},w)|^2\,\dd w\,\dd u.
\]
For a real orthonormal basis $(e_j)$ of $\mathcal H_n$, Parseval and conditional monotone convergence give
\[
 \operatorname{Tr}Q_n
 =\sum_j\mathbb E[\langle Z_n,e_j\rangle_{\mathcal H_n}^2
                                      \mid\mathcal G_n]
 =\mathbb E[\|Z_n\|_{\mathcal H_n}^2\mid\mathcal G_n].
\]
For example, a real orthonormal basis $(\varphi_j)$ of real-valued $L^2(\mathcal D_n)$ gives the real basis $(\varphi_j,\ii\varphi_j)$ of $\mathcal H_n$. The coordinates along
$\ii\varphi_j$ vanish, so this trace has the displayed normalization. On $\Omega_{n,M}$ one has $x_0\le e^{-u}\le x_1$ and
$|I_+|M^{-2}\le\int_{I_+}|\widetilde R_n(e^{-u},w)|^2\,\dd w \le |I_+|M^2$. Consequently,
\[
 \frac{x_0|I_+|}{6\beta M^2}\sum_{k\in\mathscr K_n}h_{n,k}^2
 \le\operatorname{Tr}Q_n
 \le\frac{x_1|I_+|M^2}{6\beta}\sum_{k\in\mathscr K_n}h_{n,k}^2.
\]
Substitution of \eqref{B:eq:bridge-grid-sums} gives $s_M/n$ and $S_M/n$.

For the operator norm, translate a single interval to $(0,h)$ and use the same notation $C_h$ for its covariance operator. For $f\in L^2(0,h)$,
\[
 (C_hf)(u)=\left(1-\frac uh\right)\int_0^u vf(v)\,\dd v
          +\frac uh\int_u^h(h-v)f(v)\,\dd v.
\]
The two primitives are absolutely continuous. The two terms involving $u(h-u)f(u)/h$ cancel on differentiation, leaving
\[
 (C_hf)'(u)=-\frac1h\int_0^u vf(v)\,\dd v
          +\frac1h\int_u^h(h-v)f(v)\,\dd v,
 \qquad (C_hf)''=-f.
\]
The values of $C_hf$ at $0$ and $h$ are zero. For the orthonormal sine basis $e_m(u)=\sqrt{2/h}\sin(m\pi u/h)$, the function $h^2e_m/(\pi^2m^2)$ solves this same equation with the
same boundary values. The difference has second derivative zero and vanishes at both boundary points, and hence is zero. It follows that
\begin{equation}
 C_he_m=\frac{h^2}{\pi^2m^2}e_m,\qquad
 \|C_h\|_{\mathrm{op}}=\frac{h^2}{\pi^2}.
 \label{B:eq:bridge-covariance-eigenvalues}
\end{equation}
Indeed, expansion in the sine basis bounds the operator norm by the largest displayed eigenvalue, and $e_1$ attains that bound.

For a selected cell define the real linear map $V_{n,k}:L^2(J_{n,k};\mathbb R)\to L^2(I_{n,k}\times I_+;\mathbb C)$ by $V_{n,k}f(z,w)=\beta^{-1/2}\widetilde R_n(z,w)f(-\log z)$.
The target carries its real Hilbert space structure. The norm and adjoint are
\[
 \|V_{n,k}f\|_2^2
 =\frac1\beta\int_{J_{n,k}}e^{-u}|f(u)|^2
        \int_{I_+}|\widetilde R_n(e^{-u},w)|^2\,\dd w\,\dd u
 \le\frac{x_1|I_+|M^2}{\beta}\|f\|_2^2
\]
on $\Omega_{n,M}$, and
\[
 (V_{n,k}^*F)(u)=\beta^{-1/2}e^{-u}
          \int_{I_+}\widetilde R_n(e^{-u},w)
                                  \operatorname{Re}F(e^{-u},w)\,\dd w.
\]
To verify the adjoint, substitute the definition of $V_{n,k}$ into $\operatorname{Re}\int V_{n,k}f\,\overline F$ and substitute $z=e^{-u}$. For real $f$, the resulting integral is
$\int_{J_{n,k}} f(u)(V_{n,k}^*F)(u)\,\dd u$.

The rectangle spaces for different cells are orthogonal. For $F=\bigoplus_k F_k$, conditional independence gives
\[
 \mathbb E[\langle Z_n,F\rangle_{\mathcal H_n}^2\mid\mathcal G_n]
 =\sum_{k\in\mathscr K_n}
   \langle C_{h_{n,k}}V_{n,k}^*F_k,V_{n,k}^*F_k\rangle.
\]
Polarization identifies
\[
 Q_n=\bigoplus_{k\in\mathscr K_n}
                       V_{n,k}C_{h_{n,k}}V_{n,k}^*.
\]
In particular,
\[
 \langle Q_nF,F\rangle
 \le\frac{x_1|I_+|M^2}{\beta\pi^2}
     \max_{k\in\mathscr K_n}h_{n,k}^2\sum_k\|F_k\|_2^2
 \le\frac{L_M}{n^2}\|F\|_{\mathcal H_n}^2.
\]
Positivity of $Q_n$ proves the estimate in operator norm.

The same formulas define the covariance without restricting to the localization event. For fixed $n$, all bridge intervals have positive finite length and the exterior multiplier
is bounded on the fixed rectangle. Expanding each $C_{h_{n,k}}$ in its sine basis gives finite sums of positive rank-one operators whose vectors are measurable functions of the
exterior path. Their omitted trace is bounded by
\[
 \sum_{k\in\mathscr K_n}\|V_{n,k}\|^2
          \sum_{m>J}\frac{h_{n,k}^2}{\pi^2m^2}\longrightarrow0.
\]
Thus $Q_n$ is a $\mathcal G_n$-measurable trace-class operator. For $F=\bigoplus_k F_k\in\mathcal H_n$, the Gaussian integrals for the finite sums and their mean square limits give
\[
 \mathbb E[e^{\ii\langle Z_n,F\rangle_{\mathcal H_n}}\mid\mathcal G_n]
 =\exp\!\left[-\frac12\sum_{k\in\mathscr K_n}
     \langle C_{h_{n,k}}V_{n,k}^*F_k,V_{n,k}^*F_k\rangle\right]
 =e^{-\langle Q_nF,F\rangle/2}.
\]
This verifies the conditional Gaussian distribution used above.
\end{proof}

\subsubsection{A Gaussian probability estimate}
\label{B:subsection:gaussian-small-ball-estimate}

Under the conditioning, a measurable approximation becomes a fixed shift of the remaining Gaussian fluctuation. We therefore need a probability estimate uniform in that shift.

Let $Z$ be a centered Gaussian element of a real separable Hilbert space with trace-class covariance $Q$. Put $S=\operatorname{Tr}Q>0$ and $L=\|Q\|_{\mathrm{op}}>0$. The following
estimate is uniform in a deterministic shift. The shift comparison for Gaussian measures is also a consequence of Anderson's inequality \cite{anderson1955integral}; the calculation
below gives the explicit exponent needed here.

\begin{lemma}\label{B:lem:shifted-gaussian-small-ball}
For every $h$ in the Hilbert space and $0\le r^2\le S/4$,
\begin{equation}
 \mathbb P(\|Z-h\|^2\le r^2)\le e^{-S/(8L)}.
 \label{B:eq:shifted-gaussian-small-ball}
\end{equation}
\end{lemma}
\begin{proof}
For $g\sim N(0,1)$, $t>0$, $\lambda\ge0$, and $b\in\mathbb R$, the exponent in the Gaussian integral satisfies
\[
 \frac{x^2}{2}+t(\sqrt\lambda x-b)^2
 =\frac{1+2t\lambda}{2}
     \left(x-\frac{2t\sqrt\lambda b}{1+2t\lambda}\right)^2
       +\frac{tb^2}{1+2t\lambda}.
\]
Integration therefore gives
\[
 \mathbb E e^{-t(\sqrt\lambda g-b)^2}
 =(1+2t\lambda)^{-1/2}e^{-tb^2/(1+2t\lambda)}.
\]
The formula includes $\lambda=0$.

Let $(\lambda_j,e_j)$ be the positive eigenvalues and corresponding orthonormal eigenvectors of $Q$. Its orthogonal complement is $\ker Q$. Every coordinate of $Z$ along a
countable orthonormal basis of that kernel has variance zero and hence vanishes almost surely. For the positive eigenvalues the variables
$g_j=\langle Z,e_j\rangle/\sqrt{\lambda_j}$ are independent standard Gaussians: every finite set is Gaussian with identity covariance. Moreover,
\[
 \mathbb E\sum_j\lambda_jg_j^2=\sum_j\lambda_j=S<\infty,
 \qquad Z=\sum_j\sqrt{\lambda_j}g_je_j
 \quad\text{almost surely and in }L^2.
\]
Parseval gives the equality to $Z$; the mean square error of the first $J$ terms is $\sum_{j>J}\lambda_j$. The argument also applies when there are finitely many positive
eigenvalues.

Write $h=h_0+\sum_jh_je_j$, where $h_0\in\ker Q$. Then
\[
 \|Z-h\|^2=\|h_0\|^2+\sum_j(\sqrt{\lambda_j}g_j-h_j)^2.
\]
Apply the scalar formula to the first $J$ summands and then use bounded convergence for their exponentials. This yields
\begin{equation}
 \mathbb E e^{-t\|Z-h\|^2}
 =e^{-t\|h_0\|^2}
   \prod_j(1+2t\lambda_j)^{-1/2}
              \exp\left(-\sum_j\frac{t h_j^2}{1+2t\lambda_j}\right).
 \label{B:eq:gaussian-laplace-transform}
\end{equation}
Both series determining the right side converge, since
\[
 0\le\sum_j\log(1+2t\lambda_j)\le2tS,
 \qquad
 0\le\sum_j\frac{t h_j^2}{1+2t\lambda_j}\le t\|h\|^2.
\]
The shift terms in \eqref{B:eq:gaussian-laplace-transform} are nonpositive. With $t=(2L)^{-1}$ and $0\le\lambda_j/L\le1$, use
\[
 \log(1+x)=\int_0^x\frac{\dd v}{1+v}\ge\frac x2
 \qquad(0\le x\le1)
\]
to obtain
\[
 \mathbb E e^{-t\|Z-h\|^2}
 \le\exp\!\left[-\frac12\sum_j\log(1+\lambda_j/L)\right]
 \le e^{-S/(4L)}.
\]
On $\{\|Z-h\|^2\le r^2\}$ the exponential is at least $e^{-tr^2}$. Markov's inequality therefore gives
\[
 \mathbb P(\|Z-h\|^2\le r^2)
 \le e^{tr^2}\mathbb E e^{-t\|Z-h\|^2}
 \le\exp\!\left(\frac{r^2}{2L}-\frac{S}{4L}\right)
 \le e^{-S/(8L)}.
\]
Every estimate is uniform in $h$, including its component in $\ker Q$.
\end{proof}

Let $H_n$ be any $\mathcal G_n$-measurable vector in $\mathcal H_n$ and put $b_M=s_M/(8L_M)$. On $\Omega_{n,M}$, Proposition~\ref{B:prop:conditional-covariance-bounds} and the
lemma imply
\begin{equation}
 \mathbb P\left(\|H_n-Z_n\|^2\le\frac{s_M}{4n}
                  \,\middle|\,\mathcal G_n\right)\le e^{-b_M n}.
 \label{B:eq:conditional-gaussian-small-ball}
\end{equation}
The conditional expectation can be computed by integrating the independent bridge laws from Lemma~\ref{B:lem:conditional-bridge-law}, with the exterior multiplier and $H_n$ fixed.
Under this integral the covariance is $Q_n$. On $\Omega_{n,M}$,
\[
 \frac{s_M}{4n}\le\frac{\operatorname{Tr}Q_n}{4},\qquad
 \frac{\operatorname{Tr}Q_n}{8\|Q_n\|}
 \ge\frac{s_M/n}{8L_M/n^2}=b_Mn.
\]
The inequality holds for every fixed multiplier and shift. We may therefore substitute their measurable values in this conditional integral. There is no integrability assumption on
$H_n$; the function being integrated is an indicator. In particular, its imaginary part is allowed, as included in the real Hilbert space $\mathcal H_n$.

\subsubsection{The nonlinear remainder}
\label{B:subsection:bridge-nonlinear-remainder}

We bound the nonlinear remainder by controlling the largest bridge oscillation and the sum of its squared integrals. The resulting failure probability is summable, so the remainder
does not affect the lower bound.

Define
\[
 B_n^*=\max_{k\in\mathscr K_n}\|\mathbb B_{n,k}\|_\infty,
 \qquad S_n^{\mathrm{br}}=
   \sum_{k\in\mathscr K_n}\int_{J_{n,k}}|\mathbb B_{n,k}(u)|^2\,\dd u.
\]
Set $T_0=s_+/6$, $\Lambda_0=H_0^2/\pi^2$, and $K_0=2T_0+4\Lambda_0$. Choose $L_0$ with $2L_0^2/H_0\ge6$, and put
\begin{equation}
 \mathscr E_n=\left\{B_n^*\le L_0\sqrt{\frac{\log n}{n}},\quad
                             S_n^{\mathrm{br}}\le\frac{K_0}{n}\right\}.
 \label{B:eq:bridge-size-event}
\end{equation}

\begin{lemma}\label{B:lem:bridge-size-estimate}
For $n\ge n_0\vee2$,
\begin{equation}
 \mathbb P(\mathscr E_n^c)\le 2n^{-5}+e^{-n}.
 \label{B:eq:bridge-size-probability}
\end{equation}
\end{lemma}
\begin{proof}
Let $W$ be Brownian motion on $[0,h]$, put $V=W(h)$, and set $\mathbb B(u)=W(u)-uV/h$. The covariance calculation in Lemma~\ref{B:lem:conditional-bridge-law} makes $\mathbb B$
independent of $V$. Write $p_h(y)=(2\pi h)^{-1/2}e^{-y^2/(2h)}$. For $t>0$ and $y<t$, reflection gives
\[
 \mathbb P\{\sup_{0\le u\le h}W(u)>t,\ W(h)\in\dd y\}
 =p_h(2t-y)\,\dd y.
\]
Fix $x>0$ and $0<\varepsilon<x/2$. On $\{|V|\le\varepsilon\}$, $\|W-\mathbb B\|_\infty\le\varepsilon$. Independence of $\mathbb B$ and $V$ thus gives
\[
 \frac{\int_{-\varepsilon}^{\varepsilon}
                   p_h(2(x+\varepsilon)-y)\,\dd y}
      {\int_{-\varepsilon}^{\varepsilon}p_h(y)\,\dd y}
 \le\mathbb P\{\sup\mathbb B>x\}
 \le
 \frac{\int_{-\varepsilon}^{\varepsilon}
                   p_h(2(x-\varepsilon)-y)\,\dd y}
      {\int_{-\varepsilon}^{\varepsilon}p_h(y)\,\dd y}.
\]
Both ratios tend to $p_h(2x)/p_h(0)=e^{-2x^2/h}$ as $\varepsilon\downarrow0$. Applying this identity also to $-\mathbb B$ and using the union bound proves
\[
 \mathbb P(\|\mathbb B\|_\infty>x)\le2e^{-2x^2/h}.
\]
The bound at $x=0$ follows from $1\le2$. Since $|\mathscr K_n|\le n$ and $h_{n,k}\le H_0/n$,
\[
 \mathbb P\!\left(B_n^*>L_0\sqrt{\log n/n}\right)
 \le2n\exp\!\left[-\frac{2L_0^2}{H_0}\log n\right]
 \le2n^{-5}.
\]

On each interval, expansion of the bridge in the covariance sine basis gives independent standard Gaussian coordinates $g_{k,m}$. Lemma~\ref{B:lem:conditional-bridge-law} gives
independence between distinct intervals. Parseval therefore gives
\[
 S_n^{\mathrm{br}}=\sum_{k\in\mathscr K_n}\sum_{m\ge1}
           \lambda_{k,m}g_{k,m}^2,
 \qquad \lambda_{k,m}=\frac{h_{n,k}^2}{\pi^2m^2},
\]
with almost sure convergence of the nonnegative sum, because
\[
 \sum_{k,m}\lambda_{k,m}
 =\sum_k\int_{J_{n,k}}C_{h_{n,k}}(u,u)\,\dd u
 =\frac16\sum_kh_{n,k}^2\le\frac{T_0}{n}.
\]
Every coefficient satisfies $\lambda_{k,m}\le\Lambda_0/n^2$. Set $t_n=n^2/(4\Lambda_0)$ and let $S_{n,J}^{\mathrm{br}}$ keep $m\le J$ in the sum. The scalar Gaussian integral and
independence give
\[
 \log\mathbb E e^{t_nS_{n,J}^{\mathrm{br}}}
 =-\frac12\sum_{k\in\mathscr K_n}\sum_{m=1}^J
                          \log(1-2t_n\lambda_{k,m}).
\]
Here $0\le2t_n\lambda_{k,m}\le1/2$, and $-\log(1-x)=\int_0^x(1-v)^{-1}\,\dd v\le2x$ on this range. Consequently
\[
 \log\mathbb E e^{t_nS_{n,J}^{\mathrm{br}}}
 \le2t_n\sum_{k\in\mathscr K_n}\sum_{m=1}^J\lambda_{k,m}
 \le2t_nT_0/n.
\]
The exponentials increase with $J$. Monotone convergence gives $\mathbb E e^{t_nS_n^{\mathrm{br}}}\le e^{2t_nT_0/n}<\infty$. Since $K_0-2T_0=4\Lambda_0$,
\[
 \mathbb P(S_n^{\mathrm{br}}>K_0/n)
 \le\exp[-t_n(K_0-2T_0)/n]
 =\exp\!\left[-\frac{n^2}{4\Lambda_0}
                       \frac{4\Lambda_0}{n}\right]=e^{-n}.
\]
The union bound proves \eqref{B:eq:bridge-size-probability}. All variables used in $\mathscr E_n$ are functions of the selected bridges. Thus $\mathscr E_n$ is independent of
$\mathcal G_n$ and
\[
 \mathbb P(\mathscr E_n^c\mid\mathcal G_n)
 =\mathbb P(\mathscr E_n^c)\le2n^{-5}+e^{-n}.
\]
Independence has been used only among intervals at this fixed size.
\end{proof}

\begin{proposition}\label{B:prop:bridge-remainder-bound}
For each fixed integer $M\ge1$, there is a deterministic $C_M$ such that, for all sufficiently large deterministic $n$,
\begin{equation}
 \|N_n\|_{\mathcal H_n}^2\le C_M\frac{\log n}{n^2}
       \quad\text{on }\Omega_{n,M}\cap\mathscr E_n.
 \label{B:eq:bridge-remainder-bound}
\end{equation}
In particular,
\begin{equation}
 \sum_n\mathbb P\left(\Omega_{n,M}\cap
             \left\{\|N_n\|>\frac{\sqrt{s_M}}{4\sqrt n}\right\}\right)
 <\infty.
 \label{B:eq:bridge-remainder-summability}
\end{equation}
\end{proposition}
\begin{proof}
Choose a deterministic index after which $L_0^2\log n/(\beta n)\le1$. On $\mathscr E_n$,
\[
 \|\xi_n\|_\infty^2
 \le\frac{L_0^2\log n}{\beta n}\le1.
\]
Changing variables on each selected interval gives
\[
 \begin{aligned}
 \int_{A_n}|\xi_n(z)|^2\,\dd z
 &=\frac1\beta\sum_{k\in\mathscr K_n}
       \int_{J_{n,k}}e^{-u}|\mathbb B_{n,k}(u)|^2\,\dd u\\
 &\le\frac{x_1}{\beta}S_n^{\mathrm{br}}
 \le\frac{x_1K_0}{\beta n}.
 \end{aligned}
\]
For real $x$, Taylor's integral identity is
\[
 e^x-1-x=x^2\int_0^1(1-t)e^{tx}\,\dd t,
 \qquad |e^x-1-x|^2\le\frac14e^{2|x|}|x|^4.
\]
On $\Omega_{n,M}\cap\mathscr E_n$, integration over the whole rectangle $\mathcal D_n$ therefore gives
\[
 \begin{aligned}
 \|N_n\|_{\mathcal H_n}^2
 &\le\frac{e^2M^2|I_+|}{4}\int_{A_n}|\xi_n(z)|^4\,\dd z\\
 &\le\frac{e^2M^2|I_+|}{4}\|\xi_n\|_\infty^2
                                  \int_{A_n}|\xi_n(z)|^2\,\dd z\\
 &\le\frac{e^2M^2|I_+|x_1L_0^2K_0}{4\beta^2}
                                       \frac{\log n}{n^2}.
 \end{aligned}
\]
This proves \eqref{B:eq:bridge-remainder-bound} with the deterministic constant $C_M$ in the last line. Increase the deterministic index to $N_M\ge n_0\vee2$ so that
$C_M\log n/n\le s_M/16$ for every $n\ge N_M$. For those sizes,
\[
 \Omega_{n,M}\cap
 \left\{\|N_n\|_{\mathcal H_n}>\frac{\sqrt{s_M}}{4\sqrt n}\right\}
 \subset\mathscr E_n^c.
\]
Taking probabilities and summing gives
\[
 \sum_{n\ge N_M}\mathbb P\!\left(\Omega_{n,M}\cap
       \{\|N_n\|>\sqrt{s_M}/(4\sqrt n)\}\right)
 \le\sum_{n\ge N_M}(2n^{-5}+e^{-n})<\infty.
\]
The finitely many earlier terms are bounded by one. This proves \eqref{B:eq:bridge-remainder-summability}. The probability estimate is conditional on $\mathcal G_n$. The displayed
inclusion accounts for the additional event
$\mathscr E_n$.
\end{proof}

\subsubsection{The lower bound}
\label{B:subsection:brownian-information-lower-bound}

The covariance estimate and the remainder bound now complete the proof of the lower bound. The approximation is required to use only the information in
\eqref{B:eq:brownian-sampled-information}. The following statement includes complex kernels and specifies the dependence of the event of probability one on the sequence $(T_n)$.

\begin{theorem}[Lower bound]\label{B:thm:brownian-information-lower-bound}
Fix $\beta>0$ and $a>-1$, and use the deterministic partitions of Subsection~\ref{B:subsection:brownian-data-and-grids}. Set $\mathcal N_n=\sigma\{W(-\log z_{n,k}):1\le k\le n\}$.
For each $n$, let $\mathcal U_n$ be a sigma-field independent of $\sigma\{W(u):u\ge0\}$. For every sequence $(T_n)$ of Borel $\mathfrak S_2(L^2(0,1))$-valued random variables such
that $T_n$ is $\mathcal N_n\vee\mathcal U_n$-measurable,
\begin{equation}
 \liminf_{n\to\infty}n\|T_n-\mathcal H_{\beta,a}^{-1}\|_{\mathfrak S_2}^2>0
             \qquad\text{almost surely}.
 \label{B:eq:brownian-information-lower-bound}
\end{equation}
The event of probability one may depend on $(T_n)$. Dependence between the fields $\mathcal U_n$ at different sizes is allowed.
\end{theorem}
\begin{proof}
Fix the sequence $(T_n)$. The map from a Hilbert--Schmidt operator to its $L^2$ kernel is an isometry; for an orthonormal basis $(e_j)$,
\[
 K_T(s,t)=\sum_{j\ge1}\overline{e_j(s)}(Te_j)(t),\qquad
 \|K_T\|_2^2=\sum_{j\ge1}\|Te_j\|_2^2.
\]
Restriction to $\mathcal D_n$ is continuous. Hence $H_n=K_{T_n}|_{\mathcal D_n}-\widetilde R_n|_{\mathcal D_n}$ is $\mathcal G_n$-measurable. Equations
\eqref{B:eq:bridge-exponential-factorization} and \eqref{B:eq:bridge-linear-term-and-remainder} give
\[
 D_n:=(K_{T_n}-R^B)|_{\mathcal D_n}=H_n-Z_n-N_n,
 \qquad \|T_n-R^B\|_{\mathfrak S_2}\ge\|D_n\|_{\mathcal H_n}.
\]
Fix an integer $M\ge1$ and put $\kappa_M=s_M/16$. Since $\Omega_{n,M}\in\mathcal G_n$, the conditional probability bound \eqref{B:eq:conditional-gaussian-small-ball} implies
\[
 \mathbb P\bigl(\Omega_{n,M}\cap
       \{\|H_n-Z_n\|^2\le s_M/(4n)\}\bigr)\le e^{-b_Mn}.
\]
If $\|D_n\|^2\le\kappa_M/n$ and $\|N_n\|\le\sqrt{s_M}/(4\sqrt n)$, then $\|H_n-Z_n\|\le\sqrt{s_M}/(2\sqrt n)$. Proposition~\ref{B:prop:bridge-remainder-bound} therefore gives, for
$n\ge N_M$,
\begin{equation}
 \mathbb P\left(\Omega_{n,M}\cap\{\|D_n\|^2\le\kappa_M/n\}\right)
 \le e^{-b_Mn}+2n^{-5}+e^{-n}.
 \label{B:eq:localized-lower-bound-probability}
\end{equation}
The right side is summable. Borel--Cantelli, followed by a countable intersection over $M$, gives probability one to
\[
 \Omega_*=
 \bigcap_{M\ge1}\ \bigcup_{N\ge n_0\vee N_M}\ \bigcap_{n\ge N}
 \left(\Omega_{n,M}^c\cup
                  \{n\|D_n\|_{\mathcal H_n}^2>\kappa_M\}\right).
\]
By \eqref{B:eq:conditional-kernel-covering}, almost every path belongs to $\bigcap_{n\ge n_0}\Omega_{n,M}$ for some finite integer $M$. For such a path in $\Omega_*$,
\[
 \liminf_{n\to\infty}n\|T_n-R^B\|_{\mathfrak S_2}^2
 \ge\liminf_{n\to\infty}n\|D_n\|_{\mathcal H_n}^2
 \ge\kappa_M>0.
\]
The exceptional event may depend on the fixed sequence $(T_n)$. The argument also applies to any prescribed countable collection of sequences by intersecting their events of
probability one.
\end{proof}

\subsubsection{Inverse identities and embedded resolvents}
\label{B:subsection:inverse-identities-and-embedded-resolvents}

The preceding estimates compare positive inverses on $L^2(0,1)$. We identify their bounded transforms with the resolvents on the operator domains and the finite embedded ranges.
The norm bounds on compact spectral sets are then proved separately. Throughout, the finite matrix is inverted before the identifying isometry is applied.

For a bounded operator $S$ on a Hilbert space, set
\[
 \mathscr O(S)=\{z\in\mathbb C:I-zS\text{ is boundedly invertible}\}.
\]
For bounded $R_n,R$ define
\begin{equation}
 G_n(z)=R_n(I-zR_n)^{-1},\qquad G(z)=R(I-zR)^{-1}
 \label{B:eq:resolvents-from-inverses}
\end{equation}
on their respective domains. If $z_0\in\mathscr O(S)$, put $Q_0=(I-z_0S)^{-1}$. For $|z-z_0|\|Q_0S\|<1$,
\[
 I-zS=(I-z_0S)[I-(z-z_0)Q_0S],\qquad
 (I-zS)^{-1}=\sum_{j=0}^{\infty}(z-z_0)^j(Q_0S)^jQ_0.
\]
The series converges in operator norm, uniformly on every smaller closed disk. It proves that $\mathscr O(S)$ is open and that the inverse is continuous there. When $Q_0S=0$ the
series reduces to $Q_0$.

For the finite models, $R_n=\mathscr J_n\mathcal A_n^{-1}\mathscr J_n^*$, where $\mathcal A_n>0$ and $\mathscr J_n^*\mathscr J_n=I$. Put $P_n=\mathscr J_n\mathscr J_n^*$. Then
\[
 I-zR_n=\mathscr J_n(I-z\mathcal A_n^{-1})\mathscr J_n^*+(I-P_n),
 \qquad (I-P_n)\mathscr J_n=0.
\]
Multiplication on both sides by the operator in the next display reduces to $(I-z\mathcal A_n^{-1})(I-z\mathcal A_n^{-1})^{-1}$ on $\operatorname{ran}P_n$ and to $I$ on its
orthogonal complement:
\begin{equation}
 (I-zR_n)^{-1}
 =\mathscr J_n(I-z\mathcal A_n^{-1})^{-1}\mathscr J_n^*+(I-P_n),
 \qquad z\in\rho(\mathcal A_n).
 \label{B:eq:embedded-inverse-resolvent-factor}
\end{equation}
Moreover $I-z\mathcal A_n^{-1}=\mathcal A_n^{-1}(\mathcal A_n-z)$. Thus invertibility on the finite range is equivalent to $z\in\rho(\mathcal A_n)$, and multiplication by $R_n$
gives
\begin{equation}
 \mathscr O(R_n)=\rho(\mathcal A_n),\qquad
 G_n(z)=\mathscr J_n(\mathcal A_n-z)^{-1}\mathscr J_n^*.
 \label{B:eq:embedded-resolvent-identification}
\end{equation}
The inverse factor in \eqref{B:eq:embedded-inverse-resolvent-factor} is the identity on $(\operatorname{ran}P_n)^\perp$; $G_n(z)$ vanishes there. In particular $G_n(0)=R_n$.

For $R=\mathcal H_{\beta,a}^{-1}$, the corresponding identities are
\begin{equation}
 \mathscr O(R)=\rho(\mathcal H_{\beta,a}),\qquad
 G(z)=(\mathcal H_{\beta,a}-z)^{-1}.
 \label{B:eq:bessel-resolvent-identification}
\end{equation}
To check the domain, put $H=R^{-1}$ on $\operatorname{ran}R$. Suppose $Q=(I-zR)^{-1}$ exists. The operator $RQ$ takes values in $\operatorname{dom}H$, and
\[
 (H-z)RQ=(I-zR)Q=I,\qquad
 RQ(H-z)u=RQ(I-zR)Hu=RHu=u
 \quad(u\in\operatorname{dom}H).
\]
It is a bounded two-sided inverse of $H-z$. Conversely, let $G=(H-z)^{-1}$. For $f\in L^2(0,1)$, both $Gf$ and $Rf$ belong to $\operatorname{dom}H$. Applying $R$ to $(H-z)Gf=f$,
and applying $G$ to $(H-z)Rf=f-zRf$, gives
\[
 G-R=zRG,\qquad G-R=zGR.
\]
Consequently
\[
 (I-zR)(I+zG)=I+z(G-R-zRG)=I,\qquad
 (I+zG)(I-zR)=I+z(G-R-zGR)=I.
\]
This proves \eqref{B:eq:bessel-resolvent-identification} and $(I-zR)^{-1}=I+zG(z)$. For nonreal $z$,
\[
 \|(I-zR)^{-1}\|\le1+\frac{|z|}{|\operatorname{Im}z|}.
\]
The finite factor in \eqref{B:eq:embedded-inverse-resolvent-factor} satisfies the same bound, including its identity action on the orthogonal complement. The estimates on compact
sets below also treat real points in the limiting resolvent set, including zero.

Apply Proposition~\ref{prop:shared-resolvent-transfer} with $K=R_n$ and $L=R$. It gives the resolvent difference. Multiplying that identity on the left by $I-zR_n$ and on the right
by $I-zR$ gives
\begin{equation}
 R_n-R=(I-zR_n)(G_n(z)-G(z))(I-zR).
 \label{B:eq:resolvent-inverse-difference}
\end{equation}
The mixed product cancels with its order preserved:
\[
 R_n(I-zR)-(I-zR_n)R=R_n-zR_nR-R+zR_nR=R_n-R.
\]

\subsubsection{Bounds on compact spectral sets}
\label{B:subsection:compact-spectral-sets}

We bound the inverse factors in the resolvent identity uniformly on the chosen compact set. Positivity gives a direct denominator bound outside the positive real axis; convergence
of the inverses then handles compact subsets of the full limiting resolvent set.

For a nonempty compact $K\subset\mathbb C\setminus(0,\infty)$ and $0\le M<\infty$, put
\begin{equation}
 a_K=\max_{z\in K}|z|,\qquad
 d_K(M)=\min_{z\in K,\ 0\le t\le M}|1-zt|.
 \label{B:eq:positive-resolvent-constants}
\end{equation}
The minimum is attained on $K\times[0,M]$. If $1-zt=0$, then $t>0$ and $z=t^{-1}>0$, contradicting the choice of $K$. Thus $d_K(M)>0$; at $M=0$ it equals one. Its dependence on $M$
is continuous. Indeed,
\[
 d_K(M)=\min_{z\in K,\ 0\le r\le1}|1-zMr|,\qquad
 |d_K(M)-d_K(M')|\le a_K|M-M'|.
\]
The second inequality follows by taking minima of two functions whose pointwise difference in absolute value is at most $a_K|M-M'|$. It will establish measurability when a finite
random bound is substituted for $M$.

\begin{lemma}\label{B:lem:positive-inverse-resolvent-bounds}
If $0\le R_n,R\le MI$ and $R_n-R\in\mathfrak S_2$, then
\begin{equation}
 \begin{aligned}
 \frac{\|R_n-R\|_{\mathfrak S_2}}{(1+a_KM)^2}
 &\le\inf_{z\in K}\|G_n(z)-G(z)\|_{\mathfrak S_2},\\
 \sup_{z\in K}\|G_n(z)-G(z)\|_{\mathfrak S_2}
 &\le d_K(M)^{-2}\|R_n-R\|_{\mathfrak S_2}.
 \end{aligned}
 \label{B:eq:positive-inverse-resolvent-bounds}
\end{equation}
\end{lemma}
\begin{proof}
For each $S\in\{R_n,R\}$, positivity and $S\le MI$ give $\sigma(S)\subset[0,M]$. The continuous functional calculus therefore implies, for every $z\in K$,
\[
 \|(I-zS)^{-1}\|
 =\sup_{t\in\sigma(S)}|1-zt|^{-1}\le d_K(M)^{-1},\qquad
 \|I-zS\|\le1+|z|\|S\|\le1+a_KM.
\]
Bound the two outer factors of \eqref{eq:shared-resolvent-identity} in operator norm and the middle factor in Hilbert--Schmidt norm:
\[
 \|G_n(z)-G(z)\|_{\mathfrak S_2}
 \le\|(I-zR_n)^{-1}\|\,\|R_n-R\|_{\mathfrak S_2}
                 \|(I-zR)^{-1}\|
 \le d_K(M)^{-2}\|R_n-R\|_{\mathfrak S_2}.
\]
The reverse identity gives
\[
 \|R_n-R\|_{\mathfrak S_2}
 \le\|I-zR_n\|\,\|G_n(z)-G(z)\|_{\mathfrak S_2}\|I-zR\|
 \le(1+a_KM)^2\|G_n(z)-G(z)\|_{\mathfrak S_2}.
\]
Taking the supremum in the first estimate and the infimum in the second proves \eqref{B:eq:positive-inverse-resolvent-bounds}. At $z=0$ the difference is $R_n-R$, as required by
the definition.
\end{proof}

For a general bounded $R$ and nonempty compact $K\subset\mathscr O(R)$, use $a_K=\max_K|z|$ and define
\begin{equation}
 \begin{split}
 B_K(R)&=\sup_{z\in K}\|(I-zR)^{-1}\|,\qquad
 \epsilon_K(R)=\min\left\{1,\frac1{2(1+a_K)B_K(R)}\right\},\\
 D_K(R)&=\bigl[1+a_K(\|R\|+1)\bigr]\bigl[1+a_K\|R\|\bigr].
 \end{split}
 \label{B:eq:general-resolvent-constants}
\end{equation}
Continuity of the inverse and compactness give $B_K(R)<\infty$.

\begin{lemma}\label{B:lem:compact-inverse-resolvent-bounds}
If $\|S-R\|\le\epsilon_K(R)$, then $K\subset\mathscr O(S)$ and $\sup_{z\in K}\|(I-zS)^{-1}\|\le2B_K(R)$. If also $S-R\in\mathfrak S_2$, then
\begin{align}
 D_K(R)^{-1}\|S-R\|_{\mathfrak S_2}
 &\le\inf_{z\in K}\|S(I-zS)^{-1}-G(z)\|_{\mathfrak S_2},
       \label{B:eq:compact-resolvent-lower}\\
 \sup_{z\in K}\|S(I-zS)^{-1}-G(z)\|_{\mathfrak S_2}
 &\le2B_K(R)^2\|S-R\|_{\mathfrak S_2}.
       \label{B:eq:compact-resolvent-upper}
\end{align}
\end{lemma}
\begin{proof}
Set $E=S-R$. For $z\in K$ define $F_z=z(I-zR)^{-1}E$. The definition of $\epsilon_K(R)$ gives
\[
 \sup_{z\in K}\|F_z\|
 \le a_KB_K(R)\epsilon_K(R)
 \le\frac{a_K}{2(1+a_K)}\le\frac12.
\]
The Neumann series therefore converges uniformly on $K$ and gives $\|(I-F_z)^{-1}\|\le2$. Since
\[
 I-zS=(I-zR)(I-F_z),\qquad
 (I-zS)^{-1}=(I-F_z)^{-1}(I-zR)^{-1},
\]
we have $K\subset\mathscr O(S)$ and $\sup_K\|(I-zS)^{-1}\|\le2B_K(R)$. Substitution in \eqref{eq:shared-resolvent-identity}, with $R_n=S$, yields for every $z\in K$
\[
 \|S(I-zS)^{-1}-G(z)\|_{\mathfrak S_2}
 \le2B_K(R)^2\|S-R\|_{\mathfrak S_2}.
\]
Also $\epsilon_K(R)\le1$, so $\|S\|\le\|R\|+1$. Using \eqref{B:eq:resolvent-inverse-difference},
\[
 \begin{aligned}
 \|S-R\|_{\mathfrak S_2}
 &\le\|I-zS\|\,\|S(I-zS)^{-1}-G(z)\|_{\mathfrak S_2}\|I-zR\|\\
 &\le[1+a_K(\|R\|+1)][1+a_K\|R\|]
                \|S(I-zS)^{-1}-G(z)\|_{\mathfrak S_2}.
 \end{aligned}
\]
The supremum and infimum give the two estimates. If $K=\{0\}$, then $a_K=0$, $B_K(R)=1$, $D_K(R)=1$, and $F_0=0$. Moreover $S(I-0S)^{-1}-G(0)=S-R$. Thus the lower bound is equality
and the upper bound is $\|S-R\|_{\mathfrak S_2}\le2\|S-R\|_{\mathfrak S_2}$, without division by a spectral parameter.
\end{proof}

These estimates use bounded inversion and the ideal property of the Hilbert--Schmidt norm; see Chapter~IV, \S3.3 of \cite{kato1995perturbation} and \cite{simon2005trace}
(Theorem~2.7(a)), respectively. For $R=\mathcal H_{\beta,a}^{-1}$, fix a Brownian path for which the inverse is defined. The second lemma applies to every compact subset of
$\rho(\mathcal H_{\beta,a})$, including compact subsets of positive real spectral gaps. The size of the admissible perturbation depends on that path and on the chosen compact set.

\subsubsection{Resolvent rates and tightness}
\label{B:subsection:resolvent-rates-and-tightness}

The preceding bounds give the same two-sided rate for resolvents on each compact subset of the limiting resolvent set. The comparison also gives an index $n_0(\omega,K)$ from which
the finite resolvents exist throughout the compact set
$K$.

Suppose $R_n$ and $R=\mathcal H_{\beta,a}^{-1}$ are coupled Borel Hilbert--Schmidt operators and put
$e_n=\|R_n-R\|_{\mathfrak S_2}$.

\begin{corollary}\label{B:cor:inverse-and-resolvent-rates}
Suppose that, on an event of probability one, there are $0<c_\omega\le C_\omega<\infty$ and $n_1(\omega)<\infty$ such that $c_\omega\le ne_n^2\le C_\omega$ for all
$n\ge n_1(\omega)$. On the same event, for every nonempty compact $K\subset\rho(\mathcal H_{\beta,a})$, there is $n_0(\omega,K)\ge n_1(\omega)$ such that the finite transforms are
defined on $K$ and, for every $n\ge n_0(\omega,K)$,
\begin{equation}
 \begin{aligned}
 \frac{c_{\omega,K}}n
 &\le\inf_{z\in K}\|G_n(z)-G(z)\|_{\mathfrak S_2}^2,\\
 \sup_{z\in K}\|G_n(z)-G(z)\|_{\mathfrak S_2}^2
 &\le\frac{C_{\omega,K}}n.
 \end{aligned}
 \label{B:eq:inverse-and-resolvent-rates}
\end{equation}
\end{corollary}
\begin{proof}
Fix an outcome on which the assumed bounds hold, and let $N_0$ be an index after which $c_\omega\le ne_n^2\le C_\omega$. Then $\|R_n-R\|\le e_n\le\sqrt{C_\omega/n}$ for $n\ge N_0$.
For any nonempty compact $K\subset\rho(\mathcal H_{\beta,a})$ choose an integer
\[
 N_K\ge\max\{N_0,\,C_\omega\epsilon_K(R)^{-2}\}.
\]
For $n\ge N_K$, Lemma~\ref{B:lem:compact-inverse-resolvent-bounds} gives invertibility on all of $K$ and
\[
 D_K(R)^{-2}e_n^2
 \le\inf_{z\in K}\|G_n(z)-G(z)\|_{\mathfrak S_2}^2
 \le\sup_{z\in K}\|G_n(z)-G(z)\|_{\mathfrak S_2}^2
 \le4B_K(R)^4e_n^2.
\]
Substitute the two bounds on $ne_n^2$ to obtain $c_{\omega,K}=c_\omega D_K(R)^{-2}$ and $C_{\omega,K}=4B_K(R)^4C_\omega$. The outcome was fixed before choosing $K$. Thus the same
event works for every such compact set, while $N_K$ and the constants may depend on the set and the path.
\end{proof}

For a fixed deterministic nonempty compact $K\subset\mathbb C$, set $\mathcal E_K=\{K\subset\mathscr O(R)\}$ and, on this event, define
$Y_{n,K}=n\sup_{z\in K}\|G_n(z)-G(z)\|_{\mathfrak S_2}^2$ when $K\subset\mathscr O(R_n)$, and $Y_{n,K}=+\infty$ otherwise. Set $Y_{n,K}=0$ on $\mathcal E_K^c$. This convention
records a pole of the finite resolvent on $K$, including the case where the limiting resolvent exists throughout $K$.

\begin{proposition}\label{B:prop:resolvent-tightness}
If $e_n\to0$ almost surely and $(ne_n^2)$ is tight, then
\begin{equation}
 \lim_{L\to\infty}\limsup_{n\to\infty}
       \mathbb P(\mathcal E_K\cap\{Y_{n,K}>L\})=0.
 \label{B:eq:resolvent-tightness}
\end{equation}
In particular this gives tightness when $\mathbb P(\mathcal E_K)=1$ and the transforms exist on $K$ for every $n$.
\end{proposition}
\begin{proof}
Choose a countable dense family $(v_j)$ in the unit sphere of the Hilbert space and a countable dense set $K_0\subset K$. For a bounded operator $A$, put $m(A)=\inf_j\|Av_j\|$.
Density gives $m(A)=\inf_{\|v\|=1}\|Av\|$, and
\[
 |m(A)-m(B)|\le\|A-B\|.
\]
If $m(A)>0$, then $\|Av\|\ge m(A)\|v\|$, so $A$ is injective. If $Ax_j$ converges, this inequality makes $x_j$ Cauchy and proves that the range of $A$ is closed. If in addition
$m(A^*)>0$, then $(\operatorname{ran}A)^\perp=\ker A^*=\{0\}$; the range is therefore the whole space. Conversely bounded invertibility gives $m(A)=\|A^{-1}\|^{-1}>0$ and
$m(A^*)>0$.

For bounded $S$, define
\[
 \gamma_K(S)=\inf_{z\in K_0}
    \min\{m(I-zS),\,m((I-zS)^*)\}.
\]
Both functions under the infimum are Lipschitz in $z$, with constant $\|S\|$. Hence their infima over $K_0$ equal those over $K$. The preceding invertibility criterion and
compactness show
\[
 \{K\subset\mathscr O(S)\}=\{\gamma_K(S)>0\}.
\]
In detail, a positive infimum makes both lower bounds positive at every $z\in K$. Conversely, invertibility at every point gives continuous positive lower bounds on the compact
set. Each evaluation $S\mapsto m(I-zS)$ and $S\mapsto m((I-zS)^*)$ is continuous in operator norm. The countable infimum is Borel. This proves measurability of the invertibility
events for the random Hilbert--Schmidt operators $R,R_n$. On $\mathcal E_K$,
\[
 B_K(R)=\left(\inf_{z\in K_0}m(I-zR)\right)^{-1}.
\]
It follows that $B_K(R)$, $\epsilon_K(R)$, and $D_K(R)$ are measurable on this event.

Bounded inversion is jointly continuous in $(z,S)$ on its domain by the local Neumann series. Thus $S(I-zS)^{-1}$ is Borel with values in $\mathfrak S_2$ when $S\in\mathfrak S_2$.
For fixed such $S$ write $F_S(z)=S(I-zS)^{-1}$. Multiplying by $I-zS$ and $I-wS$ gives
\[
 F_S(z)-F_S(w)=(z-w)S^2(I-zS)^{-1}(I-wS)^{-1}.
\]
If $z,w\in K\subset\mathscr O(S)$, then
\[
 \|F_S(z)-F_S(w)\|_{\mathfrak S_2}
 \le |z-w|\,\|S\|_{\mathfrak S_2}\|S\|B_K(S)^2.
\]
The supremum defining $Y_{n,K}$ therefore equals the supremum over $K_0$ whenever both transforms exist on $K$. The conventions on the complementary measurable sets prove
measurability of $Y_{n,K}$ as an extended nonnegative random variable.

On $\mathcal E_K$ let
\[
 N_K=\inf\{m\ge1:\|R_j-R\|\le\epsilon_K(R)
                         \text{ for every }j\ge m\},
\]
with value $+\infty$ for an empty set, and set $N_K=1$ on $\mathcal E_K^c$. For every integer $m\ge1$,
\[
 \mathcal E_K\cap\{N_K\le m\}
 =\mathcal E_K\cap\bigcap_{j\ge m}
       \{\|R_j-R\|\le\epsilon_K(R)\}.
\]
This proves measurability. The assumed convergence makes $N_K$ finite almost surely on $\mathcal E_K$. Set $C_K^*=4B_K(R)^4$ there, with value zero on the complement. On
$\mathcal E_K\cap\{N_K\le n\}$ Lemma~\ref{B:lem:compact-inverse-resolvent-bounds} gives $Y_{n,K}\le C_K^*ne_n^2$. For deterministic $M,L>0$, it follows that
\[
 \begin{aligned}
 \mathbb P(\mathcal E_K\cap\{Y_{n,K}>L\})
 &\le\mathbb P(\mathcal E_K\cap\{N_K>n\})
     +\mathbb P(\mathcal E_K\cap\{C_K^*>M\})\\
 &\quad+\mathbb P(ne_n^2>L/M).
 \end{aligned}
\]
Fix $\varepsilon>0$. Finiteness of $C_K^*$ on $\mathcal E_K$ allows $M$ to be chosen so that the second probability is below $\varepsilon$. Tightness of $(ne_n^2)$ then gives $L_0$
such that $\sup_n\mathbb P(ne_n^2>L/M)<\varepsilon$ for $L\ge L_0$. The first probability tends to zero as $n\to\infty$. Thus the upper limit in \eqref{B:eq:resolvent-tightness} is
at most $2\varepsilon$ for $L\ge L_0$. Letting $\varepsilon\downarrow0$ proves the double limit. Suppose $\mathbb P(\mathcal E_K)=1$ and every finite transform exists on $K$. Fix
$\varepsilon>0$. The preceding double limit gives $L_1<\infty$ with $\limsup_n\mathbb P(Y_{n,K}>L_1)<\varepsilon/2$, and hence an integer $n_1$ such that
$\sup_{n\ge n_1}\mathbb P(Y_{n,K}>L_1)<\varepsilon$. Each of the finitely many variables $Y_{n,K}$ with $n<n_1$ is finite almost surely. Therefore there is $L_2\ge L_1$ such that
\[
 \max_{n<n_1}\mathbb P(Y_{n,K}>L_2)<\varepsilon,\qquad
 \sup_{n\ge1}\mathbb P(Y_{n,K}>L_2)<\varepsilon.
\]
This includes the finitely many initial sizes and proves tightness of the whole sequence.
\end{proof}

For the two matrix couplings, use the notation $R_n$, $R$, $e_n$, and $Z=\sup_nne_n^2$ from Section~\ref{B:subsection:intro-resolvents-and-tightness}. The following proofs apply
the preceding estimates to their full embedded resolvents.

\begin{proof}[Proof of Corollary~\ref{B:cor:intro-resolvent-rates}]
Fix an outcome in the event of the two-sided inverse estimate and put $c_\omega=\tfrac12\liminf_n ne_n^2>0$ and $Z=\sup_n ne_n^2<\infty$. Choose $N_0$ such that
$ne_n^2\ge c_\omega$ for $n\ge N_0$. For a nonempty compact $K\subset\rho(\mathcal H_{\beta,a})$, set
\[
 N_K=\max\{N_0,1,\lceil Z\epsilon_K(R)^{-2}\rceil\}.
\]
For $n\ge N_K$, $\|R_n-R\|_{\rm op}\le e_n\le\sqrt{Z/n}\le\epsilon_K(R)$. Lemma~\ref{B:lem:compact-inverse-resolvent-bounds} and
\eqref{B:eq:embedded-resolvent-identification}--\eqref{B:eq:bessel-resolvent-identification} therefore place $K$ in the finite resolvent set and give
\[
 \begin{aligned}
 \inf_{z\in K}\|G_n(z)-G(z)\|_{\mathfrak S_2}^2
 &\ge\frac{e_n^2}{D_K(R)^2}\ge\frac{c_\omega}{D_K(R)^2n},\\
 \sup_{z\in K}\|G_n(z)-G(z)\|_{\mathfrak S_2}^2
 &\le4B_K(R)^4e_n^2\le\frac{4B_K(R)^4Z}{n}.
 \end{aligned}
\]
Thus $c_{\omega,K}=c_\omega D_K(R)^{-2}$ and $C_{\omega,K}=4B_K(R)^4Z$ are positive and finite. The argument is deterministic after the outcome is fixed, so its event is
independent of the subsequent choice of $K$.
\end{proof}

\begin{proof}[Proof of Corollary~\ref{B:cor:intro-tightness}]
The variable $Z=\sup_nne_n^2$ is measurable as a countable supremum and is finite almost surely. Consequently
\[
 \sup_n\mathbb P(ne_n^2>L)\le\mathbb P(Z>L)\longrightarrow0
 \qquad(L\to\infty).
\]
Also
\[
 M_\omega=\max\{\|R\|_{\rm op},\sup_n\|R_n\|_{\rm op}\}
 \le\|R\|_{\rm op}+\sqrt Z<\infty.
\]
For the compact $K\subset\mathbb C\setminus(0,\infty)$, the random variable $d_K(M_\omega)$ is measurable by the continuity proved after \eqref{B:eq:positive-resolvent-constants},
and is strictly positive. Lemma~\ref{B:lem:positive-inverse-resolvent-bounds} gives
\[
 \sup_{n\ge1}n\sup_{z\in K}\|G_n(z)-G(z)\|_{\mathfrak S_2}^2
 \le d_K(M_\omega)^{-4}Z<\infty.
\]
Each spectral supremum equals one over a countable dense subset of $K$ by the continuity proved in Proposition~\ref{B:prop:resolvent-tightness}. Its tail is therefore bounded,
uniformly in $n$, by $\mathbb P(d_K(M_\omega)^{-4}Z>L)$, which tends to zero.
\end{proof}

\subsection{Removal of the truncation and proof of the theorem for general potentials}
\label{B:sec:removal-of-the-truncation-and-proof-of-the-general-potential-theorem}

The coupled coefficient estimates control the inverse on the selected interval. To compare the full inverses, we estimate products crossing either boundary of that interval. We use
the arrays of Theorem~\ref{B:thm:simultaneous-coefficient-construction} and the estimates of Section~\ref{B:sec:normalization-and-coefficient-estimates} to bound the omitted rows
and columns, as well as the positive correction caused by restricting the triangular factor.

Throughout this subsection $V$, $\beta>0$, $a>0$, and $1<\lambda<\min\{a+1,2\}$ are fixed as in that construction. We keep its intervals $J_n,I_n,I_n^{\rm ret}$, the integers
$r_n,b_n,R_n,m_n,M_n$, and $s_n=r_n-R_*-1$ from Section~\ref{B:sec:general-potentials-and-convex-reference-measures}. Expressions involving these intervals are used from the
deterministic starting size of that construction, increased when necessary. The sampling positions, weights, and operators use the definitions \eqref{B:eq:coupling-clocks},
\eqref{B:eq:criterion-operators}, and \eqref{B:eq:normalization-restricted-inverse}, with $\theta=\theta_V$, $\alpha=4/L_V^2$, and $\mathscr J_n=\mathscr J_{n,V}$. In particular,
$\mathcal R_n=\mathcal R_{n,V,\beta,a}$ and $\mathcal R=\mathcal H_{\beta,a}^{-1}$. For estimates at a fixed size write $z_k=z_{n,k}$, $\Delta_k=\Delta_{n,k}$, and
$\phi_k=\phi(k/n)$. Constants may depend on the fixed parameters; constants carrying $\omega$ may also depend on the sampled coefficients and Brownian path.

We use $\nu_a$ from Section~\ref{B:subsection:row-bounds-near-zero-and-one}. The deterministic product $D_{n;i,j}$ of \eqref{B:eq:initial-coordinate-functionals} is defined by the
same formula for all $1\le i\le j\le n$. The reference function in \eqref{B:eq:normalization-reference-function} gives
\begin{equation}
 D_{n;i,j}=e^{\mathfrak f(i/n)-\mathfrak f(j/n)},\qquad
 D_{n;i,k}D_{n;k,j}=D_{n;i,j}\quad(i\le k\le j).
 \label{B:eq:global-reference-product}
\end{equation}
The finite products $P_{n,k}$ remain those of \eqref{B:eq:2-5}.

\subsubsection{Bounds for individual coefficients}
\label{B:subsection:full-coefficient-bounds}

The omitted part of the inverse involves coefficients outside the selected interval. The first bound combines the interior construction with the original boundary estimates,
allowing a small power loss that will remain integrable in the kernel.

\begin{lemma}\label{B:lem:full-coefficient-bounds}
Almost surely, for every $\varepsilon>0$, there are positive finite constants such that, for all sufficiently large $n$,
\begin{equation}
 C_{\omega,\varepsilon}^{-1}\phi_kz_k^{\varepsilon}
 \le Z_{n,\sigma,k}
 \le C_{\omega,\varepsilon}\phi_kz_k^{-\varepsilon},
 \qquad (\sigma,k)\in\mathscr S_n.
 \label{B:eq:full-coefficient-bounds}
\end{equation}
Here $Z_{n,x,k}=X_{n,k}$ and $Z_{n,y,k}=Y_{n,k}$.
\end{lemma}
\begin{proof}
Equation~\eqref{B:eq:simultaneous-boundary-probabilities} gives an event of probability one and a finite $n_0(\omega)$ such that $\mathcal G_n$ holds for every $n\ge n_0(\omega)$.
Intersect this event with those of Propositions~\ref{B:prop:simultaneous-local-logarithms} and \ref{B:prop:simultaneous-exterior-estimates}. For $k\le s_n$, the bound is exactly
\eqref{B:eq:original-initial-coordinate-bounds}. For $s_n+1\le k<r_n$, the coefficients are in the left boundary configuration. Its defining bound and
\eqref{B:eq:reference-minimizer-relative-error} give
\begin{equation}
 \max_{\substack{\sigma\in\{x,y\}\\s_n<k<r_n}}
 \left|\frac{Z_{n,\sigma,k}}{\phi_k}-1\right|
 \le C\bigl(\vartheta_n+r_n^{-1}+e^{-cr_n}\bigr).
 \label{B:eq:left-boundary-coefficient-bounds}
\end{equation}
Indeed, $k\asymp r_n$ on this fixed width set, and $|Z_{n,\sigma,k}-(z_n^\sharp)_{\sigma,k}|/\phi_k\le\vartheta_n$. The right side of \eqref{B:eq:left-boundary-coefficient-bounds}
tends to zero.

On $J_n$, \eqref{B:eq:conditional-minimizer-relative-centering} and \eqref{B:eq:simultaneous-local-maximum} imply
\[
 \begin{split}
 \max_{\sigma\in\{x,y\},\,k\in J_n}
 \left|\log\frac{m_{\sigma,k}^{q_n}}{\phi_k}\right|
 &\le C\bigl(r_n^{-1}+e^{-cr_n}+\vartheta_n\bigr),\\
 \max_{\sigma\in\{x,y\},\,k\in J_n}
 \left|\log\frac{Z_{n,\sigma,k}}{m_{\sigma,k}^{q_n}}\right|
 &\le C\sqrt{\ell_n/r_n}.
 \end{split}
\]
Put $t_n=C(r_n^{-1}+e^{-cr_n}+\vartheta_n+\sqrt{\ell_n/r_n})$. The two inequalities give, simultaneously for both species and all $k\in J_n$,
\[
 e^{-t_n}\le\frac{Z_{n,\sigma,k}}{\phi_k}\le e^{t_n},
 \qquad t_n\longrightarrow0.
\]
Choose $n_0$ so that the displayed errors are at most $1/2$ for $n\ge n_0$; the ratios then lie in $[1/2,2]$. The conclusion is uniform over the full conditional interval,
including its boundary points. For $k>R_n$, the estimates near the right boundary \eqref{B:eq:terminal-minimizer-bounds} and \eqref{B:eq:terminal-logarithmic-deviation}, applied by
Proposition~\ref{B:prop:simultaneous-exterior-estimates}, give $c\le Z_{n,\sigma,k}\le C$. Here $k/n\ge1/2$ whenever $R_n\ge n/2$ and $c\le\phi_k\le C$. This includes the diagonal
coefficient at $k=n$. These four index ranges cover $\mathscr S_n$. Finally intersect the events for the initial coefficients with $\varepsilon=1/h$, $h\in\mathbb N$. Since
$0<z_k\le1$, the bound with a smaller positive exponent implies the bound with any larger one. This proves the bounds simultaneously for every $\varepsilon>0$.
\end{proof}

\subsubsection{Products before the selected interval}
\label{B:subsection:products-before-retained-interval}

Products crossing from the initial coordinates toward the selected interval require more than individual coefficient bounds. The following cumulative logarithmic estimate controls
products on that intermediate range.

\begin{lemma}\label{B:lem:products-before-retained-interval}
There is a deterministic $C<\infty$ such that, almost surely for all sufficiently large $n$,
\begin{equation}
 \max_{r_n\le i<j\le m_n+1}
 \left|\log\frac{P_{n,i}}{D_{n;i,j}P_{n,j}}\right|
 \le C(1+\sqrt{\ell_n}).
 \label{B:eq:products-before-retained-interval}
\end{equation}
\end{lemma}
\begin{proof}
Condition on the exterior sigma-field $\mathscr E$ of Theorem~\ref{B:thm:simultaneous-coefficient-construction}. On $\mathcal G_n$, the full Gaussian vector of logarithmic ratios
satisfies \eqref{B:eq:simultaneous-ratio-covariance}. Therefore
\[
 \begin{split}
 \operatorname{Var}\left(\sum_{k=i}^{j-1}u_{n,k}^G
                       \,\middle|\,\mathscr E\right)
 &\le C\sum_{k=i}^{j-1}\frac1k
          \sum_{l=i}^{j-1}\sqrt{\frac{k}{l}}e^{-c|k-l|}\\
 &\le C\sum_{k=i}^{j-1}\frac1k
       \sum_{h\in\mathbb Z}(1+|h|)^{1/2}e^{-c|h|}\\
 &\le C\frac{j-i}{i}\le C_1,
        \qquad r_n\le i<j\le m_n+1.
 \end{split}
\]
We used $k/l\le1+|k-l|$ and the bounded ratio $(m_n+1)/r_n$. There are at most $(m_n-r_n+2)^2\le Cr_n^2$ intervals in this maximum. The Gaussian exponential bound and conditional
integration thus give, for deterministic $L>0$,
\begin{equation}
 \mathbb P\left(\mathcal G_n\cap
 \left\{\max_{r_n\le i<j\le m_n+1}
       \left|\sum_{k=i}^{j-1}u_{n,k}^G\right|>L\sqrt{\ell_n}\right\}\right)
 \le Cr_n^2e^{-L^2\ell_n/(2C_1)}.
 \label{B:eq:initial-interval-gaussian-probability}
\end{equation}
For $L^2/(2C_1)>4$, the right side is summable, since $r_n=O((\log n)^{A_0})$. The probabilities of $\mathcal G_n^c$ are also summable. Borel--Cantelli gives a finite $n_0(\omega)$
such that the indicated maximum bound holds for every $n\ge n_0(\omega)$.

For the conditional minimizer, use $\varepsilon_{n,k}^q$ from \eqref{B:eq:normalization-ratio-remainder}. Lemma~\ref{B:lem:normalization-minimizer-ratio} gives
\[
 \sum_{k=i}^{j-1}|\varepsilon_{n,k}^{q_n}|
 \le C\bigg[\frac1{r_n}
 +\frac{e^{-c(n-m_n)}}{1-e^{-c}} 
 +\frac{\vartheta_n}{1-e^{-c}}
       \left(e^{-c(i-r_n)}+e^{-c(R_n-m_n)}\right)\bigg]
 \le C.
\]
The first term uses $\sum_{k\ge r_n}k^{-2}\le2/r_n$; for the three exponential sums, direct evaluation gives
\[
 \begin{split}
 \sum_{k=i}^{j-1}e^{-c(n-k)}
 &=e^{-c(n-j+1)}\frac{1-e^{-c(j-i)}}{1-e^{-c}},\\
 \sum_{k=i}^{j-1}e^{-c(k-r_n)}
 &=e^{-c(i-r_n)}\frac{1-e^{-c(j-i)}}{1-e^{-c}},\\
 \sum_{k=i}^{j-1}e^{-c(R_n-k)}
 &=e^{-c(R_n-j+1)}\frac{1-e^{-c(j-i)}}{1-e^{-c}}.
 \end{split}
\]
Here $j-1\le m_n$ and each numerator $1-e^{-c(j-i)}$ is at most one, which gives precisely the preceding bounds. The identity \eqref{B:eq:global-reference-product} and the
definition of $\varepsilon_{n,k}^q$ now yield
\begin{equation}
 \sum_{k=i}^{j-1}\log\frac{m_{y,k}^{q_n}}{m_{x,k}^{q_n}}
      -\log D_{n;i,j}
       =\sum_{k=i}^{j-1}\varepsilon_{n,k}^{q_n}.
 \label{B:eq:initial-interval-minimizer-product}
\end{equation}
For the constructed logarithmic ratios, \eqref{B:eq:simultaneous-ratio-error} gives the exact relation
\[
 \sum_{k=i}^{j-1}\eta_{n,k}
 =\sum_{k=i}^{j-1}u_{n,k}^G+v_{n,i}-v_{n,j}.
\]
By \eqref{B:eq:ratio-realized-partial-sum-maximum}, $|v_{n,i}-v_{n,j}|\le C\ell_n^2/\sqrt{r_n}=o(1)$ on this range. Finally,
$\log(Y_{n,k}/X_{n,k})= \log(m_{y,k}^{q_n}/m_{x,k}^{q_n})+\eta_{n,k}$ by \eqref{B:eq:normalization-random-ratio}. Their exact sum is
\[
 \log\frac{P_{n,i}}{D_{n;i,j}P_{n,j}}
 =\sum_{k=i}^{j-1}\varepsilon_{n,k}^{q_n}
   +\sum_{k=i}^{j-1}u_{n,k}^G+v_{n,i}-v_{n,j}.
\]
On the intersection of the preceding events of probability one, its absolute value is bounded by $C+L\sqrt{\ell_n}+2C\ell_n^2/\sqrt{r_n}$ for every displayed interval. Since
$A_0>30$, the last term tends to zero; choose $n_0$ so that it is at most one for $n\ge n_0$. Enlarging the deterministic constant gives
\eqref{B:eq:products-before-retained-interval}.
\end{proof}

\subsubsection{Products on the selected interval and near one}
\label{B:subsection:retained-and-terminal-products}

On the selected interval, the coupled logarithmic representation already controls products of coefficients. Near one, the exterior estimate supplies a uniform bound instead. These
two ranges are kept separate because only the first requires a small power loss.

\begin{lemma}\label{B:lem:retained-and-terminal-products}
Almost surely, for every $\varepsilon>0$ and all sufficiently large $n$,
\begin{align}
 \left|\log\frac{P_{n,i}}{D_{n;i,j}P_{n,j}}\right|
 &\le C_{\omega,\varepsilon}
                  +\varepsilon|\log(z_iz_j)|,
                 &&m_n\le i\le j\le M_n,\label{B:eq:retained-product-bound}\\
 \left|\log\frac{P_{n,i}}{D_{n;i,j}P_{n,j}}\right|
 &\le C_\omega,
                 &&M_n\le i\le j\le n.\label{B:eq:terminal-product-bound}
\end{align}
\end{lemma}
\begin{proof}
Taking the ratio of the two identities for $p_{n,i}$ and $p_{n,j}$ in \eqref{B:eq:2-35}, with the scalar of Proposition~\ref{B:prop:coefficient-scalar-normalization}, gives
\[
 \log\frac{P_{n,i}}{P_{n,j}}
 =\frac a2\log\frac{z_i}{z_j}
  +\mathcal Y(z_i)-\mathcal Y(z_j)+\delta_{n,i}-\delta_{n,j}
  +\frac12\log\frac{\Delta_i}{\Delta_j}.
\]
Subtracting \eqref{B:eq:global-reference-product} and using \eqref{B:eq:normalization-profile-identity} yields
\begin{equation}
 \log\frac{P_{n,i}}{D_{n;i,j}P_{n,j}}
 =\mathcal Y(z_i)-\mathcal Y(z_j)+\delta_{n,i}-\delta_{n,j}
       +\frac12(\epsilon_{n,i}^{\Delta}-\epsilon_{n,j}^{\Delta}).
 \label{B:eq:retained-product-identity}
\end{equation}
Lemma~\ref{B:lem:normalization-grid-error} bounds the last difference by $C/n$, and Proposition~\ref{B:prop:normalized-coefficient-maxima} bounds the two $\delta$ terms by a
quantity tending to zero. Lemma~\ref{B:lem:bessel-brownian-growth} gives $|\mathcal Y(z)|\le C_{\omega,\varepsilon}+\varepsilon|\log z|$ on $(0,1]$. Since $z_i,z_j\le1$, the sum of
their absolute logarithms is $|\log(z_iz_j)|$. This proves \eqref{B:eq:retained-product-bound}.

For the interval near the right boundary, Proposition~\ref{B:prop:simultaneous-exterior-estimates} gives, with deterministic constants and a random starting index,
\[
 c\le (z_n^\sharp)_{\sigma,k},Z_{n,\sigma,k}\le C,\qquad
 \left|\log\frac{Z_{n,\sigma,k}}{(z_n^\sharp)_{\sigma,k}}\right|
 \le C\sqrt{\ell_n/n}\quad(k\ge M_n).
\]
Only existing coordinates are included. Consequently, for $M_n\le i<j\le n$,
\[
 \left|\sum_{k=i}^{j-1}
   \left(\log\frac{Y_{n,k}}{X_{n,k}}
       -\log\frac{(z_n^\sharp)_{y,k}}{(z_n^\sharp)_{x,k}}\right)\right|
 \le 2C(n-M_n)\sqrt{\ell_n/n}
 \le Cr_n\sqrt{\ell_n/n}=o(1).
\]
The deterministic estimate \eqref{B:eq:terminal-deterministic-product-bound} implies
\[
 \begin{split}
 \sum_{k=M_n}^{n-1}
  \left|\log\frac{(z_n^\sharp)_{y,k}}{(z_n^\sharp)_{x,k}}\right|
 &\le C\frac{n-M_n}{n}
     +C\sum_{h=1}^{n-M_n}e^{-ch}
     +\frac C{\sqrt n}\sum_{k=M_n}^{n-1}e^{-c(k-K)}\\
 &\le C\frac{r_n}{n}
       +\frac C{e^c-1}
       +\frac{Ce^{-c(M_n-K)}}{\sqrt n(1-e^{-c})}
 \le C.
 \end{split}
\]
For large $n$, $M_n/n\ge1/2$, and Lemma~\ref{B:lem:normalization-grid-error} gives
\[
 |\log D_{n;i,j}|=\left|\int_{i/n}^{j/n}\mathfrak f'(t)\,\dd t\right|
 \le2C\frac{j-i}{n}\le2C\frac{n-M_n}{n}.
\]
The logarithm of $P_{n,i}/P_{n,j}$ is the sum of the centered logarithms near the right boundary and the reference logarithms displayed above. Their two bounds, followed by this
estimate for $D_{n;i,j}$, prove \eqref{B:eq:terminal-product-bound}. The empty product $i=j$ equals one. Intersecting the growth events of Lemma~\ref{B:lem:full-coefficient-bounds}
for a countable sequence of positive exponents tending to zero gives one event for every positive $\varepsilon$.
\end{proof}

\subsubsection{Products over arbitrary index intervals}
\label{B:subsection:products-over-all-indices}

The initial, intermediate, retained and final intervals now each have a product estimate. Splitting an arbitrary product at their boundaries gives a single bound valid throughout
the row, without multiplying individual errors over all indices.

\begin{proposition}\label{B:prop:full-product-bounds}
On an event of probability one, for every $\varepsilon>0$ there are positive finite constants and $n_0(\omega,\varepsilon)<\infty$ such that, for all
$n\ge n_0(\omega,\varepsilon)$,
\begin{equation}
 C_{\omega,\varepsilon}^{-1}D_{n;i,j}(z_iz_j)^\varepsilon
 \le\frac{P_{n,i}}{P_{n,j}}
 \le C_{\omega,\varepsilon}D_{n;i,j}(z_iz_j)^{-\varepsilon},
           \qquad 1\le i\le j\le n.
 \label{B:eq:full-product-bounds}
\end{equation}
\end{proposition}
\begin{proof}
For $i\le j\le s_n+1$, use \eqref{B:eq:original-initial-coordinate-bounds}. For $s_n+1\le i\le j\le r_n$, each coefficient ratio and its reciprocal are bounded by
\eqref{B:eq:left-boundary-coefficient-bounds}. There are at most $R_*$ factors. Moreover,
\[
 |\log D_{n;i,j}|
   \le\int_{i/n}^{j/n}|\mathfrak f'(t)|\,\dd t
   \le C\log(j/i)\le C\log\frac{r_n}{r_n-R_*},
\]
by Lemma~\ref{B:lem:normalization-grid-error}. This is bounded.

For $r_n\le i\le j\le m_n+1$, the bounds on the partition intervals give $z_iz_j\le C(r_n/n)^2$. Therefore
\[
 |\log(z_iz_j)|\ge2\log(n/r_n)-C,
 \qquad
 \frac{1+\sqrt{\ell_n}}{2\log(n/r_n)-C}\longrightarrow0.
\]
Lemma~\ref{B:lem:products-before-retained-interval} now gives \eqref{B:eq:full-product-bounds} on that interval with any positive exponent. The remaining two ranges are covered by
Lemma~\ref{B:lem:retained-and-terminal-products}.

For arbitrary $i<j$, split $[i,j]$ at the points of $\{s_n+1,r_n,m_n,M_n\}$ lying strictly inside it. Write the resulting ordered indices as $i=k_0<\cdots<k_h=j$, with $h\le5$.
Both products in \eqref{B:eq:full-product-bounds} multiply across these indices, by \eqref{B:eq:2-5} and \eqref{B:eq:global-reference-product}. Let $C_{\omega,\varepsilon/10}\ge1$
dominate the finitely many constants from the five ranges, and increase the starting index so that every range estimate holds. The exact product is
\[
 \frac{P_{n,i}}{D_{n;i,j}P_{n,j}}
 =\prod_{l=1}^{h}
       \frac{P_{n,k_{l-1}}}{D_{n;k_{l-1},k_l}P_{n,k_l}}.
\]
Each factor is positive. Applying the range estimates with exponent $\varepsilon/10$ therefore gives the product of the node powers and a constant at most
$C_{\omega,\varepsilon/10}^{5}$ for the upper bound, and its reciprocal for the lower bound. If $A=z_iz_j$, then $0<A\le1$ and $z_{k_l}\ge z_i\ge A$. Thus
\begin{equation}
 \prod_{l=1}^h(z_{k_{l-1}}z_{k_l})^{-\varepsilon/10}
 =A^{-\varepsilon/10}\prod_{l=1}^{h-1}z_{k_l}^{-\varepsilon/5}
 \le A^{-(2h-1)\varepsilon/10}\le A^{-\varepsilon}.
 \label{B:eq:product-subdivision-bound}
\end{equation}
For the lower segment estimates, the reciprocal inequality is $\prod_{l=1}^h(z_{k_{l-1}}z_{k_l})^{\varepsilon/10} \ge(z_iz_j)^\varepsilon$. The product of the at most five random
constants is finite. For $i=j$, both $P_{n,i}/P_{n,j}$ and $D_{n;i,j}$ equal one, and $0<z_i^2\le1$, so both inequalities hold with any constant at least one. For each integer
$h\ge1$, let $\Omega_h$ be the probability one event proved above with exponent $1/h$. Then $\mathbb P((\bigcap_{h\ge1}\Omega_h)^c) \le\sum_{h\ge1}\mathbb P(\Omega_h^c)=0$. On this
intersection, choose $h$ with $1/h\le\varepsilon$. For $A=z_iz_j\in(0,1]$,
\[
 A^\varepsilon\le A^{1/h}\le1\le A^{-1/h}\le A^{-\varepsilon}.
\]
The two bounds at exponent $1/h$ consequently imply the bounds at exponent $\varepsilon$, with that event's finite constant and starting index.
\end{proof}

\subsubsection{The triangular factor and the full inverse}
\label{B:subsection:full-finite-factor-and-inverse}

The product bound has a direct operator interpretation through the explicit bidiagonal inverse. It gives an estimate for the triangular kernel on the whole square, which will
control the rows and columns omitted by restriction.

\begin{lemma}\label{B:lem:full-finite-factor-bound}
Almost surely, for every $\xi>0$ and all sufficiently large $n$,
\begin{equation}
 |\widetilde K_n(i,l)|
 \le C_{\omega,\xi}z_i^{a/2-\xi}z_l^{-(a+1)/2-\xi},
                    \qquad 1\le i\le l\le n.
 \label{B:eq:full-finite-factor-bound}
\end{equation}
\end{lemma}
\begin{proof}
Let $\psi_i=\Delta_i^{-1/2}\mathbf1_{[z_{i-1},z_i)}$. The factor is $\widetilde{\mathcal K}_n=(\sqrt\alpha/n) \widetilde{\mathscr J}_nB_n^{-1}\widetilde{\mathscr J}_n^*$. For
$f\in L^2(0,1)$, back substitution gives
\[
 (\widetilde{\mathcal K}_nf)(t)
 =\frac{\sqrt\alpha}{n}\sum_{l=1}^n\psi_l(t)
       \sum_{i\le l}\frac{P_{n,i}}{X_{n,l}P_{n,l}}
                          \int_0^1\psi_i(s)f(s)\,\dd s.
\]
Thus the first kernel variable is integrated. Substituting the two cell indicators gives
\begin{equation}
 \widetilde K_n(i,l)
 =\frac{\sqrt\alpha}{n\sqrt{\Delta_i\Delta_l}}
       \frac{P_{n,i}}{X_{n,l}P_{n,l}}\mathbf1_{\{i\le l\}}.
 \label{B:eq:full-finite-factor-formula}
\end{equation}
As before, the first kernel variable is integrated and the resulting function is evaluated at the second; see Section~\ref{B:subsection:embedded-kernels-and-scalar-invariance}.
Here $n\sqrt{\Delta_i\Delta_l}$ is bounded above and below uniformly. From the definition of $D_{n;i,l}$,
\[
 \frac{D_{n;i,l}}{\phi_l}
 =\left(\frac{z_i}{z_l}\right)^{a/2+1/4}
       \frac1{\sqrt{\phi_i\phi_l}}
 \le C z_l^{-1/2}\left(\frac{z_i}{z_l}\right)^{a/2},
\]
since $\phi_k\asymp\sqrt{z_k}$. Apply Lemma~\ref{B:lem:full-coefficient-bounds} and Proposition~\ref{B:prop:full-product-bounds} with exponent $\xi/2$. Their product is bounded by
\[
 C_{\omega,\xi}\frac{D_{n;i,l}}{\phi_l}
            (z_iz_l)^{-\xi/2}z_l^{-\xi/2}
 \le C_{\omega,\xi}z_l^{-1/2}
       (z_i/z_l)^{a/2}(z_iz_l)^{-\xi},
\]
where $z_l^{-\xi/2}\le(z_iz_l)^{-\xi/2}$. This proves \eqref{B:eq:full-finite-factor-bound}.
\end{proof}

\begin{proposition}\label{B:prop:full-finite-inverse-bounds}
Almost surely,
\begin{equation}
 \sup_{n\ge1}\operatorname{Tr}\widetilde{\mathcal R}_n<\infty,
 \qquad
 \sup_{n\ge1}\|\widetilde{\mathcal R}_n\|_{\mathfrak S_2}<\infty.
 \label{B:eq:full-finite-inverse-norms}
\end{equation}
For $0<\eta<\nu_a$, there is $n_0(\omega,\eta)<\infty$ such that, for all $n\ge n_0(\omega,\eta)$,
\begin{equation}
 \sum_{j=1}^n\Delta_j|\widetilde R_n(i,j)|^2
       \le C_{\omega,\eta}z_i^{\nu_a-\eta},
               \qquad 1\le i\le n.
 \label{B:eq:full-finite-row-bound}
\end{equation}
\end{proposition}
\begin{proof}
The cell representation of the inverse is
\begin{equation}
 \widetilde R_n(i,j)=\sum_{l\ge i\vee j}\Delta_l
                  \widetilde K_n(i,l)\widetilde K_n(j,l).
 \label{B:eq:full-finite-green-formula}
\end{equation}
Choose $\xi>0$ with $2\xi<a+1$ and $4\xi<1$. Lemma~\ref{B:lem:full-interval-power-integrals} and \eqref{B:eq:full-finite-factor-bound} imply
\[
 \begin{split}
 \operatorname{Tr}\widetilde{\mathcal R}_n
 &=\|\widetilde{\mathcal K}_n\|_{\mathfrak S_2}^2
   =\sum_l\Delta_l\sum_{i\le l}\Delta_i
                              |\widetilde K_n(i,l)|^2\\
 &\le C_{\omega,\xi}\sum_l\Delta_l z_l^{-a-1-2\xi}
                         \sum_{i\le l}\Delta_i z_i^{a-2\xi}\\
 &\le C_{\omega,\xi}\sum_l\Delta_l z_l^{-4\xi}
 \le C_{\omega,\xi}.
 \end{split}
\]
The trace identity is \eqref{B:eq:full-factor-norm}, applied to this finite factor. Positivity gives
$\|\widetilde{\mathcal R}_n\|_{\mathfrak S_2} \le\operatorname{Tr}\widetilde{\mathcal R}_n$. Every earlier inverse is a positive finite matrix embedded by an isometry. The maximum
of their finitely many traces is finite, which proves \eqref{B:eq:full-finite-inverse-norms}.

For $i\le j$, \eqref{B:eq:full-finite-green-formula} gives
\[
 \widetilde R_n(i,j)
 \le C_{\omega,\xi}z_i^{a/2-\xi}z_j^{a/2-\xi}
                      \sum_{l\ge j}\Delta_l z_l^{-a-1-2\xi}
 \le C_{\omega,\xi}z_i^{a/2-\xi}z_j^{-a/2-3\xi}.
\]
Taking $3\xi\le\varepsilon$ and using symmetry yields, for every $\varepsilon>0$,
\begin{equation}
 |\widetilde R_n(i,j)|
 \le C_{\omega,\varepsilon}
      (z_{i\wedge j}/z_{i\vee j})^{a/2}(z_iz_j)^{-\varepsilon}.
 \label{B:eq:full-finite-green-bound}
\end{equation}
The row norm splits into disjoint index sets,
\[
 \sum_{j=1}^n\Delta_j|\widetilde R_n(i,j)|^2
 =\sum_{j\le i}\Delta_j|\widetilde R_n(i,j)|^2
       +\sum_{j>i}\Delta_j|\widetilde R_n(i,j)|^2.
\]
The diagonal cell occurs once. To track the exponent, set $\varepsilon=\eta/16$. Since $0<\eta<\nu_a\le\min\{a,1\}$,
\[
 a-2\varepsilon>-1,\qquad
 1-4\varepsilon-(\nu_a-\eta)=1-\nu_a+3\eta/4>0.
\]
The contribution of $j\le i$ is consequently at most
\[
 C_{\omega,\eta}z_i^{-a-2\varepsilon}
           \sum_{j\le i}\Delta_j z_j^{a-2\varepsilon}
 \le C_{\omega,\eta}z_i^{1-4\varepsilon}
 \le C_{\omega,\eta}z_i^{\nu_a-\eta}.
\]
For $j>i$, write the bound as
\[
 \begin{split}
 z_i^{a-2\varepsilon}\sum_{j>i}\Delta_jz_j^{-a-2\varepsilon}
 &=z_i^{\nu_a-\eta}\sum_{j>i}\Delta_j
       \left(\frac{z_i}{z_j}\right)^{a-\nu_a+\eta-2\varepsilon}
                     z_j^{-\nu_a+\eta-4\varepsilon}\\
 &\le Cz_i^{\nu_a-\eta}.
 \end{split}
\]
The ratio exponent is positive, and $-\nu_a+\eta-4\varepsilon=-\nu_a+3\eta/4>-1$. Set $\gamma=-\nu_a+3\eta/4>-1$. The ratio $(z_i/z_j)^{a-\nu_a+7\eta/8}$ is at most one for $j>i$,
and \eqref{B:eq:full-interval-power-bounds} gives $\sum_{j=1}^n\Delta_jz_j^\gamma\le C_\gamma$ uniformly in $n$. This proves \eqref{B:eq:full-finite-row-bound}, including
$a=1$.
\end{proof}

\subsubsection{The omitted rows and columns}
\label{B:subsection:omitted-finite-rows-and-columns}

The omitted rows and columns lie outside the restricted kernel square. Their contribution can therefore be estimated directly from the full triangular factor, separately from the
correction inside that square.

Use the projections of the inverse criterion with $k_-=m_n$ and $k_+=M_n$. Thus $P_n^0$ is the projection already used in \eqref{B:eq:normalization-restricted-inverse}, while
$P_n^-$ and $P_n^+$ select the cells $k\le m_n$ and $k>M_n$. Write
\begin{equation}
 d_n=z_{n,m_n},\qquad \bar z_n=z_{n,M_n},\qquad
 \varepsilon_n^+=1-\bar z_n.
 \label{B:eq:finite-truncation-lengths}
\end{equation}
The bounds on $\theta'$ give $d_n\asymp\varepsilon_n^+\asymp r_n/n$. We increase the deterministic initial size so that
$\bar z_n\ge1/2$.

\begin{proposition}\label{B:prop:omitted-finite-rows-and-columns}
For every $0<\eta<\nu_a$, there is an almost surely finite $n_0(\omega,\eta)$ such that, for all $n\ge n_0(\omega,\eta)$,
\begin{align}
 \|P_n^-\widetilde{\mathcal R}_n\|_{\mathfrak S_2}^2
 +\|\widetilde{\mathcal R}_nP_n^-\|_{\mathfrak S_2}^2
 &\le C_{\omega,\eta}d_n^{1+\nu_a-\eta},
       \label{B:eq:finite-rows-near-zero}\\
 \|P_n^+\widetilde{\mathcal R}_n\|_{\mathfrak S_2}^2
 +\|\widetilde{\mathcal R}_nP_n^+\|_{\mathfrak S_2}^2
 &\le C_\omega(\varepsilon_n^+)^3.
       \label{B:eq:finite-rows-near-one}
\end{align}
\end{proposition}
\begin{proof}
Symmetry of the inverse kernel and \eqref{B:eq:full-finite-row-bound} give
\[
 \begin{split}
 \|P_n^-\widetilde{\mathcal R}_n\|_{\mathfrak S_2}^2
 +\|\widetilde{\mathcal R}_nP_n^-\|_{\mathfrak S_2}^2
 &=2\sum_{i\le m_n}\Delta_i\sum_j\Delta_j
                                 |\widetilde R_n(i,j)|^2\\
 &\le C_{\omega,\eta}\sum_{i\le m_n}\Delta_i z_i^{\nu_a-\eta}
 \le C_{\omega,\eta}d_n^{1+\nu_a-\eta}.
 \end{split}
\]
The last step is Lemma~\ref{B:lem:full-interval-power-integrals}.

For the other end choose $0<\xi<\min\{a/4,1/8\}$ and put $\tau_j=\sum_{l=j}^n\Delta_l=1-z_{j-1}$, with $\tau_{n+1}=0$. For $i\le j$ and $j>M_n$, the lower bound $z_j\ge1/2$ and
\eqref{B:eq:full-finite-factor-bound} imply
\begin{equation}
 |\widetilde R_n(i,j)|
 \le C_{\omega,\xi}z_i^{a/2-\xi}z_j^{a/2-\xi}
                      \sum_{l\ge j}\Delta_lz_l^{-a-1-2\xi}
 \le C_{\omega,\xi}z_i^{a/2-\xi}\tau_j.
 \label{B:eq:terminal-finite-kernel-bound}
\end{equation}
For $i>j>M_n$, symmetry instead bounds the kernel by $C_{\omega,\xi}\tau_i$, and $\tau_i\le\tau_j$. Hence
\[
 \sum_i\Delta_i|\widetilde R_n(i,j)|^2
 \le C_{\omega,\xi}\tau_j^2
       \left(\sum_{i\le j}\Delta_i z_i^{a-2\xi}
                         +\sum_{i>j}\Delta_i\right)
 \le C_{\omega,\xi}\tau_j^2.
\]
Since $\tau_j-\tau_{j+1}=\Delta_j$,
\[
 \tau_j^3-\tau_{j+1}^3-\Delta_j\tau_j^2
 =\Delta_j\tau_{j+1}(\tau_j+\tau_{j+1})\ge0.
\]
The sum from $M_n+1$ to $n$, with $\tau_{n+1}=0$, is
\[
 \sum_{j=M_n+1}^n\Delta_j\tau_j^2
 =(\varepsilon_n^+)^3
 -\sum_{j=M_n+1}^n\Delta_j\tau_{j+1}(\tau_j+\tau_{j+1}).
\]
The last sum is nonnegative. Its term at $j=n$ is zero, whereas the original sum contains $\Delta_n\tau_n^2=\Delta_n^3$. Consequently
\begin{equation}
 \sum_{j>M_n}\Delta_j(1-z_{j-1})^2
       \le(1-z_{M_n})^3=(\varepsilon_n^+)^3.
 \label{B:eq:terminal-cell-sum}
\end{equation}
Multiplying the row estimate by $\Delta_j$, summing, and doubling for the equal column norm proves \eqref{B:eq:finite-rows-near-one}.
\end{proof}

\subsubsection{Restricting the triangular factor}
\label{B:subsection:finite-factor-restriction}

Restricting the triangular factor and then forming its positive square is not the same as compressing the full inverse. The difference is the positive correction below, whose
rank-one structure makes its size explicit.

Let $\mathcal T_n$ be the correction defined before Lemma~\ref{B:lem:finite-truncation-error-equivalence}, specialized to $k_-=m_n$ and $k_+=M_n$.

\begin{lemma}\label{B:lem:finite-factor-restriction}
The operator $\mathcal T_n$ is positive and has rank one. Moreover,
\begin{equation}
 P_n^0\widetilde{\mathcal R}_nP_n^0
       =\widetilde{\mathcal R}_n^0+\mathcal T_n,
 \qquad
 \|\mathcal T_n\|_{\mathfrak S_2}^2
       \le C_\omega(\varepsilon_n^+)^2
 \label{B:eq:finite-factor-restriction-bound}
\end{equation}
almost surely for all sufficiently large $n$.
\end{lemma}
\begin{proof}
Triangularity gives $P_n^-\widetilde{\mathcal K}_nP_n^0=0$. Using $I=P_n^-+P_n^0+P_n^+$ between the two factors gives
\[
 \begin{split}
 P_n^0\widetilde{\mathcal R}_nP_n^0
 &=P_n^0\widetilde{\mathcal K}_n^*P_n^0
                         \widetilde{\mathcal K}_nP_n^0
   +P_n^0\widetilde{\mathcal K}_n^*P_n^+
                         \widetilde{\mathcal K}_nP_n^0\\
 &=\widetilde{\mathcal R}_n^0+\mathcal T_n.
 \end{split}
\]
For $i,j\in I_n^{\rm ret}$, its second summand has kernel
\begin{equation}
 T_n(i,j)=p_{n,i}p_{n,j}\sum_{l=M_n+1}^nq_{n,l},\qquad
 \|\mathcal T_n\|_{\mathfrak S_2}^2
 =\left(\sum_{l=M_n+1}^nq_{n,l}\right)^2
       \left(\sum_{i\in I_n^{\rm ret}}\Delta_i p_{n,i}^2\right)^2.
 \label{B:eq:finite-factor-correction-norm}
\end{equation}
Write $q_+=\sum_{l=M_n+1}^nq_{n,l}$ and $w=\sum_{i\in I_n^{\rm ret}}p_{n,i}\mathbf1_{[z_{i-1},z_i)}$. Then, with the first kernel variable integrated,
\[
 (\mathcal T_nf)(t)=q_+w(t)\int_0^1w(s)f(s)\,\dd s,
 \qquad
 \int_0^1\overline{f(t)}(\mathcal T_nf)(t)\,\dd t
       =q_+\left|\int_0^1w(s)f(s)\,\dd s\right|^2\ge0.
\]
All coefficients are positive and both index sets are nonempty, so $q_+>0$ and $w\ne0$. Its sole nonzero eigenvalue is $q_+\|w\|_2^2$, giving rank one and the squared norm in
\eqref{B:eq:finite-factor-correction-norm}.

To estimate the product in \eqref{B:eq:finite-factor-correction-norm}, choose $\xi$ as in the proof of Proposition~\ref{B:prop:omitted-finite-rows-and-columns}. The representation
in terms of the factor gives
\[
 \begin{split}
 0\le T_n(i,j)
 &=\sum_{l>M_n}\Delta_l
                      \widetilde K_n(i,l)\widetilde K_n(j,l)\\
 &\le C_{\omega,\xi}z_i^{a/2-\xi}z_j^{a/2-\xi}
                     \sum_{l>M_n}\Delta_lz_l^{-a-1-2\xi}\\
 &\le C_{\omega,\xi}\varepsilon_n^+
                            z_i^{a/2-\xi}z_j^{a/2-\xi}.
 \end{split}
\]
The last step uses $z_l\ge\bar z_n\ge1/2$. Squaring and summing over the two selected cell indices yields
\[
 \|\mathcal T_n\|_{\mathfrak S_2}^2
 \le C_{\omega,\xi}(\varepsilon_n^+)^2
         \left(\sum_{i\in I_n^{\rm ret}}\Delta_i z_i^{a-2\xi}\right)^2
 \le C_{\omega,\xi}(\varepsilon_n^+)^2.
\]
The coefficient product and the kernel are unchanged by the scalar normalization, by \eqref{B:eq:2-20}.
\end{proof}

\subsubsection{The truncation error}
\label{B:subsection:finite-inverse-truncation-error}

The two types of truncation error have disjoint kernel supports: one lies outside the square of the restricted kernel, and the other is the correction inside it. Their squared
Hilbert--Schmidt norms therefore add.

\begin{theorem}\label{B:thm:finite-inverse-truncation}
For every $0<\eta<\nu_a$, there is an almost surely finite $n_0(\omega,\eta)$ such that, for all $n\ge n_0(\omega,\eta)$,
\begin{equation}
 \|\widetilde{\mathcal R}_n-\widetilde{\mathcal R}_n^0\|_{\mathfrak S_2}^2
 \le C_{\omega,\eta}
     \left[d_n^{1+\nu_a-\eta}+(\varepsilon_n^+)^3
                                  +(\varepsilon_n^+)^2\right]
 =o_\omega(n^{-1}).
 \label{B:eq:finite-inverse-truncation-error}
\end{equation}
\end{theorem}
\begin{proof}
The kernel of $\widetilde{\mathcal R}_n-P_n^0\widetilde{\mathcal R}_nP_n^0$ vanishes on the square of the restricted kernel, while the kernel of $\mathcal T_n$ vanishes outside it.
Their Hilbert--Schmidt inner product is therefore zero. By Lemma~\ref{B:lem:finite-factor-restriction},
\begin{equation}
 \|\widetilde{\mathcal R}_n-\widetilde{\mathcal R}_n^0\|_{\mathfrak S_2}^2
 =\|\widetilde{\mathcal R}_n
       -P_n^0\widetilde{\mathcal R}_nP_n^0\|_{\mathfrak S_2}^2
                         +\|\mathcal T_n\|_{\mathfrak S_2}^2.
 \label{B:eq:finite-truncation-orthogonality}
\end{equation}
Put $P=P_n^0$ and $Q=I-P$. The first difference decomposes as
\[
 \widetilde{\mathcal R}_n-P\widetilde{\mathcal R}_nP
 =Q\widetilde{\mathcal R}_nQ+Q\widetilde{\mathcal R}_nP
                          +P\widetilde{\mathcal R}_nQ.
\]
These three kernels have disjoint supports. Self-adjointness gives equal norms to the last two, hence
\[
 \|\widetilde{\mathcal R}_n-P\widetilde{\mathcal R}_nP\|_{\mathfrak S_2}^2
 =\|Q\widetilde{\mathcal R}_nQ\|_{\mathfrak S_2}^2
                +2\|Q\widetilde{\mathcal R}_nP\|_{\mathfrak S_2}^2.
\]
Equivalently this is the sum of $\Delta_i\Delta_j|\widetilde R_n(i,j)|^2$ over pairs with at least one index outside $I_n^{\rm ret}$. Since
$\mathbf1_{\{i\notin I_n^{\rm ret}\ \mathrm{or}\ j\notin I_n^{\rm ret}\}} \le\mathbf1_{\{i\notin I_n^{\rm ret}\}}+ \mathbf1_{\{j\notin I_n^{\rm ret}\}}$, it is at most
\[
 \sum_{\sigma\in\{-,+\}}
  \left(\|P_n^\sigma\widetilde{\mathcal R}_n\|_{\mathfrak S_2}^2
       +\|\widetilde{\mathcal R}_nP_n^\sigma\|_{\mathfrak S_2}^2\right).
\]
Proposition~\ref{B:prop:omitted-finite-rows-and-columns} and Lemma~\ref{B:lem:finite-factor-restriction} prove the displayed upper bound. To compute its order, use
$d_n,\varepsilon_n^+\le Cr_n/n$ and $\nu_a-\eta>0$:
\begin{align}
 nd_n^{1+\nu_a-\eta}
 &\le Cr_n^{1+\nu_a-\eta}n^{-(\nu_a-\eta)}\longrightarrow0,\notag\\
 n(\varepsilon_n^+)^3&\le Cr_n^3n^{-2}\longrightarrow0,
 \qquad
 n(\varepsilon_n^+)^2\le Cr_n^2n^{-1}\longrightarrow0.
 \label{B:eq:finite-truncation-orders}
\end{align}
For $n\ge e$ one has $r_n\le2(\log n)^{A_0}$. Writing $\sigma=\nu_a-\eta>0$ bounds the three expressions by constant multiples of
\[
 (\log n)^{A_0(1+\sigma)}n^{-\sigma},\qquad
 (\log n)^{3A_0}n^{-2},\qquad
 (\log n)^{2A_0}n^{-1}.
\]
For any fixed $b,c>0$, the logarithm of $(\log n)^b n^{-c}$ is $b\log\log n-c\log n\to-\infty$. Each expression therefore tends to zero, proving
\eqref{B:eq:finite-inverse-truncation-error}.
\end{proof}

\subsubsection{Convergence of the full inverse}
\label{B:subsection:full-inverse-upper-bound}

The full inverse can now be compared through the two Brownian approximations: the sampled kernel and the continuous kernel restricted to the same cells. These are
$\mathcal R_n^{B,\mathrm{step},0}$ and $\mathcal R_n^{B,0}$ from \eqref{B:eq:criterion-brownian-cuts}, with the cells $I_n^{\rm ret}$ used by the finite restricted inverse. The
preceding truncation bound controls both the omitted kernel and the correction on the retained cells.

\begin{proposition}\label{B:prop:full-inverse-upper-bound}
Almost surely,
\begin{equation}
 Z_V:=\sup_{n\ge1}n\|\mathcal R_n-\mathcal R\|_{\mathfrak S_2}^2<\infty.
 \label{B:eq:full-inverse-error-supremum}
\end{equation}
Consequently $n\|\mathcal R_n-\mathcal R\|_{\mathfrak S_2}^2 =O_{\mathbb P}(1)$.
\end{proposition}
\begin{proof}
The isometry of the embedding by normalized cell indicators gives $\|\mathcal A_n^{-1}\|_{\mathfrak S_2} =\|\widetilde{\mathcal R}_n\|_{\mathfrak S_2}$. By
\eqref{B:eq:full-finite-inverse-norms} and Lemma~\ref{B:result:2-4},
\begin{equation}
 \|\mathcal R_n-\widetilde{\mathcal R}_n\|_{\mathfrak S_2}
 \le\frac{C_V}{n}\|\mathcal A_n^{-1}\|_{\mathfrak S_2}
 \le\frac{C_\omega}{n}.
 \label{B:eq:full-inverse-embedding-error}
\end{equation}
For the comparison on the selected interval, Corollary~\ref{B:cor:normalized-restricted-inverse} gives
\[
 \|\widetilde{\mathcal R}_n^0
          -\mathcal R_n^{B,\mathrm{step},0}\|_{\mathfrak S_2}^2
 \le C_\omega/n.
\]
The partition and cuts satisfy $m_n\asymp r_n$ and $n-M_n\asymp r_n$. Theorem~\ref{B:thm:brownian-inverse-approximation}, including its stated range of deterministic cuts, yields
\[
 \|\mathcal R_n^{B,\mathrm{step},0}
             -\mathcal R_n^{B,0}\|_{\mathfrak S_2}^2
 \le C_\omega/n.
\]
Proposition~\ref{B:prop:brownian-factor-truncation} and \eqref{B:eq:finite-truncation-orders} give
\[
 \|\mathcal R_n^{B,0}-\mathcal R\|_{\mathfrak S_2}^2
 \le C_{\omega,\eta}
     \left[d_n^{1+\nu_a-\eta}+(\varepsilon_n^+)^3
                                      +(\varepsilon_n^+)^2\right]
 =o_\omega(n^{-1}).
\]
These estimates and Theorem~\ref{B:thm:finite-inverse-truncation} apply to the exact identity
\begin{equation}
 \begin{split}
 \mathcal R_n-\mathcal R={}&
  (\mathcal R_n-\widetilde{\mathcal R}_n)
  +(\widetilde{\mathcal R}_n-\widetilde{\mathcal R}_n^0)\\
 &+(\widetilde{\mathcal R}_n^0-\mathcal R_n^{B,\mathrm{step},0})
  +(\mathcal R_n^{B,\mathrm{step},0}-\mathcal R_n^{B,0})
  +(\mathcal R_n^{B,0}-\mathcal R).
 \end{split}
 \label{B:eq:full-inverse-comparison}
\end{equation}
The squared norms of its five terms have respective orders $O_\omega(n^{-2})$, $o_\omega(n^{-1})$, $O_\omega(n^{-1})$, $O_\omega(n^{-1})$, and $o_\omega(n^{-1})$. Let $E_{n,l}$
denote the $l$th difference in \eqref{B:eq:full-inverse-comparison}. Expansion of its square gives
\[
 \begin{split}
 \|\mathcal R_n-\mathcal R\|_{\mathfrak S_2}^2
 &=\sum_{l=1}^5\|E_{n,l}\|_{\mathfrak S_2}^2
   +2\operatorname{Re}\sum_{l<h}
                  \langle E_{n,l},E_{n,h}\rangle_{\mathfrak S_2}\\
 &\le\sum_{l=1}^5\|E_{n,l}\|_{\mathfrak S_2}^2
       +\sum_{l<h}(\|E_{n,l}\|_{\mathfrak S_2}^2+
                                      \|E_{n,h}\|_{\mathfrak S_2}^2)\\
 &=5\sum_{l=1}^5\|E_{n,l}\|_{\mathfrak S_2}^2\le C_\omega/n
 \end{split}
\]
for every $n$ beyond the maximum of the five starting indices. Here Cauchy--Schwarz bounds the inner products between the differences, which are defined on the same probability
space.

Every smaller matrix has a positive finite inverse and $\mathcal R\in\mathfrak S_2$ by Proposition~\ref{B:result:2-5}. If the last estimate holds for $n\ge N(\omega)$, then
\[
 Z_V\le\max\left\{C_\omega,
       \max_{1\le n<N(\omega)}
                n(\|\mathcal R_n\|_{\mathfrak S_2}
                         +\|\mathcal R\|_{\mathfrak S_2})^2\right\}<\infty,
\]
with an empty maximum equal to zero. This proves \eqref{B:eq:full-inverse-error-supremum}. The supremum is measurable, since it is countable. Defining it arbitrarily on the null
exceptional set,
\[
 \sup_n\mathbb P\{n\|\mathcal R_n-\mathcal R\|_{\mathfrak S_2}^2>L\}
 \le\mathbb P\{Z_V>L\}\longrightarrow0.
\]
This proves tightness.
\end{proof}

\subsubsection{Proof of the theorem for general potentials}
\label{B:subsection:proof-general-potential-theorem}

\begin{proof}[Proof of Theorem~\ref{B:result:1-3}]
Theorem~\ref{B:thm:simultaneous-coefficient-construction} supplies the original coefficient law at every size and the information restriction \eqref{B:eq:1-32}.
Proposition~\ref{B:prop:full-inverse-upper-bound} gives the upper estimate. Its proof uses the scalar from Proposition~\ref{B:prop:coefficient-scalar-normalization}, the bounds in
Theorem~\ref{B:thm:normalized-coefficient-energy} and Proposition~\ref{B:prop:normalized-coefficient-maxima}, and the truncation estimate in
Theorem~\ref{B:thm:finite-inverse-truncation}.

Corollary~\ref{B:cor:simultaneous-inverse-information} applies the lower bound to these embedded inverses. Intersect its event with that of
Proposition~\ref{B:prop:full-inverse-upper-bound} and the event on which the inverse of the stochastic Bessel operator is defined. Then
\begin{equation}
 0<\liminf_{n\to\infty}
     n\|\mathcal R_n-\mathcal R\|_{\mathfrak S_2}^2
 \le\sup_{n\ge1}n\|\mathcal R_n-\mathcal R\|_{\mathfrak S_2}^2
 =Z_V<\infty.
 \label{B:eq:general-potential-two-sided-inverse}
\end{equation}
For example, if $b_\omega$ denotes the liminf, take $c_\omega^-=(1\wedge b_\omega)/2$ and $C_\omega^+=Z_V$. The definition of liminf gives the lower bound with $c_\omega^-$ for all
sufficiently large $n$. Writing $a_n(\omega)=n\|\mathcal R_n-\mathcal R\|_{\mathfrak S_2}^2$, one has $b_\omega=\sup_m\inf_{n\ge m}a_n(\omega)$, a measurable function. For the
least index $N_-$ from which $a_n\ge c_\omega^-$ holds for all $n\ge N_-$,
\[
 \{N_-\le m\}=\bigcap_{n\ge m}\{a_n\ge c_\omega^-\}
\]
on the chosen event of probability one. Extend $N_-$ by infinity on its null complement. This proves measurability and \eqref{B:eq:1-38}, with the normalization of
\eqref{B:eq:1-26} unchanged.

Fix an outcome in this event and a nonempty compact set $K\subset\rho(\mathcal H_{\beta,a})$. Use $B_K(\mathcal R)$, $D_K(\mathcal R)$, and $\epsilon_K(\mathcal R)$ from
\eqref{B:eq:general-resolvent-constants}. Since $\|\mathcal R_n-\mathcal R\|_{\rm op} \le\|\mathcal R_n-\mathcal R\|_{\mathfrak S_2}\to0$,
Lemma~\ref{B:lem:compact-inverse-resolvent-bounds} gives $n_0(\omega,K)$ such that $K\subset\rho(\mathcal A_n)$ for $n\ge n_0(\omega,K)$. To specify the constants, write
$G_n(z)=\mathcal G_{n,V,\beta,a}(z)$ and $G(z)=\mathcal G_{\beta,a}(z)$. The exact inverse identities are
\[
 \begin{split}
 G_n(z)-G(z)
   &=(I-z\mathcal R_n)^{-1}(\mathcal R_n-\mathcal R)
                                  (I-z\mathcal R)^{-1},\\
 \mathcal R_n-\mathcal R
   &=(I-z\mathcal R_n)(G_n(z)-G(z))(I-z\mathcal R).
 \end{split}
\]
The first pair of outer factors has product norm at most $2B_K(\mathcal R)^2$. Once $\|\mathcal R_n-\mathcal R\|_{\rm op}\le\epsilon_K(\mathcal R)\le1$, the second pair has product
norm at most $D_K(\mathcal R)$. Apply the Hilbert--Schmidt ideal inequality to each line, then square and use the upper and lower bounds for the error between inverses. This gives
\begin{align}
 \inf_{z\in K}
  \|\mathcal G_{n,V,\beta,a}(z)-\mathcal G_{\beta,a}(z)\|_{\mathfrak S_2}^2
 &\ge \frac{c_\omega^-}{D_K(\mathcal R)^2n},
            \label{B:eq:general-potential-resolvent-lower}\\
 \sup_{z\in K}
  \|\mathcal G_{n,V,\beta,a}(z)-\mathcal G_{\beta,a}(z)\|_{\mathfrak S_2}^2
 &\le \frac{4B_K(\mathcal R)^4C_\omega^+}{n}.
            \label{B:eq:general-potential-resolvent-upper}
\end{align}
Here the transforms of the inverses agree with the embedded resolvents by \eqref{B:eq:embedded-inverse-resolvent-factor} and \eqref{B:eq:resolvent-inverse-difference}. The lemma is
deterministic, so \eqref{B:eq:general-potential-resolvent-lower}--\eqref{B:eq:general-potential-resolvent-upper} hold for every such $K$ on the fixed event.

For a fixed nonempty compact $K\subset\mathbb C\setminus(0,\infty)$, all finite resolvents exist. The bound
\[
 M_\omega=\max\{\|\mathcal R\|_{\rm op},
                    \sup_{n\ge1}\|\mathcal R_n\|_{\rm op}\}<\infty
\]
follows from convergence and finiteness at the earlier sizes. Lemma~\ref{B:lem:positive-inverse-resolvent-bounds}, with $d_K(M_\omega)>0$ from
\eqref{B:eq:positive-resolvent-constants}, gives
\begin{equation}
 \sup_{n\ge1}n\sup_{z\in K}
  \|\mathcal G_{n,V,\beta,a}(z)-\mathcal G_{\beta,a}(z)\|_{\mathfrak S_2}^2
 \le d_K(M_\omega)^{-4}Z_V<\infty.
 \label{B:eq:general-potential-resolvent-supremum}
\end{equation}
This supremum is measurable: the resolvent difference is continuous in Hilbert--Schmidt norm on $K$, so its supremum can be taken over a fixed countable dense subset, followed by
the countable supremum in $n$. The tail of the finite random variable on the right tends to zero. This proves
\eqref{B:eq:intro-resolvent-rates}--\eqref{B:eq:intro-resolvent-supremum}.
\end{proof}

\subsubsection{Fixed compact spectral sets}
\label{B:subsection:fixed-compact-spectral-sets}

A fixed compact set may meet the random limiting spectrum, so a probabilistic statement must retain the event on which it is a subset of the resolvent set. The extended-valued
error below also records finite matrices whose resolvents are not yet defined throughout that set.

For a deterministic nonempty compact $K\subset\mathbb C$, retain $\mathcal E_K=\{K\subset\rho(\mathcal H_{\beta,a})\}$ and the extended valued variable $Y_{n,K}$ of
Proposition~\ref{B:prop:resolvent-tightness}, applied to the present coupling. Thus $Y_{n,K}$ is the normalized squared supremum of the resolvent error when the finite resolvent
exists throughout $K$, and is $+\infty$ otherwise.

\begin{corollary}\label{B:cor:general-potential-spectral-sets}
For every such $K$,
\begin{equation}
 \lim_{L\to\infty}\limsup_{n\to\infty}
       \mathbb P\bigl(\mathcal E_K\cap\{Y_{n,K}>L\}\bigr)=0.
 \label{B:eq:general-potential-spectral-set-tightness}
\end{equation}
\end{corollary}
\begin{proof}
The measurability of $\mathcal E_K$ and $Y_{n,K}$ is proved in Proposition~\ref{B:prop:resolvent-tightness}. On $\mathcal E_K$ let $N_K$ be the least integer such that
$\|\mathcal R_n-\mathcal R\|_{\rm op} \le\epsilon_K(\mathcal R)$ for every $n\ge N_K$. It is finite almost surely and measurable, since
\[
 \{N_K\le m\}\cap\mathcal E_K
 =\mathcal E_K\cap\bigcap_{n\ge m}
   \{\|\mathcal R_n-\mathcal R\|_{\rm op}
                              \le\epsilon_K(\mathcal R)\}.
\]
Set $C_K^*=4B_K(\mathcal R)^4Z_V$ on $\mathcal E_K$. Lemma~\ref{B:lem:compact-inverse-resolvent-bounds} implies $Y_{n,K}\le C_K^*$ there for $n\ge N_K$. Consequently,
\begin{equation}
 \mathbb P(\mathcal E_K\cap\{Y_{n,K}>L\})
 \le\mathbb P(\mathcal E_K\cap\{N_K>n\})
       +\mathbb P(\mathcal E_K\cap\{C_K^*>L\}).
 \label{B:eq:general-potential-spectral-set-probability}
\end{equation}
Set $N_K=\infty$ off $\mathcal E_K$ and $C_K^*=0$ there. The events $\mathcal E_K\cap\{N_K>n\}$ decrease to a null set, so continuity from above gives
\[
 \limsup_{n\to\infty}
 \mathbb P(\mathcal E_K\cap\{Y_{n,K}>L\})
 \le\mathbb P(\mathcal E_K\cap\{C_K^*>L\}).
\]
As $L\to\infty$, the events on the right decrease to $\mathcal E_K\cap\{C_K^*=\infty\}$, also null. This proves the claimed order of limits. A finite pole before $N_K$ gives
$Y_{n,K}=+\infty$ and is accounted for by the first term of \eqref{B:eq:general-potential-spectral-set-probability}. In particular it applies to $K=\{E\}$ for a fixed real $E$ on
the event $\{E\notin\sigma(\mathcal H_{\beta,a})\}$.
\end{proof}

\subsubsection{The two matrix families on one probability space}
\label{B:subsection:two-matrix-families}

The couplings for the beta-Laguerre ensemble and for general potentials may share the limiting Brownian path without sharing their auxiliary randomness. A product extension
combines the two constructions while preserving each finite row law.

\begin{corollary}\label{B:cor:two-matrix-families}
For fixed $V,\beta,a$ as above, the constructions in Theorems~\ref{B:result:1-1} and \ref{B:result:1-3} can use the same Brownian motion and limiting stochastic Bessel operator,
with the conclusions of both theorems holding almost surely.
\end{corollary}
\begin{proof}
Let $\lambda_V$ be the law of the complete auxiliary sequence used for the general potential, and let $\lambda_{\rm L}$ be the law of the auxiliary sequence used in the
beta-Laguerre construction. Section~\ref{B:subsection:simultaneous-embedded-inverses} puts both families on $\mathbb W\otimes\lambda_V\otimes\lambda_{\rm L}$. For their Borel
coefficient maps $\Phi_V$ and $\Phi_{\rm L}$ and any bounded Borel $F$,
\[
 \int F(W,\Phi_V(W,v))\,\mathbb W(\dd W)\lambda_V(\dd v)
                  \lambda_{\rm L}(\dd l)
 =\int F(W,\Phi_V(W,v))\,\mathbb W(\dd W)\lambda_V(\dd v).
\]
Integrating $v$ instead proves the corresponding identity for $\Phi_{\rm L}$. Thus each entire array has its prescribed joint distribution with $W$. Let $\Omega_V$ and
$\Omega_{\rm L}$ be the pullbacks to that space of the events in Theorems~\ref{B:result:1-3} and \ref{B:result:1-1}. Then
\[
 \mathbb P((\Omega_V\cap\Omega_{\rm L})^c)
 \le\mathbb P(\Omega_V^c)+\mathbb P(\Omega_{\rm L}^c)=0.
\]
On their intersection both limiting inverses are the same measurable function $\mathcal H_{\beta,a}(W)^{-1}$. The conclusion on compact spectral sets for the beta-Laguerre family
follows by Lemma~\ref{B:lem:compact-inverse-resolvent-bounds} and the probability estimate \eqref{B:eq:general-potential-spectral-set-probability}, using the corresponding supremum
of the inverse error.
\end{proof}

\subsection{A conditionally Gaussian random operator}
\label{sec:combined-s-conditional}
\label{S:sec:the-conditional-gaussian-operator}

The error from sampling the Brownian path has a leading term whose conditional distribution, given the sampled positions, is Gaussian in the Hilbert space of Hilbert--Schmidt
operators. For fixed $\beta>0$, Section~\ref{S:sec:the-error-from-sampling-the-brownian-path} constructs the linear interpolation terms in normal coordinates and the Brownian
bridge terms, and controls the remaining error in the weighted norm. Applying the derivative of the Green kernel introduces the scale $n^{-1/2}$. Throughout this subsection, $n$ is
large enough that $m_n<n$ and $\ell_n=m_n-R_\beta\ge1$.

\subsubsection{The expansion on the selected interval}
\label{S:sec:s10-retained-expansion}

We identify Hilbert--Schmidt operators with their square-integrable kernels and use this fixed real Hilbert space for the Gaussian expansion.

Recall the real Hilbert spaces
\[
 \mathcal H_{\mathbb R}=L^2([0,1);\mathbb R^2),\qquad
 \mathscr K=\mathfrak S_2(\mathcal H_{\mathbb R}).
\]
For operators $A,C\in\mathscr K$ with kernels $K_A,K_C$,
\begin{equation}\label{S:eq:s10-kernel-inner-product}
 \langle A,C\rangle_{\mathscr K}
 =\int_0^1\!\int_0^1
       \operatorname{tr}\bigl(K_A(s,t)^{\mathsf T}K_C(s,t)\bigr)
       \,\dd t\,\dd s.
\end{equation}
To verify this identification directly, choose a real orthonormal basis $(e_j)_{j\ge1}$ of $\mathcal H_{\mathbb R}$. If $A\in\mathscr K$, the finite kernels
\[
 K_{A,M}(s,t)=\sum_{j=1}^M(Ae_j)(s)e_j(t)^{\mathsf T}
\]
satisfy, by orthogonality in $t$,
\[
 \|K_{A,L}-K_{A,M}\|_{L^2((0,1)^2;\mathbb R^{2\times2})}^2
 =\sum_{j=M+1}^L\|Ae_j\|_{\mathcal H_{\mathbb R}}^2.
\]
They converge to an $L^2$ kernel $K_A$, with squared norm $\|A\|_{\mathfrak S_2}^2$. Integration against any finite linear combination of the $e_j$ gives the action of $A$ once $M$
includes its indices. For an arbitrary $f$, Cauchy--Schwarz gives
\[
 \left\|\int_0^1K_A(\,\cdot\,,t)f(t)\,\dd t\right\|_2
 \le\|K_A\|_{L^2}\|f\|_2,
\]
so density proves equality with $Af$ on the whole Hilbert space. Conversely, an $L^2$ matrix kernel defines a bounded operator by this inequality, and Parseval in $t$, followed by
Tonelli in $s$, gives $\sum_j\|Ae_j\|_2^2=\|K_A\|_{L^2}^2$. The correspondence is therefore onto and isometric. Polarization over the real scalars proves
\eqref{S:eq:s10-kernel-inner-product}. Rational matrix-valued step functions on rectangles with rational endpoints are a countable dense family of kernels, proving separability.

Multiplying the kernel by the indicator of a square is an orthogonal projection in this $L^2$ space. It corresponds to compressing the operator by the indicator projection on
$\mathcal H_{\mathbb R}$ and extending by zero. The reflection $f(t)\mapsto f(1-t)$ is unitary, sends $K(s,t)$ to $K(1-s,1-t)$, and preserves both the displayed inner product and
all Hilbert--Schmidt norms. These are the restriction and reflection conventions used below.

The covariance acts on $\mathscr K$. Use the rank-one convention $(A\otimes C)D=\langle D,C\rangle_{\mathscr K}A$ from Lemma~\ref{S:lem:conditional-gaussian-series}, now on this
Hilbert space. Its rank-one calculation gives $\|A\otimes C\|_{\mathfrak S_1(\mathscr K)}=\|A\|_{\mathscr K}\|C\|_{\mathscr K}$ and
$\operatorname{tr}_{\mathscr K}(A\otimes A)=\|A\|_{\mathscr K}^2$. We write $\operatorname{tr}Q$ for this covariance trace.

Condition on the sampled sigma-field $\mathscr S_n$ in \eqref{S:eq:sample-chain-sigma-field}, for one row $n$. Theorem~\ref{S:thm:joint-probability-space} gives conditionally
independent standard normals $\xi_{n,j,m,\nu}$ from the already attached bridges. The fields $\mathcal V_{n,j,m,\nu}$ and partial sums $\mathcal Z_n^{(M)}$ are those of
\eqref{S:eq:gaussian-mode-profiles}--\eqref{S:eq:gaussian-field-truncation}. The index $j$ ranges over $\mathcal J_n^{\mathrm{at}}$, including the additional intervals below $m_n$
and any zero profiles. The integer $m$ increases first; the finitely many pairs $(j,\nu)$ have the fixed order used there.

The full weighted field space is $\mathfrak H_n^+$, with norm $\|U\|_{n,+}^2=\mathcal E_n^+(U)$. After conditioning, it is a deterministic Hilbert space and
$\|\mathscr D_n\|_{\mathfrak H_n^+\to\mathscr K}\le n^{-1/2}$. Define
\begin{equation}\label{S:eq:gaussian-finite-operator-definitions}
 \begin{aligned}
 \Psi_{n,j,m,\nu}&=-\mathscr D_n\mathcal V_{n,j,m,\nu},&
 G_n^{(M)}&=-\mathscr D_n\mathcal Z_n^{(M)},\\
 G_n&=-\mathscr D_n\mathcal Z_n,&
 B_n&=-\mathscr D_n\mathcal B_n,\\
 Q_n^{(M)}&=\sum_{j\in\mathcal J_n^{\mathrm{at}}}
             \sum_{m=1}^M\sum_{\nu=1}^2
              \Psi_{n,j,m,\nu}\otimes\Psi_{n,j,m,\nu}.&&
 \end{aligned}
\end{equation}
The coefficients $\Psi_{n,j,m,\nu}$ and $B_n$ are determined by the sampled positions. Every $\Psi_{n,j,m,\nu}$ is a self-adjoint Hilbert--Schmidt operator, since its derivative
kernel has transposition symmetry. The minus signs agree with Proposition~\ref{S:prop:continuous-bridge-linearization}.

The covariance series is
\begin{equation}\label{S:eq:bridge-covariance-series}
 Q_n=\sum_{j\in\mathcal J_n^{\mathrm{at}}}
          \sum_{m\ge1}\sum_{\nu=1}^2
                   \Psi_{n,j,m,\nu}\otimes\Psi_{n,j,m,\nu}.
\end{equation}
Its convergence in trace norm is established below. Let $R_n^{\mathrm{br}}$ be the exact remainder of Proposition~\ref{S:prop:continuous-bridge-linearization}; every kernel is
extended by zero outside the square of the restricted kernel.

\begin{theorem}\label{S:thm:15.1}
On the coupling of Theorem~\ref{S:thm:joint-probability-space},
\[
 \varepsilon_{\beta,n}^{\mathrm{samp},T}
      =G_n+B_n+R_n^{\mathrm{br}},\qquad
 G_n\mid\mathscr S_n\sim\mathcal N_{\mathscr K}(0,Q_n).
\]
Almost surely, $\sqrt n\|R_n^{\mathrm{br}}\|_{\mathscr K} =O_{\beta,\omega}((\log n)^{-250})$.
\end{theorem}

The signed identity comes from the field expansion. To identify its conditional law, we first justify the Gaussian series and its covariance. After computing the conditional second
moment we finish the proof of Theorem~\ref{S:thm:15.1}.

\subsubsection{Gaussian series in Hilbert--Schmidt space}
\label{S:sec:the-conditional-gaussian-field}

Finite sums of bridge modes define the Gaussian approximation directly. Summability of their squared Hilbert--Schmidt norms gives convergence of the infinite series and a
trace-class covariance.

\begin{theorem}\label{S:thm:gaussian-finite-mode}
Conditionally on $\mathscr S_n$, $G_n^{(M)}$ converges to $G_n$ in $L^2(\mathscr K)$ and almost surely. Its limit has law $\mathcal N_{\mathscr K}(0,Q_n)$, where $Q_n$ is the
positive trace-class operator in \eqref{S:eq:bridge-covariance-series}. The convergence and conditional distributions hold simultaneously for all selected rows.
\end{theorem}
\begin{proof}
Fix a conditional law for which the sampled coefficients and the bridge construction are defined. The second sine in \eqref{S:eq:gaussian-mode-profiles} vanishes at a sampling
vertex. Consequently
\[
 |q_{n,k,j,m}(s)|\le\frac{\sqrt2}{\pi m},\qquad
 \mathcal M_{n,k}^+(\mathsf A_\nu)
 \le\|\mathsf A_\nu\|_{\mathrm{op}}^2\mathfrak G_{n,k}^+.
\]
Integrating in the normalized interval coordinate and using the bound on $\mathscr D_n$ gives
\begin{equation}\label{S:eq:gaussian-operator-single-mode}
 \|\Psi_{n,j,m,\nu}\|_{\mathscr K}^2
 \le\frac1n\|\mathcal V_{n,j,m,\nu}\|_{n,+}^2
 \le\frac{2\|\mathsf A_\nu\|_{\mathrm{op}}^2}{\pi^2}
                   \frac{e_{n,j}}{nm^2}.
\end{equation}
Here $e_{n,j}$ and $e_n^{\mathrm{KL}}=\sum_j e_{n,j}$ are the sampled quantities in \eqref{S:eq:gaussian-sine-coefficient-bound}. Since $\|\mathsf A_1\|_{\mathrm{op}}=1$ and
$\|\mathsf A_2\|_{\mathrm{op}}=1/2$, summation gives
\begin{equation}\label{S:eq:s10-summed-mode-bounds}
 \sum_{j,m,\nu}\|\Psi_{n,j,m,\nu}\|_{\mathscr K}^2
 \le\frac{5e_n^{\mathrm{KL}}}{12n},\qquad
 \sum_{j,m>M,\nu}\|\Psi_{n,j,m,\nu}\|_{\mathscr K}^2
 \le\frac{5e_n^{\mathrm{KL}}}{2\pi^2nM}.
\end{equation}
The first inequality uses $\sum_{m\ge1}m^{-2}=\pi^2/6$; the second uses $\sum_{m>M}m^{-2}\le\int_M^\infty x^{-2}\dd x=M^{-1}$. These sums are finite for each selected row, before
any uniform bound in $n$ is used.

For $L>M$, independence and centering imply exactly
\begin{equation}\label{S:eq:s10-finite-series-isometry}
 \mathbb E[\|G_n^{(L)}-G_n^{(M)}\|_{\mathscr K}^2
                           \mid\mathscr S_n]
 =\sum_j\sum_{M<m\le L}\sum_{\nu=1}^2
                        \|\Psi_{n,j,m,\nu}\|_{\mathscr K}^2.
\end{equation}
The identity uses $\mathbb E[\xi_{n,j,m,\nu}\xi_{n,j',m',\nu'}\mid\mathscr S_n] =\mathbf1_{\{(j,m,\nu)=(j',m',\nu')\}}$. We prove convergence under the fixed conditional law. Write
$t_M=\sum_{j,m>M,\nu}\|\Psi_{n,j,m,\nu}\|_{\mathscr K}^2$, so $t_M\downarrow0$. The isometry just proved makes $(G_n^{(M)})_M$ Cauchy in $L^2(\mathscr K)$.

For a finite $K>M$, put $S_L=G_n^{(L)}-G_n^{(M)}$, $M\le L\le K$, and $U=\max_{M\le L\le K}\|S_L\|$. These are Hilbert-space martingale partial sums. If $\tau$ is the first $L$
with $\|S_L\|>t$, conditional Jensen gives
$\|S_L\|\le\mathbb E[\|S_K\|\mid\mathcal F_L^{\mathrm{mode}}]$.
Multiplication by $\mathbf1_{\{\tau=L\}}$ and summation in $L$ yield
\[
 t\,\mathbb P(U>t)\le\mathbb E[\|S_K\|\mathbf1_{\{U>t\}}].
\]
Integrating in $t$ and applying Cauchy--Schwarz gives
\[
 \mathbb E U^2\le2\mathbb E[\|S_K\|U]
 \le2(\mathbb E\|S_K\|^2)^{1/2}(\mathbb E U^2)^{1/2}.
\]
All terms are finite for this finite maximum. Thus $\mathbb E U^2\le4t_M$. Letting $K\to\infty$ by monotone convergence proves
\[
 \mathbb E\left[\sup_{L\ge M}\|G_n^{(L)}-G_n^{(M)}\|_{\mathscr K}^2
       \,\middle|\,\mathscr S_n\right]\le4t_M.
\]
For the fixed sampled data choose increasing $M_r$ with $t_{M_r}\le2^{-4r}$. Markov's inequality bounds the probability that the last supremum in norm exceeds $2^{-r}$ by
$4\,2^{-2r}$. Borel--Cantelli makes the entire sequence Cauchy almost surely, since any two indices beyond $M_r$ differ by at most $2^{1-r}$ on the resulting event. Completeness of
$\mathscr K$ gives its almost sure limit. This agrees with its $L^2$ limit by uniqueness of convergence in probability. The choices $M_r$ were made only after fixing the
conditional law; the convergence event itself is the countable Cauchy event and is Borel on the original space.

Proposition~\ref{S:prop:weighted-field-covariance} identifies $\mathcal Z_n$ as the conditional $L^2$ and almost sure limit of $\mathcal Z_n^{(M)}$ on these same attached bridges.
Boundedness of $\mathscr D_n$ identifies the operator limit with the $G_n$ in \eqref{S:eq:gaussian-finite-operator-definitions}. Thus
\begin{equation}\label{S:eq:bridge-kl-expansion}
 G_n=\sum_{j\in\mathcal J_n^{\mathrm{at}}}
           \sum_{m\ge1}\sum_{\nu=1}^2
           \xi_{n,j,m,\nu}\Psi_{n,j,m,\nu}.
\end{equation}
The series is defined by the indicated order of its partial sums.

The rank-one trace norm and \eqref{S:eq:s10-summed-mode-bounds} justify absolute convergence of \eqref{S:eq:bridge-covariance-series} in $\mathfrak S_1(\mathscr K)$.

For $A\in\mathscr K$,
\[
 \langle(Q_n^{(M+1)}-Q_n^{(M)})A,A\rangle_{\mathscr K}
 =\sum_{j,\nu}\langle\Psi_{n,j,M+1,\nu},A\rangle_{\mathscr K}^2\ge0.
\]
The range of $Q_n^{(M)}$ is contained in the span of its $2M|\mathcal J_n^{\mathrm{at}}|$ defining vectors. Hence

\begin{equation}\label{S:eq:gaussian-positive-mode-order}
 0\le Q_n^{(M)}\le Q_n^{(M+1)}\le Q_n,\qquad
 \operatorname{rank}Q_n^{(M)}\le2M|\mathcal J_n^{\mathrm{at}}|.
\end{equation}
The rank bound is for $Q_n^{(M)}$ acting on $\mathscr K$; its range is spanned by the finitely many vectors
$\Psi_{n,j,m,\nu}$.

For $A_1,\ldots,A_r\in\mathscr K$ and real $u_1,\ldots,u_r$, independence of the finitely many normal variables gives
\begin{equation}\label{S:eq:gaussian-finite-characteristic}
 \mathbb E\left[\exp\left(\ii\sum_{a=1}^r
             u_a\langle G_n^{(M)},A_a\rangle_{\mathscr K}\right)
                  \,\middle|\,\mathscr S_n\right]
 =\exp\left(-\frac12\sum_{a,b=1}^r u_au_b
              \langle Q_n^{(M)}A_a,A_b\rangle_{\mathscr K}\right).
\end{equation}
Let $V=\sum_{a=1}^r u_aA_a$. Conditional $L^2$ convergence and $|\mathrm e^{\mathrm ix}-\mathrm e^{\mathrm iy}|\le|x-y|$ imply
\[
 \left|\mathbb E\!\left[\mathrm e^{\mathrm i\langle G_n^{(M)},V\rangle}
           -\mathrm e^{\mathrm i\langle G_n,V\rangle}\,\middle|\,\mathscr S_n\right]\right|
 \le\|V\|\,t_M^{1/2}\longrightarrow0.
\]
On the other side of the finite characteristic-function identity,
$|\langle(Q_n^{(M)}-Q_n)V,V\rangle|\le\|Q_n^{(M)}-Q_n\|_{\mathfrak S_1}\|V\|^2\to0$.
Consequently
\begin{equation}\label{S:eq:gaussian-limit-characteristic}
 \mathbb E\left[\exp\left(\ii\sum_{a=1}^r
                u_a\langle G_n,A_a\rangle_{\mathscr K}\right)
                  \,\middle|\,\mathscr S_n\right]
 =\exp\left(-\frac12\sum_{a,b=1}^r u_au_b
                       \langle Q_nA_a,A_b\rangle_{\mathscr K}\right).
\end{equation}
To choose one version, first take all finite tuples from a fixed countable dense subset of $\mathscr K$ and rational coefficients. There are only countably many such identities.
Intersect their sets of sampled data of full measure with the set on which the conditional covariance has finite trace. On any fixed sampled data in this intersection, the left
side is continuous in $V$, because its difference at $V,W$ is at most $(\operatorname{tr}Q_n)^{1/2}\|V-W\|$. The right side is continuous by boundedness of $Q_n$. Density therefore
extends the identity to all real coefficients and all vectors on that same set. Finally, coordinates along a countable orthonormal basis generate the Borel sigma-field of the
separable Hilbert space $\mathscr K$. Their finite-dimensional characteristic functions determine a probability measure, so the identity identifies the conditional law as
$\mathcal N_{\mathscr K}(0,Q_n)$.

Continuity of the conditional $L^2$ norm in \eqref{S:eq:s10-finite-series-isometry}, and positivity of the covariance tail, prove
\begin{equation}\label{S:eq:gaussian-operator-mode-tail}
 \begin{aligned}
 t_{n,M}
 &:=\mathbb E[\|G_n-G_n^{(M)}\|_{\mathscr K}^2\mid\mathscr S_n]\\
 &=\|Q_n-Q_n^{(M)}\|_{\mathfrak S_1(\mathscr K)}
   =\sum_{j,m>M,\nu}\|\Psi_{n,j,m,\nu}\|_{\mathscr K}^2\\
 &\le\frac{5e_n^{\mathrm{KL}}}{2\pi^2nM}.
 \end{aligned}
\end{equation}
\manualstatementlabel{S:eq:covariance-trace-norm-passage}{\theequation} For each fixed $M$, the tail is a measurable limit of sums using only coordinates with $m>M$. It is
therefore independent, under this conditional law, of the coordinates with $m\le M$ and of $G_n^{(M)}$.

The Cauchy criterion in $\mathscr K$ makes the convergence event Borel. It has conditional probability one for each $n$. Integration and countable intersection give simultaneous
convergence, with the conditional Gaussian distribution relative to each
row's sigma-field $\mathscr S_n$. \end{proof}

\subsubsection{Covariance and measurability}
\label{S:sec:s10-covariance-and-measurability}
\phantomsection\label{S:prop:gaussian-covariance-measurability}

To localize the conditional estimate, we prove that the covariance norm is measurable with respect to the sampled data. Each $\mathcal V_{n,j,m,\nu}$ is a measurable function of
those data and the deterministic profiles. The isometry $\mathcal J_n$ maps these fields into the fixed space $\mathcal H_{0,n}$. The proof in
Subsection~\ref{S:sec:s5-measurability-covariance} then shows that $\Psi_{n,j,m,\nu}=-\mathscr D_n\mathcal V_{n,j,m,\nu}$ is strongly $\mathscr S_n$-measurable. For
$v,w\in\mathscr K$,
\[
 v\otimes v-w\otimes w=(v-w)\otimes v+w\otimes(v-w)
\]
and the rank-one trace norm imply
\[
 \|v\otimes v-w\otimes w\|_{\mathfrak S_1(\mathscr K)}
 \le(\|v\|_{\mathscr K}+\|w\|_{\mathscr K})\|v-w\|_{\mathscr K}.
\]
Thus $Q_n^{(M)}$ is strongly measurable in the separable trace-class space. The event
\[
 \left\{\sum_{j,m,\nu}\|\Psi_{n,j,m,\nu}\|_{\mathscr K}^2<\infty\right\}
 \in\mathscr S_n
\]
has probability one. Define $Q_n$ there by the limit in trace norm and put it equal to zero on its complement. This gives a strongly $\mathscr S_n$-measurable version. The
inequalities
\[
 |\operatorname{tr}(Q-R)|\le\|Q-R\|_{\mathfrak S_1},\qquad
 \big|\|Q\|_{\mathrm{op}}-\|R\|_{\mathrm{op}}\big|
 \le\|Q-R\|_{\mathrm{op}}\le\|Q-R\|_{\mathfrak S_1}
\]
show that $\operatorname{tr}Q_n$, $\|Q_n\|_{\mathrm{op}}$, and $t_{n,M}$ are $\mathscr S_n$-measurable. The same construction of measurable linear maps applies to
$B_n=-\mathscr D_n\mathcal B_n$. The almost sure limit of $G_n^{(M)}$, set equal to zero off its convergence event, is strongly measurable on the common probability space; its
remaining randomness comes from the attached Gaussian coordinates.

Under the regular conditional law the exact moments are
\begin{equation}\label{S:eq:gaussian-conditional-covariance-identity}
 \begin{aligned}
 \mathbb E[G_n\mid\mathscr S_n]&=0,\\
 \mathbb E[\langle G_n,A\rangle_{\mathscr K}
              \langle G_n,C\rangle_{\mathscr K}\mid\mathscr S_n]
     &=\langle Q_nA,C\rangle_{\mathscr K},\\
 \mathbb E[\|G_n\|_{\mathscr K}^2\mid\mathscr S_n]
     &=\operatorname{tr}Q_n.
 \end{aligned}
\end{equation}
These are the moment identities in Lemma~\ref{S:lem:conditional-gaussian-series} for the ordered coefficients $\Psi_{n,j,m,\nu}$ and conditioning sigma-field $\mathscr S_n$. The
square summability required there is \eqref{S:eq:s10-summed-mode-bounds}; the covariance is the trace norm sum \eqref{S:eq:bridge-covariance-series}.

These formulas are initially integrals under a regular conditional probability measure. They do not by themselves assert that $\mathbb E\|G_n\|^2$ is finite. For an integer
$L\ge1$, put $\Lambda_{n,L}=\{\operatorname{tr}Q_n\le L\}\in\mathscr S_n$. Tonelli's theorem and the nonnegative third identity give
\[
 \mathbb E[\mathbf1_{\Lambda_{n,L}}\|G_n\|_{\mathscr K}^2]
 =\mathbb E[\mathbf1_{\Lambda_{n,L}}\operatorname{tr}Q_n]\le L.
\]
Cauchy--Schwarz makes $\mathbf1_{\Lambda_{n,L}}G_n$ Bochner integrable in $\mathscr K$. The rank-one identity $\|G_n\otimes G_n\|_{\mathfrak S_1}=\|G_n\|^2$ makes
$\mathbf1_{\Lambda_{n,L}}(G_n\otimes G_n)$ Bochner integrable in $\mathfrak S_1(\mathscr K)$. These random elements are strongly measurable by the continuous rank-one map proved
above.

Fix an orthonormal basis $(E_j)$ of $\mathscr K$ and a set $D\in\mathscr S_n$. Integration of the scalar conditional moments gives, for every $i,j$,
\[
 \begin{aligned}
 \mathbb E[\mathbf1_{D\cap\Lambda_{n,L}}\langle G_n,E_i\rangle]&=0,\\
 \mathbb E[\mathbf1_{D\cap\Lambda_{n,L}}
       \langle G_n,E_i\rangle\langle G_n,E_j\rangle]
 &=\mathbb E[\mathbf1_{D\cap\Lambda_{n,L}}\langle Q_nE_i,E_j\rangle].
 \end{aligned}
\]
Coordinates along $(E_j)$ separate vectors. Matrix entries along this basis separate bounded operators: if every entry is zero, the operator vanishes on the dense finite span and
hence everywhere. The defining property of Bochner conditional expectation, tested against these countably many continuous functionals, consequently gives
\[
 \begin{aligned}
 \mathbb E[\mathbf1_{\Lambda_{n,L}}G_n\mid\mathscr S_n]&=0,\\
 \mathbb E[\mathbf1_{\Lambda_{n,L}}(G_n\otimes G_n)\mid\mathscr S_n]
   &=\mathbf1_{\Lambda_{n,L}}Q_n.
 \end{aligned}
\]
The localized versions agree on $\Lambda_{n,L}\cap\Lambda_{n,L'}$ by uniqueness of conditional expectation. Since $\operatorname{tr}Q_n<\infty$ almost surely, their increasing
union has probability one. Thus the moments in \eqref{S:eq:gaussian-conditional-covariance-identity} have the stated localized meaning without an unconditional trace moment.

Let $\widetilde Q_n^+$ and $\widetilde Q_n^{+,(M)}$ be the field covariances of Proposition~\ref{S:prop:weighted-field-covariance}. With the adjoint taken in the full weighted
inner product,
\begin{equation}\label{S:eq:gaussian-covariance-pushforward}
 Q_n=\mathscr D_n\widetilde Q_n^+\mathscr D_n^*,\qquad
 Q_n^{(M)}=\mathscr D_n\widetilde Q_n^{+,(M)}\mathscr D_n^*.
\end{equation}
Indeed, for $v\in\mathfrak H_n^+$ and $A\in\mathscr K$,
\[
 \mathscr D_n(v\otimes_+v)\mathscr D_n^*A
 =\langle\mathscr D_n^*A,v\rangle_{n,+}\mathscr D_nv
 =\langle A,\mathscr D_nv\rangle_{\mathscr K}\mathscr D_nv.
\]
Summing proves the finite identity; the two minus signs cancel. The trace ideal estimate, Corollary~6.14 of \cite{teschl2009mathematical}, applied on the Hilbert direct sum of the
domain and target spaces, gives
\[
 \|DTD^*\|_{\mathfrak S_1}
 \le\|D\|^2\|T\|_{\mathfrak S_1}
\]
for bounded $D$ and trace-class $T$ on the indicated Hilbert spaces.

Apply this with $D=\mathscr D_n$ to the field covariance tails. Their trace norms tend to zero by \eqref{S:eq:weighted-field-mode-tail}, so the full identity follows. In
particular,
\begin{equation}\label{S:eq:s10-covariance-operator-comparison}
 \|Q_n\|_{\mathrm{op}}
 \le\frac1n\|\widetilde Q_n^+\|_{\mathrm{op}},\qquad
 \|Q_n-Q_n^{(M)}\|_{\mathrm{op}}\le t_{n,M}.
\end{equation}
The second inequality follows from positivity of the covariance tail and its trace.

\subsubsection{Conditional second moments}
\label{S:sec:s10-conditional-second-moments}

The leading conditional second moment has two contributions. The conditionally centered Gaussian random operator $G_n$ contributes its covariance trace, denoted by
$T_n^{\mathrm G}=\operatorname{tr}Q_n$ to distinguish it from the truncation time $T_n$. The term $B_n$, obtained from linear interpolation, is the conditional mean of the
linearized error and contributes its squared norm to the conditional second moment. The exact remainder is controlled separately by its norm.

\begin{corollary}\label{S:cor:gaussian-conditional-second-moment}
Under the conditional law,
\begin{equation}\label{S:eq:gaussian-trace-bias-mean}
 \mathbb E[\|G_n+B_n\|_{\mathscr K}^2\mid\mathscr S_n]
 =T_n^{\mathrm G}+\|B_n\|_{\mathscr K}^2
 \le\frac1n\left(
    \mathbb E[\mathcal E_n^+(\mathcal Z_n)\mid\mathscr S_n]
                    +\mathcal E_n^+(\mathcal B_n)\right).
\end{equation}
The same identity and bound hold for the corresponding finite sums.
\end{corollary}
\begin{proof}
Fix the sampled data. Expanding the square gives the three terms $\mathbb E\|G_n\|^2$, $2\langle\mathbb EG_n,B_n\rangle$, and $\|B_n\|^2$. The first two are $\operatorname{tr}Q_n$
and zero by \eqref{S:eq:gaussian-conditional-covariance-identity}. The derivative bound gives
\[
 \mathbb E\|G_n\|^2\le n^{-1}\mathbb E\mathcal E_n^+(\mathcal Z_n),
 \qquad \|B_n\|^2\le n^{-1}\mathcal E_n^+(\mathcal B_n).
\]
All expectations in this calculation are under that conditional law. For a finite cutoff, centering and the covariance identity give
\[
 \mathbb E[\|G_n^{(M)}+B_n\|_{\mathscr K}^2\mid\mathscr S_n]
 =\operatorname{tr}Q_n^{(M)}+\|B_n\|_{\mathscr K}^2
 \le\frac1n\left(\mathbb E[\mathcal E_n^+(\mathcal Z_n^{(M)})\mid\mathscr S_n]
                   +\mathcal E_n^+(\mathcal B_n)\right).
\]
Moreover,
\[
 \mathbb E[\mathcal E_n^+(\mathcal Z_n^{(M)})\mid\mathscr S_n]
 =\sum_{j,m\le M,\nu}\|\mathcal V_{n,j,m,\nu}\|_{n,+}^2
 \le\operatorname{tr}\widetilde Q_n^+.
\]
Thus the upper bound is uniform in $M$, with the same term from linear interpolation.

More generally, for any $\mathscr S_n$-measurable $A_n\in\mathscr K$, expansion under the fixed conditional law gives
\[
 \mathbb E[\|G_n+B_n-A_n\|_{\mathscr K}^2\mid\mathscr S_n]
 =\operatorname{tr}Q_n+\|B_n-A_n\|_{\mathscr K}^2.
\]
The conditional minimum for this linear Gaussian term is attained at $A_n=B_n$ and is $\operatorname{tr}Q_n$. The actual sampling error is not corrected in this way: its leading
conditional second moment includes $\|B_n\|^2$. This identity concerns the linear term, not the exact inverse including its nonlinear remainder, and does not establish optimality
over couplings. The sampled quantities are fixed and finite under the conditional distribution, which is the distribution used in this squared-norm calculation.
\end{proof}

We shall also use convergence of the squared norms in conditional $L^1$. Put $X=G_n+B_n$ and $Y=G_n^{(M)}+B_n$. The identity $\|X\|^2-\|Y\|^2=\langle X-Y,X+Y\rangle$ gives
\[
 \mathbb E[|\|X\|^2-\|Y\|^2|\mid\mathscr S_n]
 \le\mathbb E[\|X-Y\|^2\mid\mathscr S_n]^{1/2}
      \mathbb E[(\|X\|+\|Y\|)^2\mid\mathscr S_n]^{1/2}.
\]
The first factor is $t_{n,M}^{1/2}$. Conditional Minkowski and $\operatorname{tr}Q_n^{(M)}\le T_n^{\mathrm G}$ bound the second by $2(T_n^{\mathrm G}+\|B_n\|^2)^{1/2}$. Thus
\begin{equation}\label{S:eq:gaussian-mode-quadratic-stability}
 \mathbb E\left[ \left|\|G_n+B_n\|_{\mathscr K}^2 -\|G_n^{(M)}+B_n\|_{\mathscr K}^2\right| \,\middle|\,\mathscr S_n\right] \le 2\left((T_n^{\mathrm G}+\|B_n\|_{\mathscr K}^2)t_{n,M}\right)^{1/2}.
\end{equation}
By \eqref{S:eq:weighted-bridge-mean-bound} and \eqref{S:eq:weighted-sine-coefficient-sum}, the right sides of \eqref{S:eq:gaussian-trace-bias-mean} and
\eqref{S:eq:gaussian-mode-quadratic-stability} are at most $C_{\beta,\omega}/n$ and $C_{\beta,\omega}/(n\sqrt M)$, respectively, for all sufficiently large $n$ and every $M\ge1$.
The left sides are determined by the sampled positions.

The conditional expectation of the signed difference has the following formula. With $H=G_n-G_n^{(M)}$, expand $\|Y+H\|^2-\|Y\|^2$. The tail $H$ is centered and conditionally
independent of $G_n^{(M)}$, while $B_n$ is fixed under conditioning. Hence
\begin{equation}\label{S:eq:s10-signed-square-tail}
 \mathbb E[\|G_n+B_n\|_{\mathscr K}^2
          -\|G_n^{(M)}+B_n\|_{\mathscr K}^2\mid\mathscr S_n]
 =t_{n,M}.
\end{equation}
The mixed term has conditional mean zero.

\begin{proof}[Proof of Theorem~\ref{S:thm:15.1}]
The signed operator identity is precisely \eqref{S:eq:bridge-operator-linearization}, with the definitions \eqref{S:eq:gaussian-finite-operator-definitions}. The conditional law is
Theorem~\ref{S:thm:gaussian-finite-mode}. The estimates are
\begin{equation}\label{S:eq:gaussian-expansion-quantitative}
 \begin{aligned}
 \operatorname{tr}Q_n+\|B_n\|_{\mathscr K}^2
      &\le C_{\beta,\omega}/n,\\
 \sqrt n\|G_n+B_n\|_{\mathscr K}
      &\le C_{\beta,\omega},\\
 \sqrt n\|R_n^{\mathrm{br}}\|_{\mathscr K}
      &\le C_{\beta,\omega}(\log n)^{-250}.
 \end{aligned}
\end{equation}
The first follows from \eqref{S:eq:gaussian-trace-bias-mean} and the conditional mean bound \eqref{S:eq:weighted-bridge-mean-bound}. For the second we use the pathwise estimate
\eqref{S:eq:bridge-weighted-norm-bound}:
\[
 \sqrt n\|G_n+B_n\|_{\mathscr K}
 \le\|\mathcal Z_n+\mathcal B_n\|_{n,+}
 \le\mathcal E_n^+(\mathcal Z_n)^{1/2}
      +\mathcal E_n^+(\mathcal B_n)^{1/2}
 \le C_{\beta,\omega}.
\]
The last line of \eqref{S:eq:gaussian-expansion-quantitative} is \eqref{S:eq:continuous-bridge-remainder-rate}. Intersecting the events from
Sections~\ref{S:sec:the-error-from-sampling-the-brownian-path} and \ref{S:sec:the-conditional-gaussian-operator} gives the three bounds simultaneously.
\end{proof}

Polarization and \eqref{S:eq:gaussian-expansion-quantitative} give
\begin{equation}\label{S:eq:bridge-squared-norm-error}
 \begin{aligned}
 &n\left|\|\varepsilon_{\beta,n}^{\mathrm{samp},T}\|_{\mathscr K}^2
                 -\|G_n+B_n\|_{\mathscr K}^2\right|\\
 &\le2(\sqrt n\|G_n+B_n\|_{\mathscr K})
                (\sqrt n\|R_n^{\mathrm{br}}\|_{\mathscr K})
       +n\|R_n^{\mathrm{br}}\|_{\mathscr K}^2\\
 &\le C_{\beta,\omega}(\log n)^{-250}\longrightarrow0.
 \end{aligned}
\end{equation}
The estimate uses only the pathwise norm of $R_n^{\mathrm{br}}$.

\subsubsection{Truncation of the Gaussian series}
\label{S:sec:s10-series-truncation}

A cutoff varying with $n$ requires an estimate uniform in the row and mode indices. We obtain this joint estimate from the trace of the conditional covariance of the remainder,
without an unconditional moment assumption on the random bound for the Gaussian series.

First let $M(n)\ge1$ be deterministic integers satisfying $\sum_n M(n)^{-1}<\infty$; for example, $M(n)=n^2$. We write the order as $M(n)$ to distinguish it from the oscillation
bound $M_n$ used near the singular boundary.

\begin{corollary}\label{S:cor:gaussian-diagonal-truncation}
Almost surely,
\begin{equation}\label{S:eq:gaussian-diagonal-mode-limit}
 \sqrt n\|G_n-G_n^{(M(n))}\|_{\mathscr K}\longrightarrow0,
 \qquad
 n\|Q_n-Q_n^{(M(n))}\|_{\mathfrak S_1(\mathscr K)}\longrightarrow0.
\end{equation}
\end{corollary}
\begin{proof}
For each integer $L\ge1$ let $A_{n,L}=\{e_n^{\mathrm{KL}}\le L\}\in\mathscr S_n$. For $\delta>0$, conditional Markov and \eqref{S:eq:gaussian-operator-mode-tail} give
\[
 \mathbf1_{A_{n,L}}
 \mathbb P(\sqrt n\|G_n-G_n^{(M(n))}\|_{\mathscr K}>\delta
                       \mid\mathscr S_n)
 \le\frac{5L}{2\pi^2\delta^2M(n)}.
\]
Its sum in $n$ is finite. Lemma~\ref{S:lem:conditional-localized-borel-cantelli}, with one index in each row, applies for each integer $L$ and each positive rational $\delta$.
Intersect these countably many events. By \eqref{S:eq:weighted-sine-coefficient-sum}, almost every outcome has an integer $L$ such that $e_n^{\mathrm{KL}}\le L$ for every
sufficiently large $n$. On that outcome, for each positive rational $\delta$ there is an index $n_\delta(\omega)$ beyond which the error is at most $\delta$, proving the first
limit. The second follows pathwise from
\[
 n\|Q_n-Q_n^{(M(n))}\|_{\mathfrak S_1}
 \le\frac{5e_n^{\mathrm{KL}}}{2\pi^2M(n)}\longrightarrow0,
\]
since summability of $M(n)^{-1}$ implies $M(n)\to\infty$.
\end{proof}

An exponential estimate gives a smaller cutoff while preserving the logarithmic remainder rate. For fixed $n,M$, put $H_{n,M}=G_n-G_n^{(M)}$ and $t=t_{n,M}$. If $t=0$, the tail
vanishes almost surely under the conditional law. Suppose $t>0$. The tail is centered Gaussian with covariance $C=Q_n-Q_n^{(M)}$, trace $t$ and $\|C\|\le t$. Applying
\eqref{S:eq:covariance-gaussian-mgf} with zero shift and $\theta=(4t)^{-1}$ gives
\[
 \log\mathbb E\left[\exp\left(\frac{\|H_{n,M}\|_{\mathscr K}^2}{4t}\right)
                         \,\middle|\,\mathscr S_n\right]
 \le\frac14+
 \frac{\operatorname{tr}(C^2)}{16t^2(1-\|C\|/(2t))}
 \le\frac38\le\frac12.
\]
We used $\operatorname{tr}(C^2)\le t^2$. The infinite-dimensional exponential estimate has already been proved in Lemma~\ref{S:lem:covariance-gaussian-concentration}. Hence
\begin{equation}\label{S:eq:s10-gaussian-tail-exponential}
 \mathbb E\left[\exp\left(\frac{\|H_{n,M}\|_{\mathscr K}^2}
                                      {4t_{n,M}}\right)
                  \,\middle|\,\mathscr S_n\right]\le\me^{1/2},
 \qquad t_{n,M}>0.
\end{equation}
This bound depends only on the measurable scalar $t_{n,M}$. Markov's inequality therefore gives, for every $x\ge0$,
\begin{equation}\label{S:eq:s10-gaussian-tail-probability}
 \mathbb P\left(\|H_{n,M}\|_{\mathscr K}^2>
                 4t_{n,M}(1+x)\,\middle|\,\mathscr S_n\right)
 \le\me^{-x}.
\end{equation}
If $t_{n,M}=0$, the tail is zero almost surely, so the inequality holds in that case as well.

Take $x=4\log(n+1)+2\log(M+1)$, integrate the conditional bound, and sum over all $M\ge1$. The probability of failure for at least one $M$ in row $n$ is at most
\[
 (n+1)^{-4}\sum_{M\ge1}(M+1)^{-2}.
\]
Its sum in $n$ is finite. Borel--Cantelli and \eqref{S:eq:gaussian-operator-mode-tail} imply, on one event, for every sufficiently large $n$ and simultaneously for all $M\ge1$,
\begin{equation}\label{S:eq:s10-uniform-mode-truncation}
 \sqrt n\|G_n-G_n^{(M)}\|_{\mathscr K}
 \le C\left(
       \frac{e_n^{\mathrm{KL}}[1+\log(n+1)+\log(M+1)]}{M}
                    \right)^{1/2}.
\end{equation}
The constant is deterministic. The trace of the conditional covariance of each remainder gives a summable probability bound under the original joint law. This applies to the
dependent rows and tails with $e_n^{\mathrm{KL}}$ finite almost surely.

In particular, let
\begin{equation}\label{S:eq:s10-logarithmic-mode-cutoff}
 M(n)=\left\lceil(1+\log n)^{502}\right\rceil.
\end{equation}
For large $n$, its logarithm is at most $\log2+502\log(1+\log n)$. Thus \eqref{S:eq:s10-uniform-mode-truncation} and \eqref{S:eq:weighted-sine-coefficient-sum} give
\begin{equation}\label{S:eq:s10-logarithmic-mode-rates}
 \begin{aligned}
 \sqrt n\|G_n-G_n^{(M(n))}\|_{\mathscr K}
     &\le C_{\beta,\omega}(\log n)^{-501/2},\\
 n\|Q_n-Q_n^{(M(n))}\|_{\mathfrak S_1(\mathscr K)}
     &\le C_{\beta,\omega}(\log n)^{-502}.
 \end{aligned}
\end{equation}
The exponential estimate just proved applies to this choice of $M(n)$. Replacing $G_n$ by the finite sum adds $H_{n,M(n)}$ to the remainder. Since $501/2>250$, its normalized bound
remains $O_\omega((\log n)^{-250})$.

The uniform estimate \eqref{S:eq:s10-uniform-mode-truncation} also holds for an order selected from the sampled positions. If $M(n)$ is a finite, positive,
$\mathscr S_n$-measurable integer, then, for $A\in\mathscr K$,
\[
 \mathbb E[\me^{\ii\langle G_n^{(M(n))},A\rangle} \mid\mathscr S_n] =\sum_{l\ge1}\mathbf1_{\{M(n)=l\}} \exp\left(-\tfrac12\langle Q_n^{(l)}A,A\rangle\right) =\exp\left(-\tfrac12\langle Q_n^{(M(n))}A,A\rangle\right).
\]
The first equality follows by partitioning into the $\mathscr S_n$-measurable events $\{M(n)=l\}$; boundedness justifies summation under conditional expectation. The selection
index is $\mathscr S_n$-measurable, so the displayed conditional characteristic function remains Gaussian. This argument applies to selection from $\mathscr S_n$, with the Gaussian
coordinates integrated conditionally as above.

\subsubsection{The full inverse differences}
\label{S:sec:s10-full-inverse-differences}

The expansion into a conditionally Gaussian term and a remainder describes the error from sampling the Brownian path. Adding the smaller difference between the finite inverse and
the inverse obtained by sampling the Brownian path leaves its leading Gaussian term and conditional mean unchanged; only the remainder changes.

Recall the retained finite comparison
\[
 F_n=\mathsf R_{\beta,n}^{\mathrm{circ},T_n}
             -\mathsf R_{\beta,n}^{\mathrm{br},T_n},
 \qquad R_n^{\mathrm{fin}}=R_n^{\mathrm{br}}+F_n.
\]
The leading terms $G_n,B_n$ and covariance $Q_n$ are unchanged.

\begin{corollary}\label{S:cor:15.2}
Almost surely,
\[
 \Delta_{\beta,n}^{T}=G_n+B_n+R_n^{\mathrm{fin}},\qquad
 \sqrt n\|R_n^{\mathrm{fin}}\|_{\mathscr K}
       =O_{\beta,\omega}((\log n)^{-250}).
\]
\end{corollary}
\begin{proof}
The identity follows by adding $F_n$ to Theorem~\ref{S:thm:15.1}. Proposition~\ref{S:prop:circular-brownian-inverse-comparison} and \eqref{S:eq:circular-brownian-retained-rate}
give
\[
 \begin{aligned}
 \sqrt n\|R_n^{\mathrm{fin}}\|_{\mathscr K}
 &\le\sqrt n\|R_n^{\mathrm{br}}\|_{\mathscr K}
        +\sqrt n\|F_n\|_{\mathscr K}\\
 &\le C_{\beta,\omega}
       \left((\log n)^{-250}+(\log n)^{-501/2}\right)
 \le C'_{\beta,\omega}(\log n)^{-250}.
 \end{aligned}
\]
 \end{proof}

Let $\mathsf T_{\beta,n}^{\diamond}$ be the omitted kernels in \eqref{S:eq:boundary-tail-operator-definitions}. Keep $G_n,B_n,Q_n$ fixed, with $G_n,B_n$ extended by zero outside
the square of the restricted kernel, and define
\begin{equation}\label{S:eq:gaussian-full-remainder-definitions}
 \widehat R_n^{\mathrm{br}} =R_n^{\mathrm{br}}+\mathsf T_{\beta,n}^{\mathrm{br}} -\mathsf T_{\beta,n}^{\mathrm{sine}},\qquad \widehat R_n^{\mathrm{fin}} =R_n^{\mathrm{fin}}+\mathsf T_{\beta,n}^{\mathrm{circ}} -\mathsf T_{\beta,n}^{\mathrm{sine}}.
\end{equation}

\begin{corollary}\label{S:cor:gaussian-full-interval}
The full inverse differences satisfy
\begin{equation}\label{S:eq:gaussian-full-expansions}
 \mathsf R_{\beta,n}^{\mathrm{br}}-\mathsf R_\beta^{\mathrm{sine}} =G_n+B_n+\widehat R_n^{\mathrm{br}},\qquad \Delta_{\beta,n}=G_n+B_n+\widehat R_n^{\mathrm{fin}}.
\end{equation}
Both remainders have normalized norm $O_{\beta,\omega}((\log n)^{-250})$ almost surely.
\end{corollary}
\begin{proof}
Adding \eqref{S:eq:gaussian-full-remainder-definitions} to the retained identities gives \eqref{S:eq:gaussian-full-expansions}. Subtracting the remainder definitions gives
\[
 \widehat R_n^{\mathrm{fin}}-\widehat R_n^{\mathrm{br}} =F_n+\mathsf T_{\beta,n}^{\mathrm{circ}} -\mathsf T_{\beta,n}^{\mathrm{br}} =\mathsf R_{\beta,n}^{\mathrm{circ}} -\mathsf R_{\beta,n}^{\mathrm{br}} =F_n^{\mathrm{full}}.
\]
Theorem~\ref{S:thm:full-boundary-tails}, with $\varepsilon=a_\beta/2$, gives
\[
 n\|\mathsf T_{\beta,n}^{\alpha}
                -\mathsf T_{\beta,n}^{\chi}\|_{\mathscr K}^2
 \le2\mathfrak t_n^{\mathrm{end}}
 \le C_{\beta,\omega}n^{-a_\beta/2},
 \qquad \alpha,\chi\in\{\mathrm{circ},\mathrm{br},\mathrm{sine}\}.
\]
If $K_R$ is a retained remainder kernel and $K_T$ an omitted kernel, $\operatorname{supp}K_R\subset[0,T_n]^2$ and $K_T=0$ almost everywhere there. Thus
$\langle R,T\rangle_{\mathscr K}=\int\!\int\operatorname{tr}(K_R^{\mathsf T}K_T)=0$. Using \eqref{S:eq:gaussian-expansion-quantitative} and Corollary~\ref{S:cor:15.2}, this proves,
for
$\diamond\in\{\mathrm{br},\mathrm{fin}\}$,
\[
 n\|\widehat R_n^\diamond\|_{\mathscr K}^2
 \le C_{\beta,\omega}
             \left((\log n)^{-500}+n^{-a_\beta/2}\right).
\]
For fixed $\beta>0$,
\[
 \log\left(n^{-a_\beta/2}(\log n)^{500}\right)
 =-\frac{a_\beta}{2}\log n+500\log\log n\longrightarrow-\infty.
\]
Taking square roots gives
\begin{equation}\label{S:eq:gaussian-full-remainder-rates}
 \sqrt n\|\widehat R_n^{\mathrm{br}}\|_{\mathscr K}
 \le C_{\beta,\omega}(\log n)^{-250},\qquad
 \sqrt n\|\widehat R_n^{\mathrm{fin}}\|_{\mathscr K}
 \le C_{\beta,\omega}(\log n)^{-250}.
\end{equation}
The omitted kernels are measurable differences of the full and compressed Green kernels. The common event is the intersection with the event of
Theorem~\ref{S:thm:full-boundary-tails}; its probability remains one on the bridge extension.
\end{proof}

The bound $\sqrt n\|G_n+B_n\|_{\mathscr K}\le C_{\beta,\omega}$ and polarization give, for either full remainder,
\[
 \begin{aligned}
 &n\left|\|G_n+B_n+\widehat R_n^\diamond\|_{\mathscr K}^2
                           -\|G_n+B_n\|_{\mathscr K}^2\right|\\
 &\le2(\sqrt n\|G_n+B_n\|_{\mathscr K})
              (\sqrt n\|\widehat R_n^\diamond\|_{\mathscr K})
       +n\|\widehat R_n^\diamond\|_{\mathscr K}^2\\
 &\le C_{\beta,\omega}(\log n)^{-250}\longrightarrow0.
 \end{aligned}
\]
The cutoff \eqref{S:eq:s10-logarithmic-mode-cutoff} also permits $G_n^{(M(n))}$ in these expansions, with $G_n-G_n^{(M(n))}$ added to the remainder. The Brownian path and circular
operators are unchanged. The Gaussian distribution is that of $G_n$ conditional on $\mathscr S_n$; the remainder is a measurable error with the stated norm bound.
Section~\ref{S:sec:the-limiting-coefficient-of-the-error} evaluates the leading conditional second moment \eqref{S:eq:gaussian-trace-bias-mean}.

\subsection{Covariance estimates}
\label{sec:combined-s-covariance}
\label{S:sec:covariance-estimates}

The conditional mean-square norm of the centered Gaussian term is the covariance trace. Almost sure concentration also requires a smaller covariance operator norm. The estimates
below keep these two scales separate and are uniform over the number of retained bridge modes.

Fix $\beta>0$ and use the random operators constructed in Section~\ref{S:sec:the-conditional-gaussian-operator}. Set $\gamma=2/\beta$ and $a_\beta=\min\{1,\gamma\}$. All
covariances in this subsection act on the real Hilbert space $\mathscr K=\mathfrak S_2(\mathcal H_{\mathbb R})$. The index $j$ belongs to
$\mathcal J_n^{\mathrm{at}}=\{\ell_n,\ldots,n-1\}$, where $\ell_n=m_n-R_\beta$; the coefficient on a partition interval still has index $m_n\le k<n$. We take $n$ large enough that
$1\le\ell_n$, $m_n<n$, and $\ell_n\ge m_n/2$.

Work on the intersection of the events in Theorem~\ref{S:thm:5.1}, Propositions~\ref{S:prop:5.2}--\ref{S:prop:5.5}, Lemma~\ref{S:lem:s7-adjoint}, and
Theorem~\ref{S:thm:gaussian-finite-mode}, for positive rational losses in the boundary estimates. Since $0<k/n\le1$, a smaller rational loss gives any positive loss. Conditional
estimates will be taken under the law given $\mathscr S_n$.

The trace is of order at most $n^{-1}$. To make the normalized squared error concentrate, we need the largest covariance eigenvalue to be smaller by a positive power of $n$. The
weighted field covariance already provides this gain. We first use it to prove concentration, including uniformity over the number of Fourier modes in the partial sum. We then
compute the blocks of the Gram matrix and prove a trace lower bound from interior intervals.

\subsubsection{Upper bounds for the covariance}
\label{S:sec:s11-covariance-upper}

The trace bound follows by summing the squared norms of the modes. The operator norm requires their mutual geometry as well, and the two resulting estimates enter the concentration
argument in different ways.

\begin{proposition}\label{S:prop:covariance-upper-limit}
Almost surely, for every $\varepsilon>0$ and all sufficiently large $n$,
\begin{equation}\label{S:eq:covariance-finite-upper}
 \operatorname{tr}Q_n^{(M)}\le\frac{C_{\beta,\omega}}n,
 \qquad
 \|Q_n^{(M)}\|\le C_{\beta,\omega,\varepsilon}
                           n^{-1-a_\beta+\varepsilon},
 \qquad M\ge1.
\end{equation}
The same estimates hold for $Q_n$.
\end{proposition}
\begin{proof}
The bound for the sum of $\|\Psi_{n,j,m,\nu}\|_{\mathscr K}^2$ \eqref{S:eq:s10-summed-mode-bounds} gives
\[
 \operatorname{tr}Q_n^{(M)}
 =\sum_{j,m\le M,\nu}\|\Psi_{n,j,m,\nu}\|_{\mathscr K}^2
 \le\frac{5e_n^{\mathrm{KL}}}{12n}
 \le\frac{C_{\beta,\omega}}n.
\]
The same bound holds for $Q_n$ by summing all modes. For the operator norm, use the identity for the covariance in \eqref{S:eq:gaussian-covariance-pushforward}, with the adjoint in
the full weighted field space:
\begin{equation}\label{S:eq:covariance-weighted-space-bound}
 Q_n=\mathscr D_n\widetilde Q_n^+\mathscr D_n^*,\qquad
 \|Q_n\|\le\frac1n\|\widetilde Q_n^+\|
 \le C_{\beta,\omega,\varepsilon}n^{-1-a_\beta+\varepsilon}.
\end{equation}
Here $\|\mathscr D_n\|\le n^{-1/2}$, and \eqref{S:eq:s9-field-covariance-direct-sum} bounds the field covariance by the largest contribution of a sampling interval. Positivity of
$Q_n-Q_n^{(M)}$ gives the same bound for every $M$. This argument uses orthogonality only in the weighted space of fields, before applying $\mathscr D_n$.

The series and its limit in trace norm were constructed in Theorem~\ref{S:thm:gaussian-finite-mode}, before these uniform bounds. The identity for the positive covariance tail
\eqref{S:eq:gaussian-operator-mode-tail} gives convergence for each fixed row $n$.
\end{proof}

If $\lambda_i\ge0$ are the eigenvalues of a positive trace-class operator $Q$, then $\sum_i\lambda_i^2\le(\sup_i\lambda_i)\sum_i\lambda_i$. It follows that
\begin{equation}\label{S:eq:covariance-square-trace}
 \operatorname{tr}(Q_n^2)\le\|Q_n\|\operatorname{tr}Q_n
       \le C_{\beta,\omega,\varepsilon}n^{-2-a_\beta+\varepsilon}.
\end{equation}
The inequality $\operatorname{tr}((Q_n^{(M)})^2) \le\|Q_n^{(M)}\|\operatorname{tr}Q_n^{(M)}$ gives the same bound for each $M$.
\subsubsection{Bounds for maxima of Gaussian quadratic forms}
\label{S:sec:s11-gaussian-quadratic-forms}

We apply Lemma~\ref{S:lem:covariance-gaussian-concentration} to the exponential submartingales generated by the Gaussian partial sums. For a covariance $Q$ and a vector $b$ fixed
under the conditioning, write
\[
 T=\operatorname{tr}Q,\qquad \Lambda=\|Q\|,\qquad
 V=\operatorname{tr}(Q^2),\qquad v_b=\langle Qb,b\rangle.
\]
Fix $n$ and its sampled positions, and use the partial sums of Section~\ref{S:sec:the-conditional-gaussian-operator}. Abbreviate $G^{(M)}=G_n^{(M)}$, $Q^{(M)}=Q_n^{(M)}$, $G=G_n$,
and $Q=Q_n$. Set $G^{(0)}=0$, $Q^{(0)}=0$, and use the conditional filtration
\[
 \mathcal F_M^{\mathrm{mode}}
   =\sigma\{\xi_{n,j,m,\nu}:j\in\mathcal J_n^{\mathrm{at}},
                                  m\le M,\ \nu=1,2\}.
\]
Work under the law given $\mathscr S_n$. For a vector $b$ fixed by this conditioning, let
\[
 X_M^b=\|G^{(M)}+b\|^2-\operatorname{tr}Q^{(M)}-\|b\|^2,
 \qquad Y_M=\langle G^{(M)},b\rangle.
\]
Both are martingales. Indeed, $H_{M+1}=G^{(M+1)}-G^{(M)}$ is independent of $\mathcal F_M^{\mathrm{mode}}$, centered Gaussian, and has covariance $Q^{(M+1)}-Q^{(M)}$. Thus
\[
 \begin{aligned}
 X_{M+1}^b-X_M^b
 &=2\langle G^{(M)}+b,H_{M+1}\rangle
       +\|H_{M+1}\|^2
           -\operatorname{tr}(Q^{(M+1)}-Q^{(M)}),\\
 Y_{M+1}-Y_M&=\langle H_{M+1},b\rangle
 \end{aligned}
\]
have zero conditional expectation given $\mathcal F_M^{\mathrm{mode}}$.

For $2|\theta|\Lambda<1$, conditional Jensen's inequality makes $\me^{\theta X_M^b}$ a nonnegative submartingale. For every finite $K$, its expectation is bounded by
\[
 \mathbb E\me^{\theta X_K^b}
 \le \exp\!\left(\frac{\theta^2(V+2v_b)}{1-2|\theta|\Lambda}\right).
\]
To check the parameters, $0\le Q^{(K)}\le Q$ implies $\|Q^{(K)}\|\le\Lambda$ and $\langle Q^{(K)}b,b\rangle\le v_b$. Also, for $P=Q^{(K)}$,
\[
 \operatorname{tr}(Q^2)-\operatorname{tr}(P^2)
 =\operatorname{tr}((Q-P)Q)+\operatorname{tr}(P(Q-P))\ge0.
\]
Each trace is nonnegative: for positive operators $A,B$ of the stated classes, $\operatorname{tr}(AB) =\operatorname{tr}(B^{1/2}AB^{1/2})\ge0$. The same expansion and
$|\operatorname{tr}(AB)|\le\|A\|_{\mathfrak S_1}\|B\|$ give
\[
 |\operatorname{tr}(Q^2)-\operatorname{tr}(P^2)|
 \le(\|Q\|+\|P\|)\|Q-P\|_{\mathfrak S_1}.
\]
Together with $|\langle(Q-P)b,b\rangle|\le\|Q-P\|\|b\|^2$, this also proves convergence of the parameters in the quadratic estimates when the covariances of the finite Gaussian
sums tend to $Q$ in trace norm.

We justify the maximal inequality under this same conditional law. For $\theta>0$, convexity of the exponential and the martingale property of $X_M^b$ make
$Z_M=\mathrm e^{\theta X_M^b}$ a nonnegative submartingale. The preceding exponential-moment bound makes each $Z_M$ integrable when $2\theta\Lambda<1$. Fix a finite cutoff $K$ and
let $\tau$ be the first $M\le K$ with $X_M^b>t$, taking $\tau=\infty$ if there is none. For $M\le K$, the event $\{\tau=M\}$ belongs to $\mathcal F_M^{\mathrm{mode}}$, and
$\mathbb E[Z_K\mid\mathcal F_M^{\mathrm{mode}}]\ge Z_M$. Therefore
\[
 \begin{aligned}
 \mathrm e^{\theta t}\mathbb P(\tau\le K)
 &\le\sum_{M=0}^K\mathbb E[\mathbf1_{\{\tau=M\}}Z_M]\\
 &\le\sum_{M=0}^K\mathbb E[\mathbf1_{\{\tau=M\}}Z_K]
 \le\mathbb E Z_K.
 \end{aligned}
\]
We have proved
\[
 \mathbb P\!\left(\max_{M\le K}X_M^b>t\right)
 \le\mathrm e^{-\theta t}\mathbb E\mathrm e^{\theta X_K^b}.
\]
The calculation also applies to $-X_M^b$ and to either sign of the linear martingale $Y_M$. The argument uses the martingale property of the same Gaussian partial sums. The same
first passage bound at parameter $-\theta$ gives, for $A=V+2v_b$ and $0<2\theta\Lambda<1$,
\[
 \mathbb P\left(\max_{M\le K}|X_M^b|>t\right)
 \le2\exp\left(-\theta t+\frac{\theta^2A}{1-2\theta\Lambda}\right).
\]
For $A>0$, choose $\theta=\sqrt{x}/(\sqrt A+2\Lambda\sqrt{x})$ and $t=2\sqrt{Ax}+2\Lambda x$. The exponent is $-x$. For the linear martingale and $v_b>0$,
\[
 \mathbb P\left(\max_{M\le K}|Y_M|>t\right)
 \le2\exp(-\theta t+\theta^2v_b/2),
 \qquad \theta=t/v_b.
\]
Taking $t=\sqrt{2v_bx}$ gives the same bound $2\me^{-x}$. The zero variance cases vanish identically. Letting $K\to\infty$ by continuity from below, and taking $b=0$ for the
centered quadratic term, gives
\begin{equation}\label{S:eq:s11-mode-maximal-tails}
 \begin{aligned}
 \mathbb P\left(\sup_{M\ge0}|X_M^0|>
                             2\sqrt{Vx}+2\Lambda x\right)&\le2\me^{-x},\\
 \mathbb P\left(\sup_{M\ge0}|Y_M|>\sqrt{2v_bx}\right)&\le2\me^{-x},\\
 \mathbb P\left(\sup_{M\ge0}|X_M^b|>
                       2\sqrt{(V+2v_b)x}+2\Lambda x\right)&\le2\me^{-x}.
 \end{aligned}
\end{equation}
The same inequalities include $M=\infty$, meaning $G^{(\infty)}=G$ and $Q^{(\infty)}=Q$: the partial sums converge almost surely in $\mathscr K$ and the traces converge, so a
strict exceedance by the limit is an exceedance by some finite partial sum. Zero variance cases have identically zero martingales.

\subsubsection{Almost sure concentration}
\label{S:sec:s11-almost-sure-concentration}

The trace, covariance norm, and $\|B_n\|^2$ are $\mathscr S_n$-measurable. The events on which these quantities are bounded by a fixed integer therefore satisfy the localization
hypothesis of Lemma~\ref{S:lem:conditional-localized-borel-cantelli}.

Set $G_n^{(\infty)}=G_n$ and $Q_n^{(\infty)}=Q_n$. The supremum over $M$ below includes $\infty$.

\begin{theorem}\label{S:thm:covariance-concentration}
On one event of probability one, for every $\varepsilon>0$,
\begin{equation}\label{S:eq:covariance-concentration-rate}
 \sup_{M\in\mathbb N\cup\{\infty\}}
 \left(\left|n\|G_n^{(M)}\|_{\mathscr K}^2
                         -n\operatorname{tr}Q_n^{(M)}\right|
             +n|\langle G_n^{(M)},B_n\rangle_{\mathscr K}|\right)
       =O_{\beta,\omega,\varepsilon}(n^{-a_\beta/2+\varepsilon}).
\end{equation}
\end{theorem}
\begin{proof}
Write $\mathfrak t_n=\operatorname{tr}Q_n$ and $\mathfrak l_n=\|Q_n\|$. For fixed sampled data, diagonalize $Q_n$ and expand the centered Gaussian in its eigenbasis. If its
eigenvalues are $\lambda_j$, the independent standard normal coordinates satisfy $\mathbb E(Z_j^2-1)^2=2$ and $\mathbb E[(Z_j^2-1)Z_j]=0$. Thus
\[
 \operatorname{Var}(\|G_n+B_n\|_{\mathscr K}^2\mid\mathscr S_n)
 =2\operatorname{tr}(Q_n^2)
          +4\langle Q_nB_n,B_n\rangle_{\mathscr K}.
\]
The identities pass from finite sums in conditional $L^2$: $\sum_j\lambda_j^2<\infty$ and $\sum_j\lambda_j|\langle B_n,e_j\rangle|^2 \le\mathfrak l_n\|B_n\|^2$. The centered
quadratic part and the linear part have zero conditional covariance, not a pathwise cancellation. Fix a rational $\varepsilon_0\in(0,a_\beta)$ and set
$\kappa=a_\beta-\varepsilon_0>0$. For each deterministic integer $L\ge1$ define
\begin{equation}\label{S:eq:covariance-sa-localization}
 \Omega_{L,n}^{\mathrm{sa}}
 :=\{\mathfrak t_n\le L/n,\quad
       \mathfrak l_n\le L n^{-1-\kappa},\quad
       \|B_n\|_{\mathscr K}^2\le L/n\}.
\end{equation}
The covariance construction in Subsection~\ref{S:sec:s10-covariance-and-measurability} shows that this event belongs to $\mathscr S_n$.
Corollary~\ref{S:cor:gaussian-conditional-second-moment} and \eqref{S:eq:weighted-bridge-mean-bound} give
\[
 \mathfrak t_n+\|B_n\|_{\mathscr K}^2
 \le \frac{\mathbb E[\mathcal E_n^+(\mathcal Z_n)\mid\mathscr S_n]
                       +\mathcal E_n^+(\mathcal B_n)}{n}
 \le \frac{C_{\beta,\omega}}n.
\]
Proposition~\ref{S:prop:covariance-upper-limit}, with $\varepsilon=\varepsilon_0$, also gives $\mathfrak l_n\le C_{\beta,\omega,\varepsilon_0}n^{-1-\kappa}$. Choosing an integer
$L_\omega$ larger than these two finite constants puts $\omega$ in $\Omega_{L_\omega,n}^{\mathrm{sa}}$ for every sufficiently large $n$.

Under the conditional law given $\mathscr S_n$, put $V_n=\operatorname{tr}(Q_n^2)$ and $v_n=\langle Q_nB_n,B_n\rangle_{\mathscr K}$. On $\Omega_{L,n}^{\mathrm{sa}}$, the normalized
variance is at most $6L^2n^{-\kappa}$. The quantities in the exponential maximal estimates \eqref{S:eq:s11-mode-maximal-tails} satisfy

\[
 n^2V_n\le n^2\mathfrak l_n\mathfrak t_n
                    \le L^2n^{-\kappa},
 \quad
 n^2v_n\le n^2\mathfrak l_n\|B_n\|^2
                    \le L^2n^{-\kappa},
 \quad n\mathfrak l_n\le Ln^{-\kappa}.
\]
Apply the first two maximal inequalities in \eqref{S:eq:s11-mode-maximal-tails}, with $b=B_n$, before any localization of the conditional law. If $S_n$ denotes the supremum on the
left side of \eqref{S:eq:covariance-concentration-rate}, a union bound for these two events gives
\[
 \mathbf1_{\Omega_{L,n}^{\mathrm{sa}}}
 \mathbb P\!\left(S_n>
       6L(n^{-\kappa/2}\sqrt x+n^{-\kappa}x)
                              \,\middle|\,\mathscr S_n\right)
       \le4\me^{-x}\mathbf1_{\Omega_{L,n}^{\mathrm{sa}}}.
\]
The individual thresholds after multiplication by $n$ are $2L(n^{-\kappa/2}\sqrt x+n^{-\kappa}x)$ and $\sqrt2L n^{-\kappa/2}\sqrt x$, whose sum is at most the displayed threshold.
Multiplying by the $\mathscr S_n$-measurable indicator and using the tower property proves
\begin{equation}\label{S:eq:covariance-sa-row-tail}
 \mathbb P\!\left(\Omega_{L,n}^{\mathrm{sa}}\cap
       \{S_n>6L(n^{-\kappa/2}\sqrt x+n^{-\kappa}x)\}\right)
       \le4\me^{-x},\qquad x>0.
\end{equation}

Take $x=8\log n$, $n\ge2$, and call the localized failure event $E_{L,n}$. For each fixed $L$, $\sum_{n\ge2}\mathbb P(E_{L,n})\le4\sum_{n\ge2}n^{-8}<\infty$.
Lemma~\ref{S:lem:conditional-localized-borel-cantelli} gives
\[
 \mathbb P\left(\bigcap_{L\ge1}(\limsup_n E_{L,n})^c\right)=1.
\]
On this event choose the integer $L_\omega$ supplied by the covariance and bias bounds. For all sufficiently large $n$,
\begin{equation}\label{S:eq:s11-logarithmic-concentration-rate}
 S_n\le C_{\beta,\omega,\varepsilon_0}
       \bigl(n^{-\kappa/2}\sqrt{\log n}
                         +n^{-\kappa}\log n\bigr).
\end{equation}
Take the intersection of these events over rational $\varepsilon_0\in(0,a_\beta)$. Given $\varepsilon>0$, choose $0<\varepsilon_0<\min\{a_\beta/2,\varepsilon\}$ from this set. The
quotient of the first term in \eqref{S:eq:s11-logarithmic-concentration-rate} by $n^{-a_\beta/2+\varepsilon}$ is $n^{\varepsilon_0/2-\varepsilon}\sqrt{\log n}$, which tends to
zero. For the second term the quotient is $n^{-a_\beta/2+\varepsilon_0-\varepsilon}\log n$, which also tends to zero. This proves \eqref{S:eq:covariance-concentration-rate} on one
event for every $\varepsilon>0$.
\end{proof}

Expanding the square for every $M$ gives
\begin{equation}\label{S:eq:covariance-noncentral-rate}
 \sup_{M\in\mathbb N\cup\{\infty\}}
 n\left|\|G_n^{(M)}+B_n\|_{\mathscr K}^2
           -\operatorname{tr}Q_n^{(M)}-\|B_n\|_{\mathscr K}^2\right|
       =O_{\beta,\omega,\varepsilon}(n^{-a_\beta/2+\varepsilon}),
\end{equation}
since the expression inside the supremum is at most $2S_n$. Taking $\varepsilon<a_\beta/2$ proves convergence to zero, in particular for $M=\infty$. The direct noncentered
inequality in \eqref{S:eq:s11-mode-maximal-tails} also gives, on $\Omega_{L,n}^{\mathrm{sa}}$, the threshold $2\sqrt3 L n^{-\kappa/2}\sqrt x+2Ln^{-\kappa}x$ with probability at
most $2\me^{-x}$ after the same tower argument.

For a fixed threshold $\delta>0$, the centered maximal estimate and the scalar Gaussian maximal estimate give
\[
 \begin{aligned}
 &\mathbb P\!\left(\Omega_{L,n}^{\mathrm{sa}}\cap
    \left\{\sup_M|n\|G_n^{(M)}\|^2-n\operatorname{tr}Q_n^{(M)}|
                              >\delta\right\}\right)
 \le2\exp\!\left[-c\min\{\delta^2/L^2,\delta/L\}n^\kappa\right],\\
 &\mathbb P\!\left(\Omega_{L,n}^{\mathrm{sa}}\cap
    \left\{\sup_M n|\langle G_n^{(M)},B_n\rangle|>\delta\right\}\right)
 \le2\exp\!\left[-\delta^2 n^\kappa/(2L^2)\right].
 \end{aligned}
\]
Both probabilities are summable in $n$. This proves the conclusions from the upper bound for the covariance and the bound for the conditional mean. Since the estimates hold
simultaneously for all $M$, they apply at a random cutoff as pathwise inequalities. Conditional Gaussianity, however, has been proved only for cutoffs determined by $\mathscr S_n$.

\subsubsection{Covariances of the contributions from each interval}
\label{S:sec:s11-covariance-weights}

We compute the covariance from the images of the bridge modes. Each derivative kernel occupies two strips. Their intersections determine the rank-one off-diagonal blocks of the
Gram matrix.

The spaces of coefficients for one sampling interval are
\[
 \mathcal H_{n,j}^{(M)}=\mathbb R^{\{1,\ldots,M\}\times\{1,2\}},
 \qquad
 \mathcal H_{n,j}=\ell^2(\mathbb N\times\{1,2\}).
\]
We use the mode fields and operator coefficients already defined in \eqref{S:eq:gaussian-mode-profiles} and \eqref{S:eq:gaussian-finite-operator-definitions}. For a finite vector
$u$, put
\begin{equation}\label{S:eq:covariance-coefficient-maps}
 \begin{aligned}
 \mathcal V_{n,j}^{(M)}[u]
    &:=\sum_{m=1}^M\sum_{\nu=1}^2
                     u_{m,\nu}\mathcal V_{n,j,m,\nu},\\
 \mathcal A_{n,j}^{(M)}u
    &:=-\mathscr D_n\mathcal V_{n,j}^{(M)}[u],\\
 \mathcal A_n^{(M)}(u_j)_j
    &:=\sum_{j\in\mathcal J_n^{\mathrm{at}}}
                       \mathcal A_{n,j}^{(M)}u_j.
 \end{aligned}
\end{equation}
The domain of the last map is the finite Hilbert direct sum of the $\mathcal H_{n,j}^{(M)}$. Expanding in its coordinate basis gives
\[
 Q_n^{(M)}=\mathcal A_n^{(M)}(\mathcal A_n^{(M)})^*,
 \qquad
 \mathcal G_n^{(M)}=(\mathcal A_n^{(M)})^*\mathcal A_n^{(M)},
\]
\[
 (\mathcal G_n^{(M)})_{j l}
       =(\mathcal A_{n,j}^{(M)})^*\mathcal A_{n,l}^{(M)},
 \qquad
 \|Q_n^{(M)}\|=\|\mathcal G_n^{(M)}\|
                 =\|\mathcal A_n^{(M)}\|^2.
\]
For example, $\|A^*A\|=\sup_{\|v\|=1}\|Av\|^2=\|A\|^2$, and applying the same identity to $A^*$ gives
$\|AA^*\|=\|A\|^2$.

Write $p_{n,k}=|a_{n,k}|^2$ and $q_{n,k}=|c_{n,k}|^2$ for the sampled boundary columns. Recall the full scalar weight $\mathfrak G_{n,k}^+=\mathfrak G_{n,k}+2p_{n,k}q_{n,k}/n$.
Define
\[
 \mathcal N_{n,j}
   =\{k:m_n\le k<n,\ |k-j|\le R_\beta+2\},
\]
\begin{equation}\label{S:eq:bridge-variance-weights}
 \mathfrak v_{n,j}
     :=\frac{h_{n,j}}n
                  \sum_{k\in\mathcal N_{n,j}}\mathfrak G_{n,k}^+,
 \qquad \mathfrak v_n^*:=\max_{j\in\mathcal J_n^{\mathrm{at}}}
                                      \mathfrak v_{n,j}.
\end{equation}
For use in the ordered kernel estimate, set
\[
 \overline p_{n,j}:=\max_{k\in\mathcal N_{n,j}}p_{n,k},
 \qquad
 \overline q_{n,j}:=\max_{k\in\mathcal N_{n,j}}q_{n,k}.
\]
Empty sums and maxima are zero. All weights sum over partition indices $m_n\le k<n$.

\begin{lemma}\label{S:lem:covariance-weight-bounds}
Almost surely, for every $\varepsilon>0$ and all sufficiently large $n$,
\begin{equation}\label{S:eq:bridge-variance-sum-max}
 \sum_j\mathfrak v_{n,j}\le\frac{C_{\beta,\omega}}n,
 \qquad
 \mathfrak v_n^*\le C_{\beta,\omega,\varepsilon}
                           n^{-1-a_\beta+\varepsilon}.
\end{equation}
\end{lemma}
\begin{proof}
Put $d_\beta=2R_\beta+5$. For all sufficiently large $n$, membership $k\in\mathcal N_{n,j}$ implies $j/2\le k\le2j$ and $h_{n,j}=4/(\beta j+1)\le C_\beta/k$. Each $k$ is in at most
$d_\beta$ such sets. Interchanging the finite nonnegative sums gives
\[
 \begin{aligned}
 \sum_j\mathfrak v_{n,j}
 &\le\frac{C_\beta}{n}
       \sum_{k=m_n}^{n-1}\frac{\mathfrak G_{n,k}}k
   +\frac{C_\beta}{n^2}
       \sum_{k=m_n}^{n-1}\frac{p_{n,k}q_{n,k}}k\\
 &\le\frac{C_{\beta,\omega}}n
   +\frac{C_{\beta,\omega}(\log n)^{5/2}}{n^2}
                     (1+\log n)
 \le\frac{C_{\beta,\omega}}n.
 \end{aligned}
\]
We used Proposition~\ref{S:prop:5.5} for the first sum and \eqref{S:eq:s7-column-product-bound} for the second. The latter sum includes the contribution from the diagonal squares
to the full weight.

Fix $0<\eta<\min\{a_\beta,\varepsilon\}$. Proposition~\ref{S:prop:5.4} and the same product estimate imply
\[
 \mathfrak v_{n,j}
 \le C_{\beta,\omega,\eta}
            n^{-1-a_\beta+\eta}j^{a_\beta-1-\eta}
       +\frac{C_{\beta,\omega}(\log n)^{5/2}}{n^2j}.
\]
The exponent of $j$ in the first term is negative. In the second term, $j\ge m_n/2$ and $m_n\ge(\log n)^{512}$ give a bound $C_{\beta,\omega}n^{-2}$. Since $a_\beta\le1$ and
$\eta<\varepsilon$, both terms are at most the second bound in \eqref{S:eq:bridge-variance-sum-max}.
\end{proof}

The sampled column estimates and $j/2\le k\le2j$ also give, for every $\eta>0$,
\begin{equation}\label{S:eq:bridge-neighboring-column-bounds}
 \overline p_{n,j}\le C_{\beta,\omega,\eta}(j/n)^{-\gamma-\eta},
 \qquad
 \overline q_{n,j}\le C_{\beta,\omega,\eta}(j/n)^{\gamma-\eta}.
\end{equation}
Indeed, substitute $p_{n,k}\le C(k/n)^{-\gamma-\eta}$ and $q_{n,k}\le C(k/n)^{\gamma-\eta}$, and bound $(k/j)^b$ by $2^{|b|}$ for each of the two exponents $b$.

\subsubsection{Modes on a single sampling interval}
\label{S:sec:s11-single-interval}

A bridge mode is supported on a sampling interval in Brownian time; its operator kernel uses the original interval variable. The time change below identifies its exact support
before inner products of different modes are computed.

We use original time $t$ for kernel calculations. The time change is $\tau_\beta(t)=(4/\beta)\log(1/(1-t))$ and its inverse is $t_\beta(\theta)=1-\me^{-\beta\theta/4}$. Put
\[
 S_{n,j}:=[t_\beta(\vartheta_{n,j+1}),
                         t_\beta(\vartheta_{n,j})],
 \qquad S_{n,j}^T:=S_{n,j}\cap[0,T_n].
\]
The interiors of the $S_{n,j}$ are disjoint and are ordered in the opposite direction to $j$. By \eqref{S:eq:relative-time-at-sample}, each sine function vanishes at
$u_{n,j}(\vartheta_{n,k+1})$. Consequently
\begin{equation}\label{S:eq:s11-exact-mode-support}
 \mathcal V_{n,j,m,\nu}(t)
   =\frac{\sqrt{2h_{n,j}}}{\pi m}
     \sin\!\bigl(\pi m u_{n,j}(\tau_\beta(t))\bigr)
          \mathsf A_\nu\,\mathbf1_{[0,T_n]}(t),
\end{equation}
which vanishes outside $S_{n,j}^T$ and at its sampling boundary points. Here $u_{n,j}$ is clipped to $[0,1]$, as in \eqref{S:eq:interval-relative-time}. A nonzero value requires
$|k-j|\le R_\beta$ on the partition interval $I_{n,k}$. Thus $S_{n,j}^T$ is contained, up to boundary points, in at most $2R_\beta+1$ retained partition intervals and has measure
at most $(2R_\beta+1)/n$.

\begin{lemma}\label{S:lem:bridge-single-interval}
The maps in \eqref{S:eq:covariance-coefficient-maps} satisfy
\begin{equation}\label{S:eq:bridge-single-interval-map}
 \|\mathcal A_{n,j}^{(M)}u\|_{\mathscr K}^2
       \le\frac5{16}\mathfrak v_{n,j}\|u\|^2.
\end{equation}
They extend to Hilbert--Schmidt maps $\mathcal A_{n,j}:\mathcal H_{n,j}\to\mathscr K$, with the same bound.
\end{lemma}
\begin{proof}
For $s\in[0,1]$, let
\[
 b_m(s)=\frac{\sqrt2\sin(\pi m s)}{\pi m}.
\]
The bridge covariance identity gives $\sum_{m\ge1}b_m(s)^2=s(1-s)\le1/4$. For each $\nu$, Cauchy--Schwarz therefore gives
\[
 \left|\sqrt{h_{n,j}}\sum_{m=1}^M u_{m,\nu}b_m(s)\right|^2
 \le \frac{h_{n,j}}4\sum_{m=1}^M|u_{m,\nu}|^2.
\]
Since $\|\mathsf A_1\|_{\mathrm{op}}=1$ and $\|\mathsf A_2\|_{\mathrm{op}}=1/2$, the triangle inequality and Cauchy--Schwarz in $\nu$ yield
\[
 \|v_1\mathsf A_1+v_2\mathsf A_2\|_{\mathrm{op}}^2
       \le\frac54(|v_1|^2+|v_2|^2).
\]
Substituting $s=u_{n,j}(\tau_\beta(t))$ gives
\begin{equation}\label{S:eq:covariance-pointwise-mode-bound}
 \sup_t\|\mathcal V_{n,j}^{(M)}[u](t)\|_{\mathrm{op}}^2
               \le\frac5{16}h_{n,j}\|u\|^2.
\end{equation}
The support property and the full field norm give
\[
 \mathcal E_n^+(\mathcal V_{n,j}^{(M)}[u])
 \le\frac5{16}h_{n,j}
          \sum_{k\in\mathcal N_{n,j}}\mathfrak G_{n,k}^+\|u\|^2
 =\frac5{16}n\mathfrak v_{n,j}\|u\|^2.
\]
Applying \eqref{S:eq:continuous-frechet-bound}, whose squared operator norm bound is $1/n$, proves \eqref{S:eq:bridge-single-interval-map}.

For a single mode, $\|\mathcal V_{n,j,m,\nu}(t)\|_{\mathrm{op}}^2 \le2h_{n,j}\|\mathsf A_\nu\|_{\mathrm{op}}^2/(\pi^2m^2)$. Integration with the full weights and
\eqref{S:eq:continuous-frechet-bound} give
\[
 \|\Psi_{n,j,m,\nu}\|_{\mathscr K}^2
       \le\frac{2\|\mathsf A_\nu\|_{\mathrm{op}}^2}{\pi^2m^2}
                                                    \mathfrak v_{n,j}.
\]
Summing in $\nu$ and then in $m$, using $\sum_{m\ge1}m^{-2}=\pi^2/6$ and $\sum_{m>M}m^{-2}\le\int_M^\infty x^{-2}\dd x=1/M$, gives
\begin{equation}\label{S:eq:bridge-single-interval-trace}
 \sum_{m=1}^M\sum_{\nu=1}^2\|\Psi_{n,j,m,\nu}\|_{\mathscr K}^2 \le\frac5{12}\mathfrak v_{n,j},\qquad \sum_{m>M}\sum_{\nu=1}^2\|\Psi_{n,j,m,\nu}\|_{\mathscr K}^2 \le\frac5{2\pi^2M}\mathfrak v_{n,j}.
\end{equation}
Let $P_M$ project $\mathcal H_{n,j}$ onto the coordinates $(m,\nu)$ with $m\le M$. The finite maps agree on their common coordinates and, by
\eqref{S:eq:bridge-single-interval-map}, their norms are at most $\sqrt{5\mathfrak v_{n,j}/16}$. Their common restriction to finite sequences therefore extends uniquely to a
bounded linear map $\mathcal A_{n,j}:\mathcal H_{n,j}\to\mathscr K$, with
\[
 \mathcal A_{n,j}P_M=\mathcal A_{n,j}^{(M)}P_M,
 \qquad
 \|\mathcal A_{n,j}\|^2\le\frac5{16}\mathfrak v_{n,j}.
\]
 Its image of the coordinate vector $e_{m,\nu}$ is $\Psi_{n,j,m,\nu}$, so
\[
 \|\mathcal A_{n,j}\|_{\mathfrak S_2}^2
 =\sum_{m\ge1}\sum_{\nu=1}^2\|\Psi_{n,j,m,\nu}\|_{\mathscr K}^2
 \le\frac5{12}\mathfrak v_{n,j}.
\]
The tail bound also gives $\|\mathcal A_{n,j}(I-P_M)\|_{\mathfrak S_2}^2 \le5\mathfrak v_{n,j}/(2\pi^2M)$. Thus the finite maps, extended by zero, converge to the full map in
Hilbert--Schmidt norm.
\end{proof}

For all $j,l$, including equal or adjacent indices, the preceding lemma gives
\begin{equation}\label{S:eq:covariance-neighbor-block}
 \|(\mathcal A_{n,j}^{(M)})^*\mathcal A_{n,l}^{(M)}\|
       \le\frac5{16}\sqrt{\mathfrak v_{n,j}\mathfrak v_{n,l}}
       \le\frac5{16}\mathfrak v_n^*.
\end{equation}
The next calculation gives the inner product of the images under $\mathscr D_n$ for distinct $j,l$.

\subsubsection{Inner products of the operators $\Psi_{n,j,m,\nu}$}
\label{S:sec:s11-gram-matrix}

The order of two disjoint sampling intervals determines which strips of their derivative kernels overlap. The corresponding inner product factors into two linear functionals,
giving the rank-one block below.

Define vectors in $\mathcal H_{n,j}^{(M)}$ by their coordinates:
\[
 \begin{aligned}
 (\mathfrak a_{n,j}^{(M)})_{m,\nu}
   &:=\int_{S_{n,j}^T}
       \langle\mathcal V_{n,j,m,\nu}(t)a_n(t),a_n(t)\rangle\dd t,\\
 (\mathfrak c_{n,j}^{(M)})_{m,\nu}
   &:=\int_{S_{n,j}^T}
       \langle\mathcal V_{n,j,m,\nu}(t)c_n(t),c_n(t)\rangle\dd t.
 \end{aligned}
\]
Thus these vectors represent the two displayed linear functionals on the coefficient space. For $v,w$ in these real coefficient spaces, $v\otimes w$ denotes the map
$x\mapsto v\langle w,x\rangle$.

\begin{lemma}\label{S:lem:bridge-ordered-block}
For $j>l$, the off-diagonal block of the Gram matrix is
\begin{equation}\label{S:eq:s11-ordered-rank-one}
 (\mathcal A_{n,j}^{(M)})^*\mathcal A_{n,l}^{(M)}
       =\tfrac12\mathfrak a_{n,j}^{(M)}
                           \otimes\mathfrak c_{n,l}^{(M)},
\end{equation}
and
\begin{equation}\label{S:eq:bridge-ordered-block-forward}
 \|(\mathcal A_{n,j}^{(M)})^*\mathcal A_{n,l}^{(M)}\|
       \le\frac{C_\beta}{n^2}\sqrt{h_{n,j}h_{n,l}}
                           \overline p_{n,j}\overline q_{n,l}.
\end{equation}
For $j<l$ take adjoints. The same identities and bounds hold for the full maps $\mathcal A_{n,j}$.
\end{lemma}
\begin{proof}
Let $U=\mathcal V_{n,j}^{(M)}[u]$ and $V=\mathcal V_{n,l}^{(M)}[v]$, initially with finite coefficient vectors. If $j>l$, every interior point $x$ of $S_{n,j}^T$ precedes every
interior point $y$ of $S_{n,l}^T$ in original time. On $S_{n,j}^T\times S_{n,l}^T$ the kernels of the two first variations are
\[
 -\tfrac12(U(x)a_n(x))c_n(y)^{\mathsf T},
 \qquad
 -\tfrac12a_n(x)(V(y)c_n(y))^{\mathsf T}.
\]
On the transposed rectangle they are the transposes. The two derivative supports are contained, respectively, in
\[
 (S_{n,j}^T\times[0,T_n])\cup([0,T_n]\times S_{n,j}^T),
 \qquad
 (S_{n,l}^T\times[0,T_n])\cup([0,T_n]\times S_{n,l}^T).
\] Since the two field supports are disjoint up to boundary points, the intersection of the two derivative supports is contained
in the indicated rectangles, up to a null set. The identity $\langle ab^{\mathsf T},cd^{\mathsf T}\rangle_{\mathrm{HS}} =\langle a,c\rangle\langle b,d\rangle$ now gives
\begin{equation}\label{S:eq:covariance-ordered-bilinear}
 \langle\mathcal A_{n,j}^{(M)}u, \mathcal A_{n,l}^{(M)}v\rangle_{\mathscr K} =\frac12 \left(\int_{S_{n,j}^T}\langle U(x)a_n(x),a_n(x)\rangle\dd x\right) \left(\int_{S_{n,l}^T}\langle c_n(y),V(y)c_n(y)\rangle\dd y\right).
\end{equation}
The factor is $2(1/2)^2=1/2$. This proves \eqref{S:eq:s11-ordered-rank-one}; in particular the block has rank at most one. The calculation remains valid when both sampling supports
meet the same partition interval. Their order in time determines the two triangles.

For $\|u\|=1$, the first factor without its prefactor is bounded by
\[
 \sup_x\|U(x)\|_{\mathrm{op}}
                  |S_{n,j}^T|\overline p_{n,j}
       \le\frac{C_\beta\sqrt{h_{n,j}}}{n}\overline p_{n,j}.
\]
For $\|v\|=1$, the second factor is at most $\sup_y\|V(y)\|_{\mathrm{op}}|S_{n,l}^T|\overline q_{n,l} \le C_\beta\sqrt{h_{n,l}}\overline q_{n,l}/n$. Multiplying the two bounds
proves \eqref{S:eq:bridge-ordered-block-forward}. The two integral functionals are bounded on the dense finite sequences in $\mathcal H_{n,j}$, with norms at most
$C_\beta\sqrt{h_{n,j}}\,\overline p_{n,j}/n$ and $C_\beta\sqrt{h_{n,j}}\,\overline q_{n,j}/n$, respectively. Their Riesz representatives are
$\mathfrak a_{n,j},\mathfrak c_{n,j} \in\mathcal H_{n,j}$. In particular,
\[
 P_M\mathfrak a_{n,j}\longrightarrow\mathfrak a_{n,j},
 \qquad P_M\mathfrak c_{n,j}\longrightarrow\mathfrak c_{n,j}
 \quad\hbox{in }\mathcal H_{n,j}.
\]

Extend every finite interval map by zero on the remaining coordinates. The preceding lemma gives convergence of these maps in operator norm. Hence
\[
 \| (\mathcal A_{n,j}^{(M)})^*\mathcal A_{n,l}^{(M)} -\mathcal A_{n,j}^*\mathcal A_{n,l}\| \le\|\mathcal A_{n,j}^{(M)}-\mathcal A_{n,j}\| \|\mathcal A_{n,l}^{(M)}\| +\|\mathcal A_{n,j}\| \|\mathcal A_{n,l}^{(M)}-\mathcal A_{n,l}\| \longrightarrow0.
 \]
The rank-one side converges as well, since
\[
 \|a_M\otimes c_M-a\otimes c\|
 \le\|a_M-a\|\,\|c_M\|+\|a\|\,\|c_M-c\|\longrightarrow0
\]
for the coordinate projections just considered. Taking limits proves the identity and norm bound for the full maps. Finally, for $j<l$, taking the adjoint of the identity with
indices $l,j$ gives $\mathcal A_{n,j}^*\mathcal A_{n,l} =\tfrac12\mathfrak c_{n,j}\otimes\mathfrak a_{n,l}$. The norm bound is unchanged under this adjoint, with the two indices
exchanged.
\end{proof}

\begin{lemma}\label{S:lem:bridge-gram-row-sum}
For every $\varepsilon>0$, almost surely for all sufficiently large $n$,
\begin{equation}\label{S:eq:bridge-gram-row-sum}
 \sup_{M\ge1}\sup_j\sum_l
       \|(\mathcal A_{n,j}^{(M)})^*\mathcal A_{n,l}^{(M)}\|
       \le C_{\beta,\omega,\varepsilon}n^{-1-a_\beta+\varepsilon}.
\end{equation}
The same row bound holds for the full maps, and the column supremum has the same bound.
\end{lemma}
\begin{proof}
Choose
\[
 0<\eta<\min\{\varepsilon/8,\gamma/4,1/8\},
 \qquad \eta<(1/2-\gamma)/2\quad\hbox{when }\gamma<1/2.
\]
The diagonal block is bounded by $5\mathfrak v_n^*/16$. For $l<j$, equations \eqref{S:eq:bridge-ordered-block-forward} and \eqref{S:eq:bridge-neighboring-column-bounds} give
\[
 \begin{aligned}
 \sum_{l<j}\|(\mathcal A_{n,j}^{(M)})^*\mathcal A_{n,l}^{(M)}\|
 &\le C_{\beta,\omega,\eta}n^{-2+2\eta}
    j^{-\gamma-\eta-1/2}\sum_{l=1}^{j}l^{\gamma-\eta-1/2}\\
 &\le C_{\beta,\omega,\eta}n^{-2+2\eta}j^{-2\eta}
 \le C_{\beta,\omega,\eta}n^{-2+2\eta}.
 \end{aligned}
\]
To justify the finite power sum, for any $b>-1$ comparison of each $k^b$ with the integral on $[k,k+1]$, using a factor $2^{|b|}$, gives $\sum_{k=1}^j k^b\le C_b j^{b+1}$. Here
$b=\gamma-\eta-1/2>-1$. Extending the lower limit to $1$ only increases the nonnegative majorant.

For $l>j$ the reversed order gives
\begin{equation}\label{S:eq:bridge-regular-side-sum}
 \sum_{l>j}\|(\mathcal A_{n,j}^{(M)})^*\mathcal A_{n,l}^{(M)}\|
 \le C_{\beta,\omega,\eta}n^{-2+2\eta}j^{\gamma-\eta-1/2}
                  \sum_{l=j}^{n}l^{-\gamma-\eta-1/2}.
\end{equation}
If $\gamma\ge1/2$, set $d=\gamma+\eta-1/2>0$. Since
\[
 \sum_{l=j}^{\infty}l^{-1-d}
 \le j^{-1-d}+\int_j^\infty x^{-1-d}\dd x
 \le(1+d^{-1})j^{-d},
\]
this contribution is at most $C_{\beta,\omega,\eta}n^{-2+2\eta}j^{-2\eta}$. At $\gamma=1/2$ the same formula has $d=\eta$ and a finite constant depending on $\eta$. If
$0<\gamma<1/2$, then $d=1/2-\gamma-\eta>0$ and $\sum_{l=j}^n l^{-\gamma-\eta-1/2}\le C_d n^d$. The bound becomes
\[
 C_{\beta,\omega,\eta}n^{-3/2-\gamma+\eta}
                         j^{\gamma-\eta-1/2}
       \le C_{\beta,\omega,\eta}n^{-3/2-\gamma+\eta}.
\]
For $\gamma\ge1$, $a_\beta=1$, and $n^{-2+2\eta}\le n^{-1-a_\beta+\varepsilon}$ because $2\eta<\varepsilon$. For $1/2\le\gamma<1$, the exponent difference is
$(-2+2\eta)-(-1-\gamma+\varepsilon) =\gamma-1+2\eta-\varepsilon<0$. For $0<\gamma<1/2$, the two differences are $\gamma-1+2\eta-\varepsilon<0$ and $-1/2+\eta-\varepsilon<0$.
Together with the diagonal bound these prove \eqref{S:eq:bridge-gram-row-sum}.

The bounds for full blocks follow from Lemma~\ref{S:lem:bridge-ordered-block} applied to the extended maps and \eqref{S:eq:covariance-neighbor-block}; each row has finitely many
interval indices. Finally the block for $(l,j)$ is the adjoint of the block for $(j,l)$, so their norms form a symmetric scalar matrix. Its row and column suprema are equal.
\end{proof}

Let $T_{jl}:\mathcal H_l\to\mathcal K_j$ be bounded maps between separable Hilbert spaces, with finite or countable index sets. Put
\[
 R=\sup_j\sum_l\|T_{jl}\|,
 \qquad C=\sup_l\sum_j\|T_{jl}\|.
\]
\begin{lemma}\label{S:lem:block-schur}
If $R,C<\infty$, the blocks define a bounded map between the Hilbert direct sums, of norm at most $\sqrt{RC}$. For a self-adjoint block matrix the norm is at most $R$.
\end{lemma}
\begin{proof}
Put $b_{jl}=\|T_{jl}\|$, so $\sup_j\sum_l b_{jl}\le R$ and $\sup_l\sum_jb_{jl}\le C$. If $R=0$ or $C=0$, every block vanishes and the result is immediate. Otherwise let $x=(x_l)$
have finite support. For each $j$, the triangle inequality and weighted Cauchy--Schwarz give
\[
 \left\|\sum_lT_{jl}x_l\right\|^2
 \le\left(\sum_l b_{jl}\|x_l\|\right)^2
 \le\left(\sum_l b_{jl}\right)
                  \left(\sum_l b_{jl}\|x_l\|^2\right).
\]
A zero row contributes zero. Summing the nonnegative quantities and interchanging their sums yields
\[
 \sum_j\left\|\sum_lT_{jl}x_l\right\|^2
 \le R\sum_l\|x_l\|^2\sum_jb_{jl}
 \le RC\sum_l\|x_l\|^2.
\]
Thus the coordinate formula defines an element of the target Hilbert direct sum and has norm at most $\sqrt{RC}\|x\|$. It extends uniquely by continuity from the dense subspace of
finitely supported vectors to a bounded operator of norm at most $\sqrt{RC}$. Applying this extension to a vector supported in one coordinate shows that its blocks are exactly the
specified $T_{jl}$. If $T_{jl}^*=T_{lj}$, then $b_{jl}=b_{lj}$ and the row and column suprema agree. Moreover, the scalar product identity for vectors with finite support
$\langle Tx,y\rangle=\langle x,Ty\rangle$ extends by continuity, proving self-adjointness on the common Hilbert direct sum and the bound $\|T\|\le R$.
\end{proof}

Applying these estimates to the covariance gives a direct proof of the bound for the operator norm in Proposition~\ref{S:prop:covariance-upper-limit}:
\[
 \|Q_n^{(M)}\|=\|\mathcal G_n^{(M)}\|
 \le\sup_j\sum_l
       \|(\mathcal A_{n,j}^{(M)})^*\mathcal A_{n,l}^{(M)}\|
 \le C_{\beta,\omega,\varepsilon}n^{-1-a_\beta+\varepsilon}.
\]
The bound passes to $Q_n$ through \eqref{S:eq:covariance-trace-norm-passage}. This calculation gives the off-diagonal blocks of $\mathcal G_n^{(M)}$. The same block formulas apply
to the indices considered in the next subsubsection.

\subsubsection{A lower bound from interior intervals}
\label{S:sec:s11-covariance-lower}

The upper bounds suffice for concentration. The remaining argument shows that the trace cannot be smaller than order $n^{-1}$: a positive proportion of the first bridge modes, away
from the singular boundary, have a uniformly coercive Gram matrix after multiplication by $n^2$. Restricting the kernels to disjoint rectangles separates the modes in the
Hilbert--Schmidt norm and gives the lower bound.

Let
\[
 \mathcal J_n^\circ
   :=\{j\in\mathbb Z:\lceil2n/5\rceil\le j\le\lfloor3n/5\rfloor\},
 \qquad N_n^\circ:=|\mathcal J_n^\circ|.
\]
For large $n$, this is a subset of $\mathcal J_n^{\mathrm{at}}$ and $n/10\le N_n^\circ\le n$. Inside the sampling interval $S_{n,j}$ define
\[
 E_{n,j}:=\{t\in S_{n,j}:1/4\le u_{n,j}(\tau_\beta(t))\le3/4\}.
\]
For the mode $m=1,\nu=2$ write
\[
 \mathcal V_{n,j,1,2}(t)=\varphi_{n,j}(t)\mathsf A_2,
 \qquad
 \varphi_{n,j}(t)
   =\frac{\sqrt{2h_{n,j}}}{\pi}
          \sin\!\bigl(\pi u_{n,j}(\tau_\beta(t))\bigr)
             \mathbf1_{[0,T_n]}(t).
\]

\begin{lemma}\label{S:lem:bridge-first-mode-geometry}
For all sufficiently large $n$ and $j\in\mathcal J_n^\circ$,
\[
 S_{n,j}\subset[1/3,2/3]\subset[0,T_n],
 \qquad c_\beta/n\le|E_{n,j}|\le|S_{n,j}|\le C_\beta/n.
\]
Moreover, $|\varphi_{n,j}|\le C_\beta n^{-1/2}$ everywhere and $|\varphi_{n,j}|\ge c_\beta n^{-1/2}$ on $E_{n,j}$.
\end{lemma}
\begin{proof}
Put $r_{n,j}^+=\rho_{n,j}$ and $x_{n,j}=\beta h_{n,j}/4=(j+\beta^{-1})^{-1}$. For $0\le u\le1$ the original time corresponding to fraction $u$ of this sampling interval is
\[
 t_\beta(\vartheta_{n,j+1}+u h_{n,j})
      =1-r_{n,j}^+\me^{-u x_{n,j}}.
\]
Consequently
\[
 |S_{n,j}|=r_{n,j}^+(1-\me^{-x_{n,j}}),
 \qquad
 |E_{n,j}|=r_{n,j}^+
          (\me^{-x_{n,j}/4}-\me^{-3x_{n,j}/4}).
\]
Lemma~\ref{S:lem:sampling-position} gives $|n r_{n,j}^+-(j+1)|\le C_\beta^0$. For the chosen indices, both $r_{n,j}^+$ and $r_{n,j}^+\me^{-x_{n,j}}$ lie in $[1/3,2/3]$ for all
sufficiently large $n$, while $x_{n,j}$ and $h_{n,j}$ lie between positive deterministic multiples of $1/n$. For $0<x\le1$,
\[
 \me^{-1}x\le\int_0^x\me^{-v}\dd v\le x,
 \qquad
 \me^{-3/4}x/2\le\int_{x/4}^{3x/4}\me^{-v}\dd v\le x/2.
\]
These inequalities prove the length bounds. On $[1/4,3/4]$ the sine is at least $1/\sqrt2$, and on $[0,1]$ it is at most one. Substitution in $\varphi_{n,j}$ proves its two
estimates. Finally $m_n/n\to0$ implies that $[1/3,2/3]$ and $[3/4,7/8]$ are contained in $[0,T_n]$ for large $n$.
\end{proof}

Let $\Psi_{n,j}^\circ:=\Psi_{n,j,1,2}$. The fixed matrix $\mathsf A_2$ satisfies $\mathsf A_2^{\mathsf T}\mathsf
A_2=I_2/4$.

\begin{proposition}\label{S:prop:bridge-interior-coefficient-lower-bound}
Almost surely, for all sufficiently large $n$ and $j,l\in\mathcal J_n^\circ$,
\begin{equation}\label{S:eq:bridge-interior-coefficient-bounds}
 c_{\beta,\omega}n^{-2}\le\|\Psi_{n,j}^\circ\|_{\mathscr K}^2
                       \le C_{\beta,\omega}n^{-2}.
\end{equation}
For $j\ne l$,
\begin{equation}\label{S:eq:bridge-interior-gram-bound}
 |\langle\Psi_{n,j}^\circ,\Psi_{n,l}^\circ\rangle_{\mathscr K}|
                       \le C_{\beta,\omega}n^{-3}.
\end{equation}
\end{proposition}
\begin{proof}
The limiting columns are continuous and nonzero on compact subintervals of $[0,1)$. Their uniform sampled convergence in Lemma~\ref{S:lem:sampled-exclusive-convergence} implies
\[
 \inf_{t\in[1/3,2/3]}|a_n(t)|^2\ge c_{\beta,\omega},
 \qquad
 \int_{3/4}^{7/8}|c_n(t)|^2\dd t\ge c_{\beta,\omega},
\]
and both columns are uniformly bounded above on these intervals. The vector $u_1$ may depend on the path; after fixing the path it is nonzero, and multiplication by an invertible
matrix preserves nonvanishing.

On $E_{n,j}\times[3/4,7/8]$ the field vanishes at the second variable, so the kernel of $\Psi_{n,j}^\circ$ is exactly
\[
 -\tfrac12\varphi_{n,j}(x)
               (\mathsf A_2a_n(x))c_n(y)^{\mathsf T}.
\]
Restricting the integral of its squared matrix norm to this rectangle gives
\begin{equation}\label{S:eq:s11-interior-restricted-norm}
 \|\Psi_{n,j}^\circ\|_{\mathscr K}^2
 \ge\frac1{16}
  \left(\int_{E_{n,j}}|\varphi_{n,j}(x)|^2|a_n(x)|^2\dd x\right)
  \left(\int_{3/4}^{7/8}|c_n(y)|^2\dd y\right)
 \ge c_{\beta,\omega}n^{-2}.
\end{equation}
The factor $1/16$ is the kernel square factor $1/4$ times $|\mathsf A_2a|^2/|a|^2=1/4$. Only this term of the derivative is present on the rectangle.

For the upper bound, $j$ and every $k\in\mathcal N_{n,j}$ are between fixed positive multiples of $n$. The full scalar weights there are bounded by Proposition~\ref{S:prop:5.4} and
the sampled column bounds. Thus $\mathfrak v_{n,j}\le C_{\beta,\omega}/n^2$, and Lemma~\ref{S:lem:bridge-single-interval} proves the upper estimate. For distinct indices,
\eqref{S:eq:covariance-ordered-bilinear} applies without a separation condition. Each integral in that formula is bounded by
$C_{\beta,\omega}|S_{n,j}|\sqrt{h_{n,j}}\le C_{\beta,\omega}n^{-3/2}$ with its respective index. Their product proves \eqref{S:eq:bridge-interior-gram-bound}, also when the two
sampling intervals meet the same partition interval.
\end{proof}

Define
\[
 A_n^\circ:\mathbb R^{\mathcal J_n^\circ}\longrightarrow\mathscr K,
 \qquad A_n^\circ v=\sum_{j\in\mathcal J_n^\circ}v_j\Psi_{n,j}^\circ,
 \qquad
 Q_n^\circ:=A_n^\circ(A_n^\circ)^*
       =\sum_{j\in\mathcal J_n^\circ}
                         \Psi_{n,j}^\circ\otimes\Psi_{n,j}^\circ.
\]

\begin{theorem}\label{S:thm:bridge-effective-rank}
Almost surely, for all sufficiently large $n$,
\begin{equation}\label{S:eq:s11-interior-gram-coercivity}
 c_{\beta,\omega}n^{-2}\|v\|^2
       \le\|A_n^\circ v\|_{\mathscr K}^2
       \le C_{\beta,\omega}n^{-2}\|v\|^2,
 \qquad v\in\mathbb R^{\mathcal J_n^\circ}.
\end{equation}
In particular, $Q_n^\circ$ has rank $N_n^\circ$ and all its nonzero eigenvalues belong to
$[c_{\beta,\omega}n^{-2},C_{\beta,\omega}n^{-2}]$.
\end{theorem}
\begin{proof}
On $E_{n,j}\times[3/4,7/8]$, every mode in this sum except the $j$th vanishes in both variables. These rectangles are disjoint for different $j$, up to null sets. Therefore the
calculation in \eqref{S:eq:s11-interior-restricted-norm} yields
\[
 \|A_n^\circ v\|_{\mathscr K}^2 \ge\frac1{16} \left(\int_{3/4}^{7/8}|c_n(y)|^2\dd y\right) \sum_{j\in\mathcal J_n^\circ}|v_j|^2 \int_{E_{n,j}}|\varphi_{n,j}(x)|^2|a_n(x)|^2\dd x \ge c_{\beta,\omega}n^{-2}\sum_j|v_j|^2.
\]
For the upper bound, the diagonal entries of the Gram matrix $(A_n^\circ)^*A_n^\circ$ are at most $C_{\beta,\omega}/n^2$, and each off-diagonal entry is at most
$C_{\beta,\omega}/n^3$. Its row sums are at most $C_{\beta,\omega}(n^{-2}+N_n^\circ n^{-3})\le C_{\beta,\omega}n^{-2}$. Lemma~\ref{S:lem:block-schur} proves the upper inequality.
The lower inequality makes $A_n^\circ$ injective. Write $B=A_n^\circ$. If $B^*Bv=\lambda v$ with $\lambda>0$, then $\|Bv\|^2=\lambda\|v\|^2$ and $BB^*(Bv)=\lambda Bv$. If
$BB^*w=\lambda w$, then $B^*w\ne0$ and $B^*B(B^*w)=\lambda B^*w$. The maps $B$ and $\lambda^{-1}B^*$ are inverse on these eigenspaces. Thus the positive eigenvalues agree with
multiplicity; the preceding bounds for the quadratic form place them in the stated interval.
\end{proof}

Each selected mode occurs in $Q_n^{(M)}$ for every $M\ge1$. Subtracting its rank-one summand leaves a positive operator. The positive limit in trace norm therefore gives
\[
 0\le Q_n^\circ\le Q_n^{(M)}\le Q_n.
\]
The preceding theorem and $n/10\le N_n^\circ\le n$ imply
\begin{equation}\label{S:eq:bridge-effective-rank}
 \frac{c_{\beta,\omega}}n\le\operatorname{tr}Q_n^\circ
       \le\frac{C_{\beta,\omega}}n,
 \qquad
 \|Q_n^\circ\|\le\frac{C_{\beta,\omega}}{n^2},
 \qquad
 \frac{\operatorname{tr}Q_n^\circ}{\|Q_n^\circ\|}
       \ge c_{\beta,\omega}n.
\end{equation}
Here rank and eigenvalues refer to $Q_n^\circ$ acting on $\mathscr K$.

\begin{theorem}\label{S:thm:16.5}
On one event of probability one, for every $\varepsilon>0$ and all sufficiently large $n$,
\begin{equation}\label{S:eq:bridge-full-covariance-bounds}
 \frac{c_{\beta,\omega}}n\le\operatorname{tr}Q_n
       \le\frac{C_{\beta,\omega}}n,
 \qquad
 \|Q_n\|\le C_{\beta,\omega,\varepsilon}n^{-1-a_\beta+\varepsilon}.
\end{equation}
The same bounds hold for every $Q_n^{(M)}$, $M\ge1$, with constants independent of $M$.
\end{theorem}
\begin{proof}
For $M\ge1$, each selected first mode occurs in $Q_n^{(M)}$. The remaining rank-one summands are positive. Hence $0\le Q_n^\circ\le Q_n^{(M)}\le Q_n$.
Theorem~\ref{S:thm:bridge-effective-rank} and $N_n^\circ\ge n/10$ give
\[
 \operatorname{tr}Q_n^{(M)}\ge\operatorname{tr}Q_n^\circ
 \ge c_{\beta,\omega}n^{-2}N_n^\circ
 \ge \frac{c_{\beta,\omega}}{10n}.
\]
The same lower bound holds for $Q_n$. By Proposition~\ref{S:prop:covariance-upper-limit}, uniformly in $M$,
\[
 \operatorname{tr}Q_n^{(M)}\le\operatorname{tr}Q_n\le C_{\beta,\omega}/n,
 \qquad
 \|Q_n^{(M)}\|\le\|Q_n\|
       \le C_{\beta,\omega,\eta}n^{-1-a_\beta+\eta}.
\]
The event fixed at the beginning of the section contains these estimates for every admissible rational $\eta>0$. Given $\varepsilon>0$, choosing a rational $\eta\in(0,\varepsilon)$
and using $n^{-1-a_\beta+\eta}\le n^{-1-a_\beta+\varepsilon}$ proves the bound for every $\varepsilon>0$ on the same event.
\end{proof}

Also, for $0<\varepsilon<a_\beta$, $\operatorname{tr}Q_n/\|Q_n\| \ge c_{\beta,\omega,\varepsilon}n^{a_\beta-\varepsilon}$. The covariance norm is nonzero by the trace lower bound,
so the division is legitimate.

\subsection{The leading order coefficient of the error}
\label{sec:combined-s-coefficient}
\label{S:sec:the-limiting-coefficient-of-the-error}

The concentration estimate of Section~\ref{S:sec:covariance-estimates} leaves a conditional second moment to be computed. Fix $\beta>0$, put $\gamma=2/\beta$ and
$a_\beta=\min\{1,\gamma\}$, and write
\[
 \mathfrak C_n:=n\operatorname{tr}Q_n+n\|B_n\|_{\mathscr K}^2
       =n\mathbb E[\|G_n+B_n\|_{\mathscr K}^2\mid\mathscr S_n].
\]
We prove that
\[
 \mathfrak C_n\longrightarrow\PathCoeff(\omega)
 =\frac4\beta\int_0^1\kappa_\beta(r)
                   \frac{\Gamma_{\beta,\omega}(r)}r\,\dd r\in(0,\infty).
\]
The equality defining $\mathfrak C_n$ is Corollary~\ref{S:cor:gaussian-conditional-second-moment}; its limit is Corollary~\ref{S:cor:coefficient-conditional-mean} below. The
calculation uses the matrices $X_\theta=\mathcal A(Z_\theta)$ in the original coordinates and the fixed boundary vectors. All variables remain on the coupling of
Proposition~\ref{S:prop:7.2}, extended in Theorem~\ref{S:thm:joint-probability-space}.

In the squares of the two terms in the derivative of the Green kernel, adding the contributions from linear interpolation and the Brownian bridge gives the Brownian time elapsed
from the sample. This gives the factor $\kappa_\beta$. Products of linear interpolation terms from separated sampling intervals require a further estimate, since their increments
and coefficients are fixed under $\mathscr S_n$. We control those products by conditioning only on the terminal point and summing the increments of the squared norm of the column
associated with the terminal boundary point. The order of the variables in the Green kernel yields this cancellation before absolute values are taken. Once these products are
removed, convergence of the sampled boundary columns gives the integral above.

\subsubsection{The local second moment}
\label{S:sec:conditional-trace-and-coefficient-calculation}

The local second moment uses both operator cells and sampling intervals in Brownian time. Their indices need not coincide, so the boundary columns are recorded on the sampling
range while the kernel primitives sum only over the restricted operator cells.

Write $r_k=k/n$. For every sufficiently large row put
\[
 \mathcal K_n=\{m_n,\ldots,n-1\},\qquad
 \mathcal J_n^{\mathrm{at}}=\{\ell_n,\ldots,n-1\},\qquad
 \ell_n=m_n-R_\beta\ge1.
\]
Partition indices belong to $\mathcal K_n$ and sampling indices to $\mathcal J_n^{\mathrm{at}}$. As before, a coefficient with an invalid index is zero. The boundary columns and
their squared lengths are defined at all sample times, including the finitely many sampling indices below $m_n$:
\[
 a_{n,k}=a_\omega(\rho_{n,k}),\quad
 c_{n,k}=c_\omega(\rho_{n,k}),\quad
 p_{n,k}=|a_{n,k}|^2,\quad q_{n,k}=|c_{n,k}|^2,
 \qquad 0\le k<n.
\]
The partial sums, however, contain only partition indices:
\[
 A_{n,k}=\frac1n\sum_{\substack{l\in\mathcal K_n\\l<k}}q_{n,l},\qquad
 C_{n,k}=\frac1n\sum_{\substack{l\in\mathcal K_n\\l>k}}p_{n,l}.
\]
Thus $A_{n,m_n}=C_{n,n-1}=0$. Fix a deterministic integer $L_\beta>2R_\beta+4$. The separated sum includes all pairs with index distance at least $L_\beta$. Write $\mathbb E_n$ for
integration conditional on $\mathscr S_n$.

On $x<y$ in the original operator coordinate, the first variation is
\[
 (\mathscr D_nU)(x,y)=\frac12\bigl[
 (U(x)a_n(x))c_n(y)^{\mathsf T}
       +a_n(x)(U(y)c_n(y))^{\mathsf T}\bigr].
\]
The other triangle is its transpose. We compute the two square terms before treating their cross product. Reflect both variables by $r=1-x$ and $s=1-y$, so the triangle $x<y$
becomes $s<r$. Write $\widehat U_k(v)=U(1-(k+v)/n)$, $0<v<1$, for the field on the reflected cell with index $k$. For two different cells with indices $l<k$, the first square term
on this triangle is
\[
 \frac14|\widehat U_k(v)a_{n,k}|^2|c_{n,l}|^2.
\]
The transpose triangle has the same squared norm. Integrating over the cell of $c_{n,l}$, summing $l<k$, and then using $\dd r=\dd v/n$ gives the contribution
\[
 \frac1{2n}\sum_{k\in\mathcal K_n}A_{n,k}
                \int_0^1|\widehat U_k(v)a_{n,k}|^2\,\dd v
\]
to $\|\mathscr D_nU\|_{\mathscr K}^2$. The other square term is treated in the opposite variable and contributes
\[
 \frac1{2n}\sum_{k\in\mathcal K_n}C_{n,k}
                \int_0^1|\widehat U_k(v)c_{n,k}|^2\,\dd v.
\]
These sums exclude the square in which both variables belong to the same cell; that square must be integrated separately. The cross product is also retained, rather than absorbed
into an inequality.

For $\widehat U_k(v)=\sum_{\mu=1}^2w_{k,\mu}(v)\mathsf A_\mu$,
\[
 |\widehat U_k(v)b|^2
 =\sum_{\mu,\nu=1}^2w_{k,\mu}(v)w_{k,\nu}(v)
                  \langle\mathsf A_\mu b,\mathsf A_\nu b\rangle.
\]
Thus multiplication of the two displayed contributions by $n$ gives precisely the contraction with $\tfrac12(A_{n,k}\mathbf T(a_{n,k})+C_{n,k}\mathbf T(c_{n,k}))$ defined below.
The factor $1/2$ is the product of the kernel square factor $1/4$ and the two transpose triangles, not a covariance normalization. The integrals over squares whose two variables
lie in the same cell and the cross products are the remaining terms in the bilinear identity.

For $v\in\mathbb R^2$ define
\[
 \mathbf T(v)^{\mu\nu}
  =\langle\mathsf A_\mu v,\mathsf A_\nu v\rangle,\qquad
 \mathbf d(v)_\mu=\langle\mathsf A_\mu v,v\rangle,
 \qquad \mu,\nu\in\{1,2\}.
\]
Multiplication by $\mathsf A_1=-E_{12}$ and $\mathsf A_2=-H_0$ gives
\begin{equation}\label{S:eq:s12-affine-contractions}
 \mathbf T(v)=
 \begin{pmatrix}
 v_2^2&v_1v_2/2\\v_1v_2/2&|v|^2/4
 \end{pmatrix},\qquad
 \mathbf d(v)=\begin{pmatrix}-v_1v_2\\(v_2^2-v_1^2)/2\end{pmatrix},
 \qquad |\mathbf d(v)|^2=\frac{|v|^4}{4}.
\end{equation}
For a real matrix $H$ and a column $w$, the contraction convention is $H:(ww^{\mathsf T})=w^{\mathsf T}Hw$. In the notation of \eqref{S:eq:qv-matrix-H-definition},
\begin{equation}\label{S:eq:trace-matrix-H-definition}
 \mathbf H_{n,k}
  =\frac12\bigl(A_{n,k}\mathbf T(a_{n,k})
                  +C_{n,k}\mathbf T(c_{n,k})\bigr),\qquad
 \Gamma_{n,k}=\operatorname{tr}\mathbf H_{n,k}.
\end{equation}
Since $\operatorname{tr}\mathbf T(v)=\mathfrak h_{\mathrm{aff}}(v)$,
\[
 \tfrac18\mathfrak G_{n,k}\le\Gamma_{n,k}
          \le\tfrac58\mathfrak G_{n,k}.
\]
The entries can also be expressed in the upper half-plane coordinates of the sample. Write $d_{n,k}=q_\infty-x_{\vartheta_{n,k+1}}$ in the next display only. Here
$a=(-\sqrt p,0)^{\mathsf T}$ and $c=(d\sqrt p,p^{-1/2})^{\mathsf T}$ at the sample. Thus
\begin{equation}\label{S:eq:s12-matrix-H-entries}
 \begin{aligned}
 \mathbf H_{n,k}
 &=\begin{pmatrix}
 C_{n,k}/(2p_{n,k})&C_{n,k}d_{n,k}/4\\
 C_{n,k}d_{n,k}/4&(A_{n,k}p_{n,k}+C_{n,k}q_{n,k})/8
 \end{pmatrix},\\
 \Gamma_{n,k}
 &=\frac18A_{n,k}p_{n,k}+\frac18C_{n,k}q_{n,k}
                         +\frac{C_{n,k}}{2p_{n,k}}.
 \end{aligned}
\end{equation}
These coefficients use the fixed matrix factor $\mathcal A$. In particular, the term $C_{n,k}/(2p_{n,k})$ comes from the horizontal Brownian direction; omitting that direction
would change the leading order coefficient.

The remaining Gaussian calculation is local in the sampling index. Two operator intervals have a nonzero bridge covariance only if they contain values of the same attached bridge.
Put $u_{k,j}(s)=u_{n,j}(\theta_{n,k}(s))$ and define
\begin{equation}\label{S:eq:trace-two-point-covariance}
 \mathsf b_{n;k,l}(s,t)
 =\sum_{j\in\mathcal J_n^{\mathrm{at}}}h_{n,j}
 \bigl(\min\{u_{k,j}(s),u_{l,j}(t)\}
                       -u_{k,j}(s)u_{l,j}(t)\bigr).
\end{equation}
\begin{lemma}\label{S:lem:trace-two-point-covariance}
For every selected row,
\[
 \mathbb E_n[\widehat W_{n,k}(s)\widehat W_{n,l}(t)^{\mathsf T}]
       =\mathsf b_{n;k,l}(s,t)I_2.
\]
This covariance is zero for $|k-l|>2R_\beta$.
\end{lemma}
\begin{proof}
Work under the conditional law given $\mathscr S_n$. For each sampling interval $j$, let $\boldsymbol\beta_j$ be its two-dimensional standard bridge on $[0,1]$. Its coordinates are
independent, and different $j$ give independent bridges. The normalized coordinate of the sample $\vartheta_{n,k+1}$ in any sampling interval is an endpoint, so its bridge value is
zero. The field at an operator position therefore has the representation
\[
 \widehat W_{n,k}(s)=\sum_j\sqrt{h_{n,j}}\,
                         \boldsymbol\beta_j(u_{k,j}(s)).
\]
Coordinates outside a sampling interval are clipped to zero or one; both bridge values vanish. There are only finitely many nonzero terms at any position. For coordinate indices
$\mu,\nu\in\{1,2\}$, independence and centering give
\[
 \begin{aligned}
 \mathbb E_n[\widehat W_{n,k}^{\mu}(s)\widehat W_{n,l}^{\nu}(t)]
 &=\sum_{j,j'}\sqrt{h_{n,j}h_{n,j'}}\,
   \mathbb E_n[\beta_j^\mu(u_{k,j}(s))\beta_{j'}^\nu(u_{l,j'}(t))]\\
 &=\delta_{\mu\nu}\sum_jh_{n,j}
    \bigl(\min\{u_{k,j}(s),u_{l,j}(t)\}
                   -u_{k,j}(s)u_{l,j}(t)\bigr).
 \end{aligned}
\]
This is the stated covariance. A nonzero term requires both $|j-k|\le R_\beta$ and $|j-l|\le R_\beta$. The triangle inequality excludes such a $j$ when $|k-l|>2R_\beta$. On the
diagonal, let $u_0\in\{0,1\}$ be the normalized sample coordinate and set $\lambda=u-u_0$. For $u_0=0$, $u(1-u)=\lambda-\lambda^2$; for $u_0=1$ it is $-\lambda-\lambda^2$. Hence
$u(1-u)=|\lambda|-\lambda^2$ in both cases, giving
\begin{equation}\label{S:eq:trace-diagonal-covariance}
 \Sigma_{n,k}(s)=\mathsf b_{n;k,k}(s,s)I_2
 =\sum_jh_{n,j}\bigl(|\lambda_{n,k,j}(s)|-\lambda_{n,k,j}(s)^2\bigr)I_2.
\end{equation}

For the first $M$ Fourier modes, the profiles of \eqref{S:eq:gaussian-mode-profiles} give
\[
 \mathsf b_{n;k,l}^{(M)}(s,t)
   =\sum_jh_{n,j}\sum_{r=1}^M
                 q_{n,k,j,r}(s)q_{n,l,j,r}(t),\qquad
 \Sigma_{n,k}^{(M)}(s)=\mathsf b_{n;k,k}^{(M)}(s,s)I_2.
\]
Here $q_{n,k,j,r}(s)=\sqrt2\sin(\pi r u_{k,j}(s))/(\pi r)$, because $\sin(\pi r u_{n,j}(\vartheta_{n,k+1}))=0$. The scalar bridge covariance series is absolutely and uniformly
convergent. $|q_{n,k,j,r}(s)q_{n,l,j,r}(t)|\le2/(\pi^2r^2)$ and $\sum_{r>M}r^{-2}\le\int_M^\infty x^{-2}\dd x=M^{-1}$, so
\begin{equation}\label{S:eq:trace-covariance-tail-bound}
 |\mathsf b_{n;k,l}(s,t)-\mathsf b_{n;k,l}^{(M)}(s,t)|
 \le\frac{2}{\pi^2M}
 \sum_{\substack{j\in\mathcal J_n^{\mathrm{at}}\\
                 |j-k|\le R_\beta,\ |j-l|\le R_\beta}}h_{n,j}.
\end{equation}
On the diagonal the omitted terms are squares and hence nonnegative. For two different arguments, the omitted terms need not be nonnegative.
\end{proof}

Every product of linear interpolation terms is $\mathscr S_n$-measurable and therefore remains fixed under the conditioning. We use the bilinear forms of
Section~\ref{S:sec:stability-and-fr-echet-expansion-of-the-inverse-kernel} and decompose the linear interpolant according to its sampling interval:

\[
 \mathcal B_n=\sum_{i\in\mathcal J_n^{\mathrm{at}}}\mathcal B_n^{[i]},
 \mathcal B_{n,k}^{[i]}(s)=\lambda_{n,k,i}(s)\mathsf A(\Delta_{n,i}),
 \mathsf A(v)=\sum_{\mu=1}^2v_\mu\mathsf A_\mu.
\]
In particular,
\begin{equation}\label{S:eq:trace-interval-coefficients}
 r^a_{k,i}=\left(\int_0^1\lambda_{n,k,i}(s)\dd s\right)
                 \mathbf d(a_{n,k})\mathbin{\cdot}\Delta_{n,i},
 r^c_{k,i}=\left(\int_0^1\lambda_{n,k,i}(s)\dd s\right)
                 \mathbf d(c_{n,k})\mathbin{\cdot}\Delta_{n,i}.
\end{equation}
As in Section~\ref{S:sec:the-error-from-sampling-the-brownian-path}, the normal increment on sampling interval $i$ is expressed in the orthonormal frame $(e_1,e_2)$ at
$Z_{\vartheta_{n,i+1}}$. For a point of $I_{n,k}$, the comparison below freezes its coefficient at this starting point and estimates the resulting error.

We separate the terms containing the squared norm of each differentiated column from the cross terms and the integrals over a single partition square. The contribution from the
current interval is the integral over $I_{n,k}\times I_{n,k}$. The bridge cross term is the conditional mean of $\mathsf O_n$:

\begin{equation}\label{S:eq:trace-seven-square-bridge}
 E_n^{\mathrm{diag}}:=\sum_{k\in\mathcal K_n} \bigl(\mathsf S_{n,k}(\mathcal B_n,\mathcal B_n) +\mathbb E_n\mathsf S_{n,k}(\mathcal Z_n,\mathcal Z_n)\bigr),\qquad E_n^{\mathrm{bridge}}:=\mathbb E_n\mathsf O_n(\mathcal Z_n,\mathcal Z_n).
\end{equation}
The contribution of the linear interpolant to $\mathsf O_n$ is split according to the distance between its two sampling indices. For each index, only a bounded number of indices
contribute to the nearby terms. For the separated terms we use the martingale identity proved below under the law conditioned on the terminal point.
\begin{equation}\label{S:eq:trace-seven-chord}
 \begin{aligned}
 E_n^{\mathrm{chord,near}}
  &:=\frac1n\sum_{\substack{k,l\in\mathcal K_n\\k<l}}
           \sum_{\substack{i,j\in\mathcal J_n^{\mathrm{at}}\\|i-j|<L_\beta}}
                         r^c_{k,i}r^a_{l,j},\\
 E_n^{\mathrm{chord,far}}
  &:=\frac1n\sum_{\substack{k,l\in\mathcal K_n\\k<l}}
           \sum_{\substack{i,j\in\mathcal J_n^{\mathrm{at}}\\|i-j|\ge L_\beta}}
                         r^c_{k,i}r^a_{l,j}.
 \end{aligned}
\end{equation}
The conditional second moment of the linear interpolant plus the Brownian bridge at each point is $\mathsf V_{n,k}(s)=m_{n,k}(s)m_{n,k}(s)^{\mathsf T}+\Sigma_{n,k}(s)$. For a
traversed fraction $\lambda$, the bridge supplies $h(|\lambda|-\lambda^2)I_2$. Weighted quadratic variation gives the contribution from linear interpolation $h\lambda^2I_2$. Thus
we compare $\mathsf V_{n,k}(s)$ with $\ell_{n,k}(s)I_2$ only after contraction and summation. The last two replacements concern the deterministic time change:

\begin{equation}\label{S:eq:trace-qv-and-time-corrections}
 \begin{aligned}
 E_n^{\mathrm{qv}}
  &:=\sum_{k\in\mathcal K_n}\int_0^1\mathbf H_{n,k}:
              \bigl(\mathsf V_{n,k}(s)-\ell_{n,k}(s)I_2\bigr)\dd s,\\
 E_n^{\mathrm{interp}}
  &:=\sum_{k\in\mathcal K_n}\Gamma_{n,k}
       \left(\int_0^1\ell_{n,k}(s)\dd s-\tau_{n,k}\kappa_{n,k}\right),\\
 E_n^{\mathrm{sampling}}
  &:=\sum_{k\in\mathcal K_n}\Gamma_{n,k}
       \left(\tau_{n,k}\kappa_{n,k}-\frac4{\beta k}\kappa_\beta(r_k)\right).
 \end{aligned}
\end{equation}
The local time difference has size $4/(\beta k)$, up to its sampling offset. Since the normalized squared kernel integral assigns length $1/n$ to each remaining operator variable,
$k^{-1}=n^{-1}r_k^{-1}$ turns its main term into a Riemann sum. We keep the individual corrections to distinguish the estimates that have rates from the weighted quadratic
variation limit. The next identity is exact for each row.

\begin{proposition}\label{S:prop:trace-seven-identity}
Almost surely, for each selected row,
\begin{equation}\label{S:eq:trace-error-decomposition}
 \begin{aligned}
 \mathfrak C_n={}&n\operatorname{tr}Q_n+n\|B_n\|_{\mathscr K}^2\\
 ={}&\frac4\beta\frac1n\sum_{k=m_n}^{n-1}
       \kappa_\beta(r_k)\frac{\Gamma_{n,k}}{r_k}
    +E_n^{\mathrm{diag}}+E_n^{\mathrm{bridge}}
    +E_n^{\mathrm{chord,near}}+E_n^{\mathrm{chord,far}}\\
 & +E_n^{\mathrm{qv}}+E_n^{\mathrm{interp}}+E_n^{\mathrm{sampling}}.
 \end{aligned}
\end{equation}
\end{proposition}
\begin{proof}
We first use $\mathcal Z_n^{(M)}$ and $Q_n^{(M)}$. Only $E_n^{\mathrm{diag}}$, $E_n^{\mathrm{bridge}}$ and $E_n^{\mathrm{qv}}$ then change. Write
$\mathfrak C_n^{(M)}=n\operatorname{tr}Q_n^{(M)}+n\|B_n\|_{\mathscr K}^2$ and
\[
 t_{n,M}^{+}:=\mathbb E_n\|
       \mathcal Z_n-\mathcal Z_n^{(M)}\|_{n,+}^2
       \le Ce_n^{\mathrm{KL}}/M.
\]

Fix $\mathscr S_n$ and $M<\infty$. Every mode field belongs to $\mathfrak H_n^+$. The trace identity for the truncated Gaussian series and the common minus signs give
\[
 \mathfrak C_n^{(M)}
 =n\mathbb E_n\|\mathscr D_n\mathcal Z_n^{(M)}\|_{\mathscr K}^2
          +n\|\mathscr D_n\mathcal B_n\|_{\mathscr K}^2.
\]
Under $\mathbb E_n$, the field $\mathcal B_n$ is fixed and the Gaussian coordinates are centered. Hence the mixed bilinear terms have mean zero.

Apply Proposition~\ref{S:prop:local-trace-aux} to each quadratic expression. The squared column variations integrated over pairs of distinct cells give
\begin{equation}\label{S:eq:trace-strip-conditional-identity}
 \mathsf L_n(\mathcal B_n,\mathcal B_n)
   +\mathbb E_n\mathsf L_n(\mathcal Z_n^{(M)},\mathcal Z_n^{(M)})
 =\sum_k\int_0^1\mathbf H_{n,k}:
       \bigl(m_{n,k}m_{n,k}^{\mathsf T}+\Sigma_{n,k}^{(M)}\bigr)\dd s.
\end{equation}
The other terms in Proposition~\ref{S:prop:local-trace-aux} satisfy
\[
 \begin{aligned}
 \sum_k\bigl[\mathsf S_{n,k}(\mathcal B_n,\mathcal B_n)
       +\mathbb E_n\mathsf S_{n,k}(\mathcal Z_n^{(M)},\mathcal Z_n^{(M)})\bigr]
   &=E_n^{\mathrm{diag},(M)},\\
 \mathbb E_n\mathsf O_n(\mathcal Z_n^{(M)},\mathcal Z_n^{(M)})
   &=E_n^{\mathrm{bridge},(M)},\\
 \mathsf O_n(\mathcal B_n,\mathcal B_n)
   &=\frac1n\sum_{k<l}
       \left(\sum_i r^c_{k,i}\right)\left(\sum_j r^a_{l,j}\right)
     =E_n^{\mathrm{chord,near}}+E_n^{\mathrm{chord,far}}.
 \end{aligned}
\]
The last equality partitions every ordered pair $(i,j)$ into $|i-j|<L_\beta$ or $|i-j|\ge L_\beta$. Since $\mathbf H_{n,k}:I_2=\Gamma_{n,k}$, the same contribution equals
\[
 \begin{aligned}
 &\sum_k\int_0^1\mathbf H_{n,k}:
       (m_{n,k}m_{n,k}^{\mathsf T}+\Sigma_{n,k}^{(M)})\,\dd s\\
 &=E_n^{\mathrm{qv},(M)}+\sum_k\Gamma_{n,k}\int_0^1\ell_{n,k}(s)\,\dd s\\
 &=E_n^{\mathrm{qv},(M)}+E_n^{\mathrm{interp}}
      +\sum_k\Gamma_{n,k}\tau_{n,k}\kappa_{n,k}\\
 &=E_n^{\mathrm{qv},(M)}+E_n^{\mathrm{interp}}+E_n^{\mathrm{sampling}}
      +\frac4\beta\frac1n\sum_k\frac{\kappa_\beta(r_k)\Gamma_{n,k}}{r_k}.
 \end{aligned}
\]
Here $r_k=k/n$, and $E_n^{\mathrm{qv},(M)}$ subtracts the entire matrix $\ell_{n,k}I_2$ specified in \eqref{S:eq:trace-qv-and-time-corrections}. Adding the displayed identities
proves \eqref{S:eq:trace-error-decomposition} at finite $M$.

Let $M\to\infty$ with $n$ fixed, and put $H_M=\mathcal Z_n-\mathcal Z_n^{(M)}$. For each of the bilinear forms $\mathfrak b=\mathsf L_n,\sum_k\mathsf S_{n,k},\mathsf O_n$,
Proposition~\ref{S:prop:local-trace-aux} gives $|\mathfrak b(F,G)|\le C\|F\|_{n,+}\|G\|_{n,+}$. Under $\mathbb E_n$, $H_M$ is centered and independent of the retained modes.
Conditional Cauchy--Schwarz gives
\[
 \mathbb E_n|\mathfrak b(\mathcal Z_n^{(M)},H_M)|
 \le C\bigl(\mathbb E_n\|\mathcal Z_n^{(M)}\|_{n,+}^2\bigr)^{1/2}
          (t_{n,M}^+)^{1/2}<\infty.
\]
Conditioning further on the retained Gaussian coordinates and using $\mathbb E_n[H_M\mid\xi_{n,j,m,\nu}:m\le M]=0$ therefore yields
\[
 \mathbb E_n\mathfrak b(\mathcal Z_n^{(M)},H_M)
 =\mathbb E_n\mathfrak b(H_M,\mathcal Z_n^{(M)})=0.
\]
Expansion of the bilinear form now gives
\[
 \left|\mathbb E_n\mathfrak b(\mathcal Z_n,\mathcal Z_n)
       -\mathbb E_n\mathfrak b(\mathcal Z_n^{(M)},\mathcal Z_n^{(M)})\right|
 =|\mathbb E_n\mathfrak b(H_M,H_M)|\le Ct_{n,M}^+.
\]
Applying this calculation to $\mathfrak b(F,G)=n\langle\mathscr D_nF,\mathscr D_nG\rangle$ gives the same bound for $\mathfrak C_n-\mathfrak C_n^{(M)}$. In
$E_n^{\mathrm{qv}}-E_n^{\mathrm{qv},(M)}$, the linear interpolation term and $\ell I_2$ cancel exactly, leaving $\mathbb E_n\mathsf L_n(H_M,H_M)$. Adding the four bounds proves
\eqref{S:eq:trace-mode-conditional-tail}. The differences $\mathfrak C_n-\mathfrak C_n^{(M)}$, $E_n^{\mathrm{diag}}-E_n^{\mathrm{diag},(M)}$ and
$E_n^{\mathrm{qv}}-E_n^{\mathrm{qv},(M)}$ are nonnegative. The estimate for the off-diagonal difference allows either sign. Conditional $L^2$ convergence from
Section~\ref{S:sec:the-conditional-gaussian-operator} and these bounds identify the limits with the definitions above, proving \eqref{S:eq:trace-error-decomposition}.

For each fixed row, $\mathcal E_n^+(\mathcal B_n)<\infty$ and $\mathbb E_n\mathcal E_n^+(\mathcal Z_n)<\infty$ by Proposition~\ref{S:prop:weighted-field-covariance}. The same
bilinear bounds justify conditional Fubini and all signed conditional integrals. Measurability follows either from the expressions for the truncated Gaussian series and their
limits or from Lemma~\ref{S:lem:trace-two-point-covariance}. Taking a countable intersection over $n$ gives a common set of full measure.

The fixed-$n$ limit is controlled by
\begin{equation}\label{S:eq:trace-mode-conditional-tail}
 \begin{aligned}
 &|\mathfrak C_n-\mathfrak C_n^{(M)}|
  +|E_n^{\mathrm{diag}}-E_n^{\mathrm{diag},(M)}|
  +|E_n^{\mathrm{bridge}}-E_n^{\mathrm{bridge},(M)}|
  +|E_n^{\mathrm{qv}}-E_n^{\mathrm{qv},(M)}|\\
 &\le C t_{n,M}^{+}
                     \le Ce_n^{\mathrm{KL}}/M.
 \end{aligned}
\end{equation}
\end{proof}

\subsubsection{Terms involving nearby intervals}
\label{S:sec:s12-nearby-intervals}

The bridge covariance has finite range in the operator index. The products of linear interpolation terms with nearby indices also have only a bounded number of possible second
indices for each first index. Both can therefore be bounded by absolute values after their index ranges have been identified. The separated sum of products of linear interpolation
terms will be estimated by summation by parts.

The integrals over $I_{n,k}\times I_{n,k}$ and the covariance between distinct subintervals can also be written explicitly. Put
\[
 \mathbf V_{n;k,l}(s,t):=m_{n,k}(s)m_{n,l}(t)^{\mathsf T}
                         +\mathsf b_{n;k,l}(s,t)I_2.
\]
Then direct expansion of the two triangles in each diagonal square gives
\begin{equation}\label{S:eq:trace-square-explicit-covariance}
 \begin{aligned}
 E_n^{\mathrm{diag}}=\frac1{2n}\sum_k\int_0^1\int_0^s
 \Bigl[&q_{n,k}\mathbf T(a_{n,k}):\mathbf V_{n;k,k}(s,s)
       +p_{n,k}\mathbf T(c_{n,k}):\mathbf V_{n;k,k}(t,t)\\
       &+2\mathbf d(a_{n,k})^{\mathsf T}\mathbf V_{n;k,k}(s,t)
                              \mathbf d(c_{n,k})\Bigr]\dd t\,\dd s,
 \end{aligned}
\end{equation}
and
\begin{equation}\label{S:eq:trace-bridge-explicit-covariance}
 E_n^{\mathrm{bridge}}
 =\frac1n\sum_{k<l}\mathbf d(c_{n,k})\mathbin{\cdot}\mathbf d(a_{n,l})
             \int_0^1\!\int_0^1\mathsf b_{n;k,l}(s,t)\dd s\,\dd t.
\end{equation}
Indeed, the squared derivative on $s>t$ has the three displayed terms with coefficient $1/4$; the transposed triangle has the same integral. The Jacobian is $n^{-2}$ and the trace
normalization is $n$, giving $1/(2n)$ in \eqref{S:eq:trace-square-explicit-covariance}. The off-diagonal scalar quadratic form has coefficient $1/n$, giving
\eqref{S:eq:trace-bridge-explicit-covariance}. These formulas also follow first for finite modes and then by \eqref{S:eq:trace-mode-conditional-tail}. The cross term in each
diagonal square uses the two point covariance.

\begin{proposition}\label{S:prop:trace-corrections}
On the common event of probability one of the preceding estimates, for all sufficiently large $n$,
\begin{equation}\label{S:eq:trace-local-correction-rates}
 \begin{aligned}
 0\le E_n^{\mathrm{diag}}
    &\le C_{\beta,\omega}\frac{(\log n)^{9/2}}n,&
 |E_n^{\mathrm{bridge}}|
    &\le C_{\beta,\omega}\frac{(\log n)^{7/2}}n,\\
 |E_n^{\mathrm{chord,near}}|
    &\le C_{\beta,\omega}\frac{(\log n)^{9/2}}n,&
 E_n^{\mathrm{qv}}&\longrightarrow0.
 \end{aligned}
\end{equation}
\end{proposition}
\begin{proof}
  Corollary~\ref{S:cor:normal-increment-full-array} and the finite interaction range give, uniformly for $k\in\mathcal
  K_n$,
\begin{equation}\label{S:eq:trace-subinterval-bounds}
 \sup_s|m_{n,k}(s)|^2\le C_{\beta,\omega}\frac{\log n}{k},
 \sup_s\operatorname{tr}\Sigma_{n,k}(s)\le C_\beta/k,
   p_{n,k}q_{n,k}\le C_{\beta,\omega}(\log n)^{5/2}.
\end{equation}
For the first inequality, at most $2R_\beta+1$ increments occur, their squared magnitudes are at most $C_\beta h_{n,j}\log n$, and $h_{n,j}\le C_\beta/k$ when $|j-k|\le R_\beta$.
For the second use Lemma~\ref{S:lem:trace-two-point-covariance} and $u(1-u)\le1/4$. The last is \eqref{S:eq:s7-column-product-bound} for the sampled columns.

For any fixed integer $D\ge1$, there is $n_0(\omega,D)<\infty$ such that the sampled columns satisfy the following bounds simultaneously for all $n\ge n_0(\omega,D)$ and
$k,l\ge\ell_n$ with $|k-l|\le D$:
\begin{equation}\label{S:eq:trace-neighbor-column-products}
 p_{n,l}\le4p_{n,k},\qquad q_{n,l}\le4q_{n,k},\qquad
 q_{n,k}p_{n,l}\le C_{\beta,\omega,D}(\log n)^{5/2}.
\end{equation}

Let $s_j=\vartheta_{n,j+1}$ and $t_j=\vartheta_{n,j}$. The isometry $w\mapsto(w-\operatorname{Re}Z_{s_j})/\operatorname{Im}Z_{s_j}$ and the expansion
\eqref{S:eq:estimate-normal-affine-taylor} give, with the normal increment $\Delta_{n,j}$ expressed at $Z_{s_j}$,
\[
 X_{t_j}X_{s_j}^{-1}
   =\exp\{\mathsf A(\Delta_{n,j})+\mathcal Q_{Z_{s_j}}(\Delta_{n,j})\},
 \qquad \|\mathcal Q_{Z_{s_j}}(\Delta_{n,j})\|
                                      \le C|\Delta_{n,j}|^2.
\]
The simultaneous estimate for all increments implies $|\Delta_{n,j}|\le C_\beta\sqrt{\log n/j}$ for all these indices, with probability one for all sufficiently large $n$.
Therefore
\[
 \|X_{t_j}X_{s_j}^{-1}-I\|
 +\|X_{s_j}X_{t_j}^{-1}-I\|
       \le C_\beta\sqrt{\log n/j}.
\]
For at most $D$ consecutive factors $R_i$, use
\[
 \|R_1\cdots R_d-I\|
 \le\prod_{i=1}^d(1+\|R_i-I\|)-1,
 \qquad d\le D.
\]
The indices are comparable, and the right side is at most $C_{\beta,D}\sqrt{\log n/\ell_n}=o(1)$. Applied to either boundary vector, this proves the first two inequalities in both
directions. The identity \eqref{S:eq:s7-column-product-bound}, evaluated also at the additional sampling indices, gives the last inequality.

For the term from $I_{n,k}\times I_{n,k}$, the two transposed triangle integrals in $\mathsf S_{n,k}(F,F)$ and $|A+B|^2\le2|A|^2+2|B|^2$ give
\begin{equation}\label{S:eq:trace-same-square-majorant}
 0\le\mathsf S_{n,k}(F,F)
 \le\frac1n\int_0^1
    \bigl[sq_{n,k}|F_k(s)a_{n,k}|^2
       +(1-s)p_{n,k}|F_k(s)c_{n,k}|^2\bigr]\dd s.
\end{equation}
On the triangle $0<t<s<1$ the growing column square is integrated over a set of $t$-length $s$; the square involving $c$ gives the weight $1-s$ after interchanging the dummy
variables. For a varying field, the factors $s$ and $1-s$ remain inside the integrals. For $F=\mathsf A(y)$, the right side is at most $(C/n)p_{n,k}q_{n,k}\int_0^1|y(s)|^2\dd s$.
Apply this first to the linear interpolant and then under $\mathbb E_n$ to the Brownian bridge. Equations \eqref{S:eq:trace-subinterval-bounds} give
\[
 E_n^{\mathrm{diag}}
 \le\frac{C_{\beta,\omega}(\log n)^{7/2}}n
                 \sum_{k=m_n}^{n-1}\frac1k.
\]
Since the harmonic sum is at most $1+\log(n/m_n)\le C\log n$, the first estimate follows.

For the bridge term with distinct interval indices, the scalar contractions satisfy
\[
 \|r_k^a(\mathcal Z_n)\|_{L^2(\mathbb E_n)}\le C_\beta p_{n,k}k^{-1/2},
 \|r_k^c(\mathcal Z_n)\|_{L^2(\mathbb E_n)}\le C_\beta q_{n,k}k^{-1/2}.
\]
These follow by Minkowski in $s$ and the second inequality in \eqref{S:eq:trace-subinterval-bounds}. Lemma~\ref{S:lem:trace-two-point-covariance}, or independence of the bridges on
distinct intervals under the conditional law, makes the cross expectation zero when $|k-l|>2R_\beta$. Thus conditional Cauchy--Schwarz, with
\eqref{S:eq:trace-neighbor-column-products}, gives
\[
 |E_n^{\mathrm{bridge}}|
 \le\frac{C_\beta}n\sum_{k}\sum_{k<l\le k+2R_\beta}
       \frac{q_{n,k}p_{n,l}}{\sqrt{kl}}
 \le\frac{C_{\beta,\omega}(\log n)^{5/2}}n
                          \sum_k\frac1k.
\]

For the products of linear interpolation terms with nearby indices, \eqref{S:eq:trace-interval-coefficients} implies $|r^a_{k,i}|\le Cp_{n,k}|\Delta_{n,i}|$ and
$|r^c_{k,i}|\le Cq_{n,k}|\Delta_{n,i}|$. A nonzero term has $|i-k|,|j-l|\le R_\beta$ and $|i-j|<L_\beta$, hence $0<l-k<L_\beta+2R_\beta$. For each $k$ there are only a
deterministic number of such $(l,i,j)$, including the possibility $i=j$. The increment bound and neighboring comparison give
\[
 |E_n^{\mathrm{chord,near}}|
 \le\frac{C_{\beta,\omega}\log n}n
             \sum_k\frac{p_{n,k}q_{n,k}}k
 \le\frac{C_{\beta,\omega}(\log n)^{7/2}}n\sum_k\frac1k.
\]
This bound includes the pairs $i=j$.

Apply \eqref{S:eq:qv-sum-over-increments} with $F_{n,k}(s)=\mathbf H_{n,k}$. It gives the exact identity
\begin{equation}\label{S:eq:s12-qv-collected}
 E_n^{\mathrm{qv}} ={}\sum_{j\in\mathcal J_n^{\mathrm{at}}} \mathbf F_{n,j}^{(0)}: (\Delta_{n,j}\Delta_{n,j}^{\mathsf T}-h_{n,j}I_2) +\sum_{\substack{q\in\mathbb Z\\0<|q|\le2R_\beta}} \sum_{\substack{j,j+q\in\mathcal J_n^{\mathrm{at}}}} \mathbf F_{n,j}^{(q)}: \Delta_{n,j}\Delta_{n,j+q}^{\mathsf T},
 \end{equation}
where
\[
 \mathbf F_{n,j}^{(q)}
   =\sum_{k\in\mathcal K_n}\int_0^1
       \lambda_{n,k,j}(s)\lambda_{n,k,j+q}(s)
                             \mathbf H_{n,k}\,\dd s.
\]
Corollary~\ref{S:cor:qv-sampled-green-weights} verifies \eqref{S:eq:qv-mesh-weight-hypotheses} for $\mathbf H_{n,k}$ and every $0<\zeta<a_\beta$.
Proposition~\ref{S:prop:bridge-and-interpolation-second-moment} therefore gives
\[
 E_n^{\mathrm{qv}}=\mathcal V_n(\mathbf H)\longrightarrow0
 \quad\hbox{on }\Omega_{\mathrm{qv}}.
\]
That proposition includes both signs of the shift and the additional sampling intervals at the boundary points. The convergence is pathwise; these hypotheses do not supply a
quantitative rate for it.

The time replacements obey
\begin{equation}\label{S:eq:trace-time-correction-rates}
 |E_n^{\mathrm{interp}}| \le\frac1{3\beta}\sum_{k=m_n}^{n-1}\frac{\Gamma_{n,k}}{k^2},\qquad |E_n^{\mathrm{sampling}}| \le C_\beta\sum_{k=m_n}^{n-1}\frac{\Gamma_{n,k}}{k^2}, \sum_{k=m_n}^{n-1}\frac{\Gamma_{n,k}}{k^2}\le C_{\beta,\omega}/m_n.
\end{equation}

Finally \eqref{S:eq:time-interpolation-error} gives the first inequality in \eqref{S:eq:trace-time-correction-rates} directly. Equation~\eqref{S:eq:sampling-variance-two-errors}
gives the second bound using the absolute discrete offset. The coefficient $\Gamma_{n,k}$ is nonnegative and bounded by $(5/8)\mathfrak G_{n,k}$, so
\[
 \sum_{k=m_n}^{n-1}\frac{\Gamma_{n,k}}{k^2}
 \le\frac{5}{8m_n}\sum_{k=m_n}^{n-1}\frac{\mathfrak G_{n,k}}k
 \le C_{\beta,\omega}/m_n
\]
by Proposition~\ref{S:prop:5.5}.
\end{proof}

\subsubsection{The sum over separated sampling intervals}
\label{S:sec:s12-separated-intervals}

The separated part is the only cross term with an unbounded number of partners. The order of the two kernel variables is useful here: after collecting contributions from one
sampling interval, it turns the double sum into a product with a partial sum involving the column associated with the terminal boundary point. We establish this factorization
before estimating either factor.

Put $\mathfrak a_{n,j}=\sum_{k\in\mathcal K_n}r^a_{k,j}$ and $\mathfrak c_{n,j}=\sum_{k\in\mathcal K_n}r^c_{k,j}$.

\begin{lemma}\label{S:lem:trace-separated-pair-factorization}
For every $j\ge i+L_\beta$,
\begin{equation}\label{S:eq:trace-far-pair-factor}
 n\langle\mathscr D_n\mathcal B_n^{[j]},\mathscr D_n\mathcal B_n^{[i]}\rangle_{\mathscr K}
       =\frac1{2n}\mathfrak a_{n,j}\mathfrak c_{n,i}.
\end{equation}
\end{lemma}
\begin{proof}
Fix $j\ge i+L_\beta$ and put $F=\mathcal B_n^{[j]}$, $G=\mathcal B_n^{[i]}$. The field $F$ is supported on reflected operator cells $l$ with $|l-j|\le R_\beta$, whereas $G$ is
supported on cells $k$ with $|k-i|\le R_\beta$. Every such pair satisfies
\[
 l-k\ge j-i-2R_\beta\ge L_\beta-2R_\beta>0.
\]
Their cell supports are disjoint. In particular, every term of the bilinear derivative formula involving both fields at a single position is zero. Terms with both positions in one
cell are zero for the same reason.

It remains to calculate the ordered cross term. On a reflected rectangle with $k<l$, let $s$ belong to cell $k$ and $r$ to cell $l$. The derivative for $F$ contains only the
variation of the first column, and that for $G$ contains only the variation of the second column. The Frobenius inner product on this triangle is therefore
\[
 \frac14\langle F(1-r)a_{n,l},a_{n,l}\rangle
                \langle c_{n,k},G(1-s)c_{n,k}\rangle.
\]
The transpose triangle supplies an identical integral. By the definitions of $r^a_{l,j}$ and $r^c_{k,i}$, integration in the two normalized cell coordinates thus gives
\[
 \langle\mathscr D_nF,\mathscr D_nG\rangle_{\mathscr K}
   =\frac1{2n^2}\sum_{k<l}r^c_{k,i}r^a_{l,j}.
\]
The coefficient is $2(1/2)^2n^{-2}$: two triangles, two derivative factors, and two cell lengths. The ordering $k<l$ is automatic whenever the product is nonzero. The finite sum
consequently factors as
\[
 \sum_{k<l}r^c_{k,i}r^a_{l,j}
 =\left(\sum_k r^c_{k,i}\right)\left(\sum_l r^a_{l,j}\right)
 =\mathfrak c_{n,i}\mathfrak a_{n,j}.
\]
Multiplication by $n$ proves \eqref{S:eq:trace-far-pair-factor}. This also proves the formula when either support is empty, since both sides then vanish.

Now sum over all pairs of sampling indices. The Hilbert--Schmidt scalar product is real and symmetric, so the orderings $i<j$ and $j<i$ contribute equally. Hence
\begin{equation}\label{S:eq:trace-far-pair-collection}
 \begin{aligned}
 E_n^{\mathrm{chord,far}}
  &=n\sum_{\substack{i,j\in\mathcal J_n^{\mathrm{at}}\\|i-j|\ge L_\beta}}
       \langle\mathscr D_n\mathcal B_n^{[i]},\mathscr D_n\mathcal B_n^{[j]}\rangle_{\mathscr K}\\
  &=\frac1n\sum_{j=\ell_n+L_\beta}^{n-1}
       \mathfrak a_{n,j}\sum_{i=\ell_n}^{j-L_\beta}\mathfrak c_{n,i}.
 \end{aligned}
\end{equation}
Every sum still uses the full sampling range $\mathcal J_n^{\mathrm{at}}$, including its additional intervals near the boundary. For nearby indices, the cell supports can overlap
and the terms obtained by multiplying variations of the same column in the bilinear derivative formula need not vanish.
\end{proof}

\begin{proposition}\label{S:prop:long-range-chord}
Almost surely,
\begin{equation}\label{S:eq:long-range-chord-vanishing}
 E_n^{\mathrm{chord,far}}
 =n\sum_{|j-l|\ge L_\beta}
       \langle\mathscr D_n\mathcal B_n^{[j]},
                    \mathscr D_n\mathcal B_n^{[l]}\rangle_{\mathscr K}
 \longrightarrow0.
\end{equation}
\end{proposition}

We estimate the inner partial sum first. Smaller sampling indices correspond to later Brownian times, so the inner sum contains increments in the future of the outer one. Under
$\mathbb P^\eta$, sums of the increments of $q$ equal its boundary point difference. This replaces the inner sum by boundary values before we estimate the outer factor.
Conditioning on $\mathscr S_n$ would already determine every increment used in the linear interpolant.

\paragraph{A martingale under conditioning on the terminal point}
\label{S:sec:s12-terminal-martingale}

Fix a finite terminal point $\eta$ outside the null set in Lemma~\ref{S:lem:terminal-affine-innovations}, and work under $\mathbb P^\eta$. Put
$c_\theta=X_\theta(\eta,1)^{\mathsf T}$, $q_\theta=|c_\theta|^2$, and write $\widetilde W=\widetilde W^\eta$. The moment bound on each interval uses the conditional distribution of
the diffusion after time $s$, given $\mathcal F_s^Z$, with the specified terminal point. The reciprocal Poisson kernel is harmonic for $L^\eta$. For a fixed boundary point $\eta$,
the same function $q(z;\eta)=((x-\eta)^2+y^2)/y$ satisfies $Lq=q$, whereas $L^\eta q=0$. The choice of conditional law removes this drift.

\begin{lemma}\label{S:lem:terminal-q-martingale}
Under $\mathbb P^\eta$, the positive process $q_\theta$ satisfies
\begin{equation}\label{S:eq:terminal-q-sde}
 \dd q_\theta=2\mathbf d(c_\theta)\cdot\dd\widetilde W_\theta.
\end{equation}
It is a martingale on every finite time interval.
\end{lemma}
\begin{proof}
Here $\eta$ is finite; this has full harmonic measure. The identity from Section~\ref{S:sec:products-of-the-coupled-transfer-matrices} is
\begin{equation}\label{S:eq:terminal-q-poisson-normalization}
 q_\theta=(1+\eta^2)P_\eta(Z_\theta)^{-1}.
\end{equation}
It also follows by substituting $c=(-(x-\eta)/\sqrt y,\sqrt y)^{\mathsf T}$ and $P_\eta(z)=(1+\eta^2)y/((x-\eta)^2+y^2)$. For $L^\eta f=P_\eta^{-1}L(P_\eta f)$,
\[
 L^\eta(P_\eta^{-1})=P_\eta^{-1}L1=0.
\]
Moreover
\[
 y\partial_xq=2(x-\eta)=2\mathbf d(c)_1,
 \qquad y\partial_yq=y-(x-\eta)^2/y=2\mathbf d(c)_2.
\]
It\^o's formula under the law conditioned on $\eta$ proves \eqref{S:eq:terminal-q-sde}. The last identity in \eqref{S:eq:s12-affine-contractions} gives
\[
 \dd[q]_\theta=q_\theta^2\dd\theta,\qquad
 N_t:=\int_0^t\frac{2\mathbf d(c_\theta)}{q_\theta}
                          \cdot\dd\widetilde W_\theta,
 \qquad [N]_t=t.
\]
The predictable integrand in $N$ has norm one. We verify the conditional exponential moments that are needed below, including their integrability. Fix $s<t$ and $r\in\mathbb R$,
and put, for $s\le u\le t$,
\[
 Z_u^{(r)}=\exp\{r(N_u-N_s)-r^2(u-s)/2\}.
\]
It\^o's formula gives $\dd Z_u^{(r)}=rZ_u^{(r)}\dd N_u$. Stop at $\tau_m=\inf\{u\ge s:Z_u^{(r)}\ge m\}\wedge t$. The stopped stochastic integral is square integrable. Applying
It\^o's formula to its square, conditionally on $\mathcal F_s^Z$, and then Gronwall's inequality gives
\[
 \mathbb E^\eta[(Z_{u\wedge\tau_m}^{(r)})^2\mid\mathcal F_s^Z]
 \le1+r^2\int_s^u
       \mathbb E^\eta[(Z_{v\wedge\tau_m}^{(r)})^2\mid\mathcal F_s^Z]\,\dd v
 \le\mathrm e^{r^2(u-s)}.
\]
The stopped exponentials have conditional mean one. The uniform second-moment bound makes them uniformly integrable, and continuity of $N$ on $[s,t]$ gives $\tau_m=t$ for all
sufficiently large $m$ on each path. Passing to the limit yields
\[
 \mathbb E^\eta[\mathrm e^{r(N_t-N_s)}\mid\mathcal F_s^Z]
 =\mathrm e^{r^2(t-s)/2}.
\]
For completeness, the analogous complex exponential $\exp\{\mathrm i\theta(N_u-N_s)+\theta^2(u-s)/2\}$ has bounded modulus on $[s,t]$ and zero drift. Its conditional expectation is
one, so the conditional characteristic function of $N_t-N_s$ is $\mathrm e^{-\theta^2(t-s)/2}$. This identifies a centered normal increment independent of $\mathcal F_s^Z$.
Together with continuity and $N_0=0$, it proves that $N$ is Brownian in this filtration.

Since $q>0$, It\^o's formula applied locally to $\log q$ gives $\dd\log q=\dd N-\dd u/2$. Integrating and using the exponential moments just established gives
\begin{equation}\label{S:eq:s12-terminal-finite-moments}
 \begin{aligned}
 q_t&=q_s\exp(N_t-N_s-(t-s)/2),\\
 \mathbb E^\eta[q_t^r\mid\mathcal F_s^Z]
   &=q_s^r\exp\{r(r-1)(t-s)/2\},\qquad r\in\mathbb R.
 \end{aligned}
\end{equation}
The multiplication by the positive $\mathcal F_s^Z$-measurable factor $q_s^r$ is first justified under the conditional law, or on the events $\{q_s^r\le L\}$ and then by monotone
convergence. At $s=0$, the fixed finite terminal value gives $q_0=1+\eta^2\in(0,\infty)$, so all these moments are finite on every finite time interval. At $r=1$ the formula proves
$\mathbb E^\eta[q_t\mid\mathcal F_s^Z]=q_s$ and $\mathbb E^\eta q_t=q_0$. Thus the local martingale is a martingale on every finite time interval.

For the sampled indices, define
\begin{equation}\label{S:eq:terminal-q-sampled-increment}
 M_{n,j}^c:=\int_{\vartheta_{n,j+1}}^{\vartheta_{n,j}}
       \mathbf d(c_\theta)\cdot\dd\widetilde W_\theta
       =\tfrac12(q_{n,j-1}-q_{n,j}).
\end{equation}
The index $j-1$ appears because $q_{n,j-1}=q_{\vartheta_{n,j}}$, whereas $q_{n,j}=q_{\vartheta_{n,j+1}}$.
\end{proof}

\paragraph{Coefficients evaluated at the start of each interval}

The exact martingale increment contains the coefficient along the whole sampling interval in Brownian time. The linear interpolation term uses a coefficient evaluated at a sample.
We compare the two at the starting time $\vartheta_{n,j+1}$, where the coefficient is measurable before the increment is observed. Define
\[
 \widehat\xi_{n,j}^{c}:=
 \sum_{\mu=1}^2\langle c_{n,j},\mathsf A_\mu c_{n,j}\rangle
 \Delta\widetilde W_{n,j}^{\mu},
 \widehat\xi_{n,j}^{a}:=
 \sum_{\mu=1}^2\langle \mathsf A_\mu a_{n,j},a_{n,j}\rangle
 \Delta\widetilde W_{n,j}^{\mu},
\]
where
\[
 \Delta\widetilde W_{n,j}:=
 \widetilde W_{\vartheta_{n,j}}-\widetilde W_{\vartheta_{n,j+1}},
 a_{n,j}=a_{\vartheta_{n,j+1}},  c_{n,j}=c_{\vartheta_{n,j+1}}.
\]
Thus these coefficients are measurable with respect to the Brownian path up to $\vartheta_{n,j+1}$ together with the terminal point.

In the next lemma $p\ge2$ and $ph_{n,j}\le c_0$, where $c_0$ is a sufficiently small absolute constant. The norms are conditional on the terminal point and the past up to
$\vartheta_{n,j+1}$.

\begin{lemma}\label{S:lem:terminal-frozen-coefficients}
With an absolute constant $C$,
\begin{align}
 \left\|\widehat\xi_{n,j}^{c}-M_{n,j}^{c}\right\|_{L^p(\mathbb P^\eta\mid\mathcal F^Z_{\vartheta_{n,j+1}})}
 &\le Cp\,q_{n,j}h_{n,j},\label{S:eq:frozen-terminal-lp}\\
 \left\|\widehat\xi_{n,j}^{c}\right\|_{L^p(\mathbb P^\eta\mid\mathcal F^Z_{\vartheta_{n,j+1}})}
 &\le C\sqrt p\,q_{n,j}\sqrt{h_{n,j}},\label{S:eq:frozen-c-lp}\\
 \left\|\widehat\xi_{n,j}^{a}\right\|_{L^p(\mathbb P^\eta\mid\mathcal F^Z_{\vartheta_{n,j+1}})}
 &\le C\sqrt p\,p_{n,j}\sqrt{h_{n,j}}.\label{S:eq:frozen-a-lp}
\end{align}
The constant is independent of $\eta$, the starting state, $n$ and $j$.
\end{lemma}

\begin{proof}
Put $s=\vartheta_{n,j+1}$, $h=h_{n,j}$, and $R_t=X_{s+t}X_s^{-1}$. Under $\mathbb P^\eta$, conditional on $\mathcal F_s^Z$, It\^o{}\textquotesingle{}s formula for
$X_t=\mathcal A(Z_t)$ gives
\[
 \dd R_t=\sum_{\mu=1}^2\mathsf A_\mu R_t\,
       \dd(\widetilde W_{s+t}^\mu-\widetilde W_s^\mu)
       +\mathsf B_tR_t\,\dd t,\qquad R_0=I,
 \quad
 \mathsf B_t=\operatorname{diag}(3/8,-1/8)
       +\sum_{\mu=1}^2b_\eta^\mu(Z_{s+t})\mathsf A_\mu.
\]
The drift follows from It\^o's formula for $y^{\pm1/2}$ and $-xy^{-1/2}$, and $\|\mathsf B_t\|\le C$ because $|b_\eta|=1$. All $L^q$ norms below use this conditional law.

Stop at $\tau_N=\inf\{t\le h:\|R_t\|\ge N\}\wedge h$ and set $F_{q,N}=\|\sup_{t\le h}\|R_{t\wedge\tau_N}\|\|_q$. Lemma~\ref{S:lem:aux-martingale-moments} and Minkowski's inequality
give, for $q=p,2p$,
\[
 \left\|\sup_{t\le h}\|R_{t\wedge\tau_N}-I\|\right\|_q
 \le C(\sqrt{qh}+h)F_{q,N},\qquad
 F_{q,N}\le1+C(\sqrt{qh}+h)F_{q,N}.
\]
Choose $c_0$ so that $ph\le c_0$ makes the last coefficient at most $1/2$. Conditional Fatou as $N\to\infty$ then yields
\begin{equation}\label{S:eq:relative-affine-one-interval-moments}
 \left\|\sup_{0\le t\le h}\|R_t-I\|\right\|_{L^p(\mathbb P^\eta\mid\mathcal F_s^Z)}
 \le C\sqrt{ph},\qquad
 \left\|\sup_{0\le t\le h}\|R_t\|\right\|_{L^{2p}(\mathbb P^\eta\mid\mathcal F_s^Z)}
 \le C.
\end{equation}
The same argument on $[0,t]$ bounds $\|\sup_{u\le t}\|R_u-I\|\|_q$ by $C\sqrt{qt}$ for $q=p,2p$.

For the real quadratic forms $H_\mu(v)=v^{\mathsf T}\mathsf A_\mu v$ on $\mathbb R^2$, $c_{s+t}=R_tc_s$ gives
\[
 |H_\mu(c_{s+t})-H_\mu(c_s)|
 \le C\|R_t-I\|(1+\|R_t-I\|)|c_s|^2.
\]
H\"older's inequality therefore implies
\begin{equation}\label{S:eq:terminal-coefficient-increment-lp}
 \|H_\mu(c_{s+t})-H_\mu(c_s)\|_{L^p(\mathbb P^\eta\mid\mathcal F_s^Z)}
 \le Cq_{n,j}\sqrt{pt},\qquad 0\le t\le h.
\end{equation}
Subtracting the stochastic integrals and applying conditional Burkholder--Davis--Gundy and Minkowski gives
\[
 \|\widehat\xi_{n,j}^{c}-M_{n,j}^{c}\|_p
 \le C\sqrt p\left(\int_0^h\sum_{\mu=1}^2
       \|H_\mu(c_{s+t})-H_\mu(c_s)\|_p^2\,\dd t\right)^{1/2}
 \le Cp q_{n,j}h.
\]
This is \eqref{S:eq:frozen-terminal-lp}. Conditional on $\mathcal F_s^Z$, the increments with coefficients evaluated at the initial time are centered Gaussian; their variances, by
\eqref{S:eq:s12-affine-contractions}, are exactly
\[
 \mathbb E^\eta[|\widehat\xi_{n,j}^c|^2\mid\mathcal F_s^Z]
       =q_{n,j}^2h/4,\qquad
 \mathbb E^\eta[|\widehat\xi_{n,j}^a|^2\mid\mathcal F_s^Z]
       =p_{n,j}^2h/4.
\]
The Gaussian moment bound $\|N(0,\sigma^2)\|_p\le C\sqrt p\,\sigma$ proves \eqref{S:eq:frozen-c-lp}--\eqref{S:eq:frozen-a-lp}.

Take $p_n=8\log n$. Since $h_{n,j}\le C_\beta/\ell_n$ and $\ell_n\asymp(\log n)^{512}$, the condition $p_nh_{n,j}\le c_0$ holds uniformly for large $n$. Conditional Markov bounds
the two probabilities with thresholds $\me Cp_nq_{n,j}h_{n,j}$ and $\me C\sqrt{p_n}\,p_{n,j}\sqrt{h_{n,j}}$ by $n^{-8}$ each. Integrating over the past and terminal point and
taking a union over at most $n$ intervals gives a summable bound $2n^{-7}$. Borel--Cantelli, followed by $h_{n,j}\le C_\beta/j$, proves
\begin{equation}\label{S:eq:frozen-terminal-error}
 |\widehat\xi_{n,j}^{c}-M_{n,j}^{c}|
 \le C_{\beta,\omega}q_{n,j}\frac{\log n}{j},\qquad
 |\widehat\xi_{n,j}^{a}|
 \le C_{\beta,\omega}p_{n,j}\sqrt{\frac{\log n}{j}}
\end{equation}
for every selected interval and all sufficiently large $n$, on the original joint law and its unchanged bridge extension. The conditional thresholds include the column lengths, so
this argument requires no moments under the joint law of $q_{n,j}$, $p_{n,j}$, or $\eta$.
\end{proof}

\paragraph{Summation by parts}

The integrated fractions of the sampling intervals have bounded variation. This permits summation by parts against the exact increments of $q$, leaving boundary values and a
weighted sum of $q$ itself. For $\ell_n\le J<n$ and real deterministic weights $w_{n,j}$ put
\[
 W=\max_{\ell_n\le j\le J}|w_{n,j}|,\qquad
 V=\sum_{j=\ell_n}^{J-1}|w_{n,j+1}-w_{n,j}|.
\]

\begin{lemma}\label{S:lem:boundary-martingale-primitives}
For every $\epsilon>0$, almost surely for all sufficiently large $n$, uniformly in $J$ and the weights,
\begin{equation}\label{S:eq:weighted-terminal-primitive}
 \left|\sum_{j=\ell_n}^{J}w_{n,j}\widehat\xi_{n,j}^{c}\right|
 \le C_{\beta,\omega,\epsilon}(\log n)
 \left[Wq_{n,\ell_n-1}+(W+V)q_{n,J}\left(\frac nJ\right)^\epsilon\right].
\end{equation}
\end{lemma}
\begin{proof}
We first prove the unweighted estimate
\begin{equation}\label{S:eq:boundary-martingale-primitives}
 \left|\sum_{j=\ell_n}^{J}\widehat\xi_{n,j}^{c}\right|
 \le C_{\beta,\omega,\epsilon}(\log n)
 \left[q_{n,J}+q_{n,\ell_n-1}\right]
 \left(\frac nJ\right)^\epsilon .
\end{equation}
By Lemma~\ref{S:lem:terminal-q-martingale},
\[
 \sum_{j=\ell_n}^{J}M_{n,j}^{c}
 =\frac12(q_{n,\ell_n-1}-q_{n,J}).
\]
The bound \eqref{S:eq:frozen-terminal-error} and the bounds for the boundary solutions imply, with
$0<\eta<\min\{\gamma/2,\epsilon/4\}$,
\begin{equation}\label{S:eq:terminal-partial-sum-bound}
 \begin{aligned}
 q_{n,j}&\le C_{\omega,\eta}q_{n,J}
          (n/J)^{2\eta}(j/J)^{\gamma-\eta},  \ell_n\le j\le J,\\
 \max_{\ell_n\le j\le J}q_{n,j}
 +\sum_{j=\ell_n}^{J}\frac{q_{n,j}}j
 &\le C_{\omega,\eta,\gamma}q_{n,J}(n/J)^{2\eta}.
 \end{aligned}
\end{equation}
Use the two boundary bounds
\[
 q_{n,j}\le C(j/n)^{\gamma-\eta},
 q_{n,J}\ge C^{-1}(J/n)^{\gamma+\eta},
\]
and then the estimate for the sum of powers
\[
 J^{-(\gamma-\eta)}\sum_{j\le J}j^{\gamma-\eta-1}\le C.
\]
Therefore
\[
 \sum_{j=\ell_n}^{J}|\widehat\xi_{n,j}^c-M_{n,j}^c|
 \le C_{\beta,\omega,\epsilon}(\log n)
 q_{n,J}(n/J)^{\epsilon/2},
\]
which, together with the exact boundary point difference, proves \eqref{S:eq:boundary-martingale-primitives}.

For the weighted estimate set $e_j=\widehat\xi_{n,j}^c-M_{n,j}^c$. Exact Abel summation of the telescoping martingale increments gives
\begin{equation}\label{S:eq:terminal-weighted-exact-abel}
 \sum_{j=\ell_n}^{J}w_{n,j}\widehat\xi_{n,j}^c =\frac12\left[w_{n,\ell_n}q_{n,\ell_n-1}-w_{n,J}q_{n,J} +\sum_{j=\ell_n}^{J-1}(w_{n,j+1}-w_{n,j})q_{n,j}\right] +\sum_{j=\ell_n}^{J}w_{n,j}e_j.
\end{equation}
By \eqref{S:eq:terminal-partial-sum-bound}, the bracketed term is bounded by
\[
 C\left[Wq_{n,\ell_n-1}+(W+V)q_{n,J}(n/J)^{2\eta}\right].
\]
The error sum is bounded by
\[
 CW(\log n)q_{n,J}(n/J)^{2\eta}.
\]
Since $2\eta<\epsilon$, this proves \eqref{S:eq:weighted-terminal-primitive}. \end{proof}

On $J_{n,j}$, let
\[
 \omega_{n,j}:=\omega_{n,j}^{[m_n]}
 =\sum_{\substack{m_n\le k<n\\|k-j|\le R_\beta}}
       \int_0^1\lambda_{n,k,j}(s)\,\dd s,
   \ell_n\le j<n.
\]
Lemma~\ref{S:lem:interval-weight-variation} and its extension by zero bound give
\begin{equation}\label{S:eq:retained-chord-mass-variation}
 \max_{\ell_n\le j<n}|\omega_{n,j}|
 +\sum_{j=\ell_n}^{n-2}|\omega_{n,j+1}-\omega_{n,j}|
 \le C_\beta.
\end{equation}
These weights are deterministic functions of $(\beta,n,j,m_n)$. Applying \eqref{S:eq:weighted-terminal-primitive} with $W+V\le C_\beta$ gives
\begin{equation}\label{S:eq:interval-weight-terminal-application}
 \left|\sum_{j=\ell_n}^{J}\omega_{n,j}\widehat\xi_{n,j}^{c}\right|
 \le C_{\beta,\omega,\epsilon}(\log n)
 \bigl(q_{n,J}+q_{n,\ell_n-1}\bigr)
 \left(\frac nJ\right)^\epsilon.
\end{equation}

The contractions of Lemma~\ref{S:lem:trace-separated-pair-factorization} are
\[
 \begin{aligned}
 \mathfrak a_{n,j}
   &=\sum_{k\in\mathcal K_n}
       \left(\int_0^1\lambda_{n,k,j}(s)\dd s\right)
                       \mathbf d(a_{n,k})\cdot\Delta_{n,j},\\
 \mathfrak c_{n,j}
   &=\sum_{k\in\mathcal K_n}
       \left(\int_0^1\lambda_{n,k,j}(s)\dd s\right)
                       \mathbf d(c_{n,k})\cdot\Delta_{n,j}.
 \end{aligned}
\]
Only $|k-j|\le R_\beta$ occur.

\begin{lemma}\label{S:lem:separated-interval-coefficients}
Almost surely, uniformly over all retained sampling intervals for all sufficiently large $n$,
\begin{equation}\label{S:eq:separated-interval-coefficient-expansion}
 \mathfrak a_{n,j}=\omega_{n,j}\widehat\xi^a_{n,j}+e^a_{n,j},
 \mathfrak c_{n,j}=\omega_{n,j}\widehat\xi^c_{n,j}+e^c_{n,j},
\end{equation}
\begin{equation}\label{S:eq:separated-interval-coefficient-errors}
 |e^a_{n,j}|\le C_{\beta,\omega}\frac{\log n}{j}p_{n,j},
 |e^c_{n,j}|\le C_{\beta,\omega}\frac{\log n}{j}q_{n,j}.
\end{equation}
\end{lemma}
\begin{proof}
Set $R^{\rm nor}_{n,j}=\Delta_{n,j}-\Delta\widetilde W_{n,j}$. Corollary~\ref{S:cor:normal-increment-full-array} gives
\begin{equation}\label{S:eq:long-range-normal-innovation}
 |R^{\rm nor}_{n,j}|\le C_{\beta,\omega}h_{n,j}\log n,
 |\Delta\widetilde W_{n,j}|\le C_{\beta,\omega}\sqrt{h_{n,j}\log n}
\end{equation}
for all sufficiently large $n$, uniformly over the retained intervals. The product estimate used in \eqref{S:eq:trace-neighbor-column-products} contains at most $R_\beta$ factors
when $|k-j|\le R_\beta$, and gives
\begin{equation}\label{S:eq:long-range-neighbor-frame}
 \|X_{\vartheta_{n,k+1}}X_{\vartheta_{n,j+1}}^{-1}-I\|_{\rm op}
 \le C_{\beta,\omega}\sqrt{\frac{\log n}{j}}.
\end{equation}
Hence, after increasing the random starting index,
\[
 |a_{n,k}-a_{n,j}|\le d_{n,j}|a_{n,j}|,
 |c_{n,k}-c_{n,j}|\le d_{n,j}|c_{n,j}|,
 d_{n,j}=C_{\beta,\omega}\sqrt{\frac{\log n}{j}},
\]
and neighboring column lengths are comparable by a factor at most two.

Write $w_{n,k,j}=\int_0^1\lambda_{n,k,j}(s)\dd s$. For either choice $v_k=a_{n,k}$ or $c_{n,k}$, define the corresponding error by
\[
 e^v_{n,j}
   =\sum_{k\in\mathcal K_n}w_{n,k,j}
       (\mathbf d(v_k)-\mathbf d(v_j))\cdot\Delta_{n,j}
       +\omega_{n,j}\mathbf d(v_j)\cdot R^{\mathrm{nor}}_{n,j}.
\]
Since $\sum_kw_{n,k,j}=\omega_{n,j}$ and $\Delta_{n,j}=\Delta\widetilde W_{n,j}+R^{\mathrm{nor}}_{n,j}$, this is exactly \eqref{S:eq:separated-interval-coefficient-expansion} for
both columns. The Jacobian of $\mathbf d$ is
\[
 D\mathbf d(v)=\begin{pmatrix}-v_2&-v_1\\-v_1&v_2\end{pmatrix},
 \qquad (D\mathbf d(v))^{\mathsf T}D\mathbf d(v)=|v|^2I_2.
\]
Integration along the line segment from $w$ to $v$ therefore gives
\[
 |\mathbf d(v)-\mathbf d(w)|
      \le\tfrac12(|v|+|w|)|v-w|.
\]
Together with \eqref{S:eq:long-range-neighbor-frame} and $|\mathbf d(v_j)|=|v_j|^2/2$, this yields
\[
 |e^v_{n,j}|
 \le C_\beta|v_j|^2
       \left(\sqrt{\frac{\log n}{j}}\,|\Delta_{n,j}|
                                 +|R^{\mathrm{nor}}_{n,j}|\right)
 \le C_{\beta,\omega}|v_j|^2\frac{\log n}{j}.
\]
We used $\sum_k|w_{n,k,j}|\le2R_\beta+1$ and \eqref{S:eq:long-range-normal-innovation}. Substituting $v=a$ and $v=c$ proves the two bounds in
\eqref{S:eq:separated-interval-coefficient-errors}.

For every $\varepsilon>0$, the errors also satisfy
\begin{equation}\label{S:eq:long-range-c-error-primitive}
 \sum_{j=\ell_n}^{J}|e^c_{n,j}|
 \le C_{\beta,\omega,\varepsilon}(\log n)q_{n,J}
 \left(\frac nJ\right)^{\varepsilon}.
\end{equation}
Apply \eqref{S:eq:terminal-partial-sum-bound} with $0<\eta<\min\{\gamma/2,\varepsilon/4\}$ to the bound $|e^c_{n,j}|\le C_{\beta,\omega}q_{n,j}(\log n)/j$. Since
$2\eta<\varepsilon$, this proves \eqref{S:eq:long-range-c-error-primitive}.
\end{proof}

\begin{proof}[Proof of Proposition~\ref{S:prop:long-range-chord}]
Use the exact sum in \eqref{S:eq:trace-far-pair-collection}. With $b=\ell_n$, put
\[
 P_{n,j}=\sum_{l=b}^{j-L_\beta}\omega_{n,l}\widehat\xi_{n,l}^c,
 \qquad D_{n,j}=\sum_{l=b}^{j-L_\beta}e_{n,l}^c.
\]
Thus the inner sum in \eqref{S:eq:trace-far-pair-collection} is $P_{n,j}+D_{n,j}$. Choose once and for all
\[
 0<\epsilon<\min\{\gamma/8,1/32\},\qquad \gamma=2/\beta.
\]
On the common pathwise event the boundary estimates give
\begin{equation}\label{S:eq:separated-interval-boundary-bounds}
 p_{n,j}\le C_\omega(n/j)^{\gamma+\epsilon},\qquad
 q_{n,j}\le C_\omega(j/n)^{\gamma-\epsilon},\qquad
 p_{n,j}q_{n,j}\le C_\omega(n/j)^{2\epsilon}.
\end{equation}
The estimate \eqref{S:eq:trace-neighbor-column-products} compares indices $j$ and $j-L_\beta$, with a deterministic factor for all sufficiently large $n$.
Lemmas~\ref{S:lem:boundary-martingale-primitives} and \ref{S:lem:separated-interval-coefficients} therefore imply
\begin{equation}\label{S:eq:separated-interval-four-bounds}
 \begin{aligned}
 |P_{n,j}|&\le C_\omega(\log n)
       (q_{n,j}+q_{n,b-1})(n/j)^\epsilon,\\
 |D_{n,j}|&\le C_\omega(\log n)q_{n,j}(n/j)^\epsilon,\\
 |\widehat\xi_{n,j}^a|&\le C_\omega p_{n,j}\sqrt{\log n/j},\\
 |e_{n,j}^a|&\le C_\omega p_{n,j}(\log n)/j.
 \end{aligned}
\end{equation}
The first line uses \eqref{S:eq:weighted-terminal-primitive}. The product bound is $p_{n,j}q_{n,j}\le C_\omega(n/j)^{2\epsilon}$. Because $b\asymp(\log n)^{512}$, one has
$(\log n)/j\le1$ on the summation range for all sufficiently large $n$. Thus $|\mathfrak a_{n,j}|\le C_\omega p_{n,j}\sqrt{\log n/j}$, by
\eqref{S:eq:separated-interval-coefficient-expansion}. Substitution into \eqref{S:eq:trace-far-pair-collection} gives
\begin{equation}\label{S:eq:s12-separated-sum-majorant}
 |E_n^{\mathrm{chord,far}}|
 \le\frac{C_\omega(\log n)^{3/2}}n
       \sum_{j=b}^{n-1} j^{-1/2}p_{n,j}
                          (q_{n,j}+q_{n,b-1})(n/j)^\epsilon.
\end{equation}

The part containing $q_{n,j}$ is at most
\[
 C_\omega(\log n)^{3/2}n^{-1+3\epsilon}
       \sum_{j=b}^n j^{-1/2-3\epsilon}
 \le C_\omega(\log n)^{3/2}n^{-1/2},
\]
since $1/2+3\epsilon<1$. For the other part use $q_{n,b-1}\le C_\omega(b/n)^{\gamma-\epsilon}$ to obtain
\begin{equation}\label{S:eq:s12-separated-boundary-term}
 C_\omega(\log n)^{3/2}
   b^{\gamma-\epsilon}n^{-1+3\epsilon}
       \sum_{j=b}^n j^{-\gamma-2\epsilon-1/2}.
\end{equation}
Put $u=\gamma+2\epsilon+1/2$. Integral comparison gives
\[
 \sum_{j=b}^n j^{-u}\le
 \begin{cases}
 C_un^{1-u},&u<1,\\
 1+\log(n/b),&u=1,\\
 C_ub^{1-u},&u>1.
 \end{cases}
\]
For $u<1$, expression \eqref{S:eq:s12-separated-boundary-term} divided by $(\log n)^{3/2}n^{-1/2}$ is at most $C_\omega(b/n)^{\gamma-\epsilon}\to0$. For $u=1$, it is at most
$C_\omega b^{\gamma-\epsilon}n^{-1/2+3\epsilon}(1+\log n)\to0$. For $u>1$, it is at most $C_\omega(b/n)^{1/2-3\epsilon}\to0$. Each ratio tends to zero for the stated choice of
$\epsilon$, since $b$ is a fixed power of $\log n$.

\end{proof}

\begin{corollary}\label{S:cor:trace-far-rate}
Almost surely, for all sufficiently large $n$,
\begin{equation}\label{S:eq:trace-far-rate}
 |E_n^{\mathrm{chord,far}}|
       \le C_{\beta,\omega}(\log n)^{3/2}n^{-1/2}.
\end{equation}
\end{corollary}
\begin{proof}
Divide \eqref{S:eq:s12-separated-sum-majorant} by $(\log n)^{3/2}n^{-1/2}$. The $q_{n,j}$ contribution is bounded. The $q_{n,b-1}$ contribution tends to zero by the three ratio
estimates following \eqref{S:eq:s12-separated-boundary-term}. The quotient is therefore bounded by a finite $C_\omega$ for all $n\ge n_0(\omega)$.
\end{proof}

\subsubsection{A Riemann sum for the leading order coefficient}
\label{S:sec:s12-conditional-mean}

The preceding estimates leave only the local elapsed time in the conditional second moment. The next theorem collects those estimates; its consequence
\eqref{S:eq:bridge-coefficient-riemann-sum} is the Riemann sum needed for the coefficient. The quadratic variation argument with continuous weights proves only pathwise convergence
to zero for $E_n^{\mathrm{qv}}$, so we leave this term explicit in the quantitative bound.

\begin{theorem}\label{S:thm:trace-conditional-complete}
For each fixed $\beta>0$, almost surely,
\begin{equation}\label{S:eq:trace-conditional-absolute-vanishing}
 |E_n^{\mathrm{diag}}|+|E_n^{\mathrm{bridge}}| +|E_n^{\mathrm{chord,near}}|+|E_n^{\mathrm{chord,far}}| +|E_n^{\mathrm{qv}}|+|E_n^{\mathrm{interp}}| +|E_n^{\mathrm{sampling}}|\longrightarrow0.
 \end{equation}
\end{theorem}
\begin{proof}
Intersect the events for the conditional series, sampled boundary bounds, normal increments, weighted quadratic variation, and the preceding error estimates, including all integer
rows and admissible rational exponents. This event has probability one. Proposition~\ref{S:prop:trace-corrections} and Corollary~\ref{S:cor:trace-far-rate} give
\[
 \begin{aligned}
 &|E_n^{\mathrm{diag}}|+|E_n^{\mathrm{bridge}}|
       +|E_n^{\mathrm{chord,near}}|+|E_n^{\mathrm{chord,far}}|
       +|E_n^{\mathrm{interp}}|+|E_n^{\mathrm{sampling}}|\\
 &\qquad\le C_{\beta,\omega}
      \left(\frac{(\log n)^{9/2}}n
                  +\frac{(\log n)^{3/2}}{\sqrt n}+\frac1{m_n}\right).
 \end{aligned}
\]
Together with $E_n^{\mathrm{qv}}\to0$, from \eqref{S:eq:s12-qv-collected} and Lemma~\ref{S:lem:weighted-qv}, this proves \eqref{S:eq:trace-conditional-absolute-vanishing}. Absorb
$(\log n)^{9/2}/n$ into $(\log n)^{3/2}/\sqrt n$ in the exact decomposition \eqref{S:eq:trace-error-decomposition} to obtain
\begin{equation}\label{S:eq:trace-conditional-quantitative-remainder}
 \left|\mathfrak C_n-\frac4\beta\frac1n
   \sum_{k=m_n}^{n-1}\kappa_\beta(r_k)\frac{\Gamma_{n,k}}{r_k}
          -E_n^{\mathrm{qv}}\right|
 \le C_{\beta,\omega}
       \left(\frac{(\log n)^{3/2}}{\sqrt n}+\frac1{m_n}\right).
\end{equation}
Consequently
\begin{equation}\label{S:eq:bridge-coefficient-riemann-sum}
 n\operatorname{tr}Q_n+n\|B_n\|_{\mathscr K}^2
 =\frac4\beta\frac1n\sum_{k=m_n}^{n-1}
       \kappa_\beta(r_k)\frac{\Gamma_{n,k}}{r_k}+o_\omega(1).
\end{equation}
The quantitative bound leaves $E_n^{\mathrm{qv}}$ separate; its convergence here has no assigned rate.
\end{proof}

\subsubsection{Convergence to the leading order coefficient}
\label{S:sec:s12-sampled-coefficient}

The remaining limit is pathwise. Compact convergence of the boundary columns determines its integrand away from zero. At zero we need an integrable majorant for
$\Gamma_{\beta,\omega}(r)/r$. We prove a positive power bound for $\Gamma$ before taking the Riemann sum limit. The resulting bound also covers a growing column with infinite
squared integral.

Let $\Gamma_{\beta,\omega}$ be the continuous coefficient in \eqref{S:eq:intro-coefficient}. It is formed from the real quadratic form $\mathfrak h_{\mathrm{aff}}$ on $\mathbb R^2$
as follows:
\[
 \Gamma_{\beta,\omega}(r)=\frac12
   \bigl[A_\omega(r)\mathfrak h_{\mathrm{aff}}(a_\omega(r))
             +C_\omega(r)\mathfrak h_{\mathrm{aff}}(c_\omega(r))\bigr].
\]
The discrete coefficient is already defined in \eqref{S:eq:trace-matrix-H-definition}; its partial sums still omit the current partition interval. The continuous version of
\eqref{S:eq:s12-matrix-H-entries} is
\[
 \Gamma_{\beta,\omega}(r)=\frac18A_\omega(r)p_\omega(r)
       +\frac18C_\omega(r)q_\omega(r)
       +\frac{C_\omega(r)}{2p_\omega(r)}.
\]
We use the compact convergence and the estimates for the omitted integrals in Lemma~\ref{S:lem:sampled-exclusive-convergence}, evaluated at $r_k=k/n$.

\begin{lemma}\label{S:lem:Gamma-convergence}
On the event of Lemma~\ref{S:lem:sampled-exclusive-convergence}, for every $0<\delta<1$,
\begin{equation}\label{S:eq:Gamma-compact-convergence}
 \epsilon_n^\Gamma(\delta):=
 \max_{\substack{m_n\le k<n\\r_k\ge\delta}}
        |\Gamma_{n,k}-\Gamma_{\beta,\omega}(r_k)|\longrightarrow0.
\end{equation}
For every $0<\zeta<a_\beta$, for all sufficiently large $n$,
\begin{equation}\label{S:eq:Gamma-boundary-bound}
 0\le\Gamma_{n,k}\le K_{\beta,\omega,\zeta}r_k^\zeta,
 0\le\Gamma_{\beta,\omega}(r)\le K_{\beta,\omega,\zeta}r^\zeta.
\end{equation}
The continuous inequality holds on $0<r\le1$.
\end{lemma}

\begin{proof}
The real quadratic form $\mathfrak h_{\mathrm{aff}}(v)=(v_1^2+5v_2^2)/4$ on $\mathbb R^2$ satisfies, for $v,w\in\mathbb
R^2$,
\begin{equation}\label{S:eq:coefficient-quadratic-continuity}
 \tfrac14|v|^2\le\mathfrak h_{\mathrm{aff}}(v)\le\tfrac54|v|^2,
 \qquad
 |\mathfrak h_{\mathrm{aff}}(v)-\mathfrak h_{\mathrm{aff}}(w)|
 \le\tfrac54(|v|+|w|)|v-w|.
\end{equation}
Let $E_n^a(\delta),E_n^c(\delta),E_n^A(\delta),E_n^C(\delta)$ denote the maximal errors between sampled and continuous values at nodes $r_k\ge\delta$.
Lemma~\ref{S:lem:sampled-exclusive-convergence} makes them tend to zero and bounds all columns and their integrals on this range; only $A$ has zero as a limit of integration, and
$A_\omega(1)<\infty$. Applying
\[
 |s\mathfrak h_{\mathrm{aff}}(v)-t\mathfrak h_{\mathrm{aff}}(w)|
 \le\tfrac54\bigl(|s-t||v|^2
                +|t|(|v|+|w|)|v-w|\bigr)
\]
to the two summands of $2\Gamma$ gives
\begin{equation}\label{S:eq:Gamma-error-bound}
 \epsilon_n^\Gamma(\delta)
 \le C_{\beta,\omega,\delta}
       \bigl(E_n^A(\delta)+E_n^C(\delta)
                       +E_n^a(\delta)+E_n^c(\delta)\bigr).
\end{equation}
For entries of $\mathbf H_{n,k}$, the same calculation uses the quadratic polynomial $v^{\mathsf T}\mathsf A_\mu^{\mathsf T}\mathsf A_\nu v$, whose difference is bounded by
$\|\mathsf A_\mu\|\|\mathsf A_\nu\|(|v|+|w|)|v-w|$. At the last interval, $C_{n,n-1}=0$ and \eqref{S:eq:exclusive-C-compact-error} bounds $C_\omega(r_{n-1})$ by an integral over an
interval of length $1/n$, which tends to zero. This proves both limits on compact intervals.

Put $\gamma=2/\beta$ and choose $0<\eta<\gamma$ with $2\eta<a_\beta-\zeta$, requiring also $\eta<1-\gamma$ when $\gamma<1$. The boundary estimates and comparison of $\rho_{n,k}$
with $r_k$ give
\[
 p_{n,k}\le Kr_k^{-\gamma-\eta},\qquad
 q_{n,k}\le Kr_k^{\gamma-\eta},\qquad
 A_{n,k}\le\frac{Kr_k^{1+\gamma-\eta}}{1+\gamma-\eta},\qquad
 C_{n,k}\le K\int_{r_k}^1u^{-\gamma-\eta}\,\dd u.
\]
The last two inequalities follow by comparing monotone power sums with integrals, and hold for the continuous primitives as well. If $\gamma\ge1$, including $\gamma=1$, then
$Ap,Cq\le K'r^{1-2\eta}$, since $\int_r^1u^{-\gamma-\eta}\dd u\le r^{1-\gamma-\eta}/(\gamma+\eta-1)$. If $\gamma<1$, then $Ap\le K'r^{1-2\eta}$ and $Cq\le K'r^{\gamma-\eta}$. Every
exponent exceeds $\zeta$. For either the sampled or continuous quantities,
\[
 \tfrac18(Ap+Cq)\le\Gamma\le\tfrac58(Ap+Cq),
\]
so \eqref{S:eq:Gamma-boundary-bound} follows. In particular, $\Gamma_{\beta,\omega}$ extends continuously to zero with value zero. Only integrals of $p_\omega$ over $[r,1]$, $r>0$,
have been used.
\end{proof}

On the same pathwise event, define
\[
 f_\omega(r)=\kappa_\beta(r)\frac{\Gamma_{\beta,\omega}(r)}r,
 S_n=\frac1n\sum_{k=m_n}^{n-1}
                   \kappa_\beta(r_k)\frac{\Gamma_{n,k}}{r_k}.
\]

Extend the sampled density by zero below the selected interval and set
\[
 f_n(r)=\frac{\kappa_\beta(r_k)\Gamma_{n,k}}{r_k}
 , (r\in[k/n,(k+1)/n),\ m_n\le k<n).
\]

For fixed $\delta>0$, the convergence just proved is uniform on $[\delta,1]$. We take $n\to\infty$ there first. The bound $\Gamma(r)/r\le Kr^{\zeta-1}$ then controls the omitted
tail uniformly in the row. The proof therefore requires no auxiliary cutoff $\delta_n$ tied to a convergence rate. It applies to every deterministic cutoff with $m_n\to\infty$ and
$m_n/n\to0$.

\begin{proposition}\label{S:prop:coefficient-riemann-limit}
Almost surely, $f_\omega\in L^1(0,1)$ and $\|f_n-f_\omega\|_{L^1(0,1)}\to0$. In particular,
\begin{equation}\label{S:eq:18.81}
 S_n\longrightarrow\int_0^1 f_\omega(r)\,\dd r.
\end{equation}
\end{proposition}
\begin{proof}
Fix $0<\zeta<a_\beta\le1$, let $K=K_{\beta,\omega,\zeta}$ be the constant in \eqref{S:eq:Gamma-boundary-bound}, and put $K_\beta=\sup_{[0,1]}\kappa_\beta$. Then
$0\le f_\omega(r)\le K_\beta Kr^{\zeta-1}$. We will use the bounds
\begin{equation}\label{S:eq:coefficient-uniform-tail}
 \frac1n\sum_{\substack{m_n\le k<n\\r_k<\delta}}
       \kappa_\beta(r_k)\frac{\Gamma_{n,k}}{r_k}
 \le\frac{K_\beta K}{\zeta}\delta^\zeta,
 \qquad
 \int_0^\delta f_\omega(r)\,\dd r
 \le\frac{K_\beta K}{\zeta}\delta^\zeta
\end{equation}
and, for fixed $\delta>0$ and sufficiently large $n$,
\begin{equation}\label{S:eq:coefficient-riemann-error-conditional}
 \left|S_n-\int_0^1f_\omega(r)\,\dd r\right|
 \le\frac{2K_\beta K}{\zeta}\delta^\zeta
      +\frac{K_\beta}{\delta}\epsilon_n^\Gamma(\delta)
      +\omega_{f_\omega,[\delta,1]}(1/n)
      +\frac{\|f_\omega\|_{L^\infty([\delta,1])}}n.
\end{equation}
Indeed, $0<\zeta<1$ implies
\begin{equation}\label{S:eq:coefficient-power-sum}
 \sum_{k=1}^N k^{\zeta-1}
 \le\sum_{k=1}^N\int_{k-1}^kx^{\zeta-1}\,\dd x
 =\frac{N^\zeta}{\zeta}.
\end{equation}
Apply this with $N=\lceil n\delta\rceil-1<n\delta$ to prove the discrete part of \eqref{S:eq:coefficient-uniform-tail}; direct integration proves the other part. Empty sums are
zero.

For $k_\delta=\lceil n\delta\rceil$, choose $n_\delta$ so that $m_n\le k_\delta<n$ for every $n\ge n_\delta$. Replacing $\Gamma_{n,k}$ by $\Gamma_{\beta,\omega}(r_k)$ in the sum
over $k\ge k_\delta$ changes it by at most $K_\beta\epsilon_n^\Gamma(\delta)/\delta$. The error of the ordinary Riemann sum over $[k_\delta/n,1]$ is at most
$\omega_{f_\omega,[\delta,1]}(1/n)$; the integral over $[\delta,k_\delta/n]$ is at most $n^{-1}\|f_\omega\|_{L^\infty([\delta,1])}$. Adding the two tails proves
\eqref{S:eq:coefficient-riemann-error-conditional}. First let $n\to\infty$ with $\delta$ fixed, and then let $\delta\downarrow0$, to obtain \eqref{S:eq:18.81}.

For $L^1$ convergence, if $r\ge\delta$ then its cell has $r_k\ge\delta/2$ for large $n$, and
\[
 |f_n(r)-f_\omega(r)|
 \le\frac{2K_\beta}{\delta}\epsilon_n^\Gamma(\delta/2)
       +\omega_{f_\omega,[\delta/2,1]}(1/n).
\]
On $(0,\delta)$, bound the integral of $f_n$ by the full masses of all cells with $r_k<\delta$ and apply \eqref{S:eq:coefficient-uniform-tail}. This includes the cell meeting
$\delta$; below $m_n/n$, $f_n=0$. Hence
\[
 \limsup_{n\to\infty}\|f_n-f_\omega\|_1
 \le\frac{2K_\beta K}{\zeta}\delta^\zeta\longrightarrow0
 \quad(\delta\downarrow0).
\]
For every bounded measurable $g$ chosen after fixing the path, $|\int g(f_n-f_\omega)|\le\|g\|_\infty\|f_n-f_\omega\|_1\to0$. The sampling and boundary estimates use only
$m_n/n\to0$ and $m_n\ge1$ for all $n$ beyond a fixed index, with the omitted $A$ integral controlled by \eqref{S:eq:exclusive-A-uniform-error}. Thus every cutoff in the statement,
including the requirement $m_n\to\infty$, is covered on the same event.
\end{proof}

\begin{corollary}\label{S:cor:coefficient-time-approximations}
The two alternative sums over the deterministic sampling times satisfy
\begin{equation}\label{S:eq:coefficient-time-approximation-limits}
 \sum_{k=m_n}^{n-1}\tau_{n,k}\kappa_{n,k}\Gamma_{n,k} \longrightarrow\frac4\beta\int_0^1f_\omega(r)\,\dd r,\qquad \sum_{k=m_n}^{n-1}\Gamma_{n,k}\int_0^1\ell_{n,k}(s)\,\dd s \longrightarrow\frac4\beta\int_0^1f_\omega(r)\,\dd r.
\end{equation}
Each differs from $(4/\beta)S_n$ by at most $C_{\beta,\omega}/m_n$. Also the sum with $\kappa_{n,k}$ in place of $\kappa_\beta(r_k)$ in $S_n$ differs from $S_n$ by at most this
order.
\end{corollary}
\begin{proof}
The error $|\kappa_{n,k}-\kappa_\beta(r_k)|\le C_\beta/k$ gives the bound $C_\beta\sum_k\Gamma_{n,k}/k^2$ for the latter replacement.
Equation~\eqref{S:eq:sampling-variance-two-errors} gives the same bound for the first two replacements, using $nr_k=k$. By \eqref{S:eq:Gamma-boundary-bound} and
\eqref{S:eq:coefficient-power-sum},
\[
 \sum_{k=m_n}^{n-1}\frac{\Gamma_{n,k}}{k^2}
 \le\frac1{m_n}\sum_{k=m_n}^{n-1}\frac{\Gamma_{n,k}}k
 \le\frac{K}{m_n n^\zeta}\sum_{k=1}^{n-1}k^{\zeta-1}
 \le\frac{K}{\zeta m_n}.
\]
Proposition~\ref{S:prop:coefficient-riemann-limit} now proves both limits. The bound $O(m_n^{-1})$ controls the time replacements. Convergence in \eqref{S:eq:18.81} uses the
modulus of continuity on compact intervals and the bound near zero from the preceding proposition.
\end{proof}

\subsubsection{Positivity and convergence of the conditional mean}
\label{S:sec:s12-positive-coefficient}

Finiteness follows from the bound near zero. On any closed interval contained in $(0,1)$, both boundary columns are nonzero. Since $\mathfrak h_{\mathrm{aff}}(v)\ge |v|^2/4$, the
integrand defining the coefficient is strictly positive there, and its integral is positive.

\begin{proposition}\label{S:prop:coefficient-positive-measurable}
For each fixed $\beta>0$, the coefficient \eqref{S:eq:intro-coefficient} is measurable from the limiting Brownian path and is almost surely finite and strictly positive.
\end{proposition}
\begin{proof}
Set
\[
 A_*:=\int_0^{1/3}|c_\omega(r)|^2\,\dd r,\quad
 C_*:=\int_{2/3}^1|a_\omega(r)|^2\,\dd r,\quad
 p_*:=\min_{[1/3,2/3]}|a_\omega|^2,\quad
 q_*:=\min_{[1/3,2/3]}|c_\omega|^2.
\]
All four numbers are finite and positive: $q_\omega$ is integrable at zero, the other ranges are compact, and $[a_\omega(r),c_\omega(r)]=1$ makes both continuous columns nonzero.
For $r\in[1/3,2/3]$,
\[
 A_\omega(r)\ge A_*,\qquad C_\omega(r)\ge C_*,\qquad
 \Gamma_{\beta,\omega}(r)\ge\tfrac18(A_*p_*+C_*q_*).
\]
Use $\kappa_\beta\ge1/4$ from Proposition~\ref{S:prop:sampling-offset} and $\int_{1/3}^{2/3}\dd r/r=\log2$ for the lower bound below. Integrating \eqref{S:eq:Gamma-boundary-bound}
gives the upper bound.

The terminal point $q_\infty$ is a measurable path limit, so the columns evaluated at $\theta_\beta(r)$ are jointly measurable in $(\omega,r)$. Parameter integration of their
nonnegative squared lengths, followed by the defining integral of $\PathCoeff$, proves measurability. The integral is finite on the event where the boundary estimates hold; define
it to be zero on its null complement. The compact minima above are also measurable, since they equal infima over countable dense sets.

For $\widetilde a=\lambda a$, $\widetilde c=\lambda^{-1}c$ with $\lambda\in\mathbb R\setminus\{0\}$, homogeneity gives
\[
 \widetilde A=\lambda^{-2}A,\qquad
 \widetilde C=\lambda^2C,\qquad
 \widetilde A\mathfrak h_{\mathrm{aff}}(\widetilde a)
       =A\mathfrak h_{\mathrm{aff}}(a),\qquad
 \widetilde C\mathfrak h_{\mathrm{aff}}(\widetilde c)
       =C\mathfrak h_{\mathrm{aff}}(c).
\]
Thus $\Gamma_{\beta,\omega}$ and $\PathCoeff$ are invariant under this rescaling, with the matrix factor $\mathcal A$ held fixed. Altogether, for $0<\zeta<a_\beta$,
\begin{equation}\label{S:eq:coefficient-explicit-positivity}
 0<\frac{\log2}{8\beta}(A_*p_*+C_*q_*)
       \le\PathCoeff(\omega)
       \le\frac{4\|\kappa_\beta\|_\infty K_{\beta,\omega,\zeta}}
                    {\beta\zeta}<\infty.
\end{equation}
\end{proof}

\begin{corollary}\label{S:cor:coefficient-conditional-mean}
For the conditional covariance $Q_n$ and conditional mean $B_n$ defined above, almost surely,
\begin{equation}\label{S:eq:bridge-conditional-mean-limit}
 \mathfrak C_n=n\operatorname{tr}Q_n+n\|B_n\|_{\mathscr K}^2
       \longrightarrow\PathCoeff(\omega)\in(0,\infty).
\end{equation}
\end{corollary}
\begin{proof}
Fix a path in the intersection of the preceding events of probability one. Write $S_n$ and $f_\omega$ as in Proposition~\ref{S:prop:coefficient-riemann-limit}. The exact kernel
calculation gives $\mathfrak C_n-(4/\beta)S_n=\sum E_n$, where the sum runs over the seven corrections in \eqref{S:eq:trace-error-decomposition}. By
Theorem~\ref{S:thm:trace-conditional-complete},
\[
 |\mathfrak C_n-\PathCoeff(\omega)|
 \le\sum |E_n|+\frac4\beta\left|S_n-\int_0^1f_\omega(r)\dd r\right|
 \longrightarrow0.
\]
The two terms tend to zero by that theorem and Proposition~\ref{S:prop:coefficient-riemann-limit}, respectively. Proposition~\ref{S:prop:coefficient-positive-measurable} proves
$0<\PathCoeff(\omega)<\infty$. For each fixed $n$, Corollary~\ref{S:cor:gaussian-conditional-second-moment} gives the identity
$\mathfrak C_n=n\mathbb E[\|G_n+B_n\|_{\mathscr K}^2\mid\mathscr S_n]$. Hence $\mathfrak C_n\to\PathCoeff(\omega)$ almost surely.
\end{proof}

The estimates of Section~\ref{S:sec:covariance-estimates} also identify the leading order coefficient for finite Gaussian sums. With the convention $G_n^{(\infty)}=G_n$,
$Q_n^{(\infty)}=Q_n$, the exact covariance tail of Theorem~\ref{S:thm:gaussian-finite-mode} yields
\begin{equation}\label{S:eq:s12-finite-mode-mean}
 0\le\mathfrak C_n-
       \bigl(n\operatorname{tr}Q_n^{(M)}+n\|B_n\|_{\mathscr K}^2\bigr)
   =nt_{n,M}\le\frac{Ce_n^{\mathrm{KL}}}{M}.
\end{equation}
Let $K_n\ge1$ be deterministic integers tending to infinity. Fix $0<\epsilon<a_\beta/2$. The uniform conclusion \eqref{S:eq:covariance-noncentral-rate} and
\eqref{S:eq:s12-finite-mode-mean} imply
\begin{equation}\label{S:eq:s12-finite-mode-coefficient}
 \begin{aligned}
 &\sup_{\substack{M\in\mathbb N\cup\{\infty\}\\M\ge K_n}}
      \left|n\|G_n^{(M)}+B_n\|_{\mathscr K}^2-\PathCoeff(\omega)\right|\\
 &\quad\le |\mathfrak C_n-\PathCoeff(\omega)|
       +\frac{Ce_n^{\mathrm{KL}}}{K_n}
       +C_{\beta,\omega,\epsilon}n^{-a_\beta/2+\epsilon}
       \longrightarrow0.
 \end{aligned}
\end{equation}
All three terms tend to zero: the first by Corollary~\ref{S:cor:coefficient-conditional-mean}, the second by $e_n^{\mathrm{KL}}=O_\omega(1)$, and the third by the choice of
$\epsilon$. This is a statement about the squared norm of the finite Gaussian sum. To preserve the quantitative operator remainder, we use the specified cutoff in
\eqref{S:eq:s10-logarithmic-mode-cutoff}.

\subsubsection{Simultaneous pathwise estimates}
\label{S:sec:simultaneous-pathwise-estimates}

All stochastic identities have so far been proved under the indicated law. We now choose one Borel event for the final argument. For a nonnegative measurable sequence $(X_n)$, the
relevant sets are
\begin{equation}\label{S:eq:conditional-event-operators}
 \begin{aligned}
 \{X_n=O(1)\}
  &=\bigcup_{M,N\ge1}\bigcap_{n\ge N}\{X_n\le M\},\\
 \{X_n\to0\}
  &=\bigcap_{h\ge1}\bigcup_{N\ge1}
                              \bigcap_{n\ge N}\{X_n\le h^{-1}\},\\
 \liminf_{n\to\infty}E_n
  &=\bigcup_{N\ge1}\bigcap_{n\ge N}E_n.
 \end{aligned}
\end{equation}
Applied to $n^{a-\epsilon}X_n$, the same formulas express $X_n=O(n^{-a+\epsilon})$. Uniform bounds on a finite row use its maximum. Bounds uniform over the number of Fourier modes
in the partial sum use a supremum over $M\in\mathbb N$, which is still measurable. The Cauchy criterion for Hilbert space series is likewise countable.

For a continuous path, compact suprema equal suprema over rational times and the boundary points. For a jointly measurable matrix field,
\[
 \left\{\operatorname*{ess\,sup}_s\|U(s)\|>q\right\}
    =\left\{\int\mathbf1_{\{\|U(s)\|>q\}}\dd s>0\right\},
 \qquad q\in\mathbb Q.
\]
Both sides are measurable by parameter integration. Using $\mathcal J_n$ from Lemma~\ref{S:lem:weighted-hilbert}, compute the operator norm in the fixed Hilbert space as the
supremum over a countable dense set of unit vectors. The strong measurability statements in Sections~\ref{S:sec:stability-and-fr-echet-expansion-of-the-inverse-kernel}
and~\ref{S:sec:the-conditional-gaussian-operator} therefore make the operator bounds Borel as well. An event first obtained in a completion contains a Borel event of the same full
probability; we use such a subset below.

\phantomsection\label{S:tab:good-event-components} Choose Borel events as follows. Let $\Omega_\beta^{\mathrm{geom}}$
include continuity, the estimates for the boundary columns and their integrals, and the Green inverse identities from Theorem~\ref{S:thm:sine-selfadjoint-domain},
Proposition~\ref{S:prop:5.2}, and Lemma~\ref{S:lem:sampled-exclusive-convergence}. Let $\Omega_\beta^{\mathrm{br}}$ include the bounds for the entire increment array, bridge
expansion, and weighted quadratic variation in Section~\ref{S:sec:the-error-from-sampling-the-brownian-path}, together with the estimates through
Corollary~\ref{S:cor:coefficient-conditional-mean} proved above. These include the sum of products of linear interpolation terms over separated intervals and the finiteness and
positivity of $\PathCoeff(\omega)$.

Let $\Omega_\beta^{\mathrm G}$ include the conditional series identities of Theorem~\ref{S:thm:gaussian-finite-mode}, the pathwise energy estimate
\eqref{S:eq:bridge-weighted-norm-bound}, the quantitative expansion of Theorem~\ref{S:thm:15.1}, and the covariance and concentration estimates of Theorems~\ref{S:thm:16.5}
and~\ref{S:thm:covariance-concentration}. Include also \eqref{S:eq:s10-uniform-mode-truncation}, simultaneously for the countable pairs $(n,M)$. Finally let
$\Omega_\beta^{\mathrm{cmp}}$ include the finite comparison of Proposition~\ref{S:prop:circular-brownian-inverse-comparison} and the three inverse tail estimates of
Theorem~\ref{S:thm:full-boundary-tails}. Each of these events has probability one by the cited results. Their quantitative probability estimates were integrated under their
respective conditional laws before applying Borel--Cantelli.

For losses in exponents take the intersection over admissible positive rationals; for compact intervals take rational boundary points and integer time horizons. Each fixed choice
of parameters may have its own bound and starting index. With the event of conditional versions constructed in Section~\ref{S:sec:conditioning-fields-and-simultaneous-versions},
set
\begin{equation}\label{S:eq:conditional-initial-good-event}
 E_0=\Omega_\beta^{\mathrm{ver}}
       \cap\Omega_\beta^{\mathrm{geom}}
       \cap\Omega_\beta^{\mathrm{br}}
       \cap\Omega_\beta^{\mathrm G}
       \cap\Omega_\beta^{\mathrm{cmp}}.
\end{equation}
Thus $E_0$ is Borel and $\mathbb P(E_0)=1$. To arrange a support property for all the conditional laws used in the proof, let $(T_a)_a$ be the countable observation family
$\mathfrak C_\beta$ specified in Section~\ref{S:sec:conditioning-fields-and-simultaneous-versions}, and let $K_a(t,\dd\omega)$ be regular conditional laws of the full Borel space
given $T_a$. The next lemma is an elementary fact about these kernels and any Borel event $E_0$ of probability one.

\begin{lemma}\label{S:lem:conditional-kernel-stable-event}
There is a Borel $E\subset E_0$ of probability one such that
\begin{equation}\label{S:eq:conditional-kernel-support}
 K_a(T_a\omega,E)=1\qquad(\omega\in E,\ a\in\mathfrak C_\beta).
\end{equation}
\end{lemma}
\begin{proof}
Define recursively
\begin{equation}\label{S:eq:conditional-kernel-direct}
 E_{r+1}=E_r\cap\bigcap_a
               \{\omega:K_a(T_a\omega,E_r)=1\},\qquad
 E=\bigcap_{r\ge0}E_r.
\end{equation}
Kernel measurability makes every set Borel. If $\mathbb P(E_r)=1$, then
\[
 0\le K_a(T_a,E_r)\le1,\qquad
 \int K_a(T_a\omega,E_r)\mathbb P(\dd\omega)=\mathbb P(E_r)=1.
\]
Hence $K_a(T_a,E_r)=1$ almost surely. Induction and countability give $\mathbb P(E_r)=1$ for all $r$ and $\mathbb P(E)=1$. For $\omega\in E$, its membership in $E_{r+1}$ gives
$K_a(T_a\omega,E_r)=1$ for every $r$. Continuity from above of this fixed probability measure yields $K_a(T_a\omega,E)=\lim_rK_a(T_a\omega,E_r)=1$.
\end{proof}

\begin{theorem}\label{S:thm:common-probability-one-event}
For fixed $\beta>0$, the preceding estimates hold simultaneously on a Borel event $\Omega_\beta^{\mathrm{good}}$ of probability one satisfying
\eqref{S:eq:conditional-kernel-support}.
\end{theorem}
\begin{proof}
Apply Lemma~\ref{S:lem:conditional-kernel-stable-event} to \eqref{S:eq:conditional-initial-good-event} and set
\begin{equation}\label{S:eq:common-probability-one-event}
 \Omega_\beta^{\mathrm{good}}=\bigcap_{r\ge0}E_r.
\end{equation}
This is a Borel event of probability one contained in every event defining $E_0$. Any real loss exponent $\epsilon>0$ is covered by an admissible rational $0<\epsilon'<\epsilon$,
since
\[
 r^{a-\epsilon'}\le r^{a-\epsilon},\qquad
 r^{a+\epsilon'}\ge r^{a+\epsilon},\qquad
 n^{-a+\epsilon'}\le n^{-a+\epsilon}
 \quad(0<r\le1,\ n\ge1).
\]
Likewise, enclose a compact interval $[a,b]\subset(0,\infty)$ in $[a',T]$ with rational $0<a'<a$ and integer $T>b$. Continuity extends the bounds at rational times and their
suprema to the whole compact interval, with the same constants and starting indices. Since $E_0\subset\Omega_{\mathrm{qv}}$, Lemma~\ref{S:lem:weighted-qv} applies to every array
satisfying \eqref{S:eq:qv-random-weight-hypotheses}, including the fields $\mathbf H_{n,k}$ and induced weights checked in Proposition~\ref{S:prop:trace-corrections}.

For an observation $T_a$, let $B_a=\{t:K_a(t,\Omega_\beta^{\mathrm{good}})<1\}$. The kernel support property gives $T_a^{-1}(B_a)\cap\Omega_\beta^{\mathrm{good}}=\varnothing$,
hence $(T_a)_\#\mathbb P(B_a)=0$. Taking the countable union over $a$ proves, for the chosen conditional versions,
\begin{equation}\label{S:eq:conditional-good-conditional-probability}
 \mathbb P^\eta(\Omega_\beta^{\mathrm{good}})=1,\qquad
 \mathbb P(\Omega_\beta^{\mathrm{good}}\mid\mathscr S_n)=1,\qquad
 \mathbb P(\Omega_\beta^{\mathrm{good}}\mid\mathcal G_{n;k,j})=1.
\end{equation}
The first identity holds for almost every terminal value $\eta$. For every path in $\Omega_\beta^{\mathrm{good}}$, the latter two identities hold in the chosen versions,
simultaneously for all indicated observations. All kernels remain the chosen versions for their original sigma-fields.
\end{proof}

\phantomsection\label{S:cor:conditional-eventual-and-marginal} The constants and starting indices can be chosen measurably. For a nonnegative sequence $(X_n)$ with
$\limsup_nX_n<\infty$, define on this Borel event
\[
 \begin{aligned}
 M(X)&=\min\{m\in\mathbb N:\exists N\ \forall n\ge N,\ X_n\le m\},\\
 N(X)&=\min\{N\in\mathbb N:\forall n\ge N,\ X_n\le M(X)\}.
 \end{aligned}
\]
For integers $m,N$, their sublevel sets there are
\[
 \{M(X)\le m\}=\bigcup_N\bigcap_{n\ge N}\{X_n\le m\},\qquad
 \{N(X)\le N\}=\bigcap_{n\ge N}\{X_n\le M(X)\}.
\]
Set both variables equal to one on the Borel complement. The displayed formulas prove measurability. Applying the second construction to $X_n\le h^{-1}$ for each $h\in\mathbb N$
gives measurable starting indices for a limit. Finite maxima give common indices for finitely many estimates. The estimates then hold for every $n$ at least the chosen index.

Let $\pi$ retain the Brownian path and the entire array defining the circular operators, and discard only the auxiliary variables adjoined in
Section~\ref{S:sec:the-error-from-sampling-the-brownian-path}. Its marginal law $\mathbb P_0$ is unchanged, so for every Borel event $A_0$ in those original data,
\begin{equation}\label{S:eq:conditional-marginal-descent}
 \mathbb P_0(A_0)=\mathbb P(\pi^{-1}A_0).
\end{equation}
If $\Omega_\beta^{\mathrm{good}}\subset\pi^{-1}A_0$, then $\mathbb P_0(A_0^c)\le\mathbb P((\Omega_\beta^{\mathrm{good}})^c)=0$. The operator norm limits are Borel by
\eqref{S:eq:conditional-event-operators} and the proved measurability of the operators and coefficient. This transfers the conclusions to the original joint law.

\subsubsection{Proof of the two approximation theorems}
\label{S:sec:proof-of-main-theorems}

Fix $\beta>0$ and the joint construction used above. We now apply the estimates to the two differences defined in Section~\ref{S:sec:proof-strategy}. All norms in the following
proof are in $\mathscr K=\mathfrak S_2(\mathcal H)$, where $\mathcal H=L^2([0,1);\mathbb R^2)$. Recall that
$\Delta_{\beta,n}=F_n^{\mathrm{full}}+\varepsilon_{\beta,n}^{\mathrm{samp}}$, and that $\Pi_nD=P_nDP_n$ restricts an inverse kernel to the selected square. On that square,
\eqref{S:eq:s2-operator-expansion} reads $\Pi_n\varepsilon_{\beta,n}^{\mathrm{samp}}=G_n+B_n+R_n^{\mathrm{br}}$. The quantities $d_n=\sqrt n\,\|F_n^{\mathrm{full}}\|$ and
$u_n=\sqrt n\,\|R_n^{\mathrm{br}}\|$ tend to zero by \eqref{S:eq:s2-finite-full} and \eqref{S:eq:s2-operator-remainder}. The conditional second moment $\mathfrak C_n$ and centered
squared norm $\mathfrak s_n$ retain their definitions in \eqref{S:eq:final-mean-fluctuation}. The preceding event construction places all these estimates and the convergence to the
leading order coefficient on one event.

\begin{proof}[Proof of Theorems~\ref{S:thm:main} and~\ref{S:thm:reference-main}]
Fix a path in $\Omega_\beta^{\mathrm{good}}$ from Theorem~\ref{S:thm:common-probability-one-event}. In the real Hilbert--Schmidt space $\mathscr K$,
\[
 \mathfrak s_n=n\bigl(\|G_n\|_{\mathscr K}^2-\operatorname{tr}Q_n\bigr)
                    +2n\langle G_n,B_n\rangle_{\mathscr K}
 \longrightarrow0
\]
by the almost sure concentration estimate. With \eqref{S:eq:s2-mean-limit}, this proves
\begin{equation}\label{S:eq:final-leading-energy}
 n\|G_n+B_n\|^2=\mathfrak C_n+\mathfrak s_n
                       \longrightarrow\PathCoeff(\omega).
\end{equation}
Thus $\sqrt n\|G_n+B_n\|$ is bounded. Put $u_n=\sqrt n\|R_n^{\mathrm{br}}\|\to0$, as in \eqref{S:eq:s2-operator-remainder}. Expansion of \eqref{S:eq:s2-operator-expansion} gives
\begin{equation}\label{S:eq:final-retained-sampling-bound}
 n\bigl|\|\Pi_n\varepsilon_{\beta,n}^{\mathrm{samp}}\|^2
                       -\|G_n+B_n\|^2\bigr|
 \le2u_n\sqrt n\|G_n+B_n\|+u_n^2\longrightarrow0.
\end{equation}
Now \eqref{S:eq:s2-orthogonal-restriction} and \eqref{S:eq:s2-sampling-tail} imply
\begin{equation}\label{S:eq:s2-brownian-limit}
 n\|\varepsilon_{\beta,n}^{\mathrm{samp}}\|^2
                     \longrightarrow\PathCoeff(\omega).
\end{equation}
For the exact circular comparison $\Delta_{\beta,n}=\varepsilon_{\beta,n}^{\mathrm{samp}}+F_n^{\mathrm{full}}$, use $d_n=\sqrt n\|F_n^{\mathrm{full}}\|\to0$ and the preceding
boundedness to obtain
\[
 n\bigl|\|\Delta_{\beta,n}\|^2
                 -\|\varepsilon_{\beta,n}^{\mathrm{samp}}\|^2\bigr|
 \le2d_n\sqrt n\|\varepsilon_{\beta,n}^{\mathrm{samp}}\|+d_n^2
 \longrightarrow0.
\]
Consequently
\begin{equation}\label{S:eq:s2-circular-limit}
 n\|\Delta_{\beta,n}\|^2\longrightarrow\PathCoeff(\omega).
\end{equation}
This uses the $o(n^{-1/2})$ comparison in norm, with no independence assumption. Proposition~\ref{S:prop:coefficient-positive-measurable} provides a finite positive coefficient
measurable from $Z$; Proposition~\ref{S:prop:7.2} provides the finite walk laws and $b_n^{(n)}=q_\infty$.

Let $\pi$ retain $Z$ and the entire circular array, and let $\mathbb P_0=\pi_*\mathbb P$ be their original joint law from Theorem~\ref{S:thm:joint-probability-space}. Strong
measurability of the inverses and of $\PathCoeff$ makes the event $A_0$ of the two displayed limits Borel in these original data. Since
$\Omega_\beta^{\mathrm{good}}\subset\pi^{-1}A_0$, \eqref{S:eq:conditional-marginal-descent} gives $\mathbb P_0(A_0)=\mathbb P(\pi^{-1}A_0)=1$. This proves Theorem~\ref{S:thm:main}
on the original coupling of all rows. The sampled limit depends only on $Z$, so its probability under the Brownian marginal law is also one, proving
Theorem~\ref{S:thm:reference-main}.
\end{proof}

Set $D_{\beta,n}^{\mathrm{br}}=\|\varepsilon_{\beta,n}^{\mathrm{samp}}\|^2$ and $D_{\beta,n}^{\mathrm{circ}}=\|\Delta_{\beta,n}\|^2$. Positivity of the common limit implies, almost
surely for all sufficiently large $n$,
\begin{equation}\label{S:eq:final-two-sided-order}
 \frac{\PathCoeff(\omega)}{2n}
       \le D_{\beta,n}^{\diamond}
       \le\frac{3\PathCoeff(\omega)}{2n},\qquad
 \diamond\in\{\mathrm{br},\mathrm{circ}\}.
\end{equation}
Continuity of the square root and division by the positive limit give
\begin{equation}\label{S:eq:final-norm-and-ratio-limits}
 \sqrt{nD_{\beta,n}^{\diamond}}\longrightarrow\sqrt{\PathCoeff(\omega)},
 \qquad
 \frac{D_{\beta,n}^{\mathrm{circ}}}{D_{\beta,n}^{\mathrm{br}}}
       \longrightarrow1.
\end{equation}
If $b_n>0$ and $nb_n\to0$, then $D_{\beta,n}^{\diamond}/b_n =(nD_{\beta,n}^{\diamond})/(nb_n)\to\infty$ on this event. These are pathwise conclusions for the constructed coupling.
Passage to expectations would require uniform integrability.

\subsubsection{The remaining scalar errors}
\label{S:sec:s12-remaining-errors}

The final identity gives quantitative bounds for the individual errors in the convergence to the leading order coefficient. We retain the discrepancies $E_n^{\mathrm{qv}}$ and
$\mathfrak r_n^{\mathrm{coef}}$, and apply the rates proved for the coupling, concentration and nonlinear errors. This refines the deduction in
Section~\ref{S:sec:proof-of-main-theorems}.

Recall the retained differences $\varepsilon_{\beta,n}^{\mathrm{samp},T} =\Pi_n\varepsilon_{\beta,n}^{\mathrm{samp}}$ and $\Delta_{\beta,n}^{T}=\Pi_n\Delta_{\beta,n}$, and put
$F_n=\Pi_nF_n^{\mathrm{full}}$. Equation~\eqref{S:eq:s2-operator-expansion} and Corollary~\ref{S:cor:15.2} give
\begin{equation}\label{S:eq:final-retained-differences}
 \Delta_{\beta,n}^{T}=G_n+B_n+R_n^{\mathrm{fin}},\qquad
 R_n^{\mathrm{fin}}=R_n^{\mathrm{br}}+F_n.
\end{equation}
All these kernels are supported on the square of the restricted kernel. Since $\|\Pi_n\|=1$, the bounds \eqref{S:eq:s2-finite-full} and \eqref{S:eq:s2-operator-remainder} imply
\begin{equation}\label{S:eq:final-finite-remainder}
 \sqrt n\,\|R_n^{\mathrm{fin}}\|_{\mathscr K}
 \le \sqrt n\,\|R_n^{\mathrm{br}}\|_{\mathscr K}
       +\sqrt n\,\|F_n^{\mathrm{full}}\|_{\mathscr K}
 =O_{\beta,\omega}((\log n)^{-250}).
\end{equation}

Use $\mathfrak C_n$ and $\mathfrak s_n$ from \eqref{S:eq:final-mean-fluctuation}, and $D_{\beta,n}^{\diamond}$ from Section~\ref{S:sec:proof-of-main-theorems}. To distinguish the
error terms for which rates are available, define the two scalar discrepancies
\begin{equation}\label{S:eq:final-coefficient-discrepancies}
 \begin{aligned}
  \mathfrak r_n^{\mathrm{coef}}
   &:=\frac4\beta\left[
      \frac1n\sum_{k=m_n}^{n-1}
           \kappa_\beta(r_k)\frac{\Gamma_{n,k}}{r_k}
       -\int_0^1\kappa_\beta(r)
                    \frac{\Gamma_{\beta,\omega}(r)}r\,\dd r\right],\\
  \mathfrak e_n^{\mathrm{tr}}
   &:=E_n^{\mathrm{diag}}+E_n^{\mathrm{bridge}}
       +E_n^{\mathrm{chord,near}}+E_n^{\mathrm{chord,far}}
       +E_n^{\mathrm{qv}}+E_n^{\mathrm{interp}}+E_n^{\mathrm{sampling}}.
 \end{aligned}
\end{equation}
For $\diamond\in\{\mathrm{br},\mathrm{circ}\}$ put $R_n^{[\mathrm{br}]}=R_n^{\mathrm{br}}$ and $R_n^{[\mathrm{circ}]}=R_n^{\mathrm{fin}}$, and set
\begin{equation}\label{S:eq:final-remainder-energy-definitions}
 \mathfrak u_n^{\diamond} :=2n\langle G_n+B_n,R_n^{[\diamond]}\rangle_{\mathscr K} +n\|R_n^{[\diamond]}\|_{\mathscr K}^2,\qquad \mathfrak v_n^{\diamond} :=n\|\mathsf T_{\beta,n}^{\diamond} -\mathsf T_{\beta,n}^{\mathrm{sine}}\|_{\mathscr K}^2.
\end{equation}
Here $R_n^{[\diamond]}$ is the remainder on the square of the restricted kernel. The quantity $\mathfrak v_n^\diamond$ is the squared norm of the kernel difference on its
complement.

\phantomsection\label{S:prop:final-squared-error-identity} On $\Omega_\beta^{\mathrm{good}}$, for every selected row and
either $\diamond\in\{\mathrm{br},\mathrm{circ}\}$, one has the identity
\begin{equation}\label{S:eq:final-squared-error-identity}
 nD_{\beta,n}^{\diamond}-\PathCoeff(\omega)
   =\mathfrak r_n^{\mathrm{coef}}+\mathfrak e_n^{\mathrm{tr}}
          +\mathfrak s_n+\mathfrak u_n^{\diamond}
          +\mathfrak v_n^{\diamond}.
\end{equation}
Every term on the right tends to zero. More specifically, for each $0<\varepsilon<a_\beta/2$, for all sufficiently large
$n$,
\begin{equation}\label{S:eq:final-squared-error-bounds}
 \begin{aligned}
  |\mathfrak e_n^{\mathrm{tr}}-E_n^{\mathrm{qv}}|
   &\le C_{\beta,\omega}
             \left(\frac{(\log n)^{3/2}}{\sqrt n}+\frac1{m_n}\right),\\
  |\mathfrak s_n|&\le C_{\beta,\omega,\varepsilon}
                                          n^{-a_\beta/2+\varepsilon},\\
  |\mathfrak u_n^{\mathrm{br}}|
   &\le C_{\beta,\omega}(\log n)^{-250},
  |\mathfrak u_n^{\mathrm{circ}}|
       \le C_{\beta,\omega}(\log n)^{-250},\\
  0\le\mathfrak v_n^{\diamond}
   &\le2\mathfrak t_n^{\mathrm{end}}
        \le C_{\beta,\omega,\varepsilon}n^{-a_\beta+\varepsilon}.
 \end{aligned}
\end{equation}
After subtracting $E_n^{\mathrm{qv}}$ and $\mathfrak r_n^{\mathrm{coef}}$,
\begin{equation}\label{S:eq:final-residual-after-unrated-terms}
 \left|nD_{\beta,n}^{\diamond}-\PathCoeff(\omega)
       -\mathfrak r_n^{\mathrm{coef}}-E_n^{\mathrm{qv}}\right|
 \le C_{\beta,\omega}(\log n)^{-250},
 \qquad \diamond\in\{\mathrm{br},\mathrm{circ}\}.
\end{equation}
For these two discrepancies, the conclusion remains qualitative convergence under the stated assumptions.

The conditional second-moment identity and the definition of $\mathfrak s_n$ give
\[
 \mathfrak C_n-\PathCoeff
   =\mathfrak r_n^{\mathrm{coef}}+\mathfrak e_n^{\mathrm{tr}},
 \qquad
 n\|G_n+B_n\|^2=\mathfrak C_n+\mathfrak s_n.
\]
For either choice of $\diamond$, expand $\|G_n+B_n+R_n^{[\diamond]}\|^2$ in the real Hilbert--Schmidt space. The two extra terms are exactly $\mathfrak u_n^\diamond$. The retained
square and its complement are disjoint kernel supports; their $L^2$ inner product is zero, so restoring the complementary kernel adds exactly $\mathfrak v_n^\diamond$, with no
mixed term. These three equalities prove \eqref{S:eq:final-squared-error-identity}.

Proposition~\ref{S:prop:coefficient-riemann-limit} gives $\mathfrak r_n^{\mathrm{coef}}\to0$, and Theorem~\ref{S:thm:trace-conditional-complete} gives
$\mathfrak e_n^{\mathrm{tr}}\to0$ and its stated bound after removing $E_n^{\mathrm{qv}}$. The concentration estimate gives the bound for $\mathfrak s_n$. Put
$u_n^\diamond=\sqrt n\|R_n^{[\diamond]}\|$. Cauchy--Schwarz gives
\[
 |\mathfrak u_n^\diamond|
 \le2u_n^\diamond\sqrt{n}\|G_n+B_n\|+(u_n^\diamond)^2
 \le2u_n^\diamond\sqrt{\mathfrak C_n+|\mathfrak s_n|}+(u_n^\diamond)^2.
\]
On the fixed path the square root is bounded, because $\mathfrak C_n$ converges to a finite coefficient and $\mathfrak s_n\to0$. Equations~\eqref{S:eq:s2-operator-remainder}
and~\eqref{S:eq:final-finite-remainder} therefore give the two logarithmic bounds for $\mathfrak u_n^\diamond$. For the complement, $\|A-B\|^2\le2\|A\|^2+2\|B\|^2$ gives
$0\le\mathfrak v_n^\diamond\le2\mathfrak t_n^{\mathrm{end}}$, proving the remaining bound.

For the last displayed conclusion fix $\varepsilon=a_\beta/4>0$. The sum of all terms other than $\mathfrak r_n^{\mathrm{coef}}$ and $E_n^{\mathrm{qv}}$ is then bounded by
\[
 C_{\beta,\omega}\left[
  n^{-1/2}(\log n)^{3/2}+m_n^{-1}
  +n^{-a_\beta/4}+n^{-3a_\beta/4}+(\log n)^{-250}\right].
\]
For any fixed $c>0$ and $d\ge0$, $n^{-c}(\log n)^{d+250}\to0$: putting $t=\log n$ reduces this to
$\mathrm e^{-ct}t^{d+250}\to0$.
Also $m_n^{-1}\le(\log n)^{-512}\le(\log n)^{-250}$ for $\log n\ge1$. Choose a deterministic threshold, depending on $\beta$, beyond which the three polynomially decaying terms
above are each at most $(\log n)^{-250}$. Taking the maximum of this threshold and the finitely many random starting indices in the preceding estimates proves
\eqref{S:eq:final-residual-after-unrated-terms}.

\subsection{The critical value $\beta=2$}
\label{sec:combined-s-critical}
\label{S:sec:the-critical-value-beta-2}

The boundary columns give an operator proof of the transition at $\beta=2$. Their integrability determines whether a right boundary condition is needed. The same formula for the
Green kernel also permits a precise calculation of the dependence of the leading order coefficient on $\beta$. With one Brownian path fixed for all $\beta>0$,
Theorem~\ref{S:thm:coefficient-beta-regularity} below proves that this coefficient is $C^1$, is analytic on each side of $2$, and has a strictly positive jump in its second
derivative at $2$.

We first recall the classification at the boundary and its operator proof. We then calculate the sampling factor and prove the formula for the leading order coefficient by writing
the integrals defining the Green kernel in Brownian time. The section ends with a diagonal example and a brief account of the earlier transition results and correlation
predictions.

\subsubsection{The transition in the operator domain}
\label{S:sec:s13-operator-domain}
At spectral parameter zero the canonical equation has constant solutions. Their squared norms in the weighted space are the integrals of the squared lengths of the two boundary
columns. The solution determined by the right boundary vector is square integrable for every fixed $\beta>0$; we determine whether the other solution is square integrable. In the
original Brownian coordinates, $y_\theta=\exp(W_\theta^{(2)}-\theta/2)$ and $p_\omega(r)=y_{\theta_\beta(r)}^{-1}$. The substitution $r=\exp(-\beta\theta/4)$ gives exactly
\begin{equation}\label{S:eq:critical-endpoint-integral}
 \int_0^1p_\omega(r)\dd r
 =\frac\beta4\int_0^\infty
   \exp\left(-W_\theta^{(2)}+\frac{2-\beta}{4}\theta\right)\dd\theta.
\end{equation}
For fixed $\beta<2$, put $c=(2-\beta)/4>0$. Since $W_\theta^{(2)}/\theta\to0$, choose $T(\omega,\beta)$ so that the integrand is at least $\exp(c\theta/2)$ for
$\theta\ge T(\omega,\beta)$, and the integral diverges. For fixed $\beta>2$, put $c=(\beta-2)/4>0$; the integrand is at most $\exp(-c\theta/2)$ for $\theta\ge T(\omega,\beta)$, and
the integral is finite. At $\beta=2$, the critical case of Theorem~\ref{S:thm:5.1} gives
$\int_0^\infty\exp(-W_\theta^{(2)})\dd\theta=\infty$.

The other column satisfies $\int_0^1q_\omega(r)\dd r<\infty$ for every fixed $\beta>0$, by the same theorem. If a constant vector is written $v=\alpha u_0+\zeta u_1$, with
$\alpha\ne0$, then the first coordinate of $X_\omega(r)v$ is
\[
 \bigl(-\alpha+\zeta(q_\infty-x_{\theta_\beta(r)})\bigr)
                         y_{\theta_\beta(r)}^{-1/2}.
\]
Its squared modulus is at least $|\alpha|^2p_\omega(r)/4$ for all sufficiently small $r$. Thus in the divergent case no second independent constant solution is square integrable.
In the convergent case both basis vectors are square integrable and so is every constant vector. Proposition~\ref{S:prop:singular-weyl-alternative} and
Theorem~\ref{S:thm:sine-selfadjoint-domain} now give limit point for $\beta\le2$ and limit circle for $\beta>2$, in agreement with \cite{valko2017sine}.

In the original coordinate $t$, the domain is therefore the one in \eqref{S:eq:sine-explicit-singular-domain}: the left condition is $[f(0),u_0]=0$, and the additional right
condition $\lim_{t\uparrow1}[f(t),u_1]=0$ is imposed only when $\beta>2$. For $\beta\le2$ the latter scalar limit vanishes automatically on the maximal domain. In both cases the
terminal point determines the square integrable column of the Green kernel.

The associated stochastic sine equation has a phase transition at the same value. For $\lambda>0$, Valk\'o and Vir\'ag \cite{valko2009continuum} proved that with probability one
there is a finite time after which its phase stays above its terminal multiple of $2\pi$, exactly when $\beta\le2$; for $\beta>2$, approach from below has positive probability.

\subsubsection{The sampling factor}
\label{S:sec:the-sampling-factor-and-a-diagonal-example}
The local sampling variance is proportional to the absolute elapsed Brownian time. Its normalized factor is therefore the integral of the distance from the sampling point. Recall
$b_\beta=1/2-1/\beta$ and $\widehat\sigma_\beta(r)=b_\beta(1-r)$ from Proposition~\ref{S:prop:sampling-offset}. Write $x_+=\max\{x,0\}$. For $v\le0$, integration gives
$\int_0^1|s-v|\dd s=1/2-v$. For $0\le v\le1$, splitting at $v$ gives
\[
 \int_0^1|s-v|\dd s
 =\int_0^v(v-s)\dd s+\int_v^1(s-v)\dd s
 =\frac{v^2+(1-v)^2}{2}.
\]
Since $\widehat\sigma_\beta(r)<1/2$ for $\beta>0$ and $0\le r\le1$, these two calculations cover every value needed here. Thus
\begin{equation}\label{S:eq:critical-sampling-positive-part}
 \kappa_\beta(r)
 =\frac12-b_\beta(1-r)+(b_\beta)_+^2(1-r)^2.
\end{equation}
In particular, with $\gamma=2/\beta$,
\begin{equation}\label{S:eq:critical-sampling-gamma}
 \kappa_\beta(r)=
 \begin{cases}
 \dfrac{\gamma+(1-\gamma)r}{2},&\gamma\ge1,\\[2mm]
 \dfrac{1+\gamma^2}{4}
 +\dfrac{\gamma(1-\gamma)}2r
 +\dfrac{(1-\gamma)^2}{4}r^2,&0<\gamma<1.
 \end{cases}
\end{equation}
For $r<1$ the limiting offset is negative when $\beta<2$, zero when $\beta=2$, and strictly between zero and one when $\beta>2$. At $r=1$ it is zero for every $\beta$.

Since $b'_\beta=\beta^{-2}$ and $b''_\beta=-2\beta^{-3}$, differentiating \eqref{S:eq:critical-sampling-positive-part} gives
\[
 \partial_\beta\kappa_\beta(r)
 =-\frac{1-r}{\beta^2}
   +\frac{2(b_\beta)_+(1-r)^2}{\beta^2}.
\]
This derivative is continuous through $2$, uniformly for $0\le r\le1$. Differentiating on the two open sides gives
\[
 \partial_\beta^2\kappa_\beta(r)
 =\frac{2(1-r)}{\beta^3}
  +\mathbf1_{\{\beta>2\}}(1-r)^2
        \left(\frac2{\beta^4}-\frac{4b_\beta}{\beta^3}\right).
\]
Consequently
\begin{equation}\label{S:eq:critical-sampling-derivatives}
 \begin{aligned}
 \kappa_2(r)&=\frac12,
 &\left.\partial_\beta\kappa_\beta(r)\right|_{\beta=2}
      &=-\frac{1-r}{4},\\
 \left.\partial_\beta^2\kappa_\beta(r)\right|_{2-}
      &=\frac{1-r}{4},
 &\left.\partial_\beta^2\kappa_\beta(r)\right|_{2+}
      &=\frac{1-r}{4}+\frac{(1-r)^2}{8}.
 \end{aligned}
\end{equation}
Thus $\beta\mapsto\kappa_\beta(r)$ is $C^1$ but not twice differentiable at $2$ when $r<1$. The function at $r=1$ is constant. More explicitly, uniformly in $r\in[0,1]$, as
$h\to0$,
\begin{equation}\label{S:eq:critical-sampling-expansion}
 \kappa_{2+h}(r)
 =\frac12-\frac{1-r}{4}h+\frac{1-r}{8}h^2
       +\frac{(1-r)^2}{16}(h_+)^2+O(|h|^3).
\end{equation}
To obtain the remainder, expand $b_{2+h}=h/4-h^2/8+O(|h|^3)$. Its sign agrees with that of $h$, so $(b_{2+h})_+^2=(h_+)^2/16+O(|h|^3)$. The factors $1-r$ and $(1-r)^2$ are bounded
by one.

The finite sampling offset has a smaller correction at the critical value. For $\beta=2$ and $k+1<n$, strict convexity of $x\mapsto x^{-1}$ gives
\[
 \vartheta_{n,k+1}
 =2\sum_{j=k+1}^{n-1}\frac1{j+1/2}
 <2\sum_{j=k+1}^{n-1}\int_j^{j+1}\frac{\dd x}{x}
 =2\log\frac n{k+1}.
\]
The normalized offset $\widehat\sigma_{n,k}$ is therefore strictly negative in this range, and is zero at $k=n-1$. Formula \eqref{S:eq:sampling-critical-finite-offset} gives, for
$r_k=k/n$,
\[
 \widehat\sigma_{n,k}
 =-\frac{1-r_k^2}{24k}+O(k^{-2}).
\]
Since $\kappa_{n,k}=1/2-\widehat\sigma_{n,k}$ here, it follows that
\begin{equation}\label{S:eq:critical-finite-sampling-factor}
 \kappa_{n,k}
 =\frac12+\frac{1-r_k^2}{24k}+O(k^{-2}),\qquad
 \kappa_{n,n-1}=\frac12.
\end{equation}
The expansion is uniform for $1\le k<n$; the last equality is exact.

\subparagraph{Changing the sampling times} The jump in the second $\beta$-derivative of \eqref{S:eq:critical-sampling-positive-part} can be removed without changing the limiting
path. For this calculation only, replace the sampling time by the boundary point time
\[
 \vartheta_{n,k+1}^{0}
     :=\theta_\beta((k+1)/n),\qquad 1\le k<n.
\]
The elapsed time is now nonnegative throughout the interval. Direct integration gives
\[
 \int_0^1\left|\theta_\beta((k+s)/n) -\vartheta_{n,k+1}^{0}\right|\dd s =\frac4\beta\int_0^1\log\frac{k+1}{k+s}\dd s =\frac4\beta\left[1-k\log(1+1/k)\right].
\]
Dividing by the interval length $\tau_{n,k}=(4/\beta)\log(1+1/k)$ gives
\begin{equation}\label{S:eq:critical-alternative-sampling}
 \frac1{\tau_{n,k}}\int_0^1
  \left|\theta_\beta((k+s)/n)-\vartheta_{n,k+1}^{0}\right|\dd s
 =\frac1{\log(1+1/k)}-k
 =\frac12-\frac1{12k}+O(k^{-2}).
\end{equation}
For the last expansion use $\log(1+x)=x(1-x/2+x^2/3+O(x^3))$ and expand its reciprocal on $0\le x\le1/2$. In particular this alternative normalized factor tends to $1/2$ for every
$\beta>0$, with no change at $2$. This calculation concerns only the sampling factor; the circular increment comparison in
Sections~\ref{S:sec:coupling-the-finite-transfer-increments-to-brownian-increments}--\ref{S:sec:comparison-of-the-finite-and-sampled-brownian-operators} uses the original times.
Extending that comparison to the alternative times requires a separate estimate for the coupled increments.

\subsubsection{A second-order singularity of the leading order coefficient}
\label{S:sec:s13-coefficient-regularity}

The sampling factor isolates the possible loss of regularity at $\beta=2$. To study its effect on the full coefficient, one hyperbolic Brownian path $Z$ and its terminal point
$q_\infty$ are used for all $\beta>0$, with the time change $\theta_\beta(r)$ from the introduction. Thus \eqref{S:eq:intro-coefficient} defines $\mathcal C_\beta(\omega)$ on one
probability space, and differentiation with respect to $\beta$ keeps $Z$ fixed.

\begin{theorem}\label{S:thm:coefficient-beta-regularity}
Almost surely, there are real analytic functions $F_\omega,H_\omega$ on $(0,\infty)$, with $H_\omega>0$, such that
\begin{equation}\label{S:eq:coefficient-beta-positive-part}
 \mathcal C_\beta(\omega)=F_\omega(\beta)+(\beta-2)_+^2H_\omega(\beta).
\end{equation}
The coefficient is continuously differentiable on $(0,\infty)$ and real analytic on each side of $2$. Its one-sided second derivatives satisfy
\begin{equation}\label{S:eq:coefficient-beta-second-jump}
 \mathcal C''_{2+}(\omega)-\mathcal C''_{2-}(\omega)
 =\frac14\int_0^1(1-r)^2\frac{\Gamma_{2,\omega}(r)}r\,\dd r>0.
\end{equation}
\end{theorem}
\begin{proof}
In Brownian time put
\[
 \widetilde a_s=\mathcal A(Z_s)u_0,\quad
 \widetilde c_s=\mathcal A(Z_s)u_1,\quad
 \widetilde p_s=|\widetilde a_s|^2,\quad
 \widetilde q_s=|\widetilde c_s|^2,\quad
 h_s^a=\mathfrak h_{\mathrm{aff}}(\widetilde a_s),\quad
 h_s^c=\mathfrak h_{\mathrm{aff}}(\widetilde c_s).
\]
These functions do not depend on $\beta$. Lemma~\ref{S:lem:5.0}, applied with $0<d\le\eta/3$, bounds $\widetilde p_s=y_s^{-1}$ and $\widetilde q_s=(q_\infty-x_s)^2/y_s+y_s$ by
\begin{equation}\label{S:eq:coefficient-brownian-column-growth}
 \widetilde p_s\le K_{\eta,\omega}\me^{(1/2+\eta)s},\qquad
 \widetilde q_s\le K_{\eta,\omega}\me^{-(1/2-\eta)s}.
\end{equation}
Intersecting the events for positive rational $d$ yields one event where these bounds hold for every $\eta>0$ and $s\ge0$. The real quadratic form satisfies
$\widetilde p_s/4\le h_s^a\le5\widetilde p_s/4$ and $\widetilde q_s/4\le h_s^c\le5\widetilde q_s/4$.

For $j=0,1,2$ define
\begin{equation}\label{S:eq:coefficient-beta-moments}
 M_j(\beta)=\frac4\beta\int_0^1r^{j-1}\Gamma_{\beta,\omega}(r)\,\dd r.
\end{equation}
The Brownian path is fixed throughout this calculation. For $r=\mathrm e^{-\beta u/4}$ the columns are $\widetilde a_u,\widetilde c_u$, and substitution in the two spatial
primitives gives
\[
 A_\omega(\mathrm e^{-\beta u/4})
   =\frac\beta4\int_u^\infty\mathrm e^{-\beta v/4}\widetilde q_v\,\dd v,
 \qquad
 C_\omega(\mathrm e^{-\beta u/4})
   =\frac\beta4\int_0^u\mathrm e^{-\beta v/4}\widetilde p_v\,\dd v.
\]
Since $\Gamma_{\beta,\omega}(r)= \tfrac12(A_\omega(r)\mathfrak h_{\mathrm{aff}}(a_\omega(r)) +C_\omega(r)\mathfrak h_{\mathrm{aff}}(c_\omega(r)))$ and
$r^{j-1}\dd r=-(\beta/4)\mathrm e^{-\beta ju/4}\dd u$, the prefactor $4/\beta$ in $M_j$ cancels this Jacobian. Thus, initially as an equality of nonnegative integrals,
\[
 \begin{aligned}
 M_j(\beta)=\frac\beta8\bigg[&
 \int_0^\infty\int_u^\infty
   \mathrm e^{-\beta(ju+v)/4}h_u^a\widetilde q_v\,\dd v\,\dd u\\
 &+\int_0^\infty\int_0^u
   \mathrm e^{-\beta(ju+v)/4}\widetilde p_vh_u^c\,\dd v\,\dd u\bigg].
 \end{aligned}
\]
Tonelli permits either order before finiteness has been established. In the first integral put $(s,t)=(u,v)$; in the second put $(s,t)=(v,u)$. Both domains become
$0\le s\le t<\infty$, giving
\begin{equation}\label{S:eq:coefficient-beta-brownian-integral}
 M_j(\beta)=\frac\beta8\int_{0\le s\le t<\infty}
 \left[\me^{-\beta(js+t)/4}h_s^a\widetilde q_t
       +\me^{-\beta(s+jt)/4}\widetilde p_s h_t^c\right]\dd s\,\dd t.
\end{equation}
Only the deterministic exponential weights and the prefactor now depend on $\beta$. Fix $\beta_0>0$ and put $\beta_-=\beta_0/2$ and $\rho=\beta_0/2$. Choose
$0<\eta<\min\{1/4,\beta_-/16\}$ from the pathwise bounds already proved. Set $v=t-s$. Both nonnegative column products in the preceding integral obey
\[
 h_s^a\widetilde q_t+\widetilde p_s h_t^c
 \le C_{\eta,\omega}\mathrm e^{2\eta s}
                         \mathrm e^{-(1/2-\eta)v}.
\]
For $j=0,1,2$, the linear functions $L_1=js+t$ and $L_2=s+jt$ are nonnegative and satisfy $L_1,L_2\ge s$. If $\beta$ is complex and $|\beta-\beta_0|\le\rho$, then
$\operatorname{Re}\beta\ge\beta_-$ and $|\beta|\le3\beta_0/2$. Consequently
\begin{equation}\label{S:eq:coefficient-beta-integrable-majorant}
 \left|\beta\me^{-\beta L/4}\right|
 (h_s^a\widetilde q_t+\widetilde p_s h_t^c)
 \le K_{\beta_0,\eta,\omega}
 \me^{-(\beta_-/4-2\eta)s}\me^{-(1/2-\eta)v},
 \quad L\in\{js+t,s+jt\}.
\end{equation}
The two decay rates are positive, so the right side has finite integral over $s,v\ge0$. This proves finiteness of each real $M_j$ near $\beta_0$ and provides a common majorant for
its complex extension.

We also verify analyticity by an absolutely convergent power series. For either of the two products, write it as $f(s,t)\ge0$ and use its corresponding $L$. Define
\[
 c_m=\int_{0\le s\le t}
       \mathrm e^{-\beta_0L/4}\frac{(-L/4)^m}{m!}f(s,t)\,\dd s\,\dd t.
\]
Tonelli and the exponential series give
\[
 \begin{aligned}
 \sum_{m\ge0}\rho^m|c_m|
 &\le\int_{0\le s\le t}
       \mathrm e^{-\beta_0L/4}
       \sum_{m\ge0}\frac{(\rho L/4)^m}{m!}f(s,t)\,\dd s\,\dd t\\
 &=\int_{0\le s\le t}\mathrm e^{-(\beta_0-\rho)L/4}f(s,t)\,\dd s\,\dd t<\infty.
 \end{aligned}
\]
The last integral is bounded by the same majorant, without its bounded prefactor. For $|z|<\rho$, Fubini for this absolutely summable series therefore gives
\[
 \int_{0\le s\le t}\mathrm e^{-(\beta_0+z)L/4}f(s,t)\,\dd s\,\dd t
   =\sum_{m\ge0}c_mz^m.
\]
Multiplying by $(\beta_0+z)/8$ and adding the two terms proves that $M_j$ is holomorphic on this disk. In particular it is real analytic near $\beta_0$. For every smaller radius
$\rho'<\rho$, the differentiated series of every fixed order converges uniformly, because $m^k(\rho'/\rho)^m$ is bounded for each $k$. Equivalently, the polynomial factors in $s,t$
arising from these derivatives are integrable against a slightly weaker exponential majorant. The event was defined by the path growth bounds for rational positive exponent losses.
It was fixed before $\beta_0$ was chosen, so this argument proves real analyticity on all of $(0,\infty)$ on that one event. The derivatives act on the exponential weights, with
the time-dependent columns fixed.

Substitution of \eqref{S:eq:critical-sampling-positive-part} gives
\[
 \mathcal C_\beta=(\tfrac12-b_\beta)M_0+b_\beta M_1
           +(b_\beta)_+^2(M_0-2M_1+M_2).
\]
Thus \eqref{S:eq:coefficient-beta-positive-part} holds with the analytic functions
\begin{equation}\label{S:eq:coefficient-beta-analytic-parts}
 F_\omega(\beta)=(\tfrac12-b_\beta)M_0(\beta)+b_\beta M_1(\beta),
 \qquad
 H_\omega(\beta)=\frac{M_0(\beta)-2M_1(\beta)+M_2(\beta)}{4\beta^2}.
\end{equation}
Moreover,
\begin{equation}\label{S:eq:coefficient-beta-positive-factor}
 H_\omega(\beta)=\frac1{\beta^3}
       \int_0^1(1-r)^2\frac{\Gamma_{\beta,\omega}(r)}r\,\dd r>0.
\end{equation}
Both primitives and both columns are nonzero on $0<r<1$, so the integrand is positive there; finiteness follows from the preceding majorant.

Differentiation of the deterministic exponential factors and prefactor in \eqref{S:eq:coefficient-beta-brownian-integral} gives
\[
 \mathcal C'_2=\tfrac12M'_0(2)
                 +\tfrac14\bigl(M_1(2)-M_0(2)\bigr).
\]
The extra polynomial factors in $s,t$ are absorbed by smaller positive exponential rates in \eqref{S:eq:coefficient-beta-integrable-majorant}; the Brownian path is held fixed.
Finally $(\beta-2)_+^2H_\omega(\beta)$ has a continuous first derivative and one-sided second derivatives $0$ and $2H_\omega(2)$ at $2$. This proves
\eqref{S:eq:coefficient-beta-second-jump}, and analyticity of $F,H$ gives
\begin{equation}\label{S:eq:coefficient-beta-taylor}
 \mathcal C_{2+h}=F_\omega(2)+F_\omega'(2)h
       +\tfrac12F_\omega''(2)h^2+H_\omega(2)(h_+)^2+O_\omega(|h|^3).
\end{equation}
\end{proof}

The second-order transition in Theorem~\ref{S:thm:coefficient-beta-regularity} concerns the coefficient of the operator approximation. The proof differentiates deterministic
exponential factors in integrals of the Green columns and uses their growth bounds in Brownian time. The Brownian path $Z$ and the event on which these bounds hold are fixed before
$\beta$ is chosen. Theorems~\ref{S:thm:main} and \ref{S:thm:reference-main} give their inverse asymptotics for each fixed $\beta$. They allow intersection over a prescribed
countable set of parameters, but do not give uniform convergence in $\beta$ or differentiation of the finite error limits. For the alternative sampling factor $1/2$, evaluation of
the same functional gives $M_0(\beta)/2$, which is real analytic throughout $(0,\infty)$. An operator comparison for those alternative times would require the corresponding
estimates for the coupled increments.

\subsubsection{The diagonal example}
\label{S:sec:s13-diagonal-example}
The diagonal matrix function gives an explicit evaluation of the coefficient functional and the jump in its second derivative.

Take the boundary vectors $u_0=(-1,0)^{\mathsf T}$ and $u_1=(0,1)^{\mathsf T}$, so $[u_0,u_1]=1$, and put
\[
 X_\beta^{\mathrm e}(r)
   =\operatorname{diag}(r^{-\gamma/2},r^{\gamma/2}),
 \qquad \gamma=2/\beta.
\]
Then $a^{\mathrm e}(r)=(-r^{-\gamma/2},0)^{\mathsf T}$ and $c^{\mathrm e}(r)=(0,r^{\gamma/2})^{\mathsf T}$. Their squared lengths and primitives are
\[
 p_\beta(r)=r^{-\gamma},\quad q_\beta(r)=r^\gamma,\quad
 Q_\beta(r)=\frac{r^{1+\gamma}}{1+\gamma},\quad
 P_\beta(r)=\int_r^1u^{-\gamma}\dd u.
\]
Recall $\mathfrak h_{\mathrm{aff}}(v)=v_1^2/4+5v_2^2/4$. Hence $\mathfrak h_{\mathrm{aff}}(a^{\mathrm e})=p_\beta/4$ and $\mathfrak h_{\mathrm{aff}}(c^{\mathrm e})=5q_\beta/4$.
Substitution in the coefficient functional gives
\begin{equation}\label{S:eq:critical-diagonal-coefficient-definition}
 G_\beta^{\mathrm e}(r)
 =\frac18p_\beta(r)Q_\beta(r)
       +\frac58q_\beta(r)P_\beta(r),\qquad
 \GeoCoeff=\frac4\beta\int_0^1
           \kappa_\beta(r)\frac{G_\beta^{\mathrm e}(r)}r\dd r.
\end{equation}
The integral $P_\beta$ can diverge as $r\downarrow0$. Finiteness of the integral defining $C_\beta^{\mathrm e}$ follows by integrating first in $r$ in the positive double integral
below.

\begin{theorem}\label{S:thm:diagonal-coefficient-formula}
For every $\beta>0$,
\begin{equation}\label{S:eq:critical-diagonal-closed-form}
 \GeoCoeff=
 \begin{cases}
 \dfrac{4\beta+3}{2\beta(\beta+2)},&0<\beta\le2,\\[3mm]
 \dfrac{15\beta^4+58\beta^3+152\beta^2+112\beta+48}
 {48\beta^2(\beta+1)(\beta+2)},&\beta>2.
 \end{cases}
\end{equation}
\end{theorem}

\begin{proof}
For $l=0,1,2$, Tonelli's theorem gives
\begin{equation}\label{S:eq:estimate-geodesic-monomial-integral}
 J_l(\gamma) :=\int_0^1r^{l+\gamma-1} \left(\int_r^1u^{-\gamma}\dd u\right)\dd r =\int_0^1u^{-\gamma} \left(\int_0^u r^{l+\gamma-1}\dd r\right)\dd u =\frac1{(l+\gamma)(l+1)}.
\end{equation}
In particular, $\int_0^1q_\beta P_\beta\,\dd r/r=1/\gamma$, whereas $p_\beta Q_\beta/r=1/(1+\gamma)$. Boundedness of $\kappa_\beta$ proves finiteness of
\eqref{S:eq:critical-diagonal-coefficient-definition}, including at $\gamma=1$.

Split the sampling factor according to \eqref{S:eq:critical-sampling-positive-part}. The linear part is $\kappa_\beta^0(r)=(\gamma+(1-\gamma)r)/2$. Its two integrals are
\[
 \begin{aligned}
 \int_0^1\kappa_\beta^0(r)\dd r
   &=\frac\gamma2+\frac{1-\gamma}{4}
     =\frac{1+\gamma}{4},\\
 \int_0^1\kappa_\beta^0(r)\frac{q_\beta(r)P_\beta(r)}r\dd r
   &=\frac\gamma2J_0(\gamma)
     +\frac{1-\gamma}{2}J_1(\gamma)
     =\frac{\gamma+3}{4(\gamma+1)}.
 \end{aligned}
\]
The coefficients can have either sign, but each monomial integral is finite by \eqref{S:eq:estimate-geodesic-monomial-integral}. Using $4/\beta=2\gamma$, their contribution is
\[
 \frac\gamma4\left(
       \frac14+5\frac{\gamma+3}{4(\gamma+1)}\right)
 =\frac{\gamma(3\gamma+8)}{8(\gamma+1)}
 =\frac{3}{4\beta}+\frac{5}{4(\beta+2)}.
\]
For the remaining term set
\[
 I_2(\gamma):=\int_0^1(1-r)^2\frac{G_\beta^{\mathrm e}(r)}r\dd r.
\]
Using $(1-r)^2=1-2r+r^2$ and $J_0,J_1,J_2$ gives
\[
 I_2(\gamma) =\frac1{24(1+\gamma)} +\frac58\left( \frac1\gamma-\frac1{\gamma+1} +\frac1{3(\gamma+2)}\right) =\frac{3\gamma^2+11\gamma+15} {12\gamma(\gamma+1)(\gamma+2)}.
\]
Since $(b_\beta)_+^2=(\beta-2)_+^2/(4\beta^2)$, substitution gives
\begin{equation}\label{S:eq:critical-diagonal-positive-part}
 \GeoCoeff=\frac{3}{4\beta}+\frac{5}{4(\beta+2)}
  +(\beta-2)_+^2
       \frac{15\beta^2+22\beta+12}
             {48\beta^2(\beta+1)(\beta+2)}.
\end{equation}
For $\beta\le2$ this is the first branch in \eqref{S:eq:critical-diagonal-closed-form}. For $\beta>2$, put its terms over the common denominator $48\beta^2(\beta+1)(\beta+2)$. The
numerator is
\[
 24\beta(\beta+1)(4\beta+3)
       +(\beta-2)^2(15\beta^2+22\beta+12)
 =15\beta^4+58\beta^3+152\beta^2+112\beta+48.
\]

\end{proof}

For the derivatives at the critical value, put
\[
 F(\beta)=\frac3{4\beta}+\frac5{4(\beta+2)},\qquad
 D(\beta)=\frac{15\beta^2+22\beta+12}{48\beta^2(\beta+1)(\beta+2)}.
\]
Then $C_\beta^{\mathrm e}=F(\beta)+(\beta-2)_+^2D(\beta)$. The values $F(2)=11/16$, $F'(2)=-17/64$, $F''(2)=29/128$, and $D(2)=29/576$ give
\begin{equation}\label{S:eq:critical-diagonal-taylor}
 C_{2+h}^{\mathrm e}
 =\frac{11}{16}-\frac{17}{64}h+\frac{29}{256}h^2
          +\frac{29}{576}(h_+)^2+O(|h|^3).
\end{equation}
Indeed the remainder of the expansion of $F$ is $O(|h|^3)$, and $(h_+)^2(D(2+h)-D(2))=O(|h|^3)$. Therefore
\begin{equation}\label{S:eq:critical-diagonal-second-derivatives}
 (C^{\mathrm e})''(2-)=\frac{29}{128},\qquad
 (C^{\mathrm e})''(2+)=\frac{29}{128}+\frac{29}{288}
                         =\frac{377}{1152}.
\end{equation}
The function is $C^1$ but not twice differentiable at $2$. The jump of its second derivative is exactly twice $D(2)$; it comes from the last term of
\eqref{S:eq:critical-sampling-positive-part}.

For these columns,
\[
 \frac14\int_0^1(1-r)^2\frac{G_2^{\mathrm e}(r)}r\,\dd r
 =\frac{29}{288},
\]
which agrees with the jump in the second derivative in \eqref{S:eq:coefficient-beta-second-jump}. The value $11/16$ is the coefficient functional evaluated at these deterministic
columns.

\subsubsection{Earlier results and predictions}
\label{S:sec:s13-earlier-transitions}

Lewin's account \cite{Lewin2022} traces the log gas prediction to Forrester \cite{Forrester1984}, building on Haldane's theory of quantum fluids \cite{Haldane1981}. The predicted
competing decay powers are $x^{-2}$ and $x^{-4/\beta}$, whose relative order changes at $2$; Qu and Valk\'o record the asymptotic form in equation~(16) of \cite{QuValko2025}. Dumaz
and Malvy \cite{DumazMalvy2026} (Section~1.2) describe the expected transition as Berezinskii--Kosterlitz--Thouless. The results in Theorem~3 of \cite{QuValko2025} and
\cite{DumazMalvy2026} (Theorem~1.1) give, respectively, pointwise upper bounds and averaged decay estimates for the correlation. Their observables are spectral correlations at
large separation.

In \cite{Lewin2022}, the thermodynamic free energy per unit length of the one-dimensional log gas is an explicit real-analytic function of $(\beta,\rho)$ for $\beta,\rho>0$. Its
curvature in $\rho$ changes sign at $\beta=2$. The second-order singularity in \eqref{S:eq:coefficient-beta-second-jump} is a property of the specified operator approximation.
Equation \eqref{S:eq:critical-alternative-sampling} and the analytic functional $M_0(\beta)/2$ show its dependence on the sampling convention. The correlation prediction remains a
separate question about the limiting point process.

\clearpage

\section{Green function comparison for Wigner matrices}\label{sec:green-function-comparison-for-wigner-matrices}

The preceding convergence theorem for the GOE and GUE gives the limiting coefficient process used in the comparison for Wigner matrices. For a general Wigner matrix, however, the
coefficients themselves are dependent, so the comparison is carried out in the independent entries of the original matrix. Polynomial moments connect these entries to the Lanczos
construction. The resulting estimates, together with the estimates for the GOE and GUE, prove Theorem~\ref{W:thm:integrated-potential-and-offdiagonal-coefficients}. Convergence of
the full operator also requires localization of its low energy spectral subspace; that additional argument is given in
Section~\ref{W:sec:localization-of-the-low-energy-spectral-subspace}.

\subsection{Polynomial moments and comparison of entry distributions}\label{sec:comparison-strategy-and-prior-methods}

The observables are smooth functions of polynomial moments of degree $O(N^{1/3})$. Their comparison combines bounds for the orthogonalization map in
Section~\ref{W:sec:orthogonal-polynomials-and-spectral-statistics} with resolvent estimates for the truncated entries in Section~\ref{W:sec:truncation-and-resolvent-estimates}. The
first replacement changes the truncated entries to bounded entries; its error is controlled by the differences of their moments through order three and the fourth-order integral
remainder. The second replacement changes the bounded entries to the GOE or GUE. Here the cubic contributions must be evaluated at the same matrix before summing over the entries.
Section~\ref{W:sec:green-function-comparison} proves these two comparisons and applies them to the characteristic functions and the modulus of continuity of the integrated
potential.

This coefficient theorem controls the operator on each fixed spatial interval. The passage to the full limiting Hilbert space needs a further estimate on the entire low energy
spectral subspace, including eigenvectors with zero spectral weight at the starting vector. The localization statistic is an unnormalized trace whose spectral expansion includes
contributions from those eigenvectors. Its comparison contains an additional cubic term, estimated in Section~\ref{W:sec:localization-of-the-low-energy-spectral-subspace}; the two
probabilistic estimates are then combined through the tridiagonal convergence criterion in Section~\ref{W:sec:main-assembly}.

\subsection{Orthogonal polynomials and spectral statistics}
\label{sec:combined-w-polynomials}
\label{W:sec:orthogonal-polynomials-and-spectral-statistics}

The orthogonalization map involves about $N^{1/3}$ polynomials, and a derivative bound at fixed dimension need not suffice when the degree increases with the matrix size.
Differentiation of the integrated potential yields a cancellation. When its cumulative sum is differentiated before the summands are bounded, the summation gives a square at the
final index and differences of neighbouring off-diagonal coefficients. This is the source of the factor $\rho^{r-1}$ at derivative order $r$ in
Lemma~\ref{W:lem:tightness-directions}.

The representation through the Gram matrix makes this calculation accessible. We recall its moment formulas, also recorded in \S18.2(ix), equations~(18.2.27)--(18.2.30) of
\cite{nistDLMF18}, and develop two complementary estimates. The norm $\|\cdot\|_\diamond$ keeps the dependence on the coefficient quantities explicit and gives the bound
$\rho^{r-1}$. Hilbert--Schmidt norms of the variations of the Gram matrix instead give dimension independent estimates through order three. Both estimates concern variations of the
moments. Derivatives in standardized matrix entries require composition with the moment map and the additional factor $N^{-r/2}$.

The comparison also evaluates the observables away from well conditioned Gram matrices. Smooth spectral modifications and cutoffs extend them to the whole matrix space. Polynomial
bounds for the derivatives of these extensions control the exceptional events, while the sharper bounds through order four apply to the localized functions. The expectation of one
minus each cutoff at $H_N$ is bounded by its Gaussian expectation and the difference of the two cutoff expectations, as in Lemma~\ref{W:lem:tightness-supported-cutoff}.

\subsubsection{Polynomial moments and Jacobi coefficients}
\label{W:sec:the-gram-representation-of-the-jacobi-coefficients}
\label{W:sec:cholesky-representation-of-jacobi-coefficients}

The first Lanczos block can be recovered from the spectral measure at the starting vector. This identifies its recurrence coefficients with those of orthogonal polynomials and
gives a finite matrix formula for their dependence on the moments. Throughout this discussion $H$ is Hermitian and its deterministic starting vector has the normalization

\begin{equation}
 w\in\mathbb F_\beta^N,\qquad \|w\|=1.
 \label{W:eq:fixed-starting-vector}
\end{equation}
The spectral measure at $w$ is
\begin{equation}
 \mu_H^w=\sum_{\alpha=1}^N|\langle w,u_\alpha\rangle|^2\delta_{\lambda_\alpha},
 \qquad H=\sum_{\alpha=1}^N\lambda_\alpha u_\alpha u_\alpha^*.
 \label{W:eq:spectral-measure-at-the-starting-vector}
\end{equation}
With the scalar product convention of Section~\ref{W:sec:convergence-of-tridiagonal-operators},
\begin{equation}
 \langle f(H)w,w\rangle=w^*f(H)w
 =\sum_{\alpha=1}^N f(\lambda_\alpha)|u_\alpha^*w|^2
 =\int f\,\dd\mu_H^w.
 \label{W:eq:spectral-measure-functional-calculus}
\end{equation}
In particular, for complex-valued polynomials $p,q$,
\begin{equation}
 \langle p(H)w,q(H)w\rangle
 =w^*q(H)^*p(H)w=\int p\overline q\,\dd\mu_H^w.
 \label{W:eq:polarization-of-spectral-measures}
\end{equation}
Thus orthogonalizing polynomial vectors is the same operation as orthogonalizing polynomials in this measure. Only the first invariant block generated by $w$ is represented in this
way. The continued Lanczos construction also includes the other invariant blocks.

We use the fixed monic basis $p_j=\psi_j=U_j(\boldsymbol\cdot/2)$, with $p_{-1}=0$. It satisfies
\begin{equation}
 xp_j=p_{j+1}+p_{j-1},\qquad p_0=1,\quad p_1=x.
 \label{W:eq:chebyshev-three-term}
\end{equation}
For the semicircle measure $\dd\mu_{\rm sc}(x)=(2\pi)^{-1}\sqrt{(4-x^2)_+}\,\dd x$, this choice makes the Gram matrix the identity. Indeed, $x=2\cos\theta$ gives
\[
 p_j(2\cos\theta)=\frac{\sin((j+1)\theta)}{\sin\theta},\qquad
 \int p_ip_j\,\dd\mu_{\rm sc}
 =\frac2\pi\int_0^\pi\sin((i+1)\theta)\sin((j+1)\theta)\,\dd\theta
 =\delta_{ij}.
\]
We shall also use the product identity
\begin{equation}
 p_ip_j=\sum_{l=0}^{\min(i,j)}p_{i+j-2l}.
 \label{W:eq:moment-coordinates-Chebyshev}
\end{equation}
At $x=\zeta+\zeta^{-1}$, $p_n(x)=\sum_{u=0}^n\zeta^{n-2u}$. For $i\le j$, the coefficient of $\zeta^{i+j-2v}$ in the product is the number of pairs $0\le a\le i$, $0\le b\le j$
with $a+b=v$. It is $\min(i,v)-\max(0,v-j)+1$ when $0\le v\le i+j$. On the right of \eqref{W:eq:moment-coordinates-Chebyshev}, the same coefficient counts
$0\le l\le\min(i,v,i+j-v)$, and has the same value. Equality of these Laurent polynomials proves the identity for all $x$.

For $1\le\rho\le N$, the first $\rho$ polynomial vectors are collected in the map
\begin{equation}
 \mathcal V_\rho^{(w)}(H)c
 =\sum_{j=1}^\rho c_j\psi_{j-1}(H)w,
 \qquad \mathcal V_\rho^{(w)}(H)
 =(\psi_0(H)w,\ldots,\psi_{\rho-1}(H)w).
 \label{W:eq:krylov-chebyshev-map}
\end{equation}
We suppress the superscript $(w)$ below. Since the basis change from monomials is unitriangular,
\begin{equation}
 \operatorname{Ran}\mathcal V_\rho(H)=\mathcal K_\rho(H,w).
 \label{W:eq:chebyshev-krylov-range}
\end{equation}
The two matrices needed for orthogonalization and multiplication are
\begin{equation}
 G_\rho(H)_{ij}=\langle\psi_{j-1}(H)w,\psi_{i-1}(H)w\rangle,
 \qquad G_\rho=\mathcal V_\rho^*\mathcal V_\rho,
 \label{W:eq:polynomial-gram-matrix}
\end{equation}
\begin{equation}
 M_\rho(H)_{ij}=\langle H\psi_{j-1}(H)w,\psi_{i-1}(H)w\rangle,
 \qquad M_\rho=\mathcal V_\rho^*H\mathcal V_\rho.
 \label{W:eq:polynomial-multiplication-matrix}
\end{equation}
Equivalently,
\begin{equation}
 G_\rho(H)_{ij}=\int\psi_{i-1}\psi_{j-1}\,\dd\mu_H^w,
 \label{W:eq:gram-spectral-measure}
\end{equation}
\begin{equation}
 M_\rho(H)_{ij}=\int x\psi_{i-1}(x)\psi_{j-1}(x)\,\dd\mu_H^w(x).
 \label{W:eq:multiplication-spectral-measure}
\end{equation}
The normalization is fixed:
\begin{equation}
 G_\rho(H)_{11}=1.
 \label{W:eq:gram-first-entry}
\end{equation}
It will remain fixed when the moments are varied. Moreover,
\begin{equation}
 c^*G_\rho(H)c=\|\mathcal V_\rho(H)c\|^2.
 \label{W:eq:gram-quadratic-form}
\end{equation}
Hence $G_\rho>0$ exactly when the first $\rho$ polynomial vectors are independent. By Lemma~\ref{W:lem:lanczos-first-block}, this means that Lanczos iteration has not broken down
before these vectors have been constructed.

On the set $G_\rho>0$, orthogonalization is expressed by the factorization
\begin{equation}
 G_\rho=C_\rho^*C_\rho
 \label{W:eq:cholesky-factorization}
\end{equation}
with $C_\rho$ upper triangular and with positive diagonal. Its inverse is denoted by
\begin{equation}
 Q_\rho=C_\rho^{-1},
 \label{W:eq:inverse-cholesky}
\end{equation}
and the corresponding orthonormalized polynomial map is
\begin{equation}
 \mathcal W_\rho=\mathcal V_\rho Q_\rho.
 \label{W:eq:orthonormalized-krylov-map}
\end{equation}
Then $\mathcal W_\rho^*\mathcal W_\rho=Q_\rho^*G_\rho Q_\rho=I$. Its $k$th column is
\begin{equation}
 w_k=\mathcal W_\rho e_k
 =\sum_{j=1}^k(Q_\rho)_{jk}\psi_{j-1}(H)w,
 \qquad (Q_\rho)_{kk}>0.
 \label{W:eq:wk-polynomial}
\end{equation}
It lies in $\mathcal K_k\cap\mathcal K_{k-1}^\perp$ and has positive leading polynomial coefficient. These conditions fix its phase and are also satisfied by the Lanczos vector.
Thus
\begin{equation}
 \mathcal W_\rho e_k=w_k=q_k,\qquad 1\le k\le\rho.
 \label{W:eq:cholesky-lanczos-vectors}
\end{equation}
The matrix of multiplication in this basis is
\begin{equation}
 J_\rho=Q_\rho^*M_\rho Q_\rho,
 \label{W:eq:exact-jacobi-block}
\end{equation}
with entries
\begin{equation}
 (J_\rho)_{ij}=\langle Hq_j,q_i\rangle=q_i^*Hq_j.
 \label{W:eq:jacobi-entry-inner-product}
\end{equation}
For $i\ge j+2$, $Hq_j\in\mathcal K_{j+1}\subset\mathcal K_{i-1}$ and $q_i\perp\mathcal K_{i-1}$, so this entry is zero. Hermitian symmetry gives the same conclusion for $j\ge i+2$.
To identify the adjacent entries, write $q_j=f_{j-1}(H)w$, where $f_{j-1}$ has positive leading coefficient $\gamma_{j-1}=(Q_\rho)_{jj}$. For $j<\rho$, the expansion of $xf_{j-1}$
in the orthonormal polynomials through degree $j$ has last coefficient $\gamma_{j-1}/\gamma_j$. Thus
\[
 (J_\rho)_{j+1,j}=\langle Hq_j,q_{j+1}\rangle
 =\frac{\gamma_{j-1}}{\gamma_j}>0.
\]
These are the coefficients $b_j$ in the Lanczos recurrence, and therefore
\begin{equation}
 J_\rho=
 \begin{pmatrix}
 a_1&b_1\\
 b_1&a_2&b_2\\
 &b_2&a_3&\ddots\\
 &&\ddots&\ddots
 \end{pmatrix}_{\rho\times\rho}.
 \label{W:eq:finite-lanczos-block}
\end{equation}
Tridiagonality follows from the common moment representation of $G_\rho$ and $M_\rho$ and the fact that multiplication by the variable raises the polynomial degree by one. The
original coordinate projection is
\begin{equation}
 \mathbf P_\rho^{\mathrm{Kry}}
 =\mathcal W_\rho\mathcal W_\rho^*
 =\mathcal V_\rho G_\rho^{-1}\mathcal V_\rho^*.
 \label{W:eq:krylov-projection-formula}
\end{equation}
This is the projection in \eqref{W:eq:lanczos-gram-projection}, expressed in the polynomial basis used here.

The determinant representation will identify the boundary terms in cumulative sums. Put $\tau_k=\det G_k$, $\tau_0=1$, and let $h_k$ be the squared norm of the monic orthogonal
polynomial $\pi_k$. For $1\le k<\rho$, write
\[
 G_{k+1}=\begin{pmatrix}G_k&z_k\\z_k^*&\gamma_k\end{pmatrix}.
\]
Minimizing the norm of $\psi_k-\sum_{i<k}c_i\psi_i$ gives
\[
 \left\|\psi_k-\sum_{i<k}c_i\psi_i\right\|_{L^2(\mu_H^w)}^2
 =(c-G_k^{-1}z_k)^*G_k(c-G_k^{-1}z_k)
       +\gamma_k-z_k^*G_k^{-1}z_k.
\]
The minimum is $h_k$. Block Gaussian elimination therefore gives
\begin{equation}
 h_k=\frac{\tau_{k+1}}{\tau_k}.
 \label{W:eq:monic-norm-gram-ratio}
\end{equation}
The same formula holds at $k=0$, where $h_0=\tau_1=1$. In the monic recurrence
\begin{equation}
 x\pi_k=\pi_{k+1}+a_{k+1}\pi_k+\beta_k\pi_{k-1},\qquad \pi_{-1}=0,
 \label{W:eq:monic-three-term}
\end{equation}
self-adjointness of multiplication gives $\beta_k h_{k-1}=\langle x\pi_k,\pi_{k-1}\rangle =\langle\pi_k,x\pi_{k-1}\rangle=h_k$. After normalization, $b_k^2=\beta_k$, whence
\begin{equation}
 b_k^2=\frac{h_k}{h_{k-1}}
       =\frac{\tau_{k+1}\tau_{k-1}}{\tau_k^2}.
 \label{W:eq:offdiagonal-coefficients-gram-determinant}
\end{equation}
For the semicircle measure, $G_\rho=I$ and $M_\rho=J_{0,\rho}$, the matrix with zero diagonal and unit neighboring entries. These identities define the moment functional of the
semicircle measure used in the calculations below. The random Gram matrix is controlled by the conditioning estimates of Section~\ref{W:sec:bounds-on-gaussian-gram-matrices}.

For comparison we also need empirical moments. For a Hermitian matrix $A$, define
\begin{equation}
 G_\rho^A(H)_{ij}=\operatorname{Tr}[A\psi_{i-1}(H)\psi_{j-1}(H)],
 \label{W:eq:gram-matrix-with-insertion}
\end{equation}
\begin{equation}
 M_\rho^A(H)_{ij}=\operatorname{Tr}[A\psi_{i-1}(H)H\psi_{j-1}(H)].
 \label{W:eq:multiplication-matrix-with-insertion}
\end{equation}
These matrices are real symmetric: the spectral factors commute and are real on the real axis, and the trace of a product of two Hermitian matrices is real. Positivity of
$G_\rho^A$ is used only when $f\mapsto\operatorname{Tr}[Af(H)]$ is a positive functional. The choices $A=P=ww^*$ and $A=N^{-1}I$ give the vector and empirical moments,
respectively. Write $G_\rho^{\rm tr},M_\rho^{\rm tr}$ for the latter, and set
\begin{equation}
 P^\circ=P-N^{-1}I=ww^*-N^{-1}I.
 \label{W:eq:traceless-projection-insertion}
\end{equation}
Linearity in $A$ gives
\begin{equation}
 G_\rho^\circ=G_\rho^{P^\circ}=G_\rho-G_\rho^{\rm tr},
 \label{W:eq:centered-gram}
\end{equation}
\begin{equation}
 M_\rho^\circ=M_\rho^{P^\circ}=M_\rho-M_\rho^{\rm tr}.
 \label{W:eq:centered-multiplication}
\end{equation}
The eigenvalues of $P^\circ$ are $1-N^{-1}$ on $w$ and $-N^{-1}$ on its orthogonal complement. In particular $\|P^\circ\|\le1$ and
\begin{equation}
 \|P^\circ\|_{2,N}
 =\bigl(N^{-1}\operatorname{Tr}(P^\circ)^2\bigr)^{1/2}
 =N^{-1/2}(1-N^{-1})^{1/2}.
 \label{W:eq:traceless-projection-hs-norm}
\end{equation}
For any matrix $B$, $\operatorname{Tr}(P^\circ B)=w^*Bw-N^{-1}\operatorname{Tr}B$. In particular, $\operatorname{Tr}(P^\circ cI)=0$ for every scalar $c$. Products obtained by
differentiating a resolvent are estimated below using this trace identity. The three term recurrence also gives
\begin{equation}
 M_\rho^A(H)_{ij}=G_{\rho+1}^A(H)_{i,j+1}
       +\mathbf1_{\{j\ge2\}}G_{\rho+1}^A(H)_{i,j-1}.
 \label{W:eq:shifted-gram-identity}
\end{equation}
This is an identity of matrix functions, so it remains valid after any fixed number of entry derivatives.

\subsubsection{Variations of the integrated potential}
\label{W:sec:a-finite-sum-rule-and-potential-derivatives}

The representation through the Gram matrix describes individual recurrence coefficients. For the integrated potential, the useful quantity is their cumulative sum: its first
variation contains cancellations between successive diagonal coefficients. We first calculate the variations of the orthogonal polynomials and then sum them before estimating their
absolute values.

Let $\mathcal L$ be a real linear functional on polynomials of degree at most $2\rho-1$, with $s_j=\mathcal L(p_j)$ and $s_0=1$. The polynomials $p_j$ are fixed, whereas the
orthonormal polynomials depend on $\mathcal L$. We differentiate in fixed real linear functionals $v$ satisfying $v(1)=0$; distinct derivative labels may denote the same
functional. If $G_\rho(s)\ge\kappa I$, then for any fixed $v_1,\ldots,v_r$,
\[
 G_\rho\Bigl(s+\sum_{a=1}^rt_av_a\Bigr)
 \ge\Bigl(\kappa-\sum_{a=1}^r|t_a|\,\|\dot G_\rho[v_a]\|\Bigr)I.
\]
Thus $G_\rho\ge\kappa I/2$ when $\sum_a|t_a|\,\|\dot G_\rho[v_a]\|<\kappa/2$. The Cholesky factor with positive diagonal and the orthonormal polynomials are smooth on this
neighbourhood. The derivative bounds below use the two-sided estimate at the base functional. Set
\[
 G_n=(\mathcal L(p_ip_j))_{0\le i,j<n},\qquad
 S_n=(\mathcal L(xp_ip_j))_{0\le i,j<n},\qquad
 z_n=(\mathcal L(p_ip_n))_{0\le i<n}.
\]
Assume
\begin{equation}
 \kappa I_{\rho}\le G_{\rho}\le\kappa^{-1}I_{\rho},\qquad 0<\kappa\le1.
 \label{W:eq:projection-Gram}
\end{equation}
$\mathcal L$ defines a positive inner product on polynomials of degree at most $\rho-1$ by \eqref{W:eq:projection-Gram}; we apply Gram--Schmidt in that finite-dimensional space. If
$G_n=L_nL_n^{\mathsf T}$, with $L_n=C_n^{\mathsf T}$ the positive diagonal lower Cholesky factor, then $J_n=L_n^{-1}S_nL_n^{-\mathsf T}$ is the ordered multiplication matrix.
Degree and orthogonality make it a real Jacobi matrix, with diagonal $a_1,\ldots,a_n$ and positive off-diagonal coefficients $b_1,\ldots,b_{n-1}$. Write $D_n=\det G_n$, $D_0=1$,
and $\mathfrak h_j=D_{j+1}/D_j$. Then $b_j^2=\mathfrak h_j/\mathfrak h_{j-1}$. The squared norm $\mathfrak h_j$ of the monic orthogonal polynomial is the minimum
$\min_c(c,1)G_{j+1}(c,1)^{\mathsf T}$ over all monic polynomials of degree $j$ in the fixed monic basis. The lower bound on the Gram matrix and the choice $c=0$ give
$\kappa\le\mathfrak h_j\le\kappa^{-1}$, hence $\kappa\le b_j\le\kappa^{-1}$.

For $1\le K\le {\rho}-2$, define
\begin{equation}
 \mathcal T_K^\sigma
 =\sigma\sum_{j=1}^Ka_j
   +2\sum_{j<K}(b_j-1)+(b_K-1),\qquad
 Y_K^\sigma=-{m_N}\mathcal T_K^\sigma,
 \quad\sigma\in\{-1,1\}.
 \label{W:eq:projection-T}
\end{equation}
For $\sigma=1$, $Y_K^\sigma$ is the exact integrated potential from Section~\ref{W:sec:convergence-of-tridiagonal-operators}; $\sigma=-1$ changes the sign of the diagonal
coefficients. The quadratic form assigns weight one to the last off-diagonal coefficient and weight two to each interior coefficient. The auxiliary coefficient at zero is one and
is not a recurrence coefficient to be differentiated.

Let $q_0,\ldots,q_K$ be the orthonormal polynomials for $\mathcal L$, with positive leading coefficients. Let $\dot{\mathcal L}$ be a real moment direction annihilating constants,
and put
\[
 E_{ij}=\dot{\mathcal L}(q_iq_j),\qquad 0\le i,j\le K.
\]
Define $\mathscr L(E)$ by
\[
 \mathscr L(E)_{ij}=\begin{cases}E_{ij},&i>j,\\ E_{ii}/2,&i=j,\\0,&i<j.\end{cases}
\]
\begin{lemma}
\label{W:lem:projection-variation}
With $\mathbf q=(q_0,\ldots,q_K)^{\mathsf T}$,
\begin{align}
 D_{s}\mathbf q[\dot{\mathcal L}]&=-\mathscr L(E)\mathbf q,
 \label{W:eq:projection-qvariation}\\
 D_{s}\log b_j[\dot{\mathcal L}]&=d_j(E):=\tfrac12(E_{jj}-E_{j-1,j-1}),
 \quad 1\le j\le K,
 \label{W:eq:projection-blog}\\
 D_{s}\sum_{j\le K}a_j[\dot{\mathcal L}]
 &=b_K\dot{\mathcal L}(q_{K-1}q_K).
 \label{W:eq:projection-adiag}
\end{align}
\end{lemma}
\begin{proof}
We differentiate on the positive neighbourhood constructed above. The normalization $s_0=1$ gives $q_0=1$ throughout this neighbourhood.

Write $D_sq_j=\sum_{i=0}^j c_{ji}q_i$. For $i<j$, differentiation of $\mathcal L(q_iq_j)=0$ gives
\[
 0=\dot{\mathcal L}(q_iq_j)+\mathcal L((D_sq_i)q_j)
        +\mathcal L(q_iD_sq_j)=E_{ij}+c_{ji},
\]
since $\deg(D_sq_i)\le i<j$. Differentiation of the unit norm gives $0=E_{jj}+2c_{jj}$. Thus $c_{ji}=-E_{ji}$ for $i<j$ and $c_{jj}=-E_{jj}/2$, proving
\eqref{W:eq:projection-qvariation}.

For the squared norm of the monic orthogonal polynomial $\mathfrak h_j$, the polynomial derivative has degree less than $j$, so orthogonality gives
$D_{s}\log \mathfrak h_j=\dot{\mathcal L}(q_j^2)$. Differentiating $\log b_j=(\log \mathfrak h_j-\log \mathfrak h_{j-1})/2$ proves \eqref{W:eq:projection-blog}. Using
$xq_j=b_{j+1}q_{j+1}+a_{j+1}q_j+b_jq_{j-1}$ gives
\begin{align*}
 D_sa_{j+1}
 &=\dot{\mathcal L}(xq_j^2)+2\mathcal L(xq_jD_sq_j)\\
 &=b_{j+1}E_{j,j+1}+a_{j+1}E_{jj}+b_jE_{j,j-1}
        -a_{j+1}E_{jj}-2b_jE_{j,j-1}\\
 &=b_{j+1}\dot{\mathcal L}(q_jq_{j+1})
          -b_j\dot{\mathcal L}(q_{j-1}q_j).
\end{align*}
The second term is absent at $j=0$. Summing over $j=0,\ldots,K-1$ leaves only the last mixed product and gives \eqref{W:eq:projection-adiag}.
\end{proof}

For a vector $u=(u_1,\ldots,u_K)$ define the symmetric matrix $\mathcal M_K^\sigma(u)$ by its quadratic form
\begin{align}
 \mathbf q^{\mathsf T}\mathcal M_K^\sigma(u)\mathbf q
 ={}&\sigma u_Kq_{K-1}q_K
 +\sum_{j<K}u_j(q_j^2-q_{j-1}^2)
 +\tfrac12u_K(q_K^2-q_{K-1}^2).
 \label{W:eq:projection-Mdef}
\end{align}
This definition is linear in $u$. Write $M=\mathcal M_K^\sigma(b)$. Substituting the variations of $a_j$ and $b_j$ in \eqref{W:eq:projection-T}, and collecting adjacent squares
before estimating them, gives the derivative polynomial $\mathscr P_K^\sigma=\mathbf q^{\mathsf T}M\mathbf q$, namely
\begin{equation}
 \mathscr P_K^\sigma
 =\tfrac12b_K(q_K+\sigma q_{K-1})^2
       +\sum_{j=1}^{K-1}(b_j-b_{j+1})q_j^2-b_1q_0^2.
 \label{W:eq:potential-derivative-kernel}
\end{equation}
Writing $E(v)_{ij}=v(q_iq_j)$, we obtain
\[
 D_s\mathcal T_K^\sigma[v]=v(\mathscr P_K^\sigma)
                         =\operatorname{Tr}(ME(v)).
\]
The interior coefficients $b_j-b_{j+1}$ are controlled by their total variation. At the moment functional of the semicircle measure they vanish. For a general moment functional,
their contribution remains in the first derivative estimate.

At the moment functional of the semicircle measure $q_j=p_j$, $b_j=1$, and $\mathscr P_K^\sigma=\tfrac12(p_K+\sigma p_{K-1})^2-1$. The product identity gives
\begin{equation}
 D_s\mathcal T_K^\sigma[v]
 =\sigma\sum_{j=0}^{K-1}v(p_{2j+1})
       +\sum_{j=1}^{K-1}v(p_{2j})+\tfrac12v(p_{2K}).
 \label{W:eq:free-potential-moment-variation}
\end{equation}
In particular, the sum of the absolute coordinate derivatives is $2K-\tfrac12$. There is no dimension independent bound on this sum, even at the moment vector of the free Jacobi
matrix. The estimate below instead keeps the whole polynomial $\mathscr P_K^\sigma$ until it is paired with a moment direction.

\Needspace{6\baselineskip}

\subsubsection{Higher derivatives with respect to the moments}
\label{W:sec:derivative-growth-at-the-critical-scale}

The first variation formula suggests how the dimension enters higher derivatives. At first order, the undifferentiated column of $C$ is controlled in $\ell^1$, so the coefficient
estimate has no factor of $\rho$. Further derivatives can introduce sums over as many as $\rho$ indices. Keeping those sums explicit gives a factor $\rho^{r-1}$ at order $r$, which
is the dependence needed when these variations are later represented by resolvent integrals.

Throughout this argument $\rho\ge3$, the normalized functional satisfies \eqref{W:eq:projection-Gram}, and $\kappa$ is fixed. Matrix trace and Hilbert--Schmidt norms are denoted by
$\|\cdot\|_1$ and $\|\cdot\|_{\rm HS}$. For a normalized functional $s$ on polynomials through degree $2\rho-2$, variations satisfying $v(1)=0$ are measured by

\begin{equation}
 \norm v_\diamond=
 \max\left\{\rho^{-1}\norm{\dot G_\rho[v]}_1,
                 \max_{0\le i,j<\rho}|v(p_ip_j)|\right\}.
 \label{W:eq:tightness-norm}
\end{equation}
For $0\le n\le2\rho-2$, take $i_n=\min(n,\rho-1)$ and $j_n=n-i_n$. The polynomial $p_{i_n}p_{j_n}$ is monic of degree $n$. The transition matrix from this list to
$1,x,\ldots,x^{2\rho-2}$ is triangular with diagonal one. The list is therefore a basis. Vanishing of every $v(p_ip_j)$ forces $v=0$, proving that \eqref{W:eq:tightness-norm} is a
norm. The same formula defines a norm for complex linear functionals on these real polynomials. Only moments through degree $2\rho-2$ occur in these estimates; the extra moment
used for the last diagonal entry of $J_\rho$ is not needed for $\mathcal T_K^\sigma$, $K\le\rho-2$.

At the given positive moment functional write
\[
 q_j=\sum_{i\le j}C_{ji}p_i,\qquad CG_\rho C^{\mathsf T}=I,
\]
with positive leading coefficients. Here $C=Q_\rho^{\mathsf T}$ denotes the lower triangular coefficient matrix. The upper Cholesky factor keeps the notation $C_\rho$ of
\eqref{W:eq:cholesky-factorization}. Set
\begin{equation}
 Q_* =\max\left\{1,\max_j\sum_i|C_{ji}|,
                         \max_i\sum_j|C_{ji}|\right\},\qquad
 B=1+\sum_{j=1}^{\rho-3}|b_{j+1}-b_j|.
 \label{W:eq:comparison-QB}
\end{equation}
The row sums give $|q_j(x)|\le Q_*\max_{i<\rho}|p_i(x)|$ for $x\in\mathbb C$. The column sums control the first derivative of the coefficient matrix; the following induction gives
the dependence on $\rho$ at every fixed order.

\begin{lemma}
\label{W:lem:polynomial-coefficient-derivatives}
On $\kappa I\le G_\rho\le\kappa^{-1}I$, for every fixed $r\ge1$,
\begin{align}
 \max_{0\le i,j<\rho}|D_s^rC_{ij}[v_1,\ldots,v_r]|
 &\le C_{r,\kappa}\rho^{r-1}Q_*^{2r+1}\prod_a\|v_a\|_\diamond,\notag\\
 \max_{1\le j<\rho}|D_s^rb_j[v_1,\ldots,v_r]|
 &\le C_{r,\kappa}\rho^{r-1}Q_*^{2r}\prod_a\|v_a\|_\diamond.
 \label{W:eq:polynomial-and-jacobi-derivatives}
\end{align}
\end{lemma}
\begin{proof}
All derivatives are taken at the fixed moment functional $s$; the quantities $Q_*$ and $\|v_a\|_\diamond$ in the estimates below are evaluated there. In particular, these
quantities are not differentiated. For a symmetric matrix $E$, define $\mathscr L(E)_{ij}=E_{ij}$ for $i>j$, $\mathscr L(E)_{ii}=E_{ii}/2$, and $\mathscr L(E)_{ij}=0$ for $i<j$.
For a direction $v$, put $E(v)_{ij}=v(q_iq_j)$. Since $q_i=\sum_kC_{ik}p_k$, the definition of $\|v\|_\diamond$ gives
\begin{equation*}
 |E(v)_{ij}|
 \le \sum_{k,l}|C_{ik}|\,|C_{jl}|\,|v(p_kp_l)|
 \le Q_*^2\|v\|_\diamond.
\end{equation*}
Differentiating $\mathcal L(q_iq_j)=\delta_{ij}$ for all $0\le i,j<\rho$ gives the first variation identity
\begin{equation*}
 D_sC[v]=-\mathscr L(E(v))C.
\end{equation*}
For every $1\le j<\rho$, including $j=\rho-1$, the derivative of $b_j$ can be obtained from the monic norms. Let $\pi_j$ be the monic orthogonal polynomial and
$\mathfrak h_j=\mathcal L(\pi_j^2)$ for $0\le j<\rho$. The derivative of $\pi_j$ has degree at most $j-1$. Orthogonality therefore gives
\begin{equation*}
 D_s\mathfrak h_j[v]=v(\pi_j^2)
 +2\mathcal L(\pi_jD_s\pi_j[v])=v(\pi_j^2),
 \qquad D_s\log\mathfrak h_j[v]=v(q_j^2).
\end{equation*}
As $b_j^2=\mathfrak h_j/\mathfrak h_{j-1}$, this yields
\begin{equation}
 \begin{aligned}
 D_sC[v]&=-\mathscr L(E(v))C,\\
 D_s\log b_j[v]&=\tfrac12(E(v)_{jj}-E(v)_{j-1,j-1}),
 \quad 1\le j<\rho.
 \end{aligned}
 \label{W:eq:tightness-coefficient-recursion}
\end{equation}
using only the positivity of $G_\rho$. Consequently
\begin{equation*}
 |D_sC_{ij}[v]|
 \le \max_{i,k}|E(v)_{ik}|\sum_k|C_{kj}|
 \le Q_*^3\|v\|_\diamond,
\end{equation*}
which proves the first derivative bound.

We prove the higher derivative bounds by induction using labelled mixed derivatives. For a set of labels $A$ write $h_A=\prod_{a\in A}\|v_a\|_\diamond$ and $h_\varnothing=1$. If
$f=\sum_k\widehat f_kp_k$, set $|f|_{\rm c}=\sum_k|\widehat f_k|$. For polynomials $f,g$ of degree at most $\rho-1$, every index in their expansions in the fixed basis is less than
$\rho$. Consequently
\begin{equation*}
 |v(fg)|\le\|v\|_\diamond |f|_{\rm c}|g|_{\rm c},
\end{equation*}
because every product $p_kp_l$ has $v$ value of modulus at most $\|v\|_\diamond$. The induction hypothesis at each order $a<r$ implies
\begin{equation*}
 \max_i|D_Aq_i|_{\rm c}
 \le c_a\rho^a Q_*^{2a+1}h_A
 \quad (|A|=a\ge1),
 \qquad \max_i|q_i|_{\rm c}\le Q_*.
\end{equation*}
Indeed, $D_Aq_i=\sum_j(D_AC_{ij})p_j$ has at most $\rho$ coefficients, and the induction bound for each coefficient at order $a$ is $c_a\rho^{a-1}Q_*^{2a+1}h_A$. We take $c_0=1$ in
the corresponding bound for $a=0$.

Fix the last direction $v_r$. Differentiating its moment pairing gives the identity
\begin{equation*}
 D_AE(v_r)_{ij}
 =\sum_{U\sqcup V=A}v_r((D_Uq_i)(D_Vq_j)).
\end{equation*}
The preceding coefficient estimates, applied separately to the two factors, therefore yield
\begin{equation*}
 \max_{i,j}|D_AE(v_r)_{ij}|
 \le c'_a\rho^a Q_*^{2a+2}h_A\|v_r\|_\diamond,
 \qquad |A|=a,
\end{equation*}
where one can take
$c'_a=\sum_{u=0}^a\binom au c_uc_{a-u}$.
There are exactly $\binom au$ allocations with $|U|=u$; no constant in this estimate depends on $\rho$.

Put $S=\{1,\ldots,r-1\}$. Differentiation of the exact first variation of $C$ gives
\begin{equation*}
 D_{S\cup\{r\}}C
 =-\sum_{A\sqcup B=S}\mathscr L(D_AE(v_r))D_BC.
\end{equation*}
For $B=\varnothing$, the column sum of the undifferentiated $C$ is at most $Q_*$; hence each matrix entry of this summand is bounded by
\begin{equation*}
 c'_{r-1}\rho^{r-1}Q_*^{2r+1}h_{S\cup\{r\}}.
\end{equation*}
For $|B|=b\ge1$ and $|A|=a$, there are at most $\rho$ indices in the matrix product. The induction hypothesis for $D_BC$ gives the entrywise estimate
\begin{align*}
 &\left|\bigl[\mathscr L(D_AE(v_r))D_BC\bigr]_{ij}\right|\\
 &\qquad\le \rho
 \bigl(c'_a\rho^aQ_*^{2a+2}h_A\|v_r\|_\diamond\bigr)
 \bigl(c_b\rho^{b-1}Q_*^{2b+1}h_B\bigr)\\
 &\qquad=c'_ac_b\rho^{r-1}Q_*^{2r+1}h_{S\cup\{r\}},
\end{align*}
since $a+b=r-1$. Summing the $\binom{r-1}a$ allocations for each $a$ closes the induction with, for instance,
\begin{equation*}
 c_r=c'_{r-1}
 +\sum_{\substack{a+b=r-1\\b\ge1}}\binom{r-1}a c'_ac_b.
\end{equation*}
This constant depends on $r$ only.

For the logarithm of a neighbouring Jacobi coefficient, differentiation of the second first variation identity gives, for every nonempty label set $A$,
\begin{equation*}
 |D_A\log b_j|
 \le C_{|A|}\rho^{|A|-1}Q_*^{2|A|}h_A,
 \qquad 1\le j<\rho.
\end{equation*}
As $G_\rho$ is between $\kappa I$ and $\kappa^{-1}I$, the ratios of squared norms of successive monic orthogonal polynomials show $b_j\le\kappa^{-1}$. The ordinary multivariate
chain rule, with $S=\{1,\ldots,r\}$, reads
\begin{equation*}
 D_Sb_j=b_j\sum_{\pi\in\Pi(S)}\prod_{A\in\pi}D_A\log b_j.
\end{equation*}
If a partition has $k$ blocks, the corresponding product is bounded by
$C_r\rho^{r-k}Q_*^{2r}h_S$.
Since $1\le k\le r$, $\rho\ge1$, and there are finitely many partitions depending only on $r$, the sum is at most
$C_{r,\kappa}\rho^{r-1}Q_*^{2r}h_S$.
This proves both derivative bounds, including coinciding directions by the same labelled differentiation identities.
\end{proof}

\begin{lemma}
\label{W:lem:tightness-directions}
Under the same bounds on the Gram matrix, for $1\le j\le K\le\rho-2$ and every fixed $r\ge1$,
\begin{equation}
 |D_s^r\mathcal T_j^\sigma[v_1,\ldots,v_r]|
 \le C_{r,\kappa}B\rho^{r-1}Q_*^{2r}\prod_a\|v_a\|_\diamond.
 \label{W:eq:tightness-integrated-potential}
\end{equation}
\end{lemma}
\begin{proof}
All derivatives are evaluated at the same moment functional and the quantities $Q_*$ and $B$ are used only after differentiation. For a polynomial $f=\sum_i\widehat f_ip_i$ write
$|f|_{\rm c}=\sum_i|\widehat f_i|$. For a label set $A$ put $h_A=\prod_{a\in A}\|v_a\|_\diamond$, with $h_\varnothing=1$. Lemma~\ref{W:lem:polynomial-coefficient-derivatives} and
the fact that $q_i$ contains at most $\rho$ Chebyshev basis coefficients give
\begin{equation*}
 |D_Aq_i|_{\rm c}
 \le C_{|A|,\kappa}\rho^{|A|}Q_*^{2|A|+1}h_A
 \quad(A\ne\varnothing),
 \qquad |q_i|_{\rm c}\le Q_*.
\end{equation*}
The same estimates with an additional factor two apply to $q_j+\sigma q_{j-1}$. For every direction $v$ and two polynomials $f,g$,
\begin{equation*}
 |v(fg)|\le\|v\|_\diamond|f|_{\rm c}|g|_{\rm c}.
\end{equation*}

For $1\le j\le\rho-2$, the exact first variation is $D_s\mathcal T_j^\sigma[v]=v(\mathscr P_j^\sigma)$, with
\begin{equation*}
 \mathscr P_j^\sigma
 =\frac{b_j}{2}(q_j+\sigma q_{j-1})^2
 +\sum_{a=1}^{j-1}(b_a-b_{a+1})q_a^2-b_1q_0^2.
\end{equation*}
Write the displayed sum as $\sum_\alpha c_\alpha u_\alpha^2$. Since $q_0=1$ throughout the normalized domain and $v_r(1)=0$, the term $-b_1q_0^2$ contributes zero. For the other
terms,
\begin{equation*}
 \sum_\alpha|c_\alpha|
 \le\frac12\kappa^{-1}
 +\sum_{a=1}^{j-1}|b_a-b_{a+1}|+\kappa^{-1}
 \le C_\kappa B.
\end{equation*}
Here $j-1\le\rho-3$, as required by the definition of $B$. If $A\ne\varnothing$, then each differentiated coefficient, including one formed as the difference of two $b$
coefficients, satisfies
\begin{equation*}
 |D_Ac_\alpha|
 \le C_{|A|,\kappa}\rho^{|A|-1}Q_*^{2|A|}h_A.
\end{equation*}
There are at most $j+1\le\rho$ such coefficients. These two ways of summing the coefficients will be used separately.

Let $S=\{1,\ldots,r-1\}$. Since $v_r$ is a fixed linear functional, Leibniz' rule gives the exact formula
\begin{equation*}
 D_{S\cup\{r\}}\mathcal T_j^\sigma
 =\sum_\alpha\sum_{A\sqcup U\sqcup V=S}
   (D_Ac_\alpha)
   v_r((D_Uu_\alpha)(D_Vu_\alpha)).
\end{equation*}
For the nonconstant polynomial factors, the first estimates in the proof give, for every allocation of $U$ and $V$,
\begin{equation*}
 |v_r((D_Uu_\alpha)(D_Vu_\alpha))|
 \le C_{r,\kappa}\rho^{|U|+|V|}
 Q_*^{2(|U|+|V|)+2}h_Uh_V\|v_r\|_\diamond.
\end{equation*}
If $A=\varnothing$, multiply this by $|c_\alpha|$ and sum first in $\alpha$. The bound on $\sum_\alpha|c_\alpha|$ makes the result
\begin{equation*}
 C_{r,\kappa}B\rho^{r-1}Q_*^{2r}h_{S\cup\{r\}}.
\end{equation*}
If $|A|=k\ge1$, estimate the derivative of each coefficient and then sum over its index. Since $|U|+|V|=r-1-k$, these terms are bounded by
\begin{align*}
 & C_{r,\kappa}\rho
  \bigl(\rho^{k-1}Q_*^{2k}h_A\bigr)
  \bigl(\rho^{r-1-k}Q_*^{2(r-1-k)+2}
        h_Uh_V\|v_r\|_\diamond\bigr)\\
 &\qquad=C_{r,\kappa}\rho^{r-1}Q_*^{2r}h_{S\cup\{r\}}
 \le C_{r,\kappa}B\rho^{r-1}Q_*^{2r}h_{S\cup\{r\}},
\end{align*}
where the last inequality uses $B\ge1$. In particular, one fewer power of $\rho$ in the bound for $D_Ac_\alpha$ permits summation over up to $\rho$ coefficients.

Every one of the $r-1$ labels is assigned to exactly one of $A,U,V$; there are $3^{r-1}$ assignments. This number does not depend on $j$ or $\rho$. Summing the two estimates proves
\eqref{W:eq:tightness-integrated-potential} uniformly for $1\le j\le K\le\rho-2$ and both signs $\sigma$.
\end{proof}

The factor $\rho^{r-1}$ is needed when several derivatives of $\mathcal T_j^\sigma$ occur in a derivative of its exponential.

For a partition of $r$ derivative labels into $k$ nonempty sets, their dimension factors multiply to $m_N^k\rho^{r-k}=O(m_N^r)$ when $\rho\asymp m_N$.
Proposition~\ref{W:prop:comparison-global} applies this estimate to the exponential together with the two cutoffs.

\subsubsection{A finite sum rule and Hilbert--Schmidt estimates}
\label{W:sec:finite-coefficient-sum-rule}
\label{W:sec:derivatives-of-the-nonlinear-remainder}

The preceding estimates depend on $Q_*$ and $B$ through the coefficients of the orthonormal polynomials. Variations measured directly in Hilbert--Schmidt norm admit a different
bound. A finite trace identity controls the sum of the squared deviations of the Jacobi coefficients independently of the number of coefficients. This bound will control the
nonlinear part of the integrated potential and its derivatives.

The notation and the bounds on the Gram matrix remain those of Section~\ref{W:sec:a-finite-sum-rule-and-potential-derivatives}. The function $q(b)=b-1-\log b$, defined for $b>0$,
occurs through the identity

\[
 b^2-1-2\log b=(b-1)^2+2q(b)
\]
which separates the squared deviation from the nonlinear term $q(b)$. For $1\le n\le\rho-1$, the boundary vector is $c^{(n)}=G_n^{-1}z_n$, with $c^{(1)}_0=0$. The bound and the
nonlinear sum will be denoted by
\[
 C_\kappa=2\kappa^{-2}+\kappa^{-4}+4\log(\kappa^{-1}),\qquad
 \mathcal E_K=2\sum_{j<K}q(b_j)+q(b_K),\quad\mathcal E_0=0.
\]

\begin{theorem}
\label{W:thm:projection-energy}
Under \eqref{W:eq:projection-Gram},
\begin{align}
 \sum_{j=1}^n a_j^2+2\sum_{j<n}(b_j^2-1-2\log b_j)
 &=2c^{(n)}_{n-1}+(c^{(n)}_n)^2
     -2\log\frac{\mathfrak h_{n-1}}{\mathfrak h_0},
     \label{W:eq:projection-sumrule}\\
 \sum_{j<\rho}a_j^2+2\sum_{j<\rho-1}(b_j-1)^2
 &\le C_\kappa,
     \label{W:eq:projection-energybound}\\
 0\le\mathcal E_K&\le C_\kappa/2\qquad(0\le K\le\rho-2).
     \label{W:eq:projection-Ebound}
\end{align}
The coordinates of $c^{(n)}$ are numbered from one, with the first term on the right of \eqref{W:eq:projection-sumrule} equal to zero for $n=1$.
\end{theorem}
\begin{proof}
Fix $1\le n\le\rho-1$. All matrices in this argument are formed from the same moment functional. In the basis $p_0,\ldots,p_{n-1}$, multiplication by $x$ satisfies
\begin{equation*}
 S_n=G_nJ_{0,n}+z_ne_n^{\mathsf T}.
\end{equation*}
Indeed, for the first $n-1$ columns this follows from $xp_j=p_{j+1}+p_{j-1}$ (with $p_{-1}=0$), since both polynomials on the right then belong to the chosen basis. The last column
is $G_ne_{n-1}+z_n$ for $n\ge2$, and is $z_1$ for $n=1$. Thus, writing $c=c^{(n)}=G_n^{-1}z_n$,
\begin{equation*}
 G_n^{-1}S_n=J_{0,n}+ce_n^{\mathsf T}.
\end{equation*}
If $G_n=L_nL_n^{\mathsf T}$, the Jacobi matrix $J_n=L_n^{-1}S_nL_n^{-\mathsf T}$ is similar to $G_n^{-1}S_n$:
\begin{equation*}
 L_n^{\mathsf T}(G_n^{-1}S_n)L_n^{-\mathsf T}=J_n.
\end{equation*}
Similarity preserves the trace of the square, even though the matrix $J_{0,n}+ce_n^{\mathsf T}$ need not be symmetric. Expanding its square and using cyclicity of the finite trace
gives
\begin{align*}
 \operatorname{Tr}J_n^2
 &=\operatorname{Tr}J_{0,n}^2
   +2\operatorname{Tr}(J_{0,n}ce_n^{\mathsf T})
   +\operatorname{Tr}(ce_n^{\mathsf T}ce_n^{\mathsf T})\\
 &=2(n-1)+2e_n^{\mathsf T}J_{0,n}c
          +(e_n^{\mathsf T}c)^2\\
 &=2(n-1)+2c_{n-1}+c_n^2.
\end{align*}
Here $c_0:=0$ when $n=1$. On the other hand, $J_n$ is a Jacobi matrix, and therefore
\begin{equation*}
 \operatorname{Tr}J_n^2
 =\sum_{j=1}^na_j^2+2\sum_{j=1}^{n-1}b_j^2.
\end{equation*}
The squared norms of the monic orthogonal polynomials satisfy $b_j^2=\mathfrak h_j/\mathfrak h_{j-1}$, so their logarithms telescope:
\begin{equation*}
 4\sum_{j=1}^{n-1}\log b_j
 =2\sum_{j=1}^{n-1}
   (\log\mathfrak h_j-\log\mathfrak h_{j-1})
 =2\log\frac{\mathfrak h_{n-1}}{\mathfrak h_0}.
\end{equation*}
Subtracting $2(n-1)$ and this identity from the two expressions for $\operatorname{Tr}J_n^2$ proves \eqref{W:eq:projection-sumrule}, including $n=1$.

We now bound its right hand side without a factor depending on $n$. Because $n+1\le\rho$, the matrix $G_{n+1}$ is a principal block of $G_\rho$, and $z_n$ is obtained by deleting
the last coordinate of its last column. Consequently,
\begin{equation*}
 \|z_n\|_2\le\|G_{n+1}e_{n+1}\|_2
 \le\|G_{n+1}\|_{\rm op}\le\kappa^{-1},
 \qquad \|G_n^{-1}\|_{\rm op}\le\kappa^{-1}.
\end{equation*}
It follows that $\|c^{(n)}\|_2\le\kappa^{-2}$, and hence
\begin{equation*}
 2c_{n-1}+c_n^2\le2\kappa^{-2}+\kappa^{-4}.
\end{equation*}
Moreover, each $\mathfrak h_j$ is the minimum of the positive quadratic form on monic polynomials of degree $j$ in the basis $p_0,\ldots,p_j$. Since
$\kappa I\le G_{j+1}\le\kappa^{-1}I$, minimizing over vectors with last coordinate one and testing the vector $e_{j+1}$ gives
\begin{equation*}
 \kappa\le\mathfrak h_j\le\kappa^{-1},
 \qquad 0\le j\le\rho-1.
\end{equation*}
In particular,
\begin{equation*}
 \left|2\log\frac{\mathfrak h_{n-1}}{\mathfrak h_0}\right|
 \le4\log(\kappa^{-1}).
\end{equation*}
Combining these inequalities with \eqref{W:eq:projection-sumrule} gives an upper bound of $C_\kappa$ on its nonnegative left hand side, uniformly for $1\le n\le\rho-1$.

For $b>0$ set $g(b)=b^2-1-2\log b$. One has $g(1)=g'(1)=0$ and $g''(b)=2+2b^{-2}\ge2$; integration on the interval between $1$ and $b$ yields $g(b)\ge(b-1)^2$. Taking $n=\rho-1$ in
\eqref{W:eq:projection-sumrule} and applying this bound term by term proves \eqref{W:eq:projection-energybound}. This step uses the entire range $1\le j\le\rho-2$ of off-diagonal
terms.

Finally, $q(b)=b-1-\log b\ge0$ by $\log b\le b-1$, and
\begin{equation*}
 g(b)=(b-1)^2+2q(b)\ge2q(b).
\end{equation*}
For $1\le K\le\rho-2$ we can apply the already established sum rule with $n=K+1\le\rho-1$ and obtain
\begin{equation*}
 4\sum_{j=1}^Kq(b_j)
 \le2\sum_{j=1}^Kg(b_j)
 \le C_\kappa.
\end{equation*}
By the definition of $\mathcal E_K$,
\begin{equation*}
 0\le\mathcal E_K
 =2\sum_{j=1}^{K-1}q(b_j)+q(b_K)
 \le2\sum_{j=1}^Kq(b_j)\le\tfrac12C_\kappa.
\end{equation*}
For $K=0$ the sum is empty and $\mathcal E_0=0$. This proves \eqref{W:eq:projection-Ebound} for the entire stated range.
\end{proof}

With the same determinants $D_j$ and last columns $z_k$, define the boundary part
\[
 \mathcal C_k^\sigma=\sigma e_k^{\mathsf T}G_k^{-1}z_k
       +\tfrac12\log(D_{k+1}/D_{k-1})-\log D_1.
\]
The logarithms in this decomposition are used only where the required Gram matrices are positive definite. The inverse and Cholesky factors are later regularized on the full
coordinate space. On the nonzero set of the cutoff, they agree with the original formulas.

\begin{proposition}
\label{W:prop:moment-coordinates-boundary}
For the Jacobi coefficients determined by the positive normalized functional, the first two identities below hold for $1\le j,k\le\rho-1$. The decomposition for $Y_k^\sigma$ uses
$1\le k\le\rho-2$, and the increment identity uses
$1\le k\le\rho-3$:
\begin{align}
 b_j^2&=D_{j+1}D_{j-1}/D_j^2,\qquad
 \sum_{j=1}^ka_j=e_k^{\mathsf T}G_k^{-1}z_k,
 \label{W:eq:moment-coordinates-pivots}\\
 Y_k^\sigma&=-m_N(\mathcal C_k^\sigma+\mathcal E_k),\qquad
 \mathcal E_k\ge0,
 \label{W:eq:moment-coordinates-fullboundary}\\
 \mathcal E_{k+1}-\mathcal E_k&=q(b_k)+q(b_{k+1})\ge0.\notag
\end{align}
\end{proposition}
\begin{proof}
All identities are obtained at a fixed normalized functional for which $G_\rho>0$. Denote by $\pi_j$ the monic orthogonal polynomial of degree $j$, and put
$\eta_j=\mathcal L(\pi_j^2)>0$. The change from $p_0,\ldots,p_{j-1}$ to $\pi_0,\ldots,\pi_{j-1}$ is unitriangular; its determinant is one. Taking determinants in the identity for
the Gram matrix therefore gives
\[
 D_j=\prod_{i=0}^{j-1}\eta_i,
 \qquad \eta_j=\frac{D_{j+1}}{D_j}.
\]
For $1\le j\le\rho-1$, the polynomial $x\pi_{j-1}-\pi_j$ has degree at most $j-1$. Orthogonality therefore gives $\mathcal L(\pi_jx\pi_{j-1})=\eta_j$. The adjacent entry of the
multiplication matrix is
\[
 b_j=\mathcal L\left(\frac{\pi_j}{\sqrt{\eta_j}}\,
                   x\frac{\pi_{j-1}}{\sqrt{\eta_{j-1}}}\right)
     =\sqrt{\frac{\eta_j}{\eta_{j-1}}}>0.
\]
This calculation uses only polynomials through degree $j\le\rho-1$, so it includes $j=\rho-1$. Consequently
\[
 b_j^2=\frac{\eta_j}{\eta_{j-1}}=
 \frac{D_{j+1}D_{j-1}}{D_j^2}.
\]

For $1\le k\le\rho-1$, let $c=G_k^{-1}z_k$ and number its coordinates as $c_0,\ldots,c_{k-1}$. Orthogonal projection of $p_k$ onto the span of $p_0,\ldots,p_{k-1}$ gives
\[
 \pi_k=p_k-\sum_{i=0}^{k-1}c_i p_i.
\]
The polynomial $p_k=U_k(x/2)$ is monic and has no term of degree $k-1$. Also $p_{k-1}$ is monic and all preceding $p_i$ have degree at most $k-2$. Thus the coefficient of $x^{k-1}$
in $\pi_k$ is $-c_{k-1}$. On writing $\pi_k=x^k+v_kx^{k-1}+\cdots$, comparison of the coefficients of degree $k$ in
\[
 \pi_{k+1}=(x-a_{k+1})\pi_k-b_k^2\pi_{k-1},\qquad 1\le k\le\rho-2,
\]
gives $v_{k+1}=v_k-a_{k+1}$, with $v_1=-a_1$. It follows by induction that $-c_{k-1}=v_k=-\sum_{j=1}^k a_j$, or
$\sum_{j=1}^k a_j=e_k^{\mathsf T}G_k^{-1}z_k$.
This proves both parts of \eqref{W:eq:moment-coordinates-pivots}, including $k=1$.

Since $b_j>0$, we may take logarithms. For $k\ge1$,
\begin{align*}
 2\sum_{j=1}^{k-1}\log b_j+\log b_k
 &=\sum_{j=1}^{k-1}(\log\eta_j-\log\eta_{j-1})
       +\tfrac12(\log\eta_k-\log\eta_{k-1})\\
 &=\tfrac12(\log\eta_{k-1}+\log\eta_k)-\log\eta_0\\
 &=\tfrac12\log\frac{D_{k+1}}{D_{k-1}}-\log D_1.
\end{align*}
Here $D_0=1$ and $\eta_0=D_1=\mathcal L(1)=1$; the formula remains valid at $k=1$. Using $b-1=\log b+q(b)$ separately in the terms with weights two and one in
\eqref{W:eq:projection-T}, we get the exact identity
\[
 \mathcal T_k^\sigma
 =\sigma e_k^{\mathsf T}G_k^{-1}z_k
  +\tfrac12\log\frac{D_{k+1}}{D_{k-1}}-\log D_1
  +2\sum_{j=1}^{k-1}q(b_j)+q(b_k)
 =\mathcal C_k^\sigma+\mathcal E_k.
\]
The function $q(b)=b-1-\log b$ has $q(1)=q'(1)=0$ and $q''(b)=b^{-2}>0$, hence $q(b)\ge0$ for all $b>0$. Finally, subtracting the defining sums for $\mathcal E_{k+1}$ and
$\mathcal E_k$ gives
\[
 \mathcal E_{k+1}-\mathcal E_k
 =\left(2\sum_{j=1}^kq(b_j)+q(b_{k+1})\right)
  -\left(2\sum_{j=1}^{k-1}q(b_j)+q(b_k)\right)
 =q(b_k)+q(b_{k+1})\ge0.
\]
The initial difference is $\mathcal E_1-\mathcal E_0=q(b_1)\ge0$; the displayed two-term increment is used for $k\ge1$. All determinants appearing above are positive principal Gram
determinants, so the logarithms and quotients are defined.
\end{proof}

In particular, $\mathcal T_K^\sigma=\mathcal C_K^\sigma+\mathcal E_K$. To compute its second and third derivatives, we differentiate both $\mathcal E_K$ and the orthonormal
polynomials appearing in the first variation.

Let $v,w,z$ be fixed real linear functionals on the polynomial space, each zero on constants. In the current orthonormal basis write their matrices as $E,F,H$, respectively, and
put $A_E=\mathscr L(E)$ and $A_H=\mathscr L(H)$. Define
\begin{equation}
 M_E=\mathcal M_K^\sigma((b_jd_j(E))_{j=1}^K), \qquad N_E=M_E-A_E^{\mathsf T}M-MA_E.
 \label{W:eq:projection-N}
\end{equation}
For the variation of $E$ along $z$, set
\begin{align}
 E_{;H}&=-A_HE-EA_H^{\mathsf T},\nonumber\\
 A_{E;H}&=\mathscr L(E_{;H}),\nonumber\\
 M_{E,H}&=\mathcal M_K^\sigma
 \bigl((b_j\{d_j(E)d_j(H)+d_j(E_{;H})\})_{j=1}^K\bigr),\nonumber\\
 N_{E;H}&=M_{E,H}-A_{E;H}^{\mathsf T}M-A_E^{\mathsf T}M_H
          -M_HA_E-MA_{E;H}\nonumber\\
 &\hspace{9mm}-A_H^{\mathsf T}N_E-N_EA_H.
 \label{W:eq:projection-N3}
\end{align}

\begin{theorem}
\label{W:thm:projection-responses}
With this notation,

\begin{equation}
 D_{s}\mathcal T_K^\sigma[v]=\operatorname{Tr}(ME), \qquad D_{s}^2\mathcal T_K^\sigma[v,w]=\operatorname{Tr}(N_EF), \qquad D_{s}^3\mathcal T_K^\sigma[v,w,z]=\operatorname{Tr}(N_{E;H}F).
 \label{W:eq:integrated-potential-mixed-derivatives}
\end{equation}
The three expressions are symmetric multilinear functions of the direction labels.

\end{theorem}
\begin{proof}
The normalized orthogonal polynomials $q_j$ are smooth functions of the finitely many moments on the open set $G_\rho>0$. For each fixed functional direction $v$, define the
symmetric matrix $E(v)=(v(q_iq_j))_{0\le i,j\le K}$. By Lemma~\ref{W:lem:projection-variation}, $D_v\mathbf q=-A_E\mathbf q$ and
$D_vb_j=b_jd_j(E)$.
Differentiating the definition of $\mathcal T_K^\sigma$ and using \eqref{W:eq:projection-Mdef} therefore gives
\[
 D_v\mathcal T_K^\sigma
 =v(\mathbf q^{\mathsf T}M\mathbf q)=\operatorname{Tr}(ME).
\]
The second equality uses $v(q_iq_j)=E_{ij}$ and the symmetry of $M$ and $E$.

For the second derivative the direction $w$ is held fixed while the functional and its orthonormal polynomials vary along $v$. By the ordinary product rule,
\[
 D_v(\mathbf q^{\mathsf T}M\mathbf q)
 =\mathbf q^{\mathsf T}\bigl(M_E-A_E^{\mathsf T}M-MA_E\bigr)\mathbf q
 =\mathbf q^{\mathsf T}N_E\mathbf q.
\]
Consequently
\[
 D_vD_w\mathcal T_K^\sigma
 =w(\mathbf q^{\mathsf T}N_E\mathbf q)
 =\operatorname{Tr}(N_EF),
 \qquad F=E(w).
\]
Here $D_vM=M_E$, since $M$ is linear in $(b_j)$ and the first variation of $b_j$ is $b_jd_j(E)$.

For the third derivative, fix a third direction $z$ and write $H=E(z)$, $A_H=\mathscr L(H)$. Because the functional $v$ is fixed, the derivative of $E(v)_{ij}=v(q_iq_j)$ is due
only to the two polynomial factors:
\begin{align*}
 D_zE(v)_{ij}
 &=v((D_zq_i)q_j)+v(q_i(D_zq_j))\\
 &=-(A_HE+EA_H^{\mathsf T})_{ij}= (E_{;H})_{ij}.
\end{align*}
In particular $D_zA_E=A_{E;H}$. Differentiating $D_vb_j=b_jd_j(E)$ again gives
\[
 D_zD_vb_j
 =b_j\bigl(d_j(H)d_j(E)+d_j(E_{;H})\bigr),
\]
so $D_zM_E=M_{E,H}$. Likewise $D_zM=M_H$. The complete product derivative is
\begin{align*}
 D_zN_E
 ={}&M_{E,H}-A_{E;H}^{\mathsf T}M
       -A_E^{\mathsf T}M_H-M_HA_E-MA_{E;H}.
\end{align*}
The matrix $F=E(w)$ also depends on the polynomial basis, so $D_zF=-A_HF-FA_H^{\mathsf T}$. Therefore
\begin{align*}
 D_z\operatorname{Tr}(N_EF)
 &=\operatorname{Tr}\bigl((D_zN_E)F
                        -N_EA_HF-N_EFA_H^{\mathsf T}\bigr)\\
 &=\operatorname{Tr}\bigl((D_zN_E-A_H^{\mathsf T}N_E-N_EA_H)F\bigr)\\
 &=\operatorname{Tr}(N_{E;H}F),
\end{align*}
The middle equality uses cyclicity of the trace and gives \eqref{W:eq:projection-N3}. Smoothness on $G_\rho>0$ implies that the mixed derivatives commute, proving symmetry in the
three functionals, also when two or more coincide.
\end{proof}

\Needspace{6\baselineskip}

\begin{theorem}
\label{W:thm:projection-dimbound}
For every normalized real moment vector satisfying \eqref{W:eq:projection-Gram}, for $1\le q\le3$, $K\le {\rho}-2$, and real moment directions $v_1,\ldots,v_q$,
\begin{equation}
 |D_{s}^q\mathcal T_K^\sigma[v_1,\ldots,v_q]|
 \le C_{q,\kappa}\prod_{j=1}^q\|\dot G_{\rho}[v_j]\|_{\rm HS}.
 \label{W:eq:projection-Tbound}
\end{equation}
In particular
\begin{equation}
 \max_{\sigma=\pm}\sup_{K\le {\rho}-2}
 |D_{s}^qY_K^\sigma[v_1,\ldots,v_q]|
 \le {m_N}C_{q,\kappa}\prod_{j=1}^q\|\dot G_{\rho}[v_j]\|_{\rm HS}.
 \label{W:eq:projection-Ybound}
\end{equation}
The constants are uniform in $K$ and $\rho$, and linear interpolation obeys the same bound.
\end{theorem}
\begin{proof}
Fix $K\le\rho-2$ and put
$E=L_{K+1}^{-1}\dot G_{K+1}[v]L_{K+1}^{-\mathsf T}$.
Since $G_{K+1}=L_{K+1}L_{K+1}^{\mathsf T}\ge\kappa I$,
\[
 \|E\|_{\rm HS}
 \le\|L_{K+1}^{-1}\|^2\|\dot G_{K+1}[v]\|_{\rm HS}
 \le\kappa^{-1}\|\dot G_\rho[v]\|_{\rm HS}.
\]
The same estimate applies to the matrices $F,H$ corresponding to other directions. We prove bounds in terms of $\|E\|_{\rm HS}$, $\|F\|_{\rm HS}$, and $\|H\|_{\rm HS}$ before
returning to the original variations of the Gram matrix.

For $K\ge2$, the diagonal entries of $\mathcal M_K^\sigma(u)$, numbered from zero, are $-u_1$, $u_j-u_{j+1}$ for $1\le j\le K-2$, $u_{K-1}-u_K/2$, and $u_K/2$; its only
off-diagonal entries are the symmetric pair $\sigma u_K/2$. Thus
\[
 \|\mathcal M_K^\sigma(u)\|_{\rm HS}^2
 =u_1^2+\sum_{j=1}^{K-2}(u_j-u_{j+1})^2
     +(u_{K-1}-u_K/2)^2+\tfrac34u_K^2
 \le6\sum_{j=1}^K u_j^2.
\]
For $K=1$ the two diagonal entries are $-u_1/2,u_1/2$ and the off-diagonal entries are $\sigma u_1/2$; its squared Hilbert--Schmidt norm equals $u_1^2$, so the same upper bound
holds. Moreover, direct substitution gives $\|\mathcal M_K^\sigma(\mathbf1)\|_{\rm HS}\le\sqrt2$ for all $K\ge1$. Hence the energy estimate \eqref{W:eq:projection-energybound},
which contains every $b_j$ with $j\le K$, implies
\[
 \|M\|_{\rm HS}
 \le\sqrt2+\sqrt6\left(\sum_{j=1}^K(b_j-1)^2\right)^{1/2}
 \le\sqrt2+\sqrt{3C_\kappa}=:C_\kappa'.
\]
It follows from the first identity in \eqref{W:eq:integrated-potential-mixed-derivatives} that $|D_v\mathcal T_K^\sigma|=|\operatorname{Tr}(ME)| \le C_\kappa'\|E\|_{\rm HS}$.

The definition of $\mathscr L$ gives $\|\mathscr L(E)\|_{\rm HS}\le\|E\|_{\rm HS}$. Also, for $d_j(E)=(E_{jj}-E_{j-1,j-1})/2$,
\[
 \|d(E)\|_2
 \le\frac12(\|(E_{jj})_{j=1}^K\|_2+
                    \|(E_{j-1,j-1})_{j=1}^K\|_2)
 \le\|E\|_{\rm HS}.
\]
Since $b_j\le\kappa^{-1}$ and $\|A_E\|_{\rm HS}\le\|E\|_{\rm HS}$, the preceding matrix bound gives
\begin{align*}
 \|M_E\|_{\rm HS}&\le\sqrt6\kappa^{-1}\|d(E)\|_2
                  \le\sqrt6\kappa^{-1}\|E\|_{\rm HS},\\
 \|N_E\|_{\rm HS}
 &\le\|M_E\|_{\rm HS}
     +2\|M\|\|A_E\|_{\rm HS}
 \le (\sqrt6\kappa^{-1}+2C_\kappa')\|E\|_{\rm HS}.
\end{align*}
The second identity of \eqref{W:eq:integrated-potential-mixed-derivatives} and Hilbert--Schmidt Cauchy--Schwarz now give
$|D_vD_w\mathcal T_K^\sigma| \le C_\kappa\|E\|_{\rm HS}\|F\|_{\rm HS}$.

For the third derivative, recall $A_H=\mathscr L(H)$ and $E_{;H}=-A_HE-EA_H^{\mathsf T}$. Schatten norm inequalities yield
\[
 \|E_{;H}\|_{\rm HS}
 \le2\|A_H\|\|E\|_{\rm HS}
 \le2\|H\|_{\rm HS}\|E\|_{\rm HS}.
\]
Here $d(E)d(H)$ denotes pointwise multiplication. Thus
\[
 \|d(E)d(H)\|_2
 \le\|d(E)\|_2\|d(H)\|_\infty
 \le\|E\|_{\rm HS}\|H\|_{\rm HS},
 \quad
 \|d(E_{;H})\|_2\le2\|E\|_{\rm HS}\|H\|_{\rm HS}.
\]
Using the definition of $M_{E,H}$ and $b_j\le\kappa^{-1}$,
\[
 \|M_{E,H}\|_{\rm HS}
 \le3\sqrt6\kappa^{-1}\|E\|_{\rm HS}\|H\|_{\rm HS}.
\]
The remaining terms in $N_{E;H}$ satisfy
\begin{align*}
 \|A_{E;H}^{\mathsf T}M\|_{\rm HS}+
 \|MA_{E;H}\|_{\rm HS}
 &\le2C_\kappa'\|E_{;H}\|_{\rm HS},\\
 \|A_E^{\mathsf T}M_H\|_{\rm HS}+
 \|M_HA_E\|_{\rm HS}
 &\le2\|E\|_{\rm HS}\|M_H\|_{\rm HS},\\
 \|A_H^{\mathsf T}N_E\|_{\rm HS}+
 \|N_EA_H\|_{\rm HS}
 &\le2\|H\|_{\rm HS}\|N_E\|_{\rm HS}.
\end{align*}
Each right hand side is bounded by $C_\kappa\|E\|_{\rm HS}\|H\|_{\rm HS}$ using the preceding estimates. Together with the bound on $M_{E,H}$ this proves
\[
 \|N_{E;H}\|_{\rm HS}
 \le C_\kappa\|E\|_{\rm HS}\|H\|_{\rm HS}.
\]
The third identity of \eqref{W:eq:integrated-potential-mixed-derivatives} yields $|D_vD_wD_z\mathcal T_K^\sigma| \le C_\kappa\|E\|_{\rm HS}\|F\|_{\rm HS}\|H\|_{\rm HS}$. Applying
the first bound for the variation of the Gram matrix separately for each direction proves \eqref{W:eq:projection-Tbound} for $q=1,2,3$, with constants independent of $K,\rho$.
Since $Y_K^\sigma=-m_N\mathcal T_K^\sigma$, the same bounds multiplied by $m_N$ prove \eqref{W:eq:projection-Ybound}. For linear interpolation, the value and every moment
derivative are a convex combination of those at the two adjacent indices. The same uniform constant bounds their combination.
\end{proof}

We now bound the derivatives of $\mathcal E_K$ at each fixed order in terms of the Hilbert--Schmidt norms of $\dot G_\rho[v]$. Let $\ell_j=\log b_j$, let $w_j=2$ for $j<K$ and
$w_K=1$, and write $h_S=\prod_{a\in S}\|\dot G_\rho[v_a]\|_{\rm HS}$. For a nonempty label set $S$, differentiation of $G=LL^{\mathsf T}$ with respect to
$\mathcal L+\sum_{a\in S}t_av_a$, with the functionals $v_a$ fixed, gives
\[
 D_SL=L\mathscr L\!\left(L^{-1}
 \left[\mathbf1_{\{|S|=1\}}\dot G[v_S]
 -\sum_{\varnothing\ne A\subsetneq S}(D_AL)(D_{S\setminus A}L)^{\mathsf T}
 \right]L^{-\mathsf T}\right),
\]
where $v_S$ is used only when $S$ is a singleton. Both $L$ and its inverse have bounded operator norm in terms of $\kappa$, and $\mathscr L$ contracts Hilbert--Schmidt norm. Hence
\[
 \|D_SL\|_{\rm HS}\le C_\kappa\left(
 \mathbf1_{\{|S|=1\}}h_S+
 \sum_{\varnothing\ne A\subsetneq S}
       \|D_AL\|_{\rm HS}\|D_{S\setminus A}L\|_{\rm HS}\right).
\]
With $c_1=C_\kappa$ and $c_r=C_\kappa\sum_{a=1}^{r-1}\binom ra c_ac_{r-a}$, induction yields $\|D_SL\|_{\rm HS}\le c_{|S|}h_S$. This induction introduces no factor depending on the
matrix dimension.

The diagonal entries $x_j=L_{jj}$ lie in a fixed positive interval. The labelled chain rule reads
\[
 D_S\log x_j=\sum_{\pi\in\Pi(S)}
 (-1)^{|\pi|-1}(|\pi|-1)!x_j^{-|\pi|}
                         \prod_{A\in\pi}D_Ax_j.
\]
In each product of diagonal vectors, use one $\ell^2$ norm and $\ell^\infty$ norms for the other factors. Since $\|(D_Ax_j)_j\|_2\le\|D_AL\|_{\rm HS}$ and $\ell_j$ is the
difference of two adjacent $\log x_j$, this gives
\begin{equation}
 \|(D_S\ell_j)_j\|_2\le C_{|S|,\kappa}h_S.
 \label{W:eq:log-offdiagonal-derivative-l2}
\end{equation}
Now $q(e^\ell)=e^\ell-1-\ell$, so the exact derivative is
\[
 D_S\mathcal E_K=\sum_{j=1}^K w_j\left[
 (b_j-1)D_S\ell_j+
 b_j\sum_{\substack{\pi\in\Pi(S)\\|\pi|\ge2}}
                             \prod_{A\in\pi}D_A\ell_j\right].
\]
The energy bound gives $\|(b_j-1)_{j\le K}\|_2\le C_\kappa$. The first sum is therefore bounded by $2\|b-1\|_2\|(D_S\ell_j)_j\|_2$. For a partition with $k\ge2$ blocks
$A_1,\ldots,A_k$, use two $\ell^2$ norms:
\[
 \sum_j\prod_{i=1}^k|D_{A_i}\ell_j|
 \le\|(D_{A_1}\ell_j)_j\|_2\|(D_{A_2}\ell_j)_j\|_2
       \prod_{i=3}^k\|(D_{A_i}\ell_j)_j\|_\infty
 \le C_{|S|,\kappa}h_S.
\]
Since $\|b\|_\infty\le\kappa^{-1}$, summing the finitely many partitions proves
\[
 |D_s^r\mathcal E_K[v_1,\ldots,v_r]|
 \le C_{r,\kappa}\prod_{a=1}^r\|\dot G_\rho[v_a]\|_{\rm HS},
 \qquad K\le\rho-2.
\]

At the moment functional of the semicircle measure $D_s\mathcal E_K=0$. The second differential is
\begin{equation}
 D_s^2\mathcal E_K[v,w]
 =\frac12\sum_{j=1}^{K-1}v(p_{2j})w(p_{2j})
       +\frac14v(p_{2K})w(p_{2K}).
 \label{W:eq:second-derivative-of-the-bond-remainder}
\end{equation}
Indeed, \eqref{W:eq:projection-blog} and $p_j^2-p_{j-1}^2=p_{2j}$ give $D_s\ell_j[v]=v(p_{2j})/2$ at this point. In the second derivative only the product of the two first
variations of $\ell_j$ remains. Multiplication by $-m_N$ gives its contribution to $D_s^2Y_K^\sigma$.

\subsubsection{Weighted potentials and interpolation}
\label{W:sec:integrated-potential-boundary-weights}
\label{W:sec:boundary-terms-in-the-integrated-potential}

Values and increments of the integrated potential are linear functions of the Jacobi coefficients. Their weights come from the discrete quadratic form. We describe them for
arbitrary Hermitian matrices as well, since the smooth extensions used in the comparison need not be tridiagonal.

Let \(m_N=N^{1/3}\) and suppose \({\rho}\ge2.\) For real coefficient vectors \(c=(c_1,\ldots,c_{\rho})\in\mathbb R^{\rho}\) and
\(d=(d_1,\ldots,d_{{\rho}-1})\in\mathbb R^{{\rho}-1},\) define, for an arbitrary Hermitian matrix
\(J\in\operatorname{Herm}_{\rho}(\mathbb F_\beta)\),
\begin{equation}
 \mathscr Z_{N,\rho}(c,d;J)
 =-m_N\sum_{k=1}^{\rho}c_kJ_{kk}
       +2m_N\sum_{k=1}^{\rho-1}d_k(1-\operatorname{Re}J_{k,k+1}).
\label{W:eq:weighted-jacobi-functional}
\end{equation}
This formula is real for every Hermitian $J$. For the multiplication matrix of the polynomial basis,
\(J_{\rho}(H)_{k,k+1}=b_k>0,\) so it gives the usual off-diagonal coefficient.

Whenever \(G_{\rho}(H)>0\), set
\begin{equation}
Z_{N,{\rho}}(c,d;H)
:=
\mathscr Z_{N,{\rho}}(c,d;J_{\rho}(H)).
\label{W:eq:weighted-lanczos-observable}
\end{equation}
Then
\[
Z_{N,{\rho}}(c,d;H)
=
-m_N\sum_{k=1}^{{\rho}}c_ka_k
+
2m_N\sum_{k=1}^{{\rho}-1}d_k(1-b_k).
\]

For \(t\ge0\), let \(K_N(t):=\lfloor m_Nt\rfloor.\) If \(K_N(t)+1\le {\rho},\) define
\[
c_k^{(t)}
=
\mathbf 1_{\{k\le K_N(t)\}},
\qquad
d_k^{(t)}
=
\mathbf 1_{\{k\le K_N(t)\}}.
\]
Then the cumulative potential obtained from the tridiagonal coefficients is
\begin{equation}
Y_N(t;H)
=
Z_{N,{\rho}}(c^{(t)},d^{(t)};H).
\label{W:eq:potential-as-weighted-functional}
\end{equation}

The potential in the quadratic form instead uses $d_r=(c_r+c_{r+1})/2$ with $c_\rho=0$. Lemma~\ref{W:lem:coefficients-weighted-identity} derives these weights; the example below
explains why the difference from $Y_N$ matters for convergence.

Write $h_N=m_N^{-1}$ and $\rho=\lfloor m_NT\rfloor+r_0$, where $r_0\ge3$ is fixed. In this regime of increasing dimension we consider sufficiently large $N$ that $\rho+1\le N$,
equivalently $\rho<N$. This holds because $\rho/N\to0$ at fixed $T,r_0$. In particular, the finite spectral model in Proposition~\ref{W:prop:observable-Fourier} can be constructed
inside the original dimension $N$. We henceforth restrict to these dimensions. Let $J$ be any Hermitian ${\rho}\times {\rho}$ matrix, set $a_r=J_{rr}$ and
$b_r=\operatorname{Re}J_{r,r+1}$, and put $b_0=1$. The algebraic identities below hold for real $b_r$, including zero values. Define
\begin{equation}
 P_0(J)=0,\qquad
 P_k(J)={m_N}\sum_{r=1}^k(2-a_r-b_{r-1}-b_r),
 \qquad 1\le k<{\rho}.
 \label{W:eq:coefficients-integrated-potential-on-J}
\end{equation}
For the leading block of a Lanczos tridiagonal matrix these are the exact values $\widetilde Y_{N,k}$. Applying the same formula to the smooth matrix extension defines the extended
integrated potential. The exact identity for the boundary term \eqref{W:eq:integrated-potential-difference} is multiplied by ${m_N}$, whereas ordinary convergence of the
off-diagonal coefficients only gives $b_{k,N}-1=o(1)$. Thus changing the weight of the final off-diagonal coefficient from $1/2$ to $1$ requires a stronger estimate than the
operator criterion assumes. The following deterministic model makes the possible order one error explicit.

Take a nonzero bounded continuous $f$ of compact support with $f(0)=0$, and set, on each fixed window, $b_{k,N}=1+{h_N} f(k{h_N})$ and
$a_{k,N}=2-\bar b_{k-1,N}-b_{k,N} -{h_N}[V(k{h_N})-V((k-1){h_N})]$, where $V(0)=0$ is continuous. For large $N$ these off-diagonal coefficients are positive. Direct summation gives
\begin{equation}
 \widetilde Y_{N,k}=V(k{h_N}),\qquad
 Y_{N,k}=V(k{h_N})-f(k{h_N}),\qquad
 \max_{k\le {m_N}T}|b_{k,N}-1|\le {h_N}\|f\|_\infty.
 \label{W:eq:coefficients-boundary-identity-example}
\end{equation}
The two integrated potentials have limits $V$ and $V-f$, respectively. The last off-diagonal coefficient is therefore kept with its weight throughout the comparison.

The discrete form determines the coefficients in the following weighted sums. For $q=(q_1,\ldots,q_{\rho})\in\mathbb R^{\rho}$ with $q_{\rho}=0$, let
\begin{equation}
 (A_{\rho}q)_r=\frac{q_r+q_{r+1}}2,\qquad 1\le r<{\rho}.
 \label{W:eq:coefficients-offdiagonal-average}
\end{equation}
The averaging map $A_\rho$ acts on the deterministic weights $q$.

\Needspace{6\baselineskip}

\begin{lemma}
\label{W:lem:coefficients-weighted-identity}
For every $J$ and $q$ as above,
\begin{equation}
 \sum_{r=1}^{{\rho}-1}q_r\bigl(P_r(J)-P_{r-1}(J)\bigr)
 =\mathscr Z_{N,{\rho}}(q,A_{\rho}q;J).
 \label{W:eq:coefficients-weighted-increments}
\end{equation}
In particular, with $q_r^k=\mathbf1_{\{r\le k\}}$ and $1\le k<{\rho}$,
\begin{equation}
 P_k(J)=-{m_N}\sum_{r=1}^ka_r
       +2{m_N}\sum_{r=1}^{k-1}(1-b_r)+{m_N}(1-b_k).
 \label{W:eq:integrated-potential-weighted-sum}
\end{equation}
For $0\le\ell<k<{\rho}$, the coefficients for $P_k-P_\ell$ are $q_r=\mathbf1_{\{\ell<r\le k\}}$ and $d=A_{\rho}q$. Thus $d_k=1/2$, $d_\ell=1/2$ when $\ell\ge1$, and $d_r=1$ for
$\ell<r<k$; all other off-diagonal coefficients vanish.
\end{lemma}

\begin{proof}
Write $b_0=1$ and $q_\rho=0$. By the definition of $P_r$,
\begin{align*}
 \sum_{r=1}^{\rho-1}q_r(P_r-P_{r-1})
 &=-m_N\sum_{r=1}^{\rho-1}q_ra_r
   +m_N\sum_{r=1}^{\rho-1}q_r(1-b_{r-1})
   +m_N\sum_{r=1}^{\rho-1}q_r(1-b_r).
\end{align*}
The first of the two sums involving $b$ is, after changing $r$ to $r+1$,
\[
 \sum_{r=1}^{\rho-1}q_r(1-b_{r-1})
 =q_1(1-b_0)+\sum_{r=1}^{\rho-2}q_{r+1}(1-b_r)
 =\sum_{r=1}^{\rho-1}q_{r+1}(1-b_r).
\]
We used $b_0=1$ and the convention $q_\rho=0$ in the last equality. Therefore
\begin{align*}
 \sum_{r=1}^{\rho-1}q_r(P_r-P_{r-1})
 &=-m_N\sum_{r=1}^\rho q_ra_r
   +m_N\sum_{r=1}^{\rho-1}(q_r+q_{r+1})(1-b_r)\\
 &=\mathscr Z_{N,\rho}(q,A_\rho q;J),
\end{align*}
which proves \eqref{W:eq:coefficients-weighted-increments}.

For $1\le k<\rho$ and $q_r=\mathbf1_{\{r\le k\}}$, $(A_\rho q)_r$ is one for $r<k$, one half for $r=k$, and zero for $r>k$. Substitution gives
\[
 P_k=-m_N\sum_{r=1}^ka_r
       +2m_N\sum_{r=1}^{k-1}(1-b_r)+m_N(1-b_k),
\]
including $k=1$ (the interior sum is empty). If $q_r=\mathbf1_{\{\ell<r\le k\}}$ with $1\le\ell<k$, then $(A_\rho q)_\ell=(A_\rho q)_k=1/2$ and the interior averages equal one. For
$\ell=0$, the index $r=0$ is absent and the corresponding term would be $m_N(1-b_0)=0$. Thus the same formula represents the difference $P_k-P_\ell$ for every $0\le\ell<k<\rho$.
\end{proof}

For $t\in[0,T]$, write $k=\lfloor {m_N}t\rfloor$ and $\theta={m_N}t-k$. Define
\begin{equation}
 q_r^{\mathrm{st},t}=\mathbf1_{\{r\le k\}}, \qquad q_r^{\mathrm{lin},t}=\mathbf1_{\{r\le k\}} +\theta\mathbf1_{\{r=k+1\}}.
 \label{W:eq:coefficients-evaluation-weights}
\end{equation}
The inequality $r_0\ge3$ gives $k+1<\rho$, so the last coordinates of these vectors are zero. Lemma~\ref{W:lem:coefficients-weighted-identity} and linear interpolation give
\begin{equation}
 P^{\bullet}(t;J)
 =\mathscr Z_{N,{\rho}}(q^{\bullet,t},A_{\rho}q^{\bullet,t};J),
 \qquad \bullet\in\{\mathrm{st},\mathrm{lin}\}.
 \label{W:eq:coefficients-evaluation-identity}
\end{equation}
For $k\ge1$, the nonzero weights on the off-diagonal coefficients are $1$ for $r<k$, $(1+\theta)/2$ for $r=k$, and $\theta/2$ for $r=k+1$. For $k=0$ only $d_1=\theta/2$ remains. At
$\theta=1$ they agree with the formulas at the next node with $\theta=0$. At $t=0$ every weight is zero, so $P^{\mathrm{lin}}(0;J)=0$.

\begin{figure}[H]
\centering
\begin{tikzpicture}[x=1cm,y=1cm,every node/.style={font=\footnotesize,inner sep=2pt}]
 \node[anchor=east] at (2.1,1.9) {$P_k(J)$};
 \node[anchor=east] at (2.1,-.35) {$P^{\mathrm{lin}}(t;J)$};
 \foreach \yy in {1.9,-.35}{
  \draw[densely dashed] (2.55,\yy)--(3.15,\yy);
  \draw (3.15,\yy)--(12.45,\yy);
  \draw[densely dashed] (12.45,\yy)--(13.05,\yy);
  \foreach \x/\lab in {3.15/{k-1},6.25/k,9.35/{k+1},12.45/{k+2}}{
    \node[below=7pt] at (\x,\yy) {$\lab$};
  }
 }
 \foreach \x/\q in {3.15/1,6.25/1}{
  \fill (\x,1.9) circle (2pt); \node[above=7pt] at (\x,1.9) {$\q$};
  \fill (\x,-.35) circle (2pt); \node[above=7pt] at (\x,-.35) {$\q$};
 }
 \foreach \x in {9.35,12.45}{\draw[fill=white] (\x,1.9) circle (2pt);\node[above=7pt] at (\x,1.9) {$0$};}
 \draw[fill=black!30] (9.35,-.35) circle (2.5pt);\node[above=7pt] at (9.35,-.35) {$\theta$};
 \draw[fill=white] (12.45,-.35) circle (2pt);\node[above=7pt] at (12.45,-.35) {$0$};
 \foreach \x/\d in {4.7/1,7.8/{\frac12},10.9/0}{\node[fill=white,above=1pt] at (\x,1.9) {$\d$};}
 \foreach \x/\d in {4.7/1,7.8/{\frac{1+\theta}{2}},10.9/{\frac\theta2}}{\node[fill=white,above=1pt] at (\x,-.35) {$\d$};}
 \node at (7.55,3.05) {Vertex weights $q_r$, edge weights $d_r=(q_r+q_{r+1})/2$};
 \node at (7.55,-1.5) {$\mathscr Z_{N,\rho}(q,d;J)=-m_N\sum_rq_ra_r+2m_N\sum_rd_r(1-b_r)$};
 \node at (7.55,-2.15) {$t=(k+\theta)/m_N,\qquad 0\le\theta<1$};
\end{tikzpicture}
\caption{The weights in \eqref{W:eq:coefficients-evaluation-weights} and \eqref{W:eq:coefficients-evaluation-identity}, shown for $k\ge2$. A vertex carries the weight $q_r$ in $-m_Nq_ra_r$, and an edge carries the average $d_r$ in $2m_Nd_r(1-b_r)$. In the first row the last nonzero edge weight is $1/2$, producing $m_N(1-b_k)$ in \eqref{W:eq:integrated-potential-weighted-sum}. In the second row the weights at $k$ and $k+1$ are $(1+\theta)/2$ and $\theta/2$. At the left boundary the term with $b_0=1$ vanishes; for $k=0$ only $d_1=\theta/2$ remains.}
\label{W:fig:integrated-potential-weights}
\end{figure}
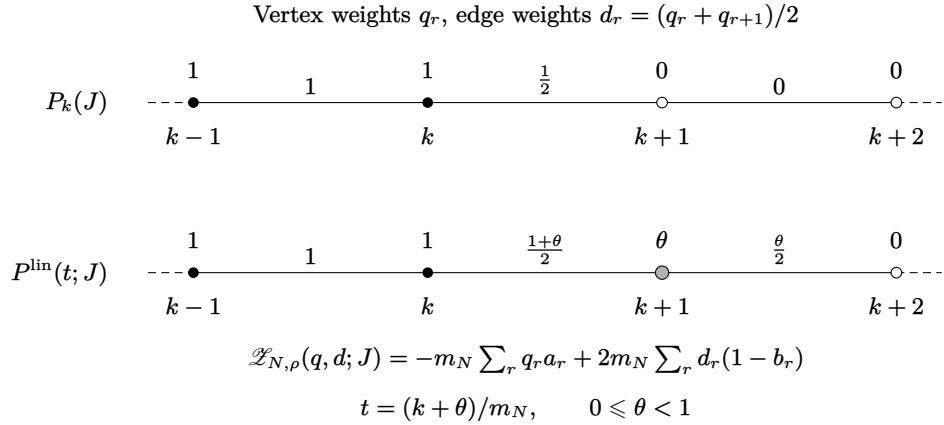

Figure~\ref{W:fig:integrated-potential-weights} shows the weights near the last included index. The averaging formula determines both the half weight in $P_k$ and the two changing
weights during linear interpolation.

Fix a positive integer $J_0$ and $\Lambda<\infty$. The coefficient class used for the following weighted potentials consists precisely of
\begin{equation}
 q=\sum_{j=1}^{J_0}\lambda_jq^{\bullet_j,t_j},\qquad d=A_{\rho}q,
 \qquad t_j\in[0,T],\quad \sum_{j=1}^{J_0}|\lambda_j|\le\Lambda.
 \label{W:eq:coefficients-coefficient-class}
\end{equation}
Each $\bullet_j$ can be step or linear, and the times and weights may depend on $N$. For these weights,
\[
 \sum_{j=1}^{J_0}\lambda_j P^{\bullet_j}(t_j;J)
       =\mathscr Z_{N,\rho}(q,A_\rho q;J).
\]
Thus the joint characteristic function is evaluated by the single exponential of this weighted statistic.

\Needspace{6\baselineskip}

For a finite vector $v$, write
\[
 \|v\|_{p,m_N}=\left(h_N\sum_r|v_r|^p\right)^{1/p}\quad(1\le p<\infty).
\]
The linear dependence on $J$ is represented by the matrix
\[
 \mathsf A_N(q,d)=-m_N\sum_rq_rE_{rr}
       -m_N\sum_{r<\rho}d_r(E_{r,r+1}+E_{r+1,r}).
\]

\begin{lemma}
\label{W:lem:integrated-potential-weight-bounds}
For $q_\rho=0$, averaging contracts $\|\cdot\|_{p,m_N}$ and $\|\cdot\|_\infty$:
\[
 \|A_\rho q\|_{p,m_N}\le\|q\|_{p,m_N},\qquad
 \|A_\rho q\|_\infty\le\|q\|_\infty.
\]
For the weights in \eqref{W:eq:coefficients-coefficient-class},
\begin{equation}
 \begin{aligned}
 \|q\|_\infty+\|A_\rho q\|_\infty&\le2\Lambda,\\
 \|q\|_{2,m_N}^2+\|A_\rho q\|_{2,m_N}^2&\le2\Lambda^2(T+h_N),\\
 \sum_{r<\rho}|q_{r+1}-q_r|&\le\Lambda,\\
 \|A_\rho q-(q_1,\ldots,q_{\rho-1})\|_{2,m_N}^2&\le h_N\Lambda^2/4.
 \end{aligned}
 \label{W:eq:coefficients-coefficient-bounds}
\end{equation}
\end{lemma}

\begin{proof}
Let $q_\rho=0$. For $p\ge1$, convexity of $x\mapsto|x|^p$ gives, for each $1\le r<\rho$,
\[
 \left|\frac{q_r+q_{r+1}}2\right|^p
 \le\frac{|q_r|^p+|q_{r+1}|^p}{2}.
\]
Summing the inequality over $r$ and including the omitted first and last contributions gives
\[
 \sum_{r=1}^{\rho-1}|(A_\rho q)_r|^p
 \le\frac12\sum_{r=1}^{\rho-1}|q_r|^p
       +\frac12\sum_{r=2}^{\rho}|q_r|^p
 \le\sum_{r=1}^{\rho}|q_r|^p.
\]
Multiplication by $h_N$ and taking the $p$th root proves
$\|A_\rho q\|_{p,m_N}\le\|q\|_{p,m_N}$.
The estimate $|(A_\rho q)_r|\le\max(|q_r|,|q_{r+1}|)$ proves contraction in the supremum norm.

Every step weight $q^{\mathrm{st},t}$ and interpolated weight $q^{\mathrm{lin},t}$ from \eqref{W:eq:coefficients-evaluation-weights} is nonnegative, nonincreasing in $r$, bounded
by one, and supported on no more than $\lfloor m_NT\rfloor+1$ indices. Because its last entry is zero, the sum of the absolute successive differences is its first entry, hence is
at most one. For $q=\sum_{j=1}^{J_0}\lambda_jq^{\bullet_j,t_j}$ and $\sum_j|\lambda_j|\le\Lambda$, the triangle inequality yields
\begin{align*}
 \|q\|_\infty&\le\Lambda,\\
 \|q\|_{2,m_N}&\le\sum_j|\lambda_j|
      \sqrt{h_N(\lfloor m_NT\rfloor+1)}
      \le\Lambda\sqrt{T+h_N},\\
 \sum_{r<\rho}|q_{r+1}-q_r|
 &\le\sum_j|\lambda_j|
      \sum_{r<\rho}|q_{r+1}^{\bullet_j,t_j}-q_r^{\bullet_j,t_j}|
 \le\Lambda.
\end{align*}
The contraction just proved supplies identical upper bounds for $A_\rho q$ in the two norms. Also,
\begin{align*}
 \|A_\rho q-(q_1,\ldots,q_{\rho-1})\|_{2,m_N}^2
 &=\frac{h_N}{4}\sum_{r<\rho}|q_{r+1}-q_r|^2\\
 &\le\frac{h_N}{4}
      \left(\sum_{r<\rho}|q_{r+1}-q_r|\right)^2
 \le\frac{h_N\Lambda^2}{4}.
\end{align*}
These estimates give all four inequalities in \eqref{W:eq:coefficients-coefficient-bounds}, uniformly in $N$, the times and the choice of step or linear interpolation.
\end{proof}

For $0\le s\le t\le T$, the identity $q_r^{\mathrm{lin},t}={m_N}\,|(0,t]\cap((r-1){h_N},r{h_N}]|$ implies that $q=q^{\mathrm{lin},t}-q^{\mathrm{lin},s}$ has nonnegative coordinates
with sum ${m_N}(t-s)$. Each is at most both one and that sum; hence ${h_N}\sum q_r^2\le {h_N}\max_rq_r\sum q_r$, proving the bounds
\begin{equation}
 0\le q_r\le\min\{1,m_N(t-s)\},\quad
 \sum_rq_r=m_N(t-s),\quad
 \|q\|_{2,m_N}^2\le(t-s)\min\{1,m_N(t-s)\}.
 \label{W:eq:coefficients-increment-coefficients}
\end{equation}
Expanding $\operatorname{Tr}(\mathsf A_N(q,d)J)$ yields
\begin{equation}
 \begin{aligned}
 \mathscr Z_{N,\rho}(q,d;J)
   &=2m_N\sum_{r<\rho}d_r+\operatorname{Tr}(\mathsf A_N(q,d)J),\\
 \|\mathsf A_N(q,A_\rho q)\|&\le3m_N\Lambda,\\
 \|\mathsf A_N(q,A_\rho q)\|_{\mathrm{HS}}^2
   &=N\bigl(\|q\|_{2,m_N}^2+2\|A_\rho q\|_{2,m_N}^2\bigr)
     \le3N\Lambda^2(T+h_N).
 \end{aligned}
 \label{W:eq:bounds-for-the-evaluation-weights}
\end{equation}
Indeed, the maximal absolute row sum of $\mathsf A_N(q,d)$ is bounded by ${m_N}(\|q\|_\infty+2\|d\|_\infty)$. Its squared Hilbert--Schmidt norm is
${m_N}^2(\sum q_r^2+2\sum d_r^2)$, which equals the displayed expression because ${m_N}^3=N$. Substitution of \eqref{W:eq:coefficients-coefficient-bounds} proves the two matrix
bounds.

\subsubsection{Reflection and the modulus of continuity}
\label{W:sec:modulus-tests-for-the-grid-values-vector}

\label{W:sec:reflection-and-the-modulus-of-continuity}

The integrated potential combines each diagonal coefficient with two neighbouring off-diagonal coefficients. The constructions for $H$ and $-H$ separate these contributions:
reflection changes the diagonal signs and preserves the neighbouring coefficients. After proving this relation, we recover the coefficients from the two potential increments and
use tightness to estimate their maxima. The reflection calculation also includes the continuation when $r_k=0$.

To transfer tightness from the Gaussian ensemble, we use a bounded smooth function of all short differences of the interpolated path. The derivatives of the logarithmic
approximation to the maximum are bounded in terms of the logarithm of the number of pairs. The modulus comparison is then combined with the coefficient bounds in
Section~\ref{W:sec:convergence-of-the-integrated-potential}.

We continue to use the scales ${m_N}=N^{1/3}$ and ${h_N}={m_N}^{-1}$, and the Lanczos continuation is the one in \eqref{W:eq:lanczos-tridiagonal-matrix}. The two paths are denoted
by $X_N^\sigma(t)=\widetilde Y_N^{\mathrm{lin}}(t;\sigma H_N,w_N)$, $\sigma\in\{+,-\}$, with $+H_N=H_N$. Both start at zero. Their continuity is measured by
\[
 \omega_T(x,\delta)=
 \sup_{\substack{s,t\in[0,T]\\|t-s|\le\delta}}|x(t)-x(s)|.
\]
All estimates involving a deterministic starting vector hold for each deterministic sequence, or uniformly in the probability sense of the convention
\eqref{W:eq:uniformity-in-the-starting-vector}. The constants in these probability bounds are uniform in the deterministic unit vector.

\Needspace{6\baselineskip}

\begin{lemma}
\label{W:lem:tightness-reflection}
For every Hermitian matrix and the specified deterministic continuation at breakdown,
\begin{equation}
 a_k(-H,w)=-a_k(H,w),\qquad b_k(-H,w)=b_k(H,w),\quad 1\le k<N.
 \label{W:eq:tightness-reflection}
\end{equation}
The diagonal identity also holds at $k=N$, and the breakdown indices agree. The diagonal sign matrix $D_{\rho}\in\mathbb R^{\rho\times\rho}$, defined by
$(D_{\rho})_{jj}=(-1)^{j-1}$, satisfies the identities for the Gram matrix and the matrix representing multiplication by the variable
\begin{equation}
 G_{\rho}(-H)=D_{\rho}G_{\rho}(H)D_{\rho},\qquad
 M_{\rho}(-H)=-D_{\rho}M_{\rho}(H)D_{\rho}.
 \label{W:eq:tightness-Gram-parity}
\end{equation}
The regularized multiplication matrix satisfies $J_{{\rho},\kappa}(-H)=-D_{\rho}J_{{\rho},\kappa}(H)D_{\rho}$. In particular, its diagonal changes sign, its real nearest-neighbour
entries do not, and cutoffs defined by functions of the eigenvalues of the Gram matrix are unchanged under reflection.
\end{lemma}
\begin{proof}
We compare the Lanczos constructions for $H$ and $-H$ with the same starting vector and the same deterministic rule for continuation. Consider first one invariant block beginning
at index $r$, and write $\varepsilon_k=(-1)^{k-r}$. At its first site the starting vector is identical for the two matrices. Suppose inductively that $q_k(-H)=\varepsilon_kq_k(H)$
and, where $k>r$, $q_{k-1}(-H)=\varepsilon_{k-1}q_{k-1}(H)$, with the neighbouring coefficient already known to agree. Then
\[
 a_k(-H)
 =q_k(-H)^*(-H)q_k(-H)=-a_k(H).
\]
The vector $r_k(-H)$ satisfies
\begin{align*}
 r_k(-H)
 &=(-H)q_k(-H)-a_k(-H)q_k(-H)
             -b_{k-1}(-H)q_{k-1}(-H)\\
 &=-\varepsilon_k(Hq_k(H)-a_k(H)q_k(H)
             -b_{k-1}(H)q_{k-1}(H))\\
 &=-\varepsilon_kr_k(H).
\end{align*}
At $k=r$ the last term in the first line is absent because the preceding coefficient vanishes. It follows that
$b_k(-H)=\|r_k(-H)\|=\|r_k(H)\|=b_k(H)$.
If this quantity is positive, normalization gives $q_{k+1}(-H)=-\varepsilon_kq_{k+1}(H) =\varepsilon_{k+1}q_{k+1}(H)$, closing the induction. If it vanishes, it vanishes for both
matrices. The accumulated subspaces are the same, since their bases differ only by signs. The continuation chooses the first nonzero projection of a coordinate basis vector onto
their common orthogonal complement. Both choices and normalizations therefore agree, and the next block starts with sign $+1$. Induction over all blocks proves the identities for
$a_k,b_k$ and the agreement of the indices for which $r_k=0$, including the final diagonal coefficient.

The Chebyshev recurrence gives $p_j(-x)=(-1)^jp_j(x)$. With $D_\rho=\operatorname{diag}(1,-1,1,-1,\ldots)$, we get $\mathcal V_\rho(-H)=\mathcal V_\rho(H)D_\rho$. Hence
\[
 G_\rho(-H)=D_\rho G_\rho(H)D_\rho,
 \qquad M_\rho(-H)=-D_\rho M_\rho(H)D_\rho.
\]
A scalar spectral function commutes with unitary conjugation, so the regularized Gram matrix transforms by the same relation. If $C$ is its positive diagonal upper Cholesky factor,
then $D_\rho CD_\rho$ is upper triangular, has the same positive diagonal and is a Cholesky factor of the conjugate matrix. Uniqueness forces it to be the factor at $-H$.
Substituting the conjugation formulas into the regularized Jacobi compression gives
\[
 J_{\rho,\kappa}(-H)
 =-D_\rho J_{\rho,\kappa}(H)D_\rho.
\]
The diagonal changes sign, whereas the real neighbouring entries are unchanged, because $(D_\rho)_{kk}(D_\rho)_{k+1,k+1}=-1$. Conjugation also preserves every eigenvalue of the
Gram matrix, so cutoffs defined solely through these eigenvalues agree for $H$ and $-H$. These identities remain valid after any of the specified deterministic continuations.
\end{proof}

Set $e_0=0$, $e_k=b_k(H,w)-1$, and $\Delta X_k^\sigma=X_N^\sigma(k{h_N})-X_N^\sigma((k-1){h_N})$. The definitions and Lemma~\ref{W:lem:tightness-reflection} give
${h_N}\Delta X_k^+=-a_k-e_{k-1}-e_k$ and ${h_N}\Delta X_k^-=a_k-e_{k-1}-e_k$. Subtraction and iteration from $e_0=0$ yield, for $1\le k\le K<N$,
\begin{align}
 a_k&=\frac {h_N}2(\Delta X_k^--\Delta X_k^+),
       \label{W:eq:tightness-diagonal-recovery}\\
 e_k&=-\frac {h_N}2\sum_{j=1}^k(-1)^{k-j}
                       (\Delta X_j^++\Delta X_j^-).
       \label{W:eq:tightness-offdiagonal-recovery}
\end{align}
Writing $J_N^\sigma(K)=\max_{j\le K}|\Delta X_j^\sigma|$, we obtain
\begin{equation}
 \max_{k\le K}|b_k-1|\le\frac{K{h_N}}{2}(J_N^+(K)+J_N^-(K)),\qquad
 \max_{k\le K}|a_k|\le\frac {h_N}2(J_N^+(K)+J_N^-(K)).
 \label{W:eq:tightness-recovery-bound}
\end{equation}
The calculation permits zero off-diagonal coefficients and therefore applies to the full tridiagonal matrix. It also applies to the diagonal and real adjacent entries of the
regularized multiplication matrix and its reflected integrated potentials, using their parity and the same convention $e_0=0$.

\Needspace{6\baselineskip}

\begin{theorem}
\label{W:thm:tightness-offdiagonal-coefficients}
Suppose both families $X_N^+$ and $X_N^-$ are tight in $C([0,T])$ for every finite $T$. Then, for every finite $L$,
\begin{equation}
 B_N(L)\xrightarrow{\mathbb P}0,\qquad
 m_N\max_{k\le\lceil m_NL\rceil+1}|a_{k,N}|
          \xrightarrow{\mathbb P}0.
 \label{W:eq:tightness-offdiagonal-coefficients-conclusion}
\end{equation}
The convergence is uniform over deterministic starting vectors when the assumed tightness is uniform.
\end{theorem}
\begin{proof}
Fix $L<\infty$ and put $K_N=\lceil m_NL\rceil+1$. For all large $N$, the points $kh_N$ with $k\le K_N$ lie in $[0,L+1]$. Since each $X_N^\sigma$ is linearly interpolated,
\[
 J_N^\sigma(K_N)
 :=\max_{1\le k\le K_N}
     |X_N^\sigma(kh_N)-X_N^\sigma((k-1)h_N)|
 \le\omega_{L+1}(X_N^\sigma,h_N).
\]
For a tight family of laws in $C([0,L+1])$, given $\varepsilon>0$ we may choose a compact set $\mathcal K$ carrying probability at least $1-\varepsilon$ for all sufficiently large
$N$. Arzel\`a--Ascoli gives
$\sup_{x\in\mathcal K}\omega_{L+1}(x,\delta)\to0$
as $\delta\downarrow0$. Since $h_N\to0$, it follows that $J_N^\sigma(K_N)\to0$ in probability for both signs.

The exact recurrences from \eqref{W:eq:tightness-diagonal-recovery} and \eqref{W:eq:tightness-offdiagonal-recovery} imply, with $e_0=0$ and $e_k=b_k-1$,
\[
 \max_{k\le K_N}|e_k|
 \le\frac{K_Nh_N}{2}\bigl(J_N^+(K_N)+J_N^-(K_N)\bigr),
 \qquad
 m_N\max_{k\le K_N}|a_k|
 \le\frac12\bigl(J_N^+(K_N)+J_N^-(K_N)\bigr).
\]
As $K_Nh_N\le L+2h_N$, these two quantities converge to zero. For $N$ sufficiently large that $K_Nh_N\le L+1$, the union bound gives the more explicit estimate
\begin{equation}
 \mathbb P(B_N(L)>\eta)
 \le\sum_{\sigma=\pm}
   \mathbb P\left(\omega_{L+1}(X_N^\sigma,h_N)>
                             \frac\eta{L+1}\right),
 \label{W:eq:tightness-offdiagonal-coefficients-probability}
\end{equation}
where $B_N(L)=\max_{k\le K_N}|b_k-1|$. Indeed, if both moduli are at most $\eta/(L+1)$, the preceding estimate is at most $\eta$.

On the event $B_N(L)<1/2$, each $b_k$ in the range is strictly positive, so there is no Lanczos breakdown there. Also
\[
 |b_k^{-1}-1|=\frac{|b_k-1|}{b_k}\le2|b_k-1|,
 \quad k\le K_N.
\]
Consequently
\[
 \mathbb P(\exists k\le K_N:b_k=0)
 \le\mathbb P(B_N(L)\ge1)\longrightarrow0
\]
and, for every $\varepsilon>0$,
\[
 \mathbb P\left(B_N(L)<\tfrac12,
          \max_{k\le K_N}|b_k^{-1}-1|>\varepsilon\right)
 \le\mathbb P(B_N(L)>\varepsilon/2)\longrightarrow0.
\]
If reciprocal coefficients are defined arbitrarily on the exceptional event, the probability of that event also tends to zero and does not change the claimed convergence. On each
interval $[(k-1)h_N,kh_N]$, a straight-line interpolation lies between its two vertex values. Therefore
\[
 \|X_N^{\sigma,\mathrm{lin}}-
           X_N^{\sigma,\mathrm{st}}\|_{L^\infty(0,L)}
 \le J_N^\sigma(K_N)\xrightarrow{\mathbb P}0.
\]

For the uniform statement fix $\eta>0$ and take a supremum over all deterministic unit vectors $w$ in the preceding union bound. Uniform tightness makes
\[
 \lim_{\delta\downarrow0}\limsup_N\sup_{\|w\|=1}
 \mathbb P\bigl(\omega_{L+1}(X_N^\sigma(\cdot;w),\delta)>u\bigr)=0
 \quad(u>0).
\]
The same bound holds with $\delta=h_N$, so all the preceding probability estimates tend to zero after the supremum over $w$. Finally the class of admissible ensembles is unchanged
by $H\mapsto-H$: centeredness, absolute tail conditions and uniform integrability of standardized squares are invariant under sign reversal, and the complex pseudovariance remains
zero. Thus the estimates may be used for both signs under the same hypotheses. This completes the proof of the theorem.
\end{proof}

\Needspace{6\baselineskip}

\begin{lemma}
\label{W:lem:tightness-dyadic}
Let $Z$ be an almost surely continuous real process on $[0,T]$ and suppose for some $p>0$, $\alpha>0$ and $A<\infty$ that
\[
 \mathbb E|Z(t)-Z(s)|^p\le A|t-s|^{1+\alpha}.
\]
Then for $0<\delta\le T$ and $u>0$,
\begin{equation}
 \mathbb P(\omega_T(Z,\delta)>u)
 \le C_{p,\alpha} A T\,u^{-p}\delta^\alpha.
 \label{W:eq:tightness-dyadic-modulus}
\end{equation}
A family with this bound uniformly and with tight initial values is tight in $C([0,T])$.
\end{lemma}
\begin{proof}
We first estimate the oscillation on an interval $I=[a,a+\ell]\subset[0,T]$. For $n\ge0$, let
\[
 M_n=\max_{0\le j<2^n}
 \left|Z\bigl(a+(j+1)\ell2^{-n}\bigr)-Z\bigl(a+j\ell2^{-n}\bigr)\right|.
\]
The increment hypothesis and the inequality $\max_j x_j^p\le\sum_jx_j^p$ for nonnegative $x_j$ give
\[
 \mathbb E M_n^p\le 2^n A(\ell2^{-n})^{1+\alpha}
 =A\ell^{1+\alpha}2^{-n\alpha}.
\]
For $t\in I$, approximate $t$ from the left by the dyadic partition points
$t_n=a+\ell2^{-n}\lfloor 2^n(t-a)/\ell\rfloor$.
Consecutive approximations are either equal or adjacent points of the finer partition. For $t<a+\ell$, $t_0=a$; for $t=a+\ell$, its difference from $a$ is at most $M_0$. Continuity
and telescoping consequently give
\[
 |Z(t)-Z(a)|\le\sum_{n\ge0}M_n,\qquad
 \operatorname{osc}_I Z\le2\sum_{n\ge0}M_n.
\]
For $0<p\le1$, subadditivity of $x^p$ and monotone convergence show
\[
 \mathbb E\left(\sum_{n\ge0}M_n\right)^p
 \le\sum_{n\ge0}\mathbb E M_n^p
 \le\frac{A\ell^{1+\alpha}}{1-2^{-\alpha}}.
\]
For $p>1$, apply Minkowski's inequality first to the finite partial sums and then use monotone convergence:
\[
 \left[\mathbb E\left(\sum_{n\ge0}M_n\right)^p\right]^{1/p}
 \le\sum_{n\ge0}(\mathbb E M_n^p)^{1/p}
 \le\frac{A^{1/p}\ell^{(1+\alpha)/p}}{1-2^{-\alpha/p}}.
\]
These estimates in particular show that the sum is finite almost surely. Thus, for every $p>0$ in the statement,
\[
 \mathbb E(\operatorname{osc}_I Z)^p\le C_{p,\alpha}A\ell^{1+\alpha}.
\]
This proves the oscillation estimate; compare the more general continuity bounds in \cite{kraetschmer2023kolmogorov}.

Every pair $0\le s\le t\le T$ with $t-s\le\delta$ belongs to one of the intervals $I_j=[j\delta,(j+2)\delta]\cap[0,T]$, $0\le j<\lceil T/\delta\rceil$. There are at most
$2T/\delta$ such intervals because $\delta\le T$, and their lengths are at most $2\delta$. Markov's inequality, followed by summation, gives
\[
 \begin{aligned}
 \mathbb P(\omega_T(Z,\delta)>u)
 &\le\sum_j\mathbb P(\operatorname{osc}_{I_j}Z>u)\\
 &\le u^{-p}\sum_j\mathbb E(\operatorname{osc}_{I_j}Z)^p
 \le C_{p,\alpha}A\frac{2T}{\delta}\,u^{-p}(2\delta)^{1+\alpha}.
 \end{aligned}
\]
Absorbing the powers of two proves \eqref{W:eq:tightness-dyadic-modulus}.

For a family with the same bound, fix $\varepsilon>0$. Tightness of the initial values supplies $R<\infty$ with $\sup_Z\mathbb P(|Z(0)|>R)<\varepsilon/2$. Choose
$\delta_j\downarrow0$ so that
\[
 C_{p,\alpha}AT\,2^{jp}\delta_j^\alpha\le\varepsilon2^{-j-1},\qquad j\ge1,
\]
and define
\[
 \mathcal K=\{x\in C([0,T]): |x(0)|\le R,\ 
                    \omega_T(x,\delta_j)\le2^{-j}\text{ for every }j\ge1\}.
\]
The modulus conditions give uniform equicontinuity. The condition with $j=1$, applied to a partition into intervals of length at most $\delta_1$, gives
$\|x\|_\infty\le R+\lceil T/\delta_1\rceil/2$ on $\mathcal K$. Each defining condition is closed for uniform convergence, since
$|\omega_T(x,\delta)-\omega_T(y,\delta)|\le2\|x-y\|_\infty$.
Arzel\`a--Ascoli makes $\mathcal K$ compact. Finally,
\[
 \sup_Z\mathbb P(Z\notin\mathcal K)
 \le\sup_Z\mathbb P(|Z(0)|>R)+
       \sum_{j\ge1}\sup_Z\mathbb P(\omega_T(Z,\delta_j)>2^{-j})
 <\varepsilon.
\]
This proves tightness in the supremum-norm topology.
\end{proof}

The compact-set argument in the proof of Theorem~\ref{W:thm:tightness-offdiagonal-coefficients}, applied on $[0,T+1]$, gives $\omega_{T+1}(\widetilde Y_N^{\rm lin},h_N)\to0$ in
probability. Each successive increment and the difference between the step and linear interpolations on the corresponding interval are bounded by this modulus. Consequently
\begin{equation}
 \max_{1\le k\le\lfloor m_NT\rfloor+1}
       |\widetilde Y_{N,k}-\widetilde Y_{N,k-1}|
       \xrightarrow{\mathbb P}0,\qquad
 \|\widetilde Y_N^{\mathrm{lin}}-\widetilde Y_N^{\mathrm{st}}\|_{C([0,T])}
       \xrightarrow{\mathbb P}0.
 \label{W:eq:coefficients-mesh-consequences}
\end{equation}
The second estimate uses the ordinary supremum norm for bounded functions on $[0,T]$. It also bounds the Skorokhod $J_1$ distance by choosing the identity time change. Thus, if the
continuous interpolations converge in distribution to a continuous process, the step paths have the same $J_1$ limit by the metric-space form of Slutsky's theorem. Uniform
tightness over deterministic unit starting vectors gives the same conclusion with the corresponding supremum of failure probabilities, because the compact set in the preceding
argument can be chosen uniformly by hypothesis.

Convergence of the joint values at finitely many fixed times would not suffice. For example, on $[0,1]$ let $f_n(t)=(1-n^2|t-1/n|)_+$. At each fixed time $t$, the functions
converge to zero, and the same is true at every fixed finite collection of times. But $f_n(1/n)=1$, whereas $f_n(1/n+n^{-2})=0$; hence $\|f_n\|_\infty=1$ and
$\omega_1(f_n,n^{-2})\ge1$ for $n\ge2$. Any uniform limit would be the pointwise zero function, in contradiction with the displayed supremum norm. The following smooth maximum
detects such failures of uniform continuity.

For two real numbers, let
$L_\lambda(y)=\lambda^{-1}\log(e^{\lambda y_1}+e^{\lambda y_2})$
and $p=e^{\lambda y_1}/(e^{\lambda y_1}+e^{\lambda y_2})$. Direct differentiation gives

\begin{equation}
 D_{y}L_\lambda(y)[v]=pv_1+(1-p)v_2, \qquad D_{y}^2L_\lambda(y)[v,w] =\lambda p(1-p)(v_1-v_2)(w_1-w_2).
 \label{W:eq:smooth-maximum-of-two-numbers}
\end{equation}
The first derivative is a weighted average and the second is a covariance. The same finite-sum calculation will control derivatives of a smooth approximation to the modulus of
continuity.

\Needspace{6\baselineskip}

\subsubsection{Smooth approximation of the modulus of continuity}
\label{W:sec:smooth-tests-for-continuity}

The supremum of short increments is the quantity needed for tightness, but it is not a smooth function of the coefficients. Linear interpolation bounds this supremum by the maximum
of finitely many differences of sampled values after the permitted separation between the sampling times is increased by two cell lengths. The elementary bound below makes that
reduction precise. It is followed by a smooth approximation whose derivative bounds depend only logarithmically on the number of pairs.

\begin{lemma}
\label{W:lem:tightness-mesh-geometry}
Let $K=\lceil T/{h_N}\rceil$ and let $y$ be linearly interpolated from $(y_0,\ldots,y_K)$ on $[0,K{h_N}]$. For $d>0$ put
\[
 \mathcal P_d=\{(i,j):0\le i<j\le K,\ (j-i){h_N}\le d+2{h_N}\},\qquad
 M_d(y)=\max_{(i,j)\in\mathcal P_d}|y_j-y_i|.
\]
Then
\begin{equation}
 \omega_T(y,d)\le M_d(y)\le\omega_{K{h_N}}(y,d+2{h_N}).
 \label{W:eq:tightness-mesh-modulus}
\end{equation}
There are at most $(K+1)^2$ pairs.
\end{lemma}
\begin{proof}
Choose cells containing $s$ and $t$, and write
\[
 y(s)=\alpha_0y_i+\alpha_1y_{i+1},\qquad
 y(t)=\beta_0y_j+\beta_1y_{j+1},
 \quad\alpha_a,\beta_b\ge0,\quad\alpha_0+\alpha_1=\beta_0+\beta_1=1.
\]
Then
\[
 y(t)-y(s)=\sum_{a,b\in\{0,1\}}\alpha_a\beta_b
                                      (y_{j+b}-y_{i+a}).
\]
For every pair in this sum, $|(j+b-i-a)h_N|\le|t-s|+2h_N\le d+2h_N$. If its two indices coincide the difference is zero; otherwise, after ordering them, it is bounded in modulus by
$M_d(y)$. As the four weights sum to one, $|y(t)-y(s)|\le M_d(y)$. Taking the supremum proves the first inequality. Every pair in $\mathcal P_d$ consists of two points in
$[0,Kh_N]$ at distance at most $d+2h_N$, giving the second. Finally $|\mathcal P_d|\le\binom{K+1}{2}\le(K+1)^2$.
\end{proof}

\Needspace{6\baselineskip}

Fix a smooth nondecreasing $\vartheta:\mathbb R\to[0,1]$ with $\vartheta=0$ on $(-\infty,1/2]$ and $\vartheta=1$ on $[1,\infty)$. Let $\mathcal P$ be a nonempty set of $D$ pairs of
coordinates, set $L_D=\log(1+2D)$, $\lambda=4L_D/u$, and define
\begin{equation}
 S_{\lambda,\mathcal P}(y) =\lambda^{-1}\log\left(1+ \sum_{(i,j)\in\mathcal P} [e^{\lambda(y_j-y_i)}+e^{-\lambda(y_j-y_i)}]\right), \qquad \Phi_{\mathcal P,u}(y)= \vartheta(S_{\lambda,\mathcal P}(y)/u).
 \label{W:eq:smooth-modulus-definition}
\end{equation}

\begin{lemma}
\label{W:lem:tightness-smooth-max}
Writing $M=\max_{\mathcal P}|y_j-y_i|$, one has
\begin{equation}
 \mathbf1_{\{M\ge u\}}\le\Phi_{\mathcal P,u}(y)
                  \le\mathbf1_{\{M>u/4\}}.
 \label{W:eq:smooth-modulus-indicators}
\end{equation}
For each fixed $r\ge1$ and all real direction vectors,
\begin{equation}
 |D_{y}^r\Phi_{\mathcal P,u}(y)[v_1,\ldots,v_r]|
 \le C_{r,\vartheta}u^{-r}(1+L_D)^{r-1}
                    \prod_{j=1}^r\|v_j\|_\infty.
 \label{W:eq:smooth-modulus-derivatives}
\end{equation}
The constants are uniform in the number and location of the coordinates, in $y$, and in the choice of the pair set.
\end{lemma}
\begin{proof}
The logarithmic maximum is between $M$ and $M+L_D/\lambda=M+u/4$. The cutoff thresholds give \eqref{W:eq:smooth-modulus-indicators}. For the derivatives, let $\mathcal L$ consist
of the zero functional and the $2D$ signed coordinate differences. For each fixed $y$, give $\ell\in\mathcal L$ probability
$e^{\lambda\ell(y)}/\sum_{\ell'\in\mathcal L}e^{\lambda\ell'(y)}$, and denote the resulting finite average by $\mathbb E_y$. Differentiating the logarithm of the sum gives
\[
 D_{y}^r S(y)[v_1,\ldots,v_r]
  =\lambda^{r-1}\sum_{\pi\in\Pi(r)}
      (-1)^{|\pi|-1}(|\pi|-1)!
      \prod_{B\in\pi}\mathbb E_y\prod_{j\in B}\ell(v_j).
\]
Here $\Pi(r)$ denotes the partitions of $\{1,\ldots,r\}$. For $Z(t)=\sum_\ell e^{\lambda\ell(y)+\lambda\sum_at_a\ell(v_a)}$, the labelled chain rule for $\lambda^{-1}\log Z$ uses
$(\log)^{(k)}(Z)=(-1)^{k-1}(k-1)!Z^{-k}$ and $\partial_BZ(0)/Z(0)=\lambda^{|B|}\mathbb E_y\prod_{a\in B}\ell(v_a)$, which proves the formula. Each $|\ell(v_j)|\le2\|v_j\|_\infty$.
The finite partition sum is bounded by a constant depending only on $r$, so $|D_{y}^rS|\le C_r\lambda^{r-1}\prod\|v_j\|_\infty$. For the scalar composition the labelled formula is
\[
 D^r\Phi(y)[v_1,\ldots,v_r]
 =\sum_{\pi\in\Pi(r)}u^{-|\pi|}\vartheta^{(|\pi|)}(S(y)/u)
             \prod_{B\in\pi}D^{|B|}S(y)[v_j:j\in B].
\]
For a partition into $k$ blocks the product of the powers of $\lambda$ is $\lambda^{\sum_{B}(|B|-1)}=\lambda^{r-k}$. Thus
\[
 |D^r\Phi(y)[v_1,\ldots,v_r]|
 \le C_{r,\vartheta}u^{-r}
       \sum_{k=1}^r(4L_D)^{r-k}\prod_{j=1}^r\|v_j\|_\infty
 \le C_{r,\vartheta}u^{-r}(1+L_D)^{r-1}
                                      \prod_{j=1}^r\|v_j\|_\infty.
\]
This proves \eqref{W:eq:smooth-modulus-derivatives}.
\end{proof}

\Needspace{6\baselineskip}

Let $\Gamma_N\in[0,1]$ be smooth and supported in a region on which the entire regularized vector of sampled values agrees with the vector obtained from the Lanczos matrix, through
index $K$. It may equal one on a smaller inner region where the eigenvalues of the Gram matrix satisfy stricter bounds. Write the cutoff expectation at $H_N^{\mathrm G}$ and the
difference of expectations as
\[
 g_N=\mathbb E(1-\Gamma_N({H_N^{\mathrm G}})),\qquad
 e_N^\Gamma=|\mathbb E\Gamma_N(H_N)-\mathbb E\Gamma_N({H_N^{\mathrm G}})|.
\]

\begin{lemma}
\label{W:lem:tightness-supported-cutoff}
Then $\mathbb E(1-\Gamma_N(H_N))\le g_N+e_N^\Gamma$. For every measurable function $f$ of the sampled values with $0\le f\le1$,
\begin{equation}
 \mathbb Ef(\mathbf X_N(H_N))
 \le\mathbb E[\Gamma_N f(\mathbf X_N^{\mathrm{reg}})](H_N)
                  +g_N+e_N^\Gamma.
 \label{W:eq:tightness-supported-removal}
\end{equation}
\end{lemma}
\begin{proof}
Write $\Gamma=\Gamma_N(H_N)$, and similarly
$\Gamma^{\mathrm G}=\Gamma_N(H_N^{\mathrm G})$.
The definitions give the exact equality
\[
 \mathbb E(1-\Gamma)
 =\mathbb E(1-\Gamma^{\mathrm G})
    +\mathbb E\Gamma^{\mathrm G}-\mathbb E\Gamma
 \le g_N+e_N^\Gamma.
\]
Let $\mathbf X_N$ be the original vector of sampled values and $\mathbf X_N^{\rm reg}$ the smooth regularization. On the event $\{\Gamma>0\}$ these two vectors coincide by the
hypothesis on the cutoff. For such a function $f$,
\begin{align*}
 \mathbb Ef(\mathbf X_N)
 &=\mathbb E[\Gamma f(\mathbf X_N)]
      +\mathbb E[(1-\Gamma)f(\mathbf X_N)]\\
 &=\mathbb E[\Gamma f(\mathbf X_N^{\rm reg})]
      +\mathbb E[(1-\Gamma)f(\mathbf X_N)]\\
 &\le\mathbb E[\Gamma f(\mathbf X_N^{\rm reg})]
      +\mathbb E(1-\Gamma)\\
 &\le\mathbb E[\Gamma f(\mathbf X_N^{\rm reg})]
      +g_N+e_N^\Gamma.
\end{align*}
The second equality is exact because the multiplier $\Gamma$ vanishes wherever equality of the original and regularized vectors is not assumed. This proves
\eqref{W:eq:tightness-supported-removal}.

When the polynomial vectors themselves are unchanged and only the formulas for their Gram matrix are regularized, the following cutoff has the required support. Choose fixed
numbers $0<a_0<a_1<a_2<a_3$ inside an interval on which the scalar modifications used for the positive Gram matrix are the identity and reciprocal on a neighbourhood of
$[a_0,a_3]$. Let $b\in C^\infty(\mathbb R;[0,1])$ satisfy $b=0$ on $(-\infty,0]$ and $b=1$ on $[1,\infty)$, and define
\[
 \gamma(x)=b\!\left(\frac{x-a_0}{a_1-a_0}\right)
            b\!\left(\frac{a_3-x}{a_3-a_2}\right),
 \qquad \Gamma_N=\det\gamma(G_\rho).
\]
The spectral theorem gives
$\Gamma_N=\prod_{j=1}^\rho\gamma(\lambda_j(G_\rho))$,
so $0\le\Gamma_N\le1$. If all eigenvalues lie in $[a_1,a_2]$, then every factor is one and $\Gamma_N=1$. If $\Gamma_N>0$, every factor is positive and every eigenvalue belongs to
$(a_0,a_3)$. On this set the regularized inverse and Cholesky formulas agree with their original versions for the same polynomial vectors. The resulting values agree at every index
appearing in $f$.

If the polynomial vectors are also changed from $p_j(H)w$ to $p_j(h_d(H))w$, this conclusion additionally uses $\Omega_d=\{h_d(H)=H\}$. The pointwise decomposition in that case is
\[
 f(\mathbf X_N)\le1-\Gamma_N+
             \Gamma_Nf(\mathbf X_N^{\rm reg})+\mathbf1_{\Omega_d^c},
 \qquad 0\le f\le1.
\]
On $\Omega_d\cap\{\Gamma_N>0\}$ the two vectors agree, and on its complement either $1-\Gamma_N=1$ or the last term is one. Taking expectations therefore adds exactly
$\mathbb P(\Omega_d^c)$ to the preceding estimate. Section~\ref{W:sec:proof-local-coefficients} uses this version when returning to the original matrix. The determinant is smooth
as a function of the entries of the Gram matrix. Because $G_\rho(-H)=D_\rho G_\rho(H)D_\rho$, the same spectral representation also proves $\Gamma_N(-H)=\Gamma_N(H)$.

\end{proof}

\Needspace{6\baselineskip}

Let $\mathbf X_N^{\sigma,\mathrm{reg}}$ denote the regularized vector of sampled values and put $F_N^\sigma=\Phi_{\mathcal P_\delta,u}(\mathbf X_N^{\sigma,\mathrm{reg}})$. Define
\begin{equation}
 e_N^\sigma(T,\delta,u)=
 \left|\mathbb E[\Gamma_NF_N^\sigma](H_N)
       -\mathbb E[\Gamma_NF_N^\sigma](H_N^{\mathrm G})\right|.
 \label{W:eq:tightness-oscillation-error}
\end{equation}
The quantities $g_N$ and $e_N^\Gamma$ are those of the preceding lemma. We use the Gaussian fourth moment modulus estimate \eqref{W:eq:tightness-gaussian-modulus}.

\begin{theorem}
\label{W:thm:modulus-comparison}
Suppose, as in Lemma~\ref{W:lem:tightness-supported-cutoff}, that the original and regularized vectors agree wherever $\Gamma_N>0$. For fixed $T,u>0$, $0<\delta\le T$ and
sufficiently large $N$,
\begin{equation}
 \mathbb P(\omega_T(X_N^\sigma,\delta)>u)
 \le C_{\beta,T}u^{-4}(\delta+2h_N)
       +e_N^\sigma(T,\delta,u)+g_N+e_N^\Gamma.
 \label{W:eq:tightness-modulus-transfer}
\end{equation}
If equality of the original and modified vectors additionally requires $\Omega_d$, the right side acquires
$\mathbb P(H_N\notin\Omega_d)+\mathbb P(H_N^{\mathrm G}\notin\Omega_d)$.
In the displayed setting, both families are tight if $g_N+e_N^\Gamma\to0$ and
\begin{equation}
 \lim_{\delta\downarrow0}\limsup_N
       \max_{\sigma=\pm}e_N^\sigma(T,\delta,u)=0
       \quad\text{for every }T,u>0.
 \label{W:eq:tightness-oscillation-input}
\end{equation}
\end{theorem}
\begin{proof}
Lemma~\ref{W:lem:tightness-mesh-geometry} and the lower bound for $\Phi_{\mathcal P_\delta,u}$ give
\[
 \mathbf1_{\{\omega_T(X_N^\sigma,\delta)>u\}}
 \le\Phi_{\mathcal P_\delta,u}(\mathbf X_N^\sigma).
\]
Multiplying by the cutoff and then comparing expectations gives
\[
 \begin{aligned}
 \mathbb P(\omega_T(X_N^\sigma,\delta)>u)
 &\le\mathbb E\Phi_{\mathcal P_\delta,u}(\mathbf X_N^\sigma)\\
 &\le\mathbb E[\Gamma_NF_N^\sigma](H_N)+g_N+e_N^\Gamma\\
 &\le\mathbb E[\Gamma_NF_N^\sigma](H_N^{\mathrm G})
                         +e_N^\sigma+g_N+e_N^\Gamma\\
 &\le\mathbb E\Phi_{\mathcal P_\delta,u}(\mathbf X_{N,G}^\sigma)
                         +e_N^\sigma+g_N+e_N^\Gamma.
 \end{aligned}
\]
The last step uses equality of the vectors of sampled values on the nonzero cutoff set and $0\le\Gamma_N\le1$. The upper indicator bound and the inequality for the modulus of
continuity on the partition therefore give
\[
 \mathbb P(\omega_T(X_N^\sigma,\delta)>u)
 \le\mathbb P\left(
   \omega_{K{h_N}}(X_{N,G}^\sigma,\delta+2{h_N})>u/4\right)
        +e_N^\sigma+g_N+e_N^\Gamma.
\]
Use \eqref{W:eq:tightness-gaussian-modulus} on $[0,T+1]$, which contains $[0,K{h_N}]$ for large $N$, and enlarge the constant. This proves \eqref{W:eq:tightness-modulus-transfer}.
If equality of the vectors is available only on $\Omega_d$, Lemma~\ref{W:lem:tightness-supported-cutoff} adds $\mathbb P(H_N\notin\Omega_d)$ in the second line. At the Gaussian
matrix, the pointwise inequality
\[
 \Gamma_NF_N^\sigma\le
 \Phi_{\mathcal P_\delta,u}(\mathbf X_{N,G}^\sigma)+\mathbf1_{\{H_N^{\mathrm G}\notin\Omega_d\}}
\]
adds the second probability. This proves the stated version with spectral modification; the same tightness conclusion holds when these additional probabilities tend to zero. For
each fixed $T,u,\delta$, $h_N\to0$ and $g_N+e_N^\Gamma\to0$ give
\[
 \limsup_N\max_{\sigma=\pm}
  \mathbb P(\omega_T(X_N^\sigma,\delta)>u)
 \le C_{\beta,T}u^{-4}\delta
       +\limsup_N\max_{\sigma=\pm}e_N^\sigma(T,\delta,u).
\]
Letting $\delta\downarrow0$ and using \eqref{W:eq:tightness-oscillation-input} gives
\[
 \lim_{\delta\downarrow0}\limsup_N\max_{\sigma=\pm}
       \mathbb P(\omega_T(X_N^\sigma,\delta)>u)=0.
\]
Here the modulus bound is an iterated limit, rather than the uniform power estimate used in Lemma~\ref{W:lem:tightness-dyadic}. We make the compactness step for this situation
explicit. For each $j\ge1$ and fixed $\varepsilon>0$, the last limit gives a number $\delta_j>0$ and an integer $N_j$ such that
\[
 \max_{\sigma=\pm}\mathbb P(\omega_T(X_N^\sigma,\delta_j)>2^{-j})
 \le\varepsilon2^{-j-1}\quad(N\ge N_j).
\]
The numbers $\delta_j$ can be decreased to make them tend to zero. For the finitely many sizes $N<N_j$, continuity of each path implies $\omega_T(X_N^\sigma,\delta)\to0$ almost
surely as $\delta\downarrow0$. Decreasing $\delta_j$ further therefore makes the displayed bound hold for those sizes too, while its validity for $N\ge N_j$ is preserved by
monotonicity in $\delta$. The set of continuous functions $x$ with $x(0)=0$ and $\omega_T(x,\delta_j)\le2^{-j}$ for every $j$ is closed, uniformly bounded and equicontinuous, as
verified in the preceding lemma. It is compact, and the sum of the failure probabilities is at most $\varepsilon/2$. This proves the stated tightness.

For its uniform version, take the supremum over deterministic unit $w$ in \eqref{W:eq:tightness-modulus-transfer}. The right side becomes
\[
 C_{\beta,T}u^{-4}(\delta+2h_N)
 +\sup_{\|w\|=1}e_N^\sigma(T,\delta,u)
 +\sup_{\|w\|=1}(g_N+e_N^\Gamma).
\]
The uniform assumptions make its iterated limit zero and give $\delta_j,N_j$ with the supremum over $w$ in the preceding probability bound. To cover the finitely many smaller sizes
uniformly as well, each Lanczos coefficient satisfies $|a_k|,b_k\le\|H_N\|$. With $b_0=1$, the interpolated potentials consequently have the deterministic Lipschitz bound
\[
 |X_N^\sigma(t;w)-X_N^\sigma(s;w)|
 \le m_N^2(3+3\|H_N\|)|t-s|,
 \qquad \sigma=\pm,
\]
for every deterministic unit $w$. For a fixed $N$, the right-hand Lipschitz constant is finite almost surely; hence
\[
 \sup_{\|w\|=1}\mathbb P(\omega_T(X_N^\sigma(\cdot;w),\delta)>2^{-j})
 \le\mathbb P\bigl(m_N^2(3+3\|H_N\|)\delta>2^{-j}\bigr)\longrightarrow0.
\]
The same choice of $\delta_j$ now covers all smaller admissible sizes and all deterministic unit $w$. Summing these uniform probability bounds gives the same compact set with
arbitrarily small complementary probability, in the sense of \eqref{W:eq:uniformity-in-the-starting-vector}.
\end{proof}

\subsubsection{Notation for comparison of entry distributions}
\label{W:sec:notation-and-scales}

The remainder of the comparison uses both normalized traces and matrix elements at the starting vector. Their different normalizations affect the powers of $N$ in the estimates, so
it is useful to specify them together with the resolvent convention. For an $N\times N$ matrix $A$, our notation is

\[
 \langle A\rangle=N^{-1}\operatorname{Tr}A,\qquad
 A^\circ=A-\langle A\rangle I,\qquad
 \|A\|_{2,N}=N^{-1/2}\|A\|_{\rm HS}.
\]
The operator, trace, and Hilbert--Schmidt norm conventions remain unchanged. The matrix units are $E_{ab}=e_ae_b^*$. The normalization of the trace and Hilbert--Schmidt norm agrees
with Section~2 of \cite{cipolloni2023rankuniform}. We specify coordinate tensor norms when defining the corresponding observables.

For expansions in the entries of $H$ we use
\[
 R_H(z)=(zI-H)^{-1},\qquad
 g(z)=\int\frac{\dd\mu_{\rm sc}(x)}{z-x}.
\]
For the scaled operators the convention is $(\mathcal A_N-z)^{-1}$; thus the alternative semicircle transform is $m_{\rm sc}=-g$. For $z=E+\ii\eta$ with $\eta>0$, set
\[
 \upsilon(z)=|\operatorname{Im}g(z)|,\quad
 \varrho(z)=\upsilon(z)/\pi,\quad
 d(z)=\operatorname{dist}(z,[-2,2]),\quad
 \Psi(z)=\sqrt{\frac{\varrho(z)}{N\eta}}+\frac1{N\eta}.
\]
We denote the density of $\mu_{\rm sc}$ by $\rho_{\rm sc}$ and its increasing quantiles by $\int_{-2}^{\gamma_j}\rho_{\rm sc}(x)\,\dd x=j/N$. Also write $\widehat j=\min(j,N+1-j)$.

For $I\subset\{1,\ldots,N\}$, $H^{(I)}$ denotes the principal minor obtained by deleting $I$, and $R_H^{(I)}(z)=(zI-H^{(I)})^{-1}$. By contrast,
$R_H(z)_{II}=\iota_I^*R_H(z)\iota_I$ is the compression of the full resolvent to the coordinate subspace indexed by $I$, where $\iota_I$ is the coordinate inclusion. The matrix
$K_{ab}$ used in entry comparison is the $N\times N$ matrix obtained by setting $h_{ab}$ and $h_{ba}$ to zero. The identities for the minor follow from Schur complementation
\cite{benaych2018lectures}.

The modulus of the resolvent is $|R_H(z)|=(R_H(z)^*R_H(z))^{1/2}$. With our sign convention and
$\eta=\operatorname{Im}z>0$,
\[
 R_H(z)-R_H(z)^*=-2\mathrm i\eta R_H(z)R_H(z)^*,\qquad
 \sum_j|R_H(z)_{ij}|^2=-\frac{\operatorname{Im}R_H(z)_{ii}}{\eta}.
\]
The second identity follows by taking the $i$th diagonal entry of the first; it is the Ward identity \cite{benaych2018lectures}.

We use stochastic domination in the sense of \cite{benaych2018lectures}. Uniform stochastic domination, denoted $X\prec Y$, means that for every $\varepsilon,D>0$ there is $N_0$
such that
\[
 \sup_{u\in\mathcal U_N}
 \mathbb P\bigl(|X_N(u)|>N^\varepsilon Y_N(u)\bigr)\le N^{-D}
 \qquad(N\ge N_0),
\]
where $u$ ranges over the specified deterministic parameters. An event holds with very high probability when, for each $D>0$, there is $N_0(D)$ such that its complement has
probability at most $N^{-D}$ for every $N\ge N_0(D)$. In an upper bound, a factor $N^{o(1)}$ denotes a deterministic sequence $q_N\ge1$ such that $\log q_N/\log N\to0$. It may
depend on the fixed parameters specified in the estimate. Moment bounds are recorded separately whenever expectations are taken.

For $a<b$, one independent entry is a real scalar when $\beta=1$ and a real vector of dimension two when $\beta=2$. In the complex case put
$x=(\sqrt2\operatorname{Re}\xi,\sqrt2\operatorname{Im}\xi)$. The corresponding matrix perturbation is
\begin{equation}
 N^{-1/2}\mathcal T_{ab}x=
 \begin{cases}
 N^{-1/2}x_1(E_{ab}+E_{ba}),&\beta=1,\\
 (2N)^{-1/2}\bigl((x_1+\ii x_2)E_{ab}
                     +(x_1-\ii x_2)E_{ba}\bigr),&\beta=2.
 \end{cases}
 \label{W:eq:standardized-entry-insertion}
\end{equation}
For a diagonal entry the perturbation is $N^{-1/2}xE_{aa}$. In the complex case integration is always with respect to the joint law of the two coordinates. In each conditional
expansion, all entries other than the selected entry are held fixed. The families of matrices, spectral parameters, and vectors used in simultaneous estimates are specified in the
corresponding results.

The lattice spacing and spectral scale are those of Section~\ref{W:sec:convergence-of-tridiagonal-operators}:
\[
 m_N=N^{1/3},\qquad h_N=m_N^{-1},\qquad
 \eta_{\rm e}=m_N^{-2}=h_N^2.
\]
Here $h_N$ is the distance between adjacent lattice points, and $\eta_{\rm e}$ is the spectral scale at the edge. The comparison arguments use the additional parameters listed
below.
\begin{center}
\small
\begin{tabularx}{\linewidth}{@{}p{.25\linewidth}X@{}}
\toprule
Parameter & Convention \\
\midrule
$L,T,E$ & Fixed intervals, time intervals, and spectral cutoffs. \\
$r\le\rho\le\rho_*$ & Dimensions of the polynomial spans and Gram matrices.
For the fixed observables at the edge, these are of order $m_N$. \\
$\kappa,\delta,u$ & Parameters controlling the Gram matrix and the modulus of continuity, fixed before the limit $N\to\infty$. \\
$e_N$, $L_N=e_N\sqrt N$ & Bounds on the matrix entries and their standardized
coordinates, with $e_N\downarrow0$ and $e_N^{-1}=N^{o(1)}$. \\
$d_N$ & Width of the spectral modification, of order
$m_N^{-2}(\log(N+2))^\gamma$, $1<\gamma<2$.
Section~\ref{W:sec:removal-of-spectral-regularization} also uses
$d_N^{(2)}$ with a squared logarithm. \\
$\mathfrak r_N$ & Threshold in the definition of $\chi_N$. \\
$\mathfrak c_N$, $\Lambda_N$ & Bounds used in the counting estimates and in
the estimates for derivatives, respectively; both are subpolynomial in $N$. \\
\bottomrule
\end{tabularx}
\end{center}
The sequences $e_N$ and $\mathfrak r_N$ are chosen from the tail and coefficient estimates for the fixed observables. A finite collection of subpolynomial bounds can be replaced by
their maximum. For a countable determining collection we choose successive thresholds in $N$ after proving the estimates for each finite subcollection. The full localization
estimate is used in the order $N\to\infty$, then $\kappa\downarrow0$ for fixed $L$, then $L\to\infty$; see Theorem~\ref{W:thm:localization-spectral-projection}. For modulus
tightness the order is $N\to\infty$, $\kappa\downarrow0$, and then $\delta\downarrow0$, with $T,u$ fixed; see Theorem~\ref{W:thm:integrated-potential-and-offdiagonal-coefficients}.

The original matrix is $H_N$. Truncation of its entries gives $H_N^{\rm tr}$. The auxiliary matrix $H_N^{\rm nor}$ is obtained by centering the truncated entries and normalizing
their covariance, as in Section~\ref{W:sec:truncation-and-resolvent-estimates}. The matrix $H_N^{\rm b}$ has uniformly bounded standardized entries with the required moments. The
matrix $H_N^{\rm G}$ is distributed according to the GOE for $\beta=1$ and the GUE for $\beta=2$. Gram matrices are denoted by $G_r$ or $G_\rho$, and a complementary principal
minor resolvent may be written $\mathcal R=R_H^{(I)}$. Entries of a minor have the normalization of their parent matrix unless stated otherwise.

For a Hermitian matrix, $H^{[d]}=h_d(H)$ denotes modification of the eigenvalues by the scalar function $h_d$. Thus $(H_N^{\rm tr})^{[d]}$ first truncates the entries and then
modifies the spectrum. The projections $P_r$ and $P_r^{[d]}$ are onto $\mathcal K_r(H,w)$ and $\mathcal K_r(H^{[d]},w)$. Their regularized counterparts, such as
$P_{r,\kappa}^{[d]}$, agree with these projections on the specified event bounding the Gram matrix and its inverse. On that event, \eqref{W:eq:lanczos-gram-projection} identifies
them with the first $r$ Lanczos coordinates when $h_d(H)=H$. The matrix $K_\alpha(t)$ has entry $\alpha$ set to zero. Its remaining entries are independent, and each distribution
is the specified convex combination, with weight $t\in[0,1]$, of the two entry laws being compared. The replacement entry is independent of $K_\alpha(t)$.

The estimates for real matrices with bounded entries are proved in \cite{lee2014necessary}; Proposition~\ref{W:prop:minor-real-input} checks their assumptions.
Theorem~\ref{W:thm:complex-main} proves the corresponding estimates for complex entries. For convex combinations of bounded and Gaussian entry distributions,
Lemmas~\ref{W:lem:cumulant-event} and~\ref{W:lem:localization-interpolation-event} verify the assumptions of the isotropic law \cite{bloemendal2014isotropic} (Theorem~2.12) and the
edge deviation estimate \cite{erdosxu2023small}. Section~\ref{W:sec:entry-truncation-estimates} treats the original truncated law.

The coefficient quantities $Q_*,B$ were defined in \eqref{W:eq:comparison-QB}. Section~\ref{W:sec:coefficient-localization} constructs their smooth cutoff $\chi_N$. We use these
functions to restrict the size of derivatives. The norms $\|\cdot\|_\diamond$ and $\|\cdot\|_{\mathscr D}$ are defined by their formulas in \eqref{W:eq:tightness-norm} and
\eqref{W:eq:comparison-direction-norm}. The nonnegative quantities $U_\alpha,V_\alpha$ are integrals of products of resolvent coordinates, defined in the corresponding estimates.

The vector $\mathbf z_N(H,w)$ consists of polynomial moments of the spectral measure at $w$, normalized traces, and moments containing $\Theta_a(H)$. The functions of $\mathbf z$
are defined on the full real coordinate space. The differentiation variables are summarized below; parameters and other arguments are held fixed.
\begin{center}
\begin{tabular}{@{}ll@{}}
\toprule
Notation & Variable \\
\midrule
$D_{\mathbf z}$ & The full real coordinate vector $\mathbf z$ \\
$D_s$ & $s=(s_0,s_1,\ldots)$ with $s_0=1$ fixed \\
$D_{\mathsf X}$ & $(G,M,\mathsf t)$ in the localization function \\
$D_G$ & The matrix argument $G$ \\
$D_y$ & The vector $y$ of sampled values \\
$D_H$ & The Hermitian matrix $H$, with $w$ fixed \\
\bottomrule
\end{tabular}
\end{center}
The coordinate $\mathsf t_E=\operatorname{Tr}\Theta_E(H)$ is unnormalized.

A matrix derivative $D_H$ is taken in the real vector space of Hermitian matrices. Its argument includes the normalization displayed in that direction. For a finite label set $S$,
$D_{H,S}$ is the mixed derivative in its labelled directions; coincident directions still have separate labels. The same convention applies to moment derivatives. For
$F(H)=A(s(H))$, repeated product and chain rules give
\[
 \begin{aligned}
 D_S(F_1\cdots F_k)
   &=\sum_{S_1\sqcup\cdots\sqcup S_k=S}(D_{S_1}F_1)\cdots(D_{S_k}F_k),\\
 D_{H,S}[A(s(H))]
   &=\sum_{\pi\in\Pi(S)}D_s^{|\pi|}A(s(H))
                         [D_{H,B}s(H):B\in\pi].
 \end{aligned}
\]
The first sum allows empty sets, with $D_\varnothing F_j=F_j$. In the second sum, $\Pi(S)$ consists of partitions into nonempty sets. To prove the second identity by induction, the
next derivative either differentiates $A$, creating a singleton set, or differentiates one argument $D_{H,B}s$, adding its label to that set. Every partition occurs once. The first
identity follows by assigning each derivative label to the factor it differentiates; the order of matrix factors is unchanged. All functions in these formulas are differentiated at
the same base point. For $f_\alpha(x)=F(K_\alpha+N^{-1/2}\mathcal T_\alpha x)$, the chain rule gives
\begin{equation}
 D_x^qf_\alpha(x)[v_1,\ldots,v_q]
 =N^{-q/2}D_H^qF(K_\alpha+N^{-1/2}\mathcal T_\alpha x)
                  [\mathcal T_\alpha v_1,\ldots,\mathcal T_\alpha v_q].
 \label{W:eq:notation-standardized-derivative}
\end{equation}
Indeed, differentiating $K_\alpha+N^{-1/2}\mathcal T_\alpha x$ with respect to $x$ in the vector $v_j$ gives $N^{-1/2}\mathcal T_\alpha v_j$, and its higher derivatives vanish. For
the unnormalized coordinate $u=N^{-1/2}x$, the factor $N^{-q/2}$ is absent. In particular,
\begin{equation}
 \begin{split}
 D_H^2F(H)[V,W]
 ={}&D_s^2A(s(H))[D_Hs(H)[V],D_Hs(H)[W]]\\
    &+D_sA(s(H))[D_H^2s(H)[V,W]].
 \end{split}
 \label{W:eq:notation-moment-matrix-chain}
\end{equation}
This follows by differentiating $D_HF(H)[V]=D_sA(s(H))[D_Hs(H)[V]]$ in direction $W$. The homogeneous term of degree $q$ in the Taylor expansion is
$D_x^qf_\alpha(0)[x,\ldots,x]/q!$.

For example, fix $v$, put $V=\mathcal T_\alpha v$, and write $s^{(j)}=D_H^js(H)[V,\ldots,V]$ at $H=K_\alpha+N^{-1/2}\mathcal T_\alpha x$. At fourth order the full composition is
\[
 \begin{aligned}
 D_x^4f_\alpha(x)[v^4]=N^{-2}\bigl\{
 &D_s^4A[s^{(1)},s^{(1)},s^{(1)},s^{(1)}]+6D_s^3A[s^{(2)},s^{(1)},s^{(1)}]\\
 &+3D_s^2A[s^{(2)},s^{(2)}]+4D_s^2A[s^{(3)},s^{(1)}]+D_sA[s^{(4)}]\bigr\}.
 \end{aligned}
\]
All derivatives of $A$ are evaluated at $s(H)$. The coefficients count partitions of four labelled positions: there are six choices of a pair and two singletons, three partitions
into two pairs, and four choices of a triple and a singleton. Repeated directions do not remove these multiplicities. In this formula $A$ is the whole function, including its
cutoffs. The bounds for $D_s^rA$ must be combined with separate estimates of the moment derivatives $s^{(2)},s^{(3)},s^{(4)}$.

For a complex off-diagonal entry let $V_j=\mathcal T_{ab}e_j$, $j=1,2$. The cubic Taylor term at $K_{ab}$ is
\[
 \frac{N^{-3/2}}6\bigl\{
 x_1^3D_H^3F[V_1^3]+3x_1^2x_2D_H^3F[V_1,V_1,V_2]
 +3x_1x_2^2D_H^3F[V_1,V_2,V_2]+x_2^3D_H^3F[V_2^3]\bigr\}.
\]
Thus taking expectation uses the joint moments $\mathbb E x_1^2x_2$ and $\mathbb E x_1x_2^2$ as well as the two pure third moments. Their vanishing or matching must follow from the
entry law under consideration. A diagonal entry has one coordinate and is represented by the matrix unit $E_{aa}$, without the off-diagonal factor $1/\sqrt2$.

For later entry formulas we write
\begin{equation}
 P_{w,N}=w_Nw_N^*,\qquad P_{w,N}^{\circ}=P_{w,N}-N^{-1}I.
 \label{W:eq:centered-starting-vector-projection}
\end{equation}

For a resolvent derivative, the trace against $P_{w,N}^{\circ}$ can be nonzero even when the undifferentiated resolvent is scalar. For example, take $H=cI$, $a<b$, and
$V=N^{-1/2}(E_{ab}+E_{ba})$. Since $R_H(z)=(z-c)^{-1}I$ and
$D_H^2R_H(z)[V,V]=2(z-c)^{-3}V^2$,
\[
 \operatorname{Tr}\!\left(P_{w,N}^{\circ}D_H^2R_H(z)[V,V]\right)
 =\frac{2}{N(z-c)^3}
       \left(|w_a|^2+|w_b|^2-\frac2N\right).
\]
Thus the second derivative in the same entry requires an estimate of the displayed coordinate weights.

\subsubsection{Spectral regularization and matrix derivatives}
\label{W:sec:spectral-regularization-and-matrix-derivatives}
The polynomial representation gives explicit matrix derivatives as long as the spectrum stays in a controlled interval. Entry comparison also requires the observables to be defined
outside that interval, where these polynomial bounds can be much larger. A smooth modification of the spectral variable provides global bounds without changing the polynomials on
the interval of interest. The estimates begin with the following scalar and matrix identities:

\begin{equation}
 p_j(2\cosh s)=\frac{\sinh((j+1)s)}{\sinh s},
 \qquad p_j(-x)=(-1)^jp_j(x),
 \label{W:eq:gram-hyperbolic}
\end{equation}
and
\begin{equation}
 D_{H}^qp_n(H)[W_1,\ldots,W_q]
 =\sum_{\sigma\in S_q}\sum_{j_0+\cdots+j_q=n-q}
 p_{j_0}(H)W_{\sigma(1)}p_{j_1}(H)\cdots
 W_{\sigma(q)}p_{j_q}(H).
 \label{W:eq:chebyshev-matrix-derivatives}
\end{equation}
For every finite matrix $H$ and sufficiently small $|t|$, the generating series converges in operator norm:
\[
 Q(t,H):=(I-tH+t^2I)^{-1}=\sum_{n\ge0}p_n(H)t^n.
\]
Differentiating $(I-tH+t^2I)Q=I$ yields $D_HQ[W]=tQWQ$. Repeated differentiation gives
\[
 D_H^qQ(t,H)[W_1,\ldots,W_q]
 =t^q\sum_{\sigma\in S_q}QW_{\sigma(1)}Q\cdots W_{\sigma(q)}Q.
\]
Multiplying the absolutely convergent series, the coefficient of $t^n$ on the right is exactly \eqref{W:eq:chebyshev-matrix-derivatives}. This proves the identity, including the
vanishing for $q>n$. The scalar hyperbolic identity follows by substituting $x=e^s+e^{-s}$ into the same recurrence.

\Needspace{6\baselineskip}

Let $H=H^*$ satisfy $\|H\|\le2+d$, and put $s=\operatorname{arcosh}(1+d/2)$. We use $F_\rho=\mathcal V_\rho(H)$, $G_\rho=F_\rho^*F_\rho$, and $P_\rho=F_\rho G_\rho^{-1}F_\rho^*$
when $G_\rho>0$. In estimates for a matrix entry, each $W_j=W_j^*$ is supported on that entry and its conjugate and satisfies
$\|W_j\|\le C N^{-1/2}$.

\begin{proposition}
\label{W:prop:gram-window}
For $0\le q\le n$,
\begin{equation}
 \|D_H^qp_n(H)[W_1,\ldots,W_q]\|
 \le q!\binom{n+q+1}{2q+1}e^{(n-q)s}\prod_j\|W_j\|.
 \label{W:eq:gram-window-poly}
\end{equation}
For normalized entry directions,
\begin{align}
 \|D_H^qF_\rho[W_1,\ldots,W_q]\|_{\mathrm{HS}}
 &\le C_qN^{-q/2}\rho^{2q+3/2}e^{\rho s},\label{W:eq:gram-window-polynomial-matrix}\\
 \|D_H^qG_\rho[W_1,\ldots,W_q]\|
 &\le C_qN^{-q/2}\rho^{2q+3}e^{2\rho s}.\label{W:eq:gram-window-Gram}
\end{align}
If $\kappa I\le G_\rho\le\kappa^{-1}I$, then, for $q\ge1$,
\begin{align}
 \|D_H^qG_\rho^{-1}[W_1,\ldots,W_q]\|
 &\le C_{q,\kappa}N^{-q/2}\rho^{5q}e^{2q\rho s},\label{W:eq:gram-window-inverse}\\
 \|D_H^qP_\rho[W_1,\ldots,W_q]\|
 &\le C_{q,\kappa}N^{-q/2}\rho^{5q+3}e^{2(q+1)\rho s}.\label{W:eq:gram-window-proj}
\end{align}
\end{proposition}
For $q>n$, $D_H^qp_n=0$. The displayed bounds are polynomial in $N$ when $\rho\le C_*N^{1/3}$ and
\begin{equation}
 d\le A_0(\log N)^2/\rho^2,\qquad A_0<\infty.
 \label{W:eq:gram-window-width}
\end{equation}
\begin{proof}
For $|x|\le2$, $|p_j(x)|\le j+1$. Outside that interval use \eqref{W:eq:gram-hyperbolic} and parity. This bounds $\|p_j(H)\|$ by $(j+1)e^{js}$. Substitute it into
\eqref{W:eq:chebyshev-matrix-derivatives} and use $\sum_{j_0+\cdots+j_q=n-q}\prod_{a=0}^q(j_a+1) =\binom{n+q+1}{2q+1}$, obtained from $(1-t)^{-2(q+1)}$. This proves
\eqref{W:eq:gram-window-poly}. For normalized matrix entry directions $\|W_j\|\le C N^{-1/2}$, the squared Hilbert--Schmidt norm of the polynomial matrix is bounded by
\[
 \begin{aligned}
 \|D_H^qF_\rho[W_1,\ldots,W_q]\|_{\rm HS}^2
 &=\sum_{n=0}^{\rho-1}\|D_H^qp_n(H)[W_1,\ldots,W_q]w\|^2\\
 &\le C_qN^{-q}e^{2\rho s}\sum_{n<\rho}(n+1)^{4q+2}
 \le C_qN^{-q}\rho^{4q+3}e^{2\rho s}.
 \end{aligned}
\]
Taking square roots proves \eqref{W:eq:gram-window-polynomial-matrix}. For a label set $S$ of size $q$,
\[
 D_{H,S}(F^*F)=\sum_{A\sqcup B=S}(D_{H,A}F)^*D_{H,B}F.
\]
Each term, with $|A|=a$ and $|B|=b$, has operator norm at most $C_qN^{-(a+b)/2}\rho^{2a+3/2+2b+3/2}e^{2\rho s}$. There are $2^q$ assignments and $a+b=q$, giving
\eqref{W:eq:gram-window-Gram}. For a nonempty label set $S$, the exact inverse formula is
\begin{equation}
 D_{H,S}G^{-1}
 =\sum_{\ell=1}^{|S|}(-1)^\ell
   \sum_{(S_1,\ldots,S_\ell)\in\operatorname{OP}_\ell(S)}
 G^{-1}(D_{H,S_1}G)G^{-1}\cdots(D_{H,S_\ell}G)G^{-1}.
 \label{W:eq:gram-inverse-labels}
\end{equation}
Differentiating $GG^{-1}=I$ in the labels $S$ gives
\[
 G D_SG^{-1}=-\sum_{\varnothing\ne A\subseteq S}(D_AG)(D_{S\setminus A}G^{-1}).
\]
Multiplying by $G^{-1}$ and inducting on $|S|$ gives \eqref{W:eq:gram-inverse-labels}, with $A$ as the first ordered block. For an ordered partition $(S_1,\ldots,S_\ell)$ of $S$,
put $q_a=|S_a|$. Since $\sum_a q_a=q$ and $q_a\ge1$, its contribution satisfies
\[
 \left\|G^{-1}(D_{H,S_1}G)G^{-1}\cdots
                  (D_{H,S_\ell}G)G^{-1}\right\|
 \le \kappa^{-\ell-1}\prod_{a=1}^{\ell}
  \left(C_{q_a}N^{-q_a/2}\rho^{2q_a+3}e^{2\rho s}\right)
 \le C_{q,\kappa}N^{-q/2}\rho^{2q+3\ell}e^{2\ell\rho s}.
\]
There are $\ell!\left\{\begin{smallmatrix}q\\\ell\end{smallmatrix}\right\}$ ordered partitions into $\ell$ nonempty
sets. Their number depends only on $q$, and $\ell\le q$ implies $2q+3\ell\le5q$. Summing gives \eqref{W:eq:gram-window-inverse}. For the projection use the exact labelled identity
\[
 D_{H,S}(FG^{-1}F^*)=
 \sum_{A\sqcup B\sqcup C=S}(D_{H,A}F)(D_{H,B}G^{-1})(D_{H,C}F)^*.
\]
Write $a=|A|$, $b=|B|$, and $c=|C|$. If $b=0$, the norm of this summand is bounded by $C_{q,\kappa}N^{-q/2}\rho^{2q+3}e^{2\rho s}$. If $b\ge1$,
\eqref{W:eq:gram-window-polynomial-matrix} and the inverse estimate give
\[
 C_{q,\kappa}N^{-(a+b+c)/2}
 \rho^{2a+3/2+5b+2c+3/2}e^{2(b+1)\rho s}
 =C_{q,\kappa}N^{-q/2}\rho^{2q+3b+3}e^{2(b+1)\rho s}.
\]
Both bounds are at most the right side of \eqref{W:eq:gram-window-proj}, since $b\le q$, $\rho\ge1$, and $s\ge0$. The sum has $3^q$ labelled terms. Finally $s\le\sqrt d$, since
$\cosh s\ge1+s^2/2$. Under \eqref{W:eq:gram-window-width},
$e^{{\rho}s}\le N^{\sqrt{A_0}}$.

For comparison with the moment estimates, take the narrower width $d\le C\rho^{-2}(\log N)^\gamma$, $\gamma<2$, and $\rho\asymp N^{1/3}$. Then $e^{\rho s}=N^{o(1)}$, and the
respective bounds for $D_{H,S}F$, $D_{H,S}G$, $D_{H,S}G^{-1}$ and $D_{H,S}P$, in the norms used above, become
\[
 N^{q/6+1/2+o(1)},\qquad N^{q/6+1+o(1)},\qquad
 N^{7q/6+o(1)},\qquad N^{7q/6+1+o(1)}.
\]
These are upper bounds in normalized entry directions, with fixed $q,\kappa$. They give polynomial bounds for individual derivatives. To make the error small after summing over the
entries, the localized comparison also uses the cancellation in the cumulative moments.
\end{proof}

\Needspace{6\baselineskip}

\begin{lemma}
\label{W:lem:gram-clipping}
For each $0<d\le1/10$ there exists a real $h_d\in C_c^\infty(\mathbb R)$ such that
\[
 h_d(x)=x\quad(|x|\le2+d),\qquad |h_d(x)|\le2+2d,
 \qquad \operatorname{supp}h_d\subset[-4,4].
\]
For every fixed $q\ge1$,
\begin{equation}
 \sup_{H=H^*}\|D_{H}^qh_d(H)[W_1,\ldots,W_q]\|
 \le C_qd^{-(q+2)}\prod_j\|W_j\|.
 \label{W:eq:spectral-modification-derivatives}
\end{equation}

\end{lemma}
\begin{proof}
Choose a fixed smooth increasing $b:\mathbb R\to[0,1]$ with $b=0$ on $(-\infty,0]$ and $b=1$ on $[1,\infty)$, and put $\zeta(t)=\int_0^t b(s)\,\mathrm ds$. This is zero for $t\le0$
and $t-c_b$ for $t\ge1$, with $c_b\in[0,1]$. For $x\ge0$ define
\[
 h_d^0(x)=x-d\zeta\bigl((x-2-d)/d\bigr),
\]
and extend oddly. It equals $x$ for $|x|\le2+d$ and is constant in absolute value between $2+d$ and $2+2d$ outside $[-2-2d,2+2d]$. Multiply by a fixed even smooth function in
$[0,1]$ that equals one on $[-3,3]$ and vanishes outside $[-4,4]$. The resulting $h_d$ has the claimed value and support bounds, and $\|h_d\|_1+\|h_d^{(q+2)}\|_1\le C_qd^{-(q+2)}$.
For $\widehat h_d(t)=\int h_d(x)e^{-\mathrm itx}\,\mathrm dx$, the estimates $|\widehat h_d(t)|\le\|h_d\|_1$ and $|t|^{q+2}|\widehat h_d(t)|\le\|h_d^{(q+2)}\|_1$ give
\[
 \begin{aligned}
 \int_{\mathbb R}|t|^q|\widehat h_d(t)|\,\dd t
 &\le\|h_d\|_1\int_{-1}^{1}|t|^q\,\dd t
   +\|h_d^{(q+2)}\|_1\int_{|t|>1}|t|^{-2}\,\dd t\\
 &=\frac2{q+1}\|h_d\|_1+2\|h_d^{(q+2)}\|_1.
 \end{aligned}
\]
Fourier inversion and the spectral theorem express $h_d(H)$ as the norm convergent integral $(2\pi)^{-1}\int\widehat h_d(t)e^{\mathrm itH}\,\mathrm dt$. Let
$\Delta_q=\{(s_0,\ldots,s_q):s_j\ge0,\ \sum_js_j=1\}$, with the $q$-dimensional coordinate measure of total mass $1/q!$. Repeated differentiation of Duhamel's formula gives
\[
 D_H^qe^{\mathrm itH}[W_1,\ldots,W_q]
 =(\mathrm it)^q\sum_{\sigma\in S_q}\int_{\Delta_q}
 e^{\mathrm its_0H}W_{\sigma(1)}e^{\mathrm its_1H}\cdots
                  W_{\sigma(q)}e^{\mathrm its_qH}\,\dd s.
\]
All exponentials are unitary. Thus its norm is at most $|t|^q q!\operatorname{vol}(\Delta_q)\prod_j\|W_j\| =|t|^q\prod_j\|W_j\|$. Consequently,
\[
 \|D_H^qh_d(H)[W_1,\ldots,W_q]\|
 \le\frac1{2\pi}\int |t|^q|\widehat h_d(t)|\,\dd t
                                    \prod_j\|W_j\|
 \le C_qd^{-(q+2)}\prod_j\|W_j\|.
\]
The same integrable majorant justifies differentiation under the integral. This proves \eqref{W:eq:spectral-modification-derivatives}.
\end{proof}

\Needspace{6\baselineskip}

Let $f\in C^q(\mathbb R;\mathbb C)$, $q\ge1$, and let $H=\sum_i\lambda_iE_i$ be its spectral decomposition. Write $f^{[q]}$ for the divided difference, extended continuously to
repeated arguments, and put $M_q=\|f^{(q)}\|_\infty$.

\begin{lemma}
\label{W:lem:smooth-functional-calculus}
For fixed Hermitian matrices $W_1,\ldots,W_q$ used as variations of $H$,
\begin{equation}
 D_H^qf(H)[W_1,\ldots,W_q]
 =\sum_{\sigma\in S_q}\sum_{i_0,\ldots,i_q}
 f^{[q]}(\lambda_{i_0},\ldots,\lambda_{i_q})
 E_{i_0}W_{\sigma(1)}E_{i_1}\cdots W_{\sigma(q)}E_{i_q}.
 \label{W:eq:smooth-matrix-derivative-formula}
\end{equation}
If $\|f^{(q)}\|_\infty<\infty$, then
\[
 \|D_H^qf(H)[W_1,\ldots,W_q]\|_{\mathrm{HS}}
 \le\|f^{(q)}\|_\infty\prod_j\|W_j\|_{\mathrm{HS}}.
\]
\end{lemma}
\begin{proof}
We first prove the formula for a polynomial by finite matrix multiplication. Let $f(x)=x^m$. If $q>m$ both the derivative and the $q$th divided difference vanish. For $q\le m$,
repeated application of the product rule to $H^m$ gives
\begin{equation*}
 D_H^qH^m[W_1,\ldots,W_q]
 =\sum_{\sigma\in\mathfrak S_q}
  \sum_{a_0+\cdots+a_q=m-q}
    H^{a_0}W_{\sigma(1)}H^{a_1}\cdots
    W_{\sigma(q)}H^{a_q}.
\end{equation*}
The divided differences of a monomial obey
\[
 (x^m)^{[q]}(x_0,\ldots,x_q)
 =\sum_{a_0+\cdots+a_q=m-q}\prod_{j=0}^qx_j^{a_j}.
\]
At $q=0$ this is the monomial itself. In the recursion $f^{[q]}(x_0,\ldots,x_q)= \bigl(f^{[q-1]}(x_1,\ldots,x_q)- f^{[q-1]}(x_0,\ldots,x_{q-1})\bigr)/(x_q-x_0)$, the terms
containing positive powers of $x_q$ and $x_0$ telescope by
$(x_q^{r+1}-x_0^{r+1})/(x_q-x_0)
  =\sum_{j=0}^r x_q^{r-j}x_0^j$.
This gives the displayed sum by induction when the nodes are
distinct, and then at repeated nodes by continuity.
Substituting $H=\sum_i\lambda_iE_i$ into each monomial product
identifies the factors of eigenvalues with this divided
difference. Linearity proves
\eqref{W:eq:smooth-matrix-derivative-formula} for all polynomials.

Next let $f$ be a polynomial, and fix a compact interval $I$ containing all the eigenvalues of the matrices under consideration. Put
$M_{q,I}=\sup_{x\in I}|f^{(q)}(x)|$.
The Genocchi--Hermite formula, in the normalization of \cite[Section~9, equation~(52)]{deboor2005divided}, is
\[
 f^{[q]}(x_0,\ldots,x_q)
 =\int_{\substack{t_1,\ldots,t_q\ge0\\t_1+\cdots+t_q\le1}}
 f^{(q)}\!\left(\left(1-\sum_{j=1}^qt_j\right)x_0+
                             \sum_{j=1}^qt_jx_j\right)\,\dd t_1\cdots\dd t_q.
\]
The integration region has volume $1/q!$. Its argument is a convex combination of the nodes and hence belongs to $I$ when all nodes belong to $I$. Consequently
\[
 |f^{[q]}(\lambda_0,\ldots,\lambda_q)|\le M_{q,I}/q!.
\]
The same integral provides the continuous extension to coinciding nodes. It also gives uniform convergence of the $q$th divided differences whenever the $q$th scalar derivatives
converge uniformly on $I$. Choose an orthonormal eigenbasis for $H$ consisting of individual eigenvectors, even when some eigenvalues coincide, and write $|W|_{\rm ent}$ for
entrywise absolute values in this basis. For each permutation $\sigma$, the corresponding summand in \eqref{W:eq:smooth-matrix-derivative-formula} satisfies
\[
 \left|\bigl(D_\sigma f(H)\bigr)_{ij}\right|
 \le\frac{M_{q,I}}{q!}
     \bigl(|W_{\sigma(1)}|_{\rm ent}\cdots
            |W_{\sigma(q)}|_{\rm ent}\bigr)_{ij}.
\]
Entrywise domination implies domination of the Hilbert--Schmidt norm. By $\|AB\|_{\rm HS}\le\|A\|_{\rm HS}\|B\|$ and $\|B\|\le\|B\|_{\rm HS}$, each product on the right has norm at
most $\prod_a\|W_a\|_{\rm HS}$; taking entrywise absolute values does not change these individual norms. There are exactly $q!$ permutations. The triangle inequality gives
\[
 \|D_H^qf(H)[W_1,\ldots,W_q]\|_{\rm HS}
 \le M_{q,I}\prod_{a=1}^q\|W_a\|_{\rm HS}.
\]
This estimate is independent of the matrix dimension and includes coinciding eigenvalues.

Finally, let $f\in C^q(\mathbb R;\mathbb C)$. Fix $H$ and a compact interval $I$ containing the spectra of all matrices in a sufficiently small neighbourhood of $H$. There are
polynomials $f_n$ such that $f_n^{(j)}\to f^{(j)}$ uniformly on $I$ for every $0\le j\le q$: approximate the continuous function $f^{(q)}$ uniformly by polynomials, integrate $q$
times, and choose the integration constants to match the first $q$ derivatives of $f$ at one point of $I$. The preceding estimate, applied with $f_n-f_m$ and each order $j\le q$,
makes $D_H^jf_n(H)$ uniformly Cauchy for matrices with spectra in $I$ and variations with bounded Hilbert--Schmidt norm. For $j=0$, uniform convergence of the functions gives
$f_n(K)\to f(K)$ in operator norm, uniformly for $\sigma(K)\subset I$. Applying the fundamental theorem of calculus on each line segment of matrices and proceeding by induction on
$j$ identifies the limits of the derivatives with the successive derivatives of $K\mapsto f(K)$. The Genocchi--Hermite integral formula also shows that $f_n^{[q]}\to f^{[q]}$
uniformly on $I^{q+1}$, including repeated arguments. Passing to the limit in the finite spectral sums proves \eqref{W:eq:smooth-matrix-derivative-formula} for $f$. Passing to the
limit gives the inequality with $\sup_{x\in I}|f^{(q)}(x)|$ for every spectral interval $I$ as above. If $M_q=\|f^{(q)}\|_\infty<\infty$, this is at most $M_q$ and proves the
global bound.
\end{proof}

In an orthonormal eigenbasis the same identity reads
\begin{equation}
 (D_H^qf(H)[W_1,\ldots,W_q])_{a_0a_q}
 =\sum_{\sigma\in S_q}\sum_{a_1,\ldots,a_{q-1}}
 f^{[q]}(\lambda_{a_0},\ldots,\lambda_{a_q})
 \prod_{j=1}^q(W_{\sigma(j)})_{a_{j-1}a_j}.
 \label{W:eq:cumulant-divided-differences}
\end{equation}

For $\mathcal D=D_H^qf(H)[W_1,\ldots,W_q]$, Hilbert--Schmidt duality gives
\begin{equation}
 \begin{aligned}
 \|\mathcal D\|_{\mathrm{HS}}&\le M_q\prod_j\|W_j\|_{\mathrm{HS}},\\
 |w^*\mathcal Dw|&\le\|\mathcal D\|_{\mathrm{HS}},\qquad
 |\operatorname{Tr}\mathcal D|\le\sqrt N\|\mathcal D\|_{\mathrm{HS}}
          \quad(\|w\|=1).
 \end{aligned}
 \label{W:eq:smooth-matrix-trace-bounds}
\end{equation}
For $Q_w=ww^*-N^{-1}I$ and $\mathscr S_f(H,w)=\operatorname{Tr}[Q_wf(H)]$, one has
\begin{equation}
 D_{H}^q\mathscr S_f(H,w)[W_1,\ldots,W_q]
 =\sum_{\sigma\in S_q}\sum_{i_0,\ldots,i_q}
 f^{[q]}(\lambda_{i_0},\ldots,\lambda_{i_q})
 \operatorname{Tr}[Q_w E_{i_0}W_{\sigma(1)}E_{i_1}
                 \cdots W_{\sigma(q)}E_{i_q}],
 \label{W:eq:centered-spectral-derivatives}
\end{equation}
and $|D_{H}^q\mathscr S_f|\le M_q\prod_j\|W_j\|_{\mathrm{HS}}$.

For a normalized trace, division by $N$ gives $N^{-1/2}\|\mathcal D\|_{\mathrm{HS}}$. In particular, $W_j=N^{-1/2}T_{\alpha,j}$ with uniformly bounded
$\|T_{\alpha,j}\|_{\mathrm{HS}}$ contributes $O_q(N^{-q/2})$.

\Needspace{6\baselineskip}

Fix $A_0>0$, $0<\kappa<1/4$, and $c_*N^{1/3}\le\rho\le C_*N^{1/3}$. For sufficiently large $N$, set
\[
 d=A_0(\log N)^2/\rho^2,\qquad H^{[d]}=h_d(H),\qquad
 F^{[d]}=(p_j(H^{[d]})w)_{0\le j<\rho},
\]
and write $G^{[d]}=(F^{[d]})^*F^{[d]}$ and $S^{[d]}=(F^{[d]})^*H^{[d]}F^{[d]}$. Choose $b_\kappa\in C_c^\infty(\mathbb R;[0,\infty))$ equal to $x-\kappa/2$ near
$[\kappa,\kappa^{-1}]$, and $\chi_\kappa\in C_c^\infty((\kappa/2,2/\kappa);[0,1])$ equal to one near that interval. Put $a_\kappa(x)=\chi_\kappa(x)/x$, extended by zero, and define
\[
 \begin{aligned}
 G_\kappa^{[d]}&=\kappa I/2+b_\kappa(G^{[d]}),&
 G_\kappa^{[d]}&=LL^*,\\
 P_\kappa^{[d]}&=F^{[d]}a_\kappa(G^{[d]})(F^{[d]})^*,&
 J_\kappa^{[d]}&=L^{-1}S^{[d]}L^{-*}.
 \end{aligned}
\]
Here $L$ is lower triangular with positive diagonal. These definitions use the full Gram matrix of the specified size. For a different size we perform the construction with its own
Gram matrix.

\begin{theorem}
\label{W:thm:gram-extension}
The six matrix functions $F^{[d]},G^{[d]},S^{[d]},G_\kappa^{[d]},P_\kappa^{[d]},J_\kappa^{[d]}$ are smooth, and
\[
 \kappa I/2\le G_\kappa^{[d]}\le C_\kappa I,\qquad
 0\le P_\kappa^{[d]}\le I.
\]
On $\|H\|\le2+d/2$ and $\kappa I\le G_\rho(H,w)\le\kappa^{-1}I$, the last three agree locally with the Gram matrix, the Krylov projection, and the Jacobi compression. For
$0\le q\le5$, their derivatives in directions $\|W_j\|_{\rm HS}\le C N^{-1/2}$ have norms at most $C_qN^{B_q}$, uniformly over Hermitian $H$, for finite exponents $B_q$ depending
only on the fixed parameters.
\end{theorem}
\begin{proof}
All derivatives below are real mixed derivatives in fixed self-adjoint matrix variations $W_1,\ldots,W_q$ satisfying $\|W_j\|_{\rm HS}\le C N^{-1/2}$. Constants may depend on
$A_0,\kappa,c_*,C_*$, the fixed derivative order, and the normalization constant for these variations, but not on the matrix or the deterministic unit vector $w$.

We begin with $F^{[d]}=(p_j(h_d(H))w)_{0\le j<\rho}$. The finite-dimensional chain rule, indexed by partitions of $\{1,\ldots,q\}$, reads
\[
 D_H^q(p_n\circ h_d)(H)[W_1,\ldots,W_q]
 =\sum_{\pi\in\Pi(q)}
 D^{|\pi|}p_n(h_d(H))
       [D_Bh_d(H)[W_a:a\in B]:B\in\pi].
\]
If $\pi$ has $k$ blocks, Proposition~\ref{W:prop:gram-window} and Lemma~\ref{W:lem:gram-clipping} bound its norm by
\begin{align*}
 C_q\rho^{2k+1}e^{\rho\sqrt{2d}}
  \prod_{B\in\pi}
    \left(d^{-(|B|+2)}\prod_{a\in B}\|W_a\|\right)
 &\le C_q\rho^{2q+1}d^{-(q+2k)}e^{\rho\sqrt{2d}}
       \prod_{a=1}^q\|W_a\|\\
 &\le C_q\rho^{2q+1}d^{-3q}e^{\rho\sqrt{2d}}
       \prod_{a=1}^q\|W_a\|.
\end{align*}
Here $1\le k\le q$, $\rho\ge1$, and $0<d\le1$. Summing at most $\rho$ squared column norms and then taking the square root gives
\begin{equation}
 \|D_H^qF^{[d]}[W_1,\ldots,W_q]\|_{\rm HS}
 \le C_qN^{-q/2}\rho^{2q+3/2}d^{-3q}
                          e^{\rho\sqrt{2d}}.
 \label{W:eq:gram-global-polynomial-matrix}
\end{equation}
At order zero, $\|F^{[d]}\|_{\rm HS}\le C\rho^{3/2}e^{\rho\sqrt{2d}}$ by the scalar Chebyshev bound. Put $\alpha=\sqrt{2A_0}$. Because $\rho\asymp N^{1/3}$ and
$d=A_0(\log N)^2/\rho^2$,
\[
 e^{\rho\sqrt{2d}}=N^\alpha,
 \quad d^{-1}\le C N^{2/3},
 \quad\|D_H^qF^{[d]}\|_{\rm HS}
       \le C_qN^{3q+1+\alpha},\qquad 0\le q\le5.
\]
This bound also holds at $q=0$.

Since $G^{[d]}=F^{[d]*}F^{[d]}$ and $S^{[d]}=F^{[d]*}h_d(H)F^{[d]}$, the labelled product rules are
\[
 D_SG^{[d]}=\sum_{A\sqcup B=S}(D_AF^{[d]})^*D_BF^{[d]},
\]
\[
 D_SS^{[d]}=\sum_{A\sqcup B\sqcup C=S}
       (D_AF^{[d]})^*D_Bh_d(H)D_CF^{[d]}.
\]
Here the empty derivative of $h_d(H)$ has norm at most three, whereas a nonempty derivative of order $b$ satisfies
$\|D_Bh_d(H)\|\le C_bN^{-b/2}d^{-(b+2)} \le C_bN^{b/6+4/3}\le C_bN^{3b}$. For each partition, the exponents of the three factors sum to at most $3q+2+2\alpha$. There are $2^q$ or
$3^q$ assignments, respectively. Therefore
\[
 \|D_H^qG^{[d]}\|+\|D_H^qS^{[d]}\|
 \le C_qN^{3q+2+2\alpha},\qquad 0\le q\le5.
\]

The scalar $b_\kappa$ has compact support and bounded derivatives of every order. Its Fourier representation and repeated Duhamel formula give
$\|D^kb_\kappa(B)[E_1,\ldots,E_k]\| \le C_{k,\kappa}\prod_a\|E_a\|$ for self-adjoint $B$. The chain rule for $b_\kappa(G^{[d]})$, with $k$ nonempty blocks, is thus bounded by
$C_{q,\kappa}N^{3q+(2+2\alpha)k}$.
As $1\le k\le q$,
\[
 \|D_H^qG_\kappa^{[d]}\|
 \le C_{q,\kappa}N^{q(5+2\alpha)},
 \qquad
 \|D_H^qG_\kappa^{[d]}\|_{\rm HS}
 \le C_{q,\kappa}N^{1/6+q(5+2\alpha)}
\]
for $q\ge1$. The second inequality uses $\rho\le C_*N^{1/3}$. Since $b_\kappa\ge0$ and is bounded,
\[
 \frac\kappa2 I\le G_\kappa^{[d]}
 \le\left(\frac\kappa2+\|b_\kappa\|_\infty\right)I.
\]
In particular its lower Cholesky factor $L$ and its inverse have operator norms bounded by a constant depending only on $\kappa$.

For a nonempty set of labels $S$, differentiate $G_\kappa^{[d]}=LL^*$ and isolate the two terms where all labels fall on one factor:
\[
 (D_SL)L^*+L(D_SL)^*
 =D_SG_\kappa^{[d]}
   -\sum_{\varnothing\ne A\subsetneq S}(D_AL)(D_{S\setminus A}L)^*
 =:Y_S.
\]
Conjugating by $L^{-1}$, and using that $L^{-1}D_SL$ is lower triangular with real diagonal, yields exactly
\begin{equation}
 D_SL=L\mathcal T(L^{-1}Y_SL^{-*}),
 \label{W:eq:gram-cholesky}
\end{equation}
where $\mathcal T$ takes the strict lower triangle and half of the diagonal. It contracts the Hilbert--Schmidt norm. Thus, writing $L_q$ for the supremum of the $q$th derivative of
$L$ in the indicated directions,
\[
 L_q\le C_{q,\kappa}\left[
    N^{1/6+q(5+2\alpha)}
   +\sum_{a=1}^{q-1}\binom qaL_aL_{q-a}\right].
\]
Set $M=31/6+2\alpha$. Since $1/6+q(5+2\alpha)\le qM$ for $q\ge1$, induction gives $L_q\le C_{q,\kappa}N^{qM}$, with an empty sum at $q=1$. The inverse identity, obtained by
differentiating $LL^{-1}=I$, is a sum over ordered partitions:
\[
 D_SL^{-1}
 =\sum_{k=1}^{|S|}(-1)^k
   \sum_{(S_1,\ldots,S_k)\in\operatorname{OP}_k(S)}
 L^{-1}(D_{S_1}L)L^{-1}\cdots(D_{S_k}L)L^{-1}.
\]
For fixed $q=|S|$, each term has exponent $M\sum_a|S_a|=qM$, and the number of ordered partitions depends only on $q$. It follows that
$\|D_SL^{-1}\|\le C_{q,\kappa}N^{qM}$.

For the regularized multiplication matrix
$J_\kappa^{[d]}=L^{-1}S^{[d]}L^{-*}$,
\[
 D_SJ_\kappa^{[d]}
 =\sum_{A\sqcup B\sqcup C=S}
  (D_AL^{-1})(D_BS^{[d]})(D_CL^{-*}).
\]
At derivative orders $a,b,c$, with $a+b+c=q$, the corresponding term has norm at most $C_{q,\kappa}N^{(a+c)M+3b+2+2\alpha} \le C_{q,\kappa}N^{qM+2+2\alpha}$, since $M\ge3$. The
$3^q$ assignments give the same bound for the sum. The scalar function $a_\kappa$ is also smooth and compactly supported. Its composition with $G^{[d]}$ has derivatives of order
$b\ge1$ bounded by $C_{b,\kappa}N^{b(5+2\alpha)}$, while the matrix itself has bounded norm. Differentiating $P_\kappa^{[d]}=F^{[d]}a_\kappa(G^{[d]})F^{[d]*}$ gives
\[
 \|D_SP_\kappa^{[d]}\|
 \le C_{q,\kappa}N^{q(5+2\alpha)+2+2\alpha}.
\]
Indeed, at allocations $(a,b,c)$ each factor contributes $3a+1+\alpha$, $(5+2\alpha)b$, and $3c+1+\alpha$, respectively, and $3\le5+2\alpha$. These estimates provide a finite
polynomial exponent for the derivatives of every one of the six matrix functions through order five, uniformly in $H$ and $w$.

It remains to verify positivity and local agreement with the unmodified matrices. Let $F^{[d]}=U\Sigma V^*$ be a singular value decomposition with singular values $s_j$. Since
$a_\kappa(x)=\chi_\kappa(x)/x$ for $x>0$ and $\chi_\kappa$ vanishes near zero,
\[
 P_\kappa^{[d]}
 =U\operatorname{diag}(\chi_\kappa(s_1^2),\ldots,
                            \chi_\kappa(s_r^2),0,\ldots,0)U^*,
\]
which implies $0\le P_\kappa^{[d]}\le I$. If $\|H\|\le2+d/2$, the scalar function $h_d$ is the identity on a neighbourhood of $\sigma(H)$. If also
$\kappa I\le G_\rho(H,w)\le\kappa^{-1}I$, the functions $b_\kappa$ and $a_\kappa$ agree locally with $x-\kappa/2$ and $x^{-1}$, respectively, on the spectrum of $G_\rho$. Hence
$G_\kappa^{[d]}=G_\rho$ and $P_\kappa^{[d]}=F G_\rho^{-1}F^*$, the orthogonal projection onto the first $\rho$ Krylov vectors. The Cholesky factor with positive diagonal is the one
used in \eqref{W:eq:exact-jacobi-block}, so $J_\kappa^{[d]}$ is the initial Lanczos Jacobi block. The strict spectral margins and the scalar identities on neighbourhoods also give
local equality of all derivatives in these regions.
\end{proof}

For the applications, choose a smooth scalar function $\gamma_\kappa$ with values in $[0,1]$ and compact support in the positive interval for the eigenvalues of the Gram matrix,
and set $\widehat\Gamma=\det\gamma_\kappa(G^{[d]})$. Let $Y$ be a fixed finite linear combination of values of the integrated potential computed from $J_\kappa^{[d]}$. Finally put
\[
 \widehat C_E=\operatorname{Tr}[(I-P_\kappa^{[d]})
                \theta_E(m_N^2(2I-H))],\qquad
 \theta_E(x)=(1+e^{x-E})^{-2}.
\]
The parameters and frequencies are fixed. We include any fixed enlargement of the Gram matrix dimension in the same convention.

\begin{corollary}
\label{W:cor:global-matrix-tests}
The functions $\widehat\Gamma$, $\widehat\Gamma e^{\mathrm iY}$, and $\widehat\Gamma\widehat C_E/(1+\widehat C_E)$ are smooth and bounded in absolute value by one. Their matrix
derivatives through order five, in normalized entry directions, are bounded by fixed powers of $N$ uniformly over all Hermitian matrices.
\end{corollary}
\begin{proof}
The scalar cutoff $\gamma_\kappa$ takes values in $[0,1]$, so $A=\gamma_\kappa(G^{[d]})$ is a positive contraction. For fixed matrices $B_1,\ldots,B_k$, multilinearity of the
determinant in its columns gives the exact expansion of $D^k\det(A)[B_1,\ldots,B_k]$: choose $k$ distinct columns of $A$, assign the labelled matrices $B_j$ to them, and replace
the selected columns in that order. Terms in which two labels are assigned to the same column vanish, since the determinant is linear in each column. There are at most $\rho^k$
ordered choices. Hadamard's inequality bounds every determinant by $\prod_j\|B_j\|$, because the unreplaced columns of the contraction $A$ have norm at most one. Thus
\[
 |D^k\det(A)[B_1,\ldots,B_k]|
 \le\rho^k\prod_j\|B_j\|.
\]
By Fourier inversion for the fixed smooth function $\gamma_\kappa$, the spectral theorem, and the chain rule, every derivative of $\gamma_\kappa(G^{[d]})$ of order $a\ge1$ has norm
at most $C_aN^{a(5+2\alpha)}$, with $\alpha=\sqrt{2A_0}$, by the proof of Theorem~\ref{W:thm:gram-extension}. Applying the labelled chain rule for the determinant and summing the
partitions of $\{1,\ldots,q\}$, each term with $k$ blocks has size
\[
 C_q\rho^kN^{(5+2\alpha)q}
 \le C_qN^{q(5+2\alpha)+q/3}.
\]
There are finitely many partitions for fixed $q$, proving the same estimate for $D^q\widehat\Gamma$.

A fixed finite linear combination of integrated potentials from $J_\kappa^{[d]}$ contains at most $2\rho$ entries and a factor $m_N$. The derivative bound for $J_\kappa^{[d]}$
gives
\[
 |D^qY|
 \le C_qm_N\rho N^{qM+2+2\alpha}
 \le C_qN^{qM+8/3+2\alpha},
 \qquad M=31/6+2\alpha.
\]
Since $J_\kappa^{[d]}$ is self-adjoint, $Y$ is real for real matrix variations. For a set $S$ of $q$ derivative labels,
\[
 D_Se^{\mathrm iY}
 =e^{\mathrm iY}
    \sum_{\pi\in\Pi(S)}\mathrm i^{|\pi|}
                         \prod_{B\in\pi}D_BY.
\]
A partition with $k$ nonempty blocks gives at most $C_qN^{qM+k(8/3+2\alpha)}\le C_qN^{q(M+8/3+2\alpha)}$, as $|e^{\mathrm iY}|=1$. The product rule, including all derivatives of
$\widehat\Gamma$, therefore bounds the derivatives of $\widehat\Gamma e^{\mathrm iY}$ by a fixed power of $N$.

Put $\Theta(H)=\theta_E(m_N^2(2I-H))$. For $t=e^{x-E}$, repeated differentiation using $\partial_x=t\partial_t$ yields the exact formula
\[
 \theta_E^{(q)}(x)
 =\sum_{j=1}^q{q\brace j}(-1)^j(j+1)!
                      \frac{t^j}{(1+t)^{j+2}},\qquad q\ge1,
\]
where ${q\brace j}$ is a Stirling number of the second kind. For $t\ge0$, each $t^j/(1+t)^{j+2}\le1$, so
$\|\theta_E^{(q)}\|_\infty\le C_q$.
Lemma~\ref{W:lem:smooth-functional-calculus}, with $D_H(m_N^2(2I-H))[W]=-m_N^2W$, gives
\[
 \|D_H^q\Theta[W_1,\ldots,W_q]\|_{\rm HS}
 \le C_qm_N^{2q}\prod_{j=1}^q\|W_j\|_{\rm HS}
 \le C_qN^{q/6}.
\]
The same matrix is a positive contraction. Since
$0\le P_\kappa^{[d]}\le I$,
\[
 0\le\widehat C_E
 =\operatorname{Tr}[(I-P_\kappa^{[d]})\Theta]\le N.
\]
For a nonempty label set $S$, differentiate both matrix factors:
\[
 D_S\widehat C_E
 =\sum_{A\sqcup B=S}
   \operatorname{Tr}\bigl(D_A(I-P_\kappa^{[d]})D_B\Theta\bigr).
\]
If $a=|A|$, $b=|B|$, $a+b=q$, use
$|\operatorname{Tr}(UV)|\le N\|U\|\|V\|$.
For $a\ge1$ the first factor is bounded by $C_aN^{a(5+2\alpha)+2+2\alpha}$, and for $a=0$ its norm is at most one. For $b\ge1$ the second factor has operator norm at most
$C_bN^{b/6}$, and at $b=0$ at most one. Thus every term is bounded by $C_qN^{qM+3+2\alpha}$, after enlarging the fixed exponent, and the $2^q$ assignments give the same bound for
their sum.

The derivatives of $h(c)=c/(1+c)$ are $h^{(k)}(c)=(-1)^{k-1}k!/(1+c)^{k+1}$ for $k\ge1$; all have modulus at most $k!$ for $c\ge0$. For every partition $\pi$ into $k$ blocks,
\[
 |h^{(k)}(\widehat C_E)|\prod_{B\in\pi}
             |D_B\widehat C_E|
 \le C_qN^{qM+k(3+2\alpha)}
 \le C_qN^{q(M+3+2\alpha)}.
\]
The product rule for $\widehat\Gamma h(\widehat C_E)$ now provides a polynomial bound including every derivative of the cutoff. Combining the exponents shows that $q(9+4\alpha)$ is
a common valid choice for these functions when $1\le q\le5$. Finally, $0\le\widehat\Gamma\le1$, $|e^{\mathrm iY}|=1$, and $0\le h(\widehat C_E)<1$; hence the three functions are
bounded in absolute value by one. Their smoothness follows from the displayed compositions of smooth finite-dimensional matrix functions and the strictly positive denominator.
\end{proof}

The polynomial estimates also control the contribution from the complement of the event used for the Schur complement. Here the polynomial exponent must be fixed before the
probability exponent is chosen. For this calculation, let $K_\alpha$ be the matrix with entry block $\alpha$ removed, and let $\mathcal F_\alpha=\sigma(K_\alpha)$. Suppose
$\mathcal E_\alpha\in\mathcal F_\alpha$ and the replacement block has conditional distribution $\nu_\alpha(\dd x)$ given $\mathcal F_\alpha$. Assume that, almost surely,
\[
 \int|\mathcal R_\alpha(K_\alpha,x)|\,\nu_\alpha(\dd x)\le CN^B.
\]
This is a bound for the absolute integrand after the replacement variable has been integrated. In particular, a factor $1+\|x\|^k$ is included using its conditional moment bound.
For a complex entry, $\nu_\alpha$ is the joint distribution of its two real coordinates. Conditional expectation and $\mathcal F_\alpha$-measurability give
\[
 \begin{aligned}
 \mathbb E\left[\mathbf1_{\mathcal E_\alpha}
       \int|\mathcal R_\alpha(K_\alpha,x)|\,\nu_\alpha(\dd x)\right]
 &\le CN^B\mathbb P(\mathcal E_\alpha),\\
 \sum_{\alpha\le N^2}\mathbb E\left[\mathbf1_{\mathcal E_\alpha}
       \int|\mathcal R_\alpha(K_\alpha,x)|\,\nu_\alpha(\dd x)\right]
 &\le C_DN^{B+2-D}
 \quad\text{if }\sup_\alpha\mathbb P(\mathcal E_\alpha)\le C_DN^{-D}.
 \end{aligned}
\]
The integrands are nonnegative, so Tonelli's theorem justifies the conditional integral and finite sum. When a single event $\mathcal E_N$ is used for every entry, this calculation
gives the bound
$CN^{B+2}\mathbb P(\mathcal E_N)$.
Section~\ref{W:sec:smooth-cutoff-for-schur-complement} supplies polynomial derivative bounds for $U_N$ whose exponents are independent of the subsequently chosen failure exponent.
After including these derivatives and the required conditional moments in $B$, one may therefore choose $D>B+2+A$ to obtain an $O(N^{-A})$ contribution. The hypothesis above is
checked for the complete expression $\mathcal R_\alpha$, including every polynomial subtracted from the observable.

For an original matrix, let $\Omega_N=\{\|H_N\|\le2+d/2\}$. If $|F_N|,|\widehat F_N|\le M_0$ and the functions agree on $\Omega_N$, then
$|F_N-\widehat F_N|\le2M_0\mathbf1_{\Omega_N^c}$, so
\begin{equation}
 |\mathbb EF_N(H_N)-\mathbb E\widehat F_N(H_N)|
 \le2M_0\mathbb P(\Omega_N^c).
 \label{W:eq:gram-ensemble-error}
\end{equation}
For functions multiplied by the same cutoff in the eigenvalues of the Gram matrix, it suffices to check equality on the nonzero cutoff set. The resulting expectation estimate is
controlled by $\mathbb P(\Omega_N^c)$.

\subsubsection{Resolvent identities and entry expansions}
\label{W:sec:multi-resolvent-identities-and-entry-expansions}

Replacing one independent entry changes the matrix only on one or two coordinate axes. This finite support is useful in the resolvent expansion: the internal factors reduce to
small matrix compressions, leaving only the two exterior vectors to be estimated in the full space. The identities below make this reduction explicit. They are the deterministic
part of Green function comparison; see \cite{erdos2012rigidity,benaych2018lectures}. Probabilistic bounds for the resulting factors will be supplied under the particular entry
distributions used in the comparison.

Here $R_j=(z_jI-H)^{-1}$, and $D_S$ denotes the mixed derivative with respect to the fixed, labelled matrix variations $(V_s)_{s\in S}$. The distinction between variation of $H$
and variation of the spectral parameter is reflected in the signs:

\begin{equation}
 D R_j[V]=R_jVR_j,\qquad
 \partial_{z_j}R_j=-R_j^2.
 \label{W:eq:resolvent-signs}
\end{equation}
Indeed, differentiating $(z_jI-H)R_j=I$ gives
\[
 -VR_j+(z_jI-H)D R_j[V]=0,\qquad
 R_j+(z_jI-H)\partial_{z_j}R_j=0,
\]
which yields the two identities after left multiplication by $R_j$.

\begin{lemma}[Derivatives of resolvents]
\label{W:lem:resolvent-derivatives}
For $q\ge1$,
\begin{equation}
 D^qR(z)[V_1,\ldots,V_q]
 =\sum_{\sigma\in\mathfrak S_q}
 R(z)V_{\sigma(1)}R(z)\cdots V_{\sigma(q)}R(z).
 \label{W:eq:resolvent-permutations}
\end{equation}
In particular, $D^qR[V,\ldots,V]/q!=R(VR)^q$. If $\mathcal W=R_0B_1R_1\cdots B_kR_k$ and the $B_j$ are independent of $H$, then
\begin{equation}
 D_S\mathcal W
 =\sum_{S_0\sqcup\cdots\sqcup S_k=S}
 (D_{S_0}R_0)B_1(D_{S_1}R_1)\cdots B_k(D_{S_k}R_k).
 \label{W:eq:resolvent-chain-labels}
\end{equation}
Empty sets are allowed and $D_\varnothing R_j=R_j$. For $|S|=q$, expanding each derivative in \eqref{W:eq:resolvent-chain-labels} gives exactly $q!\binom{k+q}{q}$ labelled terms,
counted with multiplicity.
\end{lemma}
\begin{proof}
The case $q=1$ is \eqref{W:eq:resolvent-signs}. We induct on $q$. A summand in $D_H^qR[V_1,\ldots,V_q]$ contains $q$ labelled variations and $q+1$ resolvent factors.
Differentiating in $V_{q+1}$ replaces precisely one of these $q+1$ resolvents by $RV_{q+1}R$. Every ordering of the labels $1,\ldots,q+1$ is obtained once: delete label $q+1$ from
that ordering and record its position. This proves \eqref{W:eq:resolvent-permutations}. When all variations equal $V$, every permutation gives $R(VR)^q$; division by $q!$ therefore
gives the repeated variation formula.

Applying the product rule of Section~\ref{W:sec:notation-and-scales} to the $k+1$ resolvents, with each $B_j$ fixed, gives \eqref{W:eq:resolvent-chain-labels}. To count the
products after expanding the individual resolvent derivatives, let $q_j=|S_j|$ with $\sum_{j=0}^kq_j=q$. There are $q!/(q_0!\cdots q_k!)$ allocations of labelled variations to the
factors, and $\prod_jq_j!$ orders within the factors. Thus each fixed vector $(q_0,\ldots,q_k)$ contributes $q!$ products, counted with multiplicity. There are $\binom{q+k}{k}$
nonnegative integer solutions of $q_0+\cdots+q_k=q$, giving $q!\binom{k+q}{q}$ terms. When $B_j=B_j(H)$, the product rule also differentiates these factors, as specified below.
\end{proof}

The product formula also gives a deterministic norm bound. Write $\delta_j=\operatorname{dist}(z_j,\sigma(H))>0$ and $\delta_* =\min_{0\le j\le k}\delta_j$. For each weak
composition $q_0+\cdots+q_k=q$, the corresponding sum has $q!$ labelled terms and is bounded by
\[
 q!\prod_{j=0}^k\delta_j^{-q_j-1}
       \prod_{j=1}^k\|B_j\|\prod_{a=1}^q\|V_a\|.
\]
Thus, for the matrices $B_j$ held fixed in Lemma~\ref{W:lem:resolvent-derivatives},
\[
 \|D_H^q(R_0B_1R_1\cdots B_kR_k)[V_1,\ldots,V_q]\|
 \le q!\binom{k+q}{q}\delta_*^{-q}
       \prod_{j=0}^k\delta_j^{-1}\prod_{j=1}^k\|B_j\|
       \prod_{a=1}^q\|V_a\|.
\]
The vectors in a scalar product contribute the additional factor $\|u\|\|v\|$; an unnormalized trace is bounded by $N$ times this operator norm. The coordinate compressions below
give sharper bounds for the products arising from a change of one entry.

If a matrix $V$ is supported on coordinates $J$, $|J|\le2$, write $V=\iota_JT\iota_J^*$, where $\iota_J$ is the coordinate inclusion. For a repeated direction,
\begin{equation}
 \frac1{q!}D_H^qR_H(z)[V^q]
 =R_H\iota_JT(R_{H,JJ}T)^{q-1}\iota_J^*R_H,
 \qquad q\ge1.
 \label{W:eq:intro-finite-support-derivative}
\end{equation}
This is the repeated direction case of \eqref{W:eq:resolvent-permutations}, with the coordinate inclusions contracted between consecutive factors $T$.

For matrix dependent factors, set $F_{2j}=R_j$ and $F_{2j-1}=B_j(H)$. The product rule instead reads
\[
 D_S(F_0\cdots F_{2k})
 =\sum_{S_0\sqcup\cdots\sqcup S_{2k}=S}
       (D_{S_0}F_0)\cdots(D_{S_{2k}}F_{2k}).
\]
For example $D(R_0BR_1)[V]=(DR_0[V])BR_1+R_0(DB[V])R_1+R_0B(DR_1[V])$. We use this formula for the random inverse of the Gram matrix and polynomial projection.

Let $K$ be a matrix with the selected entry set to zero, $V=\iota_JT\iota_J^*$, and $R=(zI-K)^{-1}$. Here $\iota_J:\mathbb C^{|J|}\to\mathbb C^N$ is the coordinate inclusion and
$|J|\le2$. Write $C=\iota_J^*R\iota_J$. For every integer $p\ge1$, the finite resolvent identity is
\begin{equation}
 R_{K+V}=\sum_{j=0}^{p-1}R(VR)^j+(RV)^pR_{K+V}.
 \label{W:eq:resolvent-finite-resolvent}
\end{equation}
It follows by iterating $R_{K+V}=R+RVR_{K+V}$ exactly $p$ times. On any domain where $\|CT\|\le1/2$, the same identity becomes
\begin{equation}
 R_{K+V}-R =\sum_{j=1}^{p-1}R\iota_JT(CT)^{j-1}\iota_J^*R +\mathcal D_p, \qquad \mathcal D_p =R\iota_JT(CT)^{p-1}(I-CT)^{-1}\iota_J^*R.
 \label{W:eq:resolvent-block-remainder}
\end{equation}
Indeed, the block resolvent identity is
\[
 R_{K+V}-R=R\iota_JT(I-CT)^{-1}\iota_J^*R,
\]
and expanding only the finite-dimensional inverse gives \eqref{W:eq:resolvent-block-remainder}. In particular, for $p=4$ the polynomial part is the sum of $RTR$, $RTRTR$, and
$RTRTRTR$, with the coordinate inclusions understood. The last term contains four factors representing the entry replacement and five resolvents before its bounded block inverse is
summed. Since $\|(I-CT)^{-1}\|\le2$, its contribution can be bounded before integration. Write $A=R\iota_J$ and $B=\iota_J^*R$. The trace ideal inequality
$\|AXB\|_1\le\|A\|_{\rm HS}\|X\|\|B\|_{\rm HS}$ gives
\[
 \begin{aligned}
 |u^*\mathcal D_pv|
 &\le2\|\iota_J^*R^*u\|\|\iota_J^*Rv\|
                     \|T\|^p\|C\|^{p-1},\\
 \|\mathcal D_p\|_1
 &\le2\|R\iota_J\|_{\rm HS}\|\iota_J^*R\|_{\rm HS}
                     \|T\|^p\|C\|^{p-1}.
 \end{aligned}
\]
The factors involving exterior vectors are estimated by the column bounds for the resolvent. In particular $\|C\|\le L$ and $\|T\|\le Ce_N$ imply the expansion condition once
$LCe_N\le1/2$, even when $e_N$ exceeds the imaginary part of $z$. Section~\ref{W:sec:contour-integrals-for-moment-increments} proves these compressed resolvent bounds uniformly
over the contours and permitted entry values.

For possibly different spectral parameters put
\[
 C_j=\iota_J^*R_j\iota_J,\qquad
 D=T_1C_1T_2\cdots C_{q-1}T_q,\qquad
 \mathcal W_J=R_0\iota_JD\iota_J^*R_q.
\]
Here the $T_j$ are $|J|$-dimensional matrices, while the resolvents act on the $N$-dimensional space. We will use the following contractions.

\begin{lemma}[Matrix elements and traces of resolvent products]
\label{W:lem:resolvent-contractions}
For any vectors $u,v$,
\begin{align}
 |u^*\mathcal W_Jv|
 &\le\|\iota_J^*R_0^*u\|\,\|D\|\,
                    \|\iota_J^*R_qv\|,\label{W:eq:resolvent-starting-vector-contraction}\\
 \Tr\mathcal W_J
 &=\operatorname{tr}_{|J|}\!\bigl(D\iota_J^*R_qR_0\iota_J\bigr).
 \label{W:eq:resolvent-trace-contraction}
\end{align}
If $L_j$ bounds both $\|R_je_i\|^2$ and $\|R_j^*e_i\|^2$ for $i\in J$, then
\begin{equation}
 |\Tr\mathcal W_J|\le |J|^2\|D\|\sqrt{L_qL_0}.
 \label{W:eq:resolvent-trace-energy}
\end{equation}
For the normalized trace $N^{-1}\operatorname{Tr}$ the bound is $|J|^2N^{-1}\|D\|\sqrt{L_qL_0}$.
\end{lemma}
\begin{proof}
Writing $D=(D_{ab})_{a,b\in J}$ gives
\[
 \mathcal W_J=\sum_{a,b\in J}D_{ab}(R_0e_a)(e_b^*R_q).
\]
Therefore
\[
 u^*\mathcal W_Jv=(\iota_J^*R_0^*u)^*D(\iota_J^*R_qv),\qquad
 \operatorname{Tr}\mathcal W_J=\sum_{a,b\in J}D_{ab}(R_qR_0)_{ba}.
\]
Cauchy--Schwarz proves the first inequality. The finite sum gives the compressed trace in \eqref{W:eq:resolvent-trace-contraction}. For $i,j\in J$,
\[
 |(R_qR_0)_{ij}|
 =|(R_q^*e_i)^*(R_0e_j)|\le\sqrt{L_qL_0}.
\]
Summing the at most $|J|^2$ terms proves \eqref{W:eq:resolvent-trace-energy}. When $R_q=R_0=R$, splitting the summation index between $J$ and its complement gives
\[
 (R^2)_{JJ}=(R_{JJ})^2+R_{J J^c}R_{J^cJ}.
\]
Both summands are included in the preceding bound for the squared column norms.
\end{proof}

For $z=E+\mathrm i\eta$, $\eta>0$, normality and the spectral theorem give
\begin{equation}
 S(z):=-\operatorname{Im}R(z)=\eta R(z)^*R(z),\qquad
 \sum_i|R_{ij}(z)|^2=\frac{-\operatorname{Im}R_{jj}(z)}\eta.
 \label{W:eq:resolvent-ward}
\end{equation}
The same identity holds for the squared row norm, and $\|R(z)w\|^2=w^*S(z)w/\eta$.

For the traceless matrix $P^\circ=ww^*-I/N$, the exact formula is
\begin{equation}
 \Tr(P^\circ\mathcal W_J)
 =w^*\mathcal W_Jw-\frac1N\Tr\mathcal W_J.
 \label{W:eq:resolvent-centered-contraction}
\end{equation}
The two terms have the bounds
\[
 |\operatorname{Tr}(P^\circ\mathcal W_J)|
 \le\|\iota_J^*R_0^*w\|\|D\|\|\iota_J^*R_qw\|
       +\frac{|J|^2}{N}\|D\|\sqrt{L_qL_0}.
\]
This estimates $\operatorname{Tr}(P^\circ\mathcal W_J)$ by the two exterior resolvent vectors and the column norms $L_0,L_q$.

For the interpolating matrices $H(t)$ of Section~\ref{W:sec:cumulant-cubic}, the isotropic local semicircle law gives
\begin{equation}
 \left|u^*(R_H(z)-g_t(z)I)v\right|
 \prec \left(\sqrt{\frac{|\operatorname{Im}g_t(z)|}{N\eta}}
                    +\frac1{N\eta}\right)\|u\|\|v\|,
 \qquad \eta=\operatorname{Im}z>0.
 \label{W:eq:interpolating-isotropic-law}
\end{equation}
Here $u,v\in\mathbb F_\beta^N$ are fixed before the matrix is sampled; there is no restriction other than the displayed factors $\|u\|\|v\|$, so the estimate can be read for unit
vectors and rescaled. In the comparison interpolation parameter $t$, the resolvent is evaluated at the random matrix $H(t)$, while $g_t(z)=a_t^{-1}g(z/a_t)$ is the deterministic
semicircle transform with its corresponding variance normalization. The normalization and moment assumptions needed for the isotropic local semicircle law are verified for the
bounded comparison matrices in Lemmas~\ref{W:lem:cumulant-event} and \ref{W:lem:localization-interpolation-event}. Theorem~2.12 of \cite{bloemendal2014isotropic}, after normalizing
the variance row sums, supplies \eqref{W:eq:interpolating-isotropic-law} and requires, for a fixed $\omega>0$, that
$\eta=\operatorname{Im}z\ge N^{-1+\omega}$.
For the nonreal portions at height comparable to $N^{-2/3}$, this requirement is satisfied after choosing $0<\omega<1/3$. The real portions are treated separately by the spectral
gap and Lemma~\ref{W:lem:vertical-resolvent-comparison}. To quantify the probability after a finite union, choose the exponent $D+2$ in the uniform bound for a coordinate pair
before summing over the at most $N^2$ pairs; the resulting failure probability is at most $N^{-D}$.

The density estimate at the upper edge has the following precise range. If $\eta\asymp N^{-2/3}$ and $|E-2a_t|\le C\eta$, the formula
\[
 g_t(z)=\frac{z-\sqrt{z^2-4a_t^2}}{2a_t^2},\qquad a_t=1+O(N^{-1}),
\]
with the branch asymptotic to $z$ at infinity, gives $|\operatorname{Im}g_t(z)|\le C\sqrt\eta$. Indeed, $|z-2a_t|\le C\eta$ and $|z+2a_t|\le C$ imply
$|\sqrt{z^2-4a_t^2}|\le C\sqrt\eta$. Substitution gives
\[
 \sqrt{\frac{|\operatorname{Im}g_t(z)|}{N\eta}}+\frac1{N\eta}
 \le C\bigl(N^{-1/2}\eta^{-1/4}+N^{-1}\eta^{-1}\bigr)
 \le CN^{-1/3}.
\]
On the wider portion $|E-2a_t|\le Cd_N$, the same calculation gives $|\operatorname{Im}g_t(z)|\le C\sqrt{d_N+\eta}=N^{-1/3+o(1)}$, and hence a bound $N^{-1/3+o(1)}$ at the stated
height. Elsewhere we use the density in \eqref{W:eq:interpolating-isotropic-law} itself. Uniformity here means uniform failure probabilities for fixed deterministic choices of
$t,u,v$; a finite family is handled by the union calculation above.

For an entry replacement supported on $J$, $|J|\le2$, write $V_j=\iota_JT_j\iota_J^*$ and $C_j=\iota_J^*R_j\iota_J$. The exact identity
\[
 R_0V_1R_1\cdots V_qR_q
 =R_0\iota_J(T_1C_1\cdots C_{q-1}T_q)
                         \iota_J^*R_q
\]
contains every intermediate resolvent compression; there are $q-1$ of them. If $\|C_j\|\le C$ and the squared Euclidean norm of each column of the exterior resolvent factors is at
most $Q_Nm_N$, then Lemma~\ref{W:lem:resolvent-contractions} gives
\begin{equation}
\begin{aligned}
 |u^*R_0V_1\cdots V_qR_qv|
 &\le C^{q-1}\prod_{j=1}^q\|T_j\|
     \|\iota_J^*R_0^*u\|\|\iota_J^*R_qv\|,\\
 |\operatorname{Tr}(R_0V_1\cdots V_qR_q)|
 &\le4C^{q-1}Q_Nm_N\prod_{j=1}^q\|T_j\|.
\end{aligned}
\label{W:eq:finite-support-resolvent-estimates}
\end{equation}
The factor four comes from the at most $|J|^2=4$ terms in the trace compression; the bound is not multiplied by $N$. For the centered matrix $ww^*-N^{-1}I$, the exact trace splits
into the matrix element at $w$ and $N^{-1}$ times the unnormalized trace, so the same two bounds apply separately. For the truncated heavy-tailed arrays, these pointwise contour
bounds are replaced by the conditional estimates of Sections~\ref{W:sec:row-estimates-under-gram-constraint} and \ref{W:sec:contour-integrals-for-moment-increments}.

For the deterministic matrix of third moment differences $B_{3,N}$ in \eqref{W:eq:cumulant-third-matrix}, write its spectral decomposition $B_{3,N}=\sum_j\lambda_jv_jv_j^*$ with
$\|v_j\|=1$. The vectors $v_j$, and also the starting vector $w$, are fixed when $H(t)$ is sampled. In particular,
\[
 |\operatorname{Tr}(B_{3,N}R_H(z))|
 \le\sum_j|\lambda_j|\,v_j^*|R_H(z)|v_j.
\]
On the event of Lemma~\ref{W:lem:cumulant-event}, each quadratic form is at most $q_N$. The displayed sum is therefore bounded by $q_N\sum_j|\lambda_j|=q_N\|B_{3,N}\|_1$.

The composition formula of Section~\ref{W:sec:notation-and-scales} evaluates $D_s^{|\pi|}A$ at $s(H)$ and every argument $D_{H,B}s$ at $H$. For $A=\Gamma_\rho\chi_N A_0$, each
summand of the product derivative contains $D_{S_0}\Gamma_\rho$, $D_{S_1}\chi_N$ and $D_{S_2}A_0$ at this same point. A nonzero value of $D_{S_1}\chi_N(s(H))$ restricts $Q_*(s(H))$
and $B(s(H))$ by the support conditions in Section~\ref{W:sec:coefficient-localization}. At a replacement matrix $K_\alpha+V_\alpha(x)$, these conditions are applied to
$s(K_\alpha+V_\alpha(x))$. Differentiation of the entry distributions with respect to their scalar interpolation parameter is given separately in
Lemma~\ref{W:lem:entry-distribution-interpolation}.

For full-rank deterministic matrices $B_j$, the derivative of $R_0B_1R_1\cdots B_kR_k$ is still \eqref{W:eq:resolvent-chain-labels}. If $B_j$ depends on $H$, the product formula
above assigns labels to all $2k+1$ factors. The compressed estimates for $\mathcal W_J$ apply when the intervening factors have the displayed representation as matrices supported
on the selected coordinates; products of full-rank matrices require estimates for the full factors.

The general multi-resolvent result of \cite{cipolloni2022multiresolvent} concerns full-rank products under its stronger moment hypotheses, and the refinements of
\cite{cipolloni2023rankuniform,cipolloni2024edge} impose additional assumptions involving the imaginary part of the limiting Stieltjes transform. To compare these conditions with
our contour, the semicircle transform satisfies
\[
 g(z)=\tfrac12(z-\sqrt{z^2-4}),\quad
 |\operatorname{Im}g(2+O(\eta)+\mathrm i\eta)|\asymp\sqrt\eta.
\]
At $\eta\asymp N^{-2/3}$ one consequently has
\begin{equation}
 N\eta|\operatorname{Im}g(z)|\asymp
 N\,N^{-2/3}N^{-1/3}\asymp1.
 \label{W:eq:edge-density-local-law-scale}
\end{equation}
The hypotheses involving this quantity in \cite{cipolloni2023rankuniform,cipolloni2024edge}, which require this quantity to be at least a positive power of $N$, therefore fail on
these contour segments. At this scale we use the compressed resolvent estimates \eqref{W:eq:finite-support-resolvent-estimates}, followed by the conditional integrals in
Section~\ref{W:sec:contour-integrals-for-moment-increments}.

Let $R=(zI-H)^{-1}$ with $\operatorname{Im}z=\eta>0$ and let $T=T^*$. The resolvent $R$ is normal and $S=-\operatorname{Im}R=\eta R^*R\ge0$ by \eqref{W:eq:resolvent-ward}. For
$W_q=R(TR)^q$ we have
\[
 W_q^*=(R^*T)^qR^*.
\]
Subtracting these products by replacing one of the $q+1$ resolvents at a time gives the finite algebraic identity
\[
 W_q-W_q^*
 =\sum_{j=0}^q(R^*T)^j(R-R^*)(TR)^{q-j}.
\]
Because $R-R^*=-2\mathrm iS$, this proves
\begin{equation}
 -\operatorname{Im}W_q
 =\sum_{j=0}^q(R^*T)^jS(TR)^{q-j}.
 \label{W:eq:resolvent-imaginary-word}
\end{equation}
Every summand has a positive operator in the middle, but the products as a whole need not be positive. The bound used in a quadratic form follows from Cauchy--Schwarz in the
$S$-inner product. For arbitrary $u,v$,
\[
 |u^*(R^*T)^jS(TR)^{q-j}v|
 \le\|S^{1/2}(TR)^ju\|\,
     \|S^{1/2}(TR)^{q-j}v\|.
\]
In particular, setting $u=v$ and using the arithmetic--geometric mean inequality controls the sum by the positive quadratic forms $((TR)^ju)^*S((TR)^ju)$ for $0\le j\le q$.

The same argument covers mixed derivatives after their labelled permutation sum has been taken. For self-adjoint $T_1,\ldots,T_q$ define
\[
 \mathcal W_q=\sum_{\sigma\in\mathfrak S_q}
     RT_{\sigma(1)}R\cdots T_{\sigma(q)}R.
\]
Taking the adjoint reverses the order of the $T$'s. Reversal permutes $\mathfrak S_q$, so
\[
 \mathcal W_q^*=\sum_{\sigma\in\mathfrak S_q}
   R^*T_{\sigma(1)}R^*\cdots T_{\sigma(q)}R^*.
\]
For each $\sigma$ and $0\le j\le q$, let
\[
 L_{\sigma,j}=R^*T_{\sigma(1)}R^*\cdots R^*T_{\sigma(j)},
 \qquad
 U_{\sigma,j}=T_{\sigma(j+1)}R\cdots T_{\sigma(q)}R,
\]
with the empty products interpreted as the identity. Replacing successive resolvents gives exactly
\[
 -\operatorname{Im}\mathcal W_q
   =\sum_{\sigma\in\mathfrak S_q}\sum_{j=0}^q
       L_{\sigma,j}S U_{\sigma,j},
\]
and its matrix elements satisfy
\[
 |u^*L_{\sigma,j}S U_{\sigma,j}v|
 \le\|S^{1/2}L_{\sigma,j}^*u\|
       \|S^{1/2}U_{\sigma,j}v\|.
\]
There are $q!(q+1)$ summands. Their sum is bounded by the sum of the displayed products of $S^{1/2}$-weighted norms.

Let $f$ be holomorphic on a neighbourhood of the spectra of the matrices on the finite set of entry variations under consideration, and choose a positively oriented contour
$\mathcal C$ enclosing those spectra inside this neighbourhood. For each fixed matrix,
\[
 f(H)=\frac1{2\pi\mathrm i}\int_{\mathcal C}
                 f(z)(zI-H)^{-1}\,\mathrm dz.
\]
To justify differentiation, fix $H$ and put $\delta_H=\operatorname{dist}(\mathcal C,\operatorname{spec}H)>0$. For a self-adjoint $H'$ with $\|H'-H\|<\delta_H/2$, the eigenvalue
bound gives $\operatorname{dist}(\mathcal C,\operatorname{spec}H')\ge\delta_H/2$. The permutation formula therefore yields
\[
 \sup_{z\in\mathcal C}\|D_{H'}^qR_{H'}(z)[V_1,\ldots,V_q]\|
 \le q!(2/\delta_H)^{q+1}\prod_{j=1}^q\|V_j\|.
\]
Multiplication by $|f(z)|$ gives an integrable bound on the fixed finite-length contour, uniform in this neighbourhood. The inverse identity gives the first difference quotient and
its limit there; applying the same bound successively to the differentiated products justifies every derivative under the contour integral. The argument is at fixed $N$; any
uniform estimate as $N$ grows must use the displayed dependence on $\delta_H$. The resulting formula is
\begin{equation}
 D_H^q f(H)[V_1,\ldots,V_q]
 =\frac1{2\pi\mathrm i}\int_{\mathcal C}f(z)
   \sum_{\sigma\in\mathfrak S_q}
       R V_{\sigma(1)}R\cdots V_{\sigma(q)}R\,\mathrm dz.
 \label{W:eq:resolvent-functional-calculus}
\end{equation}
The $V_j$ are fixed self-adjoint matrices, and the same $H$ occurs in every resolvent factor.

For a trace, integration by parts gives
\begin{equation}
 D_H^q\operatorname{Tr}f(H)[V_1,\ldots,V_q]
 =\frac1{2\pi\mathrm i\,q}\int_{\mathcal C}f'(z)
      \sum_{\sigma\in\mathfrak S_q}
       \operatorname{Tr}(V_{\sigma(1)}R\cdots
                              V_{\sigma(q)}R)\,\mathrm dz.
 \label{W:eq:resolvent-trace-ibp}
\end{equation}
To prove it, set $B_\sigma(z)=\operatorname{Tr}(V_{\sigma(1)}R\cdots V_{\sigma(q)}R)$ and $B(z)=\sum_\sigma B_\sigma(z)$. Differentiation of the $q$ resolvent factors in $B_\sigma$
uses $\partial_zR=-R^2$. For the factor after $V_{\sigma(j)}$, the differentiated trace is
\[
 -\operatorname{Tr}
  (V_{\sigma(1)}R\cdots V_{\sigma(j)}R^2
        V_{\sigma(j+1)}R\cdots V_{\sigma(q)}R).
\]
Cyclicity moves the resolvent portion before this squared factor to the end of the trace. The result is
\[
 -\operatorname{Tr}
  (RV_{\sigma(j+1)}R\cdots V_{\sigma(q)}R
    V_{\sigma(1)}R\cdots V_{\sigma(j)}R),
\]
where the relevant list is understood cyclically when $j=q$. For each fixed $j$, cyclic permutation of the labels of the matrices $V_j$ is a bijection of $\mathfrak S_q$. Summing
first over $\sigma$ and then over $j$ therefore yields the exact identity
\[
 B'(z)=-q\sum_{\sigma\in\mathfrak S_q}
             \operatorname{Tr}(RV_{\sigma(1)}R\cdots
                                            V_{\sigma(q)}R).
\]
The function $f(z)B(z)$ is meromorphic inside $\mathcal C$ and holomorphic near it. Its derivative has zero integral round this contour (equivalently, the residue of a derivative
at every pole is zero). Hence
\[
 0=\int_{\mathcal C}(f(z)B(z))'\,\mathrm dz
  =\int_{\mathcal C}\bigl(f'(z)B(z)+f(z)B'(z)\bigr)
                                  \,\mathrm dz.
\]
Substitute the identity for $B'$ into the trace of \eqref{W:eq:resolvent-functional-calculus}; this proves \eqref{W:eq:resolvent-trace-ibp} with exactly the factor $1/q$. For
$V_1=\cdots=V_q=V$, all $q!$ permutations agree, so division by $q!$ yields
\begin{equation}
 \frac1{q!}D_H^q\operatorname{Tr}f(H)[V^q]
 =\frac1{2\pi\mathrm i\,q}\int_{\mathcal C}
       f'(z)\operatorname{Tr}((VR)^q)\,\mathrm dz.
 \label{W:eq:resolvent-repeated-trace-ibp}
\end{equation}
For the normalized trace both sides are divided by $N$. If a scalar coefficient $a(H)$ is independent of $z$ it may be pulled through the contour integral; its derivatives in $H$
are still governed by the product rule when $H$ is varied. If $a(H,z)$ is also holomorphic in $z$, the integration by parts acts on the product and replaces $f'$ by
$\partial_z(af)=a\,f'+(\partial_z a)f$.

For the modified function $p_n\circ h_d$, smooth functional calculus gives the derivatives on the full matrix space. On a set where $h_d$ equals the identity on an open interval
containing the spectra of every matrix in an entry segment, these derivatives equal those of $p_n$. We then use the polynomial Cauchy formula. Products containing $\Theta_a$ use
the contour avoiding its poles constructed in Section~\ref{W:sec:contour-integrals-for-moment-increments}.

The projection entering the spectral localization statistic also depends on $H$. Fix the deterministic starting vector and all matrix variations. On the open set where the
$N\times r$ polynomial matrix $F(H)$ has full column rank, set
\[
 G=F^*F,\qquad A=G^{-1},\qquad P=FAF^*.
\]
For a finite set of distinct derivative labels $S$, the complete product and inverse rules are
\begin{align}
 D_SG&=\sum_{U\sqcup V=S}(D_UF)^*D_VF,\notag\\
 D_SA&=\sum_{\ell=1}^{|S|}(-1)^\ell
    \sum_{(S_1,\ldots,S_\ell)\in\operatorname{OP}(S)}
       A(D_{S_1}G)A\cdots(D_{S_\ell}G)A,
       \qquad S\ne\varnothing,
       \label{W:eq:resolvent-gram-labels}\\
 D_SP&=\sum_{U\sqcup V\sqcup W=S}
                 (D_UF)(D_VA)(D_WF)^*.
       \label{W:eq:resolvent-projection-labels}
\end{align}
In the inverse formula every ordered block $S_j$ is nonempty; empty sets are allowed in the other two formulas. The inverse identity is \eqref{W:eq:gram-inverse-labels}; the other
two follow from the product rule in Section~\ref{W:sec:notation-and-scales}. At orders two and three,
\begin{align*}
 D_{12}A={}&AG_1AG_2A+AG_2AG_1A-AG_{12}A,\\
 D_{123}A={}&-AG_{123}A
 +\sum_{\substack{U\sqcup V=\{1,2,3\}\\U,V\ne\varnothing}}
                      AG_UAG_VA\\
 &-\sum_{\sigma\in\mathfrak S_3}
     AG_{\sigma(1)}AG_{\sigma(2)}AG_{\sigma(3)}A,
\end{align*}
where $G_T=D_TG$.

There is an additional cancellation at the first order. For a single matrix variation $V$, the three terms in the product rule yield
\begin{align*}
 DP[V]
  &=(DF[V])AF^*+F(DA[V])F^*+FA(DF[V])^*,\\
 DA[V]&=-A\bigl((DF[V])^*F+F^*DF[V]\bigr)A.
\end{align*}
Substituting the second line into the first and collecting the two terms on the left and right gives
\begin{equation}
 DP[V]=(I-P)(DF[V])AF^*+FA(DF[V])^*(I-P).
 \label{W:eq:resolvent-projection-first}
\end{equation}
Put $U=FG^{-1/2}$ and $Z=(I-P)(DF[V])G^{-1/2}$. Then $U^*U=I$, $P=UU^*$, $U^*Z=0$, and
$DP[V]=ZU^*+UZ^*$.
These are off-diagonal blocks relative to $\operatorname{Ran}P\oplus\operatorname{Ran}(I-P)$, so
\[
 \|DP[V]\|_{\rm HS}^2=2\|Z\|_{\rm HS}^2,
 \qquad\|DP[V]\|=\|Z\|.
\]
If $G\ge\kappa I$, this yields the dimension-independent bounds for the first derivative
\[
 \|DP[V]\|_{\rm HS}\le\sqrt{2/\kappa}\,\|DF[V]\|_{\rm HS},
 \qquad\|DP[V]\|\le\kappa^{-1/2}\|DF[V]\|.
\]
The first-order cancellation does not eliminate the diagonal blocks at order two. Differentiating $P^2=P$ in two labels gives
\[
 P_{12}=P_{12}P+PP_{12}+P_1P_2+P_2P_1,
\]
where $P_j=D_jP$ and $P_{12}=D_{12}P$. Multiplication by $P$ and by $Q=I-P$ respectively gives
\[
 PP_{12}P=-P(P_1P_2+P_2P_1)P,
 \quad QP_{12}Q=Q(P_1P_2+P_2P_1)Q.
\]
These blocks are generally nonzero. Higher derivatives must therefore use the full labelled product and inverse formulas, not a repeated application of the norm bound for the first
derivative. The column derivatives $D_SF$ are computed by the contour functional calculus on the spectral agreement set and by the smooth scalar modification away from it. The
algebra above holds wherever $F$ has full column rank; the regularizations extend the final test functions to all Hermitian matrices.

\subsubsection{Spectral confinement and scalar regularization}
\label{W:sec:derivative-bounds-for-the-spectral-cutoff}
\label{W:sec:scalar-regularization-and-divided-differences}

The integral representations also use real spectral parameters outside the spectrum. To pass from nonreal resolvent estimates to these points, the relevant control is their
distance from the spectrum. The next lemma compares the resolvent along a short vertical segment using that distance. Spectral confinement will then place the original matrix
inside the interval where the scalar modification is the identity, allowing the modified moments and their derivatives to be identified with the original polynomial expressions.

\begin{lemma}[Resolvents away from the spectrum]
\label{W:lem:vertical-resolvent-comparison}
Let $H=H^*$, $\operatorname{dist}(E,\operatorname{spec}H)\ge\delta>0$, and $0<h\le A\delta$. For $|y|\le h$ and every vector $v$,
\begin{equation}
 \|R_H(E+\mathrm i y)v\|^2
 \le(1+A^2)\|R_H(E+\mathrm i h)v\|^2.
 \label{W:eq:vertical-energy-comparison}
\end{equation}
If $z,z'$ lie on this vertical segment, then
\begin{equation}
 |u^*(R_H(z)-R_H(z'))v|
 \le(1+A^2)|z-z'|\|R_H(E+\mathrm i h)u\|
                            \|R_H(E+\mathrm i h)v\|.
 \label{W:eq:vertical-matrix-element-comparison}
\end{equation}
On any pair of spectral points whose distance from the spectrum is at least $\delta$, the entry and quadratic form functions of $R_H$ are $\delta^{-2}$-Lipschitz for unit exterior
vectors, the squared column norms are $2\delta^{-3}$-Lipschitz, and $v^*|R_H(z)|v$ is $\delta^{-2}$-Lipschitz for $\|v\|=1$.
\end{lemma}
\begin{proof}
For the positive spectral measure $\mu_v$ of $H$ at $v$, write $r=E-\lambda$. Since $|r|\ge\delta$,
\[
 \frac1{r^2+y^2}\le\frac{1+A^2}{r^2+h^2},\qquad
 \|R_H(E+\mathrm i y)v\|^2
       =\int\frac{\mu_v(\mathrm d\lambda)}{r^2+y^2}.
\]
This proves \eqref{W:eq:vertical-energy-comparison}. The resolvent identity $R_H(z)-R_H(z')=(z'-z)R_H(z)R_H(z')$ and Cauchy--Schwarz give
\[
 |u^*(R_H(z)-R_H(z'))v|
 \le |z-z'|\|R_H(\bar z)u\|\|R_H(z')v\|.
\]
The energies at conjugate spectral points agree, so the first bound proves \eqref{W:eq:vertical-matrix-element-comparison}. For the Lipschitz statements the same identity gives
$\|R_H(z)-R_H(z')\|\le\delta^{-2}|z-z'|$. Consequently
\[
 \left|\|R_H(z)v\|^2-\|R_H(z')v\|^2\right|
 \le(\|R_H(z)v\|+\|R_H(z')v\|)
                   \|(R_H(z)-R_H(z'))v\|
 \le2\delta^{-3}|z-z'|.
\]
Finally, the reverse triangle inequality implies
\[
 \left|\frac1{|z-\lambda|}-\frac1{|z'-\lambda|}\right|
 \le\frac{|z-z'|}{|z-\lambda|\,|z'-\lambda|}
 \le\delta^{-2}|z-z'|.
\]
Integrate this scalar inequality against $\mu_v$. Normalized traces satisfy the same bounds; unnormalized traces have an additional factor $N$.
\end{proof}

The exponential factor in the polynomial bounds dictates a narrower spectral interval. With this choice it is subpolynomial in $N$. The scalar estimates below hold at every
Hermitian matrix and will be used in Sections~\ref{W:sec:green-function-comparison} and \ref{W:sec:localization-of-the-low-energy-spectral-subspace}.

Let ${m_N}=N^{1/3}$ and let $\rho_*$ be the maximal Gram matrix dimension, with $c_*{m_N}\le \rho_*\le C_*{m_N}$. Fix $A_0>0$ and $1<\gamma<2$, and set
\begin{equation}
 d_N=A_0\frac{(\log N)^\gamma}{\rho_*^2},\qquad
 \Omega_N^{(\gamma)}=\{\|H_N\|\le2+d_N/2\}.
 \label{W:eq:gram-width}
\end{equation}
All constants below may depend on these fixed parameters.

\Needspace{6\baselineskip}

\begin{proposition}
\label{W:prop:gram-window-2}
Let $H_N$ satisfy Assumption~\ref{W:ass:wigner-lee-yin}, and let $d_N,\Omega_N^{(\gamma)}$ be defined in \eqref{W:eq:gram-width}. Then
\begin{equation}
 \mathbb P((\Omega_N^{(\gamma)})^c)
 \le C\{N^{-c}+\Theta(c_0\sqrt N)+\Theta_{\mathrm d}(c_0\sqrt N)\}=o(1).
 \label{W:eq:gram-window-prob}
\end{equation}
For the Gaussian ensemble the upper bound is $C\exp[-c(\log N)^{3\gamma/2}]$. The equality of the original and modified observables under the cutoffs, and the comparison of bounded
functions in Section~\ref{W:sec:resolvents-on-exceptional-support} remain valid with this width. Moreover
\begin{equation}
 \exp(C\rho_*\sqrt{d_N})
 =\exp\{C\sqrt{A_0}(\log N)^{\gamma/2}\}=N^{o(1)}
 \label{W:eq:gram-subpoly}
\end{equation}
for every fixed $C$.
\end{proposition}
\begin{proof}
Write $\eta_{\rm e}=N^{-2/3}$ and $d=d_N=A_0\rho_*^{-2}(\log N)^\gamma$ with $1<\gamma<2$. The constants below may depend on $A_0,\gamma,c_*,C_*$ and the fixed symmetry class, but
not on $N$ or the deterministic starting vector. Since $c_*m_N\le\rho_*\le C_*m_N$,
\[
 c(\log N)^\gamma\le d/\eta_{\rm e}
       \le C(\log N)^\gamma,
 \qquad \eta_{\rm e}/d\longrightarrow0.
\]
In particular $d<1/10$ for all sufficiently large $N$.

Use the decomposition $H=W+V+E$ established in Section~\ref{W:sec:removal-of-spectral-regularization}. The normalized coarsely truncated matrix $W$ has the independent centered
entry blocks and the moment bounds required by Lemma~\ref{W:lem:resolvent-core}: off-diagonal variance $N^{-1}$, support $CN^{-\phi}$, a uniform bound on the standardized third
absolute moments, a bound $C\log(eN)$ on the standardized fourth moments, and diagonal variances at most $C/N$. Its complex entries also have zero pseudovariance. The lemma applies
with the fixed parameter $b=1/8$. More explicitly, its comparison estimate for the trace has error
\[
 \epsilon_N\le C\left(N^{-c_1}+
      e^{-d/(8\eta_{\rm e})}+e^{-1/\eta_{\rm e}}\right).
\]
The two exponential contributions remain smaller than any negative power of $N$, even after multiplication by a fixed power of $N$: for every $D>0$,
\[
 \log\bigl(N^{D+1}\log(eN)e^{-d/(16\eta_{\rm e})}\bigr)
 \le (D+1)\log N+\log\log(eN)-c(\log N)^\gamma
 \longrightarrow-\infty.
\]
Applying the moment conclusion \eqref{W:eq:truncated-core-moment-solution} with an order chosen after $D$, and then Markov's inequality, gives
\[
 \mathbb P(\|W\|>2+d/8)\le N^{-D}.
\]
The polynomial spectral window for the bounded comparison matrix is of width
\[
  C\eta_{\rm e}(\log N)^{2/3}=o(d).
\]
This verifies the width requirement in the cited moment bound. For the coordinate estimate, Lemma~\ref{W:lem:complex-weak} with $\tau=1/4$ applies at
$\eta_{\rm e}=N^{-2/3}\ge N^{-1+\tau}$. Its row estimates also hold for real entries: in the proof of Lemma~\ref{W:lem:complex-row}, conjugation then leaves each scalar unchanged,
and the bounds use centering, $\mathbb E|x_j|^2=N^{-1}$ and the same absolute moment estimates. The conditional mean of the quadratic form is $N^{-1}\operatorname{Tr}R^{(i)}$ in
either case. The subsequent resolvent identities therefore give the same coordinate estimate for real $W$.

The counting estimate is \eqref{W:eq:minor-counting}, supplied by Proposition~\ref{W:prop:minor-real-input} when $\beta=1$ and Theorem~\ref{W:thm:complex-main} when $\beta=2$. For
a fixed $D>0$, apply these estimates and the preceding bound for $\|W\|$ with failure exponent $D+1$. Let $\mathcal G_N$ be their intersection. It is measurable with respect to $W$
and, for all sufficiently large $N$,
\[
 \begin{gathered}
 \|W\|\le2+d/8,\qquad
 \sup_{|x|\le5}\max_{i,j}
 |R_W(x+\mathrm i\eta_{\rm e})_{ij}-g(x+\mathrm i\eta_{\rm e})\delta_{ij}|\le N^{-c},\\
 \sup_x|F_W(x)-F_{\rm sc}(x)|\le \mathfrak c_N/N,
 \qquad \mathfrak c_N\le N^{1/12},\\
 \mathbb P(\mathcal G_N^c)\le3N^{-D-1}\le N^{-D}.
 \end{gathered}
\]
The bounds at negative imaginary part follow by conjugation.

On this event define the two real spectral barriers $a_+=2+d/4$ and $a_-=-2-d/4$. Their distance from $\sigma(W)$ is at least $d/8$. For $\lambda\in\sigma(W)$ put
$r=a_+-\lambda\ge d/8$. The elementary inequality
\[
 \frac1r=\frac{r^2+\eta_{\rm e}^2}{r^2}\,
                      \frac r{r^2+\eta_{\rm e}^2}
 \le\left(1+\frac{64\eta_{\rm e}^2}{d^2}\right)
                      \frac r{r^2+\eta_{\rm e}^2}
\]
passes through the positive coordinate spectral measure of $W$. Thus
\[
 0\le R_W(a_+)_{ii}
 \le\left(1+\frac{64\eta_{\rm e}^2}{d^2}\right)
       \operatorname{Re}R_W(a_++\mathrm i\eta_{\rm e})_{ii}.
\]
The coarse coordinate estimate bounds the last quantity by three for sufficiently large $N$, while its prefactor tends to one; therefore $R_W(a_+)_{ii}\le4$ for every $i$. At
$a_-$, the spectral measure $\mu_i$ of $W$ at $e_i$ gives
\[
 0\le -R_W(a_-)_{ii}
 =\int\frac{\mu_i(\mathrm d\lambda)}{\lambda-a_-}
 \le\left(1+\frac{64\eta_{\rm e}^2}{d^2}\right)
       \bigl[-\operatorname{Re}R_W(a_-+\mathrm i\eta_{\rm e})_{ii}\bigr]
 \le4.
\]
Here $\lambda-a_-\ge d/8$ on the support of $\mu_i$, and the scalar inequality used at $a_+$ applies with this positive denominator. These conclusions hold simultaneously for both
barriers on the chosen high-probability event.

We estimate the off-diagonal resolvent compression at each barrier. Set $t_j=2+d/8-\lambda_j(W)\ge0$ and $n_W(t)=\#\{j:t_j\le t\}$. On $\mathcal G_N$, the eigenvalue counting
estimate gives
\[
 n_W(t)\le C(Nt^{3/2}+K_N),\quad t\ge0,
 \qquad K_N=N^{1/12}.
\]
At the lower barrier the same estimate is used for $-W$. Since $a_+-\lambda_j(W)=d/8+t_j$, integration of $2(d/8+t)^{-3}$ against this counting function gives
\begin{align*}
 \|R_W(a_+)\|_{\rm HS}^2
 &=\sum_{j=1}^N(d/8+t_j)^{-2}
   =2\int_0^\infty\frac{n_W(t)}{(d/8+t)^3}\,\mathrm dt\\
 &\le C N d^{-1/2}
         \int_0^\infty\frac{u^{3/2}}{(1+u)^3}\,\mathrm du
          +C K_Nd^{-2}\\
 &\le CN^{4/3+1/12}.
\end{align*}
The integral is finite at zero and infinity. The same bound holds for $R_W(a_-)$. The last inequality follows from $d\asymp N^{-2/3}(\log N)^\gamma$: $Nd^{-1/2}=O(N^{4/3})$ and
$K_Nd^{-2}=O(N^{4/3+1/12})$.

Let $S$ denote the support of the exceptional matrix $E$ and $\Pi_S$ the associated coordinate projection. The support construction gives, conditional on $W$,
\[
 \mathbb P(i,j\in S\mid W)\le CN^{-2+8\phi},
 \qquad i\ne j.
\]
Since $\mathcal G_N$ and $R_W$ are measurable with respect to $W$, conditional expectation gives
\begin{align*}
 &\mathbb E\left[\mathbf1_{\mathcal G_N}
   \sum_{i\ne j}\mathbf1_{\{i,j\in S\}}
                 |R_W(a_\pm)_{ij}|^2\right]\\
 &\quad=\mathbb E\left[\mathbf1_{\mathcal G_N}
   \sum_{i\ne j}\mathbb P(i,j\in S\mid W)
                       |R_W(a_\pm)_{ij}|^2\right]\\
 &\quad\le CN^{-2+8\phi}
   \mathbb E\left[\mathbf1_{\mathcal G_N}
                                \|R_W(a_\pm)\|_{\rm HS}^2\right]
 \le CN^{-2+8\phi+4/3+1/12}.
\end{align*}
For $\phi=1/100$, the exponent equals
$-2+8/100+4/3+1/12=-151/300<0$.
Markov's inequality shows that, outside an event of probability at most $CN^{-151/300}$, the matrix obtained by deleting the diagonal of $\Pi_SR_W(a_\pm)\Pi_S$ has Hilbert--Schmidt
norm at most one. Its diagonal part has norm at most four by the preceding coordinate estimate. Hence
\[
 \|R_W(a_\pm)_{SS}\|\le5
\]
at both barriers outside a set of probability at most
$\mathbb P(\mathcal G_N^c)+CN^{-151/300}$.

On $\{\|E\|\le1/20\}$, the support identity $E=\Pi_SE\Pi_S$ and the compression bound give, for
$0\le t\le1$,
\[
 \|tE_{SS}R_W(a_\pm)_{SS}\|
 \le t\,\|E\|\,\|R_W(a_\pm)_{SS}\|
 \le\frac14.
\]
The matrix determinant lemma, which is valid for these finite blocks because $a_\pm I-W$ is invertible, gives
\[
 \det(a_\pm I-W-tE)
 =\det(a_\pm I-W)
       \det(I_S-tE_{SS}R_W(a_\pm)_{SS})\ne0.
\]
Thus neither barrier is an eigenvalue of $W+tE$ for $0\le t\le1$. The eigenvalues are continuous in $t$, and at $t=0$ they all lie strictly between the barriers. No eigenvalue can
cross a barrier without meeting it, so $\|W+E\|\le2+d/4$ on the intersection of the preceding events. The remaining perturbation satisfies $\|V\|\le N^{-1+4\phi}=o(d)$ for
$\phi=1/100$; for large $N$ it contributes less than $d/4$. Consequently $\|H\|\le2+d/2$ on the same event. The union bound now gives
\[
 \mathbb P(\|H\|>2+d/2)
 \le CN^{-c}
      +C\bigl\{\Theta(c_0\sqrt N)+\Theta_{\rm d}(c_0\sqrt N)\bigr\}.
\]
Here $c>0$ can be reduced to include the fixed exponent $151/300$, the coarse coordinate and counting errors, and the normalization and support errors. The tail event
$\{\|E\|>1/20\}$ is controlled by the displayed tail functions from the truncation construction. This proves \eqref{W:eq:gram-window-prob} under the original assumptions.

For the GOE and GUE, Lemma~1.1 of \cite{erdosxu2023small} applies at $x=d/(2\eta_{\rm e})\asymp(\log N)^\gamma$, which belongs to $[1,N^{1/4}]$ for all sufficiently large $N$. It
gives $\mathbb P(\lambda_{\max}>2+d/2) \le C\exp[-c(\log N)^{3\gamma/2}]$; apply the corresponding estimate to $-H$ for the lower edge and sum the two failure probabilities. If
$\|H\|\le2+d/2$ and a perturbation $U$ satisfies $\|U\|<d/2$, then $\|H+U\|<2+d$. The identity $h_d(x)=x$ for $|x|\le2+d$ gives $h_d(H+U)=H+U$ by the spectral theorem; all
polynomial moments and functions of the Gram matrix constructed from these two matrices therefore agree with their unmodified versions on this event. Finally,
\[
 e^{C\rho_*\sqrt d}
 =\exp\bigl(C\sqrt{A_0}(\log N)^{\gamma/2}\bigr)
 =N^{C\sqrt{A_0}(\log N)^{\gamma/2-1}}=N^{o(1)}
\]
because $\gamma<2$.
\end{proof}

For the moment coordinates use the odd function $h_d$ constructed in Lemma~\ref{W:lem:gram-clipping}, with the same fixed even cutoff. That construction is valid for $0<d<1/4$; all
asymptotic uses have $d\le1/10$. Differentiating its explicit transition gives $\|h_d^{(r)}\|_\infty\le C_rd^{1-r}$ for $r\ge1$. Indeed, in the transition interval every derivative
of order $r$ contains the factor $d^{1-r}$; outside it, derivatives of the fixed cutoff are uniformly bounded. Thus the same function satisfies
\begin{equation}
 \begin{gathered}
 h_d(x)=x\quad(|x|\le2+d),\qquad |h_d(x)|\le2+2d,\\
 \operatorname{supp}h_d\subset[-4,4],\qquad
 \|h_d^{(r)}\|_\infty\le C_rd^{1-r}\quad(r\ge1).
 \end{gathered}
 \label{W:eq:gram-clip}
\end{equation}
Proposition~\ref{W:prop:gram-window-2} applies to this choice.

For the fixed polynomials $p_n$, put
\begin{equation}
 \varphi_{n,d}=p_n\circ h_d,\qquad
 s_n(H,w)=w^*\varphi_{n,d}(H)w,\qquad
 \bar s_n(H)=N^{-1}\operatorname{Tr}\varphi_{n,d}(H).
 \label{W:eq:gram-family}
\end{equation}
The function $\varphi_{n,d}-p_n(0)$ is smooth and compactly supported. Its constant part is annihilated by the trace pairing with $ww^*-N^{-1}I$ and by every matrix derivative of
positive order.

\Needspace{6\baselineskip}

\begin{lemma}
\label{W:lem:gram-reduction}
Let ${H^{[d]}}=h_d(H)$ and ${F^{[d]}}=(p_0({H^{[d]}})w,\ldots,p_{{\rho}-1}({H^{[d]}})w)$. For $0\le a,b<{\rho}$,
\begin{align}
 {G_{ab}^{[d]}}
 &=\sum_{j=0}^{\min(a,b)}s_{a+b-2j},
 \label{W:eq:gram-Gram-reduction}\\
 {S_{ab}^{[d]}}
 &=\sum_{j=0}^{\min(a,b)}
       \{s_{a+b-2j+1}+s_{a+b-2j-1}\},\qquad s_{-1}=0,
 \label{W:eq:gram-S-reduction}
\end{align}
where ${G^{[d]}}={F^{[d]}}^*{F^{[d]}}$ and ${S^{[d]}}={F^{[d]}}^*{H^{[d]}}{F^{[d]}}$. Replacing $s_n$ by $s_n-\bar s_n$ gives the corresponding centered entries. These are exact
identities of smooth functions and hence of every fixed mixed derivative. Only degrees at most $2{\rho}-1$ occur.
\end{lemma}
\begin{proof}
Write $B=h_d(H)$ and $s_n=w^*p_n(B)w$. The product identity \eqref{W:eq:moment-coordinates-Chebyshev}, applied by polynomial functional calculus, gives
\[
 \begin{aligned}
 G_{ab}^{[d]}&=w^*p_a(B)p_b(B)w
       =\sum_{j=0}^{\min(a,b)}w^*p_{a+b-2j}(B)w,\\
 S_{ab}^{[d]}&=w^*p_a(B)Bp_b(B)w
       =\sum_{j=0}^{\min(a,b)}w^*Bp_{a+b-2j}(B)w.
 \end{aligned}
\]
For $n\ge0$, \eqref{W:eq:chebyshev-three-term} reads $Bp_n(B)=p_{n+1}(B)+p_{n-1}(B)$, with $p_{-1}=0$. Substitution proves \eqref{W:eq:gram-Gram-reduction} and
\eqref{W:eq:gram-S-reduction}. The largest indices are $a+b\le2\rho-2$ and $a+b+1\le2\rho-1$, respectively.

Define $\bar s_n=N^{-1}\operatorname{Tr}p_n(B)$ and set $\bar s_{-1}=0$. Applying $N^{-1}\operatorname{Tr}$ to the same polynomial identities and subtracting gives
\[
 \begin{aligned}
 G_{ab}^{[d]}-N^{-1}\operatorname{Tr}(p_a(B)p_b(B))
   &=\sum_{j=0}^{\min(a,b)}(s_{a+b-2j}-\bar s_{a+b-2j}),\\
 S_{ab}^{[d]}-N^{-1}\operatorname{Tr}(p_a(B)Bp_b(B))
   &=\sum_{j=0}^{\min(a,b)}\sum_{\epsilon=\pm1}
                 (s_{a+b-2j+\epsilon}-\bar s_{a+b-2j+\epsilon}).
 \end{aligned}
\]
All sums are finite. If $D_{H,S}$ is a mixed derivative with respect to fixed Hermitian matrix variations, smoothness of $H\mapsto p_n(h_d(H))$ therefore gives
\[
 D_{H,S}G_{ab}^{[d]}=\sum_{j=0}^{\min(a,b)}D_{H,S}s_{a+b-2j},\qquad
 D_{H,S}S_{ab}^{[d]}=\sum_{j=0}^{\min(a,b)}\sum_{\epsilon=\pm1}
                         D_{H,S}s_{a+b-2j+\epsilon}.
\]
The centered identities have the same derivatives with $s_n-\bar s_n$ in place of $s_n$.
\end{proof}

Write $f^{[q]}(x_0,\ldots,x_q)$ for the Hermite divided difference, with $f^{[q]}(x,\ldots,x)=f^{(q)}(x)/q!$. At repeated nodes we use the continuous extension provided by the
simplex integral in Lemma~\ref{W:lem:smooth-functional-calculus}.

\Needspace{6\baselineskip}

\begin{theorem}
\label{W:thm:gram-composition}
For a smooth $h$ and a polynomial $p$, for $q\ge1$,
\begin{equation}
 (p\circ h)^{[q]}(x_0,\ldots,x_q)
 =\sum_{k=1}^q\ \sum_{0=i_0<i_1<\cdots<i_k=q}
 p^{[k]}(h(x_{i_0}),\ldots,h(x_{i_k}))
 \prod_{\ell=1}^k
 h^{[i_\ell-i_{\ell-1}]}(x_{i_{\ell-1}},\ldots,x_{i_\ell}).
 \label{W:eq:gram-composition}
\end{equation}
There are exactly $2^{q-1}$ terms. For smooth $f,g$,
\begin{equation}
 (fg)^{[q]}(x_0,\ldots,x_q)
 =\sum_{j=0}^q
 f^{[j]}(x_0,\ldots,x_j)
 g^{[q-j]}(x_j,\ldots,x_q).
 \label{W:eq:gram-product}
\end{equation}
The first identity applies to $p_n\circ h_d$, and the second to its product with $\Theta_a$.
\end{theorem}
\begin{proof}
We prove the formulas directly; they are the divided-difference composition and product formulas recorded in \cite{floater2007chain}. First suppose that $x_0,\ldots,x_q$ are
distinct. Let $X$ be the $(q+1)\times(q+1)$ upper bidiagonal matrix with diagonal $x_0,\ldots,x_q$ and ones immediately above the diagonal. For a scalar function $f$ on these
nodes, define $f(X)$ by the unique interpolation polynomial of degree at most $q$. The eigenvalues of $X$ are distinct, so evaluation of these polynomials respects products.

For a nonnegative integer $m$ and $i\le j$, multiplication of the bidiagonal factors gives
\[
 (X^m)_{ij}
 =\sum_{a_i+\cdots+a_j=m-j+i}x_i^{a_i}\cdots x_j^{a_j}
 =(x^m)^{[j-i]}(x_i,\ldots,x_j).
\]
The sum is zero when its required degree is negative. The last identity is the monomial formula proved in Lemma~\ref{W:lem:smooth-functional-calculus}; it also follows by induction
from $(x^m-y^m)/(x-y)=\sum_{a=0}^{m-1}x^ay^{m-1-a}$. By linearity and interpolation,
\[
 f(X)_{ij}=f^{[j-i]}(x_i,\ldots,x_j)\quad(i\le j),
 \qquad f(X)_{ij}=0\quad(i>j).
\]
Taking the $(0,q)$ entry of $(fg)(X)=f(X)g(X)$ now gives precisely \eqref{W:eq:gram-product}.

For the composition, set $Y=h(X)$. For $p(x)=x^m$, an individual term in $(Y^m)_{0q}$ is a product along a nondecreasing sequence of indices. Delete its repeated indices to obtain
a unique strictly increasing sequence $0=i_0<i_1<\cdots<i_k=q$. The remaining $m-k$ factors are diagonal factors, distributed among those $k+1$ indices. Consequently
\[
 (Y^m)_{0q}
 =\sum_{k=1}^q\ \sum_{0=i_0<\cdots<i_k=q}
 \left(\sum_{a_0+\cdots+a_k=m-k}\prod_{l=0}^kY_{i_li_l}^{a_l}\right)
                    \prod_{l=1}^kY_{i_{l-1}i_l}.
\]
The coefficient in parentheses is $(x^m)^{[k]}(Y_{i_0i_0},\ldots,Y_{i_ki_k})$ by the same monomial identity, which holds also when some diagonal values coincide. Linearity
therefore replaces it by $p^{[k]}(Y_{i_0i_0},\ldots,Y_{i_ki_k})$ for an arbitrary polynomial $p$. Since
\[
 Y_{ii}=h(x_i),\qquad
 Y_{ij}=h^{[j-i]}(x_i,\ldots,x_j),\qquad
 p(Y)=(p\circ h)(X),
\]
its $(0,q)$ entry is \eqref{W:eq:gram-composition}. A strictly increasing sequence is determined by the subset of $\{1,\ldots,q-1\}$ which it visits. There are $2^{q-1}$ such
subsets, proving the stated count.

Finally, the simplex integral for divided differences in Lemma~\ref{W:lem:smooth-functional-calculus} makes each side continuous in the nodes. Approximate any list with repeated
nodes by lists of distinct real nodes. The functions and their derivatives through order $q$ are bounded on a compact interval containing these lists and their images under $h$, so
every finite product above converges. This proves both formulas at repeated nodes, including the normalization $f^{[q]}(x,\ldots,x)=f^{(q)}(x)/q!$.
\end{proof}

For example the second-order kernel is
\[
 \begin{split}
 (p\circ h)^{[2]}(x_0,x_1,x_2)
 ={}&p^{[1]}(h(x_0),h(x_2))h^{[2]}(x_0,x_1,x_2)\\
 &+p^{[2]}(h(x_0),h(x_1),h(x_2))
      h^{[1]}(x_0,x_1)h^{[1]}(x_1,x_2).
 \end{split}
\]
The first term contains the second divided difference of the spectral modification and contributes on its transition region. At orders three, four, and five there are respectively
four, eight, and sixteen terms indexed by strictly increasing sequences, regardless of the number of terms in a matrix-valued chain rule.

If every $\lambda_i$ in a spectral product lies in $[-2-d,2+d]$, its $\varphi_{n,d}^{[q]}$ equals $p_n^{[q]}$. For repeated nodes the required scalar derivatives agree there as
well. Thus the divided difference kernels agree whenever all their nodes belong to this interval.

\Needspace{6\baselineskip}

\begin{proposition}
\label{W:prop:gram-bounds}
Let $d=d_N$ be \eqref{W:eq:gram-width}. Uniformly for $0\le n\le2\rho_*+1$ and each fixed $q\ge1$,
\begin{equation}
 \|\varphi_{n,d}^{(q)}\|_\infty
 \le C_q\rho_*^{2q+1}\exp(C\rho_*\sqrt d)
 =\rho_*^{2q+1}N^{o(1)}.
 \label{W:eq:gram-scalar-bound}
\end{equation}
For $q=0$ the bound is $\rho_*N^{o(1)}$. If $\|W_j\|_{\mathrm{HS}}\le C N^{-1/2}$, then
\begin{align}
 \|D_{H}^q\varphi_{n,d}(H)[W_1,\ldots,W_q]\|_{\mathrm{HS}}
 &\le C_qN^{-q/2}\rho_*^{2q+1}N^{o(1)},
 \label{W:eq:gram-matrix-bound}\\
 |D_{H}^q\mathscr S_{\varphi_{n,d}}(H,w)[W_1,\ldots,W_q]|
 &\le C_qN^{-q/2}\rho_*^{2q+1}N^{o(1)}.
 \label{W:eq:centered-moment-derivative-bound}
\end{align}
The constants are uniform over Hermitian matrices $H$.
\end{proposition}
\begin{proof}
The polynomial derivative identity of Section~\ref{W:sec:spectral-regularization-and-matrix-derivatives} in scalar directions gives, for $|y|\le2+2d$,
\[
 |p_n^{(k)}(y)|\le
 k!\binom{n+k+1}{2k+1}\exp(n\sqrt{2d})
 \le C_k(n+1)^{2k+1}\exp(n\sqrt{2d});
\]
the derivative is zero if $k>n$. The scalar partition chain rule and \eqref{W:eq:gram-clip} therefore give
\begin{equation}
 \|\varphi_{n,d}^{(q)}\|_\infty
 \le C_q\exp(n\sqrt{2d})
     \sum_{k=1}^{\min(q,n)}(n+1)^{2k+1}d^{k-q}.
 \label{W:eq:gram-sharp-scalar}
\end{equation}
Since $n+1\le 2\rho_*+2\le4\rho_*$,
\[
 (n+1)^{2k+1}d^{k-q}
 \le C_q\rho_*^{2q+1}(\rho_*^2d)^{-(q-k)}
 \le C_q\rho_*^{2q+1},
\]
using $\rho_*^2d=A_0(\log N)^\gamma\ge1$ for large $N$. For $n=0$, every derivative of positive order is zero. This proves \eqref{W:eq:gram-scalar-bound}; the zero order bound
follows from the sine or hyperbolic formula. Subpolynomial growth is \eqref{W:eq:gram-subpoly}.

Apply Lemma~\ref{W:lem:smooth-functional-calculus} to $f=\varphi_{n,d}$ and use $\|W_j\|_{\mathrm{HS}}\le CN^{-1/2}$. Its Hilbert--Schmidt and centered trace bounds give
\eqref{W:eq:gram-matrix-bound} and \eqref{W:eq:centered-moment-derivative-bound}, respectively.

\end{proof}

The product formula \eqref{W:eq:moment-coordinates-Chebyshev} writes each entry of the Gram matrix as
\[
 G_{ab}^{[d]}=\sum_{l=0}^{\min(a,b)}s_{a+b-2l},\qquad 0\le a,b<\rho_*.
\]
There are at most $\rho_*$ summands. The matrix bound \eqref{W:eq:gram-matrix-bound}, paired with $ww^*$ in Hilbert--Schmidt norm, bounds each derivative of
$s_n=w^*\varphi_{n,d}(H)w$. Pairing with $ww^*-N^{-1}I$ gives the same bound for $s_n-\bar s_n$. Differentiating the finite sum thus yields
\[
 |D_H^qG_{ab}^{[d]}[W_1,\ldots,W_q]|
 \le C_qN^{-q/2}\rho_*^{2q+2}N^{o(1)}.
\]
This estimates one entry of the Gram matrix. Summing its entries or passing to a matrix norm incurs the corresponding additional dimension factors.

The first variation of the integrated potential also depends on the coefficients of the orthonormal polynomials. We make that dependence explicit. Fix $s$ with
$\kappa I\le G_\rho(s)\le\kappa^{-1}I$ and $K\le\rho-2$, and write $\mathscr P_K^\sigma(s,x)=\sum_{n=0}^{2K}c_n(s)p_n(x)$ as in \eqref{W:eq:potential-derivative-kernel}. If
$f=\sum_{i<\rho}f_ip_i$ and $g=\sum_{j<\rho}g_jp_j$, the Chebyshev product identity gives
\[
 \sum_n|[p_n](fg)|
 \le\sum_{i,j<\rho}|f_i||g_j|(\min(i,j)+1)
 \le\rho\left(\sum_i|f_i|\right)\left(\sum_j|g_j|\right).
\]
The coefficient sums of $q_j$ and $q_K+\sigma q_{K-1}$ are at most $Q_*$ and $2Q_*$, respectively. Moreover, $b_K/2+\sum_{j<K}|b_j-b_{j+1}|+b_1\le C_\kappa B$. Applying the
preceding inequality to the displayed squares in \eqref{W:eq:potential-derivative-kernel} proves
\[
 \sum_{n=0}^{2K}|c_n(s)|\le C_\kappa B\rho Q_*^2,
 \qquad
 |v(\mathscr P_K^\sigma)|
 \le C_\kappa B\rho Q_*^2\max_{1\le n\le2K}|v(p_n)|\quad(v(1)=0).
\]
In particular, when $v(f)=D_H^q[w^*f(h_d(H))w][W_1,\ldots,W_q]$, the coefficients $c_n(s)$ are held fixed and
\[
 |v(\mathscr P_K^\sigma)|
 \le C_{q,\kappa}BQ_*^2N^{-q/2}\rho_*^{2q+2}N^{o(1)}.
\]
The integrated potential has the further factor $m_N$. Differentiating $s=s(H)$ as well also differentiates the $c_n(s)$; those contributions are governed by
Lemma~\ref{W:lem:tightness-directions}. On the support of the coefficient cutoff defined in \eqref{W:eq:comparison-cutoff}, $BQ_*^2\le8\mathfrak r_N^3=N^{o(1)}$. This is the
condition under which the last display has only one additional power of $\rho_*$ relative to the estimate for scalar moments.

For the scalar centered moments themselves, there are at most two real coordinates for each of $N(N-1)/2$ off-diagonal entries. At fixed derivative order $q$ the number of
coordinate choices within an entry is bounded by $2^q$. Thus taking absolute values and summing the individual bounds gives
\begin{equation}
 C_qN^{2-q/2}\rho_*^{2q+1}N^{o(1)}
 =N^{7/3+q/6+o(1)},
 \label{W:eq:sum-of-coordinate-derivative-bounds}
\end{equation}
where $\rho_*\asymp N^{1/3}$. The exponents are $17/6$, $3$, and $19/6$ for $q=3,4,5$. For an entry of the Gram matrix the additional factor $\rho_*$ raises these exponents by
$1/3$. For $m_Nv(\mathscr P_K^\sigma)$ the increase is $2/3$ on the support just specified; outside it the factor $BQ_*^2$ must also be included. These bounds grow with $N$. The
later comparison therefore uses the signed cubic sum and the integrated estimates before taking its final absolute value.

\subsubsection{Smooth functions of the moment coordinates}
\label{W:sec:the-moment-vector-and-its-smooth-extensions}
\label{W:sec:extensions-of-localization-and-modulus-statistics}

The spectral modification has made the basic matrix functions smooth everywhere. We now need a common set of variables for the observables built from them: the cutoff in the
eigenvalues of the Gram matrix, characteristic functions of the integrated potential, functions controlling its modulus of continuity, and the localization statistic. Moments of
the spectral measure at the starting vector, together with normalized and unnormalized traces, provide such a finite real vector. This separates the derivatives of each observable
from the changes of its moments under entry replacement.

The cutoff parameter satisfies $0<\kappa<1/4$, the scale is ${m_N}=N^{1/3}$, and the parameters of the functions under consideration are fixed. For the integrated potential or its
modulus of continuity on $[0,T]$, the required indices are

\[
 K_T=\lceil {m_N}T\rceil,\qquad r_T=K_T+2.
\]
The spectral localization statistic uses a projection of rank $r$ and an enlarged Gram matrix of dimension $\rho\ge r$. Let $\rho_*$ be the maximum of these Gram matrix dimensions,
increasing it if necessary to cover any additional leading principal block. We work with $N\ge N_0$ such that $2\le\rho_*\le N$, which is possible for the fixed windows, and
$c_*{m_N}\le \rho_*\le C_*{m_N}$. Each function uses its own Gram matrix dimension; $\rho_*$ only determines the list of available moments. Set
\begin{equation}
 M=2\rho_*-1,\qquad
 d=A_0\rho_*^{-2}(\log(N+2))^\gamma,\qquad 1<\gamma<2.
 \label{W:eq:moment-coordinates-scales}
\end{equation}
Take the odd real $C_c^\infty$ modification $h_d$ from Section~\ref{W:sec:scalar-regularization-and-divided-differences}, so that
\[
 h_d(x)=x\quad (|x|\le2+d),\qquad |h_d(x)|\le2+2d.
\]
Use the notation $H^{[d]}=h_d(H)$ and the fixed polynomial basis. We work at sufficiently large $N$ that $d<1/4$. All derivatives are real derivatives, including for Hermitian
matrices. The deterministic starting vector $w$ is fixed when differentiating, and $\|w\|=1$.

Let $\mathcal I$ be a fixed finite set of pairs $a=(\sigma,E)$, $\sigma\in\{+1,-1\}$. Define
\begin{equation}
 \Theta_a(H)=\theta_E\bigl({m_N}^2(2I-\sigma H)\bigr),\qquad
 \theta_E(x)=(1+\mathrm e^{x-E})^{-2}.
 \label{W:eq:theta-matrix-function}
\end{equation}
Including both signs for each $E$ makes the family of coordinates closed under $H\mapsto-H$; Proposition~\ref{W:prop:moment-coordinates-reflection} specifies the resulting signed
permutation. The function $\theta_E$ satisfies $\theta_E(x)\ge1/4$ for $x\le E$ and $\theta_E(x)\le e^{-2(x-E)}$ for $x\ge E$.

\Needspace{6\baselineskip}

For $1\le n\le M$ put
\begin{align}
 u_n(H)&=N^{-1}\operatorname{Tr}p_n({H^{[d]}}),&
 x_n(H,w)&={m_N}\operatorname{Tr}\bigl[(ww^*-N^{-1}I)p_n({H^{[d]}})\bigr].
 \label{W:eq:polynomial-moment-coordinates}
\end{align}
For $a\in\mathcal I$ put
\begin{align}
 t_a(H)&=\operatorname{Tr}\Theta_a(H),\nonumber\\
 v_{a,n}(H)&=N^{-1}\operatorname{Tr}
                [p_n({H^{[d]}})\Theta_a(H)],&&1\le n\le M,\nonumber\\
 y_{a,n}(H,w)&={m_N}\operatorname{Tr}
       [(ww^*-N^{-1}I)p_n({H^{[d]}})\Theta_a(H)],&&0\le n\le M.
 \label{W:eq:theta-moment-coordinates}
\end{align}

\begin{definition}
\label{W:def:moment-vector}
The moment vector $\mathbf z_N(H,w)$ consists of the coordinates in \eqref{W:eq:polynomial-moment-coordinates} and \eqref{W:eq:theta-moment-coordinates}, in a fixed order. Its real
dimension is
\begin{equation}
 D_N=2M+|\mathcal I|(2M+2)=4(|\mathcal I|+1)\rho_*-2.
 \label{W:eq:moment-coordinates-dimension}
\end{equation}
\end{definition}

For every $\mathbf z\in\mathbb R^{D_N}$, define
\begin{align}
 s_0&=1,&s_n&=u_n+x_n/{m_N}\quad(1\le n\le M),&s_{-1}&=0,\nonumber\\
 q_{a,0}&=t_a/N+y_{a,0}/{m_N},&
 q_{a,n}&=v_{a,n}+y_{a,n}/{m_N}\quad(1\le n\le M).
 \label{W:eq:moment-coordinates-reconstruction}
\end{align}
The values $s_0=1$ and $s_{-1}=0$ are held fixed under differentiation.

Every coordinate is real: the functions of $H$ appearing together in a product commute and are real on the real axis, so their product is Hermitian, also when $H$ has complex
entries. At the coordinate vector associated with a matrix,
\begin{equation}
 s_n=w^*p_n({H^{[d]}})w,\qquad
 q_{a,n}=w^*p_n({H^{[d]}})\Theta_a(H)w.
 \label{W:eq:moment-coordinates-matrixmoments}
\end{equation}
Thus $u_n$ and $v_{a,n}$ contain the factor $N^{-1}$, whereas $t_a$ is an unnormalized trace.

\Needspace{6\baselineskip} Each entry of $F_r^*F_r$, $F_r^*H^{[d]}F_r$, and $F_r^*\Theta_a(H)F_r$ is recovered from these coordinates by the following finite sums.

For $1\le r\le \rho_*$ and indices $0\le i,j<r$ define real symmetric matrices
\begin{align}
 \mathsf G_r(s)_{ij}
   &=\sum_{l=0}^{\min(i,j)}s_{i+j-2l},\label{W:eq:moment-coordinates-G}\\
 \mathsf S_r(s)_{ij}
   &=\sum_{l=0}^{\min(i,j)}
       (s_{i+j-2l+1}+s_{i+j-2l-1}),\label{W:eq:moment-coordinates-S}\\
 \mathsf B_{r,a}(q)_{ij}
   &=\sum_{l=0}^{\min(i,j)}q_{a,i+j-2l}.
 \label{W:eq:moment-coordinates-B}
\end{align}

\begin{lemma}
\label{W:lem:moment-coordinates-momentmaps}
With $F_r=(p_0({H^{[d]}})w,\ldots,p_{r-1}({H^{[d]}})w)$, their evaluations at $\mathbf z_N(H,w)$ are respectively
\[
 F_r^*F_r,\qquad F_r^*{H^{[d]}} F_r,\qquad F_r^*\Theta_a(H)F_r.
\]
All entries involve moments of degree at most $M$. The leading $r$ by $r$ compression of each larger moment matrix equals the corresponding smaller matrix before regularization.
\end{lemma}
\begin{proof}
Use the product and recurrence identities \eqref{W:eq:moment-coordinates-Chebyshev} and \eqref{W:eq:chebyshev-three-term}. For each $0\le i,j<r$, apply that product identity and
then the recurrence to obtain
\[
 \begin{aligned}
 p_i(H^{[d]})p_j(H^{[d]})
  &=\sum_{l=0}^{\min(i,j)}p_{i+j-2l}(H^{[d]}),\\
 H^{[d]}p_i(H^{[d]})p_j(H^{[d]})
  &=\sum_{l=0}^{\min(i,j)}
       (p_{i+j-2l+1}+p_{i+j-2l-1})(H^{[d]}).
 \end{aligned}
\]
Multiplication by $\Theta_a(H)$ is compatible with these identities because it commutes with every function of $H$. Taking scalar products against $w$ proves the three formulas.
The largest index in $\mathsf S_r$ is $2r-1\le M$; the other two matrices use indices at most $2r-2$. Their entries depend only on $i,j$, proving compatibility of leading principal
blocks.
\end{proof}

To compare functions along segments in $\mathbb R^{D_N}$, we extend them to coordinate vectors for which $\mathsf G_r$ need not be positive. We use the regularization method of
Theorem~\ref{W:thm:gram-extension}, with the following specific scalar functions. The cutoff is
\[
 \gamma_\kappa\in C_c^\infty((\kappa,\kappa^{-1});[0,1]),\qquad
 \gamma_\kappa=1\text{ on }[2\kappa,(2\kappa)^{-1}].
\]
Choose $\chi_\kappa\in C_c^\infty((\kappa/2,2/\kappa);[0,1])$ equal to one on $[3\kappa/4,3/(2\kappa)]$. Define, by smooth zero extension,
\begin{equation}
 a_\kappa(t)=\chi_\kappa(t)/t,\qquad
 c_\kappa(t)=\kappa/4+(t-\kappa/4)\chi_\kappa(t).
 \label{W:eq:moment-coordinates-scalars}
\end{equation}
Then
\[
 \kappa/4\le c_\kappa(t)\le2/\kappa\quad(t\in\mathbb R),\qquad
 c_\kappa(t)=t,\quad a_\kappa(t)=t^{-1}
\]
on a neighbourhood of $[\kappa,\kappa^{-1}]$. For each Gram matrix dimension $r$ set
\begin{equation}
 C_r=c_\kappa(\mathsf G_r),\quad C_r=L_rL_r^{\mathsf T},\quad
 \mathsf J_r=L_r^{-1}\mathsf S_rL_r^{-\mathsf T},\quad
 \Gamma_r=\det\gamma_\kappa(\mathsf G_r),
 \label{W:eq:moment-coordinates-J}
\end{equation}
where $L_r$ is the lower triangular Cholesky factor with positive diagonal of the real positive matrix $C_r$.

\Needspace{6\baselineskip}

\begin{lemma}
\label{W:lem:moment-coordinates-smooth}
All maps in \eqref{W:eq:moment-coordinates-J} are $C^\infty$ on $\mathbb R^{D_N}$. One has $0\le\Gamma_r\le1$ there. If $\Gamma_\rho\ne0$ and $r\le\rho$, then $C_r=\mathsf G_r>0$
and $a_\kappa(\mathsf G_r)=\mathsf G_r^{-1}$. At moment vectors obtained from a matrix,
\begin{equation}
 {P_{r,\kappa}^{[d]}}=F_ra_\kappa(\mathsf G_r)F_r^*,\qquad
 0\le{P_{r,\kappa}^{[d]}}\le I_N.
 \label{W:eq:moment-coordinates-projection}
\end{equation}
On $\{\Gamma_\rho\ne0\}$ it is the orthogonal projection onto $\operatorname{span}\{w,H^{[d]}w,\ldots,(H^{[d]})^{r-1}w\}$, and $F_rL_r^{-\mathsf T}$ is the ordered
orthonormalization of the columns of $F_r$.
\end{lemma}
\begin{proof}
Each entry of the three moment matrices is a constant plus a linear function of the coordinates, by \eqref{W:eq:moment-coordinates-G}--\eqref{W:eq:moment-coordinates-B}. For a
smooth compactly supported scalar function, Fourier inversion and the labelled Duhamel formula show that its matrix functional calculus is smooth, including at repeated
eigenvalues. This applies to $a_\kappa,\gamma_\kappa$ and $c_\kappa-\kappa/4$. Positivity bounds for $c_\kappa$ give an everywhere defined Cholesky factor and inverse. To verify
smoothness without choosing eigenvectors, use the recursion
\[
 (L_r)_{jj}=\left((C_r)_{jj}-\sum_{k<j}(L_r)_{jk}^2\right)^{1/2},\qquad
 (L_r)_{ij}=\frac{(C_r)_{ij}-\sum_{k<j}(L_r)_{ik}(L_r)_{jk}}{(L_r)_{jj}}
 \quad(i>j).
\]
For the leading $j\times j$ block $C^{(j)}$,
\[
 (L_r)_{jj}^2
 =\min_{x\in\mathbb R^{j-1}}(x,1)^{\mathsf T}C^{(j)}(x,1).
\]
The lower bound $C^{(j)}\ge(\kappa/4)I$ and the choice $x=0$ give $\kappa/4\le(L_r)_{jj}^2\le2/\kappa$. The recursion therefore uses smooth square roots and reciprocals on
intervals bounded away from zero at every step, proving global smoothness.

The eigenvalues of $\gamma_\kappa(\mathsf G_r)$ lie in $[0,1]$, so
\[
 0\le\Gamma_r=\prod_{j=1}^r
       \gamma_\kappa(\lambda_j(\mathsf G_r))\le1.
\]
If $\Gamma_\rho\ne0$, every eigenvalue of $\mathsf G_\rho$ lies in $(\kappa,\kappa^{-1})$. Interlacing then gives $\kappa I_r<\mathsf G_r<\kappa^{-1}I_r$ for $r\le\rho$.
Equation~\eqref{W:eq:moment-coordinates-scalars} yields $C_r=\mathsf G_r$ and $a_\kappa(\mathsf G_r)=\mathsf G_r^{-1}$. For an SVD $F_r=U\Sigma V^*$, functional calculus gives
\[
 F_ra_\kappa(F_r^*F_r)F_r^*
 =U\,\operatorname{diag}\bigl(s_j^2a_\kappa(s_j^2)\bigr)U^*
 =U\,\operatorname{diag}\bigl(\chi_\kappa(s_j^2)\bigr)U^*,
\]
where zero singular values contribute zero. Thus this matrix lies between zero and $I_N$. When $\Gamma_\rho\ne0$, $F_r$ has full column rank, and all nonzero eigenvalues of the
displayed matrix are one. Moreover,
\[
 (F_rL_r^{-\mathsf T})^*(F_rL_r^{-\mathsf T})
 =L_r^{-1}\mathsf G_rL_r^{-\mathsf T}=I_r,\qquad
 F_r=(F_rL_r^{-\mathsf T})L_r^{\mathsf T}.
\]
The second factor is upper triangular with positive diagonal. Uniqueness of that QR factorization identifies the columns with the orthonormal polynomials in degree order.
\end{proof}

\Needspace{6\baselineskip}

For a test function of the integrated potential with Gram matrix dimension $\rho\ge K_T+2$, write, indexing the tridiagonal matrix from one,
\[
 a_k=(\mathsf J_\rho)_{kk},\qquad
 b_k=(\mathsf J_\rho)_{k,k+1},\qquad b_0=1.
\]
For arbitrary real values of the moment coordinates these are smooth real numbers; $b_k$ need not be positive there. For $\sigma\in\{+1,-1\}$ define
\begin{equation}
 Y_0^\sigma=0,\qquad
 Y_k^\sigma={m_N}\left(2k-\sigma\sum_{j=1}^ka_j
                   -1-2\sum_{j=1}^{k-1}b_j-b_k\right),
 \quad 1\le k\le K_T.
 \label{W:eq:moment-coordinates-nodes}
\end{equation}
Interpolate linearly at sampling step ${m_N}^{-1}$, using the value at $K_T$ when ${m_N}t=K_T$ exactly.
\begin{lemma}
\label{W:lem:moment-coordinates-integrated-potential}
At moment vectors obtained from a matrix on $\{\Gamma_\rho\ne0\}$ this is precisely the integrated potential of $\sigma{H^{[d]}}$, with the auxiliary coefficient at zero
$\bar b_0=1$.
\end{lemma}
\begin{proof}
On $\{\Gamma_\rho\ne0\}$, Lemma~\ref{W:lem:moment-coordinates-smooth} identifies the Cholesky matrix with the ordered polynomial orthogonalization.
Section~\ref{W:sec:cholesky-representation-of-jacobi-coefficients} therefore gives the Lanczos coefficients with $b_j>0$. Since $k+1\le K_T+1<\rho$, every coefficient used here is
inside that block. Their sum satisfies
\[
 \sum_{j=1}^k(\bar b_{j-1}+b_j)=1+2\sum_{j<k}b_j+b_k.
\]
Substitution into \eqref{W:eq:projection-T} proves \eqref{W:eq:moment-coordinates-nodes}. Lemma~\ref{W:lem:tightness-reflection}, applied to $H^{[d]}$, changes $a_j$ to $-a_j$ and
preserves $b_j$. Thus the two signs give the integrated potentials of $H^{[d]}$ and $-H^{[d]}$, respectively. Since $h_d$ is odd, $-H^{[d]}=h_d(-H)$. Taking the same convex
combination of adjacent values proves the identity for the linear interpolants as well.
\end{proof}

\Needspace{6\baselineskip}

The logarithmic decomposition \eqref{W:eq:moment-coordinates-fullboundary} is used only where the Gram matrices are positive. The formula in \eqref{W:eq:moment-coordinates-nodes}
supplies a smooth real-valued extension everywhere; outside the set of agreement it need not come from a Jacobi matrix with positive off-diagonal entries.

For $a=(\sigma,E)\in\mathcal I$, projection rank $r$, and Gram matrix dimension $\rho\ge r$, define on $\mathbb R^{D_N}$
\begin{equation}
 C_{r,a}(\mathbf z)=t_a-
   \operatorname{Tr}\bigl[a_\kappa(\mathsf G_r(s))\mathsf B_{r,a}(q)\bigr].
 \label{W:eq:moment-coordinates-localization}
\end{equation}
At moment vectors obtained from a matrix this is
\[
 \operatorname{Tr}\bigl[(I-{P_{r,\kappa}^{[d]}})\Theta_a(H)\bigr],
 \qquad 0\le C_{r,a}\le t_a\le N.
\]
For arbitrary moment vectors we use the following smooth extension of the scalar function $c/(1+c)$.

Choose $\chi_-\in C^\infty(\mathbb R;[0,1])$ equal to zero for $c\le-1/2$ and to one for $c\ge-1/4$, and set
\begin{equation}
 \widetilde\psi(c)=
 \begin{cases}
  0,&c\le-1/2,\\
  \chi_-(c)c/(1+c),&c>-1/2.
 \end{cases}
 \label{W:eq:moment-coordinates-bounded-transform}
\end{equation}
For $g(c)=c/(1+c)$ and $j\ge1$,
\[
 g^{(j)}(c)=(-1)^{j+1}j!(1+c)^{-j-1},\qquad c>-1/2.
\]
On this interval $|g(c)|\le1$ and $|g^{(j)}(c)|\le j!2^{j+1}$. Leibniz's formula therefore gives
\[
 \widetilde\psi^{(q)}(c)=\sum_{j=0}^q\binom qj
                \chi_-^{(q-j)}(c)g^{(j)}(c),\qquad
 \|\widetilde\psi^{(q)}\|_\infty<\infty.
\]
All derivatives of $\chi_-$ vanish at $-1/2$, so the two definitions and all their derivatives agree there. Thus $\widetilde\psi\in C_b^\infty$, $|\widetilde\psi|\le1$, and
$\widetilde\psi(c)=c/(1+c)$ for $c\ge0$.

\Needspace{6\baselineskip}

\begin{lemma}
\label{W:lem:moment-coordinates-localization}
At moment vectors obtained from a matrix, $\Gamma_\rho\widetilde\psi(C_{r,a})$ is exactly the regularized localization statistic multiplied by $\Gamma_\rho$. If in addition
$h_d(H)=H$ and $\Gamma_\rho\ne0$, then, for
$\Pi_a=\mathbf1_{(-\infty,E]}({m_N}^2(2I-\sigma H))$,
\begin{equation}
 \|(I-P_r)\Pi_a\|^2\le4C_{r,a}.
 \label{W:eq:moment-coordinates-detection}
\end{equation}
The estimate holds for every real $E$ on the stated event, and $I-P_r$ projects onto the full complement of the initial Lanczos subspace.
\end{lemma}
\begin{proof}
Put $P=F_ra_\kappa(\mathsf G_r)F_r^*$. Lemma~\ref{W:lem:moment-coordinates-momentmaps} and cyclicity give
\[
 C_{r,a} =\operatorname{Tr}\Theta_a -\operatorname{Tr}\bigl(a_\kappa(\mathsf G_r)F_r^*\Theta_aF_r\bigr) =\operatorname{Tr}\bigl(\Theta_a^{1/2}(I-P)\Theta_a^{1/2}\bigr)\ge0,
\]
because $0\le P\le I$ by Lemma~\ref{W:lem:moment-coordinates-smooth}. Since the trace is nonnegative, the smooth scalar extension agrees with $c/(1+c)$ there. Under the additional
assumptions of the lemma, $P_r$ is the exact Krylov projection for $H$ and for $-H$, and $\Theta_a\ge\tfrac14\Pi_a$. More explicitly, $\|\Theta_a^{-1/2}\Pi_a\|\le2$ on that
spectral subspace, giving
\[
 \|(I-P_r)\Pi_a\|^2
 \le4\|(I-P_r)\Theta_a^{1/2}\|_{\rm HS}^2=4C_{r,a}.
\]
The complement includes all later invariant blocks of the Lanczos tridiagonal matrix. The identification of the two projections uses $h_d(H)=H$.
Sections~\ref{W:sec:resolvents-on-exceptional-support} and~\ref{W:sec:derivative-bounds-for-the-spectral-cutoff} bound the probability of its complement.
\end{proof}

Take $K=K_T$, $h_N=m_N^{-1}$ and $\mathcal P_\delta$ as in Lemma~\ref{W:lem:tightness-mesh-geometry}. For $D=|\mathcal P_\delta|>0$, use the functions in
\eqref{W:eq:smooth-modulus-definition}, with $\lambda=4\log(1+2D)/u$, and write
\begin{equation}
 \mathfrak S(y)=S_{\lambda,\mathcal P_\delta}(y),\qquad
 \mathfrak M(y)=\Phi_{\mathcal P_\delta,u}(y)
              =\vartheta(\mathfrak S(y)/u).
 \label{W:eq:moment-coordinates-modulus}
\end{equation}
For $D=0$ put $\mathfrak M=0$. The following specialization will be used with the moment coordinates.

\Needspace{6\baselineskip}

\begin{lemma}
\label{W:lem:moment-coordinates-modulus}
Let $M_\delta(y)=\max_{\mathcal P_\delta}|y_j-y_i|$, with zero for the empty pair set. Then
\begin{equation}
 \mathbf1_{\{M_\delta(y)\ge u\}}
 \le\mathfrak M(y)\le\mathbf1_{\{M_\delta(y)>u/4\}}.
 \label{W:eq:moment-coordinates-modsandwich}
\end{equation}
For the linear interpolation on $[0,K/{m_N}]$,
\[
 \omega_T(y,\delta)\le M_\delta(y)
      \le\omega_{K/{m_N}}(y,\delta+2/{m_N}).
\]
For every fixed $q\ge1$,
\begin{equation}
 |D_{y}^q\mathfrak M(y)[v_1,\ldots,v_q]|
 \le C_qu^{-q}(1+\log(1+2D))^{q-1}
                  \prod_j\|v_j\|_\infty.
 \label{W:eq:moment-coordinates-modderivative}
\end{equation}
\end{lemma}
\begin{proof}
Suppose first that $D>0$, and put $L_D=\log(1+2D)$. The largest term in the exponential sum defining $\mathfrak S$ is $\exp(\lambda M_\delta(y))$. Every term, including its
constant term, is at most this value. Hence
\[
 M_\delta(y)\le\mathfrak S(y)
 \le M_\delta(y)+\frac{L_D}{\lambda}
 =M_\delta(y)+\frac u4.
\]
If $M_\delta(y)\ge u$, then $\vartheta(\mathfrak S(y)/u)=1$. If $M_\delta(y)\le u/4$, then $\mathfrak S(y)/u\le1/2$ and this value of $\vartheta$ is zero. Since
$0\le\vartheta\le1$, these two implications prove \eqref{W:eq:moment-coordinates-modsandwich}.

For $s,t\in[0,T]$ with $|t-s|\le\delta$, write their interpolated values as convex combinations of the two endpoints of their cells. If the endpoint indices are $i,i+1$ and
$j,j+1$, respectively, then
\[
 y(t)-y(s)=\sum_{a,b\in\{0,1\}}\alpha_a\beta_b
                       (y_{j+b}-y_{i+a}),
 \qquad \alpha_a,\beta_b\ge0,\quad
 \sum_a\alpha_a=\sum_b\beta_b=1.
\]
Each pair of endpoints in this sum is at distance at most $\delta+2h_N$. Equal endpoints contribute zero; ordering the other pairs puts them in $\mathcal P_\delta$. Thus
$|y(t)-y(s)|\le M_\delta(y)$. Conversely, every pair in $\mathcal P_\delta$ consists of points of $[0,Kh_N]$ at distance at most $\delta+2h_N$. Taking the two suprema proves the
interpolation inequalities.

The logarithmic derivative calculation in Lemma~\ref{W:lem:tightness-smooth-max} gives, for $r\ge1$,
\[
 |D_y^r\mathfrak S(y)[v_1,\ldots,v_r]|
 \le C_r\lambda^{r-1}\prod_{j=1}^r\|v_j\|_\infty.
\]
Let $\Pi(q)$ denote the partitions of $\{1,\ldots,q\}$. The chain rule therefore gives
\[
 D_y^q\mathfrak M(y)[v_1,\ldots,v_q]
 =\sum_{\pi\in\Pi(q)}u^{-|\pi|}
   \vartheta^{(|\pi|)}(\mathfrak S(y)/u)
   \prod_{B\in\pi}D_y^{|B|}\mathfrak S(y)[v_j:j\in B].
\]
A partition with $k$ blocks contributes at most $C_qu^{-k}\lambda^{q-k}\prod_j\|v_j\|_\infty =C_qu^{-q}(4L_D)^{q-k}\prod_j\|v_j\|_\infty$. There are finitely many partitions for
fixed $q$, and $1\le k\le q$, so their sum is bounded by the right side of \eqref{W:eq:moment-coordinates-modderivative}.

If $D=0$, then $K=0$: for $K\ge1$ an adjacent pair would belong to $\mathcal P_\delta$. The interpolant has zero modulus, and $\mathfrak M$ and all its derivatives vanish. This
proves the remaining case.
\end{proof}

\subsubsection{Cutoffs in the moment variables}
\label{W:sec:cutoffs-in-moment-variables}

The moment representation will be used along segments of coordinate vectors, as well as at vectors obtained from matrices. A smooth coordinate cutoff gives bounded support in these
variables while leaving every vector produced by a matrix unchanged. Its coordinate scales come from the polynomial bound on the range of $h_d$. For $0\le n\le M$, these scales are

\begin{equation}
 B_n=(n+1)\exp(n\sqrt{2d}).
 \label{W:eq:moment-coordinates-Bn}
\end{equation}
The elementary trigonometric and hyperbolic formulas for $U_n$ show that $|p_n(x)|\le B_n$ on $|x|\le2+2d$. Since $\Theta_a$ is a contraction, every moment vector associated with a
matrix satisfies
\begin{align}
 |u_n|,|v_{a,n}|&\le B_n,&
 |x_n|,|y_{a,n}|&\le2{m_N}B_n,& 0\le t_a\le N.
 \label{W:eq:moment-coordinate-bounds}
\end{align}
The bounds for $y_{a,0}$ use $B_0=1$. These statements hold for every Hermitian $H$, without a norm event, because the polynomials are evaluated at $h_d(H)$.

Choose an even $\chi\in C_c^\infty(\mathbb R;[0,1])$ equal to one on $[-2,2]$ and zero outside $(-3,3)$. Define
\begin{equation}
 \begin{split}
 \Xi_N(\mathbf z)={}&
 \prod_{n=1}^{M}\chi(u_n/B_n)\chi(x_n/(2{m_N}B_n))\\
 &\times\prod_{a\in\mathcal I}
 \left\{\chi(t_a/N)
 \prod_{n=1}^{M}\chi(v_{a,n}/B_n)
 \prod_{n=0}^{M}\chi(y_{a,n}/(2{m_N}B_n))\right\}.
 \end{split}
 \label{W:eq:moment-coordinates-box}
\end{equation}
This cutoff is identically one on an open neighbourhood of every moment vector associated with a matrix, by \eqref{W:eq:moment-coordinate-bounds}. In particular
$\Xi_N(\mathbf z_N(H,w))\equiv1$ as an identity of smooth matrix functions; its matrix entry derivatives of positive order are zero. On arbitrary segments of coordinate vectors,
its derivatives occur in the product rule. The proof of Theorem~\ref{W:thm:moment-coordinates-main} also checks the coordinate bound along segments joining two moment vectors
associated with matrices.

We require cutoff bounds uniform in the Gram matrix dimension. For a coordinate box with radii $\beta_j\ge1$, write $\Xi(\mathbf z)=\prod_j\chi(z_j/\beta_j)$ and
$\mathfrak b(v)=\sum_j|v_j|/\beta_j$.

\begin{lemma}
\label{W:lem:observable-det}
If $A$ is an $r$ by $r$ complex matrix with $\|A\|\le1$, then
\begin{equation}
 |D_{A}^k\det(A)[B_1,\ldots,B_k]|
 \le\prod_{j=1}^k\|B_j\|_1,\qquad k\ge1.
 \label{W:eq:observable-detbound}
\end{equation}
For $k>r$ the derivative is zero. The bound includes singular matrices $A$.
\end{lemma}
\begin{proof}
For $k>r$ the determinant has degree $r$ in the entries of $A$, so the claimed derivative vanishes. Suppose $1\le k\le r$. We first take $B_j=u_jv_j^*$, with
$u_j,v_j\in\mathbb C^r$. Set $U=(u_1,\ldots,u_k)$, $V=(v_1,\ldots,v_k)$, and $T=\operatorname{diag}(t_1,\ldots,t_k)$. Multilinearity of the Fr\'echet derivative gives
\[
 D_A^k\det(A)[u_1v_1^*,\ldots,u_kv_k^*]
   =[t_1\cdots t_k]\det(A+UTV^*),
\]
where the bracket denotes the coefficient of the indicated monomial. Indeed $\partial_{t_1}\cdots\partial_{t_k}$ at the origin extracts that coefficient without a factorial,
because each variable has exponent one. If the columns of $U$ or $V$ are linearly dependent, that coefficient is zero by the alternating multilinearity of the determinant in the
corresponding update columns or rows.

Otherwise take thin QR factorizations $U=Q_UR_U$ and $V=Q_VR_V$ and extend $Q_U,Q_V$ to unitary matrices $\widetilde Q_U,\widetilde Q_V$. Unitary conjugations on the left and right
change the determinant only by factors of modulus one. In the transformed coordinates the update has the block form
\[
 \widetilde Q_U^*UTV^*\widetilde Q_V
 =\begin{pmatrix}R_UTR_V^*&0\\0&0\end{pmatrix}.
\]
Write $A'=\widetilde Q_U^*A\widetilde Q_V$. To obtain $t_1\cdots t_k$ from the determinant expansion one must select each of the $k$ update columns and no other update column.
Expansion along these columns gives exactly
\[
 [t_1\cdots t_k]\det(A+UTV^*)
   =\omega\det(R_U)\overline{\det(R_V)}
     \det A'_{\{k+1,\ldots,r\},\{k+1,\ldots,r\}},
 \qquad |\omega|=1.
\]
For $k=r$ the last determinant is interpreted as one. Since $\|A'\|=\|A\|\le1$, the singular values of its lower principal compression are at most one; hence its determinant has
modulus at most one, even when $A$ is singular. Hadamard's inequality gives
\[
 |\det R_U|\le\prod_{j=1}^k\|u_j\|,
 \qquad |\det R_V|\le\prod_{j=1}^k\|v_j\|.
\]
We have therefore proved $|D_A^k\det(A)[u_1v_1^*,\ldots,u_kv_k^*]| \le\prod_j\|u_j\|\|v_j\|$ without assuming that $A$ is invertible.

For arbitrary $B_j$, choose a singular value decomposition $B_j=\sum_{l=1}^{r}s_{j,l}u_{j,l}v_{j,l}^*$, where the singular vectors have unit norm and $s_{j,l}\ge0$. By
multilinearity of $D_A^k\det$ and the rank-one estimate,
\[
 \begin{aligned}
 |D_A^k\det(A)[B_1,\ldots,B_k]|
 &\le\sum_{l_1,\ldots,l_k}\prod_{j=1}^k s_{j,l_j}
 =\prod_{j=1}^k\sum_{l=1}^r s_{j,l}
 =\prod_{j=1}^k\|B_j\|_1.
 \end{aligned}
\]
All the sums are finite. This proves \eqref{W:eq:observable-detbound} in every dimension and at every contraction $A$.
\end{proof}

\begin{lemma}
\label{W:lem:gram-cutoff-trace-derivatives}
For $\Gamma(G)=\det\gamma_\kappa(G)$, Hermitian $G$, and $q\ge1$,
\begin{align}
 |D^q\Gamma(G)[E_1,\ldots,E_q]|
  &\le C_{q,\kappa}\prod_{j=1}^q\|E_j\|_1,
   \label{W:eq:gram-cutoff-trace-derivatives}\\
 |D^q\Xi(\mathbf z)[v_1,\ldots,v_q]|
  &\le C_{q,\chi}\prod_{j=1}^q\mathfrak b(v_j).
   \label{W:eq:moment-box-derivatives}
\end{align}
The constants are independent of the numbers of rows and coordinates.
\end{lemma}
\begin{proof}
Lemma~\ref{W:lem:observable-det} gives the determinant estimate. To apply it to $\gamma_\kappa(G)$, we first bound the derivatives of this matrix function.

Use the Fourier convention $\gamma_\kappa(x)=(2\pi)^{-1}\int\widehat\gamma_\kappa(t)e^{\ii tx}\,\dd t$. Apply the differentiated Duhamel formula from the proof of
Lemma~\ref{W:lem:gram-clipping} with $H=G$ and $W_i=E_i$. Place a fixed direction $E_l$ in trace norm and all other directions in operator norm. Unitarity of the exponential
factors and the simplex volume $1/j!$ give
\[
 \|D^je^{\ii tG}[E_1,\ldots,E_j]\|_1
 \le |t|^j\|E_l\|_1\prod_{i\ne l}\|E_i\|.
\]
All Fourier moments of $\widehat\gamma_\kappa$ are integrable because $\gamma_\kappa\in C_c^\infty$. Differentiation under its Fourier integral is therefore justified and yields
\[
 \|D^j\gamma_\kappa(G)[E_1,\ldots,E_j]\|_1
 \le M_j\|E_l\|_1\prod_{i\ne l}\|E_i\|
 \le M_j\prod_i\|E_i\|_1,
 \quad
 M_j=\frac1{2\pi}\int|t|^j|\widehat\gamma_\kappa(t)|\,\dd t.
\]
For $S=\{1,\ldots,q\}$, the chain rule is
\[
 D_S\Gamma(G)=\sum_{\pi\in\Pi(S)}
 D^{|\pi|}\det(\gamma_\kappa(G))
       [D_B\gamma_\kappa(G):B\in\pi].
\]
Since $\|\gamma_\kappa(G)\|\le1$, the determinant bound gives the first conclusion with
$C_{q,\kappa}=\sum_{\pi\in\Pi(S)}\prod_{B\in\pi}M_{|B|}$.

For the coordinate cutoff the exact formula is
\[
 D_S\Xi(\mathbf z)[v_1,\ldots,v_q]
 =\sum_{(j_1,\ldots,j_q)}\left(\prod_{a=1}^q\frac{v_{a,j_a}}{\beta_{j_a}}\right)
 \prod_j\chi^{(n_j)}(z_j/\beta_j),
 \quad n_j=|\{a:j_a=j\}|.
\]
At most $q$ factors have positive derivative order, and all others are at most one. The product of scalar derivative norms is bounded by a constant depending only on $q,\chi$.
Summing the absolute values of the directional coefficients gives $\prod_a\sum_j|v_{a,j}|/\beta_j$, which proves the bound for the coordinate cutoff.
\end{proof}

\Needspace{6\baselineskip}

With the fixed parameters above, define the following functions on all of $\mathbb R^{D_N}$, each using its own Gram matrix dimension:
\begin{align}
 \Phi_N^{\Gamma,\rho}(\mathbf z)
   &=\Xi_N(\mathbf z)\Gamma_\rho(s),\label{W:eq:moment-coordinates-testsG}\\
 \Phi_N^{\mathrm F,\sigma}(\mathbf z)
   &=\Xi_N(\mathbf z)\Gamma_\rho(s)
     \exp\left(\mathrm i\sum_j\lambda_jY^\sigma(t_j;\mathbf z)\right),
     \quad \rho\ge K_T+2,\label{W:eq:moment-coordinates-testsF}\\
 \Phi_N^{\mathrm M,\sigma}(\mathbf z)
   &=\Xi_N(\mathbf z)\Gamma_\rho(s)
           \mathfrak M(Y_0^\sigma,\ldots,Y_{K_T}^\sigma),
     \quad \rho\ge K_T+2,\label{W:eq:moment-coordinates-testsM}\\
 \Phi_N^{\mathrm{loc},r,\rho,a}(\mathbf z)
   &=\Xi_N(\mathbf z)\Gamma_\rho(s)
                  \widetilde\psi(C_{r,a}(\mathbf z)),
     \quad r\le\rho.\label{W:eq:moment-coordinates-testsC}
\end{align}

\begin{theorem}
\label{W:thm:moment-coordinates-main}
Each function in \eqref{W:eq:moment-coordinates-testsG}--\eqref{W:eq:moment-coordinates-testsC} belongs to $C_c^\infty(\mathbb R^{D_N})$, has absolute value at most one, and
evaluates at $\mathbf z_N(H,w)$ to its corresponding regularized matrix statistic. On $\|H\|<2+d$ these evaluations agree locally with the original statistics. For every fixed
$q\ge1$,
\begin{equation}
 \sum_{j_1,\ldots,j_q=1}^{D_N}
  \|\partial_{j_1}\cdots\partial_{j_q}\Phi_N^\nu\|_\infty
       \le C_{q,\kappa,\mathrm{parameters}}N^{C_q},
 \label{W:eq:moment-coordinates-polybound}
\end{equation}
for a finite exponent $C_q$ independent of $N$.
\end{theorem}
\begin{proof}
The finite-dimensional maps $\mathsf G_r,\mathsf S_r,\mathsf B_{r,a}$ consist of fixed matrices plus linear expressions in the coordinates introduced in
\eqref{W:eq:polynomial-moment-coordinates}; their explicit coefficient matrices are given by \eqref{W:eq:moment-coordinates-G}--\eqref{W:eq:moment-coordinates-B}. The scalar
regularizations $a_\kappa,c_\kappa,\gamma_\kappa$ and the smooth truncation of the localization statistic are defined on all real arguments. The matrix $C_r=c_\kappa(\mathsf G_r)$
is uniformly positive definite, so its Cholesky factor and inverse are smooth for every real coordinate vector.
Lemmas~\ref{W:lem:moment-coordinates-smooth}--\ref{W:lem:moment-coordinates-modulus} then show that the four observables in the statement are smooth. Their factors satisfy
\[
 0\le\Xi_N,\Gamma_\rho,\mathfrak M\le1,
 \qquad |\widetilde\psi|\le1,
 \qquad |\exp(\mathrm i\sum_j\lambda_jY_j)|=1.
\]
Thus $|\Phi_N^\nu|\le\Xi_N\Gamma_\rho\le1$. Each scalar factor $\chi(z_j/\beta_j)$ in $\Xi_N$ is supported on $[-3\beta_j,3\beta_j]$, so
$\operatorname{supp}\Phi_N^\nu\subset\prod_{j\le D_N}[-3\beta_j,3\beta_j]$; this is a compact set for each $N$.

At an actual matrix $H$, the coordinate bound \eqref{W:eq:moment-coordinates-box} places $\mathbf z_N(H,w)$ in the region where $\Xi_N=1$, with every derivative of positive order
equal to zero in a neighbourhood. Lemma~\ref{W:lem:moment-coordinates-momentmaps} gives
\[
 \mathsf G_r=F_r^*F_r,\qquad
 \mathsf S_r=F_r^*H^{[d]}F_r,\qquad
 \mathsf B_{r,a}=F_r^*\Theta_a(H)F_r.
\]
If $\Gamma_\rho(\mathbf z_N(H,w))\ne0$, every eigenvalue of the relevant Gram matrices lies in the interval on which the regularizations agree with the identity and the reciprocal.
By principal-submatrix interlacing this applies to each dimension $r\le\rho$. Cholesky orthogonalization identifies the integrated potential and projection with those constructed
from the polynomial vectors of $H^{[d]}$. The localization statistic is $\operatorname{Tr}[(I-P_r^{[d]})\Theta_a(H)]$: its polynomial projection uses $H^{[d]}$, whereas the
spectral weight is still evaluated at $H$. These identifications are exactly those of Lemmas~\ref{W:lem:moment-coordinates-integrated-potential},
\ref{W:lem:moment-coordinates-localization}, and \ref{W:lem:moment-coordinates-modulus}. If $\Gamma_\rho=0$, the corresponding functions multiplied by their cutoffs agree because
they are both zero. Finally, if $\|H\|<2+d$, take $\varepsilon=(2+d-\|H\|)/2$. Every $H'$ with $\|H'-H\|<\varepsilon$ has spectrum in $(-2-d,2+d)$, on which $h_d(x)=x$. Thus
$h_d(H')=H'$ throughout this neighbourhood, proving local agreement rather than merely equality at one matrix.

For later interpolation between two matrix coordinate vectors $\mathbf z$ and $\mathbf z'$, \eqref{W:eq:moment-coordinates-box} also gives $|(1-t)z_j+tz_j'|\le\beta_j$ for
$0\le t\le1$. The scalar cutoff is identically one on $[-2,2]$; hence all positive derivatives of $\Xi_N$ vanish on a neighbourhood of this entire segment. For the global estimate
we include every derivative of $\Xi_N$ in the product rule, including its transition regions.

We now prove a uniform polynomial estimate, keeping track of the derivative orders. Write $r\le\rho_*\le C_*N^{1/3}$, and let $q\ge1$ be fixed. The coordinate count $D_N$, the
scalar multiplier $m_N$, and the reciprocal of the smallest coordinate radius are at most a fixed power of $N$. The bounds \eqref{W:eq:moment-coordinates-box} and the choice
$1<\gamma<2$ yield
\[
 \max_{n\le2\rho_*}B_n\le (2\rho_*+1)
 \exp\{C\rho_*\sqrt d\}
 =N^{1/3+o(1)}.
\]
The estimates for the intermediate matrices will be used on the closed box
\[
 \mathcal B_N=\prod_{j=1}^{D_N}[-3\beta_j,3\beta_j].
\]
Every derivative of $\Xi_N$, including order zero, is supported in this box. The largest radius is at most $N^2$ for all sufficiently large $N$: the possible radii are $B_n$,
$2m_NB_n$, and $N$. On $\mathcal B_N$, the reconstructions of $s_n$ and $q_{a,n}$ therefore have absolute value at most $CN^2$. Each entry of
$\mathsf G_r,\mathsf S_r,\mathsf B_{r,a}$ is a sum of at most $2r$ such quantities. A derivative with respect to a single coordinate has entries bounded by $Cr$, and derivatives of
order two or higher vanish. Since $r\le C_*N^{1/3}$ and $\|M\|\le\|M\|_{\rm HS}\le r\max_{i,j}|M_{ij}|$, a fixed $A\ge4$ gives
\[
 \sup_{\mathbf z\in\mathcal B_N}
 \bigl(\|M(\mathbf z)\|_{\rm HS}+\max_j\|\partial_jM(\mathbf z)\|_{\rm HS}\bigr)
 \le N^A,
 \qquad M\in\{\mathsf G_r,\mathsf S_r,\mathsf B_{r,a}\},
\]
for all sufficiently large $N$, with $D_N,m_N,\rho_*\le N^A$ as well. The first derivatives, being constant matrices, have the same bounds outside the box. We use the bounds for
the undifferentiated matrices only on $\mathcal B_N$.

Let $\varphi$ denote any one of the fixed scalar regularizations. Fourier inversion and the iterated Duhamel formula, as in Lemma~\ref{W:lem:gram-clipping}, give for every fixed
$k$ and self-adjoint $G$
\[
 \|D^k\varphi(G)[E_1,\ldots,E_k]\|_{\rm HS}
 \le C_{k,\kappa}\|E_1\|_{\rm HS}\prod_{j=2}^k\|E_j\|,
\]
with the same bound if another factor is estimated in Hilbert--Schmidt norm. The chain rule for $C_r=c_\kappa(\mathsf G_r)$ has no derivatives of $\mathsf G_r$ of order two or
above. Consequently, for any nonempty set $S$ of coordinate derivative labels,
\[
 \|D_SC_r\|_{\rm HS}\le C_{|S|,\kappa}N^{A|S|}.
\]
The uniform lower and upper spectral bounds on $C_r=L_rL_r^{\mathsf T}$ imply $\|L_r\|+\|L_r^{-1}\|\le C_\kappa$. Applying \eqref{W:eq:gram-cholesky} with $G_\kappa^{[d]}$ replaced
by $C_r$, and using contraction of triangular truncation in Hilbert--Schmidt norm, gives
\[
 a_k:=\max_{r,\,|S|=k}\|D_SL_r\|_{\rm HS}
 \le C_{k,\kappa}\left(N^{Ak}+
       \sum_{j=1}^{k-1}\binom kj a_ja_{k-j}\right).
\]
Starting with $a_1\le C_\kappa N^A$, induction gives $a_k\le C_{k,\kappa}N^{Ak}$; the sum has finitely many terms, and every product of the induction bounds carries the power
$N^{A(j+k-j)}$. The ordered inverse formula (\ref{W:eq:gram-inverse-labels}) gives the same bound for $D_SL_r^{-1}$, with the zeroth derivative bounded by $C_\kappa$. Apply the
Leibniz rule for the product of three factors to $\mathsf J_r=L_r^{-1}\mathsf S_rL_r^{-\mathsf T}$: if $S_0\sqcup S_1\sqcup S_2=S$, then
\[
 \|(D_{S_0}L_r^{-1})(D_{S_1}\mathsf S_r)
                   (D_{S_2}L_r^{-\mathsf T})\|
 \le C_{k,\kappa}N^{A(k+1)}.
\]
Here $D_{S_1}\mathsf S_r=0$ for $|S_1|\ge2$. There are $3^k$ ordered assignments, and hence
\[
 \|D_S\mathsf J_r\|\le C_{k,\kappa}N^{A(k+1)},
 \qquad
 \max_j|D_SY_j^\sigma|\le C_{k,\kappa}N^{A(k+3)}.
\]
The second inequality follows by summing no more than $2\rho_*$ matrix entries in the integrated potential and multiplying by $m_N$; linear interpolation is a convex combination of
adjacent values.

For the determinant factor, \eqref{W:eq:gram-cutoff-trace-derivatives} gives
\[
 |D_S\Gamma_r|\le C_{k,\kappa}
        \prod_{j\in S}\|D_j\mathsf G_r\|_1
 \le C_{k,\kappa}N^{(A+1)k},
\]
using $\|E\|_1\le\sqrt r\|E\|_{\rm HS}$ and enlarging $A$ once if necessary. For the scalar localization argument
$C_{r,a}=t_a-\operatorname{Tr}[a_\kappa(\mathsf G_r)\mathsf B_{r,a}]$, use the derivative formula for matrix functions and the product rule. Since $\mathsf B_{r,a}$ is linear in
the coordinates, the same calculation yields $|D_SC_{r,a}|\le C_{k,\kappa}N^{A(k+3)}$. The estimate for $\Gamma_r$ holds for every real coordinate vector. The bounds for
$Y_j^\sigma$ and $C_{r,a}$ hold on $\mathcal B_N$, where the undifferentiated matrices $\mathsf S_r$ and $\mathsf B_{r,a}$ have just been bounded. For example, the latter
calculation is
\[
 D_SC_{r,a}=D_St_a-
 \sum_{U\sqcup V=S}\operatorname{Tr}
 \bigl(D_Ua_\kappa(\mathsf G_r)\,D_V\mathsf B_{r,a}\bigr).
\]
Only $|V|\le1$ contributes. For $U\ne\varnothing$, the functional calculus estimate bounds its first factor by $C_{|U|,\kappa}N^{A|U|}$ in operator norm; for $U=\varnothing$ its
norm is bounded by a constant. On $\mathcal B_N$ the second factor has norm at most $N^A$. Thus each trace has absolute value at most
$C_{k,\kappa}rN^{A(k+1)}\le C_{k,\kappa}N^{A(k+2)}$. There are at most $2^k$ choices of $U,V$, and $D_St_a$ is either zero or one. This proves the stated, slightly larger bound
$C_{k,\kappa}N^{A(k+3)}$.

To complete the coordinate estimates, a scalar composition at order $k$ is indexed by the partitions of its $k$ labels. For example,
\[
 D_S e^{\mathrm i\sum_j\lambda_jY_j}
  =e^{\mathrm i\sum_j\lambda_jY_j}
    \sum_{\pi\in\Pi(S)}\mathrm i^{|\pi|}
                 \prod_{B\in\pi}D_B\!\left(\sum_j\lambda_jY_j\right).
\]
For a partition with $l$ blocks the preceding bound contributes at most $C_kN^{A(k+3l)}\le C_kN^{4Ak}$. The derivatives of $\widetilde\psi$ are bounded, and
\eqref{W:eq:moment-coordinates-modderivative} gives an additional factor at most $C_k(1+\log N)^k$ for the modulus statistic. The derivatives of $\Xi_N$ satisfy the explicit
labelled formula \eqref{W:eq:moment-box-derivatives}, with $\beta_j\ge1$, so every coordinate derivative of it is bounded by $C_k$. For a selected scalar statistic $F$, the
complete product rule is
\[
 D_S(\Xi_N\Gamma_\rho F)
 =\sum_{S_0\sqcup S_1\sqcup S_2=S}
   (D_{S_0}\Xi_N)(D_{S_1}\Gamma_\rho)(D_{S_2}F).
\]
Let $r_i=|S_i|$, so $r_0+r_1+r_2=q$. A summand can be nonzero only on $\mathcal B_N$, since it contains a derivative of $\Xi_N$. Its absolute value is bounded there by
\[
 C_qN^{(A+1)r_1+4Ar_2}(1+\log N)^{r_2}
 \le C_qN^{5Aq},
\]
where $(1+\log N)^{r_2}\le N^{Ar_2}$ for sufficiently large $N$ and $A\ge1$. This includes $r_i=0$, with the corresponding undifferentiated factors bounded by one. There are $3^q$
assignments of the labels, and the partition sums within the scalar compositions have already been included in $C_q$. Every summand vanishes outside $\mathcal B_N$. Thus the larger
exponent $10Aq$ also gives the global bound
\[
 \max_{j_1,\ldots,j_q}
 \|\partial_{j_1}\cdots\partial_{j_q}\Phi_N^\nu\|_\infty
 \le C_{q,\kappa,\mathrm{parameters}}N^{10Aq}
\]
for all four families of functions and all sufficiently large $N$. Finally $D_N^q\le N^{Aq}$, so summing over the ordered choices of the coordinate indices proves
\eqref{W:eq:moment-coordinates-polybound} with $C_q=11Aq$. For each of the finitely many remaining admissible sizes, smoothness and compact support make the displayed derivative
suprema finite; increasing the constant covers those sizes. The final estimate is global because the complete functions and all their derivatives vanish outside $\mathcal B_N$.
\end{proof}

\Needspace{6\baselineskip}

If $\mathcal I$ is closed under $(\sigma,E)\mapsto(-\sigma,E)$, define a linear involution $\mathscr R$ on the coordinate vectors by
\begin{align}
 u_n&\mapsto(-1)^nu_n,& x_n&\mapsto(-1)^nx_n,\nonumber\\
 t_{\sigma,E}&\mapsto t_{-\sigma,E},&
 v_{\sigma,E,n}&\mapsto(-1)^nv_{-\sigma,E,n},&
 y_{\sigma,E,n}&\mapsto(-1)^ny_{-\sigma,E,n}.
 \label{W:eq:moment-coordinates-reflection}
\end{align}

\begin{proposition}
\label{W:prop:moment-coordinates-reflection}
The moment map and the coordinate cutoff satisfy
\[
 \mathbf z_N(-H,w)=\mathscr R\mathbf z_N(H,w),\qquad
 \Xi_N(\mathscr R\mathbf z)=\Xi_N(\mathbf z).
\]
Writing $D_r=\operatorname{diag}(1,-1,\ldots,(-1)^{r-1})$, we have
\begin{align*}
 \mathsf G_r(\mathscr R\mathbf z)&=D_r\mathsf G_r(\mathbf z)D_r,\quad
 \mathsf S_r(\mathscr R\mathbf z)=-D_r\mathsf S_r(\mathbf z)D_r,\\
 L_r(\mathscr R\mathbf z)&=D_rL_r(\mathbf z)D_r,\quad
 \mathsf J_r(\mathscr R\mathbf z)=-D_r\mathsf J_r(\mathbf z)D_r.
\end{align*}
Consequently $\Gamma_r$ is unchanged, $a_k$ changes sign, $b_k$ is unchanged, and $Y^\sigma(t;\mathscr R\mathbf z)=Y^{-\sigma}(t;\mathbf z)$. For spectral localization,
$C_{r,(\sigma,E)}(\mathscr R\mathbf z)=C_{r,(-\sigma,E)}(\mathbf z)$.

\end{proposition}
\begin{proof}
For every $n\ge0$, functional calculus gives
\[
 h_d(-H)=-h_d(H),\qquad
 p_n(h_d(-H))=(-1)^np_n(h_d(H)),\qquad
 \Theta_{\sigma,E}(-H)=\Theta_{-\sigma,E}(H).
\]
Taking normalized traces and matrix elements at $w$, including those with $\Theta_{\sigma,E}$, gives the coordinate transformation \eqref{W:eq:moment-coordinates-reflection}. For
$1\le n\le M$, substitute in $s_n=u_n+x_n/m_N$ and $q_{a,n}=v_{a,n}+y_{a,n}/m_N$. The zeroth coordinates instead satisfy $s_0=1$ and
\[
 q_{(\sigma,E),0}(\mathscr R\mathbf z)
 =t_{-\sigma,E}/N+y_{(-\sigma,E),0}/m_N
 =q_{(-\sigma,E),0}(\mathbf z).
\]
Thus, for every $0\le n\le M$,
\[
 s_n(\mathscr R\mathbf z)=(-1)^ns_n(\mathbf z),\qquad
 q_{(\sigma,E),n}(\mathscr R\mathbf z)
            =(-1)^nq_{(-\sigma,E),n}(\mathbf z).
\]
For indices $0\le i,j<r$, the formulas for the Gram matrix and the matrix representing multiplication by the variable give
\[
 \begin{aligned}
 \mathsf G_{r,ij}(\mathscr R\mathbf z)
   &=\sum_{a=0}^{\min(i,j)}(-1)^{i+j-2a}s_{i+j-2a}
     =(-1)^{i+j}\mathsf G_{r,ij}(\mathbf z),\\
 \mathsf S_{r,ij}(\mathscr R\mathbf z)
   &=\sum_{a=0}^{\min(i,j)}\sum_{\epsilon=\pm1}
               (-1)^{i+j-2a+\epsilon}s_{i+j-2a+\epsilon}
     =(-1)^{i+j+1}\mathsf S_{r,ij}(\mathbf z).
 \end{aligned}
\]
The summands with index $-1$ are zero. These are precisely the stated conjugation formulas with $D_r$.

Since $D_r^2=I$, functional calculus gives $C_r(\mathscr R\mathbf z)=D_rC_r(\mathbf z)D_r$ and
$\gamma_\kappa(\mathsf G_r(\mathscr R\mathbf z)) =D_r\gamma_\kappa(\mathsf G_r(\mathbf z))D_r$. Moreover,
\[
 (D_rL_rD_r)(D_rL_rD_r)^{\mathsf T}=D_rC_rD_r.
\]
The matrix $D_rL_rD_r$ is lower triangular and has the positive diagonal of $L_r$. Uniqueness of the Cholesky factor proves its transformation. Substitution in the definition of
$\mathsf J_r$ gives
\[
 \mathsf J_r(\mathscr R\mathbf z)
 =(D_rL_r^{-1}D_r)(-D_r\mathsf S_rD_r)(D_rL_r^{-\mathsf T}D_r)
 =-D_r\mathsf J_r(\mathbf z)D_r.
\]
Its diagonal entries change sign, and its $(k,k+1)$ entry is multiplied by $-(-1)^{k-1}(-1)^k=1$. Thus the node formula gives
\[
 Y_k^\sigma(\mathscr R\mathbf z)
 =m_N\left(2k+\sigma\sum_{j=1}^ka_j-1
                    -2\sum_{j<k}b_j-b_k\right)
 =Y_k^{-\sigma}(\mathbf z).
\]
The same equality holds after linear interpolation.

The reconstruction of $\mathsf B_{r,a}$ and the transformation of $q_{a,n}$ give, for $0\le i,j<r$,
\[
 \mathsf B_{r,(\sigma,E),ij}(\mathscr R\mathbf z)
 =\sum_{l=0}^{\min(i,j)}(-1)^{i+j-2l}q_{(-\sigma,E),i+j-2l}(\mathbf z)
 =(-1)^{i+j}\mathsf B_{r,(-\sigma,E),ij}(\mathbf z).
\]
Hence $\mathsf B_{r,(\sigma,E)}(\mathscr R\mathbf z)=D_r\mathsf B_{r,(-\sigma,E)}(\mathbf z)D_r$ for every real coordinate vector. Therefore
\[
 \begin{aligned}
 C_{r,(\sigma,E)}(\mathscr R\mathbf z)
 &=t_{-\sigma,E}-\operatorname{Tr}
       \bigl(D_ra_\kappa(\mathsf G_r)D_r
                    D_r\mathsf B_{r,(-\sigma,E)}D_r\bigr)\\
 &=t_{-\sigma,E}-\operatorname{Tr}
       \bigl(a_\kappa(\mathsf G_r)\mathsf B_{r,(-\sigma,E)}\bigr)
   =C_{r,(-\sigma,E)}(\mathbf z).
 \end{aligned}
\]
Taking determinants in the conjugation formula for $\gamma_\kappa$ proves invariance of $\Gamma_r$. Finally, $\mathscr R$ is a signed permutation of the coordinates of $\mathbf z$.
The paired coordinates have equal radii, and $\chi$ is even. Thus every factor of $\Xi_N(\mathscr R\mathbf z)=\prod_j\chi((\mathscr R\mathbf z)_j/\beta_j)$ can be paired with the
corresponding factor of $\Xi_N(\mathbf z)$, proving its invariance.
\end{proof}

For the matrix derivatives of these statistics, use $Q_w=ww^*-N^{-1}I$ as in Section~\ref{W:sec:spectral-regularization-and-matrix-derivatives}. For a finite set $S$ of labels for
fixed self-adjoint matrix increments $W_j$, put
\[
 f_n(\lambda)=p_n(h_d(\lambda)),\qquad
 f_{a,n}(\lambda)=p_n(h_d(\lambda))
                         \theta_E({m_N}^2(2-\sigma\lambda)).
\]
For nonempty $S$ the exact moment derivatives are
\begin{align}
 D_{H,S}u_n&=N^{-1}\operatorname{Tr}D_{H}^{|S|}f_n(H)[W_j:j\in S],\nonumber\\
 D_{H,S}x_n&={m_N}\operatorname{Tr}Q_wD_{H}^{|S|}f_n(H)[W_j:j\in S],\nonumber\\
 D_{H,S}t_a&=\operatorname{Tr}D_{H}^{|S|}f_{a,0}(H)[W_j:j\in S],\nonumber\\
 D_{H,S}v_{a,n}&=N^{-1}\operatorname{Tr}D_{H}^{|S|}f_{a,n}(H)[W_j:j\in S],\nonumber\\
 D_{H,S}y_{a,n}&={m_N}\operatorname{Tr}Q_wD_{H}^{|S|}f_{a,n}(H)[W_j:j\in S].
 \label{W:eq:moment-matrix-derivatives}
\end{align}
The vector $w$ is held fixed. If $W_j=N^{-1/2}T_{\alpha,j}$, each of these derivatives already contains its normalization of the matrix entry directions.

For any one of the functions of the spectral statistics, the full mixed rule is
\begin{equation}
 D_{H,S}[\Phi_N^\nu(\mathbf z_N(H,w))]
 =\sum_{\pi\in\Pi(S)}
 D_{\mathbf z}^{|\pi|}\Phi_N^\nu(\mathbf z_N(H,w))
       [D_{H,B}\mathbf z_N(H,w):B\in\pi].
 \label{W:eq:moment-coordinates-chain}
\end{equation}
This is the composition formula proved in Section~\ref{W:sec:notation-and-scales}, with $s$ replaced by $\mathbf z_N$ and $A$ by $\Phi_N^\nu$. Equations
\eqref{W:eq:moment-matrix-derivatives}--\eqref{W:eq:moment-coordinates-chain} apply through order five at every matrix $H$. The derivatives $D_{\mathbf z}^q\Phi_N^\nu$ are also
defined at arbitrary real coordinate vectors. Taylor's formula applies to every line segment in $\mathbb R^{D_N}$, since the derivatives are defined on that whole space.

The smooth cutoff $0\le U_N(H,w)\le1$ for the Schur complement estimates can be included by appending the coordinate $u_{\rm det}=U_N(H,w)$ and defining
\[
 \Phi_N^{U,\nu}(\mathbf z,u_{\rm det})
   =u_{\rm det}\chi(u_{\rm det})\Phi_N^\nu(\mathbf z).
\]
This compactly supported smooth function of the spectral statistics evaluates to $U_N\Phi_N^\nu(\mathbf z_N(H,w))$. Its modulus is at most three on the full coordinate space and at
most one on moment vectors arising from matrices. Equivalently one can keep the auxiliary cutoff function outside the moment map and use the exact rule
\begin{equation}
 D_{H,S}[U_N\Phi_N^\nu(\mathbf z_N)]
 =\sum_{A\subseteq S}(D_{H,A}U_N)
    \sum_{\pi\in\Pi(S\setminus A)}
       D_{\mathbf z}^{|\pi|}\Phi_N^\nu(\mathbf z_N)[D_{H,B}\mathbf z_N:B\in\pi].
 \label{W:eq:auxiliary-cutoff-product-rule}
\end{equation}
For an empty complement the inner expression is $\Phi_N^\nu$. All derivatives of the cutoff error remain in this formula.

The extensions above have polynomially bounded derivatives at arbitrary moment vectors and, after composition, at arbitrary Hermitian matrices.
Sections~\ref{W:sec:averaged-resolvent-and-conditional-row-estimates}--\ref{W:sec:contours-for-lanczos-comparison} estimate the moment increments under entry replacement.
Section~\ref{W:sec:removing-the-coefficient-cutoff} combines these estimates with the smooth function $\chi_N$ to obtain a summable comparison error.

\subsubsection{Derivatives of the extended functions}
\label{W:sec:critical-derivative-estimates}

The extensions just constructed can now be differentiated without restricting their arguments to moments of matrices. The following bounds are expressed in the norms of
$\dot G_\rho[v]$. Bounds on individual coordinate derivatives have an additional dependence on $\rho$. The examples following the next corollary show explicitly how that dependence
arises.

We use the extensions of Section~\ref{W:sec:cutoffs-in-moment-variables}, with Gram matrix dimension $\rho$ and potential indices $K\le\rho-2$. The variation of the Gram matrix and
the variation of the coordinate cutoff are measured together by

\[
 \mathfrak n(v)=\|\dot G_\rho[v]\|_1+
              \sum_j|v_j|/\beta_j,
\]
where $\beta_j$ are the radii of the coordinate box in $\mathbf z$. In this corollary, derivatives act on $\mathbf z$; the integrated potential is composed with the linear
reconstruction $s_n=u_n+x_n/m_N$.

\begin{corollary}
\label{W:cor:projection-fulltests}
For the function $\Phi^{\mathrm F}$ defined in \eqref{W:eq:moment-coordinates-testsF}, with fixed $\sum_j|\lambda_j|\le\Lambda$, at every real vector of moment coordinates,
\begin{equation}
 |D_{\mathbf z}^q\Phi^{\mathrm F}[v_1,\ldots,v_q]|
 \le C_{q,\kappa,\Lambda}(1+{m_N})^q\prod_j\mathfrak n(v_j),\qquad q\le3.
 \label{W:eq:projection-Fbound}
\end{equation}
For the function $\Phi^{\mathrm M}$ of the modulus of continuity, with fixed threshold $u>0$ and $D$ pairs of sampling points, one has
\begin{equation}
 |D_{\mathbf z}^q\Phi^{\mathrm M}[v_1,\ldots,v_q]|
 \le C_{q,\kappa,u}(1+{m_N})^q
       (1+\log(1+2D))^{q-1}\prod_j\mathfrak n(v_j),\qquad q\le3.
 \label{W:eq:projection-Mbound}
\end{equation}

\end{corollary}
\begin{proof}
Write $Y_\lambda=\sum_j\lambda_jY^\sigma(t_j)$ and $F=e^{\mathrm iY_\lambda}$. The exact derivatives are
\begin{align*}
 D_{\mathbf z}F[v]&=\mathrm iF\,D_{\mathbf z}Y_\lambda[v],\\
 D_{\mathbf z}^2F[v,w]&=F\{\mathrm iD_{\mathbf z}^2Y_\lambda[v,w]-D_{\mathbf z}Y_\lambda[v]D_{\mathbf z}Y_\lambda[w]\},\\
 D_{\mathbf z}^3F[v,w,z]&=F\{\mathrm iD_{\mathbf z}^3Y_\lambda[v,w,z]
    -D_{\mathbf z}^2Y_\lambda[v,w]D_{\mathbf z}Y_\lambda[z]
    -D_{\mathbf z}^2Y_\lambda[v,z]D_{\mathbf z}Y_\lambda[w]\\
 &\hspace{20mm}-D_{\mathbf z}^2Y_\lambda[w,z]D_{\mathbf z}Y_\lambda[v]
    -\mathrm iD_{\mathbf z}Y_\lambda[v]D_{\mathbf z}Y_\lambda[w]D_{\mathbf z}Y_\lambda[z]\}.
\end{align*}
Let $S$ have $q\le3$ labels. The full product rule reads
\[
 D_S(\Xi_N\Gamma_\rho F)
 =\sum_{S_0\sqcup S_1\sqcup S_2=S}
          (D_{S_0}\Xi_N)(D_{S_1}\Gamma_\rho)(D_{S_2}F).
\]
For $a=|S_2|$, the Fourier derivatives and \eqref{W:eq:projection-Ybound} give
\[
 |D_{S_2}F|\le C_{a,\kappa,\Lambda}(1+m_N)^a
          \prod_{j\in S_2}\|\dot G_\rho[v_j]\|_{\rm HS}.
\]
Lemma~\ref{W:lem:gram-cutoff-trace-derivatives} bounds the other factors as
\[
 |D_{S_0}\Xi_N|\le C_q\prod_{j\in S_0}\mathfrak b(v_j),\qquad
 |D_{S_1}\Gamma_\rho|\le C_{q,\kappa}
                       \prod_{j\in S_1}\|\dot G_\rho[v_j]\|_1.
\]
Since $\|\dot G\|_{\rm HS}\le\|\dot G\|_1\le\mathfrak n(v)$, each of the $3^q$ assignments is bounded by the right side of \eqref{W:eq:projection-Fbound}, up to a constant
depending only on $q$. Indeed, a nonzero derivative of $\Gamma_\rho$ requires the spectrum of the Gram matrix to lie in $[\kappa,\kappa^{-1}]$; otherwise the scalar cutoff vanishes
near one of its eigenvalues and $\Gamma_\rho$ is identically zero in a neighbourhood of the moment vector. In that set of moment vectors all regularizations agree locally with the
unregularized moment construction. Thus the integrated potential formula and its derivatives used above coincide with the derivatives of the extension in the moment vector. This
proves the estimate for every real coordinate vector.

For $\mathfrak M$, Lemma~\ref{W:lem:moment-coordinates-modulus} gives
\[
 |D_{y}^j\mathfrak M(y)[u_1,\ldots,u_j]|
 \le C_ju^{-j}(1+\log(1+2D))^{j-1}\prod_l\|u_l\|_\infty.
\]
For a label set $S$, the composed derivative is
\[
 D_S\mathfrak M(\mathbf Y)
 =\sum_{\pi\in\Pi(S)}
 D^{|\pi|}\mathfrak M(\mathbf Y)[D_B\mathbf Y:B\in\pi].
\]
If $|S|=r$ and $|\pi|=j$, \eqref{W:eq:projection-Ybound} bounds its summand by
\[
 C_{r,\kappa,u}m_N^j(1+\log(1+2D))^{j-1}
       \prod_{a\in S}\|\dot G_\rho[v_a]\|_{\rm HS}.
\]
For the same coordinate cutoff $\Xi_N$, the full product is
\[
 D_S(\Gamma_\rho\Xi_N\mathfrak M(\mathbf Y))
 =\sum_{S_0\sqcup S_1\sqcup S_2=S}
   (D_{S_0}\Gamma_\rho)(D_{S_1}\Xi_N)D_{S_2}\mathfrak M(\mathbf Y).
\]
The first two factors are bounded respectively by $C_r\prod_{a\in S_0}\|\dot G_\rho[v_a]\|_1$ and $C_r\prod_{a\in S_1}\sum_l|(v_a)_l|/\beta_l$. Use $j\le|S_2|\le r$, $m_N\ge1$ and
$\|\dot G\|_{\rm HS}\le\|\dot G\|_1$ to obtain \eqref{W:eq:projection-Mbound} after summing the finitely many labelled partitions. Empty label sets mean undifferentiated bounded
factors. Reflection only changes the sign in \eqref{W:eq:projection-Mdef}; all estimates hold for both signs and for exponentials of finite linear combinations of the integrated
potentials.
\end{proof}

\Needspace{6\baselineskip} A moment perturbation in one direction at the free Jacobi matrix shows why the integrated potential's scaling factor cannot be removed from these
estimates.

For $K\ge1$ and Gram matrix dimension $\rho=K+2$, consider the formal moment line
\[
 s_0=1,\qquad s_{2K}=t,\qquad s_n=0\quad(n\ne0,2K).
\]
Near $t=0$ it is positive and is represented, to every required degree, by the measure $(1+tp_{2K}(x))\mu_{\rm sc}(\mathrm dx)$ for $|t|<(2K+1)^{-1}$.
\begin{proposition}
\label{W:prop:projection-sharp}
The integrated potential is
\begin{equation}
 Y_K^+(t)=Y_K^-(t)={m_N}(1-\sqrt{1+t}).
 \label{W:eq:projection-freecurve}
\end{equation}
At this free Jacobi matrix the exponential multiplied by its cutoffs has, for a fixed nonzero frequency $\lambda$,
\begin{align}
 F'(0)&=-\mathrm i\lambda {m_N}/2,\nonumber\\
 F''(0)&=\mathrm i\lambda {m_N}/4-\lambda^2{m_N}^2/4,\nonumber\\
 F'''(0)&=-3\mathrm i\lambda {m_N}/8+3\lambda^2{m_N}^2/8
                              +\mathrm i\lambda^3{m_N}^3/8.
 \label{W:eq:free-integrated-potential-fourier-derivatives}
\end{align}
The variation of the Gram matrix has Hilbert--Schmidt norm $2$ and trace norm $1+\sqrt5$, independently of $K$.
\end{proposition}
\begin{proof}
The Chebyshev product identity \eqref{W:eq:moment-coordinates-Chebyshev} gives $\mathcal L_t(p_ip_j)=\delta_{ij}+t\, [p_{2K}](p_ip_j)$, where $[p_{2K}]$ denotes its coefficient in
the Chebyshev basis. In particular, for $0\le i,j<K$ the product has degree at most $2K-2$, so its $t$ coefficient vanishes. Thus $G_K(t)=I_K$ and all its first $K$ monic squared
norms except the degree-$K$ one equal one. The degree-$K$ polynomial has zero inner product with every polynomial of degree below $K$, because the products involved have degree at
most $2K-1$. Its squared norm is
\[
 \mathcal L_t(p_K^2)=1+t,
\]
using $p_K^2=p_{2K}+p_{2K-2}+\cdots+p_0$. The multiplication coefficients $a_j$ for $j\le K$ vanish: $xp_{j-1}^2$ has degree at most $2K-1$ and therefore its inner product is
unchanged from the semicircle measure. Hence
\[
 b_j=1\ (j<K),\qquad b_K^2=1+t,
 \qquad \mathcal T_K^\sigma(t)=\sqrt{1+t}-1,
\]
for both $\sigma=\pm1$, proving \eqref{W:eq:projection-freecurve}. The measure $(1+tp_{2K})\mu_{\rm sc}$ realizes the prescribed finite moments whenever $|t|<(2K+1)^{-1}$: on
$[-2,2]$, $|p_{2K}|\le2K+1$, so its density is positive. This also shows that all relevant Gram matrices remain positive on a neighbourhood of $t=0$.

For the realization by a finite matrix at $t=0$, take $\rho=K+2$ and the tridiagonal matrix of order $\rho$ with zero diagonal and unit neighbouring entries. Its eigenvalues are
$x_l=2\cos(\pi l/(\rho+1))$ and the squared first coordinates of its normalized eigenvectors are $\omega_l=2\sin^2(\pi l/(\rho+1))/(\rho+1)$. Discrete sine orthogonality gives
\[
 \sum_{l=1}^{\rho}\omega_l p_i(x_l)p_j(x_l)
 =\frac2{\rho+1}\sum_{l=1}^{\rho}
       \sin((i+1)\theta_l)\sin((j+1)\theta_l)
 =\delta_{ij}\quad(0\le i,j<\rho).
\]
Embedding the diagonal matrix of the $x_l$ into $N\ge\rho$ coordinates, with starting vector having entries $\sqrt{\omega_l}$ in those coordinates, gives $G_\rho=I_\rho$, $\|H\|<2$
and $h_d(H)=H$. Both cutoffs are identically one on a neighbourhood of these coordinates, since $\kappa<1/2$ and every eigenvalue of the Gram matrix equals one. Their derivatives
along $t$ consequently vanish at zero.

The two matrix norms and the derivatives of $F$ can be calculated directly. Under the variation $\dot s_{2K}=1$, all other $\dot s_j=0$, the product identity shows that the only
nonzero entries of $\dot G_{K+2}$ lie at the pairs $(K,K+2)$, $(K+2,K)$, $(K+1,K+1)$, and $(K+2,K+2)$ when the rows and columns are indexed from one. In the corresponding block of
order three this is
\[
 D=\begin{pmatrix}0&0&1\\0&1&0\\1&0&1\end{pmatrix}.
\]
Its characteristic polynomial is $(\lambda-1)(\lambda^2-\lambda-1)$, and therefore its eigenvalues are $1,(1+\sqrt5)/2,(1-\sqrt5)/2$. Consequently $\|D\|_{\rm HS}^2=4$ and
$\|D\|_1=1+\sqrt5$, with no dependence on $K$.

Set $y(t)=m_N(1-\sqrt{1+t})$ and $F(t)=\exp(\mathrm i\lambda y(t))$. Direct differentiation gives
\[
 y(0)=0,\quad y'(0)=-m_N/2,\quad
 y''(0)=m_N/4,\quad y'''(0)=-3m_N/8,
\]
and
\[
 \begin{aligned}
 F'&=\mathrm i\lambda y'F,\\
 F''&=(\mathrm i\lambda y''-\lambda^2(y')^2)F,\\
 F'''&=(\mathrm i\lambda y'''-3\lambda^2y'y''
                         -\mathrm i\lambda^3(y')^3)F.
 \end{aligned}
\]
At the moment vector of the free Jacobi matrix the cutoffs are constant one, so these are the derivatives of the full function, including its cutoffs. Substitution of the four
values of $y$ gives exactly \eqref{W:eq:free-integrated-potential-fourier-derivatives}. The calculation concerns the path $\dot x_{2K}=m_N$, $\dot u_n=0$, with all other
coordinates fixed, in the domain $\mathbb R^{D_N}$. Indeed, $s_n=u_n+x_n/m_N$ gives $\dot s_{2K}=1$ and all other $\dot s_n=0$.
\end{proof}

The preceding example shows why the factor ${m_N}$ in the integrated potential cannot be dropped. Section~\ref{W:sec:removing-the-coefficient-cutoff} uses coefficient
regularization inside the observable to obtain bounds on the sum of the absolute expectations of the Taylor remainders over matrix entries.

The next example measures the sum of absolute coordinate derivatives at $K\asymp N^{1/3}$. It motivates estimating $D_sA[v]$ as a single linear functional of $v$, using the contour
norm introduced below, rather than summing bounds for its coordinate coefficients.

\begin{proposition}
\label{W:prop:observable-Fourier}
Fix integers $K,\rho,N$ with $1\le K$ and $K+2\le\rho<N$. Take the function $\Phi^{\mathrm F}$ in \eqref{W:eq:moment-coordinates-testsF} for the upper edge integrated potential at
the single time $Kh_N$, with frequency $\lambda\in\mathbb R\setminus\{0\}$ and largest Gram matrix dimension $\rho_*=\rho$. For every unit vector $w\in\mathbb F_\beta^N$, there is
a Hermitian $N\times N$ matrix $H$ with $\|H\|<2$ and $\mathsf G_\rho=I_\rho$ for which
\begin{equation}
 \sum_{n=1}^{2K}\left|
       \frac{\partial\Phi_N^{\mathrm F}}{\partial x_n}
       (\mathbf z_N(H,w))\right|
 =|\lambda|\left(2K-\frac12\right).
 \label{W:eq:observable-Fourier-l1}
\end{equation}
In particular this quantity is of order $N^{1/3}$ when $K\asymp {m_N}$.

\end{proposition}
\begin{proof}
Choose the number of spectral nodes to be exactly
\begin{equation}
 L=\rho+1\le N,\qquad
 \theta_j=\frac{j\pi}{L+1},\qquad
 \xi_j=2\cos\theta_j,\qquad
 \omega_j=\frac{2}{L+1}\sin^2\theta_j,
 \quad 1\le j\le L.
 \label{W:eq:observable-finite-quadrature}
\end{equation}
Here $L$ denotes the number of quadrature nodes. Let $J_L^{\rm free}$ have zero diagonal and unit off-diagonal coefficients, and set $S_{rj}=\sqrt{2/(L+1)}\sin(r\theta_j)$ for
$1\le r,j\le L$. The sine recurrence, including its zero values at $r=0,L+1$, gives $J_L^{\rm free}S_{\cdot j}=\xi_jS_{\cdot j}$. The $\xi_j$ are distinct. Moreover, summing the
geometric series with ratio $\exp(2\mathrm i\theta_j)$ gives $\sum_{r=1}^L\sin^2(r\theta_j)=(L+1)/2$. Thus the columns are normalized eigenvectors of a real symmetric matrix at
distinct eigenvalues, so $S$ is orthogonal. In particular,
\begin{equation}
 \frac{2}{L+1}\sum_{j=1}^L
       \sin((i+1)\theta_j)\sin((k+1)\theta_j)
 =\delta_{ik},\qquad 0\le i,k\le L-1.
 \label{W:eq:observable-discrete-sine}
\end{equation}
The Chebyshev identity $p_n(\xi_j)=\sin((n+1)\theta_j)/\sin\theta_j$ therefore gives
\begin{equation}
 \sum_{j=1}^L\omega_jp_i(\xi_j)p_k(\xi_j)=\delta_{ik},
 \qquad 0\le i,k\le L-1.
 \label{W:eq:observable-quadrature-Gram}
\end{equation}
Taking $i=k=0$ proves $\sum_j\omega_j=1$. The recurrence $xp_k=p_{k+1}+p_{k-1}$, with $p_{-1}=0$, also gives
\begin{equation}
 \sum_{j=1}^L\omega_j\xi_jp_i(\xi_j)p_k(\xi_j)
 =\delta_{i,k+1}+\delta_{i,k-1},
 \qquad 0\le i,k<\rho.
 \label{W:eq:observable-quadrature-multiplication}
\end{equation}
Indeed $k+1\le\rho=L-1$, so every polynomial in this computation lies within the orthogonality range in \eqref{W:eq:observable-quadrature-Gram}. Hence the Gram matrix of size
$\rho$ is $I_\rho$, its Cholesky factor is $I_\rho$, and the multiplication matrix is $J_\rho^{\rm free}$. All required coefficients are therefore $a_j=0$ and $b_j=1$ for every
$j\le K$.

For the embedding into the ambient dimension, put
\begin{equation}
 D=\operatorname{diag}(\xi_1,\ldots,\xi_L,
                          \underbrace{0,\ldots,0}_{N-L}),\qquad
 w_0=(\sqrt{\omega_1},\ldots,\sqrt{\omega_L},0,\ldots,0)^{\mathsf T}.
 \label{W:eq:observable-quadrature-embedding}
\end{equation}
When $L=N$ there are no additional coordinates. Choose a deterministic orthogonal or unitary $V$ with $Vw_0=w$, and set $H=VDV^*$. Then $w^*f(H)w=\sum_j\omega_j f(\xi_j)$ for every
function on the spectrum, and $\|H\|=2\cos(\pi/(L+1))<2$. The added coordinates have zero weight in the spectral measure at $w$. Also $h_d(H)=H$, so the modified and original
polynomial moments agree. By \eqref{W:eq:moment-coordinate-bounds}--\eqref{W:eq:moment-coordinates-box}, the coordinate cutoff $\Xi_N$ is identically one on a neighbourhood of
$\mathbf z_N(H,w)$. Since $0<\kappa<1/4$ and $\mathsf G_\rho=I_\rho$, the function $\Gamma_\rho$ is also identically one on a neighbourhood of this moment vector.

Near the moment vector chosen above the scalar regularizations agree with the unmodified Gram matrix and the Jacobi matrix obtained by orthogonalization. Both cutoffs are
identically one on an open neighbourhood, so every positive-order derivative of either cutoff vanishes at the chosen point. We now differentiate with respect to the stated
coordinates rather than the matrix entries.

At $G_\rho=I_\rho$, Lemma~\ref{W:lem:projection-variation} and the Chebyshev product identity give
\[
 D_sY_K^+[v]=-m_N\left(\sum_{j=0}^{K-1}v(p_{2j+1})
             +\sum_{j=1}^{K-1}v(p_{2j})+\frac12v(p_{2K})\right).
\]
Since $v(p_0)=0$, the first two sums together contain each index $1\le n\le2K-1$ exactly once: their indices are respectively odd and even. Therefore
\[
 D_sY_K^+[v]=-m_N\left(\sum_{n=1}^{2K-1}\delta s_n
                     +\frac12\delta s_{2K}\right),
 \qquad Y_K^+=0.
\]
For a partial derivative in $x_n$, hold every $u_j$, every other $x_j$, and all coordinates defined by traces fixed. The defining reconstruction $s_j=u_j+x_j/m_N$ then gives
$\partial_{x_n}s_j=\delta_{jn}/m_N$, so
\[
 \partial_{x_n}Y_K^+
 =-\mathbf1_{\{1\le n<2K\}}-\frac12\mathbf1_{\{n=2K\}}.
\]
The function, including its cutoffs, is $\Xi_N\Gamma_\rho\exp(\mathrm i\lambda Y_K^+)$. At the moment vector of the free Jacobi matrix both cutoffs equal one and have zero first
derivative, and the exponential equals one. Hence
\begin{equation}
 \frac{\partial\Phi_N^{\mathrm F}}{\partial x_n}
       (\mathbf z_N(H,w))
 =\begin{cases}
   -\mathrm i\lambda,&1\le n<2K,\\
   -\mathrm i\lambda/2,&n=2K,\\
   0,&2K<n\le 2\rho-1.
  \end{cases}
 \label{W:eq:observable-free-coordinate-derivative}
\end{equation}
There are $2K-1$ derivatives of modulus $|\lambda|$ and one of modulus $|\lambda|/2$, so their sum is $|\lambda|(2K-1/2)$, as required. These are partial derivatives on the full
coordinate space, not derivatives of an entry distribution. Summing the absolute values proves \eqref{W:eq:observable-Fourier-l1}.
\end{proof}

The second variation of the nonlinear remainder is \eqref{W:eq:second-derivative-of-the-bond-remainder}. It remains in higher contractions even at the moment functional of the
semicircle measure.

The preceding calculation concerns individual moment derivatives. When comparing entry distributions, we estimate the polynomial representing the derivative, using the row and
column coefficient norms and the variation of the off-diagonal coefficients. The cutoff on these quantities is included in the observable; its complement is estimated by a separate
nonnegative function. Sections~\ref{W:sec:taylor-remainders-and-bounded-entries} and~\ref{W:sec:cumulant-cubic} then bound the remainder of the expansion and the signed cubic term.

\subsubsection{Contour representations of variations of the moments}
\label{W:sec:derivatives-of-the-integrated-potential-and-coefficient-localization}

The preceding examples explain why estimating every moment coordinate separately is costly. Entry variations have more structure: they admit Cauchy integral representations, which
can be paired directly with the derivative bounds of Section~\ref{W:sec:derivative-growth-at-the-critical-scale}. The norm introduced below measures these representations and
dominates $\|\cdot\|_\diamond$.

The dimensions satisfy $c_1m_N\le\rho\le\rho_*\le C_*m_N$, and $d$ is given by \eqref{W:eq:moment-coordinates-scales}. The functionals have mass one and their variations have mass
zero. For the fixed polynomials of Section~\ref{W:sec:the-gram-representation-of-the-jacobi-coefficients}, the contour is the ellipse of
Section~\ref{W:sec:contours-for-lanczos-comparison}:

\[
 z=\zeta+\zeta^{-1},\quad |\zeta|=\mathrm e^s,
 \qquad s=4\max\{\rho^{-1},\sqrt d\},\qquad z\in\mathcal C.
\]
Define
\begin{equation}
 E_\rho(z)=\frac{2\mathrm e^{\rho s}}{|z^2-4|^{1/2}},\qquad
 \nu_\rho(\mathrm dz)=E_\rho(z)^2\frac{|\mathrm dz|}{2\pi}.
 \label{W:eq:comparison-measure}
\end{equation}
The finite Chebyshev formula gives $|p_j(z)|\le E_\rho(z)$ for $j<\rho$. With $\zeta=e^{s+\mathrm i\theta}$, $\zeta-\zeta^{-1}=2(\sinh s\cos\theta+ \mathrm i\cosh s\sin\theta)$,
and hence
\[
 |\zeta-\zeta^{-1}|^2=4(\sinh^2s+\sin^2\theta),\qquad
 \int_{\mathcal C}\frac{|\mathrm dz|}{|z^2-4|}
 =\frac12\int_0^{2\pi}
                \frac{\mathrm d\theta}{(\sinh^2s+\sin^2\theta)^{1/2}}.
\]
For $0<s\le1$, symmetry and $\sin\theta\ge2\theta/\pi$ on $[0,\pi/2]$ give
\[
 \int_{\mathcal C}\frac{|\mathrm dz|}{|z^2-4|}
 \le C\int_0^{\pi/2}(s^2+\theta^2)^{-1/2}\,\mathrm d\theta
 =C\operatorname{arsinh}\!\left(\frac\pi{2s}\right)
 \le C(1+|\log s|).
\]
For $s\ge1$ the integral is at most $\pi/\sinh s\le C$.

Since $\rho s=O(1+(\log N)^{\gamma/2})$ and $\gamma<2$, it follows that $\nu_\rho(\mathcal C)=N^{o(1)}$.

For a moment variation $v$, regarded as a linear functional on the polynomials of the prescribed degree, put
\begin{equation}
 \|v\|_{\mathscr D}=\inf_k\int_{\mathcal C}|k(z)|\nu_\rho(\mathrm dz),
 \qquad
 v(f)=\frac{1}{2\pi\mathrm i}\int_{\mathcal C}f(z)k(z)\,\mathrm dz
 \quad\text{for every polynomial in that space}.
 \label{W:eq:comparison-direction-norm}
\end{equation}
The infimum is over measurable kernels with finite displayed integral and the stated moment identities. This set is nonempty. If the degree range is $0,\ldots,M$, the ellipse
encloses zero and the Laurent polynomial
\[
 k_v(z)=\sum_{j=0}^M v(x^j)z^{-j-1}
\]
satisfies, for $0\le n\le M$,
\[
 \frac1{2\pi\mathrm i}\int_{\mathcal C}z^nk_v(z)\,\dd z
 =\sum_{j=0}^Mv(x^j)\mathbf1_{\{j=n\}}=v(x^n).
\]
The kernel is continuous on the contour and therefore integrable. For any polynomial $f$ in the degree range and any representing kernel,
\[
 |v(f)|\le\sup_{z\in\mathcal C}\frac{|f(z)|}{E_\rho(z)^2}
                 \int_{\mathcal C}|k(z)|\nu_\rho(\dd z).
\]
Taking the infimum proves that $\|v\|_{\mathscr D}=0$ forces $v(x^n)=0$ for every $n$. If $k_1,k_2$ represent $v_1,v_2$, then $k_1+k_2$ represents $v_1+v_2$. Choosing their
integrals within $\varepsilon$ of the two infima and letting $\varepsilon\downarrow0$ gives the triangle inequality; scalar multiplication gives absolute homogeneity. Thus
\eqref{W:eq:comparison-direction-norm} is a norm, also on the subspace $v(1)=0$.

Let $G_\rho=(\mathcal L(p_ip_j))_{0\le i,j<\rho}$ and $\dot G_\rho[v]=(v(p_ip_j))_{0\le i,j<\rho}$. Then
\begin{equation}
 \|\dot G_\rho[v]\|_1\le\rho\|v\|_{\mathscr D}.
 \label{W:eq:comparison-Gram-direction}
\end{equation}
For $p(z)=(p_0(z),\ldots,p_{\rho-1}(z))^{\mathsf T}$,
\[
 \dot G_\rho[v]=\frac1{2\pi\mathrm i}\int_{\mathcal C}
                          p(z)p(z)^{\mathsf T}k_v(z)\,\dd z.
\]
A rank-one matrix $uv^*$ has trace norm $\|u\|_2\|v\|_2$. Apply this with $u=p(z)$ and $v=\overline{p(z)}$ to get
\[
 \|\dot G_\rho[v]\|_1
 \le\frac1{2\pi}\int_{\mathcal C}\|p(z)\|_2^2|k_v(z)|\,|\dd z|
 \le\rho\int_{\mathcal C}|k_v|\,\dd\nu_\rho.
\]
Taking the infimum gives \eqref{W:eq:comparison-Gram-direction}. Also, whenever $q_i=\sum_jC_{ij}p_j$ and every coefficient row has absolute sum at most $Q$, we have
\begin{equation}
 |v(q_iq_j)|\le Q^2\|v\|_{\mathscr D}.
 \label{W:eq:comparison-Ebound}
\end{equation}
Both conclusions follow first for a representing kernel and then by taking the infimum.

For the norm $\norm v_{\mathscr D}$ defined in \eqref{W:eq:comparison-direction-norm},
\begin{equation}
 \norm v_\diamond\le\norm v_{\mathscr D}.
 \label{W:eq:tightness-domination}
\end{equation}
Indeed $|p_i(z)p_j(z)|\le E_\rho(z)^2$, and the trace norm of the rank-one matrix $p(z)p(z)^{\mathsf T}$ is at most
$\rho E_\rho(z)^2$.

By $\|v\|_\diamond\le\|v\|_{\mathscr D}$ in \eqref{W:eq:tightness-domination}, Lemmas~\ref{W:lem:polynomial-coefficient-derivatives} and~\ref{W:lem:tightness-directions} give, for
each fixed $r\ge1$,
\begin{align}
 \max_{ij}|D_{s}^rC_{ij}[v_1,\ldots,v_r]|
 &\le C_{r,\kappa}\rho^{r-1}Q_*^{2r+1}
                            \prod_a\|v_a\|_{\mathscr D},\nonumber\\
 \max_j|D_{s}^rb_j[v_1,\ldots,v_r]|
 &\le C_{r,\kappa}\rho^{r-1}Q_*^{2r}
                            \prod_a\|v_a\|_{\mathscr D},
 \label{W:eq:comparison-coeff-derivatives}
\end{align}
for the coefficients determined by the initial principal block of the Gram matrix, and
\begin{equation}
 |D_{s}^r\mathcal T_K^\sigma[v_1,\ldots,v_r]|
 \le C_{r,\kappa}B\rho^{r-1}Q_*^{2r}
                            \prod_a\|v_a\|_{\mathscr D},
 \qquad K\le\rho-2.
 \label{W:eq:comparison-integrated-potential-derivatives}
\end{equation}
The first derivative is the polynomial pairing in \eqref{W:eq:potential-derivative-kernel}, with largest polynomial index $K$.

\subsubsection{Bounds for the polynomial coefficients}
\label{W:sec:coefficient-localization}
\label{W:sec:removing-the-coefficient-cutoff}
\label{W:sec:expectation-of-coefficient-cutoff}

The derivative bounds still contain the coefficient quantities $Q_*$ and $B$. To use them in comparison, these quantities must be restricted by a cutoff whose own derivatives can
be estimated. The bounds for bounded entries and Gaussian entries in Lemma~\ref{W:lem:comparison-ensemble-coefficients} will determine the cutoff thresholds.

Absolute values and maxima in $Q_*$ and $B$ are replaced by smooth approximations. With $\delta=\rho^{-1}$ and $a_\delta(x)=(x^2+\delta^2)^{1/2}$, the smoothed row and column sums
are

\[
 S_{j+1}=\sum_{i=0}^{\rho-1}a_\delta(C_{ji}),\qquad
 S_{\rho+i+1}=\sum_{j=0}^{\rho-1}a_\delta(C_{ji}).
\]
With $\beta_\rho=\log(2\rho)$, set
\begin{equation}
 \mathcal Q=\beta_\rho^{-1}
         \log\left(\sum_{j=1}^{2\rho}\mathrm e^{\beta_\rho S_j}\right),
 \qquad
 \mathcal B=1+\sum_{j=1}^{\rho-3}a_\delta(b_{j+1}-b_j).
 \label{W:eq:comparison-smooth-QB}
\end{equation}
The logarithm of the exponential sum defining $\mathcal Q$ is a smooth approximation to the maximum. Smoothing absolute
values adds at most one to each sum, and the smooth maximum adds at most $\log(2\rho)/\beta_\rho=1$. Since $q_0=1$,
$Q_*\ge1$ is already attained or exceeded by a coefficient row. Hence $Q_*\le\mathcal Q\le Q_*+2$ and $B\le\mathcal B\le
B+1$.

Fix $\vartheta\in C^\infty(\mathbb R;[0,1])$ equal to one on $(-\infty,1]$ and zero on $[2,\infty)$. For a deterministic sequence $\mathfrak r_N\to\infty$,
$\mathfrak r_N=N^{o(1)}$, chosen using Lemma~\ref{W:lem:comparison-ensemble-coefficients}, define
\begin{equation}
 \chi_N=\vartheta(\mathcal Q/\mathfrak r_N)\vartheta(\mathcal B/\mathfrak r_N).
 \label{W:eq:comparison-cutoff}
\end{equation}
A nonzero value or derivative of $\chi_N$ forces $Q_*\le2\mathfrak r_N$ and $B\le2\mathfrak r_N$. Its derivatives under entry replacement are included below.
\begin{lemma}
\label{W:lem:comparison-cutoff-derivatives}
On $\kappa I\le G_\rho\le\kappa^{-1}I$, for $1\le r\le4$, wherever $D_s^r\chi_N$ is nonzero,
\begin{equation}
 |D_{s}^r\chi_N[v_1,\ldots,v_r]|
 \le C_{r,\kappa}(1+\log\rho)^r\rho^r\mathfrak r_N^{3r}
                             \prod_a\|v_a\|_{\diamond}.
 \label{W:eq:comparison-cutoff-derivatives}
\end{equation}
\end{lemma}
\begin{proof}
For $k\ge1$, $\|a_\delta^{(k)}\|_\infty\le C_k\rho^{k-1}$. For a nonempty label set $S$, a smoothed row sum has derivative
\[
 D_SS_j=\sum_i\sum_{\pi\in\Pi(S)}
 a_\delta^{(|\pi|)}(C_{ji})\prod_{A\in\pi}D_AC_{ji}.
\]
For a column sum the formula is
\[
 D_SS_{\rho+i+1}=\sum_j\sum_{\pi\in\Pi(S)}
 a_\delta^{(|\pi|)}(C_{ji})\prod_{A\in\pi}D_AC_{ji}.
\] If
$|S|=r$, $|\pi|=k$, its sum over $i$ is at most
\[
 C_r\rho\,\rho^{k-1}
        \prod_{A\in\pi}(\rho^{|A|-1}Q_*^{2|A|+1})
                 \prod_{a\in S}\|v_a\|_{\diamond}
 =C_r\rho^rQ_*^{2r+k}\prod_{a\in S}\|v_a\|_{\diamond}.
\]
Since $1\le k\le r$ and $Q_*\ge1$, this is bounded by $C_r\rho^rQ_*^{3r}\prod_a\|v_a\|_{\diamond}$. For $d_j=b_{j+1}-b_j$, apply the scalar chain rule to each summand
$a_\delta(d_j)$:
\[
 D_S\mathcal B=\sum_j\sum_{\pi\in\Pi(S)}
        a_\delta^{(|\pi|)}(d_j)\prod_{A\in\pi}D_Ad_j,
 \qquad
 |D_S\mathcal B|\le C_r\rho^rQ_*^{2r}\prod_a\|v_a\|_{\diamond}.
\]
Here each derivative of an off-diagonal coefficient has the factor $\rho^{|A|-1}Q_*^{2|A|}$; summing at most $\rho$ terms gives the stated power, including all mixed derivatives.

For $L(y)=\beta^{-1}\log\sum_j e^{\beta y_j}$, let $\pi_j=e^{\beta y_j}/\sum_i e^{\beta y_i}$. Direct differentiation gives
\[
 D^kL(y)[h_1,\ldots,h_k]
 =\beta^{k-1}\sum_{\tau\in\Pi_k}
 (-1)^{|\tau|-1}(|\tau|-1)!
 \prod_{B\in\tau}\sum_j\pi_j\prod_{a\in B}h_{a,j}.
\]
Each sum over $j$ is bounded by $\prod_{a\in B}\|h_a\|_\infty$. The chain rule for $\mathcal Q=L((S_j)_j)$ therefore gives
\[
 |D_S\mathcal Q|\le C_r(1+\log\rho)^{r-1}
                \rho^rQ_*^{3r}\prod_a\|v_a\|_{\diamond}.
\]
Finally assign the labels as $A\sqcup B=S$ to the two cutoff factors. For example, for nonempty $A$,
\[
 D_A\vartheta(\mathcal Q/\mathfrak r_N)
 =\sum_{\pi\in\Pi(A)}
 \vartheta^{(|\pi|)}(\mathcal Q/\mathfrak r_N)
 \mathfrak r_N^{-|\pi|}\prod_{C\in\pi}D_C\mathcal Q,
\]
The formula for $\mathcal B$ is the same with $\mathcal Q$ replaced by $\mathcal B$; an empty set denotes the undifferentiated cutoff. Every nonzero product forces both
$\mathcal Q,\mathcal B\le2\mathfrak r_N$, since each scalar cutoff and all its derivatives vanish above two. Consequently $Q_*\le2\mathfrak r_N$ in every contributing term. The
product bounds, $\mathfrak r_N\ge1$, and the finite number of label partitions yield \eqref{W:eq:comparison-cutoff-derivatives}, including the transition regions.

\end{proof}

By \eqref{W:eq:tightness-domination}, the same bound holds with $\|v_a\|_{\mathscr D}$ in place of $\|v_a\|_\diamond$.

To remove the factor $\chi_N$, we compare the expectations of both $\chi_NF_N$ and $D_N^{\rm coef}$ under the two matrix laws. The following inequality reduces the expectation of
the cutoff error under the first matrix law to its expectation under the second.

For any test function of the matrix $F=\Gamma_\rho\Phi$ with $|\Phi|\le1$, let $F^{\rm loc}=\chi_NF$ and $D=\Gamma_\rho(1-\chi_N)$. Then $0\le D\le1$ and $|F-F^{\rm loc}|\le D$.
Let $\mathbb E_P,\mathbb E_Q$ denote expectations under two matrix distributions. If the comparisons for $F^{\rm loc}$ and $D$ have errors $\varepsilon_F$ and $\varepsilon_D$,
respectively, then
\begin{equation}
 |\mathbb E_PF-\mathbb E_QF|
 \le \varepsilon_F+\varepsilon_D+2\mathbb E_QD.
 \label{W:eq:removal-of-coefficient-cutoff}
\end{equation}
The two inequalities are
\[
 \begin{aligned}
 |\mathbb E_PF-\mathbb E_QF|
 &\le |\mathbb E_PF^{\rm loc}-\mathbb E_QF^{\rm loc}|
                +\mathbb E_P|F-F^{\rm loc}|+\mathbb E_Q|F-F^{\rm loc}|\\
 &\le\varepsilon_F+\mathbb E_PD+\mathbb E_QD,\\
 \mathbb E_PD&\le\mathbb E_QD+|\mathbb E_PD-\mathbb E_QD|
                    \le\mathbb E_QD+\varepsilon_D.
 \end{aligned}
\]
Substituting the second in the first proves \eqref{W:eq:removal-of-coefficient-cutoff}. It remains to bound $\mathbb E_QD$.

In the comparison of entry distributions we integrate the function, including its cutoffs, against the original entry laws. The cutoff is part of the differentiable function. For
example, if $\mu=(\delta_{-1}+\delta_1)/2$ and $\nu=\delta_0$, then
\[
 \int x\,(\mu-\nu)(\dd x)=0,\qquad
 \int_{(0,\infty)}x\,(\mu-\nu)(\dd x)=\tfrac12.
\]
Thus the cutoff remains in the integrand when moment cancellation is used. Its derivatives are given by the product rule of Section~\ref{W:sec:notation-and-scales}.

\subsubsection{Derivatives of the functions used in comparison}
\label{W:sec:modulus-of-continuity-tests-and-derivative-kernels}
\label{W:sec:localized-test-class}

The two cutoffs have complementary roles: $\Gamma_\rho$ keeps the Gram matrix uniformly invertible, while $\chi_N$ controls the coefficient quantities in the derivative estimate.
Both belong to the observable being compared, so their derivatives must be included even in the regions where a cutoff is changing. The bounds below apply to these complete
products and to their smooth extensions by zero.

The scale and dimensions satisfy ${m_N}=N^{1/3}$ and $c_1m_N\le\rho\le\rho_*\le C_*m_N$. We use the spectral modification $h_d$ and the cutoff in the eigenvalues of the Gram matrix
from Section~\ref{W:sec:the-moment-vector-and-its-smooth-extensions}. For finitely many indices $K_l\le\rho-2$, signs $\sigma_l\in\{-1,1\}$, and real frequencies with
$\sum_l|\lambda_l|\le\Lambda$, the function to be compared is

\begin{equation}
 F_N(H,w)=\Gamma_\rho(H,w)
   \exp\!\left(-\mathrm i {m_N}\sum_l\lambda_l
                  \mathcal T_{K_l}^{\sigma_l}(H,w)\right),
 \qquad \sigma_l\in\{-1,1\}.
 \label{W:eq:comparison-test}
\end{equation}
Convex linear interpolation between adjacent specified nodes is included. The integrated potentials are the exact cumulative sums defined in
Section~\ref{W:sec:a-finite-sum-rule-and-potential-derivatives}. At moment vectors obtained from matrices, this is exactly the exponential multiplied by $\Gamma_\rho$ in
Section~\ref{W:sec:the-moment-vector-and-its-smooth-extensions}.

Define on the set where the Gram matrix is positive definite
\begin{equation}
 F_N^{\rm loc}=\Gamma_\rho\chi_N
           \exp(-\mathrm i {m_N}\textstyle\sum_l\lambda_l
                                     \mathcal T_{K_l}^{\sigma_l}),
 \qquad D_N^{\rm coef}=\Gamma_\rho(1-\chi_N).
 \label{W:eq:comparison-two-tests}
\end{equation}
Both functions are extended by zero outside that region. The extensions are smooth on the full real coordinate space with the constant moment fixed: the cutoff in the eigenvalues
of the Gram matrix vanishes on a neighbourhood of every point where the Gram matrix is singular. On the support of its derivatives the original scalar regularizations agree with
the exact formulas on the set where the Gram matrix is positive definite.

\Needspace{6\baselineskip}

\begin{proposition}
\label{W:prop:comparison-global}
For $A_N=F_N^{\rm loc}$ or $A_N=D_N^{\rm coef}$, $|A_N|\le1$ and
\begin{equation}
 |D_{s}^rA_N(s)[v_1,\ldots,v_r]|
 \le C_r\Lambda_N{m_N}^r\prod_a\|v_a\|_{\mathscr D},
 \quad 1\le r\le4,
 \quad\Lambda_N=C\mathfrak r_N^{20}(1+\log\rho)^4=N^{o(1)}.
 \label{W:eq:comparison-global}
\end{equation}
The estimate holds for every real moment vector with $s_0=1$.
\end{proposition}
\begin{proof}
By \eqref{W:eq:gram-cutoff-trace-derivatives} and the trace norm bound \eqref{W:eq:comparison-Gram-direction},
\[
 |D_{s}^r\Gamma_\rho[v_1,\ldots,v_r]|
 \le C_r\rho^r\prod_a\|v_a\|_{\mathscr D}.
\]
Let $\theta={m_N}\sum_l\lambda_l\mathcal T_{K_l}^{\sigma_l}$. In a nonzero summand of the differentiated product $\Gamma_\rho\chi_N e^{-\mathrm i\theta}$, the factor obtained from
$\Gamma_\rho$ forces $\kappa I\le G_\rho\le\kappa^{-1}I$, and the factor obtained from $\chi_N$ forces $Q_*\le2\mathfrak r_N$ and $B\le2\mathfrak r_N$. The support conditions
include derivatives of order zero. We apply Lemma~\ref{W:lem:tightness-directions} at this same moment vector, using $\|v\|_\diamond\le\|v\|_{\mathscr D}$. For a nonempty subset
$J$ of derivative labels it gives
\[
 |D_{s,J}\theta|
 \le C_{|J|,\kappa,\Lambda}
 {m_N}\rho^{|J|-1}\mathfrak r_N^{2|J|+1}
             \prod_{a\in J}\|v_a\|_{\mathscr D}.
\]
The phase is real whenever the moment vector is real and its Gram matrix is positive definite. Hence
\begin{equation}
 D_{s,J} e^{-\mathrm i\theta}
 =e^{-\mathrm i\theta}
   \sum_{\pi\in\Pi(J)}(-\mathrm i)^{|\pi|}
                         \prod_{B\in\pi}D_{s,B}\theta,
 \label{W:eq:localized-exponential-partitions}
\end{equation}
where $\Pi(J)$ denotes the set of partitions of $J$. If $q=|J|$ and a partition has $k$ blocks, its bound is
\[
 C_q {m_N}^k\rho^{q-k}\mathfrak r_N^{2q+k}
       \prod_{a\in J}\|v_a\|_{\mathscr D}
 \le C_q{m_N}^q\mathfrak r_N^{3q}\prod_{a\in J}\|v_a\|_{\mathscr D}.
\]
This is where $\rho\asymp {m_N}$ is used. A bound with $\rho^{|J|}$ for each derivative of the integrated potential would leave an extra factor $\rho$ for every partition block.

In the full product rule for $\Gamma_\rho\chi_Ne^{-\mathrm i\theta}$, assign respectively $r_0,r_1,r_2$ derivatives to its three factors, with $r_0+r_1+r_2=r$.
Lemma~\ref{W:lem:comparison-cutoff-derivatives} and \eqref{W:eq:localized-exponential-partitions} bound this term by
\[
 C_r(1+\log\rho)^{r_1}
 \rho^{r_0+r_1}{m_N}^{r_2}\mathfrak r_N^{3r_1+3r_2}
 \prod_a\|v_a\|_{\mathscr D}
 \le C_r{m_N}^r(1+\log\rho)^r\mathfrak r_N^{3r}
 \prod_a\|v_a\|_{\mathscr D}.
\]
The bound also holds when some $r_j=0$, since each undifferentiated factor has modulus at most one and an empty product equals one. There are only finitely many terms for $r\le4$,
and $\mathfrak r_N^{20}(1+\log\rho)^4$ dominates them all. This proves \eqref{W:eq:comparison-global} throughout the set where the Gram matrix is positive definite, including the
transition region of each cutoff.

For the function measuring the coefficient cutoff error, write the identity $D_N^{\rm coef}=\Gamma_\rho-\Gamma_\rho\chi_N$. For any labelled set $S$ with $|S|=r\le4$, the
derivative of its second term is
\[
 D_S(\Gamma_\rho\chi_N)
 =\sum_{T\sqcup U=S}(D_T\Gamma_\rho)(D_U\chi_N).
\]
The determinant estimate bounds $D_T\Gamma_\rho$ by $C_r\rho^{|T|}\prod_{a\in T}\|v_a\|_{\mathscr D}$; the coefficient cutoff lemma bounds the other factor by
$C_r\rho^{|U|}\mathfrak r_N^{3|U|}(1+\log\rho)^{|U|}\prod_{a\in U}\|v_a\|_{\mathscr D}$. Their product is at most
\[
 C_r\rho^r\mathfrak r_N^{3r}(1+\log\rho)^r
       \prod_{a\in S}\|v_a\|_{\mathscr D}
 \le C_r\Lambda_N m_N^r\prod_{a\in S}\|v_a\|_{\mathscr D},
\]
since $\rho\asymp m_N$, $\mathfrak r_N\ge1$ and $r\le4$. There are precisely $2^r$ labelled allocations. The undifferentiated term satisfies
$0\le D_N^{\rm coef}=\Gamma_\rho(1-\chi_N)\le1$. Together with the bound for $D_S\Gamma_\rho$, this proves the proposition for both observables, including regions in which $\chi_N$
and all its derivatives vanish.

Finally consider a boundary point of the open set $G_\rho>0$. At least one eigenvalue of $G_\rho$ approaches zero there. The support of the scalar function $\gamma_\kappa$ is
contained in a fixed interval bounded away from zero; hence $\Gamma_\rho=\det\gamma_\kappa(G_\rho)$ is identically zero on a neighbourhood of that boundary point. Both products
$\Gamma_\rho\chi_Ne^{-\mathrm i\theta}$ and $\Gamma_\rho(1-\chi_N)$, along with all their derivatives, vanish on this neighbourhood before their unregularized factors could become
singular. Extending them by zero therefore produces smooth functions with the same bounds on the entire real coordinate space.
\end{proof}

For the modulus of continuity, fix $T,u>0$ and $0<\delta\le T$ before taking $N\to\infty$. Use $K=\lceil m_NT\rceil\le\rho-2$, the pair set $\mathcal P_\delta$ of
Lemma~\ref{W:lem:moment-coordinates-modulus}, and the exact values $Y_j^\sigma=-m_N\mathcal T_j^\sigma$, with $Y_0^\sigma=0$. Write $D=|\mathcal P_\delta|$ and
\begin{equation}
 M(y)=\Phi_{\mathcal P_\delta,u}(y),\qquad
 \widehat F_{M,N}^\sigma(H,w)=\Gamma_\rho(H,w)M(\mathbf Y^\sigma(H,w)),
 \label{W:eq:tightness-test}
\end{equation}
with the same scalar cutoff as in \eqref{W:eq:smooth-modulus-definition} and $M=0$ when $D=0$. Lemma~\ref{W:lem:moment-coordinates-modulus} gives
\begin{equation}
 |D_{y}^qM(y)[h_1,\ldots,h_q]|
 \le C_q u^{-q}(1+\log(1+2D))^{q-1}
                  \prod_a\norm{h_a}_\infty.
 \label{W:eq:tightness-soft}
\end{equation}
The same lemma gives the indicator and interpolation bounds
\begin{equation}
 \mathbf1_{\{M_\delta(y)\ge u\}}\le M(y)
       \le\mathbf1_{\{M_\delta(y)>u/4\}},\qquad
 \omega_T(y,\delta)\le M_\delta(y)\le\omega_{K/{m_N}}(y,\delta+2/{m_N}).
 \label{W:eq:tightness-sandwich}
\end{equation}
The modulus on the right refers to the same linear interpolant on $[0,K/{m_N}]$. If the pair set is empty, the convention $M=0$ gives the zero function directly.

Use the smooth function $\chi_N$ defined above, with $\mathfrak r_N=N^{o(1)}$, $\mathfrak r_N\to\infty$. On the support of it or its derivatives, $Q_*+B\le4\mathfrak r_N$. The
derivative estimate for $\chi_N$ is
\[
 |D_{s}^r\chi_N|\le C_{r,\kappa}\rho^r\mathfrak r_N^{3r}(1+\log\rho)^r
                      \prod_a\norm{v_a}_\diamond,\quad r\le4.
\]
The function $\Gamma_\rho$ satisfies $|D_{s}^r\Gamma_\rho|\le C_r\rho^r\prod\norm{v_a}_\diamond$. Set
\begin{equation}
 A_N^\sigma(s)=\Gamma_\rho(s)\chi_N(s)M(\mathbf Y^\sigma(s)),
 \qquad A_N^0(s)=\Gamma_\rho(s)\chi_N(s).
 \label{W:eq:tightness-localized}
\end{equation}
The exponential multiplied by both cutoff functions is
\begin{equation}
 A_N^{\rm F}(s)=\Gamma_\rho(s)\chi_N(s)
       \exp\!\left(-\mathrm i {m_N}\sum_l\lambda_l
                    \mathcal T_{K_l}^{\sigma_l}(s)\right),
 \qquad \sum_l|\lambda_l|\le\Lambda,
 \label{W:eq:tightness-Fourier-test}
\end{equation}
including the convex combination of adjacent sampled values in \eqref{W:eq:comparison-test}. Extend these functions by zero outside the set where the Gram matrix is positive
definite. Flatness of $\Gamma_\rho$ gives smooth functions of all real moment coordinates with $s_0=1$. The coordinate cutoff $\Xi_N$ equals one on a neighbourhood of every moment
vector arising from a matrix by \eqref{W:eq:moment-coordinate-bounds}--\eqref{W:eq:moment-coordinates-box}. Its matrix derivatives of positive order vanish there.

\Needspace{6\baselineskip}

\begin{lemma}
\label{W:lem:tightness-global}
For $A\in\{A_N^+,A_N^-,A_N^0,A_N^{\rm F}\}$, $|A|\le1$ and
\begin{equation}
 |D_{s}^rA(s)[v_1,\ldots,v_r]|
 \le \Lambda_N {m_N}^r\prod_a\norm{v_a}_\diamond,\qquad 1\le r\le4,
 \quad \Lambda_N=N^{o(1)}.
 \label{W:eq:tightness-global}
\end{equation}
The bound concerns the complete product of the function and both cutoffs.
\end{lemma}
\begin{proof}
Consider a nonzero summand in the product rule below. Its factor from $\Gamma_\rho$ implies $\kappa I\le G_\rho\le\kappa^{-1}I$, while its factor from $\chi_N$ implies
$Q_*\le2\mathfrak r_N$ and $B\le2\mathfrak r_N$. Both conditions hold at the same moment vector, including when some factors are undifferentiated.
Lemma~\ref{W:lem:tightness-directions} therefore gives, for a nonempty label set $B$,
\[
 \|D_B\mathbf Y^\sigma\|_\infty\le
 C_{|B|}m_N\rho^{|B|-1}\mathfrak r_N^{2|B|+1}
                             \prod_{a\in B}\|v_a\|_\diamond.
\]
For $S=\{1,\ldots,r\}$ the complete derivative is
\[
 D_S(\Gamma_\rho\chi_NM(\mathbf Y))
 =\sum_{S_0\sqcup S_1\sqcup S_2=S}
 (D_{S_0}\Gamma_\rho)(D_{S_1}\chi_N)
 \sum_{\pi\in\Pi(S_2)}
 D^{|\pi|}M(\mathbf Y)[D_B\mathbf Y:B\in\pi].
\]
For $S_2=\varnothing$, the inner expression means $M(\mathbf Y)$. Let $r_j=|S_j|$, and, if $r_2>0$, put $k=|\pi|$. By \eqref{W:eq:tightness-soft}, the term from differentiating the
modulus function is bounded by
\[
 C_{r,u,T}(1+\log(1+2D))^{k-1}
 m_N^k\rho^{r_2-k}\mathfrak r_N^{2r_2+k}
                         \prod_{a\in S_2}\|v_a\|_\diamond.
\]
The other factors contribute $C_r\rho^{r_0}$ and $C_r\rho^{r_1}(1+\log\rho)^{r_1}\mathfrak r_N^{3r_1}$, respectively. Since $\rho\asymp m_N$,
\[
 \rho^{r_0+r_1}m_N^k\rho^{r_2-k}\le C_rm_N^r,
 \qquad 3r_1+2r_2+k\le3r,
 \qquad D\le(K+1)^2.
\]
Thus for $r\le4$ every summand is bounded by $C_{u,T,r}m_N^r\mathfrak r_N^{12}(1+\log N)^4 \prod_a\|v_a\|_\diamond$. This estimate also covers $r_2=0$. The number of terms is at
most $3^r r^r$. It proves \eqref{W:eq:tightness-global} with a common $\Lambda_N=N^{o(1)}$. For $A_N^0$, only $S_2=\varnothing$ contributes, with last factor one. For $A_N^{\rm F}$
use $|D^k e^{\mathrm i\sum_l\lambda_ly_l}[h_1,\ldots,h_k]| \le\Lambda^k\prod_j\|h_j\|_\infty$ in the same formula. Every nonzero term lies in the support of $\chi_N$; the factor
$\Gamma_\rho$ extends the estimate across the boundary by zero.

At a fixed real moment vector, extend its real $r$th derivative by complex multilinearity in its arguments. For a complex functional $v$ on real polynomials, write
$v=v_{\mathrm R}+\mathrm i v_{\mathrm I}$. The entries of the induced variation of the Gram matrix satisfy
\[
 \dot G[v_{\mathrm R}]=\tfrac12(\dot G[v]+\overline{\dot G[v]}),\qquad
 \dot G[v_{\mathrm I}]=\tfrac1{2\mathrm i}(\dot G[v]-\overline{\dot G[v]}).
\]
Complex conjugation preserves singular values. The triangle inequality for the trace norm, together with the same inequality for each matrix entry, gives
$\|v_{\mathrm R}\|_\diamond,\|v_{\mathrm I}\|_\diamond\le\|v\|_\diamond$.
Expanding the $r$ arguments produces at most $2^r$ real multilinear expressions. Applying the proved bound to each one yields the complex multilinear estimate after multiplication
of its constant by $2^r$. The base moment vector remains real throughout this extension.
\end{proof}

Let $\mathcal C$, $E_\rho$, and $\nu_\rho$ be as in \eqref{W:eq:comparison-measure}. The contour measure has mass $\nu_\rho(\mathcal C)=N^{o(1)}$. For $z=\zeta+\zeta^{-1}$ with
$|\zeta|=e^s$, one has $|z^2-4|^{1/2}=|\zeta-\zeta^{-1}|\le2\cosh s$. Consequently $E_\rho(z)\ge e^{\rho s}/\cosh s\ge1$, since $\rho\ge1$ and $\cosh s\le e^s$. Define
$e_z(f)=f(z)-f(0)$. Since $|p_j(0)|\le1$,
\[
 \dot G[e_z]=p(z)p(z)^{\mathsf T}-p(0)p(0)^{\mathsf T},\qquad
 \norm{e_z}_\diamond\le2E_\rho(z)^2.
\]
For a function in Lemma~\ref{W:lem:tightness-global}, its first derivative has the polynomial kernel
\begin{equation}
 P_A(s,z)=D_{s}A(s)[e_z],
 \label{W:eq:tightness-response-polynomial}
\end{equation}
of degree at most $2\rho-2$. In the coordinates $s_n=\mathcal L(p_n)$, it is explicitly
\[
 P_A(s,z)=\sum_{n=1}^{2\rho-2}\frac{\partial A}{\partial s_n}(s)
                       (p_n(z)-p_n(0)).
\]
The coordinate $s_0=1$ is fixed. The derivative bound applied to $e_z$ gives the first estimate below. For the second, the fundamental theorem of calculus along the segment joining
$s$ and $s'$, both normalized by $s_0=s'_0=1$, gives
\[
 P_A(s,z)-P_A(s',z)=\int_0^1
       D_s^2A(s'+t(s-s'))[s-s',e_z]\,\dd t.
\]
Using $\|e_z\|_\diamond\le2E_\rho(z)^2$ yields
\begin{equation}
 |P_A(s,z)|\le2\Lambda_N{m_N}E_\rho(z)^2, \qquad |P_A(s,z)-P_A(s',z)|\le2\Lambda_N{m_N}^2E_\rho(z)^2\norm{s-s'}_\diamond.
 \label{W:eq:tightness-response-bounds}
\end{equation}
The derivative bound is valid along the entire segment. For a direction annihilating constants represented by a Cauchy kernel $k_v$,
\[
 \frac1{2\pi\mathrm i}\int_{\mathcal C}
       (f(z)-f(0))k_v(z)\,\mathrm dz=v(f)-f(0)v(1)=v(f).
\]
Linearity of $D_{s}A(s)$ therefore yields
\begin{equation}
 D_{s}A(s)[v]=\frac1{2\pi\mathrm i}\int_\mathcal C P_A(s,z)k_v(z)\,\mathrm d z.
 \label{W:eq:tightness-response-representation}
\end{equation}
The coefficients, including the derivatives of the sampled values and the cutoffs, remain evaluated at $s$ throughout this identity.

The same representation applies in each argument of a higher derivative. If $v_1,\ldots,v_r$ annihilate constants and $k_a$ are representing kernels, multilinearity gives, for
$1\le r\le4$,
\[
 D_s^rA(s)[v_1,\ldots,v_r]
 =\frac1{(2\pi\mathrm i)^r}
   \int_{\mathcal C^r}D_s^rA(s)[e_{z_1},\ldots,e_{z_r}]
                   \prod_{a=1}^rk_a(z_a)\,\dd z_1\cdots\dd z_r.
\]
The complex multilinear estimate just proved implies, after enlarging $\Lambda_N$ by a fixed constant,
\[
 |D_s^rA(s)[e_{z_1},\ldots,e_{z_r}]|
 \le 2^r\Lambda_Nm_N^r\prod_{a=1}^rE_\rho(z_a)^2.
\]
The absolute integral is consequently bounded by $2^r\Lambda_Nm_N^r\prod_a\int|k_a|\,\dd\nu_\rho$. This proves integrability and justifies the iterated integrals; taking the
infimum over representing kernels gives the derivative bound in the norm defined by these contour integrals. The integral representation pairs each derivative with the whole
polynomial in its corresponding variable. For fixed $r$ the mass of the product measure on the contour is $\nu_\rho(\mathcal C)^r=N^{o(1)}$. For $s=s(H)$ the chain rule gives
$D_HP_A(s(H),z)[V]=D_s^2A(s(H))[D_Hs(H)[V],e_z]$.
The bounds in this paragraph hold for $1\le r\le4$.

\subsubsection{The low energy localization statistic}
\label{W:sec:localization-statistic}

\label{W:sec:a-statistic-for-low-energy-spectral-localization}
Convergence of the coefficients on fixed intervals leaves a separate operator question: how much of the low energy spectral subspace lies beyond the polynomial vectors already
constructed? The spectral measure at $w$ cannot detect an eigenvector orthogonal to $w$. The localization statistic below is an unnormalized trace, so its spectral decomposition
includes every eigenvector in that subspace.

The scaled matrix is $A=m_N^2(2I-H)$ and the spectral weight is $\theta_E(x)=(1+\exp(x-E))^{-2}$. The projection has rank $\rho$, while the cutoff in the eigenvalues of the Gram
matrix uses the enlarged dimension $\rho_T\ge\rho$. We use the smooth contraction $P_{\rho,\kappa}^{[d]}$ of \eqref{W:eq:moment-coordinates-projection} and $\Gamma=\Gamma_{\rho_T}$
from \eqref{W:eq:moment-coordinates-J}. The former is constructed from $(p_j(h_d(H))w)_{0\le j<\rho}$; $\theta_E$ is evaluated at the original scaled matrix $A$.

On $\{\Gamma>0\}$, this contraction is the orthogonal projection $P_\rho^{[d]}$ onto the modified Krylov space. Its relation to the projection $P_\rho$ onto the first $\rho$
Lanczos vectors of $H$ requires the spectral modification to agree with $H$. Thus the two projections are identified on $\{\Gamma>0\}$ only when $H$ belongs to

\begin{equation}
 \Omega_d=\{H:h_d(H)=H\}.
 \label{W:eq:spectral-agreement-event}
\end{equation}
The norm window $\{\|H\|\le2+d/2\}$ is contained in $\Omega_d$. We use the projection identity from Lemma~\ref{W:lem:moment-coordinates-localization} only on $\Omega_d$. Define
\begin{equation}
 C_E(H)=\operatorname{Tr}[(I-P_{\rho,\kappa}^{[d]}(H))\theta_E(A)],
 \qquad Z_E(H)=\frac{C_E(H)}{1+C_E(H)},\qquad \Phi_E=\Gamma Z_E.
 \label{W:eq:gaussian-localization}
\end{equation}
All these functions are smooth for each finite $N$. Since $P_{\rho,\kappa}^{[d]}$ is a positive contraction,
\[
 C_E=\operatorname{Tr}\bigl[\theta_E(A)^{1/2}(I-P_{\rho,\kappa}^{[d]})\theta_E(A)^{1/2}\bigr],
 \qquad 0\le C_E\le\operatorname{Tr}\theta_E(A),\qquad 0\le Z_E<1.
\]
The event used for the counting bound is $\{\mathfrak K_N(H)\le K\}$ from Lemma~\ref{W:lem:gaussian-trace}. On this event, write
$n_A(t)=\#\{j:\lambda_j(A)\le t\}\le K+C(t_+)^{3/2}$. Positivity permits Tonelli's theorem, which gives
\[
 \operatorname{Tr}\theta_E(A)
 =\int_{\mathbb R}[-\theta_E'(t)]n_A(t)\,\dd t
 \le K+C\int_0^\infty[-\theta_E'(t)]t^{3/2}\,\dd t.
\]
Here $-\theta_E'(t)=2e^{t-E}(1+e^{t-E})^{-3}$. With $T_E=1+\max\{E,0\}$, the last integral is bounded by
\[
 T_E^{3/2}+2e^{2E}\int_{T_E}^\infty t^{3/2}e^{-2t}\,\dd t<\infty.
\]
This proves $C_E\le K+C_E^{\rm const}$ on the specified event, with a deterministic constant depending only on $E$ and the constant in the counting estimate. The derivatives of
$P_{\rho,\kappa}^{[d]}$ are estimated separately. On $\{\Gamma>0\}$, diagonalize the full scaled operator as $A\psi_j=\lambda_j\psi_j$ with an orthonormal eigenbasis. Then
\begin{equation}
 C_E(H)=\sum_{j=1}^N\theta_E(\lambda_j)
                   \|(I-{P_{\rho}^{[d]}})\psi_j\|^2.
 \label{W:eq:calculation-localization-eigenbasis}
\end{equation}
On $\Omega_d\cap\{\Gamma>0\}$, $P_\rho^{[d]}=P_\rho$. If $\psi_j$ belongs to an invariant block after the first, then $P_\rho\psi_j=0$ and its contribution is
$\theta_E(\lambda_j)$, while $|\langle w,\psi_j\rangle|^2=0$. Thus the trace detects this vector with its spectral weight $\theta_E(\lambda_j)$. The normalization matters: division
by $N$ would replace that contribution by $N^{-1}\theta_E(\lambda_j)$, which tends to zero for a fixed low energy eigenvalue.

The inequality $\theta_E(x)\ge1/4$ for $x\le E$ will give the operator norm estimate below. For the GOE and GUE, Lemma~\ref{W:lem:resolvent-weighted-trace} gives, on
$\Omega_N^{\mathrm{coer}}$, $\operatorname{Tr}\theta_E(A)\le C_E^{\rm const}(1+C_N)^{3/2}$, where $C_N=C_0(1+\mathfrak Z_N^4)$. On the complement,
$\operatorname{Tr}\theta_E(A)\le N$. Therefore
\[
 \mathbb E\operatorname{Tr}\theta_E(A)
 \le C_E^{\rm const}\mathbb E(1+C_N)^{3/2}
       +N\mathbb P((\Omega_N^{\mathrm{coer}})^c)
 \le C_{\beta,E}.
\]
The first expectation is bounded by \eqref{W:eq:gram-Gaussian-noise-exp-moment}, and the second term is bounded by $CN e^{-cN}$. The transformation $c\mapsto c/(1+c)$ is increasing
on $[0,\infty)$; hence $C_E>\eta^2/4$ implies $Z_E>\eta^2/(4+\eta^2)$.

Let $\Pi_E=\mathbf1_{(-\infty,E]}(A)$, and interpret $Q_{\rho}$ in the original coordinates as the complement of the first ${\rho}$ Lanczos vectors of $H$.
\begin{theorem}
\label{W:thm:gaussian-detection}
On $\{\Gamma>0\}\cap\Omega_d$,
\begin{equation}
 \|Q_{\rho}\Pi_E\|^2\le4C_E(H).
 \label{W:eq:localization-trace-bound}
\end{equation}
Consequently, for every $\eta>0$ and every measurable $0\le U\le1$,
\begin{equation}
 \mathbf1_{\{\|Q_{\rho}\Pi_E\|>\eta\}}
 \le \mathbf1_{\Omega_d^c}(H)+1-U\Gamma+
 \frac{4+\eta^2}{\eta^2}\,U\Gamma Z_E(H).
 \label{W:eq:gaussian-detection}
\end{equation}
\end{theorem}
\begin{proof}
On $\{\Gamma>0\}\cap\Omega_d$, apply Lemma~\ref{W:lem:moment-coordinates-localization} with $r=\rho$ and $a=(+1,E)$. Its polynomial projection is $P_\rho$, its trace is $C_E(H)$,
and $I-P_\rho=Q_\rho$, so \eqref{W:eq:moment-coordinates-detection} gives \eqref{W:eq:localization-trace-bound}. The complement includes every later invariant Lanczos block. Put
$\mathcal A=\{\|Q_\rho\Pi_E\|>\eta\}$. On $\Omega_d\cap\{\Gamma>0\}$, the proved norm bound gives $\mathbf1_{\mathcal A}\le(4+\eta^2)\eta^{-2}Z_E$. Since $0\le U\Gamma\le1$, and
multiplication by $U\Gamma$ annihilates the case $\Gamma=0$, it follows on $\Omega_d$ that
\[
 \begin{split}
 \mathbf1_{\mathcal A}
 &=(1-U\Gamma)\mathbf1_{\mathcal A}
                      +U\Gamma\mathbf1_{\mathcal A}\\
 &\le1-U\Gamma+\frac{4+\eta^2}{\eta^2}\,U\Gamma Z_E.
 \end{split}
\]
The right side is nonnegative everywhere. Adding $\mathbf1_{\Omega_d^c}$ therefore proves \eqref{W:eq:gaussian-detection} on the full matrix space.
\end{proof}

The expectation of $\Phi_E$ under the GOE or GUE is estimated in Section~\ref{W:sec:expectation-of-the-localization-statistic}. Its derivatives with respect to the entries of
$F_r^*\Theta_a(H)F_r$ and to $\operatorname{Tr}\Theta_a(H)$ are compared in Section~\ref{W:sec:localization-of-the-low-energy-spectral-subspace}.
Sections~\ref{W:sec:resolvents-on-exceptional-support} and~\ref{W:sec:derivative-bounds-for-the-spectral-cutoff} bound $\mathbb P(\Omega_d^c)$ for each original ensemble. Taking
expectations in \eqref{W:eq:gaussian-detection} adds this exceptional probability once to the bound for the statistic.

For $H=H_N$ and $\rho=\min\{N,\lfloor m_NL\rfloor\}$, taking expectations in \eqref{W:eq:gaussian-detection} gives
\begin{equation}
 \mathbb P(\|Q_{N,L}\Pi_N(E)\|>\eta)
 \le\mathbb P(\Omega_d^c)+\mathbb E(1-\Gamma)
                         +\frac{4+\eta^2}{\eta^2}\mathbb E\Phi_E.
 \label{W:eq:intro-localization-probability}
\end{equation}
The functions $\Gamma$ and $\Phi_E$ in this expectation are defined on every matrix. The identification with the projection onto the initial vectors of the completed Lanczos basis
was used only on $\Omega_d\cap\{\Gamma>0\}$; its two complementary cases are accounted for by $\mathbb P(\Omega_d^c)$ and $\mathbb E(1-\Gamma)$, respectively.
Section~\ref{W:sec:proof-full-localization} compares the remaining expectations and takes the limits in the stated order.

\subsection{Truncation and resolvent estimates}
\label{sec:combined-w-truncation}
\label{W:sec:truncation-and-resolvent-estimates}
\label{W:sec:entry-truncation-estimates}
\label{W:sec:comparison-resolvent-estimates}

The deterministic estimates describe how the comparison functions respond to changes of their moments. The remaining problem is to integrate those changes against the original
entry distributions and sum over all entries. This is where the tail assumption matters: the fourth moment may be infinite before truncation, and its truncated value need not
remain bounded. The corresponding moment tensors therefore stay in the estimates, while the resolvent factors are controlled conditionally.

Truncation on the original probability space gives $H_N^{\rm tr}=H_N$ outside an event of probability $o(1)$. This identity will transfer the final bounded observable after
comparison, but it supplies no prescribed polynomial rate for the exceptional event. The comparison uses the uncentered truncated entries. Centered and covariance normalized arrays
serve as auxiliary matrices for estimates requiring those assumptions; the changes in their low order moments must be accounted for when returning to $H_N^{\rm tr}$.

The needed resolvent bounds are obtained in stages because they control different parts of the expansion. Conditional Schur complement estimates first express the error in terms of
a complementary minor. Counting and entrywise estimates in both symmetry classes then bound that minor. An upper bound on the Gram matrix supplies the remaining control at the
deterministic starting vector. Integrating the resulting coordinates row by row yields the summed weights in Corollary~\ref{W:cor:integrated-comparison-weights}, which can be
paired with the derivatives from the preceding subsection.

All intervals and parameters of the functions under consideration are fixed before $N\to\infty$. Uniformity concerns the specified product distributions of the entries, principal
minors, and matrices with one entry set to zero. Bounds over each permitted segment of entry values are established before the remainder of the expansion is integrated. These
estimates are obtained from finite matrices without using a limiting point process for the eigenvalues of Wigner matrices at the spectral edge.

\subsubsection{Tail bounds and the truncated matrix}
\label{W:sec:truncation-centering-and-covariance-normalization}

The tail condition gives $\mathbb P(H_N^{\mathrm{tr}}=H_N)\to1$. Thus a bounded function of the Lanczos coefficients has the same asymptotic expectation for the two matrices. We
compare the uncentered truncated entries; the centered, covariance normalized matrix is used for the local estimates requiring exact normalization. This reduction to bounded
entries follows \cite{lee2014necessary}.

The size of the discarded moments is measured by a monotone tail function. With $\zeta=|\xi|$, it is defined by

\begin{equation}
\Theta(s)
:=
\sup_{u\ge s}
u^4\mathbb P({\zeta}>u),
\qquad
s>0.
\label{W:eq:offdiagonal-tail-function}
\end{equation}
The tail assumption gives
\begin{equation}
\Theta(s)\downarrow0
\qquad
(s\to\infty)
\label{W:eq:vanishing-offdiagonal-tail-function}
\end{equation}
by \eqref{W:eq:lee-yin-tail}.

For the diagonal array, the corresponding control is the tail of the second moment:
\begin{equation}
\Theta_{\mathrm d}(s)
:=
\sup_{N,i}
\mathbb E
\left[
d_{N,i}^2
\mathbf 1_{\{|d_{N,i}|>s\}}
\right].
\label{W:eq:diagonal-tail-function}
\end{equation}
Uniform integrability gives
\begin{equation}
\Theta_{\mathrm d}(s)\downarrow0
\qquad
(s\to\infty)
\label{W:eq:vanishing-diagonal-tail-function}
\end{equation}
by \eqref{W:eq:diagonal-uniform-integrability}.

The assumption need not give a finite fourth moment. For example, a symmetric real variable with tail $c/(s^4\log s)$ for $s\ge e$ has finite variance, whereas its fourth moment
integral contains $4c\int_e^\infty(s\log s)^{-1}\,\dd s=\infty$. Rescaling by its standard deviation preserves \eqref{W:eq:lee-yin-tail}. The bounded entry estimates of
\cite{lee2014necessary} are applied only to the truncated arrays specified in Section~\ref{W:sec:notation-and-scales}.

The following moment estimates quantify both the discarded tails and the moments that remain after truncation.

\Needspace{6\baselineskip}

\begin{lemma}
\label{W:lem:offdiag-tail-moments}
As \(L\to\infty\), the following statements hold.

For \(0<q<4\),
\begin{equation}
\mathbb E
\left[
{\zeta}^q\mathbf 1_{\{{\zeta}>L\}}
\right]
=
o(L^{q-4}).
\label{W:eq:tail-moment-below-four}
\end{equation}

In particular,
\begin{equation}
\mathbb E{\zeta}^q<\infty,
\qquad
0<q<4.
\label{W:eq:finite-below-four}
\end{equation}

For \(q=4\),
\begin{equation}
\mathbb E
\left[
{\zeta}^4\mathbf 1_{\{{\zeta}\le L\}}
\right]
=
o(\log L).
\label{W:eq:truncated-fourth-moment}
\end{equation}

For every fixed \(q>4\),
\begin{equation}
\mathbb E
\left[
{\zeta}^q\mathbf 1_{\{{\zeta}\le L\}}
\right]
=
o(L^{q-4}).
\label{W:eq:truncated-high-moment}
\end{equation}
\end{lemma}

\begin{proof}
For $a\ge0$ and $0<q<4$, the pointwise identity
\[
 a^q\mathbf1_{\{a>L\}}
 =L^q\mathbf1_{\{a>L\}}+
       q\int_L^\infty t^{q-1}\mathbf1_{\{a>t\}}\,\dd t
\]
holds: both sides vanish when $a\le L$, and otherwise the integral is $a^q-L^q$. All terms are nonnegative. Tonelli's theorem, with $a=\zeta$, therefore gives
\begin{equation}
 \mathbb E[\zeta^q\mathbf1_{\{\zeta>L\}}]
 =L^q\mathbb P(\zeta>L)+q\int_L^\infty t^{q-1}\mathbb P(\zeta>t)\,\dd t.
 \label{W:eq:tail-integration-by-parts}
\end{equation}
For $t\ge L$, the probability in the integrand is at most $\Theta(L)t^{-4}$. The integral is finite and satisfies
\[
 L^q\Theta(L)L^{-4}
   +q\Theta(L)\int_L^\infty t^{q-5}\,\dd t
 =\Theta(L)L^{q-4}\left(1+\frac q{4-q}\right).
\]
Thus
\begin{equation}
 \mathbb E[\zeta^q\mathbf1_{\{\zeta>L\}}]
 \le\frac4{4-q}\Theta(L)L^{q-4}.
 \label{W:eq:quant-tail-below-four}
\end{equation}
This tends to zero after division by $L^{q-4}$ and also proves $\mathbb E\zeta^q<\infty$, by fixing one $L$ with
$\Theta(L)<\infty$.

At the fourth moment threshold the same integral is logarithmic. For $1\le S<L$, split at $S$ and use $\min\{\zeta,L\}^4\le S^4+ 4\int_S^L t^3\mathbf1_{\{\zeta>t\}}\,\dd t$ to
obtain
\begin{equation}
 \mathbb E[\zeta^4\mathbf1_{\{\zeta\le L\}}]
 \le S^4+4\int_S^L t^3\mathbb P(\zeta>t)\,\dd t.
 \label{W:eq:fourth-split}
\end{equation}
The integral is at most $\Theta(S)\log(L/S)$. Hence
\[
 \limsup_{L\to\infty}
 \frac{\mathbb E[\zeta^4\mathbf1_{\{\zeta\le L\}}]}{\log L}
 \le4\Theta(S).
\]
Letting $S\to\infty$ proves the estimate for the truncated fourth moment.

For $q>4$, the same splitting gives
\begin{equation}
 \mathbb E[\zeta^q\mathbf1_{\{\zeta\le L\}}]
 \le S^q+q\int_S^L t^{q-1}\mathbb P(\zeta>t)\,\dd t.
 \label{W:eq:high-moment-split}
\end{equation}
The integral term is at most $q\Theta(S)L^{q-4}/(q-4)$. Divide by $L^{q-4}$, first let $L\to\infty$, and then let $S\to\infty$. This proves the higher moment estimate without
imposing a finite fourth moment.
\end{proof}

For later diagonal estimates, we need a uniform analogue at the second moment threshold.

\begin{lemma}
\label{W:lem:diagonal-truncation-moments}
Uniformly in \(N\) and \(1\le i\le N\), the following statements hold as \(L\to\infty\):
\begin{equation}
\left|
\mathbb E
\left[
d_{N,i}\mathbf 1_{\{|d_{N,i}|\le L\}}
\right]
\right|
=
o(L^{-1}),
\label{W:eq:diagonal-mean-tail}
\end{equation}
and
\begin{equation}
\mathbb E
\left[
d_{N,i}^2
\mathbf 1_{\{|d_{N,i}|>L\}}
\right]
=o(1).
\label{W:eq:diagonal-second-tail}
\end{equation}

For every fixed \(q>2\),
\begin{equation}
\sup_{N,i}
\mathbb E
\left[
|d_{N,i}|^q
\mathbf 1_{\{|d_{N,i}|\le L\}}
\right]
=
o(L^{q-2}).
\label{W:eq:diagonal-high-truncated-moment}
\end{equation}
\end{lemma}

\begin{proof}
Centering before truncation and the definition of $\Theta_{\rm d}$ give
\[
 \sup_{N,i}\left|\mathbb E[d_{N,i}\mathbf1_{\{|d_{N,i}|\le L\}}]\right|
 \le L^{-1}\Theta_{\rm d}(L)=o(L^{-1}).
\]
For the second moment estimate, subtraction of the nonnegative truncated second moment from the full second moment gives
\[
 0\le \mathbb E d_{N,i}^2-
        \mathbb E[d_{N,i}^2\mathbf1_{\{|d_{N,i}|\le L\}}]
 =\mathbb E[d_{N,i}^2\mathbf1_{\{|d_{N,i}|>L\}}]
 \le\Theta_{\rm d}(L)\longrightarrow0,
\]
uniformly in $N,i$. For $q>2$, let $M_2=\sup_{N,i}\mathbb E d_{N,i}^2$ and split the truncated moment at a fixed $S<L$. On $|d|\le S$ use $|d|^q\le S^{q-2}d^2$; on $S<|d|\le L$ use
$|d|^q\le L^{q-2}d^2$. Thus
\begin{equation}
 \sup_{N,i}L^{2-q}\mathbb E[|d_{N,i}|^q\mathbf1_{\{|d_{N,i}|\le L\}}]
 \le M_2(S/L)^{q-2}+\Theta_{\rm d}(S).
 \label{W:eq:diagonal-moment-split}
\end{equation}
Taking $L\to\infty$ and then $S\to\infty$ proves the required uniform limit.
\end{proof}

We use $L_N=e_N\sqrt N$, choosing $e_N\downarrow0$ so that the normalized tail and moment errors defined below remain small after multiplication by each fixed power of $e_N^{-1}$.
The truncated fourth moment has its own bound in Section~\ref{W:sec:moment-bounds-and-taylor-expansion-at-two-truncation-scales}.

For \(q=1,2,3\), define
\begin{equation}
\eta_q^{>}(L)
:=
L^{4-q}
\mathbb E
\left[
{\zeta}^q\mathbf 1_{\{{\zeta}>L\}}
\right].
\label{W:eq:eta-tail-q}
\end{equation}
For \(q=4\), define
\begin{equation}
\eta_4^{\le}(L)
:=
\frac{
\mathbb E[
{\zeta}^4\mathbf 1_{\{{\zeta}\le L\}}
]
}{
\log(\mathrm e+L)
}.
\label{W:eq:eta-four}
\end{equation}
For integers \(q\ge5\), define
\begin{equation}
\eta_q^{\le}(L)
:=
L^{4-q}
\mathbb E
\left[
{\zeta}^q\mathbf 1_{\{{\zeta}\le L\}}
\right].
\label{W:eq:eta-high-q}
\end{equation}
By Lemma~\ref{W:lem:offdiag-tail-moments},
\begin{equation}
\eta_q^{>}(L)\to0,
\qquad
\eta_4^{\le}(L)\to0,
\qquad
\eta_q^{\le}(L)\to0.
\label{W:eq:eta-offdiag-zero}
\end{equation}

For each integer \(q\ge3\), define
\begin{equation}
\eta_{\mathrm d,q}^{\le}(L)
:=
L^{2-q}
\sup_{N,i}
\mathbb E
\left[
|d_{N,i}|^q
\mathbf 1_{\{|d_{N,i}|\le L\}}
\right].
\label{W:eq:eta-diagonal-q}
\end{equation}
Then
\begin{equation}
\eta_{\mathrm d,q}^{\le}(L)\to0
\label{W:eq:eta-diagonal-zero}
\end{equation}
by Lemma~\ref{W:lem:diagonal-truncation-moments}.

\Needspace{6\baselineskip}

Let $(f_j)_{j\ge1}$ be a countable collection of nonnegative functions with $f_j(s)\to0$ as $s\to\infty$, containing all the tail and normalized moment bounds just listed.

\begin{lemma}
\label{W:lem:gaussian-cutoff-choice}
There are increasing integers $N_r$, $r\ge5$, and a nondecreasing sequence $\ell_N^{\mathrm{slow}}\to\infty$ such that $\ell_N^{\mathrm{slow}}=r$ for $N_r\le N<N_{r+1}$ and, with
$L_N=\sqrt N/\ell_N^{\mathrm{slow}}$,
\begin{equation}
 \ell_N^{\mathrm{slow}}\le(\log N)^{1/4},\qquad L_N\ge\ell_N^{\mathrm{slow}},
 \qquad \max_{j\le\ell_N^{\mathrm{slow}}}f_j(L_N)\le(\ell_N^{\mathrm{slow}})^{-(\ell_N^{\mathrm{slow}})^{2}}.
 \label{W:eq:gaussian-stage-bound}
\end{equation}
Consequently $(\ell_N^{\mathrm{slow}})^{A} f_j(L_N)\to0$ for every fixed $A,j$.
\end{lemma}

The same construction can enforce any countable collection of convergence requirements $g_j(N,r)\to0$ as $N\to\infty$ for each fixed integer $r$, by requiring
$g_j(N,\ell_N^{\mathrm{slow}})\le(\ell_N^{\mathrm{slow}})^{-(\ell_N^{\mathrm{slow}})^{2}}$ for $j\le\ell_N^{\mathrm{slow}}$. These requirements must be specified before the
thresholds $N_r$ are chosen.

\begin{proof}
For each integer $r\ge5$, choose $S_r\ge r$ so that $S_r\ge S_{r-1}$ and
\[
 \sup_{s\ge S_r} f_j(s)\le r^{-r^2}\qquad(1\le j\le r).
\]
Choose $N_r>N_{r-1}$ such that, for all $N\ge N_r$,
\[
 N\ge r^2 S_r^2,\qquad \log N\ge r^4,\qquad
 \max_{j\le r}g_j(N,r)\le r^{-r^2}.
\]
The last requirement is omitted when no $g_j$ have been chosen. More explicitly, for each $j\le r$ choose an integer $n_{r,j}$ such that $g_j(N,r)\le r^{-r^2}$ for all
$N\ge n_{r,j}$. One may take
\[
 N_r=1+\max\{N_{r-1},\lceil r^2S_r^2\rceil,
                  \lceil e^{r^4}\rceil,n_{r,1},\ldots,n_{r,r}\}.
\]
The terms $n_{r,j}$ are absent when no $g_j$ are specified. This is a finite integer, and the recursion makes $N_r$ strictly increasing. For $N_r\le N<N_{r+1}$ set
$\ell_N^{\mathrm{slow}}=r$. Then
\[
 L_N=\frac{\sqrt N}{r}\ge S_r\ge r,\qquad
 r\le(\log N)^{1/4},\qquad
 \max_{j\le r}f_j(L_N)\le r^{-r^2}.
\]
These are \eqref{W:eq:gaussian-stage-bound}. For fixed $A$ and $j$, once $r\ge j$,
\[
 (\ell_N^{\mathrm{slow}})^A f_j(L_N)
 \le r^{A-r^2}\longrightarrow0,
 \qquad
 \frac{\log\ell_N^{\mathrm{slow}}}{\log N}
 \le\frac{\log r}{r^4}\longrightarrow0.
\]
Also $L_N\ge r\to\infty$. The same substitution proves $g_j(N,\ell_N^{\mathrm{slow}})\le (\ell_N^{\mathrm{slow}})^{-(\ell_N^{\mathrm{slow}})^2}$ for each fixed $j$ and all
sufficiently large $N$. Assign any positive value to the finitely many terms before $N_5$.
\end{proof}

\Needspace{6\baselineskip}

We choose a nondecreasing deterministic sequence
\begin{equation}
\ell_N^{\mathrm{slow}}\longrightarrow\infty,
\qquad
\ell_N^{\mathrm{slow}}=N^{o(1)},
\label{W:eq:ell-scale}
\end{equation}
and put
\begin{equation}
{e_N}:=(\ell_N^{\mathrm{slow}})^{-1},
\qquad
L_N:={e_N}\sqrt N
=
\frac{\sqrt N}{\ell_N^{\mathrm{slow}}},
\label{W:eq:LN-definition}
\end{equation}
The sequence is chosen to satisfy the following requirements.

The scale satisfies
\begin{equation}
L_N\to\infty,
\qquad
{e_N}\downarrow0,
\qquad
{e_N}=N^{-o(1)}.
\label{W:eq:epsilon-properties}
\end{equation}

The discarded entries satisfy
\begin{equation}
N^2\mathbb P({\zeta}>L_N)\to0,
\label{W:eq:offdiag-equality-rate}
\end{equation}
and
\begin{equation}
\sum_{i=1}^{N}
\mathbb P(|d_{N,i}|>L_N)
\to0.
\label{W:eq:diagonal-equality-rate}
\end{equation}

For every fixed pair of integers \(A\ge1, \qquad q\ge1,\) the corresponding little-\(o\) factors satisfy
\begin{equation}
(\ell_N^{\mathrm{slow}})^{A}\eta_q^{>}(L_N)\to0,
\qquad
q=1,2,3,
\label{W:eq:strong-tail-absorption}
\end{equation}
\begin{equation}
(\ell_N^{\mathrm{slow}})^{A}\eta_4^{\le}(L_N)\to0,
\label{W:eq:strong-fourth-absorption}
\end{equation}
\begin{equation}
(\ell_N^{\mathrm{slow}})^{A}\eta_q^{\le}(L_N)\to0,
\qquad
q\ge5,
\label{W:eq:strong-high-absorption}
\end{equation}
and
\begin{equation}
(\ell_N^{\mathrm{slow}})^{A}\Theta_{\mathrm d}(L_N)\to0,
\label{W:eq:strong-diagonal-tail-absorption}
\end{equation}
\begin{equation}
(\ell_N^{\mathrm{slow}})^{A}
\eta_{\mathrm d,q}^{\le}(L_N)
\to0,
\qquad
q\ge3.
\label{W:eq:strong-diagonal-moment-absorption}
\end{equation}

\begin{lemma}
\label{W:lem:admissible-scale}
The sequence in \eqref{W:eq:ell-scale} can be chosen to satisfy \eqref{W:eq:epsilon-properties}--\eqref{W:eq:strong-diagonal-moment-absorption}.
\end{lemma}

\begin{proof}
Apply Lemma~\ref{W:lem:gaussian-cutoff-choice} to the countable family $\Theta,\Theta_{\rm d},\eta_q^>$ for $q=1,2,3$, $\eta_4^\le$, $\eta_q^\le$ for integers $q\ge5$, and
$\eta_{{\rm d},q}^\le$ for integers $q\ge3$. Every member tends to zero by the preceding tail estimates. The lemma gives $\ell_N^{\rm slow}\to\infty$, $\ell_N^{\rm slow}=N^{o(1)}$
and $(\ell_N^{\rm slow})^A f_j(\sqrt N/\ell_N^{\rm slow})\to0$ for every fixed $A,j$. With $e_N=(\ell_N^{\rm slow})^{-1}$, the off-diagonal members are respectively the normalized
quantities in \eqref{W:eq:strong-tail-absorption}--\eqref{W:eq:strong-diagonal-moment-absorption}. The diagonal member of order $q$ is normalized by $L_N^{q-2}$, as in its
definition. This proves all these limits for the same sequence.

Since $\mathbb P(\zeta>L)\le L^{-4}\Theta(L)$, the off-diagonal union bound is
\[
 N^2\mathbb P(\zeta>L_N)
 \le N^2L_N^{-4}\Theta(L_N)
 =(\ell_N^{\rm slow})^4\Theta(L_N)\longrightarrow0.
\]
For the diagonal entries, $L^2\mathbf1_{\{|d|>L\}}\le d^2\mathbf1_{\{|d|>L\}}$ gives
\[
 \sum_{i=1}^N\mathbb P(|d_{N,i}|>L_N)
 \le NL_N^{-2}\Theta_{\rm d}(L_N)
 =(\ell_N^{\rm slow})^2\Theta_{\rm d}(L_N)\longrightarrow0.
\]
These prove \eqref{W:eq:offdiag-equality-rate} and \eqref{W:eq:diagonal-equality-rate} for the selected sequence.
\end{proof}

The fourth-moment bound is stated after division by $\log(e+L_N)$. Proposition~\ref{W:prop:gaussian-log-obstruction} gives an example in which the unnormalized fourth moment
diverges.

\subparagraph{The uncentered truncated matrix}
\label{W:sec:moments-of-truncated-entries}

Using the scale \(L_N\) from Lemma~\ref{W:lem:admissible-scale}, define
\begin{equation}
\widehat\xi_{ij}
:=
\xi_{ij}
\mathbf 1_{\{|\xi_{ij}|\le L_N\}},
\qquad
i<j,
\label{W:eq:truncated-offdiag-truncation}
\end{equation}
and
\begin{equation}
\widehat d_{N,i}
:=
d_{N,i}
\mathbf 1_{\{|d_{N,i}|\le L_N\}}.
\label{W:eq:truncated-diagonal-truncation}
\end{equation}
Define
\begin{equation}
\widehat X_{ij}
=
\widehat\xi_{ij},
\qquad
\widehat X_{ji}
=
\overline{\widehat\xi_{ij}},
\qquad
\widehat X_{ii}
=
\widehat d_{N,i},
\label{W:eq:truncated-Xhat}
\end{equation}
and
\begin{equation}
{H_N^{\mathrm{tr}}}
:=
N^{-1/2}\widehat X_N.
\label{W:eq:truncated-Hhat}
\end{equation}

Deterministically,
\begin{equation}
\max_{1\le i,j\le N}
|H_{ij}^{\mathrm{tr}}|
\le
N^{-1/2}L_N
=
{e_N}.
\label{W:eq:entrywise-small}
\end{equation}

Define the equality event
\begin{equation}
 \Omega_N^{\mathrm{tr}} :={} \bigcap_{1\le i<j\le N} \{|\xi_{ij}|\le L_N\} \cap \bigcap_{1\le i\le N} \{|d_{N,i}|\le L_N\}.
 \label{W:eq:equality-event}
\end{equation}
On this event,
\begin{equation}
{H_N^{\mathrm{tr}}}=H_N.
\label{W:eq:matrix-equality}
\end{equation}

A union bound, together with \eqref{W:eq:offdiag-equality-rate} and \eqref{W:eq:diagonal-equality-rate}, gives $\mathbb P(\Omega_N^{\mathrm{tr}})\to1$. Since
${H_N^{\mathrm{tr}}}=H_N$ on this event, every deterministic measurable functional of the two matrices agrees with probability $1-o(1)$; in particular, any weak limit proved for
the truncated matrix transfers to the original Wigner matrix. More quantitatively, put $p_N=\mathbb P((\Omega_N^{\mathrm{tr}})^c)$. For a bounded measurable matrix statistic $F$
and a measurable event $\mathcal E$,
\[
 |\mathbb EF(H_N)-\mathbb EF(H_N^{\mathrm{tr}})|
 \le2\|F\|_\infty p_N,\qquad
 \mathbb P(H_N\in\mathcal E)
 \le\mathbb P(H_N^{\mathrm{tr}}\in\mathcal E)+p_N.
\]
Both estimates follow by splitting on $\Omega_N^{\mathrm{tr}}$ and its complement. We apply them after summing the entry comparisons, to a bounded function of the complete matrix.

With $\zeta=|\xi|$, write
\begin{equation}
L=L_N,
\qquad
\widehat\xi
=
\xi\mathbf 1_{\{|\xi|\le L\}}.
\label{W:eq:xihat-short}
\end{equation}

Since \(\mathbb E\xi=0,\) we have
\[
\mathbb E\widehat\xi
=
-
\mathbb E
\left[
\xi\mathbf 1_{\{\zeta>L\}}
\right].
\]
Therefore
\begin{equation}
|\mathbb E\widehat\xi|
\le
\mathbb E
\left[
\zeta\mathbf 1_{\{\zeta>L\}}
\right]
=
o(L^{-3}).
\label{W:eq:xihat-mean}
\end{equation}

More quantitatively, by \eqref{W:eq:quant-tail-below-four},
\begin{equation}
|\mathbb E\widehat\xi|
\le
\frac43\Theta(L)L^{-3}.
\label{W:eq:xihat-mean-quant}
\end{equation}

Next,
\[
\mathbb E|\widehat\xi|^2
=
1-
\mathbb E
\left[
\zeta^2\mathbf 1_{\{\zeta>L\}}
\right],
\]
so
\begin{equation}
\mathbb E|\widehat\xi|^2
=
1+o(L^{-2}).
\label{W:eq:xihat-second}
\end{equation}
Since \(|\mathbb E\widehat\xi|^2=o(L^{-6}),\) the centered variance satisfies
\begin{equation}
\operatorname{cum}
(\widehat\xi,\overline{\widehat\xi})
=
\mathbb E
\left|
\widehat\xi-\mathbb E\widehat\xi
\right|^2
=
1+o(L^{-2}).
\label{W:eq:xihat-variance}
\end{equation}

The pseudovariance satisfies
\begin{equation}
 \mathbb E\widehat\xi^2 = \mathbb E\xi^2 - \mathbb E \left[ \xi^2\mathbf 1_{\{\zeta>L\}} \right] = \frac2\beta-1 + o(L^{-2}),
 \label{W:eq:xihat-pseudosecond}
\end{equation}
because
\[
\left|
\mathbb E
\left[
\xi^2\mathbf 1_{\{\zeta>L\}}
\right]
\right|
\le
\mathbb E
\left[
\zeta^2\mathbf 1_{\{\zeta>L\}}
\right].
\]
After centering,
\begin{equation}
\operatorname{cum}
(\widehat\xi,\widehat\xi)
=
\frac2\beta-1
+
o(L^{-2}).
\label{W:eq:xihat-pseudovariance}
\end{equation}

In particular, in the complex Hermitian case,
\begin{equation}
\operatorname{cum}
(\widehat\xi,\widehat\xi)
=
o(L^{-2}).
\label{W:eq:complex-pseudovariance-defect}
\end{equation}

Define the normalized matrix entry
\begin{equation}
{u^{\mathrm{tr}}}
:=
N^{-1/2}\widehat\xi.
\label{W:eq:hhat-normalized}
\end{equation}
Joint cumulants are multilinear, so
\begin{equation}
\operatorname{cum}
(
{u^{\mathrm{tr}}}^{\epsilon_1},
\ldots,
{u^{\mathrm{tr}}}^{\epsilon_r}
)
=
N^{-r/2}
\operatorname{cum}
(
\widehat\xi^{\epsilon_1},
\ldots,
\widehat\xi^{\epsilon_r}
),
\label{W:eq:cumulant-homogeneity-oriented}
\end{equation}
where \(z^{+}:=z, \qquad z^{-}:=\overline z.\)

Thus
\begin{equation}
\mathbb E{u^{\mathrm{tr}}}
=
o\left(
{e_N}^{-3}N^{-2}
\right),
\label{W:eq:hhat-first-cumulant}
\end{equation}
because \(N^{-1/2}L_N^{-3} = {e_N}^{-3}N^{-2}.\) Moreover,
\begin{equation}
\operatorname{cum}
({u^{\mathrm{tr}}},\overline{{u^{\mathrm{tr}}}})
=
\frac1N
+
o\left(
{e_N}^{-2}N^{-2}
\right),
\label{W:eq:hhat-variance}
\end{equation}
and
\begin{equation}
\operatorname{cum}
({u^{\mathrm{tr}}},{u^{\mathrm{tr}}})
=
\frac1N\left(\frac2\beta-1\right)
+
o\left(
{e_N}^{-2}N^{-2}
\right).
\label{W:eq:hhat-pseudovariance}
\end{equation}

For \(\beta=2\), this becomes
\begin{equation}
\operatorname{cum}
({u^{\mathrm{tr}}},{u^{\mathrm{tr}}})
=
o\left(
{e_N}^{-2}N^{-2}
\right).
\label{W:eq:hhat-complex-pseudovariance}
\end{equation}

\subsubsection{Cumulants of the truncated entries}
\label{W:sec:cumulants-of-truncated-entries}

Truncation changes the cumulants as well as the support. In particular, an entry that was centered before truncation usually has a nonzero first cumulant afterwards. The comparison
therefore includes the resulting mean and covariance errors with their matrix normalization.

For a complex entry, the two real coordinates belong to one independent block and may be dependent within that block. A product of $r$ copies of the first coordinate and $s$ copies
of the second contributes the joint moment $\mathbb E[X_1^rX_2^s]$. Equivalently, one may use joint cumulants of the entry and its conjugate, as in equation~(4.14) of
\cite{cipolloni2022multiresolvent}. The estimates here follow directly from the moment--cumulant formula for the truncated distribution. For random variables $Z_1,\ldots,Z_p$ with
finite $p$th moments, that formula is

\begin{equation}
\operatorname{cum}(Z_1,\ldots,Z_p)
=
\sum_{\pi\in\Pi_p}
(|\pi|-1)!
(-1)^{|\pi|-1}
\prod_{B\in\pi}
\mathbb E
\left[
\prod_{j\in B}Z_j
\right],
\label{W:eq:joint-cumulant-formula}
\end{equation}
where \(\Pi_p\) denotes the set of set partitions of \(\{1,\ldots,p\}.\)

For one real random variable \(Z\), we write \(\kappa_p(Z) := \operatorname{cum}(Z,\ldots,Z).\)

Two elementary properties will be used repeatedly:

\begin{equation}
\operatorname{cum}
(c_1Z_1,\ldots,c_pZ_p)
=
\left(
\prod_{j=1}^{p}c_j
\right)
\operatorname{cum}(Z_1,\ldots,Z_p),
\label{W:eq:cumulant-multilinearity}
\end{equation}
and, for \(p\ge2\), joint cumulants are unchanged by deterministic translations of any argument.

Assume first that \(\beta=1,\) so \(\xi\) is real.

\Needspace{6\baselineskip}

\begin{proposition}
\label{W:prop:real-offdiag-cumulants}
The truncated variable \(\widehat\xi = \xi\mathbf 1_{\{|\xi|\le L_N\}}\) satisfies
\begin{align}
 \kappa_1(\widehat\xi) &= o(L_N^{-3}), \label{W:eq:real-kappa-one}\\
 \kappa_2(\widehat\xi) &= 1+o(L_N^{-2}), \label{W:eq:real-kappa-two}\\
 |\kappa_3(\widehat\xi)| &\le C, \label{W:eq:real-kappa-three}\\
 |\kappa_4(\widehat\xi)| &= o(\log L_N), \label{W:eq:real-kappa-four}\\
 \kappa_p(\widehat\xi) &= o(L_N^{p-4}). \label{W:eq:real-kappa-high}\qquad(p\ge5).
\end{align}

\end{proposition}

\begin{proof}
Equations \eqref{W:eq:xihat-mean} and \eqref{W:eq:xihat-variance} prove \eqref{W:eq:real-kappa-one} and \eqref{W:eq:real-kappa-two}. Lemma~\ref{W:lem:offdiag-tail-moments} gives
$\mathbb E|\xi|^3<\infty$; hence the first three truncated moments are bounded and
\[
 \kappa_3(\widehat\xi)=\mathbb E\widehat\xi^3
 -3\mathbb E\widehat\xi^2\mathbb E\widehat\xi
 +2(\mathbb E\widehat\xi)^3=O(1).
\]
For the fourth cumulant,
\begin{equation}
 \kappa_4(\widehat\xi)=\mathbb E\widehat\xi^4
 -4\mathbb E\widehat\xi^3\mathbb E\widehat\xi
 -3(\mathbb E\widehat\xi^2)^2
 +12\mathbb E\widehat\xi^2(\mathbb E\widehat\xi)^2
 -6(\mathbb E\widehat\xi)^4.
 \label{W:eq:fourth-cumulant-expanded}
\end{equation}
Its first term is $o(\log L_N)$ by \eqref{W:eq:truncated-fourth-moment}, and the remaining terms are $O(1)$. This proves \eqref{W:eq:real-kappa-three} and
\eqref{W:eq:real-kappa-four}.

Fix $p\ge5$ and write $m_q(L)=\mathbb E\widehat\xi^q$ for truncation at $L$. The moment estimates are
\begin{equation}
 m_1(L)=o(L^{-3}),\label{W:eq:moment-hierarchy-one}
\end{equation}
\begin{equation}
 m_2(L)=O(1),\qquad m_3(L)=O(1),
 \label{W:eq:moment-hierarchy-two-three}
\end{equation}
\begin{equation}
 m_4(L)=o(\log L),\label{W:eq:moment-hierarchy-four}
\end{equation}
and
\begin{equation}
 m_q(L)=o(L^{q-4}),\qquad q\ge5.
 \label{W:eq:moment-hierarchy-high}
\end{equation}
The partition formula reads
\begin{equation}
 \kappa_p(\widehat\xi)=\sum_{\pi\in\Pi_p}c_\pi
                  \prod_{B\in\pi}m_{|B|}(L),\qquad
 c_\pi=(-1)^{|\pi|-1}(|\pi|-1)!.
 \label{W:eq:cumulant-partition-scalar}
\end{equation}
We claim that every product satisfies
\begin{equation}
 \prod_{B\in\pi}|m_{|B|}(L)|=o(L^{p-4}).
 \label{W:eq:partition-product-goal}
\end{equation}
For the single block partition this is \eqref{W:eq:moment-hierarchy-high}. Otherwise let $\mathcal A(\pi)=\{B\in\pi:|B|\ge5\}$ and $b_4=|\{B\in\pi:|B|=4\}|$. If $\mathcal A(\pi)$
is empty, the product is $O((1+\log L)^{b_4})=o(L^{p-4})$. If $\mathcal A(\pi)$ is nonempty, then
\[
 p-4-\sum_{B\in\mathcal A(\pi)}(|B|-4)
 =\sum_{B\notin\mathcal A(\pi)}|B|+4(|\mathcal A(\pi)|-1)\ge1.
\]
In particular,
\begin{equation}
 \sum_{B\in\mathcal A(\pi)}(|B|-4)\le p-5.
 \label{W:eq:large-block-exponent}
\end{equation}
Consequently,
\[
 L^{4-p}\prod_{B\in\pi}|m_{|B|}(L)|
 \le C_p L^{-1}(1+\log L)^{b_4}\longrightarrow0.
\]
The same estimate holds for any product of numbers indexed by the sets $B$ whose absolute values are bounded by $C_{|B|}m_{|B|}(L)$, with $m_1,m_2,m_3=O(1)$ and the stated bounds
for $m_q$, $q\ge4$. Indeed, the preceding inequalities use only these absolute bounds and the sizes of the sets in the partition. The finite sum in
\eqref{W:eq:cumulant-partition-scalar} proves \eqref{W:eq:real-kappa-high}. Finally,
\[
 \kappa_p(u^{\mathrm{tr}})=N^{-p/2}\kappa_p(\widehat\xi),\qquad
 N^{-p/2}L_N^{p-4}=e_N^{p-4}N^{-2}.
\]
Together with the lower order estimates, these identities prove \eqref{W:eq:real-normalized-kappa-one}--\eqref{W:eq:real-normalized-kappa-high}.
\end{proof}

For the matrix entry \({u^{\mathrm{tr}}}=N^{-1/2}\widehat\xi,\) one has
\begin{align}
 \kappa_1({u^{\mathrm{tr}}}) &= o\left( {e_N}^{-3}N^{-2} \right), \label{W:eq:real-normalized-kappa-one}\\
 \kappa_2({u^{\mathrm{tr}}}) &= \frac1N + o\left( {e_N}^{-2}N^{-2} \right), \label{W:eq:real-normalized-kappa-two}\\
 |\kappa_3({u^{\mathrm{tr}}})| &\le CN^{-3/2}, \label{W:eq:real-normalized-kappa-three}\\
 |\kappa_4({u^{\mathrm{tr}}})| &= N^{-2}o(\log N)\le N^{-2+o(1)}, \label{W:eq:real-normalized-kappa-four}\\
 \kappa_p({u^{\mathrm{tr}}}) &= o\left( {e_N}^{p-4}N^{-2} \right). \label{W:eq:real-normalized-kappa-high}\qquad(p\ge5).
\end{align}

Assume now that \(\beta=2.\) Write
\begin{equation}
\xi
=
\frac1{\sqrt2}(X+\mathrm iY),
\label{W:eq:complex-real-coordinates}
\end{equation}
where \(X,Y\) are real.

From \(\mathbb E|\xi|^2=1 \qquad\text{and}\qquad \mathbb E\xi^2=0,\) we obtain
\begin{equation}
\mathbb EX^2=\mathbb EY^2=1,
\qquad
\mathbb EXY=0.
\label{W:eq:complex-coordinate-covariance}
\end{equation}

Define the radially truncated variables \(\widehat\xi = \xi\mathbf 1_{\{|\xi|\le L\}},\)
\begin{equation}
\widehat X
:=
\sqrt2\,\operatorname{Re}\widehat\xi,
\qquad
\widehat Y
:=
\sqrt2\,\operatorname{Im}\widehat\xi.
\label{W:eq:complex-truncated-coordinates}
\end{equation}
Then
\begin{equation}
|\widehat X|+|\widehat Y|
\le
2\zeta\mathbf 1_{\{\zeta\le L\}}.
\label{W:eq:complex-coordinate-bound}
\end{equation}

Set
\[
\widehat{\mathbf x}
=
\begin{pmatrix}
\widehat X\\
\widehat Y
\end{pmatrix}.
\]

\Needspace{6\baselineskip}

\begin{proposition}
\label{W:prop:complex-mixed-cumulants}
As $L\to\infty$, componentwise and in operator norm,
\begin{align}
\mathbb E\widehat{\mathbf x}&=o(L^{-3}),
\label{W:eq:complex-coordinate-mean}\\
\operatorname{Cov}(\widehat{\mathbf x})&=I_2+o(L^{-2}).
\label{W:eq:complex-coordinate-covariance-truncated}
\end{align}
For fixed $r\ge3$ and $\sigma_j\in\{1,2\}$, with $\widehat X^{(1)}=\widehat X$ and $\widehat X^{(2)}=\widehat Y$,
\begin{align}
|\cum(\widehat X^{(\sigma_1)},\widehat X^{(\sigma_2)},\widehat X^{(\sigma_3)})|
&\le C,
\label{W:eq:complex-third-real-cumulants}\\
\cum(\widehat X^{(\sigma_1)},\ldots,\widehat X^{(\sigma_4)})
&=o(\log L),
\label{W:eq:complex-fourth-real-cumulants}\\
\cum(\widehat X^{(\sigma_1)},\ldots,\widehat X^{(\sigma_r)})
&=o(L^{r-4}),\qquad r\ge5.
\label{W:eq:complex-high-real-cumulants}
\end{align}
\end{proposition}

\begin{proof}
Write $\zeta=|\xi|$ and let $m_q(L)=\mathbb E[\zeta^q\mathbf1_{\{\zeta\le L\}}]$. The estimates in Lemma~\ref{W:lem:offdiag-tail-moments}, uniformly in the entry distribution, give
\[
 m_q(L)\le C_q\quad(1\le q\le3),\qquad
 m_4(L)=o(\log L),\qquad m_q(L)=o(L^{q-4})\quad(q\ge5).
\]
The discarded moments also obey $\mathbb E[\zeta^q\mathbf1_{\{\zeta>L\}}]=o(L^{q-4})$ for $q=1,2$, by \eqref{W:eq:quant-tail-below-four}. Since $|X|+|Y|\le2\sqrt2\,\zeta$, the
original zero mean of each coordinate gives
\[
 |\mathbb E\widehat X|+|\mathbb E\widehat Y|
 \le2\sqrt2\,\mathbb E[\zeta\mathbf1_{\{\zeta>L\}}]
 =o(L^{-3}).
\]
For every coordinate pair $a,b\in\{1,2\}$,
\[
 \left|\mathbb E\widehat X^{(a)}\widehat X^{(b)}
             -\mathbb EX^{(a)}X^{(b)}\right|
 \le2\mathbb E[\zeta^2\mathbf1_{\{\zeta>L\}}]
 =o(L^{-2}).
\]
The original coordinate covariance is $I_2$, and subtracting the products of the truncated means adds only $o(L^{-6})$. This proves \eqref{W:eq:complex-coordinate-mean} and
\eqref{W:eq:complex-coordinate-covariance-truncated}, including the off-diagonal covariance entry.

For a labelled nonempty block $B\subset\{1,\ldots,r\}$ of size $q$, the deterministic bound $|\widehat X^{(a)}|\le\sqrt2\zeta\mathbf1_{\{\zeta\le L\}}$ gives
\begin{equation}
 \left|\mathbb E\prod_{j\in B}\widehat X^{(\sigma_j)}\right|
 \le2^{q/2}m_q(L).
 \label{W:eq:complex-block-moment-bound}
\end{equation}
We use the labelled joint cumulant formula
\[
 \operatorname{cum}(Y_1,\ldots,Y_r)
 =\sum_{\pi\in\Pi(\{1,\ldots,r\})}
       (-1)^{|\pi|-1}(|\pi|-1)!
             \prod_{B\in\pi}\mathbb E\prod_{j\in B}Y_j.
\]
For $r=3$, every block moment is bounded and there are finitely many partitions, so all third coordinate cumulants are bounded uniformly. For $r=4$, the partition consisting of one
block is $o(\log L)$ by the fourth moment estimate; every other partition consists of blocks of sizes at most three and hence contributes $O(1)=o(\log L)$.

For $r\ge5$, \eqref{W:eq:complex-block-moment-bound} bounds the product associated with a partition $\pi$ by $2^{r/2}\prod_{B\in\pi}m_{|B|}(L)$. The partition estimate in the proof
of Proposition~\ref{W:prop:real-offdiag-cumulants} gives $o(L^{r-4})$ for each product. Summing over the finitely many partitions proves \eqref{W:eq:complex-high-real-cumulants}.

Finally $\widehat\xi=(\widehat X+\mathrm i\widehat Y)/\sqrt2$, and its conjugate has the analogous expression with $-\mathrm i$. Multilinearity of a joint cumulant therefore gives,
with $c_{+,1}=c_{-,1}=1$, $c_{+,2}=\mathrm i$, and $c_{-,2}=-\mathrm i$,
\[
 \operatorname{cum}(\widehat\xi^{\epsilon_1},\ldots,\widehat\xi^{\epsilon_r})
 =2^{-r/2}\sum_{\sigma\in\{1,2\}^r}
       \left(\prod_{j=1}^r c_{\epsilon_j,\sigma_j}\right)
       \operatorname{cum}(\widehat X^{(\sigma_1)},\ldots,
                           \widehat X^{(\sigma_r)}).
\]
There are $2^r$ terms and each coefficient has modulus one. This proves the bounds for cumulants of $\widehat\xi$ and $\overline{\widehat\xi}$ displayed below. Each matrix entry
coordinate is $N^{-1/2}$ times the corresponding standardized variable, so a cumulant of order $r$ acquires the factor $N^{-r/2}$. Since $L_N=e_N\sqrt N$,
\[
 N^{-r/2}L_N^{r-4}=e_N^{r-4}N^{-2}.
\]
The normalized third, fourth, and higher cumulant bounds follow with their stated orders.
\end{proof}

Multilinearity gives the following equivalent bounds for $\epsilon_j\in\{+,-\}$:
\begin{align}
|\cum(\widehat\xi^{\epsilon_1},\widehat\xi^{\epsilon_2},\widehat\xi^{\epsilon_3})|
&\le C,
\label{W:eq:complex-third-oriented-cumulants}\\
\cum(\widehat\xi^{\epsilon_1},\ldots,\widehat\xi^{\epsilon_4})
&=o(\log L),
\label{W:eq:complex-fourth-oriented-cumulants}\\
\cum(\widehat\xi^{\epsilon_1},\ldots,\widehat\xi^{\epsilon_r})
&=o(L^{r-4}),\qquad r\ge5.
\label{W:eq:complex-high-oriented-cumulants}
\end{align}
After multiplying each argument by $N^{-1/2}$,
\begin{align}
|\cum(N^{-1/2}\widehat X^{(\sigma_1)},N^{-1/2}\widehat X^{(\sigma_2)},N^{-1/2}\widehat X^{(\sigma_3)})|
&\le CN^{-3/2},
\label{W:eq:complex-normalized-third}\\
|\cum(N^{-1/2}\widehat X^{(\sigma_1)},\ldots,N^{-1/2}\widehat X^{(\sigma_4)})|
&=N^{-2}o(\log N),
\label{W:eq:complex-normalized-fourth}\\
\cum(N^{-1/2}\widehat X^{(\sigma_1)},\ldots,N^{-1/2}\widehat X^{(\sigma_r)})
&=o({e_N}^{r-4}N^{-2}),\qquad r\ge5.
\label{W:eq:complex-normalized-high}
\end{align}

Set
\begin{equation}
\widehat d_{N,i}
=
d_{N,i}\mathbf 1_{\{|d_{N,i}|\le L_N\}},
\qquad
{u^{\mathrm{tr}}}_{N,i}^{\mathrm d}
=
N^{-1/2}\widehat d_{N,i}.
\label{W:eq:diagonal-normalized-entry}
\end{equation}

\Needspace{6\baselineskip}

\begin{proposition}
\label{W:prop:diagonal-cumulants}
Uniformly in \(N\) and \(i\),
\begin{equation}
\mathbb E\widehat d_{N,i}
=
o(L_N^{-1}),
\label{W:eq:diagonal-truncated-mean}
\end{equation}
and therefore
\begin{equation}
\kappa_1({u^{\mathrm{tr}}}_{N,i}^{\mathrm d})
=
o\left(
{e_N}^{-1}N^{-1}
\right).
\label{W:eq:diagonal-normalized-mean}
\end{equation}

Moreover,
\begin{equation}
\left|
\kappa_2({u^{\mathrm{tr}}}_{N,i}^{\mathrm d})
\right|
\le
C N^{-1}.
\label{W:eq:diagonal-second-cumulant}
\end{equation}

For every fixed \(p\ge3\),
\begin{equation}
\kappa_p(\widehat d_{N,i})
=
o(L_N^{p-2})
\label{W:eq:diagonal-truncated-high-cumulant}
\end{equation}
uniformly in \(N,i\), and hence
\begin{equation}
\kappa_p({u^{\mathrm{tr}}}_{N,i}^{\mathrm d})
=
o\left(
{e_N}^{p-2}N^{-1}
\right).
\label{W:eq:diagonal-normalized-high-cumulant}
\end{equation}
\end{proposition}

\begin{proof}
Equation \eqref{W:eq:diagonal-mean-tail} and $N^{-1/2}L_N^{-1}=e_N^{-1}N^{-1}$ prove the two bounds on the means. Also
\[
 \operatorname{Var}(\widehat d_{N,i})
 \le\mathbb E d_{N,i}^2\le M_2.
\]
For $p\ge3$, the partition formula and H\"older's inequality give
\begin{equation}
 |\kappa_p(Z)|\le C_p\mathbb E|Z|^p,
 \qquad C_p=\sum_{\pi\in\Pi_p}(|\pi|-1)!.
 \label{W:eq:cumulant-moment-bound}
\end{equation}
Indeed, each partition product satisfies
\[
 \prod_{B\in\pi}|\mathbb E Z^{|B|}|
 \le\prod_{B\in\pi}(\mathbb E|Z|^p)^{|B|/p}
 =\mathbb E|Z|^p.
\]
The constant depends only on $p$. Applying \eqref{W:eq:diagonal-high-truncated-moment}, uniformly over the diagonal labels, we obtain
\[
 |\kappa_p(\widehat d_{N,i})|=o(L_N^{p-2}),\qquad
 N^{-p/2}L_N^{p-2}=e_N^{p-2}N^{-1}.
\]
These are the claimed truncated and normalized cumulant estimates.
\end{proof}

\subsubsection{Centering and covariance normalization}
\label{W:sec:centering-and-covariance-normalization}

The truncated entries have the support needed for comparison, but their means and covariances no longer have the exact normalization required by the auxiliary local laws. Centering
and covariance normalization provide such an auxiliary array. These operations change entries even on the event $H_N^{\rm tr}=H_N$, so their effect on nonlinear observables must be
estimated rather than discarded. The formulas below identify the resulting moment changes; the actual comparison continues to use the uncentered array $H_N^{\rm tr}$.

For real symmetric matrices, $\beta=1$, the normalization involves only a scalar mean and variance:

\begin{equation}
\mu_N
:=
\mathbb E\widehat\xi,
\qquad
\sigma_N^2
:=
\operatorname{Var}(\widehat\xi).
\label{W:eq:real-standardization-parameters}
\end{equation}
The discarded moment bounds imply
\begin{equation}
\mu_N=o(L_N^{-3}),
\qquad
\sigma_N^2=1+o(L_N^{-2}).
\label{W:eq:real-standardization-parameter-bounds}
\end{equation}
In particular, \(\sigma_N>0\) for all sufficiently large \(N\), and
\begin{equation}
\sigma_N^{\pm1}
=
1+o(L_N^{-2}).
\label{W:eq:sigma-expansion}
\end{equation}

The normalized real entry is
\begin{equation}
\xi_N^{\mathrm{nor}}
:=
\frac{\widehat\xi-\mu_N}{\sigma_N}.
\label{W:eq:real-standardized-variable}
\end{equation}
Then
\begin{equation}
\mathbb E\xi_N^{\mathrm{nor}}=0,
\qquad
\mathbb E(\xi_N^{\mathrm{nor}})^2=1.
\label{W:eq:real-standardized-moments}
\end{equation}
Moreover,
\begin{equation}
|\xi_N^{\mathrm{nor}}|
\le
2L_N
\label{W:eq:real-standardized-bound}
\end{equation}
for all sufficiently large \(N\).

For every fixed \(p\ge3\),
\begin{equation}
\kappa_p(\xi_N^{\mathrm{nor}})
=
\sigma_N^{-p}\kappa_p(\widehat\xi),
\label{W:eq:standardized-cumulants-real}
\end{equation}
because cumulants of order at least two are translation invariant. Thus the bounds for cumulants of orders three and above in Proposition~\ref{W:prop:real-offdiag-cumulants} remain
unchanged up to a factor \(1+o(1)\).

For complex Hermitian matrices, we normalize the joint covariance of the two real coordinates. Suppose \(\beta=2\). Define
\begin{equation}
\boldsymbol\mu_N
:=
\mathbb E\widehat{\mathbf x},
\qquad
\Sigma_N
:=
\operatorname{Cov}(\widehat{\mathbf x}).
\label{W:eq:complex-standardization-parameters}
\end{equation}
By Proposition~\ref{W:prop:complex-mixed-cumulants},
\begin{equation}
\boldsymbol\mu_N=o(L_N^{-3}),
\qquad
\Sigma_N=I_2+o(L_N^{-2}).
\label{W:eq:complex-standardization-parameter-bounds}
\end{equation}
Hence \(\Sigma_N\) is positive definite for all sufficiently large \(N\), and
\begin{equation}
\Sigma_N^{\pm1/2}
=
I_2+o(L_N^{-2})
\label{W:eq:Sigma-expansion}
\end{equation}
in operator norm.

Define
\begin{equation}
\mathbf x_N^{\mathrm{nor}}
:=
\Sigma_N^{-1/2}
\left(
\widehat{\mathbf x}-\boldsymbol\mu_N
\right).
\label{W:eq:complex-standardized-vector}
\end{equation}
Then
\begin{equation}
\mathbb E\mathbf x_N^{\mathrm{nor}}=0,
\qquad
\operatorname{Cov}(\mathbf x_N^{\mathrm{nor}})=I_2.
\label{W:eq:complex-standardized-covariance}
\end{equation}
Writing
\[
\mathbf x_N^{\mathrm{nor}}
=
\begin{pmatrix}
X_N^{\mathrm{nor}}\\
Y_N^{\mathrm{nor}}
\end{pmatrix},
\]
define
\begin{equation}
\xi_N^{\mathrm{nor}}
:=
\frac1{\sqrt2}
(X_N^{\mathrm{nor}}+\mathrm iY_N^{\mathrm{nor}}).
\label{W:eq:complex-standardized-variable}
\end{equation}
Then
\begin{equation}
\mathbb E\xi_N^{\mathrm{nor}}=0,
\qquad
\mathbb E|\xi_N^{\mathrm{nor}}|^2=1,
\qquad
\mathbb E(\xi_N^{\mathrm{nor}})^2=0.
\label{W:eq:complex-standardized-moments}
\end{equation}
Furthermore,
\begin{equation}
|\xi_N^{\mathrm{nor}}|
\le
CL_N
\label{W:eq:complex-standardized-bound}
\end{equation}
for all sufficiently large \(N\).

Since joint cumulants are multilinear and invariant under translations at orders at least two, every fixed order cumulant tensor of \(\mathbf x_N^{\mathrm{nor}}\) is a finite
linear combination of the corresponding cumulants of
\(\widehat{\mathbf x},\) with coefficients given by entries of \(\Sigma_N^{-1/2}=I_2+o(L_N^{-2}).\) Therefore all
estimates of Proposition~\ref{W:prop:complex-mixed-cumulants} remain valid for the standardized variables.

Let \(H_N^{\mathrm{nor}}\) be the self-adjoint matrix whose off-diagonal entries are the standardized variables
\(N^{-1/2}\xi_{ij,N}^{\mathrm{nor}},\) and whose diagonal entries are the centered truncated variables
\begin{equation}
N^{-1/2}
\left(
\widehat d_{N,i}
-
\mathbb E\widehat d_{N,i}
\right).
\label{W:eq:auxiliary-standardized-matrix}
\end{equation}

The truncated and standardized off-diagonal coordinates are related exactly by
\begin{equation}
\widehat\xi
=
\mu_N+\sigma_N\xi_N^{\mathrm{nor}}
\label{W:eq:real-affine-relation}
\end{equation}
in the real case, and by
\begin{equation}
\widehat{\mathbf x}
=
\boldsymbol\mu_N
+
\Sigma_N^{1/2}\mathbf x_N^{\mathrm{nor}}
\label{W:eq:complex-affine-relation}
\end{equation}
in the complex case.

The normalized and uncentered matrices are related by the following deterministic map. Regard $\operatorname{Herm}_N(\mathbb F_\beta)$ as a real vector space, with unnormalized
Hilbert--Schmidt norm $\|H\|_{\mathrm{HS}}=(\operatorname{Tr}H^*H)^{1/2}$. For $i<j$ put
\[
 v_{ij}(H):=
 \begin{cases}
 H_{ij}\in\mathbb R,&\beta=1,\\
 (\operatorname{Re}H_{ij},\operatorname{Im}H_{ij})^{\mathsf T}
       \in\mathbb R^2,&\beta=2,
 \end{cases}
 \qquad
 S_N:=
 \begin{cases}
 \sigma_N I_1,&\beta=1,\\
 \Sigma_N^{1/2},&\beta=2.
 \end{cases}
\]
In the complex case, the passage from the standardized real coordinate vector to $v_{ij}(H)$ is multiplication by $(2N)^{-1/2}$; this scalar commutes with $\Sigma_N^{1/2}$. Define
$\mathfrak T_N$ on Hermitian matrices by
\begin{equation}
 \begin{gathered}
 (\mathfrak T_NH)_{ii}=H_{ii},\qquad
 v_{ij}(\mathfrak T_NH)=S_Nv_{ij}(H)\quad(i<j),\\
 (\mathfrak T_NH)_{ji}=\overline{(\mathfrak T_NH)_{ij}}\quad(i<j).
 \end{gathered}
 \label{W:eq:normalization-map-block-definition}
\end{equation}
This is a real-linear map preserving Hermitian matrices. Set $A_N:=\mathbb E H_N^{\mathrm{tr}}$. Its upper-triangular off-diagonal entries are $N^{-1/2}\mathbb E\widehat\xi$, and
its diagonal entries are $N^{-1/2}\mathbb E\widehat d_{N,i}$. The coordinate identities \eqref{W:eq:real-affine-relation}--\eqref{W:eq:complex-affine-relation}, together with the
unchanged centered diagonal, now prove exactly
\begin{equation}
 {H_N^{\mathrm{tr}}}=A_N+\mathfrak T_N(H_N^{\mathrm{nor}}).
 \label{W:eq:matrix-affine-relation}
\end{equation}

Let $I$ denote the identity map on this real matrix space, and define
\[
 \|\mathfrak T_N-I\|_{\mathrm{HS}\to\mathrm{HS}}
 :=\sup_{\substack{H=H^*\\H\ne0}}
       \frac{\|(\mathfrak T_N-I)H\|_{\mathrm{HS}}}
            {\|H\|_{\mathrm{HS}}}.
\]
Then, for $N\ge2$,
\begin{equation}
 \|\mathfrak T_N-I\|_{\mathrm{HS}\to\mathrm{HS}} =\delta_N, \qquad \delta_N:=\|S_N-I_\beta\|_{\mathbb R^\beta\to\mathbb R^\beta} =o(L_N^{-2})=o({e_N}^{-2}N^{-1}).
 \label{W:eq:entrywise-linear-map-close}
\end{equation}
To verify the rate, in the real case $|\sigma_N-1|=|\sigma_N^2-1|/(\sigma_N+1)=o(L_N^{-2})$ by \eqref{W:eq:real-standardization-parameter-bounds}. In the complex case the
eigenvalues of $\Sigma_N$ lie in $[1/2,3/2]$ for large $N$, and the spectral theorem gives
\[
 \|\Sigma_N^{1/2}-I_2\|
 =\max_{\lambda\in\sigma(\Sigma_N)}
                \frac{|\lambda-1|}{\sqrt\lambda+1}
 \le\frac{\|\Sigma_N-I_2\|}{1+\sqrt{\lambda_{\min}(\Sigma_N)}}
 =o(L_N^{-2}).
\]
Here the displayed norms are Euclidean operator norms on $\mathbb R^2$. For the induced norm on the space of matrices, the diagonal is unchanged, and
\begin{equation}
 \|(\mathfrak T_N-I)H\|_{\mathrm{HS}}^2 =2\sum_{i<j}\|(S_N-I_\beta)v_{ij}(H)\|_2^2 \le2\delta_N^2\sum_{i<j}\|v_{ij}(H)\|_2^2 \le\delta_N^2\|H\|_{\mathrm{HS}}^2.
 \label{W:eq:normalization-map-HS-bound}
\end{equation}
Conversely, choose a unit vector $v\in\mathbb R^\beta$ attaining $\|(S_N-I_\beta)v\|_2=\delta_N$, and a Hermitian $H$ with only the $(1,2)$ and $(2,1)$ entries nonzero and
$v_{12}(H)=v$. Then $\|H\|_{\mathrm{HS}}^2=2$ and $\|(\mathfrak T_N-I)H\|_{\mathrm{HS}}^2=2\delta_N^2$. This proves equality in \eqref{W:eq:entrywise-linear-map-close}.

The corresponding bound for the induced operator norm follows from the inequalities between the matrix norms:
\begin{equation}
 \|(\mathfrak T_N-I)H\|_{\mathrm{op}}
 \le\delta_N\|H\|_{\mathrm{HS}}
 \le\sqrt N\,\delta_N\|H\|_{\mathrm{op}}.
 \label{W:eq:normalization-map-operator-consequence}
\end{equation}
For the nonlinear observables, we use \eqref{W:eq:matrix-affine-relation} on each entry. The resulting moment differences are defined in \eqref{W:eq:comparison-defects} and
estimated in Section~\ref{W:sec:errors-from-the-first-two-truncated-moments}.

The deterministic centering matrix satisfies
\begin{equation}
\|A_N\|_{\mathrm{op}}
\le
\sqrt N\,|\mathbb E\widehat\xi|
+
N^{-1/2}
\max_i
\left|
\mathbb E\widehat d_{N,i}
\right|.
\label{W:eq:mean-matrix-norm}
\end{equation}
Indeed, the off-diagonal mean matrix has Hilbert--Schmidt norm \(\sqrt{N-1}\,|\mathbb E\widehat\xi|,\) and the diagonal mean matrix has operator norm
\(N^{-1/2} \max_i |\mathbb E\widehat d_{N,i}|.\) Consequently,
\begin{equation}
\|A_N\|_{\mathrm{op}}
=
o({e_N}^{-3}N^{-1})
+
o({e_N}^{-1}N^{-1}).
\label{W:eq:mean-matrix-small}
\end{equation}

For the scales of Lemma~\ref{W:lem:admissible-scale}, the estimates proved above give the following bounds for the matrix entries. We collect the normalizations for use in the
comparison. For a real off-diagonal coordinate ${u^{\mathrm{tr}}}=N^{-1/2}\widehat\xi$,
\begin{align}
\kappa_1({u^{\mathrm{tr}}})&=o({e_N}^{-3}N^{-2}),
\label{W:eq:truncation-output-real-one}\\
\kappa_2({u^{\mathrm{tr}}})&=N^{-1}+o({e_N}^{-2}N^{-2}),
\label{W:eq:truncation-output-real-two}\\
|\kappa_3({u^{\mathrm{tr}}})|&\le CN^{-3/2},
\label{W:eq:truncation-output-real-three}\\
|\kappa_4({u^{\mathrm{tr}}})|&=N^{-2}o(\log N)=N^{-2+o(1)},
\label{W:eq:truncation-output-real-four}\\
\kappa_p({u^{\mathrm{tr}}})&=o({e_N}^{p-4}N^{-2}),\qquad p\ge5.
\label{W:eq:truncation-output-real-high}
\end{align}
In the complex class, the mean and the errors in covariance and pseudocovariance obey
\begin{align}
\mathbb E{u^{\mathrm{tr}}}&=o({e_N}^{-3}N^{-2}),
\label{W:eq:truncation-output-complex-one}\\
\cum({u^{\mathrm{tr}}},\overline{{u^{\mathrm{tr}}}})&=N^{-1}+o({e_N}^{-2}N^{-2}),
\label{W:eq:truncation-output-complex-variance}\\
\cum({u^{\mathrm{tr}}},{u^{\mathrm{tr}}})&=o({e_N}^{-2}N^{-2}),
\label{W:eq:truncation-output-complex-pseudovariance}
\end{align}
and, with $\epsilon_j$ denoting either the entry or its conjugate, the cumulants satisfy
\begin{align}
|\cum({u^{\mathrm{tr}}}^{\epsilon_1},{u^{\mathrm{tr}}}^{\epsilon_2},{u^{\mathrm{tr}}}^{\epsilon_3})|
&\le CN^{-3/2},
\label{W:eq:truncation-output-complex-three}\\
\cum({u^{\mathrm{tr}}}^{\epsilon_1},\ldots,{u^{\mathrm{tr}}}^{\epsilon_4})
&=N^{-2}o(\log N),
\label{W:eq:truncation-output-complex-four}\\
\cum({u^{\mathrm{tr}}}^{\epsilon_1},\ldots,{u^{\mathrm{tr}}}^{\epsilon_r})
&=o({e_N}^{r-4}N^{-2}),\qquad r\ge5.
\label{W:eq:truncation-output-complex-high}
\end{align}
For every fixed $p>4$,
\begin{equation}
\mathbb E|{u^{\mathrm{tr}}}|^p=o({e_N}^{p-4}N^{-2}).
\label{W:eq:truncation-output-high-moments}
\end{equation}
For a diagonal coordinate ${u^{\mathrm{tr}}}_{N,i}^{\rm d}$, uniformly in $i$,
\begin{align}
\kappa_1({u^{\mathrm{tr}}}_{N,i}^{\rm d})&=o({e_N}^{-1}N^{-1}),
\label{W:eq:truncation-output-diagonal-one}\\
|\kappa_2({u^{\mathrm{tr}}}_{N,i}^{\rm d})|&\le CN^{-1},
\label{W:eq:truncation-output-diagonal-two}\\
\kappa_p({u^{\mathrm{tr}}}_{N,i}^{\rm d})&=o({e_N}^{p-2}N^{-1}),\qquad p\ge3.
\label{W:eq:truncation-output-diagonal-high}
\end{align}
Use the deterministic scales
\[
L_N={e_N}\sqrt N,\qquad
{e_N}\downarrow0,\qquad
{e_N}=N^{-o(1)},\qquad L_N\to\infty.
\]

Let the centered auxiliary array $H_N^{\mathrm{nor}}$ have covariance normalized off-diagonal blocks, with the diagonal only centered as in
\eqref{W:eq:auxiliary-standardized-matrix}.

\begin{theorem}
\label{W:thm:truncation-and-normalization}

The uncentered truncated matrix ${H_N^{\mathrm{tr}}}$ satisfies
\begin{align}
\mathbb P({H_N^{\mathrm{tr}}}=H_N)&\longrightarrow1,
\label{W:eq:truncation-output-equality}\\
\max_{i,j}|H_{ij}^{\mathrm{tr}}|&\le{e_N}.
\label{W:eq:truncation-output-bound}
\end{align}
For all sufficiently large $N$,
\begin{equation}
{H_N^{\mathrm{tr}}}=A_N+\mathfrak T_N(H_N^{\mathrm{nor}}),
\label{W:eq:truncation-output-affine}
\end{equation}
where $\mathfrak T_N$ is the real-linear block map \eqref{W:eq:normalization-map-block-definition}, and
\begin{align}
 \|A_N\|_{\mathrm{op}}&\le N^{-1+o(1)},
 \label{W:eq:truncation-output-mean-matrix}\\
 \|\mathfrak T_N-I\|_{\mathrm{HS}\to\mathrm{HS}}
 &=\|S_N-I_\beta\|_{\mathbb R^\beta\to\mathbb R^\beta}
   =o({e_N}^{-2}N^{-1}).
 \label{W:eq:truncation-output-linear-map}
\end{align}
and
\begin{equation}
\max_{i,j}|(H_N^{\mathrm{nor}})_{ij}|\le C{e_N}.
\label{W:eq:truncation-output-standardized-bound}
\end{equation}
\end{theorem}

\begin{proof}
Take the common cutoff from Lemma~\ref{W:lem:admissible-scale}. By the union bound,
\[
 \mathbb P(H_N^{\mathrm{tr}}\ne H_N)
 \le N^2\mathbb P(|\xi|>L_N)
       +\sum_{i=1}^N\mathbb P(|d_{N,i}|>L_N)\longrightarrow0.
\]
The entry bound follows from $L_N/\sqrt N=e_N$. The construction in \eqref{W:eq:matrix-affine-relation} gives the exact identity \eqref{W:eq:truncation-output-affine}. Its
deterministic mean matrix satisfies
\[
 \|A_N\|\le N\max_{i<j}|\mathbb EH_{ij}^{\mathrm{tr}}|
                +\max_i|\mathbb EH_{ii}^{\mathrm{tr}}|
 =o(e_N^{-3}N^{-1})+o(e_N^{-1}N^{-1})\le N^{-1+o(1)}.
\]
Equation~\eqref{W:eq:entrywise-linear-map-close} gives
$\|\mathfrak T_N-I\|_{\mathrm{HS}\to\mathrm{HS}}=\|S_N-I_\beta\|=o(e_N^{-2}N^{-1})$.
Finally $\|S_N^{-1}\|\le C$ for large $N$, so each normalized off-diagonal block has norm at most $C L_N$ before division by $\sqrt N$. For the diagonal,
\[
 |d_{N,i}\mathbf1_{\{|d_{N,i}|\le L_N\}}
       -\mathbb E[d_{N,i}\mathbf1_{\{|d_{N,i}|\le L_N\}}]|
 \le2L_N.
\]
This proves \eqref{W:eq:truncation-output-standardized-bound}.
\end{proof}

\subsubsection{The fourth moment and Taylor expansion at two truncation thresholds}
\Needspace{6\baselineskip} The preceding moment bounds leave a fourth moment factor in the error after summation over $O(N^2)$ entries. The following example shows that choosing
the truncation threshold $L_N$ more slowly does not in general remove this factor. The remainder estimates instead use two truncation thresholds, with different bounds at the
moderate radius and at the radius needed for equality with the original matrix. The signed fourth moment tensor is kept throughout this calculation.

\begin{proposition}
\label{W:prop:gaussian-log-obstruction}
Write $K_{4,N}=\mathbb E[|\xi|^4\mathbf1_{\{|\xi|\le L_N\}}]$. There is a centered, variance one real entry law satisfying the assumed tail condition such that, for every cutoff
$L_N={e_N}\sqrt N$ with ${e_N}=N^{-o(1)}\to0$ and
$N^2\mathbb P(|\xi|>L_N)\to0$,
\begin{equation}
 K_{4,N}\asymp\frac{\log N}{\log\log N},\qquad
 \log({e_N}^{-1})=O(\log\log\log N).
 \label{W:eq:gaussian-slow-tail-scales}
\end{equation}
In particular, $\log K_{4,N}/|\log{e_N}|\to\infty$ and ${e_N}^cK_{4,N}\to\infty$ for every fixed $c>0$. The same example has a complex version with zero pseudovariance.
\end{proposition}

\begin{proof}
Choose a nonnegative variable $R_0$ whose survival function is $c_*/(s^4\log\log s)$ for $s\ge e^e$, with $c_*>0$ small, and complete it below $e^e$ by putting the remaining mass
at a positive bounded point. This is a probability law because the displayed tail is decreasing and has mass less than one at $e^e$. Its second moment is finite and positive. Put
\[
 \sigma_0=(\mathbb E R_0^2)^{1/2},\qquad
 \zeta=R_0/\sigma_0,\qquad \xi=\varepsilon\zeta,
\]
where $\varepsilon$ is an independent symmetric sign. Then $\mathbb E\zeta^2=1$, $\mathbb E\xi=0$, and
\[
 \mathbb P(\zeta>s)
 =\mathbb P(R_0>\sigma_0s)
 \sim\frac{c}{s^4\log\log s}
\]
for a constant $c>0$. Thus $s^4\mathbb P(|\xi|>s)\to0$.

The untruncated fourth moment is infinite, since the layer cake formula and the last tail asymptotic give
\[
 \mathbb E|\xi|^4=\mathbb E\zeta^4
 =4\int_0^\infty t^3\mathbb P(\zeta>t)\,\dd t=\infty.
\]
It is the finite, truncated fourth moment that depends on the cutoff:
\begin{equation}
 \mathbb E[\zeta^4\mathbf1_{\{\zeta\le L\}}] =4\int_0^L t^3\mathbb P(\zeta>t)\,\dd t -L^4\mathbb P(\zeta>L) \asymp\frac{\log L}{\log\log L}.
 \label{W:eq:gaussian-slow-tail-truncated-fourth}
\end{equation}
Indeed, after $u=\log t$, the tail part of the integral is comparable to $\int^{\log L}\dd u/\log u$. For $x=\log L$, integration over $[x/2,x]$ gives a lower bound $x/(2\log x)$;
splitting the upper bound at $\sqrt x$ gives $O(\sqrt x)+2x/\log x$. The subtracted term is only $O(1/\log\log L)$. Since $\log L_N=(1/2-o(1))\log N$, this proves the first
estimate in \eqref{W:eq:gaussian-slow-tail-scales}. The equality condition becomes $e_N^{-4}/\log\log N\to0$. In particular, $e_N^{-1}\le(\log\log N)^{1/4}$ for large $N$, which
gives the second estimate. Thus $\log K_{4,N}\sim\log\log N$, whereas $|\log e_N|=O(\log\log\log N)$, and
\[
 e_N^cK_{4,N}\ge
 c'\log N/(\log\log N)^{1+c/4}\longrightarrow\infty
\]
for every fixed $c>0$. For the complex example set $\xi=\omega\zeta$, where $\omega$ is independent of $\zeta$ and uniform on $\{1,-1,\ii,-\ii\}$. The identities
\[
 \mathbb E\omega=\mathbb E\omega^2=0,\qquad |\omega|=1
\]
give $\mathbb E\xi=\mathbb E\xi^2=0$ and $\mathbb E|\xi|^2=1$. All radial tails and truncated radial moments are those of $\zeta$. For $p+q=3$, symmetry under
$\omega\mapsto-\omega$ gives
\[
 \mathbb E[(\operatorname{Re}\xi)^p(\operatorname{Im}\xi)^q]=0,
 \qquad
 \mathbb E[(\operatorname{Re}\xi)^p(\operatorname{Im}\xi)^q
                        \mathbf1_{\{\zeta\le L\}}]=0.
\]
These expectations exist because $\mathbb E\zeta^3<\infty$. To compare the uncentered truncated entry with a bounded random variable matching its first three moment tensors, put
\begin{equation}
 \widehat\xi_N=\omega\zeta\mathbf1_{\{\zeta\le L_N\}},\qquad \sigma_N^2=\mathbb E[\zeta^2\mathbf1_{\{\zeta\le L_N\}}], \qquad \xi_N^{\mathrm b}=\sigma_N\omega',
 \label{W:eq:gaussian-slow-tail-bound}
\end{equation}
where $\omega'$ is an independent copy of $\omega$. Both laws have mean zero, second radial moment $\sigma_N^2$,
nonconjugated second moment zero, and zero third real coordinate tensor. Moreover $\sigma_N\to1$ and $|\xi_N^{\mathrm
b}|\le1$.

Define their scaled vectors of real coordinates by
\[
 \mathbf u_N^{\mathrm{tr}}
   =N^{-1/2}(\sqrt2\operatorname{Re}\widehat\xi_N,
              \sqrt2\operatorname{Im}\widehat\xi_N),\qquad
 \mathbf u_N^{\mathrm b}
   =N^{-1/2}(\sqrt2\operatorname{Re}\xi_N^{\mathrm b},
              \sqrt2\operatorname{Im}\xi_N^{\mathrm b}).
\]
For $p+q=4$ write $\Delta_{pq,N}=\mathbb E[(u_{N,1}^{\mathrm{tr}})^p (u_{N,2}^{\mathrm{tr}})^q] -\mathbb E[(u_{N,1}^{\mathrm b})^p(u_{N,2}^{\mathrm b})^q]$. Since the squared
Euclidean norm of each standardized block is twice its squared radial modulus, expansion gives the exact identity
\begin{equation}
 \Delta_{40,N}+2\Delta_{22,N}+\Delta_{04,N}
   =\frac{4}{N^2}(K_{4,N}-\sigma_N^4).
 \label{W:eq:gaussian-slow-tail-fourth-defect}
\end{equation}
The three differences $\Delta_{40,N}$, $\Delta_{22,N}$, and $\Delta_{04,N}$ are the independent real-coordinate components of the fourth moment tensor that arise when one expands
the fourth-order entry replacement. Since $\sigma_N^4\le1$ and $K_{4,N}\to\infty$,
\[
 \begin{aligned}
 \frac4{N^2}(K_{4,N}-\sigma_N^4)
 &=\Delta_{40,N}+2\Delta_{22,N}+\Delta_{04,N}\\
 &\le|\Delta_{40,N}|+2|\Delta_{22,N}|+|\Delta_{04,N}|\\
 &\le4\max\{|\Delta_{40,N}|,|\Delta_{22,N}|,|\Delta_{04,N}|\}.
 \end{aligned}
\]
Consequently, for sufficiently large $N$,
\[
 \max_{(p,q)\in\{(4,0),(2,2),(0,4)\}}|\Delta_{pq,N}|
 \ge\frac{K_{4,N}-\sigma_N^4}{N^2}
 \ge\frac{K_{4,N}}{2N^2}.
\]
The factor $N^{-2}$ is exactly the normalization of four matrix entry coordinates. The discrepancy before matrix scaling is therefore at least of order $K_{4,N}$ in one of the
three real fourth-order tensor components, even though the first and third moment tensors agree and the second radial moments coincide.

For the covariance-normalized comparison, $\sigma_N^2\to1$. Choose $N_0$ such that $\sigma_N^2\ge1/2$ for $N\ge N_0$. The variable $\widehat\xi_N/\sigma_N$ has second radial moment
one, vanishing mean and complex pseudovariance, and vanishing real third moment tensor. Its fourth radial moment is
\[
 \mathbb E|\widehat\xi_N/\sigma_N|^4
    =K_{4,N}/\sigma_N^4=(1+o(1))K_{4,N},
\]
whereas $\omega'$ has fourth radial moment one. The fourth-moment difference is therefore $(1+o(1))K_{4,N}$ after exact covariance normalization as well. This proves the stated
lower bound for the fourth-moment factor when the first three coordinate moment tensors match.
\end{proof}

\subparagraph{Taylor expansion at two truncation thresholds}
\label{W:sec:moment-bounds-and-taylor-expansion-at-two-truncation-scales}
All tensors below use standardized real coordinates. A selected entry is $N^{-1/2}\mathcal T_\alpha x$; thus
\begin{equation}
 \partial_{\alpha,\sigma}F(H)=D_HF(H)[N^{-1/2}T_{\alpha,\sigma}].
 \label{W:eq:B-block-derivative}
\end{equation}
For $a<b$ in the complex case we use $T_1=(E_{ab}+E_{ba})/\sqrt2$ and $T_2=\mathrm i(E_{ab}-E_{ba})/\sqrt2$ with $x=(\sqrt2\Re\xi,\sqrt2\Im\xi)$. For an ordered moment tensor put
\begin{equation}
 M_q(X)=\mathbb EX^{\otimes q},\qquad
 \|M\|_{1,\mathrm{ord}}=\sum_{\sigma_1,\ldots,\sigma_q}
 |M_{\sigma_1\cdots\sigma_q}|.
 \label{W:eq:B-tensor-l1}
\end{equation}
The coefficient of a multiindex $\nu$ is $\Delta\mathbb EX^\nu/\nu!$; in particular the mixed quartic coefficient is
\begin{equation}
 \frac{\Delta\mathbb EX_1^2X_2^2}{2!2!}.
 \label{W:eq:B-complex-fourth-real}
\end{equation}
The expansion of a complex entry uses the joint distribution of its two real coordinates and the full coordinate tensors defined above. The formulas below are exact Taylor
identities; we estimate their derivative weights after specifying the conditioning events.

With $e_N,L_N$ as in \eqref{W:eq:LN-definition}, put
\begin{equation}
 K_N=N^\tau,\qquad \tau=\frac18.
 \label{W:eq:truncated-moment-scales}
\end{equation}
Thus $K_N<L_N$ for all sufficiently large $N$. The moderate matrix radius is $b_{*,N}=K_N/\sqrt N=N^{-3/8}$; the outer matrix radius is ${e_N}$. Constants arising from representing
a complex entry in real coordinates are fixed and have no effect on these scales.

Write $Z$ for the vector of real coordinates: $Z=\xi$ for a real off-diagonal entry and $Z=(\sqrt2\operatorname{Re}\xi,\sqrt2\operatorname{Im}\xi)$ for a complex entry. Define its
radial variable by
\[
 r_\alpha(x)=|x|\quad(d_\alpha=1),\qquad
 r_\alpha(x)=\frac{\|x\|_2}{\sqrt2}
       \quad(d_\alpha=2).
\]
For off-diagonal entries $r_\alpha(Z)=|\xi|$. In either case $\|x\|_1\le2r_\alpha(x)$. For diagonal entries use the absolute value. The truncation is
$\widehat Z=Z\mathbf1_{\{r_\alpha(Z)\le L_N\}}$, with the mass outside the truncation radius assigned to zero.

Put $R_0=|\xi|$ and use the tail functions $\Theta$ and $\Theta_{\mathrm d}$ from \eqref{W:eq:offdiagonal-tail-function} and \eqref{W:eq:diagonal-tail-function}. The deterministic
factor $\Lambda_N$ bounds the coefficient derivatives in \eqref{W:eq:comparison-global} and the contour integrals in \eqref{W:eq:tightness-response-representation}. For a fixed
number $r$ of integrations, the relevant factors are bounded by a fixed power of
\[
 1+\mathfrak r_N^{20}(1+\log\rho)^4+\nu_\rho(\mathcal C)^r.
\]
Each is $N^{o(1)}$, so their finite products can be bounded by a deterministic $\Lambda_N\ge1$ with $\log\Lambda_N/\log N\to0$. The tail threshold $e_N$ is chosen separately in
Lemma~\ref{W:lem:admissible-scale}.

\Needspace{6\baselineskip}

Set
\[
 h_4(L)=\mathbb E[R_0^4\mathbf1_{\{R_0\le L\}}],\qquad
 h_4^{>}(K,L)=\mathbb E[R_0^4\mathbf1_{\{K<R_0\le L\}}].
\]

By Lemma~\ref{W:lem:offdiag-tail-moments}, $\mathbb ER_0^3<\infty$ and $h_4(L)=o(\log L)$. The remaining estimates at the two truncation thresholds are as follows.

\begin{lemma}
\label{W:lem:truncated-moment-tail-moments}
For $1\le K<L$,
\begin{align}
 h_4^{>}(K,L)&\le\Theta(K)\{1+4\log(L/K)\},
       \label{W:eq:truncated-moment-fourth-tail}\\
 \mathbb E[R_0^5\mathbf1_{\{R_0\le K\}}]&\le K^3.
       \label{W:eq:truncated-moment-moderate-fifth}
\end{align}
\end{lemma}
\begin{proof}
For each nonnegative $a$,
\[
 a^4\mathbf1_{\{K<a\le L\}}
 \le K^4\mathbf1_{\{a>K\}}+
                     4\int_K^L t^3\mathbf1_{\{a>t\}}\,\dd t.
\]
It is an equality for $K<a\le L$, and an inequality with nonnegative right side for $a>L$. Tonelli's theorem and the tail bound give
\[
 h_4^{>}(K,L)
 \le K^4\mathbb P(R_0>K)+4\int_K^L t^3\mathbb P(R_0>t)\,\dd t
 \le\Theta(K)+4\Theta(K)\int_K^L\frac{\dd t}{t}.
\]
This proves the first inequality. The pointwise bound $R_0^5\mathbf1_{\{R_0\le K\}}\le K^3R_0^2$ and $\mathbb ER_0^2=1$ prove the second. The sharper asymptotic
$\mathbb E[R_0^5\mathbf1_{\{R_0\le K\}}]=o(K)$ is \eqref{W:eq:truncated-high-moment} with $q=5$.
\end{proof}

The threshold chosen in Lemma~\ref{W:lem:admissible-scale} also gives
\begin{equation}
 \mathbb P(H_N\ne H_N^{\mathrm{tr}})
 \le C e_N^{-4}\Theta(L_N)+e_N^{-2}\Theta_{\mathrm d}(L_N)=o(1).
 \label{W:eq:truncated-moment-equality-error}
\end{equation}
Indeed, the two terms are the off-diagonal and diagonal union bounds in that lemma.

\Needspace{6\baselineskip}

Let $M_{2,\mathrm d}=\sup_{N,i}\mathbb E d_{N,i}^2<\infty$. For $p>2$ define
\begin{equation}
 v_{p,\mathrm d}(L)
 =\sup_{N,i}L^{-(p-2)}
    \mathbb E[|d_{N,i}|^p\mathbf1_{\{|d_{N,i}|\le L\}}].
 \label{W:eq:truncated-moment-diagonal-vp}
\end{equation}

\begin{lemma}
\label{W:lem:truncated-moment-diagonal-moments}
Then $v_{p,\mathrm d}(L)\to0$. For $\widehat d=d\mathbf1_{\{|d|\le L\}}$,
\begin{align}
 |\mathbb E\widehat d|&\le L^{-1}\Theta_{\mathrm d}(L),
       \label{W:eq:truncated-moment-diagonal-mean}\\
 \mathbb E\widehat d^2&\le M_{2,\mathrm d},
       \label{W:eq:truncated-moment-diagonal-second}\\
 \mathbb E|\widehat d|^p&\le L^{p-2}v_{p,\mathrm d}(L),
       \label{W:eq:truncated-moment-diagonal-high}\\
 \mathbb E[|\widehat d|^4\mathbf1_{\{|\widehat d|>K\}}]
       &\le L^2\Theta_{\mathrm d}(K),\quad K<L.
       \label{W:eq:truncated-moment-diagonal-fourth-tail}
\end{align}
The only moment hypothesis here is uniform integrability of the diagonal squares.
\end{lemma}
\begin{proof}
Since $\mathbb Ed=0$, the first estimate follows from $\mathbb E\widehat d=-\mathbb E[d\mathbf1_{\{|d|>L\}}]$ and $|d|\mathbf1_{\{|d|>L\}}\le L^{-1}d^2\mathbf1_{\{|d|>L\}}$. For
$1\le S<L$,
\[
 L^{-(p-2)}\mathbb E[|d|^p\mathbf1_{\{|d|\le L\}}]
 \le (S/L)^{p-2}M_{2,\mathrm d}+\Theta_{\mathrm d}(S).
\]
Taking $\limsup_{L\to\infty}$ and then $S\to\infty$ gives zero uniformly over the diagonal laws. This proves \eqref{W:eq:truncated-moment-diagonal-vp} and
\eqref{W:eq:truncated-moment-diagonal-high}. Finally, $|d|^4\mathbf1_{\{K<|d|\le L\}}\le L^2d^2\mathbf1_{\{|d|>K\}}$, so taking expectations gives the last inequality. The second
moment bound follows directly from
\[
 \mathbb E\widehat d^2
 =\mathbb E[d^2\mathbf1_{\{|d|\le L\}}]
 \le\mathbb Ed^2\le M_{2,\rm d}.
\]
\end{proof}

Let $\mu_N=\mathbb E\widehat Z$, $\Sigma_N=\operatorname{Cov}(\widehat Z)$, and $Z_N^{\mathrm{nor}}=\Sigma_N^{-1/2}(\widehat Z-\mu_N)$. Here $d_\alpha=1$ or $2$, and the
untruncated real covariance is $I_{d_\alpha}$. The normalization uses the full covariance of the real coordinate vector. Use the bounded random vector $Y_N$ of
Lemma~\ref{W:lem:bounded-complex-moment-matching}: it is centered, has covariance $I_{d_\alpha}$, has the third moment tensor of $Z_N^{\mathrm{nor}}$, and has a support radius
bounded independently of $N$.

For a random vector $X$, write $M_q(X)=\mathbb E X^{\otimes q}$ and use the ordered tensor norm of \eqref{W:eq:B-tensor-l1}. Put $\Delta_q=M_q(\widehat Z)-M_q(Y_N)$.

\Needspace{6\baselineskip} We next quantify the first four moment differences between the uncentered truncated block and the bounded random vector, including the effects of
centering and covariance normalization.

\begin{lemma}
\label{W:lem:truncated-bounded-moment-differences}
For all sufficiently large $N$, with $L=L_N$,
\begin{align}
 \|\Delta_1\|_{1,\mathrm{ord}}&\le C \Theta(L)L^{-3},
       \label{W:eq:truncated-moment-off-defect-one}\\
 \|\Delta_2\|_{1,\mathrm{ord}}&\le C \Theta(L)L^{-2},
       \label{W:eq:truncated-moment-off-defect-two}\\
 \|\Delta_3\|_{1,\mathrm{ord}}&\le C \Theta(L)L^{-2},
       \label{W:eq:truncated-moment-off-defect-three}\\
 \|\Delta_4\|_{1,\mathrm{ord}}&\le C\{1+h_4(L)\}.
       \label{W:eq:truncated-moment-off-defect-four}
\end{align}
The covariance estimate also bounds the complex pseudocovariance. The cubic and fourth tensors are the ordered coordinate tensors defined above, including products containing both
real coordinates.
\end{lemma}
\begin{proof}
The first estimate follows from $\mathbb EZ=0$ and \eqref{W:eq:quant-tail-below-four}. Since $M_2(Z)=I_{d_\alpha}$,
\[
 M_2(\widehat Z)-I_{d_\alpha}
 =-\mathbb E[Z^{\otimes2}\mathbf1_{\{R_0>L\}}].
\]
This is the \emph{noncentral} second tensor in Taylor's formula, including the contribution of the mean. It proves \eqref{W:eq:truncated-moment-off-defect-two}. Consequently
\[
 \|\Sigma_N-I\|\le C \Theta(L)L^{-2}+C \Theta(L)^2L^{-6}
                   \le C \Theta(L)L^{-2}.
\]
All norms on these fixed dimensional spaces are equivalent. Thus $\Sigma_N$ is positive definite for large $N$, and its square root $S_N=\Sigma_N^{1/2}$ satisfies
$\|S_N-I\|\le C \Theta(L)L^{-2}$. In fact, for $\|\Sigma_N-I\|\le1/2$, spectral calculus gives
\[
 \|S_N-I\|=\max_{\lambda\in\sigma(\Sigma_N)}
       \frac{|\lambda-1|}{\sqrt\lambda+1}
 \le\|\Sigma_N-I\|,\qquad
 \|S_N\|\le\sqrt{3/2},\quad\|S_N^{-1}\|\le\sqrt2.
\]
These inequalities also yield $\mathbb E\|Z_N^{\rm nor}\|_2^3 \le4\|S_N^{-1}\|^3(\mathbb E\|\widehat Z\|_2^3+\|\mu_N\|_2^3)$. The right side is bounded uniformly in $N$.

The identity $\widehat Z=\mu_N+S_NZ_N^{\mathrm{nor}}$ and $\mathbb EZ_N^{\mathrm{nor}}=0$, $M_2(Z_N^{\mathrm{nor}})=I$ give
\begin{equation}
 \Delta_3=(S_N^{\otimes3}-I)M_3(Z_N^{\mathrm{nor}})
   +3\operatorname{Sym}(\mu_N\otimes\Sigma_N)
   +\mu_N^{\otimes3}.
 \label{W:eq:truncated-moment-cubic-affine}
\end{equation}
Here $\operatorname{Sym}$ is the average over the three distinct placements of the factor with one tensor index. The absolute third moment of $Z_N^{\mathrm{nor}}$ is uniformly
bounded by Lemma~\ref{W:lem:truncated-moment-tail-moments} and the uniform bounds on $S_N^{\pm1}$. The three terms on the right of \eqref{W:eq:truncated-moment-cubic-affine} are
bounded by $C \Theta(L)L^{-2}$, $C \Theta(L)L^{-3}$, and $C \Theta(L)^3L^{-9}$. This proves \eqref{W:eq:truncated-moment-off-defect-three}.

Finally, $\|M_4(\widehat Z)\|_{1,\mathrm{ord}} \le\mathbb E\|\widehat Z\|_1^4\le16h_4(L)$, whereas $Y_N$ has a bounded fourth moment. This proves the fourth estimate.

\end{proof}

We truncate the standardized Gaussian entries so that each entry replacement path lies in a compact set. If $G_\alpha$ is its standardized block, define
\begin{equation}
 G_\alpha^{[K]}=G_\alpha
                 \mathbf1_{\{r_\alpha(G_\alpha)\le K_N\}}.
 \label{W:eq:truncated-moment-Gaussian-clip}
\end{equation}
The truncation is applied separately to each entry, so the joint distribution is the product of the distributions of $G_\alpha^{[K]}$. Off diagonal,
$G_\alpha\sim N(0,I_{d_\alpha})$ in the real coordinate convention above; on the diagonal its variance is $2/\beta$. Let $H_N^{\mathrm{ub},[K]}$ have off-diagonal blocks $Y_N$ and
diagonal blocks $G_\alpha^{[K]}$. Let $H_N^{\mathrm G,[K]}$ have all blocks truncated as in \eqref{W:eq:truncated-moment-Gaussian-clip}. These two arrays have identical diagonal
laws.

\Needspace{6\baselineskip} Symmetry gives zero odd moments for $G_\alpha^{[K]}$. The next estimates calculate its covariance and the probability that at least one entry has
changed.

\begin{lemma}
\label{W:lem:truncated-moment-Gaussian-clipping}
For each fixed $0\le q\le5$, there are constants $c>0$ and $C_q$ such that
\begin{equation}
 \mathbb E[\|G_\alpha\|_1^q
             \mathbf1_{\{r_\alpha(G_\alpha)>K\}}]
 \le C_q(1+K^{q+2})\mathrm e^{-cK^2}.
 \label{W:eq:truncated-moment-Gaussian-tail}
\end{equation}
The first and third tensors of $G_\alpha^{[K]}$ vanish exactly, but its second tensor differs from that of $G_\alpha$ by the nonzero tail tensor. In particular, the difference of
the off-diagonal moment tensors between $Y_N$ and $G_\alpha^{[K]}$ satisfies
\begin{equation}
 \Delta^{\mathrm{ub},G}_1=0,\quad
 \|\Delta^{\mathrm{ub},G}_2\|_{1,\mathrm{ord}}
     \le C(1+K_N^4)\mathrm e^{-cK_N^2},\quad
 \|\Delta^{\mathrm{ub},G}_q\|_{1,\mathrm{ord}}\le C_q
       \quad(q=3,4).
 \label{W:eq:truncated-moment-bounded-Gaussian-defects}
\end{equation}
Moreover
\begin{equation}
 \mathbb P({H_N^{\mathrm G}}\ne H_N^{\mathrm G,[K]})
       \le CN^2\mathrm e^{-cK_N^2}.
 \label{W:eq:truncated-moment-Gaussian-equality}
\end{equation}
The same estimate, with $N$ in place of $N^2$, couples the matrix with bounded standardized off-diagonal entries and Gaussian diagonal entries to $H_N^{\mathrm{ub},[K]}$.
\end{lemma}
\begin{proof}
Integrating the Gaussian density in polar coordinates in dimension at most two reduces the first estimate to $C\int_{c_0K}^\infty t^{q+1}\mathrm e^{-c_1t^2}\,\dd t$. Indeed, with
$a=c_0K$ and $C_{q,c_1}=\sup_{t\ge0}t^q\mathrm e^{-c_1t^2/2}<\infty$,
\[
 \int_a^\infty t^{q+1}\mathrm e^{-c_1t^2}\,\dd t
 \le C_{q,c_1}\int_a^\infty t\mathrm e^{-c_1t^2/2}\,\dd t
 =\frac{C_{q,c_1}}{c_1}\mathrm e^{-c_1a^2/2}.
\]
Norm equivalence in dimension at most two makes the constants uniform in the real and complex cases. This proves \eqref{W:eq:truncated-moment-Gaussian-tail}. The entry truncation
set is invariant under $x\mapsto-x$, so odd tensors vanish. The identity
\[
 M_2(G_\alpha^{[K]})
   =M_2(G_\alpha)
       -\mathbb E[G_\alpha^{\otimes2}
                  \mathbf1_{\{r_\alpha(G_\alpha)>K_N\}}]
\]
proves the second moment tensor statement. For each fixed $q$, bounded support and the same Gaussian integral give
\[
 \mathbb E\|Y_N\|_1^q\le C_Y^q,\qquad
 \mathbb E\|G_\alpha^{[K]}\|_1^q
 \le\mathbb E\|G_\alpha\|_1^q\le C_q,
\]
uniformly for $d_\alpha\in\{1,2\}$. Finally, couple each truncated block to the original block by the displayed map and take a union bound. For every fixed $p$,
$\sup_{x\ge1}(1+x)^p e^{-cx^2/2}<\infty$, so $(1+K_N)^p e^{-cK_N^2}\le C_{p,c}e^{-cK_N^2/2}$.
\end{proof}

\Needspace{6\baselineskip} A centered complex block with identity covariance can still have a non-Gaussian mixed fourth moment, as the following example shows.

\begin{remark}
\label{W:rem:C-dependent-example}
Let $X$ take each of the four values
\[
 (\sqrt2,0),\qquad(-\sqrt2,0),\qquad
 (0,\sqrt2),\qquad(0,-\sqrt2)
\]
with probability $1/4$. It has mean zero and covariance $I_2$, but its coordinates are dependent. Here $\mathbb EX_1^2X_2^2=0$, while $\mathbb EG_1^2G_2^2=1$ for $G\sim N(0,I_2)$.
Thus the fourth-order joint entry comparison contains $-\tfrac14\partial_{\alpha,1}^2\partial_{\alpha,2}^2\Psi$. The coefficient is $-1/(2!2!)$, exactly as in
\eqref{W:eq:B-complex-fourth-real}; separate marginal moment matching would not justify deleting it.
\end{remark}

Write $f_{\alpha,K}(x)=\Psi(K+N^{-1/2}\mathcal T_\alpha x)$, with $K$ a matrix with the selected entry set to zero. For a fixed $x$, set $g_x(s)=f_{\alpha,K}(sx)$. All derivatives
here use the real coordinate block defined above.

\Needspace{6\baselineskip} Splitting the radial remainder of the expansion at a second radius allows different derivative estimates to be used without changing the polynomial of
degree three.

Let $f$ be $C^5$ on the segments $\{tx:0\le t\le1\}$ under consideration, and set
\[P_4f(x)=\sum_{q=0}^4D_x^q f(0)[x^{\otimes q}]/q!.\]

\begin{lemma}
\label{W:lem:truncated-moment-two-radius-exact}
One has $f=P_4f+R_{\le K}+R_{>K}$, where
\begin{align}
 R_{\le K}(f,x)
 &=\frac{\mathbf1_{\{r_\alpha(x)\le K\}}}{4!}
       \int_0^1(1-s)^4D_{x}^5f(sx)[x^{\otimes5}]\,\dd s,
       \label{W:eq:truncated-moment-small-exact}\\
 R_{>K}(f,x)
 &=\frac{\mathbf1_{\{r_\alpha(x)>K\}}}{3!}
       \int_0^1(1-s)^3
        \bigl(D_{x}^4f(sx)-D_{x}^4f(0)\bigr)[x^{\otimes4}]\,\dd s.
       \label{W:eq:truncated-moment-large-exact}
\end{align}
For two joint entry probability laws $\mu,\nu$, the signed polynomial part is exactly
\begin{equation}
 \sum_{q=1}^4\frac1{q!}
 \left\langle M_q(\mu)-M_q(\nu),\mathscr D_{\alpha,q}\Psi(K)\right\rangle.
 \label{W:eq:truncated-moment-full-polynomial}
\end{equation}
The constant term cancels because both entry laws have total mass one. The derivative norm bounds on the remainders are respectively
\begin{equation}
 \frac{\|x\|_1^5}{5!}
   \sup_{0\le s\le1}\max_{\boldsymbol\sigma\in[d_\alpha]^5}
          |\partial_\alpha^{\boldsymbol\sigma}\Psi(\iota_{\alpha,K}(sx))|,
 \qquad
 \frac{\|x\|_1^4}{12}
   \sup_{0\le s\le1}\max_{\boldsymbol\sigma\in[d_\alpha]^4}
          |\partial_\alpha^{\boldsymbol\sigma}\Psi(\iota_{\alpha,K}(sx))|,
 \label{W:eq:truncated-moment-remainder-norms}
\end{equation}
on their respective radial pieces.
\end{lemma}
\begin{proof}
Fix $x$ and put $g_x(t)=f(tx)$. The integral Taylor formulas at orders five and four give, respectively,
\[
 \begin{aligned}
 g_x(1)-\sum_{j=0}^4\frac{g_x^{(j)}(0)}{j!}
  &=\frac1{4!}\int_0^1(1-t)^4g_x^{(5)}(t)\,\dd t,\\
 g_x(1)-\sum_{j=0}^4\frac{g_x^{(j)}(0)}{j!}
  &=\frac1{3!}\int_0^1(1-t)^3
             \bigl(g_x^{(4)}(t)-g_x^{(4)}(0)\bigr)\,\dd t.
 \end{aligned}
\]
Here $g_x^{(q)}(t)=D_x^qf(tx)[x,\ldots,x]$. Use the first identity on $r_\alpha(x)\le K$ and the second on its complement. This proves \eqref{W:eq:truncated-moment-small-exact} and
\eqref{W:eq:truncated-moment-large-exact}. For $q=4,5$, the coordinate expansion gives
\[
 |D_x^qf(tx)[x,\ldots,x]|
 \le \sum_{i_1,\ldots,i_q}|x_{i_1}\cdots x_{i_q}|
                   \max_{\boldsymbol i}|\partial_{\boldsymbol i}f(tx)|
 =\|x\|_1^q\max_{\boldsymbol i}|\partial_{\boldsymbol i}f(tx)|.
\]
The integral masses are $1/5!$ and $1/4!$; the second identity contains a difference of two fourth derivatives. These observations prove
\eqref{W:eq:truncated-moment-remainder-norms}.

For $f(x)=\Psi(\iota_{\alpha,K}(x))$, integrate the displayed polynomial against $\mu-\nu$. In the notation of the statement this gives
\[
 \begin{aligned}
 \int P_4f(x)(\mu-\nu)(\dd x)
 &=\sum_{q=0}^4\frac1{q!}\sum_{i_1,\ldots,i_q}
     \partial_{i_1}\cdots\partial_{i_q}f(0)
     \int x_{i_1}\cdots x_{i_q}(\mu-\nu)(\dd x)\\
 &=\sum_{q=1}^4\frac1{q!}
      \langle M_q(\mu)-M_q(\nu),\mathscr D_{\alpha,q}\Psi(K)\rangle.
 \end{aligned}
\]
The term $q=0$ vanishes because both measures have mass one. The polynomial is integrated against the whole entry law; only the integral remainder is split at $K$.
\end{proof}

\subsubsection{Bounded entries with matching third moments}
\label{W:sec:bounded-entries-with-prescribed-moments}
A distribution with uniformly bounded entries provides an intermediate comparison between the truncated matrix and the Gaussian ensemble. It must match the third moments as well as
the first two, so that the first comparison has a small cubic error. The construction below matches the centered, covariance normalized truncation. The difference between its third
tensor and that of the uncentered truncated entry is then estimated explicitly.

All third tensors use standardized real coordinates and the ordered norm $\|M\|_{1,\rm ord}=\sum_{a,b,c}|M_{abc}|$. The original coordinate vector $X$ has dimension one or two,
$\mathbb EX=0$, and $\mathbb EXX^{\mathsf T}=I$. Its radial variable $r(X)$ is the one from Section~\ref{W:sec:moment-bounds-and-taylor-expansion-at-two-truncation-scales}, so that
$r(X)=|\xi|$, $\|X\|_1\le Cr(X)$, and its tail function is \eqref{W:eq:offdiagonal-tail-function}. At the threshold $L=L_N$ from \eqref{W:eq:LN-definition}, the truncation and its
normalization are

\begin{equation}
 \widehat X=X\mathbf1_{\{r(X)\le L\}},\qquad
 \mu=\mathbb E\widehat X,\qquad
 \Sigma=\operatorname{Cov}(\widehat X),\qquad
 X^{\mathrm{nor}}=\Sigma^{-1/2}(\widehat X-\mu).
 \label{W:eq:normalized-truncated-entry-vector}
\end{equation}
For all sufficiently large $N$, $\Sigma$ is uniformly positive definite.

\Needspace{6\baselineskip} Matching the third tensor after centering and covariance normalization leaves a small, explicitly controlled third moment difference for the original
truncated block.

\begin{lemma}
\label{W:lem:cubic-defect}
There is a uniformly bounded random vector $Y$ with the same mean, covariance, and full third tensor as $X^{\mathrm{nor}}$. Moreover,
\begin{equation}
 \tau_{3,N}:=\|M_3(\mathcal L(\widehat X))-M_3(\mathcal L(Y))\|_{1,\mathrm{ord}}
 \le C \Theta(L)L^{-2}.
 \label{W:eq:cubic-defect}
\end{equation}
The bound is uniform over the block family when its original tail and third moment bounds are uniform. Higher moments are estimated after truncation.
\end{lemma}
\begin{proof}
Apply \eqref{W:eq:quant-tail-below-four} with $p=1,2,3$. Consequently
\[
 \|\mu\|_1\le C\Theta(L)L^{-3},\qquad
 \|\Sigma-I\|\le C\Theta(L)L^{-2},
 \qquad\mathbb E\|X^{\mathrm{nor}}\|_1^3\le C.
\]
Let $S=\Sigma^{1/2}$. The relation $\widehat X=\mu+SX^{\mathrm{nor}}$ yields
\begin{equation}
 M_3(\widehat X)-M_3(X^{\mathrm{nor}}) ={}(S^{\otimes3}-I)M_3(X^{\mathrm{nor}}) +3\operatorname{Sym}(\mu\otimes\Sigma)+\mu^{\otimes3}.
 \label{W:eq:cubic-affine}
\end{equation}
In coordinates, $[3\operatorname{Sym}(\mu\otimes\Sigma)]_{ijk} =\mu_i\Sigma_{jk}+\mu_j\Sigma_{ik}+\mu_k\Sigma_{ij}$.

This is \eqref{W:eq:truncated-moment-cubic-affine} with $\widehat Z=\widehat X$ and $Z_N^{\mathrm{nor}}=X^{\mathrm{nor}}$. The estimates in that proof give
$\|M_3(\widehat X)-M_3(X^{\mathrm{nor}})\|_{1,\mathrm{ord}}\le C\Theta(L)L^{-2}$.
It remains to construct $Y$ with the same first three coordinate moment tensors as $X^{\mathrm{nor}}$.

In dimension one, put $c=\mathbb E(X^{\rm nor})^3$, set $r_\pm=(c\pm\sqrt{c^2+4})/2$, and define
\[
 \mathbb P(Y=r_+)=\frac{-r_-}{r_+-r_-},\qquad
 \mathbb P(Y=r_-)=\frac{r_+}{r_+-r_-}.
\]
Both probabilities are positive and add to one, since $r_-<0<r_+$. Their weighted mean is zero. The identity $Y^2=cY+1$ then gives $\mathbb EY^2=1$ and
$\mathbb EY^3=c\mathbb EY^2+\mathbb EY=c$. Finally $|r_\pm|\le |c|+1$, so the support is uniformly bounded.

The construction uses only zero mean, identity covariance and a uniform bound on the third absolute moment. It therefore applies to any random vector with these properties in
dimension one or two. In dimension two, set $(a,b,c,d)=(M_{111},M_{112},M_{122},M_{222})$ for $M=M_3(X^{\mathrm{nor}})$. Put
\[
 v_1=(1,0),\ v_2=(0,1),\ v_3=(1,1),\ v_4=(1,-1),
\]
\[
 c_1=a-c,\quad c_2=d-b,\quad c_3=(b+c)/2,\quad c_4=(c-b)/2.
\]
Then $M_3(X^{\mathrm{nor}})=\sum_{l=1}^4c_lv_l^{\otimes3}$ and $\sum_lv_lv_l^{\mathsf T}=3I_2$. Let $\theta=1/6$, and take independent centered two point $Z_l$ with variance
$\theta$ and third moment $c_l$. For each $l$, let $r_+>0>r_-$ be the roots of $r^2-(c_l/\theta)r-\theta=0$ and assign probabilities $-r_-/(r_+-r_-)$ and $r_+/(r_+-r_-)$ to them.
Then $\mathbb EZ_l=0$, $\mathbb EZ_l^2=-r_+r_-=\theta$, and $\mathbb EZ_l^3=\theta(r_++r_-)=c_l$. Independently let ${Z_0}=2^{-1/2}(\varepsilon_1,\varepsilon_2)$, with independent
symmetric signs, and put $Y={Z_0}+\sum_lv_lZ_l$. Its covariance is
\[
 \mathbb EYY^{\mathsf T}=\tfrac12 I_2+
           \theta\sum_{l=1}^4v_lv_l^{\mathsf T}=I_2.
\]
In a third moment, a product involving two distinct independent centered summands has expectation zero. The four surviving tensor coordinates are
\[
 \begin{aligned}
 \mathbb EY_1^3&=c_1+c_3+c_4=a,&
 \mathbb EY_1^2Y_2&=c_3-c_4=b,\\
 \mathbb EY_1Y_2^2&=c_3+c_4=c,&
 \mathbb EY_2^3&=c_2+c_3-c_4=d.
 \end{aligned}
\]
Thus $M_3(Y)=M_3(X^{\mathrm{nor}})$. To give a uniform support bound, let $M$ bound the four absolute tensor entries. Then $|c_l|\le2M$. The roots defining $Z_l$ obey
\[
 |r_\pm|\le |c_l|/\theta+\sqrt\theta
                    \le12M+1,\qquad
 \|Y\|_2\le1+(2+2\sqrt2)(12M+1).
\]
The second inequality follows from $\|Z_0\|_2=1$ and $\sum_{l=1}^4\|v_l\|_2=2+2\sqrt2$. It proves bounded support with a constant depending only on the original third moment bound.
The construction allows dependence between the two coordinates of $X$ and between those of $Y$.
\end{proof}

\subsubsection{Schur complements and conditioning on large entries}
\label{W:sec:block-schur-complements-under-the-tail-assumption}

Large off-diagonal entries determine connected components, which we treat as principal blocks together with their diagonal entries. Conditional on the labels, the entries within
these blocks, and the complementary minor, the remaining cross entries have the product distribution specified below. This permits a conditional estimate for the Schur complement
even when the entries inside a principal block are too large for a scalar estimate.

Those bounds initially depend on the complementary minor. The counting and entrywise estimates below will control that dependence; under their hypotheses,
Proposition~\ref{W:prop:schur-airy-scale} gives the scale $N^{-1/3+o(1)}$. The conditional calculation is thus available before the resolvent bounds that will be deduced from it. A
smooth cutoff will combine these estimates, with the maximum allowed component size chosen after the probability exponent.

For two coordinates, direct inversion gives

\[
 \begin{pmatrix}g^{-1}&-t\\-\overline t&g^{-1}\end{pmatrix}^{-1}
 =\frac1{g^{-2}-|t|^2}
       \begin{pmatrix}g^{-1}&t\\\overline t&g^{-1}\end{pmatrix},
 \qquad R_{11}^{-1}=g^{-1}-g|t|^2.
\]
At a nonreal spectral parameter this identity is used whenever the matrix being inverted is nonsingular. The Schur complement identity for the block inverse instead has the form
\[
 (R_{II})^{-1}+H_{II}
       =zI_I-H_{IJ}(zI_J-H_{JJ})^{-1}H_{JI}.
\]
The right side depends on $H_{IJ}$ and $H_{JJ}$; changing $H_{II}$ leaves it unchanged.

The separation scale used here is
\begin{equation}
 \ell_N^{\mathrm{Sch}}=N^{7/24},\qquad b_N=4\ell_N^{\mathrm{Sch}} N^{-1/2}=4N^{-5/24}.
 \label{W:eq:schur-scales}
\end{equation}
This scale is distinct from the slowly varying inverse cutoff $\ell_N^{\mathrm{slow}}$ in \eqref{W:eq:ell-scale}, the moderate radius $K_N=N^{1/8}$, and the equality radius $L_N$
in \eqref{W:eq:LN-definition}. In particular, for all sufficiently large $N$,
\begin{equation}
 \begin{gathered}
 e_N=(\ell_N^{\mathrm{slow}})^{-1},\qquad
 L_N=\sqrt N/\ell_N^{\mathrm{slow}},\qquad
 K_N<\ell_N^{\mathrm{Sch}},\qquad 2\ell_N^{\mathrm{Sch}}<L_N,\\
 \frac{\ell_N^{\mathrm{Sch}}}{L_N}
   =\ell_N^{\mathrm{slow}}N^{-5/24}\longrightarrow0.
 \end{gathered}
 \label{W:eq:schur-scale-hierarchy}
\end{equation}
The preceding scale hierarchy follows from $\ell_N^{\mathrm{slow}}=N^{o(1)}$; in particular, the graph threshold never replaces $L_N$. We derive the restriction on the exponent
before fixing its value. If the separating radius is $\ell=N^\gamma$ and $\sup_{t\ge\ell}t^4\mathbb P(|\xi|>t)\le C$, then
$\mathbb P(|\xi|>\ell)\le C\ell^{-4}=CN^{-4\gamma}$.
There are $N-1$ possible neighbours of any vertex, so its expected number of labelled large edges is at most $CN^{1-4\gamma}$. A polynomially small component probability therefore
requires $\gamma>1/4$.

At this hypothetical threshold, the support bound would be $b_N=4\ell N^{-1/2}=4N^{\gamma-1/2}$ and the tail integration in \eqref{W:eq:quant-tail-below-four} gives
$u_N\le CN^{-3\gamma}+\varepsilon_{1,N}$ and
$\delta_N\le CN^{-2\gamma}+C(\varepsilon_{2,N}+\varepsilon_{1,N}^2)$.
Assume temporarily the complementary minor bounds stated in \eqref{W:eq:schur-airy-hypotheses}: $\eta=N^{-2/3+o(1)}$, $R_2({\mathcal R})\le N^{1/3+o(1)}$, and
$E_\infty({\mathcal R})\le N^{o(1)}$. The row contribution and the quadratic bounded-entry contribution to \eqref{W:eq:schur-theta} are then bounded respectively by
\[
 \frac{b_N}{\sqrt N}R_2({\mathcal R})
       \le N^{\gamma-2/3+o(1)},\qquad
 b_N^2 E_\infty({\mathcal R})\le N^{2\gamma-1+o(1)}.
\]
The principal mean contribution is
\[
 \eta^{-1}u_Nb_N
 \le N^{2/3-3\gamma+\gamma-1/2+o(1)}
 =N^{1/6-2\gamma+o(1)}.
\]
The terms with $u_NN^{-1/2}$, $u_N^2$, and $\delta_N$ obey smaller bounds in the same range, while the first two terms of $\mathcal K_N$ are $N^{-1/2+o(1)}$ and $N^{-1/3+o(1)}$
under those hypotheses. To make each displayed additional contribution $o(N^{-1/3})$ one needs
\[
 2\gamma-1<-1/3,\qquad
 \gamma-2/3<-1/3,\qquad
 1/6-2\gamma<-1/3,
\]
which are jointly equivalent to $1/4<\gamma<1/3$. Taking $\gamma=7/24$ gives exponents $-5/12$, $-3/8$, and $-5/12$, respectively. To justify this choice, we first establish the
conditional moment estimates and then prove the required minor bounds independently.

Choose a fixed $C^\infty$ function $\chi:[0,\infty)\to[0,1]$ which is zero on $[0,1]$ and one on $[2,\infty)$, and put $c_N(z)=\chi(|z|/\ell_N^{\mathrm{Sch}})$. This is smooth in
the real coordinates of $z$, including at zero.

Let $\xi$ have the original off-diagonal law. Use $\Theta$ from \eqref{W:eq:offdiagonal-tail-function} and $h_4$ from
Section~\ref{W:sec:moment-bounds-and-taylor-expansion-at-two-truncation-scales}. The estimates \eqref{W:eq:quant-tail-below-four} and \eqref{W:eq:truncated-fourth-moment} control
its discarded moments and truncated fourth moment. For the other entry distributions used in the comparison, choose common upper bounds in these formulas. The centered base family
with normalized second moments consists of independent standardized off-diagonal entries satisfying
\begin{equation}
 \begin{split}
 &\mathbb E\xi_{ij}=0,\quad \mathbb E|\xi_{ij}|^2=1,\quad
 \sup_{i<j}\sup_{t\ge\ell_N^{\mathrm{Sch}}}t^4\mathbb P(|\xi_{ij}|>t)\le a_N\le C,\\
 &F_N:=C\left(1+\sup_{i<j}\mathbb E[
       |\xi_{ij}|^4\mathbf1_{\{|\xi_{ij}|\le2\ell_N^{\mathrm{Sch}}\}}]\right)=N^{o(1)}.
 \end{split}
 \label{W:eq:schur-uniform-family}
\end{equation}
The constant in $F_N$ is fixed large enough below. The uncentered truncated law carries the mean and second moment errors computed above. We therefore also allow the same tail and
truncated fourth moment bounds as in \eqref{W:eq:schur-uniform-family}, with its first two identities replaced by
\begin{equation}
 \begin{gathered}
 \sup_{i<j}|\mathbb E\xi_{ij}|\le\varepsilon_{1,N},\qquad
 \sup_{i<j}|\mathbb E|\xi_{ij}|^2-1|\le\varepsilon_{2,N},\\
 \varepsilon_{1,N}\le N^{-3/2+o(1)},\qquad
 \varepsilon_{2,N}\le N^{-1+o(1)}.
 \end{gathered}
 \label{W:eq:schur-defect-family}
\end{equation}
These are noncentral second moment bounds; the square of the mean is subtracted explicitly below. The normalized base family is the case $\varepsilon_{1,N}=\varepsilon_{2,N}=0$.
The bounds in \eqref{W:eq:schur-defect-family} allow nonzero means and covariance errors at every off-diagonal position.

For the uncentered truncated entry, the tail and fourth moment bounds can be taken as $\Theta(\ell_N^{\mathrm{Sch}})$ and $C(1+h_4(2\ell_N^{\mathrm{Sch}}))$. The mean and second
moment differences are verified explicitly in \eqref{W:eq:schur-truncated-defects} below. Centered bounded entries with the required covariance, Gaussian entries, and radially
truncated Gaussian entries are included by enlarging the common deterministic bounds as described there. For any deterministic choice of these entry laws, construct the
upper-triangular blocks independently; at an off-diagonal complex position the pair $(\operatorname{Re}\xi_{ij},\operatorname{Im}\xi_{ij})$ is sampled as a single two-dimensional
random vector. Choose the entire block law independently for each unordered pair, or choose its distribution from the stated admissible families by auxiliary selections that are
independent for distinct unordered pairs. After conditioning on those selections, the bounds on the mean, second moments, radial tail and truncated fourth moment are the same
deterministic bounds at every position. The proofs below are conditional on these selected laws and hence remain uniform when the selections are integrated out. Adjoin independent
auxiliary labels by
\begin{equation}
 \mathbb P(L_{ij}=1\mid\xi_{ij})=c_N(\xi_{ij}).
 \label{W:eq:schur-label}
\end{equation}
The pairs $(\xi_{ij},L_{ij})$ are independent over $i<j$. If $p_{ij}=\mathbb E c_N(\xi_{ij})$, the conditional law given a zero label is exactly
\begin{equation}
 \mu_{ij}^{(0)}(\dd z)=
       \frac{1-c_N(z)}{1-p_{ij}}\,\mu_{ij}(\dd z).
 \label{W:eq:schur-small-law}
\end{equation}
The conditional law given a label equal to one has density $c_N/p_{ij}$ when $p_{ij}>0$. An arbitrary choice on a zero probability label has no effect.

For a deterministic label configuration $l=(l_{ij})$ of positive probability and bounded measurable functions $f_{ij}$, independence gives
\[
 \begin{aligned}
 &\mathbb E\left[\prod_{i<j}f_{ij}(\xi_{ij})\,
                    \mathbf1_{\{L=l\}}\right]\\
 &\qquad=\prod_{i<j}\int f_{ij}(z)
             c_N(z)^{l_{ij}}(1-c_N(z))^{1-l_{ij}}\,\mu_{ij}(\dd z),\\
 &\mathbb E\left[\prod_{i<j}f_{ij}(\xi_{ij})\mid L=l\right]
       =\prod_{i<j}\int f_{ij}(z)\,\mu_{ij}^{(l_{ij})}(\dd z).
 \end{aligned}
\]
Thus the entries remain independent after the labels are fixed. If $i$ and $j$ lie in different components of $l$, then $l_{ij}=0$, so the corresponding entry has the conditional
law in \eqref{W:eq:schur-small-law}.

\Needspace{6\baselineskip}

Throughout the following estimate assume \eqref{W:eq:schur-uniform-family}, or its version allowing moment errors uniformly over the entry indices \eqref{W:eq:schur-defect-family}
with the same tail and fourth moment bounds. In the exactly normalized case put $\varepsilon_{1,N}=\varepsilon_{2,N}=0$. For all sufficiently large $N$, put
\begin{equation}
 u_N=C\{a_N(\ell_N^{\mathrm{Sch}})^{-3}+\varepsilon_{1,N}\}, \qquad \delta_N=C\{a_N(\ell_N^{\mathrm{Sch}})^{-2} +\varepsilon_{2,N}+\varepsilon_{1,N}^2\}.
 \label{W:eq:schur-mean-and-covariance-error-bounds}
\end{equation}

\begin{lemma}
\label{W:lem:schur-small-moments}
Conditional on the label being zero, let $\mu_{ij}^{(0),1}$ be the mean and let $\sigma_{ij}^2$ be the variance. Then
\begin{equation}
 \begin{split}
 &p_{ij}\le a_N(\ell_N^{\mathrm{Sch}})^{-4},\qquad
 |\mu_{ij}^{(0),1}|\le u_N,\qquad
 |\sigma_{ij}^2-1|\le\delta_N,\\
 &X_{ij}:=N^{-1/2}(\xi_{ij}-\mu_{ij}^{(0),1}),\qquad
 |X_{ij}|\le b_N,\quad
 \mathbb E|X_{ij}|^2\le 2/N,\quad
 \mathbb E|X_{ij}|^4\le F_N/N^2,\\
 &\mathbb E|X_{ij}|^p\le b_N^{p-4}F_N/N^2\quad(p\ge4),\qquad
 \frac{F_N}{N b_N^4}\le1.
 \end{split}
 \label{W:eq:schur-core-moments}
\end{equation}
The last inequality holds for all $N\ge N_0$, with $N_0$ independent of the member of the family. If the complex law also satisfies $|\mathbb E\xi_{ij}^2|\le\varepsilon_{2,N}$
(including the exactly normalized case $\mathbb E\xi_{ij}^2=0$), then $|\mathbb E[(\xi_{ij}-\mu_{ij}^{(0),1})^2\mid L_{ij}=0]| \le C\delta_N$. This quantity is included in the
conditional covariance.
\end{lemma}
\begin{proof}
The label can be nonzero only when $|\xi_{ij}|>\ell_N^{\mathrm{Sch}}$, so $p_{ij}\le a_N(\ell_N^{\mathrm{Sch}})^{-4}=o(1)$ and $1-p_{ij}\ge1/2$ once
$a_N(\ell_N^{\mathrm{Sch}})^{-4}\le1/2$, uniformly in the family. Write $m_{ij}=\mathbb E\xi_{ij}$ and $s_{ij}=\mathbb E|\xi_{ij}|^2$. The conditional law
\eqref{W:eq:schur-small-law} gives the exact formulas
\begin{equation}
 \mu_{ij}^{(0),1} =\frac{m_{ij}-\mathbb E[\xi_{ij}c_N(\xi_{ij})]}{1-p_{ij}}, \qquad \sigma_{ij}^2 =\frac{s_{ij}-\mathbb E[|\xi_{ij}|^2c_N(\xi_{ij})]}{1-p_{ij}} -|\mu_{ij}^{(0),1}|^2.
 \label{W:eq:schur-small-defect-identities}
\end{equation}
For $q=1,2$, the tail integration argument in \eqref{W:eq:quant-tail-below-four} applies uniformly to the family and gives
\[
 \mathbb E[|\xi_{ij}|^q c_N(\xi_{ij})]
 \le \frac4{4-q}a_N(\ell_N^{\mathrm{Sch}})^{q-4}.
\]
Thus $|\mu_{ij}^{(0),1}|\le2\varepsilon_{1,N} +C a_N(\ell_N^{\mathrm{Sch}})^{-3}\le u_N$. Subtracting one in the second formula of \eqref{W:eq:schur-small-defect-identities} yields
\[
 |\sigma_{ij}^2-1|
 \le 2\varepsilon_{2,N}+2p_{ij}
       +C a_N(\ell_N^{\mathrm{Sch}})^{-2}
       +|\mu_{ij}^{(0),1}|^2
 \le\delta_N.
\]
For the last inequality use $|\mu_{ij}^{(0),1}|^2\le8\varepsilon_{1,N}^2 +C a_N^2(\ell_N^{\mathrm{Sch}})^{-6}$, and then fix the constant in $\delta_N$ large enough. Also
$p_{ij}\le a_N(\ell_N^{\mathrm{Sch}})^{-2}$ and $a_N^2(\ell_N^{\mathrm{Sch}})^{-6} \le C a_N(\ell_N^{\mathrm{Sch}})^{-2}$, since $a_N\le C$ and $\ell_N^{\mathrm{Sch}}\ge1$. For the
complex pseudocovariance the exact expression is
\[
 \mathbb E[(\xi_{ij}-\mu_{ij}^{(0),1})^2\mid L_{ij}=0]
 =\frac{\mathbb E\xi_{ij}^2-\mathbb E[\xi_{ij}^2c_N(\xi_{ij})]}{1-p_{ij}}
       -(\mu_{ij}^{(0),1})^2.
\]
Its absolute value is bounded by $C\delta_N$ using the same radial second moment tail estimate. These are the conditional mean and variance of the unchanged matrix law. The
conditional law given a zero label is supported on $|\xi|\le2\ell_N^{\mathrm{Sch}}$, and its mean has modulus $o(1)$. Consequently its centered matrix entry has modulus at most
$b_N$. Write $\widetilde h_{ij}=N^{-1/2}(\xi_{ij}-\mu_{ij}^{(0),1})$ conditional on the label being zero. Then
\[
 \mathbb E_0|\widetilde h_{ij}|^4
 \le\frac8{N^2}\left(
  \frac{\mathbb E[|\xi_{ij}|^4\mathbf1_{\{|\xi_{ij}|\le2\ell_N^{\rm Sch}\}}]}
       {1-p_{ij}}+|\mu_{ij}^{(0),1}|^4\right)
 \le\frac{F_N}{N^2},
\]
with the fixed enlargement in the definition of $F_N$. Since $|\widetilde h_{ij}|\le b_N$, for every $p\ge4$,
\[
 \mathbb E_0|\widetilde h_{ij}|^p
 \le b_N^{p-4}\mathbb E_0|\widetilde h_{ij}|^4
 \le b_N^{p-4}F_NN^{-2}.
\]
Moreover $Nb_N^4=4^4N^{1/6}$. Choose $N$ large enough that $F_N\le4^4N^{1/6}$; then $F_N/(Nb_N^4)\le1$, as required in the conditional moment inequalities. Under
\eqref{W:eq:schur-defect-family} the enlarged parameters still satisfy
\begin{equation}
 u_N\le N^{-7/8+o(1)},\qquad
 \delta_N\le N^{-7/12+o(1)}.
 \label{W:eq:schur-defect-power-bounds}
\end{equation}
Thus the bounds for the differences in moment tensors, uniform in the entry indices, preserve the powers used in the subsequent conditional concentration and estimates at the
spectral edge.
\end{proof}

View $L$ as the graph on $\{1,\ldots,N\}$ with edges $L_{ij}=1$. For a fixed integer $d\ge7$, let $\mathcal A_d$ be the event that every connected component has at most $d$
vertices.

\Needspace{6\baselineskip}

\begin{lemma}
\label{W:lem:schur-components}
For a constant depending only on $d$ and the bound in \eqref{W:eq:schur-uniform-family},
\begin{equation}
 \mathbb P(\mathcal A_d^c)
 \le C_d N^{1-d/6}.
 \label{W:eq:schur-components-prob}
\end{equation}
Conditional on the graph and all entries inside its principal components, the entries between distinct components remain independent and have the conditional laws in
\eqref{W:eq:schur-small-law}. All diagonal entries are among the conditioned entries.
\end{lemma}
\begin{proof}
A connected component with at least $d+1$ vertices contains a connected set of exactly $d+1$ vertices and a spanning tree on that set. For each set there are $(d+1)^{d-1}$ labelled
trees. Each tree occurs with probability at most $(C N^{-7/6})^d$. The union bound is therefore at most
\[
 \binom N{d+1}(d+1)^{d-1}(C N^{-7/6})^d
       \le C_d N^{1-d/6}.
\]
The tree edges concern distinct unordered entries, so their probabilities multiply. For $N<d+1$ the event is empty. For a graph configuration $\ell$ of positive probability, the
product identity established before Lemma~\ref{W:lem:schur-small-moments} gives the conditional distribution $\bigotimes_{e\in\mathcal E_N}\nu_e^{\ell_e}$, where
$\mathcal E_N=\{\{i,j\}:i<j\}$ and $\nu_e^a=\mathcal L(\xi_e\mid L_e=a)$.

For this fixed graph let $J(\ell)$ be its pairs in different components and $I(\ell)$ its remaining off-diagonal pairs. Every $e\in J(\ell)$ has $\ell_e=0$. The conditional product
law gives, for bounded Borel $F$ of the internal entries and diagonal vector $d$,
\[
 \mathbb E\left[F((\xi_e)_{e\in I(\ell)},d)
                \prod_{e\in J(\ell)}f_e(\xi_e)\,\middle|\,L=\ell\right]
 =\mathbb E[F((\xi_e)_{e\in I(\ell)},d)\mid L=\ell]
                \prod_{e\in J(\ell)}\int f_e\,\dd\nu_e^0.
\]
The diagonal entries are independent of all off-diagonal pairs $(\xi_e,L_e)$ and are included in $F$. The defining property of conditional expectation now yields
\[
 \mathbb E\left[\prod_{e\in J(L)}f_e(\xi_e)\,
          \middle|\,L,(\xi_e)_{e\in I(L)},d\right]
       =\prod_{e\in J(L)}\int f_e\,\dd\nu_e^0.
\]
The conditional distribution given $L_e=0$ is the law in \eqref{W:eq:schur-small-law}. This proves the claimed conditional independence and identifies each remaining entry
distribution.
\end{proof}

This construction also gives the following decomposition of $H$. Let $A$ equal $H$ inside each component and the conditional entry means between components. Then $H=A+X$, with $X$
centered, zero on the internal blocks, and independent between blocks conditional on the internal data. On the truncation event from
Section~\ref{W:sec:moment-bounds-and-taylor-expansion-at-two-truncation-scales},
\begin{equation}
 \|A\|\le d e_N+C\sqrt N\,u_N
             \le d e_N+N^{-3/8+o(1)}=o(1)
 \quad\hbox{on }\mathcal A_d.
 \label{W:eq:schur-deformation-norm}
\end{equation}
Indeed the internal block has at most $d$ entries of modulus $e_N$ per row, and the cross component mean matrix has row sum at most $\sqrt N\,u_N$. At the $N^{-2/3}$ spectral
scale, the conditional variance correction $O(N^{-7/12+o(1)})$ must be included in the self-consistent equation at the edge.

\subsubsection{Conditional moments of the Schur complement}
\label{W:sec:conditional-moments-of-schur-complement}

Conditioning on a complementary minor expresses the error in the Schur complement through linear and quadratic expressions in independent entries. The next two lemmas estimate
those expressions with explicit dependence on the moment order. Polynomial growth of the constants in that order is essential: choosing the order logarithmically in $N$ will later
turn the conditional moment bounds into sufficiently small failure probabilities for the smooth cutoff.

\begin{lemma}
\label{W:lem:schur-rosenthal}
Let $Z_1,\ldots,Z_n$ be independent centered random vectors in a real or complex Hilbert space, and let $p\ge4$ be even. Then
\begin{equation}
 \left\|\sum_j Z_j\right\|_{L^p}
 \le Cp\left\{\left(\sum_j\mathbb E\|Z_j\|^2\right)^{1/2}
        +\left(\sum_j\mathbb E\|Z_j\|^p\right)^{1/p}\right\}.
 \label{W:eq:schur-rosenthal}
\end{equation}
Consequently, if independent centered complex scalars $X_j$ satisfy \eqref{W:eq:schur-core-moments}, then for deterministic Hilbert space vectors $v_j$,
\begin{equation}
 \left\|\sum_j X_jv_j\right\|_{L^p}
 \le Cp\left\{N^{-1/2}\Bigl(\sum_j\|v_j\|^2\Bigr)^{1/2}
                         +b_N\max_j\|v_j\|\right\}.
 \label{W:eq:schur-linear}
\end{equation}
Both statements hold with conjugated scalars as well.
\end{lemma}
\begin{proof}
We specify the classical inequality used here. For independent mean-zero vectors $Z_j$ in a real Hilbert space and $p\ge2$, Theorem~5.3, equation~(5.8), of
\cite{pinelis1994optimum} gives
\[
 \left\|\sum_jZ_j\right\|_{L^p}
 \le C\left[\sqrt p\left(\sum_j\mathbb E\|Z_j\|^2\right)^{1/2}
       +p\left(\sum_j\mathbb E\|Z_j\|^p\right)^{1/p}\right].
\]
The constant is independent of $p$, the number of summands and the Hilbert space dimension. For complex vectors we use the underlying real Hilbert space with inner product
$\operatorname{Re}\langle u,v\rangle$; its norm is unchanged. Since $p\ge4$, $\sqrt p\le p$, and the first displayed inequality implies \eqref{W:eq:schur-rosenthal}.

For the vector-valued moment bound, put $Z_j=X_jv_j$ and condition, if applicable, on the label graph and all prescribed internal entries. Under this conditional law, the $X_j$ are
independent, centered and satisfy the common bounds \eqref{W:eq:schur-core-moments}. Thus
\[
 \sum_j\mathbb E\|Z_j\|^2
 \le\frac2N\sum_j\|v_j\|^2.
\]
The number of summands is at most $N$. For every even $p\ge4$, the high moment estimate and the condition $F_N/(Nb_N^4)\le1$ give
\[
 \begin{aligned}
 \sum_j\mathbb E\|Z_j\|^p
 &\le\frac{F_N}{N^2}b_N^{p-4}\sum_j\|v_j\|^p\\
 &\le\frac{F_N}{N}b_N^{p-4}\max_j\|v_j\|^p
 \le b_N^p\max_j\|v_j\|^p.
 \end{aligned}
\]
Substitute these two estimates into the quoted Hilbert space inequality. Since $\sqrt p\le p$ and the factor $\sqrt2$ is absolute, the result is
\[
 \left\|\sum_jX_jv_j\right\|_{L^p}
 \le Cp\left\{N^{-1/2}\left(\sum_j\|v_j\|^2\right)^{1/2}
                      +b_N\max_j\|v_j\|\right\},
\]
which is \eqref{W:eq:schur-linear}. The calculation is valid under each permitted conditioning, since all parameters on its right are deterministic once that conditioning is fixed.
Conjugation changes neither centering nor the second, high moment, or support estimates and so gives the identical bound for $\overline X_j$.

Figure~\ref{W:fig:moment-hilbert-partitions} illustrates the centered index patterns associated with this classical moment inequality. \begingroup
\setlength{\intextsep}{1.2\baselineskip}
\begin{figure}[H]
\centering
\begin{tikzpicture}[x=1cm,y=1cm,every node/.style={mg text}]
 \foreach \shift/\tag in {0/a,7.0/b}{
  \begin{scope}[xshift=\shift cm]
   \foreach \k in {1,...,6}{\node[mg occurrence,label=above:{$\k$}] (\tag\k) at ({(\k-1)*1.04},0) {};}
   \foreach \u/\v in {1/2,3/4,5/6}{\draw[mg kernel] (\tag\u) to[bend left=48] (\tag\v);}
  \end{scope}
 }
 \node[anchor=west] at (0,1.0) {(a) A singleton: zero expectation};
 \node[anchor=west] at (7.0,1.0) {(b) Two blocks of size three};
 \node[mg block,label=below:{$i$}] (ai) at (0,-1.3) {};
 \node[mg block,label=below:{$j$}] (aj) at (1.6,-1.3) {};
 \node[mg block,label=below:{$k$}] (ak) at (4.2,-1.3) {};
 \draw[mg incidence] (ai)--(a1);
 \foreach \v in {2,3}{\draw[mg incidence] (aj)--(a\v);}
 \foreach \v in {4,5,6}{\draw[mg incidence] (ak)--(a\v);}
 \node[mg block,label=below:{$i$}] (bi) at (8.4,-1.3) {};
 \node[mg block,label=below:{$j$}] (bj) at (10.8,-1.3) {};
 \foreach \v in {1,3,5}{\draw[mg incidence] (bi)--(b\v);}
 \foreach \v in {2,4,6}{\draw[mg incidence] (bj)--(b\v);}
 \node at (2.6,-2.2) {$\{1\},\ \{2,3\},\ \{4,5,6\}$};
 \node at (9.6,-2.2) {$\{1,3,5\},\ \{2,4,6\}$};
\end{tikzpicture}
\caption{Index partitions in the expansion of $\|\sum_iZ_i\|^6$ used in \eqref{W:eq:schur-rosenthal}. Numbered dots denote the six factor positions; diamonds and dashed lines specify which positions use the same random vector. Solid arcs denote the Hilbert space inner products. In (a), the variable indexed by $i$ occurs once, and centering makes the expectation zero. In (b), the two indices each occur three times. Summing absolute values and allowing the two indices to coincide gives the bound $(\sum_i\mathbb E\|Z_i\|^3)^2$. Thus blocks of size three contribute to the estimate alongside the even-sized blocks.}
\label{W:fig:moment-hilbert-partitions}
\end{figure}
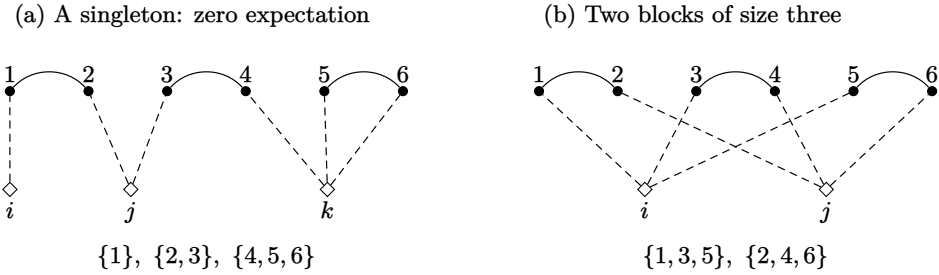
\endgroup
\end{proof}

For a deterministic $n\times n$ matrix $T$, $n\le N$, write
\[
 \begin{split}
 D_2(T)&=\Bigl(\sum_j|T_{jj}|^2\Bigr)^{1/2},\qquad
 E_\infty(T)=\max_{i,j}|T_{ij}|,\\
 R_2(T)&=\max_i\Bigl(\sum_j|T_{ij}|^2\Bigr)^{1/2}
           +\max_j\Bigl(\sum_i|T_{ij}|^2\Bigr)^{1/2},\\
 \mathcal K_N(T)&=\frac{\sqrt{F_N}}N D_2(T)
          +\frac1N\|T\|_{\mathrm{HS}}
          +\frac{b_N}{\sqrt N}R_2(T)+b_N^2E_\infty(T).
 \end{split}
\]

\Needspace{6\baselineskip} For complex vectors, we use $\langle u,v\rangle=\sum_j u_j\overline{v_j}$. Thus $(u,v)\mapsto v^*Tu=\langle Tu,v\rangle$ is a sesquilinear form on
$\mathbb C^n\times\mathbb C^n$, linear in $u$ and conjugate linear in $v$. Its diagonal value is $u^*Tu$; it can be complex when $T$ is not Hermitian. Over $\mathbb R$, the
corresponding map is bilinear. The following estimate for the diagonal values separates the fourth moments from the off-diagonal pairings and the terms controlled by bounded
support.

\begin{lemma}
\label{W:lem:schur-quadratic}
Under \eqref{W:eq:schur-core-moments}, for every even $p\ge4$,
\begin{equation}
 \left\|X^*TX-\sum_j\mathbb E|X_j|^2T_{jj}\right\|_{L^p}
       \le Cp^2\mathcal K_N(T).
 \label{W:eq:schur-quadratic}
\end{equation}
For independent vectors $X,Y$ with independent centered components satisfying the same bounds,
\begin{equation}
 \|X^*TY\|_{L^p}\le Cp^2\left\{
       N^{-1}\|T\|_{\mathrm{HS}}
       +b_NN^{-1/2}R_2(T)+b_N^2E_\infty(T)\right\}.
 \label{W:eq:schur-bilinear}
\end{equation}
The matrix $T$ need not be Hermitian.
\end{lemma}
\begin{proof}
Put $v_j=\mathbb E|X_j|^2$, $\tau_j=\mathbb EX_j^2$, $a_j=\mathbb E|X_j|^4$, and $Q=X^*TX-\sum_jv_jT_{jj}$. Independence and centering give the exact second moment identity
\[
 \mathbb E|Q|^2
 =\sum_j(a_j-v_j^2)|T_{jj}|^2
  +\sum_{i\ne j}v_iv_j|T_{ij}|^2
  +\sum_{i\ne j}\overline{\tau_i}\tau_jT_{ij}\overline{T_{ji}}.
\]
The first term comes from coincident indices; the other two come from the two pairings of distinct indices. Since $|\tau_j|\le v_j$ and
$2|T_{ij}T_{ji}|\le|T_{ij}|^2+|T_{ji}|^2$,
\[
 \mathbb E|Q|^2\le F_NN^{-2}D_2(T)^2+8N^{-2}\|T\|_{\rm HS}^2.
\]
In particular the estimate permits $\tau_j\ne0$. For independent vectors $X,Y$, the sesquilinear expression instead satisfies
\[
 \mathbb E|X^*TY|^2
 =\sum_{i,j}\mathbb E|X_i|^2\mathbb E|Y_j|^2|T_{ij}|^2.
\]

For the diagonal part of $Q$, apply \eqref{W:eq:schur-rosenthal} to $T_{jj}(|X_j|^2-v_j)$. Its variance sum is at most $F_ND_2(T)^2/N^2$, and
\[
 \sum_j\mathbb E|T_{jj}(|X_j|^2-v_j)|^p
 \le 2^pE_\infty(T)^p b_N^{2p-4}\frac{F_N}{N}
 \le(2b_N^2E_\infty(T))^p,
\]
because $F_N/(Nb_N^4)\le1$. Thus its $L^p$ norm is at most $Cp[\sqrt{F_N}D_2(T)/N+b_N^2E_\infty(T)]$.

For $X^*TY$, condition on $Y$ and apply \eqref{W:eq:schur-rosenthal}:
\[
 \|X^*TY\|_p\le Cp\left[
 N^{-1/2}\|TY\|_{L^p(\ell^2)}
 +\Bigl(\sum_i m_i\|(TY)_i\|_p^p\Bigr)^{1/p}\right],
 \qquad m_i=\mathbb E|X_i|^p.
\]
The Hilbert space and scalar versions of \eqref{W:eq:schur-linear} give
\[
 \begin{aligned}
 \|TY\|_{L^p(\ell^2)}
 &\le Cp\left[N^{-1/2}\|T\|_{\rm HS}
                   +b_N\max_j\|T_{\cdot j}\|_2\right],\\
 \|(TY)_i\|_p
 &\le Cp\left[N^{-1/2}\|T_{i\cdot}\|_2
                   +b_N\max_j|T_{ij}|\right].
 \end{aligned}
\]
Since $\sum_i m_i\le b_N^{p-4}F_N/N\le b_N^p$, Minkowski's inequality for the weights $m_i$ implies
\[
 \Bigl(\sum_i m_i\|(TY)_i\|_p^p\Bigr)^{1/p}
 \le Cp\left[b_NN^{-1/2}\max_i\|T_{i\cdot}\|_2
                       +b_N^2E_\infty(T)\right].
\]
Substitution proves \eqref{W:eq:schur-bilinear}. This calculation also holds for rectangular submatrices, with the corresponding row and column dimensions bounded by $N$.

For the off-diagonal part of $Q$, let $\epsilon_i$ be independent Bernoulli variables of mean $1/2$, independent of $X$, and put $A=\{i:\epsilon_i=1\}$. Then
\[
 \sum_{i\ne j}\overline X_iT_{ij}X_j
 =4\mathbb E_\epsilon X_A^*T_{A,A^c}X_{A^c},\qquad
 \left\|\sum_{i\ne j}\overline X_iT_{ij}X_j\right\|_p
 \le4\max_A\|X_A^*T_{A,A^c}X_{A^c}\|_p.
\]
For each fixed $A$, the two vectors are independent, and each matrix norm in \eqref{W:eq:schur-bilinear} decreases under this restriction. The bound therefore applies uniformly in
$A$. Adding the diagonal estimate proves \eqref{W:eq:schur-quadratic}.

\begingroup \setlength{\intextsep}{1.2\baselineskip}
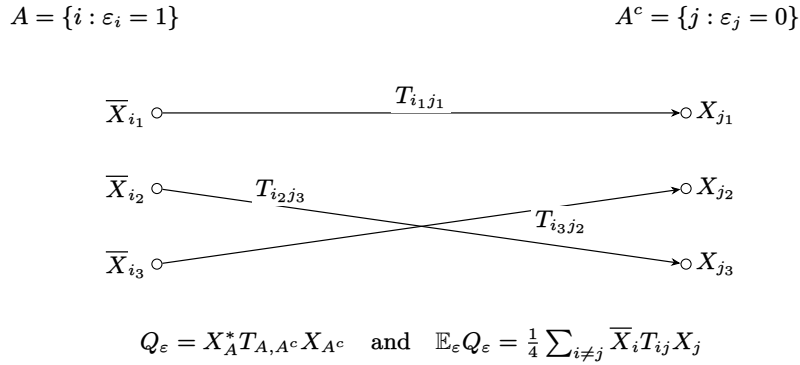
\begin{figure}[H]
\centering
\begin{tikzpicture}[x=1cm,y=1cm,every node/.style={mg text}]
 \node[anchor=west] at (0,1.25) {$A=\{i:\epsilon_i=1\}$};
 \node[anchor=west] at (8,1.25) {$A^c=\{j:\epsilon_j=0\}$};
 \foreach \k/\a in {1/i_1,2/i_2,3/i_3}{
  \node[mg site,label=left:{$\overline X_{\a}$}] (a\k) at (2,{1-\k}) {};
  \node[mg site,label=right:{$X_{j_\k}$}] (b\k) at (9,{1-\k}) {};
 }
 \draw[mg arrow] (a1)--node[mg label,above] {$T_{i_1j_1}$}(b1);
 \draw[mg arrow] (a2)--node[mg label,above,pos=.23] {$T_{i_2j_3}$}(b3);
 \draw[mg arrow] (a3)--node[mg label,below,pos=.77] {$T_{i_3j_2}$}(b2);
 \node[align=center] at (5.5,-3.05)
 {$Q_\epsilon=X_A^*T_{A,A^c}X_{A^c}$\quad and\quad
  $\mathbb E_\epsilon Q_\epsilon=\frac14\sum_{i\ne j}\overline X_iT_{ij}X_j$};
\end{tikzpicture}
\caption{The random partition used in passing from \eqref{W:eq:schur-bilinear} to \eqref{W:eq:schur-quadratic}. A fixed selector keeps the terms $\overline X_iT_{ij}X_j$ with $i\in A$ and $j\in A^c$. The vectors $X_A$ and $X_{A^c}$ are independent. Every ordered pair $i\ne j$ belongs to $A\times A^c$ with probability $1/4$, giving the factor four in the expectation over the selectors. The three arrows illustrate entries of the full matrix $T_{A,A^c}$; their incident circles specify their indices.}
\label{W:fig:schur-bernoulli-decoupling}
\end{figure}
\endgroup
\end{proof}

Fix a labelled component $I$ of size at most $d$, put $J=I^c$, and use the resolvent convention $R_H(z)=(zI-H)^{-1}$. For $\eta=|\operatorname{Im}z|>0$ define
\begin{equation}
 \begin{aligned}
 {\mathcal R}_I(z)&=(zI_J-H_{JJ})^{-1}, & B_I&=H_{IJ},\\
 m_I(z)&=N^{-1}\operatorname{Tr}{\mathcal R}_I(z), &
 Z_I(z)&=B_I{\mathcal R}_I(z)B_I^*-m_I(z)I_I.
 \end{aligned}
 \label{W:eq:schur-residual}
\end{equation}
The normalization $N^{-1}$ agrees with the variance of each matrix entry. Define the adaptive scale
\begin{equation}
 \begin{split}
 \Theta_N(T;z)={}&\mathcal K_N(T)
   +\{u_N(N^{-1/2}+b_N)+u_N^2\}\eta^{-1}
   +\frac{\delta_N}{N}\sum_j|T_{jj}|.
 \end{split}
 \label{W:eq:schur-theta}
\end{equation}
In this notation $T={\mathcal R}_I(z)$ whenever the term $\eta^{-1}$ is used.

\Needspace{6\baselineskip} Applied conditionally on a complementary minor, the preceding moment estimates control the Schur complement error with constants uniform in the
conditioning.

\begin{theorem}
\label{W:thm:schur-residual}
Fix $d\ge7$, a nonreal spectral parameter $z$, and a connected component $I$ of the label graph on $\mathcal A_d$; put $J=I^c$. Let $\mathscr F_I$ be the $\sigma$-field generated
by the label graph, all diagonal entries, the off-diagonal entries within its components, and $H_{JJ}$. Conditional on $\mathscr F_I$, the entries of $B_I=H_{IJ}$ are independent
and have the conditional laws in \eqref{W:eq:schur-small-law}. Under either family \eqref{W:eq:schur-uniform-family} or \eqref{W:eq:schur-defect-family}, for every even $p\ge4$,
\begin{equation}
 \left\|\|Z_I(z)\|_{\mathrm{HS}}\right\|_{L^p(\mathbb P(\cdot\mid\mathscr F_I))}
 \le C_dp^2\Theta_N({\mathcal R}_I(z);z).
 \label{W:eq:schur-conditional-bound}
\end{equation}
Here $\Theta_N({\mathcal R}_I(z);z)$ is $\mathscr F_I$-measurable, and the constant is independent of the conditioned values, $N$, $I$ and $z$. For each fixed $A>0$, there is
$N_0(A,d)$, uniform in the conditioned values and spectral parameter, such that for all
$N\ge N_0(A,d)$,
\begin{equation}
 \mathbb P\left(\|Z_I(z)\|_{\mathrm{HS}}>
      (\log N)^4\Theta_N({\mathcal R}_I(z);z)\mid\mathscr F_I\right)
 \le N^{-A}\quad\text{on }\mathcal A_d.
 \label{W:eq:schur-conditional-prob}
\end{equation}
The same bound holds simultaneously for any fixed polynomial number of nonreal spectral points after increasing the chosen moment order. The diagonal entries belong to
$\mathscr F_I$ and are never averaged in this conditional estimate.
\end{theorem}
\begin{proof}
All expectations and norms in this proof are conditional on the data in the statement. Put $T={\mathcal R}_I(z)$ and write $B_I=M+X$. Lemma~\ref{W:lem:schur-components} and
\eqref{W:eq:schur-core-moments} give independent centered rows $x_a$ of $X$, $\|m_a\|_2\le u_N$, and $\mathbb E|X_{aj}|^2=N^{-1}\sigma_{aj}^2$, with $|\sigma_{aj}^2-1|\le\delta_N$.
For $a,b\in I$, set
\[
 Q_{ab}=x_aTx_b^*-\mathbf1_{\{a=b\}}N^{-1}\sum_j\sigma_{aj}^2T_{jj},
 \qquad
 v_{ab}=\mathbf1_{\{a=b\}}N^{-1}\sum_j(\sigma_{aj}^2-1)T_{jj}.
\]
The $(a,b)$ entry of $Z_I$ is
\[
 (Z_I)_{ab}=Q_{ab}+x_aTm_b^*+m_aTx_b^*+m_aTm_b^*+v_{ab}.
\]
In particular its conditional mean is $m_aTm_b^*+v_{ab}$, not zero. Writing $V_0=(v_{ab})$, the two deterministic contributions satisfy
\[
 \|\mathbb E Z_I\|_{\mathrm{HS}}
 \le\|MTM^*\|_{\mathrm{HS}}+\|V_0\|_{\mathrm{HS}}
 \le d u_N^2\eta^{-1}
       +\sqrt d\,\frac{\delta_N}{N}\sum_j|T_{jj}|.
\]
Here $\|M\|_{\mathrm{HS}}^2\le d u_N^2$ and $V_0$ is diagonal. These two quantities are the corresponding deterministic terms in $\Theta_N$. For $a=b$ apply
\eqref{W:eq:schur-quadratic}, and for $a\ne b$ apply \eqref{W:eq:schur-bilinear}. Their conclusions imply
\[
 \|Q_{ab}\|_p\le Cp^2\mathcal K_N(T),\qquad
 |v_{ab}|\le\frac{\delta_N}{N}\sum_j|T_{jj}|.
\]
The linear estimate \eqref{W:eq:schur-linear}, applied with coefficient vector $Tm_b^*$, yields
\[
 \|x_aTm_b^*\|_p \le Cp\bigl(N^{-1/2}\|Tm_b^*\|_2+b_N\|Tm_b^*\|_\infty\bigr) \le Cp u_N(N^{-1/2}+b_N)\eta^{-1}.
\]
For the other linear term use $m_aTx_b^*=\overline{x_bT^*m_a^*}$ and $\|T^*\|\le\eta^{-1}$; the same bound is
\[
 \|m_aTx_b^*\|_p\le Cp u_N(N^{-1/2}+b_N)\eta^{-1}.
\]
Also $|m_aTm_b^*|\le u_N^2\eta^{-1}$. Consequently
\[
 \|(Z_I)_{ab}\|_p\le Cp^2\left\{
 \mathcal K_N(T)+[u_N(N^{-1/2}+b_N)+u_N^2]\eta^{-1}
       +\frac{\delta_N}{N}\sum_j|T_{jj}|\right\}.
\]
For $p\ge2$, Minkowski's inequality in $L^{p/2}$ gives
\[
 \bigl\|\|Z_I\|_{\mathrm{HS}}\bigr\|_p^2
 =\left\|\sum_{a,b\in I}|(Z_I)_{ab}|^2\right\|_{p/2}
 \le\sum_{a,b\in I}\|(Z_I)_{ab}\|_p^2
 \le C d^2p^4\Theta_N(T;z)^2.
\]
This proves \eqref{W:eq:schur-conditional-bound}.

If $\Theta_N(T;z)=0$, the conditional moment bound gives $Z_I=0$ almost surely. Otherwise divide by $\Theta_N(T;z)$. Choose an even integer $p$ between $(A+2)\log N$ and
$(A+2)\log N+2$. Markov's inequality now reads
\[
 \mathbb P\left(\|Z_I\|_{\mathrm{HS}}>(\log N)^4\Theta_N(T;z)
                  \mid\mathrm{conditioning}\right)
 \le\left(\frac{C_dp^2}{(\log N)^4}\right)^p
 \le\exp[-p\log\log N]\le N^{-A}
\]
for all sufficiently large $N$. For at most $N^b$ spectral points, use $A+b+1$ in this estimate and sum their failure probabilities. All bounds are uniform in the conditioned
entries and complementary minor.
\end{proof}

\Needspace{6\baselineskip} The block resolvent is reconstructed exactly from this Schur complement, whose error and trace terms are unchanged by perturbations inside the selected
principal block.

\begin{lemma}
\label{W:lem:schur-reconstruction}
For every Hermitian $H$ and every nonreal $z$,
\begin{equation}
 (R_H(z)_{II})^{-1}+H_{II}=(z-m_I(z))I_I-Z_I(z),\qquad
 R_H(z)_{II}=
       [(z-m_I(z))I_I-H_{II}-Z_I(z)]^{-1}.
 \label{W:eq:schur-reconstruction}
\end{equation}
Both $m_I$ and $Z_I$ are unchanged by every perturbation supported inside $I\times I$. Their mixed derivatives of every positive order in such directions therefore vanish.
\end{lemma}
\begin{proof}
In the decomposition $\mathbb C^N=\mathbb C^I\oplus\mathbb C^J$, put $D=zI_J-H_{JJ}$ and $S=zI_I-H_{II}-B_ID^{-1}B_I^*$. The block factorization is
\[
 zI-H=
 \begin{pmatrix}I_I&-B_ID^{-1}\\0&I_J\end{pmatrix}
 \begin{pmatrix}S&0\\0&D\end{pmatrix}
 \begin{pmatrix}I_I&0\\-D^{-1}B_I^*&I_J\end{pmatrix}.
\]
Since $z$ is nonreal, both $zI-H$ and $D$ are invertible; hence so is $S$, and inversion of these triangular factors gives
\[
 R_H(z)=
 \begin{pmatrix}
 S^{-1}&S^{-1}B_ID^{-1}\\
 D^{-1}B_I^*S^{-1}&D^{-1}+D^{-1}B_I^*S^{-1}B_ID^{-1}
 \end{pmatrix}.
\]
In particular $R_H(z)_{II}=S^{-1}$. Substituting $B_I{\mathcal R}_IB_I^*=m_I I_I+Z_I$ proves both identities in \eqref{W:eq:schur-reconstruction}. If $V$ is supported on
$I\times I$, then $B_I(H+V)=B_I(H)$ and $(H+V)_{JJ}=H_{JJ}$. Thus
\[
 m_I(H+V,z)=m_I(H,z),\qquad Z_I(H+V,z)=Z_I(H,z).
\]
Differentiation in any collection of directions with this support proves that the corresponding derivatives vanish.
\end{proof}

\subsubsection{A smooth cutoff for the Schur complement}
\label{W:sec:smooth-cutoff-for-schur-complement}

The conditional estimates are indexed by configurations of the labels that select entries for conditioning. Entry comparison requires a smooth function of the matrix, so the cutoff
averages these configurations with their smooth probabilities $\pi_L(H)$. For each fixed configuration the principal blocks remain fixed during differentiation; variation of the
probabilities is accounted for separately.

The spectral parameters belong to a deterministic set $\mathfrak Z_N\subset\mathbb C\setminus\mathbb R$ of polynomial cardinality. For each graph $L\in\mathcal A_d$, consider the
matrices $Z_I(z)$ from \eqref{W:eq:schur-residual}, with $I$ ranging over its components and $z\in\mathfrak Z_N$. The normalizing quantity $\Theta_N(\mathcal R_I(z);z)$ contains
norms and maxima; replacing these by smooth functions gives a positive scale satisfying

\begin{equation}
 \Theta_N({\mathcal R};z)\le\widetilde\Theta_N({\mathcal R};z)
      \le C\{\Theta_N({\mathcal R};z)+N^{-90}\}.
 \label{W:eq:schur-smooth-scale}
\end{equation}

To define this scale, put $\epsilon=N^{-100}$ and $s=\lceil\log(N+1)\rceil$. For a matrix $T$ of size at most $N$, write
\[
 \begin{aligned}
 D_{2,\epsilon}(T)&=
       \left(\sum_j|T_{jj}|^2+\epsilon^2\right)^{1/2},&
 H_\epsilon(T)&=(\|T\|_{\mathrm{HS}}^2+\epsilon^2)^{1/2},\\
 R_{2,\epsilon}(T)&=
   \left(\sum_i(\|T_{i\cdot}\|_2^2+\epsilon^2)^s\right)^{1/(2s)}
  +\left(\sum_j(\|T_{\cdot j}\|_2^2+\epsilon^2)^s\right)^{1/(2s)},\\
 E_{\infty,\epsilon}(T)&=
   \left(\sum_{i,j}(|T_{ij}|^2+\epsilon^2)^s\right)^{1/(2s)},&
 D_{1,\epsilon}(T)&=\sum_j(|T_{jj}|^2+\epsilon^2)^{1/2}.
 \end{aligned}
\]
Replace $D_2,\|T\|_{\mathrm{HS}},R_2,E_\infty$ and $\sum_j|T_{jj}|$ in \eqref{W:eq:schur-theta} by these functions and add $N^{-90}$. This defines $\widetilde\Theta_N$. For
nonnegative numbers $a_1,\ldots,a_m$, $m\le N^2$,
\[
 \max_j a_j\le\left(\sum_j(a_j^2+\epsilon^2)^s\right)^{1/(2s)}
 \le m^{1/(2s)}(\max_j a_j+\epsilon)
 \le\mathrm e(\max_j a_j+\epsilon).
\]
Also $0\le(a^2+\epsilon^2)^{1/2}-a\le\epsilon$. The coefficients multiplying the norm terms in \eqref{W:eq:schur-theta} are bounded for large $N$, since $F_N=N^{o(1)}$, $b_N\le1$
and $\delta_N=N^{-7/12+o(1)}$. The sum of the additive errors is therefore at most $C N\epsilon$. The term containing $\eta^{-1}$ is unchanged. These inequalities give
\eqref{W:eq:schur-smooth-scale}, and every denominator in the definition is strictly positive.

Put $\Lambda_N=(\log(N+2))^4$ and, for each component and selected point, set
\[
 q_{I,z}^L(H)=\frac{\|Z_I(z)\|_{\mathrm{HS}}^2}
                    {\Lambda_N^2\widetilde\Theta_N({\mathcal R}_I(z);z)^2}.
\]
Let $M_L$ be the number of pairs $(I,z)$ and define
\[
 S_L(H)=\frac1{\beta_L}
     \log\left(M_L^{-1}\sum_{I,z}\exp(\beta_L q_{I,z}^L(H))\right),
 \qquad \beta_L=\log(2+M_L).
\]
Choose $\chi_0\in C^\infty(\mathbb R;[0,1])$ equal to one on $(-\infty,1]$ and zero on $[2,\infty)$. Define $U_L(H)=\chi_0(S_L(H))$ for $L\in\mathcal A_d$, and $U_L=0$ otherwise.
The exact bounds
\begin{equation}
 \max q_{I,z}^L-1\le S_L\le\max q_{I,z}^L
 \label{W:eq:schur-softmax}
\end{equation}
show both that $U_L=1$ if all $q_{I,z}^L\le1$, and that every nonzero $U_L$ enforces $q_{I,z}^L<3$ for every pair
$(I,z)$.

For the original matrix, with $\xi_{ij}=\sqrt N H_{ij}$, define the conditional graph probabilities and the corresponding average by
\begin{equation}
 \pi_L(H)=\prod_{i<j}c_N(\xi_{ij})^{L_{ij}} (1-c_N(\xi_{ij}))^{1-L_{ij}}, \qquad \overline U_N(H)=\sum_L\pi_L(H)U_L(H).
 \label{W:eq:schur-posterior-average}
\end{equation}
The factors with exponents zero are interpreted as one, including under the compared ensemble laws. Finally choose a fixed increasing smooth $\vartheta_0$ which is zero on
$(-\infty,1/2]$ and one on $[3/4,\infty)$, and put
\begin{equation}
 U_N^{\mathrm{Sch}}(H)=\vartheta_0(\overline U_N(H)).
 \label{W:eq:schur-final-cutoff}
\end{equation}
The sum is finite and defines a smooth function of the matrix entries.

\Needspace{6\baselineskip}

The theorem below applies to independent upper-triangular entry blocks satisfying \eqref{W:eq:schur-uniform-family}, or \eqref{W:eq:schur-defect-family} with common deterministic
$a_N,F_N,u_N,\delta_N$. The smooth function $U_N^{\mathrm{Sch}}$ is the one defined by \eqref{W:eq:schur-final-cutoff}, for the fixed integer $d$ and polynomially sized set
$\mathfrak Z_N$. All the constants may depend on $d$ and the fixed family bounds, but not on the individual entry distributions.

\begin{theorem}
\label{W:thm:schur-cutoff}
The function $U_N^{\mathrm{Sch}}$ belongs to $C^\infty$ in the real coordinates of Hermitian matrices and satisfies $0\le U_N^{\mathrm{Sch}}\le1$. For every fixed $A>0$, uniformly
over the entry families just specified,
\begin{equation}
 \mathbb E[1-U_N^{\mathrm{Sch}}(H_N)]
  +\mathbb P(U_N^{\mathrm{Sch}}(H_N)\ne1)
 \le C_dN^{1-d/6}+N^{-A}.
 \label{W:eq:schur-cutoff-deficit}
\end{equation}
By choosing the fixed integer $d$ after $A$, either quantity can be made smaller than any prescribed power of $N^{-1}$. The construction does not depend on the deterministic
Lanczos starting vector.
\end{theorem}
\begin{proof}
For each fixed graph $L\in\mathcal A_d$, every matrix ${\mathcal R}_I(z)=(zI-H_{JJ})^{-1}$ is smooth in the real entries of $H$ at a nonreal $z$. The denominators defining
$\widetilde\Theta_N$ are strictly positive, since $\epsilon=N^{-100}$ and the final $N^{-90}$ term is positive. Thus $q_{I,z}^L$, its finite logarithmic sum $S_L$ and
$U_L=\chi_0(S_L)$ are smooth on the full matrix space. For $L\notin\mathcal A_d$, $U_L$ is the constant zero function. Each factor of
\[
 \pi_L(H)=\prod_{i<j}c_N(\sqrt N H_{ij})^{L_{ij}}
                           (1-c_N(\sqrt N H_{ij}))^{1-L_{ij}}
\]
is smooth, because $c_N$ is constant around the origin in the real plane. The identities $\pi_L\ge0$ and $\sum_L\pi_L=\prod_{i<j}[c_N+(1-c_N)]=1$ show that
$\overline U_N=\sum_L\pi_LU_L$ is a smooth function taking values in $[0,1]$. Composition with the fixed scalar $\vartheta_0$ proves smoothness and the bounds between zero and one.

On the common probability space, let $V_e$, indexed by unordered pairs, be independent uniform random variables, independent of the entries $\xi_e$. Define
$L_e=\mathbf1_{\{V_e\le c_N(\xi_e)\}}$. This realizes the conditional graph probabilities $\pi_L(H)$ exactly. Hence
\[
 \overline U_N(H)=\mathbb E[U_L(H)\mid H],\qquad
 \mathbb E(1-\overline U_N(H))=\mathbb E(1-U_L(H)).
\]
The right side splits into the event $\mathcal A_d^c$ and its complement; on the first, $0\le1-U_L\le1$. Lemma~\ref{W:lem:schur-components} gives
\[
 \mathbb E[(1-U_L)\mathbf1_{\mathcal A_d^c}]
 \le C_dN^{1-d/6}.
\]

Suppose $L\in\mathcal A_d$. It has at most $N$ components. If $|\mathfrak Z_N|\le N^b$ for all $N\ge N_0$, the set of pairs $(I,z)$ has at most $N^{b+1}$ elements. For each such
pair, condition on the $\sigma$-field $\mathscr F_I$ of Theorem~\ref{W:thm:schur-residual}. Since $\widetilde\Theta_N\ge\Theta_N$ and $\Lambda_N=(\log(N+2))^4\ge(\log N)^4$, its
conditional tail estimate with exponent $A+b+4$ gives
\[
 \begin{aligned}
 &\mathbb P\bigl(\|Z_I(z)\|_{\rm HS}>
                  \Lambda_N\widetilde\Theta_N({\mathcal R}_I(z);z),
                    \,\mathcal A_d\bigr)\\
 &\quad=\mathbb E\left[\mathbf1_{\mathcal A_d}
 \mathbb P\bigl(\|Z_I(z)\|_{\rm HS}>
            \Lambda_N\widetilde\Theta_N\mid\mathscr F_I\bigr)\right]
 \le N^{-A-b-4}.
 \end{aligned}
\]
Here $\mathcal A_d$ is $\mathscr F_I$-measurable and the conditional scale depends only on the complementary minor. The union bound over the at most $N^{b+1}$ pairs gives
\[
 \mathbb P\bigl(\mathcal A_d\cap
                   \{\max_{I,z}q_{I,z}^L>1\}\bigr)
 \le N^{b+1}N^{-A-b-4}=N^{-A-3}.
\]
When this event does not occur, $q_{I,z}^L\le1$ for every pair. By \eqref{W:eq:schur-softmax}, $S_L\le\max q_{I,z}^L\le1$, so $U_L=\chi_0(S_L)=1$. We have proved
\[
 \mathbb E(1-\overline U_N(H_N))
 \le C_dN^{1-d/6}+N^{-A-3}.
\]

For $t\in[0,1]$ the cutoff $\vartheta_0$ is one on $[3/4,1]$ and takes values in $[0,1]$. Therefore $1-\vartheta_0(t)\le4(1-t)$, and if $\vartheta_0(t)\ne1$ then $1-t>1/4$. Both
the expectation and the probability can consequently be bounded by the same quantity:
\[
 \begin{aligned}
 \mathbb E[1-U_N^{\rm Sch}(H_N)]
 &\le4\mathbb E[1-\overline U_N(H_N)],\\
 \mathbb P(U_N^{\rm Sch}(H_N)\ne1)
 &\le4\mathbb E[1-\overline U_N(H_N)].
 \end{aligned}
\]
Adding and absorbing the fixed factor eight, and applying the estimate with a larger auxiliary value of $A$, yields \eqref{W:eq:schur-cutoff-deficit} in the stated form. Finally,
to obtain a prescribed failure order $N^{-D}$, first choose the tail exponent $A>D+2$, then fix $d$ such that $d/6> D+2$, and only afterwards let $N\to\infty$. No constant in the
conditional estimate depends on the deterministic starting vector.
\end{proof}

The diagonal entries are contained in $H_{II}$ in \eqref{W:eq:schur-reconstruction}. The cutoff controls the norm $\|B_I\mathcal R_I B_I^*-m_I I_I\|_{\mathrm{HS}}$. Moreover
$\overline U_N=\sum_L\pi_LU_L$ is a convex average. If $U_N^{\mathrm{Sch}}(H)>0$, then $\overline U_N(H)>1/2$, and at least one graph with $\pi_L(H)>0$ satisfies $U_L(H)>0$. For
each such graph, the bounds for $Z_I(z)$ apply to every fixed component $I$.

\Needspace{6\baselineskip} At the edge scale, coordinate and averaged resolvent bounds reduce the conditional error in the Schur complement to the required order.

\begin{proposition}
\label{W:prop:schur-airy-scale}
Suppose, for a complementary minor and a spectral point in the upper half-plane,
\begin{equation}
 \eta=N^{-2/3+o(1)},\qquad
 \max_{i,j}|{\mathcal R}_{ij}|\le N^{o(1)},\qquad
 -N^{-1}\operatorname{Im}\operatorname{Tr}{\mathcal R}\le N^{-1/3+o(1)}.
 \label{W:eq:schur-airy-hypotheses}
\end{equation}
Then
\begin{equation}
 \Theta_N({\mathcal R};z)+\widetilde\Theta_N({\mathcal R};z)
            \le N^{-1/3+o(1)}.
 \label{W:eq:schur-airy-scale}
\end{equation}
The same statement holds below the real axis by adjunction.
\end{proposition}
\begin{proof}
The Ward identities for ${\mathcal R}=(z-H_{JJ})^{-1}$ are
\[
 \|{\mathcal R}\|_{\mathrm{HS}}^2=-\eta^{-1}\operatorname{Im}\operatorname{Tr}{\mathcal R},
 \qquad
 \sum_j|{\mathcal R}_{ij}|^2=\sum_j|{\mathcal R}_{ji}|^2
                       =-\eta^{-1}\operatorname{Im}{\mathcal R}_{ii}.
\]
The hypotheses and these identities give
\[
 D_2({\mathcal R})\le N^{1/2+o(1)},\quad
 \|{\mathcal R}\|_{\mathrm{HS}}^2\le N^{4/3+o(1)},\quad
 R_2({\mathcal R})\le N^{1/3+o(1)},\quad
 E_\infty({\mathcal R})\le N^{o(1)}.
\]
Using $F_N=N^{o(1)}$ and $b_N=4N^{-5/24}$, the four contributions are
\begin{align*}
 \frac{\sqrt{F_N}}N D_2({\mathcal R})&\le N^{-1/2+o(1)},&
 \frac1N\|{\mathcal R}\|_{\mathrm{HS}}&\le N^{-1/3+o(1)},\\
 \frac{b_N}{\sqrt N}R_2({\mathcal R})
 &\le N^{-5/24-1/2+1/3+o(1)}=N^{-3/8+o(1)},&
 b_N^2E_\infty({\mathcal R})&\le N^{-5/12+o(1)}.
\end{align*}
The remaining terms of \eqref{W:eq:schur-theta} have the following bounds, using $\eta^{-1}=N^{2/3+o(1)}$:
\[
 \begin{aligned}
 u_NN^{-1/2}\eta^{-1}
   &\le N^{-7/8-1/2+2/3+o(1)}=N^{-17/24+o(1)},\\
 u_Nb_N\eta^{-1}
   &\le N^{-7/8-5/24+2/3+o(1)}=N^{-5/12+o(1)},\\
 u_N^2\eta^{-1}
   &\le N^{-7/4+2/3+o(1)}=N^{-13/12+o(1)},\\
 \frac{\delta_N}{N}\sum_j|\mathcal R_{jj}|
   &\le\delta_N E_\infty(\mathcal R)=N^{-7/12+o(1)}.
 \end{aligned}
\]
Adding these to the four preceding contributions gives $\Theta_N(\mathcal R;z)\le N^{-1/3+o(1)}$. The defining comparison \eqref{W:eq:schur-smooth-scale} then implies
\[
 \Theta_N+\widetilde\Theta_N
 \le(1+C)\Theta_N+CN^{-90}\le N^{-1/3+o(1)}.
\]
For the lower half-plane, $\mathcal R(\bar z)=\mathcal R(z)^*$ preserves all the norms in these estimates and replaces $\eta$ by $|\operatorname{Im}z|$.
\end{proof}

To apply Proposition~\ref{W:prop:schur-airy-scale}, we still need the two bounds in \eqref{W:eq:schur-airy-hypotheses}. The deterministic estimate $\|{\mathcal R}\|\le\eta^{-1}$
would give only $N^{-1}\|{\mathcal R}\|_{\mathrm{HS}}\le N^{1/6+o(1)}$ at the edge, which is too large. Section~\ref{W:sec:entrywise-estimates-and-perturbations} obtains the
required $N^{-1/3+o(1)}$ bound from eigenvalue counting and the self-consistency equations obtained from the block Schur complement;
Section~\ref{W:sec:eigenvalue-counting-and-minor-estimates} provides the counting estimate for the Hermitian case.

\subparagraph{Differentiation with the blocks fixed}
\label{W:sec:derivatives-of-schur-complement-and-cutoff}

Use the conditional distributions of Section~\ref{W:sec:conditional-moments-of-schur-complement} and the cutoff just defined. Fix a graph and its component $I$. For a finite set of
matrix entry derivative labels $S$, let $W_j$ be the corresponding Hermitian matrices, already including the factor $N^{-1/2}$ from the standardized real coordinates of the
selected entry. Set $D_j=(W_j)_{IJ}$ and $V_j=(W_j)_{JJ}$. For $T\subseteq S$ write
\begin{equation}
 {\mathcal R}_\varnothing={\mathcal R},\qquad
 {\mathcal R}_T=\sum_{\tau\in\operatorname{Perm}(T)}
             {\mathcal R}V_{\tau(1)}{\mathcal R}\cdots V_{\tau(|T|)}{\mathcal R}\quad(T\ne\varnothing).
 \label{W:eq:schur-resolvent-derivatives}
\end{equation}
Distinct labels may denote the same matrix $W_j$. The positive sign follows from
$D_H(zI-H)^{-1}[W]=(zI-H)^{-1}W(zI-H)^{-1}$.

\Needspace{6\baselineskip}

For $1\le |S|\le5$, the labelled product rule for $Z=B{\mathcal R}B^*-N^{-1}\operatorname{Tr}{\mathcal R}\,I_I$ gives
\begin{equation}
 \begin{split}
 D_{H,S}Z={}&B{\mathcal R}_SB^*-N^{-1}\operatorname{Tr}{\mathcal R}_S\,I_I\\
 &+\sum_{j\in S}\{D_j{\mathcal R}_{S\setminus\{j\}}B^*
                       +B{\mathcal R}_{S\setminus\{j\}}D_j^*\}\\
 &+\sum_{j\ne k\in S}
          D_j{\mathcal R}_{S\setminus\{j,k\}}D_k^*.
 \end{split}
 \label{W:eq:schur-residual-derivatives}
\end{equation}
This identity holds for fixed self-adjoint $W_j$, over either scalar field.

Conditional on the minor and internal data, write $B=M+X$ as above. For an $|I|\times|J|$ matrix $D$, put $d_i=D^*e_i$ and define
\[
 \mathcal L_N(T,D)=\max_{i\in I}\left\{
 N^{-1/2}(\|Td_i\|_2+\|T^*d_i\|_2)
       +b_N(\|Td_i\|_\infty+\|T^*d_i\|_\infty)\right\}.
\]
Let $\mathcal M_S$ denote the following conditional mean:
\begin{equation}
 \begin{split}
 \mathcal M_S={}&M{\mathcal R}_SM^*
  +\operatorname{diag}_{i\in I}\left(
       N^{-1}\sum_{j\in J}(\sigma_{ij}^2-1)({\mathcal R}_S)_{jj}\right)\\
 &+\sum_{j\in S}\{D_j{\mathcal R}_{S\setminus\{j\}}M^*
                       +M{\mathcal R}_{S\setminus\{j\}}D_j^*\}\\
 &+\sum_{j\ne k\in S}D_j{\mathcal R}_{S\setminus\{j,k\}}D_k^*.
 \end{split}
 \label{W:eq:schur-derivative-mean}
\end{equation}

\begin{theorem}
\label{W:thm:schur-derivatives}
Then for each even $p\ge4$,
\begin{equation}
 \begin{split}
 \bigl\|\|D_{H,S}Z-\mathcal M_S\|_{\mathrm{HS}}\bigr\|_{L^p(\mathrm{conditional})}
 \le C_{d,|S|}p^2\biggl\{
  \mathcal K_N({\mathcal R}_S)+\mathcal L_N({\mathcal R}_S,M)
       +\sum_{j\in S}\mathcal L_N({\mathcal R}_{S\setminus\{j\}},D_j)\biggr\}.
 \end{split}
 \label{W:eq:schur-derivative-bound}
\end{equation}
For directions all supported in $I\times I$, every derivative of positive order in these identities is zero.
\end{theorem}
\begin{proof}
Apply Lemma~\ref{W:lem:resolvent-derivatives} to $H_{JJ}$ and $V_j=(W_j)_{JJ}$. It gives \eqref{W:eq:schur-resolvent-derivatives} with $q=|S|$. In the labelled Leibniz rule for
$B{\mathcal R}B^*$, at most one label can land on $B$ and at most one on $B^*$, since $B$ depends linearly on the matrix entry coordinates. Writing $D_j=D_{H,j}B$, the labelled
expansion is
\[
 \begin{split}
 D_S(B\mathcal RB^*)={}&B\mathcal R_SB^*
 +\sum_{j\in S}\bigl(D_j\mathcal R_{S\setminus\{j\}}B^*
                    +B\mathcal R_{S\setminus\{j\}}D_j^*\bigr)\\
 &+\sum_{\substack{j,k\in S\\j\ne k}}
                  D_j\mathcal R_{S\setminus\{j,k\}}D_k^*.
 \end{split}
\]
The last sum is ordered: $j$ differentiates $B$, and $k$ differentiates $B^*$. This gives $1+2q+q(q-1)$ terms before expanding the resolvent derivatives. The subtracted term is
$N^{-1}\operatorname{Tr}\mathcal R_S I_I$, since trace is linear. These are the three lines and subtraction in \eqref{W:eq:schur-residual-derivatives}.

After $B=M+X$ is substituted, the quadratic part has expectation
\[
 M{\mathcal R}_SM^*+\operatorname{diag}_{i\in I}
       \left(N^{-1}\sum_{j\in J}\sigma_{ij}^2({\mathcal R}_S)_{jj}\right).
\]
The remaining linear terms have expectation obtained by replacing $B$ by $M$, and the last line is deterministic in the conditioning. This proves
\eqref{W:eq:schur-derivative-mean}. To estimate the centered quantity, put $T={\mathcal R}_S$ and $T_j={\mathcal R}_{S\setminus\{j\}}$. For $a,b\in I$, its $(a,b)$ entry is the sum
of
\begin{align*}
 Q_{ab}&=(XTX^*)_{ab}-\mathbb E[(XTX^*)_{ab}\mid\text{minor and internal data}],\\
 L_{ab}&=(XTM^*)_{ab}+(MTX^*)_{ab},\\
 L_{ab}^{(j)}&=(D_jT_jX^*)_{ab}+(XT_jD_j^*)_{ab},\qquad j\in S.
\end{align*}
For $a=b$, \eqref{W:eq:schur-quadratic} gives $\|Q_{aa}\|_p\le Cp^2\mathcal K_N(T)$. For $a\ne b$, conditional independence and centering give
$\mathbb E[x_aTx_b^*\mid\text{minor and internal data}]=0$. Thus $Q_{ab}=x_aTx_b^*$, and \eqref{W:eq:schur-bilinear} gives $\|Q_{ab}\|_p\le Cp^2\mathcal K_N(T)$. If $d_b=D^*e_b$,
the two orientations of a linear term obey
\begin{align*}
 \|(XTD^*)_{ab}\|_p
 &\le Cp\{N^{-1/2}\|Td_b\|_2+b_N\|Td_b\|_\infty\},\\
 \|(DTX^*)_{ab}\|_p
 &\le Cp\{N^{-1/2}\|T^*d_a\|_2+b_N\|T^*d_a\|_\infty\}.
\end{align*}
Use these inequalities first with $D=M,T={\mathcal R}_S$, then with $D=D_j,T=T_j$. There is one quadratic term and $2+2|S|$ linear terms in each entry; the terms with two
differentiated row factors are already in the conditional mean. Finally, for $p\ge2$,
\[
 \big\|\|A\|_{\rm HS}\big\|_p
 \le\Big(\sum_{a,b\in I}\|A_{ab}\|_p^2\Big)^{1/2}
 \le d\max_{a,b}\|A_{ab}\|_p.
\]
Summing the displayed bounds proves \eqref{W:eq:schur-derivative-bound}. Invariance follows either from the formula, since $D_j=V_j=0$, or from
Lemma~\ref{W:lem:schur-reconstruction}.
\end{proof}

The derivative formula \eqref{W:eq:schur-residual-derivatives} distinguishes variations of the external rows from variations of the complementary minor. For instance, for labels
$j\ne k$ assigned to the two external rows, $D_j\mathcal R=D_k\mathcal R=0$, while the last line of \eqref{W:eq:schur-derivative-mean} contains the deterministic conditional
contribution
\[
 D_j\mathcal RD_k^*+D_k\mathcal RD_j^*.
\]
If instead all derivatives act within $J\times J$, then $D_j=D_k=0$ as row derivatives, but $\mathcal R_S$ remains the nonzero sum of resolvent products given in
\eqref{W:eq:schur-resolvent-derivatives}.

The resolvent divided differences used in the contour integrals have a particularly simple form. With $\mathcal R(z)=(zI_J-H_{JJ})^{-1}$, the identity
\[
 \mathcal R(z)-\mathcal R(w)=(w-z)\mathcal R(z)\mathcal R(w)
\]
gives $\mathcal R[z,w]=-\mathcal R(z)\mathcal R(w)$. Since the resolvents commute as functions of the same Hermitian minor, induction on the number of interpolation points yields
\[
 \mathcal R[z_1,\ldots,z_r]
 =(-1)^{r-1}\mathcal R(z_1)\cdots\mathcal R(z_r).
\]
At coinciding parameters this is interpreted by continuity, consistently with $\partial_z^{r-1}\mathcal R(z)/(r-1)!=(-1)^{r-1}\mathcal R(z)^r$. For a fixed complementary minor
$B=H_{JJ}$, let $\mathcal C$ be a compact rectifiable contour with $\delta=\operatorname{dist}(\mathcal C,\sigma(B))>0$. If $K=K^*$ and $\|K-B\|<\delta/2$, then
$\operatorname{dist}(\mathcal C,\sigma(K))\ge\delta/2$. The spectral theorem and Lemma~\ref{W:lem:resolvent-derivatives} give, for $z\in\mathcal C$,
\[
 \|(zI-K)^{-1}\|\le 2/\delta,\qquad
 \|D_K^q(zI-K)^{-1}[W_1,\ldots,W_q]\|
 \le q!(2/\delta)^{q+1}\prod_{j=1}^q\|W_j\|.
\]
For $\int_{\mathcal C}|f(z)|\,|\dd z|<\infty$, multiplying this bound by $|f(z)|$ gives an integrable majorant on the same neighbourhood of $B$. Differentiation therefore passes
through the contour integral, including its real points. The matrix $T=\int_{\mathcal C}f(z)\mathcal R(z)\,\dd z$ is measurable with respect to the complementary minor. Applying
the conditional row estimate to this matrix gives
\[
 \|X^*TX-\mathbb E[X^*TX\mid\mathscr F_I]\|_{L^p}
 \le Cp^2\mathcal K_N(T),\qquad
 T=\int_{\mathcal C}f(z)\mathcal R(z)\,\dd z.
\]
Here $\mathcal K_N$ is applied to the entire integrated matrix $T$; the estimate requires bounds for the Hilbert--Schmidt, row and column norms of this integrated matrix.

We next differentiate the probability weights of the graph. A derivative in an entry block is taken with respect to its standardized real coordinates, so that the change
$H_{ab}=N^{-1/2}\xi_{ab}$ is built into the derivative. For a nonempty set $S$ of derivative labels all assigned to the same off-diagonal block $\alpha=\{a,b\}$, write
$c_\alpha=c_N(\xi_{ab})$. Conditioning on all other graph labels gives
\[
 \pi_L=c_\alpha^{L_\alpha}(1-c_\alpha)^{1-L_\alpha}
                      \prod_{e\ne\alpha}c_e^{L_e}(1-c_e)^{1-L_e}.
\]
The derivative of the first factor is either $D_Sc_\alpha$ or $-D_Sc_\alpha$. Summation of its absolute value over both possible labels and then over all other labels gives the
exact identity
\begin{equation}
 \sum_L|D_{H,S}\pi_L(H)|=2|D_{H,S}c_N(\sqrt N H_{ab})|
      \le C_{|S|}(\ell_N^{\mathrm{Sch}})^{-|S|}.
 \label{W:eq:schur-posterior-derivatives}
\end{equation}
Each derivative of $c_N$ is supported in the annulus $\ell_N^{\mathrm{Sch}}<|\xi_{ab}|<2\ell_N^{\mathrm{Sch}}$, and differentiation of its fixed rescaled radial cutoff gives a
factor $(\ell_N^{\mathrm{Sch}})^{-1}$ per standardized variable. For diagonal variations the derivative of every graph probability is zero.

If the labels in $S$ are distributed among different unordered blocks, denote the set of those blocks by $\mathcal J$ and the nonempty group assigned to $\alpha$ by $S_\alpha$.
Independence of the factors in the conditional graph probabilities yields
\[
 \begin{aligned}
 \sum_L|D_{H,S}\pi_L|
 &=\prod_{\alpha\in\mathcal J}
      (|D_{H,S_\alpha}c_\alpha|
        +|D_{H,S_\alpha}(1-c_\alpha)|)\\
 &=\prod_{\alpha\in\mathcal J}2|D_{H,S_\alpha}c_\alpha|
 \le C_{|S|}(\ell_N^{\mathrm{Sch}})^{-|S|},
 \end{aligned}
\]
because the undifferentiated graph factors sum to one. For a complex entry, $S_\alpha$ includes the derivatives in both of its real coordinates, all acting on the same factor
$c_\alpha$.

Now let $s=\Lambda_N\widetilde\Theta_N({\mathcal R}_I(z);z)>0$ and $q=\|Z_I(z)\|_{\rm HS}^2s^{-2}$. For every finite label set $S$, the labelled Leibniz rule, including derivatives
of the denominator, is
\begin{equation}
 D_{H,S}q
 =\sum_{S_0\sqcup S_1\sqcup S_2=S}
       D_{H,S_0}(s^{-2})
       \operatorname{Tr}[(D_{H,S_1}Z_I)^*D_{H,S_2}Z_I],
 \label{W:eq:schur-monitor-derivatives}
\end{equation}
where $D_\varnothing Z_I=Z_I$ and, for nonempty $T$,
\[
 D_{H,T}(s^{-2})
 =\sum_{\pi\in\Pi(T)}(-1)^{|\pi|}(|\pi|+1)!
       s^{-|\pi|-2}\prod_{B\in\pi}D_{H,B}s.
\]
Indeed $(\dd^k/\dd s^k)s^{-2}=(-1)^k(k+1)!s^{-k-2}$ and the partition formula for a scalar composition applies. The sets $S_0\ne\varnothing$ correspond to differentiating the
dependence of $s$ on the complementary minor.

For the logarithmic average $S_L=\beta_L^{-1}\log(M_L^{-1}\sum_a e^{\beta_Lq_a})$, put $p_a=e^{\beta_Lq_a}/\sum_be^{\beta_Lq_b}$. The identity $\sum_ap_a=1$ gives
\[
 D_{H,j}S_L=\sum_ap_aD_{H,j}q_a,
 \qquad
 D_{H,j}p_a=\beta_Lp_a\left(D_{H,j}q_a
                   -\sum_bp_bD_{H,j}q_b\right).
\]
Iterating these two formulas proves, for every nonempty set $S$,
\[
 D_{H,S}S_L
 =\sum_{\pi\in\Pi(S)}\beta_L^{|\pi|-1}
       \operatorname{cum}_{p}
          (D_{H,B}q_a:B\in\pi),
\]
where the cumulant is taken with respect to the finite probability distribution $(p_a)$. For $\pi=\{B_1,\ldots,B_k\}$ the expression is explicitly
\[
 \operatorname{cum}_p(V_1,\ldots,V_k)
 =\sum_{\tau\in\Pi(\{1,\ldots,k\})}
          (-1)^{|\tau|-1}(|\tau|-1)!
        \prod_{C\in\tau}\left(\sum_ap_a\prod_{j\in C}V_j(a)\right),
 \quad V_j(a)=D_{H,B_j}q_a.
\]
The normalization $1/M_L$ in $S_L$ is constant and contributes no derivative.

Finally $U_L=\chi_0(S_L)$, $\overline U_N=\sum_L\pi_LU_L$, and $U_N^{\rm Sch}=\vartheta_0(\overline U_N)$. Their derivatives are
\begin{equation}
 \begin{aligned}
 D_{H,S}\overline U_N
 &=\sum_{T\subseteq S}\sum_L(D_{H,T}\pi_L)(D_{H,S\setminus T}U_L),\\
 D_{H,S}U_N^{\rm Sch}
 &=\sum_{\pi\in\Pi(S)}\vartheta_0^{(|\pi|)}(\overline U_N)
                     \prod_{B\in\pi}D_{H,B}\overline U_N.
 \end{aligned}
 \label{W:eq:derivatives-of-the-schur-cutoff}
\end{equation}
For $U_L=\chi_0(S_L)$ the same composition formula applies. Proposition~\ref{W:prop:schur-polynomial-derivatives} estimates the products in \eqref{W:eq:schur-monitor-derivatives}
and \eqref{W:eq:derivatives-of-the-schur-cutoff}.

\Needspace{6\baselineskip} The complete cutoff derivatives have polynomial growth, and each nonzero derivative is supported where one of the defining smooth cutoffs changes value.

\begin{proposition}
\label{W:prop:schur-polynomial-derivatives}
Suppose $\min_{z\in\mathfrak Z_N}|\operatorname{Im}z|\ge N^{-M}$ and $\max_{i,j}|H_{ij}|\le N^B$, for fixed $M,B$. For $1\le q\le5$ there is a finite exponent $C_q=C_q(M,B)$, also
depending on the fixed exponent specifying the discretization and smoothing convention, such that every standardized matrix entry derivative obeys
\begin{equation}
 |\partial_\alpha^{\boldsymbol\sigma}U_N^{\mathrm{Sch}}(H)|
       \le C_{q,d}N^{C_q}
                    \mathbf1_{\{U_N^{\mathrm{Sch}}(H)\ne1\}}.
 \label{W:eq:schur-polynomial-derivatives}
\end{equation}
The power $C_q$ can be chosen independently of the fixed bound on the component size $d$; only its multiplicative constant depends on $d$.
\end{proposition}
\begin{proof}
If the set of spectral parameters has at most $C N^{J}$ points, fix $J_0=J+1$. Then its cardinality is at most $N^{J_0}$ for $N\ge C$. Set
\[
 A=1000+20\bigl(\max(M,0)+\max(B,0)+J_0+1\bigr).
\]
All bounds below are for $0\le q\le5$ and sufficiently large $N$; constants may depend on $q,d$, but the displayed exponents do not. Since a normalized entry direction has norm at
most $2N^{-1/2}$, \eqref{W:eq:schur-resolvent-derivatives} and \eqref{W:eq:schur-residual-derivatives} imply
\begin{align*}
 \|D_{H,S}{\mathcal R}\|&\le q!2^qN^{(q+1)\max(M,0)},\\
 \|D_{H,S}Z\|_{\rm HS}
 &\le C_qN^{2\max(B,0)+3+(q+1)\max(M,0)}
 \le C_qN^{(q+1)A}.
\end{align*}
Here $\|B_I\|\le N^{B+1}$. For a matrix $T$ indexed by $J$,
\[
 \left|N^{-1}\operatorname{Tr}T\right|
 \le N^{-1}\sum_{j\in J}|\langle e_j,Te_j\rangle|
 \le (|J|/N)\|T\|\le\|T\|.
\]
For a matrix on $I$, $\|T\|_{\rm HS}\le|I|^{1/2}\|T\|$ and $|I|^{1/2}\le N^{1/2}$. These estimates include $q=0$.

We give the bounds for the regularized norms explicitly. Every squared quantity $t_j$ used in their construction is a sum of at most $N^2$ squared resolvent entries, or one such
squared entry, with the positive floor $N^{-200}$. The product rule and the preceding resolvent estimate give
\[
 t_j\ge N^{-200},\qquad
 |D_{H,S}t_j|\le C_qN^{(q+1)A}.
\]
For a nonempty $S$, the exact logarithmic derivative is
\[
 D_{H,S}\log t_j
 =\sum_{\pi\in\Pi(S)}(-1)^{|\pi|-1}(|\pi|-1)!
       t_j^{-|\pi|}\prod_{B\in\pi}D_{H,B}t_j.
\]
A partition into $k$ blocks contributes at most $C_qN^{200k+(q+k)A}\le C_qN^{3qA}$. If $F=(\sum_jt_j^\ell)^{1/(2\ell)}$, then
\[
 \log F=\frac1{2\ell}\log\sum_j e^{\ell\log t_j},\qquad
 |D_{H,S}\log F|
 \le C_q\sum_{\pi\in\Pi(S)}\ell^{|\pi|-1}
       \prod_{B\in\pi}\max_j|D_{H,B}\log t_j|
 \le C_qN^{4qA}.
\]
For a probability vector $\pi$ and bounded scalar functions $a_1,\ldots,a_k$ on its finite index set, the partition formula gives
\[
 |\operatorname{cum}_\pi(a_1,\ldots,a_k)|
 \le\left(\sum_{\sigma\in\Pi_k}(|\sigma|-1)!\right)
                                  \prod_{j=1}^k\|a_j\|_\infty.
\]
Indeed, each block expectation is bounded by the product of the corresponding suprema because $\sum_i\pi_i=1$. For the logarithm above take $\pi_i=t_i^\ell/\sum_jt_j^\ell$ and
$a_B=D_{H,B}\log t_i$. The chain rule supplies a factor $\ell^{k-1}$ for a partition with $k$ blocks, which proves the middle inequality independently of the number of $t_i$. Also
$F\le N^A$. The exponential chain rule therefore gives $|D_{H,S}F|\le C_qN^{(4q+1)A}$, including $q=0$ after increasing the constant. For a regularized Euclidean norm,
$F=(\sum_jt_j)^{1/2}$, this is the case $\ell=1$ of the displayed formula. The deterministic coefficients in $\widetilde\Theta_N$ and $\Lambda_N$ are bounded by $N^A$; hence, for
$s=\Lambda_N\widetilde\Theta_N$,
\[
 s\ge N^{-90},\qquad |D_{H,S}s|\le C_qN^{(4q+2)A}.
\]
In the displayed formula for $D_Ts^{-2}$, a partition into $k$ blocks has bound $C_qN^{(k+2)A+(4q+2k)A}$. Thus
\[
 |D_{H,S}s^{-2}|\le C_qN^{(7q+2)A},\qquad
 |D_{H,S}q_{I,z}|\le C_qN^{(7q+4)A}.
\]
The second estimate follows from \eqref{W:eq:schur-monitor-derivatives}: if its three label sets have sizes $a,b,c$, the exponent is at most $(7a+2+b+1+c+1)A\le(7q+4)A$.

For $q\ge1$, a $k$-block term in the cumulant formula for $D_{H,S}S_L$ is bounded by
\[
 C_q\beta^{k-1}N^{(7q+4k)A}\le C_qN^{12qA},
 \qquad k\le q,
\]
since $\beta=O(\log N)$. There are at most $N^{J_0+1}$ pairs $(I,z)$, and the cumulant uses probability weights, so this count introduces no extra factor. Applying the labelled
chain rule to the fixed function $\chi_0$ gives $|D_{H,S}U_L|\le C_{q,d}N^{12qA}$. For the graph average, \eqref{W:eq:schur-posterior-derivatives} gives
\[
 \sum_L|D_{H,T}\pi_L|\le C_qN^{|T|A},\qquad
 |D_{H,S}\overline U_N|
 \le\sum_{T\subseteq S}C_{q,d}N^{|T|A+12(q-|T|)A}
 \le C_{q,d}N^{12qA}.
\]
The order zero bounds are $\sum_L\pi_L=1$ and $|U_L|\le1$. A final application of \eqref{W:eq:derivatives-of-the-schur-cutoff} yields the same power for $D_{H,S}U_N^{\rm Sch}$. In
particular the explicit choice $C_q=20(q+1)A$ suffices. All sums used above range over at most $N^2$ matrix coordinates or $N^{J_0+1}$ pairs $(I,z)$, or use normalized graph
weights; none has exponent depending on $d$. Every derivative of positive order of $\vartheta_0$ vanishes where $\overline U_N\ge3/4$, including the boundary by smoothness. Its
labelled chain rule therefore supplies the indicator in \eqref{W:eq:schur-polynomial-derivatives}.
\end{proof}

The admissible matrices may have at most $r_0$ exceptional off-diagonal entries, where $r_0$ is fixed. Ordinary entries satisfy \eqref{W:eq:schur-uniform-family} or
\eqref{W:eq:schur-defect-family}. Exceptional entries obey the same mean, tail and truncated fourth-moment bounds, while their centered variances need only be bounded by a fixed
constant. Lemma~\ref{W:lem:schur-small-moments} gives the same conditional mean and support bounds and a bounded conditional variance for these entries. Their total contribution to
the variance bias in the Schur complement is at most
\[
 C r_0 N^{-1}E_\infty(\mathcal R)
 \le C r_0 b_N^2 E_\infty(\mathcal R)
 \le C r_0\mathcal K_N(\mathcal R),
\]
because $b_N^2\ge N^{-1}$. The ordinary entries contribute $\delta_NN^{-1}\sum_j|\mathcal R_{jj}|$ as before. The conditional linear and quadratic bounds, and hence
Theorem~\ref{W:thm:schur-cutoff}, follow with constants depending additionally on $r_0$ and the variance bound.

This extension covers a deleted entry, and an independent replacement entry $\xi$ from either family multiplied by any deterministic $s\in[0,1]$. Indeed
$|\mathbb E(s\xi)|\le\varepsilon_{1,N}$, its second moment is uniformly bounded, and for $0<s\le1$,
\begin{equation}
 \begin{aligned}
 \sup_{t\ge\ell_N^{\mathrm{Sch}}}t^4\mathbb P(|s\xi|>t)
     &\le s^4a_N,\\
 \mathbb E[|s\xi|^4\mathbf1_{\{|s\xi|\le2\ell_N^{\mathrm{Sch}}\}}]
     &\le\mathbb E[|\xi|^4\mathbf1_{\{|\xi|\le2\ell_N^{\mathrm{Sch}}\}}]
           +C a_Ns^4(1+\log(1/s))\\
     &\le C F_N.
 \end{aligned}
 \label{W:eq:schur-scaled-family-moments}
\end{equation}
For the second inequality, split at $|\xi|=2\ell_N^{\mathrm{Sch}}$. The part below this threshold only decreases under scaling; for the part between $2\ell_N^{\mathrm{Sch}}$ and
$2\ell_N^{\mathrm{Sch}}/s$, tail integration gives $C a_Ns^4(1+\int_1^{1/s}\dd t/t)$. The function $s^4(1+\log(1/s))$ is bounded on $(0,1]$, and at $s=0$ the variable vanishes.
These inequalities hold uniformly for $0\le s\le1$, with the same scalar $s$ multiplying both coordinates of a complex entry.

\Needspace{6\baselineskip}

Consider polynomially many intermediate arrays with independent upper-triangular entries, whose ordinary off-diagonal entries satisfy \eqref{W:eq:schur-uniform-family} or
\eqref{W:eq:schur-defect-family}, together with the matrices obtained by varying one entry along its radial segment. Fix common deterministic tail and fourth moment bounds, and
choose $u_N,\delta_N$ by \eqref{W:eq:schur-mean-and-covariance-error-bounds}. Allow at most $r_0$ exceptional entries as specified above. Suppose all entries are bounded by a fixed
power of $N$, with any changes in the moments caused by truncation included in these same bounds. Let $\Phi_N$ be a test function with all standardized entry derivatives through
order five bounded in modulus by $N^{B_\Phi}$ on that set.

\begin{corollary}
\label{W:cor:schur-cutoff-suppression}
For every $A>0$, the bound on the component size can be fixed sufficiently large that the sum of terms in the Taylor polynomials and integral remainders at the two truncation
thresholds that contain a derivative of positive order of $U_N^{\mathrm{Sch}}$ is at most $N^{-A}$ in absolute expectation. The bound includes the mixed real coordinates of a
complex block.
\end{corollary}
\begin{proof}
All entry derivatives in this proof are taken in standardized coordinates, including the two real coordinates of one complex off-diagonal block. Fix the finite derivative order
$r\le5$ appearing in a Taylor polynomial or its integral remainder. The derivative of the product $U_N^{\rm Sch}(H)\Phi_N(H)$ is indexed by an allocation of its $r$ derivative
labels between the cutoff and $\Phi_N$, and all terms with a nonempty set of labels on the cutoff are considered. By the labelled product rule, the number of such allocations is at
most $2^r-1\le31$; further allocations among smooth factors are already included in Proposition~\ref{W:prop:schur-polynomial-derivatives}.

Let there be at most $CN^h$ choices of compared array, replaced entry, derivative multiindex and other discrete parameters. Choose $B_\Phi$ so that every required derivative of
$\Phi_N$ is bounded by $N^{B_\Phi}$ on the stated polynomial entry bound. The comparison differentiates through order five, so each individual entry factor in its coefficients of
the expansion and integral remainders is a monomial of degree at most five in the standardized coordinates of the replaced block. If $|H_{ij}|\le N^B$, then
$|\xi_{ij}|=\sqrt N|H_{ij}|\le N^{B+1/2}$; a real coordinate of a complex block differs by at most the fixed factor $\sqrt2$. Thus the absolute value of any such monomial is at
most $C N^{5(B+1/2)}$ on the stipulated support. All fixed normalization factors in these finitely many replacement formulas are powers of $N$, so choose one exponent $b<\infty$ to
bound each complete scalar factor by $N^b$, uniformly over the two truncation radii, the compared block laws, and the interpolation parameter. The integration kernels have finite
absolute mass.

Let $C_*=5\max_{1\le j\le5}C_j$, where the constants $C_j$ are the fixed exponents in \eqref{W:eq:schur-polynomial-derivatives}. These exponents are independent of $d$. It follows
pointwise for every such product term, at every interpolation parameter $s$, that
\[
 |T_{\alpha,s}(H_{\alpha,s})|
 \le C_dN^{B_\Phi+C_*+b}
       \mathbf1_{\{U_N^{\rm Sch}(H_{\alpha,s})\ne1\}}.
\]
The indicator on the right bounds the terms where a derivative falls on $U_N^{\rm Sch}$, including its transition region. If the cutoff is identically one in a neighbourhood of the
evaluated matrix, all its positive derivatives vanish, precisely as in \eqref{W:eq:schur-polynomial-derivatives}.

The intermediate arrays have independent unordered entry blocks. The moment errors of the uncentered truncated entries are included in \eqref{W:eq:schur-defect-family}; a single
block multiplied by a number in $[0,1]$ or replaced independently is covered by the extension immediately preceding the corollary. The high-probability estimate of
Theorem~\ref{W:thm:schur-cutoff} therefore applies to each $H_{\alpha,s}$ with constants uniform in $\alpha$ and $s$. For a prescribed exponent $A>0$, fix
\[
 D>A+h+B_\Phi+C_*+b+3,
 \qquad d>6(D+2).
\]
These choices are made in this order; the cutoff derivative exponents are independent of $d$, so the second choice does not invalidate the first. Apply the cutoff theorem with
auxiliary exponent $D+2$. For all sufficiently large $N$,
\[
 \sup_{\alpha,s}\mathbb P(U_N^{\rm Sch}(H_{\alpha,s})\ne1)
 \le C_dN^{1-d/6}+N^{-D-2}
 \le C_d'N^{-D-1}.
\]
Consequently every fixed contribution satisfies
\[
 \mathbb E|T_{\alpha,s}|
 \le C_d''N^{B_\Phi+C_*+b-D-1}.
\]
For an integral remainder let $\nu$ denote its deterministic nonnegative integration measure on $[0,1]$ or $[0,1]^2$, after taking the absolute value of the scalar kernel; the
Taylor formulas used here give $\nu([0,1]^k)\le1$. Tonelli's theorem and the uniform bound at each interpolation parameter yield
\[
 \begin{aligned}
 \sum_{\alpha}\int\mathbb E|T_{\alpha,s}|\,\nu(\dd s)
 &\le C_d'''N^hN^{B_\Phi+C_*+b-D-1}\nu([0,1]^k)\\
 &\le N^{-A}
 \end{aligned}
\]
for large $N$. If the remainder has a signed deterministic kernel, its absolute value is taken before this estimate; the exceptional contribution is bounded in absolute value. This
proves the bound uniformly over the compared arrays.
\end{proof}

For the truncated entries, the parameters in Corollary~\ref{W:cor:schur-cutoff-suppression} follow from the tail estimates. Let $\xi$ with $\mathbb E\xi=0$, $\mathbb E|\xi|^2=1$,
and let $\widehat\xi=\xi\mathbf1_{\{|\xi|\le L_N\}}$. Subtraction of the discarded part gives the exact moment identities
\begin{equation}
 \begin{aligned}
 \mathbb E\widehat\xi
   &=-\mathbb E[\xi\mathbf1_{\{|\xi|>L_N\}}],&
 |\mathbb E\widehat\xi|
   &\le\tfrac43\Theta(L_N)L_N^{-3},\\
 1-\mathbb E|\widehat\xi|^2
   &=\mathbb E[|\xi|^2\mathbf1_{\{|\xi|>L_N\}}],&
 0\le1-\mathbb E|\widehat\xi|^2
   &\le2\Theta(L_N)L_N^{-2}.
 \end{aligned}
 \label{W:eq:schur-truncated-defects}
\end{equation}
The inequalities are the $q=1,2$ cases of the tail integration bound \eqref{W:eq:quant-tail-below-four}. As $L_N=N^{1/2-o(1)}$ and $\Theta(L_N)\le C$, their right sides are
respectively $N^{-3/2+o(1)}$ and $N^{-1+o(1)}$, exactly the errors allowed in \eqref{W:eq:schur-defect-family}. In the complex case the same calculation applied to $\xi^2$ and
$\mathbb E\xi^2=0$ yields $|\mathbb E\widehat\xi^2|\le2\Theta(L_N)L_N^{-2}$. The centered variance is
\[
 \operatorname{Var}(\widehat\xi)
 =\mathbb E|\widehat\xi|^2-|\mathbb E\widehat\xi|^2;
\]
the squared mean is of order $N^{-3+o(1)}$. Thus $\delta_N$ bounds both the covariance error and the error in the noncentral second moment used in the entry comparison.

For sufficiently large $N$, \eqref{W:eq:schur-scale-hierarchy} gives $2\ell_N^{\rm Sch}<L_N$. Therefore, for $t\ge\ell_N^{\rm Sch}$,
$\mathbb P(|\widehat\xi|>t)\le\mathbb P(|\xi|>t)$, whereas on $\{|\xi|\le2\ell_N^{\rm Sch}\}$ the truncation does nothing. Consequently
\[
 \sup_{t\ge\ell_N^{\rm Sch}}t^4\mathbb P(|\widehat\xi|>t)
 \le\Theta(\ell_N^{\rm Sch}),\qquad
 \mathbb E[|\widehat\xi|^4\mathbf1_{\{|\widehat\xi|\le2\ell_N^{\rm Sch}\}}]
 =h_4(2\ell_N^{\rm Sch}).
\]
These are the bounds $a_N$ and $F_N$ for the labels formed from the uncentered truncated entries.

Next let $g$ be a standardized real or circular complex Gaussian variable and truncate it radially at $K_N=N^{1/8}$. Sign symmetry in the real case and rotation symmetry in the
complex case give exactly $\mathbb E[g\mathbf1_{\{|g|\le K_N\}}]=0$; in the complex case the pseudovariance remains zero. The variance loss is
\[
 \mathbb E[|g|^2\mathbf1_{\{|g|>K_N\}}]
 \le C(1+K_N^4)e^{-cK_N^2}\le N^{-100}
\]
for large $N$. The first inequality follows by integrating the real or radial Gaussian density and integrating by parts once; its constants depend only on the fixed symmetry class.
Since all Gaussian fourth moments are finite, the truncated fourth radial moment is uniformly bounded. Add $N^{-100}$ to the common bound on the second moment error before defining
$\delta_N$. An untruncated standardized Gaussian or a uniformly bounded centered entry normalized to the prescribed covariance satisfies the corresponding exact centering and
covariance requirements. Thus the single deterministic choices of $a_N,F_N,u_N,\delta_N$ cover every entry law used in the comparison, including independent blockwise combinations.

Finally multiply one such entry by a deterministic $s\in[0,1]$ or set it to zero. The joint complex block is multiplied by one scalar $s$, and the mean changes by at most the
original mean bound. The variance is at most its original fixed upper bound. The tail and truncated fourth moment estimates remain valid by
\eqref{W:eq:schur-scaled-family-moments}, and a single block with a different variance is one of the exceptional blocks allowed before
Corollary~\ref{W:cor:schur-cutoff-suppression}. These facts verify the hypothesis uniformly for the entire interpolation segment; the probability bound can therefore be applied at
each $s$ before integration. The passage back to the untruncated arrays then adds their respective probabilities of disagreement.

The function $U_N^{\mathrm{Sch}}$ and its derivative bounds are used in Sections~\ref{W:sec:minor-estimates-by-coarse-truncation}, \ref{W:sec:resolvent-domain-for-theta}
and~\ref{W:sec:moment-increments-under-entry-replacement}. The derivatives of the functions of Gram matrices whose dimensions grow with $N$ are estimated separately in
Sections~\ref{W:sec:taylor-remainders-and-bounded-entries}, \ref{W:sec:cumulant-cubic} and~\ref{W:sec:localization-of-the-low-energy-spectral-subspace}.

\subsubsection{Eigenvalue counting estimates for matrices with bounded entries}
\label{W:sec:bounded-support-counting-inputs}

The conditional estimates for the Schur complement reduce the problem to spectral information about complementary minors. Eigenvalue counting provides that information for matrices
with a polynomial bound on their entry supports. The real symmetric statement is recorded with its precise moment hypotheses, and the complex case is proved separately. Rank and
resolvent estimates will then transfer the counting bound to the full truncated array and its minors.

This transfer must work throughout entry replacement, when different independent blocks can have different distributions. The comparison family is therefore described block by
block. Off the diagonal, $H_{ij}=N^{-1/2}\xi_{ij}$, while $H_{ii}=N^{-1/2}d_i$. The real coordinate vector is $x_{ij}=\xi_{ij}$ for real entries and
$x_{ij}=(\sqrt2\operatorname{Re}\xi_{ij},\sqrt2\operatorname{Im}\xi_{ij})$ for complex entries. The blocks are chosen deterministically from the uncentered truncated distributions,
the bounded distributions with matching third moments of Section~\ref{W:sec:moment-bounds-and-taylor-expansion-at-two-truncation-scales}, or the truncated Gaussian distributions.
Uniformly over these choices, the arrays satisfy

\begin{equation}
 \begin{split}
 &\max_{i,j}|H_{ij}|\le e_N,\qquad e_N\longrightarrow0,
       \qquad e_N^{-1}=N^{o(1)},\\
 &\sup_{i<j}\sup_{t\ge1}t^4\mathbb P(|\xi_{ij}|>t)\le C,
       \quad \mathbb E|\xi_{ij}|^3\le C,
       \quad \mathbb E|\xi_{ij}|^4\le C\log(eN),\\
 &|\mathbb E x_{ij}|\le N^{-3/2+o(1)},\qquad
       \|\operatorname{Cov}(x_{ij})-I\|\le N^{-1+o(1)},\\
 &\mathbb E d_i^2\le C,\qquad |\mathbb E d_i|\le N^{-1/2+o(1)}.
 \end{split}
 \label{W:eq:minor-family}
\end{equation}
Uniform integrability of the diagonal squares gives the truncation coupling and the estimate of its mean error. For the truncated off-diagonal entries the fourth bound follows by
integrating the fourth moment tail function to $e_N\sqrt N$; the third absolute moment is finite by the same tail integral. The bounded entries and the truncated Gaussian entries
satisfy these bounds as well. The counting estimate below uses the explicit logarithmic fourth moment bound in \eqref{W:eq:minor-family}, as required by the quoted local law
estimate.

Use the sign convention
\[
 \begin{aligned}
 R_H(z)&=(zI-H)^{-1},& m_H(z)&=N^{-1}\operatorname{Tr}R_H(z),\\
 g(z)&=\int\frac{\rho_{\rm sc}(x)}{z-x}\,\dd x,
 &\upsilon_H(z)&=-\operatorname{Im}m_H(z).
 \end{aligned}
\]
Thus $z-g(z)=g(z)^{-1}$ and $\upsilon_H\ge0$ in the upper half-plane. For a deleted coordinate set $J$, write $R_H^{(J)}=(zI-H_{J^cJ^c})^{-1}$ and
$m_H^{(J)}=N^{-1}\operatorname{Tr}R_H^{(J)}$. All minor traces use the normalization $N^{-1}$.

\Needspace{6\baselineskip} We isolate the eigenvalue counting estimate for matrices with bounded entries that will be used to control complementary minors.

\phantomsection\label{W:def:minor-counting}

Fix $\beta\in\{1,2\}$, $\phi>0$ and $C_0$. Consider matrices $W$ with independent centered real symmetric or complex Hermitian blocks, respectively, with
\[
 \mathbb E|W_{ij}|^2=N^{-1},\quad
 \mathbb E|W_{ij}|^3\le C_0N^{-3/2},\quad
 \mathbb E|W_{ij}|^4\le C_0N^{-2}\log(eN)\quad(i<j),
\]
$\mathbb E W_{ii}^2\le C_0/N$, and $|W_{ij}|\le C_0N^{-\phi}$. For $\beta=2$ also require $\mathbb E W_{ij}^2=0$ off the diagonal; the two real coordinates of a block need not be
independent. For this class we establish the following estimate, first in the real case and then in the Hermitian case: a deterministic $\mathfrak c_N=N^{o(1)}$ satisfies, with
failure at most $N^{-D}$ for every fixed $D$,
\begin{equation}
 \|W\|\le3,\qquad
 \sup_{x\in\mathbb R}|F_W(x)-F_{\rm sc}(x)|\le \mathfrak c_N/N,
 \qquad F_W(x)=N^{-1}\#\{\lambda_j(W)\le x\}.
 \label{W:eq:minor-counting}
\end{equation}
The constants and $\mathfrak c_N$ are uniform over deterministic nonidentical entry distributions satisfying the displayed bounds.

\Needspace{6\baselineskip} The real symmetric counting estimate implies an averaged resolvent bound by integration against the difference of the empirical and semicircle
distribution functions.

\begin{proposition}
\label{W:prop:minor-real-input}
For real symmetric $W$ satisfying the preceding hypotheses, \eqref{W:eq:minor-counting} holds. On the same event,
\begin{equation}
 |m_W(z)-g(z)|\le\frac{\pi \mathfrak c_N}{N\eta},\qquad
       z=E+\mathrm i\eta,\quad\eta>0.
 \label{W:eq:minor-core-trace}
\end{equation}
\end{proposition}
\begin{proof}[Application of the real counting theorem]
Apply Theorem~3.6 of \cite{lee2014necessary}. Its hypotheses are the independence and centering in Definition~2.1 of that paper, off-diagonal variance $N^{-1}$, the stated third
and fourth moment bounds, a bounded diagonal second moment after multiplication by $N$, and bounded support. The distributions may depend on $i,j,N$. For large $N$,
\[
 C_0N^{-\phi}\le N^{-\phi/2},\qquad
 C_0\log(eN)\le2C_0\log N,
\]
so the support parameter is $q=N^{\phi/2}$. The cited rigidity and counting estimates give, with $\varphi_N=(\log N)^{\log\log N}$,
\[
 \|W\|\le2+\varphi_N^{C_1}N^{-2/3}<3,\qquad
 \sup_{|x|\le5}|F_W(x)-F_{\rm sc}(x)|\le\frac{\varphi_N^{C_1}}N.
\]
Their exceptional probability is at most $N^C\exp(-c\varphi_N)$, which is smaller than $N^{-D}$ for each fixed $D$. Set $\mathfrak c_N=\varphi_N^{C_1}$; then
$\log\mathfrak c_N/\log N=C_1(\log\log N)^2/\log N\to0$. Outside $[-5,5]$ the two distribution functions agree on the norm event, proving \eqref{W:eq:minor-counting} on
$\mathbb R$. For finite positive measures $\mu,\nu$, write $F_\mu(x)=\mu(( -\infty,x])$ and $m_\mu(z)=\int(z-x)^{-1}\,\dd\mu(x)$. Their distribution functions are bounded, so
$(F_\mu(x)-F_\nu(x))/(z-x)\to0$ at both infinities. Integration by parts, followed by $\int_{\mathbb R}|z-x|^{-2}\,\dd x=\pi/|\operatorname{Im}z|$, gives
\[
 m_\mu(z)-m_\nu(z)
 =-\int_{\mathbb R}\frac{F_\mu(x)-F_\nu(x)}{(z-x)^2}\,\dd x,
 \qquad
 |m_\mu(z)-m_\nu(z)|
 \le\frac{\pi\|F_\mu-F_\nu\|_\infty}{|\operatorname{Im}z|}.
\]
This calculation applies also when the two total masses differ. With $\mu$ the normalized empirical measure of $W$ and $\nu=\mu_{\rm sc}$, it proves \eqref{W:eq:minor-core-trace}.
\end{proof}

\Needspace{6\baselineskip}
\begin{remark}
\label{W:rem:minor-complex-scope}
The cited Lee--Yin theorem concerns real symmetric matrices. The Hermitian estimate below instead follows from comparison of joint coordinate moments through order three. A
separate real matrix, whose entry moments dominate the absolute moments of the Hermitian entries, gives the preliminary spectral norm bound. The sparse local law does not give the
required precision here: its term $\min\{q^{-2}/\sqrt{\kappa+\eta},q^{-1}\}$ need not be $N^{-1/3+o(1)}$ when the coarsely truncated matrix has $q=N^\phi$ with small $\phi$.
\end{remark}

\subsubsection{Eigenvalue counting for Hermitian matrices}
\label{W:sec:eigenvalue-counting-for-complex-matrices-with-bounded-entries}

We compare high moments of the normalized trace with those of a matrix having bounded standardized entries. Matching the joint coordinate moments through order three cancels the
first three derivatives; the truncated fourth moment bounds the remaining integral. This yields the counting estimate \eqref{W:eq:minor-counting} for Hermitian matrices.

The class $\mathcal W_N(\phi,C_0)$ consists of the Hermitian matrices in the hypotheses preceding \eqref{W:eq:minor-counting}, with $\beta=2$. Its upper-triangular entries are
independent and centered, and satisfy

\begin{equation}
 \begin{split}
 &\mathbb E|H_{ij}|^2=N^{-1},\qquad \mathbb E H_{ij}^2=0,\qquad
 \mathbb E|H_{ij}|^3\le C_0N^{-3/2},\quad i<j,\\
 &\mathbb E|H_{ij}|^4\le C_0N^{-2}\log(eN),\quad i<j,\qquad
 \mathbb E H_{ii}^2\le C_0/N,\\
 &|H_{ij}|\le C_0N^{-\phi}\quad\hbox{almost surely for every }i,j.
 \end{split}
 \label{W:eq:complex-assumptions}
\end{equation}
The laws may depend on $i,j,N$. We use the resolvent and trace notation of Section~\ref{W:sec:notation-and-scales}. Stochastic domination is uniform over this class and the
indicated deterministic parameters.

\Needspace{6\baselineskip} Each off-diagonal entry is treated as one random vector $(\sqrt2\operatorname{Re}H_{ij},\sqrt2\operatorname{Im}H_{ij})$.

\begin{theorem}
\label{W:thm:complex-main}
Under \eqref{W:eq:complex-assumptions}, for each fixed $\tau>0$ and fixed $M>5$,
\begin{equation}
 \sup_{\substack{|E|\le5\\N^{-1+\tau}\le\eta\le M}}
 N\eta\,|m_H(E+\mathrm i\eta)-g(E+\mathrm i\eta)|\prec1.
 \label{W:eq:complex-averaged}
\end{equation}
Moreover, $\|H\|\le3$ with overwhelming probability and
\begin{equation}
 N\sup_{x\in\mathbb R}|F_H(x)-F_{\rm sc}(x)|\prec1.
 \label{W:eq:complex-counting}
\end{equation}
There is a deterministic $\mathfrak c_N=N^{o(1)}$, depending only on $\phi,C_0$, for which the two conclusions in \eqref{W:eq:minor-counting} hold with failure smaller than
$N^{-D}$ for every fixed $D$. This supplies the counting estimate for the Hermitian matrices in Section~\ref{W:sec:minor-estimates-by-coarse-truncation}.
\end{theorem}

For a Wigner matrix with bounded standardized entries, the averaged local semicircle law follows from \cite{benaych2018lectures}. The definition permits nonidentical complex entry
laws and bounded diagonal variances; independence is across upper-triangular entries. For the coarse norm estimate we use the estimates for real matrices with bounded entries of
\cite{lee2014necessary} on a real symmetric matrix whose entry moments dominate the absolute moments in the Hermitian trace expansion. This auxiliary matrix is constructed below.

\Needspace{6\baselineskip} Given a random vector, we construct a finitely supported vector with the same joint moments through order three. The bound on its support depends only on
the third absolute moment.

\begin{lemma}
\label{W:lem:bounded-complex-moment-matching}
Let $X\in\mathbb R^2$ satisfy $\mathbb EX=0$, $\mathbb EXX^{\mathsf T}=I_2$, and $\mathbb E\|X\|^3\le C$. There is a finitely supported $X^{\mathrm b}\in\mathbb R^2$, with support
bounded in terms of $C$ alone, whose moments of every total degree at most three agree with those of $X$. Consequently there is a Hermitian matrix $H^{\mathrm b}$ with
$|\sqrt N H^{\mathrm b}_{ij}|\le C(C_0)$, independent upper-triangular entries, zero diagonal, and the same full off-diagonal real moment tensors through degree three as $H$.
\end{lemma}
\begin{proof}
The four independent coordinates of the symmetric tensor $T=\mathbb E X^{\otimes3}$ are
\[
 a=\mathbb EX_1^3,\qquad b=\mathbb EX_1^2X_2,\qquad
 c=\mathbb EX_1X_2^2,\qquad d=\mathbb EX_2^3.
\]
Each has absolute value at most $C$. Set
\[
 v_1=(1,0)^{\mathsf T},\quad v_2=(0,1)^{\mathsf T},\quad
 v_3=(1,1)^{\mathsf T},\quad v_4=(1,-1)^{\mathsf T},
\]
\[
 c_1=a-c,\quad c_2=d-b,\quad c_3=(b+c)/2,\quad c_4=(c-b)/2.
\]
For the coordinate tuples $111,112,122,222$, the tensor $\sum_\ell c_\ell v_\ell^{\otimes3}$ has entries $c_1+c_3+c_4=a$, $c_3-c_4=b$, $c_3+c_4=c$, and $c_2+c_3-c_4=d$,
respectively. Symmetry supplies all other entries. Direct multiplication of the four column vectors also gives
\begin{equation}
 T=\sum_{\ell=1}^4c_\ell v_\ell^{\otimes3},\qquad
 \sum_{\ell=1}^4v_\ell v_\ell^{\mathsf T}=3I_2.
 \label{W:eq:complex-cubic-decomposition}
\end{equation}

Apply the two-dimensional construction in the proof of Lemma~\ref{W:lem:cubic-defect} to this tensor $T$. Its hypotheses hold because $\mathbb EX=0$, $\mathbb EXX^{\mathsf T}=I_2$
and $|T_{ijk}|\le\mathbb E\|X\|^3\le C$. With the coefficients in \eqref{W:eq:complex-cubic-decomposition}, that construction gives independent centered variables of variance $1/6$
and third moments $c_\ell$, together with the independent symmetric sign vector of covariance $I_2/2$. The resulting $X^{\mathrm b}$ has mean zero, covariance $I_2$ and third
tensor $T$. The support calculation there gives at most $2^6$ values and
\[
 \|X^{\mathrm b}\|\le1+(2+2\sqrt2)(12C+1).
\]
The constant monomial also has expectation one, proving agreement of every monomial of total degree at most three.

For the matrix application, write $H_{ij}=A+\mathrm iB$ and use $X=(\sqrt{2N}A,\sqrt{2N}B)$. The equations $\mathbb E|H_{ij}|^2=N^{-1}$ and $\mathbb EH_{ij}^2=0$ give
\[
 \mathbb EA^2=\mathbb EB^2=(2N)^{-1},\qquad
 \mathbb EAB=0,\qquad \mathbb EXX^{\mathsf T}=I_2,
\]
and $\mathbb E\|X\|^3=(2N)^{3/2}\mathbb E|H_{ij}|^3\le2^{3/2}C_0$. Apply the preceding construction to each unordered off-diagonal pair, independently across pairs, and set
$H_{ij}^{\mathrm b}=(X_1^{\mathrm b}+\mathrm iX_2^{\mathrm b})/\sqrt{2N}$ and $H_{ji}^{\mathrm b}=\overline{H_{ij}^{\mathrm b}}$. Set $H_{ii}^{\mathrm b}=0$. The off-diagonal
coordinate tensors agree through degree three with one support bound valid for every $i,j,N$. The diagonal entries are handled separately in the comparison below.
\end{proof}

We first bound the entrywise error by a negative power of the matrix size. This preliminary estimate controls the derivatives in the comparison argument for eigenvalue counting.

Absorb fixed constants by setting
\[
 \phi_0=\min\{\phi/2,1/4\},\qquad b_N=N^{-\phi_0}.
\]
All entries of $H$, $H^{\mathrm b}$, and their intermediate matrices have absolute value at most $b_N$ for large $N$.

\Needspace{6\baselineskip} We first estimate linear and quadratic forms of an independent complex row.

Let $x_1,\ldots,x_n$, $n\le N$, be independent centered complex variables with $\mathbb E|x_j|^2=N^{-1}$, $|x_j|\le b_N$, and $\mathbb E|x_j|^4\le F_N/N^2$, where
$F_N\le C\log(eN)$.

For a deterministic complex matrix $A$, put
\[
 \mathcal K(A)=\frac{\sqrt{F_N}}N\|\operatorname{diag}A\|_2
 +\frac1N\|A\|_{\mathrm{HS}}
 +\frac{b_N}{\sqrt N}\left(\max_i\|A_{i\cdot}\|_2
                          +\max_j\|A_{\cdot j}\|_2\right)
 +b_N^2\max_{i,j}|A_{ij}|.
\]

\begin{lemma}
\label{W:lem:complex-row}
For each fixed even $p\ge4$ and deterministic Hilbert space vectors $v_j$,
\begin{equation}
 \left\|\sum_jx_jv_j\right\|_{L^p}
 \le C_p\left[N^{-1/2}\Big(\sum_j\|v_j\|^2\Big)^{1/2}
                   +b_N\max_j\|v_j\|\right].
 \label{W:eq:complex-linear}
\end{equation}
\begin{equation}
 \|x^*Ax-N^{-1}\operatorname{Tr}A\|_{L^p}
 \le C_pN^{1/p}\mathcal K(A).
 \label{W:eq:complex-quadratic}
\end{equation}
Both statements hold conditionally for a random coefficient family independent of the row. In particular, for each $\delta,D>0$, a fixed sufficiently large $p$ makes their bounds,
with factor $N^\delta$ in place of $C_pN^{1/p}$, fail with probability at most $N^{-D}$.
\end{lemma}
\begin{proof}

Write $p=2k$ and $S=\sum_jx_jv_j$. Expanding gives
\[
 \mathbb E\|S\|^{2k}
 =\sum_{i_1,j_1,\ldots,i_k,j_k}
   \mathbb E\prod_{a=1}^k x_{i_a}\overline{x_{j_a}}
                  \prod_{a=1}^k\langle v_{i_a},v_{j_a}\rangle.
\]
A singleton index has zero expectation by independence and centering. For a partition $\pi$ of the $2k$ factor positions, taking absolute values and then summing over all
assignments of its blocks bounds its contribution by
\[
 \prod_{B\in\pi}\sum_j\mathbb E|x_j|^{|B|}\|v_j\|^{|B|}.
\]
This enlarges the original sum over injective assignments only after absolute values have been taken. For every remaining partition $\pi$ of the set $\{1,\ldots,2k\}$, with block
sizes $m_1,\ldots,m_\ell\ge2$, independence and $|x_j|\le b_N$ bound the absolute contribution by
\[
 \begin{aligned}
 \prod_{a=1}^\ell\sum_j\mathbb E|x_j|^{m_a}\|v_j\|^{m_a}
 &\le\prod_{a=1}^\ell
  \left(\frac1N\sum_j\|v_j\|^2\right)
               (b_N\max_j\|v_j\|)^{m_a-2}\\
 &=A_0^{2\ell}B_0^{2k-2\ell},\qquad
 A_0^2=N^{-1}\sum_j\|v_j\|^2,\quad B_0=b_N\max_j\|v_j\|.
 \end{aligned}
\]
Since $\ell\le k$, $A_0^{2\ell}B_0^{2k-2\ell}\le(A_0+B_0)^{2k}$. There are at most $(2k)^{2k}$ partitions; taking the $2k$th root proves \eqref{W:eq:complex-linear}. For
independent centered bounded scalar $Y_j$, put $V=\sum_j\mathbb E|Y_j|^2$ and $B=\max_j\|Y_j\|_\infty$. For a partition into $\ell$ blocks of sizes $m_a\ge2$, the scalar moment
contribution is at most
\[
 \prod_{a=1}^{\ell}\sum_j\mathbb E|Y_j|^{m_a}
 \le V^{\ell}B^{p-2\ell}\le(\sqrt V+B)^p.
\]
Summing over the partitions gives
\[
 \left\|\sum_jY_j\right\|_p\le C_p(\sqrt V+B).
\]

Apply that scalar bound to $Y_j=A_{jj}(|x_j|^2-N^{-1})$. Since $\mathbb E|x_j|^4\le F_N/N^2$ and $N^{-1}=\mathbb
E|x_j|^2\le b_N^2$,
\[
 \sum_j\mathbb E|Y_j|^2
 \le\frac{F_N}{N^2}\sum_j|A_{jj}|^2,
 \qquad \max_j\|Y_j\|_\infty\le b_N^2\max_j|A_{jj}|.
\]
Consequently
\[
 \left\|\sum_jA_{jj}(|x_j|^2-N^{-1})\right\|_p
 \le C_p\left[\frac{\sqrt{F_N}}N\|\operatorname{diag}A\|_2
                   +b_N^2\max_j|A_{jj}|\right].
\]
For the off-diagonal part, independent Bernoulli selectors satisfy
\[
 \sum_{i\ne j}\overline x_iA_{ij}x_j
 =4\mathbb E_\varepsilon
 \sum_{\varepsilon_i=1,\,\varepsilon_j=0}\overline x_iA_{ij}x_j.
\]
For each partition of indices the two row subfamilies are independent. Condition on the second, call it $y$, and apply \eqref{W:eq:complex-linear} to the first. This gives
\[
 C_p\left[N^{-1/2}\|\|Ay\|_2\|_{L^p}
                  +b_N\|\max_i|(Ay)_i|\|_{L^p}\right].
\]
The Hilbert space version of \eqref{W:eq:complex-linear} bounds the first expression by
\[
 C_p\left[N^{-1}\|A\|_{\mathrm{HS}}
 +b_NN^{-1/2}\max_j\|A_{\cdot j}\|_2\right].
\]

For the other term, the scalar estimate gives
\[
 \begin{aligned}
 b_N\|\max_i|(Ay)_i|\|_p
 &\le b_N\left(\sum_i\mathbb E|(Ay)_i|^p\right)^{1/p}\\
 &\le C_pN^{1/p}\left[
       b_NN^{-1/2}\max_i\|A_{i\cdot}\|_2
                 +b_N^2\max_{i,j}|A_{ij}|\right].
 \end{aligned}
\]
Add this to the Hilbert space estimate and then to the diagonal contribution. The norms of every submatrix selected by
the Bernoulli variables are bounded by the corresponding norms of $A$. Minkowski's inequality in the finite selector
average now proves \eqref{W:eq:complex-quadratic}. Finally Markov's inequality gives an upper bound $C_p^pN^{1-\delta
p}$ for the normalized quadratic deviation; choose $p>(D+2)/\delta$. For $N\ge C_p^p$, the bound is at most $N^{2-\delta
p}\le N^{-D}$.
\end{proof}

\Needspace{6\baselineskip} These row bounds first give a preliminary entrywise local law, stable under setting a single selected entry to zero and varying that entry.

\begin{lemma}
\label{W:lem:complex-weak}
Fix $\tau>0$ and fixed $M>5$. For $H$, $H^{\mathrm b}$, and every fixed family of polynomially many intermediate matrices between them, there is $c>0$ such that, with overwhelming
probability,
\begin{equation}
 \sup_{\substack{|E|\le5\\N^{-1+\tau}\le\eta\le M}}
 \max_{i,j}|R_{ij}(z)-g(z)\delta_{ij}|\le N^{-c}.
 \label{W:eq:complex-weak}
\end{equation}
In particular all these resolvent entries are bounded by a fixed constant. Deleting one selected entry and varying it arbitrarily within its support preserves a fixed entrywise
bound along the entire radial segment.
\end{lemma}
\begin{proof}
It is enough to take $0<\tau<1$; larger $\tau$ restrict the domain. Enlarge $M$, if needed, to a sufficiently large fixed number. Use a partition of spacing $N^{-K}$ in the
displayed rectangle, with $K$ a fixed large integer. Conditional on the minor $H^{(i)}$, apply Lemma~\ref{W:lem:complex-row} to the off-diagonal part of row $i$. The conditional
means are exactly zero and all variances are exactly $N^{-1}$ in every intermediate matrix. Write
\[
 R^{(i)}=(z-H^{(i)})^{-1},\qquad
 {m_H}^{(i)}=N^{-1}\operatorname{Tr}R^{(i)},\qquad
 Z_i=h_i^*R^{(i)}h_i-{m_H}^{(i)}.
\]
For each fixed minor, Lemma~\ref{W:lem:complex-row} gives
\[
 |Z_i|\le N^\delta\mathcal K(R^{(i)}),\qquad
 |h_i^*R^{(i)}e_j|
 \le N^\delta\bigl(N^{-1/2}\|R^{(i)}e_j\|_2
                        +b_N\|R^{(i)}e_j\|_\infty\bigr)
\]
except with conditional probability at most $C_p^pN^{1-\delta p}$. Let $N^J$ bound the number of choices of indices, sampling points and intermediate matrices. Integration over the
minors and a union bound give exceptional probability at most $C_p^pN^{J+1-\delta p}$. Given $D>0$, choose the even integer $p>(D+J+2)/\delta$ before taking $N\to\infty$. This
probability is then at most $N^{-D}$ for all sufficiently large $N$. The conditional estimates hold for every minor. Bounds on its resolvent will be deduced below under the
temporary bound $\Lambda\le\epsilon_0$.

Let $\Lambda(z)=\max_{ij}|R_{ij}(z)-g(z)\delta_{ij}|$. Choose $\epsilon_0=(100)^{-1}\inf_{|E|\le5,\,0\le\eta\le M}|g(E+\mathrm i\eta)|>0$. The minimum is positive by continuity of
$g$ up to the real axis and the fact that $g$ has no zero there. We first take $M$ sufficiently large and then use this fixed $\epsilon_0$. On $\Lambda\le\epsilon_0$, $|R_{ii}|$ is
bounded below. The exact minor identity
\[
 R^{(i)}_{jk}=R_{jk}-R_{ji}R_{ik}/R_{ii}
\]
then bounds all entries of $R^{(i)}$ by a fixed constant. Interlacing and integration by parts give
\begin{equation}
 |{m_H}-{m_H}^{(i)}|\le\frac\pi{N\eta}.
 \label{W:eq:complex-interlacing}
\end{equation}
The Ward identities and bounded diagonal entries give
\[
 \|R^{(i)}\|_{\mathrm{HS}}\le C\sqrt{N/\eta},\qquad
 \max_j\bigl(\|R^{(i)}e_j\|_2+\|e_j^*R^{(i)}\|_2\bigr)
 \le C\eta^{-1/2}.
\]
Consequently, under the temporary bound $\Lambda\le\epsilon_0$,
\[
 \mathcal K(R^{(i)})\le C\left[
 \sqrt{F_N/N}+(N\eta)^{-1/2}
       +b_N(N\eta)^{-1/2}+b_N^2\right].
\]
Choose $0<\delta<\min\{\phi_0,\tau\}/100$. For $\eta\ge N^{-1+\tau}$ and $b_N\le C_0N^{-\phi_0}$,
\[
 N^\delta\mathcal K(R^{(i)})
 \le N^{-1/2+\delta+o(1)}
       +C N^{-\tau/2+\delta}
       +C N^{-\phi_0-\tau/2+\delta}
       +C N^{-2\phi_0+\delta}.
\]
With $c_1=\min\{1/4,\phi_0/2,\tau/4\}>0$, this bound and $b_N+(N\eta)^{-1}$ are at most $r_N=N^{-c_1}$ for all sufficiently large $N$. The Schur complement therefore gives
\[
 R_{ii}^{-1}=z-{m_H}-\epsilon_i,\qquad
 \epsilon_i=H_{ii}+Z_i+{m_H}^{(i)}-{m_H},\qquad |\epsilon_i|\le Cr_N.
\]
The same temporary assumption bounds $z-{m_H}$ away from zero. Averaging yields
\begin{equation}
 |{m_H}-(z-{m_H})^{-1}|\le Cr_N,\qquad
 |{m_H}^2-z{m_H}+1|\le C_Mr_N.
 \label{W:eq:complex-quadratic-stability}
\end{equation}

The quadratic relation in \eqref{W:eq:complex-quadratic-stability} factors as
\[
 (m_H-g)(m_H-z+g)=m_H^2-zm_H+1.
\]
Since the product has modulus at most $C_Mr_N$, at least one factor has modulus at most $a_N=(C_Mr_N)^{1/2}$. To determine which root is relevant, fix a real part $E$ in the chosen
finite partition and vary $\eta$ from $M$ down to $N^{-1+\tau}$. Write $g_+(z)=g(z)$ and $g_-(z)=z-g(z)$. The roots are distinct when $\eta>0$, and
\[
 |g_+(z)-g_-(z)|^2=|z^2-4|,
 \qquad |z^2-4|^2=(E^2-4)^2+2(E^2+4)\eta^2+\eta^4.
\]
In particular, $|g_+-g_-|$ increases strictly with $\eta$. At $\eta=M$, $|m_H|\le M^{-1}$ while $|g_-|\ge M-M^{-1}$; choosing $M$ sufficiently large excludes $|m_H-g_-|\le a_N$.
Thus $|m_H-g_+|\le a_N$ at the initial height.

On the portion of the vertical interval where $|g_+-g_-|>4a_N$, the two closed disks of radius $a_N$ about the roots are disjoint. The continuous function $m_H(E+\mathrm i\eta)$
lies in their union by the factored equation, and cannot pass from one disk to the other while the disks remain disjoint. Since their separation is increasing in $\eta$, this
portion is an interval adjoining $M$. Hence $|m_H-g|\le a_N$ there. At every remaining height $|g_+-g_-|\le4a_N$; whichever factor is small gives
\[
 |m_H-g|\le |m_H-g_-|+|g_--g|\le5a_N.
\]
This shows $|m_H-g|\le5a_N$ throughout the connected interval on which the temporary hypothesis $\Lambda\le\epsilon_0$ holds. Using $R_{ii}^{-1}=z-m_H-\epsilon_i$, $g^{-1}=z-g$,
and the lower bound on $|R_{ii}|$, we obtain the exact identity
\[
 R_{ii}-g=R_{ii}g\,(m_H-g+\epsilon_i),
 \qquad \max_i|R_{ii}-g|\le C_M(a_N+r_N).
\]
The off-diagonal block inverse formula reads
\[
 R_{ij}=R_{ii}\,h_i^*R^{(i)}e_j,\qquad i\ne j,
\]
where $e_j$ is indexed in the complementary minor. On the same interval the conditional linear bound, the Ward estimate for $R^{(i)}$ and $\max_{k,l}|R^{(i)}_{kl}|\le C_M$ imply
\[
 \max_{i\ne j}|R_{ij}|
 \le C_MN^\delta\left((N\eta)^{-1/2}+b_N\right)
 \le C_M\left(N^{-\tau/2+\delta}+N^{-\phi_0+\delta}\right).
\]
Every exponent on the right is negative. Taking $N$ large after $\epsilon_0$ is fixed, these two estimates give $\Lambda\le\epsilon_0/2$ whenever $\Lambda\le\epsilon_0$ and the
conditional row inequalities hold.

At height $M$, $\|R^{(i)}\|\le M^{-1}$, and the deterministic bounds in the definition of $\mathcal K$ give
\[
 \mathcal K(R^{(i)})\le C_M\left(\sqrt{F_N/N}+N^{-1/2}+b_NN^{-1/2}+b_N^2\right).
\]
The conditional row inequalities therefore give $Z_i=o(1)$ and $h_i^*R^{(i)}e_j=o(1)$ uniformly over $i,j$ on the common event. Since $|H_{ii}|\le b_N=o(1)$ and
$|m_H^{(i)}|\le M^{-1}$, the Schur complement yields, uniformly for $|E|\le5$,
\[
 R_{ii}(E+\mathrm iM)=z^{-1}+O(M^{-3})+o(1),
 \qquad \max_{i\ne j}|R_{ij}(E+\mathrm iM)|=o(1).
\]
For the deterministic transform, $zg=1+g^2$ and $|g(E+\mathrm iM)|\le M^{-1}$ give
\[
 |g(E+\mathrm iM)-z^{-1}|=|g(E+\mathrm iM)^2/z|\le M^{-3}.
\]
The choice of $M$ is consistent with the definition of $\epsilon_0$. Indeed, $z=g+g^{-1}$ and $\operatorname{Im}g<0$ imply $|g|<1$ in the upper half-plane: taking imaginary parts
gives $\operatorname{Im}z=\operatorname{Im}g(1-|g|^{-2})>0$. Continuity supplies $|g|\le1$ on its boundary, and therefore
\[
 |g(z)|=|z-g(z)|^{-1}\ge(M+6)^{-1}
       \quad(|\operatorname{Re}z|\le5,\ 0\le\operatorname{Im}z\le M).
\]
Thus $\epsilon_0\ge[100(M+6)]^{-1}$. The constants in the preceding $O(M^{-3})$ bounds are absolute for $M\ge2$, so $M$ can first be fixed with their sum smaller than
$\epsilon_0/8$. Taking $N$ large afterwards makes the uniform $o(1)$ terms smaller than $\epsilon_0/8$ and gives $\Lambda(E+\mathrm iM)<\epsilon_0/4$.

To extend the conditional estimates from the finite set of spectral parameters to the rectangle, we use the following continuity bounds. For $\eta\ge N^{-1+\tau}$, the resolvent
identity gives
\[
 \|R(z)-R(z')\|\le N^{2-2\tau}|z-z'|,
 \qquad \|R^{(i)}(z)-R^{(i)}(z')\|\le N^{2-2\tau}|z-z'|
\]
when the straight segment between $z,z'$ stays in the rectangle. The scalar function $g$ has $|g'(z)|\le N^{2-2\tau}$. Using $\|h_i\|_2\le\sqrt N$, one obtains for the difference
between two spectral points
\[
 |Z_i(z)-Z_i(z')|\le2N^{3-2\tau}|z-z'|,
 \quad |h_i^*(R^{(i)}(z)-R^{(i)}(z'))e_j|
 \le N^{5/2-2\tau}|z-z'|.
\]
The functions $\mathcal K(R^{(i)})$ are sums of matrix norms multiplied by deterministic coefficients, so their differences are bounded by $C N^4|z-z'|$ on this rectangle. Take a
two-dimensional partition of spacing $N^{-K}$ with $K>6+c_1$, including its boundary heights. The row inequalities at the closest sampling point then hold everywhere with an
additional $CN^{4-K}=o(r_N)$, whereas $\Lambda$ changes by at most $CN^{2-K}$. The number of sampling points is polynomial in $N$; the conditional moment order was chosen above to
make the joint probability of all failures at most $N^{-D}$ for arbitrary fixed $D$.

Fix a real part in the partition. If the descending vertical interval on which $\Lambda\le\epsilon_0$ had a first lower boundary before $N^{-1+\tau}$, continuity would give
$\Lambda=\epsilon_0$ there. The just established estimates, including the partition approximation errors, instead give $\Lambda\le\epsilon_0/2$, a contradiction. Thus the entire
vertical interval is covered. The same resolvent Lipschitz estimate fills the gaps between successive real parts. Combining the bounds on the diagonal and off-diagonal entries, and
reducing the positive exponent if necessary, proves \eqref{W:eq:complex-weak} for the full rectangle. The event remains uniform over the fixed polynomial family of intermediate
matrices.

It remains to construct the event needed when a selected entry is replaced. Let $\alpha$ be an unordered off-diagonal pair, $K_\alpha$ the matrix with that block set to zero, and
$I$ its two coordinates. For $V=\iota_IT\iota_I^*$ with $\|T\|\le Cb_N$, the Woodbury identity is
\begin{equation}
 R_{K_\alpha+V}=R_{K_\alpha}+R_{K_\alpha,\cdot I}T(I-R_{K_\alpha,II}T)^{-1}R_{K_\alpha,I\cdot}.
 \label{W:eq:complex-woodbury}
\end{equation}
Define the event $\mathcal G_\alpha\in\sigma(K_\alpha)$ by
\[
 \mathcal G_\alpha=\left\{
 \sup_{\substack{|E|\le5\\N^{-1+\tau}\le\eta\le M}}
 \max_{a,b}|R_{K_\alpha}(E+\mathrm i\eta)_{ab}|\le3\right\}.
\]
This is measurable because resolvent entries are continuous for nonreal $z$ and the supremum can be taken on a countable dense subset. On the event in \eqref{W:eq:complex-weak},
deleting the actual sampled entry is a perturbation of size at most $Cb_N$. Apply the block identity first with base matrix $H$ and perturbation $-H_\alpha$. Its compressed
resolvent has norm at most a fixed constant since $|I|\le2$, so $(I-R_{H,II}T)^{-1}$ exists with norm at most two for all sufficiently large $N$. Each matrix entry changes by at
most $C b_N$, uniformly in $z$. Thus the resolvent estimate for the original matrix implies $\mathcal G_\alpha$; independence of $K_\alpha$ and the sampled block gives
\[
 \mathbb P_{K_\alpha}(\mathcal G_\alpha^c)
 =\mathbb P_{K_\alpha,H_\alpha}(\mathcal G_\alpha^c)
 \le\mathbb P_{K_\alpha,H_\alpha}(\text{the resolvent estimate fails})\le N^{-D}.
\]
On $\mathcal G_\alpha$, the Woodbury identity with base matrix $K_\alpha$ gives, for every permitted perturbation $V$ supported on the selected coordinates,
\[
 \sup_{\substack{|E|\le5\\N^{-1+\tau}\le\eta\le M}}
 \max_{a,b}|R_{K_\alpha+V}(E+\mathrm i\eta)_{ab}|
 \le3+Cb_N\le4.
\]
This includes every point of either replacement segment. Since $\mathcal G_\alpha$ depends only on the other entries, its indicator is constant during the conditional integration
over the replacement entry.
\end{proof}

\Needspace{6\baselineskip} \subparagraph{Comparison of resolvent trace moments}
\label{W:sec:resolvent-trace-comparison-for-complex-entries}

Bounds on resolvent entries control derivatives of normalized spectral traces and their powers when one independent entry is replaced.

\begin{lemma}
\label{W:lem:complex-derivatives}
Suppose $\max_{ij}|R_{ij}(z)|\le C$ and put $s=(N\eta)^{-1}$. Let $T_j=T_j^*$, $1\le j\le q\le4$, have uniformly bounded entries and support on the coordinates of one common matrix
entry. Here $D_H^q f[T_1,\ldots,T_q]$ denotes
$\left.\partial_{t_1}\cdots\partial_{t_q}f(H+\sum_jt_jT_j)\right|_{t=0}$.
Then
\begin{equation}
 |D_{H}^qm_H[T_1,\ldots,T_q]|\le C_qs.
 \label{W:eq:complex-trace-derivatives}
\end{equation}
For an integer $p\ge2$ and $\Phi_p(H)=|m_H(z)-g(z)|^{2p}$,
\begin{equation}
 |D_{H}^q\Phi_p[T_1,\ldots,T_q]|
 \le C_p\sum_{j=1}^{\min(q,2p)}s^j|m_H-g|^{2p-j}.
 \label{W:eq:complex-power-derivatives}
\end{equation}
For any two Hermitian matrices differing on one such block,
\begin{equation}
 |m_A(z)-m_B(z)|\le2\pi s.
 \label{W:eq:complex-rank-trace}
\end{equation}
\end{lemma}
\begin{proof}
Differentiate the inverse with respect to $H$:
\[
 D_{H}^qm_H[T_1,\ldots,T_q]
 =\frac1N\sum_{\sigma\in\mathfrak S_q}
 \operatorname{Tr}(RT_{\sigma(1)}R\cdots T_{\sigma(q)}R).
\]

Expand $T_j=\sum_{a,b}t^{(j)}_{ab}E_{ab}$ on its support. For each choice of matrix units the trace is exactly
\[
 \operatorname{Tr}(RE_{a_1b_1}R\cdots E_{a_qb_q}R)
 =R_{b_1a_2}\cdots R_{b_{q-1}a_q}(R^2)_{b_qa_1}.
\]
There are at most $4^q$ choices for each permutation, and their coefficients are bounded by a constant depending only on $q$. With $L=\max_{j,a,b}|t^{(j)}_{ab}|$ and
$C_R=\max_{a,b}|R_{ab}|$, each such trace has modulus at most $C_R^{q-1}\max_{a,b}|(R^2)_{ab}|$. Summing the $q!4^q$ possible terms therefore gives
\[
 |D_H^qm_H[T_1,\ldots,T_q]|
 \le\frac{q!4^qL^qC_R^{q-1}}N
                         \max_{a,b}|(R^2)_{ab}|.
\]
For $q=1$ the product of interior factors is empty and equals one. By the Ward identity and the Cauchy--Schwarz inequality,
\[
 |(R^2)_{ab}|\le
 \Big(\sum_j|R_{aj}|^2\Big)^{1/2}
 \Big(\sum_j|R_{jb}|^2\Big)^{1/2}
 =\frac{\sqrt{(-\operatorname{Im}R_{aa})
                    (-\operatorname{Im}R_{bb})}}\eta
 \le C/\eta.
\]

This proves \eqref{W:eq:complex-trace-derivatives}. Put $Z=m_H-g$ and $f_a=Z$ for $1\le a\le p$, $f_a=\overline Z$ for $p<a\le2p$. For the set $S=\{1,\ldots,q\}$ the full
derivative is
\[
 D_S\Phi_p=\sum_{S_1\sqcup\cdots\sqcup S_{2p}=S}
                                  \prod_{a=1}^{2p}D_{S_a}f_a.
\]
An assignment with $j$ nonempty parts has modulus at most $C_q^j s^j|Z|^{2p-j}$. For fixed $p,q$, there are at most $(2p)^q$ assignments in total. Summing for $1\le j\le\min(q,2p)$
proves \eqref{W:eq:complex-power-derivatives}. If $Z=0$, the only nonzero possibilities have $j=2p$, and the corresponding factor $|Z|^0$ is interpreted as one. For the final trace
estimate, rank interlacing gives $\sup_x|F_A(x)-F_B(x)|\le2/N$; integrate by parts and use $\int|z-x|^{-2}\,\dd x=\pi/\eta$.
\end{proof}

\Needspace{6\baselineskip} Replacing the independent entries transfers high moments of the averaged resolvent from the matrix with uniformly bounded standardized entries to the
matrix satisfying the stated support bound.

\begin{proposition}
\label{W:prop:complex-comparison}
Let $H^0=H^{\mathrm b},H^1,\ldots,H^{M_N}=H$ be any sequential replacement of upper-triangular entries, where $M_N=N(N+1)/2$. For each fixed integer $p\ge2$, fixed $z$ in the
domain of Lemma~\ref{W:lem:complex-weak}, and arbitrarily large fixed $K$, let
\[
 Y=\max_{0\le\gamma\le M_N}
       \bigl(\mathbb E|m_{H^\gamma}(z)-g(z)|^{2p}\bigr)^{1/(2p)},
 \qquad
 Y_0=\bigl(\mathbb E|m_{H^0}(z)-g(z)|^{2p}\bigr)^{1/(2p)}.
\]
Then, with $s=(N\eta)^{-1}$,
\begin{equation}
 Y^{2p}\le Y_0^{2p}
 +C_p\log(eN)\sum_{j=1}^{4}s^j(Y+C s)^{2p-j}+N^{-K},
 \label{W:eq:complex-moment-inequality}
\end{equation}
and hence
\begin{equation}
 Y\le C_p\bigl[Y_0+\log(eN)s+N^{-K/(2p)}\bigr].
 \label{W:eq:complex-moment-solution}
\end{equation}
\end{proposition}
\begin{proof}
Order the independent upper triangular blocks and write $K_\alpha$ for the matrix obtained by removing the block currently being replaced. The two possible values of this block,
with laws $\mu_\alpha$ and $\nu_\alpha$, are independent of $K_\alpha$. In an off-diagonal block $\alpha=\{a,b\}$ use the two real coordinates
\[
 u=(\sqrt2\operatorname{Re}H_{ab},\sqrt2\operatorname{Im}H_{ab}),\quad
 \mathcal T_\alpha u=u_1T_{\alpha,1}+u_2T_{\alpha,2},\quad
 T_{\alpha,1}=\frac{E_{ab}+E_{ba}}{\sqrt2},\quad
 T_{\alpha,2}=\frac{\mathrm i(E_{ab}-E_{ba})}{\sqrt2}.
\]
Thus $K_\alpha+\mathcal T_\alpha u$ is Hermitian and $\|\mathcal T_\alpha u\|\le C\|u\|_1$. We integrate against the joint distribution of the two coordinates of $u$, allowing
dependence between them. Lemma~\ref{W:lem:bounded-complex-moment-matching} gives, for every $0\le q\le3$, equality of the full symmetric tensors
\begin{equation}
 \int u^{\otimes q}\,\mu_\alpha(\dd u)=\int u^{\otimes q}\,\nu_\alpha(\dd u),\qquad
 \int\|u\|_1^4(\mu_\alpha+\nu_\alpha)(\dd u)\le C\log(eN)N^{-2}.
 \label{W:eq:complex-matching}
\end{equation}
The fourth moment estimate follows from $\|u\|_1\le2|H_{ab}|$ and the bounded standardized fourth moments of the comparison variables. In particular, matching includes each
$\mathbb E[u_1^{q-r}u_2^r]$ with $0\le r\le q\le3$.

Let $\mathcal G_\alpha\in\sigma(K_\alpha)$ be the event constructed in Lemma~\ref{W:lem:complex-weak}. For every fixed $D$ its complement has probability at most $N^{-D}$, and on
$\mathcal G_\alpha$ every matrix $K_\alpha+t\mathcal T_\alpha u$ with $0\le t\le1$ and $u$ in either entry support has uniformly bounded resolvent entries on the specified spectral
domain. Conditional on $K_\alpha$ the map
\[
 f_\alpha(u)=\Phi_p(K_\alpha+\mathcal T_\alpha u)
\]
is four times continuously differentiable in the real vector $u$. Its ordinary Taylor formula, with integral remainder, is
\[
 f_\alpha(u)=\sum_{q=0}^3\frac1{q!}D^qf_\alpha(0)[u^{\otimes q}]
 +\frac1{3!}\int_0^1(1-t)^3D_H^4\Phi_p(K_\alpha+t\mathcal T_\alpha u)[(\mathcal T_\alpha u)^4] \,\dd t.
\]
For fixed $K_\alpha$ the coefficients in the finite sum are deterministic symmetric multilinear forms. Integrating these terms against $\mu_\alpha-\nu_\alpha$ gives zero by
\eqref{W:eq:complex-matching}. The same cancellation remains true with the factor $\mathbf1_{\mathcal G_\alpha}$, because that factor is $\sigma(K_\alpha)$ measurable.

Write $Z_\alpha=m_{K_\alpha}(z)-g(z)$ and $s=(N\eta)^{-1}$. Rank interlacing, with rank at most two, gives for every $t,u$
\[
 |m_{K_\alpha+t\mathcal T_\alpha u}(z)-g(z)|\le |Z_\alpha|+2\pi s.
\]
On $\mathcal G_\alpha$, the fourth derivative estimate \eqref{W:eq:complex-power-derivatives} and multilinearity imply
\[
 |D_H^4\Phi_p(K_\alpha+t\mathcal T_\alpha u)[(\mathcal T_\alpha u)^4]|
 \le C_p\|u\|_1^4\sum_{j=1}^4s^j(|Z_\alpha|+2\pi s)^{2p-j}.
\]
Integrate $t$ using $\int_0^1(1-t)^3\,\dd t=1/4$, and then integrate $u$ under each of the two laws. Conditional on $K_\alpha$, the absolute difference of the two expectations,
restricted to $\mathcal G_\alpha$, is at most
\[
 C_p\log(eN)N^{-2}\sum_{j=1}^4s^j(|Z_\alpha|+Cs)^{2p-j}.
\]

On $\mathcal G_\alpha^c$, the estimate $|m_A(z)|\le\eta^{-1}$ for a Hermitian matrix, together with $|g(z)|\le\eta^{-1}$, gives $0\le\Phi_p(A)\le(2/\eta)^{2p}$. The sum of the
exceptional contributions from all $M_N\le N^2$ blocks is therefore at most
\[
 2N^2(2/\eta)^{2p}\sup_\alpha\mathbb P(\mathcal G_\alpha^c)
 \le C_pN^{2+2p-D}.
\]
Given the fixed error exponent $K$, choose $D>K+2p+3$ in the preceding conditional event estimate; the last quantity is then bounded by $N^{-K}$ for all large $N$.

The diagonal replacements require a separate calculation. The comparison replaces each diagonal entry by zero. The original diagonal entries are centered with
$\mathbb E|H_{aa}|^2\le C_0/N$. For $f_a(u)=\Phi_p(K_a+uE_{aa})$, the constant and linear terms cancel between the two laws because $\mathbb E H_{aa}=0$. The exact second order
remainder is
\[
 f_a(u)-f_a(0)-u f_a'(0)
 =\int_0^1(1-t) D_H^2\Phi_p(K_a+tuE_{aa})[uE_{aa},uE_{aa}]\,\dd t.
\]
The event for this single coordinate is constructed exactly as in the last paragraph of Lemma~\ref{W:lem:complex-weak}, replacing the perturbation on two coordinates by $uE_{aa}$.
The rank-one Woodbury identity bounds the resolvents of the matrix with that entry removed and of every matrix obtained by replacing the entry with a value in its support, and the
event is measurable with respect to the other entries. The second derivative estimate therefore gives a conditional difference at most
\[
 \frac{C_p}{N}\sum_{j=1}^2s^j(|m_{K_a}-g|+Cs)^{2p-j}.
\]
There are $N$ diagonal blocks; their sum contributes at most $C_p\sum_{j=1}^2s^j(Y+Cs)^{2p-j}$. Its exceptional contribution is included in the preceding choice of $D$. The
diagonal contribution is consequently controlled by $\mathbb E|H_{aa}|^2\le C_0/N$.

For every block $\alpha$, either neighboring intermediate matrix differs from $K_\alpha$ by a matrix of rank at most two. Consequently, by the deterministic trace estimate
\eqref{W:eq:complex-rank-trace} and Minkowski's inequality,
\[
 \bigl(\mathbb E|Z_\alpha|^{2p}\bigr)^{1/(2p)}\le Y+2\pi s,
 \quad
 \mathbb E(|Z_\alpha|+Cs)^{2p-j}\le(Y+C's)^{2p-j},\qquad 1\le j\le4.
\]
The second inequality is H\"older's inequality applied to the nonnegative random variable $|Z_\alpha|+Cs$. Sum the $M_N$ successive expectation differences before taking the
maximum over the intermediate matrices. The off-diagonal fourth moments contribute $N^2\cdot C\log(eN)N^{-2}=C\log(eN)$; the diagonal second moments contribute $N\cdot C/N=C$. With
the exceptional term already bounded by $N^{-K}$ this proves \eqref{W:eq:complex-moment-inequality}.

To solve it, set $B_0=Y_0+Cs+N^{-K/(2p)}$ with $C$ large enough for the shifts in the preceding bounds. If $Y\le4B_0$, the desired estimate follows. Otherwise $Y>4B_0$, and
division of \eqref{W:eq:complex-moment-inequality} by $Y^{2p}$ permits the $Y_0^{2p}$ and $N^{-K}$ terms to be absorbed into the left side. Since $Cs/Y<1/4$, for $1\le j\le4$ we
have
\[
 \frac{s^j(Y+Cs)^{2p-j}}{Y^{2p}}
 \le C_p(s/Y)^j,
 \qquad \sum_{j=1}^4(s/Y)^j\le2s/Y.
\]
It follows that $1/2\le C_p\log(eN)s/Y$, hence $Y\le C_p\log(eN)s$. Combining the two cases gives \eqref{W:eq:complex-moment-solution}. The constants are uniform over the stated
family of joint entry distributions and over the replacement order.
\end{proof}

\begin{proof}[Proof of the averaged estimate in Theorem~\ref{W:thm:complex-main}]
The matrix $H^{\mathrm b}$ of Lemma~\ref{W:lem:bounded-complex-moment-matching} has independent centered upper-triangular entries, off-diagonal variance $N^{-1}$, zero diagonal,
and uniformly bounded standardized moments of every fixed order. These are the hypotheses of Definition~2.2 in \cite{benaych2018lectures}; its diagonal variance condition permits
zero. Choose $0<\tau'\le\min\{\tau,1/5,1/M\}$. The domain in Theorem~2.6 of that paper then contains the present rectangle, since
\[
 N^{-1+\tau'}\le N^{-1+\tau},\qquad
 5\le(\tau')^{-1},\qquad M\le(\tau')^{-1}.
\]
Its averaged estimate applies to the joint complex entry law. Both Stieltjes transforms in that reference have the opposite sign to ours, so the absolute difference is unchanged.
For every fixed $p$ and $\delta>0$, this estimate and $|m_{H^{\mathrm b}}-g|\le2\eta^{-1}$ imply
\[
 Y_0\le C_{p,\delta}N^\delta s.
\]
Indeed, on the event $|m_{H^{\mathrm b}}-g|\le N^\delta s$ use this bound, and on its complement use $|m_{H^{\mathrm b}}-g|\le2\eta^{-1}$. Choosing the failure probability at most
$N^{-D_0}$ gives
\[
 Y_0^{2p}\le(N^\delta s)^{2p}+2^{2p}\eta^{-2p}N^{-D_0}
 =s^{2p}\bigl[N^{2p\delta}+2^{2p}N^{2p-D_0}\bigr].
\]
Any $D_0>2p+1$ proves the required bound, uniformly over the domain. Taking $K$ large in \eqref{W:eq:complex-moment-solution} yields
\begin{equation}
 \max_\gamma\mathbb E|m_{H^\gamma}(z)-g(z)|^{2p}
 \le C_{p,\delta}(N^\delta s)^{2p}
 \label{W:eq:complex-high-moment}
\end{equation}
For a fixed $\delta>0$, apply the preceding estimates with $\delta/2$ and use $C_p\log(eN)\le N^{\delta/2}$ for sufficiently large $N$. For $\epsilon,D>0$, take $\delta<\epsilon/2$
and then fixed $p$ sufficiently large. Markov's inequality gives $\mathbb P(|m_H-g|>N^\epsilon s)\le N^{-D}$. All constants are uniform over the given class: the support bound for
the standardized entries of $H^{\mathrm b}$, the preliminary row bounds, and the moment coefficients depend only on its fixed parameters. Use a net of spacing $N^{-6}$ in the two
real coordinates of the rectangle; it has at most $CN^{12}$ points. For a family with at most $N^a$ matrices, choose the pointwise failure exponent $D+a+14$. The union probability
is at most $CN^{-D-2}$. On the rectangle, $|\partial_zm_H|\le\eta^{-2}\le N^2$ and $|g'(z)|\le C_MN^2$, so moving to the nearest net point changes $m_H-g$ by at most $C_MN^{-4}$.
The heights of a point and its nearest net point have ratio tending to one uniformly, because $N^{-6}/\eta\le N^{-5-\tau}$; hence their values of $s=(N\eta)^{-1}$ differ by a
factor at most two. Since $s\ge (NM)^{-1}$, the interpolation error is smaller than $N^{\epsilon/2}s$ for every $N\ge N_0(\epsilon,M)$. Apply the pointwise estimate with
$\epsilon/2$ to obtain the desired $N^\epsilon s$ bound throughout the rectangle. This proves \eqref{W:eq:complex-averaged} uniformly over any fixed family of polynomially many
intermediate matrices.
\end{proof}

\Needspace{6\baselineskip} A comparison of trace moments with a real symmetric matrix gives the preliminary spectral norm bound.

\begin{lemma}
\label{W:lem:complex-norm}
Every $H\in\mathcal W_N(\phi,C_0)$ satisfies $\|H\|\le3$ with overwhelming probability.
\end{lemma}
\begin{proof}
Write $b_N=N^{-\phi_0}$ and set $A_N=2b_N$, $\delta_0=1/10$, $C_d=\max\{C_0,1\}$. For each independent unordered block construct a real centered variable $Y_{ij}$ with variance
\[
 s_{ij}=(1+\delta_0)/N\quad(i<j),\qquad
 s_{ii}=(1+\delta_0)C_d/N,
\]
by putting
\[
 Y_{ij}=\begin{cases}
 A_N,&\text{with probability }s_{ij}/(A_N^2+s_{ij}),\\
 -s_{ij}/A_N,&\text{with probability }A_N^2/(A_N^2+s_{ij}).
 \end{cases}
\]
Take these variables independently and make $Y$ real symmetric. For $r=s_{ij}/A_N^2=o(1)$, the exact moments are
\begin{equation}
 \mathbb EY_{ij}^k
 =s_{ij}A_N^{k-2}
       \frac{1+(-1)^k r^{k-1}}{1+r},\qquad k\ge2.
 \label{W:eq:complex-dominator-moments}
\end{equation}
For $r\le1/4$ the last fraction is at least $1-r\ge3/4$ when $k\ge3$; for $k=2$ it is one. Hence, uniformly for all integers $k\ge2$,
\[
 \mathbb EY_{ij}^k\ge
 \begin{cases}N^{-1}b_N^{k-2},&i<j,\\
 (C_0/N)b_N^{k-2},&i=j.
 \end{cases}
\]
All these moments are nonnegative. On the other hand, for any $u+v=k\ge2$, $|\mathbb E H_{ij}^u\overline{H_{ij}}^{\,v}|
\le\mathbb E|H_{ij}|^k\le N^{-1}b_N^{k-2}$ off the diagonal; on the diagonal the bound is $\mathbb E|H_{ii}|^k\le
C_0N^{-1}b_N^{k-2}$.

Expand $\operatorname{Tr}H^{2k}$ over closed walks $\mathbf i=(i_0,\ldots,i_{2k})$ with $i_{2k}=i_0$. For $a<b$, let $u_{ab},v_{ab}$ count traversals from $a$ to $b$ and from $b$
to $a$, respectively, and let $m_{aa}$ count loops at $a$. Independence of the unordered entry blocks gives
\begin{equation}
 \mathbb E\prod_{\ell=0}^{2k-1}H_{i_\ell i_{\ell+1}}
 =\prod_{a<b}\mathbb E[H_{ab}^{u_{ab}}
                         \overline{H_{ab}}^{\,v_{ab}}]
   \prod_a\mathbb E H_{aa}^{m_{aa}}.
 \label{W:eq:moment-walk-block-factorization}
\end{equation}
A factor of degree zero is one, and a factor of degree one vanishes. For every remaining factor, the absolute value of the moment is bounded by the nonnegative moment of $Y$ with
the same total degree. Thus, walk by walk,
\[
 \left|\mathbb E\prod_{\ell=0}^{2k-1}H_{i_\ell i_{\ell+1}}\right|
 \le\mathbb E\prod_{\ell=0}^{2k-1}Y_{i_\ell i_{\ell+1}}.
\]
Summing proves
\begin{equation}
 0\le\mathbb E\operatorname{Tr}H^{2k}
 \le\mathbb E\operatorname{Tr}Y^{2k}.
 \label{W:eq:complex-trace-domination}
\end{equation}
\begingroup
\setlength{\intextsep}{1.2\baselineskip}
\begin{figure}[H]
\centering
\begin{tikzpicture}[x=1cm,y=1cm,every node/.style={mg text}]
 \node[anchor=west] at (0,1.25) {(a) One unordered block, $a<b$};
 \node[mg site,label=left:{$a$}] (a) at (.5,0) {};
 \node[mg site,label=right:{$b$}] (b) at (4.8,0) {};
 \draw[mg arrow] (a) to[bend left=25] node[mg label,above] {$u\text{ traversals}:\ H_{ab}^{u}$} (b);
 \draw[mg arrow] (b) to[bend left=25] node[mg label,below] {$v\text{ traversals}:\ \overline{H_{ab}}^{\,v}$} (a);
 \node at (2.65,-1.65) {$|\mathbb EH_{ab}^{u}\overline{H_{ab}}^{\,v}|$};
 \node at (2.65,-2.2) {$\le\mathbb EY_{ab}^{u+v}\quad(u+v\ge2)$};
 \node[anchor=west] at (7,1.25) {(b) An allowed odd multiplicity pattern};
 \node[mg site,label=below:{$a$}] (aa) at (8.2,-1.1) {};
 \node[mg site,label=below:{$b$}] (bb) at (12.8,-1.1) {};
 \node[mg site,label=above:{$c$}] (cc) at (10.5,.7) {};
 \draw[mg arrow] (aa)--node[mg label,below] {$H_{ab}^{3}$}(bb);
 \draw[mg arrow] (bb)--node[mg label,right] {$H_{bc}^{3}$}(cc);
 \draw[mg arrow] (cc)--node[mg label,left] {$\overline{H_{ac}}^{\,3}$}(aa);
 \draw[mg arrow] (aa) .. controls +(140:1.05) and +(220:1.05) .. node[mg label,left] {$H_{aa}^{3}$}(aa);
 \node at (10.5,-2.15) {$a<b<c,\qquad 2k=12$};
\end{tikzpicture}
\caption{Grouping a closed walk by its unordered matrix entries in \eqref{W:eq:moment-walk-block-factorization}. Open circles denote visited indices, and arrows specify the direction of traversal. In (a), the two orientations contribute $H_{ab}^{u}\overline{H_{ab}}^{\,v}$ from the same entry. In (b), the walk follows $a\to b\to c\to a$ three times and then takes three loops at $a$. Its expectation is a product of four third moments, which may be nonzero for centered entries. The comparison in Lemma~\ref{W:lem:complex-norm} includes this index pattern.}
\label{W:fig:moment-closed-walk-blocks}
\end{figure}
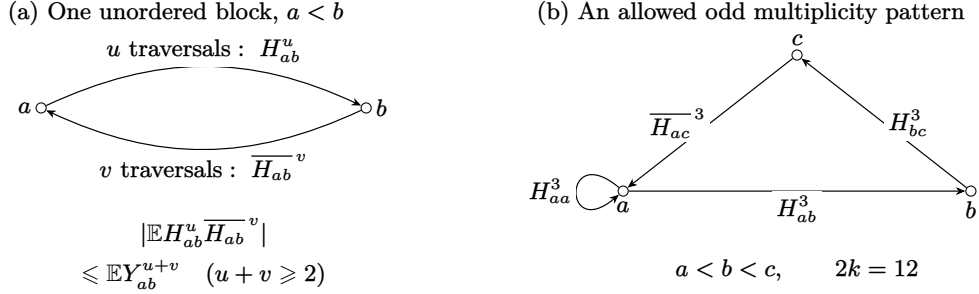
\endgroup

The real matrix $\widetilde Y=Y/\sqrt{1+\delta_0}$ has centered independent entries, off-diagonal variance $N^{-1}$, diagonal variance $C_d/N$, and support at most
$2b_N/\sqrt{1+\delta_0}$. Choose $N_0$ such that $2b_N/\sqrt{1+\delta_0}\le N^{-\phi_0/2}$ for $N\ge N_0$. Apply Lemma~3.1(2) of \cite{lee2014necessary}, with $q=N^{\phi_0/2}$ and
the diagonal variance bound in Definition~2.1 there. It gives
\[
 \mathbb P(\|Y\|>5/2)
 \le N^C\exp(-c\varphi_N^3),\qquad
 \varphi_N=(\log N)^{\log\log N},
\]
because $(2+o(1))\sqrt{1+\delta_0}<5/2$. Deterministically $\|Y\|\le NA_N$. Consequently
\[
 \mathbb E\operatorname{Tr}Y^{2k}
 \le N(5/2)^{2k}
 +N(NA_N)^{2k}N^C\exp(-c\varphi_N^3).
\]
Take $k=\lceil C_D\log N\rceil$ with $C_D\ge(D+2)/(2\log(6/5))$. Then
\[
 N(5/6)^{2k}\le N^{-D-1}.
\]
After division by $3^{2k}$, the logarithm of the second term is at most $C_D'(\log N)^2-c\varphi_N^3$. Since
\[
 \frac{\varphi_N^3}{(\log N)^2}
 =\exp\{3(\log\log N)^2-2\log\log N\}\longrightarrow\infty,
\]
that term is smaller than every fixed inverse power of $N$. Now \eqref{W:eq:complex-trace-domination} and Markov's inequality prove $\mathbb P(\|H\|>3)\le N^{-D}$.
\end{proof}

\Needspace{6\baselineskip}
\subsubsection{From averaged resolvents to counting functions}
\label{W:sec:eigenvalue-counting-and-minor-estimates}

The averaged local law controls a smoothed empirical measure, whereas the minor estimates require the number of eigenvalues in an interval. A deterministic smoothing argument
connects the two. It is formulated to allow atoms, since the empirical measure is discrete, and gives the uniform counting estimate used below.

\begin{lemma}
\label{W:lem:complex-smoothing}
Let $\mu$ be a probability measure supported in $[-3,3]$, let $\nu=\mu-\rho_{\rm sc}(x)\,\dd x$, and put $f(z)=\int(z-x)^{-1}\nu(\dd x)$. Suppose, for $0<\eta_0<1/8$ and $L\ge1$,
\[
 |f(E+\mathrm i\eta)|\le L/(N\eta),\qquad
 |E|\le4,\quad\eta_0\le\eta\le1.
\]
Then
\begin{equation}
 \sup_x|F_\mu(x)-F_{\rm sc}(x)|
 \le C\left[\eta_0+\frac L N(1+|\log\eta_0|)\right].
 \label{W:eq:complex-smoothing}
\end{equation}
Atoms of $\mu$ are allowed.
\end{lemma}
\begin{proof}
Let $D=F_\mu-F_{\rm sc}$ and $P_\eta(t)=\eta/[\pi(t^2+\eta^2)]$. Put $D_\eta=P_\eta*D$, the difference of the two smoothed distribution functions. Both measures are supported in
$[-3,3]$, so $D=0$ off that interval and $|D|\le1$. Consequently
\[
 |D_{\eta_0}(-4)|
 \le\frac{\eta_0}{\pi}\int_{-3}^3
                \frac{\dd x}{(4+x)^2+\eta_0^2}
 \le\frac{6\eta_0}{\pi}.
\]  Moreover
\[
 D_{\eta_0}'(E)=-\pi^{-1}\operatorname{Im}f(E+\mathrm i\eta_0).
\]
Analyticity of $f$ on the rectangle with vertices $-4+\mathrm i\eta_0$, $E+\mathrm i\eta_0$, $E+\mathrm i$, $-4+\mathrm i$ gives
\[
 \left|\int_{-4}^E f(x+\mathrm i\eta_0)\,\dd x\right|
 \le\left|\int_{-4}^E f(x+\mathrm i)\,\dd x\right|
 +\int_{\eta_0}^1\bigl(|f(-4+\mathrm iy)|+|f(E+\mathrm iy)|\bigr)\,\dd y.
\]
The horizontal integral has length at most eight and modulus at most $8L/N$. The two vertical integrals have total modulus at most
\[
 \frac{2L}N\int_{\eta_0}^1\frac{\dd y}{y}
                  =\frac{2L}N|\log\eta_0|.
\]
Integrating $D_{\eta_0}'$ from $-4$ and adding its initial value therefore proves
\begin{equation}
 \sup_{|E|\le4}|D_{\eta_0}(E)|
 \le C\eta_0+C L(1+|\log\eta_0|)/N.
 \label{W:eq:complex-smoothed-cdf}
\end{equation}

We now pass from the smoothed distribution function to $F_\mu$. Let $\Delta=\sup|D|$ and $M_0=\|\rho_{\rm sc}\|_\infty$. If $D(E)$ is within $\varepsilon$ of a positive supremum
$\Delta$, monotonicity of $F_\mu$ implies $D(x)\ge\Delta-\varepsilon-8M_0\eta_0$ on $[E,E+8\eta_0]$. Everywhere else $D\ge-\Delta$. The $P_{\eta_0}$-mass of that interval as seen
from $E+4\eta_0$ is $a=2\arctan(4)/\pi>1/2$. Thus
\[
 D_{\eta_0}(E+4\eta_0)
 \ge(2a-1)\Delta-a\varepsilon-8aM_0\eta_0.
\]
If $D(E)\le-\Delta+\varepsilon$, monotonicity gives $D(x)\le-\Delta+\varepsilon+8M_0\eta_0$ for $E-8\eta_0\le x\le E$. Using $D\le\Delta$ outside this interval gives
\[
 D_{\eta_0}(E-4\eta_0)
 \le-(2a-1)\Delta+a\varepsilon+8aM_0\eta_0.
\]
Since $D$ vanishes outside $[-3,3]$, all convolution centers used here lie in $[-4,4]$. Letting $\varepsilon\downarrow0$ proves $\Delta\le C[\sup_{|E|\le4}|D_{\eta_0}(E)|+\eta_0]$.
Together with \eqref{W:eq:complex-smoothed-cdf}, this proves the lemma.
\end{proof}

\begin{proof}[Completion of Theorem~\ref{W:thm:complex-main}]
Fix $\epsilon,D>0$ and take $\eta_0=N^{-1+\epsilon/4}$. On the norm event of Lemma~\ref{W:lem:complex-norm}, the spectral measure is supported in $[-3,3]$. The averaged law
\eqref{W:eq:complex-averaged}, with its domain parameter fixed at $\epsilon/4$ (or a smaller fixed number if needed), gives the hypothesis of Lemma~\ref{W:lem:complex-smoothing}
with $L=N^{\epsilon/4}$ and failure probability smaller than any prescribed inverse power of the matrix size. Then
\[
 \sup_x|F_H(x)-F_{\rm sc}(x)|
 \le C N^{-1+\epsilon/4}(1+\log N)
 \le N^{-1+\epsilon}
\]
for large $N$. This proves \eqref{W:eq:complex-counting}; large $\epsilon$ follow by using a smaller exponent first.

We now choose a deterministic sequence satisfying the hypotheses preceding \eqref{W:eq:minor-counting}. Choose
increasing deterministic integers $N_k\ge2^k$ so large that, uniformly over $\mathcal W_N(\phi,C_0)$, for every $N\ge
N_k$,
\[
 \mathbb P\left(\|H\|>3\ \text{or}\
 N\sup_x|F_H(x)-F_{\rm sc}(x)|>N^{1/k}\right)\le N^{-k}.
\]
Their existence follows from the two uniform domination statements, using an extra unit in the failure exponent before summing probabilities. Set $k_N=\max\{k:N_k\le N\}$ and
$\mathfrak c_N=N^{1/k_N}$ for $N\ge N_1$, with any positive choice for smaller $N$. For $\delta,D>0$, choose $k_0>\max\{\delta^{-1},D\}$. For $N\ge N_{k_0}$ one has $k_N\ge k_0$,
whence
\[
 \mathfrak c_N=N^{1/k_N}\le N^\delta,
 \qquad N^{-k_N}\le N^{-D}.
\]
Neither threshold depends on a selected matrix law, since $N_{k_0}$ was chosen for the supremum over the class.
\end{proof}

\Needspace{6\baselineskip}
\begin{corollary}
\label{W:cor:complex-counting-for-minors}
The counting estimate \eqref{W:eq:minor-counting} and the trace bound \eqref{W:eq:minor-core-trace} hold for $\beta=2$ under the stated joint entry hypotheses. For a principal
minor with $s$ indices deleted, with trace and counting function normalized by $N$, the corresponding bounds hold with $\mathfrak c_N+s$ in place of $\mathfrak c_N$.

\end{corollary}
\begin{proof}
Apply Theorem~\ref{W:thm:complex-main} under \eqref{W:eq:complex-assumptions}. On the event where its eigenvalue counting estimate holds,
\[
 \sup_x|F_W(x)-F_{\rm sc}(x)|\le\frac{\mathfrak c_N}{N}.
\]
The deterministic integration formula in the proof of Proposition~\ref{W:prop:minor-real-input}, applied to the empirical measure and $\mu_{\rm sc}$, gives
$|m_W(z)-g(z)|\le\pi\mathfrak c_N/(N|\operatorname{Im}z|)$ for every nonreal $z$.

For a principal minor with $s$ indices deleted, interlacing gives $n_W(x)-s\le n_{W^{(I)}}(x)\le n_W(x)$. Normalize both counting functions by $N$. Then
\[
 \sup_x|F_{W^{(I)}}(x)-F_{\rm sc}(x)|
 \le\frac{\mathfrak c_N+s}{N},\qquad
 |m_W^{(I)}(z)-g(z)|\le\frac{\pi(\mathfrak c_N+s)}{N|\operatorname{Im}z|}.
\]
The empirical measure of the minor has mass $(N-s)/N$; the same integration formula applies to these finite measures with the common denominator $N$.

The coarsely truncated and covariance normalized array of Section~\ref{W:sec:minor-estimates-by-coarse-truncation} has independent centered complex blocks and satisfies
\[
 \begin{gathered}
 \mathbb E|W_{ij}|^2=N^{-1},\quad \mathbb EW_{ij}^2=0,\quad
 \mathbb E|W_{ij}|^3\le C_0N^{-3/2},\\
 \mathbb E|W_{ij}|^4\le C_0N^{-2}\log(eN),\quad
 \mathbb EW_{ii}^2\le C_0/N,\quad |W_{ij}|\le C_0N^{-\phi}.
 \end{gathered}
\]
Thus the counting theorem applies with constants uniform over this family. For at most $N^a$ arrays, take failure exponent $D+a+1$ in Theorem~\ref{W:thm:complex-main}. Its
uniformity gives
\[
 \mathbb P(\text{some array fails the counting bound})
 \le N^aN^{-D-a-1}\le N^{-D}.
\]
The trace bounds hold on the same event.
\end{proof}

These estimates concern the empirical spectral measure. To control the spectral measure at $w$, we will also need bounds for resolvent matrix elements at that vector; the following
truncation argument supplies the minor estimates used in their Schur complement representation.

\subsubsection{Coarse truncation and estimates for minors}
\label{W:sec:minor-estimates-by-coarse-truncation}

The counting estimate is now available in both symmetry classes, but its support assumption is stronger than that of the truncated comparison family. A further, coarser truncation
produces a normalized matrix $W$ to which the estimate applies. The discarded entries are controlled by rank in the empirical counting function; the changes of means and
covariances are controlled by a small norm correction. The resulting averaged bound can then be used in the Schur complement equations for individual resolvent entries.

For a fixed $0<\phi<1/100$, the auxiliary threshold is

\begin{equation}
 T_N=N^{1/2-\phi}.
 \label{W:eq:minor-coarse-threshold}
\end{equation}
For large $N$, $N^{7/24}\ll T_N<L_N=e_N\sqrt N$. The auxiliary variables below give a decomposition of $H$ whose three parts can be estimated separately.

\Needspace{6\baselineskip} A coarse truncation decomposes the matrix into a normalized bounded support part, a small norm correction, and a low-rank remainder.

\begin{lemma}
\label{W:lem:minor-decomposition}
Under \eqref{W:eq:minor-family}, one can couple $H$ with
\begin{equation}
 H=W+V+E
 \label{W:eq:minor-decomposition}
\end{equation}
so that $W$ satisfies the hypotheses preceding \eqref{W:eq:minor-counting}, including the full real covariance normalization in the Hermitian case. For every $D>0$, with
probability at least $1-N^{-D}$,
\begin{equation}
 \operatorname{rank}E\le r_N:=\lceil3N^{5\phi}\rceil,
 \qquad \|V\|\le v_N:=N^{-1+4\phi},
 \qquad \|E\|\le C_D(e_N+N^{-\phi})=o(1).
 \label{W:eq:minor-rank-norm}
\end{equation}

\end{lemma}
\begin{proof}
We carry out the construction separately for each independent upper triangular block. Let $T_N=N^{1/2-\phi}$, $p_{ij}=\mathbb P(|\xi_{ij}|>T_N)$ off the diagonal, and use the same
threshold for the standardized diagonal variable. Independently across blocks, sample a Bernoulli label $L$ with parameter $p_{ij}$, a variable $Y$ with the conditional law given
$|\xi|\le T_N$, and a variable $Z$ with the conditional law given $|\xi|>T_N$. Choose $Y,Z$ independently of $L$ and of each other; when a conditional event has probability zero
its associated variable can be arbitrary. The variable $(1-L)Y+LZ$ has exactly the original entry law, since for every bounded measurable $f$,
\[
 \mathbb Ef((1-L)Y+LZ)=(1-p_{ij})\mathbb Ef(Y)+p_{ij}\mathbb Ef(Z)=\mathbb Ef(\xi).
\]
The construction preserves the joint real coordinate law of a complex entry. Imposing the appropriate Hermitian symmetry produces $H$. Put $B_{ij}=N^{-1/2}Y_{ij}$ and $E=H-B$; thus
$E_{ij}=0$ if the corresponding label is zero.

The fourth tail bound in \eqref{W:eq:minor-family} implies $p_{ij}\le CT_N^{-4}=CN^{-2+4\phi}$ for $i<j$. For the diagonal, the second moment gives
$p_{ii}\le CT_N^{-2}=CN^{-1+2\phi}$. If $M_o=\sum_{i<j}L_{ij}$ and $M_d=\sum_iL_{ii}$, independence yields
\[
 \mathbb EM_o\le CN^{4\phi},\qquad \mathbb EM_d\le CN^{2\phi}.
\]
For a sum $M$ of independent Bernoulli variables with mean $\mu$, the exponential Markov bound is
\[
 \mathbb P(M\ge t)\le\inf_{s>0}\exp\{-st+\mu(e^s-1)\}
 \le(e\mu/t)^t,\qquad t\ge\mu.
\]
Take $t=\lfloor N^{5\phi}\rfloor$. For every $N\ge N_0(\phi)$, the estimates $t/\mathbb EM_o\ge cN^\phi$ and $t/\mathbb EM_d\ge cN^{3\phi}$ show that each of the two count bounds
fails with probability at most $\exp(-c_\phi N^{5\phi}\log N)$. On the complementary event,
\[
 \operatorname{rank}E\le2M_o+M_d\le3N^{5\phi}\le r_N.
\]
The rank inequality uses only the fact that a matrix supported on one unordered off-diagonal pair has rank at most two and a perturbation of one diagonal entry has rank at most
one.

We next center the variables conditioned to lie below the truncation threshold and normalize the covariances of the off-diagonal variables. For $i<j$, let $X$ denote the real
standardized coordinate vector, either the scalar $\xi$ or $(\sqrt2\operatorname{Re}\xi,\sqrt2\operatorname{Im}\xi)$. Put $\mu=\mathbb E(Y)$ and $\Sigma=\operatorname{Cov}(Y)$. The
tail integrals \eqref{W:eq:quant-tail-below-four}, applied to $q=1,2$, give
\[
 \bigl\|\mathbb E[X\mathbf1_{\{|\xi|>T_N\}}]\bigr\|\le CT_N^{-3},
 \qquad \bigl\|\mathbb E[XX^{\mathsf T}\mathbf1_{\{|\xi|>T_N\}}]\bigr\|\le CT_N^{-2}.
\]
Since $p_{ij}=O(T_N^{-4})$, the exact conditional identities
\[
 \mu=\frac{\mathbb EX-\mathbb E[X\mathbf1_{\{|\xi|>T_N\}}]}{1-p_{ij}},
 \qquad
 \Sigma+\mu\mu^{\mathsf T}=
 \frac{\mathbb E[XX^{\mathsf T}]-\mathbb E[XX^{\mathsf T}\mathbf1_{\{|\xi|>T_N\}}]}{1-p_{ij}}
\]
and \eqref{W:eq:minor-family} yield
\begin{equation}
 \|\mu\|\le CN^{-3/2+3\phi},\qquad
 \|\Sigma-I\|\le CN^{-1+2\phi}.
 \label{W:eq:minor-coarse-moments}
\end{equation}
Indeed, $T_N^{-3}=N^{-3/2+3\phi}$, $T_N^{-2}=N^{-1+2\phi}$, while the specified original mean and covariance errors are of smaller order for fixed $\phi>0$. The contribution
$\mu\mu^{\mathsf T}$ is $O(T_N^{-6})$.

For large $N$ the eigenvalues of $\Sigma$ lie in $[1/2,3/2]$. Define $Y^{\mathrm{nor}}=\Sigma^{-1/2}(Y-\mu)$ for each off-diagonal block, convert back to a real or complex matrix
entry, and divide by $\sqrt N$ to obtain $W_{ij}$. Set $W_{ii}=N^{-1/2}(Y_{ii}-\mathbb EY_{ii})$ on the diagonal. Then the off-diagonal blocks are centered and satisfy exactly
\[
 \mathbb E W_{ij}\overline{W_{ij}}=N^{-1},
 \qquad \mathbb EW_{ij}^2=0\quad\text{in the complex case}.
\]
The support of $Y$ is at most $T_N$ and $\|\Sigma^{-1/2}\|\le\sqrt2$, so $|W_{ij}|\le CN^{-\phi}$. The conditional fourth moment is bounded by twice the original truncated fourth
moment. Consequently, after the bounded covariance transformation and subtraction of the small mean,
\[
 \mathbb E|W_{ij}|^3\le CN^{-3/2},\qquad
 \mathbb E|W_{ij}|^4\le CN^{-2}\log(eN),\qquad i<j.
\]
The diagonal conditional second moment is bounded by $2C/N$. These verify all hypotheses of the eigenvalue counting estimate for matrices with bounded entries, with one fixed
choice of their constants. The variables $Y$ are independent of the labels. Since $W$ is a function of $Y$, it is independent of the graph determined by those labels; $E$ also
depends on $Y$ through the removed small entries.

It remains to control $V=B-W$. Put $\delta=N^{-1+2\phi}$. The matrix $V$ is the sum of its deterministic mean matrix $V_0$ and independent centered matrices $X_{ij}$ supported on
unordered off-diagonal pairs. By \eqref{W:eq:minor-coarse-moments} and the diagonal second moment estimate,
\[
 \|V_0\|\le\max_i\sum_j|(V_0)_{ij}|\le C N^{-1+3\phi}.
\]
The second inequality includes the diagonal means: $|\mathbb EY_{ii}|\le C(N^{-1/2+o(1)}+T_N^{-1})\le CN^{-1/2+\phi}$ for sufficiently large $N$. For off-diagonal pairs,
$\|\Sigma^{-1/2}-I\|\le C\delta$ and $N^{-1/2}|Y-\mu|\le CN^{-\phi}$ show
\[
 \|X_{ij}\|\le C\delta N^{-\phi},\qquad
 \mathbb E| (X_{ij})_{ij}|^2\le C\delta^2/N.
\]
Since $X_{ij}^2=|(X_{ij})_{ij}|^2(E_{ii}+E_{jj})$, the matrix variance satisfies
\[
 \left\|\sum_{i<j}\mathbb EX_{ij}^2\right\|
 =\max_i\sum_{j\ne i}\mathbb E|(X_{ij})_{ij}|^2\le C\delta^2.
\]
The Hermitian matrix Bernstein inequality now gives, for $t>0$,
\[
 \mathbb P\!\left(\left\|\sum_{i<j}X_{ij}\right\|>t\right)
 \le2N\exp\!\left[-\frac{t^2}{C\delta^2+C\delta N^{-\phi}t}\right].
\]
For $t=\delta\log N$ the exponent is at most $-c(\log N)^2$, because $N^{-\phi}\log N=o(1)$. Thus the failure is smaller than $N^{-D}$ for each fixed $D$, and on the remaining
event
\[
 \|V\|\le CN^{-1+3\phi}+N^{-1+2\phi}\log N
 =N^{-1+4\phi}\bigl(CN^{-\phi}+N^{-2\phi}\log N\bigr)\le v_N
\]
for all sufficiently large $N$.

Finally let $\mathcal G$ be the random graph of marked off-diagonal pairs. For a fixed positive integer $d$, the occurrence of a component with at least $d+1$ vertices implies the
occurrence of a labelled tree with $d$ edges on some set of $d+1$ vertices. Independence of the labels and the bound on $p_{ij}$ imply
\[
 \mathbb P(\mathcal G\text{ has a component of size }>d)
 \le\binom N{d+1}(d+1)^{d-1}(CN^{-2+4\phi})^d
 \le C_dN^{1-d(1-4\phi)}.
\]
Choose $d$ after the prescribed failure exponent $D$, so that the last quantity is $O(N^{-D-1})$. On the event that every component has at most $d$ vertices, $E$ is block diagonal
on these components, apart from marked diagonal entries, which lie inside the components containing their vertices. Moreover $|E_{ij}|\le|H_{ij}|+|B_{ij}|\le e_N+CN^{-\phi}$, so
every row has at most $d$ nonzero entries and
\[
 \|E\|\le\max_i\sum_j|E_{ij}|\le d(e_N+CN^{-\phi})\le C_D(e_N+N^{-\phi}).
\]
Intersecting this graph event with the two count events and the Bernstein event proves all three estimates in \eqref{W:eq:minor-rank-norm}; their combined failure is at most
$N^{-D}$ after choosing the intermediate failure exponents larger than $D$. The constants are uniform over the independent entry laws satisfying \eqref{W:eq:minor-family}.
\end{proof}

\Needspace{6\baselineskip} Rank inequalities control the averaged effect of the low-rank error. The resolvent identity bounds the small norm error with the spectral density
appearing in the estimate below.

\begin{lemma}
\label{W:lem:minor-rank-resolvent}
For two Hermitian $N\times N$ matrices $A,B$ with $\operatorname{rank}(A-B)\le r$,
\begin{equation}
 |m_A(z)-m_B(z)|\le\frac{\pi r}{N\eta}.
 \label{W:eq:minor-rank-trace}
\end{equation}
For a principal minor with $s$ deleted coordinates,
\begin{equation}
 |m_H^{(J)}(z)-m_H(z)|\le\frac{\pi s}{N\eta}.
 \label{W:eq:minor-minor-trace}
\end{equation}
Finally, if $B=W+V$ and $\|V\|\le\eta/2$, then
\begin{equation}
 |m_B(z)-m_W(z)|
 \le\frac{2\|V\|}{\eta}\upsilon_W(z).
 \label{W:eq:minor-small-norm-trace}
\end{equation}
\end{lemma}
\begin{proof}

Write $n_A(x)=\dim\operatorname{Ran}\mathbf1_{(-\infty,x]}(A)$. The intersection
\[
 V_x=\operatorname{Ran}\mathbf1_{(-\infty,x]}(A)\cap\ker(A-B)
\]
has dimension at least $n_A(x)-r$. On this space $\langle Bv,v\rangle=\langle Av,v\rangle\le x\|v\|^2$. The min--max principle therefore gives $n_B(x)\ge n_A(x)-r$. Interchanging
$A$ and $B$ proves $|F_A(x)-F_B(x)|\le r/N$. For a compression to a subspace of codimension $s$, min--max gives $n_H(x)-s\le n_{H^{(J)}}(x)\le n_H(x)$. Normalize both counting
functions by $N$. The integration formula in the proof of Proposition~\ref{W:prop:minor-real-input}, with the two empirical measures, gives
\[
 |m_A(z)-m_B(z)|\le\frac{\pi r}{N|\operatorname{Im}z|}.
\]
For $H$ and $H^{(J)}$, use $s$ in place of $r$. The measures then have masses $1$ and $(N-s)/N$; boundedness of their distribution functions gives the vanishing boundary terms
required by that formula.

The identity $R_B=(I-R_WV)^{-1}R_W$ and $\|R_WV\|\le1/2$ give $\|R_B\|_{\rm HS}\le2\|R_W\|_{\rm HS}$. By the resolvent identity and the Hilbert--Schmidt trace inequality,
\[
 |m_B-m_W|\le N^{-1}\|V\|\|R_W\|_{\rm HS}\|R_B\|_{\rm HS}
 \le2\|V\|N^{-1}\|R_W\|_{\rm HS}^2.
\]
The Ward identity $N^{-1}\|R_W\|_{\rm HS}^2=\upsilon_W/\eta$ proves the result. For $\eta\ge cN^{-2/3}$, on the event $\upsilon_W\le N^{-1/3+\epsilon}$ one has
\[
 2\|V\|\upsilon_W/\eta\le (2/c)\|V\|N^{1/3+\epsilon}.
\]
The probability estimate here is used near a spectral edge. More precisely, if
\[
 cN^{-2/3}\le\eta\le CN^{-2/3},\qquad
 \min\{|E-2|,|E+2|\}\le CN^{-2/3},
\]
then the semicircle formula gives $\upsilon(z)\le C' N^{-1/3}$. On the event where \eqref{W:eq:minor-core-trace} holds,
\[
 \upsilon_W(z)\le\upsilon(z)+\frac{\pi\mathfrak c_N}{N\eta}
 \le (C'+\pi\mathfrak c_N/c)N^{-1/3}.
\]
Since $\mathfrak c_N=N^{o(1)}$, this is at most $N^{-1/3+\epsilon}$ for all sufficiently large $N$. This event has the stated overwhelming probability. At other spectral parameters
the preceding deterministic estimate uses $\upsilon(z)+\pi\mathfrak c_N/(N\eta)$ without replacing it by $N^{-1/3+\epsilon}$.
\end{proof}

\Needspace{6\baselineskip} Combining the decomposition with the eigenvalue counting estimate for matrices with bounded entries yields an averaged estimate for the array under the
stated tail assumption and for its principal minors obtained by deleting a bounded number of coordinates.

\begin{proposition}
\label{W:prop:minor-averaged}
Assume \eqref{W:eq:minor-family}. Let $\mathfrak c_N$ satisfy \eqref{W:eq:minor-counting}, and set $r_N=\lceil3N^{5\phi}\rceil$ and $v_N=N^{-1+4\phi}$. Fix $s_0$, a polynomial
family size, and $D>0$. With probability at least $1-N^{-D}$, simultaneously for every matrix in that family, every $|J|\le s_0$ and every $z=E+\mathrm i\eta$ with $\eta\ge2v_N$,
\begin{equation}
 \begin{aligned}
 |m_H^{(J)}(z)-g(z)|\le A_N(z;|J|)
 :={}&\frac{\pi(\mathfrak c_N+r_N+|J|)}{N\eta}\\
 &+\frac{2v_N}{\eta}\left(\upsilon(z)+\frac{\pi\mathfrak c_N}{N\eta}\right).
 \end{aligned}
 \label{W:eq:minor-averaged-bound}
\end{equation}
On the same event $\|H\|\le4$. The estimate also holds in the lower half-plane by conjugation. The constants are uniform over independent block distributions in the stated family.
\end{proposition}
\begin{proof}
Use $H=(W+V)+E$, where $\operatorname{rank}E\le r_N$ and $\|V\|\le v_N$. Write $\upsilon_W=-\operatorname{Im}m_W$ in the upper half-plane. The trace bound
\eqref{W:eq:minor-core-trace} implies
\[
 |m_W-g|\le\frac{\pi\mathfrak c_N}{N\eta},\qquad
 \upsilon_W\le\upsilon+\frac{\pi\mathfrak c_N}{N\eta}.
\]
For $\eta\ge2v_N$, the two perturbation estimates give, separately,
\[
 |m_H-m_{W+V}|\le\frac{\pi r_N}{N\eta},\qquad
 |m_{W+V}-m_W|\le\frac{2v_N}{\eta}\upsilon_W.
\]
Finally $|m_H^{(J)}-m_H|\le\pi s/(N\eta)$ by \eqref{W:eq:minor-minor-trace}. Adding the four differences gives
\begin{align*}
 |m_H^{(J)}-g|
 &\le |m_H^{(J)}-m_H|+|m_H-m_{W+V}|
          +|m_{W+V}-m_W|+|m_W-g|\\
 &\le\frac{\pi(s+r_N+\mathfrak c_N)}{N\eta}
       +\frac{2v_N}{\eta}
          \left(\upsilon+\frac{\pi\mathfrak c_N}{N\eta}\right).
\end{align*}
The lower half-plane follows by adjunction. All four comparisons hold on the same matrix decomposition and the same event where the eigenvalue counting estimate holds. Their bounds
are deterministic for each $z$ with $\eta\ge2v_N$; in particular $\|R_WV\|\le v_N/\eta\le1/2$. The right side is therefore the explicit function $A_N(z;s)$ in
\eqref{W:eq:minor-averaged-bound}, with no additional probability event depending on $z$ or on the deleted set $J$. Also $\|H\|\le\|W\|+\|V\|+\|E\|\le3+o(1)<4$ on the same
decomposition event. For a family of at most $N^b$ intermediate matrices, apply each probability estimate with exponent $D+b+2$. The union of its failures has probability at most
$C N^bN^{-D-b-2}\le N^{-D}$ for $N^2\ge C$.
\end{proof}

Fix $0<\tau<1/100$ and choose
\begin{equation}
 \kappa=\phi=\tau/100,\qquad
 \eta_{\mathrm c}=N^{-2/3+\kappa},\qquad \eta_*=N^{-2/3-\kappa}.
 \label{W:eq:minor-parameters}
\end{equation}
For a fixed constant $M>10$ let
\[
 \mathcal V_N=\{E+\mathrm i\eta:|E-2|\le \eta_{\mathrm c},
                             \eta_*\le\eta\le M\},
 \qquad
 \mathcal A_N=\{z\in\mathcal V_N:\eta\le \eta_{\mathrm c}\}.
\]
The negative edge follows by replacing $H$ with $-H$, and the lower half-plane by adjunction. Since $\upsilon(E+\mathrm i\eta)\le C\sqrt{|E-2|+\eta}$,
\eqref{W:eq:minor-averaged-bound} gives, for all sufficiently large $N$, uniformly over the stated family,
\begin{equation}
 \sup_{z\in\mathcal V_N}|m_H^{(J)}(z)-g(z)|=o(1),\qquad
 \sup_{z\in\mathcal A_N}\upsilon_H^{(J)}(z)
       \le N^{-1/3+8\kappa}.
 \label{W:eq:minor-averaged-critical}
\end{equation}
Take $\mathfrak c_N\le N^\kappa$. The first term of \eqref{W:eq:minor-averaged-bound} is at most $CN^{-1/3+6\kappa}$ at $\eta_*$. On $\mathcal A_N$ the small norm term is at most
$CN^{-2/3+6\kappa}+CN^{-2/3+7\kappa}$, since $v_N/\eta_* =N^{-1/3+5\kappa}$. In particular $2v_N<\eta_*$ for all sufficiently large $N$, as required when applying the small norm
comparison throughout the domain, including its lowest height. The choices $\eta_*<\eta_{\mathrm c}$ have different roles: the former bounds resolvent denominators, while the
latter bounds $|E-2|$ and the semicircle density on $\mathcal A_N$. Constants are absorbed by one further $N^\kappa$. On the rest of $\mathcal V_N$ both terms still tend to zero.
This proves \eqref{W:eq:minor-averaged-critical} from the preceding eigenvalue counting estimate.

The conditional estimates for the Schur complement will be used for matrices whose standardized entries have small, but not necessarily vanishing, means and covariance errors. In
\eqref{W:eq:minor-family} these errors are bounded by $\varepsilon_{1,N}\le N^{-3/2+o(1)}$ and $\varepsilon_{2,N}\le N^{-1+o(1)}$. Since $\ell_N^{\mathrm{Sch}}=N^{7/24}$ and the
fourth moment tail quantity satisfies $a_N\le N^{o(1)}$, the exact conditional mean and covariance formulas \eqref{W:eq:schur-small-defect-identities} give
\[
 \begin{aligned}
 u_N&\le C\bigl(a_N(\ell_N^{\mathrm{Sch}})^{-3}+\varepsilon_{1,N}\bigr)
      \le N^{-7/8+o(1)},\\
 \delta_N&\le C\bigl(a_N(\ell_N^{\mathrm{Sch}})^{-2}+\varepsilon_{2,N}+\varepsilon_{1,N}^2\bigr)
      \le N^{-7/12+o(1)}.
 \end{aligned}
\]
Indeed $3(7/24)=7/8$ and $2(7/24)=7/12$, while $\varepsilon_{1,N}^2\le N^{-3+o(1)}$. In particular the variance estimate uses the noncentral conditional second moment and subtracts
the square of the conditional mean; for a complex entry it also controls the covariance of the two real coordinates through the pseudovariance estimate in
Lemma~\ref{W:lem:schur-small-moments}. These parameters are deterministic, uniform over the compared entry distributions, and satisfy the hypotheses of
Theorem~\ref{W:thm:schur-residual}. The matrix entries still have the conditional support bound $b_N=4N^{-5/24}$. Theorem~\ref{W:thm:schur-residual} therefore applies conditional
on the graph, with each complex entry represented by the joint distribution of its two real coordinates.

\Needspace{6\baselineskip}

\begin{lemma}
\label{W:lem:minor-conditional-event}
Fix a bound on the component size $d$ and a finite covering set of polynomial cardinality in $\mathcal V_N$. On the event that every auxiliary graph component has size at most $d$,
there is an event whose complement has probability smaller than any prescribed inverse power of the matrix size and on which, for every component $I$, its complementary resolvent
${\mathcal R}=R_H^{(I)}(z)$, and all $j\notin I$,
\begin{align}
 \|B{\mathcal R}B^*-m_H^{(I)}I_I\|
   &\le C_{D,d}(\log N)^4\Theta_N({\mathcal R};z),
             \label{W:eq:minor-residual-event}\\
 \|B{\mathcal R} e_j\|_2
   &\le C_{D,d}(\log N)^2
       \left[(N^{-1/2}+u_N)\|{\mathcal R}e_j\|_2
                     +b_N\|{\mathcal R}e_j\|_\infty\right].
             \label{W:eq:minor-linear-event}
\end{align}
Here $B=H_{IJ}$ and $\Theta_N$ is the scale in Section~\ref{W:sec:smooth-cutoff-for-schur-complement}, with the enlarged parameters just stated.

\end{lemma}
\begin{proof}
Fix a realization of the label graph in $\mathcal A_d$ and one of its components $I$. Let $\mathscr F_I$ be the $\sigma$-field in Theorem~\ref{W:thm:schur-residual}. In particular,
the complementary matrix $H_{I^cI^c}$, its resolvent $\mathcal R$ and $m_H^{(I)}=N^{-1}\operatorname{Tr}\mathcal R$ are $\mathscr F_I$ measurable. Given $\mathscr F_I$, every entry
of $B=H_{I,I^c}$ has the conditional law given a zero label, and entries attached to different unordered pairs are independent by Lemma~\ref{W:lem:schur-components}. All estimates
below are conditional on this fixed data.

Theorem~\ref{W:thm:schur-residual}, with the mean and covariance parameters computed above, gives for every even $p\ge4$
\[
 \left\|\|B\mathcal RB^*-m_H^{(I)}I_I\|_{\mathrm{HS}}\right\|_{L^p(\cdot\mid\mathscr F_I)}
 \le C_dp^2\Theta_N(\mathcal R;z).
\]
Since the operator norm is bounded by the Hilbert--Schmidt norm, the same right hand side bounds its conditional $L^p$ operator norm. The scale $\Theta_N$ is nonnegative and
$\mathscr F_I$ measurable, so it can be divided out when positive.

For the linear expression fix $j\notin I$ and put $t=\mathcal Re_j$, a vector determined by $\mathscr F_I$. For $i\in I$, let $m_i$ be the conditional mean row of $B$ and
$x_i=B_{i,I^c}-m_i$. Its entries have mean zero, conditional second moments at most $2/N$, and absolute value at most $b_N$. Each conditional mean satisfies
$|(m_i)_k|\le N^{-1/2}u_N$, hence
\[
 \|m_i\|_2\le u_N,\qquad |m_i t|\le u_N\|t\|_2.
\]
The complex or real linear moment bound \eqref{W:eq:schur-linear}, applied to this fixed coefficient vector and the centered row, yields
\[
 \|(Bt)_i\|_{L^p(\cdot\mid\mathscr F_I)}
 \le C p\bigl(N^{-1/2}\|t\|_2+b_N\|t\|_\infty\bigr)+u_N\|t\|_2
 \le C'pA_j,
\]
where
\[
 A_j=(N^{-1/2}+u_N)\|\mathcal Re_j\|_2+b_N\|\mathcal Re_j\|_\infty.
\]
For $p\ge2$, applying Minkowski's inequality in conditional $L^{p/2}$ gives
\[
 \begin{aligned}
 \left\|\|B\mathcal Re_j\|_2\right\|_{L^p(\cdot\mid\mathscr F_I)}^2
 &=\left\|\sum_{i\in I}|(Bt)_i|^2\right\|_{L^{p/2}(\cdot\mid\mathscr F_I)}\\
 &\le\sum_{i\in I}\|(Bt)_i\|_{L^p(\cdot\mid\mathscr F_I)}^2
 \le d(C'pA_j)^2.
 \end{aligned}
\]
Thus the conditional $L^p$ norm of the entire vector is at most $C_dp A_j$.

Take the even integer $p_N=2\lceil\log N\rceil$ and put $L=\log N$. For positive scales, conditional Markov's inequality gives
\[
 \begin{aligned}
 &\mathbb P\!\left(\|B\mathcal RB^*-m_H^{(I)}I_I\|>L^4\Theta_N(\mathcal R;z)\,\middle|\,\mathscr F_I\right)
 \le\left(\frac{C_dp_N^2}{L^4}\right)^{p_N},\\
 &\mathbb P\!\left(\|B\mathcal Re_j\|_2>L^2 A_j\,\middle|\,\mathscr F_I\right)
 \le\left(\frac{C_dp_N}{L^2}\right)^{p_N}.
 \end{aligned}
\]
If $\Theta_N(\mathcal R;z)=0$ or $A_j=0$, the corresponding conditional moment bound forces the random norm to be zero almost surely, so the inequality holds without division.
Since $p_N\le2\log N+2$, the bases in the two displayed bounds are at most $C_d(\log N)^{-2}$ and $C_d(\log N)^{-1}$, respectively, for large $N$. Their $p_N$th powers are bounded
above by
\[
 \exp\{-2\log N\log\log N+C_d'\log N\}.
\]
For every prescribed $B_0>0$ this is at most $N^{-B_0}$ once $N$ is sufficiently large, uniformly in the realized graph, all internal entries, the complementary minor, and the
spectral point.

Let the prescribed spectral partition have at most $N^a$ points, and let the matrix family have at most $N^b$ members. On $\mathcal A_d$ there are at most $N$ components for each
matrix and at most $N$ possible indices $j$ per component. Thus there are at most $N^{a+b+2}$ conditional inequalities of the two types to check. The event $\mathcal A_d$ and the
index sets $I$ are measurable with respect to the graph, so conditional expectation and the union bound yield
\[
 \begin{aligned}
 &\mathbb P\!\left(\mathcal A_d\cap\{\text{one of the two preceding bounds fails}\}\right)\\
 &\qquad\le\mathbb E\!\left[\mathbf1_{\mathcal A_d}
 \sum_{I\in\operatorname{Comp}(L)}\sum_{z}\left(\mathbb P(\text{first bound fails}\mid\mathscr F_I)
 +\sum_{j\notin I}\mathbb P(\text{second bound fails}\mid\mathscr F_I)\right)\right]\\
 &\qquad\le2N^{a+b+2-B_0}.
 \end{aligned}
\]
Choose $B_0>D+a+b+3$. This probability is at most $N^{-D}$ for all sufficiently large $N$. The separate failure probability of $\mathcal A_d^c$ is given in
Lemma~\ref{W:lem:schur-components} and is added only when the present bounds are used without that restriction. Restoring fixed constants in the thresholds gives
\eqref{W:eq:minor-residual-event} and \eqref{W:eq:minor-linear-event} on one event, uniformly over the stated components, columns and spectral sampling points.
\end{proof}

\Needspace{6\baselineskip}
\subsubsection{Entrywise estimates and stability under perturbations}
\label{W:sec:entrywise-estimates-and-perturbations}

The averaged estimate controls the trace of the resolvent, while the conditional estimates control the error in its Schur complement equations. Together they give a continuity
argument for individual entries. The argument must be uniform over the comparison family and over all sets of deleted coordinates of bounded cardinality, since these minors
reappear when rows are conditioned upon. The following result includes this uniformity and the stability needed for entry replacements.

\begin{theorem}
\label{W:thm:minor-bootstrap}
Assume \eqref{W:eq:minor-family}, with $\beta\in\{1,2\}$. Fix $0<\tau<1/100$, a bound on the number of deleted indices $s_0$, a polynomial family size, and $D>0$. For all
sufficiently large $N$, with probability at least $1-N^{-D}$, simultaneously for every array in that family, every $|J|\le s_0$, and every $z\in\mathcal A_N$,
\begin{equation}
 \max_{i,j\notin J}|(R_H^{(J)}(z))_{ij}|\le3,
 \qquad
 -\frac1N\operatorname{Im}\operatorname{Tr}R_H^{(J)}(z)
        \le N^{-1/3+\tau}.
 \label{W:eq:minor-final-minors}
\end{equation}
The entrywise bound holds on the full vertical domain $\mathcal V_N$, with a fixed constant in place of $3$ if necessary. In fact
\begin{equation}
 \max_{i,j}|R_{H,ij}(z)-g(z)\delta_{ij}|
       \le C_{D,\tau}\bigl(e_N+N^{-1/10}\bigr)|g(z)|
 \label{W:eq:minor-coarse-entrywise}
\end{equation}
there. The counting estimates for matrices with bounded entries used in the proof are Proposition~\ref{W:prop:minor-real-input} for the real class and
Theorem~\ref{W:thm:complex-main} for the complex class.
\end{theorem}
\begin{proof}
Choose the component size bound $d$ sufficiently large that $C_dN^{1-d/6}$ is smaller than the required final failure probability, including the polynomial family union. Intersect
the averaged event of Proposition~\ref{W:prop:minor-averaged}, the norm event $\|H\|\le4$, and the conditional event of Lemma~\ref{W:lem:minor-conditional-event} on a rectangular
array of spectral parameters of spacing at most $N^{-4}$ in $(E,\eta)$. The number of points is polynomial, so the preceding concentration estimates apply simultaneously after
increasing the probability exponent.

Set $\epsilon_0=[100(d+s_0+1)]^{-1}$ and
\[
 \Lambda(z)=\max_i|R_{ii}(z)/g(z)-1|
       \mathbin{\vee}\max_{i\ne j}|R_{ij}(z)/g(z)|.
\]
Choose the fixed constant $M$ sufficiently large depending on $\epsilon_0$. Since $\|H\|\le4$, the Neumann expansion at $\eta=M$ gives
\[
 \|R_H(z)-z^{-1}I\|\le\frac4{|z|(|z|-4)},\qquad
 |g(z)-z^{-1}|\le C|z|^{-3}.
\]
Since $|g(z)|\ge c/|z|$ there, this proves $\Lambda(z)\le\epsilon_0/4$ at the row of spectral parameters with largest imaginary part.

Suppose now that $\Lambda(z)\le\epsilon_0$ at some selected point. For a component $I$ of size at most $d$, $R_{II}=g(I_I+Q_I)$ with $\|Q_I\|\le d\epsilon_0<1/10$. Thus
$\|R_{II}^{-1}\|\le2/|g|$. The principal minor identity is
\begin{equation}
 {\mathcal R}=R_{J J}-R_{J I}R_{II}^{-1}R_{I J},\qquad J=I^c.
 \label{W:eq:minor-minor-identity}
\end{equation}
For $j,k\notin I$, both vectors in its correction have Euclidean norm at most $\sqrt d\,|g|\epsilon_0$. Consequently
\begin{equation}
 |{\mathcal R}_{jk}-g\delta_{jk}|
    \le |g|(\epsilon_0+2d\epsilon_0^2),\qquad
       \max_{j,k}|{\mathcal R}_{jk}|\le2.
 \label{W:eq:minor-bootstrap-minor}
\end{equation}
The averaged estimate and interlacing control the trace of this minor.

The Ward identities now imply
\[
 N^{-1}\|{\mathcal R}\|_{\rm HS}
       =\sqrt{\upsilon_H^{(I)}/(N\eta)},\quad
 \max_j\|{\mathcal R}e_j\|_2+\max_j\|e_j^*{\mathcal R}\|_2
       \le2\sqrt{2/\eta},\quad
 \|\operatorname{diag}{\mathcal R}\|_2\le2\sqrt N.
\]
Substitute these in the definition of $\Theta_N$. Uniformly on $\mathcal V_N$ its terms are bounded by
\begin{equation}
 \begin{split}
 \Theta_N({\mathcal R};z)\le C\bigg[&
  \sqrt{F_N/N}+\sqrt{\upsilon_H^{(I)}/(N\eta)}
  +\frac{b_N}{\sqrt{N\eta}}+b_N^2\\
  &+u_N(N^{-1/2}+b_N)\eta^{-1}
      +u_N^2\eta^{-1}+\delta_N\bigg].
 \end{split}
 \label{W:eq:minor-bootstrap-theta}
\end{equation}
On $\mathcal A_N$, the bounds above and \eqref{W:eq:minor-averaged-critical} give the respective powers
\begin{equation}
 \begin{gathered}
 N^{-1/2+o(1)},\quad N^{-1/3+(9/2)\kappa},\quad
 N^{-3/8+\kappa/2},\quad N^{-5/12},\\
 N^{-17/24+\kappa+o(1)}+N^{-5/12+\kappa+o(1)},\quad
 N^{-13/12+\kappa+o(1)},\quad N^{-7/12+o(1)}.
 \end{gathered}
 \label{W:eq:minor-theta-powers}
\end{equation}
For $\eta\ge\eta_{\rm c}$, write $c_N=\mathfrak c_N$ and use $\upsilon(z)\le C\sqrt\eta$ and the displayed formula for
$A_N$:
\[
 \frac{\upsilon_H^{(I)}}{N\eta}
 \le\frac{C}{N\sqrt\eta}
 +\frac{\pi(c_N+r_N+d)}{N^2\eta^2}
 +\frac{Cv_N}{N\eta^{3/2}}
 +\frac{2\pi v_Nc_N}{N^2\eta^3}.
\]
Every term on the right decreases as $\eta$ increases. At $\eta_{\rm c}=N^{-2/3+\kappa}$, use $c_N\le N^\kappa$, $r_N\le C N^{5\kappa}$, and $v_N=N^{-1+4\kappa}$. The four terms
are bounded by
\[
 C N^{-2/3-\kappa/2},\qquad
 C N^{-2/3+3\kappa},\qquad
 C N^{-1+5\kappa/2},\qquad
 C N^{-1+2\kappa},
\]
respectively. Their square root is at most $C N^{-1/3+3\kappa/2}$, which is bounded by $N^{-1/3+(9/2)\kappa}$ for $N^{3\kappa}\ge C$. Every other term in
\eqref{W:eq:minor-bootstrap-theta} is constant or decreasing in $\eta$ once the entry and diagonal bounds have been substituted. This proves the required bounds also above
$\eta_{\rm c}$. Thus the right hand side of \eqref{W:eq:minor-residual-event} is $o(1)$ uniformly on $\mathcal V_N$; it is at most $N^{-1/3+10\kappa}$ for all $N\ge N_0(\kappa)$.
Likewise \eqref{W:eq:minor-linear-event} and \eqref{W:eq:minor-bootstrap-minor} give
\begin{equation}
 \|B{\mathcal R} e_j\|_2
 \le N^{o(1)}\left[(N^{-1/2}+u_N)\eta^{-1/2}+b_N\right]
 \le N^{-1/10}.
 \label{W:eq:minor-bootstrap-linear}
\end{equation}
The largest exponent is $-1/6+\kappa/2+o(1)<-1/10$ because $\kappa<10^{-4}$.

Put $Z_I=B{\mathcal R}B^*-m_H^{(I)}I_I$. The Schur complement formula is
\begin{equation}
 R_{II}=\left[g^{-1}I_I-\left(
        H_{II}+(m_H^{(I)}-g)I_I+Z_I\right)\right]^{-1}.
 \label{W:eq:minor-bootstrap-schur}
\end{equation}
The matrix in parentheses has norm at most $de_N+o(1)$. Since $|g|\le1$, the norm of its product with $g$ is below $1/2$ for $N\ge N_0$. The inverse formula then gives
$\|R_{II}\|\le2|g|$ and
\[
 \|R_{II}-gI_I\|\le
 2|g|^2\left(de_N+|m_H^{(I)}-g|+\|Z_I\|\right).
\]
For $j\notin I$, the other block of the inverse satisfies $R_{I j}=R_{II}B{\mathcal R} e_j$, so \eqref{W:eq:minor-bootstrap-linear} bounds its entries by $2|g|N^{-1/10}$. Every
pair of coordinates is either internal to one component or is covered by this exterior formula. We have proved the strict improvement
\begin{equation}
 \Lambda(z)\le C_d e_N+ C N^{-1/10}<\epsilon_0/2
 \label{W:eq:minor-bootstrap-improvement}
\end{equation}
for all sufficiently large $N$ whenever $\Lambda(z)\le\epsilon_0$.

To compare consecutive spectral parameters with decreasing imaginary part, $\|\partial_zR_H\|\le\eta^{-2}$ and the explicit formula for $g$ give a deterministic Lipschitz bound
$CN^2$ for $\Lambda$ on $\mathcal V_N$ (the lower bound for $|g|$ depends only on $M$). Moving one step between neighboring spectral parameters changes it by at most $CN^{-2}$.
Starting from the top row, induction keeps every new point below $\epsilon_0$ and then improves it below $\epsilon_0/2$. Interpolation in both $E$ and $\eta$ gives
\eqref{W:eq:minor-coarse-entrywise} throughout the domain.

Put $\epsilon_N=C_{D,\tau}(e_N+N^{-1/10})$. For any $J$ with $|J|\le s_0$, the proved estimate for the full matrix gives
\[
 R_{JJ}=g(I_J+D_J),\qquad \|D_J\|\le s_0\epsilon_N,
 \qquad \|R_{JJ}^{-1}\|\le 2|g|^{-1}
\]
for sufficiently large $N$. For $i,k\notin J$, the minor formula implies
\[
 |R^{(J)}_{ik}-g\delta_{ik}|
 \le |g|\epsilon_N+
     (\sqrt{s_0}|g|\epsilon_N)^2\,2|g|^{-1}
 =|g|(\epsilon_N+2s_0\epsilon_N^2).
\]
Thus every such minor has entries bounded by $3$. The estimate holds for each $J$ on the same event for the full matrix. For the density, interlacing and
\eqref{W:eq:minor-minor-trace} give
\[
 \upsilon_H^{(J)}(z)
 \le\upsilon_H(z)+\frac{\pi s_0}{N\eta}.
\]
Substitute \eqref{W:eq:minor-averaged-critical} and $8\kappa<\tau$ to obtain \eqref{W:eq:minor-final-minors}. All estimates were included in the initial event for the family of
matrices.
\end{proof}

\Needspace{6\baselineskip} The following deterministic estimate extends the event to perturbations supported on finitely many coordinates.

\begin{proposition}
\label{W:prop:minor-finite-defect}
On the event of Theorem~\ref{W:thm:minor-bootstrap}, let $\Delta=\iota_S T\iota_S^*$ be any Hermitian perturbation supported on $|S|\le s_1$ coordinates, with $\|T\|\le t_N\to0$.
The support and perturbation may depend on $H$. Uniformly over all such perturbations,
\begin{equation}
 \max_{i,j}|R_{H+\Delta,ij}-R_{H,ij}|\le C_{s_1}t_N,
 \qquad
 |m_{H+\Delta}-m_H|\le\frac{\pi s_1}{N\eta}.
 \label{W:eq:minor-finite-transfer}
\end{equation}
Every complementary minor obtained by deleting at most $s_0$ coordinates has bounded entries and satisfies the $N^{-1/3+\tau}$ density bound, with a larger fixed $\tau$ if needed.
In particular, matrices with the selected entry set to zero and their radial segments of the joint real coordinate vectors are covered simultaneously.
\end{proposition}
\begin{proof}
The exact finite-rank identity is
\[
 R_{H+\Delta}=R_H+
 R_{H,\cdot S}T(I-R_{H,SS}T)^{-1}R_{H,S\cdot}.
\]
For $N$ such that $Cs_1t_N\le1/2$, the Neumann series gives $\|(I-R_{H,SS}T)^{-1}\|\le2$. For every $i,j$, including indices in $S$, the restricted row and column have norms at
most $C\sqrt{s_1}$. Consequently
\[
 |R_{H+\Delta,ij}-R_{H,ij}|
 \le (C\sqrt{s_1})t_N\,2(C\sqrt{s_1})
 =2C^2s_1t_N.
\]
The rank bound in \eqref{W:eq:minor-rank-trace} gives the second estimate. Write $\widehat R=R_{H+\Delta}$ and let $\varepsilon_N=\sup_{i,j,z}|\widehat R_{ij}(z)-g(z)\delta_{ij}|$,
with $z$ in the domain under consideration. The original entrywise bound and the preceding inequality give $\varepsilon_N\to0$ uniformly over the perturbations. On this fixed
domain let $g_*=\inf_z|g(z)|>0$ and $g^*=\sup_z|g(z)|<\infty$. For any deletion set $J$ with $|J|\le s_0$,
\[
 \widehat R_{JJ}=g(I_J+D_J),\qquad
 \|D_J\|\le s_0\varepsilon_N/g_*,\qquad
 \|\widehat R_{JJ}^{-1}\|\le2/g_*
\]
for all sufficiently large $N$. If $i,k\notin J$, both restricted rows appearing in the minor identity have norm at most $\sqrt{s_0}\varepsilon_N$. Hence
\[
 |\widehat R^{(J)}_{ik}-g\delta_{ik}|
 \le\varepsilon_N+2s_0\varepsilon_N^2/g_*.
\]
In particular the minor entries are bounded by $g^*+1$ once the displayed error is at most one, on the same event. The full matrix rank estimate followed by interlacing for the
deletion of $J$ gives
\[
 \upsilon_{H+\Delta}^{(J)}(z)
 \le\upsilon_H(z)+\frac{\pi(s_1+s_0)}{N\eta}.
\]
On the stated domain $\eta\ge\eta_*=N^{-2/3-\kappa}$, with $\kappa=\tau/100$ fixed. Hence
\[
 \frac{\pi(s_1+s_0)}{N\eta}
 \le\pi(s_1+s_0)N^{-1/3+\kappa}
 =o(N^{-1/3+\tau'})\qquad(\tau'>\tau).
\]
This proves the claimed density bound with the permitted larger exponent, uniformly over the support, perturbation, and deletion set.

For any two permitted values of the selected entry, the difference matrix is supported on at most two coordinates and has norm at most $Ce_N$. The preceding estimate therefore
holds uniformly over the whole set of permitted values.
\end{proof}

\Needspace{6\baselineskip} We next apply the resolvent derivative formula on this event, with the normalization $N^{-1/2}$ for each entry coordinate.

\begin{corollary}
\label{W:cor:minor-derivatives}
On the same event, let ${\mathcal R}$ be a resolvent of the full matrix or of a minor obtained by deleting a fixed number of coordinates. This also includes resolvents of the
perturbed matrices described above. For $1\le q\le5$ and $W_a=N^{-1/2}T_a$, where each $T_a$ is supported on a single matrix entry and has entries bounded by a fixed constant,
\begin{align}
 \max_{i,j}|D_{H}^q{\mathcal R}[W_1,\ldots,W_q]_{ij}|
       &\le C_qN^{-q/2},\label{W:eq:minor-derivative-entry}\\
 \max_i\bigl(\|D_{H}^q{\mathcal R}\,e_i\|_2+\|e_i^*D_{H}^q{\mathcal R}\|_2\bigr)
       &\le C_qN^{-q/2}\eta^{-1/2},
                                  \label{W:eq:minor-derivative-row}\\
 \|D_{H}^q{\mathcal R}\|_{\rm HS}&\le C_qN^{-q/2}\eta^{-1},\qquad
 N^{-1}|\operatorname{Tr}D_{H}^q{\mathcal R}|
       \le C_qN^{-1-q/2}\eta^{-1}.
                                  \label{W:eq:minor-derivative-hs}
\end{align}
For a minor, restrict each $W_a$ to its surviving coordinates. The matrices $T_a$ may coincide or correspond to different real coordinates of complex entries.
\end{corollary}
\begin{proof}
Apply \eqref{W:eq:resolvent-permutations} and write $W_a=N^{-1/2}\sum_{u=1}^{4}c_{a,u}e_{i(a,u)}e_{j(a,u)}^*$, allowing zero coefficients, with $|c_{a,u}|\le C$. For each
permutation and choice of matrix units the product is
\begin{equation}
 N^{-q/2}c\,({\mathcal R}e_{a_1})
 \left(\prod_{k=1}^{q-1}{\mathcal R}_{b_k a_{k+1}}\right)
 (e_{b_q}^*{\mathcal R}),\qquad |c|\le C^q.
 \label{W:eq:matrix-unit-resolvent-product}
\end{equation}
The product in parentheses is one for $q=1$. Bounded entries and the Ward identity imply, for this individual term $A$,
\begin{align*}
 |A_{ij}|&\le C_qN^{-q/2},\\
 \|Ae_j\|_2+\|e_i^*A\|_2
 &\le C_qN^{-q/2}\max_a(\|{\mathcal R}e_a\|_2+
                                      \|e_a^*{\mathcal R}\|_2)
 \le C_qN^{-q/2}\eta^{-1/2},\\
 \|A\|_{\rm HS}
 &\le C_qN^{-q/2}\|{\mathcal R}e_{a_1}\|_2
                            \|e_{b_q}^*{\mathcal R}\|_2
 \le C_qN^{-q/2}\eta^{-1},\\
 |\operatorname{Tr}A|
 &\le C_qN^{-q/2}|e_{b_q}^*{\mathcal R}^2e_{a_1}|
 \le C_qN^{-q/2}\|e_{b_q}^*{\mathcal R}\|_2
                            \|{\mathcal R}e_{a_1}\|_2
 \le C_qN^{-q/2}\eta^{-1}.
\end{align*}
There are at most $q!4^q$ terms. Summing these inequalities, and multiplying the last by $N^{-1}$ for the normalized trace, proves
\eqref{W:eq:minor-derivative-entry}--\eqref{W:eq:minor-derivative-hs}. Restricting a direction to a minor leaves a matrix unit on the surviving coordinates or gives zero.
\end{proof}

\Needspace{6\baselineskip}

\begin{corollary}
\label{W:cor:minor-adaptive-scale}
Under the hypotheses of Theorem~\ref{W:thm:minor-bootstrap}, for every fixed bound on the number of deleted indices and every finite set of nonreal spectral parameters near the
edge with $|E-2|\le N^{-2/3+o(1)}$ and $\eta=N^{-2/3+o(1)}$, the complementary minors satisfy
\[
 \max_{i,j}|{\mathcal R}_{ij}|=O(1),\qquad
 -N^{-1}\operatorname{Im}\operatorname{Tr}{\mathcal R}\prec N^{-1/3},\qquad
 \Theta_N({\mathcal R};z)+\widetilde\Theta_N({\mathcal R};z)\prec N^{-1/3}.
\]
Here $X_N\prec N^{-1/3}$ means that, for every $\epsilon,D>0$, there is $N_0(\epsilon,D)$ such that $\mathbb P(X_N>N^{-1/3+\epsilon})\le N^{-D}$ for all $N\ge N_0(\epsilon,D)$,
uniformly over the specified families. The statement applies to the uncentered truncated arrays and to the bounded and truncated Gaussian arrays used in the comparison. Applying it
to the original ensemble adds its $o(1)$ probability that truncation changes the matrix.
\end{corollary}
\begin{proof}
Fix $\epsilon,D>0$ and a bound on the deleted indices. Choose
\[
 0<\tau\le\min\{\epsilon/8,10^{-3}\},\qquad
 \kappa=\phi=\tau/100.
\]
The meaning of the subpolynomial bounds on the spectral parameters is that, for all sufficiently large $N$, every point of the finite set in the upper half plane satisfies
\[
 |E-2|\le N^{-2/3+\kappa},\qquad
 N^{-2/3-\kappa}\le\eta\le N^{-2/3+\kappa}.
\]
Thus these spectral parameters belong to $\mathcal A_N$. Apply Theorem~\ref{W:thm:minor-bootstrap} with failure exponent $D+2$. It gives, simultaneously over the specified matrices
and deletions,
\[
 \max_{i,j}|\mathcal R_{ij}|\le3,\qquad
 \upsilon^{(I)}:=-N^{-1}\operatorname{Im}\operatorname{Tr}\mathcal R
                      \le N^{-1/3+\tau}.
\]

Proposition~\ref{W:prop:minor-finite-defect} preserves the entry bound under all allowed changes on a fixed number of coordinates. Its rank estimate and interlacing give, with
fixed $c_0$,
\[
 \upsilon_{H+\Delta}^{(I)}
 \le\upsilon_H+\frac{\pi c_0}{N\eta}
 \le N^{-1/3+\tau}+\pi c_0N^{-1/3+\kappa}
 \le N^{-1/3+2\tau}
\]
for large $N$. These deterministic estimates hold on the same event for every allowed support, including a support selected from the matrix.

Substitute the coordinate bound, $\upsilon^{(I)}\le N^{-1/3+2\tau}$, and $\eta\ge N^{-2/3-\kappa}$ into \eqref{W:eq:minor-bootstrap-theta}. The first four terms are bounded by
\[
 N^{-1/2+o(1)},\qquad
 N^{-1/3+\tau+\kappa/2},\qquad
 N^{-3/8+\kappa/2},\qquad N^{-5/12}.
\]
The mean and covariance terms are, respectively,
\[
 N^{-17/24+\kappa+o(1)}+N^{-5/12+\kappa+o(1)},\qquad
 N^{-13/12+\kappa+o(1)},\qquad N^{-7/12+o(1)}.
\]
Their exponents follow from $u_N\le N^{-7/8+o(1)}$, $b_N=4N^{-5/24}$ and $\delta_N\le N^{-7/12+o(1)}$. The largest displayed exponent is $-1/3+\tau+\kappa/2$. The choices above
give
\[
 \tau+\kappa/2\le\frac{201}{1600}\epsilon<\epsilon/4.
\]
All other exponents are at most $-1/3$ because $\kappa\le10^{-5}$. The finite sum is therefore at most $N^{-1/3+\epsilon}$ for all sufficiently large $N$, with room to absorb its
subpolynomial factors. Using $\widetilde\Theta_N\le C(\Theta_N+N^{-90})$ gives
\[
 \Theta_N(\mathcal R;z)+\widetilde\Theta_N(\mathcal R;z)
                                  \le N^{-1/3+\epsilon}
\]
on the same event for all sufficiently large $N$. Its exceptional probability is at most $N^{-D}$; $\epsilon,D$ were arbitrary. This proves both stochastic domination estimates.
Adjunction covers the lower half-plane.

The truncated, bounded, and truncated Gaussian arrays satisfy the hypotheses uniformly. For the original matrix the probability of a failure is bounded by the corresponding
truncated failure probability plus $\mathbb P(H_N\ne H_N^{\rm tr})=o(1)$.
\end{proof}

\subsubsection{The spectral domain and trace contours}
\label{W:sec:resolvent-domain-for-theta}

The localization statistic uses the matrix function $\theta_E(m_N^2(2I-H))$. Its trace detects eigenvalues below $E$ because $\theta_E(x)\ge1/4$ there. To compare this statistic,
the resolvent bounds must hold on contours adapted to that function, and its derivatives must take account of its poles. The estimates below provide this spectral domain, together
with a pole expansion and derivative bounds.

We use the comparison family \eqref{W:eq:minor-family} with $m_N=N^{1/3}$, $A=m_N^2(2I-H)$ and $R_H(z)=(zI-H)^{-1}$. Empirical measures of principal minors are normalized by the
parent dimension $N$.

Uniformity over a polynomial family means simultaneous bounds for a specified collection of at most $N^a$ intermediate matrices, with fixed $a$ and arbitrary joint coupling. The
following counting bound includes principal minors and changes to the matrix of bounded rank.

\begin{proposition}
\label{W:prop:gaussian-counting}
Let
\[
 \mathfrak K_N(H)=1+N\sup_{t\in\mathbb R}|F_H(t)-F_{\rm sc}(t)|.
\]
For the above comparison family, $\max_H\mathfrak K_N(H)\prec1$. For any fixed bounds on the number of deleted coordinates and the perturbation rank, the same conclusion holds
simultaneously for all such principal minors and perturbations of every member. The minor measures have mass $(N-|J|)/N$. More precisely, a rank-$r$ perturbation increases
$\mathfrak K_N$ by at most $r$, and deleting $s$ coordinates increases it by at most $s$. On the same event $\|H\|\le4$. If $V$ is supported on finitely many coordinates and
$\|V\|\to0$, then $\|H+V\|\le5$ whenever $\|V\|\le1$.
\end{proposition}
\begin{proof}
Fix $\epsilon,D>0$, and choose $0<\phi<\min\{1/100,\epsilon/20\}$. The cited coarse decomposition is $H=W+V+B$, with
\[
 \operatorname{rank}B\le 4N^{5\phi},\qquad
 \|V\|\le N^{-1+4\phi},\qquad \|B\|=o(1),
\]
with failure probability smaller than any prescribed inverse power of the matrix size. On the event where its eigenvalue counting estimate holds, $\|W\|\le3$ and
$\|F_W-F_{\rm sc}\|_\infty\le \mathfrak c_N/N$, with $\mathfrak c_N=N^{o(1)}$. Write $v=\|V\|$, $r=\operatorname{rank}B$. Weyl's inequality and rank interlacing give, for every
real $t$,
\[
 F_W(t-v)-r/N\le F_H(t)\le F_W(t+v)+r/N.
\]
Since $\|\rho_{\rm sc}\|_\infty=1/\pi$,
\begin{equation}
 \mathfrak K_N(H)\le1+\mathfrak c_N+r+Nv/\pi.
 \label{W:eq:gaussian-counting-transfer}
\end{equation}
Choose the factor in the counting estimate with $\mathfrak c_N\le N^{\epsilon/2}$ for large $N$ and fix $5\phi<\epsilon/2$. Then $r\le4N^{\epsilon/2}$ and $Nv\le N^{\epsilon/2}$,
so \eqref{W:eq:gaussian-counting-transfer} is at most $N^\epsilon$ for every $N\ge N_0(\epsilon)$. For each fixed $\epsilon>0$, this choice of $\phi$ gives the stated bound with
probability at least $1-N^{-D}$ for every fixed $D$ and all sufficiently large $N$. Strengthen the individual failure exponent to $D+a+2$ and take a union over the family of at
most $N^a$ matrices. The family need not be independent. The resulting event holds simultaneously over the finite family. For a principal minor with $s$ indices removed, write
$n_{H^{(J)}}(t)$ for its unnormalized counting function. The inequality $n_H(t)-s\le n_{H^{(J)}}(t)\le n_H(t)$ gives
\[
 \sup_t|n_{H^{(J)}}(t)-NF_{\rm sc}(t)|
 \le\sup_t|n_H(t)-NF_{\rm sc}(t)|+s.
\]
For a rank-$r$ perturbation the same statement holds with $r$ in place of $s$. These inequalities are for the right-continuous counting functions, so they include eigenvalues equal
to $t$. Letting $t$ increase to any boundary point proves the corresponding open boundary point inequalities as well, since $F_{\rm sc}$ is continuous. The norm bound follows from
$\|W\|\le3$, $\|V\|=o(1)$ and $\|B\|=o(1)$.
\end{proof}

There is a deterministic $K_N=N^{o(1)}$ for which these events have failure smaller than any fixed inverse polynomial, uniformly over the fixed family. Choose increasing $N_j$ so
that the bound by $N^{1/j}$ has failure at most $N^{-j}$ for all $N\ge N_j$, and set $K_N=N^{1/j_N}$ with $j_N=\max\{j:N_j\le N\}$. Add the fixed bounds on the number of deleted
indices and on the rank to $K_N$. The resulting sequence bounds the counting error on the event.

\Needspace{6\baselineskip} Integration by parts converts the counting estimate into an averaged resolvent bound at every nonreal spectral parameter.

\begin{corollary}
\label{W:cor:gaussian-averaged}
If $\mathfrak K_N(H)\le K$, then for every $\eta>0$ and real $u$,
\begin{equation}
 |m_H(u+\mathrm i\eta)-g(u+\mathrm i\eta)|
 \le\frac{\pi K}{N\eta},\qquad
 \upsilon_H(u+\mathrm i\eta)
 \le\upsilon(u+\mathrm i\eta)+\frac{\pi K}{N\eta}.
 \label{W:eq:gaussian-global-average}
\end{equation}
The same formulas hold for minors with their enlarged $K$ and common normalization. They hold for finite-rank perturbations of arbitrary amplitude with the corresponding rank
increment.
\end{corollary}
\begin{proof}
Put $D(t)=F_H(t)-F_{\rm sc}(t)$. The counting hypothesis gives $\|D\|_\infty\le K/N$. For $z=u+\mathrm i\eta$ with $\eta>0$, Stieltjes integration by parts gives
\[
 m_H(z)-g(z)
 =\int_{\mathbb R}\frac{\mathrm dD(t)}{z-t}
 =-\int_{\mathbb R}\frac{D(t)}{(z-t)^2}\,\mathrm dt.
\]
To justify the second equality, first integrate on a finite interval and then send both endpoints to infinity. The boundary terms tend to zero because $D$ is bounded and
$(z-t)^{-1}\to0$. The remaining integral converges absolutely since $|z-t|^{-2}$ is integrable. Consequently,
\[
 |m_H(z)-g(z)|
 \le\frac K N\int_{\mathbb R}
                  \frac{\mathrm dt}{(t-u)^2+\eta^2}
 =\frac{\pi K}{N\eta}.
\]
With the convention $R_H(z)=(z-H)^{-1}$, one has $\upsilon_H(z)=-\operatorname{Im}m_H(z)$ and $\upsilon(z)=-\operatorname{Im}g(z)$. Thus
\[
 \upsilon_H(z)-\upsilon(z)
 =-\operatorname{Im}(m_H(z)-g(z))
 \le|m_H(z)-g(z)|,
\]
which proves the second inequality.

For a minor, keep the normalization by $N$ in its counting function and trace. Its counting difference need not tend to zero at $+\infty$, but it remains bounded, so the boundary
term above still vanishes. Proposition~\ref{W:prop:gaussian-counting} enlarges $K$ by at most the number of deleted indices, or by the rank of a perturbation. The same integral
bound therefore applies in both cases, without a restriction on the amplitude of the perturbation.
\end{proof}

\Needspace{6\baselineskip} The coordinate estimates hold on one nonreal domain covering both edges, the bulk, and the closing portions needed for functional calculus.

Fix $b\ge6$, $c>0$, $M\ge1$, a bound on the number of deleted indices, and $0<\epsilon<1/24$.

\begin{theorem}
\label{W:thm:gaussian-domain}
On an event whose failure probability is smaller than any prescribed inverse power of the matrix size,
\begin{equation}
 \max_{i,j}|R_{H,ij}(z)-g(z)\delta_{ij}|
 \le C\bigl(e_N+N^{-1/6+\epsilon}\bigr)|g(z)|,
 \quad |\operatorname{Re}z|\le b,\quad
 cN^{-2/3}\le |\operatorname{Im}z|\le M.
 \label{W:eq:gaussian-full-domain}
\end{equation}
For each fixed bound on the number of deleted coordinates, all corresponding complementary minors have bounded entries on this domain. These conclusions are simultaneous over the
family of intermediate matrices and all perturbations supported on a fixed number of coordinates with norm $t_N\to0$, with an additional $C t_N$ in the error. On every edge window
$|u\mp2|\le C_1N^{-2/3}$, $cN^{-2/3}\le\eta\le C_1N^{-2/3}$,
\begin{equation}
 \upsilon_H^{(J)}(u+\mathrm i\eta)\prec N^{-1/3}.
 \label{W:eq:gaussian-edge-density}
\end{equation}
\end{theorem}

For $|\operatorname{Im}z|\ge M$, the deterministic bound $\|R_H(z)\|\le |\operatorname{Im}z|^{-1}$ applies. For the closing segments $z=\pm b+\mathrm iy$, $|y|\le cN^{-2/3}$, one
has $\|R_H^{(J)}(z)\|\le1$ on the norm event, also after the stated small perturbations.

\begin{proof}
It suffices to treat the upper half-plane; adjoints give the other half. The extension of Theorem~\ref{W:thm:minor-bootstrap} needed here concerns the whole compact real interval.
We carry out the continuity argument using bounds valid throughout the bulk.

Use the auxiliary graph of Section~\ref{W:sec:block-schur-complements-under-the-tail-assumption}, at separation $\ell_N^{\mathrm{Sch}}=N^{7/24}$. If the matrix family and the set
of spectral parameters have at most $N^J$ elements, choose the fixed integer $d>6(D+J+2)$. The union bound for components then gives $C_dN^{J+1-d/6}\le N^{-D}$ for sufficiently
large $N$. For a component $I$, let $G=R_H^{(I)}$, $B=H_{I,I^c}$, and $Z_I=BGB^*-m_H^{(I)}I_I$. Theorem~\ref{W:thm:schur-residual} and the conditional linear estimate give, on a
finite covering set of polynomial cardinality,
\begin{align}
 \|Z_I\|&\le C_d(\log N)^4\Theta_N(G;z),
 \label{W:eq:gaussian-residual-start}\\
 \|BGe_j\|_2&\le C_d(\log N)^2
 \left[(N^{-1/2}+u_N)\|Ge_j\|_2+b_N\|Ge_j\|_\infty\right],
 \label{W:eq:gaussian-border-start}
\end{align}
where $b_N=4N^{-5/24}$, $u_N\le N^{-7/8+o(1)}$, and $\delta_N\le N^{-7/12+o(1)}$. These events are obtained by conditioning on the graph and the minor.

A union over the finite covering set of polynomial cardinality and the at most $N$ components is valid. The adaptive normalization is exactly the one of \eqref{W:eq:schur-theta}.

In the proof we may initialize at a larger fixed height $M_*\ge M$. On this enlarged compact domain, $g$ is bounded away from zero and $|g|\le1$. For the latter, write $z=g+g^{-1}$
and use $\operatorname{Im}g<0$ and $\operatorname{Im}z>0$. Define
\[
 \Lambda(z)=\max_{i,j}|R_{H,ij}(z)/g(z)-\delta_{ij}|.
\]
Take $\epsilon_0>0$ sufficiently small compared with $d$ and the bound on the number of deleted indices. If $\Lambda\le\epsilon_0$, the exact minor identity
\begin{equation}
 G=R_{I^cI^c}-R_{I^cI}R_{II}^{-1}R_{II^c}
 \label{W:eq:gaussian-minor-identity}
\end{equation}
gives $\max_{j,k}|G_{jk}|\le C$ and $\|R_{II}^{-1}\|\le2/|g|$. Proposition~\ref{W:prop:gaussian-counting} and Corollary~\ref{W:cor:gaussian-averaged} give, on a further event,
\[
 |m_H^{(I)}-g|\le N^{-1/3+\epsilon/10},\qquad
 \upsilon_H^{(I)}\le C
\]
throughout the domain: apply the counting estimate with exponent $\epsilon/20$ and use $C N^{\epsilon/20}\le N^{\epsilon/10}$ for sufficiently large $N$. This averaged estimate
follows from the counting bound just used.

The Ward identities imply
\[
 N^{-1}\|G\|_{\rm HS}
 =\sqrt{\frac{\upsilon_H^{(I)}}{N\eta}},\qquad
 \max_j\bigl(\|Ge_j\|_2+\|e_j^*G\|_2\bigr)\le C\eta^{-1/2},
 \qquad \|\operatorname{diag}G\|_2\le C\sqrt N.
\]
Consequently the adaptive scale is at most a constant times
\begin{equation}
 \sqrt{F_N/N}+\sqrt{\upsilon_H^{(I)}/(N\eta)}
 +\frac{b_N}{\sqrt{N\eta}}+b_N^2
 +u_N(N^{-1/2}+b_N)\eta^{-1}+u_N^2\eta^{-1}+\delta_N.
 \label{W:eq:gaussian-adaptive-bound}
\end{equation}
At the smallest height the first four powers are respectively $N^{-1/2+o(1)}$, $N^{-1/6}$, $N^{-3/8}$, $N^{-5/12}$. The remaining powers are at most $N^{-17/24+o(1)}$,
$N^{-5/12+o(1)}$, $N^{-13/12+o(1)}$, and $N^{-7/12+o(1)}$. Here $F_N=N^{o(1)}$. Thus
\[
 \|Z_I\|\le N^{-1/6+\epsilon/2},\qquad
 \|BGe_j\|_2\le N^{-1/6+\epsilon/2}
\]
for all $N\ge N_0(\kappa,d)$, by \eqref{W:eq:gaussian-residual-start}--\eqref{W:eq:gaussian-border-start}.

The exact Schur complement identity is
\[
 R_{II}=
 \left[g^{-1}I_I-
 \bigl(H_{II}+(m_H^{(I)}-g)I_I+Z_I\bigr)\right]^{-1}.
\]
Put $D_I=H_{II}+(m_H^{(I)}-g)I_I+Z_I$. Then
\[
 \|D_I\|\le d e_N+N^{-1/3+\epsilon/10}+N^{-1/6+\epsilon/2}=o(1),
\]
and the Schur complement formula gives
\[
 R_{II}=g(I_I-gD_I)^{-1},\qquad
 R_{II}-gI_I=g^2(I_I-gD_I)^{-1}D_I.
\]
Choose $N_0$ such that $\|gD_I\|\le1/2$ for $N\ge N_0$. Then
\[
 \|R_{II}\|\le2|g|,\qquad
 |g|^{-1}\|R_{II}-gI_I\|\le2|g|\|D_I\|
 \le C_d(e_N+N^{-1/6+\epsilon/2}).
\]
For exterior coordinates use $R_{Ij}=R_{II}BGe_j$. This proves the strict improvement
\[
 \Lambda\le C_d(e_N+N^{-1/6+\epsilon/2})<\epsilon_0/2.
\]

On $\|H\|\le4$ and at height $M_*>8$, the resolvent identity at infinity gives
\[
 \|R_H(z)-z^{-1}I\|\le\frac4{|z|(|z|-4)},\qquad
 |g(z)-z^{-1}|\le C|z|^{-3},\qquad |g(z)|\ge c|z|^{-1}.
\]
Consequently $\Lambda(z)\le C/M_*$. Choose $M_*$ so that this is less than $\epsilon_0/4$. For the comparison between consecutive spectral parameters, the explicit equation
$g^2-zg+1=0$ gives $|g'|\le C(1+|z^2-4|^{-1/2})$ on the fixed real interval. At $\eta\ge cN^{-2/3}$ this is at most $C N^{1/3}$. Using $|R_{ij}|\le\eta^{-1}$ and the lower bound on
$|g|$ gives
\[
 |\partial_z(R_{ij}/g)|
 \le C\eta^{-2}+C\eta^{-1}N^{1/3}\le CN^2.
\]
A step between neighboring spectral parameters of size $N^{-5}$ therefore changes $\Lambda$ by at most $C N^{-3}$. Starting from the top row, the value at the next spectral
parameter is below $\epsilon_0$ and is improved by the preceding calculation using the Schur complement to $\epsilon_0/2$. Induction reaches every row. Interpolation between
selected points gives \eqref{W:eq:gaussian-full-domain} throughout the domain. For a set $J$ of bounded size, put $\epsilon_N= C(e_N+N^{-1/6+\epsilon})$. The minor identity gives
explicitly
\[
 |R^{(J)}_{ik}-g\delta_{ik}|
 \le |g|(\epsilon_N+2|J|\epsilon_N^2),\qquad i,k\notin J,
\]
once $|J|\epsilon_N\le1/2$. This controls all such minors on the same event.

For a finite coordinate perturbation $\Delta=\iota_S T\iota_S^*$ use
\[
 R_{H+\Delta}=R_H+
 R_{H,\cdot S}T(I-R_{H,SS}T)^{-1}R_{H,S\cdot}.
\]
All compression entries are bounded, and $\|T\|\le t_N\to0$. The correction is bounded entrywise by $C_{|S|}t_N$, uniformly over the support, its dependence on $H$, and every
radial parameter. For a deletion set $J$ of bounded size, write $\widetilde R_{JJ}=gI_J+E_J$. Then $\|E_J\|\le C_J(\epsilon_N+t_N)$, so $\|\widetilde R_{JJ}^{-1}\|\le2/|g|$
whenever $C_J(\epsilon_N+t_N)\le |g|/2$. The minor identity gives
\[
 |\widetilde R^{(J)}_{ik}-g\delta_{ik}|
 \le C_J\bigl(\epsilon_N+t_N+(\epsilon_N+t_N)^2\bigr),
 \qquad i,k\notin J.
\]
Thus the perturbed minor entries are uniformly bounded. Rank comparison preserves the counting bound independently of the perturbation amplitude. On the edge windows, the
square-root formula gives $\upsilon(u+\mathrm i\eta)\le C\sqrt{|u\mp2|+\eta} \le C'N^{-1/3}$. If $s$ indices are removed and the perturbation has rank at most $r$, then
\eqref{W:eq:gaussian-global-average} gives
\[
 \upsilon_{H+\Delta}^{(J)}(u+\mathrm i\eta)
 \le C'N^{-1/3}+\frac{\pi(\mathfrak K_N(H)+s+r)}{N\eta}
 \le N^{-1/3}\left[C'+\frac\pi c(\mathfrak K_N(H)+s+r)\right].
\]
For any $\varepsilon,D>0$, use the counting bound $\mathfrak K_N(H)\le N^{\varepsilon/2}$ with failure at most $N^{-D}$. The bracket is at most $N^\varepsilon$ for all
$N\ge N_0(\varepsilon)$. This proves \eqref{W:eq:gaussian-edge-density}, simultaneously over the specified families. Finally, on $\|H+\Delta\|\le5$, the distance of
$\pm b+\mathrm iy$ from the real spectrum is at least one. The two deterministic resolvent bounds follow from the spectral theorem.
\end{proof}

\Needspace{6\baselineskip} On the same domain, finite resolvent products give bounds for matrix entry derivatives of every fixed order.

\begin{corollary}
\label{W:cor:gaussian-derivatives}
Let $W_j=N^{-1/2}T_j$, where each $T_j$ is supported on one matrix entry and has bounded norm. For fixed $q\ge1$,
\[
 D_{H}^qR_H(z)[W_1,\ldots,W_q]
 =\sum_{\sigma\in S_q}R_HW_{\sigma(1)}R_H\cdots W_{\sigma(q)}R_H.
\]
On the nonreal part of Theorem~\ref{W:thm:gaussian-domain},
\begin{align*}
 \max_{i,j}|(D_{H}^qR_H)_{ij}|&\le C_qN^{-q/2},\\
 \max_i\bigl(\|D_{H}^qR_He_i\|_2+\|e_i^*D_{H}^qR_H\|_2\bigr)
 &\le C_qN^{-q/2}\eta^{-1/2},\\
 \|D_{H}^qR_H\|_{\rm HS}&\le C_qN^{-q/2}\eta^{-1},\\
 N^{-1}|\operatorname{Tr}D_{H}^qR_H|&\le C_qN^{-1-q/2}\eta^{-1}.
\end{align*}
On the vertical sides of the contour the same estimates hold with every $\eta^{-1/2}$ or $\eta^{-1}$ replaced by $1$. The statements include mixed real coordinates of a complex
block and repeated labels.
\end{corollary}
\begin{proof}
For each of the at most $q!4^q$ matrix unit terms use \eqref{W:eq:matrix-unit-resolvent-product} with ${\mathcal R}=R_H$. The common scalar coefficient is $O_q(N^{-q/2})$. If
$L(z)=\max_a\{\|R_He_a\|_2,\|e_a^*R_H\|_2\}$, the four norms of that term are bounded, respectively, by
\[
 C_qN^{-q/2},\qquad C_qN^{-q/2}L(z),\qquad
 C_qN^{-q/2}L(z)^2,\qquad C_qN^{-1-q/2}L(z)^2.
\]
The last expression is for the normalized trace and follows from $|(R_H^2)_{ba}|\le\|e_b^*R_H\|_2\|R_He_a\|_2$. On the upper half-plane the Ward identity reads
\[
 \|R_He_a\|_2^2=\|e_a^*R_H\|_2^2
       =-\eta^{-1}\operatorname{Im}(R_H)_{aa}\le C/\eta.
\]
For the lower half-plane use $R_H(\bar z)=R_H(z)^*$ and replace $\eta$ by $|\operatorname{Im}z|$. On a vertical side of the contour the spectral theorem gives
\[
 \|R_H(z)\|=\max_{\lambda\in\sigma(H)}|z-\lambda|^{-1}\le1,
\]
since $|\operatorname{Re}z|=b\ge6$ and $\|H\|\le5$ there. Both the row and column factors are therefore at most one. Substituting these two bounds and summing the finite collection
proves all four estimates on their respective domains.
\end{proof}

\subparagraph{Contour representations and poles}
\label{W:sec:contour-representations-of-theta}

Fix $0<a<\pi/4$ and let $\mathcal C_N$ be the counterclockwise boundary of
\[
 [-6,6]+\mathrm i[-a/{m_N}^2,a/{m_N}^2].
\]
On the above event, this contour encloses the whole spectrum for every permitted perturbation and stays away from it. Its horizontal sides are covered by
Theorem~\ref{W:thm:gaussian-domain}; its vertical sides are covered by the spectral distance estimate. For every polynomial $p$, Cauchy's formula gives
\begin{equation}
 D_{H}^qp(H)[W_1,\ldots,W_q]
 =\frac1{2\pi\mathrm i}\int_{\mathcal C_N}p(z)
 \sum_{\sigma\in S_q}R_H(z)W_{\sigma(1)}\cdots
 W_{\sigma(q)}R_H(z)\,\mathrm dz,
 \label{W:eq:gaussian-polynomial-contour}
\end{equation}
For $q=0$ the integrand contains one resolvent. To justify the derivatives, fix $N$ and a matrix in the norm event, and put
\[
 \delta_H=\min\{a/m_N^2,\,6-\|H\|\}>0.
\]
If $\|\sum_jt_jW_j\|<\delta_H/2$, the resolvent identity and its Neumann series give
\[
 \sup_{z\in\mathcal C_N}
 \bigl\|(zI-H-\textstyle\sum_jt_jW_j)^{-1}\bigr\|
 \le 2/\delta_H.
\]
Thus the same contour surrounds the spectra of all these nearby self-adjoint matrices. The derivative of order $q$ of its integrand has norm at most
\[
 q!\,\max_{z\in\mathcal C_N}|p(z)|
           (2/\delta_H)^{q+1}\prod_{j=1}^q\|W_j\|.
\]
This constant is finite at the fixed $N$, and the contour has finite length. Difference quotients are therefore dominated on a common neighbourhood, which proves
\eqref{W:eq:gaussian-polynomial-contour} by differentiation under the integral.

The same argument applies to the holomorphic function $z\mapsto\theta_E(m_N^2(2-z))$: its closest nonreal poles have imaginary part of absolute value $\pi/m_N^2$, greater than the
contour height $a/m_N^2$. For an arbitrary finite $H$ outside the norm event, replace the real interval $[-6,6]$ by $[-b,b]$ with $b>\|H\|+1$, keeping the same imaginary height.
This proves the functional-calculus identity at that matrix as well. Quantitative bounds uniform over such matrices are supplied by the pole expansion below; the present
differentiation argument is for a fixed $N$ and contour.

Bounds uniform in the degree also require control of $p$ on $\mathcal C_N$; the example following Proposition~\ref{W:prop:gaussian-chebyshev} exhibits its possible growth.

For fixed $E\in\mathbb R$, set
\begin{equation}
 f_E(x)=\frac1{1+\exp(x-E)},\qquad
 \theta_E(x)=f_E(x)^2.
 \label{W:eq:theta-factorization}
\end{equation}
It is holomorphic in $|\operatorname{Im}z|<\pi$, real and strictly positive on the real axis, and
\begin{equation}
 0<\theta_E\le1,\qquad
 \theta_E(x)\ge\tfrac14\quad(x\le E),\qquad
 \theta_E(x)\le\exp[-2(x-E)]\quad(x\ge E).
 \label{W:eq:theta-bounds}
\end{equation}
For each $a<\pi$, $|\theta_E(x+\mathrm iy)|\le C_a\exp[-2(x-E)_+]$ when $|y|\le a$. Indeed, $|1+r\mathrm e^{\mathrm iy}|\ge(1+r)\cos(a/2)$ for $r\ge0$. Since $\theta_E(x)\to1$ as
$x\to-\infty$, the bound $\theta_E(x)\ge1/4$ applies to every eigenvalue below $E$.

\Needspace{6\baselineskip} The counting estimate gives a bound for $\operatorname{Tr}\theta_E(A)$, including the contribution of every eigenvalue below $E$.

\begin{lemma}
\label{W:lem:gaussian-trace}
Assume $\mathfrak K_N(H)\le K$, and let $x_i$ be the eigenvalues of $A={m_N}^2(2I-H)$ with multiplicity. Then, for every real $t$,
\begin{equation}
 \#\{i:x_i\le t\}\le K+C(t_+)^{3/2},\qquad
 \operatorname{Tr}\theta_E(A)\le K+C_E.
 \label{W:eq:gaussian-full-trace}
\end{equation}
Also, for $s\ge1$ and real $x$,
\begin{equation}
 \#\{i:|x_i-x|\le s\}
 \le C\bigl(K+s\sqrt{x_+}+s^{3/2}\bigr).
 \label{W:eq:gaussian-local-count}
\end{equation}
These estimates hold for minors and for perturbations of bounded rank, with the corresponding enlarged counting constant $K$. In particular the trace is $N^{o(1)}$ on the event
where the eigenvalue counting estimate holds, even if some $x_i\to-\infty$.
\end{lemma}
\begin{proof}
Under the change of variables $u={m_N}^2(2-t)$ the unnormalized semicircle measure has density
\[
 \frac{N}{{m_N}^2}\rho_{\rm sc}(2-u/{m_N}^2)
 =\frac1{2\pi}\sqrt{4u-u^2/{m_N}^2}\,
 \mathbf1_{[0,4{m_N}^2]}(u)
 \le\frac1\pi\sqrt{u_+}.
\]
Let $n\le N$ be the dimension of $H$ and put $v=2-t/m_N^2$. With $F_H$ normalized by the parent dimension $N$,
\[
 \#\{i:x_i\le t\}=n-NF_H(v-)
 \le n-NF_{\rm sc}(v)+K
 \le N(1-F_{\rm sc}(v))+K.
\]
Continuity of $F_{\rm sc}$ justifies its value at the left limit. Integration of the density just computed gives $N(1-F_{\rm sc}(v))\le2(t_+)^{3/2}/(3\pi)$. This is also valid for
a minor, since $n-N\le0$; the mass difference contributes the nonpositive term $n-N$. For the interval $[x-s,x+s]$, the difference of the two counting functions has error at most
$2K$. Its semicircle contribution is at most
\[
 \frac1\pi\int_{x-s}^{x+s}\sqrt{u_+}\,\dd u
 \le\frac{2s}\pi\sqrt{x_+}
                  +\frac1\pi\int_{-s}^s\sqrt{|v|}\,\dd v
 =\frac{2s}\pi\sqrt{x_+}+\frac{4}{3\pi}s^{3/2}.
\]
This proves \eqref{W:eq:gaussian-local-count}, including both boundary points. Since $\theta_E$ decreases from one to zero, Tonelli's theorem gives
\[
 \operatorname{Tr}\theta_E(A)
 =\int_{\mathbb R}[-\theta_E'(t)]\#\{i:x_i\le t\}\,\mathrm dt
 \le K+C\int_0^\infty[-\theta_E'(t)]t^{3/2}\,\mathrm dt.
\]
Its integrand satisfies $-\theta_E'(t)=2e^{t-E}(1+e^{t-E})^{-3}$. Put $T=1+\max(E,0)$. Since $-\theta_E'\ge0$ and its integral on $\mathbb R$ is one,
\[
 \int_0^\infty[-\theta_E'(t)]t^{3/2}\,\dd t
 \le T^{3/2}+2e^{2E}\int_T^\infty t^{3/2}e^{-2t}\,\dd t
 \le T^{3/2}+2^{-3/2}e^{2E}\Gamma(5/2)<\infty.
\]
\end{proof}

\Needspace{6\baselineskip} We represent $\theta_E$ by paired strip integrals and a pole series, keeping the cancellation needed for convergence at negative infinity.

\begin{theorem}
\label{W:thm:gaussian-functional-calculus}
For every finite self-adjoint matrix $A$, put $S_A(w)=(wI-A)^{-1}$. If $0<a<\pi$, then
\begin{equation}
 \theta_E(A)=\frac1{2\pi\mathrm i}\int_{\mathbb R}
 \left[\theta_E(u-\mathrm ia)S_A(u-\mathrm ia)
       -\theta_E(u+\mathrm ia)S_A(u+\mathrm ia)\right] \mathrm du.
 \label{W:eq:gaussian-paired-strip}
\end{equation}
Cancellation within the bracket gives absolute integrability at minus infinity. Each undifferentiated term fails to be absolutely integrable there. At every positive derivative
order, each integral is absolutely convergent.

Let $\omega_n=\pi(2n+1)$ and $S_n^\pm=(E\pm\mathrm i\omega_n-A)^{-1}$. Then
\begin{equation}
 \theta_E(A)=\tfrac12I+
 \sum_{n=0}^{\infty}\left[S_n^++S_n^-+(S_n^+)^2+(S_n^-)^2\right].
 \label{W:eq:gaussian-poles}
\end{equation}
The series of the displayed $+$/$-$ pairs converges in operator norm.
\end{theorem}
\begin{proof}

Integrate $\theta_E(w)(wI-A)^{-1}$ around the rectangle $[-R,R]+\mathrm i[-a,a]$, with $R>2(1+\|A\|+|E|)$. The vertical side at $R$ is $O_{A,E}(e^{-2R}/R)$ and that at $-R$ is
$O_{A,E}(R^{-1})$. On the left half-line write $S_\pm=S_A(u\pm\mathrm ia)$. Then
\[
 \begin{aligned}
 \theta_E(u-\mathrm ia)S_- -\theta_E(u+\mathrm ia)S_+
 &=2\mathrm ia S_-S_+\\
 &\quad+(\theta_E(u-\mathrm ia)-1)S_-
  -(\theta_E(u+\mathrm ia)-1)S_+.
 \end{aligned}
\]
Here $\|S_\pm\|\le C_A/|u|$ for $u<-2(1+\|A\|)$ and $|\theta_E(u\pm\mathrm ia)-1|\le C_{a,E}e^u$. Thus the norm of the bracket is at most $C_A|u|^{-2}+C_{a,A,E}e^u/|u|$. On the
right it is bounded by $C_{a,A,E}e^{-2u}/u$. These integrable majorants justify the limit $R\to\infty$ and prove \eqref{W:eq:gaussian-paired-strip}. For a derivative of order
$q\ge1$, each integrand is a sum of $q!$ products with $q+1$ resolvents, bounded at the left by $C_{q,A}|u|^{-q-1}\prod_j\|V_j\|$. This also justifies differentiation under the
integral on bounded neighborhoods of
$A$.

For the pole formula, the canonical product \cite{nist2026sg58} (\S4.36, equation (4.36.2)) gives
\[
 \cosh(y/2)=\prod_{n=0}^{\infty}(1+y^2/\omega_n^2).
\]
Its logarithmic derivative, locally uniformly off its zeros, is
\[
 \tanh(y/2)=4y\sum_{n=0}^{\infty}(y^2+\omega_n^2)^{-1}.
\]
Hence
\[
 \frac1{1+\mathrm e^y}
 =\tfrac12+\sum_{n=0}^{\infty}
 \left[(\mathrm i\omega_n-y)^{-1}+(-\mathrm i\omega_n-y)^{-1}\right].
\]
On a fixed compact set of complex numbers $|y|\le R$ avoiding the poles, all terms with $\omega_n>2R$ satisfy
\[
 \left|\frac{-2y}{y^2+\omega_n^2}\right|\le\frac{8R}{3\omega_n^2},
 \qquad
 |(\mathrm i\omega_n-y)^{-2}|+|(-\mathrm i\omega_n-y)^{-2}|
                                  \le\frac8{\omega_n^2}.
\]
The sums of these majorants converge. Thus the paired series for $f(y)=(1+e^y)^{-1}$ and its first derivative converge locally uniformly off the poles. Direct differentiation gives
\[
 f'(y)=-\frac{e^y}{(1+e^y)^2},\qquad
 f(y)+f'(y)=\frac{1+e^y-e^y}{(1+e^y)^2}=f(y)^2.
\]
Adding the two convergent series therefore proves the scalar pole identity. Apply it to the finitely many real eigenvalues of $A-EI$ in its spectral decomposition to obtain
\eqref{W:eq:gaussian-poles}. Indeed, with $X=A-EI$,
\[
 S_n^++S_n^-=-2X(X^2+\omega_n^2I)^{-1},
\]
so its norm is at most $2\|X\|/\omega_n^2$, and $\|(S_n^+)^2+(S_n^-)^2\|\le2/\omega_n^2$. This proves operator norm convergence for each finite $A$ and local uniformity on bounded
sets of $A$.

\end{proof}

\begin{corollary}
\label{W:cor:theta-matrix-derivatives}
For $q\ge1$, all mixed derivatives can be taken term by term, and
\begin{equation}
 \|D_A^q\theta_E(A)[V_1,\ldots,V_q]\|
 \le C_q\prod_{j=1}^q\|V_j\|,
 \label{W:eq:gaussian-global-derivative}
\end{equation}
where $C_q$ is independent of $A,N,E$. The same bound in a Schatten norm holds with any one factor on the right in that Schatten norm.

\end{corollary}
\begin{proof}
For a resolvent $S=S_A(w)$,
\begin{equation}
 D_A^qS[V_1,\ldots,V_q]
 =\sum_{\sigma\in S_q}SV_{\sigma(1)}S\cdots V_{\sigma(q)}S.
 \label{W:eq:gaussian-pole-derivatives}
\end{equation}

Differentiating the product $S^2$ gives
\[
 D_A^q(S^2)[V_1,\ldots,V_q]
 =\sum_{U\sqcup V=\{1,\ldots,q\}}(D_US)(D_VS).
\]
For $|U|=r$, expanding both derivatives gives $r!(q-r)!$ products for each of the $\binom qr$ assignments. Concatenating the two orders produces
\[
 \sum_{\sigma\in S_q}\sum_{r=0}^q
 SV_{\sigma(1)}S\cdots V_{\sigma(r)}S^2
            V_{\sigma(r+1)}S\cdots V_{\sigma(q)}S.
\]
The empty initial or final string is omitted when $r=0$ or $r=q$. There are $(q+1)q!$ terms, each containing $q+2$ resolvents. Since $\|S_n^\pm\|\le\omega_n^{-1}$, a possible
explicit constant is
\[
 C_q=2q!\sum_{n\ge0}\omega_n^{-q-1}
       +2(q+1)!\sum_{n\ge0}\omega_n^{-q-2}<\infty.
\]
This proves uniform absolute convergence of the derivative series and \eqref{W:eq:gaussian-global-derivative}. Let $F_J$ be the finite paired partial sum. The zero order series
converges uniformly on bounded sets of self-adjoint matrices, and the displayed bounds give uniform convergence of its derivative series. Passing to the limit in
\[
 F_J(A+tV)-F_J(A)=\int_0^t DF_J(A+sV)[V]\,\dd s
\]
identifies the derivative of the limit. Repeating this identity for the multilinear derivatives justifies every fixed order. For a term containing $V_j$ and any Schatten index
$p\ge1$, write the product as $L V_j R$. The ideal inequality gives
\[
 \|LV_jR\|_p\le\|L\|\,\|V_j\|_p\,\|R\|.
\]
The product of the norms of the other matrices $V_k$ and all resolvent bounds is the same one used above; $\|V_j\|$ is simply replaced by $\|V_j\|_p$. Summing with the same
constant $C_q$ proves the Schatten norm bound.
\end{proof}

\Needspace{6\baselineskip} After edge rescaling, these representations use resolvents at the nonreal points specified below and give bounds for truncating the pole series.

For $A={m_N}^2(2I-H)$ the only resolvent points in \eqref{W:eq:gaussian-poles} are
\begin{equation}
 z_n^\pm=2-E/{m_N}^2\mp\mathrm i\omega_n/{m_N}^2,
 \qquad S_n^\pm=-{m_N}^{-2}R_H(z_n^\pm).
 \label{W:eq:gaussian-pole-domain}
\end{equation}
Every such point below a fixed height is in Theorem~\ref{W:thm:gaussian-domain}; all larger heights have deterministic resolvent bounds. The minimum height is $\pi N^{-2/3}$.

Let $\theta_E^{[J]}(A)$ denote the sum in \eqref{W:eq:gaussian-poles} with $0\le n<J$ .

\begin{corollary}
\label{W:cor:gaussian-pole-domain}
For $J\ge1$,
\begin{align}
 \|\theta_E(A)-\theta_E^{[J]}(A)\|
 &\le C(1+\|A-EI\|)/J,
 \label{W:eq:gaussian-pole-tail-zero}\\
 \|D_A^q(\theta_E-\theta_E^{[J]})(A)[V_1,\ldots,V_q]\|
 &\le C_qJ^{-q}\prod_j\|V_j\|,\qquad q\ge1.
 \label{W:eq:gaussian-pole-tail-positive}
\end{align}
The positive order estimate also holds in a Schatten norm with one factor in that norm. If $W_j=N^{-1/2}T_j$ are matrix entry directions, then
\begin{equation}
 \|D_H^q\theta_E({m_N}^2(2I-H))[W_1,\ldots,W_q]\|
 \le C_qN^{q/6}\prod_j\|T_j\|.
 \label{W:eq:gaussian-matrix-bound}
\end{equation}
\end{corollary}

On the truncated domain $\max|H_{ij}|\le1$, these values and their first five matrix entry derivatives can be approximated to any fixed polynomial accuracy by choosing a
sufficiently large power $J=N^B$.

\begin{proof}
The identity $E\pm\mathrm i\omega_n-A=m_N^2(H-z_n^\pm)$ proves \eqref{W:eq:gaussian-pole-domain}. For $r>1$, integral comparison gives
\[
 \sum_{n\ge J}\omega_n^{-r}
 \le\pi^{-r}\sum_{n\ge J}(n+1)^{-r}
 \le C_rJ^{1-r}.
\]
The paired zero order term has norm at most $2(1+\|A-EI\|)\omega_n^{-2}$, so its sum proves \eqref{W:eq:gaussian-pole-tail-zero}. At order $q\ge1$, the two types of products have
norms at most $q!\omega_n^{-q-1}\prod_j\|V_j\|$ and $(q+1)!\omega_n^{-q-2}\prod_j\|V_j\|$. Their sums are bounded by $C_qJ^{-q}\prod_j\|V_j\|$, which proves
\eqref{W:eq:gaussian-pole-tail-positive}. The Schatten version follows by assigning that norm to the chosen direction in every product.

For $A(H)=m_N^2(2I-H)$, the identities $D_HA[W_j]=-m_N^2W_j$ and $D_H^kA=0$ for $k\ge2$ leave only the partition into singleton sets in the composition formula. Thus
\[
 D_H^q[\theta_E(A(H))][W_1,\ldots,W_q]
 =D_A^q\theta_E(A)[-m_N^2W_1,\ldots,-m_N^2W_q].
\]
Since $W_j=N^{-1/2}T_j$, each matrix $-m_N^2W_j$ has norm $m_N^2N^{-1/2}\|T_j\|=N^{1/6}\|T_j\|$. This proves \eqref{W:eq:gaussian-matrix-bound}. On $\max_{ij}|H_{ij}|\le1$,
$\|A-EI\|\le m_N^2(2+N)+|E|\le C_E N^{5/3}$. The tails, after taking an unnormalized trace against a contraction, are consequently bounded by
\[
 C_E N^{8/3}J^{-1}\quad(q=0),\qquad
 C_qN^{1+q/6}J^{-q}\prod_j\|T_j\|\quad(1\le q\le5).
\]
If an additional differentiated factor has norm at most $N^{B_0}$, choose $J=\lceil N^B\rceil$ with $B>A+B_0+4$. For $q=0$ the resulting power of $N$ is $B_0+8/3-B<-A-4/3$. For
$1\le q\le5$, it is
\[
 B_0+1+q/6-qB
 \le B_0+1+q/6-q(A+B_0+4)<-A,
\]
where $A,B_0\ge0$ and $q\ge1$. The finite number of orders and Leibniz terms can be included by their fixed constants. Thus every displayed tail is $O(N^{-A})$. The exact positive
function $\theta_E$ continues to define the localization statistic.
\end{proof}

Lemma~\ref{W:lem:gaussian-trace} estimates $\operatorname{Tr}\theta_E(A)$. For products with two spectral weights we use its divided difference, as in the next calculation.

For the rest of this subsection $\theta_E^{\langle q\rangle}$ denotes the Hermite divided difference of order $q$; in particular
$\theta_E^{\langle q\rangle}(x,\ldots,x)=\theta_E^{(q)}(x)/q!$. This distinct notation avoids confusing it with a pole partial sum. Put $w_E(x)=\exp[-2(x-E)_+]$.

\Needspace{6\baselineskip} Estimating the combined divided difference kernel gives decay in the spectral variables, also when nodes coincide or tend to negative infinity.

\begin{theorem}
\label{W:thm:gaussian-divided-differences}
For every fixed $q\ge0$ and all real $x_0,\ldots,x_q$,
\begin{equation}
 |\theta_E^{\langle q\rangle}(x_0,\ldots,x_q)|
 \le C_q\sum_{j=0}^q w_E(x_j)
          \prod_{k\ne j}(1+|x_k-x_j|)^{-1}.
 \label{W:eq:gaussian-gram-cluster-bound}
\end{equation}
The estimate is uniform at repeated nodes and arbitrarily negative nodes. If $x_1,\ldots,x_n$ are the scaled eigenvalues from Lemma~\ref{W:lem:gaussian-trace}, then
\begin{equation}
 \sum_{i_0,\ldots,i_q=1}^{n}
 |\theta_E^{\langle q\rangle}(x_{i_0},\ldots,x_{i_q})|^2
 \le C_{q,E}(K+1)^{q+1}.
 \label{W:eq:gaussian-kernel-l2}
\end{equation}
Thus the sum of the squared absolute values of the divided differences is subpolynomial on the event where the eigenvalue counting estimate holds.
\end{theorem}
\begin{proof}
For $q=0$, $0<\theta_E(x)\le w_E(x)$ proves the first bound. The weighted sum estimate below, with $b=0$, proves the second. In the contour argument suppose $q\ge1$. Sort the
distinct nodes and join consecutive nodes when their gap is at most one. Each resulting cluster has diameter at most $q$; different clusters have gap greater than one. Surround a
cluster by the rectangle with imaginary sides at $\pm1/4$ and vertical sides one quarter beyond its extreme nodes. These rectangles are disjoint and stay in the holomorphy strip of
$\theta_E$. Their lengths are bounded in terms of $q$. For a cluster with extreme nodes $a,b$, the rectangle has width $b-a+1/2$, height $1/2$, and perimeter $2(b-a)+2\le2q+2$.
Every node in that cluster is at distance at least $1/4$ from its boundary. For an outside node $x_k$, let $d_k=\operatorname{dist}(x_k,[a,b])>1$. For every $z$ on the rectangle
and any node $x_j$ in that cluster,
\[
 |z-x_k|\ge d_k-1/4\ge3d_k/4,
 \qquad 1+|x_k-x_j|\le d_k+q+1\le(q+2)d_k.
\]
Thus one may take $c_q=3/[4(q+2)]$ in $|z-x_k|\ge c_q(1+|x_k-x_j|)$. The inequality $|1+re^{\mathrm iy}|\ge(1+r)\cos(1/8)$, $|y|\le1/4$, gives
\[
 |\theta_E(z)|\le\sec^2(1/8)e^{-2(\operatorname{Re}z-E)_+}
 \le\sec^2(1/8)e^{2q+1/2}w_E(x_j).
\]
The last step uses the Lipschitz bound for $x\mapsto x_+$ and $|\operatorname{Re}z-x_j|\le q+1/4$. Cauchy's divided difference formula is
\[
 \theta_E^{\langle q\rangle}(x_0,\ldots,x_q)
 =\frac1{2\pi\mathrm i}\sum_{\text{clusters}}
 \int_{\partial\mathcal R}
 \frac{\theta_E(z)}{\prod_{k=0}^q(z-x_k)}\,\mathrm dz.
\]
For distinct nodes, the residue at $x_j$ is $\theta_E(x_j)/\prod_{k\ne j}(x_j-x_k)$, and their sum is the divided difference. For repeated nodes, perturb them within the same fixed
rectangles. The integrands converge uniformly on the contours, which remain separated from every node; continuity gives the same formula, with each node contributing its
multiplicity to the denominator. Apply the preceding contour integral estimate to each term in the finite sum. For a cluster $C$ containing $j$, there are at most $q$ other nodes
and $1+|x_k-x_j|\le1+q$ for $k\in C$. Hence
\[
 \prod_{k\notin C}(1+|x_k-x_j|)^{-1}
 \le (1+q)^q\prod_{k\ne j}(1+|x_k-x_j|)^{-1}.
\]
This accounts for all factors from nodes in the same cluster with a constant depending only on $q$. Choosing one node from each cluster and then enlarging the sum to all nodes
proves \eqref{W:eq:gaussian-gram-cluster-bound}. The scalar nodes are fixed while these contours are chosen.

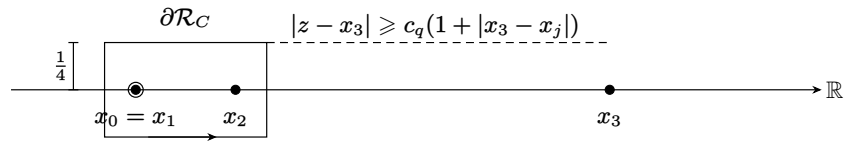
\begingroup \setlength{\intextsep}{1.2\baselineskip}
\begin{figure}[H]
\centering
\begin{tikzpicture}[x=1.65cm,y=2.5cm,every node/.style={mg text}]
 \draw[mg arrow] (0,0)--(6.5,0) node[right] {$\mathbb R$};
 \draw[mg kernel] (.75,-.25) rectangle (2.05,.25);
 \draw[mg arrow] (1.1,-.25)--(1.65,-.25);
 \node[mg occurrence] (a) at (1,0) {};
 \draw[mg kernel] (a) circle[radius=2.8pt];
 \node[mg occurrence] (b) at (1.8,0) {};
 \node[mg occurrence] (c) at (4.8,0) {};
 \node[below=7pt] at (1,0) {$x_0=x_1$};
 \node[below=7pt] at (1.8,0) {$x_2$};
 \node[below=7pt] at (4.8,0) {$x_3$};
 \node[above=4pt] at (1.4,.25) {$\partial\mathcal R_C$};
 \draw[mg incidence] (2.05,.25)--(4.8,.25);
 \node[mg label,above] at (3.4,.25)
 {$|z-x_3|\ge c_q(1+|x_3-x_j|)$};
 \draw[mg kernel] (.50,0)--(.50,.25);
 \draw[mg kernel] (.46,0)--(.54,0);
 \draw[mg kernel] (.46,.25)--(.54,.25);
 \node[left] at (.50,.125) {$\frac14$};
 \node[align=center] at (3.2,-.63)
 {$C=\{x_0,x_1,x_2\}$,\quad $x_j\in C$,\quad
   $\prod_{k=0}^3(z-x_k)=(z-x_0)^2(z-x_2)(z-x_3)$};
\end{tikzpicture}
\caption{The positively oriented contour $\partial\mathcal R_C$ in the proof of \eqref{W:eq:gaussian-gram-cluster-bound}. Filled dots denote distinct node locations. The ring at $x_0=x_1$ indicates the repeated node, whose factor in the denominator is squared. Successive distinct nodes of $C$ are at distance at most one, and the next cluster begins more than one unit away. The rectangle is at least $1/4$ from every node of $C$. The dashed guide gives the estimate for a node outside $C$ relative to the chosen $x_j\in C$.}
\label{W:fig:divided-difference-clusters}
\end{figure}
\endgroup

For the sum of squared divided differences, dyadic annuli and \eqref{W:eq:gaussian-local-count} give, for every real $x$,
\begin{equation}
 \sum_i(1+|x_i-x|)^{-2}
 \le C(K+1+\sqrt{x_+}).
 \label{W:eq:gaussian-square-convolution}
\end{equation}

For $A_0=\{i:|x_i-x|\le1\}$ and $A_r=\{i:2^{r-1}<|x_i-x|\le2^r\}$, $r\ge1$,
\[
 \begin{aligned}
 \sum_i(1+|x_i-x|)^{-2}
 &\le\#A_0+\sum_{r\ge1}2^{-2r+2}
                 \#\{i:|x_i-x|\le2^r\}\\
 &\le C\left(K+1+\sqrt{x_+}
       +\sum_{r\ge1}(K2^{-2r}+\sqrt{x_+}2^{-r}+2^{-r/2})\right)\\
 &\le C(K+1+\sqrt{x_+}).
 \end{aligned}
\]
Also, for every fixed $b\ge0$,
\begin{equation}
 \sum_i w_E(x_i)^2(1+(x_i)_+)^b\le C_{b,E}(K+1).
 \label{W:eq:gaussian-weighted-gram-cluster-sum}
\end{equation}

There are at most $K$ negative eigenvalues. The remaining indices give
\[
 \begin{aligned}
 \sum_{x_i\ge0} w_E(x_i)^2(1+x_i)^b
 &\le\sum_{l=0}^\infty
       e^{-4(l-E)_+}(l+2)^b\#\{i:l\le x_i<l+1\}\\
 &\le C\sum_{l=0}^\infty
       e^{-4(l-E)_+}(l+2)^b(K+1+\sqrt{l+1})
 \le C_{b,E}(K+1).
 \end{aligned}
\]
This proves \eqref{W:eq:gaussian-weighted-gram-cluster-sum}. Square \eqref{W:eq:gaussian-gram-cluster-bound}, use Cauchy--Schwarz on its $q+1$ terms, and sum the $q$ indices other
than the fixed one first. The result is bounded by
\[
 C_q\sum_i w_E(x_i)^2
       \left[\sum_j(1+|x_j-x_i|)^{-2}\right]^q.
\]
For each fixed index $i$, \eqref{W:eq:gaussian-square-convolution} and $K+1+\sqrt{x_+}\le2(K+1)\sqrt{1+x_+}$ give
\[
 \left[\sum_j(1+|x_j-x_i|)^{-2}\right]^q
 \le C_q(K+1)^q(1+(x_i)_+)^{q/2}.
\]
Apply \eqref{W:eq:gaussian-weighted-gram-cluster-sum} with $b=q/2$ to sum over the remaining index:
\[
 C_q(K+1)^q\sum_iw_E(x_i)^2(1+(x_i)_+)^{q/2}
 \le C_{q,E}(K+1)^{q+1}.
\]
This proves the bound on the sum of squares, including repeated spectral indices in its unrestricted sums.
\end{proof}

If $A=\sum_i x_iP_i$ is any spectral resolution, then
\begin{equation}
 D_A^q\theta_E(A)[V_1,\ldots,V_q] =\sum_{\sigma\in S_q}\sum_{i_0,\ldots,i_q} \theta_E^{\langle q\rangle}(x_{i_0},\ldots,x_{i_q}) \times P_{i_0}V_{\sigma(1)}P_{i_1}\cdots V_{\sigma(q)}P_{i_q}.
 \label{W:eq:gaussian-spectral-derivatives}
\end{equation}
This follows by expanding the resolvents in Cauchy's derivative formula and evaluating the scalar residues. Eigenvalues need not be simple; one may either use their spectral
projections or repeat an orthonormal eigenbasis. Equation~\eqref{W:eq:gaussian-kernel-l2} bounds the squared sum of these scalar coefficients; the matrix factors
$P_{i_0}V_{\sigma(1)}P_{i_1}\cdots V_{\sigma(q)}P_{i_q}$ are estimated in each application.

\Needspace{6\baselineskip} Taking absolute values of the individual pole terms too early would destroy the summability provided by the combined kernel.

\begin{remark}
\label{W:rem:absolute-pole-summation}
For $s>1$,
\[
 \sum_{n\ge0}\bigl((x-E)^2+\omega_n^2\bigr)^{-s/2}
 \asymp_s(1+|x-E|)^{1-s}.
\]
Put $r=|x-E|$ and $L=1+r$. Since $\omega_n=\pi(2n+1)$, there are constants $c_s,C_s>0$ such that
\[
 c_s(L+n)^{-s}\le(r^2+\omega_n^2)^{-s/2}
                  \le C_s(L+n)^{-s}.
\]
The decreasing function $t\mapsto(L+t)^{-s}$ satisfies
\[
 \frac{L^{1-s}}{s-1}
 \le\sum_{n=0}^\infty(L+n)^{-s}
 \le L^{-s}+\frac{L^{1-s}}{s-1}
 \le\left(1+\frac1{s-1}\right)L^{1-s}.
\]
This proves both inequalities uniformly for $r\ge0$.

For example, suppose a deterministic increasing sequence has $c(1+j^{2/3})\le1+|x_j-E|\le C(1+j^{2/3})$ for $1\le j\le N$. At $s=2$ the double sum then obeys
\[
 \sum_{j=1}^N\sum_{n=0}^\infty
       \bigl((x_j-E)^2+\omega_n^2\bigr)^{-1}
 \asymp\sum_{j=1}^N(1+j^{2/3})^{-1}\asymp N^{1/3}.
\]
The last comparison follows by integrating $t^{-2/3}$ over $[1,N]$. For this deterministic sequence, the sum of the absolute values of the individual pole contributions therefore
has an additional factor of order $N^{1/3}$. In Theorem~\ref{W:thm:gaussian-divided-differences}, the combined residues instead give the kernel bound
\eqref{W:eq:gaussian-gram-cluster-bound}, whose squared sum was estimated directly above. Thus that proof uses the decay of the combined scalar function before summing over
spectral indices.
\end{remark}

\Needspace{6\baselineskip} The polynomial identity \eqref{W:eq:chebyshev-matrix-derivatives} also applies to arbitrary finite matrices. We record it in the notation used for the
trace contours.

\begin{proposition}
\label{W:prop:gaussian-chebyshev}
Write $P_n(H)=U_n(H/2)$, with $P_0=I$, $P_1=H$ and $P_{n+1}=HP_n-P_{n-1}$. For every finite matrix $H$, every collection of fixed matrices $W_1,\ldots,W_q$, and $n\ge q$,
\begin{equation}
 D_{H}^qP_n(H)[W_1,\ldots,W_q]
 =\sum_{\sigma\in S_q}
 \sum_{j_0+\cdots+j_q=n-q}
 P_{j_0}(H)W_{\sigma(1)}P_{j_1}(H)\cdots
 W_{\sigma(q)}P_{j_q}(H).
 \label{W:eq:gaussian-chebyshev-derivatives}
\end{equation}
For $n<q$ the derivative is zero. The number of terms is $q!\binom nq$. 
\end{proposition}
\begin{proof}
Work first with formal matrix power series in a scalar variable $t$. The recurrence for $P_n$ gives
\[
 S_H(t):=\sum_{n\ge0}P_n(H)t^n,
 \qquad (I-tH+t^2I)S_H(t)=I,
 \qquad S_H(t)=(I-tH+t^2I)^{-1}.
\]
The constant term of the matrix being inverted is $I$, so its formal inverse exists. Every coefficient is a polynomial in the entries of $H$, and differentiation is defined
coefficientwise. Differentiating the inverse identity in a direction $W$ gives
\[
 D_HS_H(t)[W]=tS_H(t)WS_H(t).
\]
Induction in $q$ now yields
\[
 D_H^qS_H(t)[W_1,\ldots,W_q]
 =t^q\sum_{\sigma\in S_q}
 S_H(t)W_{\sigma(1)}S_H(t)\cdots W_{\sigma(q)}S_H(t).
\]
Indeed, differentiating one of the $q+1$ factors $S_H$ in a term of order $q$ places $W_{q+1}$ in one of the $q+1$ positions among the other direction matrices. Every permutation
of $q+1$ directions is obtained exactly once.

Expanding the $q+1$ factors $S_H(t)$ and comparing coefficients of $t^n$ proves \eqref{W:eq:gaussian-chebyshev-derivatives}. When $n\ge q$, the nonnegative indices satisfy
$j_0+\cdots+j_q=n-q$. There are
\[
 \binom{(n-q)+(q+1)-1}{(q+1)-1}=\binom nq
\]
such tuples, and $q!$ permutations of the direction matrices. This is the stated number of terms, counted with their multiplicities. For $n<q$, the factor $t^q$ makes the
coefficient zero. The coefficientwise calculation proves the identity for arbitrary finite matrices, without a spectral assumption on $H$.
\end{proof}

This proposition and \eqref{W:eq:gaussian-polynomial-contour} give two exact representations of the original matrix whose columns are polynomial functions of the matrix applied to
the starting vector, over the whole spectrum. The following calculation shows how their bounds depend on the polynomial degree.

At the scalar argument $H=2+\delta$,
\[
 P_n(2+\delta)
 =\frac{\sinh((n+1)s)}{\sinh s},\qquad
 2+\delta=2\cosh s.
\]
For $0<\delta\le1$, the equation $2\cosh s=2+\delta$ implies $c\sqrt\delta\le s\le\sqrt\delta$. Indeed, $\cosh s\ge1+s^2/2$ gives the upper bound, and $\cosh s-1\le Cs^2$ on the
resulting bounded interval gives the lower bound. With $n\asymp N^{1/3}$ and $\delta=N^{-1/4}$, therefore, $ns\asymp N^{5/24}$. For $(n+1)s\ge1$ and $s\le1$,
\[
 P_n(2+\delta)
 =\frac{\sinh((n+1)s)}{\sinh s}
 \ge \frac{1-e^{-2}}{4s}\,e^{(n+1)s}
 \ge c\exp(c'N^{5/24}).
\]
Its logarithm divided by $\log N$ tends to infinity.

To compare this with the counting estimate, let $D_N$ be diagonal with entries the semicircle quantiles of orders $(j-1/2)/N$, $1\le j\le N$. Then
$\sup_x|F_{D_N}(x)-F_{\rm sc}(x)|\le1/(2N)$. Replace its largest entry by $2+N^{-1/4}$ and call the new matrix $\widetilde D_N$. Exactly one eigenvalue is changed, so
\[
 \sup_x|F_{\widetilde D_N}(x)-F_{\rm sc}(x)|\le\frac{3}{2N},
 \qquad
 \|P_n(\widetilde D_N)\|\ge P_n(2+N^{-1/4}).
\]
Choosing the starting vector to be the corresponding coordinate vector also gives $\|P_n(\widetilde D_N)w\|=P_n(2+N^{-1/4})$. Thus the counting error is $O(N^{-1})$, while the
polynomial vector grows faster than every fixed power of $N$.

\subsubsection{Spectral confinement for the normalized matrix}
\label{W:sec:removal-of-spectral-regularization}

The example above motivates a separate estimate for the extreme eigenvalues. A sufficiently narrow spectral interval both bounds the polynomials whose degrees increase with the
matrix size and ensures agreement with the scalar modification of Section~\ref{W:sec:spectral-regularization-and-matrix-derivatives}.

The argument starts from an auxiliary ensemble with uniformly bounded standardized entries. A scalar trace comparison transfers its spectral bound to the normalized truncated
matrix. The exceptional entries are then added back using the resolvent compressed to their coordinate support and the independence built into the support construction.

For the auxiliary ensemble $W^{\rm b}$, the eigenvalue deviation theorem of \cite{erdosxu2023small} gives, for every fixed $D>0$, a constant $M_D$ such that

\begin{equation}
 \mathbb P\bigl(\|W^{\mathrm b}\|>
       2+M_DN^{-2/3}(\log N)^{2/3}\bigr)\le N^{-D}.
 \label{W:eq:bounded-entry-spectral-deviation}
\end{equation}
Both edges are covered by applying the right tail inequality to $W^{\mathrm b}$ and $-W^{\mathrm b}$. The hypotheses are independent centered Hermitian blocks, off-diagonal
variance $N^{-1}$, zero complex pseudovariance, a fixed bounded diagonal variance, and uniformly bounded standardized moments of every order. We construct a matrix with bounded
standardized entries and diagonal variance $N^{-1}$, so these hypotheses hold in both symmetry classes. We apply that theorem to this matrix and compare it with the original
truncated matrix by replacing entries. The other probabilistic estimates are the preliminary entrywise and counting results of
Section~\ref{W:sec:eigenvalue-counting-and-minor-estimates}, and matrix Bernstein's inequality in \cite{tropp2012userfriendly}. The applications of those estimates are given below.

Let ${m_N}=N^{1/3}$, ${\eta_{\mathrm e}}=N^{-2/3}$, and let $\rho_*=\rho_*(N)$ be the largest Gram matrix dimension needed for a fixed collection of functions. Fix
$0<c_*\le C_*<\infty$ and $A_0>0$, and suppose
\[
 c_*{m_N}\le \rho_*\le C_*{m_N},\qquad
 d=d_N^{(2)}=A_0\frac{(\log N)^2}{\rho_*^2},\qquad
 \Omega_N=\{\|H_N\|\le2+d/2\}.
\]
In particular $d/{\eta_{\mathrm e}}$ is bounded above and below by fixed positive multiples of $(\log N)^2$. Use the tail functions $\Theta$ and $\Theta_{\mathrm d}$ in
\eqref{W:eq:offdiagonal-tail-function} and \eqref{W:eq:diagonal-tail-function}. Resolvents have sign $R_H(z)=(zI-H)^{-1}$ and
$g(z)=\int(z-x)^{-1}\rho_{\mathrm{sc}}(x)\,\mathrm dx$.

\Needspace{6\baselineskip}

\begin{theorem}
\label{W:thm:comparison-window}
Under Assumption~\ref{W:ass:wigner-lee-yin}, for both $\beta=1,2$,
\begin{equation}
 \mathbb P(\Omega_N^c)\longrightarrow0.
 \label{W:eq:comparison-window}
\end{equation}
More precisely, there are constants $c,C,c_0>0$, depending on the fixed parameters and moment bounds but not on $N$ or a starting vector, such that, for all sufficiently large $N$,
\begin{equation}
 \mathbb P(\Omega_N^c)
 \le C\bigl(N^{-c}+\Theta(c_0\sqrt N)+\Theta_{\mathrm d}(c_0\sqrt N)\bigr).
 \label{W:eq:comparison-window-rate}
\end{equation}
For the GOE or GUE matrix ${H_N^{\mathrm G}}$,
\begin{equation}
 \mathbb P(\|{H_N^{\mathrm G}}\|>2+d/2)
 \le C\exp[-c(\log N)^3].
 \label{W:eq:comparison-Gaussian-window}
\end{equation}
The first probability estimate concerns the original and truncated matrices, with an $o(1)$ failure probability. Uniform control along entry replacements is proved separately
below.
\end{theorem}

Set, for $2<a<3$,
\[
 f_{a,{\eta_{\mathrm e}}}(x)=\frac1{1+\exp(-(x-a)/{\eta_{\mathrm e}})},\qquad
 \mathcal F_{a,{\eta_{\mathrm e}}}(B)=\operatorname{Tr}f_{a,{\eta_{\mathrm e}}}(B).
\]
This is a nonnegative function bounded by $N$, and $\mathcal F_{a,{\eta_{\mathrm e}}}(B)\ge1/2$ whenever $\lambda_{\max}(B)\ge a$. For a contour integral below, take the
counterclockwise rectangle with vertices $\pm5\pm\mathrm i {\eta_{\mathrm e}}$. The poles of $f_{a,{\eta_{\mathrm e}}}$ have imaginary parts $(2k+1)\pi {\eta_{\mathrm e}}$, so no
pole lies inside this rectangle.

\Needspace{6\baselineskip}

\begin{lemma}
\label{W:lem:comparison-exterior-derivatives}
Suppose $B=B^*$, $\|B\|\le4$, and
\begin{equation}
 \max_{i,j}|R_{B,ij}(x\pm\mathrm i{\eta_{\mathrm e}})-g(x\pm\mathrm i{\eta_{\mathrm e}})\delta_{ij}|
 \le r,\qquad -5\le x\le5,
 \label{W:eq:comparison-weak-event}
\end{equation}
where $0<r\le1$. Let $T_1,\ldots,T_q$ be Hermitian matrices supported in the same principal block of size at most two, with $\|T_j\|\le1$. For $1\le q\le4$,
\begin{equation}
 |D_{H}^q\mathcal F_{a,{\eta_{\mathrm e}}}(B)[T_1,\ldots,T_q]|
 \le C_q\bigl(r+\exp(-(a-2)/{\eta_{\mathrm e}})+\exp(-1/{\eta_{\mathrm e}})\bigr).
 \label{W:eq:comparison-exterior-derivative-bound}
\end{equation}
The derivatives are with respect to $B+\sum_jt_jT_j$. Substituting $N^{-1/2}T_j$ in every argument multiplies the derivative by $N^{-q/2}$.
\end{lemma}
\begin{proof}
The rectangle with horizontal sides at heights $\pm\eta_{\rm e}$ and vertical sides at $\pm5$ surrounds $\sigma(B)$, since $\|B\|\le4$. The nearest poles of $f_{a,\eta_{\rm e}}$
have imaginary parts $\pm\pi\eta_{\rm e}$, so this function is holomorphic on and inside the rectangle. Equation~\eqref{W:eq:resolvent-repeated-trace-ibp}, multiplied by $q!$,
gives
\begin{equation}
 D_{H}^q\mathcal F_{a,{\eta_{\mathrm e}}}(B)[T^q]
 =\frac{(q-1)!}{2\pi\mathrm i}
   \int f_{a,{\eta_{\mathrm e}}}'(z)\operatorname{Tr}(TR_B(z))^q\,\mathrm dz.
 \label{W:eq:comparison-cyclic-trace}
\end{equation}
For the fixed matrices $T_1,\ldots,T_q$, equation~\eqref{W:eq:resolvent-trace-ibp} on the same contour gives
\begin{equation}
 D_{H}^q\mathcal F_{a,{\eta_{\mathrm e}}}(B)[T_1,\ldots,T_q]
 =\frac1{2\pi\mathrm i q}\int f_{a,{\eta_{\mathrm e}}}'(z)
   \sum_{\sigma\in S_q}\operatorname{Tr}
      [T_{\sigma(1)}R_B(z)\cdots T_{\sigma(q)}R_B(z)]\,\mathrm dz.
 \label{W:eq:comparison-mixed-trace}
\end{equation}
Only the two-dimensional resolvent compression enters each trace. Let $C(z)=R_B(z)_{II}$ be that compression. For one labelled order, put $A_j=T_{\sigma(j)}C(z)$ and
$B_j=g(z)T_{\sigma(j)}$ on $\mathbb C^I$. The exact product difference is
\[
 A_1\cdots A_q-B_1\cdots B_q
 =\sum_{\ell=1}^q B_1\cdots B_{\ell-1}(A_\ell-B_\ell)
                                      A_{\ell+1}\cdots A_q.
\]
On a horizontal side, $\|C-gI\|\le2r$, $\|C\|\le C$, and $|I|\le2$. Hence the absolute trace of the difference is bounded by
\[
 \sum_{\ell=1}^q |I|\,|g|^{\ell-1}\|C-gI\|\|C\|^{q-\ell}
 \le C_qr.
\]
Furthermore
\[
 |f_{a,{\eta_{\mathrm e}}}'(x\pm\mathrm i{\eta_{\mathrm e}})|
 \le C{\eta_{\mathrm e}}^{-1}\exp(-|x-a|/{\eta_{\mathrm e}}),\qquad
 \int_{-5}^5|f_{a,{\eta_{\mathrm e}}}'(x\pm\mathrm i{\eta_{\mathrm e}})|\,\mathrm dx\le C.
\]
On the vertical sides, $\|R_B(z)\|\le1$ and $|g(z)|\le1$; their contribution to the same replacement error is at most $C_q\exp(-1/{\eta_{\mathrm e}})$. Thus the total error is
$C_q(r+\exp(-1/{\eta_{\mathrm e}}))$.

For the replaced expression the remaining scalar integral is $\int f_{a,{\eta_{\mathrm e}}}'(z)g(z)^q\,\mathrm dz$, multiplied by a bounded trace of products of the $T_j$. The
function $g$ is analytic off $[-2,2]$, has boundary values of modulus at most one, and has no other singularity in the rectangle. Deforming the contour to the two sides of this cut
therefore gives
\[
 \left|\frac1{2\pi\mathrm i}\int f_{a,{\eta_{\mathrm e}}}'(z)g(z)^q\,\mathrm dz\right|
 \le\frac1\pi\int_{-2}^2 f_{a,{\eta_{\mathrm e}}}'(x)\,\mathrm dx
 \le C\exp(-(a-2)/{\eta_{\mathrm e}}).
\]
For fixed $N,a,\eta_{\rm e}$, both $g$ and $f_{a,\eta_{\rm e}}'$ are bounded near the cut boundary points. An omitted arc of radius $\delta$ therefore contributes at most
$C_{N,a,q}\delta\to0$. On the cut the two boundary values of $g^q$ have difference of modulus at most two. Finally
\[
 \int_{-2}^2 f_{a,\eta_{\rm e}}'(x)\,\dd x
 =f_{a,\eta_{\rm e}}(2)-f_{a,\eta_{\rm e}}(-2)
 \le f_{a,\eta_{\rm e}}(2)
 \le e^{-(a-2)/\eta_{\rm e}}.
\]  These bounds
prove \eqref{W:eq:comparison-exterior-derivative-bound}. \end{proof}

We first establish confinement for the normalized matrix. Let $W$ be a matrix with bounded entries such that its upper-triangular entries are independent and centered, and for
fixed $C_0,\phi>0$ it satisfies
\begin{equation}
 \begin{split}
 &\mathbb E|W_{ij}|^2=N^{-1},\qquad
   \mathbb EW_{ij}^2=0\quad(\beta=2)\qquad(i<j),\\
 &\mathbb E|W_{ij}|^3\le C_0N^{-3/2},\qquad
   \mathbb E|W_{ij}|^4\le C_0N^{-2}\log(eN)\quad(i<j),\\
 &\mathbb EW_{ii}^2\le C_0/N,\qquad
   |W_{ij}|\le C_0N^{-\phi}\quad\text{almost surely}.
 \end{split}
 \label{W:eq:comparison-core-class}
\end{equation}
Only the diagonal second moments in this display will enter the comparison.

\Needspace{6\baselineskip}

\begin{proposition}
\label{W:prop:comparison-core-window}
For every fixed $b_1,b_2>0$, a matrix satisfying the stated entry support bounds satisfies
\begin{equation}
 \mathbb P(\|W\|>2+b_N{\eta_{\mathrm e}}(\log N)^2)\le CN^{-c},
 \qquad b_1\le b_N\le b_2,
 \label{W:eq:comparison-core-window}
\end{equation}
for some $c>0$. Constants are uniform over the class \eqref{W:eq:comparison-core-class}. This conclusion uses \eqref{W:eq:bounded-entry-spectral-deviation} only for the matrix with
bounded standardized entries constructed in the proof.
\end{proposition}
\begin{proof}
For each off-diagonal entry use the bounded moment-matching construction in the proof of Lemma~\ref{W:lem:cubic-defect}. In the real case apply it to $\sqrt N W_{ij}$; in the
complex case apply it to
\[
 X=(\sqrt{2N}\operatorname{Re}W_{ij},
       \sqrt{2N}\operatorname{Im}W_{ij}).
\]
These variables are centered, have identity covariance in their real coordinates, and have uniformly bounded third absolute moment. The construction therefore gives uniformly
bounded real coordinate vectors with the same joint moment tensors through degree three. Reversing the normalization defines $W_{ij}^{\mathrm b}$. Choose different off-diagonal
blocks independently, impose Hermitian symmetry, and take independent symmetric signs divided by $\sqrt N$ on the diagonal. Thus $W^{\mathrm b}$ has bounded standardized entries
and diagonal variance $N^{-1}$; only its off-diagonal tensors are matched with those of $W$.

Take $a=2+b_N{\eta_{\mathrm e}}(\log N)^2$ in Lemma~\ref{W:lem:comparison-exterior-derivatives}. For $\beta=2$, Lemmas~\ref{W:lem:complex-weak} and~\ref{W:lem:complex-norm}, with
for instance $\tau=1/4$, apply to the intermediate matrices of $W$ and $W^{\mathrm b}$. The height ${\eta_{\mathrm e}}$ is larger than $N^{-1+\tau}$. For $\beta=1$, apply
Lemma~3.1(1)--(2) of \cite{lee2014necessary} to these real symmetric intermediate matrices. Their means, variances, and support satisfy that lemma with $q=N^{\phi_0}$, where
$\phi_0=\min\{\phi/2,1/4\}$; its domain contains $\eta={\eta_{\mathrm e}}$ because $N\eta=N^{1/3}$ exceeds every fixed power of $\varphi_N=(\log N)^{\log\log N}$. Its entrywise
error is at most $C\varphi_N^C(N^{-\phi_0}+N^{-1/6}+N^{-1/3})$, hence a fixed negative power of $N$, and its norm bound is less than three for all $N\ge N_0(\phi_0)$. Thus in both
symmetry classes the intermediate resolvent entries differ from $g\delta_{ij}$ by a fixed negative power of $N$, and $\|H^\gamma\|\le3$, with overwhelming probability. Setting one
block to zero changes the norm by $O(N^{-\phi})$. The exact Woodbury formula on its at most two coordinates changes every resolvent entry by $O(N^{-\phi})$, since the coordinate
compressions are bounded. The same is true when any value from either entry distribution is used to replace the removed entry. After decreasing an exponent, there is $c_1>0$ such
that the event measurable with respect to the remaining entries
\[
 I_\alpha=\{\|K_\alpha\|\le7/2,\quad
  \max_{i,j,\,-5\le x\le5}
 |R_{K_\alpha,ij}(x\pm\mathrm i{\eta_{\mathrm e}})-g(x\pm\mathrm i{\eta_{\mathrm e}})\delta_{ij}|
 \le N^{-c_1}\}
\]
has complement $O_D(N^{-D})$ for every $D$. On it, every segment varying the selected entry lies in $\|B\|\le4$ and satisfies \eqref{W:eq:comparison-weak-event} with a possibly
smaller fixed negative power of $N$. Lemma~\ref{W:lem:comparison-exterior-derivatives} bounds its derivatives of orders one through four by
\begin{equation}
 \epsilon_N\le C(N^{-c_2}+\exp[-b_1(\log N)^2]).
 \label{W:eq:comparison-small-derivatives}
\end{equation}
The estimate holds on the specified intermediate matrices and interpolation segments.

For each replacement use the matrix entry coordinate $x=W_{ab}$ in the real case. In the complex case take the vector of real coordinates
\[
 x=(\sqrt2\operatorname{Re}W_{ab},
                 \sqrt2\operatorname{Im}W_{ab}),
\]
with corresponding matrix variations $(E_{ab}+E_{ba})/\sqrt2$ and $\mathrm i(E_{ab}-E_{ba})/\sqrt2$. These coordinates include the matrix normalization, so their second moments are
$O(N^{-1})$. Conditional on the remaining entries, expand $\mathcal F_{a,{\eta_{\mathrm e}}}(K_\alpha+\mathcal T_\alpha x)$ through degree three. If $\Delta\mu_\alpha$ is the
difference of these laws and $F_\alpha(x)=\mathcal F_{a,\eta_{\rm e}}(K_\alpha+\mathcal T_\alpha x)$, then, for $0\le q\le3$,
\[
 \int\frac1{q!}D^qF_\alpha(0)[x^q]\,\Delta\mu_\alpha(\dd x)
 =\frac1{q!}\sum_{r_1,\ldots,r_q}
 (\partial_{r_1}\cdots\partial_{r_q}F_\alpha)(0)
 \int x_{r_1}\cdots x_{r_q}\,\Delta\mu_\alpha(\dd x)=0.
\]
At $q=0$ this is equality of total masses. For $q\ge1$ it uses the matching joint tensors and the fact that all coefficients, including $\mathbf1_{I_\alpha}$, are fixed by
$K_\alpha$. Its exact remainder is
\[
 \frac1{3!}\int_0^1(1-t)^3
  D_{H}^4\mathcal F_{a,{\eta_{\mathrm e}}}(K_\alpha+t\mathcal T_\alpha x)
             [(\mathcal T_\alpha x)^4]\,\mathrm dt.
\]
On $I_\alpha$, its integrated absolute value is at most $C\epsilon_N N^{-2}\log(eN)$, because the fourth moment of the real coordinate vector of the replaced entry is bounded by
the preceding quantity divided by $\epsilon_N$. For a diagonal block, let $x$ be its matrix entry and $T=E_{aa}$. Taylor's formula through degree one has remainder
\[
 \int_0^1(1-t)
 D_H^2\mathcal F_{a,{\eta_{\rm e}}}(K_\alpha+txT)[xT,xT] \,\dd t.
\]
Its absolute value on $I_\alpha$ is at most $C\epsilon_N |x|^2/2$. The two centered diagonal laws have second moments bounded by $C/N$, so the constant and linear terms cancel and
the integrated absolute remainder is at most $C\epsilon_N/N$. On $I_\alpha^c$ the difference between the original function values is bounded by $2N$. Independence of the
replacement entry from $I_\alpha$ justifies both cancellations. The fourth absolute moments sum to at most $C N^2 N^{-2}\log(eN)=C\log(eN)$, and the diagonal second moments sum to
at most $C N N^{-1}=C$. If the exceptional probability at each replacement is at most $N^{-D-4}$, its contribution to the expectation of the bounded function is at most
$2N\cdot N^2N^{-D-4}=2N^{-D-1}$. Consequently telescoping the differences of expectations yields
\begin{equation}
 |\mathbb E\mathcal F_{a,{\eta_{\mathrm e}}}(W)
       -\mathbb E\mathcal F_{a,{\eta_{\mathrm e}}}(W^{\mathrm b})|
 \le C\epsilon_N\log(eN)+O_D(N^{-D}).
 \label{W:eq:comparison-core-comparison}
\end{equation}
For $W^{\mathrm b}$ and any fixed $D$, the threshold $M_D(\log N)^{2/3}$ in \eqref{W:eq:bounded-entry-spectral-deviation} is less than $(a-2)/(2{\eta_{\mathrm e}})$ for all
$N\ge N_0(D)$, since $(a-2)/\eta_{\mathrm e}$ grows faster than $(\log N)^{2/3}$. On the spectral confinement event for $W^{\mathrm b}$, every eigenvalue satisfies
$f_{a,{\eta_{\mathrm e}}}(\lambda)\le\exp(-(a-2)/(2{\eta_{\mathrm e}}))$. Therefore
\[
 \mathbb E\mathcal F_{a,{\eta_{\mathrm e}}}(W^{\mathrm b})
 \le N\exp[-b_1(\log N)^2/2]+N\,O_D(N^{-D})
 =O_D(N^{-D})
\]
by choosing the tail exponent one order larger. The external deviation estimate is used at $M_D(\log N)^{2/3}$, within its stated range; the larger enclosure follows from
$M_D(\log N)^{2/3}=o((\log N)^2)$. Since $\mathcal F_{a,{\eta_{\mathrm e}}}\ge\tfrac12\mathbf1_{\{\lambda_{\max}\ge a\}}$, the comparison estimate
\eqref{W:eq:comparison-core-comparison} bounds the probability at the upper edge. The logarithmic factors are absorbed using
$(\log N)^a e^{-c(\log N)^\gamma} \le C_{a,c}e^{-(c/2)(\log N)^\gamma}$ for any fixed $a$ and sufficiently large $N$. The matrix $-W$ satisfies the same hypotheses, and
$\lambda_{\max}(-W)=-\lambda_{\min}(W)$ gives the lower edge bound in \eqref{W:eq:comparison-core-window}.
\end{proof}

\subsubsection{Resolvents on the support of the exceptional entries}
\label{W:sec:resolvents-on-exceptional-support}

To add back the exceptional entries, we estimate the resolvent of the normalized matrix on their coordinate support. The support construction makes that matrix independent of the
labels determining the support. Throughout this proof of spectral confinement, $\phi=1/100$ and $T=N^{1/2-\phi}$ in standardized coordinates.

\begin{lemma}
\label{W:lem:comparison-decomposition}
Under Assumption~\ref{W:ass:wigner-lee-yin}, one can construct
\begin{equation}
 H=W+V+E
 \label{W:eq:comparison-decomposition}
\end{equation}
on a common probability space so that $W$ satisfies \eqref{W:eq:comparison-core-class}, $\|V\|\le N^{-1+4\phi}$ with overwhelming probability, and $E=\iota_SE_S\iota_S^*$ for a
random coordinate set $S$ independent of the entire matrix $W$. For $i\ne j$,
\begin{equation}
 \mathbb P(i\in S)\le CN^{-1+4\phi},\qquad
 \mathbb P(i,j\in S)\le CN^{-2+8\phi}.
 \label{W:eq:comparison-support-probabilities}
\end{equation}
Moreover, for a fixed sufficiently small $c_0>0$,
\begin{equation}
 \mathbb P(\|E\|>1/20)
 \le C\bigl(N^{-1+8\phi}+\Theta(c_0\sqrt N)+\Theta_{\mathrm d}(c_0\sqrt N)\bigr)
 \label{W:eq:comparison-E-small}
\end{equation}
for all sufficiently large $N$. The values of $E_S$ need not be independent of $W$.
\end{lemma}
\begin{proof}
At each independent standardized scalar or complex block generate independently a Bernoulli label $L$ of parameter $p=\mathbb P(|\xi|>T)$, a variable $Y$ with the conditional law
on $\{|\xi|\le T\}$, and a variable $Z$ with the complementary conditional law. The random variable $(1-L)Y+LZ$ has the original law. Use the same construction for every real
diagonal variable, with its own law. Labels with probability zero are simply fixed at zero. Let $B_{ij}=N^{-1/2}Y_{ij}$, and put $E_{ij}=N^{-1/2}L_{ij}(Z_{ij}-Y_{ij})$ with
Hermitian symmetry. Then $H=B+E$ exactly, and all labels are independent of $B$.

The weak fourth moment tail condition implies a global bound $\mathbb P(|\xi|>t)\le C t^{-4}$ for $t\ge1$. Tail integration gives
\[
 \mathbb E[|\xi|^r\mathbf1_{\{|\xi|>T\}}]\le CT^{r-4}
 \quad(1\le r<4),\qquad
 \mathbb E[|\xi|^4\mathbf1_{\{|\xi|\le T\}}]\le C\log(eT).
\]
For the standardized vector of real coordinates of $Y$, let $\mu$ be its mean and $\Sigma$ its covariance. For off-diagonal blocks,
\begin{equation}
 |\mu|\le CT^{-3},\qquad
 \|\Sigma-I\|\le CT^{-2},\qquad
 \mathbb E\|Y\|^3\le C,
 \label{W:eq:comparison-core-moment-corrections}
\end{equation}
To check the real coordinate normalization, take $U=\xi$ in the real case and $U=\sqrt2(\operatorname{Re}\xi,\operatorname{Im}\xi)$ in the complex case. Then $\mathbb EU=0$ and
$\mathbb EUU^{\mathsf T}=I$. For $A=\{|\xi|\le T\}$,
\[
 \mu=-\frac{\mathbb E[U\mathbf1_{A^c}]}{\mathbb P(A)},
 \qquad
 \Sigma-I=\frac{\mathbb P(A^c)I-\mathbb E[UU^{\mathsf T}\mathbf1_{A^c}]}
                    {\mathbb P(A)}-\mu\mu^{\mathsf T}.
\]
For sufficiently large $T$, $\mathbb P(A)\ge1/2$. The bounds $|U|\le\sqrt2|\xi|$ and $\|UU^{\mathsf T}\|=|U|^2$ reduce the first two estimates to the tail integrals at orders one
and two. The third follows from $\mathbb E[|U|^3\mid A]\le2\mathbb E|U|^3<\infty$. Thus all constants are fixed in both dimensions. The normalization $\Sigma^{-1/2}(Y-\mu)$
therefore has covariance $I$ and bounded third moment, logarithmic fourth moment, and support $O(T)$. Define the off-diagonal $W$ by this normalization, divided by $\sqrt N$ (and
using the complex block convention above). Center the conditional diagonal values without variance rescaling. Their second moments remain bounded and $|\mathbb E Y_{ii}|\le C/T$.
This constructs the required matrix with bounded entries. It depends only on the variables conditioned to lie below the truncation threshold and is independent of the labels.

Put $V=B-W$. Its deterministic mean has operator norm at most $C(N^{1/2}T^{-3}+N^{-1/2}T^{-1})\le CN^{-1+3\phi}$, by maximum row sums. Let $\delta=CN^{-1+2\phi}$. The centered
change due to covariance normalization is a sum of independent mean zero Hermitian matrices $V^{(ij)}$, supported on the corresponding entry pairs, satisfying
\[
 \|V^{(ij)}\|\le C\delta N^{-\phi},\qquad
 \left\|\sum_{i<j}\mathbb E(V^{(ij)})^2\right\|\le C\delta^2.
\]
The second bound follows by summing the coordinate variances, each at most $C\delta^2/N$. Matrix Bernstein gives
\[
 \mathbb P\left(\left\|\sum_{i<j}V^{(ij)}\right\|>t\right)
 \le2N\exp\left[-\frac{t^2}{C\delta^2+C\delta N^{-\phi}t}\right].
\]
At $t=\delta\log N$, the Bernstein exponent is at least $c(\log N)^2$, since $N^{-\phi}\log N\to0$. Hence the exceptional probability is smaller than every fixed inverse power of
$N$. The displayed mean bound gives
\[
 \|V\|\le\delta\log N+C N^{-1+3\phi}
 \le N^{-1+4\phi}
\]
for sufficiently large $N$ on the event in the Bernstein estimate. The covariance correction includes the complex pseudovariance.

Let $S$ be the set of vertices incident to a marked off-diagonal edge, together with marked diagonal vertices. The mark probabilities satisfy
\[
 p_{ij}\le CN^{-2+4\phi}\quad(i<j),\qquad
 p_{ii}\le CN^{-1+2\phi}.
\]
For a fixed vertex $i$, a union bound over its incident labels gives
\[
 \mathbb P(i\in S)\le p_{ii}+\sum_{j\ne i}p_{ij}
 \le C N^{-1+2\phi}+C(N-1)N^{-2+4\phi}
 \le C' N^{-1+4\phi}.
\]
The event $\{i\in S\}$ is determined by these labels only, whereas $W$ is a measurable function of the independent small conditional variables. Hence the entire random set $S$ is
independent of $W$. For $i\ne j$, either the common edge $\{i,j\}$ is marked, or both vertices are marked by their separate families of incident edges and diagonal labels. After
removing the common edge, those two families are independent. Hence
\[
 \mathbb P(i,j\in S)
 \le p_{ij}+
 \left(p_{ii}+\sum_{k\ne i,j}p_{ik}\right)
 \left(p_{jj}+\sum_{k\ne i,j}p_{jk}\right)
 \le CN^{-2+8\phi}.
\]
There are at most $N\binom{N-1}{2}$ unordered pairs of distinct off-diagonal edges with a common vertex. For each such pair their independent labels are both one with probability
at most $C^2N^{-4+8\phi}$. Thus
\[
 \mathbb P(\text{two marked off-diagonal edges share a vertex})
 \le N\binom{N-1}{2}C^2N^{-4+8\phi}
 \le C' N^{-1+8\phi}.
\]  Otherwise the marked graph
is a matching, and, allowing any diagonal marks, $E$ is block diagonal on blocks of size at most two. Consequently
\[
 \|E\|\le2\max_{ij}|E_{ij}|
 \le2\left(\max_{ij}|H_{ij}|+N^{-\phi}\right).
\]
For fixed $c_0=1/100$ and sufficiently large $N$, this is less than $1/20$ on $\max_{ij}|H_{ij}|\le c_0$. Finally,
\[
 \mathbb P\left(\max_{i<j}|H_{ij}|>c_0\right)
 \le C \Theta(c_0\sqrt N),\qquad
 \mathbb P\left(\max_i|H_{ii}|>c_0\right)
 \le C \Theta_{\mathrm d}(c_0\sqrt N).
\]
The diagonal inequality uses $\mathbb P(|d|>s)\le s^{-2}\mathbb E[d^2\mathbf1_{\{|d|>s\}}]$. This proves \eqref{W:eq:comparison-E-small}.
\end{proof}

We estimate the resolvent compressed to the support of $E$ at two real spectral parameters.

\Needspace{6\baselineskip}

\begin{lemma}
\label{W:lem:comparison-compression}
Use the coupling of Lemma~\ref{W:lem:comparison-decomposition}, and set $a_+=2+d/4$, $a_-=-2-d/4$. With probability tending to one,
\begin{equation}
 \|(R_W(a_\pm))_{SS}\|\le5.
 \label{W:eq:comparison-compressed-resolvent}
\end{equation}
The event is measurable with respect to $W$ and the support labels.
\end{lemma}
\begin{proof}
By Proposition~\ref{W:prop:comparison-core-window},
\begin{equation}
 \|W\|\le2+d/8
 \label{W:eq:comparison-inner-window}
\end{equation}
except with probability $CN^{-c}$. Intersect it with the event where the preliminary resolvent entry estimates hold at $a_\pm+\mathrm i{\eta_{\mathrm e}}$ and the event defined by

\[
 N\sup_x|F_W(x)-F_{\mathrm{sc}}(x)|\le K_N,
 \qquad K_N=N^{1/12}.
\]
For $\beta=1$, the preliminary resolvent entry estimates follow from Lemma~3.1(1) of \cite{lee2014necessary}, with $q=N^{\phi/2}$ and $\eta=\eta_{\mathrm e}=N^{-2/3}$; the support
bound $C_0N^{-\phi}\le q^{-1}$ and the domain condition $N\eta=N^{1/3}>(\log N)^{C\log\log N}$ hold for all $N\ge N_0(C_0,\phi,C)$. For $\beta=2$, use
Lemma~\ref{W:lem:complex-weak} with $\tau=1/4$. The eigenvalue counting estimate follows from Proposition~\ref{W:prop:minor-real-input} for $\beta=1$ and
Corollary~\ref{W:cor:complex-counting-for-minors} for $\beta=2$: in both cases $\mathfrak c_N=N^{o(1)}\le N^{1/12}$ for sufficiently large $N$. Thus these two events have
overwhelming probability in both symmetry classes. Their intersection with \eqref{W:eq:comparison-inner-window}, denoted $\mathcal G_W$, depends only on $W$.

The real spectral distance from either barrier to $\sigma(W)$ is at least $d/8$, which exceeds ${\eta_{\mathrm e}}$ whenever $d>8\eta_{\mathrm e}$. For the upper barrier, writing
$t=a_+-\lambda\ge d/8$ gives
\[
 \frac1t\le\left(1+\frac{{\eta_{\mathrm e}}^2}{(d/8)^2}\right)
                         \frac{t}{t^2+{\eta_{\mathrm e}}^2}.
\]
Integration against each coordinate spectral measure and the weak entrywise bound at $a_++\mathrm i{\eta_{\mathrm e}}$ yield $0<(R_W(a_+))_{ii}\le4$ for all $N\ge N_0$, uniformly
in $i$. At the lower barrier, substitute $t=\lambda-a_->0$ in the preceding scalar inequality and integrate against the coordinate spectral measures. This bounds $-R_W(a_-)_{ii}$
and gives
\begin{equation}
 \max_i|(R_W(a_\pm))_{ii}|\le4.
 \label{W:eq:comparison-real-diagonal}
\end{equation}

For the upper barrier, extend $(a_+-x)^{-2}$ constantly beyond $b=2+d/8$. This bounded monotone function has total variation $(d/8)^{-2}$ and agrees with $(a_+-x)^{-2}$ on both
relevant spectra. Stieltjes integration by parts and the counting bound give
\[
 \|R_W(a_+)\|_{\mathrm{HS}}^2
 \le N\int_{-2}^2\frac{\rho_{\mathrm{sc}}(x)}{(a_+-x)^2}\,\mathrm dx
       +K_N(d/8)^{-2}
 \le C\bigl(Nd^{-1/2}+K_Nd^{-2}\bigr).
\]
For the semicircle term, $y=2-x$ gives
\[
 \int_{-2}^2\frac{\rho_{\rm sc}(x)}{(a_+-x)^2}\,\dd x
 \le\frac1\pi\int_0^4\frac{\sqrt y}{(d/4+y)^2}\,\dd y
 \le\frac C{\sqrt d}\int_0^\infty\frac{\sqrt t}{(1+t)^2}\,\dd t
 \le Cd^{-1/2}.
\]
Reflection gives the same integral at $a_-=-a_+$, and the lower extension has the same total variation. Thus $\|R_W(a_-)\|_{\rm HS}^2\le C(Nd^{-1/2}+K_Nd^{-2})$ as well. Since
$d\asymp N^{-2/3}(\log N)^2$,
\begin{equation}
 \|R_W(a_\pm)\|_{\mathrm{HS}}^2\le C N^{4/3+1/12}
 \quad\text{on }\mathcal G_W.
 \label{W:eq:comparison-real-HS}
\end{equation}

Define the squared off-diagonal norm of the random compression by
\[
 Y_\pm=\sum_{i\ne j}\mathbf1_{\{i,j\in S\}}
                          |R_{W,ij}(a_\pm)|^2.
\]
Since $\mathcal G_W$ and every resolvent entry below are functions of $W$, independence of $S$ and $W$ gives exactly
\[
 \mathbb E[\mathbf1_{\mathcal G_W}Y_\pm]
 =\sum_{i\ne j}\mathbb P(i,j\in S)
          \mathbb E[\mathbf1_{\mathcal G_W}|R_{W,ij}(a_\pm)|^2].
\]
The support bounds in \eqref{W:eq:comparison-support-probabilities} therefore give
\[
 \mathbb E[\mathbf1_{\mathcal G_W}Y_\pm]
 \le CN^{-2+8\phi}N^{4/3+1/12}
 =CN^{-151/300},\qquad \phi=1/100.
\]
Markov's inequality implies $Y_\pm\le1$ outside an additional event of this probability. Write the compressed matrix as $D_\pm+O_\pm$, with $D_\pm$ diagonal and $O_\pm$ zero on the
diagonal. On the stated event,
\[
 \|D_\pm\|=\max_{i\in S}|R_{W,ii}(a_\pm)|\le4,\qquad
 \|O_\pm\|\le\|O_\pm\|_{\rm HS}=\sqrt{Y_\pm}\le1.
\]
Thus $\|(R_W(a_\pm))_{SS}\|\le5$. If $S=\varnothing$, the compression acts on the zero space and has norm zero. Equation~\eqref{W:eq:comparison-real-diagonal} therefore proves
\eqref{W:eq:comparison-compressed-resolvent}. Conditioning fixes $W$ and integrates the independent support labels.
\end{proof}

\begin{proof}[Proof of Theorem~\ref{W:thm:comparison-window}]
On the events above, $\|E_S\|\le1/20$ and $\|(R_W(a_\pm))_{SS}\|\le5$. For every $0\le t\le1$,
\[
 \|tE_S(R_W(a_\pm))_{SS}\|\le1/4<1.
\]
The determinant identity
\[
 \det(a_\pm I-W-tE)
 =\det(a_\pm I-W)
   \det\bigl(I-tE_S(R_W(a_\pm))_{SS}\bigr)
\]
is nonzero for both signs and every $t\in[0,1]$. If $\lambda_j(t)$ is the $j$th ordered eigenvalue of $W+tE$, then
\[
 |\lambda_j(t)-\lambda_j(t')|\le|t-t'|\|E\|,\qquad
 \lambda_j(0)\in(a_-,a_+),\qquad
 \{t\in[0,1]:\lambda_j(t)=a_\pm\}=\varnothing.
\]
The continuous image $\lambda_j([0,1])$ is connected, so it is contained in $(a_-,a_+)$. In particular $\|W+E\|\le2+d/4$. Here
\[
 \frac{N^{-1+4\phi}}d\le C N^{-1/3+4\phi}(\log N)^{-2}\longrightarrow0,
 \qquad \phi=1/100.
\]
For all sufficiently large $N$, $\|V\|\le d/4$. Adding $V$ in \eqref{W:eq:comparison-decomposition} gives $\|H\|\le2+d/2$ for these $N$. This step is deterministic and remains
valid when $E_S$ depends on $W$.

The estimates for the spectral window, coordinate resolvents, counting function, $V$, the two compression norms, and $E$ give
\[
 \begin{aligned}
 \mathbb P(\|H\|>2+d/2)
 &\le CN^{-c}+C_DN^{-D}+CN^{-151/300}+CN^{-1+8\phi}\\
 &\quad+C\Theta(c_0\sqrt N)+C\Theta_{\mathrm d}(c_0\sqrt N).
 \end{aligned}
\]
Choose $D=1$ and replace $c$ by $\min\{c,151/300,1-8\phi\}>0$. This is \eqref{W:eq:comparison-window-rate}. Each term tends to zero, proving \eqref{W:eq:comparison-window}.

For the Gaussian ensemble put $x=d/(2{\eta_{\mathrm e}})\asymp(\log N)^2$. The upper-tail estimate for the GOE and GUE in \cite{erdosxu2023small} applies for $1\le x\le N^{1/4}$
and bounds the upper edge probability by
\[
 Cx^{-3\beta/4}\exp(-2\beta x^{3/2}/3).
\]
The reflected Gaussian matrix has the same law. A union bound over the two edges proves \eqref{W:eq:comparison-Gaussian-window}. Only the GOE/GUE tail bound is used at
$x\asymp(\log N)^2$; the deviation estimate for bounded entries was used in its stated range.
\end{proof}

Choose $h_d$ and the regularized constructions exactly as in Section~\ref{W:sec:spectral-regularization-and-matrix-derivatives}; in particular $h_d(x)=x$ for $|x|\le2+d$. Use the
same scalar cutoff for the enlarged Gram matrix in the original and modified functions. On its nonzero support, every leading principal Gram matrix required by the functions lies
in the interval on which the regularizations are exact. Both original and modified functions multiplied by their cutoffs are defined to vanish when their respective cutoff for the
enlarged Gram matrix is zero.

\Needspace{6\baselineskip}

\begin{theorem}
\label{W:thm:comparison-agreement}
For every Hermitian $H$ with $\|H\|\le2+d/2$, the following hold simultaneously for all unit starting vectors and all initial blocks up to $\rho_*$. The original and modified
polynomial matrices, their Gram matrices, and the scalar cutoffs applied to the eigenvalues of their Gram matrices agree. On the nonzero common cutoff support, their ordered Jacobi
compressions, integrated potentials, and Krylov projections agree. Consequently the cutoff in the eigenvalues of the Gram matrix, the exponentials of the integrated potential, the
smooth functions of the modulus of continuity, and the bounded localization statistics agree, including the versions multiplied by the same bounded cutoff function for the Schur
complement.

These are local identities in $H$. Thus the derivatives of the smooth functions multiplied by their cutoffs agree to every order for which they are defined.
\end{theorem}
\begin{proof}
Put $\mathcal U_H=\{H+V:V=V^*,\ \|V\|<d/2\}$. For $K\in\mathcal U_H$, Weyl's inequality gives $\sigma(K)\subset(-2-d,2+d)$, and the defining property of $h_d$ gives
\[
 h_d(K)=\int h_d(x)\,\dd E_K(x)=\int x\,\dd E_K(x)=K.
\]
Consequently, for every $r\le\rho_*$ and unit vector $w$,
\[
 \widehat F_r(K,w)=F_r(K,w),\qquad
 \widehat G_r(K,w)=\widehat F_r^*\widehat F_r=F_r^*F_r=G_r(K,w).
\]
The scalar cutoffs applied to the eigenvalues of the Gram matrix therefore agree throughout $\mathcal U_H$. At every point of their common nonzero support, the enlarged Gram matrix
has the positive lower and upper bounds. For each principal Gram matrix, the variational inequalities
\[
 \lambda_{\min}(G_{\rho_*})\|v\|^2
 \le v^*G_rv\le\lambda_{\max}(G_{\rho_*})\|v\|^2
\]
follow by extending $v$ with zero coordinates. Thus all these polynomial matrices have full column rank. Uniqueness of the Cholesky factorization with positive diagonal gives
$\widehat L_r=L_r$. It follows that
\[
 \widehat Q_r=\widehat F_r\widehat L_r^{-*}=F_rL_r^{-*}=Q_r,
 \qquad \widehat J_r=Q_r^*KQ_r=J_r,
 \qquad \widehat P_r=Q_rQ_r^*=P_r.
\]
All regularized matrix functions used in these expressions equal their unregularized versions on this support, by their definitions in
Section~\ref{W:sec:spectral-regularization-and-matrix-derivatives}. The integrated potential is the specified finite sum of entries of $J_r$; hence its values agree, as do every
finite collection of sampled values and their characteristic functions and the smooth functions of their moduli of continuity.

For the localization statistic, both constructions use $\theta_E(m_N^2(2I-K))$. Their traces therefore satisfy
\[
 \operatorname{Tr}\bigl((I-\widehat P_r)
       \theta_E(m_N^2(2I-K))\bigr)
 =\operatorname{Tr}\bigl((I-P_r)
       \theta_E(m_N^2(2I-K))\bigr)
\]
on the common support. Outside it, multiplication by the common cutoff in the eigenvalues of the Gram matrix sets both functions equal to zero. Thus each pair of smooth functions
multiplied by their cutoffs $\widehat\Phi,\Phi$ agrees on the open set $\mathcal U_H$. For every order $q$ at which they are defined and all Hermitian matrix variations
$V_1,\ldots,V_q$,
\[
 D^q(\widehat\Phi-\Phi)(H)[V_1,\ldots,V_q]=0.
\]
Multiplication by the same smooth cutoff function for the Schur complement preserves the equality on $\mathcal U_H$ and all these derivative identities.
\end{proof}

\Needspace{6\baselineskip}

Let $F_N(H,w)$ and $\widehat F_N(H,w)$ be any corresponding functions multiplied by their cutoffs in Theorem~\ref{W:thm:comparison-agreement}, with absolute values at most a fixed
$M$.

\begin{corollary}
\label{W:cor:comparison-transfer}
\begin{equation}
 \sup_{\|w\|=1}
 |\mathbb EF_N(H_N,w)-\mathbb E\widehat F_N(H_N,w)|
 \le2M\mathbb P(\Omega_N^c)=o(1).
 \label{W:eq:comparison-ensemble-transfer}
\end{equation}
For a Gaussian ensemble the same error is smaller than every inverse polynomial. In particular
\begin{align}
 &\sup_{\|w\|=1}|\mathbb EF_N(H_N,w)-\mathbb EF_N({H_N^{\mathrm G}},w)|
 \nonumber\\
 &\quad\le
 \sup_{\|w\|=1}|\mathbb E\widehat F_N(H_N,w)
                      -\mathbb E\widehat F_N({H_N^{\mathrm G}},w)|
 +2M\bigl[\mathbb P(\Omega_N^c)
          +\mathbb P(\|{H_N^{\mathrm G}}\|>2+d/2)\bigr].
 \label{W:eq:comparison-two-ensembles}
\end{align}
If a truncated matrix $\bar H_N$ is used in the comparison, replace the truncated error in these formulas by $2M[\mathbb P(H_N\ne\bar H_N)+\mathbb P(\Omega_N^c)]$.

\end{corollary}
\begin{proof}
On the common event of agreement with the spectral modification $\Omega_N$, the observables satisfy $F_N(H_N,w)=\widehat F_N(H_N,w)$ for every deterministic $w$. Both have modulus
at most $M$, so
\[
 \begin{aligned}
 |\mathbb EF_N(H_N,w)-\mathbb E\widehat F_N(H_N,w)|
 &=\left|\mathbb E[(F_N(H_N,w)-\widehat F_N(H_N,w))
                                      \mathbf1_{\Omega_N^c}]\right|\\
 &\le2M\mathbb P(\Omega_N^c).
 \end{aligned}
\]
The event and its probability are independent of $w$, which permits the supremum. Apply this identity to the Gaussian ensemble with the event on which it agrees with its spectral
modification, whose failure probability has the stated inverse polynomial bounds. For each $w$, compare the original expectations through the two expectations of the modified
functions. The triangle inequality gives
\[
 \begin{aligned}
 |\mathbb EF_N(H_N,w)-\mathbb EF_N(H_N^{\rm G},w)|
 \le{}&2M\mathbb P(\Omega_N^c)
  +|\mathbb E\widehat F_N(H_N,w)-\mathbb E\widehat F_N(H_N^{\rm G},w)|\\
 &+2M\mathbb P(\|H_N^{\rm G}\|>2+d/2).
 \end{aligned}
\]
Taking the supremum proves \eqref{W:eq:comparison-two-ensembles}. For a coupled truncated matrix $\bar H_N$, on $\{H_N=\bar H_N\}\cap\Omega_N$ one also has
$F_N(H_N,w)=\widehat F_N(\bar H_N,w)$. The identical calculation gives
\[
 \sup_w|\mathbb EF_N(H_N,w)-\mathbb E\widehat F_N(\bar H_N,w)|
 \le2M\{\mathbb P(H_N\ne\bar H_N)+\mathbb P(\Omega_N^c)\}.
\]
This estimate transfers the comparison after summation over the entries to the original matrix law.
\end{proof}

The spectral window and agreement estimates hold under Assumption~\ref{W:ass:wigner-lee-yin} in both symmetry classes. Their proof combines the deviation estimate for bounded
standardized entries, a fourth-order trace comparison, and the resolvent compressed to the independent support of the exceptional entries. They are therefore available before the
comparison of growing Gram matrices in Theorem~\ref{W:thm:tightness-main}.

The confinement bound for the original law is $o(1)$, whereas the truncated interpolation has the stronger failure probabilities specified below. Equality of the original and
modified statistics is used in the final expectation comparison. Along the interpolation, derivatives are taken of their globally smooth extensions.

\subsubsection{Expected resolvent estimates}
\label{W:sec:averaged-resolvent-and-conditional-row-estimates}
The comparison requires sums of expectations of resolvent factors, so bounds holding on regularity events must be supplemented by estimates that can be integrated. The first result
below treats a resolvent matrix element for an exactly centered, variance normalized array. Conditional row estimates then supply the moment sums on the events used in
Lemma~\ref{W:lem:comparison-weights}.

The intermediate arrays in this argument have bounded diagonal entries. The passage to truncated and original diagonal distributions is made separately in
Section~\ref{W:sec:conditional-row-bounds-for-diagonal-remainders}. We keep the convention $R_H(z)=(zI-H)^{-1}$, with $\eta=\operatorname{Im}z>0$, $m_H=N^{-1}\operatorname{Tr}R_H$,
and $g$ the semicircle transform in this convention. The scales are $m_N=N^{1/3}$ and $s=N^{-1/3+o(1)}$. A perturbation in the selected entry has the form

\[
 V_\alpha(x)=N^{-1/2}\mathcal T_\alpha x,
 \qquad \mathcal T_\alpha x=\sum_{r=1}^{d_\alpha}x_rT_{\alpha,r},
 \qquad d_\alpha\le2.
\]
Distinct unordered blocks are independent. The two real coordinates of a complex block retain their joint distribution. Constants are uniform in $N$ for a fixed number of starting
vectors, fixed derivative orders, and fixed cutoff functions of the eigenvalues of the Gram matrix. Every factor denoted $N^{o(1)}$ can instead be replaced by an arbitrarily small
fixed power, with the dependence on that power retained in the estimates.

The next theorem bounds the expectation of a resolvent matrix element at a deterministic vector. Its hypotheses are the coordinate estimates and the bound on the normalized trace
stated below.

Let the centered Hermitian array have off-diagonal variance $N^{-1}$, zero pseudovariance in the complex case, and
\[
 |H_{ij}|\le e_N\longrightarrow0,\quad
 \mathbb E|H_{ij}|^3\le C N^{-3/2},\quad
 \mathbb E|H_{ij}|^4\le F_NN^{-2}\quad(i<j),
 \qquad F_N=N^{o(1)}.
\]
On the diagonal assume
\begin{equation}
 \mathbb EH_{ii}^2\le C/N,\qquad |H_{ii}|\le Ce_N,
 \qquad e_N\longrightarrow0.
 \label{W:eq:bound-on-the-intermediate-diagonal-entries}
\end{equation}
This support condition concerns the intermediate comparison arrays. The proof below establishes stability under replacing any of their entries directly from the block resolvent
identity. The error from truncating an original diagonal entry is estimated separately in Section~\ref{W:sec:conditional-row-bounds-for-diagonal-remainders}. Fix a compact real
energy interval and $cN^{-2/3}\le\eta\le M$, with $M\ge1$. Suppose, uniformly over these parameters, with failure at most $N^{-D}$ for every fixed $D$,
\begin{equation}
 \max_{i,j}|R_{H,ij}(z)|\le C_0,
 \qquad |m_H(z)-g(z)|\le \mathfrak c_N/(N\eta),\qquad \mathfrak c_N=N^{o(1)}.
 \label{W:eq:moment-weights-coordinate-trace-input}
\end{equation}
The uniform domain includes the deterministic dyadic vertical parameters between the given height and height one used in the proof of \eqref{W:eq:moment-weights-mean-absolute};
their number is $O(\log N)$. The norm below is $\|\mathbb ER\|=\sup_{\|u\|=\|v\|=1}|\mathbb Eu^*Rv|$, where $u,v$ are deterministic. The centered entry hypotheses apply to this
expectation estimate. The later conditional row argument treats uncentered truncated entries separately using their mean and covariance bounds.

\begin{theorem}
\label{W:thm:moment-weights-mean-resolvent}
Under the preceding moment and resolvent bounds,
\begin{align}
 \|\mathbb ER_H(z)\|&\le N^{o(1)},\label{W:eq:moment-weights-mean-operator}\\
 \mathbb E\,v^*|R_H(z)|v&\le N^{o(1)},\label{W:eq:moment-weights-mean-absolute}\\
 -\operatorname{Im}v^*\mathbb ER_H(z)v
 &\le N^{o(1)}\left(\upsilon(z)+(N\eta)^{-1}\right)
 \label{W:eq:moment-weights-mean-density}
\end{align}
for every deterministic unit vector $v$, uniformly in that vector. The constants are uniform over the specified spectral domain.
\end{theorem}
\begin{proof}
Fix a deterministic unit vector $v$. Set $r=R_Hv$ and $B_v=-\operatorname{Im}v^*\mathbb ER_Hv$. The exact Ward identity is
\begin{equation}
 \mathbb E\|r\|^2=B_v/\eta.
 \label{W:eq:moment-weights-expected-ward}
\end{equation}

Set the unordered block $\alpha$ to zero and denote the resulting matrix by $K_\alpha$. Let $J$ be its set of one or two coordinate indices and write $V=\iota_JT\iota_J^*$ for the
removed matrix, where $\|T\|\le Ce_N$. On the event in \eqref{W:eq:moment-weights-coordinate-trace-input},
\[
 R_{K_\alpha}=R_H-R_H\iota_JT
          (I+R_{H,JJ}T)^{-1}\iota_J^*R_H.
\]
The norm of $R_{H,JJ}$ is at most $2C_0$. For large $N$, $\|R_{H,JJ}T\|\le1/2$ and the inverse in this formula has norm at most two. For arbitrary indices $k,l$, each exterior
vector has norm at most $\sqrt2 C_0$, and therefore
\[
 |(R_{K_\alpha}-R_H)_{kl}|
 \le4C_0^2\|T\|\le C'e_N.
\]
Choose $N_0$ so that the correction in the preceding identity is at most $C_0$ for $N\ge N_0$. Then $\max_{k,l}|(R_{K_\alpha})_{kl}|\le2C_0$. Conversely, at any matrix $K$
satisfying this last bound, every permitted replacement $V(x)=\iota_JT_x\iota_J^*$, $\|T_x\|\le Ce_N$, satisfies
\[
 R_{K+V(x)}=R_K+R_K\iota_JT_x
          (I-R_{K,JJ}T_x)^{-1}\iota_J^*R_K,
 \qquad \max_{k,l}|(R_{K+V(x)})_{kl}|\le4C_0
\]
for all sufficiently large $N$, uniformly in $x$. The event that the entries of $R_K$ obey the $2C_0$ bound throughout the specified spectral domain depends only on $K$. It is
measurable by continuity of the resolvent and a countable dense subset of that domain. The original coordinate event implies this event, so its complementary probability is at most
the prescribed $N^{-D}$. The derivative estimates below are thus available throughout the full support of the conditional replacement distribution. On the exceptional set, the
deterministic bounds are
\[
 \|D_H^qR[T_1,\ldots,T_q]\|
 \le q!\eta^{-q-1}\prod_j\|T_j\|,\qquad 0\le q\le3.
\]
Each $T_j$ is a fixed self-adjoint matrix of bounded norm supported on the selected entry. Since the entry moments of orders at most four are bounded and $\eta\ge cN^{-2/3}$, each
discarded coefficient or integral remainder is bounded by $C(1+\eta^{-4})\le C N^{8/3}$. For an exceptional probability $N^{-D}$, the sum of the absolute expectations of the
$O(N^2)$ discarded terms is at most
\[
 C N^2(1+\eta^{-4})N^{-D}\le C N^{5-D}.
\]
For each $A$, choose $D>A+6$ in the probability estimate for the coordinate resolvent event. The event is $\sigma(K_\alpha)$-measurable. Conditional on $K_\alpha$, its indicator is
fixed while we expand in the one or two real coordinates of the replacement entry. The remainder in the following identity is bounded using the fourth moment of that entry.

Let $X=(X_1,\ldots,X_d)$ be centered, $d\le2$, put $M_j=\mathbb E\|X\|_1^j$, and suppose $f$ is $C^3$ on the convex hull of $0$ and the support of $X$. Write
$A_3=\max_{b,c,e}\sup|\partial_{bce}f|$ on this set. The scalar Taylor formulas, with their integral remainders, are
\begin{align*}
 f(X)&=f(0)+\sum_bX_bf_b(0)
          +\tfrac12\sum_{b,c}X_bX_cf_{bc}(0)+r_0(X),\\
 f_b(X)&=f_b(0)+\sum_cX_cf_{bc}(0)+r_b(X),\\
 f_{bc}(X)&=f_{bc}(0)+r_{bc}(X),
\end{align*}
where
\begin{align*}
 r_0(X)&=\tfrac12\int_0^1(1-t)^2
       \sum_{b,c,e}X_bX_cX_e f_{bce}(tX)\,\mathrm dt,\\
 r_b(X)&=\int_0^1(1-t)
       \sum_{c,e}X_cX_e f_{bce}(tX)\,\mathrm dt,\\
 r_{bc}(X)&=\int_0^1\sum_eX_e f_{bce}(tX)\,\mathrm dt.
\end{align*}
In particular, $|r_0(X)|\le A_3\|X\|_1^3/6$, $|r_b(X)|\le A_3\|X\|_1^2/2$, and $|r_{bc}(X)|\le A_3\|X\|_1$. With $\sigma_{ab}=\mathbb E X_aX_b$ and $\mu_{abc}=\mathbb E X_aX_bX_c$,
substitution gives the exact identity
\begin{align}
 \mathbb E[X_af(X)]
 &=\sum_b\sigma_{ab}\mathbb E f_b(X)
    +\tfrac12\sum_{b,c}\mu_{abc}\mathbb E f_{bc}(X)
    +\mathcal R_a,\notag\\
 \mathcal R_a
 &=\mathbb E[X_ar_0(X)]
   -\sum_b\sigma_{ab}\mathbb E r_b(X)
   -\tfrac12\sum_{b,c}\mu_{abc}\mathbb E r_{bc}(X).
 \label{W:eq:resolvent-cumulant-remainder}
\end{align}
The constant and linear terms in the Taylor expansion agree because $\mathbb EX=0$. Moreover,
\begin{equation}
 |\mathcal R_a|
 \le A_3\bigl(M_4/6+M_2^2/2+M_3M_1/2\bigr)
 \le \tfrac76 A_3M_4.
 \label{W:eq:resolvent-fourth-moment-only}
\end{equation}
The last inequality follows from H\"older: $M_2^2\le M_4$ and $M_3M_1\le M_4$. This real coordinate calculation applies to complex-valued $f$ and to dependent coordinates of $X$.
Centering identifies $\sigma$ and $\mu$ with the second and third cumulants; only the fourth moment is used for the remainder.

We now verify the derivative bounds used for $f=r_j=(Rv)_j$. For coordinate directions $T_1,\ldots,T_q$ supported on
$J=\{i,j\}$,
\begin{equation}
 D^qr_j[T_1,\ldots,T_q]
 =\sum_{\sigma\in\mathfrak S_q}
 e_j^*RT_{\sigma(1)}R\cdots T_{\sigma(q)}r.
 \label{W:eq:resolvent-row-derivatives}
\end{equation}
All resolvents before the final $r$ in this formula are restricted to $J$. On the event fixed by $K_\alpha$, these compressions are bounded uniformly over the permitted entry
values. Consequently, for $q\le3$,
\[
 |D^qr_j[T_1,\ldots,T_q]|
 \le C_q\prod_{\ell=1}^q\|T_\ell\|\,(|r_i|+|r_j|).
\]
Let $\mathcal A$ denote the event measurable with respect to the matrix $K$ with this block removed, and let $\mu$ be the joint law of its one or two real coordinates. Write $V(x)$
for the corresponding matrix perturbation and $r_x=(R_{K+V(x)}v)_J$. On $\mathcal A$,
\[
 r_x=(I-R_{K,JJ}V(x)_{JJ})^{-1}r_0,\qquad
 \|(I-R_{K,JJ}V(x)_{JJ})^{\pm1}\|\le C.
\]
The same bounds hold at $tx$ for $0\le t\le1$. Equivalence of the Euclidean and $\ell^1$ norms in dimension at most two gives
\[
 \sup_{x\in\operatorname{supp}\mu,\,0\le t\le1}\|r_{tx}\|_1
 \le C\|r_0\|_1
 \le C\int\|r_x\|_1\,\mu(\mathrm dx).
\]
Thus the conditional derivative constant in \eqref{W:eq:resolvent-fourth-moment-only} can be taken as $A_3(K)=C\|r_0\|_1$. Applying that identity conditionally first, and then
integrating $K$, yields
\[
 \mathbb E\bigl[\mathbf1_{\mathcal A}|\mathcal R_a(K)|\bigr] \le C M_4\mathbb E\bigl[\mathbf1_{\mathcal A}\|r_0\|_1\bigr] \le C M_4\mathbb E\left[\mathbf1_{\mathcal A} \int\|r_x\|_1\,\mu(\mathrm dx)\right] \le C M_4\mathbb E(|r_i|+|r_j|).
\]
Here $M_4=\int\|x\|_1^4\mu(\mathrm dx)$ is the deterministic block moment defined above. For the third coordinate moments, the term in \eqref{W:eq:resolvent-cumulant-remainder} is
\[
 T_a(K)=\tfrac12\sum_{b,c=1}^d\mu_{abc}
                  \int \partial_{bc}f(x)\,\mu(\dd x).
\]
The formula \eqref{W:eq:resolvent-row-derivatives} at $q=2$ has two products. Boundedness of their resolvent compressions, followed by the preceding comparison of $r_x$ with $r_0$,
gives
\[
 \max_{b,c}\sup_{x\in\operatorname{supp}\mu}
                 |\partial_{bc}f(x)|\le C\|r_0\|_1.
\]
Also $\sum_{b,c}|\mu_{abc}|\le\mathbb E(|X_a|\|X\|_1^2)\le M_3$. Consequently
\[
 \mathbb E[\mathbf1_{\mathcal A}|T_a(K)|]
 \le C M_3\mathbb E[\mathbf1_{\mathcal A}\|r_0\|_1]
 \le C M_3\mathbb E(|r_i|+|r_j|).
\]
Both estimates keep the deterministic moments of the replacement distribution outside the integral over $K$. The integral defining $T_a$ uses the joint distribution of all its real
coordinates. For a complex entry, multiplication of the two resulting identities by $1$ and $\mathrm i$ gives the identity for $H_{ij}$, with at most the sum of the two absolute
bounds. On $\mathcal A^c$, the deterministic derivative estimates and the choice $D>A+6$ above give the summed $O(N^{-A})$ error.

We write the covariance contraction in the two matrix entry orientations. In the real case, with $T_1=E_{ij}+E_{ji}$,
\[
 (RT_1r)_j=R_{ji}r_j+R_{jj}r_i.
\]
In the complex case use the matrix entry coordinates $H_{ij}=x+\mathrm i y$ and $T_2=\mathrm i(E_{ij}-E_{ji})$. The covariance assumptions give $\mathbb Ex^2=\mathbb Ey^2=1/(2N)$
and $\mathbb Exy=0$. Directly,
\[
 \partial_xr_j=R_{ji}r_j+R_{jj}r_i, \qquad \partial_yr_j=\mathrm i R_{ji}r_j-\mathrm i R_{jj}r_i, \qquad (\partial_x+\mathrm i\partial_y)r_j=2R_{jj}r_i.
\]
Thus the covariance term in $\mathbb E[H_{ij}r_j]$ is $N^{-1}\mathbb E(R_{ji}r_j+R_{jj}r_i)$ in the real case and $N^{-1}\mathbb E R_{jj}r_i$ in the complex case. The third moments
of the real coordinates contribute to \eqref{W:eq:resolvent-cumulant-remainder}. The indicator is $\sigma(K_\alpha)$-measurable and hence constant during integration over $X$.

The third and remainder terms consequently have absolute bounds
\[
 C N^{-3/2}\mathbb E(|r_i|+|r_j|),\qquad
 C F_NN^{-2}\mathbb E(|r_i|+|r_j|),
\]
respectively. For a centered diagonal entry, fix the matrix $K$ with that entry set to zero. On its regularity event the rank-one identity gives
\[
 r_i(K+xE_{ii})=\frac{r_i(K)}{1-xR_{K,ii}},\qquad
 |1-xR_{K,ii}|\ge\tfrac12.
\]
The conditional law $\mu_i$ of $x$ has mean zero. Subtracting the constant $r_i(K)$ before integration gives
\[
 \left|\int x r_i(K+xE_{ii})\,\mu_i(\dd x)\right|
 =\left|r_i(K)R_{K,ii}\int\frac{x^2}{1-xR_{K,ii}}\,\mu_i(\dd x)\right|
 \le C N^{-1}|r_i(K)|.
\]
The reverse rank-one formula gives $|r_i(K)|\le C\int|r_i(K+xE_{ii})|\,\mu_i(\dd x)$. Integrating over $K$, and adding the previously bounded exceptional contribution, proves
$|\mathbb EH_{ii}r_i|\le C N^{-1}\mathbb E|r_i|+O(N^{-A})$. Only the second moment of the replacement diagonal variable occurs.

Summing the resolvent equation $zr=v+Hr$ gives
\begin{equation}
 (z-g)\mathbb Er
 =v+\mathbb E[(m_H-g)r]+\frac{\sigma}{N}\mathbb ER_H^2v+\mathcal Z_v,
 \quad \sigma=\mathbf1_{\{\beta=1\}},
 \label{W:eq:moment-weights-schwinger}
\end{equation}
where
\begin{equation}
 \|\mathcal Z_v\|
 \le C(N^{-1/2}+F_N/N)\sqrt{B_v/\eta}+O(N^{-A}).
 \label{W:eq:moment-weights-schwinger-error}
\end{equation}
To check the vector norm, put $a_i=\mathbb E|r_i|$ and $A_1=\sum_i a_i$. Jensen's inequality and the Ward identity give
\[
 \|a\|_2^2\le\sum_i\mathbb E|r_i|^2=B_v/\eta,
 \qquad A_1\le\sqrt N\|a\|_2.
\]
The vector $z^{(3)}$ of summed third moment errors and the vector $z^{(4)}$ of fourth moment remainders satisfy componentwise
\[
 |z_i^{(3)}|\le CN^{-3/2}(Na_i+A_1),\qquad
 |z_i^{(4)}|\le CF_NN^{-2}(Na_i+A_1).
\]
Therefore
\begin{align*}
 \|z^{(3)}\|_2&\le CN^{-3/2}
       (N\|a\|_2+\sqrt N A_1)\le2CN^{-1/2}\sqrt{B_v/\eta},\\
 \|z^{(4)}\|_2&\le CF_NN^{-2}
       (N\|a\|_2+\sqrt N A_1)\le2CF_NN^{-1}\sqrt{B_v/\eta}.
\end{align*}
The covariance terms are
\[
 N^{-1}\sum_{j\ne i}R_{jj}r_i=m_Hr_i-N^{-1}R_{ii}r_i,
 \quad
 N^{-1}\sum_{j\ne i}R_{ji}r_j
   =N^{-1}(R_H^2v)_i-N^{-1}R_{ii}r_i
\]
in the real case; the second contraction is absent in the complex case by the displayed real coordinate calculation. The real case therefore contributes the full term
$N^{-1}\mathbb ER_H^2v$ in \eqref{W:eq:moment-weights-schwinger}. Their omitted diagonal part has $i$th component
\[
 d_i=\mathbb E[H_{ii}r_i]
             -\frac{1+\sigma}{N}\mathbb E[R_{ii}r_i].
\]
The centered diagonal entry bound above and $|R_{ii}|\le4C_0$ on the chosen event give $|d_i|\le CN^{-1}a_i+O(N^{-A-1})$; on its complement use the deterministic resolvent bound
and choose the failure exponent as above. Therefore
\[
 \|d\|_2\le CN^{-1}\|a\|_2+O(N^{-A})
          \le CN^{-1}\sqrt{B_v/\eta}+O(N^{-A}).
\]
Adding this inequality to the bounds for $z^{(3)}$ and $z^{(4)}$ proves \eqref{W:eq:moment-weights-schwinger-error}.

Let $M_R=\|\mathbb ER_H\|$ . On the trace event, Cauchy--Schwarz and the expected Ward identity give
\[
 \|\mathbb E[(m_H-g)r]\|
 \le\frac{\mathfrak c_N}{N\eta}\sqrt{B_v/\eta}+O(N^{-A}),\qquad
 N^{-1}\|\mathbb ER_H^2v\|
 \le\frac1{N\eta}\sqrt{B_v/\eta}.
\]
In the second inequality use $\|R_H^2v\|\le\eta^{-1}\|R_Hv\|$; in the first, the exceptional contribution is bounded by the deterministic resolvent powers and the chosen failure
probability. Since $B_v\le M_R$ and $(z-g)^{-1}=g$, with $|g|\le1$, taking norms in \eqref{W:eq:moment-weights-schwinger} and then the supremum over unit $v$ gives
\[
 M_R\le C+C A_N\sqrt{M_R},\qquad
 A_N=\frac{\mathfrak c_N+1}{N\eta^{3/2}}
       +\frac{N^{-1/2}+F_N/N}{\sqrt\eta}=N^{o(1)}.
\]

Solving this quadratic inequality in $\sqrt{M_R}$ gives $M_R\le C'(1+A_N)^2=N^{o(1)}$. For \eqref{W:eq:moment-weights-mean-absolute}, a spectral measure $\mu_v$ satisfies
\[
 \mathbb E\mu_v([u-t,u+t])
 \le2t\bigl[-\operatorname{Im}v^*\mathbb ER_H(u+\mathrm it)v\bigr]
 \le2tN^{o(1)}.
\]
For $0<\eta\le1/2$, choose $k_*\ge0$ so that $1/2<2^{k_*}\eta\le1$ and put $t_k=2^k\eta$. Let $C_N=N^{o(1)}$ be a common bound at these $O(\log N)$ heights. The central interval
and the shells satisfy
\begin{align*}
 \mathbb E\int_{|x-u|\le\eta}|z-x|^{-1}\,\mu_v(\dd x)
 &\le\eta^{-1}\mathbb E\mu_v([u-\eta,u+\eta])\le2C_N,\\
 \mathbb E\int_{t_k<|x-u|\le t_{k+1}}|z-x|^{-1}\,\mu_v(\dd x)
 &\le t_k^{-1}\mathbb E\mu_v([u-t_{k+1},u+t_{k+1}])
 \le4C_N \quad(0\le k<k_*).
\end{align*}
On $|x-u|>t_{k_*}$ the integral is at most $2\mu_v(\mathbb R)=2$. Thus
\[
 \mathbb E\,v^*|R_H(z)|v\le2C_N+4k_*C_N+2=N^{o(1)}.
\]
For $\eta>1/2$ the deterministic bound $\|R_H\|\le2$ gives the same conclusion.

Finally multiply \eqref{W:eq:moment-weights-schwinger} by $g$, take its $v$-matrix element and its negative imaginary part. On the trace event put $a_0=\mathfrak c_N/(N\eta)$.
Since $-\operatorname{Im}v^*R_Hv\ge0$,
\[
 |\operatorname{Im}\{g(m_H-g)v^*R_Hv\}|
 \le a_0\bigl[-\operatorname{Im}v^*R_Hv+v^*|R_H|v\bigr].
\]
Taking expectation in the preceding inequality and using \eqref{W:eq:moment-weights-mean-absolute} gives
\[
 \bigl|\operatorname{Im}\{g\,\mathbb E[(m_H-g)v^*R_Hv]\}\bigr|
 \le a_0B_v+a_0N^{o(1)}+O(N^{-A}).
\]
The leading term is $-\operatorname{Im}g=\upsilon$. The real covariance correction satisfies
\[
 \frac1N|\operatorname{Im}\{g\,\mathbb E(v^*R_H^2v)\}|
 \le\frac1N\mathbb E\|R_Hv\|^2=\frac{B_v}{N\eta},
\]
since the spectral theorem gives $|v^*R_H^2v|\le\int|z-x|^{-2}\mu_v(\mathrm dx)$. The remaining vector error has bound
\[
 |\operatorname{Im}\{g\,v^*\mathcal Z_v\}|
 \le C(N^{-1/2}+F_N/N)\sqrt{B_v/\eta}+O(N^{-A}).
\]
All exceptional event errors here use the deterministic bounds $\|R_H\|\le\eta^{-1}$ and $|m_H|\le\eta^{-1}$, with the failure exponent chosen after these fixed powers. Thus
\[
 B_v\le\upsilon+a_0B_v+a_0N^{o(1)}
       +\frac{B_v}{N\eta}
       +C(N^{-1/2}+F_N/N)\sqrt{B_v/\eta}+O(N^{-A}).
\]
As $a_0+(N\eta)^{-1}=o(1)$, it is at most $1/4$ for large $N$. Writing $d_N=N^{-1/2}+F_N/N$, Young's inequality gives
\[
 Cd_N\sqrt{B_v/\eta}\le\tfrac14B_v+C'\frac{d_N^2}{\eta},
 \qquad \frac{d_N^2}{\eta}
 \le\frac{2+2F_N^2/N}{N\eta}=\frac{N^{o(1)}}{N\eta}.
\]
Absorbing the two contributions of $B_v/4$ yields $B_v\le N^{o(1)}[\upsilon+(N\eta)^{-1}]$.
\end{proof}

\subsubsection{Conditional row estimates under a bound on the Gram matrix}
\label{W:sec:row-estimates-under-gram-constraint}

The empirical counting bound gives information averaged over all vectors. The comparison also needs the spectral measure at the particular starting vector $w$, and this is where
the bound on the Gram matrix enters. An upper bound on $G_\rho$ controls every polynomial of degree below $\rho$ at $w$. Approximating the resolvent by such a polynomial converts
that control into a vector estimate, with the approximation error measured relative to the semicircle resolvent. The construction uses the Chebyshev polynomial $\mathrm T_\rho$ of
the first kind, normalized by
$\mathrm T_\rho(\cosh z)=\cosh(\rho z)$.

\begin{lemma}
\label{W:lem:moment-weights-gram-resolvent}
Suppose $\sigma(H)\subset[-2\cosh t,2\cosh t]$, $0\le t<s$, and $z=2\cosh(s+\mathrm i\theta)$. Assume
\[
 G_{\rho}(H,w)=\bigl(w^*p_a(H)p_b(H)w\bigr)_{a,b<{\rho}}
 \le\kappa^{-1}I,
 \qquad \varepsilon_{\rho}:=\frac{\cosh({\rho}t)}{\sinh({\rho}s)}<1.
\]
Then
\begin{equation}
 \|R_H(z)w\|^2
 \le\kappa^{-1}
       \left(\frac{1+\sinh({\rho}s)^{-1}}{1-\varepsilon_{\rho}}\right)^2
       \int\frac{\mu_{\rm sc}(\mathrm dx)}{|z-x|^2}.
 \label{W:eq:moment-weights-gram-relative}
\end{equation}
In particular, for $t\le s/2$ and ${\rho}s\to\infty$, the bound is $C\kappa^{-1}/s$ on the ellipse. On its vertical extensions, integration of this estimate gives
$w^*|R_H(z)|w\le C\kappa^{-1}(1+|\log\eta|)$; in a fixed edge cone this last bound is $C\kappa^{-1}$.
\end{lemma}
\begin{proof}
The function
\[
 q_z(x)=\frac{1-\mathrm T_{\rho}(x/2)/\mathrm T_{\rho}(z/2)}{z-x}
\]
is a polynomial of degree ${\rho}-1$: its numerator vanishes at $x=z$. For real $x$ in the spectral interval, $|\mathrm T_{\rho}(x/2)|\le\cosh({\rho}t)$, while
$|\mathrm T_{\rho}(z/2)|\ge\sinh({\rho}s)$. Thus
\[
 |z-x|^{-1}\le(1-\varepsilon_{\rho})^{-1}|q_z(x)|.
\]
Expand $q_z=\sum_{j=0}^{\rho-1}c_jp_j$, with $c_j\in\mathbb C$. Orthogonality of the real polynomials $p_j$ for $\mu_{\rm sc}$ gives
\[
 \|q_z(H)w\|^2=c^*G_\rho(H,w)c
 \le\kappa^{-1}\|c\|_2^2
 =\kappa^{-1}\int|q_z(x)|^2\,\mu_{\rm sc}(\dd x).
\]
Integrating the squared pointwise relative error bound against the positive spectral measure at $w$ gives $\|R_H(z)w\|^2\le(1-\varepsilon_\rho)^{-2}\|q_z(H)w\|^2$. On $[-2,2]$ the
numerator of $q_z$ has modulus at most $1+\sinh({\rho}s)^{-1}$. This proves \eqref{W:eq:moment-weights-gram-relative}. For nonreal $z$ on the ellipse, $g(z)=e^{-s-\mathrm i\theta}$
on the upper half and $\eta=2\sinh s\sin\theta$, so that $\int|z-x|^{-2}\,\mathrm d\mu_{\rm sc}=e^{-s}/(2\sinh s)$. At the positive real tip $x_s=2\cosh s>2$, the distance from
$[-2,2]$ is $2\cosh s-2>0$. Differentiation of the Stieltjes integral is therefore justified by the constant bound $(2\cosh s-2)^{-2}$. Since $g(x_s)=\me^{-s}$ and
$\dd x_s/\dd s=2\sinh s$,
\[
 \int\frac{\mu_{\rm sc}(\dd x)}{(x_s-x)^2}
   =-g'(x_s)=\frac{\me^{-s}}{2\sinh s}.
\]
At the negative tip, substituting $x\mapsto-x$ in the integral gives the same value because $\mu_{\rm sc}$ is symmetric.

Increasing $|\operatorname{Im}z|$ at fixed real part increases the confocal ellipse parameter. This follows directly from the equation $u^2/(4\cosh^2 s')+y^2/(4\sinh^2s')=1$, whose
left side is strictly decreasing in $s'>0$. If $s'\ge s$, then $\cosh(\rho t)/\sinh(\rho s')\le\varepsilon_\rho<1$ and $1+\sinh(\rho s')^{-1}\le1+\sinh(\rho s)^{-1}$. Consequently
the relative approximation on each vertical extension has a constant no larger than the original one. When $t\le s/2$ and $\rho s\to\infty$, the ratio is bounded by
$2(e^{-\rho s/2}+e^{-3\rho s/2})/(1-e^{-2\rho s})\to0$, so this constant is at most a fixed multiple of $\kappa^{-1}$. For a real number $r$,
\[
 (r^2+\eta^2)^{-1/2}
 =\frac2\pi\int_0^\infty(r^2+\eta^2+t^2)^{-1}\,\mathrm dt.
\]
Let $\mu_w$ be the spectral measure of $H$ at $w$, and put $h(t)=\sqrt{\eta^2+t^2}$. The relative bound at each vertical extension has a common constant $C\kappa^{-1}$, since its
confocal parameter is at least $s$. Tonelli's theorem therefore gives
\begin{align*}
 w^*|R_H(E+\mathrm i\eta)|w
 &=\frac2\pi\int_0^\infty\int
       \frac{\mu_w(\mathrm dx)}{(E-x)^2+h(t)^2}\,\mathrm dt\\
 &\le\frac{2C}{\pi\kappa}\int_0^\infty\int_{-2}^2
       \frac{\mu_{\rm sc}(\mathrm dx)}{(E-x)^2+h(t)^2}\,\mathrm dt\\
 &=\frac C\kappa\int_{-2}^2
       \frac{\sqrt{4-x^2}}{2\pi\sqrt{(E-x)^2+\eta^2}}\,\mathrm dx.
\end{align*}
For $|E|\le3$ and $0<\eta\le1$, the last integral is at most $C\int_{-5}^{5}(u^2+\eta^2)^{-1/2}\,\mathrm du =2C\operatorname{arsinh}(5/\eta)\le C(1+|\log\eta|)$. For $|E|>3$ or
$\eta>1$ the bound by total mass is constant. In a fixed right edge cone $2-E\le C_0\eta$, write $u=2-x$. The density is at most $C\sqrt u$. If $E\ge2$, the denominator is at least
$u$ and $\int_0^4u^{-1/2}\,\mathrm du<\infty$. If $0<2-E\le C_0\eta$, split at $u_0=2(C_0+1)\eta$. On $u\le\min(u_0,4)$ the integral is at most
$C\eta^{-1}u_0^{3/2}=C_{C_0}\sqrt\eta$; on $u\ge u_0$, $|u+E-2|\ge u/2$, so its contribution is at most $2C\int_0^4u^{-1/2}\,\mathrm du$. The left edge is obtained by $x\mapsto-x$.
This proves the two bounds for the quadratic forms of $|R_H|$.
\end{proof}

Let $s=N^{-1/3+o(1)}$, $B_N=N^{o(1)}$, and let $\mathcal G$ be any event on which, at a fixed nonreal $z$,
\begin{equation}
 \begin{gathered}
 \max_{j,k}|R_{jk}|\le C,\quad \min_i|R_{ii}|\ge c>0,\quad
 \max_i\|Re_i\|^2\le B_N/s,\\
 \|Rw\|^2\le B_N/s,\quad \max_i|(Rw)_i|\le B_N.
 \end{gathered}
 \label{W:eq:moment-weights-stopping}
\end{equation}
The event $\mathcal G$ may depend on all matrix entries; below it is obtained by restricting the Gram matrix.

\Needspace{6\baselineskip} Because $\mathcal G$ may depend on the row being integrated, conditional integration uses a larger event determined by the complementary minor. The Schur
identities below construct that event after bounding $|R_{ii}|$ on $\mathcal G$.

Suppose the off-diagonal entries are independent and centered, $\mathbb E|H_{ij}|^2\le C/N$, and put
$\mu_p=\sup_{i<j}\mathbb E|H_{ij}|^p$.

\begin{proposition}
\label{W:prop:moment-weights-row-moments}
For every even $p\ge4$,
\begin{equation}
 \mathbb E\mathbf1_{\mathcal G}|(Rw)_i|^p
 \le C_pB_N^{C_p}
 \left[|w_i|^p+(Ns)^{-p/2}+\mu_p/s\right].
 \label{W:eq:moment-weights-row-p}
\end{equation}
For $p=2$ the last term is unnecessary. The bound also holds after a fixed finite set of entries has been set to zero; the proof integrates only the off-diagonal row. Under
truncation and the original tail condition,
\begin{equation}
 \mu_p\le N^{-2}F_N e_N^{p-4},\qquad p\ge4,\qquad
 F_N=N^{o(1)},\quad e_N\le1,
 \label{W:eq:moment-weights-truncated-moments}
\end{equation}
so, uniformly in the deterministic unit starting vector,
\begin{align}
 \mathbb E\mathbf1_{\mathcal G}|(Rw)_i|^4
 &\le N^{o(1)}(|w_i|+N^{-1/3})^4,
 \label{W:eq:moment-weights-row-four}\\
 \mathbb E\mathbf1_{\mathcal G}|(Rw)_i|^6
 &\le N^{o(1)}
       (|w_i|^6+N^{-2}+N^{-5/3}e_N^2),
 \label{W:eq:moment-weights-row-six}\\
 \sum_i\mathbb E\mathbf1_{\mathcal G}|(Rw)_i|^p
 &\le N^{o(1)},\qquad 3\le p\le 12.
 \label{W:eq:moment-weights-row-sums}
\end{align}
\end{proposition}
\begin{proof}
Fix the spectral parameter and the deterministic vector $w$. Write $r=Rw$, $c^{(i)}=R^{(i)}w^{(i)}$, and let $h_i$ be the off-diagonal column joining coordinate $i$ to its
complementary minor. Schur complementation gives the exact identities
\begin{equation}
 r_i=R_{ii}(w_i+h_i^*c^{(i)}),\qquad
 c_j^{(i)}=r_j-\frac{R_{ji}}{R_{ii}}r_i\quad(j\ne i).
 \label{W:eq:conditional-row-schur}
\end{equation}
We may increase $B_N$ so that $B_N\ge1$. On $\mathcal G$, \eqref{W:eq:moment-weights-stopping} implies
\[
 \|c^{(i)}\|_\infty\le CB_N,\qquad
 \|c^{(i)}\|_2^2
 \le 2\|r\|_2^2+2c^{-2}|r_i|^2\|Re_i\|_2^2
 \le CB_N^3/s.
\]
Let $\mathcal A_i$ be the event that these two bounds hold for the minor vector $R^{(i)}w^{(i)}$. The minor resolvent exists for nonreal $z$, so this defines an event in
$\sigma(H^{(i)})$, and $\mathcal G\subseteq\mathcal A_i$. For a real point of an exterior contour the same definition includes existence of the minor resolvent; spectral enclosure
and interlacing give the inclusion. We only use this extension where that enclosure has already been proved.

The multiplier $R_{ii}$ must be bounded before changing the event. Consequently, for $p\ge2$,
\begin{equation}
 \begin{split}
 \mathbb E\mathbf1_{\mathcal G}|r_i|^p
 &\le C_p|w_i|^p+
 C_p\mathbb E_{H^{(i)}}\!\left[
 \mathbf1_{\mathcal A_i}\,
       \mathbb E_i|h_i^*c^{(i)}|^p\right].
 \end{split}
 \label{W:eq:conditional-minor-domination}
\end{equation}
Here $\mathbb E_i$ integrates the off-diagonal column $h_i$, with the minor vector and $\mathbf1_{\mathcal A_i}$ fixed.
This is an enlargement of an event in the integral, not conditioning the row on $\mathcal G$. In particular, for $j\ne
k$,
\[
 \mathbb E_i[\overline{h_{i,j}}h_{i,k}]=0
\]
with $\mathbf1_{\mathcal A_i}$ and $c^{(i)}$ fixed. The inclusion $\mathcal G\subseteq\mathcal A_i$ is used before this conditional calculation.

For the fixed minor put $Y_j=\overline{h_{i,j}}c_j^{(i)}$, viewing a complex number as a vector in $\mathbb R^2$. These vectors are independent and centered under $\mathbb E_i$.
For even $p\ge4$, Lemma~\ref{W:lem:schur-rosenthal} gives
\[
 \|\textstyle\sum_jY_j\|_{L^p(\mathbb E_i)}
 \le C_p\left[(\textstyle\sum_j\mathbb E_i|Y_j|^2)^{1/2}
                      +(\textstyle\sum_j\mathbb E_i|Y_j|^p)^{1/p}\right].
\]
Use $(a+b)^p\le2^{p-1}(a^p+b^p)$ and the two estimates
\[
 \sum_j\mathbb E_i|Y_j|^2\le C N^{-1}\|c^{(i)}\|_2^2,
 \qquad
 \sum_j\mathbb E_i|Y_j|^p\le\mu_p
                   \|c^{(i)}\|_\infty^{p-2}\|c^{(i)}\|_2^2.
\]
Substitution gives
\begin{equation}
 \mathbb E_i|h_i^*c^{(i)}|^p
 \le C_p\left[
 \left(\frac{C\|c^{(i)}\|_2^2}{N}\right)^{p/2}
 +\mu_p\|c^{(i)}\|_\infty^{p-2}\|c^{(i)}\|_2^2\right].
 \label{W:eq:conditional-row-moment}
\end{equation}
For $p=2$, independence and centering give $\mathbb E_i|\sum_jY_j|^2=\sum_j\mathbb E_i|Y_j|^2$; all terms $\mathbb E_i(Y_j\overline{Y_k})$, $j\ne k$, vanish. On $\mathcal A_i$ the
preceding two bounds are at most
\[
 C_p B_N^{3p/2}(Ns)^{-p/2},\qquad
 C_p\mu_p B_N^{p+1}/s,
\]
respectively. Substituting them into \eqref{W:eq:conditional-minor-domination} proves \eqref{W:eq:moment-weights-row-p}, with a constant exponent depending only on the fixed $p$.
Here the independent summands are the complex entries of the row, each integrated with its joint real coordinate law. Restricting the sum to a deterministic set $D$ preserves
independence and decreases both nonnegative moment sums above. The constants are therefore uniform over such deletions.

For \eqref{W:eq:moment-weights-truncated-moments}, the support $|\xi_{ij}|\le e_N\sqrt N$ gives
\[
 N^{-p/2}\mathbb E|\xi_{ij}|^p
 \le N^{-2}e_N^{p-4}\mathbb E|\xi_{ij}|^4\qquad(p\ge4).
\]
The truncated fourth moment is bounded by $C\log(eN)$ under the tail assumption; a fixed centering changes this bound by at most a constant factor. With $s=N^{-1/3+o(1)}$, the
variance term for $p=4$ is $N^{-4/3+o(1)}$, and the tail term is $N^{-5/3+o(1)}$. For $p=6$ these terms are $N^{-2+o(1)}$ and $N^{-5/3+o(1)}e_N^2$. Their ratio is
$N^{1/3+o(1)}e_N^2$, which need not tend to zero. In particular $e_N^{-1}=N^{o(1)}$ does not permit the sixth moment contribution to be replaced by the Gaussian variance term. The
term involving higher moments of individual coordinates must remain in the subsequent pair sums. This proves
\eqref{W:eq:moment-weights-row-four}--\eqref{W:eq:moment-weights-row-six}. For even $p\ge4$, summing \eqref{W:eq:moment-weights-row-p} uses $\sum_i|w_i|^p\le1$ and yields
\[
 N^{o(1)}\bigl(1+N^{1-p/3}+N^{-2/3}e_N^{p-4}\bigr)=N^{o(1)}.
\]

Let $I_i$ denote the indicator in the moment bound. For $p=3$,
\[
 \mathbb E[I_i|r_i|^3]
 \le\bigl(\mathbb E[I_i|r_i|^2]\bigr)^{1/2}
      \bigl(\mathbb E[I_i|r_i|^4]\bigr)^{1/2}
 \le N^{o(1)}(|w_i|+N^{-1/3})^3.
\]
The sum is bounded by $4N^{o(1)}(\sum_i|w_i|^3+1)\le8N^{o(1)}$. For an odd $p=2k+1\ge5$, apply Cauchy--Schwarz twice:
\[
 \begin{aligned}
 \sum_i\mathbb E[I_i|r_i|^{2k+1}]
 &\le\sum_i\bigl(\mathbb E[I_i|r_i|^{2k}]\bigr)^{1/2}
                  \bigl(\mathbb E[I_i|r_i|^{2k+2}]\bigr)^{1/2}\\
 &\le\left(\sum_i\mathbb E[I_i|r_i|^{2k}]\right)^{1/2}
       \left(\sum_i\mathbb E[I_i|r_i|^{2k+2}]\right)^{1/2}
 \le N^{o(1)}.
 \end{aligned}
\]
This proves all of \eqref{W:eq:moment-weights-row-sums}.
\end{proof}

For the uncentered truncated arrays, the conditional calculation has an additional mean term. Put $\bar h_i=\mathbb Eh_i$ and $\widetilde h_i=h_i-\bar h_i$. In
\eqref{W:eq:conditional-minor-domination}, replace the inner variable by $\widetilde h_i^*c^{(i)}+\bar h_i^*c^{(i)}$ and apply \eqref{W:eq:conditional-row-moment} only to
$\widetilde h_i$. For the truncated laws, $|\mathbb EH_{ij}|\le N^{-2+o(1)}$, whence $\|\bar h_i\|_2\le N^{-3/2+o(1)}$. Thus on the minor event used at Gram matrix dimension
$\rho$,
\begin{equation}
 |\bar h_i^*c^{(i)}|
 \le \|\bar h_i\|_2\|c^{(i)}\|_2
 \le \sqrt\rho\,N^{-3/2+o(1)}
 =N^{-4/3+o(1)}\quad(\rho\asymp N^{1/3}).
 \label{W:eq:conditional-row-with-mean}
\end{equation}
This term is smaller than the $\sqrt{\rho/N}$ coordinate scale. It is included before taking any row moment. In particular the preceding estimates apply to the uncentered laws used
in the comparison, with the indicated deterministic mean term. For fixed $p$, the inequality $|x+y|^p\le2^{p-1}(|x|^p+|y|^p)$ accounts for this deterministic mean contribution in
the same row estimate.

For each pair $\alpha=(a,b)$, apply the preceding moment estimate to its matrix and regularity event. The resulting deterministic bounds involve the same coordinates of $w$, so
they can be summed over $\alpha$. Include both orientations $z,\bar z$ by setting $a_i=|(R(z)w)_i|+|(R(\bar z)w)_i|+s$ at that block. Let $U_\alpha=a_aa_b$, $V_\alpha=a_a^2+a_b^2$
on its regularity event and zero outside. Then
\begin{equation}
\begin{array}{c|c}
 \text{quantity}&\text{upper bound}\\ \hline
 \sum_{a<b}\mathbb E U_{ab}&N^{4/3+o(1)}\\
 \sum_{a<b}\mathbb E U_{ab}^2&N^{2/3+o(1)}\\
 \sum_{a<b}\mathbb E V_{ab}&N^{4/3+o(1)}\\
 \sum_{a<b}\mathbb E (U_{ab}V_{ab})&N^{2/3+o(1)}\\
 \sum_{a<b}\mathbb E V_{ab}^2&N^{1+o(1)}\\
 \sum_{a<b}\mathbb E U_{ab}^3&N^{1/3+o(1)}.
\end{array}
\label{W:eq:moment-weights-weight-table}
\end{equation}
For diagonal resolvent weights $V_a=a_a^2$ the corresponding sums are $\sum_a\mathbb EV_a\le N^{1/3+o(1)}$ and
$\sum_a\mathbb EV_a^2\le N^{o(1)}$.

\begin{proof}
Put $b_i=|w_i|+N^{-1/3}$. The second and fourth row moments imply $\mathbb Ea_i^2\le N^{o(1)}b_i^2$ and $\mathbb Ea_i^4\le N^{o(1)}b_i^4$ on each of the corresponding regularity
events. Write $S_k=\sum_i b_i^k$. Since $\|w\|_2=1$,
\[
 S_1\le\sqrt N+N^{2/3},\quad S_2\le2+2N^{1/3},\quad
 S_3\le4\sum_i|w_i|^3+4\le8,\quad
 S_4\le8\sum_i|w_i|^4+8N^{-1/3}\le16.
\]
More generally, for every fixed $p\ge3$,
$\sum_i b_i^p\le2^{p-1}(\sum_i|w_i|^p+N^{1-p/3})\le2^p$.
For each pair, H\"older with exponents $(2,2)$ for the first and second powers, and $(4/3,4)$ for the products $a_a^3a_b$ and $a_aa_b^3$, gives
\begin{align*}
 \mathbb E U_{ab}&\le N^{o(1)}b_ab_b,&
 \mathbb E U_{ab}^2&\le N^{o(1)}b_a^2b_b^2,\\
 \mathbb E V_{ab}&\le N^{o(1)}(b_a^2+b_b^2),&
 \mathbb E(U_{ab}V_{ab})&\le N^{o(1)}(b_a^3b_b+b_ab_b^3),\\
 \mathbb E V_{ab}^2&\le N^{o(1)}(b_a^4+2b_a^2b_b^2+b_b^4).
\end{align*}
In particular, summing these five inequalities gives respectively
\[
 N^{o(1)}S_1^2,\quad N^{o(1)}S_2^2,\quad
 N^{o(1)}NS_2,\quad 2N^{o(1)}S_3S_1,\quad
 2N^{o(1)}(NS_4+S_2^2).
\]
Their powers are $N^{4/3+o(1)},N^{2/3+o(1)},N^{4/3+o(1)}, N^{2/3+o(1)},N^{1+o(1)}$. The two diagonal sums are bounded by $N^{o(1)}S_2$ and $N^{o(1)}S_4$. For the last line use the
sixth moment:
\[
 (\mathbb Ea_i^6)^{1/2}
 \le N^{o(1)}(|w_i|^3+N^{-1}+N^{-5/6}e_N).
\]
Choose deterministic bounds $d_i=Q_N(|w_i|^3+N^{-1}+N^{-5/6}e_N)$, $Q_N=N^{o(1)}$, valid under the law of each matrix $K_{ab}$. Then
\[
 \sum_i d_i\le Q_N(2+N^{1/6}e_N),\qquad
 \mathbb E[I_{ab}a_a^3a_b^3]
 \le(\mathbb E[I_{ab}a_a^6])^{1/2}
      (\mathbb E[I_{ab}a_b^6])^{1/2}\le d_ad_b.
\]
Here $I_{ab}$ is that pair's regularity indicator. Thus
\[
 \sum_{a<b}\mathbb E[I_{ab}a_a^3a_b^3]
 \le\tfrac12\left(\sum_i d_i\right)^2
 \le C Q_N^2(1+N^{1/3}e_N^2)\le N^{1/3+o(1)}.
\]
The bounds for $d_i$ are uniform over the matrices $K_{ab}$; hence their deterministic sums control all pairs.
\end{proof}

\subsubsection{Spectral confinement for entry replacement}
\label{W:sec:confinement-for-entry-replacement}

\label{W:sec:moment-increments-under-entry-replacement}

A single-entry comparison integrates the replacement coordinate against its original distribution. The regularity event used in that integral must therefore be determined by the
other entries. This requirement connects the spectral bounds to the condition on the Gram matrix: confinement and coordinate estimates are proved for every permitted value of the
selected entry. An upper bound on the Gram matrix at one such value then controls the resolvent at the starting vector throughout the segment of permitted entry values. Whether
such a value exists is measurable from the matrix with the selected entry set to zero.

The moment increments below are estimated on a regularity event and combined with the observable derivatives of Section~\ref{W:sec:cutoffs-in-moment-variables} in the later
comparison. Conditional on $K_\alpha$, the replacement coordinate $x$ has its prescribed distribution. The coefficients at $x=0$ are functions of $K_\alpha$, so their conditional
integrals are the corresponding signed moment differences. We use the event, measurable with respect to $K_\alpha$, that some permitted value of entry $\alpha$ gives a matrix
satisfying the upper bound on the Gram matrix. The resolvent identity transfers the bounds at any such matrix to every permitted value of that entry. If no such matrix exists, the
function multiplied by the cutoff in the eigenvalues of the Gram matrix vanishes for every permitted value of that entry.

For a row moment restricted to the regularity event, a second conditioning step is needed: whether a permitted entry value gives the upper bound on the Gram matrix may depend on
the row. Let $q=R_{K^{(i)}}w^{(i)}$ be the vector computed from the minor. Schur complementation expresses the coordinate of the resolvent applied to the starting vector as
$R_{K,ii}(w_i+h_i^*q)$. We first bound $R_{K,ii}$ on the regularity event, then enlarge that event to one measurable with respect to the minor and controlling $q$, and only then
integrate the independent row $h_i$. The resulting bound has the form $\mathbb E[I_i|(R_Kw)_i|^4] \le C\mathbb E_{K^{(i)}}[\mathbf1_{\mathcal C_i} \mathbb E_{h_i}|w_i+h_i^*q|^4]$.
The following resolvent estimates establish this domination and the resulting conditional row moment bounds.

Write $m_N=N^{1/3}$ and $\eta_{\mathrm e}=m_N^{-2}=h_N^2$, with
\[
 d=A_0\rho_*^{-2}(\log(N+2))^\gamma,\qquad 1<\gamma<2,\qquad
 c_*{m_N}\le \rho_*\le C_*{m_N}.
\]
All cutoffs and the finite set of pairs $(\sigma,E)$ are fixed. Every Gram matrix dimension is between $c_1{m_N}$ and $\rho_*$, for a fixed $c_1>0$. The deterministic unit starting
vector $w$ is arbitrary. A selected entry is $V_\alpha(x)=N^{-1/2}\mathcal T_\alpha x$, with at most two real coordinates and with support in $I_\alpha\times I_\alpha$,
$|I_\alpha|\le2$.

Use the modified matrix $H^{[d]}$, the functions $\Theta_a$, and coordinates $(u_n,x_n,v_{a,n},y_{a,n},t_a)$ defined in
Section~\ref{W:sec:the-moment-vector-and-its-smooth-extensions}, with $M=2\rho_*-1$. In particular, $t_a=\operatorname{Tr}\Theta_a(H)$ is unnormalized, and $\Theta_a$ is evaluated
at $H$.

We work with the truncated comparison arrays and with product measures formed from convex combinations of their corresponding entry laws. Uniformly over these arrays,
\begin{align}
 |H_{ij}|&\le e_N,\quad e_N\longrightarrow0,\quad e_N^{-1}=N^{o(1)},
 &\mathbb E|\sqrt N H_{ij}|^3&\le C\quad(i<j),\nonumber\\
 \mathbb E|\sqrt N H_{ij}|^4&\le F_N=N^{o(1)}\quad(i<j),
 &\mathbb E|\sqrt N H_{ii}|^2&\le C.
 \label{W:eq:resolvent-class}
\end{align}
The fourth power tail function is uniformly bounded. In particular, for a standardized off-diagonal entry $\xi$ supported in $|\xi|\le Ce_N\sqrt N$, tail integration gives
\begin{equation}
 \mathbb E|\xi|^4
 =4\int_0^{Ce_N\sqrt N}t^3\mathbb P(|\xi|>t)\,\dd t
 \le C+ C\log(eN).
 \label{W:eq:truncated-fourth-log-from-tail}
\end{equation}
The same estimate holds for the bounded and truncated Gaussian distributions and for convex combinations of these entry laws. The coarsely truncated matrices therefore satisfy the
logarithmic fourth moment hypothesis in Lemma~\ref{W:lem:resolvent-core}. The later comparison records these bounds through the factor $F_N=N^{o(1)}$. The truncated means and
covariance errors have the bounds proved in Sections~\ref{W:sec:moment-bounds-and-taylor-expansion-at-two-truncation-scales}
and~\ref{W:sec:errors-from-the-first-two-truncated-moments}; in particular the means of the off-diagonal matrix entries are $N^{-2+o(1)}$, and the standardized covariance errors
are $N^{-1+o(1)}$. The diagonal matrix means are $N^{-1+o(1)}$. Gaussian variables may be truncated with a superpolynomially small error in returning to the original ensemble. The
conclusions below concern the truncated interpolation. The coupling with the original matrix contributes $o(1)$ to the final comparison of bounded expectations.

We use Lemma~\ref{W:lem:complex-weak}, the decomposition in Section~\ref{W:sec:resolvents-on-exceptional-support}, the coordinate and counting estimates in
Section~\ref{W:sec:resolvent-domain-for-theta}, and the conditional quadratic form bound in Lemma~\ref{W:lem:schur-quadratic}. The additional external results are the finite-$N$
upper tail bound \cite{erdosxu2023small}, applied to a matrix with bounded standardized entries, and the matrix Chernoff bound \cite{tropp2012userfriendly} (Corollary~5.2). A high
moment trace comparison will give the polynomial failure probability required to sum over the matrices with the selected entry set to zero.

\Needspace{6\baselineskip}

\begin{lemma}
\label{W:lem:resolvent-core}
Fix $\phi\in(0,1/100]$. Let $W$ be a normalized coarsely truncated matrix with independent centered blocks, off-diagonal variance $N^{-1}$, zero complex pseudovariance and support
$CN^{-\phi}$. Suppose its standardized off-diagonal entries have uniformly bounded third absolute moments and fourth moments at most $C\log(eN)$, while its diagonal entries have
second moments at most $C/N$. For every fixed $b>0$ and $D>0$,
\[
 \mathbb P(\|W\|>2+bd)\le N^{-D}
\]
for all sufficiently large $N$. The bound is uniform over the intermediate matrices used to construct the truncated matrix.
\end{lemma}
\begin{proof}
Let $a=2+bd$, $f(x)=(1+\exp(-(x-a)/{\eta_{\mathrm e}}))^{-1}$, and $Z(B)=\operatorname{Tr}f(B)$. Thus $0\le Z\le N$, and $\lambda_{\max}(B)\ge a$ implies $Z(B)\ge1/2$. For the
normalized truncated matrix and a matrix with bounded standardized entries and matching off-diagonal third moment tensors, the preliminary estimates for the spectral norm and
resolvent entries give, with failure probability smaller than any prescribed inverse power of the matrix size, $\|B\|\le4$ and
$\max_{ij}|R_{B,ij}(x\pm\mathrm i{\eta_{\mathrm e}})-g(x\pm\mathrm i{\eta_{\mathrm e}})\delta_{ij}| \le N^{-c_0}$ for $|x|\le5$. For an intermediate matrix, write $K_\alpha$ for
the matrix with block $\alpha$ set to zero. Its regularity event is $\sigma(K_\alpha)$-measurable. The Woodbury identity restricted to the perturbed coordinates extends it to every
permitted value of the selected entry, replacing $N^{-c_0}$ by
$C(N^{-c_0}+N^{-\phi})$.

We compute the derivatives of $Z(B)=\operatorname{Tr}f(B)$ on the specified exterior spectral domain. Directions $T_j$ in this paragraph are unnormalized matrix entry directions of
norm at most one on one two-dimensional block. Contour calculus on the rectangle $[-5,5]+\mathrm i[-{\eta_{\mathrm e}},{\eta_{\mathrm e}}]$, cyclicity, and integration by parts
give
\begin{equation}
 D_{H}^q Z(B)[T_1,\ldots,T_q]
 =\frac1{2\pi\mathrm i q}\int f'(z)
 \sum_{\sigma\in S_q}\operatorname{Tr}
 [T_{\sigma(1)}R_B\cdots T_{\sigma(q)}R_B] \,\mathrm dz,
 \quad 1\le q\le4.
 \label{W:eq:resolvent-exterior-derivative}
\end{equation}
This is \eqref{W:eq:comparison-mixed-trace}. Applying Lemma~\ref{W:lem:comparison-exterior-derivatives} on the regularity event bounds it by the sum of the error in the compressed
resolvent, the integral over the semicircle cut, and the vertical-side contribution. Consequently the four derivatives are bounded by
\[
 \epsilon_N=C\{N^{-c}+\exp(-bd/{\eta_{\mathrm e}})+\exp(-1/{\eta_{\mathrm e}})\}\le N^{-c/2}
\]
for all $N\ge N_0(p)$. Here the derivatives are taken with respect to unnormalized matrix entries. Allowed changes of a matrix entry have length at most $CN^{-\phi}<1$, so
integration of the first derivative bounds the change of $Z$ along any such segment by
$C\epsilon_N$.

Apply Lemma~\ref{W:lem:bounded-complex-moment-matching} to each standardized complex entry, or its real version in Lemma~\ref{W:lem:cubic-defect}. The hypotheses hold uniformly
because the standardized third absolute moments are bounded. This gives bounded centered entries with matching full off-diagonal tensors through order three. Choose bounded
centered diagonal entries of standardized variance one.

Fix $p\ge4$ and put $Y=\max_\gamma\|Z(H^\gamma)\|_{L^p}$, where $H^\gamma$ ranges over all matrices in the sequential replacement. For a set of derivative indices $S$, the
polynomial chain rule is
\[
 D_S(Z^p)=\sum_{\pi\in\Pi(S),\ |\pi|\le p}
       (p)_{|\pi|} Z^{p-|\pi|}\prod_{B\in\pi}D_BZ,
 \qquad (p)_j=p(p-1)\cdots(p-j+1).
\]
At $K_\alpha$, the preceding derivative bounds give $|D_BZ|\le\epsilon_N$ for $1\le|B|\le4$. The same bounds hold along every permitted segment of entry values, and integration of
the first derivative gives $|Z(H_t)-Z(H_{t'})|\le C\epsilon_N$. Set $Z_\alpha^0=Z(K_\alpha)$, and let $\mathcal G_\alpha$ denote the regularity event measurable in the remaining
entries. On this event,
\[
 |D_S(Z^p)(H_t)|
 \le C_p\sum_{j=1}^{|S|}\epsilon_N^j
                         (Z_\alpha^0+C\epsilon_N)^{p-j}.
\]
Condition on all entries except the selected one. In the off-diagonal Taylor formula the coefficients through order three cancel. Writing the entry coordinate as
$u=(u_1,\ldots,u_d)$, $d\le2$, its degree-$q$ coefficient has conditional integral
\[
 \frac1{q!}\sum_{r_1,\ldots,r_q}
   \partial_{r_1}\cdots\partial_{r_q}[Z^p](K_\alpha)
   \int u_{r_1}\cdots u_{r_q}(\mu_\alpha-\nu_\alpha)(\dd u)=0,
       \qquad 0\le q\le3.
\]
The derivatives and $\mathbf1_{\mathcal G_\alpha}$ are fixed by the remaining entries. At $q=0$ this uses total mass zero; the other three identities use the matched full tensors.
Its fourth-order integral remainder is bounded by
\[
 C_p\int\|u\|_1^4(\mu_\alpha+\nu_\alpha)(\dd u)
       \sum_{j=1}^4\epsilon_N^j
                        (Z_\alpha^0+C\epsilon_N)^{p-j}.
\]
Here $u$ is the unnormalized entry coordinate. The moment coefficient is at most $C N^{-2}\log(eN)$. For a diagonal coordinate, the constant and linear terms cancel, and the
second-order integral remainder has coefficient $\int |x|^2(\mu+\nu)(\mathrm dx)\le C/N$ and contains the terms $j=1,2$ in the derivative bound. There are at most $N^2$
off-diagonal coordinates and $N$ diagonal coordinates, so these coefficients sum to $C\log(eN)$ and $C$. On the exceptional event for the regularity estimates, of probability
$N^{-D}$, use $0\le Z\le N$. The sum of the discarded bounded differences is at most $2N^{p+2-D}$, which is $O(N^{-A})$ for $D>A+p+3$.

To bound the expectation over the remaining entries, introduce an independent entry $u$ with the law of either adjacent matrix. For every value in its support,
$|Z(K_\alpha+\mathcal T_\alpha u)-Z_\alpha^0|\le C\epsilon_N$ on $\mathcal G_\alpha$. Thus for $0\le j\le p$,
\[
 \begin{aligned}
 &\mathbf1_{\mathcal G_\alpha}(Z_\alpha^0+C\epsilon_N)^{p-j}\\
 &\quad\le\int\mathbf1_{\mathcal G_\alpha}
   (Z(K_\alpha+\mathcal T_\alpha u)+2C\epsilon_N)^{p-j}
                     \,\mu_\alpha(\dd u),\\
 &\mathbb E[\mathbf1_{\mathcal G_\alpha}
                       (Z_\alpha^0+C\epsilon_N)^{p-j}]
       \le(Y+2C\epsilon_N)^{p-j}.
 \end{aligned}
\]
The last inequality follows from H\"older and Minkowski for $p-j>0$, and from $\mathbb P(\mathcal G_\alpha)\le1$ for $p-j=0$. The matrix after replacement has one of the two
distributions at adjacent steps of the comparison; the entry moment was integrated while $Z_\alpha^0$ was fixed. Telescoping to each intermediate matrix and maximizing gives
\begin{equation}
 Y^p\le Y_0^p+C_p\log(eN)
       \sum_{j=1}^{4}\epsilon_N^j(Y+C\epsilon_N)^{p-j}+N^{-A}.
 \label{W:eq:resolvent-core-moment}
\end{equation}
The equality of the low order moments was used under conditioning on the remaining entries, before either integral remainder was bounded.

The deviation estimate for bounded standardized entries at $M_A{\eta_{\mathrm e}}(\log N)^{2/3}$ gives $Y_0\le N^{-A}$, with arbitrary $A$ after adjusting its constants. Indeed,
$d/\eta_{\rm e}\ge c(\log N)^\gamma$ and $\gamma>1>2/3$, choose $N_0(A,b)$ so that $M_A(\log N)^{2/3}<bd/(2\eta_{\rm e})$ for $N\ge N_0(A,b)$. On the confinement event for the
matrix with bounded entries, $Z\le N e^{-bd/(2\eta_{\rm e})}$; on its complement $Z\le N$. Thus, if that complement has probability at most $N^{-D_0}$,
\[
 Y_0\le N e^{-bd/(2\eta_{\rm e})}+N^{1-D_0/p}.
\]
Choose $D_0>p(A+1)$. Also $\log N-c(\log N)^\gamma\le-A\log N$ for all sufficiently large $N$. These choices prove $Y_0\le2N^{-A}$, and replacing $A$ by $A+1$ gives the stated
bound. If $Y\ge4(Y_0+N^{-A/p}+C\epsilon_N)$, division of \eqref{W:eq:resolvent-core-moment} by $Y^p$ gives
\[
 \frac12\le C_p\log(eN)\sum_{j=1}^4(\epsilon_N/Y)^j
 \le2C_p\log(eN)\epsilon_N/Y,
\]
after enlarging the fixed constant so that $\epsilon_N/Y\le1/2$. Hence $Y\le4C_p\log(eN)\epsilon_N$. Combining this with the complementary threshold case proves
\begin{equation}
 Y\le C_p\bigl(Y_0+N^{-A/p}+\log(\mathrm e N)\epsilon_N\bigr).
 \label{W:eq:truncated-core-moment-solution}
\end{equation}
The exponent $c>0$ in the weak entrywise estimate is fixed by the support exponent $\phi$ and can be kept fixed when the failure exponent is increased. For a fixed $D>0$, first
choose the fixed integer $p>4(D+1)/c$. Then choose $A$ large enough in the estimates for bounded standardized entries and for the exceptional events. Since
$\epsilon_N\le N^{-c/2}$, (for this fixed $p$) \eqref{W:eq:truncated-core-moment-solution} gives $Y\le N^{-c/4}$ for all $N\ge N_0(p)$. Markov's inequality now yields
$\mathbb P(Z\ge1/2)\le 2^pY^p\le N^{-D}$ for large $N$. The hypotheses also hold for $-W$; the identity $\lambda_{\max}(-W)=-\lambda_{\min}(W)$ gives the lower edge bound. All
constants are uniform in the intermediate matrices; the moment order is selected before the large-$N$ limit.
\end{proof}

\Needspace{6\baselineskip}

For the next estimate use $\phi=1/100$ and the spectral width $d$ chosen above. In particular $d\ge cN^{-2/3}$ for large $N$; hence $N^{-1+4\phi}/d\le C N^{-1/3+4\phi}\le N^{-1/4}$
once $C\le N^{1/12-4\phi}$.

\begin{lemma}
\label{W:lem:resolvent-support-cher}
Let $W$ be independent of independent Bernoulli labels $L_e$ on unordered edges and diagonal vertices. Suppose their probabilities satisfy $p_{ij}\le CN^{-2+4\phi}$ and
$p_{ii}\le CN^{-1+2\phi}$. Let $S$ be the set of vertices incident to at least one marked label. Suppose a positive matrix $R$ is measurable with respect to $W$ and satisfies
\[
 \max_iR_{ii}\le C_0,\qquad \|R\|\le C/d.
\]
For each fixed $D$ there is a fixed $K_D$ such that, conditionally on $W$,
\[
 \mathbb P(\|R_{SS}\|>K_D\mid W)\le N^{-D}.
\]
\end{lemma}
\begin{proof}
Let $\Pi_e$ be the coordinate projection on an edge or a diagonal vertex. Let $J_S$ be the coordinate inclusion from $\mathbb C^S$ and set $T=R^{1/2}J_S$. Then $T^*T=R_{SS}$ and
$TT^*=R^{1/2}\Pi_SR^{1/2}$, whence both norms are $\|T\|^2$. Moreover, on each coordinate, $\mathbf1_{\{i\in S\}}\le\sum_{e\ni i}L_e$; this gives $\Pi_S\le\sum_eL_e\Pi_e$ in the
order of positive matrices. Conjugating by $R^{1/2}$ and taking norms therefore gives
\[
 \|R_{SS}\|=\|R^{1/2}\Pi_SR^{1/2}\|
 \le\left\|\sum_eL_eR^{1/2}\Pi_eR^{1/2}\right\|.
\]
For an edge or loop $e$, positivity gives
\[
 \|R^{1/2}\Pi_eR^{1/2}\|
 =\|\Pi_eR\Pi_e\|\le\operatorname{Tr}(\Pi_eR)
 =\sum_{i\in e}R_{ii}\le2C_0.
\]
Let $D_p=\sum_ep_e\Pi_e$. This is diagonal, with
\[
 (D_p)_{ii}=p_{ii}+\sum_{j\ne i}p_{ij}
 \le C(N^{-1+2\phi}+N^{-1+4\phi})\le C'N^{-1+4\phi}.
\]
Conditional independence and $\|R\|\le C/d$ therefore give
\[
 \mu:=\left\|\mathbb E\left[\sum_eL_eR^{1/2}\Pi_eR^{1/2}\mid W\right]\right\|
 =\|R^{1/2}D_pR^{1/2}\|
 \le C N^{-1+4\phi}/d\le N^{-1/4}
\]
for all $N$ beyond a deterministic index depending on the fixed constants. If $\mu=0$, positivity and zero expected trace make the random matrix zero almost surely. Suppose
$\mu>0$. Conditional on $W$, all matrices $R^{1/2}\Pi_eR^{1/2}$ are fixed, and the independent summands have eigenvalues in $[0,2C_0]$. Matrix Chernoff gives, for
$K\ge\mathrm e\mu$,
\[
 \mathbb P\left(\left\|\sum_eL_eR^{1/2}\Pi_eR^{1/2}\right\|\ge K
       \mid W\right)
 \le N\left(\frac{\mathrm e\mu}{K}\right)^{K/(2C_0)}.
\]
Choose $K_D=1+8C_0(D+2)$. For all sufficiently large $N$ it is larger than $\mathrm e\mu$, and the displayed probability is at most
\[
 \left(\frac{\mathrm e}{K_D}\right)^{K_D/(2C_0)}
       N^{1-K_D/(8C_0)}\le N^{-D}.
\]
The conditioning fixes every coefficient matrix; the only variables in this concentration bound are the independent labels $L_e$. The vertex indicators $\mathbf1_{\{i\in S\}}$ are
not independent, since one edge can mark two vertices. In particular, the bound above is applied to the independent sum indexed by edges and loops, not to
$\sum_i\mathbf1_{\{i\in S\}}R^{1/2}e_ie_i^*R^{1/2}$. For any event $\mathcal C$ measurable with respect to $W$ on which the hypotheses hold, conditional integration yields
\[
 \mathbb P(\mathcal C\cap\{\|R_{SS}\|>K_D\})
 =\mathbb E\!\left[\mathbf1_{\mathcal C}
       \mathbb P(\|R_{SS}\|>K_D\mid W)\right]\le N^{-D}.
\]
The conditional estimate depends on the support variables and $W$; it holds for every conditional distribution of the marked entry values compatible with that construction.
\end{proof}

\Needspace{6\baselineskip}

\begin{theorem}
\label{W:thm:resolvent-window}
For every fixed $b>0,D>0$, the arrays in \eqref{W:eq:resolvent-class} satisfy
\begin{equation}
 \mathbb P(\|H\|>2+bd)\le N^{-D}
 \label{W:eq:resolvent-window}
\end{equation}
for all sufficiently large $N$, uniformly over the product laws specified above and any fixed list of polynomially many intermediate matrices.
\end{theorem}
\begin{proof}
Use the construction with independent Bernoulli labels at standardized threshold $T=N^{1/2-\phi}$, with $\phi=1/100$, as in Section~\ref{W:sec:resolvents-on-exceptional-support}:
\[
 H=W+V+E.
\]
Here $W$ is constructed from independent variables conditioned to lie below the truncation threshold. All variables are centered, and the covariance matrices of the off-diagonal
variables are normalized in their real coordinates. The labels generating $E$ are independent of $W$, though the values of $E$ need not be independent of it. For a standardized
off-diagonal block $\xi$, let $p=\mathbb P(|\xi|>T)$ and $\mu=\mathbb E\xi$. Its mean and second moment conditional on lying below the truncation threshold are
\[
 \mu_< =\frac{\mu-\mathbb E[\xi\mathbf1_{\{|\xi|>T\}}]}{1-p},
 \quad
 \mathbb E[\xi\xi^{\mathsf T}\mid |\xi|\le T]
 =\frac{\mathbb E\xi\xi^{\mathsf T}
       -\mathbb E[\xi\xi^{\mathsf T}\mathbf1_{\{|\xi|>T\}}]}{1-p}.
\]
The uniform fourth power tail function gives $p\le CT^{-4}$, $\mathbb E[|\xi|\mathbf1_{\{|\xi|>T\}}]\le CT^{-3}$, and $\mathbb E[|\xi|^2\mathbf1_{\{|\xi|>T\}}]\le CT^{-2}$.
Together with the bounds $|\mu|\le N^{-3/2+o(1)}$ and $\|\operatorname{Cov}\xi-\Sigma_0\|\le N^{-1+o(1)}$, these imply
\[
 |\mu_<|\le CN^{-3/2+3\phi},\qquad
 \|\Sigma_<-\Sigma_0\|\le CN^{-1+2\phi}.
\]
Here $\Sigma_0$ is the scalar or two-dimensional covariance matrix, whose eigenvalues are fixed positive constants. The normalized centered small block $\zeta$ can be defined by
$\xi_<=\mu_<+L\zeta$, with $\mathbb E\zeta\zeta^{\mathsf T}=\Sigma_0$ and $\|L-I\|\le CN^{-1+2\phi}$. One choice is
\[
 \zeta=\Sigma_0^{1/2}\Sigma_<^{-1/2}(\xi_<-\mu_<),\qquad
 L=\Sigma_<^{1/2}\Sigma_0^{-1/2}.
\]
Then $\mathbb E\zeta\zeta^{\mathsf T}=\Sigma_0$ and the stated reconstruction is exact. To bound $L-I$, set $A=\Sigma_<^{1/2}$, $B=\Sigma_0^{1/2}$, and $X=A-B$. They satisfy
$AX+XB=\Sigma_<-\Sigma_0$. Since the least eigenvalues of $A,B$ are bounded below by a fixed $c>0$,
\[
 X=\int_0^\infty e^{-tA}(\Sigma_<-\Sigma_0)e^{-tB}\,\dd t,
 \qquad \|X\|\le(2c)^{-1}\|\Sigma_<-\Sigma_0\|.
\]
Differentiating the integrand verifies the equation, and the exponential bound verifies convergence. Multiplication by $\Sigma_0^{-1/2}$ proves $\|L-I\|\le C N^{-1+2\phi}$. The
off-diagonal normalization correction therefore has centered independent summands $N^{-1/2}(L-I)\zeta$ of norm at most $CN^{-1+\phi}$, because $|\zeta|\le CT$. Their matrix
variance parameter obeys
\[
 v\le CN\cdot N^{-1}\|L-I\|^2\le CN^{-2+4\phi}.
\]
The bound is uniform in the pair even when $L$ varies with the pair. For the diagonal, the second moment tail estimate gives conditional mean $O(T^{-1})+N^{-1/2+o(1)}$; the
diagonal conditional distribution is centered without rescaling its variance. Hence the deterministic mean matrix $V_{\rm mean}$ satisfies
\[
 \|V_{\rm mean}\|\le
 N\max_{i<j}N^{-1/2}|\mu_{<,ij}|
       +\max_iN^{-1/2}|\mu_{<,ii}|
 \le CN^{-1+3\phi}+CN^{-1+\phi}.
\]
Matrix Bernstein, with $L_0=CN^{-1+\phi}$ and $t=\tfrac12N^{-1+4\phi}$, gives
\[
 \mathbb P(\|V-V_{\rm mean}\|>t)
 \le2N\exp\left(-\frac{t^2}{2v+2L_0t/3}\right)
 \le2N e^{-cN^{3\phi}}.
\]
Here $t^2\asymp N^{-2+8\phi}$ and $v+L_0t\le CN^{-2+5\phi}$, so the exponent is at least $cN^{3\phi}$. The mean norm is $o(t)$. It follows that $\|V\|\le N^{-1+4\phi}=o(d)$ outside
an event smaller than every fixed inverse power of $N$. This calculation includes the changes in the moments caused by truncation of \eqref{W:eq:resolvent-class} before the
conditional normalization.

For a failure exponent $D>0$, choose an integer $d_0>(D+2)/(1-4\phi)$. A component of size at least $d_0+1$ contains a tree on that many vertices. The union bound over its vertices
and labeled trees gives
\[
 \begin{aligned}
 &\mathbb P(\text{a marked off-diagonal component has size}>d_0)\\
 &\qquad\le N^{d_0+1}(d_0+1)^{d_0-1}(CN^{-2+4\phi})^{d_0}\\
 &\qquad=C_{d_0}N^{1-d_0(1-4\phi)}.
 \end{aligned}
\]
On its complement, diagonal marks do not enlarge components, and
\[
 \|E\|\le C d_0(e_N+N^{-\phi})\longrightarrow0.
\]
This bounds $E$ on the event that every marked component has at most $d_0$ vertices.

By Lemma~\ref{W:lem:resolvent-core}, arrange $\|W\|\le2+bd/8$ with failure probability smaller than any prescribed inverse power of the matrix size. At $a_\pm=\pm(2+bd/4)$ the
matrices $R_W(a_+)$ and $-R_W(a_-)$ are positive, have norm at most $C_b/d$, and have diagonal entries at most a fixed constant. For the diagonal bound compare the real resolvent
to $a_\pm+\mathrm i{\eta_{\mathrm e}}$: all spectral distances are at least $bd/8$, and
\[
 r^{-1}\le(1+(8{\eta_{\mathrm e}}/(bd))^2)\,r/(r^2+{\eta_{\mathrm e}}^2),\qquad r\ge bd/8.
\]
The coordinate estimate at the nonreal parameter and this inequality give the stated diagonal bound.

Apply Lemma~\ref{W:lem:resolvent-support-cher} at both barriers. On the event given by that lemma, choose $N_0(D)$ so that $\|E\|K_D<1/2$ for $N\ge N_0(D)$. For every
$t\in[0,1]$,
\[
 \det(a_\pm I-W-tE)
 =\det(a_\pm I-W)\det(I-tE_SR_W(a_\pm)_{SS})\ne0.
\]
For an ordered eigenvalue $\lambda_j(t)$ of $W+tE$,
\[
 |\lambda_j(t)-\lambda_j(t')|\le|t-t'|\|E\|,\qquad
 \lambda_j(0)\in(a_-,a_+),\qquad \lambda_j(t)\ne a_\pm.
\]
The intermediate value theorem gives $\lambda_j([0,1])\subset(a_-,a_+)$. The fixed compression bound and the component estimate give explicitly
$\|E_SR_W(a_\pm)_{SS}\|\le C_{d_0}K_D(e_N+N^{-\phi})\to0$, which is the bound used in the determinant. Consequently
\[
 \|H\|\le2+bd/4+N^{-1+4\phi}\le2+bd
\]
for sufficiently large $N$. For a family of at most $N^a$ matrices and a fixed exponent $D$, apply all preceding probability estimates with exponent $D+a+2$. Their constants are
uniform over the specified product laws of the entries, so at each fixed interpolation parameter
\[
 \mathbb P\!\left(\max_{H\text{ in the family}}\|H\|>2+bd\right)
 \le C N^aN^{-D-a-2}\le N^{-D}
\]
for every $N\ge N_0(D)$. This proves the stated uniform probability bound.
\end{proof}

\Needspace{6\baselineskip}

\begin{lemma}
\label{W:lem:resolvent-coordinate-energy}
For the truncated family, uniformly on compact real energies and $c{\eta_{\mathrm e}}\le\eta\le M_0$, put
\[
 X(z)=\max_i[-\operatorname{Im}R_{ii}(z)],\qquad
 \bar\varrho(z)=C\upsilon(z)+N^{o(1)}/(N\eta).
\]
With failure probability smaller than any prescribed inverse power of the matrix size,
\begin{align}
 X(z)\le N^{o(1)}\Big[&
 \eta+\bar\varrho+\sqrt{\bar\varrho/(N\eta)}
 +F_N\bar\varrho/N+b_N^2/(N\eta)\nonumber\\
 &+u_N(N^{-1/2}+b_N)/\eta+u_N^2/\eta\Big],
 \label{W:eq:resolvent-imaginary}
\end{align}
where $b_N=CN^{-5/24}$ and $u_N=N^{-7/8+o(1)}$. The factors $N^{o(1)}$ mean arbitrarily small fixed power losses, uniform over the finite spectral sets of polynomial cardinality
and the fixed parameter lists used here.
\end{lemma}
\begin{proof}
Use the auxiliary graph of Section~\ref{W:sec:block-schur-complements-under-the-tail-assumption}, at standardized threshold $N^{7/24}$. Choose $d_1>6(D+2)$. By the component bound,
the probability that some component has more than $d_1$ vertices is at most $C_{d_1}N^{1-d_1/6}\le N^{-D}$ whenever $N^{d_1/6-D-1}\ge C_{d_1}$. Condition on the labels and a
complementary minor $G=(z-H_{I^cI^c})^{-1}$. The boundary rows $B=H_{I,I^c}$ have independent centered parts, support $b_N$, variance at most $C/N$, fourth moment at most
$F_N/N^2$, and norm of the mean row at most $u_N$.

Let $S=-\operatorname{Im}G\ge0$ and $D_G=\max_j S_{jj}$. Its norm is at most $\eta^{-1}$ and $N^{-1}\operatorname{Tr}S\le\bar\varrho$, by the counting estimate on the real line and
minor interlacing. Positivity gives
\[
 \|S\|_{\rm HS}\le\sqrt{N\bar\varrho/\eta},\quad
 \max_j\|Se_j\|\le\sqrt{D_G/\eta},\quad
 \max_{ij}|S_{ij}|\le D_G,\quad
 \|\operatorname{diag}S\|_2\le\sqrt{D_GN\bar\varrho}.
\]
Apply the conditional estimate for the quadratic forms $u^*Su$, $u\in\mathbb C^{I^c}$, in Section~\ref{W:sec:conditional-moments-of-schur-complement}. One may use diagonal
quadratic forms only, since $\|BSB^*\|\le\operatorname{Tr}(BSB^*)$. For every arbitrarily small fixed $\varepsilon>0$, with a sufficiently high conditional moment order,
\begin{align*}
 \|BSB^*\|\le C_{d_1}\Big[\bar\varrho+N^\varepsilon\bigl(&\sqrt{F_ND_G\bar\varrho/N}+\sqrt{\bar\varrho/(N\eta)}\\
 &+b_N\sqrt{D_G/(N\eta)}+b_N^2D_G
 +u_N(N^{-1/2}+b_N)/\eta+u_N^2/\eta\bigr)\Big].
\end{align*}
Here the conditioning is on the graph labels and $H_{I^cI^c}$. For these fixed quantities, write $\mathcal D(S)$ for the sum inside the large parentheses in the preceding display.
The conditional estimate of order $p$ for each diagonal quadratic form has $L^p$ norm at most $C_{p,d_1}\mathcal D(S)$ after subtracting its conditional mean. The mean is bounded
by $C\bar\varrho$ plus the displayed contribution of the mean row. Markov's inequality at $N^\varepsilon\mathcal D(S)$ consequently bounds the conditional failure probability by
$C_{p,d_1}N^{-\varepsilon p}$. Positivity gives
\[
 \|BSB^*\|\le\sum_{i\in I}(BSB^*)_{ii},
\]
so at most $d_1$ such probabilities are needed for a component.

Suppose the prescribed family contains at most $N^a$ matrices and the finite set of spectral parameters has at most $N^b$ elements. Each graph has at most $N$ components. First
choose $d_1$ so that the event that a connected component exceeds this size has probability at most $N^{-D-a-2}$ for each graph. Then choose $p$ with $\varepsilon p>D+a+b+3$. The
union of the events where the conditional quadratic form estimates fail has probability at most
\[
 C_{p,d_1}d_1N^{a+b+1-\varepsilon p}\le N^{-D}
\]
for sufficiently large $N$. This bound holds for every value of the labels for which component sizes are at most $d_1$, and for every complementary minor satisfying the stated
bounds. Integrating over these variables preserves it. The full-matrix estimate for $X$ is used only in the deterministic substitution below, after this conditional calculation.

On the truncated coordinate event of Section~\ref{W:sec:resolvent-domain-for-theta}, for every set $I$ whose cardinality is bounded independently of $N$, both the compression
$R_{II}$ and its inverse have bounded operator norm: $R_{II}=gI+o(1)$ and $|g|$ is bounded away from zero on the fixed energy domain. The exact minor identity
\[
 G=R_{I^cI^c}-R_{I^cI}R_{II}^{-1}R_{II^c}
\]
and the Ward identity imply $D_G\le C_{d_1}X$. Explicitly, for $j\notin I$, $Ge_j$ is a restriction of $Re_j$ minus a linear combination of at most $d_1$ columns $Re_i$, $i\in I$.
The coefficients $R_{II}^{-1}R_{Ij}$ are bounded. Thus $\|Ge_j\|^2\le C_{d_1}X/\eta$ and $-\operatorname{Im}G_{jj}=\eta\|Ge_j\|^2\le C_{d_1}X$. Taking imaginary parts in the Schur
complement identity gives
\[
 -\operatorname{Im}R_{II}
 =R_{II}^*(\eta I+BSB^*)R_{II}.
\]
Apply the identity $\operatorname{Im}(M^{-1})=-M^{-*}(\operatorname{Im}M)M^{-1}$ with $M=R_{II}^{-1}$. In particular,
\[
 \max_{i\in I}(-\operatorname{Im}R_{ii})
 \le\|R_{II}\|^2(\eta+\|BSB^*\|).
\]
Since every coordinate belongs to one of the components, taking the maximum over the preceding identities for the imaginary parts gives
\begin{equation}
 \begin{split}
 X\le C(\eta+\bar\varrho)+C N^\varepsilon\bigg[
 &\sqrt{F_N\bar\varrho/N}\,\sqrt X
       +\sqrt{\bar\varrho/(N\eta)}
       +\frac{b_N}{\sqrt{N\eta}}\sqrt X+b_N^2X\\
 &+\frac{u_N(N^{-1/2}+b_N)}{\eta}
       +\frac{u_N^2}{\eta}\bigg].
 \end{split}
 \label{W:eq:coordinate-energy-absorption}
\end{equation}
The constant depends on the fixed bound on the component size but not on $X$. Choose $\varepsilon>0$ arbitrarily small, in particular $\varepsilon<5/24$. Since
$b_N^2=O(N^{-5/12})$, its coefficient of $X$ is at most $1/4$ for sufficiently large $N$. For each of the two terms $\alpha\sqrt X$, use $\alpha\sqrt X\le X/8+2\alpha^2$. They
contribute altogether at most $X/4$, and their remaining bounds are
\[
 C N^{2\varepsilon}\frac{F_N\bar\varrho}{N}
   +C N^{2\varepsilon}\frac{b_N^2}{N\eta}.
\]
At most $X/2$ is now on the right, so it can be moved to the left. This proves \eqref{W:eq:resolvent-imaginary}. Given $\delta,D>0$, first choose $0<\varepsilon<\delta/2$, then
choose the conditional moment order so that its exceptional probability is at most $N^{-D}$ after the finite unions. To extend the bounds from a finite covering set to the entire
domain, put $\eta_*=c\eta_{\rm e}$ and choose a rectangular array of spectral parameters of spacing at most $\eta_*/10$; its cardinality is polynomial in $N$. Increase the failure
exponent used at each point by this fixed cardinality exponent. For any $z=u+\mathrm i\eta$ choose a selected point $z'=u'+\mathrm i\eta'$ with $|z-z'|\le\eta'/4$. For real $x$,
\[
 \tfrac34|x-z'|\le|x-z|\le\tfrac54|x-z'|,
 \qquad \tfrac34\eta'\le\eta\le\tfrac54\eta'.
\]
Thus the Poisson kernels $\eta/|x-z|^2$ and $\eta'/|x-z'|^2$ have ratio between $1/8$ and $8$. Integration against a coordinate spectral measure gives $X(z)\le8X(z')$. Integration
against the semicircle measure gives $\upsilon(z')\le8\upsilon(z)$, and $(N\eta')^{-1}\le2(N\eta)^{-1}$. Every nonnegative summand on the right of \eqref{W:eq:resolvent-imaginary}
at $z'$ is consequently bounded by a fixed multiple of its value at $z$. The estimate on the selected spectral parameters therefore holds throughout the domain, with a larger
constant. Adjoints give the same estimate in the lower half-plane.
\end{proof}

\subsubsection{Contour integrals for moment increments}
\label{W:sec:contour-integrals-for-moment-increments}
\label{W:sec:moment-expansion-and-trace-normalization}

The conditional resolvent bounds are now available for the matrices along each permitted segment of entry values. Their conversion into moment increments requires a contour that
encloses these spectra and avoids the poles of $\Theta_a$. The ellipse used for polynomial moments is narrowed near the two edges for this purpose.

With the spectral and coordinate bounds of Section~\ref{W:sec:moment-increments-under-entry-replacement}, the parameters are $E_0=1+\max_{(\sigma,E)\in\mathcal I}|E|$, with $E_0=1$
for an empty label set, and $s=4\sqrt d$. A fixed $0<a_0<\pi/(16(E_0+2))$ determines the permitted height. The ellipse itself is described by

\[
 b_s=2\cosh s,\qquad
 y_s(x)=2\sinh s\sqrt{1-x^2/b_s^2},\qquad |x|\le b_s,
\]
\begin{equation}
 y(x)=\min\{y_s(x),\ a_0({\eta_{\mathrm e}}+|x-2|),\ a_0({\eta_{\mathrm e}}+|x+2|)\}.
 \label{W:eq:resolvent-contour}
\end{equation}
Let $\mathcal C_N$ be the positively oriented closed curve formed by $x\pm\mathrm i y(x)$, $-b_s\le x\le b_s$. It is rectifiable and piecewise smooth; at its two real tips it is
the original ellipse. Corners may be treated piecewise in contour integration.

\Needspace{6\baselineskip}

\begin{lemma}
\label{W:lem:resolvent-contour}
This contour encloses $[-2-d/8,2+d/8]$ and excludes every pole of every $z\mapsto\theta_E({m_N}^2(2-\sigma z))$. On an event independent of the starting vector, with failure
probability smaller than any prescribed inverse power of the matrix size, it has uniformly bounded coordinate resolvents, diagonals bounded away from zero, and
\begin{equation}
 \max_i\{\|R_H(z)e_i\|^2+\|R_H(\bar z)e_i\|^2\}
 \le {m_N}N^{o(1)},\qquad z\in\mathcal C_N.
 \label{W:eq:resolvent-column}
\end{equation}
The same estimates and $\|H+V\|<2+d/8$ hold simultaneously for every single block perturbation $V$ with $\|V\|\le C e_N$. There is a deterministic positive measure $\nu_N$ on
$\mathcal C_N$, with mass $N^{o(1)}$ and bounded below by a fixed multiple of arclength, which dominates the absolute Cauchy weights of all $p_n$ and
$p_n\theta_E({m_N}^2(2-\sigma z))$, $0\le n\le2\rho_*-1$. Moreover,
\begin{align}
 \int_{\mathcal C_N}|\theta_E({m_N}^2(2-\sigma z))|\,|\mathrm dz|
 &\le C_Ed,\label{W:eq:resolvent-soft-mass}\\
 \int_{\mathcal C_N}|p_n(z)\theta_E({m_N}^2(2-\sigma z))|\,|\mathrm dz|
 &\le {m_N}dN^{o(1)}.\label{W:eq:resolvent-mixed-mass}
\end{align}
\end{lemma}
\begin{proof}
At the real part of a pole of $\Theta_a$, $x=\sigma(2-E{\eta_{\mathrm e}})$, the height is at most $a_0{\eta_{\mathrm e}}(1+|E|)<\pi {\eta_{\mathrm e}}/2$. Thus even the first
poles at heights $\pm\pi {\eta_{\mathrm e}}$ lie outside the contour. Put $u=m_N^2(2-\sigma x)-E$ and $v=-\sigma m_N^2\operatorname{Im}z$. If $|u|\le1$, then
$|x-2\sigma|\le(1+|E|)\eta_{\rm e}$, so $|v|\le a_0(2+|E|)<\pi/2$. In this case $|1+e^{u+\mathrm iv}|\ge(1+e^u)\cos(v/2)\ge(1+e^u)/\sqrt2$. For $u\ge1$, the bound
$|1+e^{u+\mathrm iv}|\ge e^u-1\ge(1-e^{-1})e^u$ holds; for $u\le-1$ it is at least $1-e^{-1}$. These three cases give
\[
 |\Theta_a(z)|\le C e^{-2u_+}
 \le C e^{2|E|}\exp[-2(2-\sigma x)_+/\eta_{\rm e}].
\]
This proves both the uniform value bound and the decay on the interior side of the chosen edge.

Let $v=|x-2|$ near the right edge. The ellipse with parameter $s=4\sqrt d$ satisfies
\[
 y_{\rm ell}(x)^2
 =4\sinh^2s\left(1-\frac{x^2}{4\cosh^2s}\right)
 \le C d(d+v).
\]
If the linear height $a_0(\eta_{\mathrm e}+v)$ replaces the elliptic height, then $a_0^2(\eta_{\mathrm e}+v)^2\le C d(d+v)$. For $v\ge d$ this implies $v\le C_{a_0}d$; for $v<d$
the same conclusion holds. Reflection gives the same geometric bound at the left edge. Thus all modified portions lie within $C_{a_0}d$ of an edge and satisfy $|x\mp2|\le y/a_0$,
with $y\ge a_0{\eta_{\mathrm e}}$. The remaining nonreal elliptic points either have $y\ge c_{a_0}d$ or lie at the exterior real tips. For nonreal elliptic points,
$\upsilon(z)/\eta=e^{-s}/(2\sinh s)\le C/s$. For modified points $\upsilon(z)/\eta\le C\eta^{-1/2}\le C{m_N}$. At either type with $\eta\ge a_0{\eta_{\mathrm e}}$,
$\bar\varrho/\eta\le {m_N}N^{o(1)}$.

Divide \eqref{W:eq:resolvent-imaginary} by $\eta$, using $\eta\ge a_0N^{-2/3}$, $b_N=4N^{-5/24}$, $u_N\le N^{-7/8+o(1)}$ and $F_N=N^{o(1)}$. The square-root contribution satisfies
\[
 \frac1\eta\sqrt{\frac{\bar\varrho}{N\eta}}
 =\sqrt{\frac{\bar\varrho/\eta}{N\eta^2}}
 \le N^{1/3+o(1)}.
\]
The remaining terms are bounded, in their respective scales, by
\[
 \begin{aligned}
 \bar\varrho/\eta&\le N^{1/3+o(1)},&
 F_N\bar\varrho/(N\eta)&\le N^{-2/3+o(1)},\\
 b_N^2/(N\eta^2)&\le N^{-1/12+o(1)},&
 u_NN^{-1/2}/\eta^2&\le N^{-1/24+o(1)},\\
 u_Nb_N/\eta^2&\le N^{1/4+o(1)},&
 u_N^2/\eta^2&\le N^{-5/12+o(1)}.
 \end{aligned}
\]
The term $\eta$ in \eqref{W:eq:resolvent-imaginary} contributes $1$ after division by $\eta$. Each quantity is bounded by $m_NN^{o(1)}$. Finally Ward's identity
$\|R(z)e_i\|^2=-\operatorname{Im}R_{ii}(z)/\eta$ proves \eqref{W:eq:resolvent-column}.

Use Theorem~\ref{W:thm:resolvent-window} with window $d/256$. On the portion with small imaginary part next to the right real tip, write $z=x+\mathrm i y$. The ellipse formula
gives $x-(2+d/256)\ge c d$ and $|y|\le c_1d$. Choose $h=c_2d$ with $c_2\ge c_1$ fixed. At $z_h=x+\mathrm i h$, the already proved nonreal estimate is
$\|R_H(z_h)e_i\|^2\le m_NN^{o(1)}$. Equation~\eqref{W:eq:vertical-energy-comparison}, with $\delta=cd$, transfers this estimate to $z$. Moreover, writing $r=x-\lambda\ge cd$, the
scalar inequalities
\[
 \frac1r\le\left(1+\frac{h^2}{c^2d^2}\right)
                       \frac r{r^2+h^2},\qquad
 \frac r{r^2+y^2}\ge
       \left(1+\frac{c_1^2}{c^2}\right)^{-1}\frac1r
\]
show that $R_H(x)_{ii}\le C$ and $\operatorname{Re}R_H(z)_{ii}\ge c'>0$. For the lower bound, use $0<r\le6$ on this spectral window, after increasing a fixed constant. Also
$|R_H(z)|\le R_H(x)$ in spectral order, and therefore
\[
 |R_H(z)_{ij}|^2
 \le(e_i^*|R_H(z)|e_i)(e_j^*|R_H(z)|e_j)\le C^2.
\]
Reflection gives the same resolvent bounds at the left tip. On the other pieces the coordinate approximation and $|g|\ge c$ give the lower diagonal bound directly.

At $a_\pm=\pm(2+d/32)$, the distance to the unperturbed spectrum is at least $c d$. The same real comparison at height
$\eta_{\rm e}$ gives $\|R_H(a_\pm)_{II}\|\le C$ on every support $|I|\le2$. Consequently, for $\|T\|\le Ce_N$ and $0\le
t\le1$,
\[
 \|tTR_H(a_\pm)_{II}\|\le C'e_N<\tfrac12,
 \qquad \det(I-tTR_H(a_\pm)_{II})\ne0.
\]
The determinant identity and continuity of the eigenvalues place the spectrum of every matrix $H+t\iota_IT\iota_I^*$ between these barriers. For any point of the contour, the
Woodbury identity gives the exact column correction
\[
 R_{H+\iota_IT\iota_I^*}(z)e_j-R_H(z)e_j
 =R_H(z)_{\cdot I}T(I-R_H(z)_{II}T)^{-1}R_H(z)_{Ij}.
\]
Its entries are $O(e_N)$, and its Euclidean norm is at most $C e_N(\sum_{i\in I}\|R_H(z)e_i\|^2)^{1/2}$. Thus all entries stay bounded, the diagonals remain at least $c'/2$ in
modulus, and the squared column norms remain at most $m_NN^{o(1)}$, uniformly for $\|T\|\le Ce_N$. The same formula with the opposite perturbation includes deletion of the selected
entry. A finite covering set with spacing $N^{-20}$ suffices for the contour: the spectral distance there is at least $N^{-1}$ for large $N$, so the entry and squared column norm
Lipschitz constants are at most $CN^2$ and $CN^3$, respectively, by Lemma~\ref{W:lem:vertical-resolvent-comparison}. The error between a spectral parameter and its nearest point in
this set is $O(N^{-17})$. If the sets of spectral points and coordinate vector pairs have at most $N^J$ elements, use pointwise failure exponent $D+J+1$. Their union has
probability at most $N^{-D-1}$; the deterministic $O(N^{-17})$ estimate then extends the bounds from the finite covering set to the contour.

For the scalar bound, if $z$ is inside the ellipse, write $z=\zeta+\zeta^{-1}$ with $1\le|\zeta|\le e^s$. Thus
\[
 |p_n(z)|\le\frac{2e^{(2\rho_*+1)s}}{|z^2-4|^{1/2}}.
\]
On the ellipse $z=\zeta+\zeta^{-1}$, $z^2-4=(\zeta-\zeta^{-1})^2$ and $|\dd z|=|\zeta-\zeta^{-1}|\,\dd\theta$. Thus its full integral of $|z^2-4|^{-1/2}|\dd z|$ is $2\pi$. On a
clipped straight segment next to $2\sigma$, write $u=|x-2\sigma|$. Its slope has modulus $a_0$, and $|z-2\sigma|\ge c_{a_0}(\eta_{\rm e}+u)$, while $|z+2\sigma|\ge c$ for all
sufficiently large $N$. Its contribution is consequently at most
\[
 C\int_0^{Cd}(\eta_{\rm e}+u)^{-1/2}\,\dd u
 =2C\bigl(\sqrt{\eta_{\rm e}+Cd}-\sqrt{\eta_{\rm e}}\bigr)
 \le C'\sqrt d.
\]
On each of the four real intervals cut by $-2,0,2$, the smaller cap has the form $\alpha+\beta u$ with fixed coefficients. Substituting this height in the equation of the ellipse
gives a quadratic equation in $u$, hence at most two intersection points. There are therefore at most two straight subsegments in each such interval, and at most sixteen on the
upper and lower halves together. Hence a constant multiple of
\[
 \left(1+\frac{2e^{(2\rho_*+1)s}}{|z^2-4|^{1/2}}\right)
 \frac{|\mathrm dz|}{2\pi}
\]
is the required common measure, since $\rho_*s=O((\log N)^{\gamma/2})$. Its mass is uniform over the polynomial degrees.

Decompose the curve into the exterior pieces at the chosen edge, the adjacent interior pieces parameterized by $v=(2-\sigma\operatorname{Re}z)_+$, and the remaining pieces. The
first have total length at most $Cd$, the second have bounded slope, and on the third $v\ge cd$. Thus
\begin{align*}
 \int_{\mathcal C_N}|\theta_E(m_N^2(2-\sigma z))|\,|\mathrm dz|
 &\le C_Ed+C_E\int_0^{Cd}e^{-cv/\eta_{\mathrm e}}\,\mathrm dv
                                  +C_Ee^{-cd/\eta_{\mathrm e}}\\
 &\le C_E(d+\eta_{\mathrm e})+O(N^{-D})\le C_Ed
\end{align*}
for every fixed $D$, since $d/\eta_{\mathrm e}\asymp(\log N)^\gamma$ and $\gamma>1$. Finally,
\[
 \int_{\mathcal C_N}|p_n(z)\theta_E(m_N^2(2-\sigma z))|\,|\mathrm dz|
 \le (n+1)e^{ns}C_Ed\le m_NdN^{o(1)},
\]
which proves both mass estimates. The two restrictions on $\gamma$ enter separately. For fixed $A,C,c>0$, the bounds on $s,d$ give
\[
 \frac{C\rho_*s}{\log N}=O((\log N)^{\gamma/2-1})\longrightarrow0,
 \qquad A\log N-cd/\eta_{\mathrm e}\longrightarrow-\infty.
\]
The first limit uses $\gamma<2$ and permits every fixed product of the Chebyshev factors to be bounded by $N^\varepsilon$ for all $N\ge N_0(\varepsilon)$, for any fixed
$\varepsilon>0$. The second uses $\gamma>1$ and gives $N^A e^{-cd/\eta_{\mathrm e}}\to0$. These conclusions hold after all fixed derivative orders and parameters of the functions
under consideration have been chosen; neither is uniform in a moment order growing with $N$.
\end{proof}

Let $\mathcal X_\alpha$ be the closed ball of block values with $\|V_\alpha(x)\|\le 3Ce_N$, enlarged so that it contains the supports of the replacement entry distributions in its
interior. The regularity event $\mathcal B_\alpha(K)$ means that the spectral and coordinate conclusions of Lemma~\ref{W:lem:resolvent-contour} hold on this entire ball, with fixed
bounds $Q_N=N^{o(1)}$ for the column estimates and for their required vertical extensions. It is measurable with respect to $K$, and
\[
 \sup_{t,\alpha}\mathbb P(\mathcal B_\alpha(K_\alpha(t))^c)\le N^{-D}
\]
for every fixed $D$, after fixing the small exponent losses. The common event for $H(t)$ is uniform in the starting vector; the bounds for perturbations on finitely many
coordinates then give uniformity over the permitted entry values.

For a Gram matrix dimension $\rho$, define
\begin{equation}
 \mathcal A_{\alpha,\rho}(K,w)
 =\left\{\text{there is }x\in\mathcal X_\alpha\text{ with }
 G_\rho(K+V_\alpha(x),w)\le\kappa^{-1}I\right\}.
 \label{W:eq:resolvent-gram-bound}
\end{equation}
On $\mathcal B_\alpha$, these are the Gram matrices of the original polynomial vectors, since spectral truncation is inactive. For measurability use the modified Gram matrix on all
matrices and put
\[
 g(K,w)=\min_{x\in\mathcal X_\alpha}
  \lambda_{\max}\bigl(G_\rho^{[d]}(K+V_\alpha(x),w)\bigr).
\]
The ball $\mathcal X_\alpha$ is compact and the function being minimized is continuous. If $D_\alpha$ is a fixed countable dense subset of that ball, the minimum equals the infimum
over $D_\alpha$, and
\[
 \{g\le\kappa^{-1}\}
 =\bigcap_{n\ge1}\ \bigcup_{x\in D_\alpha}
 \{\lambda_{\max}(G_\rho^{[d]}(K+V_\alpha(x),w))
                                  <\kappa^{-1}+1/n\}.
\]
Every set on the right is Borel in $(K,w)$. On $\mathcal B_\alpha$ the modified and original Gram matrices agree for every $x\in\mathcal X_\alpha$, so this is precisely the event
that some permitted entry value gives the stated upper bound on the Gram matrix. Attainment of the minimum gives an admissible value also when the inequality is an equality.
Measurability follows from the continuous minimum over the compact parameter set. For a finite collection of functions take the union over its Gram matrix dimensions.

\Needspace{6\baselineskip}

\begin{lemma}
\label{W:lem:resolvent-starting-vector-gram-bound}
On $\mathcal B_\alpha\cap\mathcal A_{\alpha,\rho}$, for every entry value in the specified support and every point of the corresponding segment,
\begin{align}
 \|R_{K+V}(z)w\|^2+\|R_{K+V}(\bar z)w\|^2&\le {m_N}N^{o(1)},\nonumber\\
 \max_i|(R_{K+V}(z)w)_i|+\max_i|(R_{K+V}(\bar z)w)_i|
 &\le N^{o(1)}.\label{W:eq:resolvent-vector-on-the-regularity-event}
\end{align}
If no permitted entry value gives the stated upper bound on the Gram matrix, every composed function with the cutoff in the eigenvalues of the Gram matrix is identically zero for
all permitted entry values and in a neighbourhood of their support.
\end{lemma}
\begin{proof}
Choose a matrix $B_0$ obtained by a permitted entry replacement for which $G_\rho(B_0,w)\le\kappa^{-1}I$. On the portions of the contour that coincide with the ellipse, the
relative polynomial approximation in Lemma~\ref{W:lem:moment-weights-gram-resolvent} applies with $\rho\asymp {m_N}$, spectral parameter at most a fixed fraction of $s$, and
$\rho s\to\infty$. It gives $\|R_{B_0}(z)w\|^2\le C\kappa^{-1}/s\le C\kappa^{-1}{m_N}$. We give the spectral measure estimate implied by the bound on the Gram matrix on the
modified contour segments near the spectral edges. Let $\mu_w$ be the spectral measure of $B_0$ at $w$. Take $0<c_0\le1/4$ and $0<t_0\le c_0^2/4$. For $0<t\le t_0$, choose
$n=\min\{\rho,\lfloor c_0t^{-1/2}\rfloor\}$. Since $c_0t^{-1/2}\ge2$, this integer satisfies
\[
 \tfrac12\min\{\rho,c_0t^{-1/2}\}\le n\le
                                  \min\{\rho,c_0t^{-1/2}\}.
\]
For $t=0$ use $n=\rho$. Set
\[
 L_n(x)=\sum_{j=0}^{n-1}(j+1)p_j(x),\qquad
 S_n=\sum_{j=0}^{n-1}(j+1)^2=\frac{n(n+1)(2n+1)}6,
 \qquad n^3/3\le S_n\le n^3.
\]
If $x\in[2-t,2]$, write $x=2\cos\theta$. Since $2-x=4\sin^2(\theta/2)$ and $\sin(\theta/2)\ge\theta/\pi$ for $0\le\theta\le\pi$, $\theta\le(\pi/2)\sqrt t$ and
$n\theta\le\pi c_0/2\le\pi/8$. For $0<\theta$ and $j<n$, the sine formula therefore gives
\[
 p_j(x)=\frac{\sin((j+1)\theta)}{\sin\theta}
       \ge\frac{2(j+1)\theta/\pi}{\theta}=\frac2\pi(j+1).
\]
At $\theta=0$ use the continuous value $p_j(2)=j+1$. For $x=2\cosh v\ge2$, the function $\sinh((j+1)v)-(j+1)\sinh v$ vanishes at zero and has nonnegative derivative. Hence
$p_j(x)\ge j+1$ there as well. Consequently $L_n(x)\ge cS_n$ on all of $[2-t,\infty)$. The upper bound on the Gram matrix, applied to the coefficient vector $(j+1)_{j<n}$,
therefore implies
\begin{equation}
 c^2 S_n^2\mu_w([2-t,\infty))
 \le\int L_n(x)^2\,\mathrm d\mu_w(x)
 \le\kappa^{-1}S_n.
 \label{W:eq:gram-edge-polynomial-mass}
\end{equation}
It follows that
\[
 \mu_w([2-t,\infty))\le C\kappa^{-1}(\rho^{-3}+t_+^{3/2}).
\]
For $t>t_0$ the bound follows after increasing the constant from $\mu_w(\mathbb R)=1$; for $t<0$ it follows from the case $t=0$. Replacing $x$ by $-x$ proves the corresponding
lower edge bound.

Write $z=E+\mathrm i\eta$, with $|E-2|\le C_0\eta$ on the upper edge cone. For $r\ge\eta$, the preceding estimate gives
\[
 \mu_w([E-r,E+r])\le C_{\kappa,C_0}(\rho^{-3}+r^{3/2})
 \quad (r\le1).
\]
For $0<\eta\le1$, put $J=\lceil\log_2(1/\eta)\rceil$. Decompose the spectral integral into $|x-E|\le\eta$, the shells $2^j\eta<|x-E|\le2^{j+1}\eta$, $0\le j<J$, and the remaining
set. The mass bound holds up to radius two by the unit total mass. Thus
\[
 \|R_{B_0}(z)w\|^2
 \le C\kappa^{-1}\left[
 1+\sum_{j\ge0}\left(\rho^{-3}\eta^{-2}2^{-2j}
                          +\eta^{-1/2}2^{-j/2}\right)\right]
 \le C\kappa^{-1}(\eta^{-1/2}+\rho^{-3}\eta^{-2})
 \le C\kappa^{-1}m_N,
\]
because $\eta\ge cm_N^{-2}$ and $\rho\asymp m_N$. For $\int|z-x|^{-1}\,\mathrm d\mu_w(x)$ the same shells instead contribute
$C_{\kappa,C_0}(\rho^{-3}\eta^{-1}2^{-j} +\eta^{1/2}2^{j/2})$. For the same $J$,
\[
 \sum_{j=0}^{J-1}2^{-j}\le2,\qquad
 \eta^{1/2}\sum_{j=0}^{J-1}2^{j/2}
 =\eta^{1/2}\frac{2^{J/2}-1}{\sqrt2-1}
 \le\frac{\sqrt2}{\sqrt2-1}.
\]
Adding the central interval and the final set bounds this integral by $C_{\kappa,C_0}(1+\rho^{-3}/\eta)$. Since $\rho\asymp m_N$ and $\eta\ge c m_N^{-2}$, the last quantity is at
most $C'_{\kappa,C_0}(1+m_N^{-1})$. On each remaining elliptic piece, the relative estimate of Lemma~\ref{W:lem:moment-weights-gram-resolvent} holds on every vertical extension
with a common constant. Its spectral integral identity gives
\[
 \begin{aligned}
 w^*|R_{B_0}(E+\ii\eta)|w
 &=\frac2\pi\int_0^\infty
       \|R_{B_0}(E+\ii\sqrt{\eta^2+v^2})w\|^2\,\dd v\\
 &\le C\kappa^{-1}\int_{-2}^2
          \frac{\mu_{\rm sc}(\dd x)}{\sqrt{(E-x)^2+\eta^2}}
 \le C\kappa^{-1}(1+\log N).
 \end{aligned}
\]
For the last inequality, the nonreal pieces outside the neighbourhoods of the real contour boundary points have $|\eta|\ge cd\ge N^{-1}$, and the semicircle density is bounded. At
the exterior tips, $E$ is outside $[-2,2]$; the same integral is at most $|g(E)|\le1$. Together with the preceding estimate on the contour segments near the spectral edges, this
bounds the quadratic form at $w$ on the whole contour.

The coordinate bound follows from the same integral identity and the coordinate resolvent estimates on vertical extensions. For $0<\eta\le1$ on a nonreal piece, Ward's identity
gives
\[
 \|R_{B_0}(E+\ii y)e_i\|^2
 =\frac{-\operatorname{Im}R_{B_0,ii}(E+\ii y)}y\le C/y,
 \qquad \eta\le y\le2.
\]
For $y\ge1$, the deterministic resolvent bound gives a squared norm at most $y^{-2}$. Splitting the integral at $v=1$ therefore yields
\[
 e_i^*|R_{B_0}(E+\ii\eta)|e_i
 \le C\int_0^1\frac{\dd v}{\sqrt{\eta^2+v^2}}
          +C\int_1^\infty v^{-2}\,\dd v
 =C\operatorname{arsinh}(1/\eta)+C.
\]
The nonreal pieces have $\eta\ge c\eta_{\rm e}$, so this is $O(1+\log N)$. Near the real contour boundary points, $|R_{B_0}(E+\ii\eta)|\le|R_{B_0}(E)|$ in spectral order, and the
real coordinate resolvent bound from Lemma~\ref{W:lem:resolvent-contour} gives $e_i^*|R_{B_0}(E)|e_i\le C$. Values with $|\eta|>1$ are bounded by $1/|\eta|$, and adjunction gives
the lower half-plane. Finally polar decomposition and Cauchy--Schwarz give
\[
 |(R_{B_0}(z)w)_i|^2
 \le(e_i^*|R_{B_0}(z)|e_i)(w^*|R_{B_0}(z)|w)
 \le C\kappa^{-1}(1+\log N)^2=N^{o(1)}.
\]

For another permitted entry value, write the resulting matrix as $B_0+\iota_IT\iota_I^*$, with $|I|\le2$ and $\|T\|\le Ce_N$, and set $R_0=R_{B_0}(z)$, $r_0=R_0w$. The two exact
formulas are
\[
 (R_{B_0+\iota_IT\iota_I^*}w)_I=(I-R_{0,II}T)^{-1}(r_0)_I,
 \quad
 R_{B_0+\iota_IT\iota_I^*}w-r_0
 =R_{0,\cdot I}T(I-R_{0,II}T)^{-1}(r_0)_I.
\]
Choose $N_0$ so that the matrix subtracted from the identity has norm at most $1/2$ for $N\ge N_0$. The inverse then has norm at most two. With $\|(r_0)_I\|\le N^{o(1)}$ and the
coordinate and column bounds of Lemma~\ref{W:lem:resolvent-contour}, these identities give
\begin{align*}
 \|R_{B_0+\iota_IT\iota_I^*}w-r_0\|_2
 &\le 2Ce_N\left(\sum_{i\in I}\|R_0e_i\|^2\right)^{1/2}
                          \|(r_0)_I\|
 \le e_Nm_N^{1/2}N^{o(1)},\\
 \max_j|(R_{B_0+\iota_IT\iota_I^*}w-r_0)_j|
 &\le 2Ce_N\max_j\|R_{0,jI}\|\|(r_0)_I\|
 \le e_NN^{o(1)}.
\end{align*}
Adding the bounds for $r_0$ proves \eqref{W:eq:resolvent-vector-on-the-regularity-event} for every permitted entry value. Apply the same identities at $\bar z$ for the adjoint
resolvent.

If no permitted entry value gives the stated upper bound on the Gram matrix, compactness strengthens the strict inequality at each $x\in\mathcal X_\alpha$ to a positive gap:
\[
 \delta(K,w):=\min_{x\in\mathcal X_\alpha}
       \lambda_{\max}(G_\rho(K+V_\alpha(x),w))-\kappa^{-1}>0.
\]
Uniform continuity on this compact ball extends the lower bound $\kappa^{-1}+\delta(K,w)/2$ to a neighbourhood of the support of the replacement entry distributions, where spectral
truncation is still inactive. At least one eigenvalue remains in the zero region of $\gamma_\kappa$ there. The determinant cutoff and the composed function therefore vanish on that
neighbourhood, as do all their entry derivatives. This vanishing is deterministic for each $K,w$, with the strictly positive radius $\delta(K,w)$ just constructed. The derivatives
of the moment map itself need not vanish.
\end{proof}

We next express the moment increments as contour integrals. On the regularity event for the matrix with the selected entry set to zero put $R=R_K(z)$, ${\mathcal R_I}=R_{II}(z)$,
$T=V_\alpha(x)_{II}$. For a scalar analytic function $f$ inside the contour define
\begin{align}
 L_{f,1}&=\frac1{2\pi\mathrm i}\int_{\mathcal C_N}
       f(z)R_{\cdot I}TR_{I\cdot}\,\mathrm dz,\nonumber\\
 L_{f,2}&=\frac1{2\pi\mathrm i}\int_{\mathcal C_N}
       f(z)R_{\cdot I}T{\mathcal R_I}T R_{I\cdot}\,\mathrm dz,\nonumber\\
 L_{f,3}&=\frac1{2\pi\mathrm i}\int_{\mathcal C_N}
       f(z)R_{\cdot I}T({\mathcal R_I}T)^2R_{I\cdot}\,\mathrm dz,\nonumber\\
 L_{f,\ge4}&=\frac1{2\pi\mathrm i}\int_{\mathcal C_N}
       f(z)R_{\cdot I}T({\mathcal R_I}T)^3(I-{\mathcal R_I}T)^{-1}R_{I\cdot}\,\mathrm dz.
 \label{W:eq:resolvent-increment-formulas}
\end{align}
These are finite-dimensional matrix integrals; the first three are homogeneous in the real coordinate vector of their indicated degrees. The coefficient $L_{f,q}$ is
$D_{H}^q f(K)[V^q]/q!$ for $q\le3$. The mixed coefficients are obtained by polarization, or by summing all permutations of the labelled perturbation matrices.

For the homogeneous terms, differentiating $(zI-K-tV)R_{K+tV}(z)=I$ gives $\partial_tR_{K+tV}=R_{K+tV}VR_{K+tV}$. Induction at $t=0$ yields
\[
 \left.\partial_t^qR_{K+tV}(z)\right|_{t=0}
       =q!\,R_K(z)(V R_K(z))^q.
\]
All resolvents on the fixed contour stay bounded for $t$ in a neighborhood of zero. Hence differentiating the finite contour integral is justified by this uniform bound and gives
$L_{f,q}=D_H^q f(K)[V^q]/q!$. For $f_n=p_n$ and $f_{a,n}=p_n\theta_E({m_N}^2(2-\sigma\,\cdot))$, apply to each $L$ the following linear maps:
\begin{center}
\begin{tabular}{lll}
\toprule
statistic & scalar function & linear map on $L$\\
\midrule
$u_n$ & $f_n$ & $N^{-1}\operatorname{Tr}L$\\
$x_n$ & $f_n$ & ${m_N}\operatorname{Tr}(Q_wL)$\\
$v_{a,n}$ & $f_{a,n}$ & $N^{-1}\operatorname{Tr}L$\\
$y_{a,n}$ & $f_{a,n}$ & ${m_N}\operatorname{Tr}(Q_wL)$\\
$t_a$ & $f_{a,0}$ & $\operatorname{Tr}L$\\
\bottomrule
\end{tabular}
\end{center}
This defines vectors $\ell_\alpha,q_\alpha,c_\alpha,r_\alpha$ in precisely the coordinates of the moment vector of Section~\ref{W:sec:the-moment-vector-and-its-smooth-extensions}.

We spell out how \eqref{W:eq:resolvent-increment-formulas} produces all mixed entry derivatives. For an off-diagonal pair $a<b$ put
\[
 S_1=E_{ab}+E_{ba}\quad(\beta=1),\qquad
 S_1=\frac{E_{ab}+E_{ba}}{\sqrt2},\quad
 S_2=\frac{\mathrm i(E_{ab}-E_{ba})}{\sqrt2}\quad(\beta=2).
\]
The standardized block determines the perturbation $V(x)=N^{-1/2}\sum_rx_rS_r$. In the complex case $x_1=\sqrt2\operatorname{Re}\xi_{ab}$ and $x_2=\sqrt2\operatorname{Im}\xi_{ab}$;
a diagonal block instead has $S_1=E_{aa}$. If $\mathcal L$ is any of the five linear maps in the preceding table, its homogeneous coefficient tensor of degree $q$ is
\begin{equation}
 B_{r_1\cdots r_q}^{\mathcal L,f}
 =\frac{N^{-q/2}}{q!}\sum_{\sigma\in\mathfrak S_q}
 \mathcal L\left(\frac1{2\pi\mathrm i}
 \int_{\mathcal C_N}f(z)RS_{r_{\sigma(1)}}R\cdots
                    S_{r_{\sigma(q)}}R\,\mathrm dz\right).
 \label{W:eq:resolvent-entry-tensor}
\end{equation}
Its contribution is $\sum_{r_1,\ldots,r_q}B_{r_1\cdots r_q}^{\mathcal L,f} x_{r_1}\cdots x_{r_q}$. Formula \eqref{W:eq:resolvent-entry-tensor} gives the coefficients at orders
$q=1,2,3$ and the homogeneous degree four term. The remainder in \eqref{W:eq:resolvent-increment-formulas} includes every term of order at least four. The sum over ordered
coordinate labels supplies the multiplicity of each coordinate monomial when the tensors are contracted with the entry moments. In particular, a symmetric cubic tensor in two
coordinates contracts as
\[
 B_{111}\Delta M_{111}+3B_{112}\Delta M_{112}
       +3B_{122}\Delta M_{122}+B_{222}\Delta M_{222}.
\]
Here $\Delta M_{112}$ and $\Delta M_{122}$ are differences of the joint moments $\mathbb E X_1^2X_2$ and $\mathbb E X_1X_2^2$ of one complex entry, respectively.

Lemma~\ref{W:lem:resolvent-contractions} estimates a matrix element by its two exterior vectors and the interior resolvent compressions. For a normalized entry matrix $T$, its
trace estimate gives
\begin{equation}
 \left|\mathcal L_{\rm tr}(L_{f,q})\right|
 \le C_qN^{-q/2}\|x\|_1^q Q_Nm_N
        \int_{\mathcal C_N}|f(z)|\frac{|\mathrm dz|}{2\pi},
 \qquad \mathcal L_{\rm tr}=\Tr.
 \label{W:eq:resolvent-unnormalized-word}
\end{equation}
The same bound holds for the exact tail with $q=4$. Dividing by $N$ only for the $u$ and $v$ coordinates yields the following bounds; $m=m_N$, $m^3=N$, and $m^2d=N^{o(1)}$:
\begin{equation}
 \begin{array}{c|c|c|c}
 \text{coordinate}&\text{contour mass}&\text{trace factor}
                         &\text{result after }N^{-q/2}\|x\|_1^q\\ \hline
 u_n&N^{o(1)}&m/N&m^{-2}N^{o(1)}\\
 v_{a,n}&mdN^{o(1)}&m/N&N^{-1}N^{o(1)}\\
 t_a&Cd&m& m^{-1}N^{o(1)}
 \end{array}
 \label{W:eq:resolvent-normalizations}
\end{equation}
For the centered coordinates $x_n,y_{a,n}$ we use the exact difference \eqref{W:eq:resolvent-centered-contraction} and their additional factor $m$. The normalized trace part is
covered by $\delta_N^2=Q_N/m^2$. The unnormalized coordinate $t_a$ uses the bound for $\int_{\mathcal C_N}|\Theta_a(z)|\,|\mathrm dz|$ in \eqref{W:eq:resolvent-trace-bounds}.

\Needspace{6\baselineskip}

On $\mathcal B_\alpha\cap\mathcal A_{\alpha,\rho}$, denote the exact expansion of the moment vector by
\begin{equation}
 \mathbf z_N(K+V_\alpha(x),w)-\mathbf z_N(K,w)
 =\ell_\alpha(x)+q_\alpha(x)+c_\alpha(x)+r_\alpha(x).
 \label{W:eq:resolvent-moment-increment-identity}
\end{equation}
The coefficients are the contour integrals in \eqref{W:eq:resolvent-increment-formulas}, with the linear maps specified above.

\begin{theorem}
\label{W:thm:resolvent-increments}
There exist nonnegative weights $U_\alpha,V_\alpha$, measurable in the remaining entries, for which the following bounds hold for every admissible $x$, uniformly in the moment
indices.

\begin{align}
 \|\ell_\alpha(x)\|_\infty
 &\le C {m_N}N^{-1/2}\|x\|_1U_\alpha,\nonumber\\
 \|q_\alpha(x)\|_\infty
 &\le C {m_N}N^{-1}\|x\|_1^2V_\alpha,\nonumber\\
 \|c_\alpha(x)\|_\infty
 &\le C {m_N}N^{-3/2}\|x\|_1^3V_\alpha,\nonumber\\
 \|r_\alpha(x)\|_\infty
 &\le C {m_N}N^{-2}\|x\|_1^4V_\alpha,\label{W:eq:resolvent-moment-increment-bounds}\\
 \|\mathbf z_N(K+V_\alpha(x)) -\mathbf z_N(K)-\ell_\alpha(x)\|_\infty
 &\le C {m_N}N^{-1}\|x\|_1^2V_\alpha.\nonumber
\end{align}
For a diagonal block take $U_\alpha=V_\alpha$.
\end{theorem}
\begin{proof}
On $\mathcal B_\alpha$, every matrix $K+V_\alpha(x)$ with $x$ in the enlarged ball has spectrum inside $[-2-d/8,2+d/8]$. The modification $h_d$ is therefore the identity in a
neighbourhood of each matrix. Apply Cauchy's formula directly to the analytic functions $f_n$ and $f_{a,n}$. The contour excludes every pole of the functions $\Theta_a$. The exact
restricted resolvent identity gives
\[
 R_{K+V}-R_K=R_{\cdot I}T(I-{\mathcal R_I}T)^{-1}R_{I\cdot}.
\]
Put $Q=\mathcal R_I$. The polynomial identity
\[
 (I-QT)^{-1}=I+QT+(QT)^2+(QT)^3(I-QT)^{-1}
\]
follows on multiplying by $I-QT$, since the remaining terms cancel. Multiplication on the left by $R_{\cdot I}T$ and on the right by $R_{I\cdot}$ gives exactly the three
homogeneous integrands and the full remainder in \eqref{W:eq:resolvent-increment-formulas}. Contour integration is a linear map on these finite matrices; applying in turn the five
trace or matrix element maps in the preceding table proves \eqref{W:eq:resolvent-moment-increment-identity}. Since ${\mathcal R_I}$ is bounded and $\|T\|\le Ce_N\to0$, the inverse
has norm at most two uniformly.

A spectral moment at the starting vector and a trace contract the boundary factors of the same resolvent product differently; Figure~\ref{W:fig:moment-resolvent-words} displays
both contractions. For example, put $Q=R_{II}$ and $S=(R^2)_{II}$. The cubic integrand has the exact contractions
\begin{equation}
 \begin{aligned}
 w^*R_{\cdot I}TQTQTR_{I\cdot}w
   &=\sum_{i_1,\ldots,i_6\in I}
       (w^*R)_{i_1}T_{i_1i_2}Q_{i_2i_3}T_{i_3i_4}
          Q_{i_4i_5}T_{i_5i_6}(Rw)_{i_6},\\
 \Tr(R_{\cdot I}TQTQTR_{I\cdot})
   &=\operatorname{tr}_{I}(TQTQTS).
 \end{aligned}
 \label{W:eq:moment-word-boundary-contractions}
\end{equation}
The second identity is cyclicity of the trace, with $R_{I\cdot}R_{\cdot I}=(R^2)_{II}$. The trace then receives the normalization of its moment coordinate. \begingroup
\setlength{\intextsep}{1.2\baselineskip}
\begin{figure}[H]
\centering
\begin{tikzpicture}[x=1cm,y=1cm,every node/.style={mg text}]
 \node[anchor=west] at (0,1) {(a) Matrix element of the cubic resolvent term: two exterior factors};
 \foreach \j in {1,...,6}{\node[mg site,label=below:{$i_\j$}] (v\j) at ({1.4+(\j-1)*1.95},0) {};}
 \node[anchor=east] (wl) at (.35,0) {$w^*$};
 \node[anchor=west] (wr) at (12.2,0) {$w$};
 \draw[mg kernel] (wl)--node[mg label,above] {$R$}(v1);
 \draw[mg kernel] (v6)--node[mg label,above] {$R$}(wr);
 \foreach \a/\b in {1/2,3/4,5/6}{\draw[mg insertion] (v\a)--node[mg label,above] {$T$}(v\b);}
 \foreach \a/\b in {2/3,4/5}{\draw[mg kernel] (v\a)--node[mg label,above] {$Q$}(v\b);}
 \node[anchor=west] at (0,-1.35) {(b) Trace product: the exterior indices are summed};
 \foreach \j in {1,...,6}{\node[mg site,label=above:{$i_\j$}] (u\j) at ({1.4+(\j-1)*1.95},-2.4) {};}
 \foreach \a/\b in {1/2,3/4,5/6}{\draw[mg insertion] (u\a)--node[mg label,above] {$T$}(u\b);}
 \foreach \a/\b in {2/3,4/5}{\draw[mg kernel] (u\a)--node[mg label,above] {$Q$}(u\b);}
 \draw[mg kernel] (u6) to[out=-90,in=-90,looseness=.45]
      node[mg label,below] {$S=(R^2)_{II}$} (u1);
\end{tikzpicture}
\caption{The two cubic contractions in \eqref{W:eq:moment-word-boundary-contractions}. The vertices are summation indices in $I$, thick segments carry $T$, and the interior thin segments carry $Q=R_{II}$. In (a), the two exterior factors are the restrictions of $w^*R$ and $Rw$ used in the matrix element bound. Summation over their indices in (b) gives $S=(R^2)_{II}$. A degree-$j$ product has $j$ factors $T$ and $j-1$ interior factors $Q$; the exact error after the cubic expansion contains four factors $T$ and the bounded inverse $(I-QT)^{-1}$.}
\label{W:fig:moment-resolvent-words}
\end{figure}
\endgroup

Let $Q_N=N^{o(1)}$ dominate the constants in the coordinate norm bounds and set $\delta_N=Q_N^{1/2}/{m_N}$. At $K=K_\alpha$ define
\[
 a_i(z)=|(R_K(z)w)_i|+|(R_K(\bar z)w)_i|+\delta_N.
\]
For an off-diagonal block $\alpha=\{a,b\}$, let
\begin{equation}
 U_\alpha=\int a_a(z)a_b(z)\,\nu_N(\mathrm dz),\qquad
 V_\alpha=\int(a_a(z)^2+a_b(z)^2)\,\nu_N(\mathrm dz).
 \label{W:eq:resolvent-weights}
\end{equation}
For a diagonal block use $U_\alpha=V_\alpha=\int a_a^2\,\mathrm d\nu_N$. The deterministic term and the measure can be enlarged by a deterministic subpolynomial factor when needed.
They are functions of $K_\alpha$ and the deterministic vector $w$.

Put $l=(w^*R_K)_I$ and $r=(R_Kw)_I$. For an off-diagonal perturbation, $T_{aa}=T_{bb}=0$, and hence
\[
 |lTr|=|T_{ab}l_ar_b+T_{ba}l_br_a|
 \le2\|T\|a_a a_b.
\]
For $j\ge2$, the degree-$j$ term satisfies
\[
 |lT(QT)^{j-1}r|
 \le\|l\|_2\|r\|_2\|T\|^j\|Q\|^{j-1}
 \le C_j\|T\|^j(a_a^2+a_b^2).
\]
The bound $\|(I-QT)^{-1}\|\le2$ applies here. The exact fourth-order remainder therefore satisfies
\[
 |lT(QT)^3(I-QT)^{-1}r|
 \le C\|T\|^4(a_a^2+a_b^2).
\]
The difference minus the linear term equals $lTQT(I-QT)^{-1}r$. Its modulus is bounded by
\[
 \|l\|_2\|r\|_2\|T\|^2\|Q\|\|(I-QT)^{-1}\|
 \le C\|T\|^2(a_a^2+a_b^2).
\]
For a diagonal block $I=\{a\}$, $l,r$ are scalars with $|l|,|r|\le a_a$; the linear product is bounded by $\|T\|a_a^2$ and the degree-$j$ product by $\|T\|^j\|Q\|^{j-1}a_a^2$. The
same inverse bound handles its remainder. Integrating against the common dominating measure and using $\|T\|\le CN^{-1/2}\|x\|_1$ gives the stated homogeneous degrees. For the
normalized trace, cyclicity leaves an entry of $R_K^2$; its modulus is at most $Q_N{m_N}$ by \eqref{W:eq:resolvent-column}. After division by $N$ this is $Q_N/{m_N}^2=\delta_N^2$,
already included through the deterministic term in each factor of the weights. More explicitly, if $L_{f,j}$ is one of these matrix increments, its centered coordinate is
\[
 m_N\operatorname{Tr}(Q_wL_{f,j})
       =m_N w^*L_{f,j}w-\frac{m_N}{N}\operatorname{Tr}L_{f,j}.
\]
The first term is controlled by the two exterior resolvent factors, and the second by $m_N\delta_N^2$ under the same contour integral. The coordinate definitions give
$\Delta_j u_n+\Delta_j x_n/m_N=w^*L_{p_n,j}w$ and the corresponding identity with $p_n$ replaced by $p_n\Theta_a$, both for the homogeneous terms and for the remaining difference.

For a trace, apply \eqref{W:eq:resolvent-unnormalized-word} and the normalization table \eqref{W:eq:resolvent-normalizations}. The masses for $p_n$, $p_n\Theta_a$, and $\Theta_a$
are respectively $N^{o(1)}$, $m_NdN^{o(1)}$, and $O(d)$. Multiplying the table entries by $\|T\|^j\le C_jN^{-j/2}\|x\|_1^j$ proves \eqref{W:eq:resolvent-trace-bounds} for
$1\le j\le3$ and for the fourth-order remainder, whose block inverse has norm at most two. These sharper estimates are bounded by the common right sides in
\eqref{W:eq:resolvent-moment-increment-bounds}: the deterministic term in the resolvent weight is at least $N^{o(1)}/{m_N}^2$, so ${m_N}$ times that deterministic term also bounds
the unnormalized smoothed spectral trace. The moment vector norm is the maximum of the absolute coordinates. Since each coordinate has the same bound, with constants independent of
its index, taking this maximum preserves the bound. The subsequent summation is indexed by the selected matrix entry $\alpha$.
\end{proof}

\begin{corollary}
\label{W:cor:moment-trace-increments}
The trace coordinates have the following bounds: for $j=1,2,3$, and for the remainder with $j=4$, their degree-$j$ contributions are bounded by $N^{-j/2}\|x\|_1^j$ times
\begin{equation}
 u_n:\ N^{-2/3+o(1)},\qquad
 v_{a,n}:\ N^{-1+o(1)},\qquad
 t_a:\ N^{-1/3+o(1)}.
 \label{W:eq:resolvent-trace-bounds}
\end{equation}
For the difference minus its linear term use $j=2$ in this last statement.

\end{corollary}
\begin{proof}
For $i,j\in I$, Cauchy--Schwarz and normality of the resolvent give
\[
 |(R^2)_{ij}|\le
   \Big(\sum_k|R_{ik}|^2\Big)^{1/2}
   \Big(\sum_k|R_{kj}|^2\Big)^{1/2}\le Q_Nm_N.
\]
Thus $\|(R^2)_{II}\|\le|I|Q_Nm_N$. Cyclicity gives $\operatorname{Tr}(R_{\cdot I}D R_{I\cdot}) =\operatorname{tr}_I(D(R^2)_{II})$, whose modulus is at most
$|I|\|D\|\|(R^2)_{II}\|$. Here $|I|\le2$ and the bounded interior factors give $\|D\|\le C_jN^{-j/2}\|x\|_1^j$; the same estimate at $j=4$ includes $(I-QT)^{-1}$. For the
degree-$j$ integrand $A_j$, including the remainder at $j=4$,
\[
 |\operatorname{Tr}A_j(z)|
 \le C_jN^{-j/2}\|x\|_1^j Q_Nm_N.
\]
Multiply by the Cauchy weight and integrate. The three scalar masses are $N^{o(1)}$, $m_NdN^{o(1)}$, and $O(d)$ for $p_n$, $p_n\Theta_a$, and $\Theta_a$, respectively. Since
$m_N^3=N$ and $m_N^2d=N^{o(1)}$, their normalized bounds are
\[
 N^{-1}m_NN^{o(1)}=N^{-2/3+o(1)},\quad
 N^{-1}m_N^2dN^{o(1)}=N^{-1+o(1)},\quad
 m_NdN^{o(1)}=N^{-1/3+o(1)}.
\]
For the full quadratic remainder, put $Q=R_{II}$ and $T=V_\alpha(x)_{II}$. The exact remainder is
\[
 A_{\ge2}=R_{\cdot I}TQT(I-QT)^{-1}R_{I\cdot},\qquad
 \operatorname{Tr}A_{\ge2}
 =\operatorname{tr}_{I}\bigl(TQT(I-QT)^{-1}(R^2)_{II}\bigr).
\]
Here $\|Q\|\le C$, $\|(I-QT)^{-1}\|\le2$, and $\|(R^2)_{II}\|\le Q_Nm_N$. Thus $|\operatorname{Tr}A_{\ge2}|\le C N^{-1}\|x\|_1^2Q_Nm_N$. Multiplication by each of the three scalar
masses and the stated trace normalization gives the corresponding quadratic bound.
\end{proof}

\subsubsection{Integrated weights and conditional remainder estimates}
\label{W:sec:integrated-weights-and-conditional-remainders}

The contour formulas give bounds for the change in a moment when one entry is replaced. To use them in the interpolation identity, their weights must be integrated against the
conditional entry distributions and summed over the selected entries. The following estimates carry out those two operations together.

The contour, regularity events, and integrated weights are those of Section~\ref{W:sec:contour-integrals-for-moment-increments}. The weights are set to zero outside
$\mathcal B_\alpha\cap\mathcal A_{\alpha,\rho}$. For each $\alpha$ the expectation is taken under the distribution of its own matrix with the selected entry set to zero. The
resulting numerical bounds are then summed over $\alpha$.

\begin{proposition}
\label{W:prop:resolvent-weight-moments}
Uniformly in the deterministic unit starting vector and the interpolation parameter,
\begin{equation}
\begin{array}{c|rrrrrr}
 &\sum\mathbb EU&\sum\mathbb EU^2&\sum\mathbb EV
 &\sum\mathbb E(UV)&\sum\mathbb EV^2&\sum\mathbb EU^3\\ \hline
\text{off-diagonal}
 &N^{4/3+o(1)}&N^{2/3+o(1)}&N^{4/3+o(1)}
 &N^{2/3+o(1)}&N^{1+o(1)}&N^{1/3+o(1)}
\end{array}
 \label{W:eq:resolvent-weight-table}
\end{equation}
where each entry is an upper bound. On the diagonal, $\sum\mathbb EV\le N^{1/3+o(1)}$ and $\sum\mathbb EV^2+\sum\mathbb EV^3\le N^{o(1)}$. The replacement entries are integrated
using their truncated moments.
\end{proposition}
\begin{proof}
By Lemma~\ref{W:lem:resolvent-starting-vector-gram-bound}, $K_\alpha$ satisfies the bounds on the regularity event of Proposition~\ref{W:prop:moment-weights-row-moments}, with
$1/s$ replaced by ${m_N}N^{o(1)}$. The conditioning argument also works with a fixed zero block. Write $r=Rw$ and $c^{(i)}=R^{(i)}w^{(i)}$. On
$\mathcal B_\alpha\cap\mathcal A_{\alpha,\rho}$,
\[
 c_j^{(i)}=r_j-R_{ji}r_i/R_{ii},\quad
 \|c^{(i)}\|_\infty\le N^{o(1)},\quad
 \|c^{(i)}\|^2\le {m_N}N^{o(1)}.
\]
These two inequalities, together with existence of the minor resolvent on the real contour pieces, define an event measurable with respect to the minor that contains the original
event. First bound $|R_{ii}|$, then replace the indicator by this larger event, and finally integrate the independent row in $r_i=R_{ii}(w_i+h_i^*c^{(i)})$. For even $p\ge4$ the
resulting bound is
\[
 \mathbb E\mathbf1_{\mathcal B_\alpha\cap\mathcal A_{\alpha,\rho}}
 |r_i|^p
 \le N^{o(1)}\bigl[|w_i|^p+({m_N}/N)^{p/2}+{m_N}\mu_p\bigr],
 \qquad \mu_p=\sup_{i<j}\mathbb E|H_{ij}|^p.
\]
The deleted entry has variance zero, as allowed by the bound $C/N$. Truncated means cause an additional row scalar of size at most $N^{-3/2+o(1)}\sqrt {m_N}$, which is
$o(m_N^{-1})$. The truncation threshold gives $\mu_p\le F_NN^{-2}e_N^{p-4}$. Consequently, with $b_i=|w_i|+{m_N}^{-1}$,
\[
 \mathbb E\mathbf1_{\mathcal B_\alpha\cap\mathcal A_{\alpha,\rho}}a_i^4 \le N^{o(1)}b_i^4, \qquad \mathbb E\mathbf1_{\mathcal B_\alpha\cap\mathcal A_{\alpha,\rho}}a_i^6 \le N^{o(1)}\bigl(|w_i|^6+N^{-2}+N^{-5/3}e_N^2\bigr).
\]
Second moments have the bound $N^{o(1)}b_i^2$. The same spectral identities and row conditioning apply at every contour point, including the real exterior pieces. We use the
uniform conditional moment bounds at each fixed point before integrating over the contour.

Write $S_k=\sum_i b_i^k$. The coordinate inequalities in the proof of \eqref{W:eq:moment-weights-weight-table} give $S_1\le Cm_N^2$, $S_2\le Cm_N$, and $S_k\le C_k$ for fixed
$k\ge3$. For the selected block, write $I_\alpha=\mathbf1_{\mathcal B_\alpha\cap\mathcal A_{\alpha,\rho}}$. H\"older's inequality at two contour points gives
\[
 \begin{split}
 \mathbb E[I_\alpha a_a(z)a_b(z)a_a(z')^2]
 &\le(\mathbb E I_\alpha a_a(z)^4)^{1/4}
      (\mathbb E I_\alpha a_b(z)^4)^{1/4}
      (\mathbb E I_\alpha a_a(z')^4)^{1/2}\\
 &\le N^{o(1)}b_a^3b_b.
 \end{split}
\]
The sixth moment bound, applied at each fixed contour point, also gives
\[
 \sum_i(\mathbb E I_\alpha a_i(z)^6)^{1/2}
 \le N^{o(1)}\left(\sum_i|w_i|^3+1+N^{1/6}e_N\right)
 \le N^{o(1)}(1+N^{1/6}e_N).
\]
Let $J_N=\nu_N(\mathcal C_N)=N^{o(1)}$ and $I=\mathbf1_{\mathcal B_\alpha\cap\mathcal A_{\alpha,\rho}}$ for the specified block $\alpha$. The mixed quantity has the exact two
variable integral
\[
 \mathbb E(IU_{ab}V_{ab})=
 \iint\mathbb E\bigl[Ia_a(z)a_b(z)
                (a_a(z')^2+a_b(z')^2)\bigr]\,\nu_N(\dd z)\nu_N(\dd z').
\]
H\"older with exponents $(4,4,2)$ bounds its two terms by $N^{o(1)}b_a^3b_b$ and $N^{o(1)}b_ab_b^3$, respectively. Hence
\[
 \sum_{a<b}\mathbb E(IU_{ab}V_{ab})
 \le 2J_N^2N^{o(1)}S_3S_1=N^{2/3+o(1)}.
\]
Choose a common deterministic subpolynomial factor and set
\[
 d_i=N^{o(1)}\bigl(|w_i|^3+N^{-1}+N^{-5/6}e_N\bigr),
 \qquad (\mathbb E I a_i(z)^6)^{1/2}\le d_i
\]
uniformly in $z$. Expanding $U_{ab}^3$ as a triple integral and applying H\"older with six exponents equal to six gives
\[
 \mathbb E(IU_{ab}^3)\le J_N^3d_ad_b,
 \qquad
 \sum_{a<b}\mathbb E(IU_{ab}^3)
 \le J_N^3\left(\sum_i d_i\right)^2=N^{1/3+o(1)}.
\]
On a diagonal block $V_{ii}=\int a_i^2\,\dd\nu_N$. H\"older with three exponents equal to three gives
\[
 \sum_i\mathbb E(IV_{ii}^3)\le J_N^3\sum_i d_i^2
 \le N^{o(1)}\left(\sum_i|w_i|^6+N^{-1}+N^{-2/3}e_N^2\right)
 =N^{o(1)}.
\]
Here $S_k=\sum_i b_i^k$. For the remaining columns, with all factors evaluated at the same matrix, Tonelli and H\"older give
\begin{align*}
 \mathbb E(IU_{ab})
 &=\int\mathbb E[Ia_a(z)a_b(z)]\,\nu_N(\mathrm dz)
 \le N^{o(1)}J_N b_ab_b,\\
 \mathbb E(IU_{ab}^2)
 &=\iint\mathbb E[Ia_a(z)a_b(z)a_a(z')a_b(z')]
                     \,\nu_N(\mathrm dz)\nu_N(\mathrm dz')
 \le N^{o(1)}J_N^2b_a^2b_b^2,\\
 \mathbb E(IV_{ab})
 &\le N^{o(1)}J_N(b_a^2+b_b^2),\\
 \mathbb E(IV_{ab}^2)
 &\le N^{o(1)}J_N^2(b_a^4+2b_a^2b_b^2+b_b^4).
\end{align*}
The first line uses exponents $(2,2)$; the second uses four exponents $4$; in the last line expand the four products of squared factors and use $(2,2)$ for each. The deterministic
sums are therefore bounded by $J_NN^{o(1)}S_1^2$, $J_N^2N^{o(1)}S_2^2$, $J_NN^{o(1)}NS_2$, and $2J_N^2N^{o(1)}(NS_4+S_2^2)$, where $S_k=\sum_i b_i^k$. Since $S_1\le C m_N^2$,
$S_2\le C m_N$, $S_3+S_4\le C$, and $N=m_N^3$, they give the powers in \eqref{W:eq:resolvent-weight-table} after multiplication by $J_N$ or $J_N^2$. The diagonal quantities use
just $J_NN^{o(1)}S_2$ and $J_N^2N^{o(1)}S_4$. This proves the table under the law of each matrix $K_\alpha$ before the deterministic sums are taken.
\end{proof}

\Needspace{6\baselineskip}

\begin{theorem}
\label{W:thm:conditional-moment-expansion}
Let $F$ be one of the functions of Section~\ref{W:sec:cutoffs-in-moment-variables}, including its cutoffs, with $|F|\le1$. It may also contain the bounded cutoff for the Schur
complement. For fixed $\alpha,t$, let $\mathscr F_\alpha$ be generated by the remaining entries, let $K=K_\alpha(t)$, and set
\[
 I_\alpha=\mathbf1_{\mathcal B_\alpha\cap\mathcal A_{\alpha,\rho}},
 \qquad \delta\mu_\alpha=\mu_\alpha-\nu_\alpha.
\]
Then
\[
 \begin{split}
 \mathbb E\int F(K+V_\alpha(x))\,\delta\mu_\alpha(\dd x)
 ={}&\mathbb E I_\alpha\int
       [F(K+V_\alpha(x))-F(K)]\,\delta\mu_\alpha(\dd x)
       +\varepsilon_\alpha(t),\\
 |\varepsilon_\alpha(t)|\le{}&2\mathbb P(\mathcal B_\alpha^c).
 \end{split}
\]
In the first integral, every coefficient $C_{q,\alpha}(K)$ obtained by the finite moment expansion is paired with the full tensor $M_q(\mu_\alpha)-M_q(\nu_\alpha)$, with the factor
$I_\alpha$ still present. For every fixed $A>0$, the regularity events can be chosen so that
\[
 \sum_\alpha\int_0^1|\varepsilon_\alpha(t)|\,\dd t=O(N^{-A}),
\]
uniformly in the deterministic unit starting vector. The contribution from $\mathcal B_\alpha\cap\mathcal A_{\alpha,\rho}^c$ is zero.

\end{theorem}
\begin{proof}
Let $\mathscr F_\alpha$ be generated by the remaining entries and put $\mathcal G_\alpha=\mathcal B_\alpha\cap\mathcal A_{\alpha,\rho}$. This event and each coefficient
$C_{q,\alpha}$ at the deleted matrix are $\mathscr F_\alpha$-measurable. For the joint real coordinate laws of the selected entry,
\[
 \mathbb E\!\left[\mathbf1_{\mathcal G_\alpha}
  \int\langle C_{q,\alpha},x^{\otimes q}\rangle
           (\mu_\alpha-\nu_\alpha)(\dd x)\right]
 =\mathbb E\!\left[\mathbf1_{\mathcal G_\alpha}
  \langle C_{q,\alpha},M_q(\mu_\alpha)-M_q(\nu_\alpha)\rangle\right].
\]
For coordinate dimension $d_\alpha\le2$, the left side is the finite sum
\[
 \sum_{i_1,\ldots,i_q=1}^{d_\alpha}
 \mathbb E\!\left[\mathbf1_{\mathcal G_\alpha}
                  (C_{q,\alpha})_{i_1\cdots i_q}\right]
 \int x_{i_1}\cdots x_{i_q}(\mu_\alpha-\nu_\alpha)(\dd x).
\]
The real coordinates of a single complex entry appear jointly in this integral. Its moments are finite under the truncated supports. At each fixed $N$, the supports of all
remaining truncated entries form a bounded set. The coefficient functions are continuous on its compact closure because the extensions are smooth. They are therefore bounded there.
This finite bound and the finite entry moments justify the displayed interchanges of finite sums and expectations. Equality of the full coordinate tensors makes this coefficient
sum zero. Also $\int(\mu_\alpha-\nu_\alpha)=0$, which permits subtraction of the base value $F(K)$ on $\mathcal G_\alpha$. Theorem~\ref{W:thm:resolvent-increments} supplies the
exact increments there, with the integrated weight bounds of Proposition~\ref{W:prop:resolvent-weight-moments} and the trace normalizations in \eqref{W:eq:resolvent-trace-bounds}.

Lemma~\ref{W:lem:resolvent-starting-vector-gram-bound} gives
\[
 \mathbf1_{\mathcal B_\alpha\cap\mathcal A_{\alpha,\rho}^c}
 F(K+V_\alpha(x))=0
 \quad\text{for every allowed }x,
\]
so the conditional signed integral on this event is zero. The identity
\[
 1=\mathbf1_{\mathcal G_\alpha}
   +\mathbf1_{\mathcal B_\alpha\cap\mathcal A_{\alpha,\rho}^c}
   +\mathbf1_{\mathcal B_\alpha^c}
\]
leaves only the last event to estimate. At a fixed interpolation parameter write $K=K_\alpha(t)$. Since $|F|\le1$ and the two entry laws are probability measures, Tonelli's theorem
gives
\[
 \begin{aligned}
 &\int_0^1\sum_\alpha\mathbb E\!\left[
  \mathbf1_{\mathcal B_\alpha^c}
  \left|\int F(K+V_\alpha(x))
                    (\mu_\alpha-\nu_\alpha)(\dd x)\right|\right]\,\dd t\\
 &\qquad\le2\int_0^1\sum_\alpha
          \mathbb P(\mathcal B_\alpha(K_\alpha(t))^c)\,\dd t
 \le C_DN^{2-D}.
 \end{aligned}
\]
Taking $D>A+3$ gives $O(N^{-A})$. The supremum of the individual failure probabilities over $t$ is used inside this integral.

For the original matrix and its coupled truncation,
\[
 |\mathbb EF(H_N)-\mathbb EF(H_N^{\mathrm{tr}})|
 \le2\mathbb P(H_N\ne H_N^{\mathrm{tr}})=o(1).
\]
Thus a summed truncated comparison error $\epsilon_N$ becomes at most $\epsilon_N+2\mathbb P(H_N\ne H_N^{\mathrm{tr}})$ for the original law.
\end{proof}

\subsubsection{Taylor bounds for the cutoff and the localization statistic}
\label{W:sec:taylor-estimates-for-gram-and-localization-statistics}

The integrated weights can now be combined with derivatives of the observables. We begin with the cutoff defined as the product of a scalar function over the Gram eigenvalues and
the bounded localization statistic from Section~\ref{W:sec:extensions-of-localization-and-modulus-statistics}. Both are functions of the full vector of spectral statistics,
including the empirical trace coordinates and the coordinate cutoff. Their moment increments are controlled by Section~\ref{W:sec:integrated-weights-and-conditional-remainders}.

For the cutoff in the eigenvalues of the Gram matrix, all four contributions to the remainder of the expansion are summable along the moment segments. The localization statistic
needs a distinction between small and large increments. Local derivative bounds control the first case; boundedness of the statistic controls the complete remainder in the second.
In each case, the polynomial part of the expansion is integrated against the full distribution of the replaced entry.

The integrated potential carries an additional factor $m_N$, so its observables require the sharper derivative estimates. The calculation in
Section~\ref{W:sec:derivative-growth-at-the-critical-scale} already exhibits this growth at $G=I$. The coefficient cutoff $\chi_N$ and the estimates of
Sections~\ref{W:sec:derivatives-of-the-integrated-potential-and-coefficient-localization} and~\ref{W:sec:modulus-of-continuity-tests-and-derivative-kernels} control the
contractions with moment increments. For the localization statistic, the remainder proved here will be combined with its separate signed cubic estimate in
Section~\ref{W:sec:localization-of-the-low-energy-spectral-subspace}.

The vector $\mathbf z_N(H,w)$ and its reconstruction are those of Section~\ref{W:sec:the-moment-vector-and-its-smooth-extensions}. For fixed $r\le\rho$ and $a\in\mathcal I$, the
following notation makes their roles in the derivatives explicit:

\[
 G(\mathbf z)=\mathsf G_\rho(s),\quad
 G_r(\mathbf z)=\mathsf G_r(s),\quad
 B(\mathbf z)=\mathsf B_{r,a}(q),\quad t(\mathbf z)=t_a.
\]
Use $\dot G[v]$, $\dot B[v]$, and $\dot t[v]$ for their linear parts. Principal compression does not increase the trace norm, so $\|\dot G_r[v]\|_1\le\|\dot G[v]\|_1$. Throughout
this subsection, $\|\cdot\|_1$ denotes the matrix trace norm.

Let $\beta_j$ be the coordinate radii in \eqref{W:eq:moment-coordinates-box} and use $\mathfrak b$ from Section~\ref{W:sec:cutoffs-in-moment-variables}. Define
\begin{equation}
 \mathfrak n_\Gamma(v)=\|\dot G[v]\|_1+\mathfrak b(v),\qquad
 \mathfrak n_C(v)=\|\dot G[v]\|_1+\|\dot B[v]\|_1
                           +|\dot t[v]|+\mathfrak b(v).
 \label{W:eq:observable-norms}
\end{equation}
These are norms because all coordinate radii are positive. The full trace $t_a$ is unnormalized.

Write $\Phi_\Gamma$ and $\Phi_C$ for the cutoff and the localization statistic in \eqref{W:eq:moment-coordinates-testsG} and \eqref{W:eq:moment-coordinates-testsC}, respectively.
All cutoffs use the definitions in Section~\ref{W:sec:cutoffs-in-moment-variables}; in particular both functions have absolute value at most one.

The cutoff estimates of Lemma~\ref{W:lem:gram-cutoff-trace-derivatives} give the following bound in the norm $\mathfrak n_\Gamma$.

\begin{lemma}
\label{W:lem:observable-Gamma}
For every fixed $q\ge1$ there is $C_{q,\kappa}$, independent of the number of Lanczos steps and of $N$, such that
\begin{equation}
 |D_{\mathbf z}^q\Phi_\Gamma(\mathbf z)[v_1,\ldots,v_q]|
 \le C_{q,\kappa}\prod_{j=1}^q\mathfrak n_\Gamma(v_j)
 \qquad(\mathbf z\in\mathbb R^{D_N}).
 \label{W:eq:observable-Gammaglobal}
\end{equation}
In particular, the bound holds at every point on the moment segments in Taylor's formula.
\end{lemma}
\begin{proof}
The reconstruction satisfies $G(\mathbf z)=G(0)+\dot G[\mathbf z]$, so $DG[v]=\dot G[v]$ and $D^jG=0$ for $j\ge2$. For a set $S$ of $q$ derivative labels, the product and
composition rules consequently give
\[
 D_S\Phi_\Gamma
 =\sum_{U\sqcup V=S}(D_U\Xi_N)
             D_G^{|V|}\Gamma(G)[\dot G[v_j]:j\in V].
\]
For $U=\varnothing$ or $V=\varnothing$ the corresponding zeroth derivative is the function itself and has modulus at most one. Choose the constant as the maximum of the finitely
many derivative constants of orders $0,\ldots,q$. Every label belongs to precisely one of $U,V$, so there are $2^q$ summands, with no factor depending on the number of moment
coordinates. By Lemma~\ref{W:lem:gram-cutoff-trace-derivatives}, the absolute value of each term is at most
\[
 C_{q,\kappa}\prod_{j\in U}\mathfrak b(v_j)
                 \prod_{j\in V}\|\dot G[v_j]\|_1.
\]
Summing over the subsets $U$ bounds the derivative by $C_{q,\kappa}\prod_{j\in S}(\mathfrak b(v_j)+\|\dot G[v_j]\|_1)$, which is \eqref{W:eq:observable-Gammaglobal}.
\end{proof}

\Needspace{7\baselineskip} Use the contour $\mathcal C_N$ from Section~\ref{W:sec:contour-integrals-for-moment-increments}. Put
\[
 p^{(r)}(z)=(p_0(z),\ldots,p_{r-1}(z))^{\mathsf T},\qquad
 K_*(z)=\sum_{j=0}^{\rho_*-1}|p_j(z)|^2.
\]
For $\Theta_a(z)=\theta_E({m_N}^2(2-\sigma z))$ define
\[
 J_*(z)=\left(1+\sum_a|\Theta_a(z)|\right)
           \left(1+\sum_{n=1}^M |p_n(z)|/B_n\right).
\]
Enlarge the deterministic measure $\nu_N$ of Section~\ref{W:sec:contour-integrals-for-moment-increments} to
\begin{equation}
 \widehat\nu_N(\mathrm dz)=\nu_N(\mathrm dz)
 +C\left[
   {m_N}^{-1}K_*(z)\left(1+\sum_a|\Theta_a(z)|\right)
   +J_*(z)+{m_N}^2\sum_a|\Theta_a(z)|
 \right]\frac{|\mathrm dz|}{2\pi}.
 \label{W:eq:observable-measure}
\end{equation}
The constant may depend on the fixed cardinality of $\mathcal I$. This measure does not depend on the matrix or the starting vector.

\Needspace{6\baselineskip} The deterministic measure $\widehat\nu_N$ has total mass $N^{o(1)}$.

\begin{lemma}
\label{W:lem:observable-measure}
One has $\widehat\nu_N(\mathcal C_N)=N^{o(1)}$.
\end{lemma}
\begin{proof}
Use $s=4\sqrt d$, ${\eta_{\mathrm e}}={m_N}^{-2}$, and $z=\zeta+\zeta^{-1}$ with $1\le|\zeta|\le e^s$ on the contour. The finite Chebyshev formula gives
\[
 |p_j(z)|\le 2e^{(j+1)s}|z^2-4|^{-1/2},\qquad
 K_*(z)\le4\rho_*e^{2\rho_*s}|z^2-4|^{-1}.
\]
On an elliptic arc, $|\mathrm dz|=|\zeta-\zeta^{-1}|\,\mathrm d\theta$ and $|z^2-4|=|\zeta-\zeta^{-1}|^2$. Since $|\zeta-\zeta^{-1}|=2(\sinh^2s+\sin^2\theta)^{1/2}$, the integral
on the elliptic arcs is bounded by
\[
 \frac12\int_0^{2\pi}\frac{\dd\theta}{(\sinh^2s+\sin^2\theta)^{1/2}}
 \le C\int_0^{\pi/2}\frac{\dd\theta}{(s^2+4\theta^2/\pi^2)^{1/2}}
 \le C\operatorname{arsinh}(C/s)\le C(1+|\log s|).
\]
Here the symmetry of $\sin^2\theta$ and $\sin\theta\ge2\theta/\pi$ on $[0,\pi/2]$ justify the middle inequality. On each inward edge modification, the parameter
$u=|\operatorname{Re}z \mp2|$ obeys $|z^2-4|\ge c({\eta_{\mathrm e}}+u)$ and $|\mathrm dz|\le C\,\mathrm du$; its contribution is at most
\[
 C\int_0^{Cd}(\eta_{\rm e}+u)^{-1}\,\dd u
 =C\log(1+Cd/\eta_{\rm e}).
\]
The integral with the square-root denominator equals $2C(\sqrt{\eta_{\rm e}+Cd}-\sqrt{\eta_{\rm e}})$. Since $d/\eta_{\rm e}=N^{o(1)}$ and $d\to0$, these are bounded by
$C(1+\log N)$ and a constant, respectively. Consequently
\[
 \int_{\mathcal C_N}\frac{|\mathrm dz|}{|z^2-4|}
 \le C(1+\log N),\qquad
 \int_{\mathcal C_N}\frac{|\mathrm dz|}{|z^2-4|^{1/2}}\le C.
\]
Since $\rho_*/{m_N}=O(1)$ and $\rho_*s=O((\log N)^{\gamma/2})=o(\log N)$, the $K_*/{m_N}$ term has mass $N^{o(1)}$.

Choose
\[
 Q_N=e^{2\rho_*s}\left(1+\sum_{a\in\mathcal I}
                     \sup_{z\in\mathcal C_N}|\Theta_a(z)|\right).
\]
Then $Q_N=N^{o(1)}$. In particular it bounds the products, rather than only the individual factors, that occur below. Also $B_n=(n+1)e^{n\sqrt{2d}}\ge n+1$ and $n+1\le2\rho_*$. The
two preceding integrals therefore imply
\[
 \begin{aligned}
 \int_{\mathcal C_N}\frac{K_*(z)}{m_N}
          \left(1+\sum_a|\Theta_a(z)|\right)\frac{|\dd z|}{2\pi}
 &\le C Q_N\frac{\rho_*}{m_N}(1+\log N),\\
 \int_{\mathcal C_N}J_*(z)\frac{|\dd z|}{2\pi}
 &\le C Q_N\left(1+\sum_{n=1}^M\frac1{n+1}\right)\\
 &\le C Q_N(1+\log(M+1)).
 \end{aligned}
\]
Both bounds are subpolynomial since $\rho_*\asymp m_N$ and $M\le2\rho_*$. Finally Section~\ref{W:sec:contour-integrals-for-moment-increments} proves
$\int|\Theta_a(z)|\,|\mathrm dz|\le C_a d$; hence their contribution with coefficient ${m_N}^2$ is $O({m_N}^2d)=N^{o(1)}$.
\end{proof}

\Needspace{6\baselineskip}

On the intersection of the event where the stated spectral and resolvent bounds hold and the event that a permitted entry value gives the bound on the Gram matrix in
Section~\ref{W:sec:contour-integrals-for-moment-increments}, use the same exact decomposition $\Delta z=\ell+q+c+r$ and $b=\Delta z-\ell$.

\begin{proposition}
\label{W:prop:observable-increments}
There are nonnegative weights $U_\alpha,V_\alpha$, measurable with respect to the remaining entries, for which, with either norm in \eqref{W:eq:observable-norms},
\begin{align}
 \mathfrak n(\ell_\alpha(x))
   &\le C {m_N}N^{-1/2}\|x\|_1U_\alpha,\nonumber\\
 \mathfrak n(q_\alpha(x))+\mathfrak n(b_\alpha(x))
   &\le C {m_N}N^{-1}\|x\|_1^2V_\alpha,\nonumber\\
 \mathfrak n(c_\alpha(x))
   &\le C {m_N}N^{-3/2}\|x\|_1^3V_\alpha,\nonumber\\
 \mathfrak n(r_\alpha(x))
   &\le C {m_N}N^{-2}\|x\|_1^4V_\alpha.
 \label{W:eq:observable-increment-bounds}
\end{align}
For off-diagonal blocks their sums satisfy
\begin{equation}
\begin{array}{c|rrrrrr}
 &\sum\mathbb EU&\sum\mathbb EU^2&\sum\mathbb EV
 &\sum\mathbb E(UV)&\sum\mathbb EV^2&\sum\mathbb EU^3\\ \hline
 &N^{4/3+o(1)}&N^{2/3+o(1)}&N^{4/3+o(1)}
 &N^{2/3+o(1)}&N^{1+o(1)}&N^{1/3+o(1)}.
\end{array}
 \label{W:eq:observable-weight-sums}
\end{equation}
On the diagonal, $U_\alpha=V_\alpha$, $\sum\mathbb EV_\alpha\le N^{1/3+o(1)}$, and $\sum\mathbb EV_\alpha^2\le N^{o(1)}$. Every statement is uniform in the deterministic starting
vector and at each fixed interpolation parameter, which may subsequently be integrated.
\end{proposition}
\begin{proof}
Let $R=(zI-K_\alpha)^{-1}$, $\mathcal R_I=R_{II}$ and $T=V_\alpha(x)_{II}$. For $j=1,2,3$ write
\[
 A_j(z)=R_{\cdot I}T(\mathcal R_IT)^{j-1}R_{I\cdot},\qquad
 A_{\ge4}(z)=R_{\cdot I}T(\mathcal R_IT)^3(I-\mathcal R_IT)^{-1}R_{I\cdot}.
\]
These are exactly the finite expansion terms of Section~\ref{W:sec:contour-integrals-for-moment-increments}. The increment beyond the linear term uses
$T\mathcal R_IT(I-\mathcal R_IT)^{-1}$.

The identity $s_n=u_n+x_n/m_N$ gives $\Delta s_n=\Delta u_n+\Delta x_n/m_N$ at each pair of matrices and after every derivative of positive order. Thus the centered and empirical
trace contributions cancel within each coefficient. Its degree-$j$ Gram matrix increment is
\begin{equation}
 E_j=\frac1{2\pi\mathrm i}\int_{\mathcal C_N}
       p^{(\rho)}(z)p^{(\rho)}(z)^{\mathsf T}
       w^*A_j(z)w\,\mathrm dz.
 \label{W:eq:observable-rankone}
\end{equation}
For the mixed $r\times r$ matrix the formula is
\[
 \Delta_jB=\frac1{2\pi\mathrm i}\int_{\mathcal C_N}
 p^{(r)}(z)p^{(r)}(z)^{\mathsf T}\Theta_a(z)
                              w^*A_j(z)w\,\dd z.
\]
Here $r$ denotes the matrix dimension. Writing $p\,p^{\mathsf T}=p\,(\overline p)^*$ shows that its sole nonzero singular value is $\|p\|_2\|\overline p\|_2=\|p\|_2^2$. Formula
\eqref{W:eq:observable-rankone} also holds for the remainder and the increment beyond the linear term. The fixed moment $s_0=1$ causes no extra term because the increment of the
constant polynomial is zero.

Let $Q_N=N^{o(1)}$ dominate the coordinate constants in the column norm bounds of Section~\ref{W:sec:moment-increments-under-entry-replacement}, put $\delta_N=\sqrt{Q_N}/{m_N}$,
and set
\[
 a_i(z)=|(R(z)w)_i|+|(R(\bar z)w)_i|+\delta_N.
\]
For an off-diagonal block, $|w^*A_1w|\le C\|T\|a_aa_b$, while for $j\ge2$ it is at most $C\|T\|^j(a_a^2+a_b^2)$. The inverse in the remainder is uniformly bounded. Define
\[
 U_\alpha=\int a_aa_b\,\mathrm d\widehat\nu_N,\qquad
 V_\alpha=\int(a_a^2+a_b^2)\,\mathrm d\widehat\nu_N,
\]
with the definition using $a_a(z)^2$ for diagonal blocks. Since $\rho\le\rho_*$, we have $\|p^{(\rho)}(z)\|_2^2\le K_*(z)$, where $K_*(z)=\sum_{j=0}^{\rho_*-1}|p_j(z)|^2$. For the
homogeneous integrand $u_j(z)=w^*A_j(z)w$, the triangle inequality for the trace norm in \eqref{W:eq:observable-rankone} gives
\[
 \|E_j\|_1\le\int_{\mathcal C_N}K_*(z)|u_j(z)|\,
                                  \frac{|\dd z|}{2\pi}.
\]
The mixed matrix increment has trace norm at most
\[
 \int_{\mathcal C_N}\|p^{(r)}(z)\|_2^2
      |\Theta_a(z)|\,|u_j(z)|\,\frac{|\dd z|}{2\pi},
\]
by its Cauchy formula and the trace norm of the rank-one polynomial matrix. By \eqref{W:eq:observable-measure}, each of these measures is bounded by $Cm_N\widehat\nu_N$. Thus the
linear term is at most $Cm_N\|T\|U_\alpha$, and terms of degree $j\ge2$ are at most $C_jm_N\|T\|^jV_\alpha$. The exact remainder has the same bound with $j=4$, and the full
difference minus its linear term with $j=2$.

For the coordinate norm put $\mathrm d\ell(z)=|\mathrm dz|/(2\pi)$, $L_1(z)=a_a(z)a_b(z)$ and $L_j(z)=a_a(z)^2+a_b(z)^2$ for $j\ge2$. On a diagonal block use $L_j=a_a^2$. In all
cases $\delta_N^2\le L_j$. Let $\Delta_j$ denote the homogeneous term of degree $j$. The matrix element and normalized trace contractions give
\[
 \begin{aligned}
 \frac{|\Delta_j u_n|}{B_n}
 &\le C_j\|T\|^j\int\frac{|p_n(z)|}{B_n}\delta_N^2\,\mathrm d\ell(z),\\
 \frac{|\Delta_j x_n|}{2m_NB_n}
 &\le C_j\|T\|^j\int\frac{|p_n(z)|}{B_n}
                                  (L_j(z)+\delta_N^2)\,\mathrm d\ell(z),\\
 \frac{|\Delta_j v_{a,n}|}{B_n}
 &\le C_j\|T\|^j\int\frac{|p_n(z)\Theta_a(z)|}{B_n}
                                  \delta_N^2\,\mathrm d\ell(z),\\
 \frac{|\Delta_j y_{a,n}|}{2m_NB_n}
 &\le C_j\|T\|^j\int\frac{|p_n(z)\Theta_a(z)|}{B_n}
                             (L_j(z)+\delta_N^2)\,\mathrm d\ell(z).
 \end{aligned}
\]
The trace estimate used here is $N^{-1}|\operatorname{Tr}A_j(z)|\le C_j\|T\|^jQ_Nm_N/N =C_j\|T\|^j\delta_N^2$. The sum over the polynomial index $n$ and $a\in\mathcal I$ is bounded
by
\[
 C_j\|T\|^j\int J_*(z)L_j(z)\,\mathrm d\ell(z).
\]
The coordinate $y_{a,0}$ is covered by the constant term in $J_*$. For $t_a$ and its box radius $N$, respectively,
\[
 |\Delta_jt_a|\le C_j Q_Nm_N\|T\|^j \int|\Theta_a|\,\mathrm d\ell, \qquad |\Delta_jt_a|/N\le C_j\delta_N^2\|T\|^j \int|\Theta_a|\,\mathrm d\ell.
\]
The first is dominated by $C_jm_N\|T\|^j\int m_N^2|\Theta_a|L_j\,\mathrm d\ell$, since $m_N^3\delta_N^2=Q_Nm_N$. The measure in \eqref{W:eq:observable-measure} therefore controls
all coordinates of $\mathfrak n_\Gamma$ and $\mathfrak n_C$. The bounded inverse in $A_{\ge4}$ gives the identical estimates with $j=4$ for the remainder, and with $j=2$ for the
full nonlinear increment. Substituting $\|T\|\le C N^{-1/2}\|x\|_1$ proves \eqref{W:eq:observable-increment-bounds}.

Here $I=\mathbf1_{\mathcal B_\alpha\cap\mathcal A_{\alpha,\rho}}$ for the selected entry $\alpha$. Proposition~\ref{W:prop:resolvent-weight-moments} applies with $\widehat\nu_N$ in
place of $\nu_N$: the pointwise row estimates hold on the same contour, and $J=\widehat\nu_N(\mathcal C_N)=N^{o(1)}$. In its proof, an integral with $k$ contour variables is
bounded by $J^k$ times the same uniform row estimate. Thus every column of \eqref{W:eq:resolvent-weight-table} has the stated order for this measure as well. In particular, the
sixth-moment calculation gives
\[
 \mathbb E(IU_\alpha^3)\le J^3d_ad_b,\qquad
 \sum_i d_i\le N^{o(1)}(1+N^{1/6}e_N),
\]
with $d_i=N^{o(1)}(|w_i|^3+N^{-1}+N^{-5/6}e_N)$. The resulting sum is $J^3N^{o(1)}(1+N^{1/6}e_N)^2\le N^{1/3+o(1)}$. The remaining five sums and the two diagonal sums use the same
deterministic $S_k$ as in that proposition. All pointwise estimates are uniform over $K_\alpha$, which permits their deterministic upper bounds to be summed over $\alpha$.
\end{proof}

\Needspace{7\baselineskip} For a smooth function $\Phi$, an exact Taylor expansion of the moment vector $\Delta \mathbf z=\ell+q+c+r$, and $b=\Delta \mathbf z-\ell$, define
\begin{equation}
 P_3(\mathbf z;x)=D_{\mathbf z}\Phi(\mathbf z)\ell+D_{\mathbf z}\Phi(\mathbf z)q
            +\tfrac12D_{\mathbf z}^2\Phi(\mathbf z)[\ell,\ell]+D_{\mathbf z}\Phi(\mathbf z)c.
 \label{W:eq:observable-P3}
\end{equation}

The polynomial $P_3$ consists of the linear and quadratic moment terms and the first derivative applied to the cubic moment term. It is not the full third-order Taylor polynomial
of the composed entry function. To see the difference, set $f(t)=\Phi(\mathbf z+t\ell+t^2q+t^3c)$. Differentiating at zero gives
\[
 \frac{f^{(3)}(0)}6
 =D\Phi(\mathbf z)c+D^2\Phi(\mathbf z)[\ell,q]
                       +\frac16D^3\Phi(\mathbf z)[\ell,\ell,\ell].
\]
The remainder estimates below bound the last two terms absolutely; moment matching does not remove them. The cubic contribution in $P_3$ is $D\Phi(\mathbf z)c$. Indeed, set
$b=\Delta\mathbf z-\ell$. The three integral identities are

\begin{align*}
 \Phi(\mathbf z+\ell)-\Phi(\mathbf z)
 &=D\Phi(\mathbf z)\ell+\tfrac12D^2\Phi(\mathbf z)[\ell,\ell]\\
 &\quad+\tfrac12\int_0^1(1-u)^2D^3\Phi(\mathbf z+u\ell)[\ell^3]\,\dd u,\\
 D\Phi(\mathbf z+\ell)b-D\Phi(\mathbf z)b
 &=\int_0^1D^2\Phi(\mathbf z+u\ell)[\ell,b]\,\dd u,\\
 \Phi(\mathbf z+\ell+b)-\Phi(\mathbf z+\ell)
 &=D\Phi(\mathbf z+\ell)b\\
 &\quad+\int_0^1(1-u)D^2\Phi(\mathbf z+\ell+ub)[b,b]\,\dd u.
\end{align*}
Adding them and using $b=q+c+r$ yields
\begin{align}
 \Phi(\mathbf z+\Delta \mathbf z)-\Phi(\mathbf z)-P_3(\mathbf z;x)
 ={}&D_{\mathbf z}\Phi(\mathbf z)r\nonumber\\
 &+\tfrac12\int_0^1(1-u)^2D_{\mathbf z}^3\Phi(\mathbf z+u\ell)
                         [\ell,\ell,\ell] \,\mathrm du\nonumber\\
 &+\int_0^1D_{\mathbf z}^2\Phi(\mathbf z+u\ell)[\ell,b] \,\mathrm du\nonumber\\
 &+\int_0^1(1-u)D_{\mathbf z}^2\Phi(\mathbf z+\ell+ub)[b,b] \,\mathrm du.
 \label{W:eq:observable-exactR}
\end{align}
The coordinate derivative estimates apply on all the displayed segments.

\Needspace{6\baselineskip} The global derivative bounds for $\Phi_\Gamma$ and the increment bounds give summable absolute remainders of the expansions after integration over the
law of the selected entry.

\begin{theorem}
\label{W:thm:observable-Gram}
For $\Phi=\Phi_\Gamma$, on the regularity events of Section~\ref{W:sec:contour-integrals-for-moment-increments}, the sums of the absolute expectations of the four terms in
\eqref{W:eq:observable-exactR}, in their displayed order, are bounded by
\begin{equation}
 N^{-1/3+o(1)},\qquad N^{-1/6+o(1)},\qquad
 N^{-1/6+o(1)},\qquad N^{-1/3+o(1)}
 \label{W:eq:observable-fourrates}
\end{equation}
over off-diagonal blocks, including integration of the distributions of the entries being replaced and of the interpolation parameter. Only a bounded third absolute moment and an
$N^{o(1)}$ truncated fourth moment of the standardized replacement entry are needed. On diagonal blocks, the remainder after the linear polynomial has summed absolute expectation
$N^{-1/3+o(1)}$ using only the second moment of the entry distribution.
\end{theorem}
\begin{proof}
Lemma~\ref{W:lem:observable-Gamma} bounds the derivatives at $\mathbf z+u\ell$ and $\mathbf z+\ell+ub$, $0\le u\le1$, the points on the two segments in
\eqref{W:eq:observable-exactR}. Write
\[
 A_\alpha=C{m_N}N^{-1/2}\|x\|_1U_\alpha,\quad
 Q_\alpha=C{m_N}N^{-1}\|x\|_1^2V_\alpha,\quad
 R_\alpha=C{m_N}N^{-2}\|x\|_1^4V_\alpha.
\]
If $C_j$ is the constant in the $j$th derivative estimate, the four terms of \eqref{W:eq:observable-exactR} have the respective bounds
\[
 C_1R_\alpha,\qquad
 \frac{C_3}{2}A_\alpha^3\int_0^1(1-u)^2\,\dd u
       =\frac{C_3}{6}A_\alpha^3,\qquad
 C_2A_\alpha Q_\alpha,\qquad \frac{C_2}{2}Q_\alpha^2.
\]
For the third term the integration interval has length one; for the fourth its weight has integral $1/2$. The global estimates apply to every point of the two segments
$\mathbf z+u\ell$ and $\mathbf z+\ell+ub$, $0\le u\le1$, in the real coordinate domain. Let $I_\alpha$ be the regularity indicator and let
$M_j=\sup_\alpha\int\|x\|_1^j(\mu_\alpha+\nu_\alpha)(\dd x)$. For example, conditional integration gives exactly
\[
 \mathbb E\left[I_\alpha\int A_\alpha Q_\alpha
                         (\mu_\alpha+\nu_\alpha)(\dd x)\right]
 \le C m_N^2N^{-3/2}M_3\,\mathbb E[I_\alpha U_\alpha V_\alpha].
\]
Both weights and the indicator are fixed by the remaining entries. The other three terms use respectively $M_4$, $M_3$, and $M_4$ in the same conditional integral. Here $M_3=O(1)$
and $M_4=N^{o(1)}$. Now use \eqref{W:eq:observable-weight-sums}. The powers are explicitly
\[
\begin{array}{c|c|c}
\text{term}&\text{coefficient and resolvent weight sum}&\text{power}\\ \hline
R&{m_N}N^{-2}M_4N^{4/3+o(1)}&-1/3\\
A^3&{m_N}^3N^{-3/2}M_3N^{1/3+o(1)}&-1/6\\
AQ&{m_N}^2N^{-3/2}M_3N^{2/3+o(1)}&-1/6\\
Q^2&{m_N}^2N^{-2}M_4N^{1+o(1)}&-1/3.
\end{array}
\]
The parameter integral has length one and all estimates are uniform. For a diagonal block the exact difference has $\mathfrak n_\Gamma(\Delta z)\le Cm_NN^{-1/2}|x|V_\alpha$ and
$\mathfrak n_\Gamma(\Delta z-\ell)\le Cm_NN^{-1}|x|^2V_\alpha$. The first two global derivatives of $\Phi_\Gamma$ imply
\begin{align*}
 &|\Phi_\Gamma(z+\Delta z)-\Phi_\Gamma(z)-D\Phi_\Gamma(z)[\ell]|\\
 &\qquad\le C\mathfrak n_\Gamma(\Delta z-\ell)
                 +C\mathfrak n_\Gamma(\Delta z)^2\\
 &\qquad\le CN^{-1}|x|^2(m_NV_\alpha+m_N^2V_\alpha^2).
\end{align*}
Conditional integration uses $\mathbb E|x|^2\le C$. The diagonal sums of resolvent weights therefore give
\[
 \sum_{\alpha\ {\rm diag}}\mathbb E|\mathcal R_\alpha|
 \le N^{-1+o(1)}(m_NN^{1/3}+m_N^2)=N^{-1/3+o(1)}.
\]
\end{proof}

\Needspace{7\baselineskip} \subparagraph{The bounded localization statistic}
\label{W:sec:taylor-bounds-for-localization-statistic}

We estimate derivatives in a fixed neighborhood of a moment vector arising from a matrix. The order relation $0\le B\le G_r$ at the base vector controls the weighted matrix
$B=F_r^*\Theta(H)F_r$ and its coordinate variations throughout that neighborhood.

\begin{lemma}
\label{W:lem:local-derivative-bounds-for-localization}
Let $\mathbf z_0=\mathbf z_N(H,w)$ be any moment vector associated with a matrix. For every $\mathbf z'$ with $\mathfrak n_C(\mathbf z'-\mathbf z_0)\le1$, and $q=1,2,3$,
\begin{equation}
 |D_{\mathbf z}^q\Phi_C(\mathbf z')[v_1,\ldots,v_q]|
 \le C_{q,\kappa}\prod_j\mathfrak n_C(v_j).
 \label{W:eq:local-derivative-bounds-for-localization}
\end{equation}
The constant is independent of the number of Lanczos steps, $t_a(\mathbf z_0)$, and $N$. The estimate holds on the stated line segment in the coordinate variables, including points
which are not moments of any matrix.
\end{lemma}
\begin{proof}
At the moment vector $\mathbf z_0$ of $H$,
\begin{equation}
 0\le B(\mathbf z_0)\le G_r(\mathbf z_0),
 \label{W:eq:positivity-of-the-mixed-moment-matrix}
\end{equation}
because $0\le\Theta_a(H)\le I$. If any derivative of $\Gamma(G)$ at $G(\mathbf z')$ is nonzero, then every eigenvalue of $G(\mathbf z')$ is in $[\kappa,\kappa^{-1}]$. Indeed, an
eigenvalue outside that interval stays outside the compact support of $\gamma_\kappa$ in a neighbourhood, where the determinant is identically zero. The limiting cases are also
covered by the support gap of $\gamma_\kappa$.

It suffices, in every nonzero product rule term, to work with $\|G_r(\mathbf z')\|\le\kappa^{-1}$. The distance between the moment vectors then gives
$\|G_r(\mathbf z_0)\|\le\kappa^{-1}+1$. By \eqref{W:eq:positivity-of-the-mixed-moment-matrix},
\[
 \|B(\mathbf z')\|\le\|B(\mathbf z_0)\|+1\le\kappa^{-1}+2.
\]
Also $a_\kappa(G_r(\mathbf z'))=G_r(\mathbf z')^{-1}$ on a neighbourhood, so its derivatives are the derivatives of the inverse. With $A=G_r(\mathbf z')^{-1}$ and
$G_j=\dot G_r[v_j]$,
\[
 D_{G}^qa_\kappa(G_r)[G_1,\ldots,G_q]
 =(-1)^q\sum_{\sigma\in S_q}
 AG_{\sigma(1)}A\cdots G_{\sigma(q)}A.
\]
In particular the trace norm is at most $q!\kappa^{-q-1}\prod_j\|G_j\|_1$. The maps $G_r,B,t$ have constant first derivatives and zero derivatives of order at least two. The
labelled product rule therefore gives
\begin{align*}
 D_{\mathbf z}^qC[v_1,\ldots,v_q]={}&\mathbf1_{\{q=1\}}\dot t[v_1]
 -\operatorname{Tr}\bigl[D_{G}^qa_\kappa[G_1,\ldots,G_q]B\bigr]\\
 &-\sum_{j=1}^q\operatorname{Tr}
  \bigl[D_{G}^{q-1}a_\kappa[G_1,\ldots,\widehat G_j,\ldots,G_q]
                      \dot B[v_j]\bigr].
\end{align*}
Put $g_j=\|G_j\|_1$, $b_j=\|\dot B[v_j]\|_1$, and $\tau_j=|\dot t[v_j]|$. In the first trace term, place one $G_j$ in trace norm and the other factors in operator norm. Since
$\|G_j\|\le g_j$, this gives
\[
 \left|\operatorname{Tr}(D_G^qa_\kappa[G_1,\ldots,G_q]B)\right|
 \le q!\kappa^{-q-1}(\kappa^{-1}+2)\prod_jg_j.
\]
For $q=1$, the second trace term is $|\operatorname{Tr}(A\dot B[v_1])|\le\kappa^{-1}b_1$. For $q\ge2$, put $\dot B[v_j]$ in trace norm and all inverse and $G_k$ factors in operator
norm, obtaining
\[
 \left|\operatorname{Tr}(D_G^{q-1}a_\kappa[G_k:k\ne j]\dot B[v_j])\right|
 \le(q-1)!\kappa^{-q}b_j\prod_{k\ne j}g_k.
\]
Adding these inequalities and $\mathbf1_{\{q=1\}}|\dot t[v_1]|$ proves
\[
 |D_{\mathbf z}^q C[v_1,\ldots,v_q]|
 \le C_{q,\kappa}\prod_j(g_j+b_j+\tau_j)
 \le C_{q,\kappa}\prod_j\mathfrak n_C(v_j).
\]
For the scalar composition, its exact labelled derivative is
\[
 D_S(\widetilde\psi\circ C)
 =\sum_{\pi\in\Pi(S)}\widetilde\psi^{(|\pi|)}(C)
                              \prod_{B\in\pi}D_BC.
\]
Every direction occurs exactly once in each partition product, including repeated direction vectors. Each product has bound $C_{q,\kappa}\prod_{j\in S}\mathfrak n_C(v_j)$, by the
globally bounded scalar derivatives. Finally
\[
 D_S\Phi_C=
 \sum_{A\sqcup B\sqcup C'=S}(D_A\Xi_N)(D_B\Gamma)
                    D_{C'}(\widetilde\psi\circ C).
\]
The box and determinant estimates give, respectively,
\[
 |D_A\Xi_N|\le C_q\prod_{j\in A}\mathfrak b(v_j),\qquad
 |D_B\Gamma|\le C_{q,\kappa}\prod_{j\in B}\|\dot G[v_j]\|_1.
\]
The $3^q$ assignments and the finitely many set partitions depend only on $q$. Since every component norm is at most $\mathfrak n_C$, summing proves
\eqref{W:eq:local-derivative-bounds-for-localization}.
\end{proof}

\Needspace{6\baselineskip}

\begin{lemma}
\label{W:lem:observable-large}
Suppose $|\Phi|\le1$ and, at a base point associated with a matrix $\mathbf z_0$, its first three derivatives with respect to the spectral statistics obey a uniform multilinear
bound in a norm $\mathfrak n$ throughout $\mathfrak n(\mathbf z-\mathbf z_0)\le1$. Let $\Delta \mathbf z=\ell+q+c+r$ and $b=\Delta \mathbf z-\ell$. Suppose
\[
 \mathfrak n(\ell)\le A,\quad
 \mathfrak n(q),\mathfrak n(c),\mathfrak n(b)\le Q,\quad
 \mathfrak n(r)\le R.
\]
Then the exact remainder after the polynomial \eqref{W:eq:observable-P3} satisfies
\begin{equation}
 |\Phi(\mathbf z_0+\Delta \mathbf z)-\Phi(\mathbf z_0)-P_3(\mathbf z_0;x)|
 \le C\{R+A^3+AQ+Q^2\}.
 \label{W:eq:observable-largeremainder}
\end{equation}
On the complementary region we use the direct value bound.
\end{lemma}
\begin{proof}
Write $z=\mathbf z_0$ and $\Delta=\Delta\mathbf z=\ell+b$. The polynomial in \eqref{W:eq:observable-P3} is
\[
 P_3=D\Phi(z)[\ell+q+c]+\tfrac12D^2\Phi(z)[\ell,\ell].
\]
Suppose first that $A+Q\le1$. The two paths $z+t\ell$ and $z+\ell+tb$ stay in the unit ball about $z$. Apply the integral Taylor formula along these two paths. Since $b=q+c+r$,
subtraction of $P_3$ yields the exact remainder
\[
 \begin{aligned}
 \Phi(z+\Delta)-\Phi(z)-P_3
 ={}&D\Phi(z)[r]
  +\tfrac12\int_0^1(1-t)^2D^3\Phi(z+t\ell)[\ell,\ell,\ell]\,\dd t\\
 &+\int_0^1D^2\Phi(z+t\ell)[\ell,b]\,\dd t
  +\int_0^1(1-t)D^2\Phi(z+\ell+tb)[b,b]\,\dd t.
 \end{aligned}
\]
If $C_j$ bounds the $j$th derivative in $\mathfrak n$, the absolute value is at most $C_1R+(C_3/6)A^3+C_2AQ+(C_2/2)Q^2$.

For $A+Q>1$, the value bound and the derivatives at $z$ give
\[
 |\Phi(z+\Delta)-\Phi(z)-P_3|
 \le2+C_1(A+2Q)+(C_2/2)A^2.
\]
If $A\ge1/2$, then $1\le8A^3$, $A\le4A^3$, $A^2\le2A^3$, and $Q\le2AQ$. If $A<1/2$, then $Q>1/2$, so $1\le4Q^2$, $Q\le2Q^2$, $A\le2AQ$, and $A^2\le AQ$. In both cases the last
display is bounded by $C(A^3+AQ+Q^2)$. Together with the small increment calculation this proves \eqref{W:eq:observable-largeremainder} throughout the support of the replacement
entry distributions. The split $A+Q\le1$ is used only for estimating this same remainder. At a fixed deleted matrix, if $S=\{x:A(x)+Q(x)\le1\}$ and $\lambda$ is either entry law,
the identity remains
\[
 \int\bigl(\Phi(z+\Delta(x))-\Phi(z)\bigr)\,\lambda(\dd x)
 =\int P_3(z;x)\,\lambda(\dd x)
   +\int_S\mathcal R(x)\,\lambda(\dd x)
   +\int_{S^c}\mathcal R(x)\,\lambda(\dd x).
\]
The polynomial is integrated against the full law. Restricting its integral to $S$ would replace the prescribed moment tensors by restricted moments and destroy their exact
cancellation. On $S^c$ the displayed estimate controls the whole remainder; its individual integral terms need not obey the local derivative bounds there.
\end{proof}

\Needspace{6\baselineskip}

\begin{theorem}
\label{W:thm:localization-taylor-remainder}
For the same truncated matrices and regularity events, the exact remainder from the two successive Taylor expansions for the spectral localization function $\Phi_C$ satisfies
\begin{equation}
 \sum_{\alpha\ {\rm off}}\mathbb E|\mathcal R_{\alpha,C}|
 \le N^{-1/6+o(1)},
 \label{W:eq:localization-taylor-remainder-sum}
\end{equation}
including integration over the full entry law and integration in the interpolation parameter. On the region $A_\alpha+Q_\alpha\le1$, the four absolute integral terms of
\eqref{W:eq:observable-exactR} obey the individual rates \eqref{W:eq:observable-fourrates}. On the complementary event, apply the value bound of Lemma~\ref{W:lem:observable-large}.

For both the cutoff and the spectral localization statistic, the diagonal remainder after the linear polynomial has summed absolute expectation $N^{-1/3+o(1)}$, using only the
uniformly bounded diagonal second moment. If the stated spectral and resolvent estimates fail with probability at most $C_DN^{-D}$, the contribution of these exceptional events to
the summed absolute remainder is at most $C_D'N^{7-D}$. If the stated spectral and resolvent bounds hold and every permitted entry value fails the upper bound on the Gram matrix,
the polynomial coefficients and the exact remainder are zero.
\end{theorem}
\begin{proof}

Write $m=m_N$, $t_x=N^{-1/2}\|x\|_1$ and use the deterministic bounds
\[
 A_\alpha=Cm t_xU_\alpha,\qquad
 Q_\alpha=Cm t_x^2V_\alpha,\qquad
 R_\alpha=Cm t_x^4V_\alpha.
\]
The cubic increment has norm at most $C m t_x^3V_\alpha\le Q_\alpha$ for all sufficiently large $N$, because $t_x\le C e_N\to0$.
Lemmas~\ref{W:lem:local-derivative-bounds-for-localization} and \ref{W:lem:observable-large} give $|\mathcal R_{\alpha,C}|\le C(R_\alpha+A_\alpha^3+A_\alpha Q_\alpha+Q_\alpha^2)$.
Let $M_j$ bound the $j$th absolute moment of the standardized replacement entry under either law. Conditional integration over that entry and the summed weights in
\eqref{W:eq:observable-weight-sums} give
\[
 \begin{aligned}
 \sum_\alpha\mathbb ER_\alpha
 &\le C m N^{-2}M_4 N^{4/3+o(1)}=N^{-1/3+o(1)},\\
 \sum_\alpha\mathbb EA_\alpha^3
 &\le C m^3N^{-3/2}M_3 N^{1/3+o(1)}=N^{-1/6+o(1)},\\
 \sum_\alpha\mathbb E(A_\alpha Q_\alpha)
 &\le C m^2N^{-3/2}M_3 N^{2/3+o(1)}=N^{-1/6+o(1)},\\
 \sum_\alpha\mathbb EQ_\alpha^2
 &\le C m^2N^{-2}M_4 N^{1+o(1)}=N^{-1/3+o(1)}.
 \end{aligned}
\]
Here $M_3=O(1)$ and $M_4=N^{o(1)}$; the indicators of the conditioning events are included in every summed expectation. The estimates hold uniformly at each interpolation
parameter, so their integrals over $[0,1]$ obey the same bounds. At this stage they prove \eqref{W:eq:localization-taylor-remainder-sum} with the indicator
$\mathbf1_{\mathcal B_\alpha\cap\mathcal A_{\alpha,\rho}}$. The other base matrices are accounted for below.

For diagonal blocks, the exact resolvent difference gives
\[
 \mathfrak n(\Delta \mathbf z)\le C{m_N}N^{-1/2}|x|V_\alpha,
 \qquad
 \mathfrak n(\Delta \mathbf z-\ell)\le C{m_N}N^{-1}|x|^2V_\alpha.
\]

When $d_*=\mathfrak n(\Delta\mathbf z)\le1$, every point of $\mathbf z+u\Delta\mathbf z$, $0\le u\le1$, is in the ball where the local derivative estimate holds. The exact identity
is
\[
 \begin{aligned}
 &\Phi(\mathbf z+\Delta\mathbf z)-\Phi(\mathbf z)-D\Phi(\mathbf z)[\ell]\\
 &=D\Phi(\mathbf z)[\Delta\mathbf z-\ell]
  +\int_0^1(1-u)D^2\Phi(\mathbf z+u\Delta\mathbf z)
                      [\Delta\mathbf z,\Delta\mathbf z]\,\dd u.
 \end{aligned}
\]
Its modulus is at most $C_1\mathfrak n(\Delta\mathbf z-\ell)+(C_2/2)d_*^2$. When $d_*>1$, boundedness of $\Phi$ and the bound for $D\Phi(\mathbf z)$ give
\[
 |\Phi(\mathbf z+\Delta\mathbf z)-\Phi(\mathbf z)-D\Phi(\mathbf z)\ell|
 \le2+C\mathfrak n(\ell)
 \le2+C d_*+C\mathfrak n(\Delta\mathbf z-\ell)
 \le C'd_*^2+C\mathfrak n(\Delta\mathbf z-\ell).
\]
Thus in both cases
\[
 |\Phi(\mathbf z+\Delta \mathbf z)-\Phi(\mathbf z)-D_{\mathbf z}\Phi(\mathbf z)\ell|
 \le C\{\mathfrak n(\Delta \mathbf z-\ell)+\mathfrak n(\Delta \mathbf z)^2\}.
\]
For the cutoff one may also use its global derivative estimate directly. The diagonal sums of resolvent weights therefore yield
\[
 C{m_N}N^{-1}M_{2,\mathrm d}N^{1/3+o(1)}
 +C{m_N}^2N^{-1}M_{2,\mathrm d}N^{o(1)}
 =N^{-1/3+o(1)}.
\]
Only $M_{2,\mathrm d}=O(1)$ was used. We now include the base matrices outside the regularity event. This requires a bound for the coefficients subtracted in the definition of the
remainder as well as for the bounded values of $\Phi$.

The required coefficient bounds follow from the globally defined moment functions. To see their dependence on $N$, put
\[
 f_n(x)=p_n(h_d(x)),\qquad
 f_{a,n}(x)=p_n(h_d(x))\theta_E(m_N^2(2-\sigma x)),
 \qquad a=(\sigma,E).
\]
For $0\le j\le3$ and $n\le2\rho_*-1$, Proposition~\ref{W:prop:gram-bounds} gives
$\|f_n^{(j)}\|_\infty\le m_N^{2j+1}N^{o(1)}$.
The scalar derivatives of $\theta_E$ are bounded, so the ordinary product rule gives
\[
 \begin{split}
 \|f_{a,n}^{(j)}\|_\infty
 &\le\sum_{k=0}^j\binom jk
      \|f_n^{(k)}\|_\infty
      m_N^{2(j-k)}\|\theta_E^{(j-k)}\|_\infty\\
 &\le m_N^{2j+1}N^{o(1)}.
 \end{split}
\]
For fixed matrices $W_1,\ldots,W_j$ with $\|W_l\|_{\rm HS}\le CN^{-1/2}$, Lemma~\ref{W:lem:smooth-functional-calculus} consequently yields, at every Hermitian $K$,
\[
 \|D_K^j f_n(K)[W_1,\ldots,W_j]\|_{\rm HS}
 +\|D_K^j f_{a,n}(K)[W_1,\ldots,W_j]\|_{\rm HS}
 \le N^{-j/2}m_N^{2j+1}N^{o(1)},\qquad 1\le j\le3.
\]
Each entry of $G_\rho$ or $B$ is a sum of at most $\rho_*$ such vector moments. For a matrix $M$ of size at most $\rho_*$,
$\|M\|_1\le\sqrt{\rho_*}\|M\|_{\rm HS}\le\rho_*^{3/2}\max_{u,v}|M_{uv}|$. Since $\rho_*\asymp m_N$, it follows that
\[
 \|D_K^jG_\rho\|_1+\|D_K^jB\|_1
 \le N^{-j/2}m_N^{2j+7/2}N^{o(1)}
 =N^{j/6+7/6+o(1)}.
\]
The trace coordinate satisfies $|D_K^jt_a|\le\sqrt N\,N^{-j/2}m_N^{2j}C_j$. For the coordinate norm $\mathfrak b$, there are $O(m_N)$ coordinates, their box radii are at least one,
and the largest additional normalization in a vector moment is $m_N$. Thus
\[
 |D_K^jt_a|\le N^{j/6+1/2+o(1)},\qquad
 \mathfrak b(D_K^j\mathbf z_N)\le N^{j/6+1+o(1)}.
\]
The fixed linear reconstructions of $G_\rho,B,t_a$ commute with these derivatives. Combining the preceding inequalities proves
\[
 \mathfrak n_C(D_K^j\mathbf z_N[W_1,\ldots,W_j])
 \le N^{j/6+7/6+o(1)},\qquad j=1,2,3.
\]
The same upper bound applies to $\mathfrak n_\Gamma$.

Write $W_i=N^{-1/2}\mathcal T_\alpha e_i$ and expand the homogeneous moment coefficients as
\[
 \ell(x)=\sum_i x_i\ell_i,\quad
 q(x)=\sum_{i,j}x_ix_jq_{ij},\quad
 c(x)=\sum_{i,j,k}x_ix_jx_kc_{ijk}.
\]
Here $\ell_i=D_K\mathbf z_N[W_i]$, $q_{ij}=D_K^2\mathbf z_N[W_i,W_j]/2$, and $c_{ijk}=D_K^3\mathbf z_N[W_i,W_j,W_k]/6$. Their norms are bounded respectively by $N^{4/3+o(1)}$,
$N^{3/2+o(1)}$, and $N^{5/3+o(1)}$. The derivatives of $\Phi$ at $\mathbf z_N(K,w)$ satisfy Lemma~\ref{W:lem:observable-Gamma} or
Lemma~\ref{W:lem:local-derivative-bounds-for-localization} for every $K$. Therefore the coefficients of the polynomial \eqref{W:eq:observable-P3} satisfy
\[
 \begin{aligned}
 |D\Phi[\ell_i]|&\le C N^{4/3+o(1)},\\
 |D\Phi[q_{ij}]+\tfrac12D^2\Phi[\ell_i,\ell_j]|
       &\le C(N^{3/2+o(1)}+N^{8/3+o(1)}),\\
 |D\Phi[c_{ijk}]|&\le C N^{5/3+o(1)}.
 \end{aligned}
\]
Choose $N_0$ so that all three right sides are at most $N^5$ for $N\ge N_0$. This deterministic index is independent of the probability exponent $D$. Summing the absolute
coordinate monomials uses $\sum_{i_1,\ldots,i_j}|x_{i_1}\cdots x_{i_j}|=\|x\|_1^j$ and hence gives, on the whole matrix space,
\[
 |\mathcal R_{\alpha,C}(K,x)|
 \le C N^5(1+\|x\|_1^3)
       \quad(\alpha\text{ off the diagonal}).
\]
For a diagonal entry only the linear coefficient is subtracted, so the corresponding bound is $CN^5(1+|x|)$.

Let $\mathcal B_\alpha(t)$ be the event where the stated spectral and resolvent bounds hold, which is measurable with respect to the remaining entries. Conditional on them, the
replacement entry has its prescribed deterministic distribution. The off-diagonal third absolute moments are bounded, and the diagonal first absolute moments are bounded by the
square roots of their second moments. Thus
\[
 \begin{split}
 &\mathbb E\!\left[\mathbf1_{\mathcal B_\alpha(t)^c}
    \int|\mathcal R_{\alpha,C}(K_\alpha(t),x)|
                   (\mu_\alpha+\nu_\alpha)(\dd x)\right]\\
 &\hspace{2em}\le CN^5\mathbb P(\mathcal B_\alpha(t)^c)
 \le C_DN^{5-D}.
 \end{split}
\]
The real dimension of an entry is at most two; the integral uses its joint distribution. There are $N(N-1)/2$ off-diagonal entries and $N$ diagonal entries. Integration over
$t\in[0,1]$ and summation therefore give an additional bound $C_DN^{7-D}$. Given $A>0$, choosing $D>A+8$ makes it $O(N^{-A})$.

On $\mathcal B_\alpha\cap\mathcal A_{\alpha,\rho}^c$, the continuous function $x\mapsto\lambda_{\max}(G_\rho(K+V_\alpha(x),w))$ is greater than $\kappa^{-1}$ on the compact
enlarged ball of permitted values. Its minimum is consequently strictly greater than $\kappa^{-1}$. The determinant cutoff vanishes at every such value and on a neighbourhood of
the base moment vector. Thus $\Phi(\mathbf z_N(K+V_\alpha(x),w))=\Phi(\mathbf z_N(K,w))=0$ and all derivatives of $\Phi$ occurring in $P_3$ vanish separately. Both $P_3$ and
$\mathcal R_{\alpha,C}$ are zero there.

The three disjoint regions $\mathcal B_\alpha\cap\mathcal A_{\alpha,\rho}$, $\mathcal B_\alpha\cap\mathcal A_{\alpha,\rho}^c$, and $\mathcal B_\alpha^c$ now account for every base
matrix. Their bounds give the stated sums for the exact remainders. The replacement of the original matrix by its coupled truncation is a separate comparison of expectations of
bounded functions, with error at most $2\mathbb P(H_N\ne H_N^{\rm tr})=o(1)$.
\end{proof}

\Needspace{7\baselineskip}
\subsubsection{Contour estimates and dependence on the polynomial degree}
\label{W:sec:contours-for-lanczos-comparison}
The moment estimates must still be paired with polynomials whose degree grows with $N$. Their contour bounds therefore keep the dimension of the Gram matrix explicit until the
entry sums have been evaluated. We use the spectral statistics of Section~\ref{W:sec:the-moment-vector-and-its-smooth-extensions}, with $\rho_*\asymp m_N=N^{1/3}$ and
$d=A_0\rho_*^{-2}(\log(N+2))^\gamma$, $1<\gamma<2$. The cutoff parameter is fixed in $0<\kappa<1/4$, the dimension of the Gram matrix satisfies $8\le\rho\le\rho_*$, and the indices
of the integrated potential are at most $\rho-2$.

At a fixed interpolation parameter, $K_\alpha$ denotes the matrix with the selected entry set to zero, and the perturbation $N^{-1/2}\mathcal T_\alpha x$ represents all real
coordinates of that entry. Events used to integrate $x$ are measurable with respect to the remaining entries. This measurability is maintained in the contour estimates below.

All upper-triangular entries are independent. The matrix supports obey $\|V_\alpha(x)\|\le C e_N$, where $e_N\to0$ and $e_N^{-1}=N^{o(1)}$. For each of the two distributions of the
entries being replaced, the standardized third absolute moment is bounded and the fourth is at most $F_N=N^{o(1)}$. The same bounds hold when each entry law is replaced by a convex
combination of these two distributions. The centered off-diagonal matrix entries therefore satisfy
\begin{equation}
 \mathbb E|H_{ij}-\mathbb EH_{ij}|^2\le C/N,\quad
 \mathbb E|H_{ij}-\mathbb EH_{ij}|^4\le N^{-2+o(1)},\quad
 \mathbb E|H_{ij}-\mathbb EH_{ij}|^6
       \le N^{-2+o(1)}e_N^2.
 \label{W:eq:resolvent-entrymom}
\end{equation}
The last inequality follows by multiplying the truncated fourth moment by the square of the support bound. The truncated means are $N^{-2+o(1)}$ in matrix off-diagonal coordinates.
The diagonal entries satisfy the second moment assumptions of Section~\ref{W:sec:moment-increments-under-entry-replacement}.

Set
\[
 s=4\max\{\rho^{-1},\sqrt d\},\qquad
 \mathcal C=\{\zeta+\zeta^{-1}:\zeta=e^{s+\mathrm i\theta},\
 0\le\theta\le2\pi\}.
\]
For $\rho\ge8$ and large $N$, $s<1$, $\rho s\ge4$, and $\rho s=O(1+(\log N)^{\gamma/2})$. Thus $e^{C\rho s}=N^{o(1)}$ for every fixed $C$. The contour encloses the spectra of all
matrices obtained by the permitted entry replacements on the regularity events of Section~\ref{W:sec:contour-integrals-for-moment-increments}. On this event, the spectrum of every
such matrix lies where $h_d$ is the identity. Cauchy's formula on the ellipse therefore uses the unmodified polynomials and their matrix derivatives. The coordinates in this
calculation are the spectral moments $s_n=w^*p_n(H)w$.

\Needspace{6\baselineskip} For functions of the polynomial moments, the contour can be adapted to the dimension of the Gram matrix while retaining uniform bounds for resolvent
entries and for the resolvent applied to the starting vector.

Let $\mathcal A_{\alpha,\rho}$ be the event measurable with respect to the remaining entries that at least one matrix $B_0$ obtained by replacing the selected entry with a value in
the enlarged ball satisfies $G_\rho(B_0,w)\le\kappa^{-1}I$.
\begin{lemma}
\label{W:lem:resolvent-weight-input}
For each fixed $D$, there is a regularity event measurable with respect to the remaining entries $\mathcal B_\alpha$ with $\sup_{\alpha,t}\mathbb P(\mathcal B_\alpha^c)\le N^{-D}$.
On this event, simultaneously for every block value in the enlarged ball,
\begin{equation}
 \max_{i,j}|R_{ij}(z)|\le C,\quad
 \min_i|R_{ii}(z)|\ge c,\quad
 \max_i\big(\|R(z)e_i\|^2+\|R(\bar z)e_i\|^2\big)
       \le\rho N^{o(1)},\qquad z\in\mathcal C.
 \label{W:eq:resolvent-coordinate}
\end{equation}
On $\mathcal B_\alpha\cap\mathcal A_{\alpha,\rho}$,
\begin{equation}
 \|R(z)w\|^2+\|R(\bar z)w\|^2\le\rho N^{o(1)},\qquad
 \max_i|(R(z)w)_i|\le N^{o(1)}
 \label{W:eq:resolvent-starting-vector}
\end{equation}
for every permitted entry value. On $\mathcal B_\alpha\cap\mathcal A_{\alpha,\rho}^c$ the exponential multiplied by its cutoffs is identically zero on a neighbourhood of all
allowed entry values, and every finite entrywise remainder of the expansion of that function is zero.
\end{lemma}
\begin{proof}
We track the dependence of the coordinate and counting estimates on $\rho$. On a nonreal elliptic piece away from its small real tips, $\eta\ge cs^2\ge c\rho^{-2}$ and
\[
 \upsilon(z)/\eta=e^{-s}/(2\sinh s)\le C\rho.
\]
Counting therefore gives $\bar\varrho/\eta\le C\rho+N^{o(1)}\rho^4/N\le\rho N^{o(1)}$, since $\rho^3\le C N$. Divide the imaginary part inequality \eqref{W:eq:resolvent-imaginary}
by $\eta$. Its square-root contribution is at most
\[
 N^{o(1)}(\rho^{5/2}/\sqrt N+\rho^4/N)\le\rho N^{o(1)}.
\]
The remaining terms have the respective bounds
\[
\begin{gathered}
 N^{o(1)}F_N\rho/N,\quad N^{-17/12+o(1)}\rho^4,\quad
 N^{-11/8+o(1)}\rho^4,\\
 N^{-13/12+o(1)}\rho^4,\quad N^{-7/4+o(1)}\rho^4.
\end{gathered}
\]  After division by $\rho$ and use of $\rho^3\le C N$, these
five quantities are bounded respectively by
\[
 F_NN^{-1+o(1)},\quad N^{-5/12+o(1)},\quad
 N^{-3/8+o(1)},\quad N^{-1/12+o(1)},\quad N^{-3/4+o(1)}.
\]
All are at most $N^{o(1)}$. The square-root term divided by $\rho$ is bounded by $N^{o(1)}(\rho^{3/2}/\sqrt N+\rho^3/N)\le N^{o(1)}$. Thus every summand is at most $\rho N^{o(1)}$.
The Ward identity gives the column bound. Coordinate boundedness and nonzero diagonals follow from the same compact domain event of Sections~\ref{W:sec:resolvent-domain-for-theta}
and~\ref{W:sec:moment-increments-under-entry-replacement}.

On the right tip take $|\theta|\le s/4$. For $0<s<1$,
\[
 x=2\cosh s\cos\theta\ge2+s^2-\cosh(s)\theta^2
                    \ge2+\tfrac45s^2,
 \qquad |y|=2\sinh s|\sin\theta|\le s^2.
\]
Since $d/8\le s^2/128$, the distance from $x$ to the spectrum is at least $s^2/2$. At height $h=2s^2$ the preceding imaginary part calculation gives the column bound
$\rho N^{o(1)}$; here $\upsilon(x+\mathrm i h)/h\le C/s\le C\rho$. Lemma~\ref{W:lem:vertical-resolvent-comparison} now yields
\[
 \|R(x+\mathrm i y)e_i\|^2\le17\|R(x+2\mathrm i s^2)e_i\|^2
                              \le\rho N^{o(1)}.
\]
For each spectral distance $r=x-\lambda\ge s^2/2$, the inequalities $r^{-1}\le17r/(r^2+4s^4)$ and $r/(r^2+y^2)\ge r^{-1}/5$ show that the real resolvent is bounded above and that
$\operatorname{Re}R_{ii}(x+\mathrm i y)$ is bounded below by a fixed positive constant. The spectral distances also have an absolute upper bound on this contour. Cauchy--Schwarz
for the modulus of the resolvent bounds the off-diagonal entries. Reflection gives the same bounds at the left real boundary point of the contour.

All remaining contour points have $\eta\ge cs^2$; the estimates in the preceding paragraph already cover them. The smallest spectral distance is at least $c\rho^{-2}$, and hence at
least $N^{-1}$ for large $N$. Replacing a spectral parameter by a point of a finite covering set at distance at most $N^{-20}$ changes a coordinate matrix element by at most
$CN^{-18}$ and a column energy by at most $CN^{-17}$, by the continuity estimates in Lemma~\ref{W:lem:vertical-resolvent-comparison}. There are at most $CN^{20}$ contour points,
$N$ possible dimensions and $N^2$ coordinate pairs. Increasing the pointwise failure exponent by $25$ suffices for this union and extension.

The two real barriers of Lemma~\ref{W:lem:resolvent-contour} enclose the spectra of all matrices obtained by the permitted entry replacements. On a support $I$ with $|I|\le2$, its
bounded compression and $\|T\|\le Ce_N$ give $\|(I-R_{II}T)^{-1}\|\le2$. For a column the correction is exactly
\[
 R_{H+\iota_IT\iota_I^*}e_j-Re_j
   =R_{\cdot I}T(I-R_{II}T)^{-1}R_{Ij}.
\]
Thus its coordinate entries are $O(e_N)$ and its norm is at most $Ce_N\sqrt\rho N^{o(1)}$. This verifies all bounds in \eqref{W:eq:resolvent-coordinate} uniformly in $\rho$ and in
the selected entry value. Continuity on the compact ball of permitted entry values proves measurability of the event defined by these simultaneous bounds; the event before deletion
implies it by the same identity.

For a matrix $B_0$ satisfying the upper bound on the Gram matrix, let $t_0$ be the spectral ellipse parameter of $[-2-d/8,2+d/8]$. Then $t_0\le\sqrt d$, $s\ge4t_0$, and
$\rho s\ge4$. The relative error factor in Lemma~\ref{W:lem:moment-weights-gram-resolvent} is $\cosh(\rho t_0)/\sinh(\rho s)<1/4$. Hence
\[
 \|R_{B_0}(z)w\|^2
 \le C\kappa^{-1}\int|z-x|^{-2}\mu_{\rm sc}(\dd x)
 \le C\kappa^{-1}/s\le C\kappa^{-1}\rho.
\]
This proof also applies at the real tips. On a vertical extension the confocal parameter is at least $s$. The relative error factor therefore stays below $1/4$, and
Lemma~\ref{W:lem:moment-weights-gram-resolvent} gives
\[
 w^*|R_{B_0}(z)|w
 \le C_\kappa\int_{-2}^2|z-x|^{-1}\mu_{\rm sc}(\dd x)
 \le C_\kappa(1+\log(2+\rho)).
\]
For the second inequality, outside the tip neighbourhoods $|\operatorname{Im}z|\ge cs^2\ge c\rho^{-2}$, so integration of the bounded density gives
$C(1+|\log|\operatorname{Im}z||)$. Near a real boundary point of the contour the real part is outside $[-2,2]$, and the scalar integral is bounded by $|g(\operatorname{Re}z)|\le1$.
The coordinate spectral integral calculation in the proof of Lemma~\ref{W:lem:resolvent-starting-vector-gram-bound}, with the present lower height $c\rho^{-2}$, gives
\[
 e_i^*|R_{B_0}(z)|e_i\le C(1+\log(2+\rho)).
\]
At a real boundary point of the contour, the bound uses $|R_{B_0}(z)|\le |R_{B_0}(\operatorname{Re}z)|$ and the coordinate estimate already proved on those neighbourhoods.
Normality and Cauchy--Schwarz imply
\[
 |(R_{B_0}(z)w)_i|^2
 \le(e_i^*|R_{B_0}(z)|e_i)(w^*|R_{B_0}(z)|w)\le N^{o(1)}.
\]
Every other matrix obtained by a permitted entry replacement has the form $B_0+\iota_IT\iota_I^*$ with $\|T\|\le Ce_N$. Write $R_0=R_{B_0}$ and $r_0=R_0w$. The formulas for the
restriction to $I$ and for the whole vector are
\[
 r_I=(I-R_{0,II}T)^{-1}(r_0)_I,
 \qquad r=r_0+R_{0,\cdot I}T(I-R_{0,II}T)^{-1}(r_0)_I.
\]
The bound at $B_0$ gives $\|(r_0)_I\|\le N^{o(1)}$. Consequently
\begin{align*}
 \|r\|&\le C\sqrt\rho+
     Ce_N\Big(\sum_{i\in I}\|R_0e_i\|^2\Big)^{1/2}\|(r_0)_I\|
     \le\sqrt\rho N^{o(1)},\\
 \max_j|r_j|&\le\max_j|(r_0)_j|
            +Ce_N\max_j\|R_{0,jI}\|\|(r_0)_I\|
       \le N^{o(1)}.
\end{align*}
Apply the same identities at $\bar z$ to obtain \eqref{W:eq:resolvent-starting-vector}. All estimates hold for every matrix satisfying the stated deterministic bounds, so $B_0$ may
be chosen among the matrices obtained from permitted entry values.

Finally, if no permitted entry value gives the upper bound on the Gram matrix, every Gram matrix corresponding to an entry value in the enlarged ball has an eigenvalue greater than
$\kappa^{-1}$. The determinant cutoff is then zero in a neighbourhood of each permitted entry value and at $K_\alpha$. More precisely, for the composed observable
$F_\alpha(x)=\Phi(\mathbf z_N(K_\alpha+V_\alpha(x),w))$ which contains the determinant cutoff, $F_\alpha$ vanishes on a neighborhood of the support of the replacement entry
distributions and of $x=0$. Thus
\[
 D_x^qF_\alpha(0)=0\quad(q\ge0),\qquad
 F_\alpha(x)-\sum_{q=0}^3\frac1{q!}D_x^qF_\alpha(0)[x^q]=0
\]
for every allowed $x$. These are derivatives of the composed observable; the derivatives of the moment map itself need not vanish.
\end{proof}

Let $p^{(\rho)}(z)=(p_0(z),\ldots,p_{\rho-1}(z))^{\mathsf T}$ and put
\[
 K_\rho(z)=\sum_{j=0}^{\rho-1}|p_j(z)|^2,\qquad
 \nu_\rho^{\mathrm{pol}}(\mathrm dz)=\rho^{-1}K_\rho(z)\frac{|\mathrm dz|}{2\pi}.
\]
The Chebyshev formula gives
\[
 K_\rho(z)\le4\rho e^{2\rho s}|z^2-4|^{-1},\qquad
 \nu_\rho^{\mathrm{pol}}(\mathcal C)\le C e^{2\rho s}(1+|\log s|)=N^{o(1)}.
\]
The last integral follows from $|\mathrm dz|=|\zeta-\zeta^{-1}|\,\mathrm d\theta$.

For the matrix $K_\alpha$, write $R=R_{K_\alpha}$ and define
\[
 a_i(z)=|(R(z)w)_i|+|(R(\bar z)w)_i|+\sqrt{\rho/N},
\]
\[
 U_\alpha=\int a_a(z)a_b(z)\,\nu_\rho^{\mathrm{pol}}(\mathrm dz),\qquad
 V_\alpha=\int(a_a(z)^2+a_b(z)^2)\,\nu_\rho^{\mathrm{pol}}(\mathrm dz).
\]
Set these resolvent weights equal to zero when no permitted entry value gives the upper bound on the Gram matrix. If $E_j$ is the degree-$j$ increment of the Gram matrix and
$B_G=\Delta G-E_1$, the Woodbury identity restricted to the perturbed coordinates gives exactly
\begin{align*}
 E_j={}&\frac1{2\pi\mathrm i}\int_{\mathcal C}
 p^{(\rho)}(z)p^{(\rho)}(z)^{\mathsf T}\,
 w^*R_{\cdot I}T(G_IT)^{j-1}R_{I\cdot}w\,\mathrm dz,
 \qquad j=1,2,3.
\end{align*}
The remainder after degree three has the middle factor $T(G_IT)^3(I-G_IT)^{-1}$. The entire increment beyond $E_1$ has $T G_I T(I-G_IT)^{-1}$. Here $T=V_\alpha(x)_{II}$ and
$G_I=R_{II}$. Since the trace norm of $uv^*$ is $\|u\|\|v\|$, the matrix $p^{(\rho)}(p^{(\rho)})^{\mathsf T}$ has trace norm $K_\rho(z)$. Retaining the two resolvent factors
containing the starting vector and bounding the interior compressions proves
\begin{align}
 \|E_1\|_1&\le C\rho N^{-1/2}\|x\|_1U_\alpha,\nonumber\\
 \|B_G\|_1&\le C\rho N^{-1}\|x\|_1^2V_\alpha,\nonumber\\
 \|E_{\ge4}\|_1&\le C\rho N^{-2}\|x\|_1^4V_\alpha.
 \label{W:eq:resolvent-increments}
\end{align}
The factor $\rho$ comes from the prefactor in the definition of the measure.

\Needspace{6\baselineskip}

\begin{lemma}
\label{W:lem:resolvent-moments}
Uniformly in the deterministic starting vector and in the block set to zero,
\begin{align}
 \sum_{a<b}\mathbb EV_\alpha&\le N^{1+o(1)}\rho,\nonumber\\
 \sum_{a<b}\mathbb E(U_\alpha V_\alpha)
       &\le N^{o(1)}\sqrt{N\rho},\nonumber\\
 \sum_{a<b}\mathbb EV_\alpha^2&\le N^{o(1)}(N+\rho^2),\nonumber\\
 \sum_{a<b}\mathbb EU_\alpha^3
       &\le N^{o(1)}(1+\rho^3/N+\rho e_N^2)
        \le N^{o(1)}(1+\rho e_N^2).
 \label{W:eq:resolvent-weight-sums}
\end{align}
Different terms of a sum may be evaluated at different matrices with the selected entry set to zero.
\end{lemma}
\begin{proof}
Fix one contour point, put $r=Rw$, and delete row and column $i$. The minor vector $c^{(i)}=R^{(i)}w^{(i)}$ satisfies
\[
 c_j^{(i)}=r_j-R_{ji}r_i/R_{ii}.
\]
On the regularity event, Lemma~\ref{W:lem:resolvent-weight-input} therefore gives $\|c^{(i)}\|_\infty\le N^{o(1)}$ and $\|c^{(i)}\|^2\le\rho N^{o(1)}$. These inequalities define an
event measurable with respect to the minor containing the original event. In $r_i=R_{ii}(w_i+h_i^*c^{(i)})$, first bound $R_{ii}$, then enlarge the indicator to the minor event,
and only then condition on the independent row. The enlarged event and $c^{(i)}$ are fixed during this conditional integration.

In these expectations $\mathcal G=\mathcal B_\alpha\cap\mathcal A_{\alpha,\rho}$, with the selected entry $\alpha$ fixed. With $c=c^{(i)}$, apply
\eqref{W:eq:conditional-row-moment} to the independent centered row entries, with the minor fixed. For $p=4,6$ it gives
\[
 \mathbb E_i|h_i^*c-\mathbb E_i(h_i^*c)|^p
 \le C_p\left[(C\|c\|^2/N)^{p/2}
       +\mu_p\|c\|_\infty^{p-2}\|c\|^2\right].
\]
The estimate \eqref{W:eq:conditional-row-moment} bounds each summand using its absolute moments. Setting one entry to zero can only decrease those bounds. The row mean satisfies
$|(\mathbb Eh_i)^*c|\le N^{-2+o(1)}\sqrt N\|c\| \le\sqrt\rho N^{-3/2+o(1)}$. For $p=2,4,6$ this gives
\[
 |(\mathbb Eh_i)^*c|^p
 \le (\rho/N)^{p/2}N^{-p+o(1)}
 \le N^{o(1)}(\rho/N)^{p/2}.
\]
Consequently
\begin{align*}
 \mathbb E\mathbf1_{\mathcal G}|r_i|^2
 &\le N^{o(1)}(|w_i|^2+\rho/N),\\
 \mathbb E\mathbf1_{\mathcal G}|r_i|^4
 &\le N^{o(1)}\big(|w_i|^4+(\rho/N)^2\big),\\
 \mathbb E\mathbf1_{\mathcal G}|r_i|^6
 &\le N^{o(1)}\big(|w_i|^6+(\rho/N)^3
                         +\rho e_N^2/N^2\big).
\end{align*}
Here we used the truncated moments in \eqref{W:eq:resolvent-entrymom}. The same inequalities hold for $a_i$ and conjugate spectral points, with fixed factors.

Set $b_i=|w_i|+\sqrt{\rho/N}$. Then
\[
 \sum_i b_i\le C\sqrt{N\rho},\quad
 \sum_i b_i^2\le C\rho,\quad
 \sum_i b_i^3\le C(1+\rho^{3/2}/\sqrt N),\quad
 \sum_i b_i^4\le C(1+\rho^2/N).
\]
Write $J_\rho=\nu_\rho^{\mathrm{pol}}(\mathcal C)=N^{o(1)}$ and include the regularity indicator in every expectation. The moment bounds and H\"older's inequality give, uniformly
in the spectral arguments,
\begin{align*}
 \mathbb E[a_a(z)a_b(z)a_a(z')^2]
 &\le (\mathbb E a_a(z)^4)^{1/4}
       (\mathbb E a_b(z)^4)^{1/4}
       (\mathbb E a_a(z')^4)^{1/2}
       \le N^{o(1)}b_a^3b_b,\\
 \mathbb E[a_a(z)^2a_b(z')^2]
 &\le (\mathbb E a_a(z)^4)^{1/2}
       (\mathbb E a_b(z')^4)^{1/2}
       \le N^{o(1)}b_a^2b_b^2.
\end{align*}
For the terms with a repeated index in the second line use $b_a^4$. Expansion of the positive integrands and Tonelli's theorem now yield
\begin{align*}
 \sum_{a<b}\mathbb EV_\alpha
 &\le N^{o(1)}J_\rho(N-1)\sum_i b_i^2
       \le N^{1+o(1)}\rho,\\
 \sum_{a<b}\mathbb E(U_\alpha V_\alpha)
 &\le N^{o(1)}J_\rho^2
                  \sum_{a<b}(b_a^3b_b+b_ab_b^3)\\
 &\le N^{o(1)}\Big(\sum_i b_i^3\Big)\Big(\sum_i b_i\Big)
       \le N^{o(1)}\sqrt{N\rho},\\
 \sum_{a<b}\mathbb EV_\alpha^2
 &\le N^{o(1)}J_\rho^2
               \sum_{a<b}(b_a^4+2b_a^2b_b^2+b_b^4)\\
 &\le N^{o(1)}\left[N\sum_i b_i^4+
                           \Big(\sum_i b_i^2\Big)^2\right]
       \le N^{o(1)}(N+\rho^2).
\end{align*}
We used $\rho^3/N\le N^{o(1)}$, which follows from $\rho\le C N^{1/3}$ in the present range. For $U^3$, H\"older with six equal exponents gives
\[
 \mathbb E\prod_{j=1}^3a_a(z_j)a_b(z_j)
 \le\prod_{j=1}^3(\mathbb E a_a(z_j)^6)^{1/6}
                         (\mathbb E a_b(z_j)^6)^{1/6}.
\]
Use the deterministic bound
\[
 \big(\mathbb E\mathbf1_{\mathcal G}a_i^6\big)^{1/2}
 \le N^{o(1)}\left(|w_i|^3+
        \rho^{3/2}N^{-3/2}+\sqrt\rho e_N/N\right).
\]
Its sum is at most $N^{o(1)}(1+\rho^{3/2}/\sqrt N+\sqrt\rho e_N)$. If the last displayed bound is denoted by $d_i$, the bound for the product of six resolvent factors is
$N^{o(1)}d_ad_b$. Therefore
\[
 \sum_{a<b}\mathbb EU_\alpha^3
 \le N^{o(1)}J_\rho^3\sum_{a<b}d_ad_b
 \le N^{o(1)}\Big(\sum_i d_i\Big)^2
 \le N^{o(1)}(1+\rho^3/N+\rho e_N^2).
\]
This is the fourth line of \eqref{W:eq:resolvent-weight-sums}. For every $\alpha$, the expectation under the law of $K_\alpha$ is bounded by the displayed deterministic expression
in $w_i$. Summing these bounds proves the summed estimate for the full family $(K_\alpha)_\alpha$.
\end{proof}

\subparagraph{Summed conditional weights}
\label{W:sec:sums-of-resolvent-weights-at-the-critical-scale}
We now express the moment estimates on the regularity event in the norm \eqref{W:eq:comparison-direction-norm}, using the truncated and bounded entry distributions constructed
above. Their exact moment differences are recorded after the resolvent weight lemma. All events and expectations below precede integration over the replacement entry.

At a fixed interpolation parameter, set one entry block to zero and write
\[
 K_\alpha=H(t)^{[\alpha\leftarrow0]},\qquad
 V_\alpha(x)=N^{-1/2}\mathcal T_\alpha x,\qquad
 t_x=N^{-1/2}\|x\|_1.
\]
The variable $x$ belongs to $\mathbb R$ or $\mathbb R^2$, according to the symmetry class. For $N\ge N_0$ with $Ce_N<1$, every allowed replacement obeys $t_x\le Ce_N<1$, including
on the diagonal. The finite-dimensional moment curve is $u\mapsto s(K_\alpha+uV_\alpha(x),w)$. Write its exact decomposition as
\[
 \Delta s=\ell+q+c+r,\qquad
 \ell=s'(0),\quad q=s''(0)/2,\quad c=s'''(0)/6.
\]
Since $s_n=u_n+x_n/m_N$, the definitions give, for every matrix $K$,
\[
 s_n(K,w)=\frac1N\operatorname{Tr}p_n(h_d(K))
       +\operatorname{Tr}\bigl[(ww^*-N^{-1}I)p_n(h_d(K))\bigr]
       =w^*p_n(h_d(K))w.
\]
Differentiating gives the corresponding equalities for $\ell,q,c$; subtracting these from the increment identity gives the equality for $r$.

Let $\mathcal B_\alpha$ and $\mathcal A_{\alpha,\rho}$ be, respectively, the regularity event measurable with respect to the remaining entries and the event that a permitted entry
value gives the bound on the Gram matrix in Lemma~\ref{W:lem:resolvent-weight-input}; denote their intersection by $I_\alpha=1$. On $\mathcal B_\alpha$ the spectral truncation is
locally inactive throughout the set of permitted entry values. On $\mathcal B_\alpha\cap\mathcal A_{\alpha,\rho}^c$, the vanishing property in
Lemma~\ref{W:lem:resolvent-weight-input} makes both the observable and its expansion coefficients zero. The complementary event $\mathcal B_\alpha^c$ is estimated separately below.

\Needspace{6\baselineskip}

For the norm $\|\cdot\|_{\mathscr D}$ of \eqref{W:eq:comparison-direction-norm}, use the measure
\[
 \nu_\rho(\mathrm dz)=E_\rho(z)^2\frac{|\mathrm dz|}{2\pi},\qquad
 E_\rho(z)=\frac{2e^{\rho s}}{|z^2-4|^{1/2}}.
\]
The polynomial measure $\nu_\rho^{\mathrm{pol}}$ defined above instead averages $|p_j(z)|^2$. They satisfy $\nu_\rho^{\mathrm{pol}}\le\nu_\rho$ because
$\rho^{-1}\sum_{j<\rho}|p_j(z)|^2\le E_\rho(z)^2$. The direct contour calculation in \eqref{W:eq:comparison-measure} gives $\nu_\rho(\mathcal C)=N^{o(1)}$. In the following lemma
$U_\alpha$ and $V_\alpha$ are integrated against $\nu_\rho$, the measure in the definition of this norm.

\begin{lemma}
\label{W:lem:comparison-weights}
There are nonnegative resolvent weights $U_\alpha,V_\alpha$, measurable with respect to the remaining entries, such that on $I_\alpha=1$,
\begin{align}
 \|\ell\|_{\mathscr D}&\le Ct_xU_\alpha,&
 \|q\|_{\mathscr D}+\|\Delta s-\ell\|_{\mathscr D}
       &\le Ct_x^2V_\alpha,\nonumber\\
 \|c\|_{\mathscr D}&\le Ct_x^3V_\alpha,&
 \|r\|_{\mathscr D}&\le Ct_x^4V_\alpha.
 \label{W:eq:comparison-increments}
\end{align}
For diagonal replacements $U_\alpha=V_\alpha$ and $\|\Delta s\|_{\mathscr D}\le Ct_xV_\alpha$.
\end{lemma}
\begin{proof}

Polynomial functional calculus and \eqref{W:eq:resolvent-block-remainder}, evaluated at $K_\alpha$, give
\[
 v(f)=\frac1{2\pi\mathrm i}\int_{\mathcal C}f(z)w^*A(z)w\,\dd z
\]
for each homogeneous increment and for the remaining difference. Each such $v$ annihilates constants because $w^*Iw=1$. Using this representing kernel in the definition
\eqref{W:eq:comparison-direction-norm} gives
\[
 \|v\|_{\mathscr D}\le\int_{\mathcal C}|w^*A(z)w|\,\nu_\rho(\dd z).
\]
Write $R=(zI-K_\alpha)^{-1}$, let $J$ be the coordinate support of the selected entry, and put $G_J=R_{JJ}$ and $T=V_\alpha(x)_{JJ}$. Lemma~\ref{W:lem:resolvent-weight-input} gives
$\|G_J\|\le C$ and $\|(I-G_JT)^{-1}\|\le2$ uniformly over the permitted entry values.

Choose $Q_N=N^{o(1)}$ dominating the constants in the preceding estimates in the column norm bounds and set
\[
 a_i(z)=|(R(z)w)_i|+|(R(\bar z)w)_i|+\sqrt{Q_N\rho/N}.
\]
Put $l=(w^*R)_J$ and $r=(Rw)_J$. For an off-diagonal block $J=\{a,b\}$, the perturbation matrix satisfies $T_{aa}=T_{bb}=0$, so
\[
 |lTr|=|T_{ab}l_ar_b+T_{ba}l_br_a|\le2t_xa_a a_b.
\]
For $j=2,3$,
\[
 |lT(G_JT)^{j-1}r|
 \le\|l\|_2\|r\|_2\|T\|^j\|G_J\|^{j-1}
 \le C_jt_x^j(a_a^2+a_b^2).
\]
The remainder after degree three and the increment beyond the linear term obey
\[
 \begin{aligned}
 |lT(G_JT)^3(I-G_JT)^{-1}r|&\le C t_x^4(a_a^2+a_b^2),\\
 |lTG_JT(I-G_JT)^{-1}r|&\le C t_x^2(a_a^2+a_b^2),
 \end{aligned}
\]
because $\|G_JT\|\le C e_N\le1/2$ and the inverse has norm at most two. Define
\[
 U_\alpha=\int a_aa_b\,\mathrm d\nu_\rho,\qquad
 V_\alpha=\int(a_a^2+a_b^2)\,\mathrm d\nu_\rho.
\]
For a diagonal block $J=\{a\}$, the same products are scalars; the linear term is bounded by $t_xa_a^2$, the terms of degree $j=2,3$ by $C_jt_x^ja_a^2$, and the two tails by
$Ct_x^4a_a^2$ and $Ct_x^2a_a^2$. Use $U_\alpha=V_\alpha=\int a_a^2\,\dd\nu_\rho$. Integrating each bound against $\nu_\rho$ in the norm estimate for the moment directions proves
\eqref{W:eq:comparison-increments} and the additional diagonal estimate. Each kernel is computed from $K_\alpha$, so its weight is measurable before integration over the selected
entry. The identity $\Delta s_n=\Delta u_n+\Delta x_n/m_N$ cancels the empirical trace exactly in every term.

\end{proof}

\begin{corollary}
\label{W:cor:integrated-comparison-weights}
For off-diagonal blocks the sums restricted to the regularity event, including integration in the interpolation parameter, satisfy
\begin{equation}
\begin{gathered}
\begin{array}{c|c@{\qquad}c|c}
\text{resolvent weight}&\text{summed bound}&\text{resolvent weight}&\text{summed bound}\\\hline
 U&{m_N}^4N^{o(1)}&U^2&{m_N}^2N^{o(1)}\\
 V&{m_N}^4N^{o(1)}&UV&{m_N}^2N^{o(1)}\\
 V^2&{m_N}^3N^{o(1)}&U^3&{m_N}N^{o(1)}\\
 U^4&S_NN^{o(1)}&&
\end{array}\\[2pt]
 S_N=1+m_Ne_N^4\le Cm_N.
\end{gathered}
 \label{W:eq:comparison-weights}
\end{equation}
For diagonal replacements,
\begin{equation}
 \sum_{\alpha\ {\rm diag}}\mathbb EV_\alpha\le {m_N}N^{o(1)},\qquad
 \sum_{\alpha\ {\rm diag}}\mathbb EV_\alpha^2\le N^{o(1)}.
 \label{W:eq:comparison-diag-weights}
\end{equation}
Every expectation and sum in this corollary includes the specified restriction to the regularity event for the remaining entries. The bounds are uniform over deterministic unit
starting vectors and interpolation parameters.
\end{corollary}
\begin{proof}
On the regularity event, the bounds on resolvent entries, squared column norms, and the Gram matrix in Lemma~\ref{W:lem:resolvent-weight-input} imply the hypotheses of
Proposition~\ref{W:prop:moment-weights-row-moments}, with the squared norm of the minor coefficient vector bounded by $\rho N^{o(1)}$ and its largest coordinate bounded by
$N^{o(1)}$. The enlargement from the row dependent event to an event measurable with respect to the minor is exactly \eqref{W:eq:conditional-minor-domination}. For uncentered
entries, the row mean contributes the term bounded in \eqref{W:eq:conditional-row-with-mean}.

For each fixed block $\alpha$ and point $z$, this gives deterministic bounds, uniform in those parameters,
\begin{align}
 \mathbb E[I_\alpha a_i(z)^2]
   &\le Q_N'\bigl(|w_i|^2+\rho/N\bigr),\nonumber\\
 \mathbb E[I_\alpha a_i(z)^4]
   &\le Q_N'\bigl(|w_i|^4+(\rho/N)^2\bigr),\nonumber\\
 \mathbb E[I_\alpha a_i(z)^6]
   &\le Q_N'\bigl(|w_i|^6+(\rho/N)^3+\rho e_N^2/N^2\bigr),\nonumber\\
 \mathbb E[I_\alpha a_i(z)^8]
   &\le Q_N'\bigl(|w_i|^8+(\rho/N)^4+\rho e_N^4/N^2\bigr),
 \label{W:eq:weights-coordinate-moments}
\end{align}
where $Q_N'=N^{o(1)}$. To verify all four bounds, write $c=c^{(i)}$ for the minor coefficient vector and apply \eqref{W:eq:conditional-row-moment} with the minor fixed and the
centered row entries. For even $p\ge4$,
\[
 \mathbb E_i|h_i^*c-\mathbb E_i(h_i^*c)|^p
 \le C_p\left[(C\|c\|_2^2/N)^{p/2}
       +\mu_p\|c\|_\infty^{p-2}\|c\|_2^2\right],
 \qquad \mu_p\le CF_NN^{-2}e_N^{p-4}.
\]
Using $\|c\|_2^2\le\rho N^{o(1)}$ and $\|c\|_\infty\le N^{o(1)}$ gives the terms $(\rho/N)^{p/2}N^{o(1)}$ and $\rho e_N^{p-4}N^{-2+o(1)}$. For $p=4$, the latter is absorbed by the
former because $\rho\ge1$. For $p=6,8$, both contributions appear in \eqref{W:eq:weights-coordinate-moments}. The variance identity gives $p=2$. The row mean is at most
$\sqrt\rho N^{-3/2+o(1)}$; its $p$th power, and the deterministic additive term, are bounded by the variance term. Adding $|w_i|^p$ proves all four moment bounds. The selected
diagonal entry is absent from the off-diagonal row being integrated.

Let $J_N=\nu_\rho(\mathcal C)$; by \eqref{W:eq:comparison-measure}, $J_N=N^{o(1)}$. Replace $Q_N'$ by $C Q_N'(1+J_N)^4$, which still satisfies $Q_N'=N^{o(1)}$. Put
$b_i=|w_i|+\sqrt{\rho/N}$. The first two lines of \eqref{W:eq:weights-coordinate-moments} are bounded by $Q_N'b_i^2$ and $Q_N'b_i^4$. Keep $\alpha$ fixed when taking the joint
expectation of products of contour integrals at $K_\alpha$. For example, by Tonelli and H\"older,
\begin{equation}
 \begin{split}
 \mathbb E[I_\alpha U_\alpha V_\alpha]
 &=\iint\mathbb E\!\left[
 I_\alpha a_a(z)a_b(z)
            \bigl(a_a(z')^2+a_b(z')^2\bigr)\right]
             \,\nu_\rho(\dd z)\nu_\rho(\dd z')\\
 &\le C Q_N' J_N^2\bigl(b_a^3b_b+b_ab_b^3\bigr).
 \end{split}
 \label{W:eq:weights-cross-contour}
\end{equation}
For its first term the H\"older exponents are $4,4,2$, applied to $a_a(z)$, $a_b(z)$, and $a_a(z')^2$, with the same indicator included in each factor. Thus only fourth coordinate
moments occur, even though the contour parameters differ.

For $U$ and $U^2$, use H\"older with exponents $(2,2)$ and $(4,4,4,4)$, respectively. The second moment bound controls $V$; H\"older with exponents $(2,2)$ on the squared factors
controls $V^2$. Integrate these bounds by Tonelli's theorem. Consequently
\begin{align*}
 \mathbb E[I_\alpha U_\alpha]&\le Q_N'J_N b_ab_b,&
 \mathbb E[I_\alpha U_\alpha^2]&\le Q_N'J_N^2b_a^2b_b^2,\\
 \mathbb E[I_\alpha V_\alpha]&\le Q_N'J_N(b_a^2+b_b^2),&
 \mathbb E[I_\alpha V_\alpha^2]&\le
 Q_N'J_N^2(b_a^4+2b_a^2b_b^2+b_b^4).
\end{align*}
The two mixed choices $(i,j)=(a,b),(b,a)$ account for the coefficient $2$ in the last formula. The deterministic coordinate weights satisfy
\[
 \sum_i b_i\le C\sqrt{N\rho},\qquad
 \sum_i b_i^2\le C\rho,\qquad
 \sum_i b_i^3\le C(1+\rho^{3/2}/\sqrt N),\qquad
 \sum_i b_i^4\le C(1+\rho^2/N).
\]
Throughout the stated range $8\le\rho\le\rho_*\le C_*m_N$, the last two sums are bounded: $\rho^{3/2}/\sqrt N\le C_*^{3/2}$ and $\rho^2/N\le C_*^2m_N^{-1}$. This upper bound on
$\rho$ is the one used in the following sums. Summing the preceding pair estimates therefore gives respectively
\[
 \sum\mathbb E[I U]\le N\rho N^{o(1)},\quad
 \sum\mathbb E[I U^2]\le\rho^2N^{o(1)},\quad
 \sum\mathbb E[I V]\le N\rho N^{o(1)},
\]
\[
 \sum\mathbb E[I UV]
 \le(\sum_i b_i^3)(\sum_i b_i)N^{o(1)},\qquad
 \sum\mathbb E[I V^2]
 \le\bigl[N\sum_i b_i^4+(\sum_i b_i^2)^2\bigr]N^{o(1)}.
\]
These are the first five entries of \eqref{W:eq:comparison-weights}.

For $k=3,4$, let $M_{2k,i}$ denote the deterministic right side of the corresponding line of \eqref{W:eq:weights-coordinate-moments}. Jensen's inequality for the probability
measure $\nu_\rho/J_N$, followed by Cauchy--Schwarz inside the expectation, gives
\begin{equation}
 \begin{split}
 \mathbb E[I_\alpha U_\alpha^k]
 &\le J_N^{k-1}\int
       \mathbb E[I_\alpha a_a(z)^ka_b(z)^k]\,\nu_\rho(\dd z)\\
 &\le J_N^k(M_{2k,a}M_{2k,b})^{1/2},\\
 \sum_{a<b}\mathbb E[I_{ab}U_{ab}^k]
 &\le \tfrac12J_N^k\left(\sum_i\sqrt{M_{2k,i}}\right)^2.
 \end{split}
 \label{W:eq:weights-higher-pair-sums}
\end{equation}
The deterministic bounds $M_{2k,i}$ are uniform over pairs. We can therefore sum them even when different pairs use different matrices. Their definitions give
\begin{align*}
 \sum_i\sqrt{M_{6,i}}
 &\le N^{o(1)}(1+\rho^{3/2}/\sqrt N+\sqrt\rho\,e_N),\\
 \sum_i\sqrt{M_{8,i}}
 &\le N^{o(1)}(1+\rho^2/N+\sqrt\rho\,e_N^2).
\end{align*}
Using $(a+b+c)^2\le3(a^2+b^2+c^2)$, their squared bounds are
\[
 \begin{aligned}
 \left(\sum_i\sqrt{M_{6,i}}\right)^2
  &\le N^{o(1)}(1+\rho^3/N+\rho e_N^2),\\
 \left(\sum_i\sqrt{M_{8,i}}\right)^2
  &\le N^{o(1)}(1+\rho^4/N^2+\rho e_N^4).
 \end{aligned}
\]
Using $\rho\le C_*m_N$ gives $\rho^3/N\le C_*^3$, $\rho^4/N^2\le C_*^4m_N^{-2}$, and
\[
 1+\rho^3/N+\rho e_N^2\le C(1+m_Ne_N^2)\le Cm_N,
\]
\[
 1+\rho^4/N^2+\rho e_N^4
       \le C(1+m_Ne_N^4).
\]
Multiplication by $J_N^3$ and $J_N^4$, respectively, in \eqref{W:eq:weights-higher-pair-sums} therefore proves $\sum\mathbb E[I U^3]\le m_NN^{o(1)}$ and
$\sum\mathbb E[I U^4]\le(1+m_Ne_N^4)N^{o(1)}$ for the full permitted range of $\rho$. For a diagonal block, Jensen and the first two coordinate moments give
\[
 \sum_i\mathbb E[I_iV_i]\le
 N^{o(1)}\sum_i b_i^2\le {m_N}N^{o(1)},\qquad
 \sum_i\mathbb E[I_iV_i^2]\le
 N^{o(1)}\sum_i b_i^4\le N^{o(1)}.
\]
This proves \eqref{W:eq:comparison-diag-weights}. All inequalities are uniform in a fixed interpolation parameter. Tonelli's theorem and integration over its interval of length one
therefore preserve every bound.
\end{proof}

With the same normalization of the real coordinate vector, Lemma~\ref{W:lem:cubic-defect} constructs $H^{\rm b}_{\alpha,\rm std}$ with the first three coordinate moment tensors of
$X^{\rm nor}$. Thus its mean is zero and its covariance is the identity. The reconstruction $\widehat X=\mu+\Sigma^{1/2}X^{\rm nor}$ gives $\Delta M_1=\mu$ and
$\Delta M_2=\Sigma-I+\mu\mu^{\mathsf T}$. For $\Delta M_3$, apply \eqref{W:eq:cubic-affine}, whose three terms are bounded in the proof of that lemma. Together these give
\[
 \|\Delta M_2\|+\|\Delta M_3\|
 \le C\bigl(\|\Sigma-I\|+\|\mu\|+\|\mu\|^2+\|\mu\|^3\bigr).
\]
The truncation estimates and $L=e_N\sqrt N$ now give
\begin{equation}
 \|\Delta M_1\|\le N^{-3/2+o(1)},\qquad
 \|\Delta M_2\|+\|\Delta M_3\|\le N^{-1+o(1)}.
 \label{W:eq:comparison-defects}
\end{equation}
Indeed $\|\mu\|\le C\Theta(L)L^{-3}$ and $\|\Sigma-I\|\le C\Theta(L)L^{-2}$, while $L\to\infty$ and $e_N^{-1}=N^{o(1)}$. Thus
\[
 \Theta(L)L^{-3}\le N^{-3/2+o(1)},\qquad
 \Theta(L)(L^{-2}+L^{-3})\le N^{-1+o(1)}.
\]
Only boundedness of $\Theta(L)$ is needed for these powers. The third absolute moments are bounded, and the truncated fourth absolute moments obey
\[
 \mathbb E\|\widehat X\|^4
 \le C\left(1+\int_1^{CL}t^3\mathbb P(\|X\|>t)\,\dd t\right)
 \le C\log(e+L)=N^{o(1)}.
\]
The bounded comparison vector has the same uniform third absolute moment bound and a bounded fourth moment. For a diagonal standardized entry $\xi$ with $\mathbb E\xi=0$ and
$\mathbb E\xi^2\le C$, the pointwise inequality $|\xi|\mathbf1_{\{|\xi|>L\}}\le\xi^2/L$ gives
\[
 |\mathbb E[\xi\mathbf1_{\{|\xi|\le L\}}]|
 =|\mathbb E[\xi\mathbf1_{\{|\xi|>L\}}]|
 \le L^{-1}\mathbb E\xi^2\le C/L=N^{-1/2+o(1)}.
\]
The diagonal comparison variable is centered, so this is its first-moment difference. Both diagonal second moments remain bounded.

On the support of $\chi_N$ or a derivative of $\chi_N$, $B^4Q_*^8\le C\mathfrak r_N^{12}$. Let $\mathcal H_\alpha$ denote this support condition at the point where the derivative
is evaluated. For any nonnegative weight $W_\alpha$,
\[
 \mathbb E[I_\alpha\mathbf1_{\mathcal H_\alpha}
                 B^4Q_*^8W_\alpha]
 \le C\mathfrak r_N^{12}\mathbb E[I_\alpha W_\alpha].
\]
For each of the four choices $W_\alpha=V_\alpha,U_\alpha V_\alpha,V_\alpha^2,U_\alpha^4$, Tonelli and the pointwise bound give
\[
 \int_0^1\sum_\alpha\mathbb E[
  I_\alpha\mathbf1_{\mathcal H_\alpha}B^4Q_*^8W_\alpha]\,\dd t
 \le C\mathfrak r_N^{12}\int_0^1\sum_\alpha
                         \mathbb E[I_\alpha W_\alpha]\,\dd t.
\]
Corollary~\ref{W:cor:integrated-comparison-weights} bounds the last integrals. Since $\log(\mathfrak r_N^{12})/\log N=12\log\mathfrak r_N/\log N\to0$, this multiplication leaves
their displayed powers of $N$ unchanged. The sequence $\mathfrak r_N=N^{o(1)}$ is chosen using the bounds for $H_N^{\rm b}$ and $H_N^{\mathrm G}$ in
Lemma~\ref{W:lem:comparison-ensemble-coefficients}. All resulting losses are subpolynomial. On arbitrary line segments joining real moment vectors we use the global derivatives of
Proposition~\ref{W:prop:comparison-global}.

The entry distributions, regularity events, and integrated weights are now available for comparison. The remaining cubic contribution is summed with its signs in
Section~\ref{W:sec:green-function-comparison}. Its derivatives include the factor $m_N$ from the integrated potential. The final expectation bounds proved there are $N^{-1/3+o(1)}$
and $N^{-1/6+o(1)}$. For fixed $c>0$ and $L_N'=N^{o(1)}$, choose $N_0(c)$ so that $L_N'\le N^{c/2}$ for $N\ge N_0(c)$. Then $N^{-c}L_N'\le N^{-c/2}\to0$. The probabilities of
disagreement under truncation enter separately when returning to the original matrices.

\subsection{Green function comparison}
\label{sec:combined-w-comparison}
\label{W:sec:green-function-comparison}

We combine the deterministic derivative bounds with the conditional moment estimates to compare expectations of the observables from
Section~\ref{W:sec:orthogonal-polynomials-and-spectral-statistics}, whose derivatives have the scale $m_N=N^{1/3}$. The two comparisons require different estimates. Between
$H_N^{\rm tr}$ and $H_N^{\rm b}$, the errors are the small moment differences of the uncentered truncation and the integral remainder. Between $H_N^{\rm b}$ and the GOE or GUE, the
first two off-diagonal moment tensors agree, but the third moment difference can be of order one.

The cubic contribution in the latter comparison is treated first. Its leading term contains the first derivative of the observable, whose growth prevents a direct absolute estimate
of the entry sum. Evaluating that derivative and its resolvents at one common matrix $H(t)$ exposes the signed sum as a matrix element involving $B_{3,N}$. A trace norm bound then
controls it. The other cubic contributions admit absolute summation.

The comparison with the truncated matrix also includes the first two raw moment differences and the smaller third moment difference from the normalized matching construction.
Diagonal entries are treated using their signed mean and a remainder integrated against the sum of the two second moments. For off-diagonal entries, only the remainder is split
according to the size of the entry replacement; the polynomial part remains integrated against the full entry distributions.

Once these sums are controlled, comparison of $\Gamma_\rho(1-\chi_N)$ bounds the error from removing the coefficient cutoff. The cutoff in the eigenvalues of the Gram matrix is
removed separately, with $N\to\infty$ first at fixed cutoff parameter. Characteristic functions and modulus bounds then identify the coefficient process and its joint limit with
the off-diagonal deviations. Returning to the original ensemble adds the stated $o(1)$ errors from truncation and spectral modification. The conclusion is convergence in
distribution of the coefficient process. The comparison errors below concern bounded functions of its polynomial moments.

\subsubsection{Interpolation of the entry distributions}
\label{W:sec:interpolation-of-entry-distributions}
\label{W:sec:intro-entry-comparison}

Our comparison changes the distributions of the finitely many independent upper triangular entries. This is the discrete entry replacement argument underlying Green function
comparison \cite{erdos2012rigidity}, written below as an integral in a scalar parameter. The parameter indexes product probability measures. The relevant derivative is obtained by
differentiating those measures, so the calculation below specifies the comparison directly at each finite matrix size.

The index set $\mathcal I_N$ consists of the independent upper-triangular entries. For $\alpha\in\mathcal I_N$, the two distributions on its real coordinate space are denoted by
$\mu_\alpha$ and $\nu_\alpha$. We interpolate between these distributions and form the product measure as follows:

\[
 \lambda_{\alpha,t}=(1-t)\nu_\alpha+t\mu_\alpha,\qquad
 \lambda_t=\bigotimes_{\alpha\in\mathcal I_N}\lambda_{\alpha,t},
 \qquad 0\le t\le1.
\]
The linear map $V_\alpha(x)=N^{-1/2}\mathcal T_\alpha x$ places the real coordinate vector $x$ in entry $\alpha$ and, off the diagonal, its conjugate entry. The coordinate
dimension is one for a diagonal or real off-diagonal entry and two for a complex off-diagonal entry. Define $K_\alpha(t)$ by setting entry $\alpha$ to zero and assigning the other
entries the product distribution $\bigotimes_{\delta\ne\alpha}\lambda_{\delta,t}$. The sigma field $\mathscr F_\alpha=\sigma(K_\alpha(t))$ is fixed before integration over $x$.

For comparison with successive replacement of individual entries, choose an enumeration $\alpha_1,\ldots,\alpha_m$, put
\[
 \Lambda^{(j)}=\bigotimes_{l\le j}\mu_{\alpha_l}\otimes
                    \bigotimes_{l>j}\nu_{\alpha_l},\qquad 0\le j\le m.
\]
Telescoping gives $\mathbb E_{\Lambda^{(m)}}F-\mathbb E_{\Lambda^{(0)}}F=\sum_{j=1}^m(\mathbb E_{\Lambda^{(j)}}F-\mathbb E_{\Lambda^{(j-1)}}F)$. In each summand all entry
distributions except the $j$th agree. The following identity averages such single-entry changes over the product family $\lambda_t$.

\begin{lemma}
\label{W:lem:entry-distribution-interpolation}
For every bounded Borel function $F$ on the Hermitian matrices,
\begin{equation}
 \mathbb E_{\lambda_1}F-\mathbb E_{\lambda_0}F
 =\int_0^1\sum_{\alpha\in\mathcal I_N}\mathbb E_{K_\alpha(t)}
       \int F(K_\alpha(t)+V_\alpha(x))
                     (\mu_\alpha-\nu_\alpha)(\mathrm dx)\,\mathrm dt.
 \label{W:eq:entry-distribution-interpolation}
\end{equation}
\end{lemma}
\begin{proof}
Enumerate $\mathcal I_N=\{\alpha_1,\ldots,\alpha_m\}$. For each $J\subseteq\{1,\ldots,m\}$ define
\[
 \lambda_J=\bigotimes_{j=1}^m\lambda_{j,J},\qquad
 \lambda_{j,J}=\begin{cases}\mu_{\alpha_j},&j\in J,\\
                                  \nu_{\alpha_j},&j\notin J.
                   \end{cases}
\]
Write $\mathcal H(x_1,\ldots,x_m)=\sum_{j=1}^mV_{\alpha_j}(x_j)$. Expanding each factor $(1-t)\nu_{\alpha_j}+t\mu_{\alpha_j}$ gives $2^m$ product measures. The term indexed by $J$
chooses $\mu_{\alpha_j}$ exactly at the indices $j\in J$, and therefore has coefficient $t^{|J|}(1-t)^{m-|J|}$. Consequently
\[
 \phi(t):=\mathbb E_{\lambda_t}F
  =\sum_{J\subseteq\{1,\ldots,m\}}
      t^{|J|}(1-t)^{m-|J|}\int F(\mathcal H(x))\,\lambda_J(\mathrm dx).
\]
All coefficients are finite because $F$ is bounded. In particular $\phi$ is a polynomial, including at $t=0,1$. Differentiating the factor at position $j$ gives the signed measure
$\mu_{\alpha_j}-\nu_{\alpha_j}$; the other factors remain $\lambda_{\alpha_l,t}$. Therefore
\[
 \phi'(t)=\sum_{j=1}^m
  \int F\!\left(\sum_{l\ne j}V_{\alpha_l}(x_l)
                         +V_{\alpha_j}(x_j)\right)
       (\mu_{\alpha_j}-\nu_{\alpha_j})(\mathrm dx_j)
       \prod_{l\ne j}\lambda_{\alpha_l,t}(\mathrm dx_l).
\]
For each $j,t$, the absolute integral is at most $2\|F\|_\infty$, since the remaining factors have mass one and $\|\mu_{\alpha_j}-\nu_{\alpha_j}\|_{\rm TV}\le2$. The sum of these
bounds is integrable over $0\le t\le1$. The pushforward of $\prod_{l\ne j}\lambda_{\alpha_l,t}$ by $(x_l)_{l\ne j}\mapsto\sum_{l\ne j}V_{\alpha_l}(x_l)$ is exactly the law of
$K_{\alpha_j}(t)$. Substituting this pushforward in $\phi'(t)$ and using $\phi(1)-\phi(0)=\int_0^1\phi'(t)\,\mathrm dt$ proves \eqref{W:eq:entry-distribution-interpolation}.
\end{proof}

Conditional on $K_\alpha$, the scalar function $f_\alpha(x)=F(K_\alpha+N^{-1/2}\mathcal T_\alpha x)$ has the exact expansion
\begin{equation}
 f_\alpha(x)=\sum_{j=0}^3\frac1{j!}D^jf_\alpha(0)[x^j]
 +\frac1{3!}\int_0^1(1-t)^3D^4f_\alpha(tx)[x^4]\,\dd t.
 \label{W:eq:intro-entry-taylor}
\end{equation}
For $\psi(t)=f_\alpha(tx)$, the repeated chain rule gives $\psi^{(j)}(t)=D^jf_\alpha(tx)[x^j]$. Repeated integration of $\psi^{(4)}$ over $0\le t_4\le t_3\le t_2\le t_1\le1$ gives
a weight $(1-t_4)^3/3!$, which proves \eqref{W:eq:intro-entry-taylor}. The map $x\mapsto V_\alpha(x)$ is linear, and hence
\[
 D^jf_\alpha(x)[u_1,\ldots,u_j]
   =N^{-j/2}D_H^jF(K_\alpha+N^{-1/2}\mathcal T_\alpha x)
                    [\mathcal T_\alpha u_1,\ldots,\mathcal T_\alpha u_j].
\]
Thus every coefficient already contains its factor $N^{-j/2}$. Matching the first two entry tensors cancels the corresponding terms after integration against the signed measure
$\mu_\alpha-\nu_\alpha$. For a diagonal entry the proof instead uses
\[
 f_\alpha(x)=f_\alpha(0)+Df_\alpha(0)[x]
                 +\int_0^1(1-t)D^2f_\alpha(tx)[x,x]\,\dd t.
\]
Its quadratic remainder is bounded in Theorem~\ref{W:thm:truncation-diagonal} using the diagonal second moments. For the laws with the stated moment assumptions, the truncated
fourth moments and conditional resolvent weights are estimated before their sum over $\alpha$.

In the applications, let $\mathscr F_\alpha$ be the sigma field generated by the remaining entries. If $I_\alpha$ and a tensor $C_{q,\alpha}$ are $\mathscr F_\alpha$-measurable,
then, whenever the indicated moments exist,
\begin{equation}
 \mathbb E\left[I_\alpha\int C_{q,\alpha}:x^{\otimes q}
                         (\mu_\alpha-\nu_\alpha)(\mathrm dx)\right]
 =\mathbb E\left[I_\alpha C_{q,\alpha}:
               \bigl(M_q(\mu_\alpha)-M_q(\nu_\alpha)\bigr)\right].
 \label{W:eq:conditional-moment-contraction}
\end{equation}
Here $M_q(\mu)=\int x^{\otimes q}\mu(\mathrm dx)$, and the colon denotes summation over all ordered coordinate indices. The degree-$q$ coefficient below includes the factor $1/q!$.
The identity requires integrability of the coefficient times the indicated moment; in our applications this follows from the deterministic derivative bounds below. The error after
subtraction of the polynomial is bounded using the positive measure $\mu_\alpha+\nu_\alpha$, of total mass two.

Throughout this subsection, $c m_N\le\rho\le C m_N$, the starting vector $w$ is deterministic and has norm one, and the interval, Fourier, modulus, and cutoff parameters are fixed.
Write $\mathcal P_\rho$ for the real polynomials of degree at most $2\rho-2$. We use the normalized moment functional $s(H,w)$ from
Section~\ref{W:sec:the-moment-vector-and-its-smooth-extensions} and the norm $\|\cdot\|_\diamond$ in \eqref{W:eq:tightness-norm} on real linear functionals $v$ satisfying $v(1)=0$.

The test functions are $A_N^{\rm F},A_N^+,A_N^-,A_N^0$ from \eqref{W:eq:tightness-localized} and \eqref{W:eq:tightness-Fourier-test}, and the function $D_N^{\rm coef}$ from
Section~\ref{W:sec:removing-the-coefficient-cutoff}. For all these functions, a common deterministic $\Lambda_N=N^{o(1)}$ satisfies
\begin{equation}
 |A_N|\le1,\qquad
 |D_s^j A_N(s)[v_1,\ldots,v_j]|
 \le \Lambda_Nm_N^j\prod_{l=1}^j\|v_l\|_\diamond,
 \qquad 1\le j\le4.
 \label{W:eq:replacement-test-class}
\end{equation}
For the first four functions this is Lemma~\ref{W:lem:tightness-global}. For $D_N^{\rm coef}$ it follows by subtracting $\Gamma_\rho\chi_N$ from $\Gamma_\rho$ and applying the
derivative estimate for $\Gamma_\rho$. The constants are uniform over deterministic $w$, and the bounds hold for every real moment vector with $s_0=1$. Their complex multilinear
extensions satisfy the same estimates after a fixed enlargement of $\Lambda_N$.

The bounds in \eqref{W:eq:replacement-test-class} apply to the complete functions $A$ in Lemma~\ref{W:lem:tightness-global}. Its proof combines the derivatives of
$\Gamma_\rho\chi_N$ with those of $\Phi(\mathbf Y)$ by the product and composition formulas of Section~\ref{W:sec:notation-and-scales}. The components
$Y_j^\sigma=-m_N\mathcal T_j^\sigma$ are the exact weighted sums, with the last off-diagonal coefficient having weight one. For the modulus of continuity, the outer derivatives are
bounded in the supremum norm by \eqref{W:eq:tightness-soft}.

The exceptional event estimates require deterministic derivative bounds at every Hermitian matrix. Write $F_{A,N}(H)=A_N(s(H,w))$. For a block $\alpha$, let
$C_{r,\alpha}(K)=D_x^r[F_{A,N}(K+N^{-1/2}\mathcal T_\alpha x)]_{x=0}/r!$ and let $R_{4,\alpha}(K;x)$ be the difference between this function and its Taylor polynomial of degree
three.

\begin{lemma}
\label{W:lem:global-first-two-entry-derivatives}
\label{W:lem:global-cubic-entry-derivatives}
Let $d_\alpha\in\{1,2\}$ be the real coordinate dimension of entry $\alpha$. For $A_N$ satisfying \eqref{W:eq:replacement-test-class}, every Hermitian matrix $K$, every
$x_1,\ldots,x_r\in\mathbb R^{d_\alpha}$, and $1\le r\le4$, one has
\begin{equation}
 \left|D_x^r[F_{A,N}(K+N^{-1/2}\mathcal T_\alpha x)]_{x=0}
 [x_1,\ldots,x_r]\right|
 \le N^{4r/3+o(1)}\prod_{j=1}^r\|x_j\|_1,\qquad 1\le r\le4.
 \label{W:eq:global-cubic-entry-derivatives}
\end{equation}
The bound is uniform in $K$, $\alpha$, and the deterministic unit vector $w$. Each separate term of the third-order chain rule has a polynomial bound of the same kind. In
particular, for sufficiently large $N$,
\begin{equation}
 \|C_{r,\alpha}(K)\|_{\max}\le N^6\quad(1\le r\le3).
 \label{W:eq:global-first-two-entry-derivatives}
\end{equation}
The remainder satisfies
\begin{equation}
 |R_{4,\alpha}(K;x)|\le CN^6\|x\|_1^4.
 \label{W:eq:global-entry-derivative-bound}
\end{equation}
\end{lemma}
\begin{proof}
Put $f_{jk}=(p_j\circ h_d)(p_k\circ h_d)$, with $j,k<\rho$. The product rule and \eqref{W:eq:gram-scalar-bound} give, for $r\le4$,
\[
 \|f_{jk}^{(r)}\|_\infty
 \le\sum_{a=0}^r\binom ra
       \|\varphi_{j,d}^{(a)}\|_\infty
       \|\varphi_{k,d}^{(r-a)}\|_\infty
 \le C_r\rho_*^{2r+2}N^{o(1)}.
\]
These bounds hold on the whole real line, including the transition region of $h_d$. For a label set $S$ of size $r$ and fixed matrices $W_j=N^{-1/2}\mathcal T_\alpha x_j$, the
functional calculus estimate of Lemma~\ref{W:lem:smooth-functional-calculus} bounds each induced entry of the Gram matrix by $C_rN^{-r/2}\rho_*^{2r+2}N^{o(1)}\prod_j\|x_j\|_1$. For
a $\rho\times\rho$ matrix $B$ with singular values $s_j(B)$,
\[
 \rho^{-1}\|B\|_1
  =\rho^{-1}\sum_{j=1}^{\rho}s_j(B)
  \le\rho^{-1/2}\Bigl(\sum_{j=1}^{\rho}s_j(B)^2\Bigr)^{1/2}
  =\rho^{-1/2}\|B\|_{\mathrm{HS}}
  \le\sqrt\rho\max_{j,k}|B_{jk}|.
\]
Apply this estimate to the matrix of the differentiated functional, whose $(j,k)$ entry is the derivative of $s(H,w)(p_jp_k)$. Derivatives of positive order annihilate the constant
polynomial, because $s(H,w)(1)=w^*w=1$. The definition of $\|\cdot\|_\diamond$ and the preceding entry bound therefore give
\begin{equation}
 \|D_{H,S}s[W_j:j\in S]\|_\diamond
 \le C_r\sqrt\rho\,N^{-r/2}\rho_*^{2r+2}N^{o(1)}
       \prod_{j\in S}\|x_j\|_1
 \le N^{(r+5)/6+o(1)}\prod_{j\in S}\|x_j\|_1.
 \label{W:eq:global-moment-entry-growth}
\end{equation}
Here $\rho,\rho_*=O(m_N)$, with the constants fixed by the chosen functions. The labelled composition rule is
\[
 D_{H,S}(A_N\circ s)
   =\sum_{\pi\in\Pi(S)}D_s^{|\pi|}A_N(s(K,w))
                       [D_{H,B}s:B\in\pi].
\]
Here $D_{H,B}$ means differentiation once in each variable $t_j$, $j\in B$, of the map evaluated at $H+\sum_{j\in B}t_jW_j$, followed by $t_j=0$. The labels remain distinct also
when two of the matrices $W_j$ agree. For a partition $\pi$ of $S$ into $k$ blocks, the corresponding term is
\[
 D_s^kA_N(s(K,w))[D_{H,B}s:B\in\pi].
\]
Each derivative $D_{H,B}s$ uses the fixed matrices $N^{-1/2}\mathcal T_\alpha x_j$, $j\in B$. It annihilates constants and can therefore be used as an argument of
\eqref{W:eq:replacement-test-class}. Since $\sum_{B\in\pi}|B|=r$, multiplication of the bounds in \eqref{W:eq:global-moment-entry-growth} gives
\[
 \prod_{B\in\pi}\|D_{H,B}s\|_\diamond
 \le N^{(r+5k)/6+o(1)}\prod_{j\in S}\|x_j\|_1.
\]
The modulus of the corresponding term is thus bounded by
\[
 \Lambda_Nm_N^kN^{(r+5k)/6+o(1)}\prod_{j\in S}\|x_j\|_1
 =N^{(r+7k)/6+o(1)}\prod_{j\in S}\|x_j\|_1
 \le N^{4r/3+o(1)}\prod_{j\in S}\|x_j\|_1.
\]
There are at most $15$ set partitions at orders at most four. Summing these bounds proves the displayed bound for derivatives of orders one through four. The partition classes
$\{123\}$, $\{1\}\{23\}$ and its permutations, and $\{1\}\{2\}\{3\}$ give the three separate cubic terms, so their polynomial bounds follow without cancellation. Each tensor entry
is obtained by evaluating the multilinear derivative on the coordinate vectors $e_1,e_2$, each of $\ell^1$ norm one. Hence the multilinear bounds just proved bound the maximum
tensor norm with no additional factor depending on $N$. In particular, for $r\le3$ they give $N^{4r/3+o(1)}\le N^6$ for all sufficiently large $N$. For the remainder, define
$F_\alpha(u)=A_N(s(K+N^{-1/2}\mathcal T_\alpha u,w))$. Taylor's formula on the segment $u=tx$ gives
\[
 R_{4,\alpha}(K;x)=\frac1{3!}\int_0^1(1-t)^3
             D^4F_\alpha(tx)[x,x,x,x]\,\dd t.
\]
The fourth derivative bound is uniform on this segment. Consequently
\[
 |R_{4,\alpha}(K;x)|
 \le\frac1{24}N^{16/3+o(1)}\|x\|_1^4
 \le CN^6\|x\|_1^4
\]
for sufficiently large $N$, since the exponent denoted by $o(1)$ is then at most $2/3$.
\end{proof}

At a fixed matrix with the selected entry set to zero put $s_\alpha(x)=s(K_\alpha+N^{-1/2}\mathcal T_\alpha x,w)$ and define
\[
 \ell_r=\partial_rs_\alpha(0),\qquad
 q_{rs}=\tfrac12\partial_r\partial_ss_\alpha(0),\qquad
 c_{rst}=\tfrac16\partial_r\partial_s\partial_ts_\alpha(0).
\]
All derivatives in the following identity are evaluated at $s_\alpha(0)$:
\begin{equation}
 \begin{split}
 (C_{3,\alpha})_{rst}={}&D_{s}A[c_{rst}]
 +\tfrac13\big(D_{s}^2A[\ell_r,q_{st}]
             +D_{s}^2A[\ell_s,q_{rt}]+D_{s}^2A[\ell_t,q_{rs}]\big)\\
 &+\tfrac16D_{s}^3A[\ell_r,\ell_s,\ell_t].
 \end{split}
 \label{W:eq:cumulant-chain-tensor}
\end{equation}
This is the third-order chain rule divided by $3!$. Equivalently, with $\ell(x)=\sum_r x_r\ell_r$, $q(x)=\sum_{r,s}x_rx_sq_{rs}$ and $c(x)=\sum_{r,s,t}x_rx_sx_tc_{rst}$, its
homogeneous polynomial is
\begin{equation}
 \langle C_{3,\alpha},x^{\otimes3}\rangle
 =D_{s}A[c(x)]+D_{s}^2A[\ell(x),q(x)]
                  +\tfrac16D_{s}^3A[\ell(x),\ell(x),\ell(x)].
 \label{W:eq:cumulant-chain-polynomial}
\end{equation}
Conditional on $K_\alpha$, integration of this polynomial against the full difference of the entry distributions gives exactly the ordered tensor contraction in
\eqref{W:eq:cumulant-target-tensor}. In dimension two that contraction is
\[
 C_{111}\Delta M_{111}+3C_{112}\Delta M_{112}
       +3C_{122}\Delta M_{122}+C_{222}\Delta M_{222}.
\]
The ordered triples are $111$; $112,121,211$; $122,212,221$; and $222$. Symmetry of the third derivative makes the coefficients equal within each of these four groups, and
multiplication of scalar coordinates makes the corresponding moment components equal. This yields the multiplicities $1,3,3,1$. For $x=(x_1,x_2)$ the four moments are the integrals
of $x_1^3$, $x_1^2x_2$, $x_1x_2^2$ and $x_2^3$ against the joint signed distribution. The multiplicities above come from the ordered coordinate sum.

For the exponential observable, $A=C\exp(\mathrm i\Theta)$ with $C=\Gamma_\rho\chi_N$ and $\Theta=-m_N\sum_l\lambda_l\mathcal T_{K_l}^{\sigma_l}$. For fixed real functionals $v_a$
satisfying $v_a(1)=0$, $D_S$ denotes differentiation of $A(s+\sum_{a\in S}t_av_a)$ at $t=0$. The product rule and \eqref{W:eq:localized-exponential-partitions} give, for
$1\le|S|\le3$,
\begin{equation}
 D_SA=e^{\mathrm i\Theta}
 \sum_{U\sqcup V=S}(D_U C)
       \sum_{\pi\in\Pi(V)}\mathrm i^{|\pi|}\prod_{B\in\pi}D_B\Theta,
 \qquad
 D_UC=\sum_{L\subseteq U}(D_L\Gamma_\rho)(D_{U\setminus L}\chi_N).
 \label{W:eq:fourier-derivatives-through-third-order}
\end{equation}
For $V=\varnothing$ the inner sum is one. The starting vector and cutoff parameters are fixed throughout these derivatives.

\subsubsection{Cubic terms in the comparison with the GOE and GUE}
\label{W:sec:signed-third-moment-summation}

\label{W:sec:the-signed-third-moment-comparison}
\label{W:sec:cumulant-cubic}

For the comparison of bounded entries with the GOE or GUE, the matching first two off-diagonal moment tensors leave the cubic contribution as the main issue. The chain rule
\eqref{W:eq:cumulant-chain-polynomial} separates this contribution into three parts. Those containing $D_s^2A_N$ and $D_s^3A_N$ can be summed absolutely. The part containing
$D_sA_N$ requires cancellation in a signed matrix sum, which is available only after its coefficient and all four resolvents have been evaluated at the same matrix.

The parameters, moment functional, and class of functions are those of \eqref{W:eq:replacement-test-class}. The homogeneous bounded standardized off-diagonal distribution from
Lemma~\ref{W:lem:cubic-defect} is denoted by $\nu_N^{\rm b}$, and the corresponding real or complex standard Gaussian distribution by $\nu^G$. In complex notation an entry is
$\zeta$; it is real when $\beta=1$. Both distributions are centered and satisfy $\mathbb E|\zeta|^2=1$, with $\mathbb E\zeta^2=0$ in the complex case. Their signed difference and
the interpolating probability measure are

\[
 \Delta\nu_N=\nu_N^{\rm b}-\nu^G,
 \qquad \nu_{N,t}=t\nu_N^{\rm b}+(1-t)\nu^G.
\]
All standardized moments of each fixed order are bounded uniformly in $N,t$. The same law is used at every off-diagonal position. Each standardized diagonal entry has distribution
$t\mu_{\rm d}^{\rm b}+(1-t)\mu_{\rm d}^{G}$, where $\mu_{\rm d}^{\rm b}$ is centered with variance one and $\mu_{\rm d}^{G}$ is the corresponding Gaussian diagonal distribution.
Its variance is $2-t$ for $\beta=1$ and one for $\beta=2$. The upper triangular entry vectors are independent. Consequently $H(0)$ belongs to the GOE or GUE and $H(1)$ has the same
distribution as $H_N^{\rm b}$.

For a fixed $t$, construct one matrix $H(t)$ with these independent entries, and define $K_{ab}=H(t)^{[ab\leftarrow0]}$ for all $a<b$ on this probability space. Each $K_{ab}$ has
the marginal distribution required by the interpolation formula. For $\alpha=\{a,b\}$, the matrix corresponding to a replacement coordinate vector is $N^{-1/2}\mathcal T_{ab}x$, as
in \eqref{W:eq:standardized-entry-insertion}. Let $\mu_N^{\rm b}$ and $\mu^G$ be the corresponding real coordinate laws, so that $x=(\sqrt2\Re\zeta,\sqrt2\Im\zeta)$ when $\beta=2$.
Define
\begin{align}
 (C_{3,\alpha}(A;K))_{rst}
 &=\left.\frac1{6}\partial_{x_r}\partial_{x_s}\partial_{x_t}
   A_N\big(s(K+N^{-1/2}\mathcal T_\alpha x,w)\big)\right|_{x=0},
       \nonumber\\
 \Delta M_{rst}
 &=\int x_rx_sx_t\,(\mu_N^{\rm b}-\mu^G)(\mathrm dx),\nonumber\\
 Z_{3,N}(A;t,w)
 &=\sum_{a<b}\sum_{r,s,t'=1}^{\beta}
  (C_{3,ab}(A;K_{ab}))_{rst'}\Delta M_{rst'}.
 \label{W:eq:cumulant-target-tensor}
\end{align}
The coordinate sums run over ordered indices. The letter $t'$ is a tensor index, whereas $t\in[0,1]$ is the interpolation parameter.

For $a<b$, put $T_{ab}(\zeta)=\zeta E_{ab}+\bar\zeta E_{ba}$. Then $T_{ab}(\zeta)^2=|\zeta|^2(E_{aa}+E_{bb})$ and $T_{ab}(\zeta)^3=|\zeta|^2T_{ab}(\zeta)$. Define
\begin{equation}
 \gamma_{3,N}=\int\zeta|\zeta|^2\,\Delta\nu_N(\mathrm d\zeta),\qquad
 B_{3,N}=\sum_{a<b}\big(\gamma_{3,N}E_{ab}
                              +\overline{\gamma_{3,N}}E_{ba}\big).
 \label{W:eq:cumulant-third-matrix}
\end{equation}
The common off-diagonal distribution gives the same $\gamma_{3,N}$ at each pair. For $\beta=2$ the normalization \eqref{W:eq:standardized-entry-insertion} gives the explicit
identity
\begin{equation}
 \gamma_{3,N}=\frac1{2\sqrt2}
  \big[\Delta M_{111}+\Delta M_{122}
             +\mathrm i(\Delta M_{112}+\Delta M_{222})\big].
 \label{W:eq:cumulant-mixed-moments}
\end{equation}
For $\beta=1$ it is $\Delta M_{111}$. The second independent complex third moment, $\int\zeta^3\,\Delta\nu_N(\mathrm d\zeta)$, will occur only with off-diagonal compressed
resolvents; its full contribution is displayed in \eqref{W:eq:resolvent-complete-complex-cubic}.

\Needspace{6\baselineskip}

\begin{theorem}
\label{W:thm:cumulant-complete}
For the interpolating matrices $H(t)$ and test functions defined above, a deterministic $Q_N=N^{o(1)}$ satisfies
\begin{equation}
 \sup_{t\in[0,1]}\sup_{\|w\|=1}\mathbb E|Z_{3,N}(A;t,w)|
 \le Q_NN^{-1/6}+C_DN^{-D}
 \label{W:eq:cumulant-main}
\end{equation}
for every fixed $D>0$. Consequently,
\begin{equation}
 \sup_{\|w\|=1}\left|\int_0^1\mathbb E Z_{3,N}(A;t,w)\,\mathrm dt\right|
 \le N^{-1/6+o(1)}.
 \label{W:eq:cumulant-integrated}
\end{equation}
\end{theorem}

The relation $T_{ab}(\zeta)^3=|\zeta|^2T_{ab}(\zeta)$ identifies the useful part of the cubic resolvent product. Replacing its two interior resolvents, restricted to $\{a,b\}$, by
$g(z)I_2$ leaves $g(z)^2w^*RT_{ab}(\zeta)^3Rw$. At a common matrix, summation and signed integration therefore give $g(z)^2w^*RB_{3,N}Rw$.

The trace norm can be computed explicitly in the real case. If $\boldsymbol1$ is the all ones column, then $B_{3,N}=\gamma_{3,N}(\boldsymbol1\boldsymbol1^{\mathsf T}-I)$. Its
eigenvalues are $\gamma_{3,N}(N-1)$ and $-\gamma_{3,N}$ with multiplicity $N-1$, and hence
\begin{equation}
 \|B_{3,N}\|_1=2|\gamma_{3,N}|(N-1),
 \qquad
 \sum_{i,j}|(B_{3,N})_{ij}|=|\gamma_{3,N}|N(N-1).
 \label{W:eq:real-third-moment-matrix-norm}
\end{equation}
For complex $\gamma_{3,N}$, Lemma~\ref{W:lem:cumulant-trace-norm} proves $\|B_{3,N}\|_1\le C|\gamma_{3,N}|N\log(2N)$ by decomposing each strict triangle into dyadic rectangles.
This gives an order $N^{1+o(1)}$ trace norm in both symmetry classes.

For each pair, the polynomial $P_A(s(K_{ab}),z)$ and the resolvents are initially evaluated at $K_{ab}$. We estimate their changes when $K_{ab}$ is replaced by $H(t)$. The
resulting polynomial $P_A(s(H(t)),z)$ is common to all pairs, so it factors out of the sum containing $B_{3,N}$. Its product with the resolvents is estimated for each fixed $H(t)$.
The common off-diagonal law makes $B_{3,N}$ constant on each strict triangle. For a deterministic Hermitian matrix with unrelated bounded coefficients $\gamma_{ab}$, the general
inequality $\|B\|_1\le\sqrt N\,\|B\|_{\rm HS}$ gives only $O(N^{3/2})$. With the factor ${m_N}N^{-3/2}$ from the derivative coefficient and the matrix entry scaling, this gives the
growing upper bound $O({m_N})$. Lemma~\ref{W:lem:cumulant-trace-norm} and the estimate for the modulus of the resolvent \eqref{W:eq:cumulant-signed-sum} give $N^{1+o(1)}$ for the
homogeneous signed sum, which yields \eqref{W:eq:leading-third-moment-estimate}. This calculation uses the common off-diagonal distribution.

The chain rule \eqref{W:eq:cumulant-chain-polynomial} separates this contribution from $D_{s}^2A[\ell,q]$ and $D_{s}^3A[\ell,\ell,\ell]/6$. The latter are bounded in
\eqref{W:eq:cumulant-nonlinear}, while \eqref{W:eq:cumulant-interior-error} treats the resolvent errors.

For the cubic argument use the spectral window $d=N^{-2/3}(\log N)^\gamma$, $1<\gamma<2$, with any fixed allowed multiplicative constants. Put $s_0=4\max\{\rho^{-1},\sqrt d\}$ and
\[
 \mathcal C=\{z=\zeta+\zeta^{-1}:|\zeta|=e^{s_0}\},\qquad
 E_\rho(z)=\frac{2e^{\rho s_0}}{|z^2-4|^{1/2}},\qquad
 J_N=\int_{\mathcal C}E_\rho(z)^2\frac{|\mathrm dz|}{2\pi}.
\]
The contour is positively oriented. Since $\rho s_0=O((\log N)^{\gamma/2})=o(\log N)$, its weighted length satisfies $J_N=N^{o(1)}$. Indeed, if $\zeta=e^{s_0+\mathrm i\theta}$,
then $|\mathrm dz|/|z^2-4|=\mathrm d\theta/|\zeta-\zeta^{-1}|$, and
\[
 |\zeta-\zeta^{-1}|^2=4(\sinh^2s_0+\sin^2\theta).
\]
Using symmetry of $\sin^2\theta$ and $\sin\theta\ge2\theta/\pi$ on $[0,\pi/2]$ gives
\[
 \int_0^{2\pi}\frac{\mathrm d\theta}{|\zeta-\zeta^{-1}|}
 \le C\int_0^{\pi/2}\frac{\mathrm d\theta}{\sqrt{s_0^2+\theta^2}}
 =C\operatorname{arsinh}\!\left(\frac\pi{2s_0}\right)
 \le C(1+|\log s_0|).
\]
Consequently $J_N\le C e^{2\rho s_0}(1+|\log s_0|)$. Since $\rho s_0=o(\log N)$ and $s_0=N^{-1/3+o(1)}$, $\log J_N=o(\log N)$.

Use the polynomial $P_A$ defined in \eqref{W:eq:tightness-response-polynomial} with the present function $A_N$. The contour and the function $E_\rho$ are those used in
\eqref{W:eq:comparison-measure}. The bounds \eqref{W:eq:tightness-response-bounds} and the representation \eqref{W:eq:tightness-response-representation} apply to every function in
\eqref{W:eq:replacement-test-class}, including $D_N^{\rm coef}$: their proof uses only the first two moment derivative bounds in that class. The representation applies to linear
functionals $v$ with $v(1)=0$, and the difference estimate holds along every segment between two normalized real coordinate vectors. At the matrix $H$, the coefficient
$P_A(s(H),z)$ remains correlated with the resolvent.

\Needspace{6\baselineskip}

Let $\mathcal V_N$ be a deterministic family of at most $N+1$ unit vectors containing $w$, and put $\omega_i=|w_i|+m_N^{-1}$. We apply the next lemma with an eigenbasis of the
deterministic third moment matrix in \eqref{W:eq:cumulant-third-matrix}, together with $w$. The interpolation parameter and the starting vector are fixed when the event is
constructed; all probability estimates are uniform in those choices.

\Needspace{14\baselineskip}
\begin{lemma}
\label{W:lem:cumulant-event}
There are deterministic $q_N=N^{o(1)}$ and events $\mathcal E_{N,t,w}$ such that, for every $D>0$, there is $N_0(D)$ with $\sup_{t,w}\mathbb P(\mathcal E_{N,t,w}^c)\le N^{-D}$ for
$N\ge N_0(D)$, and on which
\begin{align}
 \max_{i,j}|R_H(z)_{ij}-g(z)\delta_{ij}|
   &\le q_Nm_N^{-1},\nonumber\\
 |(R_H(z)w)_i|+|(R_H(\bar z)w)_i|
   &\le q_N\omega_i,\nonumber\\
 v^*|R_H(z)|v&\le q_N\quad(v\in\mathcal V_N),\qquad
 \max_{a,b}|H(t)_{ab}|\le q_NN^{-1/2}
 \label{W:eq:cumulant-event-bounds}
\end{align}
uniformly for $z\in\mathcal C$. Every matrix obtained by multiplying one entry and its conjugate by a number in $[0,1]$ has spectrum in $(-2-d,2+d)$ and satisfies the first two
estimates with a larger $q_N$.
\end{lemma}
The proof is given in Section~\ref{W:sec:resolvent-bounds-for-cubic-comparison}.

On the common event and at every matrix $K_{ab}$, the contour resolvent weights can be enlarged to deterministic upper bounds
\[
 U_{ab}\le q_N'\omega_a\omega_b,\qquad
 V_{ab}\le q_N'(\omega_a^2+\omega_b^2),\qquad q_N'=N^{o(1)}.
\]
This follows by integrating the two resolvent factors containing the starting vector with the weight $E_\rho^2|\mathrm dz|/(2\pi)$; its total mass is $J_N=N^{o(1)}$. With
$\omega_i=b_i$, the coordinate inequalities proved for \eqref{W:eq:moment-weights-weight-table} give $\sum_i\omega_i\le Cm_N^2$, $\sum_i\omega_i^2\le Cm_N$, and
$\sum_i(\omega_i^3+\omega_i^4)\le C$. Consequently, pointwise on the common event,
\begin{align}
 \sum_{a<b}V_{ab}&\le {m_N}^4N^{o(1)},&
 \sum_{a<b}U_{ab}V_{ab}&\le {m_N}^2N^{o(1)},\nonumber\\
 \sum_{a<b}V_{ab}^2&\le {m_N}^3N^{o(1)},&
 \sum_{a<b}U_{ab}^3&\le N^{o(1)}.
 \label{W:eq:cumulant-integrated-weight-sums}
\end{align}
The four sums follow from the nonnegative identities and inequalities
\[
 \begin{aligned}
 \sum_{a<b}(\omega_a^2+\omega_b^2)&=(N-1)\sum_i\omega_i^2,\\
 \sum_{a<b}\omega_a\omega_b(\omega_a^2+\omega_b^2)
   &=\Bigl(\sum_i\omega_i\Bigr)\Bigl(\sum_i\omega_i^3\Bigr)
                                            -\sum_i\omega_i^4,\\
 \sum_{a<b}(\omega_a^2+\omega_b^2)^2
   &=(N-2)\sum_i\omega_i^4+\Bigl(\sum_i\omega_i^2\Bigr)^2,\\
 \sum_{a<b}\omega_a^3\omega_b^3
   &=\tfrac12\left[\Bigl(\sum_i\omega_i^3\Bigr)^2
                                      -\sum_i\omega_i^6\right]
    \le\tfrac12\Bigl(\sum_i\omega_i^3\Bigr)^2.
 \end{aligned}
\]
Substitute the preceding coordinate sums and $N=m_N^3$ into each line.

For $T=\mathcal T_{ab}x$ and $R=R_{K_{ab}}$, the exact homogeneous coefficients of the expansion of the moments on this event are
\begin{align}
 \ell(x)(f)&=\frac{N^{-1/2}}{2\pi\mathrm i}
        \int_{\mathcal C}f(z)w^*RTRw\,\mathrm dz,\nonumber\\
 q(x)(f)&=\frac{N^{-1}}{2\pi\mathrm i}
        \int_{\mathcal C}f(z)w^*RTRTRw\,\mathrm dz,\nonumber\\
 c(x)(f)&=\frac{N^{-3/2}}{2\pi\mathrm i}
        \int_{\mathcal C}f(z)w^*RTRTRTRw\,\mathrm dz.
 \label{W:eq:first-three-moment-coefficients}
\end{align}
The third matrix entry derivative is six times the last resolvent product, so its coefficient of the expansion is the displayed integral. For $f=1$, the moment functional has the
constant value one along the entry replacement, so $\ell(x)(1)=q(x)(1)=c(x)(1)=0$. Restricted coordinate resolvents are bounded, so
\begin{equation}
 \|\ell(x)\|_\diamond\le N^{-1/2}\|x\|_1U_{ab},\qquad
 \|q(x)\|_\diamond\le N^{-1}\|x\|_1^2V_{ab},\qquad
 \|c(x)\|_\diamond\le N^{-3/2}\|x\|_1^3V_{ab},
 \label{W:eq:first-three-moment-coefficient-bounds}
\end{equation}
after enlarging the deterministic constants in the resolvent weights. These homogeneous coefficient bounds hold for every $x$.

Choose the deterministic vectors in Lemma~\ref{W:lem:cumulant-event} to be $w$ and an eigenbasis of the Hermitian matrix $B_{3,N}$ in \eqref{W:eq:cumulant-third-matrix}. They
depend only on the entry laws, $N$, and $w$, and are specified before sampling $H(t)$. Let $\mathcal E$ be that one common event. By
Lemma~\ref{W:lem:global-cubic-entry-derivatives}, each of the three contracted terms in \eqref{W:eq:cumulant-chain-polynomial}, at every matrix with the selected entry set to zero,
is bounded by $N^B\int\|x\|_1^3(\mu_N^{\rm b}+\mu^G)(\mathrm dx)\le CN^B$. Fix $D>0$ and choose the common event with failure probability at most $N^{-D-B-3}$. There are at most
$N^2/2$ pairs and three contracted components, so their total exceptional contribution is bounded by $C N^{B+2}N^{-D-B-3}\le C N^{-D-1}$. In particular,
\begin{equation}
 \mathbb E\left[\mathbf1_{\mathcal E^c}
       \sum_{a<b}|\langle C_{3,ab},\Delta M_3\rangle|\right]
             \le C_DN^{-D}.
 \label{W:eq:cumulant-bad-event}
\end{equation}
Here the bound $CN^B$ holds for every base matrix after integration over $x$ against the deterministic measure $\mu_N^{\rm b}+\mu^G$. Multiplication by $\mathbf1_{\mathcal E^c}$
and expectation therefore give the displayed inequality.

Let $M_3^+=\int\|x\|_1^3(\mu_N^{\rm b}+\mu^G)(\mathrm dx)=O(1)$. Equations \eqref{W:eq:replacement-test-class}, \eqref{W:eq:first-three-moment-coefficient-bounds} and
\eqref{W:eq:cumulant-integrated-weight-sums} give on $\mathcal E$
\begin{align}
 \sum_{a<b}\left|\int D_{s}^2A[\ell(x),q(x)]\,
                    (\mu_N^{\rm b}-\mu^G)(\mathrm dx)\right|
 &\le \Lambda_NM_3^+{m_N}^2N^{-3/2}\sum_{a<b}U_{ab}V_{ab}
       \le N^{-1/6+o(1)},\nonumber\\
 \sum_{a<b}\left|\int\tfrac16D_{s}^3A[\ell(x)^3]\,
                    (\mu_N^{\rm b}-\mu^G)(\mathrm dx)\right|
 &\le C\Lambda_NM_3^+{m_N}^3N^{-3/2}\sum_{a<b}U_{ab}^3
       \le N^{-1/2+o(1)}.
 \label{W:eq:cumulant-nonlinear}
\end{align}
The last estimate uses the stronger weights available for the entry distributions of $H(t)$. Both absolute sums tend to zero. The same argument applied to $D_sA[c]$, however, gives
only
\[
 \Lambda_NM_3^+m_NN^{-3/2}\sum_{a<b}V_{ab}
 \le N^{1/6+o(1)}.
\]
This estimate does not tend to zero. The required improvement comes from summing these first derivative terms with their signs before estimating their absolute value.

\subsubsection{Evaluation of the cubic terms at a common matrix}
\label{W:sec:entry-restoration-and-cubic-contraction}
\label{W:sec:resolvent-cubic-restoration}
\label{W:sec:resolvent-complex-cubic}

The interpolation identity initially evaluates the contribution of each pair at a different matrix, with that entry set to zero. These contributions can be coupled by sampling a
single matrix $H(t)$ and defining $K_{ab}=H(t)^{[ab\leftarrow0]}$ for every $a<b$. Each $K_{ab}$ has the required marginal product distribution, so linearity of expectation
identifies the sum in \eqref{W:eq:entry-distribution-interpolation} with the expectation of the sum on this common probability space.

This construction separates two roles of the entry. The sampled value $H(t)_{ab}$ determines the change from $K_{ab}$ to $H(t)$; the variable $\zeta$ in the cubic coefficient is
integrated against the deterministic signed measure $\Delta\nu_N$. With the matrix $T_{ab}(\zeta)$ defined above, the latter coefficient contains

\begin{equation}
 W_{ab}(H,z)=\int w^*R_HT_{ab}(\zeta)R_HT_{ab}(\zeta)
                      R_HT_{ab}(\zeta)R_Hw\,\Delta\nu_N(\mathrm d\zeta).
 \label{W:eq:cumulant-word}
\end{equation}
By \eqref{W:eq:tightness-response-representation} and \eqref{W:eq:first-three-moment-coefficients}, the contracted cubic term with one derivative of the test function at $K=K_{ab}$
equals
\begin{equation}
 L_{ab}(K)=\frac{N^{-3/2}}{2\pi\mathrm i}
       \int_{\mathcal C}P_A(s(K),z)W_{ab}(K,z)\,\mathrm dz.
 \label{W:eq:cumulant-erased-linear}
\end{equation}
To sum at the common matrix $H(t)$, we estimate the changes in both factors of \eqref{W:eq:cumulant-erased-linear}.

Put $\tau_{ab}=\|H(t)-K_{ab}\|$, $u_{ab}=\omega_a\omega_b$, and $v_{ab}=\omega_a^2+\omega_b^2$. On $\mathcal E$,
$\max_{a<b}\tau_{ab}\le N^{-1/2+o(1)}$. 

\begin{lemma}
\label{W:lem:cubic-entry-restoration}
On $\mathcal E$,
\begin{equation}
 \left|\sum_{a<b}L_{ab}(K_{ab})
 -\frac{N^{-3/2}}{2\pi\mathrm i}
   \int_{\mathcal C}P_A(s(H),z)\sum_{a<b}W_{ab}(H,z)\,\mathrm dz\right|
 \le N^{-1/3+o(1)}.
 \label{W:eq:cubic-comparison-at-one-matrix}
\end{equation}
\end{lemma}
\begin{proof}
Fix $a<b$, $z\in\mathcal C$ and the integration variable $\zeta$. Write $K=K_{ab}$, $V=H-K$, $R_0=(zI-K)^{-1}$, $R_1=(zI-H)^{-1}$, and $T=T_{ab}(\zeta)$. Subtracting the products
one resolvent at a time gives
\begin{align}
 &R_0TR_0TR_0TR_0-R_1TR_1TR_1TR_1\notag\\
 &\quad=(R_0-R_1)TR_1TR_1TR_1
       +R_0T(R_0-R_1)TR_1TR_1\notag\\
 &\qquad\quad+R_0TR_0T(R_0-R_1)TR_1
       +R_0TR_0TR_0T(R_0-R_1)\notag\\
 &\quad=-\sum_{j=0}^3(R_0T)^jR_0VR_1(TR_1)^{3-j}.
 \label{W:eq:resolvent-four-restoration-terms}
\end{align}
The last equality uses $R_0-R_1=-R_0VR_1$. Each summand contains four matrices supported on $J=\{a,b\}$, namely one copy of $V$ and three copies of $T$. Every interior resolvent is
therefore restricted to $J$, and its norm is bounded on $\mathcal E$. The two exterior vectors satisfy
\[
 \|(w^*R_\epsilon)_J\|_2+
 \|(R_\epsilon w)_J\|_2\le Q_N\sqrt{v_{ab}},
 \qquad \epsilon\in\{0,1\},\quad Q_N=N^{o(1)}.
\]
Let $P_J$ restrict to these two coordinates and put $C_*=\max_{\epsilon=0,1}\|P_JR_\epsilon P_J^*\|$. The first and last factors in each scalar product below are
$w^*R_\epsilon P_J^*$ and $P_JR_{\epsilon'}w$. Between them there are three factors $P_JR_\delta P_J^*$, three copies of $T_{JJ}$, and one copy of $V_{JJ}$. Since
$\|T_{JJ}\|=|\zeta|$ and
$\|V_{JJ}\|=\tau_{ab}$,
\[
 |w^*(R_0T)^jR_0VR_1(TR_1)^{3-j}w|
 \le Q_N^2v_{ab}C_*^3\tau_{ab}|\zeta|^3,
 \qquad j=0,1,2,3.
\]
On $\mathcal E$, $C_*$ has a deterministic fixed upper bound. The triangle inequality adds precisely four such bounds. Integrating against $\Delta\nu_N$ and bounding its total
variation by $\nu_N^{\rm b}+\nu^G$ yields
\begin{equation}
 |W_{ab}(H,z)-W_{ab}(K,z)|\le Q_N'\tau_{ab}v_{ab}, \qquad |W_{ab}(K,z)|\le Q_N'v_{ab},\qquad Q_N'=N^{o(1)}.
 \label{W:eq:cumulant-word-difference}
\end{equation}
The second bound uses the preceding estimates for $w^*R_0$ and $R_0w$, with three copies of $T$. The absolute third moments are uniformly bounded. Only the compressions
$P_JR_\epsilon P_J^*$ occur between these matrices, so the full resolvent norm does not enter this estimate.

For the coefficient variation we use the second derivative of $A_N$, including its cutoffs. First estimate the change in its moment argument. The matrix in the linear term is off
diagonal, so its exterior factors involve two different coordinates. With $C_0=(R_0)_{JJ}$ and $V_J=V_{JJ}$, the restricted resolvent identity gives
\[
 R_1-R_0=R_{0,\cdot J}V_JR_{0,J\cdot}
 +R_{0,\cdot J}V_JC_0V_J(I-C_0V_J)^{-1}R_{0,J\cdot}.
\]
Write $V_J=\upsilon E_{12}+\bar\upsilon E_{21}$, so $|\upsilon|=\tau_{ab}$. The first term tested against $w$ equals
\[
 \upsilon(w^*R_0e_a)(e_b^*R_0w)
       +\bar\upsilon(w^*R_0e_b)(e_a^*R_0w).
\]
Its modulus is at most $2Q_N^2\tau_{ab}\omega_a\omega_b$. For large $N$, $\|C_0V_J\|\le1/2$, and the Neumann series gives $\|(I-C_0V_J)^{-1}\|\le2$. The second scalar product is
therefore bounded by $2C Q_N^2\tau_{ab}^2v_{ab}$. Enlarging a deterministic subpolynomial factor gives the two stated weights. For each $f\in\mathcal P_\rho$, integration of this
identity against $f(z)\,\mathrm dz/(2\pi\mathrm i)$ represents $(s(H)-s(K))(f)$. The definition of the contour norm and $\|\cdot\|_\diamond\le\|\cdot\|_{\mathscr D}$ therefore
imply
\begin{equation}
 \|s(H)-s(K)\|_\diamond
 \le Q_N''(\tau_{ab}u_{ab}+\tau_{ab}^2v_{ab}),
 \qquad Q_N''=N^{o(1)}.
 \label{W:eq:cumulant-base-increment}
\end{equation}
Here the spectral modification is the identity along the whole segment, by Lemma~\ref{W:lem:cumulant-event}.

The global extension from Section~\ref{W:sec:orthogonal-polynomials-and-spectral-statistics} is smooth along $s_u=s(K)+u(s(H)-s(K))$. Its global derivative bound therefore gives,
by the fundamental theorem of calculus,
\[
 P_A(s(H),z)-P_A(s(K),z)
 =\int_0^1D_s^2A_N(s_u)[s(H)-s(K),e_z]\,\mathrm du.
\]
The global second derivative bound and $\|e_z\|_\diamond\le2E_\rho(z)^2$ apply to every $s_u$. Subtract the products of the coefficient and resolvent integrands exactly as
\begin{align*}
 P_A(s(K),z)W_{ab}(K,z)-P_A(s(H),z)W_{ab}(H,z)
 ={}&[P_A(s(K),z)-P_A(s(H),z)]W_{ab}(K,z)\\
 &+P_A(s(H),z)[W_{ab}(K,z)-W_{ab}(H,z)].
\end{align*}
Equations~\eqref{W:eq:tightness-response-bounds}, \eqref{W:eq:cumulant-word-difference}, and \eqref{W:eq:cumulant-base-increment} bound its absolute value by
\[
 \widehat Q_NE_\rho(z)^2
 \bigl[m_N\tau_{ab}v_{ab}
       +m_N^2\tau_{ab}u_{ab}v_{ab}
       +m_N^2\tau_{ab}^2v_{ab}^2\bigr].
\]
The pair identities used in \eqref{W:eq:cumulant-integrated-weight-sums} give
\[
 \sum_{a<b}v_{ab}\le Cm_N^4,\qquad
 \sum_{a<b}u_{ab}v_{ab}\le Cm_N^2,\qquad
 \sum_{a<b}v_{ab}^2\le Cm_N^3.
\]
Multiplication by $N^{-3/2}$, integration on $\mathcal C$, and $J_N=N^{o(1)}$ now give the three bounds
\begin{align}
 N^{-3/2}m_NN^{-1/2}m_N^4&=m_N^{-1},\nonumber\\
 N^{-3/2}m_N^2N^{-1/2}m_N^2&=m_N^{-2},\nonumber\\
 N^{-3/2}m_N^2N^{-1}m_N^3&=m_N^{-5/2},
 \label{W:eq:errors-from-restoring-the-entry}
\end{align}
up to a common subpolynomial factor. Their sum is $N^{-1/3+o(1)}$, proving~\eqref{W:eq:cubic-comparison-at-one-matrix}.
\end{proof}

Let $R=R_H(z)$ and write $R_{JJ}=g(z)I_J+E_J$. On $\mathcal E$, $\|E_J\|\le {m_N}^{-1}N^{o(1)}$. For $T=T_{ab}(\zeta)_{JJ}$, the exact interior expansion is
\begin{equation}
 T(gI+E_J)T(gI+E_J)T
 =g^2T^3+gTE_JT^2+gT^2E_JT+TE_JTE_JT.
 \label{W:eq:cumulant-interior-identity}
\end{equation}
Put $\varepsilon_J=\|E_J\|$. Since $\|T\|=|\zeta|$, $\|T^2\|=|\zeta|^2$, and $|g|\le C$, the three terms obey
\[
 \|gTE_JT^2\|\le |g|\,|\zeta|^3\varepsilon_J,\qquad
 \|gT^2E_JT\|\le |g|\,|\zeta|^3\varepsilon_J,\qquad
 \|TE_JTE_JT\|\le |\zeta|^3\varepsilon_J^2.
\]
For large $N$, $\varepsilon_J\le m_N^{-1}N^{o(1)}<1$. Their sum is therefore at most $C|\zeta|^3m_N^{-1}N^{o(1)}$. The squared norms of the exterior vectors are bounded by
$N^{o(1)}(\omega_a^2+\omega_b^2)$. After integration of $|\zeta|^3$ against $\nu_N^{\rm b}+\nu^G$, whose third absolute moment is uniformly bounded, we obtain
\begin{align}
 W_{ab}(H,z)
 &=g(z)^2w^*R\left(\int T_{ab}(\zeta)^3\,
                         \Delta\nu_N(\mathrm d\zeta)\right)Rw
                    +\mathcal R_{ab}(H,z),\nonumber\\
 |\mathcal R_{ab}(H,z)|
 &\le {m_N}^{-1}N^{o(1)}(\omega_a^2+\omega_b^2),\qquad
 \sum_{a<b}|\mathcal R_{ab}(H,z)|\le N^{1+o(1)}.
 \label{W:eq:cumulant-interior-error}
\end{align}
The last sum uses $N\sum_i\omega_i^2\le C{m_N}^4$ and ${m_N}^{-1}{m_N}^4=N$. Multiplication by $N^{-3/2}$ and the bound for the derivative coefficient ${m_N}E_\rho(z)^2N^{o(1)}$,
followed by contour integration, bounds all interior errors by $N^{-1/6+o(1)}$.

We calculate the complex third moment contraction explicitly. Let
\[
 T(\zeta)=\begin{pmatrix}0&\zeta\\\bar\zeta&0\end{pmatrix},\qquad
 C=\begin{pmatrix}c_{11}&c_{12}\\c_{21}&c_{22}\end{pmatrix},\qquad
 \gamma=\int\zeta^2\bar\zeta\,\Delta\nu_N(\mathrm d\zeta),\quad
 \kappa_3=\int\zeta^3\,\Delta\nu_N(\mathrm d\zeta).
\]
The matrix $C$ is the resolvent compression at a nonreal parameter. The first two products are
\[
 TC=\begin{pmatrix}\zeta c_{21}&\zeta c_{22}\\
             \bar\zeta c_{11}&\bar\zeta c_{12}\end{pmatrix},\qquad
 TCT=\begin{pmatrix}|\zeta|^2c_{22}&\zeta^2c_{21}\\
             \bar\zeta^{\,2}c_{12}&|\zeta|^2c_{11}\end{pmatrix}.
\]
Multiplication by $CT$ then yields
\[
 TCTCT=
 \begin{pmatrix}
 c_{22}(\zeta\bar\zeta^{\,2}c_{12}+\zeta^2\bar\zeta c_{21})
 &\zeta^2\bar\zeta c_{11}c_{22}+\zeta^3c_{21}^2\\
 \zeta\bar\zeta^{\,2}c_{11}c_{22}+\bar\zeta^{\,3}c_{12}^2
 &c_{11}(\zeta\bar\zeta^{\,2}c_{12}+\zeta^2\bar\zeta c_{21})
 \end{pmatrix}.
\]
Integrating its four entries gives
\begin{equation}
 \int TCTCT\,\Delta\nu_N(\mathrm d\zeta)
 =\begin{pmatrix}
 c_{22}(\bar\gamma c_{12}+\gamma c_{21})
       &\gamma c_{11}c_{22}+\kappa_3c_{21}^2\\
 \bar\gamma c_{11}c_{22}+\bar\kappa_3c_{12}^2
       &c_{11}(\bar\gamma c_{12}+\gamma c_{21})
 \end{pmatrix}.
 \label{W:eq:resolvent-complete-complex-cubic}
\end{equation}
With the standardized real coordinates used in \eqref{W:eq:resolvent-entry-tensor},
\begin{align}
 \gamma&=\frac{\Delta M_{111}+\Delta M_{122}
                 +\mathrm i(\Delta M_{112}+\Delta M_{222})}{2\sqrt2},
                 \notag\\
 \kappa_3&=\frac{\Delta M_{111}-3\Delta M_{122}
                 +\mathrm i(3\Delta M_{112}-\Delta M_{222})}{2\sqrt2}.
 \label{W:eq:resolvent-two-complex-moments}
\end{align}
The index contractions in this calculation can be listed without any assumption on the adjoint of $C$. For $i,j\in\{1,2\}$ put $\bar i=3-i$ and $\bar j=3-j$. Because $T$ has zero
diagonal,
\[
 (TCTCT)_{ij}
   =\sum_{u=1}^2 T_{i,\bar i}C_{\bar i,u}
                         T_{u,3-u}C_{3-u,\bar j}T_{\bar j,j}.
\]
There are exactly two summands for each matrix entry, eight in all. For $(i,j)=(1,2)$ they are $\zeta\bar\zeta\zeta\,c_{22}c_{11}$ and $\zeta^3c_{21}^2$, as drawn below. Their
signed integrals are $\gamma c_{11}c_{22}$ and $\kappa_3c_{21}^2$. The other six summands are the remaining entries of \eqref{W:eq:resolvent-complete-complex-cubic}. For nonreal
$z$, the relevant adjoint identity is $C(z)^*=C(\bar z)$. The entries $c_{12}(z)$ and $c_{21}(z)$ therefore occur separately in the displayed calculation. \begingroup
\setlength{\intextsep}{1.2\baselineskip}
\begin{figure}[H]
\centering
\begin{tikzpicture}[x=1cm,y=1cm,every node/.style={mg text}]
 \foreach \k/\lab in {0/1,1/2,2/2,3/1,4/1,5/2}{
   \node[mg site,label=below:{$\lab$}] (a\k) at (1.55*\k,0) {};
 }
 \foreach \k/\lab in {0/1,1/2,2/1,3/2,4/1,5/2}{
   \node[mg site,label=below:{$\lab$}] (b\k) at (1.55*\k,-1.6) {};
 }
 \foreach \x/\y/\lab in {a0/a1/\zeta,a2/a3/\bar\zeta,a4/a5/\zeta,
                         b0/b1/\zeta,b2/b3/\zeta,b4/b5/\zeta}{
  \draw[mg insertion,-{Stealth[length=3.6pt,width=3pt]}]
      (\x)--node[above] {$\lab$}(\y);
 }
 \foreach \x/\y/\lab in {a1/a2/c_{22},a3/a4/c_{11},b1/b2/c_{21},b3/b4/c_{21}}{
  \draw[mg incidence,-{Stealth[length=3.6pt,width=3pt]}]
      (\x)--node[above] {$\lab$}(\y);
 }
 \node[anchor=west] at (8.4,0) {$\longmapsto\ \gamma c_{11}c_{22}$};
 \node[anchor=west] at (8.4,-1.6) {$\longmapsto\ \kappa_3 c_{21}^2$};
 \node[align=center] at (5.5,-2.55)
     {solid arrows: entries of $T$; dashed arrows: entries of $C$};
\end{tikzpicture}
\caption{Two terms of $(TCTCT)_{12}$ in \eqref{W:eq:resolvent-complete-complex-cubic}. Reading from left to right, solid arrows carry entries of $T$ and dashed arrows carry entries of $C$; circles specify the row and column indices in $\{1,2\}$. Integration against $\Delta\nu_N$ gives the expressions on the right. Choosing the two diagonal entries of $C$ gives $\gamma c_{11}c_{22}$, whereas the term with $\kappa_3$ contains $c_{21}^2$. With $C=gI+E$, these become $g^2\gamma$ plus terms containing $E$, and $\kappa_3E_{21}^2$, respectively.}
\label{W:fig:complex-cubic-contractions}
\end{figure}
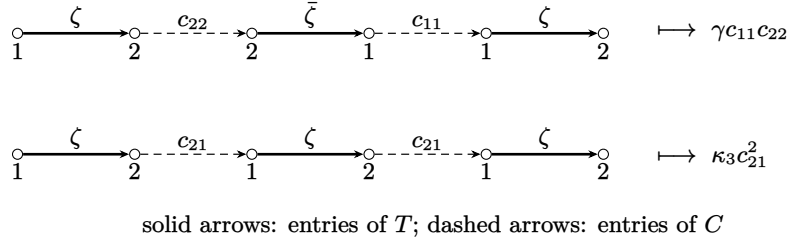
\endgroup The contributions proportional to $\kappa_3$ in \eqref{W:eq:resolvent-complete-complex-cubic} therefore
contain two off-diagonal compressed resolvent errors. If $C=gI+E$ and $\|E\|\le\varepsilon$, then
\begin{equation}
 \left\|\int TCTCT\,\Delta\nu_N
          -g^2\begin{pmatrix}0&\gamma\\\bar\gamma&0\end{pmatrix}\right\|
 \le C M_3^+\bigl(|g|\varepsilon+\varepsilon^2\bigr),\qquad
 M_3^+=\int|\zeta|^3(\nu_N^{\rm b}+\nu^G)(\mathrm d\zeta).
 \label{W:eq:resolvent-complex-cubic-error}
\end{equation}
Indeed $|c_{12}|,|c_{21}|\le\varepsilon$ and $|c_{ii}|\le|g|+\varepsilon$, while
\[
 |c_{11}c_{22}-g^2|
 =|g(E_{11}+E_{22})+E_{11}E_{22}|
 \le2|g|\varepsilon+\varepsilon^2.
\]
Each diagonal entry of the difference is at most $2M_3^+(|g|+\varepsilon)\varepsilon$ in modulus. Each off-diagonal entry is at most $M_3^+(2|g|\varepsilon+2\varepsilon^2)$. For a
$2\times2$ matrix, the operator norm is at most twice its largest entry modulus. These inequalities prove \eqref{W:eq:resolvent-complex-cubic-error}. On the common event
$\varepsilon\le Q_Nm_N^{-1}$, so this calculation gives exactly the error already bounded in \eqref{W:eq:cumulant-interior-error}, with the same two exterior factors.

The part proportional to $\kappa_3$ has the stronger bound obtained from the product of its two factors from $E$. With $\varepsilon\le Q_Nm_N^{-1}$, the sum of its absolute
contributions, including the exterior vectors, the first derivative, and the contour integral, is at most
\[
 C\Lambda_NQ_N^4 M_3^+J_N\,m_NN^{-3/2}m_N^{-2}
       \sum_{a<b}(\omega_a^2+\omega_b^2)
 \le N^{o(1)}m_N^3N^{-3/2}=N^{-1/2+o(1)}.
\]
Here $Q_N$ is increased, if necessary, to dominate the existing subpolynomial bounds. The remaining contribution at each pair is
$g(z)^2w^*R(\int T_{ab}(\zeta)^3\,\Delta\nu_N(\mathrm d\zeta))Rw$. Since every factor outside that integral is evaluated at the same $H(t)$, finite summation gives
\[
 \sum_{a<b}g(z)^2w^*R\!\left(\int T_{ab}(\zeta)^3\,
                   \Delta\nu_N(\mathrm d\zeta)\right)Rw
       =g(z)^2w^*RB_{3,N}Rw.
\]
The coefficient of the first derivative $P_A(s(H(t)),z)$ is also common to the pairs. It can be multiplied by this signed sum before an absolute value is taken. The next estimate
uses the trace norm of the deterministic matrix $B_{3,N}$.

\subsubsection{Trace norm bound for the matrix of third moment differences}
\label{W:sec:trace-norm-third-moment-matrix}
It remains to bound $w^*RB_{3,N}Rw$. The deterministic structure of $B_{3,N}$ gives a trace norm of order $N\log(2N)$, and the quadratic-form estimates for $|R|$ then control this
product.

\begin{lemma}
\label{W:lem:cumulant-trace-norm}
The deterministic Hermitian matrix in \eqref{W:eq:cumulant-third-matrix} obeys
\begin{equation}
 \|B_{3,N}\|_1\le C|\gamma_{3,N}|N\log(2N).
 \label{W:eq:cumulant-trace-norm}
\end{equation}
For real $\gamma_{3,N}$ its trace norm is exactly $2|\gamma_{3,N}|(N-1)$.
\end{lemma}
\begin{proof}
Let $U_N$ have ones strictly above the diagonal and zeros elsewhere. First suppose $N=2^k$. Partition $\{1,\ldots,N\}$ into dyadic intervals. For each interval $I$ with left and
right halves $I_L,I_R$, put $R_I=\mathbf1_{I_L}\mathbf1_{I_R}^{\mathsf T}$. For each $a<b$ there is exactly one dyadic interval whose left and right halves contain $a$ and $b$,
respectively. Hence
\[
 \sum_I\mathbf1_{I_L}(a)\mathbf1_{I_R}(b)
       =\mathbf1_{\{a<b\}},\qquad U_N=\sum_I R_I.
\]
Each $R_I$ has rank one and $\|R_I\|_1=\sqrt{|I_L||I_R|}=|I|/2$. After $\ell$ bisections, $0\le\ell<k$, the dyadic partition has $2^\ell$ intervals, each of length $N2^{-\ell}$.
Hence
\[
 \|U_N\|_1\le\sum_{\ell=0}^{k-1}\sum_{|I|=N2^{-\ell}}
       \|R_I\|_1
   =\sum_{\ell=0}^{k-1}2^\ell\frac{N2^{-\ell}}2
   =\frac{Nk}{2}.
\]
For general $N$, let $M=2^{\lceil\log_2N\rceil}<2N$ and let $J:\mathbb C^N\to\mathbb C^M$ be the coordinate inclusion. Then $U_N=J^*U_MJ$. The ideal inequality for the trace norm
gives $\|J^*U_MJ\|_1\le\|J^*\|\|U_M\|_1\|J\|$; see Theorem~2.7(a) of \cite{simon2005trace}. Consequently
\[
 \|U_N\|_1\le\|J^*\|\|U_M\|_1\|J\|
 \le\tfrac12M\log_2 M\le C N\log(2N).
\] Now
$B_{3,N}=\gamma_{3,N}U_N+\overline{\gamma_{3,N}}U_N^*$, so the triangle inequality proves \eqref{W:eq:cumulant-trace-norm}. If $\gamma_{3,N}$ is real, then
$B_{3,N}=\gamma_{3,N}(\mathbf1\mathbf1^{\mathsf T}-I)$. The vector $N^{-1/2}\mathbf1$ is an eigenvector with eigenvalue $(N-1)\gamma_{3,N}$. Every vector $v$ with
$\mathbf1^{\mathsf T}v=0$ has eigenvalue $-\gamma_{3,N}$. Since the matrix is Hermitian, its singular values are the absolute values of these eigenvalues, and their sum is
$2|\gamma_{3,N}|(N-1)$. For $N=1$, both $U_N$ and $B_{3,N}$ vanish and both identities hold with value zero.\qedhere \begingroup \setlength{\intextsep}{1.2\baselineskip}
\begin{figure}[H]
\centering
\begin{tikzpicture}[x=.47cm,y=.47cm,every node/.style={mg text}]
 \fill[black!10] (4,4) rectangle (8,8);
 \fill[black!25] (2,6) rectangle (4,8);
 \fill[black!25] (6,2) rectangle (8,4);
 \foreach \x/\y in {1/7,3/5,5/3,7/1}{\fill[black!45] (\x,\y) rectangle +(1,1);}
 \draw[step=1,black!45,line width=.3pt] (0,0) grid (8,8);
 \draw[line width=.65pt] (0,0) rectangle (8,8);
 \foreach \j in {1,...,8}{
   \node at (\j-.5,8.6) {$\j$};
   \node at (-.65,8.5-\j) {$\j$};
 }
 \node at (4,9.5) {column $b$};
 \node[rotate=90] at (-1.6,4) {row $a$};
 \node[anchor=west,align=left] at (10,7.3)
   {$|I|=8$: one $4\times4$ rectangle\\$\displaystyle\sum\|R_I\|_1=4$};
 \node[anchor=west,align=left] at (10,4.7)
   {$|I|=4$: two $2\times2$ rectangles\\$\displaystyle\sum\|R_I\|_1=2\cdot2=4$};
 \node[anchor=west,align=left] at (10,2.1)
   {$|I|=2$: four $1\times1$ rectangles\\$\displaystyle\sum\|R_I\|_1=4\cdot1=4$};
 \node[anchor=west] at (10,-.1) {$\|U_8\|_1\le4+4+4=12$};
\end{tikzpicture}
\caption{A partition of the strict upper triangle of $U_8$ into dyadic rectangles. The three shades correspond to intervals $I$ of lengths $8$, $4$, and $2$. Each rectangle represents $R_I=\mathbf1_{I_L}\mathbf1_{I_R}^{\mathsf T}$ and has trace norm $\sqrt{|I_L||I_R|}=|I|/2$. For each of the three interval lengths the sum of these trace norms is $4$, so $\|U_8\|_1\le12$. Summing the same calculation over the $\log_2N$ interval lengths proves the bound in \eqref{W:eq:cumulant-trace-norm}.}
\label{W:fig:third-moment-dyadic-rectangles}
\end{figure}
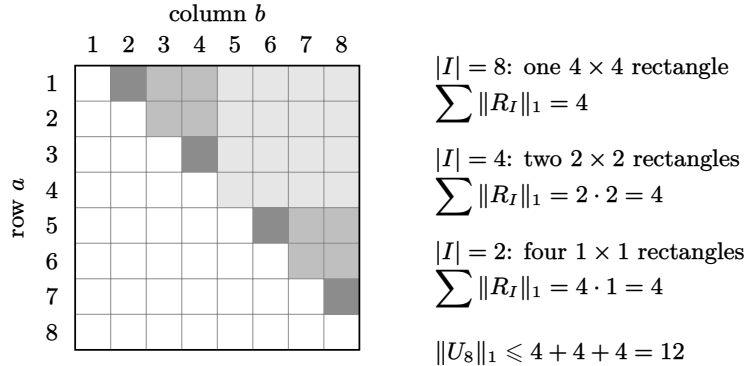
\endgroup
\end{proof}

Normality of $R$ gives the direct bound $|w^*RB_{3,N}Rw|\le\|B_{3,N}\|\,\|Rw\|^2$. The available estimates $\|B_{3,N}\|=O(N)$ and $\|Rw\|^2\le m_NN^{o(1)}$ give, after including
the derivative and entry factors, only
\[
 m_NN^{-3/2}\cdot N\cdot m_N=N^{1/6}.
\]
For a sharper estimate, write $B_{3,N}=\sum_l\lambda_lv_lv_l^*$ in a deterministic orthonormal eigenbasis and use these vectors, together with $w$, in
Lemma~\ref{W:lem:cumulant-event}. If $B_{3,N}=0$, any fixed basis suffices. The basis is fixed by $N$ and the entry distributions before $H(t)$ is sampled.

If $H=\sum_j\lambda_jP_j$, set
\[
 |R|=\sum_j|z-\lambda_j|^{-1}P_j,\qquad
 U=\sum_j\frac{|z-\lambda_j|}{z-\lambda_j}P_j.
\]
Then $U$ is unitary, commutes with $|R|$, and $R=|R|^{1/2}U|R|^{1/2}$. Cauchy--Schwarz applied to $(|R|^{1/2}u)^*U(|R|^{1/2}v)$ gives
\begin{equation}
 |u^*Rv|\le(u^*|R|u)^{1/2}(v^*|R|v)^{1/2}.
 \label{W:eq:cumulant-polar-bound}
\end{equation}
Apply this separately to $w^*Rv_l$ and $v_l^*Rw$. On the common event,
\begin{equation}
 |w^*RB_{3,N}Rw| \le(w^*|R|w)\sum_l|\lambda_l|\,v_l^*|R|v_l \le q_N^2\|B_{3,N}\|_1\le N^{1+o(1)}.
 \label{W:eq:cumulant-signed-sum}
\end{equation}
Both factors $w^*Rv_l$ and $v_l^*Rw$ are estimated by \eqref{W:eq:cumulant-polar-bound}; the eigenbasis sum is weighted by $|\lambda_l|$.

The leading sum in \eqref{W:eq:cumulant-interior-error} is $g(z)^2w^*RB_{3,N}Rw$. It follows from \eqref{W:eq:tightness-response-bounds} and \eqref{W:eq:cumulant-signed-sum} that
\begin{equation}
 \begin{split}
 &\left|\frac{N^{-3/2}}{2\pi\mathrm i}
     \int_{\mathcal C}P_A(s(H),z)g(z)^2w^*RB_{3,N}Rw\,\mathrm dz\right|\\
 &\hspace{12mm}\le
 C N^{-3/2}\Lambda_N{m_N} q_N^2\|B_{3,N}\|_1 J_N
       \le N^{-1/6+o(1)}.
 \end{split}
 \label{W:eq:leading-third-moment-estimate}
\end{equation}
Here $|g|$ is bounded on the ellipse and $|\gamma_{3,N}|=O(1)$ by the third absolute moment. Both $P_A(s(H),z)$ and $R_H(z)$ are functions of $H$. The displayed inequality holds
pointwise on $\mathcal E$ for their product, and its deterministic upper bound also bounds the absolute value of the expectation of that product multiplied by
$\mathbf1_{\mathcal E}$.

\begin{proof}[Proof of Theorem~\ref{W:thm:cumulant-complete}]
Set $K_{ab}=H(t)^{[ab\leftarrow0]}$ for every $a<b$ and write $Z_{3,N}=Z_{3,N}^{(1)}+Z_{3,N}^{(2)}+Z_{3,N}^{(3)}$ according to \eqref{W:eq:cumulant-chain-polynomial}. The
decomposition is obtained by integrating each of the three homogeneous polynomials in \eqref{W:eq:cumulant-chain-polynomial} against the full signed entry law and then summing over
$a<b$. Equation~\eqref{W:eq:cumulant-nonlinear} gives on $\mathcal E$
\[
 |Z_{3,N}^{(2)}|\le Q_Nm_N^4N^{-3/2},\qquad
 |Z_{3,N}^{(3)}|\le Q_Nm_N^3N^{-3/2}.
\]
Here and below $Q_N$ can be increased by the sum or product of finitely many deterministic subpolynomial factors. For the first, \eqref{W:eq:cubic-comparison-at-one-matrix} gives
on $\mathcal E$
\[
 Z_{3,N}^{(1)}
 =\frac{N^{-3/2}}{2\pi\mathrm i}\int_{\mathcal C}
       P_A(s(H),z)\sum_{a<b}W_{ab}(H,z)\,\mathrm dz
       +O(N^{-1/3+o(1)}).
\]
Put $E(z)=\sum_{a<b}W_{ab}(H,z)-g(z)^2w^*R_H(z)B_{3,N}R_H(z)w$. The bound \eqref{W:eq:cumulant-interior-error} and the integral estimate for $P_A$ in
\eqref{W:eq:tightness-response-bounds} give
\[
 \frac{N^{-3/2}}{2\pi}
 \int_{\mathcal C}|P_A(s(H),z)|\,|E(z)|\,|\dd z|
 \le N^{-3/2}\,m_NN^{1+o(1)}J_N
 =N^{-1/6+o(1)},
\]
where $J_N=\int_{\mathcal C}E_\rho(z)^2|\dd z|/(2\pi)=N^{o(1)}$ and $m_N=N^{1/3}$. The integral with $E$ replaced by $g^2w^*R_HB_{3,N}R_Hw$ is bounded by
\eqref{W:eq:leading-third-moment-estimate}. Adding these two bounds, the $N^{-1/3+o(1)}$ error, and the two nonlinear bounds yields
$\mathbf1_{\mathcal E}|Z_{3,N}|\le N^{-1/6+o(1)}$.

The preceding estimates have the following powers of $N$, with a deterministic factor $N^{o(1)}$ in each row:
\begin{equation}
 \begin{array}{c|c}
 \text{contribution}&\text{bound}\\ \hline
 D_{s}^2A[\ell,q]&{m_N}^4N^{-3/2}=N^{-1/6}\\
 D_{s}^3A[\ell^3]/6&{m_N}^3N^{-3/2}=N^{-1/2}\\
 \text{difference of resolvent products}&{m_N}^5N^{-2}=N^{-1/3}\\
 \text{coefficient difference: linear increment}&{m_N}^4N^{-2}=N^{-2/3}\\
 \text{coefficient difference: quadratic increment}&{m_N}^5N^{-5/2}=N^{-5/6}\\
 \text{interior errors}&{m_N}N^{-1/2}=N^{-1/6}\\
 \text{signed leading sum}&{m_N}N^{-1/2}=N^{-1/6}.
 \end{array}
 \label{W:eq:third-moment-error-sums}
\end{equation}
The last line includes $\log(2N)$ in the subpolynomial factor. The third absolute moments are uniformly bounded. For each pair, the global bound for the ordered cubic tensor is
$N^6$, and its contraction has modulus at most $N^6\|\Delta M_3\|_{1,\mathrm{ord}}\le CN^6$. There are at most $N^2/2$ pairs. Choose the common event with
$\sup_{t,w}\mathbb P(\mathcal E^c)\le N^{-D-9}$. It follows that
\[
 \sup_{t,w}\mathbb E\!\left[\mathbf1_{\mathcal E^c}|Z_{3,N}|\right]
 \le CN^8\sup_{t,w}\mathbb P(\mathcal E^c)
 \le CN^{-D-1}.
\]
Combining the two disjoint contributions gives
\[
 \mathbb E|Z_{3,N}|=
 \mathbb E[\mathbf1_{\mathcal E}|Z_{3,N}|]
 +\mathbb E[\mathbf1_{\mathcal E^c}|Z_{3,N}|]
 \le Q_NN^{-1/6}+C_DN^{-D-1}.
\]
The constants in both terms are uniform in fixed $t$ and deterministic $w$. Taking their suprema proves \eqref{W:eq:cumulant-main}. Finally,
\[
 \left|\int_0^1\mathbb E Z_{3,N}\,\mathrm dt\right|
 \le\int_0^1\mathbb E|Z_{3,N}|\,\mathrm dt
 \le N^{-1/6+o(1)},
\]
uniformly over deterministic unit starting vectors. This is \eqref{W:eq:cumulant-integrated}. To justify the integral, expand the distribution of $H(t)$ as the finite sum of
product measures used in Lemma~\ref{W:lem:entry-distribution-interpolation}. Each coefficient is a polynomial in $t$, and the integrands have the deterministic bound $CN^8$. Their
expectations are measurable and integrable on $[0,1]$; the displayed inequality is therefore the triangle inequality for an ordinary Lebesgue integral. The signed sum is formed at
each fixed $H(t)$ before taking its absolute value.

\end{proof}

\subsubsection{Resolvent estimates for the interpolating matrices}
\label{W:sec:resolvent-bounds-for-cubic-comparison}

The cubic calculation used resolvent bounds for $H(t)$ with the entry distributions defined above. We now justify those bounds for precisely that family. Every fixed standardized
moment is uniformly bounded there, which permits the isotropic local law. Its normalization has to be adjusted to the variance row sums of $H(t)$, and its nonreal estimates have to
be extended to the real spectral parameters outside the spectrum. The original truncated distributions are treated instead by the conditional estimates of
Section~\ref{W:sec:truncation-and-resolvent-estimates}.

\begin{proof}[Proof of Lemma~\ref{W:lem:cumulant-event}]
For $\beta=1$, the diagonal variance of $H(t)$ is $(2-t)/N$, so
\[
 \sum_j\mathbb E|H(t)_{ij}|^2=1+(1-t)/N=:a_t^2.
\]
For $\beta=2$, set $a_t=1$. Then $\widetilde H=H(t)/a_t$ has variance row sums one, variances between $c/N$ and $C/N$, and uniformly bounded standardized moments. These bounds hold
uniformly in $t$, since each entry distribution is a convex combination of the fixed bounded and Gaussian distributions. Theorem~2.12 of \cite{bloemendal2014isotropic} gives
\[
 |u^*R_{\widetilde H}(E+\ii\eta)v-g(E+\ii\eta)u^*v|
 \prec\sqrt{\frac{|\Im g(E+\ii\eta)|}{N\eta}}+\frac1{N\eta}
\]
for deterministic unit vectors, bounded $E$, and $N^{-1+\omega}\le\eta\le\omega^{-1}$. With $g_t(z)=a_t^{-1}g(z/a_t)$ and $R_H(z)=a_t^{-1}R_{\widetilde H}(z/a_t)$, the comparison
function for $H$ is $g_t$. On the domains below the square-root formula for $g$ gives
\[
 \partial_a\{a^{-1}g(z/a)\}=-a^{-2}g(z/a)-za^{-3}g'(z/a),
 \qquad |g_t(z)-g(z)|\le CN^{-1}s_0^{-1}=m_N^{-2}N^{o(1)}.
\]
Here $|a_t-1|=O(N^{-1})=o(s_0^2)$, $|g|\le C$ and $|g'|\le C/s_0$ on the dilation segment.

Apply the upper spectral deviation estimate of \cite{erdosxu2023small} to $H(t)$ and $-H(t)$. Its moment constants and diagonal variance parameters in $[1,2]$ are uniform in $t$.
Choose $x_D=c_D(\log N)^{2/3}$ with $(2\beta/3)c_D^{3/2}>D+2$. Then
\[
 \mathbb P\bigl(\lambda_{\max}(\pm H(t))>2+N^{-2/3}x_D\bigr)
 \le Cx_D^{-3\beta/4}e^{-(2\beta/3)x_D^{3/2}}
 \le CN^{-D-2}.
\]
Since $x_D/(N^{2/3}d/8)\to0$, with failure probability smaller than any prescribed inverse power of the matrix size,
\begin{equation}
 \sigma(H(t))\subset[-2-d/8,2+d/8].
 \label{W:eq:cumulant-confinement}
\end{equation}

On the upper ellipse, $z=2\cosh s_0\cos\theta+2\ii\sinh s_0\sin\theta$ and $g(z)=e^{-s_0-\ii\theta}$. Thus $|\Im g(z)|/\Im z=e^{-s_0}/(2\sinh s_0)$. Outside the two arcs
$|\theta|\le s_0/4$ and $|\theta-\pi|\le s_0/4$, we have $\Im z\ge cs_0^2$ and
\[
 \sqrt{\frac{|\Im g(z)|}{N\Im z}}+\frac1{N\Im z}
 \le m_N^{-1}N^{o(1)}.
\]
The local law therefore gives the required entry estimates and $|(R_H(z)w)_i|\le N^{o(1)}\omega_i$. Apply it to the ordered pairs $(e_i,e_j),(e_i,w),(w,e_i)$ and to $(v,v)$ for
$v\in\mathcal V_N$. Adjunction gives the lower half estimates.

On the right arc, $x=\Re z\ge2+4s_0^2/5$ and $|\Im z|\le s_0^2$. On the event \eqref{W:eq:cumulant-confinement}, every eigenvalue satisfies $\lambda\le2+d/8$. Since $s_0^2\ge16d$,
subtraction gives
\[
 x-\lambda\ge\frac45s_0^2-\frac d8
       \ge\left(\frac45-\frac1{128}\right)s_0^2
       \ge\frac12s_0^2.
\] At
$z_h=x+2\ii s_0^2$, the local law and Ward's identity give
\[
 \|R_H(z_h)v\|^2
 =-\frac{\Im v^*R_H(z_h)v}{2s_0^2}
 \le\frac{Cs_0+m_N^{-1}N^{o(1)}}{2s_0^2}
 \le m_NN^{o(1)}.
\]
For $r\ge s_0^2/2$ and $|y|\le2s_0^2$, $1/(r^2+y^2)\le17/(r^2+4s_0^4)$. Spectral integration and the resolvent identity consequently yield
\[
 \begin{aligned}
 \|R_H(x+\ii y)v\|^2&\le17\|R_H(z_h)v\|^2,\\
 |u^*(R_H(x+\ii y)-R_H(z_h))v|
 &\le Cs_0^2m_NN^{o(1)}=m_N^{-1}N^{o(1)},\\
 |g(x+\ii y)-g(z_h)|&\le Cs_0=m_N^{-1}N^{o(1)}.
 \end{aligned}
\]
These inequalities extend every required entry estimate to the right arc, including its intersection with the real axis. Apply the same argument to $-H$, using $R_{-H}(-z)=-R_H(z)$
and $g(-z)=-g(z)$, for the left arc. Conjugation handles the lower arcs.

It remains to bound $v^*|R_H(z)|v$, because the cubic trace estimate uses this quantity rather than a squared column norm. For each required deterministic unit vector $v$, let
$\mu_v(B)=v^*\mathbf1_B(H)v$ be its positive spectral measure. For $|E|\le4$ and $c s_0^2\le r\le16$,
\[
 -\Im v^*R_H(E+\ii r)v
   =\int\frac{r}{(E-x)^2+r^2}\,\mu_v(\mathrm dx).
\]
The integrand is at least $1/(2r)$ on $[E-r,E+r]$. The isotropic estimate, with $|g_t(E+\ii r)|\le C$, therefore implies, on its event of validity,
\[
 \mu_v([E-r,E+r])
 \le2r[-\Im v^*R_H(E+\ii r)v]
 \le C N^\varepsilon r.
\]
Here $\varepsilon>0$ may be chosen arbitrarily small; the exceptional probability may be chosen smaller than any fixed negative power of $N$. The estimates are made simultaneous in
$E,r$ below. This permits their use at the radii depending on $z$ in the next calculation.

Write $z=E+\ii\eta$. On the part of the ellipse with $|\eta|\ge c s_0^2$, set $r_0=|\eta|$ and $r_j=2^jr_0$. The central interval contributes at most
\[
 \int_{|x-E|\le r_0}|z-x|^{-1}\,\mu_v(\mathrm dx)
 \le r_0^{-1}\mu_v([E-r_0,E+r_0])\le C N^\varepsilon.
\]
For each annulus $r_j<|x-E|\le2r_j$, one has $|z-x|^{-1}\le r_j^{-1}$ and hence
\[
 \int_{r_j<|x-E|\le2r_j}|z-x|^{-1}\,\mu_v(\mathrm dx)
 \le r_j^{-1}\mu_v([E-2r_j,E+2r_j])
 \le C N^\varepsilon.
\]
Near either real point of the ellipse, the spectral distance just proved is at least $s_0^2/2$. Use $r_0=s_0^2/4$ there. The central interval has zero spectral mass, and the same
annular estimate applies. On \eqref{W:eq:cumulant-confinement}, the spectrum and $E$ lie in fixed bounded intervals. The annuli with $j\le C\log N$ therefore cover the spectrum,
giving
\[
 v^*|R_H(z)|v\le C N^\varepsilon(1+\log N).
\]
Starting with $\varepsilon/2$ in the local estimate absorbs $1+\log N$ into $N^\varepsilon$. This proves the required stochastic bound for $|R_H|$.

We specify how the estimates hold simultaneously in the spectral parameter. Take finite deterministic sets with spacing at most $N^{-20}$ in the bounded nonreal domains, on the
ellipse, and on the vertical segments just used. A two dimensional covering uses at most $C N^{40}$ points. The vector pairs needed for the coordinate estimates and the vectors in
$\mathcal V_N$ number at most $C N^2$. On the spectral event, each relevant spectral distance is at least $N^{-1}$ for large $N$. The resolvent identity and the spectral theorem
give
\[
 \begin{aligned}
 |u^*(R_H(z)-R_H(z'))v|&\le N^2|z-z'|,\\
 \big|\|R_H(z)v\|^2-\|R_H(z')v\|^2\big|&\le2N^3|z-z'|,\\
 \big|v^*(|R_H(z)|-|R_H(z')|)v\big|&\le N^2|z-z'|
 \end{aligned}
\]
for unit vectors and pairs of points in these sets. Thus the approximation error is at most $C N^{-17}$. Taking the local probability exponent to be $D+46$ bounds the union of the
$C N^{42}$ exceptional events by $C N^{-D-4}$. Continuity then extends the inequalities to the complete domains. In the entry bound, uniform moments and Markov's inequality give
\[
 \mathbb P\left(\max_{a,b}|\sqrt N H_{ab}|>N^\varepsilon\right)
 \le C_pN^{2-p\varepsilon}.
\]
Choose an even $p>(D+4)/\varepsilon$ and then take $N$ sufficiently large. This supplies the same prescribed failure probability for the maximum entry.

We also prove that the bounds hold under partial deletion of an entry. Put $z_\pm=\pm(2+d/4)$ and $K_{ab}=H-V$, where $V$ is supported on $J=\{a,b\}$. On
\eqref{W:eq:cumulant-confinement}, the upper barrier satisfies $z_+-\lambda\ge d/8$. For $r=z_+-\lambda$,
\[
 \frac1r\le65\frac{r}{r^2+d^2}.
\]
The local estimate at $z_++\ii d$ bounds its coordinate diagonal entries by a fixed constant. Integration of this inequality against the coordinate spectral measure bounds
$R_{H,ii}(z_+)$. The matrix $R_H(z_+)$ is positive, so
$|R_{H,ij}(z_+)|^2\le R_{H,ii}(z_+)R_{H,jj}(z_+)$.
Apply this calculation to $-H$ for $z_-$. It follows that both compressions $R_{H,JJ}(z_\pm)$ have uniformly bounded norms. For all $0\le u\le1$ and sufficiently large $N$,
$\|uV_{JJ}R_{H,JJ}(z_\pm)\|\le1/2$, since $\|V\|=N^{-1/2+o(1)}$. Hence
\[
 \det(z_\pm I-H+uV)
 =\det(z_\pm I-H)\det(I_J+uV_{JJ}R_{H,JJ}(z_\pm))\ne0.
\]
The ordered eigenvalues of $H-uV$ depend continuously on $u$. They begin strictly between the two barriers and cannot cross either barrier, so
\[
 \sigma(H-uV)\subset(-2-d/4,2+d/4),\qquad 0\le u\le1.
\]
This interval lies strictly inside the interval where $h_d$ equals the identity.

On $\mathcal C$, the finite matrix inverse formula gives
\[
 R_{H-uV}-R_H
 =-R_{H,\cdot J}(I_J+uV_{JJ}R_{H,JJ})^{-1}
                        uV_{JJ}R_{H,J\cdot}.
\]
The inverse in the middle has norm at most two. The bounded coordinate entries and the estimates for $R_Hw$ yield
\[
 \begin{aligned}
 |(R_{H-uV}-R_H)_{ij}|&\le N^{-1/2+o(1)},\\
 |((R_{H-uV}-R_H)w)_i|
   &\le N^{-1/2+o(1)}(\omega_a+\omega_b)
    \le m_N^{-1}N^{o(1)}.
 \end{aligned}
\]
Because $\omega_i\ge m_N^{-1}$, addition of these errors preserves the coordinate estimates in the lemma. This calculation is uniform in $u$ and in the deleted pair. All suprema
over $z$ and $u$ are measurable: on the displayed spectral intervals the matrices and resolvents depend continuously on these parameters, so countable dense subsets give the same
suprema.

Let $X_{N,t,w}$ be the maximum of the normalized suprema in \eqref{W:eq:cumulant-event-bounds}, including the first two bounds for all partially deleted matrices. Set
$X_{N,t,w}=\infty$ if $\|H(t)\|>2+d/8$ or if $\|H(t)-uV_{ab}\|>2+d/4$ for some $a<b$ and $u\in[0,1]$. The preceding estimates say that, for each fixed $\varepsilon,D>0$, there is
$N_0(\varepsilon,D)$ such that $\sup_{t,w}\mathbb P(X_{N,t,w}>N^\varepsilon)\le N^{-D}$ for all $N\ge N_0(\varepsilon,D)$. Choose increasing $N_k\ge e^{k^2}$ so that
\[
 \sup_{t,w}\mathbb P(X_{N,t,w}>N^{1/k})\le N^{-k}
 \qquad(N\ge N_k).
\]
For $k(N)=\max\{k:N_k\le N\}$, take $q_N=N^{1/k(N)}$. For $N<N_1$, take any $q_N\ge1$. Then $q_N=N^{o(1)}$ and $\mathcal E_{N,t,w}=\{X_{N,t,w}\le q_N\}$ has failure at most
$N^{-k(N)}$. A fixed polynomial loss $N^B$ on its complement has integrated expectation at most $N^{B-k(N)}$. All suprema over $t,w$ are outside the probabilities, as required.
\end{proof}

\subsubsection{Errors from the first two truncated moments}
\label{W:sec:errors-from-the-first-two-truncated-moments}

The comparison between $H_N^{\rm tr}$ and $H_N^{\rm b}$ has small errors already in its first two moments. The bounded distribution matches the centered, covariance normalized
truncation, whereas the interpolation uses the uncentered array. In particular, the quadratic coefficient is paired with the noncentral second moment. The following estimates
express these differences in terms of the tail function, using Proposition~\ref{W:prop:comparison-global}, Lemma~\ref{W:lem:tightness-global}, and the conditional moment estimates
of Lemma~\ref{W:lem:comparison-weights}.

The tensors are expressed in the standardized real coordinates $X$, with the corresponding change $N^{-1/2}\mathcal T_\alpha x$ in the matrix from
\eqref{W:eq:standardized-entry-insertion}; the factor $N^{-1/2}$ belongs to each matrix derivative. The notation $r(X)$ and $\Theta$ is that of
Section~\ref{W:sec:bounded-entries-with-prescribed-moments}, and $L=L_N$. With $\widehat X,\mu,\Sigma,X^{\rm nor}$ from \eqref{W:eq:normalized-truncated-entry-vector} and the
bounded vector $Y$ of Lemma~\ref{W:lem:cubic-defect}, we have $\mathbb EY=0$, $\mathbb EYY^{\mathsf T}=I$, and $\mathbb EY^{\otimes3}=\mathbb E(X^{\rm nor})^{\otimes3}$. The first
two moment differences are

\[
 \Delta M_1=\mathbb E\widehat X-\mathbb EY,\qquad
 \Delta M_2=\mathbb E\widehat X\widehat X^{\mathsf T}
                                      -\mathbb EYY^{\mathsf T}.
\]
They are measured in the ordered tensor norm of \eqref{W:eq:B-tensor-l1}.

The identities $\mathbb EX=\mathbb EY=0$ and $\mathbb EXX^{\mathsf T}=\mathbb EYY^{\mathsf T}=I$ give
\begin{equation}
 \Delta M_1=-\mathbb E[X\mathbf1_{\{r(X)>L\}}],\qquad
 \Delta M_2=-\mathbb E[XX^{\mathsf T}\mathbf1_{\{r(X)>L\}}].
 \label{W:eq:replacement-exact-defects}
\end{equation}
Apply \eqref{W:eq:quant-tail-below-four} to the two discarded moments. Use $\|X\|_1\le Cr(X)$ and $\|XX^{\mathsf T}\|_{1,\mathrm{ord}}=\|X\|_1^2\le Cr(X)^2$ to obtain
\begin{equation}
 \|\Delta M_1\|_{1,\mathrm{ord}}\le C\Theta(L)L^{-3},\qquad
 \|\Delta M_2\|_{1,\mathrm{ord}}\le C\Theta(L)L^{-2}.
 \label{W:eq:replacement-tail-defects}
\end{equation}

The second identity uses the noncentral second moment. The relation to the covariance is
\[
 \Sigma-I=-\mathbb E[XX^{\mathsf T}\mathbf1_{\{r(X)>L\}}]
                  -\mu\mu^{\mathsf T},\qquad
 \Delta M_2=(\Sigma-I)+\mu\mu^{\mathsf T}.
\]
In particular, $\|\Sigma-I\|\le C\Theta(L)L^{-2}+C \Theta(L)^2L^{-6}\to0$, which also justifies the inverse square root used in defining $X^{\mathrm{nor}}$. The squared mean term
in $\Sigma-I$ cancels when the raw second moment is reconstructed. It is therefore absent from \eqref{W:eq:replacement-exact-defects}. In the complex class the off-diagonal
coordinate of this tensor is $\Delta\mathbb E(X_1X_2)$, in addition to the two diagonal second moments.

The difference of the third order coordinate moment tensors is \eqref{W:eq:cubic-affine}, with $Y$ the bounded comparison vector. Applying Lemma~\ref{W:lem:cubic-defect} and the
tail estimates at $L=e_N\sqrt N$ gives
\[
 \|\Delta M_3\|_{1,\mathrm{ord}}
 \le C\Theta(L)L^{-2}
 =C\Theta(L)e_N^{-2}N^{-1}=N^{-1+o(1)}.
\]
Here $e_N^{-1}=N^{o(1)}$ and $\Theta(L)$ is bounded.

Let $\mu_\alpha$ be the law of $\widehat X$ and $\nu_\alpha$ the law of $Y$. At parameter $t$, give the independent blocks the product law
$\bigotimes_\alpha[(1-t)\nu_\alpha+t\mu_\alpha]$, with the corresponding diagonal laws, and set block $\alpha$ to zero to obtain $K_\alpha(t)$. For each fixed $t$, the replacement
entry is independent of $K_\alpha(t)$. Put
\[
 s_\alpha(x)=s(K_\alpha(t)+N^{-1/2}\mathcal T_\alpha x,w),\quad
 s_0=s_\alpha(0),\quad
 \ell_\alpha(x)=D_{x}s_\alpha(0)[x],\quad
 q_\alpha(x)=\tfrac12D_{x}^2s_\alpha(0)[x,x].
\]
Write $\ell_r=\ell_\alpha(e_r)$ and $q_{rs}=\tfrac12D_{x}^2s_\alpha(0)[e_r,e_s]$. The first two coefficient tensors in the entrywise Taylor expansion are
\begin{equation}
 (C_{1,\alpha})_r=D_{s}A_N(s_0)[\ell_r],\qquad
 (C_{2,\alpha})_{rs}
       =D_{s}A_N(s_0)[q_{rs}]
          +\tfrac12D_{s}^2A_N(s_0)[\ell_r,\ell_s].
 \label{W:eq:first-two-entry-coefficients}
\end{equation}
The associated polynomial is $\langle C_{1,\alpha},x\rangle+ \langle C_{2,\alpha},x^{\otimes2}\rangle$. The definition of $C_{2,\alpha}$ includes the factor $1/2$ from the second
derivative.

For the compact ball $\mathscr X_{\alpha,N}$ of permitted standardized entry values, enlarged as in Lemma~\ref{W:lem:resolvent-weight-input}, define
\[
 \mathcal A_{\alpha,\rho}
 =\{K_\alpha:\ \min_{x\in\mathscr X_{\alpha,N}}
       \lambda_{\max}(G_\rho(K_\alpha+V_\alpha(x),w))\le\kappa^{-1}\}.
\]
Continuity in $K_\alpha,x$ and compactness of the ball make this event measurable with respect to $\mathscr F_\alpha$. Let $\mathcal B_\alpha$ be the event of
Lemma~\ref{W:lem:resolvent-weight-input}: all permitted matrices have spectra in the stated interval and satisfy \eqref{W:eq:resolvent-coordinate} on $\mathcal C$. Put
$I_\alpha=\mathbf1_{\mathcal B_\alpha\cap\mathcal A_{\alpha,\rho}}$. Both events are fixed before the replacement variable is integrated. On $I_\alpha=1$,
Lemma~\ref{W:lem:comparison-weights} and \eqref{W:eq:tightness-domination} give
\begin{equation}
 \|\ell_\alpha(x)\|_\diamond
       \le CN^{-1/2}\|x\|_1U_\alpha,\qquad
 \|q_\alpha(x)\|_\diamond
       \le CN^{-1}\|x\|_1^2V_\alpha.
 \label{W:eq:first-two-moment-derivatives}
\end{equation}
These are the matrix Taylor coefficients of the moments, with the empirical trace in $s_n=u_n+x_n/{m_N}$ cancelling exactly under reconstruction of the moments at the starting
vector.

On this regularity event let $R=(zI-K_\alpha)^{-1}$ and $T_r=\mathcal T_\alpha e_r$. The Cauchy kernels for the linear functionals $\ell_r$ and $q_{rs}$ give
\begin{align*}
 \ell_r(f)&=\frac{N^{-1/2}}{2\pi\mathrm i}
            \int_{\mathcal C} f(z)w^*RT_rRw\,\mathrm dz,\\
 q_{rs}(f)&=\frac{N^{-1}}{4\pi\mathrm i}
       \int_{\mathcal C} f(z)w^*
          (RT_rRT_sR+RT_sRT_rR)w\,\mathrm dz.
\end{align*}
Thus the first derivative uses the two distinct support coordinates, whereas a second derivative is bounded by the sum of their squares using the bounded restricted resolvent. This
is precisely the origin of $U_\alpha$ and $V_\alpha$ in \eqref{W:eq:first-two-moment-derivatives}.

For $r\ne s$, polarization reads
\[
 q_{rs}=\tfrac14
       \bigl(q_\alpha(e_r+e_s)-q_\alpha(e_r-e_s)\bigr).
\]
Since $\|e_r+e_s\|_1=\|e_r-e_s\|_1=2$,
\[
 \|q_{rs}\|_\diamond
 \le\frac14\bigl(\|q_\alpha(e_r+e_s)\|_\diamond
                 +\|q_\alpha(e_r-e_s)\|_\diamond\bigr)
 \le 2CN^{-1}V_\alpha.
\]
For $r=s$ the same estimate follows from $q_{rr}=q_\alpha(e_r)$. Thus every ordered component has the required bound. Together with \eqref{W:eq:replacement-test-class}, this proves
\begin{equation}
 \|C_{1,\alpha}\|_{\max}
     \le C\Lambda_N{m_N}N^{-1/2}U_\alpha,\qquad
 \|C_{2,\alpha}\|_{\max}
     \le C\Lambda_NN^{-1}({m_N}V_\alpha+{m_N}^2U_\alpha^2).
 \label{W:eq:replacement-coefficient-bounds}
\end{equation}
Here $\|C\|_{\max}$ is the maximum absolute coordinate; hence $|\langle C,M\rangle|\le\|C\|_{\max}\|M\|_{1,\mathrm{ord}}$.

The first two derivative orders involve three sums of resolvent weights. Corollary~\ref{W:cor:integrated-comparison-weights}, for the weights constructed in
Lemma~\ref{W:lem:comparison-weights}, gives
\begin{equation}
 \int_0^1\sum_{\alpha\,\mathrm{off}}\mathbb E[I_\alpha U_\alpha]\,\mathrm dt
       \le Q_N{m_N}^4,\quad
 \int_0^1\sum_{\alpha\,\mathrm{off}}\mathbb E[I_\alpha V_\alpha]\,\mathrm dt
       \le Q_N{m_N}^4,\quad
 \int_0^1\sum_{\alpha\,\mathrm{off}}\mathbb E[I_\alpha U_\alpha^2]\,\mathrm dt
       \le Q_N{m_N}^2,
 \label{W:eq:replacement-three-integrated-weight-sums}
\end{equation}
where $Q_N=N^{o(1)}$ is deterministic and uniform over deterministic starting vectors. The relevant lines of that corollary are the $U$, $V$, and $U^2$ bounds, with
$\rho\asymp m_N$. Their right sides are $N\rho$, $N\rho$, and $\rho^2$, up to a common subpolynomial factor, giving the displayed powers.

On the exceptional event, the global polynomial bounds in Lemma~\ref{W:lem:global-first-two-entry-derivatives} replace the resolvent weight estimates.

At fixed $K_\alpha(t)$ the tensors and $I_\alpha$ are fixed. Thus
\begin{equation}
 I_\alpha\int\left(
 \langle C_{1,\alpha},x\rangle+
 \langle C_{2,\alpha},x^{\otimes2}\rangle\right)
             (\mu_\alpha-\nu_\alpha)(\mathrm dx)
 =I_\alpha\sum_{q=1}^2\langle C_{q,\alpha},\Delta M_{q,\alpha}\rangle.
 \label{W:eq:replacement-conditional}
\end{equation}
This integrates the entire signed joint entry distribution. In particular, in real dimension two the quadratic contraction is
\[
 C_{11}\Delta M_{11}+2C_{12}\Delta M_{12}+C_{22}\Delta M_{22}.
\]
The factorial is already included in $C_{2,\alpha}$. Both cutoffs belong to $A_N$, so the integral is over the full distribution of $x$.

\Needspace{6\baselineskip}

Define
\begin{align*}
 \mathcal E_{12,N}^{\rm good}(A;w)
 &=\int_0^1\sum_{\alpha\,\mathrm{off}}\mathbb E I_\alpha
       \sum_{q=1}^2|C_{q,\alpha}:\Delta M_{q,\alpha}|\,\mathrm dt,\\
 \mathcal E_{12,N}(A;w)
 &=\int_0^1\sum_{\alpha\,\mathrm{off}}\mathbb E
       \sum_{q=1}^2|C_{q,\alpha}:\Delta M_{q,\alpha}|\,\mathrm dt.
\end{align*}
For the bound $Q_N$ in \eqref{W:eq:replacement-three-integrated-weight-sums}, enlarged by a fixed constant as necessary, set
\[
 \varepsilon_{12,N}=C\Lambda_NQ_N\Theta(L_N)
 \left[e_N^{-3}N^{-1/3}
       +e_N^{-2}(N^{-1/3}+N^{-2/3})\right].
\]

\begin{theorem}
\label{W:thm:replacement-low-order}
There is a deterministic $Q_N=N^{o(1)}$, uniform in $t$ and in the deterministic unit vector $w$, for which
\begin{align}
 \sup_{\|w\|=1}\mathcal E_{12,N}^{\rm good}(A;w)
    &\le \varepsilon_{12,N},\label{W:eq:replacement-main}\\
 \sup_{\|w\|=1}\mathcal E_{12,N}(A;w)
    &\le \varepsilon_{12,N}+C_DN^{-D}.
       \label{W:eq:first-two-moment-errors}
\end{align}
The second estimate holds for each fixed $D>0$ after choosing the regularity event with a sufficiently small exceptional probability. In particular, the full sum is
$N^{-1/3+o(1)}$.
\end{theorem}
\begin{proof}
Write $\tau_j=\sup_\alpha\|\Delta M_{j,\alpha}\|_{1,\mathrm{ord}}$. For each ordered tensor, finite-dimensional duality gives
\[
 |C_{j,\alpha}:\Delta M_{j,\alpha}|
 \le \sum_{r_1,\ldots,r_j}|(C_{j,\alpha})_{r_1\cdots r_j}|
                         |(\Delta M_{j,\alpha})_{r_1\cdots r_j}|
 \le\|C_{j,\alpha}\|_{\max}\tau_j.
\]
The moment differences are deterministic. Multiply this inequality by the indicator $I_\alpha$, which depends only on $K_\alpha(t)$, take the expectation in $K_\alpha(t)$, and
integrate $t$. Equations \eqref{W:eq:replacement-coefficient-bounds} and \eqref{W:eq:replacement-three-integrated-weight-sums} give
\begin{align*}
 \mathcal E_{1,N}^{\rm good}
  &\le C\Lambda_NQ_N\tau_1{m_N}N^{-1/2}{m_N}^4
    =C\Lambda_NQ_N\tau_1{m_N}^5N^{-1/2},\\
 \mathcal E_{2,N}^{\rm good}
  &\le C\Lambda_NQ_N\tau_2N^{-1}({m_N}\cdot {m_N}^4+{m_N}^2\cdot {m_N}^2)
    =C\Lambda_NQ_N\tau_2N^{-1}({m_N}^5+{m_N}^4).
\end{align*}
Equation~\eqref{W:eq:replacement-tail-defects} gives $\tau_1\le C\Theta(L_N)e_N^{-3}N^{-3/2}$ and $\tau_2\le C\Theta(L_N)e_N^{-2}N^{-1}$. Since $N={m_N}^3$, the preceding bounds
become
\[
 \mathcal E_{1,N}^{\rm good}
     \le C\Lambda_NQ_N\Theta(L_N)e_N^{-3}{m_N}^{-1},\qquad
 \mathcal E_{2,N}^{\rm good}
     \le C\Lambda_NQ_N\Theta(L_N)e_N^{-2}({m_N}^{-1}+{m_N}^{-2}).
\]
This proves \eqref{W:eq:replacement-main}. Since $e_N^{-1},\Lambda_N,Q_N=N^{o(1)}$ and $\Theta(L_N)$ is bounded, the right side is $N^{-1/3+o(1)}$.

For $K_\alpha(t)$, write the $\sigma(K_\alpha(t))$-measurable regularity indicator as $I_\alpha=\mathbf1_{\mathcal B_\alpha\cap\mathcal A_{\alpha,\rho}}$. On
$\mathcal B_\alpha\setminus\mathcal A_{\alpha,\rho}$, the minimum in the definition of $\mathcal A_{\alpha,\rho}$ is strictly greater than $\kappa^{-1}$. Indeed the continuous
function attains its minimum on the compact ball. Each permitted entry value therefore has a neighbourhood where at least one Gram eigenvalue lies above the support of the scalar
determinant cutoff. The factor $\Gamma_\rho$, and hence the whole function $F_{A,N}$, vanishes in this neighbourhood. At $x=0$ every coordinate derivative is consequently zero, in
particular $C_{1,\alpha}=C_{2,\alpha}=0$. Thus the part excluded by $I_\alpha$ has support only in $\mathcal B_\alpha^c$. On $\mathcal B_\alpha^c$, use the deterministic bounds of
Lemma~\ref{W:lem:global-first-two-entry-derivatives}, with $B_1=B_2=6$. If $\sup_{\alpha,t}\mathbb P(\mathcal B_\alpha^c)\le N^{-D_*}$, then this part of the contracted coefficient
is at most
\begin{equation}
 C N^{2-D_*}\bigl(N^{B_1}\tau_1+N^{B_2}\tau_2\bigr).
 \label{W:eq:replacement-bad-coefficients}
\end{equation}
All moment differences here are deterministic moment tensors, so their contraction can be bounded before taking the probability over the remaining entries. Since $\tau_1,\tau_2$
are bounded for large $N$, choose $D_*>D+2+\max\{B_1,B_2\}+1$. Equation \eqref{W:eq:replacement-bad-coefficients} is then $O(N^{-D})$. The deterministic subpolynomial resolvent
weight bound $Q_N$ may depend on this fixed choice of failure exponent; this does not change any displayed power of $N$. Together with \eqref{W:eq:replacement-main}, this proves
\eqref{W:eq:first-two-moment-errors}.
\end{proof}

The preceding estimate controls a specific part of the exact comparison. To identify that part, write $F_{A,N}(H)=A_N(s(H,w))$ and $\Delta\mu_\alpha=\mu_\alpha-\nu_\alpha$. For an
off-diagonal entry, define
\[
 R_{4,\alpha}(K;x)
 =F_{A,N}(K+V_\alpha(x))-F_{A,N}(K)
       -\sum_{j=1}^3 C_{j,\alpha}(K):x^{\otimes j}.
\]
Substitution of this definition in \eqref{W:eq:entry-distribution-interpolation} cancels the constant term, since $\Delta\mu_\alpha$ has total mass zero. The first two
contributions to Lemma~\ref{W:lem:entry-distribution-interpolation} are therefore
\begin{equation}
 \mathcal L_{12,N}(A;w)
 =\int_0^1\sum_{\alpha\,\mathrm{off}}\mathbb E
        \sum_{q=1}^2\langle C_{q,\alpha},\Delta M_{q,\alpha}\rangle
                  \,\mathrm dt.
 \label{W:eq:replacement-low-order-functional}
\end{equation}
For each fixed $t$, $C_{j,\alpha}$ is $\mathscr F_\alpha$-measurable and $\Delta M_{j,\alpha}$ is deterministic. Thus
\[
 |\mathcal L_{12,N}(A;w)|\le\mathcal E_{12,N}(A;w)
 \le\varepsilon_{12,N}+C_DN^{-D}.
\]
The full comparison consequently reads
\[
 \begin{split}
 \mathbb E F_{A,N}(H_N^{\rm tr})-\mathbb E F_{A,N}(H_N^{\rm b})
 ={}&\mathcal L_{12,N}(A;w)
    +\int_0^1\sum_{\alpha\,\mathrm{off}}
             \mathbb E(C_{3,\alpha}:\Delta M_{3,\alpha})\,\mathrm dt\\
 &+\int_0^1\sum_{\alpha\,\mathrm{off}}\mathbb E
           \int R_{4,\alpha}(K_\alpha;x)\,\Delta\mu_\alpha(\mathrm dx)\,\mathrm dt
       +\mathcal E_{{\rm diag},N}(A;w).
 \end{split}
\]
Here $\mathcal E_{{\rm diag},N}$ denotes the sum over diagonal entries in the same interpolation identity. At fixed $N$, finite support of the distributions and the global
derivative bounds justify all integrals and sums. The subsequent estimates bound their absolute integrals uniformly in $t$ and $w$. If the first two coordinate moment tensors
agree, $\mathcal L_{12,N}=0$ by its displayed definition. For the uncentered truncation its value is bounded using \eqref{W:eq:replacement-exact-defects}.

The diagonal mean is estimated in Theorem~\ref{W:thm:truncation-diagonal} below. For centered original diagonals after truncation on the original probability space,
\[
 \sup_i|\mathbb E\widehat d_{N,i}|
       \le \Theta_{\rm d}(L_N)/L_N.
\]
The coefficient is bounded by $C\Lambda_N{m_N}N^{-1/2}V_i$, and $\int_0^1\sum_i\mathbb E[I_iV_i]\,\mathrm dt\le {m_N}N^{o(1)}$ by
Lemma~\ref{W:lem:truncation-integrated-weight-sums}. Consequently the sum of the linear contributions is bounded by $C\Lambda_NQ_N{m_N}^2N^{-1/2}\Theta_{\rm d}(L_N)/L_N$, namely
\begin{equation}
 C\Lambda_NQ_N e_N^{-1}\Theta_{\rm d}(L_N)N^{-1/3}.
 \label{W:eq:replacement-diagonal-mean}
\end{equation}
The entire quadratic contribution of the diagonal entries is included in the summable first-order remainder estimated in
Section~\ref{W:sec:conditional-row-bounds-for-diagonal-remainders}.

The first two moment contributions are now controlled. The diagonal estimate and the off-diagonal remainder of the expansion are proved next; the comparison is assembled in
Section~\ref{W:sec:removal-of-coefficient-cutoff-in-comparison}. The nonlinear localization statistic, which also contains an unnormalized trace, is reserved for
Section~\ref{W:sec:localization-of-the-low-energy-spectral-subspace}.

\subsubsection{Comparison of the diagonal entries}
\label{W:sec:diagonal-second-moment-comparison}

\label{W:sec:diagonal-comparison-under-second-moment-assumptions}

The diagonal entries are subject only to the second moment assumptions, so their comparison uses a linear expansion with a quadratic remainder. The conditional estimate of that
remainder also involves a fourth moment, but this is the fourth moment of the independent off-diagonal row. The replacement diagonal variable itself is integrated using its second
moment.

The observables and derivative bounds are those of \eqref{W:eq:replacement-test-class}. For the two standardized real diagonal distributions $\mu_i,\nu_i$, supported in
$[-CL_N,CL_N]$ with $L_N$ from \eqref{W:eq:LN-definition}, the relevant quantities are

\begin{equation}
 M_{2,N}=\sup_i\int x^2(\mu_i+\nu_i)(\mathrm dx),\qquad
 b_N=\sup_i\left|\int x(\mu_i-\nu_i)(\mathrm dx)\right|.
 \label{W:eq:truncation-moment-data}
\end{equation}
The off-diagonal entries have one of the intermediate distributions in Section~\ref{W:sec:moment-increments-under-entry-replacement}. In particular, their matrix entry means are
$N^{-2+o(1)}$, their centered variances are at most $C/N$, and their centered fourth moments are at most $F_N/N^2$, $F_N=N^{o(1)}$. These statements hold uniformly over the
interpolation parameter and remain true after deleting entries.

At parameter $t$, let $K_i(t)$ be the intermediate matrix with its $i$th diagonal entry replaced by zero. It is independent of the scalar $x$ subsequently used to replace that
entry. Let $\mathscr F_i=\sigma(K_i(t))$. The event $\mathcal B_i$ of Lemma~\ref{W:lem:resolvent-weight-input} is $\mathscr F_i$-measurable. On it, for every $x$ in the enlarged
interval of permitted values, the spectrum of $K_i(t)+N^{-1/2}xE_{ii}$ lies strictly inside both $\mathcal C$ and $(-2-d,2+d)$, and \eqref{W:eq:resolvent-coordinate} holds. With
$\mathscr X_{i,N}$ denoting that compact interval, set
\[
 \mathcal A_{i,\rho}
 =\{K_i:\ \min_{x\in\mathscr X_{i,N}}
       \lambda_{\max}(G_\rho(K_i+N^{-1/2}xE_{ii},w))\le\kappa^{-1}\},
 \qquad I_i=\mathbf1_{\mathcal B_i\cap\mathcal A_{i,\rho}}.
\]
The minimum is attained and depends measurably on $K_i$. On $\mathcal B_i\setminus\mathcal A_{i,\rho}$ the factor $\Gamma_\rho$ makes each function under consideration identically
zero near every permitted replacement and near $x=0$. This partitions the comparison into the contribution with indicator $I_i$, a zero contribution, and the exceptional event
$\mathcal B_i^c$.

\Needspace{6\baselineskip}

Write the diagonal contribution in Lemma~\ref{W:lem:entry-distribution-interpolation} as
\begin{equation}
 \mathcal E_{{\rm diag},N}(A_N;w)
 =\int_0^1\sum_{i=1}^N\mathbb E_{K_i(t)}
 \int A_N\bigl(s(K_i(t)+N^{-1/2}xE_{ii},w)\bigr)
                (\mu_i-\nu_i)(\mathrm dx)\,\mathrm dt.
 \label{W:eq:truncation-diagonal-contribution}
\end{equation}
When only the diagonal distributions vary, this is the difference of the two expectations. Let $D$ be the chosen exponent in $\sup_{i,t}\mathbb P(\mathcal B_i^c)\le N^{-D}$.

\begin{theorem}
\label{W:thm:truncation-diagonal}
For the functions satisfying \eqref{W:eq:replacement-test-class}, there is a deterministic $Q_N=N^{o(1)}$, uniform in $w$, such that
\begin{equation}
 \sup_{\|w\|=1}|\mathcal E_{{\rm diag},N}(A_N;w)|
 \le C\Lambda_NQ_N
       \left(M_{2,N}N^{-1/3}+b_NN^{1/6}\right)+CN^{1-D}.
 \label{W:eq:diagonal-comparison-bound}
\end{equation}
For bounded $M_{2,N}$, the sum of the absolute first-order remainders of the expansions is $N^{-1/3+o(1)}$. If also $b_N\le N^{-1/2+o(1)}$, the comparison in
\eqref{W:eq:diagonal-comparison-bound} has the same order.
\end{theorem}

Fix $i,t$ and write $K=K_i(t)$, $R(z)=(zI-K)^{-1}$ and $u=x/\sqrt N$. Here $u$ denotes the replacement matrix entry. On $\mathcal B_i$, $|R_{ii}(z)|\le C_0$. For every $N\ge N_0$
with $CC_0e_N\le1/2$, the bound $|u|\le Ce_N$ gives
\begin{equation}
 |1-uR_{ii}(z)|\ge\tfrac12\quad(z\in\mathcal C).
 \label{W:eq:truncation-denominator}
\end{equation}
The rank-one resolvent identity is exact:
\begin{equation}
 R_{K+uE_{ii}}(z)-R(z)
 =\frac{u}{1-uR_{ii}(z)}R(z)E_{ii}R(z).
 \label{W:eq:truncation-SM}
\end{equation}
It follows by multiplying the right side plus $R$ by $zI-K-uE_{ii}$; the coefficient of $E_{ii}R$ cancels. Let
\[
 \zeta_i(z)=(w^*R(z)e_i)(e_i^*R(z)w),\qquad
 \Delta_i(x)=s(K+uE_{ii},w)-s(K,w).
\]
Since $R(z)^*=R(\bar z)$, $|w^*R(z)e_i|=|(R(\bar z)w)_i|$. The event used here is $\mathcal B_i$. On that event $h_d(B)=B$ for $B=K$ and for every $B=K+uE_{ii}$ with a permitted
$u$, because their spectra lie in $(-2-d,2+d)$. Their spectra also lie in the interior of the positively oriented ellipse $\mathcal C$. Polynomial functional calculus therefore
gives
\[
 s(B,w)(f)=w^*f(B)w
       =\frac1{2\pi\mathrm i}\int_{\mathcal C}f(z)w^*R_B(z)w\,\mathrm dz.
\]
Subtract these two identities and use \eqref{W:eq:truncation-SM}. The resulting equality holds for every polynomial $f$ of degree at most $2\rho-2$: equivalently, it holds on a
basis of the finite dimensional polynomial space and then extends by linearity. We obtain
\begin{equation}
 \Delta_i(x)(f)=\frac1{2\pi\mathrm i}\int_{\mathcal C}
 f(z)\frac{u\zeta_i(z)}{1-uR_{ii}(z)}\,\mathrm dz.
 \label{W:eq:truncation-exact-moment}
\end{equation}
Let $v_i(f)=(2\pi\mathrm i)^{-1}\int f(z)\zeta_i(z)\,\mathrm dz$ and $\ell_i(x)=uv_i$. Differentiating $s(K+uE_{ii},w)(1)=1$ shows $v_i(1)=0$; also $\Delta_i(x)(1)=0$. The
subtracted linear term is exactly the matrix moment derivative, and
\begin{equation}
 (\Delta_i(x)-\ell_i(x))(f)=\frac1{2\pi\mathrm i}
 \int_{\mathcal C}f(z)
       \frac{u^2\zeta_i(z)R_{ii}(z)}{1-uR_{ii}(z)}\,\mathrm dz.
 \label{W:eq:diagonal-moment-remainder}
\end{equation}

Choose a deterministic $q_N=N^{o(1)}$ dominating the constants in the preceding resolvent estimates and define
\[
 a_i(z)=|(R(z)w)_i|+|(R(\bar z)w)_i|+\sqrt{q_N\rho/N},
 \qquad V_i=\int_{\mathcal C}a_i(z)^2\nu_\rho(\mathrm dz).
\]
Every $i$ uses its own matrix with the selected entry set to zero. Outside the event where the stated spectral and resolvent bounds hold define the resolvent weight to be zero. The
deterministic contour measure is \eqref{W:eq:comparison-measure} and has mass $N^{o(1)}$. If a linear functional $v$ with $v(1)=0$ has Cauchy kernel $k$, then
\[
 |v(p_jp_k)|\le\int|k(z)|\nu_\rho(\mathrm dz),\qquad
 \|\dot G_\rho[v]\|_1\le\rho\int|k(z)|\nu_\rho(\mathrm dz).
\]
The second inequality follows by integrating the rank-one matrix $p(z)p(z)^{\mathsf T}$, whose trace norm is $\|p(z)\|_2^2$ also for complex $z$. Thus
\eqref{W:eq:truncation-denominator}--\eqref{W:eq:diagonal-moment-remainder} imply on $I_i=1$
\begin{equation}
 \|v_i\|_\diamond\le CV_i,\quad
 \|\Delta_i(x)\|_\diamond\le C|u|V_i,\quad
 \|\Delta_i(x)-\ell_i(x)\|_\diamond\le Cu^2V_i.
 \label{W:eq:diagonal-moment-increment-bounds}
\end{equation}
The estimate holds for the full permitted increment.

\label{W:sec:conditional-row-bounds-for-diagonal-remainders}

To estimate $V_i^2$, we enlarge the event to one determined by the complementary minor and then integrate the off-diagonal row. \Needspace{6\baselineskip}

\begin{lemma}
\label{W:lem:truncation-integrated-weight-sums}
Uniformly in $t$ and in the deterministic starting vector, one has
\begin{equation}
 \sum_i\mathbb E[I_iV_i]\le Q_N{m_N},\qquad
 \sum_i\mathbb E[I_iV_i^2]\le Q_N,
 \qquad Q_N=N^{o(1)}.
 \label{W:eq:truncation-integrated-weight-sums}
\end{equation}
The same estimates hold after integration in $t$. They use only truncated off-diagonal fourth moments and diagonal second moments.
\end{lemma}
\begin{proof}
Fix a point $z$ of the contour. Delete row and column $i$ of $K$ and write $K^{(i)}$ for the minor, $h_i=K_{\{i\}^c,i}$ for its independent off-diagonal column, and
\[
 q=R_{K^{(i)}}(z)w^{(i)},\qquad r=R_K(z)w.
\]
Order the coordinate $i$ first and write
\[
 zI-K=\begin{pmatrix}z&-h_i^*\\-h_i&zI-K^{(i)}\end{pmatrix}.
\]
Solving $(zI-K)r=w$ in the lower block yields $r^{(i)}=R_{K^{(i)}}w^{(i)}+R_{K^{(i)}}h_i r_i=q+R_{K^{(i)}}h_i r_i$. The upper block then gives
$(z-h_i^*R_{K^{(i)}}h_i)r_i=w_i+h_i^*q$. The same equations with right side $e_i$ give $R_{K,ii}=(z-h_i^*R_{K^{(i)}}h_i)^{-1}$ and $R_{K,ji}=R_{K,ii}(R_{K^{(i)}}h_i)_j$. Therefore
\begin{equation}
 r_i=R_{K,ii}(z)(w_i+h_i^*q),\qquad
 q_j=r_j-\frac{R_{K,ji}(z)}{R_{K,ii}(z)}r_i\quad(j\ne i).
 \label{W:eq:truncation-minor-identities}
\end{equation}
On $I_i=1$, Lemma~\ref{W:lem:resolvent-weight-input} supplies $|R_{K,ii}|\in[c,C]$, bounded coordinate entries, $\|R_K e_i\|_2^2\le\rho N^{o(1)}$, and $\|r\|_2^2\le\rho N^{o(1)}$,
$\|r\|_\infty\le N^{o(1)}$. The second identity in \eqref{W:eq:truncation-minor-identities} gives
\[
 \begin{aligned}
 \|q\|_2^2&\le2\|r\|_2^2
                   +2c^{-2}|r_i|^2\|R_Ke_i\|_2^2,\\
 \|q\|_\infty&\le\|r\|_\infty
                    +c^{-1}|r_i|\max_{j\ne i}|R_{K,ji}|.
 \end{aligned}
\]
Substituting the stated bounds on the right proves
\begin{equation}
 \|q\|_2^2\le c_N\rho,\qquad
 \|q\|_\infty\le c_N,\qquad c_N=N^{o(1)}.
 \label{W:eq:truncation-minor-bounds}
\end{equation}
Define the event using only the minor by
\[
 \begin{split}
 \mathcal C_i(z)=\{K^{(i)}:\;&z\notin\sigma(K^{(i)}),\quad
     \|R_{K^{(i)}}(z)w^{(i)}\|_2^2\le c_N\rho,\\
 &\|R_{K^{(i)}}(z)w^{(i)}\|_\infty\le c_N\}.
 \end{split}
\]
At a nonreal $z$, the inverse exists for every Hermitian minor. At either real point of $\mathcal C$, the event $I_i=1$ places the spectrum of $K$ strictly between these two
points. The smallest and largest eigenvalues of a principal minor lie between those of $K$, by the variational characterization of eigenvalues. Thus $z\notin\sigma(K^{(i)})$ on
$I_i=1$ at these two real points as well. The preceding two norm calculations prove $\{I_i=1\}\subset\mathcal C_i(z)$.

For the conditional integral define $q^{\rm ext}(K^{(i)},z)=R_{K^{(i)}}(z)w^{(i)}$ on $\mathcal C_i(z)$ and zero on its complement. This is a Borel function of the minor: matrix
inversion is continuous on the open set of invertible matrices, and the two norm constraints are Borel conditions. It agrees with $q$ wherever the left side of the next inequality
is nonzero. Hence the conditional integral is well defined even at the real spectral parameters; all uses of $q$ below include $\mathbf1_{\mathcal C_i(z)}$.

Use $|R_{K,ii}|\le C$ on $I_i=1$ and $I_i\le\mathbf1_{\mathcal C_i(z)}$. Conditional expectation in $h_i$, with $K^{(i)}$ fixed, gives for $p=2,4$
\begin{equation}
 \mathbb E[I_i|r_i|^p]
 \le C_p\mathbb E_{K^{(i)}}\mathbf1_{\mathcal C_i(z)}
              \mathbb E_{h_i}|w_i+h_i^*q|^p.
 \label{W:eq:truncation-conditioning}
\end{equation}

Let $\widetilde h_i=h_i-\mathbb Eh_i$. For fixed $q$, put $Y_j=\overline{\widetilde h_{i,j}}q_j$. These are independent centered complex variables. Independence gives
$\mathbb E|\sum_jY_j|^2=\sum_j\mathbb E|Y_j|^2$ and the exact identity
\begin{equation}
 \mathbb E\left|\sum_jY_j\right|^4 ={}2\left(\sum_j\mathbb E|Y_j|^2\right)^2 +\left|\sum_j\mathbb EY_j^2\right|^2 +\sum_j\left\{\mathbb E|Y_j|^4 -2(\mathbb E|Y_j|^2)^2-|\mathbb EY_j^2|^2\right\}.
 \label{W:eq:truncation-fourth-identity}
\end{equation}
In the expansion, centering removes every assignment with a singleton index. The three pairings contribute $2(\sum_j\mathbb E|Y_j|^2)^2+|\sum_j\mathbb EY_j^2|^2$; the last sum in
\eqref{W:eq:truncation-fourth-identity} corrects the terms with all four indices equal. Figure~\ref{W:fig:moment-complex-fourth} records these four contributions, with
$v=\sum_j\mathbb E|Y_j|^2$, $s=\sum_j\mathbb EY_j^2$ and $\kappa_{4,j}=\mathbb E|Y_j|^4-2(\mathbb E|Y_j|^2)^2-|\mathbb EY_j^2|^2$. \begingroup
\setlength{\intextsep}{1.2\baselineskip}
\begin{figure}[H]
\centering
\begin{tikzpicture}[x=1cm,y=1cm,every node/.style={mg text}]
 \foreach \shift/\tag in {0/a,3.6/b,7.2/c,10.8/d}{
  \begin{scope}[xshift=\shift cm]
   \node[mg occurrence,label=above:{$1:Y$}] (\tag1) at (0,0) {};
   \node[mg occurrence,label=above:{$2:\overline Y$}] (\tag2) at (1.8,0) {};
   \node[mg occurrence,label=below:{$3:Y$}] (\tag3) at (0,-1.6) {};
   \node[mg occurrence,label=below:{$4:\overline Y$}] (\tag4) at (1.8,-1.6) {};
  \end{scope}
 }
 \foreach \a/\b in {a1/a2,a3/a4,b1/b4,b2/b3,c1/c3,c2/c4}{\draw[mg incidence] (\a)--(\b);}
 \node[mg block] (block) at (11.7,-.8) {};
 \foreach \a in {d1,d2,d3,d4}{\draw[mg incidence] (block)--(\a);}
 \node at (.9,.95) {$\{12\}\{34\}$};
 \node at (4.5,.95) {$\{14\}\{23\}$};
 \node at (8.1,.95) {$\{13\}\{24\}$};
 \node at (11.7,.95) {$\{1234\}$};
 \node at (.9,-2.55) {$v^2$};
 \node at (4.5,-2.55) {$v^2$};
 \node at (8.1,-2.55) {$|s|^2$};
 \node at (11.7,-2.55) {$\sum_j\kappa_{4,j}$};
\end{tikzpicture}
\caption{The four terms of the complex fourth moment identity \eqref{W:eq:truncation-fourth-identity}. The numbered dots are the positions in $(Y,\overline Y,Y,\overline Y)$. The three pair partitions contribute $v^2$, $v^2$, and $|s|^2$, respectively, and the diamond represents the fourth cumulant. In particular, the pairing $\{13\}\{24\}$ accounts for the pseudocovariance $s=\mathbb EY^2$. As in \eqref{W:eq:truncation-conditioning}, the minor is fixed before computing these expectations over the row entries.}
\label{W:fig:moment-complex-fourth}
\end{figure}
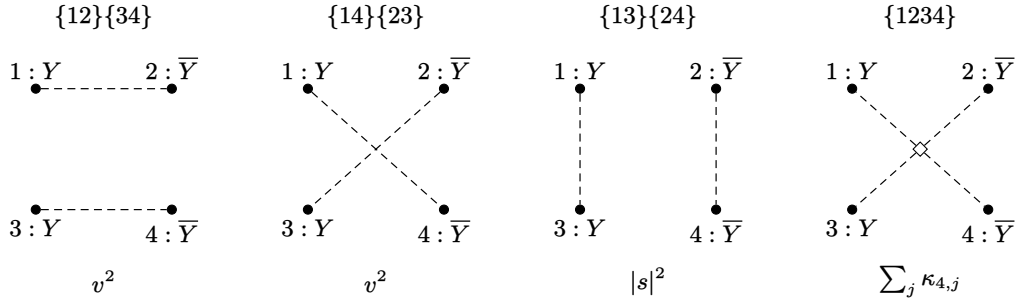
\endgroup

This expansion uses independence between distinct row entries and includes the terms $\mathbb EY_j^2$. Moreover $|\mathbb EY_j^2|\le\mathbb E|Y_j|^2$ and the two subtracted sums in
\eqref{W:eq:truncation-fourth-identity} are nonnegative. Consequently, with $v_j=\mathbb E|Y_j|^2$,
\[
 \mathbb E\Bigl|\sum_jY_j\Bigr|^4
  \le3\Bigl(\sum_jv_j\Bigr)^2+\sum_j\mathbb E|Y_j|^4.
\]
The row moment bounds and the deterministic vector $q$ give
\[
 \sum_jv_j\le\frac C N\sum_j|q_j|^2,\qquad
 \sum_j\mathbb E|Y_j|^4
   \le\frac{CF_N}{N^2}\sum_j|q_j|^4
   \le\frac{CF_N}{N^2}\|q\|_\infty^2\|q\|_2^2.
\]
These inequalities prove
\begin{equation}
 \mathbb E|\widetilde h_i^*q|^2 \le CN^{-1}\|q\|_2^2, \qquad \mathbb E|\widetilde h_i^*q|^4 \le CN^{-2}\|q\|_2^4 +CF_NN^{-2}\|q\|_\infty^2\|q\|_2^2.
 \label{W:eq:truncation-row-form}
\end{equation}
The row mean obeys
\[
 |(\mathbb Eh_i)^*q|
 \le\|\mathbb Eh_i\|_2\|q\|_2
 \le N^{-3/2+o(1)}\sqrt\rho=N^{-4/3+o(1)}.
\]
For $p=2,4$ use $|w_i+h_i^*q|^p\le3^{p-1}(|w_i|^p+|\widetilde h_i^*q|^p +|(\mathbb Eh_i)^*q|^p)$. On the minor event, \eqref{W:eq:truncation-row-form} now gives
\[
 \begin{aligned}
 \mathbb E_{h_i}|w_i+h_i^*q|^2
       &\le Q_N(|w_i|^2+\rho/N),\\
 \mathbb E_{h_i}|w_i+h_i^*q|^4
       &\le Q_N(|w_i|^4+\rho^2/N^2+F_N\rho/N^2)
        \le Q_N'(|w_i|^4+(\rho/N)^2).
 \end{aligned}
\]
For the mean term the preceding bound gives $|\mathbb Eh_i^*q|^2\le\rho N^{-3+o(1)}$ and $|\mathbb Eh_i^*q|^4\le\rho^2N^{-6+o(1)}$. These are bounded respectively by
$(\rho/N)N^{o(1)}$ and $(\rho/N)^2N^{o(1)}$. Moreover,
\[
 \frac{F_N\rho}{N^2}\le F_N\frac{\rho^2}{N^2}\qquad(\rho\ge1).
\]
Thus $Q_N'$ absorbs $F_N$ and the previously fixed subpolynomial constants. Substitute these estimates into \eqref{W:eq:truncation-conditioning}. At $\bar z$ use
$q^{\,-}=R_{K^{(i)}}(\bar z)w^{(i)}$ in the same two block equations. It satisfies the two bounds \eqref{W:eq:truncation-minor-bounds} with the same constants, and hence the two
displayed conditional moment inequalities with $q$ replaced by $q^{\,-}$. Finally
\[
 (x+y+z)^p\le3^{p-1}(x^p+y^p+z^p),\qquad x,y,z\ge0,
\]
applied to the two absolute resolvent coordinates and $\sqrt{q_N\rho/N}$ in $a_i$, gives
\begin{equation}
 \mathbb E[I_i a_i(z)^2]\le Q_N(|w_i|^2+\rho/N), \qquad \mathbb E[I_i a_i(z)^4]\le Q_N(|w_i|^4+(\rho/N)^2).
 \label{W:eq:truncation-row-moments}
\end{equation}

Let $J_N=\nu_\rho(\mathcal C)=N^{o(1)}$ and $b_i=(|w_i|^2+\rho/N)^{1/2}$. By \eqref{W:eq:truncation-row-moments},
\[
 \mathbb E[I_i a_i(z)^2]\le Q_Nb_i^2,\qquad
 \mathbb E[I_i a_i(z)^4]\le C Q_Nb_i^4.
\]
The same indicator appears at both contour points, so Cauchy--Schwarz gives
\[
 \mathbb E[I_i a_i(z)^2a_i(z')^2]
 \le\bigl(\mathbb E[I_i a_i(z)^4]\bigr)^{1/2}
     \bigl(\mathbb E[I_i a_i(z')^4]\bigr)^{1/2}
 \le C Q_Nb_i^4.
\]
Tonelli's theorem now yields
\begin{align*}
 \sum_i\mathbb E[I_iV_i]
 &\le Q_NJ_N\sum_i b_i^2=Q_NJ_N(1+\rho),\\
 \sum_i\mathbb E[I_iV_i^2]
 &\le C Q_NJ_N^2\sum_i b_i^4
 \le 2C Q_NJ_N^2(1+\rho^2/N).
\end{align*}
Here $\sum_i|w_i|^2=1$, $\sum_i|w_i|^4\le1$, and $\rho\asymp m_N$, $\rho^2/N=O(N^{-1/3})$. Increasing $Q_N$ to $C Q_N(1+J_N)^2$ gives
\eqref{W:eq:truncation-integrated-weight-sums}, since $J_N=N^{o(1)}$. The bounds are uniform at each $t$, and integration over an interval of length one preserves them.
\end{proof}

\begin{proof}[Proof of Theorem~\ref{W:thm:truncation-diagonal}]
With $s=s(K,w)$, define
\[
 \mathcal R_i(x)=A_N(s+\Delta_i(x))-A_N(s)-D_{s}A_N(s)[\ell_i(x)].
\]
Applying the fundamental theorem of calculus twice to $\theta\mapsto A_N(s+\theta\Delta_i(x))$ gives
\begin{equation}
 \mathcal R_i(x)={}D_{s}A_N(s)[\Delta_i(x)-\ell_i(x)] +\int_0^1(1-\theta) D_{s}^2A_N(s+\theta\Delta_i(x)) [\Delta_i(x),\Delta_i(x)]\,\mathrm d\theta.
 \label{W:eq:truncation-remainder-identity}
\end{equation}
The derivatives in \eqref{W:eq:truncation-remainder-identity} are with respect to the real coordinates $s_j=\mathcal L(p_j)$, with $s_0=1$ fixed. The vector $s+\theta\Delta_i(x)$
belongs to this domain for every $0\le\theta\le1$, and the globally extended $A_N$ satisfies \eqref{W:eq:replacement-test-class} there. On $I_i=1$, the two terms are bounded
respectively by
\[
 C\Lambda_Nm_Nu^2V_i,\qquad
 \int_0^1(1-\theta)\Lambda_Nm_N^2(C|u|V_i)^2\,\mathrm d\theta
     \le C\Lambda_Nm_N^2u^2V_i^2.
\]
Since $u^2=x^2/N$, this gives, for every permitted replacement value,
\begin{equation}
 |\mathcal R_i(x)|\le C\Lambda_N\frac{x^2}{N}
                         ({m_N}V_i+{m_N}^2V_i^2).
 \label{W:eq:truncation-quadratic-remainder}
\end{equation}
The bound is quadratic in the replacement entry throughout its support; the derivatives include the cutoff $\chi_N$.

Condition on $\mathscr F_i=\sigma(K_i(t))$. The indicator $I_i$ and the weight $V_i$ are then fixed, while the replacement variable has the deterministic measure $\mu_i+\nu_i$.
Thus
\[
 I_i\int|\mathcal R_i(x)|(\mu_i+\nu_i)(\mathrm dx)
 \le C\Lambda_NM_{2,N}N^{-1}I_i(m_NV_i+m_N^2V_i^2).
\]
The right side is nonnegative and integrable by Lemma~\ref{W:lem:truncation-integrated-weight-sums}. Tonelli's theorem permits the sum over $i$, the expectation over the remaining
entries, and integration in $t$ to be taken in this order. That lemma gives
\begin{align}
 \int_0^1\sum_i\mathbb E I_i\int|\mathcal R_i(x)|
                    (\mu_i+\nu_i)(\mathrm dx)\,\mathrm dt
 &\le C\Lambda_NM_{2,N}N^{-1}({m_N}Q_N{m_N}+{m_N}^2Q_N)\nonumber\\
 &\le C\Lambda_NQ_NM_{2,N}N^{-1/3}.
 \label{W:eq:truncation-remainder-summed}
\end{align}

The constant term integrates to zero. Since $I_i$ and $v_i$ are measurable with respect to the remaining entries, the linear contribution is exactly
\begin{equation}
 I_i\int D_{s}A_N(s)[\ell_i(x)](\mu_i-\nu_i)(\mathrm dx)
 =I_iN^{-1/2}D_{s}A_N(s)[v_i]\int x(\mu_i-\nu_i)(\mathrm dx).
 \label{W:eq:truncation-signed-mean}
\end{equation}
The signed first moment is a deterministic scalar of modulus at most $b_N$. Since $\|v_i\|_\diamond\le CV_i$ and $|DA_N(s)[v_i]|\le\Lambda_Nm_N\|v_i\|_\diamond$ on $I_i=1$, its
summed absolute expectation is at most
\begin{equation}
 C\Lambda_N{m_N}N^{-1/2}b_N
       \int_0^1\sum_i\mathbb E[I_iV_i]\,\mathrm dt
 \le C\Lambda_NQ_Nb_NN^{1/6}.
 \label{W:eq:truncation-mean}
\end{equation}

Lemma~\ref{W:lem:entry-distribution-interpolation}, applied with only the diagonal laws varying, identifies the sum in \eqref{W:eq:truncation-diagonal-contribution} with the
difference of the two expectations. If all entry laws vary, it is the diagonal part of that identity.

Decompose the expectation over the matrix with the selected entry set to zero according to the disjoint events $\{I_i=1\}$, $\mathcal B_i\setminus\mathcal A_{i,\rho}$, and
$\mathcal B_i^c$. On the middle event the original observable is zero for every permitted entry value, and its derivative at zero is zero, so its signed integral and its remainder
of the expansion both vanish. On the last event use the original integral, before any Taylor expansion:
\[
 \left|\int A_N(s(K_i+N^{-1/2}xE_{ii},w))
                         (\mu_i-\nu_i)(\mathrm dx)\right|
 \le\int1\,\mathrm d(\mu_i+\nu_i)=2.
\]
This bound uses no diagonal moment. Since there are $N$ diagonal blocks and $\sup_{i,t}\mathbb P(\mathcal B_i^c)\le N^{-D}$, the contribution of these exceptional events is at most
$2N^{1-D}$. Combining this with \eqref{W:eq:truncation-remainder-summed} and \eqref{W:eq:truncation-mean} proves \eqref{W:eq:diagonal-comparison-bound}.

At an arbitrary Hermitian $K$, use the functions $f_{jk}$ from Lemma~\ref{W:lem:global-first-two-entry-derivatives}. Their derivatives satisfy
$M_N^*=\max_{j,k<\rho}\|f_{jk}'\|_\infty\le N^{C_0}$, by Proposition~\ref{W:prop:gram-bounds} and the product rule. Lemma~\ref{W:lem:smooth-functional-calculus} gives
$\|D_Hf_{jk}(K)[E_{ii}]\|_{\mathrm{HS}}\le M_N^*$. Thus the functional $v(f)=D_H[w^*f(h_d(H))w]_{H=K}[E_{ii}]$ has norm $\|\cdot\|_\diamond$ at most $\sqrt\rho M_N^*$, and
\eqref{W:eq:replacement-test-class} implies
\begin{equation}
 |D_H(A_N\circ s)(K,w)[E_{ii}]|
 \le\Lambda_N{m_N}\sqrt\rho\,M_N^*\le N^{C_2}.
 \label{W:eq:global-diagonal-entry-derivative}
\end{equation}
This is a deterministic bound at every Hermitian matrix.

On the complement of the regularity event the first-order remainder of the expansion in the entry is therefore bounded by $2+N^{C_2-1/2}|x|$. Cauchy--Schwarz for the finite
positive measure $\mu_i+\nu_i$ gives $\int|x|(\mu_i+\nu_i)(\mathrm dx)\le\sqrt{2M_{2,N}}$. Measurability with respect to the remaining entries makes the sum of the integrals over
diagonal replacement values on the exceptional event at most
\begin{equation}
 CN^{1-D}\bigl(1+N^{C_2-1/2}\sqrt{M_{2,N}}\bigr).
 \label{W:eq:truncation-bad-remainder}
\end{equation}
For bounded $M_{2,N}$ and $A>0$, choose $D>A+C_2+2$. Both powers in \eqref{W:eq:truncation-bad-remainder} are then $O(N^{-A})$. With \eqref{W:eq:truncation-remainder-summed}, it
proves the theorem's bound for the absolute remainder.
\end{proof}

For the original diagonals and their truncations, use $\Theta_{\mathrm d}$ from \eqref{W:eq:diagonal-tail-function}. The mean bound is \eqref{W:eq:diagonal-mean-tail}, evaluated at
$L_N$. For comparison with any centered bounded standardized entry distribution with bounded second moment, $M_{2,N}=O(1)$ and $b_N\le e_N^{-1}N^{-1/2}\Theta_{\rm d}(L_N)$.
Theorem~\ref{W:thm:truncation-diagonal} consequently gives
\begin{equation}
 |\mathcal E_{{\rm diag},N}(A_N;w)|
 \le C\Lambda_NQ_NN^{-1/3}
       [M_{2,N}+e_N^{-1}\Theta_{\rm d}(L_N)]+CN^{1-D}
 =N^{-1/3+o(1)}.
 \label{W:eq:comparison-of-truncated-diagonals}
\end{equation}
The factor $N^{-1/3}$ gives convergence because $e_N^{-1}$, $\Lambda_N$ and $Q_N$ are subpolynomial and $\Theta_{\rm d}$ is bounded. The displayed estimate gives its dependence on
the modulus of uniform integrability. Symmetric standardized signs can be used, including when some original diagonal variances are zero.

The off-diagonal row moments required above follow from Lemma~\ref{W:lem:offdiag-tail-moments} and the centered truncated moment estimates in
Section~\ref{W:sec:moments-of-truncated-entries}.

The scale of Lemma~\ref{W:lem:admissible-scale} satisfies
\begin{equation}
 e_N^{-2}\Theta_{\mathrm d}(L_N)\longrightarrow0,\qquad
 e_N^{-4}\Theta(L_N)\longrightarrow0.
 \label{W:eq:entry-truncation-threshold-conditions}
\end{equation}
Its diagonal union bound is
\begin{equation}
 \mathbb P(\exists i:\widehat d_{N,i}\ne d_{N,i})
 \le e_N^{-2}\Theta_{\mathrm d}(L_N)=o(1).
 \label{W:eq:truncation-equality}
\end{equation}

For a bounded observable the expectation error caused by \eqref{W:eq:truncation-equality} is at most twice that probability. Under each of the two ensemble laws this contribution
is $o(1)$. Gaussian diagonals of uniformly bounded variance can be truncated symmetrically at $L_N$; they remain centered and their second moments do not increase. The union
probability is at most $CN\exp(-cL_N^2)$, smaller than every fixed inverse power of the matrix size because $L_N^2=N^{1-o(1)}$. Their mean term is zero. Truncation of the
off-diagonal Gaussian entries elsewhere in the comparison is handled by the same convention for passing to the original arrays used in
Section~\ref{W:sec:convergence-of-the-integrated-potential}.

For a function $F_N$ multiplied by its cutoff of the spectral moments at the starting vector and the function $A_N=\chi_NF_N$, one has pointwise
\begin{equation}
 |F_N-A_N|\le D_N^{\rm coef}=\Gamma_\rho(1-\chi_N).
 \label{W:eq:truncation-localization-loss}
\end{equation}
Apply Theorem~\ref{W:thm:truncation-diagonal} to $D_N^{\rm coef}$ using its global derivative bounds and Lemma~\ref{W:lem:truncation-integrated-weight-sums}. Together with the
off-diagonal comparison in Section~\ref{W:sec:taylor-remainders-and-bounded-entries}, this compares its truncated expectation with its expectation under the law of $H_N^{\rm b}$.
The bound for $H_N^{\rm b}$ then gives the corresponding bound for the truncated matrix.

For the matrices being compared $H,\widetilde H$ in that argument,
\begin{align}
 |\mathbb EF_N(H,w)-\mathbb EF_N(\widetilde H,w)|
 \le{}&|\mathbb EA_N(H,w)-\mathbb EA_N(\widetilde H,w)|\nonumber\\
 &+\mathbb ED_N^{\rm coef}(H,w)
       +\mathbb ED_N^{\rm coef}(\widetilde H,w).
 \label{W:eq:truncation-remove-cutoff}
\end{align}
Apply the estimates to the smooth functions of the moduli of continuity for both reflected processes. Equation~\eqref{W:eq:truncation-equality} bounds the error from replacing
$H_N^{\rm tr}$ by $H_N$. Section~\ref{W:sec:resolvents-on-exceptional-support} bounds the error from replacing $h_d(H_N)$ by $H_N$. Add both errors to the difference of
expectations. This proves the diagonal comparison for the exponentials and the smooth functions of the modulus of continuity, with their cutoffs, before spectral modification. The
absolute remainder bound applies to these functions multiplied by $\chi_N$.

For Theorem~\ref{W:thm:moment-weights-mean-resolvent}, let $\mathcal D_N$ be its nonreal spectral domain, including the specified dyadic vertical parameters. The assumed event
bounds all coordinates of $R_H(z)$ by $C_0$ uniformly for $z\in\mathcal D_N$. Write $H=K+u_0E_{ii}$, with $|u_0|\le Ce_N$. The rank-one identity \eqref{W:eq:truncation-SM}, applied
in the reverse direction, gives
\[
 R_K=R_H-\frac{u_0}{1+u_0R_{H,ii}}R_HE_{ii}R_H,\qquad
 \max_{a,b}|(R_K-R_H)_{ab}|\le2C_0^2|u_0|.
\]
For large $N$, both this bound and the denominator bound $|1+u_0R_{H,ii}|\ge1/2$ hold throughout $\mathcal D_N$. In particular, this bound on the entries of $R_H$ implies
\[
 \mathcal E_i=\left\{\sup_{z\in\mathcal D_N}\max_{a,b}|R_{K,ab}(z)|\le2C_0\right\}.
\]
The set $\mathcal E_i$ is measurable with respect to the remaining entries: a countable dense subset of $\mathcal D_N$ gives the same supremum by continuity of the resolvent. On
this set, $|1-uR_{K,ii}|\ge1/2$ for every $|u|\le Ce_N$, and another application of \eqref{W:eq:truncation-SM} gives
\[
 \sup_{z\in\mathcal D_N}\max_{a,b}|R_{K+uE_{ii},ab}(z)|
 \le2C_0+8C_0^2|u|\le4C_0.
\]
Thus one event controls the entire interval of replacement values and all the stated spectral parameters. For each fixed interpolation parameter, $\mathbb P(\mathcal E_i^c)$ is
bounded by the probability that the assumed entrywise bound for $R_H$ fails, uniformly in that parameter. For a specified family of at most $N^a$ matrices and all $N$ diagonal
positions, choose failure exponent $D+a+2$ there; the union has probability at most $N^{a+1}N^{-D-a-2}\le N^{-D}$. The dyadic vertical parameters are already contained in the
uniform domain.

To check the second moment step in that theorem, let its diagonal variable $u$ be centered and $f_K(u)=(R_{K+uE_{ii}}v)_i$. Then $f_K(u)=f_K(0)/(1-uR_{K,ii})$, so, conditionally on
the matrix with the selected entry set to zero,
\[
 |\mathbb E_u[u f_K(u)]|
 =\left|f_K(0)R_{K,ii}\mathbb E_u
                \frac{u^2}{1-uR_{K,ii}}\right|
 \le C\mathbb Eu^2\,|f_K(0)|.
\]
Also $|f_K(0)|\le C|f_K(u)|$ for every allowed $u$. This is the $CN^{-1}\mathbb E|f_K(u)|$ diagonal term used in the resolvent equation of that theorem; only the second moment is
used. Write $\mathcal B_i$ for this regularity event. On $\mathcal B_i^c$, $|f_K(u)|\le\eta^{-1}$ and $|u|\le Ce_N$ give
\[
 \mathbb E[\mathbf1_{\mathcal B_i^c}|u f_K(u)|]
 \le Ce_N\eta^{-1}\mathbb P(\mathcal B_i^c)
 \le C N^{1-D},
\]
using $e_N\le1$, $\eta\ge N^{-1}$ and $\mathbb P(\mathcal B_i^c)\le N^{-D}$. For an error $O(N^{-A})$ after summing $N$ diagonal entries, choose $D>A+2$.

Thus the intermediate statement uses $\mathbb EH_{ii}^2\le C/N$ and $|H_{ii}|\le Ce_N$ with $e_N\to0$. These bounds apply to the truncated matrices in the intermediate resolvent
estimate. Centering a truncated diagonal changes the bound by at most another fixed factor.

This proves the diagonal estimate with only the second moment of the replacement variable. The corresponding calculation for $\operatorname{Tr}\Theta(H)$ is separate and is proved
in Theorem~\ref{W:thm:localization-taylor-remainder}.

\subsubsection{Taylor remainders for off-diagonal entries}
\label{W:sec:taylor-expansion-in-an-entry}
\label{W:sec:resolvent-fourth-chain}

\label{W:sec:taylor-remainders-and-bounded-entries}
\label{W:sec:comparison-with-bounded-entries}
The first three coordinate moment differences control the polynomial of degree three in a selected off-diagonal entry. We estimate the difference between this polynomial and the
complete function $A_N(s(H,w))$ using the fourth absolute moment.

Here $H_N^{\rm tr}$ is the uncentered truncated array. The matrix $H_N^{\rm b}$ has independent bounded standardized off-diagonal blocks matching the first three complete tensors
of the centered, covariance normalized truncation, as in Lemma~\ref{W:lem:cubic-defect}. Its standardized diagonal entries are independent symmetric signs, so their matrix
variances are $N^{-1}$. For $A_N$ in \eqref{W:eq:replacement-test-class}, \eqref{W:eq:tightness-domination} allows the derivative bounds to be used with $\|\cdot\|_{\mathscr D}$ in
place of $\|\cdot\|_\diamond$. With the exact increments of Lemma~\ref{W:lem:comparison-weights}, the complete polynomial is

\begin{equation}
 \begin{split}
 P_{3,\alpha}(x)={}&A_N(s)+D_{s}A_N(s)[\ell+q+c]
       +\tfrac12D_{s}^2A_N(s)[\ell,\ell]\\
       &+D_{s}^2A_N(s)[\ell,q]
       +\tfrac16D_{s}^3A_N(s)[\ell,\ell,\ell].
 \end{split}
 \label{W:eq:cubic-entry-polynomial}
\end{equation}
The base vector is $s=s(K_\alpha,w)$, with the selected entry set to zero.

\Needspace{6\baselineskip}

\begin{lemma}
\label{W:lem:comparison-remainder}
For each permitted change of an off-diagonal entry on the regularity event,
\begin{equation}
 \begin{split}
 |A_N(s+\Delta s)-P_{3,\alpha}(x)|
 \le C\Lambda_Nt_x^4\big(&{m_N}V_\alpha+{m_N}^2U_\alpha V_\alpha
           +{m_N}^2V_\alpha^2\\
           &+{m_N}^3U_\alpha^2V_\alpha+{m_N}^4U_\alpha^4\big).
 \end{split}
 \label{W:eq:comparison-remainder}
\end{equation}
\end{lemma}
\begin{proof}
With $\alpha$ fixed, put $h=\Delta s-\ell=q+c+r$ and write $U=U_\alpha$, $V=V_\alpha$. Choose a fixed constant $C_0$ so that
\[
 m_N\|\ell\|_{\mathscr D}\le a=C_0m_Nt_xU,\quad
 m_N\|h\|_{\mathscr D},\ m_N\|q\|_{\mathscr D}\le b=C_0m_Nt_x^2V,
\]
\[
 m_N\|c\|_{\mathscr D}\le c_0=C_0m_Nt_x^3V,\qquad
 m_N\|r\|_{\mathscr D}\le r_0=C_0m_Nt_x^4V.
\]
We have $c_0\le b$ and $r_0\le b$ for large $N$, since $t_x\le Ce_N<1$. All derivatives without an indicated argument below are evaluated at $s$. Taylor's formula along
$u\mapsto s+u\Delta s$ gives
\[
 A_N(s+\Delta s)=\sum_{j=0}^3\frac1{j!}D_s^jA_N(s)[(\Delta s)^j]
 +\frac1{3!}\int_0^1(1-u)^3D_s^4A_N(s+u\Delta s)[(\Delta s)^4] \,\mathrm du.
\]
Subtracting \eqref{W:eq:cubic-entry-polynomial} and using symmetry of the multilinear derivatives yields the exact identity
\begin{align*}
 A_N(s+\Delta s)-P_{3,\alpha}(x)
 ={}&DA_N[r]+D^2A_N[\ell,c+r]+\tfrac12D^2A_N[h,h]\\
 &+\tfrac12D^3A_N[\ell,\ell,h]
   +\tfrac12D^3A_N[\ell,h,h]+\tfrac16D^3A_N[h,h,h]\\
 &+\tfrac16\int_0^1(1-u)^3
                      D^4A_N(s+u\Delta s)[(\Delta s)^4] \,\mathrm du.
\end{align*}
In the order displayed, the absolute value of each term is bounded by $C\Lambda_N$ times the corresponding entry in the right column:
\[
 \begin{array}{c|c}
 DA_N[r]&r_0\\
 D^2A_N[\ell,c+r]&a(c_0+r_0)\\
 D^2A_N[h,h]/2&b^2/2\\
 D^3A_N[\ell,\ell,h]/2&a^2b/2\\
 D^3A_N[\ell,h,h]/2&ab^2/2\\
 D^3A_N[h,h,h]/6&b^3/6\\
 \displaystyle\frac16\int_0^1(1-u)^3D^4A_N(s+u\Delta s)[(\Delta s)^4]\,\mathrm du
           &(a+b)^4/24.
 \end{array}
\]
The last row uses $m_N\|\Delta s\|_{\mathscr D}\le a+b$ and $\int_0^1(1-u)^3\,\mathrm du=1/4$. Thus for $a+b\le1$ the absolute value is at most
\[
 C\Lambda_N\{r_0+a(c_0+r_0)+b^2+a^2b+ab^2+b^3+(a+b)^4\}.
\]
The integral contributes at most $C_4\Lambda_N(a+b)^4/24$. On this region,
\[
 ar_0\le r_0,\qquad ab^2+b^3\le2b^2,\qquad
 (a+b)^4\le8(a^4+b^4)\le8(a^4+b^2).
\]
Thus the required bound here is $C\Lambda_N(r_0+ac_0+b^2+a^2b+a^4)$.

For $a+b>1$, use $|A_N(s+\Delta s)|,|A_N(s)|\le1$ and the bounds on the coefficients of $P_{3,\alpha}$ to obtain
\[
 |A_N(s+\Delta s)-P_{3,\alpha}(x)|
 \le2+C\Lambda_N(a+b+c_0+a^2+ab+a^3).
\]
If $a\ge1/2$, then
\[
 1\le16a^4,\quad a\le8a^4,\quad a^2\le4a^4,\quad a^3\le2a^4,
 \quad b\le b^2+4a^4,\quad ab\le2a^4+b^2/2.
\]
If $a<1/2$, then $b>1/2$, and the individual inequalities are
\[
 1\le4b^2,\quad a\le2b^2,\quad a^2\le b^2,
 \quad a^3\le\tfrac12b^2,\quad b\le2b^2,\quad ab\le b^2.
\]
Since $c_0\le b$ and $\Lambda_N\ge1$, in either case the last display is at most $C\Lambda_N(a^4+b^2)$.

The definitions give
\[
 (r_0,ac_0,b^2,a^2b,a^4)
 \le C t_x^4(m_NV,m_N^2UV,m_N^2V^2,m_N^3U^2V,m_N^4U^4)
\]
componentwise. This proves \eqref{W:eq:comparison-remainder} for every permitted entry value.
\end{proof}

The five terms in the remainder estimate also appear in the homogeneous fourth-order coefficient. Fix a coordinate vector $x$ for one entry, vary the matrix as
$K+tN^{-1/2}\mathcal T_\alpha x$, and write, for this calculation only,
\[
 s(t)=s+t\ell+t^2q+t^3c+t^4d+o(t^4).
\]
Here $d$ is the homogeneous fourth degree coefficient of the moment curve; $r_\alpha$ denotes the full remainder of the finite expansion. In the composition formula of
Section~\ref{W:sec:notation-and-scales}, the partitions of four labels have sizes $4$, $3+1$, $2+2$, $2+1+1$, and $1+1+1+1$, with counts $1,4,3,6,1$. Since $s'(0)=\ell$,
$s''(0)=2q$, $s^{(3)}(0)=6c$, and $s^{(4)}(0)=24d$, division by $4!$ gives
\begin{equation}
 DA_N[d]+D^2A_N[\ell,c]+\tfrac12D^2A_N[q,q]
       +\tfrac12D^3A_N[\ell,\ell,q]
       +\tfrac1{24}D^4A_N[\ell,\ell,\ell,\ell].
 \label{W:eq:resolvent-fourth-coefficient}
\end{equation}
To see the estimates term by term, let $t_x=N^{-1/2}\|x\|_1$ and use the same $U_\alpha,V_\alpha$ as in Lemma~\ref{W:lem:comparison-weights}. Differentiating the finite resolvent
identity at zero gives
\[
 \|\ell\|_{\mathscr D}\le Ct_xU_\alpha,\quad
 \|q\|_{\mathscr D}\le Ct_x^2V_\alpha,\quad
 \|c\|_{\mathscr D}\le Ct_x^3V_\alpha,\quad
 \|d\|_{\mathscr D}\le Ct_x^4V_\alpha.
\]
For the last estimate, the product has four factors $V_\alpha(x)$, three bounded coordinate compressions of the resolvent, and the two exterior resolvent vectors; their contour
integral is bounded by $Ct_x^4V_\alpha$. Substitution into the five multilinear derivatives yields
\[
 \begin{aligned}
 |DA_N[d]|&\le C\Lambda_Nt_x^4m_NV_\alpha,\\
 |D^2A_N[\ell,c]|&\le C\Lambda_Nt_x^4m_N^2U_\alpha V_\alpha,\\
 \tfrac12|D^2A_N[q,q]|&\le C\Lambda_Nt_x^4m_N^2V_\alpha^2,\\
 \tfrac12|D^3A_N[\ell,\ell,q]|&\le C\Lambda_Nt_x^4m_N^3U_\alpha^2V_\alpha,\\
 \tfrac1{24}|D^4A_N[\ell,\ell,\ell,\ell]|&\le C\Lambda_Nt_x^4m_N^4U_\alpha^4.
 \end{aligned}
\]
Lemma~\ref{W:lem:comparison-remainder} gives the corresponding bound for the full difference
$A_N(s+\Delta s)-P_{3,\alpha}(x)$.

\label{W:sec:summation-of-fourth-order-remainder}
We estimate the two highest products of resolvent weights in \eqref{W:eq:summed-fourth-remainder-terms}. Write $m=m_N$ and $J_N=\nu_\rho(\mathcal C)$. Keep the moment constants
separate from $J_N$. Let $C_N\ge1$ be a common deterministic upper bound for $q_N^4$ and the constants in the second, fourth and eighth row moment estimates below. Those estimates
give $C_N=N^{o(1)}$. Choose
\[
 Q_N=100C_N^8(1+J_N)^4,\qquad
 \frac{\log Q_N}{\log N}
  =\frac{\log100+8\log C_N+4\log(1+J_N)}{\log N}\longrightarrow0.
\]
This $Q_N$ bounds the finite products of moment constants and contour masses in the following calculation. The conditional row estimate gives a deterministic bound, uniform over
the pair set to zero and over the contour,
\begin{equation}
 \mathbb E[\mathbf1_{\mathcal G_\alpha}a_i^8]
 \le Q_N\left(|w_i|^8+m^{-8}+\frac{m e_N^4}{N^2}\right),
 \qquad
 \mathcal G_\alpha=\mathcal B_\alpha\cap\mathcal A_{\alpha,\rho}.
 \label{W:eq:resolvent-eighth-moment-bound}
\end{equation}
For the conditional eighth moment calculation, fix the minor obtained by removing the row being integrated, and let $q$ be its resolvent coefficient vector. The regularity event is
contained in the event measurable with respect to the minor where $\|q\|_2^2\le\rho N^{o(1)}$ and $\|q\|_\infty\le N^{o(1)}$; the selected pair is kept equal to zero throughout. To
verify this conditioning step, put $r=R_K(z)w$. The block equations give
\[
 r_i=R_{K,ii}(w_i+h_i^*q),\qquad
 q_j=r_j-\frac{R_{K,ji}}{R_{K,ii}}r_i\quad(j\ne i).
\]
These formulas hold with the diagonal entry $K_{ii}$ included in the denominator of the Schur complement formula. On $I_\alpha=1$, the bounds $|R_{K,ii}|\le C$, $|R_{K,ii}|\ge c$,
the column norm bound, and the bounds for $r$ imply the stated bounds on $q$. Let $J_i$ be the indicator of those bounds on the minor vector. It is measurable with respect to the
minor and satisfies $I_\alpha\le J_i$. Pointwise,
\[
 I_\alpha|r_i|^8
   \le C^8J_i|w_i+h_i^*q|^8.
\]
Thus conditional expectation may integrate the entire row on the right, even though $I_\alpha$ on the left can depend on that row. The diagonal entry is fixed during this
integration over the off-diagonal row. This is the order eight form of \eqref{W:eq:truncation-conditioning}. For the centered independent row variables, apply
Lemma~\ref{W:lem:schur-rosenthal} with $p=8$ conditional on the minor. With $Y_j=\overline{h_j-\mathbb Eh_j}\,q_j$ and $A_r=\sum_j\mathbb E_h|Y_j|^r$, it gives
\[
 \Bigl\|\sum_jY_j\Bigr\|_{L^8(\mathbb E_h)}
 \le8C(A_2^{1/2}+A_8^{1/8}),\qquad
 \mathbb E_h\Bigl|\sum_jY_j\Bigr|^8
 \le2^7(8C)^8(A_2^4+A_8).
\]
This is the conditional row calculation leading to \eqref{W:eq:conditional-row-moment}, with the two moment sums still in their exact form. Its constant is independent of the
number of row entries and of the fixed vector $q$. The centered row entries have support of radius at most $Ce_N$, variance at most $C/N$, and fourth moment at most $F_N/N^2$.
Therefore
\[
 \begin{aligned}
 A_2&\le CN^{-1}\|q\|_2^2,\\
 A_8&\le C F_Ne_N^4N^{-2}\sum_j|q_j|^8
     \le C F_Ne_N^4N^{-2}\|q\|_\infty^6\|q\|_2^2,\\
 \mathbb E_h\Bigl|\sum_jY_j\Bigr|^8
     &\le N^{o(1)}\bigl((\rho/N)^4+\rho e_N^4/N^2\bigr).
 \end{aligned}
\]
This uses the support bound only through $\mathbb E|h_j-\mathbb Eh_j|^8 \le Ce_N^4\mathbb E|h_j-\mathbb Eh_j|^4$. For the mean vector, each coordinate is $N^{-2+o(1)}$, and hence
$|(\mathbb Eh)^*q|\le N^{-3/2+o(1)}\sqrt\rho$. Substitute the centered and mean bounds into the Schur complement formula for $(R(z)w)_i$, using
$(x_1+\cdots+x_l)^8\le l^7\sum_jx_j^8$ for nonnegative $x_j$. Apply it at both $z$ and $\bar z$ and include the deterministic term $\sqrt{q_N\rho/N}$. After integration of the
minor this gives
\[
 \mathbb E[I_\alpha a_i(z)^8]
   \le C_N\bigl(|w_i|^8+(\rho/N)^4+\rho e_N^4/N^2\bigr).
\]
Since $\rho\asymp m$ and $N=m^3$, this is \eqref{W:eq:resolvent-eighth-moment-bound}. The conditional row inequality is uniform in the deleted pair; its constant is included in
$C_N$.

Set $\widehat\nu_\rho=J_N^{-1}\nu_\rho$, a probability measure, and $d_i=|w_i|^8+m^{-8}+me_N^4/N^2$. For each fixed pair the definition of $U_{ab}$ gives
\[
 \begin{aligned}
 \mathbb E[\mathbf1_{\mathcal G_{ab}}U_{ab}^4]
 &=J_N^4\mathbb E\!\left[\mathbf1_{\mathcal G_{ab}}
              \left(\int a_a(z)a_b(z)\,\widehat\nu_\rho(\mathrm dz)\right)^4\right]\\
 &\le J_N^3\int\mathbb E[\mathbf1_{\mathcal G_{ab}}a_a(z)^4a_b(z)^4]
                                                \,\nu_\rho(\mathrm dz)\\
 &\le J_N^3\int
   (\mathbb E[\mathbf1_{\mathcal G_{ab}}a_a(z)^8])^{1/2}
   (\mathbb E[\mathbf1_{\mathcal G_{ab}}a_b(z)^8])^{1/2}
                                                \,\nu_\rho(\mathrm dz)\\
 &\le C_NJ_N^4\sqrt{d_ad_b}.
 \end{aligned}
\]
Jensen gives the first inequality, and Cauchy--Schwarz applied to $\mathbf1_{\mathcal G_{ab}}a_a^4$ and $\mathbf1_{\mathcal G_{ab}}a_b^4$ under the same distribution of the matrix
with the selected entry set to zero gives the second. The $d_i$ are deterministic and uniform over the pair. Therefore
\[
 \sum_{a<b}\sqrt{d_ad_b}
 =\frac12\left[\left(\sum_i\sqrt{d_i}\right)^2-\sum_i d_i\right]
 \le\frac12\left(\sum_i|w_i|^4+Nm^{-4}+\sqrt m\,e_N^2\right)^2.
\]
Since $C_NJ_N^4\le Q_N/100$ and $(x+y+z)^2\le3(x^2+y^2+z^2)$ for nonnegative $x,y,z$, we obtain
\begin{align}
 \sum_{a<b}\mathbb E[\mathbf1_{\mathcal G_{ab}}U_{ab}^4]
 &\le Q_N\left(\sum_i|w_i|^4+Nm^{-4}
                +N\sqrt{m/N^2}\,e_N^2\right)^2\notag\\
 &\le Q_N(1+me_N^4)=Q_NS_N.
 \label{W:eq:resolvent-U4-sum}
\end{align}
Here $\sum_i|w_i|^4\le1$ and $Nm^{-4}=m^{-1}$. The deterministic coordinate bounds are uniform over the pair set to zero, so they can be summed even though different pairs use
different matrices. For $V_{ab}$, Jensen's inequality and $(a_a^2+a_b^2)^2\le2a_a^4+2a_b^4$ give
\[
 \begin{aligned}
 \mathbb E[\mathbf1_{\mathcal G_{ab}}V_{ab}^2]
 &\le2J_N\int\mathbb E[\mathbf1_{\mathcal G_{ab}}
                               (a_a(z)^4+a_b(z)^4)]\,\nu_\rho(\mathrm dz)\\
 &\le2C_NJ_N^2(|w_a|^4+|w_b|^4+2m^{-4}).
 \end{aligned}
\]
Summing over the pairs, using $\sum_{a<b}(|w_a|^4+|w_b|^4)=(N-1)\sum_i|w_i|^4\le N-1$, yields
\[
 \sum_{a<b}\mathbb E[\mathbf1_{\mathcal G_{ab}}V_{ab}^2]
 \le2C_NJ_N^2(N+N^2m^{-4})
 \le4C_NJ_N^2m^3\le Q_Nm^3,
\]
because $N=m^3$ and $m\ge1$. For fixed $t$, let $\mathfrak m_t$ be the finite positive measure on the disjoint union of the base probability spaces, with measure
$\mathbf1_{\mathcal G_{ab}}\,\mathbb P_{K_{ab}(t)}$ on component $(a,b)$. Cauchy--Schwarz gives
\[
 \int U^2V\,\dd\mathfrak m_t
 \le\left(\int U^4\,\dd\mathfrak m_t\right)^{1/2}
      \left(\int V^2\,\dd\mathfrak m_t\right)^{1/2}.
\]
The bounds just proved are $3C_NJ_N^4S_N$ and $4C_NJ_N^2m^3$. Their geometric mean is at most $\sqrt{12}C_NJ_N^3m^{3/2}S_N^{1/2}\le Q_Nm^{3/2}S_N^{1/2}$, by the definition of
$Q_N$. Hence
\begin{equation}
 \sum_{a<b}\mathbb E[\mathbf1_{\mathcal G_{ab}}U_{ab}^2V_{ab}]
 \le Q_Nm^{3/2}S_N^{1/2}.
 \label{W:eq:resolvent-U2V-sum}
\end{equation}
The estimate is uniform in $t$, so integration over $[0,1]$ preserves it.

Conditional on $K_\alpha(t)$, $\int t_x^4(\mu_\alpha+\nu_\alpha)(\mathrm dx) =N^{-2}\int\|x\|_1^4(\mu_\alpha+\nu_\alpha)(\mathrm dx) \le N^{-2+o(1)}$. The following table includes
all five sums; in particular the factor $S_N=1+m_Ne_N^4$ is displayed together with its power of $m_N$.

For the five terms in \eqref{W:eq:comparison-remainder}, the bounds are
\begin{equation}
\begin{array}{c|c}
\text{term after the fourth moment of the entry distribution}&
        \text{summed bound, without }\Lambda_NN^{o(1)}\\\hline
{m_N}N^{-2}V&{m_N}^{-1}\\
{m_N}^2N^{-2}UV&{m_N}^{-2}\\
{m_N}^2N^{-2}V^2&{m_N}^{-1}\\
{m_N}^3N^{-2}U^2V&{m_N}^{-3/2}S_N^{1/2}\\
{m_N}^4N^{-2}U^4&{m_N}^{-2}S_N
\end{array}
 \label{W:eq:summed-fourth-remainder-terms}
\end{equation}
Since $S_N\le Cm_N$, the two terms involving $S_N$ satisfy
\[
 m_N^{-3/2}S_N^{1/2}\le C^{1/2}m_N^{-1},\qquad
 m_N^{-2}S_N\le Cm_N^{-1}.
\]
The other three bounds are $m_N^{-1}$, $m_N^{-2}$, and $m_N^{-1}$, and $m_N^{-2}\le m_N^{-1}$ for $N\ge1$. Thus all five are $O(m_N^{-1})$. To display the dependence on the fourth
moment, write $\mathcal M_{4,N}=\sup_\alpha\int\|x\|_1^4(\mu_\alpha+\nu_\alpha)(\mathrm dx)$. The sum is bounded by
\[
 C\Lambda_NQ_N\mathcal M_{4,N}
 \bigl[m_N^{-1}+m_N^{-2}+m_N^{-1}
          +m_N^{-3/2}S_N^{1/2}+m_N^{-2}S_N\bigr]
 \le C\Lambda_NQ_N\mathcal M_{4,N}m_N^{-1}.
\]
The bound $\mathcal M_{4,N}=N^{o(1)}$ follows from the truncated fourth moment estimate. The chosen $Q_N$ bounds the five contour masses $J_N,J_N^2,J_N^2,J_N^3,J_N^4$. If
$\mathfrak I_N^{\rm good}$ denotes the full integral against the sum of the two entry laws with the indicator measurable with respect to the remaining entries, this proves
\begin{equation}
 \mathfrak I_N^{\rm good}[|R_{4,\alpha}(A_N)|]
                 \le N^{-1/3+o(1)}.
 \label{W:eq:comparison-sumR}
\end{equation}

The first three coefficient tensors of the Taylor polynomial in the entry satisfy
\begin{align}
 \|C_{1,\alpha}\|&\le C\Lambda_N{m_N}N^{-1/2}U_\alpha,\nonumber\\
 \|C_{2,\alpha}\|&\le C\Lambda_NN^{-1}
                            ({m_N}V_\alpha+{m_N}^2U_\alpha^2),\nonumber\\
 \|C_{3,\alpha}\|&\le C\Lambda_NN^{-3/2}
          ({m_N}V_\alpha+{m_N}^2U_\alpha V_\alpha+{m_N}^3U_\alpha^3).
 \label{W:eq:comparison-coefficient-bounds}
\end{align}
To verify the third line, put $c_{rst}=D_x^3s_\alpha(0)[e_r,e_s,e_t]/6$. The coefficient is the ordered tensor in \eqref{W:eq:cumulant-chain-tensor}, evaluated at this deleted
matrix. The diagonal bounds on $q(x),c(x)$ imply the mixed bounds by real polarization. For a symmetric trilinear form $T$,
\[
 T(u,v,w)=\frac1{48}\sum_{\varepsilon\in\{-1,1\}^3}
 \varepsilon_1\varepsilon_2\varepsilon_3\,
 T(\varepsilon_1u+\varepsilon_2v+\varepsilon_3w,
   \varepsilon_1u+\varepsilon_2v+\varepsilon_3w,
   \varepsilon_1u+\varepsilon_2v+\varepsilon_3w).
\]
Indeed, expansion of the right side eliminates every monomial except those containing each of $u,v,w$ once; their total coefficient is $8\cdot3!=48$. Applied to $c$ with coordinate
vectors, this gives $\|c_{rst}\|_{\mathscr D}\le CN^{-3/2}V_\alpha$. The corresponding bilinear identity was proved above. Substitution of these bounds and
$\|\ell_r\|_{\mathscr D}\le CN^{-1/2}U_\alpha$ proves all three lines of \eqref{W:eq:comparison-coefficient-bounds}. Conditioning on $K_\alpha$, each tensor is fixed, so its
integral against the full difference of the entry distributions is exactly $C_{j,\alpha}:\Delta M_j$. This remains true with the indicator measurable with respect to the remaining
entries. For the cubic term, \eqref{W:eq:comparison-defects} gives $\|\Delta M_3\|_{1,\mathrm{ord}}\le N^{-1+o(1)}$. Consequently its three contributions are bounded, respectively,
by
\[
 \Lambda_NN^{-5/2+o(1)}m_N\cdot m_N^4,
 \qquad \Lambda_NN^{-5/2+o(1)}m_N^2\cdot m_N^2,
 \qquad \Lambda_NN^{-5/2+o(1)}m_N^3\cdot m_N.
\]
Using $N=m_N^3$ gives $m_N^{-5/2},m_N^{-7/2},m_N^{-7/2}$, up to the common factor $\Lambda_NN^{o(1)}$. Together with Theorem~\ref{W:thm:replacement-low-order}, the three Taylor
orders therefore have the bounds
\begin{equation}
 {m_N}^{-1},\qquad {m_N}^{-1}+{m_N}^{-2},\qquad
                         {m_N}^{-5/2}+2{m_N}^{-7/2}.
 \label{W:eq:summed-first-three-moment-errors}
\end{equation}
For the first two orders, Theorem~\ref{W:thm:replacement-low-order} also gives the dependence on the tail function.

For diagonal replacements, Theorem~\ref{W:thm:truncation-diagonal} gives
\begin{equation}
 |A_N(s+\Delta s)-A_N(s)-D_{s}A_N(s)[\ell]|
 \le C\Lambda_NN^{-1}|x|^2({m_N}V_\alpha+{m_N}^2V_\alpha^2).
 \label{W:eq:comparison-diagonal-remainder}
\end{equation}
Its hypotheses are the global first and second derivative bounds and \eqref{W:eq:comparison-diag-weights}. Integration of the remainder and the mean term gives $N^{-1/3+o(1)}$,
including the exceptional event.

Let $\mu_\alpha^+=\mu_\alpha+\nu_\alpha$ be the sum of the two entry laws, and let $M_{j,\alpha}^+=\int\|x\|_1^j\mu_\alpha^+(\mathrm dx)$. On the off-diagonal these moments are
$N^{o(1)}$ for $j\le4$; on the diagonal $M_2^+\le C$. Since $\mathcal B_\alpha$ is measurable with respect to the remaining entries, conditional integration gives
\begin{align*}
 \mathbb E\mathbf1_{\mathcal B_\alpha^c}
   \left|\int(A(K+V(x))-A(K))\,\Delta\mu_\alpha(\mathrm dx)\right|
 &\le4\mathbb P(\mathcal B_\alpha^c),\\
 \mathbb E\mathbf1_{\mathcal B_\alpha^c}
   |C_{j,\alpha}(K):\Delta M_{j,\alpha}|
 &\le N^6M_{j,\alpha}^+\mathbb P(\mathcal B_\alpha^c),\\
 \mathbb E\mathbf1_{\mathcal B_\alpha^c}
       \int|R_{4,\alpha}(K;x)|\,\mu_\alpha^+(\mathrm dx)
 &\le C N^6M_{4,\alpha}^+\mathbb P(\mathcal B_\alpha^c).
\end{align*}
The first bound can be sharpened to $2\mathbb P(\mathcal B_\alpha^c)$ by integrating the original function values. On the diagonal, use the same argument for the linear polynomial
and its quadratic remainder, with $M_2^+$ in place of $M_4^+$. With failure probability at most $N^{-D}$, the sums over all off-diagonal blocks are respectively $O(N^{2-D})$,
$N^{8-D+o(1)}$, and $N^{8-D+o(1)}$; the diagonal bounds are at most $N^{7-D+o(1)}$. For any $A>0$, choosing $D>A+10$ makes every displayed sum $O(N^{-A})$, uniformly in the
interpolation parameter whose integration interval has length one. The additional derivatives of $U_N^{\rm Sch}$ are bounded by Corollary~\ref{W:cor:schur-cutoff-suppression} with
its explicit exponent choice. If $\mathcal B_\alpha$ holds but no permitted entry value gives the bound on the Gram matrix, the function and its derivatives vanish throughout a
neighbourhood of the support of the replacement entry distributions. Thus \eqref{W:eq:comparison-sumR} holds after all matrices with the selected entry set to zero and all
permitted replacement values are included.

\subsubsection{Comparison with bounded entries and removal of the cutoff}
\label{W:sec:removal-of-coefficient-cutoff-in-comparison}
The entrywise estimates now give a comparison of expectations. In \eqref{W:eq:entry-distribution-interpolation}, with $F(H)=A_N(s(H,w))$, the constant coefficient integrates to
zero against the difference of the truncated and bounded distributions. The cubic polynomial off the diagonal and the linear polynomial on the diagonal account for the remaining
coefficients of these expansions on the specified events. Their errors are bounded by \eqref{W:eq:summed-first-three-moment-errors}; the remainders are bounded by
\eqref{W:eq:comparison-sumR} and \eqref{W:eq:comparison-diagonal-remainder}. Adding the exceptional-event estimates yields the following result, which also provides the comparison
needed to remove the coefficient cutoff.

\begin{proposition}
\label{W:prop:comparison-local-comparison}
Let $\mathfrak r_N\to\infty$ be deterministic with $\mathfrak r_N=N^{o(1)}$, and let $A_N$ satisfy \eqref{W:eq:replacement-test-class}. Then
\begin{equation}
 \sup_{\|w\|=1}
 |\mathbb EA_N(H_N^{\mathrm{tr}},w)-\mathbb EA_N(H_N^{\rm b},w)|
 \le N^{-1/3+o(1)}.
 \label{W:eq:comparison-local-comparison}
\end{equation}
\end{proposition}
\begin{proof}
Write $m=m_N$ and $S_N=1+me_N^4\le2m$ for large $N$. Here $\mu_\alpha$ is the distribution of the uncentered truncated entry and $\nu_\alpha$ that of its bounded comparison entry.
Define $\Delta\mu_\alpha=\mu_\alpha-\nu_\alpha$ and $\mu_\alpha^+=\mu_\alpha+\nu_\alpha$. The remaining entries of $K_\alpha(t)$ have distribution
$\bigotimes_{\delta\ne\alpha}[t\mu_\delta+(1-t)\nu_\delta]$. Conditional on $\sigma(K_\alpha(t))$, the polynomial coefficients are fixed and are integrated against the signed
measure $\Delta\mu_\alpha$.

For each $\alpha$ split the expectation into $\mathcal B_\alpha\cap\mathcal A_{\alpha,\rho}$, $\mathcal B_\alpha\setminus\mathcal A_{\alpha,\rho}$, and $\mathcal B_\alpha^c$, as
above. The first part uses the explicit estimates with indicator $I_\alpha$; the second vanishes because the cutoff and its derivatives vanish; the last is bounded by the
deterministic derivative estimates. In the first part, and as an algebraic identity with the error defined by subtraction elsewhere, the zero mass of $\Delta\mu_\alpha$ gives
\[
 \begin{aligned}
 \int F(K_\alpha+V_\alpha(x))\,\Delta\mu_\alpha(\mathrm dx)
  &=\sum_{j=1}^3C_{j,\alpha}:\Delta M_{j,\alpha}
       +\int R_{4,\alpha}(K_\alpha;x)\,\Delta\mu_\alpha(\mathrm dx),
       &&\alpha\text{ off diagonal},\\
 \int F(K_i+N^{-1/2}xE_{ii})\,\Delta\mu_i(\mathrm dx)
  &=C_{1,i}\Delta M_{1,i}
       +\int\mathcal R_i(K_i;x)\,\Delta\mu_i(\mathrm dx),
       &&\alpha=\{i,i\}.
 \end{aligned}
\]
The absolute value of each remainder integral is bounded by its integral against $\mu_\alpha^+$. Thus the signed coefficient contractions and the positive remainder integrals are
distinct terms of \eqref{W:eq:entry-distribution-interpolation}. Taking expectations and integrating in $t$ gives
\begin{align*}
 |\mathbb E_\mu A_N-\mathbb E_\nu A_N|
 \le \Lambda_NN^{o(1)}\bigl[&m^{-1}+(m^{-1}+m^{-2})
               +(m^{-5/2}+2m^{-7/2})\\
 &+m^{-1}+m^{-2}+m^{-1}
           +m^{-3/2}S_N^{1/2}+m^{-2}S_N
           +m^{-1}\bigr]+O(N^{-A}).
\end{align*}
The first three terms are the signed coefficient bounds \eqref{W:eq:summed-first-three-moment-errors}. The next five are the positive fourth moment remainder bounds
\eqref{W:eq:summed-fourth-remainder-terms}, and the last term is the diagonal comparison \eqref{W:eq:comparison-diagonal-remainder}. The exceptional event calculation in
Section~\ref{W:sec:summation-of-fourth-order-remainder} permits arbitrary fixed $A$. Since $m^{-3/2}S_N^{1/2}\le\sqrt2m^{-1}$ and $m^{-2}S_N\le2m^{-1}$, the right side is
$\Lambda_Nm^{-1}N^{o(1)}+O(N^{-A})=N^{-1/3+o(1)}$. For the uniform estimate, let $\mathcal T_N(t,w)$ be the sum of the expectations of the absolute terms just bounded, before
integration in $t$. The estimates give $\mathcal T_N(t,w)\le Q_Nm^{-1}+O(N^{-A})$ for every $t$ and deterministic unit vector $w$. Integrating this pointwise deterministic bound
first and then taking the supremum gives
\[
 \sup_{\|w\|=1}\int_0^1\mathcal T_N(t,w)\,\mathrm dt
       \le Q_Nm^{-1}+O(N^{-A}).
\]
This proves \eqref{W:eq:comparison-local-comparison}. The choice of parameters has the following order. Fix the parameters of $A_N$ and any deterministic $\mathfrak r_N\to\infty$
with $\mathfrak r_N=N^{o(1)}$; these determine the polynomial derivative bounds. Next choose $D>A+10$ for the exceptional events, and finally let $N$ tend to infinity.
Lemma~\ref{W:lem:comparison-ensemble-coefficients} subsequently selects a sequence $\mathfrak r_N$ for removal of the cutoff. The comparison just proved applies to that sequence.
\end{proof}

Use $F_N$ from \eqref{W:eq:comparison-test}. Choose $\mathfrak r_N$ by Lemma~\ref{W:lem:comparison-ensemble-coefficients}, after establishing the preceding comparisons for an
arbitrary subpolynomial cutoff threshold.

\Needspace{6\baselineskip}

\begin{theorem}
\label{W:thm:fourier-comparison-with-bounded-entries}
With the coefficient cutoff chosen by Lemma~\ref{W:lem:comparison-ensemble-coefficients},
\begin{align}
 \sup_{\|w\|=1}\left|
  \mathbb EF_N(H_N^{\mathrm{tr}},w)-\mathbb EF_N(H_N^{\rm b},w)
                    \right|&\le N^{-1/3+o(1)},\label{W:eq:comparison-main}\\
 \sup_{\|w\|=1}\mathbb E D_N^{\rm coef}(H_N^{\mathrm{tr}},w)
                    &\le N^{-1/3+o(1)}.\label{W:eq:expected-coefficient-cutoff-error}
\end{align}
The localized function $F_N^{\rm loc}=\chi_NF_N$ has summed absolute fourth-order remainder $N^{-1/3+o(1)}$.
\end{theorem}

\begin{proof}[Proof of Theorem~\ref{W:thm:fourier-comparison-with-bounded-entries}]
Choose one $\mathfrak r_N$ satisfying Lemma~\ref{W:lem:comparison-ensemble-coefficients} for both $H_N^{\rm b}$ and $H_N^{\mathrm G}$. The nonnegative function
$D_N^{\rm coef}=\Gamma_\rho(1-\chi_N)$ satisfies \eqref{W:eq:replacement-test-class}. Proposition~\ref{W:prop:comparison-local-comparison} and
\eqref{W:eq:comparison-ensemble-cutoff-error} therefore give
\[
 \begin{aligned}
 \sup_{\|w\|=1}\mathbb ED_N^{\rm coef}(H_N^{\mathrm{tr}},w)
 &\le\sup_{\|w\|=1}\left|\mathbb ED_N^{\rm coef}(H_N^{\mathrm{tr}},w)
                         -\mathbb ED_N^{\rm coef}(H_N^{\rm b},w)\right|
      +\sup_{\|w\|=1}\mathbb ED_N^{\rm coef}(H_N^{\rm b},w)\\
 &\le N^{-1/3+o(1)}+N^{-D}
  \le N^{-1/3+o(1)},
 \end{aligned}
\]
where we fix, for example, $D=1$. This is \eqref{W:eq:expected-coefficient-cutoff-error}. For each $H$, the phase in $F_N(s(H,w))$ is real. Thus
\begin{equation}
 |F_N-F_N^{\rm loc}|=\Gamma_\rho(1-\chi_N)=D_N^{\rm coef}.
 \label{W:eq:comparison-ensemble-identity}
\end{equation}
Therefore
\begin{align*}
 |\mathbb EF_N(H_N^{\mathrm{tr}},w)-\mathbb EF_N(H_N^{\rm b},w)|
 &\le |\mathbb EF_N^{\rm loc}(H_N^{\mathrm{tr}},w)
                -\mathbb EF_N^{\rm loc}(H_N^{\rm b},w)|\\
 &\quad+\mathbb ED_N^{\rm coef}(H_N^{\mathrm{tr}},w)
                +\mathbb ED_N^{\rm coef}(H_N^{\rm b},w)\\
 &\le N^{-1/3+o(1)}.
\end{align*}
Proposition~\ref{W:prop:comparison-local-comparison} applies to $A_N=F_N^{\rm loc}$. Substituting that bound and \eqref{W:eq:expected-coefficient-cutoff-error} in the displayed
three term inequality proves \eqref{W:eq:comparison-main}. For the remainder bound, \eqref{W:eq:comparison-sumR} and the exceptional event bound are applied to the same function
$F_N^{\rm loc}$. In particular its Taylor derivatives still include the derivatives of $\chi_N$. Finally, on the original truncation coupling, for any function $G_N$ with
$|G_N|\le1$,
\[
 |\mathbb EG_N(H_N,w)-\mathbb EG_N(H_N^{\mathrm{tr}},w)|
 \le2\mathbb P(H_N\ne H_N^{\mathrm{tr}})=o(1).
\]
The remainder bound applies to the function, including its cutoffs, used above.

\end{proof}

\subsubsection{Comparison with the GOE and GUE}
\label{W:sec:comparison-and-convergence-of-the-integrated-potential}
The signed cubic estimate bounds the remaining off-diagonal contribution to the comparison with the GOE or GUE. The diagonal contributions are handled by the preceding second
moment comparison, which permits their variances to differ. Gaussian tails enter the integral remainders, while the polynomial coefficients are integrated against the full signed
distribution. Applying the interpolation identity to observables multiplied by $\chi_N$ now gives the comparison below; the expectation estimate for $D_N^{\rm coef}$ controls the
removal of this factor.

\begin{corollary}
\label{W:cor:cumulant-comparison}
Let $A_N$ satisfy \eqref{W:eq:replacement-test-class}, and let $H_N^{\rm b}$ and $H_N^{\mathrm G}$ be, respectively, the matrix with bounded entries and the Gaussian matrix defined
above. Then
\begin{equation}
 \sup_{\|w\|=1}\left|
  \mathbb EA_N(s(H_N^{\rm b},w))-
  \mathbb EA_N(s(H_N^{\mathrm G},w))\right|
 \le N^{-1/6+o(1)}.
 \label{W:eq:cumulant-ensemble-comparison}
\end{equation}
\end{corollary}
\begin{proof}
Fix a deterministic unit vector $w$ and write $F(H)=A_N(s(H,w))$. For each upper triangular entry $\alpha$, set
\[
 \Delta\mu_\alpha=\mu_\alpha^{\rm b}-\mu_\alpha^G,
 \qquad \lambda_{\alpha,t}=t\mu_\alpha^{\rm b}+(1-t)\mu_\alpha^G.
\]
The other entries of $K_\alpha(t)$ have distribution $\bigotimes_{\delta\ne\alpha}\lambda_{\delta,t}$. Lemma~\ref{W:lem:entry-distribution-interpolation} gives
\[
 \mathbb EF(H_N^{\rm b})-\mathbb EF(H_N^{\mathrm G})
 =\int_0^1\sum_\alpha\mathbb E_{K_\alpha(t)}
       \int F(K_\alpha(t)+V_\alpha(x))\,\Delta\mu_\alpha(\mathrm dx)\,\mathrm dt.
\]
This formula is absolutely integrable since $|F|\le1$ and each signed measure has total variation at most two.

For an off-diagonal entry, subtract its cubic polynomial. With the coefficients and normalization of \eqref{W:eq:cumulant-chain-tensor}, the resulting exact identity is
\[
 \int F(K_\alpha+V_\alpha(x))\,\Delta\mu_\alpha(\mathrm dx)
 =C_{3,\alpha}(K_\alpha):\Delta M_3
       +\int R_{4,\alpha}(K_\alpha;x)\,\Delta\mu_\alpha(\mathrm dx).
\]
Indeed, the two distributions have total mass one, zero first coordinate moments, and second coordinate tensor $I_\beta$. Thus the coefficients of degrees zero, one, and two have
zero signed integrals. In the complex case these equalities concern the joint vector of real coordinates, so every product of three coordinates remains in $\Delta M_3$. For a
diagonal entry the constant and linear coefficients integrate to zero. Its error is estimated using Theorem~\ref{W:thm:truncation-diagonal}; the standardized variances are $1$ for
the bounded signs and $2/\beta$ for the Gaussian entry.

It remains to justify the estimates when the Gaussian replacement value, or an entry of $K_\alpha(t)$, lies outside the bounded supports used in the preceding lemmas. Put
$D_N=\{x:\|x\|_1\le L_N\}$. A real Gaussian variable, or a two-dimensional Gaussian vector with the fixed covariance here, satisfies
\[
 \int_{D_N^c}(1+\|x\|_1^p)\,\mu_\alpha^G(\mathrm dx)
       \le C_pe^{-c_pL_N^2}\qquad(p\text{ fixed}).
\]
For example, choose $c>0$ so that $\mathbb E\exp(2c\|x\|_1^2)<\infty$. The inequality $(1+r^p)\mathbf1_{\{r>L\}}\le C_{p,c}e^{-cL^2}e^{2cr^2}$ proves this bound. The bounded
comparison distribution is supported in $D_N$ for large $N$.

At every Hermitian $K$, Lemma~\ref{W:lem:global-cubic-entry-derivatives} gives
\[
 |R_{4,\alpha}(K;x)|\le2+CN^6\sum_{j=1}^3\|x\|_1^j.
\]
Consequently its integral over $D_N^c$, summed over at most $N^2$ entries, is at most $CN^8e^{-cL_N^2}$. The diagonal error has the bound $2+CN^6|x|$ and satisfies the same upper
bound after summation. These estimates concern the errors only. The polynomial coefficients continue to be integrated against the full measures, so their exact moment cancellations
are unchanged.

For each remaining entry $\delta$, construct mutually independent variables $U_\delta,\xi_\delta^{\rm b},\xi_\delta^G$, with $U_\delta$ uniform on $[0,1]$ and the last two having
the indicated distributions. The triples for distinct entries are independent. At parameter $t$ use $\xi_\delta^{\rm b}$ when $U_\delta\le t$ and $\xi_\delta^G$ otherwise. In a
second matrix use instead $\widehat\xi_\delta^G=\xi_\delta^G\mathbf1_{D_N}(\xi_\delta^G)$ with the same triples. A complex entry is changed as a whole vector. This defines a
coupling of $K_\alpha(t)$ and $\widehat K_\alpha(t)$ with
\[
 \mathbb P(K_\alpha(t)\ne\widehat K_\alpha(t))
 \le\sum_{\delta\ne\alpha}\mathbb P(\xi_\delta^G\notin D_N)
 \le CN^2e^{-cL_N^2},
\]
uniformly in $t$. The replacement variable $x$ is independent of all these triples.

For any coefficient family with $\sup_K\|B_\alpha(K)\|\le N^B$, equality of the two matrices outside this event implies
\begin{equation}
 \sum_\alpha\mathbb E\|B_\alpha(K_\alpha)
                         -B_\alpha(\widehat K_\alpha)\|
 \le2N^B\sum_\alpha\mathbb P(K_\alpha\ne\widehat K_\alpha)
 \le CN^{B+4}e^{-cL_N^2}.
 \label{W:eq:gaussian-inactive-restoration}
\end{equation}
For the error itself, subtraction of the defining polynomials gives
\[
 |R_{4,\alpha}(K;x)-R_{4,\alpha}(\widehat K;x)|
 \le C\mathbf1_{\{K\ne\widehat K\}}
           \left(1+N^6\sum_{j=1}^3\|x\|_1^j\right).
\]
Integration in the independent $x$ against $\mu_\alpha^{\rm b}+\mu_\alpha^G$, followed by the sum over $\alpha$, is at most $CN^{10}e^{-cL_N^2}$. The diagonal calculation uses
degrees zero and one and is bounded by the same expression.

The auxiliary Gaussian distribution is
\[
 \widehat\mu_\alpha^G
 =\mathbf1_{D_N}\mu_\alpha^G+\mu_\alpha^G(D_N^c)\delta_0.
\]
This distribution is symmetric under $x\mapsto-x$, so its first and third coordinate tensors vanish. Write $\Sigma_{0,\alpha}$ for the original covariance and
\[
 \Delta_\alpha=\int_{D_N^c}xx^{\mathsf T}\mu_\alpha^G(\mathrm dx),
 \qquad
 \varepsilon_\alpha=
 \|\Sigma_{0,\alpha}^{-1/2}\Delta_\alpha
                         \Sigma_{0,\alpha}^{-1/2}\|.
\]
Then $0\le\Delta_\alpha\le\varepsilon_\alpha\Sigma_{0,\alpha}$ and $\varepsilon_\alpha\le Ce^{-cL_N^2}$. We restore the covariance exactly before applying the preceding estimates.
When $\varepsilon_\alpha>0$, put $\delta_\alpha=2\varepsilon_\alpha<1/2$ for large $N$ and define
\[
 C_\alpha=\Sigma_{0,\alpha}
                +\frac{1-\delta_\alpha}{\delta_\alpha}\Delta_\alpha.
\]
The matrix inequalities give
$\Sigma_{0,\alpha}\le C_\alpha\le(3/2)\Sigma_{0,\alpha}$.
For an independent vector $\epsilon_\alpha$ of independent symmetric signs, let $Z_\alpha=C_\alpha^{1/2}\epsilon_\alpha$. Its mean and third tensor are zero, its covariance is
$C_\alpha$, and its support is bounded by a fixed constant since its dimension is at most two. Set
\[
 \widetilde\mu_\alpha^G
  =(1-\delta_\alpha)\widehat\mu_\alpha^G
                       +\delta_\alpha\mathcal L(Z_\alpha).
\]
If $\varepsilon_\alpha=0$, set $\widetilde\mu_\alpha^G=\widehat\mu_\alpha^G$ and $\delta_\alpha=0$. In the nonzero case its covariance is exactly
\[
 (1-\delta_\alpha)(\Sigma_{0,\alpha}-\Delta_\alpha)
                +\delta_\alpha C_\alpha=\Sigma_{0,\alpha}.
\]
Thus $\widetilde\mu_\alpha^G$ has the same first three coordinate tensors as $\mu_\alpha^G$; it is supported in $D_N$ for large $N$ and has uniformly bounded absolute moments of
every fixed order.

For a coupling, use the variable with distribution $\widehat\mu_\alpha^G$ with probability $1-\delta_\alpha$ and the independent $Z_\alpha$ with probability $\delta_\alpha$, with
these choices independent between entries. The resulting matrix $\widetilde K_\alpha(t)$ differs from $\widehat K_\alpha(t)$ with probability at most
$\sum_{\delta\ne\alpha}\delta_\delta\le CN^2e^{-cL_N^2}$. The estimate \eqref{W:eq:gaussian-inactive-restoration} and the following bound for $R_{4,\alpha}$ therefore also apply to
this pair of matrices. For the entry being replaced, the signed difference $\widetilde\mu_\alpha^G-\mu_\alpha^G$ has, for each fixed $p$,
\[
 \int(1+\|x\|_1^p)
       |\widetilde\mu_\alpha^G-\mu_\alpha^G|(\mathrm dx)
       \le C_pe^{-c_pL_N^2}.
\]
Indeed, the difference between $\mu_\alpha^G$ and $\widehat\mu_\alpha^G$ consists of the removed tail and its mass at zero, and the remaining difference is
$\delta_\alpha(\mathcal L(Z_\alpha)-\widehat\mu_\alpha^G)$. The previously proved polynomial bound for the error thus adds at most $CN^8e^{-cL_N^2}$ after summing the integrals
over the replaced entries. The entries with restored covariance satisfy the exact covariance normalization and all moment and support bounds used by the preceding resolvent and
conditional row estimates. In particular their fourth absolute moments are uniformly bounded. For the matrices $\widetilde K_\alpha(t)$ and replacement values in $D_N$,
Lemma~\ref{W:lem:comparison-remainder} and the summed resolvent estimates bound the off-diagonal errors by $N^{-1/3+o(1)}$. Theorem~\ref{W:thm:truncation-diagonal} gives the same
bound for the diagonal errors, since the means are zero and the two second moments have bounded sum. For a prescribed $A>0$, the integrated polynomial bound on each exceptional
event is at most $CN^6$, and the sum over at most $N^2$ entries is at most $CN^{8-D}$. Choosing the failure exponent $D>A+10$ makes this contribution $O(N^{-A})$.

Since $L_N^2=N^{1-o(1)}$, for every fixed $B,A>0$,
\[
 N^Be^{-cL_N^2}
 =N^{-A}\exp\{(A+B)\log N-cL_N^2\}=O(N^{-A}).
\]
We may thus return from $\widetilde K_\alpha$ through $\widehat K_\alpha$ to $K_\alpha$, and from the distribution of the replaced entry with restored covariance to the original
one, in the error integrals. The cubic coefficients were integrated against the full signed measures throughout. Combining the entry identities gives
\begin{equation}
 \mathbb EA_N(s(H_N^{\rm b},w))-
 \mathbb EA_N(s(H_N^{\mathrm G},w))
 =\int_0^1\mathbb E Z_{3,N}(A;t,w)\,\mathrm dt
                         +O(N^{-1/3+o(1)}).
 \label{W:eq:gaussian-comparison-cubic-decomposition}
\end{equation}
Every estimate is uniform in deterministic $w$. Theorem~\ref{W:thm:cumulant-complete} bounds the integral by $N^{-1/6+o(1)}$, proving \eqref{W:eq:cumulant-ensemble-comparison}.
\end{proof}

The following identity separates the entrywise remainder of the expansion from the sum of all cubic terms, with their signs retained, estimated in
Section~\ref{W:sec:cumulant-cubic}. The tensors include every derivative of $\chi_N$. Let $K_\alpha^{\rm b,G}(t)$ have block $\alpha$ set to zero and independent remaining blocks,
with law $t\mu_\delta^{\rm b}+(1-t)\mu_\delta^{G}$ at block $\delta\ne\alpha$. Set
\begin{align}
 C_{\alpha,3}^{\rm loc}(K)
 &=\frac1{3!}D_x^3
       \left[F_N^{\rm loc}(K+N^{-1/2}\mathcal T_\alpha x,w)\right]_{x=0},
       \nonumber\\
 \mathfrak C_{3,N}^{\rm loc}(w)
 &=\int_0^1\sum_{\alpha\ {\rm off}}\mathbb E
     \left\langle C_{\alpha,3}^{\rm loc}(K_\alpha^{\rm b,G}(t)),
                         M_3^{\rm b}(\alpha)-M_3^G(\alpha)\right\rangle
                         \,\mathrm dt.
 \label{W:eq:comparison-gaussian-cubic}
\end{align}
The tensors of orders one, two, and three include the factors $1$, $1/2$, and $1/6$, respectively, as well as the entry normalization. Here $M_3^G=0$, and $M_3^{\rm b}$ is
uniformly bounded.

\Needspace{6\baselineskip}

\begin{corollary}
\label{W:cor:comparison-gaussian-reduction}
For $F_N$ in \eqref{W:eq:comparison-test} and the cutoff sequence selected in Lemma~\ref{W:lem:comparison-ensemble-coefficients}, uniformly over deterministic unit $w$,
\begin{equation}
 \mathbb EF_N(H_N^{\mathrm{tr}},w)-\mathbb EF_N({H_N^{\mathrm G}},w)
 =\mathfrak C_{3,N}^{\rm loc}(w)+O(N^{-1/3+o(1)}).
 \label{W:eq:comparison-gaussian-reduction}
\end{equation}
Moreover $\sup_{\|w\|=1}|\mathfrak C_{3,N}^{\rm loc}(w)|\le N^{-1/6+o(1)}$. Replacing $H_N^{\mathrm{tr}}$ by $H_N$ adds an error of modulus at most
$2\mathbb P(H_N\ne H_N^{\mathrm{tr}})=o(1)$.
\end{corollary}
\begin{proof}
Let $F=F_N$ and $A=F_N^{\rm loc}=\chi_NF_N$, with $w$ fixed. The exact difference of expectations is
\[
 \begin{aligned}
 \mathbb EF(H_N^{\mathrm{tr}})-\mathbb EF(H_N^{\mathrm G})
 ={}&[\mathbb EF(H_N^{\mathrm{tr}})-\mathbb EF(H_N^{\rm b})]
       +\mathbb E(F-A)(H_N^{\rm b})\\
    &+[\mathbb EA(H_N^{\rm b})-\mathbb EA(H_N^{\mathrm G})]
       +\mathbb E(A-F)(H_N^{\mathrm G}).
 \end{aligned}
\]
Theorem~\ref{W:thm:fourier-comparison-with-bounded-entries} bounds the first bracket by $N^{-1/3+o(1)}$, uniformly in $w$. For the second and fourth terms,
\eqref{W:eq:comparison-ensemble-identity} gives
\[
 \begin{aligned}
 |\mathbb E(F-A)(H_N^{\rm b},w)|
       &\le\mathbb ED_N^{\rm coef}(H_N^{\rm b},w),\\
 |\mathbb E(A-F)(H_N^{\mathrm G},w)|
       &\le\mathbb ED_N^{\rm coef}(H_N^{\mathrm G},w).
 \end{aligned}
\]
Their suprema are at most $N^{-D}$ for every fixed $D$, with the common cutoff sequence of Lemma~\ref{W:lem:comparison-ensemble-coefficients}. For the third bracket,
\eqref{W:eq:gaussian-comparison-cubic-decomposition} applied to $A$ gives
\[
 \mathbb EA(H_N^{\rm b},w)-\mathbb EA(H_N^{\mathrm G},w)
 =\int_0^1\mathbb E Z_{3,N}(A;t,w)\,\mathrm dt
                                  +O(N^{-1/3+o(1)}).
\]
The integral is exactly $\mathfrak C_{3,N}^{\rm loc}(w)$: $A=F_N^{\rm loc}$, and both definitions use $1/3!$ times the third standardized entry derivative contracted with
$M_3^{\rm b}-M_3^G$, summed over all off-diagonal blocks. Substitution in the four term identity, with $D=1$, proves \eqref{W:eq:comparison-gaussian-reduction}.
Theorem~\ref{W:thm:cumulant-complete} bounds this cubic integral by $N^{-1/6+o(1)}$. On the original coupling,
\[
 \sup_w|\mathbb EF_N(H_N,w)-\mathbb EF_N(H_N^{\mathrm{tr}},w)|
        \le2\mathbb P(H_N\ne H_N^{\mathrm{tr}})=o(1),
\]
which is the additional error in the comparison bound for the untruncated array, as in \eqref{W:eq:comparison-ensemble-transfer}.
\end{proof}

On a segment $K+uV$ obtained by changing one entry, write $f(u)=F_N(K+uV,w)$ and $c(u)=\chi_N(K+uV,w)$. The cubic coefficient in the corollary is
\begin{equation}
 \frac{(cf)'''(0)}6
 =c(0)\frac{f'''(0)}6+\frac{c'(0)f''(0)}2
       +\frac{c''(0)f'(0)}2+\frac{c'''(0)f(0)}6.
 \label{W:eq:comparison-cutoff-cubic-terms}
\end{equation}
Applying Theorem~\ref{W:thm:cumulant-complete} to $A_N^{\rm F}$ gives
\begin{equation}
 \sup_{\|w\|=1}|\mathfrak C_{3,N}^{\rm loc}(w)|\longrightarrow0.
 \label{W:eq:vanishing-fourier-cubic-term}
\end{equation}
The quantitative bound is recorded in \eqref{W:eq:fourier-third-moment-rate}.

For later use write
\begin{equation}
 \mathfrak C_{3,N}(A;w)=\int_0^1\mathbb E Z_{3,N}(A;t,w)\,\mathrm dt,
 \qquad
 \sup_{\|w\|=1}|\mathfrak C_{3,N}(A;w)|\le N^{-1/6+o(1)}.
 \label{W:eq:tightness-complete-cubic}
\end{equation}
Here $Z_{3,N}$ is the ordered contraction in \eqref{W:eq:cumulant-target-tensor}; the bound is \eqref{W:eq:cumulant-integrated}.

\Needspace{6\baselineskip}

\subsubsection{Characteristic functions and tightness}
\label{W:sec:convergence-of-the-integrated-potential}
\label{W:sec:proof-local-coefficients}

The two expectation comparisons are sufficient to identify the coefficient process once they are applied to two kinds of observables. Characteristic functions determine its
finite-dimensional distributions, and the smooth functions used to bound the modulus of continuity provide tightness. Both use the same cutoff estimates:
Proposition~\ref{W:prop:comparison-local-comparison} and Corollary~\ref{W:cor:cumulant-comparison} apply to $A_N^{\rm F}$, $A_N^+$, $A_N^-$, and $A_N^0$, while $D_N^{\rm coef}$
bounds the error from removing $\chi_N$ under the compared distributions.

The functions $\widehat F_{M,N}^\sigma$ are formed from the modified moments of Section~\ref{W:sec:the-moment-vector-and-its-smooth-extensions}. For the original matrix, truncation
and removal of the spectral modification add the $o(1)$ error in \eqref{W:eq:comparison-ensemble-transfer}. All suprema in the theorem below are over deterministic unit starting
vectors.

\begin{theorem}
\label{W:thm:tightness-main}
For every fixed $T,\delta,u$ and parameter in the bound on the Gram matrix,
\begin{equation}
 \max_{\sigma=\pm1}\sup_{\|w\|=1}
 \left|\mathbb E\widehat F_{M,N}^\sigma(H_N^{\mathrm{tr}},w)
       -\mathbb E\widehat F_{M,N}^\sigma(H_N^{\mathrm G},w)\right|
 \le N^{-1/6+o(1)}.
 \label{W:eq:tightness-main}
\end{equation}
The same estimate holds for $\Gamma_\rho$.
\end{theorem}

\begin{proof}[Proof of Theorem~\ref{W:thm:tightness-main}]
Combine Proposition~\ref{W:prop:comparison-local-comparison} with Corollary~\ref{W:cor:cumulant-comparison}, and use $D_N^{\rm coef}$ from
Section~\ref{W:sec:removing-the-coefficient-cutoff}. Equations \eqref{W:eq:expected-coefficient-cutoff-error} and \eqref{W:eq:comparison-ensemble-cutoff-error} give
\begin{equation}
 \sup_w\mathbb ED_N^{\rm coef}(H_N^{\mathrm{tr}},w)
       \le N^{-1/3+o(1)},\qquad
 \sup_w\mathbb ED_N^{\rm coef}({H_N^{\mathrm G}},w)\le N^{-D_0}
 \quad(D_0<\infty).
 \label{W:eq:expected-cutoff-errors}
\end{equation}
Since $0\le M\le1$,
\[
 |\widehat F_{M,N}^\sigma-A_N^\sigma|\le D_N^{\rm coef},
 \qquad |\Gamma_\rho-A_N^0|=D_N^{\rm coef}.
\]
For either sign put $F=\widehat F_{M,N}^\sigma$ and $A=A_N^\sigma$; for the cutoff comparison put $F=\Gamma_\rho$ and $A=A_N^0$. In both cases
\[
 \begin{aligned}
 \sup_w|\mathbb EF(H_N^{\rm tr},w)-\mathbb EF(H_N^{\rm G},w)|
 \le{}&\sup_w|\mathbb EA(H_N^{\rm tr},w)-\mathbb EA(H_N^{\rm b},w)|\\
 &+\sup_w|\mathbb EA(H_N^{\rm b},w)-\mathbb EA(H_N^{\rm G},w)|\\
 &+\sup_w\mathbb ED_N^{\rm coef}(H_N^{\rm tr},w)
       +\sup_w\mathbb ED_N^{\rm coef}(H_N^{\rm G},w)\\
 \le{}&N^{-1/3+o(1)}+N^{-1/6+o(1)}
                         +N^{-1/3+o(1)}+N^{-D_0}.
 \end{aligned}
\]
Take $D_0\ge1$. A maximum over the two signs preserves the $N^{-1/6+o(1)}$ bound, proving \eqref{W:eq:tightness-main}. The choice $F=\Gamma_\rho$ gives
\begin{equation}
 \sup_{\|w\|=1}\left|
 \mathbb E\Gamma_\rho(H_N^{\mathrm{tr}},w)
       -\mathbb E\Gamma_\rho({H_N^{\mathrm G}},w)\right|
 \le N^{-1/6+o(1)}.
 \label{W:eq:tightness-Gram-comparison}
\end{equation}
We next pass from modified moments to the original matrices. On the probability space used for entry truncation, define
\[
 \mathcal D_N^H=\{H_N=H_N^{\rm tr},\ h_d(H_N)=H_N\},
 \qquad \mathcal D_N^G=\{h_d(H_N^{\rm G})=H_N^{\rm G}\}.
\]
The truncation bounds and Proposition~\ref{W:prop:gram-window-2} give
\[
 \begin{aligned}
 p_N^H:=\mathbb P((\mathcal D_N^H)^c)
 &\le\mathbb P(H_N\ne H_N^{\rm tr})
       +\mathbb P(\|H_N\|>2+d/2)=o(1),\\
 p_N^G:=\mathbb P((\mathcal D_N^G)^c)
 &\le\mathbb P(\|H_N^{\rm G}\|>2+d/2)=o(1).
 \end{aligned}
\]
More explicitly, with $L_N=e_N\sqrt N$, the union bounds in Lemma~\ref{W:lem:admissible-scale} give
\[
 \mathbb P(H_N\ne H_N^{\rm tr})
 \le C e_N^{-4}\Theta(L_N)
                   +e_N^{-2}\Theta_{\rm d}(L_N)\longrightarrow0.
\]
The same lemma chooses $e_N$ so that both terms tend to zero. The remaining probabilities are bounded by Proposition~\ref{W:prop:gram-window-2}. These events concern the matrices
themselves, and their probabilities are independent of the choice of deterministic $w$. On $\mathcal D_N^H$ the two polynomial matrices, their Gram matrices, and their cutoff
values agree. On the subset where this common cutoff on the Gram eigenvalues is positive, the orthogonalization formulas therefore give the same integrated potential. On its zero
set both bounded functions vanish. The same reasoning applies on $\mathcal D_N^G$. Thus, for the original and modified versions $F^{\rm orig}$ and $\widehat F$, respectively,
\[
 \begin{aligned}
 |F^{\rm orig}(H_N,w)-\widehat F(H_N^{\rm tr},w)|
                      &\le2\mathbf1_{(\mathcal D_N^H)^c},\\
 |F^{\rm orig}(H_N^{\rm G},w)-\widehat F(H_N^{\rm G},w)|
                      &\le2\mathbf1_{(\mathcal D_N^G)^c}.
 \end{aligned}
\]
They agree on the respective events, and both moduli are at most one elsewhere. Thus the original expectation difference is at most
\[
 N^{-1/6+o(1)}+2\mathbb P((\mathcal D_N^H)^c)
                       +2\mathbb P((\mathcal D_N^G)^c)=o(1).
\]
Taking the positive part when required gives \eqref{W:eq:original-modulus-comparison} and \eqref{W:eq:original-gram-comparison}.
\end{proof}

\Needspace{6\baselineskip}

\begin{corollary}
\label{W:cor:tightness-Fourier}
The function $F_N$ in \eqref{W:eq:comparison-test} satisfies
\begin{align}
 \sup_{\|w\|=1}|\mathfrak C_{3,N}^{\rm loc}(w)|
 &\le N^{-1/6+o(1)},\label{W:eq:fourier-third-moment-rate}\\
 \sup_{\|w\|=1}\left|\mathbb EF_N(H_N^{\mathrm{tr}},w)
             -\mathbb EF_N(H_N^{\mathrm G},w)\right|
 &\le N^{-1/6+o(1)}.
 \label{W:eq:tightness-Fourier-comparison}
\end{align}
For the original untruncated matrix, the corresponding difference in~\eqref{W:eq:original-fourier-comparison} tends to zero.
\end{corollary}
\begin{proof}
Corollary~\ref{W:cor:comparison-gaussian-reduction} gives both $\sup_w|\mathfrak C_{3,N}^{\rm loc}(w)|\le N^{-1/6+o(1)}$ and
\[
 \sup_w|\mathbb EF_N(H_N^{\rm tr},w)-\mathbb EF_N(H_N^{\rm G},w)|
 \le\sup_w|\mathfrak C_{3,N}^{\rm loc}(w)|+N^{-1/3+o(1)}
 \le N^{-1/6+o(1)}.
\]
This proves the two displayed estimates. For the characteristic functions with their cutoffs, let $p_N^H$ be the probability that $H_N^{\rm tr}\ne H_N$ or $h_d(H_N)\ne H_N$. Let
$p_N^G=\mathbb P(h_d(H_N^{\rm G})\ne H_N^{\rm G})$. Both tend to zero. The observables have modulus at most one, so the original expectation difference is bounded by
\[
 N^{-1/6+o(1)}+2p_N^H+2p_N^G\longrightarrow0.
\]
This proves \eqref{W:eq:original-fourier-comparison} under the stated entry assumptions.
\end{proof}

\subsubsection{Tightness and identification of the coefficient limit}
\label{W:sec:tightness-and-identification-of-coefficients}
We can now complete the passage from expectation comparison to the coefficient limit. The modulus estimates give tightness for the potentials associated with both $H_N$ and $-H_N$.
Their increments control the off-diagonal coefficients through the reflection identities, and the characteristic function comparison identifies every subsequential limit of the
upper edge potential. Together these conclusions give the local coefficient theorem.

\begin{proof}[Proof of Theorem~\ref{W:thm:integrated-potential-and-offdiagonal-coefficients}]
Fix $T$ and $\kappa$, and let $\Gamma_\rho$ in this proof be computed from the modified polynomial moments, as in Theorem~\ref{W:thm:tightness-main}. Write
\[
 g_N^G(\kappa)=\sup_w\mathbb E[1-\Gamma_\rho(H_N^{\mathrm G},w)].
\]
Choose $a_\kappa,b_\kappa$ so that the scalar cutoff equals one on $[a_\kappa,b_\kappa]$, where $a_\kappa\downarrow0$ and $b_\kappa\uparrow\infty$. On $\mathcal D_N^G$ this is the
original Gram matrix. Consequently
\[
 \begin{split}
 g_N^G(\kappa)\le{}&p_N^G+
 \sup_w\mathbb P\bigl(\lambda_{\max}(G_\rho(H_N^{\mathrm G};w))>b_\kappa\bigr)\\
 &+\sup_w\mathbb P\bigl(\lambda_{\min}(G_\rho(H_N^{\mathrm G};w))^{-1}>a_\kappa^{-1}\bigr).
 \end{split}
\]
Lemma~\ref{W:lem:gaussian-gram-conditioning} proves tightness of both eigenvalue quantities, uniformly in $w$. It follows that
\[
 \lim_{\kappa\downarrow0}\limsup_{N\to\infty}g_N^G(\kappa)=0.
\]
The cutoff comparison gives
\[
 \sup_w\mathbb E[1-\Gamma_\rho(H_N^{\mathrm{tr}},w)]
 \le g_N^G(\kappa)+N^{-1/6+o(1)}.
\]

Let $\mathbf Y^{\sigma,\rm reg}$ be the values computed from the modified moments. On $\mathcal D_N^H\cap\{\Gamma_\rho>0\}$ they equal the original values of $Y_N^\sigma$. The
inequalities for the smooth function used to bound the modulus of continuity in \eqref{W:eq:tightness-sandwich} therefore give the pointwise bound
\[
 \begin{split}
 \mathbf1_{\{\omega_T(Y_N^\sigma,\delta)>u\}}
 \le{}&\mathbf1_{(\mathcal D_N^H)^c}
       +1-\Gamma_\rho(H_N^{\mathrm{tr}},w)\\
 &+\Gamma_\rho(H_N^{\mathrm{tr}},w)
               M(\mathbf Y^{\sigma,\rm reg}(H_N^{\mathrm{tr}},w)).
 \end{split}
\]
Each term on the right is nonnegative. Taking expectations, applying Theorem~\ref{W:thm:tightness-main}, and then using agreement on $\mathcal D_N^G$ yields
\[
 \begin{split}
 \sup_w\mathbb P(\omega_T(Y_N^\sigma,\delta)>u)
 \le{}&p_N^H+g_N^G(\kappa)+N^{-1/6+o(1)}\\
 &+p_N^G+
 \sup_w\mathbb P\bigl(
  \omega_{T+1}(Y_{N,G}^\sigma,\delta+2/m_N)>u/4\bigr).
 \end{split}
\]
The Gaussian modulus estimate \eqref{W:eq:tightness-gaussian-modulus} applies on $[0,T+1]$. After increasing its constant, this proves
\begin{equation}
 \begin{split}
 \sup_w\mathbb P\bigl(\omega_T(Y_N^\sigma,\delta)>u\bigr)
 \le{}& C_{\beta,T}u^{-4}(\delta+2/m_N)+g_N^G(\kappa)\\
       &+N^{-1/6+o(1)}+o(1).
 \end{split}
 \label{W:eq:tightness-modulus-probability}
\end{equation}
The last term includes truncation and spectral modification. Taking $N\to\infty$, then $\kappa\downarrow0$, and finally $\delta\downarrow0$ yields, for every $T,u>0$,
\[
 \lim_{\delta\downarrow0}\limsup_N\sup_w
 \mathbb P\bigl(\omega_T(Y_N^\sigma,\delta)>u\bigr)=0.
\]
The interpolated processes are continuous and start at zero. The tightness criterion in Theorem~7.3 of \cite{billingsley1999convergence}, after rescaling $[0,T]$ to $[0,1]$,
therefore gives tightness on every bounded interval for each deterministic sequence of starting vectors. The bounds are uniform in these vectors. If a tightness estimate is
required at all dimensions, the finitely many initial dimensions are covered by
\[
 \operatorname{Lip}(Y_N^\sigma|_{[0,T]})
 \le m_N^2(3\|H_N\|+2),
\]
whose right side is almost surely finite and independent of $w$. To pass to $C_{\rm loc}([0,\infty))$, fix $\eta>0$. For each integer $r\ge1$, tightness on $[0,r]$ provides a
compact set $\mathcal K_r\subset C([0,r])$ such that
\[
 \sup_{N,w}\mathbb P(Y_N^\sigma|_{[0,r]}\notin\mathcal K_r)
       \le\eta 2^{-r}.
\]
The finitely many initial dimensions are included by the preceding Lipschitz bound. The set
\[
 \mathcal K=\{f\in C_{\rm loc}([0,\infty)):
                    f|_{[0,r]}\in\mathcal K_r\text{ for every }r\}
\]
is compact: successively extract convergent subsequences on $[0,1],[0,2],\ldots$ and take the diagonal subsequence. The limits agree on the overlaps and define a continuous
function on $[0,\infty)$. The union bound gives $\sup_{N,w}\mathbb P(Y_N^\sigma\notin\mathcal K)\le\eta$. This proves tightness in the topology of uniform convergence on compact
intervals.

The reflection identities of Lemma~\ref{W:lem:tightness-reflection} preserve $\Gamma_\rho$ and $\chi_N$: with
$D_\rho=\operatorname{diag}(1,-1,\ldots)$,
\[
 G(-H)=D_\rho G(H)D_\rho,\qquad
 C(-H)=D_\rho C(H)D_\rho,
 \qquad b_j(-H)=b_j(H).
\]
Here $C$ is the triangular orthogonalization matrix, including its regularized version. Hence the row and column sums used in $Q_*$ and the total variation used in $\mathcal B$ are
unchanged. Negating the entries preserves their moment and tail assumptions, so the cutoff estimates apply to both signs.

Apply Theorem~\ref{W:thm:tightness-offdiagonal-coefficients} to the two reflected potential processes. The uniform modulus bounds just proved verify its hypothesis. In its
deterministic estimate, take $T=L+2$; then all indices in $B_N(L)$ are included, and
\[
 B_N(L)\le\frac{L+2+2h_N}{2}
 \sum_{\sigma=\pm}\omega_{L+3}(Y_N^\sigma,h_N)
 \xrightarrow{\mathbb P}0
\]
uniformly in deterministic unit $w$. Indeed, for every fixed $\delta>0$ and large $N$, each modulus on the right is bounded by $\omega_{L+3}(Y_N^\sigma,\delta)$; first take
$\limsup_N$ and then $\delta\downarrow0$ in the uniform modulus estimate. This proves \eqref{W:eq:tightness-offdiagonal-coefficients}.

For fixed $t_1,\ldots,t_q$ and real $\lambda_1,\ldots,\lambda_q$, let $\Phi_N(H)=\exp(\mathrm i\sum_j\lambda_jY_N^+(t_j))$. Corollary~\ref{W:cor:tightness-Fourier} implies
\begin{align*}
 |\mathbb E\Phi_N(H_N)-\mathbb E\Phi_N(H_N^{\mathrm G})|
 &\le |\mathbb E(\Gamma_\rho\Phi_N)(H_N)
           -\mathbb E(\Gamma_\rho\Phi_N)(H_N^{\mathrm G})|\\
 &\quad+\mathbb E(1-\Gamma_\rho)(H_N)
                 +\mathbb E(1-\Gamma_\rho)(H_N^{\mathrm G})\\
 &\le N^{-1/6+o(1)}+2g_N^G(\kappa)+o(1).
\end{align*}
Taking $N\to\infty$ and then $\kappa\downarrow0$, and using Theorem~\ref{W:thm:tightness-gaussian}, identifies the limiting characteristic function as
\[
 \exp\left\{\frac{\mathrm i}{2}\sum_j\lambda_jt_j^2
       -\frac2\beta\sum_{j,l}\lambda_j\lambda_l\min(t_j,t_l)\right\}.
\]
The exponent is that of $Y_\beta(t)=t^2/2+(2/\sqrt\beta)B(t)$, because its mean is $t^2/2$ and its covariance is $(4/\beta)\min(t,s)$. L\'evy's continuity theorem therefore
identifies the limiting joint distribution at every finite list of times. For any subsequence, tightness gives a further subsequence converging in $C_{\rm loc}([0,\infty))$.
Evaluation at each fixed time is continuous in this topology, so its limit has exactly these joint distributions. Values at nonnegative rational times determine a continuous
function, and their evaluation maps generate the Borel sigma field of this space. Thus every subsequential limit has the distribution of $Y_\beta$. The subsequence criterion for
weak convergence gives \eqref{W:eq:tightness-integrated-potential-limit}.

Let $d_{\rm BL}$ be the bounded Lipschitz metric on laws on $C_{\rm loc}([0,\infty))$. If convergence were not uniform, there would be $\varepsilon>0$, dimensions
$N_j\uparrow\infty$ and deterministic unit vectors $w_j$ such that
\[
 d_{\rm BL}\bigl(\mathcal L(Y_{N_j}^+(H_{N_j},w_j)),
                         \mathcal L(Y_\beta)\bigr)\ge\varepsilon.
\]
Extend $w_j$ arbitrarily to a deterministic starting sequence at all dimensions. The convergence just proved makes the displayed distances tend to zero, a contradiction.

For $0<\varepsilon<1$, choose an integer $R$ with $2^{-R}<\varepsilon/2$. The product metric satisfies
\[
 \begin{aligned}
 &\sup_w\mathbb P\left(\sum_{r\ge1}2^{-r}(1\wedge B_N(r))>\varepsilon\right)\\
 &\qquad\le\sum_{r=1}^R\sup_w\mathbb P(B_N(r)>\varepsilon/2)
 \longrightarrow0.
 \end{aligned}
\]
Write $d_0(b,0)=\sum_{r\ge1}2^{-r}(1\wedge|b_r|)$. If $\varphi$ is bounded by one and has Lipschitz constant at most one for the sum of a bounded metric on $C_{\rm loc}$ and $d_0$,
then
\[
 \begin{split}
 &\left|\mathbb E\varphi(Y_N^+,(B_N(r))_{r\ge1})
              -\mathbb E\varphi(Y_N^+,0)\right|\\
 &\qquad\le\mathbb E\min\{2,d_0((B_N(r))_{r\ge1},0)\}
 \le\varepsilon+2\mathbb P(d_0((B_N(r))_{r\ge1},0)>\varepsilon).
 \end{split}
\]
The last probability tends to zero uniformly over deterministic $w$. First let $N$ tend to infinity and then $\varepsilon$ tend to zero. The convergence of $Y_N^+$ proved above
consequently yields
\[
 \bigl(Y_N^+,(B_N(r))_{r\ge1}\bigr)
 \Longrightarrow\bigl(Y_\beta,(0)_{r\ge1}\bigr).
\]
The same inequality and the already established uniform bounded Lipschitz convergence of the first coordinate give uniform convergence of these joint distributions. This proves all
the claims of the local coefficient theorem.
\end{proof}

For the bounded functions used above, $|F_N|\le\Gamma_\rho$ and
\[
 |F_N-\chi_NF_N|\le\Gamma_\rho(1-\chi_N)=D_N^{\rm coef}.
\]
Corollary~\ref{W:cor:cumulant-comparison} and the cutoff estimate therefore give, uniformly in deterministic unit $w$,
\begin{align}
 |\mathbb EF_N(H_N^{\mathrm{tr}})-\mathbb EF_N({H_N^{\mathrm G}})|
 &\le|\mathbb E[\chi_NF_N](H_N^{\mathrm{tr}})
              -\mathbb E[\chi_NF_N]({H_N^{\mathrm G}})|\nonumber\\
 &\quad+\mathbb ED_N^{\rm coef}(H_N^{\mathrm{tr}})
              +\mathbb ED_N^{\rm coef}({H_N^{\mathrm G}})
       \le N^{-1/6+o(1)}.
 \label{W:eq:cumulant-unlocalized-ensembles}
\end{align}
For the corresponding functions of the original matrices, \eqref{W:eq:comparison-ensemble-transfer} adds at most twice each probability of disagreement under truncation or spectral
modification. Those probabilities tend to zero, yielding the comparison for the characteristic functions, $\Gamma_\rho$, and the two modulus functions.

The localization statistic additionally contains the weighted moments $F_r^*\Theta(H)F_r$ and the unnormalized trace $\operatorname{Tr}\Theta(H)$. Their entry derivatives are
estimated in Section~\ref{W:sec:localization-of-the-low-energy-spectral-subspace}.

\section{Spectral localization and the operator limit for Wigner matrices}\label{sec:spectral-localization-and-the-wigner-operator-limit}

Theorem~\ref{W:thm:integrated-potential-and-offdiagonal-coefficients} proves the local coefficient convergence required in \eqref{W:eq:operator-probabilistic-local}. The remaining
hypothesis of Theorem~\ref{W:thm:operator-probabilistic} is
\[
 \lim_{L\to\infty}\limsup_{N\to\infty}\sup_{\|w\|=1}
   \mathbb P\!\left(\left\|Q_{N,L}
                 \mathbf1_{(-\infty,E]}(\mathcal A_N)\right\|>\eta\right)=0
 \qquad(E\in\mathbb R,\ \eta>0),
\]
where $Q_{N,L}$ is defined in \eqref{W:eq:operator-coefficient-tail-notation} and $\mathcal A_N=N^{2/3}(2I-T_N)$ is the full tridiagonal matrix in
\eqref{W:eq:scaled-lanczos-operator}. Section~\ref{sec:localization-of-the-entire-low-energy-subspace} proves this estimate, including the eigenvectors in every later Lanczos
block. Section~\ref{W:sec:main-assembly} then applies Theorem~\ref{W:thm:operator-probabilistic} to obtain the operator limit for Wigner matrices.

\subsection{Localization of the low energy spectral subspace}
\label{sec:localization-of-the-entire-low-energy-subspace}
\label{W:sec:localization-of-the-low-energy-spectral-subspace}

The spectral measure at $w$ assigns zero weight to eigenvectors orthogonal to $w$, including eigenvectors in later Lanczos blocks. The unnormalized trace introduced in
Section~\ref{W:sec:localization-statistic} detects these vectors through their contribution to the full low energy spectral subspace.

We compare the bounded statistic formed from this trace and the cutoff in the eigenvalues of the Gram matrix. Section~\ref{W:sec:green-function-comparison} provides the calculation
for the terms containing $w$; the additional cubic trace terms are treated in Section~\ref{W:sec:localization-cubic-trace}. The comparison and the Gaussian trace bounds then give
the required estimate for $\|Q_{N,L}\Pi_N(E)\|$.

\subsubsection{The localization statistic}
\label{W:sec:localization-observable-and-derivatives}

We first identify the regularized polynomial projection with the projection onto the initial vectors of the completed Lanczos basis on the event where their coefficients agree. On
that event the deterministic trace inequality of Theorem~\ref{W:thm:gaussian-detection} applies to the full scaled matrix.

Fix $L>0$, $E\in\mathbb R$, $0<\kappa<1/4$, and the dimensions of the Gram matrices. For sufficiently large $N$ set
\[
 \begin{gathered}
 m_N=N^{1/3},\qquad \eta_{\mathrm e}=m_N^{-2},\qquad
 r=\lfloor m_NL\rfloor,\qquad r\le\rho\le\rho_*\asymp m_N,\\
 d=A_0\rho_*^{-2}(\log(N+2))^\gamma,\qquad 1<\gamma<2.
 \end{gathered}
\]
Constants may depend on these fixed parameters. The vector $w\in\mathbb F_\beta^N$ is deterministic and has norm one. Define
\[
 \vartheta_N(z)=\bigl(1+\exp(m_N^2(2-z)-E)\bigr)^{-2},
 \qquad \Theta(H)=\vartheta_N(H),
\]
and, for $k=r,\rho$, let
\[
 F_k(H)=(p_0(h_d(H))w,\ldots,p_{k-1}(h_d(H))w).
\]
The regularization $h_d$ and the scalar functions $a_\kappa$ and $\gamma_\kappa$ are those of Section~\ref{W:sec:spectral-regularization-and-matrix-derivatives}. In particular,
$h_d(x)=x$ for $|x|\le2+d$, $a_\kappa(x)=\chi_\kappa(x)/x$ for $x>0$, with $0\le\chi_\kappa\le1$, and $a_\kappa=0$ near zero. Both $a_\kappa(x)=x^{-1}$ and $\chi_\kappa(x)=1$ hold
on a neighbourhood of $[\kappa,\kappa^{-1}]$.

Put
\[
 G=F_\rho^*F_\rho,\qquad G_r=F_r^*F_r,\qquad
 M=F_r^*\Theta(H)F_r,\qquad \mathsf t=\Tr\Theta(H),\qquad
 \mathsf X(H,w)=(G,M,\mathsf t).
\]
The matrices $G$ and $M$ are real symmetric: their entries are spectral integrals of real functions, since $h_d(H)$ and $\Theta(H)$ commute. We regard $\mathsf X$ as a point of
$\operatorname{Sym}_\rho(\mathbb R)\times \operatorname{Sym}_r(\mathbb R)\times\mathbb R$, with $G_r$ always meaning the leading compression of $G$. On this space define
\begin{equation}
 C(\mathsf X)=\mathsf t-\Tr(a_\kappa(G_r)M),\qquad \Gamma_\rho(G)=\det\gamma_\kappa(G), \qquad \mathcal F(\mathsf X)=\Gamma_\rho(G)\widetilde\psi(C(\mathsf X)),
 \label{W:eq:localization-statistic}
\end{equation}
where $\widetilde\psi\in C_b^\infty(\mathbb R)$ agrees with $c/(1+c)$ on $[0,\infty)$ and has absolute value at most one. The matrix $\Theta(H)=\vartheta_N(H)$ is evaluated at the
original matrix $H$ throughout these definitions.

\begin{lemma}
\label{W:lem:localization-positivity}
For every Hermitian $H$, the matrix $P_\kappa=F_ra_\kappa(G_r)F_r^*$ is a positive contraction, and
\[
 0\le M\le G_r,\qquad
 C(\mathsf X(H,w))=\Tr((I-P_\kappa)\Theta(H))\ge0,
 \qquad 0\le\mathcal F(\mathsf X(H,w))\le1.
\]
On $\{\Gamma_\rho(G)>0\}$, $P_\kappa$ is the orthogonal projection onto $\operatorname{Ran}F_r$.
\end{lemma}
\begin{proof}
Take a singular value decomposition $F_r=\sum_{j=1}^{k}\sigma_j u_jv_j^*$, where $\sigma_j>0$. Functional calculus and $a_\kappa=0$ near zero give
\[
 P_\kappa=\sum_{j=1}^{k}\sigma_j^2a_\kappa(\sigma_j^2)u_ju_j^*
          =\sum_{j=1}^{k}\chi_\kappa(\sigma_j^2)u_ju_j^*.
\]
Thus $0\le P_\kappa\le I$. For every $v\in\mathbb C^r$,
\[
 v^*Mv=\|\Theta(H)^{1/2}F_rv\|^2\le\|F_rv\|^2=v^*G_rv.
\]
Cyclicity of the finite-dimensional trace yields the formula for $C$, and positivity follows from
\[
 \Tr((I-P_\kappa)\Theta)
   =\Tr\bigl(\Theta^{1/2}(I-P_\kappa)\Theta^{1/2}\bigr)\ge0.
\]
Equivalently,
\[
 C(\mathsf X(H,w))
 =\|(I-P_\kappa)^{1/2}\Theta(H)^{1/2}\|_{\mathrm{HS}}^2.
\]
The equality follows by taking the trace of the product of this operator with its adjoint. It uses only positivity of the two factors; $P_\kappa$ and $\Theta(H)$ need not commute.
The bounds on $\Gamma_\rho$ and $\widetilde\psi$ now give the bound for $\mathcal F$. If $\Gamma_\rho(G)>0$, every eigenvalue of $G$ belongs to $(\kappa,\kappa^{-1})$. Extending a
vector in $\mathbb C^r$ by zero gives $\kappa I_r<G_r<\kappa^{-1}I_r$. Hence $k=r$ and $\chi_\kappa(\sigma_j^2)=1$ for each $j$, so $P_\kappa=\sum_{j=1}^r u_ju_j^*$.
\end{proof}

The function $\mathcal F(\mathsf X(H,w))$ agrees with the localization statistic multiplied by its cutoffs in Section~\ref{W:sec:taylor-bounds-for-localization-statistic}. Its
coordinate cutoff is identically one on a neighbourhood of every moment vector obtained from a matrix. It consequently has the same matrix derivatives there. Taylor expansions
along arbitrary line segments in the full coordinate space use the smooth extension from that section.

We next relate this statistic to the tail in the tridiagonal coordinates of the spectral projection. Let $U_N$ be the unitary matrix whose columns are the Lanczos vectors of $H_N$
starting at $w$. Put
\[
 A_N^{\rm orig}=m_N^2(2I-H_N),\qquad
 \mathcal A_N=U_N^*A_N^{\rm orig}U_N,
 \qquad \Pi_N(E)=\mathbf1_{(-\infty,E]}(\mathcal A_N).
\]
For fixed $L$ and sufficiently large $N$, $r=\lfloor m_NL\rfloor<N$. Let $P_{N,L}=\sum_{j=1}^r e_je_j^*$, $Q_{N,L}=I-P_{N,L}$, and $\widehat P_{N,L}=U_NP_{N,L}U_N^*$. If $r\ge N$,
use the truncation from \eqref{W:eq:operator-coefficient-tail-notation}; then $Q_{N,L}=0$, so $\|Q_{N,L}\Pi_N(E)\|=0$.

Define $\Omega_d=\{\|H_N\|\le2+d/2\}$. On this event, $h_d(H_N)=H_N$. If also $\Gamma_\rho(G)>0$, the columns of $F_\rho$ are independent. Since $p_j$ is monic of degree $j$,
\[
 \operatorname{Ran}F_r
   =\operatorname{span}\{w,H_Nw,\ldots,H_N^{r-1}w\}.
\]
The independence of these columns rules out a Lanczos breakdown before step $r$. Section~\ref{W:sec:lanczos-tridiagonalization} then identifies this span with the first $r$ columns
of $U_N$. Hence $P_\kappa=\widehat P_{N,L}$ on this event. Functional calculus under unitary conjugation gives
\begin{align*}
 C(\mathsf X(H_N,w))
 &=\Tr\bigl((I-\widehat P_{N,L})\vartheta_N(H_N)\bigr)\\
 &=\Tr\bigl(Q_{N,L}\theta_E(\mathcal A_N)\bigr),\qquad
 \theta_E(x)=(1+e^{x-E})^{-2}.
\end{align*}
The full tridiagonal matrix is used in both expressions. If $(\lambda_{j,N},v_{j,N})_{j=1}^N$ is an orthonormal eigenbasis of $\mathcal A_N$, then on the same event
\begin{equation}
 C(\mathsf X(H_N,w))
 =\sum_{j=1}^N\theta_E(\lambda_{j,N})\|Q_{N,L}v_{j,N}\|^2.
 \label{W:eq:localization-full-eigenbasis}
\end{equation}
A unit eigenvector with $\lambda_{j,N}\le E$ and $P_{N,L}v_{j,N}=0$ contributes at least $1/4$, even if its first coordinate is zero. Replacing $C$ by $C/N$ would reduce this
contribution to $1/(4N)$. Vanishing of that normalized quantity would therefore not rule out one missing low energy direction.

Figure~\ref{W:fig:localization-later-block} illustrates this contribution for an eigenvector in a subsequent invariant block. It is the reason the trace in
\eqref{W:eq:localization-full-eigenbasis} is taken over the full tridiagonal matrix.

\begin{figure}[H]
\centering
\begin{tikzpicture}[x=1cm,y=1cm,every node/.style={font=\footnotesize,inner sep=2pt},>=Stealth]
 \node at (2.1,5.15) {The completed tridiagonal matrix $\mathcal A_N$};
 \fill[black!10] (0,1.3) rectangle (2.9,4.2);
 \fill[black!30] (0,2.75) rectangle (1.45,4.2);
 \fill[black!20] (2.9,0) rectangle (4.2,1.3);
 \draw (0,0) rectangle (4.2,4.2);
 \draw (0,1.3)--(4.2,1.3);\draw (2.9,0)--(2.9,4.2);
 \draw[densely dashed] (0,2.75)--(2.9,2.75);\draw[densely dashed] (1.45,1.3)--(1.45,4.2);
 \node at (.72,3.48) {$r\times r$};
 \node at (2.2,2.1) {$\mathcal A_N^{(1)}$};
 \node at (3.55,.65) {$\mathcal A_N^{(2)}$};
 \node at (3.55,2.75) {$0$};\node at (1.45,.65) {$0$};
 \node[above] at (1.45,4.25) {$r$};
 \node[above] at (2.9,4.25) {$D$};
 \node[above] at (4.2,4.25) {$N$};
 \node[below] at (3.55,-.08) {$D+1,\ldots,N$};
 \node[align=center] at (2.1,-.95) {$D=\dim\mathcal K_N(H_N,w)$\\$r\le\rho\le D<N$};
 \node[draw,minimum width=1.15cm,minimum height=2.9cm] at (5.6,2.75) {$0$};
 \node[draw,fill=black!20,minimum width=1.15cm,minimum height=1.3cm] at (5.6,.65) {$v^{(2)}$};
 \node[above] at (5.6,4.45) {$v$};
 \node[anchor=west] at (7.15,3.9) {$\mathcal A_Nv=\lambda v,\quad\lambda\le E,\quad\|v\|=1$};
 \node[anchor=west] at (7.15,2.85) {$e_1^*v=0,\qquad P_{N,L}v=0$};
 \node[anchor=west] at (7.15,1.8) {$Q_{N,L}v=v$};
 \node[anchor=west] at (7.15,.65) {$\theta_E(\lambda)\|Q_{N,L}v\|^2\ge\frac14$};
 \node[anchor=west] at (7.15,-.6) {$C=\displaystyle\sum_j\theta_E(\lambda_{j,N})\|Q_{N,L}v_{j,N}\|^2$};
\end{tikzpicture}
\caption{A possible later Lanczos block in the localization argument. On $\Omega_d\cap\{\Gamma_\rho>0\}$, the first $r$ coordinates lie in the first cyclic subspace, whose dimension is denoted here by $D$. The darker square marks the coordinates selected by $P_{N,L}$ within the first block. A unit eigenvector supported in a later block has zero first coordinate and belongs to $\operatorname{Ran}Q_{N,L}$. If its eigenvalue is at most $E$, it contributes at least $1/4$ to the full trace in \eqref{W:eq:localization-full-eigenbasis}. The spectral measure at the starting vector assigns it weight zero. The bound for arbitrary vectors in the low energy subspace is \eqref{W:eq:localization-detection}.}
\label{W:fig:localization-later-block}
\end{figure}

Since $\theta_E(x)\ge1/4$ for every $x\le E$, the spectral theorem implies $\theta_E(\mathcal A_N)\ge\Pi_N(E)/4$. Multiply both sides on the left and right by $Q_{N,L}$ and take
the trace:
\[
 C\ge\tfrac14\Tr(Q_{N,L}\Pi_N(E)Q_{N,L})
   =\tfrac14\|\Pi_N(E)Q_{N,L}\|_{\rm HS}^2
   \ge\tfrac14\|Q_{N,L}\Pi_N(E)\|^2.
\]
For $c\ge0$, the function $c/(1+c)$ is increasing. Therefore, on $\Omega_d\cap\{\Gamma_\rho>0\}$,
\[
 \|Q_{N,L}\Pi_N(E)\|>\eta
 \quad\Longrightarrow\quad
 \frac{C}{1+C}>\frac{\eta^2}{4+\eta^2}.
\]
Multiplying by $\Gamma_\rho$ and adding $1-\Gamma_\rho$ proves the pointwise inequality
\begin{equation}
 \mathbf1_{\{\|Q_{N,L}\Pi_N(E)\|>\eta\}}
 \le\mathbf1_{\Omega_d^c}+1-\Gamma_\rho+
       \frac{4+\eta^2}{\eta^2}\mathcal F(\mathsf X(H_N,w)).
 \label{W:eq:localization-detection}
\end{equation}
On $\{\Gamma_\rho=0\}$ its right side is at least one, so the inequality holds there as well. In particular the projection on the left contains the eigenvectors of every subsequent
Lanczos block. The trace inequality also controls every unit vector in $\operatorname{Ran}\Pi_N(E)$, including linear combinations in a multiple eigenspace, since its operator norm
conclusion is $\|Q_{N,L}\Pi_N(E)\|^2\le4C$.

For subsequent blocks, the detection bound can be stated without choosing an eigenbasis. Let $\mathcal K_N(H_N,w)=\operatorname{span}\{w,H_Nw,\ldots,H_N^{N-1}w\}$ be the first
cyclic subspace, and let $P_N^{\rm cyc}$ be its orthogonal projection in the original coordinates. This is a Borel matrix function: for the matrix of successive powers applied to
the starting vector,
$K_N=[w,H_Nw,\ldots,H_N^{N-1}w]$,
\[
 P_N^{\rm cyc}=\lim_{\varepsilon\downarrow0}
 K_N(K_N^*K_N+\varepsilon I)^{-1}K_N^*.
\]
Indeed its nonzero singular values contribute the multipliers $\sigma^2/(\sigma^2+\varepsilon)\to1$, while its orthogonal complement contributes zero. Each approximating matrix is
continuous in $(H_N,w)$. The subspace $\mathcal K_N$ and its orthogonal complement reduce $H_N$. On $\Omega_d\cap\{\Gamma_\rho>0\}$ one has
$\operatorname{Ran}\widehat P_{N,L}\subset\mathcal K_N$, and hence
\[
 C\ge\tfrac14\operatorname{Tr}\!\left[
 (I-P_N^{\rm cyc})\mathbf1_{(-\infty,E]}(A_N^{\rm orig})\right].
\]
The product in the trace is an orthogonal projection, so its trace is a nonnegative integer. It is positive exactly when the low energy spectral subspace has a component in a
subsequent block. On that event $C\ge1/4$ and $\mathcal F\ge\Gamma_\rho/5$. If the cyclic subspace has dimension less than $r$, the Gram matrix is singular on $\Omega_d$ and
$\Gamma_\rho=0$; the term $1-\Gamma_\rho$ in the detection inequality covers that case.

It therefore suffices to compare the expectation of $\mathcal F$ under the original ensemble with its expectation under the GOE or GUE, in addition to the comparison of the cutoff
expectations already proved. Corollary~\ref{W:cor:localization-comparison} gives this comparison.

\begin{corollary}
\label{W:cor:localization-comparison}
For the uncentered truncated matrix $H_N^{\mathrm{tr}}$ of Section~\ref{W:sec:moments-of-truncated-entries} and the GOE or GUE matrix ${H_N^{\mathrm G}}$, with every auxiliary
parameter fixed,
\begin{equation}
 \sup_{\|w\|=1}\left|
  \E\mathcal{F}(\mathsf{X}(H_N^{\mathrm{tr}},w))-
  \E\mathcal{F}(\mathsf{X}({H_N^{\mathrm G}},w))\right|\le N^{-1/6+o(1)}.
 \label{W:eq:localization-comparison}
\end{equation}
For the original matrix and its associated moment variables, the preceding estimate holds with an additional $o(1)$ error.
\end{corollary}

Let $\nu_N^{\rm b}$ be the homogeneous bounded entry distribution of Section~\ref{W:sec:bounded-entries-with-prescribed-moments}, and let $\nu^G$ be the corresponding standardized
Gaussian distribution. Write $\Delta\nu=\nu_N^{\rm b}-\nu^G$. For $t\in[0,1]$, let $H(t)$ have independent upper-triangular entries with standardized off-diagonal distribution
$t\nu_N^{\rm b}+(1-t)\nu^G$ and the diagonal entry distributions from Section~\ref{W:sec:taylor-remainders-and-bounded-entries}. All standardized moments of each fixed order are
uniformly bounded. The diagonal variance is $(2-t)/N$ for $\beta=1$ and $1/N$ for $\beta=2$. For $a<b$ put
\[
 K_{ab}=H(t)^{[ab\leftarrow0]},\qquad
 T_{ab}(\zeta)=\zeta E_{ab}+\overline\zeta E_{ba}.
\]
For a fixed Hermitian base $K$ define $\ell_{ab},q_{ab},c_{ab}$ by
\begin{equation}
 \mathsf X(K+uN^{-1/2}T_{ab}(\zeta),w)
 =\mathsf X(K,w)+u\ell_{ab}(K;\zeta)+u^2q_{ab}(K;\zeta)
       +u^3c_{ab}(K;\zeta)+O(u^4).
 \label{W:eq:localization-moment-coefficients}
\end{equation}
Thus $q$ and $c$ include their Taylor factorials. At $X=\mathsf X(K_{ab},w)$ the cubic coefficient of the composition is
\begin{equation}
 D\mathcal F(X)[c]+D^2\mathcal F(X)[\ell,q]
                  +\tfrac16D^3\mathcal F(X)[\ell,\ell,\ell].
 \label{W:eq:localization-complete-cubic}
\end{equation}
Let $Z_N^{(1)}(t,w)$ be the sum over $a<b$ of the integral of $D\mathcal F(X)[c]$ against $\Delta\nu$, and let $Z_N^{(3)}(t,w)$ be the sum over $a<b$ of the integral of
\eqref{W:eq:localization-complete-cubic} against $\Delta\nu$. The integration variable $\zeta$ is separate from the sampled entries of $H(t)$.

\begin{theorem}
\label{W:thm:localization-cubic}
For each fixed choice of the preceding parameters,
\begin{equation}
 \sup_{t\in[0,1]}\sup_{\|w\|=1}
 \bigl(\mathbb E|Z_N^{(1)}(t,w)|+\mathbb E|Z_N^{(3)}(t,w)|\bigr)
 \le N^{-1/6+o(1)}.
 \label{W:eq:localization-cubic-estimate}
\end{equation}
\end{theorem}

\subsubsection{Derivatives of the projection and the statistic}
\label{W:sec:projection-and-statistic-derivatives}

Differentiating the localization statistic requires derivatives of both the projection onto the polynomial subspace and the matrix function in the trace. The calculation therefore
begins with the projection derivative, before the two sources of variation are combined.

On the set where the eigenvalues of $G$ lie in $(\kappa,\kappa^{-1})$, put $F=F_r(H)$, $P=F(F^*F)^{-1}F^*$, and $\Theta=\vartheta_N(H)$. The unregularized localization statistic
there is exactly
\[
 C(H)=\Tr\bigl((I-P(H))\Theta(H)\bigr).
\]
On this open set the equality also holds after differentiation. With $A=(F^*F)^{-1}$, \eqref{W:eq:resolvent-projection-first} gives the projection derivative
\begin{equation}
 D P=(I-P)(D F)AF^*+FA(D F)^*(I-P).
 \label{W:eq:localization-projection-first-variation}
\end{equation}
For a fixed self-adjoint matrix $W$ supported on the selected entry, primes denote differentiation at $t=0$ along $H+tW$. Differentiating the two factors in
$\operatorname{Tr}((I-P)\Theta)$ gives
\begin{align}
 C'&=\Tr((I-P)\Theta')-\Tr(P'\Theta),\notag\\
 C''&=\Tr((I-P)\Theta'')-2\Tr(P'\Theta')-\Tr(P''\Theta),\notag\\
 C^{(3)}&=\Tr((I-P)\Theta^{(3)})
       -3\Tr(P'\Theta'')-3\Tr(P''\Theta')-\Tr(P^{(3)}\Theta),\notag\\
 C^{(4)}&=\Tr((I-P)\Theta^{(4)})-4\Tr(P'\Theta^{(3)})
       -6\Tr(P''\Theta'')\notag\\
 &\qquad-4\Tr(P^{(3)}\Theta')-\Tr(P^{(4)}\Theta).
 \label{W:eq:resolvent-localization-four-derivatives}
\end{align}
For mixed derivatives, use the following notation. For a label set $S$, write $F_S=D_SF$, $G_S=D_S(F^*F)$, $A=(F^*F)^{-1}$, $A_S=D_SA$, and $\Theta_S=D_S\Theta$. Empty subscripts
denote the original factors. Differentiating $G_r=F^*F$ gives
\[
 G_S=\sum_{U\sqcup V=S}F_U^*F_V.
\]
Here $G_S=D_SG_r$ is the derivative of the $r$ by $r$ Gram matrix. Its inverse has the derivatives \eqref{W:eq:gram-inverse-labels}, which in the present notation read
\begin{equation}
 A_S=\sum_{k=1}^{|S|}(-1)^k
 \sum_{\substack{S_1\dot\cup\cdots\dot\cup S_k=S\\S_j\ne\varnothing}}
 AG_{S_1}A\cdots G_{S_k}A.
 \label{W:eq:localization-inverse-partitions}
\end{equation}
The inner sum is over ordered partitions, and $A_\varnothing=A$. Substituting these derivatives in \eqref{W:eq:resolvent-projection-labels} gives, in the present notation,
\begin{equation}
 P_S=\sum_{U\dot\cup V\dot\cup W=S}F_UA_VF_W^*.
 \label{W:eq:localization-projection-partitions}
\end{equation}
Every label is assigned to one of the three factors, giving $3^{|S|}$ terms. This specifies the nine, twenty seven, and eighty one placements at orders two, three, and four,
respectively. For a repeated direction, write $G_j=\mathrm d^jG_r/\mathrm du^j$. The inverse formula at fourth order reads
\begin{align*}
 A^{(4)}={}&-AG_4A
 +4(AG_1AG_3A+AG_3AG_1A)+6AG_2AG_2A\\
 &-12(AG_1AG_1AG_2A+AG_1AG_2AG_1A+AG_2AG_1AG_1A)\\
 &+24AG_1AG_1AG_1AG_1A.
\end{align*}
The coefficients count the assignments of four labelled derivatives to ordered blocks of sizes $(4)$, $(1,3)$, $(2,2)$, $(2,1,1)$, and $(1,1,1,1)$, respectively. Substitution in
\eqref{W:eq:localization-projection-partitions} gives $P^{(4)}$. The repeated direction formulas for $C$ follow by differentiating $\Tr((I-P)\Theta)$ and using the binomial
coefficients $\binom qj$. Where $h_d$ is the identity on a neighbourhood of the spectrum of $H$, choose a contour enclosing the spectrum and avoiding the poles of $\vartheta_N$.
For nonempty $S$,
\begin{align*}
 F_Se_j&=\frac1{2\pi\mathrm i}\int p_{j-1}(z)
       \sum_{\sigma\in\operatorname{Perm}(S)}
       RW_{\sigma(1)}R\cdots W_{\sigma(|S|)}Rw\,\mathrm dz,\\
 \Theta_S&=\frac1{2\pi\mathrm i}\int\vartheta_N(z)
       \sum_{\sigma\in\operatorname{Perm}(S)}
       RW_{\sigma(1)}R\cdots W_{\sigma(|S|)}R\,\mathrm dz.
\end{align*}
The smooth modified column functions away from this region have the corresponding divided difference derivatives from Lemma~\ref{W:lem:smooth-functional-calculus}; the displayed
product identities are unchanged. In particular
\[
 \Theta_{bc}=\frac1{2\pi\mathrm i}\int\vartheta_N(z)
       (RW_bRW_cR+RW_cRW_bR)\,\mathrm dz.
\]
Combining these derivatives with the product rule for the trace gives
\begin{equation}
 D_SC=\operatorname{Tr}((I-P)\Theta_S)
       -\sum_{\varnothing\ne U\subseteq S}
                    \operatorname{Tr}(P_U\Theta_{S\setminus U}).
 \label{W:eq:resolvent-localization-labels}
\end{equation}
If the three directions coincide, the two subset sums each have three terms, giving the two coefficients $3$ in \eqref{W:eq:resolvent-localization-four-derivatives}. Division by
$3!$ then gives the homogeneous cubic Taylor coefficient. The same labelled formula for $|S|=4$ gives the coefficients $4,6,4,1$ and the fourth Taylor coefficient is its $1/4!$
multiple.

For the product involving $\Gamma_\rho$, let $\Pi(V)$ denote the set of unordered partitions of a finite label set $V$. The labelled product and chain rules give
\begin{equation}
 D_S\{\Gamma_\rho\widetilde\psi(C)\}
 =\sum_{U\sqcup V=S}(D_U\Gamma_\rho)
       \sum_{\pi\in\Pi(V)}\widetilde\psi^{(|\pi|)}(C)
                       \prod_{B\in\pi}D_BC.
 \label{W:eq:resolvent-localization-cutoff-derivatives}
\end{equation}
For $V=\varnothing$, use its single empty partition, with empty product equal to one and $\widetilde\psi^{(0)}=\widetilde\psi$. For $|S|=3$, the scalar composition contains
$\widetilde\psi'C^{(3)}+3\widetilde\psi''C'C'' +\widetilde\psi^{(3)}(C')^3$; the outer product adds all placements of one, two, and three derivatives on $\Gamma_\rho$. This is the
formula for the cubic coefficient in \eqref{W:eq:localization-complete-cubic}. The cutoff in the coordinate variables used to define the global extension is likewise differentiated
when a remainder of the expansion is estimated along a moment segment that need not arise from a matrix.

The same formulas follow by differentiating $C=\mathsf t-\Tr(G_r^{-1}M)$ and substituting $G_r=F^*F$, $M=F^*\Theta F$. Thus the moment representation and the projection
representation give identical derivatives. The smooth scalar function $a_\kappa$ continues to define the statistic at singular Gram matrices.

For $v=(\dot G,\dot M,\dot{\mathsf t})$ set
\begin{equation}
 \mathfrak{n}(v)=\|\dot G\|_1+\|\dot M\|_1+|\dot{\mathsf t}|,
 \label{W:eq:localization-norm}
\end{equation}
where matrix norms are trace norms. This is the reduced form of the spectral localization norm in Section~\ref{W:sec:taylor-bounds-for-localization-statistic}; the original norm
also includes contributions from the coordinates used in the cutoff.

\Needspace{6\baselineskip}

\begin{lemma}
\label{W:lem:localization-derivatives}
For every Hermitian matrix $H_0$, put $X_0=\mathsf X(H_0)$. For every real coordinate vector $X'$ with $\mathfrak{n}(X'-X_0)\le1$, and $j=1,2,3$,
\begin{equation}
 |D_{\mathsf X}^j\mathcal{F}(X')[v_1,\ldots,v_j]|
       \le C_{j,\kappa}\prod_l\mathfrak{n}(v_l).
 \label{W:eq:localization-local-derivatives}
\end{equation}
\end{lemma}

For two moment vectors obtained from matrices, $X=\mathsf X(H_1)$ and $Y=\mathsf X(H_2)$, the following estimate holds.
\begin{corollary}
\label{W:cor:localization-gradient-bound}
For every direction $v$,

\begin{equation}
 |(D_{\mathsf X}\mathcal{F}(X)-D_{\mathsf X}\mathcal{F}(Y))[v]|
       \le C_\kappa\mathfrak{n}(X-Y)\mathfrak{n}(v).
 \label{W:eq:localization-gradient-bound}
\end{equation}
\end{corollary}
\begin{proof}[Proof of Lemma~\ref{W:lem:localization-derivatives} and Corollary~\ref{W:cor:localization-gradient-bound}]
Consider a term in the product rule for $D^j\mathcal F$. If an eigenvalue of $G(X')$ lies outside $[\kappa,\kappa^{-1}]$, it has positive distance from the compact support of
$\gamma_\kappa$. Weyl's inequality shows that the corresponding eigenvalue stays outside that support in a neighbourhood of $G(X')$. Throughout this neighbourhood the product
$\Gamma_\rho(G)=\prod_{k=1}^{\rho}\gamma_\kappa(\lambda_k(G))$ is zero. All its derivatives there vanish. We may therefore restrict every nonzero term to
\[
 \kappa I\le G(X')\le\kappa^{-1}I,\qquad
 \|G(X_0)\|\le\|G(X')\|+\|G(X_0)-G(X')\|_1
                         \le\kappa^{-1}+1.
\]
Positivity at a matrix gives $\|M(X_0)\|\le\|G_r(X_0)\|$, and hence $\|M(X')\|\le\kappa^{-1}+2$. On $\kappa I\le G(X')\le\kappa^{-1}I$, let $A=G_r(X')^{-1}$. Its inverse
derivatives are
\[
 D_{G}^jA[B_1,\ldots,B_j]
   =(-1)^j\sum_{\pi\in S_j}AB_{\pi(1)}A\cdots B_{\pi(j)}A.
\]
For a direction $v=(\dot G,\dot M,\dot{\mathsf t})$, the first variation is
\[
 D_{\mathsf X}C[v]=\dot{\mathsf t}
       +\Tr(A\dot G_r AM)-\Tr(A\dot M).
\]
For $G_j=\dot G_r[v_j]$, $M_j=\dot M[v_j]$, the next two variations are
\begin{align*}
 D_{12}C={}&-\operatorname{Tr}(AG_1AG_2AM+AG_2AG_1AM)\\
          &+\operatorname{Tr}(AG_1AM_2+AG_2AM_1),\\
 D_{123}C={}&\sum_{\sigma\in\mathfrak S_3}
          \operatorname{Tr}(AG_{\sigma(1)}AG_{\sigma(2)}AG_{\sigma(3)}AM)\\
          &-\sum_{j=1}^3\sum_{(a,b)\in\operatorname{Perm}(\{1,2,3\}\setminus\{j\})}
          \operatorname{Tr}(AG_aAG_bAM_j).
\end{align*}
These identities follow from $D_jA=-AG_jA$; the first derivatives of $G,M,\mathsf t$ are constant in these coordinates, and their derivatives of order at least two vanish. Set
$g_j=\|G_j\|_1$, $m_j=\|M_j\|_1$, $\tau_j=|\dot{\mathsf t}[v_j]|$, and $M_*=\kappa^{-1}+2$. Putting one differentiated factor in trace norm and the other factors in operator norm
gives
\begin{align*}
 |D_1C|&\le\tau_1+\kappa^{-2}M_*g_1+\kappa^{-1}m_1,\\
 |D_{12}C|&\le2\kappa^{-3}M_*g_1g_2
                  +\kappa^{-2}(g_1m_2+g_2m_1),\\
 |D_{123}C|&\le6\kappa^{-4}M_*g_1g_2g_3
                  +2\kappa^{-3}\sum_{j=1}^3m_j\prod_{k\ne j}g_k.
\end{align*}
For example, $|\operatorname{Tr}(AG_1AG_2AM)|\le\kappa^{-3}M_*g_1g_2$ and $|\operatorname{Tr}(AG_1AG_2AM_3)|\le\kappa^{-3}g_1g_2m_3$. Each line is at most
$C_{q,\kappa}\prod_j\mathfrak n(v_j)$. The scalar derivatives through order three are
\begin{align*}
 D_1(\widetilde\psi\circ C)&=\widetilde\psi'C_1,\\
 D_{12}(\widetilde\psi\circ C)&=\widetilde\psi'C_{12}
                                      +\widetilde\psi''C_1C_2,\\
 D_{123}(\widetilde\psi\circ C)&=\widetilde\psi'C_{123}
       +\widetilde\psi''(C_{12}C_3+C_{13}C_2+C_{23}C_1)
       +\widetilde\psi^{(3)}C_1C_2C_3.
\end{align*}
Multiply by $\Gamma_\rho$ using $D_S\mathcal F=\sum_{U\sqcup V=S}(D_U\Gamma_\rho) D_V(\widetilde\psi\circ C)$. The determinant trace norm bound gives
$|D_U\Gamma_\rho|\le C_{|U|,\kappa}\prod_{j\in U}\|\dot G[v_j]\|_1$. Together with the displayed scalar formulas this proves \eqref{W:eq:localization-local-derivatives}.

For $\Delta=Y-X$ with $\mathfrak n(\Delta)\le1$, the fundamental theorem of calculus and the second derivative bound yield
\[
 |(D\mathcal F(Y)-D\mathcal F(X))[v]|
 =\left|\int_0^1D^2\mathcal F(X+t\Delta)[\Delta,v]\,\mathrm dt\right|
 \le C_\kappa\mathfrak n(\Delta)\mathfrak n(v).
\]
If $\mathfrak n(\Delta)>1$, apply the first derivative bound at $X=\mathsf X(H_1)$ and $Y=\mathsf X(H_2)$:
\[
 |(D\mathcal F(Y)-D\mathcal F(X))[v]|
 \le |D\mathcal F(Y)[v]|+|D\mathcal F(X)[v]|
 \le2C_\kappa\mathfrak n(v)
 \le2C_\kappa\mathfrak n(\Delta)\mathfrak n(v).
\] This proves
\eqref{W:eq:localization-gradient-bound} in both cases.
\end{proof}

At a point associated with a matrix where a first derivative can be nonzero, let $A=G_r^{-1}$, $a=\Gamma_\rho\widetilde\psi'(C)$, and define the real symmetric matrix $J_\Gamma$ by
$D_{G}\Gamma_\rho[B]=\Tr(J_\Gamma B)$. To identify this matrix explicitly, take $G=O\operatorname{diag}(\lambda_1,\ldots,\lambda_\rho)O^{\mathsf T}$. The derivative of the
determinant is $D\det(B)[D]=\Tr(\operatorname{adj}(B)D)$, including singular $B$. The diagonal entries of the functional calculus derivative of $\gamma_\kappa$ in this basis are
$\gamma_\kappa'(\lambda_j)(O^{\mathsf T}DO)_{jj}$. Consequently
\[
 J_\Gamma=O\operatorname{diag}\left(
 \gamma_\kappa'(\lambda_j)\prod_{k\ne j}\gamma_\kappa(\lambda_k)
                          \right)_{j=1}^{\rho}O^{\mathsf T},
 \qquad \|J_\Gamma\|\le\|\gamma_\kappa'\|_\infty.
\]
The diagonal coefficients agree on each repeated eigenspace, so the expression is independent of the choice of its orthonormal basis. The first derivative of the statistic is
\begin{equation}
 D_{\mathsf X}\mathcal{F}[v]={}\widetilde\psi(C)\Tr(J_\Gamma\dot G) +a\left\{\dot{\mathsf t} +\Tr(A\dot G_r A M)-\Tr(A\dot M)\right\}.
 \label{W:eq:localization-response}
\end{equation}
All coefficients in this formula are evaluated at the same point associated with a matrix. Outside the region in which $\Gamma_\rho$ or one of its first derivatives is nonzero, set
the following coefficients equal to zero. On that region $G_r$ has eigenvalues in $[\kappa,\kappa^{-1}]$, where the inverse formulas apply. Let
$\boldsymbol p_j(z)=(p_0(z),\ldots,p_{j-1}(z))^{\mathsf T}$ and set
\begin{align}
 P_0(X,z)&=\widetilde\psi(C)\boldsymbol p_\rho(z)^{\mathsf T}
                 J_\Gamma\boldsymbol p_\rho(z)
       +a\boldsymbol p_r(z)^{\mathsf T}AMA\boldsymbol p_r(z),
       \label{W:eq:localization-P0}\\
 P_1(X,z)&=-a\boldsymbol p_r(z)^{\mathsf T}A\boldsymbol p_r(z).
       \label{W:eq:localization-P1}
\end{align}
Here $M=F_r^*\Theta(H)F_r$ and $\mathsf t=\operatorname{Tr}\Theta(H)$, as defined above. The coefficients of $P_0$ and $P_1$ are evaluated at $X=\mathsf X(H,w)$ throughout the
following contour formulas.
\subsubsection{Contour estimates for \texorpdfstring{$\vartheta_N$}{vartheta N}}
\label{W:sec:localization-contour-estimates}

Integration by parts in the contour representation of the trace differentiates the scalar function. Its total variation on the real axis suggests the needed estimate, but the same
control must be established on the contour used by the resolvent formula.

On the real axis $\vartheta_N$ increases from zero to one, and
\[
 \vartheta_N'(x)=
 \frac{2m_N^2\exp(m_N^2(2-x)-E)}
 {(1+\exp(m_N^2(2-x)-E))^3},\qquad
 \int_{\mathbb R}\vartheta_N'(x)\,\dd x=1.
\]
The next estimate gives the corresponding bound on the contour used for the matrix functional calculus.

Let $s_*=4\sqrt d$, $b_*=2\cosh s_*$, and
\[
 y_*(x)=2\sinh s_*\sqrt{1-x^2/b_*^2},\qquad
 y(x)=\min\{y_*(x),a_0({\eta_{\mathrm e}}+|x-2|),a_0({\eta_{\mathrm e}}+|x+2|)\},
\]
with $a_0>0$ fixed sufficiently small relative to $1+|E|$. Let $\mathcal C_N$ be the positively oriented curve formed by $x\pm\ii y(x)$, $-b_*\le x\le b_*$. This is the contour of
Section~\ref{W:sec:contour-integrals-for-moment-increments}. It encloses the confined spectrum and avoids every pole of $\vartheta_N$. Put
\begin{equation}
 \mathcal P_\rho(z)=\frac{2e^{(\rho+1)s_*}}{|z^2-4|^{1/2}}.
 \label{W:eq:localization-polynomial-bound}
\end{equation}
Then $|p_j(z)|\le\mathcal P_\rho(z)$ for $j<\rho$.

\Needspace{6\baselineskip}

\begin{lemma}
\label{W:lem:localization-contour}
Uniformly for sufficiently large $N$,
\begin{align}
 \int_{\mathcal C_N}\mathcal P_\rho(z)^2|\dd z|&=N^{o(1)},
       &\sup_{\mathcal C_N}|\vartheta_N|&\le C_E,
       \label{W:eq:localization-contour-mass}\\
 \int_{\mathcal C_N}|\vartheta_N(z)|\,|\dd z|&\le C_Ed,
       &\int_{\mathcal C_N}|\vartheta_N'(z)|\,|\dd z|&\le C_E.
       \label{W:eq:localization-theta-integrals}
\end{align}
At every moment vector obtained from a matrix,
\begin{equation}
 |P_0(X,z)|+|P_1(X,z)|
      \le C_\kappa {m_N}\mathcal P_\rho(z)^2,
 \qquad |a(X)|\le1.
 \label{W:eq:localization-derivative-kernel-bound}
\end{equation}
\end{lemma}
\begin{proof}
The poles of $\vartheta_N$ have real part $2-E{\eta_{\mathrm e}}$ and imaginary parts equal to odd multiples of $\pi {\eta_{\mathrm e}}$. At that real part, the contour height is
at most $a_0{\eta_{\mathrm e}}(1+|E|)<\pi {\eta_{\mathrm e}}/2$. Thus none is enclosed. Write $u={m_N}^2(2-z)-E$. If $|\Re u|\le1$, the contour gives $|\Im u|<\pi/2$ after reducing
$a_0$. If $|\Re u|>1$, use $|1+e^u|\ge|e^{\Re u}-1|$. These observations give on the entire curve
\[
 |1+e^u|\ge c_E(1+e^{\Re u}),\qquad
 |\vartheta_N'(z)|
   =\frac{2{m_N}^2|e^u|}{|1+e^u|^3}
   \le C_E{m_N}^2e^{-c|\Re u|}.
\]
We verify which parts of the curve can meet the transition. For $x=2-v$ with $0\le v\le1$, the ellipse formula gives
\[
 y_*(x)^2=4\sinh^2s_*\left(1-\frac{(2-v)^2}{4\cosh^2s_*}\right)
          \le C s_*^2(s_*^2+v).
\]
If the straight line $a_0(\eta_{\mathrm e}+v)$ is below the ellipse, then $a_0^2v^2\le Cs_*^2(s_*^2+v)$. Dividing by $s_*^4$ and solving this quadratic inequality gives
$v\le C_{a_0}s_*^2$. The identity $y_*(-x)=y_*(x)$ gives the same bound for $x=-2+v$. For $|x-2|=O(\eta_{\mathrm e})$, one has $y_*(x)\ge c s_*^2$ and
$a_0(\eta_{\mathrm e}+|x-2|)=O(\eta_{\mathrm e})$. Since $s_*^2/\eta_{\mathrm e}\to\infty$, the contour is a straight piece there. On an exterior elliptic arc, $x>2$ implies
$|\theta|\le C s_*$ in its angular parametrization. Therefore
\[
 \operatorname{length}\{z\in\mathcal C_N:\Re z>2\}
 \le C\int_0^{Cs_*}(s_*^2+\theta^2)^{1/2}\,\dd\theta
 \le C s_*^2=O(d).
\]
The transition region and its nearby edge pieces are straight pieces with bounded slope; integration with respect to their real part yields
$\int {m_N}^2e^{-c|{m_N}^2(2-x)-E|}\,\dd x=O_E(1)$. All remaining elliptic pieces are at real distance at least $c d$ from the transition, or are on the opposite side of the
spectrum. Their contribution is bounded by $C {m_N}^2\exp(-c d/{\eta_{\mathrm e}})$ times a bounded length. Since $d/\eta_{\mathrm e}\ge c_0(\log N)^\gamma$ and $\gamma>1$, for
every fixed $D>0$ and sufficiently large $N$,
\[
 C m_N^2e^{-c d/\eta_{\mathrm e}}
 \le C\exp\bigl(\tfrac23\log N-cc_0(\log N)^\gamma\bigr)
 \le N^{-D}.
\]
This proves the last estimate of \eqref{W:eq:localization-theta-integrals}. For the mass of $\vartheta_N$, the exterior pieces have total length at most $Cd$ and
$|\vartheta_N|\le C_E$ there. On the straight interior pieces, parameterize by $v=2-\operatorname{Re}z\ge0$ and use bounded slope and exponential decay. The remaining pieces are at
distance at least $cd$ from the transition. Consequently
\[
 \int_{\mathcal C_N}|\vartheta_N(z)|\,|\mathrm dz|
 \le C_Ed+C_E\int_0^{Cd}e^{-cv/\eta_{\mathrm e}}\,\mathrm dv
                  +C_Ee^{-cd/\eta_{\mathrm e}}
 \le C_E(d+\eta_{\mathrm e})\le C_Ed.
\]

On the ellipse $z=\zeta+\zeta^{-1}$, $\zeta=e^{s_*+\ii\theta}$, one has
\[
 |\dd z|=|\zeta-\zeta^{-1}|\,\dd\theta,\qquad
 |z^2-4|=|\zeta-\zeta^{-1}|^2
       =4(\sinh^2s_*+\sin^2\theta).
\]
Symmetry and $\sin\theta\ge2\theta/\pi$ on $[0,\pi/2]$ imply
\[
 \int_{\text{ellipse}}\frac{|\dd z|}{|z^2-4|}
 \le C\int_0^{\pi/2}\frac{\dd\theta}{(s_*^2+\theta^2)^{1/2}}
 =C\operatorname{arsinh}(\pi/(2s_*))
 \le C(1+|\log s_*|).
\] On each replaced edge piece it is bounded by
$C\int_0^{Cd}({\eta_{\mathrm e}}+u)^{-1}\,\dd u$. The exponential prefactor in \eqref{W:eq:localization-polynomial-bound} is $N^{o(1)}$ because
$\rho s_*=O((\log N)^{\gamma/2})=o(\log N)$. This proves the first estimate. For $X=\mathsf X(H,w)$, the inequality $0\le\Theta(H)\le I$ gives
$0\le M=F_r^*\Theta(H)F_r\le F_r^*F_r=G_r$. Thus $0\le AMA\le AG_rA=A$ and $\|AMA\|\le\kappa^{-1}$ whenever $\kappa I\le G_r\le\kappa^{-1}I$. For any matrix $D$ occurring in the
derivative formulas,
\[
 |\boldsymbol p_j(z)^{\mathsf T}D\boldsymbol p_j(z)|
 \le\|D\|\|\boldsymbol p_j(z)\|_2^2
 \le j\|D\|\mathcal P_\rho(z)^2.
\]
The determinant derivative, $A$, and $AMA$ have operator norm bounded by $C_\kappa$ on this support. Since $j\le\rho\le C m_N$ and $|\widetilde\psi'|\le1$, substitution in
\eqref{W:eq:localization-P0}--\eqref{W:eq:localization-P1} gives $|P_0|+|P_1|\le C_\kappa m_N\mathcal P_\rho^2$ and $|a|\le1$.
\end{proof}

\subsubsection{Resolvent estimates along the interpolation}
\label{W:sec:localization-interpolation-resolvents}

The signed sum over off-diagonal entries is governed by a deterministic matrix built from their third moments. Identifying that matrix first specifies the deterministic vectors for
the resolvent bounds along the interpolation.

Define
\begin{equation}
 \gamma_3=\int\zeta|\zeta|^2\,\Delta\nu(\dd\zeta),\qquad
 B_3=\sum_{a<b}(\gamma_3E_{ab}+\bar\gamma_3E_{ba}).
 \label{W:eq:localization-third-moment-matrix}
\end{equation}
Then $T_{ab}(\zeta)^3=|\zeta|^2T_{ab}(\zeta)$, $\Tr T_{ab}(\zeta)^3=0$, and $\Tr B_3=0$. In the real coordinate convention $\zeta=(x_1+\ii x_2)/\sqrt2$,
\[
 \gamma_3=\frac1{2\sqrt2}
 \big(\Delta M_{111}+\Delta M_{122}
           +\ii(\Delta M_{112}+\Delta M_{222})\big).
\]
The leading sums contain this combination of third moments. The remaining third tensor components enter the error products and are bounded by the absolute third moment. Homogeneity
of the off-diagonal law is used in \eqref{W:eq:localization-third-moment-matrix}; it is part of Assumption~\ref{W:ass:wigner-lee-yin}.

 Fix a deterministic orthonormal eigenbasis of $B_3$ and add $w$. Denote the resulting family of at most $N+1$ vectors by $\mathcal V_N$. For a deterministic $Q_N\ge1$, let
 $\mathcal E$ be the event on which the following inequalities hold for every $z\in\mathcal C_N$:
\begin{align}
 \max_{i,j}|R_H(z)_{ij}-g(z)\delta_{ij}|&\le Q_N/{m_N},
       \label{W:eq:localization-entry-law}\\
 |(R_H(z)w)_i|+|(R_H(\bar z)w)_i|
          &\le Q_N\omega_i,
       &\omega_i&=|w_i|+{m_N}^{-1},\label{W:eq:localization-starting-vector-law}\\
 \|R_H(z)e_i\|^2+\|R_H(\bar z)e_i\|^2&\le Q_N{m_N},
       &v^*|R_H(z)|v&\le Q_N\quad(v\in\mathcal V_N),
       \label{W:eq:localization-energy-law}\\
 \max_{a<b}|H_{ab}|&\le Q_NN^{-1/2}.
       \label{W:eq:localization-entry-size}
\end{align}
Here $R_H(z)=(zI-H)^{-1}$ and $g(z)=\int(z-x)^{-1}\dd\mu_{\rm sc}(x)$. Include in $\mathcal E$ the condition
$\sigma(H(t))\subset[-2-d/8,2+d/8]$.

\begin{lemma}
\label{W:lem:localization-interpolation-event}
The sequence $Q_N$ can be chosen with $Q_N=N^{o(1)}$ so that, for each fixed $D>0$,
\[
 \sup_{t\in[0,1]}\sup_{\|w\|=1}\mathbb P(\mathcal E^c)\le N^{-D}
\]
for all sufficiently large $N$. On $\mathcal E$, every matrix obtained by partially deleting one off-diagonal entry has spectrum in $[-2-d/4,2+d/4]$ and satisfies the coordinate
and column estimates above with an enlarged deterministic
$Q_N$.
\end{lemma}
\begin{proof}
The variance normalization in Lemma~\ref{W:lem:cumulant-event} is $a_t^2=1+(1-t)/N$ in the real case and $a_t=1$ in the complex case. The normalized matrices have variance row sums
one and uniformly bounded standardized moments. The isotropic local law used there yields
\begin{equation}
 |u^*(R_H(z)-g_t(z)I)v|
 \prec \sqrt{\frac{|\Im g_t(z)|}{N|\Im z|}}
                  +\frac1{N|\Im z|},
 \qquad g_t(z)=a_t^{-1}g(z/a_t).
 \label{W:eq:localization-isotropic-input}
\end{equation}
All probability bounds below are uniform in fixed deterministic $t,w$. On the contour outside the two arcs meeting the real axis, differentiation in $a_t$ and the semicircle
formula imply
\[
 |g_t-g|\le N^{-1}m_NN^{o(1)},\qquad
 \frac{|\Im g|}{|\Im z|}\le m_NN^{o(1)},\qquad
 \sqrt{\frac{|\Im g|}{N|\Im z|}}+\frac1{N|\Im z|}
 \le m_N^{-1}N^{o(1)}.
\]
Here the modified edge pieces have $|\Im z|\ge a_0\eta_{\rm e}$; the remaining pieces outside those arcs have $|\Im z|\ge cd$ and $|\Im g|/|\Im z|\le C/\sqrt d$. Applying
\eqref{W:eq:localization-isotropic-input} to $(e_i,e_j)$, $(e_i,w)$ and $(v,v)$, $v\in\mathcal V_N$, and then using Ward's identity gives
\[
 |(R_Hw)_i|\le N^{o(1)}\omega_i,\qquad
 \|R_H(z)v\|^2
 =\frac{|\Im(v^*R_H(z)v)|}{|\Im z|}
 \le m_NN^{o(1)}
 \quad(v\in\mathcal V_N\cup\{e_i\}).
\]
The lower half-plane follows by adjunction.

The spectral deviation theorem of Erd\H{o}s and Xu \cite{erdosxu2023small}, applied to $H(t)$ and $-H(t)$, implies
\[
 \mathbb P\!\left(\sigma(H(t))\not\subset[-2-d/8,2+d/8]\right)
 \le N^{-D}
\]
for each fixed $D$ and all large $N$: choose its deviation parameter $C_D(\log N)^{2/3}$ and use $N^{2/3}d/(\log N)^{2/3}\to\infty$. At a right real tip, $s_*=4\sqrt d$ and
$|\theta|\le s_*/4$ give $r=x-\lambda\ge8d$ and $|y|\le16d$. At $z_h=x+32\ii d$, the local law and Ward's identity give
\[
 \frac{|\Im g(z_h)|}{32d}\le\frac C{\sqrt d},\qquad
 \|R_H(z_h)v\|^2\le m_NN^{o(1)}.
\]
For every point of the corresponding tip arc,
\[
 \frac1{r^2+y^2}\le\frac{17}{r^2+(32d)^2},\qquad
 \|R_H(x+\ii y)v\|^2\le17\|R_H(z_h)v\|^2.
\]
The resolvent identity and $|g'|\le C/\sqrt d$ therefore imply
\[
 |u^*(R_H(x+\ii y)-R_H(z_h))v|
 \le Cdm_NN^{o(1)}=m_N^{-1}N^{o(1)},\qquad
 |g(x+\ii y)-g(z_h)|\le C\sqrt d=m_N^{-1}N^{o(1)}.
\]
Reflection proves the same estimates at the left tip. This establishes the entry bounds, the bounds for $R_Hw$, and the column bounds on the whole contour.

For the spectral measure $\mu_v$ at a required deterministic unit vector, the local law also gives, for $u\ge c\eta_{\rm
e}$,
\[
 \mu_v([x-u,x+u])
 \le2u[-\Im v^*R_H(x+\ii u)v]\le N^{o(1)}u.
\]
Decompose into the intervals $|\lambda-x|\le\eta$ and $2^j\eta<|\lambda-x|\le2^{j+1}\eta$, where $\eta=|\Im z|$. Each contributes at most $2N^{o(1)}$ to
$\int|z-\lambda|^{-1}\mu_v(\dd\lambda)$, and only $O(\log N)$ intervals are needed. At an exterior point of spectral distance $\delta\ge cd$, start instead at $\delta/2$; the
central interval has zero mass. Thus $v^*|R_H(z)|v\le N^{o(1)}$ everywhere required.

For a partially deleted entry, write $V=H-K_{ab}=\iota_I V_I\iota_I^*$ with $I=\{a,b\}$. To obtain a consequence of the stated event $\mathcal E$, put $b_*=2\cosh(4\sqrt d)$ and
$z_\pm=\pm(2+d/4)$. Confinement and scalar spectral comparison give
\[
 0\le R_H(z_+)\le C R_H(b_*),\qquad
 0\le-R_H(z_-)\le C[-R_H(-b_*)].
\]
Indeed $(b_*-\lambda)/(2+d/4-\lambda)\le C$ on the confined spectrum, and reflection gives the bound at the left spectral barrier. The contour entry bounds at $\pm b_*$ therefore
imply $\|R_{H,II}(z_\pm)\|\le C$. Also $\|V_I\|\le N^{-1/2+o(1)}$. Hence, uniformly for $0\le u\le1$,
\[
 \det(z_\pm I-H+uV)
 =\det(z_\pm I-H)\det(I_I+uV_IR_{H,II}(z_\pm))\ne0.
\]
Continuity of the eigenvalues along $H-uV$ proves the claimed confinement. On $\mathcal C_N$,
\[
 R_{H-uV}-R_H
 =-R_{H,\cdot I}(I_I+uV_IR_{H,II})^{-1}uV_IR_{H,I\cdot},
 \qquad\|(I_I+uV_IR_{H,II})^{-1}\|\le2.
\]
Its entries and column norms are respectively bounded by $N^{-1/2+o(1)}$ and $N^{-1/2+o(1)}m_N^{1/2}$, and
\[
 |((R_{H-uV}-R_H)w)_i|
 \le N^{-1/2+o(1)}(\omega_a+\omega_b)
 \le m_N^{-1}N^{o(1)}.
\]
Adding these inequalities to the bounds for $R_H$ treats all deletions.

For simultaneous control, use deterministic $N^{-20}$ nets on the contour, the exterior vertical segments and the nonreal domains above. There are at most $CN^{40}$ net points and
$CN^2$ vector pairs. On the confinement event the relevant spectral distances exceed $N^{-1}$; the spectral theorem gives Lipschitz constants $CN^2,CN^3,CN^2$ for $u^*R(z)v$,
$\|R(z)v\|^2$ and $v^*|R(z)|v$, respectively. The error on passing from the nets is therefore at most $CN^{-17}$. Use failure exponent $D+45$ before the union bound. For the
entries, choose an even $p>(D+3)/\varepsilon$ and use
\[
 \mathbb P\!\left(\max_{a,b}|\sqrt N H_{ab}|>N^\varepsilon\right)
 \le C_pN^{2-p\varepsilon}\le N^{-D}.
\]
Finally choose $N_k\uparrow\infty$ so that, for $N\ge N_k$, all deterministic factors in these finitely many bounds are at most $N^{1/k}$ and the exceptional probability is at most
$N^{-k}$. Taking $Q_N=N^{C/k}$ on $[N_k,N_{k+1})$, with one fixed sufficiently large $C$, absorbs their fixed powers. Then $Q_N=N^{o(1)}$ and, for every fixed $D$, $k\ge D$
whenever $N\ge N_{\lceil D\rceil}$. This proves the lemma.
\end{proof}

The two cubic calculations involve $w^*RB_3Rw$ and $\Tr(B_3R)$. Both are controlled using the deterministic eigenbasis of $B_3$ and the estimates for $|R|$.

\begin{lemma}
\label{W:lem:localization-third-moment-matrix}
The deterministic Hermitian matrix $B_3$ satisfies
\begin{equation}
 \|B_3\|_1\le C|\gamma_3|N\log(2N).
 \label{W:eq:localization-third-moment-trace-norm}
\end{equation}
On $\mathcal E$,
\begin{equation}
 |w^*RB_3Rw|+|\Tr(B_3R)|\le Q_N N.
 \label{W:eq:localization-third-moment-resolvent-bounds}
\end{equation}
\end{lemma}
\begin{proof}
The matrix $B_3$ is the same homogeneous third moment matrix as in \eqref{W:eq:cumulant-third-matrix}. Lemma~\ref{W:lem:cumulant-trace-norm} therefore gives
\eqref{W:eq:localization-third-moment-trace-norm}, including the exact value in the real case $2|\gamma_3|(N-1)$. It remains to estimate its two different contractions with the
resolvent on the present contour.

Let $B_3=\sum_j\lambda_jv_jv_j^*$ in the deterministic basis used to construct $\mathcal E$. For a normal resolvent, polar decomposition and Cauchy--Schwarz give
\[
 |u^*Rv|^2\le(u^*|R|u)(v^*|R|v).
\]
Indeed, write $R=U|R|$. Normality implies that $U$ commutes with $|R|$, so Cauchy--Schwarz applied to $\langle |R|^{1/2}U^*u,|R|^{1/2}v\rangle$ proves the inequality. Applying it
separately to $w^*Rv_j$ and $v_j^*Rw$ gives their product bound, even though those two scalars need not be complex conjugates. It follows that
\[
 |w^*RB_3Rw|
 \le(w^*|R|w)\sum_j|\lambda_j|v_j^*|R|v_j,
 \qquad
 |\Tr(B_3R)|\le\sum_j|\lambda_j|v_j^*|R|v_j.
\]
By \eqref{W:eq:localization-energy-law}, the two quantities have sum at most $C|\gamma_3|(Q_N^2+Q_N)N\log(2N)$. Increasing $Q_N$ to this deterministic prefactor gives
\eqref{W:eq:localization-third-moment-resolvent-bounds}.
\end{proof}

\subsubsection{The cubic term for the unnormalized trace}
\label{W:sec:localization-cubic-trace}

The cubic expansion of the unnormalized trace contains a term linear in the resolvent error. Estimating each matrix element separately would lose the cancellation needed for this
term. We therefore compute the trace expansion before summing over entries.

Fix the common matrix $H=H(t)$ and write $R=R_H(z)$. Define
\begin{align}
 W_{ab}^{w}(z)&=\int w^*RT_{ab}(\zeta)RT_{ab}(\zeta)
                           RT_{ab}(\zeta)Rw\,\Delta\nu(\dd\zeta),
       \label{W:eq:localization-cubic-matrix-element}\\
 W_{ab}^{\rm tr}(z)&=\int\Tr(RT_{ab}(\zeta)RT_{ab}(\zeta)
                           RT_{ab}(\zeta)R)\,\Delta\nu(\dd\zeta).
       \label{W:eq:localization-cubic-trace}
\end{align}
The trace product is estimated separately. For example, with
\[
 T=\begin{pmatrix}0&1\\1&0\end{pmatrix},\qquad
 Q=gI_2+\varepsilon T,
\]
one has $TQ=gT+\varepsilon I_2$ and
\begin{equation}
 \operatorname{tr}_2((TQ)^3)=6g^2\varepsilon+2\varepsilon^3.
 \label{W:eq:trace-cubic-example}
\end{equation}
The term linear in the resolvent error therefore remains. At the bound $\|E_{ab}\|\le N^{o(1)}m_N^{-1}$, summing its absolute estimates would give the scale
$N^{-3/2}N^2m_N^{-1}=N^{1/6}$. After integration by parts, we sum this term using $B_3$ and bound the remaining terms by at least two factors of $E_{ab}$. The latter scale is
$N^{-3/2}N^2m_N^{-2}=N^{-1/6}$. This is the extra cancellation required for the localization trace, whose normalization counts each eigenstate with its full weight.

Let $a=a(H(t))$ be a scalar, independent of the contour variable, with $|a|\le1$, and define
\begin{equation}
 \mathcal J_N^{\rm tr}(H,a)
  =\frac{aN^{-3/2}}{2\pi\ii}\int_{\mathcal C_N}\vartheta_N(z)
       \sum_{a'<b'}W_{a'b'}^{\rm tr}(z)\,\dd z.
 \label{W:eq:localization-trace-contribution}
\end{equation}
\begin{lemma}
\label{W:lem:localization-trace-cubic}
On $\mathcal E$, the quantity just defined satisfies
\begin{equation}
 |\mathcal J_N^{\rm tr}(H,a)|
       \le Q_N\bigl(N^{-1/2}+N^{-1/6}\bigr).
 \label{W:eq:localization-trace-cubic-estimate}
\end{equation}
\end{lemma}
\begin{proof}
For fixed $T=T_{ab}(\zeta)$, differentiate all three resolvents in the cyclic product, with $\partial_zR=-R^2$:
\begin{align*}
 \partial_z\Tr(TRTRTR)
 ={}&-\Tr(TR^2TRTR)-\Tr(TRTR^2TR)-\Tr(TRTRTR^2)\\
 ={}&-3\Tr(RTRTRTR).
\end{align*}
Cyclically rotating the first, second, and third products by two, four, and one factors, respectively, gives $\Tr(RTRTRTR)$ in each case. Consequently,
\begin{equation}
 \Tr(RTRTRTR)=-\tfrac13\partial_z\Tr((TR)^3).
 \label{W:eq:localization-trace-differentiation}
\end{equation}
Integrating by parts on the closed piecewise smooth contour,
\begin{equation}
 \int_{\mathcal C_N}\vartheta_N(z)\Tr(RTRTRTR)\,\dd z
   =\tfrac13\int_{\mathcal C_N}\vartheta_N'(z)
                                  \Tr((TR)^3)\,\dd z.
 \label{W:eq:localization-trace-ibp}
\end{equation}
The boundary terms cancel at successive contour vertices. Analyticity near the contour and the finite third absolute moment justify integration against the entry law and the finite
pair sum.

The trace now involves only $Q=R_{II}=gI_2+E_{ab}$: $\Tr((TR)^3)=\tr_2((TQ)^3)$. Expanding this product and collecting terms under the trace gives
\begin{equation}
 \tr_2((TQ)^3)=g^3\tr_2(T^3)+3g^2\tr_2(T^3E_{ab})+r_{ab}(\zeta,z),
 \label{W:eq:localization-trace-expand}
\end{equation}
where the remainder is exactly
\begin{equation}
 r_{ab}(\zeta,z)=3g\tr_2(T^2E_{ab}TE_{ab})+\tr_2((TE_{ab})^3).
 \label{W:eq:localization-trace-quadratic-remainder}
\end{equation}
Since $T$ is off diagonal on two coordinates, $\tr_2(T^3)=0$. The coefficient of the term linear in $E_{ab}$ is the signed third moment.
Figure~\ref{W:fig:moment-trace-cubic-errors} shows the three cyclic placements producing each coefficient $3$. \begingroup \setlength{\intextsep}{1.2\baselineskip}
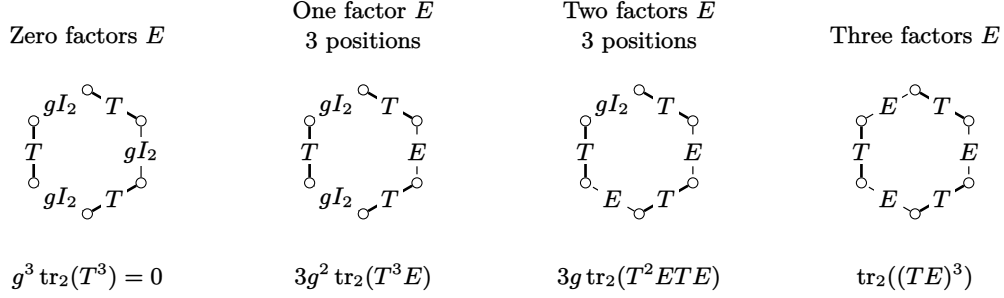
\begin{figure}[H]
\centering
\begin{tikzpicture}[x=1cm,y=1cm,every node/.style={mg text}]
 \foreach \shift/\tag in {0/a,3.65/b,7.3/c,10.95/d}{
  \begin{scope}[xshift=\shift cm]
   \foreach \k/\ang in {1/90,2/30,3/-30,4/-90,5/-150,6/150}{
    \node[mg site] (\tag\k) at ({1.25+.82*cos(\ang)},{.82*sin(\ang)}) {};
   }
   \foreach \u/\v in {1/2,3/4,5/6}{\draw[mg insertion] (\tag\u)--node[mg label] {$T$}(\tag\v);}
  \end{scope}
 }
 \foreach \u/\v in {a2/a3,a4/a5,a6/a1,b4/b5,b6/b1,c6/c1}{\draw[mg kernel] (\u)--node[mg label] {$gI_2$}(\v);}
 \foreach \u/\v in {b2/b3,c2/c3,c4/c5,d2/d3,d4/d5,d6/d1}{\draw[mg incidence] (\u)--node[mg label] {$E$}(\v);}
 \node at (1.25,1.55) {Zero factors $E$};
 \node[align=center] at (4.9,1.65) {One factor $E$\\$3$ positions};
 \node[align=center] at (8.55,1.65) {Two factors $E$\\$3$ positions};
 \node at (12.2,1.55) {Three factors $E$};
 \node at (1.25,-1.65) {$g^3\tr_2(T^3)=0$};
 \node at (4.9,-1.65) {$3g^2\tr_2(T^3E)$};
 \node at (8.55,-1.65) {$3g\tr_2(T^2ETE)$};
 \node at (12.2,-1.65) {$\tr_2((TE)^3)$};
\end{tikzpicture}
\caption{Expansion of $\tr_2(TQTQTQ)$ with $Q=gI_2+E$ in \eqref{W:eq:localization-trace-expand} and \eqref{W:eq:localization-trace-quadratic-remainder}. Vertices denote the two coordinates of the selected entry. Thick edges carry $T$; the other edges carry the indicated term from $Q$, with dashed edges denoting $E$. The factors are read clockwise. The four groups contain zero, one, two, and three factors $E$, with multiplicities $1,3,3,1$. Their contributions are the four expressions beneath the diagrams. Since $\tr_2(T^3)=0$, the first vanishes. The group with one $E$ is summed using $B_3$; the last two groups are bounded by products of at least two resolvent errors.}
\label{W:fig:moment-trace-cubic-errors}
\end{figure}
\endgroup

 Using $\|T_{ab}(\zeta)\|=|\zeta|$ and $|\tr_2 A|\le2\|A\|$, we obtain
\[
 |r_{ab}(\zeta,z)|
 \le C|\zeta|^3\big(|g|\|E_{ab}\|^2+\|E_{ab}\|^3\big)
 \le Q_N|\zeta|^3{m_N}^{-2}.
\]
After integration against the signed entry law and summation,
\begin{equation}
 \sum_{a<b}\int\Tr((T_{ab}(\zeta)R)^3)\,\Delta\nu(\dd\zeta)
     =3g^2\Tr(B_3R)+\mathcal R_{\rm tr}(z),
 \qquad |\mathcal R_{\rm tr}(z)|\le Q_NN^2{m_N}^{-2}.
 \label{W:eq:localization-trace-signed}
\end{equation}
For complex $\zeta$, the linear term can also be checked in coordinates:
\[
 \tr_2(T^3E_{ab})
 =|\zeta|^2\bigl(\zeta(E_{ab})_{ba}
                      +\bar\zeta(E_{ab})_{ab}\bigr).
\]
Integration gives $\gamma_3(E_{ab})_{ba} +\bar\gamma_3(E_{ab})_{ab}$. The entries of $E_{ab}$ are not conjugated in this identity, since a nonreal resolvent compression need not be
Hermitian. Summing these terms gives $\Tr(B_3(R-gI))=\Tr(B_3R)$ because $\Tr B_3=0$. All other components of the third tensor remain in the displayed quadratic and cubic remainder.
We estimate the resulting signed matrix sum before taking absolute values.

For fixed $H(t)$, $a(H(t))$ is constant in $z$. Thus the contour integration by parts differentiates only $\vartheta_N$, whose derivative has bounded $L^1$ norm by
\eqref{W:eq:localization-theta-integrals}. Apply \eqref{W:eq:localization-trace-ibp}--\eqref{W:eq:localization-trace-signed}, $|a|\le1$, \eqref{W:eq:localization-theta-integrals},
and \eqref{W:eq:localization-third-moment-resolvent-bounds}. The full trace contribution is at most
\begin{equation}
 \begin{split}
 C N^{-3/2}\int_{\mathcal C_N}|\vartheta_N'(z)|
       \big(Q_NN+Q_NN^2{m_N}^{-2}\big)|\dd z|
  &\le Q_N\big(N^{-1/2}+N^{1/2}{m_N}^{-2}\big)\\
  &=Q_N\big(N^{-1/2}+N^{-1/6}\big).
 \end{split}
 \label{W:eq:localization-trace-final}
\end{equation}
This proves \eqref{W:eq:localization-trace-cubic-estimate}, including the factor $1/(2\pi\ii)$ in the constant. For this unnormalized trace, the two powers come from the signed
linear term and the quadratic error, respectively.
\end{proof}

In \eqref{W:eq:localization-common-response}, the trace multiplier is $a(H)=\Gamma_\rho\widetilde\psi'(C)$, which is constant in $z$ and has absolute value at most one at the
moment vector of the matrix. The lemma therefore bounds that entire trace term. To use it in the entry comparison, we must still evaluate the derivative coefficient and all moment
increments at one common matrix. Section~\ref{W:sec:localization-moment-increments-and-kernels} makes this change of base. We then estimate the matrix element terms and the
nonlinear chain rule terms separately.

\subsubsection{Changes in moments under entry replacement}
\label{W:sec:localization-moment-increments-and-kernels}

Replacing the matrix with one entry set to zero by the matrix used in every summand changes the moment vector at which the derivative coefficient is evaluated. The weights below
keep track of that change and of the number of entry pairs, so that the estimate after summation has the required order.

On $\mathcal E$ put
\[
 u_{ab}=\omega_a\omega_b,\qquad v_{ab}=\omega_a^2+\omega_b^2.
\]
Write $S_j=\sum_i\omega_i^j$. Since $\|w\|_2=1$,
\[
 S_1\le\sqrt N+Nm_N^{-1}\le Cm_N^2,\qquad
 S_2\le2+2Nm_N^{-2}\le Cm_N,
 \qquad S_3+S_4\le C.
\]
For example, the last bound uses $\sum_i|w_i|^j\le1$ for $j\ge2$ and $Nm_N^{-3}=1$. The pair sums are
\begin{align*}
 \sum_{a<b}v_{ab}&=(N-1)S_2,\\
 \sum_{a<b}u_{ab}v_{ab}&=S_1S_3-S_4,\\
 \sum_{a<b}v_{ab}^2&=(N-2)S_4+S_2^2,\\
 \sum_{a<b}u_{ab}^3&=\tfrac12(S_3^2-S_6).
\end{align*}
Substitution gives
\begin{equation}
 \sum_{a<b}v_{ab}\le C{m_N}^4,\quad
 \sum_{a<b}u_{ab}v_{ab}\le C{m_N}^2,\quad
 \sum_{a<b}v_{ab}^2\le C{m_N}^3,\quad
 \sum_{a<b}u_{ab}^3\le C.
 \label{W:eq:localization-integrated-weight-sums}
\end{equation}
Let $K=K_{ab}$ and $\tau_{ab}=\|H-K\|\le Q_NN^{-1/2}$.
\begin{lemma}
\label{W:lem:localization-moment-increments}
On $\mathcal E$, the changes of $\mathsf X$ and its cubic coefficient satisfy
\begin{align}
 \mathfrak{n}(\mathsf{X}(H)-\mathsf{X}(K))
      &\le Q_N{m_N}(\tau_{ab}u_{ab}+\tau_{ab}^2v_{ab}),
       \label{W:eq:localization-moment-difference}\\
 \mathfrak{n}(c_{ab}(K;\zeta))
      &\le Q_N{m_N}N^{-3/2}|\zeta|^3v_{ab},
       \label{W:eq:localization-cubic-moment-bound}\\
 \mathfrak{n}(c_{ab}(H;\zeta)-c_{ab}(K;\zeta))
      &\le Q_N{m_N}N^{-3/2}\tau_{ab}|\zeta|^3v_{ab}.
       \label{W:eq:localization-cubic-difference}
\end{align}
Here $c_{ab}(H;\zeta)$ is the coefficient from \eqref{W:eq:localization-moment-coefficients} at base $H$.
\end{lemma}
\begin{proof}
Write $I=\{a,b\}$, $\iota=\iota_I$, $V=H-K=\iota B\iota^*$, $\|B\|=\tau_{ab}$. Put $T=T_{ab}(\zeta)_{II}$ and $U=\iota T\iota^*$. The auxiliary matrix $U$ is unnormalized; its norm
is $|\zeta|$. For a matrix $B_0$ on the regular event, Cauchy's formulas are
\begin{align*}
 G(B_0)&=\frac1{2\pi\mathrm i}\int_{\mathcal C_N}
       \boldsymbol p_\rho(z)\boldsymbol p_\rho(z)^{\mathsf T}
                         w^*R_{B_0}(z)w\,\mathrm dz,\\
 M(B_0)&=\frac1{2\pi\mathrm i}\int_{\mathcal C_N}
       \vartheta_N(z)\boldsymbol p_r(z)\boldsymbol p_r(z)^{\mathsf T}
                         w^*R_{B_0}(z)w\,\mathrm dz,\\
 \mathsf t(B_0)&=\frac1{2\pi\mathrm i}\int_{\mathcal C_N}
                         \vartheta_N(z)\operatorname{Tr}R_{B_0}(z)\,\mathrm dz.
\end{align*}
The two rank-one matrix weights have trace norms $\|\boldsymbol p_j\boldsymbol p_j^{\mathsf T}\|_1 =\|\boldsymbol p_j\|_2^2\le j\mathcal P_\rho^2\le Cm_N\mathcal P_\rho^2$. For the
weight in the formula for $M(B_0)$, multiply by the additional scalar $|\vartheta_N(z)|\le C_E$. The exact resolvent identity, expanded once more, is
\[
 R_H-R_K=R_KVR_K+R_KVR_KVR_H.
\]
Because $V$ is off-diagonal on $I$, its first contraction with the starting vector is bounded by $Q_N\tau_{ab}\omega_a\omega_b$. In the second, the interior compression is bounded
and each of its two exterior vectors has restricted norm at most $Q_N(\omega_a^2+\omega_b^2)^{1/2}$. Hence the second contraction is bounded by $Q_N\tau_{ab}^2v_{ab}$. Integration
gives
\[
 \|G(H)-G(K)\|_1+\|M(H)-M(K)\|_1
 \le Q_Nm_N(\tau_{ab}u_{ab}+\tau_{ab}^2v_{ab}).
\]
The cubic integrand is $N^{-3/2}w^*R_{B_0}UR_{B_0}UR_{B_0}UR_{B_0}w$. The two interior compressions are bounded, so its modulus is at most $Q_NN^{-3/2}|\zeta|^3v_{ab}$. Thus the
components corresponding to $G(B_0)$ and $M(B_0)$ satisfy \eqref{W:eq:localization-cubic-moment-bound} after integration.

For the difference of these cubic integrands put $R=R_H$ and $S=R_K$. The four term telescoping identity, with $R-S=RVS$, is
\begin{align*}
 RURURUR-SUSUSUS
 ={}&RVSUSUSUS+RURVSUSUS\\
    &+RURURVSUS+RURURURVS.
\end{align*}
In each term the left exterior vector is $w^*R\iota$ and the right one is $\iota^*Sw$. Every other resolvent lies between two matrices supported on $I$ and therefore appears as
$\iota^*R\iota$ or $\iota^*S\iota$. Their norms are bounded. In order, the four matrix factors are $(B,T,T,T)$, $(T,B,T,T)$, $(T,T,B,T)$, and $(T,T,T,B)$. Each term tested against
the starting vector consequently has modulus at most $Q_N\tau_{ab}|\zeta|^3v_{ab}$. For each of the four terms, the trace norm of its integral defining the Gram matrix is at most
\[
 CN^{-3/2}\int_{\mathcal C_N}
      m_N\mathcal P_\rho(z)^2
             Q_N\tau_{ab}|\zeta|^3v_{ab}\,|\mathrm dz|
 \le Q_Nm_NN^{-3/2}\tau_{ab}|\zeta|^3v_{ab}.
\]
For its mixed integral the same expression has the additional factor $|\vartheta_N(z)|\le C_E$, and hence the same bound. This verifies \eqref{W:eq:localization-cubic-difference}
for both matrix components.

For the unnormalized trace, cyclicity gives
\[
 \operatorname{Tr}(R_H-R_K)=\operatorname{tr}_2(B\iota^*R_KR_H\iota),
 \qquad
 \operatorname{Tr}(SUSUSUS)
 =\operatorname{tr}_2(TC_STC_STD_S),
\]
where $C_S=\iota^*S\iota$ and $D_S=\iota^*S^2\iota$. For $B_1,B_2\in\{H,K\}$,
\[
 |(\iota^*R_{B_1}R_{B_2}\iota)_{ij}|
 \le\|R_{B_1}(\bar z)e_i\|\|R_{B_2}(z)e_j\|
 \le Q_Nm_N.
\]
Thus $|\mathsf t(H)-\mathsf t(K)|\le Q_N\tau_{ab}m_Nd$, and the cubic trace coefficient is at most $Q_NN^{-3/2}|\zeta|^3m_Nd$. For each telescoped trace term the closing factor is
$\iota^*S R\iota$ and the other three compressed resolvents are bounded; its integral is at most $Q_NN^{-3/2}\tau_{ab}|\zeta|^3m_Nd$. These bounds use exactly
$\int_{\mathcal C_N}|\vartheta_N(z)|\,|\mathrm dz|\le C_Ed$. The same closing factor gives the linear and quadratic trace coefficients, including their Taylor normalizations:
\[
 |\ell_{ab}^{\mathsf t}(K;\zeta)|
       \le Q_N N^{-1/2}|\zeta|m_Nd,\qquad
 |q_{ab}^{\mathsf t}(K;\zeta)|
       \le Q_N N^{-1}|\zeta|^2m_Nd.
\]
For $G$ and $M$, the linear term for the off-diagonal perturbation is bounded by $Q_N|\zeta|u_{ab}$, and the quadratic contraction by $Q_N|\zeta|^2v_{ab}$. Their matrix weights
have integrated trace norm at most $Q_Nm_N$. Finally $m_Nd=m_N^{-1}N^{o(1)}$, $u_{ab}\ge m_N^{-2}$ and $v_{ab}\ge2m_N^{-2}$. Each trace bound is therefore included in its
corresponding right side of \eqref{W:eq:localization-moment-difference}--\eqref{W:eq:localization-cubic-difference}. Adding the three component bounds in the definition of
$\mathfrak n$ proves all three inequalities. The preceding linear and quadratic calculations also give the full $\ell_{ab}$ and $q_{ab}$ bounds used in the nonlinear cubic terms
below, because $m_Nu_{ab}\ge m_N^{-1}$ and $m_Nv_{ab}\ge2m_N^{-1}$. In particular, both norms include the scalar trace component.
\end{proof}

Let
\[
 L_{ab}(K)=\int D_{\mathsf X}\mathcal{F}(\mathsf{X}(K))[c_{ab}(K;\zeta)]\,\Delta\nu(\dd\zeta),
\]
and define $L_{ab}^{\rm com}(H)= \int D_{\mathsf X}\mathcal F(\mathsf X(H))[c_{ab}(H;\zeta)]\, \Delta\nu(\dd\zeta)$.
\begin{lemma}
\label{W:lem:localization-change-of-base}
On $\mathcal E$,
\[
 \sum_{a<b}|L_{ab}(K_{ab})-L_{ab}^{\rm com}(H)|
            \le N^{-1/3+o(1)}.
\]
\end{lemma}
\begin{proof}
For each pair the exact difference of the integrands is
\begin{equation}
 \begin{split}
 &D_{\mathsf X}\mathcal F(\mathsf X(K))[c_{ab}(K;\zeta)]
      -D_{\mathsf X}\mathcal F(\mathsf X(H))[c_{ab}(H;\zeta)]\\
 &\quad=
 \bigl(D_{\mathsf X}\mathcal F(\mathsf X(K))-D_{\mathsf X}\mathcal F(\mathsf X(H))\bigr)
                    [c_{ab}(K;\zeta)]\\
 &\qquad\quad+
 D_{\mathsf X}\mathcal F(\mathsf X(H))
                    [c_{ab}(K;\zeta)-c_{ab}(H;\zeta)].
 \end{split}
 \label{W:eq:localization-common-base-difference}
\end{equation}
Apply Corollary~\ref{W:cor:localization-gradient-bound} to $X=\mathsf X(K)$ and $Y=\mathsf X(H)$ in the first term on the right. The second line uses the first derivative bound at
$\mathsf X(H)$. The integration variable $\zeta$ has the fixed signed law $\Delta\nu$, whereas $\tau_{ab}$ is determined by $H_{ab}$. The total variation integral of $|\zeta|^3$ is
uniformly bounded. Therefore Lemma~\ref{W:lem:localization-derivatives} and \eqref{W:eq:localization-moment-difference}--\eqref{W:eq:localization-cubic-difference} give
\begin{equation}
 \begin{split}
 \sum_{a<b}|L_{ab}(K_{ab})-L_{ab}^{\rm com}(H)|
 \le Q_N\bigg(&{m_N}N^{-2}\sum v_{ab}
       +{m_N}^2N^{-2}\sum u_{ab}v_{ab}\\
       &+{m_N}^2N^{-5/2}\sum v_{ab}^2\bigg)
 \le Q_N({m_N}^{-1}+{m_N}^{-2}+{m_N}^{-5/2}).
 \end{split}
 \label{W:eq:localization-change-of-base-bound}
\end{equation}
The first term bounds the change in the cubic moment increment; the other two bound the change in the factor multiplying that increment. The constant in $Q_N$ includes
$\sup_N\int|\zeta|^3(\nu_N^{\rm b}+\nu^G)(\mathrm d\zeta)<\infty$. Indeed, substituting the pair sums in \eqref{W:eq:localization-integrated-weight-sums}, the three powers are
\[
 m_NN^{-2}m_N^4=m_N^{-1},\quad
 m_N^2N^{-2}m_N^2=m_N^{-2},\quad
 m_N^2N^{-5/2}m_N^3=m_N^{-5/2}.
\]
The first is the largest. This proves the lemma.
\end{proof}

\subsubsection{The cubic terms for matrix elements}
\label{W:sec:localization-cubic-matrix-elements}

The moment changes are now evaluated at the common matrix, where their cubic coefficients can be summed with their signs. Write $H=H(t)$ and $R=R_H(z)$, and retain the integrands
in \eqref{W:eq:localization-cubic-matrix-element}--\eqref{W:eq:localization-cubic-trace}. Their contraction with the moment derivatives is the remaining step for the matrix-element
terms.

The third resolvent derivative divided by $3!$ is the cubic resolvent product. For $k=\rho,r$ set $Q_k(z)=\boldsymbol p_k(z)\boldsymbol p_k(z)^{\mathsf T}$. The Cauchy formulas in
the proof of Lemma~\ref{W:lem:localization-moment-increments} give
\begin{align*}
 \sum_{a<b}\int c_{ab}^{G}(H;\zeta)\,\Delta\nu(\dd\zeta)
   &=\frac{N^{-3/2}}{2\pi\ii}\int Q_\rho(z)
                                   \sum_{a<b}W_{ab}^{w}(z)\,\dd z,\\
 \sum_{a<b}\int c_{ab}^{M}(H;\zeta)\,\Delta\nu(\dd\zeta)
   &=\frac{N^{-3/2}}{2\pi\ii}\int \vartheta_N(z)Q_r(z)
                                   \sum_{a<b}W_{ab}^{w}(z)\,\dd z,\\
 \sum_{a<b}\int c_{ab}^{\mathsf t}(H;\zeta)\,\Delta\nu(\dd\zeta)
   &=\frac{N^{-3/2}}{2\pi\ii}\int \vartheta_N(z)
                                   \sum_{a<b}W_{ab}^{\rm tr}(z)\,\dd z.
\end{align*}
Each integral is over $\mathcal C_N$. The finite sums and the integral over $\zeta$ may be interchanged: for fixed $H$ their integrands are bounded by a constant times $|\zeta|^3$,
which is integrable against $\nu_N^{\rm b}+\nu^G$. Contracting with \eqref{W:eq:localization-response}, and using $\Tr(DQ_k)=\boldsymbol p_k^{\mathsf T}D\boldsymbol p_k$, gives
\begin{equation}
 \begin{split}
 \sum_{a<b}L_{ab}^{\rm com}(H)
  =\frac{N^{-3/2}}{2\pi\ii}\int_{\mathcal C_N}
  \bigg[&(P_0(H,z)+\vartheta_N(z)P_1(H,z))\sum_{a<b}W_{ab}^{w}(z)\\
        &+a(H)\vartheta_N(z)\sum_{a<b}W_{ab}^{\rm tr}(z)\bigg]\dd z.
 \end{split}
 \label{W:eq:localization-common-response}
\end{equation}
Every random derivative coefficient is common across the sum over entries and remains correlated with $R$.

\Needspace{6\baselineskip}

\manualstatementlabel{W:sec:third-order-localization-comparison}{\thesubsubsection}

Use $H(t)$ and the polynomials $P_0$ and $P_1$ in Section~\ref{W:sec:localization-moment-increments-and-kernels}. For $I=\{a,b\}$ write $R_{II}=gI_2+E_{ab}$, where
$\|E_{ab}\|\le Q_N/{m_N}$ on $\mathcal E$. Expanding the two interior resolvents in \eqref{W:eq:localization-cubic-matrix-element} gives
\begin{equation}
 \sum_{a<b}W_{ab}^{w}(z)
       =g(z)^2w^*RB_3Rw+\mathcal R_w(z),
 \qquad
 |\mathcal R_w(z)|\le \frac{Q_N}{{m_N}}\sum_{a<b}v_{ab}
                     \le Q_N N.
 \label{W:eq:localization-matrix-element-sum}
\end{equation}
On $I=\{a,b\}$ put $T=T_{ab}(\zeta)_{II}$ and $R_{II}=gI+E_{ab}$. Equation~\eqref{W:eq:cumulant-interior-identity} gives
\[
 \|TR_{II}TR_{II}T-g^2T^3\|
 \le|\zeta|^3(2|g|\|E_{ab}\|+\|E_{ab}\|^2).
\]
Moreover,
\[
 \|w^*R_{\cdot I}\|_2\|R_{I\cdot}w\|_2
 \le Q_N(\omega_a^2+\omega_b^2).
\]
Since $\|E_{ab}\|\le Q_N/m_N$ and $|g|\le C$, integration against the positive measure $\nu_N^{\rm b}+\nu^G$ bounds one remainder by $Q_Nm_N^{-1}v_{ab}$. Summing the leading term
first gives $g^2w^*RB_3Rw$; summing the remainder bounds gives $Q_Nm_N^{-1}m_N^4=Q_NN$. This proves \eqref{W:eq:localization-matrix-element-sum}.

For complex entries, \eqref{W:eq:resolvent-complete-complex-cubic} applies to the complex matrix $Q=R_{II}$ and contains $\gamma_3=\int\zeta^2\bar\zeta\,\Delta\nu$ and
$\kappa_3=\int\zeta^3\,\Delta\nu$. The terms containing $\kappa_3$ have two off-diagonal entries of $E_{ab}$; the other nonleading terms have at least one error factor. The same
norm estimate therefore bounds these contributions.

By \eqref{W:eq:localization-derivative-kernel-bound}, \eqref{W:eq:localization-contour-mass}, and Lemma~\ref{W:lem:localization-third-moment-matrix}, the entire term involving the
starting vector of \eqref{W:eq:localization-common-response} is bounded by
\begin{equation}
 Q_N N^{-3/2}{m_N}N
       \int_{\mathcal C_N}\mathcal P_\rho(z)^2|\dd z|
          \le N^{-1/6+o(1)}.
 \label{W:eq:localization-matrix-element-bound}
\end{equation}
Here we used $\sup_{z\in\mathcal C_N}|\vartheta_N(z)|\le C_E$ and the local law for the deterministic vectors selected above.
\subsubsection{Expectation of the cubic terms}
\label{W:sec:localization-exceptional-events-and-comparison}

For a fixed sampled matrix $H$, expand any one of the cubic integrands in real coordinates as $\sum_{i,j,k} b_{ijk}(H)x_ix_jx_k$. Since its coefficients are fixed while the dummy
entry is integrated,
\[
 \int\sum_{i,j,k}b_{ijk}(H)x_ix_jx_k\,\Delta\mu(\dd x)
   =\sum_{i,j,k}b_{ijk}(H)\Delta M_{ijk}.
\]
Multiplication by $\mathbf1_{\mathcal E}(H)$ preserves this identity. The marginal distribution of each $K_{ab}=H(t)^{[ab\leftarrow0]}$ is the product distribution over the
remaining entries. Summing their expectations therefore gives $\mathbb E Z_N^{(1)}$ and $\mathbb E Z_N^{(3)}$ in the entry interpolation formula.

\begin{proof}[Proof of Theorem~\ref{W:thm:localization-cubic}]
We first control the complement of $\mathcal E$. Apply Lemma~\ref{W:lem:smooth-functional-calculus} to
\[
 f_{jk}(x)=p_j(h_d(x))p_k(h_d(x)),\qquad
 g_{jk}(x)=f_{jk}(x)\vartheta_N(x),\qquad \vartheta_N(x).
\]
Their derivatives through order three are polynomially bounded: Proposition~\ref{W:prop:gram-bounds} controls the polynomial factors,
$\|\vartheta_N^{(j)}\|_\infty\le C_j{m_N}^{2j}$, and the scalar product rule controls $g_{jk}$. If $M_j$ bounds these scalar derivatives, then
\[
 \begin{split}
 \mathfrak n(D_H^j\mathsf X(H)[V_1,\ldots,V_j])
 &\le(\rho^{3/2}+r^{3/2}+\sqrt N)M_j
                       \prod_l\|V_l\|_{\mathrm{HS}}.
 \end{split}
\]
Indeed, for a $k$ by $k$ matrix $B$, $\|B\|_1\le\sqrt k\|B\|_{\rm HS}\le k^{3/2}\max_{i,j}|B_{ij}|$. For the scalar trace use
$|\Tr B|\le\|I_N\|_{\rm HS}\|B\|_{\rm HS}=\sqrt N\|B\|_{\rm HS}$. These inequalities prove the three dimension factors in the display. For $0\le j\le3$, the scalar product rule
gives
\[
 M_j\le C_j\sum_{a=0}^j\rho_*^{2a+2}m_N^{2(j-a)}N^{o(1)}
       \le N^{(2j+2)/3+o(1)}.
\]
The factor $\rho^{3/2}+r^{3/2}+\sqrt N$ is $O(N^{1/2})$. For $W=N^{-1/2}T_{ab}(\zeta)$, the preceding matrix derivative bound consequently yields
\[
 \mathfrak n(D_H^j\mathsf X[W,\ldots,W])
 \le N^{-j/2+1/2+(2j+2)/3+o(1)}|\zeta|^j
 =N^{(j+7)/6+o(1)}|\zeta|^j.
\]
The Taylor factorials only decrease these fixed order bounds. In particular the exponents for $\ell,q,c$ are $4/3,3/2,5/3$. Apply Lemma~\ref{W:lem:localization-derivatives} at
$X=\mathsf X(H,w)$ to each cubic term:
\begin{align*}
 |D_{\mathsf X}\mathcal F[c]|&\le N^{5/3+o(1)}|\zeta|^3,\\
 |D_{\mathsf X}^2\mathcal F[\ell,q]|&\le N^{17/6+o(1)}|\zeta|^3,\\
 \tfrac16|D_{\mathsf X}^3\mathcal F[\ell,\ell,\ell]|
             &\le N^{4+o(1)}|\zeta|^3.
\end{align*}
Thus one common choice is $B=5$, independent of every failure exponent and valid for every Hermitian matrix. At orders one and two, the chain rule bounds are respectively
$N^{4/3+o(1)}|\zeta|$ and $(N^{3/2+o(1)}+N^{8/3+o(1)})|\zeta|^2$. Their coefficient tensors are consequently bounded by $N^5$ for sufficiently large $N$. The third moment of the
entry distribution is uniformly finite. In detail, let $M_3^+=\int|\zeta|^3(\nu_N^{\rm b}+\nu^G) (\dd\zeta)\le C$. For a failure probability bounded by $N^{-D_*}$, the discarded
sum is at most
\[
 C M_3^+N^{B+2}\mathbb P(\mathcal E^c)
       \le C N^{B+2-D_*}.
\]
The exponent $B$ in the matrix entry derivative bound is fixed independently of $D_*$. Choose $D_*>B+D+2$ after deriving that exponent and before the large-$N$ limit. Each cubic
component satisfies the displayed bound $N^B|\zeta|^3$, so this expectation estimate applies separately to all three components. It yields
\begin{equation}
 \E\mathbf1_{\mathcal E^c}\sum_{a<b}
   \left|\int D_{\mathsf X}\mathcal{F}(\mathsf{X}(K_{ab}))[c_{ab}]\,\Delta\nu\right|
      \le C_DN^{-D},
 \label{W:eq:localization-bad-event}
\end{equation}
and the same bound for either nonlinear cubic component.

On $\mathcal E$, the complete bounds on the spectral statistics are
\[
 \mathfrak{n}(\ell_{ab}(\zeta))\le Q_N{m_N}N^{-1/2}|\zeta|u_{ab},\qquad
 \mathfrak{n}(q_{ab}(\zeta))\le Q_N{m_N}N^{-1}|\zeta|^2v_{ab}.
\]
Therefore the two nonlinear terms have absolute sums over entries bounded by
\begin{align}
 \sum_{a<b}\left|\int D_{\mathsf X}^2\mathcal{F}[\ell,q]\,\Delta\nu\right|
      &\le Q_N{m_N}^2N^{-3/2}\sum u_{ab}v_{ab}
       \le N^{-1/6+o(1)},\label{W:eq:localization-second-derivative-term}\\
 \sum_{a<b}\left|\int\tfrac16D_{\mathsf X}^3\mathcal{F}[\ell^3]\,\Delta\nu\right|
      &\le Q_N{m_N}^3N^{-3/2}\sum u_{ab}^3
       \le N^{-1/2+o(1)}.\label{W:eq:localization-third-derivative-term}
\end{align}
Decompose the first derivative contribution into the difference between its values at $K_{ab}$ and $H(t)$, the matrix element integral, and the trace integral in
\eqref{W:eq:localization-common-response}. Their estimates yield
\begin{align*}
 \mathbf1_{\mathcal E}|Z_N^{(1)}|
 &\le Q_N\bigl(m_N^{-1}+m_N^{-2}+m_N^{-5/2}
                       +N^{-1/6}+N^{-1/2}+N^{-1/6}\bigr),\\
 \mathbf1_{\mathcal E}|Z_N^{(3)}-Z_N^{(1)}|
 &\le Q_N(N^{-1/6}+N^{-1/2}).
\end{align*}
Every term is bounded by $N^{-1/6+o(1)}$. Taking expectations and using \eqref{W:eq:localization-bad-event} for each cubic component proves
\eqref{W:eq:localization-cubic-estimate}. All inequalities are uniform for fixed $t$ and deterministic $w$, so the indicated suprema preserve the bounds.
\end{proof}

\begin{center}
\begin{tabular}{@{}lc@{}}
\toprule
Contribution & Bound up to $N^{o(1)}$\\
\midrule
Change in $c_{ab}$ from $K_{ab}$ to $H$ & $N^{-1/3}$\\
Variation of the derivative coefficient & $N^{-2/3}+N^{-5/6}$\\
Leading moment sum and interior error & $N^{-1/6}$\\
Signed trace terms with one error factor & $N^{-1/2}$\\
Trace terms with at least two error factors & $N^{-1/6}$\\
$D_{\mathsf X}^2\mathcal{F}[\ell,q]$ & $N^{-1/6}$\\
$D_{\mathsf X}^3\mathcal{F}[\ell^3]/6$ & $N^{-1/2}$\\
\bottomrule
\end{tabular}
\end{center}

\subsubsection{Taylor remainders for the localization statistic}
\label{W:sec:localization-taylor-remainders}

The polynomial used in the localization comparison includes the terms whose coordinate moments are summed with signs. We now express the difference from the complete statistic by
integral remainders, so that the absolute moment estimates can be applied to that difference.

For a matrix with one entry set to zero, put $X=\mathsf X(K,w)$ and write its exact moment increment as $\Delta X=\ell+q+c+r_4$. Let $h=\Delta X-\ell$ and define
\[
 P(X;x)=\mathcal F(X)+D\mathcal F(X)[\ell+q+c]
                     +\tfrac12D^2\mathcal F(X)[\ell,\ell].
\]
We use $P(X;x)$ only in this subsubsection. Its linear and quadratic coefficients are those of the expansion, and its cubic term is $D\mathcal F[c]$. The terms
$D^2\mathcal F[\ell,q]$ and $D^3\mathcal F[\ell^3]/6$ are included in the absolute remainder. We have estimated their full signed sum in Theorem~\ref{W:thm:localization-cubic} as
well, but the comparison below uses only $D\mathcal F[c]$. Thus each nonlinear cubic term is counted once, in that remainder.

Apply Lemma~\ref{W:lem:observable-large} with $\Phi=\mathcal F$, the norm $\mathfrak n$ in \eqref{W:eq:localization-norm}, and the local derivative bounds of
Lemma~\ref{W:lem:localization-derivatives}. With $a\ge\mathfrak n(\ell)$, $b\ge\max\{\mathfrak n(q),\mathfrak n(c),\mathfrak n(h)\}$, and $c_4\ge\mathfrak n(r_4)$, its conclusion
is
\[
 |\mathcal F(X+\Delta X)-P(X;x)|
 \le C_\kappa(c_4+a^3+ab+b^2).
\]
Here $P(X;x)=\mathcal F(X)+P_3(X;x)$, with $P_3$ as in \eqref{W:eq:observable-P3}.

On the regularity events measurable with respect to $K$, the moment estimates of Proposition~\ref{W:prop:observable-increments} imply
\[
 a_\alpha=Cm_NN^{-1/2}\|x\|_1U_\alpha,\qquad
 b_\alpha=Cm_NN^{-1}\|x\|_1^2V_\alpha,\qquad
 c_{4,\alpha}=Cm_NN^{-2}\|x\|_1^4V_\alpha.
\]
The bound for $c$ is included in $b_\alpha$, since every permitted entry perturbation has norm at most $Ce_N<1$. The norm in Proposition~\ref{W:prop:observable-increments}
dominates $\mathfrak n$. For the standardized off-diagonal distribution, $M_3^+=\int\|x\|_1^3(\mu_\alpha+\nu_\alpha)(\dd x)=O(1)$ and $M_4^+=N^{o(1)}$. The weights are measurable
with respect to the remaining entries. Hence their multiplication by the moments of $x$ is justified by conditional integration. The four sums in
\eqref{W:eq:observable-largeremainder} have bounds
\begin{align*}
 \sum_\alpha\mathbb E\int c_{4,\alpha}\,\dd(\mu_\alpha+\nu_\alpha)
 &\le N^{o(1)}m_NN^{-2}m_N^4=N^{-1/3+o(1)},\\
 \sum_\alpha\mathbb E\int a_\alpha^3\,\dd(\mu_\alpha+\nu_\alpha)
 &\le N^{o(1)}m_N^3N^{-3/2}m_N=N^{-1/6+o(1)},\\
 \sum_\alpha\mathbb E\int a_\alpha b_\alpha\,\dd(\mu_\alpha+\nu_\alpha)
 &\le N^{o(1)}m_N^2N^{-3/2}m_N^2=N^{-1/6+o(1)},\\
 \sum_\alpha\mathbb E\int b_\alpha^2\,\dd(\mu_\alpha+\nu_\alpha)
 &\le N^{o(1)}m_N^2N^{-2}m_N^3=N^{-1/3+o(1)}.
\end{align*}
The sums here are over off-diagonal blocks and include their regularity indicators. Each bound remains valid after integration of the interpolation parameter over $[0,1]$.

For a diagonal block, the exact resolvent identity gives
\[
 \mathfrak n(\Delta X)\le Cm_NN^{-1/2}|x|V_i,\qquad
 \mathfrak n(\Delta X-\ell)\le Cm_NN^{-1}|x|^2V_i.
\]
When $\mathfrak n(\Delta X)\le1$, the second-order integral formula and the local derivative bound imply
\[
 |\mathcal F(X+\Delta X)-\mathcal F(X)-D\mathcal F(X)[\ell]|
 \le C_\kappa\{\mathfrak n(\Delta X-\ell)+\mathfrak n(\Delta X)^2\}.
\]
When $\mathfrak n(\Delta X)>1$, the left side is at most $2+C_\kappa(\mathfrak n(\Delta X)+\mathfrak n(\Delta X-\ell))$, which has the same upper bound after increasing $C_\kappa$.
Thus, throughout the support of the replacement entry distributions, it is at most
\[
 C_\kappa N^{-1}|x|^2(m_NV_i+m_N^2V_i^2).
\]
Conditional integration and $\sum_i\mathbb EV_i\le m_NN^{o(1)}$, $\sum_i\mathbb EV_i^2\le N^{o(1)}$ bound the diagonal remainder by $N^{-1/3+o(1)}$.

To estimate the contribution from the complement of the regularity event, let $\mathcal B_\alpha$ denote the event where the stated spectral and resolvent bounds hold for the
matrix with entry $\alpha$ set to zero. On $\mathcal B_\alpha$, absence of a permitted entry value giving the upper bound on the Gram matrix makes the function, including its
cutoffs, identically zero in a neighbourhood of the support of the replacement entry distributions. On $\mathcal B_\alpha^c$, its difference between two values is bounded by $2$.
The global matrix derivative calculation in the proof of Theorem~\ref{W:thm:localization-cubic} bounds each retained coefficient by $N^5$. Therefore the difference between the
statistic and its retained polynomial is at most $C N^5(1+\|x\|_1^3)$ off the diagonal, and $C N^5(1+|x|)$ on the diagonal. Measurability of $\mathcal B_\alpha$ with respect to the
remaining entries gives a total bound $N^{7-D+o(1)}$ if $\sup_\alpha\mathbb P(\mathcal B_\alpha^c)\le N^{-D}$. Choosing $D$ after the prescribed error exponent makes this
contribution smaller than any prescribed inverse power of the matrix size. These calculations give the remainder estimates for all bases, with the coefficients integrated against
the full entry distributions.

\subsubsection{Comparison of the localization statistic}
\label{W:sec:comparison-localization-statistic}

\begin{proof}[Proof of Corollary~\ref{W:cor:localization-comparison}]
Fix a deterministic unit vector $w$ and write $m=m_N$ and $F(H)=\mathcal F(\mathsf X(H,w))$. We compare $H_N^{\rm tr}$ with $H_N^{\rm b}$, and then $H_N^{\rm b}$ with
$H_N^{\rm G}$. The bounded ensemble is the one defined in Section~\ref{W:sec:taylor-remainders-and-bounded-entries}; in particular, its standardized diagonal entries are
independent symmetric signs.

For either pair, let $\mu_{0,\alpha},\mu_{1,\alpha}$ be the standardized entry laws for the matrix perturbation $N^{-1/2}\mathcal T_\alpha x$ in
\eqref{W:eq:standardized-entry-insertion}, let $K_\alpha(t)$ have block $\alpha$ set to zero and remaining product law
$\bigotimes_{\delta\ne\alpha}[(1-t)\mu_{0,\delta}+t\mu_{1,\delta}]$, and put $\Delta M_{j,\alpha}=M_j(\mu_{1,\alpha})-M_j(\mu_{0,\alpha})$. An off-diagonal block has $\beta$ real
coordinates. At $X=\mathsf X(K_\alpha(t),w)$, use the coefficients of \eqref{W:eq:localization-moment-coefficients}, including their factorials, and write
\[
 \ell(x)=\sum_i x_i\ell_i,\qquad
 q(x)=\sum_{i,j}x_ix_jq_{ij},\qquad
 c(x)=\sum_{i,j,k}x_ix_jx_kc_{ijk}.
\]
Define symmetric coefficient tensors by
\[
 \begin{aligned}
 (C_{1,\alpha})_i&=D\mathcal F(X)[\ell_i],\\
 (C_{2,\alpha})_{ij}&=D\mathcal F(X)[q_{ij}]
                      +\tfrac12D^2\mathcal F(X)[\ell_i,\ell_j],\\
 (C_{3,\alpha})_{ijk}&=D\mathcal F(X)[c_{ijk}].
 \end{aligned}
\]
For an off-diagonal block, define $\mathcal R_\alpha(K,x)$ by the exact identity
\[
 F(K+N^{-1/2}\mathcal T_\alpha x)-F(K)
   =\sum_{j=1}^3 C_{j,\alpha}(K):x^{\otimes j}
       +\mathcal R_\alpha(K,x).
\]
For a diagonal entry the polynomial is $C_{1,i}(K)x$. Thus $\mathcal R_\alpha$ includes the two nonlinear cubic terms in \eqref{W:eq:localization-complete-cubic}. Here $:$ sums
over all ordered indices, $\|C\|_{\max}$ is the largest absolute coefficient, and $\|M\|_{1,\rm ord}$ is the sum of the absolute coefficients. In particular,
\[
 |C:M|\le\|C\|_{\max}\|M\|_{1,\rm ord}.
\]

For the remainder estimates, first replace each Gaussian distribution by
\[
 \widehat\mu^G=\mathbf1_{D_N}\mu^G+\mu^G(D_N^c)\delta_0,
 \qquad D_N=\{x:\|x\|_1\le L_N\},
\]
where $L_N=e_N\sqrt N$ is the threshold in \eqref{W:eq:LN-definition}. Leave every other law unchanged. Symmetry and the Gaussian tail integral give
\[
 \begin{gathered}
 M_1(\widehat\mu^G)=0,\qquad
 \|M_2(\widehat\mu^G)-M_2(\mu^G)\|_{1,\rm ord}
     \le C e^{-cL_N^2},\\
 \int\|x\|_1^j\,\widehat\mu^G(\dd x)
     \le\int\|x\|_1^j\,\mu^G(\dd x)
     \quad(j\ge1).
 \end{gathered}
\]
The matrix perturbation has norm $O(e_N)$, so these laws satisfy the same support and moment bounds as the truncated base laws. Denote them by $\widehat\mu_{j,\alpha}$ and their
base matrix by $\widehat K_\alpha(t)$. Use $I_\alpha=\mathbf1_{\mathcal B_\alpha\cap\mathcal A_{\alpha,\rho}}$ from Proposition~\ref{W:prop:observable-increments}; both events are
measurable with respect to this base. With $M_j^+$ bounding the corresponding $j$th absolute moments of $\widehat\mu_{0,\alpha}+\widehat\mu_{1,\alpha}$,
Lemma~\ref{W:lem:observable-large} and \eqref{W:eq:observable-weight-sums} give
\[
 \begin{aligned}
 &\sum_{\alpha\,\rm off}\mathbb E I_\alpha
        \int|\mathcal R_\alpha(\widehat K_\alpha,x)|
                   (\widehat\mu_{0,\alpha}+\widehat\mu_{1,\alpha})(\dd x)\\
 &\quad\le C Q_N\left\{
 M_4^+mN^{-2}m^4+M_3^+m^3N^{-3/2}m\right.\\
 &\hspace{7em}\left.{}+M_3^+m^2N^{-3/2}m^2
                         +M_4^+m^2N^{-2}m^3\right\}
       \le N^{-1/6+o(1)},
 \end{aligned}
\]
where $Q_N=N^{o(1)}$, $M_3^+=O(1)$, and $M_4^+=N^{o(1)}$. The diagonal calculation in Theorem~\ref{W:thm:localization-taylor-remainder} gives
\[
 \begin{aligned}
 \sum_i\mathbb E I_i\int|\mathcal R_i(\widehat K_i,x)|
                       (\widehat\mu_{0,i}+\widehat\mu_{1,i})(\dd x)
 &\le C Q_N M_{2,\rm d}^+N^{-1}(m^2+m^2)\\
 &\le N^{-1/3+o(1)}.
 \end{aligned}
\]
Both bounds are uniform in $t$ and in deterministic unit $w$. Until the Gaussian entries are restored below, evaluate $\mathcal B_\alpha$ and $\mathcal A_{\alpha,\rho}$ at
$\widehat K_\alpha(t)$.

We first include the other base matrices. If $K\in\mathcal B_\alpha\cap\mathcal A_{\alpha,\rho}^c$, the definition \eqref{W:eq:resolvent-gram-bound} and compactness of the ball
$\mathcal X_\alpha$ give
\[
 \min_{x\in\mathcal X_\alpha}
 \lambda_{\max}\bigl(G_\rho^{[d]}(K+N^{-1/2}\mathcal T_\alpha x,w)\bigr)
       >\kappa^{-1}.
\]
Here the minimum is attained because the modified Gram matrix is continuous in $x$. Every replacement entry distribution currently in use has support in the interior of
$\mathcal X_\alpha$, and $0$ also lies in that interior. The factor $\Gamma_\rho$ is therefore zero at every permitted entry value. At $x=0$, the strict inequality also makes
$\mathcal F$ identically zero on a neighbourhood of $\mathsf X(K,w)$ in the real coordinate space. Consequently
\[
 \begin{gathered}
 F(K+N^{-1/2}\mathcal T_\alpha x)=F(K)=0
       \quad(x\in\mathcal X_\alpha),\\
 D_{\mathsf X}^{j}\mathcal F(\mathsf X(K,w))=0\quad(1\le j\le3),
 \qquad C_{j,\alpha}(K)=0,\\
 \mathcal R_\alpha(K,x)=0\quad(x\in\mathcal X_\alpha).
 \end{gathered}
\]
The coefficient identity applies to each term of the polynomial defined for that entry block. In particular, $D\mathcal F[c]$ and the two nonlinear cubic contributions vanish
separately on this event.

For $K\notin\mathcal B_\alpha$, use the matrix derivative bounds at the start of the proof of Theorem~\ref{W:thm:localization-cubic}. In the norm \eqref{W:eq:localization-norm}
they state, for all Hermitian $K$,
\[
 \begin{aligned}
 \sup_{\|u\|_1\le1}\mathfrak n(\ell(u))&\le N^{4/3+o(1)},&
 \sup_{\|u\|_1\le1}\mathfrak n(q(u))&\le N^{3/2+o(1)},&
 \sup_{\|u\|_1\le1}\mathfrak n(c(u))&\le N^{5/3+o(1)}.
 \end{aligned}
\]
These bounds also hold in a diagonal direction, whose perturbation matrix has no larger Hilbert--Schmidt norm. Substitution in the definitions of the retained coefficients and
Lemma~\ref{W:lem:localization-derivatives}, followed by polarization in at most two real coordinates, gives
\[
 \begin{aligned}
 \|C_{1,\alpha}(K)\|_{\max}&\le C_\kappa N^{4/3+o(1)},\\
 \|C_{2,\alpha}(K)\|_{\max}
   &\le C_\kappa\bigl(N^{3/2+o(1)}+N^{8/3+o(1)}\bigr),\\
 \|C_{3,\alpha}(K)\|_{\max}&\le C_\kappa N^{5/3+o(1)}.
 \end{aligned}
\]
Thus all these coefficients are bounded by $N^5$ for sufficiently large $N$, uniformly in $K$ and $w$. The exponent $5$ and this dimension threshold do not depend on the failure
exponent $D$. Since $\sum_{i_1,\ldots,i_j}|x_{i_1}\cdots x_{i_j}|=\|x\|_1^j$ and $0\le F\le1$, the exact remainder identities imply
\[
 \begin{aligned}
 |\mathcal R_\alpha(K,x)|
 &\le1+N^5\sum_{j=1}^3\|x\|_1^j
       \le C N^5(1+\|x\|_1^3)
                 &&(\alpha\text{ off the diagonal}),\\
 |\mathcal R_i(K,x)|&\le1+N^5|x|
       \le C N^5(1+|x|).
 \end{aligned}
\]
All constants below may depend on the fixed cutoff parameters, the parameters of the functions under consideration and the uniform bounds on entry moments, but not on $t$, $w$ or
the block index. Constants from the regularity estimates may also depend on the fixed $D$. The off-diagonal third absolute moments are uniformly bounded. On the diagonal only the
first absolute moment is used here: for each of the two entry distributions,
$\int|x|\,\mu_{j,i}(\dd x)\le(\int|x|^2\,\mu_{j,i}(\dd x))^{1/2}\le C$.
The auxiliary laws satisfy the same bounds.

For each fixed $D>0$, the regularity estimate in Section~\ref{W:sec:contour-integrals-for-moment-increments} applies to the auxiliary product laws and gives
\[
 \sup_{t\in[0,1]}\sup_{\|w\|=1}\sup_\alpha
 \mathbb P\bigl(\mathcal B_\alpha(\widehat K_\alpha(t))^c\bigr)\le N^{-D}.
\]
Let $\mathscr F_\alpha$ be generated by the remaining auxiliary entries. The base and its regularity event are $\mathscr F_\alpha$-measurable, whereas the replacement variable has
the deterministic law $\widehat\mu_{j,\alpha}$ independently of $\mathscr F_\alpha$. Conditional integration of the preceding pointwise bounds therefore gives
\[
 \begin{aligned}
 &\mathbb E\!\left[\mathbf1_{\mathcal B_\alpha^c}
     \int|\mathcal R_\alpha(\widehat K_\alpha,x)|
                (\widehat\mu_{0,\alpha}+\widehat\mu_{1,\alpha})(\dd x)\right]\\
 &\qquad\le C N^5\mathbb P(\mathcal B_\alpha^c)
       \int(1+\|x\|_1^3)
                (\widehat\mu_{0,\alpha}+\widehat\mu_{1,\alpha})(\dd x)
       \le C N^{5-D}
              \quad(\alpha\text{ off the diagonal}),\\
 &\mathbb E\!\left[\mathbf1_{\mathcal B_i^c}
     \int|\mathcal R_i(\widehat K_i,x)|
                (\widehat\mu_{0,i}+\widehat\mu_{1,i})(\dd x)\right]
       \le C N^{5-D}.
 \end{aligned}
\]
There are $N(N-1)/2$ off-diagonal blocks and $N$ diagonal blocks. The identity
\[
 1=I_\alpha+\mathbf1_{\mathcal B_\alpha\cap\mathcal A_{\alpha,\rho}^c}
                   +\mathbf1_{\mathcal B_\alpha^c}
\]
now yields, for each $t$ and $w$,
\[
 \sum_\alpha\mathbb E\int|\mathcal R_\alpha(\widehat K_\alpha(t),x)|
              (\widehat\mu_{0,\alpha}+\widehat\mu_{1,\alpha})(\dd x)
       \le N^{-1/6+o(1)}+C N^{7-D}.
\]

We next restore the Gaussian entries. For each standardized Gaussian block $g$, its real dimension is at most two and its covariance is one of the fixed GOE or GUE covariances. The
Gaussian density and $\|g\|_1^2\le2\|g\|_2^2$ imply that, for some fixed $a>0$,
\[
 \sup_g\mathbb E\bigl[(1+\|g\|_1^3)e^{a\|g\|_1^2}\bigr]\le C.
\]
Indeed, the density is bounded by $C e^{-b\|x\|_2^2}$ with $b>0$ fixed, so any $a<b/2$ makes the integral finite. It follows that
\[
 \int_{D_N^c}(1+\|x\|_1^3)\,\mu^G(\dd x)
 \le e^{-aL_N^2}
       \mathbb E\bigl[(1+\|g\|_1^3)e^{a\|g\|_1^2}\bigr]
 \le C e^{-aL_N^2}.
\]
The same bound holds with $1+|x|$ for a Gaussian diagonal entry.

Here is a joint construction of the original and auxiliary base matrices. For every block $\delta$, independently sample $U_\delta$ uniformly on $[0,1]$ and $\xi_{j,\delta}$ with
law $\mu_{j,\delta}$, $j=0,1$. Put
\[
 \widehat\xi_{j,\delta}=
 \begin{cases}
 \xi_{j,\delta}\mathbf1_{\{\|\xi_{j,\delta}\|_1\le L_N\}},
                         &\mu_{j,\delta}\text{ Gaussian},\\
 \xi_{j,\delta},          &\text{otherwise},
 \end{cases}
\]
and define
\[
 \begin{aligned}
 K_\alpha(t)&=N^{-1/2}\sum_{\delta\ne\alpha}\mathcal T_\delta
    \bigl(\mathbf1_{\{U_\delta\le t\}}\xi_{1,\delta}
           +\mathbf1_{\{U_\delta>t\}}\xi_{0,\delta}\bigr),\\
 \widehat K_\alpha(t)&=N^{-1/2}\sum_{\delta\ne\alpha}\mathcal T_\delta
    \bigl(\mathbf1_{\{U_\delta\le t\}}\widehat\xi_{1,\delta}
           +\mathbf1_{\{U_\delta>t\}}\widehat\xi_{0,\delta}\bigr).
 \end{aligned}
\]
Their marginals are precisely the two product laws used above. A complex block is sampled as a single vector in $\mathbb R^2$ with its joint distribution. Both matrices are
measurable with respect to
$\mathscr G_\alpha=\sigma(U_\delta,\xi_{0,\delta},\xi_{1,\delta}:\delta\ne\alpha)$.
Writing $E_\alpha(t)=\{K_\alpha(t)\ne\widehat K_\alpha(t)\}$, a union bound gives
\[
 \begin{aligned}
 \mathbb P(E_\alpha(t))
 &\le\sum_{\delta\ne\alpha}\bigl[
 (1-t)\mathbb P(\xi_{0,\delta}\ne\widehat\xi_{0,\delta})
       +t\mathbb P(\xi_{1,\delta}\ne\widehat\xi_{1,\delta})\bigr]\\
 &\le C N^2e^{-aL_N^2}.
 \end{aligned}
\]
Only a selected Gaussian entry outside $D_N$ can contribute to this sum. On $E_\alpha(t)^c$ the two remainder functions agree for every replacement value. The global bounds
consequently give
\[
 \begin{aligned}
 |\mathcal R_\alpha(K_\alpha(t),x)-\mathcal R_\alpha(\widehat K_\alpha(t),x)|
 &\le C N^5\mathbf1_{E_\alpha(t)}(1+\|x\|_1^3)
                       &&(\alpha\text{ off the diagonal}),\\
 |\mathcal R_i(K_i(t),x)-\mathcal R_i(\widehat K_i(t),x)|
 &\le C N^5\mathbf1_{E_i(t)}(1+|x|).
 \end{aligned}
\]
In each integral below, $x$ is independent of $\mathscr G_\alpha$. Integrating these bounds against the original entry distributions and then taking expectations gives
\[
 \begin{aligned}
 &\sum_{\alpha\,\rm off}\mathbb E\int
  |\mathcal R_\alpha(K_\alpha(t),x)-\mathcal R_\alpha(\widehat K_\alpha(t),x)|
                   (\mu_{0,\alpha}+\mu_{1,\alpha})(\dd x)\\
 &\qquad\le C N^5\sum_{\alpha\,\rm off}\mathbb P(E_\alpha(t))
       \int(1+\|x\|_1^3)(\mu_{0,\alpha}+\mu_{1,\alpha})(\dd x)
       \le C N^9 e^{-aL_N^2},\\
 &\sum_i\mathbb E\int
  |\mathcal R_i(K_i(t),x)-\mathcal R_i(\widehat K_i(t),x)|
                   (\mu_{0,i}+\mu_{1,i})(\dd x)\\
 &\qquad\le C N^5\sum_i\mathbb P(E_i(t))
       \int(1+|x|)(\mu_{0,i}+\mu_{1,i})(\dd x)
       \le C N^8 e^{-aL_N^2}.
 \end{aligned}
\]
Thus the sum over all blocks is at most $C N^9e^{-aL_N^2}$.

It remains to restore the law of the replacement Gaussian variable. For every fixed $K$, $\mathcal R_\alpha(K,0)=0$ by its definition. Hence, for a Gaussian entry law,
\[
 \begin{aligned}
 &\int|\mathcal R_\alpha(K,x)|\,\mu^G(\dd x)
       -\int|\mathcal R_\alpha(K,x)|\,\widehat\mu^G(\dd x)\\
 &\qquad=\int_{D_N^c}|\mathcal R_\alpha(K,x)|\,\mu^G(\dd x)
       \le C N^5e^{-aL_N^2}.
 \end{aligned}
\]
The off-diagonal estimate uses $1+\|x\|_1^3$, and the diagonal estimate uses $1+|x|$. The sum over the at most two Gaussian entry distributions for each block and all blocks is
therefore at most $C N^7e^{-aL_N^2}$, uniformly in $K$. Apply this identity at $\widehat K_\alpha(t)$, and then apply the bound for replacing $\widehat K_\alpha(t)$ by
$K_\alpha(t)$. We obtain
\[
 \begin{aligned}
 &\sum_\alpha\mathbb E\int|\mathcal R_\alpha(K_\alpha(t),x)|
                  (\mu_{0,\alpha}+\mu_{1,\alpha})(\dd x)\\
 &\qquad\le N^{-1/6+o(1)}+C N^{7-D}
                       +C(N^7+N^9)e^{-aL_N^2}.
 \end{aligned}
\]
For any fixed $A>0$, choose $D>A+8$. Since $L_N^2=N^{1-o(1)}$, choose $N_0(A,a)$ so that $L_N^2\ge N^{1/2}$ and $aN^{1/2}\ge(A+9)\log N$ for every $N\ge N_0(A,a)$. Thus
\[
 N^{7-D}\le N^{-A-1},\qquad
 N^9e^{-aL_N^2}\le N^{-A}.
\]
Taking, for example, $A=1$ and integrating the uniform bound over $t\in[0,1]$ proves
\[
 \int_0^1\sum_\alpha\mathbb E\int
     |\mathcal R_\alpha(K_\alpha(t),x)|
              (\mu_{0,\alpha}+\mu_{1,\alpha})(\dd x)\,\dd t
       \le N^{-1/6+o(1)}.
\]

We now apply the exact product-measure interpolation to the original pair of entry distributions. Set $\lambda_{\alpha,t}=(1-t)\mu_{0,\alpha}+t\mu_{1,\alpha}$. The function
\[
 \Phi(t)=\int F\!\left(N^{-1/2}\sum_\alpha\mathcal T_\alpha x_\alpha\right)
                         \prod_\alpha\lambda_{\alpha,t}(\dd x_\alpha)
\]
is a polynomial in $t$: each of the finitely many measure factors is $(1-t)\mu_{0,\alpha}+t\mu_{1,\alpha}$ and $F$ is bounded. Differentiating one factor at a time, as in
Lemma~\ref{W:lem:entry-distribution-interpolation}, yields
\[
 \Phi'(t)=\sum_\alpha\mathbb E\int
       [F(K_\alpha(t)+N^{-1/2}\mathcal T_\alpha x)-F(K_\alpha(t))]
                         (\mu_{1,\alpha}-\mu_{0,\alpha})(\dd x).
\]
The subtracted term integrates to zero because the signed measure has total mass zero. Conditional on the remaining entries, each coefficient tensor is fixed, so
\[
 \int C_{j,\alpha}(K):x^{\otimes j}
                     (\mu_{1,\alpha}-\mu_{0,\alpha})(\dd x)
    =C_{j,\alpha}(K):\Delta M_{j,\alpha}.
\]
The required absolute-integrability estimate is
\[
 \begin{aligned}
 &\int_0^1\left\{
 \sum_{\alpha\,\rm off}\sum_{j=1}^3\mathbb E\int
       |C_{j,\alpha}(K_\alpha(t)):x^{\otimes j}|
                  (\mu_{0,\alpha}+\mu_{1,\alpha})(\dd x)\right.\\
 &\hspace{6em}\left.{}+\sum_i\mathbb E\int|C_{1,i}(K_i(t))x|
                  (\mu_{0,i}+\mu_{1,i})(\dd x)\right\}\,\dd t
       \le C N^5(N^2+N)\le C N^7.
 \end{aligned}
\]
This uses the off-diagonal third absolute moments and diagonal first absolute moments. Together with the preceding remainder estimate, it justifies Fubini's theorem for every term.
Substituting the exact remainder identities into $\Phi'(t)$ and integrating gives
\[
 \begin{aligned}
 \mathbb EF(H_1)-\mathbb EF(H_0)
 ={}&\int_0^1
 \sum_{\alpha\,\rm off}\sum_{j=1}^3
       \mathbb E C_{j,\alpha}(K_\alpha(t)):\Delta M_{j,\alpha}\,\dd t\\
 &+\int_0^1\sum_i\mathbb E C_{1,i}(K_i(t))\Delta M_{1,i}\,\dd t
       +\varepsilon_N,\\
 \varepsilon_N={}&\int_0^1\sum_\alpha\mathbb E\int
      \mathcal R_\alpha(K_\alpha(t),x)
                       (\mu_{1,\alpha}-\mu_{0,\alpha})(\dd x)\,\dd t,\\
 |\varepsilon_N|\le{}&\int_0^1\sum_\alpha\mathbb E\int
      |\mathcal R_\alpha(K_\alpha(t),x)|
                       (\mu_{1,\alpha}+\mu_{0,\alpha})(\dd x)\,\dd t
       \le N^{-1/6+o(1)}.
 \end{aligned}
\]
The auxiliary measures were used only to estimate this absolute remainder. Every $\Delta M_{j,\alpha}$ in the identity is the moment difference of the original pair of entry
distributions.

For $(H_0,H_1)=(H_N^{\rm G},H_N^{\rm b})$, the matched off-diagonal moments and the centered diagonal laws imply
\[
 \Delta M_{1,\alpha}=0,\qquad\Delta M_{2,\alpha}=0
       \quad(\alpha\text{ off the diagonal}),\qquad
 \Delta M_{1,i}=0.
\]
Use the same matrix $H(t)$ to construct every $K_{ab}(t)$ by setting the corresponding entry to zero. By the definition of $C_{3,ab}$,
\[
 \begin{aligned}
 \sum_{a<b}C_{3,ab}(K_{ab}(t)):\Delta M_{3,ab}
 &=\sum_{a<b}\int D\mathcal F(\mathsf X(K_{ab}(t),w))
                    [c_{ab}(K_{ab}(t);x)]
                          (\nu_N^{\rm b}-\nu^G)(\dd x)\\
 &=Z_N^{(1)}(t,w).
 \end{aligned}
\]
Here $x$ represents the whole real coordinate block of the entry. The exact interpolation identity reduces to
\[
 \mathbb EF(H_N^{\rm b})-\mathbb EF(H_N^{\rm G})
       =\int_0^1\mathbb E Z_N^{(1)}(t,w)\,\dd t+\varepsilon_N.
\]
Taking the absolute value after summing the cubic terms and applying Theorem~\ref{W:thm:localization-cubic} proves
\[
 \begin{aligned}
 \left|\mathbb EF(H_N^{\rm b})-\mathbb EF(H_N^{\rm G})\right|
 &\le\int_0^1\mathbb E|Z_N^{(1)}(t,w)|\,\dd t+|\varepsilon_N|\\
 &\le N^{-1/6+o(1)}.
 \end{aligned}
\]
The diagonal second moments enter only the remainder estimate. In fact, $\int|x|^2(\mu_{0,i}+\mu_{1,i})(\dd x)=3$ for the GOE--sign pair and $2$ for the GUE--sign pair, so their
possible difference is already covered by $M_{2,\rm d}^+$.

For $(H_0,H_1)=(H_N^{\rm b},H_N^{\rm tr})$, \eqref{W:eq:replacement-tail-defects} and Lemma~\ref{W:lem:cubic-defect} give
\[
 \|\Delta M_1\|_{1,\rm ord}\le N^{-3/2+o(1)},\qquad
 \|\Delta M_2\|_{1,\rm ord}+\|\Delta M_3\|_{1,\rm ord}
          \le N^{-1+o(1)}.
\]
Proposition~\ref{W:prop:observable-increments} and Lemma~\ref{W:lem:localization-derivatives} imply, on $I_\alpha=1$,
\[
 \begin{aligned}
 \|C_{1,\alpha}\|_{\max}&\le CmN^{-1/2}U_\alpha,\\
 \|C_{2,\alpha}\|_{\max}&\le CN^{-1}(mV_\alpha+m^2U_\alpha^2),\\
 \|C_{3,\alpha}\|_{\max}&\le CmN^{-3/2}V_\alpha.
 \end{aligned}
\]
The coefficient estimates follow componentwise from polarization:
\[
 \begin{aligned}
 q_{ij}&=\tfrac14\{q(e_i+e_j)-q(e_i-e_j)\},\\
 c_{ijk}&=\frac1{48}\sum_{\epsilon\in\{-1,1\}^3}
  \epsilon_1\epsilon_2\epsilon_3
     c(\epsilon_1e_i+\epsilon_2e_j+\epsilon_3e_k).
 \end{aligned}
\]
These formulas include the mixed components when $\beta=2$. Using \eqref{W:eq:observable-weight-sums}, the three contractions satisfy
\[
 \begin{aligned}
 \int_0^1\sum_{\alpha\,\rm off}\mathbb E I_\alpha
       |C_{1,\alpha}:\Delta M_{1,\alpha}|\,\dd t
   &\le Q_Nm^5N^{-2}=N^{-1/3+o(1)},\\
 \int_0^1\sum_{\alpha\,\rm off}\mathbb E I_\alpha
       |C_{2,\alpha}:\Delta M_{2,\alpha}|\,\dd t
   &\le Q_NN^{-2}(m^5+m^4)
                =N^{-1/3+o(1)},\\
 \int_0^1\sum_{\alpha\,\rm off}\mathbb E I_\alpha
       |C_{3,\alpha}:\Delta M_{3,\alpha}|\,\dd t
   &\le Q_Nm^5N^{-5/2}=N^{-5/6+o(1)}.
 \end{aligned}
\]
On $\mathcal B_\alpha\cap\mathcal A_{\alpha,\rho}^c$, these coefficients vanish by the argument above. For the remaining event, the global coefficient bound and the estimates for
the moment differences $\Delta M_{j,\alpha}$ give
\[
 \begin{aligned}
 &\int_0^1\sum_{\alpha\,\rm off}\sum_{j=1}^3
      \mathbb E\bigl[\mathbf1_{\mathcal B_\alpha^c}
             |C_{j,\alpha}:\Delta M_{j,\alpha}|\bigr] \,\dd t\\
 &\qquad\le N^5\int_0^1\sum_{\alpha\,\rm off}
        \mathbb P(\mathcal B_\alpha^c)
            \sum_{j=1}^3\|\Delta M_{j,\alpha}\|_{1,\rm ord}\,\dd t\\
 &\qquad\le N^5 N^2 N^{-D}N^{-1+o(1)}=N^{6-D+o(1)}.
 \end{aligned}
\]
For the diagonal, centering of the original entries and the definition of $\Theta_{\rm d}$ give
\[
 \begin{aligned}
 |\Delta M_{1,i}|&=\left|\mathbb E d_{N,i}
                       \mathbf1_{\{|d_{N,i}|>L_N\}}\right|
             \le L_N^{-1}\Theta_{\rm d}(L_N),\\
 \int_0^1\sum_i\mathbb E I_i
                  |C_{1,i}\Delta M_{1,i}|\,\dd t
   &\le Q_Nm^2N^{-1/2}L_N^{-1}\Theta_{\rm d}(L_N)
             \le N^{-1/3+o(1)}.
 \end{aligned}
\]
The diagonal contribution on the exceptional event is also controlled using only its difference of first moments:
\[
 \begin{aligned}
 &\int_0^1\sum_i\mathbb E\bigl[\mathbf1_{\mathcal B_i^c}
                       |C_{1,i}\Delta M_{1,i}|\bigr] \,\dd t\\
 &\qquad\le N^5 N N^{-D}L_N^{-1}\Theta_{\rm d}(L_N)
       \le N^{11/2-D+o(1)}.
 \end{aligned}
\]
Here $L_N^{-1}=N^{-1/2+o(1)}$ and $\Theta_{\rm d}(L_N)$ is bounded by the uniform diagonal second-moment bound. With the previously chosen $D>A+8$, both exceptional contributions
are $O(N^{-A})$. Substitution in the comparison identity proves
\[
 |\mathbb EF(H_N^{\rm tr})-\mathbb EF(H_N^{\rm b})|
       \le N^{-1/6+o(1)}.
\]
The triangle inequality proves \eqref{W:eq:localization-comparison}. Every bound used above is uniform in deterministic unit $w$.

Finally, $0\le F\le1$ and $H_N=H_N^{\rm tr}$ on $\Omega_N^{\rm tr}$ from \eqref{W:eq:equality-event}. Therefore
\[
 \sup_{\|w\|=1}|\mathbb EF(H_N)-\mathbb EF(H_N^{\rm tr})|
 \le\mathbb P((\Omega_N^{\rm tr})^c)=o(1).
\]
This proves the convergence for the original matrix.
\end{proof}

\subsubsection{Localization of the spectral projection}
\label{W:sec:proof-full-localization}

\begin{proof}[Proof of Theorem~\ref{W:thm:localization-spectral-projection}]
Use the projections and the event \(\Omega_d\) from Section~\ref{W:sec:localization-observable-and-derivatives}. The deterministic inequality \eqref{W:eq:localization-detection}
reduces the proof to three quantities: the probability of failure of spectral confinement, $\mathbb E(1-\Gamma_\rho)$, and the expectation of the bounded localization statistic.
All estimates below keep \(L,E,\eta,\kappa\) fixed until the matrix limit has been taken.

Theorem~\ref{W:thm:resolvent-window}, applied to the truncated array, and its coupling to $H_N$ give
\[
 p_{d,N}:=\mathbb P(\Omega_d^c)
 \le\mathbb P(H_N\ne H_N^{\rm tr})+
       \mathbb P(\|H_N^{\rm tr}\|>2+d/2)=o(1).
\]
This probability is independent of $w$. Define the two deterministic comparison errors for the original matrices by
\begin{align*}
 \delta_{\Gamma,N}&=
  \sup_{\|w\|=1}|\mathbb E\Gamma_\rho(H_N,w)
                    -\mathbb E\Gamma_\rho(H_N^{\rm G},w)|,\\
 \delta_{F,N}&=
  \sup_{\|w\|=1}|\mathbb E\mathcal F(\mathsf X(H_N,w))
                    -\mathbb E\mathcal F(\mathsf X(H_N^{\rm G},w))|.
\end{align*}
For fixed $L,E,\kappa$, Theorem~\ref{W:thm:tightness-main} and Corollary~\ref{W:cor:localization-comparison}, including the coupling and the spectral regularization, give
$\delta_{\Gamma,N}+\delta_{F,N}=o(1)$. Let
\[
 g_N^{\rm G}(\kappa)=\sup_{\|w\|=1}\mathbb E(1-\Gamma_\rho)(H_N^{\rm G},w),
 \qquad
 k_N^{\rm G}(L)=\sup_{\|w\|=1}\mathbb E\mathcal F(\mathsf X(H_N^{\rm G},w)).
\]
Both expectations are taken with the fixed values of $L,E,\kappa$ used above. The arguments displayed in $g_N^{\rm G}(\kappa)$ and $k_N^{\rm G}(L)$ indicate the parameters varied
in the next bounds.

We recall the estimate giving uniformity of the Gaussian bound in $\kappa$. On the coercivity and spectral events from
Section~\ref{W:sec:expectation-of-the-localization-statistic}, the weighted trace estimate gives
\[
 0\le\mathcal F(\mathsf X(H_N^{\rm G},w))
 \le \Tr\bigl(Q_{N,L}\theta_E(\mathcal A_N^{\rm G})\bigr)
 \le\frac{C_E}{L}(1+C_N)^{5/2},
 \qquad C_N=C_0(1+\mathfrak Z_N^4).
\]
The proof of Proposition~\ref{W:prop:gaussian-localization-expectation} establishes $\sup_N\mathbb E(1+C_N)^{5/2}<\infty$ and that the complement of those events has probability
tending to zero. Since $0\le\mathcal F\le1$ everywhere, taking expectations on the event and its complement separately yields
\[
 \limsup_N k_N^{\rm G}(L)\le C_{\beta,E}/L.
\]
For each fixed $L$, we also have $\lim_{\kappa\downarrow0}\limsup_Ng_N^{\rm G}(\kappa)=0$, by the bounds for the Gram matrix and its inverse under the GOE or GUE.

Taking expectations in \eqref{W:eq:localization-detection} and then taking the supremum over deterministic starting vectors gives the finite-dimensional inequality
\begin{equation}
 \sup_{\|w\|=1}\mathbb P(\|Q_{N,L}\Pi_N(E)\|>\eta) \le{}p_{d,N}+g_N^{\rm G}(\kappa)+\delta_{\Gamma,N} +\frac{4+\eta^2}{\eta^2} (k_N^{\rm G}(L)+\delta_{F,N}).
 \label{W:eq:localization-probability-bound}
\end{equation}
The constants in the comparison errors may depend on the fixed $L$ and $\kappa$, in accordance with the following order of limits. First let $N\to\infty$ with $L,E,\eta,\kappa$
fixed. Next let $\kappa\downarrow0$. The left side is independent of $\kappa$, so
\[
 \limsup_N\sup_{\|w\|=1}
       \mathbb P(\|Q_{N,L}\Pi_N(E)\|>\eta)
 \le\frac{(4+\eta^2)C_{\beta,E}}{\eta^2 L}.
\]
Finally let $L\to\infty$. This proves \eqref{W:eq:localization-spectral-projection}.

The same estimates quantify the statement concerning subsequent blocks. By the integer trace detection above, their low energy event has indicator at most
$\mathbf1_{\Omega_d^c}+1-\Gamma_\rho+5\mathcal F$. Consequently, with $\mathcal K_N=\mathcal K_N(H_N,w)$,
\[
 \sup_{\|w\|=1}\mathbb P\!\left(
 \operatorname{Ran}\mathbf1_{(-\infty,E]}(A_N^{\rm orig})
                   \not\subset\mathcal K_N\right)
 \le p_{d,N}+g_N^{\rm G}(\kappa)+\delta_{\Gamma,N}
                    +5\bigl(k_N^{\rm G}(L)+\delta_{F,N}\bigr).
\]
Taking $N\to\infty$, then $\kappa\downarrow0$, and finally $L\to\infty$ makes the right side tend to zero. This conclusion is independent of the choice of basis in the reducing
orthogonal complement. It concerns low energy states only and does not exclude Lanczos breakdown in the entire matrix.
\end{proof}

\subsection{Proof of the operator limit for Wigner matrices}
\label{W:sec:main-assembly}

The local coefficient limit and localization of the entire low energy spectral subspace are now available on the scales required by the deterministic criterion. Combining them
gives the operator limit for the continued Lanczos matrix, including the possibility of breakdown in its finite construction.

\begin{proof}[Proof of Theorem~\ref{W:thm:main-operator}]
Fix a deterministic sequence of unit vectors $w_N$. Proposition~\ref{W:prop:lanczos-construction} represents the whole matrix by its continued real symmetric Lanczos matrix, with
nonnegative off-diagonal coefficients. Proposition~\ref{W:prop:operator-equality-transfer} makes the embedded resolvent path Borel. The two hypotheses of
Theorem~\ref{W:thm:operator-probabilistic} to be verified are
\[
 \begin{aligned}
 &(\widetilde Y_N^{\rm lin},(B_N(r))_{r\ge1})
       \Longrightarrow (Y_\beta,0),\\
 &\lim_{L\to\infty}\limsup_N
   \mathbb P\bigl(\|Q_{N,L}\Pi_N(E)\|>\eta\bigr)=0
       \qquad(E\in\mathbb R,\ \eta>0).
 \end{aligned}
\]
The first is Theorem~\ref{W:thm:integrated-potential-and-offdiagonal-coefficients}, with the exact potential \eqref{W:eq:tightness-exact-grid-potential}; the second is
Theorem~\ref{W:thm:localization-spectral-projection} for the full spectral projection, including every invariant block. Both hold for the same matrix and continued basis. By
Lemma~\ref{W:lem:operator-brownian-form}, the limiting quadratic form is that of the stochastic Airy operator with Dirichlet boundary condition.
Theorem~\ref{W:thm:operator-probabilistic} therefore gives
\begin{equation}
 (\widetilde Y_N^{\rm lin},\mathscr R_N)
 \Longrightarrow(Y_\beta,\mathscr R_\beta)
 \quad\text{in }C_{\rm loc}([0,\infty))\times\mathscr Y.
 \label{W:eq:assembly-joint-convergence}
\end{equation}
Its second marginal is the required operator limit. The criterion also supplies the tight lower spectral bound.

If \eqref{W:eq:intro-uniform-distribution} failed, there would be $\varepsilon>0$, $N_j\uparrow\infty$ and deterministic unit vectors $w_j$ satisfying
\[
 d_{\rm BL}\bigl(\mathcal L(\mathscr R_{N_j}(H_{N_j},w_j)),
                         \mathcal L(\mathscr R_\beta)\bigr)\ge\varepsilon.
\]
Extend these choices to a deterministic sequence at the other dimensions. Equation~\eqref{W:eq:assembly-joint-convergence} then makes these distances tend to zero, a contradiction.

The product space is Polish by Lemma~\ref{W:lem:resolvent-path-space}, and its limiting resolvent is the Borel function of $Y_\beta$ specified in
Lemma~\ref{W:lem:operator-measurability}. Skorokhod representation provides copies converging almost surely in this product space, with each copied pair having the original joint
distribution. For its compact exhaustion $(K_j)$,
\[
 1\wedge\sup_{z\in K_j}\|\mathscr R_N(z)-\mathscr R_\beta(z)\|
 \le2^j d_{\mathscr Y}(\mathscr R_N,\mathscr R_\beta)
 \longrightarrow0.
\]
Every compact nonreal parameter set lies in some $K_j$, proving \eqref{W:eq:skorokhod-uniform-resolvent} on one event for the fixed starting sequence.
\end{proof}

\section{Spectral consequences and dependence on parameters}\label{sec:spectral-consequences-and-dependence-on-parameters}

We first extend convergence of nonreal resolvents to the limiting resolvent set and compare spectral projections associated with isolated parts of the spectrum. Identifying the
edge eigenvalues in increasing order requires a further lower spectral bound. The spectral weights of Wigner matrices in the original coordinates also require the discrete boundary
recurrence. For the positive inverse at the hard edge and the self-adjoint inverse of the canonical operator, the global comparison concerns reciprocal eigenvalues. On bounded
intervals of the original spectrum, these estimates transfer to the original eigenvalues; the eigenvector estimates depend on the corresponding spectral gaps.

The results on parameter dependence compare the stochastic Airy operators with different noise amplitudes on one Brownian path and use uniform finite size moments for deterministic
parameter sequences. A supremum of moments over the parameter does not define a pathwise supremum of independently reconstructed matrices.

\subsection{Spectral consequences for Wigner matrices}
\label{sec:combined-w-spectral}
\label{W:sec:comparison-of-the-original-ensembles-and-spectral-consequences}

Convergence in operator norm of the nonreal resolvents extends to the limiting resolvent set and yields convergence of finite spectral projections and compactly supported spectral
statistics. Ordering the eigenvalues from the bottom also requires exclusion of eigenvalues escaping to $-\infty$. Proposition~\ref{W:prop:spectral-lower-bound-in-probability}
derives the required lower bound from localization of the entire spectral subspace below a fixed energy.

To identify the spectral weights at the boundary, we use the Jacobi recurrence to control the derivative of the eigenfunction at zero. Strong local $H^1$ convergence alone does not
determine this derivative. Its distribution is then identified through the Gaussian ensemble.

\subsubsection{Comparison of the original ensembles}
\label{W:sec:comparison-of-the-original-ensembles-and-the-lower-spectral-bound}

We apply the comparison estimates to the original ensembles, including the cutoff factors and fixing the parameters of each function before taking the matrix limit.

Use the test functions of Section~\ref{W:sec:extensions-of-localization-and-modulus-statistics}, with their time intervals, evaluation times, frequencies, modulus thresholds,
energy $E$ and Gram matrix cutoffs fixed. The Gram matrix dimensions lie in the ranges specified in those definitions. The limit $N\to\infty$ is taken before any of these
parameters is varied. Evaluate the functions at the original matrix, including their factors $\Gamma$, and let $H_N^{\mathrm G}$ belong to the GOE or GUE, according to the symmetry
class. Put
\begin{align}
 d_{\Gamma,N}&=\sup_{\|w\|=1}
 \bigl[\mathbb E\Gamma(H_N^{\mathrm G},w)-\mathbb E\Gamma(H_N,w)\bigr]_+,
 \label{W:eq:original-gram-comparison}\\
 e_{{\rm F},N}&=\sup_{\|w\|=1}\left|
 \mathbb E\left[\Gamma e^{\mathrm i\sum_j\lambda_j
       \widetilde Y_N^{\rm reg}(t_j)}\right](H_N,w)
 -\mathbb E\left[\Gamma e^{\mathrm i\sum_j\lambda_j
       \widetilde Y_N^{\rm reg}(t_j)}\right](H_N^{\mathrm G},w)
 \right|,
 \label{W:eq:original-fourier-comparison}\\
 e_{{\rm M},N}^{\sigma}&=\sup_{\|w\|=1}
 \bigl[\mathbb E[\Gamma\mathcal M^\sigma](H_N,w)
       -\mathbb E[\Gamma\mathcal M^\sigma](H_N^{\mathrm G},w)\bigr]_+,
       \qquad\sigma\in\{+,-\},
 \label{W:eq:original-modulus-comparison}\\
 e_{{\rm C},N}&=\sup_{\|w\|=1}
 \bigl[\mathbb E[\Gamma Z_E](H_N,w)
       -\mathbb E[\Gamma Z_E](H_N^{\mathrm G},w)\bigr]_+.
 \label{W:eq:original-localization-comparison}
\end{align}
Every supremum is over deterministic unit vectors. Each function is evaluated with its specified Gram matrix dimension, and $Z_E=C/(1+C)$ on $\{\Gamma>0\}$, where $C$ is the
localization statistic.

\begin{theorem}
\label{W:thm:comparison-original-ensembles}
Under Assumption~\ref{W:ass:wigner-lee-yin},
\[
 d_{\Gamma,N}+e_{{\rm F},N}
       +e_{{\rm M},N}^{+}+e_{{\rm M},N}^{-}+e_{{\rm C},N}
 \longrightarrow0.
\]
\end{theorem}
\begin{proof}
Consider the five bounded functions
\[
 \Gamma,\qquad
 \Gamma\exp\!\left(\mathrm i\sum_j\lambda_j
                      \widetilde Y_N^{\rm reg}(t_j)\right),\qquad
 \Gamma\mathcal M^+,\qquad \Gamma\mathcal M^-,\qquad \Gamma Z_E.
\]
For each function let $F_N(H,w)$ denote its value at the original matrix and let $\widehat F_N(H,w)$ denote its value after the spectral modification of
Section~\ref{W:sec:extensions-of-localization-and-modulus-statistics}. The cutoff factors give $|F_N|,|\widehat F_N|\le1$. Theorem~\ref{W:thm:tightness-main} and
Corollary~\ref{W:cor:tightness-Fourier} apply to the first four functions, with the fixed parameters in the statement. For the last function,
Corollary~\ref{W:cor:localization-comparison} applies with $F_r^*\Theta(H)F_r$ and $\operatorname{Tr}\Theta(H)$ carrying the normalizations used there. Hence, for each of these
five choices,
\[
 \sup_{\|w\|=1}\left|
 \mathbb E\widehat F_N(H_N^{\rm tr},w)
       -\mathbb E\widehat F_N(H_N^{\rm G},w)\right|
 \le N^{-1/6+o(1)}.
\]
A common bound for these five errors is $r_N=N^{-1/6+o(1)}\to0$.

Construct $H_N^{\rm tr}$ from $H_N$ by the entrywise truncation of Section~\ref{W:sec:truncation-centering-and-covariance-normalization}, without recentering or rescaling the
retained entries. If $L_N=e_N\sqrt N$ is the threshold of Lemma~\ref{W:lem:admissible-scale}, the union bound gives
\[
 \begin{aligned}
 \mathbb P(H_N\ne H_N^{\rm tr})
 &\le\frac{N(N-1)}2\mathbb P(|\xi_{12}|>L_N)
          +\sum_{i=1}^N\mathbb P(|d_{N,i}|>L_N)\longrightarrow0.
 \end{aligned}
\]
The two terms tend to zero by \eqref{W:eq:offdiag-equality-rate} and \eqref{W:eq:diagonal-equality-rate}, respectively. In particular, this equality of matrices uses the diagonal
uniform integrability assumption as well as the off-diagonal tail assumption.

Let $h_d$ be the function used for the spectral modification in Lemma~\ref{W:lem:gram-clipping}, at the width $d=A_0\rho_*^{-2}(\log(N+2))^\gamma$, $1<\gamma<2$, used in the
comparison. Its constants and $c_*m_N\le\rho_*\le C_*m_N$ are fixed as above. Define
\[
 \mathcal D_N^H=\{H_N=H_N^{\rm tr},\ h_d(H_N)=H_N\},\qquad
 \mathcal D_N^G=\{h_d(H_N^{\rm G})=H_N^{\rm G}\}.
\]
Since $h_d(x)=x$ for $|x|\le2+d$, Theorem~\ref{W:thm:resolvent-window}, with $b=1/2$, implies, for every fixed $D>0$,
\[
 \begin{aligned}
 p_N^H:=\mathbb P((\mathcal D_N^H)^c)
 &\le\mathbb P(H_N\ne H_N^{\rm tr})
       +\mathbb P(\|H_N^{\rm tr}\|>2+d/2)
 \le\mathbb P(H_N\ne H_N^{\rm tr})+N^{-D}=o(1).
 \end{aligned}
\]
To apply the same theorem to the Gaussian matrix, truncate its standardized entries at $L_N$ and denote the result by $H_N^{\rm G,tr}$. Gaussian tails and $e_N^{-1}=N^{o(1)}$ give
\[
 \begin{aligned}
 \mathbb P(H_N^{\rm G}\ne H_N^{\rm G,tr})
 &\le C N^2 e^{-cL_N^2}=o(N^{-D}),\\
 p_N^G:=\mathbb P((\mathcal D_N^G)^c)
 &\le\mathbb P(H_N^{\rm G}\ne H_N^{\rm G,tr})
       +\mathbb P(\|H_N^{\rm G,tr}\|>2+d/2)=o(1).
 \end{aligned}
\]
Here $c,C$ depend only on the symmetry class, and the theorem applies to the truncated Gaussian array as specified in \eqref{W:eq:resolvent-class}. None of these events depends on
$w$. For every deterministic unit $w$,
\[
 \begin{split}
 \mathbf1_{\mathcal D_N^H}
   [F_N(H_N,w)-\widehat F_N(H_N^{\rm tr},w)]&=0,\\
 \mathbf1_{\mathcal D_N^G}
   [F_N(H_N^{\rm G},w)-\widehat F_N(H_N^{\rm G},w)]&=0.
 \end{split}
\]
On either complementary event the absolute difference is at most two. Integration and the triangle inequality therefore give
\[
 \sup_{\|w\|=1}|\mathbb EF_N(H_N,w)-\mathbb EF_N(H_N^{\rm G},w)|
 \le 2p_N^H+r_N+2p_N^G=:\delta_N\longrightarrow0.
\]
The definitions in \eqref{W:eq:original-gram-comparison}--\eqref{W:eq:original-localization-comparison} use either an absolute value or the positive part of one of these
expectation differences. Since $[x]_+\le|x|$, each of the five nonnegative quantities is at most $\delta_N$. Their sum is at most $5\delta_N$, as required.
\end{proof}

The resulting localization bound and the order of its limits are \eqref{W:eq:localization-probability-bound} and the proof in Section~\ref{W:sec:proof-full-localization}.

\subsubsection{Resolvent matrix elements in the original coordinates}
\label{W:sec:original-coordinate-resolvent-elements}

Operator-norm convergence controls matrix elements against fixed vectors in the common Hilbert space. Expressing those observables in the original matrix coordinates then requires
the Lanczos change of basis as well as the embedding.

\begin{proof}[Proof of Corollary~\ref{W:cor:intro-resolvent-matrix-elements}]
Fix the vectors $f_1,\ldots,f_p,g_1,\ldots,g_q$ in the common Hilbert space. We write the scalar products in the convention
\[
 \Phi(S)(z)=\bigl(\langle S(z)f_i,g_j\rangle\bigr)_{1\le i\le p,\,1\le j\le q}.
\]
Putting the resolvent in the other scalar-product slot gives the same estimate, with the roles of the two vector lists interchanged. For a compact set
$K\subset\mathbb C\setminus\mathbb R$, Cauchy--Schwarz gives
\[
 \begin{split}
 \sup_{z\in K}\|\Phi(S)(z)-\Phi(T)(z)\|_{\mathrm F}^2
 &=\sup_{z\in K}\sum_{i,j}|\langle(S(z)-T(z))f_i,g_j\rangle|^2\\
 &\le\left(\sum_i\|f_i\|^2\right)\left(\sum_j\|g_j\|^2\right)
       \sup_{z\in K}\|S(z)-T(z)\|_{\mathrm{op}}^2.
 \end{split}
\]
Here $\|\cdot\|_{\mathrm F}$ is the Euclidean norm of the matrix entries. In particular, every seminorm of uniform convergence on a compact set for $\Phi(S)-\Phi(T)$ is bounded by
a fixed multiple of the corresponding seminorm for $S-T$ in operator norm. Thus $\Phi$ is continuous for the locally uniform topologies. Applying the continuous mapping theorem to
the resolvent convergence gives the asserted joint convergence of all $pq$ entries, as functions of $z$.

This also proves the version joint with the other coordinates of the original convergence: apply the map which retains those coordinates and applies $\Phi$ to the resolvent
coordinate. On a representation where the resolvent paths converge almost surely locally uniformly, the displayed deterministic inequality proves the same convergence almost
surely, without a separate exceptional event for different values of $z$.

For the original matrix representation, the defining identity for the embedded resolvent gives, for every $z\notin\mathbb R$,
\[
 \langle\mathcal E_N(\mathcal A_N-z)^{-1}\mathcal E_N^*f_i,g_j\rangle
 =\langle(\mathcal A_N-z)^{-1}\mathcal E_N^*f_i,\mathcal E_N^*g_j\rangle.
\]
Conjugating the finite matrix by its Lanczos unitary gives the stated expression in the original matrix coordinates. The fixed vectors are $f_i,g_j$ in the common Hilbert space;
their coordinate representatives after this change of basis need not be deterministic. No claim about arbitrary matrix-dependent vectors follows from the continuity estimate.
\end{proof}

These observables have the following expression in the original coordinates. Direct integration gives
\begin{equation}
 (f_{i,N})_k:=(\mathcal E_N^*f_i)_k
 =h_N^{-1/2}\int_{(k-1)h_N}^{kh_N}f_i(x)\,\dd x.
 \label{W:eq:intro-test-vector-cell-integrals}
\end{equation}
Indeed $\langle\mathcal E_Nv,f_i\rangle =\sum_kv_k\overline{(f_{i,N})_k}$. Put $g_{j,N}=\mathcal E_N^*g_j$, $\widehat f_{i,N}=U_Nf_{i,N}$ and $\widehat g_{j,N}=U_Ng_{j,N}$. The
identity $\mathcal A_N=U_N^*m_N^2(2I-H_N)U_N$ gives
\begin{equation}
 (\mathcal M_N(z))_{ij} =\langle(\mathcal A_N-z)^{-1}f_{i,N},g_{j,N}\rangle =\left\langle(m_N^2(2I-H_N)-z)^{-1}\widehat f_{i,N}, \widehat g_{j,N}\right\rangle.
 \label{W:eq:intro-original-coordinate-resolvent}
\end{equation}
If the functions are supported in $[0,L]$, their coefficients in \eqref{W:eq:intro-test-vector-cell-integrals} vanish for $k>\lceil m_NL\rceil$. Their images under $U_N$
consequently lie in the span of the first $\min\{N,\lceil m_NL\rceil\}$ vectors of the continued Lanczos basis. The cell integrals are deterministic; the images are Borel functions
of $H_N$ by Proposition~\ref{W:prop:lanczos-construction} and the measurable continuation in Section~\ref{W:sec:measurability-and-truncation}. These vectors in the original matrix
coordinates depend on $H_N$ through $U_N$. Projections onto deterministic vectors in the original matrix basis are treated separately in
Section~\ref{W:sec:projections-onto-deterministic-vectors}.

\subsubsection{Continuation of the resolvent}

\label{W:sec:resolvent-continuation-and-spectral-projections}
Convergence at a nonreal parameter is the starting point, but contours around real spectral clusters require the resolvent on a larger domain. The continuation argument below
treats isometrically embedded operators and keeps the action on the orthogonal complement of the finite-dimensional range explicit.

Let $E_n:\mathscr H_n\to\mathscr H$ be isometries and let $A_n$ and $A$ be self-adjoint on the corresponding Hilbert spaces. Nonreal resolvents are taken on the complex Hilbert
spaces, using complexification when the original spaces are real. Write
\begin{equation}
 R_n(z)=E_n(A_n-z)^{-1}E_n^*,\qquad R(z)=(A-z)^{-1}.
 \label{W:eq:trace-general-resolvents}
\end{equation}
Fix $z_0\in\mathbb C\setminus\mathbb R$ and assume
\begin{equation}
 \varepsilon_n:=\|R_n(z_0)-R(z_0)\|\longrightarrow0.
 \label{W:eq:trace-nonreal-resolvent-convergence}
\end{equation}
For $K\Subset\rho(A)$ put
\begin{equation}
 d_K=\sup_{z\in K}|z-z_0|,\qquad
 M_K=\sup_{z\in K}\|(I-(z-z_0)R(z_0))^{-1}\|.
 \label{W:eq:trace-continuation-constants}
\end{equation}
Both constants are finite, by compactness of $K$ and continuity of inversion on the invertible operators.

\begin{proposition}
\label{W:prop:trace-resolvent-continuation}
If $M_Kd_K\varepsilon_n\le1/2$, then
\begin{equation}
 K\subset\rho(A_n),
 \label{W:eq:trace-eventual-resolvent-set}
\end{equation}
and
\begin{equation}
 \sup_{z\in K}\|R_n(z)-R(z)\|\le2M_K^2\varepsilon_n.
 \label{W:eq:trace-continuation-rate}
\end{equation}
In particular,
\begin{equation}
 \sup_{z\in K}\|R_n(z)-R(z)\|\longrightarrow0
 \quad(K\Subset\rho(A)).
 \label{W:eq:trace-full-resolvent-convergence}
\end{equation}
\end{proposition}

\begin{proof}
Let $P_n=E_nE_n^*$ be the orthogonal projection onto $\operatorname{Ran}E_n$, and define
\begin{equation}
 B(z):=I-(z-z_0)R(z_0),
 \qquad
 B_n(z):=I-(z-z_0)R_n(z_0).
 \label{W:eq:trace-Bz}
\end{equation}
The spectral theorem gives
\begin{equation}
 B(z)
 =(A-z)(A-z_0)^{-1}
 \label{W:eq:trace-B-limit-factorization}
\end{equation}
as a bounded operator identity. The factor is invertible exactly when $z\in\rho(A)$, since $(A-z_0)^{-1}$ is a bijection from $\mathscr H$ onto $\mathcal D(A)$. In this case
\begin{equation}
 B(z)^{-1}=(A-z_0)(A-z)^{-1},
 \qquad
 R(z)=R(z_0)B(z)^{-1}.
 \label{W:eq:trace-resolvent-from-z0}
\end{equation}

The orthogonal decomposition \( \mathscr H=\operatorname{Ran}E_n\oplus\operatorname{Ran}E_n^\perp \) gives
\begin{equation}
 B_n(z)
 =
 \left[
 E_n(A_n-z)(A_n-z_0)^{-1}E_n^*
 \right]_{\operatorname{Ran}E_n}
 \oplus I_{\operatorname{Ran}E_n^\perp}.
 \label{W:eq:trace-Bn-block-form}
\end{equation}
The first block is invertible exactly when $z\in\rho(A_n)$; the second block is the identity. For such $z$ the inverse is
\[
 B_n(z)^{-1}=E_n[I+(z-z_0)(A_n-z)^{-1}]E_n^*+(I-P_n).
\]
Multiplication by $R_n(z_0)$, using $E_n^*E_n=I$ and $R_n(z_0)(I-P_n)=0$, gives
\begin{equation}
 R_n(z)=R_n(z_0)B_n(z)^{-1}.
 \label{W:eq:trace-finite-resolvent-from-z0}
\end{equation}

For a nonempty compact $K\Subset\rho(A)$, put $\delta_K=\operatorname{dist}(K,\sigma(A))>0$. The inverse formula gives
\[
 \|B(z)^{-1}\|
 =\sup_{\lambda\in\sigma(A)}\frac{|\lambda-z_0|}{|\lambda-z|}
 \le1+\frac{|z-z_0|}{\delta_K},\qquad z\in K.
\]
Thus $M_K\le1+d_K/\delta_K<\infty$. There is nothing to check when $K$ is empty. For \(z\in K\), write
\begin{equation}
 B_n(z)
 =B(z)
 \left[
 I+B(z)^{-1}(B_n(z)-B(z))
 \right].
 \label{W:eq:trace-inversion-factorization}
\end{equation}
Since
\begin{equation}
 \sup_{z\in K}\|B_n(z)-B(z)\|
 \le
 \sup_{z\in K}|z-z_0|\,\varepsilon_n,
 \label{W:eq:trace-B-difference}
\end{equation}
the smallness assumption gives
\begin{equation}
 \sup_{z\in K}
 \|B(z)^{-1}(B_n(z)-B(z))\|
 \le\frac12.
 \label{W:eq:trace-neumann-smallness}
\end{equation}
The Neumann series in \eqref{W:eq:trace-inversion-factorization} therefore yields
\begin{equation}
 \sup_{z\in K}\|B_n(z)^{-1}\|
 \le2M_K.
 \label{W:eq:trace-Bn-inverse-bound}
\end{equation}
The resolvent identity for the inverses gives
\begin{equation}
 B_n(z)^{-1}-B(z)^{-1}
 =B_n(z)^{-1}(B(z)-B_n(z))B(z)^{-1},
 \label{W:eq:trace-inverse-difference}
\end{equation}
so
\begin{equation}
 \sup_{z\in K}\|B_n(z)^{-1}-B(z)^{-1}\|
 \le
 2M_K^2\sup_{z\in K}|z-z_0|\,\varepsilon_n.
 \label{W:eq:trace-inverse-rate}
\end{equation}
Finally, subtract \eqref{W:eq:trace-resolvent-from-z0} from \eqref{W:eq:trace-finite-resolvent-from-z0}:
\begin{equation}
 R_n(z)-R(z) ={}[R_n(z_0)-R(z_0)]B_n(z)^{-1} +R(z_0)[B_n(z)^{-1}-B(z)^{-1}].
 \label{W:eq:trace-resolvent-continuation-difference}
\end{equation}
To obtain the stated constant, multiply the difference by $B_n(z)$ on the left and $B(z)$ on the right. Since each resolvent at $z_0$ commutes with its own factor $B_n(z)$ or
$B(z)$, we obtain, with $\zeta=z-z_0$,
\[
 \begin{split}
 B_n(z)(R_n(z)-R(z))B(z)
 &=R_n(z_0)B(z)-B_n(z)R(z_0)\\
 &=R_n(z_0)-\zeta R_n(z_0)R(z_0)
       -R(z_0)+\zeta R_n(z_0)R(z_0)\\
 &=R_n(z_0)-R(z_0).
 \end{split}
\]
Thus
\[
 R_n(z)-R(z)=B_n(z)^{-1}(R_n(z_0)-R(z_0))B(z)^{-1}.
\]
The bounds $\|B_n(z)^{-1}\|\le2M_K$ and $\|B(z)^{-1}\|\le M_K$ prove \eqref{W:eq:trace-continuation-rate}. Since $\varepsilon_n\to0$, there is an index $n_K$ such that
$M_Kd_K\varepsilon_n\le1/2$ for every $n\ge n_K$. The conclusion therefore holds for every such $n$.
\end{proof}

The identity $E_n^*E_n=I$ preserves products of operators on $\mathscr H_n$ under the embedding. We shall also use a compactness criterion based on finite-dimensional compressions.

Let $\mathcal C\Subset\mathbb C\setminus\mathbb R$, and let $d=\inf_{z\in\mathcal C}|\Im z|>0$. Then
\begin{equation}
 R_n(z)^*=R_n(\bar z),\qquad
 R_n(z)-R_n(w)
       =(z-w)R_n(z)R_n(w),\qquad
 \|R_n(z)\|\le |\Im z|^{-1}.
 \label{W:eq:spectral-resolvent-algebra}
\end{equation}
In particular their restrictions to $\mathcal C$ have common Lipschitz constant $d^{-2}$. For a finite-rank orthogonal projection $F$, put
\[
 a_F(R;\mathcal C)=\sup_{z\in\mathcal C}
       \bigl(\|(I-F)R(z)\|+\|(I-F)R(z)^*\|\bigr).
\]

\begin{lemma}
\label{W:lem:compact-families-of-resolvents}
A uniformly bounded, equicontinuous family in $C(\mathcal C;\mathcal K(\mathscr H))$ is relatively compact if for every $\epsilon>0$ there is a finite-rank orthogonal projection
$F$ such that $\sup_R a_F(R;\mathcal C)<\epsilon$. Moreover,
\begin{equation}
 \sup_{\mathcal C}\|R_n-R\|
 \le a_F(R_n;\mathcal C)+a_F(R;\mathcal C)
             +\sup_{\mathcal C}\|F(R_n-R)F\|.
 \label{W:eq:spectral-compression-error}
\end{equation}
\end{lemma}
\begin{proof}
Put $G_n(z)=(A_n-z)^{-1}$. The inverse identity gives $G_n(z)-G_n(w)=(z-w)G_n(z)G_n(w)$. Thus
\[
 R_n(z)R_n(w)
 =E_nG_n(z)(E_n^*E_n)G_n(w)E_n^*
 =E_nG_n(z)G_n(w)E_n^*.
\]
Taking adjoints gives the first identity in \eqref{W:eq:spectral-resolvent-algebra}, while the spectral theorem gives $\|R_n(z)\|\le|\operatorname{Im}z|^{-1}$. Consequently
$\|R_n(z)-R_n(w)\|\le d^{-2}|z-w|$ for $z,w\in\mathcal C$. For any orthogonal projection $F$,
\[
 R-FRF=(I-F)R+FR(I-F),\qquad
 \|FR(I-F)\|\le\|(I-F)R^*\|.
\]
Apply this to $R_n$ and $R$ and subtract their two compressions to obtain \eqref{W:eq:spectral-compression-error}.

Fix $\epsilon>0$ and choose $F$ with $\sup_R a_F(R;\mathcal C)<\epsilon/6$. By compactness of $\mathcal C$, equicontinuity gives $\delta>0$ such that $\|R(z)-R(w)\|<\epsilon/6$ for
every path in the family whenever $|z-w|<\delta$. Choose a finite $\delta$-net $z_1,\ldots,z_l$ in $\mathcal C$. The compressed values $FR(z_i)F$ lie in a bounded subset of the
finite-dimensional space $\mathcal B(\operatorname{ran}F)$, which has a finite $\epsilon/6$-net. For each tuple of net points attained as approximations to $(FR(z_i)F)_{i=1}^l$,
retain one path $T_j$ from the family. There are finitely many such paths. For each $R$ choose $T_j$ with the same tuple. Then
\[
 \max_{i\le l}\|F(R(z_i)-T_j(z_i))F\|<\epsilon/3.
\]
For $z\in\mathcal C$, choose $i$ with $|z-z_i|<\delta$. The two equicontinuity bounds give
\[
 \|F(R(z)-T_j(z))F\|
 \le\|R(z)-R(z_i)\|+\|F(R(z_i)-T_j(z_i))F\|
                         +\|T_j(z_i)-T_j(z)\|<2\epsilon/3.
\]
Equation~\eqref{W:eq:spectral-compression-error} now gives $\sup_{\mathcal C}\|R-T_j\|<\epsilon$. Thus the family is totally bounded in the Banach space
$C(\mathcal C;\mathcal K(\mathscr H))$, and its closure is compact.

For applications to self-adjoint resolvents, the second term in $a_F$ is estimated on the conjugate contour. Indeed, \eqref{W:eq:spectral-resolvent-algebra} gives
\[
 \sup_{z\in\mathcal C}\|(I-F)R_n(z)^*\|
 =\sup_{z\in\overline{\mathcal C}}\|(I-F)R_n(z)\|.
\]
Applying the localization bound on $\mathcal C\cup\overline{\mathcal C}$ therefore bounds both terms defining $a_F$.
\end{proof}

\subsubsection{Spectral projections and isolated eigenvalues}
\label{W:sec:spectral-projections-and-isolated-eigenvalues}

Uniform resolvent convergence on a contour gives convergence of the enclosed spectral projection. Finite rank then makes the enclosed multiplicity stable. With the convention
$R(z)=(A-z)^{-1}$, the positively oriented Riesz integral has a minus sign.

\Needspace{6\baselineskip}

Under the assumptions of Proposition~\ref{W:prop:trace-resolvent-continuation}, let \(\Gamma\) be a finite system of pairwise disjoint positively oriented rectifiable Jordan curves
contained in \(\rho(A)\), whose bounded interiors are also pairwise disjoint. Let $\Delta$ be the union of those interiors, so the winding number of the contour system is one on
$\Delta$ and zero outside its closure. Here $\mathbf1_\Delta(A)$ means the spectral projection for $\Delta\cap\mathbb R$. Assume $\operatorname{rank}\mathbf1_\Delta(A)<\infty$.
Define $P(\Gamma)$ and, for those $n$ such that $\Gamma\subset\rho(A_n)$, $P_n(\Gamma)$ by
\begin{align}
 P(\Gamma)
 &:=-\frac1{2\pi\mathrm i}\oint_\Gamma R(z)\,\dd z,
 \label{W:eq:trace-limit-riesz}\\
 P_n(\Gamma)
 &:=-\frac1{2\pi\mathrm i}\oint_\Gamma R_n(z)\,\dd z.
 \label{W:eq:trace-finite-riesz}
\end{align}

\begin{proposition}
\label{W:prop:trace-riesz-stability}
There is an index $n_\Gamma$ such that, for every $n\ge n_\Gamma$, $\Gamma\subset\rho(A_n)$ and
\begin{equation}
 P_n(\Gamma)=E_n\mathbf1_\Delta(A_n)E_n^*,
 \qquad
 P(\Gamma)=\mathbf1_\Delta(A).
 \label{W:eq:trace-riesz-spectral-projection}
\end{equation}
Moreover,
\begin{align}
 \|P_n(\Gamma)-P(\Gamma)\|&\longrightarrow0,
 \label{W:eq:trace-riesz-convergence}\\
 \operatorname{rank}P_n(\Gamma)&=\operatorname{rank}P(\Gamma)
                   \quad(n\ge n_\Gamma),
 \label{W:eq:trace-rank-stability}\\
 \|P_n(\Gamma)-P(\Gamma)\|_1&\longrightarrow0.
 \label{W:eq:trace-riesz-trace-convergence}
\end{align}
\end{proposition}

\begin{proof}
Set $\varepsilon_n=\|R_n(z_0)-R(z_0)\|$ as in \eqref{W:eq:trace-nonreal-resolvent-convergence}. The compact contour $\Gamma\subset\rho(A)$ has finite length $\ell(\Gamma)$, and
Proposition~\ref{W:prop:trace-resolvent-continuation} gives, once
$M_\Gamma d_\Gamma\varepsilon_n\le1/2$,
\[
 \Gamma\subset\rho(A_n),\qquad
 \|P_n(\Gamma)-P(\Gamma)\|
 \le\frac{\ell(\Gamma)}{2\pi}
     \sup_{z\in\Gamma}\|R_n(z)-R(z)\|
 \le\frac{\ell(\Gamma)}{\pi}M_\Gamma^2\varepsilon_n.
\]
For the fixed operators and contour, $M_\Gamma,d_\Gamma$ and $\ell(\Gamma)$ do not depend on $n$. Choose $n_\Gamma$ so that both $M_\Gamma d_\Gamma\varepsilon_n\le1/2$ and
$\ell(\Gamma)M_\Gamma^2\varepsilon_n/\pi<1$ for every $n\ge n_\Gamma$. For real $x$ not on the contour, the scalar Cauchy formula is
\[
 -\frac1{2\pi\mathrm i}\oint_\Gamma\frac{\mathrm dz}{x-z}
 =\frac1{2\pi\mathrm i}\oint_\Gamma\frac{\mathrm dz}{z-x}
 =\mathbf1_\Delta(x).
\]
For $f,g\in\mathscr H$, the spectral measure $\mu_{f,g}$ of $A$ satisfies $|\mu_{f,g}|(\mathbb R)\le\|f\|\|g\|$. Writing $\delta=\operatorname{dist}(\Gamma,\sigma(A))>0$, we have
\[
 \int_{\mathbb R}\oint_\Gamma
     \frac{|\dd z|}{|\lambda-z|}\,\dd|\mu_{f,g}|(\lambda)
 \le\frac{\ell(\Gamma)}{\delta}\|f\|\|g\|.
\]
Fubini's theorem and the scalar Cauchy formula therefore give $P(\Gamma)=\mathbf1_\Delta(A)$. For $n$ sufficiently large the same calculation applies to $A_n$, since
$\Gamma\subset\rho(A_n)$, with its positive spectral measures on $\mathscr H_n$. Applying the isometry $E_n$ yields exactly \eqref{W:eq:trace-riesz-spectral-projection}.

Let \( r:=\operatorname{rank}P(\Gamma) \). For $\|P-Q\|<1$, the restrictions $Q|_{\operatorname{Ran}P}$ and $P|_{\operatorname{Ran}Q}$ are injective, since $Qu=0$ and $Pu=u$ imply
$\|u\|\le\|P-Q\|\|u\|$. Thus their ranks agree. Apply this with $P=P(\Gamma)$ and $Q=P_n(\Gamma)$ to obtain \eqref{W:eq:trace-rank-stability}. Finally,
\begin{equation}
 \operatorname{rank}(P_n(\Gamma)-P(\Gamma))\le2r,
 \label{W:eq:trace-riesz-difference-rank}
\end{equation}
so
\begin{equation}
 \|P_n(\Gamma)-P(\Gamma)\|_1
 \le2r\|P_n(\Gamma)-P(\Gamma)\|\longrightarrow0.
 \label{W:eq:trace-riesz-trace-bound}
\end{equation}
\end{proof}

\Needspace{6\baselineskip}

Under the assumptions of Proposition~\ref{W:prop:trace-resolvent-continuation}, let $a<b$ be real numbers in $\rho(A)$, put $I=(a,b)$, and suppose
\begin{equation}
 \sigma(A)\cap I=\{\Lambda\},
 \qquad
 r:=\operatorname{mult}_A(\Lambda)<\infty.
 \label{W:eq:trace-one-cluster-hypothesis}
\end{equation}

\begin{lemma}
\label{W:lem:trace-isolated-cluster}
For all sufficiently large \(n\), the interval \(I\) contains exactly \(r\) eigenvalues of \(A_n\), counted with multiplicity. If they are
\(\lambda_{1,n}^{I},\ldots,\lambda_{r,n}^{I}\), then
\begin{equation}
 \max_{1\le q\le r}|\lambda_{q,n}^{I}-\Lambda|
 \longrightarrow0.
 \label{W:eq:trace-cluster-eigenvalue-convergence}
\end{equation}
The embedded cluster projection converges in operator and trace norm to \(\mathbf1_{\{\Lambda\}}(A)\). In particular, if
\(r=1\), the unique eigenvalue in \(I\) is simple for all sufficiently large \(n\), and its embedded eigenspace
projection converges in operator norm to the rank-one limit projection.
\end{lemma}

\begin{proof}
Fix $\eta>0$ and take as $\Gamma$ the counterclockwise boundary of
$\{z:a<\operatorname{Re}z<b,\ |\operatorname{Im}z|<\eta\}$.
On its horizontal sides $\operatorname{Im}z=\pm\eta$, while the vertical sides meet the real axis only at the two points $a,b\in\rho(A)$. Since $A$ is self-adjoint,
$\Gamma\subset\rho(A)$. The scalar Cauchy formula identifies its Riesz projection with
\[
 P=\mathbf1_I(A)=\mathbf1_{\{\Lambda\}}(A),
 \qquad \operatorname{rank}P=r.
\]
Let $\ell_\Gamma$ be the length of $\Gamma$. With the constants $M_\Gamma,d_\Gamma$ in \eqref{W:eq:trace-continuation-constants} and $\varepsilon_n=\|R_n(z_0)-R(z_0)\|$, for all
large $n$ we have
\[
 \begin{gathered}
 P_n=E_n\mathbf1_I(A_n)E_n^*,\qquad
 \|P_n-P\|\le \frac{\ell_\Gamma}{\pi}M_\Gamma^2\varepsilon_n,\\
 \operatorname{rank}P_n=r,\qquad
 \|P_n-P\|_1\le 2r\|P_n-P\|.
 \end{gathered}
\]
The rank identity holds once $\|P_n-P\|<1$; it follows from the two injections in Proposition~\ref{W:prop:trace-riesz-stability}, and applies for arbitrary finite $r$.

Put $Q_n=\mathbf1_I(A_n)$ and $\mathscr H_n^I=\operatorname{Ran}Q_n$. The spectral theorem gives the reducing decomposition
\[
 A_n=A_n\big|_{\mathscr H_n^I}\oplus
     A_n\big|_{(\mathscr H_n^I)^\perp},\qquad
 \operatorname{dim}\mathscr H_n^I=\operatorname{rank}Q_n=r,
\]
with $\sigma(A_n|_{(\mathscr H_n^I)^\perp})\subset \mathbb R\setminus I$. On $\mathscr H_n^I$ the spectral support is contained in the bounded interval $I$, so the restriction is a
bounded self-adjoint operator on the $r$-dimensional space, hence has exactly $r$ eigenvalues with multiplicity. For each $x\in I$,
\[
 \left\|\left(A_n\big|_{(\mathscr H_n^I)^\perp}-x\right)^{-1}\right\|
 \le \frac1{\min\{x-a,b-x\}},
\]
which shows that the complementary restriction has no spectrum in $I$. In particular there are no additional eigenvalues, nor continuous spectrum, of $A_n$ in $I$.

Fix $0<\delta<\min\{\Lambda-a,b-\Lambda\}$. Set
\begin{equation}
 K_\delta
 :=[a,b]\setminus(\Lambda-\delta,\Lambda+\delta)
 \label{W:eq:trace-cluster-complement}
\end{equation}
and note that $K_\delta\Subset\rho(A)$ by \eqref{W:eq:trace-one-cluster-hypothesis} and $a,b\in\rho(A)$. Let $M_\delta,d_\delta$ be its fixed constants from
\eqref{W:eq:trace-continuation-constants}. When
$M_\delta d_\delta\varepsilon_n\le1/2$,
Proposition~\ref{W:prop:trace-resolvent-continuation} gives
\[
 K_\delta\subset\rho(A_n),\qquad
 \max_{1\le q\le r}|\lambda_{q,n}^{I}-\Lambda|<\delta.
\]
Choose $n_\delta\ge n_\Gamma$ so that $M_\delta d_\delta\varepsilon_n\le1/2$ for every $n\ge n_\delta$. The inclusion excludes spectrum from $K_\delta$, while the rank equality
gives exactly $r$ eigenvalues in $I$. Consequently, for every $n\ge n_\delta$, all of them lie in $(\Lambda-\delta,\Lambda+\delta)$. Hence
\[
 0\le\limsup_{n\to\infty}\max_{1\le q\le r}
 |\lambda_{q,n}^{I}-\Lambda|\le\delta.
\]
Letting $\delta\downarrow0$ proves \eqref{W:eq:trace-cluster-eigenvalue-convergence}; the operator- and trace-norm convergence of the projection follows from the bounds above. When
$r=1$, the finite spectral subspace is one-dimensional and therefore the eigenvalue is simple.

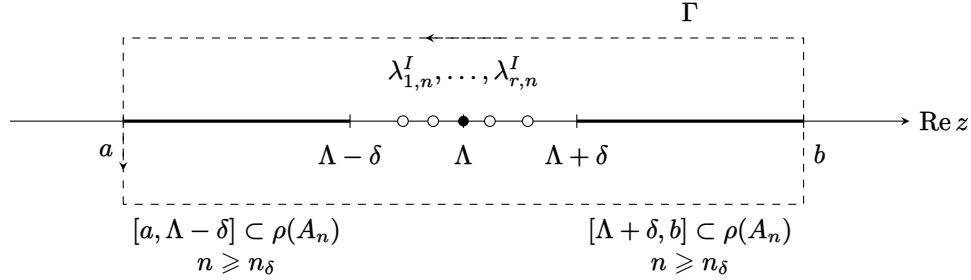
\begin{figure}[H]
\centering
\begin{tikzpicture}[x=1cm,y=1cm,>=stealth,every node/.style={font=\small}]
 \draw[->] (-.5,0)--(11.4,0) node[right] {$\operatorname{Re}z$};
 \draw[dashed] (1,-1.1) rectangle (10,1.1);
 \draw[->,dashed] (6,1.1)--(5,1.1);
 \draw[->,dashed] (1,-.2)--(1,-.65);
 \node at (8.5,1.45) {$\Gamma$};
 \draw[very thick] (1,0)--(4,0);
 \draw[very thick] (7,0)--(10,0);
 \foreach \x/\t in {1/a,4/{\Lambda-\delta},5.5/\Lambda,7/{\Lambda+\delta},10/b}{
  \draw (\x,-.09)--(\x,.09);
  \ifdim\x pt=1pt\node[below left] at (\x,-.18) {$\t$};\else
   \ifdim\x pt=10pt\node[below right] at (\x,-.18) {$\t$};\else
   \node[below] at (\x,-.18) {$\t$};\fi\fi
 }
 \fill (5.5,0) circle (2pt);
 \foreach \x in {4.7,5.1,5.85,6.35}{\draw[fill=white] (\x,0) circle (2pt);}
 \node at (5.5,.63) {$\lambda_{1,n}^{I},\ldots,\lambda_{r,n}^{I}$};
 \node[align=center] at (2.5,-1.65) {$[a,\Lambda-\delta]\subset\rho(A_n)$\\$n\ge n_\delta$};
 \node[align=center] at (8.5,-1.65) {$[\Lambda+\delta,b]\subset\rho(A_n)$\\$n\ge n_\delta$};
\end{tikzpicture}
\caption{An isolated eigenvalue $\Lambda$ of multiplicity $r$ and its projection contour $\Gamma$. The counterclockwise orientation gives $-(2\pi\mathrm i)^{-1}\oint_\Gamma(A-z)^{-1}\,\dd z$ in \eqref{W:eq:trace-limit-riesz}. Rank stability gives exactly $r$ eigenvalues of $A_n$ in $I$, counted with multiplicity. The thick intervals constitute $K_\delta$ in \eqref{W:eq:trace-cluster-complement} and lie in $\rho(A_n)$ for $n\ge n_\delta$, with $n_\delta$ chosen in the proof. The open points schematically illustrate the eigenvalues converging to $\Lambda$.}
\label{W:fig:spectral-cluster-contour}
\end{figure}
\end{proof}

\Needspace{6\baselineskip}

\begin{corollary}
\label{W:cor:trace-spectral-exactness}
Under the assumptions of Proposition~\ref{W:prop:trace-resolvent-continuation}, let $I=(a,b)$ be a bounded real interval with $a,b\in\rho(A)$. If $A$ has compact resolvent, then
\begin{equation}
 \operatorname{rank}\mathbf1_I(A_n)
 =\operatorname{rank}\mathbf1_I(A)
 \label{W:eq:trace-interval-rank}
\end{equation}
for all sufficiently large \(n\). For every compact $K\subset\rho(A)\cap\mathbb R$, there is an index $n_K$ such that $K\cap\sigma(A_n)=\varnothing$ for all $n\ge n_K$, and each
isolated eigenvalue of \(A\) attracts exactly its multiplicity of eigenvalues of \(A_n\).
\end{corollary}

\begin{proof}
The compact normal operator $(A-z_0)^{-1}$ has nonzero spectral values of finite multiplicity, with possible accumulation only at zero. A real eigenvalue $\lambda\in[a,b]$ of $A$
corresponds to the nonzero eigenvalue $(\lambda-z_0)^{-1}$ of the resolvent, and
\[
 |(\lambda-z_0)^{-1}|
 \ge \left(\max_{x\in[a,b]}|x-z_0|\right)^{-1}>0.
\]
Thus $\sigma(A)\cap[a,b]$ is finite, with finite multiplicities; this is also the compact-resolvent conclusion of Theorem~III.6.29 of \cite{kato1995perturbation}. In particular,
\[
 r:=\operatorname{rank}\mathbf1_I(A)
  =\sum_{\Lambda\in\sigma(A)\cap I}
                \operatorname{mult}_A(\Lambda)<\infty.
\]
Choose a counterclockwise rectangle $\Gamma$ with vertical sides through $a,b$ and horizontal sides at imaginary parts $\pm\eta$ for fixed $\eta>0$. Its contour is contained in
$\rho(A)$: real points on it occur only at $a$ and $b$. Proposition~\ref{W:prop:trace-riesz-stability} gives, for large $n$,
\[
 P_n=E_n\mathbf1_I(A_n)E_n^*,\quad
 P=\mathbf1_I(A),\quad
 \|P_n-P\|<1,\quad
 \operatorname{rank}P_n=r.
\]
Since $E_n$ is isometric, $\operatorname{rank}P_n$ equals the rank of $\mathbf1_I(A_n)$, proving \eqref{W:eq:trace-interval-rank}, including the case $r=0$.

For $K\Subset\rho(A)\cap\mathbb R$, define $M_K,d_K$ as in \eqref{W:eq:trace-continuation-constants}. Choose $n_K$ such that $M_Kd_K\varepsilon_n\le1/2$ for all $n\ge n_K$. Then
Proposition~\ref{W:prop:trace-resolvent-continuation} gives $K\subset\rho(A_n)$, so $\sigma(A_n)\cap K=\varnothing$.

Finally, fix an isolated eigenvalue $\lambda$ of $A$, of multiplicity $m$. Choose $\delta>0$ such that
$[\lambda-\delta,\lambda+\delta]\cap\sigma(A)=\{\lambda\}$
and $\lambda-\delta,\lambda+\delta\in\rho(A)$. Applying Lemma~\ref{W:lem:trace-isolated-cluster} to $I_\lambda=(\lambda-\delta,\lambda+\delta)$ gives
\[
 \operatorname{rank}\mathbf1_{I_\lambda}(A_n)=m,
 \qquad
 \max_{1\le j\le m}|\lambda_{j,n}^{I_\lambda}-\lambda|
 \longrightarrow0.
\]
For an interval $I$ containing several limiting eigenvalues, choose disjoint such $I_\lambda$ inside $I$, and set
$K=[a,b]\setminus\bigcup_\lambda I_\lambda$.
Then $K\Subset\rho(A)$. Take the maximum of the index $n_K$ just obtained and the finitely many indices for the intervals $I_\lambda$. Above this maximum, $K\subset\rho(A_n)$, each
$I_\lambda$ contains exactly the multiplicity of $\lambda$, and these intervals exhaust the spectrum of $A_n$ in $I$. In particular, no extra eigenvalues occur in the gaps between
these clusters.
\end{proof}

\subsubsection{Eigenvalue point measures and eigenfunctions up to phase}
\label{W:sec:spectral-measures-and-eigenfunctions}

Finite spectral projections describe isolated clusters. The vague topology combines these local conclusions into convergence of the eigenvalue point measures without imposing
control at infinity; eigenfunctions are compared through their one-dimensional projections before choosing phases.

Let $\mathcal M_{\mathrm{lf}}(\mathbb R)$ be the space of locally finite positive Radon measures with the vague topology, defined by convergence of integrals against every
$f\in C_c(\mathbb R)$. This space is Polish by Lemma~A5.5 of \cite{kallenberg2021foundations}, applied to $\mathbb R$ with its usual complete metric. In this metric bounded sets
have compact closure, so the locally finite measures and test functions in that lemma give precisely the space and topology used here. Its Borel sigma-field is generated by the
maps $\mu\mapsto\int f\,\dd\mu$ with
$f\in C_c(\mathbb R)$.

Under the deterministic assumptions of Proposition~\ref{W:prop:trace-resolvent-continuation}, suppose that $A$ and every $A_n$ have compact resolvent, as is automatic for the
finite-dimensional approximants. Their locally finite counting measures are
\begin{equation}
 \Xi_n^A:=\sum_{\lambda\in\sigma(A_n)}
       \operatorname{mult}_{A_n}(\lambda)\delta_\lambda,\qquad
 \Xi^A:=\sum_{\Lambda\in\sigma(A)}
       \operatorname{mult}_{A}(\Lambda)\delta_\Lambda.
 \label{W:eq:trace-deterministic-point-measures}
\end{equation}
They satisfy
\begin{equation}
 \Xi_n^A\longrightarrow\Xi^A\quad\text{vaguely on }\mathbb R.
 \label{W:eq:trace-deterministic-vague}
\end{equation}
Indeed, choose $a<b$ in $\rho(A)\cap\mathbb R$ with $\operatorname{supp}f\Subset(a,b)$ for $f\in C_c(\mathbb R)$. Compactness of the resolvent gives finitely many eigenvalues
$\Lambda^{(q)}$ in $(a,b)$, with multiplicities $m_q$. Choose disjoint open isolating intervals $I_q\Subset(a,b)$ around them, with $\partial I_q\subset\rho(A)$. Then
\[
 K=[a,b]\setminus\bigcup_q I_q\Subset\rho(A).
\]
Proposition~\ref{W:prop:trace-resolvent-continuation} supplies an index $n_K$ with $K\subset\rho(A_n)$ for all $n\ge n_K$. After increasing it to cover the finitely many intervals
$I_q$, Lemma~\ref{W:lem:trace-isolated-cluster} gives exactly $m_q$ eigenvalues $\lambda_{q,j}^{(n)}$ in each $I_q$, all converging to $\Lambda^{(q)}$. Eigenvalues outside $(a,b)$
contribute zero because $f$ vanishes there. Thus, including the case of an empty sum,
\[
 \langle\Xi_n^A,f\rangle
 =\sum_q\sum_{j=1}^{m_q}f(\lambda_{q,j}^{(n)})
 \longrightarrow\sum_qm_qf(\Lambda^{(q)})=\langle\Xi^A,f\rangle,
\]
which proves \eqref{W:eq:trace-deterministic-vague}.

Define
\begin{equation}
 \Xi_N:=\sum_{j=1}^N\delta_{\lambda_{j,N}}
 =\sum_{j=1}^N
  \delta_{N^{2/3}(2-\mu_j^{(N)})},
 \label{W:eq:trace-finite-point-process}
\end{equation}
and
\begin{equation}
 \Xi_\beta:=\sum_{j=1}^\infty\delta_{\Lambda_j}.
 \label{W:eq:trace-limit-point-process}
\end{equation}
The second measure is locally finite almost surely, since $\mathcal H_\beta$ is lower semibounded with compact resolvent.

We use the quotient of the Hilbert space by scalar phases. Let \(\mathbb S_\mathbb F=\{-1,1\}\) in the real case and
\(\mathbb S_\mathbb F=\{z\in\mathbb C:|z|=1\}\) in the complex case. For a separable Hilbert space \(\mathscr H\),
define
\begin{equation}
 \mathscr H_{\mathrm{ph}}:=\mathscr H/\mathbb S_\mathbb F,
 \qquad
 d_{\mathrm{ph}}([u],[v])
 :=\inf_{\omega\in\mathbb S_\mathbb F}\|u-\omega v\|.
 \label{W:eq:trace-phase-quotient-metric}
\end{equation}
The zero vector is included.

For unit vectors $u_n,u$, convergence of their rank-one projections,
\begin{equation}
 \bigl\||u_n\rangle\langle u_n|-|u\rangle\langle u|\bigr\|
 \longrightarrow0,
 \label{W:eq:trace-rank-one-projection-convergence}
\end{equation}
implies convergence up to a phase. Indeed, set $c_n=|\langle u_n,u\rangle|$ and choose $\omega_n=\overline{\langle u_n,u\rangle}/c_n$ when $c_n>0$, and $\omega_n=1$ when $c_n=0$.
Then $\langle\omega_nu_n,u\rangle=c_n$. If $c_n=1$, the two projections agree. If $c_n<1$, write $\omega_nu_n=c_nu+(1-c_n^2)^{1/2}e_n$, where $e_n\perp u$ is a unit vector. In the
basis $(u,e_n)$ the projection difference is
\[
 \begin{pmatrix}
 c_n^2-1&c_n\sqrt{1-c_n^2}\\
 c_n\sqrt{1-c_n^2}&1-c_n^2
 \end{pmatrix}.
\]
Its square is $(1-c_n^2)I_2$, and the projection difference vanishes on the orthogonal complement. In both cases,
\begin{equation}
 \bigl\||u_n\rangle\langle u_n|-|u\rangle\langle u|\bigr\|^2=1-c_n^2.
 \label{W:eq:trace-rank-one-distance}
\end{equation}
Thus $c_n\to1$, and the chosen phases satisfy
\begin{equation}
 \|\omega_nu_n-u\|^2=2-2c_n\longrightarrow0,
 \label{W:eq:trace-phase-alignment-calculation}
\end{equation}
and hence
\begin{equation}
 \|\omega_nu_n-u\|\longrightarrow0.
 \label{W:eq:trace-rank-one-vector-convergence}
\end{equation}

\Needspace{6\baselineskip}

\begin{lemma}
\label{W:lem:trace-phase-quotient}
The metric space $(\mathscr H_{\mathrm{ph}},d_{\mathrm{ph}})$ is Polish. On the quotient of the unit sphere, the map
\begin{equation}
 [u]\longmapsto P_u:=|u\rangle\langle u|
 \label{W:eq:trace-projective-map}
\end{equation}
is a homeomorphism onto the closed set of rank-one orthogonal projections in \(\mathcal K(\mathscr H)\). If \(u,v\) are unit vectors and \(c=|\langle u,v\rangle|\), then
\begin{align}
 d_{\mathrm{ph}}([u],[v])^2&=2(1-c),
 \label{W:eq:trace-phase-distance-formula}\\
 \|P_u-P_v\|^2&=1-c^2,
 \label{W:eq:trace-projection-distance-formula}
\end{align}
and hence
\begin{equation}
 \|P_u-P_v\|
 \le d_{\mathrm{ph}}([u],[v])
 \le\sqrt2\|P_u-P_v\|.
 \label{W:eq:trace-projective-equivalence}
\end{equation}
\end{lemma}

\begin{proof}
For phases $\alpha,\beta\in\mathbb S_{\mathbb F}$,
\[
 \inf_\omega\|\alpha u-\omega\beta v\|
 =\inf_\omega\|u-\overline\alpha\omega\beta v\|
 =\inf_\eta\|u-\eta v\|.
\]
Thus the distance is independent of the chosen representatives. The map $\omega\mapsto\|u-\omega v\|$ is continuous on the compact phase group, so its minimum is attained. That
minimum is zero exactly when $u=\omega v$ for one phase, which is equality of the two orbits.

For phases $\alpha,\beta$,
\[
 \|u-\alpha\beta w\|\le\|u-\alpha v\|+\|v-\beta w\|.
\]
Taking minimizing phases proves the triangle inequality. Symmetry follows by replacing a phase by its inverse. The quotient is separable because the image of a countable dense
subset of $\mathscr H$ is dense. To prove completeness, let
\([u_n]\) be Cauchy. Choose a subsequence \([u_{n_k}]\) satisfying
\begin{equation}
 d_{\mathrm{ph}}([u_{n_{k+1}}],[u_{n_k}])<2^{-k}.
 \label{W:eq:trace-phase-Cauchy-subsequence}
\end{equation}
Choose phases recursively so that representatives \(v_k\in[u_{n_k}]\) satisfy \(\|v_{k+1}-v_k\|<2^{-k}\). For $l>k$,
\[
 \|v_l-v_k\|\le\sum_{a=k}^{l-1}2^{-a}\le2^{1-k}.
\]
Completeness of $\mathscr H$ supplies $v_k\to v$, and $d_{\rm ph}([u_{n_k}],[v])\le\|v_k-v\|\to0$. Given $\epsilon>0$, choose $k$ so that this last distance is below $\epsilon/2$
and the Cauchy distances between any two terms with indices at least $n_k$ are below $\epsilon/2$. For $n\ge n_k$,
\[
 d_{\rm ph}([u_n],[v])
 \le d_{\rm ph}([u_n],[u_{n_k}])+d_{\rm ph}([u_{n_k}],[v])<\epsilon.
\]
Thus the quotient is complete; together with the countable dense set already constructed, this proves it is Polish.

For unit vectors, the phase choice above gives
\begin{equation}
 d_{\mathrm{ph}}([u],[v])^2
 =\|u-\omega v\|^2=2-2|\langle u,v\rangle|.
 \label{W:eq:trace-phase-distance-proof}
\end{equation}
Equation~\eqref{W:eq:trace-rank-one-distance} gives \eqref{W:eq:trace-projection-distance-formula}; $1-c\le1-c^2\le2(1-c)$ then proves \eqref{W:eq:trace-projective-equivalence}.

Every rank-one orthogonal projection $P$ equals $P_u$ for a unit vector $u$ in its range: self-adjointness gives $Pf=\langle f,u\rangle u$. Moreover $P_u=P_v$ implies $v=\omega u$
for a scalar $|\omega|=1$. Thus the map is a bijection onto all rank-one orthogonal projections, and \eqref{W:eq:trace-projective-equivalence} makes both it and its inverse
continuous. If $P_n\to P$ in norm with each $P_n$ rank one, then
\[
 \|P^*-P\|\le2\|P-P_n\|\to0,\qquad
 \|P^2-P\|\le(\|P\|+2)\|P-P_n\|\to0,
 \qquad |\|P\|-1|\le\|P-P_n\|\to0.
\]
Hence $P^2=P=P^*$ and $\|P\|=1$. For large \(n\), \(\|P_n-P\|<1\), so rank stability gives \(\operatorname{rank}P=1\). Thus the rank-one projection set is closed in
$\mathcal K(\mathscr H)$. This Banach space is separable: finite matrices with rational real and imaginary entries in a fixed countable orthonormal basis are dense in its
finite-rank operators, and finite-rank operators are norm dense in $\mathcal K(\mathscr H)$. In the real case only rational real entries are needed.
\end{proof}

\subsubsection{Functional calculus and lower spectral bounds}
\label{W:sec:functional-calculus-and-ordered-eigenvalues}
\label{W:sec:lower-bound-for-the-spectrum}

The functional calculus converts the resolvent limit into convergence for continuous functions vanishing at infinity. This controls bounded spectral regions. Ordering eigenvalues
from the bottom also requires a lower spectral bound.

\begin{proposition}
\label{W:prop:spectral-functional-calculus}
Under the hypotheses of Theorem~\ref{W:thm:operator-deterministic}, for every $g\in C_0(\mathbb R)$,
\begin{equation}
 \|\mathcal E_Ng(\mathcal A_N)\mathcal E_N^*-g(\mathcal H_Y)\|\to0.
 \label{W:eq:spectral-C0}
\end{equation}
The convergence is uniform on compact subsets of $C_0(\mathbb R)$. For $g\in C_c(\mathbb R)$ it also holds in trace norm.
\end{proposition}
\begin{proof}
Let $r_\pm(x)=(x\mp\ii)^{-1}$ and let $\mathcal A$ be their algebra without a constant term. Each element vanishes at infinity, $\overline{r_+}=r_-$, and $r_+(x)=r_+(y)$ implies
$x=y$. Moreover $r_+(x)\ne0$ for every real $x$. Thus $\mathcal A$ is closed under complex conjugation, separates points, and vanishes nowhere. The locally compact
Stone--Weierstrass theorem makes $\mathcal A$ dense in $C_0(\mathbb R)$. For bounded operators $B,C$ on the finite space,
\[
 (\mathcal E_NB\mathcal E_N^*)(\mathcal E_NC\mathcal E_N^*)
 =\mathcal E_NBC\mathcal E_N^*,\qquad
 (\mathcal E_NB\mathcal E_N^*)^*=\mathcal E_NB^*\mathcal E_N^*.
\]
These identities justify products and complex conjugation under embedding. For every $\epsilon>0$, choose a polynomial $p$ in the two generators, without a constant term, such that
$\|g-p\|_\infty<\epsilon$. Write $p=\sum_{a=1}^J c_a r_{a,1}\cdots r_{a,m_a}$, where $m_a\ge1$ and each $r_{a,b}(x)$ is $(x-\ii)^{-1}$ or $(x+\ii)^{-1}$. Set
\[
 \epsilon_N=\max_{z\in\{\ii,-\ii\}}
 \|\mathcal E_N(\mathcal A_N-z)^{-1}\mathcal E_N^*
                         -(\mathcal H_Y-z)^{-1}\|\longrightarrow0.
\]
Each resolvent factor has norm at most one. The identity
\[
 \prod_{b=1}^{m}A_b-\prod_{b=1}^{m}B_b
 =\sum_{b=1}^{m}A_1\cdots A_{b-1}(A_b-B_b)B_{b+1}\cdots B_m
\]
therefore bounds a monomial difference by $m\epsilon_N$, and gives
\[
 \|\mathcal E_Np(\mathcal A_N)\mathcal E_N^*-p(\mathcal H_Y)\|
 \le\epsilon_N\sum_{a=1}^J|c_a|m_a\longrightarrow0.
\]
Consequently
\begin{align*}
 \|\mathcal E_Ng(\mathcal A_N)\mathcal E_N^*-g(\mathcal H_Y)\|
 &\le2\|g-p\|_\infty+
     \|\mathcal E_Np(\mathcal A_N)\mathcal E_N^*-p(\mathcal H_Y)\|,\\
 \limsup_N\|\mathcal E_Ng(\mathcal A_N)\mathcal E_N^*-g(\mathcal H_Y)\|
 &\le2\epsilon.
\end{align*}
Letting $\epsilon\downarrow0$ proves \eqref{W:eq:spectral-C0}. For a compact set $\mathcal G\subset C_0(\mathbb R)$, choose a finite $\epsilon$-net $g_1,\ldots,g_m$ in supremum
norm. The same spectral norm bound gives
\[
 \sup_{g\in\mathcal G}\|\mathcal E_Ng(\mathcal A_N)\mathcal E_N^*
                         -g(\mathcal H_Y)\|
 \le2\epsilon+\max_{j\le m}
     \|\mathcal E_Ng_j(\mathcal A_N)\mathcal E_N^*-g_j(\mathcal H_Y)\|.
\]
Taking the two limits proves uniformity on $\mathcal G$.

For $g\in C_c(\mathbb R)$ choose $a<b$ in $\rho(\mathcal H_Y)\cap\mathbb R$ so that $\operatorname{supp}g\Subset I=(a,b)$. Such choices of $a,b$ exist because the compact resolvent
of $\mathcal H_Y$ makes its spectrum discrete. Write
\[
 Q_N=\mathbf1_I(\mathcal A_N),\quad
 P_N^I=\mathcal E_NQ_N\mathcal E_N^*,\quad
 P^I=\mathbf1_I(\mathcal H_Y).
\]
These spectral projections commute with the respective functional calculi. Since $g$ vanishes on $\mathbb R\setminus I$,
\[
 \begin{aligned}
 g(\mathcal A_N)&=Q_Ng(\mathcal A_N)Q_N,\\
 g(\mathcal H_Y)&=P^Ig(\mathcal H_Y)P^I.
 \end{aligned}
\]
By Proposition~\ref{W:prop:trace-riesz-stability}, there is a finite integer $r=\operatorname{rank}P^I$ such that $\operatorname{rank}P_N^I=r$ for all large $N$. Set
$D_N=\mathcal E_Ng(\mathcal A_N)\mathcal E_N^*-g(\mathcal H_Y)$.
The preceding identities imply
\[
 \operatorname{Ran}D_N\subset
 \operatorname{Ran}P_N^I+\operatorname{Ran}P^I,
 \qquad \operatorname{rank}D_N\le2r.
\]
If $r=0$, rank equality gives $D_N=0$ for every $N$ above the index for the fixed interval $I$. Otherwise, summing its nonzero singular values gives the explicit trace norm
estimate
\[
 \begin{aligned}
 \|D_N\|_1
 &=\sum_{j=1}^{\operatorname{rank}D_N}s_j(D_N)
 \le \operatorname{rank}D_N\,\|D_N\|\\
 &\le2r\,\|\mathcal E_Ng(\mathcal A_N)\mathcal E_N^*
                         -g(\mathcal H_Y)\|\longrightarrow0,
 \end{aligned}
\]
where the last limit is the already proved operator-norm convergence \eqref{W:eq:spectral-C0}. The finite rank $r$ depends on the chosen support interval $I$ and the limiting
operator, and is independent of $N$.
\end{proof}

To apply this functional calculus to $e^{-tx}$, we first restrict that function below a common lower spectral bound. Such a bound also allows eigenvalues to be ordered through a
shifted inverse. It follows from localization of the spectral subspace $\operatorname{ran}\Pi_N(0)$ and the local quadratic form estimate.

Let $\mathcal A_N=h_N^{-2}(2I-T_N)$ and put $Y_N^{\mathrm{st}}=\widetilde Y_N^{\mathrm{st}}$. For deterministic arrays satisfying \eqref{W:eq:operator-deterministic-local}, for
each $L<\infty$ there is $N_L$ such that the off-diagonal coefficients with indices used on $[0,L]$ lie in $[1/2,2]$ for every $N\ge N_L$. In the probabilistic argument below, this
bound is imposed through the event $B_N(L+4)\le1/2$.

\begin{proposition}
\label{W:prop:spectral-lower-bound-in-probability}
Under Assumption~\ref{W:ass:wigner-lee-yin},
\begin{equation}
 \lim_{M\to\infty}\limsup_N\sup_{\|w\|=1}
   \mathbb P\bigl(\lambda_{\min}(\mathcal A_N)<-M\bigr)=0.
 \label{W:eq:assembly-bottom-tightness}
\end{equation}
Under the deterministic hypotheses of Theorem~\ref{W:thm:operator-deterministic}, there are $B<\infty$ and an integer $N_0$ such that $\mathcal A_N\ge-BI$ for every $N\ge N_0$.
\end{proposition}
\begin{proof}
On the intersection of the events
\[
 \|Q_{N,L}\Pi_N(0)\|\le\tfrac12,\qquad B_N(L+4)\le\tfrac12,
 \qquad \|Y_N^{\rm st}\|_{\infty,[0,L+4]}\le M,
\]
the coefficients in that window lie in $[1/2,3/2]\subset[1/2,2]$. Lemma~\ref{W:lem:operator-bottom}, with $b_*=1/2$, $b^*=2$ and $\theta=1/2$, yields
\[
 \lambda_{\min}(\mathcal A_N)
 \ge-\frac{2M^2}{1/2}-\frac{8}{1-1/4}
 =-4M^2-\frac{32}{3}\ge-4M^2-16.
\]
Taking complements and then the supremum over deterministic starting vectors gives, for fixed $L,M$ and all sufficiently large $N$,
\begin{align*}
 &\sup_{\|w\|=1}\mathbb P(\lambda_{\min}(\mathcal A_N)<-4M^2-16)\\
 &\quad\le\sup_{\|w\|=1}\mathbb P(\|Q_{N,L}\Pi_N(0)\|>1/2)
       +\sup_{\|w\|=1}\mathbb P(B_N(L+4)>1/2)\\
 &\qquad\quad+\sup_{\|w\|=1}
       \mathbb P(\|Y_N^{\mathrm{st}}\|_{\infty,[0,L+4]}>M).
\end{align*}
Given $\epsilon>0$, Theorem~\ref{W:thm:localization-spectral-projection} first chooses $L$ for which the limsup of the first probability is less than $\epsilon$, uniformly over
deterministic $w$. Theorem~\ref{W:thm:integrated-potential-and-offdiagonal-coefficients} makes the second term tend to zero. Its uniform tightness on $[0,L+5]$ supplies $M$ for
which the limsup of the third supremum is less than $\epsilon$. The piecewise constant integrated potential is bounded by the linearly interpolated integrated potential on
$[0,L+5]$. For these choices, put $M_\epsilon'=4M^2+16$. The displayed probability inequality implies
\[
 \limsup_N\sup_{\|w\|=1}
 \mathbb P(\lambda_{\min}(\mathcal A_N)<-M_\epsilon')\le2\epsilon.
\]
For every $R\ge M_\epsilon'$, the event with threshold $-R$ is a subset of the event with threshold $-M_\epsilon'$. Taking $R\to\infty$ and then $\epsilon\downarrow0$ proves
\eqref{W:eq:assembly-bottom-tightness}. For deterministic data choose $L$ so that $\limsup_N\|Q_{N,L}\Pi_N(0)\|<1/2$, and put $M=1+\|Y\|_{C([0,L+5])}$. Local convergence of the
polygonal potentials bounds the step potentials on $[0,L+4]$ by $M$ for all sufficiently large $N$. At those sizes the off-diagonal coefficients are between $1/2$ and $2$, and
$\|Q_{N,L}\Pi_N(0)\|\le1/2$. The deterministic lower bound lemma therefore gives
\[
 \lambda_{\min}(\mathcal A_N)\ge-4M^2-32/3
\]
for every sufficiently large $N$, with this one fixed $M$.
\end{proof}

\begin{corollary}
\label{W:cor:semigroup-convergence}
Under the deterministic hypotheses of Theorem~\ref{W:thm:operator-deterministic}, for $0<t_0<t_1<\infty$,
\[
 \sup_{t_0\le t\le t_1}
 \|\mathcal E_Ne^{-t\mathcal A_N}\mathcal E_N^*-e^{-t\mathcal H_Y}\|
 \longrightarrow0.
\]
\end{corollary}
\begin{proof}
Proposition~\ref{W:prop:spectral-lower-bound-in-probability} gives $B<\infty$ and $N_0$ with $\mathcal A_N\ge-BI$ for $N\ge N_0$. Increase $B$ so that $\mathcal H_Y\ge-BI$ as well.
Increase $B$ to be nonnegative and set $\chi(x)=\min\{1,\max\{0,x+B+1\}\}$ and $g_t(x)=\chi(x)e^{-tx}$. Then $g_t$ vanishes for $x\le-B-1$, and for $x\ge0$,
$\sup_{t_0\le t\le t_1}|g_t(x)|\le e^{-t_0x}\to0$. On the spectra of all sufficiently large finite operators and of $\mathcal H_Y$, $\chi=1$. Integrating against their spectral
measures therefore gives $g_t(\mathcal A_N)=e^{-t\mathcal A_N}$ and $g_t(\mathcal H_Y)=e^{-t\mathcal H_Y}$. For $0<t_0\le t\le t_1$, these functions belong to $C_0$ and satisfy
\[
 \|g_t-g_s\|_\infty
 \le|t-s|\sup_{x\ge-B-1}\sup_{\tau\in[t_0,t_1]}
                                  |x|e^{-\tau x}\le C|t-s|,
\]
where one may take
\[
 C=\max\{(B+1)e^{t_1(B+1)},(e t_0)^{-1}\}.
\]
Thus $\{g_t:t_0\le t\le t_1\}$ is a compact subset of $C_0(\mathbb R)$. The compact family conclusion of Proposition~\ref{W:prop:spectral-functional-calculus} gives the claimed
uniform semigroup convergence.
\end{proof}

\subsubsection{Absence of atoms in the eigenvalue distributions of the limiting operator}
\label{W:sec:atomlessness-and-inverse-convergence}

For the stochastic Airy operator with Dirichlet boundary condition, a deterministic real number is almost surely outside the spectrum. Resolvent continuation applies pathwise to
the random set $\rho(\mathcal H_\beta)$; evaluation at a fixed real number also uses atomlessness. We prove it by conditioning on a Gaussian coordinate which adds a strictly
positive bounded function to the potential. This atomlessness argument holds for every fixed $\beta>0$.

\Needspace{6\baselineskip}
\begin{lemma}
\label{W:lem:sao-fixed-energy-atomlessness}
For every fixed $\beta>0$, let $(\Lambda_j)_{j\ge1}$ be the ordered eigenvalues of the stochastic Airy operator with Dirichlet boundary condition
\eqref{W:eq:sao-formal}--\eqref{W:eq:sao-form}. For every $j\ge1$ and every deterministic
$a\in\mathbb R$,
\begin{equation}
 \mathbb P(\Lambda_j=a)=0,\qquad
 \mathbb P\bigl(a\in\sigma(\mathcal H_\beta)\bigr)=0.
 \label{W:eq:sao-fixed-energy-atomlessness}
\end{equation}
In particular, $\mathcal H_\beta^{-1}$ is a bounded compact operator almost surely. More generally, for a deterministic finite or countable set $D\subset\mathbb R$, almost surely
every $a\in D$ belongs to $\rho(\mathcal H_\beta)$.
\end{lemma}
\begin{proof}
Fix $\beta>0$ and put $c_\beta=2/\sqrt\beta$. Let $\mathscr G$ be the continuous real paths $b$ with $b(0)=0$ and
\begin{equation}
 C_b:=\sup_{k\ge0}\frac{1}{\sqrt{\log(2+k)}}
       \sup_{s,t\in[k,k+2]}|b(t)-b(s)|<\infty,
 \label{W:eq:sao-atomless-admissible-paths}
\end{equation}
where $k$ is integer. Writing the inner supremum as $O_k(b)$, the inequality $|O_k(b)-O_k(\widetilde b)|\le 2\|b-\widetilde b\|_{C([k,k+2])}$ shows that
\[
 \mathscr G=\{b:b(0)=0\}\cap\bigcup_{m\ge1}\bigcap_{k\ge0}
       \{b:O_k(b)\le m\sqrt{\log(2+k)}\}
\]
is Borel. The reflection principle gives $\mathbb P(O_k(B)>v)\le4e^{-v^2/16}$; hence Borel--Cantelli, with $v=A\sqrt{\log(2+k)}$ and $A>\sqrt{32}$, gives
$\mathbb P(B\in\mathscr G)=1$.

For $b\in\mathscr G$, let $\mathfrak q_{\beta,b}$ be \eqref{W:eq:sao-form} with $B=b$. Its Hermitian sesquilinear form is obtained by polarization, with the scalar product linear
in the first argument. Put $\overline b(x)=\int_x^{x+1}b(s)\,\dd s$. On $[k,k+1]$,
\[
 |b-\overline b|+|\overline b'|
 \le2O_k(b)\le C_b'\sqrt{\log(2+x)},\qquad
 -\int b(|f|^2)'=-2\operatorname{Re}\int(b-\overline b)f'\overline f
                                  +\int\overline b'|f|^2.
\]
The integration by parts has no boundary term for $f\in C_c^\infty(0,\infty)$. Young's inequality gives
\[
 \begin{split}
 \left|c_\beta\int b(|f|^2)'\right|
 &\le 2c_\beta C_b'\int\sqrt{\log(2+x)}\,|f'||f|
       +c_\beta C_b'\int\sqrt{\log(2+x)}\,|f|^2\\
 &\le\frac{\varepsilon}{2}\|f'\|_2^2
   +\int\left(\frac{2c_\beta^2(C_b')^2}{\varepsilon}\log(2+x)
                    +c_\beta C_b'\sqrt{\log(2+x)}\right)|f|^2.
 \end{split}
\]
Since $\log(2+x)=o(x)$, this proves, for every $\varepsilon>0$,
\begin{equation}
 \left|c_\beta\int_0^\infty b(x)(|f(x)|^2)'\,\dd x\right|
 \le\varepsilon\left(\|f'\|_2^2+\|x^{1/2}f\|_2^2\right)
       +C_{\varepsilon,\beta,b}\|f\|_2^2,
 \qquad f\in C_c^\infty(0,\infty).
 \label{W:eq:sao-atomless-pathwise-form-bound}
\end{equation}
For $\varepsilon=1/2$, write $C=C_{1/2,\beta,b}\ge0$. The resulting bounds imply
\[
 \tfrac12\|f\|_{\mathcal V}^2
 \le\mathfrak q_{\beta,b}[f]+(C+1)\|f\|_2^2
 \le(2C+2)\|f\|_{\mathcal V}^2,
\]
where the norm is defined below. Lemma~\ref{W:lem:operator-limit-operator}, applied to $Y(x)=x^2/2+c_\beta b(x)$, therefore gives the common closed form domain
\begin{equation}
 \mathcal V=\{f\in H_0^1(0,\infty):x^{1/2}f\in L^2(0,\infty)\},
 \qquad
 \|f\|_{\mathcal V}^2
   =\|f'\|_2^2+\int_0^\infty(1+x)|f(x)|^2\,\dd x,
 \label{W:eq:sao-atomless-common-domain}
\end{equation}
and a lower semibounded self-adjoint operator $\mathcal H_{\beta,b}$ with compact resolvent. On $\mathcal V$, the Brownian term denotes the continuous extension of its value on the
form core.

Choose
\begin{equation}
 \begin{split}
 g(x)&=\sqrt2\,\mathrm e^{-x},\qquad
 h(x)=\int_0^x g(s)\,\dd s=\sqrt2(1-\mathrm e^{-x}),\\
 Z&=\int_0^\infty g(s)\,\dd B(s),\qquad
 \widehat B(x)=B(x)-Zh(x).
 \end{split}
 \label{W:eq:sao-atomless-gaussian-coordinate}
\end{equation}
The deterministic Wiener integral satisfies $Z\sim N(0,1)$ because $\|g\|_2=1$, and
\begin{equation}
 \mathbb E[B(x)Z]=h(x),\qquad
 \mathbb E[\widehat B(x)Z]=0.
 \label{W:eq:sao-atomless-gaussian-independence}
\end{equation}
Joint Gaussianity gives independence of $Z$ and every finite vector of values of $\widehat B$. Since evaluations at rational times generate the Borel sigma field of continuous
paths, $Z$ and the entire path $\widehat B$ are independent. Boundedness of $h$ gives $O_k(b+th)\le O_k(b)+2|t|\|h\|_\infty$. Thus $b+th\in\mathscr G$ for every $b\in\mathscr G$
and $t\in\mathbb R$, and $\widehat B\in\mathscr G$ almost surely.

Let $V_\beta=c_\beta g$, a bounded function positive everywhere. Integration by parts on the core and continuity on $\mathcal V$ give
\begin{equation}
 \mathfrak q_{\beta,b+t h}[f,u]
 =\mathfrak q_{\beta,b}[f,u]
        +t\int_0^\infty V_\beta(x)f(x)\overline{u(x)}\,\dd x,
 \qquad f,u\in\mathcal V.
 \label{W:eq:sao-atomless-bounded-form-perturbation}
\end{equation}
The characterization of the operator domain \eqref{W:eq:form-operator-domain}, applied in both directions, yields
\begin{equation}
 \mathcal H_{\beta,b+t h}=\mathcal H_{\beta,b}+tM_{V_\beta},\qquad
 \mathcal D(\mathcal H_{\beta,b+t h})
       =\mathcal D(\mathcal H_{\beta,b}),
 \label{W:eq:sao-atomless-bounded-operator-perturbation}
\end{equation}
since multiplication by $V_\beta$ is bounded on $L^2$.

Write $\lambda_j(b,t)$ for the ordered eigenvalues, with multiplicity. The min--max principle and $|\int V_\beta|f|^2|\le\|V_\beta\|_\infty\|f\|_2^2$ give
\begin{equation}
 |\lambda_j(b,t)-\lambda_j(b,s)|
       \le\|V_\beta\|_\infty|t-s|.
 \label{W:eq:sao-atomless-lipschitz}
\end{equation}
For $t_1<t_2$, let $S_j$ be the span of the first $j$ orthonormal eigenfunctions at $t_2$. Positivity of $V_\beta$ and compactness of the unit sphere of $S_j$ imply
\begin{equation}
 v_j(b,t_2):=\min_{\substack{f\in S_j\\\|f\|_2=1}}
           \int_0^\infty V_\beta(x)|f(x)|^2\,\dd x>0.
 \label{W:eq:sao-atomless-positive-compression}
\end{equation}
Indeed the integral is a continuous quadratic form, and its zero set consists only of the zero function because $V_\beta>0$. Using $S_j$ as a trial space at $t_1$ gives
\begin{align}
 \lambda_j(b,t_1)
 &\le\max_{\substack{f\in S_j\\\|f\|_2=1}}
                    \mathfrak q_{\beta,b+t_1h}[f]\notag\\
 &\le\lambda_j(b,t_2)-(t_2-t_1)v_j(b,t_2)
 <\lambda_j(b,t_2).
 \label{W:eq:sao-atomless-strict-monotonicity}
\end{align}

For measurability, choose an $L^2$-orthonormal sequence $(e_k)$ in $C_c^\infty(0,\infty)$ whose span is dense in $\mathcal V$. Gram--Schmidt applied to a countable dense subset of
the form core constructs such a sequence. The entries of
\[
 \bigl(\mathfrak q_{\beta,b+th}[e_s,e_r]\bigr)_{r,s=1}^m
\]
are continuous in $(b,t)$ for local uniform convergence of paths: on an interval $[0,L]$ containing their supports, their differences are bounded by
$c_\beta(\|b-\widetilde b\|_{C([0,L])} +|t-\widetilde t|\|h\|_\infty)\|(e_s\overline e_r)'\|_1$. Hence their ordered eigenvalues $\lambda_j^{(m)}(b,t)$, for $m\ge j$, are
continuous. Their trial spaces are nested, so min--max gives $\lambda_j^{(m)}\ge\lambda_j^{(m+1)}\ge\lambda_j(b,t)$.

To prove convergence, fix $(b,t)$ and approximate its first $j$ normalized eigenfunctions by $u_1^{(m)},\ldots,u_j^{(m)}\in\operatorname{span}(e_1,\ldots,e_m)$ in $\mathcal V$. Let
\[
 \begin{aligned}
 G_m&=(\langle u_s^{(m)},u_r^{(m)}\rangle)_{r,s=1}^j,\\
 Q_m&=(\mathfrak q_{\beta,b+th}[u_s^{(m)},u_r^{(m)}])_{r,s=1}^j,\\
 D_j&=\operatorname{diag}(\lambda_1(b,t),\ldots,\lambda_j(b,t)).
 \end{aligned}
\]
Continuity of the form on $\mathcal V$ gives $\delta_m:=\|G_m-I_j\|+\|Q_m-D_j\|\to0$. For $\delta_m<1$, $G_m$ is positive definite, and
\[
 \lambda_j(b,t)\le\lambda_j^{(m)}(b,t)
 \le\sup_{c\ne0}\frac{c^*Q_mc}{c^*G_mc}
 \le\lambda_j(b,t)+
       \frac{(1+|\lambda_j(b,t)|)\delta_m}{1-\delta_m}.
\]
The last inequality follows by subtracting $\lambda_j(b,t)c^*G_mc$ from the numerator and using $D_j\le\lambda_j(b,t)I_j$. Thus $\lambda_j^{(m)}\downarrow\lambda_j(b,t)$, and
$(b,t)\mapsto\lambda_j(b,t)$ is Borel on $\mathscr G\times\mathbb R$.

Let $\nu$ be the law of $\widehat B$. Independence and $B=\widehat B+Zh$ give, for deterministic $a$,
\begin{equation}
 \mathbb P(\Lambda_j=a)
 =\int_{\mathscr G}\int_{\mathbb R}
       \mathbf1_{\{\lambda_j(b,t)=a\}}
       \frac{\mathrm e^{-t^2/2}}{\sqrt{2\pi}}\,\dd t\,\nu(\dd b)
 =0,
 \label{W:eq:sao-atomless-conditional-integration}
\end{equation}
because strict monotonicity makes the set of values of the scalar variable in the inner integral empty or a singleton. Compact resolvent and lower semiboundedness give
$\sigma(\mathcal H_\beta)=\{\Lambda_j:j\ge1\}$ and $\Lambda_j\to+\infty$. A countable union therefore proves \eqref{W:eq:sao-fixed-energy-atomlessness} and the almost sure
exclusion of $D$ from the spectrum. In particular
\[
 \delta_\beta:=\operatorname{dist}(0,\sigma(\mathcal H_\beta))>0,
 \qquad \|\mathcal H_\beta^{-1}\|=\delta_\beta^{-1}<\infty
\]
almost surely, and the inverse is compact because $\Lambda_j^{-1}\to0$.
\end{proof}

\subsubsection{Inverses at real spectral parameters}
\label{W:sec:inverses-at-real-spectral-parameters}

The preceding lemma gives a bounded compact inverse $\mathscr K_\beta:=\mathcal H_\beta^{-1}$ almost surely. We now evaluate the operator limit at zero. The finite inverse must
first be defined on the exceptional singular event, since discrete entry laws are allowed.

For the finite matrix, let $\Omega_N^{\mathrm{inv}}=\{2\notin\sigma(H_N)\} =\{0\notin\sigma(\mathcal A_N)\}$, and define
\begin{equation}
 \mathscr K_N=
 \begin{cases}
  \mathcal E_N\mathcal A_N^{-1}\mathcal E_N^*,
                  &\text{on }\Omega_N^{\mathrm{inv}},\\
  0,               &\text{on }(\Omega_N^{\mathrm{inv}})^c.
 \end{cases}
 \label{W:eq:embedded-edge-inverse}
\end{equation}
On $\Omega_N^{\mathrm{inv}}$, unitary conjugation gives
\begin{equation}
 \mathscr K_N
 =N^{-2/3}\mathcal E_NU_N^*(2I-H_N)^{-1}U_N\mathcal E_N^*.
 \label{W:eq:original-wigner-edge-inverse}
\end{equation}
The value zero on the singular event makes $\mathscr K_N$ a compact operator-valued random variable for every $N$.

\Needspace{6\baselineskip}

\begin{corollary}
\label{W:cor:wigner-edge-inverse}
Under Assumption~\ref{W:ass:wigner-lee-yin}, let $w_N$ be any deterministic sequence of unit vectors satisfying \eqref{W:eq:deterministic-starting-vector}, and use its Lanczos
basis. Then
\begin{equation}
 \mathbb P\bigl(2\in\sigma(H_N)\bigr)\longrightarrow0,
 \label{W:eq:wigner-edge-invertibility}
\end{equation}
and
\begin{equation}
 \mathscr K_N\Longrightarrow\mathcal H_\beta^{-1}
 \quad\text{in }\mathcal K(L^2(\mathbb R_+))
 \quad\text{with the operator-norm topology}.
 \label{W:eq:main-inverse-convergence}
\end{equation}
The convergence is joint with the nonreal resolvent path.
\end{corollary}
There is a coupling of these operator-valued random variables with the same finite-$N$ and limiting laws on which $\Omega_N^{\mathrm{inv}}$ holds for all sufficiently large $N$ and
\begin{equation}
 \|\mathscr K_N-\mathcal H_\beta^{-1}\|\longrightarrow0
 \label{W:eq:coupled-inverse-convergence}
\end{equation}
almost surely. The event in \eqref{W:eq:wigner-edge-invertibility} is independent of the starting vector, since the unitary Lanczos transformation preserves the entire spectrum.
The value assigned on the singular event can be replaced by any measurable compact operator-valued convention without changing the distributional limit. The almost sure convergence
holds for the copies of the resolvent paths on the probability space in Theorem~\ref{W:thm:main-operator}.

\begin{proof}
Use the joint coupling of the nonreal resolvent paths in Theorem~\ref{W:thm:main-operator}, and put $T_N=\mathscr R_N(\ii)$, $T=\mathscr R_\beta(\ii)$. Almost surely $T_N\to T$ in
norm and $0\in\rho(\mathcal H_\beta)$ by Lemma~\ref{W:lem:sao-fixed-energy-atomlessness}. Thus $S=I+\ii T$ is invertible, with $S^{-1}=I-\ii\mathcal H_\beta^{-1}$. Set
$S_N=I+\ii T_N$ and $P_N=\mathcal E_N\mathcal E_N^*$. Relative to
$\operatorname{Ran}P_N\oplus\operatorname{Ran}P_N^\perp$,
\begin{equation}
 S_N=
 \left[\mathcal E_N\mathcal A_N(\mathcal A_N-\ii)^{-1}
                       \mathcal E_N^*\right]_{\operatorname{Ran}P_N}
 \oplus I_{\operatorname{Ran}P_N^\perp}.
 \label{W:eq:inverse-main-block-decomposition}
\end{equation}
Choose $N_0$ so that $\|S^{-1}\|\|S_N-S\|<1/2$ for every $N\ge N_0$. The inverse series gives invertibility of $S_N$, hence of $\mathcal A_N$, and
\[
 \|S_N^{-1}-S^{-1}\|\le2\|S^{-1}\|^2\|S_N-S\|\longrightarrow0.
\] Moreover,
\begin{equation}
 \mathscr K_N=T_NS_N^{-1},\qquad
 \mathcal H_\beta^{-1}=TS^{-1}.
 \label{W:eq:inverse-main-functional-identity}
\end{equation}
With $M=\|S^{-1}\|$, the resolvent identity yields, for large $N$,
\[
 \|\mathscr K_N-\mathcal H_\beta^{-1}\|
 =\|S_N^{-1}(T_N-T)S^{-1}\|
 \le2M^2\|T_N-T\|\longrightarrow0.
\]
The map $T\mapsto T(I+\ii T)^{-1}$, extended by zero where the denominator is singular, is Borel, since inversion is continuous on the open set of invertible operators. Thus the
coupled inverses have the original joint distributions with the resolvent paths. Their almost sure joint convergence proves \eqref{W:eq:main-inverse-convergence}; dominated
convergence for the singularity indicators proves \eqref{W:eq:wigner-edge-invertibility}. If $\widehat{\mathscr K}_N$ uses another measurable convention on that event, every
bounded continuous function $F$ of the inverse and the resolvent path satisfies
\[
 \left|\mathbb EF(\widehat{\mathscr K}_N,\mathscr R_N)
            -\mathbb EF(\mathscr K_N,\mathscr R_N)\right|
 \le2\|F\|_\infty\mathbb P((\Omega_N^{\mathrm{inv}})^c)\longrightarrow0.
\]
This also proves joint convergence for that convention.
\end{proof}

\Needspace{6\baselineskip} For finitely many deterministic real parameters, apply the inverse maps jointly to the resolvent path. At a singular finite matrix, the inverse is
assigned the value specified above.

Fix real numbers $a_1,\ldots,a_q$ and a deterministic starting sequence. For $a\in\{a_1,\ldots,a_q\}$, define
\begin{equation}
 \mathscr K_{N,a}=
 \begin{cases}
  \mathcal E_N(\mathcal A_N-a)^{-1}\mathcal E_N^*,
                         &a\notin\sigma(\mathcal A_N),\\
  0,                      &a\in\sigma(\mathcal A_N).
 \end{cases}
 \label{W:eq:real-energy-inverse-definition}
\end{equation}

\begin{corollary}
\label{W:cor:real-energy-inverse-convergence}
Under Assumption~\ref{W:ass:wigner-lee-yin}, jointly with the nonreal resolvent path,
\begin{equation}
 (\mathscr K_{N,a_\ell})_{\ell=1}^q
 \Longrightarrow
 \bigl((\mathcal H_\beta-a_\ell)^{-1}\bigr)_{\ell=1}^q
 \quad\text{in }\mathcal K(L^2(\mathbb R_+))^q.
 \label{W:eq:real-energy-inverse-convergence}
\end{equation}
Moreover,
\begin{equation}
 \mathbb P\left(
   \bigcup_{\ell=1}^q
   \{2-a_\ell N^{-2/3}\in\sigma(H_N)\}
             \right)\longrightarrow0.
 \label{W:eq:real-energy-wigner-invertibility}
\end{equation}
\end{corollary}
\begin{proof}
On the coupling of Theorem~\ref{W:thm:main-operator}, put $T_N=\mathscr R_N(\ii)$ and $T=\mathscr R_\beta(\ii)$. Then $T_N\to T$ in norm almost surely, and all $a_\ell$ belong to
$\rho(\mathcal H_\beta)$ almost surely by Lemma~\ref{W:lem:sao-fixed-energy-atomlessness}. Set $S_{N,\ell}=I-(a_\ell-\ii)T_N$ and $S_\ell=I-(a_\ell-\ii)T$. The spectral theorem
gives
\[
 S_\ell^{-1}=I+(a_\ell-\ii)(\mathcal H_\beta-a_\ell)^{-1}.
\]
For $|a_\ell-\ii|\|S_\ell^{-1}\|\|T_N-T\|<1/2$, the Neumann series gives invertibility of $S_{N,\ell}$ and $\|S_{N,\ell}^{-1}\|\le2\|S_\ell^{-1}\|$. Relative to
$\operatorname{ran}P_N\oplus\operatorname{ran}P_N^\perp$, its blocks are
\[
 S_{N,\ell}=
 \left[\mathcal E_N(\mathcal A_N-a_\ell)(\mathcal A_N-\ii)^{-1}
                   \mathcal E_N^*\right]_{\operatorname{ran}P_N}
       \oplus I_{\operatorname{ran}P_N^\perp},
 \qquad P_N=\mathcal E_N\mathcal E_N^*.
\]
Choose $N_0$ so that $\max_\ell\|S_\ell^{-1}\|\|S_{N,\ell}-S_\ell\|<1/2$ for every $N\ge N_0$. Then $a_\ell\notin\sigma(\mathcal A_N)$ for all indices $\ell$ in the finite list,
and
\[
 \begin{split}
 \mathscr K_{N,a_\ell}-(\mathcal H_\beta-a_\ell)^{-1}
 &=S_{N,\ell}^{-1}(T_N-T)S_\ell^{-1},\\
 \max_{\ell\le q}\|\mathscr K_{N,a_\ell}-(\mathcal H_\beta-a_\ell)^{-1}\|
 &\le2\max_{\ell\le q}\|S_\ell^{-1}\|^2\|T_N-T\|\longrightarrow0.
 \end{split}
\]
The inverse maps, with the stated zero convention, are Borel by continuity of inversion on its open domain. They preserve the joint law with the nonreal resolvent path. The
displayed almost sure convergence proves \eqref{W:eq:real-energy-inverse-convergence}, and dominated convergence for the finitely many singularity indicators proves
\eqref{W:eq:real-energy-wigner-invertibility}.

The inverse identity also applies to deterministic $a_N\to a$. On a path with $a\in\rho(\mathcal H_\beta)$, put $B_N=I-(a_N-\ii)T_N$, $B=I-(a-\ii)T$, and $M=\|B^{-1}\|$. Since
$\|T_N\|,\|T\|\le1$, $B_N\to B$, and multiplication by the denominators gives
\[
 \begin{split}
 T_NB_N^{-1}-TB^{-1}
 &=B_N^{-1}\{T_N-T+(a_N-a)T_NT\}B^{-1},\\
 \|T_NB_N^{-1}-TB^{-1}\|
 &\le2M^2(\|T_N-T\|+|a_N-a|)\longrightarrow0.
 \end{split}
\]
For deterministic $a$, atomlessness ensures $M<\infty$ almost surely.
\end{proof}

\begin{corollary}
\label{W:cor:real-spectral-distance}
Under the same assumptions, for each deterministic $a\in\mathbb R$,
\begin{equation}
 \operatorname{dist}(a,\sigma(\mathcal A_N))
 \Longrightarrow
 \operatorname{dist}(a,\sigma(\mathcal H_\beta)),
 \label{W:eq:real-energy-distance-convergence}
\end{equation}
and
\begin{equation}
 \lim_{\varepsilon\downarrow0}\limsup_{N\to\infty}
 \mathbb P\left(
    \min_{1\le j\le N}
    |\mu_j^{(N)}-(2-aN^{-2/3})|
        \le\varepsilon N^{-2/3}\right)=0.
 \label{W:eq:real-energy-no-concentration}
\end{equation}
\end{corollary}
\begin{proof}
On the coupling of the preceding corollary, fix an outcome on which the real inverses converge. For $N$ above the invertibility index constructed there, the spectral theorem gives
\[
 \begin{split}
 \operatorname{dist}(a,\sigma(\mathcal A_N))
 &=\|\mathscr K_{N,a}\|^{-1}
 \longrightarrow\| (\mathcal H_\beta-a)^{-1}\|^{-1}\\
 &=\operatorname{dist}(a,\sigma(\mathcal H_\beta))>0
 \quad\text{almost surely}.
 \end{split}
\]
This proves \eqref{W:eq:real-energy-distance-convergence}. Writing these distances as $D_N,D$, the Portmanteau theorem implies
\[
 \lim_{\varepsilon\downarrow0}\limsup_N\mathbb P(D_N\le\varepsilon)
 \le\lim_{\varepsilon\downarrow0}\mathbb P(D\le\varepsilon)
 =\mathbb P(D=0)=0.
\]
Finally $D_N=N^{2/3}\min_j|\mu_j^{(N)}-(2-aN^{-2/3})|$, which proves \eqref{W:eq:real-energy-no-concentration}. This is a statement about the original spectrum, independent of the
starting vector.
\end{proof}

The preceding inversion has a pathwise stability estimate. On a coupled path, put $\varepsilon_N=\|\mathscr R_N(\ii)-\mathscr R_\beta(\ii)\|$ and use $S_N,S$ from the proof of
Corollary~\ref{W:cor:wigner-edge-inverse}. Whenever $\|S^{-1}\|\varepsilon_N\le1/2$, the resolvent algebra gives
\begin{equation}
 \begin{split}
 \mathscr K_N-\mathcal H_\beta^{-1}
 &=S_N^{-1}(\mathscr R_N(\ii)-\mathscr R_\beta(\ii))S^{-1},\\
 \|\mathscr K_N-\mathcal H_\beta^{-1}\|
 &\le2(1+\delta_\beta^{-2})\varepsilon_N.
 \end{split}
 \label{W:eq:inverse-limit-stability-bound}
\end{equation}
Indeed, multiplication of the first identity on the left by $S_N$ and on the right by $S$ leaves
$\mathscr R_N(\mathrm i)S-S_N\mathscr R_\beta(\mathrm i) =\mathscr R_N(\mathrm i)-\mathscr R_\beta(\mathrm i)$. The smallness assumption gives $\|S_N^{-1}\|\le2\|S^{-1}\|$. Finally
the spectral theorem gives
\[
 \|S^{-1}\|^2
 =\sup_{\lambda\in\sigma(\mathcal H_\beta)}|1-\mathrm i/\lambda|^2
 =1+\delta_\beta^{-2}.
\]
Thus the stability constant is determined by the distance from zero to the limiting spectrum. The finite inverse has the centering and embedding specified in
\eqref{W:eq:original-wigner-edge-inverse}.

\begin{proof}[Proof of Corollary~\ref{W:cor:intro-real-inverse}]
On the joint representation used in Corollary~\ref{W:cor:real-energy-inverse-convergence}, put
\[
 T_N=\mathscr R_N(\mathrm i),\qquad T=(\mathcal H_\beta-\mathrm i)^{-1}.
\]
There is an event of probability one on which $\|T_N-T\|_{\mathrm{op}}\to0$ and $0\in\rho(\mathcal H_\beta)$. Fix a sample in this event and let $M=\|(I+\mathrm iT)^{-1}\|<\infty$.
Choose $N_0$ such that $M\|T_N-T\|_{\mathrm{op}}\le1/2$ for every $N\ge N_0$. Factoring $I+\mathrm iT_N$ about $I+\mathrm iT$ and using the Neumann series gives
\[
 \|(I+\mathrm iT_N)^{-1}\|\le2M.
\]
On $\operatorname{Ran}\mathcal E_N$, the denominator is unitarily equivalent to
\[
 I+\mathrm i(\mathcal A_N-\mathrm i)^{-1}
       =\mathcal A_N(\mathcal A_N-\mathrm i)^{-1},
\]
and it is the identity on the orthogonal complement. Its invertibility is therefore equivalent to invertibility of $\mathcal A_N$. For $N\ge N_0$,
\[
 T_N(I+\mathrm iT_N)^{-1}
   =\mathcal E_N\mathcal A_N^{-1}\mathcal E_N^*,
 \qquad T(I+\mathrm iT)^{-1}=\mathcal H_\beta^{-1}.
\]
Multiplication by the two denominators, in their displayed order, proves
\[
 \begin{split}
 T_N(I+\mathrm iT_N)^{-1}-T(I+\mathrm iT)^{-1}
   &=(I+\mathrm iT_N)^{-1}(T_N-T)(I+\mathrm iT)^{-1},\\
 \|\mathcal E_N\mathcal A_N^{-1}\mathcal E_N^*-
                 \mathcal H_\beta^{-1}\|
   &\le2M^2\|T_N-T\|\longrightarrow0.
 \end{split}
\]
The exact centering and conjugation in \eqref{W:eq:original-wigner-edge-inverse} identify this finite inverse with the one in the corollary.

The singularity indicators tend to zero on this representation. Their expectations therefore tend to zero by bounded convergence, and equality of the finite marginal laws proves
the same statement for the original matrices. Assigning zero to the inverse on the singular event makes it a Borel random compact operator. Any other Borel convention changes a
bounded continuous observable by at most twice its supremum norm times the singularity probability. Consequently the convention does not change the distributional limit. Since the
inverse converges on the same representation as the entire resolvent path, the convergence is joint with \eqref{W:eq:main-resolvent-convergence}.
\end{proof}

\subsubsection{Local weak equations and energy convergence}
\label{W:sec:local-energy-estimates-and-ordered-eigenfunctions}

We next use the Jacobi coefficients to pass to the local weak equation and strengthen convergence of the eigenfunctions. The coefficients and operators in the following arguments
are deterministic.

We use $h_N=N^{-1/3}$, $\mathcal A_N={h_N}^{-2}(2I-T_N)$, $X_N=\widetilde Y_N^{\rm lin}$, $\mathscr H=L^2(\mathbb R_+)$, and the isometric step embedding $\mathcal E_N$ and linear
embedding $\mathcal I_N$ already defined. The off-diagonal coefficients of $T_N$ are nonnegative and vanish between its invariant blocks. The auxiliary coefficient at zero is one.
Write
\[
 \mathscr R_N(z)=\mathcal E_N(\mathcal A_N-z)^{-1}\mathcal E_N^*,\qquad
 \mathscr S_N(z)=\mathcal I_N(\mathcal A_N-z)^{-1}\mathcal E_N^*,\qquad
 P_N=\mathcal E_N\mathcal E_N^*.
\]
For $z\in\mathbb C\setminus\mathbb R$, $\mathscr R_N(z)$ vanishes on $\operatorname{ran}P_N^\perp$. For the bounded operator $\mathcal E_N\mathcal A_N\mathcal E_N^*$,
\[
 (\mathcal E_N\mathcal A_N\mathcal E_N^*-z)^{-1}
       =\mathscr R_N(z)-z^{-1}(I-P_N).
\]
We use the embedded finite resolvent $\mathscr R_N$ throughout. Scalar products are linear in the first argument.

\Needspace{6\baselineskip}

Assume $X_N\to Y$ locally uniformly and $\max_{k\le N^{1/3}L}|b_{k,N}-1|\to0$ for every $L<\infty$, where $Y$ is real and continuous. Suppose $z_N\to z$, $f_N\rightharpoonup f$ in
$\mathscr H$, and
\[
 (\mathcal A_N-z_N)v_N=\mathcal E_N^*f_N,\qquad
 \sup_N\|v_N\|<\infty,\qquad \mathcal E_Nv_N\to u\quad\hbox{in }L^2.
\]

\begin{lemma}
\label{W:lem:spectral-identification}
Then $u\in H^1_{\rm loc}([0,\infty))$, $u(0)=0$, and
\begin{equation}
 \int u'\overline{\varphi'}-\int Y(u\overline\varphi)'
       -z\langle u,\varphi\rangle=\langle f,\varphi\rangle,
       \qquad\varphi\in C_c^\infty(0,\infty).
 \label{W:eq:spectral-weak-limit}
\end{equation}
If $Y$ satisfies \eqref{W:eq:operator-limit-coercivity} and $z\notin\mathbb R$, then $u=(\mathcal H_Y-z)^{-1}f$.
\end{lemma}
\begin{proof}
Put $Y_N=\widetilde Y_N^{\rm st}$ and $F_N=\mathcal I_Nv_N$. Adjacent jumps of the integrated potential are at most $2\|X_N-Y\|_{C([0,L+1])}+\omega_Y(h_N;L+1)$; hence $Y_N\to Y$
uniformly on each fixed interval. For each such interval, choose $N_L$ so that the off-diagonal coefficients lie between $b_*=1/2$ and $b^*=2$ for all $N\ge N_L$. Choose a real
cutoff $\chi\in C_c^2([0,\infty))$ with $0\le\chi\le1$, $\chi=1$ on $[0,L+1]$, support in $[0,L+2]$ and $\|\chi'\|_\infty\le2$. With $M=\sup_N\|Y_N\|_{\infty,[0,L+3]}$, the
discrete cutoff identity \eqref{W:eq:discrete-cutoff-formula} and local coercivity \eqref{W:eq:operator-local-coercive} give
\[
 \frac{b_*}{2}\|(\mathcal I_N(\chi v_N))'\|_2^2
 \le\|f_N\|\|v_N\|
       +\left(|z_N|+4b^*+\frac{2M^2}{b_*}\right)\|v_N\|^2\le C_L.
\]
Here weak convergence of $f_N$ gives $\sup_N\|f_N\|<\infty$, and $z_N\to z$ gives $\sup_N|z_N|<\infty$. Thus $C_L$ depends only on these bounds, $\sup_N\|v_N\|$, the cutoff, and
the local coefficient and potential bounds. In particular, $F_N$ is bounded in $H^1(0,L)$, including for real $z_N$. On each full cell, with $h=h_N$,
\[
 \int_{kh}^{(k+1)h}|F_N-\mathcal E_Nv_N|^2
 =\tfrac13|(v_N)_{k+1}-(v_N)_k|^2
 =\tfrac{h^2}{3}\int_{kh}^{(k+1)h}|F_N'|^2.
\]
Consequently $F_N\to u$ strongly in $L^2(0,L)$. Every weakly convergent $H^1(0,L)$ subsequence therefore has limit $u$, and the bounded sequence converges weakly to $u$ in
$H^1(0,L)$. Continuity of the trace and $F_N(0)=0$ imply
$u(0)=0$.

We record the cutoff interpolation estimate needed below. Write $u_k=h^{-1/2}(v_N)_k$, $\chi_k=\chi(kh)$ and let $\chi^{\rm lin}$ interpolate the $\chi_k$. On $x=kh+\theta h$,
\[
 \mathcal I_N(\chi v_N)-\chi^{\rm lin}F_N
 =\theta(1-\theta)(\chi_{k+1}-\chi_k)(u_{k+1}-u_k).
\]
The integrals of $\theta^2(1-\theta)^2$ and $(1-2\theta)^2$ are $1/30$ and $1/3$. Summing over the cells meeting $[0,L]$ therefore bounds the $L^2$ norm of this error by
$h^2\|\chi'\|_\infty\|F_N'\|_{L^2(0,L+1)}/\sqrt{30}$ and its derivative norm by $h\|\chi'\|_\infty\|F_N'\|_{L^2(0,L+1)}/\sqrt3$. The scalar interpolation estimates
$\|\chi^{\rm lin}-\chi\|_\infty\le h^2\|\chi''\|_\infty/8$ and $\|(\chi^{\rm lin})'-\chi'\|_\infty\le h\|\chi''\|_\infty$ then give
\begin{equation}
 \|\mathcal I_N(\chi v_N)-\chi\mathcal I_Nv_N\|_{H^1(0,L)}
 \le C_\chi h_N\|\mathcal I_Nv_N\|_{H^1(0,L+1)}.
 \label{W:eq:grid-values-cutoff-product-error}
\end{equation}
The same estimate applies to $\chi^2$.

For $\varphi\in C_c^\infty(0,\infty)$ set $w_{N,k}=\sqrt{h_N}\varphi(kh_N)$ and $g_N=\mathcal I_Nw_N$. The scalar interpolation estimates give $g_N\to\varphi$ in $H^1$ and
$\mathcal E_Nw_N\to\varphi$ in $L^2$, with all supports in one fixed interval. The exact sesquilinear weak equation is
\[
 \int\mathfrak b_NF_N'\overline{g_N'}
 -\int Y_N(F_N\overline{g_N})'
 -z_N\langle\mathcal E_Nv_N,\mathcal E_Nw_N\rangle
 =\langle f_N,\mathcal E_Nw_N\rangle.
\]
Here $\mathfrak b_Ng_N'\to\varphi'$ and $Y_Ng_N\to Y\varphi$ strongly in $L^2$, while $F_N'\rightharpoonup u'$. Also $Y_NF_N\to Yu$ strongly and $g_N'\to\varphi'$ strongly.
Expanding $(F_N\overline{g_N})'$ therefore passes both potential terms to the limit. The mass and forcing terms converge by the assumed strong convergence of $\mathcal E_Nv_N$ and
weak convergence of $f_N$. This proves \eqref{W:eq:spectral-weak-limit}.

Suppose now that \eqref{W:eq:operator-limit-coercivity} holds. Choose $0\le\chi_R\le1$, equal to one on $[0,R]$, zero after $R+1$, with $\|\chi_R'\|_\infty\le2$. On a fixed
interval $(0,M)$ containing its support, the bound
\[
 |\mathfrak q_Y[u,\psi]|
 \le\|u'\|_2\|\psi'\|_2+
 \|Y\|_\infty(\|u'\|_2\|\psi\|_2+\|u\|_2\|\psi'\|_2)
\]
extends the weak equation by density to $\chi_R^2u\in H_0^1(0,M)$. The product rule gives
\[
 \mathfrak q_Y[\chi_Ru]
 =\operatorname{Re}\mathfrak q_Y[u,\chi_R^2u]
                  +\int|\chi_R'|^2|u|^2
 =\operatorname{Re}\langle f+zu,\chi_R^2u\rangle
                  +\int|\chi_R'|^2|u|^2.
\]
Applying coercivity by the same core approximation yields
\[
 c_Y\|\chi_Ru\|_{\mathcal V}^2
 \le\|f\|\|u\|+(|z|+C_Y+4)\|u\|^2.
\]
For $R>L$ this bounds the $\mathcal V$ energy on $[0,L]$ uniformly in $L$, so $u\in\mathcal V$. Density of the core in the form norm and continuity of the form in
Lemma~\ref{W:lem:operator-limit-operator} now extend $\mathfrak q_Y[u,\psi]=\langle f+zu,\psi\rangle$ to every $\psi\in\mathcal V$. Hence $u\in\mathcal D(\mathcal H_Y)$ and
$(\mathcal H_Y-z)u=f$. Self-adjointness gives the unique solution $u=(\mathcal H_Y-z)^{-1}f$ when $z\notin\mathbb R$.
\end{proof}

\subparagraph{Convergence of the local energies}
\label{W:sec:convergence-of-local-energies}

The cutoff quadratic form identity gives convergence of the derivative norms, which strengthens the preceding weak $H^1$ convergence to strong convergence on bounded intervals.

\begin{lemma}
\label{W:lem:spectral-energy}
Under the hypotheses of Lemma~\ref{W:lem:spectral-identification}, including its $L^2$ convergence, one has
\begin{equation}
 \mathcal I_Nv_N\longrightarrow u\quad\hbox{in }H^1(0,L)
               \quad\hbox{for every }L<\infty.
 \label{W:eq:spectral-strong-energy}
\end{equation}
The conclusion includes real $z$ and zero forcing. It also holds for continuum solutions: if $u_n\in H^1_{\rm loc}([0,\infty))\cap L^2$, $u_n(0)=0$, and
\[
 \int u_n'\overline{\varphi'}-\int Y(u_n\overline\varphi)'
   =\langle f_n+z_nu_n,\varphi\rangle
       \qquad(\varphi\in C_c^\infty(0,\infty)),
\]
with $z_n\to z$, $f_n\rightharpoonup f$ in $L^2$, and $u_n\to u$ strongly in $L^2$, then $u_n\to u$ in $H^1(0,L)$ for every finite $L$.
\end{lemma}
\begin{proof}
Choose a real $\chi\in C_c^2([0,\infty))$ equal to one on $[0,L+1]$. Put $h=h_N$, $\chi_k=\chi(kh)$, $F_N=\mathcal I_Nv_N$ and $w_N=\mathcal I_N(\chi v_N)$.
Lemma~\ref{W:lem:spectral-identification} and \eqref{W:eq:grid-values-cutoff-product-error} give
\[
 w_N\rightharpoonup w:=\chi u\text{ in }H^1,
 \qquad w_N\to w\text{ in }L^2
\]
on one fixed interval containing the supports. Sampling the smooth cutoff gives
\[
 \|\mathcal E_N(\chi^2v_N)-\chi^2\mathcal E_Nv_N\|_2
 \le h\|(\chi^2)'\|_\infty\|v_N\|\to0.
\]
Thus $\mathcal E_N(\chi^2v_N)\to\chi^2u$ strongly in $L^2$. For $\chi$ itself,
\[
 \|\mathcal E_N(\chi v_N)-\chi u\|_2
 \le h\|\chi'\|_\infty\|v_N\|
      +\|\chi\|_\infty\|\mathcal E_Nv_N-u\|_2\longrightarrow0.
\]
Since $\mathcal E_N$ is an isometry, $\|\chi v_N\|\to\|\chi u\|_2$. In particular,
\[
 \begin{split}
 &|\langle f_N,\mathcal E_N(\chi^2v_N)\rangle-\langle f,\chi^2u\rangle|\\
 &\quad\le\|f_N\|\,\|\mathcal E_N(\chi^2v_N)-\chi^2u\|_2
                +|\langle f_N-f,\chi^2u\rangle|\longrightarrow0.
 \end{split}
\]
The mass term $\operatorname{Re}z_N\|\chi v_N\|^2$ converges to
$\operatorname{Re}z\|\chi u\|_2^2$.

For its remaining term put
\[
 q_{k,N}=\bar b_{k,N}\left(\frac{\chi_{k+1}-\chi_k}{h}\right)^2,
 \qquad
 \mathcal D_N=\sum_{k=0}^{N-1}q_{k,N}\operatorname{Re}
                   ((v_N)_{k+1}\overline{(v_N)_k}),
\]
using the auxiliary coefficient $\bar b_{0,N}=1$ at zero. Set $q_N^{\rm st}=q_{k,N}$ on $[kh,(k+1)h)$ for $0\le k<N$, and zero after $Nh$. Choose a fixed interval $I=[0,M]$
containing the support of $\chi$ and all cells where $q_{k,N}\ne0$, for all sufficiently large $N$. The mean value formula gives
\[
 \|q_N^{\rm st}-|\chi'|^2\|_\infty
 \le\|\mathfrak b_N-1\|_{L^\infty(I)}\|\chi'\|_\infty^2
            +2h\|\chi'\|_\infty\|\chi''\|_\infty\longrightarrow0.
\]
Only these cells contribute to $\mathcal D_N$. The local derivative bound gives
\[
 \sum_{k:q_{k,N}\ne0}|(v_N)_{k+1}-(v_N)_k|^2
 \le h^2\|F_N'\|_{L^2(I)}^2=O(h^2).
\]
On the $k$th cell $\mathcal E_Nv_N=h^{-1/2}(v_N)_{k+1}$. Therefore
\[
 \begin{split}
 &\left|\mathcal D_N-\int q_N^{\rm st}|\mathcal E_Nv_N|^2\right|\\
 &\quad\le\|q_N^{\rm st}\|_\infty
       \left(\sum_{k:q_{k,N}\ne0}|(v_N)_{k+1}-(v_N)_k|^2\right)^{1/2}\|v_N\|
 =O(h).
 \end{split}
\]
Consequently $\mathcal D_N\to\int|\chi'|^2|u|^2$. The exact cutoff identity is
\[
 \mathfrak q_N[\chi v_N]
 =\operatorname{Re}\langle f_N,\mathcal E_N(\chi^2v_N)\rangle
       +\operatorname{Re}z_N\|\chi v_N\|^2+\mathcal D_N.
\]
The preceding limits and the limiting weak equation therefore yield
\begin{equation}
 \mathfrak q_N[\chi v_N]
 \longrightarrow\operatorname{Re}\langle f+zu,\chi^2u\rangle
                         +\int|\chi'|^2|u|^2
 =\mathfrak q_Y[\chi u].
 \label{W:eq:local-cutoff-energy-limit}
\end{equation}
The last equality uses the compactly supported $H_0^1$ function $\chi^2u$, permitted by the density argument in Lemma~\ref{W:lem:spectral-identification}, and the product identity
$\mathfrak q_Y[\chi u]-\operatorname{Re}\mathfrak q_Y[u,\chi^2u] =\int|\chi'|^2|u|^2$.

Write $Y_N=\widetilde Y_N^{\rm st}$. Since $Y_Nw_N\to Yw$ strongly in $L^2$ and $w_N'\rightharpoonup w'$, the potential terms satisfy
\[
 \int Y_N(|w_N|^2)'
 =2\operatorname{Re}\langle w_N',Y_Nw_N\rangle
 \longrightarrow2\operatorname{Re}\langle w',Yw\rangle.
\]
Adding these terms to \eqref{W:eq:local-cutoff-energy-limit} gives $\int\mathfrak b_N|w_N'|^2\to\|w'\|_2^2$. On a fixed interval $I$ containing the supports,
\[
 \left|\|w_N'\|_2^2-\int\mathfrak b_N|w_N'|^2\right|
 \le\|\mathfrak b_N-1\|_{L^\infty(I)}\|w_N'\|_2^2\longrightarrow0.
\]
Hence $\|w_N'\|_2^2\to\|w'\|_2^2$, and weak convergence gives
\[
 \|w_N'-w'\|_2^2
 =\|w_N'\|_2^2+\|w'\|_2^2
           -2\operatorname{Re}\langle w_N',w'\rangle\longrightarrow0.
\]
Every cell meeting $[0,L]$ has both cutoff values equal to one for large $N$, so $w_N=F_N$ there. This proves \eqref{W:eq:spectral-strong-energy}.

For continuum solutions $u_n$, the same compactly supported test functions are admissible. Young's inequality and the product identity give, with constants depending on $Y$ and
$\chi$ on their fixed support,
\[
 \tfrac12\|(\chi u_n)'\|_2^2-C\|u_n\|_2^2
 \le\mathfrak q_Y[\chi u_n]
 =\operatorname{Re}\langle f_n+z_nu_n,\chi^2u_n\rangle
                         +\int|\chi'|^2|u_n|^2.
\]
Its right side is bounded. Thus $u_n$ is bounded in $H^1$ on each fixed interval, and the strong $L^2$ limit identifies every weak $H^1$ subsequential limit with $u$. The preceding
lemma's weak equation passes to the limit with forcing $f$ and parameter $z$. In particular,
\[
 \begin{aligned}
 \langle f_n+z_nu_n,\chi^2u_n\rangle
     &\longrightarrow\langle f+zu,\chi^2u\rangle,\\
 \int|\chi'|^2|u_n|^2&\longrightarrow\int|\chi'|^2|u|^2,\\
 \int Y(|\chi u_n|^2)'
     &=2\operatorname{Re}\langle(\chi u_n)',Y\chi u_n\rangle
       \longrightarrow2\operatorname{Re}\langle(\chi u)',Y\chi u\rangle.
 \end{aligned}
\]
Add the last line to the quadratic form identity to obtain $\|(\chi u_n)'\|_2^2\to\|(\chi u)'\|_2^2$. Weak convergence then gives
\[
 \|(\chi u_n)'-(\chi u)'\|_2^2
 =\|(\chi u_n)'\|_2^2+\|(\chi u)'\|_2^2
       -2\operatorname{Re}\langle(\chi u_n)',(\chi u)'\rangle
 \longrightarrow0.
\]
Since $\chi=1$ on $[0,L]$, this proves the continuum convergence on $[0,L]$.
\end{proof}

\subsubsection{Ordered eigenvalues and eigenfunctions}
\label{W:sec:ordered-eigenvalues-and-eigenfunctions}

The local weak equation and convergence of the local energies identify the limiting eigenfunctions. The lower bound from
Proposition~\ref{W:prop:spectral-lower-bound-in-probability} also fixes the indices of the eigenvalues: it excludes eigenvalues of $\mathcal A_N$ tending to $-\infty$ below a
convergent bounded spectral cluster.

For the limiting differential operator, simplicity follows directly from the equation with continuous integrated potential.

\begin{lemma}
\label{W:lem:continuous-integrated-potential-eigenfunctions}
Every eigenfunction $u$ of $\mathcal H_Y$ has a $C^1$ representative on $[0,\infty)$. Each eigenvalue is simple, and every nonzero eigenfunction satisfies $u'(0)\ne0$.
\end{lemma}
\begin{proof}
Let $\mathcal H_Yu=\lambda u$. The form domain gives $u\in H^1(0,L)$ and $u(0)=0$ for every finite $L$. Since $Y$ is bounded on $[0,L]$, the function $p_u=u'-Yu$ belongs to
$L^2(0,L)$. For $\psi\in C_c^\infty(0,L)$, the weak equation is
\[
 \int_0^L p_u\overline{\psi'}
 =\int_0^L(Yu'+\lambda u)\overline\psi.
\]
Thus $p_u'=-Yu'-\lambda u\in L^2(0,L)$ in the distributional sense, and $p_u\in H^1(0,L)$. Both $u$ and $p_u$ have absolutely continuous representatives up to zero. The identity
$u'=p_u+Yu$ has a continuous right side, so $u$ has a $C^1$ representative on $[0,L]$. With $U=(u,p_u)^{\mathsf T}$, the equation becomes
\begin{equation}
 U'=M_\lambda U,\qquad
 M_\lambda(x)=\begin{pmatrix}
 Y(x)&1\\-Y(x)^2-\lambda&-Y(x)
 \end{pmatrix},\qquad U(0)=(0,u'(0))^{\mathsf T}.
 \label{W:eq:continuous-integrated-potential-first-order-system}
\end{equation}
This is the system in equation~(2.4) of \cite{eckhardt2013weyl}, with $p=r=1$, $s=-Y$ and $q=-Y^2$. Its uniqueness also follows directly on $[0,L]$. If $V$ is the difference of two
solutions with the same initial value, then
$V(x)=\int_0^xM_\lambda(t)V(t)\,\dd t$.
Set $C_L=\sup_{[0,L]}\|M_\lambda\|$ and $A_L=\sup_{[0,L]}\|V\|$. Iteration of this integral identity gives, for every integer $r\ge1$,
\[
 \|V(x)\|\le A_L\frac{(C_Lx)^r}{r!},\qquad 0\le x\le L.
\]
The right side tends to zero as $r\to\infty$, so $V=0$. Zero initial values therefore give the zero solution. Since $\lambda$ and the matrix coefficients are real, a solution with
real initial values agrees with its complex conjugate and is real. If $u,v$ are nonzero eigenfunctions for the same $\lambda$, the function $u-u'(0)v/v'(0)$ has both initial values
zero and hence vanishes. Thus every eigenspace is one-dimensional. Since $L$ was arbitrary, the regularity holds on the whole half-line.
\end{proof}

Order the eigenvalues of $\mathcal A_N$ increasingly as $\lambda_{j,N}$ and those of $\mathcal H_Y$ as $\Lambda_j$. Choose normalized real limiting eigenfunctions $u_j$; their
signs will be fixed at the boundary below.

\begin{theorem}
\label{W:thm:spectral-deterministic-spectrum}
Under the hypotheses of Theorem~\ref{W:thm:operator-deterministic}, for every fixed $j$ one has $\lambda_{j,N}\to\Lambda_j$. There is an index $N_j$ such that $\lambda_{j,N}$ is
simple for $N\ge N_j$, and corresponding unit eigenvectors can be chosen so that
\begin{equation}
 \mathcal E_Nv_{j,N}\to u_j\ \text{in }L^2(\mathbb R_+),\qquad
 \mathcal I_Nv_{j,N}\to u_j\ \text{in }L^2(\mathbb R_+)
                  \cap H^1_{\mathrm{loc}}([0,\infty)).
 \label{W:eq:spectral-eigenfunction-convergence}
\end{equation}
\end{theorem}
\begin{proof}
By Proposition~\ref{W:prop:spectral-lower-bound-in-probability}, choose $B<\infty$ and $N_0$ with $\mathcal A_N\ge-BI$ for $N\ge N_0$. Choose $s>B+1$ with $s+\Lambda_1>1$; then
$-s$ lies below the finite and limiting spectra for these $N$. Proposition~\ref{W:prop:trace-resolvent-continuation} gives $\mathscr R_N(-s)\to(\mathcal H_Y+s)^{-1}$ in norm. These
are positive compact operators, with nonzero eigenvalues $(\lambda_{j,N}+s)^{-1}$ and $(\Lambda_j+s)^{-1}$ in decreasing order. The min--max inequality is
\[
 |(\lambda_{j,N}+s)^{-1}-(\Lambda_j+s)^{-1}|
 \le\|\mathscr R_N(-s)-(\mathcal H_Y+s)^{-1}\|.
\]
Since $(\Lambda_j+s)^{-1}>0$, inversion proves $\lambda_{j,N}\to\Lambda_j$ for every fixed $j$.

The limiting eigenvalues are simple by Lemma~\ref{W:lem:continuous-integrated-potential-eigenfunctions}. For a fixed $j$, choose an interval isolating $\Lambda_j$.
Proposition~\ref{W:prop:trace-riesz-stability} gives convergence of the projections. Choose $N_j$ so that their operator norm difference is less than one and $\lambda_{j,N}$ lies
in the interval for every $N\ge N_j$. The projection has rank one, so this is its only eigenvalue, with multiplicity one. Put
\[
 u_{j,N}=\frac{P_{j,N}^{\rm eig}u_j}{\|P_{j,N}^{\rm eig}u_j\|},
 \qquad v_{j,N}=\mathcal E_N^*u_{j,N}.
\]
For $\delta_N=\|P_{j,N}^{\rm eig}-P_j^{\rm eig}\|\to0$,
\[
 \|P_{j,N}^{\rm eig}u_j\|\ge1-\delta_N,\qquad
 \|u_{j,N}-u_j\|\le2\delta_N.
\]
Thus $\|v_{j,N}\|=1$, $\mathcal E_Nv_{j,N}=u_{j,N}\to u_j$ in $L^2$, and $\mathcal A_Nv_{j,N}=\lambda_{j,N}v_{j,N}$. Lemma~\ref{W:lem:spectral-energy}, with zero forcing, gives the
local strong $H^1$ convergence.

For the global linear embedding, cellwise convexity gives
\[
 \int_{kh}^{(k+1)h}|\mathcal I_Nv|^2
 \le\tfrac12(|v_k|^2+|v_{k+1}|^2),
 \qquad v_0=v_{N+1}=0.
\]
Summing the cells meeting $[L,\infty)$, for $L>2$ and $h_N<1/2$, yields
\[
 \|\mathbf1_{[L,\infty)}\mathcal I_Nv\|_2^2
 \le C\|Q_{N,L-1}v\|^2.
\]
For $x_N=h_N\lfloor(L-1)/h_N\rfloor\ge L-2$,
\[
 \|Q_{N,L-1}v_{j,N}\|
 =\|\mathbf1_{[x_N,\infty)}\mathcal E_Nv_{j,N}\|_2
 \le\|\mathcal E_Nv_{j,N}-u_j\|_2
                          +\|u_j\|_{L^2(L-2,\infty)}.
\]
Local convergence therefore gives
\[
 \limsup_N\|\mathcal I_Nv_{j,N}-u_j\|_2
 \le C\|u_j\|_{L^2(L-2,\infty)}+\|u_j\|_{L^2(L,\infty)}.
\]
Letting $L\to\infty$ proves \eqref{W:eq:spectral-eigenfunction-convergence}.
\end{proof}

\begin{corollary}
\label{W:cor:finite-spectral-sums}
Under the same hypotheses, for every $g\in C_c(\mathbb R)$, both the scalar spectral sum $\sum_jg(\lambda_{j,N})$ and the operator-valued spectral sum
$\sum_jg(\lambda_{j,N})P_{j,N}^{\rm eig}$ converge, the latter in trace norm, to their limiting counterparts.
\end{corollary}
\begin{proof}
Choose a bounded interval $I$ with $\partial I\subset\rho(\mathcal H_Y)$ and $\operatorname{supp}g\subset I$. Ordered convergence and rank stability give an index $N_I$ such that,
for every $N\ge N_I$, the eigenvalues of $\mathcal A_N$ in $I$ have exactly the index set $J=\{j:\Lambda_j\in I\}$. This set is finite. For $j\in J$,
\[
 g(\lambda_{j,N})\to g(\Lambda_j),\qquad
 \|P_{j,N}^{\rm eig}-P_j^{\rm eig}\|_1
 \le2\|P_{j,N}^{\rm eig}-P_j^{\rm eig}\|\to0,
\]
because the difference has rank at most two. The scalar sums therefore converge, and $\|P_{j,N}^{\rm eig}\|_1=1$ gives
\[
 \left\|\sum_jg(\lambda_{j,N})P_{j,N}^{\rm eig}
              -\sum_jg(\Lambda_j)P_j^{\rm eig}\right\|_1
 \le\sum_{j\in J}|g(\lambda_{j,N})-g(\Lambda_j)|
 +\sum_{j\in J}|g(\Lambda_j)|\,
                   \|P_{j,N}^{\rm eig}-P_j^{\rm eig}\|_1\to0.
\]
When $J$ is empty, both sums vanish for every $N\ge N_I$.
\end{proof}

\subsubsection{Boundary derivatives of eigenfunctions}
\label{W:sec:boundary-derivatives-of-eigenfunctions}

We next recover spectral weights in the original matrix coordinates from the boundary derivatives of the limiting eigenfunctions. Let $\xi_{j,N}$ be a unit eigenvector of $H_N$ for
its $j$th largest eigenvalue $\mu_j^{(N)}$, and put $v_{j,N}=U_N(w_N)^*\xi_{j,N}$. Since the first Lanczos vector is $w_N$,
\begin{equation}
 w_N^*\xi_{j,N}=(v_{j,N})_1,\qquad
 (\mathcal I_Nv_{j,N})'(0+)
       =h_N^{-3/2}(v_{j,N})_1
       =\sqrt N\,w_N^*\xi_{j,N}.
 \label{W:eq:ev-original-boundary-identity}
\end{equation}
The boundary derivative is not continuous in the $H^1$ topology. Indeed, choose $\phi\in C_c^\infty([0,\infty))$ with $\phi(0)=0$ and $\phi'(0)=1$, and set
$f_\varepsilon(x)=\varepsilon\phi(x/\varepsilon)$. Then
\begin{equation}
 \|f_\varepsilon\|_{H^1}^2
 =\varepsilon^3\|\phi\|_2^2+\varepsilon\|\phi'\|_2^2\longrightarrow0,
 \qquad f_\varepsilon'(0)=1.
 \label{W:eq:boundary-derivative-not-H1-continuous}
\end{equation}
We therefore use the Jacobi recurrence to control the discrete counterpart of $u'-Yu$. Its increments involve the local energy and the integrated potential, which allows us to pass
to the value at zero.

\Needspace{7\baselineskip} Let $T_N$ be real symmetric tridiagonal matrices with the conventions of Section~\ref{W:sec:local-energy-estimates-and-ordered-eigenfunctions}. Suppose
$\widetilde Y_N^{\rm lin}\to Y$ locally uniformly, where $Y$ is real and continuous, and $B_N(L)\to0$ for every $L<\infty$. Let $\lambda_N\to\lambda\in\mathbb R$ and
\[
 \mathcal A_Nv_N=\lambda_Nv_N,\qquad
 \sup_N\|v_N\|<\infty,\qquad
 \mathcal E_Nv_N\to u\quad\hbox{in }L^2(\mathbb R_+).
\]

\begin{lemma}
\label{W:lem:ev-boundary-derivative}
The limiting function $u$ belongs to $C^1$ locally, vanishes at zero, and
\begin{equation}
 h_N^{-3/2}(v_N)_1\longrightarrow u'(0).
 \label{W:eq:ev-boundary-derivative-limit}
\end{equation}
If $u\ne0$, then $u'(0)\ne0$.
\end{lemma}
\begin{proof}
Write $h=h_N$, $Y_k=\widetilde Y_{N,k}$, $u_k=h^{-1/2}(v_N)_k$, and $u_0=Y_0=0$. Use $b_0=1$ for the auxiliary coefficient and $b_k=b_{k,N}$ for $1\le k<N$. For a fixed interval
$[0,L]$, choose $N_L$ so that $\lceil L/h_N\rceil+1<N$ for $N\ge N_L$; this is possible since $h_N=N^{-1/3}$. Define
\begin{equation}
 d_k=b_k\frac{u_{k+1}-u_k}{h},\qquad
 p_k=d_k-Y_ku_k,\qquad k\ge0.
 \label{W:eq:ev-discrete-flux}
\end{equation}
Then $p_0=u_1/h=h^{-3/2}(v_N)_1$. Rearranging $(2-a_{k,N})u_k-b_{k-1}u_{k-1}-b_ku_{k+1}=h^2\lambda_Nu_k$ and using $Y_k-Y_{k-1}=h^{-1}(2-a_{k,N}-b_{k-1}-b_k)$ gives
\begin{align}
 d_k-d_{k-1}
   &=(Y_k-Y_{k-1})u_k-h\lambda_Nu_k,
       &&k\ge1,\label{W:eq:ev-flux-recurrence}\\
 p_k-p_{k-1}
   &=-Y_{k-1}(u_k-u_{k-1})-h\lambda_Nu_k.
       \label{W:eq:ev-quasiderivative-recurrence}
\end{align}
Both identities include the first cell, with $b_0=1$.

For $y=Y_{k-1}$, $b=b_{k-1}$ and $\lambda=\lambda_N$, the exact transfer matrix is
\[
 \binom{u_k}{p_k}=A_{k,N}\binom{u_{k-1}}{p_{k-1}},\qquad
 A_{k,N}
 =\begin{pmatrix}1&0\\-y-h\lambda&1\end{pmatrix}
  \begin{pmatrix}1&h/b\\0&1\end{pmatrix}
  \begin{pmatrix}1&0\\y&1\end{pmatrix}.
\]
Multiplication gives
\[
 A_{k,N}=\begin{pmatrix}
  1+hy/b&h/b\\
  -h\lambda-hy^2/b-h^2\lambda y/b&1-hy/b-h^2\lambda/b
 \end{pmatrix}.
\]
Thus $\det A_{k,N}=1$, and local boundedness of $Y_k,\lambda_N,b_k^{-1}$ gives $\|A_{k,N}-I\|\le C_Lh$ and $\|A_{k,N}\cdots A_{1,N}\|\le e^{C_Lkh}$ for $kh\le L$. The factorization
in Figure~\ref{W:fig:boundary-transfer-factors} uses values of $Y$, without its adjacent differences.
\begin{figure}[H]
\centering
\begin{tikzpicture}[x=1cm,y=1cm,>=stealth,
 every node/.style={font=\small},every path/.style={line width=.5pt}]
 \node[draw,minimum width=2.05cm,minimum height=1.12cm] (a) at (0,0)
    {$\dbinom{u_{k-1}}{p_{k-1}}$};
 \node[draw,minimum width=2.05cm,minimum height=1.12cm] (b) at (3.55,0)
    {$\dbinom{u_{k-1}}{d_{k-1}}$};
 \node[draw,minimum width=2.05cm,minimum height=1.12cm] (c) at (7.1,0)
    {$\dbinom{u_k}{d_{k-1}}$};
 \node[draw,minimum width=2.05cm,minimum height=1.12cm] (d) at (10.65,0)
    {$\dbinom{u_k}{p_k}$};
 \draw[->] (a)--node[above=18pt] {$p\mapsto p+yu$} (b);
 \draw[->] (b)--node[above=18pt] {$u\mapsto u+hd/b$} (c);
 \draw[->] (c)--node[above=18pt] {$d\mapsto d-(y+h\lambda)u$} (d);
 \node[below=10pt,align=center] at (3.55,-.55)
    {$d_{k-1}=p_{k-1}+yu_{k-1}$};
 \node[below=10pt,align=center] at (9.0,-.55)
    {$p_k=d_{k-1}-(y+h\lambda)u_k$};
\end{tikzpicture}
\caption{The three changes of variables whose product is the recurrence matrix $A_{k,N}$, with $y=Y_{k-1}$, $b=b_{k-1}$, and $\lambda=\lambda_N$. They give $p_k-p_{k-1}=-y(u_k-u_{k-1})-h\lambda u_k$. At the first cell, $u_0=Y_0=0$ and $b_0=1$, so the initial value is $p_0=u_1/h=h^{-3/2}(v_N)_1$.}
\label{W:fig:boundary-transfer-factors}
\end{figure}
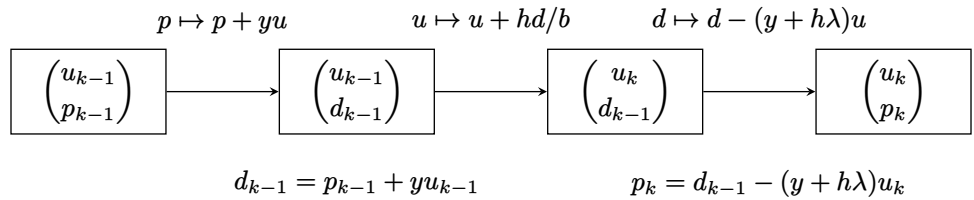

Put $f_N=\mathcal I_Nv_N$. Lemmas~\ref{W:lem:spectral-identification} and~\ref{W:lem:spectral-energy}, with zero forcing, give
\begin{equation}
 f_N\to u\quad\hbox{in }H^1(0,L+2)
 \quad\hbox{for every fixed }L>0.
 \label{W:eq:ev-local-energy-input}
\end{equation}
Also $Y_N^{\rm st}:=\widetilde Y_N^{\rm st}\to Y$ uniformly on each fixed interval, since the adjacent jumps are bounded by
$2\|\widetilde Y_N^{\rm lin}-Y\|_\infty+\omega_Y(h;L+3)$.

Let $p_N$ interpolate the $p_k$ at $kh$, and let $p_N^{\rm st}=p_k$ on $[kh,(k+1)h)$. On that cell write $f_N^-=u_k$ and $f_N^+=u_{k+1}$. Linear interpolation on that interval
gives
\begin{equation}
 \|f_N^\pm-f_N\|_{L^2(0,L)}
       \le h\|f_N'\|_{L^2(0,L+1)}.
 \label{W:eq:eigenfunction-cell-sampling}
\end{equation}
Indeed the squared cell error is $h|u_{k+1}-u_k|^2/3$; the enlarged interval includes a last incomplete cell. The two recurrences give the identities
\[
 p_N^{\rm st}=\mathfrak b_N f_N'-Y_N^{\rm st}f_N^-,\qquad
 p_N'=-Y_N^{\rm st}f_N'-\lambda_Nf_N^+.
\]
The first, together with \eqref{W:eq:ev-local-energy-input} and \eqref{W:eq:eigenfunction-cell-sampling}, proves
\begin{equation}
 \sup_N\|p_N^{\rm st}\|_{L^2(0,L+1)}<\infty,
 \qquad p_N^{\rm st}\to p:=u'-Yu\quad\hbox{in }L^2(0,L).
 \label{W:eq:ev-step-flux-limit}
\end{equation}
The second gives
\begin{equation}
 \|p_N'\|_{L^2(0,L+1)}
 \le C_L\|f_N'\|_{L^2(0,L+2)}
   +|\lambda_N|\bigl(\|f_N\|_{L^2(0,L+2)}
                   +h\|f_N'\|_{L^2(0,L+2)}\bigr)
 \le C_L,
 \label{W:eq:ev-flux-H1-bound}
\end{equation}
and, by the same convergences, $p_N'\to-Yu'-\lambda u$ in $L^2(0,L)$. Since $\|p_N-p_N^{\rm st}\|_{L^2(0,L)} \le h\|p_N'\|_{L^2(0,L+1)}$, we obtain
\[
 p_N\to p\text{ in }H^1(0,L),\qquad p'=-Yu'-\lambda u.
\]
Here the derivative identity follows by passing to the distributional derivative against smooth compactly supported test functions. The Sobolev embedding on $(0,L)$ also gives
uniform convergence of $p_N$.

Because $p,Y,u$ are continuous, $u'=p+Yu$ is continuous, so $u\in C^1([0,L])$. The Dirichlet trace is $u(0)=0$, hence $p(0)=u'(0)$. For $g\in H^1(0,L)$, integrating
$g(x)=g(0)+\int_0^xg'(t)\,\dd t$ gives
\[
 g(0)=\frac1L\int_0^Lg(x)\,\dd x
       -\frac1L\int_0^L(L-t)g'(t)\,\dd t.
\]
Cauchy--Schwarz therefore yields
\[
 |h^{-3/2}(v_N)_1-u'(0)|
 =|p_N(0)-p(0)|
 \le L^{-1/2}\|p_N-p\|_2+\sqrt{L/3}\|p_N'-p'\|_2\to0.
\]
The two recurrence identities also give explicit error bounds. On $(0,L)$,
\[
 \begin{split}
 \|p_N^{\rm st}-p\|_2
 &\le\|\mathfrak b_N\|_\infty\|f_N'-u'\|_2
       +\|\mathfrak b_N-1\|_\infty\|u'\|_2\\
 &\quad+\|Y_N^{\rm st}\|_\infty\|f_N^--u\|_2
       +\|Y_N^{\rm st}-Y\|_\infty\|u\|_2,\\
 \|p_N'-p'\|_2
 &\le\|Y_N^{\rm st}\|_\infty\|f_N'-u'\|_2
       +\|Y_N^{\rm st}-Y\|_\infty\|u'\|_2\\
 &\quad+|\lambda_N|\|f_N^+-u\|_2
       +|\lambda_N-\lambda|\|u\|_2.
 \end{split}
\]
Using \eqref{W:eq:eigenfunction-cell-sampling} and $\|p_N-p_N^{\rm st}\|_2\le h\|p_N'\|_{L^2(0,L+1)}$ in the trace inequality gives, for all large $N$,
\[
 \begin{split}
 |h^{-3/2}(v_N)_1-u'(0)|
 \le C_L\bigl(&\|f_N-u\|_{H^1(0,L+2)}
       +\|\mathfrak b_N-1\|_{L^\infty(0,L+2)}\\
       &+\|Y_N^{\rm st}-Y\|_{L^\infty(0,L+2)}
       +|\lambda_N-\lambda|+h\bigr).
 \end{split}
\]
Here $C_L$ depends on the established local bounds and on $u,Y,\lambda$, but not on $N$. This proves \eqref{W:eq:ev-boundary-derivative-limit}.

The cellwise derivatives also converge in $L^\infty(0,L)$. Indeed $\|p_N-p_N^{\rm st}\|_\infty\le\sqrt h\|p_N'\|_2$, and \eqref{W:eq:ev-local-energy-input} gives
$\|f_N^--u\|_\infty\le\|f_N-u\|_\infty+\omega_u(h;L+1)\to0$. Thus the cell identity $f_N'=\mathfrak b_N^{-1}(p_N^{\rm st}+Y_N^{\rm st}f_N^-)$, with $\mathfrak b_N^{-1}\to1$
uniformly, proves $\|f_N'-u'\|_{L^\infty(0,L)}\to0$. Finally $u'=p+Yu$, $p'=-Yp-(Y^2+\lambda)u$ is the continuous coefficient system
\eqref{W:eq:continuous-integrated-potential-first-order-system}. Its uniqueness, used in Lemma~\ref{W:lem:continuous-integrated-potential-eigenfunctions}, shows that $u'(0)=0$
would force $u=p=0$. Hence a nonzero limit has $u'(0)\ne0$.
\end{proof}

\subsubsection{Joint convergence of spectral data}
\label{W:sec:joint-convergence-of-spectral-data}

The preceding estimates compare eigenvalues, projections and boundary derivatives on a common sequence of matrices. To obtain their joint distributional limit, we apply those
deterministic estimates to copies of the pair $(\widetilde Y_N^{\rm lin},\mathscr R_N)$, together with the local coefficient bounds and a lower bound for the spectrum. Measurable
choices of eigenvectors then allow the pathwise conclusion to be transferred to the original distributions.

Fix a deterministic sequence of unit vectors $w_N$. Choose a measurable ordered orthonormal eigenbasis $\xi_{j,N}$ of $H_N$ and put
\[
 \lambda_{j,N}=N^{2/3}(2-\mu_j^{(N)}),\quad
 v_{j,N}=U_N(w_N)^*\xi_{j,N},\quad
 P_{j,N}^{\mathrm{eig}}=|\mathcal E_Nv_{j,N}\rangle
                              \langle\mathcal E_Nv_{j,N}|.
\]
On repeated eigenspaces any fixed measurable orthonormalization rule is used. Let $\varphi_j$ be the real normalized eigenfunction of $\mathcal H_\beta$ with $\varphi_j'(0)>0$, and
define
\begin{equation}
 W_j=|\varphi_j'(0)|^2.
 \label{W:eq:ev-sao-boundary-weight}
\end{equation}
The existence and uniqueness of this sign convention follow from Lemma~\ref{W:lem:continuous-integrated-potential-eigenfunctions}. For fixed $k$, set
\[
 \mathscr D_N^{(k)}=
    (\lambda_{j,N},P_{j,N}^{\mathrm{eig}},N|w_N^*\xi_{j,N}|^2)_{j\le k},
 \qquad
 \mathscr D_\beta^{(k)}=
    (\Lambda_j,|\varphi_j\rangle\langle\varphi_j|,W_j)_{j\le k}.
\]
The scalar coordinates have their usual topology, and the spectral projections have the operator norm topology in $\mathcal K(L^2(\mathbb R_+))$. For functions we use
$\mathcal X=L^2(\mathbb R_+)\cap H^1_{\mathrm{loc}}([0,\infty))$, equipped with
\[
 d_{\mathcal X}(f,g)=(1\wedge\|f-g\|_2)
       +\sum_{\ell\ge1}2^{-\ell}(1\wedge\|f-g\|_{H^1(0,\ell)}).
\]
The map $f\mapsto(f,(f|_{(0,\ell)})_{\ell\ge1})$ identifies this space with the subset of $L^2(\mathbb R_+)\times\prod_{\ell\ge1}H^1(0,\ell)$ defined by $g_\ell=f|_{(0,\ell)}$ in
$L^2(0,\ell)$ for every $\ell$. Each equality defines a closed set, since restriction on $L^2$ and the map $H^1(0,\ell)\to L^2(0,\ell)$ are continuous. The image is therefore
closed. The countable product is separable and complete, so $\mathcal X$ is a Polish space with exactly the stated global $L^2$ and local $H^1$ convergence.

\begin{proposition}
\label{W:prop:joint-spectral-data}
Under Assumption~\ref{W:ass:wigner-lee-yin}, $\mathscr D_N^{(k)}$ converges in distribution to $\mathscr D_\beta^{(k)}$, jointly with the integrated potential and the nonreal
resolvent path. After normalization to $(v_{j,N})_1>0$, the eigenfunctions $\mathcal I_Nv_{j,N}$ converge jointly in $\mathcal X$ to $\varphi_j$, $1\le j\le k$.
\end{proposition}
The normalization is made when the eigenvalue is simple and its first coordinate is nonzero; any measurable convention is used otherwise. The proof below shows that the latter
event has probability tending to zero.

\begin{proof}
Write $X_N=(\widetilde Y_N^{\rm lin},\mathscr R_N)$ and $b_N=(B_N(\ell))_{\ell\ge1}$. On $[0,\infty)^{\mathbb N}$ use the metric
\[
 d_\Pi(b,c)=\sum_{\ell\ge1}2^{-\ell}(1\wedge|b_\ell-c_\ell|).
\]
For every integer $L$ and every $\epsilon>0$,
\[
 \mathbb P(d_\Pi(b_N,0)>\epsilon+2^{-L})
 \le\sum_{\ell=1}^L\mathbb P(B_N(\ell)>\epsilon)\longrightarrow0.
\]
Letting $L\to\infty$ proves $b_N\to0$ in probability. Set $\kappa_N=(-\lambda_{1,N})_+$. Proposition~\ref{W:prop:spectral-lower-bound-in-probability} gives tightness of
$(\kappa_N)$. Its bound for all sufficiently large $N$ extends to all $N$ by enlarging the compact interval to cover the finitely many remaining distributions.
Equation~\eqref{W:eq:assembly-joint-convergence} gives tightness of $(X_N)$ and identifies its limit as
$X_\beta=(Y_\beta,\mathscr R_\beta)$.

These statements give tightness of the full tuple $(X_N,b_N,\kappa_N)$. Indeed, for $\epsilon>0$ choose compact sets $K_X,K_b,K_\kappa$ for the three marginals, each with
complementary probability at most $\epsilon/3$ uniformly in $N$. Then
\[
 \sup_N\mathbb P\bigl((X_N,b_N,\kappa_N)
       \notin K_X\times K_b\times K_\kappa\bigr)\le\epsilon.
\]
The spaces involved are separable and complete; for the resolvent coordinate this is Lemma~\ref{W:lem:resolvent-path-space}. Every subsequence therefore has a further subsequence
converging weakly to $(X_\beta,0,\kappa)$ for some finite nonnegative $\kappa$. The tightness and representation theorems for Polish spaces \cite{kallenberg2021foundations} provide
copies of this entire tuple such that
\[
 \widetilde Y_N^{\rm lin}\to Y_\beta\text{ locally uniformly},\qquad
 \mathscr R_N\to\mathscr R_\beta\text{ in }\mathscr Y,\qquad
 d_\Pi(b_N,0)\to0,\qquad \kappa_N\to\kappa
\]
almost surely along that further subsequence. The marginal distribution of $(Y_\beta,\mathscr R_\beta)$ is the one in \eqref{W:eq:assembly-joint-convergence}; in particular
$\mathscr R_\beta(z)=(\mathcal H_{Y_\beta}-z)^{-1}$.

The finite matrices can be recovered from the copied resolvents. For each $N$ put
\[
 C_N=\mathcal E_N^*\mathscr R_N(\ii)\mathcal E_N,\qquad
 \mathcal A_N=\ii I+C_N^{-1},\qquad T_N=2I-h_N^2\mathcal A_N.
\]
Evaluation at $\ii$ and compression by $\mathcal E_N$ are continuous, and inversion is continuous on the invertible finite matrices. Define the reconstruction to be zero when $C_N$
is singular. Under the original distribution, $C_N=(\mathcal A_N-\ii)^{-1}$, so this exceptional definition is used with probability zero. Real symmetry and tridiagonality of
$T_N$, nonnegativity of its adjacent coefficients, and the identities defining $\widetilde Y_N^{\rm lin}$, $B_N(\ell)$ and $\kappa_N$ are Borel conditions on the tuple. Equality of
its original and copied distributions therefore gives all these conditions with probability one for each $N$. Intersecting these events over $N$ and over integer $\ell$ gives one
event on which the reconstruction and the convergences above all hold. Below we work on this event. The limit has the form, domain, and simple discrete spectrum of
$\mathcal H_\beta$.

We first recover the order and multiplicity of the eigenvalues. Choose
\[
 s>1+\max\{\kappa,-\Lambda_1\}.
\]
There is an index $N_0$, depending on the chosen outcome and $s$, such that $\kappa_N\le\kappa+1$ and $\mathcal A_N+s>0$ for $N\ge N_0$.
Proposition~\ref{W:prop:trace-resolvent-continuation}, applied on this outcome at $z=-s$, gives
\[
 S_N:=\mathcal E_N(\mathcal A_N+s)^{-1}\mathcal E_N^*
 \longrightarrow S:=(\mathcal H_\beta+s)^{-1}
 \quad\text{in operator norm}.
\]
Both operators are compact and nonnegative. Their positive eigenvalues, in decreasing order, are $r_{j,N}=(\lambda_{j,N}+s)^{-1}$ for $j\le N$ and $r_j=(\Lambda_j+s)^{-1}$ for
$j\ge1$. Append zero eigenvalues to $S_N$. The variational characterization gives
\[
 |r_{j,N}-r_j|\le\eta_N:=\|S_N-S\|.
\]
For fixed $j$, once $N\ge j$ and $\eta_N<r_j/2$, it follows that
\[
 r_{j,N}\ge r_j/2,\qquad
 |\lambda_{j,N}-\Lambda_j|
 =\frac{|r_{j,N}-r_j|}{r_{j,N}r_j}
 \le\frac{2\eta_N}{r_j^2}\longrightarrow0.
\]
The common lower bound makes the shifted inverses positive and identifies each finite cluster with the corresponding index in the increasing spectrum.

For each $j$ choose $a_j<\Lambda_j<b_j$ in the limiting resolvent set with $[a_j,b_j]\cap\sigma(\mathcal H_\beta)=\{\Lambda_j\}$. Let $\Gamma_j$ be a rectangle enclosing this
interval, with no real points on its boundary other than $a_j,b_j$. Proposition~\ref{W:prop:trace-riesz-stability} gives, for all $N$ above an index depending on $j$ and the
outcome,
\[
 \begin{gathered}
 \Pi_{j,N}=\mathcal E_N\mathbf1_{(a_j,b_j)}(\mathcal A_N)\mathcal E_N^*,
 \qquad \Pi_j=|\varphi_j\rangle\langle\varphi_j|,\\
 \|\Pi_{j,N}-\Pi_j\|
 \le\frac{|\Gamma_j|}{2\pi}
       \sup_{z\in\Gamma_j}\|\mathscr R_N(z)-\mathscr R_\beta(z)\|
 \longrightarrow0.
 \end{gathered}
\]
When this norm is less than one the two projections have the same rank, namely one. The ordered convergence places $\lambda_{j,N}$ in $(a_j,b_j)$ for all sufficiently large $N$. It
is consequently simple and $P_{j,N}^{\mathrm{eig}}=\Pi_{j,N}$ at those sizes.

Put $\delta_{j,N}=\|\Pi_{j,N}-\Pi_j\|$. When $\delta_{j,N}<1$, define
\[
 u_{j,N}=\frac{\Pi_{j,N}\varphi_j}{\|\Pi_{j,N}\varphi_j\|},\qquad
 \widetilde v_{j,N}=\mathcal E_N^*u_{j,N}.
\]
The estimates
\[
 \|\Pi_{j,N}\varphi_j\|\ge1-\delta_{j,N},\qquad
 \|u_{j,N}-\varphi_j\|_2\le2\delta_{j,N}
\]
show that $\|\widetilde v_{j,N}\|=1$ and $\mathcal E_N\widetilde v_{j,N}\to\varphi_j$ in $L^2$. The eigenvalue equation has zero forcing. The local convergence of $Y_N$ and of the
adjacent coefficients, together with Lemma~\ref{W:lem:spectral-energy}, therefore gives
\[
 \mathcal I_N\widetilde v_{j,N}\to\varphi_j
       \quad\text{in }H^1(0,L)\quad\text{for every }L<\infty.
\]
For the global norm, the interpolation estimate in the proof of Theorem~\ref{W:thm:spectral-deterministic-spectrum} gives, for $L>2$ and
$h_N<1/2$,
\[
 \|\mathbf1_{[L,\infty)}\mathcal I_N\widetilde v_{j,N}\|_2
 \le C\|Q_{N,L-1}\widetilde v_{j,N}\|
 \le C\bigl(\|u_{j,N}-\varphi_j\|_2
              +\|\varphi_j\|_{L^2(L-2,\infty)}\bigr).
\]
Local convergence, followed by $L\to\infty$, proves convergence in $L^2$ on the half-line, and hence in $\mathcal X$.

Lemma~\ref{W:lem:ev-boundary-derivative} applies to these real eigenvectors and gives
\[
 h_N^{-3/2}(\widetilde v_{j,N})_1\longrightarrow\varphi_j'(0)>0.
\]
Thus there is a finite index $N_j$ on the chosen outcome such that the eigenvalue is simple and $(\widetilde v_{j,N})_1>0$ for every $N\ge N_j$. At those sizes this is the unique
real unit eigenvector with positive first coordinate. Since $h_N^{-3}=N$, identity~\eqref{W:eq:ev-original-boundary-identity} now yields, for every fixed $r,k$,
\begin{equation}
 \left((\lambda_{j,N})_{j=1}^r,
       (h_N^{-3}|(\widetilde v_{j,N})_1|^2)_{j=1}^k\right)
 \longrightarrow
 \left((\Lambda_j)_{j=1}^r,(|\varphi_j'(0)|^2)_{j=1}^k\right).
 \label{W:eq:ev-coupled-boundary-weights}
\end{equation}
Intersecting the events of probability one used above over $j\ge1$ gives one event on which all the preceding limits hold. For any fixed $k$, the eigenvalues are simple and the
first coordinates of the chosen vectors are positive whenever $N\ge\max_{j\le k}N_j$.

It remains to express the eigenvectors and their limiting normalizations as measurable functions of the spectral data. At finite $N$, write $a_{j,N}$ for the $j$th ordered
eigenvalue of the real symmetric matrix $\mathcal A_N$. The eigenvalues are continuous functions of its entries, and the spectral projection onto the entire corresponding
eigenspace satisfies
\[
 P_{j,N}=\lim_{m\to\infty}
       \bigl[I+m(\mathcal A_N-a_{j,N}I)^2\bigr]^{-1}.
\]
Indeed, the expression on the right acts on an eigenvector with eigenvalue $b$ as multiplication by $(1+m(b-a_{j,N})^2)^{-1}$, which tends to $\mathbf1_{\{b=a_{j,N}\}}$. Thus
$P_{j,N}$ is Borel. To select a basis when its rank exceeds one, apply Gram--Schmidt to $P_{j,N}e_1,\ldots,P_{j,N}e_N$ in that order, omitting zero remainders. If $r_l$ is the
remainder at index $l$, its selection is determined by the Borel event
\[
 \{\|r_l\|>0\}\cap\bigcap_{i<l}\{\|r_i\|=0\}.
\]
On $\{\operatorname{rank}P_{j,N}=1,\ e_1^*P_{j,N}e_1>0\}$, the choice is simply
\[
 v_{j,N}^{\mathrm{can}}=
     \frac{P_{j,N}e_1}{(e_1^*P_{j,N}e_1)^{1/2}},
 \qquad e_1^*v_{j,N}^{\mathrm{can}}=(e_1^*P_{j,N}e_1)^{1/2}>0.
\]
The remaining cases are covered by the preceding basis selection. Hence the finite eigenfunctions, their first coordinates, and the associated projections are Borel; their linear
interpolations are Borel as $\mathcal X$-valued functions because $\mathcal I_N:\mathbb R^N\to\mathcal X$ is continuous.

For the limiting operator, the proof of Lemma~\ref{W:lem:sao-fixed-energy-atomlessness} already constructs a fixed real orthonormal family $(e_m)_{m\ge1}$ in the common form core.
Let $\lambda_j^{(m)}(Y)$ be the $j$th ordered eigenvalue of the finite matrix
\[
 Q_m(Y)=(\mathfrak q_Y[e_a,e_b])_{a,b\le m}.
\]
That proof establishes continuity of $Q_m(Y)$ under locally uniform convergence of $Y$ and, by min--max and density of the form core,
\[
 \lambda_j^{(m)}(Y)\downarrow\Lambda_j(Y)\quad(m\to\infty).
\]
Thus $\Lambda_j(Y)=\inf_{m\ge j}\lambda_j^{(m)}(Y)$ is Borel. Because the limiting spectrum is simple and discrete, we may enumerate rational pairs $(a,b)$ and choose the first
satisfying
\[
 a<\Lambda_1<b<\Lambda_2 \quad(j=1),\qquad
 \Lambda_{j-1}<a<\Lambda_j<b<\Lambda_{j+1}\quad(j\ge2).
\]
These conditions are Borel and such a pair exists because the spectrum is simple and discrete. Let $\Gamma(a,b)$ be the positively oriented boundary of the rectangle with vertices
$a\pm\ii,b\pm\ii$. With $T=\mathscr R_\beta(\ii)$, the resolvent identity gives
\[
 \mathscr R_\beta(z)=T[I-(z-\ii)T]^{-1},\qquad z\in\Gamma(a,b).
\]
The inverse exists there because $\Gamma(a,b)$ is in the resolvent set. In particular,
\[
 \Pi_j=-\frac{1}{2\pi\ii}\int_{\Gamma(a,b)}
               T[I-(z-\ii)T]^{-1}\,\dd z.
\]
For a fixed rational pair the integrand depends continuously in operator norm on $T$ wherever its inverse exists, and its Riemann sums converge in operator norm to the displayed
integral. The countable choice of the pair therefore makes $\Pi_j$ Borel.

Choose a fixed real orthonormal basis $(e_m)_{m\ge1}$ of $L^2(\mathbb R_+)$. Define
\[
 m_j=\min\{m:\|\Pi_je_m\|>0\},\qquad
 f_j=\frac{\Pi_je_{m_j}}{\|\Pi_je_{m_j}\|}.
\]
The rank-one property guarantees that $m_j$ is finite. Both $m_j$ and $f_j$ are Borel, and $f_j$ is real because the limiting operator has real coefficients. Its $C^1$
representative satisfies $f_j(0)=0$ and $f_j'(0)\ne0$. For $n\ge1$ put
\[
 d_{j,n}=2n^2\int_0^{1/n}f_j(x)\,\dd x.
\]
Since $f_j(x)=xf_j'(0)+\int_0^x(f_j'(t)-f_j'(0))\,\dd t$,
\[
 |d_{j,n}-f_j'(0)|\le
    \sup_{0\le t\le1/n}|f_j'(t)-f_j'(0)|\longrightarrow0.
\]
Each $d_{j,n}$ is a continuous linear functional of the $L^2$ vector $f_j$, so $d_j=f_j'(0)$ is Borel. Set $\varphi_j=\operatorname{sgn}(d_j)f_j$; then $\varphi_j'(0)=|d_j|>0$ and
$W_j=d_j^2$ are Borel functions of the limiting resolvent and potential.

The $\mathcal X$-valued measurability can also be checked from $L^2$ data. For any integer $L\ge1$ and $0<h<1$ the map
\[
 D_{L,h}f=\frac{f(\,\cdot+h)-f}{h}\bigg|_{(0,L)}
\]
is bounded from $L^2(0,L+1)$ to $L^2(0,L)$. If $f\in H^1(0,L+1)$, the fundamental theorem for Sobolev functions and Minkowski's inequality give
\[
 \|D_{L,h}f-f'\|_{L^2(0,L)}
 \le h^{-1}\int_0^h
   \|f'(\,\cdot+t)-f'\|_{L^2(0,L)}\,\dd t\longrightarrow0.
\]
Thus $f\mapsto f'|_{(0,L)}$ is a pointwise limit of Borel maps on the functions under consideration. Together with the restrictions $f|_{(0,L)}$, these coordinates define the
closed subspace used above to realize $\mathcal X$ as a Polish space. Consequently $\varphi_j$ is Borel in $\mathcal X$.

For fixed $k$, let $E_{N,k}$ be the event that one of the first $k$ eigenspaces of $\mathcal A_N$ has dimension greater than one or is orthogonal to $e_1$. It is Borel by the
finite projection construction. On the coupled space the previously established convergence gives $\mathbf1_{E_{N,k}}\to0$ almost surely; hence
\[
 \mathbb P(E_{N,k})=\mathbb E\mathbf1_{E_{N,k}}\longrightarrow0.
\]
Outside $E_{N,k}$, the unique normalization with positive first coordinate agrees with $v_{j,N}^{\mathrm{can}}$ for every $j\le k$. Let $D_{N,k}$ denote the tuple in the
proposition using the stated convention, and $D_{N,k}^{\mathrm{can}}$ the same tuple using the Borel choices. For each bounded Borel function $F$ on the corresponding product
space,
\[
 |\mathbb EF(D_{N,k})-\mathbb EF(D_{N,k}^{\mathrm{can}})|
 \le2\|F\|_\infty\mathbb P(E_{N,k}).
\]
The copied canonical tuple has the original distribution and converges almost surely along the selected subsequence. Dominated convergence therefore gives, for every bounded
continuous $F$,
\[
 \mathbb EF(D_{N,k}^{\mathrm{can}})\longrightarrow
      \mathbb EF(D_\beta^{(k)})
\]
along that subsequence, jointly with $\widetilde Y_N^{\rm lin}$, $\mathscr R_N$, and the normalized eigenfunctions. The preceding estimate transfers this limit to $D_{N,k}$. Since
every subsequence admits a further subsequence with the same limiting distribution, the entire sequence converges as stated.
\end{proof}

\begin{proof}[Proof of Corollary~\ref{W:cor:intro-ordered-edge-spectrum}]
Fix $k$. Proposition~\ref{W:prop:joint-spectral-data} gives the joint convergence of the first $k$ increasing eigenvalues of $\mathcal A_N$ and their normalized embedded
eigenfunctions, together with $\widetilde Y_N^{\rm lin}$ and $\mathscr R_N$. The common lower spectral bound in that proposition identifies these as the first $k$ eigenvalues,
rather than merely $k$ eigenvalues in prescribed bounded intervals.

For unit vectors $u,v$ in the common Hilbert space, write $P_u$ and $P_v$ for their orthogonal rank-one projections. Expanding their difference gives
\[
 P_u-P_v=(u-v)u^*+v(u-v)^*,\qquad
 \|P_u-P_v\|_{\mathrm{op}}\le2\|u-v\|.
\]
Thus the map from a unit vector to its projection is continuous and independent of its phase. Applying this map to the first $k$ eigenfunction coordinates, and retaining the
eigenvalue coordinates, proves joint convergence of the ordered eigenvalues and projections. On the event where a finite eigenvalue is not simple, use the measurable convention
already fixed in the proposition. The probability of this event for the first $k$ labels tends to zero, so it does not change the limit.

Finally the completed Lanczos transformation is unitary and the edge scaling reverses the ordering:
\[
 \lambda_{j,N}=N^{2/3}(2-\mu_j^{(N)}),\qquad
 N^{2/3}(\mu_j^{(N)}-2)=-\lambda_{j,N},\qquad 1\le j\le k.
\]
Apply the continuous map $(x_1,\ldots,x_k)\mapsto(-x_1,\ldots,-x_k)$ to obtain the last display of the corollary. This sign change concerns the eigenvalues only; the corresponding
spectral projections are unchanged.
\end{proof}

\subsubsection{Distribution and independence of the boundary weights}
\label{W:sec:gaussian-spectral-weights-and-projections}

The weights $W_j=|\varphi_j'(0)|^2$ are fixed by the $L^2$ normalization and Dirichlet boundary condition of the limiting eigenfunctions. Their distribution follows from the
Gaussian ensembles: conjugation invariance determines the finite spectral weights, while Proposition~\ref{W:prop:joint-spectral-data} identifies their operator limit.

Set $\alpha=\beta/2$. The Gamma distribution with shape $\alpha$ and scale $1/\alpha$ has density
\begin{equation}
 \frac{\alpha^\alpha}{\Gamma(\alpha)}
        x^{\alpha-1}e^{-\alpha x}\mathbf1_{\{x>0\}}.
 \label{W:eq:ev-gamma-density}
\end{equation}
Its mean is one.

\begin{lemma}
\label{W:lem:finite-gaussian-spectral-weights}
For $H_N^{\mathrm G}$, the weights at a deterministic unit vector are independent of the ordered eigenvalues and satisfy
\begin{equation}
 \left(|(v_{j,N}^{\mathrm G})_1|^2\right)_{j=1}^N
 \stackrel{\mathrm{law}}=
 \left(\frac{X_j}{\sum_{\ell=1}^N X_\ell}\right)_{j=1}^N,
 \qquad X_j\ \text{independent with law }\operatorname{Gamma}(\alpha,1).
 \label{W:eq:ev-gaussian-dirichlet-weights}
\end{equation}
\end{lemma}
\begin{proof}
Let $\mathcal G_N=O(N)$ for $\beta=1$ and $\mathcal G_N=U(N)$ for $\beta=2$, and let $\dd U$ be normalized Haar measure. The density of $H_N^{\mathrm G}$ in its independent real
coordinates is invariant under $H\mapsto UHU^*$. Its eigenvalues are simple almost surely: the set of repeated eigenvalues is the zero set of the nonzero discriminant polynomial
and has Lebesgue measure zero.

For a matrix with simple spectrum, write $\boldsymbol\mu(H)$ for its decreasing eigenvalue vector, $P_j(H)$ for its rank-one spectral projections, and $p_j(H,w)=w^*P_j(H)w$.
Invariance gives, for bounded Borel $F,G$,
\[
 \begin{split}
 &\mathbb E\bigl[F(\boldsymbol\mu(H_N^{\mathrm G}))
                   G((p_j(H_N^{\mathrm G},w))_{j\le N})\bigr]\\
 &\quad=\mathbb E\left[F(\boldsymbol\mu(H_N^{\mathrm G}))
     \int_{\mathcal G_N}
       G((w^*UP_j(H_N^{\mathrm G})U^*w)_{j\le N})\,\dd U\right].
 \end{split}
\]
Fubini's theorem applies because the integrand is bounded. Fix $H$ with simple spectrum and choose an orthogonal or unitary matrix $Q$ whose columns are its ordered eigenvectors.
The vector inside the second function is
\[
 \bigl(|(Q^*U^*w)_j|^2\bigr)_{j\le N}.
\]
As $U$ ranges under Haar measure, $Q^*U^*w$ has the uniform distribution on the unit sphere in $\mathbb F_\beta^N$. Indeed, its distribution is invariant under $\mathcal G_N$;
averaging a continuous function over this compact group gives the same value for every unit vector because the group action is transitive.

Take independent standard real Gaussian variables $g_{j,a}$, $1\le j\le N$, $1\le a\le\beta$, and define
\[
 \zeta_j=\begin{cases}g_{j,1},&\beta=1,\\g_{j,1}+\ii g_{j,2},&\beta=2,\end{cases}
 \qquad V=\frac{\zeta}{\|\zeta\|},\qquad X_j=\frac{|\zeta_j|^2}{2}.
\]
The density of $\zeta$ is invariant under $\mathcal G_N$, so $V$ has the same sphere distribution. Hence the group integral in the first display equals
\[
 c_G:=\mathbb E\,G\left(\left(\frac{X_j}{\sum_{\ell=1}^NX_\ell}\right)_{j\le N}\right),
\]
a number independent of the fixed matrix $H$. Substitution into that display yields
\[
 \mathbb E\bigl[F(\boldsymbol\mu(H_N^{\mathrm G}))G((p_j(H_N^{\mathrm G},w))_{j\le N})\bigr]
 =\mathbb EF(\boldsymbol\mu(H_N^{\mathrm G}))\,c_G.
\]
This establishes independence from the complete ordered spectrum.

The variables $X_j$ are independent. For $\beta=1$, the two solutions $g=\pm\sqrt{2x}$ in the change of variables $x=g^2/2$ give
\[
 \mathbb P(X_j\in\dd x)=\frac{x^{-1/2}e^{-x}}{\sqrt\pi}\,\dd x.
\]
For $\beta=2$, the radial density of $(g_{j,1},g_{j,2})$ is $re^{-r^2/2}\,\dd r$; the substitution $x=r^2/2$ gives $e^{-x}\,\dd x$. Therefore, for $\alpha=\beta/2$,
\[
 \mathbb P(X_j\in\dd x)=\frac{x^{\alpha-1}e^{-x}}{\Gamma(\alpha)}\,\dd x,
 \qquad x>0.
\]
Finally, the first column of the Lanczos matrix $U_N(w)$ is $w$, so $(v_{j,N}^{\mathrm G})_1=w^*\xi_{j,N}^{\mathrm G}$. The ordered eigenvectors above thus give exactly
\eqref{W:eq:ev-gaussian-dirichlet-weights}.
\end{proof}

\begin{theorem}
\label{W:thm:ev-original-weights}
Under Assumption~\ref{W:ass:wigner-lee-yin}, for every deterministic sequence of unit vectors $w_N$ and fixed $r,k$,
\begin{equation}
 \left((\lambda_{j,N})_{j=1}^r,
        (N|w_N^*\xi_{j,N}|^2)_{j=1}^k\right)
 \Longrightarrow
 \left((\Lambda_j)_{j=1}^r,(W_j)_{j=1}^k\right).
 \label{W:eq:ev-joint-original-weights}
\end{equation}
The $W_j$ are independent with density \eqref{W:eq:ev-gamma-density}, and their entire sequence is independent of
$(\Lambda_j)_{j\ge1}$.
\end{theorem}
\begin{proof}
Proposition~\ref{W:prop:joint-spectral-data} gives \eqref{W:eq:ev-joint-original-weights} for every deterministic unit sequence $w_N$, with $W_j=|\varphi_j'(0)|^2$. The GOE and GUE
satisfy the same proposition and yield the same distribution of the limiting operator. We therefore compute the weights using Lemma~\ref{W:lem:finite-gaussian-spectral-weights}.

Let $(X_j)_{j\ge1}$ be independent with distribution $\operatorname{Gamma}(\alpha,1)$ and set $S_N=\sum_{j=1}^NX_j$. Their moment generating function is
\[
 \mathbb Ee^{tX_j}=(1-t)^{-\alpha},\qquad t<1.
\]
For $0<\varepsilon<\alpha$, applying Markov's inequality to $e^{tS_N}$ with $t=\varepsilon/(\alpha+\varepsilon)$ and then with $t=-\varepsilon/(\alpha-\varepsilon)$ gives
\[
 \begin{aligned}
 \mathbb P(S_N/N\ge\alpha+\varepsilon)
   &\le e^{-N[\varepsilon-\alpha\log(1+\varepsilon/\alpha)]},\\
 \mathbb P(S_N/N\le\alpha-\varepsilon)
   &\le e^{-N[-\varepsilon-\alpha\log(1-\varepsilon/\alpha)]}.
 \end{aligned}
\]
Both rates in square brackets are positive. Their exponential bounds are summable in $N$, so Borel--Cantelli, applied to a countable sequence $\varepsilon\downarrow0$, gives
$S_N/N\to\alpha$ almost surely. On $\{|S_N/N-\alpha|<\alpha/2\}$,
\[
 \max_{1\le j\le k}\left|\frac{NX_j}{S_N}-\frac{X_j}{\alpha}\right|
 \le\frac2{\alpha^2}\max_{1\le j\le k}X_j
                \left|\frac{S_N}{N}-\alpha\right|\longrightarrow0
\]
for every fixed $k$. Therefore
\begin{equation}
 \left(\frac{NX_j}{S_N}\right)_{j=1}^k
   \longrightarrow (X_j/\alpha)_{j=1}^k
     \quad\text{almost surely}.
 \label{W:eq:ev-dirichlet-gamma-limit}
\end{equation}
The density of $X_j/\alpha$ is \eqref{W:eq:ev-gamma-density}, by the change of variables $x=\alpha y$.

For $r,k$ fixed write $\boldsymbol\lambda_N^{\mathrm G}=(\lambda_{j,N}^{\mathrm G})_{j\le r}$ and $Q_{j,N}^{\mathrm G}=N|(v_{j,N}^{\mathrm G})_1|^2$. The finite independence in
Lemma~\ref{W:lem:finite-gaussian-spectral-weights} yields, for $s\in\mathbb R^r$, $t\in\mathbb R^k$,
\[
 \begin{split}
 &\mathbb E\exp\left\{\ii s\cdot\boldsymbol\lambda_N^{\mathrm G}
                 +\ii\sum_{j=1}^kt_jQ_{j,N}^{\mathrm G}\right\}\\
 &\qquad=\mathbb Ee^{\ii s\cdot\boldsymbol\lambda_N^{\mathrm G}}
       \mathbb E\exp\left\{\ii\sum_{j=1}^kt_j\frac{NX_j}{S_N}\right\}.
 \end{split}
\]
Proposition~\ref{W:prop:joint-spectral-data} gives convergence of the joint expression and of the eigenvalue factor. The weight factor converges by
\eqref{W:eq:ev-dirichlet-gamma-limit} and dominated convergence, since its integrand has modulus one. Moreover,
\[
 \mathbb Ee^{\ii tX_j/\alpha}=(1-\ii t/\alpha)^{-\alpha},
\]
we obtain
\[
 \mathbb E\exp\left\{\ii\sum_{j=1}^rs_j\Lambda_j
                   +\ii\sum_{j=1}^kt_jW_j\right\}
 =\mathbb E\exp\left\{\ii\sum_{j=1}^rs_j\Lambda_j\right\}
          \prod_{j=1}^k(1-\ii t_j/\alpha)^{-\alpha}.
\]
Uniqueness of characteristic functions now identifies the joint distribution for every finite $r,k$: the $W_j$ are independent with density \eqref{W:eq:ev-gamma-density}, and the
two finite vectors are independent.

For independence of the full sequences, let $\mathcal C_\Lambda$ and $\mathcal C_W$ denote the events determined by finitely many initial eigenvalues and weights, respectively,
with Borel conditions on the corresponding finite vectors. They are $\pi$-systems generating $\sigma(\Lambda_j:j\ge1)$ and $\sigma(W_j:j\ge1)$. The preceding factorization gives
\[
 \mathbb P(A\cap B)=\mathbb P(A)\mathbb P(B),
 \qquad A\in\mathcal C_\Lambda,
 \ B\in\mathcal C_W.
\]
For fixed $B\in\mathcal C_W$, the eigenvalue events satisfying this identity form a Dynkin system containing $\mathcal C_\Lambda$; the $\pi$--$\lambda$ theorem extends it to all of
$\sigma(\Lambda_j:j\ge1)$. Fixing such an $A$ and applying the same theorem to the weight events extends the identity to all of $\sigma(W_j:j\ge1)$. This establishes the
independence in the theorem. The convergence of the original Wigner matrices, already proved by Proposition~\ref{W:prop:joint-spectral-data}, has this same limiting distribution.
\end{proof}

The conclusion on simple eigenvalues and nonzero scalar products in Proposition~\ref{W:prop:joint-spectral-data} also gives, for every fixed $k$,
\[
 \mathbb P\bigl(\lambda_{1,N},\ldots,\lambda_{k,N}\text{ are simple and }
          w_N^*\xi_{j,N}\ne0\text{ for all }j\le k\bigr)\longrightarrow1.
\]
For any fixed bounded Lipschitz function $F$ of the finite tuple in \eqref{W:eq:ev-joint-original-weights}, let $X_N(w)$ denote that tuple for a deterministic unit vector $w$, and
let $X$ have the limiting law. We claim
\[
 \sup_{\|w\|=1,\,w\ {\rm deterministic}}
       |\mathbb EF(X_N(w))-\mathbb EF(X)|\longrightarrow0.
\]
If this failed, there would be $\varepsilon>0$ and a strictly increasing sequence $N_l$ for which the supremum is greater than $\varepsilon$. For each $l$, choose a deterministic
unit vector $w_{N_l}$ whose discrepancy is greater than $\varepsilon/2$. Attainment of the supremum is not required. At all other sizes choose, for example, the first coordinate
vector. This constructs one deterministic sequence $(w_N)$ to which Theorem~\ref{W:thm:ev-original-weights} applies. Its conclusion forces the discrepancy to tend to zero along
$N_l$, a contradiction. The joint convergence with the potential, resolvent, and eigenfunctions is the one in Proposition~\ref{W:prop:joint-spectral-data}. This argument does not
permit choosing $w_N$ as a function of the observed matrix.
\begin{proof}[Proof of Corollary~\ref{W:cor:intro-boundary-weights}]
Fix $r,k$ and put $m=\max\{r,k\}$. Apply Theorem~\ref{W:thm:ev-original-weights} to the first $m$ spectral labels and then project onto the first $r$ eigenvalue coordinates and the
first $k$ weight coordinates. The finite eigenvalue coordinates in this corollary are exactly
\[
 N^{2/3}(2-\mu_j^{(N)})=\lambda_{j,N}.
\]
Thus their limits are $\Lambda_j$. The opposite normalization $N^{2/3}(\mu_j^{(N)}-2)$, used in the preceding corollary, would instead have limit $-\Lambda_j$.

To identify the weights, the first Lanczos vector is $w_N$. Consequently, after conjugation to the Jacobi coordinates, the modulus of the first coordinate of the eigenvector is
$|w_N^*\xi_{j,N}|$. Equation~\eqref{W:eq:ev-original-boundary-identity} identifies its squared boundary normalization with
\[
 N|w_N^*\xi_{j,N}|^2.
\]
The boundary convergence in Lemma~\ref{W:lem:ev-boundary-derivative} therefore gives the limiting weight $|\varphi_j'(0)|^2$; its identification uses the recurrence in addition to
$L^2$ convergence. Measurable choices on the exceptional event of multiplicity or zero first coordinate have no effect, because that event has probability tending to zero for these
finitely many labels.

Writing $\alpha=\beta/2$, Theorem~\ref{W:thm:ev-original-weights} identifies the joint law of these weights as the product of the densities
\[
 \frac{\alpha^\alpha}{\Gamma(\alpha)}x^{\alpha-1}e^{-\alpha x}
       \mathbf1_{(0,\infty)}(x)\,\mathrm dx.
\]
This is the Gamma law with shape $\alpha$ and scale $1/\alpha=2/\beta$. The same theorem gives independence of the entire weight sequence from the limiting eigenvalue sequence.
Projection onto the selected finite coordinates preserves this independence and proves the corollary.
\end{proof}

\subsubsection{Projections onto deterministic vectors}
\label{W:sec:projections-onto-deterministic-vectors}

The preceding theorem identifies the squared overlaps. Multiplying each eigenvector by an independent uniform sign for $\beta=1$, or by an independent uniform complex number of
modulus one for $\beta=2$, also determines the distribution of the overlaps themselves. Let $\omega_{j,N}$ be these independent variables, independent of the matrix, and set
$\widehat\xi_{j,N}=\omega_{j,N}\xi_{j,N}$. Then
\[
 |\widehat\xi_{j,N}\rangle\langle\widehat\xi_{j,N}|
       =|\xi_{j,N}\rangle\langle\xi_{j,N}|.
\]
A standard Gaussian over $\mathbb F_\beta$ is $N(0,1)$ for $\beta=1$, and $(X+\ii Y)/\sqrt2$ for $\beta=2$, where $X,Y$ are independent $N(0,1)$ variables.

\Needspace{7\baselineskip}
\begin{corollary}
\label{W:cor:ev-one-direction-gaussian}
Under Assumption~\ref{W:ass:wigner-lee-yin}, for any deterministic unit sequence $w_N$ and fixed $r,k$,
\begin{equation}
 \left((\lambda_{j,N})_{j=1}^r,
       (\sqrt N\,w_N^*\widehat\xi_{j,N})_{j=1}^k\right)
 \Longrightarrow
 \left((\Lambda_j)_{j=1}^r,(Z_j)_{j=1}^k\right),
 \label{W:eq:ev-one-direction-gaussian}
\end{equation}
where the $Z_j$ are independent standard Gaussians over $\mathbb F_\beta$, independent of the eigenvalues. In particular this applies to $w_N=e_{i_N}$ for every deterministic
coordinate sequence $1\le i_N\le N$.
\end{corollary}
\begin{proof}
Put $Q_{j,N}=N|w_N^*\xi_{j,N}|^2$. Conditional on the matrix, multiplication by $\omega_{j,N}$ makes the argument of $w_N^*\widehat\xi_{j,N}$ uniform on the sign group or unit
circle, independently for different $j$. Thus, for any bounded continuous $F$,
\[
 \begin{split}
 &\mathbb EF\left(\boldsymbol\lambda_N,
          (\sqrt N\,w_N^*\widehat\xi_{j,N})_{j\le k}\right)\\
 &\quad=\mathbb E\int_{\mathbb T_\beta^k}
      F\left(\boldsymbol\lambda_N,(\theta_j\sqrt{Q_{j,N}})_{j\le k}\right)
                              \prod_{j=1}^k\dd\theta_j,
 \end{split}
\]
where $\mathbb T_1=\{-1,1\}$ with equal weights and $\mathbb T_2$ is the unit circle with normalized Haar measure. The integral on the right is a bounded continuous function of
$(\boldsymbol\lambda_N,(Q_{j,N})_{j\le k})$: continuity follows from that of the square root and dominated convergence in $\theta$. Theorem~\ref{W:thm:ev-original-weights}
therefore gives the limit
\[
 \bigl((\Lambda_j)_{j\le r},(\theta_j\sqrt{W_j})_{j\le k}\bigr),
\]
with the $\theta_j$ independent of each other and of the weights and eigenvalues.

For $\beta=1$, $W_j$ has Gamma shape $1/2$ and scale $2$, so
\[
 \mathbb P(\theta_j\sqrt{W_j}\in\dd x)
      =(2\pi)^{-1/2}e^{-x^2/2}\,\dd x.
\]
For $\beta=2$, $W_j$ has density $e^{-w}$ on $(0,\infty)$. With $r=\sqrt w$ and a uniform angle, the polar probability element is
\[
 \frac{1}{2\pi}\,2re^{-r^2}\,\dd r\,\dd\theta
        =\pi^{-1}e^{-|z|^2}\,\dd^2z.
\]
These are the respective standard Gaussian distributions. The product structure established in Theorem~\ref{W:thm:ev-original-weights} gives independence of the limiting
coordinates and of the eigenvalue sequence.
\end{proof}

\Needspace{6\baselineskip} Applying the limit for projection onto one deterministic vector to deterministic linear combinations gives the result for finite families of
deterministic vectors. Their limiting covariance is determined by the inner products of those vectors, including when the resulting matrix is singular.

\Needspace{7\baselineskip} Fix an eigenvalue index $j\ge1$ and a positive integer $d$. Let $q_{1,N},\ldots,q_{d,N}\in\mathbb F_\beta^N$ be deterministic and suppose their Gram
matrices converge:
\begin{equation}
 \Sigma_N=(q_{a,N}^*q_{b,N})_{a,b=1}^d\longrightarrow\Sigma.
 \label{W:eq:deterministic-vector-gram}
\end{equation}

\begin{corollary}
\label{W:cor:ev-many-directions}
Under Assumption~\ref{W:ass:wigner-lee-yin}, jointly with any fixed number of ordered edge eigenvalues,
\begin{equation}
 \left(\sqrt N\,q_{a,N}^*\widehat\xi_{j,N}\right)_{a=1}^d
       \Longrightarrow Z,
 \label{W:eq:ev-many-projections}
\end{equation}
where $Z$ is centered Gaussian with covariance $\mathbb EZZ^*=\Sigma$, independent of those eigenvalues. In the complex case $Z$ is circular, so $\mathbb EZZ^{\mathsf T}=0$. The
covariance matrix may be singular.
\end{corollary}
\begin{proof}
Write $Z_N=(\sqrt N\,q_{a,N}^*\widehat\xi_{j,N})_{a=1}^d$. For fixed $t\in\mathbb F_\beta^d$, put $q_N(t)=\sum_at_aq_{a,N}$ and $c_N=\|q_N(t)\|$. Then
\[
 t^*Z_N=\sqrt N\,q_N(t)^*\widehat\xi_{j,N},\qquad
 c_N^2=t^*\Sigma_Nt\to c^2=t^*\Sigma t.
\]
Set $w_N(t)=q_N(t)/c_N$ when $c_N>0$ and $w_N(t)=e_1$ otherwise. Corollary~\ref{W:cor:ev-one-direction-gaussian} gives $V_N=\sqrt N\,w_N(t)^*\widehat\xi_{j,N}\Longrightarrow Z_0$
jointly with any fixed eigenvalue vector, with $Z_0$ a standard Gaussian independent of that vector. Since $t^*Z_N=c_NV_N$ at every $N$, Slutsky's theorem gives the joint scalar
limit $cZ_0$, including $c=0$. Thus
\begin{equation}
 \mathbb E\exp\{\ii s\cdot\boldsymbol\lambda_N
                    +\ii\operatorname{Re}(t^*Z_N)\}
 \longrightarrow
 \mathbb E e^{\ii s\cdot\boldsymbol\Lambda}
 \begin{cases}
  \exp\{-\tfrac12 t^{\mathsf T}\Sigma t\},&\beta=1,\quad t\in\mathbb R^d,\\
  \exp\{-\tfrac14 t^*\Sigma t\},&\beta=2,\quad t\in\mathbb C^d.
 \end{cases}
 \label{W:eq:ev-joint-characteristic-functions}
\end{equation}
The nonnegative Hermitian matrix $\Sigma$ has a nonnegative square root, even if singular. Let $G$ have independent standard Gaussian coordinates over $\mathbb F_\beta$,
independently of $\boldsymbol\Lambda$, and set $Z=\Sigma^{1/2}G$. Then $\mathbb EZZ^*=\Sigma$; in the complex case $Z$ is circular and $\mathbb EZZ^{\mathsf T}=0$, with real
covariance
\[
 \operatorname{Cov}\binom{\operatorname{Re}Z}{\operatorname{Im}Z}
 =\frac12\begin{pmatrix}\operatorname{Re}\Sigma&-\operatorname{Im}\Sigma\\
                         \operatorname{Im}\Sigma&\operatorname{Re}\Sigma
          \end{pmatrix}.
\]
Its scalar characteristic functions are precisely the factors in \eqref{W:eq:ev-joint-characteristic-functions}. Since $\operatorname{Re}(t^*z)$ ranges over all real linear
functionals on $\mathbb C^d\cong\mathbb R^{2d}$, the continuity theorem for characteristic functions proves the full joint convergence and independence in the corollary.
\end{proof}

\begin{corollary}
\label{W:cor:ev-projection-products}
Under the hypotheses of Corollary~\ref{W:cor:ev-many-directions}, the phase invariant products satisfy
\begin{equation}
 \left(N(q_{a,N}^*\xi_{j,N})
                 \overline{(q_{b,N}^*\xi_{j,N})}\right)_{a,b=1}^d
          \Longrightarrow ZZ^*.
 \label{W:eq:ev-original-projection-products}
\end{equation}
In particular, for deterministic pairwise distinct coordinates $i_{1,N},\ldots,i_{d,N}$,
\begin{equation}
 \left(\sqrt N\,\widehat\xi_{j,N}(i_{a,N})\right)_{a=1}^d
         \Longrightarrow (Z_1,\ldots,Z_d),
 \label{W:eq:ev-normal-coordinate-entries}
\end{equation}
with independent standard Gaussians over $\mathbb F_\beta$ on the right. For an orthonormal family $q_{1,N},\ldots,q_{d,N}$, the mass in its span satisfies
\begin{equation}
 N\sum_{a=1}^d|q_{a,N}^*\xi_{j,N}|^2
       \Longrightarrow
       \operatorname{Gamma}(\beta d/2,2/\beta).
 \label{W:eq:ev-finite-subspace-mass}
\end{equation}
\end{corollary}
\begin{proof}
For $Z_N$ from Corollary~\ref{W:cor:ev-many-directions}, cancellation of the common phase gives
\[
 (Z_NZ_N^*)_{ab}
 =N(q_{a,N}^*\xi_{j,N})\overline{(q_{b,N}^*\xi_{j,N})}.
\]
The continuous map $z\mapsto zz^*$ proves \eqref{W:eq:ev-original-projection-products}, jointly with the finite eigenvalue vector. For distinct coordinate directions, the Gram
matrix is $I_d$, so the limit has independent standard Gaussian coordinates, proving \eqref{W:eq:ev-normal-coordinate-entries}.

For any orthonormal family, the continuous map $z\mapsto\|z\|^2$ gives the limiting mass $\sum_{a=1}^d|Z_a|^2$ with the same independent coordinates. Writing $\alpha=\beta/2$, its
Laplace transform is
\[
 \mathbb E e^{-s\sum_a|Z_a|^2}
 =\prod_{a=1}^d\left(\frac\alpha{\alpha+s}\right)^\alpha
 =\left(\frac\alpha{\alpha+s}\right)^{\alpha d},\qquad s\ge0.
\]
This is the Gamma law with shape $\alpha d=\beta d/2$ and scale $1/\alpha=2/\beta$, proving \eqref{W:eq:ev-finite-subspace-mass}.
\end{proof}

\subsubsection{Weighted spectral measures}
\label{W:sec:weighted-spectral-measures}

We now combine the eigenvalues and the squared overlaps into a measure on the rescaled spectral axis. Define the positive random measures
\begin{equation}
 \nu_{N,w}=N\sum_{j=1}^N|w_N^*\xi_{j,N}|^2\delta_{\lambda_{j,N}},
 \qquad
 \nu_{\beta}=\sum_{j\ge1}|\varphi_j'(0)|^2\delta_{\Lambda_j}.
 \label{W:eq:ev-starting-vectored-edge-measures}
\end{equation}
The first expression is independent of the basis chosen inside a repeated eigenspace, because weights at the same eigenvalue add to the corresponding spectral projection weight. We
use the vague topology on positive locally finite measures, defined by integration against $C_c(\mathbb R)$.

\Needspace{7\baselineskip}
\begin{corollary}
\label{W:cor:ev-starting-vectored-measure}
Under Assumption~\ref{W:ass:wigner-lee-yin}, for each deterministic unit sequence $w_N$,
\begin{equation}
 \nu_{N,w}\Longrightarrow\nu_\beta
       =\sum_{j\ge1}W_j\delta_{\Lambda_j}
       \quad\hbox{in the vague topology on locally finite measures}.
 \label{W:eq:ev-starting-vectored-measure-convergence}
\end{equation}
The Gamma weights are independent of the eigenvalue point process. The convergence holds jointly with the unweighted edge point process, the integrated potential and the embedded
resolvent. In particular, for every $F\in C_c(\mathbb R)$,
\begin{equation}
 Nw_N^*F\!\left(N^{2/3}(2I-H_N)\right)w_N
       \Longrightarrow \sum_{j\ge1}W_jF(\Lambda_j).
 \label{W:eq:ev-original-spectral-test}
\end{equation}
\end{corollary}
\begin{proof}
Fix the deterministic starting sequence. From any subsequence, apply Proposition~\ref{W:prop:joint-spectral-data} to the joint law of the coefficients and resolvents to obtain a
further subsequence and copies with their original joint distributions, on which
\[
 \widetilde Y_N^{\rm lin}\to Y_\beta,\qquad
 \mathscr R_N\to\mathscr R_\beta,\qquad
 B_N(\ell)\to0\ (\ell\ge1),\qquad
 \kappa_N:=(-\lambda_{1,N})_+\to\kappa<\infty
\]
almost surely. That proof applies on one event to all fixed indices and gives
\[
 \lambda_{j,N}\to\Lambda_j,\qquad
 W_{j,N}:=N|w_N^*\xi_{j,N}|^2\to W_j=|\varphi_j'(0)|^2
 \quad(j\ge1).
\]
The lower bound supplied by $\kappa_N$ is included in that argument. In particular these are ordered eigenvalues of the full finite matrices.

For a repeated eigenvalue $\lambda$ of $\mathcal A_N$, the sum of the weights is $Ne_1^*\mathbf1_{\{\lambda\}}(\mathcal A_N)e_1$, independently of the chosen eigenbasis. Thus the
copied measures have their original joint distribution with $\widetilde Y_N^{\rm lin}$ and $\mathscr R_N$, because
\[
 \int F\,\dd\nu_{N,w}=Ne_1^*F(\mathcal A_N)e_1,
 \qquad \int F\,\dd\Xi_N=\operatorname{Tr}F(\mathcal A_N)
 \quad(F\in C_c(\mathbb R))
\]
are Borel functions of the reconstructed matrix. At each fixed dimension they are continuous, by continuous functional calculus. The limiting measures are locally finite because
$\Lambda_j\to +\infty$ and every $W_j$ is finite; their integrals against continuous functions of compact support are Borel limits of finite sums of the Borel variables
$\Lambda_j,W_j$.

Fix an outcome in the common event and a compact interval $[-L,L]$. Choose $m\ge0$ with $\Lambda_{m+1}>L+1$. Ordered convergence gives an index $N_L$ such that $\lambda_{m+1,N}>L$
for every $N\ge N_L$. Thus for every $F$ supported in $[-L,L]$ the four spectral sums reduce to
\[
 \begin{aligned}
 \int F\,\dd\nu_{N,w}&=\sum_{j=1}^mW_{j,N}F(\lambda_{j,N}),&
 \int F\,\dd\nu_\beta&=\sum_{j=1}^mW_jF(\Lambda_j),\\
 \int F\,\dd\Xi_N&=\sum_{j=1}^mF(\lambda_{j,N}),&
 \int F\,\dd\Xi_\beta&=\sum_{j=1}^mF(\Lambda_j).
 \end{aligned}
\]
The empty sums cover $m=0$. For the weighted sums,
\[
 \left|\int F\,\dd\nu_{N,w}-\int F\,\dd\nu_\beta\right|
 \le\|F\|_\infty\sum_{j=1}^m|W_{j,N}-W_j|
       +\sum_{j=1}^mW_j|F(\lambda_{j,N})-F(\Lambda_j)|\to0.
\]
The unweighted sums converge by the same finite sum argument, and
\[
 \limsup_N\nu_{N,w}([-L,L])\le\sum_{j=1}^mW_j<\infty,
 \qquad \limsup_N\Xi_N([-L,L])\le m.
\]
Only local mass bounds are used.

Both measurability and convergence in the vague topology can be checked on a countable family. Take all compactly supported piecewise linear functions with rational breakpoints and
rational values, including functions equal to one on $[-\ell,\ell]$ and zero off $[-\ell-1,\ell+1]$ for each integer $\ell$. They uniformly approximate any $F\in C_c(\mathbb R)$ on
one fixed larger support $J$. Their integrals determine the vague topology and its Borel sigma field: convergence of the cutoff integrals bounds masses on $J$, and
\[
 \left|\int F\,\dd\mu_N-\int F\,\dd\mu\right|
 \le\left|\int f\,\dd\mu_N-\int f\,\dd\mu\right|
            +\|F-f\|_\infty\bigl(\mu_N(J)+\mu(J)\bigr)
\]
passes convergence from these functions to every $F$. Equivalently, their integrals embed the measure space in a countable product of copies of $\mathbb R$, with the induced
metrizable topology. This proves that the Borel measurable integrals above define random measures and that their simultaneous convergence is vague convergence. On the coupling,
therefore,
\[
 (\nu_{N,w},\Xi_N,\widetilde Y_N^{\rm lin},\mathscr R_N)
 \longrightarrow(\nu_\beta,\Xi_\beta,Y_\beta,\mathscr R_\beta)
\]
almost surely in the product of the two vague and two path topologies. Equality of the copied and original joint laws and the subsequence criterion give full sequence joint
convergence in distribution.

Theorem~\ref{W:thm:ev-original-weights} identifies $(W_j)$ as independent Gamma variables, independent of the entire eigenvalue sequence. Since $\Xi_\beta$ is a Borel function of
that sequence, it is also independent of the weights. Finally unitary functional calculus and $U_Ne_1=w_N$ give
\[
 Nw_N^*F\!\left(N^{2/3}(2I-H_N)\right)w_N
 =Ne_1^*F(\mathcal A_N)e_1=\int F\,\dd\nu_{N,w}.
\]
Evaluation at $F\in C_c(\mathbb R)$ is continuous for the vague topology, proving \eqref{W:eq:ev-original-spectral-test}. For complex $F$, apply the joint evaluation at its real
and imaginary parts.
\end{proof}

Corollary~\ref{W:cor:ev-one-direction-gaussian} gives the joint limit for one deterministic direction and finitely many eigenvectors. Corollary~\ref{W:cor:ev-many-directions} gives
the joint limit for one eigenvector and finitely many deterministic directions. These are spectral projection observables of the type studied in
\cite{knowles2013eigenvectors,bourgade2017momentflow}. The matrix $-H_N$ satisfies Assumption~\ref{W:ass:wigner-lee-yin}, and its decreasing eigenvalues are $-\mu_{N+1-j}^{(N)}$,
$1\le j\le N$. Applying the stated results to $-H_N$ gives the marginal left edge limits with eigenvalue scaling $N^{2/3}(2+\mu_{N+1-j}^{(N)})$.

\subsection{Spectral consequences of the tridiagonal approximation}
\label{sec:combined-a-spectral}
The quantitative comparison of the discrete and limiting resolvents gives rates for isolated spectral projections and eigenvalues. Increasing eigenvalue labels still require the
lower spectral bound, and the full Jacobi matrix must be distinguished from the compressed operator used in the isometric comparison.

\subsubsection{Eigenvalues, spectral projections, and point measures}\label{A:sec:spectral-consequences}

We begin with the matrices from the Gaussian ensemble and the reduced general model, together with their isometric embeddings $(A_n,J_n)$ specified in
Subsection~\ref{A:sec:atomlessness-real-resolvents-and-inverse-coupling}. The full Jacobi matrix will be treated in the enlarged space of Theorem~\ref{A:res:20.6}; compression
alone does not identify its spectrum. Set $S_n=J_n(A_n-\mathrm i)^{-1}J_n^*$ and $S=(\mathcal H_\beta-\mathrm i)^{-1}$. The nonreal comparison gives
\[
 \|S_n-S\|_{\mathrm{op}}
 \le\|S_n-S\|_{\mathfrak S_2}\le C_\omega h_n^{1/2}.
 \eqtag{21.2}
\]
The map $\lambda\mapsto(\lambda-\mathrm i)^{-1}$ sends both ends of the real line to zero. Thus resolvent convergence can control every bounded spectral interval and still miss an
eigenvalue escaping to $-\infty$. To identify the increasing labels we need a lower bound. For example, let $He_j=je_j$ on $\ell^2(\mathbb N)$ and take
$A_n=\operatorname{diag}(1,\ldots,n-1,-n)$ on its first $n$ coordinates, with inclusion $\iota_n$. The resolvent difference is diagonal. Its $n$-th entry has absolute value
$2n/(n^2+1)$, and its entries at $j>n$ have absolute value $(j^2+1)^{-1/2}$. Thus
\[
 \begin{split}
 \|\iota_n(A_n-\mathrm i)^{-1}\iota_n^*-(H-\mathrm i)^{-1}\|_{\mathrm{op}}
     &\le2/n,\\
 \|\iota_n(A_n-\mathrm i)^{-1}\iota_n^*-(H-\mathrm i)^{-1}\|_{\mathfrak S_2}^2
     &=\frac{4n^2}{(n^2+1)^2}+\sum_{j>n}\frac1{j^2+1}\longrightarrow0,
 \end{split}
\]
although $\lambda_1(A_n)=-n$. We prove the lower bound separately, with a finite size formulation for use in Section~\ref{A:sec:parameter-dependence-and-further-directions}. The
boundary estimates will show positivity after a common scalar shift, giving the required lower bound for the full spectrum.

Fix $Q>0$, $0\le R<\infty$, and $[\beta_-,\beta_+]\Subset(0,\infty)$. Choose one integer $s\ge\sup_{\beta\in[\beta_-,\beta_+]}s_\beta$ and delete the final $s$ coordinates in each
finite Jacobi matrix. In this paragraph $A_C$ is the principal matrix on $\{r_n+1,\ldots,n-s\}$ and $A_P$ is the principal matrix on $\{1,\ldots,n-s\}$. Use the coefficient events
$\mathcal L_{n,Q}$ of Proposition~\ref{A:res:left-spectral} and $\mathcal V_{n,Q}$ of Lemma~\ref{A:res:core-finite}, with the same $s$ deleted coordinates. In particular
$\mathcal V_{n,Q}$ includes $\mathcal F_{n,Q+2}$, so it controls the coefficients in the final $s$ coordinates and the off-diagonal coefficient joining them to the retained matrix.
All constants and thresholds in the next lemma are deterministic on these events and uniform in $\beta$.

\begin{lemma}\label{A:res:lower-transfer}
For $n\ge n_0(Q,R)$, on $\mathcal L_{n,Q}\cap\mathcal V_{n,Q}$, every $a\in[-R,0]$ with $A_C-a\ge I$ satisfies
\[
 A_P-a>0,\qquad \|(A_P-a)^{-1}\|_{\mathrm{op}}\le2,
 \qquad \mathcal A_{n,V,\beta}-a>0.
 \eqtag{lower-transfer}
\]
\end{lemma}

\begin{proof}
Write $r=r_n$ and $\ell_n=1+\log(n+1)$. The assumed events give $\alpha\asymp h^{-2}$, $A_L-\alpha e_Le_L^*\ge ch^{-2}r^{-2}I$, and the estimates of Lemma~\ref{A:res:core-poisson}
at $z=a+\mathrm i$, uniformly for $a\in[-R,0]$. Put $Y=A_C-a$ and $R_C=Y^{-1}$. By assumption $0<R_C\le I$. If $S_C=(Y-\mathrm i)^{-1}$, then
\[
 R_C=S_C(I+\mathrm iS_C)^{-1},\qquad
 (I+\mathrm iS_C)^{-1}=I-\mathrm iR_C.
\]
The second factor has norm at most two and commutes with $S_C$. The nonreal boundary estimates therefore give
\[
 \|R_Ce_C\|^2\le C_{Q,R}h^3\ell_n^6,
 \qquad |d_C(a)|\le C_{Q,R}h^{-1}\ell_n^6.
 \eqtag{lower-real-core}
\]
For the second inequality, the difference from $d_C(a+\mathrm i)$ is bounded by $\alpha^2\|R_Ce_C\|\|S_Ce_C\|\le Ch^{-1}\ell_n^6$; its nonreal value is bounded by
$Ch^{-1}\ell_n^4$.

Put $X=A_L-a$, $m_L=e_L^*X^{-1}e_L$ and $m_C=e_C^*Y^{-1}e_C$. Lemma~\ref{A:res:left-resolvent} gives
\[
 X>0,\quad cr^{-2}\le\kappa_L(a)/\alpha\le1,
 \quad \|X^{-1}\|\le Ch^2r^2,
 \quad \|X^{-1}e_L\|^2\le Ch^4r^2.
\]
As in the preceding real comparison,
\[
 D_n(a) =1-\alpha^2m_Lm_C =\frac{\kappa_L(a)}\alpha +\frac{d_C(a)}\alpha
 \left(1-\frac{\kappa_L(a)}\alpha\right) \ge\frac{c}{r^2}-C_{Q,R}h\ell_n^6 \ge\frac{c}{2r^2}>0
 \eqtag{positive-denominator}
\]
for all $n$ beyond a deterministic threshold. The last step uses $hr^2\ell_n^6\to0$. Positivity follows by completing the square:
\[
 \langle(A_P-a)(u,v),(u,v)\rangle ={} \|X^{1/2}(u-\alpha X^{-1}e_Le_C^*v)\|^2
 +\langle(Y-\alpha^2m_Le_Ce_C^*)v,v\rangle.
 \eqtag{schur-square}
\]
This identity holds also for complex vectors, by expansion of the real part of the cross term. With $w=Y^{-1/2}e_C$,
\[
 Y-\alpha^2m_Le_Ce_C^*
 =Y^{1/2}(I-\alpha^2m_Lww^*)Y^{1/2}.
\]
The middle operator is the identity on $w^\perp$ and is multiplication by $D_n(a)>0$ on $\operatorname{span}\{w\}$, because $\|w\|^2=m_C$. The triangular change in
\eqref{A:eq:schur-square} is invertible. Hence $A_P-a>0$.

To retain a fixed lower bound after restoring the final $s$ coordinates, compute the four blocks of the inverse. Put
\[
 a_L=X^{-1}e_L,\qquad a_C=Y^{-1}e_C,\qquad D=1-\alpha^2m_Lm_C.
\]
The off-diagonal block of $A_P-a$ is $-\alpha e_Le_C^*$. Its Schur complement with respect to $X$ is $Y-\alpha^2m_Le_Ce_C^*$. Direct multiplication, using $e_C^*Y^{-1}e_C=m_C$,
gives
\[
 (Y-\alpha^2m_Le_Ce_C^*)^{-1}
      =Y^{-1}+\frac{\alpha^2m_L}{D}a_Ca_C^*.
\]
Solving the first block equation for the left coordinate and substituting this expression into the second gives
\[
 (A_P-a)^{-1}=
 \begin{pmatrix}
 X^{-1}+\dfrac{\alpha^2m_C}{D}a_La_L^*
       &\dfrac{\alpha}{D}a_La_C^*\\[4pt]
 \dfrac{\alpha}{D}a_Ca_L^*
       &Y^{-1}+\dfrac{\alpha^2m_L}{D}a_Ca_C^*
 \end{pmatrix}.
\]
Every denominator is nonzero by \eqref{A:eq:positive-denominator}. The definitions of $\kappa_L,d_C$, together with $\kappa_L/\alpha\in[cr^{-2},1]$ and $|d_C|/\alpha=o(1)$, give
$|m_L|+|m_C|\le C\alpha^{-1}$. Hence
\[
 \left|\frac{\alpha^2m_L}{D}\right|
 +\left|\frac{\alpha^2m_C}{D}\right|
 +\left|\frac{\alpha}{D}\right|\le C\alpha r^2.
\]
For a rank-one operator $uv^*$, both its operator norm and its Hilbert--Schmidt norm equal $\|u\|\|v\|$. Using $\alpha\le Ch^{-2}$ and the column bounds already proved, the four
contributions to the operator-norm difference are bounded as follows:
\[
 \begin{aligned}
 \left\|\frac{\alpha^2m_L}{D}a_Ca_C^*\right\|
   &\le Ch^{-2}r^2(h^3\ell_n^6)=Chr^2\ell_n^6,\\
 2\left\|\frac{\alpha}{D}a_La_C^*\right\|
   &\le Ch^{-2}r^2(h^2r)(h^{3/2}\ell_n^3)
       =Ch^{3/2}r^3\ell_n^3,\\
 \left\|\frac{\alpha^2m_C}{D}a_La_L^*\right\|
   &\le Ch^{-2}r^2(h^4r^2)=Ch^2r^4,\\
 \|X^{-1}\|&\le Ch^2r^2.
 \end{aligned}
\]
The triangle inequality for these embedded blocks therefore proves
\[
 \|(A_P-a)^{-1}-P_C^*R_CP_C\|_{\mathrm{op}}
 \le C_{Q,R}\bigl(hr^2\ell_n^6+h^{3/2}r^3\ell_n^3+h^2r^4+h^2r^2\bigr)=o(1).
 \eqtag{positive-prefix-margin}
\]
Choose a deterministic size, depending on $Q,R$ and the fixed parameter interval, above which the displayed bound is at most one. Since $\|R_C\|\le1$, this gives
$\|(A_P-a)^{-1}\|\le2$. The positivity already proved then implies $A_P-a\ge I/2$ by the spectral theorem. On the same coefficient event the principal matrix on the final $s$
coordinates satisfies $A_R-a\ge ch^{-2}I$, and for a fixed exponent $p_Q$, $|\alpha_R|\le C_Qh^{-1/2}\ell_n^{p_Q}$. Thus
\[
 0\le\sigma_R(a)=\alpha_R^2e_R^*(A_R-a)^{-1}e_R
      \le C_{Q,R}h\ell_n^{2p_Q}=o(1).
\]
For $n$ sufficiently large the Schur complement satisfies
\[
 A_P-a-\sigma_R(a)e_Pe_P^*\ge\tfrac14 I.
 \eqtag{positive-terminal}
\]
Completing the square in the terminal coordinates, now with $A_R-a>0$, proves $\mathcal A_{n,V,\beta}-a>0$. The error terms in these inverse estimates are bounded by an explicit
power of $h$ times a fixed power of $r$ or $\ell_n$. Since $h\asymp n^{-1/3}$ uniformly in the specified parameter interval, one deterministic $n_0(Q,R)$ works for all the stated
matrices and real parameters.
\end{proof}

\begin{proposition}\label{A:res:common-lower}
For each fixed $\beta>0$, there are finite random $L_\omega,N_\omega$ such that
\[
 \mathcal H_\beta\ge-L_\omega I,\qquad \mathcal A_{n,\beta}^{\mathrm{DE}}\ge-L_\omega I,\qquad
 \widetilde H_n^{\mathrm{red}}\ge-L_\omega I,\qquad \mathcal A_{n,V,\beta}\ge-L_\omega I
 \qquad(n\ge N_\omega).
 \eqtag{common-lower}
\]
\end{proposition}

The statement for both models uses their conditional product coupling with the same Brownian path. Each original marginal is unchanged by taking this product.

\begin{proof}
Fix $\beta>0$ and use the conditional product coupling in the statement. Intersect the event of Brownian paths of probability one with the events in Proposition~\ref{A:res:8.2B}
and Theorem~\ref{A:res:16.2}. On this intersection there are finite shifts $\kappa_0,\kappa_1,\kappa_2\ge0$ and finite indices $N_1,N_2$ such that
\[
 \begin{aligned}
 \mathfrak q_\beta[f]+\kappa_0\|f\|^2&\ge0
       &&(f\in\mathcal D_{\mathrm A}),\\
 \langle\mathcal A_{n,\beta}^{\mathrm{DE}}u,u\rangle
           +\kappa_1\|u\|^2&\ge0
       &&(n\ge N_1),\\
 \langle\widetilde H_n^{\mathrm{red}}v,v\rangle
           +\kappa_2\|v\|^2&\ge0
       &&(n\ge N_2).
 \end{aligned}
\]
Indeed, the earlier bounds have the stronger form $c\ell_h-\kappa\|\cdot\|^2$ for the finite matrices and the analogous positive energy form for the continuum operator; dropping
the nonnegative energy gives these inequalities. For the reduced matrix, Lemma~\ref{A:res:auxiliary-boundary} identifies the complete discrete form with the form of the retained
principal matrix. The auxiliary edge coefficients cancel in that identity, so they produce no additional boundary term in the lower bound. The representation theorem for the
continuum form, and the Rayleigh quotient for the finite matrices, now give lower operator bounds with
\[
 K_\omega=\max\{0,\kappa_0,\kappa_1,\kappa_2\},\qquad
 N_* =\max\{N_1,N_2\}.
\]
Thus all three operators are bounded below by $-K_\omega I$ for every $n\ge N_*$.

For the full general matrix, put $F_n=\mathcal L_{n,3}^c\cup\mathcal V_{n,3}^c$. The finite-size estimates give $\sum_n\mathbb P(F_n)<\infty$. In particular
\[
 \mathbb P(F_n\text{ for infinitely many }n)
 \le\lim_{M\to\infty}\sum_{n\ge M}\mathbb P(F_n)=0.
\]
There is therefore a finite random $N_\partial$ such that both coefficient conditions hold for every $n\ge N_\partial$; independence across sizes is not used.

To substitute a random shift into the deterministic estimate, first take the countable intersection of the events of probability one on which Lemma~\ref{A:res:lower-transfer}
applies for $Q=3$ and $R=m$, for every integer $m\ge1$. On this intersection choose
\[
 M_\omega=\lceil K_\omega\rceil+2,\qquad a=-M_\omega.
\]
For $n\ge N_*$ the retained block satisfies $A_C-a\ge(M_\omega-K_\omega)I\ge2I$. Thus, for
\[
 N_\omega=\max\{N_*,N_\partial,n_0(3,M_\omega)\},
\]
the hypotheses of Lemma~\ref{A:res:lower-transfer} hold for every $n\ge N_\omega$, and it gives $\mathcal A_{n,V,\beta}>-M_\omega I$. The original three operators also exceed
$-M_\omega I$. Taking $L_\omega=M_\omega$ proves \eqref{A:eq:common-lower}.

All random quantities here may be chosen measurable. For example, each finite-matrix lower bound is an inequality for its smallest eigenvalue, and an index requiring it for every
subsequent size is determined by a countable intersection of such events. The integer-valued composition $n_0(3,M_\omega)$ is measurable. The conditional product coupling has the
prescribed marginal law in each model.
\end{proof}

Once the common lower bound is available, choose a real parameter below all the spectra. The shifted inverses are then positive, and their decreasing eigenvalues have the original
increasing labels. The next lemma makes the resulting estimate quantitative. Let $H$ be self-adjoint and bounded below, with compact resolvent, on a Hilbert space $\mathcal H$.
Write $\Lambda_1\le\Lambda_2\le\cdots\to\infty$ for its eigenvalues, repeated with multiplicity. Let $A_n$ be self-adjoint matrices of dimensions $m_n\to\infty$, with increasing
eigenvalues $\lambda_{k,n}$. Suppose that $J_n$ and $U$ are isometries of their respective spaces into one Hilbert space $\mathcal K$. In this paragraph put
\[
 S_n=J_n(A_n-\mathrm i)^{-1}J_n^*,\quad
 S=U(H-\mathrm i)^{-1}U^*,\quad
 \varepsilon_n=\|S_n-S\|_{\mathrm{op}}.
\]
Assume $A_n-a\ge\delta I$ and $H-a\ge\delta I$, with real $a$ and $\delta>0$, and set $b_{a,\delta}=1+|\mathrm i-a|/\delta$. The isometric images need not fill $\mathcal K$.

\begin{lemma}\label{A:res:ordered-transform}
Under these hypotheses,
\[
 \left|\frac1{\lambda_{k,n}-a}-\frac1{\Lambda_k-a}\right|
 \le b_{a,\delta}^{\,2}\varepsilon_n,
 \qquad k\le m_n.
 \eqtag{ordered-reciprocals}
\]
If $\varepsilon_n\to0$, then for each fixed $m$ and all sufficiently large $n$,
\[
 \max_{k\le m}|\lambda_{k,n}-\Lambda_k|
 \le2b_{a,\delta}^{\,2}(\Lambda_m-a)^2\varepsilon_n.
 \eqtag{ordered-rate}
\]
\end{lemma}

\begin{proof}
Set $Q_n=I+(\mathrm i-a)S_n$ and $Q=I+(\mathrm i-a)S$. On $\operatorname{Ran}J_n$, $Q_n^{-1}$ is unitarily equivalent to
\[
 (A_n-\mathrm i)(A_n-a)^{-1}
       =I+(a-\mathrm i)(A_n-a)^{-1}.
\]
It is the identity on the orthogonal complement. The same identities hold for $Q^{-1}$ with $H,U$ in place of $A_n,J_n$. The lower bounds therefore give
$\|Q_n^{-1}\|,\|Q^{-1}\|\le b_{a,\delta}$. The operators
\[
 T_n=S_nQ_n^{-1}=J_n(A_n-a)^{-1}J_n^*,\qquad
 T=SQ^{-1}=U(H-a)^{-1}U^*
\]
are compact, positive and self-adjoint. The same rational identity as in Proposition~\ref{A:res:20.2}, now with shift $\mathrm i-a$, gives
\[
 T_n-T=Q_n^{-1}(S_n-S)Q^{-1},\qquad
 \|T_n-T\|_{\mathrm{op}}
       \le b_{a,\delta}^{\,2}\varepsilon_n.
 \eqtag{positive-inverse-difference}
\]

For a positive compact operator $T$, its decreasing eigenvalues $t_k$, with zeros appended when its rank is finite, satisfy
\[
 t_k=\sup_{\substack{E\subset\mathcal K\\\dim E=k}}
        \ \inf_{\substack{u\in E\\\|u\|=1}}\langle Tu,u\rangle.
\]
The span of the first $k$ eigenvectors proves the lower bound when $t_k>0$. Every $k$-dimensional subspace intersects the orthogonal complement of the first $k-1$ eigenvectors,
proving the upper bound. If $T$ has rank less than $k$, every such subspace intersects its kernel, and positivity proves the value zero. This also proves that adjoining any zero
summand has no effect on the positive eigenvalues.

For each unit vector the two Rayleigh quotients differ by at most $\|T_n-T\|$. Taking the infimum and supremum in both directions gives $|t_{k,n}-t_k|\le\|T_n-T\|$. Positivity of
the shifted original operators identifies their decreasing inverse eigenvalues as
\[
 t_{k,n}=(\lambda_{k,n}-a)^{-1},\qquad
 t_k=(\Lambda_k-a)^{-1}.
\]
This proves \eqref{A:eq:ordered-reciprocals} for the actual increasing labels. For fixed $m$, choose $n$ so large that $b_{a,\delta}^{\,2}\varepsilon_n\le t_m/2$ and $m_n\ge m$.
Then $t_{k,n}\ge t_k/2$ for $k\le m$, and
\[
 |\lambda_{k,n}-\Lambda_k|
 =\frac{|t_{k,n}-t_k|}{t_{k,n}t_k}
 \le2t_k^{-2}b_{a,\delta}^{\,2}\varepsilon_n
 \le2(\Lambda_m-a)^2b_{a,\delta}^{\,2}\varepsilon_n.
\]
This is \eqref{A:eq:ordered-rate}. Neither simplicity nor the gap between adjacent eigenvalues is used.
\end{proof}

Ordered eigenvalues are controlled by positivity and min--max. For eigenvectors we also need to isolate the limiting spectral subspace. We do this with a contour around a nonzero
eigenvalue of the compact resolvent; zero must stay outside the contour. Return to $S=(\mathcal H_\beta-\mathrm i)^{-1}$ and $S_n=J_n(A_n-\mathrm i)^{-1}J_n^*$ for the two
isometric pairs. Their nonzero eigenvalues are $\mu_k=(\Lambda_k-\mathrm i)^{-1}$ and $\mu_{j,n}=(\lambda_{j,n}-\mathrm i)^{-1}$. Both operators are normal and compact, and
\[
 \|S_n-S\|_{\mathrm{op}}\le C_\omega h_n^{1/2}.
 \eqtag{21.5}
\]
For each fixed $k$, simplicity gives
\[
 d_k=\tfrac12\operatorname{dist}
    (\mu_k,\{0\}\cup\{\mu_j:j\ne k\})>0.
\]
Let $\Gamma_k$ be the positively oriented circle with center $\mu_k$ and radius $d_k$. Normality gives $\|(w-S)^{-1}\|\le d_k^{-1}$ on $\Gamma_k$. If $\|S_n-S\|<d_k/2$, factor
\[
 w-S_n=\bigl[I-(S_n-S)(w-S)^{-1}\bigr](w-S).
\]
The inverse of the first factor has norm at most two. Thus $\Gamma_k\subset\rho(S_n)$ and $\|(w-S_n)^{-1}\|\le2/d_k$ on the contour.

Let $\Pi_k,\Pi_{k,n}$ be the corresponding Riesz projections. The resolvent identity gives
\[
 \Pi_{k,n}-\Pi_k
 =\frac1{2\pi\mathrm i}\oint_{\Gamma_k}
      (w-S_n)^{-1}(S_n-S)(w-S)^{-1}\,\mathrm dw.
\]
It follows, both in operator norm and in Hilbert--Schmidt norm, that
\[
 \|\Pi_{k,n}-\Pi_k\|_{\mathrm{op}} \le\frac2{d_k}\|S_n-S\|_{\mathrm{op}},\qquad
 \|\Pi_{k,n}-\Pi_k\|_{\mathfrak S_2} \le\frac2{d_k}\|S_n-S\|_{\mathfrak S_2} \le
 C_{\omega,k}h_n^{1/2}.
\]
For each outcome and fixed $k$, choose $n_k$ with $2\|S_n-S\|_{\mathrm{op}}<d_k$ for every $n\ge n_k$. The projection rank argument then gives $\operatorname{rank}\Pi_{k,n}=1$.
This singles out one finite eigenvalue near $\mu_k$; the next theorem identifies its index in increasing order.

\par\Needspace{6\baselineskip}
\begin{theorem}\label{A:res:21.1}
For either isometric pair and every fixed $m$, almost surely, for all sufficiently large $n$,
\[
 \max_{k\le m}|\lambda_{k,n}-\Lambda_k|
 \le C_{\omega,m}h_n^{1/2}\le C'_{\omega,m}n^{-1/6}.
\]
\end{theorem}

\begin{proof}
Work on the intersection of the events of probability one in Proposition~\ref{A:res:common-lower} and \eqref{A:eq:21.2}. Fix a sample and $m\ge1$. There are $N_*$ and $C_*<\infty$
such that, for $n\ge N_*$,
\[
 A_n\ge-L_\omega I,\qquad \mathcal H_\beta\ge-L_\omega I,
 \qquad \|S_n-S\|_{\mathrm{op}}\le C_*h_n^{1/2}.
\]
Put $a=-L_\omega-1$ and $b=b_{a,1}$. Both operators shifted by $a$ are at least $I$. Increase $N_*$ to an index $N_m$ such that the finite dimension is at least $m$ and
\[
 b^2 C_*h_n^{1/2}\le\frac1{2(\Lambda_m-a)}
       \qquad(n\ge N_m).
\]
These requirements can be met because the dimension tends to infinity, $h_n\to0$, and $\Lambda_m-a$ is finite and positive on the fixed sample.

The positive shifted inverses in Lemma~\ref{A:res:ordered-transform} have decreasing eigenvalues $(\lambda_{k,n}-a)^{-1}$ and $(\Lambda_k-a)^{-1}$ with the same increasing energy
labels. The lemma gives
\[
 \left|(\lambda_{k,n}-a)^{-1}-(\Lambda_k-a)^{-1}\right|
      \le b^2C_*h_n^{1/2}\le\frac1{2(\Lambda_m-a)}.
\]
For $k\le m$, the finite reciprocal is consequently at least half the limiting reciprocal. Taking the difference of reciprocals therefore yields
\[
 \begin{split}
 |\lambda_{k,n}-\Lambda_k|
 &=\frac{|(\lambda_{k,n}-a)^{-1}-(\Lambda_k-a)^{-1}|}
          {(\lambda_{k,n}-a)^{-1}(\Lambda_k-a)^{-1}}\\
 &\le 2b^2(\Lambda_m-a)^2C_*h_n^{1/2}.
 \end{split}
\]
This proves the asserted rate for all $n\ge N_m$. The deterministic relation $h_n\asymp n^{-1/3}$ gives the second bound.

It remains to identify these labels with the local Riesz projections constructed above. For each fixed $k$,
\[
 \left|(\lambda_{k,n}-\mathrm i)^{-1}-(\Lambda_k-\mathrm i)^{-1}\right|
 =\frac{|\lambda_{k,n}-\Lambda_k|}
 {|\lambda_{k,n}-\mathrm i|\,|\Lambda_k-\mathrm i|}
 \le |\lambda_{k,n}-\Lambda_k|\longrightarrow0.
\]
Choose an index above which this difference is less than the radius of the disk bounded by $\Gamma_k$, and also above the index where that disk has Riesz rank one. There is exactly
one finite eigenvalue in the disk, including multiplicity, and it is the image of $\lambda_{k,n}$. Thus the ordered and local labels agree. Taking a countable intersection over $m$
gives one event for all fixed labels. The full general matrix is treated below in the enlarged Hilbert space.
\end{proof}

To obtain the total eigenvalue count on a bounded spectral interval, we must combine rank stability inside the chosen disks with exclusion of finite eigenvalues from the remaining
part of the interval. Fix $I=(a,b)$ with $a,b$ outside the limiting spectrum, and put
\[
 \Pi_I=\mathbf1_I(\mathcal H_\beta),\qquad
 \Pi_{I,n}=J_n\mathbf1_I(A_n)J_n^*.
\]
These are orthogonal projections. The following theorem, like the preceding pathwise comparisons, holds on the same event for all intervals $I=(a,b)$ satisfying
$a,b\in\rho(\mathcal H_\beta)$.

\begin{theorem}\label{A:res:21.2}
For all sufficiently large $n$,
\[
 \|\Pi_{I,n}-\Pi_I\|_{\mathfrak S_2}\le C_{\omega,I}h_n^{1/2},
 \qquad
 \#\{j:\lambda_{j,n}\in I\}=\#\{k:\Lambda_k\in I\}.
\]
\end{theorem}

\begin{proof}
Let $\gamma(t)=(t-\mathrm i)^{-1}$. There are only finitely many $\Lambda_k$ in $I$. Choose disjoint open disks about their images $\mu_k$, with closed disks avoiding zero and
every other spectral point of $S$. They can also be chosen so that
\[
 \gamma(\mathbb R)\cap\bigcup_k D_k\subset\gamma(I):
\]
the inverse map $w\mapsto\mathrm i+w^{-1}$ is continuous near each $\mu_k$, and $\Lambda_k$ has positive distance from $\{a,b\}$. Put $\Omega_I=\bigcup_kD_k$ and
$F_I=\gamma([a,b])\setminus\Omega_I$. The latter is compact and disjoint from $\sigma(S)$, hence has positive spectral distance if nonempty. For every $w$ with
$\operatorname{dist}(w,\sigma(S))>\|S_n-S\|$, the inverse series above gives $w\in\rho(S_n)$. Consequently no finite eigenvalue has image in $F_I$ for large $n$. All eigenvalues in
the disks have energies in $I$, and every finite eigenvalue in $I$ has image in the disks. In particular, $a,b\in\rho(A_n)$ for such $n$.

If $I$ contains no limiting eigenvalue, put $d_0=\operatorname{dist}(\gamma([a,b]),\sigma(S))>0$. For $n$ with
$\|S_n-S\|_{\mathrm{op}}<d_0$,
\[
 \sup_{w\in\gamma([a,b])}
       \|(S_n-S)(w-S)^{-1}\|_{\mathrm{op}}<1.
\]
The Neumann series gives $\gamma([a,b])\subset\rho(S_n)$. Hence $\Pi_{I,n}=0=\Pi_I$. Otherwise let $\Gamma_I=\partial\Omega_I$ and
$\delta_I=\operatorname{dist}(\Gamma_I,\sigma(S))>0$. The Riesz formula therefore represents exactly the two projections in the statement, and
\[
 \Pi_{I,n}-\Pi_I
 =\frac1{2\pi\mathrm i}\oint_{\Gamma_I}
      (w-S_n)^{-1}(S_n-S)(w-S)^{-1}\,\mathrm dw.
\]
For all sufficiently large $n$ the two contour resolvents have norms at most $2\delta_I^{-1}$ and $\delta_I^{-1}$. Hence
\[
 \|\Pi_{I,n}-\Pi_I\|_{\mathfrak S_2}
 \le\frac{|\Gamma_I|}{\pi\delta_I^2}\|S_n-S\|_{\mathfrak S_2}
 \le C_{\omega,I}h_n^{1/2}.
\]
Choose $n_I$ so that $C_{\omega,I}h_n^{1/2}<1$ for all $n\ge n_I$. The injectivity argument for orthogonal projections proves that their ranks are equal for these $n$, which is the
required equality of eigenvalue counts.
\end{proof}

Let $\psi_k$ be a unit eigenfunction of $\mathcal H_\beta$ and let $u_{k,n}$ be a unit eigenvector of $A_n$ for its $k$-th ordered eigenvalue. For large $n$, the eigenspace in
question is one-dimensional. Put $\phi_{k,n}=J_nu_{k,n}$ and choose its phase so that $\langle\phi_{k,n},\psi_k\rangle\ge0$. This phase is well defined for large $n$, since the
projections are close. For real eigenvectors it amounts to choosing their signs.

\begin{corollary}\label{A:res:21.3}
With this choice, for each fixed $k$, almost surely for all sufficiently large $n$,
\[
 \|J_nu_{k,n}-\psi_k\|_{L^2(\mathbb R_+)}
 \le C_{\omega,k}h_n^{1/2}.
\]
\end{corollary}

\begin{proof}
Write $c=\langle\phi_{k,n},\psi_k\rangle\in[0,1]$ and $t=(1-c^2)^{1/2}$. In the span of the two unit vectors, choosing the orthonormal basis $\psi_k,\eta$ with
$\phi_{k,n}=c\psi_k+t\eta$ when $t>0$, the difference of their projections has matrix
\[
 \begin{pmatrix}-t^2&ct\\ct&t^2\end{pmatrix}.
\]
Its eigenvalues are $t$ and $-t$. The same conclusions hold when $t=0$, with zero difference. In particular
\[
 \|\Pi_{k,n}-\Pi_k\|_{\mathrm{op}}=t,\qquad
 \|\Pi_{k,n}-\Pi_k\|_{\mathfrak S_2}^2=2t^2,\qquad
 \|\Pi_{k,n}-\Pi_k\|_{\mathfrak S_1}=2t.
\]
Also
\[
 \|\phi_{k,n}-\psi_k\|^2
 =2(1-c)\le2(1-c^2)
 =\|\Pi_{k,n}-\Pi_k\|_{\mathfrak S_2}^2.
\]
The Riesz projection estimate proves the corollary.
\end{proof}

\begin{corollary}\label{A:res:21.4}
Almost surely, for every $\varphi\in C_c(\mathbb R)$,
\[
 \|J_n\varphi(A_n)J_n^*-\varphi(\mathcal H_\beta)\|_{\mathfrak S_1}
 \longrightarrow0.
\]
Consequently $\operatorname{Tr}\varphi(A_n)\to \operatorname{Tr}\varphi(\mathcal H_\beta) =\sum_k\varphi(\Lambda_k)$.
\end{corollary}

\begin{proof}
Work on the event where Theorem~\ref{A:res:21.1} and Corollary~\ref{A:res:21.3} hold for every integer $k$; this is a countable intersection. On this fixed event, the following
argument applies to every test function in $C_c(\mathbb R)$. Choose $R$ with $\operatorname{supp}\varphi\subset[-R,R]$, and then $M$ with $\Lambda_M>R+1$. Ordered convergence gives
$\lambda_{M,n}>R$ for all large $n$. Hence all terms with index at least $M$ vanish in both spectral expansions. If $\Pi_{k,n}=\phi_{k,n}\phi_{k,n}^*$ and $\Pi_k=\psi_k\psi_k^*$,
\[
 J_n\varphi(A_n)J_n^*-\varphi(\mathcal H_\beta) ={}
 \sum_{k<M}[\varphi(\lambda_{k,n})-\varphi(\Lambda_k)]\Pi_{k,n}
 +\sum_{k<M}\varphi(\Lambda_k)(\Pi_{k,n}-\Pi_k).
\]
Each $\Pi_{k,n}$ has trace norm one. The rank-one calculation above gives $\|\Pi_{k,n}-\Pi_k\|_{\mathfrak S_1} \le C_{\omega,k}h_n^{1/2}$. With the modulus of continuity
\[
 \omega_\varphi(t)=\sup_{|x-y|\le t}|\varphi(x)-\varphi(y)|,
\]
we obtain the explicit bound
\[
 \|J_n\varphi(A_n)J_n^*-\varphi(\mathcal H_\beta)\|_{\mathfrak S_1}
 \le(M-1)\omega_\varphi(C_{\omega,M}h_n^{1/2})
       +C_{\omega,M}\|\varphi\|_\infty h_n^{1/2}\longrightarrow0.
\]

This argument applies to every compactly supported continuous function on the same fixed path. It also gives a uniform estimate for a fixed class of test functions. For $R>0$ and
$0<\alpha\le1$, choose $M$ from $R$ as above, independently of the function, and write $[\varphi]_\alpha=\sup_{x\ne y}|\varphi(x)-\varphi(y)|/|x-y|^\alpha$. For $n$ above an index
with $h_n\le1$, the preceding estimate gives
\[
 \sup_{\substack{\operatorname{supp}\varphi\subset[-R,R]\\
             \|\varphi\|_\infty+[\varphi]_\alpha\le1}}
 \|J_n\varphi(A_n)J_n^*-\varphi(\mathcal H_\beta)\|_{\mathfrak S_1}
 \le C_{\omega,R,\alpha}h_n^{\alpha/2}.
\]
The constant depends on the finitely many spectral gaps used in the projection bounds; its moments are not estimated here. The trace conclusion follows from
$|\operatorname{Tr}D|\le\|D\|_{\mathfrak S_1}$ and $J_n^*J_n=I$.
\end{proof}

Define $\xi_n=\sum_{j=1}^{m_n}\delta_{\lambda_{j,n}}$ and $\xi_\beta=\sum_{k\ge1}\delta_{\Lambda_k}$. The limiting measure is locally finite since $\Lambda_k\to\infty$.

\begin{theorem}\label{A:res:21.5}
On the constructed coupling, $\xi_n\to\xi_\beta$ almost surely in the vague topology. In particular $\xi_n\Rightarrow\xi_\beta$ in distribution.
\end{theorem}

\begin{proof}
For $\varphi\in C_c(\mathbb R)$, the integrals against the two point measures are exactly the traces in Corollary~\ref{A:res:21.4}. That corollary gives their convergence for every
such $\varphi$ on one event, which is pathwise vague convergence.

To check the measurable structure, enumerate as $(f_j)_{j\ge1}$ the continuous piecewise linear functions with compact support, rational breakpoints and rational values at those
breakpoints. This set includes, for each integer $r\ge1$, a nonnegative cutoff $\chi_r$ equal to one on $[-r,r]$ and zero outside $[-r-1,r+1]$. It is uniformly dense among
continuous functions supported in any fixed interval $(a,b)$ with $a,b\in\mathbb Z$. For locally finite positive measures, put
\[
 d(\mu,\nu)=\sum_{j\ge1}2^{-j}
          \min\{1,|\langle\mu-\nu,f_j\rangle|\}.
\]
The test functions distinguish measures, so $d$ is a metric. Vague convergence implies convergence in this metric by bounded convergence of the summands. Conversely, convergence in
$d$ implies convergence against every $f_j$, and the cutoff functions give a uniform bound on the masses of every bounded interval. If $\varphi$ is supported in $[-r,r]$,
approximate it uniformly by functions from this set supported in $[-r-1,r+1]$. For each such approximation $f_j$, the error is bounded by
\[
 |\langle\mu_n-\mu,\varphi-f_j\rangle|
 \le\|\varphi-f_j\|_\infty
       \bigl(\mu_n([-r-1,r+1])+\mu([-r-1,r+1])\bigr).
\]
First let $n\to\infty$ and then let the approximation error tend to zero. This proves that $d$ induces the vague topology. This topology is the subspace topology given by the
countably many coordinates $\langle\mu,f_j\rangle$ in $\mathbb R^{\mathbb N}$; it is second countable, and its Borel sets are generated by these coordinates. Each
$\langle\xi_n,f_j\rangle$ and $\langle\xi_\beta,f_j\rangle$ is measurable by the min--max formulas and the convergent eigenvalue sum. Thus these point measures are random elements
of this space.

For any bounded continuous functional $F$ on the space of locally finite measures with the vague topology, $F(\xi_n)\to F(\xi_\beta)$ almost surely. Bounded convergence gives
$\mathbb EF(\xi_n)\to\mathbb EF(\xi_\beta)$, which proves convergence in distribution.
\end{proof}

We now keep all eigenvalues of the general Jacobi matrix by using the isometry $\widehat{\mathcal J}_n$ from
Subsection~\ref{A:sec:atomlessness-real-resolvents-and-inverse-coupling}. Equation~\eqref{A:eq:12.43} bounds the resolvent difference on the full finite space, including the
removed coordinates and the cross terms. Let
\[
 A_n^{\mathrm{full}}=\mathcal A_{n,V,\beta},\qquad
 A_n^{\mathrm{red}}=\widetilde H_n^{\mathrm{red}},\qquad
 P_n=U_{r_n}P_{I_n}=\widetilde P_n.
\]
The nonreal boundary estimate \eqref{A:eq:12.43} gives
\[
 \begin{gathered}
 \|(A_n^{\mathrm{full}}-\mathrm i)^{-1}
       -P_n^*(A_n^{\mathrm{red}}-\mathrm i)^{-1}P_n\|_{\mathfrak S_2}^2
       =o_\omega(h_n),\\
 \|(A_n^{\mathrm{full}}-\mathrm i)^{-1}
       -P_n^*(A_n^{\mathrm{red}}-\mathrm i)^{-1}P_n\|_{\mathrm{op}}
       =o_\omega(h_n^{1/2}).
 \end{gathered}
 \eqtag{21.35}
\]
Set
\[
 \widehat S_n=\widehat{\mathcal J}_n
                (A_n^{\mathrm{full}}-\mathrm i)^{-1}
                      \widehat{\mathcal J}_n^*,\qquad
 \widehat S=(\mathcal H_\beta-\mathrm i)^{-1}\oplus0.
\]
Isometry preserves the first norm in \eqref{A:eq:21.35}. As in the real comparison, direct multiplication sends the extended reduced resolvent to
\[
 \bigl[\iota_h^{(M_n)}(A_n^{\mathrm{red}}-\mathrm i)^{-1}
                      (\iota_h^{(M_n)})^*\bigr]\oplus0.
\]
The reduced coupling and the squared triangle inequality thus give
\[
 \|\widehat S_n-\widehat S\|_{\mathfrak S_2}^2\le C_\omega h_n,
 \qquad
 \|\widehat S_n-\widehat S\|_{\mathrm{op}}\le C_\omega h_n^{1/2}.
 \eqtag{full-norm}
\]
The nonzero eigenvalues of $\widehat S_n$ are $(\lambda_{j,n}^{\mathrm{full}}-\mathrm i)^{-1}$, with multiplicity. Lemma~\ref{A:res:ordered-transform} applies to these operators:
the rational function used there vanishes at zero and therefore leaves the auxiliary zero summand unchanged.

\begin{corollary}\label{A:res:full-ordered}
For every fixed $m$, almost surely for all sufficiently large $n$,
\[
 \max_{k\le m}|\lambda_{k,n}^{(V),\mathrm{full}}-\Lambda_k|
 \le C_{\omega,V,\beta,m}h_n^{1/2}.
 \eqtag{full-ordered-rate}
\]
The projection, count, trace and point process conclusions above also hold for the full matrix, using $\widehat{\mathcal J}_n$ and the limiting inclusion $f\mapsto(f,0)$.
\end{corollary}

\begin{proof}
Let $\widehat{\mathcal H}=L^2(\mathbb R_+)\oplus\ell^2(\mathbb N)$ and let $U f=(f,0)$. Work on the event where Proposition~\ref{A:res:common-lower} and \eqref{A:eq:full-norm}
hold. Choose $a=-L_\omega-1$. For every sufficiently large $n$, both $A_n^{\mathrm{full}}-a$ and $\mathcal H_\beta-a$ are at least $I$. The rational function
\[
 F_a(t)=\frac{t}{1+(\mathrm i-a)t}
\]
vanishes at zero. On the isometric finite range it sends $(\lambda-\mathrm i)^{-1}$ to $(\lambda-a)^{-1}$, while on the orthogonal complement it leaves zero unchanged. Thus
Lemma~\ref{A:res:ordered-transform}, with $J_n=\widehat{\mathcal J}_n$ and this $U$, applies to every eigenvalue of the full matrix, not merely those of its retained principal
block. It yields \eqref{A:eq:full-ordered-rate}.

We verify the projection and count statements in this space. Fix a bounded interval $I=(a_0,b_0)$ whose endpoints lie in $\rho(\mathcal H_\beta)$, and put
\[
 \widehat\Pi_{I,n}=\widehat{\mathcal J}_n\mathbf1_I(A_n^{\mathrm{full}})
                    \widehat{\mathcal J}_n^*,\qquad
 \widehat\Pi_I=\mathbf1_I(\mathcal H_\beta)\oplus0.
\]
Under $\gamma(t)=(t-\mathrm i)^{-1}$, choose disjoint disks about the finitely many limiting eigenvalues in $I$, as in Theorem~\ref{A:res:21.2}. Their closed disks avoid zero and
all other spectral values. The compact remainder of $\gamma([a_0,b_0])$ has positive distance from $\sigma(\widehat S)$. If the operator-norm error is smaller than this distance,
the Neumann series excludes finite spectral values from the remainder. Hence the Riesz integrals over the disk boundaries represent exactly the two projections displayed above.

Let $\Gamma_I$ be the union of these boundaries and $\delta_I=\operatorname{dist}(\Gamma_I,\sigma(\widehat S))>0$. Since zero belongs to $\sigma(\widehat S)$, this distance
includes the distance to zero. For $w\in\Gamma_I$,
\[
 (w-\widehat S)^{-1}=(w-S)^{-1}\oplus w^{-1}I.
\]
The contour integral of the second term is zero, because the disks do not enclose zero. For $\|\widehat S_n-\widehat S\|<\delta_I/2$, the two contour resolvents have norms at most
$2/\delta_I$ and $1/\delta_I$, respectively. Therefore
\[
 \begin{split}
 \widehat\Pi_{I,n}-\widehat\Pi_I
 &=\frac1{2\pi\mathrm i}\oint_{\Gamma_I}
 (w-\widehat S_n)^{-1}(\widehat S_n-\widehat S)(w-\widehat S)^{-1}\,\mathrm dw,\\
 \|\widehat\Pi_{I,n}-\widehat\Pi_I\|_{\mathfrak S_2}
 &\le\frac{|\Gamma_I|}{\pi\delta_I^2}
             \|\widehat S_n-\widehat S\|_{\mathfrak S_2}
 \le C_{\omega,I}h_n^{1/2}.
 \end{split}
\]
Above an index for which the last bound is less than one, each projection is injective on the range of the other: a nonzero vector annihilated by the other projection would have
norm at most the norm difference times its own norm. Their finite ranks are therefore equal. If $I$ contains no limiting eigenvalue, the same Neumann exclusion on all of
$\gamma([a_0,b_0])$ gives $\widehat\Pi_{I,n}=0=\widehat\Pi_I$.

For a fixed simple limiting eigenvalue, put $\widehat\psi_k=U\psi_k$ and normalize $\widehat\Pi_{k,n}\widehat\psi_k$. The denominator is nonzero once the projection difference is
less than one. The rank-one calculation in Corollary~\ref{A:res:21.3} gives both
\[
 \|\widehat\phi_{k,n}-\widehat\psi_k\|\le C_{\omega,k}h_n^{1/2},
 \qquad
 \|\widehat\Pi_{k,n}-\widehat\Pi_k\|_{\mathfrak S_1}
       \le C_{\omega,k}h_n^{1/2}.
\]
For $\varphi\in C_c(\mathbb R)$ supported in $[-R,R]$, choose $M$ with $\Lambda_M>R+1$. The ordered convergence gives $\lambda_{M,n}^{\mathrm{full}}>R$ above a finite index. Thus
only $k<M$ contribute to either energy functional calculus, and
\[
 \begin{split}
 &\widehat{\mathcal J}_n\varphi(A_n^{\mathrm{full}})\widehat{\mathcal J}_n^*
                  -(\varphi(\mathcal H_\beta)\oplus0)\\
 &\quad=\sum_{k<M}[\varphi(\lambda_{k,n}^{\mathrm{full}})-\varphi(\Lambda_k)]
                        \widehat\Pi_{k,n}
       +\sum_{k<M}\varphi(\Lambda_k)(\widehat\Pi_{k,n}-\widehat\Pi_k).
 \end{split}
\]
Each finite projection has trace norm one. Consequently the trace norm of this difference is at most
\[
 (M-1)\omega_\varphi(C_{\omega,M}h_n^{1/2})
       +C_{\omega,M}\|\varphi\|_\infty h_n^{1/2},
\]
which tends to zero. The limiting expression is $U\varphi(\mathcal H_\beta)U^*=\varphi(\mathcal H_\beta)\oplus0$, even when $\varphi(0)\ne0$. The energy functional calculus is
taken on the original operator before its inclusion into the enlarged space. Taking traces proves convergence against every $\varphi\in C_c(\mathbb R)$ on the same sample. The
measurable vague-convergence argument of Theorem~\ref{A:res:21.5} then gives the point-process statement. Throughout, the auxiliary embedding is used only through its isometry
property.
\end{proof}

The preceding result has two useful consequences in the original $L^2(\mathbb R_+)$ space. Projection onto the first summand is a contraction, so
\[
 \|\mathcal J_n\mathbf1_I(\mathcal A_{n,V,\beta})\mathcal J_n^*
          -\mathbf1_I(\mathcal H_\beta)\|_{\mathfrak S_2}
 \le C_{\omega,I}h_n^{1/2}.
\]
The compressed finite operator in this display need not itself be a projection; the equality of spectral counts was proved for the orthogonal projections in the enlarged space.
Also, if $u_{k,n}^{\mathrm{full}}$ is a normalized full eigenvector, choose its phase by comparison with $(\psi_k,0)$. Then
\[
 \|\mathcal J_nu_{k,n}^{\mathrm{full}}-\psi_k\|_{L^2}^2
 +\|(I-P_n^*P_n)u_{k,n}^{\mathrm{full}}\|^2 =\|\widehat{\mathcal J}_nu_{k,n}^{\mathrm{full}}
 -(\psi_k,0)\|^2 \le C_{\omega,k}h_n.
\]
Thus the omitted coordinates carry at most $O_\omega(h_n)$ of the squared mass of any fixed normalized edge eigenvector.

For either full ensemble and fixed deterministic $a<b$, Theorem~\ref{A:res:20.1} gives $a,b\in\rho(\mathcal H_\beta)$ almost surely. The count equality consequently gives
\[
 \#\{j:\lambda_{j,n}\in(a,b)\}
       =\#\{k:\Lambda_k\in(a,b)\}
\]
for all sufficiently large $n$. These statements hold simultaneously for any specified countable collection of intervals. On the pathwise event above they also hold for every
interval whose boundary is disjoint from the spectrum for the fixed path. For each fixed $k$, the ordered estimates imply
\[
 |(\lambda_{k+1,n}-\lambda_{k,n})-(\Lambda_{k+1}-\Lambda_k)|
 \le |\lambda_{k+1,n}-\Lambda_{k+1}|
       +|\lambda_{k,n}-\Lambda_k|
 \le C_{\omega,k}h_n^{1/2}.
\]
Since $\Lambda_{k+1}-\Lambda_k>0$, the finite gap is at least half that number for large $n$. For bounded continuous $F:\mathbb R^m\to\mathbb R$, bounded convergence further gives
\[
 \mathbb EF(\lambda_{1,n},\ldots,\lambda_{m,n})
      \longrightarrow\mathbb EF(\Lambda_1,\ldots,\Lambda_m).
\]
The exact finite marginals in the coupling make this a statement about the original ensembles.

To give the normalization and sign explicitly, let $\theta_{1,n}\ge\cdots\ge\theta_{n,n}$ be the eigenvalues of the unscaled Dumitriu--Edelman matrix. Then
\[
 \lambda_{k,n}=n^{2/3}(2-\theta_{k,n}),\qquad
 \sum_{k=1}^n\delta_{n^{2/3}(2-\theta_{k,n})}
          \longrightarrow\sum_{k\ge1}\delta_{\Lambda_k}
\]
almost surely vaguely on this coupling. For the general model, write $\theta_{1,n}^{(V)}\ge\cdots\ge\theta_{n,n}^{(V)}$ for the full unscaled eigenvalues. The definition of $h$
gives
\[
 \frac{h^{-2}}{b_\star}
   =b_\star^{-1/3}\tau_\star^{-2/3}n^{2/3},\qquad
 \lambda_{k,n}^{(V),\mathrm{full}}
   =b_\star^{-1/3}\tau_\star^{-2/3}n^{2/3}
          (E_\star-\theta_{k,n}^{(V)}).
\]
The same point process limit holds for these full scaled variables, by Corollary~\ref{A:res:full-ordered}. Our limiting point measure lists the increasing eigenvalues of the
stochastic Airy operator. With the usual descending soft edge convention it is reflected, and the largest eigenvalue limit is
\[
 n^{2/3}(\theta_{1,n}-2)\ \Rightarrow\ -\Lambda_1,
 \qquad
 b_\star^{-1/3}\tau_\star^{-2/3}n^{2/3}
       (\theta_{1,n}^{(V)}-E_\star)\ \Rightarrow\ -\Lambda_1.
\]
The law of $-\Lambda_1$ is the Tracy--Widom $\beta$ law \cite{ramirez2011beta}; their lowest eigenvalue is indexed by zero. The opposite centered quantities converge to
$\Lambda_1$.

The rate $h_n^{1/2}\asymp n^{-1/6}$ for individual eigenvalues is a consequence of the global Hilbert--Schmidt estimate and the lower spectral bound. Proposition~\ref{A:res:5.1}
gives the matching lower bound for the full squared Hilbert--Schmidt error. It does not give a matching lower bound for the error of an individual eigenvalue.

\subsection{Spectral consequences of the inverse approximation and further questions}
\label{B:sec:spectral-consequences-and-further-questions}

The positive inverse permits a global comparison of reciprocal eigenvalues. Direct eigenvalues and eigenfunctions then follow locally, with constants depending on the distance from
zero and on the relevant spectral gaps. The same inverse error also controls powers and regularized determinants.

Fix either coupling in Theorems~\ref{B:result:1-1} and \ref{B:result:1-3}. The eigenvalues of the scaled finite matrix $\mathcal A_n^{\mathfrak m}$,
$\mathfrak m\in\{\mathrm L,V\}$, are written in increasing order as $\lambda_{1,n}^{\mathfrak m},\ldots,\lambda_{n,n}^{\mathfrak m}$. The limiting eigenvalues $\Lambda_j$ and their
real eigenspaces are those of Corollary~\ref{B:cor:bessel-simple-spectrum}. Use $\mathscr J_n^{\mathrm L}=\mathscr E_n$ and $\mathscr J_n^V=\mathscr J_{n,V}$ for the embeddings
already defined. In this subsection abbreviate the embedded inverses and their bounds by
\begin{equation}
 \begin{gathered}
 R_n=\mathscr J_n^{\mathfrak m}(\mathcal A_n^{\mathfrak m})^{-1}
                      (\mathscr J_n^{\mathfrak m})^*,\qquad
 R=\mathcal H_{\beta,a}^{-1},\qquad
 e_n=\|R_n-R\|_{\mathfrak S_2},\\
 Z=\sup_{n\ge1}ne_n^2,\qquad
 M=\max\{1,\|R\|_{\rm op},\sup_n\|R_n\|_{\rm op}\},\qquad
 H_2=\max\{\|R\|_{\mathfrak S_2},\sup_n\|R_n\|_{\mathfrak S_2}\}.
 \end{gathered}
 \label{B:eq:spectral-setup}
\end{equation}
The preceding theorems give $Z,M,H_2<\infty$ on an event of probability one, denoted by $\Omega_0$. On this event $R$ is positive, injective and trace class, and its positive
eigenvalues are simple. All estimates below are deterministic on $\Omega_0$. In particular,
\[
 e_n\le\sqrt{Z/n},\qquad
 H_2\le\|R\|_{\mathfrak S_2}+\sqrt Z,\qquad
 M\le1+\|R\|_{\mathfrak S_2}+\sqrt Z.
\]
Passing to a fixed direct eigenvalue introduces its positive reciprocal in the denominator, and comparing eigenfunctions also requires a spectral gap. Square summability of the
reciprocal spectrum gives the regularized determinant. Statements about expectations are treated separately at the end of this subsection.

The eigenvalue estimate below is the Hoffman--Wielandt inequality for positive compact operators \cite{bhatia1994hoffman} (Theorem~2); see also \cite{valko2020operator} for its use
in an operator limit. Spectral projections are treated by the Riesz integral \cite{kato1995perturbation} (Chapter~III, \S6.4), and determinants by the regularized product of
Theorem~6.2 of \cite{simon1977determinants}.

\subsubsection{Eigenvalues of the inverse}
\label{B:subsection:reciprocal-spectra}

Ordering the positive inverse eigenvalues decreasingly aligns them with the increasing direct spectrum. Extending the finite list by zeros includes the remaining eigenvalues of the
limiting inverse in the comparison.

Put
\begin{equation}
 \mu_j=\Lambda_j^{-1},\qquad
 \mu_{j,n}=\begin{cases}
  (\lambda_{j,n}^{\mathfrak m})^{-1},&j\le n,\\
  0,&j>n.
 \end{cases}
 \label{B:eq:spectral-reciprocal-eigenvalues}
\end{equation}
These sequences are in decreasing order.

\begin{lemma}\label{B:lem:ordered-reciprocal-eigenvalues}
For every $n$,
\begin{equation}
 \sum_{j\ge1}|\mu_{j,n}-\mu_j|^2\le e_n^2,
 \qquad
 \sup_{j\ge1}|\mu_{j,n}-\mu_j|
                    \le\|R_n-R\|_{\rm op}.
 \label{B:eq:spectral-matching}
\end{equation}
\end{lemma}
\begin{proof}
We first prove the partial trace inequality. Let $(\psi_j)$ be an orthonormal spectral basis for the positive compact operator $R$, and let $u_1,\ldots,u_m$ be orthonormal. Put
\[
 a_j=\sum_{i=1}^m|\langle u_i,\psi_j\rangle|^2.
\]
Bessel's inequality gives $0\le a_j\le1$ and $\sum_j a_j\le m$; equality holds if there is no component in the kernel. Since $\mu_j$ is decreasing and nonnegative,
\[
 \begin{split}
 \sum_j\mu_ja_j
 &\le\sum_{j\le m}\mu_ja_j+\mu_m\sum_{j>m}a_j\\
 &\le\sum_{j\le m}\mu_ja_j+\mu_m\sum_{j\le m}(1-a_j)
 \le\sum_{j\le m}\mu_j.
 \end{split}
\]
All sums are justified by nonnegativity. Expanding the Rayleigh quotients in the spectral basis proves
\begin{equation}
 \sum_{i=1}^m\langle Ru_i,u_i\rangle\le\sum_{j=1}^m\mu_j.
 \label{B:eq:spectral-partial-trace}
\end{equation}
This proof applies to any positive compact operator, with zeros appended to a finite eigenvalue list.

Let $v_1,\ldots,v_n$ be eigenvectors of $R_n$ corresponding to its decreasing positive eigenvalues $\mu_{i,n}$, and put $\mu_{n+1,n}=0$. Summation by parts gives
\[
 \begin{split}
 \operatorname{Tr}(R_nR)
 &=\sum_{i=1}^n\mu_{i,n}\langle Rv_i,v_i\rangle\\
 &=\sum_{k=1}^n(\mu_{k,n}-\mu_{k+1,n})
                         \sum_{i=1}^k\langle Rv_i,v_i\rangle\\
 &\le\sum_{k=1}^n(\mu_{k,n}-\mu_{k+1,n})\sum_{i=1}^k\mu_i
 =\sum_{i=1}^n\mu_{i,n}\mu_i.
 \end{split}
\]
The coefficients in the middle sum are nonnegative, so the partial trace bound can be applied term by term. Both $R_n$ and $R$ are Hilbert--Schmidt; their product is trace class.
Expanding the square in the Hilbert--Schmidt inner product now yields
\[
 \begin{split}
 e_n^2
 &=\sum_{i=1}^n\mu_{i,n}^2+\sum_{i\ge1}\mu_i^2
                      -2\operatorname{Tr}(R_nR)\\
 &\ge\sum_{i=1}^n(\mu_{i,n}-\mu_i)^2+\sum_{i>n}\mu_i^2
 =\sum_{i\ge1}(\mu_{i,n}-\mu_i)^2.
 \end{split}
\]
This proves the first inequality, including the reciprocal tail omitted by a comparison of only the first $n$ entries. It is the positive case of the Hoffman--Wielandt inequality
\cite{bhatia1994hoffman}; the calculation above proves the needed version directly.

For completeness, a positive compact operator $T$ satisfies
\begin{equation}
 \mu_j(T)=\inf_{\dim F=j-1}\ 
       \sup_{\substack{v\perp F\\\|v\|_2=1}}\langle Tv,v\rangle.
 \label{B:eq:spectral-minmax}
\end{equation}
If $T$ has at least $j$ positive eigenvalues, take $F$ to be the span of its first $j-1$ eigenvectors for the upper bound. For every other such $F$, its orthogonal complement
intersects the span of the first $j$ eigenvectors in a nonzero vector. The Rayleigh quotient of a unit vector in that intersection is at least $\mu_j(T)$, proving the lower bound.
If the positive rank is less than $j$, choose $F$ containing the entire range of $T$ and extend it to dimension $j-1$. The supremum is then zero; positivity gives the reverse
inequality. The infinite-dimensional ambient space supplies the required additional vectors.

For every unit vector, the Rayleigh quotients of $R_n$ and $R$ differ by at most $\|R_n-R\|_{\mathrm{op}}$. Taking the supremum and infimum in \eqref{B:eq:spectral-minmax}
preserves this bound. Interchanging $R_n$ and $R$ proves the second inequality of \eqref{B:eq:spectral-matching} for every $j$.
\end{proof}

\begin{corollary}\label{B:cor:fixed-bessel-eigenvalues}
Fix $m\ge1$. If $e_n<\mu_m/2$, then $n\ge m$ and
\begin{equation}
 \sum_{j=1}^m|\lambda_{j,n}^{\mathfrak m}-\Lambda_j|^2
          \le4\Lambda_m^4e_n^2,\qquad
 \max_{j\le m}|\lambda_{j,n}^{\mathfrak m}-\Lambda_j|
          \le2\Lambda_m^2e_n.
 \label{B:eq:spectral-fixed-eigenvalues}
\end{equation}
\end{corollary}
\begin{proof}
For $j\le m$, Lemma~\ref{B:lem:ordered-reciprocal-eigenvalues} gives $\mu_{j,n}\ge\mu_j-e_n>\mu_j/2>0$. If $n<m$, the convention $\mu_{m,n}=0$ would imply
$e_n\ge|\mu_{m,n}-\mu_m|=\mu_m$, contradicting the hypothesis. Thus all the first $m$ finite eigenvalues exist. For $j\le m$ the actual denominator satisfies
\[
 \mu_{j,n}\mu_j\ge\mu_j(\mu_j-e_n)>\tfrac12\mu_j^2
 \ge\tfrac12\mu_m^2.
\]
Moreover,
\[
 |\lambda_{j,n}^{\mathfrak m}-\Lambda_j|
 =\frac{|\mu_{j,n}-\mu_j|}{\mu_{j,n}\mu_j}
 \le2\mu_j^{-2}|\mu_{j,n}-\mu_j|
 \le2\Lambda_m^2|\mu_{j,n}-\mu_j|.
\]
Squaring the last inequality and summing, then applying \eqref{B:eq:spectral-matching}, gives
\[
 \sum_{j=1}^m|\lambda_{j,n}^{\mathfrak m}-\Lambda_j|^2
 \le4\Lambda_m^4\sum_{j=1}^m|\mu_{j,n}-\mu_j|^2
 \le4\Lambda_m^4e_n^2.
\]
Also $|\mu_{j,n}-\mu_j|\le e_n$ for every $j$, so taking the maximum in the unsquared inequality gives $\max_{j\le m}|\lambda_{j,n}^{\mathfrak m}-\Lambda_j| \le2\Lambda_m^2e_n$.
For these same labels put $d_{j,n}=\mu_{j,n}-\mu_j$. Exact reciprocation gives the more precise identity
\[
 \lambda_{j,n}^{\mathfrak m}-\Lambda_j+\Lambda_j^2d_{j,n}
 =\frac{d_{j,n}^2}{\mu_j^2\mu_{j,n}},\qquad
 \max_{j\le m}
 |\lambda_{j,n}^{\mathfrak m}-\Lambda_j+\Lambda_j^2d_{j,n}|
 \le2\Lambda_m^3 e_n^2.
\]
For fixed $m$ the scalar reciprocation remainder is therefore $O_\omega(n^{-1})$. Identifying a fluctuation law requires joint convergence of the reciprocal errors $d_{j,n}$, which
does not follow from the inverse norm bound. Both this identity and the preceding estimate for ordered eigenvalues apply without a gap between adjacent positive eigenvalues. The
eigenvector comparison below will require such a gap.
\end{proof}

\subsubsection{Eigenvalue counts and linear statistics}
\label{B:subsection:bessel-eigenvalue-counts}

A direct spectral threshold corresponds to a positive threshold for the inverse. If it lies away from the limiting spectrum, the reciprocal comparison fixes the count whenever its
operator norm error is less than half that spectral distance. The resulting finite sums determine the vague limit of the point measures.

The locally finite point measures are
\begin{equation}
 \Xi_n^{\mathfrak m}=\sum_{j=1}^n\delta_{\lambda_{j,n}^{\mathfrak m}},
 \qquad \Xi=\sum_{j\ge1}\delta_{\Lambda_j}.
 \label{B:eq:spectral-point-measures}
\end{equation}
Local finiteness of $\Xi$ follows from $\Lambda_j\to\infty$.

\begin{proposition}\label{B:prop:bessel-eigenvalue-counts}
On $\Omega_0$, for every $E>0$ outside $\sigma(\mathcal H_{\beta,a})$,
\begin{equation}
 \#\{j\le n:\lambda_{j,n}^{\mathfrak m}<E\}
       =\#\{j:\Lambda_j<E\}
 \label{B:eq:spectral-counts}
\end{equation}
for all sufficiently large $n$. On the same event, $\Xi_n^{\mathfrak m}\to\Xi$ vaguely on $(0,\infty)$.
\end{proposition}
\begin{proof}
Set $t=E^{-1}$ and $d_E=\operatorname{dist}(t,\{0\}\cup\{\mu_j:j\ge1\})$. Choose $L$ with $\mu_l<t/2$ for $l>L$. The distance from $t$ to every such eigenvalue is at least $t/2$,
while the finitely many numbers $|t-\mu_l|$, $l\le L$, are positive. Their minimum with $t/2$ is a positive lower bound for $d_E$. Also
\[
 \#\{j:\mu_j>t\}\le t^{-2}\sum_j\mu_j^2
 =E^2\|R\|_{\mathfrak S_2}^2<\infty.
\]
 If $\|R_n-R\|_{\rm op}<d_E/2$, then for every $j$, \eqref{B:eq:spectral-matching} implies
\[
 \mu_j>t\ \Longrightarrow\ \mu_{j,n}>t,
 \qquad
 \mu_j<t\ \Longrightarrow\ \mu_{j,n}<t.
\]
Hence the two sets of indices with reciprocal eigenvalue greater than $t$ agree. They are finite, and the finite set contains no index greater than $n$. This proves
\eqref{B:eq:spectral-counts}.

Let $f\in C_c((0,\infty))$ be nonzero. Choose $E>\sup\operatorname{supp}f$ outside the limiting spectrum, and let $m=\#\{j:\Lambda_j<E\}$. For all sufficiently large $n$, the
difference of the integrals in \eqref{B:eq:spectral-point-measures} is
\[
 \int f\,\dd\Xi_n^{\mathfrak m}-\int f\,\dd\Xi
       =\sum_{j=1}^m\bigl(f(\lambda_{j,n}^{\mathfrak m})-f(\Lambda_j)\bigr).
\]
If $m=0$, the sum is zero. Otherwise extend $f$ by zero to $[0,E+1]$ and let
\[
 \omega_f(\delta)=\sup\{|f(x)-f(y)|:x,y\in[0,E+1],\ |x-y|\le\delta\}.
\]
This modulus tends to zero with $\delta$. For all sufficiently large $n$, Corollary~\ref{B:cor:fixed-bessel-eigenvalues} gives
\[
 \left|\int f\,\dd\Xi_n^{\mathfrak m}-\int f\,\dd\Xi\right|
 \le m\,\omega_f(2\Lambda_m^2e_n)\longrightarrow0.
\]
The argument holds for every $f\in C_c((0,\infty))$ on the fixed path in $\Omega_0$, and therefore proves vague convergence on that event.
\end{proof}

\begin{corollary}\label{B:cor:bessel-linear-statistics}
Let $f\in C_c^1((0,\infty))$ and $B>0$ satisfy $\operatorname{supp}f\subset(0,B)$. For every $n$,
\begin{equation}
 \left|\int f\,\dd\Xi_n^{\mathfrak m}-\int f\,\dd\Xi\right|
 \le B^3\|f'\|_\infty
       (\|R_n\|_{\mathfrak S_2}^2+\|R\|_{\mathfrak S_2}^2)^{1/2}e_n.
 \label{B:eq:spectral-test-functions}
\end{equation}
\end{corollary}
\begin{proof}
Extend $f$ by zero outside $(0,\infty)$ and set $g(t)=f(t^{-1})$ for $t>0$, $g(0)=0$. The support condition makes $g$ continuously differentiable, identically zero near zero and
for all sufficiently large $t$. Where $g'$ is nonzero, $t\ge B^{-1}$, so
\[
 g'(t)=-t^{-2}f'(t^{-1}),\qquad
                  \|g'\|_\infty\le B^2\|f'\|_\infty.
\]
Let $\mathcal I_n(f)$ be the set of indices for which $g(\mu_{j,n})\ne g(\mu_j)$. For each such index at least one of $\mu_{j,n},\mu_j$ is greater than $B^{-1}$. Consequently
\[
 |\mathcal I_n(f)|\le B^2\sum_{j\in\mathcal I_n(f)}(\mu_{j,n}^2+\mu_j^2)
          \le B^2(\|R_n\|_{\mathfrak S_2}^2+\|R\|_{\mathfrak S_2}^2).
\]
The mean-value inequality and Cauchy--Schwarz now give
\[
 \begin{split}
 \left|\sum_{j\ge1}(g(\mu_{j,n})-g(\mu_j))\right|
 &\le B^2\|f'\|_\infty\sum_{j\in\mathcal I_n(f)}|\mu_{j,n}-\mu_j|\\
 &\le B^2\|f'\|_\infty |\mathcal I_n(f)|^{1/2}
                  \left(\sum_{j\ge1}|\mu_{j,n}-\mu_j|^2\right)^{1/2}.
 \end{split}
\]
The sum is finite because the preceding cardinality bound is finite; it also covers $\mathcal I_n(f)=\varnothing$, when both sides are zero. Substitution of that bound and
\eqref{B:eq:spectral-matching} gives exactly
\[
 B^2\|f'\|_\infty\, B
 (\|R_n\|_{\mathfrak S_2}^2+\|R\|_{\mathfrak S_2}^2)^{1/2}e_n.
\]
The constant applies to every $n$, including when the boundary of the support meets the limiting spectrum.
\end{proof}

\subsubsection{Eigenfunctions and spectral projections}
\label{B:subsection:bessel-eigenfunctions-and-projections}

An individual eigenfunction can be compared only after its eigenspace is isolated. The gap below includes the distance from the inverse eigenvalue to zero as well as to the other
eigenvalues, which makes the projection estimate valid for the finite-rank approximations.

Let $\psi_j$ be the real unit eigenfunction fixed in Corollary~\ref{B:cor:bessel-simple-spectrum}, and let $\Pi_j$ project onto its span. Set
\begin{equation}
 g_j=\frac14\inf\bigl(\{\mu_j\}\cup
                   \{|\mu_j-\mu_l|:l\ne j\}\bigr).
 \label{B:eq:spectral-phase-gap}
\end{equation}
Choose $L$ so that $\mu_l\le\mu_j/2$ for $l>L$. Then the infimum is bounded below by the minimum of $\mu_j/2$ and the finitely many positive gaps $|\mu_j-\mu_l|$ with $l\le L$,
$l\ne j$. In particular, $g_j>0$. For $n\ge j$, let $v_{j,n}^{\mathfrak m}$ be the real unit eigenvector of $\mathcal A_n^{\mathfrak m}$ with first coordinate positive, and set
\[
 \psi_{j,n}^{\mathfrak m}=\mathscr J_n^{\mathfrak m}v_{j,n}^{\mathfrak m},
 \qquad \Pi_{j,n}^{\mathfrak m}
     =\text{the orthogonal projection onto }\mathbb C\psi_{j,n}^{\mathfrak m}.
\]

\begin{lemma}\label{B:lem:bessel-spectral-measurability}
Every $\mathcal A_n^{\mathfrak m}$ has simple spectrum. For $n\ge j$, the maps from its coefficients to $v_{j,n}^{\mathfrak m}$ and $\Pi_{j,n}^{\mathfrak m}$ are continuous.
\end{lemma}
\begin{proof}
Write $d_k$ and $c_k$ for the diagonal and off-diagonal entries of $\mathcal A_n^{\mathfrak m}$. All $c_k$ are nonzero, since they are negative multiples of $X_{n,k}Y_{n,k}$. An
eigenvector whose first coordinate is zero has second coordinate zero by the first row equation. Induction in the remaining rows then makes every coordinate zero. Thus the map
taking an eigenvector to its first coordinate is injective on each eigenspace. The real symmetric matrix is diagonalizable, so all its eigenvalues are simple.

For an eigenvalue $\lambda=\lambda_{j,n}^{\mathfrak m}$, put $w_0=0$, $w_1=1$, and $c_0=0$. The recursion
\begin{equation}
 w_{k+1}=\frac{(\lambda-d_k)w_k-c_{k-1}w_{k-1}}{c_k},
                         \qquad 1\le k<n,
 \label{B:eq:finite-eigenvector-recursion}
\end{equation}
gives the eigenvector with first coordinate one: any eigenvector can be scaled to have this first coordinate, and the recursion then determines it uniquely. Hence its last row
equation also holds. Normalize by $(\sum_{k=1}^nw_k^2)^{1/2}$. For two real symmetric matrices $A_1,A_2$ of this fixed size, $|\langle(A_1-A_2)v,v\rangle|\le\|A_1-A_2\|_{\rm op}$
on the unit sphere. Applying the finite-dimensional min--max formula in both directions gives
\[
 |\lambda_j(A_1)-\lambda_j(A_2)|\le\|A_1-A_2\|_{\rm op}.
\]
In a neighborhood of a fixed coefficient array every denominator $c_k$ stays separated from zero. The normalized vector has denominator $\|w\|_2\ge w_1=1$. The recursion and
normalization are therefore continuous on the set of coefficient matrices under consideration. For $n=1$ the normalized vector is $(1)$.

The finite embedding is deterministic. For unit vectors $u,v$,
\[
 (P_u-P_v)h=\langle h,u\rangle(u-v)+\langle h,u-v\rangle v.
\]
For the operator $Th=\langle h,a\rangle b$, Parseval gives $\|T\|_{\mathfrak S_2}^2=\sum_k|\langle e_k,a\rangle|^2\|b\|_2^2 =\|a\|_2^2\|b\|_2^2$. The preceding decomposition
therefore yields
\[
 \|P_u-P_v\|_{\mathfrak S_2}\le2\|u-v\|_2.
\]
This proves continuity of the embedded projections.
\end{proof}

\begin{lemma}\label{B:lem:bessel-limiting-spectral-measurability}
The variables $g_j$, $\Pi_j$, and $\psi_j$ are measurable in $\mathbb R$, $\mathfrak S_2$, and $L^2(0,1)$, respectively.
\end{lemma}
\begin{proof}
Lemma~\ref{B:lem:bessel-path-measurability} makes $p^B$ and $R$ measurable in $L^2$ and $\mathfrak S_2$. The min--max estimate \eqref{B:eq:spectral-minmax} makes every $\mu_j$
continuous in operator norm on the positive compact operators. The countable infimum defining $g_j$ is therefore measurable.

Put $\zeta_j(t)=\mu_j+g_j\me^{\ii t}$ and $S_R(\zeta)=\zeta^{-1}(\zeta I-R)^{-1}R$. Since $\mu_j\ge4g_j$, the circle has distance at least $g_j$ from $\sigma(R)$ and
$|\zeta_j(t)|\ge3g_j$. The spectral theorem and the scalar contour integrals give
\begin{equation}
 \Pi_j=\frac1{2\pi}\int_0^{2\pi}
        g_j\me^{\ii t}S_R(\zeta_j(t))\,\dd t,
 \qquad
 S_R(\zeta)=(\zeta I-R)^{-1}-\zeta^{-1}I.
 \label{B:eq:bessel-measurable-projection}
\end{equation}
Indeed, on the eigenspace with eigenvalue $\mu_l$ the integral is
\[
 \frac1{2\pi\ii}\oint_{|\zeta-\mu_j|=g_j}
          \left(\frac1{\zeta-\mu_l}-\frac1\zeta\right)\,\dd\zeta
       =\mathbf1_{\{l=j\}}.
\]
The integrand takes values in $\mathfrak S_2$. Its derivative is bounded using
\[
 \begin{split}
 S_R'(\zeta)
 &=-\zeta^{-2}(\zeta I-R)^{-1}R
   -\zeta^{-1}(\zeta I-R)^{-2}R,\\
 \|S_R(\zeta_j(t))\|_{\mathfrak S_2}
 &\le\frac{\|R\|_{\mathfrak S_2}}{3g_j^2},\\
 \|S_R'(\zeta_j(t))\|_{\mathfrak S_2}
 &\le\frac{4\|R\|_{\mathfrak S_2}}{9g_j^3}.
 \end{split}
\]
The integrand subtracts the scalar resolvent because the identity operator is not Hilbert--Schmidt on $L^2(0,1)$. The subtracted term has integral zero, since the circle does not
enclose zero. In particular all the Riemann sums below belong to $\mathfrak S_2$. The chain rule gives
\[
 \Phi_j'(t)=\ii g_j\me^{\ii t}S_R(\zeta_j(t))
       +\ii g_j^2\me^{2\ii t}S_R'(\zeta_j(t)),
 \qquad \Phi_j(t)=g_j\me^{\ii t}S_R(\zeta_j(t)).
\]
Thus $\sup_t\|\Phi_j'(t)\|_{\mathfrak S_2} \le(1/3+4/9)\|R\|_{\mathfrak S_2}/g_j \le\|R\|_{\mathfrak S_2}/g_j$. In particular its left Riemann sums satisfy
\begin{equation}
 \left\|\Pi_j-\frac1N\sum_{k=0}^{N-1}\Phi_j(2\pi k/N)\right\|_{\mathfrak S_2}
 \le\frac{2\pi\|R\|_{\mathfrak S_2}}{Ng_j}.
 \label{B:eq:spectral-projection-riemann-error}
\end{equation}
For fixed $t$, each summand is measurable as an $\mathfrak S_2$-valued map. To check this in the relevant norm, use
\[
 (zI-T)^{-1}-(wI-S)^{-1}
  =(zI-T)^{-1}[(w-z)I+T-S](wI-S)^{-1}.
\]
On a neighborhood where the two resolvents and $z^{-1},w^{-1}$ are bounded, multiplication by $T,S\in\mathfrak S_2$ gives
\[
 \|z^{-1}(zI-T)^{-1}T-w^{-1}(wI-S)^{-1}S\|_{\mathfrak S_2}
       \le C\bigl(\|T-S\|_{\mathfrak S_2}+|z-w|\bigr).
\]
The constant is local in $S,w$. This proves the required continuity, and \eqref{B:eq:spectral-projection-riemann-error} proves measurability of $\Pi_j$.

The normalization from Section~\ref{B:sec:inverse-representations-and-the-stochastic-bessel-operator} and \eqref{B:eq:2-31c} give
\[
 \Pi_jp^B=\langle p^B,\psi_j\rangle\psi_j,
 \qquad \langle p^B,\psi_j\rangle>0,
 \qquad \psi_j=\frac{\Pi_jp^B}{\|\Pi_jp^B\|_2}.
\]
Evaluation is continuous because $\|Tp-Sq\|_2\le\|T-S\|_{\mathfrak S_2}\|p\|_2+ \|S\|_{\rm op}\|p-q\|_2$. Normalization is continuous away from zero. These facts prove
measurability of $\psi_j$. All maps may be assigned fixed values outside
$\Omega_0$.
\end{proof}

\begin{proposition}\label{B:prop:bessel-spectral-projections}
Fix $j\ge1$. If $\|R_n-R\|_{\rm op}<g_j/4$, then $n\ge j$ and
\begin{equation}
 \|\Pi_{j,n}^{\mathfrak m}-\Pi_j\|_{\mathfrak S_2}
                                  \le \frac{2}{g_j}e_n.
 \label{B:eq:bessel-projection-bound}
\end{equation}
\end{proposition}
\begin{proof}
Let $\Gamma_j$ be the counterclockwise circle of radius $g_j$ centered at $\mu_j$. For $\zeta\in\Gamma_j$, the distance to $\sigma(R)$ is at least $g_j$. Thus
$\|(\zeta I-R)^{-1}\|_{\rm op}\le g_j^{-1}$. Factor
\[
 \zeta I-R_n=\bigl[I-(R_n-R)(\zeta I-R)^{-1}\bigr](\zeta I-R).
\]
The norm of the perturbation in brackets is less than $1/4$; its geometric inverse series gives $\|(\zeta I-R_n)^{-1}\|_{\rm op}\le2g_j^{-1}$. The second estimate in
\eqref{B:eq:spectral-matching} places $\mu_{j,n}$ inside $\Gamma_j$. Since $\mu_j\ge4g_j$, $\mu_{j,n}>0$ and $n\ge j$. For $l\ne j$,
\[
 |\mu_{l,n}-\mu_j|
 \ge |\mu_l-\mu_j|-|\mu_{l,n}-\mu_l|
 >4g_j-g_j/4>g_j.
\]
Zero is also outside the circle. On $\ker R_n$, including the orthogonal complement of the finite embedding, the resolvent is $\zeta^{-1}I$. Its integral is zero. On an eigenvector
with reciprocal eigenvalue $\mu_{l,n}$, the scalar integral is
\[
 \frac1{2\pi\ii}\oint_{\Gamma_j}(\zeta-\mu_{l,n})^{-1}\,\dd\zeta
 =\mathbf1_{\{l=j\}}.
\]
For an eigenvector $\phi_l$ of $R$ with eigenvalue $\mu_l$,
\[
 \frac1{2\pi\ii}\oint_{\Gamma_j}(\zeta I-R)^{-1}\phi_l\,\dd\zeta
 =\left(\frac1{2\pi\ii}\oint_{\Gamma_j}
                  \frac{\dd\zeta}{\zeta-\mu_l}\right)\phi_l
 =\mathbf1_{\{l=j\}}\phi_l.
\]
The contour integral defines a bounded operator of norm at most $(2\pi g_j)g_j^{-1}/(2\pi)=1$. Since the eigenvectors form a complete orthonormal basis, equality on their finite
spans and continuity identify this operator with $\Pi_j$. The preceding finite-dimensional calculation identifies the other integral with $\Pi_{j,n}^{\mathfrak m}$. Their
difference, by the resolvent identity, is
\begin{equation}
 \Pi_{j,n}^{\mathfrak m}-\Pi_j
 =\frac1{2\pi\ii}\oint_{\Gamma_j}
   (\zeta I-R_n)^{-1}(R_n-R)(\zeta I-R)^{-1}\,\dd\zeta.
 \label{B:eq:bessel-projection-resolvent-identity}
\end{equation}
The Hilbert--Schmidt ideal inequality makes the integrand continuous in that norm. Parametrizing the circle as above gives
\[
 \begin{split}
 \|\Pi_{j,n}^{\mathfrak m}-\Pi_j\|_{\mathfrak S_2}
 &\le\frac1{2\pi}\int_0^{2\pi}
       g_j\|(\zeta_j(t)I-R_n)^{-1}\|_{\rm op}\,e_n\,
          \|(\zeta_j(t)I-R)^{-1}\|_{\rm op}\,\dd t\\
 &\le\frac{2e_n}{g_j}.
 \end{split}
\]
\end{proof}

\begin{corollary}\label{B:cor:bessel-eigenfunction-observables}
Under the hypothesis of Proposition~\ref{B:prop:bessel-spectral-projections},
\begin{equation}
 \min_{\sigma\in\{-1,1\}}
       \|\psi_{j,n}^{\mathfrak m}-\sigma\psi_j\|_2
                       \le2g_j^{-1}e_n.
 \label{B:eq:bessel-eigenfunction-bound}
\end{equation}
For every bounded operator $F$ on $L^2(0,1)$,
\begin{equation}
 |\langle\psi_{j,n}^{\mathfrak m},F\psi_{j,n}^{\mathfrak m}\rangle
                   -\langle\psi_j,F\psi_j\rangle|
                       \le4g_j^{-1}\|F\|_{\rm op}e_n.
 \label{B:eq:bessel-bounded-observable}
\end{equation}
For complex unit eigenvectors, replace the signs by complex numbers of modulus one.
\end{corollary}
\begin{proof}
Let $u,v$ be unit vectors, and change the phase of $v$ so that $c=\langle u,v\rangle$ is real and nonnegative. If $c<1$, put $s=(1-c^2)^{1/2}$ and $w=(v-cu)/s$. In the orthonormal
basis $(u,w)$, the difference of their rank-one projections is
\[
 P_v-P_u=\begin{pmatrix}-s^2&cs\\cs&s^2\end{pmatrix}.
\]
Its square is $s^2I$, and it is zero on the orthogonal complement. The same conclusions hold at $c=1$, where the difference is zero. It follows that
\begin{equation}
 \|P_v-P_u\|_{\mathfrak S_2}^2=2(1-c^2),\qquad
 \|P_v-P_u\|_{\mathfrak S_1}=2\sqrt{1-c^2},\qquad
 \min_{|\zeta|=1}\|v-\zeta u\|_2^2=2(1-c).
 \label{B:eq:rank-one-phase-identities}
\end{equation}
The phase formula also follows directly from $\|v-\zeta u\|_2^2=2-2\operatorname{Re}\langle v,\zeta u\rangle$; the maximum of the last real part over $|\zeta|=1$ is
$|\langle v,u\rangle|=c$. At $c=0$ every phase gives the same distance, and we choose the phase $1$. Since $1-c\le1-c^2$, the eigenfunction bound follows from
\eqref{B:eq:bessel-projection-bound}. For real vectors the optimizing phase is a sign. For the phase making $\|v-u\|_2$ minimal,
\[
 \begin{split}
 \langle v,Fv\rangle-\langle u,Fu\rangle
 &=\langle v-u,Fv\rangle+\langle u,F(v-u)\rangle,\\
 |\langle v,Fv\rangle-\langle u,Fu\rangle|
 &\le2\|F\|_{\rm op}\|v-u\|_2.
 \end{split}
\]
For $|\zeta|=1$, sesquilinearity gives $\langle\zeta u,F\zeta u\rangle=\langle u,Fu\rangle$. Substitution of \eqref{B:eq:bessel-eigenfunction-bound} proves
\eqref{B:eq:bessel-bounded-observable}. For multiplication by $f\in L^\infty(0,1)$ one has $\|F\|_{\rm op}=\|f\|_\infty$.
\end{proof}

The estimates concern $L^2$ eigenfunctions and bounded operators on that space. Evaluation at a point and a boundary derivative are not bounded functionals on $L^2(0,1)$; the
bounds above do not imply their convergence. For multiplication by $\mathbf1_A$ they do give, for every Borel set $A\subset(0,1)$,
\[
 \left|\int_A|\psi_{j,n}^{\mathfrak m}|^2-\int_A|\psi_j|^2\right|
 \le4g_j^{-1}e_n
\]
whenever the stated hypothesis on the spectral gap holds.

\subsubsection{Powers of the inverse and their traces}
\label{B:subsection:negative-powers-and-traces}

Powers of the inverse preserve the embedding because its adjoint composed with the embedding is the finite identity. A telescoping product then transfers the inverse error, and
products of Hilbert--Schmidt factors provide the trace-norm estimates for higher powers.

For an integer $m\ge1$, the isometry of $\mathscr J_n^{\mathfrak m}$ gives
\begin{equation}
 R_n^m=\mathscr J_n^{\mathfrak m}(\mathcal A_n^{\mathfrak m})^{-m}
                  (\mathscr J_n^{\mathfrak m})^*,\qquad
 R^m=\mathcal H_{\beta,a}^{-m}.
 \label{B:eq:spectral-embedded-powers}
\end{equation}
For $m=1$ this is the definition. Multiplication by $R_n$ uses $(\mathscr J_n^{\mathfrak m})^*\mathscr J_n^{\mathfrak m}=I_n$ and gives the identity at $m+1$. For the limiting
operator, $R^m\psi_j=\Lambda_j^{-m}\psi_j$; the eigenbasis therefore identifies $R^m$ with the bounded operator $\mathcal H_{\beta,a}^{-m}$.

\begin{proposition}\label{B:prop:bessel-inverse-powers}
For every integer $m\ge1$ and every $n$,
\begin{equation}
 \|R_n^m-R^m\|_{\mathfrak S_2}\le mM^{m-1}e_n.
 \label{B:eq:spectral-powers-HS}
\end{equation}
For $m\ge2$ one also has
\begin{equation}
 \|R_n^m-R^m\|_{\mathfrak S_1}\le mM^{m-2}H_2e_n.
 \label{B:eq:spectral-powers-trace-norm}
\end{equation}
\end{proposition}
\begin{proof}
The telescoping identity, with the factors in their displayed order, is
\begin{equation}
 R_n^m-R^m=\sum_{l=0}^{m-1}R_n^l(R_n-R)R^{m-1-l}.
 \label{B:eq:noncommuting-power-difference}
\end{equation}
The ideal inequality in Theorem~2.7(a) of \cite{simon2005trace} bounds each summand by $M^{m-1}e_n$, proving the first norm bound. For the trace norm use the Schatten product
inequality $\|AB\|_{\mathfrak S_1}\le\|A\|_{\mathfrak S_2}\|B\|_{\mathfrak S_2}$; see Theorem~2.8 of \cite{simon2005trace}, with exponents $1,2,2$. When $m\ge2$, at least one of
the powers next to $R_n-R$ supplies a second Hilbert--Schmidt factor. Specifically,
\[
 \begin{aligned}
 \|R_n^l(R_n-R)R^{m-1-l}\|_{\mathfrak S_1}
 &\le M^{l-1}H_2e_nM^{m-1-l}\quad(1\le l\le m-1),\\
 \|(R_n-R)R^{m-1}\|_{\mathfrak S_1}
 &\le e_nH_2M^{m-2}.
 \end{aligned}
\]
Summing proves \eqref{B:eq:spectral-powers-trace-norm}.
\end{proof}

\begin{corollary}\label{B:cor:bessel-negative-power-traces}
For every fixed integer $m\ge2$,
\begin{equation}
 \left|\sum_{j=1}^n(\lambda_{j,n}^{\mathfrak m})^{-m}
                 -\sum_{j\ge1}\Lambda_j^{-m}\right|
                    \le mM^{m-2}H_2e_n.
 \label{B:eq:spectral-powers-trace}
\end{equation}
\end{corollary}
\begin{proof}
For $T=R_n$ or $R$, let $(w_j)$ be a spectral basis, including zero eigenvectors where necessary. For $m\ge2$,
\[
 \sum_j\langle T^mw_j,w_j\rangle
 =\sum_j\mu_j(T)^m
 \le M^{m-2}\sum_j\mu_j(T)^2
 =M^{m-2}\|T\|_{\mathfrak S_2}^2<\infty.
\]
The positive operator $T^m$ is therefore trace class with this trace. By linearity of the trace and its bound by the trace norm,
\[
 \left|\operatorname{Tr}(R_n^m)-\operatorname{Tr}(R^m)\right|
   =\left|\operatorname{Tr}(R_n^m-R^m)\right|
   \le\|R_n^m-R^m\|_{\mathfrak S_1}.
\]
By \eqref{B:eq:spectral-embedded-powers}, the nonzero eigenvalues of $R_n^m$ are $(\lambda_{j,n}^{\mathfrak m})^{-m}$, $1\le j\le n$; the remaining eigenvalues are zero. The
corresponding eigenvalues of $R^m$ are $\Lambda_j^{-m}$. Substituting these two trace identities into the preceding inequality and using \eqref{B:eq:spectral-powers-trace-norm}
proves the stated bound.
\end{proof}

The corresponding implication for the first power fails even for positive diagonal operators with uniformly bounded traces. On $\ell^2(\mathbb N)$ let
$T=\operatorname{diag}(j^{-2})$, let $P_n$ project onto the first $n$ coordinates, and set $T_n=P_nTP_n+n^{-1}P_n$. Direct calculation gives
\begin{equation}
 \begin{split}
 \|T_n-T\|_{\mathfrak S_2}^2
       &=n^{-1}+\sum_{j>n}j^{-4},\\
 \|T_n-T\|_{\mathfrak S_1}
       &=1+\sum_{j>n}j^{-2},\qquad
 \operatorname{Tr}(T_n-T)=1-\sum_{j>n}j^{-2}.
 \end{split}
 \label{B:eq:spectral-trace-counterexample}
\end{equation}
Integral comparison gives $\sum_{j>n}j^{-4}\le(3n^3)^{-1}$ and $\sum_{j>n}j^{-2}\le n^{-1}$. Hence the squared Hilbert--Schmidt error is $n^{-1}+O(n^{-3})$, while the last two
quantities tend to one. Convergence of the traces of the first powers under the coupling of the bidiagonal matrices to the limiting inverse therefore requires an estimate beyond
\eqref{B:eq:spectral-setup}.

\subsubsection{Regularized determinants}
\label{B:subsection:regularized-determinants}

The Hilbert--Schmidt estimate controls the regularized determinant. Its zeros recover the direct spectrum, while its local uniform comparison follows from the squared summability
of the reciprocal eigenvalues.

For the present positive operators, the regularized Fredholm determinant is the following product over the eigenvalues:
\begin{equation}
 \begin{split}
 D_n^{\mathfrak m}(z)&=\det{}_2(I-zR_n)
              =\prod_{j=1}^n(1-z\mu_{j,n})\me^{z\mu_{j,n}},\\
 D(z)&=\det{}_2(I-zR)
              =\prod_{j\ge1}(1-z\mu_j)\me^{z\mu_j}.
 \end{split}
 \label{B:eq:spectral-determinants}
\end{equation}
This is the convention of Theorem~6.2 of \cite{simon1977determinants}.

\begin{lemma}\label{B:lem:bessel-regularized-products}
The functions in \eqref{B:eq:spectral-determinants} are entire. Their zeros are precisely $\lambda_{j,n}^{\mathfrak m}$ and $\Lambda_j$, respectively, and every zero is simple.
\end{lemma}
\begin{proof}
The finite product is entire. For the infinite product, fix $A>0$ and choose $J=J(A)$ such that $A\mu_j\le1/2$ for every $j>J$. On $|w|<1$ use the analytic logarithm normalized at
zero. Then
\[
 \log((1-w)e^w)=-\sum_{r=2}^{\infty}\frac{w^r}{r},\qquad
 |\log((1-w)e^w)|\le\frac{|w|^2}{2(1-|w|)}\le |w|^2
       \quad(|w|\le1/2).
\]
Consequently
\[
 \sum_{j>J}\sup_{|z|\le A}
       |\log((1-z\mu_j)e^{z\mu_j})|
 \le A^2\sum_{j>J}\mu_j^2<\infty.
\]
The logarithm series converges uniformly on the disk. Its exponential is holomorphic and has no zeros there. Multiplication by the first $J$ factors proves local uniform
convergence of the complete product. As $A$ was arbitrary, $D$ is entire, with the regularized determinant convention of \cite[Theorem~6.2]{simon1977determinants}.

On any disk, all zeros therefore come from its finitely many separated factors. The factor $(1-z\mu_j)e^{z\mu_j}$ vanishes at $z=\mu_j^{-1}$ and has derivative
$-\mu_j\mathrm e\ne0$ at that point. The other factors and the infinite exponential factor are nonzero there, since the positive eigenvalues are distinct by
Corollary~\ref{B:cor:bessel-simple-spectrum}. This proves the stated zeros and their multiplicities for $D$. Applying the same factor calculation to the finite product, and using
Lemma~\ref{B:lem:bessel-spectral-measurability}, proves the corresponding statement for $D_n^{\mathfrak m}$.
\end{proof}

For $A>0$, set
\begin{equation}
 C_A=2A^2H_2\exp(AM+A^2H_2^2).
 \label{B:eq:determinant-comparison-constant}
\end{equation}

\begin{proposition}\label{B:prop:bessel-determinant-comparison}
For every $A>0$ and every $n$,
\begin{equation}
 \sup_{|z|\le A}|D_n^{\mathfrak m}(z)-D(z)|\le C_Ae_n.
 \label{B:eq:spectral-determinant-bound}
\end{equation}
\end{proposition}
\begin{proof}
For $w\in\mathbb C$, put $x=-2\operatorname{Re}w+|w|^2$. Then $1+x=|1-w|^2\ge0$ and $1+x\le\me^x$, whence
\begin{equation}
 |(1-w)\me^w|^2
     =|1-w|^2\me^{2\operatorname{Re}w}\le\me^{|w|^2}.
 \label{B:eq:regularized-factor-modulus}
\end{equation}
For $s,t\in[0,M]$, integration of $\partial_tg_z(t)=-z^2t\me^{zt}$ gives
\begin{equation}
 |g_z(t)-g_z(s)|
 \le A^2\me^{AM}\left|\int_s^t u\,\dd u\right|
 \le A^2\me^{AM}|t-s|(t+s),\qquad |z|\le A.
 \label{B:eq:regularized-factor-difference}
\end{equation}
For $N\ge n$, continue $\mu_{j,n}$ by zero. The finite product identity is
\[
 \prod_{j=1}^Ng_z(\mu_{j,n})-\prod_{j=1}^Ng_z(\mu_j) 
 =\sum_{k=1}^N
      \left(\prod_{j<k}g_z(\mu_{j,n})\right)
      [g_z(\mu_{k,n})-g_z(\mu_k)]
      \left(\prod_{j>k}g_z(\mu_j)\right).
\]
The identity is obtained by replacing the factors in order. It uses no division by a factor and therefore holds also when a finite or limiting determinant vanishes at $z$. Since
$g_z(0)=1$, the first product for $N\ge n$ is exactly $D_n^{\mathfrak m}(z)$. Every product outside the bracket has modulus at most
\[
 \exp\!\left(\frac{A^2}{2}
    \left[\sum_j\mu_{j,n}^2+\sum_j\mu_j^2\right]\right)
                      \le\exp(A^2H_2^2)
\]
by \eqref{B:eq:regularized-factor-modulus}. Furthermore,
\begin{align}
 \sum_{j\ge1}|\mu_{j,n}-\mu_j|(\mu_{j,n}+\mu_j)
 &\le\left(\sum_{j\ge1}|\mu_{j,n}-\mu_j|^2\right)^{1/2}
       \left(\sum_{j\ge1}(\mu_{j,n}+\mu_j)^2\right)^{1/2}\notag\\
 &\le e_n(\|R_n\|_{\mathfrak S_2}+\|R\|_{\mathfrak S_2})
 \le2H_2e_n.
 \label{B:eq:determinant-spectral-sum}
\end{align}
Use \eqref{B:eq:regularized-factor-difference} in each bracket and then \eqref{B:eq:determinant-spectral-sum}. The resulting bound is $C_Ae_n$, uniformly in $N$ and $|z|\le A$.
Lemma~\ref{B:lem:bessel-regularized-products} permits passage to $N\to\infty$, proving the estimate on the entire disk.
\end{proof}

\subsubsection{Derivatives and logarithms of the determinants}
\label{B:subsection:regularized-determinant-derivatives}

Local uniform control of the entire determinants also controls their derivatives by Cauchy's formula. Statements about logarithms require a further restriction to sets without
zeros.

\begin{corollary}\label{B:cor:bessel-determinant-derivatives}
Let $K\subset\mathbb C$ be nonempty and compact, $h>0$, and $r\ge0$ an integer. With $a_K=\max_{z\in K}|z|$,
\begin{equation}
 \sup_{z\in K}|(D_n^{\mathfrak m})^{(r)}(z)-D^{(r)}(z)|
                   \le r!h^{-r}C_{a_K+h}e_n.
 \label{B:eq:spectral-determinant-derivatives}
\end{equation}
\end{corollary}
\begin{proof}
For $r\ge1$ apply Cauchy's formula to $F_n=D_n^{\mathfrak m}-D$ on the circle $|\zeta-z|=h$:
\[
 F_n^{(r)}(z)=\frac{r!}{2\pi\ii}
            \oint_{|\zeta-z|=h}\frac{F_n(\zeta)}{(\zeta-z)^{r+1}}
                                      \,\dd\zeta.
\]
The circle lies in $\{|\zeta|\le a_K+h\}$. Proposition~\ref{B:prop:bessel-determinant-comparison} and its length $2\pi h$ give
\[
 |F_n^{(r)}(z)|\le\frac{r!}{2\pi}\,(2\pi h)h^{-r-1}C_{a_K+h}e_n
       =r!h^{-r}C_{a_K+h}e_n.
\] The bound is independent of $z\in K$. The function $F_n$ is entire, so the
circle may pass through a zero of either determinant. For $r=0$, apply the same proposition on the disk of radius
$a_K+h$.
\end{proof}

For $|z|<M^{-1}$ every $1-z\mu_j$ and $1-z\mu_{j,n}$ lies in the disk of radius one about one. Let $L$ and $L_n$ denote the analytic logarithms of $D$ and $D_n^{\mathfrak m}$ on
this disk in the $z$-plane, normalized to vanish at zero.

\begin{proposition}\label{B:prop:bessel-determinant-logarithms}
For $|z|<M^{-1}$,
\begin{equation}
 L_n(z)=-\sum_{m=2}^{\infty}\frac{z^m}{m}\operatorname{Tr}(R_n^m),
 \qquad
 L(z)=-\sum_{m=2}^{\infty}\frac{z^m}{m}\operatorname{Tr}(R^m).
 \label{B:eq:regularized-determinant-trace-series}
\end{equation}
For $0<A<M^{-1}$,
\begin{equation}
 \sup_{|z|\le A}|L_n(z)-L(z)|
                 \le\frac{A^2H_2}{1-AM}e_n.
 \label{B:eq:regularized-logarithm-comparison}
\end{equation}
\end{proposition}
\begin{proof}
For $|w|<1$, integrate the geometric series along the segment from zero to $w$:
\[
 \log(1-w)+w=-\int_0^w\frac{u}{1-u}\,\dd u
                         =-\sum_{m=2}^{\infty}\frac{w^m}{m}.
\]
On $|z|\le A<M^{-1}$, absolute summability follows from
\[
 \begin{split}
 \sum_j\sum_{m=2}^{\infty}\frac{A^m\mu_j^m}{m}
 &\le\sum_j\mu_j^2
               \sum_{m=2}^{\infty}\frac{A^mM^{m-2}}m\\
 &\le\frac{A^2H_2^2}{1-AM}.
 \end{split}
\]
The same bound holds for $(\mu_{j,n})$, uniformly in $n$. It permits interchanging the sums and proves locally uniform convergence of the resulting holomorphic series.
Exponentiating the partial sums gives the partial products for $D$ and $D_n^{\mathfrak m}$. Their limits are therefore analytic logarithms of these products, with value zero at
zero. This proves \eqref{B:eq:regularized-determinant-trace-series} and the existence of the chosen branches. Uniqueness follows because the difference of two such logarithms is a
continuous function with values in $2\pi\ii\mathbb Z$, and vanishes at zero.

Corollary~\ref{B:cor:bessel-negative-power-traces} now gives
\[
 \begin{split}
 \sup_{|z|\le A}|L_n(z)-L(z)|
 &\le\sum_{m=2}^{\infty}\frac{A^m}{m}
             mM^{m-2}H_2e_n\\
 &=\frac{A^2H_2}{1-AM}e_n.
 \end{split}
\]
For an integer $q\ge1$ and $m\ge\max\{2,q\}$,
\[
 \frac{\dd^q}{\dd z^q}\frac{z^m}{m}
   =\frac{(m-1)!}{(m-q)!}z^{m-q},\qquad
 |\operatorname{Tr}(R^m)|\le H_2^2M^{m-2}.
\]
The same trace bound holds for $R_n$, uniformly in $n$. For $0<A<A'<M^{-1}$ it follows that
\[
 \begin{split}
 &\sum_{m\ge\max\{2,q\}}
 \sup_{|z|\le A'}\left|
   \frac{\dd^q}{\dd z^q}\left(\frac{z^m}{m}\operatorname{Tr}(R^m)\right)\right|\\
 &\quad\le H_2^2M^{q-2}
       \sum_{m=q}^{\infty}\frac{(m-1)!}{(m-q)!}(A'M)^{m-q}
 =\frac{H_2^2M^{q-2}(q-1)!}{(1-A'M)^q}.
 \end{split}
\]
The last equality is the binomial series, since $(m-1)!/(m-q)!=(q-1)!\binom{m-1}{q-1}$. Thus each differentiated series converges uniformly on $|z|\le A'$. Cauchy's formula on the
circle $|z|=A'$ identifies its limit with the corresponding derivative on $|z|\le A$. We also obtain the first variation along $R_t=R+tD_n$, where $D_n=R_n-R$. Let $L_t(z)$ be the
same series with $R_t$ in place of $R$. Convexity gives $\|R_t\|_{\rm op}\le M$ for $0\le t\le1$. For $m\ge2$, the product rule and cyclicity of the trace give
\[
 \frac{\dd}{\dd t}\operatorname{Tr}(R_t^m)
      =m\operatorname{Tr}(D_nR_t^{m-1}),\qquad
 \frac{\dd^2}{\dd t^2}\operatorname{Tr}(R_t^m)
      =m\sum_{l=0}^{m-2}\operatorname{Tr}
                      (D_nR_t^lD_nR_t^{m-2-l}).
\]
These products are trace class: the first has the two Hilbert--Schmidt factors $D_n,R_t$, and each term of the second has two copies of $D_n$. Convexity also gives
$\|R_t\|_{\mathfrak S_2}\le H_2$. Thus, uniformly in $t\in[0,1]$,
\[
 \left|\frac{\dd}{\dd t}\operatorname{Tr}(R_t^m)\right|
     \le mM^{m-2}H_2e_n,\qquad
 \left|\frac{\dd^2}{\dd t^2}\operatorname{Tr}(R_t^m)\right|
     \le m(m-1)M^{m-2}e_n^2.
\]
After multiplication by $A^m/m$, these majorants have sums $A^2H_2e_n/(1-AM)$ and $A^2e_n^2/(1-AM)^2$, respectively. They justify both differentiations of the series, uniformly on
$[0,1]\times\{|z|\le A\}$. Summing it and the differentiated inverse identity yields
\[
 \partial_tL_t(z)=-z^2\operatorname{Tr}((I-zR_t)^{-1}R_tD_n),\qquad
 \partial_t^2L_t(z)
 =-z^2\operatorname{Tr}((I-zR_t)^{-1}D_n(I-zR_t)^{-1}D_n).
\]
The equality
\[
 L_1(z)-L_0(z)-\partial_tL_t(z)|_{t=0}
       =\int_0^1(1-t)\partial_t^2L_t(z)\,\dd t
\]
follows by two integrations in $t$. Since $\|(I-zR_t)^{-1}\|\le(1-AM)^{-1}$, the second derivative has modulus at most $A^2(1-AM)^{-2}e_n^2$. Integration of $1-t$ gives
\[
 \sup_{|z|\le A}
 |L_n(z)-L(z)+z^2\operatorname{Tr}((I-zR)^{-1}RD_n)|
 \le\frac{A^2}{2(1-AM)^2}e_n^2.
\]
The constant depends on the admissible disk. The $O_\omega(n^{-1})$ remainder concerns the regularized logarithm; it does not provide an error bound for the trace of the first
power or the ordinary determinant. The radius $M^{-1}$ is fixed after the path is chosen. The logarithm estimates use $A<M^{-1}$, whereas the preceding determinant bounds apply on
any deterministic compact set. In particular, the linear coefficient vanishes and the quadratic coefficient of $L$ is $-\operatorname{Tr}(R^2)/2$.
\end{proof}

For a finite matrix, the ordinary determinant satisfies
\[
 \det(I-zR_n)=\me^{-z\operatorname{Tr}R_n}D_n^{\mathfrak m}(z).
\]
Since $R$ is trace class, its ordinary determinant has the same factorization into the exponential of its trace and its regularized determinant. A rate for the ordinary
determinants therefore also requires an estimate for $\operatorname{Tr}R_n-\operatorname{Tr}R$. The counterexample \eqref{B:eq:spectral-trace-counterexample} shows why this
additional estimate is needed.

\subsubsection{Measurability and rates in probability}
\label{B:subsection:bessel-spectral-probability-rates}

The pathwise estimates contain random constants, gaps and starting indices. Their conversion to rates in probability uses measurability and almost sure finiteness; no moment
assumption on these constants is needed.

\begin{lemma}\label{B:lem:random-spectral-rate}
Let $(A_n)$ be nonnegative finite measurable random variables. If $A_n\le Ce_n$ for $n\ge N$, where $C\ge0$ is measurable and finite almost surely and $N$ is a finite measurable
random positive integer, then
\begin{equation}
 \sup_{n\ge1}\sqrt n A_n<\infty\quad\text{almost surely},\qquad
 A_n=O_{\mathbb P}(n^{-1/2}).
 \label{B:eq:spectral-tightness}
\end{equation}
\end{lemma}
\begin{proof}
On $\Omega_0$ set $B=C\sqrt Z+\max_{1\le n<N}\sqrt n A_n$, with an empty maximum equal to zero. A maximum with this random finite range is measurable, since it is the countable
supremum of $\mathbf1_{\{n<N\}}\sqrt n A_n$. It is finite almost surely. For every $n$, $\sqrt n A_n\le B$, and hence
\[
 \sup_n\mathbb P\{\sqrt n A_n>L\}\le\mathbb P\{B>L\}\longrightarrow0
                            \qquad(L\to\infty).
\]
For a finite random variable $B$, the events $\{B>L\}$ decrease to the empty set up to a null set. Continuity from above of the probability measure proves the displayed limit. Thus
for every $\varepsilon>0$ a deterministic $L$ makes all the probabilities smaller than $\varepsilon$ simultaneously. This proves the $O_{\mathbb P}(n^{-1/2})$ bound from the almost
sure finiteness of $B$.
\end{proof}

We verify measurability of the errors to which the lemma is applied. Finite eigenvectors are continuous by Lemma~\ref{B:lem:bessel-spectral-measurability}, and limiting projections
and eigenfunctions are measurable by Lemma~\ref{B:lem:bessel-limiting-spectral-measurability}. The norm estimates in Proposition~\ref{B:prop:bessel-inverse-powers}, applied to two
arbitrary positive Hilbert--Schmidt operators, show that $T\mapsto T^m$ is locally continuous from $\mathfrak S_2$ to $\mathfrak S_2$, and to $\mathfrak S_1$ when $m\ge2$. Thus
both operator errors in that proposition are measurable in their stated norms, as are their traces.

For fixed $z$, $D(z)$ is the limit of measurable finite products. The Cauchy integral in Corollary~\ref{B:cor:bessel-determinant-derivatives} expresses each fixed derivative as a
limit of measurable Riemann sums. For a fixed compact $K$, take a countable dense subset $K_0$. Continuity gives, for every fixed nonnegative integer $r$,
\[
 \sup_{z\in K}|(D_n^{\mathfrak m})^{(r)}(z)-D^{(r)}(z)|
 =\sup_{z\in K_0}|(D_n^{\mathfrak m})^{(r)}(z)-D^{(r)}(z)|.
\]
Hence these suprema are measurable. Linear statistics are limits of their finite partial sums and are measurable as well.

For a positive measurable threshold $\delta$ and measurable $a_n\to0$ on $\Omega_0$, define
\[
 N=\min\{m\ge1:a_n<\delta\text{ for every }n\ge m\},
 \qquad \{N\le m\}=\bigcap_{n\ge m}\{a_n<\delta\}.
\]
Take $a_n=e_n$, $\delta=\mu_m/2$ for a fixed eigenvalue count $m$, or $a_n=\|R_n-R\|_{\rm op}$, $\delta=g_j/4$ for a fixed projection. The resulting $N$ is measurable and finite on
$\Omega_0$.

For the eigenvalue error take $A_n=(\sum_{j=1}^m|\lambda_{j,n}^{\mathfrak m}-\Lambda_j|^2)^{1/2}$ when $n\ge m$ and set $A_n=0$ otherwise. The same zero convention applies to a
fixed eigenvector or projection absent at the first few sizes. The displayed deterministic estimates and Lemma~\ref{B:lem:random-spectral-rate} give
\begin{equation}
 \begin{gathered}
 \sup_{n\ge m}n\sum_{j=1}^m
             |\lambda_{j,n}^{\mathfrak m}-\Lambda_j|^2<\infty,
 \qquad
 \sup_{n\ge j}n\|\Pi_{j,n}^{\mathfrak m}-\Pi_j\|_{\mathfrak S_2}^2<\infty,\\
 \sup_{n\ge j}n\min_{\sigma\in\{-1,1\}}
              \|\psi_{j,n}^{\mathfrak m}-\sigma\psi_j\|_2^2<\infty
                    \quad\text{almost surely}.
 \end{gathered}
 \label{B:eq:bessel-spectral-squared-rates}
\end{equation}
For a fixed linear statistic, bounded quadratic observable, inverse power, trace of a negative power, or determinant derivative on a deterministic compact set, let $A_n$ be its
absolute or norm error in the preceding estimates. In each case,
\begin{equation}
 \sup_n\sqrt n\,A_n<\infty\quad\text{almost surely},\qquad
 \lim_{L\to\infty}\sup_n\mathbb P(\sqrt n A_n>L)=0.
 \label{B:eq:bessel-fixed-observable-rates}
\end{equation}
For the projection error, $C=2/g_j$ applies for every $n\ge N$ with the measurable index $N$ constructed above. For the determinant derivative, $C=r!h^{-r}C_{a_K+h}$ applies for
every $n$. These constants are measurable and finite. The squared reciprocal error obeys the global bound $n\sum_{j\ge1}|\mu_{j,n}-\mu_j|^2\le Z$. For any fixed finite family
$A_n^{(1)},\ldots,A_n^{(s)}$,
\[
 \mathbb P\left(\max_{1\le l\le s}\sqrt n A_n^{(l)}>L\right)
 \le\sum_{l=1}^s\mathbb P(\sqrt n A_n^{(l)}>L).
\]
This proves joint tightness for that family.

\begin{proof}[Proof of Corollary~\ref{B:cor:intro-eigenvalue-rates}]
Lemma~\ref{B:lem:ordered-reciprocal-eigenvalues} proves the first inequality, including the reciprocal tail. For $j\le m$ and $e_n<(2\Lambda_m)^{-1}$, that lemma gives
$(\lambda_{j,n}^{\mathfrak m})^{-1}\ge(2\Lambda_j)^{-1}$; a reciprocal continued by zero beyond index $n$ would contradict this bound. Thus $n\ge m$ and
\[
 |\lambda_{j,n}^{\mathfrak m}-\Lambda_j|
 =\frac{|(\lambda_{j,n}^{\mathfrak m})^{-1}-\Lambda_j^{-1}|}
         {(\lambda_{j,n}^{\mathfrak m})^{-1}\Lambda_j^{-1}}
 \le2\Lambda_j^2e_n.
\]
On $e_n<\mu_m/2$, the maximum and the Euclidean norm of the first $m$ errors are at most $2\Lambda_m^2e_n$, by \eqref{B:eq:spectral-fixed-eigenvalues}. A measurable starting index
is
\[
 N_m=\min\{N:e_n<\mu_m/2\text{ for all }n\ge N\},\qquad
 \{N_m\le N\}=\bigcap_{n\ge N}\{e_n<\mu_m/2\}.
\]
For every $n\ge N_m$, the error multiplied by $\sqrt n$ is at most $2\Lambda_m^2\sqrt Z$. The finitely many earlier errors, for $n\ge m$, are finite; including them in the maximum
gives a finite measurable bound and hence tightness by its tail.
\end{proof}

\begin{proof}[Proof of Corollary~\ref{B:cor:intro-projection-rates}]
Let $g_j>0$ be as in \eqref{B:eq:spectral-phase-gap}, and put
\[
 N_j=\min\{N:\|R_n-R\|_{\rm op}<g_j/4\text{ for every }n\ge N\}.
\]
It is finite on $\Omega_0$ and measurable, since $\{N_j\le m\}=\bigcap_{n\ge m}\{\|R_n-R\|_{\rm op}<g_j/4\}$. For $n\ge N_j$, Proposition~\ref{B:prop:bessel-spectral-projections}
and \eqref{B:eq:rank-one-phase-identities} give
\[
 \begin{aligned}
 \sqrt n\min_{\sigma\in\{-1,1\}}
       \|\psi_{j,n}^{\mathfrak m}-\sigma\psi_j\|_2
 &\le\sqrt n\|\Pi_{j,n}^{\mathfrak m}-\Pi_j\|_{\mathfrak S_2}\\
 &\le\frac{2\sqrt n\,e_n}{g_j}\le\frac{2\sqrt Z}{g_j}.
 \end{aligned}
\]
The errors are measurable by Lemmas~\ref{B:lem:bessel-spectral-measurability} and \ref{B:lem:bessel-limiting-spectral-measurability}. Assign zero at $n<j$ and apply
Lemma~\ref{B:lem:random-spectral-rate} to each error with $C=2/g_j$ and $N=N_j$. This proves both rates for the projections and the eigenfunctions.
\end{proof}

\begin{proof}[Proof of Corollary~\ref{B:cor:intro-determinant-rates}]
For $A>\sup_{z\in K}|z|$, use $M,H_2$ from \eqref{B:eq:spectral-setup}. Proposition~\ref{B:prop:bessel-determinant-comparison} gives the finite bound
\[
 \sup_{|z|\le A}|D_n^{\mathfrak m}(z)-D(z)|
 \le2A^2H_2e^{AM+A^2H_2^2}e_n.
\]
Both $M$ and $H_2$ are bounded by fixed functions of $\|R\|_{\mathfrak S_2}$ and $\sqrt Z$. Cauchy's formula on circles of any fixed radius $h>0$ gives an additional factor
$r!h^{-r}$ for derivative order $r$, using the disk of radius $A+h$. With $C=r!h^{-r}C_{A+h}$ the resulting error satisfies $\sqrt n A_n\le C\sqrt Z$ for every $n$. Thus
\[
 \sup_n\mathbb P\{\sqrt n A_n>L\}
 \le\mathbb P\{C\sqrt Z>L\}\longrightarrow0.
\]
This also includes $r=0$. The compact set is deterministic and may contain zeros of either determinant, since both determinant functions are entire.
\end{proof}

\subsubsection{Spectral truncation}
\label{B:subsection:bessel-best-finite-rank-approximation}

The lower bound of Theorem~\ref{B:thm:brownian-information-lower-bound} cannot hold for arbitrary approximations of finite rank. When the full limiting inverse is available,
spectral truncation gives a smaller error, as the next proposition and corollary show.

Let $R^{[n]}=\sum_{j=1}^n\mu_j\Pi_j$ be the spectral truncation of $R$. The following proposition is the compact positive operator version of the best approximation theorem of
\cite{eckart1936approximation} (pp.~214--216).

\begin{proposition}\label{B:prop:bessel-best-finite-rank-approximation}
For every $n$,
\begin{equation}
 \inf_{\operatorname{rank}T\le n}\|R-T\|_{\mathfrak S_2}^2
       =\|R-R^{[n]}\|_{\mathfrak S_2}^2
       =\sum_{j>n}\mu_j^2.
 \label{B:eq:bessel-optimal-rank-error}
\end{equation}
The infimum is over bounded operators on $L^2(0,1)$.
\end{proposition}
\begin{proof}
Let $T$ be a bounded operator of rank $r\le n$, with no self-adjointness or positivity assumed, and let $Q$ be the orthogonal projection onto its range. A bounded finite-rank
operator is Hilbert--Schmidt. Since $(I-Q)T=0$, the two operators $Q(R-T)$ and $(I-Q)R$ have orthogonal ranges. Applying Pythagoras to their values on an orthonormal basis and
summing gives
\[
 \|R-T\|_{\mathfrak S_2}^2
 =\|Q(R-T)\|_{\mathfrak S_2}^2+\|(I-Q)R\|_{\mathfrak S_2}^2
 \ge\|R\|_{\mathfrak S_2}^2-\|QR\|_{\mathfrak S_2}^2.
\]
If $q_1,\ldots,q_r$ is an orthonormal basis of $\operatorname{Ran}Q$, then
\[
 \|QR\|_{\mathfrak S_2}^2
 =\sum_{i=1}^r\|Rq_i\|^2
 =\sum_{i=1}^r\langle R^2q_i,q_i\rangle
 \le\sum_{j=1}^r\mu_j^2.
\]
The last step is \eqref{B:eq:spectral-partial-trace} applied to the positive compact operator $R^2$, whose decreasing eigenvalues are $\mu_j^2$. It follows that
\[
 \|R-T\|_{\mathfrak S_2}^2
 \ge\sum_{j>r}\mu_j^2\ge\sum_{j>n}\mu_j^2.
\]
Conversely, the spectral truncation $R^{[n]}=\sum_{j=1}^n\mu_j\Pi_j$ has rank at most $n$ and
\[
 \|R-R^{[n]}\|_{\mathfrak S_2}^2=\sum_{j>n}\mu_j^2,
\]
because the spectral projections are orthogonal. This attains the lower bound and proves \eqref{B:eq:bessel-optimal-rank-error}. 
\end{proof}

\begin{corollary}\label{B:cor:bessel-reciprocal-tail}
One has $n\|R-R^{[n]}\|_{\mathfrak S_2}^2\longrightarrow0$.
\end{corollary}
\begin{proof}
Write $s_n=\sum_{j>n}\mu_j$. Trace-class summability gives $s_n\to0$. Since $\mu_j$ decreases,
\[
 (n+1-\lfloor n/2\rfloor)\mu_{n+1}
       \le\sum_{j=\lfloor n/2\rfloor+1}^{n+1}\mu_j
       \le s_{\lfloor n/2\rfloor}.
\]
It follows that
\begin{equation}
 n\sum_{j>n}\mu_j^2
 \le n\mu_{n+1}s_n
 \le2s_{\lfloor n/2\rfloor}s_n\longrightarrow0.
 \label{B:eq:spectral-trace-tail-proof}
\end{equation}
Here $n+1-\lfloor n/2\rfloor\ge n/2$, also for $n=1$, which gives the factor two. Both tails on the right tend to zero by trace-class summability, without any assumption of a
polynomial eigenvalue asymptotic. Proposition~\ref{B:prop:bessel-best-finite-rank-approximation} identifies the sum on the left with the squared best rank-$n$ error.
\end{proof}

The lower bound in Theorem~\ref{B:thm:brownian-information-lower-bound} concerns approximations measurable with respect to $\mathcal N_n\vee\mathcal U$, with the Brownian values
and independent auxiliary sigma-field specified there. It applies to the two constructed matrix families by their established measurability. The spectral truncation $R^{[n]}$ is
determined by the full limiting inverse. In fact it cannot satisfy that information restriction for all $n$ with one admissible auxiliary sigma-field: otherwise the lower bound
theorem would give a positive liminf for the left side of \eqref{B:eq:spectral-trace-tail-proof}, whereas that side tends to zero. The squared error of spectral truncation is
therefore smaller than the order $n^{-1}$ attained by the present coupling.

\subsubsection{Parameter restrictions in the convergence rate}
\label{B:subsection:bessel-spectral-parameter-range}

The approximation theorem is proved for fixed $a>0$ and fixed $\beta>0$. The operator construction of Section~\ref{B:sec:inverse-representations-and-the-stochastic-bessel-operator}
has the larger range $a>-1$. The two estimates that impose the smaller range in the rate proof are the choice
\[
 1<\lambda<\min\{a+1,2\}
\]
in Theorem~\ref{B:thm:weighted-inverse-stability}, and the truncation power $1+\nu_a-\eta$, where $\nu_a=\min\{a,1\}$ and $0<\eta<\nu_a$. In particular, with the cutoffs of
Section~\ref{B:sec:removal-of-the-truncation-and-proof-of-the-general-potential-theorem},
\begin{equation}
 nd_n^{1+\nu_a-\eta}
 \le Cr_n^{1+\nu_a-\eta}n^{-(\nu_a-\eta)}\longrightarrow0.
 \label{B:eq:spectral-parameter-truncation}
\end{equation}
For $r_n=\lceil(\log n)^{A_0}\rceil$, the logarithm of the upper bound is
\[
 (1+\nu_a-\eta)A_0\log\log n-(\nu_a-\eta)\log n+O(1),
\]
which tends to $-\infty$ because $\nu_a-\eta>0$. At $a=0$ neither choice of exponents is available. The same calculation explains why a theorem for fixed parameters cannot simply
be evaluated at $a=a_n\downarrow0$. For $0<a_n\le1$ and $\eta_n=a_n/2$, the scalar factor on the right of \eqref{B:eq:spectral-parameter-truncation} is
\[
 r_n^{1+a_n/2}n^{-a_n/2}
 =\exp\{(1+a_n/2)A_0\log\log n-(a_n/2)\log n+O(1)\}.
\]
It tends to zero, for example, if $a_n\log n/\log\log n\to\infty$. This only treats the displayed scalar factor: the constants in the stability and conditional moment estimates
would also have to be controlled uniformly in $a_n$. Obtaining a coupling rate for $-1<a\le0$ requires a different estimate near zero. Likewise, a limit with
$\beta=\beta_n\downarrow0$ requires estimates tracking the dependence on $\beta$ in the conditional moments and comparison constants of
Sections~\ref{B:sec:general-potentials-and-convex-reference-measures}--\ref{B:sec:transport-of-logarithmic-ratios-and-construction-of-the-coupling}.

\subsubsection{Contributions to the squared error}
\label{B:sec:an-exact-leading-error-coefficient}

The error in coupling the coefficients and the error in sampling the Brownian path may both have order $n^{-1/2}$ in norm. Their inner product can therefore contribute at order
$n^{-1}$ to the squared error, alongside the two squared norms.

For the general potential coupling, retain the operators of \eqref{B:eq:full-inverse-comparison} and put
\begin{equation}
 F_n=\widetilde{\mathcal R}_n^0-\mathcal R_n^{B,\mathrm{step},0},
 \qquad
 G_n=\mathcal R_n^{B,\mathrm{step},0}-\mathcal R_n^{B,0},
 \qquad
 E_n=R_n-R-F_n-G_n.
 \label{B:eq:leading-error-differences}
\end{equation}
The estimates for the full inverse give
\[
 \|F_n\|_{\mathfrak S_2}+\|G_n\|_{\mathfrak S_2}=O(n^{-1/2}),
 \qquad \|E_n\|_{\mathfrak S_2}=o(n^{-1/2})
 \quad\text{on }\Omega_0.
\]
Indeed, $E_n$ is the sum of the embedding error, the truncation error for the finite inverse, and the truncation error for the limiting inverse in
\eqref{B:eq:full-inverse-comparison}. Their squared norms are, respectively, $O(n^{-2})$, $o(n^{-1})$, and $o(n^{-1})$.

\begin{proposition}\label{B:prop:bessel-leading-quadratic-error}
For the general potential coupling, almost surely,
\begin{equation}
 ne_n^2={}n\|F_n\|_{\mathfrak S_2}^2+n\|G_n\|_{\mathfrak S_2}^2 
          +2n\operatorname{Re}\operatorname{Tr}(F_n^*G_n)+o(1).
 \label{B:eq:leading-error-quadratic-expansion}
\end{equation}
\end{proposition}
\begin{proof}
All products of the Hilbert--Schmidt operators in this calculation are trace class by the product inequality proved above. Expanding the squared norm of $F_n+G_n+E_n$ gives the
exact difference
\[
 \begin{split}
 ne_n^2-n\|F_n+G_n\|_{\mathfrak S_2}^2
   &=n\|E_n\|_{\mathfrak S_2}^2
       +2n\operatorname{Re}\operatorname{Tr}((F_n+G_n)^*E_n),\\
 \left|ne_n^2-n\|F_n+G_n\|_{\mathfrak S_2}^2\right|
   &\le (\sqrt n\|E_n\|_{\mathfrak S_2})^2\\
   &\quad+2\sqrt n\|E_n\|_{\mathfrak S_2}
         (\sqrt n\|F_n\|_{\mathfrak S_2}+\sqrt n\|G_n\|_{\mathfrak S_2})
    \longrightarrow0.
 \end{split}
\]
The identity $\|F_n+G_n\|_{\mathfrak S_2}^2=\|F_n\|_{\mathfrak S_2}^2+ \|G_n\|_{\mathfrak S_2}^2+2\operatorname{Re}\operatorname{Tr}(F_n^*G_n)$ proves
\eqref{B:eq:leading-error-quadratic-expansion}. If $f_n,g_n$ are the scalar integral kernels of $F_n,G_n$, the mixed product is also given by
\[
 \operatorname{Tr}(F_n^*G_n)
 =\int_0^1\!\int_0^1\overline{f_n(s,t)}g_n(s,t)\,\dd s\,\dd t.
\]
Cauchy--Schwarz in the kernel space makes this integral absolutely convergent and bounds it by $\|F_n\|_{\mathfrak S_2}\|G_n\|_{\mathfrak S_2}$. This is the mixed contribution in
the complete squared error.
\end{proof}

To identify a leading order coefficient it remains to determine the limit of $n\|F_n+G_n\|_{\mathfrak S_2}^2$. For its mixed term, the available estimates give
\[
 2n|\operatorname{Tr}(F_n^*G_n)|
     \le n\|F_n\|_{\mathfrak S_2}^2+n\|G_n\|_{\mathfrak S_2}^2=O(1).
\]
Both differences use the same Brownian path and the coupled coefficients, so this term must be included in the limit. The norm bounds by themselves allow it to remain of leading
order: for a rank-one orthogonal projection $P$, the pairs $F_n=n^{-1/2}P$, $G_n=\pm n^{-1/2}P$ satisfy
\[
 n\|F_n\|_{\mathfrak S_2}^2=n\|G_n\|_{\mathfrak S_2}^2=1,
 \qquad n\|F_n+G_n\|_{\mathfrak S_2}^2=(1\pm1)^2.
\]
Thus the separate norm bounds leave the mixed term undetermined. For the coupling of the bidiagonal matrices, $F_n$ and $G_n$ retain their definitions in
\eqref{B:eq:leading-error-differences}.

For the circular ensemble, Section~\ref{sec:circular-ensembles-canonical-systems-and-the-sine-coefficient} identifies an exact coefficient after proving that the difference between
the circular inverse and the inverse obtained by sampling the Brownian path is $o(n^{-1/2})$ in Hilbert--Schmidt norm. Its square and its inner product with the $O(n^{-1/2})$ error
from sampling the Brownian path therefore vanish after multiplication by $n$. For the coupling of the bidiagonal matrices to the limiting inverse the two terms in
\eqref{B:eq:leading-error-differences} are both currently estimated at order $n^{-1/2}$, and \eqref{B:eq:leading-error-quadratic-expansion} retains their inner product. The result
for the tridiagonal models in Section~\ref{sec:soft-edge-ensembles-and-quantitative-airy-limits} includes a contribution from the omitted spectrum. For the inverse of the
stochastic Bessel operator the squared spectral truncation error is $o(n^{-1})$ by Corollary~\ref{B:cor:bessel-reciprocal-tail}.

A sufficient condition for passing to expectations is uniform integrability of $X_n=ne_n^2$. Suppose that $n\|F_n+G_n\|_{\mathfrak S_2}^2\to X$ in probability and that
\[
 \lim_{A\to\infty}\sup_n\mathbb E[X_n\mathbf1_{\{X_n>A\}}]=0.
\]
Proposition~\ref{B:prop:bessel-leading-quadratic-error} gives $X_n\to X$ in probability. The preceding condition bounds $\sup_n\mathbb EX_n$. An almost surely convergent
subsequence and Fatou's lemma give $\mathbb EX<\infty$ and
\[
 \mathbb E(X-A)_+\le\liminf_k\mathbb E(X_{n_k}-A)_+
       \le\sup_n\mathbb E[X_n\mathbf1_{\{X_n>A\}}].
\]
For fixed $A$ and $\varepsilon>0$,
\[
 \mathbb E|\min(X_n,A)-\min(X,A)|
 \le\varepsilon+A\mathbb P(|X_n-X|>\varepsilon)\longrightarrow\varepsilon.
\]
For nonnegative $x,y$,
\[
 |x-y|\le|\min(x,A)-\min(y,A)|+(x-A)_++(y-A)_+.
\]
Applying this inequality and the two tail bounds gives
\[
 \limsup_n\mathbb E|X_n-X|
 \le\varepsilon+2\sup_n\mathbb E[X_n\mathbf1_{\{X_n>A\}}].
\]
Send $\varepsilon\downarrow0$ and then $A\to\infty$. This proves $L^1$ convergence and hence $\mathbb EX_n\to\mathbb EX$ under the stated additional hypothesis of uniform
integrability. Finiteness of the pathwise supremum does not imply that condition: for $U$ uniform on $(0,1)$, the variables $Y_n=n\mathbf1_{\{U\le1/n\}}$ satisfy
\[
 Y_n\longrightarrow0\text{ almost surely},\qquad
 \sup_nY_n\le U^{-1}<\infty\text{ almost surely},\qquad
 \mathbb EY_n=1.
\]
Passing to expectations therefore requires additional integrability estimates for this coupling.

\clearpage

\subsection{Spectral consequences for the circular ensemble}
\label{sec:combined-s-spectral}
\label{S:sec:spectral-consequences}

Fix $\beta>0$ and work on $\Omega_\beta^{\mathrm{good}}$ from Theorem~\ref{S:thm:common-probability-one-event}. The inverse estimate controls the positive and negative eigenvalues
of the compact inverse operators. To deduce convergence of the eigenvalues of the unbounded differential operators, we restrict to inverse eigenvalues separated from zero.
Comparison of a single eigenspace also requires its spectral gap. On a bounded interval of the spectrum of the differential operator, the estimate depends on the distances of the
two interval boundaries from the spectrum and is uniform over gaps between eigenvalues inside the interval.

All arguments below are deterministic on the common unweighted space. We use the real inverses and their complexifications, whose norms and eigenspaces are identified in
Lemma~\ref{S:lem:spectral-real-complex-inverse}.

Put
\[
 A_n=\mathtt r\mathrm{Circ}_{\beta,n},\qquad
 A_n^{\mathrm{br}}=\mathtt r\mathrm{Br}_{\beta,n},\qquad
 A=\mathtt r\mathrm{Sine}_\beta,
\]
and use the same symbols for the complexifications. Write
\[
 E_n=A_n-A,\qquad
 \epsilon_n^{\mathrm{sp}}=\|E_n\|_{\mathfrak S_2},\qquad
 \delta_n=\|E_n\|_{\mathrm{op}}.
\]
Theorem~\ref{S:thm:main} gives
\begin{equation}\label{S:eq:spectral-inverse-error}
 (\epsilon_n^{\mathrm{sp}})^2
 =\frac{\PathCoeff(\omega)+o_\omega(1)}n,\qquad
 0\le\delta_n\le\epsilon_n^{\mathrm{sp}}
 =O_{\beta,\omega}(n^{-1/2}).
\end{equation}
For the comparison with the sampled operator define the two errors by replacing $A_n$ by $A_n^{\mathrm{br}}$; the same formula holds with the same coefficient. Each deterministic
argument below applies to both sequences.

\subsubsection{Operator domains and reciprocal eigenvalues}
\label{S:sec:spectral-inverses-and-domains}

The comparison is between bounded inverses. We first identify the unbounded operators they invert, with their stated domains, before applying reciprocation and the resolvent
formulas. Each compact inverse is injective and has zero as a spectral accumulation point. The direct operator has zero in its resolvent set.

Use $\mathcal H_{\mathbb R}=L^2([0,1);\mathbb R^2)$ and its complexification $\mathcal H_{\mathbb C}=L^2([0,1);\mathbb C^2)$. The complex inner product is linear in its second
argument. Denote coordinatewise complex conjugation by $\mathfrak c$. If $S$ is a real operator, set
\[
 S_{\mathbb C}(x+\ii y)=Sx+\ii Sy,\qquad
 \mathcal D(S_{\mathbb C})=\mathcal D(S)+\ii\mathcal D(S).
\]
For a bounded operator the domain is the whole space. Theorem~\ref{S:thm:sine-selfadjoint-domain} gives the Green inverse for $\mathrm{Sine}_\beta$; Theorem~\ref{S:thm:3.2} gives
those for the sampled Brownian and circular operators. All three inverses are bounded, self-adjoint and injective. Write $A_X$ for the complexification of any of these inverses,
and put $\mathcal T_X=A_X^{-1}$ on $\operatorname{Ran}A_X$.

\begin{lemma}\label{S:lem:spectral-real-complex-inverse}
The operator $A_X$ is compact and self-adjoint, with dense range. The operator $\mathcal T_X$ is the complexification of the conjugated real differential operator with the domain
specified above, $0\in\rho(\mathcal T_X)$, and, for $\mu\ne0$,
\begin{equation}\label{S:eq:spectral-reciprocal-eigenspaces}
 \ker(A_X-\mu I)=\ker(\mathcal T_X-\mu^{-1}I).
\end{equation}
Real and complex eigenspace dimensions agree.
\end{lemma}

\begin{proof}
For bounded real $S$, $\|S_{\mathbb C}(x+\ii y)\|^2=\|Sx\|^2+\|Sy\|^2$. Taking $y=0$ gives the reverse operator norm inequality, and a real orthonormal basis also serves as a
complex orthonormal basis. Hence
\begin{equation}\label{S:eq:spectral-complexification-norms}
 \|S_{\mathbb C}\|_{\mathrm{op}}=\|S\|_{\mathrm{op}},\qquad
 \|S_{\mathbb C}\|_{\mathfrak S_2(\mathcal H_{\mathbb C})}
 =\|S\|_{\mathfrak S_2(\mathcal H_{\mathbb R})}.
\end{equation}
Complexification therefore preserves norm limits of finite-rank operators and compactness. The identity $(S_{\mathbb C})^*=(S^*)_{\mathbb C}$ preserves self-adjointness, and taking
real and imaginary parts preserves injectivity. Thus $A_X$ is compact, self-adjoint, and injective, with dense range.

Complexifying the two inverse identities of the already constructed real direct operator identifies its complexification with $\mathcal T_X=A_X^{-1}$ on $\operatorname{Ran}A_X$. By
the compact self-adjoint spectral theorem, $A_X$ has an orthonormal eigenbasis $(e_l)$ with nonzero real eigenvalues $\mu_l\to0$. For $v_l=\langle e_l,v\rangle$, consistently with
the inner product being linear in its second argument, $v=A_Xh$ is equivalent to $v_l=\mu_lh_l$ with $(h_l)\in\ell^2$. Therefore
\begin{equation}\label{S:eq:spectral-direct-domain-series}
 \mathcal D(\mathcal T_X)
 =\left\{\sum_lv_le_l:\sum_l|v_l/\mu_l|^2<\infty\right\},
 \qquad
 \mathcal T_X\sum_lv_le_l=\sum_l\mu_l^{-1}v_le_l.
\end{equation}
The diagonal operator with these real diagonal entries is self-adjoint. Its eigenvalues $\mu_l^{-1}$ have finite multiplicity and no finite accumulation. For $z$ outside this set,
the diagonal inverse with coefficients $(\mu_l^{-1}-z)^{-1}$ is bounded and maps into the displayed domain, because
\[
 \left|\frac{\mu_l^{-1}}{\mu_l^{-1}-z}\right|
 \le1+\frac{|z|}{\operatorname{dist}(z,\{\mu_i^{-1}\}_i)}.
\]
It is the two-sided resolvent, so the reciprocals exhaust $\sigma(\mathcal T_X)$ and $0\in\rho(\mathcal T_X)$. The same diagonal description proves
\eqref{S:eq:spectral-reciprocal-eigenspaces}.

Compactness on the infinite-dimensional space excludes a bounded inverse defined on the whole space for $A_X$, while injectivity excludes a zero eigenvector. Thus
$0\in\sigma(A_X)\setminus\sigma_{\mathrm p}(A_X)$. For real $\mu$, the equation $A_X(x+\ii y)=\mu(x+\ii y)$ holds exactly when $A_Xx=\mu x$ and $A_Xy=\mu y$; real eigenspaces
therefore complexify with their dimensions unchanged.
\end{proof}

\subsubsection{Global comparison of the reciprocal eigenvalues}
\label{S:sec:spectral-ordered-eigenvalues}

The circular inverse is not positive. Its positive and negative eigenvalue lists are therefore compared separately, with zero padding used only to match their lengths and never as
a direct eigenvalue to be reciprocated.

For a compact self-adjoint $S$, list its positive eigenvalues as $\mu_1^+(S)\ge\mu_2^+(S)\ge\cdots>0$ and its negative eigenvalues as $\mu_1^-(S)\le\mu_2^-(S)\le\cdots<0$, with
multiplicity. A finite list is continued by zeros, and an empty list is identically zero. Each infinite list tends to zero. For nonzero entries put
$\lambda_j^\pm(S)=1/\mu_j^\pm(S)$. For the inverses considered here, the positive direct eigenvalues are ordered increasingly from zero, and the negative ones decreasingly from
zero. Reciprocals and eigenvectors are assigned only to nonzero entries; padded zeros serve to compare lists of different lengths.

\begin{corollary}\label{S:cor:spectral-reciprocal-matching}
For the full circular inverses,
\begin{equation}\label{S:eq:25.6}
 \sum_{j\ge1}|\mu_j^+(A_n)-\mu_j^+(A)|^2
 +\sum_{j\ge1}|\mu_j^-(A_n)-\mu_j^-(A)|^2
 \le(\epsilon_n^{\mathrm{sp}})^2
 =\frac{\PathCoeff(\omega)+o_\omega(1)}n.
\end{equation}
The same bound holds for the inverses obtained by sampling the Brownian path, with their own Hilbert--Schmidt error and the same coefficient on the right.
\end{corollary}

We apply the min--max principle to the positive and negative eigenvalues separately. Their squared differences are then controlled by the Hoffman--Wielandt inequality for compact
self-adjoint operators.

\begin{lemma}\label{S:lem:spectral-signed-minmax}
For compact self-adjoint $S,T$ and every $j\ge1$,
\begin{equation}\label{S:eq:spectral-signed-lipschitz}
 |\mu_j^\pm(S)-\mu_j^\pm(T)|\le\|S-T\|_{\mathrm{op}}.
\end{equation}
\end{lemma}
\begin{proof}
Let $S$ be compact and self-adjoint on the infinite-dimensional Hilbert space in this subsection. Define
\[
 a_j(S)=\inf_{\dim L=j-1}\ 
       \sup_{\substack{x\perp L\\\|x\|=1}}\langle x,Sx\rangle.
\]
For each finite-dimensional $L$, choose an orthonormal sequence $(x_l)$ in $L^\perp$. It converges weakly to zero, and compactness gives $\|Sx_l\|\to0$. To see the latter
implication, any subsequence with norms bounded away from zero would have, by compactness, a further norm-convergent subsequence; its weak limit must be zero, a contradiction.
Consequently $\langle x_l,Sx_l\rangle\to0$, and the displayed supremum is at least zero.

Suppose first that $S$ has at least $j$ positive eigenvalues. The span of its first $j$ positive eigenvectors intersects $L^\perp$ nontrivially for every $L$ of dimension $j-1$. On
that span the Rayleigh quotient is at least $\mu_j^+(S)$, so $a_j(S)\ge\mu_j^+(S)$. Taking $L$ to be the span of the first $j-1$ positive eigenvectors makes every remaining
spectral value at most $\mu_j^+(S)$. The spectral expansion gives the reverse bound. Thus $a_j(S)=\mu_j^+(S)$.

If $S$ has fewer than $j$ positive eigenvalues, take $L$ containing all their eigenvectors and extend it to dimension $j-1$. The Rayleigh quotient on $L^\perp$ is nonpositive. Its
supremum is nevertheless at least zero by the preceding compactness argument, so it equals zero. This proves $a_j(S)=0=\mu_j^+(S)$, including the case where zero is not an
eigenvalue. Equivalently,
\begin{equation}\label{S:eq:spectral-signed-variational}
 \mu_j^+(S)=\max\left\{0,
 \inf_{\dim L=j-1}\ \sup_{\substack{x\perp L\\\|x\|=1}}
                         \langle x,Sx\rangle\right\}.
\end{equation}
For self-adjoint $S,T$ and unit $x$,
\[
 |\langle x,Sx\rangle-\langle x,Tx\rangle|
 \le\|S-T\|_{\mathrm{op}}.
\]
Taking the supremum over $L^\perp$ and then the infimum over $L$ gives $a_j(S)\le a_j(T)+\|S-T\|$. Reversing $S,T$ proves the estimate for the positive eigenvalue list. Finally
$\mu_j^-(S)=-\mu_j^+(-S)$ under the ordering convention in this subsection, and $\|(-S)-(-T)\|=\|S-T\|$. Applying the positive estimate to $-S,-T$ proves
\eqref{S:eq:spectral-signed-lipschitz} for the negative eigenvalue list as well. This is the min--max principle for compact operators \cite{teschl2009mathematical}.
\end{proof}

Each signed eigenvalue is continuous in operator norm, hence in Hilbert--Schmidt norm. Applied to the measurable inverse operators, this gives measurable real eigenvalues.

We use the Hoffman--Wielandt inequality for compact self-adjoint operators, with positive and negative eigenvalues ordered separately; see equation~(27), p.~490, of
\cite{bhatia1994hoffman}.

\begin{lemma}\label{S:lem:25.1}
For self-adjoint $S,T\in\mathfrak S_2$,
\begin{equation}\label{S:eq:spectral-signed-hw}
 \sum_{j\ge1}|\mu_j^+(S)-\mu_j^+(T)|^2
 +\sum_{j\ge1}|\mu_j^-(S)-\mu_j^-(T)|^2
 \le\|S-T\|_{\mathfrak S_2}^2.
\end{equation}
\end{lemma}
\begin{proof}
We first establish a trace inequality for positive Hilbert--Schmidt operators $C,D$. Write their decreasing positive eigenvalues as $c_j,d_j$, with zeros appended when necessary.
If $u_1,\ldots,u_k$ are orthonormal, the spectral expansion of $D$ gives
\[
 \sum_{i=1}^k\langle u_i,Du_i\rangle
   =\sum_j d_j a_j,\qquad
 0\le a_j\le1,\quad\sum_j a_j\le k.
\]
Splitting the sum at $k$ and using $d_j\le d_k$ for $j>k$ proves $\sum_i\langle u_i,Du_i\rangle\le\sum_{j\le k}d_j$. The argument also applies when $D$ has fewer than $k$ positive
eigenvalues, by zero padding.

Let $C^{[m]}$ be the truncation to the first $m$ positive spectral coordinates of $C$. With $c_{m+1}=0$ for this finite truncation, summation by parts in its eigenbasis gives
\[
 \begin{split}
 \operatorname{Tr}(C^{[m]}D)
 &=\sum_{k=1}^m(c_k-c_{k+1})
                  \sum_{i=1}^k\langle u_i,Du_i\rangle\\
 &\le\sum_{k=1}^m(c_k-c_{k+1})\sum_{i=1}^k d_i
 =\sum_{i=1}^m c_i d_i.
 \end{split}
\]
The trace is nonnegative because each diagonal term is nonnegative. Moreover,
\[
 |\operatorname{Tr}((C-C^{[m]})D)|
 \le\|C-C^{[m]}\|_{\mathfrak S_2}\|D\|_{\mathfrak S_2}\longrightarrow0.
\]
The scalar series $\sum_i c_i d_i$ converges by Cauchy--Schwarz. Passing to the limit therefore proves
\[
 0\le\operatorname{Tr}(CD)\le\sum_{i\ge1}c_i d_i.
\]

Now use the positive and negative parts $S=S_+-S_-$ and $T=T_+-T_-$, where $S_\pm,T_\pm$ are positive Hilbert--Schmidt operators. All four pairwise products are trace class. The
traces of the two products involving a positive part and a negative part are nonnegative, so
\[
 \begin{split}
 \operatorname{Tr}(ST)
 &=\operatorname{Tr}(S_+T_+)+\operatorname{Tr}(S_-T_-)
       -\operatorname{Tr}(S_+T_-)-\operatorname{Tr}(S_-T_+)\\
 &\le\sum_{j\ge1}\mu_j^+(S)\mu_j^+(T)
       +\sum_{j\ge1}\mu_j^-(S)\mu_j^-(T).
 \end{split}
\]
Indeed the decreasing eigenvalues of $S_-$ and $T_-$ are $-\mu_j^-(S)$ and $-\mu_j^-(T)$, whose products have the displayed sign. Since $S_+S_-=0$ and $T_+T_-=0$, the
Hilbert--Schmidt square expands as
\[
 \begin{split}
 \|S-T\|_{\mathfrak S_2}^2
 &=\|S\|_{\mathfrak S_2}^2+\|T\|_{\mathfrak S_2}^2
                 -2\operatorname{Tr}(ST)\\
 &\ge\sum_{j\ge1}|\mu_j^+(S)-\mu_j^+(T)|^2
       +\sum_{j\ge1}|\mu_j^-(S)-\mu_j^-(T)|^2.
 \end{split}
\]
This proves the estimate with constant one. The result is the ordered self-adjoint case of the infinite-dimensional Hoffman--Wielandt inequality \cite{bhatia1994hoffman}.
\end{proof}

\begin{proof}[Proof of Corollary~\ref{S:cor:spectral-reciprocal-matching}]
Fix a sample in the common event of probability one in \eqref{S:eq:spectral-inverse-error}. For the first pair take $S=A_n$ and $T=A$ in Lemma~\ref{S:lem:25.1}. Its right-hand side
is precisely $\|A_n-A\|_{\mathfrak S_2}^2$. Hence
\[
 \sum_{j\ge1}|\mu_j^+(A_n)-\mu_j^+(A)|^2
 +\sum_{j\ge1}|\mu_j^-(A_n)-\mu_j^-(A)|^2
 \le\|A_n-A\|_{\mathfrak S_2}^2.
\]
Apply the same lemma separately with $S=A_n^{\mathrm{br}}$ and $T=A$. The asymptotics for the difference of the inverses in \eqref{S:eq:spectral-inverse-error} supply the claimed
upper bound for each of these two right-hand sides. In particular, multiplication by $n$ and passage to the upper limit give the corresponding bound by the coefficient in that
equation. The comparison is an inequality; it does not identify the spectral squared-error limit with the operator squared-error coefficient.

By \eqref{S:eq:spectral-reciprocal-eigenspaces}, a nonzero inverse eigenvalue $\mu$ corresponds to the direct eigenvalue $1/\mu$, with the same multiplicity and eigenspace.
Positive and negative signs are treated separately, in the orders specified before the corollary. A zero appended to a finite list is only a comparison value for the inverse
eigenvalue sequence; it is never interpreted as the reciprocal of a finite direct eigenvalue. These observations give exactly the two conclusions of the corollary concerning
reciprocal eigenvalues.
\end{proof}

An isospectral perturbation illustrates the possible strictness of \eqref{S:eq:25.6}. Let $S=\operatorname{diag}(1,2)$ on $\mathbb R^2$, let $U_\theta$ be rotation through
$\theta$, and set $T_\theta=U_\theta S U_\theta^{\mathsf T}$. With $s=\sin\theta$, $c=\cos\theta$, direct multiplication gives
\[
 T_\theta-S=
 \begin{pmatrix}s^2&-sc\\-sc&-s^2\end{pmatrix},\qquad
 \|T_\theta-S\|_{\mathfrak S_2}^2
 =2s^4+2s^2c^2=2\sin^2\theta.
\]
Both eigenvalue lists are unchanged. Append to both operators the same diagonal operator $D e_j=(-1)^j(j+2)^{-1}e_j$ on an orthogonal infinite-dimensional space. The resulting
operators are injective, self-adjoint and Hilbert--Schmidt, and their signed lists still agree. For $\theta=n^{-1/2}$ the normalized squared operator difference tends to $2$, while
the spectral error is zero.

\subsubsection{Resolvents and isolated spectral projections}
\label{S:sec:spectral-isolated-subspaces}

For an isolated inverse eigenvalue, a contour at half the gap converts the inverse error into a projection error. The dependence on the gap comes from the two resolvents, and is
recorded explicitly below.

For a closed set $F\subset\mathbb C$, write $\operatorname{dist}(z,F)=\inf_{w\in F}|z-w|$. In the following norm estimates use $\|\cdot\|_{\mathfrak S_\infty}$ to denote the
operator norm. Let $S,T$ be bounded self-adjoint operators, put $\delta=\|T-S\|_{\mathrm{op}}$, and set $d(z)=\operatorname{dist}(z,\sigma(S))$.

\begin{lemma}\label{S:lem:spectral-neumann-inclusion}
The Hausdorff distance between $\sigma(S)$ and $\sigma(T)$ is at most $\delta$. If $d(z)>\delta$, then
\begin{equation}\label{S:eq:spectral-resolvent-perturbation}
 \begin{split}
 \|(z-T)^{-1}\|_{\mathrm{op}}&\le\frac1{d(z)-\delta},\\
 \|(z-T)^{-1}-(z-S)^{-1}\|_{\mathfrak S_p}
 &\le\frac{\|T-S\|_{\mathfrak S_p}}{d(z)(d(z)-\delta)},
 \qquad p\in\{2,\infty\}.
 \end{split}
\end{equation}
The estimate for $p=2$ assumes $T-S\in\mathfrak S_2$.
\end{lemma}
\begin{proof}
The spectral theorem gives $\|(z-S)^{-1}\|=d(z)^{-1}$, with zero included in $\sigma(S)$ when $S$ is compact on this space. If $d(z)>\delta$, the factorization
\[
 z-T=\bigl[I-(T-S)(z-S)^{-1}\bigr](z-S)
\]
and the Neumann series give $\|(z-T)^{-1}\|\le(d(z)-\delta)^{-1}$. Consequently $\sigma(T)$ lies in the closed $\delta$-neighborhood of $\sigma(S)$; interchanging $S,T$ proves the
Hausdorff bound. Finally,
\[
 (z-T)^{-1}-(z-S)^{-1}
 =(z-T)^{-1}(T-S)(z-S)^{-1}.
\]
For $p=2,\infty$, the ideal inequality $\|BEC\|_{\mathfrak S_p}\le\|B\|\|E\|_{\mathfrak S_p}\|C\|$ proves \eqref{S:eq:spectral-resolvent-perturbation}.
\end{proof}

For our inverse pair this gives, in particular,
\begin{equation}\label{S:eq:25.9}
 \sigma(A_n)\subset
 \{x\in\mathbb R:\operatorname{dist}(x,\sigma(A))\le\delta_n\},
\end{equation}
and the reverse inclusion with the roles interchanged. Both spectra contain zero.

Fix a nonzero eigenvalue $\mu$ of $A$. Let $m$ be its multiplicity and put
\[
 g=\operatorname{dist}(\mu,\sigma(A)\setminus\{\mu\}),\qquad
 \Gamma=\{\mu+(g/2)\mathrm e^{\ii t}:0\le t\le2\pi\},
\]
with positive orientation. Compactness and injectivity give $0<g\le|\mu|$ and $m<\infty$. Whenever $\delta_n<g/4$, the Riesz projections are
\begin{equation}\label{S:eq:spectral-riesz-projections}
 P=\frac1{2\pi\ii}\int_\Gamma(z-A)^{-1}\,\dd z,\qquad
 P_n=\frac1{2\pi\ii}\int_\Gamma(z-A_n)^{-1}\,\dd z.
\end{equation}
The contour method is the one used for separated spectra in Chapter~III, \S6.4, and Chapter~IV, \S3.4 of \cite{kato1995perturbation}.

\begin{theorem}\label{S:thm:spectral-isolated-projections}
If $\delta_n<g/4$, then $P,P_n$ are orthogonal projections of rank $m$, and
\begin{equation}\label{S:eq:25.13}
 \|P_n-P\|_{\mathrm{op}}
 \le\frac{\delta_n}{g/2-\delta_n}\le\frac4g\delta_n<1.
\end{equation}
Moreover,
\begin{equation}\label{S:eq:spectral-riesz-hs-bound}
 \|P_n-P\|_{\mathfrak S_2}
 \le\frac{\epsilon_n^{\mathrm{sp}}}{g/2-\delta_n}
 \le\frac4g\epsilon_n^{\mathrm{sp}}.
\end{equation}
\end{theorem}
\begin{proof}
The disk bounded by $\Gamma$ avoids zero, and $\operatorname{dist}(z,\sigma(A))\ge g/2$ on $\Gamma$. The preceding lemma gives
\begin{equation}\label{S:eq:25.10}
 \sup_{z\in\Gamma}\|(z-A)^{-1}\|_{\mathrm{op}}\le2/g,\qquad
 \sup_{z\in\Gamma}\|(z-A_n)^{-1}\|_{\mathrm{op}}
 \le\frac1{g/2-\delta_n}\le4/g.
\end{equation}
The spectral theorem and scalar contour integration identify the Riesz integrals as orthogonal projections onto the enclosed eigenspaces. These spaces have finite dimension because
the disk is separated from zero. Integrating \eqref{S:eq:spectral-resolvent-perturbation} over the contour of length $\pi g$ gives, for $p=2,\infty$,
\[
 \|P_n-P\|_{\mathfrak S_p}
 \le\frac{\|E_n\|_{\mathfrak S_p}}{g/2-\delta_n}.
\]
For $p=2$ the resolvent difference is continuous in $\mathfrak S_2$, so the integral is taken in that norm. Since
$\vartheta_n:=\|P_n-P\|_{\mathrm{op}}<1$,
\[
 v\in\operatorname{Ran}P,\quad P_nv=0
 \quad\Longrightarrow\quad
 \|v\|\le\vartheta_n\|v\|\quad\Longrightarrow\quad v=0.
\]
Thus $P_n|_{\operatorname{Ran}P}$ is injective. Interchanging the projections proves equality of their finite ranks.

Each enclosed $\nu\in\sigma(A_n)$ lies within $\delta_n$ of some $\lambda\in\sigma(A)$. If $\lambda\ne\mu$, then $|\nu-\lambda|\ge g-|\nu-\mu|>g/2>\delta_n$, a contradiction. Hence
the enclosed eigenvalues, counted with multiplicity, obey
\begin{equation}\label{S:eq:25.15}
 \max_{1\le l\le m}|\mu_{l,n}-\mu|\le\delta_n.
\end{equation}
The eigenspaces and their orthogonal complements are invariant under $\mathfrak c$, so the projections restrict to real orthogonal projections of the same ranks.
\end{proof}

\subsubsection{Individual eigenvalues of the differential operators}
\label{S:sec:spectral-fixed-direct-eigenvalues}

Reciprocation preserves each nonzero eigenspace. Its effect on the error depends on the distance of the inverse eigenvalue from zero. For the isolated eigenvalue considered above,
write $\lambda=\mu^{-1}$ and $\lambda_{l,n}=\mu_{l,n}^{-1}$. Since $\delta_n<g/4\le|\mu|/4$, equation~\eqref{S:eq:25.15} gives $|\mu_{l,n}|\ge|\mu|-\delta_n>0$. Therefore
\begin{equation}\label{S:eq:spectral-direct-eigenvalue-rate}
 \max_{1\le l\le m}|\lambda_{l,n}-\lambda|
 \le\frac{\delta_n}{|\mu|(|\mu|-\delta_n)}
 \le\frac{2\delta_n}{|\mu|^2}
 =O_{\beta,\omega,\mu}(n^{-1/2}).
\end{equation}
For each fixed nonzero signed label, choose $n_j$ so that $\delta_n<|\mu_j^\pm(A)|/2$ for every $n\ge n_j$. Equation~\eqref{S:eq:spectral-signed-lipschitz} then makes the finite
label nonzero with the same sign, and the displayed reciprocation calculation applies. The eigenspaces and projections are unchanged by reciprocation. The absolute derivative of
$\mu\mapsto\mu^{-1}$ is $|\mu|^{-2}$. This factor is bounded on sets separated from zero; it diverges as
$\mu\to0$.

There is no corresponding global unweighted $\ell^2$ estimate for differences of eigenvalues of the unbounded operators. For example, on $\ell^2$ set $Se_j=j^{-1}e_j$ and
$S_ne_j=(1+n^{-1})j^{-1}e_j$. These are compact self-adjoint injective operators, and
\[
 \|S_n-S\|_{\mathfrak S_2}^2
 =n^{-2}\sum_{j\ge1}j^{-2}\longrightarrow0,\qquad
 \sum_{j\ge1}|\lambda_j^+(S_n)-\lambda_j^+(S)|^2
 =\frac1{(n+1)^2}\sum_{j\ge1}j^2=\infty.
\]
Here the inverse operators converge in Hilbert--Schmidt norm while the unweighted squared displacement of the direct spectra is infinite.

\subsubsection{Resolvents of the differential operators}
\label{S:sec:spectral-direct-resolvents}

We express each direct resolvent as a bounded function of its inverse. This uses the controlled perturbation $E_n$ and accommodates the different domains of $\mathcal T_n$ and
$\mathcal T$. The same representation permits contours enclosing an interval of the spectrum of $\mathcal T$ that contains zero.

Let $\mathcal T=A^{-1}$ and $\mathcal T_n=A_n^{-1}$ on their ranges, as in Lemma~\ref{S:lem:spectral-real-complex-inverse}. Their domains can differ. Fix a nonempty compact
$K\subset\rho(\mathcal T)$ and set
\[
 R_K=\sup_{z\in K}|z|,\qquad
 M_K=\sup_{z\in K}\|(I-zA)^{-1}\|_{\mathrm{op}}.
\]
The proof below shows that $M_K<\infty$. For an empty $K$ the supremum statements are vacuous.

\begin{proposition}\label{S:prop:spectral-direct-resolvents}
For all sufficiently large $n$, $K\subset\rho(\mathcal T_n)$ and
\begin{equation}\label{S:eq:spectral-direct-resolvent-rate}
 \sup_{z\in K}
 \|(\mathcal T_n-z)^{-1}-(\mathcal T-z)^{-1}\|_{\mathfrak S_2}
 \le2M_K^2\epsilon_n^{\mathrm{sp}}
 =O_{\beta,\omega,K}(n^{-1/2}).
\end{equation}
The operator norm bound has $\delta_n$ in place of $\epsilon_n^{\mathrm{sp}}$.
\end{proposition}
\begin{proof}
The diagonal description \eqref{S:eq:spectral-direct-domain-series} gives, including at $z=0$,
\begin{equation}\label{S:eq:spectral-direct-bounded-formula}
 (\mathcal T-z)^{-1}=A(I-zA)^{-1},\qquad
 (I-zA)^{-1}=I+z(\mathcal T-z)^{-1}.
\end{equation}
Continuity of inversion on the open set of invertible operators and compactness of $K$ imply $M_K<\infty$. Choose $n$ so large that $R_KM_K\delta_n\le1/2$. Then
\[
 I-zA_n=\bigl[I-zE_n(I-zA)^{-1}\bigr](I-zA)
\]
is invertible throughout $K$, with inverse norm at most $2M_K$. Applying the same diagonal description to $A_n$ shows $K\subset\rho(\mathcal T_n)$ and
$(\mathcal T_n-z)^{-1}=A_n(I-zA_n)^{-1}$, with range in $\mathcal D(\mathcal T_n)=\operatorname{Ran}A_n$. Subtracting gives
\begin{equation}\label{S:eq:spectral-direct-resolvent-identity}
 (\mathcal T_n-z)^{-1}-(\mathcal T-z)^{-1}
 =(I-zA_n)^{-1}E_n(I-zA)^{-1}.
\end{equation}
Indeed multiplication by $I-zA_n$ on the left and $I-zA$ on the right reduces the identity to $A_n(I-zA)-(I-zA_n)A=E_n$. The ideal inequality proves both bounds.
\end{proof}

\subsubsection{Bounded spectral intervals}
\label{S:sec:spectral-bounded-intervals}

A contour enclosing a bounded spectral interval compares its entire invariant subspace at once. The resulting estimate depends on the distance of the two ends of the interval from
the spectrum, not on gaps between the eigenvalues inside it. It therefore remains useful when those eigenvalues approach one another.

Let $I=[a,b]$, where $a<b$ and $a,b\in\rho(\mathcal T)$. Define
\[
 L=\max\{|a|,|b|\}>0,\qquad
 d_I=\min\{1,\operatorname{dist}(a,\sigma(\mathcal T)),
                     \operatorname{dist}(b,\sigma(\mathcal T))\}>0.
\]
Let $\Gamma_I$ be the positively oriented boundary of the rectangle $[a,b]+\ii[-1,1]$, and put
\[
 R_I=\sqrt{L^2+1},\qquad
 M_I=1+R_I/d_I,\qquad
 C_I=\frac{2(b-a)+4}{\pi}M_I^2.
\]
Write $\mathcal P_I$ and $\mathcal P_{n,I}$ for the spectral projections of $\mathcal T$ and $\mathcal T_n$ on $I$. For definiteness, take $n$ large enough that
\[
 R_IM_I\delta_n\le\tfrac12,\qquad C_I\delta_n<1,\qquad
 \delta_n<(2L)^{-1},\qquad 2L^2\delta_n<d_I.
\]
They ensure, respectively, invertibility on the contour, equality of projection ranks, legitimate reciprocation, and that matched eigenvalues stay inside the interval.

\begin{theorem}\label{S:thm:spectral-fixed-window}
For all sufficiently large $n$, the two eigenvalue counts in $I$ agree, with multiplicity. The eigenvalues admit a pairing of maximal error at most $2L^2\delta_n$, and
\begin{equation}\label{S:eq:spectral-window-projections}
 \|\mathcal P_{n,I}-\mathcal P_I\|_{\mathrm{op}}
 +\|\mathcal P_{n,I}-\mathcal P_I\|_{\mathfrak S_2}
 \le C_I(\delta_n+\epsilon_n^{\mathrm{sp}})
 =O_{\beta,\omega,I}(n^{-1/2}).
\end{equation}
\end{theorem}
\begin{proof}
On a vertical side $z=c+\ii v$, $c=a,b$,
\[
 \operatorname{dist}(z,\sigma(\mathcal T))^2
 =\operatorname{dist}(c,\sigma(\mathcal T))^2+v^2;
\]
on a horizontal side, $|\operatorname{Im}z|=1$. Thus throughout $\Gamma_I$,
\[
 \operatorname{dist}(z,\sigma(\mathcal T))\ge d_I,\qquad
 |z|\le R_I,\qquad
 \|(I-zA)^{-1}\|\le1+|z|\|(\mathcal T-z)^{-1}\|\le M_I.
\]
The factorization in Proposition~\ref{S:prop:spectral-direct-resolvents} therefore puts $\Gamma_I$ in $\rho(\mathcal T_n)$ when $R_IM_I\delta_n\le1/2$ and bounds the resolvent
difference by $2M_I^2\|E_n\|_{\mathfrak S_p}$, $p=2,\infty$. The spectral theorem gives
\[
 \mathcal P_I=-\frac1{2\pi\ii}
              \int_{\Gamma_I}(\mathcal T-z)^{-1}\,\dd z,\qquad
 \mathcal P_{n,I}=-\frac1{2\pi\ii}
              \int_{\Gamma_I}(\mathcal T_n-z)^{-1}\,\dd z.
\]
Both $a$ and $b$ belong to the respective resolvent sets. The projections have finite rank because the direct spectra have no finite accumulation. Since $|\Gamma_I|=2(b-a)+4$,
\begin{equation}\label{S:eq:spectral-window-endpoint-bounds}
 \|\mathcal P_{n,I}-\mathcal P_I\|_{\mathrm{op}}\le C_I\delta_n,
 \qquad
 \|\mathcal P_{n,I}-\mathcal P_I\|_{\mathfrak S_2}
                        \le C_I\epsilon_n^{\mathrm{sp}}.
\end{equation}
If $C_I\delta_n<1$, the injectivity argument in Theorem~\ref{S:thm:spectral-isolated-projections} equates their ranks, including when the rank is zero.

Assign every limiting eigenvalue in $I$ its nonzero signed inverse label. Since $|\lambda_j^\pm(A)|\le L$, $|\mu_j^\pm(A)|\ge L^{-1}$. For $\delta_n<(2L)^{-1}$,
Lemma~\ref{S:lem:spectral-signed-minmax} makes the corresponding label of $A_n$ nonzero with the same sign and gives
\[
 |\lambda_j^\pm(A_n)-\lambda_j^\pm(A)|
 =\frac{|\mu_j^\pm(A_n)-\mu_j^\pm(A)|}
        {|\mu_j^\pm(A_n)\mu_j^\pm(A)|}
 \le2L^2\delta_n.
\]
Every limiting eigenvalue in $I$ is at least $d_I$ from $\partial I$. Thus $2L^2\delta_n<d_I$ puts all matched eigenvalues in $I$, and equality of counts exhausts the approximating
spectrum there. This proves the pairing and \eqref{S:eq:spectral-window-projections}.

For moving intervals $I_n=[a_n,b_n]$ inside a fixed bounded interval, suppose that $a_n,b_n\in\rho(\mathcal T)$. Choose $L_*\ge1$ with $|a_n|,|b_n|\le L_*$ for every $n$, and put
\[
 d_n^I=\operatorname{dist}(\{a_n,b_n\},\sigma(\mathcal T)),\qquad
 \bar d_n=\min\{1,d_n^I\},\qquad R_*=\sqrt{L_*^2+1}.
\]
The rectangles used above have length at most $4L_*+4$. Their horizontal sides have spectral distance at least one, and their vertical sides have distance at least $d_n^I$.
Therefore
\[
 M_{I_n}\le1+R_*/\bar d_n\le(1+R_*)\bar d_n^{-1},\qquad
 C_{I_n}\le C_{L_*}\bar d_n^{-2}.
\]
Assume $\epsilon_n^{\mathrm{sp}}/(d_n^I)^2\to0$. Since $\epsilon_n^{\mathrm{sp}}\to0$, this also gives
\[
 q_n:=\epsilon_n^{\mathrm{sp}}/\bar d_n^2\longrightarrow0.
\]
Indeed, $q_n$ is the maximum of $\epsilon_n^{\mathrm{sp}}$ and $\epsilon_n^{\mathrm{sp}}/(d_n^I)^2$. Using $\delta_n\le\epsilon_n^{\mathrm{sp}}$ and $\bar d_n\le1$, we obtain
\[
 \begin{aligned}
 R_{I_n}M_{I_n}\delta_n&\le C_{L_*}q_n\bar d_n\le C_{L_*}q_n,\\
 C_{I_n}\delta_n&\le C_{L_*}q_n,\\
 \delta_n&\le q_n\bar d_n^2\le q_n,\\
 2L_*^2\delta_n/\bar d_n&\le2L_*^2q_n\bar d_n\le2L_*^2q_n.
 \end{aligned}
\]
Choose $N$ so that these four upper bounds are, respectively, at most $1/2$, less than one, less than $(2L_*)^{-1}$, and less than one for $n\ge N$. They are the four conditions
needed for contour invertibility, rank equality, reciprocation and containment of the matched eigenvalues. The preceding proof thus applies separately at each such $n$ and gives
\[
 \operatorname{rank}\mathcal P_{n,I_n}=\operatorname{rank}\mathcal P_{I_n},\qquad
 \|\mathcal P_{n,I_n}-\mathcal P_{I_n}\|_{\mathfrak S_2}
 \le C_{L_*}\epsilon_n^{\mathrm{sp}}/\bar d_n^2\longrightarrow0.
\]
For this fixed limiting operator, $d_n^I$ is bounded above: choose any $\lambda_*\in\sigma(\mathcal T)$ and use $d_n^I\le L_*+|\lambda_*|$. Consequently the same estimate may be
written with $(d_n^I)^2$ in the denominator after increasing the path-dependent constant. Since $\epsilon_n^{\mathrm{sp}}=O_{\beta,\omega}(n^{-1/2})$, the condition
$n^{1/4}d_n^I\to\infty$ is sufficient. This compares two projections on the same moving interval; it does not imply that $\mathcal P_{I_n}$ converges to a fixed projection.
\end{proof}

The constant in \eqref{S:eq:spectral-window-endpoint-bounds} is determined by the interval and the distances of its two boundary points from the limiting spectrum. The interval may
contain zero. In addition to the contour argument, the bound $\|A_n\|\le\|A\|+1$ for large $n$ gives $|\lambda_n|\ge(\|A\|+1)^{-1}$ for every direct eigenvalue, so there is no
direct spectral accumulation at zero. The same interval argument applies to the operators obtained by sampling the Brownian path.

\subsubsection{Eigenvectors and first variations}
\label{S:sec:spectral-first-variations}

Return to the isolated inverse eigenvalue $\mu$, projection $P$, and gap $g$ of Subsection~\ref{S:sec:spectral-isolated-subspaces}. For a simple eigenvalue, choose the real unit
eigenvectors so that $\langle\phi_n,\phi\rangle\ge0$. For a multiple eigenvalue we compare the whole invariant subspace.

\begin{corollary}\label{S:cor:spectral-eigenvectors}
For the fixed isolated eigenvalue, both projection differences and all enclosed eigenvalue errors are $O_{\beta,\omega,\mu,g}(n^{-1/2})$. If $m=1$, each real unit eigenvector
$\phi$ admits real unit eigenvectors $\phi_n$ satisfying
\begin{equation}\label{S:eq:spectral-eigenvector-bound}
 \langle\phi_n,\phi\rangle\ge0,\qquad
 \|\phi_n-\phi\|\le\sqrt2\|P_n-P\|_{\mathrm{op}}
 \le\frac{4\sqrt2}{g}\delta_n.
\end{equation}
\end{corollary}
\begin{proof}
The rates follow from \eqref{S:eq:spectral-inverse-error} and Theorem~\ref{S:thm:spectral-isolated-projections}. For $m=1$, $P_n\phi\ne0$ by the rank argument, so put
$\phi_n=P_n\phi/\|P_n\phi\|$. It is a real unit vector spanning the enclosed one-dimensional spectral subspace, and hence an eigenvector. Let $s_n=\|(I-P_n)\phi\|$. Orthogonality
gives
\[
 \|P_n\phi\|=\sqrt{1-s_n^2},\qquad
 \langle\phi_n,\phi\rangle=\sqrt{1-s_n^2}.
\]
Thus
\[
 \|\phi_n-\phi\|^2
 =2(1-\sqrt{1-s_n^2})
 =\frac{2s_n^2}{1+\sqrt{1-s_n^2}}
 \le2s_n^2\le2\|P_n-P\|_{\mathrm{op}}^2.
\]
For complex eigenvectors, the choice $\phi_n=P_n\phi/\|P_n\phi\|$ gives $\langle\phi_n,\phi\rangle=\|P_n\phi\|>0$ and the same norm bound.
\end{proof}

Let $\mathcal H_\perp=(I-P)\mathcal H_{\mathbb C}$ and set
\[
 \mathcal S_\mu=(\mu I-A|_{\mathcal H_\perp})^{-1}(I-P),\qquad
 \mathcal S_\mu P=P\mathcal S_\mu=0,\qquad
 \|\mathcal S_\mu\|_{\mathrm{op}}\le g^{-1}.
\]
The inverse exists because $g$ is a lower bound on the distance from $\mu$ to the complementary spectrum, including zero. The sign convention is $\mu I-A$, not $A-\mu I$. We will
prove, for $\delta_n<g/4$,
\begin{equation}\label{S:eq:spectral-projection-first-variation}
 P_n-P=\mathcal S_\mu E_nP+PE_n\mathcal S_\mu+\mathcal R_{n,\mu},\qquad
 \|\mathcal R_{n,\mu}\|_{\mathfrak S_p}
 \le\frac8{g^2}\delta_n\|E_n\|_{\mathfrak S_p},\quad p\in\{2,\infty\}.
\end{equation}
Here $\mathfrak S_\infty$ denotes the operator norm. This is the isolated spectral projection expansion; compare \cite{kato1995perturbation}.

Write $R(z)=(z-A)^{-1}$ and $R_n(z)=(z-A_n)^{-1}$. The identities
\[
 R_n-R=RE_nR_n,\qquad
 R_n-R-RE_nR=RE_n(R_n-R)=RE_nRE_nR_n
\]
hold with the factors in this order. On the circle $\Gamma=\{z:|z-\mu|=g/2\}$,
\[
 \|R(z)\|\le2/g,\qquad
 \|R_n(z)\|\le(g/2-\delta_n)^{-1}.
\]
Putting either one of the two factors $E_n$ in $\mathfrak S_p$ and the other in operator norm yields
\[
 \|R(z)E_nR(z)E_nR_n(z)\|_{\mathfrak S_p}
 \le\frac{4\delta_n\|E_n\|_{\mathfrak S_p}}{g^2(g/2-\delta_n)}.
\]
The contour has length $\pi g$. After multiplication by $1/(2\pi)$ and integration, its contribution is at most
\[
 \frac{2\delta_n\|E_n\|_{\mathfrak S_p}}{g(g/2-\delta_n)}
 \le\frac8{g^2}\delta_n\|E_n\|_{\mathfrak S_p}.
\]
The integral is an integral in the relevant norm; its integrand is norm-continuous on the contour.

To compute the linear term, decompose
\[
 R(z)=\frac{P}{z-\mu}+R_\perp(z),\qquad
 R_\perp(\mu)=\mathcal S_\mu.
\]
The function $R_\perp$ is analytic in $|z-\mu|<g$. Multiplication gives
\[
 R(z)E_nR(z)=\frac{PE_nP}{(z-\mu)^2}
 +\frac{PE_nR_\perp(z)+R_\perp(z)E_nP}{z-\mu}
 +R_\perp(z)E_nR_\perp(z).
\]
The first summand has zero contour integral, as does the analytic last summand. Cauchy's formula applied to the two middle numerators gives
$PE_n\mathcal S_\mu+\mathcal S_\mu E_nP$. Integrating the resolvent difference proves \eqref{S:eq:spectral-projection-first-variation}.

The argument applies to every self-adjoint $E\in\mathfrak S_2$ with $\|E\|_{\mathrm{op}}<g/4$. The linear map $E\mapsto\mathcal S_\mu EP+PE\mathcal S_\mu$ has $\mathfrak S_2$ norm
at most $2\|E\|_{\mathfrak S_2}/g$, and the remainder divided by $\|E\|_{\mathfrak S_2}$ is at most $8\|E\|_{\mathrm{op}}/g^2$. Since $\|E\|_{\mathrm{op}}\le\|E\|_{\mathfrak S_2}$,
this quotient tends to zero as $E\to0$ in $\mathfrak S_2$. This proves Fr\'echet differentiability on the real space of self-adjoint Hilbert--Schmidt operators. Under
\eqref{S:eq:spectral-inverse-error}, the remainder is $O_{\beta,\omega,\mu,g}(n^{-1})$ in $\mathfrak S_2$. For a simple eigenvalue, the corresponding eigenvalue and normalized
eigenvector satisfy
\begin{equation}\label{S:eq:spectral-simple-first-variation}
 |\mu_n-\mu-\langle\phi,E_n\phi\rangle|
 \le\frac{\delta_n^2}{g-2\delta_n}\le\frac{2\delta_n^2}{g},\qquad
 \|\phi_n-\phi-\mathcal S_\mu E_n\phi\|
 \le\frac{8\delta_n^2}{g^2}.
\end{equation}
We prove these estimates for $\delta_n<g/4$. The isolated eigenvalue comparison gives $|\mu_n-\mu|\le\delta_n$. Choose the phase as in Corollary~\ref{S:cor:spectral-eigenvectors}.
Then $t_n=\langle\phi,\phi_n\rangle>0$, and we may write
\[
 \phi_n=t_n\phi+w_n,\qquad w_n\perp\phi,\qquad
 C_n=(I-P)A_n|_{\mathcal H_\perp},\qquad v_n=(I-P)E_n\phi.
\]
Put $B_\perp=\mu I-A|_{\mathcal H_\perp}$ and
\[
 D_n=(\mu_n-\mu)I-(I-P)E_n|_{\mathcal H_\perp}.
\]
Thus $\mu_n-C_n=B_\perp+D_n$, $\|B_\perp^{-1}\|\le g^{-1}$ and $\|D_n\|\le2\delta_n$. A Neumann series gives
\[
 \|(\mu_n-C_n)^{-1}\|\le(g-2\delta_n)^{-1}.
\]
Projecting $A_n(t_n\phi+w_n)=\mu_n(t_n\phi+w_n)$ onto $\mathcal H_\perp$ and onto $\operatorname{span}\{\phi\}$ gives, respectively,
\[
 (\mu_n-C_n)w_n=t_nv_n,\qquad
 t_n(\mu_n-\mu)=t_n\langle\phi,E_n\phi\rangle+\langle v_n,w_n\rangle.
\]
The phase convention permits division by $t_n$. Substitution yields the exact identities
\begin{equation}\label{S:eq:spectral-simple-schur-identities}
 \frac{w_n}{t_n}=(\mu_n-C_n)^{-1}v_n,\qquad
 \mu_n-\mu=\langle\phi,E_n\phi\rangle
       +\langle v_n,(\mu_n-C_n)^{-1}v_n\rangle.
\end{equation}
Since $\|v_n\|\le\delta_n$, the second identity proves the eigenvalue estimate in \eqref{S:eq:spectral-simple-first-variation}.

For the eigenvector, set $u_n=w_n/t_n$. The inverse identity on $\mathcal H_\perp$ gives
\[
 (B_\perp+D_n)^{-1}-B_\perp^{-1}
 =-(B_\perp+D_n)^{-1}D_nB_\perp^{-1}.
\]
Because $B_\perp^{-1}v_n=\mathcal S_\mu E_n\phi$, we obtain
\[
 \|u_n-\mathcal S_\mu E_n\phi\|
 \le\frac{2\delta_n^2}{g(g-2\delta_n)}\le\frac{4\delta_n^2}{g^2},\qquad
 \|u_n\|\le\frac{2\delta_n}{g}.
\]
Normalization gives $t_n=(1+\|u_n\|^2)^{-1/2}$. For $x\ge0$,
\[
 1-(1+x)^{-1/2}
 =\frac{x}{\sqrt{1+x}(1+\sqrt{1+x})}\le x/2.
\]
Consequently $0\le1-t_n\le2\delta_n^2/g^2$. Writing $v=\mathcal S_\mu E_n\phi$ and using $\phi_n=t_n(\phi+u_n)$, we find
\[
 \begin{split}
 \|\phi_n-\phi-v\|
 &\le(1-t_n)(1+\|v\|)+t_n\|u_n-v\|\\
 &\le\frac{2\delta_n^2}{g^2}(1+\delta_n/g)+\frac{4\delta_n^2}{g^2}
 \le\frac{8\delta_n^2}{g^2}.
 \end{split}
\]
This proves the second estimate, including the effect of normalization. A fluctuation theorem would additionally require a limit law for $\langle\phi,E_n\phi\rangle$.

To return to direct eigenvalues, put $\lambda=\mu^{-1}$, $\lambda_n=\mu_n^{-1}$ and $d_n=\mu_n-\mu$. The gap includes the distance to zero, so $g\le|\mu|$ and
$|\mu_n|\ge|\mu|-\delta_n\ge|\mu|/2$. The identity
\[
 \frac1{\mu+d_n}-\frac1\mu+\frac{d_n}{\mu^2}
       =\frac{d_n^2}{\mu^2(\mu+d_n)}
\]
therefore implies
\[
 \begin{split}
 \left|\lambda_n-\lambda+
             \frac{\langle\phi,E_n\phi\rangle}{\mu^2}\right|
 &\le\frac{|d_n-\langle\phi,E_n\phi\rangle|}{|\mu|^2}
       +\frac{|d_n|^2}{|\mu|^2|\mu_n|}\\
 &\le\left(\frac2{g|\mu|^2}+\frac2{|\mu|^3}\right)\delta_n^2
   =O_{\beta,\omega,\mu,g}(n^{-1}).
 \end{split}
\]
For $m>1$, choose an index beyond which $\vartheta_n=\|P_n-P\|_{\mathrm{op}}\le1/2$ and the ranks of $P_n$ and $P$ are both $m$. This is possible by
Theorem~\ref{S:thm:spectral-isolated-projections}; for example, the sufficient bound $\delta_n\le g/8$ gives $\vartheta_n\le1/2$. Regard $H_n=PP_nP$ as an operator on
$\operatorname{Ran}P$. For $v$ in this space,
\[
 \langle v,H_nv\rangle=\|P_nv\|^2
 =\|v\|^2-\|(I-P_n)v\|^2,
\]
so $(1-\vartheta_n^2)I\le H_n\le I$. Its inverse square root is defined by the finite-dimensional spectral theorem. Set
\begin{equation}\label{S:eq:spectral-subspace-isometry}
 U_n=P_nPH_n^{-1/2}:\operatorname{Ran}P\longrightarrow\operatorname{Ran}P_n.
\end{equation}
Then
\[
 U_n^*U_n=H_n^{-1/2}PP_n^2PH_n^{-1/2}=I.
\]
The map is injective between two spaces of the same finite dimension and is therefore onto. In particular $U_nU_n^*=P_n$. To estimate its distance from the identity inclusion $P$,
write
\[
 U_n-P=(P_n-P)PH_n^{-1/2}+P(H_n^{-1/2}-I).
\]
The spectral bounds on $H_n$ imply
\[
 \|U_n-P\|\le\frac{\vartheta_n}{\sqrt{1-\vartheta_n^2}}
          +\frac1{\sqrt{1-\vartheta_n^2}}-1\le2\vartheta_n.
\]
For the last inequality, divide by $\vartheta_n>0$, rationalize the second numerator, and use $\sqrt{1-\vartheta_n^2}\ge\sqrt3/2$; if $\vartheta_n=0$, then $P_n=P$ and $U_n=P$ on
$\operatorname{Ran}P$, so both sides vanish. Thus the image under $U_n$ of an orthonormal basis of $\operatorname{Ran}P$ is an orthonormal basis of the perturbed invariant
subspace.

Put $W_n=(I-P)U_n$. Since $AP=PA=\mu P$ and $U_n^*U_n=P$ on $\operatorname{Ran}P$,
\[
 U_n^*AU_n-\mu P=W_n^*(A-\mu I)W_n.
\]
For the perturbation part, direct expansion gives
\[
 U_n^*E_nU_n-PE_nP
 =(U_n^*-P)E_nU_n+PE_n(U_n-P).
\]
Adding these identities proves
\[
 \begin{split}
 U_n^*A_nU_n-\mu P-PE_nP
 ={}&W_n^*(A-\mu I)W_n\\
 &+(U_n^*-P)E_nU_n+PE_n(U_n-P).
 \end{split}
\]
Here $\|W_n\|\le\|U_n-P\|\le2\vartheta_n$, $\|U_n\|=1$, and $\|P\|=1$. Hence
\[
 \|U_n^*A_nU_n-\mu P-PE_nP\|
 \le4(\|A\|+|\mu|)\vartheta_n^2+4\delta_n\vartheta_n
 =O_{A,\mu,g}(\delta_n^2).
\]
The eigenvalues of $U_n^*A_nU_n$ are exactly the $m$ eigenvalues enclosed by the contour, with multiplicity, because $U_n$ is unitary onto that invariant subspace. For self-adjoint
matrices of dimension $m$, the ordered min--max formula changes by at most the operator norm of the perturbation. Applying it in both directions pairs these eigenvalues with those
of $\mu P+PE_nP$ with the displayed error. Thus the first-order splitting is determined by the eigenvalues of the full matrix $PE_nP$, not by its diagonal entries in an arbitrarily
chosen basis. The nearby basis supplied by $U_n$ need not diagonalize $A_n$ inside a split cluster. The vectors above lie in the common unweighted space. The canonical
eigenfunctions are $X_n^{-1}\phi_n$ and $X^{-1}\phi$. For $t<1$, whenever these multipliers are bounded on $[0,t]$,
\[
 \|X_n^{-1}\phi_n-X^{-1}\phi\|_{L^2(0,t)}
 \le\|X_n^{-1}\|_{L^\infty(0,t)}\|\phi_n-\phi\|_2
    +\|X_n^{-1}-X^{-1}\|_{L^\infty(0,t)}\|\phi\|_2.
\]
Such local comparisons therefore also use local convergence of the coordinate matrices. Bounds up to the singular boundary or for boundary evaluations require further control of
those multipliers.

\subparagraph{First variation of the direct resolvent} Retain $R_K=\sup_{z\in K}|z|$ and $M_K=\sup_{z\in K}\|(I-zA)^{-1}\|$ for a compact set $K\subset\rho(\mathcal T)$. Put
$Q(z)=I-zA$ and $Q_n(z)=I-zA_n$. Choose an index beyond which $R_KM_K\delta_n\le1/2$. The factorization of $Q_n$ about $Q$ gives $\sup_K\|Q_n(z)^{-1}\|\le2M_K$. Spectral calculus
and the bounded inverse identity imply
\[
 (\mathcal T_n-z)^{-1}-(\mathcal T-z)^{-1}
       =Q_n(z)^{-1}E_nQ(z)^{-1}.
\]
Subtract $Q(z)^{-1}E_nQ(z)^{-1}$. Since
\[
 Q_n(z)^{-1}-Q(z)^{-1}
       =zQ_n(z)^{-1}E_nQ(z)^{-1},
\]
the remainder is exactly
\[
 zQ_n(z)^{-1}E_nQ(z)^{-1}E_nQ(z)^{-1}.
\]
The factors $E_n$ remain in their displayed order. The ideal inequality, with one $E_n$ in operator norm and the other in $\mathfrak S_2$, gives
\begin{equation}\label{S:eq:spectral-direct-resolvent-first-variation}
 \begin{split}
 \sup_{z\in K}\bigl\|&(\mathcal T_n-z)^{-1}-(\mathcal T-z)^{-1}
       -(I-zA)^{-1}E_n(I-zA)^{-1}\bigr\|_{\mathfrak S_2}\\
 &\le2R_KM_K^3\delta_n\epsilon_n^{\mathrm{sp}}.
 \end{split}
\end{equation}
This is $O_{\beta,\omega,K}(n^{-1})$ under the inverse comparison. More generally, after division by $\|E_n\|_{\mathfrak S_2}$ the bound tends to zero whenever $E_n\to0$ in
$\mathfrak S_2$. Thus the displayed linear term is the derivative, uniform on $K$, with respect to the bounded inverse perturbation. It is not a derivative with respect to a
difference of the unbounded direct operators. At $z=0$ the remainder vanishes identically, as the exact inverse difference requires.
\subsubsection{Measurable projections and eigenvectors}
\label{S:sec:spectral-measurable-choices}

All preceding estimates are deterministic on each path where the inverse bound holds. They apply there to every eligible eigenvalue, compact subset of the resolvent set, and
bounded spectral interval. We now verify measurability of the choices.

Choose $\mu$ by a signed label and restrict to $\{\mu\ne0\}$. Lemma~\ref{S:lem:spectral-signed-minmax} makes the labels measurable, and hence also the positive gap and finite
multiplicity
\[
 \begin{gathered}
 d_j^\pm=
 \begin{cases}|\mu-\mu_j^\pm(A)|,&\mu_j^\pm(A)\ne\mu,\\
 +\infty,&\mu_j^\pm(A)=\mu,\end{cases}
 \qquad g=\min\{|\mu|,\inf_{j,\pm}d_j^\pm\},\\
 m=\sum_{j\ge1}\bigl(\mathbf1_{\{\mu_j^+(A)=\mu\}}
                       +\mathbf1_{\{\mu_j^-(A)=\mu\}}\bigr).
 \end{gathered}
\]
Repeated copies of $\mu$ are excluded from the gap; its distance to zero remains included.

For $z\ne0$ in the resolvent set of $A$, the identity
\[
 (z-A)^{-1}=z^{-1}I+z^{-1}(z-A)^{-1}A
\]
separates a scalar identity from a Hilbert--Schmidt operator. The circle $\Gamma$ centered at $\mu$ with radius $g/2$ does not enclose zero. Its integral of $z^{-1}I$ is therefore
zero, and
\begin{equation}\label{S:eq:spectral-measurable-riesz-integral}
 P=\frac1{2\pi\ii}\int_\Gamma z^{-1}(z-A)^{-1}A\,\dd z.
\end{equation}
Parametrize the circle by $z(t)=\mu+(g/2)\mathrm e^{\ii t}$, $0\le t\le2\pi$. The random scalars $\mu,g$ are measurable, $g\le|\mu|$, and
\[
 |z(t)|\ge|\mu|/2,\qquad
 \operatorname{dist}(z(t),\sigma(A))\ge g/2.
\]
Thus
\[
 \|z(t)^{-1}(z(t)-A)^{-1}A\|_{\mathfrak S_2}
 \le\frac4{|\mu|g}\|A\|_{\mathfrak S_2}.
\]
For each fixed $t$, the maps $(A,\mu,g)\mapsto z(t)$ and $z(t)I-A$ are measurable in their norm topologies. Inversion is continuous on the open set of invertible bounded operators.
Multiplication of a bounded operator by a Hilbert--Schmidt operator is continuous into $\mathfrak S_2$. The resulting integrand, including $z'(t)$, is therefore a measurable
$\mathfrak S_2$-valued variable. It is continuous in $t$ on each sample. Its equispaced Riemann sums are measurable finite sums and converge in $\mathfrak S_2$ to $P$. Since this
space is separable, their limit is strongly measurable.

On $\{\delta_n<g/4\}$ the same contour lies in $\rho(A_n)$, with distance at least $g/2-\delta_n$. Replacing $A$ by $A_n$ in the integral proves measurability of $P_n$ on that
measurable event; define $P_n=0$ outside it. The projection estimate gives $\vartheta_n<1$ there. On $\operatorname{Ran}P$, the spectrum of $H_n=PP_nP$ lies in
$[1-\vartheta_n^2,1]$. With $c_l=\binom{2l}{l}/4^l\in(0,1]$, the scalar power series for $(1-x)^{-1/2}$ gives
\[
 H_n^{-1/2}P=\sum_{l\ge0}c_l(P-H_n)^lP.
\]
This notation means the inverse square root on $\operatorname{Ran}P$, extended by zero on its orthogonal complement. Each term has range in an $m$-dimensional space, so its
Hilbert--Schmidt norm is at most $\sqrt m\,\vartheta_n^{2l}$. In particular,
\[
 \left\|H_n^{-1/2}P-\sum_{l=0}^Nc_l(P-H_n)^lP\right\|_{\mathfrak S_2}
 \le\frac{\sqrt m\,\vartheta_n^{2N+2}}{1-\vartheta_n^2}\longrightarrow0.
\]
The partial sums are measurable polynomials in $P,P_n$. Their limit is measurable, and multiplying it by $P_n$ proves measurability of $U_n=P_nPH_n^{-1/2}$. Set this operator to
zero off $\{\delta_n<g/4\}$. No choice of eigenvectors of $H_n$ is required for this construction. Fix a real orthonormal basis $(e_l)$ of the ambient separable Hilbert space. We
construct the basis of $\operatorname{Ran}P$ on the measurable events $\{m\ge r\}$, one vector at a time. After $\phi_1,\ldots,\phi_{r-1}$ have been selected, put
\[
 v_l^{(r)}=Pe_l-\sum_{i<r}\phi_i\langle\phi_i,Pe_l\rangle,\qquad
 k_r=\min\{l:\|v_l^{(r)}\|>0\},\qquad
 \phi_r=\frac{v_{k_r}^{(r)}}{\|v_{k_r}^{(r)}\|}.
\]
The minimum is finite on $\{m\ge r\}$. Otherwise every $Pe_l$ would lie in the span of the preceding $r-1$ vectors. Taking closed spans would force $\operatorname{Ran}P$ into that
space, contrary to its dimension. For a fixed $l$,
\[
 \{k_r=l\}=\{m\ge r\}\cap\{\|v_l^{(r)}\|>0\}
                   \cap\bigcap_{j<l}\{\|v_j^{(r)}\|=0\}.
\]
All these events are measurable by induction. The selected normalized vector is a countable measurable selection among $v_l^{(r)}/\|v_l^{(r)}\|$, on the sets where the denominator
is positive. Define it to be zero on $\{m<r\}$. This proves measurability at each step. Orthogonal subtraction and normalization show that the first $m$ selected vectors form a
real orthonormal basis of $\operatorname{Ran}P$.

On $\{\delta_n<g/4\}$, the vectors $U_n\phi_1,\ldots,U_n\phi_m$ form a measurable orthonormal basis of $\operatorname{Ran}P_n$, since $U_n$ is an onto isometry. These vectors need
not be eigenvectors individually when the cluster splits. For $m=1$, write $\phi=\phi_1$ and use instead
\[
 \phi_n=\frac{P_n\phi}{\|P_n\phi\|}.
\]
The denominator is positive: if it vanished, then $\|(P-P_n)\phi\|=1$, contradicting $\vartheta_n<1$. Its scalar product with $\phi$ equals $\|P_n\phi\|>0$, so this construction
also fixes the phase measurably. Set $\phi_n=\phi$ on the complementary event. Corollary~\ref{S:cor:spectral-eigenvectors} applies to these exact vectors and gives
\eqref{S:eq:spectral-eigenvector-bound} whenever $\delta_n<g/4$. On every sample for which $\delta_n\to0$, the gap $g$ is positive; hence there is a finite index $n_\mu$ such that
this bound holds for every $n\ge n_\mu$. The arguments use only the original inverse operators and their measurability, so \eqref{S:eq:conditional-marginal-descent} transfers the
conclusions to their original coupling. Passing to expectations requires uniform integrability.

\subsection{Dependence on the noise amplitude and further directions}
\label{sec:combined-a-parameters}
\label{A:sec:parameter-dependence-and-further-directions}
\label{A:sec:B}

For a fixed Brownian path, the dependence on $\beta$ enters the limiting quadratic form through the coefficient $c_\beta=2\beta^{-1/2}$ multiplying the Brownian term. We compare
the associated operators using these forms. The finite matrices require a different argument: their coefficient laws change with $\beta$, and the almost sure bounds already proved
have random constants and starting sizes. We return to the finite size estimates to control their moments uniformly on a compact parameter interval. Combining the two comparisons
then allows $\beta$ to vary with the matrix size.

The continuum comparison holds for every pair of parameters on one set of Brownian paths. For the finite models the supremum over $\beta$ is outside the $L^p$ norm under each
specified coupling law. It is not a pathwise supremum over sampled matrices.

Fix $I_\beta=[\beta_-,\beta_+]\Subset(0,\infty)$ and a nonempty compact set $K\Subset\mathbb C\setminus\mathbb R$. Write $\delta_K=\operatorname{dist}(K,\mathbb R)>0$. Moment
orders satisfy $1\le p<\infty$; bounds for smaller positive moments follow from monotonicity. The continuum estimates hold for every path in the set $\Omega_B$ defined in
Subsection~\ref{A:sec:finite-difference-approximation-of-the-stochastic-airy-operator}. For the finite models we use the probability kernels of
Subsection~\ref{A:sec:measurable-couplings-of-the-finite-models}, without changing them when the moment order or the required failure probability changes.

\subsubsection{Quadratic forms and resolvent continuity}\label{A:sec:continuum-forms-and-quantitative-continuity}
Use the forms \eqref{A:eq:airy-form} on $\mathcal D_{\mathrm A}$ with the norm \eqref{A:eq:airy-energy}. For $B\in\Omega_B$, put $M_B=1+\mathcal O_{1/4}(B)$ and
$P_B=1+C_{I_\beta}M_B^4$, where $C_{I_\beta}$ is fixed below. We write $\mathfrak q_{\beta,B}$ when displaying the path dependence. The identity
$\mathfrak q_{\beta,B}-\mathfrak q_{\gamma,B} =(c_\beta-c_\gamma)\mathfrak w_B$ is an identity of forms on this common domain. Since $B'$ is a distribution, we do not regard the
difference as a bounded potential on $L^2$. We first make the forms uniformly positive by one shift and compare their bounded form inverses.

\begin{proposition}\label{A:res:B.1}
One can choose $1\le\kappa_B\le P_B$ so that, uniformly for $\beta\in I_\beta$,
\[
 c\|f\|_{\mathcal E}^2
 \le\mathfrak q_{\beta,B}[f]+\kappa_B\|f\|_2^2
 \le CP_B\|f\|_{\mathcal E}^2,
 \qquad f\in\mathcal D_{\mathrm A}.
 \eqtag{B.9}
\]
The forms are closed and their associated self-adjoint operators have compact resolvent.
\end{proposition}

\begin{proof}
Apply \eqref{A:eq:form-relative} with $\alpha=1/4$. For $0<\varepsilon\le1$,
\[
 |\mathfrak w_B[f]|\le\varepsilon\mathfrak a[f]
       +C\varepsilon^{-3}M_B^4\|f\|_2^2.
\]
Put $c_*=2\beta_-^{-1/2}$ and $\varepsilon_* =\min\{1,(2c_*)^{-1}\}$. Then $c_\beta\varepsilon_*\le1/2$ for all $\beta\in I_\beta$. Choose
\[
 C_{I_\beta}=c_*C\varepsilon_*^{-3},\qquad
 \kappa_B=P_B=1+C_{I_\beta}M_B^4.
\]
Using $\mathfrak q_{\beta,B}=\mathfrak a+c_\beta\mathfrak w_B$, we obtain
\[
 \tfrac12\mathfrak a[f]+\|f\|_2^2
 \le\mathfrak q_{\beta,B}[f]+\kappa_B\|f\|_2^2
 \le\tfrac32\mathfrak a[f]+(2P_B-1)\|f\|_2^2.
\]
Since $\|f\|_{\mathcal E}^2=\mathfrak a[f]+\|f\|_2^2$, this proves \eqref{A:eq:B.9} with lower constant $1/2$ and upper constant $2P_B$.

We verify the operator consequences on the same domain. The space $\mathcal D_{\mathrm A}$, equipped with $\|\cdot\|_{\mathcal E}$, is complete: a Cauchy sequence converges in
$H^1(0,L)$ for every finite $L$, its derivatives converge in $L^2(\mathbb R_+)$, and its products with $\sqrt{1+x}$ converge in $L^2(\mathbb R_+)$. These limits agree on bounded
intervals, and the trace at zero remains zero by continuity of the $H^1(0,L)$ trace. The limit consequently belongs to $\mathcal D_{\mathrm A}$ and convergence holds in the energy
norm. The domain is dense in $L^2$, since it contains $C_c^\infty(0,\infty)$. The two bounds above show that the positive form norm of $\mathfrak q_{\beta,B}+\kappa_B\|\cdot\|_2^2$
is equivalent to this complete norm. Hence each form is densely defined and closed, and its representation gives a self-adjoint operator $\mathcal H_\beta$.

For compactness, let $(f_j)$ be bounded in the energy norm by $C$. On each interval $(0,L)$, it is bounded in $H^1(0,L)$ and therefore has a subsequence convergent in $L^2(0,L)$. A
diagonal choice gives one subsequence with this property for every positive integer $L$. Moreover,
\[
 \int_L^\infty|f_j(x)|^2\,\mathrm dx
 \le L^{-1}\int_L^\infty x|f_j(x)|^2\,\mathrm dx\le C^2/L.
\]
For two members of the chosen subsequence, the squared $L^2$ norm of their difference on $(L,\infty)$ is at most $4C^2/L$. First make this tail small by choosing $L$, and then use
convergence on $(0,L)$. The subsequence is Cauchy in $L^2(\mathbb R_+)$, proving that the inclusion $\mathcal D_{\mathrm A}\hookrightarrow L^2$ is compact.

Finally, if $u=(\mathcal H_\beta+\kappa_B)^{-1}f$, the weak equation and the lower form bound give
\[
 \tfrac12\|u\|_{\mathcal E}^2
 \le |\langle f,u\rangle|\le\|f\|_2\|u\|_{\mathcal E}.
\]
Thus the shifted inverse maps $L^2$ boundedly into the energy space. Its composition with the compact inclusion is compact. The resolvent identity then gives compactness at every
point of the resolvent set. The domain and the shift are common to all $\beta\in I_\beta$.
\end{proof}

Keep this same $\kappa_B$ for every $\beta\in I_\beta$. Let $A_\kappa=-\partial_x^2+x+\kappa_B$ with Dirichlet boundary condition. Its form domain is $\mathcal D_{\mathrm A}$, and
$A_\kappa^{1/2}$ maps this domain onto $L^2$, with inverse $A_\kappa^{-1/2}$. Since $\kappa_B\ge1$,
\[
 \mathcal L\le A_\kappa\le\kappa_B\mathcal L
\]
in the sense of quadratic forms. The relative bound for $\mathfrak w_B$, followed by polarization, makes
\[
 (u,v)\longmapsto\mathfrak w_B[A_\kappa^{-1/2}u,A_\kappa^{-1/2}v]
\]
a bounded Hermitian sesquilinear form on $L^2$. The Riesz representation theorem therefore defines a unique bounded self-adjoint $T_B$ by
\[
 \langle T_Bu,v\rangle
 =\mathfrak w_B[A_\kappa^{-1/2}u,A_\kappa^{-1/2}v].
\]
Put $C_\beta=I+c_\beta T_B$. For $f=A_\kappa^{-1/2}u$, the lower estimate in \eqref{A:eq:B.9} gives
\[
 \langle C_\beta u,u\rangle
 =\mathfrak q_{\beta,B}[f]+\kappa_B\|f\|_2^2
 \ge c\langle\mathcal L^{1/2}f,\mathcal L^{1/2}f\rangle
 \ge(c/\kappa_B)\|u\|_2^2.
\]
The last step uses $A_\kappa\le\kappa_B\mathcal L$. Thus $C_\beta$ is invertible and $\|C_\beta^{-1}\|\le\kappa_B/c$.

To identify the inverse, fix $F\in L^2$ and write $f=A_\kappa^{-1/2}u$. The weak equation for $(\mathcal H_\beta+\kappa_B)f=F$ is
\[
 \mathfrak q_{\beta,B}[f,g]+\kappa_B\langle f,g\rangle=\langle F,g\rangle,
 \qquad g\in\mathcal D_{\mathrm A}.
\]
Substituting $g=A_\kappa^{-1/2}v$, where $v$ ranges over all of $L^2$, gives
\[
 \langle C_\beta u,v\rangle=\langle A_\kappa^{-1/2}F,v\rangle.
\]
Consequently $u=C_\beta^{-1}A_\kappa^{-1/2}F$. Returning to $f$ proves the bounded operator identity
\[
 S_\beta^\kappa=(\mathcal H_\beta+\kappa_B)^{-1}
   =A_\kappa^{-1/2}(I+c_\beta T_B)^{-1}A_\kappa^{-1/2}.
\]
The argument uses the common form domain; it does not represent the distribution $B'$ by a bounded multiplication operator.

\par\Needspace{6\baselineskip}
\begin{proposition}\label{A:res:B.2}
For $\beta,\gamma\in I_\beta$,
\[
 \|S_\beta^\kappa-S_\gamma^\kappa\|_{\mathfrak S_2}
 \le C_{B,I_\beta}|c_\beta-c_\gamma|.
\]
\end{proposition}

\begin{proof}
Write $V_\theta=(I+c_\theta T_B)^{-1}$. The inverse identity gives
\[
 V_\beta-V_\gamma=(c_\gamma-c_\beta)V_\beta T_BV_\gamma,
\]
and hence
\[
 S_\beta^\kappa-S_\gamma^\kappa
 =(c_\gamma-c_\beta)A_\kappa^{-1/2}
                   V_\beta T_BV_\gamma A_\kappa^{-1/2}.
\]
The lower bound preceding the proposition implies $\sup_\theta\|V_\theta\|\le\kappa_B/c$. We estimate the two factors $A_\kappa^{-1/2}$ in $\mathfrak S_4$. The Schatten exponents
satisfy $1/2=1/4+1/4$. If $\nu_j$ are the eigenvalues of $\mathcal L$, the deterministic calculation in Lemma~\ref{A:res:weyl} gives $\nu_j\ge c j^{2/3}$. The eigenvalues of
$A_\kappa$ are $\nu_j+\kappa_B-1$, so
\[
 \|A_\kappa^{-1/2}\|_{\mathfrak S_4}^4
 =\sum_{j\ge1}(\nu_j+\kappa_B-1)^{-2}
 \le C\sum_{j\ge1}j^{-4/3}<\infty.
\]
Schatten H\"older and the ideal property now give
\[
 \|S_\beta^\kappa-S_\gamma^\kappa\|_{\mathfrak S_2}
 \le |c_\beta-c_\gamma|\,
 \|A_\kappa^{-1/2}\|_{\mathfrak S_4}^2
 \|V_\beta\|\|T_B\|\|V_\gamma\|,
\]
which proves the bound.
\end{proof}

For $R_\theta(z)=(\mathcal H_\theta-z)^{-1}$, put $A_\theta(z)=[I-(z+\kappa_B)S_\theta^\kappa]^{-1}$. For
$\lambda\in\sigma(\mathcal H_\theta)$,
\[
 \left(1-\frac{z+\kappa_B}{\lambda+\kappa_B}\right)^{-1}
 =1+\frac{z+\kappa_B}{\lambda-z}.
\]
At the accumulation point zero of $\sigma(S_\theta^\kappa)$ the inverse multiplier equals one. Thus spectral calculus proves
\[
 R_\theta(z)=S_\theta^\kappa A_\theta(z),
 \eqtag{B.20}
\]
\[
 \sup_{\theta\in I_\beta,z\in K}\|A_\theta(z)\|
 \le1+\frac{\kappa_B+\sup_K|z|}{\delta_K},
 \eqtag{B.21}
\]
and the inverse identity gives
\[
 A_\beta(z)-A_\gamma(z)
 =A_\beta(z)(z+\kappa_B)
              (S_\beta^\kappa-S_\gamma^\kappa)A_\gamma(z).
 \eqtag{B.23}
\]

\begin{theorem}\label{A:res:B.3}
For every $B\in\Omega_B$ and $\beta,\gamma\in I_\beta$,
\[
 \sup_{z\in K}\|R_\beta(z)-R_\gamma(z)\|_{\mathfrak S_2}
 \le C_{B,I_\beta,K}|\beta^{-1/2}-\gamma^{-1/2}|.
\]
\end{theorem}

\begin{proof}
Put $d=z+\kappa_B$. The shifted inverses at different parameters need not commute. Multiplying on the left by $I-dS_\beta^\kappa$ and on the right by $I-dS_\gamma^\kappa$ gives
\[
 (I-dS_\beta^\kappa)(R_\beta-R_\gamma)(I-dS_\gamma^\kappa)
 =S_\beta^\kappa(I-dS_\gamma^\kappa)
   -(I-dS_\beta^\kappa)S_\gamma^\kappa
 =S_\beta^\kappa-S_\gamma^\kappa.
\]
 Therefore $R_\beta-R_\gamma=A_\beta(S_\beta^\kappa-S_\gamma^\kappa)A_\gamma$. Proposition~\ref{A:res:B.2}, \eqref{A:eq:B.21}, and
 $|c_\beta-c_\gamma|=2|\beta^{-1/2}-\gamma^{-1/2}|$ give the result.
\end{proof}

Let $\Lambda_k(\theta)$ denote the increasing eigenvalues of $\mathcal H_\theta(B)$, including multiplicities.

\begin{proposition}\label{A:res:B.4}
For every $k\ge1$ and $\beta,\gamma\in I_\beta$,
\[
 |\Lambda_k(\beta)-\Lambda_k(\gamma)|
 \le C_{B,I_\beta,k}|c_\beta-c_\gamma|.
 \eqtag{B.26}
\]
\end{proposition}

\begin{proof}
The min--max principle and \eqref{A:eq:B.9} imply $\Lambda_k(\theta)+\kappa_B\le CP_B\nu_k$. Let $E_{k,\theta}$ be the span of the first $k$ orthonormal eigenfunctions at parameter
$\theta$. For a unit vector in this space,
\[
 c\|f\|_{\mathcal E}^2
 \le\mathfrak q_{\theta,B}[f]+\kappa_B
 \le\Lambda_k(\theta)+\kappa_B\le CP_B\nu_k.
\]
The relative form bound at a fixed value of $\varepsilon$ gives $|\mathfrak w_B[f]|\le CP_B\|f\|_{\mathcal E}^2$. On this unit sphere, the preceding energy estimate therefore gives
$|\mathfrak w_B[f]|\le C_kP_B^2$. Use $E_{k,\beta}$ as a trial subspace for the form at $\gamma$. Consequently
\[
 \begin{aligned}
\Lambda_k(\gamma)& \le\max_{f\in E_{k,\beta},\ \|f\|_2=1}\mathfrak q_{\gamma,B}[f] \le\max_{f\in E_{k,\beta},\ \|f\|_2=1} \{\mathfrak q_{\beta,B}[f] +|c_\gamma-c_\beta|\,|\mathfrak w_B[f]|\}\\
&\le\Lambda_k(\beta)+C_kP_B^2|c_\gamma-c_\beta|.
\end{aligned}
\]
Interchange $\beta$ and $\gamma$ for the reverse inequality.
\end{proof}

Equation~\eqref{A:eq:form-relative}, with $\alpha=1/4$, implies
\[
 |\mathfrak w_B[f]|
 \le\varepsilon\mathfrak a[f]
       +C(1+\varepsilon^{-4})M_B^4\|f\|_2^2,
 \qquad \varepsilon>0.
 \eqtag{B.30}
\]
The moment comparison requires an integrable upper bound for the pathwise Lipschitz constant. We obtain it in terms of $M_B$ by working on the fixed energy space of $\mathcal L$.

\par\Needspace{6\baselineskip}
\begin{proposition}\label{A:res:B.5}
The constants in Theorem~\ref{A:res:B.3} and Proposition~\ref{A:res:B.4} have versions with every finite Wiener moment. In particular,
\[
 \|\Lambda_k(\beta)-\Lambda_k(\gamma)\|_{L^p(\mathbb W)}
 \le C_{p,k,I_\beta}|c_\beta-c_\gamma|.
\]
\end{proposition}

\begin{proof}
The common comparison is
\[
 c\mathcal L\le\mathcal H_\theta+\kappa_B\le CP_B\mathcal L,
  \qquad 1\le\kappa_B\le P_B,\quad\theta\in I_\beta.
 \eqtag{continuum-polynomial-form}
\]
Define the bounded operator $T$ by the sesquilinear form
\[
 \langle Tf,g\rangle
 =\mathfrak w_B[\mathcal L^{-1/2}f,\mathcal L^{-1/2}g],
 \qquad f,g\in L^2(\mathbb R_+).
\]
Let $F_\theta(z)$ be the bounded extension of $\mathcal L^{1/2}R_\theta(z)\mathcal L^{1/2}$, initially defined on $\mathcal D_{\mathrm A}$. Its boundedness follows from the
factorization below. Equation~\eqref{A:eq:B.30}, at $\varepsilon=1$, gives $\|T\|\le CP_B$. Put $A_\theta=\mathcal H_\theta+\kappa_B$. The lower bound gives
$\|\mathcal L^{1/2}A_\theta^{-1/2}\|\le c^{-1/2}$, whereas

\[
 \|A_\theta^{1/2}R_\theta(z)A_\theta^{1/2}\|
 =\sup_{\lambda\in\sigma(\mathcal H_\theta)}
              \frac{\lambda+\kappa_B}{|\lambda-z|}
 \le1+\frac{\kappa_B+|z|}{\delta_K}.
\]
Factoring $F_\theta$ through $A_\theta^{1/2}$ proves $\sup_{\theta,z}\|F_\theta(z)\|\le C_KP_B$. To subtract the weak equations, take $f,g\in L^2$ and put
\[
 u_\gamma=R_\gamma(z)f,\qquad u_\beta=R_\beta(z)f,\qquad
 v_\beta=R_\beta(\overline z)g.
\]
All three vectors belong to the common form domain. Write
\[
 a_{\theta,z}[u,v]=\mathfrak q_{\theta,B}[u,v]-z\langle u,v\rangle.
\]
With the scalar product linear in the first variable, the resolvent and adjoint weak equations read
\[
 a_{\gamma,z}[u_\gamma,v]=\langle f,v\rangle,\qquad
 a_{\beta,z}[u_\beta,v]=\langle f,v\rangle,\qquad
 a_{\beta,z}[w,v_\beta]=\langle w,g\rangle
\]
for all vectors $v,w$ in the form domain. Since
\[
 a_{\beta,z}[u,v]-a_{\gamma,z}[u,v]
       =(c_\beta-c_\gamma)\mathfrak w_B[u,v],
\]
subtracting the first two equations and then testing with $v_\beta$ gives
\[
 \begin{split}
 \langle(R_\beta(z)-R_\gamma(z))f,g\rangle
 &=a_{\beta,z}[u_\beta-u_\gamma,v_\beta]\\
 &=(c_\gamma-c_\beta)\mathfrak w_B[u_\gamma,v_\beta].
 \end{split}
\]
This computation is valid before any multiplication by an unbounded power of $\mathcal L$.

By the definition of $T$ and the bounded extensions $F_\theta(z)$,
\[
 \mathcal L^{1/2}u_\gamma=F_\gamma(z)\mathcal L^{-1/2}f,
 \qquad
 \mathcal L^{1/2}v_\beta=F_\beta(\overline z)\mathcal L^{-1/2}g.
\]
Moreover $F_\beta(\overline z)^*=F_\beta(z)$, first by the resolvent adjoint identity on the form domain and then by bounded extension. Hence
\[
 \begin{split}
 \mathfrak w_B[u_\gamma,v_\beta]
 &=\langle T F_\gamma(z)\mathcal L^{-1/2}f,
                  F_\beta(\overline z)\mathcal L^{-1/2}g\rangle\\
 &=\langle\mathcal L^{-1/2}F_\beta(z)T F_\gamma(z)
                       \mathcal L^{-1/2}f,g\rangle.
 \end{split}
\]
Equality for all $f,g$ proves
\[
 R_\beta(z)-R_\gamma(z)
 =(c_\gamma-c_\beta)\mathcal L^{-1/2}
          F_\beta(z)T F_\gamma(z)\mathcal L^{-1/2}.
 \eqtag{continuum-polynomial-resolvent}
\]
This bounded operator identity follows from the weak equations, with the potential represented by the bounded operator $T$ between the displayed form factors. Since
$\|\mathcal L^{-1/2}\|_{\mathfrak S_4}^4 =\sum_j\nu_j^{-2}<\infty$, Schatten H\"older gives
\[
 \sup_{z\in K}\|R_\beta(z)-R_\gamma(z)\|_{\mathfrak S_2}
 \le C_{I_\beta,K}P_B^3|c_\beta-c_\gamma|.
 \eqtag{continuum-polynomial-Lipschitz}
\]
The preceding min--max calculation also gives
\[
  |\Lambda_k(\beta)-\Lambda_k(\gamma)|
 \le C_{k,I_\beta}P_B^2|c_\beta-c_\gamma|.
 \eqtag{continuum-polynomial-eigenvalue}
\]
On the span of the first $k$ eigenfunctions of $\mathcal H_\gamma$, a unit vector $u$ satisfies

$\langle u,\mathcal L u\rangle \le c^{-1}(\Lambda_k(\gamma)+\kappa_B)\le C_kP_B$, by the upper comparison with $P_B\mathcal L$ and min--max. Hence
$|\mathfrak w_B[u]|\le\|T\|\langle u,\mathcal L u\rangle \le C_kP_B^2$. Testing min--max on this span and then interchanging $\beta,\gamma$ proves the two-sided eigenvalue bound.
The same form calculation gives differentiability in the noise amplitude, with a quadratic remainder. Write $c=c_\gamma$, $s=c_\beta-c_\gamma$, and keep $T$ and $F_\theta(z)$
above. Subtracting the bounded form equations gives $F_\beta-F_\gamma=-sF_\beta T F_\gamma$. Substitution once more in the resolvent difference yields
\[
 R_\beta-R_\gamma
 =-s\mathcal L^{-1/2}F_\gamma T F_\gamma\mathcal L^{-1/2}
   +s^2\mathcal L^{-1/2}F_\beta T F_\gamma T F_\gamma
                                      \mathcal L^{-1/2}.
\]
Consequently, uniformly in $z\in K$ and $\beta,\gamma\in I_\beta$,
\[
 \|R_\beta(z)-R_\gamma(z)
       +s\mathcal L^{-1/2}F_\gamma(z)T F_\gamma(z)\mathcal L^{-1/2}
   \|_{\mathfrak S_2}
 \le C_{I_\beta,K}P_B^5s^2.
\]
There are three factors $F$ and two factors $T$ in the remainder; the two exterior factors are estimated in $\mathfrak S_4$. This proves the derivative with respect to $c$, with
one-sided interpretation at the ends of the amplitude interval. It uses the common form domain and does not differentiate the Brownian path or regard $B'$ as a bounded
multiplication operator. These polynomial upper bounds are Borel functions of $B$. By Lemma~\ref{A:res:A.6}, $\mathbb W(M_B>1+t)\le C\mathrm e^{-ct^2}$, so, for $q>0$,
\[
 \mathbb E(M_B-1)^q
 =q\int_0^\infty t^{q-1}\mathbb W(M_B>1+t)\,\mathrm dt<\infty.
\]
Thus $P_B^2$ and $P_B^3$ have every finite moment. Taking moments in the last two estimates gives the required moment bounds.
\end{proof}

\subsubsection{Spectral continuity and parameter comparison}\label{A:sec:spectral-continuity-and-parameter-comparison}
\phantomsection\label{A:para:local-spectral-counts} Theorem~\ref{A:res:B.3} and Proposition~\ref{A:res:B.4} hold for
every $B\in\Omega_B$. In particular, the resolvent and each ordered eigenvalue are locally Lipschitz in $\beta^{-1/2}$ on this single set of paths. For spectral projections, fix
$B\in\Omega_B$, $\beta_0\in I_\beta$, and $a<b$ with $a,b\notin\sigma(\mathcal H_{\beta_0})$. Write
\[
 P_\beta(a,b)=\mathbf1_{(a,b)}(\mathcal H_\beta),\qquad
 N_\beta(a,b)=\operatorname{rank}P_\beta(a,b).
\]

\par\Needspace{6\baselineskip}
\begin{theorem}\label{A:res:continuity}
There is a relative neighborhood $U$ of $\beta_0$ in $I_\beta$ such that, for $\beta,\gamma\in U$,
\[
 \|P_\beta(a,b)-P_\gamma(a,b)\|_{\mathfrak S_2}
 \le C_{B,I_\beta,a,b,\beta_0}|c_\beta-c_\gamma|,
 \qquad N_\beta(a,b)=N_{\beta_0}(a,b).
\]
\end{theorem}

\begin{proof}
Use the common shift of Proposition~\ref{A:res:B.1} and abbreviate $S_\theta=S_\theta^\kappa$. We estimate the inverse on a fixed contour in $\rho(\mathcal H_{\beta_0})$. Let
$\Gamma$ be the positively oriented rectangle with vertices $a\pm\mathrm i,b\pm\mathrm i$. Self-adjointness and the conditions $a,b\in\rho(\mathcal H_{\beta_0})$ imply
$\Gamma\subset\rho(\mathcal H_{\beta_0})$. Set
\[
 F_\theta(z)=I-(z+\kappa_B)S_\theta,\qquad
 Q_{\beta_0}(z)=F_{\beta_0}(z)^{-1}
              =I+(z+\kappa_B)R_{\beta_0}(z),
\]
and let $M=\sup_{z\in\Gamma}\|Q_{\beta_0}(z)\|<\infty$ and $Z_\Gamma=\sup_{z\in\Gamma}|z+\kappa_B|$. Choose $U$ so small that $MZ_\Gamma\|S_\theta-S_{\beta_0}\|<1/2$ on $U$; this
is possible by Proposition~\ref{A:res:B.2}. The factorization
\[
 F_\theta(z)=F_{\beta_0}(z)
 [I-(z+\kappa_B)Q_{\beta_0}(z)(S_\theta-S_{\beta_0})]
\]
and the Neumann series show that $Q_\theta(z)=F_\theta(z)^{-1}$ exists and has norm at most $2M$ on $U\times\Gamma$. On the operator domain,
$\mathcal H_\theta-z=F_\theta(z)(\mathcal H_\theta+\kappa_B)$. Consequently $\Gamma\subset\rho(\mathcal H_\theta)$ and $R_\theta(z)=S_\theta Q_\theta(z)$.

The multiplication calculation in Theorem~\ref{A:res:B.3} gives
\[
 R_\beta(z)-R_\gamma(z)
 =Q_\beta(z)(S_\beta-S_\gamma)Q_\gamma(z),
 \qquad
 \sup_{z\in\Gamma}\|R_\beta(z)-R_\gamma(z)\|_{\mathfrak S_2}
 \le4M^2C_S|c_\beta-c_\gamma|,
\]
where $C_S$ is the constant of Proposition~\ref{A:res:B.2}. The spectrum inside $\Gamma$ is finite. Integration of the scalar resolvent in its spectral expansion gives
\[
 P_\theta(a,b)=-\frac1{2\pi\mathrm i}
                    \int_\Gamma R_\theta(z)\,\mathrm dz.
\]
The sign is due to the convention $R_\theta(z)=(\mathcal H_\theta-z)^{-1}$. Integration of the difference in $\mathfrak S_2$ gives the stated bound with constant
$4M^2C_S|\Gamma|/(2\pi)$.

Shrink $U$ until $\|P_\theta-P_{\beta_0}\|_{\mathrm{op}}<1$. For $v\in\operatorname{Ran}P_{\beta_0}$ with $P_\theta
v=0$,
\[
 \|v\|=\|(P_{\beta_0}-P_\theta)v\|
 \le\|P_{\beta_0}-P_\theta\|_{\mathrm{op}}\|v\|<\|v\|
\]
unless $v=0$. Thus $P_\theta$ is injective on the range of $P_{\beta_0}$. Reversing the projections gives the opposite rank inequality. Their finite ranks are equal.
\end{proof}

The proof applies to any fixed compact contour separating a finite spectral cluster. Suppose its projection has rank one at $\beta_0$, and choose a unit vector $\psi_0$ in that
range. Shrink the relative neighborhood $U$ so that
\[
 \|P_\beta-P_{\beta_0}\|_{\mathrm{op}}\le1/2\qquad(\beta\in U).
\]
The rank remains one there. Set $u_\beta=P_\beta\psi_0$. Then
\[
 \|u_\beta\|\ge\|P_{\beta_0}\psi_0\|
      -\|(P_\beta-P_{\beta_0})\psi_0\|\ge1/2.
\]
Thus $\psi_\beta=u_\beta/\|u_\beta\|$ is well defined and spans the eigenspace. It has a consistent phase, since
\[
 \langle\psi_0,\psi_\beta\rangle
 =\frac{\langle\psi_0,P_\beta\psi_0\rangle}{\|u_\beta\|}
 =\|u_\beta\|>0.
\]
For nonzero vectors $u,v$, adding and subtracting $v/\|u\|$ and using $|\|u\|-\|v\||\le\|u-v\|$ gives
\[
 \left\|\frac{u}{\|u\|}-\frac{v}{\|v\|}\right\|
 \le\frac{2\|u-v\|}{\min\{\|u\|,\|v\|\}}.
\]
Apply this to $u_\beta,u_\gamma$ and use the lower bound $1/2$. The projection estimate yields
\[
 \|\psi_\beta-\psi_\gamma\|
 \le4\|(P_\beta-P_\gamma)\psi_0\|
 \le4\|P_\beta-P_\gamma\|_{\mathfrak S_2}
 \le C_{B,I_\beta,\Gamma}|c_\beta-c_\gamma|.
\]
This proves local Lipschitz continuity of the eigenfunction with the stated normalization. If
\[
 d_B=\inf_{\theta\in I_\beta}
       \operatorname{dist}(\{a,b\},\sigma(\mathcal H_\theta))>0,
\]
put $d_* =\min\{1,d_B\}$. For $z$ on either vertical side of the rectangle $\Gamma$, self-adjointness gives
\[
 \operatorname{dist}(z,\sigma(\mathcal H_\theta))^2
 =\operatorname{dist}(\operatorname{Re}z,\sigma(\mathcal H_\theta))^2
       +(\operatorname{Im}z)^2\ge d_B^2.
\]
On its horizontal sides the imaginary part has absolute value one. Thus the distance is at least $d_*$ throughout $I_\beta\times\Gamma$. The identity
\[
 Q_\theta(z)=[I-(z+\kappa_B)S_\theta^\kappa]^{-1}
           =I+(z+\kappa_B)R_\theta(z)
\]
gives, with $Z_\Gamma=\sup_\Gamma|z+\kappa_B|$,
\[
 \sup_{\theta\in I_\beta,z\in\Gamma}\|Q_\theta(z)\|
             \le1+Z_\Gamma/d_*=:M_*.
\]
Let $C_S$ be the constant of Proposition~\ref{A:res:B.2}. The ordered inverse identity from Theorem~\ref{A:res:B.3} and Proposition~\ref{A:res:B.2} imply
\[
 \sup_{z\in\Gamma}\|R_\beta(z)-R_\gamma(z)\|_{\mathfrak S_2}
 \le M_*^2 C_S|c_\beta-c_\gamma|.
\]
Integrating around $\Gamma$ proves
\[
 \|P_\beta(a,b)-P_\gamma(a,b)\|_{\mathfrak S_2}
 \le\frac{|\Gamma|}{2\pi}M_*^2C_S|c_\beta-c_\gamma|.
\]
For parameters sufficiently close to one another this bound is less than one in operator norm. Each projection is then injective on the range of the other, so their finite ranks
agree. The integer-valued rank function is consequently locally constant. Its level sets are relatively open; connectedness of $I_\beta$ makes the rank constant on the entire
interval.

The constants contain $d_*^{-1}$. Proposition~\ref{A:res:B.5} controls the moments of constants obtained from the form estimate for the Brownian potential, but it gives no moment
bound for $d_*^{-1}$. Also, atomlessness of an eigenvalue at each fixed parameter does not establish the displayed strictly positive infimum over a continuum of parameters. The
uniform positive lower bound for the distance to the spectrum remains a separate pathwise hypothesis. We now turn to the finite matrices. Use the laws $\Pi_{n,\beta}$ and
$\overline\Pi_{n,\beta}$ of Lemma~\ref{A:res:B.8}; their projections onto $(B^+,x)$ agree. Keeping the auxiliary variables in the latter law allows us to estimate each comparison
error without altering the matrix marginal. All suprema over $\beta$ below remain outside the corresponding expectation or $L^p$ norm.

\subsubsection{Uniform moment estimates}\label{A:sec:uniform-moment-estimates}
To obtain moment bounds, we estimate the error at each fixed size and retain the dependence of the resolvent constants on the Brownian path. On coefficient events with polynomially
small complements, these constants are bounded by a fixed polynomial in the Brownian oscillation. Off these events the original matrices remain self-adjoint, so the nonreal
resolvent bound controls the error without any assumptions on their coefficients.

Use the integer $s_*$ from \eqref{A:eq:B.45} for the number of final coordinates removed, and put
\[
 I_n^*=\{r_n+1,\ldots,n-s_*\},\quad M_n^*=|I_n^*|,\quad
 \mathcal J_{n,I_\beta}
   =\iota_{h_{n,V}}^{(M_n^*)}U_{r_n}P_{I_n^*}.
\]
All general model comparisons use this interval. We work first at sizes for which it is nonempty and the cutoffs of
Subsection~\ref{A:sec:gaussian-approximation-of-the-conditional-gibbs-measure} are ordered. If $n-s_*\le r_n$, take $I_n^*=\varnothing$, $M_n^*=0$, and the compressed embedding to
be zero. The full coefficient law at each of these finitely many sizes is still the one specified in Lemma~\ref{A:res:B.8}. Write $X_j=1+x_j$, $\ell_n=1+\log(n+1)$, and
$K^*=K\cup\{\overline z:z\in K\}$. All bounds for adjoints are made on $K^*$. Expectations involving auxiliary coordinates are under $\overline\Pi_{n,\beta}$.

We track how the constants in Subsection~\ref{A:sec:weighted-resolvent-estimates-and-comparison-principles} depend on the shift used to make the quadratic form positive.
Exponential dependence on this shift would require stronger integrability of the shift. We obtain polynomial dependence by decreasing the kernel decay rate as the shift grows and
keeping track of that rate in the weighted sums. Let $H$ be a real symmetric tridiagonal matrix on a Dirichlet interval of $N$ vertices with inherited nonnegative coordinates.
Suppose $0<h\le1$, $Nh\le C_*h^{-2}$, $X_j\le C_*h^{-2}$, and

\[
 H+\kappa I\ge c_*\mathcal L_h,\qquad
 h^2|H(j,j+1)|\le C_*,\qquad \kappa\ge1,
 \eqtag{analytic-input}
\]
where $c_*,C_*>0$ are deterministic and $\mathcal L_h$ is the matrix of the form $\ell_h$. Put $P=1+\kappa$ and $R(z)=(H-z)^{-1}$. The hypotheses do not bound the diagonal entries
in absolute value: their contribution cancels under the reciprocal exponential weights. At the left boundary of the original finite interval, the estimate below uses the number of
indices between that boundary and the evaluation point.

\begin{lemma}\label{A:res:analytic-constants}
There is $a=c_{c_*,C_*,K}/P>0$ for which, uniformly on $K$,
\[
 \begin{aligned}
\|\mathcal L_h^{1/2}R(z)\mathcal L_h^{1/2}\|_{\mathrm{op}}& \le CP,\qquad |R(z;i,j)| \le CPh(X_iX_j)^{-1/4} \mathrm e^{-a|\Phi(x_i)-\Phi(x_j)|},\\
\sum_i\mathrm e^{a|\Phi(x_i)-\Phi(x_j)|}|R(z;i,j)|^2& \le CP^4hX_j^{-3/2}.
\end{aligned}
 \eqtag{analytic-columns}
\]
The squared Hilbert--Schmidt perturbation, derivative, and boundary estimates of Subsection~\ref{A:sec:weighted-resolvent-estimates-and-comparison-principles} have constants at
most $CP^{12}$.
\end{lemma}

For a pair of matrices take the larger shift. We shall also use the conclusion for $R^{\mathrm c}(z)=(K_h^{\mathrm c}-zM_h)^{-1}$ from Proposition~\ref{A:res:4.5A}, assuming
$K_h^{\mathrm c}+\kappa M_h\ge c_*\mathcal L_h$ and the deterministic bound $h^2|(K_h^{\mathrm c})_{j,j+1}|\le C_*$. Here $M_h$ is the tridiagonal mass matrix of
Proposition~\ref{A:res:4.5A}. For the continuum forms with integrated potential $c_\beta B$ or its step approximation, the kernel, multiplication, and first-order estimates hold
with integrals instead of sums and without the factor $h$. The continuum estimate does not include a Hilbert--Schmidt bound for the differentiated half-line resolvent.

\begin{proof}
Let $A=H+\kappa I$. The lower form bound implies $\|\mathcal L_h^{1/2}A^{-1/2}\|\le c_*^{-1/2}$. Spectral calculus gives
\[
 \|A^{1/2}R(z)A^{1/2}\|
 =\sup_{\lambda\in\sigma(H)}\frac{\lambda+\kappa}{|\lambda-z|}
 \le1+\frac{\kappa+|z|}{\delta_K}\le C_KP.
\]
Factoring through $A^{-1/2}$ proves the first bound in \eqref{A:eq:analytic-columns}. Equivalently, the form inverse has norm at most $CP$ from the dual of
$(\mathbb C^N,\ell_h^{1/2})$ to that space.

Fix a vertex $i$ and a truncation threshold $R_0>0$. Set
\[
 \psi_j=\min\{|\Phi(x_j)-\Phi(x_i)|,R_0\},\quad
 \Theta=\operatorname{diag}(\mathrm e^{a\psi_j}),\quad
 d_j=a(\psi_{j+1}-\psi_j).
\]
The condition $X_j\le C_*h^{-2}$ gives $|d_j|\le Cah\sqrt{X_j}\le Ca$, also at the right boundary of the increasing interval. Diagonal entries cancel in the form of
$\Theta(H-z)\Theta^{-1}-(H-z)$. For $t_j=H(j,j+1)$ the remaining contribution is
\[
 t_j\{(\mathrm e^{d_j}-1)u_j\overline v_{j+1}
          +(\mathrm e^{-d_j}-1)u_{j+1}\overline v_j\}.
\]
Its term linear in $d_j$ contains
\[
 u_j\overline v_{j+1}-u_{j+1}\overline v_j
   =h\{u_j\overline{D_hv_j}-D_hu_j\overline v_j\}.
\]
Since $|t_j|\le C_*h^{-2}$, the sum of these linear terms is bounded in absolute value by
\[
 Ca\sum_j\sqrt{X_j}
       (|u_j||D_hv_j|+|D_hu_j||v_j|)
 \le Ca\ell_h[u]^{1/2}\ell_h[v]^{1/2}.
\]
The remainders in $\mathrm e^{\pm d_j}-1\mp d_j$ have size at most $Cd_j^2$. Their total contribution is at most
\[
 Ca^2\sum_jX_j(|u_j||v_{j+1}|+|u_{j+1}||v_j|)
 \le Ca^2\ell_h[u]^{1/2}\ell_h[v]^{1/2}.
\]
The last inequality follows from Cauchy--Schwarz and $X_{j+1}\asymp X_j$. The unperturbed form inverse has norm $CP$, whereas the conjugation error has norm $C(a+a^2)$ between the
energy space and its dual. Choose $a_0=c/P\le1$ with $c$ small enough that $CP\cdot C(a_0+a_0^2)\le1/2$. The Neumann series in the form space yields
\[
 \|\mathcal L_h^{1/2}\Theta R(z)\Theta^{-1}
                          \mathcal L_h^{1/2}\|\le CP,
 \eqtag{analytic-conjugation}
\]
uniformly in $R_0$, for $0<a\le a_0$.

Lemma~\ref{A:res:3.6} gives $\|\mathcal L_h^{-1/2}e_j\|^2\le ChX_j^{-1/2}$. Apply \eqref{A:eq:analytic-conjugation} to these two vectors at $i$ and $j$. Since $\Theta_i=1$, this
proves
\[
 |R(z;j,i)|\le CPh(X_iX_j)^{-1/4}
             \mathrm e^{-a_0\min\{|\Phi(x_j)-\Phi(x_i)|,R_0\}}.
\]
Let $R_0\to\infty$. This proves the pointwise bound.

We give the summation constants explicitly. For $0<a\le1$ and $q\in\{1/2,3/2\}$,
\[
 \sum_jX_j^{-q}\mathrm e^{-a|\Phi(x_j)-\Phi(x_i)|}
 \le Ch^{-1}X_i^{-q-1/2}a^{-b_q},
 \qquad b_{1/2}=2,\quad b_{3/2}=3.
 \eqtag{analytic-sums}
\]
To prove this, integrate on the disjoint intervals of length $h$ centered at the $x_j$. On each such interval $X(x)\asymp X_j$ and $|\Phi(x)-\Phi(x_j)|\le Ch\sqrt{X_j}\le C$. The
summand is therefore bounded by $Ch^{-1}$ times the integral of $X(x)^{-q}\mathrm e^{-a|\Phi(x)-\Phi(x_i)|}$ over that interval. The union is contained in $\mathbb R_+$, also for
an inherited coordinate. Extend the resulting integral to $\mathbb R_+$ and put $y=\Phi(x)$. Since $\mathrm dy=\sqrt{X(x)}\,\mathrm dx$ and $X(x)=(1+3y/2)^{2/3}$,
\[
 \left(\frac{X_i}{X(x)}\right)^{q+1/2}
 \le C_q(1+|y-\Phi(x_i)|)^{(2q+1)/3}.
\]
For $r\ge0$, substitution $u=at$ gives
\[
 \int_{\mathbb R}(1+|t|)^r\mathrm e^{-a|t|}\,\mathrm dt
 \le C_r a^{-1-r}.
\]
Here $r=(2q+1)/3$, so the powers of $a^{-1}$ are $5/3$ and $7/3$. Since $a\le1$, we may bound them by $a^{-2}$ and $a^{-3}$, respectively. This proves \eqref{A:eq:analytic-sums}
with constants independent of the number of vertices.

Square the kernel bound and multiply by $\mathrm e^{a_0|\Phi(x_j)-\Phi(x_i)|}$. The case $q=1/2$ gives
\[
 \sum_j\mathrm e^{a_0|\Phi(x_j)-\Phi(x_i)|}|R(z;j,i)|^2
 \le CP^2h^2X_i^{-1/2}
       (h^{-1}X_i^{-1}a_0^{-2})
 \le CP^4hX_i^{-3/2}.
\]
Decreasing $a_0$ once and calling it $a$ proves all of \eqref{A:eq:analytic-columns} with the same parameter. Cauchy--Schwarz for the weighted columns, using
$|\Phi_i-\Phi_j|\le|\Phi_i-\Phi_k|+|\Phi_k-\Phi_j|$, gives
\[
 |[R(\overline z)R(z)](i,j)|
 \le CP^4h(X_iX_j)^{-3/4}\mathrm e^{-(a/2)|\Phi_i-\Phi_j|}.
\]
Since both resolvents are functions of $H$, $R(z)R(\overline z)=R(\overline z)R(z)$. For a complex diagonal multiplier $v$, the exact trace expansion is
\[
 \|R_1M_vR_2\|_{\mathfrak S_2}^2
 =\sum_{i,j}\overline v_i v_j
       [R_1^*R_1](i,j)[R_2R_2^*](j,i).
\]
Taking absolute values, using $2|v_iv_j|\le|v_i|^2+|v_j|^2$, and then \eqref{A:eq:analytic-sums} at $q=3/2$ proves
\[
 \|R_1M_vR_2\|_{\mathfrak S_2}^2 \le CP^8h^2\sum_i |v_i|^2X_i^{-3/2} (h^{-1}X_i^{-2}a^{-3}) \le
 CP^{11}h\sum_i|v_i|^2X_i^{-7/2}.
 \eqtag{analytic-sandwich}
\]
Summing the unweighted column bound gives $\|M_vR\|_{\mathfrak S_2}^2+\|RM_v\|_{\mathfrak S_2}^2 \le CP^4h\sum_jX_j^{-3/2}|v_j|^2$. Extend vertex rows by zero at both points where
the Dirichlet condition is imposed. For $J=J_-,J_+,J_h$, the inequality $|(Ju)_j|^2\le|u_j|^2+|u_{j+1}|^2$ gives
\[
 \|M_eJR(z)\|_{\mathfrak S_2}^2
 \le\sum_j|e_j|^2\sum_k
       \bigl(|R(z;j,k)|^2+|R(z;j+1,k)|^2\bigr)
 \le CP^4h\sum_jX_j^{-3/2}|e_j|^2.
\]
The last inequality uses the row bound at $j,j+1$ and $X_{j+1}^{-3/2}\le C X_j^{-3/2}$; a row beyond the finite matrix has value zero by the boundary convention.

The first estimate in \eqref{A:eq:analytic-columns}, together with $\mathcal L_h\ge I$, gives $\|D_hR(z)\|_{\mathrm{op}}^2\le CP^2$. Consequently the two terms in
$R_1(D_h^*M_eJ+J^*M_eD_h)R_2$ have total squared Hilbert--Schmidt norm at most
\[
 CP^6\sum_j hX_j^{-3/2}|e_j|^2.
\]
Indeed, the two squared norms are bounded by
\[
 \|D_hR_1(\overline z)\|_{\mathrm{op}}^2
       \|M_eJR_2(z)\|_{\mathfrak S_2}^2,
 \qquad
 \|M_eJR_1(\overline z)\|_{\mathfrak S_2}^2
       \|D_hR_2(z)\|_{\mathrm{op}}^2,
\]
respectively. Each is at most $CP^6\sum_jhX_j^{-3/2}|e_j|^2$. At the left Dirichlet boundary, $u_j=h\sum_{k<j}D_hu_k$, so
\[
 \|R(z)e_j\|^2\le CP^2h^2j.
\]
A matrix $W$ of fixed bandwidth supported on the first $L_0$ sites has at most $C(L_0+1)$ nonzero entries. Expansion into matrix units and the last bound give
\[
 \|R_1WR_2\|_{\mathfrak S_2}^2
 \le CP^4h^4(L_0+1)^4\max_{i,j}|W_{ij}|^2.
\]
Both Dirichlet edges are included in $D_h$.

If $t_k=\lambda_k(H)+\kappa$, min--max and Lemma~\ref{A:res:3.4} give $t_k\ge ck^{2/3}$. For normalized eigenvectors $\psi_k$, $\|D_h\psi_k\|^2\le Ct_k$, and
$t_k/|\lambda_k-z|\le CP$. It follows that
\[
 \|R(z)\|_{\mathfrak S_2}^2\le CP^2\sum_{k\ge1}k^{-4/3},
 \qquad
 \|D_hR(z)\|_{\mathfrak S_2}^2
 \le CP^2\sum_{k\le N}k^{-2/3}\le CP^2h^{-1}.
\]
The factors depending on $P$ in the displayed estimates are, respectively, $P^2a^{-2}$, $P^8a^{-3}$, $P^2P^4$, $P^4$, and $P^2$. Since $a=c/P$ and $P\ge1$,
\[
 P^2a^{-2}=c^{-2}P^4,\qquad
 P^8a^{-3}=c^{-3}P^{11},\qquad
 \max\{P^{11},P^6,P^4,P^2\}\le P^{12}.
\]
The constants $c^{-2},c^{-3}$ depend only on the fixed form bounds and $K$. This proves the common polynomial bound.

For the matrices of the restricted forms we work directly with $K_h^{\mathrm c}-zM_h$, where $\frac13I\le M_h\le I$ and $K_h^{\mathrm c}+\kappa M_h\ge c_*\mathcal L_h$. Conjugation
by $M_h^{-1/2}$ is not used in this argument, since it need not preserve tridiagonality. Put $\|u\|_{\mathcal E_h}=\ell_h[u]^{1/2}$. If
\[
 (K_h^{\mathrm c}-zM_h)u=F,
\]
where $F$ is a continuous functional on the energy space, the imaginary part and then the real part of the equation tested against $u$ give
\[
 \delta_K\langle u,M_hu\rangle\le |F(u)|,\qquad
 c_*\ell_h[u]\le\operatorname{Re}F(u)
                   +(\kappa+|z|)\langle u,M_hu\rangle.
\]
Consequently
\[
 c_*\|u\|_{\mathcal E_h}^2
 \le\left(1+\frac{\kappa+|z|}{\delta_K}\right)
                      \|F\|_{\mathcal E_h^*}\|u\|_{\mathcal E_h},
 \qquad
 \|u\|_{\mathcal E_h}\le CP\|F\|_{\mathcal E_h^*}.
\]
For $u=0$ the estimate holds; for $u\ne0$, divide by its positive energy norm. The imaginary part also proves injectivity of the pencil; finite dimensionality then proves
invertibility. Thus its inverse has norm at most $CP$ from the dual energy space to the energy space.

Use the same diagonal weight $\Theta$ as above. The pencil has symmetric off-diagonal entries
\[
 t_j=(K_h^{\mathrm c})_{j,j+1}-z(M_h)_{j,j+1},\qquad |t_j|\le C_Kh^{-2}.
\]
Although these numbers may be complex, the entries in the two transposed positions coincide. Their contribution to the difference of the conjugated and original forms is
\[
 t_j\{(\mathrm e^{d_j}-1)u_j\overline v_{j+1}
           +(\mathrm e^{-d_j}-1)u_{j+1}\overline v_j\}.
\]
The diagonal terms cancel. The linear part uses
\[
 u_j\overline v_{j+1}-u_{j+1}\overline v_j
      =h(u_j\overline{D_hv_j}-D_hu_j\overline v_j),
\]
and has absolute sum at most $Ca\|u\|_{\mathcal E_h}\|v\|_{\mathcal E_h}$. The two exponential remainders are bounded by $Cd_j^2$, and their sum is at most
$Ca^2\|u\|_{\mathcal E_h}\|v\|_{\mathcal E_h}$. Choose $a=c/P$ so that the product of this perturbation norm and the inverse norm $CP$ is at most $1/2$. The Neumann series on the
energy space gives
\[
 \|\mathcal L_h^{1/2}\Theta R^{\mathrm c}(z)\Theta^{-1}
                            \mathcal L_h^{1/2}\|\le CP.
\]
Pairing with $\mathcal L_h^{-1/2}e_i$ and $\mathcal L_h^{-1/2}e_j$, and then removing the truncation of the weight, proves the same kernel bound as for the ordinary resolvent. The
summation estimate \eqref{A:eq:analytic-sums} gives the same weighted row and column bounds.

The inverse $R^{\mathrm c}(z)$ need not be normal. Its adjoint is $R^{\mathrm c}(\overline z)$, but $R^{\mathrm c}(z)$ need not commute with that adjoint. Instead, apply
Cauchy--Schwarz separately to
\[
 [(R^{\mathrm c})^*R^{\mathrm c}](i,j)
 =\sum_k\overline{R^{\mathrm c}(k,i)}R^{\mathrm c}(k,j),\qquad
 [R^{\mathrm c}(R^{\mathrm c})^*](i,j)
 =\sum_k R^{\mathrm c}(i,k)\overline{R^{\mathrm c}(j,k)}.
\]
The column estimates bound the first sum and the row estimates bound the second. Both are at most
\[
 CP^4h(X_iX_j)^{-3/4}
       \mathrm e^{-(a/2)|\Phi(x_i)-\Phi(x_j)|}
\]
in absolute value. Substituting these two distinct bounds into the exact trace expansion for $R_1M_vR_2$ proves the diagonal perturbation estimate. The same row bounds prove the
multiplier estimates for $M_eJR^{\mathrm c}$. This proves the stated kernel and perturbation estimates also when one factor is an ordinary matrix resolvent and the other is the
inverse of a matrix pencil with a nonidentity mass matrix; no commutation between them is needed. The derivative norm follows from the eigenvalue expansion alone. We may therefore
put $\widetilde H=M_h^{-1/2}K_h^{\mathrm c}M_h^{-1/2}$. If $\phi_k$ is its orthonormal eigenbasis and $u_k=M_h^{-1/2}\phi_k$, then
\[
 \langle u_k,M_hu_k\rangle=1,\quad
 \ell_h[u_k]\le C(\lambda_k+\kappa),\quad
 \lambda_k+\kappa\ge ck^{2/3}.
\]
The last bound follows from the generalized min--max principle and $M_h\le I$. Therefore
\[
 \|D_hR^{\mathrm c}(z)\|_{\mathfrak S_2}^2 \le\|M_h^{-1/2}\|^2 \sum_{k\le
 N}\frac{\|D_hu_k\|^2}{|\lambda_k-z|^2} \le CP^2\sum_{k\le N}k^{-2/3}\le CP^2h^{-1}.
\]

For the continuum, point evaluation has squared dual form norm at most $CX(x)^{-1/2}$. To see this, put $\ell=X(x)^{-1/2}$ and choose an interval $I_x\subset[0,\infty)$ containing
$x$, with length between $\ell$ and $2\ell$. The fundamental theorem of calculus, averaged over the other endpoint, gives
\[
 |u(x)|^2\le C\left(\ell^{-1}\int_{I_x}|u|^2+\ell\int_{I_x}|u'|^2\right).
\]
On this interval $X(t)$ is comparable to $X(x)$. Since $\ell^{-1}X(x)^{-1}=\ell$, the right side is at most $CX(x)^{-1/2}\|u\|_{\mathcal E}^2$. We now prove the corresponding
conjugation estimate on the form domain. Fix $y$ and $R_0>0$, put
\[
 \psi(x)=\min\{|\Phi(x)-\Phi(y)|,R_0\},\qquad \theta(x)=\mathrm e^{a\psi(x)}.
\]
This function is Lipschitz on bounded intervals, $|\psi'|\le\sqrt X$ almost everywhere, and $\theta,\theta^{-1}$ are bounded for fixed $R_0$. The product rule and
$\int X|u|^2<\infty$ show that multiplication by either function preserves the form domain and its Dirichlet trace.

In the form evaluated at $(\theta^{-1}u,\theta v)$, both the ordinary potential and the integrated random potential are unchanged: the product of these arguments is $u\overline v$,
and its weak derivative is unchanged as well. Expanding the kinetic term gives the difference
\[
 a\int\psi'u'\overline v
   -a\int\psi'u\overline{v'}
   -a^2\int(\psi')^2u\overline v.
\]
For example, the first term has absolute value at most
\[
 a\|u'\|_2\|\sqrt X\,v\|_2
 \le a\|u\|_{\mathcal E}\|v\|_{\mathcal E};
\]
the second has the same bound, and the last is bounded by $a^2\|\sqrt X\,u\|_2\|\sqrt X\,v\|_2$. Thus the perturbation norm from the energy space to its dual is at most $C(a+a^2)$.
The unperturbed form inverse has norm $CP$. Taking $a=c/P$ with $c$ small, a Neumann series bounds the conjugated form inverse by $CP$, uniformly in $R_0$.

Pairing this inverse with the point-evaluation functionals at $x,y$ gives
\[
 |R(z;x,y)|\le CP X(x)^{-1/4}X(y)^{-1/4}
      \mathrm e^{-a\min\{|\Phi(x)-\Phi(y)|,R_0\}}.
\]
Letting $R_0\to\infty$ proves the kernel estimate. For $q\in\{1/2,3/2\}$, use the normalization $\Phi(0)=0$, so $\Phi(x)=\frac23(X(x)^{3/2}-1)$. Put $u=\Phi(x)$, $v=\Phi(y)$ and
$\alpha=2(q+1/2)/3$. Then
\[
 X(x)^{-q}\,\mathrm dx=(1+3u/2)^{-\alpha}\,\mathrm du,
 \qquad
 (1+3u/2)^{-\alpha}\le C(1+3v/2)^{-\alpha}(1+|u-v|)^\alpha.
\]
For $0<a\le1$, integration over $u\ge0$ is bounded by the integral over the whole line after this estimate, and
\[
 \int_{\mathbb R}(1+|s|)^\alpha\mathrm e^{-a|s|}\,\mathrm ds
 \le C_\alpha a^{-1-\alpha}.
\]
Here $1+\alpha$ equals $5/3$ or $7/3$, respectively. Increasing these exponents to $2$ or $3$ gives
\[
 \int_0^\infty X(x)^{-q}\mathrm e^{-a|\Phi(x)-\Phi(y)|}\,\mathrm dx
 \le C X(y)^{-q-1/2}a^{-b_q},\qquad
 b_{1/2}=2,\quad b_{3/2}=3.
\]
Squaring the kernel and using $q=1/2$ therefore gives the weighted column bound
\[
 \int_0^\infty\mathrm e^{a|\Phi(x)-\Phi(y)|}|R(z;x,y)|^2\,\mathrm dx
 \le CP^2X(y)^{-1/2}X(y)^{-1}a^{-2}
 \le CP^4X(y)^{-3/2}.
\]
The adjoint resolvent gives the same row bound. Cauchy--Schwarz with these exponential weights bounds the kernels of $R(z)^*R(z)$ and $R(z)R(z)^*$ by
\[
 CP^4(X(x)X(y))^{-3/4}
       \mathrm e^{-(a/2)|\Phi(x)-\Phi(y)|}.
\]
For a bounded compactly supported multiplier $v$, the integral trace expansion, followed by $2|v(x)v(y)|\le|v(x)|^2+|v(y)|^2$, yields
\[
 \begin{split}
 \|R_1M_vR_2\|_{\mathfrak S_2}^2
 &\le CP^8\iint |v(x)v(y)|(X(x)X(y))^{-3/2}
                    \mathrm e^{-a|\Phi(x)-\Phi(y)|}\,\mathrm dx\mathrm dy\\
 &\le CP^8a^{-3}\int |v(x)|^2X(x)^{-7/2}\,\mathrm dx
 \le CP^{11}\int |v(x)|^2X(x)^{-7/2}\,\mathrm dx.
 \end{split}
\]
Approximation in the displayed weighted $L^2$ norm extends the resolvent product to the full multiplier class in this estimate.

The row bound likewise gives
\[
 \|M_eR(z)\|_{\mathfrak S_2}^2
       \le CP^4\int |e(x)|^2X(x)^{-3/2}\,\mathrm dx.
\]
The energy inverse estimate gives $\|\partial_xR(z)\|_{\mathrm{op}}^2\le CP^2$. For the form $\Gamma_e[f,g]=-\int e(f\overline g)'$, write $D=\partial_x$ on the form domain. For
bounded compactly supported $e$, the weak expression for its resolvent product is the sum
\[
 -R_1M_eDR_2-R_1D^*M_eR_2.
\]
The second term means the bounded extension $(DR_1^*)^*M_eR_2$. The first is $(M_{\overline e}R_1^*)^*DR_2$. Thus each is a product of one bounded differentiated resolvent and one
Hilbert--Schmidt multiplier product. The ideal inequality and the squared triangle inequality give
\[
 \|R_1\Gamma_e R_2\|_{\mathfrak S_2}^2
 \le CP^6\int |e(x)|^2X(x)^{-3/2}\,\mathrm dx.
\]
Initially this computation is made for bounded compactly supported $e$ on vectors in the form domain; the estimate gives the stated extension by weighted $L^2$ approximation. It
requires only the operator norm of $\partial_xR$, not a Hilbert--Schmidt norm of the differentiated half-line resolvent. For a finite Dirichlet interval, zero extension preserves
the energy and both boundary values, so the constants are independent of its right boundary. We finally record the bound needed for intervalwise regularity.

For $u=\Phi(x)$ and $v=\Phi(y)$, the row and column integrals of $X(x)|R(z;x,y)|$ are bounded by
\[
 CP\int_{\mathbb R}(1+|t|)^{1/2}\mathrm e^{-a|t|}\,\mathrm dt,
 \qquad
 CP\int_{\mathbb R}(1+|t|)^{1/6}\mathrm e^{-a|t|}\,\mathrm dt,
\]
respectively. Both are at most $CPa^{-2}\le CP^3$. Schur's test proves $\|XR(z)\|_{\mathrm{op}}\le CP^3$, including for the step approximation to $B$.
\end{proof}

We integrate the preceding bounds on the coefficient event and on its complement under the original matrix law. The following estimate allows dependence between the random
resolvent constant and the coefficient error. Let $D_{n,\beta}\ge0$ and suppose $\sup_\beta\mathbb P(E_{n,\beta,A}^c)\le C_An^{-A}$. Assume $D\le GW$ on $E=E_{n,\beta,A}$, where
\[
 \sup_{n,\beta}\|\mathbf1_EG\|_{2p}\le C_{p,A},\qquad
 \sup_\beta\|\mathbf1_EW\|_{2p}\le C_{p,A}a_n.
\]
Assume also $D\le C_Kn+Z$ everywhere, with $\sup_{n,\beta}\|Z\|_{2p}<\infty$.

\begin{lemma}\label{A:res:exceptional-moments}
Under these assumptions,
\[
 \sup_\beta\|D_{n,\beta}\|_p
 \le C_{p,A}\{a_n+n^{1-A/p}+n^{-A/(2p)}\}.
 \eqtag{exceptional-moment}
\]
\end{lemma}

For a difference of finite resolvents of rank $O(n)$, take $Z=0$. For a comparison with the continuum resolvent one may take $Z=C(1+\sup_{z\in K}\|R_\beta(z)\|_{\mathfrak S_2}^2)$.
In particular, if $a_n=h^s$, $h\asymp n^{-1/3}$, and $s>0$, it suffices that $A>p\max\{1+s/3,2s/3\}$. The same condition suffices for $a_n=h^s\ell_n^r$ with $r\ge0$.

\begin{proof}
Since $D\le GW$ on $E$, H\"older gives $\|\mathbf1_ED\|_p\le\|\mathbf1_EG\|_{2p} \|\mathbf1_EW\|_{2p}\le C_{p,A}a_n$. H\"older's inequality applies to these variables and the event
indicator with their joint distribution. On the complement,
\[
 \mathbb E[D^p;E^c] \le C_p\{n^p\mathbb P(E^c)+\mathbb E[Z^p;E^c]\} \le C_p\{n^p\mathbb
 P(E^c)+\|Z\|_{2p}^p \mathbb P(E^c)^{1/2}\} \le C_{p,A}(n^{p-A}+n^{-A/2}).
\]
Taking $p$th roots and adding the estimate on $E$ proves \eqref{A:eq:exceptional-moment}. A finite self-adjoint resolvent has norm at most $\delta_K^{-1}$ and rank at most its
dimension, so its squared Hilbert--Schmidt norm is at most $N\delta_K^{-2}$, independently of its coefficients. A contraction on either side does not increase this norm.

For the continuum, the deterministic lower coefficient in \eqref{A:eq:continuum-polynomial-form} gives $\Lambda_k(\beta)+\kappa_B\ge ck^{2/3}$. Therefore
\[
 \sup_{\beta\in I_\beta,z\in K}\|R_\beta(B,z)\|_{\mathfrak S_2}^2
 \le C_{I_\beta,K}(1+\mathcal O_{1/4}(B))^8.
 \eqtag{continuum-coarse}
\]
Indeed, each spectral summand is at most $C_K(1+\kappa_B)^2k^{-4/3}$, and the series converges. The right side has every finite Wiener moment. The argument applies also to the
continuum Dirichlet restrictions and to the step approximation. The mass bounds give the finite-rank estimate for inverses of the matrix pencils up to a deterministic factor. These
observations justify the stated choices of $Z$. Comparing the two powers of $n$ in \eqref{A:eq:exceptional-moment} with $n^{-s/3}$ proves the condition on $A$; multiplication of
the target bound by $\ell_n^r\ge1$ does not weaken it.
\end{proof}

The regular cutoffs and logarithmic separation lengths remain fixed on $I_\beta$. The Hessian of the negative logarithm of the conditional density is bounded below by
$c_V\beta_-n/2$, and
\[
 \sup_{\beta\in I_\beta}\frac1{\beta n}\le\frac1{\beta_-n}.
 \eqtag{B.53}
\]
Thus the constants in the coefficient estimates of
Sections~\ref{A:sec:the-general-jacobi-model-deterministic-reduction}--\ref{A:sec:the-general-jacobi-model-gaussian-approximation-and-comparison} are uniform on their specified
bounded sets of exterior values. Extend $M_B$ by the value one on $\Omega_B^c$; this is Borel and changes no Wiener moment. In the next lemma the finite families are the comparison
matrices of Sections~\ref{A:sec:the-gaussian-beta-ensemble} and~\ref{A:sec:the-general-jacobi-model-gaussian-approximation-and-comparison}, their specified convex interpolations,
and principal Dirichlet restrictions with inherited coordinates. The removed boundary blocks of the full general matrix are not included.

\begin{lemma}\label{A:res:B.10A}
For every $A>0$ there are Borel events $E_{n,\beta,A}^{\mathrm{DE}}$ and $E_{n,\beta,A}^V$ such that, for $n\ge N_A$,
\[
 \sup_{\beta\in I_\beta}
 \overline\Pi_{n,\beta}^{\mathrm{model}}
             ((E_{n,\beta,A}^{\mathrm{model}})^c)
 \le C_An^{-A}.
 \eqtag{B.54aa}
\]
On each event the corresponding forms satisfy
\[
 c_A\ell_h[u]\le\mathfrak h[u]+\kappa_{B,A}\|u\|^2
       \le C_A(1+\kappa_{B,A})\ell_h[u],
 \qquad \kappa_{B,A}\le C_AM_B^4,
 \eqtag{B.54ab}
\]
and
\[
 h^2|H(j,j+1)|\le C_A.
 \eqtag{B.54ac}
\]
\end{lemma}

For the fixed coupling, $A$ determines the event and the deterministic starting size $N_A$. Expectations of the Brownian variable below are taken under its original, unconditional
distribution. For every $q\ge1$,
\[
 \sup_{n,\beta}\|M_B\|_{L^q(\overline\Pi_{n,\beta})}<\infty.
 \eqtag{B.54a}
\]
In particular, the common measurable upper bound for the constants in the resolvent estimates may be chosen as
\[
 G_{n,\beta,K,A}
 :=\mathbf1_{E_{n,\beta,A}}C_{A,V,I_\beta,K}(1+M_B)^{48},
 \eqtag{B.54b}
\]
which satisfies
\[
 \sup_{n\ge N_A,\beta\in I_\beta}
 \|G_{n,\beta,K,A}\|_q
  \le C_{q,A,V,I_\beta,K}.
 \eqtag{B.54c}
\]
Lemma~\ref{A:res:B.10B} will treat the boundary error for the full matrix.

\begin{proof}
Let $\mathcal B_0=\{B^+\in\Omega_B\}\cap\mathfrak B_{\mathrm{temp}}$,

where the second event is the two-sided path event in Lemma~\ref{A:res:B.8}. Under each joint law $\overline\Pi_{n,\beta}$, the two path coordinates have Wiener law, so
$\overline\Pi_{n,\beta}(\mathcal B_0)=1$. For the general model impose the following inequalities:
\[
 \begin{gathered}
\max_{I_n^{\mathrm{mid}}}|\zeta-\widehat\zeta^{\mathrm G}| \le n^{-5/6},\qquad \max_{I_n^{\mathrm{mid}}}\frac{h^{-2}|\varepsilon_j^{\mathrm{fr}}|}{X_j} \le C_Ah^{1/14}\sqrt{\ell_n},\\
\|\widetilde e\|_\infty+\|\eta\|_\infty \le C_Ah^{1/2}\sqrt{\ell_n},\qquad \max_{I_n^{\mathrm{mid}}}\frac{h^{-2}|\rho_j|}{X_j} \le C_Ah^{3/2}\sqrt{\ell_n},\qquad \max_j(|\zeta_j|+|\zeta_j^B|) \le C_Ah^{3/2}\sqrt{\ell_n}.
\end{gathered}
\]
Let $\mathcal M_{n,A}$ denote this event and put
\[
 E_{n,\beta,A}^V=\mathcal B_0
 \cap\{\tau\in\mathcal T_n^{\mathrm{out}}(M_0)\}
 \cap\mathcal F_{n,A+3}\cap\mathcal M_{n,A}.
\]
The event $\mathcal F_{n,A+3}$ includes $\mathcal D_{n,A+3}$, on which the two conditional minimizers satisfy Proposition~\ref{A:res:14.1}. The transported fluctuation, centered at
the shorter minimizer, has coordinate $L^q$ error $C_qn^{-1}$. Proposition~\ref{A:res:14.1} adds at most $C_An^{-4}$ on $\mathcal F_{n,A+3}$ at the middle coordinates. Thus a union
bound and Markov, at fixed $q>6(A+4)$, bound failure of the first maximum by $C_qn^{1-q/6}+C_An^{-A-3}+C\mathrm e^{-cn}$. For the other maxima, Lemmas~\ref{A:res:freezing-max} and
\ref{A:res:factor-max} give the required tail bounds at the respective coordinate scales. These use conditional Gaussian tails for the Gaussian errors and concentration under the
exact coefficient law for the Gibbs fluctuations; Gaussianity of the latter is not assumed. Increasing $C_A$ makes each failure probability at most $C_An^{-A-3}+C\mathrm e^{-cn}$.
The exterior tail and \eqref{A:eq:confinement-good-event} bound the remaining failures. Adding this fixed number of bounds proves \eqref{A:eq:B.54aa}. All events are Borel by the
coefficient and kernel constructions of Subsection~\ref{A:sec:measurable-couplings-of-the-finite-models}.

We verify the form estimate on the finite size event $\mathcal F_{n,A+3}$. Lemma~\ref{A:res:12.1B} gives deterministic lower and upper constants for the form at the conditional
minimizer on $\mathcal F_{n,A+3}$. The forms at the local minimizer and at the coefficients of the deterministic difference operator with linear potential have the same type of
bounds. By \eqref{A:eq:16.33}, adding the same diagonal potential obtained from the Brownian increments changes each form by at most
$\varepsilon\ell_h+C_{\varepsilon,I_\beta}M_B^4\|\cdot\|^2$. Choose $\varepsilon=c_*/4$, where $c_*$ is the least of the three deterministic lower constants.

The remaining terms have the following absolute relative form bounds. For the transport, covariance, and factorization multipliers these are, respectively,
\[
 Cn^{-1/6},\qquad C_Ah^{1/14}\sqrt{\ell_n},
 \qquad C_Ah^{3/2}\sqrt{\ell_n}.
\]
The first follows from $h^{-2}n^{-5/6}\le Cn^{-1/6}$; the other two follow by division by $X_j$. On the far interval $X_j\ge ch^{-2}$, so replacing the scalar coefficients by their
Brownian approximations gives a relative bound of at most $C_Ah^{3/2}\sqrt{\ell_n}$. For either first-order term, $\|J\|\le1$ gives
\[
 |\langle u,\Gamma_v u\rangle|
 \le2\|v\|_\infty\|D_hu\|\|u\|
 \le\|v\|_\infty\ell_h[u]
 \le C_Ah^{1/2}\sqrt{\ell_n}\,\ell_h[u].
\]
Summing $|u_j|^2\le jh^2\|D_hu\|^2$ over the indices $j\le J_n^{\mathrm L}$ next to the left boundary gives
\[
 h^{-2}\sum_{j\le J_n^{\mathrm L}}
        |\zeta_j-\zeta_j^B||u_j|^2
 \le C_Ah^{3/2}(J_n^{\mathrm L}+1)^2\sqrt{\ell_n}\,\ell_h[u].
\]
Their sum is bounded by
\[
 a_{n,A}=C_A\bigl[ n^{-1/6} +(h^{1/14}+h^{1/2}+h^{3/2})\sqrt{\ell_n} +h^{3/2}(J_n^{\mathrm
 L}+1)^2\sqrt{\ell_n}\bigr] \longrightarrow0.
 \eqtag{localized-relative}
\]
These are precisely the seven coefficient bounds used to prove Theorem~\ref{A:res:16.2}. At the present fixed size their constants are deterministic on $E_{n,\beta,A}^V$. The
identities \eqref{A:eq:common-assembly} express every comparison form as the quadratic form at the minimizer with the Brownian potential plus a subset of these seven perturbations.
If $\mathfrak r_1,\ldots,\mathfrak r_7$ are these perturbation forms and $t_\nu\in[0,1]$, then
$|\sum_{\nu=1}^7t_\nu\mathfrak r_\nu[u]| \le\sum_{\nu=1}^7|\mathfrak r_\nu[u]|\le a_{n,A}\ell_h[u]$. Choose deterministic $N_A$ so that $a_{n,A}\le c_*/4$. Adding
$\kappa_{B,A}=\kappa_{*,A}+C_{A,I_\beta}M_B^4$ proves \eqref{A:eq:B.54ab}. The localized changes of minimizer used for \eqref{A:eq:13.67}--\eqref{A:eq:13.67a} have, in addition,
relative bound $C_Ah$ in the interior of the regular coefficient interval and $C_Ah^{3/2}\ell_n^{5/2}$ at the indices next to the left boundary. Increase $N_A$ so that
$C_Ah+C_Ah^{3/2}\ell_n^{5/2}\le c_*/8$ and $a_{n,A}\le c_*/8$ for $n\ge N_A$. The sum of all additional relative coefficients is then at most $c_*/4$, so the coefficient of
$\ell_h[u]$ in the lower bound decreases from $3c_*/4$ to at least $c_*/2$. The coefficient in the upper bound increases by at most $c_*/4$.

The kinetic coefficients at the two minimizers are nonnegative and bounded by deterministic constants on the coefficient event. The scalar changes are diagonal; the only one-sided
first-order term has off-diagonal entry $-\eta_j/h$. Thus $h^2|H(j,j+1)|\le C_A+h\|\eta\|_\infty\le C_A$, proving \eqref{A:eq:B.54ac}. Convex combinations preserve these bounds. A
vector in a principal subinterval is extended by zero; the energies of both boundary edges are then retained, proving the restriction statement with inherited coordinates.

For the Gaussian model put
\[
 E_{n,\beta,A}^{\mathrm{DE}}=\mathcal B_0\cap
 \{\max_{j\le n/2}|\Delta\zeta_j|\le h^{11/4},\quad
   \max_{j>n/2}|\Delta\zeta_j|\le h^{5/4},\quad
   \|\eta\|_\infty\le h^{1/4}\}.
\]
The coordinate estimates underlying \eqref{A:eq:8.28d}--\eqref{A:eq:8.28f} give a union bound $C_qn^{1-q/12}$. Choose $q>12(A+2)$. The diagonal relative bound is $Ch^{3/4}$ on both
halves: on the second half use $X_j\ge ch^{-2}$. The first-order relative bound is $Ch^{1/4}$. The deterministic comparisons of Proposition~\ref{A:res:7.4A} are uniform on
$I_\beta$. After the form estimate for the Brownian potential, write the deterministic lower coefficient as $c_0>0$ and the shift as $\kappa_{B,A}\le C_AM_B^4$. Choose a
deterministic size so that $C(h^{3/4}+h^{1/4})\le c_0/2$. The two perturbations then have absolute sum at most $(c_0/2)\ell_h[u]$. Subtracting this sum from the lower form bound
and adding it to the upper form bound proves \eqref{A:eq:B.54ab}. The mean coefficients are at most one and $\delta_j=h\eta_j$, proving \eqref{A:eq:B.54ac}. For the matrices
representing the discretized forms, restrict the sampled continuum form to the piecewise linear space. This gives the lower form bound. Its integrated Brownian potential is
constant between nodes, so its distributional derivative is supported at the nodes. Since only one hat function is nonzero at each interior node, the contribution of the sampled
random potential to the matrix of the restricted form is diagonal. The off-diagonal entries therefore come only from the deterministic kinetic and potential terms and satisfy the
required bound.

Finally, $B^+$ has Wiener law under every $\overline\Pi_{n,\beta}$. The tail of $\mathcal O_{1/4}(B)$ proves \eqref{A:eq:B.54a}. Lemma~\ref{A:res:analytic-constants} and
$1+\kappa_{B,A}\le C_A(1+M_B)^4$ prove \eqref{A:eq:B.54b}; raising that bound to the $q$th power and integrating proves \eqref{A:eq:B.54c}.
\end{proof}

We shall repeatedly use the Gaussian quadratic moment estimate already proved in \eqref{A:eq:gaussian-weighted-quadratic-moment}. If, conditionally on $\tau$, $Z$ is a centered
real Gaussian vector with covariance $C_\tau$ and $W\ge0$ is fixed under this conditioning, then
\[
 \|Z^{\mathsf T}WZ\|_{L^q(\mathbb P_\tau)}
 \le\|N(0,1)^2\|_{L^q}\operatorname{Tr}(WC_\tau),
 \qquad q\ge1.
 \eqtag{gaussian-quadratic-moment}
\]
The same estimate applies to a Gaussian random element whenever $Y=W^{1/2}Z$ has positive trace-class covariance $C_Y$. Choose a fixed orthonormal basis and let $P_m$ project onto
its first $m$ vectors. Then
\[
 \|P_mY\|^2\uparrow\|Y\|^2,\qquad
 \operatorname{Tr}(P_mC_YP_m)\uparrow\operatorname{Tr}C_Y.
\]
Apply the finite-dimensional estimate to $P_mY$ and pass to the limit in its $q$th moment by monotone convergence. This gives
$\|\|Y\|^2\|_{L^q}\le\|N(0,1)^2\|_{L^q}\operatorname{Tr}C_Y$.
In the applications $\operatorname{Tr}C_Y=\operatorname{Tr}(WC_\tau)$, and the estimate is applied conditionally on $\tau$ before integrating its distribution.

In the next lemma, $h$ is either model's deterministic spacing, $L=(N+1)h\asymp h^{-2}$, and $\iota_h$ is the piecewise constant isometry on $N$ coordinates. Constants may depend
on the fixed comparison constants in this relation.

\begin{lemma}\label{A:res:discretization-moments}
Under the full Wiener law,
\[
 \sup_{\beta\in I_\beta}
 \left\|\sup_{z\in K}
  \|\iota_h(\mathcal B_{N,h,\beta}^B-z)^{-1}\iota_h^*
                                  -R_\beta(z)\|_{\mathfrak S_2}^2
                      \right\|_{L^p}
 \le C_{p,I_\beta,K}h.
 \eqtag{discretization-moment}
\]
\end{lemma}

\begin{proof}
Use $P_B$ from Subsection~\ref{A:sec:continuum-forms-and-quantitative-continuity}. The forms with $e=c_\beta B$ and $e_h(x)=e(h\lfloor x/h\rfloor)$ have deterministic lower
coefficients after a shift at most $CP_B$. The estimates in Subsection~\ref{A:sec:galerkin-approximation-and-the-resolvent-error} bound the unit interval oscillations and
differences of consecutive averages of $e_h$ by $C M_B$, with the stated weight. Lemma~\ref{A:res:brownian-mesh} gives the corresponding bound for piecewise linear interpolation,
uniformly in the spacing and the length of the finite interval. The same is true of their finite difference approximations and the matrices of their restrictions to piecewise
linear functions. In this proof let $\kappa_B$ be a common shift for all these forms, with $\kappa_B\le CP_B$, increasing the deterministic prefactor if

necessary. We estimate the comparisons of Subsection~\ref{A:sec:finite-difference-approximation-of-the-stochastic-airy-operator} with the powers of $P_B$ retained.

For restriction to $(0,L)$, the cutoff identity of Lemma~\ref{A:res:4.2} has two tail terms and the two commutator terms with $\chi_L'\partial_x$ and $\chi_L''$.
Lemma~\ref{A:res:analytic-constants} bounds their squared norms by, respectively,
\[
 CP_B^4 L^{-1/2},\quad CP_B^4 L^{-3/2},\quad
 CP_B^6 L^{-3/2},\quad CP_B^4 L^{-3/2}.
\]
The third estimate factors $\partial_xR_\beta$ in operator norm; the multipliers are supported in an interval of fixed length at height $L$. Their sum is at most
$CP_B^{12}L^{-1/2}$.

The difference of the original and sampled forms on $(0,L)$ is $(h/2)\langle\cdot,\cdot\rangle+\Gamma_{e-e_h}$. The scalar resolvent identity and the first-order estimate bound its
squared resolvent error by $CP_B^{12}(h^2+\mathcal E_{\mathrm{int}}(h,L))$. The covariance trace calculated in Lemma~\ref{A:res:4.3} is at most
$C_{I_\beta}h^2\sum_jX_j^{-3/2}\le C_{I_\beta}h$. Applying \eqref{A:eq:gaussian-quadratic-moment} to the Gaussian function $X^{-3/4}(e-e_h)$ gives
\[
 \sup_{\beta\in I_\beta}
       \|\mathcal E_{\mathrm{int}}(h,L)\|_{L^q}
 \le C_{q,I_\beta}h.
 \eqtag{cell-moments}
\]

Let $R^{\mathrm{fr}}$ be the sampled continuum resolvent and $u=R^{\mathrm{fr}}(z)F$. On every open interval $I_j=(jh,(j+1)h)$ its equation is $u''=(x-h/2-z)u-F$. Thus
\[
 \left(\sum_j\|u''\|_{L^2(I_j)}^2\right)^{1/2}
 \le C_K(\|Xu\|_2+\|F\|_2)
 \le C_KP_B^3\|F\|_2.
\]
Only derivatives inside $I_j$ are used; the jumps of $u'$ at the two boundary points of $I_j$ are those of the sampled potential. The interpolation estimate \eqref{A:eq:4.33} gives
$\|u-\mathcal I_hu\|_{\mathcal E}\le CP_B^3h\|F\|_2$. The sampled sesquilinear form has continuity constant $CP_B$, and the inverse of its restriction to the piecewise linear space
has norm $CP_B$ on the energy scale, by the imaginary and real part calculation in Lemma~\ref{A:res:analytic-constants}. Orthogonality obtained by subtracting the two weak
equations therefore gives
\[
 \|u-u_h\|_{\mathcal E}
 \le(1+CP_B^2)\|u-\mathcal I_hu\|_{\mathcal E}
 \le CP_B^5h\|F\|_2.
\]
For $w=R^{\mathrm{fr}}(\overline z)(u-u_h)$, orthogonality and the adjoint equation give
\[
 \|u-u_h\|_2^2 =a_z(u-u_h,w-\mathcal I_hw) \le CP_B\|u-u_h\|_{\mathcal E} \|w-\mathcal
 I_hw\|_{\mathcal E} \le CP_B^9h^2\|F\|_2\|u-u_h\|_2.
\]
The operator norm error is thus at most $CP_B^9h^2$. The preceding bounds use $P_B^3P_B^2=P_B^5$ in the energy estimate and $P_BP_B^5P_B^3=P_B^9$ in the adjoint estimate.

To pass to Hilbert--Schmidt norm, let $P_h^{\mathrm{lin}}$ be the orthogonal projection onto the piecewise linear space. Let $\mu_k$ be the ordered eigenvalues of
$\mathcal H_{\beta,h,L}^{\mathrm{fr}}$ and put $t_k=\mu_k+\kappa_B\ge ck^{2/3}$. The projection estimate of Lemma~\ref{A:res:4.3C} gives
\[
 \|(I-P_h^{\mathrm{lin}})\psi_k\|^2
 \le C\min\{1,h^2t_k\},\qquad
 |\mu_k-z|^{-2}\le CP_B^2t_k^{-2}.
\]
The function $\min\{1,h^2t\}t^{-2}$ is bounded both by $h^2/t$ and by $t^{-2}$. The lower bound on $t_k$ is therefore enough for the next split; an upper bound on the sampled
eigenvalues is not being used. For $k_h=\lceil h^{-3}\rceil$, it follows that
\[
 \|(I-P_h^{\mathrm{lin}})R^{\mathrm{fr}}(z)\|_{\mathfrak S_2}^2
 \le CP_B^2\left(h^2\sum_{k\le k_h}k^{-2/3}
                       +\sum_{k>k_h}k^{-4/3}\right)
 \le CP_B^2h.
\]
Taking adjoints gives the other projection tail. The difference between the resolvent of the restricted form and $P_h^{\mathrm{lin}}R^{\mathrm{fr}}P_h^{\mathrm{lin}}$ has rank at
most $N\le Ch^{-3}$ and operator norm at most $CP_B^9h^2$. Its squared Hilbert--Schmidt norm is at most $CP_B^{18}h$. Adding the two projection tails proves this bound for the full
comparison with the variational restriction. Here the rank bound is applied to the operator on $\operatorname{Ran}P_h^{\mathrm{lin}}$; the two complementary components have the
squared bounds just proved.

It remains to change the mass matrix and the embedding. Use the exact identity of Proposition~\ref{A:res:4.5B}, writing $R_{\mathrm l}=(\mathcal B_{N,h,\beta}^B-z)^{-1}$ and
$R_{\mathrm c}=(K_h^{\mathrm c}-zM_h)^{-1}$:
\[
 R_{\mathrm l}-R_{\mathrm c} ={} -R_{\mathrm c}M_a-zR_{\mathrm c}M_aR_{\mathrm l} +R_{\mathrm
 c}M_{ax}R_{\mathrm l} +R_{\mathrm c}M_{a\omega^B}R_{\mathrm l} +R_{\mathrm
 c}M_bJ_-^*D_hR_{\mathrm l},
\]
where $a_j=h^2((j-1)h-z)/6$ and $b_j=h^2/6$. In particular, $a_1=-h^2z/6$ is retained. The bounds $|a_j|\le C_Kh^2X_j$, $|b_j|\le Ch^2$, and $\sum_{j\le N}X_j^{1/2}\le Ch^{-4}$
give
\[
 \|R_{\mathrm c}M_a\|_{\mathfrak S_2}^2 +\|zR_{\mathrm c}M_aR_{\mathrm l}\|_{\mathfrak S_2}^2 \le
 CP_B^4h,\quad \|R_{\mathrm c}M_{ax}R_{\mathrm l}\|_{\mathfrak S_2}^2 \le CP_B^{11}h,\quad
 \|R_{\mathrm c}M_bJ_-^*D_hR_{\mathrm l}\|_{\mathfrak S_2}^2 \le CP_B^6h^4.
\]
For the last inequality use $\|R_{\mathrm c}M_bJ_-^*\|_{\mathfrak S_2} =\|R_{\mathrm c}M_b\|_{\mathfrak S_2}$ and factor the differentiated resolvent in operator norm. The Brownian
term has squared norm at most $CP_B^{11}\mathcal Z_{\omega,h}$, where
\[
 \mathcal Z_{\omega,h}
   :=h^5\sum_{j\le N}X_j^{-3/2}|\omega_j^B|^2,
 \qquad
 \sup_{\beta\in I_\beta}\|\mathcal Z_{\omega,h}\|_{L^q}
          \le C_{q,I_\beta}h^3.
 \eqtag{mesh-gaussian-moment}
\]
Indeed $\omega_j^B=c_\beta\Delta_jB/h$ has variance $c_\beta^2/h$, so its weighted covariance trace is $C_{I_\beta}h^4\sum_jX_j^{-3/2}\le C_{I_\beta}h^3$.
Equation~\eqref{A:eq:gaussian-quadratic-moment} proves the moment bound. The squared triangle inequality applied to the five differences bounds the mass and quadrature comparison
by
$CP_B^{12}(h+\mathcal Z_{\omega,h})$.

Finally $\|I_h\|\le1$, $\|\iota_h\|=1$, and the exact interpolation identity following Proposition~\ref{A:res:4.5B} is $\|(I_h-\iota_h)u\|_2^2=(h^2/3)\|D_hu\|^2$. Writing the
change of embedding as
\[
 (I_h-\iota_h)R_{\mathrm l}I_h^*
      +\iota_hR_{\mathrm l}(I_h^*-\iota_h^*)
\]
and taking adjoints for the second term bounds its squared norm by $CP_B^2h^2h^{-1}=CP_B^2h$. Adding these bounds gives
\[
 \sup_{z\in K}
 \|\iota_hR_{\mathrm l}(z)\iota_h^*-R_\beta(z)\|_{\mathfrak S_2}^2
 \le C_KP_B^{24}
       \{h+\mathcal E_{\mathrm{int}}(h,L)+\mathcal Z_{\omega,h}\}.
 \eqtag{discretization-budget}
\]
Here $L^{-1/2}\le Ch$, and $P_B^{24}$ bounds each of the powers $P_B^{12}$ and $P_B^{18}$ already obtained. The Brownian approximation errors and $P_B$ depend on the same path.
Taking the $L^p$ norm, H\"older gives, for example,
\[
 \|P_B^{24}\mathcal E_{\mathrm{int}}\|_p
 \le\|P_B^{24}\|_{2p}\|\mathcal E_{\mathrm{int}}\|_{2p}
 \le C_ph.
\]
The other two terms satisfy
\[
 \|P_B^{24}h\|_p=h\|P_B^{24}\|_p\le C_ph,\qquad
 \|P_B^{24}\mathcal Z_{\omega,h}\|_p
 \le\|P_B^{24}\|_{2p}\|\mathcal Z_{\omega,h}\|_{2p}
 \le C_ph^3.
\]
All moments of $P_B$ are finite. Substitution into \eqref{A:eq:discretization-budget} proves \eqref{A:eq:discretization-moment}.
\end{proof}

For the Gaussian model, use the weighted errors $\mathcal E_{\eta,n}^{(\beta)}$ and $\mathcal E_{\mathrm{scal},n}^{(\beta)}$ of Propositions~\ref{A:res:A.4}--\ref{A:res:A.5}. The
quantile estimate of Lemma~\ref{A:res:A.3} is uniform for $\beta\in I_\beta$ and gives
\[
 \|\eta_j\|_{2p}^2\le C_{p,I_\beta}h,\qquad
 \|\Delta\zeta_j\|_{2p}^2
            \le C_{p,I_\beta}\frac{h^3}{1+n-j}.
\]
Minkowski therefore gives
\[
 \sup_{\beta\in I_\beta}\|\mathcal E_{\eta,n}^{(\beta)}\|_p
 \le C_ph\sum_jhX_j^{-3/2}\le C_{p,I_\beta}h.
 \eqtag{B.46}
\]
For $j\le n/2$, use $(1+n-j)^{-1}\le2/n$. For $j>n/2$, use $X_j^{-7/2}\le Ch^7$ and sum $(1+n-j)^{-1}$, including its value at $j=n$. Thus
\[
 \sup_{\beta\in I_\beta}\|\mathcal E_{\mathrm{scal},n}^{(\beta)}\|_p \le
 C_p\sum_{j=1}^n\frac{X_j^{-7/2}}{1+n-j} \le C_p(n^{-1}h^{-1}+h^7\ell_n) \le C_{p,I_\beta}h^2.
 \eqtag{B.47}
\]
Define
\[
 D_{n,\beta}^{\mathrm{DE}}(K)
 =\sup_{z\in K}\|
  \iota_h(\mathcal A_{n,\beta}^{\mathrm{DE}}-z)^{-1}\iota_h^*
                              -R_\beta(z)\|_{\mathfrak S_2}^2.
\]

\par\Needspace{6\baselineskip}
\begin{theorem}\label{A:res:B.9}
For every $p\ge1$,
\[
 \sup_{\beta\in I_\beta}
 \|D_{n,\beta}^{\mathrm{DE}}(K)\|_{L^p(\Pi_{n,\beta}^{\mathrm{DE}})}
 \le C_{p,I_\beta,K}h.
 \eqtag{B.50}
\]
\end{theorem}

\begin{proof}

The maximum estimate gives the form comparison, and weighted sums give the Hilbert--Schmidt rate. We first bound the primitives and coordinate maxima uniformly.
Lemma~\ref{A:res:A.2} gives Gaussian tails for the centered chi variables, uniformly in all degrees of freedom, including degree zero. A centered variable $Y$ with these tails
satisfies $\mathbb Ee^{tY}\le e^{Ct^2}$ for all real $t$. Let $Y'$ be an independent copy. Jensen gives $\mathbb Ee^{tY}\le\mathbb Ee^{t(Y-Y')}$; the difference is symmetric and
has moments $\mathbb E|Y-Y'|^{2k}\le C^kk!$ after increasing $C$. Symmetry and the moment bound give
\[
 \mathbb E\mathrm e^{t(Y-Y')}
 =\sum_{k=0}^\infty\frac{t^{2k}\mathbb E(Y-Y')^{2k}}{(2k)!}
 \le\sum_{k=0}^\infty\frac{(Ct^2)^k}{k!}
 =\mathrm e^{Ct^2},
\]
where $k!/(2k)!\le1/k!$. Rescaling and using independence across sites now give
\[
 \mathbb E\exp\left(t h^{-1}\sum_{j=r+1}^s\zeta_{j,n}\right)
 \le \exp\{C_{I_\beta}t^2h(s-r)\}.
\]
Optimizing exponential Markov for both signs and integrating the resulting Gaussian tail proves
\[
 \|S_s^\zeta-S_r^\zeta\|_{\psi_2}
 \le C_{I_\beta}\sqrt{h(s-r)},\qquad
 S_m^\zeta=h^{-1}\sum_{j\le m}\zeta_{j,n}.
 \eqtag{B.48}
\]
The one coordinate bound for $\eta$ gives
\[
 \mathbb P(\max_{j\le n}|\eta_{j,n}|>t)
 \le\min\!\left\{1,\sum_{j\le n}\mathbb P(|\eta_{j,n}|>t)\right\}
 \le\min\{1,2n\mathrm e^{-ct^2/h}\}.
\]
Take $t_0=(2h\log(2n)/c)^{1/2}$. For $t\ge t_0$, $2ne^{-ct^2/h}\le e^{-ct^2/(2h)}$. With $M_n=\max_{j\le n}|\eta_{j,n}|$, integration gives
\[
 \mathbb EM_n^q
 \le t_0^q+q\int_{t_0}^\infty t^{q-1}e^{-ct^2/(2h)}\,\mathrm dt
 \le t_0^q+\frac q2(2h/c)^{q/2}\Gamma(q/2).
\]
Taking $q$th roots, and using $\log(2n)\le C(1+|\log h|)$, gives
\[
 \left\|\max_{j\le n}|\eta_{j,n}|\right\|_q
 \le C_{q,I_\beta}h^{1/2}\sqrt{1+|\log h|}.
 \eqtag{B.49}
\]

Fix $A>4p+4$, set $E=E_{n,\beta,A}^{\mathrm{DE}}$, and abbreviate $G=G_{n,\beta,K,A}$, increasing its deterministic prefactor so that $G\ge1$ on $E$. Use the finite comparisons of
Subsections~\ref{A:sec:the-mean-dumitriu-edelman-coefficients} and~\ref{A:sec:gaussian-resolvent-comparison} before imposing any almost sure bound on their weighted errors. The
first-order and scalar differences contribute at most $G\mathcal E_{\eta,n}^{(\beta)}$ and $G\mathcal E_{\mathrm{scal},n}^{(\beta)}$, respectively. For the square-root coefficients
retain $R_1,R_0,a,b,r^0$ from \eqref{A:eq:7.38}--\eqref{A:eq:7.39}, and put $X_j=1+x_j$. The weighted estimates on $E$ give
\[
 \begin{aligned}
 \|R_1M_a\|_{\mathfrak S_2}^2
   +\|zR_1M_aR_0\|_{\mathfrak S_2}^2
   &\le CGh^5\sum_{j\le J}X_j^{1/2}\le CGh,\\
 \|R_1M_{ax}R_0\|_{\mathfrak S_2}^2
   +\|R_1M_{r^0}R_0\|_{\mathfrak S_2}^2
   &\le CGh^5\sum_{j\le J}X_j^{1/2}\le CGh,\\
 \|R_1M_{a\omega^B}R_0\|_{\mathfrak S_2}^2
   &\le CGh^5\sum_{j\le J}X_j^{-3/2}|\omega_j^B|^2
    \le CG\mathcal Z_{\omega,h},\\
 \|R_1M_bJ_-^*D_hR_0\|_{\mathfrak S_2}^2
   &\le CG^2h^5\sum_{j\le J}X_j^{-3/2}\le CG^2h^4.
 \end{aligned}
\]
Here $\sum_{j\le J}X_j^{1/2}\le Ch^{-4}$ and $\sum_jX_j^{-3/2}\le Ch^{-1}$. For the correction $c_{k,n}=\overline b_{k,n}-\sqrt{1-k/n}$, the two products in
\eqref{A:eq:chi-kinetic}--\eqref{A:eq:chi-diagonal} satisfy
\[
 \begin{aligned}
 \|R_2D_h^*M_cD_hR_1\|_{\mathfrak S_2}^2
 &\le CG^2\|c\|_\infty^2h^{-1}\le CG^2h^5,\\
 \|R_2M_{s^\chi}R_1\|_{\mathfrak S_2}^2
 &\le CGh\sum_j|s_j^\chi|^2X_j^{-7/2}\le CGh^2.
 \end{aligned}
\]
For either cut, its coefficient has modulus at most $h^{-2}$, whereas each resolvent column at the cut has squared norm at most $CGh^4$. The rank-two identity and the tail estimate
therefore give
\[
 \|R-R^\oplus\|_{\mathfrak S_2}^2\le CG^2h^{-4}h^8=CG^2h^4,
 \qquad
 \|R_{\mathrm{tail}}\|_{\mathfrak S_2}^2
 \le CGh\sum_{j>J}X_j^{-3/2}\le CGh.
\]
Using $G\ge1$ on $E$ and $\|\sum_{j=1}^mT_j\|_{\mathfrak S_2}^2\le m\sum_{j=1}^m\|T_j\|_{\mathfrak S_2}^2$, with a fixed number of terms, bounds their sum by
$CG^2(h+\mathcal Z_{\omega,h})$.

Consequently the squared difference $D_{n,\beta}^{\mathrm{fin}}$ between the full finite matrix and the Brownian difference resolvent satisfies
\[
 \mathbf1_E D_{n,\beta}^{\mathrm{fin}}
 \le CG^2\{h+\mathcal E_{\eta,n}^{(\beta)}
              +\mathcal E_{\mathrm{scal},n}^{(\beta)}
              +\mathcal Z_{\omega,h}\}.
\]
All terms on the right have $L^{2p}$ norm $O(h)$ or smaller, by \eqref{A:eq:B.46}--\eqref{A:eq:B.47} and \eqref{A:eq:mesh-gaussian-moment}; $\|G^2\|_{2p}\le C_{p,A}$. H\"older
therefore bounds the norm of the left side by
$C_ph$.

On $E^c$ we keep these same two matrices. Self-adjointness holds without the coefficient bounds defining $E$, and each resolvent has squared Hilbert--Schmidt norm at most
$n\delta_K^{-2}$. Hence $D_{n,\beta}^{\mathrm{fin}}\le4n\delta_K^{-2}$ there, and
\[
 \|\mathbf1_{E^c}D_{n,\beta}^{\mathrm{fin}}\|_p
 \le C_{p,A}n^{1-A/p}\le C_ph.
\]
Apply the isometry, add the error of Lemma~\ref{A:res:discretization-moments}, and use $\|T_1+T_2\|_{\mathfrak S_2}^2 \le2\|T_1\|_{\mathfrak S_2}^2+2\|T_2\|_{\mathfrak S_2}^2$.
Minkowski proves \eqref{A:eq:B.50}. The same finite-rank bound and \eqref{A:eq:continuum-coarse} cover the finitely many smaller sizes.

\end{proof}

For the general model let $A_C$ be the principal matrix on $I_n^*$ and define the actual boundary error by
\[
 \mathcal C_{\partial,n,\beta}
 =\sup_{z\in K}\| (\mathcal A_{n,V,\beta}-z)^{-1}
   -P_{I_n^*}^*(A_C-z)^{-1}P_{I_n^*}\|_{\mathfrak S_2}^2.
\]

\begin{lemma}\label{A:res:B.10B}
For every $p\ge1$,
\[
 \sup_{\beta\in I_\beta}\|\mathcal C_{\partial,n,\beta}\|_p
 \le C_{p,V,I_\beta,K}h.
 \eqtag{B.54d}
\]
\end{lemma}

\begin{proof}
Omit the final $s_*$ indices of \eqref{A:eq:B.45} and use the fixed left cut $r_n$. The coefficient events in Proposition~\ref{A:res:left-spectral} and
Lemma~\ref{A:res:core-finite} have uniform probabilities on $I_\beta$. Their constants depend on the fixed polynomial, the positive lower bound for $\beta$, and the deterministic
separation lengths. In particular the gap for the left block comes from the same deterministic half-line matrix $J^V$, independently of $\beta$. After deleting the final $s_*$
coordinates, all retained logarithmic weights are positive. Lemmas~\ref{A:res:T2}--\ref{A:res:T3} give the bounds on the final coordinate block also when $\beta<1$. Thus the event
of \eqref{A:eq:boundary-good-event} satisfies
\[
 \sup_{\beta\in I_\beta}
 \Pi_{n,\beta}^V(\mathcal G_{\partial,n,A}^c)
       \le C_{A,V,I_\beta}n^{-A}.
 \eqtag{B.54d1}
\]
On this event, the four blocks in each inverse formula obtained from a Schur complement were estimated in \eqref{A:eq:boundary-finite-rate}; both denominators are bounded there.
Substituting $r_n\le C\ell_n$ into those estimates gives
\[
 \mathcal C_{\partial,n,\beta}
 \le C_A b_n,\qquad
 b_n=h^2\ell_n^{16}+h^3\ell_n^{12}+h^4\ell_n^8.
\]
The estimate controls the full inverse before compression, including the uncorrected left inverse and both off-diagonal blocks. These are the blocks required for the isometric
spectral comparison. Since $h\asymp n^{-1/3}$, $b_n/h\to0$. On the complement, the unconditional nonreal bound gives $\mathcal C_{\partial,n,\beta}\le4n\delta_K^{-2}$. Therefore
\[
 \sup_\beta\|\mathcal C_{\partial,n,\beta}\|_p
 \le C_A b_n+C_{p,A}n^{1-A/p}.
\]
Choose $A>5p/3$. Then $n^{1-A/p}\le C h^2\le C b_n$ for large $n$. We have proved the stronger bound $\sup_\beta\|\mathcal C_{\partial,n,\beta}\|_p\le C_pb_n=o(h)$, which implies
\eqref{A:eq:B.54d}. Increasing the constant covers the smaller sizes for which the reduced interval is defined, and also the convention that the compression is zero for an empty
interval.
\end{proof}

To specify the remaining errors, for self-adjoint matrices on the retained coordinate space set
\[
 \Delta_K(A,C)=\sup_{z\in K}
       \|(A-z)^{-1}-(C-z)^{-1}\|_{\mathfrak S_2}^2.
\]
Suppressing $(n,\beta)$, use the following notation for the successive comparisons in Section~\ref{A:sec:the-general-jacobi-model-gaussian-approximation-and-comparison}:
\[
 \begin{aligned}
 \mathcal C_{\mathrm{off}}&=\Delta_K(\widetilde H_n^{\mathrm{red}},H_n^{\mathrm{diag}}),&
 \mathcal C_{\mathrm{Gauss}}&=\Delta_K(H_n^{\mathrm{diag}},H_n^{\mathrm G}),\\
 \mathcal C_{\mathrm{cov}}^{\mathrm{diag}}&=\Delta_K(H_n^{\mathrm G},\widetilde H_n),&
 \mathcal C_{\mathrm{cov}}^{\mathrm{der}}&=\Delta_K(\widetilde H_n,H_n^{B,\mathrm{mid}}),\\
 \mathcal C_{\mathrm{left}}&=\Delta_K(H_n^{B,\mathrm{mid}},H_n^{B,\mathrm{reg}}),&
 \mathcal C_{\mathrm{high}}&=\Delta_K(H_n^{B,\mathrm{reg}},H_n^{B,\circ}),\\
 \mathcal C_{\mathrm{prof}}&=\Delta_K(H_n^{B,\circ},\mathcal B_{M,h,\beta}^B),&
 \mathcal C_{\mathrm{cov}}&=\mathcal C_{\mathrm{cov}}^{\mathrm{diag}}
                              +\mathcal C_{\mathrm{cov}}^{\mathrm{der}}.
 \end{aligned}
\]
Here $M=M_n^*$, all matrices are reindexed on $\{1,\ldots,M\}$, and $\mathcal C_{\mathrm{disc}}$ is the squared error in \eqref{A:eq:discretization-moment}. Let
$\mathcal I=\{\partial,\mathrm{off},\mathrm{Gauss},\mathrm{cov}, \mathrm{left},\mathrm{high},\mathrm{prof},\mathrm{disc}\}$.

\par\Needspace{6\baselineskip}
\begin{proposition}\label{A:res:B.10}
For $p\ge1$ and $r\in\mathcal I$,
\[
 \sup_{\beta\in I_\beta}
 \|\mathcal C_{r,n,\beta}\|_{L^p(\overline\Pi_{n,\beta}^V)}
 \le C_{p,V,I_\beta,K}h.
\]
\end{proposition}

Consequently, for $L>0$,
\[
 \sup_{\beta\in I_\beta}
 \overline\Pi_{n,\beta}^V(\mathcal C_{r,n,\beta}>Lh)
 \le C_{p,V,I_\beta,K}L^{-p}.
 \eqtag{B.56}
\]
The smaller error terms satisfy the sharper uniform moment bounds
\[
 \|\mathcal C_{\mathrm{Gauss}}\|_p +\|\mathcal C_{\mathrm{cov}}^{\mathrm{diag}}\|_p \le
 C_ph^2,\qquad \|\mathcal C_{\mathrm{left}}\|_p \le C_ph^3\ell_n^5,\qquad \|\mathcal
 C_{\mathrm{high}}\|_p\le C_ph^4\ell_n.
\]
The supremum over $\beta\in I_\beta$ is understood in each of these inequalities, outside the norm.

\begin{proof}
Fix $p\ge1$, choose $A>4p+4$, and set $E=E_{n,\beta,A}^V$ and $G=G_{n,\beta,K,A}$, with $G\ge1$ on $E$. The boundary estimate is Lemma~\ref{A:res:B.10B}. For the reduced
comparisons we first estimate the errors multiplied by $\mathbf1_E$ and then integrate their exceptional complements.

The two conditional minimizers use different exterior sigma fields. Let $\tau$ denote either exterior, with marginal $\nu$, and let $T_\tau$ be its admissible bounded set. The
event $E$ implies $\tau\in T_\tau$ in both cases and may also depend on the free coordinates. Dropping its indicator inside the conditional expectation allows us to apply the
Gaussian or Gibbs estimate under the original conditional law. For every nonnegative $F$ and $q\ge1$, disintegration gives
\[
 \mathbb E[\mathbf1_EF^q]
 \le\int_{T_\tau}\mathbb E[F^q\mid\tau]\,\nu(\mathrm d\tau).
\]

For the off-diagonal fluctuation, Lemma~\ref{A:res:mean} gives $\|\eta_j\|_{2q}\le C_qh^{1/2}$ at regular coordinates and $C_qh^{1/2}\sqrt{\ell_n}$ on the far interval, after
restriction to the corresponding exterior set. The latter interval has $X_j\ge ch^{-2}$. Hence
\[
 \left\|\mathbf1_E\sum_jhX_j^{-3/2}|\eta_j|^2\right\|_q \le C_qh\sum_{\mathrm{reg}}hX_j^{-3/2}
 +C_qh\ell_n\sum_{\mathrm{far}}hX_j^{-3/2} \le C_q(h+h^2\ell_n)\le C_qh.
 \eqtag{off-local-moment}
\]
At the boundary edges, $\eta$ has the deterministic values. The first-order estimate and H\"older at order $2p$ prove
$\|\mathbf1_E\mathcal C_{\mathrm{off}}\|_p\le C_ph$.

For the transport error use the two conditional minimizers $\overline Y_n^\circ$ and $Y_n^\circ$ of Subsection~\ref{A:sec:gaussian-approximation-of-the-conditional-gibbs-measure}.
At a middle scalar coordinate,
\[
 t_j:=\zeta_j-\widehat\zeta_j^{\mathrm G} ={} -b_\star^{-1} \{(A_j-\overline A_j^\circ-g_j^A)
 +2(B_j-\overline B_j^\circ-g_j^B)\} -b_\star^{-1} \{\overline A_j^\circ-A_j^\circ +2(\overline
 B_j^\circ-B_j^\circ)\}.
\]
Proposition~\ref{A:res:14.6} bounds the first summand by $C_qn^{-1}$ in conditional $L^{2q}$, uniformly for $\tau\in\mathcal T_n^{\mathrm{out}}(M_0)$ under the original coefficient
law. The second summand is bounded by $C_An^{-4}$ on $E$, by Proposition~\ref{A:res:14.1} and the right separation length. The disintegration inequality and Minkowski's inequality
give
\[
 \left\|\mathbf1_E\sum_{I_n^{\mathrm{mid}}}X_j^{-7/2}|t_j|^2\right\|_q \le
 C_{q,A}(n^{-2}+n^{-8})\sum_jX_j^{-7/2} \le C_{q,A}h^5.
 \eqtag{transport-local-moment}
\]
The scaled multiplier is $h^{-2}t$, so the diagonal estimate supplies the factor $h^{-3}$. H\"older yields
$\|\mathbf1_E\mathcal C_{\mathrm{Gauss}}\|_p\le C_ph^2$.

The conditionally jointly Gaussian coupling constructed in Proposition~\ref{A:res:15.2} satisfies, conditionally on its admissible exterior,
\[
 \operatorname{Tr}\bigl(
    \operatorname{diag}(X_j^{-7/2})
           \operatorname{Cov}_\tau(\varepsilon^{\mathrm{fr}})\bigr)
 \le Ch^5.
\]
By \eqref{A:eq:gaussian-quadratic-moment}, its weighted square sum has conditional $L^q$ norm at most $C_qh^5$. For $\rho$, \eqref{A:eq:factor-remainder-variances} gives
$\operatorname{Var}_\tau\rho_j\le C(h^7X_j^2+h^9)$. Therefore
\[
 \left\|\sum_jX_j^{-7/2}|\rho_j|^2\right\|_{L^q(\mathbb P_\tau)}
 \le C_q\left(h^7\sum_jX_j^{-3/2}
                          +h^9\sum_jX_j^{-7/2}\right)
 \le C_qh^6.
\]
For their sum use $|u+v|^2\le2|u|^2+2|v|^2$. Apply the disintegration inequality above, the diagonal estimate, and then H\"older with $G$. This proves
$\|\mathbf1_E\mathcal C_{\mathrm{cov}}^{\mathrm{diag}}\|_p \le C_p(h^2+h^3)$.

Each $e_j$ in the spectral factorization has variance at most $Ch$, including the additional predecessor value. The continuation $\widetilde e$ equals that predecessor before the
middle interval and its final value after it. All added edges have the same variance bound, and
\[
 \left\|\mathbf1_E\sum_{j=0}^MhX_j^{-3/2}
                          |\widetilde e_j|^2\right\|_q
 \le C_qh\sum_{j=0}^MhX_j^{-3/2}\le C_qh.
 \eqtag{primitive-local-moment}
\]
At $j=0$ the convention $x_0=0$ gives $X_0=1$. Many of the continued coordinates are identical random variables. Their individual moments and Minkowski's inequality give the
estimate. The constant continuation has zero difference outside the middle interval; it therefore preserves the exact scalar identity at both cutoffs. The first-order perturbation
estimate and H\"older give $\|\mathbf1_E\mathcal C_{\mathrm{cov}}^{\mathrm{der}}\|_p\le C_ph$.

At the indices next to the left boundary, $E$ gives $\max|h^{-2}(\zeta_j-\zeta_j^B)| \le C_Ah^{-1/2}\sqrt{\ell_n}$ and $J_n^{\mathrm L}\le C\ell_n$. The boundary estimate of
Lemma~\ref{A:res:analytic-constants} implies
\[
 \mathbf1_E\mathcal C_{\mathrm{left}}
 \le C_AGh^4\ell_n^4(h^{-1/2}\sqrt{\ell_n})^2
 =C_AGh^3\ell_n^5.
\]
Taking moments proves the claimed bound. On the far interval, Lemma~\ref{A:res:mean} gives $\|\mathbf1_E\zeta_j\|_{2q}\le C_qh^{3/2}\sqrt{\ell_n}$; the Brownian marginal gives
$\|\zeta_j^B\|_{2q}\le C_qh^{3/2}$. Since $\sum_{I_n^{\mathrm{far}}}X_j^{-7/2}\le Ch^4$,
\[
 \left\|\mathbf1_E h^{-3}\sum_{I_n^{\mathrm{far}}}X_j^{-7/2}
                      |\zeta_j-\zeta_j^B|^2\right\|_q
 \le C_qh^4\ell_n.
 \eqtag{high-local-moment}
\]
The diagonal estimate and H\"older prove the corresponding bound for $\mathcal C_{\mathrm{high}}$.

Outside the indices adjacent to the left boundary, the conditional minimizer and the local profile differ by at most $C_An^{-1}$ at each regular coordinate. Let $W$ denote this
tridiagonal interior difference on the regular comparison block of $J$ vertices. Define
\[
 d_j=-h^2W_{j,j+1}\quad(1\le j<J),\qquad d_0=d_J=0,
 \qquad q_j=W_{jj}-h^{-2}(d_{j-1}+d_j).
\]
The diagonal and off-diagonal entries then give $W=D_h^*M_dD_h+M_q$. The unscaled coefficient bound and the factor $h^{-2}/b_\star$ in the matrix imply $\|d\|_\infty\le C_Ah^3$ and
$\|q\|_\infty\le C_Ah$, including the first and last entries. For the two resolvents involved,
\[
 \begin{aligned}
 \|R_1D_h^*M_dD_hR_0\|_{\mathfrak S_2}^2
 &\le\|D_hR_1(\overline z)\|_{\mathrm{op}}^2
          \|d\|_\infty^2\|D_hR_0(z)\|_{\mathfrak S_2}^2
       \le C_AG^2h^5,\\
 \|R_1M_qR_0\|_{\mathfrak S_2}^2
 &\le C_AGh\sum_jh^2X_j^{-7/2}\le C_AGh^2.
 \end{aligned}
\]
At the indices next to the left boundary its matrix entries are at most $C_Ah^{-1/2}\sqrt{\ell_n}$. Its length is at most $C\ell_n$, so Lemma~\ref{A:res:analytic-constants} bounds
the squared resolvent error by $C_AGh^4(\ell_n+1)^4(h^{-1/2}\sqrt{\ell_n})^2 \le C_A'Gh^3\ell_n^5$. The cut at the far end of the regular interval is a rank-two matrix
$t(e_Je_{J+1}^*+e_{J+1}e_J^*)$, with $|t|\le Ch^{-2}$. The column bound at $X_J\asymp h^{-2}$ is $CGh^4$. Expansion into the two rank-one resolvent products gives $CG^2h^4$ for the
squared cut error. The two remaining tails are at most
$CGh$.

For the next comparison, $d$ denotes the local profile coefficient minus one. It has the larger bound $|d_j|\le Ch^2X_j$, which need not be small at the far end of the regular
block. Compare the local coefficients with $\mathcal B_{J,h,\beta}^{B,\lambda} =T_h+X_h+M_{\omega^B}+\lambda_nI$, where $|\lambda_n|\le Ch\ell_n$. For $a_j=d_{j-1}$ and
$b_j=h^{-1}(d_j-d_{j-1})$, the product rule is
\[
 R_1D_h^*M_dD_hR_0 ={} R_1M_a+(z-\lambda_n)R_1M_aR_0-R_1M_{ax}R_0
 -R_1M_{a\omega^B}R_0-R_1M_bJ_-^*D_hR_0.
\]
All edge coefficients, including those at the two Dirichlet ends, are retained. Proposition~\ref{A:res:17.3} gives $|a_j|\le Ch^2X_j$, $|b_j|\le Ch^2$, and a deterministic
potential remainder $|r_j|\le Ch^2X_j^2$. The three weighted sums are
\[
 h\sum_j|a_j|^2X_j^{-3/2}\le Ch,\qquad
 h\sum_j|b_j|^2X_j^{-3/2}\le Ch^4,\qquad
 h\sum_j|r_j|^2X_j^{-7/2}\le Ch.
\]
The first and third use $\sum_{j\le J}X_j^{1/2}\le Ch^{-4}$. The $ax$ term has the same bound as $r$. Although the Brownian potential cancels from the difference of the two
operators, it reappears when the resolvent equation eliminates $T_hR_0$. That term is bounded by $G\mathcal Z_{\omega,h}$, whose moments are $O(h^3)$. Finally the scalar shift has
squared resolvent error at most $CG\lambda_n^2\le CGh^2\ell_n^2$. Adding all these comparisons gives
\[
 \mathbf1_E\mathcal C_{\mathrm{prof}}
 \le C_A(G+G^2)
       \{h+h^2\ell_n^2+h^3\ell_n^5+\mathcal Z_{\omega,h}\}.
\]
H\"older and \eqref{A:eq:B.54c}, at the required doubled orders, show that its $L^p$ norm is at most $C_ph$. Lemma~\ref{A:res:discretization-moments} gives the discretization bound
without a coefficient event.

Every remaining finite error on $E^c$ is a squared difference, or a fixed sum of squared differences, of resolvents of self-adjoint matrices of dimension at most $Cn$. It is at
most $C_Kn$. Thus its contribution to the $L^p$ norm is at most $C_{p,A}n^{1-A/p}$. Our choice $A>4p+4$ makes this smaller than $h^4$, and hence preserves each of the stronger
orders stated above. All the estimates are therefore under the original, unconditional coupling laws. Smaller admissible sizes are covered by increasing the constants. Markov's
inequality applied to $\mathbb E\mathcal C_r^p\le C_ph^p$ proves \eqref{A:eq:B.56}.
\end{proof}

Set
\[
 D_{n,\beta}^V(K)=\sup_{z\in K}
 \|\mathcal J_{n,I_\beta}(\mathcal A_{n,V,\beta}-z)^{-1}
                   \mathcal J_{n,I_\beta}^*-R_\beta(z)\|_{\mathfrak S_2}^2.
\]

\begin{theorem}\label{A:res:B.11}
For every $p\ge1$,
\[
 \sup_{\beta\in I_\beta}\|D_{n,\beta}^V(K)\|_{L^p(\Pi_{n,\beta}^V)}
 \le C_{p,V,I_\beta,K}h.
 \eqtag{B.64}
\]
\end{theorem}

In particular, for $L>0$,
\[
 \sup_{\beta\in I_\beta}\Pi_{n,\beta}^V(D_{n,\beta}^V(K)>Lh)
 \le C_{p,V,I_\beta,K}L^{-p}.
 \eqtag{B.65}
\]
For a fixed deterministic $\beta$, the previously constructed countable coupling also satisfies $D_{n,\beta}^V(K)\le C_{\omega,V,\beta,K}h$ almost surely for all sufficiently large
$n$.

\begin{proof}
All reduced comparison matrices are now reindexed on $\{1,\ldots,M_n^*\}$. Write $E_h=\iota_h^{(M_n^*)}$ and denote their successive values by
\[
 \begin{aligned}
 H_0&=\widetilde H_n^{\mathrm{red}},& H_1&=H_n^{\mathrm{diag}},&
 H_2&=H_n^{\mathrm G},& H_3&=\widetilde H_n,\\
 H_4&=H_n^{B,\mathrm{mid}},&H_5&=H_n^{B,\mathrm{reg}},&
 H_6&=H_n^{B,\circ},&H_7&=\mathcal B_{M_n^*,h,\beta}^{B}.
 \end{aligned}
\]
Let $S_n(z)$ be the compressed full resolvent in the definition of $D_{n,\beta}^V(K)$. Define
\[
 \begin{aligned}
 T_0(z)&=S_n(z)-E_h(H_0-z)^{-1}E_h^*,\\
 T_j(z)&=E_h\bigl[(H_{j-1}-z)^{-1}-(H_j-z)^{-1}\bigr]E_h^*,
                        &&1\le j\le7,\\
 T_8(z)&=E_h(H_7-z)^{-1}E_h^*-R_\beta(z).
 \end{aligned}
\]
These definitions give the exact telescoping identity $S_n-R_\beta=\sum_{j=0}^8T_j$. By compression of the full boundary difference,
\[
 \sup_K\|T_0\|_{\mathfrak S_2}^2\le\mathcal C_{\partial,n,\beta}.
\]
Isometric invariance bounds the seven middle squared norms by their corresponding finite errors. The two covariance differences contribute
$\mathcal C_{\mathrm{cov}}^{\mathrm{diag}}+\mathcal C_{\mathrm{cov}}^{\mathrm{der}}=\mathcal C_{\mathrm{cov}}$. The last difference is precisely the discretization error.
Cauchy--Schwarz applied to the sum of the nine resolvent differences gives
\[
 D_{n,\beta}^V(K)
 \le9\sum_{j=0}^8\sup_K\|T_j\|_{\mathfrak S_2}^2
 \le9\sum_{r\in\mathcal I}\mathcal C_{r,n,\beta}.
\]
For $p\ge1$, Minkowski's inequality and Proposition~\ref{A:res:B.10} imply
\[
 \sup_{\beta\in I_\beta}\|D_{n,\beta}^V(K)\|_p
 \le9\sum_{r\in\mathcal I}\sup_{\beta\in I_\beta}\|\mathcal C_{r,n,\beta}\|_p
 \le C_{p,V,I_\beta,K}h.
\]
The boundary contribution here is the separately proved estimate in Lemma~\ref{A:res:B.10B}. The quantity on the left depends only on the path and the original matrix, so equality
of the corresponding marginals transfers this estimate from $\overline\Pi_{n,\beta}^V$ to $\Pi_{n,\beta}^V$. This proves \eqref{A:eq:B.64}. Finally, Markov's inequality gives
\[
 \Pi_{n,\beta}^V(D_{n,\beta}^V(K)>Lh)
 \le(Lh)^{-p}\mathbb E(D_{n,\beta}^V(K))^p\le C_{p,V,I_\beta,K}L^{-p},
\]
uniformly in $\beta$, proving \eqref{A:eq:B.65}. We record the isometric versions needed below. Let $\mathcal R^C_{n,\beta}$ be the embedded resolvent of the principal block, and
let $\widehat{\mathcal J}_n$ be the isometry of the full matrix into $L^2(\mathbb R_+)\oplus\ell^2(\mathbb N)$ constructed before Theorem~\ref{A:res:20.6}, now using $I_n^*$. The
deleted coordinates are sent, in order, to the first coordinate vectors in the second summand. Direct multiplication gives
\[
 \widehat{\mathcal J}_n P_{I_n^*}^*(A_C-z)^{-1}P_{I_n^*}
              \widehat{\mathcal J}_n^*
       =\mathcal R^C_{n,\beta}(z)\oplus0.
\]
Isometric invariance of the Hilbert--Schmidt norm shows that the squared difference between this operator and the full embedded resolvent equals the actual boundary error, before
taking the supremum in $z$. Let $S_n(z)$ denote the compressed full resolvent and $\widehat S_n(z)$ the full resolvent embedded by $\widehat{\mathcal J}_n$. The ideal inequality
for compression and the preceding identity give
\[
 \sup_K\|\mathcal R^C_{n,\beta}-S_n\|_{\mathfrak S_2}^2
 \le\mathcal C_{\partial,n,\beta},\qquad
 \sup_K\|\widehat S_n-(\mathcal R^C_{n,\beta}\oplus0)\|_{\mathfrak S_2}^2
 \le\mathcal C_{\partial,n,\beta}.
\]
Consequently,
\[
 \sup_K\|\mathcal R^C_{n,\beta}-R_\beta\|_{\mathfrak S_2}^2
 \le2D_{n,\beta}^V+2\mathcal C_{\partial,n,\beta},
\]
\[
 \sup_K\|\widehat S_n-(R_\beta\oplus0)\|_{\mathfrak S_2}^2
 \le4D_{n,\beta}^V+6\mathcal C_{\partial,n,\beta}.
\]
Taking $L^p$ norms, \eqref{A:eq:B.64} and Lemma~\ref{A:res:B.10B} bound both right sides by $C_ph$. The latter version includes every finite eigenvalue.

For the almost sure estimate, fix $\beta$ and apply the kernels $\overline Q_n$ from Lemma~\ref{A:res:B.8} at parameters $(\beta,\beta,\ldots)$, using the same two-sided Brownian
path for all $n$. The complements of $\mathcal E_n^{\mathrm{KRV}}$, $\mathcal W_n$, and the boundary events have summable probabilities, as proved in
Subsection~\ref{A:sec:quadratic-form-bounds} and \eqref{A:eq:boundary-good-event}. Borel--Cantelli, without independence between sizes, gives probability one to
\[
 \Omega_{V,\beta}
 =\Omega_\partial\cap\Omega_{\mathrm{KRV}}
                         \cap\liminf_n\mathcal W_n.
\]
On $\Omega_{V,\beta}$, take the maximum of the finite indices supplied by the three coefficient events and the comparisons in Theorem~\ref{A:res:17.6}. There are
$N_{\omega,V,\beta,K}<\infty$ and $C_{\omega,V,\beta,K}<\infty$ such that $D_{n,\beta}^V(K)\le C_{\omega,V,\beta,K}h$ for every $n\ge N_{\omega,V,\beta,K}$. The same coefficient
estimates apply after deletion of the final $s_*$ coordinates.
\end{proof}

The next proposition gives the lower bound uniformly on the same parameter interval. Let $m_n\asymp n$ be deterministic.

\begin{proposition}[Lower bound]\label{A:res:B.12}
For $B\in\Omega_B$ there is $c_{B,I_\beta,K}>0$ such that
\[
 \inf_{\substack{\beta\in I_\beta,\ z\in K\\\operatorname{rank}A\le m_n}}
       \|A-R_\beta(B,z)\|_{\mathfrak S_2}^2
 \ge c_{B,I_\beta,K}m_n^{-1/3}\asymp c_{B,I_\beta,K}h.
\]
\end{proposition}

\begin{proof}
Fix $B\in\Omega_B$. The common form comparison and the eigenvalue bounds for the deterministic operator $-\partial_x^2+1+x$ give constants $c>0$ and $C_B<\infty$ such that
\[
 c j^{2/3}\le\Lambda_j(\beta)+\kappa_B\le C_Bj^{2/3},
       \qquad j\ge1,\quad\beta\in I_\beta.
\]
Let $Z_K=\sup_K|z|$ and $D_B=C_B+\kappa_B+Z_K$. Since $\Lambda_j+\kappa_B>0$, for $1\le j\le k$,
\[
 |\Lambda_j(\beta)-z|
 \le\Lambda_j(\beta)+\kappa_B+\kappa_B+|z|
 \le D_B k^{2/3}.
\]
The normal compact operator $R_\beta(z)$ therefore has at least $k$ singular values of magnitude at least $D_B^{-1}k^{-2/3}$: its eigenvectors with indices $1,\ldots,k$ already
supply $k$ mutually orthogonal singular vectors with this property. Consequently its decreasing singular values satisfy
\[
 s_k(R_\beta(z))\ge D_B^{-1}k^{-2/3},
 \qquad \beta\in I_\beta,\quad z\in K.
\]
The singular values are ordered by magnitude, independently of the increasing energy labels.

Let $T$ be any bounded operator of rank $r\le m_n$, and let $Q$ project onto its range. It is finite rank and Hilbert--Schmidt. Since $(I-Q)T=0$, orthogonal decomposition in the
range gives
\[
 \begin{split}
 \|T-R_\beta(z)\|_{\mathfrak S_2}^2
 &\ge\|(I-Q)R_\beta(z)\|_{\mathfrak S_2}^2\\
 &=\|R_\beta(z)\|_{\mathfrak S_2}^2
       -\operatorname{Tr}\bigl(QR_\beta(z)R_\beta(z)^*\bigr).
 \end{split}
\]
For a positive compact operator $C$ with decreasing eigenvalues $c_j$ and an orthogonal projection of rank $r$, expansion in an eigenbasis gives
\[
 \operatorname{Tr}(QC)=\sum_j c_j\|Qe_j\|^2\le\sum_{j=1}^r c_j.
\]
For $r=0$ the trace is zero. For $r\ge1$, the coefficients belong to $[0,1]$ and have total sum at most $r$; moving any weight from an index greater than $r$ to an unused part of
an index at most $r$ cannot decrease the sum. Equivalently, subtract $c_r$ on the first $r$ indices and use $c_j\le c_r$ on the rest. Apply this inequality to
$C=R_\beta(z)R_\beta(z)^*$. It follows that
\[
 \|T-R_\beta(z)\|_{\mathfrak S_2}^2
 \ge\sum_{k>r}s_k(R_\beta(z))^2
 \ge D_B^{-2}\sum_{k>m_n}k^{-4/3}.
\]
Since $t^{-4/3}$ decreases,
\[
 \sum_{k>m_n}k^{-4/3}\ge\int_{m_n+1}^{\infty}t^{-4/3}\,\mathrm dt
       =3(m_n+1)^{-1/3}.
\]
All constants are independent of $\beta,z,T$. Taking their infimum proves the claimed bound. Finally $m_n\asymp n$ and $h\asymp n^{-1/3}$ turn $(m_n+1)^{-1/3}$ into the stated
order $h$, after a deterministic adjustment of the comparison constants.
\end{proof}

\subsubsection{Fixed parameters and converging parameter sequences}\label{A:sec:fixed-parameters-and-converging-parameter-sequences}
We first prove a uniform lower tail bound for the smallest eigenvalue. Let $A_{n,\beta}$ denote the full Gaussian matrix, the uniformly reduced general Jacobi matrix, or the full
general Jacobi matrix. For the reduced model fix $n_0$ large enough that the retained interval is nonempty. Nonreal resolvent convergence by itself does not exclude a low
eigenvalue escaping to $-\infty$, as the example in Subsection~\ref{A:sec:spectral-consequences} shows. We therefore establish lower tightness before using the resolvent comparison
to identify increasing eigenvalue labels.

\begin{lemma}\label{A:res:lower-tightness}
For each of these models,
\[
 \lim_{R\to\infty}\sup_{n\ge n_0,\,\beta\in I_\beta}
      \mathbb P(\lambda_1(A_{n,\beta})<-R)=0.
 \eqtag{lower-tightness}
\]
The same conclusion holds for $\Lambda_1(\beta)$, uniformly in $\beta\in I_\beta$.
\end{lemma}

More precisely, for $Q>1$ and $t\ge1$ there are deterministic $R_Q(t)=C_Q(1+t)^4+2$ and $N_Q(t)$ such that
\[
 \sup_{\beta\in I_\beta}
 \mathbb P(\lambda_1(A_{n,\beta})<-R_Q(t))
 \le C\mathrm e^{-ct^2}+C_Qn^{-Q},\qquad n\ge N_Q(t).
 \eqtag{lower-finite-tail}
\]

\begin{proof}
The Brownian marginal is exact, so
\[
 \sup_{n,\beta}\overline\Pi_{n,\beta}(M_B>1+t)
 =\mathbb W(M_B>1+t)\le C\mathrm e^{-ct^2}.
 \eqtag{lower-primitive-tail}
\]
On $E_{n,\beta,Q}$, Lemma~\ref{A:res:B.10A} gives
\[
 \mathfrak a_{\mathrm{red}}[u]
 \ge c_Q\ell_h[u]-C_QM_B^4\|u\|^2.
 \eqtag{lower-core-form}
\]
Here the form is that of the full matrix in the Gaussian case. On $E_{n,\beta,Q}\cap\{M_B\le1+t\}$ this implies $A_{n,\beta}+R_Q(t)\ge I$ for the Gaussian and reduced models. For
the full Jacobi matrix, also impose the events $\mathcal L_{n,Q}\cap\mathcal V_{n,Q}$ of Lemma~\ref{A:res:lower-transfer}. Their failure probability is $C_Qn^{-Q}$, uniformly on
$I_\beta$. Apply that lemma at $a=-R_Q(t)$: its hypothesis $A_C-a\ge I$ has just been proved, and its deterministic starting size is allowed to depend on $R_Q(t)$. It yields
$\mathcal A_{n,V,\beta}-a>0$. The union bound proves \eqref{A:eq:lower-finite-tail} in all three cases. For the continuum, use $\mathcal H_\beta\ge-\kappa_B\ge-P_B$ and the same
Brownian tail.

The threshold $N_Q(t)$ depends on the proposed lower bound, so \eqref{A:eq:lower-finite-tail} is not yet a bound uniform over all sizes. We treat the finitely many smaller matrices
only after $t$ has been fixed. To do so, fix $\varepsilon>0$. Choose $t$ with $C e^{-ct^2}<\varepsilon$, and then $N\ge N_Q(t)$ with $C_QN^{-Q}<\varepsilon$. For $n\ge N$ and
$R\ge R_Q(t)$,
\[
 \sup_\beta\mathbb P(\lambda_1(A_{n,\beta})<-R)<2\varepsilon.
\]
For each of the finitely many $n_0\le n<N$, the integrable upper bound for the density \eqref{A:eq:B.44}, multiplied by any polynomial in the coefficients, is integrable uniformly
for $\beta\in I_\beta$. To bound the normalizer below, integrate the unnormalized density on $|a_i|\le1$, $1\le b_i\le2$. At this fixed size its logarithm is bounded uniformly on
the product of that box and $I_\beta$; the box has positive volume. The normalizer is consequently at least a positive constant depending on $n,I_\beta,V$. Thus every coefficient
moment is uniformly bounded at this size. Since
\[
 \|T(a,b)\|_{\mathrm{op}}
 \le\max_i|a_i|+2\max_i b_i,
\]
all moments of the fixed size scaled matrix norm are uniformly bounded as well. For a principal block its norm is no larger than that of the full matrix. Therefore
\[
 \sup_\beta\mathbb P(\lambda_1(A_{n,\beta})<-R)
 \le R^{-p}\sup_\beta\mathbb E\|A_{n,\beta}\|^p\longrightarrow0.
\]
The Gaussian and chi moment formulas give the same conclusion in the Gaussian representation. Take the maximum over these finitely many sizes and then let $\varepsilon\downarrow0$.
\end{proof}

We apply Lemma~\ref{A:res:ordered-transform} on events with a common lower spectral bound and an upper bound for the first $m$ limiting eigenvalues. On one probability space let
$A_n$ be finite self-adjoint matrices whose dimensions tend to infinity, and let $H$ be an almost surely lower semibounded operator with compact resolvent on an
infinite-dimensional space $\mathcal H_0$. Let $J_n$ and $U$ be isometries of the finite spaces and $\mathcal H_0$ into one Hilbert space $\mathcal H$. Set
\[
 S_n=J_n(A_n-\mathrm i)^{-1}J_n^*,\qquad
 S=U(H-\mathrm i)^{-1}U^*.
\]
Assume $\|S_n-S\|\to0$ in probability and $\lim_{R\to\infty}\limsup_n\mathbb P(\lambda_1(A_n)<-R)=0$.

\par\Needspace{6\baselineskip}
\begin{lemma}\label{A:res:ordered-probability}
For every fixed $m$, the increasing eigenvalues satisfy
\[
 \max_{k\le m}|\lambda_{k,n}-\Lambda_k|\longrightarrow0
 \quad\text{in probability}.
 \eqtag{ordered-probability}
\]
\end{lemma}

\begin{proof}
Write $\varepsilon_n=\|S_n-S\|$, fix $R,U_0\ge1$, and put $a=-R-1$ and $b_R=1+|\mathrm i-a|$. On the event
\[
 \mathcal A_{n,R,U_0}
 =\{\lambda_1(A_n)\ge-R,\ \Lambda_1\ge-R,\ \Lambda_m\le U_0\},
\]
both operators minus $a$ are at least $I$. Lemma~\ref{A:res:ordered-transform}, with $\delta=1$, gives
\[
 \left|\frac1{\lambda_{k,n}-a}-\frac1{\Lambda_k-a}\right|
 \le b_R^2\varepsilon_n,\qquad k\le\dim A_n.
\]
If $b_R^2\varepsilon_n\le[2(U_0-a)]^{-1}$, its reciprocal calculation gives
\[
 \max_{k\le m}|\lambda_{k,n}-\Lambda_k|
 \le2b_R^2(U_0-a)^2\varepsilon_n.
\]
The lemma includes the zero complements of both isometric ranges. For $\eta>0$, let
\[
 d_{R,U_0,\eta}=
 \min\left\{\frac1{2b_R^2(U_0-a)},
             \frac{\eta}{2b_R^2(U_0-a)^2}\right\}>0.
\]
For $n$ with $\dim A_n\ge m$ we have
\[
 \begin{split}
 \mathbb P\left(\max_{k\le m}|\lambda_{k,n}-\Lambda_k|>\eta\right)
 \le{}&\mathbb P(\lambda_1(A_n)<-R)+\mathbb P(\Lambda_1<-R)\\
 &+\mathbb P(\Lambda_m>U_0)
                  +\mathbb P(\varepsilon_n>d_{R,U_0,\eta}).
 \end{split}
 \eqtag{ordered-truncation}
\]
First take $\limsup_n$, which removes the last term. Then let $R\to\infty$ using the assumed lower tightness and lower semiboundedness of $H$, and let $U_0\to\infty$ since
$\Lambda_m$ is finite almost surely. The resulting bound is zero.
\end{proof}

Use the space $\mathcal X_K$ defined in Subsection~\ref{A:sec:measurable-couplings-of-the-finite-models}. Denote the finite resolvent in Theorems~\ref{A:res:B.9}
and~\ref{A:res:B.11} by $\mathcal R_{n,\theta}$. The corresponding moment bounds imply
\[
 \sup_{\theta\in I_\beta}
 \|\mathcal R_{n,\theta}-R_\theta\|_{L^p(\mathcal X_K)}
 \le C_{p,V,I_\beta,K}h^{1/2}.
 \eqtag{B.73}
\]
For $p\ge2$ use $\|\|T\|_{\mathcal X_K}\|_p =\|\|T\|_{\mathcal X_K}^2\|_{p/2}^{1/2}$; for $p<2$ use the $L^1$ bound for the squared norm. All variables are Borel by
Lemma~\ref{A:res:B.6}, and their $p$th moments are finite by the finite-rank bound and \eqref{A:eq:continuum-coarse}. This conversion is also where the squared error rate $h$
becomes the norm rate $h^{1/2}$ in Wasserstein distance. The Banach space $\mathcal X_K$ is separable: $K$ is compact metric and $\mathfrak S_2(L^2(\mathbb R_+))$ is separable.

For two parameters the measure
\[
 \Pi_{n,\beta,\gamma}(\mathrm dB,\mathrm dx,\mathrm dy)
 =\mathbb W(\mathrm dB)Q_n(\beta,B,\mathrm dx)
                         Q_n(\gamma,B,\mathrm dy)
\]
is a probability measure. Integrating out $y$ or $x$ recovers $\Pi_{n,\beta}$ or $\Pi_{n,\gamma}$, respectively. For probability measures on $\mathcal X_K$ with finite $p$th
moments, $W_p^{\mathcal X_K}$ is the infimum of $(\mathbb E\|X-Y\|_{\mathcal X_K}^p)^{1/p}$ over their couplings. The two sampled matrices need not approach one another when
$\gamma\to\beta$ at a fixed size: they are sampled independently conditional on $B$. We estimate their difference through the limiting operators driven by that same path. Using the
displayed joint law, Minkowski and Proposition~\ref{A:res:B.5} give
\[
 \begin{aligned}
W_p^{\mathcal X_K}(\operatorname{Law}(\mathcal R_{n,\beta}), \operatorname{Law}(\mathcal R_{n,\gamma}))& \le\|\mathcal R_{n,\beta}-R_\beta\|_{L^p(\mathcal X_K)} +\|R_\beta-R_\gamma\|_{L^p(\mathcal X_K)}\\
&+\|R_\gamma-\mathcal R_{n,\gamma}\|_{L^p(\mathcal X_K)} \le C_p\{h^{1/2}+|\beta^{-1/2}-\gamma^{-1/2}|\}.
\end{aligned}
\]

For the general model, distinguish the isometrically embedded principal block resolvent
\[
 \mathcal R^C_{n,\theta}(z)
 =\iota_h^{(M_n^*)}U_{r_n}(A_C-z)^{-1}
                    U_{r_n}^*(\iota_h^{(M_n^*)})^*
\]
from the compressed full resolvent $\mathcal R_{n,\theta}$. The boundary error gives
\[
 \|\mathcal R^C_{n,\theta}-\mathcal R_{n,\theta}\|_{\mathcal X_K}
 \le\mathcal C_{\partial,n,\theta}^{1/2}.
 \eqtag{core-spectral-transform}
\]
Hence Lemma~\ref{A:res:B.10B} transfers \eqref{A:eq:B.73} to $\mathcal R^C$. In the Gaussian case set $\mathcal R^C=\mathcal R$.

For a deterministic sequence $\beta_n\in I_\beta$, the countable coupling is
\[
 \overline\Pi^{(\beta_\bullet)}
 =\mathbb W_\pm(\mathrm dB^+,\mathrm dB^-)
   \prod_{n\ge1}\overline Q_n
     (\beta_n,B^+,B^-,\mathrm dx_n,\mathrm d\mathsf a_n).
 \eqtag{countable-lift}
\]
To specify this expression, integrate a cylinder function depending on the path and the first $m$ coordinates against the product of the first $m$ kernels. These probability laws
are consistent, since the next kernel has total mass one. Their countable extension exists on the product of the standard Borel spaces and has $(B^+,B^-,x_n,\mathsf a_n)$ marginal
$\overline\Pi_{n,\beta_n}$. In particular every preceding moment estimate applies under this single coupling. The sequence is chosen deterministically; these marginal estimates are
not being applied to a parameter selected from the observed Brownian path. Finitely many parameters outside a chosen compact interval can be treated using their own kernels.

Suppose now $\beta_n\to\beta>0$, and choose $I_\beta$ containing $\beta$ and all sufficiently large $\beta_n$.

\par\Needspace{6\baselineskip}
\begin{theorem}\label{A:res:B.13}
For every $p\ge1$,
\[
 W_p^{\mathcal X_K}(\operatorname{Law}(\mathcal R_{n,\beta_n}),
                       \operatorname{Law}(R_\beta))
 \le C_{p,V,I_\beta,K}
       \{h^{1/2}+|\beta_n^{-1/2}-\beta^{-1/2}|\}.
 \eqtag{B.78}
\]
\end{theorem}

Thus $\mathcal R_{n,\beta_n}\Rightarrow R_\beta$ in $\mathcal X_K$. Moreover, on \eqref{A:eq:countable-lift} the following convergences hold in probability, and consequently in
distribution:
\[
 (\lambda_{1,n}^{(\beta_n)},\ldots,\lambda_{m,n}^{(\beta_n)})
 \Rightarrow(\Lambda_1(\beta),\ldots,\Lambda_m(\beta)),
 \qquad m\ge1,
 \eqtag{B.79}
\]
\[
 \sum_j\delta_{\lambda_{j,n}^{(\beta_n)}}
 \Rightarrow\sum_{k\ge1}\delta_{\Lambda_k(\beta)}
 \quad\text{in the vague topology}.
 \eqtag{B.80}
\]
For the general model these statements hold for both the reduced and the full matrix spectra.

\begin{proof}
On \eqref{A:eq:countable-lift}, split the difference through the operator at the finite parameter:
\[
 \mathcal R_{n,\beta_n}-R_\beta
 =(\mathcal R_{n,\beta_n}-R_{\beta_n})+(R_{\beta_n}-R_\beta).
\]
The first term is the matrix approximation error; the second changes only the amplitude of the noise on the same path. Minkowski and \eqref{A:eq:B.73},
\eqref{A:eq:continuum-polynomial-Lipschitz} give
\[
 \|\mathcal R_{n,\beta_n}-R_\beta\|_{L^p(\mathcal X_K)}
 \le C_p\{h^{1/2}+|\beta_n^{-1/2}-\beta^{-1/2}|\}.
\]
By the definition of Wasserstein distance as the infimum over joint laws, the law on \eqref{A:eq:countable-lift} gives
\[
 W_p^{\mathcal X_K}(\operatorname{Law}(\mathcal R_{n,\beta_n}),
                    \operatorname{Law}(R_\beta))
 \le \bigl(\mathbb E\|\mathcal R_{n,\beta_n}-R_\beta\|_{\mathcal X_K}^p\bigr)^{1/p}.
\]
This proves \eqref{A:eq:B.78}. For the principal block, \eqref{A:eq:core-spectral-transform} and the continuum estimate give
\[
 \|\mathcal R^C_{n,\beta_n}-R_\beta\|_{L^p(\mathcal X_K)}
 \le\|\mathcal R^C_{n,\beta_n}-R_{\beta_n}\|_{L^p(\mathcal X_K)}
       +\|R_{\beta_n}-R_\beta\|_{L^p(\mathcal X_K)}
 \le C_p\{h^{1/2}+|\beta_n^{-1/2}-\beta^{-1/2}|\}.
\] Apply it on
$K\cup\{\mathrm i\}$. Lemmas~\ref{A:res:lower-tightness} and \ref{A:res:ordered-probability} then prove the ordered convergence for the Gaussian and reduced models on this
coupling.

For the full general matrix use $\widehat{\mathcal J}_n$ from the proof of Theorem~\ref{A:res:B.11} and put
$\widehat S_{n,\beta_n}= \widehat{\mathcal J}_n(\mathcal A_{n,V,\beta_n}-\mathrm i)^{-1} \widehat{\mathcal J}_n^*$. The isometry acts on every original coordinate. Hence
\[
 \|\widehat S_{n,\beta_n}-(R_\beta(\mathrm i)\oplus0)\| \le\mathcal
 C_{\partial,n,\beta_n}(\{\mathrm i\})^{1/2} +\|\mathcal R^C_{n,\beta_n}(\mathrm
 i)-R_\beta(\mathrm i)\| \longrightarrow0\quad\text{in probability}.
 \eqtag{full-spectral-transform}
\]
The enlarged isometry represents every finite eigenvalue. Its auxiliary zero summand adds only zero eigenvalues to the embedded resolvent. Use tightness of the lower spectral bound
for the full matrix in Lemma~\ref{A:res:ordered-probability}. This proves \eqref{A:eq:B.79} for the increasing eigenvalue labels of the full matrix.

For vague convergence, fix $f\in C_c(\mathbb R)$ supported in $[-L,L]$ and let $\varepsilon>0$. Since $\Lambda_m(\beta)\to\infty$ almost surely, choose $m$ with
$\mathbb P(\Lambda_m(\beta)\le L+1)<\varepsilon$. On
\[
 \{\Lambda_m(\beta)>L+1\}
 \cap\{\max_{j\le m}|\lambda_{j,n}^{(\beta_n)}-\Lambda_j(\beta)|<1\},
\]
all labels $j\ge m$ lie above $L$ in both spectra. If $\omega_f$ is the modulus of continuity of $f$ on $\mathbb R$, the difference of its two integrals is consequently at most
\[
 (m-1)\omega_f\!
   \left(\max_{j<m}|\lambda_{j,n}^{(\beta_n)}-\Lambda_j(\beta)|\right).
\]
This tends to zero in probability; the excluded event has limiting upper probability at most $\varepsilon$. Let
$\varepsilon\downarrow0$.

Let $(f_r)$ and $d$ be the countable family of functions and the metric constructed in the proof of Theorem~\ref{A:res:21.5}. For $\mu_n=\sum_j\delta_{\lambda_{j,n}^{(\beta_n)}}$
and $\mu=\sum_k\delta_{\Lambda_k(\beta)}$, the preceding argument gives $\langle\mu_n-\mu,f_r\rangle\to0$ in probability for every $r$. Given $\eta>0$, choose $r_0$ with
$2^{-r_0}<\eta/2$. Then
\[
 \mathbb P\{d(\mu_n,\mu)>\eta\}
 \le\sum_{r=1}^{r_0}
    \mathbb P\{|\langle\mu_n-\mu,f_r\rangle|>\eta/2\}
 \longrightarrow0.
\]
The deterministic tail of the defining metric is at most $2^{-r_0}$, which proves this bound. The measures are Borel by the same countable family of integrals. This proves
\eqref{A:eq:B.80}.

On this fixed countable coupling the estimates for every finite moment give a stronger, almost sure conclusion. For any $\varepsilon>0$, apply the norm moment bound at an integer
$p$ with $p\varepsilon>1$. Since $h_n\asymp n^{-1/3}$, Markov's inequality gives

\[
 \mathbb P\!\left(
 \|\mathcal R_{n,\beta_n}-R_{\beta_n}\|_{\mathcal X_K}
       >n^{-1/6+\varepsilon}\right)
 \le C_{p,I_\beta,K}n^{-p\varepsilon}.
\]
The same estimate holds for the principal and enlarged isometric resolvents by the preceding boundary comparison. The probability series is finite. Hence, for each fixed
$\varepsilon,K$, there is $N_{\omega,\varepsilon,K}<\infty$ such that all three errors are at most $n^{-1/6+\varepsilon}$ for $n\ge N_{\omega,\varepsilon,K}$. Intersect these
events of probability one over positive rational $\varepsilon$ and a countable compact exhaustion of $\mathbb C\setminus\mathbb R$. The pathwise continuum estimate then gives
\[
 \|\mathcal R_{n,\beta_n}-R_\beta\|_{\mathcal X_K}
 =O_\omega\!\left(n^{-1/6+\varepsilon}
                  +|\beta_n^{-1/2}-\beta^{-1/2}|\right)
 \quad(\varepsilon>0).
\]
Every compact nonreal $K$ is covered by the exhaustion. To justify ordered labels almost surely, fix $Q>1$ in the proof of Lemma~\ref{A:res:lower-tightness}. For each integer
$t\ge1$, that proof also bounds $\mathbb P(\{\lambda_1(A_{n,\beta_n})<-R_Q(t)\} \cap\{M_B\le1+t\})$ by $C_Qn^{-Q}$ for $n\ge N_Q(t)$. For each integer $t$, summability gives an
index $N_{\omega,t}$ such that $\lambda_1(A_{n,\beta_n})\ge-R_Q(t)$ for $n\ge N_{\omega,t}$ on $\{M_B\le1+t\}$. Intersect over the integer values of $t$ and choose one with
$M_B\le1+t$ on the chosen path. This gives a lower bound for all $n\ge N_{\omega,t}$. The continuum operator is bounded below by $-P_B$. Its first $m$ eigenvalues are finite. Apply
the deterministic shifted inverse estimate in Lemma~\ref{A:res:ordered-transform} on each such path to the isometric resolvents at $\mathrm i$. It gives
\[
 \max_{k\le m}|\lambda_{k,n}^{(\beta_n)}-\Lambda_k(\beta)|
 =O_\omega\!\left(n^{-1/6+\varepsilon}
                  +|\beta_n^{-1/2}-\beta^{-1/2}|\right).
\]
Taking $0<\varepsilon<1/6$ proves almost sure ordered and vague convergence by the preceding finite sum argument. This conclusion is for the one specified deterministic parameter
sequence and its countable coupling. It does not give an event for every sequence simultaneously, a bound for parameters selected from the observed path, or the fixed parameter
rate without the $n^\varepsilon$ loss.
\end{proof}

Use the bounded--Lipschitz distance
\[
 d_{\mathrm{BL}}(\mu,\nu)=
 \sup_{\|f\|_\infty+\operatorname{Lip}(f)\le1}
        \left|\int f\,\mathrm d\mu-\int f\,\mathrm d\nu\right|.
\]
We shall use the implication from weak convergence to convergence in this distance. Suppose $\mu_n\Rightarrow\mu$ and fix $\varepsilon>0$. Tightness gives $R$ such that
$\sup_n\mu_n([-R,R]^c)+\mu([-R,R]^c)\le\varepsilon$. The functions in the defining supremum are uniformly bounded and $1$-Lipschitz. On $[-R,R]$ they have a finite
$\varepsilon$-net $f_1,\ldots,f_m$ chosen from the same class. For any such $f$, choose $f_j$ with $\sup_{[-R,R]}|f-f_j|\le\varepsilon$. Then
\[
 \left|\int f\,\mathrm d(\mu_n-\mu)\right|
 \le \max_{j\le m}\left|\int f_j\,\mathrm d(\mu_n-\mu)\right|
       +2\varepsilon+2\bigl(\mu_n([-R,R]^c)+\mu([-R,R]^c)\bigr).
\]
The maximum tends to zero by weak convergence. Taking the supremum in $f$ and then letting $\varepsilon\downarrow0$ proves the implication.

\par\Needspace{6\baselineskip}
\begin{proposition}\label{A:res:B.14}
For either finite model and each fixed $k$,
\[
 \sup_{\beta\in I_\beta}
 d_{\mathrm{BL}}(\operatorname{Law}(\lambda_{k,n}^{(\beta)}),
                  \operatorname{Law}(\Lambda_k(\beta)))
 \longrightarrow0.
\]
\end{proposition}

\begin{proof}
Suppose uniform convergence fails. Then choose $\varepsilon>0$, $n_j\to\infty$, and $\beta_j\in I_\beta$ such that the displayed distance exceeds $\varepsilon$. Compactness gives a
subsequence with $\beta_j\to\beta_*$. Theorem~\ref{A:res:B.13}, applied along these sizes, gives $\lambda_{k,n_j}^{(\beta_j)}\Rightarrow\Lambda_k(\beta_*)$. Put
$\mu_j=\operatorname{Law}(\lambda_{k,n_j}^{(\beta_j)})$, $\nu_j=\operatorname{Law}(\Lambda_k(\beta_j))$, and $\nu_*=\operatorname{Law}(\Lambda_k(\beta_*))$. The bounded--Lipschitz
approximation proved before the proposition gives $d_{\mathrm{BL}}(\mu_j,\nu_*)\to0$. For every test function in the definition of $d_{\mathrm{BL}}$, its Lipschitz constant is at
most one. Applying Proposition~\ref{A:res:B.5} to the operators driven by the same Brownian path gives
\[
 d_{\mathrm{BL}}(\nu_j,\nu_*)
 \le\mathbb E|\Lambda_k(\beta_j)-\Lambda_k(\beta_*)|
 \le C_{k,I_\beta}|c_{\beta_j}-c_{\beta_*}|\longrightarrow0.
\]
Consequently $d_{\mathrm{BL}}(\mu_j,\nu_j) \le d_{\mathrm{BL}}(\mu_j,\nu_*)+d_{\mathrm{BL}}(\nu_*,\nu_j) \to0$, contradicting its lower bound $\varepsilon$. The argument applies to
the full and reduced general spectra separately, since Theorem~\ref{A:res:B.13} covers both.
\end{proof}

We next pass from bounded--Lipschitz distance to uniform convergence of distribution functions. For $\rho>0$, put
\[
 \omega_k(\rho)=\sup_{\beta\in I_\beta,\ s\in\mathbb R}
             \mathbb P(|\Lambda_k(\beta)-s|\le\rho).
\]
We first prove that $\omega_k(\rho)\to0$ as $\rho\downarrow0$. The form comparison gives
\[
 -\kappa_B\le\Lambda_k(\beta)\le CP_B\nu_k-\kappa_B,
\]
so the moment bounds for $P_B$ imply $\sup_\beta\mathbb E|\Lambda_k(\beta)|<\infty$. Markov's inequality therefore makes these laws uniformly tight.

Suppose the proposed limit is false. Choose $\varepsilon>0$, $\rho_j\downarrow0$, and $\beta_j\in I_\beta,s_j\in\mathbb R$ such that
\[
 \mathbb P(|\Lambda_k(\beta_j)-s_j|\le\rho_j)\ge\varepsilon.
\]
Uniform tightness gives $R$ with $\sup_\beta\mathbb P(|\Lambda_k(\beta)|>R)<\varepsilon/2$. For $\rho_j<1$ and $|s_j|>R+1$, the event on the left is contained in
$\{|\Lambda_k(\beta_j)|>R\}$, a contradiction. Thus the centers $s_j$ are bounded after discarding finitely many indices. Compactness permits a subsequence with $\beta_j\to\beta_*$
and $s_j\to s_*$. On the common Brownian space, Proposition~\ref{A:res:B.4} gives $\Lambda_k(\beta_j)\to\Lambda_k(\beta_*)$ almost surely. If the latter value differs from $s_*$,
choose one third of their distance; convergence of the two sequences and $\rho_j\to0$ imply that $|\Lambda_k(\beta_j)-s_j|>\rho_j$ for every sufficiently large $j$ on that sample.
Hence
\[
 \limsup_j\mathbf1_{\{|\Lambda_k(\beta_j)-s_j|\le\rho_j\}}
 \le\mathbf1_{\{\Lambda_k(\beta_*)=s_*\}}\quad\text{almost surely}.
\]
Reverse Fatou, applicable to indicators bounded by one, and Theorem~\ref{A:res:20.1} now give
\[
 \varepsilon\le\limsup_j\mathbb P(|\Lambda_k(\beta_j)-s_j|\le\rho_j)
 \le\mathbb P(\Lambda_k(\beta_*)=s_*)=0.
\]
This proves that the probabilities of intervals of vanishing length tend to zero uniformly in their centers and in the parameter.

For a fixed $s$ and $\rho>0$, define
\[
 u_\rho(x)=\bigl(1-(x-s)_+/\rho\bigr)_+,\qquad
 l_\rho(x)=\bigl(1-(x-s+\rho)_+/\rho\bigr)_+.
\]
They satisfy $l_\rho\le\mathbf1_{(-\infty,s]}\le u_\rho$, their supremum norms are at most one, and their Lipschitz constants are at most $\rho^{-1}$. Thus division by
$1+\rho^{-1}$ makes them admissible in the definition of $d_{\mathrm{BL}}$. Write $\mu_{n,\beta}$ and $\mu_\beta$ for the finite and limiting eigenvalue laws, and let
\[
 b_n=\sup_{\beta\in I_\beta}d_{\mathrm{BL}}(\mu_{n,\beta},\mu_\beta).
\]
The upper comparison is
\[
 \begin{split}
 F_{\beta,n}^{(k)}(s)-F_\beta^{(k)}(s)
 &\le\int u_\rho\,\mathrm d(\mu_{n,\beta}-\mu_\beta)
           +\int(u_\rho-\mathbf1_{(-\infty,s]})\,\mathrm d\mu_\beta\\
 &\le(1+\rho^{-1})b_n+\omega_k(\rho).
 \end{split}
\]
Using $l_\rho$ instead proves the same upper bound for $F_\beta^{(k)}(s)-F_{\beta,n}^{(k)}(s)$. Therefore
\[
 \sup_{\beta\in I_\beta,\ s\in\mathbb R}
       |F_{\beta,n}^{(k)}(s)-F_\beta^{(k)}(s)|
 \le(1+\rho^{-1})b_n+\omega_k(\rho).
\]
Proposition~\ref{A:res:B.14} gives $b_n\to0$. First take the upper limit as $n\to\infty$ with $\rho$ fixed, and then let $\rho\downarrow0$. This proves uniform convergence of the
distribution functions. The argument supplies no numerical rate, since no numerical modulus for $\omega_k(\rho)$ has been established; it also supplies no moment bound for inverse
spectral gaps. The nonreal estimates use the deterministic distance $\delta_K$ from the real axis. At real parameters it is replaced by a spectral gap, which can be small and
depends on the path. We distinguish a pathwise comparison with a gap uniform over $I_\beta$ from convergence in probability at one fixed real parameter. Fix $B\in\Omega_B$ and a
compact real set $L$ such that
\[
 d_B:=\inf_{\theta\in I_\beta}
       \operatorname{dist}(L,\sigma(\mathcal H_\theta))>0.
 \eqtag{B.90}
\]

\begin{proposition}\label{A:res:B.15}
Under \eqref{A:eq:B.90},
\[
 \sup_{a\in L}\|(\mathcal H_\beta-a)^{-1}
                    -(\mathcal H_\gamma-a)^{-1}\|_{\mathfrak S_2}
 \le C_{B,I_\beta,L}|\beta^{-1/2}-\gamma^{-1/2}|.
\]
\end{proposition}

For a fixed deterministic $a\in\mathbb R$ and $\beta_n\to\beta>0$, the isometrically embedded inverses of the Gaussian matrices and the reduced general matrices converge in
probability in $\mathfrak S_2$ to $(\mathcal H_\beta-a)^{-1}$ on \eqref{A:eq:countable-lift}. The finite operators minus $a$ are invertible with probability tending to one. The
same invertibility and convergence statements hold for the full Jacobi matrix in the enlarged isometric space, with limiting inverse $(\mathcal H_\beta-a)^{-1}\oplus0$; compression
then gives the corresponding convergence in $L^2(\mathbb R_+)$. On a singular event the inverse can be assigned the value zero.

\begin{proof}
Set $S_\theta(a)=R_\theta(a+\mathrm i)$ and $F(S)=S(I+\mathrm iS)^{-1}$. Spectral calculus gives
\[
 \sup_{\theta\in I_\beta,a\in L}
    \|(I+\mathrm iS_\theta(a))^{-1}\|
 \le1+d_B^{-1}.
\]
Multiplication by the two denominators proves
\[
 F(S)-F(T)=(I+\mathrm iS)^{-1}(S-T)(I+\mathrm iT)^{-1}.
\]
Theorem~\ref{A:res:B.3}, used on $L+\mathrm i$, proves the proposition.

For convergence in probability first use $S=R_\beta(a+\mathrm i)$ and $S_n=\mathcal R^C_{n,\beta_n}(a+\mathrm i)$. Theorem~\ref{A:res:B.13} and \eqref{A:eq:core-spectral-transform}
give $\varepsilon_n=\|S_n-S\|_{\mathfrak S_2}\to0$ in probability. Atomlessness at this fixed $a$ implies that $Z=\|(I+\mathrm iS)^{-1}\|<\infty$ almost surely. This is Borel,
since inversion is continuous on the open set of invertible bounded operators; set $Z=\infty$ elsewhere. On
\[
 E_{n,M}=\{Z\le M\}\cap\{\varepsilon_n\le(2M)^{-1}\},
 \qquad M\ge1,
\]
the Neumann series gives $\|(I+\mathrm iS_n)^{-1}\|\le2M$. On the finite isometric range this is precisely invertibility of $A_n-a$; on its orthogonal complement the denominator is
the identity. Since $F(0)=0$, $F(S_n)$ is the embedded inverse and
\[
 \|F(S_n)-F(S)\|_{\mathfrak S_2}\le2M^2\varepsilon_n.
\]
With $F(S_n)=0$ on the singular event, for $\eta>0$ we obtain
\[
 \begin{aligned}
 \mathbb P(A_n-a\text{ is singular})
 &\le\mathbb P(Z>M)+\mathbb P(\varepsilon_n>(2M)^{-1}),\\
 \mathbb P(\|F(S_n)-F(S)\|_{\mathfrak S_2}>\eta)
 &\le\mathbb P(Z>M)+
   \mathbb P\!\left(\varepsilon_n>
            \min\{(2M)^{-1},\eta/(2M^2)\}\right).
 \end{aligned}
\]
First let $n\to\infty$ with $M$ fixed, then let $M\to\infty$. The remaining probability $\mathbb P(Z>M)$ tends to zero because $Z<\infty$ almost surely. This uses no moment bound
for the inverse spectral gap. In particular, we do not infer an $L^p$ estimate for the real inverse from the nonreal moment theorem.

For the full general matrix replace $S_n$ by $\widehat{\mathcal J}_n(\mathcal A_{n,V,\beta_n}-a-\mathrm i)^{-1} \widehat{\mathcal J}_n^*$ and $S$ by $R_\beta(a+\mathrm i)\oplus0$.
The isometric moment estimate in the proof of Theorem~\ref{A:res:B.11} gives the same convergence of $\varepsilon_n$ in probability, after the continuum parameter comparison. Put
$\widehat\varepsilon_n=\|S_n-S\|_{\mathfrak S_2}$ and $\widehat Z=\|(I+\mathrm iS)^{-1}\|=\max\{1,Z\}$ on this space. On $\{\widehat Z\le M,\ \widehat\varepsilon_n\le(2M)^{-1}\}$,
\[
 \|F(S_n)-F(S)\|_{\mathfrak S_2}\le2M^2\widehat\varepsilon_n.
\]
For every $\eta>0$, assigning zero on the singular event gives
\[
 \mathbb P(\|F(S_n)-F(S)\|_{\mathfrak S_2}>\eta)
 \le\mathbb P(\widehat Z>M)
   +\mathbb P\!\left(\widehat\varepsilon_n>
             \min\{(2M)^{-1},\eta/(2M^2)\}\right).
\]
Its iterated upper limit, first as $n\to\infty$ and then as $M\to\infty$, is zero. Projection onto its first summand is a contraction and sends the full embedded inverse to
$\mathcal J_{n,I_\beta}(\mathcal A_{n,V,\beta_n}-a)^{-1} \mathcal J_{n,I_\beta}^*$. This proves convergence of the compressed resolvents.
\end{proof}

\subsubsection{Further directions}\label{A:sec:further-directions}
Two quantitative questions remain for the coupling constructed here: whether the normalized error converges, and whether the same rate holds when every matrix coordinate is
embedded in its original order.

Use the projection $P_n$ in \eqref{A:eq:23.1}. Put $Q_n=I-P_n$, $R=R_\beta(z)$ and $\mathscr R_n=\mathscr R_n(z)$. The embedded finite resolvent satisfies
$\mathscr R_n=P_n\mathscr R_nP_n$. Therefore
\[
 \mathscr R_n-R=(\mathscr R_n-P_nRP_n)-P_nRQ_n-Q_nRP_n-Q_nRQ_n.
\]
These four terms occupy distinct blocks for the orthogonal decomposition $\operatorname{Ran}P_n\oplus\operatorname{Ran}Q_n$. They are orthogonal in the Hilbert space
$\mathfrak S_2$. This follows, for example, by summing their matrix entries in an orthonormal basis formed from bases of the two subspaces; no two blocks have a common pair of
indices. Consequently
\[
 \begin{split}
 \|\mathscr R_n-R\|_{\mathfrak S_2}^2
 ={}&\|\mathscr R_n-P_nRP_n\|_{\mathfrak S_2}^2
       +\|P_nRQ_n\|_{\mathfrak S_2}^2\\
 &+\|Q_nRP_n\|_{\mathfrak S_2}^2+\|Q_nRQ_n\|_{\mathfrak S_2}^2\\
 ={}&\|\mathscr R_n-P_nRP_n\|_{\mathfrak S_2}^2
       +\|R-P_nRP_n\|_{\mathfrak S_2}^2.
 \end{split}
\]
Both terms on the last line are nonnegative and at most the full squared error. Thus the established upper bound makes both $O_\omega(h_n)$ under the coupling constructed for each
fixed parameter, uniformly on the compact nonreal sets under consideration. Moreover, $\operatorname{rank}(P_nRP_n)\le\operatorname{rank}P_n$. The latter rank is the dimension of
the retained embedding and has order $n$. Applying Proposition~\ref{A:res:B.12} to this particular bounded finite-rank operator gives
\[
 \|R_\beta(z)-P_nR_\beta(z)P_n\|_{\mathfrak S_2}^2
       \ge c_{B,I_\beta,K}h_n.
\]
This lower bound concerns approximation of the full continuum resolvent. It gives no corresponding lower bound for $\|\mathscr R_n-P_nRP_n\|_{\mathfrak S_2}^2$, whose two arguments
already act inside the same finite-dimensional range. In particular, an upper bound of order $h_n$ for that term does not identify its leading coefficient. For the specified
coupling and a fixed nonreal $z$, consider $h_n^{-1}\|\mathscr R_n(z)-R_\beta(z)\|_{\mathfrak S_2}^2$. To prove convergence, it suffices to find the limit of the normalized
projection error and of the products in the following expansion. If $T_{1,n},\ldots,T_{m,n}$ are the comparison terms after projection by $P_n$, then
\[
 h_n^{-1}\left\|\sum_{j=1}^mT_{j,n}\right\|_{\mathfrak S_2}^2
 =\sum_{j=1}^m h_n^{-1}\|T_{j,n}\|_{\mathfrak S_2}^2
   +2\operatorname{Re}\sum_{j<k}
            h_n^{-1}\langle T_{j,n},T_{k,n}\rangle_{\mathfrak S_2}.
\]
The $\mathfrak S_2$ inner products between errors constructed from the same Brownian motion can contribute at order $h_n$.

The moment bounds already show that several terms vanish after the normalization needed for a leading order coefficient. Suppose an operator error $E_n$ satisfies
$\|\|E_n\|_{\mathfrak S_2}^2\|_p\le C_ps_n$, where $s_n=o(h_n)$, and another error $T_n$ satisfies $\|\|T_n\|_{\mathfrak S_2}^2\|_p\le C_ph_n$. Cauchy--Schwarz for the
Hilbert--Schmidt inner product and H\"older in probability give
\[
 \begin{aligned}
&\left\|h_n^{-1} (\|T_n+E_n\|_{\mathfrak S_2}^2-\|T_n\|_{\mathfrak S_2}^2) \right\|_p \le h_n^{-1} \{2\|\|T_n\|_{\mathfrak S_2}\|_{2p} \|\|E_n\|_{\mathfrak S_2}\|_{2p} +\|\|E_n\|_{\mathfrak S_2}^2\|_p\}\\
&\le C_p\{(s_n/h_n)^{1/2}+s_n/h_n\}\longrightarrow0.
\end{aligned}
\]
This also holds with a supremum over $z\in K$ and uniformly in $\beta\in I_\beta$ under the corresponding moment assumptions. A fixed sum of such errors obeys the same conclusion.

In the Gaussian model this applies to the scalar quantile error and to the terms containing $c_{k,n}=\overline b_{k,n}-\sqrt{1-k/n}$ on the regular comparison block, whose squared
resolvent errors have order $h_n^2$ or smaller. In the general model it applies to the conditional Gaussian transport, the diagonal covariance errors, the direct left and far
scalar replacements, and the full boundary error; the latter has the stronger bound $b_n=o(h_n)$ proved in Lemma~\ref{A:res:B.10B}. It also applies to the $h_n^2\ell_n^2$ scalar
recentering error. The first-order coefficient terms, the remaining deterministic comparison, and the continuum projection and discretization errors still have only $O(h_n)$
squared bounds. The limits of their normalized products remain to be identified. The moment calculation removes the listed smaller terms in $L^p$. An almost sure expansion requires
additional pathwise control of the remaining error.

A second question concerns embedding all coordinates of the general Jacobi matrix in their original order. This differs from the isometry into a space with an additional orthogonal
summand used for the spectral conclusions. Let $\mathscr R_n^{\mathrm{orig}}(z)= \iota_h^{(n)}(\mathcal A_{n,V,\beta}-z)^{-1}(\iota_h^{(n)})^*$, let $\delta_n=r_nh$, and define the
right translation isometry
\[
 (S_\delta f)(x)=\mathbf1_{\{x\ge\delta\}}f(x-\delta).
\]
Directly from the definitions of the two embeddings,
\[
 \iota_h^{(n)}P_{I_n^*}^*U_{r_n}^*
       =S_{\delta_n}\iota_h^{(M_n^*)}.
\]
The boundary comparison therefore gives
\[
 \mathscr R_n^{\mathrm{orig}}(z)-R_\beta(z)
 =E_{\partial,n}(z)
   +S_{\delta_n}(\mathcal R^C_{n,\beta}(z)-R_\beta(z))S_{\delta_n}^*
   +S_{\delta_n}R_\beta(z)S_{\delta_n}^*-R_\beta(z),
\]
where $\sup_K\|E_{\partial,n}\|_{\mathfrak S_2}^2 \le\mathcal C_{\partial,n,\beta}$. The first term is $o(h)$ in every fixed moment. At a fixed $\beta$, its squared norm is also
$O_\omega(b_n)=o_\omega(h)$, by \eqref{A:eq:B.54d1} with $A>1$ and Borel--Cantelli. The second term has the established $O(h)$ squared bound, both in moments and almost surely for
fixed $\beta$. To obtain the same order for the original coordinate embedding by this argument, it would suffice to prove
\[
 \sup_{z\in K}\|S_{\delta_n}R_\beta(z)S_{\delta_n}^*-R_\beta(z)
                         \|_{\mathfrak S_2}^2=O_\omega(h),
 \qquad \delta_n=r_nh=O(h\log n).
\]
Strong continuity of $S_\delta$ gives qualitative convergence in $\mathfrak S_2$ at each fixed $z$. To see the norm estimate, write $T_m=\sum_{j=1}^m f_jg_j^*$ with
$\|R_\beta(z)-T_m\|_{\mathfrak S_2}\to0$ as $m\to\infty$. For each summand,
\[
 \|S_\delta f_j(S_\delta g_j)^*-f_jg_j^*\|_{\mathfrak S_2}
 \le\|S_\delta f_j-f_j\|_2\|g_j\|_2
       +\|f_j\|_2\|S_\delta g_j-g_j\|_2.
\]
Here $S_\delta$ is an isometry and the Hilbert--Schmidt norm of $fg^*$ is $\|f\|_2\|g\|_2$. Therefore
\[
 \begin{split}
 \|S_\delta R_\beta(z)S_\delta^*-R_\beta(z)\|_{\mathfrak S_2}
 &\le2\|R_\beta(z)-T_m\|_{\mathfrak S_2}\\
 &\quad+\sum_{j=1}^m
   \bigl(\|S_\delta f_j-f_j\|_2\|g_j\|_2
        +\|f_j\|_2\|S_\delta g_j-g_j\|_2\bigr).
 \end{split}
\]
First let $\delta\downarrow0$ for fixed $m$, and then let $m\to\infty$. The displayed estimate for simultaneous translation of the two kernel variables at $\delta_n=O(h\log n)$
remains to be proved. The scalar recentering $rhI$ estimated earlier leaves the Brownian potential and Dirichlet boundary fixed, whereas translation moves both.

\par\bigskip
\Needspace{6\baselineskip} \noindent\textbf{AI disclosure.}\enspace AI tools were used for proofreading and language editing. Given the length of the manuscript, they were also
used to help condense the presentation of some routine proofs. The author takes full responsibility for the mathematical content and the final text.
\par

\clearpage
\phantomsection 

\end{document}